\documentclass[nosumlimits,twoside]{amsart}
\usepackage{amsfonts, amsmath, amssymb}
\usepackage[bookmarksnumbered,plainpages,driverfallback=dvipdfm,backref=page]{hyperref}
\usepackage{graphicx}
\usepackage{float}
\usepackage{srcltx}
\usepackage[all]{xy}
\usepackage{version}
\usepackage[T1]{fontenc}
\usepackage{soul, color}

\newtheorem{theorem}{\sc Theorem}[section]
\newtheorem{proposition}[theorem]{\sc Proposition}

\newtheorem{lemma}[theorem]{\sc Lemma}
\newtheorem{corollary}[theorem]{\sc Corollary}
\theoremstyle{definition}

\theoremstyle{remark}
\newtheorem{remark}[theorem]{\sc Remark}

\allowdisplaybreaks
\IfFileExists{tcilatex.tex}{\input{tcilatex}}{}

\makeatletter

\begin{document}
\author{Claudia Menini}
\address{University of Ferrara, Department of Mathematics, Via Machiavelli
30, Ferrara, I-44121, Italy}
\email{men@unife.it}
\urladdr{sites.google.com/a/unife.it/claudia-menini}
\title{New examples of separable cowreaths over Clifford algebras}
\author{Blas Torrecillas}
\address{Department of Mathematics, University of Almer\'{\i}a, Almer\'{\i}%
a, Spain}
\email{btorreci@ual.es}
\urladdr{https://w3.ual.es/~btorreci}
\thanks{This paper was written while the first author was a member of the
"National Group for Algebraic and Geometric Structures and their
Applications" (GNSAGA-INdAM) and was partially supported by MIUR within the
National Research Project PRIN 2017 up to September 31, 2024. The second
author was partially supported by PID2024-158993NB-I00 and Proyecto
Lanzadera.}

\begin{abstract}
This paper continues the research we developed in \cite{MT1} and \cite%
{MT2}. The main aim of this paper is to investigate separability conditions
for a cowreath $(A\otimes H^{op},H,\psi )$ constructed by using the $8$%
-dimensional Clifford algebra $A=Cl(\alpha ,\beta _{1},\beta _{2},\gamma
_{1},\gamma _{2},\lambda )$ considered as an $H$-comodulo algebra where $H$
is the $8$-dimensional unimodular ribbon Hopf algebra $E(2)$ introduced by
Radford in \cite{R}.
\end{abstract}

\maketitle
\tableofcontents

\section{Introduction}

The categorical notion of separability was introduced in \cite{NVV} where
the notion of separable functor was given and investigated. This theory
extends the concept of separable algebra over commutative rings which was
introduced by M. Auslander and O. Goldmann in \cite{AG}. Its importance
stems also from all the applications in almost all areas of algebra (cf.
\cite{F}).

The notion of separability can be extended to (co)algebras in monoidal
categories (see \cite{BT}). In particular the study of cowreaths constructed
on two-sided Hopf modules was done in \cite[Proposition 7.4]{BCT1} and \cite%
{BCT2}. Recall that, a triple $(A,X,\psi )$is called a cowreath in a
monoidal category $\mathcal{C}$ if $X$ is an object in $\mathcal{C}$, $A$ is
an algebra in $\mathcal{C}$, $\left( X,\psi \right) $ a coalgebra in the
category $\mathcal{T}_{A}^{\#}$ and $\psi :X\otimes A\rightarrow A\otimes X $
is a morphism in $\mathcal{C}$ compatible with algebra structure of $A$ (see
section Preliminars for definitions).

\bigskip The $8$-dimensional unimodular ribbon Hopf algebra $E(2)$ was
introduced by Radford in \cite{R}. This Hopf algebra has been useful to
computing Kauffman`s invariant for knots and Henning's invariant for $3$%
-manifolds in \cite{G}. It is the unique Hopf algebra of dimension $8$ with
coradical $kC_{2}$ see \cite{S}. The general family $E(n)$ was considered in
\cite{CD}.

We will study the separability of a cowreath constructed by using the $8$%
-dimensional Clifford algebra $A=Cl(\alpha ,\beta _{1},\beta _{2},\gamma
_{1},\gamma _{2},\lambda )$ generated by $G,X_{1},X_{2}$ with relations $%
G^{2}=\alpha $, $X_{i}^{2}=\beta _{i}$ and $X_{i}G+GX_{i}=\gamma _{i}$ and $%
X_{1}X_{2}+X_{2}X_{1}=\lambda $. $A$

is an $E\left( 2\right) $-comodule algebra (in fact it is a cleft extension
if $\alpha \neq 0$ \cite{PV}). This research continues the one we did in
\cite{MT1} and \cite{MT2}. Namely in \cite{MT2} we considered the case where
$H=H_{4}$ is the Sweedler $4$-dimensional Hopf algebra over a field $k$ and $%
A=Cl(\alpha ,\beta ,\gamma )$ is the generalized Clifford algebra generated
by two elements $G,X$ and relations $G^{2}=\alpha $, $X^{2}=\beta $ and $%
XG+GX=\gamma $ which is an $H$-comodule algebra. We showed that the cowreath
$\left( A\otimes H^{op},H,\psi \right) $ is always separable. This result
needed no particular effort and thus it was natural to face the natural
subsequent step. Unfortunately here the situation appeared very much more
complex. As before, in view of \cite[Proposition 7.4]{BCT1}, we were seeking
for a $k$-linear map $B:H\otimes H\rightarrow A\otimes H^{op}$satisfying
Casimir condition $\left( \ref{Casimir Condition}\right) $, morphism
condition $\left( \ref{morphism condition}\right) $ and normalized condition
$\left( \ref{normalized condition}\right) $. We could prove that if $B$
satisfies morphism condition (actually just part of it, namely equation $%
\left( \ref{eq.h}\right) $) and $B(1_{H}\otimes 1_{H})=1_{A}\otimes 1_{H}$
then automatically normalized condition was satisfied (see \ref{MCNC}). Then
we concentrated on Casimir condition. By using equalities coming from
morphism condition we could prove (see Proposition $\left( \ref{Pro helem}%
\right) $ that all the elements of the form $B\left( h\otimes h^{\prime
}\right) $ can be obtained from elements of the form $B\left( h\otimes
1_{H}\right) $. Thus we imposed Casimir condition on the eight elements $%
B\left( h\otimes 1_{H}\right) $ for

\begin{equation*}
h\in \left\{ 1_{H},g,x_{1},x_{2},x_{1}x_{2},gx_{1},gx_{2},gx_{1}x_{2}\right\}
\end{equation*}%
and obtained as result $\left( \ref{1ot1}\right) ,\left( \ref{got1}\right)
,\left( \ref{x1}\right) ,\left( \ref{x2}\right) ,\left( \ref{x1x2}\right)
,\left( \ref{gx1}\right) ,\left( \ref{gx2}\right) $ and $\left( \ref{gx1x2}%
\right) .$Then we thought that Casimir condition was automatically satisfied
by all other elements written, as explained above, using these eight
elements . We could find some shortcuts as the one in Proposition $\left( %
\ref{Pro:gh}\right) $, but we had to verify it directly in $23$ cases. This
took us a long time and a wide space. Then we faced morphism condition. We
could prove (see $\left( \ref{mor cond}\right) $ that this condition can be
split in equality $\left( \ref{eq.a}\right) $ and $\left( \ref{eq.h}\right) $%
.We began by examing $\left( \ref{eq.a}\right) .$As we remarked in section $%
\left( \ref{mor cond}\right) $, we only need to compute equality $\left( \ref%
{eq.a}\right) $ for $G\otimes 1_{H}$, $X_{1}\otimes 1_{H}$ and $X_{2}\otimes
1_{H}$. Moreover, as we proved in $\left( \ref{REa}\right) ,$assuming that
equation $\left( \ref{eq.h}\right) $ is satisfied, then if equation $\left( %
\ref{eq.a}\right) $ holds for $h\otimes 1_{H}$ then it holds for any $%
h\otimes h^{\prime }.$ Nevertheless we had to check it for a lot of cases
but, at the end, we discovered that equation $\left( \ref{eq.a}\right) $ was
automatically satisfied. Then we faced equality $\left( \ref{eq.h}\right) $.
Here we obtained a long list of equalities $\left( \text{see \ref{LAE}}%
\right) $. At this point we decided to consider only the non degenerate case
meaning the case where all our constants $\alpha ,\beta _{1},\beta
_{2},\gamma _{1},\gamma _{2},\lambda $ are different from zero. The final
result is encoded in Theorem $\left( \ref{Theo sep}\right) $ which explains
that, in this case, we can find always a bilinear form $B$ satisfying our
three separability conditions.

\bigskip\

\section{Preliminaries}

In the following we will adopt all definitions and notations in \cite{BCT1}
except for the unit of an algebra $A,$ that we will denote by $u:\underline{1%
}\rightarrow A$.

Let $\mathcal{C}$ be a (strict) monoidal category and let $\left(
A,m,u\right) $ be an algebra in $\mathcal{C}$. Recall that a \emph{(right)
transfer morphism through} $A$ is a pair $\left( X,\psi \right) $ with $X\in
\mathcal{C}$ and%
\begin{equation*}
\psi :X\otimes A\rightarrow A\otimes X\text{ in }\mathcal{C}
\end{equation*}%
such that%
\begin{eqnarray*}
\psi \circ \left( X\otimes m\right) &=&\left( X\otimes m\right) \circ \left(
A\otimes \psi \right) \circ \left( \psi \otimes A\right) \\
\psi \circ \left( X\otimes u\right) &=&u\otimes X
\end{eqnarray*}%
The category $\mathcal{T}_{A}^{\#}$ has objects the right transfers. A
morphism $f:X\rightarrow Y$ in $\mathcal{T}_{A}^{\#}$ is a morphism $%
f:X\rightarrow A\otimes Y$ in $\mathcal{C}$ such that%
\begin{equation*}
\left( m\otimes Y\right) \circ \left( A\otimes f\right) \circ \psi =\left(
m\otimes Y\right) \circ \left( A\otimes \psi \right) \circ \left( f\otimes
A\right)
\end{equation*}%
The composition of two morphisms $f:X\rightarrow Y$ and $c:Y\rightarrow Z$
in $\mathcal{T}_{A}^{\#}$ is%
\begin{equation*}
c\odot f=\left( m\otimes Z\right) \circ \left( A\otimes c\right) \circ f
\end{equation*}%
and%
\begin{equation*}
\mathrm{Id}_{\left( X,\psi \right) }=u\otimes X.
\end{equation*}%
The tensor product of $\left( X,\psi _{X}\right) $ and $\left( Y,\psi
_{Y}\right) $ is%
\begin{equation*}
X\circledast Y=\left( X\otimes Y,\psi _{X}\odot \psi _{Y}=\left( \psi
_{X}\otimes Y\right) \circ \left( X\otimes \psi _{Y}\right) \right)
\end{equation*}%
The tensor product of $f:X\rightarrow X^{\prime }$ and $c:Y\rightarrow
Y^{\prime }$ is%
\begin{equation*}
f\circledast c=\left( m\otimes X^{\prime }\otimes Y^{\prime }\right) \circ
\left( A\otimes \psi _{X^{\prime }}\otimes Y^{\prime }\right) \circ \left(
f\otimes c\right) .
\end{equation*}%
The unit object is
\begin{equation*}
\left( \underline{1},r_{A}^{-1}\circ l_{A}:\underline{1}\otimes A\rightarrow
A\otimes \underline{1}\right)
\end{equation*}%
Recall also that a \emph{cowreath in} $\mathcal{C}$ is a triple $\left(
A,X,\psi \right) $ where $A$ is an algebra in $\mathcal{C}$ and $\left(
X,\psi \right) $ is a coalgebra in $\mathcal{T}_{A}^{\#}${\LARGE . }This
means that $\left( X,\psi \right) \in \mathcal{T}_{A}^{\#}$ and there are
morphisms%
\begin{equation*}
\delta :X\rightarrow A\otimes X\otimes X\text{ and }\varepsilon
:X\rightarrow A
\end{equation*}%
such that%
\begin{gather}
\begin{array}{c}
\left( m\otimes X\otimes X\right) \circ \left( A\otimes \psi \otimes
X\right) \circ \left( A\otimes X\otimes \psi \right) \otimes \left( \delta
\otimes A\right) =\left( m\otimes X\otimes X\right) \circ \left( A\otimes
\delta \right) \circ \psi \text{ } \\
\text{i.e. }\delta \text{ is a morphism in }\mathcal{T}_{A}^{\#}%
\end{array}
\label{cow1} \\
\begin{array}{c}
\left( m\otimes X\otimes X\otimes X\otimes X\right) \circ \left( A\otimes
\delta \otimes X\right) \circ \delta =\left( m\otimes X\otimes X\otimes
X\right) \circ \left( A\otimes \psi \otimes X\otimes X\right) \circ \left(
A\otimes X\otimes \delta \right) \circ \delta \\
\text{coassociativity}%
\end{array}%
\text{ }  \label{cow2} \\
m\circ \left( A\otimes \varepsilon \right) \circ \psi =m\circ \left(
\epsilon \otimes A\right) \text{ i.e. }\varepsilon \text{ is a morphism in }%
\mathcal{T}_{A}^{\#}  \label{cow3} \\
\left( m\otimes X\right) \circ \left( A\otimes \varepsilon \otimes X\right)
\circ \delta =u\otimes X\text{ left counit property}  \label{cow4} \\
\left( m\otimes X\right) \circ \left( A\otimes \psi \right) \circ \left(
A\otimes X\otimes \varepsilon \right) \circ \delta =u\otimes X\text{ right
counit property}  \label{cow5}
\end{gather}%
We also recall that an \emph{entwined module over a cowreath} $\left(
A,X,\psi \right) $ is a pair $\left( M,\rho :M\rightarrow M\otimes X\right) $
where $\left( M,\mu \right) \in \mathcal{C}_{A}$ satisfying%
\begin{eqnarray*}
\left( \rho \otimes X\right) \circ \rho &=&\left( \mu \otimes X\otimes
X\right) \circ \left( M\otimes \delta \right) \circ \rho \text{
coassociativity} \\
\mu \circ \left( M\otimes \varepsilon \right) \circ \rho &=&\mathrm{Id}_{M}%
\text{ counitality} \\
\rho \circ \mu &=&\left( \mu \otimes X\right) \circ \left( M\otimes \psi
\right) \circ \left( \rho \otimes A\right) \text{ A-linearity}
\end{eqnarray*}%
A morphism between entwined modules is a $A$-linear morphism $f:M\rightarrow
N$ such that $\left( f\otimes X\right) \circ \rho =\rho \circ f.$

The category of entwined modules will be denoted by $\mathcal{C}\left( \psi
\right) _{A}^{X}.$

Let $H$ be a Hopf algebra, $A$ a right $H$-comodule algebra and $A\otimes
H^{op}$ is a right $H\otimes H^{op}$-comodule algebra with
\begin{equation*}
\rho _{A\otimes H^{op}}:A\otimes H^{op}\rightarrow A\otimes H^{op}\otimes
H\otimes H^{op},
\end{equation*}%
\begin{equation*}
\rho _{A\otimes H^{op}}(d\otimes h)=(d_{0}\otimes h_{1})\otimes
(d_{1}\otimes h_{2}).
\end{equation*}%
$H$ is a right $H\otimes H^{op}$-module coalgebra via%
\begin{equation*}
h\left( h^{\prime }\otimes h^{^{\prime \prime }}\right) =h^{^{\prime \prime
}}hh^{\prime }
\end{equation*}%
and $H$ can be seen as a coalgebra in ${\mathcal{M}}_{H\otimes H^{op}}$=$_{H}%
{\mathcal{M}}_{H}.$ We consider
\begin{equation*}
\begin{array}{cccc}
\psi : & H\otimes A\otimes H^{op} & \rightarrow & A\otimes H^{op}\otimes H
\\
& h\otimes d\otimes l^{op} & \mapsto & d_{0}\otimes l_{1}\otimes l_{2}hd_{1}.%
\end{array}%
\end{equation*}

Then $\left( H,\psi \right) \in \mathcal{T}_{A\otimes H^{op}}^{\#}$ and $%
\left( H,\psi \right) $ is a coalgebra in $\mathcal{T}_{A\otimes
H^{op}}^{\#} $ via%
\begin{eqnarray*}
\delta &:&H\rightarrow A\otimes H^{op}\otimes H\otimes H \\
\delta \left( h\right) &=&1_{A}\otimes 1_{H}\otimes h_{1}\otimes h_{2}
\end{eqnarray*}%
\begin{eqnarray*}
\epsilon &:&H\rightarrow A\otimes H^{op} \\
\epsilon \left( h\right) &=&\varepsilon _{H}\left( h\right) 1_{A}\otimes
1_{H}
\end{eqnarray*}%
The category of Doi-Hopf modules ${\mathcal{M}}(H\otimes H^{op})_{A\otimes
H^{op}}^{H}$ is isomorphic to the category $_{H}{\mathcal{M}}_{A}^{H}$ of
two-sided $(H,A)$-bimodules over $H$.

A \textbf{Casimir morphism} consists of a $k$-linear map $B:H\otimes
H\rightarrow A\otimes H^{op}$ with the following properties:

1. \textbf{Casimir condition}
\begin{eqnarray*}
&&\left( m_{A\otimes H^{op}}\otimes H\right) \circ \left( A\otimes
H^{op}\otimes B\otimes H\right) \circ \left( \psi \otimes H\otimes H\right)
\circ \left( H\otimes \delta \right) \\
&=&\left( m_{A\otimes H^{op}}\otimes H\right) \circ \left( A\otimes
H^{op}\otimes \psi \right) \circ \left( A\otimes H^{op}\otimes H\otimes
B\right) \circ \left( \delta \otimes H\right)
\end{eqnarray*}%
This can be rewritten%
\begin{equation}
B^{A}(h_{2}\otimes h^{\prime })_{0}\otimes B^{H}(h_{2}\otimes h^{\prime
})_{1}\otimes B^{H}(h_{2}\otimes h^{\prime })_{2}h_{1}B^{A}(h_{2}\otimes
h^{\prime })_{1}=B^{A}(h\otimes h_{1}^{\prime })\otimes B^{H}(h\otimes
h_{1}^{\prime })\otimes h_{2}^{\prime }  \label{Casimir Condition}
\end{equation}

2. \textbf{Morphism condition} \newline
i.e. $B$ is a morphism in $\mathcal{T}_{A\otimes H^{op}}^{\#}$:
\begin{equation*}
m_{A\otimes H^{op}}\circ \left( B\otimes A\otimes H^{op}\right) =m_{A\otimes
H^{op}}\circ \left( A\otimes H^{op}\otimes B\right) \circ \left( \psi
\otimes H\right) \circ \left( H\otimes \psi \right)
\end{equation*}

This can be written as

\begin{equation}
B(h\otimes h^{\prime })\cdot (d\otimes h^{\prime \prime })=(d_{0}\otimes
h_{1}^{\prime \prime })\cdot B(h_{2}^{\prime \prime }hd_{1}\otimes
h_{3}^{\prime \prime }h^{\prime }d_{2})  \label{morphism condition}
\end{equation}%
\begin{equation*}
\begin{array}{cccc}
H\otimes \psi : & H\otimes H\otimes A\otimes H^{op} & \rightarrow & H\otimes
A\otimes H^{op}\otimes H \\
& h\otimes h^{\prime }\otimes d\otimes l^{op} & \mapsto & h\otimes
d_{0}\otimes l_{1}\otimes l_{2}h^{\prime }d_{1}.%
\end{array}%
\end{equation*}%
and we consider also the condition

3. \textbf{Normalized condition}
\begin{equation}
m_{A\otimes H^{op}}\circ \left( A\otimes H^{op}\otimes B\right) \circ \delta
=\epsilon  \label{normalized condition}
\end{equation}

\subsection{The Clifford algebra of dimension $8$ and its structure of
comodulo algebra\label{CLIFF}}

We denote by $E\left( 2\right) $ the Radford $8$ dimensional Hopf algebra.
It is generated by $g,x_{1},x_{2}$ with the relations $g^{2}=1_{H}$, $%
x_{i}^{2}=0$ and $x_{i}g=-gx_{i}$ and $x_{i}x_{j}=-x_{j}x_{i}$ for every $%
i,j\in \left\{ 1,2\right\} $. A basis for $H=E\left( 2\right) $ is given by $%
g_{H}=\left( g^{i}x_{1}^{j}x_{2}^{h}\right) _{i,j,k=0,1}.$The coalgebra
structure is given by $\Delta (g)=g\otimes g,\Delta (x_{i})=x_{i}\otimes
g+1_{H}\otimes x_{i},\epsilon (g)=1$ and $\epsilon (x_{i})=0$ and the
antipode $S(g)=g,S(x_{i})=gx_{i}$.

We consider the Clifford algebra $A=Cl(\alpha ,\beta _{i},\gamma
_{i},\lambda )$ generated by $G,X_{1},X_{2}$ with the relations $%
G^{2}=\alpha $, $X_{i}^{2}=\beta _{i}$ and $X_{i}G+GX_{i}=\gamma _{i}$ and $%
X_{1}X_{2}+X_{2}X_{1}=\lambda $. This Clifford algebra is a right $H$%
-comodule algebra via $1_{A}\rightarrow 1_{A}\otimes 1_{H},G\rightarrow
G\otimes g,X_{i}\rightarrow X_{i}\otimes g+1_{A}\otimes x_{i}$. A basis for $%
A$ is given by $G_{A}=\left( G^{i}X_{1}^{j}X_{2}^{h}\right) _{i,j,k=0,1}.$

The main aim of this paper is to find equivalent conditions in order that
the cowreath $(A\otimes H^{op},H,\psi )$ is separable with respect to a
bilinear form. In view of \cite[Proposition 7.4]{BCT1}, we have to find a
bilinear form
\begin{equation*}
B:H\otimes H\rightarrow A\otimes H^{op}
\end{equation*}%
satisfying morphism condition, normalized condition and Casimir condition.

The bilinear form $B:H\otimes H\rightarrow A\otimes H$ is given $h\otimes
h^{\prime }\rightarrow B(h,h^{\prime })$ where in term of bases%
\begin{eqnarray*}
B(h\otimes h^{\prime }) &=&\sum_{a,b\left( X_{1}\right) ,b\left(
X_{2}\right) ,d,e\left( x_{1}\right) ,e\left( x_{2}\right) =0}^{1}B(h\otimes
h^{\prime };G^{a}X_{1}^{b\left( X_{1}\right) }X_{2}^{b\left( X_{2}\right)
},g^{d}x_{1}^{e\left( x_{1}\right) }x_{2}^{e\left( x_{2}\right) }) \\
&&G^{a}X_{1}^{b\left( X_{1}\right) }X_{2}^{b\left( X_{2}\right) }\otimes
g^{d}x_{1}^{e\left( x_{1}\right) }x_{2}^{e\left( x_{2}\right) }.
\end{eqnarray*}%
For sake of simplicity, when no misunderstanding will arise, we will simply
write%
\begin{eqnarray*}
b_{1} &=&b\left( X_{1}\right) \\
b_{2} &=&b\left( X_{2}\right) \\
e_{1} &=&e\left( x_{1}\right) \\
e_{2} &=&e\left( x_{2}\right)
\end{eqnarray*}%
\begin{equation*}
g^{d}x_{1}^{e_{1}}x_{2}^{e_{2}}\cdot g^{\overline{d}}x_{1}^{\overline{e_{1}}%
}x_{2}^{\overline{e_{2}}}=\left( -1\right) ^{\overline{d}\left(
e_{1}+e_{2}\right) }\left( -1\right) ^{\overline{e_{1}}e_{2}}g^{d+\overline{d%
}}x_{1}^{e_{1}+\overline{e_{1}}}x_{2}^{e_{2}+\overline{e_{2}}}
\end{equation*}%
\begin{equation*}
g^{\overline{d}}x_{1}^{\overline{e_{1}}}x_{2}^{\overline{e_{2}}}\cdot
g^{d}x_{1}^{e_{1}}x_{2}^{e_{2}}=\left( -1\right) ^{d\left( \overline{e_{1}}+%
\overline{e_{2}}\right) }\left( -1\right) ^{e_{1}\overline{e_{2}}}g^{d+%
\overline{d}}x_{1}^{e_{1}+\overline{e_{1}}}x_{2}^{e_{2}+\overline{e_{2}}}
\end{equation*}%
\begin{equation*}
\frac{g^{d}x_{1}^{e_{1}}x_{2}^{e_{2}}\cdot g^{\overline{d}}x_{1}^{\overline{%
e_{1}}}x_{2}^{\overline{e_{2}}}}{g^{\overline{d}}x_{1}^{\overline{e_{1}}%
}x_{2}^{\overline{e_{2}}}\cdot g^{d}x_{1}^{e_{1}}x_{2}^{e_{2}}}=\frac{\left(
-1\right) ^{\overline{d}\left( e_{1}+e_{2}\right) +\overline{e_{1}}e_{2}}}{%
\left( -1\right) ^{d\left( \overline{e_{1}}+\overline{e_{2}}\right) +%
\overline{e_{1}}e_{2}+}}
\end{equation*}%
\begin{equation*}
g^{d}x_{1}^{e_{1}}x_{2}^{e_{2}}\cdot g^{\overline{d}}x_{1}^{\overline{e_{1}}%
}x_{2}^{\overline{e_{2}}}=\left( -1\right) ^{e_{1}+e_{2}}g^{d+\overline{d}%
}x_{1}^{e_{1}}x_{2}^{e_{2}}\cdot x_{1}^{\overline{e_{1}}}x_{2}^{\overline{%
e_{2}}}
\end{equation*}%
so that we get%
\begin{equation*}
B(h\otimes h^{\prime })=\sum_{a,b_{1},b_{2},d,e_{1},e_{2}=0}^{1}B(h\otimes
h^{\prime
};G^{a}X_{1}^{b_{1}}X_{2}^{b_{2}},f)G^{a}X_{1}^{b_{1}}X_{2}^{b_{2}}\otimes
g^{d}x_{1}^{e_{1}}x_{2}^{e_{2}}.
\end{equation*}%
The comultiplication on $H$ can be written:%
\begin{equation}
\Delta
(g^{d}x_{1}^{e_{1}}x_{2}^{e_{2}})=\sum_{u_{1}=0}^{e_{1}}%
\sum_{u_{2}=0}^{e_{2}}\left( -1\right) ^{\left( e_{2}+u_{2}\right)
u_{1}}g^{d}x_{1}^{\left( e_{1}-u_{1}\right) }x_{2}^{\left(
e_{2}-u_{2}\right) }\otimes
g^{d+e_{1}+e_{2}+u_{1}+u_{2}}x_{1}^{u_{1}}x_{2}^{u_{2}}  \label{form Delta}
\end{equation}

The $H$-coaction on $A$ is given by:

\begin{equation}
\rho
_{A}(G^{a}X_{1}^{b_{1}}X_{2}^{b_{2}})=\sum_{l_{1}=0}^{b_{1}}%
\sum_{l_{2}=0}^{b_{2}}\left( -1\right) ^{\left( l_{2}+b_{2}\right)
l_{1}}G^{a}X_{1}^{b_{1}-l_{1}}X_{2}^{b_{2}-l_{2}}\otimes
g^{a+b_{1}+b_{2}+l_{1}+l_{2}}x_{1}^{l_{1}}x_{2}^{l_{2}}  \label{form ro}
\end{equation}

\part{The plan}

In this part we will discuss the conditions we need for the cowreath $%
(A\otimes H^{op},H)$ is separable with respect to a bilinear form $B$ and
discover the reductions we can make.

\section{The Casimir Condition}

The Casimir condition is given by the formula in $A\otimes H^{op}\otimes H$

\begin{equation*}
B^{A}(h_{2}\otimes h^{\prime })_{0}\otimes B^{H}(h_{2}\otimes h^{\prime
})_{1}\otimes B^{H}(h_{2}\otimes h^{\prime })_{2}h_{1}B^{A}(h_{2}\otimes
h^{\prime })_{1}=B^{A}(h\otimes h_{1}^{\prime })\otimes B^{H}(h\otimes
h_{1}^{\prime })\otimes h_{2}^{\prime }
\end{equation*}

\begin{equation*}
\Delta
(g^{m}x_{1}^{n_{1}}x_{2}^{n_{2}})=\sum_{w_{1}=0}^{n_{1}}%
\sum_{w_{2}=0}^{n_{2}}\left( -1\right) ^{\left( n_{2}+w_{2}\right)
w_{1}}g^{m}x_{1}^{\left( n_{1}-w_{1}\right) }x_{2}^{\left(
n_{2}-w_{2}\right) }\otimes
g^{m+n_{1}+n_{2}+w_{1}+w_{2}}x_{1}^{w_{1}}x_{2}^{w_{2}}
\end{equation*}%
\begin{eqnarray*}
&&B^{A}(\left( g^{m}x_{1}^{n_{1}}x_{2}^{n_{2}}\right) _{2}\otimes
h)_{0}\otimes B^{H}(\left( g^{m}x_{1}^{n_{1}}x_{2}^{n_{2}}\right)
_{2}\otimes h)_{1} \\
&&\otimes B^{H}(\left( g^{m}x_{1}^{n_{1}}x_{2}^{n_{2}}\right) _{2}\otimes
h)_{2}\left( \left( g^{m}x_{1}^{n_{1}}x_{2}^{n_{2}}\right) \right)
_{1}B^{A}(\left( g^{m}x_{1}^{n_{1}}x_{2}^{n_{2}}\right) _{2}\otimes h)_{1} \\
&=&\sum_{w_{1}=0}^{n_{1}}\sum_{w_{2}=0}^{n_{2}}\left( -1\right) ^{\left(
n_{2}+w_{2}\right)
w_{1}}B^{A}(g^{m+n_{1}+n_{2}+w_{1}+w_{2}}x_{1}^{w_{1}}x_{2}^{w_{2}}\otimes
h)_{0}\otimes \\
&&\otimes
B^{H}(g^{m+n_{1}+n_{2}+w_{1}+w_{2}}x_{1}^{w_{1}}x_{2}^{w_{2}}\otimes
h)_{1}\otimes \\
&&B^{H}(g^{m+n_{1}+n_{2}+w_{1}+w_{2}}x_{1}^{w_{1}}x_{2}^{w_{2}}\otimes
h)_{2}\left( g^{m}x_{1}^{\left( n_{1}-w_{1}\right) }x_{2}^{\left(
n_{2}-w_{2}\right) }\right)
B^{A}(g^{m+n_{1}+n_{2}+w_{1}+w_{2}}x_{1}^{w_{1}}x_{2}^{w_{2}}\otimes h)_{1}
\end{eqnarray*}%
Now
\begin{eqnarray*}
&&B(g^{m+n_{1}+n_{2}+w_{1}+w_{2}}x_{1}^{w_{1}}x_{2}^{w_{2}}\otimes h) \\
&=&%
\sum_{a,b_{1},b_{2},d,e_{1},e_{2}=0}^{1}B(g^{m+n_{1}+n_{2}+w_{1}+w_{2}}x_{1}^{w_{1}}x_{2}^{w_{2}}\otimes h;G^{a}X_{1}^{b_{1}}X_{2}^{b_{2}},g^{d}x_{1}^{e_{1}}x_{2}^{e_{2}})
\\
&&G^{a}X_{1}^{b_{1}}X_{2}^{b_{2}}\otimes g^{d}x_{1}^{e_{1}}x_{2}^{e_{2}}.
\end{eqnarray*}

Thus, using $\left( \ref{form Delta}\right) $ and $\left( \ref{form ro}%
\right) $ we obtain%
\begin{eqnarray*}
&&B^{A}(\left( g^{m}x_{1}^{n_{1}}x_{2}^{n_{2}}\right) _{2}\otimes h^{\prime
})_{0}\otimes B^{H}(\left( g^{m}x_{1}^{n_{1}}x_{2}^{n_{2}}\right)
_{2}\otimes h^{\prime })_{1} \\
&&\otimes B^{H}(\left( g^{m}x_{1}^{n_{1}}x_{2}^{n_{2}}\right) _{2}\otimes
h^{\prime })_{2}\left( \left( g^{m}x_{1}^{n_{1}}x_{2}^{n_{2}}\right) \right)
_{1}B^{A}(\left( g^{m}x_{1}^{n_{1}}x_{2}^{n_{2}}\right) \otimes h^{\prime
})_{1} \\
&=&\sum_{w_{1}=0}^{n_{1}}\sum_{w_{2}=0}^{n_{2}}\left( -1\right) ^{\left(
n_{2}+w_{2}\right)
w_{1}}\sum_{a,b_{1},b_{2},d,e_{1},e_{2}=0}^{1}\sum_{l_{1}=0}^{b_{1}}%
\sum_{l_{2}=0}^{b_{2}}\sum_{u_{1}=0}^{e_{1}}\sum_{u_{2}=0}^{e_{2}}\left(
-1\right) ^{\left( b_{2}+l_{2}\right) l_{1}} \\
&&\left( -1\right) ^{\left( b_{2}+l_{2}\right) l_{1}}\left( -1\right)
^{\left( e_{2}+u_{2}\right) u_{1}}\left( -1\right) ^{\left(
b_{2}+l_{2}\right) l_{1}}\left( -1\right) ^{\left( e_{2}+u_{2}\right) u_{1}}
\\
&&B(\left( g^{m+n_{1}+n_{2}+w_{1}+w_{2}}x_{1}^{w_{1}}x_{2}^{w_{2}}\right)
\otimes h^{\prime
};G^{a}X_{1}^{b_{1}}X_{2}^{b_{2}},g^{d}x_{1}^{e_{1}}x_{2}^{e_{2}}) \\
&&G^{a}X_{1}^{b_{1}-l_{1}}X_{2}^{b_{2}-l_{2}}\otimes g^{d}x_{1}^{\left(
e_{1}-u_{1}\right) }x_{2}^{\left( e_{2}-u_{2}\right) }\otimes \\
&&\otimes
g^{d+e_{1}+e_{2}+u_{1}+u_{2}}x_{1}^{u_{1}}x_{2}^{u_{2}}g^{m}x_{1}^{\left(
n_{1}-w_{1}\right) }x_{2}^{\left( n_{2}-w_{2}\right)
}g^{a+b_{1}+b_{2}+l_{1}+l_{2}}x_{1}^{l_{1}}x_{2}^{l_{2}}
\end{eqnarray*}%
Since
\begin{eqnarray*}
&&g^{d+e_{1}+e_{2}+u_{1}+u_{2}}x_{1}^{u_{1}}x_{2}^{u_{2}}g^{m}x_{1}^{\left(
n_{1}-w_{1}\right) }x_{2}^{\left( n_{2}-w_{2}\right)
}g^{a+b_{1}+b_{2}+l_{1}+l_{2}}x_{1}^{l_{1}}x_{2}^{l_{2}} \\
&=&(-1)^{\left( a+b_{1}+b_{2}+l_{1}+l_{2}\right) \left(
u_{1}+u_{2}+n_{1}-w_{1}+n_{2}-w_{2}\right) +m(u_{1}+u_{2})+\left(
n_{1}-w_{1}\right) u_{2}+l_{1}(n_{2}-w_{2}+u_{2})} \\
&&g^{a+b_{1}+b_{2}+l_{1}+l_{2}+d+e_{1}+e_{2}+u_{1}+u_{2}+m}x_{1}^{u_{1}+l_{1}+n_{1}-w_{1}}x_{2}^{u_{2}+l_{2}+n_{2}-w_{2}}
\end{eqnarray*}%
we finally get%
\begin{eqnarray*}
&&B^{A}(\left( g^{m}x_{1}^{n_{1}}x_{2}^{n_{2}}\right) _{2}\otimes h^{\prime
})_{0}\otimes B^{H}(\left( g^{m}x_{1}^{n_{1}}x_{2}^{n_{2}}\right)
_{2}\otimes h^{\prime })_{1} \\
&&\otimes B^{H}(\left( g^{m}x_{1}^{n_{1}}x_{2}^{n_{2}}\right) _{2}\otimes
h^{\prime })_{2}\left( \left( g^{m}x_{1}^{n_{1}}x_{2}^{n_{2}}\right) \right)
_{1}B^{A}(\left( g^{m}x_{1}^{n_{1}}x_{2}^{n_{2}}\right) \otimes 1_{H})_{1} \\
&=&\sum_{w_{1}=0}^{n_{1}}\sum_{w_{2}=0}^{n_{2}}\left( -1\right) ^{\left(
n_{2}+w_{2}\right)
w_{1}}\sum_{a,b_{1},b_{2},d,e_{1},e_{2}=0}^{1}\sum_{l_{1}=0}^{b_{1}}%
\sum_{l_{2}=0}^{b_{2}}\sum_{u_{1}=0}^{e_{1}}\sum_{u_{2}=0}^{e_{2}} \\
&&\left( -1\right) ^{\alpha \left(
g^{m}x_{1}^{n_{1}-w_{1}}x_{2}^{n_{2}-w_{2}}\right) }\left( -1\right)
^{\alpha \left( g^{m}x_{1}^{n_{1}-w_{1}}x_{2}^{n_{2}-w_{2}}\right) }\left(
-1\right) ^{\alpha \left( g^{m}x_{1}^{n_{1}-w_{1}}x_{2}^{n_{2}-w_{2}}\right)
} \\
&&B(\left( g^{m+n_{1}+n_{2}+w_{1}+w_{2}}x_{1}^{w_{1}}x_{2}^{w_{2}}\right)
\otimes h^{\prime
};G^{a}X_{1}^{b_{1}}X_{2}^{b_{2}},g^{d}x_{1}^{e_{1}}x_{2}^{e_{2}}) \\
&&G^{a}X_{1}^{b_{1}-l_{1}}X_{2}^{b_{2}-l_{2}}\otimes g^{d}x_{1}^{\left(
e_{1}-u_{1}\right) }x_{2}^{\left( e_{2}-u_{2}\right) }\otimes \\
&&g^{a+b_{1}+b_{2}+l_{1}+l_{2}+d+e_{1}+e_{2}+u_{1}+u_{2}+m}x_{1}^{u_{1}+l_{1}+n_{1}-w_{1}}x_{2}^{u_{2}+l_{2}+n_{2}-w_{2}}
\end{eqnarray*}%
We set%
\begin{gather*}
\alpha \left(
g^{m}x_{1}^{n_{1}}x_{2}^{n_{2}};a,b_{1},b_{2},e_{1},e_{2};l_{1},l_{2},u_{1},u_{2}\right) =\left( u_{2}+e_{2}\right) u_{1}+\left( l_{2}+b_{2}\right) l_{1}
\\
+\left( a+b_{1}+b_{2}+l_{1}+l_{2}\right) \left(
u_{1}+u_{2}+n_{1}+n_{2}\right) +m(u_{1}+u_{2})+n_{1}u_{2}+l_{1}(n_{2}+u_{2}).
\end{gather*}%
For sake of shortness, we set, when the context will allow us,%
\begin{equation*}
\alpha \left(
g^{m}x_{1}^{n_{1}}x_{2}^{n_{2}};a,b_{1},b_{2},e_{1},e_{2};l_{1},l_{2},u_{1},u_{2}\right) =\alpha \left( g^{m}x_{1}^{n_{1}}x_{2}^{n_{2}};a,b_{1},b_{2},e_{1},e_{2}\right) .
\end{equation*}%
We deduce that%
\begin{eqnarray*}
&&\alpha \left( g^{m}x_{1}^{n_{1}-w_{1}}x_{2}^{n_{2}-w_{2}}\right) \\
&=&\left( a+b_{1}+b_{2}+l_{1}+l_{2}\right) \left(
u_{1}+u_{2}+n_{1}-w_{1}+n_{2}-w_{2}\right) +m(u_{1}+u_{2}) \\
&&+\left( n_{1}-w_{1}\right) u_{2}+l_{1}(n_{2}-w_{2}+u_{2})+\left(
l_{2}+b_{2}\right) l_{1}+\left( u_{2}+e_{2}\right) u_{1}
\end{eqnarray*}%
\begin{equation*}
B^{A}(h_{2}\otimes h^{\prime })_{0}\otimes B^{H}(h_{2}\otimes h^{\prime
})_{1}\otimes B^{H}(h_{2}\otimes h^{\prime })_{2}h_{1}B^{A}(h_{2}\otimes
h^{\prime })_{1}=B^{A}(h\otimes h_{1}^{\prime })\otimes B^{H}(h\otimes
h_{1}^{\prime })\otimes h_{2}^{\prime }
\end{equation*}%
\begin{equation*}
\Delta (g^{\mu }x_{1}^{\nu _{1}}x_{2}^{\nu _{2}})=\sum_{\omega _{1}=0}^{\nu
_{1}}\sum_{\omega _{2}=0}^{\nu _{2}}\left( -1\right) ^{\left( \nu
_{2}+\omega _{2}\right) \omega _{1}}g^{\mu }x_{1}^{\left( \nu _{1}-\omega
_{1}\right) }x_{2}^{\left( \nu _{2}-\omega _{2}\right) }\otimes g^{\mu +\nu
_{1}+\nu _{2}+\omega _{1}+\omega _{2}}x_{1}^{\omega _{1}}x_{2}^{\omega _{2}}
\end{equation*}%
We get%
\begin{eqnarray*}
&&B^{A}(h\otimes h_{1}^{\prime })\otimes B^{H}(h\otimes h_{1}^{\prime
})\otimes h_{2}^{\prime } \\
&=&\sum_{\omega _{1}=0}^{\nu _{1}}\sum_{\omega _{2}=0}^{\nu _{2}}\left(
-1\right) ^{\left( \nu _{2}+\omega _{2}\right) \omega _{1}}B^{A}(h\otimes
g^{\mu }x_{1}^{\left( \nu _{1}-\omega _{1}\right) }x_{2}^{\left( \nu
_{2}-\omega _{2}\right) })\otimes B^{H}(h\otimes g^{\mu }x_{1}^{\left( \nu
_{1}-\omega _{1}\right) }x_{2}^{\left( \nu _{2}-\omega _{2}\right) }) \\
&&\otimes g^{\mu +\nu _{1}+\nu _{2}+\omega _{1}+\omega _{2}}x_{1}^{\omega
_{1}}x_{2}^{\omega _{2}}
\end{eqnarray*}%
For $h=g^{m}x_{1}^{n_{1}}x_{2}^{n_{2}}$ and $h^{\prime }=g^{\mu }x_{1}^{\nu
_{1}}x_{2}^{\nu _{2}}$ we finally get%
\begin{gather}
\sum_{w_{1}=0}^{n_{1}}\sum_{w_{2}=0}^{n_{2}}\left( -1\right) ^{\left(
n_{2}+w_{2}\right)
w_{1}}\sum_{a,b_{1},b_{2},d,e_{1},e_{2}=0}^{1}\sum_{l_{1}=0}^{b_{1}}%
\sum_{l_{2}=0}^{b_{2}}\sum_{u_{1}=0}^{e_{1}}\sum_{u_{2}=0}^{e_{2}}
\label{MAIN FORMULA 1} \\
\left( -1\right) ^{\alpha \left(
g^{m}x_{1}^{n_{1}-w_{1}}x_{2}^{n_{2}-w_{2}};l_{1},l_{2},u_{1},u_{2}\right)
}B(g^{m+n_{1}+n_{2}+w_{1}+w_{2}}x_{1}^{w_{1}}x_{2}^{w_{2}}\otimes g^{\mu
}x_{1}^{\nu _{1}}x_{2}^{\nu
_{2}};G^{a}X_{1}^{b_{1}}X_{2}^{b_{2}},g^{d}x_{1}^{e_{1}}x_{2}^{e_{2}})
\notag \\
G^{a}X_{1}^{b_{1}-l_{1}}X_{2}^{b_{2}-l_{2}}\otimes
g^{d}x_{1}^{e_{1}-u_{1}}x_{2}^{e_{2}-u_{2}}\otimes  \notag \\
g^{a+b_{1}+b_{2}+l_{1}+l_{2}+d+e_{1}+e_{2}+u_{1}+u_{2}+m}x_{1}^{l_{1}+u_{1}+n_{1}-w_{1}}x_{2}^{l_{2}+u_{2}+n_{2}-w_{2}}
\notag \\
=\sum_{\omega _{1}=0}^{\nu _{1}}\sum_{\omega _{2}=0}^{\nu _{2}}\left(
-1\right) ^{\left( \nu _{2}+\omega _{2}\right) \omega
_{1}}B^{A}(g^{m}x_{1}^{n_{1}}x_{2}^{n_{2}}\otimes g^{\mu }x_{1}^{\nu
_{1}-\omega _{1}}x_{2}^{\nu _{2}-\omega _{2}})\otimes  \notag \\
\otimes B^{H}(g^{m}x_{1}^{n_{1}}x_{2}^{n_{2}}\otimes g^{\mu }x_{1}^{\nu
_{1}-\omega _{1}}x_{2}^{\nu _{2}-\omega _{2}})\otimes g^{\mu +\nu _{1}+\nu
_{2}+\omega _{1}+\omega _{2}}x_{1}^{\omega _{1}}x_{2}^{\omega _{2}}  \notag
\end{gather}

\section{CALCULATION of $\protect\alpha $}

For sake of simplicity we will write, whenever there is no risk of confusion,%
\begin{eqnarray*}
\alpha \left( g^{m}x_{1}^{n_{1}}x_{2}^{n_{2}}\right) &=&\alpha \left(
g^{m}x_{1}^{n_{1}}x_{2}^{n_{2}};l_{1},l_{2},u_{1},u_{2}\right) =\left(
u_{2}+e_{2}\right) u_{1}+\left( l_{2}+b_{2}\right) l_{1} \\
&+&\left( a+b_{1}+b_{2}+l_{1}+l_{2}\right) \left(
u_{1}+u_{2}+n_{1}+n_{2}\right) +m(u_{1}+u_{2}) \\
&&+n_{1}u_{2}+l_{1}(n_{2}+u_{2})
\end{eqnarray*}%
Here we list the several occurrences.

\subsubsection{$\left( l_{1},l_{2},u_{1},u_{2}\right) =\left( 0,0,0,0\right)
$%
\protect\begin{equation*}
\protect\alpha \left( g^{m}x_{1}^{n_{1}}x_{2}^{n_{2}};0,0,0,0\right) =\left(
a+b_{1}+b_{2}\right) \left( n_{1}+n_{2}\right)
\protect\end{equation*}%
}

\subsubsection{$\left( l_{1},l_{2},u_{1},u_{2}\right) =\left( 0,0,0,1\right)
$}

\begin{equation*}
\alpha \left( g^{m}x_{1}^{n_{1}}x_{2}^{n_{2}};0,0,0,1\right) =\left(
a+b_{1}+b_{2}\right) \left( 1+n_{1}+n_{2}\right) +m+n_{1}
\end{equation*}

\subsubsection{$\left( l_{1},l_{2},u_{1},u_{2}\right) =\left( 0,0,1,0\right)
$}

\begin{equation*}
\alpha \left( g^{m}x_{1}^{n_{1}}x_{2}^{n_{2}};0,0,1,0\right) =e_{2}+\left(
a+b_{1}+b_{2}\right) \left( 1+n_{1}+n_{2}\right) +m
\end{equation*}

\subsubsection{$\left( l_{1},l_{2},u_{1},u_{2}\right) =\left( 0,1,0,0\right)
$}

\begin{equation*}
\alpha \left( g^{m}x_{1}^{n_{1}}x_{2}^{n_{2}};0,1,0,0\right) =\left(
a+b_{1}+b_{2}+1\right) \left( n_{1}+n_{2}\right)
\end{equation*}

\subsubsection{$\left( l_{1},l_{2},u_{1},u_{2}\right) =\left( 1,0,0,0\right)
$}

\begin{equation*}
\alpha \left( g^{m}x_{1}^{n_{1}}x_{2}^{n_{2}};1,0,0,0\right) =b_{2}+\left(
a+b_{1}+b_{2}+1\right) \left( n_{1}+n_{2}\right) +n_{2}
\end{equation*}

\subsubsection{$\left( l_{1},l_{2},u_{1},u_{2}\right) =\left( 0,0,1,1\right)
$}

\begin{equation*}
\alpha \left( g^{m}x_{1}^{n_{1}}x_{2}^{n_{2}}\right) =\alpha \left(
g^{m}x_{1}^{n_{1}}x_{2}^{n_{2}};0,0,1,1\right) =1+e_{2}+\left(
a+b_{1}+b_{2}\right) \left( n_{1}+n_{2}\right) +n_{1}
\end{equation*}

\subsubsection{$\left( l_{1},l_{2},u_{1},u_{2}\right) =\left( 0,1,0,1\right)
$}

\begin{equation*}
\alpha \left( g^{m}x_{1}^{n_{1}}x_{2}^{n_{2}}\right) =\alpha \left(
g^{m}x_{1}^{n_{1}}x_{2}^{n_{2}};0,1,0,1\right) =\left(
a+b_{1}+b_{2}+1\right) \left( 1+n_{1}+n_{2}\right) +m+n_{1}
\end{equation*}

\subsubsection{$\left( l_{1},l_{2},u_{1},u_{2}\right) =\left( 1,0,0,1\right)
$}

\begin{equation*}
\alpha \left( g^{m}x_{1}^{n_{1}}x_{2}^{n_{2}};1,0,0,1\right) =b_{2}+\left(
a+b_{1}+b_{2}+1\right) \left( 1+n_{1}+n_{2}\right) +m+n_{1}+(n_{2}+1)
\end{equation*}

\subsubsection{$\left( l_{1},l_{2},u_{1},u_{2}\right) =\left( 0,1,1,0\right)
$}

\begin{equation*}
\alpha \left( g^{m}x_{1}^{n_{1}}x_{2}^{n_{2}};0,1,1,0\right) =e_{2}+\left(
a+b_{1}+b_{2}+1\right) \left( 1+n_{1}+n_{2}\right) +m
\end{equation*}

\subsubsection{$\left( l_{1},l_{2},u_{1},u_{2}\right) =\left( 1,0,1,0\right)
$}

\begin{equation*}
\alpha \left( g^{m}x_{1}^{n_{1}}x_{2}^{n_{2}};1,0,1,0\right)
=e_{2}+b_{2}+\left( a+b_{1}+b_{2}+1\right) \left( 1+n_{1}+n_{2}\right)
+m+n_{2}
\end{equation*}

\subsubsection{$\left( l_{1},l_{2},u_{1},u_{2}\right) =\left( 1,1,0,0\right)
$}

\begin{equation*}
\alpha \left( g^{m}x_{1}^{n_{1}}x_{2}^{n_{2}};1,1,0,0\right) =1+b_{2}+\left(
a+b_{1}+b_{2}\right) \left( n_{1}+n_{2}\right) +n_{2}
\end{equation*}

\subsubsection{$\left( l_{1},l_{2},u_{1},u_{2}\right) =\left( 0,1,1,1\right)
$}

\begin{equation*}
\alpha \left( g^{m}x_{1}^{n_{1}}x_{2}^{n_{2}};0,1,1,1\right) =\left(
1+e_{2}\right) +\left( a+b_{1}+b_{2}+1\right) \left( n_{1}+n_{2}\right)
+n_{1}
\end{equation*}

\subsubsection{$\left( l_{1},l_{2},u_{1},u_{2}\right) =\left( 1,0,1,1\right)
$}

\begin{equation*}
\alpha \left( g^{m}x_{1}^{n_{1}}x_{2}^{n_{2}};1,0,1,1\right)
=e_{2}+b_{2}+\left( a+b_{1}+b_{2}+1\right) \left( n_{1}+n_{2}\right)
+n_{1}+n_{2}
\end{equation*}

\subsubsection{$\left( l_{1},l_{2},u_{1},u_{2}\right) =\left( 1,1,0,1\right)
$}

\begin{equation*}
\alpha \left( g^{m}x_{1}^{n_{1}}x_{2}^{n_{2}};1,1,0,1\right) =b_{2}+\left(
a+b_{1}+b_{2}\right) \left( 1+n_{1}+n_{2}\right) +m+n_{1}+n_{2}
\end{equation*}

\subsubsection{$\left( l_{1},l_{2},u_{1},u_{2}\right) =\left( 1,1,1,0\right)
$}

\begin{equation*}
\alpha \left( g^{m}x_{1}^{n_{1}}x_{2}^{n_{2}};1,1,1,0\right)
=e_{2}+1+b_{2}+\left( a+b_{1}+b_{2}\right) \left( 1+n_{1}+n_{2}\right)
+m+n_{2}
\end{equation*}

\subsubsection{$\left( l_{1},l_{2},u_{1},u_{2}\right) =\left( 1,1,1,1\right)
$}

\begin{equation*}
\alpha \left( g^{m}x_{1}^{n_{1}}x_{2}^{n_{2}};1,1,1,1\right)
=e_{2}+b_{2}+\left( a+b_{1}+b_{2}\right) \left( n_{1}+n_{2}\right)
+n_{1}+n_{2}+1
\end{equation*}

\section{Morphism condition\label{mor cond}}

Now, we study when $B:H\otimes H\rightarrow A\otimes H^{op}$ is a morphism
in $\mathcal{T}_{A\otimes H^{op}}^{\sharp }$. i.e.
\begin{equation*}
m_{A\otimes H^{op}}\circ (B\otimes A\otimes H^{op})=m_{A\otimes H^{op}}\circ
(A\otimes B)\circ (\psi \otimes H)\circ (H\otimes \psi )
\end{equation*}

This can be written as, see $\left( \ref{morphism condition}\right) ,$

\begin{equation*}
B(h\otimes h^{\prime })\cdot (d\otimes h^{\prime \prime })=(d_{0}\otimes
h_{1}^{\prime \prime })\cdot B(h_{2}^{\prime \prime }hd_{1}\otimes
h_{3}^{\prime \prime }h^{\prime }d_{2})
\end{equation*}

This equality can be split in%
\begin{equation}
B(h\otimes h^{\prime })(d\otimes 1_{H})=(d_{0}\otimes 1_{H})B(hd_{1}\otimes
h^{\prime }d_{2})  \label{eq.a}
\end{equation}

\begin{equation}
B(h\otimes h^{\prime })(1_{A}\otimes h^{\prime \prime })=(1_{A}\otimes
h_{1}^{\prime \prime })B(h_{2}^{\prime \prime }h\otimes h_{3}^{\prime \prime
}h^{\prime })  \label{eq.h}
\end{equation}

In fact assume that $\left( \ref{eq.a}\right) $ and $\left( \ref{eq.h}%
\right) $ hold. Then
\begin{eqnarray*}
B(h\otimes h^{\prime })\cdot (d\otimes h^{\prime \prime }) &=&B(h\otimes
h^{\prime })\cdot \left( d\otimes 1_{H}\right) \cdot \left( 1_{A}\otimes
h^{\prime \prime }\right) \\
&&\overset{\left( \ref{eq.a}\right) }{=}(d_{0}\otimes 1_{H})\cdot
B(hd_{1}\otimes h^{\prime }d_{2})\cdot \left( 1_{A}\otimes h^{\prime \prime
}\right) \\
&&\overset{\left( \ref{eq.h}\right) }{=}(d_{0}\otimes 1_{H})\cdot
(1_{A}\otimes h_{1}^{\prime \prime })\cdot B(h_{2}^{\prime \prime
}hd_{1}\otimes h_{3}^{\prime \prime }h^{\prime }d_{2})
\end{eqnarray*}%
Now if the equality $\left( \ref{eq.a}\right) $ is true for $d,b\in A$, then
it is true for the product:

\begin{equation*}
\begin{array}{lll}
B(h\otimes h^{\prime })(db\otimes 1_{H}) & = & B(h\otimes h^{\prime
})(d\otimes 1_{H})(b\otimes 1_{H}) \\
& = & (d_{0}\otimes 1_{H})B(hd_{1}\otimes h^{\prime }d_{2})(b\otimes 1_{H})
\\
& = & (d_{0}\otimes 1_{H})(b_{0}\otimes 1_{H})B(hd_{1}b_{1}\otimes h^{\prime
}d_{2}b_{2}) \\
& = & ((db)_{0}\otimes 1_{H})B(h(db)_{1}\otimes h^{\prime }(db)_{2}).%
\end{array}%
\end{equation*}

Similarly if the equality $\left( \ref{eq.h}\right) $ is true for $s,t\in H,$
then it is true for the product:

\begin{equation*}
\begin{array}{lll}
B(h\otimes h^{\prime })(1_{A}\otimes st) & = & B(h\otimes h^{\prime
})(1_{A}\otimes t)(1\otimes s) \\
& = & (1_{A}\otimes t_{1})B(t_{2}h\otimes t_{3}h^{\prime })(1\times s) \\
& = & (1_{A}\otimes t_{1})(1_{A}\otimes s_{1})B(s_{2}t_{2}h\otimes
s_{3}t_{3}h^{\prime }) \\
& = & (1_{A}\otimes (st)_{1})B((st)_{2}h\otimes (st)_{3}h^{\prime }).%
\end{array}%
\end{equation*}

Thus we only need to compute equality $\left( \ref{eq.h}\right) $ for $%
1_{A}\otimes g$, $1_{A}\otimes x_{1}$ and $1_{A}\otimes x_{2}$ and equality $%
\left( \ref{eq.a}\right) $ for $G\otimes 1_{H}$, $X_{1}\otimes 1_{H}$ and $%
X_{2}\otimes 1_{H}$.

\section{Normalized condition}

The normalized condition is%
\begin{equation*}
B(h_{1}\otimes h_{2})=\epsilon (h)1_{A}\otimes 1_{H}\text{ for all }h\in
H_{4}.
\end{equation*}%
By using $\left( \ref{form Delta}\right) $ we deduce that%
\begin{eqnarray}
B(1_{H}\otimes 1_{H}) &=&1_{A}\otimes 1_{H},  \label{1ot1} \\
B(g\otimes g) &=&1_{A}\otimes 1_{H},  \label{gotg} \\
B(x_{i}\otimes g)+B(1_{H}\otimes x_{i}) &=&0,\text{for }i=1,2,  \label{xiotg}
\\
B(gx_{i}\otimes 1_{H})+B(g\otimes gx_{i}) &=&0,\text{for }i=1,2.
\label{gxiot1} \\
B\left( x_{1}x_{2}\otimes 1_{H}\right) +B\left( x_{1}\otimes gx_{2}\right)
-B\left( x_{2}\otimes gx_{1}\right) +B\left( 1_{H}\otimes x_{1}x_{2}\right)
&=&0  \label{x1x2ot1} \\
B\left( gx_{1}x_{2}\otimes g\right) +B\left( gx_{1}\otimes x_{2}\right)
-B\left( gx_{2}\otimes x_{1}\right) +B\left( g\otimes gx_{1}x_{2}\right) &=&0
\label{gx1x2otg}
\end{eqnarray}

\section{Morphism condition $\left( \protect\ref{eq.h}\right) $ implies
normalized condition\label{MCNC}}

\begin{proposition}
Let $B:H\otimes H\rightarrow A\otimes H^{op}.$ Then if $B$ satisfies $\left( %
\ref{eq.h}\right) $ of the morphism conditions and $\left( \ref{1ot1}\right)
$ then it satisfies normalized condition.
\end{proposition}

\begin{proof}
We apply $\mathrm{Id}_{A}\otimes \epsilon _{H}$ to $\left( \ref{eq.h}\right)
$ and obtain%
\begin{equation*}
\left( \mathrm{\mathrm{\mathrm{Id}}}\otimes \epsilon _{H}\right) \left[
B(h\otimes h^{\prime })(1_{A}\otimes h^{\prime \prime })\right] =\left(
\mathrm{\mathrm{\mathrm{Id}}}\otimes \epsilon _{H}\right) \left[ (1\otimes
h_{1}^{\prime \prime })B(h_{2}^{\prime \prime }h\otimes h_{3}^{\prime \prime
}h^{\prime })\right]
\end{equation*}

i.e.%
\begin{equation*}
\left( \mathrm{\mathrm{\mathrm{Id}}}\otimes \epsilon _{H}\right) B(h\otimes
h^{\prime })\left( 1_{A}\otimes \epsilon _{H}\left( h^{\prime \prime
}\right) \right) =B\left( h_{1}^{\prime \prime }h\otimes h_{2}^{\prime
\prime }h^{\prime }\right)
\end{equation*}%
and $h=h^{\prime }=1_{H}$%
\begin{equation*}
\left( \mathrm{\mathrm{\mathrm{Id}}}_{A}\otimes \epsilon _{H}\right)
B(1_{H}\otimes 1_{H})\left( 1_{A}\otimes \epsilon _{H}\left( h^{\prime
\prime }\right) \right) =B\left( h_{1}^{\prime \prime }\otimes h_{2}^{\prime
\prime }\right)
\end{equation*}%
i.e.%
\begin{equation*}
\left( 1_{A}\otimes \epsilon _{H}\left( h^{\prime \prime }\right) \right)
=B\left( h_{1}^{\prime \prime }\otimes h_{2}^{\prime \prime }\right)
\end{equation*}

This gives us the normalized condition.
\end{proof}

\section{Equation $B(h\otimes h^{\prime })(1_{A}\otimes h^{\prime \prime
})=(1_{A}\otimes h_{1}^{\prime \prime })B(h_{2}^{\prime \prime }h\otimes
h_{3}^{\prime \prime }h^{\prime })$}

As we remarked in section $\left( \ref{mor cond}\right) $, we only need to
compute equality $\left( \ref{eq.h}\right) $ for $1_{A}\otimes g$, $%
1_{A}\otimes x_{1}$ and $1_{A}\otimes x_{2}$.

i) For $1_{A}\otimes g$

\begin{equation}
B(h\otimes h^{\prime })\ \left( 1_{A}\otimes g\right) =(1_{A}\otimes
g)B(gh\otimes gh^{\prime })  \label{1g}
\end{equation}

ii) For $1_{A}\otimes x_{i}$

\begin{equation}
B(h\otimes h^{\prime })(1_{A}\otimes x_{i})=(1_{A}\otimes x_{i})B(gh\otimes
gh^{\prime })+B(x_{i}h\otimes gh^{\prime })+B(h\otimes x_{i}h^{\prime })
\label{1x}
\end{equation}

Moreover, we only need to compute equality $\left( \ref{eq.a}\right) $ for $%
G\otimes 1_{H}$ and $X_{i}\otimes 1_{H},$ $i=1,2$.

iii) For $G\otimes 1_{H}$

\begin{equation}
B(h\otimes h^{\prime })(G\otimes 1_{H})=(G\otimes 1_{H})B(hg\otimes
h^{\prime }g)  \label{G1}
\end{equation}

iv) For $X_{i}\otimes 1_{H}$

\begin{eqnarray}
&&B(h\otimes h^{\prime })(X_{i}\otimes 1_{H})  \label{Xi1} \\
&=&(X_{i}\otimes 1_{H})B(hg\otimes h^{\prime }g)+B(hx_{i}\otimes h^{\prime
}g)+B(h\otimes h^{\prime }x_{i})  \notag
\end{eqnarray}

\subsection{{\protect\Large Case \ }$h^{\prime \prime }=g$}

For $h^{\prime \prime }=g$ we get \ (\ref{1g}):
\begin{equation*}
B(h\otimes h^{\prime })(1_{A}\otimes g)=(1_{A}\otimes g)B(gh\otimes
gh^{\prime })
\end{equation*}

which is equivalent to

\begin{equation}
B(gh\otimes gh^{\prime })=(1_{A}\otimes g)B(h\otimes h^{\prime
})(1_{A}\otimes g)  \label{eq.10}
\end{equation}

Also equality $\left( \ref{eq.10}\right) $ rewrites as%
\begin{eqnarray*}
&&\sum_{p,q,r,s,t,u=0}^{1}B(gh\otimes gh^{\prime
};G^{p}X_{1}^{q}X_{2}^{r},g^{s}x_{1}^{t}x_{2}^{u})G^{p}X_{1}^{q}X_{2}^{r}%
\otimes g^{s}x_{1}^{t}x_{2}^{u} \\
&=&\sum_{p,q,r,s,t,u=0}^{1}B(h\otimes h^{\prime
};G^{p}X_{1}^{q}X_{2}^{r},g^{s}x_{1}^{t}x_{2}^{u})G^{p}X_{1}^{q}X_{2}^{r}%
\otimes g^{s+1}x_{1}^{t}x_{2}^{u}g. \\
&=&\sum_{p,q,r,s,t,u=0}^{1}(-1)^{t+u}B(h\otimes h^{\prime
};G^{p}X_{1}^{q}X_{2}^{r},g^{s}x_{1}^{t}x_{2}^{u})G^{p}X_{1}^{q}X_{2}^{r}%
\otimes g^{s}x_{1}^{t}x_{2}^{u}
\end{eqnarray*}%
Thus, equality (\ref{1g}) is equivalent to the following equalities.

\begin{eqnarray}
B(gh\otimes gh^{\prime };d,f) &=&B(h\otimes h^{\prime };d,f)\;\;\;\text{any }%
d\in A\text{ and }f=1_{H},g,x_{1}x_{2},gx_{1}x_{2}  \label{1gresult1} \\
B(gh\otimes gh^{\prime };d,f) &=&-B(h\otimes h^{\prime };d,f)\;\;\;\;\text{%
any }d\in A\text{ and }f=x_{i},gx_{i}  \label{1gresult2}
\end{eqnarray}

\subsection{{\protect\Large Case \ }$h^{\prime \prime }=x_{i}$}

We have%
\begin{equation*}
\left[ \left( \Delta \otimes \mathrm{Id}_{H}\right) \circ \Delta \right]
\left( x_{i}\right) =x_{i}\otimes g\otimes g+1_{H}\otimes x_{i}\otimes
g+1_{H}\otimes 1_{H}\otimes x_{i}
\end{equation*}%
\begin{equation*}
B(h\otimes h^{\prime })(1_{A}\otimes x_{i})=(1_{A}\otimes x_{i})B(gh\otimes
gh^{\prime })+(1_{A}\otimes 1_{H})B(x_{i}h\otimes gh^{\prime
})+(1_{A}\otimes 1_{H})B(h\otimes x_{i}h^{\prime })
\end{equation*}%
\begin{equation}
B(h\otimes h^{\prime })(1_{A}\otimes x_{i})=(1_{A}\otimes x_{i})B(gh\otimes
gh^{\prime })+B(x_{i}h\otimes gh^{\prime })+B(h\otimes x_{i}h^{\prime })
\label{eq:hh'}
\end{equation}%
Using $\left( \ref{eq.10}\right) ,$ we can rewrite $\left( \ref{eq:hh'}%
\right) $ in the form%
\begin{equation}
B(h\otimes h^{\prime })(1_{A}\otimes x_{i})=(1_{A}\otimes gx_{i})B(h\otimes
h^{\prime })(1_{A}\otimes g)+B(x_{i}h\otimes gh^{\prime })+B(h\otimes
x_{i}h^{\prime })  \label{eq:hh'g}
\end{equation}

Now in $\left( \ref{eq:hh'}\right) $ we consider $h=g,h^{\prime }=1_{H.}$

\begin{equation*}
B(g\otimes 1_{H})(1_{A}\otimes x_{i})=(1_{A}\otimes x_{i})B(1_{H}\otimes
g)+B(x_{i}g\otimes g)+B(g\otimes x_{i})
\end{equation*}%
\begin{equation*}
B(g\otimes 1_{H})(1_{A}\otimes x_{i})=(1_{A}\otimes x_{i})\left(
1_{H}\otimes g\right) B(g\otimes 1_{H})\left( 1_{H}\otimes g\right) +\left(
1_{H}\otimes g\right) B(x_{i}\otimes 1_{H})\left( 1_{H}\otimes g\right)
+B(g\otimes x_{i})
\end{equation*}

Moreover from $\left( \ref{eq:hh'}\right) ,$using $\left( \ref{eq.10}\right)
$ we obtain%
\begin{equation*}
B(h\otimes h^{\prime })(1_{A}\otimes x_{i})=(1_{A}\otimes gx_{i})B(h\otimes
h^{\prime })(1_{A}\otimes g)+B(x_{i}h\otimes gh^{\prime })+B(h\otimes
x_{i}h^{\prime })
\end{equation*}%
\begin{equation*}
B(g\otimes x_{i})=B(g\otimes 1_{H})(1_{A}\otimes x_{i})+\left( 1_{H}\otimes
gx_{i}\right) B(g\otimes 1_{H})\left( 1_{H}\otimes g\right) -\left(
1_{H}\otimes g\right) B(x_{i}\otimes 1_{H})\left( 1_{H}\otimes g\right)
\end{equation*}

\begin{lemma}
\label{Lem:Cas}Assume that equality $\left( \ref{eq.10}\right) $ holds. Then
Casimir condition for $B\left( gh\otimes gh^{\prime }\right) $ is equivalent
to the equality
\begin{eqnarray}
&&B^{A}(h_{2}\otimes h^{\prime })_{0}\otimes gB^{H}(h_{2}\otimes h^{\prime
})_{1}g\otimes gB^{H}(h_{2}\otimes h^{\prime })_{2}h_{1}B^{A}(h_{2}\otimes
h^{\prime })_{1}  \label{eq.03} \\
&=&B^{A}(h\otimes h_{1}^{\prime })\otimes gB^{H}(h\otimes h_{1}^{\prime
})g\otimes gh_{2}^{\prime }  \notag
\end{eqnarray}
\end{lemma}

\begin{proof}
First of all recall $\left( \ref{eq.10}\right) $%
\begin{equation*}
B(gh\otimes gh^{\prime })=(1_{A}\otimes g)B(h\otimes h^{\prime
})(1_{A}\otimes g)
\end{equation*}%
i.e.%
\begin{equation}
B^{A}(gh\otimes gh^{\prime })\otimes B^{H}(gh\otimes gh^{\prime
})=B^{A}(h\otimes h^{\prime })\otimes gB^{H}(h\otimes h^{\prime })g
\label{eq. 10bis}
\end{equation}%
Now we apply
\begin{equation*}
\rho _{A\otimes H^{op}}:A\otimes H^{op}\rightarrow A\otimes H^{op}\otimes
A\otimes H^{op},
\end{equation*}%
\begin{equation*}
\rho _{A\otimes H^{op}}(a\otimes h)=(a_{0}\otimes h_{1})\otimes
(a_{1}\otimes h_{2}).
\end{equation*}%
to both sides of $\left( \ref{eq. 10bis}\right) $ to get%
\begin{eqnarray*}
&&B^{A}(gh\otimes gh^{\prime })_{0}\otimes B^{H}(gh\otimes gh^{\prime
})_{1}\otimes B^{A}(gh\otimes gh^{\prime })_{1}\otimes B^{H}(gh\otimes
gh^{\prime })_{2} \\
&=&B^{A}(h\otimes h^{\prime })_{0}\otimes \left( gB^{H}(h\otimes h^{\prime
})g\right) _{1}\otimes B^{A}(h\otimes h^{\prime })_{1}\otimes \left(
gB^{H}(h\otimes h^{\prime })g\right) _{2}
\end{eqnarray*}%
i.e.%
\begin{eqnarray*}
&&B^{A}(gh\otimes gh^{\prime })_{0}\otimes B^{H}(gh\otimes gh^{\prime
})_{1}\otimes B^{A}(gh\otimes gh^{\prime })_{1}\otimes B^{H}(gh\otimes
gh^{\prime })_{2} \\
&=&B^{A}(h\otimes h^{\prime })_{0}\otimes gB^{H}(h\otimes h^{\prime
})_{1}g\otimes B^{A}(h\otimes h^{\prime })_{1}\otimes gB^{H}(h\otimes
h^{\prime })_{2}g
\end{eqnarray*}%
For $h=h_{2}$ we get%
\begin{eqnarray}
&&B^{A}(gh_{2}\otimes gh^{\prime })_{0}\otimes B^{H}(gh_{2}\otimes
gh^{\prime })_{1}\otimes B^{A}(gh_{2}\otimes gh^{\prime })_{1}\otimes
B^{H}(gh_{2}\otimes gh^{\prime })_{2}  \label{eq.00} \\
&=&B^{A}(h_{2}\otimes h^{\prime })_{0}\otimes gB^{H}(h_{2}\otimes h^{\prime
})_{1}g\otimes B^{A}(h_{2}\otimes h^{\prime })_{1}\otimes
gB^{H}(h_{2}\otimes h^{\prime })_{2}g  \notag
\end{eqnarray}

We write the Casimir condition
\begin{equation*}
B^{A}(h_{2}\otimes h^{\prime })_{0}\otimes B^{H}(h_{2}\otimes h^{\prime
})_{1}\otimes B^{H}(h_{2}\otimes h^{\prime })_{2}h_{1}B^{A}(h_{2}\otimes
h^{\prime })_{1}=B^{A}(h\otimes h_{1}^{\prime })\otimes B^{H}(h\otimes
h_{1}^{\prime })\otimes h_{2}^{\prime }
\end{equation*}%
for $B\left( gh\otimes gh^{\prime }\right) .$We get%
\begin{eqnarray*}
&&B^{A}(\left( gh\right) _{2}\otimes gh^{\prime })_{0}\otimes B^{H}(\left(
gh\right) _{2}\otimes gh^{\prime })_{1}\otimes B^{H}(\left( gh\right)
_{2}\otimes gh^{\prime })_{2}\left( gh\right) _{1}B^{A}(\left( gh\right)
_{2}\otimes gh^{\prime })_{1} \\
&=&B^{A}(gh\otimes \left( gh^{\prime }\right) _{1})\otimes B^{H}(gh\otimes
\left( gh^{\prime }\right) _{1})\otimes \left( gh^{\prime }\right) _{2}
\end{eqnarray*}%
i.e.%
\begin{eqnarray}
&&B^{A}(gh_{2}\otimes gh^{\prime })_{0}\otimes B^{H}(gh_{2}\otimes
gh^{\prime })_{1}\otimes B^{H}(gh_{2}\otimes gh^{\prime
})_{2}gh_{1}B^{A}(gh_{2}\otimes gh^{\prime })_{1}  \label{eq.01} \\
&=&B^{A}(gh\otimes gh_{1}^{\prime })\otimes B^{H}(gh\otimes gh_{1}^{\prime
})\otimes gh_{2}^{\prime }.  \notag
\end{eqnarray}

Now we have%
\begin{eqnarray*}
&&B^{A}(h\otimes h_{1}^{\prime })\otimes gB^{H}(h\otimes h_{1}^{\prime
})g\otimes gh_{2}^{\prime } \\
\overset{\left( \ref{eq. 10bis}\right) }{=} &&B^{A}(gh\otimes gh_{1}^{\prime
})\otimes B^{H}(gh\otimes gh_{1}^{\prime })\otimes gh_{2}^{\prime } \\
&&\overset{\left( \ref{eq.01}\right) }{=}B^{A}(gh_{2}\otimes gh^{\prime
})_{0}\otimes B^{H}(gh_{2}\otimes gh^{\prime })_{1}\otimes
B^{H}(gh_{2}\otimes gh^{\prime })_{2}gh_{1}B^{A}(gh_{2}\otimes gh^{\prime
})_{1} \\
&&\overset{\left( \ref{eq.00}\right) }{=}B^{A}(h_{2}\otimes h^{\prime
})_{0}\otimes gB^{H}(h_{2}\otimes h^{\prime })_{1}g\otimes
gB^{H}(h_{2}\otimes h^{\prime })_{2}h_{1}B^{A}(h_{2}\otimes h^{\prime })_{1}
\end{eqnarray*}

so that $\left( \ref{eq.03}\right) $ holds.

Conversely assume that $\left( \ref{eq.03}\right) $holds:%
\begin{eqnarray*}
&&B^{A}(h_{2}\otimes h^{\prime })_{0}\otimes gB^{H}(h_{2}\otimes h^{\prime
})_{1}g\otimes gB^{H}(h_{2}\otimes h^{\prime })_{2}h_{1}B^{A}(h_{2}\otimes
h^{\prime })_{1} \\
&=&B^{A}(h\otimes h_{1}^{\prime })\otimes gB^{H}(h\otimes h_{1}^{\prime
})g\otimes gh_{2}^{\prime }
\end{eqnarray*}%
We have%
\begin{eqnarray*}
&&B^{A}(gh\otimes gh_{1}^{\prime })\otimes B^{H}(gh\otimes gh_{1}^{\prime
})\otimes gh_{2}^{\prime } \\
&&\overset{\left( \ref{eq. 10bis}\right) }{=}B^{A}(h\otimes h_{1}^{\prime
})\otimes gB^{H}(h\otimes h_{1}^{\prime })g\otimes gh_{2}^{\prime } \\
&&\overset{\left( \ref{eq.03}\right) }{=}B^{A}(h_{2}\otimes h^{\prime
})_{0}\otimes gB^{H}(h_{2}\otimes h^{\prime })_{1}g\otimes
gB^{H}(h_{2}\otimes h^{\prime })_{2}h_{1}B^{A}(h_{2}\otimes h^{\prime })_{1}
\\
&&\overset{\left( \ref{eq.00}\right) }{=}B^{A}(gh_{2}\otimes gh^{\prime
})_{0}\otimes B^{H}(gh_{2}\otimes gh^{\prime })_{1}\otimes
B^{H}(gh_{2}\otimes gh^{\prime })_{2}gh_{1}\otimes B^{A}(gh_{2}\otimes
gh^{\prime })_{1}
\end{eqnarray*}%
and hence equality $\left( \ref{eq.01}\right) $ holds.
\end{proof}

\begin{proposition}
\label{Pro:gh}Assume that equality $\left( \ref{eq.10}\right) $ holds.
Assume also that, for some $h,h^{\prime }\in H$, $B\left( h\otimes h^{\prime
}\right) ,$ satisfies Casimir condition. Then also $B\left( gh\otimes
gh^{\prime }\right) $ does.
\end{proposition}

\begin{proof}
We consider Casimir condition%
\begin{equation*}
B^{A}(h_{2}\otimes h^{\prime })_{0}\otimes B^{H}(h_{2}\otimes h^{\prime
})_{1}\otimes B^{H}(h_{2}\otimes h^{\prime })_{2}h_{1}B^{A}(h_{2}\otimes
h^{\prime })_{1}=B^{A}(h\otimes h_{1}^{\prime })\otimes B^{H}(h\otimes
h_{1}^{\prime })\otimes h_{2}^{\prime }
\end{equation*}%
and multiply both sides on the left by$\left( 1_{A}\otimes g\otimes
1_{H}\right) $ and on the right by $\left( 1_{A}\otimes g\otimes g\right) $%
to get%
\begin{equation*}
B^{A}(h_{2}\otimes h^{\prime })_{0}\otimes gB^{H}(h_{2}\otimes h^{\prime
})_{1}g\otimes gB^{H}(h_{2}\otimes h^{\prime })_{2}h_{1}B^{A}(h_{2}\otimes
h^{\prime })_{1}=B^{A}(h\otimes h_{1}^{\prime })\otimes gB^{H}(h\otimes
h_{1}^{\prime })g\otimes gh_{2}^{\prime }
\end{equation*}%
which is $\left( \ref{eq.03}\right) .$Thus, by Lemma $\left( \ref{Lem:Cas}%
\right) $ we got Casimir condition for $B\left( gh\otimes gh^{\prime
}\right) $ from Casimir condition for $B\left( h\otimes h^{\prime }\right) .$
\end{proof}

\begin{corollary}
If \ $B\left( gx_{i}\otimes 1_{H}\right) $ satisfies Casimir condition, also
$B(1_{H}\otimes x_{i})$ does
\end{corollary}

\begin{proof}
By using $\left( \ref{eq:hh'g}\right) $ with $h=h^{\prime }=1_{H},$ we get%
\begin{equation*}
B(1_{H}\otimes 1_{H})(1_{A}\otimes x_{i})=(1_{A}\otimes x_{i})B(g\otimes
g^{\prime })+B(x_{i}\otimes g)+B(1_{H}\otimes x_{i})
\end{equation*}%
\begin{equation*}
(1_{A}\otimes x_{i})=(1_{A}\otimes x_{i})+B(x_{i}\otimes g)+B(1_{H}\otimes
x_{i})
\end{equation*}%
\begin{equation*}
B(1_{H}\otimes x_{i})=-B(x_{i}\otimes g)=-\left( 1\otimes g\right) B\left(
gx_{i}\otimes 1_{H}\right) \left( 1\otimes g\right)
\end{equation*}%
Apply Proposition \ref{Pro:gh}.
\end{proof}

\section{Simplification}

\begin{proposition}
\label{Pro helem}Assume that $B$ satisfies $\left( \ref{eq.h}\right) .$ Then
all the elements of the form $B\left( h\otimes h^{\prime }\right) $ can be
obtained from elements of the form $B\left( h\otimes 1_{H}\right) .$In
particular we get%
\begin{equation}
B(h\otimes x_{1}x_{2})=B(h\otimes x_{2})(1_{A}\otimes x_{1})-(1_{A}\otimes
gx_{1})B(h\otimes x_{2})(1_{A}\otimes g)-(1_{A}\otimes g)B(gx_{1}h\otimes
x_{2})(1_{A}\otimes g)  \label{simpl1}
\end{equation}%
\begin{equation}
B(h\otimes x_{i})=B(h\otimes 1_{H})(1_{A}\otimes x_{i})-(1_{A}\otimes
gx_{i})B(h\otimes 1_{H})(1_{A}\otimes g)-(1_{A}\otimes g)B(gx_{i}h\otimes
1_{H})(1_{A}\otimes g)  \label{simplx}
\end{equation}
\end{proposition}

\begin{proof}
We have that, by $\left( \ref{eq.10}\right) $,%
\begin{equation*}
B\left( h\otimes g\right) =(1_{A}\otimes g)B\left( gh\otimes 1_{H}\right)
(1_{A}\otimes g)
\end{equation*}%
Now we recall that%
\begin{equation*}
B(h\otimes h^{\prime })(1_{A}\otimes x_{i})\overset{\left( \ref{1x}\right) }{%
=}(1_{A}\otimes x_{i})B(gh\otimes gh^{\prime })+B(x_{i}h\otimes gh^{\prime
})+B(h\otimes x_{i}h^{\prime })
\end{equation*}%
so that we get%
\begin{equation}
B(h\otimes h^{\prime })(1_{A}\otimes x_{i})\overset{\left( \ref{eq.10}%
\right) }{=}(1_{A}\otimes x_{i})(1_{A}\otimes g)B(h\otimes h^{\prime
})(1_{A}\otimes g)+(1_{A}\otimes g)B(gx_{i}h\otimes h^{\prime
})(1_{A}\otimes g)+B(h\otimes x_{i}h^{\prime }).  \label{1xrew}
\end{equation}

For $h^{\prime }=1_{H},$ we get%
\begin{equation*}
B(h\otimes 1_{H})(1_{A}\otimes x_{i})=(1_{A}\otimes gx_{i})B(h\otimes
1_{H})(1_{A}\otimes g)+(1_{A}\otimes g)B(gx_{i}h\otimes 1_{H})(1_{A}\otimes
g)+B(h\otimes x_{i})
\end{equation*}%
Thus $B(h\otimes x_{i})$ depends only on terms of the form $B(h\otimes
1_{H}).$

For $h^{\prime }=x_{j}$ $j\neq i$%
\begin{equation*}
B(h\otimes x_{j})(1_{A}\otimes x_{i})=(1_{A}\otimes gx_{i})B(h\otimes
x_{j})(1_{A}\otimes g)+(1_{A}\otimes g)B(gx_{i}h\otimes x_{j})(1_{A}\otimes
g)+B(h\otimes x_{i}x_{j})
\end{equation*}%
and hence, by using the previous argument, also $B(h\otimes x_{i}x_{j})$
depends only on terms of the form $B(h\otimes 1_{H}).$

Also for $h^{\prime }=gx_{j}$ $j\neq i$%
\begin{equation*}
B(h\otimes gx_{j})(1_{A}\otimes x_{i})=-(1_{A}\otimes gx_{i})B(h\otimes
gx_{j})(1_{A}\otimes g)+(1_{A}\otimes g)B(gx_{i}h\otimes
gx_{j})(1_{A}\otimes g)+B(h\otimes x_{i}gx_{j})
\end{equation*}%
and hence, using $\left( \ref{eq.10}\right) ,$ we get%
\begin{equation*}
(1_{A}\otimes g)B(gh\otimes x_{j})(1_{A}\otimes g)(1_{A}\otimes
x_{i})=(1_{A}\otimes gx_{i})(1_{A}\otimes g)B(gh\otimes
x_{j})+B(x_{i}h\otimes x_{j})(1_{A}\otimes g)+B(h\otimes x_{i}gx_{j})
\end{equation*}%
Thus, by using the previous argument, also $B(h\otimes gx_{i}x_{j})$ this
depends only on terms of the form $B(h\otimes 1_{H}).$

In conclusion all the terms of the form $B\left( h\otimes h^{\prime }\right)
$\ can be written on elements of the form{\LARGE \ }$B\left( h\otimes
1_{H}\right) .$
\end{proof}

\begin{remark}
For $i=2$ in $\left( \ref{simplx}\right) $ we get%
\begin{equation}
B(h\otimes x_{2})=B(h\otimes 1_{H})(1_{A}\otimes x_{2})-(1_{A}\otimes
gx_{2})B(h\otimes 1_{H})(1_{A}\otimes g)-(1_{A}\otimes g)B(gx_{2}h\otimes
1_{H})(1_{A}\otimes g)  \label{simpl2i2}
\end{equation}%
so that $\left( \ref{simpl1}\right) $ becomes%
\begin{eqnarray*}
B(h\otimes x_{1}x_{2}) &=&B(h\otimes 1_{H})(1_{A}\otimes x_{2})(1_{A}\otimes
x_{1}) \\
&&-(1_{A}\otimes gx_{2})B(h\otimes 1_{H})(1_{A}\otimes g)(1_{A}\otimes x_{1})
\\
&&-(1_{A}\otimes g)B(gx_{2}h\otimes 1_{H})(1_{A}\otimes g)(1_{A}\otimes
x_{1}) \\
&&-(1_{A}\otimes gx_{1})B(h\otimes 1_{H})(1_{A}\otimes x_{2})(1_{A}\otimes g)
\\
&&+(1_{A}\otimes gx_{1})(1_{A}\otimes gx_{2})B(h\otimes 1_{H})(1_{A}\otimes
g)(1_{A}\otimes g) \\
&&+(1_{A}\otimes gx_{1})(1_{A}\otimes g)B(gx_{2}h\otimes 1_{H})(1_{A}\otimes
g)(1_{A}\otimes g) \\
&&-(1_{A}\otimes g)B(gx_{1}h\otimes 1_{H})(1_{A}\otimes x_{2})(1_{A}\otimes
g) \\
&&+(1_{A}\otimes g)(1_{A}\otimes gx_{2})B(gx_{1}h\otimes 1_{H})(1_{A}\otimes
g)(1_{A}\otimes g) \\
&&+(1_{A}\otimes g)(1_{A}\otimes g)B(gx_{2}gx_{1}h\otimes
1_{H})(1_{A}\otimes g)(1_{A}\otimes g)
\end{eqnarray*}%
i.e.%
\begin{eqnarray}
B(h\otimes x_{1}x_{2}) &=&B(h\otimes 1_{H})(1_{A}\otimes x_{1}x_{2})
\label{simplxx} \\
&&+(1_{A}\otimes gx_{2})B(h\otimes 1_{H})(1_{A}\otimes gx_{1})  \notag \\
&&+(1_{A}\otimes g)B(gx_{2}h\otimes 1_{H})(1_{A}\otimes gx_{1})  \notag \\
&&-(1_{A}\otimes gx_{1})B(h\otimes 1_{H})(1_{A}\otimes gx_{2})  \notag \\
&&+(1_{A}\otimes x_{1}x_{2})B(h\otimes 1_{H})  \notag \\
&&+(1_{A}\otimes x_{1})B(gx_{2}h\otimes 1_{H})  \notag \\
&&-(1_{A}\otimes g)B(gx_{1}h\otimes 1_{H})(1_{A}\otimes gx_{2})  \notag \\
&&-(1_{A}\otimes x_{2})B(gx_{1}h\otimes 1_{H})  \notag \\
&&+B(x_{1}x_{2}h\otimes 1_{H})  \notag
\end{eqnarray}%
From $\left( \ref{eq.h}\right) $%
\begin{equation*}
B(h\otimes h^{\prime })(1_{A}\otimes h^{\prime \prime })=(1_{A}\otimes
h_{1}^{\prime \prime })B(h_{2}^{\prime \prime }h\otimes h_{3}^{\prime \prime
}h^{\prime })
\end{equation*}%
for $h^{\prime \prime }=g$ we get and $h=gt$%
\begin{equation*}
B(gt\otimes h^{\prime })(1_{A}\otimes g)=(1_{A}\otimes g)B(ggt\otimes
gh^{\prime })
\end{equation*}%
\begin{equation*}
B(gt\otimes h^{\prime })(1_{A}\otimes g)=(1_{A}\otimes g)B(t\otimes
gh^{\prime })
\end{equation*}%
and hence%
\begin{equation}
B(h\otimes gh^{\prime })=(1_{A}\otimes g)B(gh\otimes h^{\prime
})(1_{A}\otimes g)  \label{simplgdo}
\end{equation}%
Thus for $h^{\prime }=x_{i}$ in $\left( \ref{simplgdo}\right) $%
\begin{equation*}
B(h\otimes gx_{i})=(1_{A}\otimes g)B(gh\otimes x_{i})(1_{A}\otimes g)
\end{equation*}%
and hence using $\left( \ref{simplx}\right) $

we get$\left( \ref{simplgx}\right) $%
\begin{eqnarray}
&&B(h\otimes gx_{i})  \label{simplgx} \\
&=&(1_{A}\otimes g)B(gh\otimes 1_{H})(1_{A}\otimes gx_{i})  \notag \\
&&+(1_{A}\otimes x_{i})B(gh\otimes 1_{H})  \notag \\
&&-B(gx_{i}gh\otimes 1_{H})  \notag
\end{eqnarray}%
Thus for $h^{\prime }=x_{1}x_{2}$ in $\left( \ref{simplgdo}\right) $%
\begin{equation*}
B(h\otimes gx_{1}x_{2})=(1_{A}\otimes g)B(gh\otimes x_{1}x_{2})(1_{A}\otimes
g)
\end{equation*}%
and by using $\ref{simplxx}$%
\begin{eqnarray}
B(h\otimes gx_{1}x_{2}) &=&(1_{A}\otimes g)B(gh\otimes 1_{H})(1_{A}\otimes
gx_{1}x_{2})  \label{simplgxx} \\
&&-(1_{A}\otimes x_{2})B(gh\otimes 1_{H})(1_{A}\otimes x_{1})  \notag \\
&&-B(x_{2}h\otimes 1_{H})(1_{A}\otimes x_{1})  \notag \\
&&+(1_{A}\otimes x_{1})B(gh\otimes 1_{H})(1_{A}\otimes x_{2})  \notag \\
&&+(1_{A}\otimes gx_{1}x_{2})B(gh\otimes 1_{H})(1_{A}\otimes g)  \notag \\
&&+(1_{A}\otimes gx_{1})B(x_{2}h\otimes 1_{H})(1_{A}\otimes g)  \notag \\
&&+B(x_{1}h\otimes 1_{H})(1_{A}\otimes x_{2})  \notag \\
&&-(1_{A}\otimes gx_{2})B(x_{1}h\otimes 1_{H})(1_{A}\otimes g)  \notag \\
&&+(1_{A}\otimes g)B(gx_{1}x_{2}h\otimes 1_{H})(1_{A}\otimes g)  \notag
\end{eqnarray}
\end{remark}

\section{REDUCTION of equation$\left( \protect\ref{eq.a}\right) \label{REa}$}

We want to show that, assuming that equation $\left( \ref{eq.h}\right) $ is
satisfied, if equation $\left( \ref{eq.a}\right) $ holds for $h\otimes 1_{H}$
then it holds for any $h\otimes h^{\prime }.$

\subsection{CASE\ $h^{\prime }=x_{i}$}

We want to write the left side of the formula $\left( \ref{eq.a}\right) $
for $h^{\prime }=x_{i}.$
\begin{equation*}
B(h\otimes x_{i})(d\otimes 1_{H})=(d_{0}\otimes 1_{H})B(hd_{1}\otimes
x_{i}d_{2})
\end{equation*}

We consider the left side by using $\left( \ref{simplx}\right) $
\begin{equation*}
B(h\otimes x_{i})=B(h\otimes 1_{H})(1_{A}\otimes x_{i})-(1_{A}\otimes
gx_{i})B(h\otimes 1_{H})(1_{A}\otimes g)-(1_{A}\otimes g)B(gx_{i}h\otimes
1_{H})(1_{A}\otimes g)
\end{equation*}%
We write the left side%
\begin{eqnarray*}
B(h\otimes x_{i})(d\otimes 1_{H}) &=&B(h\otimes 1_{H})(1_{A}\otimes
x_{i})(d\otimes 1_{H})-(1_{A}\otimes gx_{i})B(h\otimes 1_{H})(1_{A}\otimes
g)(d\otimes 1_{H}) \\
&&-(1_{A}\otimes g)B(gx_{i}h\otimes 1_{H})(1_{A}\otimes g)(d\otimes 1_{H})
\end{eqnarray*}%
\begin{eqnarray*}
B(h\otimes x_{i})(d\otimes 1_{H}) &=&B(h\otimes 1_{H})(d\otimes
1_{H})(1_{A}\otimes x_{i})-(1_{A}\otimes gx_{i})B(h\otimes 1_{H})(d\otimes
1_{H})(1_{A}\otimes g) \\
&&-(1_{A}\otimes g)B(gx_{i}h\otimes 1_{H})(d\otimes 1_{H})(1_{A}\otimes g)
\end{eqnarray*}%
We use here the formula for $h\otimes 1_{H}$%
\begin{eqnarray*}
B(h\otimes x_{i})(d\otimes 1_{H}) &=&(d_{0}\otimes 1_{H})B(hd_{1}\otimes
d_{2})(1_{A}\otimes x_{i})-(1_{A}\otimes gx_{i})(d_{0}\otimes
1_{H})B(hd_{1}\otimes d_{2})(1_{A}\otimes g) \\
&&-(1_{A}\otimes g)(d_{0}\otimes 1_{H})B(gx_{i}hd_{1}\otimes
d_{2})(1_{A}\otimes g)
\end{eqnarray*}%
Thus the left side of formula $\left( \ref{eq.a}\right) $for $h^{\prime
}=x_{i}$ is%
\begin{equation}
B(h\otimes x_{i})(d\otimes 1_{H})=(d_{0}\otimes 1_{H})\left[
\begin{array}{c}
B(hd_{1}\otimes d_{2})(1_{A}\otimes x_{i})-(1_{A}\otimes
gx_{i})B(hd_{1}\otimes d_{2})(1_{A}\otimes g) \\
-(1_{A}\otimes g)B(gx_{i}hd_{1}\otimes d_{2})(1_{A}\otimes g)%
\end{array}%
\right] .  \label{form axi}
\end{equation}

\subsection{Case $h^{\prime }=x_{1}x_{2}$}

We want to write the left side of the formula $\left( \ref{eq.a}\right) $
for $h^{\prime }=x_{1}x_{2}$
\begin{equation*}
B(h\otimes x_{1}x_{2})(d\otimes 1_{H})=(d_{0}\otimes 1_{H})B(hd_{1}\otimes
x_{1}x_{2}d_{2})
\end{equation*}

We consider the left side by using $\left( \ref{simplxx}\right) $

\begin{equation*}
B(h\otimes x_{1}x_{2})=B(h\otimes x_{2})(1_{A}\otimes x_{1})-(1_{A}\otimes
gx_{1})B(h\otimes x_{2})(1_{A}\otimes g)-(1_{A}\otimes g)B(gx_{1}h\otimes
x_{2})(1_{A}\otimes g)
\end{equation*}

We write the left side%
\begin{eqnarray*}
&&B(h\otimes x_{1}x_{2})(d\otimes 1_{H}) \\
&=&B(h\otimes x_{2})(1_{A}\otimes x_{1})(d\otimes 1_{H})-(1_{A}\otimes
gx_{1})B(h\otimes x_{2})(1_{A}\otimes g)(d\otimes 1_{H}) \\
&&-(1_{A}\otimes g)B(gx_{1}h\otimes x_{2})(1_{A}\otimes g)(d\otimes 1_{H})
\end{eqnarray*}%
\begin{eqnarray*}
&&B(h\otimes x_{1}x_{2})(d\otimes 1_{H}) \\
&=&B(h\otimes x_{2})(d\otimes 1_{H})(1_{A}\otimes x_{1})-(1_{A}\otimes
gx_{1})B(h\otimes x_{2})(d\otimes 1_{H})(1_{A}\otimes g) \\
&&-(1_{A}\otimes g)B(gx_{1}h\otimes x_{2})(d\otimes 1_{H})(1_{A}\otimes g)
\end{eqnarray*}%
We use here the formula $\left( \ref{form axi}\right) $
\begin{eqnarray*}
B(h\otimes x_{2})(d\otimes 1_{H}) &=&(d_{0}\otimes 1_{H})B(hd_{1}\otimes
d_{2})(1_{A}\otimes x_{2}) \\
&&-(1_{A}\otimes gx_{2})(d_{0}\otimes 1_{H})B(hd_{1}\otimes
d_{2})(1_{A}\otimes g) \\
&&-(1_{A}\otimes g)(d_{0}\otimes 1_{H})B(gx_{2}hd_{1}\otimes
d_{2})(1_{A}\otimes g)
\end{eqnarray*}%
and the formula for $gx_{1}h\otimes 1_{H}$ in $h\otimes x_{2}$
\begin{eqnarray*}
B(gx_{1}h\otimes x_{2})(d\otimes 1_{H}) &=&(d_{0}\otimes
1_{H})B(gx_{1}hd_{1}\otimes d_{2})(1_{A}\otimes x_{2}) \\
&&-(1_{A}\otimes gx_{2})(d_{0}\otimes 1_{H})B(gx_{1}hd_{1}\otimes
d_{2})(1_{A}\otimes g) \\
&&-(1_{A}\otimes g)(d_{0}\otimes 1_{H})B(gx_{2}gx_{1}hd_{1}\otimes
d_{2})(1_{A}\otimes g)
\end{eqnarray*}%
then we substitute%
\begin{eqnarray*}
&&B(h\otimes x_{1}x_{2})(d\otimes 1_{H}) \\
&=&\left[
\begin{array}{c}
(d_{0}\otimes 1_{H})B(hd_{1}\otimes d_{2})(1_{A}\otimes x_{2}) \\
-(1_{A}\otimes gx_{2})(d_{0}\otimes 1_{H})B(hd_{1}\otimes
d_{2})(1_{A}\otimes g) \\
-(1_{A}\otimes g)(d_{0}\otimes 1_{H})B(gx_{2}hd_{1}\otimes
d_{2})(1_{A}\otimes g)%
\end{array}%
\right] (1_{A}\otimes x_{1}) \\
&&-(1_{A}\otimes gx_{1})\left[
\begin{array}{c}
(d_{0}\otimes 1_{H})B(hd_{1}\otimes d_{2})(1_{A}\otimes x_{2}) \\
-(1_{A}\otimes gx_{2})(d_{0}\otimes 1_{H})B(hd_{1}\otimes
d_{2})(1_{A}\otimes g) \\
-(1_{A}\otimes g)(d_{0}\otimes 1_{H})B(gx_{2}hd_{1}\otimes
d_{2})(1_{A}\otimes g)%
\end{array}%
\right] (1_{A}\otimes g) \\
&&-(1_{A}\otimes g)\left[
\begin{array}{c}
(d_{0}\otimes 1_{H})B(gx_{1}hd_{1}\otimes d_{2})(1_{A}\otimes x_{2}) \\
-(1_{A}\otimes gx_{2})(d_{0}\otimes 1_{H})B(gx_{1}hd_{1}\otimes
d_{2})(1_{A}\otimes g) \\
-(1_{A}\otimes g)(d_{0}\otimes 1_{H})B(gx_{2}gx_{1}hd_{1}\otimes
d_{2})(1_{A}\otimes g)%
\end{array}%
\right] (1_{A}\otimes g)
\end{eqnarray*}%
We compute%
\begin{eqnarray*}
&&B(h\otimes x_{1}x_{2})(d\otimes 1_{H}) \\
&=&\left[
\begin{array}{c}
(d_{0}\otimes 1_{H})B(hd_{1}\otimes d_{2})(1_{A}\otimes x_{2})(1_{A}\otimes
x_{1}) \\
-(1_{A}\otimes gx_{i})(d_{0}\otimes 1_{H})B(hd_{1}\otimes
d_{2})(1_{A}\otimes g)(1_{A}\otimes x_{1}) \\
-(1_{A}\otimes g)(d_{0}\otimes 1_{H})B(gx_{2}hd_{1}\otimes
d_{2})(1_{A}\otimes g)(1_{A}\otimes x_{1})%
\end{array}%
\right] \\
&&-\left[
\begin{array}{c}
(1_{A}\otimes gx_{1})(d_{0}\otimes 1_{H})B(hd_{1}\otimes d_{2})(1_{A}\otimes
x_{2})(1_{A}\otimes g) \\
-(1_{A}\otimes gx_{1})(1_{A}\otimes gx_{2})(d_{0}\otimes
1_{H})B(hd_{1}\otimes d_{2})(1_{A}\otimes g)(1_{A}\otimes g) \\
-(1_{A}\otimes gx_{1})(1_{A}\otimes g)(d_{0}\otimes
1_{H})B(gx_{2}hd_{1}\otimes d_{2})(1_{A}\otimes g)(1_{A}\otimes g)%
\end{array}%
\right] \\
&&-\left[
\begin{array}{c}
(1_{A}\otimes g)(d_{0}\otimes 1_{H})B(gx_{1}hd_{1}\otimes
d_{2})(1_{A}\otimes x_{2})(1_{A}\otimes g) \\
-(1_{A}\otimes g)(1_{A}\otimes gx_{2})(d_{0}\otimes
1_{H})B(gx_{1}hd_{1}\otimes d_{2})(1_{A}\otimes g)(1_{A}\otimes g) \\
-(1_{A}\otimes g)(1_{A}\otimes g)(d_{0}\otimes
1_{H})B(gx_{2}gx_{1}hd_{1}\otimes d_{2})(1_{A}\otimes g)(1_{A}\otimes g)%
\end{array}%
\right]
\end{eqnarray*}%
We simplify%
\begin{eqnarray}
&&B(h\otimes x_{1}x_{2})(d\otimes 1_{H})  \label{form ax1x2} \\
&=&\left[
\begin{array}{c}
(d_{0}\otimes 1_{H})B(hd_{1}\otimes d_{2})(1_{A}\otimes x_{1}x_{2}) \\
+(1_{A}\otimes gx_{2})(d_{0}\otimes 1_{H})B(hd_{1}\otimes
d_{2})(1_{A}\otimes gx_{1}) \\
+(1_{A}\otimes g)(d_{0}\otimes 1_{H})B(gx_{2}hd_{1}\otimes
d_{2})(1_{A}\otimes gx_{1})%
\end{array}%
\right]  \notag \\
&&-\left[
\begin{array}{c}
(1_{A}\otimes gx_{1})(d_{0}\otimes 1_{H})B(hd_{1}\otimes d_{2})(1_{A}\otimes
gx_{2}) \\
-(1_{A}\otimes x_{1}x_{2})(d_{0}\otimes 1_{H})B(hd_{1}\otimes d_{2}) \\
-(1_{A}\otimes x_{1})(d_{0}\otimes 1_{H})B(gx_{2}hd_{1}\otimes d_{2})%
\end{array}%
\right]  \notag \\
&&-\left[
\begin{array}{c}
(1_{A}\otimes g)(d_{0}\otimes 1_{H})B(gx_{1}hd_{1}\otimes
d_{2})(1_{A}\otimes gx_{2}) \\
+(1_{A}\otimes x_{2})(d_{0}\otimes 1_{H})B(gx_{1}hd_{1}\otimes d_{2}) \\
-(d_{0}\otimes 1_{H})B(x_{1}x_{2}hd_{1}\otimes d_{2})%
\end{array}%
\right]  \notag
\end{eqnarray}

\section{CASES}

\subsubsection{Case $d=G$}

\paragraph{Case $h^{\prime }=g$}

The formula we want to prove is%
\begin{equation*}
B(h\otimes g)(G\otimes 1_{H})=(G\otimes 1_{H})B(hg\otimes 1_{H})
\end{equation*}%
By using $\left( \ref{eq.10}\right) $ this is equivalent to%
\begin{equation*}
(1_{A}\otimes g)B(gh\otimes 1_{H})(G\otimes 1_{H})(1_{A}\otimes g)=(G\otimes
1_{H})B(hg\otimes 1_{H})
\end{equation*}%
Since $\left( \ref{eq.a}\right) $ holds for any $B(h\otimes 1_{H})$ i.e.%
\begin{equation*}
B(h\otimes 1)(G\otimes 1_{H})=(G\otimes 1_{H})B(hg\otimes g)
\end{equation*}%
we deduce that%
\begin{equation*}
B(gh\otimes 1)(G\otimes 1_{H})=(G\otimes 1_{H})B(ghg\otimes gh)
\end{equation*}%
so that, by using $\left( \ref{eq.10}\right) ,$we get%
\begin{equation*}
B(gh\otimes 1)(G\otimes 1_{H})=(1_{A}\otimes g)(G\otimes 1_{H})B(hg\otimes
1_{H})(1_{A}\otimes g)
\end{equation*}%
and hence%
\begin{equation*}
(1_{A}\otimes g)B(gh\otimes 1_{H})(G\otimes 1_{H})(1_{A}\otimes g)=(G\otimes
1_{H})B(hg\otimes 1_{H}).
\end{equation*}

\paragraph{Case $h^{\prime }=x_{i}$}

The formula we want to prove is%
\begin{equation*}
B(h\otimes x_{i})(G\otimes 1_{H})=(G\otimes 1_{H})B(hg\otimes x_{i}g)
\end{equation*}

By formula $\left( \ref{form axi}\right) $we now that%
\begin{equation*}
B(h\otimes x_{i})(G\otimes 1_{H})=(G\otimes 1_{H})\left[
\begin{array}{c}
B(hg\otimes g)(1_{A}\otimes x_{i})-(1_{A}\otimes gx_{i})B(hg\otimes
g)(1_{A}\otimes g) \\
-(1_{A}\otimes g)B(gx_{i}hg\otimes g)(1_{A}\otimes g)%
\end{array}%
\right]
\end{equation*}%
Thus the left side for $d=G$ is%
\begin{equation*}
B(h\otimes x_{i})(G\otimes 1_{H})=(G\otimes 1_{H})\left[
\begin{array}{c}
B(hg\otimes g)(1_{A}\otimes x_{i})-(1_{A}\otimes gx_{i})B(hg\otimes
g)(1_{A}\otimes g) \\
-B(x_{i}hg\otimes 1_{H})%
\end{array}%
\right]
\end{equation*}%
We want to get the right side which is for $d=G$

\begin{gather*}
(d_{0}\otimes 1_{H})B(hd_{1}\otimes x_{i}d_{2})=(G\otimes 1_{H})B(hg\otimes
x_{i}g)=-(G\otimes 1_{H})B(hg\otimes gx_{i}) \\
=-(G\otimes 1_{H})(1_{A}\otimes g)B(ghg\otimes x_{i})(1_{A}\otimes g) \\
\text{we substitute the formula }\left( \ref{simplx}\right) \\
\text{ for }B\left( ghg\otimes x_{i}\right) =-(G\otimes 1_{H})(1_{A}\otimes
g)B(ghg\otimes 1_{H})(1_{A}\otimes x_{i})(1_{A}\otimes g)+ \\
(G\otimes 1_{H})(1_{A}\otimes g)(1_{A}\otimes gx_{i})B(ghg\otimes
1_{H})(1_{A}\otimes g)(1_{A}\otimes g) \\
+(G\otimes 1_{H})(1_{A}\otimes g)(1_{A}\otimes g)B(gx_{i}ghg\otimes
1_{H})(1_{A}\otimes g)(1_{A}\otimes g)
\end{gather*}%
and we obtain for the right side%
\begin{eqnarray*}
&&-(G\otimes 1_{H})(1_{A}\otimes g)B(ghg\otimes 1_{H})(1_{A}\otimes
x_{i})(1_{A}\otimes g)+ \\
&&(G\otimes 1_{H})(1_{A}\otimes g)(1_{A}\otimes gx_{i})B(ghg\otimes
1_{H})(1_{A}\otimes g)(1_{A}\otimes g) \\
&&+(G\otimes 1_{H})(1_{A}\otimes g)(1_{A}\otimes g)B(gx_{i}ghg\otimes
1_{H})(1_{A}\otimes g)(1_{A}\otimes g)
\end{eqnarray*}

Now we simplify the right side%
\begin{eqnarray*}
&&-(G\otimes 1_{H})(1_{A}\otimes g)B(ghg\otimes 1_{H})(1_{A}\otimes
x_{i})(1_{A}\otimes g)+ \\
&&(G\otimes 1_{H})(1_{A}\otimes g)(1_{A}\otimes gx_{i})B(ghg\otimes
1_{H})(1_{A}\otimes g)(1_{A}\otimes g) \\
&&-(G\otimes 1_{H})(1_{A}\otimes g)(1_{A}\otimes g)B(x_{i}hg\otimes
1_{H})(1_{A}\otimes g)(1_{A}\otimes g)
\end{eqnarray*}%
thus the right side is%
\begin{eqnarray*}
&&(G\otimes 1_{H})(1_{A}\otimes g)B(ghg\otimes 1_{H})(1_{A}\otimes
g)(1_{A}\otimes x_{i})+ \\
&&(G\otimes 1_{H})(1_{A}\otimes g)(1_{A}\otimes gx_{i})B(ghg\otimes 1_{H}) \\
&&-(G\otimes 1_{H})B(x_{i}hg\otimes 1_{H})
\end{eqnarray*}%
i.e.%
\begin{eqnarray*}
+ &&(G\otimes 1_{H})B(hg\otimes g)(1_{A}\otimes x_{i})+ \\
&&(G\otimes 1_{H})(1_{A}\otimes g)(1_{A}\otimes gx_{i})(1_{A}\otimes
g)B(hg\otimes g)(1_{A}\otimes g) \\
&&-(G\otimes 1_{H})B(x_{i}hg\otimes 1_{H})
\end{eqnarray*}%
\begin{eqnarray*}
&&\text{Thus the right side is} \\
+ &&(G\otimes 1_{H})B(hg\otimes g)(1_{A}\otimes x_{i})- \\
&&(G\otimes 1_{H})(1_{A}\otimes gx_{i})B(hg\otimes g)(1_{A}\otimes g) \\
&&-(G\otimes 1_{H})B(x_{i}hg\otimes 1_{H})
\end{eqnarray*}%
and hence it coincides with the left side.

\subsubsection{Case $h^{\prime }=x_{1}x_{2}$}

We apply the formula we know $\left( \ref{simplxx}\right) $%
\begin{equation*}
B(h\otimes x_{1}x_{2})=B(h\otimes x_{2})(1_{A}\otimes x_{1})-(1_{A}\otimes
gx_{1})B(h\otimes x_{2})(1_{A}\otimes g)-(1_{A}\otimes g)B(gx_{1}h\otimes
x_{2})(1_{A}\otimes g)
\end{equation*}

We write the left side%
\begin{eqnarray*}
&&B(h\otimes x_{1}x_{2})(G\otimes 1_{H}) \\
&=&B(h\otimes x_{2})(1_{A}\otimes x_{1})(G\otimes 1_{H})-(1_{A}\otimes
gx_{1})B(h\otimes x_{2})(1_{A}\otimes g)(G\otimes 1_{H}) \\
&&-(1_{A}\otimes g)B(gx_{1}h\otimes x_{2})(1_{A}\otimes g)(G\otimes 1_{H}) \\
&=&B(h\otimes x_{2})(G\otimes 1_{H})(1_{A}\otimes x_{1})-(1_{A}\otimes
gx_{1})B(h\otimes x_{2})(G\otimes 1_{H})(1_{A}\otimes g) \\
&&-(1_{A}\otimes g)B(gx_{1}h\otimes x_{2})(G\otimes 1_{H})(1_{A}\otimes g)
\end{eqnarray*}%
Now we can use formula a which holds for $h\otimes x_{2}$ i.e. $\left( \ref%
{form axi}\right) $%
\begin{equation*}
B(h\otimes x_{2})(G\otimes 1_{H})=(G\otimes 1_{H})B(hg\otimes x_{2}g)
\end{equation*}%
and we get

\begin{eqnarray*}
&&\text{left side} \\
&&(G\otimes 1_{H})B(hg\otimes x_{2}g)(1_{A}\otimes x_{1}) \\
&&-(1_{A}\otimes gx_{1})(G\otimes 1_{H})B(hg\otimes x_{2}g)(1_{A}\otimes g)
\\
&&-(1_{A}\otimes g)(G\otimes 1_{H})B(gx_{1}hg\otimes x_{2}g)(1_{A}\otimes g)
\end{eqnarray*}%
and we get%
\begin{eqnarray*}
&&\text{left side} \\
- &&(G\otimes 1_{H})B(hg\otimes gx_{2})(1_{A}\otimes x_{1}) \\
&&+(G\otimes 1_{H})(1_{A}\otimes gx_{1})B(hg\otimes gx_{2})(1_{A}\otimes g)
\\
&&+(G\otimes 1_{H})(1_{A}\otimes g)B(gx_{1}hg\otimes gx_{2})(1_{A}\otimes g)
\end{eqnarray*}%
\begin{eqnarray*}
&&\text{Thus the left side is} \\
&&(G\otimes 1_{H})\left[
\begin{array}{c}
-B(hg\otimes gx_{2})(1_{A}\otimes x_{1})+ \\
(1_{A}\otimes gx_{1})B(hg\otimes gx_{2})(1_{A}\otimes g)+ \\
+B(x_{1}hg\otimes x_{2})%
\end{array}%
\right]
\end{eqnarray*}

We develop the right side%
\begin{eqnarray*}
&&\text{the right side} \\
&&(G\otimes 1_{H})B(hg\otimes x_{1}x_{2}g) \\
&=&(G\otimes 1_{H})B(hg\otimes gx_{1}x_{2}) \\
&=&(G\otimes 1_{H})(1_{A}\otimes g)B(ghg\otimes x_{1}x_{2})(1_{A}\otimes g)
\end{eqnarray*}%
We apply the formula for $B\left( ghg\otimes x_{1}x_{2}\right) $ that we
know it holds $\left( \ref{simplxx}\right) $

\begin{equation*}
B(ghg\otimes x_{1}x_{2})=B(ghg\otimes x_{2})(1_{A}\otimes
x_{1})-(1_{A}\otimes gx_{1})B(ghg\otimes x_{2})(1_{A}\otimes
g)-(1_{A}\otimes g)B(gx_{1}ghg\otimes x_{2})(1_{A}\otimes g)
\end{equation*}%
which we rewrite in the following form%
\begin{equation*}
B(ghg\otimes x_{1}x_{2})=B(ghg\otimes x_{2})(1_{A}\otimes
x_{1})-(1_{A}\otimes gx_{1})B(ghg\otimes x_{2})(1_{A}\otimes
g)+B(gx_{1}hg\otimes gx_{2})
\end{equation*}

and substitute it inside%
\begin{eqnarray*}
&&\text{the right side} \\
&&(G\otimes 1_{H})(1_{A}\otimes g)B(ghg\otimes x_{2})(1_{A}\otimes
x_{1})(1_{A}\otimes g)+ \\
&&-(G\otimes 1_{H})(1_{A}\otimes g)(1_{A}\otimes gx_{1})B(ghg\otimes
x_{2})(1_{A}\otimes g)(1_{A}\otimes g)+ \\
&&+(G\otimes 1_{H})(1_{A}\otimes g)B(gx_{1}hg\otimes gx_{2})(1_{A}\otimes g)
\end{eqnarray*}%
and simplify

\begin{eqnarray*}
&&\text{the right side} \\
&=&-(G\otimes 1_{H})(1_{A}\otimes g)B(ghg\otimes x_{2})(1_{A}\otimes
g)(1_{A}\otimes x_{1})+ \\
&&-(G\otimes 1_{H})(1_{A}\otimes x_{1}g)(1_{A}\otimes g)B(ghg\otimes
x_{2})(1_{A}\otimes g)(1_{A}\otimes g)+ \\
&&+(G\otimes 1_{H})B(x_{1}hg\otimes x_{2})
\end{eqnarray*}%
we apply conjugation%
\begin{eqnarray*}
&&\text{the right side} \\
&&-(G\otimes 1_{H})B(hg\otimes gx_{2})(1_{A}\otimes x_{1})+ \\
&&+(G\otimes 1_{H})(1_{A}\otimes gx_{1})B(hg\otimes gx_{2})(1_{A}\otimes g)+
\\
&&+(G\otimes 1_{H})B(x_{1}hg\otimes x_{2})
\end{eqnarray*}%
which coincides with the left side and we conclude that formula $\left( \ref%
{eq.a}\right) $ holds.

\subsection{$d=X_{i}$}

The formula we want to prove is%
\begin{eqnarray*}
B(h\otimes h^{\prime })(X_{i}\otimes 1_{H}) &=&(X_{i}\otimes
1_{H})B(hg\otimes h^{\prime }g)+ \\
&&+(1_{A}\otimes 1_{H})B(hx_{i}\otimes h^{\prime }g) \\
&&+(1_{A}\otimes 1_{H})B(h\otimes h^{\prime }x_{i}).
\end{eqnarray*}%
Formula $\left( \ref{eq.a}\right) $
\begin{equation*}
B(h\otimes h^{\prime })(d\otimes 1_{H})=(d_{0}\otimes 1_{H})B(hd_{1}\otimes
h^{\prime }d_{2})
\end{equation*}%
in our case becomes%
\begin{equation*}
B(h\otimes h^{\prime })(X_{i}\otimes 1_{H})=(X_{i}\otimes 1_{H})B(hg\otimes
h^{\prime }g)+B(hx_{i}\otimes h^{\prime }g)+B(h\otimes h^{\prime }x_{i})
\end{equation*}

\subsubsection{Case $h^{\prime }=x_{i}$}

To write the left side we use formula $\left( \ref{form axi}\right) $%
\begin{equation*}
B(h\otimes x_{i})(d\otimes 1_{H})=(d_{0}\otimes 1_{H})\left[
\begin{array}{c}
B(hd_{1}\otimes d_{2})(1_{A}\otimes x_{i})-(1_{A}\otimes
gx_{i})B(hd_{1}\otimes d_{2})(1_{A}\otimes g) \\
-(1_{A}\otimes g)B(gx_{i}hd_{1}\otimes d_{2})(1_{A}\otimes g)%
\end{array}%
\right]
\end{equation*}%
in this case%
\begin{eqnarray*}
&&B(h\otimes x_{i})(X_{j}\otimes 1_{H}) \\
&=&(X_{j}\otimes 1_{H})\left[
\begin{array}{c}
B(hg\otimes g)(1_{A}\otimes x_{i})-(1_{A}\otimes gx_{i})B(hg\otimes
g)(1_{A}\otimes g) \\
-(1_{A}\otimes g)B(gx_{i}hg\otimes g)(1_{A}\otimes g)%
\end{array}%
\right] \\
&&+(1_{A}\otimes 1_{H})\left[
\begin{array}{c}
B(hx_{j}\otimes g)(1_{A}\otimes x_{i})-(1_{A}\otimes gx_{i})B(hx_{j}\otimes
g)(1_{A}\otimes g) \\
-(1_{A}\otimes g)B(gx_{i}hx_{j}\otimes g)(1_{A}\otimes g)%
\end{array}%
\right] \\
&&++(1_{A}\otimes 1_{H})\left[
\begin{array}{c}
B(h\otimes x_{j})(1_{A}\otimes x_{i})-(1_{A}\otimes gx_{i})B(h\otimes
x_{j})(1_{A}\otimes g) \\
-(1_{A}\otimes g)B(gx_{i}h\otimes x_{j})(1_{A}\otimes g)%
\end{array}%
\right]
\end{eqnarray*}%
Thus the left side is%
\begin{eqnarray*}
&&B(h\otimes x_{i})(X_{j}\otimes 1_{H}) \\
&=&(X_{j}\otimes 1_{H})\left[
\begin{array}{c}
B(hg\otimes g)(1_{A}\otimes x_{i})-(1_{A}\otimes gx_{i})B(hg\otimes
g)(1_{A}\otimes g) \\
-(1_{A}\otimes g)B(gx_{i}hg\otimes g)(1_{A}\otimes g)%
\end{array}%
\right] \\
&&+\left[
\begin{array}{c}
B(hx_{j}\otimes g)(1_{A}\otimes x_{i})-(1_{A}\otimes gx_{i})B(hx_{j}\otimes
g)(1_{A}\otimes g) \\
-(1_{A}\otimes g)B(gx_{i}hx_{j}\otimes g)(1_{A}\otimes g)%
\end{array}%
\right] \\
&&+\left[
\begin{array}{c}
B(h\otimes x_{j})(1_{A}\otimes x_{i})-(1_{A}\otimes gx_{i})B(h\otimes
x_{j})(1_{A}\otimes g) \\
-(1_{A}\otimes g)B(gx_{i}h\otimes x_{j})(1_{A}\otimes g)%
\end{array}%
\right]
\end{eqnarray*}%
Now we write the right side%
\begin{eqnarray*}
&&(X_{j}\otimes 1_{H})B(hg\otimes x_{i}g)+ \\
&&(1_{A}\otimes 1_{H})B(hx_{j}\otimes x_{i}g)+ \\
&&(1_{A}\otimes 1_{H})B(h\otimes x_{i}x_{j})+
\end{eqnarray*}%
which we rewrite%
\begin{eqnarray}
- &&(X_{j}\otimes 1_{H})B(hg\otimes gx_{i})+  \label{rightXjxi} \\
- &&B(hx_{j}\otimes gx_{i})+  \notag \\
&&B(h\otimes x_{i}x_{j})+  \notag
\end{eqnarray}%
we use $\left( \ref{simplgx}\right) $ for $"h"=hg,hx_{j}$

\begin{eqnarray*}
&&B(hg\otimes gx_{i}) \\
&=&(1_{A}\otimes g)B(ghg\otimes 1_{H})(1_{A}\otimes gx_{i}) \\
&&+(1_{A}\otimes x_{i})B(ghg\otimes 1_{H}) \\
&&-B(gx_{i}ghg\otimes 1_{H})
\end{eqnarray*}%
\begin{eqnarray*}
&&B(hx_{j}\otimes gx_{i}) \\
&=&(1_{A}\otimes g)B(ghx_{j}\otimes 1_{H})(1_{A}\otimes gx_{i}) \\
&&+(1_{A}\otimes x_{i})B(ghx_{j}\otimes 1_{H}) \\
&&-B(gx_{i}ghx_{j}\otimes 1_{H})
\end{eqnarray*}%
and we get that the right side is%
\begin{eqnarray*}
&&\text{right side} \\
- &&(X_{j}\otimes 1_{H})\left[
\begin{array}{c}
-(1_{A}\otimes g)B(ghg\otimes 1_{H})(1_{A}\otimes g)\left( 1_{A}\otimes
x_{i}\right) + \\
(1_{A}\otimes x_{i})B(ghg\otimes 1_{H})+ \\
-B(gx_{i}ghg\otimes 1_{H})%
\end{array}%
\right] + \\
&&-\left[
\begin{array}{c}
(1_{A}\otimes g)B(ghx_{j}\otimes 1_{H})(1_{A}\otimes gx_{i})+ \\
+(1_{A}\otimes x_{i})B(ghx_{j}\otimes 1_{H})+ \\
-B(gx_{i}ghx_{j}\otimes 1_{H})%
\end{array}%
\right] + \\
&&+B(h\otimes x_{i}x_{j})+
\end{eqnarray*}%
\begin{eqnarray*}
&&\text{right side} \\
- &&(X_{j}\otimes 1_{H})\left[
\begin{array}{c}
-B(hg\otimes g)\left( 1_{A}\otimes x_{i}\right) + \\
(1_{A}\otimes x_{i})(1_{A}\otimes g)B(hg\otimes g)(1_{A}\otimes g)+ \\
+(1_{A}\otimes g)B(gx_{i}hg\otimes g)(1_{A}\otimes g)%
\end{array}%
\right] + \\
&&-\left[
\begin{array}{c}
-B(hx_{j}\otimes g)(1_{A}\otimes x_{i})+ \\
+(1_{A}\otimes gx_{i})B(hx_{j}\otimes g)(1_{A}\otimes g)+ \\
+(1_{A}\otimes g)B(gx_{i}hx_{j}\otimes g)(1_{A}\otimes g)%
\end{array}%
\right] + \\
&&+B(h\otimes x_{i}x_{j})+
\end{eqnarray*}%
Now we consider $i=1$ and $j=2$ and we use $\left( \ref{simplxx}\right) $

\begin{equation*}
B(h\otimes x_{1}x_{2})=B(h\otimes x_{2})(1_{A}\otimes x_{1})-(1_{A}\otimes
gx_{1})B(h\otimes x_{2})(1_{A}\otimes g)-(1_{A}\otimes g)B(gx_{1}h\otimes
x_{2})(1_{A}\otimes g)
\end{equation*}%
and get%
\begin{eqnarray*}
&&\text{right side} \\
- &&(X_{j}\otimes 1_{H})\left[
\begin{array}{c}
-B(hg\otimes g)\left( 1_{A}\otimes x_{i}\right) + \\
(1_{A}\otimes x_{i})(1_{A}\otimes g)B(hg\otimes g)(1_{A}\otimes g)+ \\
+(1_{A}\otimes g)B(gx_{i}hg\otimes g)(1_{A}\otimes g)%
\end{array}%
\right] + \\
&&-\left[
\begin{array}{c}
-B(hx_{j}\otimes g)(1_{A}\otimes x_{i})+ \\
+(1_{A}\otimes gx_{i})B(hx_{j}\otimes g)(1_{A}\otimes g)+ \\
+(1_{A}\otimes g)B(gx_{i}hx_{j}\otimes g)(1_{A}\otimes g)%
\end{array}%
\right] + \\
&&+\left[
\begin{array}{c}
B(h\otimes x_{2})(1_{A}\otimes x_{1}) \\
-(1_{A}\otimes gx_{1})B(h\otimes x_{2})(1_{A}\otimes g) \\
-(1_{A}\otimes g)B(gx_{1}h\otimes x_{2})(1_{A}\otimes g)%
\end{array}%
\right] +
\end{eqnarray*}%
and we can conclude that formula $\left( \ref{eq.a}\right) $ holds.

We assume that for any $h$ and $h^{\prime }=1_{H}$ $\left( \ref{eq.a}\right)
$ holds and $d=X_{i}$ i.e.
\begin{equation*}
B(h\otimes 1_{H})(X_{i}\otimes 1_{H})=(X_{i}\otimes 1_{H})B(hg\otimes
g)+B(hx_{i}\otimes g)+B(h\otimes x_{i})
\end{equation*}%
\begin{equation*}
B(h\otimes 1_{H})(X_{1}\otimes 1_{H})=(X_{1}\otimes 1_{H})B(hg\otimes
g)+B(hx_{1}\otimes g)+B(h\otimes x_{1})
\end{equation*}

\paragraph{Case $h^{\prime }=x_{1}x_{2}$}

We have to prove that for any $h$ the formula below holds i.e.%
\begin{equation*}
B(h\otimes x_{1}x_{2})(X_{i}\otimes 1_{H})=(X_{i}\otimes 1_{H})B(hg\otimes
gx_{1}x_{2})+B(hx_{i}\otimes gx_{1}x_{2})
\end{equation*}

We use $\left( \ref{simplxx}\right) $%
\begin{eqnarray*}
B(h\otimes x_{1}x_{2}) &=&B(h\otimes 1_{H})(1_{A}\otimes x_{1}x_{2}) \\
&&+(1_{A}\otimes gx_{2})B(h\otimes 1_{H})(1_{A}\otimes gx_{1}) \\
&&+(1_{A}\otimes g)B(gx_{2}h\otimes 1_{H})(1_{A}\otimes gx_{1}) \\
&&-(1_{A}\otimes gx_{1})B(h\otimes 1_{H})(1_{A}\otimes gx_{2}) \\
&&+(1_{A}\otimes x_{1}x_{2})B(h\otimes 1_{H}) \\
&&+(1_{A}\otimes x_{1})B(gx_{2}h\otimes 1_{H}) \\
&&-(1_{A}\otimes g)B(gx_{1}h\otimes 1_{H})(1_{A}\otimes gx_{2}) \\
&&-(1_{A}\otimes x_{2})B(gx_{1}h\otimes 1_{H}) \\
&&+B(x_{1}x_{2}h\otimes 1_{H})
\end{eqnarray*}%
We rewrite this%
\begin{eqnarray*}
B(h\otimes x_{1}x_{2}) &=&B(h\otimes 1_{H})(1_{A}\otimes x_{1}x_{2}) \\
&&+(1_{A}\otimes gx_{2})B(h\otimes 1_{H})(1_{A}\otimes gx_{1}) \\
&&-(1_{A}\otimes gx_{1})B(h\otimes 1_{H})(1_{A}\otimes gx_{2}) \\
&&+(1_{A}\otimes x_{1}x_{2})B(h\otimes 1_{H}) \\
&&-(1_{A}\otimes x_{2})B(gx_{1}h\otimes 1_{H}) \\
&&-(1_{A}\otimes g)B(gx_{1}h\otimes 1_{H})(1_{A}\otimes gx_{2}) \\
&&+(1_{A}\otimes x_{1})B(gx_{2}h\otimes 1_{H}) \\
&&+(1_{A}\otimes g)B(gx_{2}h\otimes 1_{H})(1_{A}\otimes gx_{1}) \\
&&+B(x_{1}x_{2}h\otimes 1_{H})
\end{eqnarray*}

Here we do this for $X_{1}$i.e. we assume that%
\begin{equation*}
B(h\otimes 1_{H})(X_{1}\otimes 1_{H})=(X_{1}\otimes 1_{H})B(hg\otimes
g)+B(hx_{1}\otimes g)+B(h\otimes x_{1})
\end{equation*}%
Thus, using the equality above, we get%
\begin{eqnarray*}
B(h\otimes x_{1}x_{2})\left( X_{1}\otimes 1_{H}\right) &=&B(h\otimes
1_{H})\left( X_{1}\otimes 1_{H}\right) (1_{A}\otimes x_{1}x_{2}) \\
&&+(1_{A}\otimes gx_{2})B(h\otimes 1_{H})\left( X_{1}\otimes 1_{H}\right)
(1_{A}\otimes gx_{1}) \\
&&-(1_{A}\otimes gx_{1})B(h\otimes 1_{H})\left( X_{1}\otimes 1_{H}\right)
(1_{A}\otimes gx_{2}) \\
&&+(1_{A}\otimes x_{1}x_{2})B(h\otimes 1_{H})\left( X_{1}\otimes 1_{H}\right)
\\
&&-(1_{A}\otimes x_{2})B(gx_{1}h\otimes 1_{H})\left( X_{1}\otimes
1_{H}\right) \\
&&-(1_{A}\otimes g)B(gx_{1}h\otimes 1_{H})\left( X_{1}\otimes 1_{H}\right)
(1_{A}\otimes gx_{2}) \\
&&+(1_{A}\otimes x_{1})B(gx_{2}h\otimes 1_{H})\left( X_{1}\otimes
1_{H}\right) \\
&&+(1_{A}\otimes g)B(gx_{2}h\otimes 1_{H})\left( X_{1}\otimes 1_{H}\right)
(1_{A}\otimes gx_{1}) \\
&&+B(x_{1}x_{2}h\otimes 1_{H})\left( X_{1}\otimes 1_{H}\right)
\end{eqnarray*}%
\begin{equation*}
B(x_{1}x_{2}h\otimes 1_{H})\left( X_{1}\otimes 1_{H}\right) =(X_{1}\otimes
1_{H})B(x_{1}x_{2}hg\otimes g)+B(x_{1}x_{2}h\otimes x_{1})
\end{equation*}%
We split the right side in four parts and we get the following.%
\begin{eqnarray*}
&&\text{1}B(h\otimes 1_{H})\left( X_{1}\otimes 1_{H}\right) (1_{A}\otimes
x_{1}x_{2}) \\
&&+(1_{A}\otimes gx_{2})B(h\otimes 1_{H})\left( X_{1}\otimes 1_{H}\right)
(1_{A}\otimes gx_{1}) \\
&&-(1_{A}\otimes gx_{1})B(h\otimes 1_{H})\left( X_{1}\otimes 1_{H}\right)
(1_{A}\otimes gx_{2}) \\
&&+(1_{A}\otimes x_{1}x_{2})B(h\otimes 1_{H})\left( X_{1}\otimes 1_{H}\right)
\\
&=&\left[ (X_{1}\otimes 1_{H})B(hg\otimes g)+B(hx_{1}\otimes g)+B(h\otimes
x_{1})\right] (1_{A}\otimes x_{1}x_{2}) \\
&&(1_{A}\otimes gx_{2})\left[ (X_{1}\otimes 1_{H})B(hg\otimes
g)+B(hx_{1}\otimes g)+B(h\otimes x_{1})\right] (1_{A}\otimes gx_{1}) \\
&&-(1_{A}\otimes gx_{1})\left[ (X_{1}\otimes 1_{H})B(hg\otimes
g)+B(hx_{1}\otimes g)+B(h\otimes x_{1})\right] (1_{A}\otimes gx_{2}) \\
&&+(1_{A}\otimes x_{1}x_{2})\left[ (X_{1}\otimes 1_{H})B(hg\otimes
g)+B(hx_{1}\otimes g)+B(h\otimes x_{1})\right]
\end{eqnarray*}%
\begin{equation*}
B(gx_{1}h\otimes 1_{H})(X_{1}\otimes 1_{H})=(X_{1}\otimes
1_{H})B(gx_{1}hg\otimes g)+B(gx_{1}hx_{1}\otimes g)+B(gx_{1}h\otimes x_{1})
\end{equation*}%
\begin{eqnarray*}
\text{2} &&-(1_{A}\otimes x_{2})B(gx_{1}h\otimes 1_{H})\left( X_{1}\otimes
1_{H}\right) \\
&&-(1_{A}\otimes g)B(gx_{1}h\otimes 1_{H})\left( X_{1}\otimes 1_{H}\right)
(1_{A}\otimes gx_{2}) \\
&=&-(1_{A}\otimes x_{2})\left[ (X_{1}\otimes 1_{H})B(gx_{1}hg\otimes
g)+B(gx_{1}hx_{1}\otimes g)+B(gx_{1}h\otimes x_{1})\right] \\
&&-(1_{A}\otimes g)\left[ (X_{1}\otimes 1_{H})B(gx_{1}hg\otimes
g)+B(gx_{1}hx_{1}\otimes g)+B(gx_{1}h\otimes x_{1})\right] (1_{A}\otimes
gx_{2})
\end{eqnarray*}%
\begin{equation*}
B(gx_{2}h\otimes 1_{H})(X_{1}\otimes 1_{H})=(X_{1}\otimes
1_{H})B(gx_{2}hg\otimes g)+B(gx_{2}hx_{1}\otimes g)+B(gx_{2}h\otimes x_{1})
\end{equation*}%
\begin{eqnarray*}
\text{3} &&+(1_{A}\otimes x_{1})B(gx_{2}h\otimes 1_{H})\left( X_{1}\otimes
1_{H}\right) \\
&&+(1_{A}\otimes g)B(gx_{2}h\otimes 1_{H})\left( X_{1}\otimes 1_{H}\right)
(1_{A}\otimes gx_{1}) \\
&=&(1_{A}\otimes x_{1})\left[ (X_{1}\otimes 1_{H})B(gx_{2}hg\otimes
g)+B(gx_{2}hx_{1}\otimes g)+B(gx_{2}h\otimes x_{1})\right] \\
&&(1_{A}\otimes g)\left[ (X_{1}\otimes 1_{H})B(gx_{2}hg\otimes
g)+B(gx_{2}hx_{1}\otimes g)+B(gx_{2}h\otimes x_{1})\right] (1_{A}\otimes
gx_{1})
\end{eqnarray*}%
\begin{equation*}
\text{4}B(x_{1}x_{2}h\otimes 1_{H})(X_{1}\otimes 1_{H})=(X_{1}\otimes
1_{H})B(x_{1}x_{2}hg\otimes g)+B(x_{1}x_{2}h\otimes x_{1})
\end{equation*}%
Therefore the left side, using the assumption that for any $h$ $\left( \ref%
{eq.a}\right) $ holds and $d=X_{1}$ gives us%
\begin{eqnarray*}
&&B(h\otimes x_{1}x_{2})\left( X_{1}\otimes 1_{H}\right) \\
&=&\left[ (X_{1}\otimes 1_{H})B(hg\otimes g)+B(hx_{1}\otimes g)+B(h\otimes
x_{1})\right] (1_{A}\otimes x_{1}x_{2}) \\
&&(1_{A}\otimes gx_{2})\left[ (X_{1}\otimes 1_{H})B(hg\otimes
g)+B(hx_{1}\otimes g)+B(h\otimes x_{1})\right] (1_{A}\otimes gx_{1}) \\
&&-(1_{A}\otimes gx_{1})\left[ (X_{1}\otimes 1_{H})B(hg\otimes
g)+B(hx_{1}\otimes g)+B(h\otimes x_{1})\right] (1_{A}\otimes gx_{2}) \\
&&+(1_{A}\otimes x_{1}x_{2})\left[ (X_{1}\otimes 1_{H})B(hg\otimes
g)+B(hx_{1}\otimes g)+B(h\otimes x_{1})\right] + \\
&&-(1_{A}\otimes x_{2})\left[ (X_{1}\otimes 1_{H})B(gx_{1}hg\otimes
g)+B(gx_{1}hx_{1}\otimes g)+B(gx_{1}h\otimes x_{1})\right] \\
&&-(1_{A}\otimes g)\left[ (X_{1}\otimes 1_{H})B(gx_{1}hg\otimes
g)+B(gx_{1}hx_{1}\otimes g)+B(gx_{1}h\otimes x_{1})\right] (1_{A}\otimes
gx_{2})+ \\
&&(1_{A}\otimes x_{1})\left[ (X_{1}\otimes 1_{H})B(gx_{2}hg\otimes
g)+B(gx_{2}hx_{1}\otimes g)+B(gx_{2}h\otimes x_{1})\right] + \\
&&(1_{A}\otimes g)\left[ (X_{1}\otimes 1_{H})B(gx_{2}hg\otimes
g)+B(gx_{2}hx_{1}\otimes g)+B(gx_{2}h\otimes x_{1})\right] (1_{A}\otimes
gx_{1}) \\
&&(X_{1}\otimes 1_{H})B(x_{1}x_{2}hg\otimes g)+B(x_{1}x_{2}h\otimes x_{1})
\end{eqnarray*}%
Now we should compare it with the right side which is%
\begin{equation*}
(X_{1}\otimes 1_{H})B(hg\otimes gx_{1}x_{2})+B(hx_{1}\otimes gx_{1}x_{2})
\end{equation*}%
We begin by writing $B(hg\otimes gx_{1}x_{2}).$We use $\left( \ref{simplgxx}%
\right) $%
\begin{eqnarray*}
B(hg\otimes gx_{1}x_{2}) &=&(1_{A}\otimes g)B(ghg\otimes 1_{H})(1_{A}\otimes
gx_{1}x_{2}) \\
&&-(1_{A}\otimes x_{2})B(ghg\otimes 1_{H})(1_{A}\otimes x_{1}) \\
&&-B(x_{2}hg\otimes 1_{H})(1_{A}\otimes x_{1}) \\
&&+(1_{A}\otimes x_{1})B(ghg\otimes 1_{H})(1_{A}\otimes x_{2}) \\
&&+(1_{A}\otimes gx_{1}x_{2})B(ghg\otimes 1_{H})(1_{A}\otimes g) \\
&&+(1_{A}\otimes gx_{1})B(x_{2}hg\otimes 1_{H})(1_{A}\otimes g) \\
&&+B(x_{1}hg\otimes 1_{H})(1_{A}\otimes x_{2}) \\
&&-(1_{A}\otimes gx_{2})B(x_{1}hg\otimes 1_{H})(1_{A}\otimes g) \\
&&+(1_{A}\otimes g)B(gx_{1}x_{2}hg\otimes 1_{H})(1_{A}\otimes g)
\end{eqnarray*}%
and then $B(hx_{1}\otimes gx_{1}x_{2})$ always using $\left( \ref{simplgxx}%
\right) $%
\begin{eqnarray*}
B(hx_{1}\otimes gx_{1}x_{2}) &=&(1_{A}\otimes g)B(ghx_{1}\otimes
1_{H})(1_{A}\otimes gx_{1}x_{2}) \\
&&-(1_{A}\otimes x_{2})B(ghx_{1}\otimes 1_{H})(1_{A}\otimes x_{1}) \\
&&-B(x_{2}hx_{1}\otimes 1_{H})(1_{A}\otimes x_{1}) \\
&&+(1_{A}\otimes x_{1})B(ghx_{1}\otimes 1_{H})(1_{A}\otimes x_{2}) \\
&&+(1_{A}\otimes gx_{1}x_{2})B(ghx_{1}\otimes 1_{H})(1_{A}\otimes g) \\
&&+(1_{A}\otimes gx_{1})B(x_{2}hx_{1}\otimes 1_{H})(1_{A}\otimes g)
\end{eqnarray*}%
Thus the right side is

\begin{eqnarray*}
&&\text{right side}(X_{1}\otimes 1_{H})B(hg\otimes
gx_{1}x_{2})+B(hx_{1}\otimes gx_{1}x_{2}) \\
&=&(X_{1}\otimes 1_{H})(1_{A}\otimes g)B(ghg\otimes 1_{H})(1_{A}\otimes
gx_{1}x_{2}) \\
&&-(X_{1}\otimes 1_{H})(1_{A}\otimes x_{2})B(ghg\otimes 1_{H})(1_{A}\otimes
x_{1}) \\
&&-(X_{1}\otimes 1_{H})B(x_{2}hg\otimes 1_{H})(1_{A}\otimes x_{1}) \\
&&+(X_{1}\otimes 1_{H})(1_{A}\otimes x_{1})B(ghg\otimes 1_{H})(1_{A}\otimes
x_{2}) \\
&&+(X_{1}\otimes 1_{H})(1_{A}\otimes gx_{1}x_{2})B(ghg\otimes
1_{H})(1_{A}\otimes g) \\
&&+(X_{1}\otimes 1_{H})(1_{A}\otimes gx_{1})B(x_{2}hg\otimes
1_{H})(1_{A}\otimes g) \\
&&+(X_{1}\otimes 1_{H})B(x_{1}hg\otimes 1_{H})(1_{A}\otimes x_{2}) \\
&&-(X_{1}\otimes 1_{H})(1_{A}\otimes gx_{2})B(x_{1}hg\otimes
1_{H})(1_{A}\otimes g) \\
&&+(X_{1}\otimes 1_{H})(1_{A}\otimes g)B(gx_{1}x_{2}hg\otimes
1_{H})(1_{A}\otimes g) \\
&&+(1_{A}\otimes g)B(ghx_{1}\otimes 1_{H})(1_{A}\otimes gx_{1}x_{2}) \\
&&-(1_{A}\otimes x_{2})B(ghx_{1}\otimes 1_{H})(1_{A}\otimes x_{1}) \\
&&-B(x_{2}hx_{1}\otimes 1_{H})(1_{A}\otimes x_{1}) \\
&&+(1_{A}\otimes x_{1})B(ghx_{1}\otimes 1_{H})(1_{A}\otimes x_{2})+ \\
&&+(1_{A}\otimes gx_{1}x_{2})B(ghx_{1}\otimes 1_{H})(1_{A}\otimes g) \\
&&+(1_{A}\otimes gx_{1})B(x_{2}hx_{1}\otimes 1_{H})(1_{A}\otimes g)
\end{eqnarray*}%
\begin{gather*}
\text{left side + assumption}\text{ }B(h\otimes x_{1}x_{2})\left(
X_{1}\otimes 1_{H}\right) \\
=\left[ (X_{1}\otimes 1_{H})B(hg\otimes g)+B(hx_{1}\otimes g)+B(h\otimes
x_{1})\right] (1_{A}\otimes x_{1}x_{2}) \\
(1_{A}\otimes gx_{2})\left[ (X_{1}\otimes 1_{H})B(hg\otimes
g)+B(hx_{1}\otimes g)+B(h\otimes x_{1})\right] \\
(1_{A}\otimes gx_{1}) \\
-(1_{A}\otimes gx_{1})\left[ (X_{1}\otimes 1_{H})B(hg\otimes
g)+B(hx_{1}\otimes g)+B(h\otimes x_{1})\right] \\
(1_{A}\otimes gx_{2}) \\
+(1_{A}\otimes x_{1}x_{2})\left[ (X_{1}\otimes 1_{H})B(hg\otimes
g)+B(hx_{1}\otimes g)+B(h\otimes x_{1})\right] + \\
-(1_{A}\otimes x_{2})\left[ (X_{1}\otimes 1_{H})B(gx_{1}hg\otimes
g)+B(gx_{1}hx_{1}\otimes g)+B(gx_{1}h\otimes x_{1})\right] \\
-(1_{A}\otimes g)\left[ (X_{1}\otimes 1_{H})B(gx_{1}hg\otimes
g)+B(gx_{1}hx_{1}\otimes g)+B(gx_{1}h\otimes x_{1})\right] (1_{A}\otimes
gx_{2})+ \\
(1_{A}\otimes x_{1})\left[ (X_{1}\otimes 1_{H})B(gx_{2}hg\otimes
g)+B(gx_{2}hx_{1}\otimes g)+B(gx_{2}h\otimes x_{1})\right] + \\
(1_{A}\otimes g)\left[ (X_{1}\otimes 1_{H})B(gx_{2}hg\otimes
g)+B(gx_{2}hx_{1}\otimes g)+B(gx_{2}h\otimes x_{1})\right] \\
(1_{A}\otimes gx_{1}) \\
(X_{1}\otimes 1_{H})B(x_{1}x_{2}hg\otimes g)+B(x_{1}x_{2}h\otimes x_{1})
\end{gather*}%
we rewrite this{\LARGE \ }a first time by using $\left( \ref{eq.10}\right) .$%
\begin{gather*}
1\text{rewritten left side + assumption}B(h\otimes x_{1}x_{2})\left(
X_{1}\otimes 1_{H}\right) \\
=\left[ (X_{1}\otimes 1_{H})(1_{A}\otimes g)B\left( ghg\otimes 1_{H}\right)
(1_{A}\otimes g)+B(hx_{1}\otimes g)+B(h\otimes x_{1})\right] (1_{A}\otimes
x_{1}x_{2}) \\
(1_{A}\otimes gx_{2})\left[ (X_{1}\otimes 1_{H})(1_{A}\otimes g)B\left(
ghg\otimes 1_{H}\right) (1_{A}\otimes g)+B(hx_{1}\otimes g)+B(h\otimes x_{1})%
\right] (1_{A}\otimes gx_{1}) \\
-(1_{A}\otimes gx_{1})\left[ (X_{1}\otimes 1_{H})(1_{A}\otimes g)B\left(
ghg\otimes 1_{H}\right) (1_{A}\otimes g)+B(hx_{1}\otimes g)+B(h\otimes x_{1})%
\right] (1_{A}\otimes gx_{2}) \\
+(1_{A}\otimes x_{1}x_{2})\left[ (X_{1}\otimes 1_{H})(1_{A}\otimes g)B\left(
ghg\otimes 1_{H}\right) (1_{A}\otimes g)+B(hx_{1}\otimes g)+B(h\otimes x_{1})%
\right] + \\
-(1_{A}\otimes x_{2})\left[ (X_{1}\otimes 1_{H})(1_{A}\otimes g)B\left(
x_{1}hg\otimes 1_{H}\right) (1_{A}\otimes g)+B(gx_{1}h\otimes x_{1})\right]
\\
-(1_{A}\otimes g)\left[ (X_{1}\otimes 1_{H})(1_{A}\otimes g)B\left(
x_{1}hg\otimes 1_{H}\right) (1_{A}\otimes g)+B(gx_{1}h\otimes x_{1})\right]
(1_{A}\otimes gx_{2})+ \\
(1_{A}\otimes x_{1})\left[ (X_{1}\otimes 1_{H})(1_{A}\otimes g)B\left(
x_{2}hg\otimes 1_{H}\right) (1_{A}\otimes g)+B(gx_{2}hx_{1}\otimes
g)+B(gx_{2}h\otimes x_{1})\right] + \\
(1_{A}\otimes g)\left[ (X_{1}\otimes 1_{H})(1_{A}\otimes g)B\left(
x_{2}hg\otimes 1_{H}\right) (1_{A}\otimes g)+B(gx_{2}hx_{1}\otimes
g)+B(gx_{2}h\otimes x_{1})\right] (1_{A}\otimes gx_{1}) \\
(X_{1}\otimes 1_{H})B(x_{1}x_{2}hg\otimes g)+B(x_{1}x_{2}h\otimes x_{1})
\end{gather*}%
we rewrite it a second time%
\begin{gather*}
2\text{rewritten left side + assumption}B(h\otimes x_{1}x_{2})\left(
X_{1}\otimes 1_{H}\right) \\
=(X_{1}\otimes 1_{H})(1_{A}\otimes g)B\left( ghg\otimes 1_{H}\right)
(1_{A}\otimes g)(1_{A}\otimes x_{1}x_{2}) \\
+B(hx_{1}\otimes g)(1_{A}\otimes x_{1}x_{2})+B(h\otimes x_{1})(1_{A}\otimes
x_{1}x_{2}) \\
(X_{1}\otimes 1_{H})(1_{A}\otimes gx_{2})(1_{A}\otimes g)B\left( ghg\otimes
1_{H}\right) (1_{A}\otimes g)(1_{A}\otimes gx_{1}) \\
+(1_{A}\otimes gx_{2})B(hx_{1}\otimes g)(1_{A}\otimes gx_{1})+(1_{A}\otimes
gx_{2})B(h\otimes x_{1})(1_{A}\otimes gx_{1}) \\
-(X_{1}\otimes 1_{H})(1_{A}\otimes gx_{1})(1_{A}\otimes g)B\left( ghg\otimes
1_{H}\right) (1_{A}\otimes g)(1_{A}\otimes gx_{2}) \\
-(1_{A}\otimes gx_{1})B(hx_{1}\otimes g)(1_{A}\otimes gx_{2})-(1_{A}\otimes
gx_{1})B(h\otimes x_{1})(1_{A}\otimes gx_{2}) \\
+(X_{1}\otimes 1_{H})(1_{A}\otimes x_{1}x_{2})(1_{A}\otimes g)B\left(
ghg\otimes 1_{H}\right) (1_{A}\otimes g) \\
+(1_{A}\otimes x_{1}x_{2})B(hx_{1}\otimes g)+(1_{A}\otimes
x_{1}x_{2})B(h\otimes x_{1})+ \\
-(X_{1}\otimes 1_{H})(1_{A}\otimes x_{2})(1_{A}\otimes g)B\left(
x_{1}hg\otimes 1_{H}\right) (1_{A}\otimes g) \\
-(1_{A}\otimes x_{2})B(gx_{1}h\otimes x_{1}) \\
-(X_{1}\otimes 1_{H})(1_{A}\otimes g)(1_{A}\otimes g)B\left( x_{1}hg\otimes
1_{H}\right) (1_{A}\otimes g)(1_{A}\otimes gx_{2}) \\
-(1_{A}\otimes g)B(gx_{1}h\otimes x_{1})(1_{A}\otimes gx_{2})+ \\
(X_{1}\otimes 1_{H})(1_{A}\otimes x_{1})(1_{A}\otimes g)B\left(
x_{2}hg\otimes 1_{H}\right) (1_{A}\otimes g) \\
+(1_{A}\otimes x_{1})B(gx_{2}hx_{1}\otimes g)+(1_{A}\otimes
x_{1})B(gx_{2}h\otimes x_{1})+ \\
(X_{1}\otimes 1_{H})(1_{A}\otimes g)(1_{A}\otimes g)B\left( x_{2}hg\otimes
1_{H}\right) (1_{A}\otimes g)(1_{A}\otimes gx_{1}) \\
+(1_{A}\otimes g)B(gx_{2}hx_{1}\otimes g)(1_{A}\otimes gx_{1})+(1_{A}\otimes
g)B(gx_{2}h\otimes x_{1})(1_{A}\otimes gx_{1}) \\
(X_{1}\otimes 1_{H})B(x_{1}x_{2}hg\otimes g)+ \\
+B(x_{1}x_{2}h\otimes x_{1})
\end{gather*}%
we rewrite it a third time%
\begin{gather*}
3\text{rewritten left side + assumption}B(h\otimes x_{1}x_{2})\left(
X_{1}\otimes 1_{H}\right) \\
=(X_{1}\otimes 1_{H})(1_{A}\otimes g)B\left( ghg\otimes 1_{H}\right)
(1_{A}\otimes g)(1_{A}\otimes x_{1}x_{2}) \\
(X_{1}\otimes 1_{H})(1_{A}\otimes gx_{2})(1_{A}\otimes g)B\left( ghg\otimes
1_{H}\right) (1_{A}\otimes g)(1_{A}\otimes gx_{1}) \\
-(X_{1}\otimes 1_{H})(1_{A}\otimes gx_{1})(1_{A}\otimes g)B\left( ghg\otimes
1_{H}\right) (1_{A}\otimes g)(1_{A}\otimes gx_{2}) \\
+(X_{1}\otimes 1_{H})(1_{A}\otimes x_{1}x_{2})(1_{A}\otimes g)B\left(
ghg\otimes 1_{H}\right) (1_{A}\otimes g) \\
-(X_{1}\otimes 1_{H})(1_{A}\otimes x_{2})(1_{A}\otimes g)B\left(
x_{1}hg\otimes 1_{H}\right) (1_{A}\otimes g) \\
-(X_{1}\otimes 1_{H})(1_{A}\otimes g)(1_{A}\otimes g)B\left( x_{1}hg\otimes
1_{H}\right) (1_{A}\otimes g)(1_{A}\otimes gx_{2}) \\
(X_{1}\otimes 1_{H})(1_{A}\otimes x_{1})(1_{A}\otimes g)B\left(
x_{2}hg\otimes 1_{H}\right) (1_{A}\otimes g) \\
(X_{1}\otimes 1_{H})(1_{A}\otimes g)(1_{A}\otimes g)B\left( x_{2}hg\otimes
1_{H}\right) (1_{A}\otimes g)(1_{A}\otimes gx_{1}) \\
(X_{1}\otimes 1_{H})B(x_{1}x_{2}hg\otimes g) \\
+(1_{A}\otimes gx_{2})B(hx_{1}\otimes g)(1_{A}\otimes gx_{1})+(1_{A}\otimes
gx_{2})B(h\otimes x_{1})(1_{A}\otimes gx_{1}) \\
-(1_{A}\otimes gx_{1})B(hx_{1}\otimes g)(1_{A}\otimes gx_{2})-(1_{A}\otimes
gx_{1})B(h\otimes x_{1})(1_{A}\otimes gx_{2}) \\
+(1_{A}\otimes x_{1}x_{2})B(hx_{1}\otimes g)+(1_{A}\otimes
x_{1}x_{2})B(h\otimes x_{1})+ \\
-(1_{A}\otimes x_{2})B(gx_{1}h\otimes x_{1}) \\
-(1_{A}\otimes g)B(gx_{1}h\otimes x_{1})(1_{A}\otimes gx_{2})+ \\
+(1_{A}\otimes x_{1})B(gx_{2}hx_{1}\otimes g)+(1_{A}\otimes
x_{1})B(gx_{2}h\otimes x_{1})+ \\
+(1_{A}\otimes g)B(gx_{2}hx_{1}\otimes g)(1_{A}\otimes gx_{1})+(1_{A}\otimes
g)B(gx_{2}h\otimes x_{1})(1_{A}\otimes gx_{1}) \\
+B(hx_{1}\otimes g)(1_{A}\otimes x_{1}x_{2})+B(h\otimes x_{1})(1_{A}\otimes
x_{1}x_{2}) \\
+B(x_{1}x_{2}h\otimes x_{1})
\end{gather*}%
We use here%
\begin{equation*}
B\left( hg\otimes g\right) =(1_{A}\otimes g)B\left( ghg\otimes 1_{H}\right)
(1_{A}\otimes g)
\end{equation*}%
\begin{equation*}
B\left( gx_{1}hg\otimes g\right) =(1_{A}\otimes g)B\left( x_{1}hg\otimes
1_{H}\right) (1_{A}\otimes g)
\end{equation*}%
\begin{equation*}
B\left( gx_{2}hg\otimes g\right) =(1_{A}\otimes g)B\left( x_{2}hg\otimes
1_{H}\right) (1_{A}\otimes g)
\end{equation*}%
\begin{eqnarray*}
&&\text{part of right side }(X_{1}\otimes 1_{H})B(hg\otimes gx_{1}x_{2}) \\
&=&\text{1}(X_{1}\otimes 1_{H})(1_{A}\otimes g)B(ghg\otimes
1_{H})(1_{A}\otimes gx_{1}x_{2}) \\
&&\text{2}-(X_{1}\otimes 1_{H})(1_{A}\otimes x_{2})B(ghg\otimes
1_{H})(1_{A}\otimes x_{1}) \\
&&\text{7}-(X_{1}\otimes 1_{H})B(x_{2}hg\otimes 1_{H})(1_{A}\otimes x_{1}) \\
&&\text{3}+(X_{1}\otimes 1_{H})(1_{A}\otimes x_{1})B(ghg\otimes
1_{H})(1_{A}\otimes x_{2}) \\
&&\text{4}+(X_{1}\otimes 1_{H})(1_{A}\otimes gx_{1}x_{2})B(ghg\otimes
1_{H})(1_{A}\otimes g) \\
&&\text{8}+(X_{1}\otimes 1_{H})(1_{A}\otimes gx_{1})B(x_{2}hg\otimes
1_{H})(1_{A}\otimes g) \\
&&\text{5}+(X_{1}\otimes 1_{H})B(x_{1}hg\otimes 1_{H})(1_{A}\otimes x_{2}) \\
&&\text{6}-(X_{1}\otimes 1_{H})(1_{A}\otimes gx_{2})B(x_{1}hg\otimes
1_{H})(1_{A}\otimes g) \\
&&\text{9}+(X_{1}\otimes 1_{H})(1_{A}\otimes g)B(gx_{1}x_{2}hg\otimes
1_{H})(1_{A}\otimes g)
\end{eqnarray*}%
\begin{eqnarray*}
&&\text{part of left side with assumptions} \\
&=&\text{1}(X_{1}\otimes 1_{H})(1_{A}\otimes g)B\left( ghg\otimes
1_{H}\right) (1_{A}\otimes g)(1_{A}\otimes x_{1}x_{2}) \\
&&\text{2}(X_{1}\otimes 1_{H})(1_{A}\otimes gx_{2})(1_{A}\otimes g)B\left(
ghg\otimes 1_{H}\right) (1_{A}\otimes g)(1_{A}\otimes gx_{1}) \\
&&\text{3}-(X_{1}\otimes 1_{H})(1_{A}\otimes gx_{1})(1_{A}\otimes g)B\left(
ghg\otimes 1_{H}\right) (1_{A}\otimes g)(1_{A}\otimes gx_{2}) \\
&&\text{4}+(X_{1}\otimes 1_{H})(1_{A}\otimes x_{1}x_{2})(1_{A}\otimes
g)B\left( ghg\otimes 1_{H}\right) (1_{A}\otimes g) \\
&&\text{6}-(X_{1}\otimes 1_{H})(1_{A}\otimes x_{2})(1_{A}\otimes g)B\left(
x_{1}hg\otimes 1_{H}\right) (1_{A}\otimes g) \\
&&\text{5}-(X_{1}\otimes 1_{H})(1_{A}\otimes g)(1_{A}\otimes g)B\left(
x_{1}hg\otimes 1_{H}\right) (1_{A}\otimes g)(1_{A}\otimes gx_{2}) \\
&&\text{8}(X_{1}\otimes 1_{H})(1_{A}\otimes x_{1})(1_{A}\otimes g)B\left(
x_{2}hg\otimes 1_{H}\right) (1_{A}\otimes g) \\
&&\text{7}(X_{1}\otimes 1_{H})(1_{A}\otimes g)(1_{A}\otimes g)B\left(
x_{2}hg\otimes 1_{H}\right) (1_{A}\otimes g)(1_{A}\otimes gx_{1}) \\
&&\text{9}(X_{1}\otimes 1_{H})B(x_{1}x_{2}hg\otimes g)
\end{eqnarray*}%
\begin{eqnarray*}
\text{last part of right side} &&\text{6}+(1_{A}\otimes g)B(ghx_{1}\otimes
1_{H})(1_{A}\otimes gx_{1}x_{2}) \\
&&\text{1}-(1_{A}\otimes x_{2})B(ghx_{1}\otimes 1_{H})(1_{A}\otimes x_{1}) \\
&&\text{3}-B(x_{2}hx_{1}\otimes 1_{H})(1_{A}\otimes x_{1}) \\
&&\text{2}+(1_{A}\otimes x_{1})B(ghx_{1}\otimes 1_{H})(1_{A}\otimes x_{2})+
\\
&&\text{4}+(1_{A}\otimes gx_{1}x_{2})B(ghx_{1}\otimes 1_{H})(1_{A}\otimes g)
\\
&&\text{5}+(1_{A}\otimes gx_{1})B(x_{2}hx_{1}\otimes 1_{H})(1_{A}\otimes g)
\end{eqnarray*}%
\begin{eqnarray*}
\text{last part of left+assum.} &&\text{1}(1_{A}\otimes
gx_{2})B(hx_{1}\otimes g)(1_{A}\otimes gx_{1}) \\
&&+(1_{A}\otimes gx_{2})B(h\otimes x_{1})(1_{A}\otimes gx_{1}) \\
&&\text{2}-(1_{A}\otimes gx_{1})B(hx_{1}\otimes g)(1_{A}\otimes gx_{2}) \\
&&-(1_{A}\otimes gx_{1})B(h\otimes x_{1})(1_{A}\otimes gx_{2}) \\
&&\text{4}+(1_{A}\otimes x_{1}x_{2})B(hx_{1}\otimes g) \\
&&+(1_{A}\otimes x_{1}x_{2})B(h\otimes x_{1})+ \\
&&-(1_{A}\otimes x_{2})B(gx_{1}h\otimes x_{1}) \\
&&-(1_{A}\otimes g)B(gx_{1}h\otimes x_{1})(1_{A}\otimes gx_{2})+ \\
&&\text{5}+(1_{A}\otimes x_{1})B(gx_{2}hx_{1}\otimes g) \\
&&+(1_{A}\otimes x_{1})B(gx_{2}h\otimes x_{1})+ \\
&&\text{3}+(1_{A}\otimes g)B(gx_{2}hx_{1}\otimes g)(1_{A}\otimes gx_{1}) \\
&&+(1_{A}\otimes g)B(gx_{2}h\otimes x_{1})(1_{A}\otimes gx_{1}) \\
&&\text{6}+B(hx_{1}\otimes g)(1_{A}\otimes x_{1}x_{2}) \\
&&+B(h\otimes x_{1})(1_{A}\otimes x_{1}x_{2}) \\
&&+B(x_{1}x_{2}h\otimes x_{1})
\end{eqnarray*}%
Thus we have to show that the following summation is zero.%
\begin{eqnarray*}
&&\text{1}+(1_{A}\otimes gx_{2})B(h\otimes x_{1})(1_{A}\otimes gx_{1}) \\
&&\text{2}-(1_{A}\otimes gx_{1})B(h\otimes x_{1})(1_{A}\otimes gx_{2}) \\
&&\text{3}+(1_{A}\otimes x_{1}x_{2})B(h\otimes x_{1})+ \\
&&\text{5}-(1_{A}\otimes x_{2})B(gx_{1}h\otimes x_{1}) \\
&&\text{6}-(1_{A}\otimes g)B(gx_{1}h\otimes x_{1})(1_{A}\otimes gx_{2})+ \\
&&\text{7}+(1_{A}\otimes x_{1})B(gx_{2}h\otimes x_{1})+ \\
&&\text{8}+(1_{A}\otimes g)B(gx_{2}h\otimes x_{1})(1_{A}\otimes gx_{1}) \\
&&\text{4}+B(h\otimes x_{1})(1_{A}\otimes x_{1}x_{2}) \\
&&\text{9}+B(x_{1}x_{2}h\otimes x_{1})
\end{eqnarray*}%
We use $\left( \ref{simplx}\right) $%
\begin{equation*}
B(h\otimes x_{1})=B(h\otimes 1_{H})(1_{A}\otimes x_{1})-(1_{A}\otimes
gx_{1})B(h\otimes 1_{H})(1_{A}\otimes g)-(1_{A}\otimes g)B(gx_{1}h\otimes
1_{H})(1_{A}\otimes g)
\end{equation*}%
\begin{equation*}
B(gx_{1}h\otimes x_{1})=B(gx_{1}h\otimes 1_{H})(1_{A}\otimes
x_{1})-(1_{A}\otimes gx_{1})B(gx_{1}h\otimes 1_{H})(1_{A}\otimes g)
\end{equation*}%
\begin{equation*}
B(gx_{2}h\otimes x_{1})=B(gx_{2}h\otimes 1_{H})(1_{A}\otimes
x_{1})-(1_{A}\otimes gx_{1})B(gx_{2}h\otimes 1_{H})(1_{A}\otimes
g)+(1_{A}\otimes g)B(x_{1}x_{2}h\otimes 1_{H})(1_{A}\otimes g)
\end{equation*}%
\begin{equation*}
\text{1}(1_{A}\otimes gx_{2})B(h\otimes x_{1})(1_{A}\otimes gx_{1})=\text{1,1%
}-(1_{A}\otimes x_{1}x_{2})B(h\otimes 1_{H})(1_{A}\otimes x_{1})+\text{1,2}%
(1_{A}\otimes x_{2})B(gx_{1}h\otimes 1_{H})(1_{A}\otimes x_{1})
\end{equation*}%
\begin{equation*}
\text{2}-(1_{A}\otimes gx_{1})B(h\otimes x_{1})(1_{A}\otimes gx_{2})=\text{%
2,1}+(1_{A}\otimes gx_{1})B(h\otimes 1_{H})(1_{A}\otimes gx_{1}x_{2})\text{%
2,2}-(1_{A}\otimes x_{1})B(gx_{1}h\otimes 1_{H})(1_{A}\otimes x_{2})
\end{equation*}%
\begin{equation*}
\text{3}(1_{A}\otimes x_{1}x_{2})B(h\otimes x_{1})=\text{1,1}(1_{A}\otimes
x_{1}x_{2})B(h\otimes 1_{H})(1_{A}\otimes x_{1})\text{3,1}-(1_{A}\otimes
gx_{1}x_{2})B(gx_{1}h\otimes 1_{H})(1_{A}\otimes g)
\end{equation*}%
\begin{equation*}
\text{4}B(h\otimes x_{1})(1_{A}\otimes x_{1}x_{2})=\text{2,1}-(1_{A}\otimes
gx_{1})B(h\otimes 1_{H})(1_{A}\otimes gx_{1}x_{2})\text{4,1}-(1_{A}\otimes
g)B(gx_{1}h\otimes 1_{H})(1_{A}\otimes gx_{1}x_{2})
\end{equation*}%
\begin{equation*}
\text{5}-(1_{A}\otimes x_{2})B(gx_{1}h\otimes x_{1})=\text{1,2}%
-(1_{A}\otimes x_{2})B(gx_{1}h\otimes 1_{H})(1_{A}\otimes x_{1})+\text{3,1}%
(1_{A}\otimes gx_{1}x_{2})B(gx_{1}h\otimes 1_{H})(1_{A}\otimes g)
\end{equation*}%
\begin{equation*}
\text{6}-(1_{A}\otimes g)B(gx_{1}h\otimes x_{1})(1_{A}\otimes gx_{2})=\text{%
4,1}+(1_{A}\otimes g)B(gx_{1}h\otimes 1_{H})(1_{A}\otimes gx_{1}x_{2})+\text{%
2,2}(1_{A}\otimes x_{1})B(gx_{1}h\otimes 1_{H})(1_{A}\otimes x_{2})
\end{equation*}%
\begin{equation*}
\text{7}+(1_{A}\otimes x_{1})B(gx_{2}h\otimes x_{1})=\text{7,1}(1_{A}\otimes
x_{1})B(gx_{2}h\otimes 1_{H})(1_{A}\otimes x_{1})+\text{7,2}(1_{A}\otimes
gx_{1})B(x_{1}x_{2}h\otimes 1_{H})(1_{A}\otimes g)
\end{equation*}%
\begin{equation*}
\text{8}+(1_{A}\otimes g)B(gx_{2}h\otimes x_{1})(1_{A}\otimes gx_{1})=\text{%
7,1}-(1_{A}\otimes x_{1})B(gx_{2}h\otimes 1_{H})(1_{A}\otimes x_{1})\text{8,1%
}-B(x_{1}x_{2}h\otimes 1_{H})(1_{A}\otimes x_{1})
\end{equation*}%
\begin{equation*}
9+B(x_{1}x_{2}h\otimes x_{1})=\text{8,1}B(x_{1}x_{2}h\otimes
1_{H})(1_{A}\otimes x_{1})\text{7,2}-(1_{A}\otimes
gx_{1})B(x_{1}x_{2}h\otimes 1_{H})(1_{A}\otimes g)
\end{equation*}%
Therefore our claim is proved.

Now we have to do this for $X_{2}.$

We have to prove that for any $h$ the formula below holds%
\begin{equation*}
B(h\otimes x_{1}x_{2})(X_{2}\otimes 1_{H})=(X_{2}\otimes 1_{H})B(hg\otimes
gx_{1}x_{2})+B(hx_{2}\otimes gx_{1}x_{2})
\end{equation*}

\begin{eqnarray*}
B(h\otimes x_{1}x_{2})\left( X_{2}\otimes 1_{H}\right) &=&B(h\otimes
1_{H})\left( X_{2}\otimes 1_{H}\right) (1_{A}\otimes x_{1}x_{2}) \\
&&+(1_{A}\otimes gx_{2})B(h\otimes 1_{H})\left( X_{2}\otimes 1_{H}\right)
(1_{A}\otimes gx_{1}) \\
&&-(1_{A}\otimes gx_{1})B(h\otimes 1_{H})\left( X_{2}\otimes 1_{H}\right)
(1_{A}\otimes gx_{2}) \\
&&+(1_{A}\otimes x_{1}x_{2})B(h\otimes 1_{H})\left( X_{2}\otimes 1_{H}\right)
\\
&&-(1_{A}\otimes x_{2})B(gx_{1}h\otimes 1_{H})\left( X_{2}\otimes
1_{H}\right) \\
&&-(1_{A}\otimes g)B(gx_{1}h\otimes 1_{H})\left( X_{2}\otimes 1_{H}\right)
(1_{A}\otimes gx_{2}) \\
&&+(1_{A}\otimes x_{1})B(gx_{2}h\otimes 1_{H})\left( X_{2}\otimes
1_{H}\right) \\
&&+(1_{A}\otimes g)B(gx_{2}h\otimes 1_{H})\left( X_{2}\otimes 1_{H}\right)
(1_{A}\otimes gx_{1}) \\
&&+B(x_{1}x_{2}h\otimes 1_{H})\left( X_{2}\otimes 1_{H}\right)
\end{eqnarray*}%
Now we apply the assumption that the formula holds for $B(h\otimes 1_{H})$
i.e.%
\begin{equation*}
B(h\otimes 1_{H})(X_{2}\otimes 1_{H})=(X_{2}\otimes 1_{H})B(hg\otimes
g)+B(hx_{2}\otimes g)+B(h\otimes x_{2})
\end{equation*}%
\begin{eqnarray*}
&&\text{1}B(h\otimes 1_{H})\left( X_{2}\otimes 1_{H}\right) (1_{A}\otimes
x_{1}x_{2}) \\
&&+(1_{A}\otimes gx_{2})B(h\otimes 1_{H})\left( X_{2}\otimes 1_{H}\right)
(1_{A}\otimes gx_{1}) \\
&&-(1_{A}\otimes gx_{1})B(h\otimes 1_{H})\left( X_{2}\otimes 1_{H}\right)
(1_{A}\otimes gx_{2}) \\
&&+(1_{A}\otimes x_{1}x_{2})B(h\otimes 1_{H})\left( X_{2}\otimes 1_{H}\right)
\\
&=&\left[ (X_{2}\otimes 1_{H})B(hg\otimes g)+B(hx_{2}\otimes g)+B(h\otimes
x_{2})\right] (1_{A}\otimes x_{1}x_{2}) \\
&&(1_{A}\otimes gx_{2})\left[ (X_{2}\otimes 1_{H})B(hg\otimes
g)+B(hx_{2}\otimes g)+B(h\otimes x_{2})\right] (1_{A}\otimes gx_{1}) \\
&&-(1_{A}\otimes gx_{1})\left[ (X_{2}\otimes 1_{H})B(hg\otimes
g)+B(hx_{2}\otimes g)+B(h\otimes x_{2})\right] (1_{A}\otimes gx_{2}) \\
&&+(1_{A}\otimes x_{1}x_{2})\left[ (X_{2}\otimes 1_{H})B(hg\otimes
g)+B(hx_{2}\otimes g)+B(h\otimes x_{2})\right]
\end{eqnarray*}%
\begin{equation*}
B(gx_{1}h\otimes 1_{H})(X_{2}\otimes 1_{H})=(X_{2}\otimes
1_{H})B(gx_{1}hg\otimes g)+B(gx_{1}hx_{2}\otimes g)+B(gx_{1}h\otimes x_{2})
\end{equation*}%
\begin{eqnarray*}
&&\text{2}-(1_{A}\otimes x_{2})B(gx_{1}h\otimes 1_{H})\left( X_{2}\otimes
1_{H}\right) \\
&&-(1_{A}\otimes g)B(gx_{1}h\otimes 1_{H})\left( X_{2}\otimes 1_{H}\right)
(1_{A}\otimes gx_{2}) \\
&=&-(1_{A}\otimes x_{2})\left[ (X_{2}\otimes 1_{H})B(gx_{1}hg\otimes
g)+B(gx_{1}hx_{2}\otimes g)+B(gx_{1}h\otimes x_{2})\right] + \\
&&-(1_{A}\otimes g)\left[ (X_{2}\otimes 1_{H})B(gx_{1}hg\otimes
g)+B(gx_{1}hx_{2}\otimes g)+B(gx_{1}h\otimes x_{2})\right] (1_{A}\otimes
gx_{2})
\end{eqnarray*}%
\begin{equation*}
B(gx_{2}h\otimes 1_{H})(X_{2}\otimes 1_{H})=(X_{2}\otimes
1_{H})B(gx_{2}hg\otimes g)+B(gx_{2}h\otimes x_{2})
\end{equation*}%
\begin{eqnarray*}
&&\text{3}+(1_{A}\otimes x_{1})B(gx_{2}h\otimes 1_{H})\left( X_{2}\otimes
1_{H}\right) \\
&&+(1_{A}\otimes g)B(gx_{2}h\otimes 1_{H})\left( X_{2}\otimes 1_{H}\right)
(1_{A}\otimes gx_{1}) \\
&=&+(1_{A}\otimes x_{1})\left[ (X_{2}\otimes 1_{H})B(gx_{2}hg\otimes
g)+B(gx_{2}h\otimes x_{2})\right] + \\
&&+(1_{A}\otimes g)\left[ (X_{2}\otimes 1_{H})B(gx_{2}hg\otimes
g)+B(gx_{2}h\otimes x_{2})\right] (1_{A}\otimes gx_{1})
\end{eqnarray*}%
\begin{equation*}
\text{4}B(x_{1}x_{2}h\otimes 1_{H})(X_{2}\otimes 1_{H})=(X_{2}\otimes
1_{H})B(x_{1}x_{2}hg\otimes g)+B(x_{1}x_{2}h\otimes x_{2})
\end{equation*}%
Therefore the left side, using the assumption that for any $h$ $\left( \ref%
{eq.a}\right) $ holds and $d=X_{2}$ gives us%
\begin{eqnarray*}
&&\text{left side plus assumptions}\left[ (X_{2}\otimes 1_{H})B(hg\otimes
g)+B(hx_{2}\otimes g)+B(h\otimes x_{2})\right] (1_{A}\otimes x_{1}x_{2}) \\
&&(1_{A}\otimes gx_{2})\left[ (X_{2}\otimes 1_{H})B(hg\otimes
g)+B(hx_{2}\otimes g)+B(h\otimes x_{2})\right] (1_{A}\otimes gx_{1}) \\
&&-(1_{A}\otimes gx_{1})\left[ (X_{2}\otimes 1_{H})B(hg\otimes
g)+B(hx_{2}\otimes g)+B(h\otimes x_{2})\right] (1_{A}\otimes gx_{2}) \\
&&+(1_{A}\otimes x_{1}x_{2})\left[ (X_{2}\otimes 1_{H})B(hg\otimes
g)+B(hx_{2}\otimes g)+B(h\otimes x_{2})\right] + \\
&&-(1_{A}\otimes x_{2})\left[ (X_{2}\otimes 1_{H})B(gx_{1}hg\otimes
g)+B(gx_{1}hx_{2}\otimes g)+B(gx_{1}h\otimes x_{2})\right] + \\
&&-(1_{A}\otimes g)\left[ (X_{2}\otimes 1_{H})B(gx_{1}hg\otimes
g)+B(gx_{1}hx_{2}\otimes g)+B(gx_{1}h\otimes x_{2})\right] (1_{A}\otimes
gx_{2})+ \\
&&+(1_{A}\otimes x_{1})\left[ (X_{2}\otimes 1_{H})B(gx_{2}hg\otimes
g)+B(gx_{2}h\otimes x_{2})\right] + \\
&&+(1_{A}\otimes g)\left[ (X_{2}\otimes 1_{H})B(gx_{2}hg\otimes
g)+B(gx_{2}h\otimes x_{2})\right] (1_{A}\otimes gx_{1})+ \\
&&(X_{2}\otimes 1_{H})B(x_{1}x_{2}hg\otimes g)+B(x_{1}x_{2}h\otimes x_{2})
\end{eqnarray*}%
We rewrite this separating the part containing $(X_{2}\otimes 1_{H})$%
\begin{eqnarray*}
\text{new } &&\text{left side plus assumptions}(X_{2}\otimes
1_{H})B(hg\otimes g)(1_{A}\otimes x_{1}x_{2}) \\
&&(X_{2}\otimes 1_{H})(1_{A}\otimes gx_{2})B(hg\otimes g)(1_{A}\otimes
gx_{1}) \\
&&-(X_{2}\otimes 1_{H})(1_{A}\otimes gx_{1})B(hg\otimes g)(1_{A}\otimes
gx_{2}) \\
&&+(X_{2}\otimes 1_{H})(1_{A}\otimes x_{1}x_{2})B(hg\otimes g) \\
&&(-X_{2}\otimes 1_{H})(1_{A}\otimes x_{2})B(gx_{1}hg\otimes g) \\
&&-(X_{2}\otimes 1_{H})(1_{A}\otimes g)B(gx_{1}hg\otimes g)(1_{A}\otimes
gx_{2}) \\
&&+(X_{2}\otimes 1_{H})(1_{A}\otimes x_{1})B(gx_{2}hg\otimes g) \\
&&+(X_{2}\otimes 1_{H})(1_{A}\otimes g)B(gx_{2}hg\otimes g)(1_{A}\otimes
gx_{1}) \\
&&+(X_{2}\otimes 1_{H})B(x_{1}x_{2}hg\otimes g) \\
&&+B(h\otimes x_{2})(1_{A}\otimes x_{1}x_{2}) \\
&&+(1_{A}\otimes gx_{2})B(h\otimes x_{2})(1_{A}\otimes gx_{1}) \\
&&-(1_{A}\otimes gx_{1})B(h\otimes x_{2})(1_{A}\otimes gx_{2}) \\
&&+(1_{A}\otimes x_{1}x_{2})B(h\otimes x_{2}) \\
&&+B(hx_{2}\otimes g)(1_{A}\otimes x_{1}x_{2}) \\
&&+(1_{A}\otimes gx_{2})B(hx_{2}\otimes g)(1_{A}\otimes gx_{1}) \\
&&-(1_{A}\otimes gx_{1})B(hx_{2}\otimes g)(1_{A}\otimes gx_{2}) \\
&&+(1_{A}\otimes x_{1}x_{2})B(hx_{2}\otimes g)+ \\
&&-(1_{A}\otimes x_{2})B(gx_{1}hx_{2}\otimes g) \\
&&-(1_{A}\otimes g)B(gx_{1}hx_{2}\otimes g)(1_{A}\otimes gx_{2}) \\
&&-(1_{A}\otimes x_{2})B(gx_{1}h\otimes x_{2})+ \\
&&-(1_{A}\otimes g)B(gx_{1}h\otimes x_{2})(1_{A}\otimes gx_{2})+ \\
&&+(1_{A}\otimes x_{1})B(gx_{2}h\otimes x_{2})+ \\
&&+(1_{A}\otimes g)B(gx_{2}h\otimes x_{2})(1_{A}\otimes gx_{1})+ \\
&&+B(x_{1}x_{2}h\otimes x_{2})
\end{eqnarray*}%
Now we consider the right side
\begin{equation*}
\text{right side}(X_{2}\otimes 1_{H})B(hg\otimes
gx_{1}x_{2})+B(hx_{2}\otimes gx_{1}x_{2})
\end{equation*}%
We use $\left( \ref{simplgxx}\right) $%
\begin{eqnarray*}
B(h\otimes gx_{1}x_{2}) &=&(1_{A}\otimes g)B(ghg\otimes 1_{H})(1_{A}\otimes
gx_{1}x_{2}) \\
&&-(1_{A}\otimes x_{2})B(ghg\otimes 1_{H})(1_{A}\otimes x_{1}) \\
&&-B(x_{2}hg\otimes 1_{H})(1_{A}\otimes x_{1}) \\
&&+(1_{A}\otimes x_{1})B(ghg\otimes 1_{H})(1_{A}\otimes x_{2}) \\
&&+(1_{A}\otimes gx_{1}x_{2})B(ghg\otimes 1_{H})(1_{A}\otimes g) \\
&&+(1_{A}\otimes gx_{1})B(x_{2}hg\otimes 1_{H})(1_{A}\otimes g) \\
&&+B(x_{1}hg\otimes 1_{H})(1_{A}\otimes x_{2}) \\
&&-(1_{A}\otimes gx_{2})B(x_{1}hg\otimes 1_{H})(1_{A}\otimes g) \\
&&+(1_{A}\otimes g)B(gx_{1}x_{2}hg\otimes 1_{H})(1_{A}\otimes g)
\end{eqnarray*}%
Therefore%
\begin{eqnarray*}
&&\text{right side}(X_{2}\otimes 1_{H})B(hg\otimes gx_{1}x_{2}) \\
&=&(X_{2}\otimes 1_{H})(1_{A}\otimes g)B(ghg\otimes 1_{H})(1_{A}\otimes
gx_{1}x_{2}) \\
&&-(X_{2}\otimes 1_{H})(1_{A}\otimes x_{2})B(ghg\otimes 1_{H})(1_{A}\otimes
x_{1}) \\
&&-(X_{2}\otimes 1_{H})B(x_{2}hg\otimes 1_{H})(1_{A}\otimes x_{1}) \\
&&+(X_{2}\otimes 1_{H})(1_{A}\otimes x_{1})B(ghg\otimes 1_{H})(1_{A}\otimes
x_{2}) \\
&&+(X_{2}\otimes 1_{H})(1_{A}\otimes gx_{1}x_{2})B(ghg\otimes
1_{H})(1_{A}\otimes g) \\
&&+(X_{2}\otimes 1_{H})(1_{A}\otimes gx_{1})B(x_{2}hg\otimes
1_{H})(1_{A}\otimes g) \\
&&+(X_{2}\otimes 1_{H})B(x_{1}hg\otimes 1_{H})(1_{A}\otimes x_{2}) \\
&&-(X_{2}\otimes 1_{H})(1_{A}\otimes gx_{2})B(x_{1}hg\otimes
1_{H})(1_{A}\otimes g) \\
&&+(X_{2}\otimes 1_{H})(1_{A}\otimes g)B(gx_{1}x_{2}hg\otimes
1_{H})(1_{A}\otimes g)
\end{eqnarray*}%
Now we compare it with the part of the rewritten left side containing $%
X_{2}\otimes 1_{H}$%
\begin{eqnarray*}
\text{new } &&\text{left side plus assumptions containing }X_{2}\otimes 1_{H}
\\
&&\text{1}(X_{2}\otimes 1_{H})B(hg\otimes g)(1_{A}\otimes x_{1}x_{2}) \\
\text{2} &&(X_{2}\otimes 1_{H})(1_{A}\otimes gx_{2})B(hg\otimes
g)(1_{A}\otimes gx_{1}) \\
&&\text{3}-(X_{2}\otimes 1_{H})(1_{A}\otimes gx_{1})B(hg\otimes
g)(1_{A}\otimes gx_{2}) \\
&&\text{4}+(X_{2}\otimes 1_{H})(1_{A}\otimes x_{1}x_{2})B(hg\otimes g) \\
&&\text{5}(-X_{2}\otimes 1_{H})(1_{A}\otimes x_{2})B(gx_{1}hg\otimes g) \\
&&\text{6}-(X_{2}\otimes 1_{H})(1_{A}\otimes g)B(gx_{1}hg\otimes
g)(1_{A}\otimes gx_{2}) \\
&&\text{7}+(X_{2}\otimes 1_{H})(1_{A}\otimes x_{1})B(gx_{2}hg\otimes g) \\
&&\text{8}+(X_{2}\otimes 1_{H})(1_{A}\otimes g)B(gx_{2}hg\otimes
g)(1_{A}\otimes gx_{1}) \\
&&\text{9}+(X_{2}\otimes 1_{H})B(x_{1}x_{2}hg\otimes g)
\end{eqnarray*}%
\begin{eqnarray*}
&&\text{right side}(X_{2}\otimes 1_{H})B(hg\otimes gx_{1}x_{2}) \\
&=&\text{1}(X_{2}\otimes 1_{H})(1_{A}\otimes g)B(ghg\otimes
1_{H})(1_{A}\otimes gx_{1}x_{2}) \\
&&\text{2}-(X_{2}\otimes 1_{H})(1_{A}\otimes x_{2})B(ghg\otimes
1_{H})(1_{A}\otimes x_{1}) \\
&&\text{8}-(X_{2}\otimes 1_{H})B(x_{2}hg\otimes 1_{H})(1_{A}\otimes x_{1}) \\
&&\text{3}+(X_{2}\otimes 1_{H})(1_{A}\otimes x_{1})B(ghg\otimes
1_{H})(1_{A}\otimes x_{2}) \\
&&\text{4}+(X_{2}\otimes 1_{H})(1_{A}\otimes gx_{1}x_{2})B(ghg\otimes
1_{H})(1_{A}\otimes g) \\
&&\text{7}+(X_{2}\otimes 1_{H})(1_{A}\otimes gx_{1})B(x_{2}hg\otimes
1_{H})(1_{A}\otimes g) \\
&&\text{6}+(X_{2}\otimes 1_{H})B(x_{1}hg\otimes 1_{H})(1_{A}\otimes x_{2}) \\
&&\text{5}-(X_{2}\otimes 1_{H})(1_{A}\otimes gx_{2})B(x_{1}hg\otimes
1_{H})(1_{A}\otimes g) \\
&&\text{9}+(X_{2}\otimes 1_{H})(1_{A}\otimes g)B(gx_{1}x_{2}hg\otimes
1_{H})(1_{A}\otimes g)
\end{eqnarray*}%
Now we write the right part $B(hx_{2}\otimes gx_{1}x_{2})$ by using use $%
\left( \ref{simplgxx}\right) $
\begin{eqnarray*}
B(hx_{2}\otimes gx_{1}x_{2}) &=&\text{6}(1_{A}\otimes g)B(ghx_{2}\otimes
1_{H})(1_{A}\otimes gx_{1}x_{2}) \\
&&\text{1}-(1_{A}\otimes x_{2})B(ghx_{2}\otimes 1_{H})(1_{A}\otimes x_{1}) \\
&&\text{2}+(1_{A}\otimes x_{1})B(ghx_{2}\otimes 1_{H})(1_{A}\otimes x_{2}) \\
&&\text{3}+(1_{A}\otimes gx_{1}x_{2})B(ghx_{2}\otimes 1_{H})(1_{A}\otimes g)
\\
&&\text{5}+B(x_{1}hx_{2}\otimes 1_{H})(1_{A}\otimes x_{2}) \\
&&\text{4}-(1_{A}\otimes gx_{2})B(x_{1}hx_{2}\otimes 1_{H})(1_{A}\otimes g)
\end{eqnarray*}%
and we compare with what it remains from the left side%
\begin{eqnarray*}
&&+B(h\otimes x_{2})(1_{A}\otimes x_{1}x_{2}) \\
&&+(1_{A}\otimes gx_{2})B(h\otimes x_{2})(1_{A}\otimes gx_{1}) \\
&&-(1_{A}\otimes gx_{1})B(h\otimes x_{2})(1_{A}\otimes gx_{2}) \\
&&+(1_{A}\otimes x_{1}x_{2})B(h\otimes x_{2}) \\
&&\text{6}+B(hx_{2}\otimes g)(1_{A}\otimes x_{1}x_{2}) \\
&&\text{1}+(1_{A}\otimes gx_{2})B(hx_{2}\otimes g)(1_{A}\otimes gx_{1}) \\
&&\text{2}-(1_{A}\otimes gx_{1})B(hx_{2}\otimes g)(1_{A}\otimes gx_{2}) \\
&&\text{3}+(1_{A}\otimes x_{1}x_{2})B(hx_{2}\otimes g)+ \\
&&\text{4}-(1_{A}\otimes x_{2})B(gx_{1}hx_{2}\otimes g) \\
&&\text{5}-(1_{A}\otimes g)B(gx_{1}hx_{2}\otimes g)(1_{A}\otimes gx_{2}) \\
&&-(1_{A}\otimes x_{2})B(gx_{1}h\otimes x_{2})+ \\
&&-(1_{A}\otimes g)B(gx_{1}h\otimes x_{2})(1_{A}\otimes gx_{2})+ \\
&&+(1_{A}\otimes x_{1})B(gx_{2}h\otimes x_{2})+ \\
&&+(1_{A}\otimes g)B(gx_{2}h\otimes x_{2})(1_{A}\otimes gx_{1})+ \\
&&+B(x_{1}x_{2}h\otimes x_{2})
\end{eqnarray*}%
Now we have to prove the following summation is zero%
\begin{eqnarray*}
&&\text{1}+B(h\otimes x_{2})(1_{A}\otimes x_{1}x_{2}) \\
\text{2} &&+(1_{A}\otimes gx_{2})B(h\otimes x_{2})(1_{A}\otimes gx_{1}) \\
&&\text{3}-(1_{A}\otimes gx_{1})B(h\otimes x_{2})(1_{A}\otimes gx_{2}) \\
\text{4} &&+(1_{A}\otimes x_{1}x_{2})B(h\otimes x_{2}) \\
\text{5} &&-(1_{A}\otimes x_{2})B(gx_{1}h\otimes x_{2})+ \\
\text{6} &&-(1_{A}\otimes g)B(gx_{1}h\otimes x_{2})(1_{A}\otimes gx_{2})+ \\
\text{7} &&+(1_{A}\otimes x_{1})B(gx_{2}h\otimes x_{2})+ \\
\text{8} &&+(1_{A}\otimes g)B(gx_{2}h\otimes x_{2})(1_{A}\otimes gx_{1})+ \\
\text{9} &&+B(x_{1}x_{2}h\otimes x_{2})
\end{eqnarray*}%
Now we use $\left( \ref{simpl2i2}\right) $%
\begin{eqnarray*}
B(h\otimes x_{2}) &=&B(h\otimes 1_{H})(1_{A}\otimes x_{2})-(1_{A}\otimes
gx_{2})B(h\otimes 1_{H})(1_{A}\otimes g)-(1_{A}\otimes g)B(gx_{2}h\otimes
1_{H}) \\
&&(1_{A}\otimes g)
\end{eqnarray*}%
\begin{eqnarray*}
B(gx_{1}h\otimes x_{2}) &=&B(gx_{1}h\otimes 1_{H})(1_{A}\otimes
x_{2})-(1_{A}\otimes gx_{2})B(gx_{1}h\otimes 1_{H})(1_{A}\otimes
g)-(1_{A}\otimes g)B(gx_{1}gx_{2}h\otimes 1_{H}) \\
&&(1_{A}\otimes g)
\end{eqnarray*}%
\begin{equation*}
B(gx_{2}h\otimes x_{2})=B(gx_{2}h\otimes 1_{H})(1_{A}\otimes
x_{2})-(1_{A}\otimes gx_{2})B(gx_{2}h\otimes 1_{H})(1_{A}\otimes g)
\end{equation*}%
\begin{equation*}
B(x_{1}x_{2}h\otimes x_{2})=B(x_{1}x_{2}h\otimes 1_{H})(1_{A}\otimes
x_{2})-(1_{A}\otimes gx_{2})B(x_{1}x_{2}h\otimes 1_{H})(1_{A}\otimes g)
\end{equation*}%
and we get%
\begin{eqnarray*}
&&\text{1}+B(h\otimes x_{2})(1_{A}\otimes x_{1}x_{2}) \\
&=&\text{1,1}-(1_{A}\otimes gx_{2})B(h\otimes 1_{H})(1_{A}\otimes
g)(1_{A}\otimes x_{1}x_{2})+ \\
&&\text{1,2}-(1_{A}\otimes g)B(gx_{2}h\otimes 1_{H})(1_{A}\otimes
g)(1_{A}\otimes x_{1}x_{2})+
\end{eqnarray*}%
\begin{eqnarray*}
&&\text{2}+(1_{A}\otimes gx_{2})B(h\otimes x_{2})(1_{A}\otimes gx_{1}) \\
&=&\text{1,1}(1_{A}\otimes gx_{2})B(h\otimes 1_{H})(1_{A}\otimes
x_{2})(1_{A}\otimes gx_{1})+ \\
&&\text{2,1}-(1_{A}\otimes gx_{2})(1_{A}\otimes g)B(gx_{2}h\otimes
1_{H})(1_{A}\otimes g)(1_{A}\otimes gx_{1})
\end{eqnarray*}%
\begin{eqnarray*}
&&\text{3}-(1_{A}\otimes gx_{1})B(h\otimes x_{2})(1_{A}\otimes gx_{2}) \\
&=&\text{3,1}(1_{A}\otimes gx_{1})(1_{A}\otimes gx_{2})B(h\otimes
1_{H})(1_{A}\otimes g)(1_{A}\otimes gx_{2}) \\
&&+\text{3,2}(1_{A}\otimes gx_{1})(1_{A}\otimes g)B(gx_{2}h\otimes
1_{H})(1_{A}\otimes g)(1_{A}\otimes gx_{2})
\end{eqnarray*}%
\begin{eqnarray*}
&&\text{4}+(1_{A}\otimes x_{1}x_{2})B(h\otimes x_{2}) \\
&=&\text{3,1}(1_{A}\otimes x_{1}x_{2})B(h\otimes 1_{H})(1_{A}\otimes x_{2})-%
\text{4,1}(1_{A}\otimes x_{1}x_{2})(1_{A}\otimes g)B(gx_{2}h\otimes
1_{H})(1_{A}\otimes g)
\end{eqnarray*}%
\begin{eqnarray*}
&&\text{5}-(1_{A}\otimes x_{2})B(gx_{1}h\otimes x_{2})+ \\
&=&\text{5,1}-(1_{A}\otimes x_{2})B(gx_{1}h\otimes 1_{H})(1_{A}\otimes
x_{2})+\text{5,2}(1_{A}\otimes x_{2})(1_{A}\otimes g)B(x_{1}x_{2}h\otimes
1_{H})(1_{A}\otimes g)
\end{eqnarray*}%
\begin{eqnarray*}
&&\text{6}-(1_{A}\otimes g)B(gx_{1}h\otimes x_{2})(1_{A}\otimes gx_{2})+ \\
&=&+\text{5,1}(1_{A}\otimes g)(1_{A}\otimes gx_{2})B(gx_{1}h\otimes
1_{H})(1_{A}\otimes g)(1_{A}\otimes gx_{2})+\text{6,1}B(x_{1}x_{2}h\otimes
1_{H})(1_{A}\otimes g)(1_{A}\otimes gx_{2})
\end{eqnarray*}%
\begin{eqnarray*}
&&\text{7}+(1_{A}\otimes x_{1})B(gx_{2}h\otimes x_{2})+ \\
&=&\text{3,2}(1_{A}\otimes x_{1})B(gx_{2}h\otimes 1_{H})(1_{A}\otimes x_{2})-%
\text{4,1}(1_{A}\otimes x_{1})(1_{A}\otimes gx_{2})B(gx_{2}h\otimes
1_{H})(1_{A}\otimes g)
\end{eqnarray*}%
\begin{eqnarray*}
&&\text{8}+(1_{A}\otimes g)B(gx_{2}h\otimes x_{2})(1_{A}\otimes gx_{1}) \\
&=&\text{1,2}(1_{A}\otimes g)B(gx_{2}h\otimes 1_{H})(1_{A}\otimes
x_{2})(1_{A}\otimes gx_{1}) \\
&&-\text{2,1}(1_{A}\otimes g)(1_{A}\otimes gx_{2})B(gx_{2}h\otimes
1_{H})(1_{A}\otimes g)(1_{A}\otimes gx_{1})
\end{eqnarray*}%
\begin{eqnarray*}
\text{9} &&+B(x_{1}x_{2}h\otimes x_{2}) \\
&=&\text{6,1}B(x_{1}x_{2}h\otimes 1_{H})(1_{A}\otimes x_{2})-\text{5,2}%
(1_{A}\otimes gx_{2})B(x_{1}x_{2}h\otimes 1_{H})(1_{A}\otimes g)
\end{eqnarray*}%
Therefore our claim is proved.

\part{Casimir Condition\label{PART CASIMIR}}

In this part we explore the Casimir Condition. Namely we consider our seven
elements $g\otimes 1_{H},x_{1}\otimes 1_{H},x_{1}\otimes 1_{H},x_{2}\otimes
1_{H},x_{1}x_{2}\otimes 1_{H},gx_{1}\otimes 1_{H},gx_{2}\otimes
1_{H},gx_{1}x_{2}\otimes 1_{H}$ and we will assume that Casimir condition \
holds for $B$ with respect to them and deduce the particular form that we
will get in that case.

After doing this, we will consider the other elements. Because of condition $%
\left( \ref{eq.a}\right) ,$we can consider only half of them. In fact,
Proposition \ref{Pro:gh} states that, by assuming that condition $\left( \ref%
{eq.a}\right) $ holds and that, for some $h,h^{\prime }\in H$, $B\left(
h\otimes h^{\prime }\right) ,$ satisfies Casimir condition, then also $%
B\left( gh\otimes gh^{\prime }\right) $ does.

In doing this, we will extensively use Proposition \ref{Pro helem} and all
the formulas connected to it, which allow us, using $\left( \ref{eq.h}%
\right) ,$ to write all $B\left( h\otimes h^{\prime }\right) $ by using $%
B\left( h\otimes 1_{H}\right) $ where $h\otimes 1_{H}$ runs through the
seven elements.

This will be a long checking that, at the end, will show us, that no more
conditions will appear.

\section{$B\left( g\otimes 1_{H}\right) $}

We write the Casimir condition \ for $B\left( g\otimes 1_{H}\right) $

\begin{eqnarray*}
&&\sum_{a,b_{1},b_{2},d,e_{1},e_{2}=0}^{1}\sum_{l_{1}=0}^{b_{1}}%
\sum_{l_{2}=0}^{b_{2}}\sum_{u_{1}=0}^{e_{1}}\sum_{u_{2}=0}^{e_{2}}\left(
-1\right) ^{\alpha \left( g;l_{1},l_{2},u_{1},u_{2}\right) }B(g\otimes
1_{H};G^{a}X_{1}^{b_{1}}X_{2}^{b_{2}},g^{d}x_{1}^{e_{1}}x_{2}^{e_{2}}) \\
&&G^{a}X_{1}^{b_{1}-l_{1}}X_{2}^{b_{2}-l_{2}}\otimes
g^{d}x_{1}^{e_{1}-u_{1}}x_{2}^{e_{2}-u_{2}}\otimes
g^{a+b_{1}+b_{2}+l_{1}+l_{2}+d+e_{1}+e_{2}+u_{1}+u_{2}+1}x_{1}^{l_{1}+u_{1}}x_{2}^{l_{2}+u_{2}}
\\
&=&B^{A}(g\otimes 1_{H})\otimes B^{H}(g\otimes 1_{H})\otimes 1_{H}
\end{eqnarray*}%
We concentrate on the different occurrences of the third part of the tensor.

\subsection{Case $1_{H}$}

For the left side of the equation we get%
\begin{eqnarray*}
l_{1} &=&u_{1}=0 \\
l_{2} &=&u_{2}=0 \\
a+b_{1}+b_{2}+d+e_{1}+e_{2} &\equiv &1 \\
\alpha \left( g;0,0,0,0\right) &\equiv &0
\end{eqnarray*}%
and we get the equality becomes%
\begin{eqnarray*}
&&\sum_{\substack{ a,b_{1},b_{2},d,e_{1},e_{2}=0  \\ %
a+b_{1}+b_{2}+d+e_{1}+e_{2}\equiv 1}}^{1}B(g\otimes
1_{H};G^{a}X_{1}^{b_{1}}X_{2}^{b_{2}},g^{d}x_{1}^{e_{1}}x_{2}^{e_{2}})G^{a}X_{1}^{b_{1}}X_{2}^{b_{2}}\otimes g^{d}x_{1}^{e_{1}}x_{2}^{e_{2}}
\\
&=&\sum_{a,b_{1},b_{2},d,e_{1},e_{2}=0}^{1}B(g\otimes
1_{H};G^{a}X_{1}^{b_{1}}X_{2}^{b_{2}},g^{d}x_{1}^{e_{1}}x_{2}^{e_{2}})G^{a}X_{1}^{b_{1}}X_{2}^{b_{2}}\otimes g^{d}x_{1}^{e_{1}}x_{2}^{e_{2}}
\end{eqnarray*}%
Thus we obtain%
\begin{equation}
B(g\otimes
1_{H};G^{a}X_{1}^{b_{1}}X_{2}^{b_{2}},g^{d}x_{1}^{e_{1}}x_{2}^{e_{2}})=0%
\text{ whenever }a+b_{1}+b_{2}+d+e_{1}+e_{2}\equiv 0.  \label{got1 first}
\end{equation}

\subsection{Case $g$}

For the left side of the equation we get%
\begin{eqnarray*}
l_{1} &=&u_{1}=0 \\
l_{2} &=&u_{2}=0 \\
a+b_{1}+b_{2}+d+e_{1}+e_{2} &\equiv &0 \\
\alpha \left( g;0,0,0,0\right) &\equiv &0
\end{eqnarray*}%
In view of $\left( \ref{got1 first}\right) $ we get nothing new.

\subsection{Case $x_{1}$}

For the left side of the equation we get%
\begin{eqnarray*}
l_{1}+u_{1} &=&1 \\
l_{2} &=&u_{2}=0 \\
a+b_{1}+b_{2}+d+e_{1}+e_{2} &\equiv &0
\end{eqnarray*}

In view of $\left( \ref{got1 first}\right) $ we get a trivial equality.

\subsection{Case $x_{2}$}

For the left side of the equation we get%
\begin{eqnarray*}
l_{2}+u_{2} &=&1 \\
l_{1} &=&u_{1}=0 \\
a+b_{1}+b_{2}+d+e_{1}+e_{2} &\equiv &0
\end{eqnarray*}

In view of $\left( \ref{got1 first}\right) $ we get a trivial equality.

\subsection{Case $x_{1}x_{2}$}

\begin{eqnarray*}
l_{1}+u_{1} &=&1 \\
l_{2}+u_{2} &=&1 \\
a+b_{1}+b_{2}+d+e_{1}+e_{2} &\equiv &1 \\
&&\alpha \left( g;0,0,1,1\right) \equiv 1+e_{2} \\
&&\alpha \left( g;0,1,1,0\right) \equiv e_{2}+a+b_{1}+b_{2} \\
&&\alpha \left( g;1,0,0,1\right) \equiv a+b_{1} \\
&&\alpha \left( g;1,1,0,0\right) \equiv 1+b_{2}
\end{eqnarray*}%
\begin{eqnarray*}
&&\sum_{\substack{ a,b_{1},b_{2},d,e_{1},e_{2}=0  \\ %
a+b_{1}+b_{2}+d+e_{1}+e_{2}\equiv 1}}^{1}\sum_{l_{1}=0}^{b_{1}}\sum
_{\substack{ l_{2}=0  \\ l_{1}+u_{1}=1}}^{b_{2}}\sum_{u_{1}=0}^{e_{1}}\sum
_{\substack{ u_{2}=0  \\ l_{2}+u_{2}=1}}^{e_{2}}\left( -1\right) ^{\alpha
\left( g;l_{1},l_{2},u_{1},u_{2}\right) } \\
&&B(g\otimes
1_{H};G^{a}X_{1}^{b_{1}}X_{2}^{b_{2}},g^{d}x_{1}^{e_{1}}x_{2}^{e_{2}})G^{a}X_{1}^{b_{1}-l_{1}}X_{2}^{b_{2}-l_{2}}\otimes g^{d}x_{1}^{e_{1}-u_{1}}x_{2}^{e_{2}-u_{2}}
\\
&=&\sum_{\substack{ a,b_{1},b_{2},d=0  \\ a+b_{1}+b_{2}+d\equiv 1}}%
^{1}B(g\otimes
1_{H};G^{a}X_{1}^{b_{1}}X_{2}^{b_{2}},g^{d}x_{1}x_{2})G^{a}X_{1}^{b_{1}}X_{2}^{b_{2}}\otimes g^{d}+
\\
&&\sum_{\substack{ a,b_{1},d,e_{2}=0  \\ a+b_{1}d+e_{2}\equiv 1}}^{1}\left(
-1\right) ^{e_{2}+a+b_{1}}B(g\otimes
1_{H};G^{a}X_{1}^{b_{1}}X_{2},g^{d}x_{1}x_{2}^{e_{2}})G^{a}X_{1}^{b_{1}}%
\otimes g^{d}x_{2}^{e_{2}}+ \\
&&\sum_{\substack{ a,b_{2},d,e_{1}=0  \\ a+b_{2}+d+e_{1}\equiv 1}}^{1}\left(
-1\right) ^{a+1}B(g\otimes
1_{H};G^{a}X_{1}X_{2}^{b_{2}},g^{d}x_{1}^{e_{1}}x_{2})G^{a}X_{2}^{b_{2}}%
\otimes g^{d}x_{1}^{e_{1}}+ \\
&&\sum_{\substack{ a,d,e_{1},e_{2}=0  \\ a+d+e_{1}+e_{2}\equiv 1}}^{1}\left(
-1\right) ^{1+b_{2}}B(g\otimes
1_{H};G^{a}X_{1}X_{2},g^{d}x_{1}^{e_{1}}x_{2}^{e_{2}})G^{a}\otimes
g^{d}x_{1}^{e_{1}}x_{2}^{e_{2}}+
\end{eqnarray*}%
and we get%
\begin{eqnarray*}
&&\sum_{\substack{ a,d=0  \\ \overline{a+d}=1}}^{1}\left[
\begin{array}{c}
B(g\otimes 1_{H};G^{a}X_{1}X_{2},g^{d})+\left( -1\right) ^{a+1}B(g\otimes
1_{H};G^{a}X_{2},g^{d}x_{1}) \\
+\left( -1\right) ^{a}B(g\otimes 1_{H};G^{a}X_{1},g^{d}x_{2})+B(g\otimes
1_{H};G^{a},g^{d}x_{1}x_{2})%
\end{array}%
\right] G^{a}\otimes g^{d} \\
&&\sum_{\substack{ a,d=0  \\ \overline{a+d}=0}}^{1}\left[ B(g\otimes
1_{H};G^{a}X_{1}X_{2},g^{d}x_{1})+\left( -1\right) ^{a}B(g\otimes
1_{H};G^{a}X_{1},g^{d}x_{1}x_{2})\right] G^{a}\otimes g^{d}x_{1} \\
&&\sum_{\substack{ a,d=0  \\ \overline{a+d}=0}}^{1}\left[ B(g\otimes
1_{H};G^{a}X_{1}X_{2},g^{d}x_{2})+\left( -1\right) ^{a}B(g\otimes
1_{H};G^{a}X_{2},g^{d}x_{1}x_{2})\right] G^{a}\otimes g^{d}x_{2} \\
&&\sum_{\substack{ a,d=0  \\ \overline{a+d}=1}}^{1}B(g\otimes
1_{H};G^{a}X_{1}X_{2},g^{d}x_{1}x_{2})G^{a}\otimes g^{d}x_{1}x_{2} \\
&&\sum_{\substack{ a,d=0  \\ \overline{a+d}=0}}^{1}\left[ \left( -1\right)
^{a}B(g\otimes 1_{H};G^{a}X_{1}X_{2},g^{d}x_{1})+B(g\otimes
1_{H};G^{a}X_{1},g^{d}x_{1}x_{2})\right] G^{a}X_{1}\otimes g^{d}+ \\
&&+\sum_{\substack{ a,d=0  \\ \overline{a+d}=1}}^{1}\left( -1\right)
^{a+1}B(g\otimes 1_{H};G^{a}X_{1}X_{2},g^{d}x_{1}x_{2})G^{a}X_{1}\otimes
g^{d}x_{2}+ \\
&&+\sum_{\substack{ a,d=0  \\ \overline{a+d}=0}}^{1}\left[ \left( -1\right)
^{a}B(g\otimes 1_{H};G^{a}X_{1}X_{2},g^{d}x_{2})+B(g\otimes
1_{H};G^{a}X_{2},g^{d}x_{1}x_{2})\right] G^{a}X_{2}\otimes g^{d} \\
&&\sum_{\substack{ a,d=0  \\ \overline{a+d}=1}}^{1}\left( -1\right)
^{a}B(g\otimes 1_{H};G^{a}X_{1}X_{2},g^{d}x_{1}x_{2})G^{a}X_{2}\otimes
g^{d}x_{1} \\
&&\sum_{\substack{ a,d=0  \\ \overline{a+d}=1}}^{1}B(g\otimes
1_{H};G^{a}X_{1}X_{2},g^{d}x_{1}x_{2})G^{a}X_{1}X_{2}\otimes g^{d}
\end{eqnarray*}%
And we obtain%
\begin{gather}
B(g\otimes 1_{H};G^{a}X_{1}X_{2},g^{d})+\left( -1\right) ^{a+1}B(g\otimes
1_{H};G^{a}X_{2},g^{d}x_{1})+\left( -1\right) ^{a}B(g\otimes
1_{H};G^{a}X_{1},g^{d}x_{2})  \label{got1seven} \\
+B(g\otimes 1_{H};G^{a},g^{d}x_{1}x_{2})\text{ whenever }a+d\text{ }\equiv 1
\notag
\end{gather}%
\begin{equation}
B(g\otimes 1_{H};G^{a}X_{1}X_{2},g^{d}x_{1})+\left( -1\right) ^{a}B(g\otimes
1_{H};G^{a}X_{1},g^{d}x_{1}x_{2})=0\text{ whenever }a+d\text{ }\equiv 1
\label{got1nine}
\end{equation}%
\begin{equation}
B(g\otimes 1_{H};G^{a}X_{1}X_{2},g^{d}x_{2})+\left( -1\right) ^{a}B(g\otimes
1_{H};G^{a}X_{2},g^{d}x_{1}x_{2})=0\text{ whenever }a+d\text{ }\equiv 0
\label{got1ten}
\end{equation}%
\begin{equation}
B(g\otimes 1_{H};G^{a}X_{1}X_{2},g^{d}x_{1}x_{2})=0\text{ whenever }a+d\text{
}\equiv 1\text{ .}  \label{got1 thirteen}
\end{equation}

\subsection{Case $gx_{1}$}

\begin{eqnarray*}
l_{1}+u_{1} &=&1 \\
l_{2} &=&u_{2}=0 \\
a+b_{1}+b_{2}+l_{1}+l_{2}+d+e_{1}+e_{2}+u_{1}+u_{2}+1 &\equiv
&a+b_{1}+b_{2}+d+e_{1}+e_{2}\equiv 1
\end{eqnarray*}%
\begin{gather*}
\sum_{\substack{ a,b_{1},b_{2},d,e_{1},e_{2}=0  \\ %
a+b_{1}+b_{2}+d+e_{1}+e_{2}\equiv 1}}^{1}\sum_{l_{1}=0}^{b_{1}}\sum
_{\substack{ u_{1}=0  \\ l_{1}+u_{1}=1}}^{e_{1}}\left( -1\right) ^{\alpha
\left( g;l_{1},0,u_{1},0\right) }B(g\otimes
1_{H};G^{a}X_{1}^{b_{1}}X_{2}^{b_{2}},g^{d}x_{1}^{e_{1}}x_{2}^{e_{2}}) \\
G^{a}X_{1}^{b_{1}-l_{1}}X_{2}^{b_{2}}\otimes
g^{d}x_{1}^{e_{1}-u_{1}}x_{2}^{e_{2}}=0
\end{gather*}%
and we obtain%
\begin{equation}
B(g\otimes 1_{H};G^{a}X_{1}X_{2}^{b_{2}},g^{d}x_{2}^{e_{2}})+\left(
-1\right) ^{a+e_{2}+1}B(g\otimes
1_{H};G^{a}X_{2}^{b_{2}},g^{d}x_{1}x_{2}^{e_{2}})=0\text{ for }%
a+b_{2}+d+e_{2}\equiv 0  \label{got1five}
\end{equation}%
\begin{equation}
B(g\otimes 1_{H};G^{a}X_{1}X_{2}^{b_{2}},g^{d}x_{1}x_{2}^{e_{2}})=0\text{
for }a+b_{2}+d+e_{2}\equiv 1.  \label{got1six}
\end{equation}

\subsection{Case $gx_{2}$}

\begin{eqnarray*}
l_{1} &=&u_{1}=0 \\
l_{2}+u_{2} &=&1 \\
a+b_{1}+b_{2}+l_{1}+l_{2}+d+e_{1}+e_{2}+u_{1}+u_{2}+1 &\equiv
&a+b_{1}+b_{2}+d+e_{1}+e_{2}\equiv 1
\end{eqnarray*}%
\begin{gather*}
\sum_{\substack{ a,b_{1},b_{2},d,e_{1},e_{2}=0  \\ %
a+b_{1}+b_{2}+d+e_{1}+e_{2}\equiv 1}}^{1}\sum_{l_{2}=0}^{b_{2}}\sum
_{\substack{ u_{2}=0  \\ l_{2}+u_{2}=1}}^{e_{2}}\left( -1\right) ^{\alpha
\left( g;0,l_{2},0,u_{2}\right) }B(g\otimes
1_{H};G^{a}X_{1}^{b_{1}}X_{2}^{b_{2}},g^{d}x_{1}^{e_{1}}x_{2}^{e_{2}}) \\
G^{a}X_{1}^{b_{1}}X_{2}^{b_{2}-l_{2}}\otimes
g^{d}x_{1}^{e_{1}}x_{2}^{e_{2}-u_{2}}=0
\end{gather*}%
Thus we get%
\begin{equation}
B(g\otimes 1_{H};G^{a}X_{1}^{b_{1}}X_{2},g^{d}x_{1}^{e_{1}})+\left(
-1\right) ^{\left( a+b_{1}+1\right) }B(g\otimes
1_{H};G^{a}X_{1}^{b_{1}},g^{d}x_{1}^{e_{1}}x_{2})=0\text{ for }%
a+b_{1}+d+e_{1}\equiv 0  \label{got1fivebis}
\end{equation}%
\begin{equation}
B(g\otimes 1_{H};G^{a}X_{1}^{b_{1}}X_{2},g^{d}x_{1}^{e_{1}}x_{2})=0\text{
for }a+b_{1}+d+e_{1}\equiv 1.  \label{got1sixbis}
\end{equation}

\subsection{Case $gx_{1}x_{2}$}

We get%
\begin{eqnarray*}
l_{1}+u_{1} &=&1 \\
l_{2}+u_{2} &=&1 \\
a+b_{1}+b_{2}+d+e_{1}+e_{2} &\equiv &0
\end{eqnarray*}%
\begin{gather*}
\sum_{\substack{ a,b_{1},b_{2},d,e_{1},e_{2}=0  \\ %
a+b_{1}+b_{2}+d+e_{1}+e_{2}\equiv 0}}^{1}\sum_{l_{1}=0}^{b_{1}}%
\sum_{l_{2}=0}^{b_{2}}\sum_{\substack{ u_{1}=0  \\ l_{1}+u_{1}=1}}%
^{e_{1}}\sum _{\substack{ u_{2}=0  \\ l_{2}+u_{2}=1}}^{e_{2}}\left(
-1\right) ^{\alpha \left( g;l_{1},l_{2},u_{1},u_{2}\right) } \\
B(g\otimes
1_{H};G^{a}X_{1}^{b_{1}}X_{2}^{b_{2}},g^{d}x_{1}^{e_{1}}x_{2}^{e_{2}})G^{a}X_{1}^{b_{1}-l_{1}}X_{2}^{b_{2}-l_{2}}\otimes g^{d}x_{1}^{e_{1}-u_{1}}x_{2}^{e_{2}-u_{2}}=0.
\end{gather*}%
and the formula holds in view of $\left( \ref{got1 first}\right) .$

At the end we got%
\begin{equation*}
B(g\otimes
1_{H};G^{a}X_{1}^{b_{1}}X_{2}^{b_{2}},g^{d}x_{1}^{e_{1}}x_{2}^{e_{2}})=0%
\text{ whenever }a+b_{1}+b_{2}+d+e_{1}+e_{2}\equiv 0\text{ }\left( \ref{got1
first}\right)
\end{equation*}%
and%
\begin{equation*}
B(g\otimes 1_{H};G^{a}X_{1}X_{2}^{b_{2}},g^{d}x_{1}x_{2}^{e_{2}})=0\text{
for }a+b_{2}+d+e_{2}\equiv 1\text{ }\left( \ref{got1six}\right) .
\end{equation*}

\subsection{The equalities we get}

The equality $\left( \ref{got1six}\right) $ splits in%
\begin{equation}
\text{ }B(g\otimes 1_{H};GX_{1},x_{1})=0  \label{got1six,1}
\end{equation}%
\begin{equation}
\text{ }B(g\otimes 1_{H};X_{1}X_{2},x_{1})=0  \label{got1six, 2}
\end{equation}%
\begin{equation}
\text{ }B(g\otimes 1_{H};X_{1},gx_{1})=0  \label{got1six,3}
\end{equation}%
\begin{equation}
B(g\otimes 1_{H};X_{1},x_{1}x_{2})=0  \label{got1six,4}
\end{equation}%
\begin{equation}
\text{ }B(g\otimes 1_{H};GX_{1}X_{2},gx_{1})=0  \label{got1six,5}
\end{equation}%
\begin{equation}
\text{ }B(g\otimes 1_{H};GX_{1}X_{2},x_{1}x_{2})=0  \label{got1six,6}
\end{equation}%
\begin{equation}
\text{ }B(g\otimes 1_{H};GX_{1},gx_{1}x_{2})=0  \label{got1six,7}
\end{equation}%
\begin{equation}
\text{ }B(g\otimes 1_{H};X_{1}X_{2},gx_{1}x_{2})=0  \label{got1six,8}
\end{equation}%
Moreover we obtained
\begin{equation*}
B(g\otimes 1_{H};G^{a}X_{1}X_{2}^{b_{2}},g^{d}x_{2}^{e_{2}})+\left(
-1\right) ^{a+e_{2}+1}B(g\otimes
1_{H};G^{a}X_{2}^{b_{2}},g^{d}x_{1}x_{2}^{e_{2}})=0\text{ for }%
a+b_{2}+d+e_{2}\equiv 0\text{ }\left( \ref{got1five}\right)
\end{equation*}%
which splits in%
\begin{equation}
B(g\otimes 1_{H};GX_{1}X_{2},1_{H})+B(g\otimes 1_{H};GX_{2},x_{1})=0\text{ }
\label{got1five,1}
\end{equation}%
\begin{equation}
B(g\otimes 1_{H};GX_{1},g)+B(g\otimes 1_{H};G,gx_{1})=0\text{ }
\label{got1five,2}
\end{equation}%
\begin{equation}
B(g\otimes 1_{H};GX_{1},x_{2})-B(g\otimes 1_{H};G,x_{1}x_{2})=0\text{ }
\label{got1five,3}
\end{equation}%
\begin{equation}
B(g\otimes 1_{H};X_{1}X_{2},g)-B(g\otimes 1_{H};X_{2},gx_{1})=0\text{ }
\label{got1five,4}
\end{equation}%
\begin{equation}
B(g\otimes 1_{H};X_{1}X_{2},x_{2})+B(g\otimes 1_{H};X_{2},x_{1}x_{2})=0
\label{got1five,5old}
\end{equation}%
\begin{equation}
B(g\otimes 1_{H};X_{1},1_{H})-B(g\otimes 1_{H};1_{A},x_{1})=0\text{ }
\label{got1five,6}
\end{equation}

\begin{equation}
B(g\otimes 1_{H};X_{1},gx_{2})+B(g\otimes 1_{H};1_{A},gx_{1}x_{2})=0
\label{got1five,7}
\end{equation}%
\begin{equation}
B(g\otimes 1_{H};GX_{1}X_{2},gx_{2})-B(g\otimes 1_{H};GX_{2},gx_{1}x_{2})=0%
\text{ }=0\text{ }  \label{got1five,8old}
\end{equation}%
We also got%
\begin{equation*}
B(g\otimes 1_{H};G^{a}X_{1}^{b_{1}}X_{2},g^{d}x_{1}^{e_{1}})+\left(
-1\right) ^{\left( a+b_{1}+1\right) }B(g\otimes
1_{H};G^{a}X_{1}^{b_{1}},g^{d}x_{1}^{e_{1}}x_{2})=0\text{ for }%
a+b_{1}+d+e_{1}\text{ }\equiv 0\text{ }\left( \ref{got1fivebis}\right)
\end{equation*}%
which splits into%
\begin{equation}
B(g\otimes 1_{H};GX_{1}X_{2},1_{H})-B(g\otimes 1_{H};GX_{1},x_{2})=0
\label{got1fivebis,1}
\end{equation}%
\begin{equation}
B(g\otimes 1_{H};X_{1}X_{2},g)+B(g\otimes 1_{H};X_{1},gx_{2})=0
\label{got1fivebis,2}
\end{equation}%
\begin{equation}
B(g\otimes 1_{H};X_{2},gx_{1})-B(g\otimes 1_{H};1_{A},gx_{1}x_{2})=0
\label{got1fivebis,3}
\end{equation}%
\begin{equation*}
B(g\otimes 1_{H};X_{1}X_{2},x_{1})+B(g\otimes 1_{H};X_{1},x_{1}x_{2})=0
\end{equation*}%
this follows from $B(g\otimes 1_{H};X_{1}X_{2},x_{1})=0$ \ $\left( \ref%
{got1six, 2}\right) $ and $B(g\otimes 1_{H};X_{1},x_{1}x_{2})=0$ \ $\left( %
\ref{got1six,4}\right) $
\begin{equation}
B(g\otimes 1_{H};GX_{2},g)+B(g\otimes 1_{H};G,gx_{2})=0
\label{got1fivebis,5}
\end{equation}%
\begin{equation}
B(g\otimes 1_{H};GX_{2},x_{1})+B(g\otimes 1_{H};G,x_{1}x_{2})=0\text{ }
\label{got1fivebis,6}
\end{equation}%
\begin{equation}
B(g\otimes 1_{H};X_{2},1_{H})-B(g\otimes 1_{H};1_{A},x_{2})=0
\label{got1fivebis,7}
\end{equation}%
\begin{eqnarray*}
B(g\otimes 1_{H};GX_{1}X_{2},gx_{1})-B(g\otimes 1_{H};GX_{1},gx_{1}x_{2})
&=&0\text{ } \\
\text{This follows from }B(g\otimes 1_{H};GX_{1}X_{2},gx_{1}) &=&0\text{ \ }%
\left( \ref{got1six,5}\right) \\
\text{and }B(g\otimes 1_{H};GX_{1},gx_{1}x_{2}) &=&0\text{\ }\left( \ref%
{got1six,7}\right)
\end{eqnarray*}%
Moreover the obtained equality%
\begin{equation*}
B(g\otimes 1_{H};G^{a}X_{1}^{b_{1}}X_{2},g^{d}x_{1}^{e_{1}}x_{2})=0\text{
for }a+b_{1}+d+e_{1}\equiv 1\text{ }\left( \ref{got1sixbis}\right)
\end{equation*}%
splits into%
\begin{equation}
B(g\otimes 1_{H};GX_{2},x_{2})=0\text{ }  \label{got1sixbis,1}
\end{equation}%
\begin{equation*}
B(g\otimes 1_{H};X_{1}X_{2},x_{2})=0\text{ this\ is }\left( \ref{got1five,5}%
\right) \text{ .}
\end{equation*}%
\begin{equation}
B(g\otimes 1_{H};X_{2},gx_{2})=0\text{ }  \label{got1sixbis,3}
\end{equation}%
\begin{equation}
B(g\otimes 1_{H};X_{2},x_{1}x_{2})=0\text{ }  \label{got1sixbis,4}
\end{equation}%
\begin{eqnarray*}
B(g\otimes 1_{H};X_{1}X_{2},x_{2})+B(g\otimes 1_{H};X_{2},x_{1}x_{2}) &=&0%
\text{ \ }\left( \ref{got1five,5old}\right) \text{ } \\
&&\text{in view of }\left( \ref{got1sixbis,4}\right) \text{ becomes}
\end{eqnarray*}%
\begin{equation}
B(g\otimes 1_{H};X_{1}X_{2},x_{2})=0  \label{got1five,5}
\end{equation}%
\begin{equation*}
B(g\otimes 1_{H};X_{1}X_{2},gx_{1}x_{2})=0\text{ this is }\ \left( \ref%
{got1six,8}\right) \text{ }
\end{equation*}%
\begin{equation*}
B(g\otimes 1_{H};GX_{2},gx_{1}x_{2})=0\text{ \ this follows from }\left( \ref%
{got1five,8}\right)
\end{equation*}%
\begin{equation*}
B(g\otimes 1_{H};GX_{1}X_{2},x_{1}x_{2})=0\text{ this is }\left( \ref%
{got1six,6}\right) \text{ }
\end{equation*}%
\begin{equation}
B(g\otimes 1_{H};GX_{1}X_{2},gx_{2})=0\text{ }  \label{got1sixbis,8}
\end{equation}%
\begin{eqnarray*}
B(g\otimes 1_{H};GX_{1}X_{2},gx_{2})-B(g\otimes 1_{H};GX_{2},gx_{1}x_{2})
&=&0\text{ \ }\left( \ref{got1five,8old}\right) \text{ } \\
\text{in view of }B(g\otimes 1_{H};GX_{1}X_{2},gx_{2}) &=&0\text{ \ }\left( %
\ref{got1sixbis,8}\right) \\
&&\text{implies}
\end{eqnarray*}%
\begin{equation}
B(g\otimes 1_{H};GX_{1}X_{2},gx_{2})=0  \label{got1five,8}
\end{equation}%
The equality%
\begin{gather*}
B(g\otimes 1_{H};G^{a}X_{1}X_{2},g^{d})+\left( -1\right) ^{a+1}B(g\otimes
1_{H};G^{a}X_{2},g^{d}x_{1})+\left( -1\right) ^{a}B(g\otimes
1_{H};G^{a}X_{1},g^{d}x_{2}) \\
+B(g\otimes 1_{H};G^{a},g^{d}x_{1}x_{2})=0\text{ whenever }\overline{a+d}=1%
\text{ }\left( \ref{got1seven}\right)
\end{gather*}

splits into%
\begin{gather*}
B(g\otimes 1_{H};GX_{1}X_{2},1_{H})+B(g\otimes
1_{H};GX_{2},x_{1})-B(g\otimes 1_{H};GX_{1},x_{2}) \\
+B(g\otimes 1_{H};G,x_{1}x_{2})=0 \\
\text{In view of }\left( \ref{got1five,3}\right) B(g\otimes
1_{H};GX_{1},x_{2})-B(g\otimes 1_{H};G,x_{1}x_{2})=0\text{ this reduces to}
\\
B(g\otimes 1_{H};GX_{1}X_{2},1_{H})+B(g\otimes 1_{H};GX_{2},x_{1})=0 \\
\text{which is }\left( \ref{got1five,1}\right)
\end{gather*}%
and%
\begin{gather*}
+B(g\otimes 1_{H};X_{1}X_{2},g)-B(g\otimes 1_{H};X_{2},gx_{1})+B(g\otimes
1_{H};X_{1},gx_{2}) \\
+B(g\otimes 1_{H};1_{A},gx_{1}x_{2})=0\text{ } \\
\text{This follows in view of }\left( \ref{got1five,4}\right) \text{ }%
B(g\otimes 1_{H};X_{1}X_{2},g)-B(g\otimes 1_{H};X_{2},gx_{1})=0\text{ } \\
\text{and }\left( \ref{got1five,7}\right) B(g\otimes
1_{H};X_{1},gx_{2})+B(g\otimes 1_{H};1_{A},gx_{1}x_{2})=0
\end{gather*}%
The equality%
\begin{equation*}
B(g\otimes 1_{H};G^{a}X_{1}X_{2},g^{d}x_{1})+\left( -1\right) ^{a}B(g\otimes
1_{H};G^{a}X_{1},g^{d}x_{1}x_{2})=0\text{ whenever }a+d\equiv 0\left( \ref%
{got1nine}\right)
\end{equation*}%
splits in%
\begin{equation*}
B(g\otimes 1_{H};X_{1}X_{2},x_{1})+B(g\otimes 1_{H};X_{1},x_{1}x_{2})=0\text{
}\left( \ref{got1six, 2}\right) \left( \ref{got1six,4}\right)
\end{equation*}%
and%
\begin{equation*}
B(g\otimes 1_{H};GX_{1}X_{2},gx_{1})-B(g\otimes 1_{H};GX_{1},gx_{1}x_{2})=0%
\text{ }\left( \ref{got1six,5}\right) \text{ }\left( \ref{got1six,7}\right)
\end{equation*}%
The equality%
\begin{equation*}
B(g\otimes 1_{H};G^{a}X_{1}X_{2},g^{d}x_{2})+\left( -1\right) ^{a}B(g\otimes
1_{H};G^{a}X_{2},g^{d}x_{1}x_{2})=0\text{ whenever }a+d\text{ }\equiv
0\left( \ref{got1ten}\right)
\end{equation*}%
splits in
\begin{eqnarray*}
B(g\otimes 1_{H};X_{1}X_{2},x_{2})+B(g\otimes 1_{H};X_{2},x_{1}x_{2}) &=&0%
\text{ } \\
\text{this is \ }\left( \ref{got1five,5}\right) \text{ }B(g\otimes
1_{H};X_{1}X_{2},x_{2})+B(g\otimes 1_{H};X_{2},x_{1}x_{2}) &=&0
\end{eqnarray*}%
and
\begin{equation*}
B(g\otimes 1_{H};GX_{1}X_{2},gx_{2})-B(g\otimes 1_{H};GX_{2},gx_{1}x_{2})=0%
\text{ This is \ }\left( \ref{got1five,8}\right) .
\end{equation*}

Finally the equality%
\begin{equation*}
B(g\otimes 1_{H};G^{a}X_{1}X_{2},g^{d}x_{1}x_{2})=0\text{ whenever }%
a+d\equiv 1\text{ }\left( \ref{got1 thirteen}\right)
\end{equation*}%
splits into%
\begin{equation*}
B(g\otimes 1_{H};X_{1}X_{2},gx_{1}x_{2})=0\text{ }
\end{equation*}%
,which is $\left( \ref{got1six,8}\right) ,$ and
\begin{equation*}
B(g\otimes 1_{H};GX_{1}X_{2},x_{1}x_{2})=0\text{ }
\end{equation*}%
this is \ $\left( \ref{got1six,6}\right) $.

\subsection{The final form of the element $B\left( g\otimes 1_{H}\right) $}

In view of the above we get that%
\begin{eqnarray}
B\left( g\otimes 1_{H}\right) &=&B\left( g\otimes 1_{H};1_{A},g\right)
1_{A}\otimes g  \label{got1} \\
&&+B\left( g\otimes 1_{H};1_{A},x_{1}\right) 1_{A}\otimes x_{1}  \notag \\
&&+B\left( g\otimes 1_{H};1_{A},x_{2}\right) 1_{A}\otimes x_{2}  \notag \\
&&+B\left( g\otimes 1_{H};1_{A},gx_{1}x_{2}\right) 1_{A}\otimes gx_{1}x_{2}
\notag \\
&&+B\left( g\otimes 1_{H};G,1_{H}\right) G\otimes 1_{H}  \notag \\
&&+B\left( g\otimes 1_{H};G,x_{1}x_{2}\right) G\otimes x_{1}x_{2}  \notag \\
&&+B\left( g\otimes 1_{H};G,gx_{1}\right) G\otimes gx_{1}  \notag \\
&&+B\left( g\otimes 1_{H};G,gx_{2}\right) G\otimes gx_{2}  \notag \\
&&+B\left( g\otimes 1_{H};1_{A},x_{1}\right) X_{1}\otimes 1_{H}  \notag \\
&&-B\left( g\otimes 1_{H};1_{A},gx_{1}x_{2}\right) X_{1}\otimes gx_{2}
\notag \\
&&+B\left( g\otimes 1_{H};1_{A},x_{2}\right) X_{2}\otimes 1_{H}  \notag \\
&&+B\left( g\otimes 1_{H};1_{A},gx_{1}x_{2}\right) X_{2}\otimes gx_{1}
\notag \\
&&+B\left( g\otimes 1_{H};1_{A},gx_{1}x_{2}\right) X_{1}X_{2}\otimes g
\notag \\
&&-B\left( g\otimes 1_{H};G,gx_{1}\right) GX_{1}\otimes g  \notag \\
&&+B\left( g\otimes 1_{H};G,x_{1}x_{2}\right) GX_{1}\otimes x_{2}  \notag \\
&&-B\left( g\otimes 1_{H};G,gx_{2}\right) GX_{2}\otimes g  \notag \\
&&-B\left( g\otimes 1_{H};G,x_{1}x_{2}\right) GX_{2}\otimes x_{1}  \notag \\
&&+B\left( g\otimes 1_{H};G,x_{1}x_{2}\right) GX_{1}X_{2}\otimes 1_{H}
\notag
\end{eqnarray}

\section{$B\left( x_{1}\otimes 1_{H}\right) $}

Here we will assume that Casimir condition \ for $B\left( x_{1}\otimes
1_{H}\right) $ holds and we will deduce the particular form that $B\left(
x_{1}\otimes 1_{H}\right) $will get in this case.

We write the Casimir condition $\left( \ref{Casimir Condition}\right) $ for $%
B\left( x_{1}\otimes 1_{H}\right) .$%
\begin{eqnarray*}
&&\sum_{w_{1}=0}^{1}\sum_{a,b_{1},b_{2},d,e_{1},e_{2}=0}^{1}%
\sum_{l_{1}=0}^{b_{1}}\sum_{l_{2}=0}^{b_{2}}\sum_{u_{1}=0}^{e_{1}}%
\sum_{u_{2}=0}^{e_{2}}\left( -1\right) ^{\alpha \left(
x_{1}^{1-w_{1}};l_{1},l_{2},u_{1},u_{2}\right) } \\
&&B(g^{1+w_{1}}x_{1}^{w_{1}}\otimes
1_{H};G^{a}X_{1}^{b_{1}}X_{2}^{b_{2}},g^{d}x_{1}^{e_{1}}x_{2}^{e_{2}}) \\
&&G^{a}X_{1}^{b_{1}-l_{1}}X_{2}^{b_{2}-l_{2}}\otimes
g^{d}x_{1}^{e_{1}-u_{1}}x_{2}^{e_{2}-u_{2}}\otimes
g^{a+b_{1}+b_{2}+l_{1}+l_{2}+d+e_{11}+e_{12}+u_{1}+u_{2}}x_{1}^{l_{1}+u_{1}+1-w_{1}}x_{2}^{l_{2}+u_{2}}
\\
&=&B^{A}(x_{1}\otimes 1_{H})\otimes B^{H}(x_{1}\otimes 1_{H})\otimes 1_{H}
\end{eqnarray*}%
and we get%
\begin{eqnarray*}
&&\sum_{a,b_{1},b_{2},d,e_{1},e_{2}=0}^{1}\sum_{l_{1}=0}^{b_{1}}%
\sum_{l_{2}=0}^{b_{2}}\sum_{u_{1}=0}^{e_{1}}\sum_{u_{2}=0}^{e_{2}}\left(
-1\right) ^{\alpha \left( x_{1};l_{1},l_{2},u_{1},u_{2}\right) } \\
&&B(g\otimes
1_{H};G^{a}X_{1}^{b_{1}}X_{2}^{b_{2}},g^{d}x_{1}^{e_{1}}x_{2}^{e_{2}}) \\
&&G^{a}X_{1}^{b_{1}-l_{1}}X_{2}^{b_{2}-l_{2}}\otimes
g^{d}x_{1}^{e_{1}-u_{1}}x_{2}^{e_{2}-u_{2}}\otimes
g^{a+b_{1}+b_{2}+l_{1}+l_{2}+d+e_{1}+e_{2}+u_{1}+u_{2}}x_{1}^{l_{1}+u_{1}+1}x_{2}^{l_{2}+u_{2}}+
\\
&&+\sum_{a,b_{1},b_{2},d,e_{1},e_{2}=0}^{1}\sum_{l_{1}=0}^{b_{1}}%
\sum_{l_{2}=0}^{b_{2}}\sum_{u_{1}=0}^{e_{11}}\sum_{u_{2}=0}^{e_{12}}\left(
-1\right) \\
&&B(x_{1}\otimes
1_{H};G^{a}X_{1}^{b_{1}}X_{2}^{b_{2}},g^{d}x_{1}^{e_{1}}x_{2}^{e_{2}}) \\
&&G^{a}X_{1}^{b_{1}-l_{1}}X_{2}^{b_{2}-l_{2}}\otimes
g^{d}x_{1}^{e_{1}-u_{1}}x_{2}^{e_{2}-u_{2}}\otimes
g^{a+b_{1}+b_{2}+l_{1}+l_{2}+d+e_{11}+e_{12}+u_{1}+u_{2}}x_{1}^{l_{1}+u_{1}}x_{2}^{l_{2}+u_{2}}+
\\
&=&B^{A}(x_{1}\otimes 1_{H})\otimes B^{H}(x_{1}\otimes 1_{H})\otimes 1_{H}
\end{eqnarray*}%
We concentrate on the different occurrences of the third part of the tensor.

\subsection{Case $1_{H}$}

From the first summand of the left side we can not obtain $1_{H}$

Using the same argument in case $B\left( g\otimes 1_{H}\right) ,$ we obtain

\begin{equation}
B(x_{1}\otimes
1_{H};G^{a}X_{1}^{b_{1}}X_{2}^{b_{2}},g^{d}x_{1}^{e_{1}}x_{2}^{e_{2}})=0%
\text{ whenever }a+b_{1}+b_{2}+d+e_{1}+e_{2}\equiv 1  \label{x1ot1 first}
\end{equation}

We concentrate on the different occurrences of the third part of the tensor.

\subsection{Case $g$}

From the first summand of the left side we can not obtain $g$ . From the
second summand we just get $\left( \ref{x1ot1 first}\right) .$

\subsection{Case $x_{1}$}

Here we will have two summands with $x_{1}$ in the third position.

From the first summand of the left side we get $\left( \ref{got1 first}%
\right) .$From the second summand of the left side we get $\left( \ref{x1ot1
first}\right) .$

\subsection{Case $x_{2}$}

From the first summand of the left side we can not obtain $x_{2}$. From the
second summand of the left side we get $\left( \ref{x1ot1 first}\right) .$

\subsection{Case $gx_{1}$}

From the left side first summand we get%
\begin{eqnarray*}
a+b_{1}+b_{2}+d+e_{1}+e_{2} &\equiv &1 \\
l_{1} &=&u_{1}\equiv 0 \\
l_{2} &=&u_{2}\equiv 0.
\end{eqnarray*}%
From the left side second summand we get%
\begin{eqnarray*}
a+b_{1}+b_{2}+d+e_{1}+e_{2} &\equiv &0 \\
l_{1}+u_{1} &\equiv &1 \\
l_{2} &=&u_{2}\equiv 0
\end{eqnarray*}%
and, since
\begin{eqnarray*}
\alpha \left( x_{1};0,0,0,0\right) &=&a+b_{1}+b_{2} \\
\alpha \left( 1_{H};0,0,1,0\right) &=&e_{2}+a+b_{1}+b_{2} \\
\alpha \left( 1_{H};1,0,0,0\right) &=&b_{2}
\end{eqnarray*}%
we get%
\begin{eqnarray*}
&&\sum_{\substack{ a,b_{1},b_{2},d,e_{1},e_{2}=0  \\ %
a+b_{1}+b_{2}+d+e_{1}+e_{2}\equiv 1}}^{1}\left( -1\right)
^{a+b_{1}+b_{2}}B(g\otimes
1_{H};G^{a}X_{1}^{b_{1}}X_{2}^{b_{2}},g^{d}x_{1}^{e_{1}}x_{2}^{e_{2}})G^{a}X_{1}^{b_{1}}X_{2}^{b_{2}}\otimes g^{d}x_{1}^{e_{1}}x_{2}^{e_{2}}
\\
&&+\sum_{\substack{ a,b_{1},b_{2},d,e_{2}=0  \\ a+b_{1}+b_{2}+d+e_{2}\equiv
1 }}^{1}\left( -1\right) ^{e_{2}+a+b_{1}+b_{2}}B(x_{1}\otimes
1_{H};G^{a}X_{1}^{b_{1}}X_{2}^{b_{2}},g^{d}x_{1}x_{2}^{e_{2}})G^{a}X_{1}^{b_{1}}X_{2}^{b_{2}}\otimes g^{d}x_{2}^{e_{2}}
\\
&&+\sum_{\substack{ a,b_{1},b_{2},d,e_{1},e_{2}=0  \\ a+b_{2}+d+e_{1}+e_{2}%
\equiv 1}}^{1}\left( -1\right) ^{b_{2}}B(x_{1}\otimes
1_{H};G^{a}X_{1}X_{2}^{b_{2}},g^{d}x_{1}^{e_{1}}x_{2}^{e_{2}})G^{a}X_{2}^{b_{2}}\otimes g^{d}x_{1}^{e_{1}}x_{2}^{e_{2}}=0
\end{eqnarray*}

We now consider the following cases.

\subsubsection{$G^{a}\otimes g^{d}$}

\begin{gather*}
\sum_{\substack{ a,d=0  \\ a+d\equiv 1}}^{1}\left[ \left( -1\right)
^{a}B(g\otimes 1_{H};G^{a},g^{d})+\left( -1\right) ^{a}B(x_{1}\otimes
1_{H};G^{a},g^{d}x_{1})+B(x_{1}\otimes 1_{H};G^{a}X_{1},g^{d})\right] \\
G^{a}\otimes g^{d}=0
\end{gather*}%
\begin{gather}
B(g\otimes 1_{H};1_{A},g)+B(x_{1}\otimes 1_{H};1_{A},gx_{1})+B(x_{1}\otimes
1_{H};X_{1},g)=0  \label{x1ot1, second} \\
-B(g\otimes 1_{H};G,1_{H})-B(x_{1}\otimes 1_{H};G,x_{1})+B(x_{1}\otimes
1_{H};GX_{1},1_{H})=0  \label{x1ot1, third}
\end{gather}

\subsubsection{$G^{a}\otimes g^{d}x_{2}$}

\begin{gather*}
\sum_{\substack{ a,d=0  \\ a+d\equiv 0}}^{1}\left[ \left( -1\right)
^{a}B(g\otimes 1_{H};G^{a},g^{d}x_{2})+\left( -1\right) ^{a+1}B(x_{1}\otimes
1_{H};G^{a},g^{d}x_{1}x_{2})+B(x_{1}\otimes 1_{H};G^{a}X_{1},g^{d}x_{2})%
\right] \\
G^{a}\otimes g^{d}x_{2}=0
\end{gather*}%
and we get%
\begin{eqnarray}
B(g\otimes 1_{H};1_{A},x_{2})-B(x_{1}\otimes
1_{H};1_{A},x_{1}x_{2})+B(x_{1}\otimes 1_{H};X_{1},x_{2}) &=&0
\label{x1ot1, fourth} \\
-B(g\otimes 1_{H};G,gx_{2})+B(x_{1}\otimes
1_{H};G,gx_{1}x_{2})+B(x_{1}\otimes 1_{H};GX_{1},gx_{2}) &=&0
\label{x1ot1, fifth}
\end{eqnarray}

\subsubsection{$G^{a}\otimes g^{d}x_{1}$}

\begin{equation*}
\sum_{\substack{ a,d=0  \\ a+d\equiv 0}}^{1}\left[ \left( -1\right)
^{a}B(g\otimes 1_{H};G^{a},g^{d}x_{1})+B(x_{1}\otimes
1_{H};G^{a}X_{1},g^{d}x_{1})\right] G^{a}\otimes g^{d}x_{1}=0
\end{equation*}%
and we get%
\begin{eqnarray}
B(g\otimes 1_{H};1_{A},x_{1})+B(x_{1}\otimes 1_{H};X_{1},x_{1}) &=&0
\label{x1ot1, six} \\
-B(g\otimes 1_{H};G,gx_{1})+B(x_{1}\otimes 1_{H};GX_{1},gx_{1}) &=&0
\label{x1ot1, seven}
\end{eqnarray}

\subsubsection{$G^{a}X_{2}\otimes g^{d}$}

\begin{gather*}
\sum_{\substack{ a,d=0  \\ a+d\equiv 0}}^{1}\left[
\begin{array}{c}
\left( -1\right) ^{a+1}B(g\otimes 1_{H};G^{a}X_{2},g^{d})+ \\
\left( -1\right) ^{a+1}B(x_{1}\otimes
1_{H};G^{a}X_{2},g^{d}x_{1})-B(x_{1}\otimes 1_{H};G^{a}X_{1}X_{2},g^{d})%
\end{array}%
\right] \\
G^{a}X_{2}\otimes g^{d}=0
\end{gather*}%
and we get%
\begin{gather}
-B(g\otimes 1_{H};X_{2},1_{H})-B(x_{1}\otimes
1_{H};X_{2},x_{1})-B(x_{1}\otimes 1_{H};X_{1}X_{2},1_{H})=0
\label{x1ot1, eight} \\
B(g\otimes 1_{H};GX_{2},g)+B(x_{1}\otimes
1_{H};GX_{2},gx_{1})-B(x_{1}\otimes 1_{H};GX_{1}X_{2},g)=0
\label{x1ot1, nine}
\end{gather}

\subsubsection{$G^{a}X_{1}\otimes g^{d}$}

\begin{equation*}
\sum_{\substack{ a,d=0  \\ a+d\equiv 0}}^{1}\left[ \left( -1\right)
^{a+1}B(g\otimes 1_{H};G^{a}X_{1},g^{d})+\left( -1\right)
^{a+1}B(x_{1}\otimes 1_{H};G^{a}X_{1},g^{d}x_{1})\right] G^{a}X_{1}\otimes
g^{d}=0
\end{equation*}%
and we get%
\begin{eqnarray}
-B(g\otimes 1_{H};X_{1},1_{H})-B(x_{1}\otimes 1_{H};X_{1},x_{1}) &=&0
\label{x1ot1, ten} \\
B(g\otimes 1_{H};GX_{1},g)+B(x_{1}\otimes 1_{H};GX_{1},gx_{1}) &=&0
\label{x1ot1, eleven}
\end{eqnarray}

\subsubsection{$G^{a}\otimes g^{d}x_{1}x_{2}$}

\begin{equation*}
\sum_{\substack{ a,d=0  \\ a+d\equiv 1}}^{1}\left[ \left( -1\right)
^{a}B(g\otimes 1_{H};G^{a},g^{d}x_{1}x_{2})+B(x_{1}\otimes
1_{H};G^{a}X_{1},g^{d}x_{1}x_{2})\right] G^{a}\otimes g^{d}x_{1}x_{2}=0
\end{equation*}%
and we get%
\begin{eqnarray}
B(g\otimes 1_{H};1_{A},gx_{1}x_{2})+B(x_{1}\otimes 1_{H};X_{1},gx_{1}x_{2})
&=&0  \label{x1ot1, twelve} \\
-B(g\otimes 1_{H};G,x_{1}x_{2})+B(x_{1}\otimes 1_{H};GX_{1},x_{1}x_{2}) &=&0
\label{x1ot1, thirteen}
\end{eqnarray}

\subsubsection{$G^{a}X_{2}\otimes g^{d}x_{2}$}

\begin{equation*}
\sum_{\substack{ a,d=0  \\ a+d\equiv 1}}^{1}\left[
\begin{array}{c}
\left( -1\right) ^{a+1}B(g\otimes 1_{H};G^{a}X_{2},g^{d}x_{2})+\left(
-1\right) ^{a}B(x_{1}\otimes 1_{H};G^{a}X_{2},g^{d}x_{1}x_{2}) \\
-B(x_{1}\otimes 1_{H};G^{a}X_{1}X_{2},g_{1}^{d}x_{2})%
\end{array}%
\right] G^{a}X_{2}\otimes g^{d}x_{2}=0
\end{equation*}%
and we get%
\begin{eqnarray}
-B(g\otimes 1_{H};X_{2},gx_{2})+B(x_{1}\otimes
1_{H};X_{2},gx_{1}x_{2})-B(x_{1}\otimes 1_{H};X_{1}X_{2},gx_{2}) &=&0
\label{x1ot1, fourteen} \\
B(g\otimes 1_{H};GX_{2},x_{2})-B(x_{1}\otimes
1_{H};GX_{2},x_{1}x_{2})-B(x_{1}\otimes 1_{H};GX_{1}X_{2},x_{2}) &=&0
\label{x1ot1, fifteen}
\end{eqnarray}

\subsubsection{$G^{a}X_{1}\otimes g^{d}x_{2}$}

\begin{equation*}
\sum_{\substack{ a,d=0  \\ a+d\equiv 1}}^{1}\left[ \left( -1\right)
^{a+1}B(g\otimes 1_{H};G^{a}X_{1},g^{d}x_{2})+\left( -1\right)
^{a}B(x_{1}\otimes 1_{H};G^{a}X_{1},g^{d}x_{1}x_{2})\right]
G^{a}X_{1}\otimes g^{d}x_{2}=0
\end{equation*}%
and we get%
\begin{eqnarray}
-B(g\otimes 1_{H};X_{1},gx_{2})+B(x_{1}\otimes 1_{H};X_{1},gx_{1}x_{2}) &=&0
\label{x1ot1, sixteen} \\
B(g\otimes 1_{H};GX_{1},x_{2})-B(x_{1}\otimes 1_{H};GX_{1},x_{1}x_{2}) &=&0
\label{x1ot1, seventeen}
\end{eqnarray}

\subsubsection{$G^{a}X_{2}\otimes g^{d}x_{1}$}

\begin{equation*}
\sum_{\substack{ a,d=0  \\ a+d\equiv 1}}^{1}\left[ \left( -1\right)
^{a+1}B(g\otimes 1_{H};G^{a}X_{2},g^{d}x_{1})-B(x_{1}\otimes
1_{H};G^{a}X_{1}X_{2},g^{d}x_{1})\right] G^{a}X_{2}\otimes g^{d}x_{1}=0
\end{equation*}%
and we get%
\begin{eqnarray}
-B(g\otimes 1_{H};X_{2},gx_{1})-B(x_{1}\otimes 1_{H};X_{1}X_{2},gx_{1}) &=&0
\label{x1ot1, eighteen} \\
B(g\otimes 1_{H};GX_{2},x_{1})-B(x_{1}\otimes 1_{H};GX_{1}X_{2},x_{1}) &=&0
\label{x1ot1, nineteen}
\end{eqnarray}

\subsubsection{$G^{a}X_{1}\otimes g^{d}x_{1}$}

\begin{equation*}
\sum_{\substack{ a,d=0  \\ a+d\equiv 1}}^{1}\left( -1\right)
^{a+1}B(g\otimes 1_{H};G^{a}X_{1},g^{d}x_{1})G^{a}X_{1}\otimes g^{d}x_{1}=0
\end{equation*}%
and we get%
\begin{eqnarray*}
-B(g\otimes 1_{H};X_{1},gx_{1}) &=&0 \\
B(g\otimes 1_{H};GX_{1},x_{1}) &=&0
\end{eqnarray*}%
which are already known.

\subsubsection{$G^{a}X_{1}X_{2}\otimes g^{d}$}

\begin{equation*}
\sum_{\substack{ a,d=0  \\ a+d\equiv 1}}^{1}\left[ \left( -1\right)
^{a}B(g\otimes 1_{H};G^{a}X_{1}X_{2},g^{d})+\left( -1\right) ^{a}B(g\otimes
1_{H};G^{a}X_{1}X_{2},g^{d}x_{1})\right] G^{a}X_{1}X_{2}\otimes g^{d}=0
\end{equation*}%
and we get%
\begin{eqnarray}
B(g\otimes 1_{H};X_{1}X_{2},g)+B(x_{1}\otimes 1_{H};X_{1}X_{2},gx_{1}) &=&0
\label{x1ot1, twenty} \\
-B(g\otimes 1_{H};GX_{1}X_{2},1_{H})-B(x_{1}\otimes 1_{H};GX_{1}X_{2},x_{1})
&=&0  \label{x1ot1, twentyone}
\end{eqnarray}

\subsubsection{$G^{a}X_{2}\otimes g^{d}x_{1}x_{2}$}

\begin{gather*}
\sum_{\substack{ a,d=0  \\ a+d\equiv 0}}^{1}\left[ \left( -1\right)
^{a+1}B(g\otimes 1_{H};G^{a}X_{2},g^{d}x_{1}x_{2})-B(x_{1}\otimes
1_{H};G^{a}X_{1}X_{2},g^{d}x_{1}x_{2})\right] \\
G^{a}X_{1}^{b_{1}}X_{2}^{b_{2}}\otimes g^{d}x_{1}^{e_{1}}x_{2}^{e_{2}}=0
\end{gather*}%
and we get%
\begin{eqnarray}
-B(g\otimes 1_{H};X_{2},x_{1}x_{2})-B(x_{1}\otimes
1_{H};X_{1}X_{2},x_{1}x_{2}) &=&0  \label{x1ot1, twentytwo} \\
B(g\otimes 1_{H};GX_{2},gx_{1}x_{2})-B(x_{1}\otimes
1_{H};GX_{1}X_{2},gx_{1}x_{2}) &=&0  \label{x1ot1, twentythree}
\end{eqnarray}

\subsubsection{$G^{a}X_{1}\otimes g^{d}x_{1}x_{2}$}

\begin{equation*}
\sum_{\substack{ a,d=0  \\ a+d\equiv 0}}^{1}\left( -1\right)
^{a+1}B(g\otimes 1_{H};G^{a}X_{1},g^{d}x_{1}x_{2})G^{a}X_{1}\otimes
g^{d}x_{1}x_{2}=0
\end{equation*}%
and we get%
\begin{eqnarray*}
-B(g\otimes 1_{H};X_{1},x_{1}x_{2}) &=&0 \\
B(g\otimes 1_{H};GX_{1},gx_{1}x_{2}) &=&0
\end{eqnarray*}%
which we already got.

\subsubsection{$G^{a}X_{1}X_{2}\otimes g^{d}x_{2}$}

\begin{gather*}
\sum_{\substack{ a,d=0  \\ a+d\equiv 0}}^{1}\left[ \left( -1\right)
^{a}B(g\otimes 1_{H};G^{a}X_{1}X_{2},g^{d}x_{2})-B(g\otimes
1_{H};G^{a}X_{1}X_{2},g^{d}x_{1}x_{2})\right] \\
G^{a}X_{1}X_{2}\otimes g^{d}x_{2}=0
\end{gather*}%
and we get%
\begin{eqnarray}
B(g\otimes 1_{H};X_{1}X_{2},x_{2})-B(x_{1}\otimes
1_{H};X_{1}X_{2},x_{1}x_{2}) &=&0  \label{x1ot1, twentyfour} \\
-B(g\otimes 1_{H};GX_{1}X_{2},gx_{2})+B(x_{1}\otimes
1_{H};GX_{1}X_{2},gx_{1}x_{2}) &=&0  \label{x1ot1, twentyfive}
\end{eqnarray}

\subsubsection{$G^{a}X_{1}X_{2}\otimes g^{d}x_{1}$}

\begin{equation*}
\sum_{\substack{ a,d=0  \\ a+d\equiv 0}}^{1}\left( -1\right) ^{a}B(g\otimes
1_{H};G^{a}X_{1}X_{2},g^{d}x_{1})G^{a}X_{1}X_{2}\otimes g^{d}x_{1}=0
\end{equation*}%
and we get%
\begin{eqnarray*}
B(g\otimes 1_{H};X_{1}X_{2},x_{1}) &=&0 \\
-B(g\otimes 1_{H};GX_{1}X_{2},gx_{1}) &=&0
\end{eqnarray*}%
which are already known.

\subsubsection{$G^{a}X_{1}X_{2}\otimes g^{d}x_{1}x_{2}$}

\begin{equation*}
\sum_{\substack{ a,d=0  \\ a+d\equiv 1}}^{1}\left( -1\right) ^{a}B(g\otimes
1_{H};G^{a}X_{1}X_{2},g^{d}x_{1}x_{2})G^{a}X_{1}X_{2}\otimes
g^{d}x_{1}x_{2}=0
\end{equation*}%
and we get%
\begin{eqnarray*}
B(g\otimes 1_{H};X_{1}X_{2},gx_{1}x_{2}) &=&0 \\
B(g\otimes 1_{H};GX_{1}X_{2},x_{1}x_{2}) &=&0
\end{eqnarray*}%
which are already known.

\subsection{Case $gx_{2}$}

The first summand of the left side does not appear in the Casimir condition.
We consider the second summand

\begin{eqnarray*}
&&\sum_{a,b_{1},b_{2},d,e_{1},e_{2}=0}^{1}\sum_{l_{1}=0}^{b_{1}}%
\sum_{l_{2}=0}^{b_{2}}\sum_{u_{1}=0}^{e_{1}}\sum_{u_{2}=0}^{e_{2}}\left(
-1\right) ^{\alpha \left( 1_{H};l_{1},l_{2},u_{1},u_{2}\right) } \\
&&B(x_{1}\otimes
1_{H};G^{a}X_{1}^{b_{1}}X_{2}^{b_{2}},g^{d}x_{1}^{e_{1}}x_{2}^{e_{2}}) \\
&&G^{a}X_{1}^{b_{1}-l_{1}}X_{2}^{b_{2}-l_{2}}\otimes
g^{d}x_{1}^{e_{1}-u_{1}}x_{2}^{e_{2}-u_{2}}\otimes
g^{a+b_{1}+b_{2}+l_{1}+l_{2}+d+e_{1}+e_{2}+u_{1}+u_{2}}x_{1}^{l_{1}+u_{1}}x_{2}^{l_{2}+u_{2}}
\end{eqnarray*}

\begin{eqnarray*}
l_{1} &=&u_{1}=0 \\
l_{2}+u_{2} &=&1 \\
a+b_{1}+b_{2}+d+e_{1}+e_{2} &\equiv &0
\end{eqnarray*}%
so that we get%
\begin{gather*}
\sum_{\substack{ a,b_{1},b_{2},d,e_{1},e_{2}=0  \\ %
a+b_{1}+b_{2}+d+e_{1}+e_{2}\equiv 0}}^{1}\sum_{l_{2}=0}^{b_{2}}\sum
_{\substack{ u_{2}=0  \\ l_{2}+u_{2}=1}}^{e_{2}}\left( -1\right) ^{\alpha
\left( 1_{H};0,l_{2},0,u_{2}\right) }B(x_{1}\otimes
1_{H};G^{a}X_{1}^{b_{1}}X_{2}^{b_{2}},g^{d}x_{1}^{e_{1}}x_{2}^{e_{2}}) \\
G^{a}X_{1}^{b_{1}}X_{2}^{b_{2}-l_{2}}\otimes
g^{d}x_{1}^{e_{1}}x_{2}^{e_{2}-u_{2}} \\
=\sum_{\substack{ a,b_{1},b_{2},d,e_{1}=0  \\ a+b_{1}+b_{2}+d+e_{1}\equiv 1}}%
^{1}\left( -1\right) ^{\alpha \left( 1_{H};0,0,0,1\right) }B(x_{1}\otimes
1_{H};G^{a}X_{1}^{b_{1}}X_{2}^{b_{2}},g^{d}x_{1}^{e_{1}}x_{2})G^{a}X_{1}^{b_{1}}X_{2}^{b_{2}}\otimes g^{d}x_{1}^{e_{1}}+
\\
\sum_{\substack{ a,b_{1},d,e_{1},e_{2}=0  \\ a+b_{1}+d+e_{1}+e_{2}\equiv 1}}%
^{1}\left( -1\right) ^{\alpha \left( 1_{H};0,1,0,0\right) }B(x_{1}\otimes
1_{H};G^{a}X_{1}^{b_{1}}X_{2},g^{d}x_{1}^{e_{1}}x_{2}^{e_{2}})G^{a}X_{1}^{b_{1}}\otimes g^{d}x_{1}^{e_{1}}x_{2}^{e_{2}}+
\end{gather*}%
Since $\alpha \left( 1_{H};0,0,0,1\right) =a+b_{1}+b_{2}$ and $\alpha \left(
1_{H};0,1,0,0\right) =0$ we obtain%
\begin{eqnarray*}
&&\sum_{\substack{ a,b_{1},b_{2},d,e_{1}=0  \\ a+b_{1}+b_{2}+d+e_{1}\equiv 1
}}^{1}\left( -1\right) ^{a+b_{1}+b_{2}}B(x_{1}\otimes
1_{H};G^{a}X_{1}^{b_{1}}X_{2}^{b_{2}},g^{d}x_{1}^{e_{1}}x_{2})G^{a}X_{1}^{b_{1}}X_{2}^{b_{2}}\otimes g^{d}x_{1}^{e_{1}}+
\\
&&\sum_{\substack{ a,b_{1},d,e_{1},e_{2}=0  \\ a+b_{1}+d+e_{1}+e_{2}\equiv 1
}}^{1}B(x_{1}\otimes
1_{H};G^{a}X_{1}^{b_{1}}X_{2},g^{d}x_{1}^{e_{1}}x_{2}^{e_{2}})G^{a}X_{1}^{b_{1}}\otimes g^{d}x_{1}^{e_{1}}x_{2}^{e_{2}}=0
\end{eqnarray*}%
Now we proceed as in the case above.

\subsubsection{\protect\bigskip $G^{a}\otimes g^{d}$}

\begin{equation*}
\sum_{\substack{ a,d  \\ a+d=1}}^{1}\left[ \left( -1\right)
^{a}B(x_{1}\otimes 1_{H};G^{a},g^{d}x_{2})+B(x_{1}\otimes
1_{H};G^{a}X_{2},g^{d})\right] G^{a}\otimes g^{d}=0
\end{equation*}%
and we get%
\begin{eqnarray}
B(x_{1}\otimes 1_{H};1_{A},gx_{2})+B(x_{1}\otimes 1_{H};X_{2},g) &=&0
\label{x1ot1, twentysix} \\
-B(x_{1}\otimes 1_{H};G,x_{2})+B(x_{1}\otimes 1_{H};GX_{2},1_{H}) &=&0
\label{x1ot1, twentyseven}
\end{eqnarray}

\subsubsection{$G^{a}\otimes g^{d}x_{2}$}

\begin{equation*}
\sum_{\substack{ a,d  \\ a+d=0}}^{1}B(x_{1}\otimes
1_{H};G^{a}X_{2},g^{d}x_{2})G^{a}\otimes g^{d}x_{2}=0
\end{equation*}%
and we get%
\begin{eqnarray}
B(x_{1}\otimes 1_{H};X_{2},x_{2}) &=&0  \label{x1ot1, twentyeight} \\
B(x_{1}\otimes 1_{H};GX_{2},gx_{2}) &=&0  \label{x1ot1, twentynine}
\end{eqnarray}

\subsubsection{$G^{a}\otimes g^{d}x_{1}$}

\begin{equation*}
\sum_{\substack{ a,d  \\ a+d=0}}^{1}\left[ \left( -1\right)
^{a}B(x_{1}\otimes 1_{H};G^{a},g^{d}x_{1}x_{2})+B(x_{1}\otimes
1_{H};G^{a}X_{2},g^{d}x_{1})\right] G^{a}\otimes g^{d}x_{1}=0
\end{equation*}%
and we get%
\begin{eqnarray}
B(x_{1}\otimes 1_{H};1_{A},x_{1}x_{2})+B(x_{1}\otimes 1_{H};X_{2},x_{1}) &=&0
\label{x1ot1, thirty} \\
-B(x_{1}\otimes 1_{H};G,gx_{1}x_{2})+B(x_{1}\otimes 1_{H};GX_{2},gx_{1}) &=&0
\label{x1ot1, thirtyone}
\end{eqnarray}

\subsubsection{$G^{a}X_{2}\otimes g^{d}$}

\begin{equation*}
\sum_{\substack{ a,d  \\ a+d=0}}^{1}\left( -1\right) ^{a+1}B(x_{1}\otimes
1_{H};G^{a}X_{2},g^{d}x_{2})G^{a}X_{2}\otimes g^{d}=0
\end{equation*}%
and we get%
\begin{eqnarray*}
-B(x_{1}\otimes 1_{H};X_{2},x_{2}) &=&0 \\
B(x_{1}\otimes 1_{H};GX_{2},gx_{2}) &=&0
\end{eqnarray*}%
which we already obtained.

\subsubsection{$G^{a}X_{1}\otimes g^{d}$}

\begin{equation*}
\sum_{\substack{ a,d=0  \\ a+d=0}}^{1}\left[ \left( -1\right)
^{a+1}B(x_{1}\otimes 1_{H};G^{a}X_{1},g^{d}x_{2})+B(x_{1}\otimes
1_{H};G^{a}X_{1}X_{2},g^{d})\right] G^{a}X_{1}\otimes g^{d}=0
\end{equation*}%
and we get%
\begin{eqnarray}
-B(x_{1}\otimes 1_{H};X_{1},x_{2})+B(x_{1}\otimes 1_{H};X_{1}X_{2},1_{H})
&=&0  \label{x1ot1, thirtytwo} \\
B(x_{1}\otimes 1_{H};GX_{1},gx_{2})+B(x_{1}\otimes 1_{H};GX_{1}X_{2},g) &=&0
\label{x1ot1, thirtythree}
\end{eqnarray}

\subsubsection{$G^{a}\otimes g^{d}x_{1}x_{2}$}

\begin{equation*}
\sum_{\substack{ a,d=0  \\ a+d=1}}^{1}B(x_{1}\otimes
1_{H};G^{a}X_{2},g^{d}x_{1}x_{2})G^{a}\otimes g^{d}x_{1}x_{2}=0
\end{equation*}%
and we get%
\begin{eqnarray}
B(x_{1}\otimes 1_{H};X_{2},gx_{1}x_{2}) &=&0  \label{x1ot1, thirtyfour} \\
B(x_{1}\otimes 1_{H};GX_{2},x_{1}x_{2}) &=&0  \label{x1ot1, thirtyfive}
\end{eqnarray}

\subsubsection{$G^{a}X_{2}\otimes g^{d}x_{2}$}

We do not have any summand like this.

\subsubsection{$G^{a}X_{1}\otimes g^{d}x_{2}$}

\begin{equation*}
\sum_{\substack{ a,d=0  \\ a+d=1}}^{1}B(x_{1}\otimes
1_{H};G^{a}X_{1}X_{2},g^{d}x_{2})G^{a}X_{1}\otimes g^{d}x_{2}=0
\end{equation*}%
and we get%
\begin{eqnarray}
B(x_{1}\otimes 1_{H};X_{1}X_{2},gx_{2}) &=&0  \label{x1ot1, thirtysix} \\
B(x_{1}\otimes 1_{H};GX_{1}X_{2},x_{2}) &=&0  \label{x1ot1, thirtyseven}
\end{eqnarray}

\subsubsection{$G^{a}X_{2}\otimes g^{d}x_{1}$}

\begin{equation*}
\sum_{\substack{ a,d=0  \\ a+d=1}}^{1}\left( -1\right) ^{a+1}B(x_{1}\otimes
1_{H};G^{a}X_{2},g^{d}x_{1}x_{2})G^{a}X_{2}\otimes g^{d}x_{1}=0
\end{equation*}%
and we get%
\begin{eqnarray}
-B(x_{1}\otimes 1_{H};X_{2},gx_{1}x_{2}) &=&0  \label{x1ot1, thirtyeight} \\
B(x_{1}\otimes 1_{H};GX_{2},x_{1}x_{2}) &=&0  \label{x1ot1, thirtynine}
\end{eqnarray}

\subsubsection{$G^{a}X_{1}\otimes g^{d}x_{1}$}

\begin{gather*}
\sum_{\substack{ a,d=0  \\ a+d=1}}^{1}\left[ \left( -1\right)
^{a+1}B(x_{1}\otimes 1_{H};G^{a}X_{1},g^{d}x_{1}x_{2})+B(x_{1}\otimes
1_{H};G^{a}X_{1}X_{2},g^{d}x_{1})\right] \\
G^{a}X_{1}\otimes g^{d}x_{1}=0
\end{gather*}%
and we get%
\begin{eqnarray}
-B(x_{1}\otimes 1_{H};X_{1},gx_{1}x_{2})+B(x_{1}\otimes
1_{H};X_{1}X_{2},gx_{1}) &=&0  \label{x1ot1, forty} \\
B(x_{1}\otimes 1_{H};GX_{1},x_{1}x_{2})+B(x_{1}\otimes
1_{H};GX_{1}X_{2},x_{1}) &=&0  \label{x1ot1, fortyone}
\end{eqnarray}

\subsubsection{$G^{a}X_{1}X_{2}\otimes g^{d}$}

\begin{equation*}
\sum_{\substack{ a,d=0  \\ a+d=1}}^{1}\left( -1\right) ^{a}B(x_{1}\otimes
1_{H};G^{a}X_{1}X_{2},g^{d}x_{2})G^{a}X_{1}X_{2}\otimes g^{d}=0
\end{equation*}%
and we obtain%
\begin{eqnarray*}
B(x_{1}\otimes 1_{H};X_{1}X_{2},gx_{2}) &=&0 \\
-B(x_{1}\otimes 1_{H};GX_{1}X_{2},x_{2}) &=&0
\end{eqnarray*}%
which we already got.

\subsubsection{$G^{a}X_{2}\otimes g^{d}x_{1}x_{2}$}

We do not have any summand like this.

\subsubsection{$G^{a}X_{1}\otimes g^{d}x_{1}x_{2}$}

\begin{equation*}
\sum_{\substack{ a,d=0  \\ a+d=0}}^{1}B(x_{1}\otimes
1_{H};G^{a}X_{1}X_{2},g^{d}x_{1}x_{2})G^{a}X_{1}\otimes g^{d}x_{1}x_{2}=0
\end{equation*}%
and we get%
\begin{eqnarray}
B(x_{1}\otimes 1_{H};X_{1}X_{2},x_{1}x_{2}) &=&0  \label{x1ot1, fortytwo} \\
B(x_{1}\otimes 1_{H};GX_{1}X_{2},gx_{1}x_{2}) &=&0  \label{x1ot1, fortythree}
\end{eqnarray}

\subsubsection{$G^{a}X_{1}X_{2}\otimes g^{d}x_{2}$}

We do not have any summand like this.

\subsubsection{$G^{a}X_{1}X_{2}\otimes g^{d}x_{1}$}

\begin{equation*}
\sum_{\substack{ a,d=0  \\ a+d=1}}^{1}\left( -1\right) ^{a}B(x_{1}\otimes
1_{H};G^{a}X_{1}X_{2},g^{d}x_{1}x_{2})G^{a}X_{1}X_{2}\otimes g^{d}x_{1}=0
\end{equation*}%
and we get%
\begin{eqnarray*}
B(x_{1}\otimes 1_{H};X_{1}X_{2},x_{1}x_{2}) &=&0 \\
-B(x_{1}\otimes 1_{H};GX_{1}X_{2},gx_{1}x_{2}) &=&0
\end{eqnarray*}%
which we already got.

\subsubsection{$G^{a}X_{1}X_{2}\otimes g^{d}x_{1}x_{2}$}

We do not have any summand like this.

\subsection{Case $x_{1}x_{2}$}

\begin{eqnarray*}
&&\sum_{a,b_{1},b_{2},d,e_{1},e_{2}=0}^{1}\sum_{l_{1}=0}^{b_{1}}%
\sum_{l_{2}=0}^{b_{2}}\sum_{u_{1}=0}^{e_{1}}\sum_{u_{2}=0}^{e_{2}}\left(
-1\right) ^{\alpha \left( x_{1};l_{1},l_{2},u_{1},u_{2}\right) } \\
&&B(g\otimes
1_{H};G^{a}X_{1}^{b_{1}}X_{2}^{b_{2}},g^{d}x_{1}^{e_{1}}x_{2}^{e_{2}}) \\
&&G^{a}X_{1}^{b_{1}-l_{1}}X_{2}^{b_{2}-l_{2}}\otimes
g^{d}x_{1}^{e_{1}-u_{1}}x_{2}^{e_{2}-u_{2}}\otimes
g^{a+b_{1}+b_{2}+l_{1}+l_{2}+d+e_{1}+e_{2}+u_{1}+u_{2}}x_{1}^{l_{1}+u_{1}+1}x_{2}^{l_{2}+u_{2}}+
\\
&&+\sum_{a,b_{1},b_{2},d,e_{1},e_{2}=0}^{1}\sum_{l_{1}=0}^{b_{1}}%
\sum_{l_{2}=0}^{b_{2}}\sum_{u_{1}=0}^{e_{11}}\sum_{u_{2}=0}^{e_{12}}\left(
-1\right) ^{\alpha \left( 1_{H};l_{1},l_{2},u_{1},u_{2}\right) } \\
&&B(x_{1}\otimes
1_{H};G^{a}X_{1}^{b_{1}}X_{2}^{b_{2}},g^{d}x_{1}^{e_{1}}x_{2}^{e_{2}}) \\
&&G^{a}X_{1}^{b_{1}-l_{1}}X_{2}^{b_{2}-l_{2}}\otimes
g^{d}x_{1}^{e_{1}-u_{1}}x_{2}^{e_{2}-u_{2}}\otimes
g^{a+b_{1}+b_{2}+l_{1}+l_{2}+d+e_{11}+e_{12}+u_{1}+u_{2}}x_{1}^{l_{1}+u_{1}}x_{2}^{l_{2}+u_{2}}+
\\
&=&B^{A}(x_{1}\otimes 1_{H})\otimes B^{H}(x_{1}\otimes 1_{H})\otimes 1_{H}
\end{eqnarray*}%
For the first summand we get%
\begin{eqnarray*}
l_{1} &=&u_{1}=0 \\
l_{2}+u_{2} &=&1 \\
a+b_{1}+b_{2}+d+e_{1}+e_{2} &\equiv &1 \\
\alpha \left( x_{1};0,0,0,1\right) &\equiv &1 \\
\alpha \left( x_{1};0,1,0,0\right) &\equiv &a+b_{1}+b_{2}+1
\end{eqnarray*}%
and we obtain%
\begin{eqnarray*}
&&\sum_{\substack{ a,b_{1},b_{2},d,e_{1},e_{2}=0  \\ %
a+b_{1}+b_{2}+d+e_{1}+e_{2}\equiv 1}}^{1}\sum_{l_{2}=0}^{b_{2}}\sum
_{\substack{ u_{2}=0  \\ l_{2}+u_{2}=1}}^{e_{2}}\left( -1\right) ^{\alpha
\left( x_{1};0,l_{2},0,u_{2}\right) }B(g\otimes
1_{H};G^{a}X_{1}^{b_{1}}X_{2}^{b_{2}},g^{d}x_{1}^{e_{1}}x_{2}^{e_{2}}) \\
&&G^{a}X_{1}^{b_{1}}X_{2}^{b_{2}-l_{2}}\otimes
g^{d}x_{1}^{e_{1}}x_{2}^{e_{2}-u_{2}} \\
&=&\sum_{\substack{ a,b_{1},b_{2},d,e_{1}=0  \\ a+b_{1}+b_{2}+d+e_{1}\equiv
0 }}^{1}\left( -1\right) ^{\alpha \left( x_{1};0,0,0,1\right) }B(g\otimes
1_{H};G^{a}X_{1}^{b_{1}}X_{2}^{b_{2}},g^{d}x_{1}^{e_{1}}x_{2})G^{a}X_{1}^{b_{1}}X_{2}^{b_{2}}\otimes g^{d}x_{1}^{e_{1}}+
\\
&&\sum_{\substack{ a,b_{1},d,e_{1},e_{2}=0  \\ a+b_{1}+d+e_{1}+e_{2}\equiv 0
}}^{1}\left( -1\right) ^{\alpha \left( x_{1};0,1,0,0\right) }B(g\otimes
1_{H};G^{a}X_{1}^{b_{1}}X_{2},g^{d}x_{1}^{e_{1}}x_{2}^{e_{2}})G^{a}X_{1}^{b_{1}}\otimes g^{d}x_{1}^{e_{1}}x_{2}^{e_{2}-u_{2}}+
\end{eqnarray*}%
\begin{eqnarray*}
&&1\sum_{\substack{ a,d=0  \\ a+d\equiv 0}}^{1}-B(g\otimes
1_{H};G^{a},g^{d}x_{2})G^{a}\otimes g^{d} \\
&&2\sum_{\substack{ a,d=0  \\ a+d\equiv 1}}^{1}-B(g\otimes
1_{H};G^{a},g^{d}x_{1}x_{2})G^{a}\otimes g^{d}x_{1} \\
&&3\sum_{\substack{ a,d=0  \\ a+d\equiv 1}}^{1}-B(g\otimes
1_{H};G^{a}X_{2},g^{d}x_{2})G^{a}X_{2}\otimes g^{d} \\
&&4\sum_{\substack{ a,d=0  \\ a+d\equiv 0}}^{1}-B(g\otimes
1_{H};G^{a}X_{2},g^{d}x_{1}x_{2})G^{a}X_{2}\otimes g^{d}x_{1} \\
&&5\sum_{\substack{ a,d=0  \\ a+d\equiv 1}}^{1}-B(g\otimes
1_{H};G^{a}X_{1},g^{d}x_{2})G^{a}X_{1}\otimes g^{d} \\
&&6\sum_{\substack{ a,d=0  \\ a+d\equiv 0}}^{1}-B(g\otimes
1_{H};G^{a}X_{1},g^{d}x_{1}x_{2})G^{a}X_{1}\otimes g^{d}x_{1} \\
&&7\sum_{\substack{ a,d=0  \\ a+d\equiv 0}}^{1}-B(g\otimes
1_{H};G^{a}X_{1}X_{2},g^{d}x_{2})G^{a}X_{1}X_{2}\otimes g^{d} \\
&&8\sum_{\substack{ a,d=0  \\ a+d\equiv 1}}^{1}-B(g\otimes
1_{H};G^{a}X_{1}X_{2},g^{d}x_{1}x_{2})G^{a}X_{1}X_{2}\otimes g^{d}x_{1} \\
&&\sum_{\substack{ a,d=0  \\ a+d\equiv 0}}^{1}\left( -1\right)
^{a}B(g\otimes 1_{H};G^{a}X_{2},g^{d})G^{a}\otimes g^{d}+ \\
&&\sum_{\substack{ a,d=0  \\ a+d\equiv 1}}^{1}\left( -1\right)
^{a}B(g\otimes 1_{H};G^{a}X_{2},g^{d}x_{2})G^{a}\otimes g^{d}x_{2} \\
&&\sum_{\substack{ a,d=0  \\ a+d\equiv 1}}^{1}\left( -1\right)
^{a}B(g\otimes 1_{H};G^{a}X_{2},g^{d}x_{1})G^{a}\otimes g^{d}x_{1} \\
&&\sum_{\substack{ a,d=0  \\ a+d\equiv 0}}^{1}\left( -1\right)
^{a}B(g\otimes 1_{H};G^{a}X_{2},g^{d}x_{1}x_{2})G^{a}\otimes g^{d}x_{1}x_{2}
\\
&&\sum_{\substack{ a,d=0  \\ a+d\equiv 1}}^{1}\left( -1\right)
^{a+1}B(g\otimes 1_{H};G^{a}X_{1}X_{2},g^{d})G^{a}X_{1}\otimes g^{d} \\
&&\sum_{\substack{ a,d=0  \\ a+d\equiv 0}}^{1}\left( -1\right)
^{a+1}B(g\otimes 1_{H};G^{a}X_{1}X_{2},g^{d}x_{2})G^{a}X_{1}\otimes
g^{d}x_{2} \\
&&\sum_{\substack{ a,d=0  \\ a+d\equiv 0}}^{1}\left( -1\right)
^{a+1}B(g\otimes 1_{H};G^{a}X_{1}X_{2},g^{d}x_{1})G^{a}X_{1}\otimes
g^{d}x_{1} \\
&&\sum_{\substack{ a,d=0  \\ a+d\equiv 1}}^{1}\left( -1\right)
^{a+1}B(g\otimes 1_{H};G^{a}X_{1}X_{2},g^{d}x_{1}x_{2})G^{a}X_{1}\otimes
g^{d}x_{1}x_{2}
\end{eqnarray*}%
The second summand gives%
\begin{eqnarray*}
&&\sum_{\substack{ a,b_{1},b_{2},d,e_{1},e_{2}=0  \\ %
a+b_{1}+b_{2}+d+e_{1}+e_{2}\equiv 0}}^{1}\sum_{l_{1}=0}^{b_{1}}%
\sum_{l_{2}=0}^{b_{2}}\sum_{\substack{ u_{1}=0  \\ l_{1}+u_{1}=1}}%
^{e_{1}}\sum _{\substack{ u_{2}=0  \\ l_{2}+u_{2}=1}}^{e_{1}}\left(
-1\right) ^{\alpha \left( 1_{H};l_{1},l_{2},u_{1},u_{2}\right) } \\
&&B(x_{1}\otimes
1_{H};G^{a}X_{1}^{b_{1}}X_{2}^{b_{2}},g^{d}x_{1}^{e_{1}}x_{2}^{e_{2}})G^{a}X_{1}^{b_{1}-l_{1}}X_{2}^{b_{2}-l_{2}}\otimes g^{d}x_{1}^{e_{1}-u_{1}}x_{2}^{e_{2}-u_{2}}
\\
&=&\sum_{\substack{ a,d,e_{1},e_{2}=0  \\ a+d+e_{1}+e_{2}\equiv 0}}%
^{1}\left( -1\right) ^{\alpha \left( 1_{H};1,1,0,0\right) }B(x_{1}\otimes
1_{H};G^{a}X_{1}X_{2},g^{d}x_{1}^{e_{1}}x_{2}^{e_{2}})G^{a}\otimes
g^{d}x_{1}^{e_{1}}x_{2}^{e_{2}}+ \\
&&+\sum_{\substack{ a,b_{1},d,e_{2}=0  \\ a+b_{1}+d+e_{2}\equiv 0}}%
^{1}\left( -1\right) ^{\alpha \left( 1_{H};0,1,1,0\right) }B(x_{1}\otimes
1_{H};G^{a}X_{1}^{b_{1}}X_{2},g^{d}x_{1}x_{2}^{e_{2}})G^{a}X_{1}^{b_{1}}%
\otimes g^{d}x_{2}^{e_{2}}+ \\
&&+\sum_{\substack{ a,b_{2},d,e_{1}=0  \\ a+b_{2}+d+e_{1}\equiv 0}}%
^{1}\left( -1\right) ^{\alpha \left( 1_{H};1,0,0,1\right) }B(x_{1}\otimes
1_{H};G^{a}X_{1}X_{2}^{b_{2}},g^{d}x_{1}^{e_{1}}x_{2})G^{a}X_{2}^{b_{2}}%
\otimes g^{d}x_{1}^{e_{1}}+ \\
&&\sum_{\substack{ a,b_{1},b_{2},d=0  \\ a+b_{1}+b_{2}+d\equiv 0}}^{1}\left(
-1\right) ^{\alpha \left( 1_{H};0,0,1,1\right) }B(x_{1}\otimes
1_{H};G^{a}X_{1}^{b_{1}}X_{2}^{b_{2}},g^{d}x_{1}x_{2})G^{a}X_{1}^{b_{1}}X_{2}^{b_{2}}\otimes g^{d}
\end{eqnarray*}

and we get%
\begin{eqnarray*}
&&\alpha \left( 1_{H};1,1,0,0\right) \equiv 1+b_{2}\equiv 0 \\
&&\alpha \left( 1_{H};0,1,1,0\right) \equiv e_{2}+a+b_{1} \\
&&\alpha \left( 1_{H};1,0,0,1\right) \equiv a+1 \\
&&\alpha \left( 1_{H};0,0,1,1\right) \equiv 1+e_{2}\equiv 0
\end{eqnarray*}%
\begin{eqnarray*}
&&\sum_{\substack{ a,d,e_{1},e_{2}=0  \\ a+d+e_{1}+e_{2}\equiv 0}}%
^{1}B(x_{1}\otimes
1_{H};G^{a}X_{1}X_{2},g^{d}x_{1}^{e_{1}}x_{2}^{e_{2}})G^{a}\otimes
g^{d}x_{1}^{e_{1}}x_{2}^{e_{2}} \\
&&\sum_{\substack{ a,b_{1},d,e_{2}=0  \\ a+b_{1}+d+e_{2}\equiv 0}}^{1}\left(
-1\right) ^{e_{2}+a+b_{1}}B(x_{1}\otimes
1_{H};G^{a}X_{1}^{b_{1}}X_{2},g^{d}x_{1}x_{2}^{e_{2}})G^{a}X_{1}^{b_{1}}%
\otimes g^{d}x_{2}^{e_{2}} \\
&&\sum_{\substack{ a,,b_{2},d,e_{1}=0  \\ a+b_{2}+d+e_{1}\equiv 0}}%
^{1}\left( -1\right) ^{a+1}B(x_{1}\otimes
1_{H};G^{a}X_{1}X_{2}^{b_{2}},g^{d}x_{1}^{e_{1}}x_{2})G^{a}X_{2}^{b_{2}}%
\otimes g^{d}x_{1}^{e_{1}} \\
&&\sum_{\substack{ a,b_{1},b_{2},d=0  \\ a+b_{1}+b_{2}+d\equiv 0}}^{1}\left(
-1\right) ^{\alpha \left( 1_{H};0,0,1,1\right) }B(x_{1}\otimes
1_{H};G^{a}X_{1}^{b_{1}}X_{2}^{b_{2}},g^{d}x_{1}x_{2})G^{a}X_{1}^{b_{1}}X_{2}^{b_{2}}\otimes g^{d}
\end{eqnarray*}%
\begin{eqnarray*}
&&\sum_{\substack{ a,d,e_{1},e_{2}=0  \\ a+d+e_{1}+e_{2}\equiv 1}}%
^{1}B(x_{1}\otimes
1_{H};G^{a}X_{1}X_{2},g^{d}x_{1}^{e_{1}}x_{2}^{e_{2}})G^{a}\otimes
g^{d}x_{1}^{e_{1}}x_{2}^{e_{2}} \\
&&\sum_{\substack{ a,b_{1},d,e_{2}=0  \\ a+b_{1}+d+e_{2}\equiv 1}}^{1}\left(
-1\right) ^{a}B(x_{1}\otimes
1_{H};G^{a}X_{1}^{b_{1}}X_{2},g^{d}x_{1}x_{2}^{e_{2}})G^{a}X_{1}^{b_{1}}%
\otimes g^{d}x_{2}^{e_{2}} \\
&&\sum_{\substack{ a,,b_{2},d,e_{1}=0  \\ a+b_{2}+d+e_{1}\equiv 1}}%
^{1}\left( -1\right) ^{a+1}B(x_{1}\otimes
1_{H};G^{a}X_{1}X_{2}^{b_{2}},g^{d}x_{1}^{e_{1}}x_{2})G^{a}X_{2}^{b_{2}}%
\otimes g^{d}x_{1}^{e_{1}} \\
&&\sum_{\substack{ a,b_{1},b_{2},d=0  \\ a+b_{1}+b_{2}+d\equiv 1}}%
^{1}B(x_{1}\otimes
1_{H};G^{a}X_{1}^{b_{1}}X_{2}^{b_{2}},g^{d}x_{1}x_{2})G^{a}X_{1}^{b_{1}}X_{2}^{b_{2}}\otimes g^{d}
\end{eqnarray*}%
\begin{eqnarray*}
1 &&\sum_{\substack{ a,d  \\ a+d\equiv 0}}^{1}B(x_{1}\otimes
1_{H};G^{a}X_{1}X_{2},g^{d})G^{a}\otimes g^{d} \\
&&2\sum_{\substack{ a,d  \\ a+d\equiv 1}}^{1}B(x_{1}\otimes
1_{H};G^{a}X_{1}X_{2},g^{d}x_{1})G^{a}\otimes g^{d}x_{1} \\
&&3\sum_{\substack{ a,d=0  \\ a+d\equiv 1}}^{1}B(x_{1}\otimes
1_{H};G^{a}X_{1}X_{2},g^{d}x_{2})G^{a}\otimes g^{d}x_{2} \\
&&4\sum_{\substack{ a,d=0  \\ a+d\equiv 0}}^{1}B(x_{1}\otimes
1_{H};G^{a}X_{1}X_{2},g^{d}x_{1}x_{2})G^{a}\otimes g^{d}x_{1}x_{2} \\
5 &&\sum_{\substack{ a,d=0  \\ a+d\equiv 0}}^{1}\left( -1\right)
^{+a+1}B(x_{1}\otimes 1_{H};G^{a}X_{2},g^{d}x_{1})G^{a}\otimes g^{d} \\
&&6\sum_{\substack{ a,d=0  \\ a+d\equiv 1}}^{1}\left( -1\right)
^{+a+1}B(x_{1}\otimes 1_{H};G^{a}X_{1}X_{2},g^{d}x_{1})G^{a}X_{1}\otimes
g^{d} \\
&&7\sum_{\substack{ a,d=0  \\ a+d\equiv 1}}^{1}\left( -1\right)
^{+a+1}B(x_{1}\otimes 1_{H};G^{a}X_{2},g^{d}x_{1}x_{2})G^{a}\otimes
g^{d}x_{2} \\
&&8\sum_{\substack{ a,d=0  \\ a+d\equiv 0}}^{1}\left( -1\right)
^{a}B(x_{1}\otimes 1_{H};G^{a}X_{1}X_{2},g^{d}x_{1}x_{2})G^{a}X_{1}\otimes
g^{d}x_{2} \\
9 &&\sum_{\substack{ a,d=0  \\ a+d\equiv 0}}^{1}\left( -1\right)
^{a+1}B(x_{1}\otimes 1_{H};G^{a}X_{1},g^{d}x_{2})G^{a}\otimes g^{d} \\
&&10\sum_{\substack{ a,d  \\ a+d\equiv 1}}^{1}\left( -1\right)
^{a+1}B(x_{1}\otimes 1_{H};G^{a}X_{1}X_{2},g^{d}x_{2})G^{a}X_{2}\otimes g^{d}
\\
&&11\sum_{\substack{ a,d  \\ a+d\equiv 1}}^{1}\left( -1\right)
^{a+1}B(x_{1}\otimes 1_{H};G^{a}X_{1},g^{d}x_{1}x_{2})G^{a}\otimes g^{d}x_{1}
\\
&&12\sum_{\substack{ a,d=0  \\ a+d\equiv 0}}^{1}\left( -1\right)
^{a+1}B(x_{1}\otimes 1_{H};G^{a}X_{1}X_{2},g^{d}x_{1}x_{2})G^{a}X_{2}\otimes
g^{d}x_{1} \\
13 &&\sum_{\substack{ a,d=0  \\ a+d\equiv 0}}^{1}B(x_{1}\otimes
1_{H};G^{a},g^{d}x_{1}x_{2})G^{a}\otimes g^{d} \\
&&14\sum_{\substack{ a,d=0  \\ a+d\equiv 1}}^{1}B(x_{1}\otimes
1_{H};G^{a}X_{1},g^{d}x_{1}x_{2})G^{a}X_{1}\otimes g^{d} \\
&&15\sum_{\substack{ a,d=0  \\ a+d\equiv 1}}^{1}B(x_{1}\otimes
1_{H};G^{a}X_{2},g^{d}x_{1}x_{2})G^{a}X_{2}\otimes g^{d} \\
&&16\sum_{\substack{ a,d=0  \\ a+d\equiv 0}}^{1}B(x_{1}\otimes
1_{H};G^{a}X_{1}X_{2},g^{d}x_{1}x_{2})G^{a}X_{1}X_{2}\otimes g^{d}
\end{eqnarray*}%
Thus we obtain

\begin{equation*}
\sum_{\substack{ a,d  \\ a+d\equiv 0}}^{1}\left[
\begin{array}{c}
B(x_{1}\otimes 1_{H};G^{a}X_{1}X_{2},g^{d})+\left( -1\right)
^{a}B(x_{1}\otimes 1_{H};G^{a}X_{2},g^{d}x_{1}) \\
\left( -1\right) ^{a+1}B(x_{1}\otimes
1_{H};G^{a}X_{1},g^{d}x_{2})+B(x_{1}\otimes 1_{H};G^{a},g^{d}x_{1}x_{2}) \\
\left( -1\right) ^{a}B(g\otimes 1_{H};G^{a}X_{2},g^{d})-B(g\otimes
1_{H};G^{a},g^{d}x_{2})%
\end{array}%
\right] G^{a}\otimes g^{d}=0
\end{equation*}%
and we get%
\begin{equation}
\begin{array}{c}
B(x_{1}\otimes 1_{H};X_{1}X_{2},1_{H})+B(x_{1}\otimes 1_{H};X_{2},x_{1}) \\
-B(x_{1}\otimes 1_{H};X_{1},x_{2})+B(x_{1}\otimes 1_{H};1_{A},x_{1}x_{2}) \\
+B(g\otimes 1_{H};X_{2},1_{H})-B(g\otimes 1_{H};1_{A},x_{2})=0%
\end{array}
\label{x1ot1, forthyfour}
\end{equation}%
\begin{equation}
\begin{array}{c}
B(x_{1}\otimes 1_{H};GX_{1}X_{2},g)-B(x_{1}\otimes 1_{H};GX_{2},gx_{1}) \\
+B(x_{1}\otimes 1_{H};GX_{1},gx_{2})+B(x_{1}\otimes 1_{H};G,gx_{1}x_{2}) \\
-B(g\otimes 1_{H};GX_{2},g)-B(g\otimes 1_{H};G,gx_{2})=0%
\end{array}
\label{x1ot1, forthyfive}
\end{equation}%
\begin{gather*}
\sum_{\substack{ a,d  \\ a+d\equiv 1}}^{1}\left[
\begin{array}{c}
B(x_{1}\otimes 1_{H};G^{a}X_{1}X_{2},g^{d}x_{1})+\left( -1\right)
^{a+1}B(x_{1}\otimes 1_{H};G^{a}X_{1},g^{d}x_{1}x_{2})+ \\
\left( -1\right) ^{a}B(g\otimes 1_{H};G^{a}X_{2},g^{d}x_{1})-B(g\otimes
1_{H};G^{a},g^{d}x_{1}x_{2})%
\end{array}%
\right] \\
G^{a}\otimes g^{d}x_{1}=0
\end{gather*}%
and we get%
\begin{equation}
\begin{array}{c}
B(x_{1}\otimes 1_{H};X_{1}X_{2},gx_{1})-B(x_{1}\otimes
1_{H};X_{1},gx_{1}x_{2})+ \\
B(g\otimes 1_{H};X_{2},gx_{1})-B(g\otimes 1_{H};1_{A},gx_{1}x_{2})%
\end{array}%
=0  \label{x1ot1, fortysix}
\end{equation}%
\begin{equation}
\begin{array}{c}
B(x_{1}\otimes 1_{H};GX_{1}X_{2},x_{1})+B(x_{1}\otimes
1_{H};GX_{1},x_{1}x_{2})+ \\
-B(g\otimes 1_{H};GX_{2},x_{1})-B(g\otimes 1_{H};G,x_{1}x_{2})%
\end{array}%
=0  \label{x1ot1, fortyseven}
\end{equation}%
\begin{equation*}
\sum_{\substack{ a,d=0  \\ a+d\equiv 1}}^{1}\left[
\begin{array}{c}
B(x_{1}\otimes 1_{H};G^{a}X_{1}X_{2},g^{d}x_{2})+\left( -1\right)
^{a+1}B(x_{1}\otimes 1_{H};G^{a}X_{2},g^{d}x_{1}x_{2}) \\
+\left( -1\right) ^{a}B(g\otimes 1_{H};G^{a}X_{2},g^{d}x_{2})%
\end{array}%
\right] G^{a}\otimes g^{d}x_{2}=0
\end{equation*}%
and we get%
\begin{equation}
B(x_{1}\otimes 1_{H};X_{1}X_{2},gx_{2})-B(x_{1}\otimes
1_{H};X_{2},gx_{1}x_{2})+B(g\otimes 1_{H};X_{2},gx_{2})=0
\label{x1ot1, fortyeight}
\end{equation}%
\begin{equation}
B(x_{1}\otimes 1_{H};GX_{1}X_{2},x_{2})+B(x_{1}\otimes
1_{H};GX_{2},x_{1}x_{2})-B(g\otimes 1_{H};GX_{2},x_{2})=0
\label{x1ot1, fortynine}
\end{equation}%
\begin{equation*}
\sum_{\substack{ a,d=0  \\ a+d\equiv 0}}^{1}\left[ B(x_{1}\otimes
1_{H};G^{a}X_{1}X_{2},g^{d}x_{1}x_{2})+\left( -1\right) ^{a}B(g\otimes
1_{H};G^{a}X_{2},g^{d}x_{1}x_{2})\right] G^{a}\otimes g^{d}x_{1}x_{2}=0
\end{equation*}%
and we get%
\begin{equation}
B(x_{1}\otimes 1_{H};X_{1}X_{2},x_{1}x_{2})+B(g\otimes
1_{H};X_{2},x_{1}x_{2})=0  \label{x1ot1, fifty}
\end{equation}%
\begin{equation}
B(x_{1}\otimes 1_{H};GX_{1}X_{2},gx_{1}x_{2})-B(g\otimes
1_{H};GX_{2},gx_{1}x_{2})=0  \label{x1ot1, fiftyone}
\end{equation}%
and we get%
\begin{gather*}
\sum_{\substack{ a,d=0  \\ a+d\equiv 1}}^{1}\left[
\begin{array}{c}
\left( -1\right) ^{a+1}B(x_{1}\otimes
1_{H};G^{a}X_{1}X_{2},g^{d}x_{1})+B(x_{1}\otimes
1_{H};G^{a}X_{1},g^{d}x_{1}x_{2})+ \\
\left( -1\right) ^{a+1}B(g\otimes 1_{H};G^{a}X_{1}X_{2},g^{d})-B(g\otimes
1_{H};G^{a}X_{1},g^{d}x_{2})%
\end{array}%
\right] \\
G^{a}X_{1}\otimes g^{d}=0
\end{gather*}%
and we get%
\begin{equation}
-B(x_{1}\otimes 1_{H};X_{1}X_{2},gx_{1})+B(x_{1}\otimes
1_{H};X_{1},gx_{1}x_{2})-B(g\otimes 1_{H};X_{1}X_{2},g)-B(g\otimes
1_{H};X_{1},gx_{2})=0  \label{x1ot1, fiftytwo}
\end{equation}%
\begin{equation}
B(x_{1}\otimes 1_{H};GX_{1}X_{2},x_{1})+B(x_{1}\otimes
1_{H};GX_{1},x_{1}x_{2})+B(g\otimes 1_{H};GX_{1}X_{2},1_{H})-B(g\otimes
1_{H};GX_{1},x_{2})=0  \label{x1ot1, fiftythree}
\end{equation}%
\begin{gather*}
\sum_{\substack{ a,d=0  \\ a+d\equiv 0}}^{1}\left[ \left( -1\right)
^{a}B(x_{1}\otimes 1_{H};G^{a}X_{1}X_{2},g^{d}x_{1}x_{2})+\left( -1\right)
^{a+1}B(g\otimes 1_{H};G^{a}X_{1}X_{2},g^{d}x_{2})\right] \\
G^{a}X_{1}\otimes g^{d}x_{2}=0
\end{gather*}%
and we get%
\begin{equation}
B(x_{1}\otimes 1_{H};X_{1}X_{2},x_{1}x_{2})-B(g\otimes
1_{H};X_{1}X_{2},x_{2})=0  \label{x1ot1, fiftyfour}
\end{equation}%
\begin{equation}
-B(x_{1}\otimes 1_{H};GX_{1}X_{2},gx_{1}x_{2})+B(g\otimes
1_{H};GX_{1}X_{2},gx_{2})=0  \label{x1ot1, fiftyfive}
\end{equation}%
\begin{gather*}
\sum_{\substack{ a,d  \\ a+d\equiv 1}}^{1}\left[
\begin{array}{c}
\left( -1\right) ^{a+1}B(x_{1}\otimes 1_{H};G^{a}X_{1}X_{2},g^{d}x_{2}) \\
+B(x_{1}\otimes 1_{H};G^{a}X_{2},g^{d}x_{1}x_{2})-B(g\otimes
1_{H};G^{a}X_{2},g^{d}x_{2})%
\end{array}%
\right] \\
G^{a}X_{2}\otimes g^{d}=0
\end{gather*}%
and we get%
\begin{equation}
-B(x_{1}\otimes 1_{H};X_{1}X_{2},gx_{2})+B(x_{1}\otimes
1_{H};X_{2},gx_{1}x_{2})-B(g\otimes 1_{H};X_{2},gx_{2})=0
\label{x1ot1, fiftysix}
\end{equation}%
\begin{equation}
B(x_{1}\otimes 1_{H};GX_{1}X_{2},x_{2})+B(x_{1}\otimes
1_{H};GX_{2},x_{1}x_{2})-B(g\otimes 1_{H};GX_{2},x_{2})=0
\label{x1ot1, fiftyseven}
\end{equation}%
\begin{equation*}
\sum_{\substack{ a,d=0  \\ a+d\equiv 0}}^{1}\left[ \left( -1\right)
^{a+1}B(x_{1}\otimes 1_{H};G^{a}X_{1}X_{2},g^{d}x_{1}x_{2})-B(g\otimes
1_{H};G^{a}X_{2},g^{d}x_{1}x_{2})\right] G^{a}X_{2}\otimes g^{d}x_{1}=0
\end{equation*}%
and we get%
\begin{equation*}
-B(x_{1}\otimes 1_{H};X_{1}X_{2},x_{1}x_{2})-B(g\otimes
1_{H};X_{2},x_{1}x_{2})=0
\end{equation*}%
already got and%
\begin{equation}
B(x_{1}\otimes 1_{H};GX_{1}X_{2},gx_{1}x_{2})-B(g\otimes
1_{H};GX_{2},gx_{1}x_{2})=0  \label{x1ot1, fiftyeight}
\end{equation}%
\begin{equation*}
\sum_{\substack{ a,d=0  \\ a+d\equiv 0}}^{1}\left[ B(x_{1}\otimes
1_{H};G^{a}X_{1}X_{2},g^{d}x_{1}x_{2})-B(g\otimes
1_{H};G^{a}X_{1}X_{2},g^{d}x_{2})\right] G^{a}X_{1}X_{2}\otimes g^{d}=0
\end{equation*}%
and we get%
\begin{equation*}
B(x_{1}\otimes 1_{H};X_{1}X_{2},x_{1}x_{2})-B(g\otimes
1_{H};X_{1}X_{2},x_{2})=0
\end{equation*}%
already got%
\begin{equation}
B(x_{1}\otimes 1_{H};GX_{1}X_{2},gx_{1}x_{2})-B(g\otimes
1_{H};GX_{1}X_{2},gx_{2})=0  \label{x1ot1, fiftynine}
\end{equation}

\subsection{Case $gx_{1}x_{2}$}

The first summand of the left side

\begin{eqnarray*}
&&\sum_{\substack{ a,b_{1},b_{2},d,e_{1},e_{2}=0  \\ \overline{%
a+b_{1}+b_{2}+d+e_{1}+e_{2}}=1}}^{1}\sum_{l_{1}=0}^{b_{1}}%
\sum_{l_{2}=0}^{b_{2}}\sum_{u_{1}=0}^{e_{1}}\sum_{u_{2}=0}^{e_{2}}\left(
-1\right) ^{\alpha }B(g\otimes
1_{H};G^{a}X_{1}^{b_{1}}X_{2}^{b_{2}},g^{d}x_{1}^{e_{1}}x_{2}^{e_{2}}) \\
&&G^{a}X_{1}^{b_{1}-l_{1}}X_{2}^{b_{2}-l_{2}}\otimes g^{d}x_{1}^{\left(
e_{1}-u_{1}\right) }x_{2}^{\left( e_{2}-u_{2}\right) }\otimes
g^{a+b_{1}+b_{2}+d+e_{1}+e_{2}+l_{1}+l_{2}+u_{1}+u_{2}}x_{1}^{u_{1}+l_{1}+1}x_{2}^{l_{2}+u_{2}}
\end{eqnarray*}%
where%
\begin{eqnarray*}
\alpha \left( g^{m}x_{1}^{n_{1}}x_{2}^{n_{2}}\right) &=&\alpha \left(
g^{m}x_{1}^{n_{1}}x_{2}^{n_{2}};l_{1},l_{2},u_{1},u_{2}\right) =\left(
u_{2}+e_{2}\right) u_{1}+\left( l_{2}+b_{2}\right) l_{1} \\
&+&\left( a+b_{1}+b_{2}+l_{1}+l_{2}\right) \left(
u_{1}+u_{2}+n_{1}+n_{2}\right) +m(u_{1}+u_{2})+n_{1}u_{2}+l_{1}(n_{2}+u_{2})
\end{eqnarray*}%
\begin{eqnarray*}
\alpha &=&\alpha \left( a,b_{1},b_{2},e_{2},l_{1},l_{2},u_{1},u_{2}\right)
=\left( l_{2}+b_{2}\right) l_{1}+\left( u_{2}+e_{2}\right) u_{1}+ \\
&&+\left( a+b_{1}+b_{2}+l_{1}+l_{2}\right) \left( u_{1}+u_{2}+1\right)
+l_{1}u_{2}+u_{2}
\end{eqnarray*}%
we get%
\begin{eqnarray*}
a+b_{1}+b_{2}+d+e_{1}+e_{2}+l_{1}+l_{2}+u_{1}+u_{2} &\equiv &1 \\
u_{1}+l_{1} &\equiv &0 \\
l_{2}+u_{2} &\equiv &1
\end{eqnarray*}%
i.e.%
\begin{eqnarray*}
\overline{a+b_{1}+b_{2}+d+e_{1}+e_{2}} &\equiv &0 \\
u_{1}+l_{1} &\equiv &0 \\
l_{2}+u_{2} &\equiv &1
\end{eqnarray*}%
In view of $\left( \ref{got1 first}\right) $, we already know that $%
B(g\otimes
1_{H};G^{a}X_{1}^{b_{1}}X_{2}^{b_{2}},g^{d}x_{1}^{e_{1}}x_{2}^{e_{2}})=0$
whenever

$\overline{a+b_{1}+b_{2}+d+e_{1}+e_{2}}=0.$

The second summand of the left side is%
\begin{gather*}
\sum_{a,b_{1},b_{2},d,e_{1},e_{2}=0}^{1}\sum_{l_{1}=0}^{b_{1}}%
\sum_{l_{2}=0}^{b_{2}}\sum_{u_{1}=0}^{e_{1}}\sum_{u_{2}=0}^{e_{2}}(-1)^{%
\alpha \left( 1_{H},l_{1},l_{2},u_{1},u_{2}\right) } \\
B(x_{1}\otimes
1_{H};G^{a}X_{1}^{b_{1}}X_{2}^{b_{2}},g^{d}x_{1}^{e_{1}}x_{2}^{e_{2}}) \\
G^{a}X_{1}^{b_{1}-l_{1}}X_{2}^{b_{2}-l_{2}}\otimes g^{d}x_{1}^{\left(
e_{1}-u_{1}\right) }x_{2}^{\left( e_{2}-u_{2}\right) }\otimes
g^{a+b_{1}+b_{2}+d+e_{1}+e_{2}+l_{1}+l_{2}+u_{1}+u_{2}}x_{1}^{u_{1}+l_{1}}x_{2}^{u_{2}+l_{2}}
\end{gather*}

thus we get%
\begin{eqnarray*}
a+b_{1}+b_{2}+d+e_{1}+e_{2}+l_{1}+l_{2}+u_{1}+u_{2} &=&1 \\
l_{1}+u_{1} &=&1 \\
l_{2}+u_{2} &=&1
\end{eqnarray*}%
i.e.%
\begin{eqnarray*}
a+b_{1}+b_{2}+d+e_{1}+e_{2} &\equiv &1 \\
l_{1}+u_{1} &\equiv &1 \\
l_{2}+u_{2} &\equiv &1
\end{eqnarray*}

Since $\alpha \left( 1_{H},0,0,1,1\right) =1+e_{2},\alpha \left(
1_{H},1,0,0,1\right) =a+b_{1},$

$\alpha \left( 1_{H},0,1,1,0\right) =e_{2}+a+b_{1}+b_{2}+1,\alpha \left(
1_{H},1,1,0,0\right) =1+b_{2}$

we get%
\begin{gather*}
\sum_{\substack{ a,b_{1},b_{2},d,e_{1},e_{2}=0  \\ %
a+b_{1}+b_{2}+d+e_{1}+e_{2}\equiv 1}}^{1}\sum_{l_{1}=0}^{b_{1}}%
\sum_{l_{2}=0}^{b_{2}}\sum_{\substack{ u_{1}=0  \\ l_{1}+u_{1}\equiv 1}}%
^{e_{1}}\sum_{\substack{ u_{2}=0  \\ l_{2}+u_{2}\equiv 1}}%
^{e_{2}}(-1)^{\alpha \left( 1_{H},l_{1},l_{2},u_{1},u_{2}\right) } \\
B(x_{1}\otimes
1_{H};G^{a}X_{1}^{b_{1}}X_{2}^{b_{2}},g^{d}x_{1}^{e_{1}}x_{2}^{e_{2}})G^{a}X_{1}^{b_{1}-l_{1}}X_{2}^{b_{2}-l_{2}}\otimes g^{d}x_{1}^{\left( e_{1}-u_{1}\right) }x_{2}^{\left( e_{2}-u_{2}\right) }=
\\
\sum_{\substack{ a,b_{1},b_{2},d=0  \\ a+b_{1}+b_{2}+d\equiv 1}}^{1}\left(
-1\right) ^{1+e_{2}}B(x_{1}\otimes
1_{H};G^{a}X_{1}^{b_{1}}X_{2}^{b_{2}},g^{d}x_{1}x_{2})G^{a}X_{1}^{b_{1}}X_{2}^{b_{2}}\otimes g^{d}+
\\
\sum_{\substack{ a,b_{2},d,e_{1}=0  \\ a+b_{2}+d+e_{1}\equiv 1}}%
^{1}(-1)^{a+b_{1}}B(x_{1}\otimes
1_{H};G^{a}X_{1}X_{2}^{b_{2}},g^{d}x_{1}^{e_{1}}x_{2})G^{a}X_{2}^{b_{2}}%
\otimes g^{d}x_{1}^{e_{1}}+ \\
\sum_{\substack{ a,b_{1},d,e_{2}=0  \\ a+b_{1}+d+e_{2}\equiv 1}}%
^{1}(-1)^{e_{2}+a+b_{1}+b_{2}+1}B(x_{1}\otimes
1_{H};G^{a}X_{1}^{b_{1}}X_{2},g^{d}x_{1}x_{2}^{e_{2}})G^{a}X_{1}^{b_{1}}%
\otimes g^{d}x_{2}^{e_{2}}+ \\
\sum_{\substack{ a,d,e_{1},e_{2}=0  \\ a+d+e_{1}+e_{2}\equiv 1}}%
^{1}(-1)^{1+b_{2}}B(x_{1}\otimes
1_{H};G^{a}X_{1}X_{2},g^{d}x_{1}^{e_{1}}x_{2}^{e_{2}})G^{a}\otimes
g^{d}x_{1}^{e_{1}}x_{2}^{e_{2}}
\end{gather*}%
Thus we obtain

\begin{eqnarray*}
&&\sum_{\substack{ a,b_{1},b_{2},d=0  \\ a+b_{1}+b_{2}+d\equiv 1}}%
^{1}B(x_{1}\otimes
1_{H};G^{a}X_{1}^{b_{1}}X_{2}^{b_{2}},g^{d}x_{1}x_{2})G^{a}X_{1}^{b_{1}}X_{2}^{b_{2}}\otimes g^{d}+
\\
&&\sum_{\substack{ a,b_{2},d,e_{1}=0  \\ a+b_{2}+d+e_{1}\equiv 1}}%
^{1}(-1)^{a+1}B(x_{1}\otimes
1_{H};G^{a}X_{1}X_{2}^{b_{2}},g^{d}x_{1}^{e_{1}}x_{2})G^{a}X_{2}^{b_{2}}%
\otimes g^{d}x_{1}^{e_{1}}+ \\
&&\sum_{\substack{ a,b_{1},d,e_{2}=0  \\ a+b_{1}+d+e_{2}\equiv 1}}%
^{1}(-1)^{e_{2}+a+b_{1}}B(x_{1}\otimes
1_{H};G^{a}X_{1}^{b_{1}}X_{2},g^{d}x_{1}x_{2}^{e_{2}})G^{a}X_{1}^{b_{1}}%
\otimes g^{d}x_{2}^{e_{2}}+ \\
&&\sum_{\substack{ a,d,e_{1},e_{2}=0  \\ a+d+e_{1}+e_{2}\equiv 1}}%
^{1}B(x_{1}\otimes
1_{H};G^{a}X_{1}X_{2},g^{d}x_{1}^{e_{1}}x_{2}^{e_{2}})G^{a}\otimes
g^{d}x_{1}^{e_{1}}x_{2}^{e_{2}}
\end{eqnarray*}

Here \ we consider all cases.

\subsubsection{$G^{a}\otimes g^{d}$}

\begin{eqnarray*}
&&\sum_{\substack{ a,b_{1},b_{2},d=0  \\ a+b_{1}+b_{2}+d\equiv 1}}%
^{1}B(x_{1}\otimes
1_{H};G^{a}X_{1}^{b_{1}}X_{2}^{b_{2}},g^{d}x_{1}x_{2})G^{a}X_{1}^{b_{1}}X_{2}^{b_{2}}\otimes g^{d}+
\\
&&\sum_{\substack{ a,b_{2},d,e_{1}=0  \\ a+b_{2}+d+e_{1}\equiv 1}}%
^{1}(-1)^{a+1}B(x_{1}\otimes
1_{H};G^{a}X_{1}X_{2}^{b_{2}},g^{d}x_{1}^{e_{1}}x_{2})G^{a}X_{2}^{b_{2}}%
\otimes g^{d}x_{1}^{e_{1}}+ \\
&&\sum_{\substack{ a,b_{1},d,e_{2}=0  \\ a+b_{1}+d+e_{2}\equiv 1}}%
^{1}(-1)^{e_{2}+a+b_{1}}B(x_{1}\otimes
1_{H};G^{a}X_{1}^{b_{1}}X_{2},g^{d}x_{1}x_{2}^{e_{2}})G^{a}X_{1}^{b_{1}}%
\otimes g^{d}x_{2}^{e_{2}}+ \\
&&\sum_{\substack{ a,d,e_{1},e_{2}=0  \\ a+d+e_{1}+e_{2}\equiv 1}}%
^{1}B(x_{1}\otimes
1_{H};G^{a}X_{1}X_{2},g^{d}x_{1}^{e_{1}}x_{2}^{e_{2}})G^{a}\otimes
g^{d}x_{1}^{e_{1}}x_{2}^{e_{2}}
\end{eqnarray*}%
\begin{equation*}
\sum_{\substack{ a,d=0  \\ a+d\equiv 1}}^{1}\left[
\begin{array}{c}
B(x_{1}\otimes 1_{H};G^{a},g^{d}x_{1}x_{2})+(-1)^{a+1}B(x_{1}\otimes
1_{H};G^{a}X_{1},g^{d}x_{2})+ \\
+(-1)^{a}B(x_{1}\otimes 1_{H};G^{a}X_{2},g^{d}x_{1})+B(x_{1}\otimes
1_{H};G^{a}X_{1}X_{2},g^{d})%
\end{array}%
\right] G^{a}\otimes g^{d}+
\end{equation*}%
and we get%
\begin{equation}
\begin{array}{c}
+B(x_{1}\otimes 1_{H};1_{A},gx_{1}x_{2})-B(x_{1}\otimes 1_{H};X_{1},gx_{2})+
\\
+B(x_{1}\otimes 1_{H};X_{2},gx_{1})+B(x_{1}\otimes 1_{H};X_{1}X_{2},g)=0%
\end{array}
\label{x1ot1,sixty}
\end{equation}%
\begin{equation}
\begin{array}{c}
B(x_{1}\otimes 1_{H};G,x_{1}x_{2})+B(x_{1}\otimes 1_{H};GX_{1},x_{2})+ \\
-B(x_{1}\otimes 1_{H};GX_{2},x_{1})+B(x_{1}\otimes 1_{H};GX_{1}X_{2},1_{H})=0%
\end{array}
\label{x1ot1,sixtyone}
\end{equation}

\subsubsection{$G^{a}\otimes g^{d}x_{2}$}

\begin{equation*}
\sum_{\substack{ a,d=0  \\ a+d\equiv 0}}^{1}\left[ (-1)^{1+a}B(x_{1}\otimes
1_{H};G^{a}X_{2},g^{d}x_{1}x_{2})+B(x_{1}\otimes
1_{H};G^{a}X_{1}X_{2},g^{d}x_{2})\right] G^{a}\otimes g^{d}x_{2}=0
\end{equation*}%
and we get%
\begin{eqnarray}
-B(x_{1}\otimes 1_{H};X_{2},x_{1}x_{2})+B(x_{1}\otimes
1_{H};X_{1}X_{2},x_{2}) &=&0  \label{x1ot1, sixtytwo} \\
B(x_{1}\otimes 1_{H};GX_{2},gx_{1}x_{2})+B(x_{1}\otimes
1_{H};GX_{1}X_{2},gx_{2}) &=&0  \label{x1ot1, sixtythree}
\end{eqnarray}

\subsubsection{$G^{a}\otimes g^{d}x_{1}$}

\begin{equation*}
\sum_{\substack{ a,d=0  \\ a+d\equiv 0}}^{1}\left[ (-1)^{a+1}B(x_{1}\otimes
1_{H};G^{a}X_{1},g^{d}x_{1}x_{2})+B(x_{1}\otimes
1_{H};G^{a}X_{1}X_{2},g^{d}x_{1})\right] G^{a}\otimes g^{d}x_{1}
\end{equation*}%
and we get%
\begin{eqnarray}
-B(x_{1}\otimes 1_{H};X_{1},x_{1}x_{2})+B(x_{1}\otimes
1_{H};X_{1}X_{2},x_{1}) &=&0  \label{x1ot1, sixtyfour} \\
(B(x_{1}\otimes 1_{H};GX_{1},gx_{1}x_{2})+B(x_{1}\otimes
1_{H};GX_{1}X_{2},gx_{1}) &=&0  \label{x1ot1, sixtyfive}
\end{eqnarray}

\subsubsection{$G^{a}X_{2}\otimes g^{d}$}

\begin{equation*}
\sum_{\substack{ a,d=0  \\ a+d\equiv 0}}^{1}\left[ B(x_{1}\otimes
1_{H};G^{a}X_{2},g^{d}x_{1}x_{2})+(-1)^{a+1}B(x_{1}\otimes
1_{H};G^{a}X_{1}X_{2},g^{d}x_{2})\right] G^{a}X_{2}\otimes g^{d}=0
\end{equation*}%
and we get%
\begin{eqnarray*}
B(x_{1}\otimes 1_{H};X_{2},x_{1}x_{2})-B(x_{1}\otimes
1_{H};X_{1}X_{2},x_{2}) &=&0 \\
B(x_{1}\otimes 1_{H};GX_{2},gx_{1}x_{2})+B(x_{1}\otimes
1_{H};GX_{1}X_{2},gx_{2}) &=&0
\end{eqnarray*}%
which we already obtained.

\subsubsection{$G^{a}X_{1}\otimes g^{d}$}

\begin{equation*}
\sum_{\substack{ a,d=0  \\ a+d\equiv 0}}^{1}\left[ B(x_{1}\otimes
1_{H};G^{a}X_{1},g^{d}x_{1}x_{2})+(-1)^{a+1}B(x_{1}\otimes
1_{H};G^{a}X_{1}X_{2},g^{d}x_{1})\right] G^{a}X_{1}\otimes g^{d}=0
\end{equation*}%
and we get%
\begin{eqnarray*}
B(x_{1}\otimes 1_{H};X_{1},x_{1}x_{2})-B(x_{1}\otimes
1_{H};X_{1}X_{2},x_{1}) &=&0 \\
B(x_{1}\otimes 1_{H};GX_{1},gx_{1}x_{2})+B(x_{1}\otimes
1_{H};GX_{1}X_{2},gx_{1}) &=&0
\end{eqnarray*}%
which we already obtained.

\subsubsection{$G^{a}\otimes g^{d}x_{1}x_{2}$}

\begin{equation*}
\sum_{\substack{ a,d=0  \\ a+d\equiv 1}}^{1}B(x_{1}\otimes
1_{H};G^{a}X_{1}X_{2},g^{d}x_{1}x_{2})G^{a}\otimes g^{d}x_{1}x_{2}=0
\end{equation*}%
and we get
\begin{eqnarray}
B(x_{1}\otimes 1_{H};X_{1}X_{2},gx_{1}x_{2}) &=&0  \label{x1ot1, sixtysix} \\
B(x_{1}\otimes 1_{H};GX_{1}X_{2},x_{1}x_{2}) &=&0  \label{x1ot1, sixtyseven}
\end{eqnarray}

\subsubsection{$G^{a}X_{2}\otimes g^{d}x_{2}$}

There no summands like this.

\subsubsection{$G^{a}X_{1}\otimes g^{d}x_{2}$}

\begin{equation*}
\sum_{\substack{ a,d=0  \\ a+d\equiv 1}}^{1}(-1)^{a}B(x_{1}\otimes
1_{H};G^{a}X_{1}X_{2},g^{d}x_{1}x_{2})G^{a}X_{1}\otimes g^{d}x_{2}=0
\end{equation*}%
and we get%
\begin{eqnarray*}
B(x_{1}\otimes 1_{H};X_{1}X_{2},gx_{1}x_{2}) &=&0 \\
B(x_{1}\otimes 1_{H};GX_{1}X_{2},x_{1}x_{2}) &=&0
\end{eqnarray*}%
which we already obtained

\subsubsection{$G^{a}X_{2}\otimes g^{d}x_{1}$}

\begin{equation*}
\sum_{\substack{ a,d=0  \\ a+d\equiv 1}}^{1}(-1)^{a+1}B(x_{1}\otimes
1_{H};G^{a}X_{1}X_{2},g^{d}x_{1}x_{2})G^{a}X_{2}\otimes g^{d}x_{1}=0
\end{equation*}%
and we get%
\begin{eqnarray*}
-B(x_{1}\otimes 1_{H};X_{1}X_{2},gx_{1}x_{2}) &=&0 \\
B(x_{1}\otimes 1_{H};GX_{1}X_{2},x_{1}x_{2}) &=&0
\end{eqnarray*}%
which we already obtained

\subsubsection{$G^{a}X_{1}\otimes g^{d}x_{1}$}

There are no terms like this.

\subsubsection{$G^{a}X_{1}X_{2}\otimes g^{d}$}

\begin{equation*}
\sum_{\substack{ a,d=0  \\ a+d\equiv 1}}^{1}B(x_{1}\otimes
1_{H};G^{a}X_{1}X_{2},g^{d}x_{1}x_{2})G^{a}X_{1}X_{2}\otimes g^{d}=0
\end{equation*}%
and we get%
\begin{eqnarray*}
B(x_{1}\otimes 1_{H};X_{1}X_{2},gx_{1}x_{2}) &=&0 \\
B(x_{1}\otimes 1_{H};GX_{1}X_{2},x_{1}x_{2}) &=&0
\end{eqnarray*}%
which we already obtained.

\subsubsection{$G^{a}X_{2}\otimes g^{d}x_{1}x_{2}$}

There is no term like this.

\subsubsection{$G^{a}X_{1}\otimes g^{d}x_{1}x_{2}$}

There is no term like this.

\subsubsection{$G^{a}X_{1}X_{2}\otimes g^{d}x_{2}$}

There is no term like this.

\subsubsection{$G^{a}X_{1}X_{2}\otimes g^{d}x_{1}$}

There is no term like this.

\subsubsection{$G^{a}X_{1}X_{2}\otimes g^{d}x_{1}x_{2}$}

There is no term like this.

By using all the equalities we got above from $\left( \ref{x1ot1 first}%
\right) $ to $\left( \ref{x1ot1, sixtyseven}\right) $ we obtain the
following form of $B\left( x_{1}\otimes 1_{H}\right) .$

\subsection{The final form of the element $B\left( x_{1}\otimes 1_{H}\right)
$}

Here we use also the final form of $B\left( g\otimes 1_{H}\right) .$
\begin{eqnarray}
&&B\left( x_{1}\otimes 1_{H}\right)  \label{x1} \\
&=&B\left( x_{1}\otimes 1_{H};1_{A},1_{H}\right) 1_{A}\otimes 1_{H}+  \notag
\\
&&B\left( x_{1}\otimes 1_{H};1_{A},x_{1}x_{2}\right) 1_{A}\otimes x_{1}x_{2}+
\notag \\
&&B\left( x_{1}\otimes 1_{H};1_{A},gx_{1}\right) 1_{A}\otimes gx_{1}+  \notag
\\
&&B\left( x_{1}\otimes 1_{H};1_{A},gx_{2}\right) 1_{A}\otimes gx_{2}+  \notag
\\
&&B\left( x_{1}\otimes 1_{H};G,g\right) G\otimes g+  \notag \\
&&B\left( x_{1}\otimes 1_{H};G,x_{1}\right) G\otimes x_{1}+  \notag \\
&&B\left( x_{1}\otimes 1_{H};G,x_{2}\right) G\otimes x_{2}+  \notag \\
&&B\left( x_{1}\otimes 1_{H};G,gx_{1}x_{2}\right) G\otimes gx_{1}x_{2}+
\notag \\
&&\left[ -B(g\otimes 1_{H};1_{A},g)-B(x_{1}\otimes 1_{H};1_{A},gx_{1})\right]
X_{1}\otimes g+  \notag \\
&&-B(g\otimes 1_{H};1_{A},x_{1})X_{1}\otimes x_{1}+  \notag \\
&&\left[ -B(g\otimes 1_{H};1_{A},x_{2})+B(x_{1}\otimes
1_{H};1_{A},x_{1}x_{2})\right] X_{1}\otimes x_{2}+  \notag \\
&&-B(g\otimes 1_{H};1_{A},gx_{1}x_{2})X_{1}\otimes gx_{1}x_{2}+  \notag \\
&&-B(x_{1}\otimes 1_{H};1_{A},gx_{2})X_{2}\otimes g+  \notag \\
&&-B(x_{1}\otimes 1_{H};1_{A},x_{1}x_{2})X_{2}\otimes x_{1}+  \notag \\
&&\left[ -B(g\otimes 1_{H};X_{2},1_{H})+B(x_{1}\otimes
1_{H};1_{A},x_{1}x_{2})\right] X_{1}X_{2}\otimes 1_{H}+  \notag \\
&&-B(g\otimes 1_{H};X_{1}X_{2},g)X_{1}X_{2}\otimes gx_{1}+  \notag \\
&&\left[ B(g\otimes 1_{H};G,1_{H})+B(x_{1}\otimes 1_{H};G,x_{1})\right]
GX_{1}\otimes 1_{H}+  \notag \\
&&B(g\otimes 1_{H};GX_{1}X_{2},1_{H})GX_{1}\otimes x_{1}x_{2}+  \notag \\
&&-B(g\otimes 1_{H};GX_{1},g)GX_{1}\otimes gx_{1}+  \notag \\
&&\left[ -B(g\otimes 1_{H};GX_{2},g)-B(x_{1}\otimes 1_{H};G,gx_{1}x_{2})%
\right] GX_{1}\otimes gx_{2}+  \notag \\
&&B(x_{1}\otimes 1_{H};G,x_{2})GX_{2}\otimes 1_{H}+  \notag \\
&&B(x_{1}\otimes 1_{H};G,gx_{1}x_{2})GX_{2}\otimes gx_{1}+  \notag \\
&&\left[ -B(g\otimes 1_{H};G,gx_{2})+B(x_{1}\otimes 1_{H};G,gx_{1}x_{2})%
\right] GX_{1}X_{2}\otimes g+  \notag \\
&&-B(g\otimes 1_{H};G,x_{1}x_{2})GX_{1}X_{2}\otimes x_{1}  \notag
\end{eqnarray}

\section{$B(x_{2}\otimes 1_{H})$}

We write the Casimir formula $\left( \ref{MAIN FORMULA 1}\right) $ for $%
B(x_{2}\otimes 1_{H})$
\begin{eqnarray*}
&&\sum_{w_{2}=0}^{1}\sum_{a,b_{1},b_{2},d,e_{1},e_{2}=0}^{1}%
\sum_{l_{1}=0}^{b_{1}}\sum_{l_{2}=0}^{b_{2}}\sum_{u_{1}=0}^{e_{1}}%
\sum_{u_{2}=0}^{e_{2}}\left( -1\right) ^{\alpha \left(
x_{2}^{1-w_{2}};l_{1},l_{2},u_{1},u_{2}\right) } \\
&&B(g^{1+w_{2}}x_{2}^{w_{2}}\otimes g^{\mu }x_{1}^{\nu _{1}}x_{2}^{\nu
_{2}};G^{a}X_{1}^{b_{1}}X_{2}^{b_{2}},g^{d}x_{1}^{e_{1}}x_{2}^{e_{2}}) \\
&&G^{a}X_{1}^{b_{1}-l_{1}}X_{2}^{b_{2}-l_{2}}\otimes
g^{d}x_{1}^{e_{1}-u_{1}}x_{2}^{e_{2}-u_{2}}\otimes \\
&&g^{a+b_{1}+b_{2}+l_{1}+l_{2}+d+e_{1}+e_{2}+u_{1}+u_{2}}x_{1}^{l_{1}+u_{1}}x_{2}^{l_{2}+u_{2}+1-w_{2}}
\\
&=&B^{A}(x_{2}\otimes 1_{H})\otimes B^{H}(x_{2}\otimes 1_{H})\otimes 1_{H}
\end{eqnarray*}%
and we obtain

\begin{eqnarray*}
&&\sum_{a,b_{1},b_{2},d,e_{1},e_{2}=0}^{1}\sum_{l_{1}=0}^{b_{1}}%
\sum_{l_{2}=0}^{b_{2}}\sum_{u_{1}=0}^{e_{1}}\sum_{u_{2}=0}^{e_{2}}\left(
-1\right) ^{\alpha \left( x_{2};l_{1},l_{2},u_{1},u_{2}\right) } \\
&&B(g\otimes
1_{H};G^{a}X_{1}^{b_{1}}X_{2}^{b_{2}},g^{d}x_{1}^{e_{1}}x_{2}^{e_{2}})G^{a}X_{1}^{b_{1}-l_{1}}X_{2}^{b_{2}-l_{2}}\otimes g^{d}x_{1}^{e_{1}-u_{1}}x_{2}^{e_{2}-u_{2}}
\\
&&\otimes
g^{a+b_{1}+b_{2}+l_{1}+l_{2}+d+e_{1}+e_{2}+u_{1}+u_{2}}x_{1}^{l_{1}+u_{1}}x_{2}^{l_{2}+u_{2}+1}
\\
&&+\sum_{a,b_{1},b_{2},d,e_{1},e_{2}=0}^{1}\sum_{l_{1}=0}^{b_{1}}%
\sum_{l_{2}=0}^{b_{2}}\sum_{u_{1}=0}^{e_{1}}\sum_{u_{2}=0}^{e_{2}}\left(
-1\right) ^{\alpha \left( 1_{H};l_{1},l_{2},u_{1},u_{2}\right) } \\
&&B(x_{2}\otimes
1_{H};G^{a}X_{1}^{b_{1}}X_{2}^{b_{2}},g^{d}x_{1}^{e_{1}}x_{2}^{e_{2}})G^{a}X_{1}^{b_{1}-l_{1}}X_{2}^{b_{2}-l_{2}}\otimes g^{d}x_{1}^{e_{1}-u_{1}}x_{2}^{e_{2}-u_{2}}
\\
&&\otimes
g^{a+b_{1}+b_{2}+l_{1}+l_{2}+d+e_{1}+e_{2}+u_{1}+u_{2}}x_{1}^{l_{1}+u_{1}}x_{2}^{l_{2}+u_{2}}
\\
&=&B^{A}(x_{2}\otimes 1_{H})\otimes B^{H}(x_{2}\otimes 1_{H})\otimes 1_{H}
\end{eqnarray*}

$\newline
$

We concentrate on the different occurrences of the third part of the tensor.

\subsection{Case $1_{H}$}

Using the same argument in case $B\left( x_{1}\otimes 1_{H}\right) ,$ we
obtain

\begin{equation}
B(x_{2}\otimes
1_{H};G^{a}X_{1}^{b_{1}}X_{2}^{b_{2}},g^{d}x_{1}^{e_{1}}x_{2}^{e_{2}})=0%
\text{ whenever }a+b_{1}+b_{2}+d+e_{1}+e_{2}\text{ }\equiv 1
\label{xot2 first}
\end{equation}

\subsection{Case $g$}

From the first summand of the left side we can not obtain $g$ . Therefore we
get

\begin{eqnarray*}
g^{a+b_{1}+b_{2}+l_{1}+l_{2}+d+e_{1}+e_{2}+u_{1}+u_{2}}x_{1}^{l_{1}+u_{1}}x_{2}^{l_{2}+u_{2}} &=&g\Rightarrow
\\
u_{1} &=&l_{1}=l_{2}=u_{2}=0 \\
a+b_{1}+b_{2}+l_{1}+l_{2}+d+e_{1}+e_{2}+u_{1}+u_{2} &\equiv &1\text{ so
that, } \\
B(x_{2}\otimes
1_{H};G^{a}X_{1}^{b_{1}}X_{2}^{b_{2}},g^{d}x_{1}^{e_{1}}x_{2}^{e_{2}}) &=&0%
\text{ whenever }a+b_{1}+b_{2}+d+e_{1}+e_{2}\equiv 1. \\
&&\text{This is }\left( \ref{xot2 first}\right)
\end{eqnarray*}

\subsection{Case $x_{2}$}

Here we will have two summands with $x_{2}$ in the third position.

From the first summand of the left side we get\qquad
\begin{eqnarray*}
g^{a+b_{1}+b_{2}+d+e_{1}+e_{2}+l_{1}+l_{2}+u_{1}+u_{2}}x_{1}^{u_{1}+l_{1}}x_{2}^{l_{2}+u_{2}+1} &=&x_{2}
\\
a+b_{1}+b_{2}+d+e_{1}+e_{2}+l_{1}+l_{2}+u_{1}+u_{2} &\equiv &0 \\
u_{1}+l_{1} &=&0\text{ so that }u_{1}=0=l_{1} \\
l_{2}+u_{2} &=&0\text{ so that }l_{2}=0=u_{2}\text{ and} \\
a+b_{1}+b_{2}+d+e_{1}+e_{2} &\equiv &0\text{ but we already know} \\
B(g\otimes
1_{H};G^{a}X_{1}^{b_{1}}X_{2}^{b_{2}},g^{d}x_{1}^{e_{1}}x_{2}^{e_{2}}) &=&0%
\text{ whenever }a+b_{1}+b_{2}+d+e_{1}+e_{2}\equiv 0\text{ }\left( \ref{got1
first}\right)
\end{eqnarray*}%
From the second summand of the left side we get%
\begin{eqnarray*}
g^{a+b_{1}+b_{2}+l_{1}+l_{2}+d+e_{1}+e_{2}+u_{1}+u_{2}}x_{1}^{l_{1}+u_{1}}x_{2}^{l_{2}+u_{2}} &=&x_{2}\Rightarrow
\\
l_{1}+u_{1} &=&0 \\
l_{2}+u_{2} &=&1 \\
a+b_{1}+b_{2}+l_{1}+l_{2}+d+e_{1}+e_{2}+u_{1}+u_{2} &\equiv &0 \\
a+b_{1}+b_{2}+d+e_{1}+e_{2} &\equiv &1\text{ We already know, in view of }%
\left( \ref{xot2 first}\right) \text{ that} \\
B(x_{2}\otimes
1_{H};G^{a}X_{1}^{b_{1}}X_{2}^{b_{2}},g^{d}x_{1}^{e_{1}}x_{2}^{e_{2}}) &=&0%
\text{ whenever }a+b_{1}+b_{2}+d+e_{1}+e_{2}\equiv 1.
\end{eqnarray*}

\subsection{Case $x_{1}$}

From the first summand of the left side we can not obtain $x_{1}.$ From the
second summand we get%
\begin{eqnarray*}
g^{a+b_{1}+b_{2}+l_{1}+l_{2}+d+e_{1}+e_{2}+u_{1}+u_{2}}x_{1}^{l_{1}+u_{1}}x_{2}^{l_{2}+u_{2}} &=&x_{1}\Rightarrow
\\
l_{1}+u_{1} &=&0 \\
l_{2} &=&u_{2}=0 \\
a+b_{1}+b_{2}+l_{1}+l_{2}+d+e_{1}+e_{2}+u_{1}+u_{2} &\equiv &0\text{ so
that, } \\
a+b_{1}+b_{2}+d+e_{1}+e_{2} &\equiv &1 \\
B(x_{2}\otimes
1_{H};G^{a}X_{1}^{b_{1}}X_{2}^{b_{2}},g^{d}x_{1}^{e_{1}}x_{2}^{e_{2}}) &=&0%
\text{ whenever }a+b_{1}+b_{2}+d+e_{1}+e_{2}\equiv 1. \\
&&\text{This is }\left( \ref{xot2 first}\right)
\end{eqnarray*}

\subsection{Case $x_{1}x_{2}$}

\begin{eqnarray*}
&&\sum_{a,b_{1},b_{2},d,e_{1},e_{2}=0}^{1}\sum_{l_{1}=0}^{b_{1}}%
\sum_{l_{2}=0}^{b_{2}}\sum_{u_{1}=0}^{e_{1}}\sum_{u_{2}=0}^{e_{2}}\left(
-1\right) ^{\alpha \left( x_{2};l_{1},l_{2},u_{1},u_{2}\right) } \\
&&B(g\otimes
1_{H};G^{a}X_{1}^{b_{1}}X_{2}^{b_{2}},g^{d}x_{1}^{e_{1}}x_{2}^{e_{2}}) \\
&&G^{a}X_{1}^{b_{1}-l_{1}}X_{2}^{b_{2}-l_{2}}\otimes
g^{d}x_{1}^{e_{1}-u_{1}}x_{2}^{e_{2}-u_{2}}\otimes
g^{a+b_{1}+b_{2}+l_{1}+l_{2}+d+e_{1}+e_{2}+u_{1}+u_{2}}x_{1}^{l_{1}+u_{1}}x_{2}^{l_{2}+u_{2}+1}
\\
&&+\sum_{a,b_{1},b_{2},d,e_{1},e_{2}=0}^{1}\sum_{l_{1}=0}^{b_{1}}%
\sum_{l_{2}=0}^{b_{2}}\sum_{u_{1}=0}^{e_{1}}\sum_{u_{2}=0}^{e_{2}}\left(
-1\right) ^{\alpha \left( 1_{H};l_{1},l_{2},u_{1},u_{2}\right) } \\
&&B(x_{2}\otimes
1_{H};G^{a}X_{1}^{b_{1}}X_{2}^{b_{2}},g^{d}x_{1}^{e_{1}}x_{2}^{e_{2}}) \\
&&G^{a}X_{1}^{b_{1}-l_{1}}X_{2}^{b_{2}-l_{2}}\otimes
g^{d}x_{1}^{e_{1}-u_{1}}x_{2}^{e_{2}-u_{2}}\otimes
g^{a+b_{1}+b_{2}+l_{1}+l_{2}+d+e_{1}+e_{2}+u_{1}+u_{2}}x_{1}^{l_{1}+u_{1}}x_{2}^{l_{2}+u_{2}}
\\
&=&B^{A}(x_{2}\otimes 1_{H})\otimes B^{H}(x_{2}\otimes 1_{H})\otimes 1_{H}
\end{eqnarray*}

From the first summand we get%
\begin{eqnarray*}
g^{a+b_{1}+b_{2}+d+e_{1}+e_{2}+l_{1}+l_{2}+u_{1}+u_{2}}x_{1}^{u_{1}+l_{1}}x_{2}^{l_{2}+u_{2}+1} &=&x_{1}x_{2}
\\
a+b_{1}+b_{2}+d+e_{1}+e_{2}+l_{1}+l_{2}+u_{1}+u_{2} &\equiv &0 \\
u_{1}+l_{1} &=&1\text{ } \\
l_{2}+u_{2} &=&0\text{ so that }l_{2}=0=u_{2}\text{ and} \\
a+b_{1}+b_{2}+d+e_{1}+e_{2} &\equiv &1
\end{eqnarray*}%
and we get%
\begin{eqnarray*}
&&\sum_{\substack{ a,b_{1},b_{2},d,e_{1},e_{2}=0  \\ %
a+b_{1}+b_{2}+d+e_{1}+e_{2}\equiv 1}}^{1}\sum_{l_{1}=0}^{b_{1}}%
\sum_{u_{1}=0}^{e_{1}}\left( -1\right) ^{\alpha \left(
x_{2};l_{1},0,u_{1},0\right) }B(g\otimes
1_{H};G^{a}X_{1}^{b_{1}}X_{2}^{b_{2}},g^{d}x_{1}^{e_{1}}x_{2}^{e_{2}}) \\
&&G^{a}X_{1}^{b_{1}-l_{1}}X_{2}^{b_{2}}\otimes
g^{d}x_{1}^{e_{1}-u_{1}}x_{2}^{e_{2}} \\
&=&\sum_{\substack{ a,b_{1},b_{2},d,e_{1},e_{2}=0  \\ a+b_{1}+b_{2}+d+e_{2}%
\equiv 0}}^{1}\left( -1\right) ^{\alpha \left( x_{2};0,0,1,0\right)
}B(g\otimes
1_{H};G^{a}X_{1}^{b_{1}}X_{2}^{b_{2}},g^{d}x_{1}x_{2}^{e_{2}})G^{a}X_{1}^{b_{1}}X_{2}^{b_{2}}\otimes g^{d}x_{2}^{e_{2}}+
\\
&&\sum_{\substack{ a,b_{1},b_{2},d,e_{1},e_{2}=0  \\ a+b_{2}+d+e_{1}+e_{2}%
\equiv 0}}^{1}\left( -1\right) ^{\alpha \left( x_{2};1,0,0,0\right)
}B(g\otimes
1_{H};G^{a}X_{1}X_{2}^{b_{2}},g^{d}x_{1}^{e_{1}}x_{2}^{e_{2}})G^{a}X_{2}^{b_{2}}\otimes g^{d}x_{1}^{e_{1}}x_{2}^{e_{2}}
\end{eqnarray*}

From the second summand we get%
\begin{eqnarray*}
g^{a+b_{1}+b_{2}+l_{1}+l_{2}+d+e_{1}+e_{2}+u_{1}+u_{2}}x_{1}^{l_{1}+u_{1}}x_{2}^{l_{2}+u_{2}} &=&x_{1}x_{2}\Rightarrow
\\
l_{1}+u_{1} &=&1 \\
l_{2}+u_{2} &=&1 \\
a+b_{1}+b_{2}+l_{1}+l_{2}+d+e_{1}+e_{2}+u_{1}+u_{2} &\equiv &0\text{ so
that, } \\
a+b_{1}+b_{2}+d+e_{1}+e_{2} &\equiv &0
\end{eqnarray*}%
and we get

\begin{eqnarray*}
&&\sum_{\substack{ a,b_{1},b_{2},d,e_{1},e_{2}=0  \\ %
a+b_{1}+b_{2}+d+e_{1}+e_{2}\equiv 0}}^{1}\sum_{l_{1}=0}^{b_{1}}%
\sum_{l_{2}=0}^{b_{2}}\sum_{\substack{ u_{1}=0  \\ l_{1}+u_{1}=1}}%
^{e_{1}}\sum _{\substack{ u_{2}=0  \\ l_{2}+u_{2}=1}}^{e_{2}}\left(
-1\right) ^{\alpha \left( 1_{H};l_{1},l_{2},u_{1},u_{2}\right) } \\
&&B(x_{2}\otimes
1_{H};G^{a}X_{1}^{b_{1}}X_{2}^{b_{2}},g^{d}x_{1}^{e_{1}}x_{2}^{e_{2}})G^{a}X_{1}^{b_{1}-l_{1}}X_{2}^{b_{2}-l_{2}}\otimes g^{d}x_{1}^{e_{1}-u_{1}}x_{2}^{e_{2}-u_{2}}
\\
&=&\sum_{\substack{ a,b_{1},b_{2},d=0  \\ a+b_{1}+b_{2}+d\equiv 0}}%
^{1}\left( -1\right) ^{\alpha \left( 1_{H};0,0,1,1\right) }B(x_{2}\otimes
1_{H};G^{a}X_{1}^{b_{1}}X_{2}^{b_{2}},g^{d}x_{1}x_{2})G^{a}X_{1}^{b_{1}}X_{2}^{b_{2}}\otimes g^{d}+
\\
&&\sum_{\substack{ a,b_{1},d,e_{2}=0  \\ a+b_{1}+d+e_{2}\equiv 0}}^{1}\left(
-1\right) ^{\alpha \left( 1_{H};0,1,1,0\right) }B(x_{2}\otimes
1_{H};G^{a}X_{1}^{b_{1}}X_{2},g^{d}x_{1}x_{2}^{e_{2}})G^{a}X_{1}^{b_{1}}%
\otimes g^{d}x_{2}^{e_{2}}+ \\
&&\sum_{\substack{ a,b_{2},d,e_{1}=0  \\ a+b_{2}+d+e_{1}\equiv 0}}^{1}\left(
-1\right) ^{\alpha \left( 1_{H};1,0,0,1\right) }B(x_{2}\otimes
1_{H};G^{a}X_{1}X_{2}^{b_{2}},g^{d}x_{1}^{e_{1}}x_{2})G^{a}X_{2}^{b_{2}}%
\otimes g^{d}x_{1}^{e_{1}}+ \\
&&\sum_{\substack{ a,d,e_{1},e_{2}=0  \\ a+d+e_{1}+e_{2}\equiv 0}}^{1}\left(
-1\right) ^{\alpha \left( 1_{H};1,1,0,0\right) }B(x_{2}\otimes
1_{H};G^{a}X_{1}X_{2},g^{d}x_{1}^{e_{1}}x_{2}^{e_{2}})G^{a}\otimes
g^{d}x_{1}^{e_{1}}x_{2}^{e_{2}}
\end{eqnarray*}

Since $\alpha \left( x_{2};0,0,1,0\right) =e_{2},\alpha \left(
x_{2};1,0,0,0\right) =a+b_{1},\alpha \left( 1_{H};0,0,1,1\right) =1+e_{2},$

$\alpha \left( 1_{H};0,1,1,0\right) =e_{2}+a+b_{1}+b_{2}+1,$

$\alpha \left( 1_{H};1,0,0,1\right) =a+b_{1},\alpha \left(
1_{H};1,1,0,0\right) =1+b_{2}$ we obtain

\begin{eqnarray*}
&&\sum_{\substack{ a,b_{1},b_{2},d,e_{1},e_{2}=0  \\ a+b_{1}+b_{2}+d+e_{2}%
\equiv 0}}^{1}\left( -1\right) ^{e_{2}}B(g\otimes
1_{H};G^{a}X_{1}^{b_{1}}X_{2}^{b_{2}},g^{d}x_{1}x_{2}^{e_{2}})G^{a}X_{1}^{b_{1}}X_{2}^{b_{2}}\otimes g^{d}x_{2}^{e_{2}}+
\\
&&\sum_{\substack{ a,b_{2},d,e_{1},e_{2}=0  \\ a+b_{2}+d+e_{1}+e_{2}\equiv 0
}}^{1}\left( -1\right) ^{a+1}B(g\otimes
1_{H};G^{a}X_{1}X_{2}^{b_{2}},g^{d}x_{1}^{e_{1}}x_{2}^{e_{2}})G^{a}X_{2}^{b_{2}}\otimes g^{d}x_{1}^{e_{1}}x_{2}^{e_{2}}+
\\
&&\sum_{\substack{ a,b_{1},b_{2},d=0  \\ a+b_{1}+b_{2}+d\equiv 0}}%
^{1}B(x_{2}\otimes
1_{H};G^{a}X_{1}^{b_{1}}X_{2}^{b_{2}},g^{d}x_{1}x_{2})G^{a}X_{1}^{b_{1}}X_{2}^{b_{2}}\otimes g^{d}+
\\
&&\sum_{\substack{ a,b_{1},d,e_{2}=0  \\ a+b_{1}+d+e_{2}\equiv 0}}^{1}\left(
-1\right) ^{e_{2}+a+b_{1}}B(x_{2}\otimes
1_{H};G^{a}X_{1}^{b_{1}}X_{2},g^{d}x_{1}x_{2}^{e_{2}})G^{a}X_{1}^{b_{1}}%
\otimes g^{d}x_{2}^{e_{2}}+ \\
&&\sum_{\substack{ a,b_{2},d,e_{1}=0  \\ a+b_{2}+d+e_{1}\equiv 0}}^{1}\left(
-1\right) ^{a+1}B(x_{2}\otimes
1_{H};G^{a}X_{1}X_{2}^{b_{2}},g^{d}x_{1}^{e_{1}}x_{2})G^{a}X_{2}^{b_{2}}%
\otimes g^{d}x_{1}^{e_{1}}+ \\
&&\sum_{\substack{ a,d,e_{1},e_{2}=0  \\ a+d+e_{1}+e_{2}\equiv 0}}%
^{1}B(x_{2}\otimes
1_{H};G^{a}X_{1}X_{2},g^{d}x_{1}^{e_{1}}x_{2}^{e_{2}})G^{a}\otimes
g^{d}x_{1}^{e_{1}}x_{2}^{e_{2}}
\end{eqnarray*}

\subsubsection{$G^{a}\otimes g^{d}$}

\begin{eqnarray*}
&&\sum_{\substack{ a,d=0  \\ a+d\equiv 0}}^{1}\left[
\begin{array}{c}
B(g\otimes 1_{H};G^{a},g^{d}x_{1})+\left( -1\right) ^{a+1}B(g\otimes
1_{H};G^{a}X_{1},g^{d})+ \\
+B(x_{2}\otimes 1_{H};G^{a},g^{d}x_{1}x_{2}) \\
+\left( -1\right) ^{a_{1}}B(x_{2}\otimes 1_{H};G^{a}X_{2},g^{d}x_{1})+\left(
-1\right) ^{a+1}B(x_{2}\otimes 1_{H};G^{a}X_{1},g^{d}x_{2}) \\
+B(x_{2}\otimes 1_{H};G^{a}X_{1}X_{2},g^{d})%
\end{array}%
\right] \\
G^{a}\otimes g^{d} &=&0
\end{eqnarray*}%
and we obtain%
\begin{equation}
\begin{array}{c}
B(g\otimes 1_{H};1_{A},x_{1})-B(g\otimes 1_{H};X_{1},1_{H})+B(x_{2}\otimes
1_{H};1_{A},x_{1}x_{2}) \\
+B(x_{2}\otimes 1_{H};X_{2},x_{1})-B(x_{2}\otimes
1_{H};X_{1},x_{2})+B(x_{2}\otimes 1_{H};X_{1}X_{2},1_{H})%
\end{array}%
=0  \label{x2ot1, second}
\end{equation}%
and
\begin{equation}
\begin{array}{c}
B(g\otimes 1_{H};G,gx_{1})+B(g\otimes 1_{H};GX_{1},g)+B(x_{2}\otimes
1_{H};G,gx_{1}x_{2}) \\
-B(x_{2}\otimes 1_{H};GX_{2},gx_{1})+B(x_{2}\otimes
1_{H};GX_{1},gx_{2})+B(x_{2}\otimes 1_{H};GX_{1}X_{2},g)%
\end{array}%
=0  \label{x2ot1, third}
\end{equation}

\subsubsection{$G^{a}\otimes g^{d}x_{2}$}

\begin{equation*}
\sum_{\substack{ a,d=0  \\ a+d\equiv 1}}^{1}\left[
\begin{array}{c}
-B(g\otimes 1_{H};G^{a},g^{d}x_{1}x_{2})+\left( -1\right) ^{a+1}B(g\otimes
1_{H};G^{a}X_{1},g^{d}x_{2})+ \\
\left( -1\right) ^{a+1}B(x_{2}\otimes
1_{H};G^{a}X_{2},g^{d}x_{1}x_{2})+B(x_{2}\otimes
1_{H};G^{a}X_{1}X_{2},g^{d}x_{2})%
\end{array}%
\right] G^{a}\otimes g^{d}x_{2}=0
\end{equation*}%
and we get%
\begin{equation}
\begin{array}{c}
-B(g\otimes 1_{H};1_{A},gx_{1}x_{2})-B(g\otimes 1_{H};X_{1},gx_{2})+ \\
-B(x_{2}\otimes 1_{H};X_{2},gx_{1}x_{2})+B(x_{2}\otimes
1_{H};X_{1}X_{2},gx_{2})%
\end{array}%
=0  \label{x2ot1, fourth}
\end{equation}%
\begin{equation}
\begin{array}{c}
-B(g\otimes 1_{H};G,x_{1}x_{2})+B(g\otimes 1_{H};GX_{1},x_{2})+ \\
+B(x_{2}\otimes 1_{H};GX_{2},x_{1}x_{2})+B(x_{2}\otimes
1_{H};GX_{1}X_{2},x_{2})%
\end{array}%
=0  \label{x2ot1, fifth}
\end{equation}

\subsubsection{$G^{a}\otimes g^{d}x_{1}$}

\begin{gather*}
\sum_{\substack{ a,d=0  \\ a+d\equiv 1}}^{1}\left[
\begin{array}{c}
\left( -1\right) ^{a+1}B(g\otimes 1_{H};G^{a}X_{1},g^{d}x_{1})+\left(
-1\right) ^{a+1}B(x_{2}\otimes 1_{H};G^{a}X_{1},g^{d}x_{1}x_{2})+ \\
B(x_{2}\otimes 1_{H};G^{a}X_{1}X_{2},g^{d}x_{1})%
\end{array}%
\right] \\
G^{a}\otimes g^{d}x_{1}=0
\end{gather*}%
and we get%
\begin{equation}
-B(g\otimes 1_{H};X_{1},gx_{1})-B(x_{2}\otimes
1_{H};X_{1},gx_{1}x_{2})+B(x_{2}\otimes 1_{H};X_{1}X_{2},gx_{1})=0
\label{x2ot1, sixth}
\end{equation}%
\begin{equation}
B(g\otimes 1_{H};GX_{1},x_{1})+B(x_{2}\otimes
1_{H};GX_{1},x_{1}x_{2})+B(x_{2}\otimes 1_{H};GX_{1}X_{2},x_{1})=0
\label{x2ot1, seventh}
\end{equation}

\subsubsection{$G^{a}X_{2}\otimes g^{d}$}

\begin{equation*}
\sum_{\substack{ a,d=0  \\ a+d\equiv 1}}^{1}\left[
\begin{array}{c}
B(g\otimes 1_{H};G^{a}X_{2},g^{d}x_{1})+\left( -1\right) ^{a+1}B(g\otimes
1_{H};G^{a}X_{1}X_{2},g^{d})+ \\
B(x_{2}\otimes 1_{H};G^{a}X_{2},g^{d}x_{1}x_{2})+\left( -1\right)
^{a+1}B(x_{2}\otimes 1_{H};G^{a}X_{1}X_{2},g^{d}x_{2})%
\end{array}%
\right] G^{a}X_{2}\otimes g^{d}=0
\end{equation*}%
and we get%
\begin{equation}
\begin{array}{c}
B(g\otimes 1_{H};X_{2},gx_{1})-B(g\otimes 1_{H};X_{1}X_{2},g)+ \\
B(x_{2}\otimes 1_{H};X_{2},gx_{1}x_{2})-B(x_{2}\otimes
1_{H};X_{1}X_{2},gx_{2})%
\end{array}%
=0  \label{x2ot1, eighth}
\end{equation}%
and%
\begin{equation}
\begin{array}{c}
B(g\otimes 1_{H};GX_{2},x_{1})+B(g\otimes 1_{H};GX_{1}X_{2},1_{H})+ \\
B(x_{2}\otimes 1_{H};GX_{2},x_{1}x_{2})+B(x_{2}\otimes
1_{H};GX_{1}X_{2},x_{2})%
\end{array}%
=0  \label{x2ot1, ninth}
\end{equation}

\subsubsection{$G^{a}X_{1}\otimes g^{d}$}

\begin{gather*}
\sum_{\substack{ a,d=0  \\ a+d\equiv 1}}^{1}\left[
\begin{array}{c}
B(g\otimes 1_{H};G^{a}X_{1},g^{d}x_{1})+B(x_{2}\otimes
1_{H};G^{a}X_{1},g^{d}x_{1}x_{2}) \\
+\left( -1\right) ^{a+1}B(x_{2}\otimes 1_{H};G^{a}X_{1}X_{2},g^{d}x_{1})%
\end{array}%
\right] \\
G^{a}X_{1}\otimes g^{d}=0
\end{gather*}%
and we get%
\begin{equation}
B(g\otimes 1_{H};X_{1},gx_{1})+B(x_{2}\otimes
1_{H};X_{1},gx_{1}x_{2})-B(x_{2}\otimes 1_{H};X_{1}X_{2},gx_{1})=0
\label{x2ot1, tenth}
\end{equation}%
\begin{equation}
B(g\otimes 1_{H};GX_{1},x_{1})+B(x_{2}\otimes
1_{H};GX_{1},x_{1}x_{2})+B(x_{2}\otimes 1_{H};GX_{1}X_{2},x_{1})=0
\label{x2ot1, eleventh}
\end{equation}

\subsubsection{$G^{a}\otimes g^{d}x_{1}x_{2}$}

\begin{equation*}
\sum_{\substack{ a,d=0  \\ a+d\equiv 0}}^{1}\left[ \left( -1\right)
^{a+1}B(g\otimes 1_{H};G^{a}X_{1},g^{d}x_{1}x_{2})+B(x_{2}\otimes
1_{H};G^{a}X_{1}X_{2},g^{d}x_{1}x_{2})\right] G^{a}\otimes g^{d}x_{1}x_{2}=0
\end{equation*}%
and we get%
\begin{equation}
-B(g\otimes 1_{H};X_{1},x_{1}x_{2})+B(x_{2}\otimes
1_{H};X_{1}X_{2},x_{1}x_{2})=0  \label{x2ot1, twelfth}
\end{equation}%
\begin{equation}
B(g\otimes 1_{H};GX_{1},gx_{1}x_{2})+B(x_{2}\otimes
1_{H};GX_{1}X_{2},gx_{1}x_{2})=0  \label{x2ot1, thirteenth}
\end{equation}

\subsubsection{$G^{a}X_{2}\otimes g^{d}x_{2}$}

\begin{equation*}
\sum_{\substack{ a,d=0  \\ a+d\equiv 0}}^{1}\left[ -B(g\otimes
1_{H};G^{a}X_{2},g^{d}x_{1}x_{2})+\left( -1\right) ^{a+1}B(g\otimes
1_{H};G^{a}X_{1}X_{2},g^{d}x_{2})\right] G^{a}X_{2}\otimes g^{d}x_{2}=0
\end{equation*}%
and we get%
\begin{equation*}
-B(g\otimes 1_{H};X_{2},x_{1}x_{2})-B(g\otimes 1_{H};X_{1}X_{2},x_{2})=0
\end{equation*}%
\begin{equation*}
-B(g\otimes 1_{H};GX_{2},gx_{1}x_{2})+B(g\otimes 1_{H};GX_{1}X_{2},gx_{2})=0
\end{equation*}%
which are already known.

\subsubsection{$G^{a}X_{1}\otimes g^{d}x_{2}$}

\begin{equation*}
\sum_{\substack{ a,d=0  \\ a+d\equiv 0}}^{1}\left[ -B(g\otimes
1_{H};G^{a}X_{1},g^{d}x_{1}x_{2})+\left( -1\right) ^{a}B(x_{2}\otimes
1_{H};G^{a}X_{1}X_{2},g^{d}x_{1}x_{2})\right] G^{a}X_{1}\otimes g^{d}x_{2}=0
\end{equation*}%
and we get%
\begin{equation*}
-B(g\otimes 1_{H};X_{1},x_{1}x_{2})+B(x_{2}\otimes
1_{H};X_{1}X_{2},x_{1}x_{2})=0
\end{equation*}%
\begin{equation*}
-B(g\otimes 1_{H};GX_{1},gx_{1}x_{2})-B(x_{2}\otimes
1_{H};GX_{1}X_{2},gx_{1}x_{2})=0
\end{equation*}%
which are already known

\subsubsection{$G^{a}X_{2}\otimes g^{d}x_{1}$}

\begin{gather*}
\sum_{\substack{ a,d=0  \\ a+d\equiv 0}}^{1}\left[ \left( -1\right)
^{a+1}B(g\otimes 1_{H};G^{a}X_{1}X_{2},g^{d}x_{1})+\left( -1\right)
^{a+1}B(x_{2}\otimes 1_{H};G^{a}X_{1}X_{2},g^{d}x_{1}x_{2})\right] \\
G^{a}X_{2}\otimes g^{d}x_{1}=0
\end{gather*}%
and we get%
\begin{equation}
-B(g\otimes 1_{H};X_{1}X_{2},x_{1})-B(x_{2}\otimes
1_{H};X_{1}X_{2},x_{1}x_{2})=0  \label{x2ot1, fourteen}
\end{equation}%
\begin{equation}
B(g\otimes 1_{H};GX_{1}X_{2},gx_{1})+B(x_{2}\otimes
1_{H};GX_{1}X_{2},gx_{1}x_{2})=0  \label{x2ot1, fifteen}
\end{equation}

\subsubsection{$G^{a}X_{1}\otimes g^{d}x_{1}$}

There are no terms like this.

\subsubsection{$G^{a}X_{1}X_{2}\otimes g^{d}$}

\begin{equation*}
\sum_{\substack{ a,d=0  \\ a+d\equiv 0}}^{1}\left[ B(g\otimes
1_{H};G^{a}X_{1}X_{2},g^{d}x_{1})+B(x_{2}\otimes
1_{H};G^{a}X_{1}X_{2},g^{d}x_{1}x_{2})\right] G^{a}X_{1}X_{2}\otimes g^{d}=0
\end{equation*}%
and we get%
\begin{equation*}
B(g\otimes 1_{H};X_{1}X_{2},x_{1})+B(x_{2}\otimes
1_{H};X_{1}X_{2},x_{1}x_{2})=0
\end{equation*}%
\begin{equation*}
B(g\otimes 1_{H};GX_{1}X_{2},gx_{1})+B(x_{2}\otimes
1_{H};GX_{1}X_{2},gx_{1}x_{2})=0
\end{equation*}%
which are already known.

\subsubsection{$G^{a}X_{2}\otimes g^{d}x_{1}x_{2}$}

\begin{equation*}
\sum_{\substack{ a,d=0  \\ a+d\equiv 1}}^{1}\left( -1\right)
^{a+1}B(g\otimes 1_{H};G^{a}X_{1}X_{2},g^{d}x_{1}x_{2})G^{a}X_{2}\otimes
g^{d}x_{1}x_{2}=0
\end{equation*}%
and we get%
\begin{equation*}
-B(g\otimes 1_{H};X_{1}X_{2},gx_{1}x_{2})=0
\end{equation*}%
\begin{equation*}
B(g\otimes 1_{H};GX_{1}X_{2},x_{1}x_{2})=0
\end{equation*}%
which are already known.

\subsubsection{$G^{a}X_{1}\otimes g^{d}x_{1}x_{2}$}

There is no term like this.

\subsubsection{$G^{a}X_{1}X_{2}\otimes g^{d}x_{2}$}

\begin{equation*}
\sum_{\substack{ a,d=0  \\ a+d\equiv 1}}^{1}-B(g\otimes
1_{H};G^{a}X_{1}X_{2},g^{d}x_{1}x_{2})G^{a}X_{1}X_{2}\otimes g^{d}x_{2}=0
\end{equation*}%
and we get%
\begin{equation*}
-B(g\otimes 1_{H};X_{1}X_{2},gx_{1}x_{2})=0
\end{equation*}%
\begin{equation*}
-B(g\otimes 1_{H};GX_{1}X_{2},x_{1}x_{2})=0
\end{equation*}%
which are already known.

\subsubsection{$G^{a}X_{1}X_{2}\otimes g^{d}x_{1}$}

There is no term like this.

\paragraph{$G^{a}X_{1}X_{2}\otimes g^{d}x_{1}x_{2}$}

There is no term like this.

\subsection{Case $gx_{1}$}

From the first summand of the left side we can not obtain $gx_{1}.$From the
second summand of the left side we get%
\begin{eqnarray*}
g^{a+b_{1}+b_{2}+l_{1}+l_{2}+d+e_{1}+e_{2}+u_{1}+u_{2}}x_{1}^{l_{1}+u_{1}}x_{2}^{l_{2}+u_{2}} &=&gx_{1}\Rightarrow
\\
l_{1}+u_{1} &=&1 \\
l_{2}+u_{2} &=&0\Rightarrow l_{2}=u_{2}=0 \\
a+b_{1}+b_{2}+l_{1}+l_{2}+d+e_{1}+e_{2}+u_{1}+u_{2} &\equiv &1 \\
a+b_{1}+b_{2}+d+e_{1}+e_{2} &\equiv &0\text{ }
\end{eqnarray*}%
\begin{gather*}
\sum_{\substack{ a,b_{1},b_{2},d,e_{1},e_{2}=0  \\ %
a+b_{1}+b_{2}+d+e_{1}+e_{2}\equiv 0\text{ }}}^{1}\sum_{l_{1}=0}^{b_{1}}\sum
_{\substack{ u_{1}=0  \\ l_{1}+u_{1}=1}}^{e_{1}}\left( -1\right) ^{\alpha
\left( 1_{H};l_{1},0,u_{1},0\right) }B(x_{2}\otimes
1_{H};G^{a}X_{1}^{b_{1}}X_{2}^{b_{2}},g^{d}x_{1}^{e_{1}}x_{2}^{e_{2}}) \\
G^{a}X_{1}^{b_{1}-l_{1}}X_{2}^{b_{2}}\otimes
g^{d}x_{1}^{e_{1}-u_{1}}x_{2}^{e_{2}}=0
\end{gather*}%
Since%
\begin{equation*}
\alpha \left( g^{m}x_{1}^{n_{1}}x_{2}^{n_{2}};0,0,1,0\right) =e_{2}+\left(
a+b_{1}+b_{2}\right)
\end{equation*}%
and%
\begin{equation*}
\alpha \left( g^{m}x_{1}^{n_{1}}x_{2}^{n_{2}};1,0,0,0\right) =b_{2}
\end{equation*}%
we obtain%
\begin{gather*}
\sum_{\substack{ a,b_{1},b_{2},d,e_{2}=0  \\ a+b_{1}+b_{2}+d+e_{2}\equiv 1%
\text{ }}}^{1}\left( -1\right) ^{e_{2}+\left( a+b_{1}+b_{2}\right)
}B(x_{2}\otimes
1_{H};G^{a}X_{1}^{b_{1}}X_{2}^{b_{2}},g^{d}x_{1}x_{2}^{e_{2}})G^{a}X_{1}^{b_{1}}X_{2}^{b_{2}}\otimes g^{d}x_{2}^{e_{2}}+
\\
\sum_{\substack{ a,b_{2},d,e_{1},e_{2}=0  \\ a+b_{2}+d+e_{1}+e_{2}\equiv 1%
\text{ }}}^{1}\left( -1\right) ^{b_{2}}B(x_{2}\otimes
1_{H};G^{a}X_{1}X_{2}^{b_{2}},g^{d}x_{1}^{e_{1}}x_{2}^{e_{2}})G^{a}X_{2}^{b_{2}}\otimes g^{d}x_{1}^{e_{1}}x_{2}^{e_{2}}=0
\end{gather*}

\subsubsection{$G^{a}\otimes g^{d}$}

\begin{equation*}
\sum_{\substack{ a,d=0  \\ a+d\equiv 1\text{ }}}^{1}\left[ \left( -1\right)
^{a}B(x_{2}\otimes 1_{H};G^{a},g^{d}x_{1})+B(x_{2}\otimes
1_{H};G^{a}X_{1},g^{d})\right] G^{a}\otimes g^{d}=0
\end{equation*}%
and we get%
\begin{equation}
B(x_{2}\otimes 1_{H};1_{A},gx_{1})+B(x_{2}\otimes 1_{H};X_{1},g)=0
\label{x2ot1, sixteen}
\end{equation}%
\begin{equation}
-B(x_{2}\otimes 1_{H};G,x_{1})+B(x_{2}\otimes 1_{H};GX_{1},1_{H})=0
\label{x2ot1, seventeen}
\end{equation}

\subsubsection{$G^{a}\otimes g^{d}x_{2}$}

\begin{equation*}
\sum_{\substack{ a,d=0  \\ a+d\equiv 0\text{ }}}^{1}\left[ \left( -1\right)
^{a+1}B(x_{2}\otimes 1_{H};G^{a},g^{d}x_{1}x_{2})+B(x_{2}\otimes
1_{H};G^{a}X_{1},g^{d}x_{2})\right] G^{a}\otimes g^{d}x_{2}=0
\end{equation*}%
and we get%
\begin{equation}
-B(x_{2}\otimes 1_{H};1_{A},x_{1}x_{2})+B(x_{2}\otimes 1_{H};X_{1},x_{2})=0
\label{x2ot1, eighteen}
\end{equation}%
\begin{equation}
B(x_{2}\otimes 1_{H};G,gx_{1}x_{2})+B(x_{2}\otimes 1_{H};GX_{1},gx_{2})=0
\label{x2ot1, nineteen}
\end{equation}

\subsubsection{$G^{a}\otimes g^{d}x_{1}$}

\begin{equation*}
\sum_{\substack{ a,d=0  \\ a+d\equiv 0\text{ }}}^{1}B(x_{2}\otimes
1_{H};G^{a}X_{1},g^{d}x_{1})G^{a}\otimes g^{d}x_{1}=0
\end{equation*}%
and we get%
\begin{equation}
B(x_{2}\otimes 1_{H};X_{1},x_{1})=0  \label{x2ot1, twenty}
\end{equation}%
\begin{equation}
B(x_{2}\otimes 1_{H};GX_{1},gx_{1})=0  \label{x2ot1,twentyone}
\end{equation}

\subsubsection{$G^{a}X_{2}\otimes g^{d}$}

\begin{equation*}
\sum_{\substack{ a,d=0  \\ a+d\equiv 0\text{ }}}^{1}\left[ \left( -1\right)
^{a+1}B(x_{2}\otimes 1_{H};G^{a}X_{2},g^{d}x_{1})-B(x_{2}\otimes
1_{H};G^{a}X_{1}X_{2},g^{d})\right] G^{a}X_{2}\otimes g^{d}=0
\end{equation*}%
and we get%
\begin{equation}
-B(x_{2}\otimes 1_{H};X_{2},x_{1})-B(x_{2}\otimes 1_{H};X_{1}X_{2},1_{H})=0
\label{x2ot1,twentytwo}
\end{equation}%
\begin{equation}
B(x_{2}\otimes 1_{H};GX_{2},gx_{1})-B(x_{2}\otimes 1_{H};GX_{1}X_{2},g)=0
\label{x2ot1,twentythree}
\end{equation}

\subsubsection{$G^{a}X_{1}\otimes g^{d}$}

\begin{equation*}
\sum_{\substack{ a,d=0  \\ a+d\equiv 0\text{ }}}^{1}\left( -1\right)
^{a+1}B(x_{2}\otimes 1_{H};G^{a}X_{1},g^{d}x_{1})G^{a}X_{1}\otimes g^{d}=0
\end{equation*}%
and we get%
\begin{equation*}
-B(x_{2}\otimes 1_{H};X_{1},x_{1})=0
\end{equation*}%
\begin{equation*}
B(x_{2}\otimes 1_{H};GX_{1},gx_{1})=0
\end{equation*}%
which we already got.

\subsubsection{$G^{a}\otimes g^{d}x_{1}x_{2}$}

\begin{equation*}
\sum_{\substack{ a,d=0  \\ a+d\equiv 1\text{ }}}^{1}B(x_{2}\otimes
1_{H};G^{a}X_{1},g^{d}x_{1}x_{2})G^{a}\otimes g^{d}x_{1}x_{2}=0
\end{equation*}%
and we get%
\begin{equation}
B(x_{2}\otimes 1_{H};X_{1},gx_{1}x_{2})=0  \label{x2ot1,twentyfour}
\end{equation}%
\begin{equation}
B(x_{2}\otimes 1_{H};GX_{1},x_{1}x_{2})=0  \label{x2ot1,twentyfive}
\end{equation}

\subsubsection{$G^{a}X_{2}\otimes g^{d}x_{2}$}

\begin{equation*}
\sum_{\substack{ a,d=0  \\ a+d\equiv 1\text{ }}}^{1}\left[ \left( -1\right)
^{a}B(x_{2}\otimes 1_{H};G^{a}X_{2},g^{d}x_{1}x_{2})-B(x_{2}\otimes
1_{H};G^{a}X_{1}X_{2},g^{d}x_{2})\right] G^{a}X_{2}\otimes g^{d}x_{2}=0
\end{equation*}%
and we get%
\begin{equation}
B(x_{2}\otimes 1_{H};X_{2},gx_{1}x_{2})-B(x_{2}\otimes
1_{H};X_{1}X_{2},gx_{2})=0  \label{x2ot1,twentysix}
\end{equation}%
\begin{equation}
-B(x_{2}\otimes 1_{H};GX_{2},x_{1}x_{2})-B(x_{2}\otimes
1_{H};GX_{1}X_{2},x_{2})=0  \label{x2ot1,twentyseven}
\end{equation}

\subsubsection{$G^{a}X_{1}\otimes g^{d}x_{2}$}

\begin{equation*}
\sum_{\substack{ a,d=0  \\ a+d\equiv 1\text{ }}}^{1}\left( -1\right)
^{a}B(x_{2}\otimes 1_{H};G^{a}X_{1},g^{d}x_{1}x_{2})G^{a}X_{1}\otimes
g^{d}x_{2}=0
\end{equation*}%
and we get%
\begin{equation*}
B(x_{2}\otimes 1_{H};X_{1},gx_{1}x_{2})=0
\end{equation*}%
\begin{equation*}
-B(x_{2}\otimes 1_{H};GX_{1},x_{1}x_{2})=0
\end{equation*}%
which we already got.

\subsubsection{$G^{a}X_{2}\otimes g^{d}x_{1}$}

\begin{equation*}
\sum_{\substack{ a,d=0  \\ a+d\equiv 1\text{ }}}^{1}-B(x_{2}\otimes
1_{H};G^{a}X_{1}X_{2},g^{d}x_{1})G^{a}X_{2}\otimes g^{d}x_{1}=0
\end{equation*}%
and we get%
\begin{equation}
-B(x_{2}\otimes 1_{H};X_{1}X_{2},gx_{1})=0  \label{x2ot1,twentyeight}
\end{equation}%
\begin{equation}
-B(x_{2}\otimes 1_{H};GX_{1}X_{2},x_{1})=0  \label{x2ot1,twentynine}
\end{equation}

\subsubsection{$G^{a}X_{1}\otimes g^{d}x_{1}$}

There is no term like this.

\subsubsection{$G^{a}X_{1}X_{2}\otimes g^{d}$}

\begin{equation*}
\sum_{\substack{ a,d=0  \\ a+d\equiv 1\text{ }}}^{1}\left( -1\right)
^{a}B(x_{2}\otimes 1_{H};G^{a}X_{1}X_{2},g^{d}x_{1})G^{a}X_{1}X_{2}\otimes
g^{d}=0
\end{equation*}%
and we get%
\begin{equation*}
B(x_{2}\otimes 1_{H};X_{1}X_{2},gx_{1})=0
\end{equation*}%
\begin{equation*}
-B(x_{2}\otimes 1_{H};GX_{1}X_{2},x_{1})=0
\end{equation*}%
which we already got.

\subsubsection{$G^{a}X_{2}\otimes g^{d}x_{1}x_{2}$}

\begin{equation*}
\sum_{\substack{ a,d=0  \\ a+d\equiv 0\text{ }}}^{1}-B(x_{2}\otimes
1_{H};G^{a}X_{1}X_{2},g^{d}x_{1}x_{2})G^{a}X_{2}\otimes g^{d}x_{1}x_{2}=0
\end{equation*}%
and we get%
\begin{equation}
-B(x_{2}\otimes 1_{H};X_{1}X_{2},x_{1}x_{2})=0  \label{x2ot1,thirty}
\end{equation}%
\begin{equation}
-B(x_{2}\otimes 1_{H};GX_{1}X_{2},gx_{1}x_{2})=0  \label{x2ot1,thirtyone}
\end{equation}

\subsubsection{$G^{a}X_{1}\otimes g^{d}x_{1}x_{2}$}

There is no term like this.

\subsubsection{$G^{a}X_{1}X_{2}\otimes g^{d}x_{2}$}

\begin{equation*}
\sum_{\substack{ a,d=0  \\ a+d\equiv 0\text{ }}}^{1}\left( -1\right)
^{a+1}B(x_{2}\otimes
1_{H};G^{a}X_{1}X_{2},g^{d}x_{1}x_{2})G^{a}X_{1}X_{2}\otimes g^{d}x_{2}=0
\end{equation*}%
and we get%
\begin{equation*}
-B(x_{2}\otimes 1_{H};X_{1}X_{2},x_{1}x_{2})=0
\end{equation*}%
\begin{equation*}
B(x_{2}\otimes 1_{H};GX_{1}X_{2},gx_{1}x_{2})=0
\end{equation*}%
which we already got.

\subsubsection{$G^{a}X_{1}X_{2}\otimes g^{d}x_{1}$}

There is no term like this.

\subsubsection{$G^{a}X_{1}X_{2}\otimes g^{d}x_{1}x_{2}$}

There is no term like this.

\subsection{Case $gx_{2}$}

From the first summand of the left side we get%
\begin{eqnarray*}
g^{a+b_{1}+b_{2}+d+e_{1}+e_{2}+l_{1}+l_{2}+u_{1}+u_{2}}x_{1}^{u_{1}+l_{1}}x_{2}^{l_{2}+u_{2}+1} &\equiv &gx_{2}
\\
a+b_{1}+b_{2}+d+e_{1}+e_{2}+l_{1}+l_{2}+u_{1}+u_{2} &\equiv &1 \\
l_{1}+u_{1} &=&0\text{ so that }l_{1}=0=u_{1} \\
l_{2}+u_{2} &=&0\text{ so that }l_{2}=0=u_{2}\text{ and} \\
a+b_{1}+b_{2}+d+e_{1}+e_{2} &\equiv &1\text{ }
\end{eqnarray*}%
and%
\begin{equation*}
\alpha \left( g^{m}x_{1}^{n_{1}}x_{2}^{n_{2}};0,0,0,0\right) =\left(
a+b_{1}+b_{2}\right) \left( n_{1}+n_{2}\right)
\end{equation*}
\begin{equation*}
\alpha =\alpha \left( x_{2},0,0,0,0\right) =a+b_{1}+b_{2}
\end{equation*}%
From the second summand of the left side we get%
\begin{eqnarray*}
g^{a+b_{1}+b_{2}+d+e_{1}+e_{2}+l_{1}+l_{2}+u_{1}+u_{2}}x_{1}^{u_{1}+l_{1}}x_{2}^{l_{2}+u_{2}} &=&gx_{2}\Rightarrow
\\
u_{1} &=&l_{1}=0 \\
l_{2}+u_{2} &=&1 \\
a+b_{1}+b_{2}+d+e_{1}+e_{2}+l_{1}+l_{2}+u_{1}+u_{2} &\equiv &1\text{ so
that, } \\
a+b_{1}+b_{2}+d+e_{1}+e_{2} &\equiv &0\text{ }
\end{eqnarray*}%
Since $\alpha \left( 1_{H};0,0,0,1\right) =a+b_{1}+b_{2}$ and $\alpha \left(
1_{H};0,1,0,0\right) =0$

\begin{eqnarray*}
0 &=&\sum_{\substack{ a,b_{1},b_{2},d,e_{1},e_{2}=0  \\ %
a+b_{1}+b_{2}+d+e_{1}+e_{2}\equiv 1\text{ }}}^{1}\left( -1\right) ^{\left(
a+b_{1}+b_{2}\right) }B(g\otimes
1_{H};G^{a}X_{1}^{b_{1}}X_{2}^{b_{2}},g^{d}x_{1}^{e_{1}}x_{2}^{e_{2}})G^{a}X_{1}^{b_{1}}X_{2}^{b_{2}}\otimes g^{d}x_{1}^{e_{1}}x_{2}^{e_{2}}
\\
&&+\sum_{\substack{ a,b_{1},b_{2},d,e_{1},e_{2}=0  \\ %
a+b_{1}+b_{2}+d+e_{1}+e_{2}\equiv 0\text{ }}}^{1}\sum_{l_{2}=0}^{b_{2}}\sum
_{\substack{ u_{2}=0  \\ l_{2}+u_{2}=1}}^{e_{2}}(-1)^{\alpha \left(
1_{H};0,l_{2},0,u_{2}\right) } \\
&&B(x_{2}\otimes \
1_{H};G^{a}X_{1}^{b_{1}}X_{2}^{b_{2}},g^{d}x_{1}^{e_{1}}x_{2}^{e_{2}})G^{a}X_{1}^{b_{1}}X_{2}^{b_{2}-l_{2}}\otimes g^{d}x_{1}^{e_{1}}x_{2}^{\left( e_{2}-u_{2}\right) }
\\
&=&\sum_{\substack{ a,b_{1},b_{2},d,e_{1},e_{2}=0  \\ %
a+b_{1}+b_{2}+d+e_{1}+e_{2}\equiv 1\text{ }}}^{1}\left( -1\right) ^{\left(
a+b_{1}+b_{2}\right) }B(g\otimes
1_{H};G^{a}X_{1}^{b_{1}}X_{2}^{b_{2}},g^{d}x_{1}^{e_{1}}x_{2}^{e_{2}}) \\
&&G^{a}X_{1}^{b_{1}}X_{2}^{b_{2}}\otimes g^{d}x_{1}^{e_{1}}x_{2}^{e_{2}}+ \\
&&+\sum_{\substack{ a,b_{1},b_{2},d,e_{1},e_{2}=0  \\ a+b_{1}+b_{2}+d+e_{1}%
\equiv 1\text{ }}}^{1}(-1)^{a+b_{1}+b_{2}}B(x_{2}\otimes \
1_{H};G^{a}X_{1}^{b_{1}}X_{2}^{b_{2}},g^{d}x_{1}^{e_{1}}x_{2})G^{a}X_{1}^{b_{1}}X_{2}^{b_{2}}\otimes g^{d}x_{1}^{e_{1}}
\\
&&+\sum_{\substack{ a,b_{1},d,e_{1},e_{2}=0  \\ a+b_{1}+d+e_{1}+e_{2}\equiv 1%
\text{ }}}^{1}B(x_{2}\otimes \
1_{H};G^{a}X_{1}^{b_{1}}X_{2},g^{d}x_{1}^{e_{1}}x_{2}^{e_{2}})G^{a}X_{1}^{b_{1}}\otimes g^{d}x_{1}^{e_{1}}x_{2}^{e_{2}}
\end{eqnarray*}%
and we obtain%
\begin{eqnarray*}
&&\sum_{\substack{ a,b_{1},b_{2},d,e_{1},e_{2}=0  \\ %
a+b_{1}+b_{2}+d+e_{1}+e_{2}\equiv 1\text{ }}}^{1}\left( -1\right) ^{\left(
a+b_{1}+b_{2}\right) }B(g\otimes
1_{H};G^{a}X_{1}^{b_{1}}X_{2}^{b_{2}},g^{d}x_{1}^{e_{1}}x_{2}^{e_{2}})G^{a}X_{1}^{b_{1}}X_{2}^{b_{2}}\otimes g^{d}x_{1}^{e_{1}}x_{2}^{e_{2}}+
\\
&&+\sum_{\substack{ a,b_{1},b_{2},d,e_{1},e_{2}=0  \\ a+b_{1}+b_{2}+d+e_{1}%
\equiv 1\text{ }}}^{1}(-1)^{a+b_{1}+b_{2}}B(x_{2}\otimes \
1_{H};G^{a}X_{1}^{b_{1}}X_{2}^{b_{2}},g^{d}x_{1}^{e_{1}}x_{2})G^{a}X_{1}^{b_{1}}X_{2}^{b_{2}}\otimes g^{d}x_{1}^{e_{1}}
\\
&&+\sum_{\substack{ a,b_{1},d,e_{1},e_{2}=0  \\ a+b_{1}+d+e_{1}+e_{2}\equiv 1%
\text{ }}}^{1}B(x_{2}\otimes \
1_{H};G^{a}X_{1}^{b_{1}}X_{2},g^{d}x_{1}^{e_{1}}x_{2}^{e_{2}})G^{a}X_{1}^{b_{1}}\otimes g^{d}x_{1}^{e_{1}}x_{2}^{e_{2}}.
\end{eqnarray*}

\subsubsection{$G^{a}\otimes g^{d}$}

\begin{equation*}
\sum_{\substack{ a,d=0  \\ a+d\equiv 1\text{ }}}^{1}\left[
\begin{array}{c}
\left( -1\right) ^{a}B(g\otimes 1_{H};G^{a},g^{d})+(-1)^{a}B(x_{2}\otimes \
1_{H};G^{a},g^{d}x_{2})+ \\
+B(x_{2}\otimes \ 1_{H};G^{a}X_{2},g^{d})%
\end{array}%
\right] G^{a}\otimes g^{d}=0
\end{equation*}%
and we get%
\begin{equation}
B(g\otimes 1_{H};1_{A},g)+B(x_{2}\otimes \
1_{H};1_{A},gx_{2})+B(x_{2}\otimes \ 1_{H};X_{2},g)=0
\label{x2ot1,thirtytwo}
\end{equation}%
\begin{equation}
-B(g\otimes 1_{H};G,1_{H})-B(x_{2}\otimes \ 1_{H};G,x_{2})+B(x_{2}\otimes \
1_{H};GX_{2},1_{H})=0  \label{x2ot1,thirtythree}
\end{equation}

\subsubsection{$G^{a}\otimes g^{d}x_{2}$}

\begin{equation*}
\sum_{\substack{ a,d=0  \\ a+d\equiv 0\text{ }}}^{1}\left[ \left( -1\right)
^{a}B(g\otimes 1_{H};G^{a},g^{d}x_{2})+B(x_{2}\otimes \
1_{H};G^{a}X_{2},g^{d}x_{2})\right] G^{a}\otimes g^{d}x_{2}=0
\end{equation*}%
and we get%
\begin{equation}
B(g\otimes 1_{H};1_{A},x_{2})+B(x_{2}\otimes \ 1_{H};X_{2},x_{2})=0
\label{x2ot1,thirtyfour}
\end{equation}%
\begin{equation}
-B(g\otimes 1_{H};G,gx_{2})+B(x_{2}\otimes \ 1_{H};GX_{2},gx_{2})=0
\label{x2ot1,thirtyfive}
\end{equation}

\subsubsection{$G^{a}\otimes g^{d}x_{1}$}

\begin{equation*}
\sum_{\substack{ a,d=0  \\ a+d\equiv 0\text{ }}}^{1}\left[
\begin{array}{c}
\left( -1\right) ^{a}B(g\otimes
1_{H};G^{a},g^{d}x_{1})+(-1)^{a}B(x_{2}\otimes \ 1_{H};G^{a},g^{d}x_{1}x_{2})
\\
+B(x_{2}\otimes \ 1_{H};G^{a}X_{2},g^{d}x_{1})%
\end{array}%
\right] G^{a}\otimes g^{d}x_{1}=0
\end{equation*}%
and we get%
\begin{equation}
B(g\otimes 1_{H};1_{A},x_{1})+B(x_{2}\otimes \
1_{H};1_{A},x_{1}x_{2})+B(x_{2}\otimes \ 1_{H};X_{2},x_{1})=0
\label{x2ot1,thirtysix}
\end{equation}%
\begin{equation}
-B(g\otimes 1_{H};G,gx_{1})-B(x_{2}\otimes \
1_{H};G,gx_{1}x_{2})+B(x_{2}\otimes \ 1_{H};GX_{2},gx_{1})=0
\label{x2ot1,thirtyseven}
\end{equation}

\subsubsection{$G^{a}X_{2}\otimes g^{d}$}

\begin{equation*}
\sum_{\substack{ a,d=0  \\ a+d\equiv 0\text{ }}}^{1}\left[ \left( -1\right)
^{a+1}B(g\otimes 1_{H};G^{a}X_{2},g^{d})+(-1)^{a+1}B(x_{2}\otimes \
1_{H};G^{a}X_{2},g^{d}x_{2})\right] G^{a}X_{2}\otimes g^{d}=0
\end{equation*}%
and we get%
\begin{equation}
-B(g\otimes 1_{H};X_{2},1_{H})-B(x_{2}\otimes \ 1_{H};X_{2},x_{2})=0
\label{x2ot1,thirtyeight}
\end{equation}%
\begin{equation}
B(g\otimes 1_{H};GX_{2},g)+B(x_{2}\otimes \ 1_{H};GX_{2},gx_{2})=0
\label{x2ot1,thirtynine}
\end{equation}

\subsubsection{$G^{a}X_{1}\otimes g^{d}$}

\begin{equation*}
\sum_{\substack{ a,d=0  \\ a+d\equiv 0\text{ }}}^{1}\left[
\begin{array}{c}
\left( -1\right) ^{\left( a+1\right) }B(g\otimes
1_{H};G^{a}X_{1},g^{d})+(-1)^{a+1}B(x_{2}\otimes \
1_{H};G^{a}X_{1},g^{d}x_{2}) \\
+B(x_{2}\otimes \ 1_{H};G^{a}X_{1}X_{2},g^{d})%
\end{array}%
\right] G^{a}X_{1}\otimes g^{d}=0
\end{equation*}%
and we get%
\begin{equation*}
-B(g\otimes 1_{H};X_{1},1_{H})-B(x_{2}\otimes \
1_{H};X_{1},x_{2})+B(x_{2}\otimes \ 1_{H};X_{1}X_{2},1_{H})=0
\end{equation*}%
\begin{equation*}
+B(g\otimes 1_{H};GX_{1},g)+B(x_{2}\otimes \
1_{H};GX_{1},gx_{2})+B(x_{2}\otimes \ 1_{H};GX_{1}X_{2},g)=0
\end{equation*}

\subsubsection{$G^{a}\otimes g^{d}x_{1}x_{2}$}

\begin{equation*}
\sum_{\substack{ a,d=0  \\ a+d\equiv 1\text{ }}}^{1}\left[ \left( -1\right)
^{a}B(g\otimes 1_{H};G^{a},g^{d}x_{1}x_{2})+B(x_{2}\otimes \
1_{H};G^{a}X_{2},g^{d}x_{1}x_{2})\right] G^{a}\otimes g^{d}x_{1}x_{2}=0
\end{equation*}%
and we get%
\begin{equation}
B(g\otimes 1_{H};1_{A},gx_{1}x_{2})+B(x_{2}\otimes \
1_{H};X_{2},gx_{1}x_{2})=0  \label{x2ot1,fourthy}
\end{equation}%
\begin{equation}
-B(g\otimes 1_{H};G,x_{1}x_{2})+B(x_{2}\otimes \ 1_{H};GX_{2},x_{1}x_{2})=0
\label{x2ot1,fourthyone}
\end{equation}

\subsubsection{$G^{a}X_{2}\otimes g^{d}x_{2}$}

\begin{equation*}
\sum_{\substack{ a,d=0  \\ a+d\equiv 1\text{ }}}^{1}\left( -1\right)
^{a+1}B(g\otimes 1_{H};G^{a}X_{2},g^{d}x_{2})G^{a}X_{2}\otimes g^{d}x_{2}=0
\end{equation*}

and we get%
\begin{equation}
-B(g\otimes 1_{H};X_{2},gx_{2})=0  \label{x2ot1,fourthytwo}
\end{equation}%
\begin{equation}
B(g\otimes 1_{H};GX_{2},x_{2})=0  \label{x2ot1,fourthythree}
\end{equation}

\subsubsection{$G^{a}X_{1}\otimes g^{d}x_{2}$}

\begin{equation*}
\sum_{\substack{ a,d=0  \\ a+d\equiv 1\text{ }}}^{1}\left[ \left( -1\right)
^{a+1}B(g\otimes 1_{H};G^{a}X_{1},g^{d}x_{2})+B(x_{2}\otimes \
1_{H};G^{a}X_{1}X_{2},g^{d}x_{2})\right] G^{a}X_{1}\otimes g^{d}x_{2}=0
\end{equation*}%
and we get%
\begin{equation}
-B(g\otimes 1_{H};X_{1},gx_{2})+B(x_{2}\otimes \ 1_{H};X_{1}X_{2},gx_{2})=0
\label{x2ot1,fourthyfour}
\end{equation}%
\begin{equation}
B(g\otimes 1_{H};GX_{1},x_{2})+B(x_{2}\otimes \ 1_{H};GX_{1}X_{2},x_{2})=0
\label{x2ot1,fourthyfive}
\end{equation}

\subsubsection{$G^{a}X_{2}\otimes g^{d}x_{1}$}

\begin{gather*}
\sum_{\substack{ a,d=0  \\ a+d\equiv 1\text{ }}}^{1}\left[ \left( -1\right)
^{a+1}B(g\otimes 1_{H};G^{a}X_{2},g^{d}x_{1})+(-1)^{a+1}B(x_{2}\otimes \
1_{H};G^{a}X_{2},g^{d}x_{1}x_{2})\right] \\
G^{a}X_{2}\otimes g^{d}x_{1}=0
\end{gather*}%
and we get%
\begin{equation}
-B(g\otimes 1_{H};X_{2},gx_{1})-B(x_{2}\otimes \ 1_{H};X_{2},gx_{1}x_{2})=0
\label{x2ot1,fourthysix}
\end{equation}%
\begin{equation}
B(g\otimes 1_{H};GX_{2},x_{1})+B(x_{2}\otimes \ 1_{H};GX_{2},x_{1}x_{2})=0
\label{x2ot1,fourthyseven}
\end{equation}

\subsubsection{$G^{a}X_{1}\otimes g^{d}x_{1}$}

\begin{gather*}
\sum_{\substack{ a,d=0  \\ a+d\equiv 1\text{ }}}^{1}\left[
\begin{array}{c}
\left( -1\right) ^{a+1}B(g\otimes
1_{H};G^{a}X_{1},g^{d}x_{1})+(-1)^{a+1}B(x_{2}\otimes \
1_{H};G^{a}X_{1},g^{d}x_{1}x_{2}) \\
+B(x_{2}\otimes \ 1_{H};G^{a}X_{1}X_{2},g^{d}x_{1})%
\end{array}%
\right] \\
G^{a}X_{1}\otimes g^{d}x_{1}=0
\end{gather*}%
and we get%
\begin{equation}
-B(g\otimes 1_{H};X_{1},gx_{1})-B(x_{2}\otimes \
1_{H};X_{1},gx_{1}x_{2})+B(x_{2}\otimes \ 1_{H};X_{1}X_{2},gx_{1})=0
\label{x2ot1,fourthyeight}
\end{equation}%
\begin{equation}
B(g\otimes 1_{H};GX_{1},x_{1})+B(x_{2}\otimes \
1_{H};GX_{1},x_{1}x_{2})+B(x_{2}\otimes \ 1_{H};GX_{1}X_{2},x_{1})=0
\label{x2ot1,fourthynine}
\end{equation}

\subsubsection{$G^{a}X_{1}X_{2}\otimes g^{d}$}

\begin{equation*}
\sum_{\substack{ a,d=0  \\ a+d\equiv 1\text{ }}}^{1}\left[ \left( -1\right)
^{a}B(g\otimes 1_{H};G^{a}X_{1}X_{2},g^{d})+\left( -1\right)
^{a}B(x_{2}\otimes \ 1_{H};G^{a}X_{1}X_{2},g^{d}x_{2})\right]
G^{a}X_{1}X_{2}\otimes g^{d}=0
\end{equation*}%
and we get%
\begin{equation}
B(g\otimes 1_{H};X_{1}X_{2},g)+B(x_{2}\otimes \ 1_{H};X_{1}X_{2},gx_{2})=0
\label{x2ot1,fifthy}
\end{equation}%
\begin{equation}
-B(g\otimes 1_{H};GX_{1}X_{2},1_{H})-B(x_{2}\otimes \
1_{H};GX_{1}X_{2},x_{2})=0  \label{x2ot1,fifthyone}
\end{equation}

\subsubsection{$G^{a}X_{2}\otimes g^{d}x_{1}x_{2}$}

\begin{equation*}
\sum_{\substack{ a,d=0  \\ a+d\equiv 0\text{ }}}^{1}\left( -1\right)
^{a+1}B(g\otimes 1_{H};G^{a}X_{2},g^{d}x_{1}x_{2})G^{a}X_{2}\otimes
g^{d}x_{1}x_{2}=0
\end{equation*}%
and we get%
\begin{equation*}
-B(g\otimes 1_{H};X_{2},x_{1}x_{2})=0
\end{equation*}%
\begin{equation*}
B(g\otimes 1_{H};GX_{2},gx_{1}x_{2})=0
\end{equation*}%
which are already known.

\subsubsection{$G^{a}X_{1}\otimes g^{d}x_{1}x_{2}$}

\begin{equation*}
\sum_{\substack{ a,d=0  \\ a+d\equiv 0\text{ }}}^{1}\left[ \left( -1\right)
^{a+1}B(g\otimes 1_{H};G^{a}X_{1},g^{d}x_{1}x_{2})+B(x_{2}\otimes \
1_{H};G^{a}X_{1}X_{2},g^{d}x_{1}x_{2})\right] G^{a}X_{1}\otimes
g^{d}x_{1}x_{2}=0
\end{equation*}%
and we get%
\begin{equation}
-B(g\otimes 1_{H};X_{1},x_{1}x_{2})+B(x_{2}\otimes \
1_{H};X_{1}X_{2},x_{1}x_{2})=0  \label{x2ot1,fifthytwo}
\end{equation}%
\begin{equation}
B(g\otimes 1_{H};GX_{1},gx_{1}x_{2})+B(x_{2}\otimes \
1_{H};GX_{1}X_{2},gx_{1}x_{2})=0  \label{x2ot1,fifthythree}
\end{equation}

\subsubsection{$G^{a}X_{1}X_{2}\otimes g^{d}x_{2}$}

\begin{equation*}
\sum_{\substack{ a,d=0  \\ a+d\equiv 0\text{ }}}^{1}\left( -1\right)
^{a}B(g\otimes 1_{H};G^{a}X_{1}X_{2},g^{d}x_{2})G^{a}X_{1}X_{2}\otimes
g^{d}x_{2}=0
\end{equation*}%
and we get%
\begin{equation*}
B(g\otimes 1_{H};X_{1}X_{2},x_{2})=0
\end{equation*}%
\begin{equation*}
-B(g\otimes 1_{H};GX_{1}X_{2},gx_{2})=0
\end{equation*}%
which are already known.

\subsubsection{$G^{a}X_{1}X_{2}\otimes g^{d}x_{1}$}

\begin{gather*}
\sum_{\substack{ a,d=0  \\ a+d\equiv 0\text{ }}}^{1}\left[ \left( -1\right)
^{a}B(g\otimes 1_{H};G^{a}X_{1}X_{2},g^{d}x_{1})+\left( -1\right)
^{a}B(x_{2}\otimes \ 1_{H};G^{a}X_{1}X_{2},g^{d}x_{1}x_{2})\right] \\
G^{a}X_{1}X_{2}\otimes g^{d}x_{1}=0
\end{gather*}%
and we get%
\begin{equation}
B(g\otimes 1_{H};X_{1}X_{2},x_{1})+B(x_{2}\otimes \
1_{H};X_{1}X_{2},x_{1}x_{2})=0  \label{x2ot1,fifthyfour}
\end{equation}%
\begin{equation}
-B(g\otimes 1_{H};GX_{1}X_{2},gx_{1})-B(x_{2}\otimes \
1_{H};GX_{1}X_{2},gx_{1}x_{2})=0  \label{x2ot1,fifthyfive}
\end{equation}

\subsubsection{$G^{a}X_{1}X_{2}\otimes g^{d}x_{1}x_{2}$}

\begin{equation*}
\sum_{\substack{ a,d=0  \\ a+d\equiv 1\text{ }}}^{1}\left( -1\right)
^{a}B(g\otimes 1_{H};G^{a}X_{1}X_{2},g^{d}x_{1}x_{2})G^{a}X_{1}X_{2}\otimes
g^{d}x_{1}x_{2}=0
\end{equation*}%
and we get%
\begin{equation*}
B(g\otimes 1_{H};X_{1}X_{2},gx_{1}x_{2})=0
\end{equation*}%
which are already known.%
\begin{equation*}
-B(g\otimes 1_{H};GX_{1}X_{2},x_{1}x_{2})=0
\end{equation*}

\subsection{Case $gx_{1}x_{2}$}

For the first summand we get{\Huge \ }%
\begin{eqnarray*}
&&B^{A}(g\otimes 1_{H})_{0}\otimes B^{H}(g\otimes 1_{H})_{1}\otimes
B^{H}(g\otimes 1_{H})_{2}x_{2}B^{A}(g\otimes 1_{H})_{1} \\
&=&\sum_{a,b_{1},b_{2},d,e_{1},e_{2}=0}^{1}\sum_{l_{1}=0}^{b_{1}}%
\sum_{l_{2}=0}^{b_{2}}\sum_{u_{1}=0}^{e_{1}}\sum_{u_{2}=0}^{e_{2}}\left(
-1\right) ^{\alpha }B(g\otimes
1_{H};G^{a}X_{1}^{b_{1}}X_{2}^{b_{2}},g^{d}x_{1}^{e_{1}}x_{2}^{e_{2}}) \\
&&G^{a}X_{1}^{b_{1}-l_{1}}X_{2}^{b_{2}-l_{2}}\otimes g^{d}x_{1}^{\left(
e_{1}-u_{1}\right) }x_{2}^{\left( e_{2}-u_{2}\right) }\otimes
g^{a+b_{1}+b_{2}+d+e_{1}+e_{2}+l_{1}+l_{2}+u_{1}+u_{2}}x_{1}^{u_{1}+l_{1}}x_{2}^{l_{2}+u_{2}+1}
\end{eqnarray*}%
\begin{eqnarray*}
a+b_{1}+b_{2}+d+e_{1}+e_{2}+l_{1}+l_{2}+u_{1}+u_{2} &\equiv &1 \\
u_{1}+l_{1} &=&1 \\
l_{2} &=&u_{2}=0
\end{eqnarray*}%
so that%
\begin{eqnarray*}
a+b_{1}+b_{2}+d+e_{1}+e_{2} &\equiv &0 \\
u_{1}+l_{1} &=&1 \\
l_{2} &=&u_{2}=0
\end{eqnarray*}

where%
\begin{equation*}
\alpha =(a+b_{1}+b_{2}+l_{1})(u_{1}+1)+l_{1}\left( b_{2}+1\right) +e_{2}u_{1}
\end{equation*}%
so that we get%
\begin{eqnarray*}
&&\sum_{\substack{ a,b_{1},b_{2},d,e_{1},e_{2}=0  \\ %
a+b_{1}+b_{2}+d+e_{1}+e_{2}\equiv 0}}^{1}\sum_{l_{1}=0}^{b_{1}}\sum
_{\substack{ u_{1}=0  \\ l_{1}+u_{1}=1}}^{e_{1}}\left( -1\right)
^{(a+b_{1}+b_{2}+l_{1})(u_{1}+1)+l_{1}\left( b_{2}+1\right) +e_{2}u_{1}} \\
&&B(g\otimes
1_{H};G^{a}X_{1}^{b_{1}}X_{2}^{b_{2}},g^{d}x_{1}^{e_{1}}x_{2}^{e_{2}})G^{a}X_{1}^{b_{1}-l_{1}}X_{2}^{b_{2}}\otimes g^{d}x_{1}^{\left( e_{1}-u_{1}\right) }x_{2}^{e_{2}}
\end{eqnarray*}%
i.e.%
\begin{eqnarray*}
&&\sum_{_{\substack{ a,b_{1},b_{2},d,e_{1},e_{2}=0  \\ %
a+b_{1}+b_{2}+d+e_{1}+e_{2}\equiv 0}}}^{1}\sum_{l_{1}=0}^{b_{1}}\sum
_{\substack{ u_{1}=0  \\ l_{1}+u_{1}=1}}^{e_{1}}\left( -1\right)
^{(a+b_{1}+b_{2}+l_{1})(u_{1}+1)+l_{1}\left( b_{2}+1\right) +e_{2}u_{1}} \\
&&B(g\otimes
1_{H};G^{a}X_{1}^{b_{1}}X_{2}^{b_{2}},g^{d}x_{1}^{e_{1}}x_{2}^{e_{2}})G^{a}X_{1}^{b_{1}-l_{1}}X_{2}^{b_{2}}\otimes g^{d}x_{1}^{\left( e_{1}-u_{1}\right) }x_{2}^{e_{2}}
\\
&=&\sum_{_{\substack{ a,b_{1},b_{2},d,e_{2}=0  \\ a+b_{1}+b_{2}+d+e_{2}%
\equiv 1}}}^{1}\left( -1\right) ^{e_{2}}B(g\otimes
1_{H};G^{a}X_{1}^{b_{1}}X_{2}^{b_{2}},g^{d}x_{1}x_{2}^{e_{2}})G^{a}X_{1}^{b_{1}}X_{2}^{b_{2}}\otimes g^{d}x_{2}^{e_{2}}+
\\
&&\sum_{_{\substack{ a,b_{2},d,e_{1},e_{2}=0  \\ a+b_{2}+d+e_{1}+e_{2}\equiv
1 }}}^{1}\left( -1\right) ^{(a+1)}B(g\otimes
1_{H};G^{a}X_{1}X_{2}^{b_{2}},g^{d}x_{1}^{e_{1}}x_{2}^{e_{2}})G^{a}X_{2}^{b_{2}}\otimes g^{d}x_{1}^{e_{1}}x_{2}^{e_{2}}
\end{eqnarray*}

For the second summand, we have%
\begin{eqnarray*}
&&B^{A}(x_{2}\otimes \ 1_{H})_{0}\otimes B^{H}(x_{2}\otimes \
1_{H})_{1}\otimes B^{H}(x_{2}\otimes \ 1_{H})_{2}B^{A}(x_{2}\otimes \
1_{H})_{1} \\
&=&\sum_{a,b_{1},b_{2},d,e_{1},e_{2}=0}^{1}\sum_{l_{1}=0}^{b_{1}}%
\sum_{l_{2}=0}^{b_{2}}\sum_{u_{1}=0}^{e_{1}}\sum_{u_{2}=0}^{e_{2}}(-1)^{%
\beta }B(x_{2}\otimes \
1_{H};G^{a}X_{1}^{b_{1}}X_{2}^{b_{2}},g^{d}x_{1}^{e_{1}}x_{2}^{e_{2}}) \\
&&G^{a}X_{1}^{b_{1}-l_{1}}X_{2}^{b_{2}-l_{2}}\otimes g^{d}x_{1}^{\left(
e_{1}-u_{1}\right) }x_{2}^{\left( e_{2}-u_{2}\right) }\otimes
g^{a+b_{1}+b_{2}+d+e_{1}+e_{2}+l_{1}+l_{2}+u_{1}+u_{2}}x_{1}^{u_{1}+l_{1}}x_{2}^{u_{2}+l_{2}}
\end{eqnarray*}

where $\beta =\left( l_{2}+b_{2}\right) l_{1}+\left( u_{2}+e_{2}\right)
u_{1}+\left( a+b_{1}+b_{2}+l_{1}+l_{2}\right) \left( u_{1}+u_{2}\right)
+l_{1}u_{2}$

\begin{eqnarray*}
a+b_{1}+b_{2}+d+e_{1}+e_{2}+l_{1}+l_{2}+u_{1}+u_{2} &\equiv &1 \\
l_{1}+u_{1} &=&1 \\
l_{2}+u_{2} &=&1
\end{eqnarray*}%
so that%
\begin{equation*}
a+b_{1}+b_{2}+d+e_{1}+e_{2}\equiv 1
\end{equation*}%
\begin{gather*}
\sum_{\substack{ a,b_{1},b_{2},d,e_{1},e_{2}=0  \\ %
a+b_{1}+b_{2}+d+e_{1}+e_{2}\equiv 1}}^{1}\sum_{l_{1}=0}^{b_{1}}%
\sum_{l_{2}=0}^{b_{2}}\sum_{\substack{ u_{1}=0  \\ l_{1}+u_{1}=1}}%
^{e_{1}}\sum _{\substack{ u_{2}=0  \\ l_{2}+u_{2}=1}}^{e_{2}} \\
(-1)^{\left( l_{2}+b_{2}\right) l_{1}+\left( u_{2}+e_{2}\right) u_{1}+\left(
a+b_{1}+b_{2}+l_{1}+l_{2}\right) \left( u_{1}+u_{2}\right) +l_{1}u_{2}} \\
B(x_{2}\otimes \
1_{H};G^{a}X_{1}^{b_{1}}X_{2}^{b_{2}},g^{d}x_{1}^{e_{1}}x_{2}^{e_{2}})G^{a}X_{1}^{b_{1}-l_{1}}X_{2}^{b_{2}-l_{2}}\otimes g^{d}x_{1}^{\left( e_{1}-u_{1}\right) }x_{2}^{\left( e_{2}-u_{2}\right) }=
\\
\sum_{\substack{ a,b_{1},b_{2},d=0  \\ a+b_{1}+b_{2}+d\equiv 1}}%
^{1}B(x_{2}\otimes \
1_{H};G^{a}X_{1}^{b_{1}}X_{2}^{b_{2}},g^{d}x_{1}x_{2})G^{a}X_{1}^{b_{1}}X_{2}^{b_{2}}\otimes g^{d}+
\\
\sum_{\substack{ a,b_{1}d,e_{2}=0  \\ a+b_{1}+d+e_{2}\equiv 1}}%
^{1}(-1)^{e_{2}+\left( a+b_{1}\right) }B(x_{2}\otimes \
1_{H};G^{a}X_{1}^{b_{1}}X_{2},g^{d}x_{1}x_{2}^{e_{2}})G^{a}X_{1}^{b_{1}}%
\otimes g^{d}x_{2}^{e_{2}}+ \\
\sum_{\substack{ a,b_{2},d,e_{1}=0  \\ a+b_{2}+d+e_{1}\equiv 1}}%
^{1}(-1)^{a+1}B(x_{2}\otimes \
1_{H};G^{a}X_{1}X_{2}^{b_{2}},g^{d}x_{1}^{e_{1}}x_{2})G^{a}X_{2}^{b_{2}}%
\otimes g^{d}x_{1}^{e_{1}}+ \\
+\sum_{\substack{ a,d,e_{1},e_{2}=0  \\ a+d+e_{1}+e_{2}\equiv 1}}%
^{1}B(x_{2}\otimes \
1_{H};G^{a}X_{1}X_{2},g^{d}x_{1}^{e_{1}}x_{2}^{e_{2}})G^{a}\otimes
g^{d}x_{1}^{e_{1}}x_{2}^{e_{2}}+
\end{gather*}%
thus we get%
\begin{gather*}
\sum_{_{\substack{ a,b_{1},b_{2},d,e_{2}=0  \\ a+b_{1}+b_{2}+d+e_{2}\equiv 1
}}}^{1}\left( -1\right) ^{e_{2}}B(g\otimes
1_{H};G^{a}X_{1}^{b_{1}}X_{2}^{b_{2}},g^{d}x_{1}x_{2}^{e_{2}})G^{a}X_{1}^{b_{1}}X_{2}^{b_{2}}\otimes g^{d}x_{2}^{e_{2}}+
\\
\sum_{_{\substack{ a,b_{2},d,e_{1},e_{2}=0  \\ a+b_{2}+d+e_{1}+e_{2}\equiv 1
}}}^{1}\left( -1\right) ^{(a+1)}B(g\otimes
1_{H};G^{a}X_{1}X_{2}^{b_{2}},g^{d}x_{1}^{e_{1}}x_{2}^{e_{2}})G^{a}X_{2}^{b_{2}}\otimes g^{d}x_{1}^{e_{1}}x_{2}^{e_{2}}
\\
+\sum_{\substack{ a,b_{1},b_{2},d=0  \\ a+b_{1}+b_{2}+d\equiv 1}}%
^{1}B(x_{2}\otimes \
1_{H};G^{a}X_{1}^{b_{1}}X_{2}^{b_{2}},g^{d}x_{1}x_{2})G^{a}X_{1}^{b_{1}}X_{2}^{b_{2}}\otimes g^{d}+
\\
\sum_{\substack{ a,b_{1}d,e_{2}=0  \\ a+b_{1}+d+e_{2}\equiv 1}}%
^{1}(-1)^{+e_{2}+\left( a+b_{1}\right) }B(x_{2}\otimes \
1_{H};G^{a}X_{1}^{b_{1}}X_{2},g^{d}x_{1}x_{2}^{e_{2}})G^{a}X_{1}^{b_{1}}%
\otimes g^{d}x_{2}^{e_{2}}+ \\
\sum_{\substack{ a,b_{2},d,e_{1}=0  \\ a+b_{2}+d+e_{1}\equiv 1}}%
^{1}(-1)^{a+1}B(x_{2}\otimes \
1_{H};G^{a}X_{1}X_{2}^{b_{2}},g^{d}x_{1}^{e_{1}}x_{2})G^{a}X_{2}^{b_{2}}%
\otimes g^{d}x_{1}^{e_{1}}+ \\
\sum_{\substack{ a,d,e_{1},e_{2}=0  \\ a+d+e_{1}+e_{2}\equiv 1}}%
^{1}B(x_{2}\otimes \
1_{H};G^{a}X_{1}X_{2},g^{d}x_{1}^{e_{1}}x_{2}^{e_{2}})G^{a}\otimes
g^{d}x_{1}^{e_{1}}x_{2}^{e_{2}}+
\end{gather*}

\subsubsection{$G^{a}\otimes g^{d}$}

\begin{equation*}
\sum_{\substack{ a,d  \\ a+d\equiv 1}}^{1}\left[
\begin{array}{c}
B(g\otimes 1_{H};G^{a},g^{d}x_{1})+\left( -1\right) ^{(a+1)}B(g\otimes
1_{H};G^{a}X_{1},g^{d})+ \\
B(x_{2}\otimes \ 1_{H};G^{a},g^{d}x_{1}x_{2})+(-1)^{a}B(x_{2}\otimes \
1_{H};G^{a}X_{2},g^{d}x_{1})+ \\
(-1)^{a+1}B(x_{2}\otimes \ 1_{H};G^{a}X_{1},g^{d}x_{2})+B(x_{2}\otimes \
1_{H};G^{a}X_{1}X_{2},g^{d})%
\end{array}%
\right] G^{a}\otimes g^{d}=0
\end{equation*}%
and we get%
\begin{equation}
\begin{array}{c}
B(g\otimes 1_{H};1_{A},gx_{1})-B(g\otimes 1_{H};X_{1},g)+ \\
B(x_{2}\otimes \ 1_{H};1_{A},gx_{1}x_{2})+B(x_{2}\otimes \
1_{H};X_{2},gx_{1})+ \\
-B(x_{2}\otimes \ 1_{H};X_{1},gx_{2})+B(x_{2}\otimes \ 1_{H};X_{1}X_{2},g)%
\end{array}%
=0  \label{x2ot1,fifthysix}
\end{equation}%
\begin{equation}
\begin{array}{c}
B(g\otimes 1_{H};G,x_{1})+B(g\otimes 1_{H};GX_{1},1_{H})+ \\
B(x_{2}\otimes \ 1_{H};G,x_{1}x_{2})-B(x_{2}\otimes \ 1_{H};GX_{2},x_{1})+
\\
+B(x_{2}\otimes \ 1_{H};GX_{1},x_{2})+B(x_{2}\otimes \
1_{H};GX_{1}X_{2},1_{H})%
\end{array}%
=0  \label{x2ot1,fifthyseven}
\end{equation}

\subsubsection{$G^{a}\otimes g^{d}x_{2}$}

\begin{equation*}
\sum_{\substack{ a,d  \\ a+d\equiv 0}}^{1}\left[
\begin{array}{c}
-B(g\otimes 1_{H};G^{a},g^{d}x_{1}x_{2})+\left( -1\right) ^{(a+1)}B(g\otimes
1_{H};G^{a}X_{1},g^{d}x_{2})+ \\
+(-1)^{a+1}B(x_{2}\otimes \ 1_{H};G^{a}X_{2},g^{d}x_{1}x_{2})+B(x_{2}\otimes
\ 1_{H};G^{a}X_{1}X_{2},g^{d}x_{2})%
\end{array}%
\right] G^{a}\otimes g^{d}x_{2}+
\end{equation*}%
and we get%
\begin{equation}
\begin{array}{c}
-B(g\otimes 1_{H};1_{A},x_{1}x_{2})-B(g\otimes 1_{H};X_{1},x_{2})+ \\
-B(x_{2}\otimes \ 1_{H};X_{2},x_{1}x_{2})+B(x_{2}\otimes \
1_{H};X_{1}X_{2},x_{2})%
\end{array}%
=0  \label{x2ot1,fifthyeight}
\end{equation}%
\begin{equation}
\begin{array}{c}
-B(g\otimes 1_{H};G,gx_{1}x_{2})+B(g\otimes 1_{H};GX_{1},gx_{2})+ \\
+B(x_{2}\otimes \ 1_{H};GX_{2},gx_{1}x_{2})+B(x_{2}\otimes \
1_{H};GX_{1}X_{2},gx_{2})%
\end{array}%
=0  \label{x2ot1,fifthynine}
\end{equation}

\subsubsection{$G^{a}\otimes g^{d}x_{1}$}

\begin{gather*}
\sum_{\substack{ a,d  \\ a+d\equiv 0}}^{1}\left[
\begin{array}{c}
\left( -1\right) ^{(a+1)}B(g\otimes
1_{H};G^{a}X_{1},g^{d}x_{1})+(-1)^{a+1}B(x_{2}\otimes \
1_{H};G^{a}X_{1},g^{d}x_{1}x_{2})+ \\
B(x_{2}\otimes \ 1_{H};G^{a}X_{1}X_{2},g^{d}x_{1})%
\end{array}%
\right] \\
G^{a}\otimes g^{d}x_{1}=0
\end{gather*}%
and we get%
\begin{equation}
-B(g\otimes 1_{H};X_{1},x_{1})-B(x_{2}\otimes \
1_{H};X_{1},x_{1}x_{2})+B(x_{2}\otimes \ 1_{H};X_{1}X_{2},x_{1})=0
\label{x2ot1,sixty}
\end{equation}%
\begin{equation}
B(g\otimes 1_{H};GX_{1},gx_{1})+B(x_{2}\otimes \
1_{H};GX_{1},gx_{1}x_{2})+B(x_{2}\otimes \ 1_{H};GX_{1}X_{2},gx_{1})=0
\label{x2ot1,sixtyone}
\end{equation}

\subsubsection{$G^{a}X_{2}\otimes g^{d}$}

\begin{equation*}
\sum_{\substack{ a,d  \\ a+d\equiv 0}}^{1}\left[
\begin{array}{c}
\left( -1\right) ^{e_{2}}B(g\otimes 1_{H};G^{a}X_{2},g^{d}x_{1})+\left(
-1\right) ^{(a+1)}B(g\otimes 1_{H};G^{a}X_{1}X_{2},g^{d})+ \\
+B(x_{2}\otimes \ 1_{H};G^{a}X_{2},g^{d}x_{1}x_{2})+(-1)^{a+1}B(x_{2}\otimes
\ 1_{H};G^{a}X_{1}X_{2},g^{d}x_{2})%
\end{array}%
\right] G^{a}X_{2}\otimes g^{d}=0
\end{equation*}

and we get%
\begin{equation}
\begin{array}{c}
B(g\otimes 1_{H};X_{2},x_{1})-B(g\otimes 1_{H};X_{1}X_{2},1_{H})+ \\
+B(x_{2}\otimes \ 1_{H};X_{2},x_{1}x_{2})-B(x_{2}\otimes \
1_{H};X_{1}X_{2},x_{2})%
\end{array}%
=0  \label{x2ot1,sixtytwo}
\end{equation}%
\begin{equation}
\begin{array}{c}
B(g\otimes 1_{H};GX_{2},gx_{1})+B(g\otimes 1_{H};GX_{1}X_{2},g)+ \\
+B(x_{2}\otimes \ 1_{H};GX_{2},gx_{1}x_{2})+B(x_{2}\otimes \
1_{H};GX_{1}X_{2},gx_{2})%
\end{array}%
=0  \label{x2ot1,sixtythree}
\end{equation}

\subsubsection{$G^{a}X_{1}\otimes g^{d}$}

\begin{equation*}
\sum_{\substack{ a,d  \\ a+d\equiv 0}}^{1}\left[
\begin{array}{c}
B(g\otimes 1_{H};G^{a}X_{1},g^{d}x_{1})+B(x_{2}\otimes \
1_{H};G^{a}X_{1},g^{d}x_{1}x_{2})+ \\
+(-1)^{a+1}B(x_{2}\otimes \ 1_{H};G^{a}X_{1}X_{2},g^{d}x_{1})%
\end{array}%
\right] G^{a}X_{1}\otimes g^{d}=0
\end{equation*}%
and we get

\subsubsection{$G^{a}\otimes g^{d}x_{1}x_{2}$}

\begin{equation*}
\sum_{\substack{ a,d  \\ a+d\equiv 1}}^{1}\left[ \left( -1\right)
^{(a+1)}B(g\otimes 1_{H};G^{a}X_{1},g^{d}x_{1}x_{2})+B(x_{2}\otimes \
1_{H};G^{a}X_{1}X_{2},g^{d}x_{1}x_{2})\right] G^{a}\otimes g^{d}x_{1}x_{2}=0
\end{equation*}%
and we get%
\begin{equation}
-B(g\otimes 1_{H};X_{1},gx_{1}x_{2})+B(x_{2}\otimes \
1_{H};X_{1}X_{2},gx_{1}x_{2})=0  \label{x2ot1,sixtysix}
\end{equation}%
\begin{equation}
B(g\otimes 1_{H};GX_{1},x_{1}x_{2})+B(x_{2}\otimes \
1_{H};GX_{1}X_{2},x_{1}x_{2})=0  \label{x2ot1,sixtyseven}
\end{equation}

\subsubsection{$G^{a}X_{2}\otimes g^{d}x_{2}$}

\begin{equation*}
\sum_{\substack{ a,d  \\ a+d\equiv 1}}^{1}\left[ -B(g\otimes
1_{H};G^{a}X_{2},g^{d}x_{1}x_{2})+\left( -1\right) ^{(a+1)}B(g\otimes
1_{H};G^{a}X_{1}X_{2},g^{d}x_{2})\right] G^{a}X_{2}\otimes g^{d}x_{2}=0
\end{equation*}%
and we get%
\begin{equation*}
-B(g\otimes 1_{H};X_{2},gx_{1}x_{2})-B(g\otimes 1_{H};X_{1}X_{2},gx_{2})=0
\end{equation*}%
\begin{equation*}
-B(g\otimes 1_{H};GX_{2},x_{1}x_{2})+B(g\otimes 1_{H};GX_{1}X_{2},x_{2})=0
\end{equation*}%
which are already known

\subsubsection{$G^{a}X_{1}\otimes g^{d}x_{2}$}

\begin{equation*}
\sum_{\substack{ a,d  \\ a+d\equiv 1}}^{1}\left[ -B(g\otimes
1_{H};G^{a}X_{1},g^{d}x_{1}x_{2})+(-1)^{a}B(x_{2}\otimes \
1_{H};G^{a}X_{1}X_{2},g^{d}x_{1}x_{2})\right] G^{a}X_{1}\otimes g^{d}x_{2}=0
\end{equation*}%
and we get%
\begin{equation*}
-B(g\otimes 1_{H};X_{1},gx_{1}x_{2})+B(x_{2}\otimes \
1_{H};X_{1}X_{2},gx_{1}x_{2})=0
\end{equation*}%
\begin{equation*}
-B(g\otimes 1_{H};GX_{1},x_{1}x_{2})-B(x_{2}\otimes \
1_{H};GX_{1}X_{2},x_{1}x_{2})=0
\end{equation*}%
which are already known.

\subsubsection{$G^{a}X_{2}\otimes g^{d}x_{1}$}

\begin{gather*}
\sum_{\substack{ a,d  \\ a+d\equiv 1}}^{1}\left[ \left( -1\right)
^{(a+1)}B(g\otimes
1_{H};G^{a}X_{1}X_{2},g^{d}x_{1})+(-1)^{a+1}B(x_{2}\otimes \
1_{H};G^{a}X_{1}X_{2},g^{d}x_{1}x_{2})\right] \\
G^{a}X_{2}\otimes g^{d}x_{1}=0
\end{gather*}%
and we get%
\begin{equation}
-B(g\otimes 1_{H};X_{1}X_{2},gx_{1})-B(x_{2}\otimes \
1_{H};X_{1}X_{2},gx_{1}x_{2})=0  \label{x2ot1,sixtyeight}
\end{equation}%
\begin{equation}
B(g\otimes 1_{H};GX_{1}X_{2},x_{1})+B(x_{2}\otimes \
1_{H};GX_{1}X_{2},x_{1}x_{2})=0  \label{x2ot1,sixtysnine}
\end{equation}

\paragraph{$G^{a}X_{1}\otimes g^{d}x_{1}$}

There is no term like this.

\subsubsection{$G^{a}X_{1}X_{2}\otimes g^{d}$}

\begin{equation*}
\sum_{\substack{ a,d  \\ a+d\equiv 1}}^{1}\left[ B(g\otimes
1_{H};G^{a}X_{1}X_{2},g^{d}x_{1})+B(x_{2}\otimes \
1_{H};G^{a}X_{1}X_{2},g^{d}x_{1}x_{2})\right] G^{a}X_{1}X_{2}\otimes g^{d}=0
\end{equation*}%
and we get%
\begin{equation*}
B(g\otimes 1_{H};X_{1}X_{2},gx_{1})+B(x_{2}\otimes \
1_{H};X_{1}X_{2},gx_{1}x_{2})=0
\end{equation*}%
\begin{equation*}
B(g\otimes 1_{H};GX_{1}X_{2},x_{1})+B(x_{2}\otimes \
1_{H};GX_{1}X_{2},x_{1}x_{2})=0
\end{equation*}%
which are already known.

\subsubsection{$G^{a}X_{2}\otimes g^{d}x_{1}x_{2}$}

\begin{equation*}
\sum_{\substack{ a,d  \\ a+d\equiv 0}}^{1}\left( -1\right)
^{(a+1)}B(g\otimes
1_{H};G^{a}X_{1}X_{2}^{b_{2}},g^{d}x_{1}^{e_{1}}x_{2}^{e_{2}})G^{a}X_{2}%
\otimes g^{d}x_{1}x_{2}=0
\end{equation*}%
and we get%
\begin{equation*}
-B(g\otimes 1_{H};X_{1}X_{2},x_{1}x_{2})=0
\end{equation*}%
\begin{equation*}
B(g\otimes 1_{H};GX_{1}X_{2},gx_{1}x_{2})=0
\end{equation*}%
which are already known

\subsubsection{$G^{a}X_{1}\otimes g^{d}x_{1}x_{2}$}

There is no term like this.

\subsubsection{$G^{a}X_{1}X_{2}\otimes g^{d}x_{2}$}

\begin{equation*}
\sum_{\substack{ a,d  \\ a+d\equiv 0}}^{1}-B(g\otimes
1_{H};G^{a}X_{1}X_{2},g^{d}x_{1}x_{2})G^{a}X_{1}X_{2}\otimes g^{d}x_{2}=0
\end{equation*}%
and we get%
\begin{equation*}
-B(g\otimes 1_{H};X_{1}X_{2},x_{1}x_{2})=0
\end{equation*}%
\begin{equation*}
-B(g\otimes 1_{H};GX_{1}X_{2},gx_{1}x_{2})=0
\end{equation*}%
which are already known.

\subsubsection{$G^{a}X_{1}X_{2}\otimes g^{d}x_{1}$}

There is no term like this.

\subsubsection{$G^{a}X_{1}X_{2}\otimes g^{d}x_{1}x_{2}$}

There is no term like this.

\bigskip By using all the equalities we got above from $\left( \ref{xot2
first}\right) $ to $\left( \ref{x2ot1,sixtysnine}\right) $ we obtain the
following form of $B\left( x_{2}\otimes 1_{H}\right) .$

\subsection{The final form of the element $B\left( x_{2}\otimes 1_{H}\right)
$}

\begin{eqnarray}
B\left( x_{2}\otimes 1_{H}\right) &=&B\left( x_{2}\otimes
1_{H};1_{A},1_{H}\right) 1_{A}\otimes 1_{H}  \label{x2} \\
&&+B\left( x_{2}\otimes 1_{H};1_{A},x_{1}x_{2}\right) 1_{A}\otimes x_{1}x_{2}
\notag \\
&&+B\left( x_{2}\otimes 1_{H};1_{A},gx_{1}\right) 1_{A}\otimes gx_{1}+
\notag \\
&&+B\left( x_{2}\otimes 1_{H};1_{A},gx_{2}\right) 1_{A}\otimes gx_{2}++
\notag \\
&&+B\left( x_{2}\otimes 1_{H};G,g\right) G\otimes g  \notag \\
&&+B\left( x_{2}\otimes 1_{H};G,x_{1}\right) G\otimes x_{1}+  \notag \\
&&+B\left( x_{2}\otimes 1_{H};G,x_{2}\right) G\otimes x_{2}+  \notag \\
&&+B\left( x_{2}\otimes 1_{H};G,gx_{1}x_{2}\right) G\otimes gx_{1}x_{2}+
\notag \\
&&-B(x_{2}\otimes 1_{H};1_{A},gx_{1})X_{1}\otimes g+  \notag \\
&&+B(x_{2}\otimes 1_{H};1_{A},x_{1}x_{2})X_{1}\otimes x_{2}+  \notag \\
&&+\left[ -B(g\otimes 1_{H};1_{A},g)-B(x_{2}\otimes \ 1_{H};1_{A},gx_{2})%
\right] X_{2}\otimes g+  \notag \\
&&+\left[ -B(g\otimes 1_{H};1_{A},x_{1})-B(x_{2}\otimes \
1_{H};1_{A},x_{1}x_{2})\right] X_{2}\otimes x_{1}  \notag \\
&&-B(g\otimes 1_{H};1_{A},x_{2})X_{2}\otimes x_{2}+  \notag \\
&&-B(g\otimes 1_{H};X_{1}X_{2},g)X_{2}\otimes gx_{1}x_{2}+  \notag \\
&&+\left[ B(g\otimes 1_{H};1_{A},x_{1})+B(x_{2}\otimes \
1_{H};1_{A},x_{1}x_{2})\right] X_{1}X_{2}\otimes 1_{H}+  \notag \\
&&-B(g\otimes 1_{H};X_{1}X_{2},g)X_{1}X_{2}\otimes gx_{2}+  \notag \\
&&+B(x_{2}\otimes 1_{H};G,x_{1})GX_{1}\otimes 1_{H}+  \notag \\
&&-B(x_{2}\otimes 1_{H};G,gx_{1}x_{2})GX_{1}\otimes gx_{2}+  \notag \\
&&+\left[ B(g\otimes 1_{H};G,1_{H})+B(x_{2}\otimes \ 1_{H};G,x_{2}\right]
GX_{2}\otimes 1_{H}+  \notag \\
&&+B(g\otimes 1_{H};GX_{1}X_{2},1_{H})GX_{2}\otimes x_{1}x_{2}+  \notag \\
&&+\left[ B(g\otimes 1_{H};G,gx_{1})+B(x_{2}\otimes \ 1_{H};G,gx_{1}x_{2})%
\right] GX_{2}\otimes gx_{1}+  \notag \\
&&+B(g\otimes 1_{H};G,gx_{2})GX_{2}\otimes gx_{2}+  \notag \\
&&+\left[ B(g\otimes 1_{H};G,gx_{1})+B(x_{2}\otimes \ 1_{H};G,gx_{1}x_{2})%
\right] GX_{1}X_{2}\otimes g+  \notag \\
&&-B(g\otimes 1_{H};GX_{1}X_{2},1_{H})GX_{1}X_{2}\otimes x_{2}  \notag
\end{eqnarray}

\section{$B$$\left( gx_{1}\otimes 1\right) $}

We write the Casimir formula $\left( \ref{MAIN FORMULA 1}\right) $ for $B$%
\bigskip $\left( gx_{1}\otimes 1_{H}\right) .$%
\begin{eqnarray*}
&&\sum_{w_{1}=0}^{1}\sum_{a,b_{1},b_{2},d,e_{1},e_{2}=0}^{1}%
\sum_{l_{1}=0}^{b_{1}}\sum_{l_{2}=0}^{b_{2}}\sum_{u_{1}=0}^{e_{1}}%
\sum_{u_{2}=0}^{e_{2}}\left( -1\right) ^{\alpha \left(
gx_{1}^{1-w_{1}};l_{1},l_{2},u_{1},u_{2}\right) } \\
&&B(g^{w_{1}}x_{1}^{w_{1}}\otimes
1_{H};G^{a}X_{1}^{b_{1}}X_{2}^{b_{2}},g^{d}x_{1}^{e_{1}}x_{2}^{e_{2}}) \\
&&G^{a}X_{1}^{b_{1}-l_{1}}X_{2}^{b_{2}-l_{2}}\otimes
g^{d}x_{1}^{e_{1}-u_{1}}x_{2}^{e_{2}-u_{2}}\otimes \\
&&g^{a+b_{1}+b_{2}+l_{1}+l_{2}+d+e_{1}+e_{2}+u_{1}+u_{2}+1}x_{1}^{l_{1}+u_{1}+1-w_{1}}x_{2}^{l_{2}+u_{2}}
\\
&=&B\bigskip \left( gx_{1}\otimes 1_{H}\right) \otimes 1_{H}
\end{eqnarray*}%
and we get%
\begin{eqnarray*}
&&\sum_{a,b_{1},b_{2},d,e_{1},e_{2}=0}^{1}\sum_{l_{1}=0}^{b_{1}}%
\sum_{l_{2}=0}^{b_{2}}\sum_{u_{1}=0}^{e_{1}}\sum_{u_{2}=0}^{e_{2}}\left(
-1\right) ^{\alpha \left( gx_{1};l_{1},l_{2},u_{1},u_{2}\right) } \\
&&B(1_{H}\otimes
1_{H};G^{a}X_{1}^{b_{1}}X_{2}^{b_{2}},g^{d}x_{1}^{e_{1}}x_{2}^{e_{2}}) \\
&&G^{a}X_{1}^{b_{1}-l_{1}}X_{2}^{b_{2}-l_{2}}\otimes
g^{d}x_{1}^{e_{1}-u_{1}}x_{2}^{e_{2}-u_{2}}\otimes
g^{a+b_{1}+b_{2}+l_{1}+l_{2}+d+e_{1}+e_{2}+u_{1}+u_{2}+1}x_{1}^{l_{1}+u_{1}+1}x_{2}^{l_{2}+u_{2}}
\\
&&+\sum_{a,b_{1},b_{2},d,e_{1},e_{2}=0}^{1}\sum_{l_{1}=0}^{b_{1}}%
\sum_{l_{2}=0}^{b_{2}}\sum_{u_{1}=0}^{e_{1}}\sum_{u_{2}=0}^{e_{2}}\left(
-1\right) ^{\alpha \left( g;l_{1},l_{2},u_{1},u_{2}\right) } \\
&&B(gx_{1}\otimes
1_{H};G^{a}X_{1}^{b_{1}}X_{2}^{b_{2}},g^{d}x_{1}^{e_{1}}x_{2}^{e_{2}}) \\
&&G^{a}X_{1}^{b_{1}-l_{1}}X_{2}^{b_{2}-l_{2}}\otimes
g^{d}x_{1}^{e_{1}-u_{1}}x_{2}^{e_{2}-u_{2}}\otimes
g^{a+b_{1}+b_{2}+l_{1}+l_{2}+d+e_{1}+e_{2}+u_{1}+u_{2}+1}x_{1}^{l_{1}+u_{1}}x_{2}^{l_{2}+u_{2}}
\\
&=&B\bigskip \left( gx_{1}\otimes 1_{H}\right) \otimes 1_{H}.
\end{eqnarray*}%
Since $B(1_{H}\otimes 1_{H})=1_{A}\otimes 1_{H}$ we obtain

\begin{eqnarray*}
&&1_{A}\otimes 1_{H}\otimes gx_{1}+ \\
&&+\sum_{a,b_{1},b_{2},d,e_{1},e_{2}=0}^{1}\sum_{l_{1}=0}^{b_{1}}%
\sum_{l_{2}=0}^{b_{2}}\sum_{u_{1}=0}^{e_{1}}\sum_{u_{2}=0}^{e_{2}}\left(
-1\right) ^{\alpha \left( g;l_{1},l_{2},u_{1},u_{2}\right) } \\
&&B(gx_{1}\otimes
1_{H};G^{a}X_{1}^{b_{1}}X_{2}^{b_{2}},g^{d}x_{1}^{e_{1}}x_{2}^{e_{2}}) \\
&&G^{a}X_{1}^{b_{1}-l_{1}}X_{2}^{b_{2}-l_{2}}\otimes
g^{d}x_{1}^{e_{1}-u_{1}}x_{2}^{e_{2}-u_{2}}\otimes
g^{a+b_{1}+b_{2}+l_{1}+l_{2}+d+e_{1}+e_{2}+u_{1}+u_{2}+1}x_{1}^{l_{1}+u_{1}}x_{2}^{l_{2}+u_{2}}
\\
&=&B\bigskip \left( gx_{1}\otimes 1_{H}\right) \otimes 1_{H}.
\end{eqnarray*}

\subsection{Case $1_{H}$}

\begin{equation*}
g^{d+e_{1}+e_{2}+u_{1}+u_{2}+a+b_{1}+b_{2}+l_{1}+l_{2}+1}x_{1}^{u_{1}+l_{1}}x_{2}^{u_{2}+l_{2}}=1_{H}
\end{equation*}%
and hence%
\begin{eqnarray*}
u_{1}+l_{1} &=&0 \\
u_{2}+l_{2} &=&0 \\
d+e_{1}+e_{2}+u_{1}+u_{2}+a+b_{1}+b_{2}+l_{1}+l_{2}+1 &\equiv &0
\end{eqnarray*}%
i.e.%
\begin{eqnarray*}
u_{1} &=&l_{1}=0 \\
u_{2} &=&l_{2}=0 \\
a+b_{1}+b_{2}+d+e_{1}+e_{2} &\equiv &1
\end{eqnarray*}%
\begin{equation*}
\alpha \left( g;0,0,0,0\right) =0
\end{equation*}%
\begin{eqnarray*}
&&\sum_{\substack{ a,b_{1},b_{2},d,e_{1},e_{2}=0  \\ %
a+b_{1}+b_{2}+d+e_{1}+e_{2}\equiv 1}}^{1}B(gx_{1}\otimes
1_{H};G^{a}X_{1}^{b_{1}}X_{2}^{b_{2}},g^{d}x_{1}^{e_{1}}x_{2}^{e_{2}})G^{a}X_{1}^{b_{1}}X_{2}^{b_{2}}\otimes g^{d}x_{1}^{e_{1}}x_{2}^{e_{2}}
\\
&=&\sum_{a,b_{1},b_{2},d,e_{1},e_{2}=0}^{1}B(gx_{1}\otimes
1_{H};G^{a}X_{1}^{b_{1}}X_{2}^{b_{2}},g^{d}x_{1}^{e_{1}}x_{2}^{e_{2}})G^{a}X_{1}^{b_{1}}X_{2}^{b_{2}}\otimes g^{d}x_{1}^{e_{1}}x_{2}^{e_{2}}
\end{eqnarray*}%
which implies that%
\begin{equation}
B(gx_{1}\otimes
1_{H};G^{a}X_{1}^{b_{1}}X_{2}^{b_{2}},g^{d}x_{1}^{e_{1}}x_{2}^{e_{2}})=0%
\text{ for }a+b_{1}+b_{2}+d+e_{1}+e_{2}\equiv 0  \label{gx1ot1, first}
\end{equation}

\subsection{Case $g$}

\begin{equation*}
g^{a+b_{1}+b_{2}+d+e_{1}+e_{2}+u_{1}+u_{2}+l_{1}+l_{2}+1}x_{1}^{u_{1}+l_{1}}x_{2}^{u_{2}+l_{2}}=g
\end{equation*}%
this means%
\begin{eqnarray*}
u_{1}+l_{1} &=&0 \\
u_{2}+l_{2} &=&0 \\
d+e_{1}+e_{2}+u_{1}+u_{2}+a+b_{1}+b_{2}+l_{1}+l_{2}+1 &\equiv &1
\end{eqnarray*}%
\begin{eqnarray*}
u_{1} &=&l_{1}=0 \\
u_{2} &=&l_{2}=0 \\
a+b_{1}+b_{2}+d+e_{1}+e_{2} &\equiv &0
\end{eqnarray*}%
which follows from case $1_{H}$.

\subsection{Case $x_{1}$}

\begin{equation*}
g^{d+e_{1}+e_{2}+u_{1}+u_{2}+a+b_{1}+b_{2}+l_{1}+l_{2}+1}x_{1}^{u_{1}+l_{1}}x_{2}^{u_{2}+l_{2}}=x_{1}
\end{equation*}%
\begin{eqnarray*}
g^{a+b_{1}+b_{2}+d+e_{1}+e_{2}+l_{1}+l_{2}+u_{1}+u_{2}+1}x_{1}^{u_{1}+l_{1}}x_{2}^{l_{2}+u_{2}} &=&x_{1}\Rightarrow
\\
l_{2} &=&u_{2}=0 \\
u_{1}+l_{1} &=&1 \\
d+e_{1}+e_{2}+u_{1}+u_{2}+a+b_{1}+b_{2}+l_{1}+l_{2}+1 &\equiv &0 \\
a+b_{1}+b_{2}+d+e_{1}+e_{2} &\equiv &0\text{ }
\end{eqnarray*}%
which follows from case $1_{H}$.

\subsection{Case $x_{2}$}

\begin{equation*}
g^{d+e_{1}+e_{2}+u_{1}+u_{2}+a+b_{1}+b_{2}+l_{1}+l_{2}+1}x_{1}^{u_{1}+l_{1}}x_{2}^{u_{2}+l_{2}}=x_{2}
\end{equation*}%
\begin{eqnarray*}
u_{1} &=&l_{1}=0 \\
u_{2}+l_{2} &=&1 \\
d+e_{1}+e_{2}+u_{1}+u_{2}+a+b_{1}+b_{2}+l_{1}+l_{2}+1 &\equiv &0
\end{eqnarray*}%
\begin{eqnarray*}
u_{1} &=&l_{1}=0 \\
u_{2}+l_{2} &=&1 \\
a+b_{1}+b_{2}+d+e_{1}+e_{2} &\equiv &0\text{ }
\end{eqnarray*}%
which follows from case $1_{H}$.

\subsection{Case $x_{1}x_{2}$}

\begin{eqnarray*}
&&\sum_{a,b_{1},b_{2},d,e_{1},e_{2}=0}^{1}\sum_{l_{1}=0}^{b_{1}}%
\sum_{l_{2}=0}^{b_{2}}\sum_{u_{1}=0}^{e_{1}}\sum_{u_{2}=0}^{e_{2}}\left(
-1\right) ^{\alpha \left( g;l_{1},l_{2},u_{1},u_{2}\right) } \\
&&B(gx_{1}\otimes
1_{H};G^{a}X_{1}^{b_{1}}X_{2}^{b_{2}},g^{d}x_{1}^{e_{1}}x_{2}^{e_{2}}) \\
&&G^{a}X_{1}^{b_{1}-l_{1}}X_{2}^{b_{2}-l_{2}}\otimes
g^{d}x_{1}^{e_{1}-u_{1}}x_{2}^{e_{2}-u_{2}}\otimes
g^{a+b_{1}+b_{2}+l_{1}+l_{2}+d+e_{1}+e_{2}+u_{1}+u_{2}+1}x_{1}^{l_{1}+u_{1}}x_{2}^{l_{2}+u_{2}}
\end{eqnarray*}

\begin{equation*}
g^{d+e_{1}+e_{2}+u_{1}+u_{2}+a+b_{1}+b_{2}+l_{1}+l_{2}+1}x_{1}^{u_{1}+l_{1}}x_{2}^{u_{2}+l_{2}}=x_{1}x_{2}
\end{equation*}%
\begin{eqnarray*}
u_{1}+l_{1} &=&1 \\
u_{2}+l_{2} &=&1 \\
d+e_{1}+e_{2}+u_{1}+u_{2}+a+b_{1}+b_{2}+l_{1}+l_{2}+1 &\equiv &0
\end{eqnarray*}%
\begin{eqnarray*}
u_{1}+l_{1} &=&1 \\
u_{2}+l_{2} &=&1 \\
a+b_{1}+b_{2}+d+e_{1}+e_{2} &\equiv &1
\end{eqnarray*}%
we get

\begin{eqnarray*}
&&\sum_{\substack{ a,b_{1},b_{2},d=0  \\ a+b_{1}+b_{2}+d\equiv 1}}^{1}\left(
-1\right) ^{\alpha \left( g;0,0,1,1\right) }B(gx_{1}\otimes
1_{H};G^{a}X_{1}^{b_{1}}X_{2}^{b_{2}},g^{d}x_{1}x_{2})G^{a}X_{1}^{b_{1}}X_{2}^{b_{2}}\otimes g^{d}
\\
&&\sum_{\substack{ a,b_{1},d,e_{2}=0  \\ a+b_{1}+d+e_{2}\equiv 1}}^{1}\left(
-1\right) ^{\alpha \left( g;0,1,1,0\right) }B(gx_{1}\otimes
1_{H};G^{a}X_{1}^{b_{1}}X_{2},g^{d}x_{1}x_{2}^{e_{2}})G^{a}X_{1}^{b_{1}}%
\otimes g^{d}x_{2}^{e_{2}} \\
&&\sum_{\substack{ a,,d,e_{1},e_{2}=0  \\ a+d+e_{1}+e_{2}\equiv 1}}%
^{1}\left( -1\right) ^{\alpha \left( g;1,1,0,0\right) }B(gx_{1}\otimes
1_{H};G^{a}X_{1}X_{2},g^{d}x_{1}^{e_{1}}x_{2}^{e_{2}})G\otimes
g^{d}x_{1}^{e_{1}}x_{2}^{e_{2}} \\
&&\sum_{\substack{ a,b_{2},d,e_{1}=0  \\ a+b_{2}+d+e_{1}\equiv 1}}^{1}\left(
-1\right) ^{\alpha \left( g;1,0,0,1\right) }B(gx_{1}\otimes
1_{H};G^{a}X_{1}X_{2}^{b_{2}},g^{d}x_{1}^{e_{1}}x_{2})G^{a}X_{2}^{b_{2}}%
\otimes g^{d}x_{1}^{e_{1}}
\end{eqnarray*}%
since $\alpha \left( g;0,0,1,1\right) =1+e_{2}$, $\alpha \left(
g;0,1,1,0\right) =e_{2}+a+b_{1}+b_{2}$,

$\alpha \left( g;1,1,0,0\right) =1+b_{2}$, $\alpha \left( g;1,0,0,1\right)
=\alpha \left( g^{m}x_{1}^{n_{1}}x_{2}^{n_{2}};1,0,0,1\right) =a+b_{1}+1$

\begin{eqnarray*}
&&\sum_{\substack{ a,b_{1},b_{2},d=0  \\ a+b_{1}+b_{2}+d\equiv 1}}%
^{1}B(gx_{1}\otimes
1_{H};G^{a}X_{1}^{b_{1}}X_{2}^{b_{2}},g^{d}x_{1}x_{2})G^{a}X_{1}^{b_{1}}X_{2}^{b_{2}}\otimes g^{d}
\\
&&\sum_{\substack{ a,b_{1},d,e_{2}=0  \\ a+b_{1}+d+e_{2}\equiv 1}}^{1}\left(
-1\right) ^{e_{2}+a+b_{1}+1}B(gx_{1}\otimes
1_{H};G^{a}X_{1}^{b_{1}}X_{2},g^{d}x_{1}x_{2}^{e_{2}})G^{a}X_{1}^{b_{1}}%
\otimes g^{d}x_{2}^{e_{2}} \\
&&\sum_{\substack{ a,,d,e_{1},e_{2}=0  \\ a+d+e_{1}+e_{2}\equiv 1}}%
^{1}B(gx_{1}\otimes
1_{H};G^{a}X_{1}X_{2},g^{d}x_{1}^{e_{1}}x_{2}^{e_{2}})G\otimes
g^{d}x_{1}^{e_{1}}x_{2}^{e_{2}} \\
&&\sum_{\substack{ a,b_{2},d,e_{1}=0  \\ a+b_{2}+d+e_{1}\equiv 1}}^{1}\left(
-1\right) ^{a}B(gx_{1}\otimes
1_{H};G^{a}X_{1}X_{2}^{b_{2}},g^{d}x_{1}^{e_{1}}x_{2})G^{a}X_{2}^{b_{2}}%
\otimes g^{d}x_{1}^{e_{1}}
\end{eqnarray*}

\subsubsection{$G^{a}\otimes g^{d}$}

\begin{equation*}
\sum_{\substack{ a,d=0  \\ a+d\equiv 1}}^{1}\left[
\begin{array}{c}
B(gx_{1}\otimes 1_{H};G^{a},g^{d}x_{1}x_{2})+\left( -1\right)
^{a}B(gx_{1}\otimes 1_{H};G^{a}X_{2},g^{d}x_{1})+ \\
+B(gx_{1}\otimes 1_{H};G^{a}X_{1}X_{2},g^{d})+\left( -1\right)
^{a}B(gx_{1}\otimes 1_{H};G^{a}X_{1},g^{d}x_{2})%
\end{array}%
\right] G^{a}\otimes g^{d}=0
\end{equation*}%
and we get%
\begin{equation}
\begin{array}{c}
B(gx_{1}\otimes 1_{H};,gx_{1}x_{2})-B(gx_{1}\otimes 1_{H};X_{2},gx_{1})+ \\
+B(gx_{1}\otimes 1_{H};X_{1}X_{2},g)+B(gx_{1}\otimes 1_{H};X_{1},gx_{2})%
\end{array}%
=0  \label{gx1ot1, second}
\end{equation}%
\begin{equation}
\begin{array}{c}
B(gx_{1}\otimes 1_{H};G,x_{1}x_{2})+B(gx_{1}\otimes 1_{H};GX_{2},x_{1})+ \\
+B(gx_{1}\otimes 1_{H};GX_{1}X_{2},1_{H})-B(gx_{1}\otimes 1_{H};GX_{1},x_{2})%
\end{array}%
=0  \label{gx1ot1, third}
\end{equation}

\subsubsection{$G^{a}\otimes g^{d}x_{2}$}

\begin{equation*}
\sum_{\substack{ a,d=0  \\ a+d\equiv 0}}^{1}\left[ \left( -1\right)
^{a}B(gx_{1}\otimes 1_{H};G^{a}X_{2},g^{d}x_{1}x_{2})+B(gx_{1}\otimes
1_{H};G^{a}X_{1}X_{2},g^{d}x_{2})\right] G^{a}\otimes g^{d}x_{2}=0
\end{equation*}%
and we get%
\begin{equation}
B(gx_{1}\otimes 1_{H};X_{2},x_{1}x_{2})+B(gx_{1}\otimes
1_{H};X_{1}X_{2},x_{2})=0  \label{gx1ot1, four}
\end{equation}%
\begin{equation}
-B(gx_{1}\otimes 1_{H};GX_{2},gx_{1}x_{2})+B(gx_{1}\otimes
1_{H};GX_{1}X_{2},gx_{2})=0  \label{gx1ot1, five}
\end{equation}

\subsubsection{$G^{a}\otimes g^{d}x_{1}$}

\begin{gather*}
\sum_{\substack{ a,d=0  \\ a+d\equiv 0}}^{1}\left[ B(gx_{1}\otimes
1_{H};G^{a}X_{1}X_{2},g^{d}x_{1})+\left( -1\right) ^{a}B(gx_{1}\otimes
1_{H};G^{a}X_{1},g^{d}x_{1}x_{2})\right] \\
G^{a}\otimes g^{d}x_{1}=0
\end{gather*}%
and we get%
\begin{equation}
B(gx_{1}\otimes 1_{H};X_{1}X_{2},x_{1})+B(gx_{1}\otimes
1_{H};X_{1},x_{1}x_{2})=0  \label{gx1ot1, six}
\end{equation}%
\begin{equation}
B(gx_{1}\otimes 1_{H};GX_{1}X_{2},gx_{1})-B(gx_{1}\otimes
1_{H};GX_{1},gx_{1}x_{2})=0  \label{gx1ot1, seven}
\end{equation}

\subsubsection{$G^{a}X_{2}\otimes g^{d}$}

\begin{gather*}
\sum_{\substack{ a,d=0  \\ a+d\equiv 0}}^{1}\left[ B(gx_{1}\otimes
1_{H};G^{a}X_{2},g^{d}x_{1}x_{2})+\left( -1\right) ^{a}B(gx_{1}\otimes
1_{H};G^{a}X_{1}X_{2},g^{d}x_{2})\right] \\
G^{a}X_{2}\otimes g^{d}=0
\end{gather*}%
and we get%
\begin{equation*}
B(gx_{1}\otimes 1_{H};X_{2},x_{1}x_{2})+B(gx_{1}\otimes
1_{H};X_{1}X_{2},x_{2})=0
\end{equation*}%
\begin{equation*}
B(gx_{1}\otimes 1_{H};GX_{2},gx_{1}x_{2})-B(gx_{1}\otimes
1_{H};GX_{1}X_{2},gx_{2})=0
\end{equation*}%
which we already got.

\subsubsection{$G^{a}X_{1}\otimes g^{d}$}

\begin{gather*}
\sum_{\substack{ a,d=0  \\ a+d\equiv 0}}^{1}\left[ B(gx_{1}\otimes
1_{H};G^{a}X_{1},g^{d}x_{1}x_{2})+\left( -1\right) ^{a}B(gx_{1}\otimes
1_{H};G^{a}X_{1}X_{2},g^{d}x_{1})\right] \\
G^{a}X_{1}\otimes g^{d}=0
\end{gather*}%
and we get%
\begin{equation*}
B(gx_{1}\otimes 1_{H};X_{1},x_{1}x_{2})+B(gx_{1}\otimes
1_{H};X_{1}X_{2},x_{1})=0
\end{equation*}%
\begin{equation*}
B(gx_{1}\otimes 1_{H};GX_{1},gx_{1}x_{2})-B(gx_{1}\otimes
1_{H};GX_{1}X_{2},gx_{1})=0
\end{equation*}%
which we already got.

\subsubsection{$G^{a}\otimes g^{d}x_{1}x_{2}$}

\begin{equation*}
\sum_{\substack{ a,d=0  \\ a+d\equiv 1}}^{1}B(gx_{1}\otimes
1_{H};G^{a}X_{1}X_{2},g^{d}x_{1}x_{2})G^{a}\otimes g^{d}x_{1}x_{2}=0
\end{equation*}%
and we get%
\begin{equation}
B(gx_{1}\otimes 1_{H};X_{1}X_{2},gx_{1}x_{2})=0  \label{gx1ot1, eigth}
\end{equation}%
\begin{equation}
B(gx_{1}\otimes 1_{H};GX_{1}X_{2},x_{1}x_{2})=0  \label{gx1ot1, nine}
\end{equation}

\subsubsection{$G^{a}X_{2}\otimes g^{d}x_{2}$}

We do not have any equality like this.

\subsubsection{$G^{a}X_{1}\otimes g^{d}x_{2}$}

\begin{equation*}
\sum_{\substack{ a,d=0  \\ a+d\equiv 1}}^{1}\left( -1\right)
^{a+1}B(gx_{1}\otimes
1_{H};G^{a}X_{1}X_{2},g^{d}x_{1}x_{2})G^{a}X_{1}\otimes g^{d}x_{2}=0
\end{equation*}%
and we get%
\begin{equation*}
-B(gx_{1}\otimes 1_{H};X_{1}X_{2},gx_{1}x_{2})=0
\end{equation*}%
\begin{equation*}
B(gx_{1}\otimes 1_{H};GX_{1}X_{2},x_{1}x_{2})=0
\end{equation*}%
which are already known.

Going on this way, we get no new information.

\subsubsection{Case $gx_{1}$}

\begin{eqnarray*}
&&1_{A}\otimes 1_{H}\otimes gx_{1} \\
&&+\sum_{a,b_{1},b_{2},d,e_{1},e_{2}=0}^{1}\sum_{l_{1}=0}^{b_{1}}%
\sum_{l_{2}=0}^{b_{2}}\sum_{u_{1}=0}^{e_{1}}\sum_{u_{2}=0}^{e_{2}}\left(
-1\right) ^{\alpha \left( g;l_{1},l_{2},u_{1},u_{2}\right) } \\
&&B(gx_{1}\otimes
1_{H};G^{a}X_{1}^{b_{1}}X_{2}^{b_{2}},g^{d}x_{1}^{e_{1}}x_{2}^{e_{2}}) \\
&&G^{a}X_{1}^{b_{1}-l_{1}}X_{2}^{b_{2}-l_{2}}\otimes
g^{d}x_{1}^{e_{1}-u_{1}}x_{2}^{e_{2}-u_{2}}\otimes
g^{a+b_{1}+b_{2}+l_{1}+l_{2}+d+e_{1}+e_{2}+u_{1}+u_{2}+1}x_{1}^{l_{1}+u_{1}}x_{2}^{l_{2}+u_{2}}
\\
&=&B\bigskip \left( gx_{1}\otimes 1_{H}\right) \otimes 1_{H}.
\end{eqnarray*}%
\begin{eqnarray*}
l_{1}+u_{1} &=&1 \\
l_{2} &=&u_{2}=0 \\
d+e_{1}+e_{2}+u_{1}+u_{2}+a+b_{1}+b_{2}+l_{1}+l_{2}+1 &\equiv &1
\end{eqnarray*}%
\begin{eqnarray*}
l_{1}+u_{1} &=&1 \\
l_{2} &=&u_{2}=0 \\
a+b_{1}+b_{2}+d+e_{1}+e_{2} &\equiv &1\text{ }
\end{eqnarray*}%
and we get

\begin{eqnarray*}
&&1_{A}\otimes 1_{H}+ \\
&&+\sum_{\substack{ a,b_{1},b_{2},d,e_{1},e_{2}=0  \\ %
a+b_{1}+b_{2}+d+e_{1}+e_{2}\equiv 1\text{ }}}^{1}\sum_{l_{1}=0}^{b_{1}}\sum
_{\substack{ u_{1}=0  \\ l_{1}+u_{1}=1}}^{e_{1}}\left( -1\right) ^{\alpha
\left( g;l_{1},0,u_{1},0\right) }B(gx_{1}\otimes
1_{H};G^{a}X_{1}^{b_{1}}X_{2}^{b_{2}},g^{d}x_{1}^{e_{1}}x_{2}^{e_{2}}) \\
&&G^{a}X_{1}^{b_{1}-l_{1}}X_{2}^{b_{2}}\otimes
g^{d}x_{1}^{e_{1}-u_{1}}x_{2}^{e_{2}}=0
\end{eqnarray*}%
and, since $\alpha \left( g;0,0,1,0\right) =e_{2}+a+b_{1}+b_{2}+1$, $\alpha
\left( g;1,0,0,0\right) =b_{2}$ we obtain%
\begin{eqnarray*}
&&1_{A}\otimes 1_{H} \\
&&+\sum_{\substack{ a,b_{1},b_{2},d,,e_{2}=0  \\ a+b_{1}+b_{2}+d+e_{2}\equiv
0\text{ }}}^{1}\left( -1\right) ^{e_{2}+a+b_{1}+b_{2}+1}B(gx_{1}\otimes
1_{H};G^{a}X_{1}^{b_{1}}X_{2}^{b_{2}},g^{d}x_{1}x_{2}^{e_{2}})G^{a}X_{1}^{b_{1}}X_{2}^{b_{2}}\otimes g^{d}x_{2}^{e_{2}}+
\\
&&\sum_{\substack{ a,b_{1},b_{2},d,e_{1},e_{2}=0  \\ a+b_{2}+d+e_{1}+e_{2}%
\equiv 0\text{ }}}^{1}\left( -1\right) ^{b_{2}}B(gx_{1}\otimes
1_{H};G^{a}X_{1}X_{2}^{b_{2}},g^{d}x_{1}^{e_{1}}x_{2}^{e_{2}})G^{a}X_{2}^{b_{2}}\otimes g^{d}x_{1}^{e_{1}}x_{2}^{e_{2}}=0
\end{eqnarray*}

\subsubsection{$G^{a}\otimes g^{d}$}

\begin{eqnarray*}
&&1_{A}\otimes 1_{H}+ \\
&&+\sum_{\substack{ a,d=0  \\ a+d\equiv 0\text{ }}}^{1}\left[ \left(
-1\right) ^{a+1}B(gx_{1}\otimes 1_{H};G^{a},g^{d}x_{1})+B(gx_{1}\otimes
1_{H};G^{a}X_{1},g^{d})\right] G^{a}\otimes g^{d}=0
\end{eqnarray*}%
and we get%
\begin{equation}
1+-B(gx_{1}\otimes 1_{H};1_{H},x_{1})+B(gx_{1}\otimes 1_{H};X_{1},1_{H})=0
\label{gx1ot1,ten}
\end{equation}%
\begin{equation}
\left[ B(gx_{1}\otimes 1_{H};G,gx_{1})+B(gx_{1}\otimes 1_{H};GX_{1},g)\right]
=0  \label{gx1ot1, eleven}
\end{equation}

\subsubsection{$G^{a}\otimes g^{d}x_{2}$}

\begin{equation*}
\sum_{\substack{ a,d=0  \\ a+d\equiv 1\text{ }}}^{1}\left[ \left( -1\right)
^{a}B(gx_{1}\otimes 1_{H};G^{a},g^{d}x_{1}x_{2})+B(gx_{1}\otimes
1_{H};G^{a}X_{1},g^{d}x_{2})\right] G^{a}\otimes g^{d}x_{2}=0
\end{equation*}%
and we get%
\begin{equation}
B(gx_{1}\otimes 1_{H};1_{A},gx_{1}x_{2})+B(gx_{1}\otimes
1_{H};X_{1},gx_{2})=0  \label{gx1ot1, twelve}
\end{equation}%
\begin{equation}
-B(gx_{1}\otimes 1_{H};G,x_{1}x_{2})+B(gx_{1}\otimes 1_{H};GX_{1},x_{2})=0
\label{gx1ot1, thirteen}
\end{equation}

\subsubsection{$G^{a}\otimes g^{d}x_{1}$}

\begin{equation*}
\sum_{\substack{ a,d=0  \\ a+d\equiv 1\text{ }}}^{1}B(gx_{1}\otimes
1_{H};G^{a}X_{1},g^{d}x_{1})G^{a}\otimes g^{d}x_{1}=0
\end{equation*}%
and we get%
\begin{equation}
B(gx_{1}\otimes 1_{H};X_{1},gx_{1})=0  \label{gx1ot1,fourteen}
\end{equation}%
\begin{equation}
B(gx_{1}\otimes 1_{H};GX_{1},x_{1})=0  \label{gx1ot1, fiveteen}
\end{equation}

\subsubsection{$G^{a}X_{2}\otimes g^{d}$}

\begin{equation*}
\sum_{\substack{ a,d=0  \\ a+d\equiv 1\text{ }}}^{1}\left[ \left( -1\right)
^{a}B(gx_{1}\otimes 1_{H};G^{a}X_{2},g^{d}x_{1})-B(gx_{1}\otimes
1_{H};G^{a}X_{1}X_{2},g^{d})\right] G^{a}X_{2}\otimes g^{d}=0
\end{equation*}%
and we get%
\begin{equation}
B(gx_{1}\otimes 1_{H};X_{2},gx_{1})-B(gx_{1}\otimes 1_{H};X_{1}X_{2},g)=0
\label{gx1ot1, sixteen}
\end{equation}%
\begin{equation}
-B(gx_{1}\otimes 1_{H};GX_{2},x_{1})-B(gx_{1}\otimes
1_{H};GX_{1}X_{2},1_{H})=0  \label{gx1ot1,seventeen}
\end{equation}

\subsubsection{$G^{a}X_{1}\otimes g^{d}$}

\begin{equation*}
\sum_{\substack{ a,d=0  \\ a+d\equiv 1\text{ }}}^{1}\left( -1\right)
^{a}B(gx_{1}\otimes 1_{H};G^{a}X_{1},g^{d}x_{1})G^{a}X_{1}\otimes g^{d}=0
\end{equation*}%
and we get%
\begin{equation*}
B(gx_{1}\otimes 1_{H};X_{1},gx_{1})=0
\end{equation*}%
\begin{equation*}
-B(gx_{1}\otimes 1_{H};GX_{1},x_{1})=0
\end{equation*}%
which are already known.

\subsubsection{$G^{a}\otimes g^{d}x_{1}x_{2}$}

\begin{equation*}
\sum_{\substack{ a,d=0  \\ a+d\equiv 0\text{ }}}^{1}B(gx_{1}\otimes
1_{H};G^{a}X_{1},g^{d}x_{1}x_{2})G^{a}\otimes g^{d}x_{1}x_{2}=0
\end{equation*}%
and we get%
\begin{equation}
B(gx_{1}\otimes 1_{H};X_{1},gx_{1}x_{2})=0  \label{gx1ot1,eigtheen}
\end{equation}%
\begin{equation}
B(gx_{1}\otimes 1_{H};GX_{1},gx_{1}x_{2})=0  \label{gx1ot1,nineteen}
\end{equation}

\subsubsection{$G^{a}X_{2}\otimes g^{d}x_{2}$}

\begin{gather*}
\sum_{\substack{ a,d=0  \\ a+d\equiv 0\text{ }}}^{1}\left[ \left( -1\right)
^{a+1}B(gx_{1}\otimes 1_{H};G^{a}X_{2},g^{d}x_{1}x_{2})-B(gx_{1}\otimes
1_{H};G^{a}X_{1}X_{2},g^{d}x_{2})\right] \\
G^{a}X_{2}\otimes g^{d}x_{2}=0
\end{gather*}%
and we get%
\begin{equation}
-B(gx_{1}\otimes 1_{H};X_{2},x_{1}x_{2})-B(gx_{1}\otimes
1_{H};X_{1}X_{2},x_{2})=0  \label{gx1ot1,twenty}
\end{equation}%
\begin{equation}
B(gx_{1}\otimes 1_{H};GX_{2},gx_{1}x_{2})-B(gx_{1}\otimes
1_{H};GX_{1}X_{2},gx_{2})=0  \label{gx1ot1,twentyone}
\end{equation}

\subsubsection{$G^{a}X_{1}\otimes g^{d}x_{2}$}

\begin{equation*}
\sum_{\substack{ a,d=0  \\ a+d\equiv 0\text{ }}}^{1}\left( -1\right)
^{a+1}B(gx_{1}\otimes 1_{H};G^{a}X_{1},g^{d}x_{1}x_{2})G^{a}X_{1}\otimes
g^{d}x_{2}=0
\end{equation*}%
and we get%
\begin{equation*}
-B(gx_{1}\otimes 1_{H};X_{1},x_{1}x_{2})=0
\end{equation*}%
\begin{equation*}
B(gx_{1}\otimes 1_{H};GX_{1},gx_{1}x_{2})=0
\end{equation*}%
which are already known.

\subsubsection{$G^{a}X_{2}\otimes g^{d}x_{1}$}

\begin{equation*}
\sum_{\substack{ a,d=0  \\ a+d\equiv 0\text{ }}}^{1}-B(gx_{1}\otimes
1_{H};G^{a}X_{1}X_{2},g^{d}x_{1})G^{a}X_{2}\otimes g^{d}x_{1}=0
\end{equation*}%
and we get%
\begin{equation}
-B(gx_{1}\otimes 1_{H};X_{1}X_{2},x_{1})=0  \label{gx1ot1,twentytwo}
\end{equation}%
\begin{equation}
-B(gx_{1}\otimes 1_{H};GX_{1}X_{2},gx_{1})=0  \label{gx1ot1,twentythree}
\end{equation}

\subsubsection{$G^{a}X_{1}\otimes g^{d}x_{1}$}

We do not have a term like this.

\subsubsection{$G^{a}X_{1}X_{2}\otimes g^{d}$}

\begin{equation*}
\sum_{\substack{ a,d=0  \\ a+d\equiv 0\text{ }}}^{1}\left( -1\right)
^{a+1}B(gx_{1}\otimes
1_{H};G^{a}X_{1}X_{2},g^{d}x_{1})G^{a}X_{1}X_{2}\otimes g^{d}=0
\end{equation*}%
and we get%
\begin{equation*}
-B(gx_{1}\otimes 1_{H};X_{1}X_{2},x_{1})=0
\end{equation*}%
\begin{equation*}
B(gx_{1}\otimes 1_{H};GX_{1}X_{2},gx_{1})=0
\end{equation*}%
which are already known.

\subsubsection{$G^{a}X_{2}\otimes g^{d}x_{1}x_{2}$}

\begin{equation*}
\sum_{\substack{ a,d=0  \\ a+d\equiv 1\text{ }}}^{1}-B(gx_{1}\otimes
1_{H};G^{a}X_{1}X_{2},g^{d}x_{1}x_{2})G^{a}X_{2}\otimes g^{d}x_{1}x_{2}=0
\end{equation*}%
and we get%
\begin{equation}
-B(gx_{1}\otimes 1_{H};X_{1}X_{2},gx_{1}x_{2})=0  \label{gx1ot1,twentyfour}
\end{equation}%
\begin{equation}
-B(gx_{1}\otimes 1_{H};GX_{1}X_{2},x_{1}x_{2})=0  \label{gx1ot1,twentyfive.}
\end{equation}

By going on this way, we get no further information

\subsection{Case $gx_{2}$}

Since $\alpha \left( g;0,0,0,1\right) \equiv a+b_{1}+b_{2}+1,\alpha \left(
g;0,1,0,0\right) \equiv 0$ we get%
\begin{eqnarray*}
&&\sum_{\substack{ a,b_{1},b_{2},d,e_{1}=0  \\ a+b_{1}+b_{2}+e_{1}+d\equiv 1%
\text{ }}}^{1}\left( -1\right) ^{a+b_{1}+b_{2}+1}B(gx_{1}\otimes
1_{H};G^{a}X_{1}^{b_{1}}X_{2}^{b_{2}},g^{d}x_{1}^{e_{1}}x_{2})G^{a}X_{1}^{b_{1}}X_{2}^{b_{2}}\otimes g^{d}x_{1}^{e_{1}}
\\
&&\sum_{\substack{ a,b_{1},d,e_{1},e_{2}=0  \\ a+b_{1}+d+e_{1}+e_{2}\equiv 1%
\text{ }}}^{1}B(gx_{1}\otimes
1_{H};G^{a}X_{1}^{b_{1}}X_{2},g^{d}x_{1}^{e_{1}}x_{2}^{e_{2}})G^{a}X_{1}^{b_{1}}\otimes g^{d}x_{1}^{e_{1}}x_{2}^{e_{2}}
\end{eqnarray*}

\subsubsection{$G^{a}\otimes g^{d}$}

\begin{equation*}
\sum_{\substack{ a,d=0  \\ a+d\equiv 0\text{ }}}^{1}\left[ \left( -1\right)
^{a+1}B(gx_{1}\otimes 1_{H};G^{a},g^{d}x_{2})+B(gx_{1}\otimes
1_{H};G^{a}X_{2},g^{d})\right] G^{a}\otimes g^{d}=0
\end{equation*}%
and we get%
\begin{equation}
-B(gx_{1}\otimes 1_{H};1_{A},x_{2})+B(gx_{1}\otimes 1_{H};X_{2},1_{H})=0
\label{gx1ot1,twentysix}
\end{equation}%
\begin{equation}
B(gx_{1}\otimes 1_{H};G,gx_{2})+B(gx_{1}\otimes 1_{H};GX_{2},g)=0
\label{gx1ot1,twentyseven}
\end{equation}

\paragraph{$G^{a}\otimes g^{d}x_{2}$}

\begin{equation*}
\sum_{\substack{ a,d=0  \\ a+d\equiv 1\text{ }}}^{1}B(gx_{1}\otimes
1_{H};G^{a}X_{2},g^{d}x_{2})G^{a}\otimes g^{d}x_{2}=0
\end{equation*}%
and we get%
\begin{equation}
B(gx_{1}\otimes 1_{H};X_{2},gx_{2})=0  \label{gx1ot1,twentyeigth}
\end{equation}%
\begin{equation}
B(gx_{1}\otimes 1_{H};GX_{2},x_{2})=0  \label{gx1ot1,twentynine}
\end{equation}

\subsubsection{$G^{a}\otimes g^{d}x_{1}$}

\begin{equation*}
\sum_{\substack{ a,d=0  \\ a+d\equiv 1\text{ }}}^{1}\left[ \left( -1\right)
^{a+1}B(gx_{1}\otimes 1_{H};G^{a},g^{d}x_{1}x_{2})+B(gx_{1}\otimes
1_{H};G^{a}X_{2},g^{d}x_{1})\right] G^{a}\otimes g^{d}x_{1}=0
\end{equation*}%
and we get%
\begin{equation}
-B(gx_{1}\otimes 1_{H};1_{A},gx_{1}x_{2})+B(gx_{1}\otimes
1_{H};X_{2},gx_{1})=0  \label{gx1ot1,thirty}
\end{equation}%
\begin{equation}
B(gx_{1}\otimes 1_{H};G,x_{1}x_{2})+B(gx_{1}\otimes 1_{H};GX_{2},x_{1})=0
\label{gx1ot1,thirtyone}
\end{equation}

\subsubsection{$G^{a}X_{2}\otimes g^{d}$}

\begin{equation*}
\sum_{\substack{ a,d=0  \\ a+d\equiv 1\text{ }}}^{1}\left( -1\right)
^{a}B(gx_{1}\otimes 1_{H};G^{a}X_{2},g^{d}x_{2})G^{a}X_{2}\otimes g^{d}=0
\end{equation*}%
and we get%
\begin{equation*}
B(gx_{1}\otimes 1_{H};X_{2},gx_{2})=0
\end{equation*}%
\begin{equation*}
-B(gx_{1}\otimes 1_{H};GX_{2},x_{2})=0
\end{equation*}%
which are already known.

\subsubsection{$G^{a}X_{1}\otimes g^{d}$}

\begin{equation*}
\sum_{\substack{ a,d=0  \\ a+d\equiv 1\text{ }}}^{1}\left[ \left( -1\right)
^{a}B(gx_{1}\otimes 1_{H};G^{a}X_{1},g^{d}x_{2})+B(gx_{1}\otimes
1_{H};G^{a}X_{1}X_{2},g^{d})\right] G^{a}X_{1}\otimes g^{d}=0
\end{equation*}%
and we get%
\begin{equation}
B(gx_{1}\otimes 1_{H};X_{1},gx_{2})+B(gx_{1}\otimes 1_{H};X_{1}X_{2},g)=0
\label{gx1ot1,thirtytwo}
\end{equation}%
\begin{equation}
-B(gx_{1}\otimes 1_{H};GX_{1},x_{2})+B(gx_{1}\otimes
1_{H};GX_{1}X_{2},1_{H})=0  \label{gx1ot1,thirtythree}
\end{equation}

\subsubsection{$G^{a}\otimes g^{d}x_{1}x_{2}$}

\begin{equation*}
\sum_{\substack{ a,d=0  \\ a+d\equiv 0\text{ }}}^{1}B(gx_{1}\otimes
1_{H};G^{a}X_{2},g^{d}x_{1}x_{2})G^{a}\otimes g^{d}x_{1}x_{2}=0
\end{equation*}%
and we get%
\begin{equation}
B(gx_{1}\otimes 1_{H};X_{2},x_{1}x_{2})=0  \label{gx1ot1,thirtyfour}
\end{equation}%
\begin{equation}
B(gx_{1}\otimes 1_{H};GX_{2},gx_{1}x_{2})=0  \label{gx1ot1,thirtyfive}
\end{equation}

\subsubsection{$G^{a}X_{2}\otimes g^{d}x_{2}$}

We do not have a term like this.

\subsubsection{$G^{a}X_{1}\otimes g^{d}x_{2}$}

\begin{equation*}
\sum_{\substack{ a,d=0  \\ a+d\equiv 0\text{ }}}^{1}B(gx_{1}\otimes
1_{H};G^{a}X_{1}X_{2},g^{d}x_{2})G^{a}X_{1}\otimes g^{d}x_{2}=0
\end{equation*}%
and we get%
\begin{equation}
B(gx_{1}\otimes 1_{H};X_{1}X_{2},x_{2})=0  \label{gx1ot1,thirtysix}
\end{equation}%
\begin{equation}
B(gx_{1}\otimes 1_{H};GX_{1}X_{2},gx_{2})=0  \label{gx1ot1,thirtyseven}
\end{equation}

\subsubsection{$G^{a}X_{2}\otimes g^{d}x_{1}$}

\begin{equation*}
\sum_{\substack{ a,d=0  \\ a+d\equiv 0\text{ }}}^{1}\left( -1\right)
^{a}B(gx_{1}\otimes 1_{H};G^{a}X_{2},g^{d}x_{1}x_{2})G^{a}X_{2}\otimes
g^{d}x_{1}=0
\end{equation*}%
and we get%
\begin{equation*}
B(gx_{1}\otimes 1_{H};X_{2},x_{1}x_{2})=0
\end{equation*}%
\begin{equation*}
-B(gx_{1}\otimes 1_{H};GX_{2},gx_{1}x_{2})=0
\end{equation*}%
which we already got.

\subsubsection{$G^{a}X_{1}\otimes g^{d}x_{1}$}

\begin{equation*}
\sum_{\substack{ a,d=0  \\ a+d\equiv 0\text{ }}}^{1}\left[ \left( -1\right)
^{a}B(gx_{1}\otimes 1_{H};G^{a}X_{1},g^{d}x_{1}x_{2})+B(gx_{1}\otimes
1_{H};G^{a}X_{1}X_{2},g^{d}x_{1})\right] G^{a}X_{1}\otimes g^{d}x_{1}=0
\end{equation*}%
and we get%
\begin{equation}
B(gx_{1}\otimes 1_{H};X_{1},x_{1}x_{2})+B(gx_{1}\otimes
1_{H};X_{1}X_{2},x_{1})=0  \label{gx1ot1,thirtyseigth}
\end{equation}%
\begin{equation}
-B(gx_{1}\otimes 1_{H};GX_{1},gx_{1}x_{2})+B(gx_{1}\otimes
1_{H};GX_{1}X_{2},gx_{1})=0  \label{gx1ot1,thirtynine}
\end{equation}

\subsubsection{$G^{a}X_{1}X_{2}\otimes g^{d}$}

\begin{equation*}
\sum_{\substack{ a,d=0  \\ a+d\equiv 0\text{ }}}^{1}\left( -1\right)
^{a+1}B(gx_{1}\otimes
1_{H};G^{a}X_{1}X_{2},g^{d}x_{2})G^{a}X_{1}X_{2}\otimes g^{d}=0
\end{equation*}%
and we get%
\begin{equation}
-B(gx_{1}\otimes 1_{H};X_{1}X_{2},x_{2})=0  \label{gx1ot1,fourthy}
\end{equation}%
\begin{equation}
B(gx_{1}\otimes 1_{H};GX_{1}X_{2},gx_{2})=0  \label{gx1ot1,fourthyone}
\end{equation}%
which are known.

\subsubsection{$G^{a}X_{2}\otimes g^{d}x_{1}x_{2}$}

We do not have any term like this

\subsubsection{$G^{a}X_{1}\otimes g^{d}x_{1}x_{2}$}

\begin{equation*}
\sum_{\substack{ a,d=0  \\ a+d\equiv 1\text{ }}}^{1}B(gx_{1}\otimes
1_{H};G^{a}X_{1}X_{2},g^{d}x_{1}x_{2})G^{a}X_{1}\otimes g^{d}x_{1}x_{2}=0
\end{equation*}%
and we get%
\begin{equation*}
B(gx_{1}\otimes 1_{H};X_{1}X_{2},gx_{1}x_{2})=0
\end{equation*}%
\begin{equation*}
B(gx_{1}\otimes 1_{H};GX_{1}X_{2},x_{1}x_{2})=0
\end{equation*}

Going on this way, we do not get anything new

\subsection{Case $gx_{1}x_{2}$}

\begin{eqnarray*}
&&+\sum_{a,b_{1},b_{2},d,e_{1},e_{2}=0}^{1}\sum_{l_{1}=0}^{b_{1}}%
\sum_{l_{2}=0}^{b_{2}}\sum_{u_{1}=0}^{e_{1}}\sum_{u_{2}=0}^{e_{2}}\left(
-1\right) ^{\alpha \left( g;l_{1},l_{2},u_{1},u_{2}\right) } \\
&&B(gx_{1}\otimes
1_{H};G^{a}X_{1}^{b_{1}}X_{2}^{b_{2}},g^{d}x_{1}^{e_{1}}x_{2}^{e_{2}}) \\
&&G^{a}X_{1}^{b_{1}-l_{1}}X_{2}^{b_{2}-l_{2}}\otimes
g^{d}x_{1}^{e_{1}-u_{1}}x_{2}^{e_{2}-u_{2}}\otimes
g^{a+b_{1}+b_{2}+l_{1}+l_{2}+d+e_{1}+e_{2}+u_{1}+u_{2}+1}x_{1}^{l_{1}+u_{1}}x_{2}^{l_{2}+u_{2}}
\\
&=&0
\end{eqnarray*}%
\begin{eqnarray*}
a+b_{1}+b_{2}+l_{1}+l_{2}+d+e_{1}+e_{2}+u_{1}+u_{2}+1 &\equiv &1 \\
l_{1}+u_{1} &\equiv &1 \\
l_{2}+u_{2} &\equiv &1
\end{eqnarray*}

and we get%
\begin{eqnarray*}
a+b_{1}+b_{2}+l_{1}+l_{2}+d+e_{1}+e_{2}+u_{1}+u_{2} &\equiv &0 \\
l_{1}+u_{1} &\equiv &1 \\
l_{2}+u_{2} &\equiv &1
\end{eqnarray*}%
Since, by $\left( \ref{gx1ot1, first}\right) $%
\begin{equation*}
B(gx_{1}\otimes
1_{H};G^{a}X_{1}^{b_{1}}X_{2}^{b_{2}},g^{d}x_{1}^{e_{1}}x_{2}^{e_{2}})=0%
\text{ for }a+b_{1}+b_{2}+d+e_{1}+e_{2}\equiv 0
\end{equation*}%
we get nothing new.

By using all the equalities we got above from $\left( \ref{gx1ot1, first}%
\right) $ to $\left( \ref{gx1ot1,fourthyone}\right) $ we obtain the
following form of $B\left( gx_{1}\otimes 1_{H}\right) .$

\subsection{The final form of the element $B\left( gx_{1}\otimes
1_{H}\right) $}

\begin{eqnarray*}
B\left( gx_{1}\otimes 1_{H}\right) &=&B(gx_{1}\otimes
1_{H};1_{A},g)1_{A}\otimes g+ \\
&&+B(gx_{1}\otimes 1_{H};1_{A},x_{1})1_{A}\otimes x_{1} \\
&&+B(gx_{1}\otimes 1_{H};1_{A},x_{2})1_{A}\otimes x_{2} \\
&&+B(gx_{1}\otimes 1_{H};1_{A},gx_{1}x_{2})1_{A}\otimes gx_{1}x_{2} \\
&&+B(gx_{1}\otimes 1_{H};G,1_{H})G\otimes 1_{H}+ \\
&&+B(gx_{1}\otimes 1_{H};G,x_{1}x_{2})G\otimes x_{1}x_{2}+ \\
&&+B(gx_{1}\otimes 1_{H};G,gx_{1})G\otimes gx_{1}+ \\
&&+B(gx_{1}\otimes 1_{H};G,gx_{2})G\otimes gx_{2}+ \\
&&+\left[ -1+B(gx_{1}\otimes 1_{H};1_{H},x_{1})\right] X_{1}\otimes 1_{H}+ \\
&&-B(gx_{1}\otimes 1_{H};1_{A},gx_{1}x_{2})X_{1}\otimes gx_{2}+ \\
&&+B(gx_{1}\otimes 1_{H};1_{A},x_{2})X_{2}\otimes 1_{H}+ \\
&&+B(gx_{1}\otimes 1_{H};1_{A},gx_{1}x_{2})X_{2}\otimes gx_{1}+ \\
&&+B(gx_{1}\otimes 1_{H};1_{A},gx_{1}x_{2})X_{1}X_{2}\otimes g+ \\
&&-B(gx_{1}\otimes 1_{H};G,gx_{1})GX_{1}\otimes g+ \\
&&+B(gx_{1}\otimes 1_{H};G,x_{1}x_{2})GX_{1}\otimes x_{2}+ \\
&&-B(gx_{1}\otimes 1_{H};G,gx_{2})GX_{2}\otimes g+ \\
&&-B(gx_{1}\otimes 1_{H};G,x_{1}x_{2})GX_{2}\otimes x_{1}+ \\
&&+B(gx_{1}\otimes 1_{H};G,x_{1}x_{2})GX_{1}X_{2}\otimes 1_{H}
\end{eqnarray*}

\section{$B\left( gx_{2}\otimes 1_{H}\right) $}

We write the Casimir formula $\left( \ref{MAIN FORMULA 1}\right) $ for $B$ $%
\left( gx_{2}\otimes 1_{H}\right) .$

For $h=g^{m}x_{1}^{n_{1}}x_{2}^{n_{2}}$ and $h^{\prime }=g^{\mu }x_{1}^{\nu
_{1}}x_{2}^{\nu _{2}}$ we finally get%
\begin{eqnarray*}
&&\sum_{a,b_{1},b_{2},d,e_{1},e_{2}=0}^{1}\sum_{l_{1}=0}^{b_{1}}%
\sum_{l_{2}=0}^{b_{2}}\sum_{u_{1}=0}^{e_{1}}\sum_{u_{2}=0}^{e_{2}}\left(
-1\right) ^{\alpha \left( gx_{2};l_{1},l_{2},u_{1},u_{2}\right) } \\
&&B(1_{H}\otimes
1_{H};G^{a}X_{1}^{b_{1}}X_{2}^{b_{2}},g^{d}x_{1}^{e_{1}}x_{2}^{e_{2}}) \\
&&G^{a}X_{1}^{b_{1}-l_{1}}X_{2}^{b_{2}-l_{2}}\otimes
g^{d}x_{1}^{e_{1}-u_{1}}x_{2}^{e_{2}-u_{2}}\otimes
g^{a+b_{1}+b_{2}+l_{1}+l_{2}+d+e_{1}+e_{2}+u_{1}+u_{2}+1}x_{1}^{l_{1}+u_{1}}x_{2}^{l_{2}+u_{2}+1}
\\
&&\sum_{a,b_{1},b_{2},d,e_{1},e_{2}=0}^{1}\sum_{l_{1}=0}^{b_{1}}%
\sum_{l_{2}=0}^{b_{2}}\sum_{u_{1}=0}^{e_{1}}\sum_{u_{2}=0}^{e_{2}}\left(
-1\right) ^{\alpha \left( g;l_{1},l_{2},u_{1},u_{2}\right) } \\
&&B(gx_{2}\otimes
1_{H};G^{a}X_{1}^{b_{1}}X_{2}^{b_{2}},g^{d}x_{1}^{e_{1}}x_{2}^{e_{2}}) \\
&&G^{a}X_{1}^{b_{1}-l_{1}}X_{2}^{b_{2}-l_{2}}\otimes
g^{d}x_{1}^{e_{1}-u_{1}}x_{2}^{e_{2}-u_{2}}\otimes
g^{a+b_{1}+b_{2}+l_{1}+l_{2}+d+e_{1}+e_{2}+u_{1}+u_{2}+1}x_{1}^{l_{1}+u_{1}}x_{2}^{l_{2}+u_{2}}
\\
&=&B\bigskip \left( gx_{2}\otimes 1_{H}\right) \otimes 1_{H}
\end{eqnarray*}%
Since $B(1_{H}\otimes 1_{H})=1_{A}\otimes 1_{H}$ we obtain%
\begin{eqnarray*}
&&1_{A}\otimes 1_{H}\otimes gx_{2} \\
+
&&\sum_{a,b_{1},b_{2},d,e_{1},e_{2}=0}^{1}\sum_{l_{1}=0}^{b_{1}}%
\sum_{l_{2}=0}^{b_{2}}\sum_{u_{1}=0}^{e_{1}}\sum_{u_{2}=0}^{e_{2}}\left(
-1\right) ^{\alpha \left( g;l_{1},l_{2},u_{1},u_{2}\right) } \\
&&B(gx_{2}\otimes
1_{H};G^{a}X_{1}^{b_{1}}X_{2}^{b_{2}},g^{d}x_{1}^{e_{1}}x_{2}^{e_{2}}) \\
&&G^{a}X_{1}^{b_{1}-l_{1}}X_{2}^{b_{2}-l_{2}}\otimes
g^{d}x_{1}^{e_{1}-u_{1}}x_{2}^{e_{2}-u_{2}}\otimes
g^{a+b_{1}+b_{2}+l_{1}+l_{2}+d+e_{1}+e_{2}+u_{1}+u_{2}+1}x_{1}^{l_{1}+u_{1}}x_{2}^{l_{2}+u_{2}}
\\
&=&B\bigskip \left( gx_{2}\otimes 1_{H}\right) \otimes 1_{H}
\end{eqnarray*}

\subsection{Case $1_{H}$}

The first summand of the left side has no term like this. For the second
summand, we get

\begin{equation*}
g^{d+e_{1}+e_{2}+u_{1}+u_{2}+a+b_{1}+b_{2}+l_{1}+l_{2}+1}x_{1}^{u_{1}+l_{1}}x_{2}^{u_{2}+l_{2}}=1_{H}
\end{equation*}%
and hence%
\begin{eqnarray*}
u_{1}+l_{1} &=&0 \\
u_{2}+l_{2} &=&0 \\
d+e_{1}+e_{2}+u_{1}+u_{2}+a+b_{1}+b_{2}+l_{1}+l_{2}+1 &\equiv &0
\end{eqnarray*}%
i.e.%
\begin{eqnarray*}
u_{1} &=&l_{1}=0 \\
u_{2} &=&l_{2}=0 \\
a+b_{1}+b_{2}+d+e_{1}+e_{2} &\equiv &1
\end{eqnarray*}%
\begin{eqnarray*}
&&\sum_{\substack{ a,b_{1},b_{2},d,e_{1},e_{2}=0  \\ %
a+b_{1}+b_{2}+d+e_{1}+e_{2}\equiv 1}}^{1}\left( -1\right) ^{\alpha \left(
g;l_{1},l_{2},u_{1},u_{2}\right) }B(gx_{2}\otimes \
1_{H};G^{a}X_{1}^{b_{1}}X_{2}^{b_{2}},g^{d}x_{1}^{e_{1}}x_{2}^{e_{2}}) \\
&&G^{a}X_{1}^{b_{1}}X_{2}^{b_{2}}\otimes g^{d}x_{1}^{e_{1}}x_{2}^{e_{2}} \\
&=&\sum_{a,b_{1},b_{2},d,e_{1},e_{2}=0}^{1}B(gx_{2}\otimes \
1_{H};G^{a}X_{1}^{b_{1}}X_{2}^{b_{2}},g^{d}x_{1}^{e_{1}}x_{2}^{e_{2}})G^{a}X_{1}^{b_{1}}X_{2}^{b_{2}}\otimes g^{d}x_{1}^{e_{1}}x_{2}^{e_{2}}
\end{eqnarray*}%
which implies that%
\begin{equation}
B(gx_{2}\otimes \
1_{H};G^{a}X_{1}^{b_{1}}X_{2}^{b_{2}},g^{d}x_{1}^{e_{1}}x_{2}^{e_{2}})=0%
\text{ for }a+b_{1}+b_{2}+d+e_{1}+e_{2}\equiv 0  \label{gx2ot1, first}
\end{equation}

\subsection{Case $g$}

\begin{equation*}
g^{a+b_{1}+b_{2}+d+e_{1}+e_{2}+u_{1}+u_{2}+l_{1}+l_{2}+1}x_{1}^{u_{1}+l_{1}}x_{2}^{u_{2}+l_{2}}=g
\end{equation*}%
this means%
\begin{eqnarray*}
u_{1}+l_{1} &=&0 \\
u_{2}+l_{2} &=&0 \\
a+b_{1}+b_{2}+d+e_{1}+e_{2}+u_{1}+u_{2}+l_{1}+l_{2}+1 &\equiv &1
\end{eqnarray*}%
\begin{eqnarray*}
u_{1} &=&l_{1}=0 \\
u_{2} &=&l_{2}=0 \\
a+b_{1}+b_{2}+d+e_{1}+e_{2} &\equiv &0
\end{eqnarray*}%
which follows from case $1_{H}$.

\subsection{Case $x_{1}$}

\begin{equation*}
g^{d+e_{1}+e_{2}+u_{1}+u_{2}+a+b_{1}+b_{2}+l_{1}+l_{2}+1}x_{1}^{u_{1}+l_{1}}x_{2}^{u_{2}+l_{2}}=x_{1}
\end{equation*}%
\begin{eqnarray*}
g^{a+b_{1}+b_{2}+d+e_{1}+e_{2}+l_{1}+l_{2}+u_{1}+u_{2}+1}x_{1}^{u_{1}+l_{1}}x_{2}^{l_{2}+u_{2}} &=&x_{1}\Rightarrow
\\
l_{2} &=&u_{2}=0 \\
u_{1}+l_{1} &=&1 \\
a+b_{1}+b_{2}+d+e_{1}+e_{2}+l_{1}+l_{2}+u_{1}+u_{2}+1 &\equiv &0 \\
a+b_{1}+b_{2}+d+e_{1}+e_{2} &\equiv &0\text{ }
\end{eqnarray*}%
which follows from case $1_{H}$.

\subsection{Case $x_{2}$}

\begin{equation*}
g^{d+e_{1}+e_{2}+u_{1}+u_{2}+a+b_{1}+b_{2}+l_{1}+l_{2}+1}x_{1}^{u_{1}+l_{1}}x_{2}^{u_{2}+l_{2}}=x_{2}
\end{equation*}%
\begin{eqnarray*}
u_{1} &=&l_{1}=0 \\
u_{2}+l_{2} &=&1 \\
d+e_{1}+e_{2}+u_{1}+u_{2}+a+b_{1}+b_{2}+l_{1}+l_{2}+1 &\equiv &0
\end{eqnarray*}%
\begin{eqnarray*}
u_{1} &=&l_{1}=0 \\
u_{2}+l_{2} &=&1 \\
d+e_{1}+e_{2}+a+b_{1}+b_{2} &\equiv &0
\end{eqnarray*}%
which follows from case $1_{H}$.

\subsection{Case $x_{1}x_{2}$}

\begin{eqnarray*}
&&\sum_{a,b_{1},b_{2},d,e_{1},e_{2}=0}^{1}\sum_{l_{1}=0}^{b_{1}}%
\sum_{l_{2}=0}^{b_{2}}\sum_{u_{1}=0}^{e_{1}}\sum_{u_{2}=0}^{e_{2}}\left(
-1\right) ^{\alpha \left( g;l_{1},l_{2},u_{1},u_{2}\right) } \\
&&B(gx_{2}\otimes
1_{H};G^{a}X_{1}^{b_{1}}X_{2}^{b_{2}},g^{d}x_{1}^{e_{1}}x_{2}^{e_{2}}) \\
&&G^{a}X_{1}^{b_{1}-l_{1}}X_{2}^{b_{2}-l_{2}}\otimes
g^{d}x_{1}^{e_{1}-u_{1}}x_{2}^{e_{2}-u_{2}}\otimes
g^{a+b_{1}+b_{2}+l_{1}+l_{2}+d+e_{1}+e_{2}+u_{1}+u_{2}+1}x_{1}^{l_{1}+u_{1}}x_{2}^{l_{2}+u_{2}}
\\
&=&0
\end{eqnarray*}

We get%
\begin{eqnarray*}
l_{1}+u_{1} &=&1 \\
l_{2}+u_{2} &=&1 \\
a+b_{1}+b_{2}+d+e_{1}+e_{2} &\equiv &1
\end{eqnarray*}%
and hence we obtain, since $\alpha \left( g;0,0,1,1\right) =1+e_{2},$ $%
\alpha \left( g;0,1,1,0\right) =e_{2}+a+b_{1}+b_{2}$, $\alpha \left(
g;1,0,0,1\right) =a+b_{1}+1$ and

$\alpha \left( g;1,1,0,0\right) =1+b_{2}$

\subsubsection{\protect\bigskip}

\begin{eqnarray*}
&&\sum_{\substack{ a,b_{1},b_{2},d=0  \\ a+b_{1}+b_{2}+d\equiv 1}}%
^{1}B(gx_{2}\otimes
1_{H};G^{a}X_{1}^{b_{1}}X_{2}^{b_{2}},g^{d}x_{1}x_{2})G^{a}X_{1}^{b_{1}}X_{2}^{b_{2}}\otimes g^{d}+
\\
&&+\sum_{\substack{ a,b_{1},d,e_{2}=0  \\ a+b_{1}+d+e_{2}\equiv 1}}%
^{1}\left( -1\right) ^{e_{2}+a+b_{1}+1}B(gx_{2}\otimes
1_{H};G^{a}X_{1}^{b_{1}}X_{2},g^{d}x_{1}x_{2}^{e_{2}})G^{a}X_{1}^{b_{1}}%
\otimes g^{d}x_{2}^{e_{2}}+ \\
&&+\sum_{\substack{ a,b_{2},d,e_{1}=0  \\ a+b_{2}+d+e_{1}\equiv 1}}%
^{1}\left( -1\right) ^{a}B(gx_{2}\otimes
1_{H};G^{a}X_{1}X_{2}^{b_{2}},g^{d}x_{1}^{e_{1}}x_{2})G^{a}X_{2}^{b_{2}}%
\otimes g^{d}x_{1}^{e_{1}}+ \\
&&+\sum_{\substack{ a,d,e_{1},e_{2}=0  \\ a+d+e_{1}+e_{2}\equiv 1}}%
^{1}B(gx_{2}\otimes
1_{H};G^{a}X_{1}X_{2},g^{d}x_{1}^{e_{1}}x_{2}^{e_{2}})G^{a}\otimes
g^{d}x_{1}^{e_{1}}x_{2}^{e_{2}} \\
&=&0
\end{eqnarray*}%
$\left( 0,0,0,0\right) G^{a}\otimes g^{d}$

\begin{equation*}
\sum_{\substack{ a,d=0  \\ a+d\equiv 1}}^{1}\left[
\begin{array}{c}
B(gx_{2}\otimes 1_{H};G^{a},g^{d}x_{1}x_{2})+\left( -1\right)
^{a+1}B(gx_{2}\otimes 1_{H};G^{a}X_{2},g^{d}x_{1})+ \\
+\left( -1\right) ^{a}B(gx_{2}\otimes
1_{H};G^{a}X_{1},g^{d}x_{2})+B(gx_{2}\otimes 1_{H};G^{a}X_{1}X_{2},g^{d})%
\end{array}%
\right] G^{a}\otimes g^{d}=0
\end{equation*}%
and we get%
\begin{equation}
\begin{array}{c}
B(gx_{2}\otimes 1_{H};1_{A},gx_{1}x_{2})-B(gx_{2}\otimes 1_{H};X_{2},gx_{1})+
\\
+B(gx_{2}\otimes 1_{H};X_{1},gx_{2})+B(gx_{2}\otimes 1_{H};X_{1}X_{2},g)%
\end{array}%
=0  \label{gx2ot1, second}
\end{equation}%
and%
\begin{equation}
\begin{array}{c}
B(gx_{2}\otimes 1_{H};G,x_{1}x_{2})+B(gx_{2}\otimes 1_{H};GX_{2},x_{1})+ \\
-B(gx_{2}\otimes 1_{H};GX_{1},x_{2})+B(gx_{2}\otimes 1_{H};GX_{1}X_{2},1_{H})%
\end{array}%
=0  \label{gx2ot1, third}
\end{equation}

\subsubsection{$G^{a}\otimes g^{d}x_{2}$}

\begin{equation*}
\sum_{\substack{ a,d=0  \\ a+d\equiv 0}}^{1}\left[ \left( -1\right)
^{a}B(gx_{2}\otimes 1_{H};G^{a}X_{2},g^{d}x_{1}x_{2})+B(gx_{2}\otimes
1_{H};G^{a}X_{1}X_{2},g^{d}x_{2})\right] G^{a}\otimes g^{d}x_{2}=0
\end{equation*}%
and we get%
\begin{equation}
B(gx_{2}\otimes 1_{H};X_{2},x_{1}x_{2})+B(gx_{2}\otimes
1_{H};X_{1}X_{2},x_{2})=0  \label{gx2ot1, four}
\end{equation}%
and%
\begin{equation}
-B(gx_{2}\otimes 1_{H};GX_{2},gx_{1}x_{2})+B(gx_{2}\otimes
1_{H};GX_{1}X_{2},gx_{2})=0  \label{gx2ot1, five}
\end{equation}

\subsubsection{$G^{a}\otimes g^{d}x_{1}$}

\begin{equation*}
\sum_{\substack{ a,d=0  \\ a+d\equiv 0}}^{1}\left[ \left( -1\right)
^{a}B(gx_{2}\otimes 1_{H};G^{a}X_{1},g^{d}x_{1}x_{2})+B(gx_{2}\otimes
1_{H};G^{a}X_{1}X_{2},g^{d}x_{1})\right] G^{a}\otimes g^{d}x_{1}=0
\end{equation*}%
and we get%
\begin{equation}
B(gx_{2}\otimes 1_{H};X_{1},x_{1}x_{2})+B(gx_{2}\otimes
1_{H};X_{1}X_{2},x_{1})=0  \label{gx2ot1, six}
\end{equation}%
and%
\begin{equation}
-B(gx_{2}\otimes 1_{H};GX_{1},gx_{1}x_{2})+B(gx_{2}\otimes
1_{H};GX_{1}X_{2},gx_{1})=0  \label{gx2ot1, seven}
\end{equation}

\subsubsection{$G^{a}X_{2}\otimes g^{d}$}

\begin{equation*}
\sum_{\substack{ a,d=0  \\ a+d\equiv 0}}^{1}\left[ B(gx_{2}\otimes
1_{H};G^{a}X_{2},g^{d}x_{1}x_{2})+\left( -1\right) ^{a}B(gx_{2}\otimes
1_{H};G^{a}X_{1}X_{2},g^{d}x_{2})\right] G^{a}X_{2}\otimes g^{d}=0
\end{equation*}%
And from this we obtain the same equalities we got in case $G^{a}\otimes
g^{d}x_{2}.$

\subsubsection{$G^{a}X_{1}\otimes g^{d}$}

\begin{equation*}
\sum_{\substack{ a,d=0  \\ a+d\equiv 0}}^{1}\left[ B(gx_{2}\otimes
1_{H};G^{a}X_{1},g^{d}x_{1}x_{2})+\left( -1\right) ^{a}B(gx_{2}\otimes
1_{H};G^{a}X_{1}X_{2},g^{d}x_{1})\right] G^{a}X_{1}\otimes g^{d}=0
\end{equation*}%
And from this we obtain the same equalities we got in case $G^{a}\otimes
g^{d}x_{1}.$

\subsubsection{$G^{a}\otimes g^{d}x_{1}x_{2}$}

\begin{equation*}
\sum_{\substack{ a,d=0  \\ a+d\equiv 1}}^{1}B(gx_{2}\otimes
1_{H};G^{a}X_{1}X_{2},g^{d}x_{1}x_{2})G^{a}\otimes g^{d}x_{1}x_{2}=0
\end{equation*}%
and we get%
\begin{equation}
B(gx_{2}\otimes 1_{H};X_{1}X_{2},gx_{1}x_{2})=0  \label{gx2ot1, eight}
\end{equation}%
\begin{equation}
B(gx_{2}\otimes 1_{H};GX_{1}X_{2},x_{1}x_{2})=0  \label{gx2ot1, nine}
\end{equation}

\subsubsection{$G^{a}X_{2}\otimes g^{d}x_{2}$}

There is no term like this.

\subsubsection{$G^{a}X_{1}\otimes g^{d}x_{2}$}

\begin{equation*}
\sum_{\substack{ a,d=0  \\ a+d\equiv 1}}^{1}\left( -1\right)
^{a+1}B(gx_{2}\otimes
1_{H};G^{a}X_{1}X_{2},g^{d}x_{1}x_{2})G^{a}X_{1}\otimes g^{d}x_{2}=0
\end{equation*}%
and we get the same equations as in case $G^{a}\otimes g^{d}x_{1}x_{2}.$

\subsubsection{$G^{a}X_{2}\otimes g^{d}x_{1}$}

\begin{equation*}
\sum_{\substack{ a,d=0  \\ a+d\equiv 1}}^{1}\left( -1\right)
^{a}B(gx_{2}\otimes 1_{H};G^{a}X_{1}X_{2},g^{d}x_{1}x_{2})G^{a}X_{2}\otimes
g^{d}x_{1}=0
\end{equation*}%
and we get the same equations as in case $G^{a}\otimes g^{d}x_{1}x_{2}.$

\subsubsection{$G^{a}X_{1}\otimes g^{d}x_{1}$}

There is no term like this.

\subsubsection{$G^{a}X_{1}X_{2}\otimes g^{d}$}

\begin{equation*}
\sum_{\substack{ a,d=0  \\ a+d\equiv 1}}^{1}B(gx_{2}\otimes
1_{H};G^{a}X_{1}X_{2},g^{d}x_{1}x_{2})G^{a}X_{1}X_{2}\otimes g^{d}=0
\end{equation*}%
and we get the same equations as in case $G^{a}\otimes g^{d}x_{1}x_{2}.$

\subsubsection{$G^{a}X_{2}\otimes g^{d}x_{1}x_{2}$}

There is no term like this.

\subsubsection{$G^{a}X_{1}\otimes g^{d}x_{1}x_{2}$}

\subsubsection{$G^{a}X_{1}X_{2}\otimes g^{d}x_{2}$}

There is no term like this.

\subsubsection{$G^{a}X_{1}X_{2}\otimes g^{d}x_{1}$}

There is no term like this.

\subsubsection{$G^{a}X_{1}X_{2}\otimes g^{d}x_{1}x_{2}$}

There is no term like this.

\subsection{Case $gx_{1}$}

\begin{eqnarray*}
&&\sum_{a,b_{1},b_{2},d,e_{1},e_{2}=0}^{1}\sum_{l_{1}=0}^{b_{1}}%
\sum_{l_{2}=0}^{b_{2}}\sum_{u_{1}=0}^{e_{1}}\sum_{u_{2}=0}^{e_{2}}\left(
-1\right) ^{\alpha \left( g;l_{1},l_{2},u_{1},u_{2}\right) } \\
&&B(gx_{2}\otimes
1_{H};G^{a}X_{1}^{b_{1}}X_{2}^{b_{2}},g^{d}x_{1}^{e_{1}}x_{2}^{e_{2}}) \\
&&G^{a}X_{1}^{b_{1}-l_{1}}X_{2}^{b_{2}-l_{2}}\otimes
g^{d}x_{1}^{e_{1}-u_{1}}x_{2}^{e_{2}-u_{2}}\otimes
g^{a+b_{1}+b_{2}+l_{1}+l_{2}+d+e_{1}+e_{2}+u_{1}+u_{2}+1}x_{1}^{l_{1}+u_{1}}x_{2}^{l_{2}+u_{2}}
\end{eqnarray*}

\begin{eqnarray*}
l_{1}+u_{1} &=&1 \\
l_{2} &=&u_{2}=0 \\
a+b_{1}+b_{2}+d+e_{1}+e_{2} &\equiv &1
\end{eqnarray*}%
and, since $\alpha \left( g;0,0,1,0\right) =e_{2}+\left(
a+b_{1}+b_{2}\right) +1$ and $\alpha \left( g;1,0,0,0\right) =b_{2},$ we get%
\begin{eqnarray*}
&&\sum_{\substack{ a,b_{1},b_{2},d,e_{2}=0  \\ a+b_{1}+b_{2}+d+e_{2}\equiv 0
}}^{1}\left( -1\right) ^{e_{2}+\left( a+b_{1}+b_{2}\right)
+1}B(gx_{2}\otimes
1_{H};G^{a}X_{1}^{b_{1}}X_{2}^{b_{2}},g^{d}x_{1}x_{2}^{e_{2}})G^{a}X_{1}^{b_{1}}X_{2}^{b_{2}}\otimes g^{d}x_{2}^{e_{2}}
\\
&&+\sum_{\substack{ a,b_{2},d,e_{1},e_{2}=0  \\ a+b_{2}+d+e_{1}+e_{2}\equiv
0 }}^{1}\left( -1\right) ^{b_{2}}B(gx_{2}\otimes
1_{H};G^{a}X_{1}X_{2}^{b_{2}},g^{d}x_{1}^{e_{1}}x_{2}^{e_{2}})G^{a}X_{2}^{b_{2}}\otimes g^{d}x_{1}^{e_{1}}x_{2}^{e_{2}}.
\end{eqnarray*}

\paragraph{$G^{a}\otimes g^{d}$}

\begin{equation*}
\sum_{\substack{ a,d=0  \\ a+d\equiv 0}}^{1}\left[ \left( -1\right)
^{a+1}B(gx_{2}\otimes 1_{H};G^{a},g^{d}x_{1})+B(gx_{2}\otimes
1_{H};G^{a}X_{1},g^{d})\right] G^{a}\otimes g^{d}=0
\end{equation*}%
and we get%
\begin{equation}
-B(gx_{2}\otimes 1_{H};1_{A},x_{1})+B(gx_{2}\otimes 1_{H};X_{1},1_{H})=0
\label{gx2ot1,ten}
\end{equation}%
\begin{equation}
B(gx_{2}\otimes 1_{H};G,gx_{1})+B(gx_{2}\otimes 1_{H};GX_{1},g)=0
\label{gx2ot1, eleven}
\end{equation}

\subsubsection{$G^{a}\otimes g^{d}x_{2}$}

\begin{equation*}
\sum_{\substack{ a,d=0  \\ a+d\equiv 1}}^{1}\left[ \left( -1\right)
^{a}B(gx_{2}\otimes 1_{H};G^{a},g^{d}x_{1}x_{2})+B(gx_{2}\otimes
1_{H};G^{a}X_{1},g^{d}x_{2})\right] G^{a}\otimes g^{d}x_{2}=0
\end{equation*}%
and we get%
\begin{equation}
B(gx_{2}\otimes 1_{H};1_{A},gx_{1}x_{2})+B(gx_{2}\otimes
1_{H};X_{1},gx_{2})=0  \label{gx2ot1, twelve}
\end{equation}%
\begin{equation}
-B(gx_{2}\otimes 1_{H};G,x_{1}x_{2})+B(gx_{2}\otimes 1_{H};GX_{1},x_{2})=0
\label{gx2ot1, thirteen}
\end{equation}

\subsubsection{$G^{a}\otimes g^{d}x_{1}$}

\begin{equation*}
\sum_{\substack{ a,d=0  \\ a+d\equiv 1}}^{1}B(gx_{2}\otimes
1_{H};G^{a}X_{1},g^{d}x_{1})G^{a}\otimes g^{d}x_{1}=0
\end{equation*}%
and we get%
\begin{equation}
B(gx_{2}\otimes 1_{H};X_{1},gx_{1})=0  \label{gx2ot1, fourteen}
\end{equation}%
and%
\begin{equation}
B(gx_{2}\otimes 1_{H};GX_{1},x_{1})=0  \label{gx2ot1, fiveteen}
\end{equation}

\subsubsection{$G^{a}X_{2}\otimes g^{d}$}

\begin{equation*}
\sum_{\substack{ a,d=0  \\ a+d\equiv 1}}^{1}\left[ \left( -1\right)
^{a}B(gx_{2}\otimes 1_{H};G^{a}X_{2},g^{d}x_{1})-B(gx_{2}\otimes
1_{H};G^{a}X_{1}X_{2},g^{d})\right] G^{a}X_{2}\otimes g^{d}=0
\end{equation*}%
and we get%
\begin{equation}
B(gx_{2}\otimes 1_{H};X_{2},gx_{1})-B(gx_{2}\otimes 1_{H};X_{1}X_{2},g)=0
\label{gx2ot1, sixteen}
\end{equation}%
\begin{equation}
-B(gx_{2}\otimes 1_{H};GX_{2},x_{1})-B(gx_{2}\otimes
1_{H};GX_{1}X_{2},1_{H})=0  \label{gx2ot1, seventeen}
\end{equation}

\subsubsection{$G^{a}X_{1}\otimes g^{d}$}

\begin{equation*}
\sum_{\substack{ a,d=0  \\ a+d\equiv 1}}^{1}\left( -1\right)
^{a}B(gx_{2}\otimes 1_{H};G^{a}X_{1},g^{d}x_{1})G^{a}X_{1}\otimes g^{d}
\end{equation*}%
and we get%
\begin{equation*}
B(gx_{2}\otimes 1_{H};X_{1},gx_{1})=0
\end{equation*}%
and%
\begin{equation*}
-B(gx_{2}\otimes 1_{H};GX_{1},x_{1})=0
\end{equation*}%
which we already know.

\subsubsection{$G^{a}\otimes g^{d}x_{1}x_{2}$}

\begin{equation*}
\sum_{\substack{ a,d=0  \\ a+d\equiv 0}}^{1}B(gx_{2}\otimes
1_{H};G^{a}X_{1},g^{d}x_{1}x_{2})G^{a}\otimes g^{d}x_{1}x_{2}=0
\end{equation*}%
and we get%
\begin{equation}
B(gx_{2}\otimes 1_{H};X_{1},x_{1}x_{2})=0  \label{gx2ot1, eighteen}
\end{equation}%
\begin{equation}
B(gx_{2}\otimes 1_{H};GX_{1},gx_{1}x_{2})=0  \label{gx2ot1, nineteen}
\end{equation}

\subsubsection{$G^{a}X_{2}\otimes g^{d}x_{2}$}

\begin{equation*}
\sum_{\substack{ a,d=0  \\ a+d\equiv 0}}^{1}\left[ \left( -1\right)
^{a+1}B(gx_{2}\otimes 1_{H};G^{a}X_{2},g^{d}x_{1}x_{2})-B(gx_{2}\otimes
1_{H};G^{a}X_{1}X_{2},g^{d}x_{2})\right] G^{a}X_{2}\otimes g^{d}x_{2}=0
\end{equation*}%
and we get%
\begin{equation}
-B(gx_{2}\otimes 1_{H};X_{2},x_{1}x_{2})-B(gx_{2}\otimes
1_{H};X_{1}X_{2},x_{2})=0  \label{gx2ot1, twenty}
\end{equation}%
\begin{equation}
B(gx_{2}\otimes 1_{H};GX_{2},gx_{1}x_{2})-B(gx_{2}\otimes
1_{H};GX_{1}X_{2},gx_{2})=0  \label{gx2ot1, twentyone}
\end{equation}

\subsubsection{$G^{a}X_{1}\otimes g^{d}x_{2}$}

\begin{equation*}
\sum_{\substack{ a,d=0  \\ a+d\equiv 0}}^{1}\left( -1\right)
^{a+1}B(gx_{2}\otimes 1_{H};G^{a}X_{1},g^{d}x_{1}x_{2})G^{a}X_{1}\otimes
g^{d}x_{2}=0
\end{equation*}%
and we get the same equations as in case $G^{a}\otimes g^{d}x_{1}x_{2}.$

\subsubsection{$G^{a}X_{2}\otimes g^{d}x_{1}$}

\begin{equation*}
\sum_{\substack{ a,d=0  \\ a+d\equiv 0}}^{1}-B(gx_{2}\otimes
1_{H};G^{a}X_{1}X_{2},g^{d}x_{1})G^{a}X_{2}\otimes g^{d}x_{1}=0
\end{equation*}%
and we get%
\begin{equation}
-B(gx_{2}\otimes 1_{H};X_{1}X_{2},x_{1})=0  \label{gx2ot1, twentytwo}
\end{equation}%
\begin{equation}
-B(gx_{2}\otimes 1_{H};GX_{1}X_{2},gx_{1})=0  \label{gx2ot1, twentythree}
\end{equation}

\subsubsection{$G^{a}X_{1}\otimes g^{d}x_{1}$}

We do not have any term like this.

\subsubsection{$G^{a}X_{1}X_{2}\otimes g^{d}$}

\begin{equation*}
\sum_{\substack{ a,d=0  \\ a+d\equiv 0}}^{1}\left( -1\right)
^{a+1}B(gx_{2}\otimes
1_{H};G^{a}X_{1}X_{2},g^{d}x_{1})G^{a}X_{1}X_{2}\otimes g^{d}=0
\end{equation*}%
and we get the same equations as in case $G^{a}X_{2}\otimes g^{d}x_{1}.$

\subsubsection{$G^{a}X_{2}\otimes g^{d}x_{1}x_{2}$}

\begin{equation*}
\sum_{\substack{ a,d=0  \\ a+d\equiv 1}}^{1}-B(gx_{2}\otimes
1_{H};G^{a}X_{1}X_{2},g^{d}x_{1}x_{2})G^{a}X_{2}\otimes g^{d}x_{1}x_{2}=0
\end{equation*}%
and we get the same equations as in case $G^{a}\otimes g^{d}x_{1}x_{2}$ of
the subsection $x_{1}x_{2.}$

\subsubsection{$G^{a}X_{1}\otimes g^{d}x_{1}x_{2}$}

There is no term like this.

\subsubsection{$G^{a}X_{1}X_{2}\otimes g^{d}x_{2}$}

\begin{equation*}
\sum_{\substack{ a,d=0  \\ a+d\equiv 1}}^{1}\left( -1\right)
^{a}B(gx_{2}\otimes
1_{H};G^{a}X_{1}X_{2},g^{d}x_{1}x_{2})G^{a}X_{1}X_{2}\otimes g^{d}x_{2}=0
\end{equation*}%
and we get the same equations as in case $G^{a}\otimes g^{d}x_{1}x_{2}$ of
the subsection $x_{1}x_{2.}$

\subsubsection{$G^{a}X_{1}X_{2}\otimes g^{d}x_{1}$}

We do not have any term like this.

\subsubsection{$G^{a}X_{1}X_{2}\otimes g^{d}x_{1}x_{2}$}

We do not have any term like this.

\subsection{Case $gx_{2}$}

We have%
\begin{eqnarray*}
&&1_{A}\otimes 1_{H}\otimes gx_{2} \\
+
&&\sum_{a,b_{1},b_{2},d,e_{1},e_{2}=0}^{1}\sum_{l_{1}=0}^{b_{1}}%
\sum_{l_{2}=0}^{b_{2}}\sum_{u_{1}=0}^{e_{1}}\sum_{u_{2}=0}^{e_{2}}\left(
-1\right) ^{\alpha \left( g;l_{1},l_{2},u_{1},u_{2}\right) } \\
&&B(gx_{2}\otimes
1_{H};G^{a}X_{1}^{b_{1}}X_{2}^{b_{2}},g^{d}x_{1}^{e_{1}}x_{2}^{e_{2}}) \\
&&G^{a}X_{1}^{b_{1}-l_{1}}X_{2}^{b_{2}-l_{2}}\otimes
g^{d}x_{1}^{e_{1}-u_{1}}x_{2}^{e_{2}-u_{2}}\otimes
g^{a+b_{1}+b_{2}+l_{1}+l_{2}+d+e_{1}+e_{2}+u_{1}+u_{2}+1}x_{1}^{l_{1}+u_{1}}x_{2}^{l_{2}+u_{2}}
\\
&=&0
\end{eqnarray*}%
\begin{eqnarray*}
l_{1} &=&u_{1}=0 \\
l_{2}+u_{2} &=&1 \\
a+b_{1}+b_{2}+d+e_{1}+e_{2} &\equiv &1
\end{eqnarray*}%
Since $\alpha \left( g;0,0,0,1\right) =\left( a+b_{1}+b_{2}\right) +1$ and $%
\alpha \left( g;0,1,0,0\right) =0$ we get%
\begin{gather*}
1_{A}\otimes 1_{H}+ \\
+\sum_{\substack{ a,b_{1},b_{2},d,e_{1}=0  \\ a+b_{1}+b_{2}+d+e_{1}\equiv 0}}%
^{1}B(gx_{2}\otimes
1_{H};G^{a}X_{1}^{b_{1}}X_{2}^{b_{2}},g^{d}x_{1}^{e_{1}}x_{2})G^{a}X_{1}^{b_{1}}X_{2}^{b_{2}}\otimes g^{d}x_{1}^{e_{1}}+
\\
+\sum_{\substack{ a,b_{1},d,e_{1},e_{2}=0  \\ a+b_{1}+d+e_{1}+e_{2}\equiv 0}}%
^{1}B(gx_{2}\otimes
1_{H};G^{a}X_{1}^{b_{1}}X_{2},g^{d}x_{1}^{e_{1}}x_{2}^{e_{2}})G^{a}X_{1}^{b_{1}}\otimes g^{d}x_{1}^{e_{1}}x_{2}^{e_{2}}=0
\end{gather*}

\subsubsection{$G^{a}\otimes g^{d}$}

\begin{gather*}
1_{A}\otimes 1_{H}+ \\
\sum_{\substack{ a,d=0  \\ a+d\equiv 0}}^{1}\left[ \left( -1\right)
^{a+1}B(gx_{2}\otimes 1_{H};G^{a},g^{d}x_{2})+B(gx_{2}\otimes
1_{H};G^{a}X_{2},g^{d})\right] G^{a}\otimes g^{d}=0
\end{gather*}%
and we get%
\begin{equation}
1-B(gx_{2}\otimes 1_{H};1_{A},x_{2})+B(gx_{2}\otimes 1_{H};X_{2},1_{H})=0
\label{gx2ot1, twentyfour}
\end{equation}%
and%
\begin{equation}
B(gx_{2}\otimes 1_{H};G,gx_{2})+B(gx_{2}\otimes 1_{H};GX_{2},g)=0
\label{gx2ot1, twentyfive}
\end{equation}

\subsubsection{$G^{a}\otimes g^{d}x_{2}$}

\begin{equation*}
\sum_{\substack{ a,d=0  \\ a+d\equiv 1}}^{1}B(gx_{2}\otimes
1_{H};G^{a}X_{2},g^{d}x_{2})G^{a}\otimes g^{d}x_{2}=0
\end{equation*}%
and we get%
\begin{equation}
B(gx_{2}\otimes 1_{H};X_{2},gx_{2})=0  \label{gx2ot1, twentysix}
\end{equation}%
and%
\begin{equation}
B(gx_{2}\otimes 1_{H};GX_{2},x_{2})=0  \label{gx2ot1, twentyseven}
\end{equation}

\subsubsection{$G^{a}\otimes g^{d}x_{1}$}

\begin{equation*}
\sum_{\substack{ a,d=0  \\ a+d\equiv 1}}^{1}\left[ \left( -1\right)
^{a+1}B(gx_{2}\otimes 1_{H};G^{a},g^{d}x_{1}x_{2})+B(gx_{2}\otimes
1_{H};G^{a}X_{2},g^{d}x_{1})\right] G^{a}\otimes g^{d}x_{1}=0
\end{equation*}%
and we get%
\begin{equation}
-B(gx_{2}\otimes 1_{H};1_{A},gx_{1}x_{2})+B(gx_{2}\otimes
1_{H};X_{2},gx_{1})=0  \label{gx2ot1, twentyeight}
\end{equation}%
\begin{equation}
B(gx_{2}\otimes 1_{H};G,x_{1}x_{2})+B(gx_{2}\otimes 1_{H};GX_{2},x_{1})=0
\label{gx2ot1, twentynine}
\end{equation}

\subsubsection{$G^{a}X_{2}\otimes g^{d}$}

\begin{equation*}
\sum_{\substack{ a,d=0  \\ a+d\equiv 1}}^{1}\left( -1\right)
^{a}B(gx_{2}\otimes 1_{H};G^{a}X_{2},g^{d}x_{2})G^{a}X_{2}\otimes g^{d}=0
\end{equation*}%
and we get%
\begin{equation*}
B(gx_{2}\otimes 1_{H};X_{2},gx_{2})=0
\end{equation*}%
and%
\begin{equation*}
-B(gx_{2}\otimes 1_{H};GX_{2},x_{2})=0
\end{equation*}%
that we already obtained.

\subsubsection{$G^{a}X_{1}\otimes g^{d}$}

\begin{equation*}
\sum_{\substack{ a,d=0  \\ a+d\equiv 1}}^{1}\left[ \left( -1\right)
^{a}B(gx_{2}\otimes 1_{H};G^{a}X_{1},g^{d}x_{2})+B(gx_{2}\otimes
1_{H};G^{a}X_{1}X_{2},g^{d})\right] G^{a}X_{1}\otimes g^{d}=0
\end{equation*}%
and we get%
\begin{equation}
B(gx_{2}\otimes 1_{H};X_{1},gx_{2})+B(gx_{2}\otimes 1_{H};X_{1}X_{2},g)=0
\label{gx2ot1, thirty}
\end{equation}%
and%
\begin{equation}
-B(gx_{2}\otimes 1_{H};GX_{1},x_{2})+B(gx_{2}\otimes
1_{H};GX_{1}X_{2},1_{H})=0.  \label{gx2ot1, thirtyone}
\end{equation}

\subsubsection{$G^{a}\otimes g^{d}x_{1}x_{2}$}

\begin{equation*}
\sum_{\substack{ a,d=0  \\ a+d\equiv 0}}^{1}B(gx_{2}\otimes
1_{H};G^{a}X_{2},g^{d}x_{1}x_{2})G^{a}\otimes g^{d}x_{1}x_{2}=0
\end{equation*}%
and we get%
\begin{equation}
B(gx_{2}\otimes 1_{H};X_{2},x_{1}x_{2})=0  \label{gx2ot1, thirtytwo}
\end{equation}%
and%
\begin{equation}
B(gx_{2}\otimes 1_{H};GX_{2},gx_{1}x_{2})=0  \label{gx2ot1, thirtythree}
\end{equation}

\subsubsection{$G^{a}X_{2}\otimes g^{d}x_{2}$}

We do not have any term like this.

\subsubsection{$G^{a}X_{1}\otimes g^{d}x_{2}$}

\begin{equation*}
\sum_{\substack{ a,d=0  \\ a+d\equiv 0}}^{1}B(gx_{2}\otimes
1_{H};G^{a}X_{1}X_{2},g^{d}x_{2})G^{a}X_{1}\otimes g^{d}x_{2}=0
\end{equation*}%
and we get%
\begin{equation}
B(gx_{2}\otimes 1_{H};X_{1}X_{2},x_{2})=0  \label{gx2ot1, thirtyfour}
\end{equation}%
and%
\begin{equation}
B(gx_{2}\otimes 1_{H};GX_{1}X_{2},gx_{2})=0  \label{gx2ot1, thirtyfive}
\end{equation}

\subsubsection{$G^{a}X_{2}\otimes g^{d}x_{1}$}

\begin{equation*}
\sum_{\substack{ a,d=0  \\ a+d\equiv 0}}^{1}\left( -1\right)
^{a}B(gx_{2}\otimes 1_{H};G^{a}X_{2},g^{d}x_{1}x_{2})G^{a}X_{2}\otimes
g^{d}x_{1}=0
\end{equation*}%
and we get%
\begin{equation*}
B(gx_{2}\otimes 1_{H};X_{2},x_{1}x_{2})=0
\end{equation*}%
and%
\begin{equation*}
-B(gx_{2}\otimes 1_{H};GX_{2},gx_{1}x_{2})=0
\end{equation*}%
which we already got.

\subsubsection{$G^{a}X_{1}\otimes g^{d}x_{1}$}

\begin{equation*}
\sum_{\substack{ a,d=0  \\ a+d\equiv 0}}^{1}\left[ \left( -1\right)
^{a}B(gx_{2}\otimes 1_{H};G^{a}X_{1},g^{d}x_{1}x_{2})+B(gx_{2}\otimes
1_{H};G^{a}X_{1}X_{2},g^{d}x_{1})\right] G^{a}X_{1}\otimes g^{d}x_{1}=0
\end{equation*}%
and we get%
\begin{equation}
B(gx_{2}\otimes 1_{H};X_{1},x_{1}x_{2})+B(gx_{2}\otimes
1_{H};X_{1}X_{2},x_{1})=0  \label{gx2ot1, thirtysix}
\end{equation}%
\begin{equation}
-B(gx_{2}\otimes 1_{H};GX_{1},gx_{1}x_{2})+B(gx_{2}\otimes
1_{H};GX_{1}X_{2},gx_{1})=0  \label{gx2ot1, thirtyseven}
\end{equation}

\subsubsection{$G^{a}X_{1}X_{2}\otimes g^{d}$}

\begin{equation*}
\sum_{\substack{ a,d=0  \\ a+d\equiv 0}}^{1}\left( -1\right)
^{a+1}B(gx_{2}\otimes
1_{H};G^{a}X_{1}X_{2},g^{d}x_{2})G^{a}X_{1}^{b_{1}}X_{2}^{b_{2}}\otimes
g^{d}x_{1}^{e_{1}}=0
\end{equation*}%
and we get%
\begin{equation*}
-B(gx_{2}\otimes 1_{H};X_{1}X_{2},x_{2})=0
\end{equation*}%
and%
\begin{equation*}
B(gx_{2}\otimes 1_{H};GX_{1}X_{2},gx_{2})=0
\end{equation*}%
which we already got.

\subsubsection{$G^{a}X_{2}\otimes g^{d}x_{1}x_{2}$}

There is no term like this.

\subsubsection{$G^{a}X_{1}\otimes g^{d}x_{1}x_{2}$}

\begin{equation*}
\sum_{\substack{ a,d=0  \\ a+d\equiv 1}}^{1}B(gx_{2}\otimes
1_{H};G^{a}X_{1}X_{2},g^{d}x_{1}x_{2})G^{a}X_{1}\otimes g^{d}x_{1}x_{2}=0
\end{equation*}%
and we get%
\begin{equation*}
B(gx_{2}\otimes 1_{H};X_{1}X_{2},gx_{1}x_{2})=0
\end{equation*}%
and%
\begin{equation*}
B(gx_{2}\otimes 1_{H};GX_{1}X_{2},x_{1}x_{2})=0
\end{equation*}%
which we already know from case $gx_{1}.$

\subsubsection{$G^{a}X_{1}X_{2}\otimes g^{d}x_{2}$}

There is no term like this.

\subsubsection{$G^{a}X_{1}X_{2}\otimes g^{d}x_{1}$}

\begin{equation*}
\sum_{\substack{ a,d=0  \\ a+d\equiv 1}}^{1}\left( -1\right)
^{a+1}B(gx_{2}\otimes
1_{H};G^{a}X_{1}X_{2},g^{d}x_{1}x_{2})G^{a}X_{1}X_{2}\otimes g^{d}x_{1}=0
\end{equation*}%
and we get%
\begin{equation*}
-B(gx_{2}\otimes 1_{H};X_{1}X_{2},gx_{1}x_{2})=0
\end{equation*}%
and%
\begin{equation*}
B(gx_{2}\otimes 1_{H};GX_{1}X_{2},x_{1}x_{2})=0
\end{equation*}%
which we already got

\subsubsection{$G^{a}X_{1}X_{2}\otimes g^{d}x_{1}x_{2}$}

There is no term like this.

\subsection{Case $gx_{1}x_{2}$}

\begin{eqnarray*}
&&\sum_{a,b_{1},b_{2},d,e_{1},e_{2}=0}^{1}\sum_{l_{1}=0}^{b_{1}}%
\sum_{l_{2}=0}^{b_{2}}\sum_{u_{1}=0}^{e_{1}}\sum_{u_{2}=0}^{e_{2}}\left(
-1\right) ^{\alpha \left( g;l_{1},l_{2},u_{1},u_{2}\right) } \\
&&B(gx_{2}\otimes
1_{H};G^{a}X_{1}^{b_{1}}X_{2}^{b_{2}},g^{d}x_{1}^{e_{1}}x_{2}^{e_{2}}) \\
&&G^{a}X_{1}^{b_{1}-l_{1}}X_{2}^{b_{2}-l_{2}}\otimes
g^{d}x_{1}^{e_{1}-u_{1}}x_{2}^{e_{2}-u_{2}}\otimes
g^{a+b_{1}+b_{2}+l_{1}+l_{2}+d+e_{1}+e_{2}+u_{1}+u_{2}+1}x_{1}^{l_{1}+u_{1}}x_{2}^{l_{2}+u_{2}}=0
\end{eqnarray*}%
\begin{eqnarray*}
l_{1}+u_{1} &=&1 \\
l_{2}+u_{2} &=&1 \\
d+e_{1}+e_{2}+a+b_{1}+b_{2} &\equiv &0
\end{eqnarray*}%
By $\left( \ref{gx2ot1, first}\right) $ we already know that $%
B(gx_{2}\otimes
1_{H};G^{a}X_{1}^{b_{1}}X_{2}^{b_{2}},g^{d}x_{1}^{e_{1}}x_{2}^{e_{2}})=0$
whenever $d+e_{1}+e_{2}+a+b_{1}+b_{2}\equiv 0.$

By using all the equalities we got above from $\left( \ref{gx2ot1, first}%
\right) $ to $\left( \ref{gx2ot1, thirtyseven}\right) $ we obtain the
following form of $B\left( gx_{2}\otimes 1_{H}\right) .$

\subsection{The final form of the element $B\left( gx_{2}\otimes
1_{H}\right) $}

\begin{eqnarray*}
B\left( gx_{2}\otimes 1_{H}\right) &=&B(gx_{2}\otimes
1_{H};1_{A},g)1_{A}\otimes g+ \\
&&+B(gx_{2}\otimes 1_{H};1_{A},x_{1})1_{A}\otimes x_{1} \\
&&+B(gx_{2}\otimes 1_{H};1_{A},x_{2})1_{A}\otimes x_{2} \\
&&+B(gx_{2}\otimes 1_{H};1_{A},gx_{1}x_{2})1_{A}\otimes gx_{1}x_{2} \\
&&+B(gx_{2}\otimes 1_{H};G,1_{H})G\otimes 1_{H}+ \\
&&+B(gx_{2}\otimes 1_{H};G,x_{1}x_{2})G\otimes x_{1}x_{2}+ \\
&&+B(gx_{2}\otimes 1_{H};G,gx_{1})G\otimes gx_{1}+ \\
&&+B(gx_{2}\otimes 1_{H};G,gx_{2})G\otimes gx_{2}+ \\
&&+B(gx_{2}\otimes 1_{H};1_{A},x_{1})X_{1}\otimes 1_{H}+ \\
&&-B(gx_{2}\otimes 1_{H};1_{A},gx_{1}x_{2})X_{1}\otimes gx_{2}+ \\
&&\left[ -1+B(gx_{2}\otimes 1_{H};1_{A},x_{2})\right] X_{2}\otimes 1_{H}+ \\
&&+B(gx_{2}\otimes 1_{H};1_{A},gx_{1}x_{2})X_{2}\otimes gx_{1}+ \\
&&+B(gx_{2}\otimes 1_{H};1_{A},gx_{1}x_{2})X_{1}X_{2}\otimes g+ \\
&&-B(gx_{2}\otimes 1_{H};G,gx_{1})GX_{1}\otimes g+ \\
&&+B(gx_{2}\otimes 1_{H};G,x_{1}x_{2})GX_{1}\otimes x_{2}+ \\
&&-B(gx_{2}\otimes 1_{H};G,gx_{2})GX_{2}\otimes g+ \\
&&-B(gx_{2}\otimes 1_{H};G,x_{1}x_{2})GX_{2}\otimes x_{1}+ \\
&&+B(gx_{2}\otimes 1_{H};G,x_{1}x_{2})GX_{1}X_{2}\otimes 1_{H}
\end{eqnarray*}

\section{$B(x_{1}x_{2}\otimes 1_{H})$}

We write the Casimir formula $\left( \ref{MAIN FORMULA 1}\right) $ for $%
B(x_{1}x_{2}\otimes 1_{H}).$
\begin{eqnarray*}
&&\sum_{a,b_{1},b_{2},d,e_{1},e_{2}=0}^{1}\sum_{l_{1}=0}^{b_{1}}%
\sum_{l_{2}=0}^{b_{2}}\sum_{u_{1}=0}^{e_{1}}\sum_{u_{2}=0}^{e_{2}}\left(
-1\right) ^{\alpha \left( x_{1}x_{2};l_{1},l_{2},u_{1},u_{2}\right) } \\
&&B(1_{H}\otimes
1_{H};G^{a}X_{1}^{b_{1}}X_{2}^{b_{2}},g^{d}x_{1}^{e_{1}}x_{2}^{e_{2}}) \\
&&G^{a}X_{1}^{b_{1}-l_{1}}X_{2}^{b_{2}-l_{2}}\otimes
g^{d}x_{1}^{e_{1}-u_{1}}x_{2}^{e_{2}-u_{2}}\otimes
g^{a+b_{1}+b_{2}+l_{1}+l_{2}+d+e_{1}+e_{2}+u_{1}+u_{2}}x_{1}^{l_{1}+u_{1}+1}x_{2}^{l_{2}+u_{2}+1}
\\
&&\sum_{a,b_{1},b_{2},d,e_{1},e_{2}=0}^{1}\sum_{l_{1}=0}^{b_{1}}%
\sum_{l_{2}=0}^{b_{2}}\sum_{u_{1}=0}^{e_{1}}\sum_{u_{2}=0}^{e_{2}}\left(
-1\right) ^{\alpha \left( x_{1};l_{1},l_{2},u_{1},u_{2}\right) } \\
&&B(gx_{2}\otimes
1_{H};G^{a}X_{1}^{b_{1}}X_{2}^{b_{2}},g^{d}x_{1}^{e_{1}}x_{2}^{e_{2}}) \\
&&G^{a}X_{1}^{b_{1}-l_{1}}X_{2}^{b_{2}-l_{2}}\otimes
g^{d}x_{1}^{e_{1}-u_{1}}x_{2}^{e_{2}-u_{2}}\otimes
g^{a+b_{1}+b_{2}+l_{1}+l_{2}+d+e_{1}+e_{2}+u_{1}+u_{2}}x_{1}^{l_{1}+u_{1}+1}x_{2}^{l_{2}+u_{2}}
\\
&&-\sum_{a,b_{1},b_{2},d,e_{1},e_{2}=0}^{1}\sum_{l_{1}=0}^{b_{1}}%
\sum_{l_{2}=0}^{b_{2}}\sum_{u_{1}=0}^{e_{1}}\sum_{u_{2}=0}^{e_{2}}\left(
-1\right) ^{\alpha \left( x_{2};l_{1},l_{2},u_{1},u_{2}\right) } \\
&&B(gx_{1}\otimes
1_{H};G^{a}X_{1}^{b_{1}}X_{2}^{b_{2}},g^{d}x_{1}^{e_{1}}x_{2}^{e_{2}}) \\
&&G^{a}X_{1}^{b_{1}-l_{1}}X_{2}^{b_{2}-l_{2}}\otimes
g^{d}x_{1}^{e_{1}-u_{1}}x_{2}^{e_{2}-u_{2}}\otimes
g^{a+b_{1}+b_{2}+l_{1}+l_{2}+d+e_{1}+e_{2}+u_{1}+u_{2}}x_{1}^{l_{1}+u_{1}}x_{2}^{l_{2}+u_{2}+1}
\\
&&+\sum_{a,b_{1},b_{2},d,e_{1},e_{2}=0}^{1}\sum_{l_{1}=0}^{b_{1}}%
\sum_{l_{2}=0}^{b_{2}}\sum_{u_{1}=0}^{e_{1}}\sum_{u_{2}=0}^{e_{2}}\left(
-1\right) ^{\alpha \left( 1_{H};l_{1},l_{2},u_{1},u_{2}\right) } \\
&&B(x_{1}x_{2}\otimes
1_{H};G^{a}X_{1}^{b_{1}}X_{2}^{b_{2}},g^{d}x_{1}^{e_{1}}x_{2}^{e_{2}}) \\
&&G^{a}X_{1}^{b_{1}-l_{1}}X_{2}^{b_{2}-l_{2}}\otimes
g^{d}x_{1}^{e_{1}-u_{1}}x_{2}^{e_{2}-u_{2}}\otimes
g^{a+b_{1}+b_{2}+l_{1}+l_{2}+d+e_{1}+e_{2}+u_{1}+u_{2}}x_{1}^{l_{1}+u_{1}}x_{2}^{l_{2}+u_{2}}
\\
&=&B^{A}(x_{1}x_{2}\otimes 1_{H})\otimes B^{H}(x_{1}x_{2}\otimes
1_{H})\otimes 1_{H}
\end{eqnarray*}

\subsection{Case $1_{H}$}

\begin{eqnarray}
&&+\sum_{a,b_{1},b_{2},d,e_{1},e_{2}=0}^{1}\sum_{l_{1}=0}^{b_{1}}%
\sum_{l_{2}=0}^{b_{2}}\sum_{u_{1}=0}^{e_{1}}\sum_{u_{2}=0}^{e_{2}}\left(
-1\right) ^{\alpha \left( 1_{H};l_{1},l_{2},u_{1},u_{2}\right) }
\label{eqX1X2} \\
&&B(x_{1}x_{2}\otimes
1_{H};G^{a}X_{1}^{b_{1}}X_{2}^{b_{2}},g^{d}x_{1}^{e_{1}}x_{2}^{e_{2}})
\notag \\
&&G^{a}X_{1}^{b_{1}-l_{1}}X_{2}^{b_{2}-l_{2}}\otimes
g^{d}x_{1}^{e_{1}-u_{1}}x_{2}^{e_{2}-u_{2}}\otimes
g^{a+b_{1}+b_{2}+l_{1}+l_{2}+d+e_{1}+e_{2}+u_{1}+u_{2}}x_{1}^{l_{1}+u_{1}}x_{2}^{l_{2}+u_{2}}
\notag \\
&=&B^{A}(x_{1}x_{2}\otimes 1_{H})\otimes B^{H}(x_{1}x_{2}\otimes
1_{H})\otimes 1_{H}  \notag
\end{eqnarray}%
\begin{eqnarray*}
a+b_{1}+b_{2}+l_{1}+l_{2}+d+e_{1}+e_{2}+u_{1}+u_{2} &\equiv &0 \\
l_{1}+u_{1} &\equiv &0 \\
l_{2}+u_{2} &\equiv &0
\end{eqnarray*}%
and we get%
\begin{equation}
B(x_{1}x_{2}\otimes
1_{H};G^{a}X_{1}^{b_{1}}X_{2}^{b_{2}},g^{d}x_{1}^{e_{1}}x_{2}^{e_{2}})=0%
\text{ whenever }a+b_{1}+b_{2}+d+e_{1}+e_{2}\equiv 1  \label{x1otx2, first}
\end{equation}

\subsection{Case $g$}

\begin{eqnarray*}
&&+\sum_{a,b_{1},b_{2},d,e_{1},e_{2}=0}^{1}\sum_{l_{1}=0}^{b_{1}}%
\sum_{l_{2}=0}^{b_{2}}\sum_{u_{1}=0}^{e_{1}}\sum_{u_{2}=0}^{e_{2}}\left(
-1\right) ^{\alpha \left( 1_{H};l_{1},l_{2},u_{1},u_{2}\right) } \\
&&B(x_{1}x_{2}\otimes
1_{H};G^{a}X_{1}^{b_{1}}X_{2}^{b_{2}},g^{d}x_{1}^{e_{1}}x_{2}^{e_{2}}) \\
&&G^{a}X_{1}^{b_{1}-l_{1}}X_{2}^{b_{2}-l_{2}}\otimes
g^{d}x_{1}^{e_{1}-u_{1}}x_{2}^{e_{2}-u_{2}}\otimes
g^{a+b_{1}+b_{2}+l_{1}+l_{2}+d+e_{1}+e_{2}+u_{1}+u_{2}}x_{1}^{l_{1}+u_{1}}x_{2}^{l_{2}+u_{2}}=0
\end{eqnarray*}%
\begin{eqnarray*}
a+b_{1}+b_{2}+d+e_{1}+e_{2} &\equiv &1 \\
l_{1}+u_{1} &\equiv &0 \\
l_{2}+u_{2} &\equiv &0
\end{eqnarray*}%
and the equality holds in view of $\left( \ref{x1x2ot1}\right) .$

\subsection{Case $x_{1}$}

\begin{eqnarray*}
&&\sum_{a,b_{1},b_{2},d,e_{1},e_{2}=0}^{1}\sum_{l_{1}=0}^{b_{1}}%
\sum_{l_{2}=0}^{b_{2}}\sum_{u_{1}=0}^{e_{1}}\sum_{u_{2}=0}^{e_{2}}\left(
-1\right) ^{\alpha \left( x_{1};l_{1},l_{2},u_{1},u_{2}\right) } \\
&&B(gx_{2}\otimes
1_{H};G^{a}X_{1}^{b_{1}}X_{2}^{b_{2}},g^{d}x_{1}^{e_{1}}x_{2}^{e_{2}}) \\
&&G^{a}X_{1}^{b_{1}-l_{1}}X_{2}^{b_{2}-l_{2}}\otimes
g^{d}x_{1}^{e_{1}-u_{1}}x_{2}^{e_{2}-u_{2}}\otimes
g^{a+b_{1}+b_{2}+l_{1}+l_{2}+d+e_{1}+e_{2}+u_{1}+u_{2}}x_{1}^{l_{1}+u_{1}+1}x_{2}^{l_{2}+u_{2}}
\\
&&+\sum_{a,b_{1},b_{2},d,e_{1},e_{2}=0}^{1}\sum_{l_{1}=0}^{b_{1}}%
\sum_{l_{2}=0}^{b_{2}}\sum_{u_{1}=0}^{e_{1}}\sum_{u_{2}=0}^{e_{2}}\left(
-1\right) ^{\alpha \left( 1_{H};l_{1},l_{2},u_{1},u_{2}\right) } \\
&&B(x_{1}x_{2}\otimes
1_{H};G^{a}X_{1}^{b_{1}}X_{2}^{b_{2}},g^{d}x_{1}^{e_{1}}x_{2}^{e_{2}}) \\
&&G^{a}X_{1}^{b_{1}-l_{1}}X_{2}^{b_{2}-l_{2}}\otimes
g^{d}x_{1}^{e_{1}-u_{1}}x_{2}^{e_{2}-u_{2}}\otimes
g^{a+b_{1}+b_{2}+l_{1}+l_{2}+d+e_{1}+e_{2}+u_{1}+u_{2}}x_{1}^{l_{1}+u_{1}}x_{2}^{l_{2}+u_{2}}
\\
&=&0
\end{eqnarray*}%
For the first summand we get

\begin{eqnarray*}
a+b_{1}+b_{2}+l_{1}+l_{2}+d+e_{1}+e_{2}+u_{1}+u_{2} &\equiv &0 \\
l_{1}+u_{1} &=&0 \\
l_{2}+u_{2} &=&0
\end{eqnarray*}%
and we get%
\begin{eqnarray*}
a+b_{1}+b_{2}+d+e_{1}+e_{2} &\equiv &0 \\
l_{1} &=&u_{1}=0 \\
l_{2} &=&u_{2}=0
\end{eqnarray*}%
By $\left( \ref{gx2ot1, first}\right) ,$ we have%
\begin{equation*}
B(gx_{2}\otimes \
1_{H};G^{a}X_{1}^{b_{1}}X_{2}^{b_{2}},g^{d}x_{1}^{e_{1}}x_{2}^{e_{2}})=0%
\text{ for }a+b_{1}+b_{2}+d+e_{1}+e_{2}\equiv 0
\end{equation*}%
so the first summand is zero. For the second summand we have

\begin{eqnarray*}
a+b_{1}+b_{2}+l_{1}+l_{2}+d+e_{1}+e_{2}+u_{1}+u_{2} &\equiv &0 \\
l_{1}+u_{1} &=&1 \\
l_{2}+u_{2} &=&0
\end{eqnarray*}%
and we get%
\begin{eqnarray*}
a+b_{1}+b_{2}+d+e_{1}+e_{2} &\equiv &1 \\
l_{1}+u_{1} &=&1 \\
l_{2} &=&u_{2}=0
\end{eqnarray*}

and the equality holds in view of $\left( \ref{x1x2ot1}\right) .$

\subsection{Case $x_{2}$}

\begin{eqnarray*}
&&-\sum_{a,b_{1},b_{2},d,e_{1},e_{2}=0}^{1}\sum_{l_{1}=0}^{b_{1}}%
\sum_{l_{2}=0}^{b_{2}}\sum_{u_{1}=0}^{e_{1}}\sum_{u_{2}=0}^{e_{2}}\left(
-1\right) ^{\alpha \left( x_{2};l_{1},l_{2},u_{1},u_{2}\right) } \\
&&B(gx_{1}\otimes
1_{H};G^{a}X_{1}^{b_{1}}X_{2}^{b_{2}},g^{d}x_{1}^{e_{1}}x_{2}^{e_{2}}) \\
&&G^{a}X_{1}^{b_{1}-l_{1}}X_{2}^{b_{2}-l_{2}}\otimes
g^{d}x_{1}^{e_{1}-u_{1}}x_{2}^{e_{2}-u_{2}}\otimes
g^{a+b_{1}+b_{2}+l_{1}+l_{2}+d+e_{1}+e_{2}+u_{1}+u_{2}}x_{1}^{l_{1}+u_{1}}x_{2}^{l_{2}+u_{2}+1}
\\
&&+\sum_{a,b_{1},b_{2},d,e_{1},e_{2}=0}^{1}\sum_{l_{1}=0}^{b_{1}}%
\sum_{l_{2}=0}^{b_{2}}\sum_{u_{1}=0}^{e_{1}}\sum_{u_{2}=0}^{e_{2}}\left(
-1\right) ^{\alpha \left( 1_{H};l_{1},l_{2},u_{1},u_{2}\right) } \\
&&B(x_{1}x_{2}\otimes
1_{H};G^{a}X_{1}^{b_{1}}X_{2}^{b_{2}},g^{d}x_{1}^{e_{1}}x_{2}^{e_{2}}) \\
&&G^{a}X_{1}^{b_{1}-l_{1}}X_{2}^{b_{2}-l_{2}}\otimes
g^{d}x_{1}^{e_{1}-u_{1}}x_{2}^{e_{2}-u_{2}}\otimes
g^{a+b_{1}+b_{2}+l_{1}+l_{2}+d+e_{1}+e_{2}+u_{1}+u_{2}}x_{1}^{l_{1}+u_{1}}x_{2}^{l_{2}+u_{2}}
\\
&=&0
\end{eqnarray*}%
For the first summand we have%
\begin{eqnarray*}
a+b_{1}+b_{2}+l_{1}+l_{2}+d+e_{1}+e_{2}+u_{1}+u_{2} &\equiv &0 \\
l_{1}+u_{1} &=&0 \\
l_{2}+u_{2} &=&0
\end{eqnarray*}%
so that we get

\begin{eqnarray*}
a+b_{1}+b_{2}+d+e_{1}+e_{2} &\equiv &0 \\
l_{1} &=&u_{1}=0 \\
l_{2} &=&u_{2}=0
\end{eqnarray*}%
In view of $\left( \ref{gx1ot1, first}\right) ,$ we know that $%
B(gx_{1}\otimes
1_{H};G^{a}X_{1}^{b_{1}}X_{2}^{b_{2}},g^{d}x_{1}^{e_{1}}x_{2}^{e_{2}})=0$
whenever $a+b_{1}+b_{2}+d+e_{1}+e_{2}\equiv 0$ so that the first summand is
zero. Let us consider the second summand.

\begin{eqnarray*}
a+b_{1}+b_{2}+l_{1}+l_{2}+d+e_{1}+e_{2}+u_{1}+u_{2} &\equiv &0 \\
l_{1}+u_{1} &=&0 \\
l_{2}+u_{2} &=&1
\end{eqnarray*}%
and we get%
\begin{eqnarray*}
a+b_{1}+b_{2}+d+e_{1}+e_{2} &\equiv &1 \\
l_{1} &=&u_{1}=0 \\
l_{2}+u_{2} &=&1
\end{eqnarray*}%
and the equality holds in view of $\left( \ref{x1x2ot1}\right) .$

\subsection{Case $x_{1}x_{2}$}

\begin{eqnarray*}
&&\sum_{a,b_{1},b_{2},d,e_{1},e_{2}=0}^{1}\sum_{l_{1}=0}^{b_{1}}%
\sum_{l_{2}=0}^{b_{2}}\sum_{u_{1}=0}^{e_{1}}\sum_{u_{2}=0}^{e_{2}}\left(
-1\right) ^{\alpha \left( x_{1}x_{2};l_{1},l_{2},u_{1},u_{2}\right) } \\
&&B(1_{H}\otimes
1_{H};G^{a}X_{1}^{b_{1}}X_{2}^{b_{2}},g^{d}x_{1}^{e_{1}}x_{2}^{e_{2}}) \\
&&G^{a}X_{1}^{b_{1}-l_{1}}X_{2}^{b_{2}-l_{2}}\otimes
g^{d}x_{1}^{e_{1}-u_{1}}x_{2}^{e_{2}-u_{2}}\otimes
g^{a+b_{1}+b_{2}+l_{1}+l_{2}+d+e_{1}+e_{2}+u_{1}+u_{2}}x_{1}^{l_{1}+u_{1}+1}x_{2}^{l_{2}+u_{2}+1}
\\
&&\sum_{a,b_{1},b_{2},d,e_{1},e_{2}=0}^{1}\sum_{l_{1}=0}^{b_{1}}%
\sum_{l_{2}=0}^{b_{2}}\sum_{u_{1}=0}^{e_{1}}\sum_{u_{2}=0}^{e_{2}}\left(
-1\right) ^{\alpha \left( x_{1};l_{1},l_{2},u_{1},u_{2}\right) } \\
&&B(gx_{2}\otimes
1_{H};G^{a}X_{1}^{b_{1}}X_{2}^{b_{2}},g^{d}x_{1}^{e_{1}}x_{2}^{e_{2}}) \\
&&G^{a}X_{1}^{b_{1}-l_{1}}X_{2}^{b_{2}-l_{2}}\otimes
g^{d}x_{1}^{e_{1}-u_{1}}x_{2}^{e_{2}-u_{2}}\otimes
g^{a+b_{1}+b_{2}+l_{1}+l_{2}+d+e_{1}+e_{2}+u_{1}+u_{2}}x_{1}^{l_{1}+u_{1}+1}x_{2}^{l_{2}+u_{2}}
\\
&&-\sum_{a,b_{1},b_{2},d,e_{1},e_{2}=0}^{1}\sum_{l_{1}=0}^{b_{1}}%
\sum_{l_{2}=0}^{b_{2}}\sum_{u_{1}=0}^{e_{1}}\sum_{u_{2}=0}^{e_{2}}\left(
-1\right) ^{\alpha \left( x_{2};l_{1},l_{2},u_{1},u_{2}\right) } \\
&&B(gx_{1}\otimes
1_{H};G^{a}X_{1}^{b_{1}}X_{2}^{b_{2}},g^{d}x_{1}^{e_{1}}x_{2}^{e_{2}}) \\
&&G^{a}X_{1}^{b_{1}-l_{1}}X_{2}^{b_{2}-l_{2}}\otimes
g^{d}x_{1}^{e_{1}-u_{1}}x_{2}^{e_{2}-u_{2}}\otimes
g^{a+b_{1}+b_{2}+l_{1}+l_{2}+d+e_{1}+e_{2}+u_{1}+u_{2}}x_{1}^{l_{1}+u_{1}}x_{2}^{l_{2}+u_{2}+1}
\\
&&+\sum_{a,b_{1},b_{2},d,e_{1},e_{2}=0}^{1}\sum_{l_{1}=0}^{b_{1}}%
\sum_{l_{2}=0}^{b_{2}}\sum_{u_{1}=0}^{e_{1}}\sum_{u_{2}=0}^{e_{2}}\left(
-1\right) ^{\alpha \left( 1_{H};l_{1},l_{2},u_{1},u_{2}\right) } \\
&&B(x_{1}x_{2}\otimes
1_{H};G^{a}X_{1}^{b_{1}}X_{2}^{b_{2}},g^{d}x_{1}^{e_{1}}x_{2}^{e_{2}}) \\
&&G^{a}X_{1}^{b_{1}-l_{1}}X_{2}^{b_{2}-l_{2}}\otimes
g^{d}x_{1}^{e_{1}-u_{1}}x_{2}^{e_{2}-u_{2}}\otimes
g^{a+b_{1}+b_{2}+l_{1}+l_{2}+d+e_{1}+e_{2}+u_{1}+u_{2}}x_{1}^{l_{1}+u_{1}}x_{2}^{l_{2}+u_{2}}
\\
&=&B^{A}(x_{1}x_{2}\otimes 1_{H})\otimes B^{H}(x_{1}x_{2}\otimes
1_{H})\otimes 1_{H}
\end{eqnarray*}%
The first summand gives us $1_{A}\otimes 1_{H}\otimes x_{1}x_{2}.$Let us
consider the second summand.
\begin{eqnarray*}
&&\sum_{a,b_{1},b_{2},d,e_{1},e_{2}=0}^{1}\sum_{l_{1}=0}^{b_{1}}%
\sum_{l_{2}=0}^{b_{2}}\sum_{u_{1}=0}^{e_{1}}\sum_{u_{2}=0}^{e_{2}}\left(
-1\right) ^{\alpha \left( x_{1};l_{1},l_{2},u_{1},u_{2}\right) } \\
&&B(gx_{2}\otimes
1_{H};G^{a}X_{1}^{b_{1}}X_{2}^{b_{2}},g^{d}x_{1}^{e_{1}}x_{2}^{e_{2}}) \\
&&G^{a}X_{1}^{b_{1}-l_{1}}X_{2}^{b_{2}-l_{2}}\otimes
g^{d}x_{1}^{e_{1}-u_{1}}x_{2}^{e_{2}-u_{2}}\otimes
g^{a+b_{1}+b_{2}+l_{1}+l_{2}+d+e_{1}+e_{2}+u_{1}+u_{2}}x_{1}^{l_{1}+u_{1}+1}x_{2}^{l_{2}+u_{2}}
\end{eqnarray*}%
We get%
\begin{eqnarray*}
a+b_{1}+b_{2}+l_{1}+l_{2}+d+e_{1}+e_{2}+u_{1}+u_{2} &\equiv &0 \\
l_{1}+u_{1}+1 &=&1 \\
l_{2}+u_{2} &=&1
\end{eqnarray*}%
i.e.

\begin{eqnarray*}
a+b_{1}+b_{2}+d+e_{1}+e_{2} &\equiv &1 \\
l_{1} &=&u_{1}=0 \\
l_{2}+u_{2} &=&1
\end{eqnarray*}%
and we obtain%
\begin{eqnarray*}
&&\sum_{\substack{ a,b_{1},b_{2},d,e_{1},e_{2}=0  \\ %
a+b_{1}+b_{2}+d+e_{1}+e_{2}\equiv 1}}^{1}\sum_{l_{2}=0}^{b_{2}}\sum
_{\substack{ u_{2}=0  \\ l_{2}+u_{2}=1}}^{e_{2}}\left( -1\right) ^{\alpha
\left( x_{1};0,l_{2},0,u_{2}\right) }B(gx_{2}\otimes
1_{H};G^{a}X_{1}^{b_{1}}X_{2}^{b_{2}},g^{d}x_{1}^{e_{1}}x_{2}^{e_{2}}) \\
&&G^{a}X_{1}^{b_{1}}X_{2}^{b_{2}-l_{2}}\otimes
g^{d}x_{1}^{e_{1}}x_{2}^{e_{2}-u_{2}}.
\end{eqnarray*}

\begin{eqnarray*}
&&\sum_{\substack{ a,b_{1},b_{2},d,e_{1}=0  \\ a+b_{1}+b_{2}+d+e_{1}\equiv 0
}}^{1}-B(gx_{2}\otimes
1_{H};G^{a}X_{1}^{b_{1}}X_{2}^{b_{2}},g^{d}x_{1}^{e_{1}}x_{2})G^{a}X_{1}^{b_{1}}X_{2}^{b_{2}}\otimes g^{d}x_{1}^{e_{1}}+
\\
&&\sum_{\substack{ a,b_{1},d,e_{1},e_{2}=0  \\ a+b_{1}+d+e_{1}+e_{2}\equiv 0
}}^{1}\left( -1\right) ^{a+b_{1}}B(gx_{2}\otimes
1_{H};G^{a}X_{1}^{b_{1}}X_{2},g^{d}x_{1}^{e_{1}}x_{2}^{e_{2}})G^{a}X_{1}^{b_{1}}\otimes g^{d}x_{1}^{e_{1}}x_{2}^{e_{2}}.
\end{eqnarray*}%
Let us consider the third summand.%
\begin{eqnarray*}
&&-\sum_{a,b_{1},b_{2},d,e_{1},e_{2}=0}^{1}\sum_{l_{1}=0}^{b_{1}}%
\sum_{l_{2}=0}^{b_{2}}\sum_{u_{1}=0}^{e_{1}}\sum_{u_{2}=0}^{e_{2}}\left(
-1\right) ^{\alpha \left( x_{2};l_{1},l_{2},u_{1},u_{2}\right) } \\
&&B(gx_{1}\otimes
1_{H};G^{a}X_{1}^{b_{1}}X_{2}^{b_{2}},g^{d}x_{1}^{e_{1}}x_{2}^{e_{2}}) \\
&&G^{a}X_{1}^{b_{1}-l_{1}}X_{2}^{b_{2}-l_{2}}\otimes
g^{d}x_{1}^{e_{1}-u_{1}}x_{2}^{e_{2}-u_{2}}\otimes
g^{a+b_{1}+b_{2}+l_{1}+l_{2}+d+e_{1}+e_{2}+u_{1}+u_{2}}x_{1}^{l_{1}+u_{1}}x_{2}^{l_{2}+u_{2}+1}
\end{eqnarray*}

\begin{eqnarray*}
a+b_{1}+b_{2}+l_{1}+l_{2}+d+e_{1}+e_{2}+u_{1}+u_{2} &\equiv &0 \\
l_{1}+u_{1} &=&1 \\
l_{2}+u_{2}+1 &=&1
\end{eqnarray*}%
i.e.

\begin{eqnarray*}
a+b_{1}+b_{2}+d+e_{1}+e_{2} &\equiv &1 \\
l_{1}+u_{1} &=&1 \\
l_{2} &=&u_{2}=0
\end{eqnarray*}%
\begin{eqnarray*}
&&-\sum_{\substack{ a,b_{1},b_{2},d,e_{1},e_{2}=0  \\ %
a+b_{1}+b_{2}+d+e_{1}+e_{2}\equiv 1}}^{1}\sum_{l_{1}=0}^{b_{1}}\sum
_{\substack{ u_{1}=0  \\ l_{1}+u_{1}=1}}^{e_{1}}\left( -1\right) ^{\alpha
\left( x_{2};l_{1},0,u_{1},0\right) }B(gx_{1}\otimes
1_{H};G^{a}X_{1}^{b_{1}}X_{2}^{b_{2}},g^{d}x_{1}^{e_{1}}x_{2}^{e_{2}}) \\
&&G^{a}X_{1}^{b_{1}-l_{1}}X_{2}^{b_{2}}\otimes
g^{d}x_{1}^{e_{1}-u_{1}}x_{2}^{e_{2}}
\end{eqnarray*}

Since $\alpha \left( x_{2};0,0,1,0\right) =e_{2}$ and $\alpha \left(
x_{2};1,0,0,0\right) =a+b_{1}$ we obtain%
\begin{eqnarray*}
&&-\sum_{\substack{ a,b_{1},b_{2},d,e_{2}=0  \\ a+b_{1}+b_{2}+d+e_{2}\equiv
0 }}^{1}\left( -1\right) ^{e_{2}}B(gx_{1}\otimes
1_{H};G^{a}X_{1}^{b_{1}}X_{2}^{b_{2}},g^{d}x_{1}x_{2}^{e_{2}})G^{a}X_{1}^{b_{1}}X_{2}^{b_{2}}\otimes g^{d}x_{2}^{e_{2}}+
\\
&&-\sum_{\substack{ a,b_{2},d,e_{1},e_{2}=0  \\ a+b_{2}+d+e_{1}+e_{2}\equiv
0 }}^{1}\left( -1\right) ^{a+1}B(gx_{1}\otimes
1_{H};G^{a}X_{1}X_{2}^{b_{2}},g^{d}x_{1}^{e_{1}}x_{2}^{e_{2}})G^{a}X_{2}^{b_{2}}\otimes g^{d}x_{1}^{e_{1}}x_{2}^{e_{2}}.
\end{eqnarray*}%
Let us consider the fourth term.

\begin{eqnarray*}
&&+\sum_{a,b_{1},b_{2},d,e_{1},e_{2}=0}^{1}\sum_{l_{1}=0}^{b_{1}}%
\sum_{l_{2}=0}^{b_{2}}\sum_{u_{1}=0}^{e_{1}}\sum_{u_{2}=0}^{e_{2}}\left(
-1\right) ^{\alpha \left( 1_{H};l_{1},l_{2},u_{1},u_{2}\right) } \\
&&B(x_{1}x_{2}\otimes
1_{H};G^{a}X_{1}^{b_{1}}X_{2}^{b_{2}},g^{d}x_{1}^{e_{1}}x_{2}^{e_{2}}) \\
&&G^{a}X_{1}^{b_{1}-l_{1}}X_{2}^{b_{2}-l_{2}}\otimes
g^{d}x_{1}^{e_{1}-u_{1}}x_{2}^{e_{2}-u_{2}}\otimes
g^{a+b_{1}+b_{2}+l_{1}+l_{2}+d+e_{1}+e_{2}+u_{1}+u_{2}}x_{1}^{l_{1}+u_{1}}x_{2}^{l_{2}+u_{2}}
\end{eqnarray*}%
\begin{eqnarray*}
a+b_{1}+b_{2}+d+e_{1}+e_{2} &\equiv &0 \\
l_{1}+u_{1} &=&1 \\
l_{2}+u_{2} &=&1
\end{eqnarray*}

\bigskip

\begin{eqnarray*}
&&\sum_{\substack{ a,b_{1},b_{2},d,e_{1},e_{2}=0  \\ %
a+b_{1}+b_{2}+d+e_{1}+e_{2}\equiv 0}}^{1}\sum_{l_{1}=0}^{b_{1}}%
\sum_{l_{2}=0}^{b_{2}}\sum_{\substack{ u_{1}=0  \\ l_{1}+u_{1}=1}}%
^{e_{1}}\sum _{\substack{ u_{2}=0  \\ l_{2}+u_{2}=1}}^{e_{2}}\left(
-1\right) ^{\alpha \left( 1_{H};l_{1},l_{2},u_{1},u_{2}\right) } \\
&&B(x_{1}x_{2}\otimes
1_{H};G^{a}X_{1}^{b_{1}}X_{2}^{b_{2}},g^{d}x_{1}^{e_{1}}x_{2}^{e_{2}})G^{a}X_{1}^{b_{1}-l_{1}}X_{2}^{b_{2}-l_{2}}\otimes g^{d}x_{1}^{e_{1}-u_{1}}x_{2}^{e_{2}-u_{2}}
\end{eqnarray*}%
and, since $\alpha \left( 1_{H};0,0,1,1\right) =1+e_{2},$ $\alpha \left(
1_{H};0,1,1,0\right) =e_{2}+\left( a+b_{1}+b_{2}+1\right) ,$ $\alpha \left(
1_{H};1,0,0,1\right) =a+b_{1}$ and $\alpha \left( 1_{H};1,1,0,0\right)
=1+b_{2},$ we get%
\begin{eqnarray*}
&&\sum_{\substack{ a,b_{1},b_{2},d,=0  \\ a+b_{1}+b_{2}+d\equiv 0}}%
^{1}B(x_{1}x_{2}\otimes
1_{H};G^{a}X_{1}^{b_{1}}X_{2}^{b_{2}},g^{d}x_{1}x_{2})G^{a}X_{1}^{b_{1}}X_{2}^{b_{2}}\otimes g^{d}+
\\
&&\sum_{\substack{ a,b_{1},d,e_{2}=0  \\ a+b_{1}+d+e_{2}\equiv 0}}^{1}\left(
-1\right) ^{e_{2}+a+b_{1}}B(x_{1}x_{2}\otimes
1_{H};G^{a}X_{1}^{b_{1}}X_{2},g^{d}x_{1}x_{2}^{e_{2}})G^{a}X_{1}^{b_{1}}%
\otimes g^{d}x_{2}^{e_{2}}+ \\
&&\sum_{\substack{ a,b_{2},d,e_{1}=0  \\ a+b_{2}+d+e_{1}\equiv 0}}^{1}\left(
-1\right) ^{a+1}B(x_{1}x_{2}\otimes
1_{H};G^{a}X_{1}X_{2}^{b_{2}},g^{d}x_{1}^{e_{1}}x_{2})G^{a}X_{2}^{b_{2}}%
\otimes g^{d}x_{1}^{e_{1}}+ \\
&&\sum_{\substack{ a,d,e_{1},e_{2}=0  \\ a+d+e_{1}+e_{2}\equiv 0}}%
^{1}B(x_{1}x_{2}\otimes
1_{H};G^{a}X_{1}X_{2},g^{d}x_{1}^{e_{1}}x_{2}^{e_{2}})G^{a}\otimes
g^{d}x_{1}^{e_{1}}x_{2}^{e_{2}}.
\end{eqnarray*}%
Summing up we get

\begin{gather*}
1_{A}\otimes 1_{H}+ \\
\sum_{\substack{ a,b_{1},b_{2},d,e_{1}=0  \\ a+b_{1}+b_{2}+d+e_{1}\equiv 0}}%
^{1}-B(gx_{2}\otimes
1_{H};G^{a}X_{1}^{b_{1}}X_{2}^{b_{2}},g^{d}x_{1}^{e_{1}}x_{2})G^{a}X_{1}^{b_{1}}X_{2}^{b_{2}}\otimes g^{d}x_{1}^{e_{1}}+
\\
\sum_{\substack{ a,b_{1},d,e_{1},e_{2}=0  \\ a+b_{1}+d+e_{1}+e_{2}\equiv 0}}%
^{1}\left( -1\right) ^{a+b_{1}}B(gx_{2}\otimes
1_{H};G^{a}X_{1}^{b_{1}}X_{2},g^{d}x_{1}^{e_{1}}x_{2}^{e_{2}})G^{a}X_{1}^{b_{1}}\otimes g^{d}x_{1}^{e_{1}}x_{2}^{e_{2}}+
\\
-\sum_{\substack{ a,b_{1},b_{2},d,e_{2}=0  \\ a+b_{1}+b_{2}+d+e_{2}\equiv 0}}%
^{1}\left( -1\right) ^{e_{2}}B(gx_{1}\otimes
1_{H};G^{a}X_{1}^{b_{1}}X_{2}^{b_{2}},g^{d}x_{1}x_{2}^{e_{2}})G^{a}X_{1}^{b_{1}}X_{2}^{b_{2}}\otimes g^{d}x_{2}^{e_{2}}+
\\
-\sum_{\substack{ a,b_{2},d,e_{1},e_{2}=0  \\ a+b_{2}+d+e_{1}+e_{2}\equiv 0}}%
^{1}\left( -1\right) ^{a+1}B(gx_{1}\otimes
1_{H};G^{a}X_{1}X_{2}^{b_{2}},g^{d}x_{1}^{e_{1}}x_{2}^{e_{2}})G^{a}X_{2}^{b_{2}}\otimes g^{d}x_{1}^{e_{1}}x_{2}^{e_{2}}+
\\
\sum_{\substack{ a,b_{1},b_{2},d,=0  \\ a+b_{1}+b_{2}+d\equiv 0}}%
^{1}B(x_{1}x_{2}\otimes
1_{H};G^{a}X_{1}^{b_{1}}X_{2}^{b_{2}},g^{d}x_{1}x_{2})G^{a}X_{1}^{b_{1}}X_{2}^{b_{2}}\otimes g^{d}+
\\
\sum_{\substack{ a,b_{1},d,e_{2}=0  \\ a+b_{1}+d+e_{2}\equiv 0}}^{1}\left(
-1\right) ^{e_{2}+a+b_{1}}B(x_{1}x_{2}\otimes
1_{H};G^{a}X_{1}^{b_{1}}X_{2},g^{d}x_{1}x_{2}^{e_{2}})G^{a}X_{1}^{b_{1}}%
\otimes g^{d}x_{2}^{e_{2}}+ \\
\sum_{\substack{ a,b_{2},d,e_{1}=0  \\ a+b_{2}+d+e_{1}\equiv 0}}^{1}\left(
-1\right) ^{a+1}B(x_{1}x_{2}\otimes
1_{H};G^{a}X_{1}X_{2}^{b_{2}},g^{d}x_{1}^{e_{1}}x_{2})G^{a}X_{2}^{b_{2}}%
\otimes g^{d}x_{1}^{e_{1}}+ \\
\sum_{\substack{ a,d,e_{1},e_{2}=0  \\ a+d+e_{1}+e_{2}\equiv 0}}%
^{1}B(x_{1}x_{2}\otimes
1_{H};G^{a}X_{1}X_{2},g^{d}x_{1}^{e_{1}}x_{2}^{e_{2}})G^{a}\otimes
g^{d}x_{1}^{e_{1}}x_{2}^{e_{2}}=0
\end{gather*}

\subsubsection{$G^{a}\otimes g^{d}$}

\begin{gather*}
1_{A}\otimes 1_{H}+ \\
\sum_{\substack{ a,d=0  \\ a+d\equiv 0}}^{1}\left[
\begin{array}{c}
-B(gx_{2}\otimes 1_{H};G^{a},g^{d}x_{2})+\left( -1\right)
^{a}B(gx_{2}\otimes 1_{H};G^{a}X_{2},g^{d})+ \\
-B(gx_{1}\otimes 1_{H};G^{a},g^{d}x_{1})+\left( -1\right)
^{a}B(gx_{1}\otimes 1_{H};G^{a}X_{1},g^{d})+ \\
+B(x_{1}x_{2}\otimes 1_{H};G^{a},g^{d}x_{1}x_{2})+\left( -1\right)
^{a}B(x_{1}x_{2}\otimes 1_{H};G^{a}X_{2},g^{d}x_{1})+ \\
+\left( -1\right) ^{a+1}B(x_{1}x_{2}\otimes
1_{H};G^{a}X_{1},g^{d}x_{2})+B(x_{1}x_{2}\otimes 1_{H};G^{a}X_{1}X_{2},g^{d})%
\end{array}%
\right] G^{a}\otimes g^{d}=0
\end{gather*}%
and we get%
\begin{equation*}
\begin{array}{c}
1-B(gx_{2}\otimes 1_{H};1_{A},x_{2})+B(gx_{2}\otimes 1_{H};X_{2},1_{H})+ \\
-B(gx_{1}\otimes 1_{H};1_{A},x_{1})+B(gx_{1}\otimes 1_{H};X_{1},1_{H})+ \\
+B(x_{1}x_{2}\otimes 1_{H};1_{A},x_{1}x_{2})+B(x_{1}x_{2}\otimes
1_{H};X_{2},x_{1})+ \\
-B(x_{1}x_{2}\otimes 1_{H};X_{1},x_{2})+B(x_{1}x_{2}\otimes
1_{H};X_{1}X_{2},1_{H})%
\end{array}%
=0
\end{equation*}%
and%
\begin{equation*}
\begin{array}{c}
-B(gx_{2}\otimes 1_{H};G,gx_{2})-B(gx_{2}\otimes 1_{H};GX_{2},g)+ \\
-B(gx_{1}\otimes 1_{H};G,gx_{1})-B(gx_{1}\otimes 1_{H};GX_{1},g)+ \\
+B(x_{1}x_{2}\otimes 1_{H};G,gx_{1}x_{2})-B(x_{1}x_{2}\otimes
1_{H};GX_{2},gx_{1})+ \\
+B(x_{1}x_{2}\otimes 1_{H};GX_{1},gx_{2})+B(x_{1}x_{2}\otimes
1_{H};GX_{1}X_{2},g)%
\end{array}%
=0.
\end{equation*}%
By using the form of the elements $B\left( gx_{1}\otimes 1_{H}\right) $ and $%
B\left( gx_{2}\otimes 1_{H}\right) $ we get%
\begin{equation*}
\begin{array}{c}
-B(gx_{2}\otimes 1_{H};1_{A},x_{2})+B(gx_{2}\otimes 1_{H};1_{A},x_{2})+ \\
-B(gx_{1}\otimes 1_{H};1_{A},x_{1})-1+B(gx_{1}\otimes 1_{H};1_{H},x_{1})+ \\
+B(x_{1}x_{2}\otimes 1_{H};1_{A},x_{1}x_{2})+B(x_{1}x_{2}\otimes
1_{H};X_{2},x_{1})+ \\
-B(x_{1}x_{2}\otimes 1_{H};X_{1},x_{2})+B(x_{1}x_{2}\otimes
1_{H};X_{1}X_{2},1_{H})%
\end{array}%
=0
\end{equation*}%
i.e.%
\begin{equation}
\begin{array}{c}
-1+B(x_{1}x_{2}\otimes 1_{H};1_{A},x_{1}x_{2})+B(x_{1}x_{2}\otimes
1_{H};X_{2},x_{1})+ \\
-B(x_{1}x_{2}\otimes 1_{H};X_{1},x_{2})+B(x_{1}x_{2}\otimes
1_{H};X_{1}X_{2},1_{H})%
\end{array}%
=0  \label{x1otx2, second}
\end{equation}%
and%
\begin{equation*}
\begin{array}{c}
-B(gx_{2}\otimes 1_{H};G,gx_{2})+B(gx_{2}\otimes 1_{H};G,gx_{2})+ \\
-B(gx_{1}\otimes 1_{H};G,gx_{1})+B(gx_{1}\otimes 1_{H};G,gx_{1})+ \\
+B(x_{1}x_{2}\otimes 1_{H};G,gx_{1}x_{2})-B(x_{1}x_{2}\otimes
1_{H};GX_{2},gx_{1})+ \\
+B(x_{1}x_{2}\otimes 1_{H};GX_{1},gx_{2})+B(x_{1}x_{2}\otimes
1_{H};GX_{1}X_{2},g)%
\end{array}%
=0
\end{equation*}%
i.e.%
\begin{equation}
\begin{array}{c}
+B(x_{1}x_{2}\otimes 1_{H};G,gx_{1}x_{2})-B(x_{1}x_{2}\otimes
1_{H};GX_{2},gx_{1})+ \\
+B(x_{1}x_{2}\otimes 1_{H};GX_{1},gx_{2})+B(x_{1}x_{2}\otimes
1_{H};GX_{1}X_{2},g)%
\end{array}%
=0  \label{x1otx2, third}
\end{equation}

\subsubsection{$G^{a}\otimes g^{d}x_{2}$}

\begin{gather*}
\sum_{\substack{ a,d=0  \\ a+d\equiv 1}}^{1}\left[
\begin{array}{c}
\left( -1\right) ^{a}B(gx_{2}\otimes
1_{H};G^{a}X_{2},g^{d}x_{2})+B(gx_{1}\otimes 1_{H};G^{a},g^{d}x_{1}x_{2})+
\\
+\left( -1\right) ^{a}B(gx_{1}\otimes 1_{H};G^{a}X_{1},g^{d}x_{2})+\left(
-1\right) ^{a+1}B(x_{1}x_{2}\otimes 1_{H};G^{a}X_{2},g^{d}x_{1}x_{2})+ \\
+B(x_{1}x_{2}\otimes 1_{H};G^{a}X_{1}X_{2},g^{d}x_{2}%
\end{array}%
\right] \\
G^{a}\otimes g^{d}x_{2}=0
\end{gather*}%
and we get%
\begin{equation*}
\begin{array}{c}
B(gx_{2}\otimes 1_{H};X_{2},gx_{2})+B(gx_{1}\otimes 1_{H};1_{A},gx_{1}x_{2})+
\\
+B(gx_{1}\otimes 1_{H};X_{1},gx_{2})-B(x_{1}x_{2}\otimes
1_{H};X_{2},gx_{1}x_{2})+ \\
+B(x_{1}x_{2}\otimes 1_{H};X_{1}X_{2},gx_{2})%
\end{array}%
=0
\end{equation*}%
and%
\begin{equation*}
\begin{array}{c}
-B(gx_{2}\otimes 1_{H};GX_{2},x_{2})+B(gx_{1}\otimes 1_{H};G,x_{1}x_{2})+ \\
-B(gx_{1}\otimes 1_{H};GX_{1},x_{2})+B(x_{1}x_{2}\otimes
1_{H};GX_{2},x_{1}x_{2})+ \\
+B(x_{1}x_{2}\otimes 1_{H};GX_{1}X_{2},x_{2})%
\end{array}%
=0
\end{equation*}%
By using the form of the elements $B\left( gx_{1}\otimes 1_{H}\right) $ and $%
B\left( gx_{2}\otimes 1_{H}\right) $ we get%
\begin{equation*}
\begin{array}{c}
B(gx_{1}\otimes 1_{H};1_{A},gx_{1}x_{2})+ \\
-B(gx_{1}\otimes 1_{H};1_{A},gx_{1}x_{2})-B(x_{1}x_{2}\otimes
1_{H};X_{2},gx_{1}x_{2})+ \\
+B(x_{1}x_{2}\otimes 1_{H};X_{1}X_{2},gx_{2})%
\end{array}%
=0
\end{equation*}%
i.e.%
\begin{equation}
-B(x_{1}x_{2}\otimes 1_{H};X_{2},gx_{1}x_{2})+B(x_{1}x_{2}\otimes
1_{H};X_{1}X_{2},gx_{2})=0  \label{x1otx2, four}
\end{equation}%
and%
\begin{equation*}
\begin{array}{c}
+B(gx_{1}\otimes 1_{H};G,x_{1}x_{2})+ \\
-B(gx_{1}\otimes 1_{H};G,x_{1}x_{2})+B(x_{1}x_{2}\otimes
1_{H};GX_{2},x_{1}x_{2})+ \\
+B(x_{1}x_{2}\otimes 1_{H};GX_{1}X_{2},x_{2})%
\end{array}%
=0
\end{equation*}%
i.e.%
\begin{equation}
B(x_{1}x_{2}\otimes 1_{H};GX_{2},x_{1}x_{2})+B(x_{1}x_{2}\otimes
1_{H};GX_{1}X_{2},x_{2})=0  \label{x1otx2, five}
\end{equation}

\subsubsection{$G^{a}\otimes g^{d}x_{1}$}

\begin{gather*}
\sum_{\substack{ a,d=0  \\ a+d\equiv 1}}^{1}\left[
\begin{array}{c}
-B(gx_{2}\otimes 1_{H};G^{a},g^{d}x_{1}x_{2})+\left( -1\right)
^{a}B(gx_{2}\otimes 1_{H};G^{a}X_{2},g^{d}x_{1})+ \\
+\left( -1\right) ^{a}B(gx_{1}\otimes 1_{H};G^{a}X_{1},g^{d}x_{1})+\left(
-1\right) ^{a+1}B(x_{1}x_{2}\otimes 1_{H};G^{a}X_{1},g^{d}x_{1}x_{2})+ \\
+B(x_{1}x_{2}\otimes 1_{H};G^{a}X_{1}X_{2},g^{d}x_{1})%
\end{array}%
\right] \\
G^{a}\otimes g^{d}x_{1}=0
\end{gather*}%
and we get%
\begin{equation*}
\begin{array}{c}
-B(gx_{2}\otimes 1_{H};1_{A},gx_{1}x_{2})+B(gx_{2}\otimes
1_{H};X_{2},gx_{1})+ \\
+B(gx_{1}\otimes 1_{H};X_{1},gx_{1})-B(x_{1}x_{2}\otimes
1_{H};X_{1},gx_{1}x_{2})+ \\
+B(x_{1}x_{2}\otimes 1_{H};X_{1}X_{2},gx_{1})=0%
\end{array}%
\end{equation*}%
and%
\begin{equation*}
\begin{array}{c}
-B(gx_{2}\otimes 1_{H};G,x_{1}x_{2})-B(gx_{2}\otimes 1_{H};GX_{2},x_{1})+ \\
-B(gx_{1}\otimes 1_{H};GX_{1},x_{1})+B(x_{1}x_{2}\otimes
1_{H};GX_{1},x_{1}x_{2})+ \\
+B(x_{1}x_{2}\otimes 1_{H};GX_{1}X_{2},x_{1})=0.%
\end{array}%
\end{equation*}%
By using the form of the elements $B\left( gx_{1}\otimes 1_{H}\right) $ and $%
B\left( gx_{2}\otimes 1_{H}\right) $ we get%
\begin{equation*}
\begin{array}{c}
-B(gx_{2}\otimes 1_{H};1_{A},gx_{1}x_{2})+B(gx_{2}\otimes
1_{H};1_{A},gx_{1}x_{2})+ \\
-B(x_{1}x_{2}\otimes 1_{H};X_{1},gx_{1}x_{2})+ \\
+B(x_{1}x_{2}\otimes 1_{H};X_{1}X_{2},gx_{1})=0%
\end{array}%
\end{equation*}%
i.e.%
\begin{equation}
-B(x_{1}x_{2}\otimes 1_{H};X_{1},gx_{1}x_{2})+B(x_{1}x_{2}\otimes
1_{H};X_{1}X_{2},gx_{1})=0  \label{x1otx2, six}
\end{equation}%
and%
\begin{equation*}
\begin{array}{c}
-B(gx_{2}\otimes 1_{H};G,x_{1}x_{2})+B(gx_{2}\otimes 1_{H};G,x_{1}x_{2})+ \\
+B(x_{1}x_{2}\otimes 1_{H};GX_{1},x_{1}x_{2})+ \\
+B(x_{1}x_{2}\otimes 1_{H};GX_{1}X_{2},x_{1})=0.%
\end{array}%
\end{equation*}%
i.e.%
\begin{equation}
+B(x_{1}x_{2}\otimes 1_{H};GX_{1},x_{1}x_{2})+B(x_{1}x_{2}\otimes
1_{H};GX_{1}X_{2},x_{1})=0.  \label{x1otx2, seven}
\end{equation}

\subsubsection{$G^{a}X_{2}\otimes g^{d}$}

\begin{gather*}
\sum_{\substack{ a,d=0  \\ a+d\equiv 1}}^{1}\left[
\begin{array}{c}
-B(gx_{2}\otimes 1_{H};G^{a}X_{2},g^{d}x_{2})-B(gx_{1}\otimes
1_{H};G^{a}X_{2},g^{d}x_{1})+ \\
+\left( -1\right) ^{a}B(gx_{1}\otimes
1_{H};G^{a}X_{1}X_{2},g^{d})+B(x_{1}x_{2}\otimes
1_{H};G^{a}X_{2},g^{d}x_{1}x_{2})+ \\
+\left( -1\right) ^{a+1}B(x_{1}x_{2}\otimes 1_{H};G^{a}X_{1}X_{2},g^{d}x_{2})%
\end{array}%
\right] \\
G^{a}X_{2}\otimes g^{d}=0
\end{gather*}%
and we get%
\begin{equation*}
\begin{array}{c}
-B(gx_{2}\otimes 1_{H};X_{2},gx_{2})-B(gx_{1}\otimes 1_{H};X_{2},gx_{1})+ \\
+B(gx_{1}\otimes 1_{H};X_{1}X_{2},g)+B(x_{1}x_{2}\otimes
1_{H};X_{2},gx_{1}x_{2})+ \\
-B(x_{1}x_{2}\otimes 1_{H};X_{1}X_{2},gx_{2})%
\end{array}%
=0
\end{equation*}%
and%
\begin{equation*}
\begin{array}{c}
-B(gx_{2}\otimes 1_{H};GX_{2},x_{2})-B(gx_{1}\otimes 1_{H};GX_{2},x_{1})+ \\
-B(gx_{1}\otimes 1_{H};GX_{1}X_{2},1_{H})+B(x_{1}x_{2}\otimes
1_{H};GX_{2},x_{1}x_{2})+ \\
+B(x_{1}x_{2}\otimes 1_{H};GX_{1}X_{2},x_{2})%
\end{array}%
=0.
\end{equation*}%
By using the form of the elements $B\left( gx_{1}\otimes 1_{H}\right) $ and $%
B\left( gx_{2}\otimes 1_{H}\right) $ we get%
\begin{equation*}
\begin{array}{c}
-B(gx_{1}\otimes 1_{H};1_{A},gx_{1}x_{2})+ \\
+B(gx_{2}\otimes 1_{H};1_{A},gx_{1}x_{2})+B(x_{1}x_{2}\otimes
1_{H};X_{2},gx_{1}x_{2})+ \\
-B(x_{1}x_{2}\otimes 1_{H};X_{1}X_{2},gx_{2})%
\end{array}%
=0
\end{equation*}%
i.e.%
\begin{equation*}
B(x_{1}x_{2}\otimes 1_{H};X_{2},gx_{1}x_{2})-B(x_{1}x_{2}\otimes
1_{H};X_{1}X_{2},gx_{2})=0
\end{equation*}%
which is $\left( \ref{x1otx2, four}\right) $ and%
\begin{equation*}
\begin{array}{c}
B(gx_{1}\otimes 1_{H};G,x_{1}x_{2})+ \\
-B(gx_{1}\otimes 1_{H};G,x_{1}x_{2})+B(x_{1}x_{2}\otimes
1_{H};GX_{2},x_{1}x_{2})+ \\
+B(x_{1}x_{2}\otimes 1_{H};GX_{1}X_{2},x_{2})%
\end{array}%
=0.
\end{equation*}%
i.e.%
\begin{equation*}
B(x_{1}x_{2}\otimes 1_{H};GX_{2},x_{1}x_{2})+B(x_{1}x_{2}\otimes
1_{H};GX_{1}X_{2},x_{2})=0
\end{equation*}%
this is $\left( \ref{x1otx2, five}\right) .$

\subsubsection{$G^{a}X_{1}\otimes g^{d}$}

\begin{gather*}
\sum_{\substack{ a,d=0  \\ a+d\equiv 1}}^{1}\left[
\begin{array}{c}
-B(gx_{2}\otimes 1_{H};G^{a}X_{1},g^{d}x_{2})+\left( -1\right)
^{a+1}B(gx_{2}\otimes 1_{H};G^{a}X_{1}X_{2},g^{d})+ \\
-B(gx_{1}\otimes 1_{H};G^{a}X_{1},g^{d}x_{1})+B(x_{1}x_{2}\otimes
1_{H};G^{a}X_{1},g^{d}x_{1}x_{2})+ \\
\left( -1\right) ^{a+1}B(x_{1}x_{2}\otimes 1_{H};G^{a}X_{1}X_{2},g^{d}x_{1})%
\end{array}%
\right] \\
G^{a}X_{1}\otimes g^{d}=0
\end{gather*}%
and we get%
\begin{equation*}
\begin{array}{c}
-B(gx_{2}\otimes 1_{H};X_{1},gx_{2})-B(gx_{2}\otimes 1_{H};X_{1}X_{2},g)+ \\
-B(gx_{1}\otimes 1_{H};X_{1},gx_{1})+B(x_{1}x_{2}\otimes
1_{H};X_{1},gx_{1}x_{2})+ \\
-B(x_{1}x_{2}\otimes 1_{H};X_{1}X_{2},gx_{1})%
\end{array}%
=0
\end{equation*}%
and%
\begin{equation*}
\begin{array}{c}
-B(gx_{2}\otimes 1_{H};GX_{1},x_{2})+B(gx_{2}\otimes
1_{H};GX_{1}X_{2},1_{H})+ \\
-B(gx_{1}\otimes 1_{H};GX_{1},x_{1})+B(x_{1}x_{2}\otimes
1_{H};GX_{1},x_{1}x_{2})+ \\
+B(x_{1}x_{2}\otimes 1_{H};GX_{1}X_{2},x_{1})%
\end{array}%
=0.
\end{equation*}%
By using the form of the elements $B\left( gx_{1}\otimes 1_{H}\right) $ and $%
B\left( gx_{2}\otimes 1_{H}\right) $ we get%
\begin{equation*}
\begin{array}{c}
B(gx_{2}\otimes 1_{H};1_{A},gx_{1}x_{2})-B(gx_{2}\otimes
1_{H};1_{A},gx_{1}x_{2})+ \\
+B(x_{1}x_{2}\otimes 1_{H};X_{1},gx_{1}x_{2})+ \\
-B(x_{1}x_{2}\otimes 1_{H};X_{1}X_{2},gx_{1})%
\end{array}%
=0
\end{equation*}%
i.e.%
\begin{equation*}
+B(x_{1}x_{2}\otimes 1_{H};X_{1},gx_{1}x_{2})-B(x_{1}x_{2}\otimes
1_{H};X_{1}X_{2},gx_{1})=0
\end{equation*}%
which is $\left( \ref{x1otx2, six}\right) $ and%
\begin{equation*}
\begin{array}{c}
-B(gx_{2}\otimes 1_{H};G,x_{1}x_{2})GX_{1}+B(gx_{2}\otimes
1_{H};G,x_{1}x_{2})+ \\
+B(x_{1}x_{2}\otimes 1_{H};GX_{1},x_{1}x_{2})+B(x_{1}x_{2}\otimes
1_{H};GX_{1}X_{2},x_{1})%
\end{array}%
=0.
\end{equation*}%
i.e.%
\begin{equation*}
B(x_{1}x_{2}\otimes 1_{H};GX_{1},x_{1}x_{2})+B(x_{1}x_{2}\otimes
1_{H};GX_{1}X_{2},x_{1})=0
\end{equation*}%
which is $\left( \ref{x1otx2, seven}\right) .$

\subsubsection{$G^{a}\otimes g^{d}x_{1}x_{2}$}

\begin{gather*}
\sum_{\substack{ a,d=0  \\ a+d\equiv 0}}^{1}\left[
\begin{array}{c}
\left( -1\right) ^{a}B(gx_{2}\otimes
1_{H};G^{a}X_{2},g^{d}x_{1}x_{2})+\left( -1\right) ^{a}B(gx_{1}\otimes
1_{H};G^{a}X_{1},g^{d}x_{1}x_{2})+ \\
+B(x_{1}x_{2}\otimes 1_{H};G^{a}X_{1}X_{2},g^{d}x_{1}x_{2})%
\end{array}%
\right] \\
G^{a}\otimes g^{d}x_{1}x_{2}=0
\end{gather*}%
and we get%
\begin{equation*}
\begin{array}{c}
B(gx_{2}\otimes 1_{H};X_{2},x_{1}x_{2})+B(gx_{1}\otimes
1_{H};X_{1},x_{1}x_{2})+ \\
+B(x_{1}x_{2}\otimes 1_{H};X_{1}X_{2},x_{1}x_{2})%
\end{array}%
=0
\end{equation*}%
and%
\begin{equation*}
\begin{array}{c}
-B(gx_{2}\otimes 1_{H};GX_{2},gx_{1}x_{2})-B(gx_{1}\otimes
1_{H};GX_{1},gx_{1}x_{2})+ \\
+B(x_{1}x_{2}\otimes 1_{H};GX_{1}X_{2},gx_{1}x_{2})%
\end{array}%
=0.
\end{equation*}%
By using the form of the elements $B\left( gx_{1}\otimes 1_{H}\right) $ and $%
B\left( gx_{2}\otimes 1_{H}\right) $ we get%
\begin{equation}
B(x_{1}x_{2}\otimes 1_{H};X_{1}X_{2},x_{1}x_{2})=0  \label{x1otx2, eight}
\end{equation}%
and%
\begin{equation}
B(x_{1}x_{2}\otimes 1_{H};GX_{1}X_{2},gx_{1}x_{2})=0.  \label{x1otx2, nine}
\end{equation}

\subsubsection{$G^{a}X_{2}\otimes g^{d}x_{2}$}

\begin{equation*}
\sum_{\substack{ a,d=0  \\ a+d\equiv 0}}^{1}\left[ B(gx_{1}\otimes
1_{H};G^{a}X_{2},g^{d}x_{1}x_{2})+\left( -1\right) ^{a}B(gx_{1}\otimes
1_{H};G^{a}X_{1}X_{2},g^{d}x_{2})\right] G^{a}X_{2}\otimes g^{d}x_{2}=0
\end{equation*}%
and we get%
\begin{equation*}
B(gx_{1}\otimes 1_{H};X_{2},x_{1}x_{2})+B(gx_{1}\otimes
1_{H};X_{1}X_{2},x_{2})=0
\end{equation*}%
\begin{equation*}
B(gx_{1}\otimes 1_{H};GX_{2},gx_{1}x_{2})-B(gx_{1}\otimes
1_{H};GX_{1}X_{2},gx_{2})=0.
\end{equation*}%
By the form of the element $B(gx_{1}\otimes 1_{H})$ we know that each
summand in equalities above is zero.

\subsubsection{$G^{a}X_{1}\otimes g^{d}x_{2}$}

\begin{gather*}
\sum_{\substack{ a,d=0  \\ a+d\equiv 0}}^{1}\left[
\begin{array}{c}
\left( -1\right) ^{a+1}B(gx_{2}\otimes
1_{H};G^{a}X_{1}X_{2},g^{d}x_{2})+B(gx_{1}\otimes
1_{H};G^{a}X_{1},g^{d}x_{1}x_{2})+ \\
+\left( -1\right) ^{a}B(x_{1}x_{2}\otimes
1_{H};G^{a}X_{1}X_{2},g^{d}x_{1}x_{2})%
\end{array}%
\right] \\
G^{a}X_{1}\otimes g^{d}x_{2}=0
\end{gather*}%
and we get%
\begin{equation*}
\begin{array}{c}
-B(gx_{2}\otimes 1_{H};X_{1}X_{2},x_{2})+B(gx_{1}\otimes
1_{H};X_{1},x_{1}x_{2})+ \\
+B(x_{1}x_{2}\otimes 1_{H};X_{1}X_{2},x_{1}x_{2})%
\end{array}%
=0
\end{equation*}%
\begin{equation*}
\begin{array}{c}
B(gx_{2}\otimes 1_{H};GX_{1}X_{2},gx_{2})+B(gx_{1}\otimes
1_{H};GX_{1},gx_{1}x_{2})+ \\
-B(x_{1}x_{2}\otimes 1_{H};GX_{1}X_{2},gx_{1}x_{2})%
\end{array}%
=0.
\end{equation*}%
By using the form of the elements $B\left( gx_{1}\otimes 1_{H}\right) $ and $%
B\left( gx_{2}\otimes 1_{H}\right) $ we get%
\begin{equation*}
B(x_{1}x_{2}\otimes 1_{H};X_{1}X_{2},x_{1}x_{2})=0
\end{equation*}%
and%
\begin{equation*}
B(x_{1}x_{2}\otimes 1_{H};GX_{1}X_{2},gx_{1}x_{2})=0
\end{equation*}%
which we already got.

\subsubsection{$G^{a}X_{2}\otimes g^{d}x_{1}$}

\begin{gather*}
\sum_{\substack{ a,d=0  \\ a+d\equiv 0}}^{1}\left[
\begin{array}{c}
-B(gx_{2}\otimes 1_{H};G^{a}X_{2},g^{d}x_{1}x_{2})+\left( -1\right)
^{a}B(gx_{1}\otimes 1_{H};G^{a}X_{1}X_{2},g^{d}x_{1})+ \\
+\left( -1\right) ^{a+1}B(x_{1}x_{2}\otimes
1_{H};G^{a}X_{1}X_{2},g^{d}x_{1}x_{2})%
\end{array}%
\right] \\
G^{a}X_{2}\otimes g^{d}x_{1}=0
\end{gather*}%
and we get%
\begin{equation*}
\begin{array}{c}
-B(gx_{2}\otimes 1_{H};X_{2},x_{1}x_{2})+B(gx_{1}\otimes
1_{H};X_{1}X_{2},x_{1})+ \\
-B(x_{1}x_{2}\otimes 1_{H};X_{1}X_{2},x_{1}x_{2})%
\end{array}%
=0
\end{equation*}%
and%
\begin{equation*}
\begin{array}{c}
-B(gx_{2}\otimes 1_{H};GX_{2},gx_{1}x_{2})-B(gx_{1}\otimes
1_{H};GX_{1}X_{2},gx_{1})+ \\
+B(x_{1}x_{2}\otimes 1_{H};GX_{1}X_{2},gx_{1}x_{2})%
\end{array}%
=0.
\end{equation*}%
By using the form of the elements $B\left( gx_{1}\otimes 1_{H}\right) $ and $%
B\left( gx_{2}\otimes 1_{H}\right) $ we get%
\begin{equation*}
B(x_{1}x_{2}\otimes 1_{H};X_{1}X_{2},x_{1}x_{2})=0
\end{equation*}%
and%
\begin{equation*}
B(x_{1}x_{2}\otimes 1_{H};GX_{1}X_{2},gx_{1}x_{2})=0
\end{equation*}%
which we already got.

\subsubsection{$G^{a}X_{1}\otimes g^{d}x_{1}$}

\begin{gather*}
\sum_{\substack{ a,d=0  \\ a+d\equiv 0}}^{1}\left[ -B(gx_{2}\otimes
1_{H};G^{a}X_{1},g^{d}x_{1}x_{2})+\left( -1\right) ^{a+1}B(gx_{2}\otimes
1_{H};G^{a}X_{1}X_{2},g^{d}x_{1})\right] \\
G^{a}X_{1}\otimes g^{d}x_{1}=0
\end{gather*}%
and we get%
\begin{equation*}
-B(gx_{2}\otimes 1_{H};X_{1},x_{1}x_{2})-B(gx_{2}\otimes
1_{H};X_{1}X_{2},x_{1})=0
\end{equation*}%
and%
\begin{equation*}
-B(gx_{2}\otimes 1_{H};GX_{1},gx_{1}x_{2})+B(gx_{2}\otimes
1_{H};GX_{1}X_{2},gx_{1})=0.
\end{equation*}%
By the form of the element $B(gx_{2}\otimes 1_{H})$ we know that each
summand in equalities above is zero.

\subsubsection{$G^{a}X_{1}X_{2}\otimes g^{d}$}

\begin{gather*}
\sum_{\substack{ a,d=0  \\ a+d\equiv 0}}^{1}\left[
\begin{array}{c}
-B(gx_{2}\otimes 1_{H};G^{a}X_{1}X_{2},g^{d}x_{2})-B(gx_{1}\otimes
1_{H};G^{a}X_{1}X_{2},g^{d}x_{1})+ \\
+B(x_{1}x_{2}\otimes 1_{H};G^{a}X_{1}X_{2},g^{d}x_{1}x_{2})%
\end{array}%
\right] \\
G^{a}X_{1}X_{2}\otimes g^{d}=0
\end{gather*}%
and we get%
\begin{equation*}
\begin{array}{c}
-B(gx_{2}\otimes 1_{H};X_{1}X_{2},x_{2})-B(gx_{1}\otimes
1_{H};X_{1}X_{2},x_{1})+ \\
+B(x_{1}x_{2}\otimes 1_{H};X_{1}X_{2},x_{1}x_{2})%
\end{array}%
=0
\end{equation*}%
and%
\begin{equation*}
\begin{array}{c}
-B(gx_{2}\otimes 1_{H};GX_{1}X_{2},gx_{2})-B(gx_{1}\otimes
1_{H};GX_{1}X_{2},gx_{1})+ \\
+B(x_{1}x_{2}\otimes 1_{H};GX_{1}X_{2},gx_{1}x_{2})%
\end{array}%
=0.
\end{equation*}%
By using the form of the elements $B\left( gx_{1}\otimes 1_{H}\right) $ and $%
B\left( gx_{2}\otimes 1_{H}\right) $ we get%
\begin{equation*}
B(x_{1}x_{2}\otimes 1_{H};X_{1}X_{2},x_{1}x_{2})=0
\end{equation*}%
and%
\begin{equation*}
B(x_{1}x_{2}\otimes 1_{H};GX_{1}X_{2},gx_{1}x_{2})=0
\end{equation*}%
which we already got.

\subsubsection{$G^{a}X_{2}\otimes g^{d}x_{1}x_{2}$}

\begin{equation*}
\sum_{\substack{ a,d=0  \\ a+d\equiv 1}}^{1}\left( -1\right)
^{a}B(gx_{1}\otimes 1_{H};G^{a}X_{1}X_{2},g^{d}x_{1}x_{2})G^{a}X_{2}\otimes
g^{d}x_{1}x_{2}=0
\end{equation*}%
and we get%
\begin{equation*}
B(gx_{1}\otimes 1_{H};X_{1}X_{2},gx_{1}x_{2})=0
\end{equation*}%
and%
\begin{equation*}
-B(gx_{1}\otimes 1_{H};GX_{1}X_{2},x_{1}x_{2})=0
\end{equation*}%
By the form of the element $B(gx_{1}\otimes 1_{H})$ we know that these
equalities are already known.

\subsubsection{$G^{a}X_{1}\otimes g^{d}x_{1}x_{2}$}

\begin{equation*}
\sum_{\substack{ a,d=0  \\ a+d\equiv 1}}^{1}\left( -1\right)
^{a+1}B(gx_{2}\otimes
1_{H};G^{a}X_{1}X_{2},g^{d}x_{1}x_{2})G^{a}X_{1}\otimes g^{d}x_{1}x_{2}
\end{equation*}%
and we get%
\begin{equation*}
-B(gx_{2}\otimes 1_{H};X_{1}X_{2},gx_{1}x_{2})=0
\end{equation*}%
and%
\begin{equation*}
B(gx_{2}\otimes 1_{H};GX_{1}X_{2},x_{1}x_{2})=0
\end{equation*}%
By the form of the element $B(gx_{2}\otimes 1_{H})$ we know that these
equalities are already known.

\subsubsection{$G^{a}X_{1}X_{2}\otimes g^{d}x_{2}$}

\begin{equation*}
\sum_{\substack{ a,d=0  \\ a+d\equiv 1}}^{1}B(gx_{1}\otimes
1_{H};G^{a}X_{1}X_{2},g^{d}x_{1}x_{2})G^{a}X_{1}X_{2}\otimes g^{d}x_{2}=0
\end{equation*}%
and we get%
\begin{equation*}
B(gx_{1}\otimes 1_{H};X_{1}X_{2},gx_{1}x_{2})=0
\end{equation*}%
and%
\begin{equation*}
B(gx_{1}\otimes 1_{H};GX_{1}X_{2},x_{1}x_{2})=0
\end{equation*}%
By the form of the element $B(gx_{1}\otimes 1_{H})$ we know that these
equalities are already known.

\subsubsection{$G^{a}X_{1}X_{2}\otimes g^{d}x_{1}$}

\begin{equation*}
\sum_{\substack{ a,d=0  \\ a+d\equiv 1}}^{1}-B(gx_{2}\otimes
1_{H};G^{a}X_{1}X_{2},g^{d}x_{1}x_{2})G^{a}X_{1}X_{2}\otimes g^{d}x_{1}=0
\end{equation*}%
and we get%
\begin{equation*}
-B(gx_{2}\otimes 1_{H};X_{1}X_{2},gx_{1}x_{2})=0
\end{equation*}%
and%
\begin{equation*}
-B(gx_{2}\otimes 1_{H};GX_{1}X_{2},x_{1}x_{2})=0
\end{equation*}%
which we already got.

\subsubsection{$G^{a}X_{1}X_{2}\otimes g^{d}x_{1}x_{2}$}

There is no term like this.

\subsection{Case $gx_{1}$}

\begin{eqnarray*}
&&\sum_{a,b_{1},b_{2},d,e_{1},e_{2}=0}^{1}\sum_{l_{1}=0}^{b_{1}}%
\sum_{l_{2}=0}^{b_{2}}\sum_{u_{1}=0}^{e_{1}}\sum_{u_{2}=0}^{e_{2}}\left(
-1\right) ^{\alpha \left( x_{1}x_{2};l_{1},l_{2},u_{1},u_{2}\right) } \\
&&B(1_{H}\otimes
1_{H};G^{a}X_{1}^{b_{1}}X_{2}^{b_{2}},g^{d}x_{1}^{e_{1}}x_{2}^{e_{2}}) \\
&&G^{a}X_{1}^{b_{1}-l_{1}}X_{2}^{b_{2}-l_{2}}\otimes
g^{d}x_{1}^{e_{1}-u_{1}}x_{2}^{e_{2}-u_{2}}\otimes
g^{a+b_{1}+b_{2}+l_{1}+l_{2}+d+e_{1}+e_{2}+u_{1}+u_{2}}x_{1}^{l_{1}+u_{1}+1}x_{2}^{l_{2}+u_{2}+1}
\\
+
&&\sum_{a,b_{1},b_{2},d,e_{1},e_{2}=0}^{1}\sum_{l_{1}=0}^{b_{1}}%
\sum_{l_{2}=0}^{b_{2}}\sum_{u_{1}=0}^{e_{1}}\sum_{u_{2}=0}^{e_{2}}\left(
-1\right) ^{\alpha \left( x_{1};l_{1},l_{2},u_{1},u_{2}\right) } \\
&&B(gx_{2}\otimes
1_{H};G^{a}X_{1}^{b_{1}}X_{2}^{b_{2}},g^{d}x_{1}^{e_{1}}x_{2}^{e_{2}}) \\
&&G^{a}X_{1}^{b_{1}-l_{1}}X_{2}^{b_{2}-l_{2}}\otimes
g^{d}x_{1}^{e_{1}-u_{1}}x_{2}^{e_{2}-u_{2}}\otimes
g^{a+b_{1}+b_{2}+l_{1}+l_{2}+d+e_{1}+e_{2}+u_{1}+u_{2}}x_{1}^{l_{1}+u_{1}+1}x_{2}^{l_{2}+u_{2}}
\\
&&-\sum_{a,b_{1},b_{2},d,e_{1},e_{2}=0}^{1}\sum_{l_{1}=0}^{b_{1}}%
\sum_{l_{2}=0}^{b_{2}}\sum_{u_{1}=0}^{e_{1}}\sum_{u_{2}=0}^{e_{2}}\left(
-1\right) ^{\alpha \left( x_{2};l_{1},l_{2},u_{1},u_{2}\right) } \\
&&B(gx_{1}\otimes
1_{H};G^{a}X_{1}^{b_{1}}X_{2}^{b_{2}},g^{d}x_{1}^{e_{1}}x_{2}^{e_{2}}) \\
&&G^{a}X_{1}^{b_{1}-l_{1}}X_{2}^{b_{2}-l_{2}}\otimes
g^{d}x_{1}^{e_{1}-u_{1}}x_{2}^{e_{2}-u_{2}}\otimes
g^{a+b_{1}+b_{2}+l_{1}+l_{2}+d+e_{1}+e_{2}+u_{1}+u_{2}}x_{1}^{l_{1}+u_{1}}x_{2}^{l_{2}+u_{2}+1}
\\
&&+\sum_{a,b_{1},b_{2},d,e_{1},e_{2}=0}^{1}\sum_{l_{1}=0}^{b_{1}}%
\sum_{l_{2}=0}^{b_{2}}\sum_{u_{1}=0}^{e_{1}}\sum_{u_{2}=0}^{e_{2}}\left(
-1\right) ^{\alpha \left( 1_{H};l_{1},l_{2},u_{1},u_{2}\right) } \\
&&B(x_{1}x_{2}\otimes
1_{H};G^{a}X_{1}^{b_{1}}X_{2}^{b_{2}},g^{d}x_{1}^{e_{1}}x_{2}^{e_{2}}) \\
&&G^{a}X_{1}^{b_{1}-l_{1}}X_{2}^{b_{2}-l_{2}}\otimes
g^{d}x_{1}^{e_{1}-u_{1}}x_{2}^{e_{2}-u_{2}}\otimes
g^{a+b_{1}+b_{2}+l_{1}+l_{2}+d+e_{1}+e_{2}+u_{1}+u_{2}}x_{1}^{l_{1}+u_{1}}x_{2}^{l_{2}+u_{2}}
\\
&=&B^{A}(x_{1}x_{2}\otimes 1_{H})\otimes B^{H}(x_{1}x_{2}\otimes
1_{H})\otimes 1_{H}
\end{eqnarray*}%
We have to consider only the second and the fourth summand.

The second summand gives us

\begin{eqnarray*}
&&\sum_{a,b_{1},b_{2},d,e_{1},e_{2}=0}^{1}\sum_{l_{1}=0}^{b_{1}}%
\sum_{l_{2}=0}^{b_{2}}\sum_{u_{1}=0}^{e_{1}}\sum_{u_{2}=0}^{e_{2}}\left(
-1\right) ^{\alpha \left( x_{1};l_{1},l_{2},u_{1},u_{2}\right) } \\
&&B(gx_{2}\otimes
1_{H};G^{a}X_{1}^{b_{1}}X_{2}^{b_{2}},g^{d}x_{1}^{e_{1}}x_{2}^{e_{2}}) \\
&&G^{a}X_{1}^{b_{1}-l_{1}}X_{2}^{b_{2}-l_{2}}\otimes
g^{d}x_{1}^{e_{1}-u_{1}}x_{2}^{e_{2}-u_{2}}\otimes
g^{a+b_{1}+b_{2}+l_{1}+l_{2}+d+e_{1}+e_{2}+u_{1}+u_{2}}x_{1}^{l_{1}+u_{1}+1}x_{2}^{l_{2}+u_{2}}
\end{eqnarray*}%
We get%
\begin{eqnarray*}
a+b_{1}+b_{2}+l_{1}+l_{2}+d+e_{1}+e_{2}+u_{1}+u_{2} &\equiv &1 \\
l_{1}+u_{1}+1 &=&1 \\
l_{2}+u_{2} &=&0
\end{eqnarray*}%
i.e.%
\begin{eqnarray*}
a+b_{1}+b_{2}+d+e_{1}+e_{2} &\equiv &1 \\
l_{1} &=&u_{1}=0 \\
l_{2} &=&u_{2}=0
\end{eqnarray*}

and, since $\alpha \left( x_{1};0,0,0,0\right) =a+b_{1}+b_{2},$ we obtain%
\begin{equation*}
\sum_{\substack{ a,b_{1},b_{2},d,e_{1},e_{2}=0  \\ %
a+b_{1}+b_{2}+d+e_{1}+e_{2}\equiv 1}}^{1}\left( -1\right)
^{a+b_{1}+b_{2}}B(gx_{2}\otimes
1_{H};G^{a}X_{1}^{b_{1}}X_{2}^{b_{2}},g^{d}x_{1}^{e_{1}}x_{2}^{e_{2}})G^{a}X_{1}^{b_{1}}X_{2}^{b_{2}}\otimes g^{d}x_{1}^{e_{1}}x_{2}^{e_{2}}.
\end{equation*}%
Let us consider the fourth summand%
\begin{eqnarray*}
&&+\sum_{a,b_{1},b_{2},d,e_{1},e_{2}=0}^{1}\sum_{l_{1}=0}^{b_{1}}%
\sum_{l_{2}=0}^{b_{2}}\sum_{u_{1}=0}^{e_{1}}\sum_{u_{2}=0}^{e_{2}}\left(
-1\right) ^{\alpha \left( 1_{H};l_{1},l_{2},u_{1},u_{2}\right) } \\
&&B(x_{1}x_{2}\otimes
1_{H};G^{a}X_{1}^{b_{1}}X_{2}^{b_{2}},g^{d}x_{1}^{e_{1}}x_{2}^{e_{2}}) \\
&&G^{a}X_{1}^{b_{1}-l_{1}}X_{2}^{b_{2}-l_{2}}\otimes
g^{d}x_{1}^{e_{1}-u_{1}}x_{2}^{e_{2}-u_{2}}\otimes
g^{a+b_{1}+b_{2}+l_{1}+l_{2}+d+e_{1}+e_{2}+u_{1}+u_{2}}x_{1}^{l_{1}+u_{1}}x_{2}^{l_{2}+u_{2}}.
\end{eqnarray*}%
We get%
\begin{eqnarray*}
a+b_{1}+b_{2}+l_{1}+l_{2}+d+e_{1}+e_{2}+u_{1}+u_{2} &\equiv &1 \\
l_{1}+u_{1} &=&1 \\
l_{2}+u_{2} &=&0
\end{eqnarray*}%
i.e.%
\begin{eqnarray*}
a+b_{1}+b_{2}+d+e_{1}+e_{2} &\equiv &0 \\
l_{1}+u_{1} &=&1 \\
l_{2} &=&u_{2}=0
\end{eqnarray*}%
and we get%
\begin{eqnarray*}
&&\sum_{\substack{ a,b_{1},b_{2},d,e_{1},e_{2}=0  \\ %
a+b_{1}+b_{2}+d+e_{1}+e_{2}\equiv 0}}^{1}\sum_{l_{1}=0}^{b_{1}}\sum
_{\substack{ u_{1}=0  \\ l_{1}+u_{1}=1}}^{e_{1}}\left( -1\right) ^{\alpha
\left( 1_{H};l_{1},0,u_{1},0\right) } \\
&&B(x_{1}x_{2}\otimes
1_{H};G^{a}X_{1}^{b_{1}}X_{2}^{b_{2}},g^{d}x_{1}^{e_{1}}x_{2}^{e_{2}})G^{a}X_{1}^{b_{1}-l_{1}}X_{2}^{b_{2}}\otimes g^{d}x_{1}^{e_{1}-u_{1}}x_{2}^{e_{2}}.
\end{eqnarray*}%
Since $\alpha \left( 1_{H};0,0,1,0\right) =e_{2}+\left( a+b_{1}+b_{2}\right)
$ and $\alpha \left( 1_{H};1,0,0,0\right) =b_{2},$ we obtain
\begin{eqnarray*}
&&\sum_{\substack{ a,b_{1},b_{2},d,e_{2}=0  \\ a+b_{1}+b_{2}+d+e_{2}\equiv 1
}}^{1}\left( -1\right) ^{e_{2}+\left( a+b_{1}+b_{2}\right)
}B(x_{1}x_{2}\otimes
1_{H};G^{a}X_{1}^{b_{1}}X_{2}^{b_{2}},g^{d}x_{1}x_{2}^{e_{2}})G^{a}X_{1}^{b_{1}}X_{2}^{b_{2}}\otimes g^{d}x_{2}^{e_{2}}+
\\
&&+\sum_{\substack{ a,b_{2},d,e_{1},e_{2}=0  \\ a+b_{2}+d+e_{1}+e_{2}\equiv
1 }}^{1}\left( -1\right) ^{b_{2}}B(x_{1}x_{2}\otimes
1_{H};G^{a}X_{1}X_{2}^{b_{2}},g^{d}x_{1}^{e_{1}}x_{2}^{e_{2}})G^{a}X_{2}^{b_{2}}\otimes g^{d}x_{1}^{e_{1}}x_{2}^{e_{2}}.
\end{eqnarray*}%
Summing up, we obtain:%
\begin{eqnarray*}
&&\sum_{\substack{ a,b_{1},b_{2},d,e_{1},e_{2}=0  \\ %
a+b_{1}+b_{2}+d+e_{1}+e_{2}\equiv 1}}^{1}\left( -1\right)
^{a+b_{1}+b_{2}}B(gx_{2}\otimes
1_{H};G^{a}X_{1}^{b_{1}}X_{2}^{b_{2}},g^{d}x_{1}^{e_{1}}x_{2}^{e_{2}})G^{a}X_{1}^{b_{1}}X_{2}^{b_{2}}\otimes g^{d}x_{1}^{e_{1}}x_{2}^{e_{2}}+
\\
&&\sum_{\substack{ a,b_{1},b_{2},d,e_{2}=0  \\ a+b_{1}+b_{2}+d+e_{2}\equiv 1
}}^{1}\left( -1\right) ^{e_{2}+\left( a+b_{1}+b_{2}\right)
}B(x_{1}x_{2}\otimes
1_{H};G^{a}X_{1}^{b_{1}}X_{2}^{b_{2}},g^{d}x_{1}x_{2}^{e_{2}})G^{a}X_{1}^{b_{1}}X_{2}^{b_{2}}\otimes g^{d}x_{2}^{e_{2}}+
\\
&&+\sum_{\substack{ a,b_{2},d,e_{1},e_{2}=0  \\ a+b_{2}+d+e_{1}+e_{2}\equiv
1 }}^{1}\left( -1\right) ^{b_{2}}B(x_{1}x_{2}\otimes
1_{H};G^{a}X_{1}X_{2}^{b_{2}},g^{d}x_{1}^{e_{1}}x_{2}^{e_{2}})G^{a}X_{2}^{b_{2}}\otimes g^{d}x_{1}^{e_{1}}x_{2}^{e_{2}}.
\end{eqnarray*}

\subsubsection{$G^{a}\otimes g^{d}$}

\begin{equation*}
\sum_{\substack{ a,d=0  \\ a+d\equiv 1}}^{1}\left[
\begin{array}{c}
\left( -1\right) ^{a}B(gx_{2}\otimes 1_{H};G^{a},g^{d})+\left( -1\right)
^{a}B(x_{1}x_{2}\otimes 1_{H};G^{a},g^{d}x_{1})+ \\
+B(x_{1}x_{2}\otimes 1_{H};G^{a}X_{1},g^{d})%
\end{array}%
\right] G^{a}\otimes g^{d}=0
\end{equation*}%
and we get%
\begin{equation}
B(gx_{2}\otimes 1_{H};1_{A},g)+B(x_{1}x_{2}\otimes
1_{H};1_{A},gx_{1})+B(x_{1}x_{2}\otimes 1_{H};X_{1},g)=0  \label{x1otx2,ten}
\end{equation}%
and%
\begin{equation}
-B(gx_{2}\otimes 1_{H};G,1_{H})-B(x_{1}x_{2}\otimes
1_{H};G,x_{1})+B(x_{1}x_{2}\otimes 1_{H};GX_{1},1_{H})=0
\label{x1otx2, eleven}
\end{equation}

\subsubsection{$G^{a}\otimes g^{d}x_{2}$}

\begin{gather*}
\sum_{\substack{ a,d=0  \\ a+d\equiv 0}}^{1}\left[
\begin{array}{c}
\left( -1\right) ^{a}B(gx_{2}\otimes 1_{H};G^{a},g^{d}x_{2})+\left(
-1\right) ^{a+1}B(x_{1}x_{2}\otimes 1_{H};G^{a},g^{d}x_{1}x_{2})+ \\
+B(x_{1}x_{2}\otimes 1_{H};G^{a}X_{1},g^{d}x_{2})%
\end{array}%
\right] \\
G^{a}\otimes g^{d}x_{2}=0
\end{gather*}%
and we get%
\begin{equation}
\begin{array}{c}
B(gx_{2}\otimes 1_{H};1_{A},x_{2})-B(x_{1}x_{2}\otimes
1_{H};1_{A},x_{1}x_{2})+ \\
+B(x_{1}x_{2}\otimes 1_{H};X_{1},x_{2})%
\end{array}%
=0  \label{x1otx2, twelve}
\end{equation}%
and%
\begin{equation}
\begin{array}{c}
-B(gx_{2}\otimes 1_{H};G,gx_{2})+B(x_{1}x_{2}\otimes 1_{H};G,gx_{1}x_{2})+
\\
+B(x_{1}x_{2}\otimes 1_{H};GX_{1},gx_{2})%
\end{array}%
=0.  \label{x1otx2, thirteen}
\end{equation}

\subsubsection{$G^{a}\otimes g^{d}x_{1}$}

\begin{equation*}
\sum_{\substack{ a,d=0  \\ a+d\equiv 0}}^{1}\left[ \left( -1\right)
^{a}B(gx_{2}\otimes 1_{H};G^{a},g^{d}x_{1})+B(x_{1}x_{2}\otimes
1_{H};G^{a}X_{1},g^{d}x_{1})\right] G^{a}\otimes g^{d}x_{1}=0
\end{equation*}%
and we get%
\begin{equation}
B(gx_{2}\otimes 1_{H};1_{A},x_{1})+B(x_{1}x_{2}\otimes 1_{H};X_{1},x_{1})=0
\label{x1otx2, fourteen}
\end{equation}%
and%
\begin{equation}
-B(gx_{2}\otimes 1_{H};G,gx_{1})+B(x_{1}x_{2}\otimes 1_{H};GX_{1},gx_{1})=0
\label{x1otx2, fifteen}
\end{equation}

\subsubsection{$G^{a}X_{2}\otimes g^{d}$}

\begin{gather*}
\sum_{\substack{ a,d=0  \\ a+d\equiv 0}}^{1}\left[
\begin{array}{c}
\left( -1\right) ^{a+1}B(gx_{2}\otimes 1_{H};G^{a}X_{2},g^{d})+\left(
-1\right) ^{a+1}B(x_{1}x_{2}\otimes 1_{H};G^{a}X_{2},g^{d}x_{1})+ \\
-B(x_{1}x_{2}\otimes 1_{H};G^{a}X_{1}X_{2},g^{d})%
\end{array}%
\right] \\
G^{a}X_{2}\otimes g^{d}=0
\end{gather*}%
and we get%
\begin{equation*}
-B(gx_{2}\otimes 1_{H};X_{2},1_{H})-B(x_{1}x_{2}\otimes
1_{H};X_{2},x_{1})-B(x_{1}x_{2}\otimes 1_{H};X_{1}X_{2},1_{H})=0
\end{equation*}%
and%
\begin{equation*}
B(gx_{2}\otimes 1_{H};GX_{2},g)+B(x_{1}x_{2}\otimes
1_{H};GX_{2},gx_{1})-B(x_{1}x_{2}\otimes 1_{H};GX_{1}X_{2},g)=0.
\end{equation*}%
By using the form of the element $B\left( gx_{2}\otimes 1_{H}\right) $ we get%
\begin{equation}
\left[ 1-B(gx_{2}\otimes 1_{H};1_{A},x_{2})\right] -B(x_{1}x_{2}\otimes
1_{H};X_{2},x_{1})-B(x_{1}x_{2}\otimes 1_{H};X_{1}X_{2},1_{H})=0
\label{x1otx2, sixteen}
\end{equation}%
and%
\begin{equation}
-B(gx_{2}\otimes 1_{H};G,gx_{2})+B(x_{1}x_{2}\otimes
1_{H};GX_{2},gx_{1})-B(x_{1}x_{2}\otimes 1_{H};GX_{1}X_{2},g)=0.
\label{x1otx2, seventeen}
\end{equation}

\subsubsection{$G^{a}X_{1}\otimes g^{d}$}

\begin{gather*}
\sum_{\substack{ a,d=0  \\ a+d\equiv 0}}^{1}\left[ \left( -1\right)
^{a+1}B(gx_{2}\otimes 1_{H};G^{a}X_{1},g^{d})+\left( -1\right)
^{a+1}B(x_{1}x_{2}\otimes 1_{H};G^{a}X_{1},g^{d}x_{1})\right] \\
G^{a}X_{1}\otimes g^{d}=0
\end{gather*}%
and we get%
\begin{equation*}
-B(gx_{2}\otimes 1_{H};X_{1},1_{H})-B(x_{1}x_{2}\otimes 1_{H};X_{1},x_{1})=0
\end{equation*}%
and%
\begin{equation*}
B(gx_{2}\otimes 1_{H};GX_{1},g)+B(x_{1}x_{2}\otimes 1_{H};GX_{1},gx_{1})=0.
\end{equation*}%
By using the form of the element $B\left( gx_{2}\otimes 1_{H}\right) $ we get%
\begin{equation}
-B(gx_{2}\otimes 1_{H};1_{A},x_{1})-B(x_{1}x_{2}\otimes 1_{H};X_{1},x_{1})=0
\label{x1otx2, eighteen}
\end{equation}%
and%
\begin{equation}
-B(gx_{2}\otimes 1_{H};G,gx_{1})+B(x_{1}x_{2}\otimes 1_{H};GX_{1},gx_{1})=0.
\label{x1otx2, nineteen}
\end{equation}

\subsubsection{$G^{a}\otimes g^{d}x_{1}x_{2}$}

\begin{equation*}
\sum_{\substack{ a,d=0  \\ a+d\equiv 1}}^{1}\left[ \left( -1\right)
^{a}B(gx_{2}\otimes 1_{H};G^{a},g^{d}x_{1}x_{2})+B(x_{1}x_{2}\otimes
1_{H};G^{a}X_{1},g^{d}x_{1}x_{2})\right] G^{a}\otimes g^{d}x_{1}x_{2}=0
\end{equation*}%
and we get%
\begin{equation}
B(gx_{2}\otimes 1_{H};1_{A},gx_{1}x_{2})+B(x_{1}x_{2}\otimes
1_{H};X_{1},gx_{1}x_{2})=0  \label{x1otx2, twenty}
\end{equation}%
and%
\begin{equation}
-B(gx_{2}\otimes 1_{H};G,x_{1}x_{2})+B(x_{1}x_{2}\otimes
1_{H};GX_{1},x_{1}x_{2})=0  \label{x1otx2, twentyone}
\end{equation}

\subsubsection{$G^{a}X_{2}\otimes g^{d}x_{2}$}

\begin{equation*}
\sum_{\substack{ a,d=0  \\ a+d\equiv 1}}^{1}\left[
\begin{array}{c}
\left( -1\right) ^{a+1}B(gx_{2}\otimes 1_{H};G^{a}X_{2},g^{d}x_{2})+\left(
-1\right) ^{a}B(x_{1}x_{2}\otimes 1_{H};G^{a}X_{2},g^{d}x_{1}x_{2})+ \\
-B(x_{1}x_{2}\otimes 1_{H};G^{a}X_{1}X_{2},g^{d}x_{2})%
\end{array}%
\right] G^{a}X_{2}\otimes g^{d}x_{2}
\end{equation*}%
and we get%
\begin{equation*}
\begin{array}{c}
-B(gx_{2}\otimes 1_{H};X_{2},gx_{2})+B(x_{1}x_{2}\otimes
1_{H};X_{2},gx_{1}x_{2})+ \\
-B(x_{1}x_{2}\otimes 1_{H};X_{1}X_{2},gx_{2})%
\end{array}%
=0
\end{equation*}%
and%
\begin{equation*}
\begin{array}{c}
B(gx_{2}\otimes 1_{H};GX_{2},x_{2})-B(x_{1}x_{2}\otimes
1_{H};GX_{2},x_{1}x_{2})+ \\
-B(x_{1}x_{2}\otimes 1_{H};GX_{1}X_{2},x_{2})%
\end{array}%
=0.
\end{equation*}%
By using the form of the element $B\left( gx_{2}\otimes 1_{H}\right) $ we get%
\begin{equation}
B(x_{1}x_{2}\otimes 1_{H};X_{2},gx_{1}x_{2})-B(x_{1}x_{2}\otimes
1_{H};X_{1}X_{2},gx_{2})=0  \label{x1otx2, twentytwo}
\end{equation}%
and%
\begin{equation}
-B(x_{1}x_{2}\otimes 1_{H};GX_{2},x_{1}x_{2})-B(x_{1}x_{2}\otimes
1_{H};GX_{1}X_{2},x_{2})=0.  \label{x1otx2, twentythree}
\end{equation}

\subsubsection{$G^{a}X_{1}\otimes g^{d}x_{2}$}

\begin{equation*}
\sum_{\substack{ a,d=0  \\ a+d\equiv 1}}^{1}\left[ \left( -1\right)
^{a+1}B(gx_{2}\otimes 1_{H};G^{a}X_{1},g^{d}x_{2})+\left( -1\right)
^{a}B(x_{1}x_{2}\otimes 1_{H};G^{a}X_{1},g^{d}x_{1}x_{2})\right]
G^{a}X_{1}\otimes g^{d}x_{2}=0
\end{equation*}%
and we get%
\begin{equation*}
-B(gx_{2}\otimes 1_{H};X_{1},gx_{2})+B(x_{1}x_{2}\otimes
1_{H};X_{1},gx_{1}x_{2})=0
\end{equation*}%
and%
\begin{equation*}
B(gx_{2}\otimes 1_{H};GX_{1},x_{2})-B(x_{1}x_{2}\otimes
1_{H};GX_{1},x_{1}x_{2})=0.
\end{equation*}%
By using the form of the element $B\left( gx_{2}\otimes 1_{H}\right) $ we get%
\begin{equation}
B(gx_{2}\otimes 1_{H};1_{A},gx_{1}x_{2})+B(x_{1}x_{2}\otimes
1_{H};X_{1},gx_{1}x_{2})=0  \label{x1otx2, twentyfour}
\end{equation}%
and%
\begin{equation}
B(gx_{2}\otimes 1_{H};G,x_{1}x_{2})-B(x_{1}x_{2}\otimes
1_{H};GX_{1},x_{1}x_{2})=0.  \label{x1otx2, twentyfive}
\end{equation}

\subsubsection{$G^{a}X_{2}\otimes g^{d}x_{1}$}

\begin{equation*}
\sum_{\substack{ a,d=0  \\ a+d\equiv 1}}^{1}\left[ \left( -1\right)
^{a+1}B(gx_{2}\otimes 1_{H};G^{a}X_{2},g^{d}x_{1})-B(x_{1}x_{2}\otimes
1_{H};G^{a}X_{1}X_{2},g^{d}x_{1})\right] G^{a}X_{2}\otimes g^{d}x_{1}=0
\end{equation*}%
and we get%
\begin{equation*}
-B(gx_{2}\otimes 1_{H};X_{2},gx_{1})-B(x_{1}x_{2}\otimes
1_{H};X_{1}X_{2},gx_{1})=0
\end{equation*}%
and%
\begin{equation*}
B(gx_{2}\otimes 1_{H};GX_{2},x_{1})-B(x_{1}x_{2}\otimes
1_{H};GX_{1}X_{2},x_{1})=0.
\end{equation*}%
By using the form of the element $B\left( gx_{2}\otimes 1_{H}\right) $ we get%
\begin{equation}
-B(gx_{2}\otimes 1_{H};1_{A},gx_{1}x_{2})-B(x_{1}x_{2}\otimes
1_{H};X_{1}X_{2},gx_{1})=0  \label{x1otx2, twentysix}
\end{equation}%
and%
\begin{equation}
-B(gx_{2}\otimes 1_{H};G,x_{1}x_{2})-B(x_{1}x_{2}\otimes
1_{H};GX_{1}X_{2},x_{1})=0.  \label{x1otx2, twentyseven}
\end{equation}

\subsubsection{$G^{a}X_{1}\otimes g^{d}x_{1}$}

\begin{equation*}
\sum_{\substack{ a,d=0  \\ a+d\equiv 1}}^{1}\left( -1\right)
^{a+1}B(gx_{2}\otimes 1_{H};G^{a}X_{1},g^{d}x_{1})G^{a}X_{1}\otimes
g^{d}x_{1}=0
\end{equation*}%
and we get%
\begin{equation*}
-B(gx_{2}\otimes 1_{H};X_{1},gx_{1})=0
\end{equation*}%
and%
\begin{equation*}
B(gx_{2}\otimes 1_{H};GX_{1},x_{1})=0.
\end{equation*}%
By the form of the element $B(gx_{2}\otimes 1_{H})$ we know that these
equalities are already known.

\subsubsection{$G^{a}X_{1}X_{2}\otimes g^{d}$}

\begin{equation*}
\sum_{\substack{ a,d=0  \\ a+d\equiv 1}}^{1}\left[ \left( -1\right)
^{a}B(gx_{2}\otimes 1_{H};G^{a}X_{1}X_{2},g^{d})+\left( -1\right)
^{a}B(x_{1}x_{2}\otimes 1_{H};G^{a}X_{1}X_{2},g^{d}x_{1})\right]
G^{a}X_{1}X_{2}\otimes g^{d}=0
\end{equation*}%
and we get%
\begin{equation*}
B(gx_{2}\otimes 1_{H};X_{1}X_{2},g)+B(x_{1}x_{2}\otimes
1_{H};X_{1}X_{2},gx_{1})=0
\end{equation*}%
and%
\begin{equation*}
-B(gx_{2}\otimes 1_{H};GX_{1}X_{2},1_{H})-B(x_{1}x_{2}\otimes
1_{H};GX_{1}X_{2},x_{1})=0.
\end{equation*}%
By using the form of the element $B\left( gx_{2}\otimes 1_{H}\right) $ we get%
\begin{equation*}
B(gx_{2}\otimes 1_{H};1_{A},gx_{1}x_{2})+B(x_{1}x_{2}\otimes
1_{H};X_{1}X_{2},gx_{1})=0
\end{equation*}%
and%
\begin{equation*}
-B(gx_{2}\otimes 1_{H};G,x_{1}x_{2})-B(x_{1}x_{2}\otimes
1_{H};GX_{1}X_{2},x_{1})=0.
\end{equation*}%
These are respectively $\left( \ref{x1otx2, twentysix}\right) $ and $\left( %
\ref{x1otx2, twentyseven}\right) .$

\subsubsection{$G^{a}X_{2}\otimes g^{d}x_{1}x_{2}$}

\begin{equation*}
\sum_{\substack{ a,d=0  \\ a+d\equiv 0}}^{1}\left[ \left( -1\right)
^{a+1}B(gx_{2}\otimes 1_{H};G^{a}X_{2},g^{d}x_{1}x_{2})-B(x_{1}x_{2}\otimes
1_{H};G^{a}X_{1}X_{2},g^{d}x_{1}x_{2})\right] G^{a}X_{2}\otimes
g^{d}x_{1}x_{2}=0
\end{equation*}%
and we get%
\begin{equation*}
-B(gx_{2}\otimes 1_{H};X_{2},x_{1}x_{2})-B(x_{1}x_{2}\otimes
1_{H};X_{1}X_{2},x_{1}x_{2})=0
\end{equation*}%
and%
\begin{equation*}
B(gx_{2}\otimes 1_{H};GX_{2},gx_{1}x_{2})-B(x_{1}x_{2}\otimes
1_{H};GX_{1}X_{2},gx_{1}x_{2})=0.
\end{equation*}%
By using the form of the element $B\left( gx_{2}\otimes 1_{H}\right) $ we get%
\begin{equation*}
B(x_{1}x_{2}\otimes 1_{H};X_{1}X_{2},x_{1}x_{2})=0
\end{equation*}%
and%
\begin{equation*}
B(x_{1}x_{2}\otimes 1_{H};GX_{1}X_{2},gx_{1}x_{2})=0.
\end{equation*}%
which we already got.

\subsubsection{$G^{a}X_{1}\otimes g^{d}x_{1}x_{2}$}

\begin{equation*}
\sum_{\substack{ a,d=0  \\ a+d\equiv 0}}^{1}\left( -1\right)
^{a+1}B(gx_{2}\otimes 1_{H};G^{a}X_{1},g^{d}x_{1}x_{2})G^{a}X_{1}\otimes
g^{d}x_{1}x_{2}=0
\end{equation*}%
and we get%
\begin{equation*}
-B(gx_{2}\otimes 1_{H};X_{1},x_{1}x_{2})=0
\end{equation*}%
and%
\begin{equation*}
B(gx_{2}\otimes 1_{H};GX_{1},gx_{1}x_{2})=0
\end{equation*}%
which are already known because of the form of the element $B\left(
gx_{2}\otimes 1_{H}\right) .$

\subsubsection{$G^{a}X_{1}X_{2}\otimes g^{d}x_{2}$}

\begin{gather*}
\sum_{\substack{ a,d=0  \\ a+d\equiv 0}}^{1}\left[ \left( -1\right)
^{a}B(gx_{2}\otimes 1_{H};G^{a}X_{1}X_{2},g^{d}x_{2})+\left( -1\right)
^{a+1}B(x_{1}x_{2}\otimes 1_{H};G^{a}X_{1}X_{2},g^{d}x_{1}x_{2})\right] \\
G^{a}X_{1}X_{2}\otimes g^{d}x_{2}=0
\end{gather*}%
and we get%
\begin{equation*}
B(gx_{2}\otimes 1_{H};X_{1}X_{2},x_{2})-B(x_{1}x_{2}\otimes
1_{H};X_{1}X_{2},x_{1}x_{2})=0
\end{equation*}%
and%
\begin{equation*}
-B(gx_{2}\otimes 1_{H};GX_{1}X_{2},gx_{2})+B(x_{1}x_{2}\otimes
1_{H};GX_{1}X_{2},gx_{1}x_{2})=0.
\end{equation*}%
By using the form of the element $B\left( gx_{2}\otimes 1_{H}\right) $ we get%
\begin{equation*}
B(x_{1}x_{2}\otimes 1_{H};X_{1}X_{2},x_{1}x_{2})=0
\end{equation*}%
and%
\begin{equation*}
B(x_{1}x_{2}\otimes 1_{H};GX_{1}X_{2},gx_{1}x_{2})=0.
\end{equation*}%
which we already got.

\subsubsection{$G^{a}X_{1}X_{2}\otimes g^{d}x_{1}$}

\begin{equation*}
\sum_{\substack{ a,d=0  \\ a+d\equiv 0}}^{1}\left( -1\right)
^{a}B(gx_{2}\otimes 1_{H};G^{a}X_{1}X_{2},g^{d}x_{1})G^{a}X_{1}X_{2}\otimes
g^{d}x_{1}=0
\end{equation*}%
and we get%
\begin{equation*}
B(gx_{2}\otimes 1_{H};X_{1}X_{2},x_{1})=0
\end{equation*}%
and

\begin{equation*}
-B(gx_{2}\otimes 1_{H};GX_{1}X_{2},gx_{1})=0
\end{equation*}%
which are already known because of the form of the element $B\left(
gx_{2}\otimes 1_{H}\right) .$

\subsubsection{$G^{a}X_{1}X_{2}\otimes g^{d}x_{1}x_{2}$}

\begin{equation*}
\sum_{\substack{ a,d=0  \\ a+d\equiv 1}}^{1}\left( -1\right)
^{a}B(gx_{2}\otimes
1_{H};G^{a}X_{1}X_{2},g^{d}x_{1}x_{2})G^{a}X_{1}X_{2}\otimes
g^{d}x_{1}x_{2}=0
\end{equation*}%
and we get%
\begin{equation*}
B(gx_{2}\otimes 1_{H};X_{1}X_{2},gx_{1}x_{2})=0
\end{equation*}%
and%
\begin{equation*}
-B(gx_{2}\otimes 1_{H};GX_{1}X_{2},x_{1}x_{2})=0
\end{equation*}%
which are already known because of the form of the element $B\left(
gx_{2}\otimes 1_{H}\right) .$

\subsection{Case $gx_{2}$}

\begin{eqnarray*}
&&\sum_{a,b_{1},b_{2},d,e_{1},e_{2}=0}^{1}\sum_{l_{1}=0}^{b_{1}}%
\sum_{l_{2}=0}^{b_{2}}\sum_{u_{1}=0}^{e_{1}}\sum_{u_{2}=0}^{e_{2}}\left(
-1\right) ^{\alpha \left( x_{1}x_{2};l_{1},l_{2},u_{1},u_{2}\right) } \\
&&B(1_{H}\otimes
1_{H};G^{a}X_{1}^{b_{1}}X_{2}^{b_{2}},g^{d}x_{1}^{e_{1}}x_{2}^{e_{2}}) \\
&&G^{a}X_{1}^{b_{1}-l_{1}}X_{2}^{b_{2}-l_{2}}\otimes
g^{d}x_{1}^{e_{1}-u_{1}}x_{2}^{e_{2}-u_{2}}\otimes
g^{a+b_{1}+b_{2}+l_{1}+l_{2}+d+e_{1}+e_{2}+u_{1}+u_{2}}x_{1}^{l_{1}+u_{1}+1}x_{2}^{l_{2}+u_{2}+1}
\\
+
&&\sum_{a,b_{1},b_{2},d,e_{1},e_{2}=0}^{1}\sum_{l_{1}=0}^{b_{1}}%
\sum_{l_{2}=0}^{b_{2}}\sum_{u_{1}=0}^{e_{1}}\sum_{u_{2}=0}^{e_{2}}\left(
-1\right) ^{\alpha \left( x_{1};l_{1},l_{2},u_{1},u_{2}\right) } \\
&&B(gx_{2}\otimes
1_{H};G^{a}X_{1}^{b_{1}}X_{2}^{b_{2}},g^{d}x_{1}^{e_{1}}x_{2}^{e_{2}}) \\
&&G^{a}X_{1}^{b_{1}-l_{1}}X_{2}^{b_{2}-l_{2}}\otimes
g^{d}x_{1}^{e_{1}-u_{1}}x_{2}^{e_{2}-u_{2}}\otimes
g^{a+b_{1}+b_{2}+l_{1}+l_{2}+d+e_{1}+e_{2}+u_{1}+u_{2}}x_{1}^{l_{1}+u_{1}+1}x_{2}^{l_{2}+u_{2}}
\\
&&-\sum_{a,b_{1},b_{2},d,e_{1},e_{2}=0}^{1}\sum_{l_{1}=0}^{b_{1}}%
\sum_{l_{2}=0}^{b_{2}}\sum_{u_{1}=0}^{e_{1}}\sum_{u_{2}=0}^{e_{2}}\left(
-1\right) ^{\alpha \left( x_{2};l_{1},l_{2},u_{1},u_{2}\right) } \\
&&B(gx_{1}\otimes
1_{H};G^{a}X_{1}^{b_{1}}X_{2}^{b_{2}},g^{d}x_{1}^{e_{1}}x_{2}^{e_{2}}) \\
&&G^{a}X_{1}^{b_{1}-l_{1}}X_{2}^{b_{2}-l_{2}}\otimes
g^{d}x_{1}^{e_{1}-u_{1}}x_{2}^{e_{2}-u_{2}}\otimes
g^{a+b_{1}+b_{2}+l_{1}+l_{2}+d+e_{1}+e_{2}+u_{1}+u_{2}}x_{1}^{l_{1}+u_{1}}x_{2}^{l_{2}+u_{2}+1}
\\
&&+\sum_{a,b_{1},b_{2},d,e_{1},e_{2}=0}^{1}\sum_{l_{1}=0}^{b_{1}}%
\sum_{l_{2}=0}^{b_{2}}\sum_{u_{1}=0}^{e_{1}}\sum_{u_{2}=0}^{e_{2}}\left(
-1\right) ^{\alpha \left( 1_{H};l_{1},l_{2},u_{1},u_{2}\right)
}B(x_{1}x_{2}\otimes
1_{H};G^{a}X_{1}^{b_{1}}X_{2}^{b_{2}},g^{d}x_{1}^{e_{1}}x_{2}^{e_{2}}) \\
&&G^{a}X_{1}^{b_{1}-l_{1}}X_{2}^{b_{2}-l_{2}}\otimes
g^{d}x_{1}^{e_{1}-u_{1}}x_{2}^{e_{2}-u_{2}}\otimes
g^{a+b_{1}+b_{2}+l_{1}+l_{2}+d+e_{1}+e_{2}+u_{1}+u_{2}}x_{1}^{l_{1}+u_{1}}x_{2}^{l_{2}+u_{2}}
\\
&=&B^{A}(x_{1}x_{2}\otimes 1_{H})\otimes B^{H}(x_{1}x_{2}\otimes
1_{H})\otimes 1_{H}
\end{eqnarray*}%
\begin{eqnarray*}
&&-\sum_{a,b_{1},b_{2},d,e_{1},e_{2}=0}^{1}\sum_{l_{1}=0}^{b_{1}}%
\sum_{l_{2}=0}^{b_{2}}\sum_{u_{1}=0}^{e_{1}}\sum_{u_{2}=0}^{e_{2}}\left(
-1\right) ^{\alpha \left( x_{2};l_{1},l_{2},u_{1},u_{2}\right) } \\
&&B(gx_{1}\otimes
1_{H};G^{a}X_{1}^{b_{1}}X_{2}^{b_{2}},g^{d}x_{1}^{e_{1}}x_{2}^{e_{2}}) \\
&&G^{a}X_{1}^{b_{1}-l_{1}}X_{2}^{b_{2}-l_{2}}\otimes
g^{d}x_{1}^{e_{1}-u_{1}}x_{2}^{e_{2}-u_{2}}\otimes
g^{a+b_{1}+b_{2}+l_{1}+l_{2}+d+e_{1}+e_{2}+u_{1}+u_{2}}x_{1}^{l_{1}+u_{1}}x_{2}^{l_{2}+u_{2}+1}
\end{eqnarray*}%
We have to consider only the third and the fourth summand.

We examine the third summand.%
\begin{eqnarray*}
&&-\sum_{a,b_{1},b_{2},d,e_{1},e_{2}=0}^{1}\sum_{l_{1}=0}^{b_{1}}%
\sum_{l_{2}=0}^{b_{2}}\sum_{u_{1}=0}^{e_{1}}\sum_{u_{2}=0}^{e_{2}}\left(
-1\right) ^{\alpha \left( x_{2};l_{1},l_{2},u_{1},u_{2}\right) } \\
&&B(gx_{1}\otimes
1_{H};G^{a}X_{1}^{b_{1}}X_{2}^{b_{2}},g^{d}x_{1}^{e_{1}}x_{2}^{e_{2}}) \\
&&G^{a}X_{1}^{b_{1}-l_{1}}X_{2}^{b_{2}-l_{2}}\otimes
g^{d}x_{1}^{e_{1}-u_{1}}x_{2}^{e_{2}-u_{2}}\otimes
g^{a+b_{1}+b_{2}+l_{1}+l_{2}+d+e_{1}+e_{2}+u_{1}+u_{2}}x_{1}^{l_{1}+u_{1}}x_{2}^{l_{2}+u_{2}+1}
\end{eqnarray*}%
We get%
\begin{eqnarray*}
a+b_{1}+b_{2}+l_{1}+l_{2}+d+e_{1}+e_{2}+u_{1}+u_{2} &\equiv &1 \\
l_{1}+u_{1} &=&0 \\
l_{2}+u_{2}+1 &=&1
\end{eqnarray*}%
i.e.%
\begin{eqnarray*}
a+b_{1}+b_{2}+d+e_{1}+e_{2} &\equiv &1 \\
l_{1} &=&u_{1}=0 \\
l_{2} &=&u_{2}=0.
\end{eqnarray*}%
Since $\alpha \left( x_{2};0,0,0,0\right) =a+b_{1}+b_{2},$ we obtain%
\begin{equation*}
\sum_{\substack{ a,b_{1},b_{2},d,e_{1},e_{2}=0  \\ %
a+b_{1}+b_{2}+d+e_{1}+e_{2}\equiv 1}}^{1}\left( -1\right)
^{a+b_{1}+b_{2}+1}B(gx_{1}\otimes
1_{H};G^{a}X_{1}^{b_{1}}X_{2}^{b_{2}},g^{d}x_{1}^{e_{1}}x_{2}^{e_{2}})G^{a}X_{1}^{b_{1}}X_{2}^{b_{2}}\otimes g^{d}x_{1}^{e_{1}}x_{2}^{e_{2}}.
\end{equation*}%
We examine the fourth summand.%
\begin{eqnarray*}
&&\sum_{a,b_{1},b_{2},d,e_{1},e_{2}=0}^{1}\sum_{l_{1}=0}^{b_{1}}%
\sum_{l_{2}=0}^{b_{2}}\sum_{u_{1}=0}^{e_{1}}\sum_{u_{2}=0}^{e_{2}}\left(
-1\right) ^{\alpha \left( 1_{H};l_{1},l_{2},u_{1},u_{2}\right) } \\
&&B(x_{1}x_{2}\otimes
1_{H};G^{a}X_{1}^{b_{1}}X_{2}^{b_{2}},g^{d}x_{1}^{e_{1}}x_{2}^{e_{2}}) \\
&&G^{a}X_{1}^{b_{1}-l_{1}}X_{2}^{b_{2}-l_{2}}\otimes
g^{d}x_{1}^{e_{1}-u_{1}}x_{2}^{e_{2}-u_{2}}\otimes
g^{a+b_{1}+b_{2}+l_{1}+l_{2}+d+e_{1}+e_{2}+u_{1}+u_{2}}x_{1}^{l_{1}+u_{1}}x_{2}^{l_{2}+u_{2}}
\end{eqnarray*}%
We get%
\begin{eqnarray*}
a+b_{1}+b_{2}+l_{1}+l_{2}+d+e_{1}+e_{2}+u_{1}+u_{2} &\equiv &1 \\
l_{1}+u_{1} &=&0 \\
l_{2}+u_{2} &=&1
\end{eqnarray*}%
i.e.%
\begin{eqnarray*}
a+b_{1}+b_{2}+d+e_{1}+e_{2} &\equiv &0 \\
l_{1} &=&u_{1}=0 \\
l_{2}+u_{2} &=&1.
\end{eqnarray*}%
We obtain%
\begin{eqnarray*}
&&\sum_{\substack{ a,b_{1},b_{2},d,e_{1},e_{2}=0  \\ %
a+b_{1}+b_{2}+d+e_{1}+e_{2}\equiv 0}}^{1}\sum_{l_{2}=0}^{b_{2}}\sum
_{\substack{ u_{2}=0  \\ l_{2}+u_{2}=1}}^{e_{2}}\left( -1\right) ^{\alpha
\left( 1_{H};0,l_{2},0,u_{2}\right) }B(x_{1}x_{2}\otimes
1_{H};G^{a}X_{1}^{b_{1}}X_{2}^{b_{2}},g^{d}x_{1}^{e_{1}}x_{2}^{e_{2}}) \\
&&G^{a}X_{1}^{b_{1}}X_{2}^{b_{2}-l_{2}}\otimes
g^{d}x_{1}^{e_{1}}x_{2}^{e_{2}-u_{2}}.
\end{eqnarray*}

Since $\alpha \left( 1_{H};0,0,0,1\right) =a+b_{1}+b_{2}$ and $\alpha \left(
1_{H};0,1,0,0\right) =0,$ we get%
\begin{eqnarray*}
&&\sum_{\substack{ a,b_{1},b_{2},d,e_{1}=0  \\ a+b_{1}+b_{2}+d+e_{1}\equiv 1
}}^{1}\left( -1\right) ^{a+b_{1}+b_{2}}B(x_{1}x_{2}\otimes
1_{H};G^{a}X_{1}^{b_{1}}X_{2}^{b_{2}},g^{d}x_{1}^{e_{1}}x_{2})G^{a}X_{1}^{b_{1}}X_{2}^{b_{2}}\otimes g^{d}x_{1}^{e_{1}}+
\\
&&+\sum_{\substack{ a,b_{1},d,e_{1},e_{2}=0  \\ a+b_{1}+d+e_{1}+e_{2}\equiv
1 }}^{1}B(x_{1}x_{2}\otimes
1_{H};G^{a}X_{1}^{b_{1}}X_{2},g^{d}x_{1}^{e_{1}}x_{2}^{e_{2}})G^{a}X_{1}^{b_{1}}\otimes g^{d}x_{1}^{e_{1}}x_{2}^{e_{2}}.
\end{eqnarray*}%
Summing up we obtain%
\begin{eqnarray*}
&&\sum_{\substack{ a,b_{1},b_{2},d,e_{1},e_{2}=0  \\ %
a+b_{1}+b_{2}+d+e_{1}+e_{2}\equiv 1}}^{1}\left( -1\right)
^{a+b_{1}+b_{2}+1}B(gx_{1}\otimes
1_{H};G^{a}X_{1}^{b_{1}}X_{2}^{b_{2}},g^{d}x_{1}^{e_{1}}x_{2}^{e_{2}})G^{a}X_{1}^{b_{1}}X_{2}^{b_{2}}\otimes g^{d}x_{1}^{e_{1}}x_{2}^{e_{2}}
\\
&&\sum_{\substack{ a,b_{1},b_{2},d,e_{1}=0  \\ a+b_{1}+b_{2}+d+e_{1}\equiv 1
}}^{1}\left( -1\right) ^{a+b_{1}+b_{2}}B(x_{1}x_{2}\otimes
1_{H};G^{a}X_{1}^{b_{1}}X_{2}^{b_{2}},g^{d}x_{1}^{e_{1}}x_{2})G^{a}X_{1}^{b_{1}}X_{2}^{b_{2}}\otimes g^{d}x_{1}^{e_{1}}+
\\
&&+\sum_{\substack{ a,b_{1},d,e_{1},e_{2}=0  \\ a+b_{1}+d+e_{1}+e_{2}\equiv
1 }}^{1}B(x_{1}x_{2}\otimes
1_{H};G^{a}X_{1}^{b_{1}}X_{2},g^{d}x_{1}^{e_{1}}x_{2}^{e_{2}})G^{a}X_{1}^{b_{1}}\otimes g^{d}x_{1}^{e_{1}}x_{2}^{e_{2}}=0
\end{eqnarray*}

\subsubsection{$G^{a}\otimes g^{d}$}

\begin{equation*}
\sum_{\substack{ a,d=0  \\ a+d\equiv 1}}^{1}\left[
\begin{array}{c}
\left( -1\right) ^{a+1}B(gx_{1}\otimes 1_{H};G^{a},g^{d})+\left( -1\right)
^{a}B(x_{1}x_{2}\otimes 1_{H};G^{a},g^{d}x_{2})+ \\
B(x_{1}x_{2}\otimes 1_{H};G^{a}X_{2},g^{d})%
\end{array}%
\right] G^{a}\otimes g^{d}=0
\end{equation*}%
and we get%
\begin{equation}
-B(gx_{1}\otimes 1_{H};1_{A},g)+B(x_{1}x_{2}\otimes
1_{H};1_{A},gx_{2})+B(x_{1}x_{2}\otimes 1_{H};X_{2},g)=0
\label{x1otx2, twentyeight}
\end{equation}%
and%
\begin{equation}
B(gx_{1}\otimes 1_{H};G,1_{H})-B(x_{1}x_{2}\otimes
1_{H};G,x_{2})+B(x_{1}x_{2}\otimes 1_{H};GX_{2},1_{H})=0.
\label{x1otx2, twentynine}
\end{equation}

\subsubsection{$G^{a}\otimes g^{d}x_{2}$}

\begin{equation*}
\sum_{\substack{ a,d=0  \\ a+d\equiv 0}}^{1}\left[ \left( -1\right)
^{a+1}B(gx_{1}\otimes 1_{H};G^{a},g^{d}x_{2})+B(x_{1}x_{2}\otimes
1_{H};G^{a}X_{2},g^{d}x_{2})\right] G^{a}\otimes g^{d}x_{2}=0
\end{equation*}%
and we get%
\begin{equation}
-B(gx_{1}\otimes 1_{H};1_{A},x_{2})+B(x_{1}x_{2}\otimes 1_{H};X_{2},x_{2})=0
\label{x1otx2, thirty}
\end{equation}%
and%
\begin{equation}
B(gx_{1}\otimes 1_{H};G,gx_{2})+B(x_{1}x_{2}\otimes 1_{H};GX_{2},gx_{2})=0.
\label{x1otx2, thirtyone}
\end{equation}

\subsubsection{$G^{a}\otimes g^{d}x_{1}$}

\begin{equation*}
\sum_{\substack{ a,d=0  \\ a+d\equiv 0}}^{1}\left[
\begin{array}{c}
\left( -1\right) ^{a+1}B(gx_{1}\otimes 1_{H};G^{a},g^{d}x_{1})+\left(
-1\right) ^{a}B(x_{1}x_{2}\otimes 1_{H};G^{a},g^{d}x_{1}x_{2})+ \\
+B(x_{1}x_{2}\otimes 1_{H};G^{a}X_{2},g^{d}x_{1})%
\end{array}%
\right] G^{a}\otimes g^{d}x_{1}
\end{equation*}%
and we get%
\begin{equation}
-B(gx_{1}\otimes 1_{H};1_{A},x_{1})+B(x_{1}x_{2}\otimes
1_{H};1_{A},x_{1}x_{2})+B(x_{1}x_{2}\otimes 1_{H};X_{2},x_{1})=0
\label{x1otx2, thirtytwo}
\end{equation}%
and%
\begin{equation}
B(gx_{1}\otimes 1_{H};G,gx_{1})-B(x_{1}x_{2}\otimes
1_{H};G,gx_{1}x_{2})+B(x_{1}x_{2}\otimes 1_{H};GX_{2},gx_{1})=0.
\label{x1otx2, thirtythree}
\end{equation}

\subsubsection{$G^{a}X_{2}\otimes g^{d}$}

\begin{equation*}
\sum_{\substack{ a,d=0  \\ a+d\equiv 0}}^{1}\left[ \left( -1\right)
^{a}B(gx_{1}\otimes 1_{H};G^{a}X_{2},g^{d})+\left( -1\right)
^{a+1}B(x_{1}x_{2}\otimes 1_{H};G^{a}X_{2},g^{d}x_{2})\right]
G^{a}X_{2}\otimes g^{d}=0
\end{equation*}%
and we get%
\begin{equation*}
B(gx_{1}\otimes 1_{H};X_{2},1_{H})-B(x_{1}x_{2}\otimes 1_{H};X_{2},x_{2})=0
\end{equation*}%
and%
\begin{equation*}
-B(gx_{1}\otimes 1_{H};GX_{2},g)+B(x_{1}x_{2}\otimes 1_{H};GX_{2},gx_{2})=0.
\end{equation*}%
By using the form of the element $B\left( gx_{1}\otimes 1_{H}\right) $ we get%
\begin{equation*}
B(gx_{1}\otimes 1_{H};1_{A},x_{2})-B(x_{1}x_{2}\otimes 1_{H};X_{2},x_{2})=0
\end{equation*}%
and%
\begin{equation*}
B(gx_{2}\otimes 1_{H};G,gx_{2})+B(x_{1}x_{2}\otimes 1_{H};GX_{2},gx_{2})=0
\end{equation*}%
which are respectively $\left( \ref{x1otx2, thirty}\right) $ and $\left( \ref%
{x1otx2, thirtyone}\right) .$

\subsubsection{$G^{a}X_{1}\otimes g^{d}$}

\begin{equation*}
\sum_{\substack{ a,d=0  \\ a+d\equiv 0}}^{1}\left[
\begin{array}{c}
\left( -1\right) ^{a}B(gx_{1}\otimes 1_{H};G^{a}X_{1},g^{d})+\left(
-1\right) ^{a+1}B(x_{1}x_{2}\otimes 1_{H};G^{a}X_{1},g^{d}x_{2})+ \\
+B(x_{1}x_{2}\otimes 1_{H};G^{a}X_{1}X_{2},g^{d})%
\end{array}%
\right] G^{a}X_{1}\otimes g^{d}=0
\end{equation*}%
and we get%
\begin{equation*}
\begin{array}{c}
B(gx_{1}\otimes 1_{H};X_{1},1_{H})-B(x_{1}x_{2}\otimes 1_{H};X_{1},x_{2})+
\\
+B(x_{1}x_{2}\otimes 1_{H};X_{1}X_{2},1_{H})%
\end{array}%
=0
\end{equation*}%
and%
\begin{equation*}
\begin{array}{c}
-B(gx_{1}\otimes 1_{H};GX_{1},g)+B(x_{1}x_{2}\otimes 1_{H};GX_{1},gx_{2})+
\\
+B(x_{1}x_{2}\otimes 1_{H};GX_{1}X_{2},g)%
\end{array}%
=0.
\end{equation*}%
By using the form of the element $B\left( gx_{1}\otimes 1_{H}\right) $ we get%
\begin{equation}
\begin{array}{c}
\left[ -1+B(gx_{1}\otimes 1_{H};1_{H},x_{1})\right] -B(x_{1}x_{2}\otimes
1_{H};X_{1},x_{2})+ \\
+B(x_{1}x_{2}\otimes 1_{H};X_{1}X_{2},1_{H})%
\end{array}%
=0  \label{x1otx2, thirtyfour}
\end{equation}%
and%
\begin{equation}
B(gx_{2}\otimes 1_{H};G,gx_{1})+B(x_{1}x_{2}\otimes
1_{H};GX_{1},gx_{2})+B(x_{1}x_{2}\otimes 1_{H};GX_{1}X_{2},g)=0.
\label{x1otx2, thirtyfive}
\end{equation}

\subsubsection{$G^{a}\otimes g^{d}x_{1}x_{2}$}

\begin{equation*}
\sum_{\substack{ a,d=0  \\ a+d\equiv 1}}^{1}\left[ \left( -1\right)
^{a+1}B(gx_{1}\otimes 1_{H};G^{a},g^{d}x_{1}x_{2})+B(x_{1}x_{2}\otimes
1_{H};G^{a}X_{2},g^{d}x_{1}x_{2})\right] G^{a}\otimes g^{d}x_{1}x_{2}=0
\end{equation*}%
and we get%
\begin{equation}
-B(gx_{1}\otimes 1_{H};1_{A},gx_{1}x_{2})+B(x_{1}x_{2}\otimes
1_{H};X_{2},gx_{1}x_{2})=0  \label{x1otx2, thirtysix}
\end{equation}%
and%
\begin{equation}
B(gx_{1}\otimes 1_{H};G,x_{1}x_{2})+B(x_{1}x_{2}\otimes
1_{H};GX_{2},x_{1}x_{2})=0  \label{x1otx2, thirtyseven}
\end{equation}

\subsubsection{$G^{a}X_{2}\otimes g^{d}x_{2}$}

\begin{equation*}
\sum_{\substack{ a,d=0  \\ a+d\equiv 1}}^{1}\left( -1\right)
^{a}B(gx_{1}\otimes 1_{H};G^{a}X_{2},g^{d}x_{2})G^{a}X_{2}\otimes
g^{d}x_{2}=0
\end{equation*}%
and we get%
\begin{equation*}
B(gx_{1}\otimes 1_{H};X_{2},gx_{2})=0
\end{equation*}%
and%
\begin{equation*}
-B(gx_{1}\otimes 1_{H};GX_{2},x_{2})=0
\end{equation*}%
which are already known because of the form of the element $B\left(
gx_{1}\otimes 1_{H}\right) .$

\subsubsection{$G^{a}X_{1}\otimes g^{d}x_{2}$}

\begin{equation*}
\sum_{\substack{ a,d=0  \\ a+d\equiv 1}}^{1}\left[ \left( -1\right)
^{a}B(gx_{1}\otimes 1_{H};G^{a}X_{1},g^{d}x_{2})+B(x_{1}x_{2}\otimes
1_{H};G^{a}X_{1}X_{2},g^{d}x_{2})\right] G^{a}X_{1}\otimes g^{d}x_{2}=0
\end{equation*}%
and we get%
\begin{equation*}
B(gx_{1}\otimes 1_{H};X_{1},gx_{2})+B(x_{1}x_{2}\otimes
1_{H};X_{1}X_{2},gx_{2})=0
\end{equation*}%
and%
\begin{equation*}
-B(gx_{1}\otimes 1_{H};GX_{1},x_{2})+B(x_{1}x_{2}\otimes
1_{H};GX_{1}X_{2},x_{2})=0.
\end{equation*}%
By using the form of the element $B\left( gx_{1}\otimes 1_{H}\right) $ we get%
\begin{equation}
-B(gx_{1}\otimes 1_{H};1_{A},gx_{1}x_{2})+B(x_{1}x_{2}\otimes
1_{H};X_{1}X_{2},gx_{2})=0  \label{x1otx2, thirtyeight}
\end{equation}%
and%
\begin{equation}
-B(gx_{2}\otimes 1_{H};G,x_{1}x_{2})+B(x_{1}x_{2}\otimes
1_{H};GX_{1}X_{2},x_{2})=0  \label{x1otx2, thirtynine}
\end{equation}

\subsubsection{$G^{a}X_{2}\otimes g^{d}x_{1}$}

\begin{equation*}
\sum_{\substack{ a,d=0  \\ a+d\equiv 1}}^{1}\left[ \left( -1\right)
^{a}B(gx_{1}\otimes 1_{H};G^{a}X_{2},g^{d}x_{1})+\left( -1\right)
^{a+1}B(x_{1}x_{2}\otimes 1_{H};G^{a}X_{2},g^{d}x_{1}x_{2})\right]
G^{a}X_{2}\otimes g^{d}x_{1}=0
\end{equation*}%
and we get%
\begin{equation*}
B(gx_{1}\otimes 1_{H};X_{2},gx_{1})-B(x_{1}x_{2}\otimes
1_{H};X_{2},gx_{1}x_{2})=0
\end{equation*}%
and%
\begin{equation*}
-B(gx_{1}\otimes 1_{H};GX_{2},x_{1})+B(x_{1}x_{2}\otimes
1_{H};GX_{2},x_{1}x_{2})=0.
\end{equation*}%
By using the form of the element $B\left( gx_{1}\otimes 1_{H}\right) $ we get%
\begin{equation*}
B(gx_{1}\otimes 1_{H};1_{A},gx_{1}x_{2})-B(x_{1}x_{2}\otimes
1_{H};X_{2},gx_{1}x_{2})=0
\end{equation*}%
and%
\begin{equation*}
B(gx_{2}\otimes 1_{H};G,x_{1}x_{2})+B(x_{1}x_{2}\otimes
1_{H};GX_{2},x_{1}x_{2})=0.
\end{equation*}%
which are respectively $\left( \ref{x1otx2, thirtysix}\right) $ and $\left( %
\ref{x1otx2, thirtyseven}\right) .$

\subsubsection{$G^{a}X_{1}\otimes g^{d}x_{1}$}

\begin{equation*}
\sum_{\substack{ a,d=0  \\ a+d\equiv 1}}^{1}\left[
\begin{array}{c}
\left( -1\right) ^{a}B(gx_{1}\otimes 1_{H};G^{a}X_{1},g^{d}x_{1})+\left(
-1\right) ^{a+1}B(x_{1}x_{2}\otimes 1_{H};G^{a}X_{1},g^{d}x_{1}x_{2})+ \\
B(x_{1}x_{2}\otimes 1_{H};G^{a}X_{1}X_{2},g^{d}x_{1})%
\end{array}%
\right] G^{a}X_{1}\otimes g^{d}x_{1}=0
\end{equation*}%
and we get%
\begin{equation*}
\begin{array}{c}
B(gx_{1}\otimes 1_{H};X_{1},gx_{1})-B(x_{1}x_{2}\otimes
1_{H};X_{1},gx_{1}x_{2})+ \\
B(x_{1}x_{2}\otimes 1_{H};X_{1}X_{2},gx_{1})%
\end{array}%
=0
\end{equation*}%
and%
\begin{equation*}
\begin{array}{c}
-B(gx_{1}\otimes 1_{H};GX_{1},x_{1})+B(x_{1}x_{2}\otimes
1_{H};GX_{1},x_{1}x_{2})+ \\
B(x_{1}x_{2}\otimes 1_{H};GX_{1}X_{2},x_{1})%
\end{array}%
=0.
\end{equation*}%
By using the form of the element $B\left( gx_{1}\otimes 1_{H}\right) $ we get%
\begin{equation}
-B(x_{1}x_{2}\otimes 1_{H};X_{1},gx_{1}x_{2})+B(x_{1}x_{2}\otimes
1_{H};X_{1}X_{2},gx_{1})=0  \label{x1otx2, forty}
\end{equation}%
and%
\begin{equation}
B(x_{1}x_{2}\otimes 1_{H};GX_{1},x_{1}x_{2})+B(x_{1}x_{2}\otimes
1_{H};GX_{1}X_{2},x_{1})=0  \label{x1otx2, fortyone}
\end{equation}

\subsubsection{$G^{a}X_{1}X_{2}\otimes g^{d}$}

\begin{equation*}
\sum_{\substack{ a,d=0  \\ a+d\equiv 1}}^{1}\left[ \left( -1\right)
^{a+1}B(gx_{1}\otimes 1_{H};G^{a}X_{1}X_{2},g^{d})+\left( -1\right)
^{a}B(x_{1}x_{2}\otimes 1_{H};G^{a}X_{1}X_{2},g^{d}x_{2})\right]
G^{a}X_{1}X_{2}\otimes g^{d}=0
\end{equation*}%
and we get%
\begin{equation*}
-B(gx_{1}\otimes 1_{H};X_{1}X_{2},g)+B(x_{1}x_{2}\otimes
1_{H};X_{1}X_{2},gx_{2})=0
\end{equation*}%
and%
\begin{equation*}
B(gx_{1}\otimes 1_{H};GX_{1}X_{2},1_{H})-B(x_{1}x_{2}\otimes
1_{H};GX_{1}X_{2},x_{2})=0.
\end{equation*}%
By using the form of the element $B\left( gx_{1}\otimes 1_{H}\right) $ we get%
\begin{equation*}
-B(gx_{1}\otimes 1_{H};1_{A},gx_{1}x_{2})+B(x_{1}x_{2}\otimes
1_{H};X_{1}X_{2},gx_{2})=0
\end{equation*}%
and%
\begin{equation*}
B(gx_{1}\otimes 1_{H};G,x_{1}x_{2})-B(x_{1}x_{2}\otimes
1_{H};GX_{1}X_{2},x_{2})=0.
\end{equation*}%
These are respectively $\left( \ref{x1otx2, thirtyeight}\right) $ and $%
\left( \ref{x1otx2, thirtynine}\right) .$

\subsubsection{$G^{a}X_{2}\otimes g^{d}x_{1}x_{2}$}

\begin{equation*}
\sum_{\substack{ a,d=0  \\ a+d\equiv 0}}^{1}\left( -1\right)
^{a}B(gx_{1}\otimes 1_{H};G^{a}X_{2},g^{d}x_{1}x_{2})G^{a}X_{2}\otimes
g^{d}x_{1}x_{2}=0
\end{equation*}%
and we get%
\begin{equation*}
B(gx_{1}\otimes 1_{H};X_{2},x_{1}x_{2})=0
\end{equation*}%
and%
\begin{equation*}
-B(gx_{1}\otimes 1_{H};GX_{2},gx_{1}x_{2})=0.
\end{equation*}%
which are already known because of the form of the element $B\left(
gx_{1}\otimes 1_{H}\right) .$

\subsubsection{$G^{a}X_{1}\otimes g^{d}x_{1}x_{2}$}

\begin{equation*}
\sum_{\substack{ a,d=0  \\ a+d\equiv 0}}^{1}\left[ \left( -1\right)
^{a}B(gx_{1}\otimes 1_{H};G^{a}X_{1},g^{d}x_{1}x_{2})+B(x_{1}x_{2}\otimes
1_{H};G^{a}X_{1}X_{2},g^{d}x_{1}x_{2})\right] G^{a}X_{1}\otimes
g^{d}x_{1}x_{2}=0
\end{equation*}%
and we get%
\begin{equation*}
B(gx_{1}\otimes 1_{H};X_{1},x_{1}x_{2})+B(x_{1}x_{2}\otimes
1_{H};X_{1}X_{2},x_{1}x_{2})=0
\end{equation*}%
and%
\begin{equation*}
-B(gx_{1}\otimes 1_{H};GX_{1},gx_{1}x_{2})+B(x_{1}x_{2}\otimes
1_{H};GX_{1}X_{2},gx_{1}x_{2})=0.
\end{equation*}%
By using the form of the element $B\left( gx_{1}\otimes 1_{H}\right) $ we get%
\begin{equation*}
B(x_{1}x_{2}\otimes 1_{H};X_{1}X_{2},x_{1}x_{2})=0
\end{equation*}%
and%
\begin{equation*}
B(x_{1}x_{2}\otimes 1_{H};GX_{1}X_{2},gx_{1}x_{2})=0.
\end{equation*}%
which we already got.

\subsubsection{$G^{a}X_{1}X_{2}\otimes g^{d}x_{2}$}

\begin{equation*}
\sum_{\substack{ a,d=0  \\ a+d\equiv 0}}^{1}\left( -1\right)
^{a+1}B(gx_{1}\otimes
1_{H};G^{a}X_{1}X_{2},g^{d}x_{2})G^{a}X_{1}X_{2}\otimes g^{d}x_{2}=0
\end{equation*}%
and we get%
\begin{equation*}
-B(gx_{1}\otimes 1_{H};X_{1}X_{2},x_{2})=0
\end{equation*}%
and%
\begin{equation*}
B(gx_{1}\otimes 1_{H};GX_{1}X_{2},gx_{2})=0
\end{equation*}%
which are already known because of the form of the element $B\left(
gx_{1}\otimes 1_{H}\right) .$

\subsubsection{$G^{a}X_{1}X_{2}\otimes g^{d}x_{1}$}

\begin{gather*}
\sum_{\substack{ a,d=0  \\ a+d\equiv 0}}^{1}\left[ \left( -1\right)
^{a+1}B(gx_{1}\otimes 1_{H};G^{a}X_{1}X_{2},g^{d}x_{1})+\left( -1\right)
^{a}B(x_{1}x_{2}\otimes 1_{H};G^{a}X_{1}X_{2},g^{d}x_{1}x_{2})\right] \\
G^{a}X_{1}X_{2}\otimes g^{d}x_{1}=0
\end{gather*}%
and we get%
\begin{equation*}
-B(gx_{1}\otimes 1_{H};X_{1}X_{2},x_{1})+B(x_{1}x_{2}\otimes
1_{H};X_{1}X_{2},x_{1}x_{2})=0
\end{equation*}%
and%
\begin{equation*}
B(gx_{1}\otimes 1_{H};GX_{1}X_{2},gx_{1})-B(x_{1}x_{2}\otimes
1_{H};GX_{1}X_{2},gx_{1}x_{2})=0.
\end{equation*}%
By using the form of the element $B\left( gx_{1}\otimes 1_{H}\right) $ we get%
\begin{equation*}
B(x_{1}x_{2}\otimes 1_{H};X_{1}X_{2},x_{1}x_{2})=0
\end{equation*}%
and%
\begin{equation*}
B(x_{1}x_{2}\otimes 1_{H};GX_{1}X_{2},gx_{1}x_{2})=0.
\end{equation*}%
which we already got.

\subsubsection{$G^{a}X_{1}X_{2}\otimes g^{d}x_{1}x_{2}$}

\begin{equation*}
\sum_{\substack{ a,d=0  \\ a+d\equiv 1}}^{1}\left( -1\right)
^{a+1}B(gx_{1}\otimes
1_{H};G^{a}X_{1}X_{2},g^{d}x_{1}x_{2})G^{a}X_{1}X_{2}\otimes
g^{d}x_{1}x_{2}=0
\end{equation*}%
and we get%
\begin{equation*}
-B(gx_{1}\otimes 1_{H};X_{1}X_{2},gx_{1}x_{2})=0
\end{equation*}%
and%
\begin{equation*}
B(gx_{1}\otimes 1_{H};GX_{1}X_{2},x_{1}x_{2})=0
\end{equation*}%
which are already known because of the form of the element $B\left(
gx_{1}\otimes 1_{H}\right) .$

\subsection{Case $gx_{1}x_{2}$}

\begin{eqnarray*}
&&\sum_{a,b_{1},b_{2},d,e_{1},e_{2}=0}^{1}\sum_{l_{1}=0}^{b_{1}}%
\sum_{l_{2}=0}^{b_{2}}\sum_{u_{1}=0}^{e_{1}}\sum_{u_{2}=0}^{e_{2}}\left(
-1\right) ^{\alpha \left( x_{1}x_{2};l_{1},l_{2},u_{1},u_{2}\right) } \\
&&B(1_{H}\otimes
1_{H};G^{a}X_{1}^{b_{1}}X_{2}^{b_{2}},g^{d}x_{1}^{e_{1}}x_{2}^{e_{2}}) \\
&&G^{a}X_{1}^{b_{1}-l_{1}}X_{2}^{b_{2}-l_{2}}\otimes
g^{d}x_{1}^{e_{1}-u_{1}}x_{2}^{e_{2}-u_{2}}\otimes
g^{a+b_{1}+b_{2}+l_{1}+l_{2}+d+e_{1}+e_{2}+u_{1}+u_{2}}x_{1}^{l_{1}+u_{1}+1}x_{2}^{l_{2}+u_{2}+1}
\\
+
&&\sum_{a,b_{1},b_{2},d,e_{1},e_{2}=0}^{1}\sum_{l_{1}=0}^{b_{1}}%
\sum_{l_{2}=0}^{b_{2}}\sum_{u_{1}=0}^{e_{1}}\sum_{u_{2}=0}^{e_{2}}\left(
-1\right) ^{\alpha \left( x_{1};l_{1},l_{2},u_{1},u_{2}\right) } \\
&&B(gx_{2}\otimes
1_{H};G^{a}X_{1}^{b_{1}}X_{2}^{b_{2}},g^{d}x_{1}^{e_{1}}x_{2}^{e_{2}}) \\
&&G^{a}X_{1}^{b_{1}-l_{1}}X_{2}^{b_{2}-l_{2}}\otimes
g^{d}x_{1}^{e_{1}-u_{1}}x_{2}^{e_{2}-u_{2}}\otimes
g^{a+b_{1}+b_{2}+l_{1}+l_{2}+d+e_{1}+e_{2}+u_{1}+u_{2}}x_{1}^{l_{1}+u_{1}+1}x_{2}^{l_{2}+u_{2}}
\\
&&-\sum_{a,b_{1},b_{2},d,e_{1},e_{2}=0}^{1}\sum_{l_{1}=0}^{b_{1}}%
\sum_{l_{2}=0}^{b_{2}}\sum_{u_{1}=0}^{e_{1}}\sum_{u_{2}=0}^{e_{2}}\left(
-1\right) ^{\alpha \left( x_{2};l_{1},l_{2},u_{1},u_{2}\right) } \\
&&B(gx_{1}\otimes
1_{H};G^{a}X_{1}^{b_{1}}X_{2}^{b_{2}},g^{d}x_{1}^{e_{1}}x_{2}^{e_{2}}) \\
&&G^{a}X_{1}^{b_{1}-l_{1}}X_{2}^{b_{2}-l_{2}}\otimes
g^{d}x_{1}^{e_{1}-u_{1}}x_{2}^{e_{2}-u_{2}}\otimes
g^{a+b_{1}+b_{2}+l_{1}+l_{2}+d+e_{1}+e_{2}+u_{1}+u_{2}}x_{1}^{l_{1}+u_{1}}x_{2}^{l_{2}+u_{2}+1}
\\
&&+\sum_{a,b_{1},b_{2},d,e_{1},e_{2}=0}^{1}\sum_{l_{1}=0}^{b_{1}}%
\sum_{l_{2}=0}^{b_{2}}\sum_{u_{1}=0}^{e_{1}}\sum_{u_{2}=0}^{e_{2}}\left(
-1\right) ^{\alpha \left( 1_{H};l_{1},l_{2},u_{1},u_{2}\right)
}B(x_{1}x_{2}\otimes
1_{H};G^{a}X_{1}^{b_{1}}X_{2}^{b_{2}},g^{d}x_{1}^{e_{1}}x_{2}^{e_{2}}) \\
&&G^{a}X_{1}^{b_{1}-l_{1}}X_{2}^{b_{2}-l_{2}}\otimes
g^{d}x_{1}^{e_{1}-u_{1}}x_{2}^{e_{2}-u_{2}}\otimes
g^{a+b_{1}+b_{2}+l_{1}+l_{2}+d+e_{1}+e_{2}+u_{1}+u_{2}}x_{1}^{l_{1}+u_{1}}x_{2}^{l_{2}+u_{2}}
\\
&=&B^{A}(x_{1}x_{2}\otimes 1_{H})\otimes B^{H}(x_{1}x_{2}\otimes
1_{H})\otimes 1_{H}
\end{eqnarray*}%
The first summand gives us $1_{A}\otimes 1_{H}\otimes x_{1}x_{2.}$ The
second summand gives us%
\begin{eqnarray*}
&&\sum_{a,b_{1},b_{2},d,e_{1},e_{2}=0}^{1}\sum_{l_{1}=0}^{b_{1}}%
\sum_{l_{2}=0}^{b_{2}}\sum_{u_{1}=0}^{e_{1}}\sum_{u_{2}=0}^{e_{2}}\left(
-1\right) ^{\alpha \left( x_{1};l_{1},l_{2},u_{1},u_{2}\right)
}B(gx_{2}\otimes
1_{H};G^{a}X_{1}^{b_{1}}X_{2}^{b_{2}},g^{d}x_{1}^{e_{1}}x_{2}^{e_{2}}) \\
&&G^{a}X_{1}^{b_{1}-l_{1}}X_{2}^{b_{2}-l_{2}}\otimes
g^{d}x_{1}^{e_{1}-u_{1}}x_{2}^{e_{2}-u_{2}}\otimes
g^{a+b_{1}+b_{2}+l_{1}+l_{2}+d+e_{1}+e_{2}+u_{1}+u_{2}}x_{1}^{l_{1}+u_{1}+1}x_{2}^{l_{2}+u_{2}}
\end{eqnarray*}%
\begin{eqnarray*}
a+b_{1}+b_{2}+l_{1}+l_{2}+d+e_{1}+e_{2}+u_{1}+u_{2} &\equiv &1 \\
l_{1}+u_{1}+1 &=&1 \\
l_{2}+u_{2} &=&1
\end{eqnarray*}%
i.e.%
\begin{eqnarray*}
a+b_{1}+b_{2}+d+e_{1}+e_{2} &\equiv &0 \\
l_{1} &=&u_{1}=1 \\
l_{2}+u_{2} &=&1
\end{eqnarray*}%
This summand is zero since, by $\left( \ref{gx2ot1, first}\right) ,$ $%
B(gx_{2}\otimes \
1_{H};G^{a}X_{1}^{b_{1}}X_{2}^{b_{2}},g^{d}x_{1}^{e_{1}}x_{2}^{e_{2}})=0$
whenever $a+b_{1}+b_{2}+d+e_{1}+e_{2}\equiv 0.$

An analogous statement holds for the third summand.

Let us consider the third summand%
\begin{eqnarray*}
&&\sum_{a,b_{1},b_{2},d,e_{1},e_{2}=0}^{1}\sum_{l_{1}=0}^{b_{1}}%
\sum_{l_{2}=0}^{b_{2}}\sum_{u_{1}=0}^{e_{1}}\sum_{u_{2}=0}^{e_{2}}\left(
-1\right) ^{\alpha \left( 1_{H};l_{1},l_{2},u_{1},u_{2}\right) } \\
&&B(x_{1}x_{2}\otimes
1_{H};G^{a}X_{1}^{b_{1}}X_{2}^{b_{2}},g^{d}x_{1}^{e_{1}}x_{2}^{e_{2}}) \\
&&G^{a}X_{1}^{b_{1}-l_{1}}X_{2}^{b_{2}-l_{2}}\otimes
g^{d}x_{1}^{e_{1}-u_{1}}x_{2}^{e_{2}-u_{2}}\otimes
g^{a+b_{1}+b_{2}+l_{1}+l_{2}+d+e_{1}+e_{2}+u_{1}+u_{2}}x_{1}^{l_{1}+u_{1}}x_{2}^{l_{2}+u_{2}}.
\end{eqnarray*}%
We get%
\begin{eqnarray*}
a+b_{1}+b_{2}+l_{1}+l_{2}+d+e_{1}+e_{2}+u_{1}+u_{2} &\equiv &1 \\
l_{1}+u_{1} &=&1 \\
l_{2}+u_{2} &=&1
\end{eqnarray*}%
i.e.%
\begin{eqnarray*}
a+b_{1}+b_{2}+d+e_{1}+e_{2} &\equiv &1 \\
l_{1}+u_{1} &=&1 \\
l_{2}+u_{2} &=&1
\end{eqnarray*}%
and we obtain%
\begin{eqnarray*}
&&\sum_{\substack{ a,b_{1},b_{2},d,e_{1},e_{2}=0  \\ %
a+b_{1}+b_{2}+d+e_{1}+e_{2}\equiv 1}}^{1}\sum_{l_{1}=0}^{b_{1}}%
\sum_{l_{2}=0}^{b_{2}}\sum_{\substack{ u_{1}=0  \\ l_{1}+u_{1}=1}}%
^{e_{1}}\sum _{\substack{ u_{2}=0  \\ l_{2}+u_{2}=1}}^{e_{2}}\left(
-1\right) ^{\alpha \left( 1_{H};l_{1},l_{2},u_{1},u_{2}\right) } \\
&&B(x_{1}x_{2}\otimes
1_{H};G^{a}X_{1}^{b_{1}}X_{2}^{b_{2}},g^{d}x_{1}^{e_{1}}x_{2}^{e_{2}})G^{a}X_{1}^{b_{1}-l_{1}}X_{2}^{b_{2}-l_{2}}\otimes g^{d}x_{1}^{e_{1}-u_{1}}x_{2}^{e_{2}-u_{2}}
\end{eqnarray*}%
i.e.%
\begin{eqnarray*}
&&\sum_{\substack{ a,b_{1},b_{2},d=0  \\ a+b_{1}+b_{2}+d\equiv 1}}^{1}\left(
-1\right) ^{\alpha \left( 1_{H};0,0,1,1\right) }B(x_{1}x_{2}\otimes
1_{H};G^{a}X_{1}^{b_{1}}X_{2}^{b_{2}},g^{d}x_{1}x_{2})G^{a}X_{1}^{b_{1}}X_{2}^{b_{2}}\otimes g^{d}+
\\
&&+\sum_{\substack{ a,b_{1},d,e_{2}=0  \\ a+b_{1}+d+e_{2}\equiv 1}}%
^{1}\left( -1\right) ^{\alpha \left( 1_{H};0,1,1,0\right)
}B(x_{1}x_{2}\otimes
1_{H};G^{a}X_{1}^{b_{1}}X_{2},g^{d}x_{1}x_{2}^{e_{2}})G^{a}X_{1}^{b_{1}}%
\otimes g^{d}x_{2}^{e_{2}}+ \\
&&+\sum_{\substack{ a,b_{2},d,e_{1}=0  \\ a+b_{2}+d+e_{1}\equiv 1}}%
^{1}\left( -1\right) ^{\alpha \left( 1_{H};1,0,0,1\right)
}B(x_{1}x_{2}\otimes
1_{H};G^{a}X_{1}X_{2}^{b_{2}},g^{d}x_{1}^{e_{1}}x_{2})G^{a}X_{2}^{b_{2}}%
\otimes g^{d}x_{1}^{e_{1}}+ \\
&&+\sum_{\substack{ a,d,e_{1},e_{2}=0  \\ a+d+e_{1}+e_{2}\equiv 1}}%
^{1}\left( -1\right) ^{\alpha \left( 1_{H};1,1,0,0\right)
}B(x_{1}x_{2}\otimes
1_{H};G^{a}X_{1}X_{2},g^{d}x_{1}^{e_{1}}x_{2}^{e_{2}})G^{a}\otimes
g^{d}x_{1}^{e_{1}}x_{2}^{e_{2}}.
\end{eqnarray*}%
By taking in account that $\alpha \left( 1_{H};0,0,1,1\right) =1+e_{2}$, $%
\alpha \left( 1_{H};0,1,1,0\right) =e_{2}+\left( a+b_{1}+b_{2}+1\right) $, $%
\alpha \left( 1_{H};1,0,0,1\right) =a+b_{1}$ and $\alpha \left(
1_{H};1,1,0,0\right) =1+b_{2},$ we finally get%
\begin{eqnarray*}
&&\sum_{\substack{ a,b_{1},b_{2},d=0  \\ a+b_{1}+b_{2}+d\equiv 1}}%
^{1}B(x_{1}x_{2}\otimes
1_{H};G^{a}X_{1}^{b_{1}}X_{2}^{b_{2}},g^{d}x_{1}x_{2})G^{a}X_{1}^{b_{1}}X_{2}^{b_{2}}\otimes g^{d}+
\\
&&+\sum_{\substack{ a,b_{1},d,e_{2}=0  \\ a+b_{1}+d+e_{2}\equiv 1}}%
^{1}\left( -1\right) ^{e_{2}+a+b_{1}}B(x_{1}x_{2}\otimes
1_{H};G^{a}X_{1}^{b_{1}}X_{2},g^{d}x_{1}x_{2}^{e_{2}})G^{a}X_{1}^{b_{1}}%
\otimes g^{d}x_{2}^{e_{2}}+ \\
&&+\sum_{\substack{ a,b_{2},d,e_{1}=0  \\ a+b_{2}+d+e_{1}\equiv 1}}%
^{1}\left( -1\right) ^{a+1}B(x_{1}x_{2}\otimes
1_{H};G^{a}X_{1}X_{2}^{b_{2}},g^{d}x_{1}^{e_{1}}x_{2})G^{a}X_{2}^{b_{2}}%
\otimes g^{d}x_{1}^{e_{1}}+ \\
&&+\sum_{\substack{ a,d,e_{1},e_{2}=0  \\ a+d+e_{1}+e_{2}\equiv 1}}%
^{1}B(x_{1}x_{2}\otimes
1_{H};G^{a}X_{1}X_{2},g^{d}x_{1}^{e_{1}}x_{2}^{e_{2}})G^{a}\otimes
g^{d}x_{1}^{e_{1}}x_{2}^{e_{2}}=0
\end{eqnarray*}

\subsubsection{$G^{a}\otimes g^{d}$}

\begin{equation*}
\sum_{\substack{ a,d=0  \\ a+d\equiv 1}}^{1}\left[
\begin{array}{c}
B(x_{1}x_{2}\otimes 1_{H};G^{a},g^{d}x_{1}x_{2})+\left( -1\right)
^{a}B(x_{1}x_{2}\otimes 1_{H};G^{a}X_{2},g^{d}x_{1})+ \\
+\left( -1\right) ^{a+1}B(x_{1}x_{2}\otimes
1_{H};G^{a}X_{1},g^{d}x_{2})+B(x_{1}x_{2}\otimes 1_{H};G^{a}X_{1}X_{2},g^{d})%
\end{array}%
\right] G^{a}\otimes g^{d}=0
\end{equation*}%
and we get%
\begin{equation}
\begin{array}{c}
B(x_{1}x_{2}\otimes 1_{H};1_{A},gx_{1}x_{2})+B(x_{1}x_{2}\otimes
1_{H};X_{2},gx_{1})+ \\
-B(x_{1}x_{2}\otimes 1_{H};X_{1},gx_{2})+B(x_{1}x_{2}\otimes
1_{H};X_{1}X_{2},g)%
\end{array}%
=0  \label{x1otx2, fortytwo}
\end{equation}%
and%
\begin{equation}
\begin{array}{c}
B(x_{1}x_{2}\otimes 1_{H};G,x_{1}x_{2})-B(x_{1}x_{2}\otimes
1_{H};GX_{2},x_{1})+ \\
+B(x_{1}x_{2}\otimes 1_{H};GX_{1},x_{2})+B(x_{1}x_{2}\otimes
1_{H};GX_{1}X_{2},1_{H})%
\end{array}%
=0.  \label{x1otx2, fortythree}
\end{equation}

\subsubsection{$G^{a}\otimes g^{d}x_{2}$}

\begin{equation*}
\sum_{\substack{ a,d=0  \\ a+d\equiv 0}}^{1}\left[ \left( -1\right)
^{a+1}B(x_{1}x_{2}\otimes
1_{H};G^{a}X_{2},g^{d}x_{1}x_{2})+B(x_{1}x_{2}\otimes
1_{H};G^{a}X_{1}X_{2},g^{d}x_{2})\right] G^{a}\otimes g^{d}x_{2}=0
\end{equation*}%
and we get%
\begin{equation}
-B(x_{1}x_{2}\otimes 1_{H};X_{2},x_{1}x_{2})+B(x_{1}x_{2}\otimes
1_{H};X_{1}X_{2},x_{2})=0  \label{x1otx2, fortyfour}
\end{equation}%
and%
\begin{equation}
B(x_{1}x_{2}\otimes 1_{H};GX_{2},gx_{1}x_{2})+B(x_{1}x_{2}\otimes
1_{H};GX_{1}X_{2},gx_{2})=0.  \label{x1otx2, fortyfive}
\end{equation}

\subsubsection{$G^{a}\otimes g^{d}x_{1}$}

\begin{equation*}
\sum_{\substack{ a,d=0  \\ a+d\equiv 0}}^{1}\left[ \left( -1\right)
^{a+1}B(x_{1}x_{2}\otimes
1_{H};G^{a}X_{1},g^{d}x_{1}x_{2})+B(x_{1}x_{2}\otimes
1_{H};G^{a}X_{1}X_{2},g^{d}x_{1})\right] G^{a}\otimes g^{d}x_{1}=0
\end{equation*}%
and we get%
\begin{equation}
-B(x_{1}x_{2}\otimes 1_{H};X_{1},x_{1}x_{2})+B(x_{1}x_{2}\otimes
1_{H};X_{1}X_{2},x_{1})=0  \label{x1otx2, fortysix}
\end{equation}%
and%
\begin{equation}
B(x_{1}x_{2}\otimes 1_{H};GX_{1},gx_{1}x_{2})+B(x_{1}x_{2}\otimes
1_{H};GX_{1}X_{2},gx_{1})=0.  \label{x1otx2, fortyseven}
\end{equation}

\subsubsection{$G^{a}X_{2}\otimes g^{d}$}

\begin{equation*}
\sum_{\substack{ a,d=0  \\ a+d\equiv 0}}^{1}\left[ B(x_{1}x_{2}\otimes
1_{H};G^{a}X_{2},g^{d}x_{1}x_{2})+\left( -1\right) ^{a+1}B(x_{1}x_{2}\otimes
1_{H};G^{a}X_{1}X_{2},g^{d}x_{2})\right] G^{a}X_{2}\otimes g^{d}=0
\end{equation*}%
and we get%
\begin{equation*}
B(x_{1}x_{2}\otimes 1_{H};X_{2},x_{1}x_{2})-B(x_{1}x_{2}\otimes
1_{H};X_{1}X_{2},x_{2})=0
\end{equation*}%
and%
\begin{equation*}
B(x_{1}x_{2}\otimes 1_{H};GX_{2},gx_{1}x_{2})+B(x_{1}x_{2}\otimes
1_{H};GX_{1}X_{2},gx_{2})=0.
\end{equation*}%
These are respectively $\left( \ref{x1otx2, fortyfour}\right) $ and $\left( %
\ref{x1otx2, fortyfive}\right) .$

\subsubsection{$G^{a}X_{1}\otimes g^{d}$}

\begin{equation*}
\sum_{\substack{ a,d=0  \\ a+d\equiv 0}}^{1}\left[ B(x_{1}x_{2}\otimes
1_{H};G^{a}X_{1},g^{d}x_{1}x_{2})+\left( -1\right) ^{a+1}B(x_{1}x_{2}\otimes
1_{H};G^{a}X_{1}X_{2},g^{d}x_{1})\right] G^{a}X_{1}\otimes g^{d}=0
\end{equation*}%
and we get%
\begin{equation*}
B(x_{1}x_{2}\otimes 1_{H};X_{1},x_{1}x_{2})-B(x_{1}x_{2}\otimes
1_{H};X_{1}X_{2},x_{1})=0
\end{equation*}%
and%
\begin{equation*}
B(x_{1}x_{2}\otimes 1_{H};GX_{1},gx_{1}x_{2})+B(x_{1}x_{2}\otimes
1_{H};GX_{1}X_{2},gx_{1})=0
\end{equation*}%
these are respectively $\left( \ref{x1otx2, fortysix}\right) $ and $\left( %
\ref{x1otx2, fortyseven}\right) .$

\subsubsection{$G^{a}\otimes g^{d}x_{1}x_{2}$}

\begin{equation*}
\sum_{\substack{ a,d=0  \\ a+d\equiv 1}}^{1}B(x_{1}x_{2}\otimes
1_{H};G^{a}X_{1}X_{2},g^{d}x_{1}x_{2})G^{a}\otimes g^{d}x_{1}x_{2}=0
\end{equation*}%
and we get%
\begin{equation}
B(x_{1}x_{2}\otimes 1_{H};X_{1}X_{2},gx_{1}x_{2})=0
\label{x1otx2, fortyeight}
\end{equation}%
and%
\begin{equation}
B(x_{1}x_{2}\otimes 1_{H};GX_{1}X_{2},x_{1}x_{2})=0.
\label{x1otx2, fortynine}
\end{equation}

\subsubsection{$G^{a}X_{2}\otimes g^{d}x_{2}$}

We do not have any summand like this.

\subsubsection{$G^{a}X_{1}\otimes g^{d}x_{2}$}

\begin{equation*}
\sum_{\substack{ a,d=0  \\ a+d\equiv 1}}^{1}\left( -1\right)
^{a}B(x_{1}x_{2}\otimes
1_{H};G^{a}X_{1}X_{2},g^{d}x_{1}x_{2})G^{a}X_{1}\otimes g^{d}x_{2}=0
\end{equation*}%
and we get%
\begin{equation*}
B(x_{1}x_{2}\otimes 1_{H};X_{1}X_{2},gx_{1}x_{2})=0
\end{equation*}%
and%
\begin{equation*}
-B(x_{1}x_{2}\otimes 1_{H};GX_{1}X_{2},x_{1}x_{2})=0
\end{equation*}%
which we just got.

\subsubsection{$G^{a}X_{2}\otimes g^{d}x_{1}$}

\begin{equation*}
\sum_{\substack{ a,d=0  \\ a+d\equiv 1}}^{1}\left( -1\right)
^{a+1}B(x_{1}x_{2}\otimes
1_{H};G^{a}X_{1}X_{2},g^{d}x_{1}x_{2})G^{a}X_{2}\otimes g^{d}x_{1}=0
\end{equation*}%
and we get%
\begin{equation*}
B(x_{1}x_{2}\otimes 1_{H};X_{1}X_{2},gx_{1}x_{2})=0
\end{equation*}%
and%
\begin{equation*}
-B(x_{1}x_{2}\otimes 1_{H};GX_{1}X_{2},x_{1}x_{2})=0
\end{equation*}%
which we just got.

\subsubsection{$G^{a}X_{1}\otimes g^{d}x_{1}$}

There is no term like this.

\subsubsection{$G^{a}X_{1}X_{2}\otimes g^{d}$}

\begin{equation*}
\sum_{\substack{ a,d=0  \\ a+d\equiv 1}}^{1}B(x_{1}x_{2}\otimes
1_{H};G^{a}X_{1}X_{2},g^{d}x_{1}x_{2})G^{a}X_{1}X_{2}\otimes g^{d}=0
\end{equation*}%
and we get%
\begin{equation*}
B(x_{1}x_{2}\otimes 1_{H};X_{1}X_{2},gx_{1}x_{2})=0
\end{equation*}%
and%
\begin{equation*}
B(x_{1}x_{2}\otimes 1_{H};GX_{1}X_{2},x_{1}x_{2})=0
\end{equation*}%
which we just got.

\subsubsection{$G^{a}X_{2}\otimes g^{d}x_{1}x_{2}$}

There is no term like this.

\subsubsection{$G^{a}X_{1}\otimes g^{d}x_{1}x_{2}$}

There is no term like this.

\subsubsection{$G^{a}X_{1}X_{2}\otimes g^{d}x_{2}$}

There is no term like this.

\subsubsection{$G^{a}X_{1}X_{2}\otimes g^{d}x_{1}$}

There is no term like this.

\subsubsection{$G^{a}X_{1}X_{2}\otimes g^{d}x_{1}x_{2}$}

There is no term like this.

By using equalities from $\left( \ref{x1otx2, first}\right) $ to $\left( \ref%
{x1otx2, fortynine}\right) $ we obtain the following form of $B\left(
x_{1}x_{2}\otimes 1_{H}\right) .$

\subsection{The final form of the element $B(x_{1}x_{2}\otimes 1_{H})$}

\begin{eqnarray}
B(x_{1}x_{2}\otimes 1_{H}) &=&B(x_{1}x_{2}\otimes
1_{H};1_{A},1_{H})1_{A}\otimes 1_{H}+  \label{x1x2} \\
&&+B(x_{1}x_{2}\otimes 1_{H};1_{A},x_{1}x_{2})1_{A}\otimes x_{1}x_{2}+
\notag \\
&&+B(x_{1}x_{2}\otimes 1_{H};1_{A},gx_{1})1_{A}\otimes gx_{1}+  \notag \\
&&+B(x_{1}x_{2}\otimes 1_{H};1_{A},gx_{2})1_{A}\otimes gx_{2}+  \notag \\
&&+B(x_{1}x_{2}\otimes 1_{H};G,g)G\otimes g+  \notag \\
&&+B(x_{1}x_{2}\otimes 1_{H};G,x_{1})G\otimes x_{1}+  \notag \\
&&+B(x_{1}x_{2}\otimes 1_{H};G,x_{2})G\otimes x_{2}+  \notag \\
&&+B(x_{1}x_{2}\otimes 1_{H};G,gx_{1}x_{2})G\otimes gx_{1}x_{2}+  \notag \\
&&+B(x_{1}x_{2}\otimes 1_{H};X_{1},g)X_{1}\otimes g+  \notag \\
&&+B(x_{1}x_{2}\otimes 1_{H};X_{1},x_{1})X_{1}\otimes x_{1}+  \notag \\
&&+B(x_{1}x_{2}\otimes 1_{H};X_{1},x_{2})X_{1}\otimes x_{2}+  \notag \\
&&+B(x_{1}x_{2}\otimes 1_{H};X_{1},gx_{1}x_{2})X_{1}\otimes gx_{1}x_{2}+
\notag \\
&&+B(x_{1}x_{2}\otimes 1_{H};X_{2},g)X_{2}\otimes g+  \notag \\
&&+B(x_{1}x_{2}\otimes 1_{H};X_{2},x_{1})X_{2}\otimes x_{1}+  \notag \\
&&+B(x_{1}x_{2}\otimes 1_{H};X_{2},x_{2})X_{2}\otimes x_{2}+  \notag \\
&&+B(x_{1}x_{2}\otimes 1_{H};X_{2},gx_{1}x_{2})X_{2}\otimes gx_{1}x_{2}+
\notag \\
&&+\left[ +1-B(x_{1}x_{2}\otimes 1_{H};1_{A},x_{1}x_{2})-B(x_{1}x_{2}\otimes
1_{H};X_{2},x_{1})+B(x_{1}x_{2}\otimes 1_{H};X_{1},x_{2})\right]
X_{1}X_{2}\otimes 1_{H}+  \notag \\
&&+B(x_{1}x_{2}\otimes 1_{H};X_{1},gx_{1}x_{2})X_{1}X_{2}\otimes gx_{1}
\notag \\
&&+B(x_{1}x_{2}\otimes 1_{H};X_{2},gx_{1}x_{2})X_{1}X_{2}\otimes gx_{2}
\notag \\
&&+B(x_{1}x_{2}\otimes 1_{H};GX_{1},1_{H})GX_{1}\otimes 1_{H}+  \notag \\
&&+B(x_{1}x_{2}\otimes 1_{H};GX_{1},x_{1}x_{2})GX_{1}\otimes x_{1}x_{2}
\notag \\
&&+B(x_{1}x_{2}\otimes 1_{H};GX_{1},gx_{1})GX_{1}\otimes gx_{1}  \notag \\
&&+B(x_{1}x_{2}\otimes 1_{H};GX_{1},gx_{2})GX_{1}\otimes gx_{2}+  \notag \\
&&+B(x_{1}x_{2}\otimes 1_{H};GX_{2},1_{H})GX_{2}\otimes 1_{H}+  \notag \\
&&+B(x_{1}x_{2}\otimes 1_{H};GX_{2},x_{1}x_{2})GX_{2}\otimes x_{1}x_{2}+
\notag \\
&&+B(x_{1}x_{2}\otimes 1_{H};GX_{2},gx_{1})GX_{2}\otimes gx_{1}+  \notag \\
&&+B(x_{1}x_{2}\otimes 1_{H};GX_{2},gx_{2})GX_{2}\otimes gx_{2}+  \notag \\
&&+\left[ -B(x_{1}x_{2}\otimes 1_{H};G,gx_{1}x_{2})+B(x_{1}x_{2}\otimes
1_{H};GX_{2},gx_{1})-B(x_{1}x_{2}\otimes 1_{H};GX_{1},gx_{2})\right]
GX_{1}X_{2}\otimes g  \notag \\
&&-B(x_{1}x_{2}\otimes 1_{H};GX_{1},x_{1}x_{2})GX_{1}X_{2}\otimes x_{1}+
\notag \\
&&-B(x_{1}x_{2}\otimes 1_{H};GX_{2},x_{1}x_{2})GX_{1}X_{2}\otimes x_{2}
\notag
\end{eqnarray}

\section{$B\left( gx_{1}\otimes 1_{H}\right) $ and $B\left( gx_{2}\otimes
1_{H}\right) $ rewritten}

Now we use equalities from $\left( \ref{x1otx2, first}\right) $ to $\left( %
\ref{x1otx2, fortynine}\right) $ to rewrite $B\left( gx_{1}\otimes
1_{H}\right) $ and $B\left( gx_{2}\otimes 1_{H}\right) .$

\begin{eqnarray}
B\left( gx_{1}\otimes 1_{H}\right) &=&\left[ B(x_{1}x_{2}\otimes
1_{H};1_{A},gx_{2})+B(x_{1}x_{2}\otimes 1_{H};X_{2},g)\right] 1_{A}\otimes g+
\label{gx1} \\
&&+\left[ B(x_{1}x_{2}\otimes 1_{H};1_{A},x_{1}x_{2})+B(x_{1}x_{2}\otimes
1_{H};X_{2},x_{1})\right] 1_{A}\otimes x_{1}  \notag \\
&&+B(x_{1}x_{2}\otimes 1_{H};X_{2},x_{2})1_{A}\otimes x_{2}  \notag \\
&&+B(x_{1}x_{2}\otimes 1_{H};X_{2},gx_{1}x_{2})1_{A}\otimes gx_{1}x_{2}
\notag \\
&&+\left[ B(x_{1}x_{2}\otimes 1_{H};G,x_{2})-B(x_{1}x_{2}\otimes
1_{H};GX_{2},1_{H})\right] G\otimes 1_{H}+  \notag \\
&&-B(x_{1}x_{2}\otimes 1_{H};GX_{2},x_{1}x_{2})G\otimes x_{1}x_{2}+  \notag
\\
&&+\left[ B(x_{1}x_{2}\otimes 1_{H};G,gx_{1}x_{2})-B(x_{1}x_{2}\otimes
1_{H};GX_{2},gx_{1})\right] G\otimes gx_{1}+  \notag \\
&&-B(x_{1}x_{2}\otimes 1_{H};GX_{2},gx_{2})G\otimes gx_{2}+  \notag \\
&&+\left[ -1+B(x_{1}x_{2}\otimes 1_{H};1_{A},x_{1}x_{2})+B(x_{1}x_{2}\otimes
1_{H};X_{2},x_{1})\right] X_{1}\otimes 1_{H}+  \notag \\
&&-B(x_{1}x_{2}\otimes 1_{H};X_{2},gx_{1}x_{2})X_{1}\otimes gx_{2}+  \notag
\\
&&+B(x_{1}x_{2}\otimes 1_{H};X_{2},x_{2})X_{2}\otimes 1_{H}+  \notag \\
&&+B(x_{1}x_{2}\otimes 1_{H};X_{2},gx_{1}x_{2})X_{2}\otimes gx_{1}+  \notag
\\
&&+B(x_{1}x_{2}\otimes 1_{H};X_{2},gx_{1}x_{2})X_{1}X_{2}\otimes g+  \notag
\\
&&\left[ -B(x_{1}x_{2}\otimes 1_{H};G,gx_{1}x_{2})+B(x_{1}x_{2}\otimes
1_{H};GX_{2},gx_{1})\right] GX_{1}\otimes g+  \notag \\
&&-B(x_{1}x_{2}\otimes 1_{H};GX_{2},x_{1}x_{2})GX_{1}\otimes x_{2}+  \notag
\\
&&B(x_{1}x_{2}\otimes 1_{H};GX_{2},gx_{2})GX_{2}\otimes g+  \notag \\
&&+B(x_{1}x_{2}\otimes 1_{H};GX_{2},x_{1}x_{2})GX_{2}\otimes x_{1}+  \notag
\\
&&-B(x_{1}x_{2}\otimes 1_{H};GX_{2},x_{1}x_{2})GX_{1}X_{2}\otimes 1_{H}
\notag
\end{eqnarray}%
\begin{eqnarray}
B\left( gx_{2}\otimes 1_{H}\right) &=&\left[ -B(x_{1}x_{2}\otimes
1_{H};1_{A},gx_{1})-B(x_{1}x_{2}\otimes 1_{H};X_{1},g)\right] 1_{A}\otimes g+
\label{gx2} \\
&&-B(x_{1}x_{2}\otimes 1_{H};X_{1},x_{1})1_{A}\otimes x_{1}  \notag \\
&&+\left[ B(x_{1}x_{2}\otimes 1_{H};1_{A},x_{1}x_{2})-B(x_{1}x_{2}\otimes
1_{H};X_{1},x_{2})\right] 1_{A}\otimes x_{2}  \notag \\
&&-B(x_{1}x_{2}\otimes 1_{H};X_{1},gx_{1}x_{2})1_{A}\otimes gx_{1}x_{2}
\notag \\
&&+\left[ -B(x_{1}x_{2}\otimes 1_{H};G,x_{1})+B(x_{1}x_{2}\otimes
1_{H};GX_{1},1_{H})\right] G\otimes 1_{H}+  \notag \\
&&+B(x_{1}x_{2}\otimes 1_{H};GX_{1},x_{1}x_{2})G\otimes x_{1}x_{2}+  \notag
\\
&&+B(x_{1}x_{2}\otimes 1_{H};GX_{1},gx_{1})G\otimes gx_{1}+  \notag \\
&&+\left[ B(x_{1}x_{2}\otimes 1_{H};G,gx_{1}x_{2})+B(x_{1}x_{2}\otimes
1_{H};GX_{1},gx_{2})\right] G\otimes gx_{2}+  \notag \\
&&-B(x_{1}x_{2}\otimes 1_{H};X_{1},x_{1})X_{1}\otimes 1_{H}+  \notag \\
&&+B(x_{1}x_{2}\otimes 1_{H};X_{1},gx_{1}x_{2})X_{1}\otimes gx_{2}+  \notag
\\
&&\left[ -1+B(x_{1}x_{2}\otimes 1_{H};1_{A},x_{1}x_{2})-B(x_{1}x_{2}\otimes
1_{H};X_{1},x_{2})\right] X_{2}\otimes 1_{H}+  \notag \\
&&-B(x_{1}x_{2}\otimes 1_{H};X_{1},gx_{1}x_{2})X_{2}\otimes gx_{1}+  \notag
\\
&&-B(x_{1}x_{2}\otimes 1_{H};X_{1},gx_{1}x_{2})X_{1}X_{2}\otimes g+  \notag
\\
&&-B(x_{1}x_{2}\otimes 1_{H};GX_{1},gx_{1})GX_{1}\otimes g+  \notag \\
&&+B(x_{1}x_{2}\otimes 1_{H};GX_{1},x_{1}x_{2})GX_{1}\otimes x_{2}+  \notag
\\
&&+\left[ -B(x_{1}x_{2}\otimes 1_{H};G,gx_{1}x_{2})-B(x_{1}x_{2}\otimes
1_{H};GX_{1},gx_{2})\right] GX_{2}\otimes g+  \notag \\
&&-B(x_{1}x_{2}\otimes 1_{H};GX_{1},x_{1}x_{2})GX_{2}\otimes x_{1}+  \notag
\\
&&+B(x_{1}x_{2}\otimes 1_{H};GX_{1},x_{1}x_{2})GX_{1}X_{2}\otimes 1_{H}
\notag
\end{eqnarray}

\section{$B(gx_{1}x_{2}\otimes 1_{H})$}

We write the Casimir formula $\left( \ref{MAIN FORMULA 1}\right) $ for $%
B(gx_{1}x_{2}\otimes 1_{H}).$

\begin{eqnarray*}
&&\sum_{a,b_{1},b_{2},d,e_{1},e_{2}=0}^{1}\sum_{l_{1}=0}^{b_{1}}%
\sum_{l_{2}=0}^{b_{2}}\sum_{u_{1}=0}^{e_{1}}\sum_{u_{2}=0}^{e_{2}}\left(
-1\right) ^{\alpha \left( gx_{1}x_{2};l_{1},l_{2},u_{1},u_{2}\right) } \\
&&B(g\otimes
1_{H};G^{a}X_{1}^{b_{1}}X_{2}^{b_{2}},g^{d}x_{1}^{e_{1}}x_{2}^{e_{2}}) \\
&&G^{a}X_{1}^{b_{1}-l_{1}}X_{2}^{b_{2}-l_{2}}\otimes
g^{d}x_{1}^{e_{1}-u_{1}}x_{2}^{e_{2}-u_{2}}\otimes
g^{a+b_{1}+b_{2}+l_{1}+l_{2}+d+e_{1}+e_{2}+u_{1}+u_{2}+1}x_{1}^{l_{1}+u_{1}+1}x_{2}^{l_{2}+u_{2}+1}
\\
&&\sum_{a,b_{1},b_{2},d,e_{1},e_{2}=0}^{1}\sum_{l_{1}=0}^{b_{1}}%
\sum_{l_{2}=0}^{b_{2}}\sum_{u_{1}=0}^{e_{1}}\sum_{u_{2}=0}^{e_{2}}\left(
-1\right) ^{\alpha \left( gx_{1};l_{1},l_{2},u_{1},u_{2}\right) } \\
&&B(x_{2}\otimes
1_{H};G^{a}X_{1}^{b_{1}}X_{2}^{b_{2}},g^{d}x_{1}^{e_{1}}x_{2}^{e_{2}}) \\
&&G^{a}X_{1}^{b_{1}-l_{1}}X_{2}^{b_{2}-l_{2}}\otimes
g^{d}x_{1}^{e_{1}-u_{1}}x_{2}^{e_{2}-u_{2}}\otimes
g^{a+b_{1}+b_{2}+l_{1}+l_{2}+d+e_{1}+e_{2}+u_{1}+u_{2}+1}x_{1}^{l_{1}+u_{1}+1}x_{2}^{l_{2}+u_{2}}
\\
&&-\sum_{a,b_{1},b_{2},d,e_{1},e_{2}=0}^{1}\sum_{l_{1}=0}^{b_{1}}%
\sum_{l_{2}=0}^{b_{2}}\sum_{u_{1}=0}^{e_{1}}\sum_{u_{2}=0}^{e_{2}}\left(
-1\right) ^{\alpha \left( gx_{2};l_{1},l_{2},u_{1},u_{2}\right) } \\
&&B(x_{1}\otimes
1_{H};G^{a}X_{1}^{b_{1}}X_{2}^{b_{2}},g^{d}x_{1}^{e_{1}}x_{2}^{e_{2}}) \\
&&G^{a}X_{1}^{b_{1}-l_{1}}X_{2}^{b_{2}-l_{2}}\otimes
g^{d}x_{1}^{e_{1}-u_{1}}x_{2}^{e_{2}-u_{2}}\otimes
g^{a+b_{1}+b_{2}+l_{1}+l_{2}+d+e_{1}+e_{2}+u_{1}+u_{2}+1}x_{1}^{l_{1}+u_{1}}x_{2}^{l_{2}+u_{2}+1}
\\
&&+\sum_{a,b_{1},b_{2},d,e_{1},e_{2}=0}^{1}\sum_{l_{1}=0}^{b_{1}}%
\sum_{l_{2}=0}^{b_{2}}\sum_{u_{1}=0}^{e_{1}}\sum_{u_{2}=0}^{e_{2}}\left(
-1\right) ^{\alpha \left( g;l_{1},l_{2},u_{1},u_{2}\right)
}B(gx_{1}x_{2}\otimes
1_{H};G^{a}X_{1}^{b_{1}}X_{2}^{b_{2}},g^{d}x_{1}^{e_{1}}x_{2}^{e_{2}}) \\
&&G^{a}X_{1}^{b_{1}-l_{1}}X_{2}^{b_{2}-l_{2}}\otimes
g^{d}x_{1}^{e_{1}-u_{1}}x_{2}^{e_{2}-u_{2}}\otimes
g^{a+b_{1}+b_{2}+l_{1}+l_{2}+d+e_{1}+e_{2}+u_{1}+u_{2}+1}x_{1}^{l_{1}+u_{1}}x_{2}^{l_{2}+u_{2}}
\\
&=&B^{A}(gx_{1}x_{2}\otimes 1_{H})\otimes B^{H}(gx_{1}x_{2}\otimes
1_{H})\otimes 1_{H}
\end{eqnarray*}

\subsection{Case $1_{H}$}

We have to consider only the last summand.%
\begin{eqnarray*}
&&\sum_{a,b_{1},b_{2},d,e_{1},e_{2}=0}^{1}\sum_{l_{1}=0}^{b_{1}}%
\sum_{l_{2}=0}^{b_{2}}\sum_{u_{1}=0}^{e_{1}}\sum_{u_{2}=0}^{e_{2}}\left(
-1\right) ^{\alpha \left( g;l_{1},l_{2},u_{1},u_{2}\right) } \\
&&B(gx_{1}x_{2}\otimes
1_{H};G^{a}X_{1}^{b_{1}}X_{2}^{b_{2}},g^{d}x_{1}^{e_{1}}x_{2}^{e_{2}}) \\
&&G^{a}X_{1}^{b_{1}-l_{1}}X_{2}^{b_{2}-l_{2}}\otimes
g^{d}x_{1}^{e_{1}-u_{1}}x_{2}^{e_{2}-u_{2}}\otimes
g^{a+b_{1}+b_{2}+l_{1}+l_{2}+d+e_{1}+e_{2}+u_{1}+u_{2}+1}x_{1}^{l_{1}+u_{1}}x_{2}^{l_{2}+u_{2}}
\end{eqnarray*}%
and we get%
\begin{eqnarray*}
a+b_{1}+b_{2}+l_{1}+l_{2}+d+e_{1}+e_{2}+u_{1}+u_{2}+1 &\equiv &0 \\
l_{1}+u_{1} &=&0 \\
l_{2}+u_{2} &=&0
\end{eqnarray*}%
i.e.%
\begin{eqnarray*}
a+b_{1}+b_{2}+d+e_{1}+e_{2} &\equiv &1 \\
l_{1} &=&u_{1} \\
l_{2} &=&u_{2}
\end{eqnarray*}%
and we get%
\begin{eqnarray*}
&&\sum_{\substack{ a,b_{1},b_{2},d,e_{1},e_{2}=0  \\ %
a+b_{1}+b_{2}+d+e_{1}+e_{2}\equiv 1}}^{1}\left( -1\right) ^{\alpha \left(
g;0,0,0,0\right) }B(gx_{1}x_{2}\otimes
1_{H};G^{a}X_{1}^{b_{1}}X_{2}^{b_{2}},g^{d}x_{1}^{e_{1}}x_{2}^{e_{2}})G^{a}X_{1}^{b_{1}}X_{2}^{b_{2}}\otimes g^{d}x_{1}^{e_{1}}x_{2}^{e_{2}}
\\
&=&B^{A}(gx_{1}x_{2}\otimes 1_{H})\otimes B^{H}(gx_{1}x_{2}\otimes 1_{H})
\end{eqnarray*}%
so that we get%
\begin{equation}
B(gx_{1}x_{2}\otimes
1_{H};G^{a}X_{1}^{b_{1}}X_{2}^{b_{2}},g^{d}x_{1}^{e_{1}}x_{2}^{e_{2}})=0%
\text{ whenever }a+b_{1}+b_{2}+d+e_{1}+e_{2}\equiv 0  \label{gx1x2ot1, first}
\end{equation}

\subsection{Case $g$}

We have to consider only the last summand.%
\begin{eqnarray*}
&&+\sum_{a,b_{1},b_{2},d,e_{1},e_{2}=0}^{1}\sum_{l_{1}=0}^{b_{1}}%
\sum_{l_{2}=0}^{b_{2}}\sum_{u_{1}=0}^{e_{1}}\sum_{u_{2}=0}^{e_{2}}\left(
-1\right) ^{\alpha \left( g;l_{1},l_{2},u_{1},u_{2}\right) } \\
&&B(gx_{1}x_{2}\otimes
1_{H};G^{a}X_{1}^{b_{1}}X_{2}^{b_{2}},g^{d}x_{1}^{e_{1}}x_{2}^{e_{2}}) \\
&&G^{a}X_{1}^{b_{1}-l_{1}}X_{2}^{b_{2}-l_{2}}\otimes
g^{d}x_{1}^{e_{1}-u_{1}}x_{2}^{e_{2}-u_{2}}\otimes
g^{a+b_{1}+b_{2}+l_{1}+l_{2}+d+e_{1}+e_{2}+u_{1}+u_{2}+1}x_{1}^{l_{1}+u_{1}}x_{2}^{l_{2}+u_{2}}
\end{eqnarray*}%
and we get%
\begin{eqnarray*}
a+b_{1}+b_{2}+l_{1}+l_{2}+d+e_{1}+e_{2}+u_{1}+u_{2}+1 &\equiv &1 \\
l_{1}+u_{1} &=&0 \\
l_{2}+u_{2} &=&0
\end{eqnarray*}%
i.e.%
\begin{eqnarray*}
a+b_{1}+b_{2}+d+e_{1}+e_{2} &\equiv &0 \\
l_{1} &=&u_{1} \\
l_{2} &=&u_{2}
\end{eqnarray*}%
By $\left( \ref{gx1x2ot1, first}\right) $ we know that $B(gx_{1}x_{2}\otimes
1_{H};G^{a}X_{1}^{b_{1}}X_{2}^{b_{2}},g^{d}x_{1}^{e_{1}}x_{2}^{e_{2}})=0$
whenever $a+b_{1}+b_{2}+d+e_{1}+e_{2}\equiv 0$. Thus the equality is
satisfied.

\subsection{Case $x_{1}$}

We have to consider only the second and the fourth one. Let us consider the
second one%
\begin{eqnarray*}
&&\sum_{a,b_{1},b_{2},d,e_{1},e_{2}=0}^{1}\sum_{l_{1}=0}^{b_{1}}%
\sum_{l_{2}=0}^{b_{2}}\sum_{u_{1}=0}^{e_{1}}\sum_{u_{2}=0}^{e_{2}}\left(
-1\right) ^{\alpha \left( gx_{1};l_{1},l_{2},u_{1},u_{2}\right)
}B(x_{2}\otimes
1_{H};G^{a}X_{1}^{b_{1}}X_{2}^{b_{2}},g^{d}x_{1}^{e_{1}}x_{2}^{e_{2}}) \\
&&G^{a}X_{1}^{b_{1}-l_{1}}X_{2}^{b_{2}-l_{2}}\otimes
g^{d}x_{1}^{e_{1}-u_{1}}x_{2}^{e_{2}-u_{2}}\otimes
g^{a+b_{1}+b_{2}+l_{1}+l_{2}+d+e_{1}+e_{2}+u_{1}+u_{2}+1}x_{1}^{l_{1}+u_{1}+1}x_{2}^{l_{2}+u_{2}}
\end{eqnarray*}%
\begin{eqnarray*}
a+b_{1}+b_{2}+l_{1}+l_{2}+d+e_{1}+e_{2}+u_{1}+u_{2}+1 &\equiv &0 \\
l_{1}+u_{1} &=&0 \\
l_{2}+u_{2} &=&0
\end{eqnarray*}%
i.e.%
\begin{eqnarray*}
a+b_{1}+b_{2}+d+e_{1}+e_{2} &\equiv &1 \\
l_{1} &=&u_{1} \\
l_{2} &=&u_{2}
\end{eqnarray*}%
and we get%
\begin{equation*}
\sum_{\substack{ a,b_{1},b_{2},d,e_{1},e_{2}=0  \\ %
a+b_{1}+b_{2}+d+e_{1}+e_{2}\equiv 1}}^{1}\left( -1\right) ^{\alpha \left(
gx_{1};0,0,0,0\right) }B(x_{2}\otimes
1_{H};G^{a}X_{1}^{b_{1}}X_{2}^{b_{2}},g^{d}x_{1}^{e_{1}}x_{2}^{e_{2}})G^{a}X_{1}^{b_{1}}X_{2}^{b_{2}}\otimes g^{d}x_{1}^{e_{1}}x_{2}^{e_{2}}.
\end{equation*}%
Since, by \ref{xot2 first} $B(x_{2}\otimes
1_{H};G^{a}X_{1}^{b_{1}}X_{2}^{b_{2}},g^{d}x_{1}^{e_{1}}x_{2}^{e_{2}})=0$
whenever $a+b_{1}+b_{2}+d+e_{1}+e_{2}$ $\equiv 1$, this summand is always
zero. Let us consider the fourth summand.%
\begin{eqnarray*}
&&+\sum_{a,b_{1},b_{2},d,e_{1},e_{2}=0}^{1}\sum_{l_{1}=0}^{b_{1}}%
\sum_{l_{2}=0}^{b_{2}}\sum_{u_{1}=0}^{e_{1}}\sum_{u_{2}=0}^{e_{2}}\left(
-1\right) ^{\alpha \left( g;l_{1},l_{2},u_{1},u_{2}\right) } \\
&&B(gx_{1}x_{2}\otimes
1_{H};G^{a}X_{1}^{b_{1}}X_{2}^{b_{2}},g^{d}x_{1}^{e_{1}}x_{2}^{e_{2}}) \\
&&G^{a}X_{1}^{b_{1}-l_{1}}X_{2}^{b_{2}-l_{2}}\otimes
g^{d}x_{1}^{e_{1}-u_{1}}x_{2}^{e_{2}-u_{2}}\otimes
g^{a+b_{1}+b_{2}+l_{1}+l_{2}+d+e_{1}+e_{2}+u_{1}+u_{2}+1}x_{1}^{l_{1}+u_{1}}x_{2}^{l_{2}+u_{2}}
\end{eqnarray*}%
\begin{eqnarray*}
a+b_{1}+b_{2}+l_{1}+l_{2}+d+e_{1}+e_{2}+u_{1}+u_{2}+1 &\equiv &0 \\
l_{1}+u_{1} &=&1 \\
l_{2}+u_{2} &=&0
\end{eqnarray*}%
i.e.%
\begin{eqnarray*}
a+b_{1}+b_{2}+d+e_{1}+e_{2} &\equiv &0 \\
l_{1}+u_{1} &=&1 \\
l_{2} &=&u_{2}
\end{eqnarray*}%
By $\left( \ref{gx1x2ot1, first}\right) $ we know that $B(gx_{1}x_{2}\otimes
1_{H};G^{a}X_{1}^{b_{1}}X_{2}^{b_{2}},g^{d}x_{1}^{e_{1}}x_{2}^{e_{2}})=0$
whenever $a+b_{1}+b_{2}+d+e_{1}+e_{2}\equiv 0,$ thus this summand is zero.
In conclusion, since $\alpha \left( gx_{1};0,0,0,0\right) =a+b_{1}+b_{2},$
we just get%
\begin{equation*}
\sum_{\substack{ a,b_{1},b_{2},d,e_{1},e_{2}=0  \\ %
a+b_{1}+b_{2}+d+e_{1}+e_{2}\equiv 1}}^{1}\left( -1\right)
^{a+b_{1}+b_{2}}B(x_{2}\otimes
1_{H};G^{a}X_{1}^{b_{1}}X_{2}^{b_{2}},g^{d}x_{1}^{e_{1}}x_{2}^{e_{2}})G^{a}X_{1}^{b_{1}}X_{2}^{b_{2}}\otimes g^{d}x_{1}^{e_{1}}x_{2}^{e_{2}}=0
\end{equation*}

\subsection{Case $x_{2}$}

We have to consider only the third and the fourth summand. Let us examine
the third summand.

\begin{eqnarray*}
&&-\sum_{a,b_{1},b_{2},d,e_{1},e_{2}=0}^{1}\sum_{l_{1}=0}^{b_{1}}%
\sum_{l_{2}=0}^{b_{2}}\sum_{u_{1}=0}^{e_{1}}\sum_{u_{2}=0}^{e_{2}}\left(
-1\right) ^{\alpha \left( gx_{2};l_{1},l_{2},u_{1},u_{2}\right) } \\
&&B(x_{1}\otimes
1_{H};G^{a}X_{1}^{b_{1}}X_{2}^{b_{2}},g^{d}x_{1}^{e_{1}}x_{2}^{e_{2}}) \\
&&G^{a}X_{1}^{b_{1}-l_{1}}X_{2}^{b_{2}-l_{2}}\otimes
g^{d}x_{1}^{e_{1}-u_{1}}x_{2}^{e_{2}-u_{2}}\otimes
g^{a+b_{1}+b_{2}+l_{1}+l_{2}+d+e_{1}+e_{2}+u_{1}+u_{2}+1}x_{1}^{l_{1}+u_{1}}x_{2}^{l_{2}+u_{2}+1}
\end{eqnarray*}%
We get%
\begin{eqnarray*}
a+b_{1}+b_{2}+l_{1}+l_{2}+d+e_{1}+e_{2}+u_{1}+u_{2}+1 &\equiv &0 \\
l_{1}+u_{1} &=&0 \\
l_{2}+u_{2}+1 &=&1
\end{eqnarray*}%
Thus we deduce that $a+b_{1}+b_{2}+d+e_{1}+e_{2}\equiv 1.$ Since, by $\left( %
\ref{x1ot1 first}\right) ,$ $B(x_{1}\otimes
1_{H};G^{a}X_{1}^{b_{1}}X_{2}^{b_{2}},g^{d}x_{1}^{e_{1}}x_{2}^{e_{2}})=0$
whenever $a+b_{1}+b_{2}+d+e_{1}+e_{2}\equiv 1,$ this summand is always zero.

Let us consider the last summand.

\begin{eqnarray*}
&&+\sum_{a,b_{1},b_{2},d,e_{1},e_{2}=0}^{1}\sum_{l_{1}=0}^{b_{1}}%
\sum_{l_{2}=0}^{b_{2}}\sum_{u_{1}=0}^{e_{1}}\sum_{u_{2}=0}^{e_{2}}\left(
-1\right) ^{\alpha \left( g;l_{1},l_{2},u_{1},u_{2}\right) } \\
&&B(gx_{1}x_{2}\otimes
1_{H};G^{a}X_{1}^{b_{1}}X_{2}^{b_{2}},g^{d}x_{1}^{e_{1}}x_{2}^{e_{2}}) \\
&&G^{a}X_{1}^{b_{1}-l_{1}}X_{2}^{b_{2}-l_{2}}\otimes
g^{d}x_{1}^{e_{1}-u_{1}}x_{2}^{e_{2}-u_{2}}\otimes
g^{a+b_{1}+b_{2}+l_{1}+l_{2}+d+e_{1}+e_{2}+u_{1}+u_{2}+1}x_{1}^{l_{1}+u_{1}}x_{2}^{l_{2}+u_{2}}
\end{eqnarray*}%
\begin{eqnarray*}
a+b_{1}+b_{2}+l_{1}+l_{2}+d+e_{1}+e_{2}+u_{1}+u_{2}+1 &\equiv &0 \\
l_{1}+u_{1} &=&0 \\
l_{2}+u_{2} &=&1
\end{eqnarray*}%
i.e.%
\begin{eqnarray*}
a+b_{1}+b_{2}+d+e_{1}+e_{2} &\equiv &0 \\
l_{1}+u_{1} &=&1 \\
l_{2} &=&u_{2}
\end{eqnarray*}%
By $\left( \ref{gx1x2ot1, first}\right) $ we know that $B(gx_{1}x_{2}\otimes
1_{H};G^{a}X_{1}^{b_{1}}X_{2}^{b_{2}},g^{d}x_{1}^{e_{1}}x_{2}^{e_{2}})=0$
whenever $a+b_{1}+b_{2}+d+e_{1}+e_{2}\equiv 0,$ thus this summand is zero.

\subsection{Case $x_{1}x_{2}$}

Let us consider the first summand%
\begin{eqnarray*}
&&\sum_{a,b_{1},b_{2},d,e_{1},e_{2}=0}^{1}\sum_{l_{1}=0}^{b_{1}}%
\sum_{l_{2}=0}^{b_{2}}\sum_{u_{1}=0}^{e_{1}}\sum_{u_{2}=0}^{e_{2}}\left(
-1\right) ^{\alpha \left( gx_{1}x_{2};l_{1},l_{2},u_{1},u_{2}\right) } \\
&&B(g\otimes
1_{H};G^{a}X_{1}^{b_{1}}X_{2}^{b_{2}},g^{d}x_{1}^{e_{1}}x_{2}^{e_{2}}) \\
&&G^{a}X_{1}^{b_{1}-l_{1}}X_{2}^{b_{2}-l_{2}}\otimes
g^{d}x_{1}^{e_{1}-u_{1}}x_{2}^{e_{2}-u_{2}}\otimes
g^{a+b_{1}+b_{2}+l_{1}+l_{2}+d+e_{1}+e_{2}+u_{1}+u_{2}+1}x_{1}^{l_{1}+u_{1}+1}x_{2}^{l_{2}+u_{2}+1}
\end{eqnarray*}%
We get%
\begin{eqnarray*}
a+b_{1}+b_{2}+l_{1}+l_{2}+d+e_{1}+e_{2}+u_{1}+u_{2}+1 &\equiv &0 \\
l_{1}+u_{1}+1 &=&1 \\
l_{2}+u_{2}+1 &=&1
\end{eqnarray*}%
and hence we deduce that%
\begin{eqnarray*}
a+b_{1}+b_{2}+d+e_{1}+e_{2} &\equiv &1 \\
l_{1}+u_{1} &=&0 \\
l_{2}+u_{2} &=&0
\end{eqnarray*}%
and, since $\alpha \left( gx_{1}x_{2};0,0,0,0\right) =0,$ we get%
\begin{equation*}
\sum_{\substack{ a,b_{1},b_{2},d,e_{1},e_{2}=0  \\ %
a+b_{1}+b_{2}+d+e_{1}+e_{2}\equiv 1}}^{1}B(g\otimes
1_{H};G^{a}X_{1}^{b_{1}}X_{2}^{b_{2}},g^{d}x_{1}^{e_{1}}x_{2}^{e_{2}})G^{a}X_{1}^{b_{1}}X_{2}^{b_{2}}\otimes g^{d}x_{1}^{e_{1}}x_{2}^{e_{2}}.
\end{equation*}%
Let us consider the second summand.%
\begin{eqnarray*}
&&\sum_{a,b_{1},b_{2},d,e_{1},e_{2}=0}^{1}\sum_{l_{1}=0}^{b_{1}}%
\sum_{l_{2}=0}^{b_{2}}\sum_{u_{1}=0}^{e_{1}}\sum_{u_{2}=0}^{e_{2}}\left(
-1\right) ^{\alpha \left( gx_{1};l_{1},l_{2},u_{1},u_{2}\right) } \\
&&B(x_{2}\otimes
1_{H};G^{a}X_{1}^{b_{1}}X_{2}^{b_{2}},g^{d}x_{1}^{e_{1}}x_{2}^{e_{2}}) \\
&&G^{a}X_{1}^{b_{1}-l_{1}}X_{2}^{b_{2}-l_{2}}\otimes
g^{d}x_{1}^{e_{1}-u_{1}}x_{2}^{e_{2}-u_{2}}\otimes
g^{a+b_{1}+b_{2}+l_{1}+l_{2}+d+e_{1}+e_{2}+u_{1}+u_{2}+1}x_{1}^{l_{1}+u_{1}+1}x_{2}^{l_{2}+u_{2}}
\end{eqnarray*}%
We get%
\begin{eqnarray*}
a+b_{1}+b_{2}+l_{1}+l_{2}+d+e_{1}+e_{2}+u_{1}+u_{2}+1 &\equiv &0 \\
l_{1}+u_{1}+1 &=&1 \\
l_{2}+u_{2} &=&1
\end{eqnarray*}%
and hence we deduce that%
\begin{eqnarray*}
a+b_{1}+b_{2}+d+e_{1}+e_{2} &\equiv &0 \\
l_{1}+u_{1} &=&0 \\
l_{2}+u_{2} &=&1
\end{eqnarray*}%
and we get%
\begin{eqnarray*}
&&\sum_{\substack{ a,b_{1},b_{2},d,e_{1},e_{2}=0  \\ %
a+b_{1}+b_{2}+d+e_{1}+e_{2}\equiv 0}}^{1}\sum_{l_{2}=0}^{b_{2}}\sum
_{\substack{ u_{2}=0  \\ l_{2}+u_{2}=1}}^{e_{2}}\left( -1\right) ^{\alpha
\left( gx_{1};0,l_{2},0,u_{2}\right) } \\
&&B(x_{2}\otimes
1_{H};G^{a}X_{1}^{b_{1}}X_{2}^{b_{2}},g^{d}x_{1}^{e_{1}}x_{2}^{e_{2}})G^{a}X_{1}^{b_{1}}X_{2}^{b_{2}-l_{2}}\otimes g^{d}x_{1}^{e_{1}}x_{2}^{e_{2}-u_{2}}.
\end{eqnarray*}%
Since $\alpha \left( gx_{1};0,0,0,1\right) =0$ and $\alpha \left(
gx_{1};0,1,0,0\right) =a+b_{1}+b_{2}+1$ we obtain%
\begin{eqnarray*}
&&\sum_{\substack{ a,b_{1},b_{2},d,e_{1}=0  \\ a+b_{1}+b_{2}+d+e_{1}\equiv 1
}}^{1}B(x_{2}\otimes
1_{H};G^{a}X_{1}^{b_{1}}X_{2}^{b_{2}},g^{d}x_{1}^{e_{1}}x_{2})G^{a}X_{1}^{b_{1}}X_{2}^{b_{2}}\otimes g^{d}x_{1}^{e_{1}}+
\\
&&+\sum_{\substack{ a,b_{1},d,e_{1},e_{2}=0  \\ a+b_{1}+d+e_{1}+e_{2}\equiv
1 }}^{1}\left( -1\right) ^{a+b_{1}}B(x_{2}\otimes
1_{H};G^{a}X_{1}^{b_{1}}X_{2},g^{d}x_{1}^{e_{1}}x_{2}^{e_{2}})G^{a}X_{1}^{b_{1}}\otimes g^{d}x_{1}^{e_{1}}x_{2}^{e_{2}}.
\end{eqnarray*}%
Let us consider the third summand.

\begin{eqnarray*}
&&-\sum_{a,b_{1},b_{2},d,e_{1},e_{2}=0}^{1}\sum_{l_{1}=0}^{b_{1}}%
\sum_{l_{2}=0}^{b_{2}}\sum_{u_{1}=0}^{e_{1}}\sum_{u_{2}=0}^{e_{2}}\left(
-1\right) ^{\alpha \left( gx_{2};l_{1},l_{2},u_{1},u_{2}\right) } \\
&&B(x_{1}\otimes
1_{H};G^{a}X_{1}^{b_{1}}X_{2}^{b_{2}},g^{d}x_{1}^{e_{1}}x_{2}^{e_{2}}) \\
&&G^{a}X_{1}^{b_{1}-l_{1}}X_{2}^{b_{2}-l_{2}}\otimes
g^{d}x_{1}^{e_{1}-u_{1}}x_{2}^{e_{2}-u_{2}}\otimes
g^{a+b_{1}+b_{2}+l_{1}+l_{2}+d+e_{1}+e_{2}+u_{1}+u_{2}+1}x_{1}^{l_{1}+u_{1}}x_{2}^{l_{2}+u_{2}+1}
\end{eqnarray*}%
We get%
\begin{eqnarray*}
a+b_{1}+b_{2}+l_{1}+l_{2}+d+e_{1}+e_{2}+u_{1}+u_{2}+1 &\equiv &0 \\
l_{1}+u_{1} &=&1 \\
l_{2}+u_{2}+1 &=&1
\end{eqnarray*}%
and hence we deduce that%
\begin{eqnarray*}
a+b_{1}+b_{2}+d+e_{1}+e_{2} &\equiv &0 \\
l_{1}+u_{1} &=&1 \\
l_{2}+u_{2} &=&0
\end{eqnarray*}%
and we obtain%
\begin{eqnarray*}
&&\sum_{\substack{ a,b_{1},b_{2},d,e_{1},e_{2}=0  \\ %
a+b_{1}+b_{2}+d+e_{1}+e_{2}\equiv 0}}^{1}\sum_{l_{1}=0}^{b_{1}}\sum
_{\substack{ u_{1}=0  \\ l_{1}+u_{1}=1}}^{e_{1}}\left( -1\right) ^{\alpha
\left( gx_{2};l_{1},0,u_{1},0\right) +1}B(x_{1}\otimes
1_{H};G^{a}X_{1}^{b_{1}}X_{2}^{b_{2}},g^{d}x_{1}^{e_{1}}x_{2}^{e_{2}}) \\
&&G^{a}X_{1}^{b_{1}-l_{1}}X_{2}^{b_{2}}\otimes
g^{d}x_{1}^{e_{1}-u_{1}}x_{2}^{e_{2}}.
\end{eqnarray*}%
Since $\alpha \left( gx_{2};0,0,1,0\right) =e_{2}+1$ and $\alpha \left(
gx_{2};1,0,0,0\right) =a+b_{1},$ we get
\begin{eqnarray*}
&&\sum_{\substack{ a,b_{1},b_{2},d,e_{2}=0  \\ a+b_{1}+b_{2}+d+e_{2}\equiv 1
}}^{1}\left( -1\right) ^{e_{2}}B(x_{1}\otimes
1_{H};G^{a}X_{1}^{b_{1}}X_{2}^{b_{2}},g^{d}x_{1}x_{2}^{e_{2}})G^{a}X_{1}^{b_{1}}X_{2}^{b_{2}}\otimes g^{d}x_{2}^{e_{2}}+
\\
&&+\sum_{\substack{ a,b_{2},d,e_{1},e_{2}=0  \\ a+b_{2}+d+e_{1}+e_{2}\equiv
1 }}^{1}\left( -1\right) ^{a}B(x_{1}\otimes
1_{H};G^{a}X_{1}X_{2}^{b_{2}},g^{d}x_{1}^{e_{1}}x_{2}^{e_{2}})G^{a}X_{2}^{b_{2}}\otimes g^{d}x_{1}^{e_{1}}x_{2}^{e_{2}}.
\end{eqnarray*}%
Let us consider the fourth summand.%
\begin{eqnarray*}
&&+\sum_{a,b_{1},b_{2},d,e_{1},e_{2}=0}^{1}\sum_{l_{1}=0}^{b_{1}}%
\sum_{l_{2}=0}^{b_{2}}\sum_{u_{1}=0}^{e_{1}}\sum_{u_{2}=0}^{e_{2}}\left(
-1\right) ^{\alpha \left( g;l_{1},l_{2},u_{1},u_{2}\right) } \\
&&B(gx_{1}x_{2}\otimes
1_{H};G^{a}X_{1}^{b_{1}}X_{2}^{b_{2}},g^{d}x_{1}^{e_{1}}x_{2}^{e_{2}}) \\
&&G^{a}X_{1}^{b_{1}-l_{1}}X_{2}^{b_{2}-l_{2}}\otimes
g^{d}x_{1}^{e_{1}-u_{1}}x_{2}^{e_{2}-u_{2}}\otimes
g^{a+b_{1}+b_{2}+l_{1}+l_{2}+d+e_{1}+e_{2}+u_{1}+u_{2}+1}x_{1}^{l_{1}+u_{1}}x_{2}^{l_{2}+u_{2}}
\end{eqnarray*}

We get%
\begin{eqnarray*}
a+b_{1}+b_{2}+l_{1}+l_{2}+d+e_{1}+e_{2}+u_{1}+u_{2}+1 &\equiv &0 \\
l_{1}+u_{1} &=&1 \\
l_{2}+u_{2} &=&1
\end{eqnarray*}%
and hence we deduce that%
\begin{eqnarray*}
a+b_{1}+b_{2}+d+e_{1}+e_{2} &\equiv &1 \\
l_{1}+u_{1} &=&1 \\
l_{2}+u_{2} &=&1
\end{eqnarray*}%
Since, $\alpha \left( g;0,0,1,1\right) =1+e_{2},$ $\alpha \left(
g;0,1,1,0\right) =e_{2}+a+b_{1}+b_{2},$ $\alpha \left( g;1,0,0,1\right)
=a+b_{1}+1$ and $\alpha \left( g;1,1,0,0\right) =1+b_{2},$ we obtain%
\begin{eqnarray*}
&&\sum_{\substack{ a,b_{1},b_{2},d,e_{1},e_{2}=0  \\ %
a+b_{1}+b_{2}+d+e_{1}+e_{2}\equiv 1}}^{1}B(g\otimes
1_{H};G^{a}X_{1}^{b_{1}}X_{2}^{b_{2}},g^{d}x_{1}^{e_{1}}x_{2}^{e_{2}})G^{a}X_{1}^{b_{1}}X_{2}^{b_{2}}\otimes g^{d}x_{1}^{e_{1}}x_{2}^{e_{2}}
\\
&&\sum_{\substack{ a,b_{1},b_{2},d,e_{1}=0  \\ a+b_{1}+b_{2}+d+e_{1}\equiv 1
}}^{1}B(x_{2}\otimes
1_{H};G^{a}X_{1}^{b_{1}}X_{2}^{b_{2}},g^{d}x_{1}^{e_{1}}x_{2})G^{a}X_{1}^{b_{1}}X_{2}^{b_{2}}\otimes g^{d}x_{1}^{e_{1}}+
\\
&&+\sum_{\substack{ a,b_{1},d,e_{1},e_{2}=0  \\ a+b_{1}+d+e_{1}+e_{2}\equiv
1 }}^{1}\left( -1\right) ^{a+b_{1}}B(x_{2}\otimes
1_{H};G^{a}X_{1}^{b_{1}}X_{2},g^{d}x_{1}^{e_{1}}x_{2}^{e_{2}})G^{a}X_{1}^{b_{1}}\otimes g^{d}x_{1}^{e_{1}}x_{2}^{e_{2}}
\\
&&\sum_{\substack{ a,b_{1},b_{2},d,e_{2}=0  \\ a+b_{1}+b_{2}+d+e_{2}\equiv 1
}}^{1}\left( -1\right) ^{e_{2}}B(x_{1}\otimes
1_{H};G^{a}X_{1}^{b_{1}}X_{2}^{b_{2}},g^{d}x_{1}x_{2}^{e_{2}})G^{a}X_{1}^{b_{1}}X_{2}^{b_{2}}\otimes g^{d}x_{2}^{e_{2}}+
\\
&&+\sum_{\substack{ a,b_{2},d,e_{1},e_{2}=0  \\ a+b_{2}+d+e_{1}+e_{2}\equiv
1 }}^{1}\left( -1\right) ^{a}B(x_{1}\otimes
1_{H};G^{a}X_{1}X_{2}^{b_{2}},g^{d}x_{1}^{e_{1}}x_{2}^{e_{2}})G^{a}X_{2}^{b_{2}}\otimes g^{d}x_{1}^{e_{1}}x_{2}^{e_{2}}+
\\
&&\sum_{\substack{ a,b_{1},b_{2},d=0  \\ a+b_{1}+b_{2}+d\equiv 1}}%
^{1}B(gx_{1}x_{2}\otimes
1_{H};G^{a}X_{1}^{b_{1}}X_{2}^{b_{2}},g^{d}x_{1}x_{2})G^{a}X_{1}^{b_{1}}X_{2}^{b_{2}}\otimes g^{d}+
\\
&&+\sum_{\substack{ a,b_{1},d,e_{2}=0  \\ a+b_{1}+d+e_{2}\equiv 1}}%
^{1}\left( -1\right) ^{e_{2}+a+b_{1}+1}B(gx_{1}x_{2}\otimes
1_{H};G^{a}X_{1}^{b_{1}}X_{2},g^{d}x_{1}x_{2}^{e_{2}})G^{a}X_{1}^{b_{1}}%
\otimes g^{d}x_{2}^{e_{2}}+ \\
&&+\sum_{\substack{ a,b_{2},d,e_{1}=0  \\ a+b_{2}+d+e_{1}\equiv 1}}%
^{1}\left( -1\right) ^{a}B(gx_{1}x_{2}\otimes
1_{H};G^{a}X_{1}X_{2}^{b_{2}},g^{d}x_{1}^{e_{1}}x_{2})G^{a}X_{2}^{b_{2}}%
\otimes g^{d}x_{1}^{e_{1}}+ \\
&&+\sum_{\substack{ a,d,e_{1},e_{2}=0  \\ a+d+e_{1}+e_{2}\equiv 1}}%
^{1}B(gx_{1}x_{2}\otimes
1_{H};G^{a}X_{1}X_{2},g^{d}x_{1}^{e_{1}}x_{2}^{e_{2}})G^{a}\otimes
g^{d}x_{1}^{e_{1}}x_{2}^{e_{2}}.
\end{eqnarray*}%
In conclusion we get

\begin{eqnarray*}
&&\sum_{\substack{ a,b_{1},b_{2},d=0  \\ a+b_{1}+b_{2}+d\equiv 1}}%
^{1}B(gx_{1}x_{2}\otimes
1_{H};G^{a}X_{1}^{b_{1}}X_{2}^{b_{2}},g^{d}x_{1}x_{2})G^{a}X_{1}^{b_{1}}X_{2}^{b_{2}}\otimes g^{d}+
\\
&&+\sum_{\substack{ a,b_{1},d,e_{2}=0  \\ a+b_{1}+d+e_{2}\equiv 1}}%
^{1}\left( -1\right) ^{e_{2}+a+b_{1}+1}B(gx_{1}x_{2}\otimes
1_{H};G^{a}X_{1}^{b_{1}}X_{2},g^{d}x_{1}x_{2}^{e_{2}})G^{a}X_{1}^{b_{1}}%
\otimes g^{d}x_{2}^{e_{2}}+ \\
&&+\sum_{\substack{ a,b_{2},d,e_{1}=0  \\ a+b_{2}+d+e_{1}\equiv 1}}%
^{1}\left( -1\right) ^{a}B(gx_{1}x_{2}\otimes
1_{H};G^{a}X_{1}X_{2}^{b_{2}},g^{d}x_{1}^{e_{1}}x_{2})G^{a}X_{2}^{b_{2}}%
\otimes g^{d}x_{1}^{e_{1}}+ \\
&&+\sum_{\substack{ a,d,e_{1},e_{2}=0  \\ a+d+e_{1}+e_{2}\equiv 1}}%
^{1}B(gx_{1}x_{2}\otimes
1_{H};G^{a}X_{1}X_{2},g^{d}x_{1}^{e_{1}}x_{2}^{e_{2}})G^{a}\otimes
g^{d}x_{1}^{e_{1}}x_{2}^{e_{2}}.
\end{eqnarray*}

\subsubsection{$G^{a}\otimes g^{d}$}

\begin{equation*}
\sum_{\substack{ a,d=0  \\ a+d\equiv 1}}^{1}\left[
\begin{array}{c}
B(g\otimes 1_{H};G^{a},g^{d})+B(x_{2}\otimes 1_{H};G^{a},g^{d}x_{2})+\left(
-1\right) ^{a}B(x_{2}\otimes 1_{H};G^{a}X_{2},g^{d})+ \\
+B(x_{1}\otimes 1_{H};G^{a},g^{d}x_{1})+\left( -1\right) ^{a}B(x_{1}\otimes
1_{H};G^{a}X_{1},g^{d}) \\
+B(gx_{1}x_{2}\otimes 1_{H};G^{a},g^{d}x_{1}x_{2})+\left( -1\right)
^{a+1}B(gx_{1}x_{2}\otimes 1_{H};G^{a}X_{2},g^{d}x_{1}) \\
+\left( -1\right) ^{a}B(gx_{1}x_{2}\otimes
1_{H};G^{a}X_{1},g^{d}x_{2})+B(gx_{1}x_{2}\otimes
1_{H};G^{a}X_{1}X_{2},g^{d})%
\end{array}%
\right] G^{a}\otimes g^{d}=0
\end{equation*}

and we get%
\begin{equation*}
\begin{array}{c}
B(g\otimes 1_{H};1_{A},g)+B(x_{2}\otimes 1_{H};1_{A},gx_{2})+B(x_{2}\otimes
1_{H};X_{2},g)+ \\
+B(x_{1}\otimes 1_{H};1_{A},gx_{1})+B(x_{1}\otimes
1_{H};X_{1},g)+B(gx_{1}x_{2}\otimes 1_{H};1_{A},gx_{1}x_{2})+ \\
-B(gx_{1}x_{2}\otimes 1_{H};X_{2},gx_{1})+B(gx_{1}x_{2}\otimes
1_{H};X_{1},gx_{2})+B(gx_{1}x_{2}\otimes 1_{H};X_{1}X_{2},g)%
\end{array}%
=0
\end{equation*}%
and%
\begin{equation*}
\begin{array}{c}
B(g\otimes 1_{H};G,1_{H})+B(x_{2}\otimes 1_{H};G,x_{2})-B(x_{2}\otimes
1_{H};GX_{2},1_{H})+ \\
+B(x_{1}\otimes 1_{H};G,x_{1})-B(x_{1}\otimes
1_{H};GX_{1},1_{H})+B(gx_{1}x_{2}\otimes 1_{H};G,x_{1}x_{2})+ \\
+B(gx_{1}x_{2}\otimes 1_{H};GX_{2},x_{1})-B(gx_{1}x_{2}\otimes
1_{H};GX_{1},x_{2})+B(gx_{1}x_{2}\otimes 1_{H};GX_{1}X_{2},1_{H})%
\end{array}%
=0
\end{equation*}%
By taking in account the form of the elements $B(g\otimes 1_{H})$, $%
B(x_{1}\otimes 1_{H})$ and $B(x_{2}\otimes \ 1_{H}),$ we get
\begin{equation*}
\begin{array}{c}
B(g\otimes 1_{H};1_{A},g)+B(x_{2}\otimes 1_{H};1_{A},gx_{2})+\left[
-B(g\otimes 1_{H};1_{A},g)-B(x_{2}\otimes \ 1_{H};1_{A},gx_{2})\right] + \\
+B(x_{1}\otimes 1_{H};1_{A},gx_{1})+\left[ -B(g\otimes
1_{H};1_{A},g)-B(x_{1}\otimes 1_{H};1_{A},gx_{1})\right] +B(gx_{1}x_{2}%
\otimes 1_{H};1_{A},gx_{1}x_{2})+ \\
-B(gx_{1}x_{2}\otimes 1_{H};X_{2},gx_{1})+B(gx_{1}x_{2}\otimes
1_{H};X_{1},gx_{2})+B(gx_{1}x_{2}\otimes 1_{H};X_{1}X_{2},g)%
\end{array}%
=0
\end{equation*}%
i.e.%
\begin{equation}
\begin{array}{c}
-B(g\otimes 1_{H};1_{A},g))+B(gx_{1}x_{2}\otimes 1_{H};1_{A},gx_{1}x_{2})+
\\
-B(gx_{1}x_{2}\otimes 1_{H};X_{2},gx_{1})+B(gx_{1}x_{2}\otimes
1_{H};X_{1},gx_{2})+B(gx_{1}x_{2}\otimes 1_{H};X_{1}X_{2},g)%
\end{array}%
=0  \label{gx1x2ot1, second}
\end{equation}%
and%
\begin{equation*}
\begin{array}{c}
B(g\otimes 1_{H};G,1_{H})+B(x_{2}\otimes 1_{H};G,x_{2})-\left[ B(g\otimes
1_{H};G,1_{H})+B(x_{2}\otimes \ 1_{H};G,x_{2}\right] + \\
+B(x_{1}\otimes 1_{H};G,x_{1})-\left[ B(g\otimes
1_{H};G,1_{H})+B(x_{1}\otimes 1_{H};G,x_{1})\right] +B(gx_{1}x_{2}\otimes
1_{H};G,x_{1}x_{2})+ \\
+B(gx_{1}x_{2}\otimes 1_{H};GX_{2},x_{1})-B(gx_{1}x_{2}\otimes
1_{H};GX_{1},x_{2})+B(gx_{1}x_{2}\otimes 1_{H};GX_{1}X_{2},1_{H})%
\end{array}%
=0
\end{equation*}%
i.e.%
\begin{equation}
\begin{array}{c}
-B(g\otimes 1_{H};G,1_{H})+B(gx_{1}x_{2}\otimes 1_{H};G,x_{1}x_{2})+ \\
+B(gx_{1}x_{2}\otimes 1_{H};GX_{2},x_{1})-B(gx_{1}x_{2}\otimes
1_{H};GX_{1},x_{2})+B(gx_{1}x_{2}\otimes 1_{H};GX_{1}X_{2},1_{H})%
\end{array}%
=0  \label{gx1x2ot1, third}
\end{equation}

\subsubsection{$G^{a}\otimes g^{d}x_{2}$}

\begin{equation*}
\sum_{\substack{ a,d=0  \\ a+d\equiv 0}}^{1}\left[
\begin{array}{c}
B(g\otimes 1_{H};G^{a},g^{d}x_{2})+\left( -1\right) ^{a}B(x_{2}\otimes
1_{H};G^{a}X_{2},g^{d}x_{2}) \\
-B(x_{1}\otimes 1_{H};G^{a},g^{d}x_{1}x_{2})+\left( -1\right)
^{a}B(x_{1}\otimes 1_{H};G^{a}X_{1},g^{d}x_{2})+ \\
+\left( -1\right) ^{a}B(gx_{1}x_{2}\otimes
1_{H};G^{a}X_{2},g^{d}x_{1}x_{2})+B(gx_{1}x_{2}\otimes
1_{H};G^{a}X_{1}X_{2},g^{d}x_{2})%
\end{array}%
\right] G^{a}\otimes g^{d}x_{2}=0
\end{equation*}%
and we get%
\begin{equation*}
\begin{array}{c}
B(g\otimes 1_{H};1_{A},x_{2})+B(x_{2}\otimes 1_{H};X_{2},x_{2}) \\
-B(x_{1}\otimes 1_{H};1_{A},x_{1}x_{2})+B(x_{1}\otimes 1_{H};X_{1},x_{2})+
\\
+B(gx_{1}x_{2}\otimes 1_{H};X_{2},x_{1}x_{2})+B(gx_{1}x_{2}\otimes
1_{H};X_{1}X_{2},x_{2})%
\end{array}%
=0
\end{equation*}%
and%
\begin{equation*}
\begin{array}{c}
B(g\otimes 1_{H};G,gx_{2})-B(x_{2}\otimes 1_{H};GX_{2},gx_{2}) \\
-B(x_{1}\otimes 1_{H};G,gx_{1}x_{2})-B(x_{1}\otimes 1_{H};GX_{1},gx_{2})+ \\
-B(gx_{1}x_{2}\otimes 1_{H};GX_{2},gx_{1}x_{2})+B(gx_{1}x_{2}\otimes
1_{H};GX_{1}X_{2},gx_{2})%
\end{array}%
=0.
\end{equation*}%
By taking in account the form of the elements $B(g\otimes 1_{H})$, $%
B(x_{1}\otimes 1_{H})$ and $B(x_{2}\otimes \ 1_{H}),$ we get%
\begin{equation*}
\begin{array}{c}
B(g\otimes 1_{H};1_{A},x_{2})-B(g\otimes 1_{H};1_{A},x_{2}) \\
-B(x_{1}\otimes 1_{H};1_{A},x_{1}x_{2})+\left[ -B(g\otimes
1_{H};X_{2},1_{H})+B(x_{1}\otimes 1_{H};1_{A},x_{1}x_{2})\right] + \\
+B(gx_{1}x_{2}\otimes 1_{H};X_{2},x_{1}x_{2})+B(gx_{1}x_{2}\otimes
1_{H};X_{1}X_{2},x_{2})%
\end{array}%
=0
\end{equation*}%
i.e.%
\begin{equation}
-B(g\otimes 1_{H};X_{2},1_{H})+B(gx_{1}x_{2}\otimes
1_{H};X_{2},x_{1}x_{2})+B(gx_{1}x_{2}\otimes 1_{H};X_{1}X_{2},x_{2})=0
\label{gx1x2ot1, fourth}
\end{equation}%
and

\begin{equation*}
\begin{array}{c}
B\left( g\otimes 1_{H};G,gx_{2}\right) -B(g\otimes 1_{H};G,gx_{2}) \\
-B(x_{1}\otimes 1_{H};G,gx_{1}x_{2})-\left[ -B(g\otimes
1_{H};GX_{2},g)-B(x_{1}\otimes 1_{H};G,gx_{1}x_{2})\right] \\
-B(gx_{1}x_{2}\otimes 1_{H};GX_{2},gx_{1}x_{2})+B(gx_{1}x_{2}\otimes
1_{H};GX_{1}X_{2},gx_{2})%
\end{array}%
=0.
\end{equation*}%
i.e.%
\begin{equation}
B(g\otimes 1_{H};GX_{2},g)-B(gx_{1}x_{2}\otimes
1_{H};GX_{2},gx_{1}x_{2})+B(gx_{1}x_{2}\otimes 1_{H};GX_{1}X_{2},gx_{2})=0.
\label{gx1x2ot1, fifth}
\end{equation}

\subsubsection{$G^{a}\otimes g^{d}x_{1}$}

\begin{equation*}
\sum_{\substack{ a,d=0  \\ a+d\equiv 0}}^{1}\left[
\begin{array}{c}
B(g\otimes 1_{H};G^{a},g^{d}x_{1})+B(x_{2}\otimes
1_{H};G^{a},g^{d}x_{1}x_{2})+ \\
\left( -1\right) ^{a}B(x_{2}\otimes 1_{H};G^{a}X_{2},g^{d}x_{1})+\left(
-1\right) ^{a}B(x_{1}\otimes 1_{H};G^{a}X_{1},g^{d}x_{1})+ \\
+\left( -1\right) ^{a}B(gx_{1}x_{2}\otimes
1_{H};G^{a}X_{1},g^{d}x_{1}x_{2})+B(gx_{1}x_{2}\otimes
1_{H};G^{a}X_{1}X_{2},g^{d}x_{1})%
\end{array}%
\right] G^{a}\otimes g^{d}x_{1}=0
\end{equation*}%
and we get%
\begin{equation*}
\begin{array}{c}
B(g\otimes 1_{H};1_{A},x_{1})+B(x_{2}\otimes 1_{H};1_{A},x_{1}x_{2})+ \\
+B(x_{2}\otimes 1_{H};X_{2},x_{1})+B(x_{1}\otimes 1_{H};X_{1},x_{1})+ \\
+B(gx_{1}x_{2}\otimes 1_{H};X_{1},x_{1}x_{2})+B(gx_{1}x_{2}\otimes
1_{H};X_{1}X_{2},x_{1})%
\end{array}%
=0
\end{equation*}%
and%
\begin{equation*}
\begin{array}{c}
B(g\otimes 1_{H};G,gx_{1})+B(x_{2}\otimes 1_{H};G,gx_{1}x_{2})+ \\
-B(x_{2}\otimes 1_{H};GX_{2},gx_{1})-B(x_{1}\otimes 1_{H};GX_{1},gx_{1})+ \\
-B(gx_{1}x_{2}\otimes 1_{H};GX_{1},gx_{1}x_{2})+B(gx_{1}x_{2}\otimes
1_{H};GX_{1}X_{2},gx_{1})%
\end{array}%
=0.
\end{equation*}%
By taking in account the form of the elements $B(g\otimes 1_{H})$, $%
B(x_{1}\otimes 1_{H})$ and $B(x_{2}\otimes \ 1_{H}),$ we get%
\begin{equation*}
\begin{array}{c}
B(g\otimes 1_{H};1_{A},x_{1})+B(x_{2}\otimes 1_{H};1_{A},x_{1}x_{2})+ \\
+\left[ -B(g\otimes 1_{H};1_{A},x_{1})-B(x_{2}\otimes \
1_{H};1_{A},x_{1}x_{2})\right] -B(g\otimes 1_{H};X_{1},1_{H})+ \\
+B(gx_{1}x_{2}\otimes 1_{H};X_{1},x_{1}x_{2})+B(gx_{1}x_{2}\otimes
1_{H};X_{1}X_{2},x_{1})%
\end{array}%
=0
\end{equation*}%
i.e%
\begin{equation}
-B(g\otimes 1_{H};X_{1},1_{H})+B(gx_{1}x_{2}\otimes
1_{H};X_{1},x_{1}x_{2})+B(gx_{1}x_{2}\otimes 1_{H};X_{1}X_{2},x_{1})=0
\label{gx1x2ot1, sixth}
\end{equation}%
and%
\begin{equation*}
\begin{array}{c}
B(g\otimes 1_{H};G,gx_{1})+B(x_{2}\otimes 1_{H};G,gx_{1}x_{2})+ \\
-\left[ B(g\otimes 1_{H};G,gx_{1})+B(x_{2}\otimes \ 1_{H};G,gx_{1}x_{2})%
\right] +B(g\otimes 1_{H};GX_{1},g)+ \\
-B(gx_{1}x_{2}\otimes 1_{H};GX_{1},gx_{1}x_{2})+B(gx_{1}x_{2}\otimes
1_{H};GX_{1}X_{2},gx_{1})%
\end{array}%
=0.
\end{equation*}%
i.e.%
\begin{equation}
B(g\otimes 1_{H};GX_{1},g)-B(gx_{1}x_{2}\otimes
1_{H};GX_{1},gx_{1}x_{2})+B(gx_{1}x_{2}\otimes 1_{H};GX_{1}X_{2},gx_{1})=0.
\label{gx1x2ot1, seventh}
\end{equation}

\subsubsection{$G^{a}X_{2}\otimes g^{d}$}

\begin{equation*}
\sum_{\substack{ a,d=0  \\ a+d\equiv 0}}^{1}\left[
\begin{array}{c}
B(g\otimes 1_{H};G^{a}X_{2},g^{d})+B(x_{2}\otimes
1_{H};G^{a}X_{2},g^{d}x_{2})+ \\
+B(x_{1}\otimes 1_{H};G^{a}X_{2},g^{d}x_{1})+\left( -1\right)
^{a}B(x_{1}\otimes 1_{H};G^{a}X_{1}X_{2},g^{d})+ \\
+B(gx_{1}x_{2}\otimes 1_{H};G^{a}X_{2},g^{d}x_{1}x_{2})+\left( -1\right)
^{a}B(gx_{1}x_{2}\otimes 1_{H};G^{a}X_{1}X_{2},g^{d}x_{2})%
\end{array}%
\right] G^{a}X_{2}\otimes g^{d}=0
\end{equation*}%
and we get%
\begin{equation*}
\begin{array}{c}
B(g\otimes 1_{H};X_{2},1_{H})+B(x_{2}\otimes 1_{H};X_{2},x_{2})+ \\
+B(x_{1}\otimes 1_{H};X_{2},x_{1})+B(x_{1}\otimes 1_{H};X_{1}X_{2},1_{H})+
\\
B(gx_{1}x_{2}\otimes 1_{H};X_{2},x_{1}x_{2})+B(gx_{1}x_{2}\otimes
1_{H};X_{1}X_{2},x_{2})%
\end{array}%
=0
\end{equation*}%
and%
\begin{equation*}
\begin{array}{c}
B(g\otimes 1_{H};GX_{2},g)+B(x_{2}\otimes 1_{H};GX_{2},gx_{2})+ \\
+B(x_{1}\otimes 1_{H};GX_{2},gx_{1})-B(x_{1}\otimes 1_{H};GX_{1}X_{2},G)+ \\
+B(gx_{1}x_{2}\otimes 1_{H};GX_{2},gx_{1}x_{2})-B(gx_{1}x_{2}\otimes
1_{H};GX_{1}X_{2},gx_{2})%
\end{array}%
=0.
\end{equation*}%
By taking in account the form of the elements $B(g\otimes 1_{H})$, $%
B(x_{1}\otimes 1_{H})$ and $B(x_{2}\otimes \ 1_{H}),$ we get%
\begin{equation*}
\begin{array}{c}
B\left( g\otimes 1_{H};1_{A},x_{2}\right) -B(g\otimes 1_{H};1_{A},x_{2})+ \\
-B(x_{1}\otimes 1_{H};1_{A},x_{1}x_{2})+\left[ -B(g\otimes
1_{H};X_{2},1_{H})+B(x_{1}\otimes 1_{H};1_{A},x_{1}x_{2})\right] + \\
B(gx_{1}x_{2}\otimes 1_{H};X_{2},x_{1}x_{2})+B(gx_{1}x_{2}\otimes
1_{H};X_{1}X_{2},x_{2})%
\end{array}%
=0
\end{equation*}%
i.e.%
\begin{equation}
-B(g\otimes 1_{H};X_{2},1_{H})+B(gx_{1}x_{2}\otimes
1_{H};X_{2},x_{1}x_{2})+B(gx_{1}x_{2}\otimes 1_{H};X_{1}X_{2},x_{2})=0
\label{gx1x2ot1, eight}
\end{equation}%
and%
\begin{equation*}
\begin{array}{c}
-B\left( g\otimes 1_{H};G,gx_{2}\right) +B(g\otimes 1_{H};G,gx_{2})+ \\
+B(x_{1}\otimes 1_{H};G,gx_{1}x_{2})-\left[ B(g\otimes
1_{H};GX_{2},g)+B(x_{1}\otimes 1_{H};G,gx_{1}x_{2})\right] \\
+B(gx_{1}x_{2}\otimes 1_{H};GX_{2},gx_{1}x_{2})-B(gx_{1}x_{2}\otimes
1_{H};GX_{1}X_{2},gx_{2})%
\end{array}%
=0.
\end{equation*}%
i.e.%
\begin{equation}
-B(g\otimes 1_{H};GX_{2},g)+B(gx_{1}x_{2}\otimes
1_{H};GX_{2},gx_{1}x_{2})-B(gx_{1}x_{2}\otimes 1_{H};GX_{1}X_{2},gx_{2})=0.
\label{gx1x2ot1, nine}
\end{equation}

\subsubsection{$G^{a}X_{1}\otimes g^{d}$}

\begin{equation*}
\sum_{\substack{ a,d=0  \\ a+d\equiv 0}}^{1}\left[
\begin{array}{c}
B(g\otimes 1_{H};G^{a}X_{1},g^{d})+B(x_{2}\otimes
1_{H};G^{a}X_{1},g^{d}x_{2})+ \\
\left( -1\right) ^{a+1}B(x_{2}\otimes
1_{H};G^{a}X_{1}X_{2},g^{d})+B(x_{1}\otimes 1_{H};G^{a}X_{1},g^{d}x_{1})+ \\
+B(gx_{1}x_{2}\otimes 1_{H};G^{a}X_{1},g^{d}x_{1}x_{2})+\left( -1\right)
^{a}B(gx_{1}x_{2}\otimes 1_{H};G^{a}X_{1}X_{2},g^{d}x_{1})%
\end{array}%
\right] G^{a}X_{1}\otimes g^{d}=0
\end{equation*}%
and we get%
\begin{equation*}
\begin{array}{c}
B(g\otimes 1_{H};X_{1},1_{H})+B(x_{2}\otimes 1_{H};X_{1},x_{2})+ \\
-B(x_{2}\otimes 1_{H};X_{1}X_{2},1_{H})+B(x_{1}\otimes 1_{H};X_{1},x_{1})+
\\
+B(gx_{1}x_{2}\otimes 1_{H};X_{1},x_{1}x_{2})+B(gx_{1}x_{2}\otimes
1_{H};X_{1}X_{2},x_{1})%
\end{array}%
=0
\end{equation*}%
and%
\begin{equation*}
\begin{array}{c}
B(g\otimes 1_{H};GX_{1},g)+B(x_{2}\otimes 1_{H};GX_{1},gx_{2})+ \\
+B(x_{2}\otimes 1_{H};GX_{1}X_{2},g)+B(x_{1}\otimes 1_{H};GX_{1},gx_{1})+ \\
+B(gx_{1}x_{2}\otimes 1_{H};GX_{1},gx_{1}x_{2})-B(gx_{1}x_{2}\otimes
1_{H};GX_{1}X_{2},gx_{1})%
\end{array}%
=0.
\end{equation*}%
By taking in account the form of the elements $B(g\otimes 1_{H})$, $%
B(x_{1}\otimes 1_{H})$ and $B(x_{2}\otimes \ 1_{H}),$ we get%
\begin{equation*}
\begin{array}{c}
B\left( g\otimes 1_{H};1_{A},x_{1}\right) +B(x_{2}\otimes
1_{H};1_{A},x_{1}x_{2})+ \\
-\left[ B(g\otimes 1_{H};1_{A},x_{1})+B(x_{2}\otimes \
1_{H};1_{A},x_{1}x_{2})\right] -B\left( g\otimes 1_{H};1_{A},x_{1}\right) +
\\
+B(gx_{1}x_{2}\otimes 1_{H};X_{1},x_{1}x_{2})+B(gx_{1}x_{2}\otimes
1_{H};X_{1}X_{2},x_{1})%
\end{array}%
=0
\end{equation*}%
i.e.

\begin{equation}
-B(g\otimes 1_{H};1_{A},x_{1})+B(gx_{1}x_{2}\otimes
1_{H};X_{1},x_{1}x_{2})+B(gx_{1}x_{2}\otimes 1_{H};X_{1}X_{2},x_{1})=0
\label{gx1x2ot1, ten}
\end{equation}%
and%
\begin{equation*}
\begin{array}{c}
B(g\otimes 1_{H};GX_{1},g)-B(x_{2}\otimes 1_{H};G,gx_{1}x_{2})+ \\
+\left[ B(g\otimes 1_{H};G,gx_{1})+B(x_{2}\otimes \ 1_{H};G,gx_{1}x_{2})%
\right] -B(g\otimes 1_{H};GX_{1},g)+ \\
+B(gx_{1}x_{2}\otimes 1_{H};GX_{1},gx_{1}x_{2})-B(gx_{1}x_{2}\otimes
1_{H};GX_{1}X_{2},gx_{1})%
\end{array}%
=0.
\end{equation*}%
i.e.

\begin{equation}
B(g\otimes 1_{H};G,gx_{1})+B(gx_{1}x_{2}\otimes
1_{H};GX_{1},gx_{1}x_{2})-B(gx_{1}x_{2}\otimes 1_{H};GX_{1}X_{2},gx_{1})=0.
\label{gx1x2ot1, eleven}
\end{equation}

\subsubsection{$G^{a}\otimes g^{d}x_{1}x_{2}$}

\begin{equation*}
\sum_{\substack{ a,d=0  \\ a+d\equiv 1}}^{1}\left[
\begin{array}{c}
B(g\otimes 1_{H};G^{a},g^{d}x_{1}x_{2})+\left( -1\right) ^{a}B(x_{2}\otimes
1_{H};G^{a}X_{2},g^{d}x_{1}x_{2})+ \\
+\left( -1\right) ^{a}B(x_{1}\otimes
1_{H};G^{a}X_{1},g^{d}x_{1}x_{2})+B(gx_{1}x_{2}\otimes
1_{H};G^{a}X_{1}X_{2},g^{d}x_{1}x_{2})%
\end{array}%
\right] G^{a}\otimes g^{d}x_{1}x_{2}=0
\end{equation*}%
and we get%
\begin{equation*}
\begin{array}{c}
B(g\otimes 1_{H};1_{A},gx_{1}x_{2})+B(x_{2}\otimes 1_{H};X_{2},gx_{1}x_{2})+
\\
+B(x_{1}\otimes 1_{H};X_{1},gx_{1}x_{2})+B(gx_{1}x_{2}\otimes
1_{H};X_{1}X_{2},gx_{1}x_{2})%
\end{array}%
=0
\end{equation*}%
and%
\begin{equation*}
\begin{array}{c}
B(g\otimes 1_{H};G,x_{1}x_{2})-B(x_{2}\otimes 1_{H};GX_{2},x_{1}x_{2})+ \\
-B(x_{1}\otimes 1_{H};GX_{1},x_{1}x_{2})+B(gx_{1}x_{2}\otimes
1_{H};GX_{1}X_{2},x_{1}x_{2})%
\end{array}%
=0.
\end{equation*}%
By taking in account the form of the elements $B(g\otimes 1_{H})$, $%
B(x_{1}\otimes 1_{H})$ and $B(x_{2}\otimes \ 1_{H}),$ we get%
\begin{equation*}
\begin{array}{c}
B\left( g\otimes 1_{H};1_{A},gx_{1}x_{2}\right) -B\left( g\otimes
1_{H};1_{A},gx_{1}x_{2}\right) + \\
-B\left( g\otimes 1_{H};1_{A},gx_{1}x_{2}\right) +B(gx_{1}x_{2}\otimes
1_{H};X_{1}X_{2},gx_{1}x_{2})%
\end{array}%
=0
\end{equation*}%
i.e.%
\begin{equation}
-B\left( g\otimes 1_{H};1_{A},gx_{1}x_{2}\right) +B(gx_{1}x_{2}\otimes
1_{H};X_{1}X_{2},gx_{1}x_{2})=0  \label{gx1x2ot1, twelve}
\end{equation}%
and%
\begin{equation*}
\begin{array}{c}
B(g\otimes 1_{H};G,x_{1}x_{2})-B(g\otimes 1_{H};G,x_{1}x_{2})+ \\
-B(g\otimes 1_{H};G,x_{1}x_{2})+B(gx_{1}x_{2}\otimes
1_{H};GX_{1}X_{2},x_{1}x_{2})%
\end{array}%
=0.
\end{equation*}%
i.e.%
\begin{equation}
-B(g\otimes 1_{H};G,x_{1}x_{2})+B(gx_{1}x_{2}\otimes
1_{H};GX_{1}X_{2},x_{1}x_{2})=0.  \label{gx1x2ot1, thirteen}
\end{equation}

\subsubsection{$G^{a}X_{2}\otimes g^{d}x_{2}$}

\begin{equation*}
\sum_{\substack{ a,d=0  \\ a+d\equiv 1}}^{1}\left[
\begin{array}{c}
B(g\otimes 1_{H};G^{a}X_{2},g^{d}x_{2})-B(x_{1}\otimes
1_{H};G^{a}X_{2},g^{d}x_{1}x_{2}) \\
+\left( -1\right) ^{a}B(x_{1}\otimes 1_{H};G^{a}X_{1}X_{2},g^{d}x_{2})%
\end{array}%
\right] G^{a}X_{2}\otimes g^{d}x_{2}=0
\end{equation*}%
and we get%
\begin{equation*}
\begin{array}{c}
B(g\otimes 1_{H};X_{2},gx_{2})-B(x_{1}\otimes 1_{H};X_{2},gx_{1}x_{2}) \\
+B(x_{1}\otimes 1_{H};X_{1}X_{2},gx_{2})%
\end{array}%
=0
\end{equation*}%
and%
\begin{equation*}
\begin{array}{c}
B(g\otimes 1_{H};GX_{2},x_{2})-B(x_{1}\otimes 1_{H};GX_{2},x_{1}x_{2}) \\
-B(x_{1}\otimes 1_{H};GX_{1}X_{2},x_{2})%
\end{array}%
=0.
\end{equation*}%
By taking in account the form of the elements $B(g\otimes 1_{H})$, $%
B(x_{1}\otimes 1_{H})$ and $B(x_{2}\otimes \ 1_{H}),$ we know that these
equalities are satisfied since all their summands are zero.

\subsubsection{$G^{a}X_{1}\otimes g^{d}x_{2}$}

\begin{equation*}
\sum_{\substack{ a,d=0  \\ a+d\equiv 1}}^{1}\left[
\begin{array}{c}
B(g\otimes 1_{H};G^{a}X_{1},g^{d}x_{2})+\left( -1\right)
^{a+1}B(x_{2}\otimes 1_{H};G^{a}X_{1}X_{2},g^{d}x_{2})+ \\
-B(x_{1}\otimes 1_{H};G^{a}X_{1},g^{d}x_{1}x_{2})+\left( -1\right)
^{a+1}B(gx_{1}x_{2}\otimes 1_{H};G^{a}X_{1}X_{2},g^{d}x_{1}x_{2})%
\end{array}%
\right] G^{a}X_{1}\otimes g^{d}x_{2}=0
\end{equation*}%
and we get%
\begin{equation*}
\begin{array}{c}
B(g\otimes 1_{H};X_{1},gx_{2})-B(x_{2}\otimes 1_{H};X_{1}X_{2},gx_{2})+ \\
-B(x_{1}\otimes 1_{H};X_{1},gx_{1}x_{2})-B(gx_{1}x_{2}\otimes
1_{H};X_{1}X_{2},gx_{1}x_{2})%
\end{array}%
=0
\end{equation*}%
and%
\begin{equation*}
\begin{array}{c}
B(g\otimes 1_{H};GX_{1},x_{2})+B(x_{2}\otimes 1_{H};GX_{1}X_{2},x_{2})+ \\
-B(x_{1}\otimes 1_{H};GX_{1},x_{1}x_{2})+B(gx_{1}x_{2}\otimes
1_{H};GX_{1}X_{2},x_{1}x_{2})%
\end{array}%
=0.
\end{equation*}%
By taking in account the form of the elements $B(g\otimes 1_{H})$, $%
B(x_{1}\otimes 1_{H})$ and $B(x_{2}\otimes \ 1_{H}),$ we get%
\begin{equation*}
\begin{array}{c}
-B\left( g\otimes 1_{H};1_{A},gx_{1}x_{2}\right) +B\left( g\otimes
1_{H};1_{A},gx_{1}x_{2}\right) + \\
+B\left( g\otimes 1_{H};1_{A},gx_{1}x_{2}\right) -B(gx_{1}x_{2}\otimes
1_{H};X_{1}X_{2},gx_{1}x_{2})%
\end{array}%
=0
\end{equation*}%
i.e.%
\begin{equation*}
B\left( g\otimes 1_{H};1_{A},gx_{1}x_{2}\right) -B(gx_{1}x_{2}\otimes
1_{H};X_{1}X_{2},gx_{1}x_{2})=0
\end{equation*}%
and%
\begin{equation*}
\begin{array}{c}
B\left( g\otimes 1_{H};G,x_{1}x_{2}\right) -B\left( g\otimes
1_{H};G,x_{1}x_{2}\right) + \\
-B\left( g\otimes 1_{H};G,x_{1}x_{2}\right) +B(gx_{1}x_{2}\otimes
1_{H};GX_{1}X_{2},x_{1}x_{2})%
\end{array}%
=0.
\end{equation*}%
i.e.%
\begin{equation*}
-B\left( g\otimes 1_{H};G,x_{1}x_{2}\right) +B(gx_{1}x_{2}\otimes
1_{H};GX_{1}X_{2},x_{1}x_{2})=0.
\end{equation*}%
These equalities are just $\left( \ref{gx1x2ot1, twelve}\right) $ and $%
\left( \ref{gx1x2ot1, thirteen}\right) $ we got above.

\subsubsection{$G^{a}X_{2}\otimes g^{d}x_{1}$}

\begin{equation*}
\sum_{\substack{ a,d=0  \\ a+d\equiv 1}}^{1}\left[
\begin{array}{c}
B(g\otimes 1_{H};G^{a}X_{2},g^{d}x_{1})+B(x_{2}\otimes
1_{H};G^{a}X_{2},g^{d}x_{1}x_{2})+ \\
+\left( -1\right) ^{a}B(x_{1}\otimes
1_{H};G^{a}X_{1}X_{2},g^{d}x_{1})+\left( -1\right) ^{a}B(gx_{1}x_{2}\otimes
1_{H};G^{a}X_{1}X_{2},g^{d}x_{1}x_{2})%
\end{array}%
\right] G^{a}X_{2}\otimes g^{d}x_{1}=0
\end{equation*}%
and we get%
\begin{equation*}
\begin{array}{c}
B(g\otimes 1_{H};X_{2},gx_{1})+B(x_{2}\otimes 1_{H};X_{2},gx_{1}x_{2})+ \\
+B(x_{1}\otimes 1_{H};X_{1}X_{2},gx_{1})+B(gx_{1}x_{2}\otimes
1_{H};X_{1}X_{2},gx_{1}x_{2})%
\end{array}%
=0
\end{equation*}%
and%
\begin{equation*}
\begin{array}{c}
B(g\otimes 1_{H};GX_{2},x_{1})+B(x_{2}\otimes 1_{H};GX_{2},x_{1}x_{2})+ \\
-B(x_{1}\otimes 1_{H};GX_{1}X_{2},x_{1})-B(gx_{1}x_{2}\otimes
1_{H};GX_{1}X_{2},x_{1}x_{2})%
\end{array}%
=0.
\end{equation*}%
By taking in account the form of the elements $B(g\otimes 1_{H})$, $%
B(x_{1}\otimes 1_{H})$ and $B(x_{2}\otimes \ 1_{H}),$ we get%
\begin{equation*}
\begin{array}{c}
B\left( g\otimes 1_{H};1_{A},gx_{1}x_{2}\right) -B(g\otimes
1_{H};X_{1}X_{2},g)+ \\
-B(g\otimes 1_{H};X_{1}X_{2},g)+B(gx_{1}x_{2}\otimes
1_{H};X_{1}X_{2},gx_{1}x_{2})%
\end{array}%
=0
\end{equation*}%
i.e.%
\begin{equation*}
-B\left( g\otimes 1_{H};1_{A},gx_{1}x_{2}\right) +B(gx_{1}x_{2}\otimes
1_{H};X_{1}X_{2},gx_{1}x_{2})=0
\end{equation*}%
and%
\begin{equation*}
\begin{array}{c}
-B\left( g\otimes 1_{H};G,x_{1}x_{2}\right) +B\left( g\otimes
1_{H};G,x_{1}x_{2}\right) + \\
+B\left( g\otimes 1_{H};G,x_{1}x_{2}\right) -B(gx_{1}x_{2}\otimes
1_{H};GX_{1}X_{2},x_{1}x_{2})%
\end{array}%
=0.
\end{equation*}%
i.e.%
\begin{equation*}
B\left( g\otimes 1_{H};G,x_{1}x_{2}\right) -B(gx_{1}x_{2}\otimes
1_{H};GX_{1}X_{2},x_{1}x_{2})=0.
\end{equation*}%
These equalities are just $\left( \ref{gx1x2ot1, twelve}\right) $ and $%
\left( \ref{gx1x2ot1, thirteen}\right) $ we got above.

\subsubsection{$G^{a}X_{1}\otimes g^{d}x_{1}$}

\begin{equation*}
\sum_{\substack{ a,d=0  \\ a+d\equiv 1}}^{1}\left[
\begin{array}{c}
B(g\otimes 1_{H};G^{a}X_{1},g^{d}x_{1})+B(x_{2}\otimes
1_{H};G^{a}X_{1},g^{d}x_{1}x_{2})+ \\
\left( -1\right) ^{a+1}B(x_{2}\otimes 1_{H};G^{a}X_{1}X_{2},g^{d}x_{1})%
\end{array}%
\right] G^{a}X_{1}\otimes g^{d}x_{1}=0
\end{equation*}%
and we get%
\begin{equation*}
\begin{array}{c}
B(g\otimes 1_{H};X_{1},gx_{1})+B(x_{2}\otimes 1_{H};X_{1},gx_{1}x_{2})+ \\
-B(x_{2}\otimes 1_{H};X_{1}X_{2},gx_{1})%
\end{array}%
=0
\end{equation*}%
and%
\begin{equation*}
\begin{array}{c}
B(g\otimes 1_{H};GX_{1},x_{1})+B(x_{2}\otimes 1_{H};GX_{1},x_{1}x_{2})+ \\
+B(x_{2}\otimes 1_{H};GX_{1}X_{2},x_{1})%
\end{array}%
=0.
\end{equation*}%
By taking in account the form of the elements $B(g\otimes 1_{H})$, $%
B(x_{1}\otimes 1_{H})$ and $B(x_{2}\otimes \ 1_{H}),$ we know that these
equalities are satisfied since all their summands are zero.

\subsubsection{$G^{a}X_{1}X_{2}\otimes g^{d}$}

\begin{equation*}
\sum_{\substack{ a,d=0  \\ a+d\equiv 1}}^{1}\left[
\begin{array}{c}
B(g\otimes 1_{H};G^{a}X_{1}X_{2},g^{d})+B(x_{2}\otimes
1_{H};G^{a}X_{1}X_{2},g^{d}x_{2})+ \\
+B(x_{1}\otimes 1_{H};G^{a}X_{1}X_{2},g^{d}x_{1})+B(gx_{1}x_{2}\otimes
1_{H};G^{a}X_{1}X_{2},g^{d}x_{1}x_{2})%
\end{array}%
\right] G^{a}X_{1}X_{2}\otimes g^{d}=0
\end{equation*}%
and we get%
\begin{equation*}
\begin{array}{c}
B(g\otimes 1_{H};X_{1}X_{2},g)+B(x_{2}\otimes 1_{H};X_{1}X_{2},gx_{2})+ \\
+B(x_{1}\otimes 1_{H};X_{1}X_{2},gx_{1})+B(gx_{1}x_{2}\otimes
1_{H};X_{1}X_{2},gx_{1}x_{2})%
\end{array}%
=0
\end{equation*}%
and%
\begin{equation*}
\begin{array}{c}
B(g\otimes 1_{H};GX_{1}X_{2},1_{H})+B(x_{2}\otimes 1_{H};GX_{1}X_{2},x_{2})+
\\
+B(x_{1}\otimes 1_{H};GX_{1}X_{2},x_{1})+B(gx_{1}x_{2}\otimes
1_{H};GX_{1}X_{2},x_{1}x_{2})%
\end{array}%
=0.
\end{equation*}%
By taking in account the form of the elements $B(g\otimes 1_{H})$, $%
B(x_{1}\otimes 1_{H})$ and $B(x_{2}\otimes \ 1_{H}),$ we get%
\begin{equation*}
\begin{array}{c}
B(g\otimes 1_{H};X_{1}X_{2},g)-B\left( g\otimes
1_{H};1_{A},gx_{1}x_{2}\right) \\
-B(g\otimes 1_{H};X_{1}X_{2},g)+B(gx_{1}x_{2}\otimes
1_{H};X_{1}X_{2},gx_{1}x_{2})%
\end{array}%
=0
\end{equation*}%
i.e%
\begin{equation*}
-B\left( g\otimes 1_{H};1_{A},gx_{1}x_{2}\right) +B(gx_{1}x_{2}\otimes
1_{H};X_{1}X_{2},gx_{1}x_{2})=0
\end{equation*}%
and%
\begin{equation*}
\begin{array}{c}
B(g\otimes 1_{H};GX_{1}X_{2},1_{H})-B(g\otimes 1_{H};GX_{1}X_{2},1_{H})+ \\
-B\left( g\otimes 1_{H};G,x_{1}x_{2}\right) +B(gx_{1}x_{2}\otimes
1_{H};GX_{1}X_{2},x_{1}x_{2})%
\end{array}%
=0
\end{equation*}%
i.e.%
\begin{equation*}
-B\left( g\otimes 1_{H};G,x_{1}x_{2}\right) +B(gx_{1}x_{2}\otimes
1_{H};GX_{1}X_{2},x_{1}x_{2})=0.
\end{equation*}%
These equalities are just $\left( \ref{gx1x2ot1, twelve}\right) $ and $%
\left( \ref{gx1x2ot1, thirteen}\right) $ we got above.

\subsubsection{$G^{a}X_{2}\otimes g^{d}x_{1}x_{2}$}

\begin{equation*}
\sum_{\substack{ a,d=0  \\ a+d\equiv 0}}^{1}\left[ B(g\otimes
1_{H};G^{a}X_{2},g^{d}x_{1}x_{2})+\left( -1\right) ^{a}B(x_{1}\otimes
1_{H};G^{a}X_{1}X_{2},g^{d}x_{1}x_{2})\right] G^{a}X_{2}\otimes
g^{d}x_{1}x_{2}=0
\end{equation*}%
and we get%
\begin{equation*}
B(g\otimes 1_{H};X_{2},x_{1}x_{2})+B(x_{1}\otimes
1_{H};X_{1}X_{2},x_{1}x_{2})=0
\end{equation*}%
and%
\begin{equation*}
B(g\otimes 1_{H};GX_{2},gx_{1}x_{2})-B(x_{1}\otimes
1_{H};GX_{1}X_{2},gx_{1}x_{2})=0
\end{equation*}%
By taking in account the form of the elements $B(g\otimes 1_{H})$ and $%
B(x_{1}\otimes 1_{H}),$ we know that these equalities are satisfied since
all their summands are zero.

\subsubsection{$G^{a}X_{1}\otimes g^{d}x_{1}x_{2}$}

\begin{equation*}
\sum_{\substack{ a,d=0  \\ a+d\equiv 0}}^{1}\left[ B(g\otimes
1_{H};G^{a}X_{1},g^{d}x_{1}x_{2})+\left( -1\right) ^{a+1}B(x_{2}\otimes
1_{H};G^{a}X_{1}X_{2},g^{d}x_{1}x_{2})\right] G^{a}X_{1}\otimes
g^{d}x_{1}x_{2}=0
\end{equation*}%
It is easily checked that by the form of the elements $B(g\otimes 1_{H})$
and $B(x_{2}\otimes 1_{H}),$ we know that these equalities are satisfied
since all their summands are zero.

\subsubsection{$G^{a}X_{1}X_{2}\otimes g^{d}x_{2}$}

\begin{equation*}
\sum_{\substack{ a,d=0  \\ a+d\equiv 0}}^{1}\left[ B(g\otimes
1_{H};G^{a}X_{1}X_{2},g^{d}x_{2})-B(x_{1}\otimes
1_{H};G^{a}X_{1}X_{2},g^{d}x_{1}x_{2})\right] G^{a}X_{1}X_{2}\otimes
g^{d}x_{2}=0
\end{equation*}%
Same conclusions as in case above.

\subsubsection{$G^{a}X_{1}X_{2}\otimes g^{d}x_{1}$}

\begin{equation*}
\sum_{\substack{ a,d=0  \\ a+d\equiv 0}}^{1}\left[ B(g\otimes
1_{H};G^{a}X_{1}X_{2},g^{d}x_{1})+B(x_{2}\otimes
1_{H};G^{a}X_{1}X_{2},g^{d}x_{1}x_{2})\right] G^{a}X_{1}X_{2}\otimes
g^{d}x_{1}=0.
\end{equation*}%
Same conclusions as in case above.

\subsubsection{$G^{a}X_{1}X_{2}\otimes g^{d}x_{1}x_{2}$}

\begin{equation*}
\sum_{\substack{ a,d=0  \\ a+d\equiv 0}}^{1}B(g\otimes
1_{H};G^{a}X_{1}X_{2},g^{d}x_{1}x_{2})G^{a}X_{1}X_{2}\otimes
g^{d}x_{1}x_{2}=0
\end{equation*}%
Same conclusions as in case above.

\subsection{Case $gx_{1}$}

\begin{eqnarray*}
&&\sum_{a,b_{1},b_{2},d,e_{1},e_{2}=0}^{1}\sum_{l_{1}=0}^{b_{1}}%
\sum_{l_{2}=0}^{b_{2}}\sum_{u_{1}=0}^{e_{1}}\sum_{u_{2}=0}^{e_{2}}\left(
-1\right) ^{\alpha \left( gx_{1};l_{1},l_{2},u_{1},u_{2}\right) } \\
&&B(x_{2}\otimes
1_{H};G^{a}X_{1}^{b_{1}}X_{2}^{b_{2}},g^{d}x_{1}^{e_{1}}x_{2}^{e_{2}}) \\
&&G^{a}X_{1}^{b_{1}-l_{1}}X_{2}^{b_{2}-l_{2}}\otimes
g^{d}x_{1}^{e_{1}-u_{1}}x_{2}^{e_{2}-u_{2}}\otimes
g^{a+b_{1}+b_{2}+l_{1}+l_{2}+d+e_{1}+e_{2}+u_{1}+u_{2}+1}x_{1}^{l_{1}+u_{1}+1}x_{2}^{l_{2}+u_{2}}
\\
&&+\sum_{a,b_{1},b_{2},d,e_{1},e_{2}=0}^{1}\sum_{l_{1}=0}^{b_{1}}%
\sum_{l_{2}=0}^{b_{2}}\sum_{u_{1}=0}^{e_{1}}\sum_{u_{2}=0}^{e_{2}}\left(
-1\right) ^{\alpha \left( g;l_{1},l_{2},u_{1},u_{2}\right) } \\
&&B(gx_{1}x_{2}\otimes
1_{H};G^{a}X_{1}^{b_{1}}X_{2}^{b_{2}},g^{d}x_{1}^{e_{1}}x_{2}^{e_{2}}) \\
&&G^{a}X_{1}^{b_{1}-l_{1}}X_{2}^{b_{2}-l_{2}}\otimes
g^{d}x_{1}^{e_{1}-u_{1}}x_{2}^{e_{2}-u_{2}}\otimes
g^{a+b_{1}+b_{2}+l_{1}+l_{2}+d+e_{1}+e_{2}+u_{1}+u_{2}+1}x_{1}^{l_{1}+u_{1}}x_{2}^{l_{2}+u_{2}}=0
\end{eqnarray*}%
We examine the first summand.%
\begin{eqnarray*}
&&\sum_{a,b_{1},b_{2},d,e_{1},e_{2}=0}^{1}\sum_{l_{1}=0}^{b_{1}}%
\sum_{l_{2}=0}^{b_{2}}\sum_{u_{1}=0}^{e_{1}}\sum_{u_{2}=0}^{e_{2}}\left(
-1\right) ^{\alpha \left( gx_{1};l_{1},l_{2},u_{1},u_{2}\right) } \\
&&B(x_{2}\otimes
1_{H};G^{a}X_{1}^{b_{1}}X_{2}^{b_{2}},g^{d}x_{1}^{e_{1}}x_{2}^{e_{2}}) \\
&&G^{a}X_{1}^{b_{1}-l_{1}}X_{2}^{b_{2}-l_{2}}\otimes
g^{d}x_{1}^{e_{1}-u_{1}}x_{2}^{e_{2}-u_{2}}\otimes
g^{a+b_{1}+b_{2}+l_{1}+l_{2}+d+e_{1}+e_{2}+u_{1}+u_{2}+1}x_{1}^{l_{1}+u_{1}+1}x_{2}^{l_{2}+u_{2}}
\end{eqnarray*}%
We get%
\begin{eqnarray*}
a+b_{1}+b_{2}+l_{1}+l_{2}+d+e_{1}+e_{2}+u_{1}+u_{2}+1 &\equiv &1 \\
l_{1}+u_{1}+1 &=&1 \\
l_{2}+u_{2} &=&0
\end{eqnarray*}%
i.e.%
\begin{eqnarray*}
a+b_{1}+b_{2}+d+e_{1}+e_{2} &\equiv &0 \\
l_{1}+u_{1} &=&0 \\
l_{2}+u_{2} &=&0
\end{eqnarray*}%
and we obtain, since $\alpha \left( gx_{1};0,0,0,0\right) =a+b_{1}+b_{2}$%
\begin{equation*}
\sum_{\substack{ a,b_{1},b_{2},d,e_{1},e_{2}=0  \\ %
a+b_{1}+b_{2}+d+e_{1}+e_{2}\equiv 0}}^{1}\left( -1\right)
^{a+b_{1}+b_{2}}B(x_{2}\otimes
1_{H};G^{a}X_{1}^{b_{1}}X_{2}^{b_{2}},g^{d}x_{1}^{e_{1}}x_{2}^{e_{2}})G^{a}X_{1}^{b_{1}}X_{2}^{b_{2}}\otimes g^{d}x_{1}^{e_{1}}x_{2}^{e_{2}}.
\end{equation*}%
We examine the second summand.

\begin{eqnarray*}
&&+\sum_{a,b_{1},b_{2},d,e_{1},e_{2}=0}^{1}\sum_{l_{1}=0}^{b_{1}}%
\sum_{l_{2}=0}^{b_{2}}\sum_{u_{1}=0}^{e_{1}}\sum_{u_{2}=0}^{e_{2}}\left(
-1\right) ^{\alpha \left( g;l_{1},l_{2},u_{1},u_{2}\right) } \\
&&B(gx_{1}x_{2}\otimes
1_{H};G^{a}X_{1}^{b_{1}}X_{2}^{b_{2}},g^{d}x_{1}^{e_{1}}x_{2}^{e_{2}}) \\
&&G^{a}X_{1}^{b_{1}-l_{1}}X_{2}^{b_{2}-l_{2}}\otimes
g^{d}x_{1}^{e_{1}-u_{1}}x_{2}^{e_{2}-u_{2}}\otimes
g^{a+b_{1}+b_{2}+l_{1}+l_{2}+d+e_{1}+e_{2}+u_{1}+u_{2}+1}x_{1}^{l_{1}+u_{1}}x_{2}^{l_{2}+u_{2}}
\end{eqnarray*}%
We get%
\begin{eqnarray*}
a+b_{1}+b_{2}+l_{1}+l_{2}+d+e_{1}+e_{2}+u_{1}+u_{2}+1 &\equiv &1 \\
l_{1}+u_{1} &=&1 \\
l_{2}+u_{2} &=&0
\end{eqnarray*}%
i.e.%
\begin{eqnarray*}
a+b_{1}+b_{2}+d+e_{1}+e_{2} &\equiv &1 \\
l_{1}+u_{1} &=&1 \\
l_{2} &=&u_{2}=0
\end{eqnarray*}%
and we obtain%
\begin{eqnarray*}
&&\sum_{\substack{ a,b_{1},b_{2},d,e_{1},e_{2}=0  \\ %
a+b_{1}+b_{2}+d+e_{1}+e_{2}\equiv 1}}^{1}\sum_{l_{1}=0}^{b_{1}}\sum
_{\substack{ u_{1}=0  \\ l_{1}+u_{1}=1}}^{e_{1}}\left( -1\right) ^{\alpha
\left( g;l_{1},0,u_{1},0\right) } \\
&&B(gx_{1}x_{2}\otimes
1_{H};G^{a}X_{1}^{b_{1}}X_{2}^{b_{2}},g^{d}x_{1}^{e_{1}}x_{2}^{e_{2}})G^{a}X_{1}^{b_{1}-l_{1}}X_{2}^{b_{2}}\otimes g^{d}x_{1}^{e_{1}-u_{1}}x_{2}^{e_{2}}.
\end{eqnarray*}%
Since $\alpha \left( g;0,0,1,0\right) =e_{2}+a+b_{1}+b_{2}+1$ and $\alpha
\left( g;1,0,0,0\right) =b_{2},$ we finally get
\begin{eqnarray*}
&&\sum_{\substack{ a,b_{1},b_{2},d,e_{2}=0  \\ a+b_{1}+b_{2}+d+e_{2}\equiv 0
}}^{1}\left( -1\right) ^{e_{2}+a+b_{1}+b_{2}+1}B(gx_{1}x_{2}\otimes
1_{H};G^{a}X_{1}^{b_{1}}X_{2}^{b_{2}},g^{d}x_{1}x_{2}^{e_{2}})G^{a}X_{1}^{b_{1}}X_{2}^{b_{2}}\otimes g^{d}x_{2}^{e_{2}}+
\\
&&\sum_{\substack{ a,b_{2},d,e_{1},e_{2}=0  \\ a+b_{2}+d+e_{1}+e_{2}\equiv 0
}}^{1}\left( -1\right) ^{b_{2}}B(gx_{1}x_{2}\otimes
1_{H};G^{a}X_{1}X_{2}^{b_{2}},g^{d}x_{1}^{e_{1}}x_{2}^{e_{2}})G^{a}X_{2}^{b_{2}}\otimes g^{d}x_{1}^{e_{1}}x_{2}^{e_{2}}.
\end{eqnarray*}%
In conclusion we obtain%
\begin{eqnarray*}
&&\sum_{\substack{ a,b_{1},b_{2},d,e_{1},e_{2}=0  \\ %
a+b_{1}+b_{2}+d+e_{1}+e_{2}\equiv 0}}^{1}\left( -1\right)
^{a+b_{1}+b_{2}}B(x_{2}\otimes
1_{H};G^{a}X_{1}^{b_{1}}X_{2}^{b_{2}},g^{d}x_{1}^{e_{1}}x_{2}^{e_{2}})G^{a}X_{1}^{b_{1}}X_{2}^{b_{2}}\otimes g^{d}x_{1}^{e_{1}}x_{2}^{e_{2}}
\\
&&\sum_{\substack{ a,b_{1},b_{2},d,e_{2}=0  \\ a+b_{1}+b_{2}+d+e_{2}\equiv 0
}}^{1}\left( -1\right) ^{e_{2}+a+b_{1}+b_{2}+1}B(gx_{1}x_{2}\otimes
1_{H};G^{a}X_{1}^{b_{1}}X_{2}^{b_{2}},g^{d}x_{1}x_{2}^{e_{2}})G^{a}X_{1}^{b_{1}}X_{2}^{b_{2}}\otimes g^{d}x_{2}^{e_{2}}+
\\
&&\sum_{\substack{ a,b_{2},d,e_{1},e_{2}=0  \\ a+b_{2}+d+e_{1}+e_{2}\equiv 0
}}^{1}\left( -1\right) ^{b_{2}}B(gx_{1}x_{2}\otimes
1_{H};G^{a}X_{1}X_{2}^{b_{2}},g^{d}x_{1}^{e_{1}}x_{2}^{e_{2}})G^{a}X_{2}^{b_{2}}\otimes g^{d}x_{1}^{e_{1}}x_{2}^{e_{2}}=0
\end{eqnarray*}

\subsubsection{$G^{a}\otimes g^{d}$}

\begin{equation*}
\sum_{\substack{ a,d=0  \\ a+d\equiv 0}}^{1}\left[
\begin{array}{c}
\left( -1\right) ^{a}B(x_{2}\otimes 1_{H};G^{a},g^{d})+\left( -1\right)
^{a+1}B(gx_{1}x_{2}\otimes 1_{H};G^{a},g^{d}x_{1})+ \\
+B(gx_{1}x_{2}\otimes 1_{H};G^{a}X_{1},g^{d})%
\end{array}%
\right] G^{a}\otimes g^{d}=0
\end{equation*}%
and we get%
\begin{equation}
B(x_{2}\otimes 1_{H};1_{A},1_{H})-B(gx_{1}x_{2}\otimes
1_{H};1_{A},x_{1})+B(gx_{1}x_{2}\otimes 1_{H};X_{1},1_{H})=0
\label{gx1x2ot1, fourteen}
\end{equation}%
and%
\begin{equation}
-B(x_{2}\otimes 1_{H};G,g)+B(gx_{1}x_{2}\otimes
1_{H};G,gx_{1})+B(gx_{1}x_{2}\otimes 1_{H};GX_{1},g)=0.
\label{gx1x2ot1, fifteen}
\end{equation}

\subsubsection{$G^{a}\otimes g^{d}x_{2}$}

\begin{equation*}
\sum_{\substack{ a,d=0  \\ a+d\equiv 1}}^{1}\left[
\begin{array}{c}
\left( -1\right) ^{a}B(x_{2}\otimes 1_{H};G^{a},g^{d}x_{2})+\left( -1\right)
^{a}B(gx_{1}x_{2}\otimes 1_{H};G^{a},g^{d}x_{1}x_{2})+ \\
+B(gx_{1}x_{2}\otimes 1_{H};G^{a}X_{1},g^{d}x_{2})%
\end{array}%
\right] G^{a}\otimes g^{d}x_{2}=0
\end{equation*}%
and we get%
\begin{equation}
B(x_{2}\otimes 1_{H};1_{A},gx_{2})+B(gx_{1}x_{2}\otimes
1_{H};1_{A},gx_{1}x_{2})+B(gx_{1}x_{2}\otimes 1_{H};X_{1},gx_{2})=0
\label{gx1x2ot1, sixteen}
\end{equation}%
and%
\begin{equation}
-B(x_{2}\otimes 1_{H};G,x_{2})-B(gx_{1}x_{2}\otimes
1_{H};G,x_{1}x_{2})+B(gx_{1}x_{2}\otimes 1_{H};GX_{1},x_{2})=0.
\label{gx1x2ot1, seventeen}
\end{equation}

\subsubsection{$G^{a}\otimes g^{d}x_{1}$}

\begin{equation*}
\sum_{\substack{ a,d=0  \\ a+d\equiv 1}}^{1}\left[ \left( -1\right)
^{a}B(x_{2}\otimes 1_{H};G^{a},g^{d}x_{1})+B(gx_{1}x_{2}\otimes
1_{H};G^{a}X_{1},g^{d}x_{1})\right] G^{a}\otimes g^{d}x_{1}=0
\end{equation*}%
and we get%
\begin{equation}
B(x_{2}\otimes 1_{H};1_{A},gx_{1})+B(gx_{1}x_{2}\otimes 1_{H};X_{1},gx_{1})=0
\label{gx1x2ot1, eighteen}
\end{equation}%
and%
\begin{equation}
-B(x_{2}\otimes 1_{H};G,x_{1})+B(gx_{1}x_{2}\otimes 1_{H};GX_{1},x_{1})=0.
\label{gx1x2ot1, nineteen}
\end{equation}

\subsubsection{$G^{a}X_{2}\otimes g^{d}$}

\begin{equation*}
\sum_{\substack{ a,d=0  \\ a+d\equiv 1}}^{1}\left[
\begin{array}{c}
\left( -1\right) ^{a+1}B(x_{2}\otimes 1_{H};G^{a}X_{2},g^{d})+\left(
-1\right) ^{a}B(gx_{1}x_{2}\otimes 1_{H};G^{a}X_{2},g^{d}x_{1})+ \\
-B(gx_{1}x_{2}\otimes 1_{H};G^{a}X_{1}X_{2},g^{d})%
\end{array}%
\right] G^{a}X_{2}\otimes g^{d}=0
\end{equation*}%
and we get%
\begin{equation*}
-B(x_{2}\otimes 1_{H};X_{2},g)+B(gx_{1}x_{2}\otimes
1_{H};X_{2},gx_{1})-B(gx_{1}x_{2}\otimes 1_{H};X_{1}X_{2},g)=0
\end{equation*}%
and%
\begin{equation*}
B(x_{2}\otimes 1_{H};GX_{2},1_{H})-B(gx_{1}x_{2}\otimes
1_{H};GX_{2},x_{1})-B(gx_{1}x_{2}\otimes 1_{H};GX_{1}X_{2},1_{H})=0.
\end{equation*}%
By taking in account the form of the element $B(x_{2}\otimes \ 1_{H}),$ we
get%
\begin{equation}
B(g\otimes 1_{H};1_{A},g)+B(x_{2}\otimes \
1_{H};1_{A},gx_{2})+B(gx_{1}x_{2}\otimes
1_{H};X_{2},gx_{1})-B(gx_{1}x_{2}\otimes 1_{H};X_{1}X_{2},g)=0
\label{gx1x2ot1, twenty}
\end{equation}%
and%
\begin{equation}
B(g\otimes 1_{H};G,1_{H})+B(x_{2}\otimes \
1_{H};G,x_{2})-B(gx_{1}x_{2}\otimes 1_{H};GX_{2},x_{1})-B(gx_{1}x_{2}\otimes
1_{H};GX_{1}X_{2},1_{H})=0.  \label{gx1x2ot1, twentyone}
\end{equation}

\subsubsection{$G^{a}X_{1}\otimes g^{d}$}

\begin{equation*}
\sum_{\substack{ a,d=0  \\ a+d\equiv 1}}^{1}\left[ \left( -1\right)
^{a+1}B(x_{2}\otimes 1_{H};G^{a}X_{1},g^{d})+\left( -1\right)
^{a}B(gx_{1}x_{2}\otimes 1_{H};G^{a}X_{1},g^{d}x_{1})\right]
G^{a}X_{1}\otimes g^{d}=0
\end{equation*}%
and we get%
\begin{equation*}
-B(x_{2}\otimes 1_{H};X_{1},g)+B(gx_{1}x_{2}\otimes 1_{H};X_{1},gx_{1})=0
\end{equation*}%
and%
\begin{equation*}
B(x_{2}\otimes 1_{H};GX_{1},1_{H})-B(gx_{1}x_{2}\otimes
1_{H};GX_{1},x_{1})=0.
\end{equation*}%
By taking in account the form of the element $B(x_{1}\otimes 1_{H}),$ we
obtain%
\begin{equation}
B(x_{2}\otimes 1_{H};1_{A},gx_{1})+B(gx_{1}x_{2}\otimes 1_{H};X_{1},gx_{1})=0
\label{gx1x2ot1, twentytwo}
\end{equation}%
and%
\begin{equation}
B(x_{2}\otimes 1_{H};G,x_{1})-B(gx_{1}x_{2}\otimes 1_{H};GX_{1},x_{1})=0
\label{gx1x2ot1, twentythree}
\end{equation}

\subsubsection{$G^{a}\otimes g^{d}x_{1}x_{2}$}

\begin{equation*}
\sum_{\substack{ a,d=0  \\ a+d\equiv 0}}^{1}\left[ \left( -1\right)
^{a}B(x_{2}\otimes 1_{H};G^{a},g^{d}x_{1}x_{2})+B(gx_{1}x_{2}\otimes
1_{H};G^{a}X_{1},g^{d}x_{1}x_{2})\right] G^{a}\otimes g^{d}x_{1}x_{2}=0
\end{equation*}%
and we get%
\begin{equation}
B(x_{2}\otimes 1_{H};1_{A},x_{1}x_{2})+B(gx_{1}x_{2}\otimes
1_{H};X_{1},x_{1}x_{2})=0  \label{gx1x2ot1, twentyfour}
\end{equation}%
and%
\begin{equation}
-B(x_{2}\otimes 1_{H};G,gx_{1}x_{2})+B(gx_{1}x_{2}\otimes
1_{H};GX_{1},gx_{1}x_{2})=0.  \label{gx1x2ot1, twentyfive}
\end{equation}

\subsubsection{$G^{a}X_{2}\otimes g^{d}x_{2}$}

\begin{equation*}
\sum_{\substack{ a,d=0  \\ a+d\equiv 0}}^{1}\left[
\begin{array}{c}
\left( -1\right) ^{a+1}B(x_{2}\otimes 1_{H};G^{a}X_{2},g^{d}x_{2})+\left(
-1\right) ^{a+1}B(gx_{1}x_{2}\otimes 1_{H};G^{a}X_{2},g^{d}x_{1}x_{2})+ \\
-B(gx_{1}x_{2}\otimes 1_{H};G^{a}X_{1}X_{2},g^{d}x_{2})%
\end{array}%
\right] G^{a}X_{2}\otimes g^{d}x_{2}
\end{equation*}%
and we get%
\begin{equation*}
-B(x_{2}\otimes 1_{H};X_{2},x_{2})-B(gx_{1}x_{2}\otimes
1_{H};X_{2},x_{1}x_{2})-B(gx_{1}x_{2}\otimes 1_{H};X_{1}X_{2},x_{2})=0
\end{equation*}%
and%
\begin{equation*}
B(x_{2}\otimes 1_{H};GX_{2},gx_{2})+B(gx_{1}x_{2}\otimes
1_{H};GX_{2},gx_{1}x_{2})-B(gx_{1}x_{2}\otimes 1_{H};GX_{1}X_{2},gx_{2})=0.
\end{equation*}%
By taking in account the form of the element $B(x_{2}\otimes 1_{H}),$ we
obtain%
\begin{equation}
B(g\otimes 1_{H};1_{A},x_{2})-B(gx_{1}x_{2}\otimes
1_{H};X_{2},x_{1}x_{2})-B(gx_{1}x_{2}\otimes 1_{H};X_{1}X_{2},x_{2})=0
\label{gx1x2ot1, twentysix}
\end{equation}%
and%
\begin{equation}
B(g\otimes 1_{H};G,gx_{2})+B(gx_{1}x_{2}\otimes
1_{H};GX_{2},gx_{1}x_{2})-B(gx_{1}x_{2}\otimes 1_{H};GX_{1}X_{2},gx_{2})=0.
\label{gx1x2ot1, twentyseven}
\end{equation}

\subsubsection{$G^{a}X_{1}\otimes g^{d}x_{2}$}

\begin{equation*}
\sum_{\substack{ a,d=0  \\ a+d\equiv 0}}^{1}\left[ \left( -1\right)
^{a+1}B(x_{2}\otimes 1_{H};G^{a}X_{1},g^{d}x_{2})+\left( -1\right)
^{a+1}B(gx_{1}x_{2}\otimes 1_{H};G^{a}X_{1},g^{d}x_{1}x_{2})\right]
G^{a}X_{1}\otimes g^{d}x_{2}=0
\end{equation*}%
and we get%
\begin{equation*}
-B(x_{2}\otimes 1_{H};X_{1},x_{2})-B(gx_{1}x_{2}\otimes
1_{H};X_{1},x_{1}x_{2})=0
\end{equation*}%
and%
\begin{equation*}
B(x_{2}\otimes 1_{H};GX_{1},gx_{2})+B(gx_{1}x_{2}\otimes
1_{H};GX_{1},gx_{1}x_{2})=0.
\end{equation*}%
By taking in account the form of the element $B(x_{2}\otimes 1_{H}),$ we
obtain%
\begin{equation*}
-B(x_{2}\otimes 1_{H};1_{A},x_{1}x_{2})-B(gx_{1}x_{2}\otimes
1_{H};X_{1},x_{1}x_{2})=0
\end{equation*}%
and%
\begin{equation*}
-B(x_{2}\otimes 1_{H};G,gx_{1}x_{2})+B(gx_{1}x_{2}\otimes
1_{H};GX_{1},gx_{1}x_{2})=0.
\end{equation*}%
These are respectively $\left( \ref{gx1x2ot1, twentyfour}\right) $ and $%
\left( \ref{gx1x2ot1, twentyfive}\right) .$

\subsubsection{$G^{a}X_{2}\otimes g^{d}x_{1}$}

\begin{equation*}
\sum_{\substack{ a,d=0  \\ a+d\equiv 0}}^{1}\left[ \left( -1\right)
^{a+1}B(x_{2}\otimes 1_{H};G^{a}X_{2},g^{d}x_{1})-B(gx_{1}x_{2}\otimes
1_{H};G^{a}X_{1}X_{2},g^{d}x_{1})\right] G^{a}X_{2}\otimes g^{d}x_{1}=0
\end{equation*}%
and we get%
\begin{equation*}
-B(x_{2}\otimes 1_{H};X_{2},x_{1})-B(gx_{1}x_{2}\otimes
1_{H};X_{1}X_{2},x_{1})=0
\end{equation*}%
and%
\begin{equation*}
B(x_{2}\otimes 1_{H};GX_{2},gx_{1})-B(gx_{1}x_{2}\otimes
1_{H};GX_{1}X_{2},gx_{1})=0.
\end{equation*}%
By taking in account the form of the element $B(x_{2}\otimes 1_{H}),$ we
obtain%
\begin{equation}
B(g\otimes 1_{H};1_{A},x_{1})+B(x_{2}\otimes \
1_{H};1_{A},x_{1}x_{2})-B(gx_{1}x_{2}\otimes 1_{H};X_{1}X_{2},x_{1})=0
\label{gx1x2ot1, twentyeight}
\end{equation}%
and%
\begin{equation}
B(g\otimes 1_{H};G,gx_{1})+B(x_{2}\otimes \
1_{H};G,gx_{1}x_{2})-B(gx_{1}x_{2}\otimes 1_{H};GX_{1}X_{2},gx_{1})=0.
\label{gx1x2ot1, twentynine}
\end{equation}

\subsubsection{$G^{a}X_{1}\otimes g^{d}x_{1}$}

\begin{equation*}
\sum_{\substack{ a,d=0  \\ a+d\equiv 0}}^{1}\left( -1\right)
^{a+1}B(x_{2}\otimes 1_{H};G^{a}X_{1},g^{d}x_{1})G^{a}X_{1}\otimes
g^{d}x_{1}=0
\end{equation*}%
and we get%
\begin{equation*}
-B(x_{2}\otimes 1_{H};X_{1},x_{1})=0
\end{equation*}%
and%
\begin{equation*}
B(x_{2}\otimes 1_{H};GX_{1},gx_{1})=0.
\end{equation*}%
These are already known as it is easily checked by the form of the element $%
B(x_{2}\otimes 1_{H}).$

\subsubsection{$G^{a}X_{1}X_{2}\otimes g^{d}$}

\begin{equation*}
\sum_{\substack{ a,d=0  \\ a+d\equiv 0}}^{1}\left[ \left( -1\right)
^{a}B(x_{2}\otimes 1_{H};G^{a}X_{1}X_{2},g^{d})+\left( -1\right)
^{a+1}B(gx_{1}x_{2}\otimes 1_{H};G^{a}X_{1}X_{2},g^{d}x_{1})\right]
G^{a}X_{1}X_{2}\otimes g^{d}=0
\end{equation*}%
and we get%
\begin{equation*}
B(x_{2}\otimes 1_{H};X_{1}X_{2},1_{H})-B(gx_{1}x_{2}\otimes
1_{H};X_{1}X_{2},x_{1})=0
\end{equation*}%
and%
\begin{equation*}
-B(x_{2}\otimes 1_{H};GX_{1}X_{2},g)+B(gx_{1}x_{2}\otimes
1_{H};GX_{1}X_{2},gx_{1})=0.
\end{equation*}%
By taking in account the form of the element $B(x_{2}\otimes 1_{H}),$ we
obtain%
\begin{equation*}
B(g\otimes 1_{H};1_{A},x_{1})+B(x_{2}\otimes \
1_{H};1_{A},x_{1}x_{2})-B(gx_{1}x_{2}\otimes 1_{H};X_{1}X_{2},x_{1})=0
\end{equation*}%
and%
\begin{equation*}
-B(g\otimes 1_{H};G,gx_{1})-B(x_{2}\otimes \
1_{H};G,gx_{1}x_{2})+B(gx_{1}x_{2}\otimes 1_{H};GX_{1}X_{2},gx_{1})=0.
\end{equation*}%
These are respectively $\left( \ref{gx1x2ot1, twentyeight}\right) $ and $%
\left( \ref{gx1x2ot1, twentynine}\right) .$

\subsubsection{$G^{a}X_{2}\otimes g^{d}x_{1}x_{2}$}

\begin{equation*}
\sum_{\substack{ a,d=0  \\ a+d\equiv 1}}^{1}\left[ \left( -1\right)
^{a+1}B(x_{2}\otimes 1_{H};G^{a}X_{2},g^{d}x_{1}x_{2})-B(gx_{1}x_{2}\otimes
1_{H};G^{a}X_{1}X_{2},g^{d}x_{1}x_{2})\right] G^{a}X_{2}\otimes
g^{d}x_{1}x_{2}=0
\end{equation*}%
and we get%
\begin{equation*}
-B(x_{2}\otimes 1_{H};X_{2},gx_{1}x_{2})-B(gx_{1}x_{2}\otimes
1_{H};X_{1}X_{2},gx_{1}x_{2})=0
\end{equation*}%
and%
\begin{equation*}
B(x_{2}\otimes 1_{H};GX_{2},x_{1}x_{2})-B(gx_{1}x_{2}\otimes
1_{H};GX_{1}X_{2},x_{1}x_{2})=0.
\end{equation*}%
By taking in account the form of the element $B(x_{2}\otimes 1_{H}),$ we
obtain%
\begin{equation}
B(g\otimes 1_{H};X_{1}X_{2},g)-B(gx_{1}x_{2}\otimes
1_{H};X_{1}X_{2},gx_{1}x_{2})=0  \label{gx1x2ot1, thirty}
\end{equation}%
and%
\begin{equation}
B(g\otimes 1_{H};GX_{1}X_{2},1_{H})-B(gx_{1}x_{2}\otimes
1_{H};GX_{1}X_{2},x_{1}x_{2})=0.  \label{gx1x2ot1, thirtyone}
\end{equation}

\subsubsection{$G^{a}X_{1}\otimes g^{d}x_{1}x_{2}$}

\begin{equation*}
\sum_{\substack{ a,d=0  \\ a+d\equiv 1}}^{1}\left( -1\right)
^{a+1}B(x_{2}\otimes 1_{H};G^{a}X_{1},g^{d}x_{1}x_{2})G^{a}X_{1}\otimes
g^{d}x_{1}x_{2}=0
\end{equation*}%
and we get%
\begin{equation*}
-B(x_{2}\otimes 1_{H};X_{1},gx_{1}x_{2})=0
\end{equation*}%
and%
\begin{equation*}
B(x_{2}\otimes 1_{H};GX_{1},x_{1}x_{2})=0.
\end{equation*}%
These are already known as it is easily checked by the form of the element $%
B(x_{2}\otimes 1_{H}).$

\subsubsection{$G^{a}X_{1}X_{2}\otimes g^{d}x_{2}$}

\begin{equation*}
\sum_{\substack{ a,d=0  \\ a+d\equiv 1}}^{1}\left[ \left( -1\right)
^{a}B(x_{2}\otimes 1_{H};G^{a}X_{1}X_{2},g^{d}x_{2})+\left( -1\right)
^{a}B(gx_{1}x_{2}\otimes 1_{H};G^{a}X_{1}X_{2},g^{d}x_{1}x_{2})\right]
G^{a}X_{1}X_{2}\otimes g^{d}x_{2}=0
\end{equation*}%
and we get%
\begin{equation*}
B(x_{2}\otimes 1_{H};X_{1}X_{2},gx_{2})+B(gx_{1}x_{2}\otimes
1_{H};X_{1}X_{2},gx_{1}x_{2})=0
\end{equation*}%
and%
\begin{equation*}
-B(x_{2}\otimes 1_{H};GX_{1}X_{2},x_{2})-B(gx_{1}x_{2}\otimes
1_{H};GX_{1}X_{2},x_{1}x_{2})=0.
\end{equation*}%
By taking in account the form of the element $B(x_{2}\otimes 1_{H}),$ we
obtain%
\begin{equation*}
-B(g\otimes 1_{H};X_{1}X_{2},g)+B(gx_{1}x_{2}\otimes
1_{H};X_{1}X_{2},gx_{1}x_{2})=0
\end{equation*}%
and%
\begin{equation*}
B(g\otimes 1_{H};GX_{1}X_{2},1_{H})-B(gx_{1}x_{2}\otimes
1_{H};GX_{1}X_{2},x_{1}x_{2})=0.
\end{equation*}%
These are respectively $\left( \ref{gx1x2ot1, thirty}\right) $ and $\left( %
\ref{gx1x2ot1, thirtyone}\right) .$

\subsubsection{$G^{a}X_{1}X_{2}\otimes g^{d}x_{1}$}

\begin{equation*}
\sum_{\substack{ a,d=0  \\ a+d\equiv 1}}^{1}\left( -1\right)
^{a}B(x_{2}\otimes 1_{H};G^{a}X_{1}X_{2},g^{d}x_{1})G^{a}X_{1}X_{2}\otimes
g^{d}x_{1}=0
\end{equation*}%
and we get%
\begin{equation*}
B(x_{2}\otimes 1_{H};X_{1}X_{2},gx_{1})=0
\end{equation*}%
and%
\begin{equation*}
-B(x_{2}\otimes 1_{H};GX_{1}X_{2},x_{1})=0.
\end{equation*}%
These are already known as it is easily checked by the form of the element $%
B(x_{2}\otimes 1_{H}).$

\subsubsection{$G^{a}X_{1}X_{2}\otimes g^{d}x_{1}x_{2}$}

\begin{equation*}
\sum_{\substack{ a,d=0  \\ a+d\equiv 0}}^{1}\left( -1\right)
^{a}B(x_{2}\otimes
1_{H};G^{a}X_{1}X_{2},g^{d}x_{1}x_{2})G^{a}X_{1}X_{2}\otimes
g^{d}x_{1}x_{2}=0
\end{equation*}%
and we get%
\begin{equation*}
B(x_{2}\otimes 1_{H};X_{1}X_{2},x_{1}x_{2})=0
\end{equation*}%
and%
\begin{equation*}
-B(x_{2}\otimes 1_{H};GX_{1}X_{2},g.x_{1}x_{2})=0
\end{equation*}%
These are already known as it is easily checked by the form of the element $%
B(x_{2}\otimes 1_{H}).$

\subsection{Case $gx_{2}$}

\begin{eqnarray*}
&&-\sum_{a,b_{1},b_{2},d,e_{1},e_{2}=0}^{1}\sum_{l_{1}=0}^{b_{1}}%
\sum_{l_{2}=0}^{b_{2}}\sum_{u_{1}=0}^{e_{1}}\sum_{u_{2}=0}^{e_{2}}\left(
-1\right) ^{\alpha \left( gx_{2};l_{1},l_{2},u_{1},u_{2}\right) } \\
&&B(x_{1}\otimes
1_{H};G^{a}X_{1}^{b_{1}}X_{2}^{b_{2}},g^{d}x_{1}^{e_{1}}x_{2}^{e_{2}}) \\
&&G^{a}X_{1}^{b_{1}-l_{1}}X_{2}^{b_{2}-l_{2}}\otimes
g^{d}x_{1}^{e_{1}-u_{1}}x_{2}^{e_{2}-u_{2}}\otimes
g^{a+b_{1}+b_{2}+l_{1}+l_{2}+d+e_{1}+e_{2}+u_{1}+u_{2}+1}x_{1}^{l_{1}+u_{1}}x_{2}^{l_{2}+u_{2}+1}
\\
&&+\sum_{a,b_{1},b_{2},d,e_{1},e_{2}=0}^{1}\sum_{l_{1}=0}^{b_{1}}%
\sum_{l_{2}=0}^{b_{2}}\sum_{u_{1}=0}^{e_{1}}\sum_{u_{2}=0}^{e_{2}}\left(
-1\right) ^{\alpha \left( g;l_{1},l_{2},u_{1},u_{2}\right) } \\
&&B(gx_{1}x_{2}\otimes
1_{H};G^{a}X_{1}^{b_{1}}X_{2}^{b_{2}},g^{d}x_{1}^{e_{1}}x_{2}^{e_{2}}) \\
&&G^{a}X_{1}^{b_{1}-l_{1}}X_{2}^{b_{2}-l_{2}}\otimes
g^{d}x_{1}^{e_{1}-u_{1}}x_{2}^{e_{2}-u_{2}}\otimes
g^{a+b_{1}+b_{2}+l_{1}+l_{2}+d+e_{1}+e_{2}+u_{1}+u_{2}+1}x_{1}^{l_{1}+u_{1}}x_{2}^{l_{2}+u_{2}}=0
\end{eqnarray*}%
We consider the first summand.

\begin{eqnarray*}
&&-\sum_{a,b_{1},b_{2},d,e_{1},e_{2}=0}^{1}\sum_{l_{1}=0}^{b_{1}}%
\sum_{l_{2}=0}^{b_{2}}\sum_{u_{1}=0}^{e_{1}}\sum_{u_{2}=0}^{e_{2}}\left(
-1\right) ^{\alpha \left( gx_{2};l_{1},l_{2},u_{1},u_{2}\right) } \\
&&B(x_{1}\otimes
1_{H};G^{a}X_{1}^{b_{1}}X_{2}^{b_{2}},g^{d}x_{1}^{e_{1}}x_{2}^{e_{2}}) \\
&&G^{a}X_{1}^{b_{1}-l_{1}}X_{2}^{b_{2}-l_{2}}\otimes
g^{d}x_{1}^{e_{1}-u_{1}}x_{2}^{e_{2}-u_{2}}\otimes
g^{a+b_{1}+b_{2}+l_{1}+l_{2}+d+e_{1}+e_{2}+u_{1}+u_{2}+1}x_{1}^{l_{1}+u_{1}}x_{2}^{l_{2}+u_{2}+1}.
\end{eqnarray*}%
\newline
We get%
\begin{eqnarray*}
a+b_{1}+b_{2}+l_{1}+l_{2}+d+e_{1}+e_{2}+u_{1}+u_{2}+1 &\equiv &1 \\
l_{1}+u_{1} &=&0 \\
l_{2}+u_{2}+1 &=&1
\end{eqnarray*}%
i.e.

\begin{eqnarray*}
a+b_{1}+b_{2}+d+e_{1}+e_{2} &\equiv &0 \\
l_{1} &=&u_{1}=0 \\
l_{2} &=&u_{2}=0.
\end{eqnarray*}%
Since $\alpha \left( gx_{2};0,0,0,0\right) =a+b_{1}+b_{2}$ we obtain%
\begin{equation*}
\sum_{\substack{ a,b_{1},b_{2},d,e_{1},e_{2}=0  \\ %
a+b_{1}+b_{2}+d+e_{1}+e_{2}\equiv 0}}^{1}\left( -1\right)
^{a+b_{1}+b_{2}+1}B(x_{1}\otimes
1_{H};G^{a}X_{1}^{b_{1}}X_{2}^{b_{2}},g^{d}x_{1}^{e_{1}}x_{2}^{e_{2}})G^{a}X_{1}^{b_{1}}X_{2}^{b_{2}}\otimes g^{d}x_{1}^{e_{1}}x_{2}^{e}.
\end{equation*}%
Let us consider the second summand

\begin{eqnarray*}
&&+\sum_{a,b_{1},b_{2},d,e_{1},e_{2}=0}^{1}\sum_{l_{1}=0}^{b_{1}}%
\sum_{l_{2}=0}^{b_{2}}\sum_{u_{1}=0}^{e_{1}}\sum_{u_{2}=0}^{e_{2}}\left(
-1\right) ^{\alpha \left( g;l_{1},l_{2},u_{1},u_{2}\right) } \\
&&B(gx_{1}x_{2}\otimes
1_{H};G^{a}X_{1}^{b_{1}}X_{2}^{b_{2}},g^{d}x_{1}^{e_{1}}x_{2}^{e_{2}}) \\
&&G^{a}X_{1}^{b_{1}-l_{1}}X_{2}^{b_{2}-l_{2}}\otimes
g^{d}x_{1}^{e_{1}-u_{1}}x_{2}^{e_{2}-u_{2}}\otimes
g^{a+b_{1}+b_{2}+l_{1}+l_{2}+d+e_{1}+e_{2}+u_{1}+u_{2}+1}x_{1}^{l_{1}+u_{1}}x_{2}^{l_{2}+u_{2}}=0.
\end{eqnarray*}%
We get%
\begin{eqnarray*}
a+b_{1}+b_{2}+l_{1}+l_{2}+d+e_{1}+e_{2}+u_{1}+u_{2}+1 &\equiv &1 \\
l_{1}+u_{1} &=&0 \\
l_{2}+u_{2} &=&1
\end{eqnarray*}%
i.e.%
\begin{eqnarray*}
a+b_{1}+b_{2}+d+e_{1}+e_{2} &\equiv &1 \\
l_{1} &=&u_{1}=0 \\
l_{2}+u_{2} &=&1.
\end{eqnarray*}%
We obtain%
\begin{eqnarray*}
&&\sum_{\substack{ a,b_{1},b_{2},d,e_{1},e_{2}=0  \\ %
a+b_{1}+b_{2}+d+e_{1}+e_{2}\equiv 1}}^{1}\sum_{l_{2}=0}^{b_{2}}\sum
_{\substack{ u_{2}=0  \\ l_{2}+u_{2}=1}}^{e_{2}}\left( -1\right) ^{\alpha
\left( g;0,l_{2},0,u_{2}\right) } \\
&&B(gx_{1}x_{2}\otimes
1_{H};G^{a}X_{1}^{b_{1}}X_{2}^{b_{2}},g^{d}x_{1}^{e_{1}}x_{2}^{e_{2}})G^{a}X_{1}^{b_{1}}X_{2}^{b_{2}-l_{2}}\otimes g^{d}x_{1}^{e_{1}}x_{2}^{e_{2}-u_{2}}.
\end{eqnarray*}%
Since $\alpha \left( g;0,0,0,1\right) =a+b_{1}+b_{2}+1$ and $\alpha \left(
g;0,1,0,0\right) =0,$ we thus get
\begin{eqnarray*}
&&\sum_{\substack{ a,b_{1},b_{2},d,e_{1}=0  \\ a+b_{1}+b_{2}+d+e_{1}\equiv 0
}}^{1}\left( -1\right) ^{a+b_{1}+b_{2}+1}B(gx_{1}x_{2}\otimes
1_{H};G^{a}X_{1}^{b_{1}}X_{2}^{b_{2}},g^{d}x_{1}^{e_{1}}x_{2})G^{a}X_{1}^{b_{1}}X_{2}^{b_{2}}\otimes g^{d}x_{1}^{e_{1}}+
\\
&&\sum_{\substack{ a,b_{1},d,e_{1},e_{2}=0  \\ a+b_{1}+d+e_{1}+e_{2}\equiv 0
}}^{1}B(gx_{1}x_{2}\otimes
1_{H};G^{a}X_{1}^{b_{1}}X_{2},g^{d}x_{1}^{e_{1}}x_{2}^{e_{2}})G^{a}X_{1}^{b_{1}}\otimes g^{d}x_{1}^{e_{1}}x_{2}^{e_{2}}.
\end{eqnarray*}%
In conclusion we get

\begin{eqnarray*}
&&\sum_{\substack{ a,b_{1},b_{2},d,e_{1},e_{2}=0  \\ %
a+b_{1}+b_{2}+d+e_{1}+e_{2}\equiv 0}}^{1}\left( -1\right)
^{a+b_{1}+b_{2}+1}B(x_{1}\otimes
1_{H};G^{a}X_{1}^{b_{1}}X_{2}^{b_{2}},g^{d}x_{1}^{e_{1}}x_{2}^{e_{2}})G^{a}X_{1}^{b_{1}}X_{2}^{b_{2}}\otimes g^{d}x_{1}^{e_{1}}x_{2}^{e}+
\\
+ &&\sum_{\substack{ a,b_{1},b_{2},d,e_{1}=0  \\ a+b_{1}+b_{2}+d+e_{1}\equiv
0 }}^{1}\left( -1\right) ^{a+b_{1}+b_{2}+1}B(gx_{1}x_{2}\otimes
1_{H};G^{a}X_{1}^{b_{1}}X_{2}^{b_{2}},g^{d}x_{1}^{e_{1}}x_{2})G^{a}X_{1}^{b_{1}}X_{2}^{b_{2}}\otimes g^{d}x_{1}^{e_{1}}+
\\
+ &&\sum_{\substack{ a,b_{1},d,e_{1},e_{2}=0  \\ a+b_{1}+d+e_{1}+e_{2}\equiv
0 }}^{1}B(gx_{1}x_{2}\otimes
1_{H};G^{a}X_{1}^{b_{1}}X_{2},g^{d}x_{1}^{e_{1}}x_{2}^{e_{2}})G^{a}X_{1}^{b_{1}}\otimes g^{d}x_{1}^{e_{1}}x_{2}^{e_{2}}=0.
\end{eqnarray*}

\subsubsection{$G^{a}\otimes g^{d}$}

\begin{equation*}
\sum_{\substack{ a,d=0  \\ a+d\equiv 0}}^{1}\left[
\begin{array}{c}
\left( -1\right) ^{a+1}B(x_{1}\otimes 1_{H};G^{a},g^{d})+\left( -1\right)
^{a+1}B(gx_{1}x_{2}\otimes 1_{H};G^{a},g^{d}x_{2})+ \\
+B(gx_{1}x_{2}\otimes 1_{H};G^{a}X_{2},g^{d})%
\end{array}%
\right] G^{a}\otimes g^{d}=0
\end{equation*}%
and we get%
\begin{equation}
-B(x_{1}\otimes 1_{H};1_{A},1_{H})-B(gx_{1}x_{2}\otimes
1_{H};1_{A},x_{2})+B(gx_{1}x_{2}\otimes 1_{H};X_{2},1_{H})=0
\label{gx1x2ot1, thirtytwo}
\end{equation}%
and%
\begin{equation}
B(x_{1}\otimes 1_{H};G,g)+B(gx_{1}x_{2}\otimes
1_{H};G,gx_{2})+B(gx_{1}x_{2}\otimes 1_{H};GX_{2},g)=0.
\label{gx1x2ot1, thirtythree}
\end{equation}

\subsubsection{$G^{a}\otimes g^{d}x_{2}$}

\begin{equation*}
\sum_{\substack{ a,d=0  \\ a+d\equiv 1}}^{1}\left[ \left( -1\right)
^{a+1}B(x_{1}\otimes 1_{H};G^{a},g^{d}x_{2})+B(gx_{1}x_{2}\otimes
1_{H};G^{a}X_{2},g^{d}x_{2})\right] G^{a}\otimes g^{d}x_{2}=0
\end{equation*}%
and we get%
\begin{equation}
-B(x_{1}\otimes 1_{H};1_{A},gx_{2})+B(gx_{1}x_{2}\otimes
1_{H};X_{2},gx_{2})=0  \label{gx1x2ot1, thirtyfour}
\end{equation}%
and%
\begin{equation}
B(x_{1}\otimes 1_{H};G,x_{2})+B(gx_{1}x_{2}\otimes 1_{H};GX_{2},x_{2})=0.
\label{gx1x2ot1, thirtyfive}
\end{equation}

\subsubsection{$G^{a}\otimes g^{d}x_{1}$}

\begin{equation*}
\sum_{\substack{ a,d=0  \\ a+d\equiv 1}}^{1}\left[
\begin{array}{c}
\left( -1\right) ^{a+1}B(x_{1}\otimes 1_{H};G^{a},g^{d}x_{1})+\left(
-1\right) ^{a+1}B(gx_{1}x_{2}\otimes 1_{H};G^{a},g^{d}x_{1}x_{2})+ \\
+B(gx_{1}x_{2}\otimes 1_{H};G^{a}X_{2},g^{d}x_{1})%
\end{array}%
\right] G^{a}\otimes g^{d}x_{1}=0
\end{equation*}%
and we get%
\begin{equation}
-B(x_{1}\otimes 1_{H};1_{A},gx_{1})-B(gx_{1}x_{2}\otimes
1_{H};1_{A},gx_{1}x_{2})+B(gx_{1}x_{2}\otimes 1_{H};X_{2},gx_{1})=0
\label{gx1x2ot1, thirtysix}
\end{equation}%
and%
\begin{equation}
B(x_{1}\otimes 1_{H};G,x_{1})+B(gx_{1}x_{2}\otimes
1_{H};G,x_{1}x_{2})+B(gx_{1}x_{2}\otimes 1_{H};GX_{2},x_{1})=0.
\label{gx1x2ot1, thirtyseven}
\end{equation}

\subsubsection{$G^{a}X_{2}\otimes g^{d}$}

\begin{equation*}
\sum_{\substack{ a,d=0  \\ a+d\equiv 1}}^{1}\left[ \left( -1\right)
^{a}B(x_{1}\otimes 1_{H};G^{a}X_{2},g^{d}+\left( -1\right)
^{a}B(gx_{1}x_{2}\otimes 1_{H};G^{a}X_{2},g^{d}x_{2})\right]
G^{a}X_{2}\otimes g^{d}=0
\end{equation*}%
and we get%
\begin{equation*}
B(x_{1}\otimes 1_{H};X_{2},g)+B(gx_{1}x_{2}\otimes 1_{H};X_{2},gx_{2})=0
\end{equation*}%
and%
\begin{equation*}
-B(x_{1}\otimes 1_{H};GX_{2},1_{H})-B(gx_{1}x_{2}\otimes
1_{H};GX_{2},x_{2})=0.
\end{equation*}%
By taking in account the form of the element $B(x_{1}\otimes 1_{H}),$ we
obtain%
\begin{equation*}
-B(x_{1}\otimes 1_{H};1_{A},gx_{2}))+B(gx_{1}x_{2}\otimes
1_{H};X_{2},gx_{2})=0
\end{equation*}%
and%
\begin{equation*}
-B(x_{1}\otimes 1_{H};G,x_{2})-B(gx_{1}x_{2}\otimes 1_{H};GX_{2},x_{2})=0.
\end{equation*}%
These are respectively $\left( \ref{gx1x2ot1, thirtyfour}\right) $ and $%
\left( \ref{gx1x2ot1, thirtyfive}\right) .$

\subsubsection{$G^{a}X_{1}\otimes g^{d}$}

\begin{equation*}
\sum_{\substack{ a,d=0  \\ a+d\equiv 1}}^{1}\left[
\begin{array}{c}
\left( -1\right) ^{a}B(x_{1}\otimes 1_{H};G^{a}X_{1},g^{d})+\left( -1\right)
^{a}B(gx_{1}x_{2}\otimes 1_{H};G^{a}X_{1},g^{d}x_{2})+ \\
+B(gx_{1}x_{2}\otimes 1_{H};G^{a}X_{1}X_{2},g^{d})%
\end{array}%
\right] G^{a}X_{1}\otimes g^{d}=0
\end{equation*}%
and we get%
\begin{equation*}
B(x_{1}\otimes 1_{H};X_{1},g)+B(gx_{1}x_{2}\otimes
1_{H};X_{1},gx_{2})+B(gx_{1}x_{2}\otimes 1_{H};X_{1}X_{2},g)=0
\end{equation*}%
and%
\begin{equation*}
-B(x_{1}\otimes 1_{H};GX_{1},1_{H})-B(gx_{1}x_{2}\otimes
1_{H};GX_{1},x_{2})+B(gx_{1}x_{2}\otimes 1_{H};GX_{1}X_{2},1_{H})=0.
\end{equation*}%
By taking in account the form of the element $B(x_{1}\otimes 1_{H}),$ we
obtain%
\begin{equation}
\left[ -B(g\otimes 1_{H};1_{A},g)-B(x_{1}\otimes 1_{H};1_{A},gx_{1})\right]
+B(gx_{1}x_{2}\otimes 1_{H};X_{1},gx_{2})+B(gx_{1}x_{2}\otimes
1_{H};X_{1}X_{2},g)=0  \label{gx1x2ot1, thirtyeight}
\end{equation}%
and%
\begin{equation}
\left[ -B(g\otimes 1_{H};G,1_{H})-B(x_{1}\otimes 1_{H};G,x_{1})\right]
-B(gx_{1}x_{2}\otimes 1_{H};GX_{1},x_{2})+B(gx_{1}x_{2}\otimes
1_{H};GX_{1}X_{2},1_{H})=0.  \label{gx1x2ot1, thirtynine}
\end{equation}

\subsubsection{$G^{a}\otimes g^{d}x_{1}x_{2}$}

\begin{equation*}
\sum_{\substack{ a,d=0  \\ a+d\equiv 0}}^{1}\left[ \left( -1\right)
^{a+1}B(x_{1}\otimes 1_{H};G^{a},g^{d}x_{1}x_{2})+B(gx_{1}x_{2}\otimes
1_{H};G^{a}X_{2},g^{d}x_{1}x_{2})\right] G^{a}\otimes g^{d}x_{1}x_{2}=0
\end{equation*}%
and we get%
\begin{equation}
-B(x_{1}\otimes 1_{H};1_{A},x_{1}x_{2})+B(gx_{1}x_{2}\otimes
1_{H};X_{2},x_{1}x_{2})=0  \label{gx1x2ot1, forty}
\end{equation}%
and%
\begin{equation}
B(x_{1}\otimes 1_{H};G,gx_{1}x_{2})+B(gx_{1}x_{2}\otimes
1_{H};GX_{2},gx_{1}x_{2})=0  \label{gx1x2ot1, fortyone}
\end{equation}

\subsubsection{$G^{a}X_{2}\otimes g^{d}x_{2}$}

\begin{equation*}
\sum_{\substack{ a,d=0  \\ a+d\equiv 0}}^{1}\left( -1\right)
^{a}B(x_{1}\otimes 1_{H};G^{a}X_{2},g^{d}x_{2})G^{a}X_{2}\otimes g^{d}x_{2}=0
\end{equation*}%
and we get%
\begin{equation*}
B(x_{1}\otimes 1_{H};X_{2},x_{2})=0
\end{equation*}%
and%
\begin{equation*}
-B(x_{1}\otimes 1_{H};GX_{2},gx_{2})=0.
\end{equation*}%
These are already known as it is easily checked by the form of the element $%
B(x_{1}\otimes 1_{H}).$

\subsubsection{$G^{a}X_{1}\otimes g^{d}x_{2}$}

\begin{equation*}
\sum_{\substack{ a,d=0  \\ a+d\equiv 0}}^{1}\left[ \left( -1\right)
^{a}B(x_{1}\otimes 1_{H};G^{a}X_{1},g^{d}x_{2})+B(gx_{1}x_{2}\otimes
1_{H};G^{a}X_{1}X_{2},g^{d}x_{2})\right] G^{a}X_{1}\otimes g^{d}x_{2}=0
\end{equation*}%
and we get%
\begin{equation*}
B(x_{1}\otimes 1_{H};X_{1},x_{2})+B(gx_{1}x_{2}\otimes
1_{H};X_{1}X_{2},x_{2})=0
\end{equation*}%
and%
\begin{equation*}
-B(x_{1}\otimes 1_{H};GX_{1},gx_{2})+B(gx_{1}x_{2}\otimes
1_{H};GX_{1}X_{2},gx_{2})=0.
\end{equation*}%
By taking in account the form of the element $B(x_{1}\otimes 1_{H}),$ we
obtain
\begin{equation}
\left[ -B(g\otimes 1_{H};X_{2},1_{H})+B(x_{1}\otimes 1_{H};1_{A},x_{1}x_{2})%
\right] +B(gx_{1}x_{2}\otimes 1_{H};X_{1}X_{2},x_{2})=0
\label{gx1x2ot1, fortytwo}
\end{equation}%
and
\begin{equation}
\left[ B(g\otimes 1_{H};GX_{2},g)+B(x_{1}\otimes 1_{H};G,gx_{1}x_{2})\right]
+B(gx_{1}x_{2}\otimes 1_{H};GX_{1}X_{2},gx_{2})=0.
\label{gx1x2ot1, fortythree}
\end{equation}

\subsubsection{$G^{a}X_{2}\otimes g^{d}x_{1}$}

\begin{equation*}
\sum_{\substack{ a,d=0  \\ a+d\equiv 0}}^{1}\left[ \left( -1\right)
^{a}B(x_{1}\otimes 1_{H};G^{a}X_{2},g^{d}x_{1})+\left( -1\right)
^{a}B(gx_{1}x_{2}\otimes 1_{H};G^{a}X_{2},g^{d}x_{1}x_{2})\right]
G^{a}X_{2}\otimes g^{d}x_{1}=0
\end{equation*}%
and we get%
\begin{equation*}
B(x_{1}\otimes 1_{H};X_{2},x_{1})+B(gx_{1}x_{2}\otimes
1_{H};X_{2},x_{1}x_{2})=0
\end{equation*}%
and%
\begin{equation*}
-B(x_{1}\otimes 1_{H};GX_{2},gx_{1})-B(gx_{1}x_{2}\otimes
1_{H};GX_{2},gx_{1}x_{2})=0.
\end{equation*}%
By taking in account the form of the element $B(x_{1}\otimes 1_{H}),$ we
obtain%
\begin{equation*}
-B(x_{1}\otimes 1_{H};1_{A},x_{1}x_{2})+B(gx_{1}x_{2}\otimes
1_{H};X_{2},x_{1}x_{2})=0
\end{equation*}%
and%
\begin{equation*}
-B(x_{1}\otimes 1_{H};G,gx_{1}x_{2})-B(gx_{1}x_{2}\otimes
1_{H};GX_{2},gx_{1}x_{2})=0.
\end{equation*}%
These are respectively $\left( \ref{gx1x2ot1, forty}\right) $ and $\left( %
\ref{gx1x2ot1, fortyone}\right) .$

\subsubsection{$G^{a}X_{1}\otimes g^{d}x_{1}$}

\begin{equation*}
\sum_{\substack{ a,d=0  \\ a+d\equiv 0}}^{1}\left[
\begin{array}{c}
\left( -1\right) ^{a}B(x_{1}\otimes 1_{H};G^{a}X_{1},g^{d}x_{1})+\left(
-1\right) ^{a}B(gx_{1}x_{2}\otimes 1_{H};G^{a}X_{1},g^{d}x_{1}x_{2})+ \\
B(gx_{1}x_{2}\otimes 1_{H};G^{a}X_{1}X_{2},g^{d}x_{1})%
\end{array}%
\right] G^{a}X_{1}\otimes g^{d}x_{1}=0
\end{equation*}%
and we get%
\begin{equation*}
B(x_{1}\otimes 1_{H};X_{1},x_{1})+B(gx_{1}x_{2}\otimes
1_{H};X_{1},x_{1}x_{2})+B(gx_{1}x_{2}\otimes 1_{H};X_{1}X_{2},x_{1})=0
\end{equation*}%
and%
\begin{equation*}
-B(x_{1}\otimes 1_{H};GX_{1},gx_{1})-B(gx_{1}x_{2}\otimes
1_{H};GX_{1},gx_{1}x_{2})+B(gx_{1}x_{2}\otimes 1_{H};GX_{1}X_{2},gx_{1})=0.
\end{equation*}%
By taking in account the form of the element $B(x_{1}\otimes 1_{H}),$ we
obtain%
\begin{equation}
-B(g\otimes 1_{H};X_{1},1_{H})+B(gx_{1}x_{2}\otimes
1_{H};X_{1},x_{1}x_{2})+B(gx_{1}x_{2}\otimes 1_{H};X_{1}X_{2},x_{1})=0
\label{gx1x2ot1, fortyfour}
\end{equation}%
and%
\begin{equation}
B(g\otimes 1_{H};GX_{1},g)-B(gx_{1}x_{2}\otimes
1_{H};GX_{1},gx_{1}x_{2})+B(gx_{1}x_{2}\otimes 1_{H};GX_{1}X_{2},gx_{1})=0.
\label{gx1x2ot1, fortyfive}
\end{equation}

\subsubsection{$G^{a}X_{1}X_{2}\otimes g^{d}$}

\begin{equation*}
\sum_{\substack{ a,d=0  \\ a+d\equiv 0}}^{1}\left[ \left( -1\right)
^{a+1}B(x_{1}\otimes 1_{H};G^{a}X_{1}X_{2},g^{d})+\left( -1\right)
^{a+1}B(gx_{1}x_{2}\otimes 1_{H};G^{a}X_{1}X_{2},g^{d}x_{2})\right]
G^{a}X_{1}X_{2}\otimes g^{d}=0
\end{equation*}%
and we get%
\begin{equation*}
-B(x_{1}\otimes 1_{H};X_{1}X_{2},1_{H})-B(gx_{1}x_{2}\otimes
1_{H};X_{1}X_{2},x_{2})=0
\end{equation*}%
and%
\begin{equation*}
B(x_{1}\otimes 1_{H};GX_{1}X_{2},g)+B(gx_{1}x_{2}\otimes
1_{H};GX_{1}X_{2},gx_{2})=0.
\end{equation*}%
By taking in account the form of the element $B(x_{1}\otimes 1_{H}),$ we
obtain%
\begin{equation*}
\left[ B(g\otimes 1_{H};X_{2},1_{H})-B(x_{1}\otimes 1_{H};1_{A},x_{1}x_{2})%
\right] -B(gx_{1}x_{2}\otimes 1_{H};X_{1}X_{2},x_{2})=0
\end{equation*}%
and%
\begin{equation*}
\left[ B(g\otimes 1_{H};GX_{2},g)+B(x_{1}\otimes 1_{H};G,gx_{1}x_{2})\right]
+B(gx_{1}x_{2}\otimes 1_{H};GX_{1}X_{2},gx_{2})=0.
\end{equation*}%
These are respectively $\left( \ref{gx1x2ot1, fortytwo}\right) $ and $\left( %
\ref{gx1x2ot1, fortythree}\right) .$

\subsubsection{$G^{a}X_{2}\otimes g^{d}x_{1}x_{2}$}

\begin{equation*}
\sum_{\substack{ a,d=0  \\ a+d\equiv 1}}^{1}\left( -1\right)
^{a}B(x_{1}\otimes 1_{H};G^{a}X_{2},g^{d}x_{1}x_{2})G^{a}X_{2}\otimes
g^{d}x_{1}x_{2}=0
\end{equation*}%
and we get%
\begin{equation*}
B(x_{1}\otimes 1_{H};X_{2},gx_{1}x_{2})=0
\end{equation*}%
and%
\begin{equation*}
-B(x_{1}\otimes 1_{H};GX_{2},x_{1}x_{2})=0.
\end{equation*}%
These are already known as it is easily checked by the form of the element $%
B(x_{1}\otimes 1_{H}).$

\subsubsection{$G^{a}X_{1}\otimes g^{d}x_{1}x_{2}$}

\begin{equation*}
\sum_{\substack{ a,d=0  \\ a+d\equiv 1}}^{1}\left[ \left( -1\right)
^{a}B(x_{1}\otimes 1_{H};G^{a}X_{1},g^{d}x_{1}x_{2})+B(gx_{1}x_{2}\otimes
1_{H};G^{a}X_{1}X_{2},g^{d}x_{1}x_{2})\right] G^{a}X_{1}\otimes
g^{d}x_{1}x_{2}=0
\end{equation*}%
and we get%
\begin{equation*}
B(x_{1}\otimes 1_{H};X_{1},gx_{1}x_{2})+B(gx_{1}x_{2}\otimes
1_{H};X_{1}X_{2},gx_{1}x_{2})=0
\end{equation*}%
and%
\begin{equation*}
-B(x_{1}\otimes 1_{H};GX_{1},x_{1}x_{2})+B(gx_{1}x_{2}\otimes
1_{H};GX_{1}X_{2},x_{1}x_{2})=0.
\end{equation*}%
By taking in account the form of the element $B(x_{1}\otimes 1_{H}),$ we
obtain%
\begin{equation}
-B(g\otimes 1_{H};X_{1}X_{2},g)+B(gx_{1}x_{2}\otimes
1_{H};X_{1}X_{2},gx_{1}x_{2})=0  \label{gx1x2ot1, fortysix}
\end{equation}%
and%
\begin{equation}
-B(g\otimes 1_{H};GX_{1}X_{2},1_{H})+B(gx_{1}x_{2}\otimes
1_{H};GX_{1}X_{2},x_{1}x_{2})=0.  \label{gx1x2ot1, fortyseven}
\end{equation}

\subsubsection{$G^{a}X_{1}X_{2}\otimes g^{d}x_{2}$}

\begin{equation*}
\sum_{\substack{ a,d=0  \\ a+d\equiv 1}}^{1}\left( -1\right)
^{a+1}B(x_{1}\otimes 1_{H};G^{a}X_{1}X_{2},g^{d}x_{2})G^{a}X_{1}X_{2}\otimes
g^{d}x_{2}=0
\end{equation*}%
and we get%
\begin{equation*}
-B(x_{1}\otimes 1_{H};X_{1}X_{2},gx_{2})=0
\end{equation*}%
and%
\begin{equation*}
B(x_{1}\otimes 1_{H};GX_{1}X_{2},x_{2})=0.
\end{equation*}%
These are already known as it is easily checked by the form of the element $%
B(x_{1}\otimes 1_{H}).$

\subsubsection{$G^{a}X_{1}X_{2}\otimes g^{d}x_{1}$}

\begin{gather*}
\sum_{\substack{ a,d=0  \\ a+d\equiv 1}}^{1}\left[ \left( -1\right)
^{a+1}B(x_{1}\otimes 1_{H};G^{a}X_{1}X_{2},g^{d}x_{1})+\left( -1\right)
^{a+1}B(gx_{1}x_{2}\otimes 1_{H};G^{a}X_{1}X_{2},g^{d}x_{1}x_{2})\right] \\
G^{a}X_{1}X_{2}\otimes g^{d}x_{1}=0
\end{gather*}%
and we get
\begin{equation*}
-B(x_{1}\otimes 1_{H};X_{1}X_{2},gx_{1})-B(gx_{1}x_{2}\otimes
1_{H};X_{1}X_{2},gx_{1}x_{2})=0.
\end{equation*}%
and%
\begin{equation*}
B(x_{1}\otimes 1_{H};GX_{1}X_{2},x_{1})+B(gx_{1}x_{2}\otimes
1_{H};GX_{1}X_{2},x_{1}x_{2})=0.
\end{equation*}%
By taking in account the form of the element $B(x_{1}\otimes 1_{H}),$ we
obtain%
\begin{equation*}
B(g\otimes 1_{H};X_{1}X_{2},g)-B(gx_{1}x_{2}\otimes
1_{H};X_{1}X_{2},gx_{1}x_{2})=0
\end{equation*}%
and%
\begin{equation*}
-B(g\otimes 1_{H};GX_{1}X_{2},1_{H})+B(gx_{1}x_{2}\otimes
1_{H};GX_{1}X_{2},x_{1}x_{2})=0.
\end{equation*}%
These are respectively $\left( \ref{gx1x2ot1, fortysix}\right) $ and $\left( %
\ref{gx1x2ot1, fortyseven}\right) .$

\subsubsection{$G^{a}X_{1}X_{2}\otimes g^{d}x_{1}x_{2}$}

\begin{equation*}
\sum_{\substack{ a,d=0  \\ a+d\equiv 0}}^{1}\left( -1\right)
^{a+1}B(x_{1}\otimes
1_{H};G^{a}X_{1}X_{2},g^{d}x_{1}x_{2})G^{a}X_{1}X_{2}\otimes
g^{d}x_{1}x_{2}=0
\end{equation*}%
and we get%
\begin{equation*}
-B(x_{1}\otimes 1_{H};X_{1}X_{2},x_{1}x_{2})=0
\end{equation*}%
and%
\begin{equation*}
B(x_{1}\otimes 1_{H};GX_{1}X_{2},gx_{1}x_{2})=0
\end{equation*}%
These are already known as it is easily checked by the form of the element $%
B(x_{1}\otimes 1_{H}).$

\subsection{Case $gx_{1}x_{2}$}

\begin{eqnarray*}
&&\sum_{a,b_{1},b_{2},d,e_{1},e_{2}=0}^{1}\sum_{l_{1}=0}^{b_{1}}%
\sum_{l_{2}=0}^{b_{2}}\sum_{u_{1}=0}^{e_{1}}\sum_{u_{2}=0}^{e_{2}}\left(
-1\right) ^{\alpha \left( gx_{1}x_{2};l_{1},l_{2},u_{1},u_{2}\right) } \\
&&B(g\otimes
1_{H};G^{a}X_{1}^{b_{1}}X_{2}^{b_{2}},g^{d}x_{1}^{e_{1}}x_{2}^{e_{2}}) \\
&&G^{a}X_{1}^{b_{1}-l_{1}}X_{2}^{b_{2}-l_{2}}\otimes
g^{d}x_{1}^{e_{1}-u_{1}}x_{2}^{e_{2}-u_{2}}\otimes
g^{a+b_{1}+b_{2}+l_{1}+l_{2}+d+e_{1}+e_{2}+u_{1}+u_{2}+1}x_{1}^{l_{1}+u_{1}+1}x_{2}^{l_{2}+u_{2}+1}
\\
+
&&\sum_{a,b_{1},b_{2},d,e_{1},e_{2}=0}^{1}\sum_{l_{1}=0}^{b_{1}}%
\sum_{l_{2}=0}^{b_{2}}\sum_{u_{1}=0}^{e_{1}}\sum_{u_{2}=0}^{e_{2}}\left(
-1\right) ^{\alpha \left( gx_{1};l_{1},l_{2},u_{1},u_{2}\right) } \\
&&B(x_{2}\otimes
1_{H};G^{a}X_{1}^{b_{1}}X_{2}^{b_{2}},g^{d}x_{1}^{e_{1}}x_{2}^{e_{2}}) \\
&&G^{a}X_{1}^{b_{1}-l_{1}}X_{2}^{b_{2}-l_{2}}\otimes
g^{d}x_{1}^{e_{1}-u_{1}}x_{2}^{e_{2}-u_{2}}\otimes
g^{a+b_{1}+b_{2}+l_{1}+l_{2}+d+e_{1}+e_{2}+u_{1}+u_{2}+1}x_{1}^{l_{1}+u_{1}+1}x_{2}^{l_{2}+u_{2}}
\\
&&-\sum_{a,b_{1},b_{2},d,e_{1},e_{2}=0}^{1}\sum_{l_{1}=0}^{b_{1}}%
\sum_{l_{2}=0}^{b_{2}}\sum_{u_{1}=0}^{e_{1}}\sum_{u_{2}=0}^{e_{2}}\left(
-1\right) ^{\alpha \left( gx_{2};l_{1},l_{2},u_{1},u_{2}\right) } \\
&&B(x_{1}\otimes
1_{H};G^{a}X_{1}^{b_{1}}X_{2}^{b_{2}},g^{d}x_{1}^{e_{1}}x_{2}^{e_{2}}) \\
&&G^{a}X_{1}^{b_{1}-l_{1}}X_{2}^{b_{2}-l_{2}}\otimes
g^{d}x_{1}^{e_{1}-u_{1}}x_{2}^{e_{2}-u_{2}}\otimes
g^{a+b_{1}+b_{2}+l_{1}+l_{2}+d+e_{1}+e_{2}+u_{1}+u_{2}+1}x_{1}^{l_{1}+u_{1}}x_{2}^{l_{2}+u_{2}+1}
\\
&&+\sum_{a,b_{1},b_{2},d,e_{1},e_{2}=0}^{1}\sum_{l_{1}=0}^{b_{1}}%
\sum_{l_{2}=0}^{b_{2}}\sum_{u_{1}=0}^{e_{1}}\sum_{u_{2}=0}^{e_{2}}\left(
-1\right) ^{\alpha \left( g;l_{1},l_{2},u_{1},u_{2}\right) } \\
&&B(gx_{1}x_{2}\otimes
1_{H};G^{a}X_{1}^{b_{1}}X_{2}^{b_{2}},g^{d}x_{1}^{e_{1}}x_{2}^{e_{2}}) \\
&&G^{a}X_{1}^{b_{1}-l_{1}}X_{2}^{b_{2}-l_{2}}\otimes
g^{d}x_{1}^{e_{1}-u_{1}}x_{2}^{e_{2}-u_{2}}\otimes
g^{a+b_{1}+b_{2}+l_{1}+l_{2}+d+e_{1}+e_{2}+u_{1}+u_{2}+1}x_{1}^{l_{1}+u_{1}}x_{2}^{l_{2}+u_{2}}
\\
&=&0
\end{eqnarray*}%
We consider the first summand.%
\begin{eqnarray*}
&&\sum_{a,b_{1},b_{2},d,e_{1},e_{2}=0}^{1}\sum_{l_{1}=0}^{b_{1}}%
\sum_{l_{2}=0}^{b_{2}}\sum_{u_{1}=0}^{e_{1}}\sum_{u_{2}=0}^{e_{2}}\left(
-1\right) ^{\alpha \left( gx_{1}x_{2};l_{1},l_{2},u_{1},u_{2}\right) } \\
&&B(g\otimes
1_{H};G^{a}X_{1}^{b_{1}}X_{2}^{b_{2}},g^{d}x_{1}^{e_{1}}x_{2}^{e_{2}}) \\
&&G^{a}X_{1}^{b_{1}-l_{1}}X_{2}^{b_{2}-l_{2}}\otimes
g^{d}x_{1}^{e_{1}-u_{1}}x_{2}^{e_{2}-u_{2}}\otimes
g^{a+b_{1}+b_{2}+l_{1}+l_{2}+d+e_{1}+e_{2}+u_{1}+u_{2}+1}x_{1}^{l_{1}+u_{1}+1}x_{2}^{l_{2}+u_{2}+1}
\end{eqnarray*}%
We get%
\begin{eqnarray*}
a+b_{1}+b_{2}+l_{1}+l_{2}+d+e_{1}+e_{2}+u_{1}+u_{2}+1 &\equiv &1 \\
l_{1}+u_{1}+1 &=&1 \\
l_{2}+u_{2}+1 &=&1
\end{eqnarray*}%
i.e.%
\begin{eqnarray*}
a+b_{1}+b_{2}+d+e_{1}+e_{2} &\equiv &0 \\
l_{1} &=&u_{1}=0 \\
l_{2} &=&u_{2}=0
\end{eqnarray*}%
Since, by $\left( \ref{got1 first}\right) ,B(g\otimes
1_{H};G^{a}X_{1}^{b_{1}}X_{2}^{b_{2}},g^{d}x_{1}^{e_{1}}x_{2}^{e_{2}})=0$
whenever $a+b_{1}+b_{2}+d+e_{1}+e_{2}\equiv 0$ this summand is zero.

We consider the second summand%
\begin{eqnarray*}
+
&&\sum_{a,b_{1},b_{2},d,e_{1},e_{2}=0}^{1}\sum_{l_{1}=0}^{b_{1}}%
\sum_{l_{2}=0}^{b_{2}}\sum_{u_{1}=0}^{e_{1}}\sum_{u_{2}=0}^{e_{2}}\left(
-1\right) ^{\alpha \left( gx_{1};l_{1},l_{2},u_{1},u_{2}\right) } \\
&&B(x_{2}\otimes
1_{H};G^{a}X_{1}^{b_{1}}X_{2}^{b_{2}},g^{d}x_{1}^{e_{1}}x_{2}^{e_{2}}) \\
&&G^{a}X_{1}^{b_{1}-l_{1}}X_{2}^{b_{2}-l_{2}}\otimes
g^{d}x_{1}^{e_{1}-u_{1}}x_{2}^{e_{2}-u_{2}}\otimes
g^{a+b_{1}+b_{2}+l_{1}+l_{2}+d+e_{1}+e_{2}+u_{1}+u_{2}+1}x_{1}^{l_{1}+u_{1}+1}x_{2}^{l_{2}+u_{2}}
\end{eqnarray*}%
We get%
\begin{eqnarray*}
a+b_{1}+b_{2}+l_{1}+l_{2}+d+e_{1}+e_{2}+u_{1}+u_{2}+1 &\equiv &1 \\
l_{1}+u_{1}+1 &=&1 \\
l_{2}+u_{2} &=&1
\end{eqnarray*}%
i.e.%
\begin{eqnarray*}
a+b_{1}+b_{2}+d+e_{1}+e_{2} &\equiv &1 \\
l_{1} &=&u_{1}=0 \\
l_{2}+u_{2} &=&1
\end{eqnarray*}%
Since, by $\left( \ref{x1ot1 first}\right) $ we know that $B(x_{1}\otimes
1_{H};G^{a}X_{1}^{b_{1}}X_{2}^{b_{2}},g^{d}x_{1}^{e_{1}}x_{2}^{e_{2}})=0$
whenever $a+b_{1}+b_{2}+d+e_{1}+e_{2}\equiv 1,$ this summand is zero.

Let us consider the third summand%
\begin{eqnarray*}
&&-\sum_{a,b_{1},b_{2},d,e_{1},e_{2}=0}^{1}\sum_{l_{1}=0}^{b_{1}}%
\sum_{l_{2}=0}^{b_{2}}\sum_{u_{1}=0}^{e_{1}}\sum_{u_{2}=0}^{e_{2}}\left(
-1\right) ^{\alpha \left( gx_{2};l_{1},l_{2},u_{1},u_{2}\right) } \\
&&B(x_{1}\otimes
1_{H};G^{a}X_{1}^{b_{1}}X_{2}^{b_{2}},g^{d}x_{1}^{e_{1}}x_{2}^{e_{2}}) \\
&&G^{a}X_{1}^{b_{1}-l_{1}}X_{2}^{b_{2}-l_{2}}\otimes
g^{d}x_{1}^{e_{1}-u_{1}}x_{2}^{e_{2}-u_{2}}\otimes
g^{a+b_{1}+b_{2}+l_{1}+l_{2}+d+e_{1}+e_{2}+u_{1}+u_{2}+1}x_{1}^{l_{1}+u_{1}}x_{2}^{l_{2}+u_{2}+1}
\end{eqnarray*}%
We get%
\begin{eqnarray*}
a+b_{1}+b_{2}+l_{1}+l_{2}+d+e_{1}+e_{2}+u_{1}+u_{2}+1 &\equiv &1 \\
l_{1}+u_{1} &=&1 \\
l_{2}+u_{2}+1 &=&1
\end{eqnarray*}%
i.e.%
\begin{eqnarray*}
a+b_{1}+b_{2}+d+e_{1}+e_{2} &\equiv &1 \\
l_{1}+u_{1} &=&1 \\
l_{2} &=&u_{2}=0
\end{eqnarray*}%
Since, by $\left( \ref{x1ot1 first}\right) ,$ $B(x_{1}\otimes
1_{H};G^{a}X_{1}^{b_{1}}X_{2}^{b_{2}},g^{d}x_{1}^{e_{1}}x_{2}^{e_{2}})=0$
whenever $a+b_{1}+b_{2}+d+e_{1}+e_{2}\equiv 1,$ this summand is zero. We now
consider the last summand.%
\begin{eqnarray*}
&&+\sum_{a,b_{1},b_{2},d,e_{1},e_{2}=0}^{1}\sum_{l_{1}=0}^{b_{1}}%
\sum_{l_{2}=0}^{b_{2}}\sum_{u_{1}=0}^{e_{1}}\sum_{u_{2}=0}^{e_{2}}\left(
-1\right) ^{\alpha \left( g;l_{1},l_{2},u_{1},u_{2}\right) } \\
&&B(gx_{1}x_{2}\otimes
1_{H};G^{a}X_{1}^{b_{1}}X_{2}^{b_{2}},g^{d}x_{1}^{e_{1}}x_{2}^{e_{2}}) \\
&&G^{a}X_{1}^{b_{1}-l_{1}}X_{2}^{b_{2}-l_{2}}\otimes
g^{d}x_{1}^{e_{1}-u_{1}}x_{2}^{e_{2}-u_{2}}\otimes
g^{a+b_{1}+b_{2}+l_{1}+l_{2}+d+e_{1}+e_{2}+u_{1}+u_{2}+1}x_{1}^{l_{1}+u_{1}}x_{2}^{l_{2}+u_{2}}
\end{eqnarray*}%
We get%
\begin{eqnarray*}
a+b_{1}+b_{2}+l_{1}+l_{2}+d+e_{1}+e_{2}+u_{1}+u_{2}+1 &\equiv &1 \\
l_{1}+u_{1} &=&1 \\
l_{2}+u_{2} &=&1
\end{eqnarray*}%
i.e.%
\begin{eqnarray*}
a+b_{1}+b_{2}+d+e_{1}+e_{2} &\equiv &0 \\
l_{1}+u_{1} &=&1 \\
l_{2}+u_{2} &=&1
\end{eqnarray*}%
Since, by $\left( \ref{gx1x2ot1, first}\right) ,$ $B(gx_{1}x_{2}\otimes
1_{H};G^{a}X_{1}^{b_{1}}X_{2}^{b_{2}},g^{d}x_{1}^{e_{1}}x_{2}^{e_{2}})=0$
whenever $a+b_{1}+b_{2}+d+e_{1}+e_{2}\equiv 0,$ even this summand is zero.

By using equalities from $\left( \ref{gx1x2ot1, first}\right) $ to $\left( %
\ref{gx1x2ot1, fortyseven}\right) $ we obtain the following form of $B\left(
gx_{1}x_{2}\otimes 1_{H}\right) .$

\subsection{The final form of the element $B(gx_{1}x_{2}\otimes 1_{H})$}

\begin{eqnarray}
B\left( gx_{1}x_{2}\otimes 1_{H}\right) &=&B(gx_{1}x_{2}\otimes
1_{H};1_{A},g)1_{A}\otimes g+  \label{gx1x2} \\
&&+B(gx_{1}x_{2}\otimes 1_{H};1_{A},x_{1})1_{A}\otimes x_{1}  \notag \\
&&+B(gx_{1}x_{2}\otimes 1_{H};1_{A},x_{2})1_{A}\otimes x_{2}  \notag \\
&&+B(gx_{1}x_{2}\otimes 1_{H};1_{A},gx_{1}x_{2})1_{A}\otimes gx_{1}x_{2}
\notag \\
&&+B(gx_{1}x_{2}\otimes 1_{H};G,1_{H})G\otimes 1_{H}+  \notag \\
&&+B(gx_{1}x_{2}\otimes 1_{H};G,x_{1}x_{2})G\otimes x_{1}x_{2}+  \notag \\
&&+B(gx_{1}x_{2}\otimes 1_{H};G,gx_{1})G\otimes gx_{1}+  \notag \\
&&+B(gx_{1}x_{2}\otimes 1_{H};G,gx_{2})G\otimes gx_{2}+  \notag \\
&&+\left[ -B(x_{2}\otimes 1_{H};1_{A},1_{H})+B(gx_{1}x_{2}\otimes
1_{H};1_{A},x_{1})\right] X_{1}\otimes 1_{H}+  \notag \\
&&-B(x_{2}\otimes 1_{H};1_{A},x_{1}x_{2})X_{1}\otimes x_{1}x_{2}+  \notag \\
&&-B(x_{2}\otimes 1_{H};1_{A},gx_{1})X_{1}\otimes gx_{1}+  \notag \\
&&+\left[ -B(x_{2}\otimes 1_{H};1_{A},gx_{2})-B(gx_{1}x_{2}\otimes
1_{H};1_{A},gx_{1}x_{2})\right] X_{1}\otimes gx_{2}+  \notag \\
&&+\left[ B(x_{1}\otimes 1_{H};1_{A},1_{H})+B(gx_{1}x_{2}\otimes
1_{H};1_{A},x_{2})\right] X_{2}\otimes 1_{H}+B(x_{1}\otimes
1_{H};1_{A},x_{1}x_{2})X_{2}\otimes x_{1}x_{2}+  \notag \\
&&+\left[ B(x_{1}\otimes 1_{H};1_{A},gx_{1})+B(gx_{1}x_{2}\otimes
1_{H};1_{A},gx_{1}x_{2})\right] X_{2}\otimes gx_{1}+  \notag \\
&&+B(x_{1}\otimes 1_{H};1_{A},gx_{2})X_{2}\otimes gx_{2}+  \notag \\
&&+\left[
\begin{array}{c}
B(g\otimes 1_{H};1_{A},g)+B(x_{2}\otimes \ 1_{H};1_{A},gx_{2})+ \\
+B(x_{1}\otimes 1_{H};1_{A},gx_{1})+B(gx_{1}x_{2}\otimes
1_{H};1_{A},gx_{1}x_{2})%
\end{array}%
\right] X_{1}X_{2}\otimes g+  \notag \\
&&+\left[ B(g\otimes 1_{H};1_{A},x_{1})+B(x_{2}\otimes \
1_{H};1_{A},x_{1}x_{2})\right] X_{1}X_{2}\otimes x_{1}+  \notag \\
&&+\left[ B(g\otimes 1_{H};X_{2},1_{H})-B(x_{1}\otimes
1_{H};1_{A},x_{1}x_{2})\right] X_{1}X_{2}\otimes x_{2}+  \notag \\
&&+B\left( g\otimes 1_{H};1_{A},gx_{1}x_{2}\right) X_{1}X_{2}\otimes
gx_{1}x_{2}  \notag \\
&&+\left[ B(x_{2}\otimes 1_{H};G,g)-B(gx_{1}x_{2}\otimes 1_{H};G,gx_{1})%
\right] GX_{1}\otimes g+  \notag \\
&&+B(x_{2}\otimes 1_{H};G,x_{1})GX_{1}\otimes x_{1}+  \notag \\
&&+\left[ B(x_{2}\otimes 1_{H};G,x_{2})+B(gx_{1}x_{2}\otimes
1_{H};G,x_{1}x_{2})\right] GX_{1}\otimes x_{2}+  \notag \\
&&+B(x_{2}\otimes 1_{H};G,gx_{1}x_{2})GX_{1}\otimes gx_{1}x_{2}+  \notag \\
&&+\left[ -B(x_{1}\otimes 1_{H};G,g)-B(gx_{1}x_{2}\otimes 1_{H};G,gx_{2})%
\right] GX_{2}\otimes g+  \notag \\
&&+\left[ -B(x_{1}\otimes 1_{H};G,x_{1})-B(gx_{1}x_{2}\otimes
1_{H};G,x_{1}x_{2})\right] GX_{2}\otimes x_{1}+  \notag \\
&&-B(x_{1}\otimes 1_{H};G,x_{2})GX_{2}\otimes x_{2}-B(x_{1}\otimes
1_{H};G,gx_{1}x_{2})GX_{2}\otimes gx_{1}x_{2}+  \notag \\
&&+\left[
\begin{array}{c}
B(g\otimes 1_{H};G,1_{H})+B(x_{2}\otimes \ 1_{H};G,x_{2})+ \\
+B(x_{1}\otimes 1_{H};G,x_{1})+B(gx_{1}x_{2}\otimes 1_{H};G,x_{1}x_{2})%
\end{array}%
\right] GX_{1}X_{2}\otimes 1_{H}  \notag \\
&&+B(g\otimes 1_{H};G,x_{1}x_{2})GX_{1}X_{2}\otimes x_{1}x_{2}+  \notag \\
&&\left[ B(g\otimes 1_{H};G,gx_{1})+B(x_{2}\otimes \ 1_{H};G,gx_{1}x_{2})%
\right] GX_{1}X_{2}\otimes gx_{1}+  \notag \\
&&+\left[ B(g\otimes 1_{H};G,gx_{2})-B(x_{1}\otimes 1_{H};G,gx_{1}x_{2})%
\right] GX_{1}X_{2}\otimes gx_{2}+  \notag
\end{eqnarray}

\section{$B\left( 1\otimes x_{1}\right) $ and $B\left( 1\otimes x_{2}\right)
$}

From $\left( \ref{simplx}\right) $ we get%
\begin{equation*}
B(1_{H}\otimes x_{i})=(1_{A}\otimes x_{i})-(1_{A}\otimes
gx_{i})(1_{A}\otimes g)-(1_{A}\otimes g)B(gx_{i}\otimes 1_{H})(1_{A}\otimes
g)
\end{equation*}%
\begin{equation*}
B(1_{H}\otimes x_{i})=(1_{A}\otimes x_{i})-(1_{A}\otimes
x_{i})-(1_{H}\otimes g)B(gx_{i}\otimes 1_{H})(1_{A}\otimes g)
\end{equation*}%
so that%
\begin{equation}
B(1_{H}\otimes x_{i})=-(1_{A}\otimes g)B(gx_{i}\otimes 1_{H})(1_{A}\otimes
g).  \label{form:1otxi}
\end{equation}%
In view of Proposition \ref{Pro:gh}, we do not have to check Casimir
condition for both $B(1_{H}\otimes x_{1})$ and $B(1_{H}\otimes x_{2}).$

We get%
\begin{eqnarray*}
B\left( 1\otimes x_{1}\right) &=&\left[ -B(x_{1}x_{2}\otimes
1_{H};1_{A},gx_{2})-B(x_{1}x_{2}\otimes 1_{H};X_{2},g)\right] 1_{A}\otimes g+
\\
&&+\left[ B(x_{1}x_{2}\otimes 1_{H};1_{A},x_{1}x_{2})+B(x_{1}x_{2}\otimes
1_{H};X_{2},x_{1})\right] 1_{A}\otimes x_{1} \\
&&+B(x_{1}x_{2}\otimes 1_{H};X_{2},x_{2})1_{A}\otimes x_{2} \\
&&-B(x_{1}x_{2}\otimes 1_{H};X_{2},gx_{1}x_{2})1_{A}\otimes gx_{1}x_{2} \\
&&\left[ -B(x_{1}x_{2}\otimes 1_{H};G,x_{2})+B(x_{1}x_{2}\otimes
1_{H};GX_{2},1_{H})\right] G\otimes 1_{H}+ \\
&&+B(x_{1}x_{2}\otimes 1_{H};GX_{2},x_{1}x_{2})G\otimes x_{1}x_{2}+ \\
&&+\left[ B(x_{1}x_{2}\otimes 1_{H};G,gx_{1}x_{2})-B(x_{1}x_{2}\otimes
1_{H};GX_{2},gx_{1})\right] G\otimes gx_{1}+ \\
&&-B(x_{1}x_{2}\otimes 1_{H};GX_{2},gx_{2})G\otimes gx_{2}+ \\
&&+\left[ +1-B(x_{1}x_{2}\otimes 1_{H};1_{A},x_{1}x_{2})-B(x_{1}x_{2}\otimes
1_{H};X_{2},x_{1})\right] X_{1}\otimes 1_{H}+ \\
&&-B(x_{1}x_{2}\otimes 1_{H};X_{2},gx_{1}x_{2})X_{1}\otimes gx_{2}+ \\
&&-B(x_{1}x_{2}\otimes 1_{H};X_{2},x_{2})X_{2}\otimes 1_{H}+ \\
&&+B(x_{1}x_{2}\otimes 1_{H};X_{2},gx_{1}x_{2})X_{2}\otimes gx_{1}+ \\
&&-B(x_{1}x_{2}\otimes 1_{H};X_{2},gx_{1}x_{2})X_{1}X_{2}\otimes g+ \\
&&+\left[ +B(x_{1}x_{2}\otimes 1_{H};G,gx_{1}x_{2})-B(x_{1}x_{2}\otimes
1_{H};GX_{2},gx_{1})\right] GX_{1}\otimes g+ \\
&&-B(x_{1}x_{2}\otimes 1_{H};GX_{2},x_{1}x_{2})GX_{1}\otimes x_{2}+ \\
&&-B(x_{1}x_{2}\otimes 1_{H};GX_{2},gx_{2})GX_{2}\otimes g+ \\
&&+B(x_{1}x_{2}\otimes 1_{H};GX_{2},x_{1}x_{2})GX_{2}\otimes x_{1}+ \\
&&+B(x_{1}x_{2}\otimes 1_{H};GX_{2},x_{1}x_{2})GX_{1}X_{2}\otimes 1_{H}
\end{eqnarray*}%
and%
\begin{eqnarray*}
B(1_{H}\otimes x_{2}) &=&\left[ +B(x_{1}x_{2}\otimes
1_{H};1_{A},gx_{1})+B(x_{1}x_{2}\otimes 1_{H};X_{1},g)\right] 1_{A}\otimes g+
\\
&&-B(x_{1}x_{2}\otimes 1_{H};X_{1},x_{1})1_{A}\otimes x_{1} \\
&&+\left[ B(x_{1}x_{2}\otimes 1_{H};1_{A},x_{1}x_{2})-B(x_{1}x_{2}\otimes
1_{H};X_{1},x_{2})\right] 1_{A}\otimes x_{2} \\
&&+B(x_{1}x_{2}\otimes 1_{H};X_{1},gx_{1}x_{2})1_{A}\otimes gx_{1}x_{2} \\
&&-\left[ -B(x_{1}x_{2}\otimes 1_{H};G,x_{1})+B(x_{1}x_{2}\otimes
1_{H};GX_{1},1_{H})\right] G\otimes 1_{H}+ \\
&&-B(x_{1}x_{2}\otimes 1_{H};GX_{1},x_{1}x_{2})G\otimes x_{1}x_{2}+ \\
&&+B(x_{1}x_{2}\otimes 1_{H};GX_{1},gx_{1})G\otimes gx_{1}+ \\
&&+\left[ B(x_{1}x_{2}\otimes 1_{H};G,gx_{1}x_{2})+B(x_{1}x_{2}\otimes
1_{H};GX_{1},gx_{2})\right] G\otimes gx_{2}+ \\
&&+B(x_{1}x_{2}\otimes 1_{H};X_{1},x_{1})X_{1}\otimes 1_{H}+ \\
&&+B(x_{1}x_{2}\otimes 1_{H};X_{1},gx_{1}x_{2})X_{1}\otimes gx_{2}+ \\
&&-\left[ -1+B(x_{1}x_{2}\otimes 1_{H};1_{A},x_{1}x_{2})-B(x_{1}x_{2}\otimes
1_{H};X_{1},x_{2})\right] X_{2}\otimes 1_{H}+ \\
&&-B(x_{1}x_{2}\otimes 1_{H};X_{1},gx_{1}x_{2})X_{2}\otimes gx_{1}+ \\
&&+B(x_{1}x_{2}\otimes 1_{H};X_{1},gx_{1}x_{2})X_{1}X_{2}\otimes g+ \\
&&+B(x_{1}x_{2}\otimes 1_{H};GX_{1},gx_{1})GX_{1}\otimes g+ \\
&&+B(x_{1}x_{2}\otimes 1_{H};GX_{1},x_{1}x_{2})GX_{1}\otimes x_{2}+ \\
&&-\left[ -B(x_{1}x_{2}\otimes 1_{H};G,gx_{1}x_{2})-B(x_{1}x_{2}\otimes
1_{H};GX_{1},gx_{2})\right] GX_{2}\otimes g+ \\
&&-B(x_{1}x_{2}\otimes 1_{H};GX_{1},x_{1}x_{2})GX_{2}\otimes x_{1}+ \\
&&-B(x_{1}x_{2}\otimes 1_{H};GX_{1},x_{1}x_{2})GX_{1}X_{2}\otimes 1_{H}.
\end{eqnarray*}

\section{$B\left( 1\otimes gx_{1}\right) $}

From $\left( \ref{simplgx}\right) $ we get%
\begin{eqnarray}
&&B(1_{H}\otimes gx_{1})  \label{form 1otgx1} \\
&=&(1_{A}\otimes g)B(g\otimes 1_{H})(1_{A}\otimes gx_{1})  \notag \\
&&+(1_{A}\otimes x_{1})B(g\otimes 1_{H})+B(x_{1}\otimes 1_{H})  \notag
\end{eqnarray}%
and we get%
\begin{eqnarray*}
&&B(1_{H}\otimes gx_{1}) \\
&=&B\left( x_{1}\otimes 1_{H};1_{A},1_{H}\right) 1_{A}\otimes 1_{H}+ \\
&&+B\left( x_{1}\otimes 1_{H};1_{A},x_{1}x_{2}\right) 1_{A}\otimes
x_{1}x_{2}+ \\
&&+\left[ 2B\left( g\otimes 1_{H};1_{A},g\right) +B\left( x_{1}\otimes
1_{H};1_{A},gx_{1}\right) \right] 1_{A}\otimes gx_{1}+ \\
&&+B\left( x_{1}\otimes 1_{H};1_{A},gx_{2}\right) 1_{A}\otimes gx_{2}+ \\
&&+B\left( x_{1}\otimes 1_{H};G,g\right) G\otimes g+ \\
&&+B\left( x_{1}\otimes 1_{H};G,x_{1}\right) G\otimes x_{1}+ \\
&&+B\left( x_{1}\otimes 1_{H};G,x_{2}\right) G\otimes x_{2}+ \\
&&+\left[ -2B\left( g\otimes 1_{H};G,gx_{2}\right) +B\left( x_{1}\otimes
1_{H};G,gx_{1}x_{2}\right) \right] G\otimes gx_{1}x_{2} \\
&&+\left[ -B(g\otimes 1_{H};1_{A},g)-B(x_{1}\otimes 1_{H};1_{A},gx_{1})%
\right] X_{1}\otimes g+ \\
&&-B(g\otimes 1_{H};X_{1},1_{H})X_{1}\otimes x_{1}+ \\
&&+\left[ -B(g\otimes 1_{H};X_{2},1_{H})+B(x_{1}\otimes
1_{H};1_{A},x_{1}x_{2})\right] X_{1}\otimes x_{2}+ \\
&&+B\left( g\otimes 1_{H};1_{A},gx_{1}x_{2}\right) X_{1}\otimes gx_{1}x_{2}
\\
&&-B(x_{1}\otimes 1_{H};1_{A},gx_{2})X_{2}\otimes g+ \\
&&-B(x_{1}\otimes 1_{H};1_{A},x_{1}x_{2})X_{2}\otimes x_{1}+ \\
&&+\left[ -B(g\otimes 1_{H};1_{A},x_{2})+B(x_{1}\otimes
1_{H};1_{A},x_{1}x_{2})\right] X_{1}X_{2}\otimes 1_{H} \\
&&+B\left( g\otimes 1_{H};1_{A},gx_{1}x_{2}\right) X_{1}X_{2}\otimes gx_{1}
\\
&&+\left[ B(g\otimes 1_{H};G,1_{H})+B(x_{1}\otimes 1_{H};G,x_{1})\right]
GX_{1}\otimes 1_{H}+ \\
&&B(g\otimes 1_{H};G,x_{1}x_{2})GX_{1}\otimes x_{1}x_{2}+ \\
&&\left[ -B\left( g\otimes 1_{H};G,gx_{1}\right) \right] GX_{1}\otimes gx_{1}
\\
&&\left[ B(g\otimes 1_{H};G,gx_{2})-B(x_{1}\otimes 1_{H};G,gx_{1}x_{2})%
\right] GX_{1}\otimes gx_{2}+ \\
&&+B(x_{1}\otimes 1_{H};G,x_{2})GX_{2}\otimes 1_{H}+ \\
&&+\left[ -2B\left( g\otimes 1_{H};G,gx_{2}\right) +B(x_{1}\otimes
1_{H};G,gx_{1}x_{2})\right] GX_{2}\otimes gx_{1}+ \\
&&+\left[ -B(g\otimes 1_{H};G,gx_{2})+B(x_{1}\otimes 1_{H};G,gx_{1}x_{2})%
\right] GX_{1}X_{2}\otimes g \\
&&-B\left( g\otimes 1_{H};G,x_{1}x_{2}\right) GX_{1}X_{2}\otimes x_{1}
\end{eqnarray*}%
We write Casimir condition $\left( \ref{MAIN FORMULA 1}\right) $ for $%
B(1_{H}\otimes gx_{1})$
\begin{eqnarray*}
&&\sum_{a,b_{1},b_{2},d,e_{1},e_{2}=0}^{1}\sum_{l_{1}=0}^{b_{1}}%
\sum_{l_{2}=0}^{b_{2}}\sum_{u_{1}=0}^{e_{1}}\sum_{u_{2}=0}^{e_{2}}\left(
-1\right) ^{\alpha \left( 1_{H};l_{1},l_{2},u_{1},u_{2}\right) } \\
&&B(1_{H}\otimes
gx_{1};G^{a}X_{1}^{b_{1}}X_{2}^{b_{2}},g^{d}x_{1}^{e_{1}}x_{2}^{e_{2}}) \\
&&G^{a}X_{1}^{b_{1}-l_{1}}X_{2}^{b_{2}-l_{2}}\otimes
g^{d}x_{1}^{e_{1}-u_{1}}x_{2}^{e_{2}-u_{2}}\otimes \\
&&g^{a+b_{1}+b_{2}+l_{1}+l_{2}+d+e_{1}+e_{2}+u_{1}+u_{2}}x_{1}^{l_{1}+u_{1}}x_{2}^{l_{2}+u_{2}}
\\
&=&\sum_{\omega _{1}=0}^{1}B^{A}(1\otimes gx_{1}^{1-\omega _{1}})\otimes
B^{H}(1\otimes gx_{1}^{1-\omega _{1}})\otimes g^{\omega _{1}}x_{1}^{\omega
_{1}} \\
&=&B^{A}(1_{H}\otimes gx_{1})\otimes B^{H}(1_{H}\otimes gx_{1})\otimes 1_{H}+
\\
+ &&B^{A}(1_{H}\otimes g)\otimes B^{H}(1_{H}\otimes g)\otimes gx_{1}
\end{eqnarray*}

\subsection{case $1_{H}$}

We get

\begin{eqnarray*}
a+b_{1}+b_{2}+l_{1}+l_{2}+d+e_{1}+e_{2}+u_{1}+u_{2} &\equiv &0 \\
l_{1}+u_{1} &\equiv &0 \\
l_{2}+u_{2} &\equiv &0
\end{eqnarray*}%
i.e.%
\begin{eqnarray*}
a+b_{1}+b_{2}+d+e_{1}+e_{2} &\equiv &0 \\
l_{1} &=&u_{1}=0 \\
l_{2} &=&u_{2}\equiv 0
\end{eqnarray*}%
and we obtain

\begin{eqnarray*}
&&\sum_{\substack{ a,b_{1},b_{2},d,e_{1},e_{2}=0  \\ %
a+b_{1}+b_{2}+d+e_{1}+e_{2}\equiv 0}}^{1}\left( -1\right) ^{\alpha \left(
1_{H};0,0,0,0\right) }B(1_{H}\otimes
gx_{1};G^{a}X_{1}^{b_{1}}X_{2}^{b_{2}},g^{d}x_{1}^{e_{1}}x_{2}^{e_{2}})G^{a}X_{1}^{b_{1}}X_{2}^{b_{2}}\otimes g^{d}x_{1}^{e_{1}}x_{2}^{e_{2}}
\\
&=&=B^{A}(1_{H}\otimes gx_{1})\otimes B^{H}(1_{H}\otimes gx_{1})\otimes 1_{H}
\end{eqnarray*}%
and hence
\begin{equation}
B(1_{H}\otimes
gx_{1};G^{a}X_{1}^{b_{1}}X_{2}^{b_{2}},g^{d}x_{1}^{e_{1}}x_{2}^{e_{2}})=0%
\text{ whenever }a+b_{1}+b_{2}+d+e_{1}+e_{2}\equiv 1  \label{1otgx1, first}
\end{equation}

\subsection{case $g,x_{1}$ and $x_{2}$}

all give us%
\begin{equation*}
a+b_{1}+b_{2}+d+e_{1}+e_{2}\equiv 1
\end{equation*}%
and we get no new information.

\subsection{case $x_{1}x_{2}$}

\begin{eqnarray*}
&&\sum_{a,b_{1},b_{2},d,e_{1},e_{2}=0}^{1}\sum_{l_{1}=0}^{b_{1}}%
\sum_{l_{2}=0}^{b_{2}}\sum_{u_{1}=0}^{e_{1}}\sum_{u_{2}=0}^{e_{2}}\left(
-1\right) ^{\alpha \left( 1_{H};l_{1},l_{2},u_{1},u_{2}\right) } \\
&&B(1_{H}\otimes
gx_{1};G^{a}X_{1}^{b_{1}}X_{2}^{b_{2}},g^{d}x_{1}^{e_{1}}x_{2}^{e_{2}}) \\
&&G^{a}X_{1}^{b_{1}-l_{1}}X_{2}^{b_{2}-l_{2}}\otimes
g^{d}x_{1}^{e_{1}-u_{1}}x_{2}^{e_{2}-u_{2}}\otimes \\
&&g^{a+b_{1}+b_{2}+l_{1}+l_{2}+d+e_{1}+e_{2}+u_{1}+u_{2}}x_{1}^{l_{1}+u_{1}}x_{2}^{l_{2}+u_{2}}
\\
&=&\sum_{\omega _{1}=0}^{1}B^{A}(1\otimes gx_{1}^{1-\omega _{1}})\otimes
B^{H}(1\otimes gx_{1}^{1-\omega _{1}})\otimes g^{\omega _{1}}x_{1}^{\omega
_{1}} \\
&=&B^{A}(1_{H}\otimes gx_{1})\otimes B^{H}(1_{H}\otimes gx_{1})\otimes 1_{H}
\\
&&B^{A}(1_{H}\otimes g)\otimes B^{H}(1_{H}\otimes g)\otimes gx_{1}
\end{eqnarray*}

\begin{eqnarray*}
a+b_{1}+b_{2}+l_{1}+l_{2}+d+e_{1}+e_{2}+u_{1}+u_{2} &\equiv &0 \\
l_{1}+u_{1} &\equiv &1 \\
l_{2}+u_{2} &\equiv &1
\end{eqnarray*}%
\begin{eqnarray*}
a+b_{1}+b_{2}+d+e_{1}+e_{2} &\equiv &0 \\
l_{1}+u_{1} &\equiv &1 \\
l_{2}+u_{2} &\equiv &1
\end{eqnarray*}%
\begin{eqnarray*}
&&\sum_{\substack{ a,b_{1},b_{2},d,e_{1},e_{2}=0  \\ %
a+b_{1}+b_{2}+d+e_{1}+e_{2}\equiv 0}}^{1}\sum_{l_{1}=0}^{b_{1}}%
\sum_{l_{2}=0}^{b_{2}}\sum_{\substack{ u_{1}=0  \\ l_{1}+u_{1}\equiv 1}}%
^{e_{1}}\sum_{\substack{ u_{2}=0  \\ l_{2}+u_{2}\equiv 1}}^{e_{2}}\left(
-1\right) ^{\alpha \left( 1_{H};l_{1},l_{2},u_{1},u_{2}\right) } \\
&&B(1_{H}\otimes
gx_{1};G^{a}X_{1}^{b_{1}}X_{2}^{b_{2}},g^{d}x_{1}^{e_{1}}x_{2}^{e_{2}})G^{a}X_{1}^{b_{1}-l_{1}}X_{2}^{b_{2}-l_{2}}\otimes g^{d}x_{1}^{e_{1}-u_{1}}x_{2}^{e_{2}-u_{2}}\otimes x_{1}x_{2}
\\
&&=0
\end{eqnarray*}%
Since
\begin{eqnarray*}
\alpha \left( 1_{H};0,0,1,1\right) &=&1+e_{2}\text{, }\alpha \left(
1_{H};0,1,1,0\right) =e_{2}+a+b_{1}+b_{2}+1 \\
,\alpha \left( 1_{H};1,0,0,1\right) &=&a+b_{1}\text{and }\alpha \left(
1_{H};1,1,0,0\right) =1+b_{2}
\end{eqnarray*}%
we obtain%
\begin{eqnarray*}
&&\sum_{\substack{ a,b_{1},b_{2},d=0  \\ a+b_{1}+b_{2}+d\equiv 0}}%
^{1}B(1_{H}\otimes
gx_{1};G^{a}X_{1}^{b_{1}}X_{2}^{b_{2}},g^{d}x_{1}x_{2})G^{a}X_{1}^{b_{1}}X_{2}^{b_{2}}\otimes g^{d}
\\
&&+\sum_{\substack{ a,b_{1},d,e_{2}=0  \\ a+b_{1}+d+e_{2}\equiv 0}}%
^{1}\left( -1\right) ^{e_{2}+a+b_{1}}B(1_{H}\otimes
gx_{1};G^{a}X_{1}^{b_{1}}X_{2},g^{d}x_{1}x_{2}^{e_{2}})G^{a}X_{1}^{b_{1}}%
\otimes g^{d}x_{2}^{e_{2}} \\
&&+\sum_{\substack{ a,b_{2},d,e_{1},=0  \\ a+b_{2}+d+e_{1}\equiv 0}}%
^{1}\left( -1\right) ^{a+1}B(1_{H}\otimes
gx_{1};G^{a}X_{1}X_{2}^{b_{2}},g^{d}x_{1}^{e_{1}}x_{2})G^{a}X_{2}^{b_{2}}%
\otimes g^{d}x_{1}^{e_{1}} \\
&&+\sum_{\substack{ a,d,e_{1},e_{2}=0  \\ a+d+e_{1}+e_{2}\equiv 0}}%
^{1}B(1_{H}\otimes
gx_{1};G^{a}X_{1}X_{2},g^{d}x_{1}^{e_{1}}x_{2}^{e_{2}})G^{a}\otimes
g^{d}x_{1}^{e_{1}}x_{2}^{e_{2}} \\
&=&0
\end{eqnarray*}

\subsubsection{ $\left( 0,0,0,0\right) G^{a}\otimes g^{d}$}

We
\begin{equation*}
\sum_{\substack{ a,d=0  \\ a+d\equiv 0}}^{1}\left[
\begin{array}{c}
B(1_{H}\otimes gx_{1};G^{a},g^{d}x_{1}x_{2})+\left( -1\right)
^{a}B(1_{H}\otimes gx_{1};G^{a}X_{2},g^{d}x_{1})+ \\
+\left( -1\right) ^{a+1}B(1_{H}\otimes
gx_{1};G^{a}X_{1},g^{d}x_{2})+B(1_{H}\otimes gx_{1};G^{a}X_{1}X_{2},g^{d})%
\end{array}%
\right] G^{a}\otimes g^{d}
\end{equation*}%
and we get
\begin{eqnarray*}
&&+B(1_{H}\otimes gx_{1};1_{A},x_{1}x_{2})+B(1_{H}\otimes gx_{1};X_{2},x_{1})
\\
&&-B(1_{H}\otimes gx_{1};X_{1},x_{2})+B(1_{H}\otimes
gx_{1};X_{1}X_{2},1_{H})=0
\end{eqnarray*}%
and%
\begin{eqnarray*}
&&+B(1_{H}\otimes gx_{1};G,gx_{1}x_{2})-B(1_{H}\otimes gx_{1};GX_{2},gx_{1})
\\
&&+B(1_{H}\otimes gx_{1};GX_{1},gx_{2})+B(1_{H}\otimes
gx_{1};GX_{1}X_{2},g)=0
\end{eqnarray*}%
By the form of the element $B(1_{H}\otimes gx_{1}),$ both equalities give us
no new information.

\subsubsection{$G^{a}\otimes g^{d}x_{2}$}

\begin{equation*}
\sum_{\substack{ a,d=0  \\ a+d\equiv 1}}^{1}\left[ \left( -1\right)
^{a+1}B(1_{H}\otimes gx_{1};G^{a}X_{2},g^{d}x_{1}x_{2})+B(1_{H}\otimes
gx_{1};G^{a}X_{1}X_{2},g^{d}x_{2})\right] G^{a}\otimes g^{d}x_{2}=0
\end{equation*}%
and we get

\begin{equation*}
-B(1_{H}\otimes gx_{1};X_{2},gx_{1}x_{2})+B(1_{H}\otimes
gx_{1};X_{1}X_{2},gx_{2})=0
\end{equation*}%
and%
\begin{equation*}
B(1_{H}\otimes gx_{1};GX_{2},x_{1}x_{2})+B(1_{H}\otimes
gx_{1};GX_{1}X_{2},x_{2})=0
\end{equation*}%
By the form of the element $B(1_{H}\otimes gx_{1}),$ both equalities give us
no new information.

\subsubsection{$G^{a}\otimes g^{d}x_{1}$}

\begin{equation*}
\sum_{\substack{ a,d=0  \\ a+d\equiv 1}}^{1}\left[ \left( -1\right)
^{a+1}B(1_{H}\otimes gx_{1};G^{a}X_{1},g^{d}x_{1}x_{2})+B(1_{H}\otimes
gx_{1};G^{a}X_{1}X_{2},g^{d}x_{1})\right] G^{a}\otimes g^{d}x_{1}=0
\end{equation*}

and we get%
\begin{equation*}
-B(1_{H}\otimes gx_{1};X_{1},gx_{1}x_{2})+B(1_{H}\otimes
gx_{1};X_{1}X_{2},gx_{1})=0
\end{equation*}%
and%
\begin{equation*}
B(1_{H}\otimes gx_{1};GX_{1},x_{1}x_{2})+B(1_{H}\otimes
gx_{1};GX_{1}X_{2},x_{1})=0.
\end{equation*}%
By the form of the element $B(1_{H}\otimes gx_{1}),$ both equalities give us
no new information.

\subsubsection{$G^{a}X_{2}\otimes g^{d}$}

\begin{equation*}
\sum_{\substack{ a,d=0  \\ a+d\equiv 1}}^{1}\left[ B(1_{H}\otimes
gx_{1};G^{a}X_{2},g^{d}x_{1}x_{2})+\left( -1\right) ^{a+1}B(1_{H}\otimes
gx_{1};G^{a}X_{1}X_{2},g^{d}x_{2})\right] G^{a}X_{2}\otimes g^{d}=0
\end{equation*}%
and we get%
\begin{equation*}
B(1_{H}\otimes gx_{1};X_{2},gx_{1}x_{2})-B(1_{H}\otimes
gx_{1};X_{1}X_{2},gx_{2})=0
\end{equation*}%
and%
\begin{equation*}
B(1_{H}\otimes gx_{1};GX_{2},x_{1}x_{2})+B(1_{H}\otimes
gx_{1};GX_{1}X_{2},x_{2})=0.
\end{equation*}%
By the form of the element $B(1_{H}\otimes gx_{1}),$ both equalities give us
no new information.

\subsubsection{$G^{a}X_{1}\otimes g^{d}$}

\begin{equation*}
\sum_{\substack{ a,d=0  \\ a+d\equiv 1}}^{1}\left[ B(1_{H}\otimes
gx_{1};G^{a}X_{1},g^{d}x_{1}x_{2})+\left( -1\right) ^{a+1}B(1_{H}\otimes
gx_{1};G^{a}X_{1}X_{2},g^{d}x_{1})\right] G^{a}X_{1}\otimes g^{d}=0
\end{equation*}

and we get%
\begin{equation*}
B(1_{H}\otimes gx_{1};X_{1},gx_{1}x_{2})-B(1_{H}\otimes
gx_{1};X_{1}X_{2},gx_{1})=0
\end{equation*}%
\begin{equation*}
B(1_{H}\otimes gx_{1};GX_{1},x_{1}x_{2})+B(1_{H}\otimes
gx_{1};GX_{1}X_{2},x_{1})=0
\end{equation*}%
both of them we already got.

\subsubsection{$G^{a}\otimes g^{d}x_{1}x_{2}$}

\begin{equation*}
\sum_{\substack{ a,d=0  \\ a+d\equiv 0}}^{1}B(1_{H}\otimes
gx_{1};G^{a}X_{1}X_{2},g^{d}x_{1}x_{2})G^{a}\otimes g^{d}x_{1}x_{2}=0
\end{equation*}%
and we get%
\begin{equation*}
B(1_{H}\otimes gx_{1};X_{1}X_{2},x_{1}x_{2})=0
\end{equation*}%
\begin{equation*}
B(1_{H}\otimes gx_{1};GX_{1}X_{2},gx_{1}x_{2})=0.
\end{equation*}%
By the form of the element $B(1_{H}\otimes gx_{1}),$ both equalities give us
no new information.

\subsubsection{$G^{a}X_{2}\otimes g^{d}x_{2}$}

We do not have any summand like this.

\subsubsection{$G^{a}X_{1}\otimes g^{d}x_{2}$}

\begin{equation*}
\sum_{\substack{ a,d=0  \\ a+d\equiv 0}}^{1}\left( -1\right)
^{a}B(1_{H}\otimes gx_{1};G^{a}X_{1}X_{2},g^{d}x_{1}x_{2})G^{a}X_{1}\otimes
g^{d}x_{2}=0
\end{equation*}%
\begin{equation*}
B(1_{H}\otimes gx_{1};X_{1}X_{2},x_{1}x_{2})=0
\end{equation*}%
\begin{equation*}
B(1_{H}\otimes gx_{1};GX_{1}X_{2},gx_{1}x_{2})=0
\end{equation*}%
already got

\subsubsection{$G^{a}X_{2}\otimes g^{d}x_{1}$}

\begin{equation*}
\sum_{\substack{ a,d=0  \\ a+d\equiv 0}}^{1}\left( -1\right)
^{a+1}B(1_{H}\otimes
gx_{1};G^{a}X_{1}X_{2},g^{d}x_{1}x_{2})G^{a}X_{2}\otimes g^{d}x_{1}=0
\end{equation*}%
\begin{equation*}
B(1_{H}\otimes gx_{1};X_{1}X_{2},x_{1}x_{2})=0
\end{equation*}%
\begin{equation*}
B(1_{H}\otimes gx_{1};GX_{1}X_{2},gx_{1}x_{2})=0
\end{equation*}%
already got

\subsubsection{$G^{a}X_{1}\otimes g^{d}x_{1}$}

We do not have any summand like this.

\subsubsection{$G^{a}X_{1}X_{2}\otimes g^{d}$}

\begin{equation*}
\sum_{\substack{ a,d=0  \\ a+d\equiv 0}}^{1}B(1_{H}\otimes
gx_{1};G^{a}X_{1}X_{2},g^{d}x_{1}x_{2})G^{a}X_{1}X_{2}\otimes g^{d}=0
\end{equation*}%
and we get%
\begin{equation*}
B(1_{H}\otimes gx_{1};X_{1}X_{2},x_{1}x_{2})=0
\end{equation*}%
and%
\begin{equation*}
B(1_{H}\otimes gx_{1};GX_{1}X_{2},gx_{1}x_{2})=0
\end{equation*}%
which we already got

\subsubsection{$G^{a}X_{2}\otimes g^{d}x_{1}x_{2}$}

We do not have any summand like this.

\subsubsection{$G^{a}X_{1}\otimes g^{d}x_{1}x_{2}$}

We do not have any summand like this.

\subsubsection{$G^{a}X_{1}X_{2}\otimes g^{d}x_{2}$}

We do not have any summand like this.

\subsubsection{$G^{a}X_{1}X_{2}\otimes g^{d}x_{1}$}

We do not have any summand like this.

\subsubsection{$G^{a}X_{1}X_{2}\otimes g^{d}x_{1}x_{2}$}

We do not have any summand like this.

\subsection{case $gx_{1}$}

\begin{eqnarray*}
&&\sum_{a,b_{1},b_{2},d,e_{1},e_{2}=0}^{1}\sum_{l_{1}=0}^{b_{1}}%
\sum_{l_{2}=0}^{b_{2}}\sum_{u_{1}=0}^{e_{1}}\sum_{u_{2}=0}^{e_{2}}\left(
-1\right) ^{\alpha \left( 1_{H};l_{1},l_{2},u_{1},u_{2}\right) } \\
&&B(1_{H}\otimes
gx_{1};G^{a}X_{1}^{b_{1}}X_{2}^{b_{2}},g^{d}x_{1}^{e_{1}}x_{2}^{e_{2}}) \\
&&G^{a}X_{1}^{b_{1}-l_{1}}X_{2}^{b_{2}-l_{2}}\otimes
g^{d}x_{1}^{e_{1}-u_{1}}x_{2}^{e_{2}-u_{2}}\otimes \\
&&g^{a+b_{1}+b_{2}+l_{1}+l_{2}+d+e_{1}+e_{2}+u_{1}+u_{2}}x_{1}^{l_{1}+u_{1}}x_{2}^{l_{2}+u_{2}}
\\
&=&\sum_{\omega _{1}=0}^{1}B^{A}(1\otimes gx_{1}^{1-\omega _{1}})\otimes
B^{H}(1\otimes gx_{1}^{1-\omega _{1}})\otimes g^{\omega _{1}}x_{1}^{\omega
_{1}} \\
&=&B^{A}(1_{H}\otimes gx_{1})\otimes B^{H}(1_{H}\otimes gx_{1})\otimes 1_{H}
\\
&&B^{A}(1_{H}\otimes g)\otimes B^{H}(1_{H}\otimes g)\otimes gx_{1}
\end{eqnarray*}

\begin{eqnarray*}
a+b_{1}+b_{2}+l_{1}+l_{2}+d+e_{1}+e_{2}+u_{1}+u_{2} &\equiv &1 \\
l_{1}+u_{1} &\equiv &1 \\
l_{2}+u_{2} &\equiv &0
\end{eqnarray*}%
and we get%
\begin{eqnarray*}
a+b_{1}+b_{2}+d+e_{1}+e_{2} &\equiv &0 \\
l_{1}+u_{1} &\equiv &1 \\
l_{2} &=&u_{2}\equiv 0
\end{eqnarray*}%
\begin{gather*}
\sum_{\substack{ a,b_{1},b_{2},d,e_{1},e_{2}=0  \\ %
a+b_{1}+b_{2}+d+e_{1}+e_{2}\equiv 0}}^{1}\sum_{l_{1}=0}^{b_{1}}\sum
_{\substack{ u_{1}=0  \\ l_{1}+u_{1}\equiv 1}}^{e_{1}}\left( -1\right)
^{\alpha \left( 1_{H};l_{1},0,u_{1},0\right) }B(1_{H}\otimes
gx_{1};G^{a}X_{1}^{b_{1}}X_{2}^{b_{2}},g^{d}x_{1}^{e_{1}}x_{2}^{e_{2}}) \\
G^{a}X_{1}^{b_{1}-l_{1}}X_{2}^{b_{2}}\otimes
g^{d}x_{1}^{e_{1}-u_{1}}x_{2}^{e_{2}}=B^{A}(1_{H}\otimes g)\otimes
B^{H}(1_{H}\otimes g)
\end{gather*}%
Since
\begin{eqnarray*}
\alpha \left( 1_{H};0,0,1,0\right) &=&e_{2}+\left( a+b_{1}+b_{2}\right) \\
\alpha \left( 1_{H};1,0,0,0\right) &=&b_{2}
\end{eqnarray*}%
we get%
\begin{eqnarray*}
&&\sum_{\substack{ a,b_{1},b_{2},d,e_{2}=0  \\ a+b_{1}+b_{2}+d+e_{2}\equiv 1
}}^{1}\left( -1\right) ^{e_{2}+\left( a+b_{1}+b_{2}\right) }B(1_{H}\otimes
gx_{1};G^{a}X_{1}^{b_{1}}X_{2}^{b_{2}},g^{d}x_{1}x_{2}^{e_{2}})G^{a}X_{1}^{b_{1}}X_{2}^{b_{2}}\otimes g^{d}x_{2}^{e_{2}}+
\\
&&+\sum_{\substack{ a,b_{2},d,e_{1},e_{2}=0  \\ a+b_{2}+d+e_{1}+e_{2}\equiv
1 }}^{1}\left( -1\right) ^{b_{2}}B(1_{H}\otimes
gx_{1};G^{a}X_{1}X_{2}^{b_{2}},g^{d}x_{1}^{e_{1}}x_{2}^{e_{2}})G^{a}X_{2}^{b_{2}}\otimes g^{d}x_{1}^{e_{1}}x_{2}^{e_{2}}
\\
&=&B^{A}(1_{H}\otimes g)\otimes B^{H}(1_{H}\otimes g)
\end{eqnarray*}

Since%
\begin{equation*}
B(1_{H}\otimes
g,G^{a}X_{1}^{b_{1}}X_{2}^{b_{2}},g^{d}x_{1}^{e_{1}}x_{2}^{e_{2}})=\pm
B(g1_{H}\otimes
1_{H},G^{a}X_{1}^{b_{1}}X_{2}^{b_{2}},g^{d}x_{1}^{e_{1}}x_{2}^{e_{2}})
\end{equation*}%
we also know that
\begin{equation*}
\sum_{\substack{ a,b_{1},b_{2},d,e_{1},e_{2}=0  \\ %
a+b_{1}+b_{2}+d+e_{1}+e_{2}\equiv 1}}^{1}B(1_{H}\otimes
g,G^{a}X_{1}^{b_{1}}X_{2}^{b_{2}},g^{d}x_{1}^{e_{1}}x_{2}^{e_{2}})G^{a}X_{1}^{b_{1}}X_{2}^{b_{2}}\otimes g^{d}x_{1}^{e_{1}}x_{2}^{e_{2}}=0
\end{equation*}%
i.e.%
\begin{equation*}
B(1_{H}\otimes
g,G^{a}X_{1}^{b_{1}}X_{2}^{b_{2}},g^{d}x_{1}^{e_{1}}x_{2}^{e_{2}})=0\text{
whenever }a+b_{1}+b_{2}+d+e_{1}+e_{2}\equiv 1
\end{equation*}

\subsubsection{$G^{a}\otimes g^{d}$}

\begin{equation*}
\sum_{\substack{ a,d=0  \\ a+d\equiv 1}}^{1}\left[ \left( -1\right)
^{a}B(1_{H}\otimes gx_{1};G^{a},g^{d}x_{1})+B(1_{H}\otimes
gx_{1};G^{a}X_{1},g^{d})-B(1_{H}\otimes g;G^{a},g^{d})\right] G^{a}\otimes
g^{d}=0
\end{equation*}%
and we get%
\begin{equation*}
B(1_{H}\otimes gx_{1};X_{1},g)+B(1_{H}\otimes
gx_{1};1_{A},gx_{1})-B(1_{H}\otimes g;1_{A},g)=0
\end{equation*}

\begin{equation*}
B(1_{H}\otimes gx_{1};GX_{1},1_{H})-B(1_{H}\otimes
gx_{1};G,x_{1})-B(1_{H}\otimes g;G,1_{H})=0.
\end{equation*}%
By the form of the element $B(1_{H}\otimes gx_{1}),$ both equalities give us
no new information.

\subsubsection{$G^{a}\otimes g^{d}x_{2}$}

\begin{gather*}
\sum_{\substack{ a,d=0  \\ a+d\equiv 0}}^{1}\left[ \left( -1\right)
^{a+1}B(1_{H}\otimes gx_{1};G^{a},g^{d}x_{1}x_{2})+B(1_{H}\otimes
gx_{1};G^{a}X_{1},g^{d}x_{2})-B(1_{H}\otimes g;G^{a},g^{d}x_{2})\right] \\
G^{a}\otimes g^{d}x_{2}=0
\end{gather*}%
and we obtain%
\begin{equation*}
B(1_{H}\otimes gx_{1};X_{1},x_{2})-B(1_{H}\otimes
gx_{1};1_{A},x_{1}x_{2})-B(1_{H}\otimes g;1_{A},x_{2})=0
\end{equation*}%
\begin{equation*}
B(1_{H}\otimes gx_{1};GX_{1},gx_{2})+B(1_{H}\otimes
gx_{1};G,gx_{1}x_{2})-B(1_{H}\otimes g;G,gx_{2})=0
\end{equation*}%
By the form of the element $B(1_{H}\otimes gx_{1}),$ both equalities give us
no new information.

\subsubsection{$G^{a}\otimes g^{d}x_{1}$}

\begin{equation*}
\sum_{\substack{ a,d=0  \\ a+d\equiv 0}}^{1}B(1_{H}\otimes
gx_{1};G^{a}X_{1},g^{d}x_{1})G^{a}\otimes g^{d}x_{1}=B(1_{H}\otimes
g;G^{a},g^{d}x_{1})G^{a}\otimes g^{d}x_{1}
\end{equation*}%
and we get%
\begin{equation*}
B(1_{H}\otimes gx_{1};X_{1},x_{1})-B(1_{H}\otimes g;1_{A},x_{1})=0
\end{equation*}%
\begin{equation*}
B(1_{H}\otimes gx_{1};GX_{1},gx_{1})-B(1_{H}\otimes g;G,gx_{1})=0.
\end{equation*}%
By the form of the element $B(1_{H}\otimes gx_{1}),$ both equalities give us
no new information.

\subsubsection{$G^{a}X_{2}\otimes g^{d}$}

\begin{gather*}
\sum_{\substack{ a,d=0  \\ a+d\equiv 0}}^{1}\left[ \left( -1\right)
^{a+1}B(1_{H}\otimes gx_{1};G^{a}X_{2},g^{d}x_{1})-B(1_{H}\otimes
gx_{1};G^{a}X_{1}X_{2},g^{d})-B(1_{H}\otimes g;G^{a}X_{2},g^{d})\right] \\
G^{a}X_{2}\otimes g^{d}=0
\end{gather*}%
and we get%
\begin{equation*}
-B(1_{H}\otimes gx_{1};X_{1}X_{2},1_{H})-B(1_{H}\otimes
gx_{1};X_{2},x_{1})-B(1_{H}\otimes g;X_{2},1_{H})=0
\end{equation*}%
\begin{equation*}
-B(1_{H}\otimes gx_{1};GX_{1}X_{2},g)+B(1_{H}\otimes
gx_{1};GX_{2},gx_{1})-B(1_{H}\otimes g;GX_{2},g)=0.
\end{equation*}%
By the form of the element $B(1_{H}\otimes gx_{1}),$ both equalities give us
no new information.

\subsubsection{$G^{a}X_{1}\otimes g^{d}$}

\begin{equation*}
\sum_{\substack{ a,d=0  \\ a+d\equiv 0}}^{1}\left[ \left( -1\right)
^{a+1}B(1_{H}\otimes gx_{1};G^{a}X_{1},g^{d}x_{1})-B(1_{H}\otimes
g;G^{a}X_{1},g^{d})\right] G^{a}X_{1}\otimes g^{d}=0
\end{equation*}%
and we get%
\begin{equation*}
-B(1_{H}\otimes gx_{1};X_{1},x_{1})-B(1_{H}\otimes g;X_{1},1_{H})=0
\end{equation*}%
\begin{equation*}
B(1_{H}\otimes gx_{1};GX_{1},gx_{1})-B(1_{H}\otimes g;GX_{1},g)=0
\end{equation*}%
which we already got.

\subsubsection{$G^{a}\otimes g^{d}x_{1}x_{2}$}

\begin{equation*}
\sum_{\substack{ a,d=0  \\ a+d\equiv 1}}^{1}\left[ B(1_{H}\otimes
gx_{1};G^{a}X_{1},g^{d}x_{1}x_{2})-B(1_{H}\otimes g;G^{a},g^{d}x_{1}x_{2})%
\right] G^{a}\otimes g^{d}x_{1}x_{2}
\end{equation*}%
and we get%
\begin{equation*}
B(1_{H}\otimes gx_{1};X_{1},gx_{1}x_{2})-B(1_{H}\otimes
g;1_{A},gx_{1}x_{2})=0
\end{equation*}%
\begin{equation*}
B(1_{H}\otimes gx_{1};GX_{1},x_{1}x_{2})-B(1_{H}\otimes g;G,x_{1}x_{2})=0.
\end{equation*}%
By the form of the element $B(1_{H}\otimes gx_{1}),$ both equalities give us
no new information.

\subsubsection{$G^{a}X_{2}\otimes g^{d}x_{2}$}

\begin{gather*}
\sum_{\substack{ a,d=0  \\ a+d\equiv 1}}^{1}\left[ \left( -1\right)
^{a}B(1_{H}\otimes gx_{1};G^{a}X_{2},g^{d}x_{1}x_{2})-B(1_{H}\otimes
gx_{1};G^{a}X_{1}X_{2},g^{d}x_{2})-B(1_{H}\otimes g;G^{a}X_{2},g^{d}x_{2})%
\right] \\
G^{a}X_{2}\otimes g^{d}x_{2}=0
\end{gather*}%
and we get%
\begin{equation*}
-B(1_{H}\otimes gx_{1};X_{1}X_{2},gx_{2})+B(1_{H}\otimes
gx_{1};X_{2},gx_{1}x_{2})-B(1_{H}\otimes g;X_{2}\otimes gx_{2})=0
\end{equation*}%
\begin{equation*}
-B(1_{H}\otimes gx_{1};GX_{1}X_{2},x_{2})-B(1_{H}\otimes
gx_{1};GX_{2},x_{1}x_{2})-B(1_{H}\otimes g;GX_{2}\otimes x_{2})=0
\end{equation*}%
By the form of the element $B(1_{H}\otimes gx_{1}),$ both equalities give us
no new information.

\subsubsection{$G^{a}X_{1}\otimes g^{d}x_{2}$}

\begin{equation*}
\sum_{\substack{ a,d=0  \\ a+d\equiv 1}}^{1}\left[ \left( -1\right)
^{a}B(1_{H}\otimes gx_{1};G^{a}X_{1},g^{d}x_{1}x_{2})-B(1_{H}\otimes
g;G^{a}X_{1},g^{d}x_{2})\right] G^{a}X_{1}\otimes g^{d}x_{2}=0
\end{equation*}%
and we get%
\begin{equation*}
B(1_{H}\otimes gx_{1};X_{1},gx_{1}x_{2})-B(1_{H}\otimes g;X_{1},gx_{2})=0
\end{equation*}%
\begin{equation*}
-B(1_{H}\otimes gx_{1};GX_{1},x_{1}x_{2})-B(1_{H}\otimes g;GX_{1},x_{2})=0
\end{equation*}%
By the form of the element $B(1_{H}\otimes gx_{1}),$ both equalities give us
no new information.

\subsubsection{$G^{a}X_{2}\otimes g^{d}x_{1}$}

\begin{equation*}
\sum_{\substack{ a,d=0  \\ a+d\equiv 1}}^{1}\left[ -B(1_{H}\otimes
gx_{1};G^{a}X_{1}X_{2},g^{d}x_{1})-B(1_{H}\otimes g;G^{a}X_{2},g^{d}x_{1})%
\right] G^{a}X_{2}\otimes g^{d}x_{1}=0
\end{equation*}%
and we get%
\begin{equation*}
-B(1_{H}\otimes gx_{1};X_{1}X_{2},gx_{1})-B(1_{H}\otimes g;X_{2},gx_{1})=0.
\end{equation*}%
\begin{equation*}
-B(1_{H}\otimes gx_{1};GX_{1}X_{2},x_{1})-B(1_{H}\otimes g;GX_{2},x_{1})=0.
\end{equation*}%
By the form of the element $B(1_{H}\otimes gx_{1}),$ both equalities give us
no new information.

\subsubsection{$G^{a}X_{1}\otimes g^{d}x_{1}$}

\begin{equation*}
-\sum_{\substack{ a,d=0  \\ a+d\equiv 1}}^{1}B(1_{H}\otimes
g;G^{a}X_{1},g^{d}x_{1})=0
\end{equation*}%
and we get%
\begin{equation*}
-B(1_{H}\otimes g;X_{1},gx_{1})=0
\end{equation*}%
and%
\begin{equation*}
-B(1_{H}\otimes g;GX_{1},x_{1})=0
\end{equation*}%
By the form of the element $B(g\otimes 1_{H}),$ both equalities give us no
new information.

\subsubsection{$G^{a}X_{1}X_{2}\otimes g^{d}$}

\begin{equation*}
\sum_{\substack{ a,d=0  \\ a+d\equiv 1}}^{1}\left[ \left( -1\right)
^{a}B(1_{H}\otimes gx_{1};G^{a}X_{1}X_{2},g^{d}x_{1})-B(1_{H}\otimes
g;G^{a}X_{1}X_{2},g^{d})\right] G^{a}X_{1}X_{2}\otimes g^{d}=0
\end{equation*}%
and we get%
\begin{equation*}
B(1_{H}\otimes gx_{1};X_{1}X_{2},gx_{1})-B(1_{H}\otimes g;X_{1}X_{2},g)=0
\end{equation*}%
\begin{equation*}
-B(1_{H}\otimes gx_{1};GX_{1}X_{2},x_{1})-B(1_{H}\otimes
g;GX_{1}X_{2},1_{H})=0.
\end{equation*}%
By the form of the element $B(1_{H}\otimes gx_{1}),$ both equalities give us
no new information.

\subsubsection{$G^{a}X_{2}\otimes g^{d}x_{1}x_{2}$}

\begin{equation*}
\sum_{\substack{ a,d=0  \\ a+d\equiv 0}}^{1}\left[ -B(1_{H}\otimes
gx_{1};G^{a}X_{1}X_{2},g^{d}x_{1}x_{2})-B(1_{H}\otimes
g;G^{a}X_{2},g^{d}x_{1}x_{2})\right] G^{a}X_{2}\otimes g^{d}x_{1}x_{2}=0
\end{equation*}

and we get%
\begin{equation*}
-B(1_{H}\otimes gx_{1};X_{1}X_{2},x_{1}x_{2})-B(1_{H}\otimes
g;X_{2},x_{1}x_{2})=0
\end{equation*}%
\begin{equation*}
-B(1_{H}\otimes gx_{1};GX_{1}X_{2},gx_{1}x_{2})-B(1_{H}\otimes
g;GX_{2},gx_{1}x_{2})=0.
\end{equation*}%
By the form of the elements $B(1_{H}\otimes gx_{1})$ and $B(g\otimes 1_{H}),$
both equalities give us no new information.

\subsubsection{$G^{a}X_{1}\otimes g^{d}x_{1}x_{2}$}

\begin{equation*}
\sum_{\substack{ a,d=0  \\ a+d\equiv 0}}^{1}B(1_{H}\otimes
g;G^{a}X_{1},g^{d}x_{1}x_{2})=0
\end{equation*}%
and we get%
\begin{equation*}
-B(1_{H}\otimes g;X_{1},x_{1}x_{2})=0
\end{equation*}%
\begin{equation*}
-B(1_{H}\otimes g;GX_{1},gx_{1}x_{2})=0
\end{equation*}%
By the form of the element $B(g\otimes 1_{H}),$ we do not get any new
information.

\subsubsection{$G^{a}X_{1}X_{2}\otimes g^{d}x_{2}$}

\begin{equation*}
\sum_{\substack{ a,d=0  \\ a+d\equiv 0}}^{1}\left[ \left( -1\right)
^{a+1}B(1_{H}\otimes gx_{1};G^{a}X_{1}X_{2},g^{d}x_{1}x_{2})-B(1_{H}\otimes
g;G^{a}X_{1}X_{2},g^{d}x_{2})\right] G^{a}X_{1}X_{2}\otimes g^{d}x_{2}=0
\end{equation*}%
and we get%
\begin{equation*}
B(1_{H}\otimes gx_{1};GX_{1}X_{2},gx_{1}x_{2})-B(1_{H}\otimes
g;GX_{1}X_{2},gx_{2})=0
\end{equation*}%
\begin{equation*}
-B(1_{H}\otimes gx_{1};X_{1}X_{2},x_{1}x_{2})-B(1_{H}\otimes
g;X_{1}X_{2},x_{2})=0.
\end{equation*}%
By the form of the elements $B(1_{H}\otimes gx_{1})$ and $B(g\otimes 1_{H}),$
both equalities give us no new information.

\subsubsection{$G^{a}X_{1}X_{2}\otimes g^{d}x_{1}$}

\begin{equation*}
\sum_{\substack{ a,d=0  \\ a+d\equiv 0}}^{1}B(1_{H}\otimes
g;G^{a}X_{1}X_{2},g^{d}x_{1})G^{a}X_{1}X_{2}\otimes g^{d}x_{1}=0
\end{equation*}%
and we get%
\begin{equation*}
-B(1_{H}\otimes g;X_{1}X_{2},x_{1})=0
\end{equation*}%
\begin{equation*}
-B(1_{H}\otimes g;GX_{1}X_{2},gx_{1})=0
\end{equation*}%
By the form of the element $B(g\otimes 1_{H}),$ we do not get any new
information.

\subsubsection{$G^{a}X_{1}X_{2}\otimes g^{d}x_{1}x_{2}$}

\begin{equation*}
\sum_{\substack{ a,d=0  \\ a+d\equiv 0}}^{1}B(1_{H}\otimes
g;G^{a}X_{1}X_{2},g^{d}x_{1}x_{2})G^{a}X_{1}X_{2}\otimes g^{d}x_{1}x_{2}=0
\end{equation*}

and we get

\begin{equation*}
-B(1_{H}\otimes g;X_{1}X_{2},x_{1}x_{2})=0
\end{equation*}%
\begin{equation*}
-B(1_{H}\otimes g;GX_{1}X_{2},gx_{1}x_{2})=0
\end{equation*}

By the form of the element $B(g\otimes 1_{H}),$ we do not get any new
information.

\subsection{case $gx_{2}$}

\begin{eqnarray*}
&&\sum_{a,b_{1},b_{2},d,e_{1},e_{2}=0}^{1}\sum_{l_{1}=0}^{b_{1}}%
\sum_{l_{2}=0}^{b_{2}}\sum_{u_{1}=0}^{e_{1}}\sum_{u_{2}=0}^{e_{2}}\left(
-1\right) ^{\alpha \left( 1_{H};l_{1},l_{2},u_{1},u_{2}\right) } \\
&&B(1_{H}\otimes
gx_{1};G^{a}X_{1}^{b_{1}}X_{2}^{b_{2}},g^{d}x_{1}^{e_{1}}x_{2}^{e_{2}}) \\
&&G^{a}X_{1}^{b_{1}-l_{1}}X_{2}^{b_{2}-l_{2}}\otimes
g^{d}x_{1}^{e_{1}-u_{1}}x_{2}^{e_{2}-u_{2}}\otimes \\
&&g^{a+b_{1}+b_{2}+l_{1}+l_{2}+d+e_{1}+e_{2}+u_{1}+u_{2}}x_{1}^{l_{1}+u_{1}}x_{2}^{l_{2}+u_{2}}
\\
&=&0
\end{eqnarray*}

\begin{eqnarray*}
a+b_{1}+b_{2}+l_{1}+l_{2}+d+e_{1}+e_{2}+u_{1}+u_{2} &\equiv &1 \\
l_{1}+u_{1} &\equiv &0 \\
l_{2}+u_{2} &\equiv &1
\end{eqnarray*}%
\begin{eqnarray*}
a+b_{1}+b_{2}+d+e_{1}+e_{2} &\equiv &0 \\
l_{1} &=&u_{1}\equiv 0 \\
l_{2}+u_{2} &\equiv &1
\end{eqnarray*}%
\begin{eqnarray*}
&&\sum_{\substack{ a,b_{1},b_{2},d,e_{1},e_{2}=0  \\ %
a+b_{1}+b_{2}+d+e_{1}+e_{2}\equiv 0}}^{1}\sum_{l_{2}=0}^{b_{2}}\sum
_{\substack{ u_{2}=0  \\ l_{2}+u_{2}\equiv 1}}^{e_{2}}\left( -1\right)
^{\alpha \left( 1_{H};0,l_{2},0,u_{2}\right) } \\
&&B(1_{H}\otimes
gx_{1};G^{a}X_{1}^{b_{1}}X_{2}^{b_{2}},g^{d}x_{1}^{e_{1}}x_{2}^{e_{2}})G^{a}X_{1}^{b_{1}}X_{2}^{b_{2}-l_{2}}\otimes g^{d}x_{1}^{e_{1}}x_{2}^{e_{2}-u_{2}}
\\
&=&0
\end{eqnarray*}%
Since $\alpha \left( 1_{H};0,0,0,1\right) =a+b_{1}+b_{2}$ and $\alpha \left(
1_{H};0,1,0,0\right) =0,$ we obtain%
\begin{eqnarray*}
&&\sum_{\substack{ a,b_{1},b_{2},d,e_{1}=0  \\ a+b_{1}+b_{2}+d+e_{1}\equiv 1
}}^{1}\left( -1\right) ^{a+b_{1}+b_{2}}B(1_{H}\otimes
gx_{1};G^{a}X_{1}^{b_{1}}X_{2}^{b_{2}},g^{d}x_{1}^{e_{1}}x_{2})G^{a}X_{1}^{b_{1}}X_{2}^{b_{2}}\otimes g^{d}x_{1}^{e_{1}}+
\\
&&+\sum_{\substack{ a,b_{1},d,e_{1},e_{2}=0  \\ a+b_{1}+d+e_{1}+e_{2}\equiv
1 }}^{1}B(1_{H}\otimes
gx_{1};G^{a}X_{1}^{b_{1}}X_{2},g^{d}x_{1}^{e_{1}}x_{2}^{e_{2}})G^{a}X_{1}^{b_{1}}\otimes g^{d}x_{1}^{e_{1}}x_{2}^{e_{2}}
\\
&=&0.
\end{eqnarray*}

\subsubsection{$\left( 0,0,0,0\right) G^{a}\otimes g^{d}$}

\begin{equation*}
\sum_{\substack{ a,d=0  \\ a+d\equiv 1}}^{1}\left[ \left( -1\right)
^{a}B(1_{H}\otimes gx_{1};G^{a},g^{d}x_{2})+B(1_{H}\otimes
gx_{1};G^{a}X_{2},g^{d})\right] G^{a}\otimes g^{d}=0
\end{equation*}%
and we get%
\begin{equation*}
B(1_{H}\otimes gx_{1};1_{A},gx_{2})+B(1_{H}\otimes gx_{1};X_{2},g)=0
\end{equation*}%
\begin{equation*}
-B(1_{H}\otimes gx_{1};G,x_{2})+B(1_{H}\otimes gx_{1};GX_{2},1_{H})=0.
\end{equation*}%
By the form of the element $B(1_{H}\otimes gx_{1}),$ both equalities give us
no new information.

\subsubsection{$G^{a}\otimes g^{d}x_{2}$}

\begin{equation*}
\sum_{\substack{ a,d=0  \\ a+d\equiv 0}}^{1}-B(1_{H}\otimes
gx_{1};G^{a}X_{2},g^{d}x_{2})G^{a}\otimes g^{d}x_{2}=0.
\end{equation*}%
and we get%
\begin{equation*}
B(1_{H}\otimes gx_{1};X_{2},x_{2})=0
\end{equation*}%
\begin{equation*}
B(1_{H}\otimes gx_{1};GX_{2},gx_{2})=0.
\end{equation*}%
By the form of the element $B(1_{H}\otimes gx_{1}),$ both equalities give us
no new information.

\subsubsection{$G^{a}\otimes g^{d}x_{1}$}

\begin{equation*}
\sum_{\substack{ a,d=0  \\ a+d\equiv 0}}^{1}\left[ \left( -1\right)
^{a}B(1_{H}\otimes gx_{1};G^{a},g^{d}x_{1}x_{2})+B(1_{H}\otimes
gx_{1};G^{a}X_{2},g^{d}x_{1})\right] G^{a}\otimes g^{d}x_{1}=0
\end{equation*}%
and we get%
\begin{equation*}
B(1_{H}\otimes gx_{1};1_{A},x_{1}x_{2})+B(1_{H}\otimes gx_{1};X_{2},x_{1})=0
\end{equation*}%
\begin{equation*}
-B(1_{H}\otimes gx_{1};G,gx_{1}x_{2})+B(1_{H}\otimes gx_{1};GX_{2},gx_{1})=0.
\end{equation*}%
By the form of the element $B(1_{H}\otimes gx_{1}),$ both equalities give us
no new information.

\subsubsection{$G^{a}X_{2}\otimes g^{d}$}

\begin{equation*}
\sum_{\substack{ a,d=0  \\ a+d\equiv 0}}^{1}\left( -1\right)
^{a+1}B(1_{H}\otimes gx_{1};G^{a}X_{2},g^{d}x_{2})G^{a}X_{2}\otimes g^{d}=0
\end{equation*}%
and we get%
\begin{equation*}
-B(1_{H}\otimes gx_{1};X_{2},x_{2})=0
\end{equation*}%
\begin{equation*}
B(1_{H}\otimes gx_{1};GX_{2},gx_{2})=0.
\end{equation*}%
By the form of the element $B(1_{H}\otimes gx_{1}),$ both equalities give us
no new information.

\subsubsection{$G^{a}X_{1}\otimes g^{d}$}

\begin{equation*}
\sum_{\substack{ a,d=0  \\ a+d\equiv 0}}^{1}\left[ \left( -1\right)
^{a+1}B(1_{H}\otimes gx_{1};G^{a}X_{1},g^{d}x_{2})+B(1_{H}\otimes
gx_{1};G^{a}X_{1}X_{2},g^{d})\right] G^{a}X_{1}\otimes g^{d}=0
\end{equation*}%
and we get%
\begin{equation*}
-B(1_{H}\otimes gx_{1};X_{1},x_{2})+B(1_{H}\otimes gx_{1};X_{1}X_{2},1_{H})=0
\end{equation*}%
\begin{equation*}
B(1_{H}\otimes gx_{1};GX_{1},gx_{2})+B(1_{H}\otimes gx_{1};GX_{1}X_{2},g)=0
\end{equation*}%
By the form of the element $B(1_{H}\otimes gx_{1}),$ both equalities give us
no new information.

\subsubsection{$G^{a}\otimes g^{d}x_{1}x_{2}$}

\begin{equation*}
\sum_{\substack{ a,d=0  \\ a+d\equiv 1}}^{1}B(1_{H}\otimes
gx_{1};G^{a}X_{2},g^{d}x_{1}x_{2})G^{a}\otimes g^{d}x_{1}x_{2}=0.
\end{equation*}%
and we get%
\begin{equation*}
B(1_{H}\otimes gx_{1};X_{2},gx_{1}x_{2})=0
\end{equation*}%
\begin{equation*}
B(1_{H}\otimes gx_{1};GX_{2},x_{1}x_{2})=0.
\end{equation*}%
By the form of the element $B(1_{H}\otimes gx_{1}),$ both equalities give us
no new information.

\subsubsection{$G^{a}X_{2}\otimes g^{d}x_{2}$}

We do not have any term like this

\subsubsection{$G^{a}X_{1}\otimes g^{d}x_{2}$}

\begin{equation*}
\sum_{\substack{ a,d=0  \\ a+d\equiv 1}}^{1}B(1_{H}\otimes
gx_{1};G^{a}X_{1}X_{2},g^{d}x_{2})G^{a}X_{1}\otimes g^{d}x_{2}=0
\end{equation*}%
and we get%
\begin{equation*}
B(1_{H}\otimes gx_{1};X_{1}X_{2},gx_{2})=0
\end{equation*}%
and%
\begin{equation*}
B(1_{H}\otimes gx_{1};GX_{1}X_{2},x_{2})=0.
\end{equation*}%
By the form of the element $B(1_{H}\otimes gx_{1}),$ both equalities give us
no new information.

\subsubsection{$G^{a}X_{2}\otimes g^{d}x_{1}$}

\begin{equation*}
\sum_{\substack{ a,d=0  \\ a+d\equiv 1}}^{1}\left( -1\right)
^{a+1}B(1_{H}\otimes gx_{1};G^{a}X_{2},g^{d}x_{1}x_{2})G^{a}X_{2}\otimes
g^{d}x_{1}=0
\end{equation*}%
and we get%
\begin{equation*}
-B(1_{H}\otimes gx_{1};X_{2},gx_{1}x_{2})=0
\end{equation*}%
and%
\begin{equation*}
B(1_{H}\otimes gx_{1};GX_{2},x_{1}x_{2})=0.
\end{equation*}%
By the form of the element $B(1_{H}\otimes gx_{1}),$ both equalities give us
no new information.

\subsubsection{$G^{a}X_{1}\otimes g^{d}x_{1}$}

\begin{equation*}
\sum_{\substack{ a,d=0  \\ a+d\equiv 1}}^{1}\left[ \left( -1\right)
^{a+1}B(1_{H}\otimes gx_{1};G^{a}X_{1},g^{d}x_{1}x_{2})+B(1_{H}\otimes
gx_{1};G^{a}X_{1}X_{2},g^{d}x_{1})\right] G^{a}X_{1}\otimes g^{d}x_{1}=0
\end{equation*}%
and we get%
\begin{equation*}
-B(1_{H}\otimes gx_{1};X_{1},gx_{1}x_{2})+B(1_{H}\otimes
gx_{1};X_{1}X_{2},gx_{1})=0
\end{equation*}%
\begin{equation*}
B(1_{H}\otimes gx_{1};GX_{1},x_{1}x_{2})+B(1_{H}\otimes
gx_{1};GX_{1}X_{2},x_{1})=0
\end{equation*}%
By the form of the element $B(1_{H}\otimes gx_{1}),$ both equalities give us
no new information.

\subsubsection{$G^{a}X_{1}X_{2}\otimes g^{d}$}

\begin{equation*}
\sum_{\substack{ a,d=0  \\ a+d\equiv 1}}^{1}\left( -1\right)
^{a}B(1_{H}\otimes gx_{1};G^{a}X_{1}X_{2},g^{d}x_{2})G^{a}X_{1}X_{2}\otimes
g^{d}=0
\end{equation*}%
and we get%
\begin{equation*}
B(1_{H}\otimes gx_{1};X_{1}X_{2},gx_{2})=0
\end{equation*}%
and%
\begin{equation*}
-B(1_{H}\otimes gx_{1};GX_{1}X_{2},x_{2})=0.
\end{equation*}%
By the form of the element $B(1_{H}\otimes gx_{1}),$ both equalities give us
no new information.

\subsubsection{$G^{a}X_{2}\otimes g^{d}x_{1}x_{2}$}

We do not have any term like this.

\subsubsection{$G^{a}X_{1}\otimes g^{d}x_{1}x_{2}$}

\begin{equation*}
\sum_{\substack{ a,d=0  \\ a+d\equiv 0}}^{1}B(1_{H}\otimes
gx_{1};G^{a}X_{1}X_{2},g^{d}x_{1}x_{2})G^{a}X_{1}\otimes g^{d}x_{1}x_{2}=0
\end{equation*}%
and we get%
\begin{equation*}
B(1_{H}\otimes gx_{1};X_{1}X_{2},x_{1}x_{2})=0
\end{equation*}%
and%
\begin{equation*}
B(1_{H}\otimes gx_{1};GX_{1}X_{2},gx_{1}x_{2})=0
\end{equation*}%
which we already got.

\subsubsection{$G^{a}X_{1}X_{2}\otimes g^{d}x_{2}$}

We do not have any term like this

\subsubsection{$G^{a}X_{1}X_{2}\otimes g^{d}x_{1}$}

\begin{equation*}
\sum_{\substack{ a,d=0  \\ a+d\equiv 0}}^{1}\left( -1\right)
^{a_{2}}B(1_{H}\otimes
gx_{1};G^{a}X_{1}X_{2},g^{d}x_{1}x_{2})G^{a}X_{1}X_{2}\otimes g^{d}x_{1}=0
\end{equation*}%
and we get%
\begin{equation*}
B(1_{H}\otimes gx_{1};X_{1}X_{2},x_{1}x_{2})=0
\end{equation*}%
and%
\begin{equation*}
B(1_{H}\otimes gx_{1};GX_{1}X_{2},gx_{1}x_{2})=0
\end{equation*}%
which we already got.

\subsubsection{$G^{a}X_{1}X_{2}\otimes g^{d}x_{1}x_{2}$}

We do not have any term like this.

\subsection{case $gx_{1}x_{2}$}

\begin{eqnarray*}
&&\sum_{a,b_{1},b_{2},d,e_{1},e_{2}=0}^{1}\sum_{l_{1}=0}^{b_{1}}%
\sum_{l_{2}=0}^{b_{2}}\sum_{u_{1}=0}^{e_{1}}\sum_{u_{2}=0}^{e_{2}}\left(
-1\right) ^{\alpha \left( 1_{H};l_{1},l_{2},u_{1},u_{2}\right) } \\
&&B(1_{H}\otimes
gx_{1};G^{a}X_{1}^{b_{1}}X_{2}^{b_{2}},g^{d}x_{1}^{e_{1}}x_{2}^{e_{2}}) \\
&&G^{a}X_{1}^{b_{1}-l_{1}}X_{2}^{b_{2}-l_{2}}\otimes
g^{d}x_{1}^{e_{1}-u_{1}}x_{2}^{e_{2}-u_{2}}\otimes \\
&&g^{a+b_{1}+b_{2}+l_{1}+l_{2}+d+e_{1}+e_{2}+u_{1}+u_{2}}x_{1}^{l_{1}+u_{1}}x_{2}^{l_{2}+u_{2}}
\\
&=&\sum_{\omega _{1}=0}^{1}B^{A}(1\otimes gx_{1}^{1-\omega _{1}})\otimes
B^{H}(1\otimes gx_{1}^{1-\omega _{1}})\otimes g^{\omega _{1}}x_{1}^{\omega
_{1}} \\
&=&B^{A}(1_{H}\otimes gx_{1})\otimes B^{H}(1_{H}\otimes gx_{1})\otimes 1_{H}
\\
&&B^{A}(1_{H}\otimes g)\otimes B^{H}(1_{H}\otimes g)\otimes gx_{1}
\end{eqnarray*}%
\begin{eqnarray*}
&&\sum_{a,b_{1},b_{2},d,e_{1},e_{2}=0}^{1}\sum_{l_{1}=0}^{b_{1}}%
\sum_{l_{2}=0}^{b_{2}}\sum_{u_{1}=0}^{e_{1}}\sum_{u_{2}=0}^{e_{2}}\left(
-1\right) ^{\alpha \left( 1_{H};l_{1},l_{2},u_{1},u_{2}\right) } \\
&&B(1_{H}\otimes
gx_{1};G^{a}X_{1}^{b_{1}}X_{2}^{b_{2}},g^{d}x_{1}^{e_{1}}x_{2}^{e_{2}}) \\
&&G^{a}X_{1}^{b_{1}-l_{1}}X_{2}^{b_{2}-l_{2}}\otimes
g^{d}x_{1}^{e_{1}-u_{1}}x_{2}^{e_{2}-u_{2}}\otimes \\
&&g^{a+b_{1}+b_{2}+l_{1}+l_{2}+d+e_{1}+e_{2}+u_{1}+u_{2}}x_{1}^{l_{1}+u_{1}}x_{2}^{l_{2}+u_{2}}
\\
&=&\sum_{\omega _{1}=0}^{1}B^{A}(1_{H}\otimes gx_{1}^{1-\omega _{1}})\otimes
B^{H}(1_{H}\otimes gx_{1}^{1-\omega _{1}})\otimes gx_{1}^{\omega _{1}}
\end{eqnarray*}

\begin{eqnarray*}
a+b_{1}+b_{2}+l_{1}+l_{2}+d+e_{1}+e_{2}+u_{1}+u_{2} &\equiv &1 \\
l_{1}+u_{1} &=&1 \\
l_{2}+u_{2} &=&1
\end{eqnarray*}%
so that%
\begin{equation*}
a+b_{1}+b_{2}+d+e_{1}+e_{2}\equiv 1
\end{equation*}%
and we know that, in view of $\left( \ref{1otgx1, first}\right) ,$%
\begin{equation*}
B(1_{H}\otimes
gx_{1};G^{a}X_{1}^{b_{1}}X_{2}^{b_{2}},g^{d}x_{1}^{e_{1}}x_{2}^{e_{2}})=0%
\text{ whenever }a+b_{1}+b_{2}+d+e_{1}+e_{2}\equiv 1.
\end{equation*}

\section{$B\left( 1\otimes gx_{2}\right) $}

From $\left( \ref{simplgx}\right) $ we get%
\begin{eqnarray}
&&B(1_{H}\otimes gx_{2})  \label{form 1otgx2} \\
&=&(1_{A}\otimes g)B(g\otimes 1_{H})(1_{A}\otimes gx_{2})  \notag \\
&&+(1_{A}\otimes x_{2})B(g\otimes 1_{H})  \notag \\
&&+B(x_{2}\otimes 1_{H})  \notag
\end{eqnarray}%
so that we obtain%
\begin{eqnarray*}
B\left( 1_{H}\otimes gx_{2}\right) &=&B\left( x_{2}\otimes
1_{H};1_{A},1_{H}\right) 1_{A}\otimes 1_{H} \\
&&+B\left( x_{2}\otimes 1_{H};1_{A},x_{1}x_{2}\right) 1_{A}\otimes x_{1}x_{2}
\\
&&+B\left( x_{2}\otimes 1_{H};1_{A},gx_{1}\right) 1_{A}\otimes gx_{1}+ \\
&&+\left[ 2B\left( g\otimes 1_{H};1_{A},g\right) +B\left( x_{2}\otimes
1_{H};1_{A},gx_{2}\right) \right] 1_{A}\otimes gx_{2}+ \\
&&+B\left( x_{2}\otimes 1_{H};G,g\right) G\otimes g \\
&&+B\left( x_{2}\otimes 1_{H};G,x_{1}\right) G\otimes x_{1}+ \\
&&+B\left( x_{2}\otimes 1_{H};G,x_{2}\right) G\otimes x_{2}+ \\
&&+\left[ 2B\left( g\otimes 1_{H};G,gx_{1}\right) +B\left( x_{2}\otimes
1_{H};G,gx_{1}x_{2}\right) \right] G\otimes gx_{1}x_{2}+ \\
&&-B(x_{2}\otimes 1_{H};1_{A},gx_{1})X_{1}\otimes g+ \\
&&+B(x_{2}\otimes 1_{H};1_{A},x_{1}x_{2})X_{1}\otimes x_{2}+ \\
&&+\left[ -B(g\otimes 1_{H};1_{A},g)-B(x_{2}\otimes \ 1_{H};1_{A},gx_{2})%
\right] X_{2}\otimes g+ \\
&&+\left[ -B(g\otimes 1_{H};1_{A},x_{1})-B(x_{2}\otimes \
1_{H};1_{A},x_{1}x_{2})\right] X_{2}\otimes x_{1} \\
&&-B(g\otimes 1_{H};1_{A},x_{2})X_{2}\otimes x_{2}+ \\
&&+B\left( g\otimes 1_{H};1_{A},gx_{1}x_{2}\right) X_{2}\otimes gx_{1}x_{2}+
\\
&&+\left[ B(g\otimes 1_{H};1_{A},x_{1})+B(x_{2}\otimes \
1_{H};1_{A},x_{1}x_{2})\right] X_{1}X_{2}\otimes 1_{H}+ \\
&&+B\left( g\otimes 1_{H};1_{A},gx_{1}x_{2}\right) X_{1}X_{2}\otimes gx_{2}+
\\
&&+B(x_{2}\otimes 1_{H};G,x_{1})GX_{1}\otimes 1_{H}+ \\
&&\left[ -2B\left( g\otimes 1_{H};G,gx_{1}\right) -B(x_{2}\otimes
1_{H};G,gx_{1}x_{2})\right] GX_{1}\otimes gx_{2}+ \\
&&+\left[ B(g\otimes 1_{H};G,1_{H})+B(x_{2}\otimes \ 1_{H};G,x_{2}\right]
GX_{2}\otimes 1_{H}+ \\
&&+B(g\otimes 1_{H};GX_{1}X_{2},1_{H})GX_{2}\otimes x_{1}x_{2}+ \\
&&+\left[ B(g\otimes 1_{H};G,gx_{1})+B(x_{2}\otimes \ 1_{H};G,gx_{1}x_{2})%
\right] GX_{2}\otimes gx_{1}+ \\
&&-B(g\otimes 1_{H};G,gx_{2})GX_{2}\otimes gx_{2} \\
&&+\left[ B(g\otimes 1_{H};G,gx_{1})+B(x_{2}\otimes \ 1_{H};G,gx_{1}x_{2})%
\right] GX_{1}X_{2}\otimes g+ \\
&&-B(g\otimes 1_{H};GX_{1}X_{2},1_{H})GX_{1}X_{2}\otimes x_{2}
\end{eqnarray*}

From $\left( \ref{MAIN FORMULA 1}\right) $ we get%
\begin{eqnarray*}
&&\sum_{a,b_{1},b_{2},d,e_{1},e_{2}=0}^{1}\sum_{l_{1}=0}^{b_{1}}%
\sum_{l_{2}=0}^{b_{2}}\sum_{u_{1}=0}^{e_{1}}\sum_{u_{2}=0}^{e_{2}}\left(
-1\right) ^{\alpha \left( 1_{H};l_{1},l_{2},u_{1},u_{2}\right) } \\
&&B(1_{H}\otimes
gx_{2};G^{a}X_{1}^{b_{1}}X_{2}^{b_{2}},g^{d}x_{1}^{e_{1}}x_{2}^{e_{2}}) \\
&&G^{a}X_{1}^{b_{1}-l_{1}}X_{2}^{b_{2}-l_{2}}\otimes
g^{d}x_{1}^{e_{1}-u_{1}}x_{2}^{e_{2}-u_{2}}\otimes \\
&&g^{a+b_{1}+b_{2}+l_{1}+l_{2}+d+e_{1}+e_{2}+u_{1}+u_{2}}x_{1}^{l_{1}+u_{1}}x_{2}^{l_{2}+u_{2}}
\\
&=&\sum_{\omega _{2}=0}^{1}B^{A}(1_{H}\otimes gx_{2}^{1-\omega _{2}})\otimes
B^{H}(1_{H}\otimes gx_{2}^{1-\omega _{2}})\otimes g^{\omega
_{2}}x_{2}^{\omega _{2}}
\end{eqnarray*}%
i.e.
\begin{eqnarray*}
&&\sum_{a,b_{1},b_{2},d,e_{1},e_{2}=0}^{1}\sum_{l_{1}=0}^{b_{1}}%
\sum_{l_{2}=0}^{b_{2}}\sum_{u_{1}=0}^{e_{1}}\sum_{u_{2}=0}^{e_{2}}\left(
-1\right) ^{\alpha \left( 1_{H};l_{1},l_{2},u_{1},u_{2}\right) } \\
&&B(1_{H}\otimes
gx_{2};G^{a}X_{1}^{b_{1}}X_{2}^{b_{2}},g^{d}x_{1}^{e_{1}}x_{2}^{e_{2}}) \\
&&G^{a}X_{1}^{b_{1}-l_{1}}X_{2}^{b_{2}-l_{2}}\otimes
g^{d}x_{1}^{e_{1}-u_{1}}x_{2}^{e_{2}-u_{2}}\otimes \\
&&g^{a+b_{1}+b_{2}+l_{1}+l_{2}+d+e_{1}+e_{2}+u_{1}+u_{2}}x_{1}^{l_{1}+u_{1}}x_{2}^{l_{2}+u_{2}}
\\
&=&B^{A}(1_{H}\otimes gx_{2})\otimes B^{H}(1_{H}\otimes gx_{2}^{{}})\otimes
1_{H}+ \\
&&B^{A}(1_{H}\otimes g)\otimes B^{H}(1_{H}\otimes g)\otimes gx_{2}
\end{eqnarray*}

\subsection{case $1_{H}$}

\begin{eqnarray*}
a+b_{1}+b_{2}+l_{1}+l_{2}+d+e_{1}+e_{2}+u_{1}+u_{2} &=&0 \\
l_{1}+u_{1} &=&0 \\
l_{2}+u_{2} &=&0
\end{eqnarray*}%
so that%
\begin{eqnarray*}
a+b_{1}+b_{2}+d+e_{1}+e_{2} &=&0 \\
l_{1} &=&u_{1}=0 \\
l_{2} &=&u_{2}=0
\end{eqnarray*}%
the right side of the equality gives us
\begin{equation*}
B^{A}(1_{H}\otimes gx_{2})\otimes B^{H}(1_{H}\otimes gx_{2})
\end{equation*}%
and the equality becomes, as
\begin{equation*}
\alpha \left( 1_{H};0,0,0,0\right) =\left( a+b_{1}+b_{2}\right) \left(
n_{1}+n_{2}\right) =0
\end{equation*}%
\begin{eqnarray*}
&&\sum_{\substack{ a,b_{1},b_{2},d,e_{1},e_{2}=0  \\ %
a+b_{1}+b_{2}+d+e_{1}+e_{2}\equiv 0}}^{1}B(1_{H}\otimes
gx_{2};G^{a}X_{1}^{b_{1}}X_{2}^{b_{2}},g^{d}x_{1}^{e_{1}}x_{2}^{e_{2}}) \\
&&G^{a}X_{1}^{b_{1}}X_{2}^{b_{2}}\otimes
g^{d}x_{1}^{e_{1}}x_{2}^{e_{2}}\otimes 1_{H} \\
&=&B^{A}(1_{H}\otimes gx_{2})\otimes B^{H}(1_{H}\otimes gx_{2})\otimes 1_{H}
\end{eqnarray*}%
so that we deduce that%
\begin{equation}
B(1_{H}\otimes
gx_{2};G^{a}X_{1}^{b_{1}}X_{2}^{b_{2}},g^{d}x_{1}^{e_{1}}x_{2}^{e_{2}})=0%
\text{ whenever }a+b_{1}+b_{2}+d+e_{1}+e_{2}\equiv 1  \label{1otgx2, first}
\end{equation}

\subsection{case $g$}

\begin{eqnarray*}
&&\sum_{a,b_{1},b_{2},d,e_{1},e_{2}=0}^{1}\sum_{l_{1}=0}^{b_{1}}%
\sum_{l_{2}=0}^{b_{2}}\sum_{u_{1}=0}^{e_{1}}\sum_{u_{2}=0}^{e_{2}}\left(
-1\right) ^{\alpha \left( 1_{H};l_{1},l_{2},u_{1},u_{2}\right) } \\
&&B(1_{H}\otimes
gx_{2};G^{a}X_{1}^{b_{1}}X_{2}^{b_{2}},g^{d}x_{1}^{e_{1}}x_{2}^{e_{2}}) \\
&&G^{a}X_{1}^{b_{1}-l_{1}}X_{2}^{b_{2}-l_{2}}\otimes
g^{d}x_{1}^{e_{1}-u_{1}}x_{2}^{e_{2}-u_{2}}\otimes \\
&&g^{a+b_{1}+b_{2}+l_{1}+l_{2}+d+e_{1}+e_{2}+u_{1}+u_{2}}x_{1}^{l_{1}+u_{1}}x_{2}^{l_{2}+u_{2}}
\\
&=&\sum_{\omega _{2}=0}^{1}B^{A}(1_{H}\otimes gx_{2}^{1-\omega _{2}})\otimes
B^{H}(1_{H}\otimes gx_{2}^{1-\omega _{2}})\otimes g^{\omega
_{2}}x_{2}^{\omega _{2}}
\end{eqnarray*}%
From the left side%
\begin{eqnarray*}
a+b_{1}+b_{2}+l_{1}+l_{2}+d+e_{1}+e_{2}+u_{1}+u_{2} &\equiv &1 \\
l_{1}+u_{1} &=&0 \\
l_{2}+u_{2} &\equiv &0
\end{eqnarray*}%
\begin{eqnarray*}
a+b_{1}+b_{2}+l_{1}+l_{2}+d+e_{1}+e_{2}+u_{1}+u_{2} &\equiv &1 \\
l_{1}+u_{1} &=&0 \\
l_{2}+u_{2} &\equiv &0
\end{eqnarray*}%
i.e.

\begin{eqnarray*}
a+b_{1}+b_{2}+d+e_{1}+e_{2} &=&1 \\
l_{1} &=&u_{1}=0 \\
l_{2} &=&u_{2}=0
\end{eqnarray*}

From the right side we get just zero. Thus we get
\begin{equation*}
B(1_{H}\otimes
gx_{2};G^{a}X_{1}^{b_{1}}X_{2}^{b_{2}},g^{d}x_{1}^{e_{1}}x_{2}^{e_{2}})=0%
\text{ whenever }a+b_{1}+b_{2}+d+e_{1}+e_{2}\equiv 1
\end{equation*}%
which we already know.

\subsection{case $x_{1}$}

The right side of the equality is $0.$ As far as the left side is concerned,
we get

\begin{equation*}
B(1_{H}\otimes
gx_{2};G^{a}X_{1}^{b_{1}}X_{2}^{b_{2}},g^{d}x_{1}^{e_{1}}x_{2}^{e_{2}})=0%
\text{ whenever }a+b_{1}+b_{2}+d+e_{1}+e_{2}\equiv 1
\end{equation*}%
which we already know.

\subsection{case $x_{2}$}

\begin{equation*}
B(1_{H}\otimes
gx_{2};G^{a}X_{1}^{b_{1}}X_{2}^{b_{2}},g^{d}x_{1}^{e_{1}}x_{2}^{e_{2}})=0%
\text{ whenever }a+b_{1}+b_{2}+d+e_{1}+e_{2}\equiv 1
\end{equation*}%
which we already know.

\subsection{case $x_{1}x_{2}$}

The right side of the equality is $0.$ As far as the left side is concerned,
we get

\begin{eqnarray*}
a+b_{1}+b_{2}+l_{1}+l_{2}+d+e_{1}+e_{2}+u_{1}+u_{2} &\equiv &0 \\
l_{1}+u_{1} &=&1 \\
l_{2}+u_{2} &\equiv &1
\end{eqnarray*}%
and hence%
\begin{eqnarray*}
a+b_{1}+b_{2}+d+e_{1}+e_{2} &\equiv &0 \\
l_{1}+u_{1} &=&1 \\
l_{2}+u_{2} &\equiv &1
\end{eqnarray*}%
Since
\begin{eqnarray*}
\alpha \left( 1_{H};0,0,1,1\right) &=&1+e_{2}\text{, }\alpha \left(
1_{H};0,1,1,0\right) =e_{2}+a+b_{1}+b_{2}+1, \\
\alpha \left( 1_{H};1,0,0,1\right) &=&a+b_{1}\text{and }\alpha \left(
1_{H};1,1,0,0\right) =1+b_{2},
\end{eqnarray*}%
we get%
\begin{eqnarray*}
&&\sum_{\substack{ a,b_{1},b_{2},d=0  \\ a+b_{1}+b_{2}+d\equiv 0}}%
^{1}B(1_{H}\otimes
gx_{2};G^{a}X_{1}^{b_{1}}X_{2}^{b_{2}},g^{d}x_{1}x_{2})G^{a}X_{1}^{b_{1}}X_{2}^{b_{2}}\otimes g^{d}
\\
&&\sum_{\substack{ a,b_{1},d,e_{2}=0  \\ a+b_{1}+d+e_{2}\equiv 0}}^{1}\left(
-1\right) ^{e_{2}+a+b_{1}}B(1_{H}\otimes
gx_{2};G^{a}X_{1}^{b_{1}}X_{2},g^{d}x_{1}x_{2}^{e_{2}})G^{a}X_{1}^{b_{1}}%
\otimes g^{d}x_{2}^{e_{2}} \\
&&\sum_{\substack{ a,b_{2},d,e_{1}=0  \\ a+b_{2}+d+e_{1}\equiv 0}}^{1}\left(
-1\right) ^{a+1}B(1_{H}\otimes
gx_{2};G^{a}X_{1}X_{2}^{b_{2}},g^{d}x_{1}^{e_{1}}x_{2})G^{a}X_{2}^{b_{2}}%
\otimes g^{d}x_{1}^{e_{1}} \\
&&\sum_{\substack{ a,d,e_{1},e_{2}=0  \\ a+d+e_{1}+e_{2}\equiv 0}}%
^{1}B(1_{H}\otimes
gx_{2};G^{a}X_{1}X_{2},g^{d}x_{1}^{e_{1}}x_{2}^{e_{2}})G^{a}\otimes
g^{d}x_{1}^{e_{1}}x_{2}^{e_{2}} \\
&=&0
\end{eqnarray*}

\subsubsection{$G^{a}\otimes g^{d}$}

\begin{eqnarray*}
&&\sum_{\substack{ a,d=0  \\ a+d\equiv 0}}^{1}B(1_{H}\otimes
gx_{2};G^{a},g^{d}x_{1}x_{2})G^{a}\otimes g^{d}+ \\
&&\sum_{\substack{ a,d=0  \\ a+d\equiv 0}}^{1}\left( -1\right)
^{e_{2}+a+b_{1}}B(1_{H}\otimes gx_{2};G^{a}X_{2},g^{d}x_{1})G^{a}\otimes
g^{d}+ \\
&&\sum_{\substack{ a,d=0  \\ a+d\equiv 0}}^{1}\left( -1\right)
^{a+1}B(1_{H}\otimes gx_{2};G^{a}X_{1},g^{d}x_{2})G^{a}\otimes g^{d}+ \\
&&\sum_{\substack{ a,d=0  \\ a+d\equiv 0}}^{1}B(1_{H}\otimes
gx_{2};G^{a}X_{1}X_{2},g^{d})G^{a}\otimes g^{d}
\end{eqnarray*}%
and we get%
\begin{equation*}
\sum_{\substack{ a,d=0  \\ a+d\equiv 0}}^{1}\left[
\begin{array}{c}
B(1_{H}\otimes gx_{2};G^{a},g^{d}x_{1}x_{2})+\left( -1\right)
^{a}B(1_{H}\otimes gx_{2};G^{a}X_{2},g^{d}x_{1}) \\
\left( -1\right) ^{a+1}B(1_{H}\otimes
gx_{2};G^{a}X_{1},g^{d}x_{2})+B(1_{H}\otimes gx_{2};G^{a}X_{1}X_{2},g^{d})%
\end{array}%
\right] G^{a}\otimes g^{d}=0
\end{equation*}%
i.e.%
\begin{eqnarray*}
&&B(1_{H}\otimes gx_{2};1_{A},x_{1}x_{2})+B(1_{H}\otimes gx_{2};X_{2},x_{1})+
\\
&&-B(1_{H}\otimes gx_{2};X_{1},x_{2})+B(1_{H}\otimes gx_{2};X_{1}X_{2},1_{H})
\\
&=&0
\end{eqnarray*}%
and%
\begin{eqnarray*}
&&B(1_{H}\otimes gx_{2};G,gx_{1}x_{2})-B(1_{H}\otimes gx_{2};GX_{2},gx_{1})+
\\
&&+B(1_{H}\otimes gx_{2};GX_{1},gx_{2})+B(1_{H}\otimes gx_{2};GX_{1}X_{2},g)
\\
&=&0
\end{eqnarray*}%
By the form of the element, these equalities are satisfied.

\subsubsection{$G^{a}\otimes g^{d}x_{2}$}

\begin{equation*}
\sum_{\substack{ a,d=0  \\ a+d\equiv 1}}^{1}\left[ \left( -1\right)
^{1+a}B(1_{H}\otimes gx_{2};G^{a}X_{2},g^{d}x_{1}x_{2})+B(1_{H}\otimes
gx_{2};G^{a}X_{1}X_{2},g^{d}x_{1}^{e_{1}}x_{2}^{e_{2}})\right] G^{a}\otimes
g^{d}x_{2}=0
\end{equation*}%
and we get%
\begin{equation*}
-B(1_{H}\otimes gx_{2};X_{2},gx_{1}x_{2})+B(1_{H}\otimes
gx_{2};X_{1}X_{2},gx_{2})=0
\end{equation*}%
and%
\begin{equation*}
B(1_{H}\otimes gx_{2};GX_{2},x_{1}x_{2})+B(1_{H}\otimes
gx_{2};GX_{1}X_{2},x_{2})=0.
\end{equation*}%
By the form of the element, these equalities are satisfied.

\subsubsection{$G^{a}\otimes g^{d}x_{1}$}

\begin{equation*}
\sum_{\substack{ a,d=0  \\ a+d\equiv 1}}^{1}\left[ \left( -1\right)
^{a+1}B(1_{H}\otimes gx_{2};G^{a}X_{1},g^{d}x_{1}x_{2})+B(1_{H}\otimes
gx_{2};G^{a}X_{1}X_{2},g^{d}x_{1})\right] G^{a}\otimes g^{d}x_{1}=0
\end{equation*}%
and we get%
\begin{equation*}
-B(1_{H}\otimes gx_{2};X_{1},gx_{1}x_{2})+B(1_{H}\otimes
gx_{2};X_{1}X_{2},gx_{1})=0,
\end{equation*}%
\begin{equation*}
B(1_{H}\otimes gx_{2};GX_{1},x_{1}x_{2})+B(1_{H}\otimes
gx_{2};GX_{1}X_{2},x_{1})=0\text{ .}
\end{equation*}%
By the form of the element, these equalities are satisfied.

\subsubsection{$G^{a}X_{2}\otimes g^{d}$}

\begin{equation*}
\sum_{\substack{ a,d=0  \\ a+d\equiv 1}}^{1}\left[ B(1_{H}\otimes
gx_{2};G^{a}X_{2},g^{d}x_{1}x_{2})+\left( -1\right) ^{a+1}B(1_{H}\otimes
gx_{2};G^{a}X_{1}X_{2},g^{d}x_{2})\right] G^{a}X_{2}\otimes g^{d}=0
\end{equation*}%
We get%
\begin{equation*}
B(1_{H}\otimes gx_{2};X_{2},gx_{1}x_{2})-B(1_{H}\otimes
gx_{2};X_{1}X_{2},gx_{2})=0\text{ }
\end{equation*}%
and%
\begin{equation*}
B(1_{H}\otimes gx_{2};GX_{2},x_{1}x_{2})+B(1_{H}\otimes
gx_{2};GX_{1}X_{2},x_{2})=0\text{ }.
\end{equation*}%
Both of them were already found.

\subsubsection{$G^{a}X_{1}\otimes g^{d}$}

\begin{equation*}
\sum_{\substack{ a,d=0  \\ a+d\equiv 1}}^{1}\left[ B(1_{H}\otimes
gx_{2};G^{a}X_{1},g^{d}x_{1}x_{2})+\left( -1\right) ^{a+1}B(1_{H}\otimes
gx_{2};G^{a}X_{1}X_{2},g^{d}x_{1})\right] G^{a}X_{1}\otimes g^{d}=0
\end{equation*}%
and we get%
\begin{equation*}
B(1_{H}\otimes gx_{2};X_{1},gx_{1}x_{2})-B(1_{H}\otimes
gx_{2};X_{1}X_{2},gx_{1})=0\text{ }
\end{equation*}%
and%
\begin{equation*}
B(1_{H}\otimes gx_{2};GX_{1},x_{1}x_{2})+B(1_{H}\otimes
gx_{2};GX_{1}X_{2},x_{1})=0\text{ .}
\end{equation*}%
Both of them were already found.

\subsubsection{$G^{a}\otimes g^{d}x_{1}x_{2}$}

\begin{equation*}
\sum_{\substack{ a,d=0  \\ a+d\equiv 0}}^{1}B(1_{H}\otimes
gx_{2};G^{a}X_{1}X_{2},g^{d}x_{1}x_{2})G^{a}\otimes g^{d}x_{1}x_{2}=0
\end{equation*}%
We get%
\begin{equation*}
B(1_{H}\otimes gx_{2};X_{1}X_{2},x_{1}x_{2})=0
\end{equation*}%
and%
\begin{equation*}
B(1_{H}\otimes gx_{2};GX_{1}X_{2},gx_{1}x_{2})=0.
\end{equation*}%
By the form of the element, these equalities are satisfied.

\subsubsection{$G^{a}X_{2}\otimes g^{d}x_{2}$}

There is no term like this.

\subsubsection{$G^{a}X_{1}\otimes g^{d}x_{2}$}

We get%
\begin{equation*}
B(1_{H}\otimes gx_{2};X_{1}X_{2},x_{1}x_{2})=0
\end{equation*}%
and%
\begin{equation*}
B(1_{H}\otimes gx_{2};GX_{1}X_{2},gx_{1}x_{2})=0\text{ .}
\end{equation*}%
Both of them were already found

\subsubsection{$G^{a}X_{2}\otimes g^{d}x_{1}$}

\begin{equation*}
\sum_{\substack{ a,d=0  \\ a+d\equiv 0}}^{1}B(1_{H}\otimes
gx_{2};G^{a}X_{1}X_{2},g^{d}x_{1}x_{2})G^{a}X_{2}\otimes g^{d}x_{1}=0
\end{equation*}%
and we obtain%
\begin{equation*}
B(1_{H}\otimes gx_{2};X_{1}X_{2},x_{1}x_{2})=0
\end{equation*}%
and%
\begin{equation*}
B(1_{H}\otimes gx_{2};GX_{1}X_{2},gx_{1}x_{2})=0
\end{equation*}%
which we already got.

\subsubsection{$G^{a}X_{1}\otimes g^{d}x_{1}$}

There is no term like this.

\subsubsection{$G^{a}X_{1}X_{2}\otimes g^{d}$}

\begin{equation*}
\sum_{\substack{ a,d=0  \\ a+d\equiv 0}}^{1}B(1_{H}\otimes
gx_{2};G^{a}X_{1}X_{2},g^{d}x_{1}x_{2})G^{a}X_{1}X_{2}\otimes g^{d}=0
\end{equation*}%
and we get%
\begin{equation*}
B(1_{H}\otimes gx_{2};X_{1}X_{2},x_{1}x_{2})=0
\end{equation*}%
and%
\begin{equation*}
B(1_{H}\otimes gx_{2};GX_{1}X_{2},gx_{1}x_{2})=0
\end{equation*}%
which we already got.

\subsubsection{$G^{a}X_{2}\otimes g^{d}x_{1}x_{2}$}

There is no term like this.

\subsubsection{$G^{a}X_{1}\otimes g^{d}x_{1}x_{2}$}

There is no term like this.

\subsubsection{$G^{a}X_{1}X_{2}\otimes g^{d}x_{2}$}

There is no term like this.

\subsubsection{$G^{a}X_{1}X_{2}\otimes g^{d}x_{1}$}

There is no term like this.

\subsubsection{$G^{a}X_{1}X_{2}\otimes g^{d}x_{1}x_{2}$}

There is no term like this.

\subsection{case $gx_{1}$}

The right side of the equality is zero.%
\begin{eqnarray*}
a+b_{1}+b_{2}+l_{1}+l_{2}+d+e_{1}+e_{2}+u_{1}+u_{2} &=&1 \\
l_{1}+u_{1} &=&1 \\
l_{2}+u_{2} &=&0
\end{eqnarray*}%
\begin{eqnarray*}
a+b_{1}+b_{2}+d+e_{1}+e_{2} &=&0 \\
l_{1}+u_{1} &=&1 \\
l_{2} &=&u_{2}=0
\end{eqnarray*}%
We get%
\begin{eqnarray*}
&&\sum_{\substack{ a,b_{1},b_{2},d,e_{1},e_{2}=0  \\ %
a+b_{1}+b_{2}+d+e_{1}+e_{2}=0}}^{1}\sum_{l_{1}=0}^{b_{1}}\sum_{\substack{ %
u_{1}=0  \\ l_{1}+u_{1}=1}}^{e_{1}}\left( -1\right) ^{\alpha \left(
1_{H};l_{1},0,u_{1},0\right) }B(1_{H}\otimes
gx_{2};G^{a}X_{1}^{b_{1}}X_{2}^{b_{2}},g^{d}x_{1}^{e_{1}}x_{2}^{e_{2}}) \\
&&G^{a}X_{1}^{b_{1}-l_{1}}X_{2}^{b_{2}}\otimes
g^{d}x_{1}^{e_{1}-u_{1}}x_{2}^{e_{2}}=0
\end{eqnarray*}%
Since $\alpha \left( 1_{H};0,0,1,0\right) =e_{2}+\left( a+b_{1}+b_{2}\right)
$ and $\alpha \left( 1_{H};1,0,0,0\right) =b_{2}$ we get

\begin{eqnarray*}
&&\sum_{\substack{ a,b_{1},b_{2},d,e_{2}=0  \\ a+b_{1}+b_{2}+d+e_{2}=1}}%
^{1}\left( -1\right) ^{e_{2}+\left( a+b_{1}+b_{2}\right) }B(1_{H}\otimes
gx_{2};G^{a}X_{1}^{b_{1}}X_{2}^{b_{2}},g^{d}x_{1}x_{2}^{e_{2}})G^{a}X_{1}^{b_{1}}X_{2}^{b_{2}}\otimes g^{d}x_{2}^{e_{2}}+
\\
&&\sum_{\substack{ a,b_{2},d,e_{1},e_{2}=0  \\ a+b_{2}+d+e_{1}+e_{2}=1}}%
^{1}\left( -1\right) ^{b_{2}}B(1_{H}\otimes
gx_{2};G^{a}X_{1}X_{2}^{b_{2}},g^{d}x_{1}^{e_{1}}x_{2}^{e_{2}})G^{a}X_{2}^{b_{2}}\otimes g^{d}x_{1}^{e_{1}}x_{2}^{e_{2}}
\end{eqnarray*}

\subsubsection{$G^{a}\otimes g^{d}$}

\begin{equation*}
\sum_{\substack{ a,d=0  \\ a+d\equiv 1}}^{1}\left[ \left( -1\right)
^{a}B(1_{H}\otimes gx_{2};G^{a},g^{d}x_{1})+B(1_{H}\otimes
gx_{2};G^{a}X_{1},g^{d})\right] G^{a}\otimes g^{d}=0
\end{equation*}%
and we get%
\begin{equation*}
B(1_{H}\otimes gx_{2};1_{A},gx_{1})+B(1_{H}\otimes gx_{2};X_{1},g)=0
\end{equation*}%
and%
\begin{equation*}
-B(1_{H}\otimes gx_{2};G,x_{1})+B(1_{H}\otimes gx_{2};GX_{1},1_{H})=0.
\end{equation*}%
By the form of the element, these equalities are satisfied.

\subsubsection{$G^{a}\otimes g^{d}x_{2}$}

\begin{equation*}
\sum_{\substack{ a,d=0  \\ a+d\equiv 0}}^{1}\left[ \left( -1\right)
^{1+a}B(1_{H}\otimes gx_{2};G^{a},g^{d}x_{1}x_{2})+B(1_{H}\otimes
gx_{2};G^{a}X_{1},g^{d}x_{2})\right] G^{a}\otimes g^{d}x_{2}=0
\end{equation*}%
and we get%
\begin{equation*}
-B(1_{H}\otimes gx_{2};1_{A},x_{1}x_{2})+B(1_{H}\otimes gx_{2};X_{1},x_{2})=0
\end{equation*}%
and%
\begin{equation*}
B(1_{H}\otimes gx_{2};G,gx_{1}x_{2})+B(1_{H}\otimes gx_{2};GX_{1},gx_{2})=0.
\end{equation*}%
By the form of the element, these equalities are satisfied.

\subsubsection{$G^{a}\otimes g^{d}x_{1}$}

\begin{equation*}
\sum_{\substack{ a,d=0  \\ a+d\equiv 0}}^{1}B(1_{H}\otimes
gx_{2};G^{a}X_{1},g^{d}x_{1})G^{a}\otimes g^{d}x_{1}
\end{equation*}%
and we get%
\begin{equation*}
B(1_{H}\otimes gx_{2};X_{1},x_{1})=0
\end{equation*}%
and%
\begin{equation*}
B(1_{H}\otimes gx_{2};GX_{1},gx_{1})=0.
\end{equation*}%
By the form of the element, these equalities are satisfied.

\subsubsection{$G^{a}X_{2}\otimes g^{d}$}

\begin{equation*}
\sum_{\substack{ a,d=0  \\ a+d\equiv 0}}^{1}\left[ \left( -1\right)
^{a+1}B(1_{H}\otimes gx_{2};G^{a}X_{2},g^{d}x_{1})-B(1_{H}\otimes
gx_{2};G^{a}X_{1}X_{2},g^{d})\right] G^{a}X_{2}\otimes g^{d}=0
\end{equation*}%
and we get%
\begin{equation*}
-B(1_{H}\otimes gx_{2};X_{2},x_{1})-B(1_{H}\otimes gx_{2};X_{1}X_{2},1_{H})=0
\end{equation*}%
and%
\begin{equation*}
B(1_{H}\otimes gx_{2};GX_{2},gx_{1})-B(1_{H}\otimes gx_{2};GX_{1}X_{2},g)=0.
\end{equation*}%
By the form of the element, these equalities are satisfied.

\subsubsection{$G^{a}X_{1}\otimes g^{d}$}

\begin{equation*}
\sum_{\substack{ a,d=0  \\ a+d\equiv 0}}^{1}\left( -1\right)
^{a+1}B(1_{H}\otimes gx_{2};G^{a}X_{1},g^{d}x_{1})G^{a}X_{1}\otimes g^{d}=0
\end{equation*}%
We obtain%
\begin{equation*}
-B(1_{H}\otimes gx_{2};X_{1},x_{1})=0\text{ }
\end{equation*}%
and%
\begin{equation*}
B(1_{H}\otimes gx_{2};GX_{1},gx_{1})=0\text{ .}
\end{equation*}%
Both of them were already got.

\subsubsection{$G^{a}\otimes g^{d}x_{1}x_{2}$}

\begin{equation*}
\sum_{\substack{ a,d=0  \\ a+d\equiv 1}}^{1}B(1_{H}\otimes
gx_{2};G^{a}X_{1},g^{d}x_{1}x_{2})G^{a}\otimes g^{d}x_{1}x_{2}=0
\end{equation*}%
and we get%
\begin{equation*}
B(1_{H}\otimes gx_{2};X_{1},gx_{1}x_{2})=0
\end{equation*}%
and%
\begin{equation*}
B(1_{H}\otimes gx_{2};GX_{1},x_{1}x_{2})=0.
\end{equation*}%
By the form of the element, these equalities are satisfied.

\subsubsection{$G^{a}X_{2}\otimes g^{d}x_{2}$}

\begin{equation*}
\sum_{\substack{ a,d=0  \\ a+d\equiv 1}}^{1}\left[ \left( -1\right)
^{a}B(1_{H}\otimes gx_{2};G^{a}X_{2},g^{d}x_{1}x_{2})-B(1_{H}\otimes
gx_{2};G^{a}X_{1}X_{2},g^{d}x_{2})\right] G^{a}X_{2}\otimes g^{d}x_{2}=0
\end{equation*}%
and we obtained%
\begin{equation*}
B(1_{H}\otimes gx_{2};X_{2},gx_{1}x_{2})-B(1_{H}\otimes
gx_{2};X_{1}X_{2},gx_{2})=0
\end{equation*}%
\begin{equation*}
-B(1_{H}\otimes gx_{2};GX_{2},x_{1}x_{2})-B(1_{H}\otimes
gx_{2};GX_{1}X_{2},x_{2})=0.
\end{equation*}%
By the form of the element, these equalities are satisfied.

\subsubsection{$G^{a}X_{1}\otimes g^{d}x_{2}$}

\begin{equation*}
\sum_{\substack{ a,d=0  \\ a+d\equiv 1}}^{1}\left( -1\right)
^{a}B(1_{H}\otimes gx_{2};G^{a}X_{1},g^{d}x_{1}x_{2})G^{a}X_{1}\otimes
g^{d}x_{2}=0
\end{equation*}%
and we get the equalities below which we already got.%
\begin{equation*}
B(1_{H}\otimes gx_{2};X_{1},gx_{1}x_{2})=0\text{ ,}
\end{equation*}%
\begin{equation*}
-B(1_{H}\otimes gx_{2};GX_{1},x_{1}x_{2})=0\text{ .}
\end{equation*}

\subsubsection{$G^{a}X_{2}\otimes g^{d}x_{1}$}

\begin{equation*}
\sum_{\substack{ a,d=0  \\ a+d\equiv 1}}^{1}-B(1_{H}\otimes
gx_{2};G^{a}X_{1}X_{2},g^{d}x_{1})G^{a}X_{2}\otimes g^{d}x_{1}=0
\end{equation*}%
and we get%
\begin{equation*}
-B(1_{H}\otimes gx_{2};X_{1}X_{2},gx_{1})=0
\end{equation*}%
and%
\begin{equation*}
-B(1_{H}\otimes gx_{2};GX_{1}X_{2},x_{1})=0.
\end{equation*}%
By the form of the element, these equalities are satisfied.

\subsubsection{$G^{a}X_{1}\otimes g^{d}x_{1}$}

There is no term like this

\subsubsection{$G^{a}X_{1}X_{2}\otimes g^{d}$}

\begin{equation*}
\sum_{\substack{ a,d=0  \\ a+d\equiv 1}}^{1}\left( -1\right)
^{a}B(1_{H}\otimes gx_{2};G^{a}X_{1}X_{2},g^{d}x_{1})G^{a}X_{1}X_{2}\otimes
g^{d}=0
\end{equation*}%
and we get the equalities below which we already got.%
\begin{equation*}
B(1_{H}\otimes gx_{2};X_{1}X_{2},gx_{1})=0\text{ ,}
\end{equation*}%
\begin{equation*}
-B(1_{H}\otimes gx_{2};GX_{1}X_{2},x_{1})=0\text{.}
\end{equation*}

\subsubsection{$G^{a}X_{2}\otimes g^{d}x_{1}x_{2}$}

\begin{equation*}
\sum_{\substack{ a,d=0  \\ a+d\equiv 0}}^{1}-B(1_{H}\otimes
gx_{2};G^{a}X_{1}X_{2},g^{d}x_{1}x_{2})G^{a}X_{2}\otimes g^{d}x_{1}x_{2}
\end{equation*}%
and we get the equalities below which we already got.%
\begin{equation*}
-B(1_{H}\otimes gx_{2};X_{1}X_{2},x_{1}x_{2})=0,\text{ }
\end{equation*}%
\begin{equation*}
-B(1_{H}\otimes gx_{2};GX_{1}X_{2},gx_{1}x_{2})=0.\text{ }
\end{equation*}

\subsubsection{$G^{a}X_{1}\otimes g^{d}x_{1}x_{2}$}

There is no term like this

\subsubsection{$G^{a}X_{1}X_{2}\otimes g^{d}x_{2}$}

\begin{equation*}
\sum_{\substack{ a,d=0  \\ a+d\equiv 0}}^{1}\left( -1\right)
^{a+1}B(1_{H}\otimes
gx_{2};G^{a}X_{1}X_{2},g^{d}x_{1}x_{2})G^{a}X_{1}X_{2}\otimes g^{d}x_{2}=0
\end{equation*}%
and we get the equalities below which we already got.%
\begin{equation*}
-B(1_{H}\otimes gx_{2};X_{1}X_{2},x_{1}x_{2})\text{ ,}
\end{equation*}%
\begin{equation*}
B(1_{H}\otimes gx_{2};GX_{1}X_{2},gx_{1}x_{2})\text{ .}
\end{equation*}

\subsubsection{$G^{a}X_{1}X_{2}\otimes g^{d}x_{1}$}

There is no term like this.

\subsubsection{$G^{a}X_{1}X_{2}\otimes g^{d}x_{1}x_{2}$}

There is no term like this.

\subsection{case $gx_{2}$}

\begin{eqnarray*}
&&\sum_{a,b_{1},b_{2},d,e_{1},e_{2}=0}^{1}\sum_{l_{1}=0}^{b_{1}}%
\sum_{l_{2}=0}^{b_{2}}\sum_{u_{1}=0}^{e_{1}}\sum_{u_{2}=0}^{e_{2}}\left(
-1\right) ^{\alpha \left( 1_{H};l_{1},l_{2},u_{1},u_{2}\right) } \\
&&B(1_{H}\otimes
gx_{2};G^{a}X_{1}^{b_{1}}X_{2}^{b_{2}},g^{d}x_{1}^{e_{1}}x_{2}^{e_{2}}) \\
&&G^{a}X_{1}^{b_{1}-l_{1}}X_{2}^{b_{2}-l_{2}}\otimes
g^{d}x_{1}^{e_{1}-u_{1}}x_{2}^{e_{2}-u_{2}}\otimes \\
&&g^{a+b_{1}+b_{2}+l_{1}+l_{2}+d+e_{1}+e_{2}+u_{1}+u_{2}}x_{1}^{l_{1}+u_{1}}x_{2}^{l_{2}+u_{2}}
\\
&=&\sum_{\omega _{2}=0}^{1}B^{A}(1_{H}\otimes gx_{2}^{1-\omega _{2}})\otimes
B^{H}(1_{H}\otimes gx_{2}^{1-\omega _{2}})\otimes g^{\omega
_{2}}x_{2}^{\omega _{2}}
\end{eqnarray*}%
\begin{eqnarray*}
a+b_{1}+b_{2}+l_{1}+l_{2}+d+e_{1}+e_{2}+u_{1}+u_{2} &=&1 \\
l_{1}+u_{1} &=&0\Rightarrow l_{1}=u_{1}=0 \\
l_{2}+u_{2} &=&1
\end{eqnarray*}%
\begin{eqnarray*}
a+b_{1}+b_{2}+d+e_{1}+e_{2} &=&0 \\
l_{1}+u_{1} &=&0\Rightarrow l_{1}=u_{1}=0 \\
l_{2}+u_{2} &=&1
\end{eqnarray*}

\begin{eqnarray*}
&&\sum_{\substack{ a,b_{1},b_{2},d,e_{1},e_{2}=0  \\ %
a+b_{1}+b_{2}+d+e_{1}+e_{2}=0}}^{1}\sum_{l_{2}=0}^{b_{2}}\sum_{\substack{ %
u_{2}=0  \\ l_{2}+u_{2}=1}}^{e_{2}}\left( -1\right) ^{\alpha \left(
1_{H};0,l_{2},0,u_{2}\right) }B(1_{H}\otimes
gx_{2};G^{a}X_{1}^{b_{1}}X_{2}^{b_{2}},g^{d}x_{1}^{e_{1}}x_{2}^{e_{2}}) \\
&&G^{a}X_{1}^{b_{1}}X_{2}^{b_{2}-l_{2}}\otimes
g^{d}x_{1}^{e_{1}}x_{2}^{e_{2}-u_{2}}\otimes gx_{2}=B^{A}(1_{H}\otimes
g)\otimes B^{H}(1_{H}\otimes g)\otimes gx_{2}
\end{eqnarray*}%
\begin{eqnarray*}
&&\sum_{\substack{ a,b_{1},b_{2},d,e_{1},e_{2}=0  \\ a+b_{1}+b_{2}+d+e_{1}=1
}}^{1}\left( -1\right) ^{\left( a+b_{1}+b_{2}\right) }B(1_{H}\otimes
gx_{2};G^{a}X_{1}^{b_{1}}X_{2}^{b_{2}},g^{d}x_{1}^{e_{1}}x_{2})G^{a}X_{1}^{b_{1}}X_{2}^{b_{2}}\otimes g^{d}x_{1}^{e_{1}}+
\\
&&\sum_{\substack{ a,b_{1},b_{2},d,e_{1},e_{2}=0  \\ a+b_{1}+d+e_{1}+e_{2}=1
}}^{1}B(1_{H}\otimes
gx_{2};G^{a}X_{1}^{b_{1}}X_{2},g^{d}x_{1}^{e_{1}}x_{2}^{e_{2}})G^{a}X_{1}^{b_{1}}\otimes g^{d}x_{1}^{e_{1}}x_{2}^{e_{2}}=
\\
&=&B^{A}(1_{H}\otimes g)\otimes B^{H}(1_{H}\otimes g)
\end{eqnarray*}

\subsubsection{$G^{a}\otimes g^{d}$}

The left side is%
\begin{equation*}
\sum_{\substack{ a,d=0  \\ a+d=1}}^{1}\left[ \left( -1\right)
^{a}B(1_{H}\otimes gx_{2};G^{a},g^{d}x_{2})+B(1_{H}\otimes
gx_{2};G^{a}X_{2},g^{d})\right] G^{a}\otimes g^{d}
\end{equation*}%
so that we obtain

\begin{eqnarray*}
B(1_{H}\otimes gx_{2};1_{A},gx_{2})+B(1_{H}\otimes gx_{2};X_{2},g)
&=&B\left( 1_{H}\otimes g;1_{A},g\right) \\
-B(1_{H}\otimes gx_{2};G,x_{2})+B(1_{H}\otimes gx_{2};GX_{2},1_{H})
&=&B\left( 1_{H}\otimes g;G,1_{H}\right)
\end{eqnarray*}%
By the form of the elements $B(1_{H}\otimes gx_{2})$ and $B\left( g\otimes
1_{H}\right) ,$ these equalities are satisfied.

\subsubsection{$G^{a}\otimes g^{d}x_{2}$}

The left side is%
\begin{equation*}
\sum_{\substack{ a,d=0  \\ a+d=0}}^{1}B(1_{H}\otimes
gx_{2};G^{a}X_{2},g^{d}x_{2})G^{a}\otimes g^{d}x_{2}
\end{equation*}

so that we obtain%
\begin{equation*}
B(1_{H}\otimes gx_{2};X_{2},x_{2})=B\left( 1_{H}\otimes g;1_{A},x_{2}\right)
\end{equation*}%
and%
\begin{equation*}
B(1_{H}\otimes gx_{2};GX_{2},gx_{2})=B\left( 1_{H}\otimes g;G,gx_{2}\right) .
\end{equation*}%
By the form of the elements $B(1_{H}\otimes gx_{2})$ and $B\left( g\otimes
1_{H}\right) ,$ these equalities are satisfied.

\subsubsection{$G^{a}\otimes g^{d}x_{1}$}

The left side is%
\begin{equation*}
\sum_{\substack{ a,d=0  \\ a+d=0}}^{1}\left[ \left( -1\right)
^{a}B(1_{H}\otimes gx_{2};G^{a},g^{d}x_{1}x_{2})+B(1_{H}\otimes
gx_{2};G^{a}X_{2},g^{d}x_{1})\right] G^{a}\otimes g^{d}x_{1}
\end{equation*}%
and we get%
\begin{equation*}
B(1_{H}\otimes gx_{2};1_{A},x_{1}x_{2})+B(1_{H}\otimes
gx_{2};X_{2},x_{1})=B\left( 1_{H}\otimes g;1_{A},x_{1}\right) ,
\end{equation*}%
\begin{equation*}
-B(1_{H}\otimes gx_{2};G,gx_{1}x_{2})+B(1_{H}\otimes
gx_{2};GX_{2},gx_{1})=B\left( 1_{H}\otimes g;G,gx_{1}\right) .
\end{equation*}

By the form of the elements $B(1_{H}\otimes gx_{2})$ and $B\left( g\otimes
1_{H}\right) ,$ these equalities are satisfied.

\subsubsection{$G^{a}X_{2}\otimes g^{d}$}

The left side is%
\begin{equation*}
\sum_{\substack{ a,d=0  \\ a+d=0}}^{1}\left( -1\right) ^{\left( a+1\right)
}B(1_{H}\otimes gx_{2};G^{a}X_{2},g^{d}x_{2})G^{a}X_{2}\otimes g^{d}
\end{equation*}%
and we get%
\begin{equation*}
-B(1_{H}\otimes gx_{2};X_{2},x_{2})=B\left( 1_{H}\otimes g;X_{2},1_{H}\right)
\end{equation*}%
and%
\begin{equation*}
+B(1_{H}\otimes gx_{2};GX_{2},gx_{2})=B\left( 1_{H}\otimes g;GX_{2},g\right)
.
\end{equation*}%
By the form of the elements $B(1_{H}\otimes gx_{2})$ and $B\left( g\otimes
1_{H}\right) ,$ these equalities are satisfied.

\subsubsection{$G^{a}X_{1}\otimes g^{d}$}

The left side is%
\begin{equation*}
\sum_{\substack{ a,d=0  \\ a+d=0}}^{1}\left[ \left( -1\right) ^{\left(
a+1\right) }B(1_{H}\otimes gx_{2};G^{a}X_{1},g^{d}x_{2})+B(1_{H}\otimes
gx_{2};G^{a}X_{1}X_{2},g^{d})\right] G^{a}X_{1}\otimes g^{d}+
\end{equation*}%
and we get%
\begin{equation*}
-B(1_{H}\otimes gx_{2};X_{1},x_{2})+B(1_{H}\otimes
gx_{2};X_{1}X_{2},1_{H})=B\left( 1_{H}\otimes g;X_{1},1_{H}\right)
\end{equation*}%
and%
\begin{equation*}
B(1_{H}\otimes gx_{2};GX_{1},gx_{2})+B(1_{H}\otimes
gx_{2};GX_{1}X_{2},g)=B\left( 1_{H}\otimes g;GX_{1},g\right) .
\end{equation*}%
By the form of the elements $B(1_{H}\otimes gx_{2})$ and $B\left( g\otimes
1_{H}\right) ,$ these equalities are satisfied.

\subsubsection{$G^{a}\otimes g^{d}x_{1}x_{2}$}

\begin{equation*}
\sum_{\substack{ a,d=0  \\ a+d=1}}^{1}B(1_{H}\otimes
gx_{2};G^{a}X_{2},g^{d}x_{1}x_{2})G^{a}\otimes g^{d}x_{1}x_{2}
\end{equation*}%
We get%
\begin{equation*}
B(1_{H}\otimes gx_{2};X_{2},gx_{1}x_{2})=B\left( 1_{H}\otimes
g,1_{H},gx_{1}x_{2}\right)
\end{equation*}%
and%
\begin{equation*}
B(1_{H}\otimes gx_{2};GX_{2},x_{1}x_{2})=B\left( 1_{H}\otimes
g,G,x_{1}x_{2}\right) .
\end{equation*}%
By the form of the elements $B(1_{H}\otimes gx_{2})$ and $B\left( g\otimes
1_{H}\right) ,$ these equalities are satisfied.

\subsubsection{$G^{a}X_{2}\otimes g^{d}x_{2}$}

There is no term in the first side.

\subsubsection{$G^{a}X_{1}\otimes g^{d}x_{2}$}

\begin{equation*}
\sum_{\substack{ a,d=0  \\ a+d=1}}^{1}B(1_{H}\otimes
gx_{2};G^{a}X_{1}X_{2},g^{d}x_{2})G^{a}X_{1}\otimes g^{d}x_{2}
\end{equation*}%
We get%
\begin{equation*}
B(1_{H}\otimes gx_{2};X_{1}X_{2},gx_{2})=B\left( 1_{H}\otimes
g;X_{1},gx_{2}\right)
\end{equation*}%
and%
\begin{equation*}
B(1_{H}\otimes gx_{2};GX_{1}X_{2},x_{2})=B\left( 1_{H}\otimes
g;GX_{1},x_{2}\right) .
\end{equation*}%
By the form of the elements $B(1_{H}\otimes gx_{2})$ and $B\left( g\otimes
1_{H}\right) ,$ these equalities are satisfied.

\subsubsection{$G^{a}X_{2}\otimes g^{d}x_{1}$}

\begin{equation*}
\sum_{\substack{ a,d=0  \\ a+d=1}}^{1}\left( -1\right) ^{\left( a+1\right)
}B(1_{H}\otimes gx_{2};G^{a}X_{2},g^{d}x_{1}x_{2})G^{a}X_{2}\otimes
g^{d}x_{1}
\end{equation*}%
We get%
\begin{equation*}
-B(1_{H}\otimes gx_{2};X_{2},gx_{1}x_{2})=B\left( 1_{H}\otimes
g;X_{2},gx_{1}\right)
\end{equation*}%
and%
\begin{equation*}
B(1_{H}\otimes gx_{2};GX_{2},x_{1}x_{2})=B\left( 1_{H}\otimes
g;GX_{2},x_{1}\right) .
\end{equation*}%
By the form of the elements $B(1_{H}\otimes gx_{2})$ and $B\left( g\otimes
1_{H}\right) ,$ these equalities are satisfied.

\subsubsection{$G^{a}X_{1}\otimes g^{d}x_{1}$}

\begin{equation*}
\sum_{\substack{ a,d=0  \\ a+d=1}}^{1}\left[ \left( -1\right) ^{\left(
a+1\right) }B(1_{H}\otimes gx_{2};G^{a}X_{1},g^{d}x_{1}x_{2})+B(1_{H}\otimes
gx_{2};G^{a}X_{1}X_{2},g^{d}x_{1})\right] G^{a}X_{1}\otimes g^{d}x_{1}
\end{equation*}%
and we get%
\begin{equation*}
-B(1_{H}\otimes gx_{2};X_{1},gx_{1}x_{2})+B(1_{H}\otimes
gx_{2};X_{1}X_{2},gx_{1})=B\left( 1_{H}\otimes g;X_{1},gx_{1}\right)
\end{equation*}%
and
\begin{equation*}
B(1_{H}\otimes gx_{2};GX_{1},x_{1}x_{2})+B(1_{H}\otimes
gx_{2};GX_{1}X_{2},x_{1})=B\left( 1_{H}\otimes g;GX_{1},x_{1}\right) .
\end{equation*}%
By the form of the elements $B(1_{H}\otimes gx_{2})$ and $B\left( g\otimes
1_{H}\right) ,$ these equalities are satisfied.

\subsubsection{$G^{a}X_{1}X_{2}\otimes g^{d}$}

\begin{equation*}
\sum_{\substack{ a,d=0  \\ a+d=1}}^{1}\left( -1\right) ^{a}B(1_{H}\otimes
gx_{2};G^{a}X_{1}X_{2},g^{d}x_{2})G^{a}X_{1}X_{2}\otimes g^{d}
\end{equation*}%
and we get%
\begin{equation*}
B(1_{H}\otimes gx_{2};X_{1}X_{2},gx_{2})=B\left( 1_{H}\otimes
g;X_{1}X_{2},g\right)
\end{equation*}%
and%
\begin{equation*}
-B(1_{H}\otimes gx_{2};GX_{1}X_{2},x_{2})=B\left( 1_{H}\otimes
g;GX_{1}X_{2},1_{H}\right)
\end{equation*}%
By the form of the elements $B(1_{H}\otimes gx_{2})$ and $B\left( g\otimes
1_{H}\right) ,$ these equalities are satisfied.

\subsubsection{$G^{a}X_{2}\otimes g^{d}x_{1}x_{2}$}

There is no term like this in the first side.

\subsubsection{$G^{a}X_{1}\otimes g^{d}x_{1}x_{2}$}

\begin{equation*}
\sum_{\substack{ a,d=0  \\ a+d=0}}^{1}B(1_{H}\otimes
gx_{2};G^{a}X_{1}X_{2},g^{d}x_{1}x_{2})G^{a}X_{1}\otimes g^{d}x_{1}x_{2}
\end{equation*}%
and we get%
\begin{equation*}
B(1_{H}\otimes gx_{2};X_{1}X_{2},x_{1}x_{2})=B\left( 1_{H}\otimes
g;X_{1},x_{1}x_{2}\right)
\end{equation*}%
and%
\begin{equation*}
B(1_{H}\otimes gx_{2};GX_{1}X_{2},gx_{1}x_{2})=B\left( 1_{H}\otimes
g;GX_{1},gx_{1}x_{2}\right) .
\end{equation*}%
By the form of the elements $B(1_{H}\otimes gx_{2})$ and $B\left( g\otimes
1_{H}\right) ,$ these equalities are satisfied.

\subsubsection{$G^{a}X_{1}X_{2}\otimes g^{d}x_{2}$}

There is no term like this in the left side.

\subsubsection{$G^{a}X_{1}X_{2}\otimes g^{d}x_{1}$}

\begin{equation*}
\sum_{\substack{ a,d=0  \\ a+d=0}}^{1}\left( -1\right) ^{a}B(1_{H}\otimes
gx_{2};G^{a}X_{1}X_{2},g^{d}x_{1}x_{2})G^{a}X_{1}X_{2}\otimes g^{d}x_{1}
\end{equation*}%
we get%
\begin{equation*}
B(1_{H}\otimes gx_{2};X_{1}X_{2},x_{1}x_{2})=B\left( 1_{H}\otimes
g;X_{1}X_{2},x_{1}\right)
\end{equation*}%
and%
\begin{equation*}
-B(1_{H}\otimes gx_{2};GX_{1}X_{2},gx_{1}x_{2})=B\left( 1_{H}\otimes
g;GX_{1}X_{2},gx_{1}\right) .
\end{equation*}%
By the form of the elements $B(1_{H}\otimes gx_{2})$ and $B\left( g\otimes
1_{H}\right) ,$ these equalities are satisfied.

\subsubsection{$G^{a}X_{1}X_{2}\otimes g^{d}x_{1}x_{2}$}

There is no term like this in the first side.

\subsection{case $gx_{1}x_{2}$}

\begin{eqnarray*}
&&\sum_{\substack{ a,b_{1},b_{2},d,e_{1},e_{2}=0  \\ %
a+b_{1}+b_{2}+d+e_{1}+e_{2}=1}}^{1}\sum_{l_{1}=0}^{b_{1}}%
\sum_{l_{2}=0}^{b_{2}}\sum_{\substack{ u_{1}=0  \\ l_{1}+u_{1}=1}}%
^{e_{1}}\sum _{\substack{ u_{2}=0  \\ l_{2}+u_{2}=1}}^{e_{2}}\left(
-1\right) ^{\alpha \left( 1_{H};l_{1},l_{2},u_{1},u_{2}\right) } \\
&&B(1_{H}\otimes
gx_{2};G^{a}X_{1}^{b_{1}}X_{2}^{b_{2}},g^{d}x_{1}^{e_{1}}x_{2}^{e_{2}}) \\
&&G^{a}X_{1}^{b_{1}-l_{1}}X_{2}^{b_{2}-l_{2}}\otimes
g^{d}x_{1}^{e_{1}-u_{1}}x_{2}^{e_{2}-u_{2}}\otimes \\
&&g^{a+b_{1}+b_{2}+d+e_{1}+e_{2}}x_{1}^{l_{1}+u_{1}}x_{2}^{l_{2}+u_{2}} \\
&=&\sum_{\omega _{2}=0}^{1}B^{A}(1_{H}\otimes gx_{2}^{1-\omega _{2}})\otimes
B^{H}(1_{H}\otimes gx_{2}^{1-\omega _{2}})\otimes g^{\omega
_{2}}x_{2}^{\omega _{2}}
\end{eqnarray*}%
Since%
\begin{equation*}
B(1_{H}\otimes
gx_{2};G^{a}X_{1}^{b_{1}}X_{2}^{b_{2}},g^{d}x_{1}^{e_{1}}x_{2}^{e_{2}})=0%
\text{ whenever }a+b_{1}+b_{2}+d+e_{1}+e_{2}\equiv 1
\end{equation*}%
we have nothing to check.

\section{$B\left( 1_{H}\otimes x_{1}x_{2}\right) $}

For $h=1_{H}$ in $\left( \ref{simpl1}\right) $ and by using $\left( \ref%
{form:1otxi}\right) $ for $i=2,$ we get%
\begin{equation*}
B(1_{H}\otimes x_{1}x_{2})=B(1_{H}\otimes x_{2})(1_{A}\otimes
x_{1})-(1_{A}\otimes gx_{1})B(1_{H}\otimes x_{2})(1_{A}\otimes
g)-(1_{A}\otimes g)B(gx_{1}\otimes x_{2})(1_{A}\otimes g)
\end{equation*}

and%
\begin{equation*}
B(1_{H}\otimes x_{2})=-(1_{A}\otimes g)B(gx_{2}\otimes 1_{H})(1_{A}\otimes
g).
\end{equation*}%
For $h=gx_{1}$ in $\left( \ref{simpl2i2}\right) $ we get%
\begin{equation*}
B(gx_{1}\otimes x_{2})=B(gx_{1}\otimes 1_{H})(1_{A}\otimes
x_{2})-(1_{A}\otimes gx_{2})B(gx_{1}\otimes 1_{H})(1_{A}\otimes
g)-(1_{A}\otimes g)B(x_{1}x_{2}\otimes 1_{H})(1_{A}\otimes g)
\end{equation*}%
Finally we get%
\begin{eqnarray}
&&B(1_{H}\otimes x_{1}x_{2})  \label{form 1otx1x2} \\
&=&-(1_{A}\otimes g)B(gx_{2}\otimes 1_{H})(1_{A}\otimes g)(1_{A}\otimes
x_{1})  \notag \\
&&+(1_{A}\otimes gx_{1})(1_{A}\otimes g)B(gx_{2}\otimes 1_{H})(1_{A}\otimes
g)(1_{A}\otimes g)  \notag \\
&&-(1_{A}\otimes g)B(gx_{1}\otimes 1_{H})(1_{A}\otimes x_{2})(1_{A}\otimes g)
\notag \\
&&+(1_{A}\otimes g)(1_{A}\otimes gx_{2})B(gx_{1}\otimes 1_{H})(1_{A}\otimes
g)(1_{A}\otimes g)  \notag \\
&&+(1_{A}\otimes g)(1_{A}\otimes g)B(x_{1}x_{2}\otimes 1_{H})(1_{A}\otimes
g)(1_{A}\otimes g)  \notag
\end{eqnarray}

and we obtain%
\begin{eqnarray*}
&&B(1_{H}\otimes x_{1}x_{2}) \\
&=&B(x_{1}x_{2}\otimes 1_{H};1_{A},1_{H})1_{A}\otimes 1_{H}+ \\
&&+B(x_{1}x_{2}\otimes 1_{H};1_{A},x_{1}x_{2})1_{A}\otimes x_{1}x_{2}+ \\
&&+\left[ -B(x_{1}x_{2}\otimes 1_{H};1_{A},gx_{1})-2B(x_{1}x_{2}\otimes
1_{H};X_{1},g)\right] 1_{A}\otimes gx_{1} \\
&&+\left[ -B(x_{1}x_{2}\otimes 1_{H};1_{A},gx_{2})-2B(x_{1}x_{2}\otimes
1_{H};X_{2},g)\right] \left( 1_{A}\otimes gx_{2}\right) + \\
&&+B(x_{1}x_{2}\otimes 1_{H};G,g)G\otimes g+ \\
&&+B(x_{1}x_{2}\otimes 1_{H};G,x_{1})G\otimes x_{1}+B(x_{1}x_{2}\otimes
1_{H};G,x_{2})G\otimes x_{2}+ \\
&&+\left[
\begin{array}{c}
-3B(x_{1}x_{2}\otimes 1_{H};G,gx_{1}x_{2})- \\
2B(x_{1}x_{2}\otimes 1_{H};GX_{1},gx_{2})+2B(x_{1}x_{2}\otimes
1_{H};GX_{2},gx_{1})%
\end{array}%
\right] G\otimes gx_{1}x_{2} \\
&&+B(x_{1}x_{2}\otimes 1_{H};X_{1},g)X_{1}\otimes g+ \\
&&+B(x_{1}x_{2}\otimes 1_{H};X_{1},x_{1})X_{1}\otimes
x_{1}+B(x_{1}x_{2}\otimes 1_{H};X_{1},x_{2})X_{1}\otimes x_{2}+ \\
&&-B(x_{1}x_{2}\otimes 1_{H};X_{1},gx_{1}x_{2})X_{1}\otimes
gx_{1}x_{2}+B(x_{1}x_{2}\otimes 1_{H};X_{2},g)X_{2}\otimes g+ \\
&&+B(x_{1}x_{2}\otimes 1_{H};X_{2},x_{1})X_{2}\otimes
x_{1}+B(x_{1}x_{2}\otimes 1_{H};X_{2},x_{2})X_{2}\otimes x_{2}+ \\
&&-B(x_{1}x_{2}\otimes 1_{H};X_{2},gx_{1}x_{2})X_{2}\otimes gx_{1}x_{2}+ \\
&&+\left[
\begin{array}{c}
+1-B(x_{1}x_{2}\otimes 1_{H};1_{A},x_{1}x_{2}) \\
-B(x_{1}x_{2}\otimes 1_{H};X_{2},x_{1})+B(x_{1}x_{2}\otimes
1_{H};X_{1},x_{2})%
\end{array}%
\right] \\
&&X_{1}X_{2}\otimes 1_{H}+ \\
&&-B(x_{1}x_{2}\otimes 1_{H};X_{1},gx_{1}x_{2})X_{1}X_{2}\otimes gx_{1} \\
&&-B(x_{1}x_{2}\otimes 1_{H};X_{2},gx_{1}x_{2})X_{1}X_{2}\otimes gx_{2} \\
&&+B(x_{1}x_{2}\otimes 1_{H};GX_{1},1_{H})GX_{1}\otimes 1_{H}+ \\
&&+B(x_{1}x_{2}\otimes 1_{H};GX_{1},x_{1}x_{2})GX_{1}\otimes x_{1}x_{2} \\
&&-B(x_{1}x_{2}\otimes 1_{H};GX_{1},gx_{1})GX_{1}\otimes gx_{1} \\
&&\left[ B(x_{1}x_{2}\otimes 1_{H};GX_{1},gx_{2})+2B(x_{1}x_{2}\otimes
1_{H};G,gx_{1}x_{2})-2B(x_{1}x_{2}\otimes 1_{H};GX_{2},gx_{1})\right]
GX_{1}\otimes gx_{2} \\
&&+B(x_{1}x_{2}\otimes 1_{H};GX_{2},1_{H})GX_{2}\otimes 1_{H}+ \\
&&+B(x_{1}x_{2}\otimes 1_{H};GX_{2},x_{1}x_{2})GX_{2}\otimes x_{1}x_{2}+ \\
&&+\left[
\begin{array}{c}
B(x_{1}x_{2}\otimes 1_{H};GX_{2},gx_{1})-2B(x_{1}x_{2}\otimes
1_{H};G,gx_{1}x_{2}) \\
-2B(x_{1}x_{2}\otimes 1_{H};GX_{1},gx_{2})%
\end{array}%
\right] GX_{2}\otimes gx_{1} \\
&&-B(x_{1}x_{2}\otimes 1_{H};GX_{2},gx_{2})GX_{2}\otimes gx_{2}GX_{2}\otimes
gx_{1} \\
&&+\left[
\begin{array}{c}
-B(x_{1}x_{2}\otimes 1_{H};G,gx_{1}x_{2})+B(x_{1}x_{2}\otimes
1_{H};GX_{2},gx_{1}) \\
-B(x_{1}x_{2}\otimes 1_{H};GX_{1},gx_{2})%
\end{array}%
\right] GX_{1}X_{2}\otimes g \\
&&-B(x_{1}x_{2}\otimes 1_{H};GX_{1},x_{1}x_{2})GX_{1}X_{2}\otimes x_{1}+ \\
&&-B(x_{1}x_{2}\otimes 1_{H};GX_{2},x_{1}x_{2})GX_{1}X_{2}\otimes x_{2}
\end{eqnarray*}%
We write Casimir condition $\left( \ref{MAIN FORMULA 1}\right) $ for $%
B(1_{A}\otimes x_{1}x_{2}).$
\begin{eqnarray*}
&&\sum_{a,b_{1},b_{2},d,e_{1},e_{2}=0}^{1}\sum_{l_{1}=0}^{b_{1}}%
\sum_{l_{2}=0}^{b_{2}}\sum_{u_{1}=0}^{e_{1}}\sum_{u_{2}=0}^{e_{2}}\left(
-1\right) ^{\alpha \left( 1_{H};l_{1},l_{2},u_{1},u_{2}\right) } \\
&&B(1_{H}\otimes
x_{1}x_{2};G^{a}X_{1}^{b_{1}}X_{2}^{b_{2}},g^{d}x_{1}^{e_{1}}x_{2}^{e_{2}})G^{a}X_{1}^{b_{1}-l_{1}}X_{2}^{b_{2}-l_{2}}\otimes g^{d}x_{1}^{e_{1}-u_{1}}x_{2}^{e_{2}-u_{2}}\otimes
\\
&&g^{a+b_{1}+b_{2}+l_{1}+l_{2}+d+e_{1}+e_{2}+u_{1}+u_{2}}x_{1}^{l_{1}+u_{1}}x_{2}^{l_{2}+u_{2}}
\\
&=&B^{A}(1_{A}\otimes x_{1}x_{2})\otimes B^{H}(1_{A}\otimes
x_{1}x_{2})\otimes 1_{H} \\
&&B^{A}(1_{A}\otimes x_{1})\otimes B^{H}(1_{A}\otimes x_{1})\otimes gx_{2} \\
&&-B^{A}(1_{A}\otimes x_{2})\otimes B^{H}1_{A}\otimes x_{2})\otimes gx_{1} \\
&&+B^{A}(1_{A}\otimes 1_{H})\otimes B^{H}(1_{A}\otimes 1_{H})\otimes
x_{1}x_{2}
\end{eqnarray*}

\subsection{Case $1_{H}$}

Considering the left side we get
\begin{eqnarray*}
l_{1} &=&u_{1}=0 \\
l_{2} &=&u_{2}=0 \\
a+b_{1}+b_{2}+d+e_{1}+e_{2} &\equiv &0
\end{eqnarray*}%
and we get that%
\begin{equation}
B(1_{H}\otimes
x_{1}x_{2};G^{a}X_{1}^{b_{1}}X_{2}^{b_{2}},g^{d}x_{1}^{e_{1}}x_{2}^{e_{2}})=0%
\text{ whenever }a+b_{1}+b_{2}+d+e_{1}+e_{2}\equiv 1  \label{1otx1x2, first}
\end{equation}%
Thus we obtain%
\begin{eqnarray*}
&&\sum_{\substack{ a,b_{1},b_{2},d,e_{1},e_{2}=0  \\ %
a+b_{1}+b_{2}+d+e_{1}+e_{2}\equiv 0}}^{1}\sum_{l_{1}=0}^{b_{1}}%
\sum_{l_{2}=0}^{b_{2}}\sum_{u_{1}=0}^{e_{1}}\sum_{u_{2}=0}^{e_{2}}\left(
-1\right) ^{\alpha \left( 1_{H};l_{1},l_{2},u_{1},u_{2}\right) } \\
&&B(1_{H}\otimes
x_{1}x_{2};G^{a}X_{1}^{b_{1}}X_{2}^{b_{2}},g^{d}x_{1}^{e_{1}}x_{2}^{e_{2}})G^{a}X_{1}^{b_{1}-l_{1}}X_{2}^{b_{2}-l_{2}}\otimes g^{d}x_{1}^{e_{1}-u_{1}}x_{2}^{e_{2}-u_{2}}\otimes
\\
&&g^{a+b_{1}+b_{2}+l_{1}+l_{2}+d+e_{1}+e_{2}+u_{1}+u_{2}}x_{1}^{l_{1}+u_{1}}x_{2}^{l_{2}+u_{2}}
\\
&=&B^{A}(1_{A}\otimes x_{1}x_{2})\otimes B^{H}(1_{A}\otimes
x_{1}x_{2})\otimes 1_{H} \\
&&B^{A}(1_{A}\otimes x_{1})\otimes B^{H}(1_{A}\otimes x_{1})\otimes gx_{2} \\
&&-B^{A}(1_{A}\otimes x_{2})\otimes B^{H}1_{A}\otimes x_{2})\otimes gx_{1} \\
&&+B^{A}(1_{A}\otimes 1_{H})\otimes B^{H}(1_{A}\otimes 1_{H})\otimes
x_{1}x_{2}
\end{eqnarray*}

\subsection{Case $g$}

Gives no other information

\subsection{Case $x_{1}$}

By considering the left side we get%
\begin{eqnarray*}
&&\sum_{a,b_{1},b_{2},d,e_{1},e_{2}=0}^{1}\sum_{l_{1}=0}^{b_{1}}%
\sum_{l_{2}=0}^{b_{2}}\sum_{u_{1}=0}^{e_{1}}\sum_{u_{2}=0}^{e_{2}}\left(
-1\right) ^{\alpha \left( 1_{H};l_{1},l_{2},u_{1},u_{2}\right) } \\
&&B(1_{H}\otimes
x_{1}x_{2};G^{a}X_{1}^{b_{1}}X_{2}^{b_{2}},g^{d}x_{1}^{e_{1}}x_{2}^{e_{2}})G^{a}X_{1}^{b_{1}-l_{1}}X_{2}^{b_{2}-l_{2}}\otimes g^{d}x_{1}^{e_{1}-u_{1}}x_{2}^{e_{2}-u_{2}}\otimes
\\
&&g^{a+b_{1}+b_{2}+l_{1}+l_{2}+d+e_{1}+e_{2}+u_{1}+u_{2}}x_{1}^{l_{1}+u_{1}}x_{2}^{l_{2}+u_{2}}
\end{eqnarray*}%
\begin{eqnarray*}
l_{1}+u_{1} &=&1 \\
l_{2} &=&u_{2}=0 \\
a+b_{1}+b_{2}+d+e_{1}+e_{2} &\equiv &1
\end{eqnarray*}

\subsection{Case $x_{2}$}

By considering the left side we get%
\begin{eqnarray*}
l_{1} &=&u_{1}=0 \\
l_{2}+u_{2} &=&1 \\
a+b_{1}+b_{2}+d+e_{1}+e_{2} &\equiv &1
\end{eqnarray*}%
{\Huge \ }

Since%
\begin{equation*}
B(1_{H}\otimes
x_{1}x_{2};G^{a}X_{1}^{b_{1}}X_{2}^{b_{2}},g^{d}x_{1}^{e_{1}}x_{2}^{e_{2}})=0%
\text{ whenever }a+b_{1}+b_{2}+d+e_{1}+e_{2}\equiv 1
\end{equation*}%
we get no new information.

\subsection{Case $x_{1}x_{2}$}

By considering the left side we get%
\begin{eqnarray*}
l_{1}+u_{1} &=&1 \\
l_{2}+u_{2} &=&1 \\
a+b_{1}+b_{2}+d+e_{1}+e_{2} &\equiv &0
\end{eqnarray*}%
\begin{gather*}
\sum_{\substack{ a,b_{1},b_{2},d,e_{1},e_{2}=0  \\ %
a+b_{1}+b_{2}+d+e_{1}+e_{2}\equiv 0}}^{1}\sum_{l_{1}=0}^{b_{1}}%
\sum_{l_{2}=0}^{b_{2}}\sum_{\substack{ u_{1}=0  \\ l_{1}+u_{1}=1}}%
^{e_{1}}\sum _{\substack{ u_{2}=0  \\ l_{2}+u_{2}=1}}^{e_{2}}\left(
-1\right) ^{\alpha \left( 1_{H};l_{1},l_{2},u_{1},u_{2}\right) } \\
B(1_{H}\otimes
x_{1}x_{2};G^{a}X_{1}^{b_{1}}X_{2}^{b_{2}},g^{d}x_{1}^{e_{1}}x_{2}^{e_{2}})G^{a}X_{1}^{b_{1}-l_{1}}X_{2}^{b_{2}-l_{2}}\otimes g^{d}x_{1}^{e_{1}-u_{1}}x_{2}^{e_{2}-u_{2}}
\end{gather*}%
and hence%
\begin{eqnarray*}
&&\sum_{\substack{ a,b_{1},b_{2},d=0  \\ a+b_{1}+b_{2}+d\equiv 0}}^{1}\left(
-1\right) ^{\alpha \left( 1_{H};0,0,1,1\right) }B(1_{H}\otimes
x_{1}x_{2};G^{a}X_{1}^{b_{1}}X_{2}^{b_{2}},g^{d}x_{1}x_{2}) \\
&&G^{a}X_{1}^{b_{1}}X_{2}^{b_{2}}\otimes g^{d}x_{1}^{e_{1}-1}x_{2}^{e_{2}-1}+
\\
&&\sum_{\substack{ a,b_{1},b_{2},d,e_{1},e_{2}=0  \\ %
a+b_{1}+b_{2}+d+e_{1}+e_{2}\equiv 0}}^{1}\left( -1\right) ^{\alpha \left(
1_{H};0,1,1,0\right) }B(1_{H}\otimes
x_{1}x_{2};G^{a}X_{1}^{b_{1}}X_{2},g^{d}x_{1}x_{2}^{e_{2}})G^{a}X_{1}^{b_{1}}X_{2}^{b_{2}-1}\otimes g^{d}x_{1}^{e_{1}-1}x_{2}^{e_{2}}+
\\
&&\sum_{\substack{ a,b_{1},b_{2},d,e_{1},e_{2}=0  \\ %
a+b_{1}+b_{2}+d+e_{1}+e_{2}\equiv 0}}^{1}\left( -1\right) ^{\alpha \left(
1_{H};1,0,0,1\right) }B(1_{H}\otimes
x_{1}x_{2};G^{a}X_{1}X_{2}^{b_{2}},g^{d}x_{1}x_{2})G^{a}X_{1}^{b_{1}-1}X_{2}^{b_{2}}\otimes g^{d}x_{1}^{e_{1}}x_{2}^{e_{2}-1}+
\\
&&\sum_{\substack{ a,b_{1},b_{2},d,e_{1},e_{2}=0  \\ %
a+b_{1}+b_{2}+d+e_{1}+e_{2}\equiv 0}}^{1}\left( -1\right) ^{\alpha \left(
1_{H};1,1,0,0\right) }B(1_{H}\otimes
x_{1}x_{2};G^{a}X_{1}X_{2},g^{d}x_{1}^{e_{1}}x_{2}^{e_{2}})G^{a}X_{1}^{b_{1}-1}X_{2}^{b_{2}-1}\otimes g^{d}x_{1}^{e_{1}}x_{2}^{e_{2}}
\end{eqnarray*}%
since%
\begin{eqnarray*}
&&\alpha \left( 1_{H};0,0,1,1\right) \equiv 1+e_{2}\equiv 0 \\
&&\alpha \left( 1_{H};0,1,1,0\right) \equiv e_{2}+a+b_{1}+b_{2}+1\equiv
e_{2}+a+b_{1}\equiv d \\
&&\alpha \left( 1_{H};1,0,0,1\right) \equiv b_{2}+\left(
a+b_{1}+b_{2}+1\right) +1\equiv a+b_{1}\equiv a+1 \\
&&\alpha \left( 1_{H};1,1,0,0\right) \equiv 1+b_{2}\equiv 0
\end{eqnarray*}%
By considering also the right side we get

\begin{gather*}
\sum_{\substack{ a,b_{1},b_{2},d,=0  \\ a+b_{1}+b_{2}+d\equiv 0}}%
^{1}B(1_{H}\otimes
x_{1}x_{2};G^{a}X_{1}^{b_{1}}X_{2}^{b_{2}},g^{d}x_{1}x_{2})G^{a}X_{1}^{b_{1}}X_{2}^{b_{2}}\otimes g^{d}+
\\
\sum_{\substack{ a,b_{1},d,e_{2}=0  \\ a+b_{1}+d+e_{2}\equiv 0}}^{1}\left(
-1\right) ^{d}B(1_{H}\otimes
x_{1}x_{2};G^{a}X_{1}^{b_{1}}X_{2},g^{d}x_{1}x_{2}^{e_{2}})G^{a}X_{1}^{b_{1}}\otimes g^{d}x_{2}^{e_{2}}+
\\
\sum_{\substack{ a,b_{2},d,e_{1}=0  \\ a+b_{2}+d+e_{1}\equiv 0}}^{1}\left(
-1\right) ^{a+1}B(1_{H}\otimes
x_{1}x_{2};G^{a}X_{1}X_{2}^{b_{2}},g^{d}x_{1}^{e_{1}}x_{2})G^{a}X_{2}^{b_{2}}\otimes g^{d}x_{1}^{e_{1}}+
\\
\sum_{\substack{ a,d,e_{1},e_{2}=0  \\ a+d+e_{1}+e_{2}\equiv 0}}%
^{1}B(1_{H}\otimes
x_{1}x_{2};G^{a}X_{1}X_{2},g^{d}x_{1}^{e_{1}}x_{2}^{e_{2}})G^{a}\otimes
g^{d}x_{1}^{e_{1}}x_{2}^{e_{2}} \\
-B(1_{A}\otimes
1_{H};G^{a}X_{1}^{b_{1}}X_{2}^{b_{2}},g^{d}x_{1}^{e_{1}}x_{2}^{e_{2}})G^{a}X_{1}^{b_{1}}X_{2}^{b_{2}}\otimes g^{d}x_{1}^{e_{1}}x_{2}^{e_{2}}=0
\end{gather*}

\subsubsection{$G^{a}\otimes g^{d}$}

\begin{equation*}
+\sum_{\substack{ a,d,=0  \\ a+d\equiv 0}}^{1}+\left[
\begin{array}{c}
-B(1_{A}\otimes 1_{H};G^{a},g^{d})+B(1_{H}\otimes
x_{1}x_{2};G^{a},g^{d}x_{1}x_{2})+ \\
+\left( -1\right) ^{d}B(1_{H}\otimes x_{1}x_{2};G^{a}X_{2},g^{d}x_{1})+ \\
\left( -1\right) ^{a+1}B(1_{H}\otimes
x_{1}x_{2};G^{a}X_{1},g^{d}x_{2})+B(1_{H}\otimes
x_{1}x_{2};G^{a}X_{1}X_{2},g^{d})%
\end{array}%
\right] G^{a}\otimes g^{d}
\end{equation*}

and we get%
\begin{gather*}
-B(1_{A}\otimes 1_{H};1_{A},1_{H})+B(1_{H}\otimes
x_{1}x_{2};1_{A},x_{1}x_{2})+B(1_{H}\otimes x_{1}x_{2};X_{2},x_{1}) \\
-B(1_{H}\otimes x_{1}x_{2};X_{1},x_{2})+B(1_{H}\otimes
x_{1}x_{2};X_{1}X_{2},1_{H})=0
\end{gather*}%
and%
\begin{gather*}
-B(1_{A}\otimes 1_{H};G,g)+B(1_{H}\otimes
x_{1}x_{2};G,gx_{1}x_{2})-B(1_{H}\otimes x_{1}x_{2};GX_{2},gx_{1}) \\
+B(1_{H}\otimes x_{1}x_{2};GX_{1},gx_{2})+B(1_{H}\otimes
x_{1}x_{2};GX_{1}X_{2},g)=0.
\end{gather*}%
By using the forms of $B(1_{H}\otimes x_{1}x_{2})$ we get no new information.

\subsubsection{$G^{a}\otimes g^{d}x_{2}$}

\begin{gather*}
\sum_{\substack{ a,d=0  \\ a+d\equiv 1}}^{1}\left[ \left( -1\right)
^{d}B(1_{H}\otimes x_{1}x_{2};G^{a}X_{2},g^{d}x_{1}x_{2})+B(1_{H}\otimes
x_{1}x_{2};G^{a}X_{1}X_{2},g^{d}x_{2})\right] \\
G^{a}\otimes g^{d}x_{2}=0
\end{gather*}%
and we get%
\begin{equation*}
-B(1_{H}\otimes x_{1}x_{2};X_{2},gx_{1}x_{2})+B(1_{H}\otimes
x_{1}x_{2};X_{1}X_{2},gx_{2})=0
\end{equation*}%
and%
\begin{equation*}
B(1_{H}\otimes x_{1}x_{2};GX_{2},x_{1}x_{2})+B(1_{H}\otimes
x_{1}x_{2};GX_{1}X_{2},x_{2})=0.
\end{equation*}%
By using the forms of $B(1_{H}\otimes x_{1}x_{2})$ we get no new information.

\subsubsection{$G^{a}\otimes g^{d}x_{1}$}

\begin{equation*}
\sum_{\substack{ a,d=0  \\ a+d\equiv 1}}^{1}\left[ \left( -1\right)
^{a+1}B(1_{H}\otimes x_{1}x_{2};G^{a}X_{1},g^{d}x_{1}x_{2})+B(1_{H}\otimes
x_{1}x_{2};G^{a}X_{1}X_{2},g^{d}x_{1})\right] G^{a}\otimes g^{d}x_{1}=0
\end{equation*}%
We get%
\begin{equation*}
-B(1_{H}\otimes x_{1}x_{2};X_{1},gx_{1}x_{2})+B(1_{H}\otimes
x_{1}x_{2};X_{1}X_{2},gx_{1})=0
\end{equation*}%
and%
\begin{equation*}
B(1_{H}\otimes x_{1}x_{2};GX_{1},x_{1}x_{2})+B(1_{H}\otimes
x_{1}x_{2};GX_{1}X_{2},x_{1})=0.
\end{equation*}%
By using the forms of $B(1_{H}\otimes x_{1}x_{2})$ we get no new information.

\subsubsection{$G^{a}X_{2}\otimes g^{d}$}

\begin{gather*}
\sum_{\substack{ a,d=0  \\ a+d\equiv 1}}^{1}B(1_{H}\otimes
x_{1}x_{2};G^{a}X_{2},g^{d}x_{1}x_{2})+\left( -1\right) ^{a+1}B(1_{H}\otimes
x_{1}x_{2};G^{a}X_{1}X_{2},g^{d}x_{2})] \\
G^{a}X_{2}\otimes g^{d}=0
\end{gather*}%
and we get%
\begin{equation*}
B(1_{H}\otimes x_{1}x_{2};X_{2},gx_{1}x_{2})-B(1_{H}\otimes
x_{1}x_{2};X_{1}X_{2},gx_{2})=0
\end{equation*}%
\begin{equation*}
B(1_{H}\otimes x_{1}x_{2};GX_{2},x_{1}x_{2})+B(1_{H}\otimes
x_{1}x_{2};GX_{1}X_{2},gx_{2})=0
\end{equation*}%
which we already obtained.

\subsubsection{$G^{a}X_{1}\otimes g^{d}$}

\begin{equation*}
\sum_{\substack{ a,d=0  \\ a+d\equiv 1}}^{1}\left[ B(1_{H}\otimes
x_{1}x_{2};G^{a}X_{1},g^{d}x_{1}x_{2})+\left( -1\right) ^{d}B(1_{H}\otimes
x_{1}x_{2};G^{a}X_{1}X_{2},g^{d}x_{1})\right] G^{a}X_{1}\otimes g^{d}=0
\end{equation*}%
and we get%
\begin{equation*}
B(1_{H}\otimes x_{1}x_{2};X_{1},gx_{1}x_{2})-B(1_{H}\otimes
x_{1}x_{2};X_{1}X_{2},gx_{1})=0
\end{equation*}%
\begin{equation*}
B(1_{H}\otimes x_{1}x_{2};GX_{1},x_{1}x_{2})+B(1_{H}\otimes
x_{1}x_{2};GX_{1}X_{2},x_{1})=0
\end{equation*}%
which we already obtained.

\subsubsection{$G^{a}\otimes g^{d}x_{1}x_{2}$}

\begin{equation*}
\sum_{\substack{ a,d=0  \\ a+d\equiv 0}}^{1}B(1_{H}\otimes
x_{1}x_{2};G^{a}X_{1}X_{2},g^{d}x_{1}x_{2})G^{a}\otimes g^{d}x_{1}x_{2}=0
\end{equation*}%
We get%
\begin{equation*}
B(1_{H}\otimes x_{1}x_{2};X_{1}X_{2},x_{1}x_{2})=0
\end{equation*}%
\begin{equation*}
B(1_{H}\otimes x_{1}x_{2};GX_{1}X_{2},gx_{1}x_{2})=0
\end{equation*}%
which are already known

\subsubsection{$G^{a}X_{2}\otimes g^{d}x_{2}$}

There are no terms.

\subsubsection{$G^{a}X_{1}\otimes g^{d}x_{2}$}

\begin{equation*}
\sum_{\substack{ a,d=0  \\ a+d\equiv 0}}^{1}\left( -1\right)
^{d}B(1_{H}\otimes
x_{1}x_{2};G^{a}X_{1}X_{2},g^{d}x_{1}x_{2})G^{a}X_{1}\otimes g^{d}x_{2}=0
\end{equation*}%
\begin{equation*}
B(1_{H}\otimes x_{1}x_{2};X_{1}X_{2},x_{1}x_{2})=0
\end{equation*}%
\begin{equation*}
B(1_{H}\otimes x_{1}x_{2};GX_{1}X_{2},gx_{1}x_{2})=0
\end{equation*}%
already got

\subsubsection{$G^{a}X_{2}\otimes g^{d}x_{1}$}

\begin{equation*}
\sum_{\substack{ a,d=0  \\ a+d\equiv 0}}^{1}\left( -1\right)
^{a+1}B(1_{H}\otimes
x_{1}x_{2};G^{a}X_{1}X_{2},g^{d}x_{1}x_{2})G^{a}X_{2}\otimes g^{d}x_{1}=0
\end{equation*}%
We obtain%
\begin{equation*}
B(1_{H}\otimes x_{1}x_{2};X_{1}X_{2},x_{1}x_{2})=0
\end{equation*}%
\begin{equation*}
B(1_{H}\otimes x_{1}x_{2};GX_{1}X_{2},gx_{1}x_{2})=0
\end{equation*}%
already got

\subsubsection{$G^{a}X_{1}\otimes g^{d}x_{1}$}

We get nothing

\subsubsection{$G^{a}X_{1}X_{2}\otimes g^{d}$}

\begin{equation*}
\sum_{\substack{ a,d=0  \\ a+d\equiv 0}}^{1}B(1_{H}\otimes
x_{1}x_{2};G^{a}X_{1}X_{2},g^{d}x_{1}x_{2})G^{a}X_{1}X_{2}\otimes g^{d}
\end{equation*}%
We obtain%
\begin{equation*}
B(1_{H}\otimes x_{1}x_{2};X_{1}X_{2},x_{1}x_{2})=0
\end{equation*}%
\begin{equation*}
B(1_{H}\otimes x_{1}x_{2};GX_{1}X_{2},gx_{1}x_{2})=0
\end{equation*}%
already got

\subsubsection{$G^{a}X_{2}\otimes g^{d}x_{1}x_{2}$}

No terms

\subsubsection{$G^{a}X_{1}\otimes g^{d}x_{1}x_{2}$}

No terms

\subsubsection{$G^{a}X_{1}X_{2}\otimes g^{d}x_{2}$}

No terms

\subsubsection{$G^{a}X_{1}X_{2}\otimes g^{d}x_{1}$}

No terms

\subsubsection{$G^{a}X_{1}X_{2}\otimes g^{d}x_{1}x_{2}$}

No terms

\subsection{Case $gx_{1}$}

\begin{eqnarray*}
l_{1}+u_{1} &=&1 \\
l_{2} &=&u_{2}=0 \\
a+b_{1}+b_{2}+d+e_{1}+e_{2} &\equiv &0
\end{eqnarray*}

\begin{eqnarray*}
&&\sum_{\substack{ a,b_{1},b_{2},d,e_{1},e_{2}=0  \\ %
a+b_{1}+b_{2}+d+e_{1}+e_{2}\equiv 0}}^{1}\sum_{l_{1}=0}^{b_{1}}%
\sum_{u_{1}=0}^{e_{1}}\left( -1\right) ^{\alpha \left(
1_{H};l_{1},0,u_{1},0\right) }B(1_{H}\otimes
x_{1}x_{2};G^{a}X_{1}^{b_{1}}X_{2}^{b_{2}},g^{d}x_{1}^{e_{1}}x_{2}^{e_{2}})
\\
&&G^{a}X_{1}^{b_{1}-l_{1}}X_{2}^{b_{2}}\otimes
g^{d}x_{1}^{e_{1}-u_{1}}x_{2}^{e_{2}} \\
&=&-B^{A}(1_{A}\otimes x_{2})\otimes B^{H}(1_{A}\otimes x_{2})
\end{eqnarray*}%
\begin{eqnarray*}
&&\sum_{\substack{ a,b_{1},b_{2},d,e_{2}=0  \\ a+b_{1}+b_{2}+d+e_{2}\equiv 1
}}^{1}\left( -1\right) ^{\alpha \left( 1_{H};0,0,1,0\right) }B(1_{H}\otimes
x_{1}x_{2};G^{a}X_{1}^{b_{1}}X_{2}^{b_{2}},g^{d}x_{1}x_{2}^{e_{2}})G^{a}X_{1}^{b_{1}}X_{2}^{b_{2}}\otimes g^{d}x_{2}^{e_{2}}+
\\
&&\sum_{\substack{ a,b_{2},d,e_{1},e_{2}=0  \\ a+b_{2}+d+e_{1}+e_{2}\equiv 1
}}^{1}\left( -1\right) ^{\alpha \left( 1_{H};1,0,0,0\right) }B(1_{H}\otimes
x_{1}x_{2};G^{a}X_{1}X_{2}^{b_{2}},g^{d}x_{1}^{e_{1}}x_{2}^{e_{2}})G^{a}X_{2}^{b_{2}}\otimes g^{d}x_{1}^{e_{1}}x_{2}^{e_{2}}
\\
&=&-B^{A}(1_{A}\otimes x_{2})\otimes B^{H}1_{A}\otimes x_{2})
\end{eqnarray*}%
Since%
\begin{eqnarray*}
&&\alpha \left( 1_{H};0,0,1,0\right) \equiv e_{2}+a+b_{1}+b_{2} \\
&&\alpha \left( 1_{H};1,0,0,0\right) \equiv b_{2}
\end{eqnarray*}%
we finally obtain%
\begin{eqnarray*}
&&\sum_{\substack{ a,b_{1},b_{2},d,e_{2}=0  \\ a+b_{1}+b_{2}+d+e_{2}\equiv 1
}}^{1}\left( -1\right) ^{e_{2}+a+b_{1}+b_{2}}B(1_{H}\otimes
x_{1}x_{2};G^{a}X_{1}^{b_{1}}X_{2}^{b_{2}},g^{d}x_{1}x_{2}^{e_{2}})G^{a}X_{1}^{b_{1}}X_{2}^{b_{2}}\otimes g^{d}x_{2}^{e_{2}}+
\\
&&\sum_{\substack{ a,b_{2},d,e_{1},e_{2}=0  \\ a+b_{2}+d+e_{1}+e_{2}\equiv 1
}}^{1}\left( -1\right) ^{b_{2}}B(1_{H}\otimes
x_{1}x_{2};G^{a}X_{1}X_{2}^{b_{2}},g^{d}x_{1}^{e_{1}}x_{2}^{e_{2}})G^{a}X_{2}^{b_{2}}\otimes g^{d}x_{1}^{e_{1}}x_{2}^{e_{2}}
\\
&=&-\sum_{a,b_{1},b_{2},d,e_{1},e_{2}=0}^{1}B(1_{A}\otimes
x_{2};G^{a}X_{1}^{b_{1}}X_{2}^{b_{2}},g^{d}x_{1}^{e_{1}}x_{2}^{e_{2}})G^{a}X_{1}^{b_{1}}X_{2}^{b_{2}}\otimes g^{d}x_{1}^{e_{1}}x_{2}^{e_{2}}
\end{eqnarray*}

\subsubsection{$G^{a}\otimes g^{d}$}

\begin{gather*}
\sum_{\substack{ a,d=0  \\ a+d\equiv 1}}^{1}\left[ B(1_{H}\otimes
x_{2};G^{a},g^{d})+\left( -1\right) ^{a}B(1_{H}\otimes
x_{1}x_{2};G^{a},g^{d}x_{1})+B(1_{H}\otimes x_{1}x_{2};G^{a}X_{1},g^{d})%
\right] \\
G^{a}\otimes g^{d}=0
\end{gather*}%
We get%
\begin{equation*}
B(1_{A}\otimes x_{2};1_{A},g)+B(1_{H}\otimes
x_{1}x_{2};1_{A},gx_{1})+B(1_{H}\otimes x_{1}x_{2};X_{1},g)=0
\end{equation*}%
and%
\begin{equation*}
B(1_{A}\otimes x_{2};G,1_{H})-B(1_{H}\otimes
x_{1}x_{2};G,x_{1})+B(1_{H}\otimes x_{1}x_{2};GX_{1},1_{H})=0.
\end{equation*}%
By using the forms of $B(1_{H}\otimes x_{1}x_{2})$ and of $B(1_{A}\otimes
x_{2}),$ we get no new information.

\subsubsection{$G^{a}\otimes g^{d}x_{2}${\protect\Huge \ }}

\begin{gather*}
\sum_{\substack{ a,d=0  \\ a+d\equiv 0}}^{1}\left[ \left( -1\right)
^{a+1}B(1_{H}\otimes x_{1}x_{2};G^{a},g^{d}x_{1}x_{2})+B(1_{H}\otimes
x_{1}x_{2};G^{a}X_{1},g^{d}x_{2})+B(1_{H}\otimes x_{2};G^{a},g^{d}x_{2})%
\right] \\
G^{a}\otimes g^{d}x_{2}=0
\end{gather*}%
We get%
\begin{equation*}
B(1_{H}\otimes x_{2};1_{A},x_{2})-B(1_{H}\otimes
x_{1}x_{2};1_{A},x_{1}x_{2})+B(1_{H}\otimes x_{1}x_{2};X_{1},x_{2})=0
\end{equation*}%
and%
\begin{equation*}
B(1_{H}\otimes x_{2};G,gx_{2})+B(1_{H}\otimes
x_{1}x_{2};G,gx_{1}x_{2})+B(1_{H}\otimes x_{1}x_{2};GX_{1},gx_{2})=0
\end{equation*}%
By using the form of $B(1_{H}\otimes x_{1}x_{2})$ and of $B(1_{A}\otimes
x_{2}),$ we get no new information.

\subsubsection{$G^{a}\otimes g^{d}x_{1}$}

\begin{equation*}
\sum_{\substack{ a,d=0  \\ a+d\equiv 0}}^{1}\left[ B(1_{H}\otimes
x_{2};G^{a},g^{d}x_{1})+B(1_{H}\otimes x_{1}x_{2};G^{a}X_{1},g^{d}x_{1})%
\right] G^{a}\otimes g^{d}x_{1}=0
\end{equation*}%
We get%
\begin{equation*}
B(1_{H}\otimes x_{2};1_{A},x_{1})+B(1_{H}\otimes x_{1}x_{2};X_{1},x_{1})=0
\end{equation*}%
and%
\begin{equation*}
B(1_{H}\otimes x_{2};G,gx_{1})+B(1_{H}\otimes x_{1}x_{2};GX_{1},gx_{1})=0.
\end{equation*}%
By using the form of $B(1_{H}\otimes x_{1}x_{2})$ and of $B(1_{A}\otimes
x_{2}),$ we get no new information.

\subsubsection{$G^{a}X_{2}\otimes g^{d}$}

\begin{equation*}
\sum_{\substack{ a,d=0  \\ a+d\equiv 0}}^{1}\left[ B(1_{H}\otimes
x_{2};G^{a}X_{2},g^{d})+\left( -1\right) ^{a+1}B(1_{H}\otimes
x_{1}x_{2};G^{a}X_{2},g^{d}x_{1})-B(1_{H}\otimes
x_{1}x_{2};G^{a}X_{1}X_{2},g^{d})\right] G^{a}X_{2}\otimes g^{d}=0
\end{equation*}%
We get%
\begin{equation*}
B(1_{H}\otimes x_{2};X_{2},1_{H})-B(1_{H}\otimes
x_{1}x_{2};X_{2},x_{1})-B(1_{H}\otimes x_{1}x_{2};X_{1}X_{2},1_{H})=0
\end{equation*}%
and%
\begin{equation*}
B(1_{H}\otimes x_{2};GX_{2},g)+B(1_{H}\otimes
x_{1}x_{2};GX_{2},gx_{1})-B(1_{H}\otimes x_{1}x_{2};GX_{1}X_{2},g)=0.
\end{equation*}%
By using the forms of $B(1_{H}\otimes x_{1}x_{2})$ and of $B(1_{A}\otimes
x_{2})$ $,$ we get no new information.

\subsubsection{$G^{a}X_{1}\otimes g^{d}$}

\begin{equation*}
\sum_{\substack{ a,d=0  \\ a+d\equiv 0}}^{1}\left[ B(1_{H}\otimes
x_{2};G^{a}X_{1},g^{d})+\left( -1\right) ^{a+1}B(1_{H}\otimes
x_{1}x_{2};G^{a}X_{1},g^{d}x_{1})\right] G^{a}X_{1}\otimes g^{d}=0
\end{equation*}%
We get%
\begin{equation*}
B(1_{H}\otimes x_{2};X_{1},1_{H})-B(1_{H}\otimes x_{1}x_{2};X_{1},x_{1})=0
\end{equation*}%
and%
\begin{equation*}
B(1_{H}\otimes x_{2};GX_{1},g)+B(1_{H}\otimes x_{1}x_{2};GX_{1},gx_{1})=0.
\end{equation*}%
By using the form of $B(1_{H}\otimes x_{1}x_{2})$ and of $B(1_{A}\otimes
x_{2})$ $,$ we get no new information.

\subsubsection{$G^{a}\otimes g^{d}x_{1}x_{2}$}

\begin{equation*}
\sum_{\substack{ a,d=0  \\ a+d\equiv 1}}^{1}\left[ B(1_{H}\otimes
x_{2};G^{a},g^{d}x_{1}x_{2})+B(1_{H}\otimes
x_{1}x_{2};G^{a}X_{1},g^{d}x_{1}x_{2})\right] G^{a}\otimes g^{d}x_{1}x_{2}=0
\end{equation*}%
We get%
\begin{equation*}
B(1_{H}\otimes x_{2};1_{A},gx_{1}x_{2})+B(1_{H}\otimes
x_{1}x_{2};X_{1},gx_{1}x_{2})=0
\end{equation*}%
and%
\begin{equation*}
B(1_{H}\otimes x_{2};G,x_{1}x_{2})+B(1_{H}\otimes
x_{1}x_{2};GX_{1},x_{1}x_{2})=0
\end{equation*}%
By using the forms of $B(1_{H}\otimes x_{1}x_{2})$ and of $B(1_{A}\otimes
x_{2})$ $,$ we get no new information.

\subsubsection{$G^{a}X_{2}\otimes g^{d}x_{2}$}

\begin{equation*}
\sum_{\substack{ a,d=0  \\ a+d\equiv 1}}^{1}\left[
\begin{array}{c}
B(1_{H}\otimes x_{2};G^{a}X_{2},g^{d}x_{2})+\left( -1\right)
^{a}B(1_{H}\otimes x_{1}x_{2};G^{a}X_{2},g^{d}x_{1}x_{2}) \\
-B(1_{H}\otimes x_{1}x_{2};G^{a}X_{1}X_{2},g^{d}x_{2})%
\end{array}%
\right] G^{a}X_{2}\otimes g^{d}x_{2}=0.
\end{equation*}%
We get%
\begin{equation*}
B(1_{H}\otimes x_{2};X_{2},gx_{2})+B(1_{H}\otimes
x_{1}x_{2};X_{2},gx_{1}x_{2})-B(1_{H}\otimes x_{1}x_{2};X_{1}X_{2},gx_{2})=0
\end{equation*}%
and%
\begin{equation*}
B(1_{H}\otimes x_{2};GX_{2},x_{2})-B(1_{H}\otimes
x_{1}x_{2};GX_{2},x_{1}x_{2})-B(1_{H}\otimes x_{1}x_{2};GX_{1}X_{2},x_{2})=0.
\end{equation*}%
By using the form of $B(1_{H}\otimes x_{1}x_{2})$ and of $B(1_{A}\otimes
x_{2})$ $,$ we get no new information.

\subsubsection{$G^{a}X_{1}\otimes g^{d}x_{2}$}

\begin{equation*}
\sum_{\substack{ a,d=0  \\ a+d\equiv 1}}^{1}\left[ B(1_{H}\otimes
x_{2};G^{a}X_{1},g^{d}x_{2})+\left( -1\right) ^{a}B(1_{H}\otimes
x_{1}x_{2};G^{a}X_{1},g^{d}x_{1}x_{2})\right] G^{a}X_{1}\otimes g^{d}x_{2}=0
\end{equation*}%
We get%
\begin{equation*}
B(1_{H}\otimes x_{2};X_{1},gx_{2})+B(1_{H}\otimes
x_{1}x_{2};X_{1},gx_{1}x_{2})=0
\end{equation*}%
and%
\begin{equation*}
B(1_{H}\otimes x_{2};GX_{1},x_{2})-B(1_{H}\otimes
x_{1}x_{2};GX_{1},x_{1}x_{2})=0.
\end{equation*}%
By using the form of $B(1_{H}\otimes x_{1}x_{2})$ and of $B(1_{A}\otimes
x_{2})$ $,$ we get no new information.

\subsubsection{$G^{a}X_{2}\otimes g^{d}x_{1}$}

\begin{equation*}
\sum_{\substack{ a,d=0  \\ a+d\equiv 1}}^{1}\left[ B(1_{H}\otimes
x_{2};G^{a}X_{2},g^{d}x_{1})-B(1_{H}\otimes
x_{1}x_{2};G^{a}X_{1}X_{2},g^{d}x_{1})\right] G^{a}X_{2}\otimes g^{d}x_{1}=0
\end{equation*}%
We get%
\begin{equation*}
B(1_{H}\otimes x_{2};X_{2},gx_{1})-B(1_{H}\otimes
x_{1}x_{2};X_{1}X_{2},gx_{1})=0
\end{equation*}%
and
\begin{equation*}
B(1_{H}\otimes x_{2};GX_{2},x_{1})-B(1_{H}\otimes
x_{1}x_{2};GX_{1}X_{2},x_{1})=0.
\end{equation*}%
By using the forms of $B(1_{H}\otimes x_{1}x_{2})$ and of $B(1_{A}\otimes
x_{2})$ $,$ we get no new information.

\subsubsection{$G^{a}X_{1}\otimes g^{d}x_{1}$}

\begin{equation*}
\sum_{\substack{ a,d=0  \\ a+d\equiv 1}}^{1}B(1_{H}\otimes
x_{2};G^{a}X_{1},g^{d}x_{1})G^{a}X_{1}\otimes g^{d}x_{1}=0
\end{equation*}%
we get%
\begin{equation*}
B(1_{H}\otimes x_{2};X_{1},gx_{1})=0
\end{equation*}%
\begin{equation*}
B(1_{H}\otimes x_{2};GX_{1},x_{1})=0
\end{equation*}%
which are already known.

\subsubsection{$G^{a}X_{1}X_{2}\otimes g^{d}$}

\begin{equation*}
\sum_{\substack{ a,d=0  \\ a+d\equiv 1}}^{1}\left[ B(1_{H}\otimes
x_{2};G^{a}X_{1}X_{2},g^{d})+\left( -1\right) ^{a}B(1_{H}\otimes
x_{1}x_{2};G^{a}X_{1}X_{2},g^{d}x_{1})\right] G^{a}X_{1}X_{2}\otimes g^{d}=0
\end{equation*}%
We get%
\begin{equation*}
B(1_{H}\otimes x_{2};X_{1}X_{2},g)+B(1_{H}\otimes
x_{1}x_{2};X_{1}X_{2},gx_{1})=0
\end{equation*}%
and%
\begin{equation*}
B(1_{H}\otimes x_{2};GX_{1}X_{2},1_{H})-B(1_{H}\otimes
x_{1}x_{2};GX_{1}X_{2},x_{1})=0
\end{equation*}%
By using the forms of $B(1_{H}\otimes x_{1}x_{2})$ and of $B(1_{A}\otimes
x_{2})$ $,$ we get no new information.

\subsubsection{$G^{a}X_{2}\otimes g^{d}x_{1}x_{2}$}

\begin{equation*}
\sum_{\substack{ a,d=0  \\ a+d\equiv 0}}^{1}\left[ B(1_{H}\otimes
x_{2};G^{a}X_{2},g^{d}x_{1}x_{2})-B(1_{H}\otimes
x_{1}x_{2};G^{a}X_{1}X_{2},g^{d}x_{1}x_{2})\right] G^{a}X_{2}\otimes
g^{d}x_{1}x_{2}=0
\end{equation*}%
We get%
\begin{equation*}
B(1_{H}\otimes x_{2};X_{2},x_{1}x_{2})-B(1_{H}\otimes
x_{1}x_{2};X_{1}X_{2},x_{1}x_{2})=0
\end{equation*}%
and%
\begin{equation*}
B(1_{H}\otimes x_{2};GX_{2},gx_{1}x_{2})-B(1_{H}\otimes
x_{1}x_{2};GX_{1}X_{2},gx_{1}x_{2})=0
\end{equation*}%
By using the forms of $B(1_{H}\otimes x_{1}x_{2})$ and of $B(1_{A}\otimes
x_{2})$ $,$ we get no new information.

\subsubsection{$G^{a}X_{1}\otimes g^{d}x_{1}x_{2}$}

\begin{equation*}
\sum_{\substack{ a,d=0  \\ a+d\equiv 0}}^{1}B(1_{H}\otimes
x_{2};G^{a}X_{1},g^{d}x_{1}x_{2})G^{a}X_{1}\otimes g^{d}x_{1}x_{2}=0
\end{equation*}%
We get%
\begin{equation*}
B(1_{H}\otimes x_{2};X_{1},x_{1}x_{2})=0
\end{equation*}%
and%
\begin{equation*}
B(1_{H}\otimes x_{2};GX_{1},gx_{1}x_{2})=0
\end{equation*}%
These are already known.

\subsubsection{$G^{a}X_{1}X_{2}\otimes g^{d}x_{2}$}

\begin{equation*}
\sum_{\substack{ a,d=0  \\ a+d\equiv 0}}^{1}\left[ B(1_{H}\otimes
x_{2};G^{a}X_{1}X_{2},g^{d}x_{2})+\left( -1\right) ^{a}B(1_{H}\otimes
x_{1}x_{2};G^{a}X_{1}X_{2},g^{d}x_{1}x_{2})\right] G^{a}X_{1}X_{2}\otimes
g^{d}x_{2}=0
\end{equation*}%
We get%
\begin{equation*}
B(1_{H}\otimes x_{2};X_{1}X_{2},x_{2})-B(1_{H}\otimes
x_{1}x_{2};X_{1}X_{2},x_{1}x_{2})=0
\end{equation*}%
\begin{equation*}
B(1_{H}\otimes x_{2};GX_{1}X_{2},gx_{2})+B(1_{H}\otimes
x_{1}x_{2};GX_{1}X_{2},gx_{1}x_{2})=0
\end{equation*}%
Both of them are already known.

\subsubsection{$G^{a}X_{1}X_{2}\otimes g^{d}x_{1}$}

\begin{equation*}
\sum_{\substack{ a,d=0  \\ a+d\equiv 0}}^{1}B(1_{H}\otimes
x_{2};G^{a}X_{1}X_{2},g^{d}x_{1})G^{a}X_{1}X_{2}\otimes g^{d}x_{1}=0
\end{equation*}%
We get%
\begin{equation*}
B(1_{H}\otimes x_{2};X_{1}X_{2},x_{1})=0
\end{equation*}%
\begin{equation*}
B(1_{H}\otimes x_{2};GX_{1}X_{2},gx_{1})=0
\end{equation*}%
Both of them were already known.

\subsubsection{$G^{a}X_{1}X_{2}\otimes g^{d}x_{1}x_{2}$}

\begin{equation*}
\sum_{\substack{ a,d=0  \\ a+d\equiv 0}}^{1}B(1_{H}\otimes
x_{2};G^{a}X_{1}X_{2},g^{d}x_{1}x_{2})G^{a}X_{1}X_{2}\otimes
g^{d}x_{1}x_{2}=0
\end{equation*}%
We get%
\begin{equation*}
B(1_{H}\otimes x_{2};X_{1}X_{2},x_{1}x_{2})=0
\end{equation*}%
and%
\begin{equation*}
B(1_{H}\otimes x_{2};GX_{1}X_{2},gx_{1}x_{2})=0
\end{equation*}%
Both of them were already known.

\subsection{Case $gx_{2}$}

\begin{eqnarray*}
l_{1} &=&u_{1}=0 \\
l_{2}+u_{2} &=&1 \\
a+b_{1}+b_{2}+d+e_{1}+e_{2} &\equiv &0
\end{eqnarray*}%
\begin{eqnarray*}
&&\sum_{\substack{ a,b_{1},b_{2},d,e_{1},e_{2}=0  \\ %
a+b_{1}+b_{2}+d+e_{1}+e_{2}\equiv 0}}^{1}\sum_{l_{2}=0}^{b_{2}}\sum
_{\substack{ u_{2}=0  \\ l_{2}+u_{2}=1}}^{e_{2}}\left( -1\right) ^{\alpha
\left( 1_{H};0,l_{2},0,u_{2}\right) }B(1_{H}\otimes
x_{1}x_{2};G^{a}X_{1}^{b_{1}}X_{2}^{b_{2}},g^{d}x_{1}^{e_{1}}x_{2}^{e_{2}})
\\
&&G^{a}X_{1}^{b_{1}}X_{2}^{b_{2}-l_{2}}\otimes
g^{d}x_{1}^{e_{1}}x_{2}^{e_{2}-u_{2}}\otimes \\
&=&B^{A}(1_{A}\otimes x_{1})\otimes B^{H}(1_{A}\otimes x_{1})\otimes gx_{2}
\end{eqnarray*}%
\begin{eqnarray*}
&&\sum_{\substack{ a,b_{1},b_{2},d,e_{1}=0  \\ a+b_{1}+b_{2}+d+e_{1}\equiv 1
}}^{1}\left( -1\right) ^{\alpha \left( 1_{H};0,0,0,1\right) }B(1_{H}\otimes
x_{1}x_{2};G^{a}X_{1}^{b_{1}}X_{2}^{b_{2}},g^{d}x_{1}^{e_{1}}x_{2})G^{a}X_{1}^{b_{1}}X_{2}^{b_{2}}\otimes g^{d}x_{1}^{e_{1}}+
\\
&&\sum_{\substack{ a,b_{1},d,e_{1},e_{2}=0  \\ a+b_{1}+d+e_{1}+e_{2}\equiv 1
}}^{1}\left( -1\right) ^{\alpha \left( 1_{H};0,1,0,0\right) }B(1_{H}\otimes
x_{1}x_{2};G^{a}X_{1}^{b_{1}}X_{2},g^{d}x_{1}^{e_{1}}x_{2}^{e_{2}})G^{a}X_{1}^{b_{1}}\otimes g^{d}x_{1}^{e_{1}}x_{2}^{e_{2}}+
\\
&=&B^{A}(1_{A}\otimes x_{1})\otimes B^{H}(1_{A}\otimes x_{1})
\end{eqnarray*}%
\begin{eqnarray*}
\alpha \left( 1_{H};0,0,0,1\right) &\equiv &a+b_{1}+b_{2} \\
\alpha \left( 1_{H};0,1,0,0\right) &\equiv &0
\end{eqnarray*}%
we get%
\begin{gather*}
-\sum_{a,b_{1},b_{2},d,e_{1},e_{2}=0}^{1}B(1_{A}\otimes
x_{1},G^{a}X_{1}^{b_{1}}X_{2}^{b_{2}},g^{d}x_{1}^{e_{1}}x_{2}^{e_{2}})G^{a}X_{1}^{b_{1}}X_{2}^{b_{2}}\otimes g^{d}x_{1}^{e_{1}}x_{2}^{e_{2}}
\\
\sum_{\substack{ a,b_{1},b_{2},d,e_{1}=0  \\ a+b_{1}+b_{2}+d+e_{1}\equiv 1}}%
^{1}\left( -1\right) ^{a+b_{1}+b_{2}}B(1_{H}\otimes
x_{1}x_{2};G^{a}X_{1}^{b_{1}}X_{2}^{b_{2}},g^{d}x_{1}^{e_{1}}x_{2})G^{a}X_{1}^{b_{1}}X_{2}^{b_{2}}\otimes g^{d}x_{1}^{e_{1}}+
\\
+\sum_{\substack{ a,b_{1},d,e_{1},e_{2}=0  \\ a+b_{1}+d+e_{1}+e_{2}\equiv 1}}%
^{1}B(1_{H}\otimes
x_{1}x_{2};G^{a}X_{1}^{b_{1}}X_{2},g^{d}x_{1}^{e_{1}}x_{2}^{e_{2}})G^{a}X_{1}^{b_{1}}\otimes g^{d}x_{1}^{e_{1}}x_{2}^{e_{2}}=0
\end{gather*}

\subsubsection{$G^{a}\otimes g^{d}$}

\begin{equation*}
\sum_{\substack{ a,d=0  \\ a+d\equiv 1}}^{1}\left[ -B(1_{A}\otimes
x_{1},G^{a},g^{d})+\left( -1\right) ^{a}B(1_{H}\otimes
x_{1}x_{2};G^{a},g^{d}x_{2})+B(1_{H}\otimes x_{1}x_{2};G^{a}X_{2},g^{d})%
\right] G^{a}\otimes g^{d}=0.
\end{equation*}%
We get%
\begin{equation*}
-B(1_{A}\otimes x_{1},1_{A},g)+B(1_{H}\otimes
x_{1}x_{2};1_{A},gx_{2})+B(1_{H}\otimes x_{1}x_{2};X_{2},g)=0
\end{equation*}%
and%
\begin{equation*}
-B(1_{A}\otimes x_{1},G,1_{H})-B(1_{H}\otimes
x_{1}x_{2};G,x_{2})+B(1_{H}\otimes x_{1}x_{2};GX_{2},1_{H})=0.
\end{equation*}%
By using the forms of $B(1_{H}\otimes x_{1}x_{2})$ and of $B(1_{A}\otimes
x_{1})$ $,$ we get no new information.

\subsubsection{$G^{a}\otimes g^{d}x_{2}$}

\begin{equation*}
\sum_{\substack{ a,d=0  \\ a+d\equiv 0}}^{1}\left[ -B(1_{A}\otimes
x_{1},G^{a},g^{d}x_{2})+B(1_{H}\otimes x_{1}x_{2};G^{a}X_{2},g^{d}x_{2})%
\right] G^{a}\otimes g^{d}x_{2}=0
\end{equation*}%
We get%
\begin{equation*}
-B(1_{A}\otimes x_{1},1_{A},x_{2})+B(1_{H}\otimes x_{1}x_{2};X_{2},x_{2})=0
\end{equation*}%
and%
\begin{equation*}
-B(1_{A}\otimes x_{1},G,gx_{2})+B(1_{H}\otimes x_{1}x_{2};GX_{2},gx_{2})=0.
\end{equation*}%
By using the forms of $B(1_{H}\otimes x_{1}x_{2})$ and of $B(1_{A}\otimes
x_{1})$ $,$ we get no new information.

\subsubsection{$G^{a}\otimes g^{d}x_{1}$}

\begin{gather*}
\sum_{\substack{ a,d=0  \\ a+d\equiv 0}}^{1}\left[ -B(1_{A}\otimes
x_{1},G^{a},g^{d}x_{1})+\left( -1\right) ^{a}B(1_{H}\otimes
x_{1}x_{2};G^{a},g^{d}x_{1}x_{2})+B(1_{H}\otimes
x_{1}x_{2};G^{a}X_{2},g^{d}x_{1})\right] \\
G^{a}\otimes g^{d}x_{1}=0
\end{gather*}%
We get%
\begin{equation*}
-B(1_{A}\otimes x_{1},1_{A},x_{1})+B(1_{H}\otimes
x_{1}x_{2};1_{A},x_{1}x_{2})+B(1_{H}\otimes x_{1}x_{2};X_{2},x_{1})=0
\end{equation*}%
and%
\begin{equation*}
-B(1_{A}\otimes x_{1},G,gx_{1})-B(1_{H}\otimes
x_{1}x_{2};G,gx_{1}x_{2})+B(1_{H}\otimes x_{1}x_{2};GX_{2},gx_{1})=0.
\end{equation*}%
By using the forms of $B(1_{H}\otimes x_{1}x_{2})$ and of $B(1_{A}\otimes
x_{1})$ $,$ we get no new information.

\subsubsection{$G^{a}X_{2}\otimes g^{d}$}

\begin{equation*}
\sum_{\substack{ a,d=0  \\ a+d\equiv 0}}^{1}\left[ -B(1_{A}\otimes
x_{1},G^{a}X_{2},g^{d})+\left( -1\right) ^{a+1}B(1_{H}\otimes
x_{1}x_{2};G^{a}X_{2},g^{d}x_{2})\right] G^{a}X_{2}\otimes g^{d}=0
\end{equation*}%
We get%
\begin{equation*}
-B(1_{A}\otimes x_{1},X_{2},1_{H})-B(1_{H}\otimes x_{1}x_{2};X_{2},x_{2})=0
\end{equation*}%
and%
\begin{equation*}
-B(1_{A}\otimes x_{1},GX_{2},g)+B(1_{H}\otimes x_{1}x_{2};GX_{2},gx_{2})=0.
\end{equation*}%
By using the forms of $B(1_{H}\otimes x_{1}x_{2})$ and of $B(1_{A}\otimes
x_{1})$ $,$ we get no new information.

\subsubsection{$G^{a}X_{1}\otimes g^{d}$}

\begin{gather*}
\sum_{\substack{ a,d=0  \\ a+d\equiv 0}}^{1}\left[ -B(1_{A}\otimes
x_{1},G^{a}X_{1},g^{d})+\left( -1\right) ^{a+1}B(1_{H}\otimes
x_{1}x_{2};G^{a}X_{1},g^{d}x_{2})+B(1_{H}\otimes
x_{1}x_{2};G^{a}X_{1}X_{2},g^{d})\right] \\
G^{a}X_{1}\otimes g^{d}=0
\end{gather*}%
We get%
\begin{equation*}
-B(1_{A}\otimes x_{1},X_{1},1_{H})-B(1_{H}\otimes
x_{1}x_{2};X_{1},x_{2})+B(1_{H}\otimes x_{1}x_{2};X_{1}X_{2},1_{H})=0
\end{equation*}%
and%
\begin{equation*}
-B(1_{A}\otimes x_{1},GX_{1},g)+B(1_{H}\otimes
x_{1}x_{2};GX_{1},gx_{2})+B(1_{H}\otimes x_{1}x_{2};GX_{1}X_{2},g)=0.
\end{equation*}%
By using the forms of $B(1_{H}\otimes x_{1}x_{2})$ and of $B(1_{A}\otimes
x_{1})$ $,$ we get no new information.

\subsubsection{$G^{a}\otimes g^{d}x_{1}x_{2}$}

\begin{equation*}
\sum_{\substack{ a,d=0  \\ a+d\equiv 1}}^{1}\left[ -B(1_{A}\otimes
x_{1},G^{a},g^{d}x_{1}x_{2})+B(1_{H}\otimes
x_{1}x_{2};G^{a}X_{2},g^{d}x_{1}x_{2})\right] G^{a}\otimes g^{d}x_{1}x_{2}=0
\end{equation*}%
We get%
\begin{equation*}
-B(1_{A}\otimes x_{1},1_{A},gx_{1}x_{2})+B(1_{H}\otimes
x_{1}x_{2};X_{2},gx_{1}x_{2})=0
\end{equation*}%
and%
\begin{equation*}
-B(1_{A}\otimes x_{1},G,x_{1}x_{2})+B(1_{H}\otimes
x_{1}x_{2};GX_{2},x_{1}x_{2})=0.
\end{equation*}%
By using the forms of $B(1_{H}\otimes x_{1}x_{2})$ and of $B(1_{A}\otimes
x_{1})$ $,$ we get no new information.

\subsubsection{$G^{a}X_{2}\otimes g^{d}x_{2}$}

\begin{equation*}
\sum_{\substack{ a,d=0  \\ a+d\equiv 1}}^{1}\left[ -B(1_{A}\otimes
x_{1},G^{a}X_{2},g^{d}x_{2})+B(1_{H}\otimes x_{1}x_{2};G^{a}X_{2},g^{d}x_{2})%
\right] G^{a}X_{2}\otimes g^{d}x_{2}=0
\end{equation*}%
We get%
\begin{equation*}
-B(1_{A}\otimes x_{1},X_{2},gx_{2})=0
\end{equation*}%
\begin{equation*}
-B(1_{A}\otimes x_{1},GX_{2},x_{2})=0
\end{equation*}%
already known

\subsubsection{$G^{a}X_{1}\otimes g^{d}x_{2}$}

\begin{equation*}
\sum_{\substack{ a,d=0  \\ a+d\equiv 1}}^{1}\left[ -B(1_{A}\otimes
x_{1},G^{a}X_{1},g^{d}x_{2})+B(1_{H}\otimes
x_{1}x_{2};G^{a}X_{1}X_{2},g^{d}x_{2})\right] G^{a}X_{1}\otimes g^{d}x_{2}=0
\end{equation*}%
We get%
\begin{equation*}
-B(1_{A}\otimes x_{1},X_{1},gx_{2})+B(1_{H}\otimes
x_{1}x_{2};X_{1}X_{2},gx_{2})=0
\end{equation*}%
and%
\begin{equation*}
-B(1_{A}\otimes x_{1},GX_{1},x_{2})+B(1_{H}\otimes
x_{1}x_{2};GX_{1}X_{2},x_{2})=0.
\end{equation*}%
By using the forms of $B(1_{H}\otimes x_{1}x_{2})$ and of $B(1_{A}\otimes
x_{1})$ $,$ we get no new information.

\subsubsection{$G^{a}X_{2}\otimes g^{d}x_{1}$}

\begin{equation*}
\sum_{\substack{ a,d=0  \\ a+d\equiv 1}}^{1}\left[ -B(1_{A}\otimes
x_{1},G^{a}X_{2},g^{d}x_{1})+\left( -1\right) ^{a+1}B(1_{H}\otimes
x_{1}x_{2};G^{a}X_{2},g^{d}x_{1}x_{2})\right] G^{a}X_{2}\otimes g^{d}x_{1}=0
\end{equation*}%
We get%
\begin{equation*}
-B(1_{A}\otimes x_{1},X_{2},gx_{1})-B(1_{H}\otimes
x_{1}x_{2};X_{2},gx_{1}x_{2})=0
\end{equation*}%
and%
\begin{equation*}
-B(1_{A}\otimes x_{1},GX_{2},x_{1})+B(1_{H}\otimes
x_{1}x_{2};GX_{2},x_{1}x_{2})=0.
\end{equation*}%
By using the forms of $B(1_{H}\otimes x_{1}x_{2})$ and of $B(1_{A}\otimes
x_{1})$ $,$ we get no new information.

\subsubsection{$G^{a}X_{1}\otimes g^{d}x_{1}$}

\begin{equation*}
\sum_{\substack{ a,d=0  \\ a+d\equiv 1}}^{1}\left[
\begin{array}{c}
-B(1_{A}\otimes x_{1},G^{a}X_{1},g^{d}x_{1})+\left( -1\right)
^{a+1}B(1_{H}\otimes x_{1}x_{2};G^{a}X_{1},g^{d}x_{1}x_{2}) \\
+B(1_{H}\otimes x_{1}x_{2};G^{a}X_{1}X_{2},g^{d}x_{1})%
\end{array}%
\right] G^{a}X_{1}\otimes g^{d}x_{1}=0
\end{equation*}%
We get%
\begin{equation*}
-B(1_{A}\otimes x_{1},X_{1},gx_{1})-B(1_{H}\otimes
x_{1}x_{2};X_{1},gx_{1}x_{2})+B(1_{H}\otimes x_{1}x_{2};X_{1}X_{2},gx_{1})=0
\end{equation*}%
and%
\begin{equation*}
-B(1_{A}\otimes x_{1},GX_{1},x_{1})+B(1_{H}\otimes
x_{1}x_{2};GX_{1},x_{1}x_{2})+B(1_{H}\otimes x_{1}x_{2};GX_{1}X_{2},x_{1})=0.
\end{equation*}%
By using the forms of $B(1_{H}\otimes x_{1}x_{2})$ and of $B(1_{A}\otimes
x_{1})$ $,$ we get no new information.

\subsubsection{$G^{a}X_{1}X_{2}\otimes g^{d}$}

\begin{equation*}
\sum_{\substack{ a,d=0  \\ a+d\equiv 1}}^{1}\left[
\begin{array}{c}
-B(1_{A}\otimes x_{1},G^{a}X_{1}X_{2},g^{d})+\left( -1\right)
^{a}B(1_{H}\otimes x_{1}x_{2};G^{a}X_{1}X_{2},g^{d}x_{2}) \\
+B(1_{H}\otimes x_{1}x_{2};G^{a}X_{1}X_{2},g^{d})%
\end{array}%
\right] G^{a}X_{1}X_{2}\otimes g^{d}=0
\end{equation*}%
We get%
\begin{equation*}
-B(1_{A}\otimes x_{1},X_{1}X_{2},g)+B(1_{H}\otimes
x_{1}x_{2};X_{1}X_{2},gx_{2})=0
\end{equation*}%
and%
\begin{equation*}
-B(1_{A}\otimes x_{1},GX_{1}X_{2},1_{H})-B(1_{H}\otimes
x_{1}x_{2};GX_{1}X_{2},x_{2})=0.
\end{equation*}%
By using the forms of $B(1_{H}\otimes x_{1}x_{2})$ and of $B(1_{A}\otimes
x_{1})$ $,$ we get no new information.

\subsubsection{$G^{a}X_{2}\otimes g^{d}x_{1}x_{2}$}

\begin{equation*}
\sum_{\substack{ a,d=0  \\ a+d\equiv 0}}^{1}-B(1_{A}\otimes
x_{1},G^{a}X_{2},g^{d}x_{1}x_{2})G^{a}X_{2}\otimes g^{d}x_{1}x_{2}=0
\end{equation*}%
We get%
\begin{equation*}
-B(1_{A}\otimes x_{1},X_{2},x_{1}x_{2})=0
\end{equation*}%
and%
\begin{equation*}
-B(1_{A}\otimes x_{1},GX_{2},gx_{1}x_{2})=0
\end{equation*}%
which are already known.

\subsubsection{$G^{a}X_{1}\otimes g^{d}x_{1}x_{2}$}

\begin{equation*}
\sum_{\substack{ a,d=0  \\ a+d\equiv 0}}^{1}\left[ -B(1_{A}\otimes
x_{1},G^{a}X_{1},g^{d}x_{1}x_{2})+B(1_{H}\otimes
x_{1}x_{2};G^{a}X_{1}X_{2},g^{d}x_{1}x_{2})\right] G^{a}X_{1}\otimes
g^{d}x_{1}x_{2}=0
\end{equation*}%
We get%
\begin{equation*}
-B(1_{A}\otimes x_{1},X_{1},x_{1}x_{2})+B(1_{H}\otimes
x_{1}x_{2};X_{1}X_{2},x_{1}x_{2})=0
\end{equation*}%
and%
\begin{equation*}
-B(1_{A}\otimes x_{1},GX_{1},gx_{1}x_{2})+B(1_{H}\otimes
x_{1}x_{2};GX_{1}X_{2},gx_{1}x_{2})=0.
\end{equation*}%
By using the forms of $B(1_{H}\otimes x_{1}x_{2})$ and of $B(1_{A}\otimes
x_{1})$ $,$ we get no new information.

\subsubsection{$G^{a}X_{1}X_{2}\otimes g^{d}x_{2}$}

\begin{equation*}
\sum_{\substack{ a,d=0  \\ a+d\equiv 0}}^{1}-B(1_{A}\otimes
x_{1},G^{a}X_{1}X_{2},g^{d}x_{2})G^{a}X_{1}X_{2}\otimes g^{d}x_{2}=0
\end{equation*}%
We get%
\begin{equation*}
-B(1_{A}\otimes x_{1},X_{1}X_{2},x_{2})=0
\end{equation*}%
and%
\begin{equation*}
-B(1_{A}\otimes x_{1},GX_{1}X_{2},gx_{2})=0
\end{equation*}%
which are already known.

\subsubsection{$G^{a}X_{1}X_{2}\otimes g^{d}x_{1}$}

\begin{equation*}
\sum_{\substack{ a,d=0  \\ a+d\equiv 0}}^{1}\left[ -B(1_{A}\otimes
x_{1},G^{a}X_{1}X_{2},g^{d}x_{1})+\left( -1\right) ^{a}B(1_{H}\otimes
x_{1}x_{2};G^{a}X_{1}X_{2},g^{d}x_{1}x_{2})\right] G^{a}X_{1}X_{2}\otimes
g^{d}x_{1}=0
\end{equation*}%
We get%
\begin{equation*}
-B(1_{A}\otimes x_{1},X_{1}X_{2},x_{1})+B(1_{H}\otimes
x_{1}x_{2};X_{1}X_{2},x_{1}x_{2})=0
\end{equation*}%
and%
\begin{equation*}
-B(1_{A}\otimes x_{1},GX_{1}X_{2},gx_{1})-B(1_{H}\otimes
x_{1}x_{2};GX_{1}X_{2},gx_{1}x_{2})=0
\end{equation*}%
By using the forms of $B(1_{H}\otimes x_{1}x_{2})$ and of $B(1_{A}\otimes
x_{1})$ $,$ we get no new information.

\subsubsection{$G^{a}X_{1}X_{2}\otimes g^{d}x_{1}x_{2}$}

\begin{equation*}
\sum_{\substack{ a,d=0  \\ a+d\equiv 1}}^{1}-B(1_{A}\otimes
x_{1},G^{a}X_{1}X_{2},g^{d}x_{1}x_{2})G^{a}X_{1}X_{2}\otimes
g^{d}x_{1}x_{2}=0
\end{equation*}%
We get%
\begin{equation*}
-B(1_{A}\otimes x_{1},X_{1}X_{2},gx_{1}x_{2})=0
\end{equation*}%
and%
\begin{equation*}
-B(1_{A}\otimes x_{1},GX_{1}X_{2},x_{1}x_{2})=0
\end{equation*}%
which are already known.

\subsection{Case $gx_{1}x_{2}$}

Nothing from the right side. From the left side we get%
\begin{eqnarray*}
l_{1}+u_{1} &=&1 \\
l_{2}+u_{2} &=&1 \\
a+b_{1}+b_{2}+d+e_{1}+e_{2} &\equiv &1
\end{eqnarray*}%
Since%
\begin{equation*}
B(1_{H}\otimes
x_{1}x_{2};G^{a}X_{1}^{b_{1}}X_{2}^{b_{2}},g^{d}x_{1}^{e_{1}}x_{2}^{e_{2}})=0%
\text{ whenever }a+b_{1}+b_{2}+d+e_{1}+e_{2}\equiv 1
\end{equation*}%
nothing to check.

\section{$B\left( x_{1}\otimes x_{1}\right) $}

From $\left( \ref{simplx}\right) $%
\begin{equation*}
B(x_{1}\otimes x_{1})=B(x_{1}\otimes 1_{H})(1_{A}\otimes
x_{1})-(1_{A}\otimes gx_{1})B(x_{1}\otimes 1_{H})(1_{A}\otimes
g)-(1_{A}\otimes g)B(gx_{1}x_{1}\otimes 1_{H})(1_{A}\otimes g)
\end{equation*}%
we get
\begin{equation}
B(x_{1}\otimes x_{1})=B(x_{1}\otimes 1_{H})(1_{A}\otimes
x_{1})-(1_{A}\otimes gx_{1})B(x_{1}\otimes 1_{H})(1_{A}\otimes g)
\label{form x1otx1}
\end{equation}%
Therefore we obtain%
\begin{eqnarray*}
&&B\left( x_{1}\otimes x_{1}\right) \\
&=&-2B\left( x_{1}\otimes 1_{H};1_{A},gx_{2}\right) 1_{A}\otimes gx_{1}x_{2}+
\\
&-&2B\left( x_{1}\otimes 1_{H};G,g\right) G\otimes gx_{1}+ \\
&&+2\left[ B(g\otimes 1_{H};1_{A},g)+B(x_{1}\otimes 1_{H};1_{A},gx_{1})%
\right] X_{1}\otimes gx_{1}+ \\
&&+2B(x_{1}\otimes 1_{H};1_{A},gx_{2})X_{2}\otimes gx_{1}+ \\
&&+2\left[ B(g\otimes 1_{H};GX_{2},g)+B(x_{1}\otimes 1_{H};G,gx_{1}x_{2})%
\right] GX_{1}\otimes gx_{1}x_{2}+ \\
&&+2\left[ +B(g\otimes 1_{H};G,gx_{2})-B(x_{1}\otimes 1_{H};G,gx_{1}x_{2})%
\right] GX_{1}X_{2}\otimes gx_{1}.
\end{eqnarray*}%
We write Casimir formula for $B\left( x_{1}\otimes x_{1}\right) $ and we get
\begin{eqnarray*}
&&\sum_{w_{1}=0}^{1}\sum_{a,b_{1},b_{2},d,e_{1},e_{2}=0}^{1}%
\sum_{l_{1}=0}^{b_{1}}\sum_{l_{2}=0}^{b_{2}}\sum_{u_{1}=0}^{e_{1}}%
\sum_{u_{2}=0}^{e_{2}}\left( -1\right) ^{\alpha \left(
x_{1}^{1-w_{1}};l_{1},l_{2},u_{1},u_{2}\right) } \\
&&B(g^{1+w_{1}}x_{1}^{w_{1}}\otimes
x_{1};G^{a}X_{1}^{b_{1}}X_{2}^{b_{2}},g^{d}x_{1}^{e_{1}}x_{2}^{e_{2}}) \\
&&G^{a}X_{1}^{b_{1}-l_{1}}X_{2}^{b_{2}-l_{2}}\otimes
g^{d}x_{1}^{e_{1}-u_{1}}x_{2}^{e_{2}-u_{2}}\otimes \\
&&g^{a+b_{1}+b_{2}+l_{1}+l_{2}+d+e_{1}+e_{2}+u_{1}+u_{2}}x_{1}^{l_{1}+u_{1}+1-w_{1}}x_{2}^{l_{2}+u_{2}}
\\
&=&\sum_{\omega _{1}=0}^{\nu _{1}}B^{A}(x_{1}\otimes x_{1}^{\nu _{1}-\omega
_{1}})\otimes B^{H}(x_{1}\otimes x_{1}^{\nu _{1}-\omega _{1}})\otimes
g^{1+\omega _{1}}x_{1}^{\omega _{1}}.
\end{eqnarray*}

\subsubsection{$B\left( x_{1}\otimes x_{1};1_{A},gx_{1}x_{2}\right) $}

\begin{eqnarray*}
a &=&b_{1}=b_{2}=0 \\
d &=&e_{1}=e_{2}=1
\end{eqnarray*}%
\begin{eqnarray*}
a+b_{1}+b_{2}+d+e_{1}+e &\equiv &1 \\
l_{1}+u_{1}+1-w_{1} &=&0\Rightarrow l_{1}=u_{1}=0,w_{1}=1 \\
l_{2}+u_{2} &=&1
\end{eqnarray*}%
\begin{eqnarray*}
&&\sum_{w_{1}=0}^{1}\sum_{u_{1}=0}^{1}\sum_{u_{2}=0}^{1}\left( -1\right)
^{\alpha \left( x_{1}^{1-w_{1}};0,0,u_{1},u_{2}\right) } \\
&&B(g^{1+w_{1}}x_{1}^{w_{1}}\otimes x_{1};1_{A},gx_{1}x_{2})1_{A}\otimes
gx_{1}^{1-u_{1}}x_{2}^{1-u_{2}}\otimes
g^{1+u_{1}+u_{2}}x_{1}^{u_{1}+1-w_{1}}x_{2}^{u_{2}} \\
&=&\sum_{\omega _{1}=0}^{1}B^{A}(x_{1}\otimes x_{1}^{1-\omega _{1}})\otimes
B^{H}(x_{1}\otimes x_{1}^{1-\omega _{1}})\otimes g^{1+\omega
_{1}}x_{1}^{\omega _{1}}.
\end{eqnarray*}%
\begin{eqnarray*}
&&\sum_{u_{1}=0}^{1}\sum_{u_{2}=0}^{1}\left( -1\right) ^{\alpha \left(
x_{1};0,0,u_{1},u_{2}\right) } \\
&&B(g\otimes x_{1};1_{A},gx_{1}x_{2})1_{A}\otimes
gx_{1}^{1-u_{1}}x_{2}^{1-u_{2}}\otimes
g^{1+u_{1}+u_{2}}x_{1}^{u_{1}+1}x_{2}^{u_{2}} \\
&&+\sum_{u_{1}=0}^{1}\sum_{u_{2}=0}^{1}\left( -1\right) ^{\alpha \left(
1_{H};0,0,u_{1},u_{2}\right) } \\
&&B(x_{1}^{{}}\otimes x_{1};1_{A},gx_{1}x_{2})1_{A}\otimes
gx_{1}^{1-u_{1}}x_{2}^{1-u_{2}}\otimes
g^{1+u_{1}+u_{2}}x_{1}^{u_{1}}x_{2}^{u_{2}} \\
&=&\sum_{\omega _{1}=0}^{\nu _{1}}B^{A}(x_{1}\otimes x_{1}^{\nu _{1}-\omega
_{1}})\otimes B^{H}(x_{1}\otimes x_{1}^{\nu _{1}-\omega _{1}})\otimes
g^{1+\omega _{1}}x_{1}^{\omega _{1}}.
\end{eqnarray*}

\begin{eqnarray*}
&&\left( -1\right) ^{\alpha \left( x_{1};0,0,0,0\right) }B(g\otimes
x_{1};1_{A},gx_{1}x_{2})1_{A}\otimes gx_{1}x_{2}\otimes gx_{1}+ \\
&&+\left( -1\right) ^{\alpha \left( x_{1};0,0,0,1\right) }B(g\otimes
x_{1};1_{A},gx_{1}x_{2})1_{A}\otimes gx_{1}^{{}}\otimes x_{1}x_{2} \\
&&+\left( -1\right) ^{\alpha \left( x_{1};0,0,1,0\right) }B(g\otimes
x_{1};1_{A},gx_{1}x_{2})1_{A}\otimes gx_{2}^{{}}\otimes x_{2} \\
&&+\left( -1\right) ^{\alpha \left( x_{1};0,0,1,1\right) }B(g\otimes
x_{1};1_{A},gx_{1}x_{2})1_{A}\otimes g\otimes x_{2}^{{}} \\
&&+\left( -1\right) ^{\alpha \left( 1_{H};0,0,0,0\right) }B(x_{1}\otimes
x_{1};1_{A},gx_{1}x_{2})1_{A}\otimes gx_{1}^{{}}x_{2}^{{}}\otimes g \\
&&+\left( -1\right) ^{\alpha \left( 1_{H};0,0,0,1\right) }B(x_{1}\otimes
x_{1};1_{A},gx_{1}x_{2})1_{A}\otimes gx_{1}^{{}}\otimes x_{2}^{{}} \\
&&+\left( -1\right) ^{\alpha \left( 1_{H};0,0,1,0\right) }B(x_{1}\otimes
x_{1};1_{A},gx_{1}x_{2})1_{A}\otimes gx_{2}\otimes x_{1}^{{}} \\
&&+\left( -1\right) ^{\alpha \left( 1_{H};0,0,1,1\right) }B(x_{1}\otimes
x_{1};1_{A},gx_{1}x_{2})1_{A}\otimes g\otimes gx_{1}x_{2}
\end{eqnarray*}

\paragraph{case $1_{A}\otimes gx_{1}x_{2}\otimes gx_{1}$}

\begin{eqnarray*}
&&\sum_{u_{1}=0}^{1}\sum_{u_{2}=0}^{1}\left( -1\right) ^{\alpha \left(
x_{1};0,0,u_{1},u_{2}\right) } \\
&&B(g\otimes x_{1};1_{A},gx_{1}x_{2})1_{A}\otimes
gx_{1}^{1-u_{1}}x_{2}^{1-u_{2}}\otimes
g^{1+u_{1}+u_{2}}x_{1}^{u_{1}+1}x_{2}^{u_{2}} \\
&&+\sum_{u_{1}=0}^{1}\sum_{u_{2}=0}^{1}\left( -1\right) ^{\alpha \left(
1_{H};0,0,u_{1},u_{2}\right) } \\
&&B(x_{1}\otimes x_{1};1_{A},gx_{1}x_{2})1_{A}\otimes
gx_{1}^{1-u_{1}}x_{2}^{1-u_{2}}\otimes
g^{1+u_{1}+u_{2}}x_{1}^{u_{1}}x_{2}^{u_{2}} \\
&=&B^{A}(x_{1}\otimes x_{1}^{{}})\otimes B^{H}(x_{1}\otimes
x_{1}^{{}})\otimes g+ \\
&&B^{A}(x_{1}\otimes 1_{H})\otimes B^{H}(x_{1}\otimes 1_{H})\otimes x_{1}
\end{eqnarray*}%
First summand left side gives us%
\begin{equation*}
u_{1}=u_{2}=0
\end{equation*}%
Second summand is impossible. Right side is zero.

\subsection{case $1_{H}$}

Since%
\begin{eqnarray*}
1+\omega _{1} &\equiv &0 \\
\omega _{1} &=&0
\end{eqnarray*}%
is impossible, the right side of the equality is zero. From the left side we
get%
\begin{eqnarray*}
a+b_{1}+b_{2}+l_{1}+l_{2}+d+e_{1}+e_{2}+u_{1}+u_{2} &\equiv &0 \\
l_{1}+u_{1}+1-w_{1} &=&0 \\
l_{2}+u_{2} &=&0
\end{eqnarray*}%
\begin{eqnarray*}
a+b_{1}+b_{2}+d+e_{1}+e_{2} &\equiv &0 \\
l_{1}+u_{1}+1-w_{1} &=&0\Rightarrow w_{1}=1,l_{1}=u_{1}=0 \\
l_{2}+u_{2} &=&0\Rightarrow l_{2}=u_{2}=0
\end{eqnarray*}%
and we get, since $\alpha \left( 1_{H};0,0,0,0\right) =0,$%
\begin{equation*}
\sum_{\substack{ a,b_{1},b_{2},d,e_{1},e_{2}=0  \\ %
a+b_{1}+b_{2}+d+e_{1}+e_{2}\equiv 0}}^{1}B(x_{1}\otimes
x_{1};G^{a}X_{1}^{b_{1}}X_{2}^{b_{2}},g^{d}x_{1}^{e_{1}}x_{2}^{e_{2}})G^{a}X_{1}^{b_{1}}X_{2}^{b_{2}}\otimes g^{d}x_{1}^{e_{1}}x_{2}^{e_{2}}=0
\end{equation*}%
Thus we get that%
\begin{equation}
B(x_{1}\otimes
x_{1};G^{a}X_{1}^{b_{1}}X_{2}^{b_{2}},g^{d}x_{1}^{e_{1}}x_{2}^{e_{2}})=0%
\text{ whenever }a+b_{1}+b_{2}+d+e_{1}+e_{2}\equiv 0.  \label{x1otx1, first}
\end{equation}

This is satisfied by the form of the element.

\subsection{case $g$}

From the right side we get%
\begin{eqnarray*}
1+\omega _{1} &\equiv &1 \\
\omega _{1} &=&0
\end{eqnarray*}%
i.e. $\omega _{1}=0.$From the left side we get%
\begin{eqnarray*}
a+b_{1}+b_{2}+l_{1}+l_{2}+d+e_{1}+e_{2}+u_{1}+u_{2} &\equiv &1 \\
l_{1}+u_{1}+1-w_{1} &=&0 \\
l_{2}+u_{2} &=&0
\end{eqnarray*}%
\begin{eqnarray*}
a+b_{1}+b_{2}+d+e_{1}+e_{2} &\equiv &1 \\
l_{1}+u_{1}+1-w_{1} &=&0\Rightarrow w_{1}=1,l_{1}=u_{1}=0 \\
l_{2}+u_{2} &=&0\Rightarrow l_{2}=u_{2}=0
\end{eqnarray*}%
Thus we obtain%
\begin{eqnarray*}
&&\sum_{\substack{ a,b_{1},b_{2},d,e_{1},e_{2}=0  \\ %
a+b_{1}+b_{2}+d+e_{1}+e_{2}\equiv 1}}^{1}B(x_{1}\otimes
x_{1};G^{a}X_{1}^{b_{1}}X_{2}^{b_{2}},g^{d}x_{1}^{e_{1}}x_{2}^{e_{2}})G^{a}X_{1}^{b_{1}}X_{2}^{b_{2}}\otimes g^{d}x_{1}^{e_{1}}x_{2}^{e_{2}}
\\
&=&B^{A}(x_{1}\otimes x_{1})\otimes B^{H}(x_{1}\otimes x_{1})
\end{eqnarray*}%
which already holds in view of $\left( \ref{x1otx1, first}\right) .$

\subsection{case $x_{1}$}

From the right side we get%
\begin{eqnarray*}
1+\omega _{1} &\equiv &0 \\
\omega _{1} &=&1
\end{eqnarray*}%
i.e. $\omega _{1}=1.$From the left side we get%
\begin{eqnarray*}
a+b_{1}+b_{2}+d+e_{1}+e_{2} &\equiv &w_{1} \\
l_{1}+u_{1}+1-w_{1} &=&1\Rightarrow l_{1}+u_{1}=w_{1} \\
l_{2}+u_{2} &=&0\Rightarrow l_{2}=u_{2}=0
\end{eqnarray*}%
\begin{eqnarray*}
&&\sum_{w_{1}=0}^{1}\sum_{a,b_{1},b_{2},d,e_{1},e_{2}=0}^{1}%
\sum_{l_{1}=0}^{b_{1}}\sum_{\substack{ u_{1}=0  \\ l_{1}+u_{1}=w_{1}}}%
^{e_{1}}\left( -1\right) ^{\alpha \left(
x_{1}^{1-w_{1}};l_{1},0,u_{1},0\right) } \\
&&B(g^{1+w_{1}}x_{1}^{w_{1}}\otimes
x_{1};G^{a}X_{1}^{b_{1}}X_{2}^{b_{2}},g^{d}x_{1}^{e_{1}}x_{2}^{e_{2}}) \\
&&G^{a}X_{1}^{b_{1}-l_{1}}X_{2}^{b_{2}}\otimes
g^{d}x_{1}^{e_{1}-u_{1}}x_{2}^{e_{2}} \\
&=&B^{A}(x_{1}\otimes 1_{H})\otimes B^{H}(x_{1}\otimes 1_{H})
\end{eqnarray*}%
\begin{equation*}
\alpha \left( x_{1}^{1-w_{1}};l_{1},0,u_{1},0\right)
=e_{2}u_{1}+b_{2}l_{1}+\left( a+b_{1}+b_{2}+l_{1}\right) \left(
u_{1}+1-w_{1}\right)
\end{equation*}%
We consider the left side.

\begin{eqnarray*}
&&\sum_{w_{1}=0}^{1}\sum_{\substack{ a,b_{1},b_{2},d,e_{1},e_{2}=0  \\ %
a+b_{1}+b_{2}+d+e_{1}+e_{2}\equiv w_{1}}}^{1}\sum_{l_{1}=0}^{b_{1}}\sum
_{\substack{ u_{1}=0  \\ l_{1}+u_{1}=w_{1}}}^{e_{1}}\left( -1\right)
^{\alpha \left( x_{1}^{1-w_{1}};l_{1},0,u_{1},0\right) } \\
&&B(g^{1+w_{1}}x_{1}^{w_{1}}\otimes
x_{1};G^{a}X_{1}^{b_{1}}X_{2}^{b_{2}},g^{d}x_{1}^{e_{1}}x_{2}^{e_{2}})G^{a}X_{1}^{b_{1}-l_{1}}X_{2}^{b_{2}}\otimes g^{d}x_{1}^{e_{1}-u_{1}}x_{2}^{e_{2}}
\end{eqnarray*}%
Since%
\begin{eqnarray*}
\alpha \left( x_{1};0,0,0,0\right) &=&a+b_{1}+b_{2} \\
\alpha \left( 1_{H};0,0,1,0\right) &=&e_{2}+a+b_{1}+b_{2} \\
\alpha \left( 1_{H};1,0,0,0\right) &=&b_{2}
\end{eqnarray*}%
we obtain%
\begin{eqnarray*}
&&\sum_{\substack{ a,b_{1},b_{2},d,e_{1},e_{2}=0  \\ %
a+b_{1}+b_{2}+d+e_{1}+e_{2}\equiv 0}}^{1}\left( -1\right)
^{a+b_{1}+b_{2}}B(g\otimes
x_{1};G^{a}X_{1}^{b_{1}}X_{2}^{b_{2}},g^{d}x_{1}^{e_{1}}x_{2}^{e_{2}})G^{a}X_{1}^{b_{1}}X_{2}^{b_{2}}\otimes g^{d}x_{1}^{e_{1}}x_{2}^{e_{2}}
\\
&&\sum_{\substack{ a,b_{1},b_{2},d,e_{2}=0  \\ a+b_{1}+b_{2}+d+e_{2}\equiv 0
}}^{1}\left( -1\right) ^{e_{2}+a+b_{1}+b_{2}}B(x_{1}\otimes
x_{1};G^{a}X_{1}^{b_{1}}X_{2}^{b_{2}},g^{d}x_{1}x_{2}^{e_{2}})G^{a}X_{1}^{b_{1}}X_{2}^{b_{2}}\otimes g^{d}x_{2}^{e_{2}}
\\
&&+\sum_{\substack{ a,b_{2},d,e_{1},e_{2}=0  \\ a+b_{2}+d+e_{1}+e_{2}\equiv
0 }}^{1}\left( -1\right) ^{b_{2}}B(x_{1}\otimes
x_{1};G^{a}X_{1}X_{2}^{b_{2}},g^{d}x_{1}^{e_{1}}x_{2}^{e_{2}})G^{a}X_{2}^{b_{2}}\otimes g^{d}x_{1}^{e_{1}}x_{2}^{e_{2}}
\end{eqnarray*}

\subsubsection{$G^{a}\otimes g^{d}$}

\begin{equation*}
\sum_{\substack{ a,d=0  \\ a+d\equiv 0}}^{1}\left[ \left( -1\right)
^{a}B(g\otimes x_{1};G^{a},g^{d})+\left( -1\right) ^{a}B(x_{1}\otimes
x_{1};G^{a},g^{d}x_{1})+B(x_{1}\otimes x_{1};G^{a}X_{1},g^{d})\right]
G^{a}\otimes g^{d}
\end{equation*}%
and we get%
\begin{equation*}
B(g\otimes x_{1};1_{A},1_{H})+B(x_{1}\otimes
x_{1};1_{A},x_{1})+B(x_{1}\otimes x_{1};X_{1},1_{H})-B(x_{1}\otimes
1_{H};1_{A},1_{H})=0
\end{equation*}

\begin{equation*}
-B(g\otimes x_{1};G,g)-B(x_{1}\otimes x_{1};G,gx_{1})+B(x_{1}\otimes
x_{1};GX_{1},g)-B(x_{1}\otimes 1_{H};G,g)=0
\end{equation*}

\subsubsection{$G^{a}\otimes g^{d}x_{2}$}

We obtain%
\begin{gather*}
\sum_{\substack{ a,d=0  \\ a+d\equiv 1}}^{1}\left[ \left( -1\right)
^{a}B(g\otimes x_{1};G^{a},g^{d}x_{2})+\left( -1\right) ^{1+a}B(x_{1}\otimes
x_{1};G^{a},g^{d}x_{1}x_{2})+B(x_{1}\otimes x_{1};G^{a}X_{1},g^{d}x_{2})%
\right] \\
G^{a}\otimes g^{d}x_{2}
\end{gather*}%
and we get%
\begin{equation*}
B(g\otimes x_{1};1_{A},gx_{2})-B(x_{1}\otimes
x_{1};1_{A},gx_{1}x_{2})+B(x_{1}\otimes x_{1};X_{1},gx_{2})-B(x_{1}\otimes
1_{H};1_{A},gx_{2})=0
\end{equation*}%
and%
\begin{equation*}
-B(g\otimes x_{1};G,x_{2})+B(x_{1}\otimes x_{1};G,x_{1}x_{2})+B(x_{1}\otimes
x_{1};GX_{1},x_{2})-B(x_{1}\otimes 1_{H};G,x_{2})=0
\end{equation*}%
which are satisfied in view of the form of the elements.

\subsubsection{$G^{a}\otimes g^{d}x_{1}$}

\begin{equation*}
\sum_{\substack{ a,d=0  \\ a+d\equiv 1}}^{1}\left[ \left( -1\right)
^{a}B(g\otimes x_{1};G^{a},g^{d}x_{1})+B(x_{1}\otimes
x_{1};G^{a}X_{1},g^{d}x_{1})\right] G^{a}\otimes g^{d}x_{1}
\end{equation*}%
and we get%
\begin{equation*}
B(g\otimes x_{1};1_{A},gx_{1})+B(x_{1}\otimes
x_{1};X_{1},gx_{1})-B(x_{1}\otimes 1_{H};1_{A},gx_{1})=0
\end{equation*}%
and%
\begin{equation*}
-B(g\otimes x_{1};G,x_{1})+B(x_{1}\otimes x_{1};GX_{1},x_{1})-B(x_{1}\otimes
1_{H};G,x_{1})=0
\end{equation*}%
which are satisfied in view of the form of the elements.

\subsubsection{$G^{a}X_{2}\otimes g^{d}$}

We obtain%
\begin{gather*}
\sum_{\substack{ a,d=0  \\ a+d\equiv 1}}^{1}\left[ \left( -1\right)
^{a+1}B(g\otimes x_{1};G^{a}X_{2},g^{d})+\left( -1\right)
^{a+1}B(x_{1}\otimes x_{1};G^{a}X_{2},g^{d}x_{1})-B(x_{1}\otimes
x_{1};G^{a}X_{1}X_{2},g^{d})\right] \\
G^{a}X_{2}\otimes g^{d}
\end{gather*}%
and we get%
\begin{equation*}
-B(g\otimes x_{1};X_{2},g)-B(x_{1}\otimes x_{1};X_{2},gx_{1})-B(x_{1}\otimes
x_{1};X_{1}X_{2},g)-B(x_{1}\otimes 1_{H};X_{2},g)=0
\end{equation*}%
and%
\begin{equation*}
B(g\otimes x_{1};GX_{2},1_{H})+B(x_{1}\otimes
x_{1};GX_{2},x_{1})-B(x_{1}\otimes x_{1};GX_{1}X_{2},1_{H})-B(x_{1}\otimes
1_{H};GX_{2},1_{H})=0
\end{equation*}%
which are satisfied in view of the form of the elements.

\subsubsection{$G^{a}X_{1}\otimes g^{d}$}

We obtain%
\begin{equation*}
\sum_{\substack{ a,d=0  \\ a+d\equiv 1}}^{1}\left[ \left( -1\right)
^{a+1}B(g\otimes x_{1};G^{a}X_{1},g^{d})+\left( -1\right)
^{a+1}B(x_{1}\otimes x_{1};G^{a}X_{1},g^{d}x_{1})\right] G^{a}X_{1}\otimes
g^{d}
\end{equation*}%
and we get%
\begin{equation*}
-B(g\otimes x_{1};X_{1},g)-B(x_{1}\otimes x_{1};X_{1},gx_{1})-B(x_{1}\otimes
1_{H};X_{1},g)=0
\end{equation*}%
and%
\begin{equation*}
B(g\otimes x_{1};GX_{1},1_{H})+B(x_{1}\otimes
x_{1};GX_{1},x_{1})-B(x_{1}\otimes 1_{H};GX_{1},1_{H})=0
\end{equation*}%
which are satisfied in view of the form of the elements.

\subsubsection{$G^{a}\otimes g^{d}x_{1}x_{2}$}

We obtain%
\begin{equation*}
\sum_{\substack{ a,d=0  \\ a+d\equiv 0}}^{1}\left[ \left( -1\right)
^{a}B(g\otimes x_{1};G^{a},g^{d}x_{1}x_{2})+B(x_{1}\otimes
x_{1};G^{a}X_{1},g^{d}x_{1}x_{2})\right] G^{a}\otimes g^{d}x_{1}x_{2}
\end{equation*}%
and we get%
\begin{equation*}
B(g\otimes x_{1};1_{A},x_{1}x_{2})+B(x_{1}\otimes
x_{1};X_{1},x_{1}x_{2})-B(x_{1}\otimes 1_{H};1_{A},x_{1}x_{2})=0
\end{equation*}%
and%
\begin{equation*}
-B(g\otimes x_{1};G,gx_{1}x_{2})+B(x_{1}\otimes
x_{1};GX_{1},gx_{1}x_{2})-B(x_{1}\otimes 1_{H};G,gx_{1}x_{2})=0
\end{equation*}%
which are satisfied in view of the form of the elements.

\subsubsection{$G^{a}X_{2}\otimes g^{d}x_{2}$}

We obtain%
\begin{gather*}
\sum_{\substack{ a,d=0  \\ a+d\equiv 0}}^{1}\left[ \left( -1\right)
^{a+1}B(g\otimes x_{1};G^{a}X_{2},g^{d}x_{2})+\left( -1\right)
^{a}B(x_{1}\otimes x_{1};G^{a}X_{2},g^{d}x_{1}x_{2})-B(x_{1}\otimes
x_{1};G^{a}X_{1}X_{2},g^{d}x_{2})\right] \\
G^{a}X_{2}\otimes g^{d}x_{2}
\end{gather*}%
and we get%
\begin{equation*}
-B(g\otimes x_{1};X_{2},x_{2})+B(x_{1}\otimes
x_{1};X_{2},x_{1}x_{2})-B(x_{1}\otimes
x_{1};X_{1}X_{2},x_{2})-B(x_{1}\otimes 1_{H};X_{2},x_{2})=0
\end{equation*}%
and%
\begin{equation*}
B(g\otimes x_{1};GX_{2},gx_{2})-B(x_{1}\otimes
x_{1};GX_{2},gx_{1}x_{2})-B(x_{1}\otimes
x_{1};GX_{1}X_{2},gx_{2})-B(x_{1}\otimes 1_{H};GX_{2},gx_{2})=0
\end{equation*}%
which are satisfied in view of the form of the elements.

\subsubsection{$G^{a}X_{1}\otimes g^{d}x_{2}$}

We obtain%
\begin{equation*}
\sum_{\substack{ a,d=0  \\ a+d\equiv 0}}^{1}\left[ \left( -1\right)
^{a+1}B(g\otimes x_{1};G^{a}X_{1},g^{d}x_{2})+\left( -1\right)
^{a}B(x_{1}\otimes x_{1};G^{a}X_{1},g^{d}x_{1}x_{2})\right]
G^{a}X_{1}\otimes g^{d}x_{2}
\end{equation*}%
and we get%
\begin{equation*}
-B(g\otimes x_{1};X_{1},x_{2})+B(x_{1}\otimes
x_{1};X_{1},x_{1}x_{2})-B(x_{1}\otimes 1_{H};X_{1},x_{2})=0
\end{equation*}%
and%
\begin{equation*}
B(g\otimes x_{1};GX_{1},gx_{2})-B(x_{1}\otimes
x_{1};GX_{1},gx_{1}x_{2})-B(x_{1}\otimes 1_{H};GX_{1},gx_{2})=0
\end{equation*}%
which are satisfied in view of the form of the elements.

\subsubsection{$G^{a}X_{2}\otimes g^{d}x_{1}$}

We obtain%
\begin{equation*}
\sum_{\substack{ a,d=0  \\ a+d\equiv 0}}^{1}\left[ \left( -1\right)
^{a+1}B(g\otimes x_{1};G^{a}X_{2},g^{d}x_{1})-B(x_{1}\otimes
x_{1};G^{a}X_{1}X_{2},g^{d}x_{1})\right] G^{a}X_{2}\otimes g^{d}x_{1}
\end{equation*}%
and we get%
\begin{equation*}
-B(g\otimes x_{1};X_{2},x_{1})-B(x_{1}\otimes
x_{1};X_{1}X_{2},x_{1})-B(x_{1}\otimes 1_{H};X_{2},x_{1})=0
\end{equation*}%
and%
\begin{equation*}
B(g\otimes x_{1};GX_{2},gx_{1})-B(x_{1}\otimes
x_{1};GX_{1}X_{2},gx_{1})-B(x_{1}\otimes 1_{H};GX_{2},gx_{1})=0
\end{equation*}%
which are satisfied in view of the form of the elements.

\subsubsection{$G^{a}X_{1}\otimes g^{d}x_{1}$}

We obtain%
\begin{equation*}
\sum_{\substack{ a,d=0  \\ a+d\equiv 0}}^{1}\left( -1\right)
^{a+1}B(g\otimes x_{1};G^{a}X_{1},g^{d}x_{1})G^{a}X_{1}\otimes g^{d}x_{1}
\end{equation*}%
and we get%
\begin{equation*}
-B(g\otimes x_{1};X_{1},x_{1})-B(x_{1}\otimes 1_{H};X_{1},x_{1})=0
\end{equation*}%
and%
\begin{equation*}
B(g\otimes x_{1};GX_{1},gx_{1})-B(x_{1}\otimes 1_{H};GX_{1},gx_{1})=0
\end{equation*}%
which are satisfied in view of the form of the elements.

\subsubsection{$G^{a}X_{1}X_{2}\otimes g^{d}$}

We obtain%
\begin{equation*}
\sum_{\substack{ a,d=0  \\ a+d\equiv 0}}^{1}\left[ \left( -1\right)
^{a}B(g\otimes x_{1};G^{a}X_{1}X_{2},g^{d})+\left( -1\right)
^{a}B(x_{1}\otimes x_{1};G^{a}X_{1}X_{2},g^{d}x_{1})\right]
G^{a}X_{1}X_{2}\otimes g^{d}
\end{equation*}%
and we get%
\begin{equation*}
B(g\otimes x_{1};X_{1}X_{2},1_{H})+B(x_{1}\otimes
x_{1};X_{1}X_{2},x_{1})-B(x_{1}\otimes 1_{H};X_{1}X_{2},1_{H})=0
\end{equation*}%
and%
\begin{equation*}
-B(g\otimes x_{1};GX_{1}X_{2},g)-B(x_{1}\otimes
x_{1};GX_{1}X_{2},gx_{1})-B(x_{1}\otimes 1_{H};GX_{1}X_{2},g)=0
\end{equation*}%
which are satisfied in view of the form of the elements.

\subsubsection{$G^{a}X_{2}\otimes g^{d}x_{1}x_{2}$}

\begin{equation*}
\sum_{\substack{ a,d=0  \\ a+d\equiv 1}}^{1}\left[ \left( -1\right)
^{a+1}B(g\otimes x_{1};G^{a}X_{2},g^{d}x_{1}x_{2})-B(x_{1}\otimes
x_{1};G^{a}X_{1}X_{2},g^{d}x_{1}x_{2})\right] G^{a}X_{2}\otimes
g^{d}x_{1}x_{2}
\end{equation*}%
and we get%
\begin{equation*}
-B(g\otimes x_{1};X_{2},gx_{1}x_{2})-B(x_{1}\otimes
x_{1};X_{1}X_{2},gx_{1}x_{2})-B(x_{1}\otimes 1_{H};X_{2},gx_{1}x_{2})=0
\end{equation*}%
and%
\begin{equation*}
B(g\otimes x_{1};GX_{2},x_{1}x_{2})-B(x_{1}\otimes
x_{1};GX_{1}X_{2},x_{1}x_{2})-B(x_{1}\otimes 1_{H};GX_{2},x_{1}x_{2})=0
\end{equation*}%
which are satisfied in view of the form of the elements.

\subsubsection{$G^{a}X_{1}\otimes g^{d}x_{1}x_{2}$}

We obtain%
\begin{equation*}
\sum_{\substack{ a,d=0  \\ a+d\equiv 1}}^{1}\left( -1\right)
^{a+1}B(g\otimes x_{1};G^{a}X_{1},g^{d}x_{1}x_{2})G^{a}X_{1}\otimes
g^{d}x_{1}x_{2}
\end{equation*}%
and we get%
\begin{equation*}
-B(g\otimes x_{1};X_{1},gx_{1}x_{2})-B(x_{1}\otimes
1_{H};X_{1},gx_{1}x_{2})=0
\end{equation*}%
and%
\begin{equation*}
B(g\otimes x_{1};GX_{1},x_{1}x_{2})-B(x_{1}\otimes 1_{H};GX_{1},x_{1}x_{2})=0
\end{equation*}%
which are satisfied in view of the form of the elements.

\subsubsection{$G^{a}X_{1}X_{2}\otimes g^{d}x_{2}$}

We obtain%
\begin{equation*}
\sum_{\substack{ a,d=0  \\ a+d\equiv 1}}^{1}\left[ \left( -1\right)
^{a}B(g\otimes x_{1};G^{a}X_{1}X_{2},g^{d}x_{2})+\left( -1\right)
^{a+1}B(x_{1}\otimes x_{1};G^{a}X_{1}X_{2},g^{d}x_{1}x_{2})\right]
G^{a}X_{1}X_{2}\otimes g^{d}x_{2}
\end{equation*}%
and we get%
\begin{equation*}
B(g\otimes x_{1};X_{1}X_{2},gx_{2})-B(x_{1}\otimes
x_{1};X_{1}X_{2},gx_{1}x_{2})-B(x_{1}\otimes 1_{H};X_{1}X_{2},gx_{2})=0
\end{equation*}%
and%
\begin{equation*}
-B(g\otimes x_{1};GX_{1}X_{2},x_{2})+B(x_{1}\otimes
x_{1};GX_{1}X_{2},x_{1}x_{2})-B(x_{1}\otimes 1_{H};GX_{1}X_{2},x_{2})=0
\end{equation*}%
which are satisfied in view of the form of the elements.

\subsubsection{$G^{a}X_{1}X_{2}\otimes g^{d}x_{1}$}

We obtain%
\begin{equation*}
\sum_{\substack{ a,d=0  \\ a+d\equiv 1}}^{1}\left( -1\right) ^{a}B(g\otimes
x_{1};G^{a}X_{1}X_{2},g^{d}x_{1})G^{a}X_{1}X_{2}\otimes g^{d}x_{1}
\end{equation*}%
and we get%
\begin{equation*}
B(g\otimes x_{1};X_{1}X_{2},gx_{1})-B(x_{1}\otimes 1_{H};X_{1}X_{2},gx_{1})=0
\end{equation*}%
and%
\begin{equation*}
-B(g\otimes x_{1};GX_{1}X_{2},x_{1})-B(x_{1}\otimes
1_{H};GX_{1}X_{2},x_{1})=0
\end{equation*}%
which are satisfied in view of the form of the elements.

\subsubsection{$G^{a}X_{1}X_{2}\otimes g^{d}x_{1}x_{2}$}

We obtain%
\begin{equation*}
\sum_{\substack{ a,d=0  \\ a+d\equiv 1}}^{1}\left( -1\right) ^{a}B(g\otimes
x_{1};G^{a}X_{1}X_{2},g^{d}x_{1}x_{2})G^{a}X_{1}X_{2}\otimes g^{d}x_{1}x_{2}
\end{equation*}%
and we get%
\begin{equation*}
B(g\otimes x_{1};X_{1}X_{2},x_{1}x_{2})-B(x_{1}\otimes
1_{H};X_{1}X_{2},x_{1}x_{2})=0
\end{equation*}

and%
\begin{equation*}
-B(g\otimes x_{1};GX_{1}X_{2},gx_{1}x_{2})-B(x_{1}\otimes
1_{H};GX_{1}X_{2},gx_{1}x_{2})=0
\end{equation*}%
which are satisfied in view of the form of the elements.

\subsection{case $x_{2}$}

The right side of the equality gives us zero.

For the left side of the equality we get%
\begin{eqnarray*}
a+b_{1}+b_{2}+d+e_{1}+e &\equiv &1 \\
l_{1}+u_{1}+1-w_{1} &=&0\Rightarrow l_{1}=u_{1}=0,w_{1}=1 \\
l_{2}+u_{2} &=&1
\end{eqnarray*}%
Since%
\begin{equation*}
\alpha \left( 1_{H};0,l_{2},0,u_{2}\right) =\left(
a+b_{1}+b_{2}+l_{2}\right) u_{2}.
\end{equation*}%
we get%
\begin{eqnarray*}
&&\sum_{\substack{ a,b_{1},b_{2},d,e_{1}=0  \\ a+b_{1}+b_{2}+d+e_{1}\equiv 0
}}^{1}\left( -1\right) ^{a+b_{1}+b_{2}}B(x_{1}\otimes
x_{1};G^{a}X_{1}^{b_{1}}X_{2}^{b_{2}},g^{d}x_{1}^{e_{1}}x_{2})G^{a}X_{1}^{b_{1}}X_{2}^{b_{2}}\otimes g^{d}x_{1}^{e_{1}}
\\
+ &&\sum_{\substack{ a,b_{1},d,e_{1},e_{2}=0  \\ a+b_{1}+d+e_{1}+e_{2}\equiv
0 }}^{1}B(x_{1}\otimes
x_{1};G^{a}X_{1}^{b_{1}}X_{2},g^{d}x_{1}^{e_{1}}x_{2}^{e_{2}})G^{a}X_{1}^{b_{1}}\otimes g^{d}x_{1}^{e_{1}}x_{2}^{e_{2}}
\end{eqnarray*}%
Since%
\begin{eqnarray*}
&&B\left( x_{1}\otimes x_{1}\right) \\
&=&-2B\left( x_{1}\otimes 1_{H};1_{A},gx_{2}\right) 1_{A}\otimes gx_{1}x_{2}+
\\
&-&2B\left( x_{1}\otimes 1_{H};G,g\right) G\otimes gx_{1}+ \\
&&+2\left[ B(g\otimes 1_{H};1_{A},g)+B(x_{1}\otimes 1_{H};1_{A},gx_{1})%
\right] X_{1}\otimes gx_{1}+ \\
&&+2B(x_{1}\otimes 1_{H};1_{A},gx_{2})X_{2}\otimes gx_{1}+ \\
&&+2\left[ B(g\otimes 1_{H};GX_{2},g)+B(x_{1}\otimes 1_{H};G,gx_{1}x_{2})%
\right] GX_{1},gx_{1}x_{2}+ \\
&&+2\left[ +B(g\otimes 1_{H};G,gx_{2})-B(x_{1}\otimes 1_{H};G,gx_{1}x_{2})%
\right] GX_{1}X_{2}\otimes gx_{1}.
\end{eqnarray*}%
we have to consider only the cases where the above coefficients of the
elements appear.
\begin{eqnarray*}
&&\sum_{\substack{ a,b_{1},b_{2},d,e_{1}=0  \\ a+b_{1}+b_{2}+d+e_{1}\equiv 0
}}^{1}\left( -1\right) ^{a+b_{1}+b_{2}}B(x_{1}\otimes
x_{1};G^{a}X_{1}^{b_{1}}X_{2}^{b_{2}},g^{d}x_{1}^{e_{1}}x_{2})G^{a}X_{1}^{b_{1}}X_{2}^{b_{2}}\otimes g^{d}x_{1}^{e_{1}}
\\
+ &&\sum_{\substack{ a,b_{1},d,e_{1},e_{2}=0  \\ a+b_{1}+d+e_{1}+e_{2}\equiv
0 }}^{1}B(x_{1}\otimes
x_{1};G^{a}X_{1}^{b_{1}}X_{2},g^{d}x_{1}^{e_{1}}x_{2}^{e_{2}})G^{a}X_{1}^{b_{1}}\otimes g^{d}x_{1}^{e_{1}}x_{2}^{e_{2}}=0
\end{eqnarray*}

\begin{eqnarray*}
&&\sum_{\substack{ a,b_{1},b_{2},d,e_{1}=0  \\ a+b_{1}+b_{2}+d+e_{1}\equiv 0
}}^{1}\left( -1\right) ^{a+b_{1}+b_{2}}B(x_{1}\otimes
x_{1};G^{a}X_{1}^{b_{1}}X_{2}^{b_{2}},g^{d}x_{1}^{e_{1}}x_{2})G^{a}X_{1}^{b_{1}}X_{2}^{b_{2}}\otimes g^{d}x_{1}^{e_{1}}
\\
+ &&\sum_{\substack{ a,b_{1},d,e_{1},e_{2}=0  \\ a+b_{1}+d+e_{1}+e_{2}\equiv
0 }}^{1}B(x_{1}\otimes
x_{1};G^{a}X_{1}^{b_{1}}X_{2},g^{d}x_{1}^{e_{1}}x_{2}^{e_{2}})G^{a}X_{1}^{b_{1}}\otimes g^{d}x_{1}^{e_{1}}x_{2}^{e_{2}}=0
\end{eqnarray*}

Now we look for term in $B(x_{1}\otimes x_{1};1_{A},gx_{1}x_{2})$ in the
equality above. It appears only in the first summand with these condition%
\begin{eqnarray*}
b_{1} &=&b_{2}=0 \\
e_{1} &=&1
\end{eqnarray*}%
and we obtain%
\begin{equation*}
\sum_{\substack{ a,d,=0  \\ a+d\equiv 1}}^{1}\left( -1\right)
^{a}B(x_{1}\otimes x_{1};G^{a},g^{d}x_{1}x_{2})G^{a}\otimes g^{d}x_{1}
\end{equation*}%
therefore we need to consider when in the second summand it appears the
element of the basis $G^{a}\otimes g^{d}x_{1}.$ This implies that
\begin{equation*}
e_{2}=0.
\end{equation*}%
Therefore we get
\begin{equation*}
\sum_{\substack{ a,d,=0  \\ a+d\equiv 1}}^{1}\left[ \left( -1\right)
^{a}B(x_{1}\otimes x_{1};G^{a},g^{d}x_{1}x_{2})+B(x_{1}\otimes
x_{1};G^{a}X_{2},g^{d}x_{1})\right] G^{a}\otimes g^{d}x_{1}=0
\end{equation*}%
and we get%
\begin{equation*}
B(x_{1}\otimes x_{1};1_{A},gx_{1}x_{2})+B(x_{1}\otimes x_{1};X_{2},gx_{1})=0
\end{equation*}%
and%
\begin{equation*}
-B(x_{1}\otimes x_{1};G,x_{1}x_{2})+B(x_{1}\otimes x_{1};GX_{2},x_{1})
\end{equation*}%
which are satisfied in view of the form of the element.

\paragraph{$B(x_{1}\otimes x_{1};G,gx_{1})$}

This term cannot appear in the equality above.

\paragraph{$B(x_{1}\otimes x_{1};X_{1},gx_{1})$}

This term cannot appear in the equality above.

\paragraph{$B(x_{1}\otimes x_{1};X_{2},gx_{1})$}

This term already appear in case $B(x_{1}\otimes x_{1};1_{A},gx_{1}x_{2})$%
which are satisfied in view of the form of the element.

\paragraph{$B(x_{1}\otimes x_{1};GX_{1},gx_{1}x_{2})$}

We proceed as above. From first summand we get%
\begin{equation*}
b_{1}=e_{1}=1,b_{2}=0\Rightarrow G^{a}X_{1}\otimes g^{d}x_{1}
\end{equation*}%
From the second summand we get%
\begin{equation*}
e_{2}=0
\end{equation*}

\begin{equation*}
\sum_{\substack{ a,d,=0  \\ a+d\equiv 0}}^{1}\left[ \left( -1\right)
^{a+1}B(x_{1}\otimes x_{1};G^{a}X_{1},g^{d}x_{1}x_{2})+B(x_{1}\otimes
x_{1};G^{a}X_{1}X_{2},g^{d}x_{1})\right] G^{a}X_{1}\otimes g^{d}x_{1}=0
\end{equation*}%
and we get%
\begin{equation*}
B(x_{1}\otimes x_{1};GX_{1},gx_{1}x_{2})+B(x_{1}\otimes
x_{1};GX_{1}X_{2},gx_{1})=0
\end{equation*}%
which is already known.

\paragraph{$B(x_{1}\otimes x_{1};GX_{1}X_{2},gx_{1})$}

This term already appear in the case above.

\subsection{case $x_{1}x_{2}$}

The right side is zero. The left side gives us%
\begin{eqnarray*}
&&\sum_{w_{1}=0}^{1}\sum_{a,b_{1},b_{2},d,e_{1},e_{2}=0}^{1}%
\sum_{l_{1}=0}^{b_{1}}\sum_{l_{2}=0}^{b_{2}}\sum_{u_{1}=0}^{e_{1}}%
\sum_{u_{2}=0}^{e_{2}}\left( -1\right) ^{\alpha \left(
x_{1}^{1-w_{1}};l_{1},l_{2},u_{1},u_{2}\right) } \\
&&B(g^{n_{1}+w_{1}}x_{1}^{w_{1}}\otimes
x_{1};G^{a}X_{1}^{b_{1}}X_{2}^{b_{2}},g^{d}x_{1}^{e_{1}}x_{2}^{e_{2}}) \\
&&G^{a}X_{1}^{b_{1}-l_{1}}X_{2}^{b_{2}-l_{2}}\otimes
g^{d}x_{1}^{e_{1}-u_{1}}x_{2}^{e_{2}-u_{2}}\otimes \\
&&g^{a+b_{1}+b_{2}+l_{1}+l_{2}+d+e_{1}+e_{2}+u_{1}+u_{2}}x_{1}^{l_{1}+u_{1}+1-w_{1}}x_{2}^{l_{2}+u_{2}}
\\
&=&\sum_{\omega _{1}=0}^{\nu _{1}}B^{A}(x_{1}\otimes x_{1}^{\nu _{1}-\omega
_{1}})\otimes B^{H}(x_{1}\otimes x_{1}^{\nu _{1}-\omega _{1}})\otimes
g^{1+\omega _{1}}x_{1}^{\omega _{1}}
\end{eqnarray*}%
\begin{eqnarray*}
a+b_{1}+b_{2}+l_{1}+l_{2}+d+e_{1}+e_{2}+u_{1}+u_{2} &\equiv &0 \\
l_{1}+u_{1}+1-w_{1} &=&1\Rightarrow l_{1}+u_{1}=w_{1} \\
l_{2}+u_{2} &=&1
\end{eqnarray*}%
so that%
\begin{eqnarray*}
a+b_{1}+b_{2}+d+e_{1}+e_{2} &\equiv &1+w_{1} \\
l_{1}+u_{1}+1-w_{1} &=&1\Rightarrow l_{1}+u_{1}=w_{1} \\
l_{2}+u_{2} &=&1
\end{eqnarray*}

Since%
\begin{eqnarray*}
\alpha \left( x_{1};0,0,0,1\right) &=&\left( a+b_{1}+b_{2}\right) \left(
1+1\right) +1\equiv 1 \\
\alpha \left( x_{1};0,1,0,0\right) &=&\left( a+b_{1}+b_{2}+1\right) \\
\alpha \left( 1_{H};0,0,1,1\right) &=&1+e_{2} \\
\alpha \left( 1_{H};0,1,1,0\right) &=&e_{2}+a+b_{1}+b_{2}+1 \\
\alpha \left( 1_{H};1,0,0,1\right) &=&a+b_{1} \\
\alpha \left( 1_{H};1,1,0,0\right) &=&1+b_{2}
\end{eqnarray*}%
we obtain that the left side is%
\begin{eqnarray*}
&&\sum_{\substack{ a,b_{1},b_{2},d,e_{1}=0  \\ a+b_{1}+b_{2}+d+e_{1}\equiv 0
}}^{1}-B(g\otimes
x_{1};G^{a}X_{1}^{b_{1}}X_{2}^{b_{2}},g^{d}x_{1}^{e_{1}}x_{2})G^{a}X_{1}^{b_{1}}X_{2}^{b_{2}}\otimes g^{d}x_{1}^{e_{1}}
\\
&&+\sum_{\substack{ a,b_{1},d,e_{1},e_{2}=0  \\ a+b_{1}+d+e_{1}+e_{2}\equiv
0 }}^{1}\left( -1\right) ^{a+b_{1}}B(g\otimes
x_{1};G^{a}X_{1}^{b_{1}}X_{2},g^{d}x_{1}^{e_{1}}x_{2}^{e_{2}})G^{a}X_{1}^{b_{1}}\otimes g^{d}x_{1}^{e_{1}}x_{2}^{e_{2}}
\\
+ &&\sum_{\substack{ a,b_{1},b_{2},d=0  \\ a+b_{1}+b_{2}+d\equiv 0}}%
^{1}B(x_{1}\otimes
x_{1};G^{a}X_{1}^{b_{1}}X_{2}^{b_{2}},g^{d}x_{1}x_{2})G^{a}X_{1}^{b_{1}}X_{2}^{b_{2}}\otimes g^{d}
\\
&&+\sum_{\substack{ a,b_{1},d,e_{2}=0  \\ a+b_{1}+d+e_{2}\equiv 0}}%
^{1}\left( -1\right) ^{e_{2}+a+b_{1}}B(x_{1}\otimes
x_{1};G^{a}X_{1}^{b_{1}}X_{2},g^{d}x_{1}x_{2}^{e_{2}})G^{a}X_{1}^{b_{1}}%
\otimes g^{d}x_{2}^{e_{2}} \\
&&+\sum_{\substack{ a,b_{2},d,e_{1},e_{2}=0  \\ a+b_{2}+d+e_{1}\equiv 0}}%
^{1}\left( -1\right) ^{a+1}B(x_{1}\otimes
x_{1};G^{a}X_{1}X_{2}^{b_{2}},g^{d}x_{1}^{e_{1}}x_{2})G^{a}X_{2}^{b_{2}}%
\otimes g^{d}x_{1}^{e_{1}} \\
&&+\sum_{\substack{ a,d,e_{1},e_{2}=0  \\ a+d+e_{1}+e_{2}\equiv 0}}%
^{1}B(x_{1}\otimes
x_{1};G^{a}X_{1}X_{2},g^{d}x_{1}^{e_{1}}x_{2}^{e_{2}})G^{a}\otimes
g^{d}x_{1}^{e_{1}}x_{2}^{e_{2}}
\end{eqnarray*}

\subsubsection{$G^{a}\otimes g^{d}$}

We obtain%
\begin{gather*}
\sum_{\substack{ a,d=0  \\ a+d\equiv 0}}^{1}\left[
\begin{array}{c}
-B(g\otimes x_{1};G^{a},g^{d}x_{2})+\left( -1\right) ^{a}B(g\otimes
x_{1};G^{a}X_{2},g^{d})+B(x_{1}\otimes x_{1};G^{a},g^{d}x_{1}x_{2}) \\
+\left( -1\right) ^{a}B(x_{1}\otimes x_{1};G^{a}X_{2},g^{d}x_{1})+\left(
-1\right) ^{a+1}B(x_{1}\otimes x_{1};G^{a}X_{1},g^{d}x_{2})+B(x_{1}\otimes
x_{1};G^{a}X_{1}X_{2},g^{d})%
\end{array}%
\right] \\
G^{a}\otimes g^{d}=0
\end{gather*}%
and we get%
\begin{gather*}
-B(g\otimes x_{1};1_{A},x_{2})+B(g\otimes x_{1};X_{2},1_{H})+B(x_{1}\otimes
x_{1};1_{A},x_{1}x_{2}) \\
+B(x_{1}\otimes x_{1};X_{2},x_{1})-B(x_{1}\otimes
x_{1};X_{1},x_{2})+B(x_{1}\otimes x_{1};X_{1}X_{2},1_{H})=0
\end{gather*}%
and%
\begin{gather*}
-B(g\otimes x_{1};G,gx_{2})-B(g\otimes x_{1};GX_{2},g)+B(x_{1}\otimes
x_{1};G,gx_{1}x_{2}) \\
-B(x_{1}\otimes x_{1};GX_{2},gx_{1})+B(x_{1}\otimes
x_{1};GX_{1},gx_{2})+B(x_{1}\otimes x_{1};GX_{1}X_{2},g)=0
\end{gather*}%
which are satisfied in view of the form of the elements.

\subsubsection{$G^{a}\otimes g^{d}x_{2}$}

We obtain%
\begin{gather*}
\sum_{\substack{ a,d=0  \\ a+d\equiv 1}}^{1}\left[ \left( -1\right)
^{a}B(g\otimes x_{1};G^{a}X_{2},g^{d}x_{2})+\left( -1\right)
^{1+a}B(x_{1}\otimes x_{1};G^{a}X_{2},g^{d}x_{1}x_{2})+B(x_{1}\otimes
x_{1};G^{a}X_{1}X_{2},g^{d}x_{2})\right] \\
G^{a}\otimes g^{d}x_{2}=0
\end{gather*}%
and we get%
\begin{equation*}
B(g\otimes x_{1};X_{2},gx_{2})-B(x_{1}\otimes
x_{1};X_{2},gx_{1}x_{2})+B(x_{1}\otimes x_{1};X_{1}X_{2},gx_{2})=0
\end{equation*}%
and%
\begin{equation*}
-B(g\otimes x_{1};GX_{2},x_{2})+B(x_{1}\otimes
x_{1};GX_{2},x_{1}x_{2})+B(x_{1}\otimes x_{1};GX_{1}X_{2},x_{2})=0
\end{equation*}%
which are satisfied in view of the form of the elements.

\subsubsection{$G^{a}\otimes g^{d}x_{1}$}

We obtain%
\begin{equation*}
\sum_{\substack{ a,d=0  \\ a+d\equiv 1}}^{1}\left[
\begin{array}{c}
-B(g\otimes x_{1};G^{a},g^{d}x_{1}x_{2})+\left( -1\right) ^{a}B(g\otimes
x_{1};G^{a}X_{2},g^{d}x_{1})+ \\
\left( -1\right) ^{a+1}B(x_{1}\otimes
x_{1};G^{a}X_{1},g^{d}x_{1}x_{2})+B(x_{1}\otimes
x_{1};G^{a}X_{1}X_{2},g^{d}x_{1})%
\end{array}%
\right] G^{a}\otimes g^{d}x_{1}=0
\end{equation*}%
and we get%
\begin{gather*}
-B(g\otimes x_{1};1_{A},gx_{1}x_{2})+B(g\otimes x_{1};X_{2},gx_{1})+ \\
-B(x_{1}\otimes x_{1};X_{1},gx_{1}x_{2})+B(x_{1}\otimes
x_{1};X_{1}X_{2},gx_{1})=0
\end{gather*}%
and%
\begin{gather*}
-B(g\otimes x_{1};G,x_{1}x_{2})-B(g\otimes x_{1};GX_{2},x_{1})+ \\
+B(x_{1}\otimes x_{1};GX_{1},x_{1}x_{2})+B(x_{1}\otimes
x_{1};GX_{1}X_{2},x_{1})=0
\end{gather*}%
which are satisfied in view of the form of the elements.

\subsubsection{$G^{a}X_{2}\otimes g^{d}$}

We obtain%
\begin{gather*}
\sum_{\substack{ a,d=0  \\ a+d\equiv 1}}^{1}\left[ -B(g\otimes
x_{1};G^{a}X_{2},g^{d}x_{2})+B(x_{1}\otimes
x_{1};G^{a}X_{2},g^{d}x_{1}x_{2})+\left( -1\right) ^{a+1}B(x_{1}\otimes
x_{1};G^{a}X_{1}X_{2},g^{d}x_{2})\right] \\
G^{a}X_{2}\otimes g^{d}=0
\end{gather*}%
and we get%
\begin{equation*}
-B(g\otimes x_{1};X_{2},gx_{2})+B(x_{1}\otimes
x_{1};X_{2},gx_{1}x_{2})-B(x_{1}\otimes x_{1};X_{1}X_{2},gx_{2})=0\text{ }
\end{equation*}%
and%
\begin{equation*}
-B(g\otimes x_{1};GX_{2},x_{2})+B(x_{1}\otimes
x_{1};GX_{2},x_{1}x_{2})+B(x_{1}\otimes x_{1};GX_{1}X_{2},x_{2})=0.
\end{equation*}%
Both of them were already found.

\subsubsection{$G^{a}X_{1}\otimes g^{d}$}

We obtain%
\begin{equation*}
\sum_{\substack{ a,d=0  \\ a+d\equiv 1}}^{1}\left[
\begin{array}{c}
-B(g\otimes x_{1};G^{a}X_{1},g^{d}x_{2})+\left( -1\right) ^{a+1}B(g\otimes
x_{1};G^{a}X_{1}X_{2},g^{d})+ \\
B(x_{1}\otimes x_{1};G^{a}X_{1},g^{d}x_{1}x_{2})+\left( -1\right)
^{a+1}B(x_{1}\otimes x_{1};G^{a}X_{1}X_{2},g^{d}x_{1})%
\end{array}%
\right] G^{a}X_{1}\otimes g^{d}=0
\end{equation*}%
and we get%
\begin{gather*}
-B(g\otimes x_{1};X_{1},gx_{2})-B(g\otimes x_{1};X_{1}X_{2},g)+ \\
+B(x_{1}\otimes x_{1};X_{1},gx_{1}x_{2})-B(x_{1}\otimes
x_{1};X_{1}X_{2},gx_{1})=0
\end{gather*}%
and%
\begin{gather*}
-B(g\otimes x_{1};GX_{1},x_{2})+B(g\otimes x_{1};GX_{1}X_{2},1_{H})+ \\
B(x_{1}\otimes x_{1};GX_{1},x_{1}x_{2})+B(x_{1}\otimes
x_{1};GX_{1}X_{2},x_{1})=0
\end{gather*}%
which are satisfied in view of the form of the elements.

\subsubsection{$G^{a}\otimes g^{d}x_{1}x_{2}$}

We obtain%
\begin{equation*}
\sum_{\substack{ a,d=0  \\ a+d\equiv 0}}^{1}\left[ \left( -1\right)
^{a}B(g\otimes x_{1};G^{a}X_{2},g^{d}x_{1}x_{2})+B(x_{1}\otimes
x_{1};G^{a}X_{1}X_{2},g^{d}x_{1}x_{2})\right] G^{a}\otimes g^{d}x_{1}x_{2}=0
\end{equation*}%
and we get%
\begin{equation*}
B(g\otimes x_{1};X_{2},x_{1}x_{2})+B(x_{1}\otimes
x_{1};X_{1}X_{2},x_{1}x_{2})=0
\end{equation*}%
and%
\begin{equation*}
-B(g\otimes x_{1};GX_{2},gx_{1}x_{2})+B(x_{1}\otimes
x_{1};GX_{1}X_{2},gx_{1}x_{2})=0
\end{equation*}%
which are satisfied in view of the form of the elements.

\subsubsection{$G^{a}X_{2}\otimes g^{d}x_{2}$}

This is not possible

\subsubsection{$G^{a}X_{1}\otimes g^{d}x_{2}$}

We get%
\begin{equation*}
-B(g\otimes x_{1};X_{1}X_{2},x_{2})+B(x_{1}\otimes
x_{1};X_{1}X_{2},x_{1}x_{2})=0\text{ }
\end{equation*}%
and%
\begin{equation*}
B(g\otimes x_{1};GX_{1}X_{2},gx_{2})-B(x_{1}\otimes
x_{1};GX_{1}X_{2},gx_{1}x_{2})=0.
\end{equation*}%
Both of them were already found.

\subsubsection{$G^{a}X_{2}\otimes g^{d}x_{1}$}

We obtain%
\begin{equation*}
\sum_{\substack{ a,d=0  \\ a+d\equiv 0}}^{1}\left[ -B(g\otimes
x_{1};G^{a}X_{2},g^{d}x_{1}x_{2})+\left( -1\right) ^{a+1}B(x_{1}\otimes
x_{1};G^{a}X_{1}X_{2},g^{d}x_{1}x_{2})\right] G^{a}X_{2}\otimes g^{d}x_{1}=0
\end{equation*}%
and we get%
\begin{equation*}
-B(g\otimes x_{1};X_{2},x_{1}x_{2})-B(x_{1}\otimes
x_{1};X_{1}X_{2},x_{1}x_{2})=0\text{ }
\end{equation*}%
and%
\begin{equation*}
-B(g\otimes x_{1};GX_{2},gx_{1}x_{2})+B(x_{1}\otimes
x_{1};GX_{1}X_{2},gx_{1}x_{2})=0.\text{ }
\end{equation*}%
Both of them were already found.

\subsubsection{$G^{a}X_{1}\otimes g^{d}x_{1}$}

We obtain%
\begin{equation*}
\sum_{\substack{ a,d=0  \\ a+d\equiv 0}}^{1}\left[ -B(g\otimes
x_{1};G^{a}X_{1},g^{d}x_{1}x_{2})+\left( -1\right) ^{a+1}B(g\otimes
x_{1};G^{a}X_{1}X_{2},g^{d}x_{1})\right] G^{a}X_{1}\otimes g^{d}x_{1}=0
\end{equation*}%
and we get%
\begin{equation*}
-B(g\otimes x_{1};X_{1},x_{1}x_{2})-B(g\otimes x_{1};X_{1}X_{2},x_{1})=0
\end{equation*}%
and%
\begin{equation*}
-B(g\otimes x_{1};GX_{1},gx_{1}x_{2})+B(g\otimes x_{1};GX_{1}X_{2},gx_{1})=0
\end{equation*}%
which are satisfied in view of the form of the element.

\subsubsection{$G^{a}X_{1}X_{2}\otimes g^{d}$}

We obtain%
\begin{equation*}
\sum_{\substack{ a,d=0  \\ a+d\equiv 0}}^{1}\left[ -B(g\otimes
x_{1};G^{a}X_{1}X_{2},g^{d}x_{2})+B(x_{1}\otimes
x_{1};G^{a}X_{1}X_{2},g^{d}x_{1}x_{2})\right] G^{a}X_{1}X_{2}\otimes g^{d}=0
\end{equation*}%
and we get%
\begin{equation*}
-B(g\otimes x_{1};X_{1}X_{2},x_{2})+B(x_{1}\otimes
x_{1};X_{1}X_{2},x_{1}x_{2})=0
\end{equation*}%
and%
\begin{equation*}
-B(g\otimes x_{1};GX_{1}X_{2},gx_{2})+B(x_{1}\otimes
x_{1};GX_{1}X_{2},gx_{1}x_{2})=0\text{ .}
\end{equation*}%
Both of them were already found.

\subsubsection{$G^{a}X_{2}\otimes g^{d}x_{1}x_{2}$}

This is impossible

\subsubsection{$G^{a}X_{1}\otimes g^{d}x_{1}x_{2}$}

We obtain%
\begin{equation*}
\sum_{\substack{ a,d=0  \\ a+d\equiv 1}}^{1}\left( -1\right)
^{a+1}B(g\otimes x_{1};G^{a}X_{1}X_{2},g^{d}x_{1}x_{2})G^{a}X_{1}\otimes
g^{d}x_{1}x_{2}=0
\end{equation*}%
and we get%
\begin{equation*}
-B(g\otimes x_{1};X_{1}X_{2},gx_{1}x_{2})=0\text{ }
\end{equation*}%
and%
\begin{equation*}
B(g\otimes x_{1};GX_{1}X_{2},x_{1}x_{2})=0
\end{equation*}%
Both of them were already found.

\subsubsection{$G^{a}X_{1}X_{2}\otimes g^{d}x_{2}$}

This is impossible

\subsubsection{$G^{a}X_{1}X_{2}\otimes g^{d}x_{1}$}

We obtain%
\begin{equation*}
\sum_{\substack{ a,d=0  \\ a+d\equiv 1}}^{1}-B(g\otimes
x_{1};G^{a}X_{1}X_{2},g^{d}x_{1}x_{2})G^{a}X_{1}X_{2}\otimes g^{d}x_{1}=0
\end{equation*}%
and we get%
\begin{equation*}
-B(g\otimes x_{1};X_{1}X_{2},gx_{1}x_{2})=0\text{ }
\end{equation*}%
and%
\begin{equation*}
-B(g\otimes x_{1};GX_{1}X_{2},x_{1}x_{2})=0\text{ .}
\end{equation*}%
Both of them were already found.

\subsubsection{$G^{a}X_{1}X_{2}\otimes g^{d}x_{1}x_{2}$}

This is impossible

\subsection{Case $gx_{1}$}

The right term of Casimir condition is zero. The left term gives us%
\begin{eqnarray*}
a+b_{1}+b_{2}+d+e_{1}+e_{2} &\equiv &1+w_{1} \\
l_{1}+u_{1}+1-w_{1} &=&1\Rightarrow l_{1}+u_{1}=w_{1} \\
l_{2}+u_{2} &=&0\Rightarrow l_{2}=u_{2}=0
\end{eqnarray*}%
and hence%
\begin{eqnarray*}
&&\sum_{w_{1}=0}^{1}\sum_{\substack{ a,b_{1},b_{2},d,e_{1},e_{2}=0  \\ %
a+b_{1}+b_{2}+d+e_{1}+e_{2}\equiv 1+w_{1}}}^{1}\sum_{l_{1}=0}^{b_{1}}\sum
_{\substack{ u_{1}=0  \\ l_{1}+u_{1}=w_{1}}}^{e_{1}}\left( -1\right)
^{\alpha \left( x_{1}^{1-w_{1}};l_{1},0,u_{1},0\right) } \\
&&B(g^{n_{1}+w_{1}}x_{1}^{w_{1}}\otimes
x_{1};G^{a}X_{1}^{b_{1}}X_{2}^{b_{2}},g^{d}x_{1}^{e_{1}}x_{2}^{e_{2}}) \\
&&G^{a}X_{1}^{b_{1}-l_{1}}X_{2}^{b_{2}}\otimes
g^{d}x_{1}^{e_{1}-u_{1}}x_{2}^{e_{2}}.
\end{eqnarray*}%
Thus we get%
\begin{eqnarray*}
B\left( x_{1}\otimes x_{1};gx_{1}\right) &=&\sum_{\substack{ %
a,b_{1},b_{2},d,e_{1},e_{2}=0  \\ a+b_{1}+b_{2}+d+e_{1}+e_{2}\equiv 1}}%
^{1}\left( -1\right) ^{\alpha \left( x_{1};0,0,0,0\right) }B(g\otimes
x_{1};G^{a}X_{1}^{b_{1}}X_{2}^{b_{2}},g^{d}x_{1}^{e_{1}}x_{2}^{e_{2}}) \\
&&G^{a}X_{1}^{b_{1}}X_{2}^{b_{2}}\otimes g^{d}x_{1}^{e_{1}}x_{2}^{e_{2}}
\end{eqnarray*}%
Thus we get

\begin{equation*}
B(g\otimes
x_{1};G^{a}X_{1}^{b_{1}}X_{2}^{b_{2}},g^{d}x_{1}^{e_{1}}x_{2}^{e_{2}})=0%
\text{ whenever }a+b_{1}+b_{2}+d+e_{1}+e_{2}\equiv 1
\end{equation*}%
Since by $\left( \ref{eq.10}\right) $%
\begin{equation*}
B(g\otimes x_{1})=(1_{A}\otimes g)B(1_{H}\otimes gx_{1})(1_{A}\otimes g)
\end{equation*}%
This is already known by $\left( \ref{1otgx1, first}\right) .$

\subsection{case $gx_{2}$}

The right side is zero.

The left side gives us%
\begin{eqnarray*}
a+b_{1}+b_{2}+l_{1}+l_{2}+d+e_{1}+e_{2}+u_{1}+u_{2} &\equiv &1 \\
l_{1}+u_{1}+1-w_{1} &=&0 \\
l_{2}+u_{2} &=&0
\end{eqnarray*}%
so that%
\begin{eqnarray*}
a+b_{1}+b_{2}+d+e_{1}+e_{2} &\equiv &0 \\
l_{1}+u_{1}+1-w_{1} &=&0\Rightarrow w_{1}=1,l_{1}=u_{1}=0 \\
l_{2}+u_{2} &=&1.
\end{eqnarray*}%
We get%
\begin{gather*}
\sum_{\substack{ a,b_{1},b_{2},d,e_{1},e_{2}=0  \\ %
a+b_{1}+b_{2}+d+e_{1}+e_{2}\equiv 0}}^{1}\sum_{l_{2}=0}^{b_{2}}\sum
_{\substack{ u_{2}=0  \\ l_{2}+u_{2}=1}}^{e_{2}}\left( -1\right) ^{\alpha
\left( 1_{H};0,l_{2},0,u_{2}\right) }B(x_{1}\otimes
x_{1};G^{a}X_{1}^{b_{1}}X_{2}^{b_{2}},g^{d}x_{1}^{e_{1}}x_{2}^{e_{2}}) \\
G^{a}X_{1}^{b_{1}}X_{2}^{b_{2}-l_{2}}\otimes
g^{d}x_{1}^{e_{1}}x_{2}^{e_{2}-u_{2}}=0
\end{gather*}

Since, by $\left( \ref{x1otx1, first}\right) ,B(x_{1}\otimes
x_{1};G^{a}X_{1}^{b_{1}}X_{2}^{b_{2}},g^{d}x_{1}^{e_{1}}x_{2}^{e_{2}})=0$
whenever $a+b_{1}+b_{2}+d+e_{1}+e_{2}\equiv 0,$ we get that this is already
zero.

\subsection{case $gx_{1}x_{2}$}

The right side is zero. The left side gives us%
\begin{eqnarray*}
a+b_{1}+b_{2}+l_{1}+l_{2}+d+e_{1}+e_{2}+u_{1}+u_{2} &\equiv &1 \\
l_{1}+u_{1}+1-w_{1} &\equiv &1 \\
l_{2}+u_{2} &\equiv &1.
\end{eqnarray*}%
Thus we get%
\begin{eqnarray*}
a+b_{1}+b_{2}+d+e_{1}+e_{2} &\equiv &w_{1} \\
l_{1}+u_{1}+1-w_{1} &=&1\Rightarrow l_{1}+u_{1}=w_{1} \\
l_{2}+u_{2} &=&1.
\end{eqnarray*}

As%
\begin{eqnarray*}
\alpha \left( x_{1};0,0,0,1\right) &=&\alpha \left(
g^{m}x_{1}^{n_{1}}x_{2}^{n_{2}};0,0,0,1\right) \left( a+b_{1}+b_{2}\right)
\left( 1+n_{1}+n_{2}\right) +m+n_{1} \\
&=&\left( a+b_{1}+b_{2}\right) \left( 1+1\right) +1\equiv 1 \\
\alpha \left( x_{1};0,1,0,0\right) &=&\alpha \left(
g^{m}x_{1}^{n_{1}}x_{2}^{n_{2}};0,1,0,0\right) =\left(
a+b_{1}+b_{2}+1\right) \left( n_{1}+n_{2}\right) =\left(
a+b_{1}+b_{2}+1\right) \\
\alpha \left( 1_{H};0,0,1,1\right) &=&\alpha \left(
g^{m}x_{1}^{n_{1}}x_{2}^{n_{2}};0,0,1,1\right) =1+e_{2}+\left(
a+b_{1}+b_{2}\right) \left( n_{1}+n_{2}\right) +n_{1}=1+e_{2} \\
\alpha \left( 1_{H};0,1,1,0\right) &=&\alpha \left(
g^{m}x_{1}^{n_{1}}x_{2}^{n_{2}};0,1,1,0\right) =e_{2}+\left(
a+b_{1}+b_{2}+1\right) \left( 1+n_{1}+n_{2}\right) +m \\
&=&e_{2}+a+b_{1}+b_{2}+1 \\
\alpha \left( 1_{H};1,0,0,1\right) &=&\alpha \left(
g^{m}x_{1}^{n_{1}}x_{2}^{n_{2}};1,0,0,1\right) =b_{2}+\left(
a+b_{1}+b_{2}+1\right) \left( 1+n_{1}+n_{2}\right) +m+n_{1}+(n_{2}+1) \\
&\equiv &b_{2}+\left( +b_{2}+1\right) +1\equiv a+b_{1} \\
\alpha \left( 1_{H};1,1,0,0\right) &=&\alpha \left(
g^{m}x_{1}^{n_{1}}x_{2}^{n_{2}};1,1,0,0\right) =1+b_{2}+\left(
a+b_{1}+b_{2}\right) \left( n_{1}+n_{2}\right) +n_{2}=1+b_{2}
\end{eqnarray*}%
we obtain%
\begin{gather*}
\sum_{\substack{ a,b_{1},b_{2},d,e_{1}=0  \\ a+b_{1}+b_{2}+d+e_{1}\equiv 1}}%
^{1}-B(g\otimes
x_{1};G^{a}X_{1}^{b_{1}}X_{2}^{b_{2}},g^{d}x_{1}^{e_{1}}x_{2})G^{a}X_{1}^{b_{1}}X_{2}^{b_{2}}\otimes g^{d}x_{1}^{e_{1}}+
\\
+\sum_{\substack{ a,b_{1},d,e_{1},e_{2}=0  \\ a+b_{1}+d+e_{1}+e_{2}\equiv 1}}%
^{1}\left( -1\right) ^{a+b_{1}}B(g\otimes
x_{1};G^{a}X_{1}^{b_{1}}X_{2},g^{d}x_{1}^{e_{1}}x_{2}^{e_{2}})G^{a}X_{1}^{b_{1}}\otimes g^{d}x_{1}^{e_{1}}x_{2}^{e_{2}}
\\
+\sum_{\substack{ a,b_{1},b_{2},d=0  \\ a+b_{1}+b_{2}+d\equiv 1}}%
^{1}B(x_{1}\otimes
x_{1};G^{a}X_{1}^{b_{1}}X_{2}^{b_{2}},g^{d}x_{1}x_{2})G^{a}X_{1}^{b_{1}}X_{2}^{b_{2}}\otimes g^{d}
\\
+\sum_{\substack{ a,b_{1},d,e_{2}=0  \\ a+b_{1}+d+e_{2}\equiv 1}}^{1}\left(
-1\right) ^{e_{2}+a+b_{1}}B(x_{1}\otimes
x_{1};G^{a}X_{1}^{b_{1}}X_{2},g^{d}x_{1}x_{2}^{e_{2}})G^{a}X_{1}^{b_{1}}%
\otimes g^{d}x_{2}^{e_{2}} \\
+\sum_{\substack{ a,b_{2},d,e_{1}=0  \\ a+b_{2}+d+e_{1}\equiv 1}}^{1}\left(
-1\right) ^{a+1}B(x_{1}\otimes
x_{1};G^{a}X_{1}X_{2}^{b_{2}},g^{d}x_{1}^{e_{1}}x_{2})G^{a}X_{2}^{b_{2}}%
\otimes g^{d}x_{1}^{e_{1}} \\
+\sum_{\substack{ a,d,e_{1},e_{2}=0  \\ a+d+e_{1}+e_{2}\equiv 1}}%
^{1}B(x_{1}\otimes
x_{1};G^{a}X_{1}X_{2},g^{d}x_{1}^{e_{1}}x_{2}^{e_{2}})G^{a}\otimes
g^{d}x_{1}^{e_{1}}x_{2}^{e_{2}}=0
\end{gather*}

\subsubsection{$G^{a}\otimes g^{d}$}

We obtain%
\begin{gather*}
\sum_{\substack{ a,d=0  \\ a+d\equiv 1}}^{1}\left[
\begin{array}{c}
-B(g\otimes x_{1};G^{a},g^{d}x_{2})+\left( -1\right) ^{a}B(g\otimes
x_{1};G^{a}X_{2},g^{d})+B(x_{1}\otimes x_{1};G^{a},g^{d}x_{1}x_{2}) \\
+\left( -1\right) ^{a}B(x_{1}\otimes x_{1};G^{a}X_{2},g^{d}x_{1})+\left(
-1\right) ^{a+1}B(x_{1}\otimes x_{1};G^{a}X_{1},g^{d}x_{2})+B(x_{1}\otimes
x_{1};G^{a}X_{1}X_{2},g^{d})%
\end{array}%
\right] \\
G^{a}\otimes g^{d}=0
\end{gather*}%
and we get%
\begin{gather*}
-B(g\otimes x_{1};1_{A},gx_{2})+B(g\otimes x_{1};X_{2},g)+B(x_{1}\otimes
x_{1};1_{A},gx_{1}x_{2})+ \\
+B(x_{1}\otimes x_{1};X_{2},gx_{1})-B(x_{1}\otimes
x_{1};X_{1},gx_{2})+B(x_{1}\otimes x_{1};X_{1}X_{2},g)=0
\end{gather*}%
and%
\begin{gather*}
-B(g\otimes x_{1};G,x_{2})-B(g\otimes x_{1};GX_{2},1_{H})+B(x_{1}\otimes
x_{1};G,x_{1}x_{2})+ \\
-B(x_{1}\otimes x_{1};GX_{2},x_{1})+B(x_{1}\otimes
x_{1};GX_{1},x_{2})+B(x_{1}\otimes x_{1};GX_{1}X_{2},1_{H})=0
\end{gather*}%
which are satisfied in view of the form of the element.

\subsubsection{$G^{a}\otimes g^{d}x_{2}$}

We obtain%
\begin{gather*}
\sum_{\substack{ a,d=0  \\ a+d\equiv 0}}^{1}\left[ \left( -1\right)
^{a}B(g\otimes x_{1};G^{a}X_{2},g^{d}x_{2})+\left( -1\right)
^{1+a}B(x_{1}\otimes x_{1};G^{a}X_{2},g^{d}x_{1}x_{2})+B(x_{1}\otimes
x_{1};G^{a}X_{1}X_{2},g^{d}x_{2})\right] \\
G^{a}\otimes g^{d}x_{2}=0
\end{gather*}%
and we get%
\begin{equation*}
B(g\otimes x_{1};X_{2},x_{2})-B(x_{1}\otimes
x_{1};X_{2},x_{1}x_{2})+B(x_{1}\otimes x_{1};X_{1}X_{2},x_{2})=0
\end{equation*}%
and%
\begin{equation*}
-B(g\otimes x_{1};GX_{2},gx_{2})+B(x_{1}\otimes
x_{1};GX_{2},gx_{1}x_{2})+B(x_{1}\otimes x_{1};GX_{1}X_{2},gx_{2})=0
\end{equation*}%
which are satisfied in view of the form of the element.

\subsubsection{$G^{a}\otimes g^{d}x_{1}$}

We obtain%
\begin{gather*}
\sum_{\substack{ a,d=0  \\ a+d\equiv 0}}^{1}\left[
\begin{array}{c}
-B(g\otimes x_{1};G^{a},g^{d}x_{1}x_{2})+\left( -1\right) ^{a}B(g\otimes
x_{1};G^{a}X_{2},g^{d}x_{1})+ \\
\left( -1\right) ^{a+1}B(x_{1}\otimes
x_{1};G^{a}X_{1},g^{d}x_{1}x_{2})+B(x_{1}\otimes
x_{1};G^{a}X_{1}X_{2},g^{d}x_{1})%
\end{array}%
\right] \\
G^{a}\otimes g^{d}x_{1}=0
\end{gather*}%
and we get%
\begin{gather*}
-B(g\otimes x_{1};1_{A},x_{1}x_{2})+B(g\otimes x_{1};X_{2},x_{1})+ \\
-B(x_{1}\otimes x_{1};X_{1},x_{1}x_{2})+B(x_{1}\otimes
x_{1};X_{1}X_{2},x_{1})=0
\end{gather*}%
and%
\begin{gather*}
-B(g\otimes x_{1};G,gx_{1}x_{2})-B(g\otimes x_{1};GX_{2},gx_{1})+ \\
+B(x_{1}\otimes x_{1};GX_{1},gx_{1}x_{2})+B(x_{1}\otimes
x_{1};GX_{1}X_{2},gx_{1})=0
\end{gather*}%
which are satisfied in view of the form of the elements.

\subsubsection{$G^{a}X_{2}\otimes g^{d}$}

We obtain%
\begin{gather*}
\sum_{\substack{ a,d=0  \\ a+d\equiv 0}}^{1}\left[ -B(g\otimes
x_{1};G^{a}X_{2},g^{d}x_{2})+B(x_{1}\otimes
x_{1};G^{a}X_{2},g^{d}x_{1}x_{2})+\left( -1\right) ^{a+1}B(x_{1}\otimes
x_{1};G^{a}X_{1}X_{2},g^{d}x_{2})\right] \\
G^{a}X_{2}\otimes g^{d}=0
\end{gather*}%
and we get%
\begin{equation*}
-B(g\otimes x_{1};X_{2},x_{2})+B(x_{1}\otimes
x_{1};X_{2},x_{1}x_{2})-B(x_{1}\otimes x_{1};X_{1}X_{2},x_{2})=0\text{ }
\end{equation*}%
and%
\begin{equation*}
-B(g\otimes x_{1};GX_{2},gx_{2})+B(x_{1}\otimes
x_{1};GX_{2},gx_{1}x_{2})+B(x_{1}\otimes x_{1};GX_{1}X_{2},gx_{2})=0.\text{ }
\end{equation*}%
Both of them were already found.

\subsubsection{$G^{a}X_{1}\otimes g^{d}$}

We obtain%
\begin{equation*}
\sum_{\substack{ a,d=0  \\ a+d\equiv 0}}^{1}\left[
\begin{array}{c}
-B(g\otimes x_{1};G^{a}X_{1},g^{d}x_{2})+\left( -1\right) ^{a+1}B(g\otimes
x_{1};G^{a}X_{1}X_{2},g^{d}) \\
+B(x_{1}\otimes x_{1};G^{a}X_{1},g^{d}x_{1}x_{2})+\left( -1\right)
^{a+1}B(x_{1}\otimes x_{1};G^{a}X_{1}X_{2},g^{d}x_{1})+%
\end{array}%
\right] G^{a}X_{1}\otimes g^{d}=0
\end{equation*}%
and we get%
\begin{gather*}
-B(g\otimes x_{1};X_{1},x_{2})-B(g\otimes x_{1};X_{1}X_{2},1_{H})+ \\
+B(x_{1}\otimes x_{1};X_{1},x_{1}x_{2})-B(x_{1}\otimes
x_{1};X_{1}X_{2},x_{1})=0
\end{gather*}%
and%
\begin{gather*}
-B(g\otimes x_{1};GX_{1},gx_{2})+B(g\otimes x_{1};GX_{1}X_{2},g)+ \\
+B(x_{1}\otimes x_{1};GX_{1},gx_{1}x_{2})+B(x_{1}\otimes
x_{1};GX_{1}X_{2},gx_{1})=0
\end{gather*}%
which are satisfied in view of the form of the elements.

\subsubsection{$G^{a}\otimes g^{d}x_{1}x_{2}$}

We obtain%
\begin{equation*}
\sum_{\substack{ a,d=0  \\ a+d\equiv 1}}^{1}\left[ \left( -1\right)
^{a}B(g\otimes x_{1};G^{a}X_{2},g^{d}x_{1}x_{2})+B(x_{1}\otimes
x_{1};G^{a}X_{1}X_{2},g^{d}x_{1}x_{2})\right] G^{a}\otimes g^{d}x_{1}x_{2}=0
\end{equation*}%
and we get%
\begin{equation*}
B(g\otimes x_{1};X_{2},gx_{1}x_{2})+B(x_{1}\otimes
x_{1};X_{1}X_{2},gx_{1}x_{2})=0
\end{equation*}%
and%
\begin{equation*}
-B(g\otimes x_{1};GX_{2},x_{1}x_{2})+B(x_{1}\otimes
x_{1};GX_{1}X_{2},x_{1}x_{2})=0
\end{equation*}%
which are satisfied in view of the form of the elements.

\subsubsection{$G^{a}X_{2}\otimes g^{d}x_{2}$}

This is impossible

\subsubsection{$G^{a}X_{1}\otimes g^{d}x_{2}$}

We obtain%
\begin{equation*}
\sum_{\substack{ a,d=0  \\ a+d\equiv 1}}^{1}\left[ \left( -1\right)
^{a+1}B(g\otimes x_{1};G^{a}X_{1}X_{2},g^{d}x_{2})+\left( -1\right)
^{a}B(x_{1}\otimes x_{1};G^{a}X_{1}X_{2},g^{d}x_{1}x_{2})\right]
G^{a}X_{1}\otimes g^{d}x_{2}=0
\end{equation*}%
and we get%
\begin{equation*}
-B(g\otimes x_{1};X_{1}X_{2},gx_{2})+B(x_{1}\otimes
x_{1};X_{1}X_{2},gx_{1}x_{2})=0
\end{equation*}%
and%
\begin{equation*}
B(g\otimes x_{1};GX_{1}X_{2},x_{2})-B(x_{1}\otimes
x_{1};GX_{1}X_{2},x_{1}x_{2})=0
\end{equation*}%
which are satisfied in view of the form of the elements.

\subsubsection{$G^{a}X_{2}\otimes g^{d}x_{1}$}

We obtain%
\begin{equation*}
\sum_{\substack{ a,d=0  \\ a+d\equiv 1}}^{1}\left[ -B(g\otimes
x_{1};G^{a}X_{2},g^{d}x_{1}x_{2})+\left( -1\right) ^{a+1}B(x_{1}\otimes
x_{1};G^{a}X_{1}X_{2},g^{d}x_{1}x_{2})\right] G^{a}X_{2}\otimes g^{d}x_{1}=0
\end{equation*}%
and we get%
\begin{equation*}
-B(g\otimes x_{1};X_{2},gx_{1}x_{2})-B(x_{1}\otimes
x_{1};X_{1}X_{2},gx_{1}x_{2})=0\text{ }
\end{equation*}%
and%
\begin{equation*}
-B(g\otimes x_{1};GX_{2},x_{1}x_{2})+B(x_{1}\otimes
x_{1};GX_{1}X_{2},x_{1}x_{2})=0\text{.}
\end{equation*}%
Both of them were already found.

\subsubsection{$G^{a}X_{1}\otimes g^{d}x_{1}$}

We obtain%
\begin{equation*}
\sum_{\substack{ a,d=0  \\ a+d\equiv 1}}^{1}\left[ -B(g\otimes
x_{1};G^{a}X_{1},g^{d}x_{1}x_{2})+\left( -1\right) ^{a+1}B(g\otimes
x_{1};G^{a}X_{1}X_{2},g^{d}x_{1})\right] G^{a}X_{1}\otimes g^{d}x_{1}=0
\end{equation*}%
and we get%
\begin{equation*}
-B(g\otimes x_{1};X_{1},gx_{1}x_{2})-B(g\otimes x_{1};X_{1}X_{2},gx_{1})=0
\end{equation*}%
and%
\begin{equation*}
-B(g\otimes x_{1};GX_{1},x_{1}x_{2})+B(g\otimes x_{1};GX_{1}X_{2},x_{1})=0%
\text{ }
\end{equation*}%
which are satisfied in view of the form of the element.

\subsubsection{$G^{a}X_{1}X_{2}\otimes g^{d}$}

We obtain%
\begin{equation*}
\sum_{\substack{ a,d=0  \\ a+d\equiv 1}}^{1}\left[ -B(g\otimes
x_{1};G^{a}X_{1}X_{2},g^{d}x_{2})+B(x_{1}\otimes
x_{1};G^{a}X_{1}X_{2},g^{d}x_{1}x_{2})\right] G^{a}X_{1}X_{2}\otimes g^{d}=0
\end{equation*}%
and we get%
\begin{equation*}
-B(g\otimes x_{1};X_{1}X_{2},gx_{2})+B(x_{1}\otimes
x_{1};X_{1}X_{2},gx_{1}x_{2})=0\text{ }
\end{equation*}%
and%
\begin{equation*}
-B(g\otimes x_{1};GX_{1}X_{2},x_{2})+B(x_{1}\otimes
x_{1};GX_{1}X_{2},x_{1}x_{2})=0\text{.}
\end{equation*}%
Both of them were already found.

\subsubsection{$G^{a}X_{2}\otimes g^{d}x_{1}x_{2}$}

This is impossible

\subsubsection{$G^{a}X_{1}\otimes g^{d}x_{1}x_{2}$}

We obtain%
\begin{equation*}
\sum_{\substack{ a,d=0  \\ a+d\equiv 0}}^{1}\left( -1\right)
^{a+1}B(g\otimes x_{1};G^{a}X_{1}X_{2},g^{d}x_{1}x_{2})G^{a}X_{1}\otimes
g^{d}x_{1}x_{2}=0
\end{equation*}%
and we get%
\begin{equation*}
-B(g\otimes x_{1};X_{1}X_{2},x_{1}x_{2})=0
\end{equation*}%
and%
\begin{equation*}
B(g\otimes x_{1};GX_{1}X_{2},gx_{1}x_{2})=0.
\end{equation*}%
which are satisfied in view of the form of the element.

\subsubsection{$G^{a}X_{1}X_{2}\otimes g^{d}x_{2}$}

This is impossible

\subsubsection{$G^{a}X_{1}X_{2}\otimes g^{d}x_{1}$}

We obtain%
\begin{equation*}
\sum_{\substack{ a,d=0  \\ a+d\equiv 0}}^{1}-B(g\otimes
x_{1};G^{a}X_{1}X_{2},g^{d}x_{1}x_{2})G^{a}X_{1}X_{2}\otimes g^{d}x_{1}=0
\end{equation*}%
and we get%
\begin{equation*}
-B(g\otimes x_{1};X_{1}X_{2},x_{1}x_{2})=0\text{ }
\end{equation*}%
and%
\begin{equation*}
B(g\otimes x_{1};GX_{1}X_{2},gx_{1}x_{2})=0\text{.}
\end{equation*}%
Both of them were already found.

\subsubsection{$G^{a}X_{1}X_{2}\otimes g^{d}x_{1}x_{2}$}

This is impossible.

\section{$B\left( x_{1}\otimes x_{2}\right) $}

From formula $\left( \ref{simpl2i2}\right) $ we get%
\begin{eqnarray}
B(x_{1}\otimes x_{2}) &=&B(x_{1}\otimes 1_{H})(1_{A}\otimes
x_{2})-(1_{A}\otimes gx_{2})B(x_{1}\otimes 1_{H})(1_{A}\otimes g)
\label{form x1otx2} \\
&&+(1_{A}\otimes g)B(gx_{1}x_{2}\otimes 1_{H})(1_{A}\otimes g)  \notag
\end{eqnarray}%
and we obtain

We write the Casimir condition for $x_{1}\otimes x_{2}$

\begin{eqnarray*}
&&\sum_{w_{1}=0}^{1}\sum_{a,b_{1},b_{2},d,e_{1},e_{2}=0}^{1}%
\sum_{l_{1}=0}^{b_{1}}\sum_{l_{2}=0}^{b_{2}}\sum_{u_{1}=0}^{e_{1}}%
\sum_{u_{2}=0}^{e_{2}}\left( -1\right) ^{\alpha \left(
x_{1}^{1-w_{1}};l_{1},l_{2},u_{1},u_{2}\right) } \\
&&B(g^{1+w_{1}}x_{1}^{w_{1}}\otimes
x_{2};G^{a}X_{1}^{b_{1}}X_{2}^{b_{2}},g^{d}x_{1}^{e_{1}}x_{2}^{e_{2}})G^{a}X_{1}^{b_{1}-l_{1}}X_{2}^{b_{2}-l_{2}}\otimes g^{d}x_{1}^{e_{1}-u_{1}}x_{2}^{e_{2}-u_{2}}\otimes
\\
&&g^{a+b_{1}+b_{2}+l_{1}+l_{2}+d+e_{1}+e_{2}+u_{1}+u_{2}}x_{1}^{l_{1}+u_{1}+1-w_{1}}x_{2}^{l_{2}+u_{2}}
\\
&=&\sum_{\omega _{2}=0}^{1}B^{A}(x_{1}\otimes x_{2}^{1-\omega _{2}})\otimes
B^{H}(x_{1}\otimes x_{2}^{1-\omega _{2}})\otimes g^{1+\omega
_{2}}x_{2}^{\omega _{2}}.
\end{eqnarray*}

\subsection{Case $1_{H}$}

We get from the left side

\begin{eqnarray*}
l_{1}+u_{1}+1-w_{1} &=&0\Rightarrow w_{1}=1,l_{1}=u_{1}=0 \\
l_{2} &=&u_{2}=0 \\
a+b_{1}+b_{2}+d+e_{1}+e_{2} &\equiv &0
\end{eqnarray*}%
and we obtain%
\begin{gather*}
\sum_{\substack{ a,b_{1},b_{2},d,e_{1},e_{2}=0  \\ %
a+b_{1}+b_{2}+d+e_{1}+e_{2}\equiv 0}}^{1}\left( -1\right) ^{\alpha \left(
1_{H};0,0,0,0\right) }B(x_{1}\otimes
x_{2};G^{a}X_{1}^{b_{1}}X_{2}^{b_{2}},g^{d}x_{1}^{e_{1}}x_{2}^{e_{2}})G^{a}X_{1}^{b_{1}}X_{2}^{b_{2}}\otimes g^{d}x_{1}^{e_{1}}x_{2}^{e_{2}}
\\
\otimes 1_{H}=0
\end{gather*}%
Since there is nothing on the right side we get%
\begin{equation*}
B(x_{1}\otimes
x_{2};G^{a}X_{1}^{b_{1}}X_{2}^{b_{2}},g^{d}x_{1}^{e_{1}}x_{2}^{e_{2}})=0%
\text{ for every }a+b_{1}+b_{2}+d+e_{1}+e_{2}\equiv 0
\end{equation*}%
and this condition is satisfied in view of the form of the element.

\subsection{Case $g$}

As usual this gives no further information.

\subsection{Case $x_{1}$}

Nothing from right side. The left side gives us%
\begin{eqnarray*}
l_{1}+u_{1}+1-w_{1} &=&1\Rightarrow l_{1}+u_{1}=w_{1} \\
l_{2} &=&u_{2}=0 \\
a+b_{1}+b_{2}+d+e_{1}+e_{2} &\equiv &w_{1}
\end{eqnarray*}%
Since%
\begin{eqnarray*}
&&\alpha \left( x_{1};0,0,0,0\right) =a+b_{1}+b_{2} \\
&&\alpha \left( 1_{H};0,0,1,0\right) =e_{2}+\left( a+b_{1}+b_{2}\right) \\
&&\alpha \left( 1_{H};1,0,0,0\right) =b_{2}
\end{eqnarray*}%
we get that the equality is%
\begin{eqnarray*}
&&\sum_{\substack{ a,b_{1},b_{2},d,e_{1},e_{2}=0  \\ %
a+b_{1}+b_{2}+d+e_{1}+e_{2}\equiv 0}}^{1}\left( -1\right)
^{a+b_{1}+b_{2}}B(g\otimes
x_{2};G^{a}X_{1}^{b_{1}}X_{2}^{b_{2}},g^{d}x_{1}^{e_{1}}x_{2}^{e_{2}})G^{a}X_{1}^{b_{1}}X_{2}^{b_{2}}\otimes g^{d}x_{1}^{e_{1}}x_{2}^{e_{2}}
\\
&&\sum_{\substack{ a,b_{1},b_{2},d,e_{2}=0  \\ a+b_{1}+b_{2}+d+e_{2}\equiv 0
}}^{1}\left( -1\right) ^{e_{2}+\left( a+b_{1}+b_{2}\right) }B(x_{1}\otimes
x_{2};G^{a}X_{1}^{b_{1}}X_{2}^{b_{2}},g^{d}x_{1}x_{2}^{e_{2}})G^{a}X_{1}^{b_{1}}X_{2}^{b_{2}}\otimes g^{d}x_{2}^{e_{2}}+
\\
&&\sum_{\substack{ a,b_{2},d,e_{1},e_{2}=0  \\ a+b_{2}+d+e_{1}+e_{2}\equiv 0
}}^{1}\left( -1\right) ^{b_{2}}B(x_{1}\otimes
x_{2};G^{a}X_{1}X_{2}^{b_{2}},g^{d}x_{1}^{e_{1}}x_{2}^{e_{2}})G^{a}X_{2}^{b_{2}}\otimes g^{d}x_{1}^{e_{1}}x_{2}^{e_{2}}=0
\end{eqnarray*}

\subsubsection{$G^{a}\otimes g^{d}$}

We obtain%
\begin{equation*}
\sum_{\substack{ a,d=0  \\ a+d\equiv 0}}^{1}\left[
\begin{array}{c}
\left( -1\right) ^{a}B(g\otimes x_{2};G^{a},g^{d})+\left( -1\right)
^{a}B(x_{1}\otimes x_{2};G^{a},g^{d}x_{1}) \\
+B(x_{1}\otimes x_{2};G^{a}X_{1},g^{d})%
\end{array}%
\right] G^{a}\otimes g^{d}=0
\end{equation*}%
and we get%
\begin{equation*}
B(g\otimes x_{2};1_{A},1_{H})+B(x_{1}\otimes
x_{2};1_{A},x_{1})+B(x_{1}\otimes x_{2};X_{1},1_{H})=0
\end{equation*}

and%
\begin{equation*}
-B(g\otimes x_{2};G,g)-B(x_{1}\otimes x_{2};G,gx_{1})+B(x_{1}\otimes
x_{2};GX_{1},g)=0
\end{equation*}%
which are satisfied in view of the form of the elements.

\subsubsection{$G^{a}\otimes g^{d}x_{2}$}

We obtain%
\begin{equation*}
\sum_{\substack{ a,d=0  \\ a+d\equiv 1}}^{1}\left[
\begin{array}{c}
\left( -1\right) ^{a}B(g\otimes x_{2};G^{a},g^{d}x_{2})+\left( -1\right)
^{1+a}B(x_{1}\otimes x_{2};G^{a},g^{d}x_{1}x_{2}) \\
+B(x_{1}\otimes x_{2};G^{a}X_{1},g^{d}x_{2})%
\end{array}%
\right] G^{a}\otimes g^{d}x_{2}=0
\end{equation*}%
and we get%
\begin{equation*}
B(g\otimes x_{2};1_{A},gx_{2})-B(x_{1}\otimes
x_{2};1_{A},gx_{1}x_{2})+B(x_{1}\otimes x_{2};X_{1},gx_{2})=0
\end{equation*}%
and%
\begin{equation*}
-B(g\otimes x_{2};G,x_{2})+B(x_{1}\otimes x_{2};G,x_{1}x_{2})+B(x_{1}\otimes
x_{2};GX_{1},x_{2})=0
\end{equation*}%
which are satisfied in view of the form of the elements.

\subsubsection{$G^{a}\otimes g^{d}x_{1}$}

We obtain%
\begin{equation*}
\sum_{\substack{ a,d=0  \\ a+d\equiv 1}}^{1}\left[ \left( -1\right)
^{a}B(g\otimes x_{2};G^{a},g^{d}x_{1})+B(x_{1}\otimes
x_{2};G^{a}X_{1},g^{d}x_{1})\right] G^{a}\otimes g^{d}x_{1}=0
\end{equation*}%
Thus we get%
\begin{equation*}
B(g\otimes x_{2};1_{A},gx_{1})+B(x_{1}\otimes x_{2};X_{1},gx_{1})=0
\end{equation*}%
and%
\begin{equation*}
-B(g\otimes x_{2};G,x_{1})+B(x_{1}\otimes x_{2};GX_{1},x_{1})=0
\end{equation*}%
which hold in view of the form of the elements.

\subsubsection{$G^{a}X_{2}\otimes g^{d}$}

We obtain%
\begin{equation*}
\sum_{\substack{ a,d=0  \\ a+d\equiv 1}}^{1}\left[
\begin{array}{c}
\left( -1\right) ^{a+1}B(g\otimes x_{2};G^{a}X_{2},g^{d})+ \\
\left( -1\right) ^{a+1}B(x_{1}\otimes
x_{2};G^{a}X_{2},g^{d}x_{1})-B(x_{1}\otimes x_{2};G^{a}X_{1}X_{2},g^{d})%
\end{array}%
\right] G^{a}X_{2}\otimes g^{d}=0
\end{equation*}%
and we get%
\begin{equation*}
-B(g\otimes x_{2};X_{2},g)-B(x_{1}\otimes x_{2};X_{2},gx_{1})-B(x_{1}\otimes
x_{2};X_{1}X_{2},g)=0
\end{equation*}%
and%
\begin{equation*}
B(g\otimes x_{2};GX_{2},1_{H})+B(x_{1}\otimes
x_{2};GX_{2},x_{1})-B(x_{1}\otimes x_{2};GX_{1}X_{2},1_{H})=0
\end{equation*}%
which hold in view of the form of the elements.

\subsubsection{$G^{a}X_{1}\otimes g^{d}$}

\begin{equation*}
\sum_{\substack{ a,d=0  \\ a+d\equiv 1}}^{1}\left[
\begin{array}{c}
\left( -1\right) ^{a+1}B(g\otimes x_{2};G^{a}X_{1},g^{d}) \\
+\left( -1\right) ^{a+1}B(x_{1}\otimes x_{2};G^{a}X_{1},g^{d}x_{1})%
\end{array}%
\right] G^{a}X_{1}\otimes g^{d}=0
\end{equation*}%
We get%
\begin{equation*}
-B(g\otimes x_{2};X_{1},g)-B(x_{1}\otimes x_{2};X_{1},gx_{1})=0
\end{equation*}%
and%
\begin{equation*}
B(g\otimes x_{2};GX_{1},1_{H})+B(x_{1}\otimes x_{2};GX_{1},x_{1})=0
\end{equation*}%
which hold in view of the form of the elements.

\subsubsection{$G^{a}\otimes g^{d}x_{1}x_{2}$}

We obtain%
\begin{equation*}
\sum_{\substack{ a,d=0  \\ a+d\equiv 0}}^{1}\left[
\begin{array}{c}
\left( -1\right) ^{a}B(g\otimes x_{2};G^{a},g^{d}x_{1}x_{2}) \\
+B(x_{1}\otimes x_{2};G^{a}X_{1},g^{d}x_{1}x_{2})%
\end{array}%
\right] G^{a}\otimes g^{d}x_{1}x_{2}=0
\end{equation*}%
and we get%
\begin{equation*}
B(g\otimes x_{2};1_{A},x_{1}x_{2})+B(x_{1}\otimes x_{2};X_{1},x_{1}x_{2})=0
\end{equation*}%
and%
\begin{equation*}
-B(g\otimes x_{2};G,gx_{1}x_{2})+B(x_{1}\otimes x_{2};GX_{1},gx_{1}x_{2})=0
\end{equation*}%
which hold in view of the form of the elements.

\subsubsection{$G^{a}X_{2}\otimes g^{d}x_{2}$}

We obtain%
\begin{gather*}
\sum_{\substack{ a,d=0  \\ a+d\equiv 0}}^{1}\left[
\begin{array}{c}
\left( -1\right) ^{a+1}B(g\otimes x_{2};G^{a}X_{2},g^{d}x_{2})+ \\
\left( -1\right) ^{a}B(x_{1}\otimes
x_{2};G^{a}X_{2},g^{d}x_{1}x_{2})-B(x_{1}\otimes
x_{2};G^{a}X_{1}X_{2},g^{d}x_{2})%
\end{array}%
\right] \\
G^{a}X_{2}\otimes g^{d}x_{2}=0
\end{gather*}%
and we get%
\begin{equation*}
-B(g\otimes x_{2};X_{2},x_{2})+B(x_{1}\otimes
x_{2};X_{2},x_{1}x_{2})-B(x_{1}\otimes x_{2};X_{1}X_{2},x_{2})=0
\end{equation*}%
and%
\begin{equation*}
B(g\otimes x_{2};GX_{2},gx_{2})-B(x_{1}\otimes
x_{2};GX_{2},gx_{1}x_{2})-B(x_{1}\otimes x_{2};GX_{1}X_{2},gx_{2})=0
\end{equation*}%
which hold in view of the form of the elements.

\subsubsection{$G^{a}X_{1}\otimes g^{d}x_{2}$}

\begin{equation*}
\sum_{\substack{ a,d=0  \\ a+d\equiv 0}}^{1}\left[ B(g\otimes
x_{2};G^{a}X_{1},g^{d}x_{2})+\left( -1\right) ^{a}B(x_{1}\otimes
x_{2};G^{a}X_{1},g^{d}x_{1})\right] G^{a}X_{1}\otimes g^{d}x_{2}=0
\end{equation*}%
We get%
\begin{equation*}
-B(g\otimes x_{2};X_{1},x_{2})+B(x_{1}\otimes x_{2};X_{1},x_{1}x_{2})=0
\end{equation*}%
and%
\begin{equation*}
B(g\otimes x_{2};GX_{1},gx_{2})-B(x_{1}\otimes x_{2};GX_{1},gx_{1}x_{2})=0
\end{equation*}%
which hold in view of the form of the elements.

\subsubsection{$G^{a}X_{2}\otimes g^{d}x_{1}$}

We obtain%
\begin{equation*}
\sum_{\substack{ a,d=0  \\ a+d\equiv 0}}^{1}\left[ \left( -1\right)
^{a+1}B(g\otimes x_{2};G^{a}X_{2},g^{d}x_{1})-B(x_{1}\otimes
x_{2};G^{a}X_{1}X_{2},g^{d}x_{1})\right] G^{a}X_{2}\otimes g^{d}x_{1}=0
\end{equation*}%
and we get%
\begin{equation*}
-B(g\otimes x_{2};X_{2},x_{1})-B(x_{1}\otimes x_{2};X_{1}X_{2},x_{1})=0
\end{equation*}%
and%
\begin{equation*}
+B(g\otimes x_{2};GX_{2},gx_{1})-B(x_{1}\otimes x_{2};GX_{1}X_{2},gx_{1})=0
\end{equation*}%
which hold in view of the form of the elements.

\subsubsection{$G^{a}X_{1}\otimes g^{d}x_{1}$}

\begin{equation*}
\sum_{\substack{ a,d=0  \\ a+d\equiv 0}}^{1}\left( -1\right)
^{a+1}B(g\otimes x_{2};G^{a}X_{1},g^{d}x_{1})G^{a}X_{1}\otimes g^{d}x_{1}=0
\end{equation*}%
We get%
\begin{equation*}
-B(g\otimes x_{2};X_{1},x_{1})=0
\end{equation*}%
and%
\begin{equation*}
B(g\otimes x_{2};GX_{1},gx_{1})=0
\end{equation*}%
which hold in view of the form of the element.

\subsubsection{$G^{a}X_{1}X_{2}\otimes g^{d}$}

We obtain%
\begin{equation*}
\sum_{\substack{ a,d=0  \\ a+d\equiv 0}}^{1}\left[
\begin{array}{c}
\left( -1\right) ^{a}B(g\otimes x_{2};G^{a}X_{1}X_{2})+ \\
\left( -1\right) ^{a}B(x_{1}\otimes x_{2};G^{a}X_{1}X_{2},g^{d}x_{1})%
\end{array}%
\right] G^{a}X_{1}X_{2}\otimes g^{d}=0
\end{equation*}%
and we get%
\begin{equation*}
B(g\otimes x_{2};X_{1}X_{2},1_{H})+B(x_{1}\otimes x_{2};X_{1}X_{2},x_{1})=0
\end{equation*}%
and%
\begin{equation*}
-B(g\otimes x_{2};GX_{1}X_{2},g)-B(x_{1}\otimes x_{2};GX_{1}X_{2},gx_{1})=0
\end{equation*}%
which hold in view of the form of the elements.

\subsubsection{$G^{a}X_{2}\otimes g^{d}x_{1}x_{2}$}

We obtain%
\begin{equation*}
\sum_{\substack{ a,d=0  \\ a+d\equiv 1}}^{1}\left[
\begin{array}{c}
\left( -1\right) ^{a+1}B(g\otimes x_{2};G^{a}X_{2},g^{d}x_{1}x_{2}) \\
-B(x_{1}\otimes x_{2};G^{a}X_{1}X_{2},g^{d}x_{1}x_{2})%
\end{array}%
\right] G^{a}X_{2}\otimes g^{d}x_{1}x_{2}=0
\end{equation*}%
and we get%
\begin{equation*}
-B(g\otimes x_{2};X_{2},gx_{1}x_{2})-B(x_{1}\otimes
x_{2};X_{1}X_{2},gx_{1}x_{2})=0
\end{equation*}%
and%
\begin{equation*}
B(g\otimes x_{2};GX_{2},x_{1}x_{2})-B(x_{1}\otimes
x_{2};GX_{1}X_{2},x_{1}x_{2})=0
\end{equation*}%
which hold in view of the form of the elements.

\subsubsection{$G^{a}X_{1}\otimes g^{d}x_{1}x_{2}$}

We obtain%
\begin{equation*}
\sum_{\substack{ a,d=0  \\ a+d\equiv 1}}^{1}\left( -1\right)
^{a+1}B(g\otimes x_{2};G^{a}X_{1},g^{d}x_{1}x_{2})G^{a}X_{1}\otimes
g^{d}x_{1}x_{2}=0
\end{equation*}%
and we get%
\begin{equation*}
-B(g\otimes x_{2};X_{1},gx_{1}x_{2})=0
\end{equation*}%
\begin{equation*}
B(g\otimes x_{2};GX_{1},x_{1}x_{2})=0
\end{equation*}%
which hold in view of the form of the element.

\subsubsection{$G^{a}X_{1}X_{2}\otimes g^{d}x_{2}$}

We obtain%
\begin{equation*}
\sum_{\substack{ a,d=0  \\ a+d\equiv 1}}^{1}\left[
\begin{array}{c}
\left( -1\right) ^{a}B(g\otimes x_{2};G^{a}X_{1}X_{2},g^{d}x_{2}) \\
+\left( -1\right) ^{a+1}B(x_{1}\otimes x_{2};G^{a}X_{1}X_{2},g^{d}x_{1}x_{2})%
\end{array}%
\right] G^{a}X_{1}X_{2}\otimes g^{d}x_{2}=0
\end{equation*}%
and we get%
\begin{equation*}
B(g\otimes x_{2};X_{1}X_{2},gx_{2})-B(x_{1}\otimes
x_{2};X_{1}X_{2},gx_{1}x_{2})=0
\end{equation*}%
and%
\begin{equation*}
-B(g\otimes x_{2};GX_{1}X_{2},x_{2})+B(x_{1}\otimes
x_{2};GX_{1}X_{2},x_{1}x_{2})=0
\end{equation*}%
which hold in view of the form of the elements.

\subsubsection{$G^{a}X_{1}X_{2}\otimes g^{d}x_{1}$}

We obtain%
\begin{equation*}
\sum_{\substack{ a,d=0  \\ a+d\equiv 1}}^{1}\left( -1\right) ^{a}B(g\otimes
x_{2};G^{a}X_{1}X_{2},g^{d}x_{1})G^{a}X_{1}X_{2}\otimes g^{d}x_{1}=0
\end{equation*}%
We get%
\begin{equation*}
B(g\otimes x_{2};X_{1}X_{2},gx_{1})=0
\end{equation*}%
and%
\begin{equation*}
-B(g\otimes x_{2};GX_{1}X_{2},x_{1})=0
\end{equation*}%
which hold in view of the form of the element.

\subsubsection{$G^{a}X_{1}X_{2}\otimes g^{d}x_{1}x_{2}$}

We obtain%
\begin{equation*}
\sum_{\substack{ a,d=0  \\ a+d\equiv 0}}^{1}\left( -1\right) ^{a}B(g\otimes
x_{2};G^{a}X_{1}X_{2},g^{d}x_{1}x_{2})G^{a}X_{1}X_{2}\otimes
g^{d}x_{1}x_{2}=0
\end{equation*}%
and we get%
\begin{equation*}
B(g\otimes x_{2};X_{1}X_{2},x_{1}x_{2})=0
\end{equation*}%
and%
\begin{equation*}
-B(g\otimes x_{2};GX_{1}X_{2},gx_{1}x_{2})=0
\end{equation*}%
which hold in view of the form of the element.

\subsection{Case $x_{2}$}

Left side we get%
\begin{eqnarray*}
l_{1}+u_{1}+1-w_{1} &=&0\Rightarrow w_{1}=1\text{ and }l_{1}=u_{1}=0 \\
l_{2}+u_{2} &=&1 \\
a+b_{1}+b_{2}+d+e_{1}+e_{2} &\equiv &1
\end{eqnarray*}%
Right side we get%
\begin{equation*}
\omega _{2}=1
\end{equation*}%
Thus we obtain%
\begin{eqnarray*}
&&\sum_{\substack{ a,b_{1},b_{2},d,e_{1},e_{2}=0  \\ %
a+b_{1}+b_{2}+d+e_{1}+e_{2}\equiv 1}}^{1}\sum_{l_{2}=0}^{b_{2}}\sum
_{\substack{ u_{2}=0  \\ l_{2}+u_{2}=1}}^{e_{2}}\left( -1\right) ^{\alpha
\left( 1_{H};0,l_{2},0,u_{2}\right) } \\
&&B(x_{1}\otimes
x_{2};G^{a}X_{1}^{b_{1}}X_{2}^{b_{2}},g^{d}x_{1}^{e_{1}}x_{2}^{e_{2}})G^{a}X_{1}^{b_{1}}X_{2}^{b_{2}-l_{2}}\otimes g^{d}x_{1}^{e_{1}}x_{2}^{e_{2}-u_{2}}\otimes x_{2}
\\
&=&\sum_{\omega _{2}=0}^{1}B^{A}(x_{1}\otimes 1_{H})\otimes
B^{H}(x_{1}\otimes 1_{H})\otimes x_{2}
\end{eqnarray*}%
i.e.

as%
\begin{eqnarray*}
&&\alpha \left( 1_{H};0,0,0,1\right) \equiv a+b_{1}+b_{2} \\
&&\alpha \left( 1_{H};0,1,0,0\right) \equiv 0
\end{eqnarray*}%
we obtain%
\begin{eqnarray*}
&&\sum_{\substack{ a,b_{1},b_{2},d,e_{1}=0  \\ a+b_{1}+b_{2}+d+e_{1}\equiv 0
}}^{1}\left( -1\right) ^{a+b_{1}+b_{2}}B(x_{1}\otimes
x_{2};G^{a}X_{1}^{b_{1}}X_{2}^{b_{2}},g^{d}x_{1}^{e_{1}}x_{2}) \\
&&G^{a}X_{1}^{b_{1}}X_{2}^{b_{2}}\otimes g^{d}x_{1}^{e_{1}}\otimes x_{2}+ \\
&&+\sum_{\substack{ a,b_{1},d,e_{1},e_{2}=0  \\ a+b_{1}+d+e_{1}+e_{2}\equiv
0 }}^{1}B(x_{1}\otimes
x_{2};G^{a}X_{1}^{b_{1}}X_{2},g^{d}x_{1}^{e_{1}}x_{2}^{e_{2}})G^{a}X_{1}^{b_{1}}\otimes g^{d}x_{1}^{e_{1}}x_{2}^{e_{2}}\otimes x_{2}
\\
&=&B^{A}(x_{1}\otimes 1_{H})\otimes B^{H}(x_{1}\otimes 1_{H})\otimes x_{2}.
\end{eqnarray*}

\subsubsection{$G^{a}\otimes g^{d}$}

We obtain%
\begin{eqnarray*}
&&\sum_{\substack{ a,d=0  \\ a+d\equiv 0}}^{1}\left[ \left( -1\right)
^{a}B(x_{1}\otimes x_{2};G^{a},g^{d}x_{2})+B(x_{1}\otimes
x_{2};G^{a}X_{2},g^{d})\right] G^{a}\otimes g^{d} \\
&=&B(x_{1}\otimes 1_{H};G^{a},g^{d})G^{a}\otimes g^{d}
\end{eqnarray*}%
We get%
\begin{equation*}
B(x_{1}\otimes x_{2};1_{A},x_{2})+B(x_{1}\otimes
x_{2};X_{2},1_{H})-B(x_{1}\otimes 1_{H};1_{A},1_{H})=0
\end{equation*}%
and%
\begin{equation*}
-B(x_{1}\otimes x_{2};G,gx_{2})+B(x_{1}\otimes
x_{2};GX_{2},g)-B(x_{1}\otimes 1_{H};G,g)=0
\end{equation*}%
which hold in view of the form of the element.

\subsubsection{$G^{a}\otimes g^{d}x_{2}$}

\begin{eqnarray*}
&&\sum_{\substack{ a,d=0  \\ a+d\equiv 1}}^{1}B(x_{1}\otimes
x_{2};G^{a}X_{2},g^{d}x_{2})G^{a}\otimes g^{d}x_{2} \\
&=&B(x_{1}\otimes 1_{H};G^{a},g^{d}x_{2})G^{a}\otimes g^{d}x_{2}
\end{eqnarray*}%
We get%
\begin{equation*}
B(x_{1}\otimes x_{2};X_{2},gx_{2})-B(x_{1}\otimes 1_{H};1_{A},gx_{2})=0
\end{equation*}%
and%
\begin{equation*}
B(x_{1}\otimes x_{2};GX_{2},x_{2})-B(x_{1}\otimes 1_{H};G,x_{2})=0
\end{equation*}%
which hold in view of the form of the element.

\subsubsection{$G^{a}\otimes g^{d}x_{1}$}

We obtain%
\begin{eqnarray*}
&&\sum_{\substack{ a,d=0  \\ a+d\equiv 1}}^{1}\left[ \left( -1\right)
^{a}B(x_{1}\otimes x_{2};G^{a},g^{d}x_{1}x_{2})+B(x_{1}\otimes
x_{2};G^{a}X_{2},g^{d}x_{1})\right] G^{a}\otimes g^{d}x_{1} \\
&=&B(x_{1}\otimes 1_{H};G^{a},g^{d}x_{1})G^{a}\otimes g^{d}x_{1}
\end{eqnarray*}%
We get%
\begin{equation*}
B(x_{1}\otimes x_{2};1_{A},gx_{1}x_{2})+B(x_{1}\otimes
x_{2};X_{2},gx_{1})-B(x_{1}\otimes 1_{H};1_{A},gx_{1})=0
\end{equation*}%
and%
\begin{equation*}
-B(x_{1}\otimes x_{2};G,x_{1}x_{2})+B(x_{1}\otimes
x_{2};GX_{2},x_{1})-B(x_{1}\otimes 1_{H};G,x_{1})=0
\end{equation*}%
which hold in view of the form of the element.

\subsubsection{$G^{a}X_{2}\otimes g^{d}$}

We obtain%
\begin{eqnarray*}
&&\sum_{\substack{ a,d=0  \\ a+d\equiv 1}}^{1}\left( -1\right)
^{a+1}B(x_{1}\otimes x_{2};G^{a}X_{2},g^{d}x_{2})G^{a}X_{2}\otimes g^{d} \\
&=&B(x_{1}\otimes 1_{H};G^{a}X_{2},g^{d})G^{a}X_{2}\otimes g^{d}
\end{eqnarray*}%
and we get%
\begin{equation*}
B(x_{1}\otimes x_{2};X_{2},gx_{2})+B(x_{1}\otimes 1_{H};X_{2},g)=0
\end{equation*}%
and%
\begin{equation*}
B(x_{1}\otimes x_{2};GX_{2},x_{2})-B(x_{1}\otimes 1_{H};GX_{2},1_{H})=0
\end{equation*}%
which hold in view of the form of the element.

\subsubsection{$G^{a}X_{1}\otimes g^{d}$}

\begin{eqnarray*}
&&\sum_{\substack{ a,d=0  \\ a+d\equiv 1}}^{1}\left[ \left( -1\right)
^{a+1}B(x_{1}\otimes x_{2};G^{a}X_{1},g^{d}x_{2})+B(x_{1}\otimes
x_{2};G^{a}X_{1}X_{2},g^{d})\right] G^{a}X_{1}\otimes g^{d}+ \\
&=&B(x_{1}\otimes 1_{H};G^{a}X_{1},g^{d})G^{a}X_{1}\otimes g^{d}
\end{eqnarray*}%
We get%
\begin{equation*}
-B(x_{1}\otimes x_{2};X_{1},gx_{2})+B(x_{1}\otimes
x_{2};X_{1}X_{2},g)-B(x_{1}\otimes 1_{H};X_{1},g)=0.
\end{equation*}%
and%
\begin{equation*}
B(x_{1}\otimes x_{2};GX_{1},x_{2})+B(x_{1}\otimes
x_{2};GX_{1}X_{2},1_{H})-B(x_{1}\otimes 1_{H};GX_{1},1_{H})=0
\end{equation*}%
which hold in view of the form of the element.

\subsubsection{$G^{a}\otimes g^{d}x_{1}x_{2}$}

We obtain

\begin{equation*}
\sum_{\substack{ a,d=0  \\ a+d\equiv 0}}^{1}\left[ B(x_{1}\otimes
x_{2};G^{a}X_{2},g^{d}x_{1}x_{2})-B(x_{1}\otimes 1_{H};G^{a},g^{d}x_{1}x_{2})%
\right] G^{a}\otimes g^{d}x_{1}x_{2}=0.
\end{equation*}%
and we get%
\begin{equation*}
B(x_{1}\otimes x_{2};X_{2},x_{1}x_{2})-B(x_{1}\otimes
1_{H};1_{A},x_{1}x_{2})=0
\end{equation*}%
and%
\begin{equation*}
B(x_{1}\otimes x_{2};GX_{2},gx_{1}x_{2})-B(x_{1}\otimes
1_{H};G,gx_{1}x_{2})=0
\end{equation*}%
which hold in view of the form of the element.

\subsubsection{$G^{a}X_{2}\otimes g^{d}x_{2}$}

Nothing on the left side. Thus we get%
\begin{equation*}
\sum_{\substack{ a,d=0  \\ a+d\equiv 0}}^{1}B(x_{1}\otimes
1_{H};G^{a}X_{2},g^{d}x_{2})G^{a}X_{2}\otimes g^{d}x_{2}=0
\end{equation*}%
so that we get%
\begin{eqnarray*}
B(x_{1}\otimes 1_{H};X_{2},x_{2}) &=&0 \\
B(x_{1}\otimes 1_{H};GX_{2},gx_{2}) &=&0
\end{eqnarray*}%
which hold in view of the form of the elements.

\subsubsection{$G^{a}X_{1}\otimes g^{d}x_{2}$}

We obtain%
\begin{eqnarray*}
&&\sum_{\substack{ a,d=0  \\ a+d\equiv 0}}^{1}B(x_{1}\otimes
x_{2};G^{a}X_{1}X_{2},g^{d}x_{2})G^{a}X_{1}\otimes g^{d}x_{2} \\
&=&\sum_{\substack{ a,d=0  \\ a+d\equiv 0}}^{1}B(x_{1}\otimes
1_{H};G^{a}X_{1},g^{d}x_{2})G^{a}X_{1}\otimes g^{d}x_{2}
\end{eqnarray*}%
and we get%
\begin{equation*}
B(x_{1}\otimes x_{2};X_{1}X_{2},x_{2})-B(x_{1}\otimes 1_{H};X_{1},x_{2})=0.
\end{equation*}%
and%
\begin{equation*}
B(x_{1}\otimes x_{2};GX_{1}X_{2},gx_{2})-B(x_{1}\otimes
1_{H};GX_{1},gx_{2})=0
\end{equation*}%
which hold in view of the form of the elements.

\subsubsection{$G^{a}X_{2}\otimes g^{d}x_{1}$}

We obtain%
\begin{equation*}
\sum_{\substack{ a,d=0  \\ a+d\equiv 0}}^{1}\left[ \left( -1\right)
^{a+1}B(x_{1}\otimes x_{2};G^{a}X_{2},g^{d}x_{1}x_{2})-B(x_{1}\otimes
1_{H};G^{a}X_{2},g^{d}x_{1})\right] G^{a}X_{2}\otimes g^{d}x_{1}=0
\end{equation*}%
and we get%
\begin{equation*}
-B(x_{1}\otimes x_{2};X_{2},x_{1}x_{2})-B(x_{1}\otimes 1_{H};X_{2},x_{1})=0
\end{equation*}%
and%
\begin{equation*}
B(x_{1}\otimes x_{2};GX_{2},gx_{1}x_{2})-B(x_{1}\otimes
1_{H};GX_{2},gx_{1})=0
\end{equation*}%
which hold in view of the form of the elements.

\subsubsection{$G^{a}X_{1}\otimes g^{d}x_{1}$}

We obtain%
\begin{eqnarray*}
&&\sum_{\substack{ a,d=0  \\ a+d\equiv 0}}^{1}\left[ \left( -1\right)
^{a+1}B(x_{1}\otimes x_{2};G^{a}X_{1},g^{d}x_{1}x_{2})+B(x_{1}\otimes
x_{2};G^{a}X_{1}X_{2},g^{d}x_{1})\right] G^{a}X_{1}\otimes g^{d}x_{1}+ \\
&=&\sum_{\substack{ a,d=0  \\ a+d\equiv 0}}^{1}B(x_{1}\otimes
1_{H};G^{a}X_{1},g^{d}x_{1})G^{a}X_{1}\otimes g^{d}x_{1}
\end{eqnarray*}%
and we get%
\begin{equation*}
-B(x_{1}\otimes x_{2};X_{1},x_{1}x_{2})+B(x_{1}\otimes
x_{2};X_{1}X_{2},x_{1})-B(x_{1}\otimes 1_{H};X_{1},x_{1})=0
\end{equation*}%
and%
\begin{equation*}
B(x_{1}\otimes x_{2};GX_{1},gx_{1}x_{2})+B(x_{1}\otimes
x_{2};GX_{1}X_{2},gx_{1})-B(x_{1}\otimes 1_{H};GX_{1},gx_{1})=0
\end{equation*}%
which hold in view of the form of the elements.

\subsubsection{$G^{a}X_{1}X_{2}\otimes g^{d}$}

We obtain%
\begin{equation*}
\sum_{\substack{ a,d=0  \\ a+d\equiv 0}}^{1}\left[ \left( -1\right)
^{a}B(x_{1}\otimes x_{2};G^{a}X_{1}X_{2},g^{d}x_{2})-B(x_{1}\otimes
1_{H};G^{a}X_{1}X_{2},g^{d})\right] G^{a}X_{1}X_{2}\otimes g^{d}=0
\end{equation*}%
and we get%
\begin{equation*}
B(x_{1}\otimes x_{2};X_{1}X_{2},x_{2})-B(x_{1}\otimes
1_{H};X_{1}X_{2},1_{H})=0
\end{equation*}%
and%
\begin{equation*}
-B(x_{1}\otimes x_{2};GX_{1}X_{2},gx_{2})-B(x_{1}\otimes
1_{H};GX_{1}X_{2},g)=0
\end{equation*}%
which hold in view of the form of the elements.

\subsubsection{$G^{a}X_{2}\otimes g^{d}x_{1}x_{2}$}

We obtain%
\begin{equation*}
\sum_{\substack{ a,d=0  \\ a+d\equiv 1}}^{1}B(x_{1}\otimes
1_{H};G^{a}X_{2}\otimes g^{d}x_{1}x_{2})G^{a}X_{2}\otimes g^{d}x_{1}x_{2}=0
\end{equation*}%
and we get%
\begin{equation*}
B(x_{1}\otimes 1_{H};X_{2},gx_{1}x_{2})=0
\end{equation*}%
\begin{equation*}
B(x_{1}\otimes 1_{H};GX_{2}\otimes gx_{1}x_{2})=0
\end{equation*}%
which hold in view of the form of the element.

\subsubsection{$G^{a}X_{1}\otimes g^{d}x_{1}x_{2}$}

We obtain%
\begin{equation*}
\sum_{\substack{ a,d=0  \\ a+d\equiv 1}}^{1}\left[ B(x_{1}\otimes
x_{2};G^{a}X_{1}X_{2},g^{d}x_{1}x_{2})-B(x_{1}\otimes
1_{H};G^{a}X_{1},g^{d}x_{1}x_{2})\right] G^{a}X_{1}\otimes g^{d}x_{1}x_{2}=0
\end{equation*}%
and we get%
\begin{equation*}
B(x_{1}\otimes x_{2};X_{1}X_{2},gx_{1}x_{2})-B(x_{1}\otimes
1_{H};X_{1},gx_{1}x_{2})=0
\end{equation*}%
and%
\begin{equation*}
B(x_{1}\otimes x_{2};GX_{1}X_{2},x_{1}x_{2})-B(x_{1}\otimes
1_{H};GX_{1},x_{1}x_{2})=0
\end{equation*}%
which hold in view of the form of the element.

\subsubsection{$G^{a}X_{1}X_{2}\otimes g^{d}x_{2}$}

We obtain%
\begin{equation*}
\sum_{\substack{ a,d=0  \\ a+d\equiv 1}}^{1}B(x_{1}\otimes
1_{H};G^{a}X_{1}X_{2},g^{d}x_{2})G^{a}X_{1}X_{2}\otimes g^{d}x_{2}=0
\end{equation*}%
and we get%
\begin{equation*}
B(x_{1}\otimes 1_{H};X_{1}X_{2},x_{2})=0
\end{equation*}%
and%
\begin{equation*}
B(x_{1}\otimes 1_{H};GX_{1}X_{2},x_{2})=0
\end{equation*}%
which hold in view of the form of the element.

\subsubsection{$G^{a}X_{1}X_{2}\otimes g^{d}x_{1}$}

We obtain%
\begin{equation*}
\sum_{\substack{ a,d=0  \\ a+d\equiv 1}}^{1}\left[ \left( -1\right)
^{a}B(x_{1}\otimes x_{2};G^{a}X_{1}X_{2},g^{d}x_{1}x_{2})-B(x_{1}\otimes
1_{H};G^{a}X_{1}X_{2},g^{d}x_{1})\right] G^{a}X_{1}X_{2}\otimes g^{d}x_{1}=0
\end{equation*}%
and we get%
\begin{equation*}
+B(x_{1}\otimes x_{2};X_{1}X_{2},gx_{1}x_{2})-B(x_{1}\otimes
1_{H};X_{1}X_{2},gx_{1})=0
\end{equation*}%
and%
\begin{equation*}
-B(x_{1}\otimes x_{2};GX_{1}X_{2},x_{1}x_{2})-B(x_{1}\otimes
1_{H};GX_{1}X_{2},x_{1})=0
\end{equation*}%
which hold in view of the form of the elements.

\subsubsection{$G^{a}X_{1}X_{2}\otimes g^{d}x_{1}x_{2}$}

We obtain%
\begin{equation*}
\sum_{\substack{ a,d=0  \\ a+d\equiv 0}}^{1}B(x_{1}\otimes
1_{H};G^{a}X_{1}X_{2},g^{d}x_{1}x_{2})G^{a}X_{1}X_{2}\otimes
g^{d}x_{1}x_{2}=0
\end{equation*}%
and we get%
\begin{equation*}
B(x_{1}\otimes 1_{H};X_{1}X_{2},x_{1}x_{2})=0
\end{equation*}

and%
\begin{equation*}
B(x_{1}\otimes 1_{H};GX_{1}X_{2},gx_{1}x_{2})=0
\end{equation*}%
which hold in view of the form of the element.

\subsection{Case $x_{1}x_{2}$}

We get from the left side%
\begin{eqnarray*}
&&\sum_{w_{1}=0}^{1}\sum_{a,b_{1},b_{2},d,e_{1},e_{2}=0}^{1}%
\sum_{l_{1}=0}^{b_{1}}\sum_{l_{2}=0}^{b_{2}}\sum_{u_{1}=0}^{e_{1}}%
\sum_{u_{2}=0}^{e_{2}}\left( -1\right) ^{\alpha \left(
x_{1}^{1-w_{1}};l_{1},l_{2},u_{1},u_{2}\right) } \\
&&B(g^{1+w_{1}}x_{1}^{w_{1}}\otimes
x_{2};G^{a}X_{1}^{b_{1}}X_{2}^{b_{2}},g^{d}x_{1}^{e_{1}}x_{2}^{e_{2}})G^{a}X_{1}^{b_{1}-l_{1}}X_{2}^{b_{2}-l_{2}}\otimes g^{d}x_{1}^{e_{1}-u_{1}}x_{2}^{e_{2}-u_{2}}\otimes
\\
&&g^{a+b_{1}+b_{2}+l_{1}+l_{2}+d+e_{1}+e_{2}+u_{1}+u_{2}}x_{1}^{l_{1}+u_{1}+1-w_{1}}x_{2}^{l_{2}+u_{2}}
\\
&=&\sum_{\omega _{2}=0}^{1}B^{A}(x_{1}\otimes x_{2}^{1-\omega _{2}})\otimes
B^{H}(x_{1}\otimes x_{2}^{1-\omega _{2}})\otimes g^{1+\omega
_{2}}x_{2}^{\omega _{2}}
\end{eqnarray*}%
\begin{eqnarray*}
l_{1}+u_{1}+1-w_{1} &=&1\text{ i.e.}l_{1}+u_{1}=w_{1} \\
l_{2}+u_{2} &=&1 \\
a+b_{1}+b_{2}+d+e_{1}+e_{2} &\equiv &1+w_{1}
\end{eqnarray*}%
Right side is zero. So that we get%
\begin{gather*}
\sum_{\substack{ a,b_{1},b_{2},d,e_{1},e_{2}=0  \\ %
a+b_{1}+b_{2}+d+e_{1}+e_{2}\equiv 1}}^{1}\sum_{l_{2}=0}^{b_{2}}\sum
_{\substack{ u_{2}=0  \\ l_{2}+u_{2}=1}}^{e_{2}}\left( -1\right) ^{\alpha
\left( x_{1};0,l_{2},0,u_{2}\right) }B(g\otimes
x_{2};G^{a}X_{1}^{b_{1}}X_{2}^{b_{2}},g^{d}x_{1}^{e_{1}}x_{2}^{e_{2}}) \\
G^{a}X_{1}^{b_{1}}X_{2}^{b_{2}-l_{2}}\otimes
g^{d}x_{1}^{e_{1}}x_{2}^{e_{2}-u_{2}} \\
+\sum_{\substack{ a,b_{1},b_{2},d,e_{1},e_{2}=0  \\ %
a+b_{1}+b_{2}+d+e_{1}+e_{2}\equiv 0}}^{1}\sum_{l_{1}=0}^{b_{1}}%
\sum_{l_{2}=0}^{b_{2}}\sum_{\substack{ u_{1}=0  \\ l_{1}+u_{1}=1}}%
^{e_{1}}\sum _{\substack{ u_{2}=0  \\ l_{2}+u_{2}=1}}^{e_{2}}\left(
-1\right) ^{\alpha \left( 1_{H};l_{1},l_{2},u_{1},u_{2}\right) } \\
B(x_{1}\otimes
x_{2};G^{a}X_{1}^{b_{1}}X_{2}^{b_{2}},g^{d}x_{1}^{e_{1}}x_{2}^{e_{2}})G^{a}X_{1}^{b_{1}-l_{1}}X_{2}^{b_{2}-l_{2}}\otimes g^{d}x_{1}^{e_{1}-u_{1}}x_{2}^{e_{2}-u_{2}}=0
\end{gather*}%
i.e.%
\begin{eqnarray*}
&&\sum_{\substack{ a,b_{1},b_{2},d,e=0  \\ a+b_{1}+b_{2}+d+e_{1}\equiv 0}}%
^{1}\left( -1\right) ^{\alpha \left( x_{1};0,0,0,1\right) }B(g\otimes
x_{2};G^{a}X_{1}^{b_{1}}X_{2}^{b_{2}},g^{d}x_{1}^{e_{1}}x_{2})G^{a}X_{1}^{b_{1}}X_{2}^{b_{2}}\otimes g^{d}x_{1}^{e_{1}}+
\\
&&\sum_{\substack{ a,b_{1},,d,e_{1},e_{2}=0  \\ a+b_{1}+d+e_{1}+e_{2}\equiv
0 }}^{1}\left( -1\right) ^{\alpha \left( x_{1};0,1,0,0\right) }B(g\otimes
x_{2};G^{a}X_{1}^{b_{1}}X_{2},g^{d}x_{1}^{e_{1}}x_{2}^{e_{2}})G^{a}X_{1}^{b_{1}}\otimes g^{d}x_{1}^{e_{1}}x_{2}^{e_{2}}
\\
&&+\sum_{\substack{ a,b_{1},b_{2},d=0  \\ a+b_{1}+b_{2}+d\equiv 0}}%
^{1}\left( -1\right) ^{\alpha \left( 1_{H};0,0,1,1\right) }B(x_{1}\otimes
x_{2};G^{a}X_{1}^{b_{1}}X_{2}^{b_{2}},g^{d}x_{1}x_{2})G^{a}X_{1}^{b_{1}}X_{2}^{b_{2}}\otimes g^{d}
\\
&&+\sum_{\substack{ a,b_{1},d,e_{2}=0  \\ a+b_{1}+d+e_{2}\equiv 0}}%
^{1}\left( -1\right) ^{\alpha \left( 1_{H};0,1,1,0\right) }B(x_{1}\otimes
x_{2};G^{a}X_{1}^{b_{1}}X_{2},g^{d}x_{1}x_{2}^{e_{2}})G^{a}X_{1}^{b_{1}}%
\otimes g^{d}x_{2}^{e_{2}} \\
&&+\sum_{\substack{ a,b_{2},d,e_{1}=0  \\ a+b_{2}+d+e_{1}\equiv 0}}%
^{1}\left( -1\right) ^{\alpha \left( 1_{H};1,0,0,1\right) }B(x_{1}\otimes
x_{2};G^{a}X_{1}X_{2}^{b_{2}},g^{d}x_{1}^{e_{1}}x_{2})G^{a}X_{2}^{b_{2}}%
\otimes g^{d}x_{1}^{e_{1}} \\
&&+\sum_{\substack{ a,d,e_{1},e_{2}=0  \\ a+d+e_{1}+e_{2}\equiv 0}}%
^{1}\left( -1\right) ^{\alpha \left( 1_{H};1,1,0,0\right) }B(x_{1}\otimes
x_{2};G^{a}X_{1}X_{2},g^{d}x_{1}^{e_{1}}x_{2}^{e_{2}})G^{a}\otimes
g^{d}x_{1}^{e_{1}}x_{2}^{e_{2}}
\end{eqnarray*}%
Since
\begin{eqnarray*}
B(g\otimes
x_{2};G^{a}X_{1}^{b_{1}}X_{2}^{b_{2}},g^{d}x_{1}^{e_{1}}x_{2}^{e_{2}}) &=&0%
\text{ for }a+b_{1}+b_{2}+d+e_{1}+e_{2}\equiv 1\text{ and} \\
B(x_{1}\otimes
x_{2};G^{a}X_{1}^{b_{1}}X_{2}^{b_{2}},g^{d}x_{1}^{e_{1}}x_{2}^{e_{2}}) &=&0%
\text{ for }a+b_{1}+b_{2}+d+e_{1}+e_{2}\equiv 0
\end{eqnarray*}%
all summands are equal to zero and the equality holds.

\subsection{Case $gx_{1}$}

From the left side%
\begin{eqnarray*}
l_{1}+u_{1} &=&w_{1} \\
l_{2} &=&u_{2}=0 \\
a+b_{1}+b_{2}+d+e_{1}+e_{2} &\equiv &1+w_{1}
\end{eqnarray*}%
Right side is zero. We get%
\begin{gather*}
\sum_{\substack{ a,b_{1},b_{2},d,e_{1},e_{2}=0  \\ %
a+b_{1}+b_{2}+d+e_{1}+e_{2}\equiv 1}}^{1}\left( -1\right) ^{\alpha \left(
x_{1};0,0,0,0\right) }B(g\otimes
x_{2};G^{a}X_{1}^{b_{1}}X_{2}^{b_{2}},g^{d}x_{1}^{e_{1}}x_{2}^{e_{2}})G^{a}X_{1}^{b_{1}}X_{2}^{b_{2}}\otimes g^{d}x_{1}^{e_{1}}x_{2}^{e_{2}}+
\\
+\sum_{\substack{ a,b_{1},b_{2},d,e_{2}=0  \\ a+b_{1}+b_{2}+d+e_{2}\equiv 1}}%
^{1}\left( -1\right) ^{\alpha \left( 1_{H};0,0,1,0\right) }B(x_{1}\otimes
x_{2};G^{a}X_{1}^{b_{1}}X_{2}^{b_{2}},g^{d}x_{1}x_{2}^{e_{2}})G^{a}X_{1}^{b_{1}}X_{2}^{b_{2}}\otimes g^{d}x_{2}^{e_{2}}
\\
+\sum_{\substack{ a,b_{2},d,e_{1},e_{2}=0  \\ a+b_{2}+d+e_{1}+e_{2}\equiv 1}}%
^{1}\left( -1\right) ^{\alpha \left( 1_{H};1,0,0,0\right) }B(x_{1}\otimes
x_{2};G^{a}X_{1}X_{2}^{b_{2}},g^{d}x_{1}^{e_{1}}x_{2}^{e_{2}})G^{a}X_{2}^{b_{2}}\otimes g^{d}x_{1}^{e_{1}}x_{2}^{e_{2}}=0
\end{gather*}%
Since
\begin{eqnarray*}
B(g\otimes
x_{2};G^{a}X_{1}^{b_{1}}X_{2}^{b_{2}},g^{d}x_{1}^{e_{1}}x_{2}^{e_{2}}) &=&0%
\text{ for }a+b_{1}+b_{2}+d+e_{1}+e_{2}\equiv 1\text{ and} \\
B(x_{1}\otimes
x_{2};G^{a}X_{1}^{b_{1}}X_{2}^{b_{2}},g^{d}x_{1}^{e_{1}}x_{2}^{e_{2}}) &=&0%
\text{ for }a+b_{1}+b_{2}+d+e_{1}+e_{2}\equiv 0
\end{eqnarray*}%
all summands are equal to zero and the equality holds.

\subsection{Case $gx_{2}$}

From the left side we get%
\begin{eqnarray*}
l_{1}+u_{1}+1-w_{1} &=&0\Rightarrow w_{1}=1,l_{1}=u_{1}=0 \\
l_{2}+u_{2} &=&1 \\
a+b_{1}+b_{2}+d+e_{1}+e_{2} &\equiv &0
\end{eqnarray*}%
Right side nothing. Thus we get%
\begin{gather*}
\sum_{\substack{ a,b_{1},b_{2},d,e_{1},e_{2}=0  \\ %
a+b_{1}+b_{2}+d+e_{1}+e_{2}\equiv 0}}^{1}\sum_{l_{2}=0}^{b_{2}}\sum
_{\substack{ u_{2}=0  \\ l_{2}+u_{2}=1}}^{e_{2}}\left( -1\right) ^{\alpha
\left( 1_{H};0,l_{2},0,u_{2}\right) } \\
B(x_{1}\otimes
x_{2};G^{a}X_{1}^{b_{1}}X_{2}^{b_{2}},g^{d}x_{1}^{e_{1}}x_{2}^{e_{2}})G^{a}X_{1}^{b_{1}}X_{2}^{b_{2}-l_{2}}\otimes g^{d}x_{1}^{e_{1}}x_{2}^{e_{2}-u_{2}}=0
\end{gather*}%
i.e.%
\begin{eqnarray*}
&&\sum_{\substack{ a,b_{1},b_{2},d,e_{1}=0  \\ a+b_{1}+b_{2}+d+e_{1}\equiv 1
}}^{1}\left( -1\right) ^{\alpha \left( 1_{H};0,0,0,1\right) }B(x_{1}\otimes
x_{2};G^{a}X_{1}^{b_{1}}X_{2}^{b_{2}},g^{d}x_{1}^{e_{1}}x_{2})G^{a}X_{1}^{b_{1}}X_{2}^{b_{2}}\otimes g^{d}x_{1}^{e_{1}}
\\
+ &&\sum_{\substack{ a,b_{1},d,e_{1},e_{2}=0  \\ a+b_{1}+d+e_{1}+e_{2}\equiv
1 }}^{1}\left( -1\right) ^{\alpha \left( 1_{H};0,1,0,0\right)
}B(x_{1}\otimes
x_{2};G^{a}X_{1}^{b_{1}}X_{2},g^{d}x_{1}^{e_{1}}x_{2}^{e_{2}})G^{a}X_{1}^{b_{1}}\otimes g^{d}x_{1}^{e_{1}}x_{2}^{e_{2}}
\end{eqnarray*}%
Since
\begin{equation*}
B(x_{1}\otimes
x_{2};G^{a}X_{1}^{b_{1}}X_{2}^{b_{2}},g^{d}x_{1}^{e_{1}}x_{2}^{e_{2}})=0%
\text{ for }a+b_{1}+b_{2}+d+e_{1}+e_{2}\equiv 0
\end{equation*}%
all summands are equal to zero and the equality holds.

\subsection{Case $gx_{1}x_{2}$}

Nothing from right side. From the left side we get%
\begin{eqnarray*}
l_{1}+u_{1} &=&w_{1} \\
l_{2}+u_{2} &=&1 \\
a+b_{1}+b_{2}+d+e_{1}+e_{2} &\equiv &w_{1}
\end{eqnarray*}%
and hence%
\begin{eqnarray*}
&&\sum_{w_{1}=0}^{1}\sum_{\substack{ a,b_{1},b_{2},d,e_{1},e_{2}=0  \\ %
a+b_{1}+b_{2}+d+e_{1}+e_{2}\equiv w_{1}}}^{1}\sum_{l_{1}=0}^{b_{1}}%
\sum_{l_{2}=0}^{b_{2}}\sum_{\substack{ u_{1}=0  \\ l_{1}+u_{1}=w_{1}}}%
^{e_{1}}\sum_{\substack{ u_{2}=0  \\ l_{2}+u_{2}=1}}^{e_{2}}\left( -1\right)
^{\alpha \left( x_{1}^{1-w_{1}};l_{1},l_{2},u_{1},u_{2}\right) } \\
&&B(g^{1+w_{1}}x_{1}^{w_{1}}\otimes
x_{2};G^{a}X_{1}^{b_{1}}X_{2}^{b_{2}},g^{d}x_{1}^{e_{1}}x_{2}^{e_{2}})G^{a}X_{1}^{b_{1}-l_{1}}X_{2}^{b_{2}-l_{2}}\otimes g^{d}x_{1}^{e_{1}-u_{1}}x_{2}^{e_{2}-u_{2}}
\end{eqnarray*}%
Thus we get%
\begin{eqnarray*}
&&\sum_{w_{1}=0}^{1}\sum_{\substack{ a,b_{1},b_{2},d,e_{1}=0  \\ %
a+b_{1}+b_{2}+d+e_{1}\equiv w_{1}+1}}^{1}\sum_{l_{1}=0}^{b_{1}}\sum
_{\substack{ u_{1}=0  \\ l_{1}+u_{1}=w_{1}}}^{e_{1}}\left( -1\right)
^{\alpha \left( x_{1}^{1-w_{1}};l_{1},0,u_{1},1\right) } \\
&&B(g^{1+w_{1}}x_{1}^{w_{1}}\otimes
x_{2};G^{a}X_{1}^{b_{1}}X_{2}^{b_{2}},g^{d}x_{1}^{e_{1}}x_{2})G^{a}X_{1}^{b_{1}-l_{1}}X_{2}^{b_{2}}\otimes g^{d}x_{1}^{e_{1}-u_{1}}+
\\
&&\sum_{w_{1}=0}^{1}\sum_{\substack{ a,b_{1},d,e_{1},e_{2}=0  \\ %
a+b_{1}+d+e_{1}+e_{2}\equiv w_{1}+1}}^{1}\sum_{l_{1}=0}^{b_{1}}\sum
_{\substack{ u_{1}=0  \\ l_{1}+u_{1}=w_{1}}}^{e_{1}}\left( -1\right)
^{\alpha \left( x_{1}^{1-w_{1}};l_{1},1,u_{1},0\right) } \\
&&B(g^{1+w_{1}}x_{1}^{w_{1}}\otimes
x_{2};G^{a}X_{1}^{b_{1}}X_{2},g^{d}x_{1}^{e_{1}}x_{2}^{e_{2}})G^{a}X_{1}^{b_{1}-l_{1}}\otimes g^{d}x_{1}^{e_{1}-u_{1}}x_{2}^{e_{2}}
\end{eqnarray*}%
and hence%
\begin{eqnarray*}
&&\sum_{\substack{ a,b_{1},b_{2},d,e_{1}=0  \\ a+b_{1}+b_{2}+d+e_{1}\equiv 1
}}^{1}\left( -1\right) ^{\alpha \left( x_{1};0,0,0,1\right) } \\
&&B(g\otimes
x_{2};G^{a}X_{1}^{b_{1}}X_{2}^{b_{2}},g^{d}x_{1}^{e_{1}}x_{2})G^{a}X_{1}^{b_{1}}X_{2}^{b_{2}}\otimes g^{d}x_{1}^{e_{1}}+
\\
&&\sum_{\substack{ a,b_{1},b_{2},d,e_{1}=0  \\ a+b_{1}+b_{2}+d+e_{1}\equiv 0
}}^{1}\sum_{l_{1}=0}^{b_{1}}\sum_{\substack{ u_{1}=0  \\ l_{1}+u_{1}=1}}%
^{e_{1}}\left( -1\right) ^{\alpha \left( 1_{H};l_{1},0,u_{1},1\right) } \\
&&B(x_{1}\otimes
x_{2};G^{a}X_{1}^{b_{1}}X_{2}^{b_{2}},g^{d}x_{1}^{e_{1}}x_{2})G^{a}X_{1}^{b_{1}-l_{1}}X_{2}^{b_{2}}\otimes g^{d}x_{1}^{e_{1}-u_{1}}+
\\
&&\sum_{\substack{ a,b_{1},d,e_{1},e_{2}=0  \\ a+b_{1}+d+e_{1}+e_{2}\equiv 1
}}^{1}\left( -1\right) ^{\alpha \left( x_{1};0,1,0,0\right) } \\
&&B(g\otimes
x_{2};G^{a}X_{1}^{b_{1}}X_{2},g^{d}x_{1}^{e_{1}}x_{2}^{e_{2}})G^{a}X_{1}^{b_{1}}\otimes g^{d}x_{1}^{e_{1}}x_{2}^{e_{2}}+
\\
&&\sum_{\substack{ a,b_{1},d,e_{1},e_{2}=0  \\ a+b_{1}+d+e_{1}+e_{2}\equiv 0
}}^{1}\sum_{l_{1}=0}^{b_{1}}\sum_{\substack{ u_{1}=0  \\ l_{1}+u_{1}=1}}%
^{e_{1}}\left( -1\right) ^{\alpha \left( 1_{H};l_{1},1,u_{1},0\right) } \\
&&B(x_{1}\otimes
x_{2};G^{a}X_{1}^{b_{1}}X_{2},g^{d}x_{1}^{e_{1}}x_{2}^{e_{2}})G^{a}X_{1}^{b_{1}-l_{1}}\otimes g^{d}x_{1}^{e_{1}-u_{1}}x_{2}^{e_{2}}
\end{eqnarray*}%
Finally we get%
\begin{eqnarray*}
&&\sum_{\substack{ a,b_{1},b_{2},d,e_{1}=0  \\ a+b_{1}+b_{2}+d+e_{1}\equiv 1
}}^{1}\left( -1\right) ^{\alpha \left( x_{1};0,0,0,1\right) }B(g\otimes
x_{2};G^{a}X_{1}^{b_{1}}X_{2}^{b_{2}},g^{d}x_{1}^{e_{1}}x_{2})G^{a}X_{1}^{b_{1}}X_{2}^{b_{2}}\otimes g^{d}x_{1}^{e_{1}}
\\
&&+\sum_{\substack{ a,b_{1},d,e_{1},e_{2}=0  \\ a+b_{1}+d+e_{1}+e_{2}\equiv
1 }}^{1}\left( -1\right) ^{\alpha \left( x_{1};0,1,0,0\right) }B(g\otimes
x_{2};G^{a}X_{1}^{b_{1}}X_{2},g^{d}x_{1}^{e_{1}}x_{2}^{e_{2}})G^{a}X_{1}^{b_{1}}\otimes g^{d}x_{1}^{e_{1}}x_{2}^{e_{2}}+
\\
&&\sum_{\substack{ a,b_{1},b_{2},d=0  \\ a+b_{1}+b_{2}+d\equiv 1}}^{1}\left(
-1\right) ^{\alpha \left( 1_{H};0,0,1,1\right) }B(x_{1}\otimes
x_{2};G^{a}X_{1}^{b_{1}}X_{2}^{b_{2}},g^{d}x_{1}x_{2})G^{a}X_{1}^{b_{1}}X_{2}^{b_{2}}\otimes g^{d}+
\\
&&\sum_{\substack{ a,b_{2},d,e_{1}=0  \\ a+b_{2}+d+e_{1}\equiv 1}}^{1}\left(
-1\right) ^{\alpha \left( 1_{H};1,0,0,1\right) }B(x_{1}\otimes
x_{2};G^{a}X_{1}X_{2}^{b_{2}},g^{d}x_{1}^{e_{1}}x_{2})G^{a}X_{2}^{b_{2}}%
\otimes g^{d}x_{1}^{e_{1}}+ \\
&&\sum_{\substack{ a,b_{1},d,e_{2}=0  \\ a+b_{1}+d+e_{2}\equiv 1}}^{1}\left(
-1\right) ^{\alpha \left( 1_{H};0,1,1,0\right) }B(x_{1}\otimes
x_{2};G^{a}X_{1}^{b_{1}}X_{2},g^{d}x_{1}x_{2}^{e_{2}})G^{a}X_{1}^{b_{1}}%
\otimes g^{d}x_{2}^{e_{2}} \\
&&\sum_{\substack{ a,d,e_{1},e_{2}=0  \\ a+d+e_{1}+e_{2}\equiv 1}}^{1}\left(
-1\right) ^{\alpha \left( 1_{H};1,1,0,0\right) }B(x_{1}\otimes
x_{2};G^{a}X_{1}X_{2},g^{d}x_{1}^{e_{1}}x_{2}^{e_{2}})G^{a}\otimes
g^{d}x_{1}^{e_{1}}x_{2}^{e_{2}}
\end{eqnarray*}%
Since%
\begin{eqnarray*}
\alpha \left( x_{1};0,0,0,1\right) &\equiv &1 \\
\alpha \left( x_{1};0,1,0,0\right) &\equiv &a+b_{1}+b_{2}+1 \\
\alpha \left( 1_{H};0,0,1,1\right) &\equiv &1+e_{2} \\
\alpha \left( 1_{H};1,0,0,1\right) &\equiv &a+b_{1} \\
\alpha \left( 1_{H};0,1,1,0\right) &\equiv &e_{2}+a+b_{1}+b_{2}+1 \\
\alpha \left( 1_{H};1,1,0,0\right) &\equiv &1+b_{2}
\end{eqnarray*}%
we obtain

\begin{eqnarray*}
&&\sum_{\substack{ a,b_{1},b_{2},d,e_{1}=0  \\ a+b_{1}+b_{2}+d+e_{1}\equiv 1
}}^{1}-B(g\otimes
x_{2};G^{a}X_{1}^{b_{1}}X_{2}^{b_{2}},g^{d}x_{1}^{e_{1}}x_{2})G^{a}X_{1}^{b_{1}}X_{2}^{b_{2}}\otimes g^{d}x_{1}^{e_{1}}
\\
&&+\sum_{\substack{ a,b_{1},d,e_{1},e_{2}=0  \\ a+b_{1}+d+e_{1}+e_{2}\equiv
1 }}^{1}\left( -1\right) ^{a+b_{1}}B(g\otimes
x_{2};G^{a}X_{1}^{b_{1}}X_{2},g^{d}x_{1}^{e_{1}}x_{2}^{e_{2}})G^{a}X_{1}^{b_{1}}\otimes g^{d}x_{1}^{e_{1}}x_{2}^{e_{2}}+
\\
&&\sum_{\substack{ a,b_{1},b_{2},d=0  \\ a+b_{1}+b_{2}+d\equiv 1}}%
^{1}B(x_{1}\otimes
x_{2};G^{a}X_{1}^{b_{1}}X_{2}^{b_{2}},g^{d}x_{1}x_{2})G^{a}X_{1}^{b_{1}}X_{2}^{b_{2}}\otimes g^{d}+
\\
&&\sum_{\substack{ a,b_{2},d,e_{1}=0  \\ a+b_{2}+d+e_{1}\equiv 1}}^{1}\left(
-1\right) ^{a+1}B(x_{1}\otimes
x_{2};G^{a}X_{1}X_{2}^{b_{2}},g^{d}x_{1}^{e_{1}}x_{2})G^{a}X_{2}^{b_{2}}%
\otimes g^{d}x_{1}^{e_{1}}+ \\
&&\sum_{\substack{ a,b_{1},d,e_{2}=0  \\ a+b_{1}+d+e_{2}\equiv 1}}^{1}\left(
-1\right) ^{e_{2}+a+b_{1}}B(x_{1}\otimes
x_{2};G^{a}X_{1}^{b_{1}}X_{2},g^{d}x_{1}x_{2}^{e_{2}})G^{a}X_{1}^{b_{1}}%
\otimes g^{d}x_{2}^{e_{2}} \\
&&\sum_{\substack{ a,d,e_{1},e_{2}=0  \\ a+d+e_{1}+e_{2}\equiv 1}}%
^{1}B(x_{1}\otimes
x_{2};G^{a}X_{1}X_{2},g^{d}x_{1}^{e_{1}}x_{2}^{e_{2}})G^{a}\otimes
g^{d}x_{1}^{e_{1}}x_{2}^{e_{2}}
\end{eqnarray*}

\subsubsection{$G^{a}\otimes g^{d}$}

We obtain%
\begin{equation*}
\sum_{\substack{ a,d=0  \\ a+d\equiv 1}}^{1}\left[
\begin{array}{c}
-B(g\otimes x_{2};G^{a},g^{d}x_{2})+\left( -1\right) ^{a}B(g\otimes
x_{2};G^{a}X_{2},g^{d}) \\
+B(x_{1}\otimes x_{2};G^{a},g^{d}x_{1}x_{2})+\left( -1\right)
^{a+1}B(x_{1}\otimes x_{2};G^{a}X_{1},g^{d}x_{2}) \\
+\left( -1\right) ^{a}B(x_{1}\otimes
x_{2};G^{a}X_{2},g^{d}x_{1})+B(x_{1}\otimes x_{2};G^{a}X_{1}X_{2},g^{d})%
\end{array}%
\right] G^{a}\otimes g^{d}=0
\end{equation*}%
and we get%
\begin{gather*}
-B(g\otimes x_{2};1_{A},gx_{2})+B(g\otimes x_{2};X_{2},g)+B(x_{1}\otimes
x_{2};1_{A},gx_{1}x_{2})+ \\
-B(x_{1}\otimes x_{2};X_{1},gx_{2})+B(x_{1}\otimes
x_{2};X_{2},gx_{1})+B(x_{1}\otimes x_{2};X_{1}X_{2},g)=0
\end{gather*}%
and%
\begin{gather*}
-B(g\otimes x_{2};G,x_{2})-B(g\otimes x_{2};GX_{2},1_{H})+B(x_{1}\otimes
x_{2};G,x_{1}x_{2})+ \\
+B(x_{1}\otimes x_{2};GX_{1},x_{2})-B(x_{1}\otimes
x_{2};GX_{2},x_{1})+B(x_{1}\otimes x_{2};GX_{1}X_{2},1_{H})=0
\end{gather*}%
which hold in view of the form of the elements.

\subsubsection{$G^{a}\otimes g^{d}x_{2}$}

We obtain%
\begin{equation*}
\sum_{\substack{ a,d=0  \\ a+d\equiv 0}}^{1}\left[
\begin{array}{c}
\left( -1\right) ^{a}B(g\otimes x_{2};G^{a}X_{2},g^{d}x_{2})+ \\
\left( -1\right) ^{a+1}B(x_{1}\otimes
x_{2};G^{a}X_{2},g^{d}x_{1}x_{2})+B(x_{1}\otimes
x_{2};G^{a}X_{1}X_{2},g^{d}x_{2})%
\end{array}%
\right] G^{a}\otimes g^{d}x_{2}=0
\end{equation*}%
and we get%
\begin{equation*}
B(g\otimes x_{2};X_{2},x_{2})-B(x_{1}\otimes
x_{2};X_{2},x_{1}x_{2})+B(x_{1}\otimes x_{2};X_{1}X_{2},x_{2})=0
\end{equation*}%
and%
\begin{equation*}
-B(g\otimes x_{2};GX_{2},gx_{2})+B(x_{1}\otimes
x_{2};GX_{2},gx_{1}x_{2})+B(x_{1}\otimes x_{2};GX_{1}X_{2},gx_{2})=0
\end{equation*}%
which hold in view of the form of the elements.

\subsubsection{$G^{a}\otimes g^{d}x_{1}$}

We obtain%
\begin{equation*}
\sum_{\substack{ a,d=0  \\ a+d\equiv 0}}^{1}\left[
\begin{array}{c}
-B(g\otimes x_{2};G^{a},g^{d}x_{1}x_{2})+\left( -1\right) ^{a}B(g\otimes
x_{2};G^{a}X_{2},g^{d}x_{1})+ \\
\left( -1\right) ^{a+1}B(x_{1}\otimes
x_{2};G^{a}X_{1},g^{d}x_{1}x_{2})+B(x_{1}\otimes
x_{2};G^{a}X_{1}X_{2},g^{d}x_{1})%
\end{array}%
\right] G^{a}\otimes g^{d}x_{1}=0
\end{equation*}%
and we get%
\begin{equation*}
-B(g\otimes x_{2};1_{A},x_{1}x_{2})+B(g\otimes
x_{2};X_{2},x_{1})-B(x_{1}\otimes x_{2};X_{1},x_{1}x_{2})+B(x_{1}\otimes
x_{2};X_{1}X_{2},x_{1})=0
\end{equation*}%
and%
\begin{gather*}
-B(g\otimes x_{2};G,gx_{1}x_{2})-B(g\otimes
x_{2};GX_{2},gx_{1})+B(x_{1}\otimes x_{2};GX_{1},gx_{1}x_{2}) \\
+B(x_{1}\otimes x_{2};GX_{1}X_{2},gx_{1})=0
\end{gather*}%
which hold in view of the form of the elements.

\subsubsection{$G^{a}X_{2}\otimes g^{d}$}

We obtain%
\begin{equation*}
\sum_{\substack{ a,d=0  \\ a+d\equiv 0}}^{1}\left[
\begin{array}{c}
-B(g\otimes x_{2};G^{a}X_{2},g^{d}x_{2})+B(x_{1}\otimes
x_{2};G^{a}X_{2},g^{d}x_{1}x_{2}) \\
+\left( -1\right) ^{a+1}B(x_{1}\otimes x_{2};G^{a}X_{1}X_{2},g^{d}x_{2})%
\end{array}%
\right] G^{a}X_{2}\otimes g^{d}=0
\end{equation*}%
and we get%
\begin{equation*}
-B(g\otimes x_{2};X_{2},x_{2})+B(x_{1}\otimes
x_{2};X_{2},x_{1}x_{2})-B(x_{1}\otimes x_{2};X_{1}X_{2},x_{2})=0
\end{equation*}%
and%
\begin{equation*}
-B(g\otimes x_{2};GX_{2},gx_{2})+B(x_{1}\otimes
x_{2};GX_{2},gx_{1}x_{2})+B(x_{1}\otimes x_{2};GX_{1}X_{2},gx_{2})=0
\end{equation*}%
Both of them, already got.

\subsubsection{$G^{a}X_{1}\otimes g^{d}$}

We obtain%
\begin{equation*}
\sum_{\substack{ a,d=0  \\ a+d\equiv 0}}^{1}\left[
\begin{array}{c}
-B(g\otimes x_{2};G^{a}X_{1},g^{d}x_{2})+\left( -1\right) ^{a+1}B(g\otimes
x_{2};G^{a}X_{1}X_{2},g^{d}) \\
+B(x_{1}\otimes x_{2};G^{a}X_{1},g^{d}x_{1}x_{2})+\left( -1\right)
^{a+1}B(x_{1}\otimes x_{2};G^{a}X_{1}X_{2},g^{d}x_{1})%
\end{array}%
\right] G^{a}X_{1}\otimes g^{d}=0.
\end{equation*}%
and we get%
\begin{equation*}
-B(g\otimes x_{2};X_{1},x_{2})-B(g\otimes
x_{2};X_{1}X_{2},1_{H})+B(x_{1}\otimes
x_{2};X_{1},x_{1}x_{2})-B(x_{1}\otimes x_{2};X_{1}X_{2},x_{1})=0
\end{equation*}%
and%
\begin{equation*}
-B(g\otimes x_{2};GX_{1},gx_{2})+B(g\otimes
x_{2};GX_{1}X_{2},g)+B(x_{1}\otimes x_{2};GX_{1},gx_{1}x_{2})+B(x_{1}\otimes
x_{2};GX_{1}X_{2},gx_{1})=0
\end{equation*}%
Both of them, already got.

\subsubsection{$G^{a}\otimes g^{d}x_{1}x_{2}$}

We obtain%
\begin{equation*}
\sum_{\substack{ a,d=0  \\ a+d\equiv 1}}^{1}\left[ \left( -1\right)
^{a}B(g\otimes x_{2};G^{a}X_{2},g^{d}x_{1}x_{2})+B(x_{1}\otimes
x_{2};G^{a}X_{1}X_{2},g^{d}x_{1}x_{2})\right] G^{a}\otimes g^{d}x_{1}x_{2}=0
\end{equation*}%
and we get%
\begin{equation*}
B(g\otimes x_{2};X_{2},gx_{1}x_{2})+B(x_{1}\otimes
x_{2};X_{1}X_{2},gx_{1}x_{2})=0
\end{equation*}%
and%
\begin{equation*}
-B(g\otimes x_{2};GX_{2},x_{1}x_{2})+B(x_{1}\otimes
x_{2};GX_{1}X_{2},x_{1}x_{2})=0
\end{equation*}%
which hold in view of the form of the elements.

\subsubsection{$G^{a}X_{2}\otimes g^{d}x_{2}$}

There are no summands like this

\subsubsection{$G^{a}X_{1}\otimes g^{d}x_{2}$}

We obtain%
\begin{equation*}
\sum_{\substack{ a,d=0  \\ a+d\equiv 1}}^{1}\left[ \left( -1\right)
^{a+1}B(g\otimes x_{2};G^{a}X_{1}X_{2},g^{d}x_{2})+\left( -1\right)
^{a}B(x_{1}\otimes x_{2};G^{a}X_{1}X_{2},g^{d}x_{1}x_{2})\right]
G^{a}X_{1}\otimes g^{d}x_{2}=0
\end{equation*}%
and we get%
\begin{equation*}
-B(g\otimes x_{2};X_{1}X_{2},gx_{2})+B(x_{1}\otimes
x_{2};X_{1}X_{2},gx_{1}x_{2})=0
\end{equation*}

and%
\begin{equation*}
B(g\otimes x_{2};GX_{1}X_{2},x_{2})-B(x_{1}\otimes
x_{2};GX_{1}X_{2},x_{1}x_{2})=0
\end{equation*}%
which hold in view of the form of the elements.

\subsubsection{$G^{a}X_{2}\otimes g^{d}x_{1}$}

We obtain%
\begin{equation*}
\sum_{\substack{ a,d=0  \\ a+d\equiv 1}}^{1}\left[ -B(g\otimes
x_{2};G^{a}X_{2},g^{d}x_{1}x_{2})+\left( -1\right) ^{a+1}B(x_{1}\otimes
x_{2};G^{a}X_{1}X_{2},g^{d}x_{1}x_{2})\right] G^{a}X_{2}\otimes g^{d}x_{1}=0
\end{equation*}%
and we get%
\begin{equation*}
-B(g\otimes x_{2};X_{2},gx_{1}x_{2})-B(x_{1}\otimes
x_{2};X_{1}X_{2},gx_{1}x_{2})=0
\end{equation*}%
and%
\begin{equation*}
-B(g\otimes x_{2};GX_{2},x_{1}x_{2})+B(x_{1}\otimes
x_{2};GX_{1}X_{2},x_{1}x_{2})=0
\end{equation*}%
which hold in view of the form of the elements.

\subsubsection{$G^{a}X_{1}\otimes g^{d}x_{1}$}

We obtain%
\begin{equation*}
\sum_{\substack{ a,d=0  \\ a+d\equiv 1}}^{1}\left[ -B(g\otimes
x_{2};G^{a}X_{1},g^{d}x_{1}x_{2})+\left( -1\right) ^{a+1}B(g\otimes
x_{2};G^{a}X_{1}X_{2},g^{d}x_{1})\right] G^{a}X_{1}\otimes g^{d}x_{1}=0
\end{equation*}%
and we get%
\begin{equation*}
-B(g\otimes x_{2};X_{1},gx_{1}x_{2})-B(g\otimes x_{2};X_{1}X_{2},gx_{1})=0
\end{equation*}%
and%
\begin{equation*}
-B(g\otimes x_{2};GX_{1},x_{1}x_{2})+B(g\otimes x_{2};GX_{1}X_{2},x_{1})=0
\end{equation*}%
which hold in view of the form of the element.

\subsubsection{$G^{a}X_{1}X_{2}\otimes g^{d}$}

We obtain%
\begin{equation*}
\sum_{\substack{ a,d=0  \\ a+d\equiv 1}}^{1}\left[ -B(g\otimes
x_{2};G^{a}X_{1}X_{2},g^{d}x_{2})+B(x_{1}\otimes
x_{2};G^{a}X_{1}X_{2},g^{d}x_{1}x_{2})\right] G^{a}X_{1}X_{2}\otimes g^{d}=0
\end{equation*}%
and we get%
\begin{equation*}
-B(g\otimes x_{2};X_{1}X_{2},gx_{2})+B(x_{1}\otimes
x_{2};X_{1}X_{2},gx_{1}x_{2})=0
\end{equation*}%
and%
\begin{equation*}
-B(g\otimes x_{2};GX_{1}X_{2},x_{2})+B(x_{1}\otimes
x_{2};GX_{1}X_{2},x_{1}x_{2})=0
\end{equation*}%
which we already got.

\subsubsection{$G^{a}X_{2}\otimes g^{d}x_{1}x_{2}$}

There are no summands like this.

\subsubsection{$G^{a}X_{1}\otimes g^{d}x_{1}x_{2}$}

We obtain%
\begin{equation*}
\sum_{\substack{ a,d=0  \\ a+d\equiv 0}}^{1}\left( -1\right)
^{a+1}B(g\otimes x_{2};G^{a}X_{1}X_{2},g^{d}x_{1}x_{2})G^{a}X_{1}\otimes
g^{d}x_{1}x_{2}=0
\end{equation*}%
and we get%
\begin{equation*}
-B(g\otimes x_{2};X_{1}X_{2},x_{1}x_{2})=0
\end{equation*}%
and%
\begin{equation*}
B(g\otimes x_{2};GX_{1}X_{2},gx_{1}x_{2})=0
\end{equation*}%
which hold in view of the form of the element.

\subsubsection{$G^{a}X_{1}X_{2}\otimes g^{d}x_{2}$}

No summands

\subsubsection{$G^{a}X_{1}X_{2}\otimes g^{d}x_{1}$}

We obtain%
\begin{equation*}
\sum_{\substack{ a,d=0  \\ a+d\equiv 0}}^{1}-B(g\otimes
x_{2};G^{a}X_{1}X_{2},g^{d}x_{1}x_{2})G^{a}X_{1}X_{2}\otimes g^{d}x_{1}=0
\end{equation*}%
and we get%
\begin{equation*}
-B(g\otimes x_{2};X_{1}X_{2},x_{1}x_{2})=0
\end{equation*}%
\begin{equation*}
B(g\otimes x_{2};GX_{1}X_{2},gx_{1}x_{2})=0
\end{equation*}%
which appear just above.

\subsubsection{$G^{a}X_{1}X_{2}\otimes g^{d}x_{1}x_{2}$}

There are no summands like this.

\section{$B\left( 1\otimes gx_{1}x_{2}\right) $}

By using $\left( \ref{eq:hh'}\right) ,$ for $i=1$ we get%
\begin{equation*}
B(h\otimes h^{\prime })(1_{A}\otimes x_{1})=(1_{A}\otimes x_{1})B(gh\otimes
gh^{\prime })+B(x_{1}h\otimes gh^{\prime })+B(h\otimes x_{1}h^{\prime })
\end{equation*}%
Then for $h=1_{H},h^{\prime }=gx_{2}$ and by applying $\left( \ref{eq.10}%
\right) $, we obtain%
\begin{equation}
B(1_{H}\otimes gx_{1}x_{2})=(1_{A}\otimes gx_{1})B(1_{H}\otimes
gx_{2})\left( 1_{A}\otimes g\right) +B(x_{1}\otimes x_{2})-B(1_{H}\otimes
gx_{2})(1_{A}\otimes x_{1})  \label{form 1otgx1x2}
\end{equation}%
Since%
\begin{eqnarray*}
&&B(x_{1}\otimes x_{2})\overset{\left( \ref{form x1otx2}\right) }{=}%
B(x_{1}\otimes 1_{H})(1_{A}\otimes x_{2})-(1_{A}\otimes
gx_{2})B(x_{1}\otimes 1_{H})(1_{A}\otimes g) \\
&&+(1_{A}\otimes g)B(gx_{1}x_{2}\otimes 1_{H})(1_{A}\otimes g)
\end{eqnarray*}%
we finally get%
\begin{eqnarray*}
&&B(1_{H}\otimes gx_{1}x_{2})\overset{\left( \ref{form 1otgx1x2}\right) }{=}%
(1_{A}\otimes gx_{1})B(1_{H}\otimes gx_{2})\left( 1_{A}\otimes g\right) \\
&&+B(x_{1}\otimes x_{2})-B(1_{H}\otimes gx_{2})(1_{A}\otimes x_{1}) \\
&=&+(1_{A}\otimes x_{1})B(g\otimes 1_{H})(1_{A}\otimes x_{2}) \\
&&+(1_{A}\otimes gx_{1}x_{2})B(g\otimes 1_{H})\left( 1_{A}\otimes g\right) \\
&&+(1_{A}\otimes gx_{1})B(x_{2}\otimes 1_{H})\left( 1_{A}\otimes g\right) \\
&&+B(x_{1}\otimes 1_{H})(1_{A}\otimes x_{2}) \\
&&-(1_{A}\otimes gx_{2})B(x_{1}\otimes 1_{H})(1_{A}\otimes g) \\
&&+(1_{A}\otimes g)B(gx_{1}x_{2}\otimes 1_{H})(1_{A}\otimes g) \\
&&+(1_{A}\otimes g)B(g\otimes 1_{H})(1_{A}\otimes gx_{1}x_{2})+ \\
&&-(1_{A}\otimes x_{2})B(g\otimes 1_{H})(1_{A}\otimes x_{1})+ \\
&&-B(x_{2}\otimes 1_{H})(1_{A}\otimes x_{1})
\end{eqnarray*}

By using $\left( \ref{form x1otx2}\right) $ we obtain%
\begin{eqnarray*}
B\left( 1\otimes gx_{1}x_{2}\right) &=&B(gx_{1}x_{2}\otimes
1_{H};1_{A},g)1_{A}\otimes g+ \\
&&-B(gx_{1}x_{2}\otimes 1_{H};1_{A},x_{1})1_{A}\otimes x_{1} \\
&&-B(gx_{1}x_{2}\otimes 1_{H};1_{A},x_{2})1_{A}\otimes x_{2} \\
&&+\left[
\begin{array}{c}
4B\left( g\otimes 1_{H};1_{A},g\right) +2B\left( x_{1}\otimes
1_{H};1_{A},gx_{1}\right) \\
+2B\left( x_{2}\otimes 1_{H};1_{A},gx_{2}\right) +B(gx_{1}x_{2}\otimes
1_{H};1_{A},gx_{1}x_{2})%
\end{array}%
\right] 1_{A}\otimes gx_{1}x_{2}+ \\
&&+B(gx_{1}x_{2}\otimes 1_{H};G,1_{H})G\otimes 1_{H}+B(gx_{1}x_{2}\otimes
1_{H};G,x_{1}x_{2})G\otimes x_{1}x_{2}+ \\
&&\left[ 2B\left( x_{2}\otimes 1_{H};G,g\right) -B(gx_{1}x_{2}\otimes
1_{H};G,gx_{1})\right] G\otimes gx_{1} \\
&&\left[ -2B\left( x_{1}\otimes 1_{H};G,g\right) -B(gx_{1}x_{2}\otimes
1_{H};G,gx_{2})\right] G\otimes gx_{2}+ \\
&&+\left[ -B(x_{2}\otimes 1_{H};1_{A},1_{H})+B(gx_{1}x_{2}\otimes
1_{H};1_{A},x_{1})\right] X_{1}\otimes 1_{H}+ \\
&&-B(x_{2}\otimes 1_{H};1_{A},x_{1}x_{2})X_{1}\otimes x_{1}x_{2}+ \\
&&-B(x_{2}\otimes 1_{H};1_{A},gx_{1})X_{1}\otimes gx_{1} \\
&&+\left[
\begin{array}{c}
2B(g\otimes 1_{H};1_{A},g)+2B(x_{1}\otimes 1_{H};1_{A},gx_{1}) \\
+B(x_{2}\otimes 1_{H};1_{A},gx_{2})+B(gx_{1}x_{2}\otimes
1_{H};1_{A},gx_{1}x_{2})%
\end{array}%
\right] X_{1}\otimes gx_{2}+ \\
&&+\left[ B(x_{1}\otimes 1_{H};1_{A},1_{H})+B(gx_{1}x_{2}\otimes
1_{H};1_{A},x_{2})\right] X_{2}\otimes 1_{H} \\
&&+B(x_{1}\otimes 1_{H};1_{A},x_{1}x_{2})X_{2}\otimes x_{1}x_{2}+ \\
&&+\left[
\begin{array}{c}
-2B(g\otimes 1_{H};1_{A},g)-B(x_{1}\otimes 1_{H};1_{A},gx_{1}) \\
-2B(x_{2}\otimes \ 1_{H};1_{A},gx_{2})-B(gx_{1}x_{2}\otimes
1_{H};1_{A},gx_{1}x_{2})%
\end{array}%
\right] X_{2}\otimes gx_{1}+ \\
&&+B(x_{1}\otimes 1_{H};1_{A},gx_{2})X_{2}\otimes gx_{2}+ \\
&&+\left[
\begin{array}{c}
B(g\otimes 1_{H};1_{A},g)+B(x_{2}\otimes \ 1_{H};1_{A},gx_{2}) \\
+B(x_{1}\otimes 1_{H};1_{A},gx_{1})+B(gx_{1}x_{2}\otimes
1_{H};1_{A},gx_{1}x_{2})%
\end{array}%
\right] X_{1}X_{2}\otimes g+ \\
&&-\left[ B(g\otimes 1_{H};1_{A},x_{1})+B(x_{2}\otimes \
1_{H};1_{A},x_{1}x_{2})\right] X_{1}X_{2}\otimes x_{1}+ \\
&&-\left[ B(g\otimes 1_{H};X_{2},1_{H})-B(x_{1}\otimes
1_{H};1_{A},x_{1}x_{2})\right] X_{1}X_{2}\otimes x_{2}+ \\
&&+B\left( g\otimes 1_{H};1_{A},gx_{1}x_{2}\right) X_{1}X_{2}\otimes
gx_{1}x_{2}+ \\
&&+\left[ B(x_{2}\otimes 1_{H};G,g)-B(gx_{1}x_{2}\otimes 1_{H};G,gx_{1})%
\right] GX_{1}\otimes g+ \\
&&-B(x_{2}\otimes 1_{H};G,x_{1})GX_{1}\otimes x_{1}+ \\
&&-\left[ B(x_{2}\otimes 1_{H};G,x_{2})+B(gx_{1}x_{2}\otimes
1_{H};G,x_{1}x_{2})\right] GX_{1}\otimes x_{2}+ \\
&&\left[ -2B\left( g\otimes 1_{H};G,gx_{1}\right) -B(x_{2}\otimes
1_{H};G,gx_{1}x_{2})\right] GX_{1}\otimes gx_{1}x_{2}+ \\
&&+\left[ -B(x_{1}\otimes 1_{H};G,g)-B(gx_{1}x_{2}\otimes 1_{H};G,gx_{2})%
\right] GX_{2}\otimes g+ \\
&&+\left[ +B(x_{1}\otimes 1_{H};G,x_{1})+B(gx_{1}x_{2}\otimes
1_{H};G,x_{1}x_{2})\right] GX_{2}\otimes x_{1}+ \\
&&+B(x_{1}\otimes 1_{H};G,x_{2})GX_{2}\otimes x_{2}+ \\
&&+\left[ -2B(g\otimes 1_{H};G,gx_{2})+B(x_{1}\otimes 1_{H};G,gx_{1}x_{2})%
\right] GX_{2}\otimes gx_{1}x_{2}+ \\
&&+\left[
\begin{array}{c}
B(g\otimes 1_{H};G,1_{H})+B(x_{2}\otimes \ 1_{H};G,x_{2}) \\
+B(x_{1}\otimes 1_{H};G,x_{1})+B(gx_{1}x_{2}\otimes 1_{H};G,x_{1}x_{2})%
\end{array}%
\right] GX_{1}X_{2}\otimes 1_{H} \\
&&+B(g\otimes 1_{H};G,x_{1}x_{2})GX_{1}X_{2}\otimes x_{1}x_{2}+\left[
\begin{array}{c}
B(g\otimes 1_{H};G,gx_{1}) \\
+B(x_{2}\otimes \ 1_{H};G,gx_{1}x_{2})%
\end{array}%
\right] GX_{1}X_{2}\otimes gx_{1}+ \\
&&+\left[ +B(g\otimes 1_{H};G,gx_{2})-B(x_{1}\otimes 1_{H};G,gx_{1}x_{2})%
\right] GX_{1}X_{2}\otimes gx_{2}
\end{eqnarray*}

We write the Casimir condition for $B\left( 1\otimes gx_{1}x_{2}\right) $
and we get%
\begin{eqnarray*}
&&\sum_{a,b_{1},b_{2},d,e_{1},e_{2}=0}^{1}\sum_{l_{1}=0}^{b_{1}}%
\sum_{l_{2}=0}^{b_{2}}\sum_{u_{1}=0}^{e_{1}}\sum_{u_{2}=0}^{e_{2}}\left(
-1\right) ^{\alpha \left( 1_{H};l_{1},l_{2},u_{1},u_{2}\right) } \\
&&B(1_{H}\otimes
gx_{1}x_{2};G^{a}X_{1}^{b_{1}}X_{2}^{b_{2}},g^{d}x_{1}^{e_{1}}x_{2}^{e_{2}})
\\
&&G^{a}X_{1}^{b_{1}-l_{1}}X_{2}^{b_{2}-l_{2}}\otimes
g^{d}x_{1}^{e_{1}-u_{1}}x_{2}^{e_{2}-u_{2}}\otimes \\
&&g^{a+b_{1}+b_{2}+l_{1}+l_{2}+d+e_{1}+e_{2}+u_{1}+u_{2}}x_{1}^{l_{1}+u_{1}}x_{2}^{l_{2}+u_{2}}
\\
&=&\sum_{\omega _{1}=0}^{1}\sum_{\omega _{2}=0}^{1}\left( -1\right) ^{\left(
1+\omega _{2}\right) \omega _{1}}B^{A}(1_{H}\otimes gx_{1}^{1-\omega
_{1}}x_{2}^{1-\omega _{2}})\otimes B^{H}(1_{H}\otimes gx_{1}^{1-\omega
_{1}}x_{2}^{1-\omega _{2}})\otimes g^{1+\omega _{1}+\omega
_{2}}x_{1}^{\omega _{1}}x_{2}^{\omega _{2}}
\end{eqnarray*}

\subsection{case $1_{H}$}

We note that the right side of the equality is $0$. The left side gives us%
\begin{eqnarray*}
a+b_{1}+b_{2}+l_{1}+l_{2}+d+e_{1}+e_{2}+u_{1}+u_{2} &\equiv &0 \\
l_{1}+u_{1} &=&0 \\
l_{2}+u_{2} &=&0
\end{eqnarray*}%
so that%
\begin{eqnarray*}
a+b_{1}+b_{2}+d+e_{1}+e_{2} &\equiv &0 \\
l_{1} &=&u_{1}=0 \\
l_{2} &=&u_{2}=0
\end{eqnarray*}%
and we get since%
\begin{equation*}
\alpha \left( 1_{H};0,0,0,0\right) =0
\end{equation*}%
\begin{equation*}
\sum_{\substack{ a,b_{1},b_{2},d,e_{1},e_{2}=0  \\ %
a+b_{1}+b_{2}+d+e_{1}+e_{2}\equiv 0}}^{1}B(1_{H}\otimes
gx_{1}x_{2};G^{a}X_{1}^{b_{1}}X_{2}^{b_{2}},g^{d}x_{1}^{e_{1}}x_{2}^{e_{2}})G^{a}X_{1}^{b_{1}}X_{2}^{b_{2}}\otimes g^{d}x_{1}^{e_{1}}x_{2}^{e_{2}}=0
\end{equation*}%
and we deduce that%
\begin{equation*}
B(1_{H}\otimes
gx_{1}x_{2};G^{a}X_{1}^{b_{1}}X_{2}^{b_{2}},g^{d}x_{1}^{e_{1}}x_{2}^{e_{2}})=0%
\text{ whenever }a+b_{1}+b_{2}+d+e_{1}+e_{2}\equiv 0
\end{equation*}

\subsection{case $g$}

The right side is%
\begin{equation*}
B^{A}(1_{H}\otimes gx_{1}x_{2})\otimes B^{H}(1_{H}\otimes
gx_{1}x_{2})\otimes g
\end{equation*}%
The left side gives us%
\begin{eqnarray*}
a+b_{1}+b_{2}+l_{1}+l_{2}+d+e_{1}+e_{2}+u_{1}+u_{2} &\equiv &1 \\
l_{1}+u_{1} &=&0 \\
l_{2}+u_{2} &=&0
\end{eqnarray*}%
so that%
\begin{eqnarray*}
a+b_{1}+b_{2}+d+e_{1}+e_{2} &\equiv &1 \\
l_{1} &=&u_{1}=0 \\
l_{2} &=&u_{2}=0
\end{eqnarray*}%
and we get%
\begin{eqnarray*}
&&\sum_{\substack{ a,b_{1},b_{2},d,e_{1},e_{2}=0  \\ %
a+b_{1}+b_{2}+d+e_{1}+e_{2}\equiv 1}}^{1}B(1_{H}\otimes
gx_{1}x_{2};G^{a}X_{1}^{b_{1}}X_{2}^{b_{2}},g^{d}x_{1}^{e_{1}}x_{2}^{e_{2}})G^{a}X_{1}^{b_{1}}X_{2}^{b_{2}}\otimes g^{d}x_{1}^{e_{1}}x_{2}^{e_{2}}
\\
&&=B^{A}(1_{H}\otimes gx_{1}x_{2})\otimes B^{H}(1_{H}\otimes gx_{1}x_{2})
\end{eqnarray*}%
which holds as

\begin{equation*}
B(1_{H}\otimes
gx_{1}x_{2};G^{a}X_{1}^{b_{1}}X_{2}^{b_{2}},g^{d}x_{1}^{e_{1}}x_{2}^{e_{2}})=0%
\text{ whenever }a+b_{1}+b_{2}+d+e_{1}+e_{2}\equiv 0
\end{equation*}

\subsection{case $x_{1}$}

The right side is%
\begin{equation*}
-B^{A}(1_{H}\otimes gx_{2})\otimes B^{H}(1_{H}\otimes gx_{2})\otimes x_{1}
\end{equation*}%
The left side gives us%
\begin{eqnarray*}
a+b_{1}+b_{2}+l_{1}+l_{2}+d+e_{1}+e_{2}+u_{1}+u_{2} &\equiv &0 \\
l_{1}+u_{1} &=&1 \\
l_{2}+u_{2} &=&0
\end{eqnarray*}%
so that%
\begin{eqnarray*}
a+b_{1}+b_{2}+d+e_{1}+e_{2} &\equiv &1 \\
l_{1}+u_{1} &=&1 \\
l_{2} &=&u_{2}=0
\end{eqnarray*}%
Since%
\begin{eqnarray*}
\alpha \left( 1_{H};0,0,1,0\right) &=&e_{2}+\left( a+b_{1}+b_{2}\right)
\left( 1+n_{1}+n_{2}\right) +m \\
&=&e_{2}+\left( a+b_{1}+b_{2}\right)
\end{eqnarray*}%
and%
\begin{eqnarray*}
\alpha \left( 1_{H};1,0,0,0\right) &=&b_{2}+\left( a+b_{1}+b_{2}+1\right)
\left( n_{1}+n_{2}\right) +n_{2}) \\
&=&b_{2}
\end{eqnarray*}%
the equality becomes

\begin{eqnarray*}
&&\sum_{\substack{ a,b_{1},b_{2},d,e_{2}=0  \\ a+b_{1}+b_{2}+d+e_{2}\equiv 0
}}^{1}\left( -1\right) ^{e_{2}+a+b_{1}+b_{2}}B(1_{H}\otimes
gx_{1}x_{2};G^{a}X_{1}^{b_{1}}X_{2}^{b_{2}},g^{d}x_{1}x_{2}^{e_{2}})G^{a}X_{1}^{b_{1}}X_{2}^{b_{2}}\otimes g^{d}x_{2}^{e_{2}}+
\\
&&\sum_{\substack{ a,b_{2},d,e_{1},e_{2}=0  \\ a+b_{2}+d+e_{1}+e_{2}\equiv 0
}}^{1}\left( -1\right) ^{b_{2}}B(1_{H}\otimes
gx_{1}x_{2};G^{a}X_{1}X_{2}^{b_{2}},g^{d}x_{1}^{e_{1}}x_{2}^{e_{2}})G^{a}X_{2}^{b_{2}}\otimes g^{d}x_{1}^{e_{1}}x_{2}^{e_{2}}
\\
&=&-B^{A}(1_{H}\otimes gx_{2})\otimes B^{H}(1_{H}\otimes gx_{2})
\end{eqnarray*}

\subsubsection{$G^{a}\otimes g^{d}$}

We obtain%
\begin{eqnarray*}
&&\sum_{\substack{ a,d=0  \\ a+d\equiv 0}}^{1}\left[ \left( -1\right)
^{a}B(1_{H}\otimes gx_{1}x_{2};G^{a},g^{d}x_{1})+B(1_{H}\otimes
gx_{1}x_{2};G^{a}X_{1},g^{d})\right] G^{a}\otimes g^{d} \\
&&+B(1_{H}\otimes gx_{2};G^{a},g^{d})G^{a}\otimes g^{d} \\
&=&0
\end{eqnarray*}%
and we get%
\begin{gather*}
B(1_{H}\otimes gx_{1}x_{2};1_{A},x_{1})+B(1_{H}\otimes
gx_{1}x_{2};X_{1},1_{H})+ \\
+B(1_{H}\otimes gx_{2};1_{A},1_{H})=0.
\end{gather*}%
and%
\begin{gather*}
-B(1_{H}\otimes gx_{1}x_{2};G,gx_{1})+B(1_{H}\otimes gx_{1}x_{2};GX_{1},g)+
\\
+B(1_{H}\otimes gx_{2};G,g)=0.
\end{gather*}%
By

\subsubsection{$G^{a}\otimes g^{d}x_{2}$}

We obtain%
\begin{gather*}
\sum_{\substack{ a,d=0  \\ a+d\equiv 1}}^{1}\left( -1\right)
^{a+1}B(1_{H}\otimes gx_{1}x_{2};G^{a},g^{d}x_{1}x_{2})G^{a}\otimes
g^{d}x_{2}+ \\
\sum_{\substack{ a,d=0  \\ a+d\equiv 1}}^{1}B(1_{H}\otimes
gx_{1}x_{2};G^{a}X_{1},g^{d}x_{2})G^{a}\otimes g^{d}x_{2} \\
+B(1_{H}\otimes gx_{2};G^{a},g^{d}x_{2})G^{a}\otimes g^{d}x_{2}=0
\end{gather*}%
and we get%
\begin{gather*}
-B(1_{H}\otimes gx_{1}x_{2};1_{A},gx_{1}x_{2})+B(1_{H}\otimes
gx_{1}x_{2};X_{1},gx_{2})+ \\
+B(1_{H}\otimes gx_{2};1_{A},gx_{2})=0.
\end{gather*}%
and%
\begin{gather*}
B(1_{H}\otimes gx_{1}x_{2};G,x_{1}x_{2})+B(1_{H}\otimes
gx_{1}x_{2};GX_{1},x_{2})+ \\
+B(1_{H}\otimes gx_{2};G,x_{2})=0.
\end{gather*}%
which are satisfied by using the form of the elements.

\subsubsection{$G^{a}\otimes g^{d}x_{1}$}

We obtain%
\begin{equation*}
\sum_{\substack{ a,d=0  \\ a+d\equiv 1}}^{1}B(1_{H}\otimes
gx_{1}x_{2};G^{a}X_{1},g^{d}x_{1})G^{a}\otimes g^{d}x_{1}=0
\end{equation*}%
and we get%
\begin{equation*}
B(1_{H}\otimes gx_{1}x_{2};X_{1},gx_{1})+B(1_{H}\otimes
gx_{2};1_{A},gx_{1})=0
\end{equation*}%
and%
\begin{equation*}
B(1_{H}\otimes gx_{1}x_{2};GX_{1},x_{1})+B(1_{H}\otimes gx_{2};G,x_{1})=0
\end{equation*}%
which are satisfied by using the form of the elements.

\subsubsection{$G^{a}X_{2}\otimes g^{d}$}

The left side is%
\begin{equation*}
\sum_{\substack{ a,d=0  \\ a+d\equiv 1}}^{1}\left[ \left( -1\right)
^{a+1}B(1_{H}\otimes gx_{1}x_{2};G^{a}X_{2},g^{d}x_{1})-B(1_{H}\otimes
gx_{1}x_{2};G^{a}X_{1}X_{2},g^{d})\right] G^{a}X_{2}\otimes g^{d}+
\end{equation*}%
and we get%
\begin{equation*}
-B(1_{H}\otimes gx_{1}x_{2};X_{2},gx_{1})-B(1_{H}\otimes
gx_{1}x_{2};X_{1}X_{2},g)+B(1_{H}\otimes gx_{2};X_{2},g)=0
\end{equation*}%
and%
\begin{equation*}
+B(1_{H}\otimes gx_{1}x_{2};GX_{2},x_{1})-B(1_{H}\otimes
gx_{1}x_{2};GX_{1}X_{2},1_{H})+B(1_{H}\otimes gx_{2};GX_{2},1_{H})=0
\end{equation*}%
which are satisfied in view of the form of the elements.

\subsubsection{$G^{a}X_{1}\otimes g^{d}$}

The left side gives us%
\begin{equation*}
\sum_{\substack{ a,d=0  \\ a+d\equiv 1}}^{1}\left( -1\right)
^{a+1}B(1_{H}\otimes
gx_{1}x_{2};G^{a}X_{1},g^{d}x_{1})G^{a}X_{1}^{b_{1}}X_{2}^{b_{2}}\otimes
g^{d}x_{2}^{e_{2}}
\end{equation*}%
and we get%
\begin{equation*}
-B(1_{H}\otimes gx_{1}x_{2};X_{1},gx_{1})+B(1_{H}\otimes gx_{2};X_{1},g)=0
\end{equation*}%
and%
\begin{equation*}
+B(1_{H}\otimes gx_{1}x_{2};GX_{1},x_{1})+B(1_{H}\otimes
gx_{2};GX_{1},1_{H})=0
\end{equation*}%
which are satisfied in view of the form of the elements.

\subsubsection{$G^{a}\otimes g^{d}x_{1}x_{2}$}

The left side gives us%
\begin{equation*}
\sum_{\substack{ a,d=0  \\ a+d\equiv 0}}^{1}B(1_{H}\otimes
gx_{1}x_{2};G^{a}X_{1},g^{d}x_{1}x_{2})G^{a}\otimes g^{d}x_{1}x_{2}
\end{equation*}%
and we obtain%
\begin{equation*}
B(1_{H}\otimes gx_{1}x_{2};X_{1},x_{1}x_{2})+B(1_{H}\otimes
gx_{2};1_{A},x_{1}x_{2})=0
\end{equation*}%
and%
\begin{equation*}
B(1_{H}\otimes gx_{1}x_{2};GX_{1},gx_{1}x_{2})+B(1_{H}\otimes
gx_{2};G,gx_{1}x_{2})=0
\end{equation*}%
which are satisfied in view of the form of the elements.

\subsubsection{$G^{a}X_{2}\otimes g^{d}x_{2}$}

The left side gives us%
\begin{equation*}
\sum_{\substack{ a,d=0  \\ a+d\equiv 0}}^{1}\left[ \left( -1\right)
^{a}B(1_{H}\otimes gx_{1}x_{2};G^{a}X_{2},g^{d}x_{1}x_{2})-B(1_{H}\otimes
gx_{1}x_{2};G^{a}X_{1}X_{2},g^{d}x_{2})\right] G^{a}X_{2}\otimes g^{d}x_{2}
\end{equation*}%
and we get%
\begin{equation*}
B(1_{H}\otimes gx_{1}x_{2};X_{2},x_{1}x_{2})-B(1_{H}\otimes
gx_{1}x_{2};X_{1}X_{2},x_{2})+B(1_{H}\otimes gx_{2};X_{2},x_{2})=0
\end{equation*}%
and%
\begin{equation*}
-B(1_{H}\otimes gx_{1}x_{2};GX_{2},gx_{1}x_{2})-B(1_{H}\otimes
gx_{1}x_{2};GX_{1}X_{2},gx_{2})+B(1_{H}\otimes gx_{2};GX_{2},gx_{2})=0
\end{equation*}%
which are satisfied in view of the form of the elements.

\subsubsection{$G^{a}X_{1}\otimes g^{d}x_{2}$}

The left side gives us%
\begin{equation*}
\sum_{\substack{ a,d=0  \\ a+d\equiv 0}}^{1}\left( -1\right)
^{a}B(1_{H}\otimes gx_{1}x_{2};G^{a}X_{1},g^{d}x_{1}x_{2})G^{a}X_{1}\otimes
g^{d}x_{2}
\end{equation*}%
and we get%
\begin{equation*}
B(1_{H}\otimes gx_{1}x_{2};X_{1},x_{1}x_{2})+B(1_{H}\otimes
gx_{2};X_{1},x_{2})=0
\end{equation*}%
and%
\begin{equation*}
-B(1_{H}\otimes gx_{1}x_{2};GX_{1},gx_{1}x_{2})+B(1_{H}\otimes
gx_{2};GX_{1},gx_{2})=0
\end{equation*}%
which are satisfied in view of the form of the elements.

\subsubsection{$G^{a}X_{2}\otimes g^{d}x_{1}$}

The left side gives us%
\begin{equation*}
\sum_{\substack{ a,d=0  \\ a+d\equiv 0}}^{1}-B(1_{H}\otimes
gx_{1}x_{2};G^{a}X_{1}X_{2},g^{d}x_{1})G^{a}X_{2}\otimes g^{d}x_{1}
\end{equation*}%
and we get%
\begin{equation*}
-B(1_{H}\otimes gx_{1}x_{2};X_{1}X_{2},x_{1})+B(1_{H}\otimes
gx_{2};X_{2},x_{1})=0
\end{equation*}%
and%
\begin{equation*}
-B(1_{H}\otimes gx_{1}x_{2};GX_{1}X_{2},gx_{1})+B(1_{H}\otimes
gx_{2};GX_{2},gx_{1})=0
\end{equation*}%
which are satisfied in view of the form of the elements.

\subsubsection{$G^{a}X_{1}\otimes g^{d}x_{1}$}

There are no terms like this.

\subsubsection{$G^{a}X_{1}X_{2}\otimes g^{d}$}

The left side gives us%
\begin{equation*}
\sum_{\substack{ a,d=0  \\ a+d\equiv 0}}^{1}\left( -1\right)
^{a}B(1_{H}\otimes
gx_{1}x_{2};G^{a}X_{1}X_{2},g^{d}x_{1})G^{a}X_{1}X_{2}\otimes g^{d}
\end{equation*}%
and we get%
\begin{equation*}
B(1_{H}\otimes gx_{1}x_{2};X_{1}X_{2},x_{1})+B(1_{H}\otimes
gx_{2};X_{1}X_{2},1_{H})=0
\end{equation*}%
and%
\begin{equation*}
-B(1_{H}\otimes gx_{1}x_{2};GX_{1}X_{2},gx_{1})+B(1_{H}\otimes
gx_{2};GX_{1}X_{2},g)=0
\end{equation*}%
which are satisfied in view of the form of the elements.

\subsubsection{$G^{a}X_{2}\otimes g^{d}x_{1}x_{2}$}

The left side gives us%
\begin{equation*}
\sum_{\substack{ a,d=0  \\ a+d\equiv 1}}^{1}-B(1_{H}\otimes
gx_{1}x_{2};G^{a}X_{1}X_{2},g^{d}x_{1}x_{2})G^{a}X_{2}\otimes g^{d}x_{1}x_{2}
\end{equation*}%
and we get%
\begin{equation*}
-B(1_{H}\otimes gx_{1}x_{2};X_{1}X_{2},gx_{1}x_{2})+B(1_{H}\otimes
gx_{2};X_{2},gx_{1}x_{2})=0
\end{equation*}%
and%
\begin{equation*}
-B(1_{H}\otimes gx_{1}x_{2};GX_{1}X_{2},x_{1}x_{2})+B(1_{H}\otimes
gx_{2};GX_{2},x_{1}x_{2})=0
\end{equation*}%
which are satisfied in view of the form of the elements.

\subsubsection{$G^{a}X_{1}\otimes g^{d}x_{1}x_{2}$}

There are no terms like this.

\subsubsection{$G^{a}X_{1}X_{2}\otimes g^{d}x_{2}$}

The left side gives us%
\begin{equation*}
\sum_{\substack{ a,d=0  \\ a+d\equiv 1}}^{1}\left( -1\right)
^{a+1}B(1_{H}\otimes
gx_{1}x_{2};G^{a}X_{1}X_{2},g^{d}x_{1}x_{2})G^{a}X_{1}X_{2}\otimes g^{d}x_{2}
\end{equation*}%
and we get%
\begin{equation*}
-B(1_{H}\otimes gx_{1}x_{2};X_{1}X_{2},gx_{1}x_{2})+B(1_{H}\otimes
gx_{2};X_{1}X_{2},gx_{2})=0
\end{equation*}%
and%
\begin{equation*}
+B(1_{H}\otimes gx_{1}x_{2};GX_{1}X_{2},x_{1}x_{2})+B(1_{H}\otimes
gx_{2};GX_{1}X_{2},x_{2})=0
\end{equation*}%
which are satisfied in view of the form of the elements.

\subsubsection{$G^{a}X_{1}X_{2}\otimes g^{d}x_{1}$}

There are no terms like this.

\subsubsection{$G^{a}X_{1}X_{2}\otimes g^{d}x_{1}x_{2}$}

There are no terms like this.

\subsection{case $x_{2}$}

The right side is%
\begin{equation*}
B^{A}(1_{H}\otimes gx_{1})\otimes B^{H}(1_{H}\otimes gx_{1})\otimes x_{2}
\end{equation*}%
The left side gives us%
\begin{eqnarray*}
a+b_{1}+b_{2}+l_{1}+l_{2}+d+e_{1}+e_{2}+u_{1}+u_{2} &\equiv &0 \\
l_{1}+u_{1} &=&0 \\
l_{2}+u_{2} &=&1
\end{eqnarray*}%
so that%
\begin{eqnarray*}
a+b_{1}+b_{2}+d+e_{1}+e_{2} &\equiv &1 \\
l_{1} &=&u_{1}=0 \\
l_{2}+u_{2} &=&1
\end{eqnarray*}

Since%
\begin{eqnarray*}
\alpha \left( 1_{H};0,0,0,1\right) &=&+\left( a+b_{1}+b_{2}\right) \left(
1+n_{1}+n_{2}\right) +m+n_{1} \\
&=&a+b_{1}+b_{2}
\end{eqnarray*}%
and%
\begin{equation*}
\alpha \left( 1_{H};0,1,0,0\right) =\left( a+b_{1}+b_{2}+1\right) \left(
n_{1}+n_{2}\right) =0
\end{equation*}%
we obtain that the equality is%
\begin{eqnarray*}
&&\sum_{\substack{ a,b_{1},b_{2},d,e_{1}=0  \\ a+b_{1}+b_{2}+d+e_{1}\equiv 0
}}^{1}\left( -1\right) ^{a+b_{1}+b_{2}}B(1_{H}\otimes
gx_{1}x_{2};G^{a}X_{1}^{b_{1}}X_{2}^{b_{2}},g^{d}x_{1}^{e_{1}}x_{2})G^{a}X_{1}^{b_{1}}X_{2}^{b_{2}}\otimes g^{d}x_{1}^{e_{1}}+
\\
&&+\sum_{\substack{ a,b_{1},d,e_{1},e_{2}=0  \\ a+b_{1}+d+e_{1}+e_{2}\equiv
1 }}^{1}B(1_{H}\otimes
gx_{1}x_{2};G^{a}X_{1}^{b_{1}}X_{2},g^{d}x_{1}^{e_{1}}x_{2}^{e_{2}})G^{a}X_{1}^{b_{1}}\otimes g^{d}x_{1}^{e_{1}}x_{2}^{e_{2}}
\\
&=&B^{A}(1_{H}\otimes gx_{1})\otimes B^{H}(1_{H}\otimes gx_{1})
\end{eqnarray*}

\subsection{case $x_{1}x_{2}$}

The right side is $0.$ The left side gives us

\begin{eqnarray*}
a+b_{1}+b_{2}+l_{1}+l_{2}+d+e_{1}+e_{2}+u_{1}+u_{2} &\equiv &0 \\
l_{1}+u_{1} &=&1 \\
l_{2}+u_{2} &=&1
\end{eqnarray*}%
so that%
\begin{eqnarray*}
a+b_{1}+b_{2}+d+e_{1}+e_{2} &\equiv &0 \\
l_{1}+u_{1} &=&1 \\
l_{2}+u_{2} &=&1
\end{eqnarray*}%
this is $0$ as%
\begin{equation*}
B(1_{H}\otimes
gx_{1}x_{2};G^{a}X_{1}^{b_{1}}X_{2}^{b_{2}},g^{d}x_{1}^{e_{1}}x_{2}^{e_{2}})=0%
\text{ whenever }a+b_{1}+b_{2}+d+e_{1}+e_{2}\equiv 0
\end{equation*}

\subsection{case $gx_{1}$ and $gx_{2}$}

The right side is 0. The left side gives us%
\begin{eqnarray*}
a+b_{1}+b_{2}+l_{1}+l_{2}+d+e_{1}+e_{2}+u_{1}+u_{2} &\equiv &1 \\
l_{1}+u_{1} &=&1 \\
l_{2}+u_{2} &=&0
\end{eqnarray*}%
so that%
\begin{eqnarray*}
a+b_{1}+b_{2}+d+e_{1}+e_{2} &\equiv &0 \\
l_{1}+u_{1} &=&1 \\
l_{2}+u_{2} &=&0
\end{eqnarray*}

this is $0$ as%
\begin{equation*}
B(1_{H}\otimes
gx_{1}x_{2};G^{a}X_{1}^{b_{1}}X_{2}^{b_{2}},g^{d}x_{1}^{e_{1}}x_{2}^{e_{2}})=0%
\text{ whenever }a+b_{1}+b_{2}+d+e_{1}+e_{2}\equiv 0
\end{equation*}

The case $gx_{2}$ is completely symmetric.

\subsection{case $gx_{1}x_{2}$}

The right side is%
\begin{equation*}
B^{A}(1_{H}\otimes g)\otimes B^{H}(1_{H}\otimes g)\otimes gx_{1}x_{2}
\end{equation*}%
The left side gives us%
\begin{eqnarray*}
a+b_{1}+b_{2}+l_{1}+l_{2}+d+e_{1}+e_{2}+u_{1}+u_{2} &\equiv &1 \\
l_{1}+u_{1} &=&1 \\
l_{2}+u_{2} &=&1
\end{eqnarray*}%
so that%
\begin{eqnarray*}
a+b_{1}+b_{2}+d+e_{1}+e_{2} &\equiv &1 \\
l_{1}+u_{1} &=&1 \\
l_{2}+u_{2} &=&1.
\end{eqnarray*}%
As%
\begin{equation*}
\alpha \left( 1_{H};0,0,1,1\right) =1+e_{2}
\end{equation*}%
\begin{equation*}
\alpha \left( 1_{H};0,1,1,0\right) =e_{2}+a+b_{1}+b_{2}+1
\end{equation*}

we obtain that
\begin{eqnarray*}
&&\sum_{\substack{ a,b_{1},b_{2},d,=0  \\ a+b_{1}+b_{2}+d\equiv 1}}%
^{1}B(1_{H}\otimes
gx_{1}x_{2};G^{a}X_{1}^{b_{1}}X_{2}^{b_{2}},g^{d}x_{1}x_{2})G^{a}X_{1}^{b_{1}}X_{2}^{b_{2}}\otimes g^{d}+
\\
&&\sum_{\substack{ a,b_{1},d,e_{2}=0  \\ a+b_{1}+d+e_{2}\equiv 1}}^{1}\left(
-1\right) ^{e_{2}+a+b_{1}}B(1_{H}\otimes
gx_{1}x_{2};G^{a}X_{1}^{b_{1}}X_{2},g^{d}x_{1}x_{2}^{e_{2}})G^{a}X_{1}^{b_{1}}\otimes g^{d}x_{2}^{e_{2}}
\\
&&\sum_{\substack{ a,b_{2},d,e_{1}=0  \\ a+b_{2}+d+e_{1}\equiv 1}}^{1}\left(
-1\right) ^{a+1}B(1_{H}\otimes
gx_{1}x_{2};G^{a}X_{1}X_{2}^{b_{2}},g^{d}x_{1}^{e_{1}}x_{2})G^{a}X_{2}^{b_{2}}\otimes g^{d}x_{1}^{e_{1}}
\\
&&\sum_{\substack{ a,d,e_{1},e_{2}=0  \\ a+d+e_{1}+e_{2}\equiv 1}}%
^{1}B(1_{H}\otimes
gx_{1}x_{2};G^{a}X_{1}X_{2},g^{d}x_{1}^{e_{1}}x_{2}^{e_{2}})G^{a}\otimes
g^{d}x_{1}^{e_{1}}x_{2}^{e_{2}}= \\
&&=B^{A}(1_{H}\otimes g)
\end{eqnarray*}

\subsubsection{$G^{a}\otimes g^{d}$}

We obtain%
\begin{equation*}
=\sum_{\substack{ a,d,=0  \\ a+d\equiv 1}}^{1}\left[
\begin{array}{c}
B(1_{H}\otimes gx_{1}x_{2};G^{a},g^{d}x_{1}x_{2})+\left( -1\right)
^{a}B(1_{H}\otimes gx_{1}x_{2};G^{a}X_{2},g^{d}x_{1})+ \\
\left( -1\right) ^{a+1}B(1_{H}\otimes
gx_{1}x_{2};G^{a}X_{1},g^{d}x_{2})+B(1_{H}\otimes
gx_{1}x_{2};G^{a}X_{1}X_{2},g^{d})%
\end{array}%
\right] G^{a}\otimes g^{d}+
\end{equation*}%
and we get%
\begin{gather*}
B(1_{H}\otimes gx_{1}x_{2};1_{A},gx_{1}x_{2})+B(1_{H}\otimes
gx_{1}x_{2};X_{2},gx_{1}) \\
-B(1_{H}\otimes gx_{1}x_{2};X_{1},gx_{2})+B(1_{H}\otimes
gx_{1}x_{2};X_{1}X_{2},g)=B(1_{H}\otimes g;1_{A},g)
\end{gather*}%
and%
\begin{gather*}
B(1_{H}\otimes gx_{1}x_{2};G,x_{1}x_{2})-B(1_{H}\otimes
gx_{1}x_{2};GX_{2},x_{1}) \\
+B(1_{H}\otimes gx_{1}x_{2};GX_{1},x_{2})+B(1_{H}\otimes
gx_{1}x_{2};GX_{1}X_{2},1_{H})=B(1_{H}\otimes g;G,1_{H}).
\end{gather*}%
which are satisfied in view of the form of the elements.

\subsubsection{$G^{a}\otimes g^{d}x_{2}$}

We obtain%
\begin{eqnarray*}
&&\sum_{\substack{ a,d=0  \\ a+d\equiv 0}}^{1}\left[ \left( -1\right)
^{1+a}B(1_{H}\otimes gx_{1}x_{2};G^{a}X_{2},g^{d}x_{1}x_{2})+B(1_{H}\otimes
gx_{1}x_{2};G^{a}X_{1}X_{2},g^{d}x_{2})\right] G^{a}\otimes g^{d}x_{2} \\
&=&\sum_{\substack{ a,d=0  \\ a+d\equiv 0}}^{1}B(1_{H}\otimes
g;G^{a},g^{d}x_{2})G^{a}\otimes g^{d}x_{2}
\end{eqnarray*}%
and we get
\begin{equation*}
-B(1_{H}\otimes gx_{1}x_{2};X_{2},x_{1}x_{2})+B(1_{H}\otimes
gx_{1}x_{2};X_{1}X_{2},x_{2})=-B(g\otimes 1_{H};1_{A},x_{2})
\end{equation*}%
and%
\begin{equation*}
B(1_{H}\otimes gx_{1}x_{2};GX_{2},gx_{1}x_{2})+B(1_{H}\otimes
gx_{1}x_{2};GX_{1}X_{2},gx_{2})=-B(g\otimes 1_{H};G,gx_{2})
\end{equation*}%
which are satisfied in view of the form of the elements.

\subsubsection{$G^{a}\otimes g^{d}x_{1}$}

\begin{eqnarray*}
&&\sum_{\substack{ a,d=0  \\ a+d\equiv 0}}^{1}\left[ \left( -1\right)
^{a+1}B(1_{H}\otimes gx_{1}x_{2};G^{a}X_{1},g^{d}x_{1}x_{2})+B(1_{H}\otimes
gx_{1}x_{2};G^{a}X_{1}X_{2},g^{d}x_{1})\right] G^{a}\otimes g^{d}x_{1} \\
= &&\sum_{\substack{ a,d=0  \\ a+d\equiv 0}}^{1}B(1_{H}\otimes
g;G^{a},g^{d}x_{1})G^{a}\otimes g^{d}x_{1}
\end{eqnarray*}%
and we get%
\begin{equation*}
-B(1_{H}\otimes gx_{1}x_{2};X_{1},x_{1}x_{2})+B(1_{H}\otimes
gx_{1}x_{2};X_{1}X_{2},x_{1})=-B(g\otimes 1_{H};1_{A},x_{1})
\end{equation*}%
and%
\begin{equation*}
B(1_{H}\otimes gx_{1}x_{2};GX_{1},gx_{1}x_{2})+B(1_{H}\otimes
gx_{1}x_{2};GX_{1}X_{2},gx_{1})=-B(g\otimes 1_{H};G,gx_{1})
\end{equation*}%
which are satisfied in view of the form of the elements.

\subsubsection{$G^{a}X_{2}\otimes g^{d}$}

We obtain%
\begin{gather*}
\sum_{\substack{ a,d=0  \\ a+d\equiv 0}}^{1}\left[ B(1_{H}\otimes
gx_{1}x_{2};G^{a}X_{2},g^{d}x_{1}x_{2})+\left( -1\right)
^{a+1}B(1_{H}\otimes gx_{1}x_{2};G^{a}X_{1}X_{2},g^{d}x_{2})\right] \\
G^{a}X_{2}\otimes g^{d}=\sum_{\substack{ a,d=0  \\ a+d\equiv 0}}%
^{1}B(1_{H}\otimes g;G^{a}X_{2},g^{d})G^{a}X_{2}\otimes g^{d}
\end{gather*}%
and we get%
\begin{equation*}
B(1_{H}\otimes gx_{1}x_{2};X_{2},x_{1}x_{2})-B(1_{H}\otimes
gx_{1}x_{2};X_{1}X_{2},x_{2})=B(g\otimes 1_{H};X_{2},1_{H})
\end{equation*}%
and%
\begin{equation*}
B(1_{H}\otimes gx_{1}x_{2};GX_{2},gx_{1}x_{2})+B(1_{H}\otimes
gx_{1}x_{2};GX_{1}X_{2},gx_{2})=B(g\otimes 1_{H};GX_{2},g)
\end{equation*}%
which are satisfied in view of the form of the elements.

\subsubsection{$G^{a}X_{1}\otimes g^{d}$}

We obtain%
\begin{eqnarray*}
&&\sum_{\substack{ a,d=0  \\ a+d\equiv 0}}^{1}\left[
\begin{array}{c}
B(1_{H}\otimes gx_{1}x_{2};G^{a}X_{1},g^{d}x_{1}x_{2}) \\
+\left( -1\right) ^{a+1}B(1_{H}\otimes
gx_{1}x_{2};G^{a}X_{1}X_{2},g^{d}x_{1})%
\end{array}%
\right] G^{a}X_{1}\otimes g^{d} \\
&=&\sum_{\substack{ a,d=0  \\ a+d\equiv 0}}^{1}B(g\otimes
1_{H};G^{a}X_{1},g^{d})G^{a}X_{1}\otimes g^{d}
\end{eqnarray*}%
and we get%
\begin{equation*}
B(1_{H}\otimes gx_{1}x_{2};X_{1},x_{1}x_{2})-B(1_{H}\otimes
gx_{1}x_{2};X_{1}X_{2},x_{1})=B(1_{H}\otimes g;X_{1},1_{H})
\end{equation*}%
and

\begin{equation*}
B(1_{H}\otimes gx_{1}x_{2};GX_{1},gx_{1}x_{2})+B(1_{H}\otimes
gx_{1}x_{2};GX_{1}X_{2},gx_{1})=B(1_{H}\otimes g;GX_{1},g)
\end{equation*}%
which we already got in case $G^{a}\otimes g^{d}x_{1}.$

\subsubsection{$G^{a}\otimes g^{d}x_{1}x_{2}$}

We obtain%
\begin{equation*}
\sum_{\substack{ a,d=0  \\ a+d\equiv 1}}^{1}B(1_{H}\otimes
gx_{1}x_{2};G^{a}X_{1}X_{2},g^{d}x_{1}x_{2})G^{a}\otimes
g^{d}x_{1}x_{2}=\sum _{\substack{ a,d=0  \\ a+d\equiv 1}}^{1}B(1_{H}\otimes
g;G^{a},g^{d}x_{1}x_{2})G^{a}\otimes g^{d}x_{1}x_{2}
\end{equation*}%
and we get
\begin{equation*}
B(1_{H}\otimes gx_{1}x_{2};X_{1}X_{2},gx_{1}x_{2})=B(1_{H}\otimes
g;1_{A},gx_{1}x_{2})
\end{equation*}%
and%
\begin{equation*}
B(1_{H}\otimes gx_{1}x_{2};GX_{1}X_{2},x_{1}x_{2})=B(1_{H}\otimes
g;G,x_{1}x_{2})
\end{equation*}%
which hold in view of the form of the elements.

\subsubsection{$G^{a}X_{2}\otimes g^{d}x_{2}$}

There is no term like this in the left side. Thus we get%
\begin{equation*}
B(1_{H}\otimes g;X_{2},gx_{2})=0
\end{equation*}%
and%
\begin{equation*}
B(1_{H}\otimes g;GX_{2},x_{2})=0
\end{equation*}%
which hold in view of the form of the element.

\subsubsection{$G^{a}X_{1}\otimes g^{d}x_{2}$}

We obtain%
\begin{eqnarray*}
&&\sum_{\substack{ a,d=0  \\ a+d\equiv 1}}^{1}\left( -1\right)
^{a}B(1_{H}\otimes
gx_{1}x_{2};G^{a}X_{1}X_{2},g^{d}x_{1}x_{2})G^{a}X_{1}\otimes g^{d}x_{2} \\
&=&\sum_{\substack{ a,d=0  \\ a+d\equiv 1}}^{1}B(1_{H}\otimes
g;G^{a}X_{1},g^{d}x_{2})G^{a}X_{1}\otimes g^{d}x_{2}
\end{eqnarray*}%
and we get%
\begin{equation*}
B(1_{H}\otimes gx_{1}x_{2};X_{1}X_{2},gx_{1}x_{2})=B(1_{H}\otimes
g;X_{1},gx_{2})\text{ }
\end{equation*}%
and%
\begin{equation*}
-B(1_{H}\otimes gx_{1}x_{2};GX_{1}X_{2},x_{1}x_{2})=B(1_{H}\otimes
g;GX_{1},x_{2})
\end{equation*}%
which hold in view of the form of the elements.

\subsubsection{$G^{a}X_{2}\otimes g^{d}x_{1}$}

We obtain%
\begin{equation*}
\sum_{\substack{ a,d=0  \\ a+d\equiv 1}}^{1}\left( -1\right)
^{a+1}B(1_{H}\otimes
gx_{1}x_{2};G^{a}X_{1}X_{2},g^{d}x_{1}x_{2})G^{a}X_{2}\otimes g^{d}x_{1}
\end{equation*}%
and we get%
\begin{equation*}
-B(1_{H}\otimes gx_{1}x_{2};X_{1}X_{2},gx_{1}x_{2})=B(1_{H}\otimes
g;X_{2},gx_{1})
\end{equation*}%
and%
\begin{equation*}
B(1_{H}\otimes gx_{1}x_{2};GX_{1}X_{2},x_{1}x_{2})=B(1_{H}\otimes
g;GX_{2},x_{1})
\end{equation*}%
which hold in view of the form of the elements.

\subsubsection{$G^{a}X_{1}\otimes g^{d}x_{1}$}

There is no term like this in the left side. Thus we obtain%
\begin{equation*}
B(1_{H}\otimes g;X_{1},gx_{1})=0
\end{equation*}%
and%
\begin{equation*}
B(1_{H}\otimes g;GX_{1},x_{1})=0
\end{equation*}%
which hold in view of the form of the element.

\subsubsection{$G^{a}X_{1}X_{2}\otimes g^{d}$}

We obtain%
\begin{equation*}
\sum_{\substack{ a,d=0  \\ a+d\equiv 1}}^{1}B(1_{H}\otimes
gx_{1}x_{2};G^{a}X_{1}X_{2},g^{d}x_{1}x_{2})G^{a}X_{1}X_{2}\otimes g^{d}+
\end{equation*}%
\begin{equation*}
B(1_{H}\otimes gx_{1}x_{2};X_{1}X_{2},gx_{1}x_{2})=B\left( 1_{H}\otimes
g;X_{1}X_{2},g\right)
\end{equation*}%
and%
\begin{equation*}
B(1_{H}\otimes gx_{1}x_{2};GX_{1}X_{2},x_{1}x_{2})=B\left( 1_{H}\otimes
g;GX_{1}X_{2},1_{H}\right)
\end{equation*}%
which hold in view of the form of the elements.

\subsubsection{$G^{a}X_{2}\otimes g^{d}x_{1}x_{2}$}

Nothing on the left side. We get%
\begin{equation*}
B\left( 1_{H}\otimes g;X_{2},x_{1}x_{2}\right) =0
\end{equation*}%
and%
\begin{equation*}
B\left( 1_{H}\otimes g;GX_{2},gx_{1}x_{2}\right) =0
\end{equation*}%
which hold in view of the form of the element.

\subsubsection{$G^{a}X_{1}\otimes g^{d}x_{1}x_{2}$}

Nothing on the left side. We get%
\begin{equation*}
B\left( 1_{H}\otimes g;X_{1},x_{1}x_{2}\right) =0
\end{equation*}%
and%
\begin{equation*}
B\left( 1_{H}\otimes g;GX_{1},gx_{1}x_{2}\right) =0
\end{equation*}%
which hold in view of the form of the element.

\subsubsection{$G^{a}X_{1}X_{2}\otimes g^{d}x_{2}$}

Nothing on the left side. We get%
\begin{equation*}
B\left( 1_{H}\otimes g;X_{1}X_{2},x_{2}\right) =0
\end{equation*}%
\begin{equation*}
B\left( 1_{H}\otimes g;GX_{1}X_{2},gx_{2}\right) =0
\end{equation*}%
which hold in view of the form of the element.

\subsubsection{$G^{a}X_{1}X_{2}\otimes g^{d}x_{1}$}

Nothing on the left side. We get%
\begin{equation*}
B\left( 1_{H}\otimes g;X_{1}X_{2},x_{1}\right) =0
\end{equation*}%
and%
\begin{equation*}
B\left( 1_{H}\otimes g;GX_{1}X_{2},gx_{1}\right) =0
\end{equation*}%
which hold in view of the form of the element.

\subsubsection{$G^{a}X_{1}X_{2}\otimes g^{d}x_{1}x_{2}$}

Nothing on the left side. We get%
\begin{equation*}
B\left( 1_{H}\otimes g;X_{1}X_{2},x_{1}x_{2}\right) =0
\end{equation*}

and%
\begin{equation*}
B\left( 1_{H}\otimes g;GX_{1}X_{2},gx_{1}x_{2}\right) =0
\end{equation*}%
which hold in view of the form of the element.

\section{$B\left( x_{1}\otimes x_{1}x_{2}\right) $}

By applying $\left( \ref{simpl1}\right) ,$ we get%
\begin{eqnarray*}
B(x_{1}\otimes x_{1}x_{2}) &=&B(x_{1}\otimes x_{2})(1_{A}\otimes
x_{1})-(1_{A}\otimes gx_{1})B(x_{1}\otimes x_{2})(1_{A}\otimes g) \\
&&-(1_{A}\otimes g)B(gx_{1}x_{1}\otimes x_{2})(1_{A}\otimes g)
\end{eqnarray*}%
i.e.%
\begin{equation*}
B(x_{1}\otimes x_{1}x_{2})=B(x_{1}\otimes x_{2})(1_{A}\otimes
x_{1})-(1_{A}\otimes gx_{1})B(x_{1}\otimes x_{2})(1_{A}\otimes g)
\end{equation*}%
By using $\left( \ref{form x1otx2}\right) ,$ we get%
\begin{eqnarray}
B(x_{1}\otimes x_{1}x_{2}) &=&B(x_{1}\otimes 1_{H})(1_{A}\otimes x_{1}x_{2})
\label{form x1otx1x2} \\
&&+(1_{A}\otimes gx_{2})B(x_{1}\otimes 1_{H})(1_{A}\otimes gx_{1})  \notag \\
&&-(1_{A}\otimes g)B(gx_{1}x_{2}\otimes 1_{H})(1_{A}\otimes gx_{1})  \notag
\\
&&-(1_{A}\otimes gx_{1})B(x_{1}\otimes 1_{H})(1_{A}\otimes gx_{2})  \notag \\
&&+(1_{A}\otimes x_{1}x_{2})B(x_{1}\otimes 1_{H})  \notag \\
&&-(1_{A}\otimes x_{1})B(gx_{1}x_{2}\otimes 1_{H})  \notag
\end{eqnarray}%
and hence%
\begin{eqnarray*}
B(x_{1}\otimes x_{1}x_{2}) &=&-2B(gx_{1}x_{2}\otimes
1_{H};1_{A},g)1_{A}\otimes gx_{1}+ \\
&&\left[ 4B\left( x_{1}\otimes 1_{H};G,g\right) +2B(gx_{1}x_{2}\otimes
1_{H};G,gx_{2})\right] G\otimes gx_{1}x_{2}+ \\
&&+\left[
\begin{array}{c}
-4B(g\otimes 1_{H};1_{A},g)-4B(x_{1}\otimes 1_{H};1_{A},gx_{1}) \\
-2B(x_{2}\otimes 1_{H};1_{A},gx_{2})-2B(gx_{1}x_{2}\otimes
1_{H};1_{A},gx_{1}x_{2})%
\end{array}%
\right] X_{1}\otimes gx_{1}x_{2}+ \\
&&-2B(x_{1}\otimes 1_{H};1_{A},gx_{2})X_{2}\otimes gx_{1}x_{2}+ \\
&&+\left[
\begin{array}{c}
-2B(g\otimes 1_{H};1_{A},g)-2B(x_{2}\otimes \ 1_{H};1_{A},gx_{2}) \\
-2B(x_{1}\otimes 1_{H};1_{A},gx_{1})-2B(gx_{1}x_{2}\otimes
1_{H};1_{A},gx_{1}x_{2})%
\end{array}%
\right] X_{1}X_{2}\otimes gx_{1}+ \\
&&+\left[ -2B(x_{2}\otimes 1_{H};G,g)+2B(gx_{1}x_{2}\otimes 1_{H};G,gx_{1})%
\right] GX_{1}\otimes gx_{1}+ \\
&&+\left[ +2B(x_{1}\otimes 1_{H};G,g)+2B(gx_{1}x_{2}\otimes 1_{H};G,gx_{2})%
\right] GX_{2}\otimes gx_{1}+ \\
&&+\left[ -2B(g\otimes 1_{H};G,gx_{2})+2B(x_{1}\otimes 1_{H};G,gx_{1}x_{2})%
\right] GX_{1}X_{2}\otimes gx_{1}x_{2}
\end{eqnarray*}

We write the Casimir formula for $B\left( x_{1}\otimes x_{1}x_{2}\right) :$%
\begin{eqnarray*}
&&\sum_{a,b_{1},b_{2},d,e_{1},e_{2}=0}^{1}\sum_{l_{1}=0}^{b_{1}}%
\sum_{l_{2}=0}^{b_{2}}\sum_{u_{1}=0}^{e_{1}}\sum_{u_{2}=0}^{e_{2}}\left(
-1\right) ^{\alpha \left( x_{1};l_{1},l_{2},u_{1},u_{2}\right) } \\
&&B(g\otimes
x_{1}x_{2};G^{a}X_{1}^{b_{1}}X_{2}^{b_{2}},g^{d}x_{1}^{e_{1}}x_{2}^{e_{2}})G^{a}X_{1}^{b_{1}-l_{1}}X_{2}^{b_{2}-l_{2}}\otimes g^{d}x_{1}^{e_{1}-u_{1}}x_{2}^{e_{2}-u_{2}}\otimes
\\
&&g^{a+b_{1}+b_{2}+l_{1}+l_{2}+d+e_{1}+e_{2}+u_{1}+u_{2}}x_{1}^{l_{1}+u_{1}+1}x_{2}^{l_{2}+u_{2}}
\\
&&+\sum_{a,b_{1},b_{2},d,e_{1},e_{2}=0}^{1}\sum_{l_{1}=0}^{b_{1}}%
\sum_{l_{2}=0}^{b_{2}}\sum_{u_{1}=0}^{e_{1}}\sum_{u_{2}=0}^{e_{2}}\left(
-1\right) ^{\alpha \left( 1_{H};l_{1},l_{2},u_{1},u_{2}\right) } \\
&&B(x_{1}\otimes
x_{1}x_{2};G^{a}X_{1}^{b_{1}}X_{2}^{b_{2}},g^{d}x_{1}^{e_{1}}x_{2}^{e_{2}})G^{a}X_{1}^{b_{1}-l_{1}}X_{2}^{b_{2}-l_{2}}\otimes g^{d}x_{1}^{e_{1}-u_{1}}x_{2}^{e_{2}-u_{2}}\otimes
\\
&&g^{a+b_{1}+b_{2}+l_{1}+l_{2}+d+e_{1}+e_{2}+u_{1}+u_{2}}x_{1}^{l_{1}+u_{1}}x_{2}^{l_{2}+u_{2}}
\\
&=&B^{A}(x_{1}\otimes x_{1}x_{2})\otimes B^{H}(x_{1}\otimes
x_{1}x_{2})\otimes 1_{H} \\
&&B^{A}(x_{1}\otimes x_{1})\otimes B^{H}(x_{1}\otimes x_{1})\otimes gx_{2} \\
&&\left( -1\right) B^{A}(x_{1}\otimes x_{2})\otimes B^{H}(x_{1}\otimes
x_{2})\otimes gx_{1} \\
&&B^{A}(x_{1}\otimes 1_{H})\otimes B^{H}(x_{1}\otimes 1_{H})\otimes
x_{1}x_{2}
\end{eqnarray*}%
and we obtain%
\begin{eqnarray*}
&&\sum_{a,b_{1},b_{2},d,e_{1},e_{2}=0}^{1}\sum_{l_{1}=0}^{b_{1}}%
\sum_{l_{2}=0}^{b_{2}}\sum_{u_{1}=0}^{e_{1}}\sum_{u_{2}=0}^{e_{2}}\left(
-1\right) ^{\alpha \left( x_{1};l_{1},l_{2},u_{1},u_{2}\right) } \\
&&B(g\otimes
x_{1}x_{2};G^{a}X_{1}^{b_{1}}X_{2}^{b_{2}},g^{d}x_{1}^{e_{1}}x_{2}^{e_{2}})G^{a}X_{1}^{b_{1}-l_{1}}X_{2}^{b_{2}-l_{2}}\otimes g^{d}x_{1}^{e_{1}-u_{1}}x_{2}^{e_{2}-u_{2}}\otimes
\\
&&g^{a+b_{1}+b_{2}+l_{1}+l_{2}+d+e_{1}+e_{2}+u_{1}+u_{2}}x_{1}^{l_{1}+u_{1}+1}x_{2}^{l_{2}+u_{2}}
\\
&&+\sum_{a,b_{1},b_{2},d,e_{1},e_{2}=0}^{1}\sum_{l_{1}=0}^{b_{1}}%
\sum_{l_{2}=0}^{b_{2}}\sum_{u_{1}=0}^{e_{1}}\sum_{u_{2}=0}^{e_{2}}\left(
-1\right) ^{\alpha \left( 1_{H};l_{1},l_{2},u_{1},u_{2}\right) } \\
&&B(x_{1}\otimes
x_{1}x_{2};G^{a}X_{1}^{b_{1}}X_{2}^{b_{2}},g^{d}x_{1}^{e_{1}}x_{2}^{e_{2}})G^{a}X_{1}^{b_{1}-l_{1}}X_{2}^{b_{2}-l_{2}}\otimes g^{d}x_{1}^{e_{1}-u_{1}}x_{2}^{e_{2}-u_{2}}\otimes
\\
&&g^{a+b_{1}+b_{2}+l_{1}+l_{2}+d+e_{1}+e_{2}+u_{1}+u_{2}}x_{1}^{l_{1}+u_{1}}x_{2}^{l_{2}+u_{2}}
\\
&=&B^{A}(x_{1}\otimes x_{1}x_{2})\otimes B^{H}(x_{1}\otimes
x_{1}x_{2})\otimes 1_{H} \\
&&B^{A}(x_{1}\otimes x_{1})\otimes B^{H}(x_{1}\otimes x_{1})\otimes gx_{2} \\
&&\left( -1\right) B^{A}(x_{1}\otimes x_{2})\otimes B^{H}(x_{1}\otimes
x_{2})\otimes gx_{1} \\
&&B^{A}(x_{1}\otimes 1_{H})\otimes B^{H}(x_{1}\otimes 1_{H})\otimes
x_{1}x_{2}
\end{eqnarray*}

\subsection{$B\left( x_{1}\otimes x_{1}x_{2},1_{A},gx_{1}\right) $}

Now we look for the terms in $B(x_{1}\otimes x_{1}x_{2};1_{A},gx_{1})$ in
the equality above. We get%
\begin{eqnarray*}
a &=&b_{1}=b_{2}=0 \\
d &=&e_{1}=1,e_{2}=0
\end{eqnarray*}%
and we get, since $\alpha \left( 1_{H};0,0,0,0\right) =0$ and $\alpha \left(
1_{H};0,0,1,0\right) =e_{2}+\left( a+b_{1}+b_{2}\right) =0$
\begin{eqnarray*}
&&+B(x_{1}\otimes x_{1}x_{2};1_{A},gx_{1})1_{A}\otimes gx_{1}\otimes 1_{H} \\
&&B(x_{1}\otimes x_{1}x_{2};1_{A},gx_{1})1_{A}\otimes g\otimes gx_{1}
\end{eqnarray*}

\subsubsection{Case $1_{A}\otimes gx_{1}\otimes 1_{H}$}

\begin{eqnarray*}
&&\sum_{a,b_{1},b_{2},d,e_{1},e_{2}=0}^{1}\sum_{l_{1}=0}^{b_{1}}%
\sum_{l_{2}=0}^{b_{2}}\sum_{u_{1}=0}^{e_{1}}\sum_{u_{2}=0}^{e_{2}}\left(
-1\right) ^{\alpha \left( x_{1};l_{1},l_{2},u_{1},u_{2}\right) } \\
&&B(g\otimes
x_{1}x_{2};G^{a}X_{1}^{b_{1}}X_{2}^{b_{2}},g^{d}x_{1}^{e_{1}}x_{2}^{e_{2}})G^{a}X_{1}^{b_{1}-l_{1}}X_{2}^{b_{2}-l_{2}}\otimes g^{d}x_{1}^{e_{1}-u_{1}}x_{2}^{e_{2}-u_{2}}\otimes
\\
&&g^{a+b_{1}+b_{2}+l_{1}+l_{2}+d+e_{1}+e_{2}+u_{1}+u_{2}}x_{1}^{l_{1}+u_{1}+1}x_{2}^{l_{2}+u_{2}}
\\
&&+\sum_{a,b_{1},b_{2},d,e_{1},e_{2}=0}^{1}\sum_{l_{1}=0}^{b_{1}}%
\sum_{l_{2}=0}^{b_{2}}\sum_{u_{1}=0}^{e_{1}}\sum_{u_{2}=0}^{e_{2}}\left(
-1\right) ^{\alpha \left( 1_{H};l_{1},l_{2},u_{1},u_{2}\right) } \\
&&B(x_{1}\otimes
x_{1}x_{2};G^{a}X_{1}^{b_{1}}X_{2}^{b_{2}},g^{d}x_{1}^{e_{1}}x_{2}^{e_{2}})G^{a}X_{1}^{b_{1}-l_{1}}X_{2}^{b_{2}-l_{2}}\otimes g^{d}x_{1}^{e_{1}-u_{1}}x_{2}^{e_{2}-u_{2}}\otimes
\\
&&g^{a+b_{1}+b_{2}+l_{1}+l_{2}+d+e_{1}+e_{2}+u_{1}+u_{2}}x_{1}^{l_{1}+u_{1}}x_{2}^{l_{2}+u_{2}}
\\
&=&B^{A}(x_{1}\otimes x_{1}x_{2})\otimes B^{H}(x_{1}\otimes
x_{1}x_{2})\otimes 1_{H} \\
&&B^{A}(x_{1}\otimes x_{1})\otimes B^{H}(x_{1}\otimes x_{1})\otimes gx_{2} \\
&&-B^{A}(x_{1}\otimes x_{2})\otimes B^{H}(x_{1}\otimes x_{2})\otimes gx_{1}
\\
&&B^{A}(x_{1}\otimes 1_{H})\otimes B^{H}(x_{1}\otimes 1_{H})\otimes
x_{1}x_{2}
\end{eqnarray*}

First summand left side is not possible. Second summand
\begin{equation*}
\sum_{u_{1}=0}^{1}B(x_{1}\otimes x_{1}x_{2};1_{A},gx_{1})1_{A}\otimes
gx_{1}^{1-u_{1}}\otimes 1_{H}
\end{equation*}%
\begin{eqnarray*}
l_{1} &=&u_{1}=0 \\
l_{2} &=&u_{2}=0 \\
a+b_{1}+b_{2}+d+e_{1}+e_{2} &\equiv &0 \\
d &=&1 \\
e_{1} &=&1 \\
e_{2} &=&0 \\
a &=&0 \\
b_{1} &=&0 \\
b_{2} &=&0
\end{eqnarray*}%
and we get%
\begin{equation*}
+B(x_{1}\otimes x_{1}x_{2};1_{A},gx_{1})1_{A}\otimes gx_{1}\otimes 1_{H}
\end{equation*}%
On the right side the only contribution is from the first summand%
\begin{equation*}
B(x_{1}\otimes x_{1}x_{2};1_{A},gx_{1})1_{A}\otimes gx_{1}\otimes 1_{H}
\end{equation*}%
and we get nothing new.

\subsubsection{Case $1_{A}\otimes g\otimes gx_{1}$}

First summand of the left side%
\begin{eqnarray*}
&&\sum_{a,b_{1},b_{2},d,e_{1},e_{2}=0}^{1}\sum_{l_{1}=0}^{b_{1}}%
\sum_{l_{2}=0}^{b_{2}}\sum_{u_{1}=0}^{e_{1}}\sum_{u_{2}=0}^{e_{2}}\left(
-1\right) ^{\alpha \left( x_{1};l_{1},l_{2},u_{1},u_{2}\right) } \\
&&B(g\otimes
x_{1}x_{2};G^{a}X_{1}^{b_{1}}X_{2}^{b_{2}},g^{d}x_{1}^{e_{1}}x_{2}^{e_{2}})G^{a}X_{1}^{b_{1}-l_{1}}X_{2}^{b_{2}-l_{2}}\otimes g^{d}x_{1}^{e_{1}-u_{1}}x_{2}^{e_{2}-u_{2}}\otimes
\\
&&g^{a+b_{1}+b_{2}+l_{1}+l_{2}+d+e_{1}+e_{2}+u_{1}+u_{2}}x_{1}^{l_{1}+u_{1}+1}x_{2}^{l_{2}+u_{2}}
\end{eqnarray*}%
\begin{eqnarray*}
l_{1} &=&u_{1}=l_{2}=u_{2}=0 \\
a+b_{1}+b_{2}+d+e_{1}+e_{2} &\equiv &1 \\
d &=&1,e_{1}=e_{2}=0 \\
a &=&b_{1}=b_{2}=0
\end{eqnarray*}%
and we get%
\begin{equation*}
B(g\otimes x_{1}x_{2};1_{A},g)1_{A}\otimes g\otimes gx_{1}
\end{equation*}%
Second summand of the left side%
\begin{eqnarray*}
&&+\sum_{a,b_{1},b_{2},d,e_{1},e_{2}=0}^{1}\sum_{l_{1}=0}^{b_{1}}%
\sum_{l_{2}=0}^{b_{2}}\sum_{u_{1}=0}^{e_{1}}\sum_{u_{2}=0}^{e_{2}}\left(
-1\right) ^{\alpha \left( 1_{H};l_{1},l_{2},u_{1},u_{2}\right) } \\
&&B(x_{1}\otimes
x_{1}x_{2};G^{a}X_{1}^{b_{1}}X_{2}^{b_{2}},g^{d}x_{1}^{e_{1}}x_{2}^{e_{2}})G^{a}X_{1}^{b_{1}-l_{1}}X_{2}^{b_{2}-l_{2}}\otimes g^{d}x_{1}^{e_{1}-u_{1}}x_{2}^{e_{2}-u_{2}}\otimes
\\
&&g^{a+b_{1}+b_{2}+l_{1}+l_{2}+d+e_{1}+e_{2}+u_{1}+u_{2}}x_{1}^{l_{1}+u_{1}}x_{2}^{l_{2}+u_{2}}
\end{eqnarray*}

\begin{eqnarray*}
l_{1}+u_{1} &=&1,l_{2}=u_{2}=0 \\
a+b_{1}+b_{2}+d+e_{1}+e_{2} &\equiv &0 \\
d &=&1,e_{1}=u_{1},e_{2}=u_{2}=0 \\
a &=&b_{2}=0,b_{1}=l_{1}
\end{eqnarray*}%
since $\alpha \left( 1_{H};0,0,1,0\right) =e_{2}+\left( a+b_{1}+b_{2}\right)
=0$ and $\alpha \left( 1_{H};1,0,0,0\right) =b_{2}=0,$ we get%
\begin{equation*}
\left[ B(x_{1}\otimes x_{1}x_{2};1_{A},gx_{1})+B(x_{1}\otimes
x_{1}x_{2};X_{1},g)\right] 1_{A}\otimes g\otimes gx_{1}
\end{equation*}

Considering also the right side $-B^{A}(x_{1}\otimes x_{2})\otimes
B^{H}(x_{1}\otimes x_{2})\otimes gx_{1}$ we get%
\begin{equation*}
B(g\otimes x_{1}x_{2};1_{A},g)+B(x_{1}\otimes
x_{1}x_{2};1_{A},gx_{1})+B(x_{1}\otimes x_{1}x_{2};X_{1},g)=-B(x_{1}\otimes
x_{2};1_{A},g)
\end{equation*}

which is satisfied in view of the form of the elements.

\subsection{$B\left( x_{1}\otimes x_{1}x_{2},G,gx_{1}x_{2}\right) $}

Now we look for the terms in $B\left( x_{1}\otimes
x_{1}x_{2},G,gx_{1}x_{2}\right) $ in the Casimir equality and we get

\begin{equation*}
a=1,b_{1}=b_{2}=0,d=e_{1}=e_{2}=1
\end{equation*}%
and we get%
\begin{eqnarray*}
&&+\sum_{u_{1}=0}^{1}\sum_{u_{2}=0}^{1}\left( -1\right) ^{\alpha \left(
1_{H};0,0,u_{1},u_{2}\right) }B(x_{1}\otimes
x_{1}x_{2};G,gx_{1}x_{2})G\otimes gx_{1}^{1-u_{1}}x_{2}^{1-u_{2}}\otimes
g^{u_{1}+u_{2}}x_{1}^{u_{1}}x_{2}^{u_{2}} \\
&=&\left( -1\right) ^{\alpha \left( 1_{H};0,0,0,0\right) }B(x_{1}\otimes
x_{1}x_{2};G,gx_{1}x_{2})G\otimes gx_{1}x_{2}\otimes 1_{H}+ \\
&&+\left( -1\right) ^{\alpha \left( 1_{H};0,0,0,1\right) }B(x_{1}\otimes
x_{1}x_{2};G,gx_{1}x_{2})G\otimes gx_{1}\otimes gx_{2}+ \\
&&+\left( -1\right) ^{\alpha \left( 1_{H};0,0,1,0\right) }B(x_{1}\otimes
x_{1}x_{2};G,gx_{1}x_{2})G\otimes gx_{2}\otimes gx_{1}+ \\
&&+\left( -1\right) ^{\alpha \left( 1_{H};0,0,1,1\right) }B(x_{1}\otimes
x_{1}x_{2};G,gx_{1}x_{2})G\otimes g\otimes x_{1}x_{2}
\end{eqnarray*}

\subsubsection{Case $G\otimes gx_{1}x_{2}\otimes 1_{H}$}

First summand of the left side is not possible. Second summand of the left
side%
\begin{eqnarray*}
&&+\sum_{a,b_{1},b_{2},d,e_{1},e_{2}=0}^{1}\sum_{l_{1}=0}^{b_{1}}%
\sum_{l_{2}=0}^{b_{2}}\sum_{u_{1}=0}^{e_{1}}\sum_{u_{2}=0}^{e_{2}}\left(
-1\right) ^{\alpha \left( 1_{H};l_{1},l_{2},u_{1},u_{2}\right) } \\
&&B(x_{1}\otimes
x_{1}x_{2};G^{a}X_{1}^{b_{1}}X_{2}^{b_{2}},g^{d}x_{1}^{e_{1}}x_{2}^{e_{2}})G^{a}X_{1}^{b_{1}-l_{1}}X_{2}^{b_{2}-l_{2}}\otimes g^{d}x_{1}^{e_{1}-u_{1}}x_{2}^{e_{2}-u_{2}}\otimes
\\
&&g^{a+b_{1}+b_{2}+l_{1}+l_{2}+d+e_{1}+e_{2}+u_{1}+u_{2}}x_{1}^{l_{1}+u_{1}}x_{2}^{l_{2}+u_{2}}
\end{eqnarray*}

we get%
\begin{eqnarray*}
l_{1} &=&u_{1}=l_{2}=u_{2}=0 \\
d &=&e_{1}=e_{2}=1 \\
a &=&1,b_{1}=b_{2}=0
\end{eqnarray*}%
and we obtain%
\begin{equation*}
B(x_{1}\otimes x_{1}x_{2};G,gx_{1}x_{2})G\otimes gx_{1}x_{2}\otimes 1_{H}
\end{equation*}%
Considering also the right side, we get\newline
\begin{equation*}
B(x_{1}\otimes x_{1}x_{2};G,gx_{1}x_{2})G\otimes gx_{1}x_{2}\otimes
1_{H}=B(x_{1}\otimes x_{1}x_{2};G,gx_{1}x_{2})G\otimes gx_{1}x_{2}\otimes
1_{H}
\end{equation*}%
and we obtain no new information.

\subsubsection{Case $G\otimes gx_{1}\otimes gx_{2}$}

First summand of the left side is not possible. Second summand of the left
side%
\begin{eqnarray*}
&&+\sum_{a,b_{1},b_{2},d,e_{1},e_{2}=0}^{1}\sum_{l_{1}=0}^{b_{1}}%
\sum_{l_{2}=0}^{b_{2}}\sum_{u_{1}=0}^{e_{1}}\sum_{u_{2}=0}^{e_{2}}\left(
-1\right) ^{\alpha \left( 1_{H};l_{1},l_{2},u_{1},u_{2}\right) } \\
&&B(x_{1}\otimes
x_{1}x_{2};G^{a}X_{1}^{b_{1}}X_{2}^{b_{2}},g^{d}x_{1}^{e_{1}}x_{2}^{e_{2}})G^{a}X_{1}^{b_{1}-l_{1}}X_{2}^{b_{2}-l_{2}}\otimes g^{d}x_{1}^{e_{1}-u_{1}}x_{2}^{e_{2}-u_{2}}\otimes
\\
&&g^{a+b_{1}+b_{2}+l_{1}+l_{2}+d+e_{1}+e_{2}+u_{1}+u_{2}}x_{1}^{l_{1}+u_{1}}x_{2}^{l_{2}+u_{2}}
\end{eqnarray*}%
\begin{eqnarray*}
l_{1} &=&u_{1}=0,l_{2}+u_{2}=1 \\
a+b_{1}+b_{2}+d+e_{1}+e_{2} &=&0 \\
d &=&1,e_{1}=1,e_{2}=u_{2} \\
a &=&1,b_{1}=0,b_{2}=l_{2}.
\end{eqnarray*}%
Since $\alpha \left( 1_{H};0,0,0,1\right) =a+b_{1}+b_{2}=1$ and $\alpha
\left( 1_{H};0,1,0,0\right) =0$ we obtain%
\begin{equation*}
\left[ -B(x_{1}\otimes x_{1}x_{2};G,gx_{1}x_{2})+B(x_{1}\otimes
x_{1}x_{2};GX_{2},gx_{1})\right] G\otimes gx_{1}\otimes gx_{2}.
\end{equation*}%
Considering also the right side, we obtain%
\begin{equation*}
\left[ -B(x_{1}\otimes x_{1}x_{2};G,gx_{1}x_{2})+B(x_{1}\otimes
x_{1}x_{2};GX_{2},gx_{1})\right] G\otimes gx_{1}\otimes
gx_{2}=B(x_{1}\otimes x_{1};G,gx_{1})G\otimes gx_{1}\otimes gx_{2}
\end{equation*}%
which is satisfied in view of the form of the elements.

\subsubsection{Case $G\otimes gx_{2}\otimes gx_{1}$}

First summand left side gives us
\begin{eqnarray*}
l_{1} &=&u_{1}=0=l_{2}=u_{2} \\
a+b_{1}+b_{2}+d+e_{1}+e_{2} &\equiv &1 \\
d &=&1,e_{1}=0,e_{2}=1 \\
a &=&1,b_{1}=0,b_{2}=0.
\end{eqnarray*}%
Since $\alpha \left( x_{1};0,0,0,0\right) =a+b_{1}+b_{2}=1$ we get%
\begin{equation*}
-B(g\otimes x_{1}x_{2};G,gx_{2})G\otimes gx_{2}\otimes gx_{1}.
\end{equation*}%
Second summand left side gives us%
\begin{eqnarray*}
&&+\sum_{a,b_{1},b_{2},d,e_{1},e_{2}=0}^{1}\sum_{l_{1}=0}^{b_{1}}%
\sum_{l_{2}=0}^{b_{2}}\sum_{u_{1}=0}^{e_{1}}\sum_{u_{2}=0}^{e_{2}}\left(
-1\right) ^{\alpha \left( 1_{H};l_{1},l_{2},u_{1},u_{2}\right) } \\
&&B(x_{1}\otimes
x_{1}x_{2};G^{a}X_{1}^{b_{1}}X_{2}^{b_{2}},g^{d}x_{1}^{e_{1}}x_{2}^{e_{2}})G^{a}X_{1}^{b_{1}-l_{1}}X_{2}^{b_{2}-l_{2}}\otimes g^{d}x_{1}^{e_{1}-u_{1}}x_{2}^{e_{2}-u_{2}}\otimes
\\
&&g^{a+b_{1}+b_{2}+l_{1}+l_{2}+d+e_{1}+e_{2}+u_{1}+u_{2}}x_{1}^{l_{1}+u_{1}}x_{2}^{l_{2}+u_{2}}
\end{eqnarray*}%
\begin{eqnarray*}
l_{1}+u_{1} &=&1,l_{2}=u_{2}=0 \\
a+b_{1}+b_{2}+d+e_{1}+e_{2} &\equiv &1 \\
d &=&1,e_{1}=u_{1},e_{2}=1 \\
a &=&1,b_{1}=l_{1},b_{2}=0
\end{eqnarray*}%
Since $\alpha \left( 1_{H};0,0,1,0\right) =e_{2}+\left( a+b_{1}+b_{2}\right)
\equiv 0$ and $\alpha \left( 1_{H};1,0,0,0\right) =0$ we get%
\begin{equation*}
\left[ B(x_{1}\otimes x_{1}x_{2};G,gx_{1}x_{2})+B(x_{1}\otimes
x_{1}x_{2};GX_{1},gx_{2})\right] G\otimes gx_{2}\otimes gx_{1}
\end{equation*}%
By considering also the right side, we obtain%
\begin{equation*}
\left[
\begin{array}{c}
-B(g\otimes x_{1}x_{2};G,gx_{2})+B(x_{1}\otimes x_{1}x_{2};G,gx_{1}x_{2}) \\
+B(x_{1}\otimes x_{1}x_{2};GX_{1},gx_{2})+B(x_{1}\otimes x_{2};G,gx_{2})%
\end{array}%
\right] G\otimes gx_{2}\otimes gx_{1}=0
\end{equation*}%
which is satisfied in view of the form of the element.

\subsubsection{Case $G\otimes g\otimes x_{1}x_{2}$}

\begin{eqnarray*}
&&\sum_{a,b_{1},b_{2},d,e_{1},e_{2}=0}^{1}\sum_{l_{1}=0}^{b_{1}}%
\sum_{l_{2}=0}^{b_{2}}\sum_{u_{1}=0}^{e_{1}}\sum_{u_{2}=0}^{e_{2}}\left(
-1\right) ^{\alpha \left( x_{1};l_{1},l_{2},u_{1},u_{2}\right) } \\
&&B(g\otimes
x_{1}x_{2};G^{a}X_{1}^{b_{1}}X_{2}^{b_{2}},g^{d}x_{1}^{e_{1}}x_{2}^{e_{2}})G^{a}X_{1}^{b_{1}-l_{1}}X_{2}^{b_{2}-l_{2}}\otimes g^{d}x_{1}^{e_{1}-u_{1}}x_{2}^{e_{2}-u_{2}}\otimes
\\
&&g^{a+b_{1}+b_{2}+l_{1}+l_{2}+d+e_{1}+e_{2}+u_{1}+u_{2}}x_{1}^{l_{1}+u_{1}+1}x_{2}^{l_{2}+u_{2}}
\end{eqnarray*}

First summand of the left side%
\begin{eqnarray*}
l_{1} &=&u_{1}=0,l_{2}+u_{2}=1 \\
a+b_{1}+b_{2}+d+e_{1}+e_{2} &\equiv &1 \\
d &=&1,e_{1}=0,e_{2}=u_{2} \\
a &=&1,b_{1}=0,b_{2}=l_{2}
\end{eqnarray*}%
Since $\alpha \left( x_{1};0,0,0,1\right) =1$ and $\alpha \left(
x_{1};0,1,0,0\right) =a+b_{1}+b_{2}+1=1$%
\begin{eqnarray*}
\alpha \left( x_{1};0,0,0,1\right) &=&1 \\
\alpha \left( x_{1};0,1,0,0\right) &=&a+b_{1}+b_{2}+1
\end{eqnarray*}%
\begin{equation*}
\left[ -B(g\otimes x_{1}x_{2};G,gx_{2})-B(g\otimes x_{1}x_{2};GX_{2},g)%
\right] G\otimes g\otimes x_{1}x_{2}
\end{equation*}%
Second summand of the left side%
\begin{eqnarray*}
l_{1}+u_{1} &=&1,l_{2}+u_{2}=1 \\
a+b_{1}+b_{2}+d+e_{1}+e_{2} &=&0 \\
d &=&1,e_{1}=u_{1},e_{2}=u_{2} \\
a &=&1,b_{1}=l_{1},b_{2}=l_{2}.
\end{eqnarray*}%
\begin{eqnarray*}
\alpha \left( 1_{H};0,0,1,1\right) &=&1+e_{2}\text{, }\alpha \left(
1_{H};0,1,1,0\right) =e_{2}+a+b_{1}+b_{2}+1, \\
\alpha \left( 1_{H};1,0,0,1\right) &=&a+b_{1}\text{ and }\alpha \left(
1_{H};1,1,0,0\right) =1+b_{2}
\end{eqnarray*}%
\begin{eqnarray*}
\alpha \left( 1_{H};0,0,1,1\right) &=&1+e_{2}\equiv 0\text{, }\alpha \left(
1_{H};0,1,1,0\right) =e_{2}+a+b_{1}+b_{2}+1\equiv 1, \\
\alpha \left( 1_{H};1,0,0,1\right) &=&a+b_{1}\equiv 0\text{ and }\alpha
\left( 1_{H};1,1,0,0\right) =1+b_{2}\equiv 0
\end{eqnarray*}%
\begin{equation*}
\left[
\begin{array}{c}
B(x_{1}\otimes x_{1}x_{2};G,gx_{1}x_{2})-B(x_{1}\otimes
x_{1}x_{2};GX_{2},gx_{1})+ \\
B(x_{1}\otimes x_{1}x_{2};GX_{1},gx_{2})+B(x_{1}\otimes
x_{1}x_{2};GX_{1}X_{2},g)%
\end{array}%
\right] G\otimes g\otimes x_{1}x_{2}
\end{equation*}%
By taking in account also the left side, we get%
\begin{eqnarray*}
&&-B(g\otimes x_{1}x_{2};G,gx_{2})-B(g\otimes
x_{1}x_{2};GX_{2},g)+B(x_{1}\otimes x_{1}x_{2};G,gx_{1}x_{2})+ \\
&&-B(x_{1}\otimes x_{1}x_{2};GX_{2},gx_{1})+B(x_{1}\otimes
x_{1}x_{2};GX_{1},gx_{2})+B(x_{1}\otimes x_{1}x_{2};GX_{1}X_{2},g) \\
&=&B(x_{1}\otimes 1_{H};G,g)
\end{eqnarray*}%
which is satisfied in view of the form of the elements.

\subsection{$B\left( x_{1}\otimes x_{1}x_{2};X_{1},gx_{1}x_{2}\right) $}

Now we look for the terms in $B\left( x_{1}\otimes
x_{1}x_{2};X_{1},gx_{1}x_{2}\right) $ in the Casimir equality and we get

\begin{equation*}
a=0,b_{1}=1,b_{2}=0,d=e_{1}=e_{2}=1
\end{equation*}%
\begin{eqnarray*}
&&B(x_{1}\otimes
x_{1}x_{2};G^{a}X_{1}^{b_{1}}X_{2}^{b_{2}},g^{d}x_{1}^{e_{1}}x_{2}^{e_{2}})G^{a}X_{1}^{b_{1}-l_{1}}X_{2}^{b_{2}-l_{2}}\otimes g^{d}x_{1}^{e_{1}-u_{1}}x_{2}^{e_{2}-u_{2}}\otimes
\\
&&g^{a+b_{1}+b_{2}+l_{1}+l_{2}+d+e_{1}+e_{2}+u_{1}+u_{2}}x_{1}^{l_{1}+u_{1}}x_{2}^{l_{2}+u_{2}}
\end{eqnarray*}

\begin{gather*}
\sum_{l_{1}=0}^{1}\sum_{u_{1}=0}^{1}\sum_{u_{2}=0}^{1}\left( -1\right)
^{\alpha \left( 1_{H};l_{1},0,u_{1},u_{2}\right) }B(x_{1}\otimes
x_{1}x_{2};X_{1},gx_{1}x_{2})X_{1}^{1-l_{1}} \\
\otimes gx_{1}^{1-u_{1}}x_{2}^{1-u_{2}}\otimes
g^{l_{1}+u_{1}+u_{2}}x_{1}^{l_{1}+u_{1}}x_{2}^{l_{2}+u_{2}} \\
=\left( -1\right) ^{\alpha \left( 1_{H};0,0,0,0\right) }B(x_{1}\otimes
x_{1}x_{2};X_{1},gx_{1}x_{2})X_{1}\otimes gx_{1}x_{2}\otimes 1_{H} \\
+\left( -1\right) ^{\alpha \left( 1_{H};0,0,0,1\right) }B(x_{1}\otimes
x_{1}x_{2};X_{1},gx_{1}x_{2})X_{1}\otimes gx_{1}\otimes gx_{2} \\
+\left( -1\right) ^{\alpha \left( 1_{H};0,0,1,0\right) }B(x_{1}\otimes
x_{1}x_{2};X_{1},gx_{1}x_{2})X_{1}\otimes gx_{2}\otimes gx_{1} \\
+\left( -1\right) ^{\alpha \left( 1_{H};0,0,1,1\right) }B(x_{1}\otimes
x_{1}x_{2};X_{1},gx_{1}x_{2})X_{1}\otimes g\otimes x_{1}x_{2} \\
+\left( -1\right) ^{\alpha \left( 1_{H};1,0,0,0\right) }B(x_{1}\otimes
x_{1}x_{2};X_{1},gx_{1}x_{2})1_{A}\otimes gx_{1}x_{2}\otimes gx_{1} \\
+\left( -1\right) ^{\alpha \left( 1_{H};1,0,1,0\right) }B(x_{1}\otimes
x_{1}x_{2};X_{1},gx_{1}x_{2})1_{A}\otimes
gx_{1}^{1-u_{1}}x_{2}^{1-u_{2}}\otimes
g^{u_{1}+u_{2}}x_{1}^{1+1}x_{2}^{u_{2}}=0 \\
+\left( -1\right) ^{\alpha \left( 1_{H};1,0,0,1\right) }B(x_{1}\otimes
x_{1}x_{2};X_{1},gx_{1}x_{2})1_{A}\otimes gx_{1}\otimes x_{1}x_{2} \\
+\left( -1\right) ^{\alpha \left( 1_{H};1,0,1,1\right) }B(x_{1}\otimes
x_{1}x_{2};X_{1},gx_{1}x_{2})1_{A}\otimes
gx_{1}^{1-u_{1}}x_{2}^{1-u_{2}}\otimes
g^{u_{1}+u_{2}}x_{1}^{1+u_{1}}x_{2}^{u_{2}}=0
\end{gather*}

\subsubsection{Case $X_{1}\otimes gx_{1}x_{2}\otimes 1_{H}$}

\begin{eqnarray*}
&&+\sum_{a,b_{1},b_{2},d,e_{1},e_{2}=0}^{1}\sum_{l_{1}=0}^{b_{1}}%
\sum_{l_{2}=0}^{b_{2}}\sum_{u_{1}=0}^{e_{1}}\sum_{u_{2}=0}^{e_{2}}\left(
-1\right) ^{\alpha \left( 1_{H};l_{1},l_{2},u_{1},u_{2}\right) } \\
&&B(x_{1}\otimes
x_{1}x_{2};G^{a}X_{1}^{b_{1}}X_{2}^{b_{2}},g^{d}x_{1}^{e_{1}}x_{2}^{e_{2}})G^{a}X_{1}^{b_{1}-l_{1}}X_{2}^{b_{2}-l_{2}}\otimes g^{d}x_{1}^{e_{1}-u_{1}}x_{2}^{e_{2}-u_{2}}\otimes
\\
&&g^{a+b_{1}+b_{2}+l_{1}+l_{2}+d+e_{1}+e_{2}+u_{1}+u_{2}}x_{1}^{l_{1}+u_{1}}x_{2}^{l_{2}+u_{2}}
\end{eqnarray*}

Nothing from the first summand of the left side. The second summand of the
left side gives us%
\begin{eqnarray*}
l_{1} &=&u_{1}=l_{2}=u_{2}=0 \\
a+b_{1}+b_{2}+d+e_{1}+e_{2} &\equiv &0 \\
d &=&1,e_{1}=1,e_{2}=1 \\
a &=&0,b_{1}=1,b_{2}=0
\end{eqnarray*}%
and we get%
\begin{equation*}
B(x_{1}\otimes x_{1}x_{2};X_{1},gx_{1}x_{2})X_{1}\otimes gx_{1}x_{2}\otimes
1_{H}
\end{equation*}%
By taking care also of the left side we get%
\begin{equation*}
B(x_{1}\otimes x_{1}x_{2};X_{1},gx_{1}x_{2})=B(x_{1}\otimes
x_{1}x_{2};X_{1},gx_{1}x_{2})
\end{equation*}%
so that we get no new information.

\subsubsection{Case $X_{1}\otimes gx_{1}\otimes gx_{2}$}

\bigskip

There is nothing from the first summand of the left side while the second
summand gives us%
\begin{eqnarray*}
l_{1} &=&u_{1}=0,l_{2}+u_{2}=1 \\
a+b_{1}+b_{2}+d+e_{1}+e_{2} &\equiv &0 \\
d &=&1,e_{1}=1,e_{2}=u_{2} \\
a &=&0,b_{1}=1,b_{2}=l_{2}
\end{eqnarray*}%
Since $\alpha \left( 1_{H};0,0,0,1\right) =1$ and $\alpha \left(
1_{H};0,1,0,0\right) =0$ we get%
\begin{equation*}
\left[ -B(x_{1}\otimes x_{1}x_{2};X_{1},gx_{1}x_{2})+B(x_{1}\otimes
x_{1}x_{2};X_{1}X_{2},gx_{1})\right] X_{1}\otimes gx_{1}\otimes gx_{2}
\end{equation*}

By taking in account also the right side we get%
\begin{equation*}
-B(x_{1}\otimes x_{1}x_{2};X_{1},gx_{1}x_{2})+B(x_{1}\otimes
x_{1}x_{2};X_{1}X_{2},gx_{1})=B(x_{1}\otimes x_{1};X_{1},gx_{1})
\end{equation*}%
which holds in view of the form of the elements.

\subsubsection{Case $X_{1}\otimes gx_{2}\otimes gx_{1}$}

From the fist summand of the left side we get%
\begin{eqnarray*}
l_{1} &=&u_{1}=l_{2}=u_{2}=0 \\
a+b_{1}+b_{2}+d+e_{1}+e_{2} &\equiv &1 \\
d &=&1,e_{1}=0,e_{2}=1 \\
a &=&0,b_{1}=1,b_{2}=0
\end{eqnarray*}%
Since $\alpha \left( x_{1};0,0,0,0\right) =a+b_{1}+b_{2}=1$we obtain%
\begin{equation*}
-B(g\otimes x_{1}x_{2};X_{1},gx_{2})X_{1}\otimes gx_{2}\otimes gx_{1}.
\end{equation*}%
From the second summand of the left side we get

\begin{eqnarray*}
l_{1}+u_{1} &=&1,l_{2}=u_{2}=0 \\
a+b_{1}+b_{2}+d+e_{1}+e_{2} &\equiv &0 \\
d &=&1,e_{1}=u_{1},e_{2}=1 \\
a &=&0,b_{1}-l_{1}=1\Rightarrow b_{1}=1\text{ and }l_{1}=0\Rightarrow
u_{1}=1=e_{1},b_{2}=0
\end{eqnarray*}%
Since $\alpha \left( 1_{H};0,0,1,0\right) =e_{2}+\left( a+b_{1}+b_{2}\right)
=1+1\equiv 0$ we get%
\begin{equation*}
B(x_{1}\otimes x_{1}x_{2};X_{1},gx_{1}x_{2})X_{1}\otimes gx_{2}\otimes
gx_{1}.
\end{equation*}%
By taking in account also the right side, we get%
\begin{equation*}
-B(g\otimes x_{1}x_{2};X_{1},gx_{2})+B(x_{1}\otimes
x_{1}x_{2};X_{1},gx_{1}x_{2})+B(x_{1}\otimes x_{2},X_{1},gx_{2})=0
\end{equation*}

which holds in view of the form of the elements.

\subsubsection{Case $X_{1}\otimes g\otimes x_{1}x_{2}$}

From the fist summand of the left side we get%
\begin{eqnarray*}
l_{1} &=&u_{1}=0,l_{2}+u_{2}=1 \\
a+b_{1}+b_{2}+d+e_{1}+e_{2} &\equiv &1 \\
d &=&1,e_{1}=0,e_{2}=u_{2} \\
a &=&0,b_{1}=1,b_{2}=l_{2}
\end{eqnarray*}%
Since $\alpha \left( x_{1};0,0,0,1\right) =1$ and $\alpha \left(
x_{1};0,1,0,0\right) =\left( a+b_{1}+b_{2}+1\right) \equiv 1$ we obtain
\begin{equation*}
\left[ -B(g\otimes x_{1}x_{2};X_{1},gx_{2})-B(g\otimes
x_{1}x_{2};X_{1}X_{2},g)\right] X_{1}\otimes g\otimes x_{1}x_{2}
\end{equation*}%
From the second summand of the left side we get%
\begin{eqnarray*}
l_{1}+u_{1} &=&1,l_{2}+u_{2}=1 \\
a+b_{1}+b_{2}+d+e_{1}+e_{2} &\equiv &0 \\
d &=&1,e_{1}=u_{1},e_{2}=u_{2} \\
a &=&0,b_{1}-l_{1}=1\Rightarrow b_{1}=1,l_{1}=0,u_{1}=1,e_{1}=1,b_{2}=l_{2}
\end{eqnarray*}%
Since $\alpha \left( 1_{H};0,0,1,1\right) =$ $1+e_{2}\equiv 0$ and $\alpha
\left( 1_{H};0,1,1,0\right) =0+\left( 0+1+1+1\right) \equiv 1,$ we obtain%
\begin{equation*}
\left[ B(x_{1}\otimes x_{1}x_{2};X_{1},gx_{1}x_{2})-B(x_{1}\otimes
x_{1}x_{2};X_{1}X_{2},gx_{1})\right] X_{1}\otimes g\otimes x_{1}x_{2}
\end{equation*}

By considering also the right side, we get%
\begin{gather*}
-B(g\otimes x_{1}x_{2};X_{1},gx_{2})-B(g\otimes
x_{1}x_{2};X_{1}X_{2},g)+B(x_{1}\otimes x_{1}x_{2};X_{1},gx_{1}x_{2}) \\
-B(x_{1}\otimes x_{1}x_{2};X_{1}X_{2},gx_{1})=B(x_{1}\otimes 1_{H};X_{1},g)
\end{gather*}%
which holds in view of the form of the elements.

\subsubsection{Case $1_{A}\otimes gx_{1}x_{2}\otimes gx_{1}$}

From the first summand of the left side we get%
\begin{eqnarray*}
l_{1} &=&u_{1}=l_{2}=u_{2}=0 \\
a+b_{1}+b_{2}+d+e_{1}+e_{2} &\equiv &1 \\
d &=&1,e_{1}=1,e_{2}=1 \\
a &=&0,b_{1}=0,b_{2}=0
\end{eqnarray*}%
Since $\alpha \left( x_{1};0,0,0,0\right) $ $=a+b_{1}+b_{2}=0,$ we obtain%
\begin{equation*}
B(g\otimes x_{1}x_{2};1_{A},gx_{1}x_{2})1_{A}\otimes gx_{1}x_{2}\otimes
gx_{1}.
\end{equation*}%
From the second summand of the left side we get%
\begin{eqnarray*}
l_{1}+u_{1} &=&1,l_{2}=u_{2}=0 \\
a+b_{1}+b_{2}+d+e_{1}+e_{2} &\equiv &0 \\
d &=&1,e_{1}-u_{1}=1\Rightarrow e_{1}=1,u_{1}=0,l_{1}=1=b_{1},e_{2}=1 \\
a &=&0,b_{2}=0
\end{eqnarray*}%
Since $\alpha \left( 1_{H};1,0,0,0\right) =b_{2}=0,$ we obtain%
\begin{equation*}
B(x_{1}\otimes x_{1}x_{2};X_{1},gx_{1}x_{2})1_{A}\otimes gx_{1}x_{2}\otimes
gx_{1}
\end{equation*}%
By considering also the left side, we get%
\begin{equation*}
B(g\otimes x_{1}x_{2};1_{A},gx_{1}x_{2})+B(x_{1}\otimes
x_{1}x_{2};X_{1},gx_{1}x_{2})+B(x_{1}\otimes x_{2};1_{A},gx_{1}x_{2})=0
\end{equation*}%
which holds in view of the form of the elements.

\subsubsection{Case $1_{A}\otimes gx_{1}\otimes x_{1}x_{2}$}

From the first summand of the left side we get%
\begin{eqnarray*}
l_{1} &=&u_{1}=0,l_{2}+u_{2}=1 \\
a+b_{1}+b_{2}+d+e_{1}+e_{2} &\equiv &1 \\
d &=&1,e_{1}=1,e_{2}=u_{2} \\
a &=&0,b_{1}=0,b_{2}=l_{2}
\end{eqnarray*}%
Since $\alpha \left( x_{1};0,0,0,1\right) =1$ and $\alpha \left(
x_{1};0,1,0,0\right) =a+b_{1}+b_{2}+1=0$%
\begin{equation*}
\left[ -B(g\otimes x_{1}x_{2};1_{A},gx_{1}x_{2})+B(g\otimes
x_{1}x_{2};X_{2},gx_{1})\right] 1_{A}\otimes gx_{1}\otimes x_{1}x_{2}
\end{equation*}%
From the second summand of the left side we get%
\begin{eqnarray*}
l_{1}+u_{1} &=&1,l_{2}+u_{2}=1 \\
a+b_{1}+b_{2}+d+e_{1}+e_{2} &\equiv &1 \\
d &=&1,e_{1}-u_{1}=1\Rightarrow e_{1}=1,u_{1}=0,l_{1}=1,e_{2}=u_{2} \\
a &=&0,b_{1}=1,b_{2}=l_{2}
\end{eqnarray*}%
Since $\alpha \left( 1_{H};1,0,0,1\right) $ $=a+b_{1}=1$ and $\alpha \left(
1_{H};1,1,0,0\right) =1+b_{2}=0,$ we obtain%
\begin{equation*}
\left[ -B(x_{1}\otimes x_{1}x_{2};X_{1},gx_{1}x_{2})+B(x_{1}\otimes
x_{1}x_{2};X_{1}X_{2},gx_{1})\right] 1_{A}\otimes gx_{1}\otimes x_{1}x_{2}.
\end{equation*}%
By considering also the right side, we get%
\begin{gather*}
-B(g\otimes x_{1}x_{2};1_{A},gx_{1}x_{2})+B(g\otimes
x_{1}x_{2};X_{2},gx_{1})-B(x_{1}\otimes x_{1}x_{2};X_{1},gx_{1}x_{2}) \\
+B(x_{1}\otimes x_{1}x_{2};X_{1}X_{2},gx_{1})=B(x_{1}\otimes
1_{H};1_{A},gx_{1})
\end{gather*}%
which holds in view of the form of the elements.

\subsection{$B(x_{1}\otimes x_{1}x_{2};X_{2},gx_{1}x_{2})$}

\begin{equation*}
a=0,b_{1}=0,b_{2}=1,d=1,e_{1}=1,e_{2}=1
\end{equation*}

and we get%
\begin{multline*}
+\sum_{l_{2}=0}^{1}\sum_{u_{1}=0}^{1}\sum_{u_{2}=0}^{1}\left( -1\right)
^{\alpha \left( 1_{H};0,l_{2},u_{1},u_{2}\right) }B(x_{1}\otimes
x_{1}x_{2};X_{2},gx_{1}x_{2})X_{2}^{1-l_{2}} \\
\otimes gx_{1}^{1-u_{1}}x_{2}^{1-u_{2}}\otimes
g^{l_{2}+u_{1}+u_{2}}x_{1}^{u_{1}}x_{2}^{l_{2}+u_{2}}= \\
+\left( -1\right) ^{\alpha \left( 1_{H};0,0,0,0\right) }B(x_{1}\otimes
x_{1}x_{2};X_{2},gx_{1}x_{2})X_{2}\otimes gx_{1}x_{2}\otimes 1_{H}+ \\
+\left( -1\right) ^{\alpha \left( 1_{H};0,0,0,1\right) }B(x_{1}\otimes
x_{1}x_{2};X_{2},gx_{1}x_{2})X_{2}\otimes gx_{1}\otimes gx_{2}+ \\
+\left( -1\right) ^{\alpha \left( 1_{H};0,0,1,0\right) }B(x_{1}\otimes
x_{1}x_{2};X_{2},gx_{1}x_{2})X_{2}\otimes gx_{2}\otimes gx_{1}+ \\
+\left( -1\right) ^{\alpha \left( 1_{H};0,0,1,1\right) }B(x_{1}\otimes
x_{1}x_{2};X_{2},gx_{1}x_{2})X_{2}\otimes g\otimes x_{1}x_{2}+ \\
+\left( -1\right) ^{\alpha \left( 1_{H};0,1,0,0\right) }B(x_{1}\otimes
x_{1}x_{2};X_{2},gx_{1}x_{2})1_{A}\otimes gx_{1}x_{2}\otimes gx_{2}+ \\
+\left( -1\right) ^{\alpha \left( 1_{H};0,1,0,1\right) }B(x_{1}\otimes
x_{1}x_{2};X_{2},gx_{1}x_{2})1_{A}\otimes
gx_{1}^{1-u_{1}}x_{2}^{1-u_{2}}\otimes
g^{l_{2}+u_{1}+u_{2}}x_{1}^{u_{1}}x_{2}^{1+1}+=0 \\
+\left( -1\right) ^{\alpha \left( 1_{H};0,1,1,0\right) }B(x_{1}\otimes
x_{1}x_{2};X_{2},gx_{1}x_{2})1_{A}\otimes gx_{2}\otimes x_{1}x_{2}+\text{%
doing} \\
+\left( -1\right) ^{\alpha \left( 1_{H};0,11,1\right) }B(x_{1}\otimes
x_{1}x_{2};X_{2},gx_{1}x_{2})1_{A}\otimes
gx_{1}^{1-u_{1}}x_{2}^{1-u_{2}}\otimes
g^{l_{2}+u_{1}+u_{2}}x_{1}^{u_{1}}x_{2}^{1+1}=0
\end{multline*}

\subsubsection{Case $X_{2}\otimes gx_{1}x_{2}\otimes 1_{H}$}

Nothing from the first summand of the left side. From the second summand of
the left side we get

\begin{eqnarray*}
l_{1} &=&u_{1}=l_{2}=u_{2}=0 \\
a+b_{1}+b_{2}+d+e_{1}+e_{2} &\equiv &0 \\
d &=&e_{1}=e_{2}=1 \\
a &=&b_{1}=0,b_{2}=1
\end{eqnarray*}%
and we get%
\begin{equation*}
B(x_{1}\otimes x_{1}x_{2};X_{2},gx_{1}x_{2})X_{2}\otimes gx_{1}x_{2}\otimes
1_{H}.
\end{equation*}%
By considering also the right side we obtain%
\begin{equation*}
B(x_{1}\otimes x_{1}x_{2};X_{2},gx_{1}x_{2})=B(x_{1}\otimes
x_{1}x_{2};X_{2},gx_{1}x_{2})
\end{equation*}%
and thus we get nothing new.

\subsubsection{Case $X_{2}\otimes gx_{1}\otimes gx_{2}$}

Nothing from the first summand of the left side. Second summand of the left
side gives us

\begin{eqnarray*}
l_{1} &=&u_{1}=0,l_{2}+u_{2}=1 \\
a+b_{1}+b_{2}+d+e_{1}+e_{2} &\equiv &0 \\
d &=&e_{1}=1,e_{2}=u_{2} \\
a &=&b_{1}=0,b_{2}-l_{2}=1\Rightarrow b_{2}=1,l_{2}=0,e_{2}=u_{2}=1
\end{eqnarray*}%
Since $\alpha \left( 1_{H};0,0,0,1\right) =a+b_{1}+b_{2}=1$ we get%
\begin{equation*}
-B(x_{1}\otimes x_{1}x_{2};X_{2},gx_{1}x_{2})X_{2}\otimes gx_{1}\otimes
gx_{2}.
\end{equation*}%
By considering also the right side we get%
\begin{equation*}
-B(x_{1}\otimes x_{1}x_{2};X_{2},gx_{1}x_{2})=B(x_{1}\otimes
x_{1};X_{2},gx_{1})
\end{equation*}%
\begin{eqnarray*}
-B(x_{1}\otimes x_{1}x_{2};X_{2},gx_{1}x_{2}) &=&-2B(x_{1}\otimes
1_{H};1_{A},gx_{2}) \\
&=&B(x_{1}\otimes x_{1};X_{2},gx_{1})=2B(x_{1}\otimes 1_{H};1_{A},gx_{2})
\end{eqnarray*}

\subsubsection{Case $X_{2}\otimes gx_{2}\otimes gx_{1}$}

From the first summand of the left side we get%
\begin{eqnarray*}
l_{1} &=&u_{1}=0=l_{2}=u_{2}=0 \\
a+b_{1}+b_{2}+d+e_{1}+e_{2} &\equiv &1 \\
d &=&e_{2}=1,e_{1}=0 \\
a &=&b_{1}=0,b_{2}=1.
\end{eqnarray*}%
Since $\alpha \left( x_{1};0,0,0,0\right) =a+b_{1}+b_{2}=1,$ we get%
\begin{equation*}
-B(g\otimes x_{1}x_{2};X_{2},gx_{2})X_{2}\otimes gx_{2}\otimes gx_{1}.
\end{equation*}%
From the first summand of the left side we get%
\begin{eqnarray*}
l_{1}+u_{1} &=&1,l_{2}=u_{2}=0 \\
a+b_{1}+b_{2}+d+e_{1}+e_{2} &\equiv &0 \\
d &=&e_{2}=1,e_{1}=u_{1} \\
a &=&0,b_{1}=l_{1},b_{2}=1
\end{eqnarray*}%
Since $\alpha \left( 1_{H};0,0,1,0\right) =e_{2}+\left( a+b_{1}+b_{2}\right)
\equiv 0$ and $\alpha \left( 1_{H};1,0,0,0\right) =b_{2}=1$%
\begin{equation*}
\left[ B(x_{1}\otimes x_{1}x_{2};X_{2},gx_{1}x_{2})-B(x_{1}\otimes
x_{1}x_{2};X_{1}X_{2},gx_{2})\right] X_{2}\otimes gx_{2}\otimes gx_{1}
\end{equation*}%
By considering also the right side, we obtain%
\begin{gather*}
-B(g\otimes x_{1}x_{2};X_{2},gx_{2})+B(x_{1}\otimes
x_{1}x_{2};X_{2},gx_{1}x_{2}) \\
-B(x_{1}\otimes x_{1}x_{2};X_{1}X_{2},gx_{2})+B(x_{1}\otimes
x_{2};X_{2},gx_{2})=0
\end{gather*}%
which holds in view of the form of the elements.

\subsubsection{Case $X_{2}\otimes g\otimes x_{1}x_{2}$}

From the first summand of the left side we get%
\begin{eqnarray*}
l_{1} &=&u_{1}=0,l_{2}+u_{2}=1 \\
a+b_{1}+b_{2}+d+e_{1}+e_{2} &\equiv &1 \\
d &=&1,e_{1}=0,e_{2}=u_{2} \\
a &=&b_{1}=0,b_{2}-l_{2}=1\Rightarrow b_{2}=1,l_{2}=0\Rightarrow
e_{2}=u_{2}=1.
\end{eqnarray*}%
Since $\alpha \left( x_{1};0,0,0,1\right) =1,$ we get%
\begin{equation*}
-B(g\otimes x_{1}x_{2};X_{2},gx_{2})X_{2}\otimes g\otimes x_{1}x_{2}.
\end{equation*}%
From the second summand of the left side we get%
\begin{eqnarray*}
l_{1}+u_{1} &=&1,l_{2}+u_{2}=1 \\
a+b_{1}+b_{2}+d+e_{1}+e_{2} &\equiv &0 \\
d &=&1,e_{1}=u_{1},e_{2}=u_{2} \\
a &=&0,b_{1}=l_{1},b_{2}-l_{2}=1\Rightarrow b_{2}=1,l_{2}=0\Rightarrow
e_{2}=u_{2}=1.
\end{eqnarray*}%
Since $\alpha \left( 1_{H};0,0,1,1\right) =1+e_{2}\equiv 0$ and $\alpha
\left( 1_{H};1,0,0,1\right) =a+b_{1}=1$%
\begin{equation*}
\left[ B(x_{1}\otimes x_{1}x_{2};X_{2},gx_{1}x_{2})-B(x_{1}\otimes
x_{1}x_{2};X_{1}X_{2},gx_{2})\right] )X_{2}\otimes g\otimes x_{1}x_{2}.
\end{equation*}%
By considering also the right side, we obtain%
\begin{equation*}
-B(g\otimes x_{1}x_{2};X_{2},gx_{2})+B(x_{1}\otimes
x_{1}x_{2};X_{2},gx_{1}x_{2})-B(x_{1}\otimes
x_{1}x_{2};X_{1}X_{2},gx_{2})-B(x_{1}\otimes 1_{H},X_{2},g)=0
\end{equation*}%
which holds in view of the form of the element.

\subsubsection{Case $1_{A}\otimes gx_{1}x_{2}\otimes gx_{2}$}

Nothing from the first summand of the left side. From the second summand we
get%
\begin{eqnarray*}
l_{1} &=&u_{1}=0,l_{2}+u_{2}=1 \\
a+b_{1}+b_{2}+d+e_{1}+e_{2} &\equiv &0 \\
d &=&1,e_{1}=1,e_{2}-u_{2}=1\Rightarrow e_{2}=1,u_{2}=0,l_{2}=1 \\
a &=&b_{1}=0,b_{2}=l_{2}=1.
\end{eqnarray*}%
Since $\alpha \left( 1_{H};0,1,0,0\right) =0,$ we obtain%
\begin{equation*}
B(x_{1}\otimes x_{1}x_{2};X_{2},gx_{1}x_{2})1_{A}\otimes gx_{1}x_{2}\otimes
gx_{2}.
\end{equation*}%
By considering also the right side, we obtain%
\begin{equation*}
B(x_{1}\otimes x_{1}x_{2};X_{2},gx_{1}x_{2})-B(x_{1}\otimes
x_{1};1_{A},gx_{1}x_{2})=0
\end{equation*}%
\begin{eqnarray*}
B(x_{1}\otimes x_{1}x_{2};X_{2},gx_{1}x_{2}) &=&-2B(x_{1}\otimes
1_{H};1_{A},gx_{2}) \\
-B(x_{1}\otimes x_{1};1_{A},gx_{1}x_{2}) &=&+2B\left( x_{1}\otimes
1_{H};1_{A},gx_{2}\right)
\end{eqnarray*}%
which holds in view of the form of the elements.

\subsubsection{Case $1_{A}\otimes gx_{2}\otimes x_{1}x_{2}$}

From the first summand of the left side we get%
\begin{eqnarray*}
&&\sum_{a,b_{1},b_{2},d,e_{1},e_{2}=0}^{1}\sum_{l_{1}=0}^{b_{1}}%
\sum_{l_{2}=0}^{b_{2}}\sum_{u_{1}=0}^{e_{1}}\sum_{u_{2}=0}^{e_{2}}\left(
-1\right) ^{\alpha \left( x_{1};l_{1},l_{2},u_{1},u_{2}\right) } \\
&&B(g\otimes
x_{1}x_{2};G^{a}X_{1}^{b_{1}}X_{2}^{b_{2}},g^{d}x_{1}^{e_{1}}x_{2}^{e_{2}})G^{a}X_{1}^{b_{1}-l_{1}}X_{2}^{b_{2}-l_{2}}\otimes g^{d}x_{1}^{e_{1}-u_{1}}x_{2}^{e_{2}-u_{2}}\otimes
\\
&&g^{a+b_{1}+b_{2}+l_{1}+l_{2}+d+e_{1}+e_{2}+u_{1}+u_{2}}x_{1}^{l_{1}+u_{1}+1}x_{2}^{l_{2}+u_{2}}
\end{eqnarray*}%
\begin{eqnarray*}
l_{1} &=&u_{1}=0,l_{2}+u_{2}=1 \\
a+b_{1}+b_{2}+d+e_{1}+e_{2} &\equiv &1 \\
d &=&1,e_{1}=0,e_{2}-u_{2}=1\Rightarrow e_{2}=1,u_{2}=0,l_{2}=1 \\
a &=&b_{1}=0,b_{2}=l_{2}=1.
\end{eqnarray*}%
Since $\alpha \left( x_{1};0,1,0,0\right) =a+b_{1}+b_{2}+1\equiv 0,$ we get%
\begin{equation*}
B(g\otimes x_{1}x_{2};X_{2},gx_{2})1_{A}\otimes gx_{2}\otimes x_{1}x_{2}.
\end{equation*}%
From the second summand, we obtain%
\begin{eqnarray*}
&&+\sum_{a,b_{1},b_{2},d,e_{1},e_{2}=0}^{1}\sum_{l_{1}=0}^{b_{1}}%
\sum_{l_{2}=0}^{b_{2}}\sum_{u_{1}=0}^{e_{1}}\sum_{u_{2}=0}^{e_{2}}\left(
-1\right) ^{\alpha \left( 1_{H};l_{1},l_{2},u_{1},u_{2}\right) } \\
&&B(x_{1}\otimes
x_{1}x_{2};G^{a}X_{1}^{b_{1}}X_{2}^{b_{2}},g^{d}x_{1}^{e_{1}}x_{2}^{e_{2}})G^{a}X_{1}^{b_{1}-l_{1}}X_{2}^{b_{2}-l_{2}}\otimes g^{d}x_{1}^{e_{1}-u_{1}}x_{2}^{e_{2}-u_{2}}\otimes
\\
&&g^{a+b_{1}+b_{2}+l_{1}+l_{2}+d+e_{1}+e_{2}+u_{1}+u_{2}}x_{1}^{l_{1}+u_{1}}x_{2}^{l_{2}+u_{2}}
\end{eqnarray*}%
\begin{eqnarray*}
l_{1}+u_{1} &=&1,l_{2}+u_{2}=1 \\
a+b_{1}+b_{2}+d+e_{1}+e_{2} &\equiv &0 \\
d &=&1,e_{1}=u_{1},e_{2}-u_{2}=1\Rightarrow e_{2}=1,u_{2}=0,l_{2}=1 \\
a &=&0,b_{1}=l_{1},b_{2}=l_{2}=1.
\end{eqnarray*}%
Since $\alpha \left( 1_{H};0,1,1,0\right) =1+0+0+1+1\equiv 1$ and $\alpha
\left( 1_{H};1,1,0,0\right) =1+b_{2}=1+1\equiv 0,$ w obtain%
\begin{equation*}
\left[ -B(x_{1}\otimes x_{1}x_{2};X_{2},gx_{1}x_{2})+B(x_{1}\otimes
x_{1}x_{2};X_{1}X_{2},gx_{2})\right] 1_{A}\otimes gx_{2}\otimes x_{1}x_{2}
\end{equation*}%
By considering also the right side, we get%
\begin{equation*}
B(g\otimes x_{1}x_{2};X_{2},gx_{2})-B(x_{1}\otimes
x_{1}x_{2};X_{2},gx_{1}x_{2})+B(x_{1}\otimes
x_{1}x_{2};X_{1}X_{2},gx_{2})-B(x_{1}\otimes 1_{H};1_{A},gx_{2})=0
\end{equation*}%
\begin{eqnarray*}
B(g\otimes x_{1}x_{2};X_{2},gx_{2}) &=&-B(x_{1}\otimes 1_{H};1_{A},gx_{2}) \\
-B(x_{1}\otimes x_{1}x_{2};X_{2},gx_{1}x_{2} &=&+2B(x_{1}\otimes
1_{H};1_{A},gx_{2}) \\
+B(x_{1}\otimes x_{1}x_{2};X_{1}X_{2},gx_{2}) &=&0 \\
-B(x_{1}\otimes 1_{H};1_{A},gx_{2}) &=&-B(x_{1}\otimes 1_{H};1_{A},gx_{2})
\end{eqnarray*}%
which holds in view of the form of the elements.

\subsection{$B(x_{1}\otimes x_{1}x_{2};X_{1}X_{2},gx_{1})$}

\begin{eqnarray*}
&&+\sum_{a,b_{1},b_{2},d,e_{1},e_{2}=0}^{1}\sum_{l_{1}=0}^{b_{1}}%
\sum_{l_{2}=0}^{b_{2}}\sum_{u_{1}=0}^{e_{1}}\sum_{u_{2}=0}^{e_{2}}\left(
-1\right) ^{\alpha \left( 1_{H};l_{1},l_{2},u_{1},u_{2}\right) } \\
&&B(x_{1}\otimes
x_{1}x_{2};G^{a}X_{1}^{b_{1}}X_{2}^{b_{2}},g^{d}x_{1}^{e_{1}}x_{2}^{e_{2}})G^{a}X_{1}^{b_{1}-l_{1}}X_{2}^{b_{2}-l_{2}}\otimes g^{d}x_{1}^{e_{1}-u_{1}}x_{2}^{e_{2}-u_{2}}\otimes
\\
&&g^{a+b_{1}+b_{2}+l_{1}+l_{2}+d+e_{1}+e_{2}+u_{1}+u_{2}}x_{1}^{l_{1}+u_{1}}x_{2}^{l_{2}+u_{2}}
\end{eqnarray*}

\begin{equation*}
a=0,b_{1}=1,b_{2}=1,d=1,e_{1}=1,e_{2}=0
\end{equation*}%
and we get%
\begin{multline*}
\sum_{l_{1}=0}^{1}\sum_{l_{2}=0}^{1}\sum_{u_{1}=0}^{1}\left( -1\right)
^{\alpha \left( 1_{H};l_{1},l_{2},u_{1},0\right) }B(x_{1}\otimes
x_{1}x_{2};X_{1}X_{2},gx_{1})X_{1}^{1-l_{1}}X_{2}^{1-l_{2}} \\
\otimes gx_{1}^{1-u_{1}}\otimes
g^{l_{1}+l_{2}+u_{1}}x_{1}^{l_{1}+u_{1}}x_{2}^{l_{2}} \\
=\left( -1\right) ^{\alpha \left( 1_{H};0,0,0,0\right) }B(x_{1}\otimes
x_{1}x_{2};X_{1}X_{2},gx_{1})X_{1}X_{2}\otimes gx_{1}\otimes 1_{H} \\
+\left( -1\right) ^{\alpha \left( 1_{H};0,0,1,0\right) }B(x_{1}\otimes
x_{1}x_{2};X_{1}X_{2},gx_{1})X_{1}X_{2}\otimes g\otimes gx_{1} \\
+\left( -1\right) ^{\alpha \left( 1_{H};0,1,0,0\right) }B(x_{1}\otimes
x_{1}x_{2};X_{1}X_{2},gx_{1})X_{1}\otimes gx_{1}\otimes gx_{2} \\
+\left( -1\right) ^{\alpha \left( 1_{H};0,1,1,0\right) }B(x_{1}\otimes
x_{1}x_{2};X_{1}X_{2},gx_{1})X_{1}\otimes g\otimes x_{1}x_{2} \\
+\left( -1\right) ^{\alpha \left( 1_{H};1,0,0,0\right) }B(x_{1}\otimes
x_{1}x_{2};X_{1}X_{2},gx_{1})X_{2}\otimes gx_{1}\otimes gx_{1} \\
+\left( -1\right) ^{\alpha \left( 1_{H};1,0,1,0\right) }B(x_{1}\otimes
x_{1}x_{2};X_{1}X_{2},gx_{1})X_{2}\otimes gx_{1}^{1-u_{1}}\otimes
g^{1++u_{1}}x_{1}^{1+1}=0 \\
+\left( -1\right) ^{\alpha \left( 1_{H};1,1,0,0\right) }B(x_{1}\otimes
x_{1}x_{2};X_{1}X_{2},gx_{1})1_{H}\otimes gx_{1}\otimes x_{1}x_{2} \\
+\left( -1\right) ^{\alpha \left( 1_{H};1,1,1,0\right) }B(x_{1}\otimes
x_{1}x_{2};X_{1}X_{2},gx_{1})X_{2}^{1-l_{2}}\otimes gx_{1}^{1-u_{1}}\otimes
g^{1+l_{2}+u_{1}}x_{1}^{1+1}x_{2}^{l_{2}}=0
\end{multline*}

\subsubsection{Case $X_{1}X_{2}\otimes gx_{1}\otimes 1_{H}$}

Nothing from the first summand of the left side. The second summand gives us%
\begin{eqnarray*}
&&+\sum_{a,b_{1},b_{2},d,e_{1},e_{2}=0}^{1}\sum_{l_{1}=0}^{b_{1}}%
\sum_{l_{2}=0}^{b_{2}}\sum_{u_{1}=0}^{e_{1}}\sum_{u_{2}=0}^{e_{2}}\left(
-1\right) ^{\alpha \left( 1_{H};l_{1},l_{2},u_{1},u_{2}\right) } \\
&&B(x_{1}\otimes
x_{1}x_{2};G^{a}X_{1}^{b_{1}}X_{2}^{b_{2}},g^{d}x_{1}^{e_{1}}x_{2}^{e_{2}})G^{a}X_{1}^{b_{1}-l_{1}}X_{2}^{b_{2}-l_{2}}\otimes g^{d}x_{1}^{e_{1}-u_{1}}x_{2}^{e_{2}-u_{2}}\otimes
\\
&&g^{a+b_{1}+b_{2}+l_{1}+l_{2}+d+e_{1}+e_{2}+u_{1}+u_{2}}x_{1}^{l_{1}+u_{1}}x_{2}^{l_{2}+u_{2}}
\end{eqnarray*}

\begin{eqnarray*}
l_{1} &=&u_{1}=l_{2}=u_{2}=0 \\
a+b_{1}+b_{2}+d+e_{1}+e_{2} &\equiv &0 \\
d &=&1,e_{1}=1,e_{2}=0 \\
a &=&0,b_{1}=1,b_{2}=1
\end{eqnarray*}%
Since $\alpha \left( 1_{H};0,0,0,0\right) =0,$ we obtain
\begin{equation*}
B(x_{1}\otimes x_{1}x_{2};X_{1}X_{2},gx_{1})X_{1}X_{2}\otimes gx_{1}\otimes
1_{H}.
\end{equation*}%
By considering also the right side we get%
\begin{equation*}
B(x_{1}\otimes x_{1}x_{2};X_{1}X_{2},gx_{1})=B(x_{1}\otimes
x_{1}x_{2};X_{1}X_{2},gx_{1}).
\end{equation*}

\subsubsection{Case $X_{1}X_{2}\otimes g\otimes gx_{1}$}

From the first summand of the left side we get%
\begin{eqnarray*}
l_{1} &=&u_{1}=l_{2}=u_{2}=0 \\
a+b_{1}+b_{2}+d+e_{1}+e_{2} &\equiv &1 \\
d &=&1,e_{1}=e_{2}=0, \\
a &=&0,b_{1}=b_{2}=1.
\end{eqnarray*}%
Since $\alpha \left( x_{1};0,0,0,0\right) =a+b_{1}+b_{2}\equiv 0,$ we get%
\begin{equation*}
B(g\otimes x_{1}x_{2};X_{1}X_{2},g)X_{1}X_{2}\otimes g\otimes gx_{1}
\end{equation*}%
From the second summand of the left side we get%
\begin{eqnarray*}
l_{1}+u_{1} &=&1,l_{2}=u_{2}=0 \\
a+b_{1}+b_{2}+d+e_{1}+e_{2} &\equiv &1 \\
d &=&1,e_{1}=u_{1},e_{2}=0 \\
a &=&0,b_{1}-l_{1}=1\Rightarrow b_{1}=1,l_{1}=0,e_{1}=u_{1}=1,b_{2}=1.
\end{eqnarray*}%
Since $\alpha \left( 1_{H};0,0,1,0\right) =e_{2}+\left( a+b_{1}+b_{2}\right)
=1+1\equiv 0,$ we obtain%
\begin{equation*}
B(x_{1}\otimes x_{1}x_{2};X_{1}X_{2},gx_{1})X_{1}X_{2}\otimes g\otimes
gx_{1}.
\end{equation*}%
By considering also the right side, we get%
\begin{equation*}
B(g\otimes x_{1}x_{2};X_{1}X_{2},g)+B(x_{1}\otimes
x_{1}x_{2};X_{1}X_{2},gx_{1})+B(x_{1}\otimes x_{2};X_{1}X_{2},g)=0
\end{equation*}%
\begin{eqnarray*}
B(g\otimes x_{1}x_{2};X_{1}X_{2},g) &=&\left[
\begin{array}{c}
B(g\otimes 1_{H};1_{A},g)+B(x_{2}\otimes \ 1_{H};1_{A},gx_{2}) \\
+B(x_{1}\otimes 1_{H};1_{A},gx_{1})+B(gx_{1}x_{2}\otimes
1_{H};1_{A},gx_{1}x_{2})%
\end{array}%
\right] \\
+B(x_{1}\otimes x_{1}x_{2};X_{1}X_{2},gx_{1}) &=&\left[
\begin{array}{c}
-2B(g\otimes 1_{H};1_{A},g)-2B(x_{2}\otimes \ 1_{H};1_{A},gx_{2}) \\
-2B(x_{1}\otimes 1_{H};1_{A},gx_{1})-2B(gx_{1}x_{2}\otimes
1_{H};1_{A},gx_{1}x_{2})%
\end{array}%
\right] \\
+B(x_{1}\otimes x_{2};X_{1}X_{2},g) &=&\left[
\begin{array}{c}
B(g\otimes 1_{H};1_{A},g)+B(x_{2}\otimes \ 1_{H};1_{A},gx_{2}) \\
+B(x_{1}\otimes 1_{H};1_{A},gx_{1})+B(gx_{1}x_{2}\otimes
1_{H};1_{A},gx_{1}x_{2})%
\end{array}%
\right]
\end{eqnarray*}%
which holds in view of the form of the elements.

\subsubsection{Case $X_{1}\otimes gx_{1}\otimes gx_{2}$}

Nothing from the first summand of the left side. From the second summand we
get%
\begin{eqnarray*}
&&+\sum_{a,b_{1},b_{2},d,e_{1},e_{2}=0}^{1}\sum_{l_{1}=0}^{b_{1}}%
\sum_{l_{2}=0}^{b_{2}}\sum_{u_{1}=0}^{e_{1}}\sum_{u_{2}=0}^{e_{2}}\left(
-1\right) ^{\alpha \left( 1_{H};l_{1},l_{2},u_{1},u_{2}\right) } \\
&&B(x_{1}\otimes
x_{1}x_{2};G^{a}X_{1}^{b_{1}}X_{2}^{b_{2}},g^{d}x_{1}^{e_{1}}x_{2}^{e_{2}})G^{a}X_{1}^{b_{1}-l_{1}}X_{2}^{b_{2}-l_{2}}\otimes g^{d}x_{1}^{e_{1}-u_{1}}x_{2}^{e_{2}-u_{2}}\otimes
\\
&&g^{a+b_{1}+b_{2}+l_{1}+l_{2}+d+e_{1}+e_{2}+u_{1}+u_{2}}x_{1}^{l_{1}+u_{1}}x_{2}^{l_{2}+u_{2}}
\end{eqnarray*}%
\begin{eqnarray*}
l_{1} &=&u_{1}=0,l_{2}+u_{2}=1 \\
a+b_{1}+b_{2}+d+e_{1}+e_{2} &\equiv &0 \\
d &=&1,e_{1}=1,e_{2}=u_{2} \\
a &=&0,b_{1}=1,b_{2}=l_{2}
\end{eqnarray*}%
Since $\alpha \left( 1_{H};0,0,0,1\right) =$ $a+b_{1}+b_{2}=1$ and $\alpha
\left( 1_{H};0,1,0,0\right) =0,$ we get%
\begin{equation*}
\left[ -B(x_{1}\otimes x_{1}x_{2};X_{1},gx_{1}x_{2})+B(x_{1}\otimes
x_{1}x_{2};X_{1}X_{2},gx_{1})\right] X_{1}\otimes gx_{1}\otimes gx_{2}.
\end{equation*}%
By considering also the right side, we get%
\begin{equation*}
-B(x_{1}\otimes x_{1}x_{2};X_{1},gx_{1}x_{2})+B(x_{1}\otimes
x_{1}x_{2};X_{1}X_{2},gx_{1})-B(x_{1}\otimes x_{1};X_{1},gx_{1})=0
\end{equation*}%
\begin{eqnarray*}
-B(x_{1}\otimes x_{1}x_{2};X_{1},gx_{1}x_{2}) &=&\left[
\begin{array}{c}
+4B(g\otimes 1_{H};1_{A},g)+4B(x_{1}\otimes 1_{H};1_{A},gx_{1}) \\
+2B(x_{2}\otimes 1_{H};1_{A},gx_{2})+2B(gx_{1}x_{2}\otimes
1_{H};1_{A},gx_{1}x_{2})%
\end{array}%
\right] \\
+B(x_{1}\otimes x_{1}x_{2};X_{1}X_{2},gx_{1}) &=&+\left[
\begin{array}{c}
-2B(g\otimes 1_{H};1_{A},g)-2B(x_{2}\otimes \ 1_{H};1_{A},gx_{2}) \\
-2B(x_{1}\otimes 1_{H};1_{A},gx_{1})-2B(gx_{1}x_{2}\otimes
1_{H};1_{A},gx_{1}x_{2})%
\end{array}%
\right] \\
-B(x_{1}\otimes x_{1};X_{1},gx_{1}) &=&\left[ -2B(g\otimes
1_{H};1_{A},g)-2B(x_{1}\otimes 1_{H};1_{A},gx_{1})\right]
\end{eqnarray*}%
which holds in view of the form of the elements.

\subsubsection{Case $X_{1}\otimes g\otimes x_{1}x_{2}$}

From the first summand of the left side we get%
\begin{eqnarray*}
l_{1} &=&u_{1}=0,l_{2}+u_{2}=1 \\
a+b_{1}+b_{2}+d+e_{1}+e_{2} &\equiv &1 \\
d &=&1,e_{1}=0,e_{2}=u_{2}, \\
a &=&0,b_{1}=1,b_{2}=l_{2}.
\end{eqnarray*}%
Since $\alpha \left( x_{1};0,0,0,1\right) \equiv 1$ and $\alpha \left(
x_{1};0,1,0,0\right) \equiv a+b_{1}+b_{2}+1\equiv 0+1+1+1\equiv 1$
\begin{equation*}
\left[ -B(g\otimes x_{1}x_{2};X_{1},gx_{2})-B(g\otimes
x_{1}x_{2};X_{1}X_{2},g)\right] X_{1}\otimes g\otimes x_{1}x_{2}
\end{equation*}%
From the second summand of the left side we get%
\begin{eqnarray*}
l_{1}+u_{1} &=&1,l_{2}+u_{2}=1 \\
a+b_{1}+b_{2}+d+e_{1}+e_{2} &\equiv &0 \\
d &=&1,e_{1}=u_{1},e_{2}=u_{2}, \\
a &=&0,b_{1}-l_{1}=1\Rightarrow b_{1}=1,l_{1}=0,e_{1}=u_{1}=1,b_{2}=l_{2}.
\end{eqnarray*}%
Since $\alpha \left( 1_{H};0,0,1,1\right) \equiv 1+e_{2}\equiv 0$ and $%
\alpha \left( 1_{H};0,1,1,0\right) =e_{2}+a+b_{1}+b_{2}+1\equiv
0+0+1+1+1\equiv 1,$ we get%
\begin{equation*}
\left[ B(x_{1}\otimes x_{1}x_{2};X_{1},gx_{1}x_{2})+B(x_{1}\otimes
x_{1}x_{2};X_{1}X_{2},gx_{1})\right] X_{1}\otimes g\otimes x_{1}x_{2}.
\end{equation*}%
By taking in account also the right side, we get%
\begin{gather*}
-B(g\otimes x_{1}x_{2};X_{1},gx_{2})-B(g\otimes
x_{1}x_{2};X_{1}X_{2},g)+B(x_{1}\otimes x_{1}x_{2};X_{1},gx_{1}x_{2}) \\
-B(x_{1}\otimes x_{1}x_{2};X_{1}X_{2},gx_{1})-B(x_{1}\otimes 1_{H};X_{1},g)=0
\end{gather*}%
which holds in view of the form of the elements.

\subsubsection{Case $X_{2}\otimes gx_{1}\otimes gx_{1}$}

From the first summand of the left side we get%
\begin{eqnarray*}
&&\sum_{a,b_{1},b_{2},d,e_{1},e_{2}=0}^{1}\sum_{l_{1}=0}^{b_{1}}%
\sum_{l_{2}=0}^{b_{2}}\sum_{u_{1}=0}^{e_{1}}\sum_{u_{2}=0}^{e_{2}}\left(
-1\right) ^{\alpha \left( x_{1};l_{1},l_{2},u_{1},u_{2}\right) } \\
&&B(g\otimes
x_{1}x_{2};G^{a}X_{1}^{b_{1}}X_{2}^{b_{2}},g^{d}x_{1}^{e_{1}}x_{2}^{e_{2}})G^{a}X_{1}^{b_{1}-l_{1}}X_{2}^{b_{2}-l_{2}}\otimes g^{d}x_{1}^{e_{1}-u_{1}}x_{2}^{e_{2}-u_{2}}\otimes
\\
&&g^{a+b_{1}+b_{2}+l_{1}+l_{2}+d+e_{1}+e_{2}+u_{1}+u_{2}}x_{1}^{l_{1}+u_{1}+1}x_{2}^{l_{2}+u_{2}}
\end{eqnarray*}%
\begin{eqnarray*}
l_{1} &=&u_{1}=l_{2}=u_{2}=0 \\
a+b_{1}+b_{2}+d+e_{1}+e_{2} &\equiv &1 \\
d &=&1,e_{1}=1,e_{2}=0, \\
a &=&0,b_{1}=0,b_{2}=1.
\end{eqnarray*}%
Since $\alpha \left( x_{1};0,0,0,0\right) $ $=a+b_{1}+b_{2}=1$ we get%
\begin{equation*}
-B(g\otimes x_{1}x_{2};X_{2},gx_{1})X_{2}\otimes gx_{1}\otimes gx_{1}.
\end{equation*}%
From the second summand of the left side we get%
\begin{eqnarray*}
l_{1}+u_{1} &=&1,l_{2}=u_{2}=0 \\
a+b_{1}+b_{2}+d+e_{1}+e_{2} &\equiv &0 \\
d &=&1,e_{1}-u_{1}=1\Rightarrow e_{1}=1,u_{1}=0,l_{1}=1,e_{2}=0, \\
a &=&0,b_{1}=l_{1}=1,b_{2}=1.
\end{eqnarray*}%
Since $\alpha \left( 1_{H};1,0,0,0\right) =b_{2}=1,$ we get%
\begin{equation*}
-B(x_{1}\otimes x_{1}x_{2};X_{1}X_{2},gx_{1})X_{2}\otimes gx_{1}\otimes
gx_{1}.
\end{equation*}%
By taking in account also the right side, we get%
\begin{equation*}
-B(g\otimes x_{1}x_{2};X_{2},gx_{1})-B(x_{1}\otimes
x_{1}x_{2};X_{1}X_{2},gx_{1})+B(x_{1}\otimes x_{2};X_{2},gx_{1})=0
\end{equation*}%
which holds in view of the form of the element.

\subsubsection{Case $1_{H}\otimes gx_{1}\otimes x_{1}x_{2}$}

This case already appeared in $B\left( x_{1}\otimes
x_{1}x_{2};X_{1},gx_{1}x_{2}\right) .$

\subsection{$B(x_{1}\otimes x_{1}x_{2};GX_{1},gx_{1})$}

We get%
\begin{equation*}
a=1,b_{1}=1,b_{2}=0,d=1,e_{1}=1,e_{2}=0
\end{equation*}%
and then we obtain%
\begin{gather*}
+\sum_{l_{1}=0}^{b_{1}}\sum_{u_{1}=0}^{e_{1}}\left( -1\right) ^{\alpha
\left( 1_{H};l_{1},0,u_{1},0\right) }B(x_{1}\otimes
x_{1}x_{2};GX_{1},gx_{1})GX_{1}^{1-l_{1}}\otimes gx_{1}^{1-u_{1}}\otimes
g^{l_{1}+u_{1}}x_{1}^{l_{1}+u_{1}} \\
=\left( -1\right) ^{\alpha \left( 1_{H};0,0,0,0\right) }B(x_{1}\otimes
x_{1}x_{2};GX_{1},gx_{1})GX_{1}\otimes gx_{1}\otimes 1_{H} \\
+\left( -1\right) ^{\alpha \left( 1_{H};0,0,1,0\right) }B(x_{1}\otimes
x_{1}x_{2};GX_{1},gx_{1})GX_{1}\otimes g\otimes gx_{1} \\
+\left( -1\right) ^{\alpha \left( 1_{H};1,0,0,0\right) }B(x_{1}\otimes
x_{1}x_{2};GX_{1},gx_{1})G\otimes gx_{1}\otimes gx_{1} \\
+\left( -1\right) ^{\alpha \left( 1_{H};1,0,1,0\right) }B(x_{1}\otimes
x_{1}x_{2};GX_{1},gx_{1})G\otimes gx_{1}^{1-u_{1}}\otimes
g^{1+u_{1}}x_{1}^{1+1}=0
\end{gather*}

\subsubsection{Case $GX_{1}\otimes gx_{1}\otimes 1_{H}$}

From the first summand of the left side, we get nothing. From the second
summand we get%
\begin{eqnarray*}
&&+\sum_{a,b_{1},b_{2},d,e_{1},e_{2}=0}^{1}\sum_{l_{1}=0}^{b_{1}}%
\sum_{l_{2}=0}^{b_{2}}\sum_{u_{1}=0}^{e_{1}}\sum_{u_{2}=0}^{e_{2}}\left(
-1\right) ^{\alpha \left( 1_{H};l_{1},l_{2},u_{1},u_{2}\right) } \\
&&B(x_{1}\otimes
x_{1}x_{2};G^{a}X_{1}^{b_{1}}X_{2}^{b_{2}},g^{d}x_{1}^{e_{1}}x_{2}^{e_{2}})G^{a}X_{1}^{b_{1}-l_{1}}X_{2}^{b_{2}-l_{2}}\otimes g^{d}x_{1}^{e_{1}-u_{1}}x_{2}^{e_{2}-u_{2}}\otimes
\\
&&g^{a+b_{1}+b_{2}+l_{1}+l_{2}+d+e_{1}+e_{2}+u_{1}+u_{2}}x_{1}^{l_{1}+u_{1}}x_{2}^{l_{2}+u_{2}}
\end{eqnarray*}%
\begin{eqnarray*}
l_{1} &=&u_{1}=l_{2}=u_{2}=0 \\
a+b_{1}+b_{2}+d+e_{1}+e_{2} &\equiv &0 \\
d &=&1,e_{1}=1,e_{2}=0, \\
a &=&1,b_{1}=1,b_{2}=0
\end{eqnarray*}%
Since $\alpha \left( 1_{H};0,0,0,0\right) =0,$ we obtain
\begin{equation*}
B(x_{1}\otimes x_{1}x_{2};GX_{1},gx_{1})GX_{1}\otimes gx_{1}\otimes 1_{H}
\end{equation*}%
Considering also the right side, we obtain%
\begin{equation*}
B(x_{1}\otimes x_{1}x_{2};GX_{1},gx_{1})=B(x_{1}\otimes
x_{1}x_{2};GX_{1},gx_{1})
\end{equation*}%
which trivially holds.

\subsubsection{Case $GX_{1}\otimes g\otimes gx_{1}$}

From the first summand of the left side we get%
\begin{eqnarray*}
l_{1} &=&u_{1}=l_{2}=u_{2}=0 \\
a+b_{1}+b_{2}+d+e_{1}+e_{2} &\equiv &1 \\
d &=&1,e_{1}=0,e_{2}=0, \\
a &=&1,b_{1}=1,b_{2}=0
\end{eqnarray*}%
Since $\alpha \left( x_{1};0,0,0,0\right) =a+b_{1}+b_{2}\equiv 0,$ we obtain%
\begin{equation*}
B(g\otimes x_{1}x_{2};GX_{1},g)GX_{1}\otimes g\otimes gx_{1}.
\end{equation*}%
From the second summand of the left side we get%
\begin{eqnarray*}
l_{1}+u_{1} &=&1,l_{2}=u_{2}=0 \\
a+b_{1}+b_{2}+d+e_{1}+e_{2} &\equiv &1 \\
d &=&1,e_{1}=u_{1},e_{2}=0, \\
a &=&1,b_{1}=1,l_{1}=0\Rightarrow e_{1}=u_{1}=1,b_{2}=0
\end{eqnarray*}%
Since $\alpha \left( 1_{H};0,0,1,0\right) =e_{2}+\left( a+b_{1}+b_{2}\right)
\equiv 0$ we get%
\begin{equation*}
B(x_{1}\otimes x_{1}x_{2};GX_{1},gx_{1})GX_{1}\otimes g\otimes gx_{1}
\end{equation*}%
Thus, by considering also the right side, we obtain
\begin{equation*}
B(g\otimes x_{1}x_{2};GX_{1},g)+B(x_{1}\otimes
x_{1}x_{2};GX_{1},gx_{1})+B(x_{1}\otimes x_{2};GX_{1},g)=0
\end{equation*}%
which holds in view of the form of the elements.

\subsubsection{Case $G\otimes gx_{1}\otimes gx_{1}$}

From the first summand of the left side we get%
\begin{eqnarray*}
l_{1} &=&u_{1}=l_{2}=u_{2}=0 \\
a+b_{1}+b_{2}+d+e_{1}+e_{2} &\equiv &1 \\
d &=&1,e_{1}=1,e_{2}=0, \\
a &=&1,b_{1}=0,b_{2}=0
\end{eqnarray*}%
Since $\alpha \left( x_{1};0,0,0,0\right) =a+b_{1}+b_{2}\equiv 1$, we get%
\begin{equation*}
-B(g\otimes x_{1}x_{2};G,gx_{1})G\otimes gx_{1}\otimes gx_{1}
\end{equation*}%
From the second summand of the left side we get%
\begin{eqnarray*}
l_{1}+u_{1} &=&1,l_{2}=u_{2}=0 \\
a+b_{1}+b_{2}+d+e_{1}+e_{2} &\equiv &1 \\
d &=&1,e_{1}=1,u_{1}=0\Rightarrow b_{1}=l_{1}=1,e_{2}=0, \\
a &=&1,1,b_{2}=0
\end{eqnarray*}%
Since $\alpha \left( 1_{H};1,0,0,0\right) =b_{2}\equiv 0$ we get%
\begin{equation*}
B(x_{1}\otimes x_{1}x_{2};GX_{1},gx_{1})G\otimes gx_{1}\otimes gx_{1.}
\end{equation*}%
By considering also the right side we obtain%
\begin{equation*}
-B(g\otimes x_{1}x_{2};G,gx_{1})+B(x_{1}\otimes
x_{1}x_{2};GX_{1},gx_{1})+B(x_{1}\otimes x_{2};G,gx_{1})=0
\end{equation*}%
which holds in view of the form of the elements.

\subsection{$B(x_{1}\otimes x_{1}x_{2};GX_{2},gx_{1})$}

\begin{eqnarray*}
a &=&1,b_{1}=0,b_{2}=1, \\
d &=&1,e_{1}=1,e_{2}=0 \\
a+b_{1}+b_{2}+l_{1}+l_{2}+d+e_{1}+e_{2}+u_{1}+u_{2}
&=&1+1+l_{2}+1+1+u_{1}\equiv l_{2}+u_{1}
\end{eqnarray*}%
and we get%
\begin{eqnarray*}
&&\sum_{l_{2}=0}^{1}\sum_{u_{1}=0}^{1}\left( -1\right) ^{\alpha \left(
1_{H};0,l_{2},u_{1},0\right) }B(x_{1}\otimes
x_{1}x_{2};GX_{2},gx_{1})GX_{2}^{1-l_{2}}\otimes gx_{1}^{1-u_{1}}\otimes
g^{l_{2}+u_{1}}x_{1}^{u_{1}}x_{2}^{l_{2}} \\
&=&\left( -1\right) ^{\alpha \left( 1_{H};0,0,0,0\right) }B(x_{1}\otimes
x_{1}x_{2};GX_{2},gx_{1})GX_{2}\otimes gx_{1}\otimes 1_{H}+ \\
&&+\left( -1\right) ^{\alpha \left( 1_{H};0,0,1,0\right) }B(x_{1}\otimes
x_{1}x_{2};GX_{2},gx_{1})GX_{2}\otimes g\otimes gx_{1}+ \\
&&+\left( -1\right) ^{\alpha \left( 1_{H};0,1,0,0\right) }B(x_{1}\otimes
x_{1}x_{2};GX_{2},gx_{1})G\otimes gx_{1}\otimes gx_{2}+ \\
&&+\left( -1\right) ^{\alpha \left( 1_{H};0,1,1,0\right) }B(x_{1}\otimes
x_{1}x_{2};GX_{2},gx_{1})G\otimes g\otimes x_{1}x_{2}
\end{eqnarray*}

\subsubsection{Case $GX_{2}\otimes gx_{1}\otimes 1_{H}$}

First term of the left side gives us nothing. Second term of the left side

\begin{eqnarray*}
l_{1} &=&u_{1}=l_{2}=u_{2}=0 \\
a+b_{1}+b_{2}+d+e_{1}+e_{2} &\equiv &0 \\
d &=&1,e_{1}=1,e_{2}=0, \\
a &=&1,b_{1}=0,b_{2}=1
\end{eqnarray*}%
since $\alpha \left( 1_{H};0,0,0,0\right) =0$ we get%
\begin{equation*}
B(x_{1}\otimes x_{1}x_{2};GX_{2},gx_{1})GX_{2}\otimes gx_{1}\otimes 1_{H}
\end{equation*}%
Considering also the right side we get

\begin{equation*}
B(x_{1}\otimes x_{1}x_{2};GX_{2},gx_{1})=B(x_{1}\otimes
x_{1}x_{2};GX_{2},gx_{1})
\end{equation*}%
so that we obtain nothing new.

\subsubsection{Case $GX_{2}\otimes g\otimes gx_{1}$}

From the first summand of the left side we get
\begin{eqnarray*}
l_{1} &=&u_{1}=l_{2}=u_{2}=0 \\
a+b_{1}+b_{2}+d+e_{1}+e_{2} &\equiv &1 \\
d &=&1,e_{1}=0,e_{2}=0, \\
a &=&1,b_{1}=0,b_{2}=1
\end{eqnarray*}%
Since $\alpha \left( x_{1};0,0,0,0\right) =a+b_{1}+b_{2}\equiv 0$, we obtain%
\begin{equation*}
B(g\otimes x_{1}x_{2};GX_{2},g)GX_{2}\otimes g\otimes gx_{1.}
\end{equation*}%
From the second summand of the left side we get
\begin{eqnarray*}
l_{1}+u_{1} &=&1 \\
l_{2} &=&u_{2}=0 \\
a+b_{1}+b_{2}+d+e_{1}+e_{2} &\equiv &1 \\
d &=&1,e_{1}=u_{1},e_{2}=0, \\
a &=&1,b_{1}=l_{1},b_{2}=1\text{ }
\end{eqnarray*}%
Since $\alpha \left( 1_{H};0,0,1,0\right) =e_{2}+\left( a+b_{1}+b_{2}\right)
\equiv 0$ and $\alpha \left( 1_{H};1,0,0,0\right) =b_{2}=1$%
\begin{equation*}
\left[ B(x_{1}\otimes x_{1}x_{2};GX_{2},gx_{1})-B(x_{1}\otimes
x_{1}x_{2};GX_{1}X_{2},g)\right] GX_{2}\otimes g\otimes gx_{1}.
\end{equation*}%
By considering also the right side we get%
\begin{gather*}
B(g\otimes x_{1}x_{2};GX_{2},g)+B(x_{1}\otimes x_{1}x_{2};GX_{2},gx_{1}) \\
-B(x_{1}\otimes x_{1}x_{2};GX_{1}X_{2},g)+B(x_{1}\otimes x_{2};GX_{2},g)=0.
\end{gather*}%
which holds in view of the form of the elements.

\subsubsection{Case $G\otimes gx_{1}\otimes gx_{2}$}

Nothing from the first summand of the left side.

For the second summand we get%
\begin{eqnarray*}
&&+\sum_{a,b_{1},b_{2},d,e_{1},e_{2}=0}^{1}\sum_{l_{1}=0}^{b_{1}}%
\sum_{l_{2}=0}^{b_{2}}\sum_{u_{1}=0}^{e_{1}}\sum_{u_{2}=0}^{e_{2}}\left(
-1\right) ^{\alpha \left( 1_{H};l_{1},l_{2},u_{1},u_{2}\right) } \\
&&B(x_{1}\otimes
x_{1}x_{2};G^{a}X_{1}^{b_{1}}X_{2}^{b_{2}},g^{d}x_{1}^{e_{1}}x_{2}^{e_{2}})G^{a}X_{1}^{b_{1}-l_{1}}X_{2}^{b_{2}-l_{2}}\otimes g^{d}x_{1}^{e_{1}-u_{1}}x_{2}^{e_{2}-u_{2}}\otimes
\\
&&g^{a+b_{1}+b_{2}+l_{1}+l_{2}+d+e_{1}+e_{2}+u_{1}+u_{2}}x_{1}^{l_{1}+u_{1}}x_{2}^{l_{2}+u_{2}}
\end{eqnarray*}%
\begin{eqnarray*}
l_{1} &=&u_{1}=0 \\
l_{2}+u_{2} &=&1 \\
a+b_{1}+b_{2}+d+e_{1}+e_{2} &\equiv &0 \\
d &=&1,e_{1}=1,e_{2}=u_{2}, \\
a &=&1,b_{1}=0,b_{2}=l_{2}\text{ }
\end{eqnarray*}%
Since $\alpha \left( 1_{H};0,0,0,1\right) =a+b_{1}+b_{2}\equiv 1$ and $%
\alpha \left( 1_{H};0,1,0,0\right) =0,$ we get

\begin{equation*}
\left[ -B(x_{1}\otimes x_{1}x_{2};G,gx_{1}x_{2})+B(x_{1}\otimes
x_{1}x_{2};GX_{2},gx_{1})\right] G\otimes gx_{1}\otimes gx_{2}
\end{equation*}%
By taking in account also the right side we get%
\begin{equation*}
-B(x_{1}\otimes x_{1}x_{2};G,gx_{1}x_{2})+B(x_{1}\otimes
x_{1}x_{2};GX_{2},gx_{1})-B(x_{1}\otimes x_{1};G,gx_{1})=0
\end{equation*}%
which holds in view of the form of the elements.

\subsubsection{Case $G\otimes g\otimes x_{1}x_{2}$}

From the first summand of the left side we get
\begin{eqnarray*}
l_{1} &=&u_{1}=0 \\
l_{2}+u_{2} &=&1 \\
a+b_{1}+b_{2}+d+e_{1}+e_{2} &\equiv &1 \\
d &=&1,e_{1}=0,e_{2}=u_{2}, \\
a &=&1,b_{1}=0,b_{2}=l_{2}
\end{eqnarray*}%
Since $\alpha \left( x_{1};0,0,0,1\right) \equiv 1$ and $\alpha \left(
x_{1};0,1,0,0\right) =a+b_{1}+b_{2}+1\equiv 1$ we obtain%
\begin{equation*}
\left[ -B(g\otimes x_{1}x_{2};G,gx_{2})-B(g\otimes x_{1}x_{2};GX_{2},g)%
\right] G\otimes g\otimes x_{1}x_{2}
\end{equation*}%
From the second summand of the left side we get
\begin{eqnarray*}
l_{1}+u_{1} &=&1 \\
l_{2}+u_{2} &=&1 \\
a+b_{1}+b_{2}+d+e_{1}+e_{2} &\equiv &0 \\
d &=&1,e_{1}=u_{1},e_{2}=u_{2}, \\
a &=&1,b_{1}=l_{1},b_{2}=l_{2}
\end{eqnarray*}%
Since $\alpha \left( 1_{H};0,0,1,1\right) =1+e_{2}\equiv 0$, $\alpha \left(
1_{H};0,1,1,0\right) \equiv 1,\alpha \left( 1_{H};1,0,0,1\right) \equiv 0$
and $\alpha \left( 1_{H};1,1,0,0\right) \equiv 0$%
\begin{equation*}
\left[
\begin{array}{c}
B(x_{1}\otimes x_{1}x_{2};G,gx_{1}x_{2})-B(x_{1}\otimes
x_{1}x_{2};GX_{2},gx_{1}) \\
+B(x_{1}\otimes x_{1}x_{2};GX_{1},gx_{2})+B(x_{1}\otimes
x_{1}x_{2};GX_{1}X_{2},g)%
\end{array}%
\right] G\otimes g\otimes x_{1}x_{2}
\end{equation*}%
By considering also the right side, finally we get%
\begin{gather*}
-B(g\otimes x_{1}x_{2};G,gx_{2})-B(g\otimes
x_{1}x_{2};GX_{2},g)+B(x_{1}\otimes x_{1}x_{2};G,gx_{1}x_{2}) \\
-B(x_{1}\otimes x_{1}x_{2};GX_{2},gx_{1})+B(x_{1}\otimes
x_{1}x_{2};GX_{1},gx_{2})+ \\
+B(x_{1}\otimes x_{1}x_{2};GX_{1}X_{2},g)+-B(x_{1}\otimes 1_{H};G,g)=0
\end{gather*}%
which holds in view of the form of the elements.

\begin{eqnarray*}
&&\sum_{a,b_{1},b_{2},d,e_{1},e_{2}=0}^{1}\sum_{l_{1}=0}^{b_{1}}%
\sum_{l_{2}=0}^{b_{2}}\sum_{u_{1}=0}^{e_{1}}\sum_{u_{2}=0}^{e_{2}}\left(
-1\right) ^{\alpha \left( x_{1};l_{1},l_{2},u_{1},u_{2}\right) } \\
&&B(g\otimes
x_{1}x_{2};G^{a}X_{1}^{b_{1}}X_{2}^{b_{2}},g^{d}x_{1}^{e_{1}}x_{2}^{e_{2}})G^{a}X_{1}^{b_{1}-l_{1}}X_{2}^{b_{2}-l_{2}}\otimes g^{d}x_{1}^{e_{1}-u_{1}}x_{2}^{e_{2}-u_{2}}\otimes
\\
&&g^{a+b_{1}+b_{2}+l_{1}+l_{2}+d+e_{1}+e_{2}+u_{1}+u_{2}}x_{1}^{l_{1}+u_{1}+1}x_{2}^{l_{2}+u_{2}}
\\
&&+\sum_{a,b_{1},b_{2},d,e_{1},e_{2}=0}^{1}\sum_{l_{1}=0}^{b_{1}}%
\sum_{l_{2}=0}^{b_{2}}\sum_{u_{1}=0}^{e_{1}}\sum_{u_{2}=0}^{e_{2}}\left(
-1\right) ^{\alpha \left( 1_{H};l_{1},l_{2},u_{1},u_{2}\right) } \\
&&B(x_{1}\otimes
x_{1}x_{2};G^{a}X_{1}^{b_{1}}X_{2}^{b_{2}},g^{d}x_{1}^{e_{1}}x_{2}^{e_{2}})G^{a}X_{1}^{b_{1}-l_{1}}X_{2}^{b_{2}-l_{2}}\otimes g^{d}x_{1}^{e_{1}-u_{1}}x_{2}^{e_{2}-u_{2}}\otimes
\\
&&g^{a+b_{1}+b_{2}+l_{1}+l_{2}+d+e_{1}+e_{2}+u_{1}+u_{2}}x_{1}^{l_{1}+u_{1}}x_{2}^{l_{2}+u_{2}}
\\
&=&B^{A}(x_{1}\otimes x_{1}x_{2})\otimes B^{H}(x_{1}\otimes
x_{1}x_{2})\otimes 1_{H} \\
&&B^{A}(x_{1}\otimes x_{1})\otimes B^{H}(x_{1}\otimes x_{1})\otimes gx_{2} \\
&&\left( -1\right) B^{A}(x_{1}\otimes x_{2})\otimes B^{H}(x_{1}\otimes
x_{2})\otimes gx_{1} \\
&&B^{A}(x_{1}\otimes 1_{H})\otimes B^{H}(x_{1}\otimes 1_{H})\otimes
x_{1}x_{2}
\end{eqnarray*}

\subsection{$B(x_{1}\otimes x_{1}x_{2};GX_{1}X_{2},gx_{1}x_{2})$}

We deduce that%
\begin{equation*}
a=b_{1}=b_{2}=d=e_{1}=e_{2}=1
\end{equation*}%
and we obtain

\begin{gather*}
\sum_{l_{1}=0}^{1}\sum_{l_{2}=0}^{1}\sum_{u_{1}=0}^{1}\sum_{u_{2}=0}^{1}%
\left( -1\right) ^{\alpha \left( 1_{H};l_{1},l_{2},u_{1},u_{2}\right)
}B(x_{1}\otimes x_{1}x_{2};GX_{1}X_{2},gx_{1}x_{2}) \\
G^{a}X_{1}^{1-l_{1}}X_{2}^{1-l_{2}}\otimes
gx_{1}^{1-u_{1}}x_{2}^{1-u_{2}}\otimes
g^{l_{1}+l_{2}+u_{1}+u_{2}}x_{1}^{l_{1}+u_{1}}x_{2}^{l_{2}+u_{2}} \\
=\left( -1\right) ^{\alpha \left( 1_{H};0,0,0,0\right) }B(x_{1}\otimes
x_{1}x_{2};GX_{1}X_{2},gx_{1}x_{2})GX_{1}X_{2}\otimes gx_{1}x_{2}\otimes
1_{H}+ \\
+\left( -1\right) ^{\alpha \left( 1_{H};0,0,0,1\right) }B(x_{1}\otimes
x_{1}x_{2};GX_{1}X_{2},gx_{1}x_{2})GX_{1}X_{2}\otimes gx_{1}\otimes gx_{2}%
\text{ } \\
+\left( -1\right) ^{\alpha \left( 1_{H};0,0,1,0\right) }B(x_{1}\otimes
x_{1}x_{2};GX_{1}X_{2},gx_{1}x_{2})GX_{1}X_{2}\otimes gx_{2}\otimes gx_{1}+
\\
+\left( -1\right) ^{\alpha \left( 1_{H};0,0,1,1\right) }B(x_{1}\otimes
x_{1}x_{2};GX_{1}X_{2},gx_{1}x_{2})GX_{1}X_{2}\otimes g\otimes x_{1}x_{2}+ \\
+\left( -1\right) ^{\alpha \left( 1_{H};0,1,0,0\right) }B(x_{1}\otimes
x_{1}x_{2};GX_{1}X_{2},gx_{1}x_{2})GX_{1}\otimes gx_{1}x_{2}\otimes gx_{2}+
\\
+\left( -1\right) ^{\alpha \left( 1_{H};0,1,0,1\right) }B(x_{1}\otimes
x_{1}x_{2};GX_{1}X_{2},gx_{1}x_{2})GX_{1}\otimes
gx_{1}x_{2}^{1-u_{2}}\otimes g^{u_{2}}x_{2}^{1+1}+=0 \\
+\left( -1\right) ^{\alpha \left( 1_{H};0,1,1,0\right) }B(x_{1}\otimes
x_{1}x_{2};GX_{1}X_{2},gx_{1}x_{2})GX_{1}\otimes gx_{2}\otimes x_{1}x_{2}+ \\
+\left( -1\right) ^{\alpha \left( 1_{H};0,1,1,1\right) }B(x_{1}\otimes
x_{1}x_{2};GX_{1}X_{2},gx_{1}x_{2})GX_{1}\otimes gx_{2}^{1-u_{2}}\otimes
g^{1+u_{2}}x_{1}x_{2}^{1+1}+=0 \\
+\left( -1\right) ^{\alpha \left( 1_{H};1,0,0,0\right) }B(x_{1}\otimes
x_{1}x_{2};GX_{1}X_{2},gx_{1}x_{2})GX_{2}\otimes gx_{1}x_{2}\otimes gx_{1}+
\\
+\left( -1\right) ^{\alpha \left( 1_{H};1,0,0,1\right) }B(x_{1}\otimes
x_{1}x_{2};GX_{1}X_{2},gx_{1}x_{2})GX_{2}\otimes
gx_{1}x_{2}^{1-u_{2}}\otimes g^{1+u_{2}}x_{1}^{1+0}x_{2}^{1+1}+=0 \\
+\left( -1\right) ^{\alpha \left( 1_{H};1,0,1,u_{2}\right) }B(x_{1}\otimes
x_{1}x_{2};GX_{1}X_{2},gx_{1}x_{2})GX_{2}\otimes
gx_{1}^{1-u_{1}}x_{2}^{1-u_{2}}\otimes
g^{1+u_{1}+u_{2}}x_{1}^{1+1}x_{2}^{1+u_{2}}+=0 \\
+\left( -1\right) ^{\alpha \left( 1_{H};1,1,0,0\right) }B(x_{1}\otimes
x_{1}x_{2};GX_{1}X_{2},gx_{1}x_{2})G\otimes gx_{1}x_{2}\otimes x_{1}x_{2}+ \\
+\left( -1\right) ^{\alpha \left( 1_{H};1,1,0,1\right) }B(x_{1}\otimes
x_{1}x_{2};GX_{1}X_{2},gx_{1}x_{2})G\otimes gx_{1}x_{2}^{1-u_{2}}\otimes
g^{u_{2}}x_{1}x_{2}^{1+1}+=0 \\
+\left( -1\right) ^{\alpha \left( 1_{H};1,1,1,u_{2}\right) }B(x_{1}\otimes
x_{1}x_{2};GX_{1}X_{2},gx_{1}x_{2})G\otimes
gx_{1}^{1-u_{1}}x_{2}^{1-u_{2}}\otimes
g^{u_{1}+u_{2}}x_{1}^{1+1}x_{2}^{1+u_{2}}=0.
\end{gather*}%
thus we obtain%
\begin{gather*}
\sum_{l_{1}=0}^{1}\sum_{l_{2}=0}^{1}\sum_{u_{1}=0}^{1}\sum_{u_{2}=0}^{1}%
\left( -1\right) ^{\alpha \left( 1_{H};l_{1},l_{2},u_{1},u_{2}\right)
}B(x_{1}\otimes x_{1}x_{2};GX_{1}X_{2},gx_{1}x_{2}) \\
G^{a}X_{1}^{1-l_{1}}X_{2}^{1-l_{2}}\otimes
gx_{1}^{1-u_{1}}x_{2}^{1-u_{2}}\otimes
g^{l_{1}+l_{2}+u_{1}+u_{2}}x_{1}^{l_{1}+u_{1}}x_{2}^{l_{2}+u_{2}} \\
=\left( -1\right) ^{\alpha \left( 1_{H};0,0,0,0\right) }B(x_{1}\otimes
x_{1}x_{2};GX_{1}X_{2},gx_{1}x_{2})GX_{1}X_{2}\otimes gx_{1}x_{2}\otimes
1_{H}+ \\
+\left( -1\right) ^{\alpha \left( 1_{H};0,0,0,1\right) }B(x_{1}\otimes
x_{1}x_{2};GX_{1}X_{2},gx_{1}x_{2})GX_{1}X_{2}\otimes gx_{1}\otimes gx_{2} \\
+\left( -1\right) ^{\alpha \left( 1_{H};0,0,1,0\right) }B(x_{1}\otimes
x_{1}x_{2};GX_{1}X_{2},gx_{1}x_{2})GX_{1}X_{2}\otimes gx_{2}\otimes gx_{1}+
\\
+\left( -1\right) ^{\alpha \left( 1_{H};0,0,1,1\right) }B(x_{1}\otimes
x_{1}x_{2};GX_{1}X_{2},gx_{1}x_{2})GX_{1}X_{2}\otimes g\otimes x_{1}x_{2}+ \\
+\left( -1\right) ^{\alpha \left( 1_{H};0,1,0,0\right) }B(x_{1}\otimes
x_{1}x_{2};GX_{1}X_{2},gx_{1}x_{2})GX_{1}\otimes gx_{1}x_{2}\otimes x_{2}+ \\
+\left( -1\right) ^{\alpha \left( 1_{H};0,1,1,0\right) }B(x_{1}\otimes
x_{1}x_{2};GX_{1}X_{2},gx_{1}x_{2})GX_{1}\otimes gx_{2}\otimes gx_{1}x_{2}+
\\
+\left( -1\right) ^{\alpha \left( 1_{H};1,0,0,0\right) }B(x_{1}\otimes
x_{1}x_{2};GX_{1}X_{2},gx_{1}x_{2})GX_{2}\otimes gx_{1}x_{2}\otimes
gx_{1}x_{2}+ \\
+\left( -1\right) ^{\alpha \left( 1_{H};1,1,0,0\right) }B(x_{1}\otimes
x_{1}x_{2};GX_{1}X_{2},gx_{1}x_{2})G\otimes gx_{1}x_{2}\otimes
g^{u_{2}}x_{1}x_{2}+
\end{gather*}

\subsubsection{Case $GX_{1}X_{2}\otimes gx_{1}x_{2}\otimes 1_{H}$}

We obtain%
\begin{equation*}
B(x_{1}\otimes x_{1}x_{2};GX_{1}X_{2},gx_{1}x_{2})=B(x_{1}\otimes
x_{1}x_{2};GX_{1}X_{2},gx_{1}x_{2})
\end{equation*}%
and we get nothing new.

\subsubsection{Case $GX_{1}X_{2}\otimes gx_{1}\otimes gx_{2}$}

Nothing from the first summand of the left side. From the second summand of
the left side we get%
\begin{eqnarray*}
&&+\sum_{a,b_{1},b_{2},d,e_{1},e_{2}=0}^{1}\sum_{l_{1}=0}^{b_{1}}%
\sum_{l_{2}=0}^{b_{2}}\sum_{u_{1}=0}^{e_{1}}\sum_{u_{2}=0}^{e_{2}}\left(
-1\right) ^{\alpha \left( 1_{H};l_{1},l_{2},u_{1},u_{2}\right) } \\
&&B(x_{1}\otimes
x_{1}x_{2};G^{a}X_{1}^{b_{1}}X_{2}^{b_{2}},g^{d}x_{1}^{e_{1}}x_{2}^{e_{2}})G^{a}X_{1}^{b_{1}-l_{1}}X_{2}^{b_{2}-l_{2}}\otimes g^{d}x_{1}^{e_{1}-u_{1}}x_{2}^{e_{2}-u_{2}}\otimes
\\
&&g^{a+b_{1}+b_{2}+l_{1}+l_{2}+d+e_{1}+e_{2}+u_{1}+u_{2}}x_{1}^{l_{1}+u_{1}}x_{2}^{l_{2}+u_{2}}
\end{eqnarray*}%
\begin{eqnarray*}
l_{1} &=&u_{1}=0 \\
l_{2}+u_{2} &=&1 \\
a+b_{1}+b_{2}+d+e_{1}+e_{2} &\equiv &0 \\
d &=&1,e_{1}=1,e_{2}=u_{2}, \\
a &=&1,b_{1}=1,b_{2}-l_{2}=1\Rightarrow b_{2}=1,l_{2}=0,u_{2}=e_{2}=1
\end{eqnarray*}%
Since $\alpha \left( 1_{H};0,0,0,1\right) =a+b_{1}+b_{2}\equiv 1,$ we obtain%
\begin{equation*}
-B(x_{1}\otimes x_{1}x_{2};GX_{1}X_{2},gx_{1}x_{2})GX_{1}X_{2}\otimes
gx_{1}\otimes gx_{2}.
\end{equation*}%
By considering also the right side, we obtain%
\begin{equation*}
-B(x_{1}\otimes x_{1}x_{2};GX_{1}X_{2},gx_{1}x_{2})-B(x_{1}\otimes
x_{1};GX_{1}X_{2},gx_{1})=0
\end{equation*}%
which holds in view of the form of the elements.

\subsubsection{Case $GX_{1}X_{2}\otimes gx_{2}\otimes gx_{1}$}

From the fist summand of the left side we get%
\begin{eqnarray*}
l_{1} &=&u_{1}=l_{2}=u_{2}=0 \\
a+b_{1}+b_{2}+d+e_{1}+e_{2} &\equiv &1 \\
d &=&1,e_{1}=0,e_{2}=1, \\
a &=&1,b_{1}=1,b_{2}=1
\end{eqnarray*}%
Since $\alpha \left( x_{1};0,0,0,0\right) =a+b_{1}+b_{2}\equiv 1$, we obtain%
\begin{equation*}
-B(g\otimes x_{1}x_{2};GX_{1}X_{2},gx_{2})GX_{1}X_{2}\otimes gx_{2}\otimes
gx_{1}
\end{equation*}%
From the second summand of the left side we get%
\begin{eqnarray*}
l_{1}+u_{1} &=&1 \\
l_{2} &=&u_{2}=0 \\
a+b_{1}+b_{2}+d+e_{1}+e_{2} &\equiv &0 \\
d &=&1,e_{1}=u_{1},e_{2}=1, \\
a &=&1,b_{1}-l_{1}=1\Rightarrow b_{1}=1,l_{1}=0,e_{1}=u_{1}=1,b_{2}=1
\end{eqnarray*}%
Since $\alpha \left( 1_{H};0,0,1,0\right) =e_{2}+\left( a+b_{1}+b_{2}\right)
\equiv 0,$ we obtain%
\begin{equation*}
B(x_{1}\otimes x_{1}x_{2};GX_{1}X_{2},gx_{1}x_{2})GX_{1}X_{2}\otimes
gx_{2}\otimes gx_{1}.
\end{equation*}%
By considering also the right side we get%
\begin{equation*}
-B(g\otimes x_{1}x_{2};GX_{1}X_{2},gx_{2})+B(x_{1}\otimes
x_{1}x_{2};GX_{1}X_{2},gx_{1}x_{2})+B(x_{1}\otimes
x_{2};GX_{1}X_{2},gx_{2})=0
\end{equation*}%
which holds in view of the form of the elements.

\subsubsection{Case $GX_{1}X_{2}\otimes g\otimes x_{1}x_{2}$}

From the fist summand of the left side we get%
\begin{eqnarray*}
l_{1} &=&u_{1}=0 \\
l_{2}+u_{2} &=&1 \\
a+b_{1}+b_{2}+d+e_{1}+e_{2} &\equiv &1 \\
d &=&1,e_{1}=0,e_{2}=u_{2}, \\
a &=&1,b_{1}=1,b_{2}-l_{2}=1\Rightarrow b_{2}=1,l_{2}=0,e_{2}=u_{2}=1
\end{eqnarray*}%
Since $\alpha \left( x_{1};0,0,0,1\right) \equiv 1$, we obtain%
\begin{equation*}
-B(g\otimes x_{1}x_{2};GX_{1}X_{2},gx_{2})GX_{1}X_{2}\otimes g\otimes
x_{1}x_{2}.
\end{equation*}%
From the second summand of the left side we get%
\begin{gather*}
l_{1}+u_{1}=1 \\
l_{2}+u_{2}=1 \\
a+b_{1}+b_{2}+d+e_{1}+e_{2}\equiv 1 \\
d=1,e_{1}=u_{1},e_{2}=u_{2}, \\
a=1,b_{1}-l_{1}=1\Rightarrow
b_{1}=1,l_{1}=0,e_{1}=u_{1}=1,b_{2}-l_{2}=1\Rightarrow
b_{2}=1,l_{2}=0,e_{2}=u_{2}=1
\end{gather*}%
Since $\alpha \left( 1_{H};0,0,1,1\right) =1+e_{2}\equiv 0,$ we obtain%
\begin{equation*}
B(x_{1}\otimes x_{1}x_{2};GX_{1}X_{2},gx_{1}x_{2})GX_{1}X_{2}\otimes
g\otimes x_{1}x_{2}.
\end{equation*}%
By considering also the right side, we get%
\begin{equation*}
-B(g\otimes x_{1}x_{2};GX_{1}X_{2},gx_{2})+B(x_{1}\otimes
x_{1}x_{2};GX_{1}X_{2},gx_{1}x_{2})-B(x_{1}\otimes 1_{H};GX_{1}X_{2},g)=0
\end{equation*}%
which holds in view of the form of the elements.

\subsubsection{Case $GX_{1}\otimes gx_{1}x_{2}\otimes gx_{2}$}

Nothing from the first summand of the left side. From the second summand of
the left side we get%
\begin{eqnarray*}
l_{1} &=&u_{1}=0,l_{2}+u_{2}=1 \\
a &=&1,b_{1}=1,b_{2}=l_{2}, \\
d &=&1,e_{1}=1,e_{2}-u_{2}=1\Rightarrow u_{2}=0,e_{2}=1,b_{2}=l_{2}=1.
\end{eqnarray*}%
Since $\alpha \left( 1_{H};0,1,0,0\right) \equiv 0,$ we obtain%
\begin{equation*}
+B(x_{1}\otimes x_{1}x_{2};GX_{1}X_{2},gx_{1}x_{2})GX_{1}\otimes
gx_{1}x_{2}\otimes gx_{2}.
\end{equation*}%
By considering also the right side, we obtain%
\begin{equation*}
+B(x_{1}\otimes x_{1}x_{2};GX_{1}X_{2},gx_{1}x_{2})-B(x_{1}\otimes
x_{1};GX_{1},gx_{1}x_{2})=0
\end{equation*}%
which holds in view of the form of the elements.

\subsubsection{Case $GX_{1}\otimes gx_{2}\otimes x_{1}x_{2}$}

From the first summand of the left side we get%
\begin{gather*}
l_{1}=u_{1}=0,l_{2}+u_{2}=1 \\
a=1,b_{1}=1,b_{2}=l_{2} \\
d=1,e_{1}=0,e_{2}-u_{2}=1\Rightarrow e_{2}=1,u_{2}=0,b_{2}=l_{2}=1
\end{gather*}%
Since $\alpha \left( x_{1};0,1,0,0\right) \equiv a+b_{1}+b_{2}+1\equiv 0$,
we get%
\begin{equation*}
+B(g\otimes x_{1}x_{2};GX_{1}X_{2},gx_{2})GX_{1}\otimes gx_{2}\otimes
x_{1}x_{2}
\end{equation*}%
From the second term of the left side, we get%
\begin{gather*}
l_{1}+u_{1}=1,l_{2}+u_{2}=1 \\
a=1,b_{1}-l_{1}=1\Rightarrow b_{1}=1,l_{1}=0,u_{1}=1,b_{2}=l_{2} \\
d=1,e_{1}=u_{1}=1,e_{2}-u_{2}=1\Rightarrow e_{2}=1,u_{2}=0,b_{2}=l_{2}=1
\end{gather*}%
Since $\alpha \left( 1_{H};0,1,1,0\right) \equiv e_{2}+a+b_{1}+b_{2}+1\equiv
1,$we get

\begin{equation*}
-B(x_{1}\otimes x_{1}x_{2};GX_{1}X_{2},gx_{1}x_{2})GX_{1}\otimes
gx_{2}\otimes x_{1}x_{2}.
\end{equation*}%
By considering also the right side, we get%
\begin{equation*}
+B(g\otimes x_{1}x_{2};GX_{1}X_{2},gx_{2})-B(x_{1}\otimes
x_{1}x_{2};GX_{1}X_{2},gx_{1}x_{2})-B(x_{1}\otimes 1_{H};GX_{1},gx_{2})=0
\end{equation*}%
which holds in view of the form of the elements.

\subsubsection{Case $GX_{2}\otimes gx_{1}x_{2}\otimes gx_{1}$}

From the fist summand of the left side we get%
\begin{eqnarray*}
l_{1} &=&u_{1}=l_{2}=u_{2}=0, \\
a &=&b_{2}=1,b_{1}=0, \\
d &=&1,e_{1}=0,e_{2}=1,
\end{eqnarray*}%
Since $\alpha \left( x_{1};0,0,0,0\right) =a+b_{1}+b_{2}\equiv 0$, we obtain%
\begin{equation*}
+B(g\otimes x_{1}x_{2};GX_{2},gx_{1}x_{2})GX_{2}\otimes gx_{1}x_{2}\otimes
gx_{1}.
\end{equation*}%
From the second summand of the left side we get%
\begin{eqnarray*}
l_{1}+u_{1} &=&1,l_{2}=u_{2}=0 \\
a &=&b_{2}=1,b_{1}=l_{1} \\
d &=&e_{2}=1,e_{1}-u_{1}=1\Rightarrow e_{1}=1,u_{1}=0,b_{1}=l_{1}=1.
\end{eqnarray*}%
Since $\alpha \left( 1_{H};1,0,0,0\right) \equiv b_{2}\equiv 1,$ we obtain%
\begin{equation*}
+B(x_{1}\otimes x_{1}x_{2};GX_{1}X_{2},gx_{1}x_{2})GX_{2}\otimes
gx_{1}x_{2}\otimes gx_{1}.
\end{equation*}%
By considering also the right side we get%
\begin{equation*}
+B(g\otimes x_{1}x_{2};GX_{2},gx_{1}x_{2})-B(x_{1}\otimes
x_{1}x_{2};GX_{1}X_{2},gx_{1}x_{2})+B(x_{1}\otimes
x_{2};GX_{2},gx_{1}x_{2})=0
\end{equation*}%
which holds in view of the form of the elements.

\subsubsection{Case $GX_{2}\otimes gx_{1}x_{2}\otimes gx_{1}x_{2}$}

From the first term of the left side we get%
\begin{gather*}
l_{1}=u_{1}=0 \\
l_{2}+u_{2}=1 \\
a+b_{1}+b_{2}+d+e_{1}+e_{2}\equiv 0 \\
d=1,e_{1}=1,e_{2}-u_{2}=1\Rightarrow e_{2}=1,u_{2}=0,l_{2}=1 \\
a=1,b_{1}=0, \\
b_{2}-l_{2}=1
\end{gather*}%
This has no solution. From the second term of the left side we get$%
GX_{2}\otimes gx_{1}x_{2}\otimes gx_{1}x_{2}$%
\begin{gather*}
l_{1}+u_{1}=1 \\
l_{2}+u_{2}=1 \\
a+b_{1}+b_{2}+d+e_{1}+e_{2}\equiv 1 \\
d=1,e_{1}-u_{1}=1\Rightarrow e_{1}=1,u_{1}=0,l_{1}=1 \\
e_{2}-u_{2}=1\Rightarrow e_{2}=1,u_{2}=0,l_{2}=1, \\
a=1,b_{2}-l_{2}=1
\end{gather*}%
This has no solution.

\subsubsection{Case $G\otimes gx_{1}x_{2}\otimes x_{1}x_{2}$}

From the first term of the left side we get%
\begin{gather*}
l_{1}=u_{1}=0 \\
l_{2}+u_{2}=1 \\
a+b_{1}+b_{2}+d+e_{1}+e_{2}\equiv 1 \\
d=1,e_{1}=1,e_{2}-u_{2}=1\Rightarrow e_{2}=1,u_{2}=0,l_{2}=1 \\
a=1,b_{1}=0,b_{2}=l_{2}=1
\end{gather*}%
Since $\alpha \left( x_{1};0,1,0,0\right) =a+b_{1}+b_{2}+1\equiv 1,$ we
obtain%
\begin{equation*}
-B(g\otimes x_{1}x_{2};GX_{2},gx_{1}x_{2})G\otimes gx_{1}x_{2}\otimes
x_{1}x_{2}
\end{equation*}%
From the second term of the left side we get%
\begin{gather*}
l_{1}+u_{1}=1 \\
l_{2}+u_{2}=1 \\
a+b_{1}+b_{2}+d+e_{1}+e_{2}\equiv 0 \\
d=1,e_{1}-u_{1}=1\Rightarrow e_{1}=1,u_{1}=0,l_{1}=1 \\
e_{2}-u_{2}=1\Rightarrow e_{2}=1,u_{2}=0,l_{2}=1 \\
a=1,b_{1}=l_{1}=1,b_{2}=l_{2}=1
\end{gather*}%
Since $\alpha \left( 1_{H};1,1,0,0\right) \equiv 1+b_{2}\equiv 0,$ we get%
\begin{equation*}
B(x_{1}\otimes x_{1}x_{2};GX_{1}X_{2},gx_{1}x_{2})G\otimes
gx_{1}x_{2}\otimes x_{1}x_{2}
\end{equation*}

By considering also the right side, we obtain%
\begin{equation*}
-B(g\otimes x_{1}x_{2};GX_{2},gx_{1}x_{2})+B(x_{1}\otimes
x_{1}x_{2};GX_{1}X_{2},gx_{1}x_{2})-B(x_{1}\otimes 1_{H};G,gx_{1}x_{2})=0
\end{equation*}%
which holds in view of the form of the elements.

\section{$B$$\left( x_{1}\otimes gx_{1}\right) $}

Using $\left( \ref{simplgx}\right) ,$ we get%
\begin{eqnarray}
&&B(x_{1}\otimes gx_{1})  \label{form x1otgx1} \\
&=&(1_{A}\otimes g)B(gx_{1}\otimes 1_{H})(1_{A}\otimes gx_{1})  \notag \\
&&+(1_{A}\otimes x_{1})B(gx_{1}\otimes 1_{H})  \notag \\
&&-B(gx_{1}gx_{1}\otimes 1_{H})  \notag
\end{eqnarray}%
We deduce that%
\begin{eqnarray*}
B\left( x_{1}\otimes gx_{1}\right) &=&2\left[ B(x_{1}x_{2}\otimes
1_{H};1_{A},gx_{2})+B(x_{1}x_{2}\otimes 1_{H};X_{2},g)\right] 1_{A}\otimes
gx_{1}+ \\
&&+2B(x_{1}x_{2}\otimes 1_{H};GX_{2},gx_{2})G\otimes gx_{1}x_{2}+ \\
&&+2B(x_{1}x_{2}\otimes 1_{H};X_{2},gx_{1}x_{2})X_{1}\otimes gx_{1}x_{2}+ \\
&&+2B(x_{1}x_{2}\otimes 1_{H};X_{2},gx_{1}x_{2})X_{1}X_{2}\otimes gx_{1}+ \\
&&+2\left[ -B(x_{1}x_{2}\otimes 1_{H};G,gx_{1}x_{2})+B(x_{1}x_{2}\otimes
1_{H};GX_{2},gx_{1})\right] GX_{1}\otimes gx_{1}+ \\
&&+2B(x_{1}x_{2}\otimes 1_{H};GX_{2},gx_{2})GX_{2}\otimes gx_{1}+
\end{eqnarray*}%
We write the Casimir formula for $B(x_{1}\otimes gx_{1})$:
\begin{eqnarray*}
&&\sum_{w_{1}=0}^{1}\sum_{a,b_{1},b_{2},d,e_{1},e_{2}=0}^{1}%
\sum_{l_{1}=0}^{b_{1}}\sum_{l_{2}=0}^{b_{2}}\sum_{u_{1}=0}^{e_{1}}%
\sum_{u_{2}=0}^{e_{2}}\left( -1\right) ^{\alpha \left(
x_{1}^{1-w_{1}};l_{1},l_{2},u_{1},u_{2}\right) } \\
&&B(g^{1+w_{1}}x_{1}^{w_{1}}\otimes
gx_{1};G^{a}X_{1}^{b_{1}}X_{2}^{b_{2}},g^{d}x_{1}^{e_{1}}x_{2}^{e_{2}}) \\
&&G^{a}X_{1}^{b_{1}-l_{1}}X_{2}^{b_{2}-l_{2}}\otimes
g^{d}x_{1}^{e_{1}-u_{1}}x_{2}^{e_{2}-u_{2}}\otimes \\
&&g^{a+b_{1}+b_{2}+l_{1}+l_{2}+d+e_{1}+e_{2}+u_{1}+u_{2}}x_{1}^{l_{1}+u_{1}+1-w_{1}}x_{2}^{l_{2}+u_{2}}
\\
&=&\sum_{\omega _{1}=0}^{1}B^{A}(x_{1}\otimes gx_{1}^{1-\omega _{1}})\otimes
B^{H}(x_{1}\otimes gx_{1}^{1-\omega _{1}})\otimes g^{\omega
_{1}}x_{1}^{\omega _{1}}
\end{eqnarray*}%
Thus we get%
\begin{eqnarray*}
&&\sum_{a,b_{1},b_{2},d,e_{1},e_{2}=0}^{1}\sum_{l_{1}=0}^{b_{1}}%
\sum_{l_{2}=0}^{b_{2}}\sum_{u_{1}=0}^{e_{1}}\sum_{u_{2}=0}^{e_{2}}\left(
-1\right) ^{\alpha \left( x_{1};l_{1},l_{2},u_{1},u_{2}\right) } \\
&&B(g\otimes
gx_{1};G^{a}X_{1}^{b_{1}}X_{2}^{b_{2}},g^{d}x_{1}^{e_{1}}x_{2}^{e_{2}})G^{a}X_{1}^{b_{1}-l_{1}}X_{2}^{b_{2}-l_{2}}\otimes g^{d}x_{1}^{e_{1}-u_{1}}x_{2}^{e_{2}-u_{2}}\otimes
\\
&&g^{a+b_{1}+b_{2}+l_{1}+l_{2}+d+e_{1}+e_{2}+u_{1}+u_{2}}x_{1}^{l_{1}+u_{1}+1}x_{2}^{l_{2}+u_{2}}
\\
&&+\sum_{a,b_{1},b_{2},d,e_{1},e_{2}=0}^{1}\sum_{l_{1}=0}^{b_{1}}%
\sum_{l_{2}=0}^{b_{2}}\sum_{u_{1}=0}^{e_{1}}\sum_{u_{2}=0}^{e_{2}}\left(
-1\right) ^{\alpha \left( 1_{H};l_{1},l_{2},u_{1},u_{2}\right) } \\
&&B(x_{1}\otimes
gx_{1};G^{a}X_{1}^{b_{1}}X_{2}^{b_{2}},g^{d}x_{1}^{e_{1}}x_{2}^{e_{2}})G^{a}X_{1}^{b_{1}-l_{1}}X_{2}^{b_{2}-l_{2}}\otimes g^{d}x_{1}^{e_{1}-u_{1}}x_{2}^{e_{2}-u_{2}}\otimes
\\
&&g^{a+b_{1}+b_{2}+l_{1}+l_{2}+d+e_{1}+e_{2}+u_{1}+u_{2}}x_{1}^{l_{1}+u_{1}}x_{2}^{l_{2}+u_{2}}
\\
&=&B^{A}(x_{1}\otimes gx_{1})\otimes B^{H}(x_{1}\otimes gx_{1})\otimes 1_{H}+
\\
&&+B^{A}(x_{1}\otimes g)\otimes B^{H}(x_{1}\otimes g)\otimes gx_{1}+
\end{eqnarray*}

\subsection{$B(x_{1}\otimes gx_{1};1_{A},gx_{1})$}

Now we look for the terms in $B(x_{1}\otimes gx_{1};1_{A},gx_{1})$ in the
equality above. We get%
\begin{eqnarray*}
a &=&b_{1}=b_{2}=0 \\
d &=&e_{1}=1 \\
e_{2} &=&0
\end{eqnarray*}%
and hence%
\begin{eqnarray*}
&&\sum_{u_{1}=0}^{1}\left( -1\right) ^{\alpha \left(
1_{H};0,0,u_{1},0\right) }B(x_{1}\otimes gx_{1};1_{A},gx_{1})1_{A}\otimes
gx_{1}^{1-u_{1}}\otimes g^{u_{1}}x_{1}^{u_{1}} \\
= &&\left( -1\right) ^{\alpha \left( 1_{H};0,0,0,0\right) }B(x_{1}\otimes
gx_{1};1_{A},gx_{1})1_{A}\otimes gx_{1}\otimes 1_{H} \\
&&+\left( -1\right) ^{\alpha \left( 1_{H};0,0,1,0\right) }B(x_{1}\otimes
gx_{1};1_{A},gx_{1})1_{A}\otimes g\otimes gx_{1}
\end{eqnarray*}%
Thus we have to look for all the terms containing $1_{A}\otimes
gx_{1}\otimes 1_{H}$ and for those containing $1_{A}\otimes g\otimes gx_{1}.$

\subsubsection{Case $1_{A}\otimes gx_{1}\otimes 1_{H}$}

We look for terms $1_{A}\otimes gx_{1}\otimes 1_{H}$ in%
\begin{eqnarray*}
&&\sum_{a,b_{1},b_{2},d,e_{1},e_{2}=0}^{1}\sum_{l_{1}=0}^{b_{1}}%
\sum_{l_{2}=0}^{b_{2}}\sum_{u_{1}=0}^{e_{1}}\sum_{u_{2}=0}^{e_{2}}\left(
-1\right) ^{\alpha \left( x_{1};l_{1},l_{2},u_{1},u_{2}\right) } \\
&&B(g\otimes
gx_{1};G^{a}X_{1}^{b_{1}}X_{2}^{b_{2}},g^{d}x_{1}^{e_{1}}x_{2}^{e_{2}})G^{a}X_{1}^{b_{1}-l_{1}}X_{2}^{b_{2}-l_{2}}\otimes g^{d}x_{1}^{e_{1}-u_{1}}x_{2}^{e_{2}-u_{2}}\otimes
\\
&&g^{a+b_{1}+b_{2}+l_{1}+l_{2}+d+e_{1}+e_{2}+u_{1}+u_{2}}x_{1}^{l_{1}+u_{1}+1}x_{2}^{l_{2}+u_{2}}
\\
&&+\sum_{a,b_{1},b_{2},d,e_{1},e_{2}=0}^{1}\sum_{l_{1}=0}^{b_{1}}%
\sum_{l_{2}=0}^{b_{2}}\sum_{u_{1}=0}^{e_{1}}\sum_{u_{2}=0}^{e_{2}}\left(
-1\right) ^{\alpha \left( 1_{H};l_{1},l_{2},u_{1},u_{2}\right) } \\
&&B(x_{1}\otimes
gx_{1};G^{a}X_{1}^{b_{1}}X_{2}^{b_{2}},g^{d}x_{1}^{e_{1}}x_{2}^{e_{2}})G^{a}X_{1}^{b_{1}-l_{1}}X_{2}^{b_{2}-l_{2}}\otimes g^{d}x_{1}^{e_{1}-u_{1}}x_{2}^{e_{2}-u_{2}}\otimes
\\
&&g^{a+b_{1}+b_{2}+l_{1}+l_{2}+d+e_{1}+e_{2}+u_{1}+u_{2}}x_{1}^{l_{1}+u_{1}}x_{2}^{l_{2}+u_{2}}
\\
&=&B^{A}(x_{1}\otimes gx_{1})\otimes B^{H}(x_{1}\otimes gx_{1})\otimes 1_{H}
\\
&&+B^{A}(x_{1}\otimes g)\otimes B^{H}(x_{1}\otimes g)\otimes gx_{1}+
\end{eqnarray*}%
In the first summand of the left side of the equation there is no such a
term.

Looking to the second summand of the left side of the equation we get

\begin{eqnarray*}
l_{1} &=&u_{1}=0 \\
l_{2} &=&u_{2}=0 \\
a+b_{1}+b_{2}+d+e_{1}+e_{2} &=&0 \\
d &=&1 \\
e_{1} &=&1 \\
e_{2} &=&u_{2}=0 \\
a &=&0 \\
b_{1} &=&l_{1}=0 \\
b_{2} &=&l_{2}=0
\end{eqnarray*}%
and we get, since $\alpha \left( 1_{H};0,0,0,0\right) \equiv 0$%
\begin{equation*}
B(x_{1}\otimes gx_{1};1_{A},gx_{1})1_{A}\otimes gx_{1}\otimes 1_{H}
\end{equation*}%
Looking to the right side of the equation we get%
\begin{equation*}
B\left( x_{1}\otimes gx_{1},1_{A},gx_{1}\right) 1_{A}\otimes gx_{1}\otimes
1_{H}
\end{equation*}%
Thus we get%
\begin{equation*}
B(x_{1}\otimes gx_{1};1_{A},gx_{1})-B\left( x_{1}\otimes
gx_{1},1_{A},gx_{1}\right) =0
\end{equation*}%
which is trivial.

\subsubsection{Case $1_{A}\otimes g\otimes gx_{1}$}

Proceeding in a similar way as above we get, first summand left side of the
equation%
\begin{eqnarray*}
l_{2} &=&u_{2}=0 \\
l_{1} &=&u_{1}=0 \\
a+b_{1}+b_{2}+d+e_{1}+e_{2} &\equiv &1 \\
d &=&1 \\
e_{1} &=&u_{1}=0 \\
e_{2} &=&u_{2}=0 \\
a &=&0 \\
b_{1} &=&u_{1}=0 \\
b_{2} &=&u_{2}=0
\end{eqnarray*}%
Thus we get
\begin{equation*}
B(g\otimes gx_{1};1_{A},g)1_{A}\otimes g\otimes gx_{1}
\end{equation*}

Second summand of the left side of the equation

\begin{eqnarray*}
l_{1}+u_{1} &=&1 \\
l_{2} &=&u_{2}=0 \\
a+b_{1}+b_{2}+d+e_{1}+e_{2} &\equiv &0 \\
d &=&1 \\
e_{1} &=&u_{1} \\
e_{2} &=&u_{2}=0 \\
a &=&0 \\
b_{1} &=&l_{1} \\
b_{2} &=&l_{2}=0
\end{eqnarray*}%
Since%
\begin{eqnarray*}
\alpha \left( 1_{H};0,0,1,0\right) &\equiv &0 \\
\alpha \left( 1_{H};1,0,0,0\right) &\equiv &0
\end{eqnarray*}%
and we get%
\begin{eqnarray*}
&&+\sum_{\substack{ b_{1},e_{1}=0  \\ b_{1}+e_{1}\equiv 1}}^{1}\left(
-1\right) ^{\alpha \left( 1_{H};b_{1},0,e_{1},0\right) }B(x_{1}\otimes
gx_{1};X_{1}^{b_{1}},gx_{1}^{e_{1}})1_{A}\otimes g\otimes gx_{1} \\
&=&\left[ B(x_{1}\otimes gx_{1};1_{A},gx_{1})+B(x_{1}\otimes gx_{1};X_{1},g)%
\right] 1_{A}\otimes g\otimes gx_{1}
\end{eqnarray*}%
In the second side we get%
\begin{equation*}
B(x_{1}\otimes g;1_{A},g)1_{A}\otimes g\otimes gx_{1}
\end{equation*}%
Thus we obtain the equality%
\begin{equation*}
B(g\otimes gx_{1};1_{A},g)+B(x_{1}\otimes
gx_{1};1_{A},gx_{1})+B(x_{1}\otimes gx_{1};X_{1},g)-B(x_{1}\otimes
g;1_{A},g)=0
\end{equation*}%
By the form of the elements this equality holds.

\subsection{$B(x_{1}\otimes gx_{1};G,gx_{1}x_{2})$}

Now we look for the terms in $B(x_{1}\otimes gx_{1};G,gx_{1}x_{2})$ in the
equality
\begin{eqnarray*}
&&\sum_{a,b_{1},b_{2},d,e_{1},e_{2}=0}^{1}\sum_{l_{1}=0}^{b_{1}}%
\sum_{l_{2}=0}^{b_{2}}\sum_{u_{1}=0}^{e_{1}}\sum_{u_{2}=0}^{e_{2}}\left(
-1\right) ^{\alpha \left( x_{1};l_{1},l_{2},u_{1},u_{2}\right) } \\
&&B(g\otimes
gx_{1};G^{a}X_{1}^{b_{1}}X_{2}^{b_{2}},g^{d}x_{1}^{e_{1}}x_{2}^{e_{2}})G^{a}X_{1}^{b_{1}-l_{1}}X_{2}^{b_{2}-l_{2}}\otimes g^{d}x_{1}^{e_{1}-u_{1}}x_{2}^{e_{2}-u_{2}}\otimes
\\
&&g^{a+b_{1}+b_{2}+l_{1}+l_{2}+d+e_{1}+e_{2}+u_{1}+u_{2}}x_{1}^{l_{1}+u_{1}+1}x_{2}^{l_{2}+u_{2}}
\\
&&+\sum_{a,b_{1},b_{2},d,e_{1},e_{2}=0}^{1}\sum_{l_{1}=0}^{b_{1}}%
\sum_{l_{2}=0}^{b_{2}}\sum_{u_{1}=0}^{e_{1}}\sum_{u_{2}=0}^{e_{2}}\left(
-1\right) ^{\alpha \left( 1_{H};l_{1},l_{2},u_{1},u_{2}\right) } \\
&&B(x_{1}\otimes
gx_{1};G^{a}X_{1}^{b_{1}}X_{2}^{b_{2}},g^{d}x_{1}^{e_{1}}x_{2}^{e_{2}})G^{a}X_{1}^{b_{1}-l_{1}}X_{2}^{b_{2}-l_{2}}\otimes g^{d}x_{1}^{e_{1}-u_{1}}x_{2}^{e_{2}-u_{2}}\otimes
\\
&&g^{a+b_{1}+b_{2}+l_{1}+l_{2}+d+e_{1}+e_{2}+u_{1}+u_{2}}x_{1}^{l_{1}+u_{1}}x_{2}^{l_{2}+u_{2}}
\\
&=&B^{A}(x_{1}\otimes gx_{1})\otimes B^{H}(x_{1}\otimes gx_{1})\otimes 1_{H}+
\\
&&+B^{A}(x_{1}\otimes g)\otimes B^{H}(x_{1}\otimes g)\otimes gx_{1}+
\end{eqnarray*}%
We get%
\begin{gather*}
a=1,b_{1}=0=b_{2},d=e_{1}=e_{2}=1 \\
a+b_{1}+b_{2}+l_{1}+l_{2}+d+e_{1}+e_{2}+u_{1}+u_{2}\equiv u_{1}+u_{2}
\end{gather*}%
and hence we get%
\begin{eqnarray*}
&&\sum_{u_{1}=0}^{1}\sum_{u_{2}=0}^{1}\left( -1\right) ^{\alpha \left(
1_{H};0,0,u_{1},u_{2}\right) }B(x_{1}\otimes gx_{1};G,gx_{1}x_{2})G\otimes
gx_{1}^{1-u_{1}}x_{2}^{1-u_{2}}\otimes
g^{u_{1}+u_{2}}x_{1}^{u_{1}}x_{2}^{u_{2}} \\
&=&\left( -1\right) ^{\alpha \left( 1_{H};0,0,0,0\right) }B(x_{1}\otimes
gx_{1};G,gx_{1}x_{2})G\otimes gx_{1}x_{2}\otimes 1_{H} \\
&&\left( -1\right) ^{\alpha \left( 1_{H};0,0,0,1\right) }B(x_{1}\otimes
gx_{1};G,gx_{1}x_{2})G\otimes gx_{1}\otimes gx_{2} \\
&&\left( -1\right) ^{\alpha \left( 1_{H};0,0,1,0\right) }B(x_{1}\otimes
gx_{1};G,gx_{1}x_{2})G\otimes gx_{2}\otimes gx_{1} \\
&&\left( -1\right) ^{\alpha \left( 1_{H};0,0,1,1\right) }B(x_{1}\otimes
gx_{1};G,gx_{1}x_{2})G\otimes g\otimes x_{1}x_{2}
\end{eqnarray*}

\subsubsection{Case $G\otimes gx_{1}x_{2}\otimes 1_{H}$}

By looking for terms in $G\otimes gx_{1}x_{2}\otimes 1_{H}$ in the first
summand of our equality we realize that \ there are no such a term.

By looking for terms in $G\otimes gx_{1}x_{2}\otimes 1_{H}$ in the second
summand of our equality we get
\begin{eqnarray*}
&&+\sum_{a,b_{1},b_{2},d,e_{1},e_{2}=0}^{1}\sum_{l_{1}=0}^{b_{1}}%
\sum_{l_{2}=0}^{b_{2}}\sum_{u_{1}=0}^{e_{1}}\sum_{u_{2}=0}^{e_{2}}\left(
-1\right) ^{\alpha \left( 1_{H};l_{1},l_{2},u_{1},u_{2}\right) } \\
&&B(x_{1}\otimes
gx_{1};G^{a}X_{1}^{b_{1}}X_{2}^{b_{2}},g^{d}x_{1}^{e_{1}}x_{2}^{e_{2}})G^{a}X_{1}^{b_{1}-l_{1}}X_{2}^{b_{2}-l_{2}}\otimes g^{d}x_{1}^{e_{1}-u_{1}}x_{2}^{e_{2}-u_{2}}\otimes
\\
&&g^{a+b_{1}+b_{2}+l_{1}+l_{2}+d+e_{1}+e_{2}+u_{1}+u_{2}}x_{1}^{l_{1}+u_{1}}x_{2}^{l_{2}+u_{2}}
\end{eqnarray*}%
\begin{eqnarray*}
l_{1} &=&u_{1}=0 \\
l_{2} &=&u_{2}=0 \\
a+b_{1}+b_{2}+d+e_{1}+e_{2} &\equiv &0 \\
d &=&1,e_{1}=1,e_{2}=1 \\
a &=&1,b_{1}=l_{1}=0,b_{2}=l_{2}=0
\end{eqnarray*}%
Since $\alpha \left( 1_{H};0,0,0,0\right) =0$ we get%
\begin{equation*}
B(x_{1}\otimes gx_{1};G,gx_{1}x_{2})G\otimes gx_{1}x_{2}\otimes 1_{H}
\end{equation*}%
Doing the same for the right side of our equality we get%
\begin{equation*}
B(x_{1}\otimes gx_{1};G,gx_{1}x_{2})G\otimes gx_{1}x_{2}\otimes
1_{H}=B(x_{1}\otimes gx_{1};G,gx_{1}x_{2})G\otimes gx_{1}x_{2}\otimes 1_{H}
\end{equation*}%
so we get no information.

\subsubsection{Case $G\otimes gx_{1}\otimes gx_{2}$}

Proceeding as above, in the first summand there is no term. In the second
summand%
\begin{eqnarray*}
&&+\sum_{a,b_{1},b_{2},d,e_{1},e_{2}=0}^{1}\sum_{l_{1}=0}^{b_{1}}%
\sum_{l_{2}=0}^{b_{2}}\sum_{u_{1}=0}^{e_{1}}\sum_{u_{2}=0}^{e_{2}}\left(
-1\right) ^{\alpha \left( 1_{H};l_{1},l_{2},u_{1},u_{2}\right) } \\
&&B(x_{1}\otimes
gx_{1};G^{a}X_{1}^{b_{1}}X_{2}^{b_{2}},g^{d}x_{1}^{e_{1}}x_{2}^{e_{2}})G^{a}X_{1}^{b_{1}-l_{1}}X_{2}^{b_{2}-l_{2}}\otimes g^{d}x_{1}^{e_{1}-u_{1}}x_{2}^{e_{2}-u_{2}}\otimes
\\
&&g^{a+b_{1}+b_{2}+l_{1}+l_{2}+d+e_{1}+e_{2}+u_{1}+u_{2}}x_{1}^{l_{1}+u_{1}}x_{2}^{l_{2}+u_{2}}
\end{eqnarray*}%
\begin{eqnarray*}
l_{1} &=&u_{1}=0 \\
l_{2}+u_{2} &=&1 \\
a+b_{1}+b_{2}+d+e_{1}+e_{2} &\equiv &0 \\
d &=&1,e_{1}=1,e_{2}=u_{2} \\
a &=&1,b_{1}=l_{1}=0,b_{2}=l_{2}
\end{eqnarray*}%
so we get%
\begin{eqnarray*}
&&+\sum_{\substack{ b_{2},e_{2}  \\ b_{2}+e_{2}=1}}^{1}%
\sum_{l_{2}=0}^{b_{2}}\sum_{u_{2}=0}^{e_{2}}\left( -1\right) ^{\alpha \left(
1_{H};0,b_{2},0,e_{2}\right) }B(x_{1}\otimes
gx_{1};GX_{2}^{b_{2}},gx_{1}x_{2}^{e_{2}})G\otimes gx_{1}\otimes gx_{2} \\
&=&\left[
\begin{array}{c}
\left( -1\right) ^{\alpha \left( 1_{H};0,0,0,1\right) }B(x_{1}\otimes
gx_{1};G,gx_{1}x_{2}) \\
+\left( -1\right) ^{\alpha \left( 1_{H};0,1,0,0\right) }B(x_{1}\otimes
gx_{1};GX_{2},gx_{1})%
\end{array}%
\right] G\otimes gx_{1}\otimes gx_{2}
\end{eqnarray*}

As

\begin{eqnarray*}
&&\alpha \left( 1_{H};0,0,0,1\right) =1 \\
&&\alpha \left( 1_{H};0,1,0,0\right) =0
\end{eqnarray*}%
we get%
\begin{equation*}
\left[ -B(x_{1}\otimes gx_{1};G,gx_{1}x_{2})+B(x_{1}\otimes
gx_{1};GX_{2},gx_{1})\right] G\otimes gx_{1}\otimes gx_{2}
\end{equation*}%
Since there is no term on the right side, we obtain%
\begin{equation*}
-B(x_{1}\otimes gx_{1};G,gx_{1}x_{2})+B(x_{1}\otimes gx_{1};GX_{2},gx_{1})=0
\end{equation*}%
which holds in view of the form of the element.

\subsubsection{Case $G\otimes gx_{2}\otimes gx_{1}$}

First summand%
\begin{eqnarray*}
&&\sum_{a,b_{1},b_{2},d,e_{1},e_{2}=0}^{1}\sum_{l_{1}=0}^{b_{1}}%
\sum_{l_{2}=0}^{b_{2}}\sum_{u_{1}=0}^{e_{1}}\sum_{u_{2}=0}^{e_{2}}\left(
-1\right) ^{\alpha \left( x_{1};l_{1},l_{2},u_{1},u_{2}\right) } \\
&&B(g\otimes
gx_{1};G^{a}X_{1}^{b_{1}}X_{2}^{b_{2}},g^{d}x_{1}^{e_{1}}x_{2}^{e_{2}})G^{a}X_{1}^{b_{1}-l_{1}}X_{2}^{b_{2}-l_{2}}\otimes g^{d}x_{1}^{e_{1}-u_{1}}x_{2}^{e_{2}-u_{2}}\otimes
\\
&&g^{a+b_{1}+b_{2}+l_{1}+l_{2}+d+e_{1}+e_{2}+u_{1}+u_{2}}x_{1}^{l_{1}+u_{1}+1}x_{2}^{l_{2}+u_{2}}
\end{eqnarray*}

\begin{eqnarray*}
l_{2} &=&u_{2}=0 \\
l_{1} &=&u_{1}=0 \\
a+b_{1}+b_{2}+d+e_{1}+e_{2} &\equiv &1 \\
d &=&1,e_{1}=u_{1}=0,e_{2}=1 \\
a &=&1,b_{1}=u_{1}=0,b_{2}=u_{2}=0
\end{eqnarray*}%
\begin{equation*}
\alpha \left( x_{1};0,0,0,0\right) \equiv \left( 1+0+0\right)
\end{equation*}%
\begin{equation*}
\left( -1\right) B(g\otimes gx_{1};G,gx_{2})G\otimes gx_{2}\otimes gx_{1}
\end{equation*}%
Second summand

\begin{eqnarray*}
&&+\sum_{a,b_{1},b_{2},d,e_{1},e_{2}=0}^{1}\sum_{l_{1}=0}^{b_{1}}%
\sum_{l_{2}=0}^{b_{2}}\sum_{u_{1}=0}^{e_{1}}\sum_{u_{2}=0}^{e_{2}}\left(
-1\right) ^{\alpha \left( 1_{H};l_{1},l_{2},u_{1},u_{2}\right) } \\
&&B(x_{1}\otimes
gx_{1};G^{a}X_{1}^{b_{1}}X_{2}^{b_{2}},g^{d}x_{1}^{e_{1}}x_{2}^{e_{2}})G^{a}X_{1}^{b_{1}-l_{1}}X_{2}^{b_{2}-l_{2}}\otimes g^{d}x_{1}^{e_{1}-u_{1}}x_{2}^{e_{2}-u_{2}}\otimes
\\
&&g^{a+b_{1}+b_{2}+l_{1}+l_{2}+d+e_{1}+e_{2}+u_{1}+u_{2}}x_{1}^{l_{1}+u_{1}}x_{2}^{l_{2}+u_{2}}
\end{eqnarray*}%
\begin{eqnarray*}
l_{1}+u_{1} &=&1 \\
l_{2} &=&u_{2}=0 \\
a+b_{1}+b_{2}+d+e_{1}+e_{2} &\equiv &0 \\
d &=&1,e_{1}=u_{1},e_{2}=1 \\
a &=&1,b_{1}=l_{1},b_{2}=l_{2}=0
\end{eqnarray*}

\begin{eqnarray*}
&&+\sum_{\substack{ b_{1},e_{1}=0  \\ b_{1}+e_{1}=1}}^{1}\left( -1\right)
^{\alpha \left( 1_{H};b_{1},0,e_{1},0\right) }B(x_{1}\otimes
gx_{1};GX_{1}^{b_{1}},gx_{1}^{e_{1}}x_{2})G\otimes gx_{2}\otimes gx_{1} \\
&=&\left[ \left( -1\right) ^{\alpha \left( 1_{H};0,0,1,0\right)
}B(x_{1}\otimes gx_{1};G,gx_{1}x_{2})+\left( -1\right) ^{\alpha \left(
1_{H};1,0,0,0\right) }B(x_{1}\otimes gx_{1};GX_{1},gx_{2})\right] G\otimes
gx_{2}\otimes gx_{1}
\end{eqnarray*}%
We get%
\begin{eqnarray*}
&&\alpha \left( 1_{H};0,0,1,0\right) =e_{2}+\left( a+b_{1}+b_{2}\right)
=1+\left( 1+0+0\right) =0 \\
&&\alpha \left( 1_{H};1,0,0,0\right) =b_{2}=0
\end{eqnarray*}%
\begin{equation*}
\left[ B(x_{1}\otimes gx_{1};G,gx_{1}x_{2})+B(x_{1}\otimes
gx_{1};GX_{1},gx_{2})\right] G\otimes gx_{2}\otimes gx_{1}
\end{equation*}%
Thus, by considering also the right side of the equation, we obtain%
\begin{equation*}
B(x_{1}\otimes gx_{1};G,gx_{1}x_{2})+B(x_{1}\otimes
gx_{1};GX_{1},gx_{2})-B(g\otimes gx_{1};G,gx_{2})-B(x_{1}\otimes
g;G,gx_{2})=0
\end{equation*}

which holds in view of the form of the elements.

\subsubsection{Case $G\otimes g\otimes x_{1}x_{2}$}

First summand%
\begin{eqnarray*}
&&\sum_{a,b_{1},b_{2},d,e_{1},e_{2}=0}^{1}\sum_{l_{1}=0}^{b_{1}}%
\sum_{l_{2}=0}^{b_{2}}\sum_{u_{1}=0}^{e_{1}}\sum_{u_{2}=0}^{e_{2}}\left(
-1\right) ^{\alpha \left( x_{1};l_{1},l_{2},u_{1},u_{2}\right) } \\
&&B(g\otimes
gx_{1};G^{a}X_{1}^{b_{1}}X_{2}^{b_{2}},g^{d}x_{1}^{e_{1}}x_{2}^{e_{2}})G^{a}X_{1}^{b_{1}-l_{1}}X_{2}^{b_{2}-l_{2}}\otimes g^{d}x_{1}^{e_{1}-u_{1}}x_{2}^{e_{2}-u_{2}}\otimes
\\
&&g^{a+b_{1}+b_{2}+l_{1}+l_{2}+d+e_{1}+e_{2}+u_{1}+u_{2}}x_{1}^{l_{1}+u_{1}+1}x_{2}^{l_{2}+u_{2}}
\end{eqnarray*}

\begin{eqnarray*}
l_{2}+u_{2} &=&1 \\
l_{1} &=&u_{1}=0 \\
a+b_{1}+b_{2}+d+e_{1}+e_{2} &\equiv &1 \\
d &=&1,e_{1}=u_{1}=0,e_{2}=u_{2} \\
a &=&1,b_{1}=l_{1}=0,b_{2}=l_{2}
\end{eqnarray*}%
\begin{eqnarray*}
&&\sum_{\substack{ b_{2},e_{2}=0  \\ b_{2}+e_{2}=1}}^{1}\left( -1\right)
^{\alpha \left( x_{1};0,b_{2},0,e_{2}\right) }B(g\otimes
gx_{1};GX_{2}^{b_{2}},gx_{2}^{e_{2}})G\otimes g\otimes x_{1}x_{2} \\
&=&\left[
\begin{array}{c}
\left( -1\right) ^{\alpha \left( x_{1};0,0,0,1\right) }B(g\otimes
gx_{1};G,gx_{2}) \\
+\left( -1\right) ^{\alpha \left( x_{1};0,1,0,0\right) }B(g\otimes
gx_{1};GX_{2},g)%
\end{array}%
\right] G\otimes g\otimes x_{1}x_{2}
\end{eqnarray*}%
\begin{eqnarray*}
\alpha \left( x_{1};0,0,0,1\right) &\equiv &1 \\
\alpha \left( x_{1};0,1,0,0\right) &=&a+b_{1}+b_{2}+1\equiv 1+0+1+1\equiv 1
\end{eqnarray*}%
\begin{equation*}
\left[ -B(g\otimes gx_{1};G,gx_{2})-B(g\otimes gx_{1};GX_{2},g)\right]
G\otimes g\otimes x_{1}x_{2}
\end{equation*}%
Second summand%
\begin{eqnarray*}
&&+\sum_{a,b_{1},b_{2},d,e_{1},e_{2}=0}^{1}\sum_{l_{1}=0}^{b_{1}}%
\sum_{l_{2}=0}^{b_{2}}\sum_{u_{1}=0}^{e_{1}}\sum_{u_{2}=0}^{e_{2}}\left(
-1\right) ^{\alpha \left( 1_{H};l_{1},l_{2},u_{1},u_{2}\right) } \\
&&B(x_{1}\otimes
gx_{1};G^{a}X_{1}^{b_{1}}X_{2}^{b_{2}},g^{d}x_{1}^{e_{1}}x_{2}^{e_{2}})G^{a}X_{1}^{b_{1}-l_{1}}X_{2}^{b_{2}-l_{2}}\otimes g^{d}x_{1}^{e_{1}-u_{1}}x_{2}^{e_{2}-u_{2}}\otimes
\\
&&g^{a+b_{1}+b_{2}+l_{1}+l_{2}+d+e_{1}+e_{2}+u_{1}+u_{2}}x_{1}^{l_{1}+u_{1}}x_{2}^{l_{2}+u_{2}}
\end{eqnarray*}%
\begin{eqnarray*}
l_{1}+u_{1} &=&1 \\
l_{2}+u_{2} &=&1 \\
a+b_{1}+b_{2}+d+e_{1}+e_{2} &\equiv &0 \\
d &=&1 \\
e_{1} &=&u_{1} \\
e_{2} &=&u_{2} \\
a &=&1 \\
b_{1} &=&l_{1} \\
b_{2} &=&l_{2}
\end{eqnarray*}%
\begin{eqnarray*}
&&+\sum_{\substack{ b_{1},b_{2},e_{1},e_{2}=0  \\ b_{1}+e_{1}=1  \\ %
b_{2}+e_{2}=1}}^{1}\left( -1\right) ^{\alpha \left(
1_{H};b_{1},b_{2},e_{1},e_{2}\right) }B(x_{1}\otimes
gx_{1};GX_{1}^{b_{1}}X_{2}^{b_{2}},gx_{1}^{e_{1}}x_{2}^{e_{2}})G\otimes
g\otimes x_{1}x_{2} \\
&&\left( -1\right) ^{\alpha \left( 1_{H};0,0,1,1\right) }B(x_{1}\otimes
gx_{1};G,gx_{1}x_{2})G\otimes g\otimes x_{1}x_{2} \\
&&\left( -1\right) ^{\alpha \left( 1_{H};1,0,0,1\right) }B(x_{1}\otimes
gx_{1};GX_{1},gx_{2})G\otimes g\otimes x_{1}x_{2} \\
&&\left( -1\right) ^{\alpha \left( 1_{H};0,1,1,0\right) }B(x_{1}\otimes
gx_{1};GX_{2},gx_{1})G\otimes g\otimes x_{1}x_{2} \\
&&\left( -1\right) ^{\alpha \left( 1_{H};1,1,0,0\right) }B(x_{1}\otimes
gx_{1};GX_{1}X_{2},g)G\otimes g\otimes x_{1}x_{2}
\end{eqnarray*}%
\begin{eqnarray*}
&&\alpha \left( 1_{H};0,0,1,1\right) \equiv 1+e_{2}\equiv 1+1\equiv 0 \\
&&\alpha \left( 1_{H};1,0,0,1\right) \equiv a+b_{1}\equiv 1+1\equiv 0 \\
&&\alpha \left( 1_{H};0,1,1,0\right) \equiv e_{2}+\left(
a+b_{1}+b_{2}+1\right) \equiv 0+1+0+1+1\equiv 1 \\
&&\alpha \left( 1_{H};1,1,0,0\right) \equiv 1+b_{2}\equiv 0
\end{eqnarray*}%
and we get%
\begin{equation*}
\left[
\begin{array}{c}
B(x_{1}\otimes gx_{1};G,gx_{1}x_{2})+B(x_{1}\otimes gx_{1};GX_{1},gx_{2}) \\
-B(x_{1}\otimes gx_{1};GX_{2},gx_{1})+B(x_{1}\otimes gx_{1};GX_{1}X_{2},g)%
\end{array}%
\right] G\otimes g\otimes x_{1}x_{2}
\end{equation*}%
Since there is nothing on the right side we get%
\begin{gather*}
-B(g\otimes gx_{1};G,gx_{2})-B(g\otimes gx_{1};GX_{2},g)+B(x_{1}\otimes
gx_{1};G,gx_{1}x_{2})+ \\
B(x_{1}\otimes gx_{1};GX_{1},gx_{2})-B(x_{1}\otimes
gx_{1};GX_{2},gx_{1})+B(x_{1}\otimes gx_{1};GX_{1}X_{2},g)=0
\end{gather*}

which holds in view of the form of the elements.

\subsection{$B(x_{1}\otimes gx_{1};X_{1},gx_{1}x_{2})$}

We deduce that%
\begin{gather*}
a=0,b_{1}=1,b_{2}=0,d=1,e_{1}=e_{2}=1. \\
a+b_{1}+b_{2}+l_{1}+l_{2}+d+e_{1}+e_{2}+u_{1}+u_{2}\equiv l_{1}+u_{1}+u_{2}
\end{gather*}%
and we get%
\begin{gather*}
\sum_{a,b_{1},b_{2},d,e_{1},e_{2}=0}^{1}\sum_{l_{1}=0}^{b_{1}}%
\sum_{u_{1}=0}^{e_{1}}\sum_{u_{2}=0}^{e_{2}}\left( -1\right) ^{\alpha \left(
1_{H};l_{1},0,u_{1},u_{2}\right) }B(x_{1}\otimes gx_{1};X_{1},gx_{1}x_{2}) \\
X_{1}^{1-l_{1}}\otimes gx_{1}^{1-u_{1}}x_{2}^{1-u_{2}}\otimes
g^{l_{1}+u_{1}+u_{2}}x_{1}^{l_{1}+u_{1}}x_{2}^{u_{2}}= \\
=\left( -1\right) ^{\alpha \left( 1_{H};0,0,0,0\right) }B(x_{1}\otimes
gx_{1};X_{1},gx_{1}x_{2})X_{1}\otimes gx_{1}x_{2}\otimes 1_{H}+ \\
\left( -1\right) ^{\alpha \left( 1_{H};0,0,0,1\right) }B(x_{1}\otimes
gx_{1};X_{1},gx_{1}x_{2})X_{1}\otimes gx_{1}\otimes gx_{2}+ \\
\left( -1\right) ^{\alpha \left( 1_{H};0,0,1,0\right) }B(x_{1}\otimes
gx_{1};X_{1},gx_{1}x_{2})X_{1}\otimes gx_{2}\otimes gx_{1}+ \\
\left( -1\right) ^{\alpha \left( 1_{H};0,0,1,1\right) }B(x_{1}\otimes
gx_{1};X_{1},gx_{1}x_{2})X_{1}\otimes g\otimes x_{1}x_{2}+ \\
\left( -1\right) ^{\alpha \left( 1_{H};1,0,0,0\right) }B(x_{1}\otimes
gx_{1};X_{1},gx_{1}x_{2})1_{A}\otimes gx_{1}x_{2}\otimes gx_{1}+ \\
\left( -1\right) ^{\alpha \left( 1_{H};1,0,0,1\right) }B(x_{1}\otimes
gx_{1};X_{1},gx_{1}x_{2})1_{A}\otimes gx_{1}\otimes x_{1}x_{2}+ \\
\left( -1\right) ^{\alpha \left( 1_{H};1,0,1,0\right) }B(x_{1}\otimes
gx_{1};X_{1},gx_{1}x_{2})1_{A}\otimes gx_{2}^{1-u_{2}}\otimes 0+ \\
\left( -1\right) ^{\alpha \left( 1_{H};1,0,1,1\right) }B(x_{1}\otimes
gx_{1};X_{1},gx_{1}x_{2})1_{A}\otimes gx_{1}^{1-u_{1}}x_{2}^{1-u_{2}}\otimes
0+
\end{gather*}

\subsubsection{Case $X_{1}\otimes gx_{1}x_{2}\otimes 1_{H}$}

First summand nothing.

Second summand%
\begin{eqnarray*}
&&+\sum_{a,b_{1},b_{2},d,e_{1},e_{2}=0}^{1}\sum_{l_{1}=0}^{b_{1}}%
\sum_{l_{2}=0}^{b_{2}}\sum_{u_{1}=0}^{e_{1}}\sum_{u_{2}=0}^{e_{2}}\left(
-1\right) ^{\alpha \left( 1_{H};l_{1},l_{2},u_{1},u_{2}\right) } \\
&&B(x_{1}\otimes
gx_{1};G^{a}X_{1}^{b_{1}}X_{2}^{b_{2}},g^{d}x_{1}^{e_{1}}x_{2}^{e_{2}})G^{a}X_{1}^{b_{1}-l_{1}}X_{2}^{b_{2}-l_{2}}\otimes g^{d}x_{1}^{e_{1}-u_{1}}x_{2}^{e_{2}-u_{2}}\otimes
\\
&&g^{a+b_{1}+b_{2}+l_{1}+l_{2}+d+e_{1}+e_{2}+u_{1}+u_{2}}x_{1}^{l_{1}+u_{1}}x_{2}^{l_{2}+u_{2}}
\end{eqnarray*}%
\begin{eqnarray*}
l_{1} &=&u_{1}=0 \\
l_{2} &=&u_{2}=0 \\
a+b_{1}+b_{2}+d+e_{1}+e_{2} &\equiv &0 \\
d &=&1 \\
e_{1} &=&1 \\
e_{2} &=&1 \\
a &=&0 \\
b_{1} &=&1 \\
b_{2} &=&l_{2}=0
\end{eqnarray*}%
Considering the right side we get an obvious equality.%
\begin{equation*}
B(x_{1}\otimes gx_{1};X_{1},gx_{1}x_{2})X_{1}\otimes gx_{1}x_{2}\otimes
1_{H}=B(x_{1}\otimes gx_{1};X_{1},gx_{1}x_{2})X_{1}\otimes
gx_{1}x_{2}\otimes 1_{H}
\end{equation*}

\subsubsection{Case $X_{1}\otimes gx_{1}\otimes gx_{2}$}

Nothing in the first summand. Second summand%
\begin{eqnarray*}
&&+\sum_{\substack{ b_{2},e_{2}=0  \\ b_{2}+e_{2}=1}}^{1}\left( -1\right)
^{\alpha \left( 1_{H};0,b_{2},0,e_{2}\right) }B(x_{1}\otimes
gx_{1};X_{1}X_{2}^{b_{2}},gx_{1}x_{2}^{e_{2}})X_{1}\otimes gx_{1}\otimes
gx_{2} \\
&=&\left[
\begin{array}{c}
\left( -1\right) ^{\alpha \left( 1_{H};0,0,0,1\right) }B(x_{1}\otimes
gx_{1};X_{1},gx_{1}x_{2}) \\
+\left( -1\right) ^{\alpha \left( 1_{H};0,1,0,0\right) }B(x_{1}\otimes
gx_{1};X_{1}X_{2},gx_{1})%
\end{array}%
\right] \\
&&X_{1}\otimes gx_{1}\otimes gx_{2}
\end{eqnarray*}%
\begin{eqnarray*}
&&\alpha \left( 1_{H};0,0,0,1\right) =a+b_{1}+b_{2}=0+1+0=1 \\
&&\alpha \left( 1_{H};0,1,0,0\right) =0
\end{eqnarray*}%
\begin{equation*}
-B(x_{1}\otimes gx_{1};X_{1},gx_{1}x_{2})+B(x_{1}\otimes
gx_{1};X_{1}X_{2},gx_{1})=0
\end{equation*}%
which holds in view of the form of the element.

\subsubsection{Case $X_{1}\otimes gx_{2}\otimes gx_{1}$}

\begin{eqnarray*}
&&\sum_{a,b_{1},b_{2},d,e_{1},e_{2}=0}^{1}\sum_{l_{1}=0}^{b_{1}}%
\sum_{l_{2}=0}^{b_{2}}\sum_{u_{1}=0}^{e_{1}}\sum_{u_{2}=0}^{e_{2}}\left(
-1\right) ^{\alpha \left( x_{1};l_{1},l_{2},u_{1},u_{2}\right) } \\
&&B(g\otimes
gx_{1};G^{a}X_{1}^{b_{1}}X_{2}^{b_{2}},g^{d}x_{1}^{e_{1}}x_{2}^{e_{2}})G^{a}X_{1}^{b_{1}-l_{1}}X_{2}^{b_{2}-l_{2}}\otimes g^{d}x_{1}^{e_{1}-u_{1}}x_{2}^{e_{2}-u_{2}}\otimes
\\
&&g^{a+b_{1}+b_{2}+l_{1}+l_{2}+d+e_{1}+e_{2}+u_{1}+u_{2}}x_{1}^{l_{1}+u_{1}+1}x_{2}^{l_{2}+u_{2}}
\end{eqnarray*}

First summand%
\begin{eqnarray*}
l_{1} &=&u_{1}=0 \\
l_{2} &=&u_{2}=0 \\
a+b_{1}+b_{2}+d+e_{1}+e_{2} &\equiv &1 \\
d &=&1 \\
e_{1} &=&u_{1}=0 \\
e_{2} &=&1 \\
a &=&0 \\
b_{1} &=&1 \\
b_{2} &=&l_{2}=0
\end{eqnarray*}%
and we get%
\begin{equation*}
\left( -1\right) ^{\alpha \left( x_{1};0,0,0,0\right) }B(g\otimes
gx_{1};X_{1},gx_{2})X_{1}\otimes gx_{2}\otimes gx_{1}
\end{equation*}%
\begin{equation*}
\alpha \left( x_{1};0,0,0,0\right) =a+b_{1}+b_{2}=1
\end{equation*}%
so that we get%
\begin{equation*}
-B(g\otimes gx_{1};X_{1},gx_{2})X_{1}\otimes gx_{2}\otimes gx_{1}
\end{equation*}%
Second summand%
\begin{eqnarray*}
&&+\sum_{a,b_{1},b_{2},d,e_{1},e_{2}=0}^{1}\sum_{l_{1}=0}^{b_{1}}%
\sum_{l_{2}=0}^{b_{2}}\sum_{u_{1}=0}^{e_{1}}\sum_{u_{2}=0}^{e_{2}}\left(
-1\right) ^{\alpha \left( 1_{H};l_{1},l_{2},u_{1},u_{2}\right) } \\
&&B(x_{1}\otimes
gx_{1};G^{a}X_{1}^{b_{1}}X_{2}^{b_{2}},g^{d}x_{1}^{e_{1}}x_{2}^{e_{2}})G^{a}X_{1}^{b_{1}-l_{1}}X_{2}^{b_{2}-l_{2}}\otimes g^{d}x_{1}^{e_{1}-u_{1}}x_{2}^{e_{2}-u_{2}}\otimes
\\
&&g^{a+b_{1}+b_{2}+l_{1}+l_{2}+d+e_{1}+e_{2}+u_{1}+u_{2}}x_{1}^{l_{1}+u_{1}}x_{2}^{l_{2}+u_{2}}
\end{eqnarray*}%
\begin{eqnarray*}
l_{1}+u_{1} &=&1 \\
l_{2} &=&u_{2}=0 \\
a+b_{1}+b_{2}+d+e_{1}+e_{2} &\equiv &0 \\
d &=&1 \\
e_{1} &=&u_{1}=1 \\
e_{2} &=&1 \\
a &=&0 \\
b_{1} &=&1\Rightarrow l_{1}=0,u_{1}=1 \\
b_{2} &=&l_{2}=0
\end{eqnarray*}%
\begin{equation*}
\left[ \left( -1\right) ^{\alpha \left( 1_{H};0,0,1,0\right) }B(x_{1}\otimes
gx_{1};X_{1},gx_{1}x_{2})\right] X_{1}\otimes gx_{2}\otimes gx_{1}
\end{equation*}%
\begin{equation*}
\alpha \left( 1_{H};0,0,1,0\right) \equiv e_{2}+a+b_{1}+b_{2}=1+1+0\equiv 0
\end{equation*}%
and we get%
\begin{equation*}
-B(g\otimes gx_{1};X_{1},gx_{2})X_{1}\otimes gx_{2}\otimes gx_{1}
\end{equation*}%
Finally, considering also the right side we get%
\begin{equation*}
-B(g\otimes gx_{1};X_{1},gx_{2})+B(x_{1}\otimes
gx_{1};X_{1},gx_{1}x_{2})-B(x_{1}\otimes g;X_{1},gx_{2})=0
\end{equation*}%
which holds in view of the form of the elements.

\subsubsection{Case $X_{1}\otimes g\otimes x_{1}x_{2}$}

First summand%
\begin{eqnarray*}
&&\sum_{a,b_{1},b_{2},d,e_{1},e_{2}=0}^{1}\sum_{l_{1}=0}^{b_{1}}%
\sum_{l_{2}=0}^{b_{2}}\sum_{u_{1}=0}^{e_{1}}\sum_{u_{2}=0}^{e_{2}}\left(
-1\right) ^{\alpha \left( x_{1};l_{1},l_{2},u_{1},u_{2}\right) } \\
&&B(g\otimes
gx_{1};G^{a}X_{1}^{b_{1}}X_{2}^{b_{2}},g^{d}x_{1}^{e_{1}}x_{2}^{e_{2}})G^{a}X_{1}^{b_{1}-l_{1}}X_{2}^{b_{2}-l_{2}}\otimes g^{d}x_{1}^{e_{1}-u_{1}}x_{2}^{e_{2}-u_{2}}\otimes
\\
&&g^{a+b_{1}+b_{2}+l_{1}+l_{2}+d+e_{1}+e_{2}+u_{1}+u_{2}}x_{1}^{l_{1}+u_{1}+1}x_{2}^{l_{2}+u_{2}}
\end{eqnarray*}%
\begin{eqnarray*}
l_{1} &=&u_{1}=0 \\
l_{2}+u_{2} &=&1 \\
a+b_{1}+b_{2}+d+e_{1}+e_{2} &\equiv &1 \\
d &=&1 \\
e_{1} &=&0 \\
e_{2} &=&u_{2} \\
a &=&0 \\
b_{1} &=&1 \\
b_{2} &=&l_{2}
\end{eqnarray*}%
\begin{eqnarray*}
&&\sum_{\substack{ b_{2},e_{2}=0  \\ b_{2}+e_{2}=1}}^{1}\left( -1\right)
^{\alpha \left( x_{1};0,b_{2},0,e_{2}\right) }B(g\otimes
gx_{1};X_{1}X_{2}^{b_{2}},gx_{2}^{e_{2}})X_{1}\otimes g\otimes x_{1}x_{2} \\
&=&\left[ \left( -1\right) ^{\alpha \left( x_{1};0,0,0,1\right) }B(g\otimes
gx_{1};X_{1},gx_{2})+\left( -1\right) ^{\alpha \left( x_{1};0,1,0,0\right)
}B(g\otimes gx_{1};X_{1}X_{2},g)\right] X_{1}\otimes g\otimes x_{1}x_{2}
\end{eqnarray*}%
\begin{eqnarray*}
&&\alpha \left( x_{1};0,0,0,1\right) \equiv 1 \\
&&\alpha \left( x_{1};0,1,0,0\right) \equiv 0+1+1+1\equiv 1
\end{eqnarray*}%
\begin{equation*}
-B(g\otimes gx_{1};X_{1},gx_{2})-B(g\otimes gx_{1};X_{1}X_{2},g)
\end{equation*}

Second summand%
\begin{eqnarray*}
l_{1}+u_{1} &=&1 \\
l_{2}+u_{2} &=&1 \\
a+b_{1}+b_{2}+d+e_{1}+e_{2} &\equiv &0 \\
d &=&1 \\
e_{1} &=&u_{1}=1 \\
e_{2} &=&u_{2} \\
a &=&0 \\
b_{1} &=&1\Rightarrow l_{1}=0,u_{1}=1 \\
b_{2} &=&l_{2}
\end{eqnarray*}

\begin{eqnarray*}
&&+\sum_{\substack{ b_{2},e_{2}=0  \\ b_{2}+e_{2}=1}}^{1}\left( -1\right)
^{\alpha \left( 1_{H};0,b_{2},1,e_{2}\right) }B(x_{1}\otimes
gx_{1};X_{1}X_{2}^{b_{2}},gx_{1}x_{2}^{e_{2}})X_{1}\otimes g\otimes
x_{1}x_{2} \\
&=&\left[
\begin{array}{c}
\left( -1\right) ^{\alpha \left( 1_{H};0,0,1,1\right) }B(x_{1}\otimes
gx_{1};X_{1},gx_{1}x_{2}) \\
+\left( -1\right) ^{\alpha \left( 1_{H};0,1,1,0\right) }B(x_{1}\otimes
gx_{1};X_{1}X_{2},gx_{1})%
\end{array}%
\right] X_{1}\otimes g\otimes x_{1}x_{2}
\end{eqnarray*}%
\begin{eqnarray*}
&&\alpha \left( 1_{H};0,0,1,1\right) \equiv 1+e_{2}=1+1\equiv 0 \\
&&\alpha \left( 1_{H};0,1,1,0\right) \equiv 0+0+1+1+1\equiv 1
\end{eqnarray*}%
and we get%
\begin{gather*}
-B(g\otimes gx_{1};X_{1},gx_{2})-B(g\otimes gx_{1};X_{1}X_{2},g) \\
+B(x_{1}\otimes gx_{1};X_{1},gx_{1}x_{2})-B(x_{1}\otimes
gx_{1};X_{1}X_{2},gx_{1})=0
\end{gather*}%
which hold in view of the form of the element.

\subsubsection{Case $1_{A}\otimes gx_{1}x_{2}\otimes gx_{1}$}

\begin{eqnarray*}
&&\sum_{a,b_{1},b_{2},d,e_{1},e_{2}=0}^{1}\sum_{l_{1}=0}^{b_{1}}%
\sum_{l_{2}=0}^{b_{2}}\sum_{u_{1}=0}^{e_{1}}\sum_{u_{2}=0}^{e_{2}}\left(
-1\right) ^{\alpha \left( x_{1};l_{1},l_{2},u_{1},u_{2}\right) } \\
&&B(g\otimes
gx_{1};G^{a}X_{1}^{b_{1}}X_{2}^{b_{2}},g^{d}x_{1}^{e_{1}}x_{2}^{e_{2}})G^{a}X_{1}^{b_{1}-l_{1}}X_{2}^{b_{2}-l_{2}}\otimes g^{d}x_{1}^{e_{1}-u_{1}}x_{2}^{e_{2}-u_{2}}\otimes
\\
&&g^{a+b_{1}+b_{2}+l_{1}+l_{2}+d+e_{1}+e_{2}+u_{1}+u_{2}}x_{1}^{l_{1}+u_{1}+1}x_{2}^{l_{2}+u_{2}}
\end{eqnarray*}

First summand%
\begin{eqnarray*}
l_{1} &=&u_{1}=0 \\
l_{2} &=&u_{2}=0 \\
a+b_{1}+b_{2}+d+e_{1}+e_{2} &\equiv &1 \\
d &=&1 \\
e_{1} &=&1 \\
e_{2} &=&1 \\
a &=&0 \\
b_{1} &=&0 \\
b_{2} &=&0
\end{eqnarray*}

\begin{equation*}
\left( -1\right) ^{\alpha \left( x_{1};0,0,0,0\right) }B(g\otimes
gx_{1};1_{H},gx_{1}x_{2})1_{A}\otimes gx_{1}x_{2}\otimes gx_{1}
\end{equation*}

\begin{equation*}
\alpha \left( x_{1};0,0,0,0\right) \equiv a+b_{1}+b_{2}=0
\end{equation*}%
Thus we obtain%
\begin{equation*}
B(g\otimes gx_{1};1_{H},gx_{1}x_{2})
\end{equation*}

Second summand%
\begin{eqnarray*}
&&+\sum_{a,b_{1},b_{2},d,e_{1},e_{2}=0}^{1}\sum_{l_{1}=0}^{b_{1}}%
\sum_{l_{2}=0}^{b_{2}}\sum_{u_{1}=0}^{e_{1}}\sum_{u_{2}=0}^{e_{2}}\left(
-1\right) ^{\alpha \left( 1_{H};l_{1},l_{2},u_{1},u_{2}\right) } \\
&&B(x_{1}\otimes
gx_{1};G^{a}X_{1}^{b_{1}}X_{2}^{b_{2}},g^{d}x_{1}^{e_{1}}x_{2}^{e_{2}})G^{a}X_{1}^{b_{1}-l_{1}}X_{2}^{b_{2}-l_{2}}\otimes g^{d}x_{1}^{e_{1}-u_{1}}x_{2}^{e_{2}-u_{2}}\otimes
\\
&&g^{a+b_{1}+b_{2}+l_{1}+l_{2}+d+e_{1}+e_{2}+u_{1}+u_{2}}x_{1}^{l_{1}+u_{1}}x_{2}^{l_{2}+u_{2}}
\end{eqnarray*}

\begin{eqnarray*}
l_{1}+u_{1} &=&1 \\
l_{2} &=&u_{2}=0 \\
a+b_{1}+b_{2}+d+e_{1}+e_{2} &\equiv &0 \\
d &=&1 \\
e_{1}-u_{1} &=&1\Rightarrow e_{1}=1,u_{1}=0,l_{1}=1 \\
e_{2} &=&1 \\
a &=&0 \\
b_{1} &=&l_{1}=1 \\
b_{2} &=&l_{2}=0
\end{eqnarray*}

\begin{equation*}
\left( -1\right) ^{\alpha \left( 1_{H};1,0,0,0\right) }B(x_{1}\otimes
gx_{1};X_{1},gx_{1}x_{2})1_{A}\otimes gx_{1}x_{2}\otimes gx_{1}
\end{equation*}%
\begin{equation*}
\alpha \left( 1_{H};1,0,0,0\right) \equiv b_{2}=0
\end{equation*}%
Considering also the right side, we finally get

\begin{equation*}
B(g\otimes gx_{1};1_{H},gx_{1}x_{2})+B(x_{1}\otimes
gx_{1};X_{1},gx_{1}x_{2})-B(x_{1}\otimes g;1_{H},gx_{1}x_{2})=0
\end{equation*}%
which hold in view of the form of the element.

\subsubsection{Case $1_{A}\otimes gx_{1}\otimes x_{1}x_{2}$}

First summand%
\begin{eqnarray*}
l_{1} &=&u_{1}=0 \\
l_{2}+u_{2} &=&1 \\
a+b_{1}+b_{2}+d+e_{1}+e_{2} &\equiv &1 \\
d &=&1 \\
e_{1} &=&1 \\
e_{2} &=&u_{1} \\
a &=&0 \\
b_{1} &=&0 \\
b_{2} &=&l_{2}
\end{eqnarray*}

\begin{eqnarray*}
&&\sum_{\substack{ b_{2},e_{2}=0  \\ b_{2}+e_{2}=1}}^{1}\left( -1\right)
^{\alpha \left( x_{1};0,b_{2},0,e_{2}\right) }B(g\otimes
gx_{1};X_{2}^{b_{2}},gx_{1}x_{2}^{e_{2}})1_{A}\otimes gx_{1}\otimes
x_{1}x_{2} \\
&=&\left[
\begin{array}{c}
\left( -1\right) ^{\alpha \left( x_{1};0,0,0,1\right) }B(g\otimes
gx_{1};1_{A},gx_{1}x_{2}) \\
+\left( -1\right) ^{\alpha \left( x_{1};0,1,0,0\right) }B(g\otimes
gx_{1};X_{2},gx_{1})%
\end{array}%
\right] 1_{A}\otimes gx_{1}\otimes x_{1}x_{2}
\end{eqnarray*}

\begin{eqnarray*}
&&\alpha \left( x_{1};0,0,0,1\right) \equiv 1 \\
&&\alpha \left( x_{1};0,1,0,0\right) \equiv a+b_{1}+b_{2}+1=0+0+1+1\equiv 0
\end{eqnarray*}%
\begin{equation*}
-B(g\otimes gx_{1};1_{A},gx_{1}x_{2})+B(g\otimes gx_{1};X_{2},gx_{1})
\end{equation*}

Second summand%
\begin{eqnarray*}
&&+\sum_{a,b_{1},b_{2},d,e_{1},e_{2}=0}^{1}\sum_{l_{1}=0}^{b_{1}}%
\sum_{l_{2}=0}^{b_{2}}\sum_{u_{1}=0}^{e_{1}}\sum_{u_{2}=0}^{e_{2}}\left(
-1\right) ^{\alpha \left( 1_{H};l_{1},l_{2},u_{1},u_{2}\right) } \\
&&B(x_{1}\otimes
gx_{1};G^{a}X_{1}^{b_{1}}X_{2}^{b_{2}},g^{d}x_{1}^{e_{1}}x_{2}^{e_{2}})G^{a}X_{1}^{b_{1}-l_{1}}X_{2}^{b_{2}-l_{2}}\otimes g^{d}x_{1}^{e_{1}-u_{1}}x_{2}^{e_{2}-u_{2}}\otimes
\\
&&g^{a+b_{1}+b_{2}+l_{1}+l_{2}+d+e_{1}+e_{2}+u_{1}+u_{2}}x_{1}^{l_{1}+u_{1}}x_{2}^{l_{2}+u_{2}}
\end{eqnarray*}

\begin{eqnarray*}
l_{1}+u_{1} &=&1 \\
l_{2}+u_{2} &=&1 \\
a+b_{1}+b_{2}+d+e_{1}+e_{2} &\equiv &0 \\
d &=&1 \\
e_{1}-u_{1} &=&1\Rightarrow e_{1}=1,u_{1}=0,l_{1}=1 \\
e_{2} &=&u_{2} \\
a &=&0 \\
b_{1} &=&l_{1}=1 \\
b_{2} &=&l_{2}
\end{eqnarray*}

\begin{eqnarray*}
&&\sum_{\substack{ b_{2},e_{2}=0  \\ b_{2}+e_{2}=1}}^{1}\left( -1\right)
^{\alpha \left( 1_{H};1,b_{2},0,e_{2}\right) }B(x_{1}\otimes
gx_{1};X_{1}X_{2}^{b_{2}},gx_{1}x_{2}^{e_{2}})1_{A}\otimes gx_{1}\otimes
x_{1}x_{2} \\
&=&\left[
\begin{array}{c}
\left( -1\right) ^{\alpha \left( 1_{H};1,0,0,1\right) }B(x_{1}\otimes
gx_{1};X_{1},gx_{1}x_{2}) \\
+\left( -1\right) ^{\alpha \left( 1_{H};1,1,0,0\right) }B(x_{1}\otimes
gx_{1};X_{1}X_{2},gx_{1})%
\end{array}%
\right] 1_{A}\otimes gx_{1}\otimes x_{1}x_{2}
\end{eqnarray*}

\begin{eqnarray*}
&&\alpha \left( 1_{H};1,0,0,1\right) \equiv a+b_{1}\equiv 1 \\
&&\alpha \left( 1_{H};1,1,0,0\right) \equiv 1+b_{2}\equiv 0
\end{eqnarray*}%
\begin{equation*}
-B(x_{1}\otimes gx_{1};X_{1},gx_{1}x_{2})+B(x_{1}\otimes
gx_{1};X_{1}X_{2},gx_{1})
\end{equation*}%
Thus finally we get

\begin{equation*}
-B(g\otimes gx_{1};1_{A},gx_{1}x_{2})+B(g\otimes
gx_{1};X_{2},gx_{1})-B(x_{1}\otimes gx_{1};X_{1},gx_{1}x_{2})+B(x_{1}\otimes
gx_{1};X_{1}X_{2},gx_{1})=0
\end{equation*}%
which holds in view of the form of the elements.

\subsection{$B\left( x_{1}\otimes gx_{1};X_{1}X_{2},gx_{1}\right) $}

\begin{equation*}
a=0,b_{1}=b_{2}=1,d=1,e_{1}=1,e_{2}=0
\end{equation*}%
\begin{eqnarray*}
&&\sum_{l_{1}=0}^{1}\sum_{l_{2}=0}^{1}\sum_{u_{1}=0}^{1}\left( -1\right)
^{\alpha \left( 1_{H};l_{1},l_{2},u_{1},0\right) }B(x_{1}\otimes
gx_{1};X_{1}X_{2},gx_{1}) \\
&&X_{1}^{1-l_{1}}X_{2}^{1-l_{2}}\otimes gx_{1}^{1-u_{1}}\otimes
g^{l_{1}+l_{2}+u_{1}}x_{1}^{l_{1}+u_{1}}x_{2}^{l_{2}} \\
&=&\left( -1\right) ^{\alpha \left( 1_{H};0,0,00\right) }B(x_{1}\otimes
gx_{1};X_{1}X_{2},gx_{1})X_{1}X_{2}\otimes gx_{1}\otimes 1_{H}+ \\
&&+\left( -1\right) ^{\alpha \left( 1_{H};0,0,1,0\right) }B(x_{1}\otimes
gx_{1};X_{1}X_{2},gx_{1})X_{1}X_{2}\otimes g\otimes gx_{1}+ \\
&&+\left( -1\right) ^{\alpha \left( 1_{H};0,1,0,0\right) }B(x_{1}\otimes
gx_{1};X_{1}X_{2},gx_{1})X_{1}\otimes gx_{1}\otimes gx_{2}+ \\
&&+\left( -1\right) ^{\alpha \left( 1_{H};1,0,0,0\right) }B(x_{1}\otimes
gx_{1};X_{1}X_{2},gx_{1})X_{2}\otimes gx_{1}\otimes gx_{1}+ \\
&&+\left( -1\right) ^{\alpha \left( 1_{H};1,1,0,0\right) }B(x_{1}\otimes
gx_{1};X_{1}X_{2},gx_{1})1_{A}\otimes gx_{1}\otimes x_{1}x_{2}+ \\
&&+\left( -1\right) ^{\alpha \left( 1_{H};0,1,1,0\right) }B(x_{1}\otimes
gx_{1};X_{1}X_{2},gx_{1})X_{1}\otimes g\otimes x_{1}x_{2}+ \\
&&+\left( -1\right) ^{\alpha \left( 1_{H};1,0,1,0\right) }B(x_{1}\otimes
gx_{1};X_{1}X_{2},gx_{1})X_{1}^{1-l_{1}}X_{2}^{1-l_{2}}\otimes
gx_{1}^{1-u_{1}}\otimes g^{l_{1}+l_{2}+u_{1}}x_{1}^{1+1}x_{2}^{l_{2}}+=0 \\
&&+\left( -1\right) ^{\alpha \left( 1_{H};1,1,1,0\right) }B(x_{1}\otimes
gx_{1};X_{1}X_{2},gx_{1})X_{1}^{1-l_{1}}X_{2}^{1-l_{2}}\otimes
gx_{1}^{1-u_{1}}\otimes g^{l_{1}+l_{2}+u_{1}}x_{1}^{1+1}x_{2}^{l_{2}}+=0
\end{eqnarray*}

\subsubsection{Case $X_{1}X_{2}\otimes gx_{1}\otimes 1_{H}$}

First summand

\begin{eqnarray*}
&&\sum_{a,b_{1},b_{2},d,e_{1},e_{2}=0}^{1}\sum_{l_{1}=0}^{b_{1}}%
\sum_{l_{2}=0}^{b_{2}}\sum_{u_{1}=0}^{e_{1}}\sum_{u_{2}=0}^{e_{2}}\left(
-1\right) ^{\alpha \left( x_{1};l_{1},l_{2},u_{1},u_{2}\right) } \\
&&B(g\otimes
gx_{1};G^{a}X_{1}^{b_{1}}X_{2}^{b_{2}},g^{d}x_{1}^{e_{1}}x_{2}^{e_{2}})G^{a}X_{1}^{b_{1}-l_{1}}X_{2}^{b_{2}-l_{2}}\otimes g^{d}x_{1}^{e_{1}-u_{1}}x_{2}^{e_{2}-u_{2}}\otimes
\\
&&g^{a+b_{1}+b_{2}+l_{1}+l_{2}+d+e_{1}+e_{2}+u_{1}+u_{2}}x_{1}^{l_{1}+u_{1}+1}x_{2}^{l_{2}+u_{2}}
\end{eqnarray*}%
This cannot occur.

Second summand%
\begin{eqnarray*}
&&+\sum_{a,b_{1},b_{2},d,e_{1},e_{2}=0}^{1}\sum_{l_{1}=0}^{b_{1}}%
\sum_{l_{2}=0}^{b_{2}}\sum_{u_{1}=0}^{e_{1}}\sum_{u_{2}=0}^{e_{2}}\left(
-1\right) ^{\alpha \left( 1_{H};l_{1},l_{2},u_{1},u_{2}\right) } \\
&&B(x_{1}\otimes
gx_{1};G^{a}X_{1}^{b_{1}}X_{2}^{b_{2}},g^{d}x_{1}^{e_{1}}x_{2}^{e_{2}})G^{a}X_{1}^{b_{1}-l_{1}}X_{2}^{b_{2}-l_{2}}\otimes g^{d}x_{1}^{e_{1}-u_{1}}x_{2}^{e_{2}-u_{2}}\otimes
\\
&&g^{a+b_{1}+b_{2}+l_{1}+l_{2}+d+e_{1}+e_{2}+u_{1}+u_{2}}x_{1}^{l_{1}+u_{1}}x_{2}^{l_{2}+u_{2}}
\end{eqnarray*}

\begin{eqnarray*}
l_{1} &=&u_{1}=0 \\
l_{2} &=&u_{2}=0 \\
a+b_{1}+b_{2}+d+e_{1}+e_{2} &\equiv &0 \\
d &=&1 \\
e_{1}-u_{1} &=&1\Rightarrow e_{1}=1 \\
e_{2} &=&u_{2}=0 \\
a &=&0 \\
b_{1} &=&1 \\
b_{2} &=&1
\end{eqnarray*}

\begin{equation*}
\left( -1\right) ^{\alpha \left( 1_{H};0,0,0,0\right) }B(x_{1}\otimes
gx_{1};X_{1}X_{2},gx_{1})X_{1}X_{2}\otimes gx_{1}\otimes 1_{H}
\end{equation*}%
\begin{equation*}
\alpha \left( 1_{H};0,0,0,0\right) \equiv 0
\end{equation*}%
Considering also the right side, we finally get%
\begin{equation*}
B(x_{1}\otimes gx_{1};X_{1}X_{2},gx_{1})=B(x_{1}\otimes
gx_{1};X_{1}X_{2},gx_{1})
\end{equation*}

which is obvious.

\subsubsection{Case $X_{1}X_{2}\otimes g\otimes gx_{1}$}

First summand%
\begin{eqnarray*}
&&\sum_{a,b_{1},b_{2},d,e_{1},e_{2}=0}^{1}\sum_{l_{1}=0}^{b_{1}}%
\sum_{l_{2}=0}^{b_{2}}\sum_{u_{1}=0}^{e_{1}}\sum_{u_{2}=0}^{e_{2}}\left(
-1\right) ^{\alpha \left( x_{1};l_{1},l_{2},u_{1},u_{2}\right) } \\
&&B(g\otimes
gx_{1};G^{a}X_{1}^{b_{1}}X_{2}^{b_{2}},g^{d}x_{1}^{e_{1}}x_{2}^{e_{2}})G^{a}X_{1}^{b_{1}-l_{1}}X_{2}^{b_{2}-l_{2}}\otimes g^{d}x_{1}^{e_{1}-u_{1}}x_{2}^{e_{2}-u_{2}}\otimes
\\
&&g^{a+b_{1}+b_{2}+l_{1}+l_{2}+d+e_{1}+e_{2}+u_{1}+u_{2}}x_{1}^{l_{1}+u_{1}+1}x_{2}^{l_{2}+u_{2}}
\end{eqnarray*}

\begin{eqnarray*}
l_{1} &=&u_{1}=0 \\
l_{2} &=&u_{2}=0 \\
a+b_{1}+b_{2}+d+e_{1}+e_{2} &\equiv &1 \\
d &=&1 \\
e_{1} &=&0 \\
e_{2} &=&0 \\
a &=&0 \\
b_{1} &=&1 \\
b_{2} &=&1
\end{eqnarray*}%
\begin{equation*}
\left( -1\right) ^{\alpha \left( x_{1};0,0,0,0\right) }B(g\otimes
gx_{1};X_{1}X_{2},g)X_{1}X_{2}\otimes g\otimes gx_{1}
\end{equation*}

\begin{equation*}
\alpha \left( x_{1};0,0,0,0\right) \equiv a+b_{1}+b_{2}=0+1+1\equiv 0
\end{equation*}%
and we get
\begin{equation*}
B(g\otimes gx_{1};X_{1}X_{2},g)X_{1}X_{2}\otimes g\otimes gx_{1}
\end{equation*}

Second summand%
\begin{eqnarray*}
&&+\sum_{a,b_{1},b_{2},d,e_{1},e_{2}=0}^{1}\sum_{l_{1}=0}^{b_{1}}%
\sum_{l_{2}=0}^{b_{2}}\sum_{u_{1}=0}^{e_{1}}\sum_{u_{2}=0}^{e_{2}}\left(
-1\right) ^{\alpha \left( 1_{H};l_{1},l_{2},u_{1},u_{2}\right) } \\
&&B(x_{1}\otimes
gx_{1};G^{a}X_{1}^{b_{1}}X_{2}^{b_{2}},g^{d}x_{1}^{e_{1}}x_{2}^{e_{2}})G^{a}X_{1}^{b_{1}-l_{1}}X_{2}^{b_{2}-l_{2}}\otimes g^{d}x_{1}^{e_{1}-u_{1}}x_{2}^{e_{2}-u_{2}}\otimes
\\
&&g^{a+b_{1}+b_{2}+l_{1}+l_{2}+d+e_{1}+e_{2}+u_{1}+u_{2}}x_{1}^{l_{1}+u_{1}}x_{2}^{l_{2}+u_{2}}
\end{eqnarray*}%
\begin{eqnarray*}
l_{1}+u_{1} &=&1 \\
l_{2} &=&u_{2}=0 \\
a+b_{1}+b_{2}+d+e_{1}+e_{2} &\equiv &0 \\
d &=&1 \\
e_{1} &=&u_{1}=1 \\
e_{2} &=&0 \\
a &=&0 \\
b_{1}-l_{1} &=&1\Rightarrow b_{1}=1,l_{1}=0,u_{1}=1 \\
b_{2} &=&1
\end{eqnarray*}%
\begin{equation*}
\left( -1\right) ^{\alpha \left( 1_{H};0,0,1,0\right) }B(x_{1}\otimes
gx_{1};X_{1}X_{2},gx_{1})X_{1}X_{2}\otimes g\otimes gx_{1}
\end{equation*}

\begin{equation*}
\alpha \left( 1_{H};0,0,1,0\right) \equiv e_{2}+a+b_{1}+b_{2}=0+0+1+1\equiv 0
\end{equation*}%
Considering also the right side of the equality, we obtain%
\begin{equation*}
B(g\otimes gx_{1};X_{1}X_{2},g)+B(x_{1}\otimes
gx_{1};X_{1}X_{2},gx_{1})-B(x_{1}\otimes g;X_{1}X_{2},g)=0
\end{equation*}%
which holds in view of the form of the elements.

\subsubsection{Case $X_{1}\otimes gx_{1}\otimes gx_{2}$}

Already considered in $B(x_{1}\otimes gx_{1};X_{1},gx_{1}x_{2}).$

\subsubsection{Case $X_{2}\otimes gx_{1}\otimes gx_{1}$}

First summand%
\begin{eqnarray*}
&&\sum_{a,b_{1},b_{2},d,e_{1},e_{2}=0}^{1}\sum_{l_{1}=0}^{b_{1}}%
\sum_{l_{2}=0}^{b_{2}}\sum_{u_{1}=0}^{e_{1}}\sum_{u_{2}=0}^{e_{2}}\left(
-1\right) ^{\alpha \left( x_{1};l_{1},l_{2},u_{1},u_{2}\right) } \\
&&B(g\otimes
gx_{1};G^{a}X_{1}^{b_{1}}X_{2}^{b_{2}},g^{d}x_{1}^{e_{1}}x_{2}^{e_{2}})G^{a}X_{1}^{b_{1}-l_{1}}X_{2}^{b_{2}-l_{2}}\otimes g^{d}x_{1}^{e_{1}-u_{1}}x_{2}^{e_{2}-u_{2}}\otimes
\\
&&g^{a+b_{1}+b_{2}+l_{1}+l_{2}+d+e_{1}+e_{2}+u_{1}+u_{2}}x_{1}^{l_{1}+u_{1}+1}x_{2}^{l_{2}+u_{2}}
\end{eqnarray*}

\begin{eqnarray*}
l_{1} &=&u_{1}=0 \\
l_{2} &=&u_{2}=0 \\
a+b_{1}+b_{2}+d+e_{1}+e_{2} &\equiv &1 \\
d &=&1 \\
e_{1} &=&1 \\
e_{2} &=&0 \\
a &=&0 \\
b_{1} &=&0 \\
b_{2} &=&1
\end{eqnarray*}%
\begin{equation*}
\left( -1\right) ^{\alpha \left( x_{1};0,0,0,0\right) }B(g\otimes
gx_{1};X_{2},gx_{1})X_{2}\otimes gx_{1}\otimes gx_{1}
\end{equation*}

and we get%
\begin{equation*}
\alpha \left( x_{1};0,0,0,0\right) \equiv a+b_{1}+b_{2}=0+0+1\equiv 1
\end{equation*}%
\begin{equation*}
-B(g\otimes gx_{1};X_{2},gx_{1})
\end{equation*}

Second summand%
\begin{eqnarray*}
l_{1}+u_{1} &=&1 \\
l_{2} &=&u_{2}=0 \\
a+b_{1}+b_{2}+d+e_{1}+e_{2} &\equiv &0 \\
d &=&1 \\
e_{1}-u_{1} &=&1\Rightarrow e_{1}=1,u_{1}=0,l_{1}=1 \\
e_{2} &=&0 \\
a &=&0 \\
b_{1} &=&l_{1}=1 \\
b_{2} &=&1
\end{eqnarray*}%
\begin{equation*}
\left( -1\right) ^{\alpha \left( 1_{H};1,0,0,0\right) }B(x_{1}\otimes
gx_{1};X_{1}X_{2},gx_{1})X_{2}\otimes gx_{1}\otimes gx_{1}
\end{equation*}%
\begin{equation*}
\alpha \left( 1_{H};1,0,0,0\right) \equiv b_{2}=1
\end{equation*}%
Considering also the right side, we finally get%
\begin{equation*}
-B(g\otimes gx_{1};X_{2},gx_{1})-B(x_{1}\otimes
gx_{1};X_{1}X_{2},gx_{1})-B(x_{1}\otimes g;X_{2},gx_{1})=0
\end{equation*}%
which holds in view of the form of the elements.

\subsubsection{Case $1_{A}\otimes gx_{1}\otimes x_{1}x_{2}$}

Already considered in $B(x_{1}\otimes gx_{1};X_{1},gx_{1}x_{2}).$

\subsubsection{Case $X_{1}\otimes g\otimes x_{1}x_{2}$}

Already considered in $B(x_{1}\otimes gx_{1};X_{1},gx_{1}x_{2}).$

\subsection{$B(x_{1}\otimes gx_{1};GX_{1},gx_{1})$}

\begin{equation*}
a=1,b_{1}=1,b_{2}=0,d=1,e_{1}=1,e_{2}=0.
\end{equation*}%
and we obtain

thus we get

\begin{gather*}
\sum_{l_{1}=0}^{1}\sum_{u_{1}=0}^{1}\left( -1\right) ^{\alpha \left(
1_{H};l_{1},0,u_{1},0\right) }B(x_{1}\otimes
gx_{1};GX_{1},gx_{1})GX_{1}^{1-l_{1}}\otimes gx_{1}^{1-u_{1}}\otimes
g^{l_{1}+u_{1}}x_{1}^{l_{1}+u_{1}} \\
=\left( -1\right) ^{\alpha \left( 1_{H};0,0,0,0\right) }B(x_{1}\otimes
gx_{1};GX_{1},gx_{1})GX_{1}\otimes gx_{1}\otimes 1_{H}+ \\
+\left( -1\right) ^{\alpha \left( 1_{H};0,0,1,0\right) }B(x_{1}\otimes
gx_{1};GX_{1},gx_{1})GX_{1}\otimes g\otimes gx_{1}+ \\
+\left( -1\right) ^{\alpha \left( 1_{H};1,0,0,0\right) }B(x_{1}\otimes
gx_{1};GX_{1},gx_{1})G\otimes gx_{1}\otimes gx_{1}+ \\
+\left( -1\right) ^{\alpha \left( 1_{H};1,0,1,0\right) }B(x_{1}\otimes
gx_{1};GX_{1},gx_{1})G\otimes g\otimes g^{l_{1}+u_{1}}x_{1}^{1+1}+=0
\end{gather*}

\subsubsection{Case $GX_{1}\otimes gx_{1}\otimes 1_{H}$}

First summand is zero. Second summand
\begin{eqnarray*}
l_{1} &=&u_{1}=0 \\
l_{2} &=&u_{2}=0 \\
a+b_{1}+b_{2}+d+e_{1}+e_{2} &\equiv &0 \\
d &=&1 \\
e_{1} &=&1 \\
e_{2} &=&0 \\
a &=&1 \\
b_{1} &=&1 \\
b_{2} &=&0
\end{eqnarray*}%
and we get%
\begin{equation*}
B(x_{1}\otimes gx_{1};GX_{1},gx_{1})GX_{1}\otimes gx_{1}\otimes 1_{H}
\end{equation*}%
From the right part of the equality we get the same, so there is nothing new.

\subsubsection{Case $GX_{1}\otimes g\otimes gx_{1}$}

First summand%
\begin{eqnarray*}
&&\sum_{a,b_{1},b_{2},d,e_{1},e_{2}=0}^{1}\sum_{l_{1}=0}^{b_{1}}%
\sum_{l_{2}=0}^{b_{2}}\sum_{u_{1}=0}^{e_{1}}\sum_{u_{2}=0}^{e_{2}}\left(
-1\right) ^{\alpha \left( x_{1};l_{1},l_{2},u_{1},u_{2}\right) } \\
&&B(g\otimes
gx_{1};G^{a}X_{1}^{b_{1}}X_{2}^{b_{2}},g^{d}x_{1}^{e_{1}}x_{2}^{e_{2}})G^{a}X_{1}^{b_{1}-l_{1}}X_{2}^{b_{2}-l_{2}}\otimes g^{d}x_{1}^{e_{1}-u_{1}}x_{2}^{e_{2}-u_{2}}\otimes
\\
&&g^{a+b_{1}+b_{2}+l_{1}+l_{2}+d+e_{1}+e_{2}+u_{1}+u_{2}}x_{1}^{l_{1}+u_{1}+1}x_{2}^{l_{2}+u_{2}}
\end{eqnarray*}%
\begin{eqnarray*}
l_{1} &=&u_{1}=0 \\
l_{2} &=&u_{2}=0 \\
a+b_{1}+b_{2}+d+e_{1}+e_{2} &\equiv &1 \\
d &=&1 \\
e_{1} &=&0 \\
e_{2} &=&0 \\
a &=&1 \\
b_{1} &=&1 \\
b_{2} &=&0
\end{eqnarray*}%
\begin{equation*}
\alpha \left( x_{1};0,0,0,0\right) \equiv 1+1+0\equiv 0
\end{equation*}%
\begin{equation*}
B(g\otimes gx_{1};GX_{1},g)GX_{1}\otimes g\otimes gx_{1}
\end{equation*}

Second summand

\begin{eqnarray*}
l_{1}+u_{1} &=&1 \\
l_{2} &=&u_{2}=0 \\
a+b_{1}+b_{2}+d+e_{1}+e_{2} &\equiv &0 \\
d &=&1 \\
e_{1} &=&u_{1}=1 \\
e_{2} &=&0 \\
a &=&1 \\
b_{1} &=&1\Rightarrow l_{1}=0,u_{1}=1 \\
b_{2} &=&0
\end{eqnarray*}%
and we get%
\begin{equation*}
\alpha \left( 1_{H};0,0,1,0\right) \equiv e_{2}+\left( a+b_{1}+b_{2}\right)
\equiv 0+1+1+0\equiv 0
\end{equation*}%
\begin{equation*}
B(x_{1}\otimes gx_{1};GX_{1},gx_{1})GX_{1}\otimes g\otimes gx_{1}
\end{equation*}%
Considering also the right side of the equality we get

\begin{equation*}
B(g\otimes gx_{1};GX_{1},g)+B(x_{1}\otimes
gx_{1};GX_{1},gx_{1})-B(x_{1}\otimes g;GX_{1},g)=0
\end{equation*}%
which holds in view of the form of the elements.

\subsubsection{Case $G\otimes gx_{1}\otimes gx_{1}$}

First summand%
\begin{eqnarray*}
l_{1} &=&u_{1}=0 \\
l_{2} &=&u_{2}=0 \\
a+b_{1}+b_{2}+d+e_{1}+e_{2} &\equiv &1 \\
d &=&1,e_{1}=1,e_{2}=0,a=1,b_{1}=0,b_{2}=0.
\end{eqnarray*}%
\begin{equation*}
\alpha \left( x_{1};0,0,0,0\right) \equiv a+b_{1}+b_{2}\equiv 1+0+0=1
\end{equation*}%
\begin{equation*}
-B(g\otimes gx_{1};G,gx_{1})G\otimes gx_{1}\otimes gx_{1}
\end{equation*}

Second summand

\begin{eqnarray*}
l_{1}+u_{1} &=&1 \\
l_{2} &=&u_{2}=0 \\
a+b_{1}+b_{2}+d+e_{1}+e_{2} &\equiv &0 \\
d &=&1 \\
e_{1} &=&1\Rightarrow u_{1}=0,l_{1}=1 \\
e_{2} &=&0,a=1,b_{1}=1,b_{2}=0,
\end{eqnarray*}%
\begin{equation*}
\alpha \left( 1_{H};1,0,0,0\right) \equiv b_{2}=0
\end{equation*}%
\begin{equation*}
B(x_{1}\otimes gx_{1};GX_{1},gx_{1})G\otimes gx_{1}\otimes gx_{1}
\end{equation*}%
Considering also the right side we get%
\begin{equation*}
-B(g\otimes gx_{1};G,gx_{1})+B(x_{1}\otimes
gx_{1};GX_{1},gx_{1})-B(x_{1}\otimes g;G,gx_{1})=0
\end{equation*}%
which holds in view of the form of the elements.

\subsection{$B(x_{1}\otimes gx_{1};GX_{2},gx_{1})$}

\begin{equation*}
a=1,b_{1}=0,b_{2}=1,d=1,e_{1}=1,e_{2}=0,
\end{equation*}%
and we get%
\begin{eqnarray*}
&&+\sum_{l_{2}=0}^{1}\sum_{u_{1}=0}^{1}\left( -1\right) ^{\alpha \left(
1_{H};0,l_{2},u_{1},0\right) }B(x_{1}\otimes
gx_{1};GX_{2},gx_{1})GX_{2}^{1-l_{2}}\otimes gx_{1}^{1-u_{1}}\otimes
g^{l_{2}+u_{1}}x_{1}^{u_{1}}x_{2}^{l_{2}} \\
= &&\left( -1\right) ^{\alpha \left( 1_{H};0,0,0,0\right) }B(x_{1}\otimes
gx_{1};GX_{2},gx_{1})GX_{2}\otimes gx_{1}\otimes 1_{H}+ \\
+ &&\left( -1\right) ^{\alpha \left( 1_{H};0,0,1,0\right) }B(x_{1}\otimes
gx_{1};GX_{2},gx_{1})GX_{2}\otimes g\otimes gx_{1}+ \\
+ &&\left( -1\right) ^{\alpha \left( 1_{H};0,1,0,0\right) }B(x_{1}\otimes
gx_{1};GX_{2},gx_{1})G\otimes gx_{1}\otimes gx_{2}+ \\
&&+\left( -1\right) ^{\alpha \left( 1_{H};0,1,1,0\right) }B(x_{1}\otimes
gx_{1};GX_{2},gx_{1})G\otimes g\otimes x_{1}x_{2}+
\end{eqnarray*}

\subsubsection{Case $GX_{2}\otimes gx_{1}\otimes 1_{H}$}

First summand nothing. Second summand

\begin{eqnarray*}
&&+\sum_{a,b_{1},b_{2},d,e_{1},e_{2}=0}^{1}\sum_{l_{1}=0}^{b_{1}}%
\sum_{l_{2}=0}^{b_{2}}\sum_{u_{1}=0}^{e_{1}}\sum_{u_{2}=0}^{e_{2}}\left(
-1\right) ^{\alpha \left( 1_{H};l_{1},l_{2},u_{1},u_{2}\right) } \\
&&B(x_{1}\otimes
gx_{1};G^{a}X_{1}^{b_{1}}X_{2}^{b_{2}},g^{d}x_{1}^{e_{1}}x_{2}^{e_{2}})G^{a}X_{1}^{b_{1}-l_{1}}X_{2}^{b_{2}-l_{2}}\otimes g^{d}x_{1}^{e_{1}-u_{1}}x_{2}^{e_{2}-u_{2}}\otimes
\\
&&g^{a+b_{1}+b_{2}+l_{1}+l_{2}+d+e_{1}+e_{2}+u_{1}+u_{2}}x_{1}^{l_{1}+u_{1}}x_{2}^{l_{2}+u_{2}}
\end{eqnarray*}%
\begin{eqnarray*}
l_{1} &=&u_{1}=0 \\
l_{2} &=&u_{2}=0 \\
a+b_{1}+b_{2}+d+e_{1}+e_{2} &\equiv &0 \\
d &=&1,e_{1}=1,e_{2}=0,a=1,b_{1}=0,b_{2}=1.
\end{eqnarray*}%
\begin{equation*}
\alpha \left( 1_{H};l_{1},l_{2},u_{1},u_{2}\right) \equiv 0
\end{equation*}

\begin{equation*}
B(x_{1}\otimes gx_{1};GX_{2}\otimes gx_{1})GX_{2}\otimes gx_{1}\otimes 1_{H}
\end{equation*}%
Considering also the right side, we get an obvious equality.

\subsubsection{$GX_{2}\otimes g\otimes gx_{1}$}

First summand%
\begin{equation*}
B(g\otimes gx_{1};GX_{2},g)GX_{2}\otimes g\otimes gx_{1}
\end{equation*}%
\begin{eqnarray*}
l_{1} &=&u_{1}=0 \\
l_{2} &=&u_{2}=0 \\
a+b_{1}+b_{2}+d+e_{1}+e_{2} &\equiv &1 \\
d &=&1,e_{1}=0,e_{2}=0,a=1,b_{1}=0,b_{2}=1.
\end{eqnarray*}%
\begin{equation*}
\alpha \left( x_{1};0,0,0,0\right) \equiv a+b_{1}+b_{2}\equiv 1+0+1\equiv 0
\end{equation*}%
and we get%
\begin{equation*}
B(g\otimes gx_{1};GX_{2},g)GX_{2}\otimes g\otimes gx_{1}
\end{equation*}

Second summand%
\begin{eqnarray*}
l_{1}+u_{1} &=&1 \\
l_{2} &=&u_{2}=0 \\
a+b_{1}+b_{2}+d+e_{1}+e_{2} &\equiv &0 \\
d &=&1,e_{1}=u_{1},e_{2}=0,a=1,b_{1}=l_{1},b_{2}=1.
\end{eqnarray*}%
\begin{eqnarray*}
&&+\sum_{\substack{ b_{1},e_{1}=0  \\ b_{1}+e_{1}\equiv 1}}^{1}\left(
-1\right) ^{\alpha \left( 1_{H};b_{1},0,e_{1},0\right) }B(x_{1}\otimes
gx_{1};GX_{1}^{b_{1}}X_{2},gx_{1}^{e_{1}})GX_{2}\otimes g\otimes gx_{1} \\
&=&\left( -1\right) ^{\alpha \left( 1_{H};0,0,1,0\right) }B(x_{1}\otimes
gx_{1};GX_{2},gx_{1})GX_{2}\otimes g\otimes gx_{1} \\
&&\left( -1\right) ^{\alpha \left( 1_{H};1,0,0,0\right) }B(x_{1}\otimes
gx_{1};GX_{1}X_{2},g)GX_{2}\otimes g\otimes gx_{1}
\end{eqnarray*}%
\begin{eqnarray*}
\alpha \left( 1_{H};0,0,1,0\right) &\equiv &e_{2}+a+b_{1}+b_{2}\equiv
1+1\equiv 0 \\
\alpha \left( 1_{H};1,0,0,0\right) &\equiv &b_{2}=1
\end{eqnarray*}%
\begin{equation*}
\left[ B(x_{1}\otimes gx_{1};GX_{2},gx_{1})-B(x_{1}\otimes
gx_{1};GX_{1}X_{2},g)\right] GX_{2}\otimes g\otimes gx_{1}
\end{equation*}%
Considering also the right side we get%
\begin{equation*}
B(g\otimes gx_{1};GX_{2},g)+B(x_{1}\otimes
gx_{1};GX_{2},gx_{1})-B(x_{1}\otimes gx_{1};GX_{1}X_{2},g)-B(x_{1}\otimes
g;GX_{2},g)=0
\end{equation*}%
which holds in view of the form of the elements.

\subsubsection{Case $G\otimes gx_{1}\otimes gx_{2}$}

Already considered in $B(x_{1}\otimes gx_{1};G,gx_{1}x_{2}).$

\subsubsection{Case $G\otimes g\otimes x_{1}x_{2}$}

First summand%
\begin{eqnarray*}
&&\sum_{a,b_{1},b_{2},d,e_{1},e_{2}=0}^{1}\sum_{l_{1}=0}^{b_{1}}%
\sum_{l_{2}=0}^{b_{2}}\sum_{u_{1}=0}^{e_{1}}\sum_{u_{2}=0}^{e_{2}}\left(
-1\right) ^{\alpha \left( x_{1};l_{1},l_{2},u_{1},u_{2}\right) } \\
&&B(g\otimes
gx_{1};G^{a}X_{1}^{b_{1}}X_{2}^{b_{2}},g^{d}x_{1}^{e_{1}}x_{2}^{e_{2}})G^{a}X_{1}^{b_{1}-l_{1}}X_{2}^{b_{2}-l_{2}}\otimes g^{d}x_{1}^{e_{1}-u_{1}}x_{2}^{e_{2}-u_{2}}\otimes
\\
&&g^{a+b_{1}+b_{2}+l_{1}+l_{2}+d+e_{1}+e_{2}+u_{1}+u_{2}}x_{1}^{l_{1}+u_{1}+1}x_{2}^{l_{2}+u_{2}}
\end{eqnarray*}%
\begin{eqnarray*}
l_{1} &=&u_{1}=0 \\
l_{2}+u_{2} &=&1 \\
a+b_{1}+b_{2}+d+e_{1}+e_{2} &\equiv &1 \\
d &=&1,e_{1}=0,e_{2}=u_{2},a=1,b_{1}=0,b_{2}=l_{2}.
\end{eqnarray*}%
\begin{eqnarray*}
\alpha \left( x_{1};0,0,0,1\right) &\equiv &a+b_{1}+b_{2}\equiv 1+0+0\equiv 1
\\
\alpha \left( x_{1};0,1,0,0\right) &\equiv &a+b_{1}+b_{2}+1\equiv
1+0+1+1\equiv 1
\end{eqnarray*}%
\begin{eqnarray*}
&&\sum_{\substack{ b_{2},e_{2}=0  \\ b_{2}+u_{2}\equiv 1}}^{1}\left(
-1\right) ^{\alpha \left( x_{1};0,b_{2},0,e_{2}\right) }B(g\otimes
gx_{1};GX_{2}^{b_{2}},gx_{2}^{e_{2}})G\otimes g\otimes x_{1}x_{2} \\
&=&\left[ -B(g\otimes gx_{1};G,gx_{2})-B(g\otimes gx_{1};GX_{2},g)\right]
G\otimes g\otimes x_{1}x_{2}
\end{eqnarray*}

\begin{equation*}
-B(g\otimes gx_{1};G,gx_{2})-B(g\otimes gx_{1};GX_{2},g)
\end{equation*}%
Second summand%
\begin{eqnarray*}
l_{1}+u_{1} &=&1 \\
l_{2}+u_{2} &=&1 \\
a+b_{1}+b_{2}+d+e_{1}+e_{2} &\equiv &0 \\
d &=&1,e_{1}=u_{1},e_{2}=u_{2},a=1,b_{1}=l_{1},b_{2}=l_{2}=0
\end{eqnarray*}%
\begin{eqnarray*}
&&\alpha \left( 1_{H};0,0,1,1\right) \equiv 1+1\equiv 0 \\
&&\alpha \left( 1_{H};1,0,0,1\right) \equiv b_{2}+a+b_{1}+b_{2}\equiv
1+1\equiv 0 \\
&&\alpha \left( 1_{H};0,1,1,0\right) \equiv e_{2}+a+b_{1}+b_{2}+1\equiv
0+0+1\equiv 1 \\
&&\alpha \left( 1_{H};1,1,0,0\right) \equiv 1+b_{2}\equiv 0
\end{eqnarray*}%
\begin{eqnarray*}
&&\sum_{\substack{ b_{1},e_{1},b_{2},e_{2}  \\ b_{1}+e_{1}\equiv 1  \\ %
b_{2}+e_{2}\equiv 1}}\left( -1\right) ^{\alpha \left(
1_{H};b_{1},b_{2},e_{1},e_{2}\right) }B(x_{1}\otimes
gx_{1};GX_{1}^{b_{1}}X_{2}^{b_{2}},gx_{1}^{e_{1}}x_{2}^{e_{2}})G\otimes
g\otimes x_{1}x_{2} \\
&&B(x_{1}\otimes gx_{1};G,gx_{1}x_{2})G\otimes g\otimes x_{1}x_{2}+ \\
&&B(x_{1}\otimes gx_{1};GX_{1},gx_{2})G\otimes g\otimes x_{1}x_{2}+ \\
&&-B(x_{1}\otimes gx_{1};GX_{2},gx_{1})G\otimes g\otimes x_{1}x_{2}+ \\
&&B(x_{1}\otimes gx_{1};GX_{1}X_{2},g)G\otimes g\otimes x_{1}x_{2}+
\end{eqnarray*}

Since there is nothing on the right side we get%
\begin{gather*}
-B(g\otimes gx_{1};G,gx_{2})-B(g\otimes gx_{1};GX_{2},g)+B(x_{1}\otimes
gx_{1};G,gx_{1}x_{2})+ \\
+B(x_{1}\otimes gx_{1};GX_{1},gx_{2})-B(x_{1}\otimes
gx_{1};GX_{2},gx_{1})+B(x_{1}\otimes gx_{1};GX_{1}X_{2},g)=0
\end{gather*}%
which we already considered in Case $G\otimes g\otimes x_{1}x_{2}$ of $%
B(x_{1}\otimes gx_{1};G,gx_{1}x_{2}).$

\section{$B \left( x_{1}\otimes gx_{2}\right) $}

By using $\left( \ref{simplgx}\right) ,$ we get%
\begin{eqnarray}
&&B(x_{1}\otimes gx_{2})  \label{form x1otgx2} \\
&=&(1_{A}\otimes g)B(gx_{1}\otimes 1_{H})(1_{A}\otimes gx_{2})  \notag \\
&&+(1_{A}\otimes x_{2})B(gx_{1}\otimes 1_{H})  \notag \\
&&-B(gx_{2}gx_{1}\otimes 1_{H})  \notag
\end{eqnarray}%
so that we obtain

\begin{eqnarray*}
B(x_{1}\otimes gx_{2}) &=&-B(x_{1}x_{2}\otimes
1_{H};1_{A},1_{H})1_{A}\otimes 1_{H}+ \\
&&-B(x_{1}x_{2}\otimes 1_{H};1_{A},x_{1}x_{2})1_{A}\otimes x_{1}x_{2}+ \\
&&-B(x_{1}x_{2}\otimes 1_{H};1_{A},gx_{1})1_{A}\otimes gx_{1}+ \\
&&\left[ B(x_{1}x_{2}\otimes 1_{H};1_{A},gx_{2})+2B(x_{1}x_{2}\otimes
1_{H};X_{2},g)\right] 1_{A}\otimes gx_{2}+ \\
&&-B(x_{1}x_{2}\otimes 1_{H};G,g)G\otimes g+ \\
&&-B(x_{1}x_{2}\otimes 1_{H};G,x_{1})G\otimes x_{1}+ \\
&&-B(x_{1}x_{2}\otimes 1_{H};G,x_{2})G\otimes x_{2}+ \\
&&\left[ B(x_{1}x_{2}\otimes 1_{H};G,gx_{1}x_{2})-2B(x_{1}x_{2}\otimes
1_{H};GX_{2},gx_{1})\right] G\otimes gx_{1}x_{2}+ \\
&&-B(x_{1}x_{2}\otimes 1_{H};X_{1},g)X_{1}\otimes g+ \\
&&-B(x_{1}x_{2}\otimes 1_{H};X_{1},x_{1})X_{1}\otimes x_{1}+ \\
&&-B(x_{1}x_{2}\otimes 1_{H};X_{1},x_{2})X_{1}\otimes x_{2}+ \\
&&-B(x_{1}x_{2}\otimes 1_{H};X_{1},gx_{1}x_{2})X_{1}\otimes gx_{1}x_{2}+ \\
&&-B(x_{1}x_{2}\otimes 1_{H};X_{2},g)X_{2}\otimes g+ \\
&&-B(x_{1}x_{2}\otimes 1_{H};X_{2},x_{1})X_{2}\otimes x_{1}+ \\
&&-B(x_{1}x_{2}\otimes 1_{H};X_{2},x_{2})X_{2}\otimes x_{2}+ \\
&&+B(x_{1}x_{2}\otimes 1_{H};X_{2},gx_{1}x_{2})X_{2}\otimes gx_{1}x_{2}+ \\
&&-\left[
\begin{array}{c}
+1-B(x_{1}x_{2}\otimes 1_{H};1_{A},x_{1}x_{2})- \\
B(x_{1}x_{2}\otimes 1_{H};X_{2},x_{1})+B(x_{1}x_{2}\otimes 1_{H};X_{1},x_{2})%
\end{array}%
\right] X_{1}X_{2}\otimes 1_{H}+ \\
&&-B(x_{1}x_{2}\otimes 1_{H};X_{1},gx_{1}x_{2})X_{1}X_{2}\otimes gx_{1} \\
&&+B(x_{1}x_{2}\otimes 1_{H};X_{2},gx_{1}x_{2})X_{1}X_{2}\otimes gx_{2} \\
&&-B(x_{1}x_{2}\otimes 1_{H};GX_{1},1_{H})GX_{1}\otimes 1_{H}+ \\
&&-B(x_{1}x_{2}\otimes 1_{H};GX_{1},x_{1}x_{2})GX_{1}\otimes x_{1}x_{2} \\
&&-B(x_{1}x_{2}\otimes 1_{H};GX_{1},gx_{1})GX_{1}\otimes gx_{1} \\
&&\left[
\begin{array}{c}
-2B(x_{1}x_{2}\otimes 1_{H};G,gx_{1}x_{2})+ \\
+2B(x_{1}x_{2}\otimes 1_{H};GX_{2},gx_{1})-B(x_{1}x_{2}\otimes
1_{H};GX_{1},gx_{2})%
\end{array}%
\right] GX_{1}\otimes gx_{2}+ \\
&&-B(x_{1}x_{2}\otimes 1_{H};GX_{2},1_{H})GX_{2}\otimes 1_{H}+ \\
&&-B(x_{1}x_{2}\otimes 1_{H};GX_{2},x_{1}x_{2})GX_{2}\otimes x_{1}x_{2}+ \\
&&-B(x_{1}x_{2}\otimes 1_{H};GX_{2},gx_{1})GX_{2}\otimes gx_{1}+ \\
&&+B(x_{1}x_{2}\otimes 1_{H};GX_{2},gx_{2})GX_{2}\otimes gx_{2}+ \\
&&-\left[
\begin{array}{c}
-B(x_{1}x_{2}\otimes 1_{H};G,gx_{1}x_{2})+ \\
+B(x_{1}x_{2}\otimes 1_{H};GX_{2},gx_{1})-B(x_{1}x_{2}\otimes
1_{H};GX_{1},gx_{2})%
\end{array}%
\right] GX_{1}X_{2}\otimes g \\
&&+B(x_{1}x_{2}\otimes 1_{H};GX_{1},x_{1}x_{2})GX_{1}X_{2}\otimes x_{1}+ \\
&&+B(x_{1}x_{2}\otimes 1_{H};GX_{2},x_{1}x_{2})GX_{1}X_{2}\otimes x_{2}
\end{eqnarray*}%
We write the Casimir formula for $B(x_{1}\otimes gx_{2})$%
\begin{eqnarray*}
&&\sum_{a,b_{1},b_{2},d,e_{1},e_{2}=0}^{1}\sum_{l_{1}=0}^{b_{1}}%
\sum_{l_{2}=0}^{b_{2}}\sum_{u_{1}=0}^{e_{1}}\sum_{u_{2}=0}^{e_{2}}\left(
-1\right) ^{\alpha \left( x_{1};l_{1},l_{2},u_{1},u_{2}\right) } \\
&&B(g\otimes
gx_{2};G^{a}X_{1}^{b_{1}}X_{2}^{b_{2}},g^{d}x_{1}^{e_{1}}x_{2}^{e_{2}})G^{a}X_{1}^{b_{1}-l_{1}}X_{2}^{b_{2}-l_{2}}
\\
&&\otimes g^{d}x_{1}^{e_{1}-u_{1}}x_{2}^{e_{2}-u_{2}}\otimes
g^{a+b_{1}+b_{2}+l_{1}+l_{2}+d+e_{1}+e_{2}+u_{1}+u_{2}}x_{1}^{l_{1}+u_{1}+1}x_{2}^{l_{2}+u_{2}}
\\
&&\sum_{a,b_{1},b_{2},d,e_{1},e_{2}=0}^{1}\sum_{l_{1}=0}^{b_{1}}%
\sum_{l_{2}=0}^{b_{2}}\sum_{u_{1}=0}^{e_{1}}\sum_{u_{2}=0}^{e_{2}}\left(
-1\right) ^{\alpha \left( 1_{H};l_{1},l_{2},u_{1},u_{2}\right) } \\
&&B(x_{1}\otimes
gx_{2};G^{a}X_{1}^{b_{1}}X_{2}^{b_{2}},g^{d}x_{1}^{e_{1}}x_{2}^{e_{2}})G^{a}X_{1}^{b_{1}-l_{1}}X_{2}^{b_{2}-l_{2}}
\\
&&\otimes g^{d}x_{1}^{e_{1}-u_{1}}x_{2}^{e_{2}-u_{2}}\otimes
g^{a+b_{1}+b_{2}+l_{1}+l_{2}+d+e_{1}+e_{2}+u_{1}+u_{2}}x_{1}^{l_{1}+u_{1}}x_{2}^{l_{2}+u_{2}}
\\
&=&B^{A}(x_{1}\otimes gx_{2})\otimes B^{H}(x_{1}\otimes gx_{2})\otimes 1_{H}+
\\
&&B^{A}(x_{1}\otimes g)\otimes B^{H}(x_{1}\otimes g)\otimes gx_{2}
\end{eqnarray*}

\subsection{Cases $1_{H}$}

give us%
\begin{equation*}
B(x_{1}\otimes
gx_{2};G^{a}X_{1}^{b_{1}}X_{2}^{b_{2}},g^{d}x_{1}^{e_{1}}x_{2}^{e_{2}})=0%
\text{ whenever }a+b_{1}+b_{2}+d+e_{1}+e_{2}\equiv 1.
\end{equation*}

\subsection{Case $g$}

\begin{eqnarray*}
a+b_{1}+b_{2}+d+e_{1}+e_{2} &\equiv &1 \\
l_{1}+u_{1} &=&l_{2}+u_{2}=0
\end{eqnarray*}%
Since%
\begin{equation*}
B(x_{1}\otimes
gx_{2};G^{a}X_{1}^{b_{1}}X_{2}^{b_{2}},g^{d}x_{1}^{e_{1}}x_{2}^{e_{2}})=0%
\text{ whenever }a+b_{1}+b_{2}+d+e_{1}+e_{2}\equiv 1.
\end{equation*}
the Casimir condition is satisfied.

\subsection{Case $x_{1}$}

\begin{eqnarray*}
&&\sum_{a,b_{1},b_{2},d,e_{1},e_{2}=0}^{1}\sum_{l_{1}=0}^{b_{1}}%
\sum_{l_{2}=0}^{b_{2}}\sum_{u_{1}=0}^{e_{1}}\sum_{u_{2}=0}^{e_{2}}\left(
-1\right) ^{\alpha \left( x_{1};l_{1},l_{2},u_{1},u_{2}\right) } \\
&&B(g\otimes
gx_{2};G^{a}X_{1}^{b_{1}}X_{2}^{b_{2}},g^{d}x_{1}^{e_{1}}x_{2}^{e_{2}}) \\
&&G^{a}X_{1}^{b_{1}-l_{1}}X_{2}^{b_{2}-l_{2}}\otimes
g^{d}x_{1}^{e_{1}-u_{1}}x_{2}^{e_{2}-u_{2}}\otimes
g^{a+b_{1}+b_{2}+l_{1}+l_{2}+d+e_{1}+e_{2}+u_{1}+u_{2}}x_{1}^{l_{1}+u_{1}+1}x_{2}^{l_{2}+u_{2}}
\\
&&\sum_{a,b_{1},b_{2},d,e_{1},e_{2}=0}^{1}\sum_{l_{1}=0}^{b_{1}}%
\sum_{l_{2}=0}^{b_{2}}\sum_{u_{1}=0}^{e_{1}}\sum_{u_{2}=0}^{e_{2}}\left(
-1\right) ^{\alpha \left( 1_{H};l_{1},l_{2},u_{1},u_{2}\right) } \\
&&B(x_{1}\otimes
gx_{2};G^{a}X_{1}^{b_{1}}X_{2}^{b_{2}},g^{d}x_{1}^{e_{1}}x_{2}^{e_{2}}) \\
&&G^{a}X_{1}^{b_{1}-l_{1}}X_{2}^{b_{2}-l_{2}}\otimes
g^{d}x_{1}^{e_{1}-u_{1}}x_{2}^{e_{2}-u_{2}}\otimes
g^{a+b_{1}+b_{2}+l_{1}+l_{2}+d+e_{1}+e_{2}+u_{1}+u_{2}}x_{1}^{l_{1}+u_{1}}x_{2}^{l_{2}+u_{2}}
\\
&=&B^{A}(x_{1}\otimes gx_{2})\otimes B^{H}(x_{1}\otimes gx_{2})\otimes 1_{H}+
\\
&&B^{A}(x_{1}\otimes g)\otimes B^{H}(x_{1}\otimes g)\otimes gx_{2}
\end{eqnarray*}

First summand is zero since%
\begin{equation*}
B(1_{H}\otimes
x_{2};G^{a}X_{1}^{b_{1}}X_{2}^{b_{2}},g^{d}x_{1}^{e_{1}}x_{2}^{e_{2}})=0%
\text{ whenever }a+b_{1}+b_{2}+d+e_{1}+e_{2}\equiv 0
\end{equation*}%
The second summand is going to be also zero as%
\begin{eqnarray*}
l_{1}+u_{1} &=&1 \\
l_{2}+u_{2} &=&0 \\
a+b_{1}+b_{2}+d+e_{1}+e_{2} &\equiv &1
\end{eqnarray*}%
and we know that these are zero from above.

\subsection{Case $x_{2}$}

The first summand is zero. The second summand is also zero as%
\begin{eqnarray*}
l_{1}+u_{1} &=&0 \\
l_{2}+u_{2} &=&1 \\
a+b_{1}+b_{2}+d+e_{1}+e_{2} &\equiv &1
\end{eqnarray*}%
and we know that these are zero from above.

\subsection{Case $x_{1}x_{2}$}

From the first summand we get%
\begin{eqnarray*}
l_{1}+u_{1} &=&0 \\
l_{2}+u_{2} &=&1 \\
a+b_{1}+b_{2}+d+e_{1}+e_{2} &\equiv &1
\end{eqnarray*}%
\begin{gather*}
\sum_{\substack{ a,b_{1},b_{2},d,e_{1},e_{2}=0  \\ %
a+b_{1}+b_{2}+d+e_{1}+e_{2}\equiv 1}}^{1}\sum_{l_{2}=0}^{b_{2}}%
\sum_{u_{2}=0}^{e_{2}}\left( -1\right) ^{\alpha \left(
x_{1};l_{1},l_{2},u_{1},u_{2}\right) } \\
B(g\otimes
gx_{2};G^{a}X_{1}^{b_{1}}X_{2}^{b_{2}},g^{d}x_{1}^{e_{1}}x_{2}^{e_{2}})G^{a}X_{1}^{b_{1}-l_{1}}X_{2}^{b_{2}-l_{2}}\otimes g^{d}x_{1}^{e_{1}-u_{1}}x_{2}^{e_{2}-u_{2}}
\end{gather*}%
Since we have $\alpha \left( x_{1};0,0,0,1\right) \equiv 1$ and $\alpha
\left( x_{1};0,1,0,0\right) \equiv a+b_{1}+b_{2}+1$ we obtain%
\begin{gather*}
\sum_{\substack{ a,b_{1},b_{2},d,e_{1},e_{2}=0  \\ %
a+b_{1}+b_{2}+d+e_{1}+e_{2}\equiv 1}}^{1}\sum_{l_{2}=0}^{b_{2}}\sum
_{\substack{ u_{2}=0  \\ l_{2}+u_{2}=1}}^{e_{2}}\left( -1\right) ^{\alpha
\left( x_{1};0,l_{2},0,u_{2}\right) }B(g\otimes
gx_{2};G^{a}X_{1}^{b_{1}}X_{2}^{b_{2}},g^{d}x_{1}^{e_{1}}x_{2}^{e_{2}}) \\
G^{a}X_{1}^{b_{1}}X_{2}^{b_{2}-l_{2}}\otimes
g^{d}x_{1}^{e_{1}}x_{2}^{e_{2}-u_{2}}= \\
=\sum_{\substack{ a,b_{1},b_{2},d,e_{1}=0  \\ a+b_{1}+b_{2}+d+e_{1}\equiv 0}}%
^{1}-B(g\otimes
gx_{2};G^{a}X_{1}^{b_{1}}X_{2}^{b_{2}},g^{d}x_{1}^{e_{1}}x_{2})G^{a}X_{1}^{b_{1}}X_{2}^{b_{2}}\otimes g^{d}x_{1}^{e_{1}}+
\\
+\sum_{\substack{ a,b_{1},d,e_{1},e_{2}=0  \\ a+b_{1}+d+e_{1}+e_{2}\equiv 0}}%
^{1}\left( -1\right) ^{a+b_{1}}B(g\otimes
gx_{2};G^{a}X_{1}^{b_{1}}X_{2},g^{d}x_{1}^{e_{1}}x_{2}^{e_{2}})G^{a}X_{1}^{b_{1}}\otimes g^{d}x_{1}^{e_{1}}x_{2}^{e_{2}}
\end{gather*}%
From the second summand of the left side we get%
\begin{eqnarray*}
l_{1}+u_{1} &=&1 \\
l_{2}+u_{2} &=&1 \\
a+b_{1}+b_{2}+d+e_{1}+e_{2} &\equiv &0
\end{eqnarray*}%
and we get
\begin{eqnarray*}
&&\sum_{\substack{ a,b_{1},b_{2},d,e_{1},e_{2}=0  \\ %
a+b_{1}+b_{2}+d+e_{1}+e_{2}\equiv 0}}^{1}\sum_{l_{1}=0}^{b_{1}}%
\sum_{l_{2}=0}^{b_{2}}\sum_{\substack{ u_{1}=0  \\ l_{1}+u_{1}=1}}%
^{e_{1}}\sum _{\substack{ u_{2}=0  \\ l_{2}+u_{2}=1}}^{e_{2}}\left(
-1\right) ^{\alpha \left( 1_{H};l_{1},l_{2},u_{1},u_{2}\right) } \\
&&B(x_{1}\otimes
gx_{2};G^{a}X_{1}^{b_{1}}X_{2}^{b_{2}},g^{d}x_{1}^{e_{1}}x_{2}^{e_{2}})G^{a}X_{1}^{b_{1}-l_{1}}X_{2}^{b_{2}-l_{2}}\otimes g^{d}x_{1}^{e_{1}-u_{1}}x_{2}^{e_{2}-u_{2}}\otimes x_{1}x_{2}
\end{eqnarray*}%
Since%
\begin{eqnarray*}
&&\alpha \left( 1_{H};0,0,1,1\right) \equiv 1+e_{2}\equiv 0 \\
&&\alpha \left( 1_{H};0,1,1,0\right) \equiv e_{2}+a+b_{1} \\
&&\alpha \left( 1_{H};1,0,0,1\right) \equiv a+1 \\
&&\alpha \left( 1_{H};1,1,0,0\right) \equiv 1+1\equiv 0
\end{eqnarray*}%
\begin{eqnarray*}
&&\sum_{\substack{ a,b_{1},b_{2},d=0  \\ a+b_{1}+b_{2}+d\equiv 0}}%
^{1}B(x_{1}\otimes
gx_{2};G^{a}X_{1}^{b_{1}}X_{2}^{b_{2}},g^{d}x_{1}x_{2})G^{a}X_{1}^{b_{1}}X_{2}^{b_{2}}\otimes g^{d}\otimes x_{1}x_{2}
\\
&&\sum_{\substack{ a,b_{1},d,e_{2}=0  \\ a+b_{1}+d+e_{2}\equiv 0}}^{1}\left(
-1\right) ^{e_{2}+a+b_{1}}B(x_{1}\otimes
gx_{2};G^{a}X_{1}^{b_{1}}X_{2},g^{d}x_{1}x_{2}^{e_{2}})G^{a}X_{1}^{b_{1}}%
\otimes g^{d}x_{2}^{e_{2}}\otimes x_{1}x_{2} \\
&&\sum_{\substack{ a,b_{2},d,e_{1}=0  \\ a+b_{2}+d+e_{1}\equiv 0}}^{1}\left(
-1\right) ^{a+1}B(x_{1}\otimes
gx_{2};G^{a}X_{1}X_{2}^{b_{2}},g^{d}x_{1}^{e_{1}}x_{2})G^{a}X_{2}^{b_{2}}%
\otimes g^{d}x_{1}^{e_{1}}\otimes x_{1}x_{2} \\
&&\sum_{\substack{ a,d,e_{1},e_{2}=0  \\ a+d+e_{1}+e_{2}\equiv 0}}%
^{1}B(x_{1}\otimes
gx_{2};G^{a}X_{1}X_{2},g^{d}x_{1}^{e_{1}}x_{2}^{e_{2}})G^{a}\otimes
g^{d}x_{1}^{e_{1}}x_{2}^{e_{2}}\otimes x_{1}x_{2}
\end{eqnarray*}%
and we finally get%
\begin{eqnarray*}
&&\sum_{\substack{ a,b_{1},b_{2},d,e_{1}=0  \\ a+b_{1}+b_{2}+d+e_{1}\equiv 0
}}^{1}-B(g\otimes
gx_{2};G^{a}X_{1}^{b_{1}}X_{2}^{b_{2}},g^{d}x_{1}^{e_{1}}x_{2})G^{a}X_{1}^{b_{1}}X_{2}^{b_{2}}\otimes g^{d}x_{1}^{e_{1}}
\\
&&+\sum_{\substack{ a,b_{1},d,e_{1},e_{2}=0  \\ a+b_{1}+d+e_{1}+e_{2}\equiv
0 }}^{1}\left( -1\right) ^{a+b_{1}}B(g\otimes
gx_{2};G^{a}X_{1}^{b_{1}}X_{2},g^{d}x_{1}^{e_{1}}x_{2}^{e_{2}})G^{a}X_{1}^{b_{1}}\otimes g^{d}x_{1}^{e_{1}}x_{2}^{e_{2}}
\\
+ &&\sum_{\substack{ a,b_{1},b_{2},d=0  \\ a+b_{1}+b_{2}+d\equiv 0}}%
^{1}B(x_{1}\otimes
gx_{2};G^{a}X_{1}^{b_{1}}X_{2}^{b_{2}},g^{d}x_{1}x_{2})G^{a}X_{1}^{b_{1}}X_{2}^{b_{2}}\otimes g^{d}
\\
+ &&\sum_{\substack{ a,b_{1},d,e_{2}=0  \\ a+b_{1}+d+e_{2}\equiv 0}}%
^{1}\left( -1\right) ^{e_{2}+a+b_{1}}B(x_{1}\otimes
gx_{2};G^{a}X_{1}^{b_{1}}X_{2},g^{d}x_{1}x_{2}^{e_{2}})G^{a}X_{1}^{b_{1}}%
\otimes g^{d}x_{2}^{e_{2}} \\
+ &&\sum_{\substack{ a,b_{2},d,e_{1}=0  \\ a+b_{2}+d+e_{1}\equiv 0}}%
^{1}\left( -1\right) ^{a+1}B(x_{1}\otimes
gx_{2};G^{a}X_{1}X_{2}^{b_{2}},g^{d}x_{1}^{e_{1}}x_{2})G^{a}X_{2}^{b_{2}}%
\otimes g^{d}x_{1}^{e_{1}} \\
+ &&\sum_{\substack{ a,d,e_{1},e_{2}=0  \\ a+d+e_{1}+e_{2}\equiv 0}}%
^{1}B(x_{1}\otimes
gx_{2};G^{a}X_{1}X_{2},g^{d}x_{1}^{e_{1}}x_{2}^{e_{2}})G^{a}\otimes
g^{d}x_{1}^{e_{1}}x_{2}^{e_{2}}=0
\end{eqnarray*}

\subsubsection{$G^{a}\otimes g^{d}$}

\begin{equation*}
\sum_{\substack{ a,d  \\ a+d=0}}\left[
\begin{array}{c}
-B(g\otimes gx_{2};G^{a},g^{d}x_{2})+\left( -1\right) ^{a}B(g\otimes
gx_{2};G^{a}X_{2},g^{d}) \\
+B(x_{1}\otimes gx_{2};G^{a},g^{d}x_{1}x_{2})+\left( -1\right)
^{a}B(x_{1}\otimes gx_{2};G^{a}X_{2},g^{d}x_{1}) \\
+\left( -1\right) ^{a+1}B(x_{1}\otimes
gx_{2};G^{a}X_{1},g^{d}x_{2})+B(x_{1}\otimes gx_{2};G^{a}X_{1}X_{2},g^{d})%
\end{array}%
\right] G^{a}\otimes g^{d}=0
\end{equation*}%
and we get%
\begin{equation*}
\begin{array}{c}
-B(g\otimes gx_{2};1_{A},x_{2})+B(g\otimes
gx_{2};X_{2},1_{H})+B(x_{1}\otimes gx_{2};1_{A},x_{1}x_{2})+ \\
+B(x_{1}\otimes gx_{2};X_{2},x_{1})-B(x_{1}\otimes
gx_{2};X_{1},x_{2})+B(x_{1}\otimes gx_{2};X_{1}X_{2},1_{H})=0%
\end{array}%
\end{equation*}%
and%
\begin{equation*}
\begin{array}{c}
-B(g\otimes gx_{2};G,gx_{2})-B(g\otimes gx_{2};GX_{2},g)+B(x_{1}\otimes
gx_{2};G,gx_{1}x_{2})+ \\
-B(x_{1}\otimes gx_{2};GX_{2},gx_{1})+B(x_{1}\otimes
gx_{2};GX_{1},gx_{2})+B(x_{1}\otimes gx_{2};GX_{1}X_{2},g)=0%
\end{array}%
\end{equation*}%
which are satisfied in view of the form of the elements.

\subsubsection{$G^{a}\otimes g^{d}x_{2}$}

We get%
\begin{equation*}
\sum_{\substack{ a,d  \\ a+d=1}}\left[
\begin{array}{c}
\left( -1\right) ^{a}B(g\otimes gx_{2};G^{a}X_{2},g^{d}x_{2})+ \\
\left( -1\right) ^{1+a}B(x_{1}\otimes
gx_{2};G^{a}X_{2},g^{d}x_{1}x_{2})+B(x_{1}\otimes
gx_{2};G^{a}X_{1}X_{2},g^{d}x_{2})%
\end{array}%
\right] G^{a}\otimes g^{d}x_{2}=0
\end{equation*}%
and we obtain%
\begin{equation*}
B(g\otimes gx_{2};X_{2},gx_{2})-B(x_{1}\otimes
gx_{2};X_{2},gx_{1}x_{2})+B(x_{1}\otimes gx_{2};X_{1}X_{2},gx_{2})=0
\end{equation*}%
and%
\begin{equation*}
-B(g\otimes gx_{2};GX_{2},x_{2})+B(x_{1}\otimes
gx_{2};GX_{2},x_{1}x_{2})+B(x_{1}\otimes gx_{2};GX_{1}X_{2},x_{2})=0
\end{equation*}%
which hold in view of the form of the elements.

\subsubsection{$G^{a}\otimes g^{d}x_{1}$}

We obtain%
\begin{equation*}
\sum_{\substack{ a,d  \\ a+d=1}}\left[
\begin{array}{c}
-B(g\otimes gx_{2};G^{a},g^{d}x_{1}x_{2})+\left( -1\right) ^{a}B(g\otimes
gx_{2};G^{a}X_{2},g^{d}x_{1})+ \\
\left( -1\right) ^{a+1}B(x_{1}\otimes
gx_{2};G^{a}X_{1},g^{d}x_{1}x_{2})+B(x_{1}\otimes
gx_{2};G^{a}X_{1}X_{2},g^{d}x_{1})%
\end{array}%
\right] G^{a}\otimes g^{d}x_{1}=0
\end{equation*}%
and we get%
\begin{equation*}
\begin{array}{c}
-B(g\otimes gx_{2};1_{A},gx_{1}x_{2})+B(g\otimes gx_{2};X_{2},gx_{1})+ \\
-B(x_{1}\otimes gx_{2};X_{1},gx_{1}x_{2})+B(x_{1}\otimes
gx_{2};X_{1}X_{2},gx_{1})=0%
\end{array}%
\end{equation*}%
and%
\begin{equation*}
\begin{array}{c}
-B(g\otimes gx_{2};G,x_{1}x_{2})-B(g\otimes gx_{2};GX_{2},x_{1})+ \\
+B(x_{1}\otimes gx_{2};GX_{1},x_{1}x_{2})+B(x_{1}\otimes
gx_{2};GX_{1}X_{2},x_{1})=0%
\end{array}%
\end{equation*}%
which hold in view of the form of the elements

\subsubsection{$G^{a}X_{2}\otimes g^{d}$}

We obtain
\begin{equation*}
\sum_{\substack{ a,d  \\ a+d=1}}\left[
\begin{array}{c}
-B(g\otimes gx_{2};G^{a}X_{2},g^{d}x_{2})+B(x_{1}\otimes
gx_{2};G^{a}X_{2},g^{d}x_{1}x_{2}) \\
+\left( -1\right) ^{a+1}B(x_{1}\otimes gx_{2};G^{a}X_{1}X_{2},g^{d}x_{2})%
\end{array}%
\right] G^{a}X_{2}\otimes g^{d}=0
\end{equation*}%
and we get%
\begin{equation*}
-B(g\otimes gx_{2};X_{2},gx_{2})+B(x_{1}\otimes
gx_{2};X_{2},gx_{1}x_{2})-B(x_{1}\otimes gx_{2};X_{1}X_{2},gx_{2})=0
\end{equation*}%
and%
\begin{equation*}
-B(g\otimes gx_{2};GX_{2},x_{2})+B(x_{1}\otimes
gx_{2};GX_{2},x_{1}x_{2})+B(x_{1}\otimes gx_{2};GX_{1}X_{2},x_{2})=0
\end{equation*}%
which hold in view of the form of the elements.

\subsubsection{$G^{a}X_{1}\otimes g^{d}$}

We obtain%
\begin{equation*}
\sum_{\substack{ a,d  \\ a+d=1}}\left[
\begin{array}{c}
-B(g\otimes gx_{2};G^{a}X_{1},g^{d}x_{2})+\left( -1\right) ^{a+1}B(g\otimes
gx_{2};G^{a}X_{1}X_{2},g^{d}) \\
+B(x_{1}\otimes gx_{2};G^{a}X_{1},g^{d}x_{1}x_{2})+\left( -1\right)
^{a+1}B(x_{1}\otimes gx_{2};G^{a}X_{1}X_{2},g^{d}x_{1})%
\end{array}%
\right] G^{a}X_{1}\otimes g^{d}=0
\end{equation*}%
and we get
\begin{equation*}
\begin{array}{c}
-B(g\otimes gx_{2};X_{1},gx_{2})-B(g\otimes gx_{2};X_{1}X_{2},g) \\
+B(x_{1}\otimes gx_{2};X_{1},gx_{1}x_{2})-B(x_{1}\otimes
gx_{2};X_{1}X_{2},gx_{1})=0%
\end{array}%
\end{equation*}%
and which hold in view of the form of the element.

\subsubsection{$G^{a}\otimes g^{d}x_{1}x_{2}$}

We obtain%
\begin{equation*}
\sum_{\substack{ a,d  \\ a+d=0}}\left[ \left( -1\right) ^{a}B(g\otimes
gx_{2};G^{a}X_{2},g^{d}x_{1}x_{2})+B(x_{1}\otimes
gx_{2};G^{a}X_{1}X_{2},g^{d}x_{1}x_{2})\right] G^{a}\otimes g^{d}x_{1}x_{2}
\end{equation*}%
and we get%
\begin{equation*}
B(g\otimes gx_{2};X_{2},x_{1}x_{2})+B(x_{1}\otimes
gx_{2};X_{1}X_{2},x_{1}x_{2})=0
\end{equation*}%
and%
\begin{equation*}
-B(g\otimes gx_{2};GX_{2},gx_{1}x_{2})+B(x_{1}\otimes
gx_{2};GX_{1}X_{2},gx_{1}x_{2})=0
\end{equation*}%
which hold in view of the form of the elements.

\subsubsection{$G^{a}X_{2}\otimes g^{d}x_{2}$}

There is no term like this.

\subsubsection{$G^{a}X_{1}\otimes g^{d}x_{2}$}

We obtain%
\begin{equation*}
\sum_{\substack{ a,d  \\ a+d=0}}\left[ \left( -1\right) ^{a+1}B(g\otimes
gx_{2};G^{a}X_{1}X_{2},g^{d}x_{2})+\left( -1\right) ^{a}B(x_{1}\otimes
gx_{2};G^{a}X_{1}X_{2},g^{d}x_{1}x_{2})\right] G^{a}X_{1}\otimes g^{d}x_{2}=0
\end{equation*}%
and we get%
\begin{equation*}
-B(g\otimes gx_{2};X_{1}X_{2},x_{2})+B(x_{1}\otimes
gx_{2};X_{1}X_{2},x_{1}x_{2})=0
\end{equation*}%
and%
\begin{equation*}
B(g\otimes gx_{2};GX_{1}X_{2},gx_{2})-B(x_{1}\otimes
gx_{2};GX_{1}X_{2},gx_{1}x_{2})=0
\end{equation*}%
\begin{equation*}
B(g\otimes gx_{2};GX_{1}X_{2},gx_{2})-B(x_{1}\otimes
gx_{2};GX_{1}X_{2},gx_{1}x_{2})=0
\end{equation*}%
which hold in view of the form of the elements.

\subsubsection{$G^{a}X_{2}\otimes g^{d}x_{1}$}

We obtain%
\begin{equation*}
\sum_{\substack{ a,d  \\ a+d=0}}\left[ -B(g\otimes
gx_{2};G^{a}X_{2},g^{d}x_{1}x_{2})+\left( -1\right) ^{a+1}B(x_{1}\otimes
gx_{2};G^{a}X_{1}X_{2},g^{d}x_{1}x_{2})\right] G^{a}X_{2}\otimes g^{d}x_{1}=0
\end{equation*}%
and we get%
\begin{equation*}
-B(g\otimes gx_{2};X_{2},x_{1}x_{2})-B(x_{1}\otimes
gx_{2};X_{1}X_{2},x_{1}x_{2})=0
\end{equation*}%
and%
\begin{equation*}
-B(g\otimes gx_{2};GX_{2},gx_{1}x_{2})+B(x_{1}\otimes
gx_{2};GX_{1}X_{2},gx_{1}x_{2})=0
\end{equation*}%
which hold in view of the form of the elements.

\subsubsection{$G^{a}X_{1}\otimes g^{d}x_{1}$}

We obtain
\begin{equation*}
\sum_{\substack{ a,d  \\ a+d=0}}\left[ -B(g\otimes
gx_{2};G^{a}X_{1},g^{d}x_{1}x_{2})+\left( -1\right) ^{a+1}B(g\otimes
gx_{2};G^{a}X_{1}X_{2},g^{d}x_{1})\right] G^{a}X_{1}\otimes g^{d}x_{1}=0
\end{equation*}%
and we get%
\begin{equation*}
-B(g\otimes gx_{2};X_{1},x_{1}x_{2})-B(g\otimes gx_{2};X_{1}X_{2},x_{1})=0
\end{equation*}%
and%
\begin{equation*}
-B(g\otimes gx_{2};GX_{1},gx_{1}x_{2})+B(g\otimes
gx_{2};GX_{1}X_{2},gx_{1})=0
\end{equation*}%
which hold in view of the form of the elements.

\subsubsection{$G^{a}X_{1}X_{2}\otimes g^{d}$}

We obtain%
\begin{equation*}
\sum_{\substack{ a,d  \\ a+d=0}}\left[ -B(g\otimes
gx_{2};G^{a}X_{1}X_{2},g^{d}x_{2})+B(x_{1}\otimes
gx_{2};G^{a}X_{1}X_{2},g^{d}x_{1}x_{2})\right] G^{a}X_{1}X_{2}\otimes g^{d}=0
\end{equation*}%
and we get%
\begin{equation*}
-B(g\otimes gx_{2};X_{1}X_{2},x_{2})+B(x_{1}\otimes
gx_{2};X_{1}X_{2},x_{1}x_{2})=0
\end{equation*}%
and%
\begin{equation*}
-B(g\otimes gx_{2};GX_{1}X_{2},gx_{2})+B(x_{1}\otimes
gx_{2};GX_{1}X_{2},gx_{1}x_{2}))=0
\end{equation*}%
which hold in view of the form of the elements.

\subsubsection{$G^{a}X_{2}\otimes g^{d}x_{1}x_{2}$}

There is no term like this

\subsubsection{$G^{a}X_{1}\otimes g^{d}x_{1}x_{2}$}

We obtain%
\begin{equation*}
\sum_{\substack{ a,d  \\ a+d=1}}\left( -1\right) ^{a+1}B(g\otimes
gx_{2};G^{a}X_{1}X_{2},g^{d}x_{1}x_{2})G^{a}X_{1}\otimes g^{d}x_{1}x_{2}=0
\end{equation*}%
and we get%
\begin{equation*}
-B(g\otimes gx_{2};X_{1}X_{2},gx_{1}x_{2})=0
\end{equation*}%
and%
\begin{equation*}
+B(g\otimes gx_{2};GX_{1}X_{2},x_{1}x_{2})=0
\end{equation*}%
which hold in view of the form of the element.

\subsubsection{$G^{a}X_{1}X_{2}\otimes g^{d}x_{2}$}

There is no term like this.

\subsubsection{$G^{a}X_{1}X_{2}\otimes g^{d}x_{1}$}

We obtain%
\begin{equation*}
\sum_{\substack{ a,d  \\ a+d=1}}-B(g\otimes
gx_{2};G^{a}X_{1}X_{2},g^{d}x_{1}x_{2})G^{a}X_{1}X_{2}\otimes g^{d}x_{1}=0
\end{equation*}%
and we get%
\begin{equation*}
-B(g\otimes gx_{2};X_{1}X_{2},gx_{1}x_{2})=0
\end{equation*}%
and%
\begin{equation*}
-B(g\otimes gx_{2};GX_{1}X_{2},gx_{1}x_{2})=0
\end{equation*}%
which we already obtain in case $G^{a}X_{1}\otimes g^{d}x_{1}x_{2}$ above.

\subsubsection{$G^{a}X_{1}X_{2}\otimes g^{d}x_{1}x_{2}$}

There is no term like this.

\subsection{Case $gx_{1}$}

First summand of the left side gives us%
\begin{eqnarray*}
l_{1} &=&u_{1}=l_{2}=u_{2}=0 \\
a+b_{1}+b_{2}+d+e_{1}+e_{2} &\equiv &1
\end{eqnarray*}%
Since $\alpha \left( x_{1};0,0,0,0\right) =a+b_{1}+b_{2},$ we get%
\begin{equation*}
\sum_{\substack{ a,b_{1},b_{2},d,e_{1},e_{2}=0  \\ %
a+b_{1}+b_{2}+d+e_{1}+e_{2}\equiv 1}}^{1}\left( -1\right)
^{a+b_{1}+b_{2}}B(g\otimes
gx_{2};G^{a}X_{1}^{b_{1}}X_{2}^{b_{2}},g^{d}x_{1}^{e_{1}}x_{2}^{e_{2}})G^{a}X_{1}^{b_{1}}X_{2}^{b_{2}}\otimes g^{d}x_{1}^{e_{1}}x_{2}^{e_{2}}
\end{equation*}%
The second summand of the left side gives us%
\begin{eqnarray*}
l_{1}+u_{1} &=&1 \\
l_{2} &=&u_{2}=0 \\
a+b_{1}+b_{2}+d+e_{1}+e_{2} &\equiv &0 \\
\alpha \left( 1_{H};0,0,1,0\right) &\equiv &e_{2}+\left( a+b_{1}+b_{2}\right)
\\
\alpha \left( 1_{H};1,0,0,0\right) &\equiv &b_{2}
\end{eqnarray*}%
and we get
\begin{eqnarray*}
&&\sum_{\substack{ a,b_{1},b_{2},d,e_{1},e_{2}=0  \\ %
a+b_{1}+b_{2}+d+e_{1}+e_{2}\equiv 0}}^{1}\sum_{l_{1}=0}^{b_{1}}\sum
_{\substack{ u_{1}=0  \\ l_{1}+u_{1}=1}}^{e_{1}}\left( -1\right) ^{\alpha
\left( 1_{H};l_{1},0,u_{1},0\right) }B(x_{1}\otimes
gx_{2};G^{a}X_{1}^{b_{1}}X_{2}^{b_{2}},g^{d}x_{1}^{e_{1}}x_{2}^{e_{2}}) \\
&&G^{a}X_{1}^{b_{1}-l_{1}}X_{2}^{b_{2}}\otimes
g^{d}x_{1}^{e_{1}-u_{1}}x_{2}^{e_{2}} \\
&=&\sum_{\substack{ a,b_{1},b_{2},d,e_{2}=0  \\ a+b_{1}+b_{2}+d+e_{2}\equiv
1 }}^{1}\left( -1\right) ^{e_{2}+\left( a+b_{1}+b_{2}\right) }B(x_{1}\otimes
gx_{2};G^{a}X_{1}^{b_{1}}X_{2}^{b_{2}},g^{d}x_{1}x_{2}^{e_{2}})G^{a}X_{1}^{b_{1}}X_{2}^{b_{2}}\otimes g^{d}x_{2}^{e_{2}}+
\\
&&+\sum_{\substack{ a,b_{2},d,e_{1},e_{2}=0  \\ a+b_{2}+d+e_{1}+e_{2}\equiv
1 }}^{1}\left( -1\right) ^{b_{2}}B(x_{1}\otimes
gx_{2};G^{a}X_{1}X_{2}^{b_{2}},g^{d}x_{1}^{e_{1}}x_{2}^{e_{2}})G^{a}X_{2}^{b_{2}}\otimes g^{d}x_{1}^{e_{1}}x_{2}^{e_{2}}
\end{eqnarray*}%
By considering also the right side, we get

\begin{eqnarray*}
&&\sum_{\substack{ a,b_{1},b_{2},d,e_{1},e_{2}=0  \\ %
a+b_{1}+b_{2}+d+e_{1}+e_{2}\equiv 1}}^{1}\left( -1\right)
^{a+b_{1}+b_{2}}B(g\otimes
gx_{2};G^{a}X_{1}^{b_{1}}X_{2}^{b_{2}},g^{d}x_{1}^{e_{1}}x_{2}^{e_{2}})G^{a}X_{1}^{b_{1}}X_{2}^{b_{2}}\otimes g^{d}x_{1}^{e_{1}}x_{2}^{e_{2}}
\\
&&+\sum_{\substack{ a,b_{1},b_{2},d,e_{2}=0  \\ a+b_{1}+b_{2}+d+e_{2}\equiv
1 }}^{1}\left( -1\right) ^{e_{2}+\left( a+b_{1}+b_{2}\right) }B(x_{1}\otimes
gx_{2};G^{a}X_{1}^{b_{1}}X_{2}^{b_{2}},g^{d}x_{1}x_{2}^{e_{2}})G^{a}X_{1}^{b_{1}}X_{2}^{b_{2}}\otimes g^{d}x_{2}^{e_{2}}+
\\
&&\sum_{\substack{ a,b_{2},d,e_{1},e_{2}=0  \\ a+b_{2}+d+e_{1}+e_{2}\equiv 1
}}^{1}\left( -1\right) ^{b_{2}}B(x_{1}\otimes
gx_{2};G^{a}X_{1}X_{2}^{b_{2}},g^{d}x_{1}^{e_{1}}x_{2}^{e_{2}})G^{a}X_{2}^{b_{2}}\otimes g^{d}x_{1}^{e_{1}}x_{2}^{e_{2}}=0
\end{eqnarray*}

We obtain%
\begin{equation*}
\sum_{\substack{ a,d  \\ a+d\equiv 1}}^{1}\left[
\begin{array}{c}
\left( -1\right) ^{a}B(g\otimes gx_{2};G^{a},g^{d})+\left( -1\right)
^{a}B(x_{1}\otimes gx_{2};G^{a},g^{d}x_{1})+ \\
+B(x_{1}\otimes gx_{2};G^{a}X_{1},g^{d})%
\end{array}%
\right] G^{a}\otimes g^{d}=0
\end{equation*}%
and we get%
\begin{equation*}
B(g\otimes gx_{2};1_{A},g)+B(x_{1}\otimes
gx_{2};1_{A},gx_{1})+B(x_{1}\otimes gx_{2};X_{1},g)=0
\end{equation*}%
and%
\begin{equation*}
-B(g\otimes gx_{2};G,1_{H})-B(x_{1}\otimes gx_{2};G,x_{1})+B(x_{1}\otimes
gx_{2};GX_{1},1_{H})=0
\end{equation*}%
which hold in view of the form of the elements.

\subsubsection{$G^{a}\otimes g^{d}x_{2}$}

We obtain%
\begin{equation*}
\sum_{\substack{ a,d  \\ a+d\equiv 0}}^{1}\left[
\begin{array}{c}
\left( -1\right) ^{a}B(g\otimes gx_{2};G^{a},g^{d}x_{2}) \\
+\left( -1\right) ^{a+1}B(x_{1}\otimes
gx_{2};G^{a},g^{d}x_{1}x_{2})+B(x_{1}\otimes gx_{2};G^{a}X_{1},g^{d}x_{2})%
\end{array}%
\right] G^{a}\otimes g^{d}x_{2}=0
\end{equation*}%
and we get%
\begin{equation*}
B(g\otimes gx_{2};1_{A},x_{2})-B(x_{1}\otimes
gx_{2};1_{A},x_{1}x_{2})+B(x_{1}\otimes gx_{2};X_{1},x_{2})=0
\end{equation*}%
and%
\begin{equation*}
-B(g\otimes gx_{2};G,gx_{2})+B(x_{1}\otimes
gx_{2};G,gx_{1}x_{2})+B(x_{1}\otimes gx_{2};GX_{1},gx_{2})=0
\end{equation*}%
which hold in view of the form of the elements.

\subsubsection{$G^{a}\otimes g^{d}x_{1}$}

We obtain%
\begin{equation*}
\sum_{\substack{ a,d  \\ a+d\equiv 0}}^{1}\left[ \left( -1\right)
^{a}B(g\otimes gx_{2};G^{a},g^{d}x_{1})+B(x_{1}\otimes
gx_{2};G^{a}X_{1},g^{d}x_{1})\right] G^{a}\otimes g^{d}x_{1}=0
\end{equation*}%
and we get%
\begin{equation*}
B(g\otimes gx_{2};1_{A},x_{1})+B(x_{1}\otimes gx_{2};X_{1},x_{1})=0
\end{equation*}%
and%
\begin{equation*}
-B(g\otimes gx_{2};G,gx_{1})+B(x_{1}\otimes gx_{2};GX_{1},gx_{1})=0
\end{equation*}%
which hold in view of the form of the elements.

\subsubsection{$G^{a}X_{2}\otimes g^{d}$}

We obtain%
\begin{equation*}
\sum_{\substack{ a,d  \\ a+d\equiv 0}}^{1}\left[
\begin{array}{c}
\left( -1\right) ^{a+1}B(g\otimes gx_{2};G^{a}X_{2},g^{d})+\left( -1\right)
^{a+1}B(x_{1}\otimes gx_{2};G^{a}X_{2},g^{d}x_{1})+ \\
-B(x_{1}\otimes gx_{2};G^{a}X_{1}X_{2},g^{d})%
\end{array}%
\right] G^{a}X_{2}\otimes g^{d}=0
\end{equation*}%
and we get%
\begin{equation*}
-B(g\otimes gx_{2};X_{2},1_{H})-B(x_{1}\otimes
gx_{2};X_{2},x_{1})-B(x_{1}\otimes gx_{2};X_{1}X_{2},1_{H})=0
\end{equation*}%
and%
\begin{equation*}
B(g\otimes gx_{2};GX_{2},g)+B(x_{1}\otimes
gx_{2};GX_{2},gx_{1})-B(x_{1}\otimes gx_{2};GX_{1}X_{2},g)=0
\end{equation*}%
which hold in view of the form of the elements.

\subsubsection{$G^{a}X_{1}\otimes g^{d}$}

We obtain%
\begin{equation*}
\sum_{\substack{ a,d  \\ a+d\equiv 0}}^{1}\left[ \left( -1\right)
^{a+1}B(g\otimes gx_{2};G^{a}X_{1},g^{d})+\left( -1\right)
^{a+1}B(x_{1}\otimes gx_{2};G^{a}X_{1},g^{d}x_{1})\right] G^{a}X_{1}\otimes
g^{d}=0
\end{equation*}%
and we get%
\begin{equation*}
-B(g\otimes gx_{2};X_{1},1_{H})-B(x_{1}\otimes gx_{2};X_{1},x_{1})=0
\end{equation*}%
and%
\begin{equation*}
B(g\otimes gx_{2};GX_{1},g)+B(x_{1}\otimes gx_{2};GX_{1},gx_{1})=0
\end{equation*}%
which hold in view of the form of the elements.

\subsubsection{$G^{a}\otimes g^{d}x_{1}x_{2}$}

We obtain%
\begin{equation*}
\sum_{\substack{ a,d  \\ a+d\equiv 1}}^{1}\left[ \left( -1\right)
^{a}B(g\otimes gx_{2};G^{a},g^{d}x_{1}x_{2})+B(x_{1}\otimes
gx_{2};G^{a}X_{1},g^{d}x_{1}x_{2})\right] G^{a}\otimes g^{d}x_{1}x_{2}=0
\end{equation*}%
and we get%
\begin{equation*}
B(g\otimes gx_{2};1_{A},gx_{1}x_{2})+B(x_{1}\otimes
gx_{2};X_{1},gx_{1}x_{2})=0
\end{equation*}%
and%
\begin{equation*}
-B(g\otimes gx_{2};G,x_{1}x_{2})+B(x_{1}\otimes gx_{2};GX_{1},x_{1}x_{2})=0
\end{equation*}%
which hold in view of the form of the elements.

\subsubsection{$G^{a}X_{2}\otimes g^{d}x_{2}$}

We obtain%
\begin{equation*}
\sum_{\substack{ a,d  \\ a+d\equiv 1}}^{1}\left[
\begin{array}{c}
\left( -1\right) ^{a+1}B(g\otimes gx_{2};G^{a}X_{2},g^{d}x_{2})+\left(
-1\right) ^{a}B(x_{1}\otimes gx_{2};G^{a}X_{2},g^{d}x_{1}x_{2})+ \\
-B(x_{1}\otimes gx_{2};G^{a}X_{1}X_{2},g^{d}x_{2})%
\end{array}%
\right] G^{a}X_{2}\otimes g^{d}x_{2}=0
\end{equation*}%
and we get%
\begin{equation*}
-B(g\otimes gx_{2};X_{2},gx_{2})+B(x_{1}\otimes
gx_{2};X_{2},gx_{1}x_{2})-B(x_{1}\otimes gx_{2};X_{1}X_{2},gx_{2})=0
\end{equation*}%
and%
\begin{equation*}
B(g\otimes gx_{2};GX_{2},x_{2})-B(x_{1}\otimes
gx_{2};GX_{2},x_{1}x_{2})-B(x_{1}\otimes gx_{2};GX_{1}X_{2},x_{2})=0
\end{equation*}%
which hold in view of the form of the elements.

\subsubsection{$G^{a}X_{1}\otimes g^{d}x_{2}$}

We obtain%
\begin{equation*}
\sum_{\substack{ a,d  \\ a+d\equiv 1}}^{1}\left[ \left( -1\right)
^{a+1}B(g\otimes gx_{2};G^{a}X_{1},g^{d}x_{2})+\left( -1\right)
^{a}B(x_{1}\otimes gx_{2};G^{a}X_{1},g^{d}x_{1}x_{2})\right]
G^{a}X_{1}\otimes g^{d}x_{2}=0
\end{equation*}%
and we get%
\begin{equation*}
-B(g\otimes gx_{2};X_{1},gx_{2})+B(x_{1}\otimes gx_{2};X_{1},gx_{1}x_{2})=0
\end{equation*}%
and%
\begin{equation*}
+B(g\otimes gx_{2};GX_{1},x_{2})-B(x_{1}\otimes gx_{2};GX_{1},x_{1}x_{2})=0
\end{equation*}%
which hold in view of the form of the elements.

\subsubsection{$G^{a}X_{2}\otimes g^{d}x_{1}$}

We obtain%
\begin{equation*}
\sum_{\substack{ a,d  \\ a+d\equiv 1}}^{1}\left[ \left( -1\right)
^{a+1}B(g\otimes gx_{2};G^{a}X_{2},g^{d}x_{1})-B(x_{1}\otimes
gx_{2};G^{a}X_{1}X_{2},g^{d}x_{1})\right] G^{a}X_{2}\otimes g^{d}x_{1}=0
\end{equation*}%
and we get%
\begin{equation*}
-B(g\otimes gx_{2};X_{2},gx_{1})-B(x_{1}\otimes gx_{2};X_{1}X_{2},gx_{1})=0.
\end{equation*}%
and%
\begin{equation*}
B(g\otimes gx_{2};GX_{2},x_{1})-B(x_{1}\otimes gx_{2};GX_{1}X_{2},x_{1})=0
\end{equation*}%
which hold in view of the form of the elements.

\subsubsection{$G^{a}X_{1}\otimes g^{d}x_{1}$}

We obtain%
\begin{equation*}
\sum_{\substack{ a,d  \\ a+d\equiv 1}}^{1}\left( -1\right) ^{a+1}B(g\otimes
gx_{2};G^{a}X_{1},g^{d}x_{1})G^{a}X_{1}\otimes g^{d}x_{1}=0
\end{equation*}%
and we get%
\begin{equation*}
-B(g\otimes gx_{2};X_{1},gx_{1})=0
\end{equation*}%
and%
\begin{equation*}
B(g\otimes gx_{2};GX_{1},x_{1})=0
\end{equation*}%
which hold in view of the form of the element.

\subsubsection{$G^{a}X_{1}X_{2}\otimes g^{d}$}

We obtain%
\begin{equation*}
\sum_{\substack{ a,d  \\ a+d\equiv 1}}^{1}\left[ \left( -1\right)
^{a}B(g\otimes gx_{2};G^{a}X_{1}X_{2},g^{d})+\left( -1\right)
^{a}B(x_{1}\otimes gx_{2};G^{a}X_{1}X_{2},g^{d}x_{1})\right]
G^{a}X_{1}X_{2}\otimes g^{d}=0
\end{equation*}%
and we get%
\begin{equation*}
B(g\otimes gx_{2};X_{1}X_{2},g)+B(x_{1}\otimes gx_{2};X_{1}X_{2},gx_{1})=0
\end{equation*}%
and%
\begin{equation*}
-B(g\otimes gx_{2};GX_{1}X_{2},1_{H})-B(x_{1}\otimes
gx_{2};GX_{1}X_{2},x_{1})=0
\end{equation*}%
which hold in view of the form of the elements.

\subsubsection{$G^{a}X_{2}\otimes g^{d}x_{1}x_{2}$}

We obtain%
\begin{equation*}
\sum_{\substack{ a,d  \\ a+d\equiv 0}}^{1}\left[ \left( -1\right)
^{a+1}B(g\otimes gx_{2};G^{a}X_{2},g^{d}x_{1}x_{2})-B(x_{1}\otimes
gx_{2};G^{a}X_{1}X_{2},g^{d}x_{1}x_{2})\right] G^{a}X_{2}\otimes
g^{d}x_{1}x_{2}=0
\end{equation*}%
and we get%
\begin{equation*}
-B(g\otimes gx_{2};X_{2},x_{1}x_{2})-B(x_{1}\otimes
gx_{2};X_{1}X_{2},x_{1}x_{2})=0
\end{equation*}%
and%
\begin{equation*}
B(g\otimes gx_{2};GX_{2},gx_{1}x_{2})-B(x_{1}\otimes
gx_{2};GX_{1}X_{2},gx_{1}x_{2})=0
\end{equation*}%
which hold in view of the form of the elements.

\subsubsection{$G^{a}X_{1}\otimes g^{d}x_{1}x_{2}$}

We obtain%
\begin{equation*}
\sum_{\substack{ a,d  \\ a+d\equiv 0}}^{1}\left( -1\right) ^{a+1}B(g\otimes
gx_{2};G^{a}X_{1},g^{d}x_{1}x_{2})G^{a}X_{1}\otimes g^{d}x_{1}x_{2}=0
\end{equation*}%
and we get%
\begin{equation*}
-B(g\otimes gx_{2};X_{1},x_{1}x_{2})=0
\end{equation*}%
and%
\begin{equation*}
B(g\otimes gx_{2};GX_{1},gx_{1}x_{2})=0
\end{equation*}%
which hold in view of the form of the element.

\subsubsection{$G^{a}X_{1}X_{2}\otimes g^{d}x_{2}$}

We obtain%
\begin{equation*}
\sum_{\substack{ a,d  \\ a+d\equiv 0}}^{1}\left[
\begin{array}{c}
\left( -1\right) ^{a}B(g\otimes gx_{2};G^{a}X_{1}X_{2},g^{d}x_{2}) \\
+\left( -1\right) ^{a+1}B(x_{1}\otimes
gx_{2};G^{a}X_{1}X_{2},g^{d}x_{1}x_{2})%
\end{array}%
\right] G^{a}X_{1}X_{2}\otimes g^{d}x_{2}=0
\end{equation*}%
and we get%
\begin{equation*}
B(g\otimes gx_{2};X_{1}X_{2},x_{2})-B(x_{1}\otimes
gx_{2};X_{1}X_{2},x_{1}x_{2})=0
\end{equation*}%
and%
\begin{equation*}
-B(g\otimes gx_{2};GX_{1}X_{2},gx_{2})+B(x_{1}\otimes
gx_{2};GX_{1}X_{2},gx_{1}x_{2})=0
\end{equation*}%
which hold in view of the form of the elements.

\subsubsection{$G^{a}X_{1}X_{2}\otimes g^{d}x_{1}$}

We obtain%
\begin{equation*}
\sum_{\substack{ a,d  \\ a+d\equiv 0}}^{1}\left( -1\right) ^{a}B(g\otimes
gx_{2};G^{a}X_{1}X_{2},g^{d}x_{1})G^{a}X_{1}X_{2}\otimes g^{d}x_{1}=0
\end{equation*}%
and we get%
\begin{equation*}
B(g\otimes gx_{2};X_{1}X_{2},x_{1})=0
\end{equation*}%
and%
\begin{equation*}
-B(g\otimes gx_{2};GX_{1}X_{2},gx_{1})=0
\end{equation*}%
which hold in view of the form of the element.

\subsubsection{$G^{a}X_{1}X_{2}\otimes g^{d}x_{1}x_{2}$}

We obtain%
\begin{equation*}
\sum_{\substack{ a,d  \\ a+d\equiv 0}}^{1}\left( -1\right) ^{a}B(g\otimes
gx_{2};G^{a}X_{1}X_{2},g^{d}x_{1})G^{a}X_{1}X_{2}\otimes g^{d}x_{1}=0
\end{equation*}%
and we get%
\begin{equation*}
B(g\otimes gx_{2};X_{1}X_{2},x_{1})=0
\end{equation*}%
and%
\begin{equation*}
-B(g\otimes gx_{2};GX_{1}X_{2},gx_{1})=0
\end{equation*}%
which hold in view of the form of the element.

\subsection{Case $gx_{2}$}

\begin{eqnarray*}
l_{1} &=&u_{1}=0 \\
l_{2}+u_{2} &=&1 \\
a+b_{1}+b_{2}+d+e_{1}+e_{2} &\equiv &0 \\
\alpha \left( 1_{H};0,0,0,1\right) &\equiv &a+b_{1}+b_{2} \\
\alpha \left( 1_{H};0,1,0,0\right) &\equiv &0
\end{eqnarray*}%
\begin{eqnarray*}
&&\sum_{\substack{ a,b_{1},b_{2},d,e_{1}=0  \\ a+b_{1}+b_{2}+d+e_{1}\equiv 1
}}^{1}\left( -1\right) ^{a+b_{1}+b_{2}}B(x_{1}\otimes
gx_{2};G^{a}X_{1}^{b_{1}}X_{2}^{b_{2}},g^{d}x_{1}^{e_{1}}x_{2})G^{a}X_{1}^{b_{1}}X_{2}^{b_{2}}\otimes g^{d}x_{1}^{e_{1}}+
\\
&&+\sum_{\substack{ a,b_{1},d,e_{1},e_{2}=0  \\ a+b_{1}+d+e_{1}+e_{2}\equiv
1 }}^{1}B(x_{1}\otimes
gx_{2};G^{a}X_{1}^{b_{1}}X_{2},g^{d}x_{1}^{e_{1}}x_{2}^{e_{2}})G^{a}X_{1}^{b_{1}}\otimes g^{d}x_{1}^{e_{1}}x_{2}^{e_{2}}
\\
&=&B^{A}(x_{1}\otimes g)\otimes B^{H}(x_{1}\otimes g)
\end{eqnarray*}

\subsubsection{$G^{a}\otimes g^{d}$}

\begin{equation*}
\sum_{\substack{ a,d,=0  \\ a+d\equiv 1}}^{1}-B(x_{1}\otimes g;G^{a},g^{d})+%
\left[ \left( -1\right) ^{a}B(x_{1}\otimes
gx_{2};G^{a},g^{d}x_{2})+B(x_{1}\otimes gx_{2};G^{a}X_{2},g^{d})\right]
G^{a}\otimes g^{d}=0
\end{equation*}%
We get%
\begin{equation*}
-B(x_{1}\otimes g;1_{A},g)+B(x_{1}\otimes
gx_{2};1_{A},gx_{2})+B(x_{1}\otimes gx_{2};X_{2},g)=0
\end{equation*}%
and%
\begin{equation*}
-B(x_{1}\otimes g;G,1_{H})-B(x_{1}\otimes gx_{2};G,x_{2})+B(x_{1}\otimes
gx_{2};GX_{2},1_{H})=0
\end{equation*}%
which hold in view of the form of the elements.

\subsubsection{$G^{a}\otimes g^{d}x_{2}$}

\begin{equation*}
\sum_{\substack{ a,d,=0  \\ a+d\equiv 0}}^{1}\left[ -B(x_{1}\otimes
g;G^{a},g^{d}x_{2})+B(x_{1}\otimes gx_{2};G^{a}X_{2},g^{d}x_{2})\right]
G^{a}\otimes g^{d}x_{2}=0
\end{equation*}%
We get%
\begin{equation*}
-B(x_{1}\otimes g;1_{A},x_{2})+B(x_{1}\otimes gx_{2};X_{2},x_{2})=0
\end{equation*}%
and%
\begin{equation*}
-B(x_{1}\otimes g;G,gx_{2})+B(x_{1}\otimes gx_{2};GX_{2},gx_{2})=0
\end{equation*}%
which hold in view of the form of the elements.$G^{a}\otimes g^{d}x_{1}$

\begin{equation*}
\sum_{\substack{ a,d,=0  \\ a+d\equiv 0}}^{1}\left[
\begin{array}{c}
-B(x_{1}\otimes g;G^{a},g^{d}x_{1})+ \\
\left( -1\right) ^{a}B(x_{1}\otimes
gx_{2};G^{a},g^{d}x_{1}x_{2})+B(x_{1}\otimes gx_{2};G^{a}X_{2},g^{d}x_{1})%
\end{array}%
\right] G^{a}\otimes g^{d}x_{1}=0
\end{equation*}%
We get%
\begin{equation*}
-B(x_{1}\otimes g;1_{A},x_{1})+B(x_{1}\otimes
gx_{2};1_{A},x_{1}x_{2})+B(x_{1}\otimes gx_{2};X_{2},x_{1})=0
\end{equation*}%
and%
\begin{equation*}
-B(x_{1}\otimes g;G,gx_{1})-B(x_{1}\otimes
gx_{2};G,gx_{1}x_{2})+B(x_{1}\otimes gx_{2};GX_{2},gx_{1})=0
\end{equation*}%
which hold in view of the form of the elements.

\subsubsection{$G^{a}X_{2}\otimes g^{d}$}

\begin{equation*}
\sum_{\substack{ a,d,=0  \\ a+d\equiv 0}}^{1}\left[
\begin{array}{c}
-B(x_{1}\otimes g;G^{a}X_{2},g^{d}) \\
+\left( -1\right) ^{a+1}B(x_{1}\otimes gx_{2};G^{a}X_{2},g^{d}x_{2})%
\end{array}%
\right] G^{a}X_{2}\otimes g^{d}=0
\end{equation*}%
We get%
\begin{equation*}
-B(x_{1}\otimes g;X_{2},1_{H})-B(x_{1}\otimes gx_{2};X_{2},x_{2})=0
\end{equation*}%
and%
\begin{equation*}
-B(x_{1}\otimes g;GX_{2},g)+B(x_{1}\otimes gx_{2};GX_{2},gx_{2})=0
\end{equation*}%
which hold in view of the form of the elements.

\subsubsection{$G^{a}X_{1}\otimes g^{d}$}

\begin{equation*}
\sum_{\substack{ a,d,=0  \\ a+d\equiv 0}}^{1}\left[
\begin{array}{c}
-B(x_{1}\otimes g;G^{a}X_{1},g^{d}) \\
+\left( -1\right) ^{a+1}B(x_{1}\otimes
gx_{2};G^{a}X_{1},g^{d}x_{2})+B(x_{1}\otimes gx_{2};G^{a}X_{1}X_{2},g^{d})%
\end{array}%
\right] G^{a}X_{1}\otimes g^{d}=0
\end{equation*}%
We get%
\begin{equation*}
-B(x_{1}\otimes g;X_{1},1_{H})-B(x_{1}\otimes
gx_{2};X_{1},x_{2})+B(x_{1}\otimes gx_{2};X_{1}X_{2},1_{H})=0
\end{equation*}%
and%
\begin{equation*}
-B(x_{1}\otimes g;GX_{1},g)+B(x_{1}\otimes
gx_{2};GX_{1},gx_{2})+B(x_{1}\otimes gx_{2};GX_{1}X_{2},g)=0
\end{equation*}%
which hold in view of the form of the elements.

\subsubsection{$G^{a}\otimes g^{d}x_{1}x_{2}$}

\begin{equation*}
\sum_{\substack{ a,d,=0  \\ a+d\equiv 1}}^{1}\left[ -B(x_{1}\otimes
g;G^{a},g^{d}x_{1}x_{2})+B(x_{1}\otimes gx_{2};G^{a}X_{2},g^{d}x_{1}x_{2})%
\right] G^{a}\otimes g^{d}x_{1}x_{2}=0
\end{equation*}%
We get%
\begin{equation*}
-B(x_{1}\otimes g;1_{A},gx_{1}x_{2})+B(x_{1}\otimes
gx_{2};X_{2},gx_{1}x_{2})=0
\end{equation*}%
and%
\begin{equation*}
-B(x_{1}\otimes g;G,x_{1}x_{2})+B(x_{1}\otimes gx_{2};GX_{2},x_{1}x_{2})=0
\end{equation*}%
which hold in view of the form of the elements.

\subsubsection{$G^{a}X_{2}\otimes g^{d}x_{2}$}

Nothing

\subsubsection{$G^{a}X_{1}\otimes g^{d}x_{2}$}

\begin{equation*}
\sum_{\substack{ a,d,=0  \\ a+d\equiv 1}}^{1}\left[ -B(x_{1}\otimes
g;G^{a}X_{1},g^{d}x_{2})+B(x_{1}\otimes gx_{2};G^{a}X_{1}X_{2},g^{d}x_{2})%
\right] G^{a}X_{1}\otimes g^{d}x_{2}=0
\end{equation*}%
We get%
\begin{equation*}
-B(x_{1}\otimes g;X_{1},gx_{2})+B(x_{1}\otimes gx_{2};X_{1}X_{2},gx_{2})=0
\end{equation*}%
and%
\begin{equation*}
-B(x_{1}\otimes g;GX_{1},x_{2})+B(x_{1}\otimes gx_{2};GX_{1}X_{2},x_{2})=0
\end{equation*}%
which hold in view of the form of the elements.

\subsubsection{$G^{a}X_{2}\otimes g^{d}x_{1}$}

\begin{equation*}
\sum_{\substack{ a,d,=0  \\ a+d\equiv 1}}^{1}\left[ -B(x_{1}\otimes
g;G^{a}X_{2}\otimes g^{d}x_{1})+\left( -1\right) ^{a+1}B(x_{1}\otimes
gx_{2};G^{a}X_{2},g^{d}x_{1}x_{2})\right] G^{a}X_{2}\otimes g^{d}x_{1}=0
\end{equation*}%
We get%
\begin{equation*}
-B(x_{1}\otimes g;X_{2},gx_{1})-B(x_{1}\otimes gx_{2};X_{2},gx_{1}x_{2})=0
\end{equation*}%
and%
\begin{equation*}
-B(x_{1}\otimes g;GX_{2},x_{1})+B(x_{1}\otimes gx_{2};GX_{2}x_{1}x_{2})=0
\end{equation*}%
which hold in view of the form of the elements.

\subsubsection{$G^{a}X_{1}\otimes g^{d}x_{1}$}

\begin{equation*}
\sum_{\substack{ a,d,=0  \\ a+d\equiv 1}}^{1}\left[
\begin{array}{c}
-B(x_{1}\otimes g;G^{a}X_{1},g^{d}x_{1})+ \\
\left( -1\right) ^{a+1}B(x_{1}\otimes
gx_{2};G^{a}X_{1},g^{d}x_{1}x_{2})+B(x_{1}\otimes
gx_{2};G^{a}X_{1}X_{2},g^{d}x_{1})%
\end{array}%
\right] G^{a}X_{1}\otimes g^{d}x_{1}=0
\end{equation*}%
and we get%
\begin{equation*}
-B(x_{1}\otimes g;X_{1},gx_{1})-B(x_{1}\otimes
gx_{2};X_{1},gx_{1}x_{2})+B(x_{1}\otimes gx_{2};X_{1}X_{2},gx_{1})=0
\end{equation*}%
and%
\begin{equation*}
-B(x_{1}\otimes g;GX_{1},x_{1})+B(x_{1}\otimes
gx_{2};GX_{1},x_{1}x_{2})+B(x_{1}\otimes gx_{2};GX_{1}X_{2},x_{1})=0
\end{equation*}%
which hold in view of the form of the elements.

\subsubsection{$G^{a}X_{1}X_{2}\otimes g^{d}$}

\begin{equation*}
\sum_{\substack{ a,d,=0  \\ a+d\equiv 1}}^{1}\left[ -B(x_{1}\otimes
g;G^{a}X_{1}X_{2},g^{d})+\left( -1\right) ^{a}B(x_{1}\otimes
gx_{2};G^{a}X_{1}X_{2},g^{d}x_{2})\right] G^{a}X_{1}X_{2}\otimes g^{d}=0
\end{equation*}%
We get%
\begin{equation*}
-B(x_{1}\otimes g;X_{1}X_{2},g)+B(x_{1}\otimes gx_{2};X_{1}X_{2},gx_{2})=0
\end{equation*}%
and%
\begin{equation*}
-B(x_{1}\otimes g;GX_{1}X_{2},1_{H})-B(x_{1}\otimes
gx_{2};GX_{1}X_{2},x_{2})=0
\end{equation*}%
which hold in view of the form of the elements.

\subsubsection{$G^{a}X_{2}\otimes g^{d}x_{1}x_{2}$}

Nothing

\subsubsection{$G^{a}X_{1}\otimes g^{d}x_{1}x_{2}$}

\begin{equation*}
\sum_{\substack{ a,d,=0  \\ a+d\equiv 0}}^{1}\left[ -B(x_{1}\otimes
g;G^{a}X_{1},g^{d}x_{1}x_{2})+B(x_{1}\otimes
gx_{2};G^{a}X_{1}X_{2},g^{d}x_{1}x_{2})\right] G^{a}X_{1}\otimes
g^{d}x_{1}x_{2}=0
\end{equation*}%
We get%
\begin{equation*}
-B(x_{1}\otimes g;X_{1},x_{1}x_{2})+B(x_{1}\otimes
gx_{2};X_{1}X_{2},x_{1}x_{2})=0
\end{equation*}%
and%
\begin{equation*}
-B(x_{1}\otimes g;GX_{1},gx_{1}x_{2})+B(x_{1}\otimes
gx_{2};GX_{1}X_{2},gx_{1}x_{2})=0
\end{equation*}%
which hold in view of the form of the elements.

\subsubsection{$G^{a}X_{1}X_{2}\otimes g^{d}x_{2}$}

Nothing

\subsubsection{$G^{a}X_{1}X_{2}\otimes g^{d}x_{1}$}

\begin{equation*}
\sum_{\substack{ a,d,=0  \\ a+d\equiv 0}}^{1}\left[ -B(x_{1}\otimes
g;G^{a}X_{1}X_{2}\otimes g^{d}x_{1})+\left( -1\right) ^{a}B(x_{1}\otimes
gx_{2};G^{a}X_{1}X_{2},g^{d}x_{1}x_{2})\right] G^{a}X_{1}X_{2}\otimes
g^{d}x_{1}=0
\end{equation*}%
We get%
\begin{equation*}
-B(x_{1}\otimes g;X_{1}X_{2},x_{1})+B(x_{1}\otimes
gx_{2};X_{1}X_{2},x_{1}x_{2})=0
\end{equation*}%
and%
\begin{equation*}
-B(x_{1}\otimes g;GX_{1}X_{2},gx_{1})-B(x_{1}\otimes
gx_{2};GX_{1}X_{2},gx_{1}x_{2})=0
\end{equation*}%
which hold in view of the form of the elements.

\subsubsection{$G^{a}X_{1}X_{2}\otimes g^{d}x_{1}x_{2}$}

Nothing.

\subsection{Case $gx_{1}x_{2}$}

\begin{eqnarray*}
&&\sum_{a,b_{1},b_{2},d,e_{1},e_{2}=0}^{1}\sum_{l_{1}=0}^{b_{1}}%
\sum_{l_{2}=0}^{b_{2}}\sum_{u_{1}=0}^{e_{1}}\sum_{u_{2}=0}^{e_{2}}\left(
-1\right) ^{\alpha \left( 1_{H};l_{1},l_{2},u_{1},u_{2}\right) } \\
&&B(x_{1}\otimes
gx_{2};G^{a}X_{1}^{b_{1}}X_{2}^{b_{2}},g^{d}x_{1}^{e_{1}}x_{2}^{e_{2}})G^{a}X_{1}^{b_{1}-l_{1}}X_{2}^{b_{2}-l_{2}}
\\
&&\otimes g^{d}x_{1}^{e_{1}-u_{1}}x_{2}^{e_{2}-u_{2}}\otimes
g^{a+b_{1}+b_{2}+l_{1}+l_{2}+d+e_{1}+e_{2}+u_{1}+u_{2}}x_{1}^{l_{1}+u_{1}}x_{2}^{l_{2}+u_{2}}
\end{eqnarray*}%
\begin{eqnarray*}
l_{1}+u_{1} &=&1 \\
l_{2}+u_{2} &=&1 \\
a+b_{1}+b_{2}+d+e_{1}+e_{2} &\equiv &1
\end{eqnarray*}%
Since%
\begin{equation*}
B(x_{1}\otimes
gx_{2};G^{a}X_{1}^{b_{1}}X_{2}^{b_{2}},g^{d}x_{1}^{e_{1}}x_{2}^{e_{2}})=0%
\text{ whenever }a+b_{1}+b_{2}+d+e_{1}+e_{2}\equiv 1
\end{equation*}%
this is trivial.

\section{$B\left( x_{1}\otimes gx_{1}x_{2}\right) $}

By using $\left( \ref{simpl1}\right) $ we have
\begin{equation*}
B(h\otimes x_{1}x_{2})=B(h\otimes x_{2})(1_{A}\otimes x_{1})-(1_{A}\otimes
gx_{1})B(h\otimes x_{2})(1_{A}\otimes g)-(1_{A}\otimes g)B(gx_{1}h\otimes
x_{2})(1_{A}\otimes g)
\end{equation*}%
so that, for $h=gx_{1}$ we obtain%
\begin{equation*}
B(gx_{1}\otimes x_{1}x_{2})=B(gx_{1}\otimes x_{2})(1_{A}\otimes
x_{1})-(1_{A}\otimes gx_{1})B(gx_{1}\otimes x_{2})(1_{A}\otimes
g)-(1_{A}\otimes g)B(gx_{1}gx_{1}\otimes x_{2})(1_{A}\otimes g)
\end{equation*}%
i.e.%
\begin{equation*}
B(gx_{1}\otimes x_{1}x_{2})=B(gx_{1}\otimes x_{2})(1_{A}\otimes
x_{1})-(1_{A}\otimes gx_{1})B(gx_{1}\otimes x_{2})(1_{A}\otimes g)
\end{equation*}%
and hence%
\begin{eqnarray*}
B(x_{1}\otimes gx_{1}x_{2}) &=&(1_{A}\otimes g)B(gx_{1}\otimes
x_{2})(1_{A}\otimes x_{1})(1_{A}\otimes g) \\
&&-(1_{A}\otimes g)(1_{A}\otimes gx_{1})B(gx_{1}\otimes x_{2})(1_{A}\otimes
g)(1_{A}\otimes g)
\end{eqnarray*}%
so that%
\begin{equation*}
B(x_{1}\otimes gx_{1}x_{2})=(1_{A}\otimes g)B(gx_{1}\otimes
x_{2})(1_{A}\otimes gx_{1})+(1_{A}\otimes x_{1})B(gx_{1}\otimes x_{2})
\end{equation*}%
i.e.%
\begin{equation*}
B(x_{1}\otimes gx_{1}x_{2})=-B(x_{1}\otimes gx_{2})(1_{A}\otimes
x_{1})+(1_{A}\otimes gx_{1})B(x_{1}\otimes gx_{2})(1_{A}\otimes g)
\end{equation*}%
Since%
\begin{equation*}
B(x_{1}\otimes gx_{2})\overset{\left( \ref{form x1otgx2}\right) }{=}%
(1_{A}\otimes g)B(gx_{1}\otimes 1_{H})(1_{A}\otimes gx_{2})+(1_{A}\otimes
x_{2})B(gx_{1}\otimes 1_{H})-B(x_{1}x_{2}\otimes 1_{H})
\end{equation*}%
\begin{eqnarray}
B(x_{1}\otimes gx_{1}x_{2}) &=&(1_{A}\otimes g)B(gx_{1}\otimes
1_{H})(1_{A}\otimes gx_{1}x_{2})  \label{form x1otgx1x2} \\
&&-(1_{A}\otimes x_{2})B(gx_{1}\otimes 1_{H})(1_{A}\otimes x_{1})  \notag \\
&&+B(x_{1}x_{2}\otimes 1_{H})(1_{A}\otimes x_{1})  \notag \\
&&+(1_{A}\otimes x_{1})B(gx_{1}\otimes 1_{H})(1_{A}\otimes x_{2})  \notag \\
&&+(1_{A}\otimes gx_{1}x_{2})B(gx_{1}\otimes 1_{H})(1_{A}\otimes g)  \notag
\\
&&-(1_{A}\otimes gx_{1})B(x_{1}x_{2}\otimes 1_{H})(1_{A}\otimes g)  \notag
\end{eqnarray}%
Finally we obtain%
\begin{eqnarray*}
B\left( x_{1}\otimes gx_{1}x_{2}\right) &=&\left[ 2B(x_{1}x_{2}\otimes
1_{H};1_{A},gx_{2})+4B(x_{1}x_{2}\otimes 1_{H};X_{2},g)\right] 1_{A}\otimes
gx_{1}x_{2}+ \\
&&-2B(x_{1}x_{2}\otimes 1_{H};G,g)G\otimes gx_{1}+ \\
&&-2B(x_{1}x_{2}\otimes 1_{H};X_{1},g)X_{1}\otimes gx_{1}+ \\
&&-2B(x_{1}x_{2}\otimes 1_{H};X_{2},g)X_{2}\otimes gx_{1}+ \\
&&+2B(x_{1}x_{2}\otimes 1_{H};X_{2},gx_{1}x_{2})X_{1}X_{2}\otimes gx_{1}x_{2}
\\
&&+\left[
\begin{array}{c}
-4B(x_{1}x_{2}\otimes 1_{H};G,gx_{1}x_{2})+ \\
4B(x_{1}x_{2}\otimes 1_{H};GX_{2},gx_{1})-2B(x_{1}x_{2}\otimes
1_{H};GX_{1},gx_{2})%
\end{array}%
\right] GX_{1}\otimes gx_{1}x_{2}+ \\
&&+2B(x_{1}x_{2}\otimes 1_{H};GX_{2},gx_{2})GX_{2}\otimes gx_{1}x_{2}+ \\
&&-2\left[
\begin{array}{c}
-B(x_{1}x_{2}\otimes 1_{H};G,gx_{1}x_{2}) \\
+B(x_{1}x_{2}\otimes 1_{H};GX_{2},gx_{1})-B(x_{1}x_{2}\otimes
1_{H};GX_{1},gx_{2})%
\end{array}%
\right] GX_{1}X_{2}\otimes gx_{1}
\end{eqnarray*}

The Casimir formula for $B\left( x_{1}\otimes gx_{1}x_{2}\right) $
\begin{eqnarray*}
&&\sum_{a,b_{1},b_{2},d,e_{1},e_{2}=0}^{1}\sum_{l_{1}=0}^{b_{1}}%
\sum_{l_{2}=0}^{b_{2}}\sum_{u_{1}=0}^{e_{1}}\sum_{u_{2}=0}^{e_{2}}\left(
-1\right) ^{\alpha \left( x_{1};l_{1},l_{2},u_{1},u_{2}\right) } \\
&&B(g\otimes
gx_{1}x_{2};G^{a}X_{1}^{b_{1}}X_{2}^{b_{2}},g^{d}x_{1}^{e_{1}}x_{2}^{e_{2}})G^{a}X_{1}^{b_{1}-l_{1}}X_{2}^{b_{2}-l_{2}}\otimes
\\
&&g^{d}x_{1}^{e_{1}-u_{1}}x_{2}^{e_{2}-u_{2}}\otimes
g^{a+b_{1}+b_{2}+l_{1}+l_{2}+d+e_{1}+e_{2}+u_{1}+u_{2}}x_{1}^{l_{1}+u_{1}+1}x_{2}^{l_{2}+u_{2}}
\\
&&\sum_{a,b_{1},b_{2},d,e_{1},e_{2}=0}^{1}\sum_{l_{1}=0}^{b_{1}}%
\sum_{l_{2}=0}^{b_{2}}\sum_{u_{1}=0}^{e_{1}}\sum_{u_{2}=0}^{e_{2}}\left(
-1\right) ^{\alpha \left( 1_{H};l_{1},l_{2},u_{1},u_{2}\right) } \\
&&B(x_{1}\otimes
gx_{1}x_{2};G^{a}X_{1}^{b_{1}}X_{2}^{b_{2}},g^{d}x_{1}^{e_{1}}x_{2}^{e_{2}})G^{a}X_{1}^{b_{1}-l_{1}}X_{2}^{b_{2}-l_{2}}\otimes
\\
&&g^{d}x_{1}^{e_{1}-u_{1}}x_{2}^{e_{2}-u_{2}}\otimes
g^{a+b_{1}+b_{2}+l_{1}+l_{2}+d+e_{1}+e_{2}+u_{1}+u_{2}}x_{1}^{l_{1}+u_{1}}x_{2}^{l_{2}+u_{2}}
\\
&=&B^{A}(x_{1}\otimes gx_{1}x_{2})\otimes B^{H}(x_{1}\otimes
gx_{1}x_{2})\otimes g \\
&&+B^{A}(x_{1}\otimes gx_{1})\otimes B^{H}(x_{1}\otimes gx_{1})\otimes x_{2}
\\
&&-B^{A}(x_{1}\otimes gx_{2})\otimes B^{H}(x_{1}\otimes gx_{2})\otimes x_{1}
\\
&&+B^{A}(x_{1}\otimes g)\otimes B^{H}(x_{1}\otimes g)\otimes gx_{1}x_{2}
\end{eqnarray*}

\subsection{$B\left( x_{1}\otimes gx_{1}x_{2};1_{A},gx_{1}x_{2}\right) $}

We get%
\begin{eqnarray*}
a &=&b_{1}=b_{2}=0 \\
d &=&e_{1}=e_{2}=1
\end{eqnarray*}%
and we obtain%
\begin{eqnarray*}
&&\sum_{u_{1}=0}^{1}\sum_{u_{2}=0}^{1}\left( -1\right) ^{\alpha \left(
1_{H};0,0,u_{1},u_{2}\right) }B(x_{1}\otimes
gx_{1}x_{2};1_{A},gx_{1}x_{2})1_{A}\otimes
gx_{1}^{1-u_{1}}x_{2}^{1-u_{2}}\otimes
g^{1+u_{1}+u_{2}}x_{1}^{u_{1}}x_{2}^{u_{2}} \\
&&=\left( -1\right) ^{\alpha \left( 1_{H};0,0,0,0\right) }B(x_{1}\otimes
gx_{1}x_{2};1_{A},gx_{1}x_{2})1_{A}\otimes gx_{1}x_{2}\otimes g+ \\
&&+\left( -1\right) ^{\alpha \left( 1_{H};0,0,0,1\right) }B(x_{1}\otimes
gx_{1}x_{2};1_{A},gx_{1}x_{2})1_{A}\otimes gx_{1}\otimes x_{2}+ \\
&&+\left( -1\right) ^{\alpha \left( 1_{H};0,0,1,0\right) }B(x_{1}\otimes
gx_{1}x_{2};1_{A},gx_{1}x_{2})1_{A}\otimes gx_{2}\otimes x_{1}+ \\
&&+\left( -1\right) ^{\alpha \left( 1_{H};0,0,1,1\right) }B(x_{1}\otimes
gx_{1}x_{2};1_{A},gx_{1}x_{2})1_{A}\otimes g\otimes gx_{1}x_{2}
\end{eqnarray*}

\subsubsection{Case $1_{A}\otimes gx_{1}x_{2}\otimes g$}

Nothing from the first summand of the left side. Second summand of the left
side gives us%
\begin{eqnarray*}
l_{1} &=&u_{1}=l_{2}=u_{2}=0 \\
a &=&b_{1}=b_{2}=0 \\
d &=&e_{1}=e_{2}=1.
\end{eqnarray*}%
By considering also the right side, we obtain a trivial equality.

\subsubsection{Case $1_{A}\otimes gx_{1}\otimes x_{2}$}

Nothing from the first summand of the left side. Second summand of the left
side gives us%
\begin{eqnarray*}
l_{1} &=&u_{1}=0 \\
l_{2}+u_{2} &=&1 \\
a &=&b_{1}=0 \\
b_{2} &=&l_{2} \\
d &=&e_{1}=1 \\
e_{2} &=&u_{2}.
\end{eqnarray*}%
Since $\alpha \left( 1_{H};0,0,0,1\right) \equiv a+b_{1}+b_{2}\equiv 0$ and $%
\alpha \left( 1_{H};0,1,0,0\right) \equiv 0$, we get%
\begin{equation*}
\left[ B(x_{1}\otimes gx_{1}x_{2};1_{A},gx_{1}x_{2})+B(x_{1}\otimes
gx_{1}x_{2};X_{2},gx_{1})\right] 1_{A}\otimes gx_{1}\otimes x_{2}.
\end{equation*}%
By considering also the right side, we get%
\begin{equation*}
B(x_{1}\otimes gx_{1}x_{2};1_{A},gx_{1}x_{2})+B(x_{1}\otimes
gx_{1}x_{2};X_{2},gx_{1})-B(x_{1}\otimes gx_{1};1_{A},gx_{1})=0.
\end{equation*}%
which holds in view of the form of the elements.

\subsubsection{Case $1_{A}\otimes gx_{2}\otimes x_{1}$}

From the first summand of the left side we get%
\begin{eqnarray*}
l_{1} &=&u_{1}=l_{2}=u_{2}=0 \\
a &=&b_{1}=b_{2}=0 \\
d &=&e_{2}=1 \\
e_{1} &=&0
\end{eqnarray*}%
Since $\alpha \left( x_{1};0,0,0,0\right) \equiv a+b_{1}+b_{2}=0$ and we get%
\begin{equation*}
B(g\otimes gx_{1}x_{2};1_{A},gx_{2})1_{A}\otimes gx_{2}\otimes x_{1}.
\end{equation*}%
For the second summand of the left side we get%
\begin{eqnarray*}
l_{1}+u_{1} &=&1,l_{2}=u_{2}=0 \\
a &=&b_{2}=0,b_{1}=l_{1} \\
d &=&e_{2}=1 \\
e_{1} &=&u_{1}
\end{eqnarray*}%
Since $\alpha \left( 1_{H};0,0,1,0\right) \equiv e_{2}+\left(
a+b_{1}+b_{2}\right) \equiv 1$ and $\alpha \left( 1_{H};1,0,0,0\right)
\equiv b_{2}=0,$ we get

\begin{equation*}
\left[ -B(x_{1}\otimes gx_{1}x_{2};1_{A},gx_{1}x_{2})+B(x_{1}\otimes
gx_{1}x_{2};X_{1},gx_{2})\right] 1_{A}\otimes gx_{2}\otimes x_{1}.
\end{equation*}%
By considering also the right side of the equation, we obtain%
\begin{equation*}
+B(g\otimes gx_{1}x_{2};1_{A},gx_{2})-B(x_{1}\otimes
gx_{1}x_{2};1_{A},gx_{1}x_{2})+B(x_{1}\otimes
gx_{1}x_{2};X_{1},gx_{2})+B(x_{1}\otimes gx_{2};1_{A},gx_{2})=0
\end{equation*}

which holds in view of the form of the elements.

\subsubsection{Case $1_{A}\otimes g\otimes gx_{1}x_{2}$}

From the first summand of the left side we get%
\begin{eqnarray*}
l_{1} &=&u_{1}=0 \\
l_{2}+u_{2} &=&1 \\
a &=&b_{1}=0 \\
b_{2} &=&l_{2} \\
d &=&1,e_{1}=0 \\
e_{2} &=&u_{2}
\end{eqnarray*}%
Since $\alpha \left( x_{1};0,0,0,1\right) \equiv 1$ and $\alpha \left(
x_{1};0,1,0,0\right) \equiv 0$, we get

\begin{equation*}
\left[ -B(g\otimes gx_{1}x_{2};1_{A},gx_{2})+B(g\otimes gx_{1}x_{2};X_{2},g)%
\right] 1_{A}\otimes g\otimes gx_{1}x_{2}.
\end{equation*}%
From the second summand of the left side we get%
\begin{eqnarray*}
l_{1}+u_{1} &=&1 \\
l_{2}+u_{2} &=&1 \\
a &=&0 \\
b_{1} &=&l_{1},b_{2}=l_{2} \\
e_{1} &=&u_{1},e_{2}=u_{2} \\
d &=&1.
\end{eqnarray*}%
Since $\alpha \left( 1_{H};0,0,1,1\right) \equiv 1+e_{2}\equiv 0$, $\alpha
\left( 1_{H};0,1,1,0\right) \equiv e_{2}+a+b_{1}+b_{2}+1\equiv 0$, $\alpha
\left( 1_{H};1,0,0,1\right) \equiv 1$ and $\alpha \left(
1_{H};1,1,0,0\right) \equiv 0,$ we obtain%
\begin{equation*}
\left[
\begin{array}{c}
B(x_{1}\otimes gx_{1}x_{2};1_{A},gx_{1}x_{2})+B(x_{1}\otimes
gx_{1}x_{2};X_{2},gx_{1})+ \\
-B(x_{1}\otimes gx_{1}x_{2};X_{1},gx_{2})+B(x_{1}\otimes
gx_{1}x_{2};X_{1}X_{2},g)%
\end{array}%
\right] 1_{A}\otimes g\otimes gx_{1}x_{2}.
\end{equation*}%
By considering also the right side, we obtain%
\begin{equation*}
\begin{array}{c}
-B(g\otimes gx_{1}x_{2};1_{A},gx_{2})+B(g\otimes gx_{1}x_{2};X_{2},g) \\
B(x_{1}\otimes gx_{1}x_{2};1_{A},gx_{1}x_{2})+B(x_{1}\otimes
gx_{1}x_{2};X_{2},gx_{1})+ \\
-B(x_{1}\otimes gx_{1}x_{2};X_{1},gx_{2})+B(x_{1}\otimes
gx_{1}x_{2};X_{1}X_{2},g)+ \\
-B(x_{1}\otimes g;1_{A},g)=0.%
\end{array}%
\end{equation*}%
which holds in view of the form of the elements.

\subsection{$B\left( x_{1}\otimes gx_{1}x_{2};G,gx_{1}\right) $}

We get%
\begin{eqnarray*}
a &=&1,b_{1}=b_{2}=0 \\
d &=&e_{1}=1,e_{2}=0
\end{eqnarray*}

and we obtain%
\begin{eqnarray*}
&&\left( -1\right) ^{\alpha \left( 1_{H};0,0,0,0\right) }B(x_{1}\otimes
gx_{1}x_{2};G,gx_{1})G\otimes gx_{1}\otimes g+ \\
&&\left( -1\right) ^{\alpha \left( 1_{H};0,0,1,0\right) }B(x_{1}\otimes
gx_{1}x_{2};G,gx_{1})G\otimes g\otimes x_{1}
\end{eqnarray*}

\subsubsection{Case $G\otimes gx_{1}\otimes g$}

Nothing from the first summand of the left side. Second summand of the left
side gives us%
\begin{eqnarray*}
l_{1} &=&u_{1}=l_{2}=u_{2}=0 \\
a &=&1,b_{1}=b_{2}=0 \\
d &=&e_{1}=1,e_{2}=0.
\end{eqnarray*}%
Since $\alpha \left( 1_{H};0,0,0,0\right) =0,$ we get%
\begin{equation*}
B(x_{1}\otimes gx_{1}x_{2};G,gx_{1})G\otimes gx_{1}\otimes g.
\end{equation*}%
By considering also the right side, we obtain a trivial equality.

\subsubsection{Case $G\otimes g\otimes x_{1}$}

From the first summand of the left side we get%
\begin{eqnarray*}
l_{1} &=&u_{1}=l_{2}=u_{2}=0 \\
a &=&1,b_{1}=b_{2}=0 \\
d &=&1,e_{1}=e_{2}=0.
\end{eqnarray*}%
Since $\alpha \left( x_{1};0,0,0,0\right) \equiv a+b_{1}+b_{2}\equiv 1,$ we
get%
\begin{equation*}
-B(g\otimes gx_{1}x_{2};G,g)G\otimes g\otimes x_{1}.
\end{equation*}%
By considering the second summand of the left side, we get%
\begin{eqnarray*}
l_{1}+u_{1} &=&1,l_{2}=u_{2}=0 \\
a &=&1,b_{1}=l_{1},b_{2}=0 \\
d &=&1,e_{1}=u_{1},e_{2}=0.
\end{eqnarray*}%
Since $\alpha \left( 1_{H};0,0,1,0\right) \equiv e_{2}+\left(
a+b_{1}+b_{2}\right) \equiv 1$ and $\alpha \left( 1_{H};1,0,0,0\right)
\equiv b_{2}=0,$ we obtain%
\begin{equation*}
\left[ -B(x_{1}\otimes gx_{1}x_{2};G,gx_{1})+B(x_{1}\otimes
gx_{1}x_{2};GX_{1},g)\right] G\otimes g\otimes x_{1}.
\end{equation*}%
By considering also the right side, we get%
\begin{gather*}
-B(g\otimes gx_{1}x_{2};G,g)-B(x_{1}\otimes gx_{1}x_{2};G,gx_{1}) \\
+B(x_{1}\otimes gx_{1}x_{2};GX_{1},g)+B(x_{1}\otimes gx_{2};G,g)=0
\end{gather*}%
which holds in view of the form of the elements.

\subsection{$B\left( x_{1}\otimes gx_{1}x_{2};X_{1},gx_{1}\right) $}

We get%
\begin{eqnarray*}
b_{1} &=&1,a=b_{2}=0 \\
d &=&e_{1}=1,e_{2}=0
\end{eqnarray*}

and we obtain%
\begin{gather*}
\left( -1\right) ^{\alpha \left( 1_{H};0,0,0,0\right) }B(x_{1}\otimes
gx_{1}x_{2};X_{1},gx_{1})X_{1}\otimes gx_{1}\otimes g+ \\
+\left( -1\right) ^{\alpha \left( 1_{H};0,0,1,0\right) }B(x_{1}\otimes
gx_{1}x_{2};X_{1},gx_{1})X_{1}\otimes g\otimes x_{1}+ \\
+\left( -1\right) ^{\alpha \left( 1_{H};1,0,0,0\right) }B(x_{1}\otimes
gx_{1}x_{2};X_{1},gx_{1})1_{A}\otimes gx_{1}\otimes x_{1}+ \\
+\left( -1\right) ^{\alpha \left( 1_{H};1,0,1,0\right) }B(x_{1}\otimes
gx_{1}x_{2};X_{1},gx_{1})1_{A}\otimes g\otimes gx_{1}^{2}+=0
\end{gather*}

\subsubsection{Case $X_{1}\otimes gx_{1}\otimes g$}

We get a trivial equality.

\subsubsection{Case $X_{1}\otimes g\otimes x_{1}$}

From the first summand of the left side we get%
\begin{eqnarray*}
l_{1} &=&u_{1}=l_{2}=u_{2}=0 \\
a &=&b_{2}=0,b_{1}=1 \\
d &=&1,e_{1}=e_{2}=0
\end{eqnarray*}%
Since $\alpha \left( x_{1};0,0,0,0\right) \equiv a+b_{1}+b_{2}=1,$ we get%
\begin{equation*}
-B(g\otimes gx_{1}x_{2};X_{1},g)X_{1}\otimes g\otimes x_{1}.
\end{equation*}%
By considering the second summand of the left side, we get%
\begin{eqnarray*}
l_{1}+u_{1} &=&1,l_{2}=u_{2}=0 \\
a &=&0,b_{1}-l_{1}=1\Rightarrow b_{1}=1,l_{1}=0,u_{1}=1,b_{2}=0 \\
d &=&1,e_{1}=u_{1}=1,e_{2}=0.
\end{eqnarray*}%
Since $\alpha \left( 1_{H};0,0,1,0\right) \equiv e_{2}+\left(
a+b_{1}+b_{2}\right) \equiv -1$ we obtain%
\begin{equation*}
-B(x_{1}\otimes gx_{1}x_{2};X_{1},gx_{1})X_{1}\otimes g\otimes x_{1}.
\end{equation*}%
By considering also the right side, we get%
\begin{equation*}
-B(g\otimes gx_{1}x_{2};X_{1},g)-B(x_{1}\otimes
gx_{1}x_{2};X_{1},gx_{1})+B(x_{1}\otimes gx_{2};X_{1},g)=0
\end{equation*}%
which hold in view of the form of the elements.

\subsubsection{Case $1_{A}\otimes gx_{1}\otimes x_{1}$}

From the first summand of the left side we get%
\begin{eqnarray*}
l_{1} &=&u_{1}=l_{2}=u_{2}=0 \\
a &=&b_{2}=b_{1}=0 \\
d &=&e_{1}=1,e_{2}=0
\end{eqnarray*}%
Since $\alpha \left( x_{1};0,0,0,0\right) \equiv a+b_{1}+b_{2}=0,$ we get%
\begin{equation*}
B(g\otimes gx_{1}x_{2};1_{A},gx_{1})1_{A}\otimes gx_{1}\otimes x_{1}.
\end{equation*}%
By considering the second summand of the left side, we get%
\begin{eqnarray*}
l_{1}+u_{1} &=&1,l_{2}=u_{2}=0 \\
a &=&b_{2}=0,b_{1}=l_{1}, \\
d &=&1,e_{1}-u_{1}=1\Rightarrow e_{1}=1,u_{1}=0,l_{1}=b_{1}=1,e_{2}=0.
\end{eqnarray*}%
Since $\alpha \left( 1_{H};1,0,0,0\right) =b_{2}=0,$ we obtain%
\begin{equation*}
B(x_{1}\otimes gx_{1}x_{2};X_{1},gx_{1})1_{A}\otimes gx_{1}\otimes x_{1}.
\end{equation*}%
By considering also the right side we get%
\begin{equation*}
B(g\otimes gx_{1}x_{2};1_{A},gx_{1})+B(x_{1}\otimes
gx_{1}x_{2};X_{1},gx_{1})+B(x_{1}\otimes gx_{2};1_{A},gx_{1})=0
\end{equation*}

which holds in view of the form of the element.

\subsection{$B\left( x_{1}\otimes gx_{1}x_{2};X_{2},gx_{1}\right) $}

We get%
\begin{eqnarray*}
b_{2} &=&1,a=b_{1}=0 \\
d &=&e_{1}=1,e_{2}=0
\end{eqnarray*}

and we obtain%
\begin{eqnarray*}
&&\left( -1\right) ^{\alpha \left( 1_{H};0,0,0,0\right) }B(x_{1}\otimes
gx_{1}x_{2};X_{2},gx_{1})X_{2}\otimes gx_{1}\otimes g+ \\
&&\left( -1\right) ^{\alpha \left( 1_{H};0,0,1,0\right) }B(x_{1}\otimes
gx_{1}x_{2};X_{2},gx_{1})X_{2}\otimes g\otimes x_{1} \\
&&\left( -1\right) ^{\alpha \left( 1_{H};0,1,0,0\right) }B(x_{1}\otimes
gx_{1}x_{2};X_{2},gx_{1})1_{A}\otimes gx_{1}\otimes x_{2} \\
&&\left( -1\right) ^{\alpha \left( 1_{H};0,1,1,0\right) }B(x_{1}\otimes
gx_{1}x_{2};X_{2},gx_{1})1_{A}\otimes g\otimes gx_{1}x_{2}
\end{eqnarray*}

\subsubsection{Case $X_{2}\otimes gx_{1}\otimes g$}

Nothing new.

\subsubsection{Case $X_{2}\otimes g\otimes x_{1}$}

From the first summand of the left side we get%
\begin{eqnarray*}
l_{1} &=&u_{1}=l_{2}=u_{2}=0 \\
a &=&b_{1}=0,b_{2}=1 \\
d &=&1,e_{1}=e_{2}=0
\end{eqnarray*}%
Since $\alpha \left( x_{1};0,0,0,0\right) \equiv a+b_{1}+b_{2}=1,$ we get%
\begin{equation*}
-B(g\otimes gx_{1}x_{2};X_{2},g)X_{2}\otimes g\otimes x_{1}.
\end{equation*}%
By considering the second summand of the left side, we get%
\begin{eqnarray*}
l_{1}+u_{1} &=&1,l_{2}=u_{2}=0 \\
a &=&0,b_{1}=l_{1},b_{2}=1 \\
d &=&1,e_{1}=u_{1},e_{2}=0.
\end{eqnarray*}%
Since $\alpha \left( 1_{H};0,0,1,0\right) =e_{2}+\left( a+b_{1}+b_{2}\right)
=1$ and $\alpha \left( 1_{H};1,0,0,0\right) =b_{2}=1$ we obtain%
\begin{equation*}
\left[ -B(x_{1}\otimes gx_{1}x_{2};X_{2},gx_{1})-B(x_{1}\otimes
gx_{1}x_{2};X_{1}X_{2},g)\right] X_{2}\otimes g\otimes x_{1}.
\end{equation*}%
By considering also the right side, we obtain%
\begin{gather*}
-B(g\otimes gx_{1}x_{2};X_{2},g)-B(x_{1}\otimes gx_{1}x_{2};X_{2},gx_{1}) \\
-B(x_{1}\otimes gx_{1}x_{2};X_{1}X_{2},g)+B(x_{1}\otimes gx_{2};X_{2},g)=0
\end{gather*}%
which holds in view of the form of the elements.

\subsubsection{Case $1_{A}\otimes gx_{1}\otimes x_{2}$}

We already considered this case in $B\left( x_{1}\otimes
gx_{1}x_{2};1_{A},gx_{1}x_{2}\right) .$

\subsubsection{Case $1_{A}\otimes g\otimes gx_{1}x_{2}$}

We already considered this case in $B\left( x_{1}\otimes
gx_{1}x_{2};1_{A},gx_{1}x_{2}\right) .$

\subsection{$B\left( x_{1}\otimes gx_{1}x_{2};X_{1}X_{2},gx_{1}x_{2}\right) $%
}

We have%
\begin{equation*}
a=0,b_{1}=b_{2}=1,d=e_{1}=e_{2}=1
\end{equation*}%
\begin{equation*}
X_{1}^{1-l_{1}}X_{2}^{1-l_{2}}\otimes gx_{1}^{1-u_{1}}x_{2}^{1-u}\otimes
\otimes g^{1+l_{1}+l_{2}+u_{1}+u_{2}}x_{1}^{l_{1}+u_{1}}x_{2}^{l_{2}+u_{2}}
\end{equation*}%
\begin{gather*}
\left( -1\right) ^{\alpha \left( 1_{H};0,0,1,0\right) }B(x_{1}\otimes
gx_{1}x_{2};X_{1}X_{2},gx_{1}x_{2})X_{1}X_{2}\otimes gx_{2}\otimes x_{1}+ \\
\left( -1\right) ^{\alpha \left( 1_{H};0,1,1,0\right) }B(x_{1}\otimes
gx_{1}x_{2};X_{1}X_{2},gx_{1}x_{2})X_{1}\otimes gx_{2}\otimes gx_{1}x_{2} \\
+\left( -1\right) ^{\alpha \left( 1_{H};1,0,1,0\right) }B(x_{1}\otimes
gx_{1}x_{2};X_{1}X_{2},gx_{1}x_{2})X_{1}^{1-l_{1}}X_{2}^{1-l_{2}}\otimes
gx_{2}\otimes g^{+l_{1}+l_{2}}x_{1}^{1+1}x_{2}^{l_{2}}=0 \\
+\left( -1\right) ^{\alpha \left( 1_{H};1,1,1,0\right) }B(x_{1}\otimes
gx_{1}x_{2};X_{1}X_{2},gx_{1}x_{2})X_{1}^{1-l_{1}}X_{2}^{1-l_{2}}\otimes
gx_{2}\otimes g^{+l_{1}+l_{2}}x_{1}^{1+1}x_{2}^{l_{2}}=0 \\
+\left( -1\right) ^{\alpha \left( 1_{H};0,0,0,0\right) }B(x_{1}\otimes
gx_{1}x_{2};X_{1}X_{2},gx_{1}x_{2})X_{1}X_{2}\otimes gx_{1}x_{2}\otimes g+ \\
+\left( -1\right) ^{\alpha \left( 1_{H};0,1,0,0\right) }B(x_{1}\otimes
gx_{1}x_{2};X_{1}X_{2},gx_{1}x_{2})X_{1}\otimes gx_{1}x_{2}\otimes x_{2}+ \\
+\left( -1\right) ^{\alpha \left( 1_{H};1,0,0,0\right) }B(x_{1}\otimes
gx_{1}x_{2};X_{1}X_{2},gx_{1}x_{2})X_{2}\otimes gx_{1}x_{2}\otimes x_{1}+ \\
+\left( -1\right) ^{\alpha \left( 1_{H};1,1,0,0\right) }B(x_{1}\otimes
gx_{1}x_{2};X_{1}X_{2},gx_{1}x_{2})1_{A}\otimes gx_{1}x_{2}\otimes
gx_{1}x_{2}+ \\
+\left( -1\right) ^{\alpha \left( 1_{H};0,0,0,1\right) }B(x_{1}\otimes
gx_{1}x_{2};X_{1}X_{2},gx_{1}x_{2})X_{1}X_{2}\otimes gx_{1}\otimes x_{2}+ \\
+\left( -1\right) ^{\alpha \left( 1_{H};0,1,0,1\right) }B(x_{1}\otimes
gx_{1}x_{2};X_{1}X_{2},gx_{1}x_{2})X_{1}^{1-l_{1}}X_{2}^{1-l_{2}}\otimes
gx_{1}\otimes g^{1+l_{1}+l_{2}+1}x_{1}^{l_{1}}x_{2}^{1+1}=0 \\
+\left( -1\right) ^{\alpha \left( 1_{H};1,0,0,1\right) }B(x_{1}\otimes
gx_{1}x_{2};X_{1}X_{2},gx_{1}x_{2})X_{2}\otimes gx_{1}\otimes gx_{1}x_{2}+ \\
+\left( -1\right) ^{\alpha \left( 1_{H};1,1,0,1\right) }B(x_{1}\otimes
gx_{1}x_{2};X_{1}X_{2},gx_{1}x_{2})X_{1}^{1-l_{1}}X_{2}^{1-l_{2}}\otimes
gx_{1}\otimes g^{1+l_{1}+l_{2}+1}x_{1}^{l_{1}}x_{2}^{1+1}=0 \\
+\left( -1\right) ^{\alpha \left( 1_{H};0,0,1,1\right) }B(x_{1}\otimes
gx_{1}x_{2};X_{1}X_{2},gx_{1}x_{2})X_{1}X_{2}\otimes g\otimes gx_{1}x_{2}+ \\
+\left( -1\right) ^{\alpha \left( 1_{H};0,1,1,1\right) }B(x_{1}\otimes
gx_{1}x_{2};X_{1}X_{2},gx_{1}x_{2})X_{2}^{1-l_{2}}\otimes g\otimes
g^{1+l_{1}+l_{2}}x_{1}^{l_{1}+1}x_{2}^{1+1}=0 \\
+\left( -1\right) ^{\alpha \left( 1_{H};1,0,1,1\right) }B(x_{1}\otimes
gx_{1}x_{2};X_{1}X_{2},gx_{1}x_{2})X_{2}^{1-l_{2}}\otimes g\otimes
g^{1+l_{1}+l_{2}}x_{1}^{1+1}x_{2}^{l_{2}+1}=0 \\
+\left( -1\right) ^{\alpha \left( 1_{H};1,1,1,1\right) }B(x_{1}\otimes
gx_{1}x_{2};X_{1}X_{2},gx_{1}x_{2})X_{2}^{1-l_{2}}\otimes g\otimes
g^{1+l_{1}+l_{2}}x_{1}^{1+1}x_{2}^{1+1}=0
\end{gather*}

\subsubsection{Case $X_{1}X_{2}\otimes gx_{2}\otimes x_{1}$}

From the first summand of the left side we get%
\begin{eqnarray*}
l_{1} &=&u_{1}=l_{2}=u_{2}=0 \\
a &=&0,b_{1}=b_{2}=1 \\
d &=&e_{2}=1,e_{1}=0.
\end{eqnarray*}%
Since $\alpha \left( x_{1};0,0,0,0\right) \equiv a+b_{1}+b_{2}=1,$ we get%
\begin{equation*}
B(g\otimes gx_{1}x_{2};X_{1}X_{2},gx_{2})X_{1}X_{2}\otimes gx_{2}\otimes
x_{1}.
\end{equation*}%
By considering the second summand of the left side, we get%
\begin{eqnarray*}
l_{1}+u_{1} &=&1,l_{2}=u_{2}=0 \\
a &=&0,b_{1}=1,l_{1}=0\Rightarrow u_{1}=1,b_{2}=1 \\
d &=&1,e_{1}=u_{1}=1,e_{2}=1.
\end{eqnarray*}%
Since $\alpha \left( 1_{H};0,0,1,0\right) =e_{2}+\left( a+b_{1}+b_{2}\right)
\equiv 1,$ we get%
\begin{equation*}
-B(x_{1}\otimes gx_{1}x_{2};X_{1}X_{2},gx_{1}x_{2})X_{1}X_{2}\otimes
gx_{2}\otimes x_{1}.
\end{equation*}%
By considering also the right side, we get%
\begin{equation*}
B(g\otimes gx_{1}x_{2};X_{1}X_{2},gx_{2})-B(x_{1}\otimes
gx_{1}x_{2};X_{1}X_{2},gx_{1}x_{2})+B(x_{1}\otimes
gx_{2};X_{1}X_{2},gx_{2})=0
\end{equation*}%
which holds in view of the form of the elements.

\subsubsection{Case $X_{1}\otimes gx_{2}\otimes gx_{1}x_{2}$}

From the first summand of the left side we get%
\begin{eqnarray*}
l_{1} &=&u_{1}=0,l_{2}+u_{2}=1 \\
a &=&0,b_{1}=1,b_{2}=l_{2} \\
d &=&1,e_{1}=0,e_{2}=1,u_{2}=0\Rightarrow b_{2}=l_{2}=1.
\end{eqnarray*}

Since $\alpha \left( x_{1};0,1,0,0\right) =a+b_{1}+b_{2}+1\equiv 1,$ we get%
\begin{equation*}
-B(g\otimes gx_{1}x_{2};X_{1}X_{2},gx_{2})X_{1}\otimes gx_{2}\otimes
gx_{1}x_{2}.
\end{equation*}%
By considering the second summand of the left side, we get%
\begin{eqnarray*}
l_{1}+u_{1} &=&1,l_{2}+u_{2}=1 \\
a &=&0,b_{1}=1,l_{1}=0\Rightarrow u_{1}=1,b_{2}=l_{2} \\
d &=&1,e_{1}=u_{1}=1,e_{2}=1,u_{2}=0\Rightarrow b_{2}=l_{2}=1.
\end{eqnarray*}%
Since $\alpha \left( 1_{H};0,1,1,0\right) =e_{2}+a+b_{1}+b_{2}+1\equiv 0,$
we obtain%
\begin{equation*}
B(x_{1}\otimes gx_{1}x_{2};X_{1}X_{2},gx_{1}x_{2})X_{1}\otimes gx_{2}\otimes
gx_{1}x_{2}.
\end{equation*}%
By considering also the right side, we get%
\begin{equation*}
-B(g\otimes gx_{1}x_{2};X_{1}X_{2},gx_{2})+B(x_{1}\otimes
gx_{1}x_{2};X_{1}X_{2},gx_{1}x_{2})-B(x_{1}\otimes g;X_{1},gx_{2})=0
\end{equation*}%
which holds in view of the form of the element.

\subsubsection{Case $X_{1}X_{2}\otimes gx_{1}x_{2}\otimes g$}

We get a trivial equality.

\subsubsection{Case $X_{1}\otimes gx_{1}x_{2}\otimes x_{2}$}

Nothing from the first summand of the left side. From the second summand of
the left side we get%
\begin{eqnarray*}
l_{1} &=&u_{1}=0,l_{2}+u_{2}=1 \\
a &=&0,b_{1}=1,b_{2}=l_{2} \\
d &=&1,e_{1}=1,e_{2}=1,u_{2}=0\Rightarrow b_{2}=l_{2}=1.
\end{eqnarray*}%
Since $\alpha \left( 1_{H},0,1,0,0\right) =0,$ we get%
\begin{equation*}
B(x_{1}\otimes gx_{1}x_{2};X_{1}X_{2},gx_{1}x_{2})X_{1}\otimes
gx_{1}x_{2}\otimes x_{2}.
\end{equation*}%
By considering also the right side we get%
\begin{equation*}
B(x_{1}\otimes gx_{1}x_{2};X_{1}X_{2},gx_{1}x_{2})-B(x_{1}\otimes
gx_{1};X_{1},gx_{1}x_{2})=0
\end{equation*}%
which holds in view of the form of the elements.

\subsubsection{Case $X_{2}\otimes gx_{1}x_{2}\otimes x_{1}$}

From the first summand of the left side we get%
\begin{eqnarray*}
l_{1} &=&u_{1}=l_{2}=u_{2}=0, \\
a &=&0,b_{1}=0,b_{2}=1 \\
d &=&1,e_{1}=1,e_{2}=1.
\end{eqnarray*}%
Since $\alpha \left( x_{1};0,0,0,0\right) =1$ we get%
\begin{equation*}
-B(g\otimes gx_{1}x_{2};X_{2},gx_{1}x_{2})X_{2}\otimes gx_{1}x_{2}\otimes
x_{1}
\end{equation*}%
From the second summand of the left side we get%
\begin{eqnarray*}
l_{1}+u_{1} &=&1,l_{2}=u_{2}=0 \\
a &=&0,b_{1}=l_{1},b_{2}=1 \\
d &=&1,e_{1}-u_{1}=1\Rightarrow e_{1}=1,u_{1}=0,b_{1}=l_{1}=1,e_{2}=1.
\end{eqnarray*}%
Since $\alpha \left( 1_{H};1,0,0,0\right) =b_{2}=1$ we get

\begin{equation*}
-B(x_{1}\otimes gx_{1}x_{2};X_{1}X_{2},gx_{1}x_{2})X_{2}\otimes
gx_{1}x_{2}\otimes x_{1}.
\end{equation*}

By considering also the right side we get%
\begin{equation*}
-B(g\otimes gx_{1}x_{2};X_{2},gx_{1}x_{2})-B(x_{1}\otimes
gx_{1}x_{2};X_{1}X_{2},gx_{1}x_{2})+B(x_{1}\otimes
gx_{2};X_{2},gx_{1}x_{2})=0
\end{equation*}%
which holds in view of the form of the elements.

\subsubsection{Case $1_{A}\otimes gx_{1}x_{2}\otimes gx_{1}x_{2}$}

From the first summand of the left side we get%
\begin{eqnarray*}
l_{1} &=&u_{1}=0,l_{2}+u_{2}=1, \\
a &=&0,b_{1}=0,b_{2}=l_{2} \\
d &=&1,e_{1}=1,e_{2}-u_{2}=1\Rightarrow e_{2}=1,u_{2}=0,b_{2}=l_{2}=1.
\end{eqnarray*}%
Since $\alpha \left( x_{1};0,1,0,0\right) =a+b_{1}+b_{2}+1\equiv 0,$ we get%
\begin{equation*}
B(g\otimes gx_{1}x_{2};X_{2},gx_{1}x_{2})1_{A}\otimes gx_{1}x_{2}\otimes
gx_{1}x_{2}.
\end{equation*}%
From the second summand of the left side we get%
\begin{eqnarray*}
l_{1}+u_{1} &=&1,l_{2}+u_{2}=1 \\
a &=&0,b_{1}=l_{1},b_{2}=l_{2} \\
d &=&1,e_{1}-u_{1}=1\Rightarrow e_{1}=1,u_{1}=0,b_{1}=l_{1}=1 \\
e_{2}-u_{2} &=&1.\Rightarrow e_{2}=1,u_{2}=0,b_{2}=l_{2}=1
\end{eqnarray*}%
Since $\alpha \left( 1_{H};1,1,0,0\right) \equiv 1+b_{2}\equiv 0,$ we get%
\begin{equation*}
B(x_{1}\otimes gx_{1}x_{2};X_{1}X_{2},gx_{1}x_{2})1_{A}\otimes
gx_{1}x_{2}\otimes gx_{1}x_{2}
\end{equation*}%
By considering also the right side we get%
\begin{equation*}
B(g\otimes gx_{1}x_{2};X_{2},gx_{1}x_{2})+B(x_{1}\otimes
gx_{1}x_{2};X_{1}X_{2},gx_{1}x_{2})-B(x_{1}\otimes g;1_{A},gx_{1}x_{2})=0
\end{equation*}%
which holds in view of the form of the elements.

\subsubsection{Case $X_{1}X_{2}\otimes gx_{1}\otimes x_{2}$}

Nothing from the first term of the left side. The second term of the left
side gives us

\begin{eqnarray*}
l_{1} &=&u_{1}=0,l_{2}+u_{2}=1 \\
a &=&0,b_{1}=1,b_{2}-l_{2}=1\Rightarrow b_{2}=1,l_{2}=0,u_{2}=1 \\
d &=&1,e_{1}=1,e_{2}=u_{2}=1.
\end{eqnarray*}%
Since $\alpha \left( 1_{H};0,0,0,1\right) =a+b_{1}+b_{2}\equiv 0,$ we obtain%
\begin{equation*}
B(x_{1}\otimes gx_{1}x_{2};X_{1}X_{2},gx_{1}x_{2})X_{1}X_{2}\otimes
gx_{1}\otimes x_{2}.
\end{equation*}%
By considering also the right side we get%
\begin{equation*}
B(x_{1}\otimes gx_{1}x_{2};X_{1}X_{2},gx_{1}x_{2})-B(x_{1}\otimes
gx_{1};X_{1}X_{2},gx_{1})=0
\end{equation*}%
which holds in view of the form of the elements.

\subsubsection{Case $X_{2}\otimes gx_{1}\otimes gx_{1}x_{2}$}

From the first summand of the left side we get%
\begin{eqnarray*}
&&\sum_{a,b_{1},b_{2},d,e_{1},e_{2}=0}^{1}\sum_{l_{1}=0}^{b_{1}}%
\sum_{l_{2}=0}^{b_{2}}\sum_{u_{1}=0}^{e_{1}}\sum_{u_{2}=0}^{e_{2}}\left(
-1\right) ^{\alpha \left( x_{1};l_{1},l_{2},u_{1},u_{2}\right) } \\
&&B(g\otimes
gx_{1}x_{2};G^{a}X_{1}^{b_{1}}X_{2}^{b_{2}},g^{d}x_{1}^{e_{1}}x_{2}^{e_{2}})G^{a}X_{1}^{b_{1}-l_{1}}X_{2}^{b_{2}-l_{2}}\otimes
\\
&&g^{d}x_{1}^{e_{1}-u_{1}}x_{2}^{e_{2}-u_{2}}\otimes
g^{a+b_{1}+b_{2}+l_{1}+l_{2}+d+e_{1}+e_{2}+u_{1}+u_{2}}x_{1}^{l_{1}+u_{1}+1}x_{2}^{l_{2}+u_{2}}
\end{eqnarray*}%
\begin{eqnarray*}
l_{1} &=&u_{1}=0,l_{2}+u_{2}=1, \\
a &=&0,b_{1}=0,b_{2}-l_{2}=1\Rightarrow b_{2}=1,l_{2}=0,u_{2}=1. \\
d &=&1,e_{1}=1,e_{2}=u_{2}=1.
\end{eqnarray*}%
Since $\alpha \left( x_{1};0,0,0,1\right) =1,$ we obtain%
\begin{equation*}
-B(g\otimes gx_{1}x_{2};X_{2},gx_{1}x_{2})X_{2}\otimes gx_{1}\otimes
gx_{1}x_{2}.
\end{equation*}%
From the second summand of the left side we get%
\begin{eqnarray*}
l_{1}+u_{1} &=&1,l_{2}+u_{2}=1 \\
a &=&0,b_{1}=l_{1},b_{2}-l_{2}=1\Rightarrow b_{2}=1,l_{2}=0,u_{2}=1 \\
d &=&1,e_{1}-u_{1}=1\Rightarrow e_{1}=1,u_{1}=0,b_{1}=l_{1}=1,e_{2}=u_{2}=1
\end{eqnarray*}%
Since $\alpha \left( 1_{H};1,0,0,1\right) \equiv a+b_{1}=1,$ we get%
\begin{equation*}
-B(x_{1}\otimes gx_{1}x_{2};X_{1}X_{2},gx_{1}x_{2})X_{2}\otimes
gx_{1}\otimes gx_{1}x_{2}.
\end{equation*}%
By considering also the right side, we get%
\begin{equation*}
-B(g\otimes gx_{1}x_{2};X_{2},gx_{1}x_{2})-B(x_{1}\otimes
gx_{1}x_{2};X_{1}X_{2},gx_{1}x_{2})-B(x_{1}\otimes g;X_{2},gx_{1})=0
\end{equation*}%
which holds in view of the form of the elements.

\subsubsection{Case $X_{1}X_{2}\otimes g\otimes gx_{1}x_{2}$}

From the first summand of the left side we get%
\begin{eqnarray*}
&&\sum_{a,b_{1},b_{2},d,e_{1},e_{2}=0}^{1}\sum_{l_{1}=0}^{b_{1}}%
\sum_{l_{2}=0}^{b_{2}}\sum_{u_{1}=0}^{e_{1}}\sum_{u_{2}=0}^{e_{2}}\left(
-1\right) ^{\alpha \left( x_{1};l_{1},l_{2},u_{1},u_{2}\right) }B(g\otimes
gx_{1}x_{2};G^{a}X_{1}^{b_{1}}X_{2}^{b_{2}},g^{d}x_{1}^{e_{1}}x_{2}^{e_{2}})
\\
&&G^{a}X_{1}^{b_{1}-l_{1}}X_{2}^{b_{2}-l_{2}}\otimes
g^{d}x_{1}^{e_{1}-u_{1}}x_{2}^{e_{2}-u_{2}}\otimes
g^{a+b_{1}+b_{2}+l_{1}+l_{2}+d+e_{1}+e_{2}+u_{1}+u_{2}}x_{1}^{l_{1}+u_{1}+1}x_{2}^{l_{2}+u_{2}}
\end{eqnarray*}%
\begin{eqnarray*}
l_{1} &=&u_{1}=0,l_{2}+u_{2}=1, \\
a &=&0,b_{1}=1,b_{2}-l_{2}=1\Rightarrow b_{2}=1,l_{2}=0,u_{2}=1. \\
d &=&1,e_{1}=0,e_{2}=u_{2}=1.
\end{eqnarray*}%
Since $\alpha \left( x_{1};0,0,0,1\right) =1,$ we obtain%
\begin{equation*}
-B(g\otimes gx_{1}x_{2};X_{1}X_{2},gx_{2})X_{1}X_{2}\otimes g\otimes
gx_{1}x_{2}.
\end{equation*}%
From the second summand of the left side we get%
\begin{eqnarray*}
l_{1}+u_{1} &=&1,l_{2}+u_{2}=1 \\
a &=&0,b_{1}-l_{1}=1\Rightarrow b_{1}=1,l_{1}=0,u_{1}=1 \\
b_{2}-l_{2} &=&1\Rightarrow b_{2}=1,l_{2}=0,u_{2}=1 \\
d &=&1,e_{1}=u_{1}=1,e_{2}=u_{2}=1
\end{eqnarray*}%
Since $\alpha \left( 1_{H};0,0,1,1\right) =1+e_{2}\equiv 0$ we obtain%
\begin{equation*}
B(x_{1}\otimes gx_{1}x_{2};X_{1}X_{2},gx_{1}x_{2})X_{1}X_{2}\otimes g\otimes
gx_{1}x_{2}.
\end{equation*}%
By considering also the right side we obtain%
\begin{equation*}
-B(g\otimes gx_{1}x_{2};X_{1}X_{2},gx_{2})+B(x_{1}\otimes
gx_{1}x_{2};X_{1}X_{2},gx_{1}x_{2})-B(x_{1}\otimes g;X_{1}X_{2},g)=0
\end{equation*}%
which holds in view of the form of the elements.

\subsection{$B\left( x_{1}\otimes gx_{1}x_{2};GX_{1},gx_{1}x_{2}\right) $}

We deduce that%
\begin{eqnarray*}
a &=&1,b_{1}=1,b_{2}=0 \\
d &=&e_{1}=e_{2}=1
\end{eqnarray*}%
and we get%
\begin{gather*}
\left( -1\right) ^{\alpha \left( 1_{H};0,0,0,0\right) }B(x_{1}\otimes
gx_{1}x_{2};GX_{1},gx_{1}x_{2})GX_{1}\otimes gx_{1}x_{2}\otimes g \\
\left( -1\right) ^{\alpha \left( 1_{H};1,0,0,0\right) }B(x_{1}\otimes
gx_{1}x_{2};GX_{1},gx_{1}x_{2})G\otimes gx_{1}x_{2}\otimes x_{1} \\
\left( -1\right) ^{\alpha \left( 1_{H};0,0,1,0\right) }B(x_{1}\otimes
gx_{1}x_{2};GX_{1},gx_{1}x_{2})GX_{1}\otimes gx_{2}\otimes x_{1} \\
\left( -1\right) ^{\alpha \left( 1_{H};1,0,1,0\right) }B(x_{1}\otimes
gx_{1}x_{2};GX_{1},gx_{1}x_{2})GX_{1}^{1-l_{1}}\otimes gx_{2}\otimes
g^{+l_{11}}x_{1}^{1+1}=0 \\
\left( -1\right) ^{\alpha \left( 1_{H};0,0,0,1\right) }B(x_{1}\otimes
gx_{1}x_{2};GX_{1},gx_{1}x_{2})GX_{1}\otimes gx_{1}\otimes x_{2} \\
\left( -1\right) ^{\alpha \left( 1_{H};1,0,0,1\right) }B(x_{1}\otimes
gx_{1}x_{2};GX_{1},gx_{1}x_{2})G\otimes gx_{1}\otimes gx_{1}x_{2} \\
\left( -1\right) ^{\alpha \left( 1_{H};0,0,1,1\right) }B(x_{1}\otimes
gx_{1}x_{2};GX_{1},gx_{1}x_{2})GX_{1}\otimes g\otimes gx_{1}x_{2} \\
\left( -1\right) ^{\alpha \left( 1_{H};1,0,1,1\right) }B(x_{1}\otimes
gx_{1}x_{2};GX_{1},gx_{1}x_{2})GX_{1}^{1-l_{1}}\otimes g\otimes
g^{+l_{1}+1}x_{1}^{1+1}x_{2}=0
\end{gather*}

\subsubsection{Case $G\otimes gx_{1}x_{2}\otimes x_{1}$}

From the first summand of the left side we get%
\begin{eqnarray*}
l_{1} &=&u_{1}=l_{2}=u_{2}=0 \\
a &=&1,b_{1}=b_{2}=0 \\
d &=&e_{1}=e_{2}=1.
\end{eqnarray*}%
Since $\alpha \left( x_{1};0,0,0,0\right) \equiv a+b_{1}+b_{2}=1,$ we get%
\begin{equation*}
-B(g\otimes gx_{1}x_{2};G,gx_{1}x_{2})G\otimes gx_{1}x_{2}\otimes x_{1}.
\end{equation*}%
By considering the second summand of the left side, we get%
\begin{eqnarray*}
l_{1}+u_{1} &=&1,l_{2}=u_{2}=0 \\
a &=&1,b_{1}=l_{1},b_{2}=l_{2}=0 \\
d &=&1,e_{1}-u_{1}=1\Rightarrow e_{1}=1,u_{1}=0,b_{1}=l_{1}=1. \\
e_{2} &=&1
\end{eqnarray*}%
Since $\alpha \left( 1_{H};1,0,0,0\right) =b_{2}=0,$ we obtain%
\begin{equation*}
B(x_{1}\otimes gx_{1}x_{2};GX_{1},gx_{1}x_{2})G\otimes gx_{1}x_{2}\otimes
x_{1}
\end{equation*}%
By considering also the right side we obtain%
\begin{equation*}
-B(g\otimes gx_{1}x_{2};G,gx_{1}x_{2})+B(x_{1}\otimes
gx_{1}x_{2};GX_{1},gx_{1}x_{2})+B(x_{1}\otimes gx_{2};G,gx_{1}x_{2})=0
\end{equation*}%
which holds in view of the form of the elements.

\subsubsection{Case $GX_{1}\otimes gx_{2}\otimes x_{1}$}

From the first summand of the left side we get%
\begin{eqnarray*}
l_{1} &=&u_{1}=l_{2}=u_{2}=0 \\
a &=&1,b_{1}=1,b_{2}=0 \\
d &=&e_{2}=1,e_{1}=0
\end{eqnarray*}%
Since $\alpha \left( x_{1};0,0,0,0\right) \equiv a+b_{1}+b_{2}=0,$ we get%
\begin{equation*}
B(g\otimes gx_{1}x_{2};GX_{1},gx_{2})GX_{1}\otimes gx_{2}\otimes x_{1}.
\end{equation*}%
By considering the second summand of the left side, we get%
\begin{eqnarray*}
l_{1}+u_{1} &=&1,l_{2}=u_{2}=0 \\
a &=&1,b_{1}-l_{1}=1,\Rightarrow b_{1}=1,l_{1}=0,u_{1}=1,b_{2}=l_{2}=0 \\
d &=&1,e_{1}=u_{1}=1,e_{2}=1.
\end{eqnarray*}%
Since $\alpha \left( 1_{H};0,0,1,0\right) =e_{2}+\left( a+b_{1}+b_{2}\right)
\equiv 1,$ we obtain%
\begin{equation*}
-B(x_{1}\otimes gx_{1}x_{2};GX_{1},gx_{1}x_{2})GX_{1}\otimes gx_{2}\otimes
x_{1}.
\end{equation*}%
By considering also the right side we obtain%
\begin{equation*}
B(g\otimes gx_{1}x_{2};GX_{1},gx_{2})-B(x_{1}\otimes
gx_{1}x_{2};GX_{1},gx_{1}x_{2})+B(x_{1}\otimes gx_{2};GX_{1},gx_{2})=0
\end{equation*}%
which holds in view of the form of the elements.

\subsubsection{Case $GX_{1}\otimes gx_{1}\otimes x_{2}$}

Nothing from the first term of the left side. The second term of the left
side gives us

\begin{eqnarray*}
l_{1} &=&u_{1}=0,l_{2}+u_{2}=1 \\
a &=&1,b_{1}=1,b_{2}=l_{2} \\
d &=&1,e_{1}=1,e_{2}=u_{2}.
\end{eqnarray*}%
Since $\alpha \left( 1_{H};0,0,0,1\right) =a+b_{1}+b_{2}\equiv 0$ and $%
\alpha \left( 1_{H};0,1,0,0\right) =0$ we get%
\begin{equation*}
\left[ B(x_{1}\otimes gx_{1}x_{2};GX_{1},gx_{1}x_{2})+B(x_{1}\otimes
gx_{1}x_{2};GX_{1}X_{2},gx_{1})\right] GX_{1}\otimes gx_{1}\otimes x_{2}
\end{equation*}%
By considering also the right side, we get

\begin{equation*}
B(x_{1}\otimes gx_{1}x_{2};GX_{1},gx_{1}x_{2})+B(x_{1}\otimes
gx_{1}x_{2};GX_{1}X_{2},gx_{1})-B(x_{1}\otimes gx_{1};GX_{1},gx_{1})=0
\end{equation*}%
which holds in view of the form of the elements.

\subsubsection{Case $G\otimes gx_{1}\otimes gx_{1}x_{2}$}

From the first summand of the left side we get%
\begin{eqnarray*}
l_{1} &=&u_{1}=0,l_{2}+u_{2}=1 \\
a &=&1,b_{1}=0,b_{2}=l_{2} \\
d &=&e_{1}=1,e_{2}=u_{2}
\end{eqnarray*}%
Since $\alpha \left( x_{1};0,0,0,1\right) \equiv 1$ and $\alpha \left(
x_{1};0,1,0,0\right) \equiv a+b_{1}+b_{2}+1\equiv 1$ we get%
\begin{equation*}
\left[ -B(g\otimes gx_{1}x_{2};G,gx_{1}x_{2})-B(g\otimes
gx_{1}x_{2};GX_{2},gx_{1})\right] G\otimes gx_{1}\otimes gx_{1}x_{2}
\end{equation*}%
By considering the second summand of the left side, we get%
\begin{eqnarray*}
l_{1}+u_{1} &=&1,l_{2}+u_{2}=1 \\
a &=&1,b_{1}=l_{1},b_{2}=l_{2} \\
d &=&1,e_{1}-u_{1}=1\Rightarrow e_{1}=1,u_{1}=0,b_{1}=l_{1}=1,e_{2}=u_{2}.
\end{eqnarray*}%
Since $\alpha \left( 1_{H};1,0,0,1\right) \equiv 0$ and $\alpha \left(
1_{H};1,1,0,0\right) \equiv 0$%
\begin{equation*}
\left[ B(x_{1}\otimes gx_{1}x_{2};GX_{1},gx_{1}x_{2})+B(x_{1}\otimes
gx_{1}x_{2};GX_{1}X_{2},gx_{1})\right] G\otimes gx_{1}\otimes gx_{1}x_{2}.
\end{equation*}%
By considering also the right side we get%
\begin{gather*}
-B(g\otimes gx_{1}x_{2};G,gx_{1}x_{2})-B(g\otimes gx_{1}x_{2};GX_{2},gx_{1})+
\\
+B(x_{1}\otimes gx_{1}x_{2};GX_{1},gx_{1}x_{2})+B(x_{1}\otimes
gx_{1}x_{2};GX_{1}X_{2},gx_{1})-B(x_{1}\otimes g;G,gx_{1})=0
\end{gather*}

which holds in view of the form of the elements.

\subsubsection{Case $GX_{1}\otimes g\otimes gx_{1}x_{2}$}

From the first summand of the left side we get%
\begin{eqnarray*}
l_{1} &=&u_{1}=0,l_{2}+u_{2}=1 \\
a &=&1,b_{1}=1,b_{2}=l_{2} \\
d &=&1,e_{1}=0,e_{2}=u_{2}
\end{eqnarray*}%
Since $\alpha \left( x_{1};0,0,0,1\right) \equiv 1$ and $\alpha \left(
x_{1};0,1,0,0\right) \equiv a+b_{1}+b_{2}+1\equiv 0$ we get%
\begin{equation*}
\left[ -B(g\otimes gx_{1}x_{2};GX_{1},gx_{2})+B(g\otimes
gx_{1}x_{2};GX_{1}X_{2},g)\right] GX_{1}\otimes g\otimes gx_{1}x_{2}
\end{equation*}%
By considering the second summand of the left side, we get%
\begin{eqnarray*}
l_{1}+u_{1} &=&1,l_{2}+u_{2}=1 \\
a &=&1,b_{1}-l_{1}=1\Rightarrow b_{1}=1,l_{1}=0,u_{1}=1,b_{2}=l_{2} \\
d &=&1,e_{1}=u_{1}=1,e_{2}=u_{2}.
\end{eqnarray*}%
Since $\alpha \left( 1_{H};0,0,0,1\right) \equiv a+b_{1}+b_{2}\equiv 0$ and $%
\alpha \left( 1_{H};0,1,0,0\right) \equiv 0,$ we obtain%
\begin{equation*}
\left[ B(x_{1}\otimes gx_{1}x_{2};GX_{1},gx_{1}x_{2})+B(x_{1}\otimes
gx_{1}x_{2};GX_{1}X_{2},gx_{1})\right] GX_{1}\otimes g\otimes gx_{1}x_{2}.
\end{equation*}%
By considering also the right side, we get%
\begin{gather*}
-B(g\otimes gx_{1}x_{2};GX_{1},gx_{2})+B(g\otimes gx_{1}x_{2};GX_{1}X_{2},g)+
\\
B(x_{1}\otimes gx_{1}x_{2};GX_{1},gx_{1}x_{2})+B(x_{1}\otimes
gx_{1}x_{2};GX_{1}X_{2},gx_{1})-B(x_{1}\otimes g;GX_{1},g)=0
\end{gather*}%
which holds in view of the form of the elements.

\subsection{$B\left( x_{1}\otimes gx_{1}x_{2};GX_{2},gx_{1}x_{2}\,\right) $}

We have%
\begin{eqnarray*}
a &=&1,b_{1}=0,b_{2}=1 \\
d &=&e_{1}=e_{2}=1
\end{eqnarray*}%
and we get%
\begin{gather*}
\left( -1\right) ^{\alpha \left( 1_{H};0,0,0,0\right) }B\left( x_{1}\otimes
gx_{1}x_{2};GX_{2},gx_{1}x_{2}\,\right) GX_{2}\otimes gx_{1}x_{2}\otimes g+
\\
\left( -1\right) ^{\alpha \left( 1_{H};0,1,0,0\right) }B\left( x_{1}\otimes
gx_{1}x_{2};GX_{2},gx_{1}x_{2}\,\right) G\otimes gx_{1}x_{2}\otimes x_{2}+ \\
\left( -1\right) ^{\alpha \left( 1_{H};0,0,1,0\right) }B\left( x_{1}\otimes
gx_{1}x_{2};GX_{2},gx_{1}x_{2}\,\right) GX_{2}\otimes gx_{2}\otimes x_{1}+ \\
\left( -1\right) ^{\alpha \left( 1_{H};0,1,1,0\right) }B\left( x_{1}\otimes
gx_{1}x_{2};GX_{2},gx_{1}x_{2}\,\right) G\otimes gx_{2}\otimes gx_{1}x_{2}+
\\
\left( -1\right) ^{\alpha \left( 1_{H};0,0,0,1\right) }B\left( x_{1}\otimes
gx_{1}x_{2};GX_{2},gx_{1}x_{2}\,\right) GX_{2}\otimes gx_{1}\otimes x_{2}+ \\
\left( -1\right) ^{\alpha \left( 1_{H};0,1,0,1\right) }B\left( x_{1}\otimes
gx_{1}x_{2};GX_{2},gx_{1}x_{2}\,\right) G\otimes gx_{1}\otimes
g^{l_{2}}x_{2}^{1+1}=0 \\
\left( -1\right) ^{\alpha \left( 1_{H};0,0,1,1\right) }B\left( x_{1}\otimes
gx_{1}x_{2};GX_{2},gx_{1}x_{2}\,\right) GX_{2}\otimes g\otimes gx_{1}x_{2} \\
\left( -1\right) ^{\alpha \left( 1_{H};0,1,1,1\right) }B\left( x_{1}\otimes
gx_{1}x_{2};GX_{2},gx_{1}x_{2}\,\right) GX_{2}^{1-l_{2}}\otimes g\otimes
g^{l_{2}+1}x_{1}x_{2}^{1+1}=0
\end{gather*}

\subsubsection{Case $G\otimes gx_{1}x_{2}\otimes x_{2}$}

Nothing from the first term of the left side. The second term of the left
side gives us

\begin{eqnarray*}
l_{1} &=&u_{1}=0,l_{2}+u_{2}=1 \\
a &=&1,b_{1}=0,b_{2}=l_{2} \\
d &=&1,e_{1}=1,e_{2}-u_{2}=1\Rightarrow e_{2}=1,u_{2}=0,b_{2}=l_{2}=1.
\end{eqnarray*}%
Since $\alpha \left( 1_{H};0,1,0,0\right) =0$ we get%
\begin{equation*}
B(x_{1}\otimes gx_{1}x_{2};GX_{2},gx_{1}x_{2})G\otimes gx_{1}x_{2}\otimes
x_{2}.
\end{equation*}%
By considering also the right side, we get%
\begin{equation*}
B(x_{1}\otimes gx_{1}x_{2};GX_{2},gx_{1}x_{2})-B(x_{1}\otimes
gx_{1};G,gx_{1}x_{2})=0
\end{equation*}%
which holds in view of the form of the elements.

\subsubsection{Case $GX_{2}\otimes gx_{2}\otimes x_{1}$}

From the first summand of the left side we get%
\begin{eqnarray*}
l_{1} &=&u_{1}=l_{2}=u_{2}=0 \\
a &=&1,b_{1}=0,b_{2}=1 \\
d &=&e_{2}=1,e_{1}=0
\end{eqnarray*}%
Since $\alpha \left( x_{1};0,0,0,0\right) \equiv a+b_{1}+b_{2}\equiv 0,$ we
get%
\begin{equation*}
B(g\otimes gx_{1}x_{2};GX_{2},gx_{2})GX_{2}\otimes gx_{2}\otimes x_{1}.
\end{equation*}%
By considering the second summand of the left side, we get%
\begin{eqnarray*}
l_{1}+u_{1} &=&1,l_{2}=u_{2}=0 \\
a &=&1,b_{1}=l_{1},b_{2}=1 \\
d &=&1,e_{1}=u_{1},e_{2}=1
\end{eqnarray*}%
Since $\alpha \left( 1_{H};0,0,1,0\right) =e_{2}+\left( a+b_{1}+b_{2}\right)
\equiv 1$ and $\alpha \left( 1_{H};1,0,0,0\right) =b_{2}=1,$ we obtain%
\begin{equation*}
\left[ -B(x_{1}\otimes gx_{1}x_{2};GX_{2},gx_{1}x_{2})-B(x_{1}\otimes
gx_{1}x_{2};GX_{1}X_{2},gx_{2})\right] GX_{2}\otimes gx_{2}\otimes x_{1}.
\end{equation*}%
By considering also the right side, we get%
\begin{gather*}
B(g\otimes gx_{1}x_{2};GX_{2},gx_{2})-B(x_{1}\otimes
gx_{1}x_{2};GX_{2},gx_{1}x_{2}) \\
-B(x_{1}\otimes gx_{1}x_{2};GX_{1}X_{2},gx_{2})+B(x_{1}\otimes
gx_{2};GX_{2},gx_{2})=0
\end{gather*}%
which holds in view of the form of the elements.

\subsubsection{Case $G\otimes gx_{2}\otimes gx_{1}x_{2}$}

From the first summand of the left side we get%
\begin{eqnarray*}
l_{1} &=&u_{1}=0,l_{2}+u_{2}=1 \\
a &=&1,b_{1}=0,b_{2}=l_{2} \\
d &=&1,e_{1}=0,e_{2}-u_{2}=1\Rightarrow e_{2}=1,u_{2}=0,b_{2}=l_{2}=1
\end{eqnarray*}%
Since $\alpha \left( x_{1};0,1,0,0\right) \equiv a+b_{1}+b_{2}+1\equiv 1$ we
get%
\begin{equation*}
-B(g\otimes gx_{1}x_{2};GX_{2},gx_{2})G\otimes gx_{2}\otimes gx_{1}x_{2}
\end{equation*}%
By considering the second summand of the left side, we get%
\begin{eqnarray*}
l_{1}+u_{1} &=&1,l_{2}+u_{2}=1 \\
a &=&1,b_{1}=l_{1},b_{2}=l_{2} \\
d &=&1,e_{1}=u_{1},e_{2}-u_{2}=1\Rightarrow e_{2}=1,u_{2}=0,b_{2}=l_{2}=1.
\end{eqnarray*}%
Since $\alpha \left( 1_{H};0,1,1,0\right) \equiv e_{2}+a+b_{1}+b_{2}+1\equiv
0$ and $\alpha \left( 1_{H};1,1,0,0\right) \equiv 1+b_{2}\equiv 0,$ we obtain%
\begin{equation*}
\left[ B(x_{1}\otimes gx_{1}x_{2};GX_{2},gx_{1}x_{2})+B(x_{1}\otimes
gx_{1}x_{2};GX_{1}X_{2},gx_{2})\right] G\otimes gx_{2}\otimes gx_{1}x_{2}.
\end{equation*}%
By considering also the right side, we get%
\begin{gather*}
-B(g\otimes gx_{1}x_{2};GX_{2},gx_{2})+B(x_{1}\otimes
gx_{1}x_{2};GX_{2},gx_{1}x_{2})+ \\
+B(x_{1}\otimes gx_{1}x_{2};GX_{1}X_{2},gx_{2})-B(x_{1}\otimes g;G,gx_{2})=0
\end{gather*}%
which holds in view of the form of the elements.

\subsubsection{Case $GX_{2}\otimes gx_{1}\otimes x_{2}$}

Nothing from the first term of the left side. The second term of the left
side gives us

\begin{eqnarray*}
l_{1} &=&u_{1}=0,l_{2}+u_{2}=1 \\
a &=&1,b_{1}=0,b_{2}-l_{2}=1\Rightarrow b_{2}=1,l_{2}=0,u_{2}=1 \\
d &=&1,e_{1}=1,e_{2}=u_{2}=1.
\end{eqnarray*}%
Since $\alpha \left( 1_{H};0,0,0,1\right) =a+b_{1}+b_{2}\equiv 0,$ we get%
\begin{equation*}
B(x_{1}\otimes gx_{1}x_{2};GX_{2},gx_{1}x_{2})GX_{2}\otimes gx_{1}\otimes
x_{2}.
\end{equation*}%
By considering also the right side, we get%
\begin{equation*}
B(x_{1}\otimes gx_{1}x_{2};GX_{2},gx_{1}x_{2})-B(x_{1}\otimes
gx_{1};GX_{2},gx_{1})=0
\end{equation*}%
which holds in view of the form of the elements.

\subsubsection{Case $GX_{2}\otimes g\otimes gx_{1}x_{2}$}

From the first summand of the left side we get%
\begin{eqnarray*}
l_{1} &=&u_{1}=0,l_{2}+u_{2}=1 \\
a &=&1,b_{1}=0,b_{2}-l_{2}=1\Rightarrow b_{2}=1,l_{2}=0,u_{2}=1 \\
d &=&1,e_{1}=0,e_{2}=u_{2}=1
\end{eqnarray*}%
Since $\alpha \left( x_{1};0,0,0,1\right) \equiv 1$we get%
\begin{equation*}
-B(g\otimes gx_{1}x_{2};GX_{2},gx_{2})GX_{2}\otimes g\otimes gx_{1}x_{2}
\end{equation*}%
By considering the second summand of the left side, we get%
\begin{eqnarray*}
l_{1}+u_{1} &=&1,l_{2}+u_{2}=1 \\
a &=&1,b_{1}=l_{1},b_{2}-l_{2}=1\Rightarrow b_{2}=1,l_{2}=0,u_{2}=1 \\
d &=&1,e_{1}=u_{1},e_{2}=u_{2}=1.
\end{eqnarray*}%
Since $\alpha \left( 1_{H};0,0,1,1\right) \equiv 1+e_{2}\equiv 0$ and $%
\alpha \left( 1_{H};1,0,0,1\right) \equiv a+b_{1}\equiv 0$ we obtain%
\begin{equation*}
\left[ B(x_{1}\otimes gx_{1}x_{2};GX_{2},gx_{1}x_{2})+B(x_{1}\otimes
gx_{1}x_{2};GX_{1}X_{2},gx_{2})\right] GX_{2}\otimes g\otimes gx_{1}x_{2}.
\end{equation*}%
By considering also the right side, we get%
\begin{gather*}
-B(g\otimes gx_{1}x_{2};GX_{2},gx_{2})+B(x_{1}\otimes
gx_{1}x_{2};GX_{2},gx_{1}x_{2})+ \\
+B(x_{1}\otimes gx_{1}x_{2};GX_{1}X_{2},gx_{2})-B(x_{1}\otimes g;GX_{2},g)=0
\end{gather*}%
which holds in view of the form of the elements.

\subsection{$B\left( x_{1}\otimes gx_{1}x_{2};GX_{1}X_{2},gx_{1}\right) $}

We get%
\begin{eqnarray*}
a &=&b_{1}=b_{2}=1 \\
d &=&e_{1}=1,e_{2}=0
\end{eqnarray*}%
so that we obtain%
\begin{gather*}
\left( -1\right) ^{\alpha \left( 1_{H};0,0,0,0\right) }B(x_{1}\otimes
gx_{1}x_{2};GX_{1}X_{2},gx_{1})GX_{1}X_{2}\otimes gx_{1}\otimes g+ \\
\left( -1\right) ^{\alpha \left( 1_{H};1,0,0,0\right) }B(x_{1}\otimes
gx_{1}x_{2};GX_{1}X_{2},gx_{1})GX_{2}\otimes gx_{1}\otimes x_{1}+ \\
\left( -1\right) ^{\alpha \left( 1_{H};0,1,0,0\right) }B(x_{1}\otimes
gx_{1}x_{2};GX_{1}X_{2},gx_{1})GX_{1}\otimes gx_{1}\otimes x_{2}+ \\
\left( -1\right) ^{\alpha \left( 1_{H};1,1,0,0\right) }B(x_{1}\otimes
gx_{1}x_{2};GX_{1}X_{2},gx_{1})G\otimes gx_{1}\otimes gx_{1}x_{2}+ \\
\left( -1\right) ^{\alpha \left( 1_{H};0,0,1,0\right) }B(x_{1}\otimes
gx_{1}x_{2};GX_{1}X_{2},gx_{1})GX_{1}X_{2}\otimes g\otimes x_{1}+ \\
\left( -1\right) ^{\alpha \left( 1_{H};1,0,1,0\right) }B(x_{1}\otimes
gx_{1}x_{2};GX_{1}X_{2},gx_{1})GX_{1}^{1-l_{1}}X_{2}\otimes g\otimes
g^{l_{1}}x_{1}^{1+1}=0 \\
\left( -1\right) ^{\alpha \left( 1_{H};0,1,1,0\right) }B(x_{1}\otimes
gx_{1}x_{2};GX_{1}X_{2},gx_{1})GX_{1}\otimes g\otimes gx_{1}x_{2}+ \\
\left( -1\right) ^{\alpha \left( 1_{H};1,1,1,0\right) }B(x_{1}\otimes
gx_{1}x_{2};GX_{1}X_{2},gx_{1})GX_{1}^{1-l_{1}}\otimes g\otimes
g^{l_{1}+1}x_{1}^{1+1}x_{2}=0
\end{gather*}

\subsubsection{Case $GX_{2}\otimes gx_{1}\otimes x_{1}$}

From the first summand of the left side we get%
\begin{eqnarray*}
l_{1} &=&u_{1}=l_{2}=u_{2}=0 \\
a &=&1,b_{1}=0,b_{2}=1 \\
d &=&e_{1}=1,e_{2}=0
\end{eqnarray*}%
Since $\alpha \left( x_{1};0,0,0,0\right) \equiv a+b_{1}+b_{2}\equiv 0,$ we
get%
\begin{equation*}
B(g\otimes gx_{1}x_{2};GX_{2},gx_{1})GX_{2}\otimes gx_{1}\otimes x_{1}.
\end{equation*}%
By considering the second summand of the left side, we get%
\begin{eqnarray*}
l_{1}+u_{1} &=&1,l_{2}=u_{2}=0 \\
a &=&1,b_{1}=l_{1},b_{2}=1 \\
d &=&1,e_{1}-u_{1}=1\Rightarrow e_{1}=1,u_{1}=0,b_{1}=l_{1}=1,e_{2}=0.
\end{eqnarray*}%
Since $\alpha \left( 1_{H};1,0,0,0\right) =b_{2}=1,$ we obtain%
\begin{equation*}
-B(x_{1}\otimes gx_{1}x_{2};GX_{1}X_{2},gx_{1})GX_{2}\otimes gx_{1}\otimes
x_{1}.
\end{equation*}%
By considering also the right side, we get%
\begin{equation*}
B(g\otimes gx_{1}x_{2};GX_{2},gx_{1})-B(x_{1}\otimes
gx_{1}x_{2};GX_{1}X_{2},gx_{1})+B(x_{1}\otimes gx_{2};GX_{2},gx_{1})=0
\end{equation*}%
which holds in view of the form of the elements.

\subsubsection{Case $GX_{1}\otimes gx_{1}\otimes x_{2}$}

Nothing from the first term of the left side. The second term of the left
side gives us

\begin{eqnarray*}
l_{1} &=&u_{1}=0,l_{2}+u_{2}=1 \\
a &=&1,b_{1}=1,b_{2}=l_{2} \\
d &=&1,e_{1}=1,e_{2}=u_{2}
\end{eqnarray*}%
Since $\alpha \left( 1_{H};0,0,0,1\right) =a+b_{1}+b_{2}\equiv 0$ and $%
\alpha \left( 1_{H};0,1,0,0\right) =0,$ we get%
\begin{equation*}
\left[ B(x_{1}\otimes gx_{1}x_{2};GX_{1},gx_{1}x_{2})+B(x_{1}\otimes
gx_{1}x_{2};GX_{1}X_{2},gx_{1})\right] GX_{1}\otimes gx_{1}\otimes x_{2}.
\end{equation*}%
By considering also the right side, we get%
\begin{equation*}
B(x_{1}\otimes gx_{1}x_{2};GX_{1},gx_{1}x_{2})+B(x_{1}\otimes
gx_{1}x_{2};GX_{1}X_{2},gx_{1})-B(x_{1}\otimes gx_{1};GX_{1},gx_{1})=0
\end{equation*}%
which holds in view of the form of the elements.

\subsubsection{Case $G\otimes gx_{1}\otimes gx_{1}x_{2}$}

Already done in $B\left( x_{1}\otimes gx_{1}x_{2};GX_{1},gx_{1}x_{2}\right) $%
.

\subsubsection{Case $GX_{1}X_{2}\otimes g\otimes x_{1}$}

From the first summand of the left side we get%
\begin{eqnarray*}
l_{1} &=&u_{1}=l_{2}=u_{2}=0 \\
a &=&b_{1}=b_{2}=1 \\
d &=&1,e_{1}=e_{2}=0
\end{eqnarray*}%
Since $\alpha \left( x_{1};0,0,0,0\right) \equiv a+b_{1}+b_{2}\equiv 1,$ we
get%
\begin{equation*}
-B(g\otimes gx_{1}x_{2};GX_{1}X_{2},g)GX_{1}X_{2}\otimes g\otimes x_{1}.
\end{equation*}%
By considering the second summand of the left side, we get%
\begin{eqnarray*}
l_{1}+u_{1} &=&1,l_{2}=u_{2}=0 \\
a &=&1,b_{1}-l_{1}=1\Rightarrow b_{1}=1,l_{1}=0,u_{1}=1,b_{2}=1 \\
d &=&1,e_{1}=u_{1}=1,e_{2}=0
\end{eqnarray*}%
Since $\alpha \left( 1_{H},0,0,1,0\right) =e_{2}+\left( a+b_{1}+b_{2}\right)
\equiv 1,$ we obtain%
\begin{equation*}
-B(x_{1}\otimes gx_{1}x_{2};GX_{1}X_{2},gx_{1})GX_{1}X_{2}\otimes g\otimes
x_{1}.
\end{equation*}%
By considering also the right side, we get%
\begin{equation*}
-B(g\otimes gx_{1}x_{2};GX_{1}X_{2},g)-B(x_{1}\otimes
gx_{1}x_{2};GX_{1}X_{2},gx_{1})+B(x_{1}\otimes gx_{2};GX_{1}X_{2},g)=0
\end{equation*}%
which holds in view of the form of the elements.

\subsubsection{Case $GX_{1}\otimes g\otimes gx_{1}x_{2}$}

Already done in $B\left( x_{1}\otimes gx_{1}x_{2};GX_{1},gx_{1}x_{2}\right)
. $

\section{$B\left( x_{2}\otimes x_{1}\right) $}

From $\left( \ref{simplx}\right) $ we get%
\begin{eqnarray}
B(x_{2}\otimes x_{1}) &=&B(x_{2}\otimes 1_{H})(1_{A}\otimes
x_{1})-(1_{A}\otimes gx_{1})B(x_{2}\otimes 1_{H})(1_{A}\otimes g)
\label{form x2otx1} \\
&&-(1_{A}\otimes g)B(gx_{1}x_{2}\otimes 1_{H})(1_{A}\otimes g)  \notag
\end{eqnarray}%
We obtain%
\begin{eqnarray*}
B(x_{2}\otimes x_{1}) &=&-B(gx_{1}x_{2}\otimes 1_{H};1_{A},g)1_{A}\otimes g+
\\
&&+B(gx_{1}x_{2}\otimes 1_{H};1_{A},x_{1})1_{A}\otimes
x_{1}+B(gx_{1}x_{2}\otimes 1_{H};1_{A},x_{2})1_{A}\otimes x_{2} \\
&&+\left[ -2B\left( x_{2}\otimes 1_{H};1_{A},gx_{2}\right)
-B(gx_{1}x_{2}\otimes 1_{H};1_{A},gx_{1}x_{2})\right] 1_{A}\otimes
gx_{1}x_{2} \\
&&-B(gx_{1}x_{2}\otimes 1_{H};G,1_{H})G\otimes 1_{H}+ \\
&&-B(gx_{1}x_{2}\otimes 1_{H};G,x_{1}x_{2})G\otimes x_{1}x_{2}+ \\
&&+\left[ -2B\left( x_{2}\otimes 1_{H};G,g\right) +B(gx_{1}x_{2}\otimes
1_{H};G,gx_{1})\right] G\otimes gx_{1}+ \\
&&+B(gx_{1}x_{2}\otimes 1_{H};G,gx_{2})G\otimes gx_{2}+ \\
&&+\left[ +B(x_{2}\otimes 1_{H};1_{A},1_{H})-B(gx_{1}x_{2}\otimes
1_{H};1_{A},x_{1})\right] X_{1}\otimes 1_{H}+ \\
&&+B(x_{2}\otimes 1_{H};1_{A},x_{1}x_{2})X_{1}\otimes
x_{1}x_{2}+B(x_{2}\otimes 1_{H};1_{A},gx_{1})X_{1}\otimes gx_{1}+ \\
&&+\left[ -B(x_{2}\otimes 1_{H};1_{A},gx_{2})-B(gx_{1}x_{2}\otimes
1_{H};1_{A},gx_{1}x_{2})\right] X_{1}\otimes gx_{2}+ \\
&&+\left[ -B(x_{1}\otimes 1_{H};1_{A},1_{H})-B(gx_{1}x_{2}\otimes
1_{H};1_{A},x_{2})\right] X_{2}\otimes 1_{H} \\
&&-B(x_{1}\otimes 1_{H};1_{A},x_{1}x_{2})X_{2}\otimes x_{1}x_{2}+ \\
&&+\left[
\begin{array}{c}
+2B(g\otimes 1_{H};1_{A},g)+2B(x_{2}\otimes \ 1_{H};1_{A},gx_{2}) \\
+B(x_{1}\otimes 1_{H};1_{A},gx_{1})+B(gx_{1}x_{2}\otimes
1_{H};1_{A},gx_{1}x_{2})%
\end{array}%
\right] X_{2}\otimes gx_{1}+ \\
&&+B(x_{1}\otimes 1_{H};1_{A},gx_{2})X_{2}\otimes gx_{2}+ \\
&&-\left[
\begin{array}{c}
+B(g\otimes 1_{H};1_{A},g)+B(x_{2}\otimes \ 1_{H};1_{A},gx_{2}) \\
+B(x_{1}\otimes 1_{H};1_{A},gx_{1})+B(gx_{1}x_{2}\otimes
1_{H};1_{A},gx_{1}x_{2})%
\end{array}%
\right] X_{1}X_{2}\otimes g+ \\
&&+\left[ B(g\otimes 1_{H};1_{A},x_{1})+B(x_{2}\otimes \
1_{H};1_{A},x_{1}x_{2})\right] X_{1}X_{2}\otimes x_{1}+ \\
&&+\left[ B(g\otimes 1_{H};X_{2},1_{H})-B(x_{1}\otimes
1_{H};1_{A},x_{1}x_{2})\right] X_{1}X_{2}\otimes x_{2}+ \\
&&+B\left( g\otimes 1_{H};1_{A},gx_{1}x_{2}\right) X_{1}X_{2}\otimes
gx_{1}x_{2} \\
&&-\left[ B(x_{2}\otimes 1_{H};G,g)-B(gx_{1}x_{2}\otimes 1_{H};G,gx_{1})%
\right] GX_{1}\otimes g+ \\
&&+B(x_{2}\otimes 1_{H};G,x_{1})GX_{1}\otimes x_{1}+ \\
&&+\left[ B(x_{2}\otimes 1_{H};G,x_{2})+B(gx_{1}x_{2}\otimes
1_{H};G,x_{1}x_{2})\right] GX_{1}\otimes x_{2}+ \\
&&+B(x_{2}\otimes 1_{H};G,gx_{1}x_{2})GX_{1}\otimes gx_{1}x_{2}+ \\
&&-\left[ -B(x_{1}\otimes 1_{H};G,g)-B(gx_{1}x_{2}\otimes 1_{H};G,gx_{2})%
\right] GX_{2}\otimes g+ \\
&&+\left[ -B(x_{1}\otimes 1_{H};G,x_{1})-B(gx_{1}x_{2}\otimes
1_{H};G,x_{1}x_{2})\right] GX_{2}\otimes x_{1}+ \\
&&-B(x_{1}\otimes 1_{H};G,x_{2})GX_{2}\otimes x_{2} \\
&&+\left[ -2B(g\otimes 1_{H};G,gx_{2})+B(x_{1}\otimes 1_{H};G,gx_{1}x_{2})%
\right] GX_{2}\otimes gx_{1}x_{2}+ \\
&&-\left[
\begin{array}{c}
B(g\otimes 1_{H};G,1_{H})+B(x_{2}\otimes \ 1_{H};G,x_{2}) \\
+B(x_{1}\otimes 1_{H};G,x_{1})+B(gx_{1}x_{2}\otimes 1_{H};G,x_{1}x_{2})%
\end{array}%
\right] GX_{1}X_{2}\otimes 1_{H} \\
&&-B(g\otimes 1_{H};G,x_{1}x_{2})GX_{1}X_{2}\otimes x_{1}x_{2} \\
&&+\left[ -B(g\otimes 1_{H};G,gx_{1})-B(x_{2}\otimes \ 1_{H};G,gx_{1}x_{2})%
\right] GX_{1}X_{2}\otimes gx_{1}+ \\
&&+\left[ B(g\otimes 1_{H};G,gx_{2})-B(x_{1}\otimes 1_{H};G,gx_{1}x_{2})%
\right] GX_{1}X_{2}\otimes gx_{2}+
\end{eqnarray*}

We write the Casimir formula for $B\left( x_{2}\otimes x_{1}\right) $%
\begin{eqnarray*}
&&\sum_{a,b_{1},b_{2},d,e_{1},e_{2}=0}^{1}\sum_{l_{1}=0}^{b_{1}}%
\sum_{l_{2}=0}^{b_{2}}\sum_{u_{1}=0}^{e_{1}}\sum_{u_{2}=0}^{e_{2}}\left(
-1\right) ^{\alpha \left( x_{2};l_{1},l_{2},u_{1},u_{2}\right) } \\
&&B(g\otimes
x_{1};G^{a}X_{1}^{b_{1}}X_{2}^{b_{2}},g^{d}x_{1}^{e_{1}}x_{2}^{e_{2}}) \\
&&G^{a}X_{1}^{b_{1}-l_{1}}X_{2}^{b_{2}-l_{2}}\otimes
g^{d}x_{1}^{e_{1}-u_{1}}x_{2}^{e_{2}-u_{2}}\otimes
g^{a+b_{1}+b_{2}+l_{1}+l_{2}+d+e_{1}+e_{2}+u_{1}+u_{2}}x_{1}^{l_{1}+u_{1}}x_{2}^{l_{2}+u_{2}+1}+
\\
&&+\sum_{a,b_{1},b_{2},d,e_{1},e_{2}=0}^{1}\sum_{l_{1}=0}^{b_{1}}%
\sum_{l_{2}=0}^{b_{2}}\sum_{u_{1}=0}^{e_{1}}\sum_{u_{2}=0}^{e_{2}}\left(
-1\right) ^{\alpha \left( 1_{H};l_{1},l_{2},u_{1},u_{2}\right) } \\
&&B(x_{2}\otimes
x_{1};G^{a}X_{1}^{b_{1}}X_{2}^{b_{2}},g^{d}x_{1}^{e_{1}}x_{2}^{e_{2}}) \\
&&G^{a}X_{1}^{b_{1}-l_{1}}X_{2}^{b_{2}-l_{2}}\otimes
g^{d}x_{1}^{e_{1}-u_{1}}x_{2}^{e_{2}-u_{2}}\otimes
g^{a+b_{1}+b_{2}+l_{1}+l_{2}+d+e_{1}+e_{2}+u_{1}+u_{2}}x_{1}^{l_{1}+u_{1}}x_{2}^{l_{2}+u_{2}}
\\
&=&B^{A}(x_{2}\otimes x_{1})\otimes B^{H}(x_{2}\otimes x_{1})\otimes g+ \\
&&B^{A}(x_{2}\otimes 1_{H})\otimes B^{H}(x_{2}\otimes 1_{H})\otimes x_{1}
\end{eqnarray*}

\subsection{\ $B\left( x_{2}\otimes x_{1};1_{A},g\right) $}

We get%
\begin{equation*}
B(x_{2}\otimes x_{1};1_{A},g)1_{A}\otimes g\otimes g
\end{equation*}

\subsubsection{Case $1_{A}\otimes g\otimes g$}

This case is trivial

\subsubsection{\ $B\left( x_{2}\otimes x_{1};1_{A},x_{1}\right) $}

We get%
\begin{eqnarray*}
&&+\left( -1\right) ^{\alpha \left( 1_{H};0,0,0,0\right) }B\left(
x_{2}\otimes x_{1};1_{A},x_{1}\right) 1_{A}\otimes x_{1}\otimes g+ \\
&&+\left( -1\right) ^{\alpha \left( 1_{H};0,0,1,0\right) }B\left(
x_{2}\otimes x_{1};1_{A},x_{1}\right) 1_{A}\otimes 1_{H}\otimes x_{1}
\end{eqnarray*}

\subsubsection{Case $1_{A}\otimes x_{1}\otimes g$}

it is trivial.

\subsubsection{Case $1_{A}\otimes 1_{H}\otimes x_{1}$}

Nothing from the first summand of the left side. The second summand gives us%
\begin{eqnarray*}
l_{1}+u_{1} &=&1 \\
l_{2} &=&u_{2}=0 \\
a &=&b_{2}=0,b_{1}=l_{1} \\
d &=&e_{2}=0,e_{1}=u_{1}
\end{eqnarray*}%
since $\alpha \left( 1_{H};0,0,1,0\right) =e_{2}+\left( a+b_{1}+b_{2}\right)
\equiv 0$ and $\alpha \left( 1_{H};1,0,0,0\right) =b_{2}=0,$ we get%
\begin{equation*}
\left[ B(x_{2}\otimes x_{1};1_{A},x_{1})+B(x_{2}\otimes x_{1};X_{1},1_{H})%
\right] 1_{A}\otimes 1_{H}\otimes x_{1}
\end{equation*}%
By considering also the right side, we obtain%
\begin{equation*}
B(x_{2}\otimes x_{1};1_{A},x_{1})+B(x_{2}\otimes
x_{1};X_{1},1_{H})-B(x_{2}\otimes 1_{H};1_{A},1_{H})=0
\end{equation*}

which holds in view of the form of the elements.

\subsection{\ $B\left( x_{2}\otimes x_{1};1_{A},x_{2}\right) $}

We get%
\begin{eqnarray*}
&&\left( -1\right) ^{\alpha \left( 1_{H};0,0,0,0\right) }B(x_{2}\otimes
x_{1};1_{A},x_{2})1_{A}\otimes x_{2}\otimes g+ \\
&&\left( -1\right) ^{\alpha \left( 1_{H};0,0,0,1\right) }B(x_{2}\otimes
x_{1};1_{A},x_{2})1_{A}\otimes 1_{H}\otimes x_{2}
\end{eqnarray*}

\subsubsection{Case $1_{A}\otimes x_{2}\otimes g$}

This case is trivial.

\subsubsection{Case $1_{A}\otimes 1_{H}\otimes x_{2}$}

From the first summand of the left side we get

~%
\begin{eqnarray*}
l_{1} &=&u_{1}=0 \\
l_{2} &=&u_{2}=0 \\
a &=&b_{1}=b_{2}=0, \\
d &=&e_{1}=e_{2}=0
\end{eqnarray*}%
Since $\alpha \left( x_{2};0,0,0,0\right) =a+b_{1}+b_{2}=0$ we get

\begin{equation*}
B(g\otimes x_{1};1_{A},1_{H})1_{A}\otimes 1_{H}\otimes x_{2}
\end{equation*}%
From the second summand of the left side we get

~%
\begin{eqnarray*}
l_{1} &=&u_{1}=0 \\
l_{2}+u_{2} &=&1 \\
a &=&b_{1}=0,b_{2}=l_{2}, \\
d &=&e_{1}=0,e_{2}=u_{2}
\end{eqnarray*}%
since $\alpha \left( 1_{H};0,0,0,1\right) \equiv a+b_{1}+b_{2}=0$ and $%
\alpha \left( 1_{H};0,1,0,0\right) \equiv 0$ we obtain%
\begin{equation*}
\left[ B(x_{2}\otimes x_{1};1_{A},x_{2})-B(x_{2}\otimes x_{1};X_{2},1_{H})%
\right] 1_{A}\otimes 1_{H}\otimes x_{2}.
\end{equation*}%
By considering also the right side of the equation, we get%
\begin{equation*}
B(g\otimes x_{1};1_{A},1_{H})+B(x_{2}\otimes
x_{1};1_{A},x_{2})+B(x_{2}\otimes x_{1};X_{2},1_{H})=0
\end{equation*}%
which holds in view of the form of the elements.

\subsection{ $B\left( x_{2}\otimes x_{1};1_{A},gx_{1}x_{2}\right) $}

We obtain%
\begin{eqnarray*}
&&\left( -1\right) ^{\alpha \left( 1_{H};0,0,0,0\right) }B(x_{2}\otimes
x_{1};1_{A},gx_{1}x_{2})1_{A}\otimes gx_{1}x_{2}\otimes g+ \\
&&\left( -1\right) ^{\alpha \left( 1_{H};0,0,0,1\right) }B(x_{2}\otimes
x_{1};1_{A},gx_{1}x_{2})1_{A}\otimes gx_{1}\otimes x_{2}+ \\
&&\left( -1\right) ^{\alpha \left( 1_{H};0,0,1,0\right) }B(x_{2}\otimes
x_{1};1_{A},gx_{1}x_{2})1_{A}\otimes gx_{2}\otimes x_{1}+ \\
&&\left( -1\right) ^{\alpha \left( 1_{H};0,0,1,1\right) }B(x_{2}\otimes
x_{1};1_{A},gx_{1}x_{2})1_{A}\otimes g\otimes gx_{1}x_{2}
\end{eqnarray*}

\subsubsection{Case $1_{A}\otimes gx_{1}x_{2}\otimes g$}

this is trivial.

\subsubsection{Case $1_{A}\otimes gx_{1}\otimes x_{2}$}

From the first summand of the left side we get

~%
\begin{eqnarray*}
l_{1} &=&u_{1}=0 \\
l_{2} &=&u_{2}=0 \\
a &=&b_{1}=b_{2}=0, \\
d &=&1,e_{1}=1,e_{2}=0
\end{eqnarray*}%
Since $\alpha \left( x_{2};0,0,0,0\right) =a+b_{1}+b_{2}=0$ we get

\begin{equation*}
B(g\otimes x_{1};1_{A},gx_{1})1_{A}\otimes gx_{1}\otimes x_{2}
\end{equation*}%
From the second summand of the left side we get

~%
\begin{eqnarray*}
l_{1} &=&u_{1}=0 \\
l_{2}+u_{2} &=&1 \\
a &=&b_{1}=0,b_{2}=l_{2}, \\
d &=&e_{1}=1,e_{2}=u_{2}
\end{eqnarray*}%
since $\alpha \left( 1_{H};0,0,0,1\right) \equiv a+b_{1}+b_{2}=0$ and $%
\alpha \left( 1_{H};0,1,0,0\right) \equiv 0$ we obtain%
\begin{equation*}
\left[ B(x_{2}\otimes x_{1};1_{A},gx_{1}x_{2})+B(x_{2}\otimes
x_{1};X_{2},gx_{1})\right] 1_{A}\otimes gx_{1}\otimes x_{2}.
\end{equation*}%
By considering also the right side, we get%
\begin{equation*}
B(g\otimes x_{1};1_{A},gx_{1})+B(x_{2}\otimes
x_{1};1_{A},gx_{1}x_{2})+B(x_{2}\otimes x_{1};X_{2},gx_{1})=0
\end{equation*}%
which holds in view of the form of the elements.

\subsubsection{Case $1_{A}\otimes gx_{1}\otimes x_{1}$}

Nothing from the first summand of the left side. The second summand gives us%
\begin{eqnarray*}
l_{1}+u_{1} &=&1 \\
l_{2} &=&u_{2}=0 \\
a &=&b_{2}=0,b_{1}=l_{1} \\
d &=&1,e_{2}=0,e_{1}-u_{1}=1\Rightarrow e_{1}=1,u_{1}=0,b_{1}=l_{1}=1
\end{eqnarray*}%
since $\alpha \left( 1_{H};1,0,0,0\right) =b_{2}=0,$ we get%
\begin{equation*}
B(x_{2}\otimes x_{1};X_{1},gx_{1})1_{A}\otimes gx_{1}\otimes x_{1}
\end{equation*}%
By considering also the right side, we get%
\begin{equation*}
B(x_{2}\otimes x_{1};X_{1},gx_{1})-B(x_{2}\otimes 1_{H};1_{A},gx_{1})=0
\end{equation*}%
which holds in view of the form of the elements.

\subsubsection{Case $1_{A}\otimes g\otimes gx_{1}x_{2}$}

From the first summand of the left side we get

~%
\begin{eqnarray*}
l_{1}+u_{1} &=&1 \\
l_{2} &=&u_{2}=0 \\
a &=&b_{2}=0,b_{1}=l_{1} \\
d &=&1,e_{2}=0,e_{1}=u_{1}
\end{eqnarray*}%
Since $\alpha \left( x_{2};0,0,1,0\right) $ $\equiv e_{2}=0$ and $\alpha
\left( x_{2};1,0,0,0\right) \equiv a+b_{1}=1$ we get

\begin{equation*}
\left[ B(g\otimes x_{1};1_{A},gx_{1})-B(g\otimes x_{1};X_{1},g)\right]
1_{A}\otimes g\otimes gx_{1}x_{2}.
\end{equation*}%
From the second summand of the left side we get

~%
\begin{eqnarray*}
l_{1}+u_{1} &=&1 \\
l_{2}+u_{2} &=&1 \\
a &=&0,b_{1}=l_{1},b_{2}=l_{2}, \\
d &=&1,e_{1}=u_{1},e_{2}=u_{2}
\end{eqnarray*}%
Since $\alpha \left( 1_{H};0,0,1,1\right) \equiv 1+1\equiv 0$,$\alpha \left(
1_{H};0,1,1,0\right) \equiv e_{2}+a+b_{1}+b_{2}+1\equiv 0,$

$\alpha \left( 1_{H};1,0,0,1\right) \equiv a+b_{1}=1$ and $\alpha \left(
1_{H};1,1,0,0\right) \equiv 1+b_{2}\equiv 0,$ we get

\begin{equation*}
\left[
\begin{array}{c}
B(x_{2}\otimes x_{1};1_{A},gx_{1}x_{2})+B(x_{2}\otimes x_{1};X_{2},gx_{1})+
\\
-B(x_{2}\otimes x_{1};X_{1},gx_{2})+B(x_{2}\otimes x_{1};X_{1}X_{2},g)%
\end{array}%
\right] 1_{A}\otimes g\otimes gx_{1}x_{2}.
\end{equation*}%
By considering also the right side, we get%
\begin{gather*}
B(g\otimes x_{1};1_{A},gx_{1})-B(g\otimes x_{1};X_{1},g)+B(x_{2}\otimes
x_{1};1_{A},gx_{1}x_{2})+ \\
+B(x_{2}\otimes x_{1};X_{2},gx_{1})-B(x_{2}\otimes
x_{1};X_{1},gx_{2})+B(x_{2}\otimes x_{1};X_{1}X_{2},g)=0
\end{gather*}%
which holds in view of the form of the elements.

\subsection{\ $B\left( x_{2}\otimes x_{1};G,1_{H}\right) $}

We deduce that
\begin{equation*}
a=1,b_{1}=b_{2}=0,d=e_{1}=e_{2}=0
\end{equation*}%
and we get%
\begin{equation*}
B\left( x_{2}\otimes x_{1};G,1_{H}\right) G\otimes 1_{H}\otimes g
\end{equation*}

\subsubsection{Case $G\otimes 1_{H}\otimes g$}

this is trivial.

\subsection{ $B\left( x_{2}\otimes x_{1};G,x_{1}x_{2}\right) $}

\begin{equation*}
a=1,b_{1}=b_{2}=0,d=0,e_{1}=e_{2}=1
\end{equation*}%
We deduce that%
\begin{eqnarray*}
&&\left( -1\right) ^{\alpha \left( 1_{H};0,0,0,0\right) }B(x_{2}\otimes
x_{1};G,x_{1}x_{2})G\otimes x_{1}x_{2}\otimes g+ \\
&&\left( -1\right) ^{\alpha \left( 1_{H};0,0,0,1\right) }B(x_{2}\otimes
x_{1};G,x_{1}x_{2})G\otimes x_{1}\otimes x_{2}+ \\
&&\left( -1\right) ^{\alpha \left( 1_{H};0,0,1,0\right) }B(x_{2}\otimes
x_{1};G,x_{1}x_{2})G\otimes x_{2}\otimes x_{1}+ \\
&&\left( -1\right) ^{\alpha \left( 1_{H};0,0,1,1\right) }B(x_{2}\otimes
x_{1};G,x_{1}x_{2})G\otimes 1_{H}\otimes gx_{1}x_{2}+
\end{eqnarray*}

\subsubsection{Case $G\otimes x_{1}x_{2}\otimes g$}

This is trivial.

\subsubsection{Case $G\otimes x_{1}\otimes x_{2}$}

From the first summand of the left side we get

~%
\begin{eqnarray*}
l_{1} &=&u_{1}=0 \\
l_{2} &=&u_{2}=0 \\
a &=&1,b_{1}=b_{2}=0, \\
d &=&0,e_{1}=1,e_{2}=0
\end{eqnarray*}%
Since $\alpha \left( x_{2};0,0,0,0\right) =a+b_{1}+b_{2}=1$ we get%
\begin{equation*}
-B(g\otimes x_{1};G,x_{1})G\otimes x_{1}\otimes x_{2}
\end{equation*}%
From the second summand of the left side we get

~%
\begin{eqnarray*}
l_{1} &=&u_{1}=0 \\
l_{2}+u_{2} &=&1 \\
a &=&1,b_{1}=0,b_{2}=l_{2}, \\
d &=&0,e_{1}=1,e_{2}=u_{2}
\end{eqnarray*}%
since $\alpha \left( 1_{H};0,0,0,1\right) \equiv a+b_{1}+b_{2}=1$ and $%
\alpha \left( 1_{H};0,1,0,0\right) \equiv 0$ we obtain%
\begin{equation*}
\left[ -B(x_{2}\otimes x_{1};G,x_{1}x_{2})+B(x_{2}\otimes x_{1};GX_{2},x_{1})%
\right] G\otimes x_{1}\otimes x_{2}.
\end{equation*}%
Since there is nothing from the right ide, we obtain%
\begin{equation*}
-B(g\otimes x_{1};G,x_{1})-B(x_{2}\otimes x_{1};G,x_{1}x_{2})+B(x_{2}\otimes
x_{1};GX_{2},x_{1})=0
\end{equation*}%
which holds in view of the form of the elements.

\subsubsection{Case $G\otimes x_{2}\otimes x_{1}$}

Nothing from the first summand of the left side. From the second summand of
the left side, we deduce that

\begin{eqnarray*}
l_{1}+u_{1} &=&1 \\
l_{2} &=&u_{2}=0 \\
a &=&1,b_{1}=l_{1},b_{2}=0, \\
d &=&0,e_{1}=u_{1},e_{2}=1.
\end{eqnarray*}%
Since $\alpha \left( 1_{H};0,0,1,0\right) \equiv e_{2}+\left(
a+b_{1}+b_{2}\right) \equiv 0$ and $\alpha \left( 1_{H};1,0,0,0\right)
\equiv b_{2}=0,$ we get%
\begin{equation*}
\left[ B(x_{2}\otimes x_{1};G,x_{1}x_{2})+B(x_{2}\otimes x_{1};GX_{1},x_{2})%
\right] G\otimes x_{2}\otimes x_{1}.
\end{equation*}%
By considering also the right side we get%
\begin{equation*}
B(x_{2}\otimes x_{1};G,x_{1}x_{2})+B(x_{2}\otimes
x_{1};GX_{1},x_{2})-B(x_{2}\otimes 1_{H};G,x_{2})=0
\end{equation*}%
which holds in view of the form of the element.

\subsubsection{Case $G\otimes 1_{H}\otimes gx_{1}x_{2}$}

From the first summand of the left side we get

~%
\begin{eqnarray*}
l_{1}+u_{1} &=&1 \\
l_{2} &=&u_{2}=0 \\
a &=&1,b_{2}=0,b_{1}=l_{1} \\
d &=&0,e_{2}=0,e_{1}=u_{1}
\end{eqnarray*}%
Since $\alpha \left( x_{2};0,0,1,0\right) $ $\equiv e_{2}=0$ and $\alpha
\left( x_{2};1,0,0,0\right) \equiv a+b_{1}=0$ we get

\begin{equation*}
\left[ B(g\otimes x_{1};G,x_{1})+B(g\otimes x_{1};GX_{1},1_{H})\right]
G\otimes 1_{H}\otimes gx_{1}x_{2}.
\end{equation*}%
From the second summand of the left side we get

~%
\begin{eqnarray*}
l_{1}+u_{1} &=&1 \\
l_{2}+u_{2} &=&1 \\
a &=&1,b_{1}=l_{1},b_{2}=l_{2}, \\
d &=&0,e_{1}=u_{1},e_{2}=u_{2}
\end{eqnarray*}%
Since $\alpha \left( 1_{H};0,0,1,1\right) =1+e_{2}\equiv 1+1\equiv 0$, $%
\alpha \left( 1_{H};0,1,1,0\right) \equiv e_{2}+a+b_{1}+b_{2}+1\equiv
1,\alpha \left( 1_{H};1,0,0,1\right) \equiv a+b_{1}=0$ and $\alpha \left(
1_{H};1,1,0,0\right) \equiv 1+b_{2}\equiv 0,$ we get%
\begin{equation*}
\left[
\begin{array}{c}
B(x_{2}\otimes x_{1};G,x_{1}x_{2})+B(x_{2}\otimes x_{1};GX_{2},x_{1})+ \\
-B(x_{2}\otimes x_{1};GX_{1},x_{2})+B(x_{2}\otimes x_{1};GX_{1}X_{2},1_{H})%
\end{array}%
\right] G\otimes 1_{H}\otimes gx_{1}x_{2}.
\end{equation*}%
Since there is nothing from the right side, we get%
\begin{equation*}
\begin{array}{c}
B(g\otimes x_{1};G,x_{1})+B(g\otimes x_{1};GX_{1},1_{H}) \\
+B(x_{2}\otimes x_{1};G,x_{1}x_{2})-B(x_{2}\otimes x_{1};GX_{2},x_{1})+ \\
+B(x_{2}\otimes x_{1};GX_{1},x_{2})+B(x_{2}\otimes x_{1};GX_{1}X_{2},1_{H})%
\end{array}%
=0
\end{equation*}%
which holds in view of the form of the elements.

\subsection{$B\left( x_{2}\otimes x_{1};G,gx_{1}\right) $}

We obtain%
\begin{eqnarray*}
&&\left( -1\right) ^{\alpha \left( 1_{H};0,0,0,0\right) }B(x_{2}\otimes
x_{1};G,gx_{1})G\otimes gx_{1}\otimes g+ \\
&&\left( -1\right) ^{\alpha \left( 1_{H};0,0,1,0\right) }B(x_{2}\otimes
x_{1};G,gx_{1})G\otimes g\otimes x_{1}
\end{eqnarray*}

\subsubsection{Case $G\otimes gx_{1}\otimes g$}

This is trivial.

\subsubsection{Case $G\otimes g\otimes x_{1}$}

Nothing from the first summand of the left side. From the second summand of
the left side, we deduce that

\begin{eqnarray*}
l_{1}+u_{1} &=&1 \\
l_{2} &=&u_{2}=0 \\
a &=&1,b_{1}=l_{1},b_{2}=0, \\
d &=&1,e_{1}=u_{1},e_{2}=0.
\end{eqnarray*}%
Since $\alpha \left( 1_{H};0,0,1,0\right) \equiv e_{2}+\left(
a+b_{1}+b_{2}\right) \equiv 1$ and $\alpha \left( 1_{H};1,0,0,0\right)
\equiv b_{2}=0,$ we get%
\begin{equation*}
\left[ -B(x_{2}\otimes x_{1};G,gx_{1})+B(x_{2}\otimes x_{1};GX_{1},g)\right]
G\otimes g\otimes x_{1}.
\end{equation*}%
By considering also the right side we get%
\begin{equation*}
-B(x_{2}\otimes x_{1};G,gx_{1})+B(x_{2}\otimes
x_{1};GX_{1},g)-B(x_{2}\otimes 1_{H};G,g)=0
\end{equation*}%
which holds in view of the form of the elements.

\subsection{$B\left( x_{2}\otimes x_{1};G,gx_{2}\right) $}

We deduce that%
\begin{eqnarray*}
&&\left( -1\right) ^{\alpha \left( 1_{H};0,0,0,0\right) }B(x_{2}\otimes
x_{1};G,gx_{2})G\otimes gx_{2}\otimes g+ \\
&&\left( -1\right) ^{\alpha \left( 1_{H};0,0,0,u_{2}\right) }B(x_{2}\otimes
x_{1};G,gx_{2})G\otimes g\otimes x_{2}+
\end{eqnarray*}

\subsubsection{Case $G\otimes gx_{2}\otimes g$}

This is trivial.

\subsubsection{Case $G\otimes g\otimes x_{2}$}

From the first summand of the left side we get

~%
\begin{eqnarray*}
l_{1} &=&u_{1}=0 \\
l_{2} &=&u_{2}=0 \\
a &=&1,b_{1}=b_{2}=0, \\
d &=&1,e_{1}=0,e_{2}=0
\end{eqnarray*}%
Since $\alpha \left( x_{2};0,0,0,0\right) =a+b_{1}+b_{2}=1$ we get%
\begin{equation*}
-B(g\otimes x_{1};G,g)G\otimes g\otimes x_{2}
\end{equation*}%
From the second summand of the left side we get

~%
\begin{eqnarray*}
l_{1} &=&u_{1}=0 \\
l_{2}+u_{2} &=&1 \\
a &=&1,b_{1}=0,b_{2}=l_{2}, \\
d &=&1,e_{1}=0,e_{2}=u_{2}
\end{eqnarray*}%
since $\alpha \left( 1_{H};0,0,0,1\right) \equiv a+b_{1}+b_{2}=1$ and $%
\alpha \left( 1_{H};0,1,0,0\right) \equiv 0$ we obtain%
\begin{equation*}
\left[ -B(x_{2}\otimes x_{1};G,gx_{2})+B(x_{2}\otimes x_{1};GX_{2},g)\right]
G\otimes g\otimes x_{2}.
\end{equation*}%
Since there is nothing from the right side, we get%
\begin{equation*}
-B(g\otimes x_{1};G,g)-B(x_{2}\otimes x_{1};G,gx_{2})+B(x_{2}\otimes
x_{1};GX_{2},g)=0
\end{equation*}%
which holds in view of the form of the elements.

\subsection{$B\left( x_{2}\otimes x_{1};X_{1},1_{H}\right) $}

We deduce that
\begin{equation*}
a=0,b_{1}=1,b_{2}=0,d=e_{1}=e_{2}=0
\end{equation*}%
and we get%
\begin{eqnarray*}
&&\left( -1\right) ^{\alpha \left( 1_{H};0,0,0,0\right) }B(x_{2}\otimes
x_{1};X_{1},1_{H})X_{1}\otimes 1_{H}\otimes g+ \\
&&\left( -1\right) ^{\alpha \left( 1_{H};1,0,0,0\right) }B(x_{2}\otimes
x_{1};X_{1},1_{H})1_{A}\otimes 1_{H}\otimes x_{1}
\end{eqnarray*}

\subsubsection{Case $X_{1}\otimes 1_{H}\otimes g$}

It is trivial.

\subsubsection{Case $1_{A}\otimes 1_{H}\otimes x_{1}$}

Nothing from the first summand of the left side. From the second summand of
the left side, we deduce that

\begin{eqnarray*}
l_{1}+u_{1} &=&1 \\
l_{2} &=&u_{2}=0 \\
a &=&0,b_{1}=l_{1},b_{2}=0, \\
d &=&0,e_{1}=u_{1},e_{2}=0.
\end{eqnarray*}%
Since $\alpha \left( 1_{H};0,0,1,0\right) \equiv e_{2}+\left(
a+b_{1}+b_{2}\right) \equiv 0$ and $\alpha \left( 1_{H};1,0,0,0\right)
\equiv b_{2}=0,$ we get%
\begin{equation*}
\left[ B(x_{2}\otimes x_{1};1_{A},x_{1})+B(x_{2}\otimes x_{1};X_{1},1_{H})%
\right] 1_{A}\otimes 1_{H}\otimes x_{1}.
\end{equation*}%
By considering also the right side we get%
\begin{equation*}
B(x_{2}\otimes x_{1};1_{A},x_{1})+B(x_{2}\otimes
x_{1};X_{1},1_{H})-B(x_{2}\otimes 1_{H};1_{A},1_{H})=0
\end{equation*}%
which holds in view of the form of the elements.

\subsection{$B\left( x_{2}\otimes x_{1};X_{1},x_{1}x_{2}\right) $}

We deduce that%
\begin{eqnarray*}
a &=&0,b_{1}=1,b_{2}=0 \\
d &=&0,e_{1}=e_{2}=1
\end{eqnarray*}%
and we get%
\begin{gather*}
\left( -1\right) ^{\alpha \left( 1_{H};0,0,0,0\right) }B(x_{2}\otimes
x_{1};X_{1},x_{1}x_{2})X_{1}\otimes x_{1}x_{2}\otimes g+ \\
\left( -1\right) ^{\alpha \left( 1_{H};1,0,0,0\right) }B(x_{2}\otimes
x_{1};X_{1},x_{1}x_{2})1_{A}\otimes x_{1}x_{2}\otimes x_{1} \\
\left( -1\right) ^{\alpha \left( 1_{H};0,0,1,0\right) }B(x_{2}\otimes
x_{1};X_{1},x_{1}x_{2})X_{1}\otimes x_{2}\otimes x_{1}+ \\
\left( -1\right) ^{\alpha \left( 1_{H};1,0,1,0\right) }B(x_{2}\otimes
x_{1};X_{1},x_{1}x_{2})X_{1}^{1-l_{1}}\otimes x_{2}\otimes
g^{l_{1}}x_{1}^{1+1}=0 \\
\left( -1\right) ^{\alpha \left( 1_{H};0,0,0,1\right) }B(x_{2}\otimes
x_{1};X_{1},x_{1}x_{2})X_{1}\otimes x_{1}\otimes x_{2}+ \\
\left( -1\right) ^{\alpha \left( 1_{H};1,0,0,1\right) }B(x_{2}\otimes
x_{1};X_{1},x_{1}x_{2})1_{A}\otimes x_{1}\otimes gx_{1}x_{2}+ \\
\left( -1\right) ^{\alpha \left( 1_{H};0,0,1,1\right) }B(x_{2}\otimes
x_{1};X_{1},x_{1}x_{2})X_{1}\otimes 1_{H}\otimes gx_{1}x_{2}+ \\
\left( -1\right) ^{\alpha \left( 1_{H};1,0,1,1\right) }B(x_{2}\otimes
x_{1};X_{1},x_{1}x_{2})X_{1}^{1-l_{1}}\otimes 1_{H}\otimes
g^{l_{1}+1}x_{1}^{1+1}x_{2}=0
\end{gather*}

\subsubsection{Case $X_{1}\otimes x_{1}x_{2}\otimes g$}

This is trivial.

\subsubsection{Case $1_{A}\otimes x_{1}x_{2}\otimes x_{1}$}

Nothing from the first summand of the left side. From the second summand of
the left side, we deduce that

\begin{eqnarray*}
l_{1}+u_{1} &=&1 \\
l_{2} &=&u_{2}=0 \\
a &=&0,b_{1}=l_{1},b_{2}=0, \\
d &=&0,e_{1}-u_{1}=1\Rightarrow e_{1}=1,u_{1}=0,b_{1}=l_{1}=1,e_{2}=1.
\end{eqnarray*}%
Since $\alpha \left( 1_{H};1,0,0,0\right) \equiv b_{2}=0,$we get%
\begin{equation*}
B(x_{2}\otimes x_{1};X_{1},x_{1}x_{2})1_{A}\otimes x_{1}x_{2}\otimes x_{1}.
\end{equation*}%
By considering also the right side, we obtain%
\begin{equation*}
B(x_{2}\otimes x_{1};X_{1},x_{1}x_{2})-B(x_{2}\otimes
1_{H};1_{A},x_{1}x_{2})=0
\end{equation*}%
which holds in view of the form of the elements.

\subsubsection{Case $X_{1}\otimes x_{2}\otimes x_{1}$}

Nothing from the first summand of the left side. From the second summand of
the left side, we deduce that

\begin{eqnarray*}
l_{1}+u_{1} &=&1 \\
l_{2} &=&u_{2}=0 \\
a &=&0,b_{1}-l_{1}=1\Rightarrow b_{1}=1,l_{1}=0,u_{1}=1,b_{2}=0, \\
d &=&0,e_{1}=u_{1}=1,e_{2}=1.
\end{eqnarray*}%
Since $\alpha \left( 1_{H};0,0,1,0\right) =e_{2}+\left( a+b_{1}+b_{2}\right)
\equiv 0,$ we get%
\begin{equation*}
B(x_{2}\otimes x_{1};X_{1},x_{1}x_{2})X_{1}\otimes x_{2}\otimes x_{1}.
\end{equation*}%
By considering also the right side, we get%
\begin{equation*}
B(x_{2}\otimes x_{1};X_{1},x_{1}x_{2})-B(x_{2}\otimes 1_{H};X_{1},x_{2})=0
\end{equation*}%
which holds in view of the form of the elements.

\subsubsection{Case $X_{1}\otimes x_{1}\otimes x_{2}$}

From the first summand of the left side we get

\begin{eqnarray*}
l_{1} &=&u_{1}=0 \\
l_{2} &=&u_{2}=0 \\
a &=&0,b_{1}=1,b_{2}=0, \\
d &=&0,e_{1}=1,e_{2}=0
\end{eqnarray*}%
Since $\alpha \left( x_{2};0,0,0,0\right) =a+b_{1}+b_{2}=1$ we get%
\begin{equation*}
-B(g\otimes x_{1};X_{1},x_{1})X_{1}\otimes x_{1}\otimes x_{2}
\end{equation*}%
From the second summand of the left side we get

~%
\begin{eqnarray*}
l_{1} &=&u_{1}=0 \\
l_{2}+u_{2} &=&1 \\
a &=&0,b_{1}=1,b_{2}=l_{2}, \\
d &=&0,e_{1}=1,e_{2}=u_{2}
\end{eqnarray*}%
since $\alpha \left( 1_{H};0,0,0,1\right) \equiv a+b_{1}+b_{2}=1$ and $%
\alpha \left( 1_{H};0,1,0,0\right) \equiv 0$ we obtain%
\begin{equation*}
\left[ -B(x_{2}\otimes x_{1};X_{1},x_{1}x_{2})+B(x_{2}\otimes
x_{1};X_{1}X_{2},x_{1})\right] X_{1}\otimes x_{1}\otimes x_{2}.
\end{equation*}%
Since there is nothing in the right side, we obtain

\begin{equation*}
-B(g\otimes x_{1};X_{1},x_{1})-B(x_{2}\otimes
x_{1};X_{1},x_{1}x_{2})+B(x_{2}\otimes x_{1};X_{1}X_{2},x_{1})=0
\end{equation*}%
which holds in view of the form of the elements.

\subsubsection{Case $1_{A}\otimes x_{1}\otimes gx_{1}x_{2}$}

From the first summand of the left side we get

\begin{eqnarray*}
l_{1}+u_{1} &=&1 \\
l_{2} &=&u_{2}=0 \\
a &=&0,b_{2}=0,b_{1}=l_{1} \\
d &=&0,e_{2}=0,e_{1}-u_{1}=1\Rightarrow e_{1}=1,u_{1}=0,b_{1}=l_{1}=1.
\end{eqnarray*}%
Since $\alpha \left( x_{2};1,0,0,0\right) \equiv a+b_{1}\equiv 1$ we get

\begin{equation*}
\left[ -B(g\otimes x_{1};X_{1},x_{1})\right] 1_{A}\otimes x_{1}\otimes
gx_{1}x_{2}.
\end{equation*}%
From the second summand of the left side we get

\begin{eqnarray*}
l_{1}+u_{1} &=&1 \\
l_{2}+u_{2} &=&1 \\
a &=&0,b_{1}=l_{1},b_{2}=l_{2}, \\
d &=&0,e_{1}-u_{1}=1\Rightarrow e_{1}=1,u_{1}=0,b_{1}=l_{1}=1,e_{2}=u_{2}
\end{eqnarray*}%
Since $\alpha \left( 1_{H};1,0,0,1\right) \equiv a+b_{1}=1$ and $\alpha
\left( 1_{H};1,1,0,0\right) \equiv 1+b_{2}\equiv 0,$ we get%
\begin{equation*}
\left[ -B(x_{2}\otimes x_{1};X_{1},x_{1}x_{2})+B(x_{2}\otimes
x_{1};X_{1}X_{2},x_{1})\right] 1_{A}\otimes x_{1}\otimes gx_{1}x_{2}.
\end{equation*}%
Since there is nothing from the right side, we obtain%
\begin{equation*}
-B(g\otimes x_{1};X_{1},x_{1})-B(x_{2}\otimes
x_{1};X_{1},x_{1}x_{2})+B(x_{2}\otimes x_{1};X_{1}X_{2},x_{1})=0
\end{equation*}%
which we just got above.

\subsubsection{Case $X_{1}\otimes 1_{H}\otimes gx_{1}x_{2}$}

From the first summand of the left side we get

~%
\begin{eqnarray*}
l_{1}+u_{1} &=&1 \\
l_{2} &=&u_{2}=0 \\
a &=&0,b_{2}=0,b_{1}-l_{1}=1\Rightarrow b_{1}=1,l_{1}=0,u_{1}=1 \\
d &=&0,e_{2}=0,e_{1}=u_{1}=1.
\end{eqnarray*}%
Since $\alpha \left( x_{2};0,0,1,0\right) \equiv e_{2}=0,$ we get

\begin{equation*}
\left[ B(g\otimes x_{1};X_{1},x_{1})\right] X_{1}\otimes 1_{H}\otimes
gx_{1}x_{2}.
\end{equation*}%
From the second summand of the left side we get

~%
\begin{eqnarray*}
l_{1}+u_{1} &=&1 \\
l_{2}+u_{2} &=&1 \\
a &=&0,b_{1}-l_{1}=1\Rightarrow b_{1}=1,l_{1}=0,u_{1}=1,b_{2}=l_{2}, \\
d &=&0,e_{1}=u_{1}=1,e_{2}=u_{2}
\end{eqnarray*}%
Since $\alpha \left( 1_{H};0,0,1,1\right) \equiv $ $1+e_{2}\equiv 0$ and $%
\alpha \left( 1_{H};0,1,1,0\right) \equiv e_{2}+a+b_{1}+b_{2}+1\equiv 1,$ we
get%
\begin{equation*}
\left[ B(x_{2}\otimes x_{1};X_{1},x_{1}x_{2})-B(x_{2}\otimes
x_{1};X_{1}X_{2},x_{1})\right] X_{1}\otimes 1_{H}\otimes gx_{1}x_{2}.
\end{equation*}%
Since there is nothing from the right side, we obtain%
\begin{equation*}
+B(g\otimes x_{1};X_{1},x_{1})+B(x_{2}\otimes
x_{1};X_{1},x_{1}x_{2})-B(x_{2}\otimes x_{1};X_{1}X_{2},x_{1})=0
\end{equation*}%
which we just got above.

\subsection{$B\left( x_{2}\otimes x_{1};X_{1},gx_{1}\right) $}

We deduce that%
\begin{eqnarray*}
a &=&0,b_{1}=1,b_{2}=0 \\
d &=&1,e_{1}=1,e_{2}=0
\end{eqnarray*}%
and we get%
\begin{eqnarray*}
&&\left( -1\right) ^{\alpha \left( 1_{H};0,0,1,0\right) }B(x_{2}\otimes
x_{1};X_{1},gx_{1})X_{1}\otimes g\otimes x_{1}+ \\
&&+\left( -1\right) ^{\alpha \left( 1_{H};1,0,0,0\right) }B(x_{2}\otimes
x_{1};X_{1},gx_{1})1_{A}\otimes gx_{1}\otimes x_{1}
\end{eqnarray*}

\subsubsection{Case $X_{1}\otimes g\otimes x_{1}$}

Nothing from the first summand of the left side. From the second summand of
the left side, we deduce that

\begin{eqnarray*}
l_{1}+u_{1} &=&1 \\
l_{2} &=&u_{2}=0 \\
a &=&0,b_{1}-l_{1}=1\Rightarrow b_{1}=1,l_{1}=0,u_{1}=1,b_{2}=0, \\
d &=&1,e_{1}=u_{1}=1,e_{2}=0.
\end{eqnarray*}%
Since $\alpha \left( 1_{H};0,0,1,0\right) \equiv e_{2}+\left(
a+b_{1}+b_{2}\right) \equiv 1$ $,$ we get%
\begin{equation*}
-B(x_{2}\otimes x_{1};X_{1},gx_{1})X_{1}\otimes g\otimes x_{1}.
\end{equation*}%
By considering also the right side, we obtain%
\begin{equation*}
-B(x_{2}\otimes x_{1};X_{1},gx_{1})-B(x_{2}\otimes 1_{H};X_{1},g)=0
\end{equation*}%
which holds in view of the form of the elements.

\subsubsection{Case $1_{A}\otimes gx_{1}\otimes x_{1}$}

Nothing from the first summand of the left side. From the second summand of
the left side, we deduce that

\begin{eqnarray*}
l_{1}+u_{1} &=&1 \\
l_{2} &=&u_{2}=0 \\
a &=&0,b_{1}=l_{1},b_{2}=0, \\
d &=&1,e_{1}-u_{1}=1\Rightarrow e_{1}=1,u_{1}=0,b_{1}=l_{1}=1,e_{2}=0.
\end{eqnarray*}%
Since $\alpha \left( 1_{H};1,0,0,0\right) \equiv b_{2}=0,$ we obtain%
\begin{equation*}
B(x_{2}\otimes x_{1};X_{1},gx_{1})1_{A}\otimes gx_{1}\otimes x_{1}.
\end{equation*}%
By considering also the right side, we obtain%
\begin{equation*}
B(x_{2}\otimes x_{1};X_{1},gx_{1})-B(x_{2}\otimes 1_{H};1_{A},gx_{1})=0
\end{equation*}%
which holds in view of the form of the elements.

\subsection{$B\left( x_{2}\otimes x_{1};X_{1},gx_{2}\right) $}

We deduce that%
\begin{eqnarray*}
a &=&0,b_{1}=1,b_{2}=0 \\
d &=&1,e_{1}=0,e_{2}=1
\end{eqnarray*}%
and we get%
\begin{eqnarray*}
&&\left( -1\right) ^{\alpha \left( 1_{H};0,0,0,1\right) }B(x_{2}\otimes
x_{1};X_{1},gx_{2})X_{1}\otimes g\otimes x_{2}+ \\
&&+\left( -1\right) ^{\alpha \left( 1_{H};1,0,0,0\right) }B(x_{2}\otimes
x_{1};X_{1},gx_{2})1_{A}\otimes gx_{2}\otimes x_{1}
\end{eqnarray*}

\subsubsection{Case $X_{1}\otimes g\otimes x_{2}$}

From the first summand of the left side we get

\begin{eqnarray*}
l_{1} &=&u_{1}=0 \\
l_{2} &=&u_{2}=0 \\
a &=&0,b_{1}=1,b_{2}=0, \\
d &=&1,e_{1}=0,e_{2}=0
\end{eqnarray*}%
Since $\alpha \left( x_{2};0,0,0,0\right) =a+b_{1}+b_{2}=1$ we get%
\begin{equation*}
-B(g\otimes x_{1};X_{1},g)X_{1}\otimes g\otimes x_{2}
\end{equation*}%
From the second summand of the left side we get

~%
\begin{eqnarray*}
l_{1} &=&u_{1}=0 \\
l_{2}+u_{2} &=&1 \\
a &=&0,b_{1}=1,b_{2}=l_{2}, \\
d &=&1,e_{1}=0,e_{2}=u_{2}
\end{eqnarray*}%
since $\alpha \left( 1_{H};0,0,0,1\right) \equiv a+b_{1}+b_{2}=1$ and $%
\alpha \left( 1_{H};0,1,0,0\right) \equiv 0$ we obtain%
\begin{equation*}
\left[ -B(x_{2}\otimes x_{1};X_{1},gx_{2})+B(x_{2}\otimes x_{1};X_{1}X_{2},g)%
\right] X_{1}\otimes g\otimes x_{2}.
\end{equation*}%
Since there is nothing in the right side, we obtain%
\begin{equation*}
-B(g\otimes x_{1};X_{1},g)-B(x_{2}\otimes x_{1};X_{1},gx_{2})+B(x_{2}\otimes
x_{1};X_{1}X_{2},g)=0
\end{equation*}%
which holds in view of the form of the elements.

\subsubsection{Case $1_{A}\otimes gx_{2}\otimes x_{1}$}

Nothing from the first summand of the left side. From the second summand of
the left side, we deduce that

\begin{eqnarray*}
l_{1}+u_{1} &=&1 \\
l_{2} &=&u_{2}=0 \\
a &=&0,b_{1}=l_{1},b_{2}=0, \\
d &=&1,e_{2}=1\Rightarrow e_{1}-u_{1}=0,.
\end{eqnarray*}%
Since $\alpha \left( 1_{H};0,0,1,0\right) =e_{2}+\left( a+b_{1}+b_{2}\right)
=1$ and $\alpha \left( 1_{H};1,0,0,0\right) \equiv b_{2}=0$ we obtain%
\begin{equation*}
\left[ -B(x_{2}\otimes x_{1};1_{A},gx_{1}x_{2})+B(x_{2}\otimes
x_{1};X_{1},gx_{2})\right] 1_{A}\otimes gx_{2}\otimes x_{1}.
\end{equation*}%
By considering also the right side, we obtain%
\begin{equation*}
-B(x_{2}\otimes x_{1};1_{A},gx_{1}x_{2})+B(x_{2}\otimes
x_{1};X_{1},gx_{2})-B(x_{2}\otimes 1_{H};1_{A},gx_{2})=0
\end{equation*}%
which holds in view of the form of the elements.

\subsection{$B\left( x_{2}\otimes x_{1};X_{2},1_{H}\right) $}

We deduce that
\begin{equation*}
a=b_{1}=0,b_{2}=1,d=e_{1}=e_{2}=0
\end{equation*}%
and we get%
\begin{eqnarray*}
&&\left( -1\right) ^{\alpha \left( 1_{H};0,0,0,0\right) }B(x_{2}\otimes
x_{1};X_{2},1_{H})X_{2}\otimes 1_{H}\otimes g+ \\
&&\left( -1\right) ^{\alpha \left( 1_{H};0,1,0,0\right) }B(x_{2}\otimes
x_{1};X_{2},1_{H})1_{A}\otimes 1_{H}\otimes x_{2}
\end{eqnarray*}

\subsubsection{Case $X_{2}\otimes 1_{H}\otimes g$}

This is trivial.

\subsubsection{Case $1_{A}\otimes 1_{H}\otimes x_{2}$}

From the first summand of the left side we get

\begin{eqnarray*}
l_{1} &=&u_{1}=0 \\
l_{2} &=&u_{2}=0 \\
a &=&0,b_{1}=0,b_{2}=0, \\
d &=&0,e_{1}=0,e_{2}=0
\end{eqnarray*}%
Since $\alpha \left( x_{2};0,0,0,0\right) =a+b_{1}+b_{2}=0$ we get%
\begin{equation*}
B(g\otimes x_{1};1_{A},1_{H})1_{A}\otimes 1_{H}\otimes x_{2}
\end{equation*}%
From the second summand of the left side we get

~%
\begin{eqnarray*}
l_{1} &=&u_{1}=0 \\
l_{2}+u_{2} &=&1 \\
a &=&0,b_{1}=0,b_{2}=l_{2}, \\
d &=&0,e_{1}=0,e_{2}=u_{2}
\end{eqnarray*}%
since $\alpha \left( 1_{H};0,0,0,1\right) \equiv a+b_{1}+b_{2}=0$ and $%
\alpha \left( 1_{H};0,1,0,0\right) \equiv 0$ we obtain

\begin{equation*}
\left[ B(x_{2}\otimes x_{1};1_{A},x_{2})+B(x_{2}\otimes x_{1};X_{2},1_{H})%
\right] 1_{A}\otimes 1_{H}\otimes x_{2}.
\end{equation*}%
Since there is nothing from the right side, we obtain%
\begin{equation*}
B(g\otimes x_{1};1_{A},1_{H})+B(x_{2}\otimes
x_{1};1_{A},x_{2})+B(x_{2}\otimes x_{1};X_{2},1_{H})=0
\end{equation*}%
which holds in view of the form of the elements.

\subsection{$B\left( x_{2}\otimes x_{1};X_{2},x_{1}x_{2}\right) $}

We deduce that%
\begin{eqnarray*}
a &=&b_{1}=0,b_{2}=1 \\
d &=&0,e_{1}=e_{2}=1
\end{eqnarray*}%
and we get%
\begin{equation*}
\left( -1\right) ^{\alpha \left( 1_{H};0,l_{2},u_{1},u_{2}\right)
}X_{2}^{1-l_{2}}\otimes x_{1}^{1-u_{1}}x_{2}^{1-u_{2}}\otimes
g^{1+l_{2}+u_{1}+u_{2}}x_{1}^{u_{1}}x_{2}^{l_{2}+u_{2}}
\end{equation*}%
\begin{gather*}
\left( -1\right) ^{\alpha \left( 1_{H};0,0,0,0\right) }B\left( x_{2}\otimes
x_{1};X_{2},x_{1}x_{2}\right) X_{2}\otimes x_{1}x_{2}\otimes g+ \\
\left( -1\right) ^{\alpha \left( 1_{H};0,1,0,0\right) }B\left( x_{2}\otimes
x_{1};X_{2},x_{1}x_{2}\right) 1_{A}\otimes x_{1}x_{2}\otimes x_{2} \\
\left( -1\right) ^{\alpha \left( 1_{H};0,0,1,0\right) }B\left( x_{2}\otimes
x_{1};X_{2},x_{1}x_{2}\right) X_{2}\otimes x_{2}\otimes x_{1} \\
\left( -1\right) ^{\alpha \left( 1_{H};0,1,1,0\right) }B\left( x_{2}\otimes
x_{1};X_{2},x_{1}x_{2}\right) 1_{A}\otimes x_{2}\otimes gx_{1}x_{2} \\
\left( -1\right) ^{\alpha \left( 1_{H};0,0,0,1\right) }B\left( x_{2}\otimes
x_{1};X_{2},x_{1}x_{2}\right) X_{2}\otimes x_{1}\otimes x_{2} \\
\left( -1\right) ^{\alpha \left( 1_{H};0,1,0,1\right) }B\left( x_{2}\otimes
x_{1};X_{2},x_{1}x_{2}\right) X_{2}^{1-l_{2}}\otimes x_{1}\otimes
g^{l_{2}+1}x_{2}^{1+1}=0 \\
\left( -1\right) ^{\alpha \left( 1_{H};0,0,1,1\right) }B\left( x_{2}\otimes
x_{1};X_{2},x_{1}x_{2}\right) X_{2}\otimes 1_{H}\otimes gx_{1}x_{2} \\
\left( -1\right) ^{\alpha \left( 1_{H};0,1,1,1\right) }B\left( x_{2}\otimes
x_{1};X_{2},x_{1}x_{2}\right) X_{2}^{1-l_{2}}\otimes 1_{H}\otimes
g^{l_{2}+1}x_{1}x_{2}^{1+1}=0
\end{gather*}

\subsubsection{Case $X_{2}\otimes x_{1}x_{2}\otimes g$}

This is trivial.

\subsubsection{Case $1_{A}\otimes x_{1}x_{2}\otimes x_{2}$}

From the first summand of the left side we get

\begin{eqnarray*}
l_{1} &=&u_{1}=0 \\
l_{2} &=&u_{2}=0 \\
a &=&b_{1}=b_{2}=0, \\
d &=&0,e_{1}=e_{2}=1
\end{eqnarray*}%
Since $\alpha \left( x_{2};0,0,0,0\right) =a+b_{1}+b_{2}=0$ we get%
\begin{equation*}
B(g\otimes x_{1};1_{A},x_{1}x_{2})1_{A}\otimes x_{1}x_{2}\otimes x_{2}
\end{equation*}%
From the second summand of the left side we get

~%
\begin{eqnarray*}
l_{1} &=&u_{1}=0 \\
l_{2}+u_{2} &=&1 \\
a &=&b_{1}=0,b_{2}=l_{2}, \\
d &=&0,e_{1}=1,e_{2}-u_{2}=1\Rightarrow e_{2}=1,u_{2}=0,b_{2}=l_{2}=1
\end{eqnarray*}%
since $\alpha \left( 1_{H};0,1,0,0\right) \equiv 0$ we obtain%
\begin{equation*}
B(x_{2}\otimes x_{1};X_{2},x_{1}x_{2})1_{A}\otimes x_{1}x_{2}\otimes x_{2}.
\end{equation*}%
Since there is nothing from the right side, we get%
\begin{equation*}
B(g\otimes x_{1};1_{A},x_{1}x_{2})+B(x_{2}\otimes x_{1};X_{2},x_{1}x_{2})=0
\end{equation*}%
which holds in view of the form of the elements.

\subsubsection{Case $X_{2}\otimes x_{2}\otimes x_{1}$}

Nothing from the first summand of the left side. From the second summand of
the left side, we deduce that

\begin{eqnarray*}
l_{1}+u_{1} &=&1 \\
l_{2} &=&u_{2}=0 \\
a &=&0,b_{1}=l_{1},b_{2}=1, \\
d &=&0,e_{1}=u_{1},e_{2}=1.
\end{eqnarray*}%
Since $\alpha \left( 1_{H};0,0,1,0\right) =e_{2}+\left( a+b_{1}+b_{2}\right)
\equiv 0,$and $\alpha \left( 1_{H};1,0,0,0\right) =b_{2}=1$ we get%
\begin{equation*}
\left[ B(x_{2}\otimes x_{1};X_{2},x_{1}x_{2})-B(x_{2}\otimes
x_{1};X_{1}X_{2},x_{2})\right] X_{2}\otimes x_{2}\otimes x_{1}.
\end{equation*}%
By considering also the right side, we get%
\begin{equation*}
B(x_{2}\otimes x_{1};X_{2},x_{1}x_{2})-B(x_{2}\otimes
x_{1};X_{1}X_{2},x_{2})-B(x_{2}\otimes 1_{H};X_{2},x_{2})=0
\end{equation*}%
which holds in view of the form of the elements.

\subsubsection{Case $1_{A}\otimes x_{2}\otimes gx_{1}x_{2}$}

From the first summand of the left side we get

\begin{eqnarray*}
l_{1}+u_{1} &=&1 \\
l_{2} &=&u_{2}=0 \\
a &=&0,b_{2}=0,b_{1}=l_{1} \\
d &=&0,e_{2}=1,e_{1}=u_{1}.
\end{eqnarray*}%
Since $\alpha \left( x_{2};0,0,1,0\right) =e_{2}=1$ $\ $and $\alpha \left(
x_{2};1,0,0,0\right) \equiv a+b_{1}\equiv 1$ we get

\begin{equation*}
\left[ -B(g\otimes x_{1};1_{A},x_{1}x_{2})-B(g\otimes x_{1};X_{1},x_{2})%
\right] 1_{A}\otimes x_{2}\otimes gx_{1}x_{2}.
\end{equation*}%
From the second summand of the left side we get

\begin{eqnarray*}
l_{1}+u_{1} &=&1 \\
l_{2}+u_{2} &=&1 \\
a &=&0,b_{1}=l_{1},b_{2}=l_{2}, \\
d &=&0,e_{1}=u_{1},e_{2}-u_{2}=1\Rightarrow e_{2}=1,u_{2}=0,b_{2}=l_{2}=1
\end{eqnarray*}%
Since $\alpha \left( 1_{H};0,1,1,0\right) \equiv e_{2}+a+b_{1}+b_{2}+1\equiv
1$ and $\alpha \left( 1_{H};1,1,0,0\right) \equiv 1+b_{2}\equiv 0,$ we get%
\begin{equation*}
\left[ -B(x_{2}\otimes x_{1};X_{2},x_{1}x_{2})+B(x_{2}\otimes
x_{1};X_{1}X_{2},x_{2})\right] 1_{A}\otimes x_{2}\otimes gx_{1}x_{2}.
\end{equation*}%
Since there is nothing from the right side, we obtain%
\begin{gather*}
-B(g\otimes x_{1};1_{A},x_{1}x_{2})-B(g\otimes x_{1};X_{1},x_{2})+ \\
-B(x_{2}\otimes x_{1};X_{2},x_{1}x_{2})+B(x_{2}\otimes
x_{1};X_{1}X_{2},x_{2})=0
\end{gather*}%
which holds in view of the form of the elements.

\subsubsection{Case $X_{2}\otimes x_{1}\otimes x_{2}$}

From the first summand of the left side we get

\begin{eqnarray*}
l_{1} &=&u_{1}=0 \\
l_{2} &=&u_{2}=0 \\
a &=&0,b_{1}=0,b_{2}=1, \\
d &=&e_{2}=0,e_{1}=1,
\end{eqnarray*}%
Since $\alpha \left( x_{2};0,0,0,0\right) =a+b_{1}+b_{2}=1$ we get%
\begin{equation*}
-B(g\otimes x_{1};X_{2},x_{1})X_{2}\otimes x_{1}\otimes x_{2}
\end{equation*}%
From the second summand of the left side we get

~%
\begin{eqnarray*}
l_{1} &=&u_{1}=0 \\
l_{2}+u_{2} &=&1 \\
a &=&b_{1}=0,b_{2}-l_{2}=1\Rightarrow b_{2}=1,l_{2}=0,u_{2}=1, \\
d &=&0,e_{1}=1,e_{2}=u_{2}=1
\end{eqnarray*}%
since $\alpha \left( 1_{H};0,0,0,1\right) \equiv a+b_{1}+b_{2}=1$ we obtain%
\begin{equation*}
-B(x_{2}\otimes x_{1};X_{2},x_{1}x_{2})X_{2}\otimes x_{1}\otimes x_{2}.
\end{equation*}%
Since there is nothing in the right side, we obtain%
\begin{equation*}
-B(g\otimes x_{1};X_{2},x_{1})-B(x_{2}\otimes x_{1};X_{2},x_{1}x_{2})=0
\end{equation*}%
which holds in view of the form of the elements.

\subsubsection{Case $X_{2}\otimes 1_{H}\otimes gx_{1}x_{2}$}

From the first summand of the left side we get

\begin{eqnarray*}
l_{1}+u_{1} &=&1 \\
l_{2} &=&u_{2}=0 \\
a &=&0,b_{2}=1,b_{1}=l_{1} \\
d &=&0,e_{2}=0,e_{1}=u_{1}.
\end{eqnarray*}%
Since $\alpha \left( x_{2};0,0,1,0\right) =e_{2}=0$ $\ $and $\alpha \left(
x_{2};1,0,0,0\right) \equiv a+b_{1}\equiv 1$ we get

\begin{equation*}
\left[ B(g\otimes x_{1};X_{2},x_{1})-B(g\otimes x_{1};X_{1}X_{2},1_{H})%
\right] X_{2}\otimes 1_{H}\otimes gx_{1}x_{2}.
\end{equation*}%
From the second summand of the left side we get

\begin{eqnarray*}
l_{1}+u_{1} &=&1 \\
l_{2}+u_{2} &=&1 \\
a &=&0,b_{1}=l_{1},b_{2}-l_{2}=1\Rightarrow b_{2}=1,l_{2}=0,u_{2}=1, \\
d &=&0,e_{1}=u_{1},e_{2}=u_{2}=1
\end{eqnarray*}%
Since $\alpha \left( 1_{H};0,0,1,1\right) \equiv 1+e_{2}\equiv 0$ and $%
\alpha \left( 1_{H};1,0,0,1\right) \equiv a+b_{1}\equiv 1$ we get

\begin{equation*}
\left[ B(x_{2}\otimes x_{1};X_{2},x_{1}x_{2})-B(x_{2}\otimes
x_{1};X_{1}X_{2},x_{2})\right] X_{2}\otimes 1_{H}\otimes gx_{1}x_{2}.
\end{equation*}%
Since there is nothing in the right side, we obtain%
\begin{gather*}
B(g\otimes x_{1};X_{2},x_{1})-B(g\otimes x_{1};X_{1}X_{2},1_{H})+ \\
B(x_{2}\otimes x_{1};X_{2},x_{1}x_{2})-B(x_{2}\otimes
x_{1};X_{1}X_{2},x_{2})=0
\end{gather*}%
which holds in view of the form of the elements.

\subsection{$B(x_{2}\otimes x_{1};X_{2},gx_{1})$}

We deduce that
\begin{equation*}
a=b_{1}=0,b_{2}=1,d=e_{1}=1,e_{2}=0
\end{equation*}%
and we get%
\begin{eqnarray*}
&&\left( -1\right) ^{\alpha \left( 1_{H};0,0,0,0\right) }B(x_{2}\otimes
x_{1};X_{2},gx_{1})X_{2}\otimes gx_{1}\otimes g+ \\
&&\left( -1\right) ^{\alpha \left( 1_{H};0,1,0,0\right) }B(x_{2}\otimes
x_{1};X_{2},gx_{1})1_{A}\otimes gx_{1}\otimes x_{2}+ \\
&&\left( -1\right) ^{\alpha \left( 1_{H};0,0,1,0\right) }B(x_{2}\otimes
x_{1};X_{2},gx_{1})X_{2}\otimes g\otimes x_{1}+ \\
&&\left( -1\right) ^{\alpha \left( 1_{H};0,1,1,0\right) }B(x_{2}\otimes
x_{1};X_{2},gx_{1})1_{A}\otimes g\otimes gx_{1}x_{2}
\end{eqnarray*}

\subsubsection{Case $X_{2}\otimes gx_{1}\otimes g$}

This case is trivial.

\subsubsection{Case $1_{A}\otimes gx_{1}\otimes x_{2}$}

From the first summand of the left side we get

\begin{eqnarray*}
l_{1} &=&u_{1}=0 \\
l_{2} &=&u_{2}=0 \\
a &=&0,b_{1}=0,b_{2}=0, \\
d &=&e_{1}=1,e_{2}=1,
\end{eqnarray*}%
Since $\alpha \left( x_{2};0,0,0,0\right) =a+b_{1}+b_{2}=0$ we get%
\begin{equation*}
+B(g\otimes x_{1};1_{A},gx_{1})1_{A}\otimes gx_{1}\otimes x_{2}
\end{equation*}%
From the second summand of the left side we get

~%
\begin{eqnarray*}
l_{1} &=&u_{1}=0 \\
l_{2}+u_{2} &=&1 \\
a &=&b_{1}=0,b_{2}=l_{2}, \\
d &=&e_{1}=1,e_{2}=u_{2}
\end{eqnarray*}%
since $\alpha \left( 1_{H};0,0,0,1\right) \equiv a+b_{1}+b_{2}\equiv 0$ $%
\alpha \left( 1_{H};0,1,0,0\right) =0$ we obtain%
\begin{equation*}
\left[ B(x_{2}\otimes x_{1};1_{A},gx_{1}x_{2})+B(x_{2}\otimes
x_{1};X_{2},gx_{1})\right] 1_{A}\otimes gx_{1}\otimes x_{2}.
\end{equation*}

Since there is nothing in the right side, we obtain%
\begin{equation*}
B(g\otimes x_{1};1_{A},gx_{1})+B(x_{2}\otimes
x_{1};1_{A},gx_{1}x_{2})+B(x_{2}\otimes x_{1};X_{2},gx_{1})=0
\end{equation*}%
which holds in view of the form of the elements.

\subsubsection{Case $X_{2}\otimes g\otimes x_{1}$}

Nothing from the first summand of the left side. From the second summand of
the left side, we deduce that

\begin{eqnarray*}
l_{1}+u_{1} &=&1 \\
l_{2} &=&u_{2}=0 \\
a &=&0,b_{1}=l_{1},b_{2}=1, \\
d &=&1,e_{1}=u_{1},e_{2}=0.
\end{eqnarray*}%
Since $\alpha \left( 1_{H};0,0,1,0\right) =e_{2}+\left( a+b_{1}+b_{2}\right)
\equiv 1,$and $\alpha \left( 1_{H};1,0,0,0\right) =b_{2}=1$ we get%
\begin{equation*}
\left[ -B(x_{2}\otimes x_{1};X_{2},gx_{1})-B(x_{2}\otimes x_{1};X_{1}X_{2},g)%
\right] X_{2}\otimes g\otimes x_{1}.
\end{equation*}%
By considering also the right side, we get

\begin{equation*}
-B(x_{2}\otimes x_{1};X_{2},gx_{1})-B(x_{2}\otimes
x_{1};X_{1}X_{2},g)-B(x_{2}\otimes 1_{H};X_{2},g)=0
\end{equation*}%
which holds in view of the form of the elements.

\subsubsection{Case $1_{A}\otimes g\otimes gx_{1}x_{2}$}

From the first summand of the left side we get

\begin{eqnarray*}
l_{1}+u_{1} &=&1 \\
l_{2} &=&u_{2}=0 \\
a &=&0,b_{2}=0,b_{1}=l_{1} \\
d &=&1,e_{2}=0,e_{1}=u_{1}.
\end{eqnarray*}%
Since $\alpha \left( x_{2};0,0,1,0\right) =e_{2}=0$ $\ $and $\alpha \left(
x_{2};1,0,0,0\right) \equiv a+b_{1}\equiv 1$ we get

\begin{equation*}
\left[ B(g\otimes x_{1};1_{A},gx_{1})-B(g\otimes x_{1};X_{1},g)\right]
1_{A}\otimes g\otimes gx_{1}x_{2}.
\end{equation*}%
From the second summand of the left side we get

\begin{eqnarray*}
l_{1}+u_{1} &=&1 \\
l_{2}+u_{2} &=&1 \\
a &=&0,b_{1}=l_{1},b_{2}=l_{2}, \\
d &=&1,e_{1}=u_{1},e_{2}=u_{2}
\end{eqnarray*}%
\begin{eqnarray*}
\alpha \left( 1_{H};0,0,1,1\right) &=&1+e_{2}\equiv 0\text{, }\alpha \left(
1_{H};0,1,1,0\right) =e_{2}+a+b_{1}+b_{2}+1\equiv 0, \\
\alpha \left( 1_{H};1,0,0,1\right) &\equiv &a+b_{1}=1\text{and }\alpha
\left( 1_{H};1,1,0,0\right) \equiv 1+b_{2}\equiv 0
\end{eqnarray*}%
$,$ we get%
\begin{equation*}
\left[
\begin{array}{c}
B(x_{2}\otimes x_{1};1_{A},gx_{1}x_{2})+B(x_{2}\otimes x_{1};X_{2},gx_{1})+
\\
-B(x_{2}\otimes x_{1};X_{1},gx_{2})+B(x_{2}\otimes x_{1};X_{1}X_{2},1_{H})%
\end{array}%
\right] 1_{A}\otimes x_{2}\otimes gx_{1}x_{2}.
\end{equation*}%
Since there is nothing from the right side, we obtain%
\begin{eqnarray*}
&&B(g\otimes x_{1};1_{A},gx_{1})-B(g\otimes x_{1};X_{1},g)+ \\
&&%
\begin{array}{c}
B(x_{2}\otimes x_{1};1_{A},gx_{1}x_{2})+B(x_{2}\otimes x_{1};X_{2},gx_{1})+
\\
-B(x_{2}\otimes x_{1};X_{1},gx_{2})+B(x_{2}\otimes x_{1};X_{1}X_{2},1_{H})%
\end{array}%
\end{eqnarray*}%
which we already got.

\subsection{$B(x_{2}\otimes x_{1};X_{2},gx_{2})$}

We deduce that
\begin{equation*}
a=b_{1}=0,b_{2}=1,d=e_{2}=1,e_{1}=0
\end{equation*}%
and we get%
\begin{eqnarray*}
&&\left( -1\right) ^{\alpha \left( 1_{H};0,0,0,1\right) }B(x_{2}\otimes
x_{1};X_{2},gx_{2})X_{2}\otimes g\otimes x_{2}+ \\
&&\left( -1\right) ^{\alpha \left( 1_{H};0,1,0,0\right) }B(x_{2}\otimes
x_{1};X_{2},gx_{2})1_{A}\otimes gx_{2}\otimes x_{2}+
\end{eqnarray*}

\subsubsection{Case $X_{2}\otimes g\otimes x_{2}$}

From the first summand of the left side we get

\begin{eqnarray*}
l_{1} &=&u_{1}=0 \\
l_{2} &=&u_{2}=0 \\
a &=&0,b_{1}=0,b_{2}=1, \\
d &=&1,e_{1}=0,e_{2}=0
\end{eqnarray*}%
Since $\alpha \left( x_{2};0,0,0,0\right) =a+b_{1}+b_{2}=1$ we get%
\begin{equation*}
-B(g\otimes x_{1};X_{2},g)X_{2}\otimes g\otimes x_{2}
\end{equation*}%
From the second summand of the left side we get

~%
\begin{eqnarray*}
l_{1} &=&u_{1}=0 \\
l_{2}+u_{2} &=&1 \\
a &=&0,b_{1}=0,b_{2}-l_{2}=1\Rightarrow b_{2}=1,l_{2}=0,u_{2}=1 \\
d &=&1,e_{1}=0,e_{2}=u_{2}=1
\end{eqnarray*}%
since $\alpha \left( 1_{H};0,0,0,1\right) \equiv a+b_{1}+b_{2}=1$ we obtain

\begin{equation*}
-B(x_{2}\otimes x_{1};X_{2},gx_{2})X_{2}\otimes g\otimes x_{2}.
\end{equation*}%
Since there is nothing from the right side, we obtain%
\begin{equation*}
-B(g\otimes x_{1};X_{2},g)-B(x_{2}\otimes x_{1};X_{2},gx_{2})=0
\end{equation*}%
which holds in view of the form of the elements.

\subsubsection{Case $1_{A}\otimes gx_{2}\otimes x_{2}$}

From the first summand of the left side we get

\begin{eqnarray*}
l_{1} &=&u_{1}=0 \\
l_{2} &=&u_{2}=0 \\
a &=&b_{1}=b_{2}=0, \\
d &=&1,e_{1}=0,e_{2}=1
\end{eqnarray*}%
Since $\alpha \left( x_{2};0,0,0,0\right) =a+b_{1}+b_{2}=0$ we get%
\begin{equation*}
+B(g\otimes x_{1};1_{A},gx_{2})1_{A}\otimes gx_{2}\otimes x_{2}
\end{equation*}%
From the second summand of the left side we get

~%
\begin{eqnarray*}
l_{1} &=&u_{1}=0 \\
l_{2}+u_{2} &=&1 \\
a &=&0,b_{1}=0,b_{2}=l_{2}, \\
d &=&1,e_{1}=0,e_{2}-u_{2}=1\Rightarrow e_{2}=1,u_{2}=0,b_{2}=l_{2}=1
\end{eqnarray*}%
since $\alpha \left( 1_{H};0,1,0,0\right) \equiv 0$ we obtain

\begin{equation*}
+B(x_{2}\otimes x_{1};X_{2},gx_{2})1_{A}\otimes gx_{2}\otimes x_{2}.
\end{equation*}%
Since there is nothing from the right side, we obtain%
\begin{equation*}
+B(g\otimes x_{1};1_{A},gx_{2})+B(x_{2}\otimes x_{1};X_{2},gx_{2})=0
\end{equation*}%
which holds in view of the form of the elements.

\subsection{$B(x_{2}\otimes x_{1};X_{1}X_{2},g)$}

We deduce that
\begin{equation*}
a=0,b_{1}=b_{2}=1,d=1,e_{1}=e_{2}=0
\end{equation*}%
and we get%
\begin{eqnarray*}
&&\left( -1\right) ^{\alpha \left( 1_{H};0,0,0,0\right) }B(x_{2}\otimes
x_{1};X_{1}X_{2},g)X_{1}X_{2}\otimes g\otimes g \\
&&\left( -1\right) ^{\alpha \left( 1_{H};1,0,0,0\right) }B(x_{2}\otimes
x_{1};X_{1}X_{2},g)X_{2}\otimes g\otimes x_{1} \\
&&\left( -1\right) ^{\alpha \left( 1_{H};0,1,0,0\right) }B(x_{2}\otimes
x_{1};X_{1}X_{2},g)X_{1}\otimes g\otimes x_{2} \\
&&\left( -1\right) ^{\alpha \left( 1_{H};1,1,0,0\right) }B(x_{2}\otimes
x_{1};X_{1}X_{2},g)1_{A}\otimes g\otimes gx_{1}x_{2}
\end{eqnarray*}

\subsubsection{Case $X_{1}X_{2}\otimes g\otimes g$}

This case is trivial

\subsubsection{Case $X_{2}\otimes g\otimes x_{1}$}

Already considered in subsection $B(x_{2}\otimes x_{1};X_{2},gx_{1}).$

\subsubsection{Case $X_{1}\otimes g\otimes x_{2}$}

Already considered in subsection $B\left( x_{2}\otimes
x_{1};X_{1},gx_{2}\right) .$

\subsubsection{Case $1_{A}\otimes g\otimes gx_{1}x_{2}$}

\ Already considered in subsection $B\left( x_{2}\otimes
x_{1};1_{A},gx_{1}x_{2}\right) .$

\subsection{$B(x_{2}\otimes x_{1};X_{1}X_{2},x_{1})$}

We deduce that
\begin{equation*}
a=0,b_{1}=b_{2}=1,d=e_{2}=0,e_{1}=1
\end{equation*}%
and we get%
\begin{equation*}
\left( -1\right) ^{\alpha \left( 1_{H};l_{1},l_{2},u_{1},0\right)
}B(x_{2}\otimes x_{1};X_{1}X_{2},x_{1})X_{1}^{1-l_{1}}X_{2}^{1-l_{2}}\otimes
x_{1}^{1-u_{1}}\otimes
g^{1+l_{1}+l_{2}+u_{1}}x_{1}^{l_{1}+u_{1}}x_{2}^{l_{2}}
\end{equation*}%
\begin{gather*}
\left( -1\right) ^{\alpha \left( 1_{H};0,0,0,0\right) }B(x_{2}\otimes
x_{1};X_{1}X_{2},x_{1})X_{1}X_{2}\otimes x_{1}\otimes g \\
+\left( -1\right) ^{\alpha \left( 1_{H};1,0,0,0\right) }B(x_{2}\otimes
x_{1};X_{1}X_{2},x_{1})X_{2}\otimes x_{1}\otimes x_{1} \\
+\left( -1\right) ^{\alpha \left( 1_{H};0,1,0,0\right) }B(x_{2}\otimes
x_{1};X_{1}X_{2},x_{1})X_{1}\otimes x_{1}\otimes x_{2} \\
+\left( -1\right) ^{\alpha \left( 1_{H};1,1,0,0\right) }B(x_{2}\otimes
x_{1};X_{1}X_{2},x_{1})1_{A}\otimes x_{1}\otimes gx_{1}x_{2} \\
+\left( -1\right) ^{\alpha \left( 1_{H};0,0,1,0\right) }B(x_{2}\otimes
x_{1};X_{1}X_{2},x_{1})X_{1}X_{2}\otimes 1_{H}\otimes x_{1} \\
+\left( -1\right) ^{\alpha \left( 1_{H};1,0,1,0\right) }B(x_{2}\otimes
x_{1};X_{1}X_{2},x_{1})X_{1}^{1-l_{1}}X_{2}\otimes 1_{H}\otimes
g^{l_{1}}x_{1}^{1+1}=0 \\
+\left( -1\right) ^{\alpha \left( 1_{H};0,1,1,0\right) }B(x_{2}\otimes
x_{1};X_{1}X_{2},x_{1})X_{1}\otimes 1_{H}\otimes gx_{1}x_{2} \\
++\left( -1\right) ^{\alpha \left( 1_{H};1,1,1,0\right) }B(x_{2}\otimes
x_{1};X_{1}X_{2},x_{1})1_{A}\otimes 1_{H}\otimes x_{1}^{1+1}x_{2}=0
\end{gather*}

\subsubsection{Case $X_{1}X_{2}\otimes x_{1}\otimes g$}

This is trivial

\subsubsection{Case $X_{2}\otimes x_{1}\otimes x_{1}$}

Nothing from the first summand of the left side. From the second summand of
the left side, we deduce that

\begin{eqnarray*}
l_{1}+u_{1} &=&1 \\
l_{2} &=&u_{2}=0 \\
a &=&0,b_{1}=l_{1},b_{2}=1, \\
d &=&0,e_{1}-u_{1}=1\Rightarrow e_{1}=1,u_{1}=0,b_{1}=l_{1}=1,e_{2}=0.
\end{eqnarray*}%
Since $\alpha \left( 1_{H};1,0,0,0\right) =b_{2}=1$ we get%
\begin{equation*}
-B(x_{2}\otimes x_{1};X_{1}X_{2},x_{1})X_{2}\otimes x_{1}\otimes x_{1}.
\end{equation*}%
By considering also the right side, we get%
\begin{equation*}
-B(x_{2}\otimes x_{1};X_{1}X_{2},x_{1})-B(x_{2}\otimes 1_{H};X_{2},x_{1})=0
\end{equation*}%
which holds in view of the form of the elements.

\subsubsection{Case $X_{1}\otimes x_{1}\otimes x_{2}$}

Already considered in subsection $B\left( x_{2}\otimes
x_{1};X_{1},x_{1}x_{2}\right) .$

\subsubsection{Case $1_{A}\otimes x_{1}\otimes gx_{1}x_{2}$}

Already considered in subsection $B\left( x_{2}\otimes
x_{1};X_{1},x_{1}x_{2}\right) .$

\subsubsection{Case $X_{1}X_{2}\otimes 1_{H}\otimes x_{1}$}

Nothing from the first summand of the left side. From the second summand of
the left side, we deduce that

\begin{eqnarray*}
l_{1}+u_{1} &=&1 \\
l_{2} &=&u_{2}=0 \\
a &=&0,b_{1}-l_{1}=1\Rightarrow b_{1}=1,l_{1}=0,u_{1}=1,b_{2}=1, \\
d &=&0,e_{1}=u_{1}=1,e_{2}=0.
\end{eqnarray*}%
Since $\alpha \left( 1_{H};0,0,1,0\right) =e_{2}+\left( a+b_{1}+b_{2}\right)
\equiv 0$ we get%
\begin{equation*}
B(x_{2}\otimes x_{1};X_{1}X_{2},x_{1})X_{1}X_{2}\otimes 1_{H}\otimes x_{1}.
\end{equation*}%
By considering also the right side, we get%
\begin{equation*}
B(x_{2}\otimes x_{1};X_{1}X_{2},x_{1})-B(x_{2}\otimes
1_{H};X_{1}X_{2},1_{H})=0
\end{equation*}%
which holds in view of the form of the elements.

\subsubsection{Case $X_{1}\otimes 1_{H}\otimes gx_{1}x_{2}$}

Already considered in subsection $B\left( x_{2}\otimes
x_{1};X_{1},x_{1}x_{2}\right) .$

\subsection{$B\left( x_{2}\otimes x_{1};X_{1}X_{2},x_{2}\right) $}

We deduce that
\begin{equation*}
a=0,b_{1}=b_{2}=1,d=e_{1}=0,e_{2}=1
\end{equation*}%
and we get%
\begin{equation*}
\left( -1\right) ^{\alpha \left( 1_{H};l_{1},l_{2},0,u_{2}\right) }B\left(
x_{2}\otimes x_{1};X_{1}X_{2},x_{2}\right)
X_{1}^{1-l_{1}}X_{2}^{1-l_{2}}\otimes x_{2}^{1-u_{2}}\otimes
g^{1+l_{1}+l_{2}+u_{2}}x_{1}^{l_{1}}x_{2}^{l_{2}+u_{2}}
\end{equation*}%
\begin{gather*}
\left( -1\right) ^{\alpha \left( 1_{H};0,0,0,0\right) }B\left( x_{2}\otimes
x_{1};X_{1}X_{2},x_{2}\right) X_{1}X_{2}\otimes x_{2}\otimes g \\
\left( -1\right) ^{\alpha \left( 1_{H};1,0,0,0\right) }B\left( x_{2}\otimes
x_{1};X_{1}X_{2},x_{2}\right) X_{2}\otimes x_{2}\otimes x_{1} \\
\left( -1\right) ^{\alpha \left( 1_{H};0,1,0,0\right) }B\left( x_{2}\otimes
x_{1};X_{1}X_{2},x_{2}\right) X_{1}\otimes x_{2}\otimes x_{2} \\
\left( -1\right) ^{\alpha \left( 1_{H};1,1,0,0\right) }B\left( x_{2}\otimes
x_{1};X_{1}X_{2},x_{2}\right) 1_{A}\otimes x_{2}\otimes gx_{1}x_{2} \\
\left( -1\right) ^{\alpha \left( 1_{H};0,0,0,1\right) }B\left( x_{2}\otimes
x_{1};X_{1}X_{2},x_{2}\right) X_{1}X_{2}\otimes 1_{H}\otimes x_{2} \\
\left( -1\right) ^{\alpha \left( 1_{H};1,0,0,1\right) }B\left( x_{2}\otimes
x_{1};X_{1}X_{2},x_{2}\right) X_{2}\otimes 1_{H}\otimes gx_{1}x_{2} \\
\left( -1\right) ^{\alpha \left( 1_{H};l_{1},1,0,1\right) }B\left(
x_{2}\otimes x_{1};X_{1}X_{2},x_{2}\right) X_{1}^{1-l_{1}}\otimes
1_{H}\otimes g^{l_{1}}x_{1}^{l_{1}}x_{2}^{1+1}=0
\end{gather*}

\subsubsection{Case $X_{1}X_{2}\otimes x_{2}\otimes g$}

This case is trivial.

\subsubsection{Case $X_{2}\otimes x_{2}\otimes x_{1}$}

Already done in subsection $B\left( x_{2}\otimes
x_{1};X_{2},x_{1}x_{2}\right) .$

\subsubsection{Case $X_{1}\otimes x_{2}\otimes x_{2}$}

From the first summand of the left side we get

\begin{eqnarray*}
l_{1} &=&u_{1}=0 \\
l_{2} &=&u_{2}=0 \\
a &=&b_{2}=0,b_{1}=1 \\
d &=&e_{1}=0,e_{2}=1
\end{eqnarray*}%
Since $\alpha \left( x_{2};0,0,0,0\right) =a+b_{1}+b_{2}\equiv 1$ we get%
\begin{equation*}
-B(g\otimes x_{1};X_{1},x_{2})X_{1}\otimes x_{2}\otimes x_{2}
\end{equation*}%
From the second summand of the left side we get

~%
\begin{eqnarray*}
l_{1} &=&u_{1}=0 \\
l_{2}+u_{2} &=&1 \\
a &=&0,b_{1}=1,b_{2}=l_{2}, \\
d &=&e_{1}=0,e_{2}-u_{2}=1\Rightarrow e_{2}=1,u_{2}=0,b_{2}=l_{2}=1
\end{eqnarray*}%
since $\alpha \left( 1_{H};0,1,0,0\right) \equiv 0$, we get%
\begin{equation*}
B(x_{2}\otimes x_{1};X_{1}X_{2},x_{2})X_{1}\otimes x_{2}\otimes x_{2}.
\end{equation*}%
Since there is nothing from the right side, we get%
\begin{equation*}
-B(g\otimes x_{1};X_{1},x_{2})+B(x_{2}\otimes x_{1};X_{1}X_{2},x_{2})=0
\end{equation*}%
which holds in view of the form of the elements.

\subsubsection{Case $1_{A}\otimes x_{2}\otimes gx_{1}x_{2}$}

This case already appeared in subsection $B\left( x_{2}\otimes
x_{1};X_{2},x_{1}x_{2}\right) .$

\subsubsection{Case $X_{1}X_{2}\otimes 1_{H}\otimes x_{2}$}

From the first summand of the left side we get

\begin{eqnarray*}
l_{1} &=&u_{1}=0 \\
l_{2} &=&u_{2}=0 \\
a &=&0,b_{1}=b_{2}=1, \\
d &=&e_{1}=e_{2}=0
\end{eqnarray*}%
Since $\alpha \left( x_{2};0,0,0,0\right) =a+b_{1}+b_{2}\equiv 0$ we get%
\begin{equation*}
+B(g\otimes x_{1};X_{1}X_{2},1_{H})X_{1}X_{2}\otimes 1_{H}\otimes x_{2}
\end{equation*}%
From the second summand of the left side we get

~%
\begin{eqnarray*}
l_{1} &=&u_{1}=0 \\
l_{2}+u_{2} &=&1 \\
a &=&0,b_{1}=1,b_{2}-l_{2}=1\Rightarrow b_{2}=1,l_{2}=0,u_{2}=1 \\
d &=&e_{1}=0,e_{2}=u_{2}=1.
\end{eqnarray*}%
Since $\alpha \left( 1_{H};0,0,0,1\right) \equiv a+b_{1}+b_{2}\equiv 0,$ we
obtain

\begin{equation*}
+B(x_{2}\otimes x_{1};X_{1}X_{2},x_{2})X_{1}X_{2}\otimes 1_{H}\otimes x_{2}.
\end{equation*}%
Since there is nothing from the right side, we get%
\begin{equation*}
B(g\otimes x_{1};X_{1}X_{2},1_{H})+B(x_{2}\otimes x_{1};X_{1}X_{2},x_{2})=0
\end{equation*}%
which holds in view of the form of the elements.

\subsubsection{Case $X_{2}\otimes 1_{H}\otimes gx_{1}x_{2}$}

This was already considered in subsection $B\left( x_{2}\otimes
x_{1};X_{2},x_{1}x_{2}\right) .$

\subsection{$B\left( x_{2}\otimes x_{1};X_{1}X_{2},gx_{1}x_{2}\right) $}

We deduce that
\begin{equation*}
a=0,b_{1}=b_{2}=1,d=e_{1}=e_{2}=1
\end{equation*}%
and we get%
\begin{gather*}
\left( -1\right) ^{\alpha \left( 1_{H};0,0,0,0\right) }B\left( x_{2}\otimes
x_{1};X_{1}X_{2},gx_{1}x_{2}\right) X_{1}X_{2}\otimes gx_{1}x_{2}\otimes g \\
+\left( -1\right) ^{\alpha \left( 1_{H};1,0,0,0\right) }B\left( x_{2}\otimes
x_{1};X_{1}X_{2},gx_{1}x_{2}\right) X_{2}\otimes gx_{1}x_{2}\otimes x_{1} \\
+\left( -1\right) ^{\alpha \left( 1_{H};0,1,0,0\right) }B\left( x_{2}\otimes
x_{1};X_{1}X_{2},gx_{1}x_{2}\right) X_{1}\otimes gx_{1}x_{2}\otimes x_{2} \\
+\left( -1\right) ^{\alpha \left( 1_{H};1,1,0,0\right) }B\left( x_{2}\otimes
x_{1};X_{1}X_{2},gx_{1}x_{2}\right) 1_{A}\otimes gx_{1}x_{2}\otimes
gx_{1}x_{2} \\
+\left( -1\right) ^{\alpha \left( 1_{H};0,0,1,0\right) }B\left( x_{2}\otimes
x_{1};X_{1}X_{2},gx_{1}x_{2}\right) X_{1}X_{2}\otimes gx_{2}\otimes x_{1} \\
+\left( -1\right) ^{\alpha \left( 1_{H};1,0,1,0\right) }B\left( x_{2}\otimes
x_{1};X_{1}X_{2},gx_{1}x_{2}\right) X_{1}^{1-l_{1}}X_{2}\otimes
gx_{2}\otimes g^{l_{1}}x_{1}^{1+1}=0 \\
+\left( -1\right) ^{\alpha \left( 1_{H};0,1,1,0\right) }B\left( x_{2}\otimes
x_{1};X_{1}X_{2},gx_{1}x_{2}\right) X_{1}\otimes gx_{2}\otimes gx_{1}x_{2} \\
+\left( -1\right) ^{\alpha \left( 1_{H};1,1,1,0\right) }B\left( x_{2}\otimes
x_{1};X_{1}X_{2},gx_{1}x_{2}\right) X_{1}^{1-l_{1}}\otimes gx_{2}\otimes
g^{l_{1}+1}x_{1}^{1+1}x_{2}=0 \\
+\left( -1\right) ^{\alpha \left( 1_{H};0,0,0,1\right) }B\left( x_{2}\otimes
x_{1};X_{1}X_{2},gx_{1}x_{2}\right) X_{1}X_{2}\otimes gx_{1}\otimes x_{2} \\
+\left( -1\right) ^{\alpha \left( 1_{H};1,0,0,1\right) }B\left( x_{2}\otimes
x_{1};X_{1}X_{2},gx_{1}x_{2}\right) X_{2}\otimes gx_{1}\otimes gx_{1}x_{2} \\
+\left( -1\right) ^{\alpha \left( 1_{H};l_{1},1,0,1\right) }B\left(
x_{2}\otimes x_{1};X_{1}X_{2},gx_{1}x_{2}\right)
X_{1}^{1-l_{1}}X_{2}^{1-l_{2}}\otimes gx_{1}\otimes
g^{l_{1}+l_{2}}x_{1}^{l_{1}}x_{2}^{1+1}=0 \\
+\left( -1\right) ^{\alpha \left( 1_{H};0,0,1,1\right) }B\left( x_{2}\otimes
x_{1};X_{1}X_{2},gx_{1}x_{2}\right) X_{1}X_{2}\otimes g\otimes gx_{1}x_{2} \\
+\left( -1\right) ^{\alpha \left( 1_{H};1,0,1,1\right) }B\left( x_{2}\otimes
x_{1};X_{1}X_{2},gx_{1}x_{2}\right) X_{1}^{1-l_{1}}X_{2}\otimes g\otimes
g^{l_{1}+1}x_{1}^{1+1}x_{2}=0 \\
+\left( -1\right) ^{\alpha \left( 1_{H};1,1,1,1\right) }B\left( x_{2}\otimes
x_{1};X_{1}X_{2},gx_{1}x_{2}\right) X_{1}^{1-l_{1}}X_{2}^{1-l_{2}}\otimes
g\otimes g^{l_{1}+l_{2}+1}x_{1}^{l_{1}+1}x_{2}^{1+1}=0
\end{gather*}

\subsubsection{Case $X_{1}X_{2}\otimes gx_{1}x_{2}\otimes g$}

This case is trivial.

\subsubsection{Case $X_{2}\otimes gx_{1}x_{2}\otimes x_{1}$}

Nothing from the first summand of the left side. From the second summand of
the left side, we deduce that

\begin{eqnarray*}
l_{1}+u_{1} &=&1 \\
l_{2} &=&u_{2}=0 \\
a &=&0,b_{1}=l_{1},b_{2}=1, \\
d &=&1,e_{1}-u_{1}=1\Rightarrow e_{1}=1,u_{1}=0,b_{1}=l_{1}=1,e_{2}=1.
\end{eqnarray*}%
Since $\alpha \left( 1_{H};1,0,0,0\right) =b_{2}=1$ we get%
\begin{equation*}
-B(x_{2}\otimes x_{1};X_{1}X_{2},gx_{1}x_{2})X_{2}\otimes gx_{1}x_{2}\otimes
x_{1}.
\end{equation*}%
By considering also the right side, we get%
\begin{equation*}
-B(x_{2}\otimes x_{1};X_{1}X_{2},gx_{1}x_{2})-B(x_{2}\otimes
1_{H};X_{2},gx_{1}x_{2})=0
\end{equation*}%
which holds in view of the form of the elements.

\subsubsection{Case $X_{1}\otimes gx_{1}x_{2}\otimes x_{2}$}

From the first summand of the left side we get

\begin{eqnarray*}
l_{1} &=&u_{1}=0 \\
l_{2} &=&u_{2}=0 \\
a &=&b_{2}=0,b_{1}=1 \\
d &=&e_{1}=e_{2}=1
\end{eqnarray*}%
Since $\alpha \left( x_{2};0,0,0,0\right) =a+b_{1}+b_{2}=1$ we get%
\begin{equation*}
-B(g\otimes x_{1};X_{1},gx_{1}x_{2})X_{1}\otimes gx_{1}x_{2}\otimes x_{2}
\end{equation*}%
From the second summand of the left side we get

~%
\begin{eqnarray*}
l_{1} &=&u_{1}=0 \\
l_{2}+u_{2} &=&1 \\
a &=&0,b_{1}=1,b_{2}=l_{2}, \\
d &=&1,e_{1}=1,e_{2}-u_{2}=1\Rightarrow e_{2}=1,u_{2}=0,b_{2}=l_{2}=1
\end{eqnarray*}%
since $\alpha \left( 1_{H};0,1,0,0\right) \equiv 0$ we obtain

\begin{equation*}
+B(x_{2}\otimes x_{1};X_{1}X_{2},gx_{1}x_{2})X_{1}\otimes gx_{1}x_{2}\otimes
x_{2}.
\end{equation*}%
Since there is nothing from the right side, we obtain%
\begin{equation*}
-B(g\otimes x_{1};X_{1},gx_{1}x_{2})+B(x_{2}\otimes
x_{1};X_{1}X_{2},gx_{1}x_{2})=0
\end{equation*}%
which holds in view of the form of the elements.

\subsubsection{Case $1_{A}\otimes gx_{1}x_{2}\otimes gx_{1}x_{2}$}

From the first summand of the left side we get

\begin{eqnarray*}
l_{1}+u_{1} &=&1 \\
l_{2} &=&u_{2}=0 \\
a &=&0,b_{2}=0,b_{1}=l_{1} \\
d &=&e_{2}=1,e_{1}-u_{1}=1\Rightarrow e_{1}=1,u_{1}=0,b_{1}=l_{1}=1.
\end{eqnarray*}%
Since $\alpha \left( x_{2};1,0,0,0\right) \equiv a+b_{1}\equiv 1$ we get

\begin{equation*}
-B(g\otimes x_{1};X_{1},gx_{1}x_{2})1_{A}\otimes gx_{1}x_{2}\otimes
gx_{1}x_{2}.
\end{equation*}%
From the second summand of the left side we get

\begin{eqnarray*}
l_{1}+u_{1} &=&1 \\
l_{2}+u_{2} &=&1 \\
a &=&0,b_{1}=l_{1},b_{2}=l_{2}, \\
d &=&1,e_{1}-u_{1}=1\Rightarrow e_{1}=1,u_{1}=0,b_{1}=l_{1}=1, \\
e_{2}-u_{2} &=&1\Rightarrow e_{2}=1,u_{2}=0,b_{2}=l_{2}=1
\end{eqnarray*}%
Since $\alpha \left( 1_{H};1,1,0,0\right) \equiv 1+b_{2}\equiv 0,$ we get%
\begin{equation*}
B(x_{1}\otimes x_{2};X_{1}X_{2},gx_{1}x_{2})1_{A}\otimes gx_{1}x_{2}\otimes
gx_{1}x_{2}.
\end{equation*}%
Since there is nothing from the right side, we obtain%
\begin{equation*}
-B(g\otimes x_{1};X_{1},gx_{1}x_{2})+B(x_{1}\otimes
x_{2};X_{1}X_{2},gx_{1}x_{2})=0
\end{equation*}%
which holds in view of the form of the elements.

\subsubsection{Case $X_{1}X_{2}\otimes gx_{2}\otimes x_{1}$}

Nothing from the first summand of the left side. From the second summand of
the left side, we deduce that

\begin{eqnarray*}
l_{1}+u_{1} &=&1 \\
l_{2} &=&u_{2}=0 \\
a &=&0,b_{1}-l_{1}=1\Rightarrow b_{1}=1,l_{1}=0,u_{1}=1,b_{2}=1, \\
d &=&1,e_{1}=u_{1}=1,e_{2}=1.
\end{eqnarray*}%
Since $\alpha \left( 1_{H};0,0,1,0\right) \equiv e_{2}+\left(
a+b_{1}+b_{2}\right) \equiv 1$ we get%
\begin{equation*}
-B(x_{2}\otimes x_{1};X_{1}X_{2},gx_{1}x_{2})X_{1}X_{2}\otimes gx_{2}\otimes
x_{1}.
\end{equation*}%
By considering also the right side, we get

\begin{equation*}
-B(x_{2}\otimes x_{1};X_{1}X_{2},gx_{1}x_{2})-B(x_{2}\otimes
1_{H};X_{1}X_{2},gx_{2})=0
\end{equation*}%
which holds in view of the form of the elements.

\subsubsection{Case $X_{1}\otimes gx_{2}\otimes gx_{1}x_{2}$}

From the first summand of the left side we get

\begin{eqnarray*}
l_{1}+u_{1} &=&1 \\
l_{2} &=&u_{2}=0 \\
a &=&0,b_{2}=0,b_{1}-l_{1}=1\Rightarrow b_{1}=1,l_{1}=0,u_{1}=1 \\
d &=&e_{2}=1,e_{1}=u_{1}=1.
\end{eqnarray*}%
Since $\alpha \left( x_{2};0,0,1,0\right) \equiv e_{2}=1,$ we get

\begin{equation*}
-B(g\otimes x_{1};X_{1},gx_{1}x_{2})X_{1}\otimes gx_{2}\otimes gx_{1}x_{2}.
\end{equation*}%
From the second summand of the left side we get

\begin{eqnarray*}
l_{1}+u_{1} &=&1 \\
l_{2}+u_{2} &=&1 \\
a &=&0,b_{1}-l_{1}=1\Rightarrow b_{1}=1,l_{1}=0,u_{1}=1,b_{2}=l_{2}, \\
d &=&1,e_{1}=u_{1}=1, \\
e_{2}-u_{2} &=&1\Rightarrow e_{2}=1,u_{2}=0,b_{2}=l_{2}=1
\end{eqnarray*}%
Since $\alpha \left( 1_{H};0,1,1,0\right) \equiv e_{2}+a+b_{1}+b_{2}+1\equiv
0,$ we get%
\begin{equation*}
B(x_{1}\otimes x_{2};X_{1}X_{2},gx_{1}x_{2})X_{1}\otimes gx_{2}\otimes
gx_{1}x_{2}.
\end{equation*}%
Since there is nothing from the right side, we obtain%
\begin{equation*}
-B(g\otimes x_{1};X_{1},gx_{1}x_{2})+B(x_{1}\otimes
x_{2};X_{1}X_{2},gx_{1}x_{2})=0
\end{equation*}%
which holds in view of the form of the elements.

\subsubsection{Case $X_{1}X_{2}\otimes gx_{1}\otimes x_{2}$}

From the first summand of the left side we get

\begin{eqnarray*}
l_{1} &=&u_{1}=0 \\
l_{2} &=&u_{2}=0 \\
a &=&0,b_{1}=b_{2}=1, \\
d &=&e_{1}=1,e_{2}=0
\end{eqnarray*}%
Since $\alpha \left( x_{2};0,0,0,0\right) =a+b_{1}+b_{2}\equiv 0$ we get%
\begin{equation*}
+B(g\otimes x_{1};X_{1}X_{2},gx_{1})X_{1}X_{2}\otimes gx_{1}\otimes x_{2}
\end{equation*}%
From the second summand of the left side we get

~%
\begin{eqnarray*}
l_{1} &=&u_{1}=0 \\
l_{2}+u_{2} &=&1 \\
a &=&0,b_{1}=1,b_{2}-l_{2}=1\Rightarrow b_{2}=1,l_{2}=0,u_{2}=1 \\
d &=&e_{1}=1,e_{2}=u_{2}=1
\end{eqnarray*}%
since $\alpha \left( 1_{H};0,0,0,1\right) \equiv a+b_{1}+b_{2}\equiv 0$ we
obtain

\begin{equation*}
+B(x_{2}\otimes x_{1};X_{1}X_{2},gx_{1}x_{2})X_{2}\otimes g\otimes x_{2}.
\end{equation*}%
Since there is nothing from the right side, we obtain%
\begin{equation*}
B(g\otimes x_{1};X_{1}X_{2},gx_{1})+B(x_{2}\otimes
x_{1};X_{1}X_{2},gx_{1}x_{2})=0
\end{equation*}%
which holds in view of the form of the elements.

\subsubsection{Case $X_{2}\otimes gx_{1}\otimes gx_{1}x_{2}$}

From the first summand of the left side we get

\begin{eqnarray*}
l_{1}+u_{1} &=&1 \\
l_{2} &=&u_{2}=0 \\
a &=&0,b_{2}=1,b_{1}=l_{1} \\
d &=&1,e_{1}-u_{1}=1\Rightarrow e_{1}=1,u_{1}=0,b_{1}=l_{1}=1,e_{2}=0
\end{eqnarray*}%
Since $\alpha \left( x_{2};1,0,0,0\right) \equiv a+b_{1}+1+1\equiv 1,$ we get

\begin{equation*}
-B(g\otimes x_{1};X_{1}X_{2},gx_{1})X_{2}\otimes gx_{1}\otimes gx_{1}x_{2}.
\end{equation*}%
From the second summand of the left side we get

\begin{eqnarray*}
l_{1}+u_{1} &=&1 \\
l_{2}+u_{2} &=&1 \\
a &=&0,b_{1}=l_{1}, \\
b_{2}-l_{2} &=&1\Rightarrow b_{2}=1,l_{2}=0,u_{2}=1 \\
d &=&1,e_{1}-u_{1}=1\Rightarrow e_{1}=1,u_{1}=0,b_{1}=l_{1}=1 \\
e_{2} &=&u_{2}=1
\end{eqnarray*}%
Since $\alpha \left( 1_{H};1,0,0,1\right) \equiv a+b_{1}\equiv 1,$ we get%
\begin{equation*}
-B(x_{1}\otimes x_{2};X_{1}X_{2},gx_{1}x_{2})X_{1}\otimes gx_{2}\otimes
gx_{1}x_{2}.
\end{equation*}%
Since there is nothing from the right side, we obtain%
\begin{equation*}
-B(g\otimes x_{1};X_{1}X_{2},gx_{1})-B(x_{1}\otimes
x_{2};X_{1}X_{2},gx_{1}x_{2})=0
\end{equation*}%
which we got in the previous case.

\subsubsection{Case $X_{1}X_{2}\otimes g\otimes gx_{1}x_{2}$}

From the first summand of the left side we get

\begin{eqnarray*}
l_{1}+u_{1} &=&1 \\
l_{2} &=&u_{2}=0 \\
a &=&0,b_{1}-l_{1}=1\Rightarrow b_{1}=1,l_{1}=0,u_{1}=1,b_{2}=1, \\
d &=&1,e_{1}=u_{1}=1,e_{2}=0
\end{eqnarray*}%
Since $\alpha \left( x_{2};0,0,1,0\right) \equiv e_{2}=0,$ we get

\begin{equation*}
B(g\otimes x_{1};X_{1}X_{2},gx_{1})X_{1}X_{2}\otimes g\otimes gx_{1}x_{2}.
\end{equation*}%
From the second summand of the left side we get

\begin{eqnarray*}
l_{1}+u_{1} &=&1 \\
l_{2}+u_{2} &=&1 \\
a &=&0,b_{1}-l_{1}=1\Rightarrow b_{1}=1,l_{1}=0,u_{1}=1, \\
b_{2}-l_{2} &=&1\Rightarrow b_{2}=1,l_{2}=0,u_{2}=1, \\
d &=&1,e_{1}=u_{1}=1,e_{2}=u_{2}=1
\end{eqnarray*}%
Since $\alpha \left( 1_{H};0,0,1,1\right) \equiv 1+e_{2}\equiv 0,$ we get%
\begin{equation*}
B(x_{1}\otimes x_{2};X_{1}X_{2},gx_{1}x_{2})X_{1}X_{2}\otimes g\otimes
gx_{1}x_{2}.
\end{equation*}%
Since there is nothing from the right side, we obtain

\begin{equation*}
+B(g\otimes x_{1};X_{1}X_{2},gx_{1})+B(x_{1}\otimes
x_{2};X_{1}X_{2},gx_{1}x_{2})=0
\end{equation*}%
which is the previous equality.

\subsection{$B(x_{1}\otimes x_{2};GX_{1},g)$}

We deduce that%
\begin{equation*}
a=b_{1}=1,d=1
\end{equation*}%
and we get%
\begin{eqnarray*}
&&\left( -1\right) ^{\alpha \left( 1_{H};0,0,0,0\right) }B(x_{2}\otimes
x_{1};GX_{1},g)GX_{1}\otimes g\otimes g+ \\
&&\left( -1\right) ^{\alpha \left( 1_{H};1,0,0,0\right) }B(x_{2}\otimes
x_{1};GX_{1},g)G\otimes g\otimes x_{1}
\end{eqnarray*}

\subsubsection{Case $GX_{1}\otimes g\otimes g$}

This is trivial.

\subsubsection{Case $G\otimes g\otimes x_{1}$}

This case was already considered in subsection $B\left( x_{2}\otimes
x_{1};G,gx_{1}\right) .$

\subsection{$B\left( x_{2}\otimes x_{1};GX_{1},x_{1}\right) $}

We deduce that%
\begin{equation*}
a=b_{1}=1,e_{1}=1
\end{equation*}%
and we get%
\begin{gather*}
\left( -1\right) ^{\alpha \left( 1_{H};0,0,0,0\right) }B\left( x_{2}\otimes
x_{1};GX_{1},x_{1}\right) GX_{1}\otimes x_{1}\otimes g+ \\
+\left( -1\right) ^{\alpha \left( 1_{H};1,0,0,0\right) }B\left( x_{2}\otimes
x_{1};GX_{1},x_{1}\right) G\otimes x_{1}\otimes x_{1}+ \\
+\left( -1\right) ^{\alpha \left( 1_{H};0,0,1,0\right) }B\left( x_{2}\otimes
x_{1};GX_{1},x_{1}\right) GX_{1}\otimes 1_{H}\otimes x_{1}+ \\
+\left( -1\right) ^{\alpha \left( 1_{H};1,0,1,0\right) }B\left( x_{2}\otimes
x_{1};GX_{1},x_{1}\right) GX_{1}^{1-l_{1}}\otimes 1_{H}\otimes
g^{l_{1}}x_{1}^{1+1}=0
\end{gather*}

\subsubsection{Case $G\otimes x_{1}\otimes x_{1}$}

Nothing from the first summand of the left side. From the second summand of
the left side, we deduce that

\begin{eqnarray*}
l_{1}+u_{1} &=&1 \\
l_{2} &=&u_{2}=0 \\
a &=&1,b_{1}=l_{1},b_{2}=0, \\
d &=&0,e_{1}-u_{1}=1\Rightarrow e_{1}=1,u_{1}=0,b_{1}=l_{1}=1,e_{2}=0.
\end{eqnarray*}%
Since $\alpha \left( 1_{H};1,0,0,0\right) =b_{2}=0$ we get%
\begin{equation*}
B(x_{2}\otimes x_{1};GX_{1},x_{1})G\otimes x_{1}\otimes x_{1}.
\end{equation*}%
By considering also the right side, we get%
\begin{equation*}
B(x_{2}\otimes x_{1};GX_{1},x_{1})-B(x_{2}\otimes 1_{H};G,x_{1})=0
\end{equation*}

which holds in view of the form of the elements.

\subsubsection{Case $GX_{1}\otimes 1_{H}\otimes x_{1}$}

Nothing from the first summand of the left side. From the second summand of
the left side, we deduce that

\begin{eqnarray*}
l_{1}+u_{1} &=&1 \\
l_{2} &=&u_{2}=0 \\
a &=&1,b_{1}-l_{1}=1\Rightarrow b_{1}=1,l_{1}=0,u_{1}=1,b_{2}=0, \\
d &=&0,e_{1}=u_{1}=1,e_{2}=0.
\end{eqnarray*}%
Since $\alpha \left( 1_{H};0,0,1,0\right) \equiv e_{2}+\left(
a+b_{1}+b_{2}\right) \equiv 0$ we get%
\begin{equation*}
B(x_{2}\otimes x_{1};GX_{1},x_{1})GX_{1}\otimes 1_{H}\otimes x_{1}.
\end{equation*}%
By considering also the right side, we get%
\begin{equation*}
B(x_{2}\otimes x_{1};GX_{1},x_{1})-B(x_{2}\otimes 1_{H};GX_{1},1_{H})=0
\end{equation*}%
which holds in view of the form of the elements.

\subsection{$B\left( x_{2}\otimes x_{1};GX_{1},x_{2}\right) $}

We deduce that%
\begin{equation*}
a=b_{1}=1,e_{2}=1
\end{equation*}%
and we get%
\begin{eqnarray*}
&&\left( -1\right) ^{\alpha \left( 1_{H};0,0,0,0\right) }B\left(
x_{2}\otimes x_{1};GX_{1},x_{2}\right) GX_{1}^{{}}\otimes x_{2}\otimes g+ \\
&&\left( -1\right) ^{\alpha \left( 1_{H};1,0,0,0\right) }B\left(
x_{2}\otimes x_{1};GX_{1},x_{2}\right) G\otimes x_{2}\otimes x_{1}+ \\
&&\left( -1\right) ^{\alpha \left( 1_{H};0,0,0,1\right) }B\left(
x_{2}\otimes x_{1};GX_{1},x_{2}\right) GX_{1}\otimes 1_{H}\otimes x_{2}+ \\
&&\left( -1\right) ^{\alpha \left( 1_{H};1,0,0,1\right) }B\left(
x_{2}\otimes x_{1};GX_{1},x_{2}\right) G\otimes 1_{H}\otimes gx_{1}x_{2}
\end{eqnarray*}

\subsubsection{Case $G\otimes x_{2}\otimes x_{1}$}

This case was already considered in subsection $B\left( x_{2}\otimes
x_{1};G,x_{1}x_{2}\right) .$

\subsubsection{Case $GX_{1}\otimes 1_{H}\otimes x_{2}$}

From the first summand of the left side we get

\begin{eqnarray*}
l_{1} &=&u_{1}=0 \\
l_{2} &=&u_{2}=0 \\
a &=&b_{1}=1,b_{2}=0, \\
d &=&e_{1}=e_{2}=0
\end{eqnarray*}%
Since $\alpha \left( x_{2};0,0,0,0\right) =a+b_{1}+b_{2}\equiv 0$ we get%
\begin{equation*}
+B(g\otimes x_{1};GX_{1},1_{H})GX_{1}\otimes 1_{H}\otimes x_{2}
\end{equation*}%
From the second summand of the left side we get

~%
\begin{eqnarray*}
l_{1} &=&u_{1}=0 \\
l_{2}+u_{2} &=&1 \\
a &=&b_{1}=1,b_{2}=l_{2} \\
d &=&e_{1}=0,e_{2}=u_{2}
\end{eqnarray*}%
since $\alpha \left( 1_{H};0,0,0,1\right) \equiv a+b_{1}+b_{2}\equiv 0$ and $%
\alpha \left( 1_{H};0,1,0,0\right) \equiv 0$ we obtain

\begin{equation*}
\left[ B(x_{2}\otimes x_{1};GX_{1},x_{2})+B(x_{2}\otimes
x_{1};GX_{1}X_{2},1_{H})\right] GX_{1}\otimes 1_{H}\otimes x_{2}.
\end{equation*}%
Since there is nothing from the right side, we obtain\qquad
\begin{equation*}
B(g\otimes x_{1};GX_{1},1_{H})+B(x_{2}\otimes
x_{1};GX_{1},x_{2})+B(x_{2}\otimes x_{1};GX_{1}X_{2},1_{H})=0
\end{equation*}%
which holds in view of the form of the elements.

\subsubsection{Case $G\otimes 1_{H}\otimes gx_{1}x_{2}$}

This case was already considered in subsection $B\left( x_{2}\otimes
x_{1};G,x_{1}x_{2}\right) .$

\subsection{$B\left( x_{2}\otimes x_{1};GX_{1},gx_{1}x_{2}\right) $}

We deduce that%
\begin{eqnarray*}
a &=&b_{1}=1,b_{2}=0, \\
d &=&e_{1}=e_{2}=1
\end{eqnarray*}%
and we get%
\begin{gather*}
\left( -1\right) ^{\alpha \left( 1_{H};0,0,0,0\right) }B\left( x_{2}\otimes
x_{1};GX_{1},gx_{1}x_{2}\right) GX_{1}\otimes gx_{1}x_{2}\otimes g+ \\
\left( -1\right) ^{\alpha \left( 1_{H};1,0,0,0\right) }B\left( x_{2}\otimes
x_{1};GX_{1},gx_{1}x_{2}\right) G\otimes gx_{1}x_{2}\otimes x_{1} \\
\left( -1\right) ^{\alpha \left( 1_{H};0,0,1,0\right) }B\left( x_{2}\otimes
x_{1};GX_{1},gx_{1}x_{2}\right) GX_{1}\otimes gx_{2}\otimes x_{1} \\
\left( -1\right) ^{\alpha \left( 1_{H};1,0,1,0\right) }B\left( x_{2}\otimes
x_{1};GX_{1},gx_{1}x_{2}\right) GX_{1}^{1-l_{1}}\otimes gx_{2}\otimes
g^{l_{1}}x_{1}^{1+1}=0 \\
+\left( -1\right) ^{\alpha \left( 1_{H};0,0,0,1\right) }B\left( x_{2}\otimes
x_{1};GX_{1},gx_{1}x_{2}\right) GX_{1}\otimes gx_{1}\otimes x_{2}+ \\
+\left( -1\right) ^{\alpha \left( 1_{H};1,0,0,1\right) }B\left( x_{2}\otimes
x_{1};GX_{1},gx_{1}x_{2}\right) G\otimes gx_{1}\otimes gx_{1}x_{2}\text{ } \\
+\left( -1\right) ^{\alpha \left( 1_{H};0,0,1,1\right) }B\left( x_{2}\otimes
x_{1};GX_{1},gx_{1}x_{2}\right) GX_{1}\otimes g\otimes gx_{1}x_{2} \\
+\left( -1\right) ^{\alpha \left( 1_{H};1,0,1,1\right) }B\left( x_{2}\otimes
x_{1};GX_{1},gx_{1}x_{2}\right) GX_{1}^{1-l_{1}}\otimes g\otimes
g^{l_{1}+1}x_{1}^{1+1}x_{2}=0
\end{gather*}

\subsubsection{Case $G\otimes gx_{1}x_{2}\otimes x_{1}$}

Nothing from the first summand of the left side. From the second summand of
the left side, we deduce that

\begin{eqnarray*}
l_{1}+u_{1} &=&1 \\
l_{2} &=&u_{2}=0 \\
a &=&1,b_{1}=l_{1},b_{2}=0, \\
d &=&1,e_{1}-u_{1}=1\Rightarrow e_{1}=1,u_{1}=0,b_{1}=l_{1}=1,e_{2}=1.
\end{eqnarray*}%
Since $\alpha \left( 1_{H};1,0,0,0\right) \equiv b_{2}=0$ we get%
\begin{equation*}
B(x_{2}\otimes x_{1};GX_{1},gx_{1}x_{2})G\otimes gx_{1}x_{2}\otimes x_{1}.
\end{equation*}%
By considering also the right side, we get%
\begin{equation*}
B(x_{2}\otimes x_{1};GX_{1},gx_{1}x_{2})-B(x_{2}\otimes
1_{H};G,gx_{1}x_{2})=0
\end{equation*}%
which holds in view of the form of the elements.

\subsubsection{Case $GX_{1}\otimes gx_{2}\otimes x_{1}$}

Nothing from the first summand of the left side. From the second summand of
the left side, we deduce that

\begin{eqnarray*}
l_{1}+u_{1} &=&1 \\
l_{2} &=&u_{2}=0 \\
a &=&1,b_{1}-l_{1}=1\Rightarrow b_{1}=1,l_{1}=0,u_{1}=1,b_{2}=0, \\
d &=&1,e_{1}=u_{1}=1,e_{2}=1.
\end{eqnarray*}%
Since $\alpha \left( 1_{H};0,0,1,0\right) \equiv e_{2}+\left(
a+b_{1}+b_{2}\right) \equiv 1$ we get%
\begin{equation*}
-B(x_{2}\otimes x_{1};GX_{1},gx_{1}x_{2})GX_{1}\otimes gx_{2}\otimes x_{1}.
\end{equation*}%
By considering also the right side, we get

\begin{equation*}
-B(x_{2}\otimes x_{1};GX_{1},gx_{1}x_{2})-B(x_{2}\otimes
1_{H};GX_{1},gx_{2})=0
\end{equation*}%
which holds in view of the form of the elements.

\subsubsection{Case $GX_{1}\otimes gx_{1}\otimes x_{2}$}

From the first summand of the left side we get

\begin{eqnarray*}
l_{1} &=&u_{1}=0 \\
l_{2} &=&u_{2}=0 \\
a &=&b_{1}=1,b_{2}=0, \\
d &=&e_{1}=1,e_{2}=0
\end{eqnarray*}%
Since $\alpha \left( x_{2};0,0,0,0\right) =a+b_{1}+b_{2}\equiv 0$ we get%
\begin{equation*}
+B(g\otimes x_{1};GX_{1},gx_{1})GX_{1}\otimes gx_{1}\otimes x_{2}
\end{equation*}%
From the second summand of the left side we get

~%
\begin{eqnarray*}
l_{1} &=&u_{1}=0 \\
l_{2}+u_{2} &=&1 \\
a &=&b_{1}=1,b_{2}=l_{2} \\
d &=&e_{1}=1,e_{2}=u_{2}
\end{eqnarray*}%
since $\alpha \left( 1_{H};0,0,0,1\right) \equiv a+b_{1}+b_{2}\equiv 0$ and $%
\alpha \left( 1_{H};0,1,0,0\right) \equiv 0$ we obtain

\begin{equation*}
\left[ B(x_{2}\otimes x_{1};GX_{1},gx_{1}x_{2})+B(x_{2}\otimes
x_{1};GX_{1}X_{2},gx_{1})\right] GX_{1}\otimes 1_{H}\otimes x_{2}.
\end{equation*}%
Since there is nothing from the right side, we obtain

\begin{equation*}
B(g\otimes x_{1};GX_{1},gx_{1})+B(x_{2}\otimes
x_{1};GX_{1},gx_{1}x_{2})+B(x_{2}\otimes x_{1};GX_{1}X_{2},gx_{1})=0
\end{equation*}%
which holds in view of the form of the elements.

\subsubsection{Case $G\otimes gx_{1}\otimes gx_{1}x_{2}$}

From the first summand of the left side we get

\begin{eqnarray*}
l_{1}+u_{1} &=&1 \\
l_{2} &=&u_{2}=0 \\
a &=&1,b_{1}=l_{1},b_{2}=0, \\
d &=&1,e_{1}-u_{1}=1\Rightarrow e_{1}=1,u_{1}=0,b_{1}=l_{1}=1,e_{2}=0
\end{eqnarray*}%
Since $\alpha \left( x_{2};1,0,0,0\right) \equiv a+b_{1}\equiv 0,$ we get

\begin{equation*}
B(g\otimes x_{1};GX_{1},gx_{1})G\otimes gx_{1}\otimes gx_{1}x_{2}.
\end{equation*}%
From the second summand of the left side we get

\begin{eqnarray*}
l_{1}+u_{1} &=&1 \\
l_{2}+u_{2} &=&1 \\
a &=&1,b_{1}=l_{1},b_{2}=l_{2} \\
d &=&1,e_{1}-u_{1}=1\Rightarrow e_{1}=1,u_{1}=0,b_{1}=l_{1}=1,e_{2}=u_{2}
\end{eqnarray*}%
Since $\alpha \left( 1_{H};1,0,0,1\right) \equiv a+b_{1}\equiv 0$ and $%
\alpha \left( 1_{H};1,1,0,0\right) \equiv 1+b_{2}\equiv 0,$ we get%
\begin{equation*}
\left[ B(x_{1}\otimes x_{2};GX_{1},gx_{1}x_{2})+B(x_{1}\otimes
x_{2};GX_{1}X_{2},gx_{1})\right] X_{1}X_{2}\otimes g\otimes gx_{1}x_{2}.
\end{equation*}%
Since there is nothing from the right side, we obtain%
\begin{equation*}
B(g\otimes x_{1};GX_{1},gx_{1})+B(x_{1}\otimes
x_{2};GX_{1},gx_{1}x_{2})+B(x_{1}\otimes x_{2};GX_{1}X_{2},gx_{1})=0
\end{equation*}%
which holds in view of the form of the elements.

\subsubsection{Case $GX_{1}\otimes g\otimes gx_{1}x_{2}$}

From the first summand of the left side we get

\begin{eqnarray*}
l_{1}+u_{1} &=&1 \\
l_{2} &=&u_{2}=0 \\
a &=&1,b_{1}-l_{1}=1\Rightarrow b_{1}=1,l_{1}=0,u_{1}=1,b_{2}=0, \\
d &=&1,e_{1}=u_{1}=1,e_{2}=0
\end{eqnarray*}%
Since $\alpha \left( x_{2};0,0,1,0\right) \equiv e_{2}=0,$ we get

\begin{equation*}
B(g\otimes x_{1};GX_{1},gx_{1})GX_{1}\otimes g\otimes gx_{1}x_{2}.
\end{equation*}%
From the second summand of the left side we get

\begin{eqnarray*}
l_{1}+u_{1} &=&1 \\
l_{2}+u_{2} &=&1 \\
a &=&1,b_{1}-l_{1}=1\Rightarrow b_{1}=1,l_{1}=0,u_{1}=1,b_{2}=l_{2} \\
d &=&1,e_{1}=u_{1}=1,e_{2}=u_{2}
\end{eqnarray*}%
Since $\alpha \left( 1_{H};0,0,1,1\right) \equiv 1+e_{2}\equiv 0$ and $%
\alpha \left( 1_{H};0,1,1,0\right) \equiv e_{2}+a+b_{1}+b_{2}+1\equiv 0,$ we
get%
\begin{equation*}
\left[ B(x_{1}\otimes x_{2};GX_{1},gx_{1}x_{2})+B(x_{1}\otimes
x_{2};GX_{1}X_{2},gx_{1})\right] X_{1}X_{2}\otimes g\otimes gx_{1}x_{2}.
\end{equation*}%
Since there is nothing from the right side, we obtain%
\begin{equation*}
B(g\otimes x_{1};GX_{1},gx_{1})+B(x_{1}\otimes
x_{2};GX_{1},gx_{1}x_{2})+B(x_{1}\otimes x_{2};GX_{1}X_{2},gx_{1})=0
\end{equation*}%
which we just got in the previous case.

\subsection{$B(x_{2}\otimes x_{1};GX_{2},g)$}

We deduce that
\begin{equation*}
a=1,b_{2}=1,d=1
\end{equation*}%
and we obtain%
\begin{eqnarray*}
&&\left( -1\right) ^{\alpha \left( 1_{H};0,0,0,0\right) }B(x_{2}\otimes
x_{1};GX_{2},g)GX_{2}\otimes g\otimes g+ \\
&&\left( -1\right) ^{\alpha \left( 1_{H};0,1,0,0\right) }B(x_{2}\otimes
x_{1};GX_{2},g)G\otimes g\otimes x_{2}
\end{eqnarray*}

\subsubsection{Case $G\otimes g\otimes x_{2}$}

This case was already considered in subsection $B\left( x_{2}\otimes
x_{1};G,gx_{2}\right) .$

\subsection{$B(x_{2}\otimes x_{1};GX_{2},x_{1})$}

We deduce that%
\begin{equation*}
a=1,b_{2}=1,e_{1}=1
\end{equation*}%
and get%
\begin{eqnarray*}
&&\left( -1\right) ^{\alpha \left( 1_{H};0,0,0,0\right) }B(x_{2}\otimes
x_{1};GX_{2},x_{1})GX_{2}\otimes x_{1}\otimes g+ \\
&&+\left( -1\right) ^{\alpha \left( 1_{H};0,1,0,0\right) }B(x_{2}\otimes
x_{1};GX_{2},x_{1})G\otimes x_{1}\otimes x_{2}+ \\
&&+\left( -1\right) ^{\alpha \left( 1_{H};0,0,1,0\right) }B(x_{2}\otimes
x_{1};GX_{2},x_{1})GX_{2}\otimes 1_{H}\otimes x_{1}+ \\
&&+\left( -1\right) ^{\alpha \left( 1_{H};0,1,1,0\right) }B(x_{2}\otimes
x_{1};GX_{2},x_{1})G\otimes 1_{H}\otimes gx_{1}x_{2}
\end{eqnarray*}

\subsubsection{Case $G\otimes x_{1}\otimes x_{2}$}

This case was already considered in subsection $B\left( x_{2}\otimes
x_{1};G,x_{1}x_{2}\right) .$

\subsubsection{Case $GX_{2}\otimes 1_{H}\otimes x_{1}$}

Nothing from the first summand of the left side. From the second summand of
the left side, we deduce that

\begin{eqnarray*}
l_{1}+u_{1} &=&1 \\
l_{2} &=&u_{2}=0 \\
a &=&1,b_{1}=l_{1},b_{2}=1, \\
d &=&0,e_{1}=u_{1},e_{2}=0.
\end{eqnarray*}%
Since $\alpha \left( 1_{H};0,0,1,0\right) \equiv e_{2}+\left(
a+b_{1}+b_{2}\right) \equiv 0$ and $\alpha \left( 1_{H};1,0,0,0\right)
\equiv b_{2}=1$ we get%
\begin{equation*}
\left[ B(x_{2}\otimes x_{1};GX_{2},x_{1})-B(x_{2}\otimes
x_{1};GX_{1}X_{2},1_{H})\right] GX_{2}\otimes 1_{H}\otimes x_{1}.
\end{equation*}%
By considering also the right side, we get%
\begin{equation*}
B(x_{2}\otimes x_{1};GX_{2},x_{1})-B(x_{2}\otimes
x_{1};GX_{1}X_{2},1_{H})-B(x_{2}\otimes 1_{H};GX_{2},1_{H})=0
\end{equation*}%
which holds in view of the form of the elements.

\subsubsection{Case $G\otimes 1_{H}\otimes gx_{1}x_{2}$}

This case was already considered in subsection $B\left( x_{2}\otimes
x_{1};G,x_{1}x_{2}\right) .$

\subsection{$B(x_{2}\otimes x_{1};GX_{2},x_{2})$}

We deduce that%
\begin{equation*}
a=b_{2}=1,e_{2}=1
\end{equation*}%
and we get
\begin{gather*}
\left( -1\right) ^{\alpha \left( 1_{H};0,0,0,0\right) }B(x_{2}\otimes
x_{1};GX_{2},x_{2})GX_{2}\otimes x_{2}\otimes g+ \\
\left( -1\right) ^{\alpha \left( 1_{H};0,1,0,0\right) }B(x_{2}\otimes
x_{1};GX_{2},x_{2})G\otimes x_{2}\otimes x_{2}+ \\
\left( -1\right) ^{\alpha \left( 1_{H};0,0,0,1\right) }B(x_{2}\otimes
x_{1};GX_{2},x_{2})GX_{2}\otimes 1_{H}\otimes x_{2}+ \\
\left( -1\right) ^{\alpha \left( 1_{H};0,1,0,1\right) }B(x_{2}\otimes
x_{1};GX_{2},x_{2})GX_{2}^{1-l_{2}}\otimes 1_{H}\otimes
g^{l_{2}}x_{2}^{1+1}=0.
\end{gather*}

\subsubsection{Case $G\otimes x_{2}\otimes x_{2}$}

From the first summand of the left side we get

\begin{eqnarray*}
l_{1} &=&u_{1}=0 \\
l_{2} &=&u_{2}=0 \\
a &=&1,b_{1}=b_{2}=0, \\
d &=&e_{1}=0,e_{2}=1
\end{eqnarray*}%
Since $\alpha \left( x_{2};0,0,0,0\right) =a+b_{1}+b_{2}\equiv 1$ we get%
\begin{equation*}
-B(g\otimes x_{1};G,x_{2})G\otimes x_{2}\otimes x_{2}.
\end{equation*}%
From the second summand of the left side we get

~%
\begin{eqnarray*}
l_{1} &=&u_{1}=0 \\
l_{2}+u_{2} &=&1 \\
a &=&1,b_{1}=0,b_{2}=l_{2} \\
d &=&e_{1}=0,e_{2}-u_{2}=1\Rightarrow e_{2}=1,u_{2}=0,b_{2}=l_{2}=1.
\end{eqnarray*}%
Since $\alpha \left( 1_{H};0,1,0,0\right) \equiv 0$ we obtain

\begin{equation*}
B(x_{2}\otimes x_{1};GX_{2},x_{2})G\otimes x_{2}\otimes x_{2}.
\end{equation*}%
Since there is nothing from the right side, we obtain

\begin{equation*}
-B(g\otimes x_{1};G,x_{2})+B(x_{2}\otimes x_{1};GX_{2},x_{2})=0
\end{equation*}%
which holds in view of the form of the elements.

\subsection{$B(x_{2}\otimes x_{1};GX_{2},gx_{1}x_{2})$}

We deduce that%
\begin{equation*}
a=b_{2}=1,d=e_{1}=e_{2}=1
\end{equation*}%
and we get%
\begin{gather*}
\left( -1\right) ^{\alpha \left( 1_{H};0,0,0,0\right) }B(x_{2}\otimes
x_{1};GX_{2},gx_{1}x_{2})GX_{2}\otimes gx_{1}x_{2}\otimes g \\
\left( -1\right) ^{\alpha \left( 1_{H};0,1,0,0\right) }B(x_{2}\otimes
x_{1};GX_{2},gx_{1}x_{2})G\otimes gx_{1}x_{2}\otimes x_{2} \\
\left( -1\right) ^{\alpha \left( 1_{H};0,0,1,0\right) }B(x_{2}\otimes
x_{1};GX_{2},gx_{1}x_{2})GX_{2}\otimes gx_{2}\otimes x_{1} \\
\left( -1\right) ^{\alpha \left( 1_{H};0,1,1,0\right) }B(x_{2}\otimes
x_{1};GX_{2},gx_{1}x_{2})G\otimes gx_{2}\otimes gx_{1}x_{2} \\
\left( -1\right) ^{\alpha \left( 1_{H};0,0,0,1\right) }B(x_{2}\otimes
x_{1};GX_{2},gx_{1}x_{2})GX_{2}\otimes gx_{1}\otimes x_{2} \\
\left( -1\right) ^{\alpha \left( 1_{H};0,1,0,1\right) }B(x_{2}\otimes
x_{1};GX_{2},gx_{1}x_{2})GX_{2}^{1-l_{2}}\otimes gx_{1}\otimes
g^{l_{2}}x_{2}^{1+1}=0 \\
\left( -1\right) ^{\alpha \left( 1_{H};0,0,1,1\right) }B(x_{2}\otimes
x_{1};GX_{2},gx_{1}x_{2})GX_{2}\otimes g\otimes gx_{1}x_{2} \\
\left( -1\right) ^{\alpha \left( 1_{H};0,1,1,1\right) }B(x_{2}\otimes
x_{1};GX_{2},gx_{1}x_{2})GX_{2}^{1-l_{2}}\otimes g\otimes
g^{l_{2}+1}x_{1}x_{2}^{1+1}=0.
\end{gather*}

\subsubsection{Case $G\otimes gx_{1}x_{2}\otimes x_{2}$}

From the first summand of the left side we get

\begin{eqnarray*}
l_{1} &=&u_{1}=0 \\
l_{2} &=&u_{2}=0 \\
a &=&1,b_{1}=b_{2}=0, \\
d &=&e_{1}=e_{2}=1
\end{eqnarray*}%
Since $\alpha \left( x_{2};0,0,0,0\right) =a+b_{1}+b_{2}\equiv 1$ we get%
\begin{equation*}
-B(g\otimes x_{1};G,gx_{1}x_{2})G\otimes gx_{1}x_{2}\otimes x_{2}.
\end{equation*}%
From the second summand of the left side we get

~%
\begin{eqnarray*}
l_{1} &=&u_{1}=0 \\
l_{2}+u_{2} &=&1 \\
a &=&1,b_{1}=0,b_{2}=l_{2} \\
d &=&e_{1}=1,e_{2}-u_{2}=1\Rightarrow e_{2}=1,u_{2}=0,b_{2}=l_{2}=1.
\end{eqnarray*}%
Since $\alpha \left( 1_{H};0,1,0,0\right) \equiv 0$ we obtain

\begin{equation*}
B(x_{2}\otimes x_{1};GX_{2},gx_{1}x_{2})G\otimes gx_{1}x_{2}\otimes x_{2}.
\end{equation*}%
Since there is nothing from the right side, we obtain%
\begin{equation*}
-B(g\otimes x_{1};G,gx_{1}x_{2})+B(x_{2}\otimes x_{1};GX_{2},gx_{1}x_{2})=0
\end{equation*}

which holds in view of the form of the elements.

\subsubsection{Case $GX_{2}\otimes gx_{2}\otimes x_{1}$}

Nothing from the first summand of the left side. From the second summand of
the left side, we deduce that

\begin{eqnarray*}
l_{1}+u_{1} &=&1 \\
l_{2} &=&u_{2}=0 \\
a &=&1,b_{1}=l_{1},b_{2}=1, \\
d &=&e_{2}=1,e_{1}=u_{1}.
\end{eqnarray*}%
Since $\alpha \left( 1_{H};0,0,1,0\right) \equiv e_{2}+\left(
a+b_{1}+b_{2}\right) \equiv 1$ and $\alpha \left( 1_{H};1,0,0,0\right)
\equiv b_{2}=1$ we get%
\begin{equation*}
\left[ -B(x_{2}\otimes x_{1};GX_{2},gx_{1}x_{2})-B(x_{2}\otimes
x_{1};GX_{1}X_{2},gx_{2})\right] GX_{2}\otimes gx_{2}\otimes x_{1}.
\end{equation*}%
By considering also the right side, we get%
\begin{equation*}
-B(x_{2}\otimes x_{1};GX_{2},gx_{1}x_{2})-B(x_{2}\otimes
x_{1};GX_{1}X_{2},gx_{2})-B(x_{2}\otimes 1_{H};GX_{2},gx_{2})=0
\end{equation*}%
which holds in view of the form of the elements.

\subsubsection{Case $G\otimes gx_{2}\otimes gx_{1}x_{2}$}

From the first summand of the left side we get

\begin{eqnarray*}
l_{1}+u_{1} &=&1 \\
l_{2} &=&u_{2}=0 \\
a &=&1,b_{1}=l_{1},b_{2}=0, \\
d &=&1,e_{1}=u_{1},e_{2}=1
\end{eqnarray*}%
Since $\alpha \left( x_{2};0,0,1,0\right) \equiv e_{2}=1$ and $\alpha \left(
x_{2};1,0,0,0\right) \equiv a+b_{1}\equiv 0,$ we get

\begin{equation*}
\left[ -B(g\otimes x_{1};G,gx_{1}x_{2})+B(g\otimes x_{1};GX_{1},gx_{2})%
\right] G\otimes gx_{2}\otimes gx_{1}x_{2}.
\end{equation*}%
From the second summand of the left side we get

\begin{eqnarray*}
l_{1}+u_{1} &=&1 \\
l_{2}+u_{2} &=&1 \\
a &=&1,b_{1}=l_{1},b_{2}=l_{2} \\
d &=&1,e_{1}=u_{1},e_{2}-u_{2}=1\Rightarrow e_{2}=1,u_{2}=0,b_{2}=l_{2}=1
\end{eqnarray*}%
Since $\alpha \left( 1_{H};0,1,1,0\right) \equiv e_{2}+a+b_{1}+b_{2}+1\equiv
0$ and $\alpha \left( 1_{H};1,1,0,0\right) \equiv 1+b_{2}\equiv 0,$ we get%
\begin{equation*}
\left[ B(x_{2}\otimes x_{1};GX_{2},gx_{1}x_{2})+B(x_{2}\otimes
x_{1};GX_{1}X_{2},gx_{2})\right] X_{1}X_{2}\otimes g\otimes gx_{1}x_{2}.
\end{equation*}%
Since there is nothing from the right side, we obtain%
\begin{gather*}
-B(g\otimes x_{1};G,gx_{1}x_{2})+B(g\otimes x_{1};GX_{1},gx_{2})+ \\
B(x_{2}\otimes x_{1};GX_{2},gx_{1}x_{2})+B(x_{2}\otimes
x_{1};GX_{1}X_{2},gx_{2})=0
\end{gather*}%
which holds in view of the form of the elements.

\subsubsection{Case $GX_{2}\otimes gx_{1}\otimes x_{2}$}

From the first summand of the left side we get

\begin{eqnarray*}
l_{1} &=&u_{1}=0 \\
l_{2} &=&u_{2}=0 \\
a &=&b_{2}=1,b_{1}=0, \\
d &=&e_{1}=1,e_{2}=0
\end{eqnarray*}%
Since $\alpha \left( x_{2};0,0,0,0\right) =a+b_{1}+b_{2}\equiv 0$ we get%
\begin{equation*}
B(g\otimes x_{1};GX_{2},gx_{1})GX_{2}\otimes gx_{1}\otimes x_{2}.
\end{equation*}%
From the second summand of the left side we get

~%
\begin{eqnarray*}
l_{1} &=&u_{1}=0 \\
l_{2}+u_{2} &=&1 \\
a &=&1,b_{1}=0,b_{2}-l_{2}=1\Rightarrow b_{2}=1,l_{2}=0,u_{2}=1 \\
d &=&e_{1}=1,e_{2}=u_{2}=1.
\end{eqnarray*}%
Since $\alpha \left( 1_{H};0,0,0,1\right) \equiv a+b_{1}+b_{2}\equiv 0$ we
obtain

\begin{equation*}
B(x_{2}\otimes x_{1};GX_{2},gx_{1}x_{2})GX_{2}\otimes gx_{1}\otimes x_{2}.
\end{equation*}%
Since there is nothing from the right side, we obtain

\begin{equation*}
B(g\otimes x_{1};GX_{2},gx_{1})+B(x_{2}\otimes x_{1};GX_{2},gx_{1}x_{2})=0
\end{equation*}%
which holds in view of the form of the elements.

\subsubsection{Case $GX_{2}\otimes g\otimes gx_{1}x_{2}$}

From the first summand of the left side we get

\begin{eqnarray*}
l_{1}+u_{1} &=&1 \\
l_{2} &=&u_{2}=0 \\
a &=&1,b_{1}=l_{1},b_{2}=1, \\
d &=&1,e_{1}=u_{1},e_{2}=0
\end{eqnarray*}%
Since $\alpha \left( x_{2};0,0,1,0\right) \equiv e_{2}=0$ and $\alpha \left(
x_{2};1,0,0,0\right) \equiv a+b_{1}\equiv 0,$ we get

\begin{equation*}
\left[ B(g\otimes x_{1};GX_{2},gx_{1})+B(g\otimes x_{1};GX_{1}X_{2},g)\right]
GX_{2}\otimes g\otimes gx_{1}x_{2}.
\end{equation*}%
From the second summand of the left side we get

\begin{eqnarray*}
l_{1}+u_{1} &=&1 \\
l_{2}+u_{2} &=&1 \\
a &=&1,b_{1}=l_{1},b_{2}-l_{2}=1\Rightarrow b_{2}=1,l_{2}=0,u_{2}=1, \\
d &=&1,e_{1}=u_{1},e_{2}=u_{2}=1.
\end{eqnarray*}%
Since $\alpha \left( 1_{H};0,0,1,1\right) \equiv 1+e_{2}\equiv 0$ and $%
\alpha \left( 1_{H};1,0,0,1\right) \equiv a+b_{1}\equiv 0,$ we get%
\begin{equation*}
\left[ B(x_{2}\otimes x_{1};GX_{2},gx_{1}x_{2})+B(x_{2}\otimes
x_{1};GX_{1}X_{2},gx_{2})\right] X_{1}X_{2}\otimes g\otimes gx_{1}x_{2}.
\end{equation*}%
Since there is nothing from the right side, we obtain%
\begin{gather*}
B(g\otimes x_{1};GX_{2},gx_{1})+B(g\otimes x_{1};GX_{1}X_{2},g)+ \\
B(x_{2}\otimes x_{1};GX_{2},gx_{1}x_{2})+B(x_{2}\otimes
x_{1};GX_{1}X_{2},gx_{2})=0
\end{gather*}%
which holds in view of the form of the elements.

\subsection{$B(x_{2}\otimes x_{1};GX_{1}X_{2},1_{H})$}

We deduce that%
\begin{equation*}
a=b_{1}=b_{2}=1,d=e_{1}=e_{2}=0
\end{equation*}%
and we get%
\begin{eqnarray*}
&&\left( -1\right) ^{\alpha \left( 1_{H};0,0,0,0\right) }B(x_{2}\otimes
x_{1};GX_{1}X_{2},1_{H})GX_{1}X_{2}\otimes 1_{H}\otimes g+ \\
&&\left( -1\right) ^{\alpha \left( 1_{H};1,0,0,0\right) }B(x_{2}\otimes
x_{1};GX_{1}X_{2},1_{H})GX_{2}\otimes 1_{H}\otimes x_{1}+ \\
&&\left( -1\right) ^{\alpha \left( 1_{H};0,1,0,0\right) }B(x_{2}\otimes
x_{1};GX_{1}X_{2},1_{H})GX_{1}\otimes 1_{H}\otimes x_{2}+ \\
&&\left( -1\right) ^{\alpha \left( 1_{H};1,1,0,0\right) }B(x_{2}\otimes
x_{1};GX_{1}X_{2},1_{H})G\otimes 1_{H}\otimes gx_{1}x_{2}
\end{eqnarray*}

\subsubsection{Case $GX_{2}\otimes 1_{H}\otimes x_{1}$}

This case was already considered in section $B(x_{2}\otimes
x_{1};GX_{2},x_{1}).$

\subsubsection{Case $GX_{1}\otimes 1_{H}\otimes x_{2}$}

This case was already considered in section $B\left( x_{2}\otimes
x_{1};GX_{1},x_{2}\right) .$

\subsubsection{Case $G\otimes 1_{H}\otimes gx_{1}x_{2}$}

This case was already considered in section $B\left( x_{2}\otimes
x_{1};G,x_{1}x_{2}\right) .$

\subsection{$B(x_{2}\otimes x_{1};GX_{1}X_{2},x_{1}x_{2})$}

We deduce that%
\begin{eqnarray*}
a &=&b_{1}=b_{2}=1 \\
d &=&0,e_{1}=e_{2}=1
\end{eqnarray*}%
and we get%
\begin{eqnarray*}
&&\left( -1\right) ^{\alpha \left( 1_{H};l_{1},l_{2},u_{1},u_{2}\right)
}B(x_{2}\otimes
x_{1};GX_{1}X_{2},x_{1}x_{2})GX_{1}^{1-l_{1}}X_{2}^{1-l_{2}}\otimes
x_{1}^{1-u_{1}}x_{2}^{1-u_{2}} \\
&&\otimes g^{1+l_{1}+l_{2}+u_{1}+u_{2}}x_{1}^{l_{1}+u_{1}}x_{2}^{l_{2}+u_{2}}
\end{eqnarray*}%
\begin{gather*}
\left( -1\right) ^{\alpha \left( 1_{H};0,0,0,0\right) }B(x_{2}\otimes
x_{1};GX_{1}X_{2},x_{1}x_{2})GX_{1}X_{2}\otimes x_{1}x_{2}\otimes g \\
+\left( -1\right) ^{\alpha \left( 1_{H};1,0,0,0\right) }B(x_{2}\otimes
x_{1};GX_{1}X_{2},x_{1}x_{2})GX_{2}\otimes x_{1}x_{2}\otimes x_{1} \\
+\left( -1\right) ^{\alpha \left( 1_{H};0,1,0,0\right) }B(x_{2}\otimes
x_{1};GX_{1}X_{2},x_{1}x_{2})GX_{1}\otimes x_{1}x_{2}\otimes x_{2} \\
+\left( -1\right) ^{\alpha \left( 1_{H};1,1,0,0\right) }B(x_{2}\otimes
x_{1};GX_{1}X_{2},x_{1}x_{2})G\otimes x_{1}x_{2}\otimes gx_{1}x_{2} \\
+\left( -1\right) ^{\alpha \left( 1_{H};0,0,1,0\right) }B(x_{2}\otimes
x_{1};GX_{1}X_{2},x_{1}x_{2})GX_{1}X_{2}\otimes x_{2}\otimes x_{1} \\
+\left( -1\right) ^{\alpha \left( 1_{H};1,0,1,0\right) }B(x_{2}\otimes
x_{1};GX_{1}X_{2},x_{1}x_{2})GX_{1}^{1-l_{1}}X_{2}\otimes x_{2}\otimes
g^{l_{1}}x_{1}^{1+1}=0 \\
+\left( -1\right) ^{\alpha \left( 1_{H};0,1,1,0\right) }B(x_{2}\otimes
x_{1};GX_{1}X_{2},x_{1}x_{2})GX_{1}\otimes x_{2}\otimes gx_{1}x_{2} \\
+\left( -1\right) ^{\alpha \left( 1_{H};1,1,1,0\right) }B(x_{2}\otimes
x_{1};GX_{1}X_{2},x_{1}x_{2})GX_{1}^{1-l_{1}}\otimes x_{2}\otimes
g^{l_{1}+1}x_{1}^{1+1}x_{2}=0 \\
+\left( -1\right) ^{\alpha \left( 1_{H};0,0,0,1\right) }B(x_{2}\otimes
x_{1};GX_{1}X_{2},x_{1}x_{2})GX_{1}X_{2}\otimes x_{1}\otimes x_{2} \\
+\left( -1\right) ^{\alpha \left( 1_{H};1,0,0,1\right) }B(x_{2}\otimes
x_{1};GX_{1}X_{2},x_{1}x_{2})GX_{2}\otimes x_{1}\otimes gx_{1}x_{2} \\
+\left( -1\right) ^{\alpha \left( 1_{H};l_{1},1,0,1\right) }B(x_{2}\otimes
x_{1};GX_{1}X_{2},x_{1}x_{2})GX_{1}^{1-l_{1}}X_{2}^{1-l_{2}}\otimes
x_{1}\otimes g^{l_{1}+l_{2}}x_{1}^{l_{1}}x_{2}^{1+1}=0 \\
+\left( -1\right) ^{\alpha \left( 1_{H};0,0,1,1\right) }B(x_{2}\otimes
x_{1};GX_{1}X_{2},x_{1}x_{2})GX_{1}X_{2}\otimes 1_{H}\otimes gx_{1}x_{2} \\
+\left( -1\right) ^{\alpha \left( 1_{H};1,0,1,1\right) }B(x_{2}\otimes
x_{1};GX_{1}X_{2},x_{1}x_{2})GX_{1}^{1-l_{1}}X_{2}\otimes 1_{H}\otimes
g^{l_{1}+1}x_{1}^{1+1}x_{2}=0 \\
+\left( -1\right) ^{\alpha \left( 1_{H};1,1,1,1\right) }B(x_{2}\otimes
x_{1};GX_{1}X_{2},x_{1}x_{2})GX_{1}^{1-l_{1}}X_{2}^{1-l_{2}}\otimes
1_{H}\otimes g^{l_{1}+l_{2}+1}x_{1}^{l_{1}+1}x_{2}^{1+1}=0
\end{gather*}

\subsubsection{Case $GX_{2}\otimes x_{1}x_{2}\otimes x_{1}$}

Nothing from the first summand of the left side. From the second summand of
the left side, we deduce that

\begin{eqnarray*}
l_{1}+u_{1} &=&1 \\
l_{2} &=&u_{2}=0 \\
a &=&1,b_{1}=l_{1},b_{2}=1, \\
d &=&0,e_{2}=1,e_{1}-u_{1}=1\Rightarrow e_{1}=1,u_{1}=0,b_{1}=l_{1}=1.
\end{eqnarray*}%
Since $\alpha \left( 1_{H};1,0,0,0\right) \equiv b_{2}=1$ we get%
\begin{equation*}
-B(x_{2}\otimes x_{1};GX_{1}X_{2},x_{1}x_{2})GX_{2}\otimes x_{1}x_{2}\otimes
x_{1}.
\end{equation*}%
By considering also the right side, we get%
\begin{equation*}
-B(x_{2}\otimes x_{1};GX_{1}X_{2},x_{1}x_{2})-B(x_{2}\otimes
1_{H};GX_{2},x_{1}x_{2})=0
\end{equation*}%
which holds in view of the form of the elements.

\subsubsection{Case $GX_{1}\otimes x_{1}x_{2}\otimes x_{2}$}

From the first summand of the left side we get

\begin{eqnarray*}
l_{1} &=&u_{1}=0 \\
l_{2} &=&u_{2}=0 \\
a &=&b_{1}=1,b_{2}=0, \\
d &=&0,e_{1}=e_{2}=1
\end{eqnarray*}%
Since $\alpha \left( x_{2};0,0,0,0\right) =a+b_{1}+b_{2}\equiv 0$ we get%
\begin{equation*}
B(g\otimes x_{1};GX_{1},x_{1}x_{2})GX_{1}\otimes x_{1}x_{2}\otimes x_{2}.
\end{equation*}%
From the second summand of the left side we get

~%
\begin{eqnarray*}
l_{1} &=&u_{1}=0 \\
l_{2}+u_{2} &=&1 \\
a &=&1,b_{1}=1,b_{2}=l_{2} \\
d &=&0,e_{1}=1,e_{2}-u_{2}=1\Rightarrow e_{2}=1,u_{2}=0,b_{2}=l_{2}=1
\end{eqnarray*}%
Since $\alpha \left( 1_{H};0,1,0,0\right) \equiv $ 0 we obtain

\begin{equation*}
B(x_{2}\otimes x_{1};GX_{1}X_{2},x_{1}x_{2})GX_{1}\otimes x_{1}x_{2}\otimes
x_{2}.
\end{equation*}%
Since there is nothing from the right side, we obtain

\begin{equation*}
B(g\otimes x_{1};GX_{1},x_{1}x_{2})+B(x_{2}\otimes
x_{1};GX_{1}X_{2},x_{1}x_{2})=0
\end{equation*}%
which holds in view of the form of the elements.

\subsubsection{Case $G\otimes x_{1}x_{2}\otimes gx_{1}x_{2}$}

From the first summand of the left side we get

\begin{eqnarray*}
l_{1}+u_{1} &=&1 \\
l_{2} &=&u_{2}=0 \\
a &=&1,b_{1}=l_{1},b_{2}=0, \\
d &=&0,e_{1}-u_{1}=1\Rightarrow e_{1}=1,u_{1}=0,b_{1}=l_{1}=1,e_{2}=1
\end{eqnarray*}%
Since $\alpha \left( x_{2};1,0,0,0\right) \equiv a+b_{1}\equiv 0,$ we get

\begin{equation*}
B(g\otimes x_{1};GX_{1},x_{1}x_{2})G\otimes x_{1}x_{2}\otimes gx_{1}x_{2}.
\end{equation*}%
From the second summand of the left side we get

\begin{eqnarray*}
l_{1}+u_{1} &=&1 \\
l_{2}+u_{2} &=&1 \\
a &=&1,b_{1}=l_{1},b_{2}=l_{2}, \\
d &=&0,e_{1}-u_{1}=1\Rightarrow e_{1}=1,u_{1}=0,b_{1}=l_{1}=1, \\
e_{2}-u_{2} &=&1\Rightarrow e_{2}=1,u_{2}=0,b_{2}=l_{2}=1.
\end{eqnarray*}%
Since $\alpha \left( 1_{H};1,1,0,0\right) \equiv 1+b_{2}\equiv 0,$ we get%
\begin{equation*}
B(x_{2}\otimes x_{1};GX_{1}X_{2},x_{1}x_{2})G\otimes x_{1}x_{2}\otimes
gx_{1}x_{2}.
\end{equation*}%
Since there is nothing from the right side, we obtain%
\begin{equation*}
B(g\otimes x_{1};GX_{1},x_{1}x_{2})+B(x_{2}\otimes
x_{1};GX_{1}X_{2},x_{1}x_{2})=0
\end{equation*}%
which we just got above.

\subsubsection{Case $GX_{1}X_{2}\otimes x_{2}\otimes x_{1}$}

Nothing from the first summand of the left side. From the second summand of
the left side, we deduce that

\begin{eqnarray*}
l_{1}+u_{1} &=&1 \\
l_{2} &=&u_{2}=0 \\
a &=&1,b_{1}-l_{1}=1\Rightarrow b_{1}=1,l_{1}=0,u_{1}=1,b_{2}=1, \\
d &=&0,e_{2}=1,e_{1}=u_{1}=1.
\end{eqnarray*}%
Since $\alpha \left( 1_{H};0,0,1,0\right) \equiv e_{2}+\left(
a+b_{1}+b_{2}\right) \equiv 0$ we get%
\begin{equation*}
+B(x_{2}\otimes x_{1};GX_{1}X_{2},x_{1}x_{2})GX_{1}X_{2}\otimes x_{2}\otimes
x_{1}.
\end{equation*}%
By considering also the right side, we get%
\begin{equation*}
+B(x_{2}\otimes x_{1};GX_{1}X_{2},x_{1}x_{2})-B(x_{2}\otimes
1_{H};GX_{1}X_{2},x_{2})=0
\end{equation*}%
which holds in view of the form of the elements.

\subsubsection{Case $GX_{1}\otimes x_{2}\otimes gx_{1}x_{2}$}

From the first summand of the left side we get

\begin{eqnarray*}
l_{1}+u_{1} &=&1 \\
l_{2} &=&u_{2}=0 \\
a &=&1,b_{1}-l_{1}=1\Rightarrow b_{1}=1,l_{1}=0,u_{1}=1,b_{2}=0, \\
d &=&0,e_{1}=u_{1}=1,e_{2}=1
\end{eqnarray*}%
Since $\alpha \left( x_{2};0,0,1,0\right) \equiv e_{2}=1,$ we get

\begin{equation*}
-B(g\otimes x_{1};GX_{1},x_{1}x_{2})GX_{1}\otimes x_{2}\otimes gx_{1}x_{2}.
\end{equation*}%
From the second summand of the left side we get

\begin{eqnarray*}
l_{1}+u_{1} &=&1 \\
l_{2}+u_{2} &=&1 \\
a &=&1,b_{1}-l_{1}=1\Rightarrow b_{1}=1,l_{1}=0,u_{1}=1,b_{2}=l_{2}, \\
d &=&0,e_{1}=u_{1}=1, \\
e_{2}-u_{2} &=&1\Rightarrow e_{2}=1,u_{2}=0,b_{2}=l_{2}=1.
\end{eqnarray*}%
Since $\alpha \left( 1_{H};0,1,1,0\right) \equiv e_{2}+a+b_{1}+b_{2}+1\equiv
1,$ we get%
\begin{equation*}
-B(x_{2}\otimes x_{1};GX_{1}X_{2},x_{1}x_{2})G\otimes x_{1}x_{2}\otimes
gx_{1}x_{2}.
\end{equation*}%
Since there is nothing from the right side, we obtain%
\begin{equation*}
-B(g\otimes x_{1};GX_{1},x_{1}x_{2})-B(x_{2}\otimes
x_{1};GX_{1}X_{2},x_{1}x_{2})=0
\end{equation*}%
which we just got in Case $G\otimes x_{1}x_{2}\otimes gx_{1}x_{2}.$

\subsubsection{Case $GX_{1}X_{2}\otimes x_{1}\otimes x_{2}$}

From the first summand of the left side we get

\begin{eqnarray*}
l_{1} &=&u_{1}=0 \\
l_{2} &=&u_{2}=0 \\
a &=&b_{1}=b_{2}=1, \\
d &=&0,e_{1}=1,e_{2}=0
\end{eqnarray*}%
Since $\alpha \left( x_{2};0,0,0,0\right) =a+b_{1}+b_{2}\equiv 1$ we get%
\begin{equation*}
-B(g\otimes x_{1};GX_{1}X_{2},x_{1})GX_{1}X_{2}\otimes x_{1}\otimes x_{2}.
\end{equation*}%
From the second summand of the left side we get

~%
\begin{eqnarray*}
l_{1} &=&u_{1}=0 \\
l_{2}+u_{2} &=&1 \\
a &=&1,b_{1}=1,b_{2}-l_{2}=1\Rightarrow b_{2}=1,l_{2}=0,u_{2}=1 \\
d &=&0,e_{1}=1,e_{2}=u_{2}=1
\end{eqnarray*}%
Since $\alpha \left( 1_{H};0,0,0,1\right) \equiv a+b_{1}+b_{2}\equiv 1,$ we
obtain

\begin{equation*}
-B(x_{2}\otimes x_{1};GX_{1}X_{2},x_{1}x_{2})GX_{1}\otimes x_{1}x_{2}\otimes
x_{2}.
\end{equation*}%
Since there is nothing from the right side, we obtain%
\begin{equation*}
-B(g\otimes x_{1};GX_{1}X_{2},x_{1})-B(x_{2}\otimes
x_{1};GX_{1}X_{2},x_{1}x_{2})=0
\end{equation*}%
which holds in view of the form of the elements.

\subsubsection{Case $GX_{2}\otimes x_{1}\otimes gx_{1}x_{2}$}

From the first summand of the left side we get

\begin{eqnarray*}
l_{1}+u_{1} &=&1 \\
l_{2} &=&u_{2}=0 \\
a &=&1,b_{1}=l_{1},b_{2}=1, \\
d &=&0,e_{1}-u_{1}=1\Rightarrow e_{1}=1,u_{1}=0,b_{1}=l_{1}=1,e_{2}=0
\end{eqnarray*}%
Since $\alpha \left( x_{2};1,0,0,0\right) \equiv a+b_{1}\equiv 0$ we get

\begin{equation*}
+B(g\otimes x_{1};GX_{1}X_{2},x_{1})GX_{2}\otimes x_{1}\otimes gx_{1}x_{2}.
\end{equation*}%
From the second summand of the left side we get

\begin{eqnarray*}
l_{1}+u_{1} &=&1 \\
l_{2}+u_{2} &=&1 \\
a &=&1,b_{1}=l_{1},b_{2}-l_{2}=1\Rightarrow b_{2}=1,l_{2}=0,u_{2}=1, \\
d &=&0,e_{1}-u_{1}=1\Rightarrow e_{1}=1,u_{1}=0,b_{1}=l_{1}=1, \\
e_{2} &=&u_{2}=1.
\end{eqnarray*}%
Since $\alpha \left( 1_{H};1,0,0,1\right) \equiv a+b_{1}\equiv 0,$ we get%
\begin{equation*}
+B(x_{2}\otimes x_{1};GX_{1}X_{2},x_{1}x_{2})GX_{2}\otimes x_{1}\otimes
gx_{1}x_{2}.
\end{equation*}%
Since there is nothing from the right side, we obtain

\begin{equation*}
+B(g\otimes x_{1};GX_{1}X_{2},x_{1})+B(x_{2}\otimes
x_{1};GX_{1}X_{2},x_{1}x_{2})=0
\end{equation*}%
which we just got in Case $GX_{1}X_{2}\otimes x_{1}\otimes x_{2}.$

\subsubsection{Case $GX_{1}X_{2}\otimes 1_{H}\otimes gx_{1}x_{2}$}

From the first summand of the left side we get

\begin{eqnarray*}
l_{1}+u_{1} &=&1 \\
l_{2} &=&u_{2}=0 \\
a &=&1,b_{1}-l_{1}=1\Rightarrow b_{1}=1,l_{1}=0,u_{1}=1,b_{2}=1, \\
d &=&0,e_{1}=u_{1}=1,e_{2}=0
\end{eqnarray*}%
Since $\alpha \left( x_{2};1,0,0,0\right) \equiv a+b_{1}\equiv 0,$ we get

\begin{equation*}
B(g\otimes x_{1};GX_{1}X_{2},x_{1})GX_{1}X_{2}\otimes 1_{H}\otimes
gx_{1}x_{2}.
\end{equation*}%
From the second summand of the left side we get

\begin{eqnarray*}
l_{1}+u_{1} &=&1 \\
l_{2}+u_{2} &=&1 \\
a &=&1,b_{1}-l_{1}=1\Rightarrow b_{1}=1,l_{1}=0,u_{1}=1, \\
b_{2}-l_{2} &=&1\Rightarrow b_{2}=1,l_{2}=0,u_{2}=1 \\
d &=&0,e_{1}=u_{1}=1,e_{2}=u_{2}=1.
\end{eqnarray*}%
Since $\alpha \left( 1_{H};0,0,1,1\right) \equiv 1+e_{2}\equiv 0,$ we get%
\begin{equation*}
B(x_{2}\otimes x_{1};GX_{1}X_{2},x_{1}x_{2})GX_{1}X_{2}\otimes 1_{H}\otimes
gx_{1}x_{2}.
\end{equation*}%
Since there is nothing from the right side, we obtain

\begin{equation*}
+B(g\otimes x_{1};GX_{1}X_{2},x_{1})+B(x_{2}\otimes
x_{1};GX_{1}X_{2},x_{1}x_{2})=0
\end{equation*}%
which we just got in Case $GX_{1}X_{2}\otimes x_{1}\otimes x_{2}.$

\subsection{$B(x_{2}\otimes x_{1};GX_{1}X_{2},gx_{1})$}

We deduce that%
\begin{equation*}
a=b_{1}=b_{2}=1,d=e_{1}=1
\end{equation*}%
and we get%
\begin{equation*}
\left( -1\right) ^{\alpha \left( 1_{H};l_{1},l_{2},u_{1},0\right)
}B(x_{2}\otimes
x_{1};GX_{1}X_{2},gx_{1})GX_{1}^{1-l_{1}}X_{2}^{1-l_{2}}\otimes
gx_{1}^{1-u_{1}}\otimes
g^{1+l_{1}+l_{2}+u_{1}}x_{1}^{l_{1}+u_{1}}x_{2}^{l_{2}}
\end{equation*}%
\begin{gather*}
\left( -1\right) ^{\alpha \left( 1_{H};0,0,0,0\right) }B(x_{2}\otimes
x_{1};GX_{1}X_{2},gx_{1})GX_{1}X_{2}\otimes gx_{1}\otimes g \\
\left( -1\right) ^{\alpha \left( 1_{H};1,0,0,0\right) }B(x_{2}\otimes
x_{1};GX_{1}X_{2},gx_{1})GX_{2}\otimes gx_{1}\otimes x_{1} \\
\left( -1\right) ^{\alpha \left( 1_{H};0,1,0,0\right) }B(x_{2}\otimes
x_{1};GX_{1}X_{2},gx_{1})GX_{1}\otimes gx_{1}\otimes x_{2} \\
\left( -1\right) ^{\alpha \left( 1_{H};1,1,0,0\right) }B(x_{2}\otimes
x_{1};GX_{1}X_{2},gx_{1})G\otimes gx_{1}\otimes gx_{1}x_{2} \\
\left( -1\right) ^{\alpha \left( 1_{H};0,0,1,0\right) }B(x_{2}\otimes
x_{1};GX_{1}X_{2},gx_{1})GX_{1}X_{2}\otimes g\otimes x_{1} \\
\left( -1\right) ^{\alpha \left( 1_{H};1,0,1,0\right) }B(x_{2}\otimes
x_{1};GX_{1}X_{2},gx_{1})GX_{1}^{1-l_{1}}X_{2}^{1-l_{2}}\otimes
gx_{1}^{1-u_{1}}\otimes g^{1+l_{1}+l_{2}+u_{1}}x_{1}^{1+1}x_{2}^{l_{2}}=0 \\
\left( -1\right) ^{\alpha \left( 1_{H};0,1,1,0\right) }B(x_{2}\otimes
x_{1};GX_{1}X_{2},gx_{1})GX_{1}\otimes g\otimes gx_{1}x_{2} \\
\left( -1\right) ^{\alpha \left( 1_{H};1,1,1,0\right) }B(x_{2}\otimes
x_{1};GX_{1}X_{2},gx_{1})GX_{1}^{1-l_{1}}X_{2}^{1-l_{2}}\otimes
gx_{1}^{1-u_{1}}\otimes g^{1+l_{1}+l_{2}+u_{1}}x_{1}^{1+1}x_{2}^{l_{2}}=0.
\end{gather*}

\subsubsection{Case $GX_{2}\otimes gx_{1}\otimes x_{1}$}

Nothing from the first summand of the left side. From the second summand of
the left side, we deduce that

\begin{eqnarray*}
l_{1}+u_{1} &=&1 \\
l_{2} &=&u_{2}=0 \\
a &=&1,b_{1}=l_{1},b_{2}=1, \\
d &=&1,e_{1}-u_{1}=1\Rightarrow e_{1}=1,u_{1}=0,b_{1}=l_{1}=1,e_{2}=0.
\end{eqnarray*}%
Since $\alpha \left( 1_{H};1,0,0,0\right) \equiv b_{2}\equiv 1$ we get%
\begin{equation*}
-B(x_{2}\otimes x_{1};GX_{1}X_{2},gx_{1})GX_{2}\otimes gx_{1}\otimes x_{1}.
\end{equation*}%
By considering also the right side, we get%
\begin{equation*}
-B(x_{2}\otimes x_{1};GX_{1}X_{2},gx_{1})-B(x_{2}\otimes
1_{H};GX_{2},gx_{1})=0
\end{equation*}%
which holds in view of the form of the elements.

\subsubsection{Case $GX_{1}\otimes gx_{1}\otimes x_{2}$}

This case was already considered in subsection $B\left( x_{2}\otimes
x_{1};GX_{1},gx_{1}x_{2}\right) .$

\subsubsection{Case $G\otimes gx_{1}\otimes gx_{1}x_{2}$}

This case was already considered in subsection $B\left( x_{2}\otimes
x_{1};GX_{1},gx_{1}x_{2}\right) .$

\subsubsection{Case $GX_{1}X_{2}\otimes g\otimes x_{1}$}

Nothing from the first summand of the left side. From the second summand of
the left side, we deduce that

\begin{eqnarray*}
l_{1}+u_{1} &=&1 \\
l_{2} &=&u_{2}=0 \\
a &=&1,b_{1}-l_{1}=1\Rightarrow b_{1}=1,l_{1}=0,u_{1}=1,b_{2}=1, \\
d &=&1,e_{1}=u_{1}=1,e_{2}=0.
\end{eqnarray*}%
Since $\alpha \left( 1_{H};0,0,1,0\right) \equiv e_{2}+\left(
a+b_{1}+b_{2}\right) \equiv 1$ we get%
\begin{equation*}
-B(x_{2}\otimes x_{1};GX_{1}X_{2},gx_{1})GX_{2}\otimes gx_{1}\otimes x_{1}.
\end{equation*}%
By considering also the right side, we get%
\begin{equation*}
-B(x_{2}\otimes x_{1};GX_{1}X_{2},gx_{1})-B(x_{2}\otimes
1_{H};GX_{2},gx_{1})=0
\end{equation*}%
which we already got in case $GX_{2}\otimes gx_{1}\otimes x_{1}.$

\subsubsection{Case $GX_{1}\otimes g\otimes gx_{1}x_{2}$}

This case was already considered in case $B\left( x_{2}\otimes
x_{1};GX_{1},gx_{1}x_{2}\right) .$

\subsection{$B(x_{2}\otimes x_{1};GX_{1}X_{2},gx_{2})$}

We deduce that
\begin{equation*}
a=b_{1}=b_{2}=1,d=e_{2}=1
\end{equation*}%
and we get%
\begin{equation*}
\left( -1\right) ^{\alpha \left( 1_{H};l_{1},l_{2},0,u_{2}\right)
}B(x_{2}\otimes
x_{1};GX_{1}X_{2},gx_{2})GX_{1}^{1-l_{1}}X_{2}^{1-l_{2}}\otimes
gx_{2}^{1-u_{2}}\otimes
g^{1+l_{1}+l_{2}+u_{2}}x_{1}^{l_{1}}x_{2}^{l_{2}+u_{2}}
\end{equation*}%
\begin{gather*}
\left( -1\right) ^{\alpha \left( 1_{H};0,0,0,0\right) }B(x_{2}\otimes
x_{1};GX_{1}X_{2},gx_{2})GX_{1}X_{2}\otimes gx_{2}\otimes g \\
\left( -1\right) ^{\alpha \left( 1_{H};1,0,0,0\right) }B(x_{2}\otimes
x_{1};GX_{1}X_{2},gx_{2})GX_{2}\otimes gx_{2}\otimes x_{1} \\
\left( -1\right) ^{\alpha \left( 1_{H};0,1,0,0\right) }B(x_{2}\otimes
x_{1};GX_{1}X_{2},gx_{2})GX_{1}\otimes gx_{2}\otimes x_{2} \\
\left( -1\right) ^{\alpha \left( 1_{H};1,1,0,0\right) }B(x_{2}\otimes
x_{1};GX_{1}X_{2},gx_{2})G\otimes gx_{2}\otimes gx_{1}x_{2} \\
\left( -1\right) ^{\alpha \left( 1_{H};0,0,0,1\right) }B(x_{2}\otimes
x_{1};GX_{1}X_{2},gx_{2})GX_{1}X_{2}\otimes g\otimes x_{2} \\
\left( -1\right) ^{\alpha \left( 1_{H};1,0,0,1\right) }B(x_{2}\otimes
x_{1};GX_{1}X_{2},gx_{2})GX_{2}\otimes g\otimes gx_{1}x_{2} \\
\left( -1\right) ^{\alpha \left( 1_{H};0,1,0,1\right) }B(x_{2}\otimes
x_{1};GX_{1}X_{2},gx_{2})GX_{1}\otimes g\otimes gx_{2}^{1+1}=0 \\
\left( -1\right) ^{\alpha \left( 1_{H};1,1,0,1\right) }B(x_{2}\otimes
x_{1};GX_{1}X_{2},gx_{2})GX_{1}^{1-l_{1}}X_{2}^{1-l_{2}}\otimes
gx_{2}^{1-u_{2}}\otimes g^{1+l_{1}+l_{2}+u_{2}}x_{1}x_{2}^{1+1}=0.
\end{gather*}

\subsubsection{Case $GX_{2}\otimes gx_{2}\otimes x_{1}$}

This Case already appeared in subsection $B(x_{2}\otimes
x_{1};GX_{2},gx_{1}x_{2}).$

\subsubsection{Case $GX_{1}\otimes gx_{2}\otimes x_{2}$}

From the first summand of the left side we get

\begin{eqnarray*}
l_{1} &=&u_{1}=0 \\
l_{2} &=&u_{2}=0 \\
a &=&b_{1}=1,b_{2}=0 \\
d &=&1,e_{1}=0,e_{2}=1
\end{eqnarray*}%
Since $\alpha \left( x_{2};0,0,0,0\right) =a+b_{1}+b_{2}\equiv 0$ we get%
\begin{equation*}
B(g\otimes x_{1};GX_{1},gx_{2})GX_{1}\otimes gx_{2}\otimes x_{2}.
\end{equation*}%
From the second summand of the left side we get

~%
\begin{eqnarray*}
l_{1} &=&u_{1}=0 \\
l_{2}+u_{2} &=&1 \\
a &=&1,b_{1}=1,b_{2}=l_{2} \\
d &=&1,e_{1}=0,e_{2}-u_{2}=1\Rightarrow e_{2}=1,u_{2}=0,b_{2}=l_{2}=1
\end{eqnarray*}%
Since \ $\alpha \left( 1_{H};0,1,0,0\right) \equiv 0,$ we obtain

\begin{equation*}
B(x_{2}\otimes x_{1};GX_{1}X_{2},gx_{2})GX_{1}\otimes gx_{2}\otimes x_{2}.
\end{equation*}%
Since there is nothing from the right side, we obtain%
\begin{equation*}
B(g\otimes x_{1};GX_{1},gx_{2})+B(x_{2}\otimes x_{1};GX_{1}X_{2},gx_{2})=0
\end{equation*}%
which holds in view of the form of the elements.

\subsubsection{Case $G\otimes gx_{2}\otimes gx_{1}x_{2}$}

This case was already considered in subsection $B(x_{2}\otimes
x_{1};GX_{2},gx_{1}x_{2}).$

\subsubsection{Case $GX_{1}X_{2}\otimes g\otimes x_{2}$}

From the first summand of the left side we get

\begin{eqnarray*}
l_{1} &=&u_{1}=0 \\
l_{2} &=&u_{2}=0 \\
a &=&b_{1}=b_{2}=1, \\
d &=&1,e_{1}=0,e_{2}=0
\end{eqnarray*}%
Since $\alpha \left( x_{2};0,0,0,0\right) =a+b_{1}+b_{2}\equiv 1$ we get%
\begin{equation*}
-B(g\otimes x_{1};GX_{1}X_{2},g)GX_{1}X_{2}\otimes g\otimes x_{2}.
\end{equation*}%
From the second summand of the left side we get

~%
\begin{eqnarray*}
l_{1} &=&u_{1}=0 \\
l_{2}+u_{2} &=&1 \\
a &=&1,b_{1}=1,b_{2}-l_{2}=1\Rightarrow b_{2}=1,l_{2}=0,u_{2}=1 \\
d &=&1,e_{1}=0,e_{2}=u_{2}=1.
\end{eqnarray*}%
Since $\alpha \left( 1_{H};0,0,0,1\right) \equiv a+b_{1}+b_{2}\equiv 1,$ we
obtain

\begin{equation*}
-B(x_{2}\otimes x_{1};GX_{1}X_{2},gx_{2})GX_{1}X_{2}\otimes g\otimes x_{2}.
\end{equation*}%
Since there is nothing from the right side, we obtain%
\begin{equation*}
-B(g\otimes x_{1};GX_{1}X_{2},g)-B(x_{2}\otimes x_{1};GX_{1}X_{2},gx_{2})=0
\end{equation*}%
which holds in view of the form of the elements.

\subsubsection{Case $GX_{2}\otimes g\otimes gx_{1}x_{2}$}

This Case was already considered in subsection $B(x_{2}\otimes
x_{1};GX_{2},gx_{1}x_{2}).$

\section{$B\left( x_{2}\otimes x_{2}\right) $}

From $\left( \ref{simplx}\right) $ we get%
\begin{equation}
B(x_{2}\otimes x_{2})=B(x_{2}\otimes 1_{H})(1_{A}\otimes
x_{2})-(1_{A}\otimes gx_{2})B(x_{2}\otimes 1_{H})(1_{A}\otimes g).
\label{form x2otx2}
\end{equation}%
and we obtain%
\begin{eqnarray*}
B\left( x_{2}\otimes x_{2}\right) &=&2B\left( x_{2}\otimes
1_{H};1_{A},gx_{1}\right) 1_{A}\otimes gx_{1}x_{2}+ \\
&&-2B\left( x_{2}\otimes 1_{H};G,g\right) G\otimes gx_{2} \\
&&+2B(x_{2}\otimes 1_{H};1_{A},gx_{1})X_{1}\otimes gx_{2}+ \\
&&+2\left[ +B(g\otimes 1_{H};1_{A},g)+B(x_{2}\otimes \ 1_{H};1_{A},gx_{2})%
\right] X_{2}\otimes gx_{2}+ \\
&&+2\left[ B(g\otimes 1_{H};G,gx_{1})+B(x_{2}\otimes \ 1_{H};G,gx_{1}x_{2})%
\right] GX_{2}\otimes gx_{1}x_{2}+ \\
&&+2\left[ -B(g\otimes 1_{H};G,gx_{1})-B(x_{2}\otimes \ 1_{H};G,gx_{1}x_{2})%
\right] GX_{1}X_{2}\otimes gx_{2}+
\end{eqnarray*}

\subsection{$B\left( x_{2}\otimes x_{2};1_{A},gx_{1}x_{2}\right) $}

We deduce that%
\begin{equation*}
a=0,d=e_{1}=e_{2}=1
\end{equation*}%
and we get%
\begin{eqnarray*}
&&\left( -1\right) ^{\alpha \left( 1_{H};0,0,0,0\right) }B\left(
x_{2}\otimes x_{2};1_{A},gx_{1}x_{2}\right) 1_{A}\otimes gx_{1}x_{2}\otimes g
\\
&&\left( -1\right) ^{\alpha \left( 1_{H};0,0,1,0\right) }B\left(
x_{2}\otimes x_{2};1_{A},gx_{1}x_{2}\right) 1_{A}\otimes gx_{2}\otimes x_{1}
\\
&&\left( -1\right) ^{\alpha \left( 1_{H};0,0,0,1\right) }B\left(
x_{2}\otimes x_{2};1_{A},gx_{1}x_{2}\right) 1_{A}\otimes gx_{1}\otimes x_{2}
\\
&&\left( -1\right) ^{\alpha \left( 1_{H};0,0,1,1\right) }B\left(
x_{2}\otimes x_{2};1_{A},gx_{1}x_{2}\right) 1_{A}\otimes g\otimes gx_{1}x_{2}
\end{eqnarray*}

\subsubsection{Case $1_{A}\otimes gx_{1}x_{2}\otimes g$}

This case is trivial.

\subsubsection{Case $1_{A}\otimes gx_{2}\otimes x_{1}$}

Nothing from the first summand of the left side. The second summand of the
left side gives us%
\begin{eqnarray*}
l_{1}+u_{1} &=&1,l_{2}=u_{2}=0 \\
a &=&0,b_{1}=l_{1},b_{2}=0 \\
d &=&1,e_{2}=1,e_{1}=u_{1}
\end{eqnarray*}%
Since $\alpha \left( 1_{H};0,0,1,0\right) \equiv e_{2}+\left(
a+b_{1}+b_{2}\right) \equiv 0$ and $\alpha \left( 1_{H};1,0,0,0\right)
\equiv b_{2}=0,$ we get%
\begin{equation*}
\left[ -B(x_{2}\otimes x_{2};1_{A},gx_{1}x_{2})+B(x_{2}\otimes
x_{2};X_{1},gx_{2})\right] 1_{A}\otimes gx_{2}\otimes x_{1}.
\end{equation*}%
Since there is nothing in the right side, we obtain%
\begin{equation*}
-B(x_{2}\otimes x_{2};1_{A},gx_{1}x_{2})+B(x_{2}\otimes x_{2};X_{1},gx_{2})=0
\end{equation*}%
which holds in view of the form of the element.

\subsubsection{Case $1_{A}\otimes gx_{1}\otimes x_{2}$}

From the first summand of the left side we get

\begin{eqnarray*}
l_{1} &=&u_{1}=0 \\
l_{2} &=&u_{2}=0 \\
a &=&b_{1}=b_{2}=0 \\
d &=&1,e_{1}=1,e_{2}=0
\end{eqnarray*}%
Since $\alpha \left( x_{2};0,0,0,0\right) =a+b_{1}+b_{2}\equiv 0$ we get%
\begin{equation*}
B(g\otimes x_{2};1_{A},gx_{1})1_{A}\otimes gx_{1}\otimes x_{2}.
\end{equation*}%
From the second summand of the left side we get

~%
\begin{eqnarray*}
l_{1} &=&u_{1}=0 \\
l_{2}+u_{2} &=&1 \\
a &=&b_{1}=0,b_{2}=l_{2} \\
d &=&1,e_{1}=1,e_{2}=u_{2}
\end{eqnarray*}%
Since $\alpha \left( 1_{H};0,0,0,1\right) \equiv a+b_{1}+b_{2}\equiv 0$\ and
$\alpha \left( 1_{H};0,1,0,0\right) \equiv 0$ and we obtain

\begin{equation*}
\left[ B(x_{2}\otimes x_{2};X_{2},gx_{1})+B(x_{2}\otimes
x_{2};1_{A},gx_{1}x_{2})\right] 1_{A}\otimes gx_{1}\otimes x_{2}.
\end{equation*}%
Considering also the right side, we get%
\begin{gather*}
B(g\otimes x_{2};1_{A},gx_{1})+B(x_{2}\otimes x_{2};X_{2},gx_{1}) \\
+B(x_{2}\otimes x_{2};1_{A},gx_{1}x_{2})-B(x_{2}\otimes 1_{H};1_{A},gx_{1})=0
\end{gather*}%
which holds in view of the form of the elements.

\subsubsection{Case $1_{A}\otimes g\otimes gx_{1}x_{2}$}

From the first summand of the left side we get

\begin{eqnarray*}
l_{1}+u_{1} &=&1 \\
l_{2} &=&u_{2}=0 \\
a &=&0,b_{1}=l_{1},b_{2}=0, \\
d &=&1,e_{1}=u_{1},e_{2}=0
\end{eqnarray*}%
Since $\alpha \left( x_{2};0,0,1,0\right) =e_{2}=0$ and $\alpha \left(
x_{2};1,0,0,0\right) \equiv a+b_{1}\equiv 1,$ we get

\begin{equation*}
\left[ B(g\otimes x_{2};1_{A},gx_{1})-B(g\otimes x_{2};X_{1},g)\right]
1_{A}\otimes g\otimes gx_{1}x_{2}.
\end{equation*}%
From the second summand of the left side we get

\begin{eqnarray*}
l_{1}+u_{1} &=&1 \\
l_{2}+u_{2} &=&1 \\
a &=&0,b_{1}=l_{1},b_{2}=l_{2}, \\
d &=&1,e_{1}=u_{1},e_{2}=u_{2}.
\end{eqnarray*}%
Since
\begin{eqnarray*}
\alpha \left( 1_{H};0,0,1,1\right) &\equiv &1+e_{2}\equiv 0\text{, }\alpha
\left( 1_{H};0,1,1,0\right) \equiv e_{2}+a+b_{1}+b_{2}+1\equiv 0, \\
\alpha \left( 1_{H};1,0,0,1\right) &\equiv &a+b_{1}\equiv 1\text{and }\alpha
\left( 1_{H};1,1,0,0\right) \equiv 1+b_{2}\equiv 0
\end{eqnarray*}%
$,$ we get%
\begin{equation*}
\left[
\begin{array}{c}
B(x_{2}\otimes x_{2};1_{A},gx_{1}x_{2})+B(x_{2}\otimes x_{2};X_{2},gx_{1})+
\\
-B(x_{2}\otimes x_{2};X_{1},gx_{2})+B(x_{2}\otimes x_{2};X_{1}X_{2},g)%
\end{array}%
\right] 1_{A}\otimes g\otimes gx_{1}x_{2}.
\end{equation*}%
Since there is nothing from the right side, we obtain%
\begin{gather*}
B(g\otimes x_{2};1_{A},gx_{1})-B(g\otimes x_{2};X_{1},g)+ \\
B(x_{2}\otimes x_{2};1_{A},gx_{1}x_{2})+B(x_{2}\otimes x_{2};X_{2},gx_{1})+
\\
-+B(x_{2}\otimes x_{2};X_{1}X_{2},g)=0
\end{gather*}%
which holds in view of the form of the elements.

\subsection{$B(x_{2}\otimes x_{2};G,gx_{2})$}

We deduce that
\begin{equation*}
a=1,d=e_{2}=1
\end{equation*}%
and we get%
\begin{eqnarray*}
&&\left( -1\right) ^{\alpha \left( 1_{H};0,0,0,0\right) }B(x_{2}\otimes
x_{2};G,gx_{2})G\otimes gx_{2}\otimes g \\
&&+\left( -1\right) ^{\alpha \left( 1_{H};0,0,0,1\right) }B(x_{2}\otimes
x_{2};G,gx_{2})G\otimes g\otimes x_{2}.
\end{eqnarray*}

\subsubsection{Case $G\otimes g\otimes x_{2}$}

From the first summand of the left side we get

\begin{eqnarray*}
l_{1} &=&u_{1}=0 \\
l_{2} &=&u_{2}=0 \\
a &=&1,b_{1}=b_{2}=0 \\
d &=&1,e_{1}=e_{2}=0
\end{eqnarray*}%
Since $\alpha \left( x_{2};0,0,0,0\right) =a+b_{1}+b_{2}\equiv 1$ we get%
\begin{equation*}
-B(g\otimes x_{2};G,g)G\otimes g\otimes x_{2}.
\end{equation*}%
From the second summand of the left side we get

~%
\begin{eqnarray*}
l_{1} &=&u_{1}=0 \\
l_{2}+u_{2} &=&1 \\
a &=&1,b_{1}=0,b_{2}=l_{2} \\
d &=&1,e_{1}=0,e_{2}=u_{2}
\end{eqnarray*}%
Since $\alpha \left( 1_{H};0,0,0,1\right) \equiv a+b_{1}+b_{2}\equiv 1$\ and
$\alpha \left( 1_{H};0,1,0,0\right) \equiv 0$ and we obtain

\begin{equation*}
\left[ -B(x_{2}\otimes x_{2};G,gx_{2})+B(x_{2}\otimes x_{2};GX_{2},g)\right]
G\otimes g\otimes x_{2}.
\end{equation*}%
Considering also the right side, we get%
\begin{gather*}
-B(g\otimes x_{2};G,g)-B(x_{2}\otimes x_{2};G,gx_{2})+ \\
+B(x_{2}\otimes x_{2};GX_{2},g)-B(x_{2}\otimes 1_{H};G,g)=0
\end{gather*}

which holds in view of the form of the elements.

\subsection{$B(x_{2}\otimes x_{2};X_{1},gx_{2})$}

We deduce that
\begin{equation*}
b_{1}=1,d=1,e_{2}=1
\end{equation*}%
and we get%
\begin{eqnarray*}
&&\left( -1\right) ^{\alpha \left( 1_{H};0,0,0,0\right) }B(x_{2}\otimes
x_{2};X_{1},gx_{2})X_{1}\otimes gx_{2}\otimes g+ \\
&&\left( -1\right) ^{\alpha \left( 1_{H};1,0,0,0\right) }B(x_{2}\otimes
x_{2};X_{1},gx_{2})1_{A}\otimes gx_{2}\otimes x_{1} \\
&&\left( -1\right) ^{\alpha \left( 1_{H};0,0,0,1\right) }B(x_{2}\otimes
x_{2};X_{1},gx_{2})X_{1}\otimes g\otimes x_{2} \\
&&\left( -1\right) ^{\alpha \left( 1_{H};1,0,0,1\right) }B(x_{2}\otimes
x_{2};X_{1},gx_{2})1_{A}\otimes g\otimes gx_{1}x_{2}
\end{eqnarray*}

\subsubsection{Case $1_{A}\otimes gx_{2}\otimes x_{1}$}

This was already considered in subsection $B\left( x_{2}\otimes
x_{2};1_{A},gx_{1}x_{2}\right) .$

\subsubsection{Case $X_{1}\otimes g\otimes x_{2}$}

From the first summand of the left side we get

\begin{eqnarray*}
l_{1} &=&u_{1}=0 \\
l_{2} &=&u_{2}=0 \\
a &=&0,b_{1}=1,b_{2}=0 \\
d &=&1,e_{1}=e_{2}=0
\end{eqnarray*}%
Since $\alpha \left( x_{2};0,0,0,0\right) =a+b_{1}+b_{2}\equiv 1$ we get%
\begin{equation*}
-B(g\otimes x_{2};X_{1},g)X_{1}\otimes g\otimes x_{2}.
\end{equation*}%
From the second summand of the left side we get

~%
\begin{eqnarray*}
l_{1} &=&u_{1}=0 \\
l_{2}+u_{2} &=&1 \\
a &=&0,b_{1}=1,b_{2}=l_{2} \\
d &=&1,e_{1}=0,e_{2}=u_{2}
\end{eqnarray*}%
Since $\alpha \left( 1_{H};0,0,0,1\right) \equiv a+b_{1}+b_{2}\equiv 1$\ and
$\alpha \left( 1_{H};0,1,0,0\right) \equiv 0$ and we obtain

\begin{equation*}
\left[ -B(x_{2}\otimes x_{2};X_{1},gx_{2})+B(x_{2}\otimes x_{2};X_{1}X_{2},g)%
\right] X_{1}\otimes g\otimes x_{2}.
\end{equation*}%
Considering also the right side, we get%
\begin{gather*}
-B(g\otimes x_{2};X_{1},g)-B(x_{2}\otimes x_{2};X_{1},gx_{2})+ \\
+B(x_{2}\otimes x_{2};X_{1}X_{2},g)-B(x_{2}\otimes 1_{H};X_{1},g)=0
\end{gather*}%
which holds in view of the form of the elements.

\subsubsection{Case $1_{A}\otimes g\otimes gx_{1}x_{2}$}

This case was already considered in subsection $B\left( x_{2}\otimes
x_{2};1_{A},gx_{1}x_{2}\right) .$

\subsection{$B(x_{2}\otimes x_{2};X_{2},gx_{2})$}

We deduce that
\begin{equation*}
b_{2}=1,d=e_{2}=1
\end{equation*}%
and we get%
\begin{gather*}
\left( -1\right) ^{\alpha \left( 1_{H};0,0,0,0\right) }B(x_{2}\otimes
x_{2};X_{2},gx_{2})X_{2}\otimes gx_{2}\otimes g+ \\
+\left( -1\right) ^{\alpha \left( 1_{H};0,1,0,0\right) }B(x_{2}\otimes
x_{2};X_{2},gx_{2})1_{A}\otimes gx_{2}\otimes x_{2}+ \\
+\left( -1\right) ^{\alpha \left( 1_{H};0,0,0,1\right) }B(x_{2}\otimes
x_{2};X_{2},gx_{2})X_{2}\otimes g\otimes x_{2}+ \\
++\left( -1\right) ^{\alpha \left( 1_{H};0,1,0,1\right) }B(x_{2}\otimes
x_{2};X_{2},gx_{2})X_{2}^{1-l_{2}}\otimes g\otimes g^{l_{2}}x_{2}^{1+1}=0
\end{gather*}

\subsubsection{Case $1_{A}\otimes gx_{2}\otimes x_{2}$}

From the first summand of the left side we get

\begin{eqnarray*}
l_{1} &=&u_{1}=0 \\
l_{2} &=&u_{2}=0 \\
a &=&0,b_{1}=b_{2}=0 \\
d &=&1,e_{1}=0,e_{2}=1
\end{eqnarray*}%
Since $\alpha \left( x_{2};0,0,0,0\right) =a+b_{1}+b_{2}\equiv 0$ we get%
\begin{equation*}
B(g\otimes x_{2};1_{A},gx_{2})1_{A}\otimes gx_{2}\otimes x_{2}.
\end{equation*}%
From the second summand of the left side we get

~%
\begin{eqnarray*}
l_{1} &=&u_{1}=0 \\
l_{2}+u_{2} &=&1 \\
a &=&0,b_{1}=0,b_{2}=l_{2} \\
d &=&1,e_{1}=0,e_{2}-u_{2}=1\Rightarrow e_{2}=1,u_{2}=0,b_{2}=l_{2}=1
\end{eqnarray*}%
Since $\alpha \left( 1_{H};0,1,0,0\right) \equiv 0$ and we obtain

\begin{equation*}
B(x_{2}\otimes x_{2};X_{2},gx_{2})1_{A}\otimes gx_{2}\otimes x_{2}.
\end{equation*}%
Considering also the right side, we get%
\begin{equation*}
B(g\otimes x_{2};1_{A},gx_{2})+B(x_{2}\otimes
x_{2};X_{2},gx_{2})-B(x_{2}\otimes 1_{H};1_{A},gx_{2})=0
\end{equation*}

which holds in view of the form of the elements.

\subsubsection{Case $X_{2}\otimes g\otimes x_{2}$}

From the first summand of the left side we get

\begin{eqnarray*}
l_{1} &=&u_{1}=0 \\
l_{2} &=&u_{2}=0 \\
a &=&b_{1}=0,b_{2}=1 \\
d &=&1,e_{1}=e_{2}=
\end{eqnarray*}%
Since $\alpha \left( x_{2};0,0,0,0\right) =a+b_{1}+b_{2}\equiv 1$ we get%
\begin{equation*}
-B(g\otimes x_{2};X_{2},g)X_{2}\otimes g\otimes x_{2}.
\end{equation*}%
From the second summand of the left side we get

~%
\begin{eqnarray*}
l_{1} &=&u_{1}=0 \\
l_{2}+u_{2} &=&1 \\
a &=&b_{1}=0,b_{2}-l_{2}=1\Rightarrow b_{2}=1,l_{2}=0,u_{2}=1 \\
d &=&1,e_{1}=0,e_{2}=u_{2}=1
\end{eqnarray*}%
Since $\alpha \left( 1_{H};0,0,0,1\right) \equiv a+b_{1}+b_{2}\equiv 1$ and
we obtain

\begin{equation*}
-B(x_{2}\otimes x_{2};X_{2},gx_{2})X_{2}\otimes g\otimes x_{2}.
\end{equation*}%
Considering also the right side, we get%
\begin{equation*}
-B(g\otimes x_{2};X_{2},g)-B(x_{2}\otimes x_{2};X_{2},gx_{2})-B(x_{2}\otimes
1_{H};X_{2},g)=0
\end{equation*}%
which holds in view of the form of the elements.

\subsection{$B(x_{2}\otimes x_{2};GX_{2},gx_{1}x_{2})$}

We deduce that
\begin{equation*}
a=b_{2}=1,d=e_{1}=e_{2}=1
\end{equation*}%
and we get

\begin{equation*}
\left( -1\right) ^{\alpha \left( 1_{H};0,l_{2},u_{1},u_{2}\right)
}B(x_{2}\otimes x_{1};GX_{1}X_{2},gx_{2})GX_{2}^{1-l_{2}}\otimes
gx_{1}^{1-u_{1}}x_{2}^{1-u_{2}}\otimes
g^{1+l_{2}+u_{1}+u_{2}}x_{1}^{u_{1}}x_{2}^{l_{2}+u_{2}}
\end{equation*}%
\begin{gather*}
\left( -1\right) ^{\alpha \left( 1_{H};0,0,0,0\right) }B(x_{2}\otimes
x_{1};GX_{1}X_{2},gx_{2})GX_{2}\otimes gx_{1}x_{2}\otimes g \\
\left( -1\right) ^{\alpha \left( 1_{H};0,1,0,0\right) }B(x_{2}\otimes
x_{1};GX_{1}X_{2},gx_{2})G\otimes gx_{1}x_{2}\otimes x_{2} \\
\left( -1\right) ^{\alpha \left( 1_{H};0,0,1,0\right) }B(x_{2}\otimes
x_{1};GX_{1}X_{2},gx_{2})GX_{2}\otimes gx_{2}\otimes x_{1} \\
\left( -1\right) ^{\alpha \left( 1_{H};0,1,1,0\right) }B(x_{2}\otimes
x_{1};GX_{1}X_{2},gx_{2})G\otimes gx_{2}\otimes gx_{1}x_{2} \\
\left( -1\right) ^{\alpha \left( 1_{H};0,0,0,1\right) }B(x_{2}\otimes
x_{1};GX_{1}X_{2},gx_{2})GX_{2}\otimes gx_{1}\otimes x_{2} \\
\left( -1\right) ^{\alpha \left( 1_{H};0,1,0,1\right) }B(x_{2}\otimes
x_{1};GX_{1}X_{2},gx_{2})GX_{2}^{1-l_{2}}\otimes gx_{1}\otimes
g^{+l_{2}}x_{2}^{1+1}=0 \\
\left( -1\right) ^{\alpha \left( 1_{H};0,0,1,1\right) }B(x_{2}\otimes
x_{1};GX_{1}X_{2},gx_{2})GX_{2}\otimes g\otimes gx_{1}x_{2} \\
\left( -1\right) ^{\alpha \left( 1_{H};0,1,1,1\right) }B(x_{2}\otimes
x_{1};GX_{1}X_{2},gx_{2})GX_{2}^{1-l_{2}}\otimes g\otimes
g^{+l_{2}+1}x_{1}x_{2}^{1+1}=0
\end{gather*}

\subsubsection{Case $G\otimes gx_{1}x_{2}\otimes x_{2}$}

From the first summand of the left side we get

\begin{eqnarray*}
l_{1} &=&u_{1}=0 \\
l_{2} &=&u_{2}=0 \\
a &=&1,b_{1}=b_{2}=0 \\
d &=&e_{1}=e_{2}=1
\end{eqnarray*}%
Since $\alpha \left( x_{2};0,0,0,0\right) =a+b_{1}+b_{2}\equiv 1$ we get%
\begin{equation*}
-B(g\otimes x_{2};G,gx_{1}x_{2})G\otimes gx_{1}x_{2}\otimes x_{2}.
\end{equation*}%
From the second summand of the left side we get

~%
\begin{eqnarray*}
l_{1} &=&u_{1}=0 \\
l_{2}+u_{2} &=&1 \\
a &=&1,b_{1}=0,b_{2}=l_{2} \\
d &=&1,e_{1}=1,e_{2}-u_{2}=1\Rightarrow e_{2}=1,u_{2}=0,b_{2}=l_{2}=1
\end{eqnarray*}%
Since $\alpha \left( 1_{H};0,1,0,0\right) \equiv 0$ and we obtain

\begin{equation*}
B(x_{2}\otimes x_{2};GX_{2},gx_{1}x_{2})G\otimes gx_{1}x_{2}\otimes x_{2}.
\end{equation*}%
Considering also the right side, we get%
\begin{equation*}
-B(g\otimes x_{2};G,gx_{1}x_{2})+B(x_{2}\otimes
x_{2};GX_{2},gx_{1}x_{2})-B(x_{2}\otimes 1_{H};G,gx_{1}x_{2})=0
\end{equation*}%
which holds in view of the form of the elements.

\subsubsection{Case $GX_{2}\otimes gx_{2}\otimes x_{1}$}

Nothing from the first summand of the left side. The second summand of the
left side gives us%
\begin{eqnarray*}
l_{1}+u_{1} &=&1,l_{2}=u_{2}=0 \\
a &=&1,b_{1}=l_{1},b_{2}=1 \\
d &=&1,e_{2}=1,e_{1}=u_{1}
\end{eqnarray*}%
Since $\alpha \left( 1_{H};0,0,1,0\right) \equiv e_{2}+\left(
a+b_{1}+b_{2}\right) \equiv 1$ and $\alpha \left( 1_{H};1,0,0,0\right)
\equiv b_{2}=0,$ we get%
\begin{equation*}
\left[ -B(x_{2}\otimes x_{2};GX_{2},gx_{1}x_{2})+B(x_{2}\otimes
x_{2};GX_{1}X_{2},gx_{2})\right] GX_{2}\otimes gx_{2}\otimes x_{1}.
\end{equation*}%
Since there is nothing in the right side, we obtain%
\begin{equation*}
-B(x_{2}\otimes x_{2};GX_{2},gx_{1}x_{2})+B(x_{2}\otimes
x_{2};GX_{1}X_{2},gx_{2})=0
\end{equation*}%
which holds in view of the form of the elements.

\subsubsection{Case $G\otimes gx_{2}\otimes gx_{1}x_{2}$}

From the first summand of the left side we get

\begin{eqnarray*}
l_{1}+u_{1} &=&1 \\
l_{2} &=&u_{2}=0 \\
a &=&1,b_{1}=l_{1},b_{2}=0, \\
d &=&1,e_{1}=u_{1},e_{2}=1
\end{eqnarray*}%
Since $\alpha \left( x_{2};0,0,1,0\right) =e_{2}=1$ and $\alpha \left(
x_{2};1,0,0,0\right) \equiv a+b_{1}\equiv 0,$ we get

\begin{equation*}
\left[ -B(g\otimes x_{2};G,gx_{1}x_{2})+B(g\otimes x_{2};GX_{1},gx_{2})%
\right] G\otimes gx_{2}\otimes gx_{1}x_{2}.
\end{equation*}%
From the second summand of the left side we get

\begin{eqnarray*}
l_{1}+u_{1} &=&1 \\
l_{2}+u_{2} &=&1 \\
a &=&1,b_{1}=l_{1},b_{2}=l_{2}, \\
d &=&1,e_{1}=u_{1},e_{2}-u_{2}=1\Rightarrow e_{2}=1,u_{2}=0,b_{2}=l_{2}=1.
\end{eqnarray*}%
Since $\alpha \left( 1_{H};0,1,1,0\right) \equiv e_{2}+a+b_{1}+b_{2}+1\equiv
0$ and $\alpha \left( 1_{H};1,1,0,0\right) \equiv 1+b_{2}\equiv 0$, we get%
\begin{equation*}
\left[ B(x_{2}\otimes x_{2};GX_{2},gx_{1}x_{2})+B(x_{2}\otimes
x_{2};GX_{1}X_{2},gx_{2})+\right] G\otimes gx_{2}\otimes gx_{1}x_{2}.
\end{equation*}%
Since there is nothing from the right side, we obtain%
\begin{gather*}
-B(g\otimes x_{2};G,gx_{1}x_{2})+B(g\otimes x_{2};GX_{1},gx_{2})+ \\
B(x_{2}\otimes x_{2};GX_{2},gx_{1}x_{2})+B(x_{2}\otimes
x_{2};GX_{1}X_{2},gx_{2})=0
\end{gather*}%
which holds in view of the form of the elements.

\subsubsection{Case $GX_{2}\otimes gx_{1}\otimes x_{2}$}

From the first summand of the left side we get

\begin{eqnarray*}
l_{1} &=&u_{1}=0 \\
l_{2} &=&u_{2}=0 \\
a &=&1,b_{1}=0,b_{2}=1 \\
d &=&e_{1}=1,e_{2}=0
\end{eqnarray*}%
Since $\alpha \left( x_{2};0,0,0,0\right) =a+b_{1}+b_{2}\equiv 0$ we get%
\begin{equation*}
+B(g\otimes x_{2};GX_{2},gx_{1})GX_{2}\otimes gx_{1}\otimes x_{2}.
\end{equation*}%
From the second summand of the left side we get

~%
\begin{eqnarray*}
l_{1} &=&u_{1}=0 \\
l_{2}+u_{2} &=&1 \\
a &=&1,b_{1}=0,b_{2}-l_{2}=1\Rightarrow b_{2}=1,l_{2}=0,u_{2}=1 \\
d &=&1,e_{1}=1,e_{2}=u_{2}=1
\end{eqnarray*}%
Since $\alpha \left( 1_{H};0,0,0,1\right) \equiv a+b_{1}+b_{2}\equiv 0$ and
we obtain

\begin{equation*}
B(x_{2}\otimes x_{2};GX_{2},gx_{1}x_{2})GX_{2}\otimes gx_{1}\otimes x_{2}.
\end{equation*}%
Considering also the right side, we get%
\begin{equation*}
B(g\otimes x_{2};GX_{2},gx_{1})+B(x_{2}\otimes
x_{2};GX_{2},gx_{1}x_{2})-B(x_{2}\otimes 1_{H};GX_{2},gx_{1})=0
\end{equation*}%
which holds in view of the form of the elements.

\subsubsection{Case $GX_{2}\otimes g\otimes gx_{1}x_{2}$}

From the first summand of the left side we get

\begin{eqnarray*}
l_{1}+u_{1} &=&1 \\
l_{2} &=&u_{2}=0 \\
a &=&1,b_{1}=l_{1},b_{2}=1, \\
d &=&1,e_{1}=u_{1},e_{2}=0
\end{eqnarray*}%
Since $\alpha \left( x_{2};0,0,1,0\right) =e_{2}=0$ and $\alpha \left(
x_{2};1,0,0,0\right) \equiv a+b_{1}\equiv 0,$ we get

\begin{equation*}
\left[ +B(g\otimes x_{2};GX_{2},gx_{1})+B(g\otimes x_{2};GX_{1}X_{2},g)%
\right] GX_{2}\otimes g\otimes gx_{1}x_{2}.
\end{equation*}%
From the second summand of the left side we get

\begin{eqnarray*}
l_{1}+u_{1} &=&1 \\
l_{2}+u_{2} &=&1 \\
a &=&1,b_{1}=l_{1},b_{2}-l_{2}=1\Rightarrow b_{2}=1,l_{2}=0,u_{2}=1 \\
d &=&1,e_{1}=u_{1},e_{2}=u_{2}=1.
\end{eqnarray*}%
Since $\alpha \left( 1_{H};0,0,1,1\right) \equiv 1+e_{2}\equiv 0$ and $%
\alpha \left( 1_{H};1,0,0,1\right) \equiv a+b_{1}\equiv 0$, we get%
\begin{equation*}
\left[ B(x_{2}\otimes x_{2};GX_{2},gx_{1}x_{2})+B(x_{2}\otimes
x_{2};GX_{1}X_{2},gx_{2})+\right] GX_{2}\otimes g\otimes gx_{1}x_{2}.
\end{equation*}%
Since there is nothing from the right side, we get%
\begin{gather*}
+B(g\otimes x_{2};GX_{2},gx_{1})+B(g\otimes x_{2};GX_{1}X_{2},g)+ \\
B(x_{2}\otimes x_{2};GX_{2},gx_{1}x_{2})+B(x_{2}\otimes
x_{2};GX_{1}X_{2},gx_{2})=0
\end{gather*}%
which holds in view of the form of the elements.

\subsection{$B(x_{2}\otimes x_{2};GX_{1}X_{2},gx_{2})$}

We deduce that%
\begin{equation*}
a=b_{1}=b_{2}=1,d=e_{2}=1
\end{equation*}%
and we get%
\begin{equation*}
\left( -1\right) ^{\alpha \left( 1_{H};l_{1},l_{2},0,u_{2}\right)
}B(x_{2}\otimes
x_{2};GX_{1}X_{2},gx_{2})GX_{1}^{1-l_{1}}X_{2}^{1-l_{2}}\otimes
gx_{2}^{1-u_{2}}\otimes
g^{1+l_{1}+l_{2}+u_{2}}x_{1}^{l_{1}}x_{2}^{l_{2}+u_{2}}
\end{equation*}%
\begin{gather*}
\left( -1\right) ^{\alpha \left( 1_{H};0,0,0,0\right) }B(x_{2}\otimes
x_{2};GX_{1}X_{2},gx_{2})GX_{1}X_{2}\otimes gx_{2}\otimes g \\
\left( -1\right) ^{\alpha \left( 1_{H};1,0,0,0\right) }B(x_{2}\otimes
x_{2};GX_{1}X_{2},gx_{2})GX_{2}\otimes gx_{2}\otimes x_{1} \\
\left( -1\right) ^{\alpha \left( 1_{H};0,1,0,0\right) }B(x_{2}\otimes
x_{2};GX_{1}X_{2},gx_{2})GX_{1}\otimes gx_{2}\otimes x_{2} \\
\left( -1\right) ^{\alpha \left( 1_{H};1,1,0,0\right) }B(x_{2}\otimes
x_{2};GX_{1}X_{2},gx_{2})G\otimes gx_{2}\otimes gx_{1}x_{2} \\
\left( -1\right) ^{\alpha \left( 1_{H};0,0,0,1\right) }B(x_{2}\otimes
x_{2};GX_{1}X_{2},gx_{2})GX_{1}X_{2}\otimes g\otimes x_{2} \\
\left( -1\right) ^{\alpha \left( 1_{H};1,0,0,1\right) }B(x_{2}\otimes
x_{2};GX_{1}X_{2},gx_{2})GX_{2}\otimes g\otimes gx_{1}x_{2} \\
\left( -1\right) ^{\alpha \left( 1_{H};l_{1},1,0,1\right) }B(x_{2}\otimes
x_{2};GX_{1}X_{2},gx_{2})GX_{1}^{1-l_{1}}X_{2}^{1-l_{2}}\otimes g\otimes
g^{l_{1}+l_{2}}x_{1}^{l_{1}}x_{2}^{1+1}=0
\end{gather*}

\subsubsection{Case $GX_{2}\otimes gx_{2}\otimes x_{1}$}

Already considered in subsection $B(x_{2}\otimes x_{2};GX_{2},gx_{1}x_{2}).$

\subsubsection{Case $GX_{1}\otimes gx_{2}\otimes x_{2}$}

From the first summand of the left side we get

\begin{eqnarray*}
l_{1} &=&u_{1}=0 \\
l_{2} &=&u_{2}=0 \\
a &=&b_{1}=1,b_{2}=0 \\
d &=&e_{2}=1,e_{1}=0
\end{eqnarray*}%
Since $\alpha \left( x_{2};0,0,0,0\right) =a+b_{1}+b_{2}\equiv 0$ we get%
\begin{equation*}
+B(g\otimes x_{2};GX_{1},gx_{2})GX_{1}\otimes gx_{2}\otimes x_{2}.
\end{equation*}%
From the second summand of the left side we get

~%
\begin{eqnarray*}
l_{1} &=&u_{1}=0 \\
l_{2}+u_{2} &=&1 \\
a &=&b_{1}=1,b_{2}=l_{2} \\
d &=&1,e_{1}=0,e_{2}-u_{2}=1\Rightarrow e_{2}=1,u_{2}=0,b_{2}=l_{2}=1
\end{eqnarray*}%
Since $\alpha \left( 1_{H};0,1,0,0\right) \equiv 0$ and we obtain

\begin{equation*}
B(x_{2}\otimes x_{2};GX_{1}X_{2},gx_{2})GX_{1}\otimes gx_{2}\otimes x_{2}.
\end{equation*}%
Considering also the right side, we get%
\begin{equation*}
B(g\otimes x_{2};GX_{1},gx_{2})+B(x_{2}\otimes
x_{2};GX_{1}X_{2},gx_{2})-B(x_{2}\otimes 1_{H};GX_{1},gx_{2})=0
\end{equation*}%
which holds in view of the form of the elements.

\subsubsection{Case $G\otimes gx_{2}\otimes gx_{1}x_{2}$}

This case was already considered in subsection $B(x_{2}\otimes
x_{2};GX_{2},gx_{1}x_{2}).$

\subsubsection{Case $GX_{1}X_{2}\otimes g\otimes x_{2}$}

From the first summand of the left side we get

\begin{eqnarray*}
l_{1} &=&u_{1}=0 \\
l_{2} &=&u_{2}=0 \\
a &=&b_{1}=b_{2}=1 \\
d &=&1,e_{1}=e_{2}=0.
\end{eqnarray*}%
Since $\alpha \left( x_{2};0,0,0,0\right) =a+b_{1}+b_{2}\equiv 1$ we get%
\begin{equation*}
-B(g\otimes x_{2};GX_{1}X_{2},g)GX_{1}X_{2}\otimes g\otimes x_{2}.
\end{equation*}%
From the second summand of the left side we get

~%
\begin{eqnarray*}
l_{1} &=&u_{1}=0 \\
l_{2}+u_{2} &=&1 \\
a &=&b_{1}=1,b_{2}-l_{2}=1\Rightarrow b_{2}=1,l_{2}=0,u_{2}=1 \\
d &=&1,e_{1}=0,e_{2}=u_{2}=1
\end{eqnarray*}%
Since $\alpha \left( 1_{H};0,0,0,1\right) \equiv a+b_{1}+b_{2}\equiv 1$ and
we obtain

\begin{equation*}
-B(x_{2}\otimes x_{2};GX_{1}X_{2},gx_{2})GX_{1}X_{2}\otimes g\otimes x_{2}.
\end{equation*}%
Considering also the right side, we get%
\begin{equation*}
-B(g\otimes x_{2};GX_{1}X_{2},g)-B(x_{2}\otimes
x_{2};GX_{1}X_{2},gx_{2})-B(x_{2}\otimes 1_{H};GX_{1}X_{2},g)=0
\end{equation*}%
which holds in view of the form of the elements.

\subsubsection{Case $GX_{2}\otimes g\otimes gx_{1}x_{2}$}

This case was already considered in subsection $B(x_{2}\otimes
x_{2};GX_{2},gx_{1}x_{2}).$

\section{$B\left( x_{2}\otimes x_{1}x_{2}\right) $}

From $\left( \ref{simplxx}\right) $ we get%
\begin{eqnarray}
B(x_{2}\otimes x_{1}x_{2}) &=&B(x_{2}\otimes 1_{H})(1_{A}\otimes x_{1}x_{2})
\label{form x2otx1x2} \\
&&+(1_{A}\otimes gx_{2})B(x_{2}\otimes 1_{H})(1_{A}\otimes gx_{1})  \notag \\
&&-(1_{A}\otimes gx_{1})B(x_{2}\otimes 1_{H})(1_{A}\otimes gx_{2})  \notag \\
&&+(1_{A}\otimes x_{1}x_{2})B(x_{2}\otimes 1_{H})  \notag \\
&&-(1_{A}\otimes g)B(gx_{1}x_{2}\otimes 1_{H})(1_{A}\otimes gx_{2})  \notag
\\
&&-(1_{A}\otimes x_{2})B(gx_{1}x_{2}\otimes 1_{H}).  \notag
\end{eqnarray}

Thus we obtain%
\begin{eqnarray*}
&&B\left( x_{2}\otimes x_{1}x_{2}\right) =-2B(gx_{1}x_{2}\otimes
1_{H};1_{A},g)1_{A}\otimes gx_{2}+ \\
&&\left[ 4B\left( x_{2}\otimes 1_{H};G,g\right) -2B(gx_{1}x_{2}\otimes
1_{H};G,gx_{1})\right] G\otimes gx_{1}x_{2}+ \\
&&-2B(x_{2}\otimes 1_{H};1_{A},gx_{1})X_{1}\otimes gx_{1}x_{2}+ \\
&&+\left[
\begin{array}{c}
-4B(g\otimes 1_{H};1_{A},g)-4B(x_{2}\otimes \ 1_{H};1_{A},gx_{2}) \\
-2B(x_{1}\otimes 1_{H};1_{A},gx_{1})-2B(gx_{1}x_{2}\otimes
1_{H};1_{A},gx_{1}x_{2})%
\end{array}%
\right] X_{2}\otimes gx_{1}x_{2}+ \\
&&+\left[
\begin{array}{c}
-2B(g\otimes 1_{H};1_{A},g)-2B(x_{2}\otimes \ 1_{H};1_{A},gx_{2}) \\
-2B(x_{1}\otimes 1_{H};1_{A},gx_{1})-2B(gx_{1}x_{2}\otimes
1_{H};1_{A},gx_{1}x_{2})%
\end{array}%
\right] X_{1}X_{2}\otimes gx_{2}+ \\
&&+\left[ -2B(x_{2}\otimes 1_{H};G,g)+2B(gx_{1}x_{2}\otimes 1_{H};G,gx_{1})%
\right] GX_{1}\otimes gx_{2}+ \\
&&+\left[ +2B(x_{1}\otimes 1_{H};G,g)+2B(gx_{1}x_{2}\otimes 1_{H};G,gx_{2})%
\right] GX_{2}\otimes gx_{2}+ \\
&&\left[ 2B(g\otimes 1_{H};G,gx_{1})+2B(x_{2}\otimes \ 1_{H};G,gx_{1}x_{2})%
\right] GX_{1}X_{2}\otimes gx_{1}x_{2}+
\end{eqnarray*}%
We write the Casimir formula for $B\left( x_{2}\otimes x_{1}x_{2}\right) .$%
\begin{eqnarray*}
&&\sum_{a,b_{1},b_{2},d,e_{1},e_{2}=0}^{1}\sum_{l_{1}=0}^{b_{1}}%
\sum_{l_{2}=0}^{b_{2}}\sum_{u_{1}=0}^{e_{1}}\sum_{u_{2}=0}^{e_{2}}\left(
-1\right) ^{\alpha \left( x_{2};l_{1},l_{2},u_{1},u_{2}\right) } \\
&&B(g\otimes
x_{1}x_{2};G^{a}X_{1}^{b_{1}}X_{2}^{b_{2}},g^{d}x_{1}^{e_{1}}x_{2}^{e_{2}})
\\
&&G^{a}X_{1}^{b_{1}-l_{1}}X_{2}^{b_{2}-l_{2}}\otimes
g^{d}x_{1}^{e_{1}-u_{1}}x_{2}^{e_{2}-u_{2}}\otimes
g^{a+b_{1}+b_{2}+l_{1}+l_{2}+d+e_{1}+e_{2}+u_{1}+u_{2}}x_{1}^{l_{1}+u_{1}}x_{2}^{l_{2}+u_{2}+1}
\\
&&\sum_{a,b_{1},b_{2},d,e_{1},e_{2}=0}^{1}\sum_{l_{1}=0}^{b_{1}}%
\sum_{l_{2}=0}^{b_{2}}\sum_{u_{1}=0}^{e_{1}}\sum_{u_{2}=0}^{e_{2}}\left(
-1\right) ^{\alpha \left( 1_{H};l_{1},l_{2},u_{1},u_{2}\right) } \\
&&B(x_{2}\otimes
x_{1}x_{2};G^{a}X_{1}^{b_{1}}X_{2}^{b_{2}},g^{d}x_{1}^{e_{1}}x_{2}^{e_{2}})
\\
&&G^{a}X_{1}^{b_{1}-l_{1}}X_{2}^{b_{2}-l_{2}}\otimes
g^{d}x_{1}^{e_{1}-u_{1}}x_{2}^{e_{2}-u_{2}}\otimes
g^{a+b_{1}+b_{2}+l_{1}+l_{2}+d+e_{1}+e_{2}+u_{1}+u_{2}}x_{1}^{l_{1}+u_{1}}x_{2}^{l_{2}+u_{2}}
\\
&=&B^{A}(x_{2}\otimes x_{1}x_{2})\otimes B^{H}(x_{2}\otimes
x_{1}x_{2})\otimes 1_{H} \\
&&B^{A}(x_{2}\otimes x_{1})\otimes B^{H}(x_{2}\otimes x_{1})\otimes gx_{2} \\
&&-B^{A}(x_{2}\otimes x_{2})\otimes B^{H}(x_{2}\otimes x_{2})\otimes gx_{1}
\\
&&B^{A}(x_{2}\otimes 1_{H})\otimes B^{H}(x_{2}\otimes 1_{H})\otimes
x_{1}x_{2}
\end{eqnarray*}

\subsection{$B(x_{2}\otimes x_{1}x_{2};1_{A},gx_{2})$}

We deduce that

\begin{eqnarray*}
a &=&b_{1}=b_{2}=0 \\
d &=&e_{2}=1,e_{1}=0
\end{eqnarray*}%
\begin{eqnarray*}
&&\left( -1\right) ^{\alpha \left( 1_{H};0,0,0,0\right) }B(x_{2}\otimes
x_{1}x_{2};1_{A},gx_{2})1_{A}\otimes gx_{2}\otimes 1_{H} \\
&&+\left( -1\right) ^{\alpha \left( 1_{H};0,0,0,1\right) }B(x_{2}\otimes
x_{1}x_{2};1_{A},gx_{2})1_{A}\otimes g\otimes gx_{2}
\end{eqnarray*}

\subsubsection{Case $1_{A}\otimes g\otimes gx_{2}$}

First summand%
\begin{eqnarray*}
l_{1} &=&u_{1}=0,l_{2}=u_{2}=0 \\
a &=&b_{1}=b_{2}=0,d=1,e_{1}=e_{2}=0
\end{eqnarray*}%
and we get

\begin{equation*}
B(g\otimes x_{1}x_{2};1_{A},g)1_{A}\otimes g\otimes gx_{2}
\end{equation*}%
second summand
\begin{eqnarray*}
l_{1} &=&u_{1}=0,l_{2}+u_{2}=1 \\
a &=&b_{1}=0,b_{2}=l_{2} \\
d &=&1,e_{1}=0,e_{2}=u_{2}
\end{eqnarray*}%
and we get%
\begin{eqnarray*}
&&B(g\otimes x_{1}x_{2};1_{A},g)+B(x_{2}\otimes x_{1}x_{2};1_{A},gx_{2}) \\
&&+B(x_{2}\otimes x_{1}x_{2};X_{2},g)-B(x_{2}\otimes x_{1};1_{A},g)=0
\end{eqnarray*}%
which holds in view of the form of the elements

\subsection{$B\left( x_{2}\otimes x_{1}x_{2};G,gx_{1}x_{2}\right) $}

\begin{eqnarray*}
a &=&1,b_{1}=b_{2}=0 \\
d &=&e_{1}=e_{2}=1
\end{eqnarray*}%
\begin{eqnarray*}
&&\sum_{u_{1}=0}^{1}\sum_{u_{2}=0}^{1}\left( -1\right) ^{\alpha \left(
1_{H};0,0,u_{1},u_{2}\right) }B(x_{2}\otimes
x_{1}x_{2};G,gx_{1}x_{2})G\otimes gx_{1}^{1-u_{1}}x_{2}^{1-u_{2}} \\
&&\otimes g^{a+1+u_{1}+u_{2}}x_{1}^{u_{1}}x_{2}^{u_{2}} \\
&=&\left( -1\right) ^{\alpha \left( 1_{H};0,0,0,0\right) }B(x_{2}\otimes
x_{1}x_{2};G,gx_{1}x_{2})G\otimes gx_{1}x_{2}\otimes 1_{H} \\
&&\left( -1\right) ^{\alpha \left( 1_{H};0,0,0,1\right) }B(x_{2}\otimes
x_{1}x_{2};G,gx_{1}x_{2})G\otimes gx_{1}\otimes gx_{2} \\
&&\left( -1\right) ^{\alpha \left( 1_{H};0,0,1,0\right) }B(x_{2}\otimes
x_{1}x_{2};G,gx_{1}x_{2})G\otimes gx_{2}\otimes gx_{1} \\
&&\left( -1\right) ^{\alpha \left( 1_{H};0,0,1,1\right) }B(x_{2}\otimes
x_{1}x_{2};G,gx_{1}x_{2})G\otimes g\otimes x_{1}x_{2}
\end{eqnarray*}

\subsubsection{Case $G\otimes gx_{1}x_{2}\otimes 1_{H}$}

This gives us nothing new.

\subsubsection{Case $G\otimes gx_{1}\otimes gx_{2}$}

The first summand of the left side of the equality gives us%
\begin{eqnarray*}
l_{1} &=&u_{1}=0,l_{2}=u_{2}=0 \\
a &=&1,b_{1}=b_{2}=0 \\
d &=&e_{1}=1,e_{2}=0 \\
&&
\end{eqnarray*}%
Since $\alpha \left( x_{2};0,0,0,0\right) \equiv a+b_{1}+b_{2}\equiv 1$ we
get%
\begin{equation*}
-B(g\otimes x_{1}x_{2};G,gx_{1})G\otimes gx_{1}\otimes gx_{2}
\end{equation*}

Second summand%
\begin{eqnarray*}
l_{1} &=&u_{1}=0,l_{2}+u_{2}=1 \\
a &=&1,b_{1}=0,b_{2}=l_{2} \\
d &=&e_{1}=1,e_{2}=u_{2}
\end{eqnarray*}%
Since $\alpha \left( 1_{H};0,0,0,1\right) \equiv 1+0+0=1$ and $\alpha \left(
1_{H};0,1,0,0\right) \equiv 0,$ we get%
\begin{equation*}
\left[ -B(x_{2}\otimes x_{1}x_{2};G,gx_{1}x_{2})+B(x_{2}\otimes
x_{1}x_{2};GX_{2},gx_{1})\right] G\otimes gx_{1}\otimes gx_{2}
\end{equation*}%
By taking in account also the right side of the equality we get

\begin{gather*}
-B(g\otimes x_{1}x_{2};G,gx_{1})-B(x_{2}\otimes
x_{1};G,gx_{1})-B(x_{2}\otimes x_{1}x_{2};G,gx_{1}x_{2}) \\
+B(x_{2}\otimes x_{1}x_{2};GX_{2},gx_{1})=0
\end{gather*}%
which holds in view of the form of the elements.

\subsubsection{Case $G\otimes gx_{2}\otimes gx_{1}$}

The first summand of the left side gives us nothing. The second summand of
the left side gives us%
\begin{eqnarray*}
l_{1}+u_{1} &=&1,l_{2}=u_{2}=0 \\
a &=&1,b_{1}=l_{1},b_{2}=0 \\
d &=&1,e_{1}=u_{1},e_{2}=1.
\end{eqnarray*}%
Since $\alpha \left( 1_{H};0,0,1,1\right) \equiv 1+e_{2}\equiv 0$ and $%
\alpha \left( 1_{H};1,0,0,1\right) \equiv a+b_{1}\equiv 0,$ we obtain
\begin{equation*}
\left[ B(x_{2}\otimes x_{1}x_{2};G,gx_{1}x_{2})+B(x_{2}\otimes
x_{1}x_{2};GX_{1},gx_{2})\right] G\otimes gx_{2}\otimes gx_{1}
\end{equation*}%
By taking in account also the right side we get%
\begin{equation*}
B(x_{2}\otimes x_{2};G,gx_{2})+B(x_{2}\otimes
x_{1}x_{2};G,gx_{1}x_{2})+B(x_{2}\otimes x_{1}x_{2};GX_{1},gx_{2})=0
\end{equation*}%
which holds in view of the form of the elements.

\subsubsection{Case $G\otimes g\otimes x_{1}x_{2}$}

First summand of the left side of the equality gives us%
\begin{eqnarray*}
l_{1}+u_{1} &=&1,l_{2}=u_{2}=0 \\
a &=&1,b_{1}=l_{1},b_{2}=0 \\
d &=&1,e_{1}=u_{1},e_{2}=0.
\end{eqnarray*}%
Since $\alpha \left( x_{2};0,0,1,0\right) \equiv e_{2}=0$ and $\alpha \left(
x_{2};1,0,0,0\right) \equiv b_{2}+a+b_{1}+b_{2}+1+1\equiv 0,$ we obtain
\begin{equation*}
\left[ B(g\otimes x_{1}x_{2};G,gx_{1})+B(g\otimes x_{1}x_{2};GX_{1},g)\right]
G\otimes g\otimes x_{1}x_{2}.
\end{equation*}%
Second summand of the left side gives us%
\begin{eqnarray*}
l_{1}+u_{1} &=&1,l_{2}+u_{2}=1 \\
a &=&1,b_{1}=l_{1},b_{2}=l_{2} \\
d &=&1,e_{1}=u_{1},e_{2}=u_{2}
\end{eqnarray*}%
Since%
\begin{eqnarray*}
&&\alpha \left( 1_{H};0,0,1,1\right) \equiv 1+e_{2}\equiv 0 \\
&&\alpha \left( 1_{H};0,1,1,0\right) \equiv e_{2}+\left(
a+b_{1}+b_{2}+1\right) \equiv 1 \\
&&\alpha \left( 1_{H};1,0,0,1\right) \equiv a+b_{1}\equiv 0 \\
&&\alpha \left( 1_{H};1,1,0,0\right) \equiv 1+b_{2}\equiv 0
\end{eqnarray*}%
we obtain%
\begin{equation*}
\left[
\begin{array}{c}
B(x_{2}\otimes x_{1}x_{2};G,gx_{1}x_{2})-B(x_{2}\otimes
x_{1}x_{2};GX_{2},gx_{1})+ \\
+B(x_{2}\otimes x_{1}x_{2};GX_{1},gx_{2})+B(x_{2}\otimes
x_{1}x_{2};GX_{1}X_{2},g)%
\end{array}%
\right] G\otimes g\otimes x_{1}x_{2}.
\end{equation*}%
By taking in account also the right side, we get%
\begin{gather*}
B(g\otimes x_{1}x_{2};G,gx_{1})+B(g\otimes
x_{1}x_{2};GX_{1},g)-B(x_{2}\otimes 1_{H};G,g)+ \\
+B(x_{2}\otimes x_{1}x_{2};G,gx_{1}x_{2})+B(x_{2}\otimes
x_{1}x_{2};GX_{1},gx_{2}) \\
-B(x_{2}\otimes x_{1}x_{2};GX_{2},gx_{1})+B(x_{2}\otimes
x_{1}x_{2};GX_{1}X_{2},g)=0
\end{gather*}%
which holds in view of the form of the elements.

\subsection{$B(x_{2}\otimes x_{1}x_{2};X_{2},gx_{1}x_{2})$}

We deduce that%
\begin{eqnarray*}
a &=&0,b_{1}=0,b_{2}=1 \\
d &=&e_{1}=e_{2}=1
\end{eqnarray*}%
and we get%
\begin{eqnarray*}
&&\left( -1\right) ^{\alpha \left( 1_{H};0,0,0,0\right) }B(x_{2}\otimes
x_{1}x_{2};X_{2},gx_{1}x_{2})X_{2}\otimes gx_{1}x_{2}\otimes 1_{H} \\
&&\left( -1\right) ^{\alpha \left( 1_{H};0,0,0,1\right) }B(x_{2}\otimes
x_{1}x_{2};X_{2},gx_{1}x_{2})X_{2}\otimes gx_{1}\otimes gx_{2} \\
&&\left( -1\right) ^{\alpha \left( 1_{H};0,1,0,0\right) }B(x_{2}\otimes
x_{1}x_{2};X_{2},gx_{1}x_{2})1_{A}\otimes gx_{1}x_{2}\otimes gx_{2} \\
&&\left( -1\right) ^{\alpha \left( 1_{H};0,1,0,1\right) }B(x_{2}\otimes
x_{1}x_{2};X_{2},gx_{1}x_{2})X_{2}\otimes gx_{1}\otimes
g^{l_{2}+u_{1}+u_{2}}x_{1}^{u_{1}}x_{2}^{1+1}=0 \\
&&\left( -1\right) ^{\alpha \left( 1_{H};0,0,1,0\right) }B(x_{2}\otimes
x_{1}x_{2};X_{2},gx_{1}x_{2})X_{2}\otimes gx_{2}\otimes gx_{1} \\
&&\left( -1\right) ^{\alpha \left( 1_{H};0,0,1,1\right) }B(x_{2}\otimes
x_{1}x_{2};X_{2},gx_{1}x_{2})X_{2}\otimes g\otimes x_{1}x_{2} \\
&&\left( -1\right) ^{\alpha \left( 1_{H};0,1,1,0\right) }B(x_{2}\otimes
x_{1}x_{2};X_{2},gx_{1}x_{2})1_{A}\otimes gx_{2}\otimes x_{1}x_{2} \\
&&\left( -1\right) ^{\alpha \left( 1_{H};0,1,1,1\right) }B(x_{2}\otimes
x_{1}x_{2};X_{2},gx_{1}x_{2})X_{2}^{1-l_{2}}\otimes
gx_{1}^{1-u_{1}}x_{2}^{1-u_{2}}\otimes
g^{l_{2}+u_{1}+u_{2}}x_{1}^{u_{1}}x_{2}^{1+1}=0
\end{eqnarray*}

\subsubsection{Case $X_{2}\otimes gx_{1}x_{2}\otimes 1_{H}$}

Nothing new.

\subsubsection{Case $X_{2}\otimes gx_{1}\otimes gx_{2}$}

First summand of the left side gives us%
\begin{eqnarray*}
l_{1} &=&u_{1}=0,l_{2}=u_{2}=0 \\
a &=&b_{1}=0,b_{2}=1 \\
d &=&e_{1}=1,e_{2}=0.
\end{eqnarray*}%
Since $\alpha \left( x_{1};0,0,0,0\right) \equiv a+b_{1}+b_{2}\equiv 1,$ we
get%
\begin{equation*}
-B(g\otimes x_{1}x_{2};X_{2},gx_{1})X_{2}\otimes gx_{1}\otimes gx_{2}
\end{equation*}%
Second summand of the left side gives us%
\begin{eqnarray*}
l_{1} &=&u_{1}=0,l_{2}+u_{2}=1 \\
a &=&b_{1}=0,b_{2}-l_{2}=1\Rightarrow b_{2}=1,l_{2}=0\Rightarrow u_{2}=1 \\
d &=&e_{1}=1,e_{2}=u_{2}=1.
\end{eqnarray*}%
Since $\alpha \left( 1_{H};0,0,0,1\right) \equiv a+b_{1}+b_{2}\equiv 1,$ we
obtain%
\begin{equation*}
-B(x_{2}\otimes x_{1}x_{2};X_{2},gx_{1}x_{2})X_{2}\otimes gx_{1}\otimes
gx_{2}.
\end{equation*}%
By taking in account the right side we get

\begin{equation*}
-B(g\otimes x_{1}x_{2};X_{2},gx_{1})-B(x_{2}\otimes
x_{1},X_{2},gx_{1})-B(x_{2}\otimes x_{1}x_{2};X_{2},gx_{1}x_{2})=0
\end{equation*}%
which holds in view of the form of the elements.

\subsubsection{Case $1_{A}\otimes gx_{1}x_{2}\otimes gx_{2}$}

First summand of the left side gives us%
\begin{eqnarray*}
l_{1} &=&u_{1}=0,l_{2}=u_{2}=0 \\
a &=&b_{1}=b_{2}=0 \\
d &=&e_{1}=e_{2}=1
\end{eqnarray*}%
Since $\alpha \left( x_{2};0,0,0,0\right) \equiv a+b_{1}+b_{2}\equiv 0,$ we
get%
\begin{equation*}
B(g\otimes x_{1}x_{2};1_{A},gx_{1}x_{2})1_{A}\otimes gx_{1}x_{2}\otimes
gx_{2}
\end{equation*}%
Second summand of the left side gives us%
\begin{eqnarray*}
l_{1} &=&u_{1}=0,l_{2}+u_{2}=1 \\
a &=&b_{1}=0,b_{2}=l_{2} \\
d &=&1,e_{1}=1,e_{2}-u_{2}=1\Rightarrow e_{2}=1,u_{2}=0,b_{2}=l_{2}=1
\end{eqnarray*}%
Since $\alpha \left( 1_{H};0,1,0,0\right) \equiv 0$, we get%
\begin{equation*}
B(x_{2}\otimes x_{1}x_{2};X_{2},gx_{1}x_{2})1_{A}\otimes gx_{1}x_{2}\otimes
gx_{2}
\end{equation*}%
By taking in account also the right side we obtain%
\begin{equation*}
B(g\otimes x_{1}x_{2};1_{A},gx_{1}x_{2})-B(x_{2}\otimes
x_{1};1_{A},gx_{1}x_{2})+B(x_{2}\otimes x_{1}x_{2};X_{2},gx_{1}x_{2})=0
\end{equation*}%
which holds in view of the form of the elements.

\subsubsection{Case $X_{2}\otimes gx_{2}\otimes gx_{1}$}

First summand of the left side, nothing. Second summand of the left side
gives us%
\begin{eqnarray*}
l_{1}+u_{1} &=&1,l_{2}=u_{2}=0 \\
a &=&0,b_{1}=l_{1},b_{2}=1 \\
d &=&1,e_{1}=u_{1},e_{2}=1
\end{eqnarray*}%
Since $\alpha \left( 1_{H};0,0,1,0\right) \equiv e_{2}+a+b_{1}+b_{2}\equiv 0$
and $\alpha \left( 1_{H};1,0,0,0\right) \equiv b_{2}=1$, we get%
\begin{equation*}
\left[ B(x_{2}\otimes x_{1}x_{2};X_{2},gx_{1}x_{2})-B(x_{2}\otimes
x_{1}x_{2};X_{1}X_{2},gx_{2})\right] X_{2}\otimes gx_{2}\otimes gx_{1}=0
\end{equation*}%
By taking in account also the right side we get%
\begin{equation*}
+B(x_{2}\otimes x_{1}x_{2};X_{2},gx_{1}x_{2})+B(x_{2}\otimes
x_{2};X_{2},gx_{2})-B(x_{2}\otimes x_{1}x_{2};X_{1}X_{2},gx_{2})=0
\end{equation*}%
which holds in view of the form of the elements.

\subsubsection{Case $X_{2}\otimes g\otimes x_{1}x_{2}$}

First summand of the left side gives us%
\begin{eqnarray*}
l_{1}+u_{1} &=&1,l_{2}=u_{2}=0 \\
a &=&0,b_{1}=l_{1},b_{2}=1 \\
d &=&1,e_{1}=u_{1},e_{2}=0.
\end{eqnarray*}%
Since $\alpha \left( x_{2};0,0,1,0\right) \equiv e_{2}=0$ and $\alpha \left(
x_{2};1,0,0,0\right) \equiv b_{2}+\left( a+b_{1}+b_{2}+1\right) +1\equiv 1$
we get%
\begin{equation*}
\left[ B(g\otimes x_{1}x_{2};X_{2},gx_{1})-B(g\otimes
x_{1}x_{2};X_{1}X_{2},g)\right] X_{2}\otimes g\otimes x_{1}x_{2}
\end{equation*}%
Second summand of the left side gives us%
\begin{eqnarray*}
l_{1}+u_{1} &=&1,l_{2}+u_{2}=1 \\
a &=&0,b_{1}=l_{1},b_{2}-l_{2}=1\Rightarrow b_{2}=1,l_{2}=0\Rightarrow
u_{2}=1 \\
d &=&1,e_{1}=u_{1},e_{2}=u_{2}=1.
\end{eqnarray*}%
Since $\alpha \left( 1_{H};0,0,1,1\right) \equiv 1+e_{2}\equiv 0$ and $%
\alpha \left( 1_{H};1,0,0,1\right) \equiv a+b_{1}=1,$ we get
\begin{equation*}
\left[ B(x_{2}\otimes x_{1}x_{2};X_{2},gx_{1}x_{2})-B(x_{2}\otimes
x_{1}x_{2};X_{1}X_{2},gx_{2})\right] X_{2}\otimes g\otimes x_{1}x_{2}
\end{equation*}%
By taking in account also the right side we obtain%
\begin{gather*}
B(g\otimes x_{1}x_{2};X_{2},gx_{1})-B(g\otimes
x_{1}x_{2};X_{1}X_{2},g)-B(x_{2}\otimes 1_{H};X_{2},g) \\
+B(x_{2}\otimes x_{1}x_{2};X_{2},gx_{1}x_{2})-B(x_{2}\otimes
x_{1}x_{2};X_{1}X_{2},gx_{2})=0
\end{gather*}%
which holds in view of the form of the elements.

\subsubsection{Case $1_{A}\otimes gx_{2}\otimes x_{1}x_{2}$}

First summand of the left side gives us%
\begin{eqnarray*}
l_{1}+u_{1} &=&1,l_{2}=u_{2}=0, \\
a &=&0,b_{1}=l_{1},b_{2}=0, \\
d &=&1,e_{1}=u_{1},e_{2}=1.
\end{eqnarray*}%
Since $\alpha \left( x_{2};0,0,1,0\right) \equiv e_{2}=1$ and $\alpha \left(
x_{2};1,0,0,0\right) \equiv a+b_{1}+1+1\equiv a+b_{1}\equiv 1,$ we obtain
\begin{equation*}
\left[ -B(g\otimes x_{1}x_{2};1_{A},gx_{1}x_{2})-B(g\otimes
x_{1}x_{2};X_{1},gx_{2})\right] 1_{A}\otimes gx_{2}\otimes x_{1}x_{2}
\end{equation*}%
Second summand of the left side gives us%
\begin{eqnarray*}
l_{1}+u_{1} &=&1,l_{2}+u_{2}=1 \\
a &=&0,b_{1}=l_{1}, \\
d &=&1,e_{1}=u_{1},e_{2}-u_{2}=1\Rightarrow e_{2}=1,u_{2}=0\Rightarrow
b_{2}=l_{2}=1
\end{eqnarray*}%
Since $\alpha \left( 1_{H};0,1,1,0\right) \equiv e_{2}+\left(
a+b_{1}+b_{2}+1\right) \equiv 1+0+0+1+1\equiv 1$

and $\alpha \left( 1_{H};1,1,0,0\right) \equiv 1+b_{2}+\left(
a+b_{1}+b_{2}\right) \equiv 1+1\equiv 0,$ we get%
\begin{equation*}
\left[ -B(x_{2}\otimes x_{1}x_{2};X_{2},gx_{1}x_{2})+B(x_{2}\otimes
x_{1}x_{2};X_{1}X_{2},gx_{2})\right] 1_{A}\otimes gx_{2}\otimes x_{1}x_{2}
\end{equation*}%
By taking in account also the right side we obtain%
\begin{gather*}
-B(g\otimes x_{1}x_{2};1_{A},gx_{1}x_{2})-B(g\otimes x_{1}x_{2};X_{1},gx_{2})
\\
-B(x_{2}\otimes 1_{H};1_{A},gx_{2})-B(x_{2}\otimes
x_{1}x_{2};X_{2},gx_{1}x_{2})+B(x_{2}\otimes x_{1}x_{2};X_{1}X_{2},gx_{2})=0
\end{gather*}%
which holds in view of the form of the elements.

\subsection{$B\left( x_{2}\otimes x_{1}x_{2};GX_{1},gx_{2}\right) $}

We deduce that%
\begin{equation*}
a=1,b_{1}=1,b_{2}=0,d=1,e_{1}=0,e_{2}=1
\end{equation*}

and we get

\begin{eqnarray*}
&&\left( -1\right) ^{\alpha \left( 1_{H};0,0,0,0\right) }B(x_{2}\otimes
x_{1}x_{2};GX_{1},gx_{2})GX_{1}\otimes gx_{2}\otimes 1_{H} \\
&&\left( -1\right) ^{\alpha \left( 1_{H};0,0,0,1\right) }B(x_{2}\otimes
x_{1}x_{2};GX_{1},gx_{2})GX_{1}\otimes g\otimes gx_{2} \\
&&\left( -1\right) ^{\alpha \left( 1_{H};1,0,0,0\right) }B(x_{2}\otimes
x_{1}x_{2};GX_{1},gx_{2})G\otimes gx_{2}\otimes gx_{1} \\
&&\left( -1\right) ^{\alpha \left( 1_{H};1,0,0,1\right) }B(x_{2}\otimes
x_{1}x_{2};GX_{1},gx_{2})G\otimes g\otimes x_{1}x_{2}
\end{eqnarray*}

\subsubsection{Case $GX_{1}\otimes g\otimes gx_{2}$}

First summand of the left side gives us%
\begin{eqnarray*}
l_{1} &=&u_{1}=0,l_{2}=u_{2}=0 \\
a &=&1,b_{1}=1,b_{2}=0, \\
d &=&1,e_{1}=0,e_{2}=0.
\end{eqnarray*}%
Since $\alpha \left( x_{2};0,0,0,0\right) \equiv a+b_{1}+b_{2}\equiv 0,$ we
get%
\begin{equation*}
B(g\otimes x_{1}x_{2};GX_{1},g)GX_{1}\otimes g\otimes gx_{2}
\end{equation*}%
Second summand of the left side gives us

\begin{eqnarray*}
l_{1} &=&u_{1}=0,l_{2}+u_{2}=1 \\
a &=&1,b_{1}=1,b_{2}=l_{2} \\
d &=&1,e_{1}=0,e_{2}=u_{2}
\end{eqnarray*}%
Since $\alpha \left( 1_{H};0,0,0,1\right) \equiv a+b_{1}+b_{2}\equiv
1+0+1\equiv 0$ and $\alpha \left( 1_{H};0,1,0,0\right) \equiv 0,$ we obtain
\begin{equation*}
\left[ B(x_{2}\otimes x_{1}x_{2};GX_{1},gx_{2})+B(x_{2}\otimes
x_{1}x_{2};GX_{1}X_{2},g)\right] GX_{1}\otimes g\otimes gx_{2}.
\end{equation*}%
Since there is nothing on the right side, we get%
\begin{equation*}
B(g\otimes x_{1}x_{2};GX_{1},g)-B(x_{2}\otimes
x_{1};GX_{1},g)+B(x_{2}\otimes x_{1}x_{2};GX_{1},gx_{2})+B(x_{2}\otimes
x_{1}x_{2};GX_{1}X_{2},g)=0
\end{equation*}

which holds in view of the form of the elements.

\subsubsection{Case $G\otimes gx_{2}\otimes gx_{1}$}

This case was already considered in subsection $B\left( x_{2}\otimes
x_{1}x_{2};G,gx_{1}x_{2}\right) .$

\subsubsection{Case $G\otimes g\otimes x_{1}x_{2}$}

This case was already considered in subsection $B\left( x_{2}\otimes
x_{1}x_{2};G,gx_{1}x_{2}\right) $

\subsection{$B\left( x_{2}\otimes x_{1}x_{2};GX_{2},gx_{2}\right) $}

We deduce that%
\begin{eqnarray*}
a &=&1,b_{1}=0,b_{2}=1 \\
d &=&1,e_{1}=0,e_{2}=1
\end{eqnarray*}%
and we get%
\begin{eqnarray*}
&&\left( -1\right) ^{\alpha \left( 1_{H};0,0,0,0\right) }B(x_{2}\otimes
x_{1}x_{2};G^{a}X_{2},g^{d}x_{2})GX_{2}\otimes gx_{2}\otimes 1_{H}+ \\
&&\left( -1\right) ^{\alpha \left( 1_{H};0,0,0,1\right) }B(x_{2}\otimes
x_{1}x_{2};G^{a}X_{2},g^{d}x_{2})GX_{2}\otimes g\otimes gx_{2}+ \\
&&\left( -1\right) ^{\alpha \left( 1_{H};0,1,0,0\right) }B(x_{2}\otimes
x_{1}x_{2};G^{a}X_{2},g^{d}x_{2})G\otimes gx_{2}\otimes gx_{2}+ \\
&&\left( -1\right) ^{\alpha \left( 1_{H};0,1,0,1\right) }B(x_{2}\otimes
x_{1}x_{2};G^{a}X_{2},g^{d}x_{2})GX_{2}^{1-l_{2}}\otimes
gx_{2}^{1-u_{2}}\otimes g^{l_{2}+u_{2}}x_{2}^{1+1}=0
\end{eqnarray*}

\subsubsection{Case $GX_{2}\otimes g\otimes gx_{2}$}

First summand of the left side of the equality gives us

\begin{eqnarray*}
l_{1} &=&u_{1}=0,l_{2}=u_{2}=0 \\
a &=&1,b_{1}=0,b_{2}=1 \\
d &=&1,e_{1}=0,e_{2}=u_{2},
\end{eqnarray*}%
Since $\alpha \left( x_{2};0,0,0,0\right) \equiv 0,$ we get%
\begin{equation*}
B(g\otimes x_{1}x_{2};GX_{2},g)GX_{2}\otimes g\otimes gx_{2}
\end{equation*}%
Second summand of the left side gives us%
\begin{eqnarray*}
l_{1} &=&u_{1}=0,l_{2}+u_{2}=1 \\
a &=&1,b_{1}=0,b_{2}-l_{2}=1\Rightarrow b_{2}=1,l_{2}=0\Rightarrow u_{2}=1 \\
d &=&1,e_{1}=0,e_{2}=u_{2}=1.
\end{eqnarray*}%
Since $\alpha \left( 1_{H};0,0,0,1\right) \equiv a+b_{1}+b_{2}\equiv 0,$ we
obtain%
\begin{equation*}
B(x_{2}\otimes x_{1}x_{2};GX_{2},gx_{2})GX_{2}\otimes g\otimes gx_{2}
\end{equation*}%
By considering also the right side we get

\begin{equation*}
B(g\otimes x_{1}x_{2};GX_{2},g)-B(x_{2}\otimes
x_{1};GX_{2},g)+B(x_{2}\otimes x_{1}x_{2};GX_{2},gx_{2})=0
\end{equation*}%
which holds in view of the form of the elements.

\subsubsection{Case $G\otimes gx_{2}\otimes gx_{2}$}

First summand of the left side of the equality gives us

\begin{eqnarray*}
l_{1} &=&u_{1}=0,l_{2}=u_{2}=0 \\
a &=&1,b_{1}=0,b_{2}=0 \\
d &=&1,e_{1}=0,e_{2}=1.
\end{eqnarray*}%
Since $\alpha \left( x_{2};0,0,0,0\right) \equiv a+b_{1}+b_{2}\equiv 1,$ we
obtain%
\begin{equation*}
-B(g\otimes x_{1}x_{2};G,gx_{2})G\otimes gx_{2}\otimes gx_{2}.
\end{equation*}

Second summand of the left side gives us%
\begin{eqnarray*}
l_{1} &=&u_{1}=0,l_{2}+u_{2}=1 \\
a &=&1,b_{1}=0,b_{2}=l_{2} \\
d &=&1,e_{1}=0,e_{2}-u_{2}=1\Rightarrow e_{2}=1,u_{2}=0\Rightarrow
b_{2}=l_{2}=1.
\end{eqnarray*}%
Since $\alpha \left( 1_{H};0,1,0,0\right) \equiv 0,$ we obtain%
\begin{equation*}
B(x_{2}\otimes x_{1}x_{2};GX_{2},gx_{2})G\otimes gx_{2}\otimes gx_{2}.
\end{equation*}%
By taking in account also the right side, we get%
\begin{equation*}
-B(g\otimes x_{1}x_{2};G,gx_{2})-B(x_{2}\otimes
x_{1};G,gx_{2})+B(x_{2}\otimes x_{1}x_{2};GX_{2},gx_{2})=0
\end{equation*}

which holds in view of the form of the elements.

\subsection{$B(x_{2}\otimes x_{1}x_{2};GX_{1}X_{2},gx_{1}x_{2})$}

We deduce that%
\begin{equation*}
a=b_{1}=b_{2}=1,d=e_{1}=e_{2}=1
\end{equation*}%
and we get%
\begin{gather*}
\left( -1\right) ^{\alpha \left( 1_{H};0,0,0,0\right) }B(x_{2}\otimes
x_{1}x_{2};GX_{1}X_{2},gx_{1}x_{2})GX_{1}X_{2}\otimes gx_{1}x_{2}\otimes
1_{H} \\
+\left( -1\right) ^{\alpha \left( 1_{H};1,0,0,0\right) }B(x_{2}\otimes
x_{1}x_{2};GX_{1}X_{2},gx_{1}x_{2})GX_{2}\otimes gx_{1}x_{2}\otimes gx_{1} \\
+\left( -1\right) ^{\alpha \left( 1_{H};0,1,0,0\right) }B(x_{2}\otimes
x_{1}x_{2};GX_{1}X_{2},gx_{1}x_{2})GX_{1}\otimes gx_{1}x_{2}\otimes gx_{2} \\
+\left( -1\right) ^{\alpha \left( 1_{H};1,1,0,0\right) }B(x_{2}\otimes
x_{1}x_{2};GX_{1}X_{2},gx_{1}x_{2})G\otimes gx_{1}x_{2}\otimes x_{1}x_{2} \\
+\left( -1\right) ^{\alpha \left( 1_{H};0,0,1,0\right) }B(x_{2}\otimes
x_{1}x_{2};GX_{1}X_{2},gx_{1}x_{2})GX_{1}X_{2}\otimes gx_{2}\otimes gx_{1} \\
+\left( -1\right) ^{\alpha \left( 1_{H};1,0,1,0\right) }B(x_{2}\otimes
x_{1}x_{2};GX_{1}X_{2},gx_{1}x_{2})GX_{1}^{1-l_{1}}X_{2}\otimes
gx_{2}\otimes g^{l_{1}+1}x_{1}^{1+1}=0 \\
+\left( -1\right) ^{\alpha \left( 1_{H};0,1,1,0\right) }B(x_{2}\otimes
x_{1}x_{2};GX_{1}X_{2},gx_{1}x_{2})GX_{1}\otimes gx_{2}\otimes x_{1}x_{2} \\
+\left( -1\right) ^{\alpha \left( 1_{H};1,1,1,0\right) }B(x_{2}\otimes
x_{1}x_{2};GX_{1}X_{2},gx_{1}x_{2})GX_{1}^{1-l_{1}}\otimes gx_{2}\otimes
g^{l_{1}}x_{1}^{1+1}x_{2}=0 \\
+\left( -1\right) ^{\alpha \left( 1_{H};0,0,0,1\right) }B(x_{2}\otimes
x_{1}x_{2};GX_{1}X_{2},gx_{1}x_{2})GX_{1}X_{2}\otimes gx_{1}\otimes gx_{2} \\
+\left( -1\right) ^{\alpha \left( 1_{H};1,0,0,1\right) }B(x_{2}\otimes
x_{1}x_{2};GX_{1}X_{2},gx_{1}x_{2})GX_{2}\otimes gx_{1}\otimes x_{1}x_{2} \\
+\left( -1\right) ^{\alpha \left( 1_{H};l_{1},1,0,1\right) }B(x_{2}\otimes
x_{1}x_{2};GX_{1}X_{2},gx_{1}x_{2})GX_{1}^{1-l_{1}}X_{2}^{1-l_{2}}\otimes
gx_{1}\otimes g^{l_{1}+l_{2}+1}x_{1}^{l_{1}}x_{2}^{1+1}=0 \\
+\left( -1\right) ^{\alpha \left( 1_{H};0,0,1,1\right) }B(x_{2}\otimes
x_{1}x_{2};GX_{1}X_{2},gx_{1}x_{2})GX_{1}X_{2}\otimes g\otimes x_{1}x_{2} \\
+\left( -1\right) ^{\alpha \left( 1_{H};1,0,1,1\right) }B(x_{2}\otimes
x_{1}x_{2};GX_{1}X_{2},gx_{1}x_{2})GX_{1}^{1-l_{1}}X_{2}\otimes g\otimes
g^{l_{1}+}x_{1}^{1+1}x_{2}=0 \\
+\left( -1\right) ^{\alpha \left( 1_{H};1,1,1,1\right) }B(x_{2}\otimes
x_{1}x_{2};GX_{1}X_{2},gx_{1}x_{2})GX_{1}^{1-l_{1}}X_{2}^{1-l_{2}}\otimes
g\otimes g^{l_{1}+l_{2}+}x_{1}^{l_{1}+1}x_{2}^{1+1}=0
\end{gather*}

\subsubsection{Case $GX_{2}\otimes gx_{1}x_{2}\otimes gx_{1}$}

First summand of the left side, nothing. Second summand of the left side
gives us%
\begin{eqnarray*}
l_{1}+u_{1} &=&1,l_{2}=u_{2}=0 \\
a &=&1,b_{1}=l_{1},b_{2}=1 \\
d &=&1,e_{1}-u_{1}=1\Rightarrow e_{1}=1,u_{1}=0,b_{1}=l_{1}=1,e_{2}=1.
\end{eqnarray*}%
Since $\alpha \left( 1_{H};1,0,0,0\right) \equiv b_{2}=1$, we get%
\begin{equation*}
-B(x_{2}\otimes x_{1}x_{2};GX_{1}X_{2},gx_{1}x_{2})GX_{2}\otimes
gx_{1}x_{2}\otimes gx_{1}
\end{equation*}%
By taking in account also the right side we get%
\begin{equation*}
-B(x_{2}\otimes x_{1}x_{2};GX_{1}X_{2},gx_{1}x_{2})+B(x_{2}\otimes
x_{2};GX_{2},gx_{1}x_{2})=0
\end{equation*}%
which holds in view of the form of the elements.

\subsubsection{Case $GX_{1}\otimes gx_{1}x_{2}\otimes gx_{2} $}

First summand of the left side of the equality gives us

\begin{eqnarray*}
l_{1} &=&u_{1}=0,l_{2}=u_{2}=0 \\
a &=&b_{1}=1,b_{2}=0 \\
d &=&e_{1}=e_{2}=1.
\end{eqnarray*}%
Since $\alpha \left( x_{2};0,0,0,0\right) \equiv a+b_{1}+b_{2}\equiv 0,$ we
obtain%
\begin{equation*}
B(g\otimes x_{1}x_{2};GX_{1},gx_{1}x_{2})GX_{1}\otimes gx_{1}x_{2}\otimes
gx_{2}.
\end{equation*}

Second summand of the left side gives us%
\begin{eqnarray*}
l_{1} &=&u_{1}=0,l_{2}+u_{2}=1 \\
a &=&1,b_{1}=1,b_{2}=l_{2} \\
d &=&1,e_{1}=1,e_{2}-u_{2}=1\Rightarrow e_{2}=1,u_{2}=0\Rightarrow
b_{2}=l_{2}=1.
\end{eqnarray*}%
Since $\alpha \left( 1_{H};0,1,0,0\right) \equiv 0,$ we obtain%
\begin{equation*}
B(x_{2}\otimes x_{1}x_{2};GX_{1}X_{2},gx_{1}x_{2})GX_{1}\otimes
gx_{1}x_{2}\otimes gx_{2}..
\end{equation*}%
By taking in account also the right side, we get

\begin{equation*}
B(g\otimes x_{1}x_{2};GX_{1},gx_{1}x_{2})+B(x_{2}\otimes
x_{1}x_{2};GX_{1}X_{2},gx_{1}x_{2})-B(x_{2}\otimes
x_{1};GX_{1},gx_{1}x_{2})=0
\end{equation*}%
which holds in view of the form of the elements.

\subsubsection{Case $G\otimes gx_{1}x_{2}\otimes x_{1}x_{2}$}

First summand of the left side gives us%
\begin{eqnarray*}
l_{1}+u_{1} &=&1,l_{2}=u_{2}=0, \\
a &=&1,b_{1}=l_{1},b_{2}=0, \\
d &=&1,e_{1}-u_{1}=1\Rightarrow e_{1}=1,u_{1}=0,b_{1}=l_{1}=1,e_{2}=1.
\end{eqnarray*}%
Since $\alpha \left( x_{2};1,0,0,0\right) \equiv a+b_{1}+1+1\equiv
a+b_{1}\equiv 0,$ we obtain
\begin{equation*}
B(g\otimes x_{1}x_{2};GX_{1},gx_{1}x_{2})G\otimes gx_{1}x_{2}\otimes
x_{1}x_{2}
\end{equation*}%
Second summand of the left side gives us%
\begin{eqnarray*}
l_{1}+u_{1} &=&1,l_{2}+u_{2}=1 \\
a &=&1,b_{1}=l_{1},b_{2}=l_{2} \\
d &=&1,e_{1}-u_{1}=1\Rightarrow e_{1}=1,u_{1}=0,b_{1}=l_{1}=1, \\
e_{2}-u_{2} &=&1\Rightarrow e_{2}=1,u_{2}=0\Rightarrow b_{2}=l_{2}=1
\end{eqnarray*}%
Since $\alpha \left( 1_{H};1,1,0,0\right) \equiv 1+b_{2}\equiv 0,$ we get%
\begin{equation*}
+B(x_{2}\otimes x_{1}x_{2};GX_{1}X_{2},gx_{1}x_{2})G\otimes
gx_{1}x_{2}\otimes x_{1}x_{2}
\end{equation*}%
By taking in account also the right side we obtain%
\begin{equation*}
B(g\otimes x_{1}x_{2};GX_{1},gx_{1}x_{2})+B(x_{2}\otimes
x_{1}x_{2};GX_{1}X_{2},gx_{1}x_{2})-B(x_{2}\otimes 1_{H};G,gx_{1}x_{2})=0
\end{equation*}%
which holds in view of the form of the elements.

\subsubsection{Case $GX_{1}X_{2}\otimes gx_{2}\otimes gx_{1}$}

First summand of the left side, nothing. Second summand of the left side
gives us%
\begin{eqnarray*}
l_{1}+u_{1} &=&1,l_{2}=u_{2}=0 \\
a &=&1,b_{1}-l_{1}=1\Rightarrow b_{1}=1,l_{1}=0,u_{1}=1,b_{2}=1 \\
d &=&1,e_{1}=u_{1}=1,e_{2}=1.
\end{eqnarray*}%
Since $\alpha \left( 1_{H};0,0,1,0\right) \equiv e_{2}+\left(
a+b_{1}+b_{2}\right) \equiv 0$, we get%
\begin{equation*}
B(x_{2}\otimes x_{1}x_{2};GX_{1}X_{2},gx_{1}x_{2})GX_{1}X_{2}\otimes
gx_{2}\otimes gx_{1}
\end{equation*}%
By taking in account also the right side we get%
\begin{equation*}
B(x_{2}\otimes x_{1}x_{2};GX_{1}X_{2},gx_{1}x_{2})+B(x_{2}\otimes
x_{2};GX_{1}X_{2},gx_{2})=0
\end{equation*}%
which holds in view of the form of the elements.

\subsubsection{Case $GX_{1}\otimes gx_{2}\otimes x_{1}x_{2}$}

First summand of the left side gives us%
\begin{eqnarray*}
l_{1}+u_{1} &=&1,l_{2}=u_{2}=0, \\
a &=&1,b_{1}-l_{1}=1\Rightarrow b_{1}=1,l_{1}=0,u_{1}=1,b_{2}=0, \\
d &=&1,e_{1}=u_{1}=1,e_{2}=1.
\end{eqnarray*}%
Since $\alpha \left( x_{2};0,0,1,0\right) \equiv e_{2}\equiv 1,$ we obtain
\begin{equation*}
-B(g\otimes x_{1}x_{2};GX_{1},gx_{1}x_{2})GX_{1}\otimes gx_{2}\otimes
x_{1}x_{2}
\end{equation*}%
Second summand of the left side gives us%
\begin{eqnarray*}
l_{1}+u_{1} &=&1,l_{2}+u_{2}=1 \\
a &=&1,b_{1}-l_{1}=1\Rightarrow b_{1}=1,l_{1}=0,u_{1}=1,b_{2}=l_{2} \\
d &=&1,e_{1}=u_{1}=1,e_{2}-u_{2}=1\Rightarrow e_{2}=1,u_{2}=0\Rightarrow
b_{2}=l_{2}=1
\end{eqnarray*}%
Since $\alpha \left( 1_{H};0,1,1,0\right) \equiv e_{2}+a+b_{1}+b_{2}+1\equiv
1,$ we get%
\begin{equation*}
-B(x_{2}\otimes x_{1}x_{2};GX_{1}X_{2},gx_{1}x_{2})GX_{1}\otimes
gx_{2}\otimes x_{1}x_{2}
\end{equation*}%
By taking in account also the right side we obtain%
\begin{equation*}
-B(g\otimes x_{1}x_{2};GX_{1},gx_{1}x_{2})-B(x_{2}\otimes
x_{1}x_{2};GX_{1}X_{2},gx_{1}x_{2})-B(x_{2}\otimes 1_{H};GX_{1},gx_{2})=0
\end{equation*}%
which holds in view of the form of the elements.

\subsubsection{Case $GX_{1}X_{2}\otimes gx_{1}\otimes gx_{2}$}

First summand of the left side of the equality gives us

\begin{eqnarray*}
l_{1} &=&u_{1}=0,l_{2}=u_{2}=0 \\
a &=&b_{1}=b_{2}=1 \\
d &=&e_{1}=1,e_{2}=0.
\end{eqnarray*}%
Since $\alpha \left( x_{2};0,0,0,0\right) \equiv a+b_{1}+b_{2}\equiv 1,$ we
obtain%
\begin{equation*}
-B(g\otimes x_{1}x_{2};GX_{1}X_{2},gx_{1})GX_{1}X_{2}\otimes gx_{1}\otimes
gx_{2}.
\end{equation*}

Second summand of the left side gives us%
\begin{eqnarray*}
l_{1} &=&u_{1}=0,l_{2}+u_{2}=1 \\
a &=&1,b_{1}=1,b_{2}-l_{2}=1\Rightarrow b_{2}=1,l_{2}=0,u_{2}=1 \\
d &=&1,e_{1}=1,e_{2}=u_{2}=1.
\end{eqnarray*}%
Since $\alpha \left( 1_{H};0,0,0,1\right) \equiv a+b_{1}+b_{2}\equiv 1,$ we
obtain%
\begin{equation*}
-B(x_{2}\otimes x_{1}x_{2};GX_{1}X_{2},gx_{1}x_{2})GX_{1}X_{2}\otimes
gx_{1}\otimes gx_{2}.
\end{equation*}%
By taking in account also the right side, we get

\begin{equation*}
-B(g\otimes x_{1}x_{2};GX_{1}X_{2},gx_{1})-B(x_{2}\otimes
x_{1}x_{2};GX_{1}X_{2},gx_{1}x_{2})-B(x_{2}\otimes
x_{1};GX_{1}X_{2},gx_{1})=0
\end{equation*}%
which holds in view of the form of the elements.

\subsubsection{Case $GX_{2}\otimes gx_{1}\otimes x_{1}x_{2}$}

First summand of the left side gives us%
\begin{eqnarray*}
l_{1}+u_{1} &=&1,l_{2}=u_{2}=0, \\
a &=&1,b_{1}=l_{1},b_{2}=1, \\
d &=&1,e_{1}-u_{1}=1\Rightarrow e_{1}=1,u_{1}=0,b_{1}=l_{1}=1,e_{2}=0.
\end{eqnarray*}%
Since $\alpha \left( x_{2};1,0,0,0\right) \equiv a+b_{1}+1+1\equiv 0,$ we
obtain
\begin{equation*}
B(g\otimes x_{1}x_{2};GX_{1}X_{2},gx_{1})GX_{2}\otimes gx_{1}\otimes
x_{1}x_{2}
\end{equation*}%
Second summand of the left side gives us%
\begin{eqnarray*}
l_{1}+u_{1} &=&1,l_{2}+u_{2}=1 \\
a &=&1,b_{1}=l_{1},b_{2}-l_{2}=1\Rightarrow b_{2}=1,l_{2}=0,u_{2}=1 \\
d &=&1,e_{1}-u_{1}=1\Rightarrow e_{1}=1,u_{1}=0,b_{1}=l_{1}=1,e_{2}=u_{2}=1
\end{eqnarray*}%
Since $\alpha \left( 1_{H};1,0,0,1\right) \equiv a+b_{1}\equiv 0,$ we get%
\begin{equation*}
B(x_{2}\otimes x_{1}x_{2};GX_{1}X_{2},gx_{1}x_{2})GX_{2}\otimes
gx_{1}\otimes x_{1}x_{2}
\end{equation*}%
By taking in account also the right side we obtain%
\begin{equation*}
B(g\otimes x_{1}x_{2};GX_{1}X_{2},gx_{1})+B(x_{2}\otimes
x_{1}x_{2};GX_{1}X_{2},gx_{1}x_{2})-B(x_{2}\otimes 1_{H};GX_{2},gx_{1})=0
\end{equation*}%
which holds in view of the form of the elements.

\subsubsection{Case $GX_{1}X_{2}\otimes g\otimes x_{1}x_{2}$}

First summand of the left side gives us%
\begin{eqnarray*}
l_{1}+u_{1} &=&1,l_{2}=u_{2}=0, \\
a &=&1,b_{1}-l_{1}=1\Rightarrow b_{1}=1,l_{1}=0,u_{1}=1,b_{2}=1, \\
d &=&1,e_{1}=u_{1}=1,e_{2}=0.
\end{eqnarray*}%
$\alpha \left( x_{2};0,1,0,0\right) \equiv a+b_{1}+b_{2}+1\equiv 0,$ we
obtain
\begin{equation*}
B(g\otimes x_{1}x_{2};GX_{1}X_{2},gx_{1})GX_{1}X_{2}\otimes g\otimes
x_{1}x_{2}
\end{equation*}%
Second summand of the left side gives us%
\begin{eqnarray*}
l_{1}+u_{1} &=&1,l_{2}+u_{2}=1 \\
a &=&1,b_{1}-l_{1}=1\Rightarrow b_{1}=1,l_{1}=0,u_{1}=1 \\
b_{2}-l_{2} &=&1\Rightarrow b_{2}=1,l_{2}=0,u_{2}=1 \\
d &=&1,e_{1}=u_{1}=1,e_{2}=u_{2}=1
\end{eqnarray*}%
Since $\alpha \left( 1_{H};0,0,1,1\right) \equiv 1+e_{2}\equiv 0,$ we get%
\begin{equation*}
B(x_{2}\otimes x_{1}x_{2};GX_{1}X_{2},gx_{1}x_{2})GX_{1}X_{2}\otimes
g\otimes x_{1}x_{2}
\end{equation*}%
By taking in account also the right side we obtain%
\begin{equation*}
B(g\otimes x_{1}x_{2};GX_{1}X_{2},gx_{1})+B(x_{2}\otimes
x_{1}x_{2};GX_{1}X_{2},gx_{1}x_{2})-B(x_{2}\otimes 1_{H};GX_{1}X_{2},g)=0
\end{equation*}%
which holds in view of the form of the elements.

\section{$B\left( x_{2}\otimes gx_{1}\right) $}

By using $\left( \ref{simplgx}\right) $ we obtain%
\begin{eqnarray}
&&B(x_{2}\otimes gx_{1})  \label{form x2otgx1} \\
&=&(1_{A}\otimes g)B(gx_{2}\otimes 1_{H})(1_{A}\otimes gx_{1})+(1_{A}\otimes
x_{1})B(gx_{2}\otimes 1_{H})  \notag \\
&&+B(x_{1}x_{2}\otimes 1_{H}).  \notag
\end{eqnarray}%
Thus we obtain%
\begin{eqnarray*}
B(x_{2}\otimes gx_{1}) &=&B(x_{1}x_{2}\otimes 1_{H};1_{A},1_{H})1_{A}\otimes
1_{H}+ \\
&&+B(x_{1}x_{2}\otimes 1_{H};1_{A},x_{1}x_{2})1_{A}\otimes x_{1}x_{2}+ \\
&&+\left[ -B(x_{1}x_{2}\otimes 1_{H};1_{A},gx_{1})-2B(x_{1}x_{2}\otimes
1_{H};X_{1},g)\right] 1_{A}\otimes gx_{1} \\
&&+B(x_{1}x_{2}\otimes 1_{H};1_{A},gx_{2})1_{A}\otimes gx_{2}+ \\
&&+B(x_{1}x_{2}\otimes 1_{H};G,g)G\otimes g+ \\
&&+B(x_{1}x_{2}\otimes 1_{H};G,x_{1})G\otimes x_{1}+ \\
&&+B(x_{1}x_{2}\otimes 1_{H};G,x_{2})G\otimes x_{2}+ \\
&&+\left[ -B(x_{1}x_{2}\otimes 1_{H};G,gx_{1}x_{2})-2B(x_{1}x_{2}\otimes
1_{H};GX_{1},gx_{2})\right] G\otimes gx_{1}x_{2} \\
&&+B(x_{1}x_{2}\otimes 1_{H};X_{1},g)X_{1}\otimes g+ \\
&&+B(x_{1}x_{2}\otimes 1_{H};X_{1},x_{1})X_{1}\otimes x_{1}+ \\
&&+B(x_{1}x_{2}\otimes 1_{H};X_{1},x_{2})X_{1}\otimes x_{2}+ \\
&&+B(x_{1}x_{2}\otimes 1_{H};X_{1},gx_{1}x_{2})X_{1}\otimes gx_{1}x_{2}+ \\
&&+B(x_{1}x_{2}\otimes 1_{H};X_{2},g)X_{2}\otimes g+ \\
&&+B(x_{1}x_{2}\otimes 1_{H};X_{2},x_{1})X_{2}\otimes x_{1}+ \\
&&+B(x_{1}x_{2}\otimes 1_{H};X_{2},x_{2})X_{2}\otimes x_{2}+ \\
&&+B(x_{1}x_{2}\otimes 1_{H};X_{2},gx_{1}x_{2})X_{2}\otimes gx_{1}x_{2}+ \\
&&+\left[
\begin{array}{c}
+1-B(x_{1}x_{2}\otimes 1_{H};1_{A},x_{1}x_{2}) \\
-B(x_{1}x_{2}\otimes 1_{H};X_{2},x_{1})+B(x_{1}x_{2}\otimes
1_{H};X_{1},x_{2})%
\end{array}%
\right] X_{1}X_{2}\otimes 1_{H}+ \\
&&-3B(x_{1}x_{2}\otimes 1_{H};X_{1},gx_{1}x_{2})X_{1}X_{2}\otimes gx_{1} \\
&&+B(x_{1}x_{2}\otimes 1_{H};X_{2},gx_{1}x_{2})X_{1}X_{2}\otimes gx_{2} \\
&&+B(x_{1}x_{2}\otimes 1_{H};GX_{1},1_{H})GX_{1}\otimes 1_{H}+ \\
&&+B(x_{1}x_{2}\otimes 1_{H};GX_{1},x_{1}x_{2})GX_{1}\otimes x_{1}x_{2} \\
&&-B(x_{1}x_{2}\otimes 1_{H};GX_{1},gx_{1})GX_{1}\otimes gx_{1} \\
&&+B(x_{1}x_{2}\otimes 1_{H};GX_{1},gx_{2})GX_{1}\otimes gx_{2}+ \\
&&+B(x_{1}x_{2}\otimes 1_{H};GX_{2},1_{H})GX_{2}\otimes 1_{H}+ \\
&&+B(x_{1}x_{2}\otimes 1_{H};GX_{2},x_{1}x_{2})GX_{2}\otimes x_{1}x_{2}+ \\
&&+\left[
\begin{array}{c}
B(x_{1}x_{2}\otimes 1_{H};GX_{2},gx_{1})- \\
2B(x_{1}x_{2}\otimes 1_{H};G,gx_{1}x_{2})-2B(x_{1}x_{2}\otimes
1_{H};GX_{1},gx_{2})%
\end{array}%
\right] GX_{2}\otimes gx_{1}+ \\
&&+B(x_{1}x_{2}\otimes 1_{H};GX_{2},gx_{2})GX_{2}\otimes gx_{2}+ \\
&&+\left[
\begin{array}{c}
-B(x_{1}x_{2}\otimes 1_{H};G,gx_{1}x_{2})+ \\
B(x_{1}x_{2}\otimes 1_{H};GX_{2},gx_{1})-B(x_{1}x_{2}\otimes
1_{H};GX_{1},gx_{2})%
\end{array}%
\right] GX_{1}X_{2}\otimes g \\
&&-B(x_{1}x_{2}\otimes 1_{H};GX_{1},x_{1}x_{2})GX_{1}X_{2}\otimes x_{1}+ \\
&&-B(x_{1}x_{2}\otimes 1_{H};GX_{2},x_{1}x_{2})GX_{1}X_{2}\otimes x_{2}
\end{eqnarray*}%
We set%
\begin{gather*}
C=\left[ -2B(x_{1}x_{2}\otimes 1_{H};1_{A},gx_{1})-2B(x_{1}x_{2}\otimes
1_{H};X_{1},g)\right] 1_{A}\otimes gx_{1}+ \\
+\left[ -2B(x_{1}x_{2}\otimes 1_{H};G,gx_{1}x_{2})-2B(x_{1}x_{2}\otimes
1_{H};GX_{1},gx_{2})\right] G\otimes gx_{1}x_{2}+ \\
-2B(x_{1}x_{2}\otimes 1_{H};X_{1},gx_{1}x_{2})X_{1}\otimes gx_{1}x_{2}+ \\
-2B(x_{1}x_{2}\otimes 1_{H};X_{1},gx_{1}x_{2})X_{1}X_{2}\otimes gx_{1}+ \\
-2B(x_{1}x_{2}\otimes 1_{H};GX_{1},gx_{1})GX_{1}\otimes gx_{1}+ \\
+\left[ -2B(x_{1}x_{2}\otimes 1_{H};G,gx_{1}x_{2})-2B(x_{1}x_{2}\otimes
1_{H};GX_{1},gx_{2})\right] GX_{2}\otimes gx_{1}
\end{gather*}%
Then%
\begin{gather*}
B(x_{2}\otimes gx_{1})=C+B(x_{1}x_{2}\otimes 1_{H}). \\
B(x_{2}\otimes
gx_{1};G^{a}X_{1}^{b_{1}}X_{2}^{b_{2}},g^{d}x_{1}^{e_{1}}x_{2}^{e_{2}})= \\
=C(G^{a}X_{1}^{b_{1}}X_{2}^{b_{2}},g^{d}x_{1}^{e_{1}}x_{2}^{e_{2}})+B(x_{1}x_{2}\otimes 1_{H};G^{a}X_{1}^{b_{1}}X_{2}^{b_{2}},g^{d}x_{1}^{e_{1}}x_{2}^{e_{2}})
\end{gather*}%
so that%
\begin{gather*}
\sum_{a,b_{1},b_{2},d,e_{1},e_{2}=0}^{1}\sum_{l_{1}=0}^{b_{1}}%
\sum_{l_{2}=0}^{b_{2}}\sum_{u_{1}=0}^{e_{1}}\sum_{u_{2}=0}^{e_{2}}\left(
-1\right) ^{\alpha \left( 1_{H};l_{1},l_{2},u_{1},u_{2}\right) } \\
B(x_{2}\otimes
gx_{1};G^{a}X_{1}^{b_{1}}X_{2}^{b_{2}},g^{d}x_{1}^{e_{1}}x_{2}^{e_{2}}) \\
G^{a}X_{1}^{b_{1}-l_{1}}X_{2}^{b_{2}-l_{2}}\otimes
g^{d}x_{1}^{e_{1}-u_{1}}x_{2}^{e_{2}-u_{2}}\otimes
g^{a+b_{1}+b_{2}+l_{1}+l_{2}+d+e_{1}+e_{2}+u_{1}+u_{2}}x_{1}^{l_{1}+u_{1}}x_{2}^{l_{2}+u_{2}}=
\\
\sum_{a,b_{1},b_{2},d,e_{1},e_{2}=0}^{1}\sum_{l_{1}=0}^{b_{1}}%
\sum_{l_{2}=0}^{b_{2}}\sum_{u_{1}=0}^{e_{1}}\sum_{u_{2}=0}^{e_{2}}\left(
-1\right) ^{\alpha \left( 1_{H};l_{1},l_{2},u_{1},u_{2}\right) } \\
\left[
\begin{array}{c}
C(G^{a}X_{1}^{b_{1}}X_{2}^{b_{2}},g^{d}x_{1}^{e_{1}}x_{2}^{e_{2}}) \\
+B(x_{1}x_{2}\otimes
1_{H};G^{a}X_{1}^{b_{1}}X_{2}^{b_{2}},g^{d}x_{1}^{e_{1}}x_{2}^{e_{2}})%
\end{array}%
\right] \\
G^{a}X_{1}^{b_{1}-l_{1}}X_{2}^{b_{2}-l_{2}}\otimes
g^{d}x_{1}^{e_{1}-u_{1}}x_{2}^{e_{2}-u_{2}}\otimes
g^{a+b_{1}+b_{2}+l_{1}+l_{2}+d+e_{1}+e_{2}+u_{1}+u_{2}}x_{1}^{l_{1}+u_{1}}x_{2}^{l_{2}+u_{2}}=
\\
\sum_{a,b_{1},b_{2},d,e_{1},e_{2}=0}^{1}\sum_{l_{1}=0}^{b_{1}}%
\sum_{l_{2}=0}^{b_{2}}\sum_{u_{1}=0}^{e_{1}}\sum_{u_{2}=0}^{e_{2}}\left(
-1\right) ^{\alpha \left( 1_{H};l_{1},l_{2},u_{1},u_{2}\right)
}C(G^{a}X_{1}^{b_{1}}X_{2}^{b_{2}},g^{d}x_{1}^{e_{1}}x_{2}^{e_{2}}) \\
G^{a}X_{1}^{b_{1}-l_{1}}X_{2}^{b_{2}-l_{2}}\otimes
g^{d}x_{1}^{e_{1}-u_{1}}x_{2}^{e_{2}-u_{2}}\otimes
g^{a+b_{1}+b_{2}+l_{1}+l_{2}+d+e_{1}+e_{2}+u_{1}+u_{2}}x_{1}^{l_{1}+u_{1}}x_{2}^{l_{2}+u_{2}}+
\\
\sum_{a,b_{1},b_{2},d,e_{1},e_{2}=0}^{1}\sum_{l_{1}=0}^{b_{1}}%
\sum_{l_{2}=0}^{b_{2}}\sum_{u_{1}=0}^{e_{1}}\sum_{u_{2}=0}^{e_{2}}\left(
-1\right) ^{\alpha \left( 1_{H};l_{1},l_{2},u_{1},u_{2}\right)
}B(x_{1}x_{2}\otimes
1_{H};G^{a}X_{1}^{b_{1}}X_{2}^{b_{2}},g^{d}x_{1}^{e_{1}}x_{2}^{e_{2}}) \\
G^{a}X_{1}^{b_{1}-l_{1}}X_{2}^{b_{2}-l_{2}}\otimes
g^{d}x_{1}^{e_{1}-u_{1}}x_{2}^{e_{2}-u_{2}}\otimes
g^{a+b_{1}+b_{2}+l_{1}+l_{2}+d+e_{1}+e_{2}+u_{1}+u_{2}}x_{1}^{l_{1}+u_{1}}x_{2}^{l_{2}+u_{2}}
\end{gather*}%
We now consider the Casimir formula for $B(x_{2}\otimes gx_{1})$
\begin{eqnarray*}
&&\sum_{a,b_{1},b_{2},d,e_{1},e_{2}=0}^{1}\sum_{l_{1}=0}^{b_{1}}%
\sum_{l_{2}=0}^{b_{2}}\sum_{u_{1}=0}^{e_{1}}\sum_{u_{2}=0}^{e_{2}}\left(
-1\right) ^{\alpha \left( x_{2};l_{1},l_{2},u_{1},u_{2}\right) } \\
&&B(g\otimes
gx_{1};G^{a}X_{1}^{b_{1}}X_{2}^{b_{2}},g^{d}x_{1}^{e_{1}}x_{2}^{e_{2}}) \\
&&G^{a}X_{1}^{b_{1}-l_{1}}X_{2}^{b_{2}-l_{2}}\otimes
g^{d}x_{1}^{e_{1}-u_{1}}x_{2}^{e_{2}-u_{2}}\otimes
g^{a+b_{1}+b_{2}+l_{1}+l_{2}+d+e_{1}+e_{2}+u_{1}+u_{2}}x_{1}^{l_{1}+u_{1}}x_{2}^{l_{2}+u_{2}+1}+
\\
&&+\sum_{a,b_{1},b_{2},d,e_{1},e_{2}=0}^{1}\sum_{l_{1}=0}^{b_{1}}%
\sum_{l_{2}=0}^{b_{2}}\sum_{u_{1}=0}^{e_{1}}\sum_{u_{2}=0}^{e_{2}}\left(
-1\right) ^{\alpha \left( 1_{H};l_{1},l_{2},u_{1},u_{2}\right) } \\
&&B(x_{2}\otimes
gx_{1};G^{a}X_{1}^{b_{1}}X_{2}^{b_{2}},g^{d}x_{1}^{e_{1}}x_{2}^{e_{2}}) \\
&&G^{a}X_{1}^{b_{1}-l_{1}}X_{2}^{b_{2}-l_{2}}\otimes
g^{d}x_{1}^{e_{1}-u_{1}}x_{2}^{e_{2}-u_{2}}\otimes
g^{a+b_{1}+b_{2}+l_{1}+l_{2}+d+e_{1}+e_{2}+u_{1}+u_{2}}x_{1}^{l_{1}+u_{1}}x_{2}^{l_{2}+u_{2}}
\\
&=&B^{A}(x_{2}\otimes gx_{1})\otimes B^{H}(x_{2}\otimes gx_{1})\otimes 1_{H}+
\\
&&B^{A}(x_{2}\otimes g)\otimes B^{H}(x_{2}\otimes g)\otimes gx_{1}+
\end{eqnarray*}%
which we can rewrite%
\begin{gather*}
\sum_{a,b_{1},b_{2},d,e_{1},e_{2}=0}^{1}\sum_{l_{1}=0}^{b_{1}}%
\sum_{l_{2}=0}^{b_{2}}\sum_{u_{1}=0}^{e_{1}}\sum_{u_{2}=0}^{e_{2}}\left(
-1\right) ^{\alpha \left( x_{2};l_{1},l_{2},u_{1},u_{2}\right) }B(g\otimes
gx_{1};G^{a}X_{1}^{b_{1}}X_{2}^{b_{2}},g^{d}x_{1}^{e_{1}}x_{2}^{e_{2}}) \\
G^{a}X_{1}^{b_{1}-l_{1}}X_{2}^{b_{2}-l_{2}}\otimes
g^{d}x_{1}^{e_{1}-u_{1}}x_{2}^{e_{2}-u_{2}}\otimes
g^{a+b_{1}+b_{2}+l_{1}+l_{2}+d+e_{1}+e_{2}+u_{1}+u_{2}}x_{1}^{l_{1}+u_{1}}x_{2}^{l_{2}+u_{2}+1}+
\\
\sum_{a,b_{1},b_{2},d,e_{1},e_{2}=0}^{1}\sum_{l_{1}=0}^{b_{1}}%
\sum_{l_{2}=0}^{b_{2}}\sum_{u_{1}=0}^{e_{1}}\sum_{u_{2}=0}^{e_{2}}\left(
-1\right) ^{\alpha \left( 1_{H};l_{1},l_{2},u_{1},u_{2}\right)
}C(G^{a}X_{1}^{b_{1}}X_{2}^{b_{2}},g^{d}x_{1}^{e_{1}}x_{2}^{e_{2}}) \\
G^{a}X_{1}^{b_{1}-l_{1}}X_{2}^{b_{2}-l_{2}}\otimes
g^{d}x_{1}^{e_{1}-u_{1}}x_{2}^{e_{2}-u_{2}}\otimes
g^{a+b_{1}+b_{2}+l_{1}+l_{2}+d+e_{1}+e_{2}+u_{1}+u_{2}}x_{1}^{l_{1}+u_{1}}x_{2}^{l_{2}+u_{2}}+
\\
\sum_{a,b_{1},b_{2},d,e_{1},e_{2}=0}^{1}\sum_{l_{1}=0}^{b_{1}}%
\sum_{l_{2}=0}^{b_{2}}\sum_{u_{1}=0}^{e_{1}}\sum_{u_{2}=0}^{e_{2}}\left(
-1\right) ^{\alpha \left( 1_{H};l_{1},l_{2},u_{1},u_{2}\right)
}B(x_{1}x_{2}\otimes
1_{H};G^{a}X_{1}^{b_{1}}X_{2}^{b_{2}},g^{d}x_{1}^{e_{1}}x_{2}^{e_{2}}) \\
G^{a}X_{1}^{b_{1}-l_{1}}X_{2}^{b_{2}-l_{2}}\otimes
g^{d}x_{1}^{e_{1}-u_{1}}x_{2}^{e_{2}-u_{2}}\otimes
g^{a+b_{1}+b_{2}+l_{1}+l_{2}+d+e_{1}+e_{2}+u_{1}+u_{2}}x_{1}^{l_{1}+u_{1}}x_{2}^{l_{2}+u_{2}}
\\
=C\otimes 1_{H}+B^{A}(x_{1}x_{2}\otimes 1_{H})\otimes
B^{H}(x_{1}x_{2}\otimes 1_{H})\otimes 1_{H} \\
B^{A}(x_{2}\otimes g)\otimes B^{H}(x_{2}\otimes g)\otimes gx_{1}+
\end{gather*}%
Now by formula \ref{eqX1X2}, we have that%
\begin{eqnarray*}
&&+\sum_{a,b_{1},b_{2},d,e_{1},e_{2}=0}^{1}\sum_{l_{1}=0}^{b_{1}}%
\sum_{l_{2}=0}^{b_{2}}\sum_{u_{1}=0}^{e_{1}}\sum_{u_{2}=0}^{e_{2}}\left(
-1\right) ^{\alpha \left( 1_{H};l_{1},l_{2},u_{1},u_{2}\right) } \\
&&B(x_{1}x_{2}\otimes
1_{H};G^{a}X_{1}^{b_{1}}X_{2}^{b_{2}},g^{d}x_{1}^{e_{1}}x_{2}^{e_{2}}) \\
&&G^{a}X_{1}^{b_{1}-l_{1}}X_{2}^{b_{2}-l_{2}}\otimes
g^{d}x_{1}^{e_{1}-u_{1}}x_{2}^{e_{2}-u_{2}}\otimes
g^{a+b_{1}+b_{2}+l_{1}+l_{2}+d+e_{1}+e_{2}+u_{1}+u_{2}}x_{1}^{l_{1}+u_{1}}x_{2}^{l_{2}+u_{2}}
\\
&=&B^{A}(x_{1}x_{2}\otimes 1_{H})\otimes B^{H}(x_{1}x_{2}\otimes
1_{H})\otimes 1_{H}
\end{eqnarray*}%
and hence what we have to consider is%
\begin{gather*}
\sum_{a,b_{1},b_{2},d,e_{1},e_{2}=0}^{1}\sum_{l_{1}=0}^{b_{1}}%
\sum_{l_{2}=0}^{b_{2}}\sum_{u_{1}=0}^{e_{1}}\sum_{u_{2}=0}^{e_{2}}\left(
-1\right) ^{\alpha \left( x_{2};l_{1},l_{2},u_{1},u_{2}\right) } \\
B(g\otimes
gx_{1};G^{a}X_{1}^{b_{1}}X_{2}^{b_{2}},g^{d}x_{1}^{e_{1}}x_{2}^{e_{2}}) \\
G^{a}X_{1}^{b_{1}-l_{1}}X_{2}^{b_{2}-l_{2}}\otimes
g^{d}x_{1}^{e_{1}-u_{1}}x_{2}^{e_{2}-u_{2}}\otimes
g^{a+b_{1}+b_{2}+l_{1}+l_{2}+d+e_{1}+e_{2}+u_{1}+u_{2}}x_{1}^{l_{1}+u_{1}}x_{2}^{l_{2}+u_{2}+1}+
\\
\sum_{a,b_{1},b_{2},d,e_{1},e_{2}=0}^{1}\sum_{l_{1}=0}^{b_{1}}%
\sum_{l_{2}=0}^{b_{2}}\sum_{u_{1}=0}^{e_{1}}\sum_{u_{2}=0}^{e_{2}}\left(
-1\right) ^{\alpha \left( 1_{H};l_{1},l_{2},u_{1},u_{2}\right)
}C(G^{a}X_{1}^{b_{1}}X_{2}^{b_{2}},g^{d}x_{1}^{e_{1}}x_{2}^{e_{2}}) \\
G^{a}X_{1}^{b_{1}-l_{1}}X_{2}^{b_{2}-l_{2}}\otimes
g^{d}x_{1}^{e_{1}-u_{1}}x_{2}^{e_{2}-u_{2}}\otimes
g^{a+b_{1}+b_{2}+l_{1}+l_{2}+d+e_{1}+e_{2}+u_{1}+u_{2}}x_{1}^{l_{1}+u_{1}}x_{2}^{l_{2}+u_{2}}+
\\
=C\otimes 1_{H} \\
B^{A}(x_{2}\otimes g)\otimes B^{H}(x_{2}\otimes g)\otimes gx_{1}
\end{gather*}%
Now
\begin{gather*}
C=\left[ -2B(x_{1}x_{2}\otimes 1_{H};1_{A},gx_{1})-2B(x_{1}x_{2}\otimes
1_{H};X_{1},g)\right] 1_{A}\otimes gx_{1}+ \\
+\left[ -2B(x_{1}x_{2}\otimes 1_{H};G,gx_{1}x_{2})-2B(x_{1}x_{2}\otimes
1_{H};GX_{1},gx_{2})\right] G\otimes gx_{1}x_{2}+ \\
-2B(x_{1}x_{2}\otimes 1_{H};X_{1},gx_{1}x_{2})X_{1}\otimes gx_{1}x_{2}+ \\
-2B(x_{1}x_{2}\otimes 1_{H};X_{1},gx_{1}x_{2})X_{1}X_{2}\otimes gx_{1}+ \\
-2B(x_{1}x_{2}\otimes 1_{H};GX_{1},gx_{1})GX_{1}\otimes gx_{1}+ \\
+\left[ -2B(x_{1}x_{2}\otimes 1_{H};G,gx_{1}x_{2})-2B(x_{1}x_{2}\otimes
1_{H};GX_{1},gx_{2})\right] GX_{2}\otimes gx_{1}
\end{gather*}%
and we proceed in the usual way.

\subsection{$B(x_{2}\otimes gx_{1};1_{A},gx_{1})$}

\begin{eqnarray*}
a &=&b_{1}=b_{2}=0 \\
d &=&e_{1}=1,e_{2}=0
\end{eqnarray*}%
\begin{eqnarray*}
&&\left( -1\right) ^{\alpha \left( 1_{H};0,0,0,0\right) }B(x_{2}\otimes
gx_{1};1_{A},gx_{1})1_{A}\otimes gx_{1}\otimes 1_{H}+ \\
&&\left( -1\right) ^{\alpha \left( 1_{H};0,0,1,0\right) }B(x_{2}\otimes
gx_{1};1_{A},gx_{1})1_{A}\otimes g\otimes gx_{1}+
\end{eqnarray*}

\subsubsection{Case $1_{A}\otimes g\otimes gx_{1}$}

First summand of the left side, nothing. Second summand of the left side
gives us%
\begin{eqnarray*}
l_{1}+u_{1} &=&1,l_{2}=u_{2}=0 \\
a &=&0,b_{1}=l_{1},b_{2}=0 \\
d &=&1,e_{1}=u_{1},e_{2}=0.
\end{eqnarray*}%
Since $\alpha \left( 1_{H};0,0,1,0\right) \equiv e_{2}+\left(
a+b_{1}+b_{2}\right) =0$ and $\alpha \left( 1_{H};1,0,0,0\right) \equiv
b_{2}=0$, we get%
\begin{equation*}
\left[ B(x_{2}\otimes gx_{1};1_{A},gx_{1})+B(x_{2}\otimes gx_{1};X_{1},g)%
\right] 1_{A}\otimes g\otimes gx_{1}
\end{equation*}%
By taking in account also the right side we get%
\begin{equation*}
B(x_{2}\otimes gx_{1};1_{A},gx_{1})+B(x_{2}\otimes
gx_{1};X_{1},g)-B(x_{2}\otimes g;1_{A},g)=0
\end{equation*}%
which holds in view of the form of the elements.

\subsection{$B(x_{2}\otimes gx_{1};G,gx_{1}x_{2})$}

\begin{eqnarray*}
a &=&1,b_{1}=b_{2}=0 \\
d &=&e_{1}=e_{2}=1
\end{eqnarray*}%
\begin{eqnarray*}
&&+\left( -1\right) ^{\alpha \left( 1_{H};0,0,0,0\right) }B(x_{2}\otimes
gx_{1};G,gx_{1}x_{2})G\otimes gx_{1}x_{2}\otimes 1_{H} \\
&&+\left( -1\right) ^{\alpha \left( 1_{H};0,0,1,0\right) }B(x_{2}\otimes
gx_{1};G,gx_{1}x_{2})G\otimes gx_{2}\otimes gx_{1} \\
&&+\left( -1\right) ^{\alpha \left( 1_{H};0,0,0,1\right) }B(x_{2}\otimes
gx_{1};G,gx_{1}x_{2})G\otimes gx_{1}\otimes gx_{2} \\
&&+\left( -1\right) ^{\alpha \left( 1_{H};0,0,1,1\right) }B(x_{2}\otimes
gx_{1};G,gx_{1}x_{2})G\otimes g\otimes x_{1}x_{2}
\end{eqnarray*}

\subsubsection{Case $G\otimes gx_{2}\otimes gx_{1}$}

First summand of the left side, nothing. Second summand of the left side
gives us%
\begin{eqnarray*}
l_{1}+u_{1} &=&1,l_{2}=u_{2}=0 \\
a &=&1,b_{1}=l_{1},b_{2}=0 \\
d &=&1,e_{1}=u_{1},e_{2}=1.
\end{eqnarray*}%
Since $\alpha \left( 1_{H};0,0,1,0\right) \equiv e_{2}+\left(
a+b_{1}+b_{2}\right) =0$ and $\alpha \left( 1_{H};1,0,0,0\right) \equiv
b_{2}=0$, we get%
\begin{equation*}
\left[ B(x_{2}\otimes gx_{1};G,gx_{1}x_{2})+B(x_{2}\otimes
gx_{1};GX_{1},gx_{2})\right] G\otimes gx_{2}\otimes gx_{1}
\end{equation*}%
By taking in account also the right side we get%
\begin{equation*}
B(x_{2}\otimes gx_{1};G,gx_{1}x_{2})+B(x_{2}\otimes
gx_{1};GX_{1},gx_{2})-B(x_{2}\otimes g;G,gx_{2})=0
\end{equation*}%
which holds in view of the form of the elements.

\subsubsection{Case $G\otimes gx_{1}\otimes gx_{2}$}

First summand of the left side of the equality gives us

\begin{eqnarray*}
l_{1} &=&u_{1}=0,l_{2}=u_{2}=0 \\
a &=&1 \\
d &=&e_{1}=1,e_{2}=0.
\end{eqnarray*}%
Since $\alpha \left( x_{2};0,0,0,0\right) \equiv a+b_{1}+b_{2}\equiv 1,$ we
obtain%
\begin{equation*}
-B(g\otimes gx_{1};G,gx_{1})G\otimes gx_{1}\otimes gx_{2}.
\end{equation*}

Second summand of the left side gives us%
\begin{eqnarray*}
l_{1} &=&u_{1}=0,l_{2}+u_{2}=1 \\
a &=&1,b_{1}=0,b_{2}=l_{2} \\
d &=&1,e_{1}=1,e_{2}=u_{2}.
\end{eqnarray*}%
Since $\alpha \left( 1_{H};0,0,0,1\right) \equiv a+b_{1}+b_{2}\equiv 1$ and $%
\alpha \left( 1_{H};0,1,0,0\right) \equiv 0,$ we obtain%
\begin{equation*}
\left[ -B(x_{2}\otimes gx_{1};G,gx_{1}x_{2})+B(x_{2}\otimes
gx_{1};GX_{2},gx_{1})\right] G\otimes gx_{1}\otimes gx_{2}.
\end{equation*}%
Since there is nothing in the right side, we get%
\begin{equation*}
-B(g\otimes gx_{1};G,gx_{1})-B(x_{2}\otimes
gx_{1};G,gx_{1}x_{2})+B(x_{2}\otimes gx_{1};GX_{2},gx_{1})=0
\end{equation*}%
which holds in view of the form of the elements.

\subsubsection{Case $G\otimes g\otimes x_{1}x_{2}$}

First summand of the left side gives us%
\begin{eqnarray*}
l_{1}+u_{1} &=&1,l_{2}=u_{2}=0, \\
a &=&1,b_{1}=l_{1},b_{2}=0, \\
d &=&1,e_{1}=u_{1},e_{2}=0.
\end{eqnarray*}%
$\alpha \left( x_{2};0,0,1,0\right) \equiv $ $e_{2}=0$ and $\alpha \left(
x_{2};1,0,0,0\right) \equiv b_{2}+\left( a+b_{1}+b_{2}+1\right) +1\equiv
a+b_{1}\equiv 0,$ we obtain
\begin{equation*}
\left[ B(g\otimes gx_{1};G,gx_{1})+B(g\otimes gx_{1};GX_{1},g)\right]
G\otimes g\otimes x_{1}x_{2}
\end{equation*}%
Second summand of the left side gives us%
\begin{eqnarray*}
l_{1}+u_{1} &=&1,l_{2}+u_{2}=1 \\
a &=&1,b_{1}=l_{1},b_{2}=l_{2} \\
d &=&1,e_{1}=u_{1},e_{2}=u_{2}
\end{eqnarray*}%
Since
\begin{eqnarray*}
\alpha \left( 1_{H};0,0,1,1\right) &=&1+e_{2}\equiv 0\text{, }\alpha \left(
1_{H};0,1,1,0\right) =e_{2}+a+b_{1}+b_{2}+1\equiv 1, \\
\alpha \left( 1_{H};1,0,0,1\right) &\equiv &a+b_{1}\equiv 0\text{and }\alpha
\left( 1_{H};1,1,0,0\right) \equiv 1+b_{2}\equiv 0
\end{eqnarray*}%
we get%
\begin{equation*}
\left[
\begin{array}{c}
B(x_{2}\otimes gx_{1};G,gx_{1}x_{2})-B(x_{2}\otimes gx_{1};GX_{2},gx_{1})+
\\
B(x_{2}\otimes gx_{1};GX_{1},gx_{2})+B(x_{2}\otimes gx_{1};GX_{1}X_{2},g)%
\end{array}%
\right] G\otimes g\otimes x_{1}x_{2}
\end{equation*}%
Since there is nothing in the right side, we obtain%
\begin{gather*}
B(g\otimes gx_{1};G,gx_{1})+B(g\otimes gx_{1};GX_{1},g) \\
B(x_{2}\otimes gx_{1};G,gx_{1}x_{2})-B(x_{2}\otimes gx_{1};GX_{2},gx_{1})+ \\
B(x_{2}\otimes gx_{1};GX_{1},gx_{2})+B(x_{2}\otimes gx_{1};GX_{1}X_{2},g)=0
\end{gather*}%
which holds in view of the form of the elements.

\subsection{$B(x_{2}\otimes gx_{1};X_{1},gx_{1}x_{2})$}

\begin{eqnarray*}
b_{1} &=&1,a=b_{2}=0 \\
d &=&e_{1}=e_{2}=1
\end{eqnarray*}%
\begin{gather*}
\left( -1\right) ^{\alpha \left( 1_{H};0,0,0,0\right) }B(x_{2}\otimes
gx_{1};X_{1},gx_{1}x_{2})X_{1}\otimes gx_{1}x_{2}\otimes 1_{H}+ \\
\left( -1\right) ^{\alpha \left( 1_{H};1,0,0,0\right) }B(x_{2}\otimes
gx_{1};X_{1},gx_{1}x_{2})1_{A}\otimes gx_{1}x_{2}\otimes gx_{1}+ \\
\left( -1\right) ^{\alpha \left( 1_{H};0,0,1,0\right) }B(x_{2}\otimes
gx_{1};X_{1},gx_{1}x_{2})X_{1}\otimes gx_{2}\otimes gx_{1}+ \\
\left( -1\right) ^{\alpha \left( 1_{H};1,0,1,0\right) }B(x_{2}\otimes
gx_{1};X_{1},gx_{1}x_{2})X_{1}^{1-l_{1}}\otimes gx_{2}\otimes
g^{+l_{1}+1}x_{1}^{1+1}+=0 \\
\left( -1\right) ^{\alpha \left( 1_{H};0,0,0,1\right) }B(x_{2}\otimes
gx_{1};X_{1},gx_{1}x_{2})X_{1}\otimes gx_{1}\otimes gx_{2} \\
\left( -1\right) ^{\alpha \left( 1_{H};1,0,0,1\right) }B(x_{2}\otimes
gx_{1};X_{1},gx_{1}x_{2})1_{A}\otimes gx_{1}\otimes x_{1}x_{2} \\
\left( -1\right) ^{\alpha \left( 1_{H};0,0,1,1\right) }B(x_{2}\otimes
gx_{1};X_{1},gx_{1}x_{2})X_{1}\otimes g\otimes x_{1}x_{2} \\
\left( -1\right) ^{\alpha \left( 1_{H};1,0,1,1\right) }B(x_{2}\otimes
gx_{1};X_{1},gx_{1}x_{2})X_{1}^{1-l_{1}}\otimes g\otimes
g^{+l_{1}+}x_{1}^{1+1}x_{2}=0
\end{gather*}

\subsubsection{Case $1_{A}\otimes gx_{1}x_{2}\otimes gx_{1}$}

First summand of the left side, nothing. Second summand of the left side
gives us%
\begin{eqnarray*}
l_{1}+u_{1} &=&1,l_{2}=u_{2}=0 \\
a &=&0,b_{1}=l_{1},b_{2}=0 \\
d &=&1,e_{1}-u_{1}=1\Rightarrow e_{1}=1,u_{1}=0,b_{1}=l_{1}=1,e_{2}=1.
\end{eqnarray*}%
Since $\alpha \left( 1_{H};1,0,0,0\right) \equiv b_{2}=0$, we get%
\begin{equation*}
B(x_{2}\otimes gx_{1};X_{1},gx_{1}x_{2})1_{A}\otimes gx_{1}x_{2}\otimes
gx_{1}
\end{equation*}%
By taking in account also the right side we get%
\begin{equation*}
B(x_{2}\otimes gx_{1};X_{1},gx_{1}x_{2})-B(x_{2}\otimes
g;1_{A},gx_{1}x_{2})=0
\end{equation*}%
which holds in view of the form of the elements.

\subsubsection{Case $X_{1}\otimes gx_{2}\otimes gx_{1}$}

First summand of the left side, nothing. Second summand of the left side
gives us%
\begin{eqnarray*}
l_{1}+u_{1} &=&1,l_{2}=u_{2}=0 \\
a &=&0,b_{1}-l_{1}=1\Rightarrow b_{1}=1,l_{1}=0,u_{1}=1,b_{2}=0 \\
d &=&1,e_{1}=u_{1}=1,e_{2}=1.
\end{eqnarray*}%
Since $\alpha \left( 1_{H};0,0,1,0\right) \equiv e_{2}+\left(
a+b_{1}+b_{2}\right) \equiv 0$, we get%
\begin{equation*}
B(x_{2}\otimes gx_{1};X_{1},gx_{1}x_{2})X_{1}\otimes gx_{2}\otimes gx_{1}
\end{equation*}%
By taking in account also the right side we get%
\begin{equation*}
B(x_{2}\otimes gx_{1};X_{1},gx_{1}x_{2})-B(x_{2}\otimes g;X_{1},gx_{2})=0
\end{equation*}%
which holds in view of the form of the elements.

\subsubsection{Case $X_{1}\otimes gx_{1}\otimes gx_{2}$}

First summand of the left side of the equality gives us

\begin{eqnarray*}
l_{1} &=&u_{1}=0,l_{2}=u_{2}=0 \\
a &=&0,b_{1}=1,b_{2}=0 \\
d &=&e_{1}=1,e_{2}=0.
\end{eqnarray*}%
Since $\alpha \left( x_{2};0,0,0,0\right) \equiv a+b_{1}+b_{2}\equiv 1,$ we
obtain%
\begin{equation*}
-B(g\otimes gx_{1};X_{1},gx_{1})X_{1}\otimes gx_{1}\otimes gx_{2}.
\end{equation*}

Second summand of the left side gives us%
\begin{eqnarray*}
l_{1} &=&u_{1}=0,l_{2}+u_{2}=1 \\
a &=&0,b_{1}=1,b_{2}=l_{2} \\
d &=&1,e_{1}=1,e_{2}=u_{2}.
\end{eqnarray*}%
Since $\alpha \left( 1_{H};0,0,0,1\right) \equiv a+b_{1}+b_{2}\equiv 1$ and $%
\alpha \left( 1_{H};0,1,0,0\right) \equiv 0,$ we obtain%
\begin{equation*}
\left[ -B(x_{2}\otimes gx_{1};X_{1},gx_{1}x_{2})+B(x_{2}\otimes
gx_{1};X_{1}X_{2},gx_{1})\right] X_{1}\otimes gx_{1}\otimes gx_{2}.
\end{equation*}%
Since there is nothing in the right side, we get%
\begin{equation*}
-B(g\otimes gx_{1};X_{1},gx_{1})-B(x_{2}\otimes
gx_{1};X_{1},gx_{1}x_{2})+B(x_{2}\otimes gx_{1};X_{1}X_{2},gx_{1})=0
\end{equation*}%
which holds in view of the form of the elements.

\subsubsection{Case $1_{A}\otimes gx_{1}\otimes x_{1}x_{2}$}

First summand of the left side gives us%
\begin{eqnarray*}
l_{1}+u_{1} &=&1,l_{2}=u_{2}=0, \\
a &=&0,b_{1}=l_{1},b_{2}=0, \\
d &=&1,e_{1}-u_{1}=1\Rightarrow e_{1}=1,u_{1}=0,b_{1}=l_{1}=1,e_{2}=0.
\end{eqnarray*}%
Since $\alpha \left( x_{2};1,0,0,0\right) \equiv b_{2}+\left(
a+b_{1}+b_{2}+1\right) +1\equiv a+b_{1}\equiv 1,$ we obtain
\begin{equation*}
-B(g\otimes gx_{1};X_{1},gx_{1})1_{A}\otimes gx_{1}\otimes x_{1}x_{2}
\end{equation*}%
Second summand of the left side gives us%
\begin{eqnarray*}
l_{1}+u_{1} &=&1,l_{2}+u_{2}=1 \\
a &=&0,b_{1}=l_{1},b_{2}=l_{2} \\
d &=&1,e_{1}-u_{1}=1\Rightarrow e_{1}=1,u_{1}=0,b_{1}=l_{1}=1,e_{2}=u_{2}
\end{eqnarray*}%
Since $\alpha \left( 1_{H};1,0,0,1\right) \equiv a+b_{1}\equiv 1$ and $%
\alpha \left( 1_{H};1,1,0,0\right) \equiv 1+b_{2}\equiv 0$ we get%
\begin{equation*}
\left[ -B(x_{2}\otimes gx_{1};X_{1},gx_{1}x_{2})+B(x_{2}\otimes
gx_{1};X_{1}X_{2},gx_{1})\right] 1_{A}\otimes gx_{1}\otimes x_{1}x_{2}
\end{equation*}%
Since there is nothing in the right side, we obtain%
\begin{equation*}
-B(g\otimes gx_{1};X_{1},gx_{1})-B(x_{2}\otimes
gx_{1};X_{1},gx_{1}x_{2})+B(x_{2}\otimes gx_{1};X_{1}X_{2},gx_{1})=0
\end{equation*}%
which we already got in case $X_{1}\otimes gx_{1}\otimes gx_{2}.$

\subsubsection{Case $X_{1}\otimes g\otimes x_{1}x_{2}$}

First summand of the left side gives us%
\begin{eqnarray*}
l_{1}+u_{1} &=&1,l_{2}=u_{2}=0, \\
a &=&0,b_{1}-l_{1}=1\Rightarrow b_{1}=1,l_{1}=0,u_{1}=1,b_{2}=0, \\
d &=&1,e_{1}=u_{1}=1,e_{2}=0.
\end{eqnarray*}%
Since $\alpha \left( x_{2};0,0,1,0\right) \equiv e_{2}=0,$ we obtain
\begin{equation*}
B(g\otimes gx_{1};X_{1},gx_{1})X_{1}\otimes g\otimes x_{1}x_{2}
\end{equation*}%
Second summand of the left side gives us%
\begin{eqnarray*}
l_{1}+u_{1} &=&1,l_{2}+u_{2}=1 \\
a &=&0,b_{1}-l_{1}=1\Rightarrow b_{1}=1,l_{1}=0,u_{1}=1,b_{2}=l_{2} \\
d &=&1,e_{1}=u_{1}=1,e_{2}=u_{2}
\end{eqnarray*}%
Since $\alpha \left( 1_{H};0,0,1,1\right) \equiv 1+e_{2}\equiv 0$ and $%
\alpha \left( 1_{H};0,1,1,0\right) \equiv e_{2}+a+b_{1}+b_{2}+1\equiv 1$ we
get%
\begin{equation*}
\left[ B(x_{2}\otimes gx_{1};X_{1},gx_{1}x_{2})-B(x_{2}\otimes
gx_{1};X_{1}X_{2},gx_{1})\right] X_{1}\otimes g\otimes x_{1}x_{2}
\end{equation*}%
Since there is nothing in the right side, we obtain%
\begin{equation*}
+B(g\otimes gx_{1};X_{1},gx_{1})+B(x_{2}\otimes
gx_{1};X_{1},gx_{1}x_{2})-B(x_{2}\otimes gx_{1};X_{1}X_{2},gx_{1})=0
\end{equation*}%
which we already got in case $X_{1}\otimes gx_{1}\otimes gx_{2}.$

\subsection{$B(x_{2}\otimes gx_{1};X_{1}X_{2},gx_{1})$}

\begin{eqnarray*}
a &=&0,b_{1}=b_{2}=1 \\
d &=&e_{1}=1,e_{2}=0
\end{eqnarray*}

\begin{gather*}
\left( -1\right) ^{\alpha \left( 1_{H};0,0,0,0\right) }B(x_{2}\otimes
gx_{1};X_{1}X_{2},gx_{1})X_{1}X_{2}\otimes gx_{1}\otimes 1_{H}+ \\
\left( -1\right) ^{\alpha \left( 1_{H};1,0,0,0\right) }B(x_{2}\otimes
gx_{1};X_{1}X_{2},gx_{1})X_{2}\otimes gx_{1}\otimes gx_{1}+ \\
\left( -1\right) ^{\alpha \left( 1_{H};0,1,0,0\right) }B(x_{2}\otimes
gx_{1};X_{1}X_{2},gx_{1})X_{1}\otimes gx_{1}\otimes gx_{2}+ \\
\left( -1\right) ^{\alpha \left( 1_{H};1,1,0,0\right) }B(x_{2}\otimes
gx_{1};X_{1}X_{2},gx_{1})1_{A}\otimes gx_{1}\otimes x_{1}x_{2}+ \\
\left( -1\right) ^{\alpha \left( 1_{H};0,0,1,0\right) }B(x_{2}\otimes
gx_{1};X_{1}X_{2},gx_{1})X_{1}X_{2}\otimes g\otimes gx_{1}+ \\
\left( -1\right) ^{\alpha \left( 1_{H};1,0,1,0\right) }B(x_{2}\otimes
gx_{1};X_{1}X_{2},gx_{1})X_{1}^{1-l_{1}}X_{2}\otimes g\otimes
g^{l_{1}+1}x_{1}^{1+1}=0 \\
\left( -1\right) ^{\alpha \left( 1_{H};0,1,1,0\right) }B(x_{2}\otimes
gx_{1};X_{1}X_{2},gx_{1})X_{1}\otimes g\otimes x_{1}x_{2}+ \\
\left( -1\right) ^{\alpha \left( 1_{H};1,1,1,0\right) }B(x_{2}\otimes
gx_{1};X_{1}X_{2},gx_{1})X_{1}^{1-l_{1}}\otimes g\otimes
g^{l_{1}+}x_{1}^{1+1}x_{2}=0
\end{gather*}

\subsubsection{Case $X_{2}\otimes gx_{1}\otimes gx_{1}$}

First summand of the left side, nothing. Second summand of the left side
gives us%
\begin{eqnarray*}
l_{1}+u_{1} &=&1,l_{2}=u_{2}=0 \\
a &=&0,b_{1}=l_{1},b_{2}=1 \\
d &=&1,e_{1}-u_{1}=1\Rightarrow e_{1}=1,u_{1}=0,b_{1}=l_{1}=1,e_{2}=0.
\end{eqnarray*}%
Since $\alpha \left( 1_{H};1,0,0,0\right) \equiv b_{2}=1$, we get%
\begin{equation*}
-B(x_{2}\otimes gx_{1};X_{1}X_{2},gx_{1})X_{2}\otimes gx_{1}\otimes gx_{1}
\end{equation*}%
By taking in account also the right side we get%
\begin{equation*}
-B(x_{2}\otimes gx_{1};X_{1}X_{2},gx_{1})-B(x_{2}\otimes g;X_{2},gx_{1})=0
\end{equation*}%
which holds in view of the form of the elements.

\subsubsection{Case $X_{1}\otimes gx_{1}\otimes gx_{2}$}

This case was already considered in the subsection $B(x_{2}\otimes
gx_{1};X_{1},gx_{1}x_{2}).$

\subsubsection{Case $1_{A}\otimes gx_{1}\otimes x_{1}x_{2}$}

This case was already considered in the subsection $B(x_{2}\otimes
gx_{1};X_{1},gx_{1}x_{2}).$

\subsubsection{Case $X_{1}X_{2}\otimes g\otimes gx_{1}$}

First summand of the left side, nothing. Second summand of the left side
gives us%
\begin{eqnarray*}
l_{1}+u_{1} &=&1,l_{2}=u_{2}=0 \\
a &=&0,b_{1}-l_{1}=1\Rightarrow b_{1}=1,l_{1}=0,u_{1}=1,b_{2}=1 \\
d &=&1,e_{1}=u_{1}=1,e_{2}=0.
\end{eqnarray*}%
Since $\alpha \left( 1_{H};0,0,1,0\right) \equiv e_{2}+\left(
a+b_{1}+b_{2}\right) \equiv 0$, we get%
\begin{equation*}
B(x_{2}\otimes gx_{1};X_{1}X_{2},gx_{1})X_{1}X_{2}\otimes g\otimes gx_{1}
\end{equation*}%
By taking in account also the right side we get%
\begin{equation*}
B(x_{2}\otimes gx_{1};X_{1}X_{2},gx_{1})-B(x_{2}\otimes g;X_{1}X_{2},g)=0
\end{equation*}%
which holds in view of the form of the elements.

\subsubsection{Case $X_{1}\otimes g\otimes x_{1}x_{2}$}

This case was already considered in the subsection $B(x_{2}\otimes
gx_{1};X_{1},gx_{1}x_{2}).$

\subsection{$B(x_{2}\otimes gx_{1};GX_{1},gx_{1})$}

We deduce that%
\begin{eqnarray*}
a &=&b_{1}=1,b_{2}=0 \\
d &=&e_{1}=1,e_{2}=0
\end{eqnarray*}%
and we get%
\begin{gather*}
\left( -1\right) ^{\alpha \left( 1_{H};0,0,0,0\right) }B(x_{2}\otimes
gx_{1};GX_{1},gx_{1})GX_{1}\otimes gx_{1}\otimes 1_{H}+ \\
\left( -1\right) ^{\alpha \left( 1_{H};1,0,0,0\right) }B(x_{2}\otimes
gx_{1};GX_{1},gx_{1})G\otimes gx_{1}\otimes gx_{1}+ \\
\left( -1\right) ^{\alpha \left( 1_{H};0,0,1,0\right) }B(x_{2}\otimes
gx_{1};GX_{1},gx_{1})GX_{1}\otimes g\otimes gx_{1}+ \\
\left( -1\right) ^{\alpha \left( 1_{H};1,0,1,0\right) }B(x_{2}\otimes
gx_{1};GX_{1},gx_{1})GX_{1}^{1-l_{1}}\otimes g\otimes
g^{l_{1}+1}x_{1}^{1+1}=0.
\end{gather*}

\subsubsection{Case $G\otimes gx_{1}\otimes gx_{1}$}

First summand of the left side, nothing. Second summand of the left side
gives us%
\begin{eqnarray*}
l_{1}+u_{1} &=&1,l_{2}=u_{2}=0 \\
a &=&1,b_{1}=l_{1},b_{2}=0 \\
d &=&1,e_{1}-u_{1}=1\Rightarrow e_{1}=1,u_{1}=0,b_{1}=l_{1}=1,e_{2}=0.
\end{eqnarray*}%
Since $\alpha \left( 1_{H};1,0,0,0\right) \equiv b_{2}=0$, we get%
\begin{equation*}
+B(x_{2}\otimes gx_{1};GX_{1},gx_{1})G\otimes gx_{1}\otimes gx_{1}
\end{equation*}%
By taking in account also the right side we get%
\begin{equation*}
B(x_{2}\otimes gx_{1};GX_{1},gx_{1})-B(x_{2}\otimes g;G,gx_{1})=0
\end{equation*}%
which holds in view of the form of the elements.

\subsubsection{Case $GX_{1}\otimes g\otimes gx_{1}$}

First summand of the left side, nothing. Second summand of the left side
gives us%
\begin{eqnarray*}
l_{1}+u_{1} &=&1,l_{2}=u_{2}=0 \\
a &=&1,b_{1}-l_{1}=1\Rightarrow b_{1}=1,l_{1}=0,u_{1}=1,b_{2}=0 \\
d &=&1,e_{1}=u_{1}=1,e_{2}=0.
\end{eqnarray*}%
Since $\alpha \left( 1_{H};0,0,1,0\right) \equiv e_{2}+\left(
a+b_{1}+b_{2}\right) \equiv 0$, we get%
\begin{equation*}
+B(x_{2}\otimes gx_{1};GX_{1},gx_{1})GX_{1}\otimes g\otimes gx_{1}
\end{equation*}%
By taking in account also the right side we get%
\begin{equation*}
B(x_{2}\otimes gx_{1};GX_{1},gx_{1})-B(x_{2}\otimes g;GX_{1},g)=0
\end{equation*}%
which holds in view of the form of the elements.

\subsection{$B(x_{2}\otimes gx_{1};GX_{2},gx_{1})$}

\begin{eqnarray*}
a &=&b_{2}=1,b_{1}=0 \\
d &=&e_{1}=1,e_{2}=0
\end{eqnarray*}%
\begin{eqnarray*}
&&\left( -1\right) ^{\alpha \left( 1_{H};0,0,0,0\right) }B(x_{2}\otimes
gx_{1};GX_{2},gx_{1})GX_{2}\otimes gx_{1}\otimes 1_{H}+ \\
&&\left( -1\right) ^{\alpha \left( 1_{H};0,1,0,0\right) }B(x_{2}\otimes
gx_{1};GX_{2},gx_{1})G\otimes gx_{1}\otimes gx_{2}+ \\
&&\left( -1\right) ^{\alpha \left( 1_{H};0,0,1,0\right) }B(x_{2}\otimes
gx_{1};GX_{2},gx_{1})GX_{2}\otimes g\otimes gx_{1}+ \\
&&\left( -1\right) ^{\alpha \left( 1_{H};0,1,1,0\right) }B(x_{2}\otimes
gx_{1};GX_{2},gx_{1})G\otimes g\otimes x_{1}x_{2}
\end{eqnarray*}

\subsubsection{Case $G\otimes gx_{1}\otimes gx_{2}$}

This case was already considered in subsection $B(x_{2}\otimes
gx_{1};G,gx_{1}x_{2}).$

\subsubsection{Case $GX_{2}\otimes g\otimes gx_{1}$}

First summand of the left side, nothing. Second summand of the left side
gives us%
\begin{eqnarray*}
l_{1}+u_{1} &=&1,l_{2}=u_{2}=0 \\
a &=&1,b_{1}=l_{1},b_{2}=1 \\
d &=&1,e_{1}=u_{1},e_{2}=0.
\end{eqnarray*}%
Since $\alpha \left( 1_{H};0,0,1,0\right) \equiv e_{2}+\left(
a+b_{1}+b_{2}\right) \equiv 0$ and $\alpha \left( 1_{H};1,0,0,0\right)
\equiv 1$we get%
\begin{equation*}
\left[ +B(x_{2}\otimes gx_{1};GX_{2},gx_{1})-B(x_{2}\otimes
gx_{1};GX_{1}X_{2},g)\right] GX_{2}\otimes g\otimes gx_{1}
\end{equation*}%
By taking in account also the right side we get%
\begin{equation*}
B(x_{2}\otimes gx_{1};GX_{2},gx_{1})-B(x_{2}\otimes
gx_{1};GX_{1}X_{2},g)-B(x_{2}\otimes g;GX_{2},g)=0
\end{equation*}%
which holds in view of the form of the elements.

\subsubsection{Case $G\otimes g\otimes x_{1}x_{2}$}

This case was already considered in subsection $B(x_{2}\otimes
gx_{1};G,gx_{1}x_{2})$

.

\section{$B\left( x_{2}\otimes gx_{2}\right) $}

By using $\left( \ref{simplgx}\right) $ we get
\begin{eqnarray}
&&B(x_{2}\otimes gx_{2})  \label{form x2otgx2} \\
&=&(1_{A}\otimes g)B(gx_{2}\otimes 1_{H})(1_{A}\otimes gx_{2})  \notag \\
&&+(1_{A}\otimes x_{2})B(gx_{2}\otimes 1_{H})  \notag
\end{eqnarray}%
and we obtain
\begin{eqnarray*}
B(x_{2}\otimes gx_{2}) &=&\left[ -2B(x_{1}x_{2}\otimes
1_{H};1_{A},gx_{1})-2B(x_{1}x_{2}\otimes 1_{H};X_{1},g)\right] 1_{A}\otimes
gx_{2}+ \\
&&+2B(x_{1}x_{2}\otimes 1_{H};GX_{1},gx_{1})G\otimes gx_{1}x_{2}+ \\
&&-2B(x_{1}x_{2}\otimes 1_{H};X_{1},gx_{1}x_{2})X_{2}\otimes gx_{1}x_{2}+ \\
&&-2B(x_{1}x_{2}\otimes 1_{H};X_{1},gx_{1}x_{2})X_{1}X_{2}\otimes gx_{2}+ \\
&&-2B(x_{1}x_{2}\otimes 1_{H};GX_{1},gx_{1})GX_{1}\otimes gx_{2}+ \\
&&+\left[ -2B(x_{1}x_{2}\otimes 1_{H};G,gx_{1}x_{2})-2B(x_{1}x_{2}\otimes
1_{H};GX_{1},gx_{2})\right] GX_{2}\otimes gx_{2}+
\end{eqnarray*}

We write Casimir condition for $B(x_{2}\otimes gx_{2})$
\begin{eqnarray*}
&&\sum_{a,b_{1},b_{2},d,e_{1},e_{2}=0}^{1}\sum_{l_{1}=0}^{b_{1}}%
\sum_{l_{2}=0}^{b_{2}}\sum_{u_{1}=0}^{e_{1}}\sum_{u_{2}=0}^{e_{2}}\left(
-1\right) ^{\alpha \left( x_{2};l_{1},l_{2},u_{1},u_{2}\right) } \\
&&B(g\otimes
gx_{2};G^{a}X_{1}^{b_{1}}X_{2}^{b_{2}},g^{d}x_{1}^{e_{1}}x_{2}^{e_{2}}) \\
&&G^{a}X_{1}^{b_{1}-l_{1}}X_{2}^{b_{2}-l_{2}}\otimes
g^{d}x_{1}^{e_{1}-u_{1}}x_{2}^{e_{2}-u_{2}}\otimes
g^{a+b_{1}+b_{2}+l_{1}+l_{2}+d+e_{1}+e_{2}+u_{1}+u_{2}}x_{1}^{l_{1}+u_{1}}x_{2}^{l_{2}+u_{2}+1}+
\\
&&+\sum_{a,b_{1},b_{2},d,e_{1},e_{2}=0}^{1}\sum_{l_{1}=0}^{b_{1}}%
\sum_{l_{2}=0}^{b_{2}}\sum_{u_{1}=0}^{e_{1}}\sum_{u_{2}=0}^{e_{2}}\left(
-1\right) ^{\alpha \left( 1_{H};l_{1},l_{2},u_{1},u_{2}\right) } \\
&&B(x_{2}\otimes
gx_{2};G^{a}X_{1}^{b_{1}}X_{2}^{b_{2}},g^{d}x_{1}^{e_{1}}x_{2}^{e_{2}}) \\
&&G^{a}X_{1}^{b_{1}-l_{1}}X_{2}^{b_{2}-l_{2}}\otimes
g^{d}x_{1}^{e_{1}-u_{1}}x_{2}^{e_{2}-u_{2}}\otimes
g^{a+b_{1}+b_{2}+l_{1}+l_{2}+d+e_{1}+e_{2}+u_{1}+u_{2}}x_{1}^{l_{1}+u_{1}}x_{2}^{l_{2}+u_{2}}
\\
&=&B^{A}(x_{2}\otimes gx_{2})\otimes B^{H}(x_{2}\otimes gx_{2})\otimes 1_{H}+
\\
&&B^{A}(x_{2}\otimes g)\otimes B^{H}(x_{2}\otimes g)\otimes gx_{2}+
\end{eqnarray*}

\subsection{$B(x_{2}\otimes gx_{2};1_{A},gx_{2})$}

\begin{eqnarray*}
a &=&b_{1}=b_{2}=0 \\
d &=&e_{2}=1,e_{1}=0
\end{eqnarray*}%
\begin{eqnarray*}
&&\left( -1\right) ^{\alpha \left( 1_{H};0,0,0,0\right) }B(x_{2}\otimes
gx_{2};1_{A},gx_{2})1_{A}\otimes gx_{2}\otimes 1_{H}+ \\
&&\left( -1\right) ^{\alpha \left( 1_{H};0,0,0,1\right) }B(x_{2}\otimes
gx_{2};1_{A},gx_{2})1_{A}\otimes g\otimes gx_{2}
\end{eqnarray*}

\subsubsection{Case $1_{A}\otimes g\otimes gx_{2}$}

First summand of the left side of the equality gives us

\begin{eqnarray*}
l_{1} &=&u_{1}=0,l_{2}=u_{2}=0 \\
a &=&b_{1}=b_{2}=0 \\
d &=&1,e_{1}=e_{2}=0.
\end{eqnarray*}%
Since $\alpha \left( x_{2};0,0,0,0\right) \equiv a+b_{1}+b_{2}\equiv 0,$ we
obtain%
\begin{equation*}
+B(g\otimes gx_{2};1_{A},g)1_{A}\otimes g\otimes gx_{2}.
\end{equation*}

Second summand of the left side gives us%
\begin{eqnarray*}
l_{1} &=&u_{1}=0,l_{2}+u_{2}=1 \\
a &=&0,b_{1}=0,b_{2}=l_{2} \\
d &=&1,e_{1}=0,e_{2}=u_{2}.
\end{eqnarray*}%
Since $\alpha \left( 1_{H};0,0,0,1\right) \equiv a+b_{1}+b_{2}\equiv 0$ and $%
\alpha \left( 1_{H};0,1,0,0\right) \equiv 0,$ we obtain%
\begin{equation*}
\left[ B(x_{2}\otimes gx_{2};1_{A},gx_{2})+B(x_{2}\otimes gx_{2};X_{2},g)%
\right] 1_{A}\otimes g\otimes gx_{2}.
\end{equation*}%
By taking in account also the right side we get, we get%
\begin{equation*}
B(g\otimes gx_{2};1_{A},g)+B(x_{2}\otimes
gx_{2};1_{A},gx_{2})+B(x_{2}\otimes gx_{2};X_{2},g)-B(x_{2}\otimes
g;1_{A},g)=0
\end{equation*}%
which holds in view of the form of the elements.

\subsection{$B(x_{2}\otimes gx_{2};G,gx_{1}x_{2})$}

We deduce that%
\begin{equation*}
a=1,d=e_{1}=e_{2}=1
\end{equation*}%
and we get%
\begin{eqnarray*}
&&\left( -1\right) ^{\alpha \left( 1_{H};0,0,0,0\right) }B\left(
x_{2}\otimes gx_{2};G,gx_{1}x_{2}\right) G\otimes gx_{1}x_{2}\otimes 1_{H} \\
&&\left( -1\right) ^{\alpha \left( 1_{H};0,0,1,0\right) }B\left(
x_{2}\otimes gx_{2};G,gx_{1}x_{2}\right) G\otimes gx_{2}\otimes gx_{1} \\
&&\left( -1\right) ^{\alpha \left( 1_{H};0,0,0,1\right) }B\left(
x_{2}\otimes gx_{2};G,gx_{1}x_{2}\right) G\otimes gx_{1}\otimes gx_{2} \\
&&\left( -1\right) ^{\alpha \left( 1_{H};0,0,1,1\right) }B\left(
x_{2}\otimes gx_{2};G,gx_{1}x_{2}\right) G\otimes g\otimes x_{1}x_{2}
\end{eqnarray*}

\subsubsection{Case $G\otimes gx_{2}\otimes gx_{1}$}

First summand of the left side, nothing. Second summand of the left side
gives us%
\begin{eqnarray*}
l_{1}+u_{1} &=&1,l_{2}=u_{2}=0 \\
a &=&1,b_{1}=l_{1},b_{2}=0 \\
d &=&1,e_{1}=u_{1},e_{2}=1.
\end{eqnarray*}%
Since $\alpha \left( 1_{H};0,0,1,0\right) \equiv e_{2}+\left(
a+b_{1}+b_{2}\right) \equiv 0$ and $\alpha \left( 1_{H};1,0,0,0\right)
\equiv b_{2}\equiv 0$, we get%
\begin{equation*}
\left[ +B(x_{2}\otimes gx_{2};G,gx_{1}x_{2})+B(x_{2}\otimes
gx_{2};GX_{1},gx_{2})\right] G\otimes gx_{2}\otimes gx_{1}
\end{equation*}%
Since there is nothing in the right side we get%
\begin{equation*}
B(x_{2}\otimes gx_{2};G,gx_{1}x_{2})+B(x_{2}\otimes gx_{2};GX_{1},gx_{2})=0
\end{equation*}%
which holds in view of the form of the elements.

\subsubsection{Case $G\otimes gx_{1}\otimes gx_{2}$}

First summand of the left side of the equality gives us

\begin{eqnarray*}
l_{1} &=&u_{1}=0,l_{2}=u_{2}=0 \\
a &=&1,b_{1}=b_{2}=0 \\
d &=&1,e_{1}=1,e_{2}=0.
\end{eqnarray*}%
Since $\alpha \left( x_{2};0,0,0,0\right) \equiv a+b_{1}+b_{2}\equiv 1,$ we
obtain%
\begin{equation*}
-B(g\otimes gx_{2};G,gx_{1})G\otimes gx_{1}\otimes gx_{2}.
\end{equation*}

Second summand of the left side gives us%
\begin{eqnarray*}
l_{1} &=&u_{1}=0,l_{2}+u_{2}=1 \\
a &=&1,b_{1}=0,b_{2}=l_{2} \\
d &=&1,e_{1}=1,e_{2}=u_{2}.
\end{eqnarray*}%
Since $\alpha \left( 1_{H};0,0,0,1\right) \equiv a+b_{1}+b_{2}\equiv 1$ and $%
\alpha \left( 1_{H};0,1,0,0\right) \equiv 0,$ we obtain%
\begin{equation*}
\left[ -B(x_{2}\otimes gx_{2};G,gx_{1}x_{2})+B(x_{2}\otimes
gx_{2};GX_{2},gx_{1})\right] G\otimes gx_{1}\otimes gx_{2}.
\end{equation*}%
By taking in account also the right side we get, we get%
\begin{equation*}
-B(g\otimes gx_{2};G,gx_{1})-B(x_{2}\otimes
gx_{2};G,gx_{1}x_{2})+B(x_{2}\otimes gx_{2};GX_{2},gx_{1})-B(x_{2}\otimes
g;G,gx_{1})=0
\end{equation*}%
\begin{eqnarray*}
-B(g\otimes gx_{2};G,gx_{1}) &=&B(x_{1}x_{2}\otimes 1_{H};GX_{1},gx_{1}) \\
-B(x_{2}\otimes gx_{2};G,gx_{1}x_{2}) &=&-2B(x_{1}x_{2}\otimes
1_{H};GX_{1},gx_{1}) \\
+B(x_{2}\otimes gx_{2};GX_{2},gx_{1}) &=&0 \\
-B(x_{2}\otimes g;G,gx_{1}) &=&B(x_{1}x_{2}\otimes 1_{H};GX_{1},gx_{1})
\end{eqnarray*}%
which holds in view of the form of the elements.

\subsubsection{Case $G\otimes g\otimes x_{1}x_{2}$}

First summand of the left side of the equality gives us%
\begin{eqnarray*}
l_{1}+u_{1} &=&1,l_{2}=u_{2}=0 \\
a &=&1,b_{1}=l_{1},b_{2}=0 \\
d &=&1,e_{1}=u_{1},e_{2}=0.
\end{eqnarray*}%
Since $\alpha \left( x_{2};0,0,1,0\right) \equiv e_{2}=0$ and $\alpha \left(
x_{2};1,0,0,0\right) \equiv b_{2}+a+b_{1}+b_{2}+1+1\equiv 0,$ we obtain
\begin{equation*}
\left[ B(g\otimes gx_{2};G,gx_{1})+B(g\otimes gx_{2};GX_{1},g)\right]
G\otimes g\otimes x_{1}x_{2}.
\end{equation*}%
Second summand of the left side gives us%
\begin{eqnarray*}
l_{1}+u_{1} &=&1,l_{2}+u_{2}=1 \\
a &=&1,b_{1}=l_{1},b_{2}=l_{2} \\
d &=&1,e_{1}=u_{1},e_{2}=u_{2}
\end{eqnarray*}%
Since%
\begin{eqnarray*}
&&\alpha \left( 1_{H};0,0,1,1\right) \equiv 1+e_{2}\equiv 0 \\
&&\alpha \left( 1_{H};0,1,1,0\right) \equiv e_{2}+\left(
a+b_{1}+b_{2}+1\right) \equiv 1 \\
&&\alpha \left( 1_{H};1,0,0,1\right) \equiv a+b_{1}\equiv 0 \\
&&\alpha \left( 1_{H};1,1,0,0\right) \equiv 1+b_{2}\equiv 0
\end{eqnarray*}%
we obtain%
\begin{equation*}
\left[
\begin{array}{c}
B(x_{2}\otimes gx_{2};G,gx_{1}x_{2})-B(x_{2}\otimes gx_{2};GX_{2},gx_{1})+
\\
+B(x_{2}\otimes gx_{2};GX_{1},gx_{2})+B(x_{2}\otimes gx_{2};GX_{1}X_{2},g)%
\end{array}%
\right] G\otimes g\otimes x_{1}x_{2}.
\end{equation*}%
By taking in account also the right side, we get%
\begin{gather*}
B(g\otimes gx_{2};G,gx_{1})+B(g\otimes gx_{2};GX_{1},g)+ \\
+B(x_{2}\otimes gx_{2};G,gx_{1}x_{2})+B(x_{2}\otimes
gx_{2};GX_{1},gx_{2})-B(x_{2}\otimes gx_{2};GX_{2},gx_{1})+B(x_{2}\otimes
gx_{2};GX_{1}X_{2},g)=0
\end{gather*}%
which holds in view of the form of the elements.

\subsection{$B(x_{2}\otimes gx_{2};X_{2},gx_{1}x_{2})$}

We deduce that
\begin{equation*}
b_{2}=1,d=e_{1}=e_{2}=1
\end{equation*}%
and we get

\begin{equation*}
\left( -1\right) ^{\alpha \left( 1_{H};0,l_{2},u_{1},u_{2}\right)
}B(x_{2}\otimes gx_{2};X_{2},gx_{1}x_{2})X_{2}^{1-l_{2}}\otimes
gx_{1}^{1-u_{1}}x_{2}^{1-u_{2}}\otimes
g^{l_{2}+u_{1}+u_{2}}x_{1}^{u_{1}}x_{2}^{l_{2}+u_{2}}
\end{equation*}%
\begin{gather*}
\left( -1\right) ^{\alpha \left( 1_{H};0,0,0,0\right) }B(x_{2}\otimes
gx_{2};X_{2},gx_{1}x_{2})X_{2}\otimes gx_{1}x_{2}\otimes 1_{A} \\
\left( -1\right) ^{\alpha \left( 1_{H};0,1,0,0\right) }B(x_{2}\otimes
gx_{2};X_{2},gx_{1}x_{2})1_{A}\otimes gx_{1}x_{2}\otimes gx_{2} \\
\left( -1\right) ^{\alpha \left( 1_{H};0,0,1,0\right) }B(x_{2}\otimes
gx_{2};X_{2},gx_{1}x_{2})X_{2}\otimes gx_{2}\otimes gx_{1} \\
\left( -1\right) ^{\alpha \left( 1_{H};0,1,1,0\right) }B(x_{2}\otimes
gx_{2};X_{2},gx_{1}x_{2})1_{A}\otimes gx_{2}\otimes x_{1}x_{2} \\
\left( -1\right) ^{\alpha \left( 1_{H};0,0,0,1\right) }B(x_{2}\otimes
gx_{2};X_{2},gx_{1}x_{2})X_{2}\otimes gx_{1}\otimes gx_{2} \\
\left( -1\right) ^{\alpha \left( 1_{H};0,1,0,1\right) }B(x_{2}\otimes
gx_{2};X_{2},gx_{1}x_{2})X_{2}^{1-l_{2}}\otimes gx_{1}\otimes
g^{+l_{2}+1}x_{2}^{1+1}=0 \\
\left( -1\right) ^{\alpha \left( 1_{H};0,0,1,1\right) }B(x_{2}\otimes
gx_{2};X_{2},gx_{1}x_{2})X_{2}\otimes g\otimes x_{1}x_{2} \\
\left( -1\right) ^{\alpha \left( 1_{H};0,1,1,1\right) }B(x_{2}\otimes
gx_{2};X_{2},gx_{1}x_{2})X_{2}^{1-l_{2}}\otimes g\otimes
g^{+l_{2}}x_{1}x_{2}^{1+1}=0
\end{gather*}

\subsubsection{Case $1_{A}\otimes gx_{1}x_{2}\otimes gx_{2}$}

First summand of the left side of the equality gives us

\begin{eqnarray*}
l_{1} &=&u_{1}=0,l_{2}=u_{2}=0 \\
a &=&0 \\
d &=&e_{1}=e_{2}=1.
\end{eqnarray*}%
Since $\alpha \left( x_{2};0,0,0,0\right) \equiv a+b_{1}+b_{2}\equiv 0,$ we
obtain%
\begin{equation*}
B(g\otimes gx_{2};1_{A},gx_{1}x_{2})1_{A}\otimes gx_{1}x_{2}\otimes gx_{2}.
\end{equation*}

Second summand of the left side gives us%
\begin{eqnarray*}
l_{1} &=&u_{1}=0,l_{2}+u_{2}=1 \\
a &=&b_{1}=0,b_{2}=l_{2} \\
d &=&1,e_{1}=1,e_{2}-u_{2}=1\Rightarrow e_{2}=1,u_{2}=0,b_{2}=l_{2}=1.
\end{eqnarray*}%
Since $\alpha \left( 1_{H};0,1,0,0\right) \equiv 0,$ we obtain%
\begin{equation*}
+B(x_{2}\otimes gx_{2};X_{2},gx_{1}x_{2})1_{A}\otimes gx_{1}x_{2}\otimes
gx_{2}.
\end{equation*}%
By taking in account also the right side, we get%
\begin{equation*}
B(g\otimes gx_{2};1_{A},gx_{1}x_{2})+B(x_{2}\otimes
gx_{2};X_{2},gx_{1}x_{2})-B(x_{2}\otimes g;1_{A},gx_{1}x_{2})=0
\end{equation*}%
which holds in view of the form of the elements.

\subsubsection{Case $X_{2}\otimes gx_{2}\otimes gx_{1}$}

First summand of the left side, nothing. Second summand of the left side
gives us%
\begin{eqnarray*}
l_{1}+u_{1} &=&1,l_{2}=u_{2}=0 \\
a &=&0,b_{1}=l_{1},b_{2}=1 \\
d &=&1,e_{1}=u_{1},e_{2}=1.
\end{eqnarray*}%
Since $\alpha \left( 1_{H};0,0,1,0\right) \equiv e_{2}+\left(
a+b_{1}+b_{2}\right) \equiv 0$ and $\alpha \left( 1_{H};1,0,0,0\right)
\equiv b_{2}\equiv 1$, we get%
\begin{equation*}
\left[ +B(x_{2}\otimes gx_{2};X_{2},gx_{1}x_{2})-B(x_{2}\otimes
gx_{2};X_{1}X_{2},gx_{2})\right] X_{2}\otimes gx_{2}\otimes gx_{1}
\end{equation*}%
Since there is nothing in the right side we get%
\begin{equation*}
+B(x_{2}\otimes gx_{2};X_{2},gx_{1}x_{2})-B(x_{2}\otimes
gx_{2};X_{1}X_{2},gx_{2})=0
\end{equation*}%
which holds in view of the form of the elements.

\subsubsection{Case $1_{A}\otimes gx_{2}\otimes x_{1}x_{2}$}

First summand of the left side of the equality gives us%
\begin{eqnarray*}
l_{1}+u_{1} &=&1,l_{2}=u_{2}=0 \\
a &=&0,b_{1}=l_{1},b_{2}=0 \\
d &=&1,e_{1}=u_{1},e_{2}=1.
\end{eqnarray*}%
Since $\alpha \left( x_{2};0,0,1,0\right) \equiv e_{2}=1$ and $\alpha \left(
x_{2};1,0,0,0\right) \equiv b_{2}+a+b_{1}+b_{2}+1+1\equiv 1,$ we obtain
\begin{equation*}
\left[ -B(g\otimes gx_{2};1_{A},gx_{1}x_{2})-B(g\otimes gx_{2};X_{1},gx_{2})%
\right] 1_{A}\otimes gx_{2}\otimes x_{1}x_{2}.
\end{equation*}%
Second summand of the left side gives us%
\begin{eqnarray*}
l_{1}+u_{1} &=&1,l_{2}+u_{2}=1 \\
a &=&0,b_{1}=l_{1},b_{2}=l_{2} \\
d &=&1,e_{1}=u_{1},e_{2}-u_{2}=1\Rightarrow e_{2}=1,u_{2}=0,b_{2}=l_{2}=1
\end{eqnarray*}%
Since $\alpha \left( 1_{H};0,1,1,0\right) \equiv e_{2}+\left(
a+b_{1}+b_{2}+1\right) \equiv 1$ and $\alpha \left( 1_{H};1,1,0,0\right)
\equiv 1+b_{2}\equiv 0$ we obtain%
\begin{equation*}
\left[ -B(x_{2}\otimes gx_{2};X_{2},gx_{1}x_{2})+B(x_{2}\otimes
gx_{2};X_{1}X_{2},gx_{2})\right] G\otimes g\otimes x_{1}x_{2}.
\end{equation*}%
Since there is nothing in the right side, we get%
\begin{gather*}
-B(g\otimes gx_{2};1_{A},gx_{1}x_{2})-B(g\otimes gx_{2};X_{1},gx_{2})+ \\
-B(x_{2}\otimes gx_{2};X_{2},gx_{1}x_{2})+B(x_{2}\otimes
gx_{2};X_{1}X_{2},gx_{2})=0
\end{gather*}%
which holds in view of the form of the elements.

\subsubsection{Case $X_{2}\otimes gx_{1}\otimes gx_{2}$}

First summand of the left side of the equality gives us

\begin{eqnarray*}
l_{1} &=&u_{1}=0,l_{2}=u_{2}=0 \\
a &=&0,b_{2}=1,b_{1}=0 \\
d &=&e_{1}=1,e_{2}=0
\end{eqnarray*}%
Since $\alpha \left( x_{2};0,0,0,0\right) \equiv a+b_{1}+b_{2}\equiv 1,$ we
obtain%
\begin{equation*}
-B(g\otimes gx_{2};X_{2},gx_{1})X_{2}\otimes gx_{1}\otimes gx_{2}.
\end{equation*}

Second summand of the left side gives us%
\begin{eqnarray*}
l_{1} &=&u_{1}=0,l_{2}+u_{2}=1 \\
a &=&b_{1}=0,b_{2}-l_{2}=1\Rightarrow b_{2}=1,l_{2}=0,u_{2}=1 \\
d &=&1,e_{1}=1,e_{2}=u_{2}=1
\end{eqnarray*}%
Since $\alpha \left( 1_{H};0,0,0,1\right) \equiv a+b_{1}+b_{2}\equiv 1,$ we
obtain%
\begin{equation*}
-B(x_{2}\otimes gx_{2};X_{2},gx_{1}x_{2})X_{2}\otimes gx_{1}\otimes gx_{2}.
\end{equation*}%
By taking in account also the right side, we get%
\begin{equation*}
-B(g\otimes gx_{2};X_{2},gx_{1})-B(x_{2}\otimes
gx_{2};X_{2},gx_{1}x_{2})-B(x_{2}\otimes g;X_{2},gx_{1})=0
\end{equation*}%
which holds in view of the form of the elements.

\subsubsection{Case $X_{2}\otimes g\otimes x_{1}x_{2}$}

First summand of the left side of the equality gives us%
\begin{eqnarray*}
l_{1}+u_{1} &=&1,l_{2}=u_{2}=0 \\
a &=&0,b_{1}=l_{1},b_{2}=1 \\
d &=&1,e_{1}=u_{1},e_{2}=0.
\end{eqnarray*}%
Since $\alpha \left( x_{2};0,0,1,0\right) \equiv e_{2}=0$ and $\alpha \left(
x_{2};1,0,0,0\right) \equiv b_{2}+a+b_{1}+b_{2}+1+1\equiv 1,$ we obtain
\begin{equation*}
\left[ B(g\otimes gx_{2};X_{2},gx_{1})-B(g\otimes gx_{2};X_{1}X_{2},g)\right]
X_{2}\otimes g\otimes x_{1}x_{2}.
\end{equation*}%
Second summand of the left side gives us%
\begin{eqnarray*}
l_{1}+u_{1} &=&1,l_{2}+u_{2}=1 \\
a &=&0,b_{1}=l_{1},b_{2}-l_{2}=1\Rightarrow b_{2}=1,l_{2}=0,u_{2}=1 \\
d &=&1,e_{1}=u_{1},e_{2}=u_{2}=1
\end{eqnarray*}%
Since $\alpha \left( 1_{H};0,0,1,1\right) \equiv 1+e_{2}\equiv 0$ and $%
\alpha \left( 1_{H};1,0,0,1\right) \equiv a+b_{1}\equiv 1$ we obtain%
\begin{equation*}
\left[ +B(x_{2}\otimes gx_{2};X_{2},gx_{1}x_{2})-B(x_{2}\otimes
gx_{2};X_{1}X_{2},gx_{2})\right] X_{2}\otimes g\otimes x_{1}x_{2}.
\end{equation*}%
Since there is nothing in the right side, we get%
\begin{gather*}
B(g\otimes gx_{2};X_{2},gx_{1})-B(g\otimes gx_{2};X_{1}X_{2},g)+ \\
+B(x_{2}\otimes gx_{2};X_{2},gx_{1}x_{2})-B(x_{2}\otimes
gx_{2};X_{1}X_{2},gx_{2})=0
\end{gather*}%
which holds in view of the form of the elements.

\subsection{$B(x_{2}\otimes gx_{2};X_{1}X_{2},gx_{2})$}

We have%
\begin{equation*}
a=0,b_{1}=b_{2}=1,d=e_{2}=1
\end{equation*}%
and we get%
\begin{gather*}
+\left( -1\right) ^{\alpha \left( 1_{H};0,0,0,0\right) }B(x_{2}\otimes
gx_{2};X_{1}X_{2},gx_{2})X_{1}X_{2}\otimes gx_{2}\otimes 1_{A}+ \\
+\left( -1\right) ^{\alpha \left( 1_{H};0,1,0,0\right) }B(x_{2}\otimes
gx_{2};X_{1}X_{2},gx_{2})X_{1}\otimes gx_{2}\otimes gx_{2}+ \\
+\left( -1\right) ^{\alpha \left( 1_{H};1,0,0,0\right) }B(x_{2}\otimes
gx_{2};X_{1}X_{2},gx_{2})X_{2}\otimes gx_{2}\otimes gx_{1}+ \\
+\left( -1\right) ^{\alpha \left( 1_{H};1,1,0,0\right) }B(x_{2}\otimes
gx_{2};X_{1}X_{2},gx_{2})1_{A}\otimes gx_{2}\otimes x_{1}x_{2}+ \\
+\left( -1\right) ^{\alpha \left( 1_{H};0,0,0,1\right) }B(x_{2}\otimes
gx_{2};X_{1}X_{2},gx_{2})X_{1}X_{2}\otimes g\otimes gx_{2}+ \\
+\left( -1\right) ^{\alpha \left( 1_{H};0,1,0,1\right) }B(x_{2}\otimes
gx_{2};X_{1}X_{2},gx_{2})X_{1}^{1-l_{1}}X_{2}^{1-l_{2}}\otimes g\otimes
g^{l_{1}+l_{2}+1}x_{1}^{l_{1}}x_{2}^{1+1}=0 \\
+\left( -1\right) ^{\alpha \left( 1_{H};1,0,0,1\right) }B(x_{2}\otimes
gx_{2};X_{1}X_{2},gx_{2})X_{2}\otimes g\otimes x_{1}x_{2}+ \\
+\left( -1\right) ^{\alpha \left( 1_{H};1,1,0,1\right) }B(x_{2}\otimes
gx_{2};X_{1}X_{2},gx_{2})X_{1}^{1-l_{1}}X_{2}^{1-l_{2}}\otimes g\otimes
g^{1+l_{1}+l_{2}}x_{1}^{l_{1}}x_{2}^{1+1}=0.
\end{gather*}

\subsubsection{Case $X_{1}\otimes gx_{2}\otimes gx_{2}$}

First summand of the left side of the equality gives us

\begin{eqnarray*}
l_{1} &=&u_{1}=0,l_{2}=u_{2}=0 \\
a &=&0,b_{1}=1,b_{2}=0 \\
d &=&e_{2}=1,e_{1}=0
\end{eqnarray*}%
Since $\alpha \left( x_{2};0,0,0,0\right) \equiv a+b_{1}+b_{2}\equiv 1,$ we
obtain%
\begin{equation*}
-B(g\otimes gx_{2};X_{1},gx_{2})X_{1}\otimes gx_{2}\otimes gx_{2}.
\end{equation*}

Second summand of the left side gives us%
\begin{eqnarray*}
l_{1} &=&u_{1}=0,l_{2}+u_{2}=1 \\
a &=&0,b_{1}=1,b_{2}=l_{2} \\
d &=&1,e_{1}=0,e_{2}-u_{2}=1\Rightarrow e_{2}=1,u_{2}=0,b_{2}=l_{2}=1.
\end{eqnarray*}%
Since $\alpha \left( 1_{H};0,1,0,0\right) \equiv 0,$ we obtain%
\begin{equation*}
+B(x_{2}\otimes gx_{2};X_{1}X_{2},gx_{2})X_{1}\otimes gx_{2}\otimes gx_{2}.
\end{equation*}%
By taking in account also the right side, we get%
\begin{equation*}
-B(g\otimes gx_{2};X_{1},gx_{2})+B(x_{2}\otimes
gx_{2};X_{1}X_{2},gx_{2})-B(x_{2}\otimes g;X_{1},gx_{2})=0
\end{equation*}%
which holds in view of the form of the elements.

\subsubsection{Case $X_{2}\otimes gx_{2}\otimes gx_{1}$}

Already considered in subsection $B(x_{2}\otimes gx_{2};X_{2},gx_{1}x_{2}).$

\subsubsection{Case $1_{A}\otimes gx_{2}\otimes x_{1}x_{2}$}

Already considered in subsection $B(x_{2}\otimes gx_{2};X_{2},gx_{1}x_{2}).$

\subsubsection{Case $X_{1}X_{2}\otimes g\otimes gx_{2}$}

First summand of the left side of the equality gives us

\begin{eqnarray*}
l_{1} &=&u_{1}=0,l_{2}=u_{2}=0 \\
a &=&0,b_{1}=b_{2}=1, \\
d &=&1,e_{1}=e_{2}=0
\end{eqnarray*}%
Since $\alpha \left( x_{2};0,0,0,0\right) \equiv a+b_{1}+b_{2}\equiv 0,$ we
obtain%
\begin{equation*}
B(g\otimes gx_{2};X_{1}X_{2},g)X_{1}X_{2}\otimes g\otimes gx_{2}.
\end{equation*}

Second summand of the left side gives us%
\begin{eqnarray*}
l_{1} &=&u_{1}=0,l_{2}+u_{2}=1 \\
a &=&0,b_{1}=1,b_{2}-l_{2}=1\Rightarrow b_{2}=1,l_{2}=0,u_{2}=1 \\
d &=&1,e_{1}=0,e_{2}=u_{2}=1
\end{eqnarray*}%
Since $\alpha \left( 1_{H};0,0,0,1\right) \equiv a+b_{1}+b_{2}\equiv 0,$ we
obtain%
\begin{equation*}
B(x_{2}\otimes gx_{2};X_{1}X_{2},gx_{2})X_{1}X_{2}\otimes g\otimes gx_{2}.
\end{equation*}%
By taking in account also the right side, we get%
\begin{equation*}
B(g\otimes gx_{2};X_{1}X_{2},g)+B(x_{2}\otimes
gx_{2};X_{1}X_{2},gx_{2})-B(x_{2}\otimes g;X_{1}X_{2},g)=0
\end{equation*}%
which holds in view of the form of the elements.

\subsubsection{Case $X_{2}\otimes g\otimes x_{1}x_{2}$}

This case was already considered in subsection $B(x_{2}\otimes
gx_{2};X_{2},gx_{1}x_{2}).$

\subsection{$B(x_{2}\otimes gx_{2};GX_{1},gx_{2})$}

We deduce that
\begin{equation*}
a=b_{1}=1,d=1,e_{2}=1
\end{equation*}%
and we get%
\begin{equation*}
GX_{1}^{1-l_{1}}\otimes gx_{2}^{1-u_{2}}\otimes
g^{l_{1}+u_{2}}x_{1}^{l_{1}}x_{2}^{+u_{2}}
\end{equation*}%
\begin{eqnarray*}
&&\left( -1\right) ^{\alpha \left( 1_{H};0,0,0,0\right) }B(x_{2}\otimes
x_{2};X_{1},gx_{2})GX_{1}\otimes gx_{2}\otimes 1_{H}+ \\
&&\left( -1\right) ^{\alpha \left( 1_{H};1,0,0,0\right) }B(x_{2}\otimes
x_{2};X_{1},gx_{2})G\otimes gx_{2}\otimes gx_{1} \\
&&\left( -1\right) ^{\alpha \left( 1_{H};0,0,0,1\right) }B(x_{2}\otimes
x_{2};X_{1},gx_{2})GX_{1}\otimes g\otimes gx_{2} \\
&&\left( -1\right) ^{\alpha \left( 1_{H};1,0,0,1\right) }B(x_{2}\otimes
x_{2};X_{1},gx_{2})G\otimes g\otimes x_{1}x_{2}
\end{eqnarray*}

\subsubsection{Case $G\otimes gx_{2}\otimes gx_{1}$}

This case was already considered in subsection $B(x_{2}\otimes
gx_{2};G,gx_{1}x_{2}).$

\subsubsection{Case $GX_{1}\otimes g\otimes gx_{2}$}

First summand of the left side of the equality gives us

\begin{eqnarray*}
l_{1} &=&u_{1}=0,l_{2}=u_{2}=0 \\
a &=&b_{1}=1,b_{2}=0 \\
d &=&1,e_{1}=e_{2}=0
\end{eqnarray*}%
Since $\alpha \left( x_{2};0,0,0,0\right) \equiv a+b_{1}+b_{2}\equiv 0,$ we
obtain%
\begin{equation*}
B(g\otimes gx_{2};GX_{1},g)GX_{1}\otimes g\otimes gx_{2}.
\end{equation*}

Second summand of the left side gives us%
\begin{eqnarray*}
l_{1} &=&u_{1}=0,l_{2}+u_{2}=1 \\
a &=&1,b_{1}=1,b_{2}=l_{2} \\
d &=&1,e_{1}=0,e_{2}=u_{2}
\end{eqnarray*}%
Since $\alpha \left( 1_{H};0,0,0,1\right) \equiv a+b_{1}+b_{2}\equiv 0$ and $%
\alpha \left( 1_{H};0,1,0,0\right) \equiv 0$ we obtain%
\begin{equation*}
\left[ +B(x_{2}\otimes gx_{2};GX_{1},gx_{2})+B(x_{2}\otimes
gx_{2};GX_{1}X_{2},g)\right] GX_{1}\otimes g\otimes gx_{2}.
\end{equation*}%
By taking in account also the right side, we get%
\begin{equation*}
B(g\otimes gx_{2};GX_{1},g)+B(x_{2}\otimes
gx_{2};GX_{1},gx_{2})+B(x_{2}\otimes gx_{2};GX_{1}X_{2},g)-B(x_{2}\otimes
g;GX_{1},g)=0
\end{equation*}%
which holds in view of the form of the elements.

\subsubsection{Case $G\otimes g\otimes x_{1}x_{2}$}

This case was already considered in subsection $B(x_{2}\otimes
gx_{2};G,gx_{1}x_{2}).$

\subsection{$B(x_{2}\otimes gx_{2};GX_{2},gx_{2})$}

We have%
\begin{equation*}
a=1,b_{2}=1,b_{1}=0,d=e_{2}=1,e_{1}=0
\end{equation*}%
\begin{equation*}
GX_{2}^{1-l_{2}}\otimes gx_{2}^{1-u_{2}}\otimes
g^{l_{2}+u_{2}}x_{2}^{l_{2}+u_{2}}
\end{equation*}%
\begin{gather*}
+\left( -1\right) ^{\alpha \left( 1_{H};0,0,0,0\right) }B(x_{1}\otimes
gx_{1}x_{2};X_{1}X_{2},gx_{1}x_{2})GX_{2}\otimes gx_{2}\otimes 1_{A}+ \\
+\left( -1\right) ^{\alpha \left( 1_{H};0,1,0,0\right) }B(x_{1}\otimes
gx_{1}x_{2};X_{1}X_{2},gx_{1}x_{2})G\otimes gx_{2}\otimes gx_{2}+ \\
+\left( -1\right) ^{\alpha \left( 1_{H};0,0,0,1\right) }B(x_{1}\otimes
gx_{1}x_{2};X_{1}X_{2},gx_{1}x_{2})GX_{2}\otimes g\otimes gx_{2}+ \\
+\left( -1\right) ^{\alpha \left( 1_{H};0,1,0,1\right) }B(x_{1}\otimes
gx_{1}x_{2};X_{1}X_{2},gx_{1}x_{2})GX_{2}^{1-l_{2}}\otimes g\otimes
g^{+l_{2}}x_{2}^{1+1}=0
\end{gather*}

\subsubsection{Case $G\otimes gx_{2}\otimes gx_{2}$}

First summand of the left side of the equality gives us

\begin{eqnarray*}
l_{1} &=&u_{1}=0,l_{2}=u_{2}=0 \\
a &=&1,b_{1}=b_{2}=0 \\
d &=&1,e_{1}=0,e_{2}=1
\end{eqnarray*}%
Since $\alpha \left( x_{2};0,0,0,0\right) \equiv a+b_{1}+b_{2}\equiv 1,$ we
obtain%
\begin{equation*}
-B(g\otimes gx_{2};G,gx_{2})G\otimes gx_{2}\otimes gx_{2}.
\end{equation*}

Second summand of the left side gives us%
\begin{eqnarray*}
l_{1} &=&u_{1}=0,l_{2}+u_{2}=1 \\
a &=&1,b_{1}=0,b_{2}=l_{2} \\
d &=&1,e_{1}=0,e_{2}-u_{2}=1\Rightarrow e_{2}=1,u_{2}=0,b_{2}=l_{2}=1
\end{eqnarray*}%
Since $\alpha \left( 1_{H};0,1,0,0\right) \equiv 0$ we obtain%
\begin{equation*}
B(x_{2}\otimes gx_{2};GX_{2},gx_{2})G\otimes gx_{2}\otimes gx_{2}..
\end{equation*}%
By taking in account also the right side, we get%
\begin{equation*}
-B(g\otimes gx_{2};G,gx_{2})+B(x_{2}\otimes
gx_{2};GX_{2},gx_{2})-B(x_{2}\otimes g;G,gx_{2})=0
\end{equation*}%
which holds in view of the form of the elements.

\subsubsection{Case $GX_{2}\otimes g\otimes gx_{2}$}

First summand of the left side of the equality gives us

\begin{eqnarray*}
l_{1} &=&u_{1}=0,l_{2}=u_{2}=0 \\
a &=&b_{2}=1,b_{1}=0 \\
d &=&1,e_{1}=0,e_{2}=0
\end{eqnarray*}%
Since $\alpha \left( x_{2};0,0,0,0\right) \equiv a+b_{1}+b_{2}\equiv 0,$ we
obtain%
\begin{equation*}
B(g\otimes gx_{2};GX_{2},g)GX_{2}\otimes g\otimes gx_{2}.
\end{equation*}

Second summand of the left side gives us%
\begin{eqnarray*}
l_{1} &=&u_{1}=0,l_{2}+u_{2}=1 \\
a &=&1,b_{1}=0,b_{2}-l_{2}=1\Rightarrow b_{2}=1,l_{2}=0,u_{2}=1 \\
d &=&1,e_{1}=0,e_{2}=u_{2}=1
\end{eqnarray*}%
Since $\alpha \left( 1_{H};0,0,0,1\right) \equiv a+b_{1}+b_{2}\equiv 0$ we
obtain%
\begin{equation*}
B(x_{2}\otimes gx_{2};GX_{2},gx_{2})GX_{2}\otimes g\otimes gx_{2}.
\end{equation*}%
By taking in account also the right side, we get%
\begin{equation*}
B(g\otimes gx_{2};GX_{2},g)+B(x_{2}\otimes
gx_{2};GX_{2},gx_{2})-B(x_{2}\otimes g;GX_{2},g)=0
\end{equation*}%
which holds in view of the form of the elements.

\section{$B\left( x_{2}\otimes gx_{1}x_{2}\right) $}

By using $\ref{simplgxx}$ and we get%
\begin{eqnarray}
B(x_{2}\otimes gx_{1}x_{2}) &=&(1_{A}\otimes g)B(gx_{2}\otimes
1_{H})(1_{A}\otimes gx_{1}x_{2})  \label{form x2otgx1x2} \\
&&-(1_{A}\otimes x_{2})B(gx_{2}\otimes 1_{H})(1_{A}\otimes x_{1})  \notag \\
&&+(1_{A}\otimes x_{1})B(gx_{2}\otimes 1_{H})(1_{A}\otimes x_{2})  \notag \\
&&+(1_{A}\otimes gx_{1}x_{2})B(gx_{2}\otimes 1_{H})(1_{A}\otimes g)  \notag
\\
&&+B(x_{1}x_{2}\otimes 1_{H})(1_{A}\otimes x_{2})  \notag \\
&&-(1_{A}\otimes gx_{2})B(x_{1}x_{2}\otimes 1_{H})(1_{A}\otimes g)  \notag
\end{eqnarray}%
Therefore we obtain%
\begin{gather*}
B(x_{2}\otimes gx_{1}x_{2})= \\
\left[ -2B(x_{1}x_{2}\otimes 1_{H};1_{A},gx_{1})-4B(x_{1}x_{2}\otimes
1_{H};X_{1},g)\right] 1_{A}\otimes gx_{1}x_{2} \\
-2B(x_{1}x_{2}\otimes 1_{H};G,g)G\otimes gx_{2}+ \\
-2B(x_{1}x_{2}\otimes 1_{H};X_{1},g)X_{1}\otimes gx_{2}+ \\
-2B(x_{1}x_{2}\otimes 1_{H};X_{2},g)X_{2}\otimes gx_{2}+ \\
-2B(x_{1}x_{2}\otimes 1_{H};X_{1},gx_{1}x_{2})X_{1}X_{2}\otimes gx_{1}x_{2}
\\
-2B(x_{1}x_{2}\otimes 1_{H};GX_{1},gx_{1})GX_{1}\otimes gx_{1}x_{2} \\
\left[ -4B(x_{1}x_{2}\otimes 1_{H};G,gx_{1}x_{2})-4B(x_{1}x_{2}\otimes
1_{H};GX_{1},gx_{2})+2B(x_{1}x_{2}\otimes 1_{H};GX_{2},gx_{1})\right]
GX_{2}\otimes gx_{1}x_{2}+ \\
+\left[ +2B(x_{1}x_{2}\otimes 1_{H};G,gx_{1}x_{2})-2B(x_{1}x_{2}\otimes
1_{H};GX_{2},gx_{1})+2B(x_{1}x_{2}\otimes 1_{H};GX_{1},gx_{2})\right]
GX_{1}X_{2}\otimes gx_{2}
\end{gather*}

We write the Casimir formula for $B\left( x_{2}\otimes gx_{1}x_{2}\right) $%
\begin{eqnarray*}
&&\sum_{a,b_{1},b_{2},d,e_{1},e_{2}=0}^{1}\sum_{l_{1}=0}^{b_{1}}%
\sum_{l_{2}=0}^{b_{2}}\sum_{u_{1}=0}^{e_{1}}\sum_{u_{2}=0}^{e_{2}}\left(
-1\right) ^{\alpha \left( x_{2};l_{1},l_{2},u_{1},u_{2}\right) } \\
&&B(g\otimes
gx_{1}x_{2};G^{a}X_{1}^{b_{1}}X_{2}^{b_{2}},g^{d}x_{1}^{e_{1}}x_{2}^{e_{2}})
\\
&&G^{a}X_{1}^{b_{1}-l_{1}}X_{2}^{b_{2}-l_{2}}\otimes
g^{d}x_{1}^{e_{1}-u_{1}}x_{2}^{e_{2}-u_{2}}\otimes
g^{a+b_{1}+b_{2}+l_{1}+l_{2}+d+e_{1}+e_{2}+u_{1}+u_{2}}x_{1}^{l_{1}+u_{1}}x_{2}^{l_{2}+u_{2}+1}+
\\
&&+\sum_{a,b_{1},b_{2},d,e_{1},e_{2}=0}^{1}\sum_{l_{1}=0}^{b_{1}}%
\sum_{l_{2}=0}^{b_{2}}\sum_{u_{1}=0}^{e_{1}}\sum_{u_{2}=0}^{e_{2}}\left(
-1\right) ^{\alpha \left( 1_{H};l_{1},l_{2},u_{1},u_{2}\right) } \\
&&B(x_{2}\otimes
gx_{1}x_{2};G^{a}X_{1}^{b_{1}}X_{2}^{b_{2}},g^{d}x_{1}^{e_{1}}x_{2}^{e_{2}})
\\
&&G^{a}X_{1}^{b_{1}-l_{1}}X_{2}^{b_{2}-l_{2}}\otimes
g^{d}x_{1}^{e_{1}-u_{1}}x_{2}^{e_{2}-u_{2}}\otimes
g^{a+b_{1}+b_{2}+l_{1}+l_{2}+d+e_{1}+e_{2}+u_{1}+u_{2}}x_{1}^{l_{1}+u_{1}}x_{2}^{l_{2}+u_{2}}
\\
&=&B^{A}(x_{2}\otimes gx_{1}x_{2})\otimes B^{H}(x_{2}\otimes
gx_{1}x_{2})\otimes g+\text{ } \\
&&B^{A}(x_{2}\otimes gx_{1})\otimes B^{H}(x_{2}\otimes gx_{1})\otimes x_{2}%
\text{ } \\
&&-B^{A}(x_{2}\otimes gx_{2})\otimes B^{H}(x_{2}\otimes gx_{2})\otimes x_{1}
\\
&&B^{A}(x_{2}\otimes g)\otimes B^{H}(x_{2}\otimes g)\otimes gx_{1}x_{2}
\end{eqnarray*}

\subsection{$B\left( x_{2}\otimes gx_{1}x_{2};1_{A},gx_{1}x_{2}\right) $}

We deduce that%
\begin{eqnarray*}
a &=&b_{1}=b_{2}=0 \\
d &=&e_{1}=e_{2}=1 \\
a+b_{1}+b_{2}+l_{1}+l_{2}+d+e_{1}+e_{2}+u_{1}+u_{2} &\equiv &u_{1}+u_{2}+1
\end{eqnarray*}%
and get%
\begin{eqnarray*}
&&\sum_{u_{1}=0}^{e_{1}}\sum_{u_{2}=0}^{e_{2}}\left( -1\right) ^{\alpha
\left( 1_{H};0,0,u_{1},u_{2}\right) }B(x_{2}\otimes
gx_{1}x_{2};1_{A},gx_{1}x_{2})1_{A}\otimes \\
&&\otimes gx_{1}^{1-u_{1}}x_{2}^{1-u_{2}}\otimes
g^{u_{1}+u_{2}+1}x_{1}^{u_{1}}x_{2}^{u_{2}} \\
&=&\left( -1\right) ^{\alpha \left( 1_{H};0,0,0,0\right) }B(x_{2}\otimes
gx_{1}x_{2};1_{A},gx_{1}x_{2})1_{A}\otimes gx_{1}x_{2}\otimes g+ \\
&&\left( -1\right) ^{\alpha \left( 1_{H};0,0,0,1\right) }B(x_{2}\otimes
gx_{1}x_{2};1_{A},gx_{1}x_{2})1_{A}\otimes gx_{1}\otimes x_{2}+\text{ } \\
&&\left( -1\right) ^{\alpha \left( 1_{H};0,0,1,0\right) }B(x_{2}\otimes
gx_{1}x_{2};1_{A},gx_{1}x_{2})1_{A}\otimes gx_{2}\otimes x_{1}+ \\
&&\left( -1\right) ^{\alpha \left( 1_{H};0,0,1,1\right) }B(x_{2}\otimes
gx_{1}x_{2};1_{A},gx_{1}x_{2})1_{A}\otimes g\otimes gx_{1}x_{2}+
\end{eqnarray*}

\subsubsection{Case $1_{A}\otimes gx_{1}\otimes x_{2}$}

First summand of the left side of the equality gives us%
\begin{eqnarray*}
l_{1} &=&u_{1}=0,l_{2}=u_{2}=0 \\
a+b_{1}+b_{2}+d+e_{1}+e_{2} &\equiv &0 \\
a &=&b_{1}=b_{2}=0 \\
d &=&1,e_{1}=1,e_{2}=0
\end{eqnarray*}%
Since $\alpha \left( x_{2};0,0,0,0\right) \equiv a+b_{1}+b_{2}\equiv 0,$ we
get%
\begin{equation*}
B(g\otimes gx_{1}x_{2};1_{A},gx_{1})1_{A}\otimes gx_{1}\otimes x_{2}
\end{equation*}%
Second summand of the left side of the equality gives us%
\begin{eqnarray*}
l_{1} &=&u_{1}=0,l_{2}+u_{2}=1 \\
a &=&b_{1}=0,b_{2}=l_{2} \\
d &=&1,e_{1}=1,e_{2}=u_{2}
\end{eqnarray*}%
Since $\alpha \left( 1_{H};0,0,0,1\right) \equiv 0$ and $\alpha \left(
1_{H};0,1,0,0\right) \equiv 0$ we obtain
\begin{equation*}
\left[ B(x_{2}\otimes gx_{1}x_{2};1_{A},gx_{1}x_{2})+B(x_{2}\otimes
gx_{1}x_{2};X_{2},gx_{1})\right] 1_{A}\otimes gx_{1}\otimes x_{2}
\end{equation*}

Considering also the right side we get%
\begin{equation*}
B(g\otimes gx_{1}x_{2};1_{A},gx_{1})-B(x_{2}\otimes
gx_{1};1_{A},gx_{1})+B(x_{2}\otimes
gx_{1}x_{2};1_{A},gx_{1}x_{2})+B(x_{2}\otimes gx_{1}x_{2};X_{2},gx_{1})=0.
\end{equation*}%
which holds in view of the form of the elements.

\subsubsection{Case $1_{A}\otimes gx_{2}\otimes x_{1}$}

First summand of the left side of the equality gives us nothing. Second
summand of the left side of the equality gives us%
\begin{eqnarray*}
l_{1}+u_{1} &=&1,l_{2}=u_{2}=0 \\
a &=&b_{2}=0,b_{1}=l_{1} \\
d &=&e_{2}=1,e_{1}=u_{1}
\end{eqnarray*}%
Since $\alpha \left( 1_{H};0,0,1,0\right) \equiv e_{2}+\left(
a+b_{1}+b_{2}\right) \equiv 1$ and $\alpha \left( 1_{H};1,0,0,0\right)
\equiv b_{2}=0$ we get
\begin{equation*}
\left[ -B(x_{2}\otimes gx_{1}x_{2};1_{A},gx_{1}x_{2})+B(x_{2}\otimes
gx_{1}x_{2};X_{1},gx_{2})\right] 1_{A}\otimes gx_{2}\otimes x_{1}
\end{equation*}%
By considering also the right side of the equality we obtain%
\begin{equation*}
B(x_{2}\otimes gx_{2};1_{A},gx_{2})+-B(x_{2}\otimes
gx_{1}x_{2};1_{A},gx_{1}x_{2})+B(x_{2}\otimes gx_{1}x_{2};X_{1},gx_{2})=0
\end{equation*}%
which holds in view of the form of the elements.

\subsubsection{Case $1_{A}\otimes g\otimes gx_{1}x_{2}$}

First summand of the left side of the equality gives us%
\begin{eqnarray*}
l_{1}+u_{1} &=&1,l_{2}=u_{2}=0 \\
a &=&b_{2}=0,b_{1}=l_{1} \\
d &=&1,e_{1}=u_{1},e_{2}=0
\end{eqnarray*}%
Since $\alpha \left( x_{2};0,0,1,0\right) \equiv e_{2}=0$ and $\alpha \left(
x_{2};1,0,0,0\right) \equiv a+b_{1}=1$ we obtain
\begin{equation*}
\left[ B(g\otimes gx_{1}x_{2};1_{A},gx_{1})-B(g\otimes gx_{1}x_{2};X_{1},g)%
\right] 1_{A}\otimes g\otimes gx_{1}x_{2}.
\end{equation*}%
Second summand of the left side of the equality gives us%
\begin{eqnarray*}
l_{1}+u_{1} &=&1,l_{2}+u_{2}=1 \\
a &=&0,b_{1}=l_{1},b_{2}=l_{2}, \\
d &=&1,e_{1}=u_{1},e_{2}=u_{2}.
\end{eqnarray*}%
Since
\begin{eqnarray*}
\alpha \left( 1_{H};1,0,0,1\right) &\equiv &a+b_{1}=1 \\
\alpha \left( 1_{H};1,1,0,0\right) &\equiv &1+b_{2}\equiv 0 \\
\alpha \left( 1_{H};0,0,1,1\right) &\equiv &1+e_{2}\equiv 0 \\
\alpha \left( 1_{H};0,1,1,0\right) &\equiv &e_{2}+\left(
a+b_{1}+b_{2}+1\right) \equiv 0
\end{eqnarray*}%
we obtain%
\begin{equation*}
\left[
\begin{array}{c}
-B(x_{2}\otimes gx_{1}x_{2};X_{1},gx_{2})+B(x_{2}\otimes
gx_{1}x_{2};X_{1}X_{2},g) \\
+B(x_{2}\otimes gx_{1}x_{2};X_{2},gx_{1})+B(x_{2}\otimes
gx_{1}x_{2};1_{A},gx_{1}x_{2})%
\end{array}%
\right] 1_{A}\otimes g\otimes gx_{1}x_{2}.
\end{equation*}

By considering also the right side of the equality we get%
\begin{gather*}
B(g\otimes gx_{1}x_{2};1_{A},gx_{1})-B(g\otimes
gx_{1}x_{2};X_{1},g)-B(x_{2}\otimes g;1_{A},g)-B(x_{2}\otimes
gx_{1}x_{2};X_{1},gx_{2}) \\
+B(x_{2}\otimes gx_{1}x_{2};X_{1}X_{2},g)+B(x_{2}\otimes
gx_{1}x_{2};X_{2},gx_{1})+B(x_{2}\otimes gx_{1}x_{2};1_{A},gx_{1}x_{2})=0
\end{gather*}%
which holds in view of the form of the elements.

\subsection{$B\left( x_{2}\otimes gx_{1}x_{2};G,gx_{2}\right) $}

We deduce that%
\begin{equation*}
a=1,b_{1}=b_{2}=0,d=1,e_{1}=0,e_{2}=1
\end{equation*}%
and we get

\begin{eqnarray*}
&&\sum_{u_{2}=0}^{1}\left( -1\right) ^{\alpha \left(
1_{H};0,0,0,u_{2}\right) }B(x_{2}\otimes gx_{1}x_{2};G,gx_{2})G\otimes
gx_{2}^{1-u_{2}}\otimes g^{u_{2}+1}x_{2}^{u_{2}} \\
&=&\left( -1\right) ^{\alpha \left( 1_{H};0,0,0,0\right) }B(x_{2}\otimes
gx_{1}x_{2};G,gx_{2})G\otimes gx_{2}\otimes g+\text{ } \\
&&\left( -1\right) ^{\alpha \left( 1_{H};0,0,0,1\right) }B(x_{2}\otimes
gx_{1}x_{2};G,gx_{2})G\otimes g\otimes x_{2}\text{ }
\end{eqnarray*}

\subsubsection{Case $G\otimes g\otimes x_{2}$}

First summand of the left side of the equality gives us%
\begin{eqnarray*}
l_{1} &=&u_{1}=0,l_{2}=u_{2}=0 \\
a &=&1,b_{1}=b_{2}=0
\end{eqnarray*}%
Since $\alpha \left( x_{2};0,0,0,0\right) \equiv a+b_{1}+b_{2}=1$ we get
\begin{equation*}
-B(g\otimes gx_{1}x_{2};G,g)G\otimes g\otimes x_{2}.
\end{equation*}%
Second summand of the left side of the equality gives us

\begin{eqnarray*}
l_{1} &=&u_{1}=0,l_{2}+u_{2}=1 \\
a &=&1,b_{1}=0,b_{2}=l_{2}, \\
d &=&1,e_{1}=0,e_{2}=u_{2}.
\end{eqnarray*}%
Since $\alpha \left( 1_{H};0,0,0,1\right) \equiv a+b_{1}+b_{2}\equiv 1$ and $%
\alpha \left( 1_{H};0,1,0,0\right) \equiv 0$ we get
\begin{equation*}
\left[ -B(x_{2}\otimes gx_{1}x_{2};G,gx_{2})+B(x_{2}\otimes
gx_{1}x_{2};GX_{2},g)\right] G\otimes g\otimes x_{2}
\end{equation*}

By considering also the right side of the equality we get

\begin{equation*}
-B(g\otimes gx_{1}x_{2};G,g)-B(x_{2}\otimes gx_{1};G,g)-B(x_{2}\otimes
gx_{1}x_{2};G,gx_{2})+B(x_{2}\otimes gx_{1}x_{2};GX_{2},g)=0
\end{equation*}%
which holds in view of the form of the elements.

\subsection{$B\left( x_{2}\otimes gx_{1}x_{2};X_{1},gx_{2}\right) $}

We deduce that%
\begin{equation*}
b_{1}=1,d=e_{2}=1
\end{equation*}%
and we get%
\begin{eqnarray*}
&&\sum_{l_{1}=0}^{1}\sum_{u_{2}=0}^{1}\left( -1\right) ^{\alpha \left(
1_{H};l_{1},0,0,u_{2}\right) }B(x_{2}\otimes
gx_{1}x_{2};X_{1},gx_{2})X_{1}^{1-l_{1}}\otimes \\
&&gx_{2}^{1-u_{2}}\otimes g^{l_{1}+u_{2}+1}x_{1}^{l_{1}}x_{2}^{u_{2}} \\
&=&\left( -1\right) ^{\alpha \left( 1_{H};0,0,0,0\right) }B(x_{2}\otimes
gx_{1}x_{2};X_{1},gx_{2})X_{1}\otimes gx_{2}\otimes g+ \\
&&\left( -1\right) ^{\alpha \left( 1_{H};0,0,0,1\right) }B(x_{2}\otimes
gx_{1}x_{2};X_{1},gx_{2})X_{1}\otimes g\otimes x_{2}+ \\
&&\left( -1\right) ^{\alpha \left( 1_{H};1,0,0,0\right) }B(x_{2}\otimes
gx_{1}x_{2};X_{1},gx_{2})1_{A}\otimes gx_{2}\otimes x_{1} \\
&&\left( -1\right) ^{\alpha \left( 1_{H};1,0,0,1\right) }B(x_{2}\otimes
gx_{1}x_{2};X_{1},gx_{2})1_{A}\otimes g\otimes gx_{1}x_{2}.
\end{eqnarray*}

\subsubsection{Case $X_{1}\otimes g\otimes x_{2}$}

First summand of the left side of the equality gives us%
\begin{eqnarray*}
l_{1} &=&u_{1}=0,l_{2}=u_{2}=0 \\
a &=&b_{2}=0,b_{1}=1, \\
d &=&1,e_{1}=e_{2}=0.
\end{eqnarray*}%
Since $\alpha \left( x_{2};0,0,0,0\right) \equiv a+b_{1}+b_{2}=1,$ we get
\begin{equation*}
-B(g\otimes gx_{1}x_{2};X_{1},g)X_{1}\otimes g\otimes x_{2}
\end{equation*}%
Second summand of the left side of the equality gives us%
\begin{eqnarray*}
l_{1} &=&u_{1}=0,l_{2}+u_{2}=1 \\
a &=&0,b_{1}=1,b_{2}=l_{2} \\
d &=&1,e_{1}=0,e_{2}=u_{2}
\end{eqnarray*}%
Since $\alpha \left( 1_{H};0,0,0,1\right) \equiv a+b_{1}+b_{2}\equiv 1$ and $%
\alpha \left( 1_{H};0,1,0,0\right) \equiv 0,$ we obtain
\begin{equation*}
\left[ -B(x_{2}\otimes gx_{1}x_{2};X_{1},gx_{2})+B(x_{2}\otimes
gx_{1}x_{2};X_{1}X_{2},g)\right] X_{1}\otimes g\otimes x_{2}
\end{equation*}

By considering also the right side of the equality we get

\begin{eqnarray*}
&&-B(g\otimes gx_{1}x_{2};X_{1},g)-B(x_{2}\otimes gx_{1};X_{1},g) \\
&&-B(x_{2}\otimes gx_{1}x_{2};X_{1},gx_{2})+B(x_{2}\otimes
gx_{1}x_{2};X_{1}X_{2},g).
\end{eqnarray*}%
which holds in view of the form of the elements.

\subsubsection{Case $1_{A}\otimes gx_{2}\otimes x_{1}$}

This case was already considered in subsection $B\left( x_{2}\otimes
gx_{1}x_{2};1_{A},gx_{1}x_{2}\right) .$

\subsubsection{Case $1_{A}\otimes g\otimes gx_{1}x_{2}$}

This case was already considered in subsection $B\left( x_{2}\otimes
gx_{1}x_{2};1_{A},gx_{1}x_{2}\right) .$

\subsection{$B\left( x_{2}\otimes gx_{1}x_{2};X_{2},gx_{2}\right) $}

We deduce that
\begin{equation*}
b_{2}=1,d=e_{2}=1
\end{equation*}%
and we get%
\begin{gather*}
\left( -1\right) ^{\alpha \left( 1_{H};0,0,0,0\right) }X_{2}\otimes
gx_{2}\otimes g+ \\
\left( -1\right) ^{\alpha \left( 1_{H};0,1,0,0\right) }1_{A}\otimes
gx_{2}\otimes x_{2}+ \\
\left( -1\right) ^{\alpha \left( 1_{H};0,0,0,1\right) }X_{2}\otimes g\otimes
x_{2} \\
\left( -1\right) ^{\alpha \left( 1_{H};0,1,0,1\right)
}X_{2}^{1-l_{2}}\otimes g\otimes g^{l_{2}}x_{2}^{1+1}=0
\end{gather*}

\subsubsection{Case $1_{A}\otimes gx_{2}\otimes x_{2}$}

First summand of the left side of the equality gives us%
\begin{eqnarray*}
l_{1} &=&u_{1}=0,l_{2}=u_{2}=0 \\
a &=&b_{1}=b_{2}=0 \\
d &=&e_{2}=1,e_{1}=0
\end{eqnarray*}%
Since $\alpha \left( x_{2};0,0,0,0\right) \equiv a+b_{1}+b_{2}\equiv 0,$ we
get%
\begin{equation*}
B(g\otimes gx_{1}x_{2};1_{A},gx_{2})1_{A}\otimes gx_{2}\otimes x_{2}
\end{equation*}%
Second summand of the left side of the equality gives us%
\begin{eqnarray*}
l_{1} &=&u_{1}=0,l_{2}+u_{2}=1 \\
a &=&b_{1}=0,b_{2}=l_{2} \\
d &=&1,e_{1}=0,e_{2}-u_{2}=1\Rightarrow e_{2}=1,u_{2}=0,b_{2}=l_{2}=1
\end{eqnarray*}%
Since $\alpha \left( 1_{H};0,1,0,0\right) \equiv 0$ we obtain
\begin{equation*}
B(x_{2}\otimes gx_{1}x_{2};X_{2},gx_{2})1_{A}\otimes gx_{2}\otimes x_{2}
\end{equation*}

Considering also the right side we get%
\begin{equation*}
B(g\otimes gx_{1}x_{2};1_{A},gx_{2})-B(x_{2}\otimes
gx_{1};1_{A},gx_{2})+B(x_{2}\otimes gx_{1}x_{2};X_{2},gx_{2})=0.
\end{equation*}%
which holds in view of the form of the elements.

\subsubsection{$X_{2}\otimes g\otimes x_{2}$}

First summand of the left side of the equality gives us%
\begin{eqnarray*}
l_{1} &=&u_{1}=0,l_{2}=u_{2}=0 \\
a &=&b_{1}=0,b_{2}=1 \\
d &=&1,e_{1}=e_{2}=0
\end{eqnarray*}%
Since $\alpha \left( x_{2};0,0,0,0\right) \equiv a+b_{1}+b_{2}\equiv 1,$ we
get%
\begin{equation*}
-B(g\otimes gx_{1}x_{2};X_{2},g)X_{2}\otimes g\otimes x_{2}
\end{equation*}%
Second summand of the left side of the equality gives us%
\begin{eqnarray*}
l_{1} &=&u_{1}=0,l_{2}+u_{2}=1 \\
a &=&b_{1}=0,b_{2}-l_{2}=1\Rightarrow b_{2}=1,l_{2}=0,u_{2}=1 \\
d &=&1,e_{1}=0,e_{2}=u_{2}=1
\end{eqnarray*}%
Since $\alpha \left( 1_{H};0,0,0,1\right) \equiv a+b_{1}+b_{2}\equiv 1$ we
obtain
\begin{equation*}
-B(x_{2}\otimes gx_{1}x_{2};X_{2},gx_{2})X_{2}\otimes g\otimes x_{2}
\end{equation*}

Considering also the right side we get%
\begin{equation*}
-B(g\otimes gx_{1}x_{2};X_{2},g)-B(x_{2}\otimes
gx_{1};X_{2},g)-B(x_{2}\otimes gx_{1}x_{2};X_{2},gx_{2})=0.
\end{equation*}%
which holds in view of the form of the elements.

\subsection{$B\left( x_{2}\otimes gx_{1}x_{2};X_{1}X_{2},gx_{1}x_{2}\right) $%
}

We deduce that%
\begin{equation*}
a=0,b_{1}=b_{2}=1,d=e_{1}=e_{2}=1
\end{equation*}%
\begin{equation*}
X_{1}^{1-l_{1}}X_{2}^{1-l_{2}}\otimes gx_{1}^{1-u_{1}}x_{2}^{1-u}\otimes
g^{1+l_{1}+l_{2}+u_{1}+u_{2}}x_{1}^{l_{1}+u_{1}}x_{2}^{l_{2}+u_{2}}
\end{equation*}%
and we get

\begin{gather*}
\left( -1\right) ^{\alpha \left( 1_{H};0,0,1,0\right) }B(x_{2}\otimes
gx_{1}x_{2};X_{1}X_{2},gx_{1}x_{2})X_{1}X_{2}\otimes gx_{2}\otimes x_{1}+ \\
\left( -1\right) ^{\alpha \left( 1_{H};0,1,1,0\right) }B(x_{2}\otimes
gx_{1}x_{2};X_{1}X_{2},gx_{1}x_{2})X_{1}\otimes gx_{2}\otimes gx_{1}x_{2} \\
+\left( -1\right) ^{\alpha \left( 1_{H};1,0,1,0\right) }B(x_{2}\otimes
gx_{1}x_{2};X_{1}X_{2},gx_{1}x_{2})X_{1}^{1-l_{1}}X_{2}^{1-l_{2}}\otimes
gx_{2}\otimes g^{+l_{1}+l_{2}}x_{1}^{1+1}x_{2}^{l_{2}}=0 \\
+\left( -1\right) ^{\alpha \left( 1_{H};1,1,1,0\right) }B(x_{2}\otimes
gx_{1}x_{2};X_{1}X_{2},gx_{1}x_{2})X_{1}^{1-l_{1}}X_{2}^{1-l_{2}}\otimes
gx_{2}\otimes g^{+l_{1}+l_{2}}x_{1}^{1+1}x_{2}^{l_{2}}=0 \\
+\left( -1\right) ^{\alpha \left( 1_{H};0,0,0,0\right) }B(x_{2}\otimes
gx_{1}x_{2};X_{1}X_{2},gx_{1}x_{2})X_{1}X_{2}\otimes gx_{1}x_{2}\otimes g+ \\
+\left( -1\right) ^{\alpha \left( 1_{H};0,1,0,0\right) }B(x_{2}\otimes
gx_{1}x_{2};X_{1}X_{2},gx_{1}x_{2})X_{1}\otimes gx_{1}x_{2}\otimes x_{2}+ \\
+\left( -1\right) ^{\alpha \left( 1_{H};1,0,0,0\right) }B(x_{2}\otimes
gx_{1}x_{2};X_{1}X_{2},gx_{1}x_{2})X_{2}\otimes gx_{1}x_{2}\otimes x_{1}+ \\
+\left( -1\right) ^{\alpha \left( 1_{H};1,1,0,0\right) }B(x_{2}\otimes
gx_{1}x_{2};X_{1}X_{2},gx_{1}x_{2})1_{A}\otimes gx_{1}x_{2}\otimes
gx_{1}x_{2}+ \\
+\left( -1\right) ^{\alpha \left( 1_{H};0,0,0,1\right) }B(x_{2}\otimes
gx_{1}x_{2};X_{1}X_{2},gx_{1}x_{2})X_{1}X_{2}\otimes gx_{1}\otimes x_{2}+ \\
+\left( -1\right) ^{\alpha \left( 1_{H};0,1,0,1\right) }B(x_{2}\otimes
gx_{1}x_{2};X_{1}X_{2},gx_{1}x_{2})X_{1}^{1-l_{1}}X_{2}^{1-l_{2}}\otimes
gx_{1}\otimes g^{1+l_{1}+l_{2}+1}x_{1}^{l_{1}}x_{2}^{1+1}=0 \\
+\left( -1\right) ^{\alpha \left( 1_{H};1,0,0,1\right) }B(x_{2}\otimes
gx_{1}x_{2};X_{1}X_{2},gx_{1}x_{2})X_{2}\otimes gx_{1}\otimes gx_{1}x_{2}+ \\
+\left( -1\right) ^{\alpha \left( 1_{H};1,1,0,1\right) }B(x_{2}\otimes
gx_{1}x_{2};X_{1}X_{2},gx_{1}x_{2})X_{1}^{1-l_{1}}X_{2}^{1-l_{2}}\otimes
gx_{1}\otimes g^{1+l_{1}+l_{2}+1}x_{1}^{l_{1}}x_{2}^{1+1}=0 \\
+\left( -1\right) ^{\alpha \left( 1_{H};0,0,1,1\right) }B(x_{2}\otimes
gx_{1}x_{2};X_{1}X_{2},gx_{1}x_{2})X_{1}X_{2}\otimes g\otimes gx_{1}x_{2}+ \\
+\left( -1\right) ^{\alpha \left( 1_{H};0,1,1,1\right) }B(x_{2}\otimes
gx_{1}x_{2};X_{1}X_{2},gx_{1}x_{2})X_{2}^{1-l_{2}}\otimes g\otimes
g^{1+l_{1}+l_{2}}x_{1}^{l_{1}+1}x_{2}^{1+1}=0 \\
+\left( -1\right) ^{\alpha \left( 1_{H};1,0,1,1\right) }B(x_{2}\otimes
gx_{1}x_{2};X_{1}X_{2},gx_{1}x_{2})X_{2}^{1-l_{2}}\otimes g\otimes
g^{1+l_{1}+l_{2}}x_{1}^{1+1}x_{2}^{l_{2}+1}=0 \\
+\left( -1\right) ^{\alpha \left( 1_{H};1,1,1,1\right) }B(x_{2}\otimes
gx_{1}x_{2};X_{1}X_{2},gx_{1}x_{2})X_{2}^{1-l_{2}}\otimes g\otimes
g^{1+l_{1}+l_{2}}x_{1}^{1+1}x_{2}^{1+1}=0
\end{gather*}

\subsubsection{Case $X_{1}X_{2}\otimes gx_{2}\otimes x_{1}$}

First summand of the left side of the equality gives us nothing. Second
summand of the left side of the equality gives us%
\begin{eqnarray*}
l_{1}+u_{1} &=&1,l_{2}=u_{2}=0 \\
a &=&0,b_{1}-l_{1}=1\Rightarrow b_{1}=1,l_{1}=0,u_{1}=1,b_{2}=1 \\
d &=&e_{2}=1,e_{1}=u_{1}=1
\end{eqnarray*}%
Since $\alpha \left( 1_{H};0,0,1,0\right) \equiv e_{2}+\left(
a+b_{1}+b_{2}\right) \equiv 1$ we get
\begin{equation*}
-B(x_{2}\otimes gx_{1}x_{2};X_{1}X_{2},gx_{1}x_{2})X_{1}X_{2}\otimes
gx_{2}\otimes x_{1}
\end{equation*}%
By considering also the right side of the equality we obtain%
\begin{equation*}
B(x_{2}\otimes gx_{2};X_{1}X_{2},gx_{2})-B(x_{2}\otimes
gx_{1}x_{2};X_{1}X_{2},gx_{1}x_{2})=0
\end{equation*}%
which holds in view of the form of the elements.

\subsubsection{Case $X_{1}\otimes gx_{2}\otimes gx_{1}x_{2}$}

First summand of the left side of the equality gives us%
\begin{eqnarray*}
l_{1}+u_{1} &=&1,l_{2}=u_{2}=0 \\
a &=&b_{2}=0,b_{1}-l_{1}=1\Rightarrow b_{1}=1,l_{1}=0,u_{1}=1 \\
d &=&1,e_{1}=u_{1}=1,e_{2}=1
\end{eqnarray*}%
Since $\alpha \left( x_{2};1,0,0,0\right) \equiv a+b_{1}=1$ we obtain
\begin{equation*}
-B(g\otimes gx_{1}x_{2};X_{1},gx_{1}x_{2})X_{1}\otimes gx_{2}\otimes
gx_{1}x_{2}.
\end{equation*}%
Second summand of the left side of the equality gives us%
\begin{eqnarray*}
l_{1}+u_{1} &=&1,l_{2}+u_{2}=1 \\
a &=&0,b_{1}-l_{1}=1\Rightarrow b_{1}=1,l_{1}=0,u_{1}=1,b_{2}=l_{2}, \\
d &=&1,e_{1}=u_{1}=1,e_{2}-u_{2}=1\Rightarrow e_{2}=1,u_{2}=0,b_{2}=l_{2}=1.
\end{eqnarray*}%
Since $\alpha \left( 1_{H};0,1,1,0\right) \equiv e_{2}+\left(
a+b_{1}+b_{2}+1\right) \equiv 0$ we obtain%
\begin{equation*}
B(x_{2}\otimes gx_{1}x_{2};X_{1}X_{2},x_{1}x_{2})X_{1}\otimes gx_{2}\otimes
gx_{1}x_{2}.
\end{equation*}

By considering also the right side of the equality we get%
\begin{equation*}
-B(x_{2}\otimes g;X_{1},gx_{2})-B(g\otimes
gx_{1}x_{2};X_{1},gx_{1}x_{2})+B(x_{2}\otimes
gx_{1}x_{2};X_{1}X_{2},x_{1}x_{2})=0
\end{equation*}%
which holds in view of the form of the elements.

\subsubsection{Case $X_{1}\otimes gx_{1}x_{2}\otimes x_{2}$}

First summand of the left side of the equality gives us%
\begin{eqnarray*}
l_{1} &=&u_{1}=0,l_{2}=u_{2}=0 \\
a &=&b_{2}=0,b_{1}=1
\end{eqnarray*}%
Since $\alpha \left( x_{2};0,0,0,0\right) \equiv a+b_{1}+b_{2}=1$ we get
\begin{equation*}
-B(g\otimes gx_{1}x_{2};X_{1},gx_{1}x_{2})X_{1}\otimes gx_{1}x_{2}\otimes
x_{2}.
\end{equation*}%
Second summand of the left side of the equality gives us

\begin{eqnarray*}
l_{1} &=&u_{1}=0,l_{2}+u_{2}=1 \\
a &=&0,b_{1}=1,b_{2}=l_{2}, \\
d &=&1,e_{1}=1,e_{2}-u_{2}=1\Rightarrow e_{2}=1,u_{2}=0,b_{2}=l_{2}=1.
\end{eqnarray*}%
Since and $\alpha \left( 1_{H};0,1,0,0\right) \equiv 0$ we get
\begin{equation*}
+B(x_{2}\otimes gx_{1}x_{2};X_{1}X_{2},gx_{1}x_{2})X_{1}\otimes
gx_{1}x_{2}\otimes x_{2}
\end{equation*}

By considering also the right side of the equality we get

\begin{equation*}
-B(x_{2}\otimes gx_{1};X_{1},gx_{1}x_{2})-B(g\otimes
gx_{1}x_{2};X_{1},gx_{1}x_{2})+B(x_{2}\otimes
gx_{1}x_{2};X_{1}X_{2},gx_{1}x_{2})=0
\end{equation*}%
which holds in view of the form of the elements.

\subsubsection{Case $X_{2}\otimes gx_{1}x_{2}\otimes x_{1}$}

First summand of the left side of the equality gives us nothing. Second
summand of the left side of the equality gives us%
\begin{eqnarray*}
l_{1}+u_{1} &=&1,l_{2}=u_{2}=0 \\
a &=&0,b_{1}-l_{1}=1\Rightarrow b_{1}=1,l_{1}=0,u_{1}=1,b_{2}=1 \\
d &=&e_{2}=1,e_{1}=u_{1}=1
\end{eqnarray*}%
Since $\alpha \left( 1_{H};0,0,1,0\right) \equiv e_{2}+\left(
a+b_{1}+b_{2}\right) \equiv 1$ we get
\begin{equation*}
-B(x_{2}\otimes gx_{1}x_{2};X_{1}X_{2},gx_{1}x_{2})X_{1}X_{2}\otimes
gx_{2}\otimes x_{1}
\end{equation*}%
By considering also the right side of the equality we obtain%
\begin{equation*}
B(x_{2}\otimes gx_{2};X_{1}X_{2},gx_{2})-B(x_{2}\otimes
gx_{1}x_{2};X_{1}X_{2},gx_{1}x_{2})=0
\end{equation*}%
which holds in view of the form of the elements.

\subsubsection{Case $1_{A}\otimes gx_{1}x_{2}\otimes gx_{1}x_{2}$}

First summand of the left side of the equality gives us%
\begin{eqnarray*}
l_{1}+u_{1} &=&1,l_{2}=u_{2}=0 \\
a &=&b_{2}=0,b_{1}=l_{1} \\
d &=&1,e_{1}-u_{1}=1\Rightarrow e_{1}=1,u_{1}=0,b_{1}=l_{1}=1,e_{2}=1
\end{eqnarray*}%
Since $\alpha \left( x_{2};1,0,0,0\right) \equiv a+b_{1}=1$ we obtain
\begin{equation*}
-B(g\otimes gx_{1}x_{2};X_{1},gx_{1}x_{2})1_{A}\otimes gx_{1}x_{2}\otimes
gx_{1}x_{2}.
\end{equation*}%
Second summand of the left side of the equality gives us%
\begin{eqnarray*}
l_{1}+u_{1} &=&1,l_{2}+u_{2}=1 \\
a &=&0,b_{1}=l_{1}=1,b_{2}=l_{2}, \\
d &=&1,e_{1}-u_{1}=1\Rightarrow e_{1}=1,u_{1}=0,b_{1}=l_{1}=1 \\
e_{2}-u_{2} &=&1\Rightarrow e_{2}=1,u_{2}=0,b_{2}=l_{2}=1.
\end{eqnarray*}%
Since $\alpha \left( 1_{H};1,1,1,1\right) \equiv e_{2}+\left(
a+b_{1}+b_{2}+1\right) \equiv 0$ we obtain%
\begin{equation*}
B(x_{2}\otimes gx_{1}x_{2};X_{1}X_{2},gx_{1}x_{2})1_{A}\otimes
gx_{1}x_{2}\otimes gx_{1}x_{2}.
\end{equation*}

By considering also the right side of the equality we get%
\begin{equation*}
-B(x_{2}\otimes g;1_{A},gx_{1}x_{2})-B(g\otimes
gx_{1}x_{2};X_{1},gx_{1}x_{2})+B(x_{2}\otimes
gx_{1}x_{2};X_{1}X_{2},gx_{1}x_{2})=0
\end{equation*}%
which holds in view of the form of the elements.

\subsubsection{Case $X_{1}X_{2}\otimes gx_{1}\otimes x_{2}$}

First summand of the left side of the equality gives us%
\begin{eqnarray*}
l_{1} &=&u_{1}=0,l_{2}=u_{2}=0 \\
a &=&0,b_{2}=b_{1}=1
\end{eqnarray*}%
Since $\alpha \left( x_{2};0,0,0,0\right) \equiv a+b_{1}+b_{2}=0$ we get
\begin{equation*}
B(g\otimes gx_{1}x_{2};X_{1}X_{2},gx_{1})X_{1}X_{2}\otimes gx_{1}\otimes
x_{2}.
\end{equation*}%
Second summand of the left side of the equality gives us

\begin{eqnarray*}
l_{1} &=&u_{1}=0,l_{2}+u_{2}=1 \\
a &=&0,b_{1}=1,b_{2}-l_{2}=1\Rightarrow b_{2}=1,l_{2}=0,u_{2}=1 \\
d &=&1,e_{1}=1,e_{2}=u_{2}=1.
\end{eqnarray*}%
Since and $\alpha \left( 1_{H};0,0,0,1\right) \equiv a+b_{1}+b_{2}\equiv 0$
we get
\begin{equation*}
+B(x_{2}\otimes gx_{1}x_{2};X_{1}X_{2},gx_{1}x_{2})X_{1}X_{2}\otimes
gx_{1}\otimes x_{2}
\end{equation*}

By considering also the right side of the equality we get

\begin{equation*}
-B(x_{2}\otimes gx_{1};X_{1}X_{2},gx_{1})+B(g\otimes
gx_{1}x_{2};X_{1}X_{2},gx_{1})+B(x_{2}\otimes
gx_{1}x_{2};X_{1}X_{2},gx_{1}x_{2})=0
\end{equation*}%
which holds in view of the form of the elements.

\subsubsection{Case $X_{2}\otimes gx_{1}\otimes gx_{1}x_{2}$}

First summand of the left side of the equality gives us%
\begin{eqnarray*}
l_{1}+u_{1} &=&1,l_{2}=u_{2}=0 \\
a &=&0,b_{1}=l_{1},b_{2}=1 \\
d &=&1,e_{1}-u_{1}=1\Rightarrow e_{1}=1,u_{1}=0,b_{1}=l_{1}=1,e_{2}=0
\end{eqnarray*}%
Since $\alpha \left( x_{2};1,0,0,0\right) \equiv a+b_{1}=1$ we obtain
\begin{equation*}
-B(g\otimes gx_{1}x_{2};X_{1}X_{2},gx_{1})X_{2}\otimes gx_{1}\otimes
gx_{1}x_{2}.
\end{equation*}%
Second summand of the left side of the equality gives us%
\begin{eqnarray*}
l_{1}+u_{1} &=&1,l_{2}+u_{2}=1 \\
a &=&0,b_{1}=l_{1},b_{2}-l_{2}=1\Rightarrow b_{2}=1,l_{2}=0,u_{2}=1, \\
d &=&1,e_{1}-u_{1}=1\Rightarrow e_{1}=1,u_{1}=0,b_{1}=l_{1}=1,e_{2}=u_{2}=1
\end{eqnarray*}%
Since $\alpha \left( 1_{H};1,0,0,1\right) \equiv a+b_{1}\equiv 1$ we obtain%
\begin{equation*}
-B(x_{2}\otimes gx_{1}x_{2};X_{1}X_{2},gx_{1}x_{2})X_{2}\otimes
gx_{1}\otimes gx_{1}x_{2}.
\end{equation*}

By considering also the right side of the equality we get%
\begin{equation*}
-B(x_{2}\otimes g;X_{2},gx_{1})-B(g\otimes
gx_{1}x_{2};X_{1}X_{2},gx_{1})-B(x_{2}\otimes
gx_{1}x_{2};X_{1}X_{2},gx_{1}x_{2})=0
\end{equation*}%
which holds in view of the form of the elements.

\subsubsection{Case $X_{1}X_{2}\otimes g\otimes gx_{1}x_{2}$}

First summand of the left side of the equality gives us%
\begin{eqnarray*}
l_{1}+u_{1} &=&1,l_{2}=u_{2}=0 \\
a &=&0,b_{1}-l_{1}=1\Rightarrow b_{1}=1,l_{1}=0,u_{1}=1,b_{2}=1 \\
d &=&1,e_{1}=u_{1}=1,e_{2}=0
\end{eqnarray*}%
Since $\alpha \left( x_{2};0,0,1,0\right) \equiv e_{2}=0$ we obtain
\begin{equation*}
B(g\otimes gx_{1}x_{2};X_{1}X_{2},gx_{1})X_{1}X_{2}\otimes g\otimes
gx_{1}x_{2}.
\end{equation*}%
Second summand of the left side of the equality gives us%
\begin{eqnarray*}
l_{1}+u_{1} &=&1,l_{2}+u_{2}=1 \\
a &=&0,b_{1}-l_{1}=1\Rightarrow b_{1}=1,l_{1}=0,u_{1}=1, \\
b_{2}-l_{2} &=&1\Rightarrow b_{2}=1,l_{2}=0,u_{2}=1, \\
d &=&1,e_{1}=u_{1}=1,e_{2}=u_{2}=1
\end{eqnarray*}%
Since $\alpha \left( 1_{H};0,0,1,1\right) \equiv 1+e_{2}\equiv 0$ we obtain%
\begin{equation*}
+B(x_{2}\otimes gx_{1}x_{2};X_{1}X_{2},gx_{1}x_{2})X_{1}X_{2}\otimes
g\otimes gx_{1}x_{2}.
\end{equation*}

By considering also the right side of the equality we get%
\begin{equation*}
-B(x_{2}\otimes g;X_{1}X_{2},g)+B(g\otimes
gx_{1}x_{2};X_{1}X_{2},gx_{1})+B(x_{2}\otimes
gx_{1}x_{2};X_{1}X_{2},gx_{1}x_{2})=0
\end{equation*}%
which holds in view of the form of the elements.

\subsection{$B\left( x_{2}\otimes gx_{1}x_{2};GX_{1},gx_{1}x_{2}\right) $}

We deduce that%
\begin{eqnarray*}
a &=&1,b_{1}=1,b_{2}=0 \\
d &=&e_{1}=e_{2}=1
\end{eqnarray*}%
and we get%
\begin{gather*}
\left( -1\right) ^{\alpha \left( 1_{H};0,0,0,0\right) }B(x_{2}\otimes
gx_{1}x_{2};GX_{1},gx_{1}x_{2})GX_{1}\otimes gx_{1}x_{2}\otimes g \\
\left( -1\right) ^{\alpha \left( 1_{H};1,0,0,0\right) }B(x_{2}\otimes
gx_{1}x_{2};GX_{1},gx_{1}x_{2})G\otimes gx_{1}x_{2}\otimes x_{1} \\
\left( -1\right) ^{\alpha \left( 1_{H};0,0,1,0\right) }B(x_{2}\otimes
gx_{1}x_{2};GX_{1},gx_{1}x_{2})GX_{1}\otimes gx_{2}\otimes x_{1} \\
\left( -1\right) ^{\alpha \left( 1_{H};1,0,1,0\right) }B(x_{2}\otimes
gx_{1}x_{2};GX_{1},gx_{1}x_{2})GX_{1}^{1-l_{1}}\otimes gx_{2}\otimes
g^{+l_{11}}x_{1}^{1+1}=0 \\
\left( -1\right) ^{\alpha \left( 1_{H};0,0,0,1\right) }B(x_{2}\otimes
gx_{1}x_{2};GX_{1},gx_{1}x_{2})GX_{1}\otimes gx_{1}\otimes x_{2} \\
\left( -1\right) ^{\alpha \left( 1_{H};1,0,0,1\right) }B(x_{2}\otimes
gx_{1}x_{2};GX_{1},gx_{1}x_{2})G\otimes gx_{1}\otimes gx_{1}x_{2} \\
\left( -1\right) ^{\alpha \left( 1_{H};0,0,1,1\right) }B(x_{2}\otimes
gx_{1}x_{2};GX_{1},gx_{1}x_{2})GX_{1}\otimes g\otimes gx_{1}x_{2} \\
\left( -1\right) ^{\alpha \left( 1_{H};1,0,1,1\right) }B(x_{2}\otimes
gx_{1}x_{2};GX_{1},gx_{1}x_{2})GX_{1}^{1-l_{1}}\otimes g\otimes
g^{+l_{1}+1}x_{1}^{1+1}x_{2}=0
\end{gather*}

\subsubsection{Case $G\otimes gx_{1}x_{2}\otimes x_{1}$}

First summand of the left side of the equality gives us nothing. Second
summand of the left side of the equality gives us%
\begin{eqnarray*}
l_{1}+u_{1} &=&1,l_{2}=u_{2}=0 \\
a &=&1,b_{1}=l_{1},b_{2}=0 \\
d &=&e_{2}=1,e_{1}-u_{1}=1\Rightarrow e_{1}=1,u_{1}=0,b_{1}=l_{1}=1
\end{eqnarray*}%
Since $\alpha \left( 1_{H};1,0,0,0\right) \equiv b_{2}\equiv 0$ we get
\begin{equation*}
B(x_{2}\otimes gx_{1}x_{2};GX_{1},gx_{1}x_{2})G\otimes gx_{1}x_{2}\otimes
x_{1}
\end{equation*}%
By considering also the right side of the equality we obtain%
\begin{equation*}
B(x_{2}\otimes gx_{2};G,gx_{1}x_{2})+B(x_{2}\otimes
gx_{1}x_{2};GX_{1},gx_{1}x_{2})=0
\end{equation*}%
which holds in view of the form of the elements.

\subsubsection{Case $GX_{1}\otimes gx_{2}\otimes x_{1}$}

First summand of the left side of the equality gives us nothing. Second
summand of the left side of the equality gives us%
\begin{eqnarray*}
l_{1}+u_{1} &=&1,l_{2}=u_{2}=0 \\
a &=&1,b_{1}-l_{1}=1\Rightarrow b_{1}=1,l_{1}=0,u_{1}=1,b_{2}=0 \\
d &=&e_{2}=1,e_{1}=u_{1}=1
\end{eqnarray*}%
Since $\alpha \left( 1_{H};0,0,1,0\right) \equiv e_{2}\equiv 1$ we get
\begin{equation*}
-B(x_{2}\otimes gx_{1}x_{2};GX_{1},gx_{1}x_{2})GX_{1}\otimes gx_{2}\otimes
x_{1}
\end{equation*}%
By considering also the right side of the equality we obtain%
\begin{equation*}
B(x_{2}\otimes gx_{2};GX_{1},gx_{2})-B(x_{2}\otimes
gx_{1}x_{2};GX_{1},gx_{1}x_{2})=0
\end{equation*}%
which holds in view of the form of the elements.

\subsubsection{Case $GX_{1}\otimes gx_{1}\otimes x_{2}$}

First summand of the left side of the equality gives us%
\begin{eqnarray*}
l_{1} &=&u_{1}=0,l_{2}=u_{2}=0 \\
a &=&b_{1}=1,b_{2}=0 \\
d &=&e_{1}=1,e_{2}=0
\end{eqnarray*}%
Since $\alpha \left( x_{2};0,0,0,0\right) \equiv a+b_{1}+b_{2}\equiv 0$ we
get
\begin{equation*}
B(g\otimes gx_{1}x_{2};GX_{1},gx_{1})GX_{1}\otimes gx_{1}\otimes x_{2}.
\end{equation*}%
Second summand of the left side of the equality gives us

\begin{eqnarray*}
l_{1} &=&u_{1}=0,l_{2}+u_{2}=1 \\
a &=&1,b_{1}=1,b_{2}=l_{2} \\
d &=&1,e_{1}=1,e_{2}=u_{2}.
\end{eqnarray*}%
Since and $\alpha \left( 1_{H};0,0,0,1\right) \equiv a+b_{1}+b_{2}\equiv 0$
and $\alpha \left( 1_{H};0,1,0,0\right) \equiv 0$ we get
\begin{equation*}
\left[ +B(x_{2}\otimes gx_{1}x_{2};GX_{1},gx_{1}x_{2})+B(x_{2}\otimes
gx_{1}x_{2};GX_{1}X_{2},gx_{1})\right] GX_{1}\otimes gx_{1}\otimes x_{2}.
\end{equation*}

By considering also the right side of the equality we get

\begin{gather*}
-B(x_{2}\otimes gx_{1};GX_{1},gx_{1})+B(g\otimes gx_{1}x_{2};GX_{1},gx_{1})
\\
+B(x_{2}\otimes gx_{1}x_{2};GX_{1},gx_{1}x_{2})+B(x_{2}\otimes
gx_{1}x_{2};GX_{1}X_{2},gx_{1})=0
\end{gather*}%
which holds in view of the form of the elements.

\subsubsection{Case $G\otimes gx_{1}\otimes gx_{1}x_{2}$}

First summand of the left side of the equality gives us%
\begin{eqnarray*}
l_{1}+u_{1} &=&1,l_{2}=u_{2}=0 \\
a &=&1,b_{1}=l_{1},b_{2}=0 \\
d &=&1,e_{1}-u_{1}=1\Rightarrow e_{1}=1,u_{1}=0,b_{1}=l_{1}=1,e_{2}=0
\end{eqnarray*}%
Since $\alpha \left( x_{2};1,0,0,0\right) \equiv a+b_{1}\equiv 0$ we obtain
\begin{equation*}
B(g\otimes gx_{1}x_{2};GX_{1},gx_{1})G\otimes gx_{1}\otimes gx_{1}x_{2}.
\end{equation*}%
Second summand of the left side of the equality gives us%
\begin{eqnarray*}
l_{1}+u_{1} &=&1,l_{2}+u_{2}=1 \\
a &=&1,b_{1}=l_{1},b_{2}=l_{2} \\
d &=&1,e_{1}-u_{1}=1\Rightarrow e_{1}=1,u_{1}=0,b_{1}=l_{1}=1,e_{2}=u_{2}
\end{eqnarray*}%
Since $\alpha \left( 1_{H};1,0,0,1\right) \equiv a+b_{1}\equiv 0$ and $%
\alpha \left( 1_{H};1,1,0,0\right) \equiv 1+b_{2}\equiv 0$ we obtain%
\begin{equation*}
\left[ B(x_{2}\otimes gx_{1}x_{2};GX_{1},gx_{1}x_{2})+B(x_{2}\otimes
gx_{1}x_{2};GX_{1}X_{2},gx_{1})\right] G\otimes gx_{1}\otimes gx_{1}x_{2}.
\end{equation*}

By considering also the right side of the equality we get%
\begin{gather*}
-B(x_{2}\otimes g;G,gx_{1})+B(g\otimes gx_{1}x_{2};GX_{1},gx_{1}) \\
+B(x_{2}\otimes gx_{1}x_{2};GX_{1},gx_{1}x_{2})+B(x_{2}\otimes
gx_{1}x_{2};GX_{1}X_{2},gx_{1})=0
\end{gather*}%
which holds in view of the form of the elements.

\subsubsection{Case $GX_{1}\otimes g\otimes gx_{1}x_{2}$}

First summand of the left side of the equality gives us%
\begin{eqnarray*}
l_{1}+u_{1} &=&1,l_{2}=u_{2}=0 \\
a &=&1,b_{1}-l_{1}=1\Rightarrow b_{1}=1,l_{1}=0,u_{1}=1,b_{2}=0 \\
d &=&1,e_{1}=u_{1}=1,e_{2}=0
\end{eqnarray*}%
Since $\alpha \left( x_{2};1,0,0,0\right) \equiv a+b_{1}\equiv 0$ we obtain
\begin{equation*}
B(g\otimes gx_{1}x_{2};GX_{1},gx_{1})GX_{1}\otimes g\otimes gx_{1}x_{2}.
\end{equation*}%
Second summand of the left side of the equality gives us%
\begin{eqnarray*}
l_{1}+u_{1} &=&1,l_{2}+u_{2}=1 \\
a &=&1,b_{1}-l_{1}=1\Rightarrow b_{1}=1,l_{1}=0,u_{1}=1,b_{2}=l_{2} \\
d &=&1,e_{1}=u_{1}=1,e_{2}=u_{2}
\end{eqnarray*}%
Since $\alpha \left( 1_{H};0,0,1,1\right) \equiv 1+e_{2}\equiv 0$ and $%
\alpha \left( 1_{H};0,1,1,0\right) \equiv e_{2}+a+b_{1}+b_{2}+1\equiv 0$ we
obtain%
\begin{equation*}
\left[ B(x_{2}\otimes gx_{1}x_{2};GX_{1},gx_{1}x_{2})+B(x_{2}\otimes
gx_{1}x_{2};GX_{1}X_{2},gx_{1})\right] GX_{1}\otimes g\otimes gx_{1}x_{2}.
\end{equation*}

By considering also the right side of the equality we get%
\begin{gather*}
-B(x_{2}\otimes g;GX_{1},g)+B(g\otimes gx_{1}x_{2};GX_{1},gx_{1}) \\
+B(x_{2}\otimes gx_{1}x_{2};GX_{1},gx_{1}x_{2})+B(x_{2}\otimes
gx_{1}x_{2};GX_{1}X_{2},gx_{1})=0
\end{gather*}%
which we already got above.

\subsection{$B\left( x_{2}\otimes gx_{1}x_{2};GX_{2},gx_{1}x_{2}\right) $}

We have%
\begin{eqnarray*}
a &=&1,b_{1}=0,b_{2}=1 \\
d &=&e_{1}=e_{2}=1
\end{eqnarray*}%
and we get%
\begin{gather*}
\left( -1\right) ^{\alpha \left( 1_{H};0,0,0,0\right) }B\left( x_{2}\otimes
gx_{1}x_{2};GX_{2},gx_{1}x_{2}\,\right) GX_{2}\otimes gx_{1}x_{2}\otimes g+
\\
\left( -1\right) ^{\alpha \left( 1_{H};0,1,0,0\right) }B\left( x_{2}\otimes
gx_{1}x_{2};GX_{2},gx_{1}x_{2}\,\right) G\otimes gx_{1}x_{2}\otimes x_{2}+ \\
\left( -1\right) ^{\alpha \left( 1_{H};0,0,1,0\right) }B\left( x_{2}\otimes
gx_{1}x_{2};GX_{2},gx_{1}x_{2}\,\right) GX_{2}\otimes gx_{2}\otimes x_{1}+ \\
\left( -1\right) ^{\alpha \left( 1_{H};0,1,1,0\right) }B\left( x_{2}\otimes
gx_{1}x_{2};GX_{2},gx_{1}x_{2}\,\right) G\otimes gx_{2}\otimes gx_{1}x_{2}+
\\
\left( -1\right) ^{\alpha \left( 1_{H};0,0,0,1\right) }B\left( x_{2}\otimes
gx_{1}x_{2};GX_{2},gx_{1}x_{2}\,\right) GX_{2}\otimes gx_{1}\otimes x_{2}+ \\
\left( -1\right) ^{\alpha \left( 1_{H};0,1,0,1\right) }B\left( x_{2}\otimes
gx_{1}x_{2};GX_{2},gx_{1}x_{2}\,\right) G\otimes gx_{1}\otimes
g^{l_{2}}x_{2}^{1+1}=0 \\
\left( -1\right) ^{\alpha \left( 1_{H};0,0,1,1\right) }B\left( x_{2}\otimes
gx_{1}x_{2};GX_{2},gx_{1}x_{2}\,\right) GX_{2}\otimes g\otimes gx_{1}x_{2} \\
\left( -1\right) ^{\alpha \left( 1_{H};0,1,1,1\right) }B\left( x_{2}\otimes
gx_{1}x_{2};GX_{2},gx_{1}x_{2}\,\right) GX_{2}^{1-l_{2}}\otimes g\otimes
g^{l_{2}+1}x_{1}x_{2}^{1+1}=0
\end{gather*}

\subsubsection{Case $G\otimes gx_{1}x_{2}\otimes x_{2}$}

First summand of the left side of the equality gives us%
\begin{eqnarray*}
l_{1} &=&u_{1}=0,l_{2}=u_{2}=0 \\
a &=&1,b_{1}=b_{2}=0 \\
d &=&e_{1}=e_{2}=1
\end{eqnarray*}%
Since $\alpha \left( x_{2};0,0,0,0\right) \equiv a+b_{1}+b_{2}\equiv 1$ we
get
\begin{equation*}
-B(g\otimes gx_{1}x_{2};G,gx_{1}x_{2})G\otimes gx_{1}x_{2}\otimes x_{2}.
\end{equation*}%
Second summand of the left side of the equality gives us

\begin{eqnarray*}
l_{1} &=&u_{1}=0,l_{2}+u_{2}=1 \\
a &=&1,b_{1}=0,b_{2}=l_{2} \\
d &=&1,e_{1}=1,e_{2}-u_{2}=1\Rightarrow e_{2}=1,u_{2}=0,b_{2}=l_{2}=1.
\end{eqnarray*}%
Since $\alpha \left( 1_{H};0,1,0,0\right) \equiv 0$ we get
\begin{equation*}
+B(x_{2}\otimes gx_{1}x_{2};GX_{2},gx_{1}x_{2})G\otimes gx_{1}x_{2}\otimes
x_{2}.
\end{equation*}

By considering also the right side of the equality we get

\begin{equation*}
-B(x_{2}\otimes gx_{1};G,gx_{1}x_{2})-B(g\otimes
gx_{1}x_{2};G,gx_{1}x_{2})+B(x_{2}\otimes gx_{1}x_{2};GX_{2},gx_{1}x_{2})=0
\end{equation*}%
which holds in view of the form of the elements.

\subsubsection{Case $GX_{2}\otimes gx_{2}\otimes x_{1}$}

First summand of the left side of the equality gives us nothing. Second
summand of the left side of the equality gives us%
\begin{eqnarray*}
l_{1}+u_{1} &=&1,l_{2}=u_{2}=0 \\
a &=&1,b_{1}=l_{1},b_{2}=1 \\
d &=&e_{2}=1,e_{1}=u_{1}
\end{eqnarray*}%
Since $\alpha \left( 1_{H};0,0,1,0\right) \equiv e_{2}+\left(
a+b_{1}+b_{2}\right) \equiv 1$ and $\alpha \left( 1_{H};1,0,0,0\right)
\equiv b_{2}=1$ we get
\begin{equation*}
\left[ -B(x_{2}\otimes gx_{1}x_{2};GX_{2},gx_{1}x_{2})-B(x_{2}\otimes
gx_{1}x_{2};GX_{1}X_{2},gx_{2})\right] GX_{2}\otimes gx_{2}\otimes x_{1}
\end{equation*}%
By considering also the right side of the equality we obtain%
\begin{equation*}
B(x_{2}\otimes gx_{2};GX_{2},gx_{2})-B(x_{2}\otimes
gx_{1}x_{2};GX_{2},gx_{1}x_{2})-B(x_{2}\otimes
gx_{1}x_{2};GX_{1}X_{2},gx_{2})=0
\end{equation*}%
which holds in view of the form of the elements.

\subsubsection{Case $G\otimes gx_{2}\otimes gx_{1}x_{2}$}

First summand of the left side of the equality gives us%
\begin{eqnarray*}
l_{1}+u_{1} &=&1,l_{2}=u_{2}=0 \\
a &=&1,b_{1}=l_{1},b_{2}=0 \\
d &=&1,e_{1}=u_{1},e_{2}=1
\end{eqnarray*}%
Since $\alpha \left( x_{2};0,0,1,0\right) \equiv e_{2}=1$ and $\alpha \left(
x_{2};1,0,0,0\right) \equiv a+b_{1}\equiv 0$ we obtain
\begin{equation*}
\left[ -B(g\otimes gx_{1}x_{2};G,gx_{1}x_{2})+B(g\otimes
gx_{1}x_{2};GX_{1},gx_{2})\right] G\otimes gx_{2}\otimes gx_{1}x_{2}.
\end{equation*}%
Second summand of the left side of the equality gives us%
\begin{eqnarray*}
l_{1}+u_{1} &=&1,l_{2}+u_{2}=1 \\
a &=&1,b_{1}=l_{1},b_{2}=l_{2} \\
d &=&1,e_{1}=u_{1},e_{2}-u_{2}=1\Rightarrow e_{2}=1,u_{2}=0,b_{2}=l_{2}=1
\end{eqnarray*}%
Since $\alpha \left( 1_{H};0,1,1,1\right) \equiv 1+e_{2}\equiv 0$ and $%
\alpha \left( 1_{H};1,1,0,1\right) \equiv a+b_{1}\equiv 0$ we obtain%
\begin{equation*}
\left[ B(x_{2}\otimes gx_{1}x_{2};GX_{2},gx_{1}x_{2})+B(x_{2}\otimes
gx_{1}x_{2};GX_{1}X_{2},gx_{2})\right] G\otimes gx_{2}\otimes gx_{1}x_{2}.
\end{equation*}

By considering also the right side of the equality we get%
\begin{gather*}
-B(x_{2}\otimes g;G,gx_{2})-B(g\otimes gx_{1}x_{2};G,gx_{1}x_{2})+B(g\otimes
gx_{1}x_{2};GX_{1},gx_{2}) \\
+B(x_{2}\otimes gx_{1}x_{2};GX_{2},gx_{1}x_{2})+B(x_{2}\otimes
gx_{1}x_{2};GX_{1}X_{2},gx_{2})=0
\end{gather*}%
which holds in view of the form of the elements.

\subsubsection{Case $GX_{2}\otimes gx_{1}\otimes x_{2}$}

First summand of the left side of the equality gives us%
\begin{eqnarray*}
l_{1} &=&u_{1}=0,l_{2}=u_{2}=0 \\
a &=&b_{2}=1,b_{1}=0 \\
d &=&e_{1}=1,e_{2}=0
\end{eqnarray*}%
Since $\alpha \left( x_{2};0,0,0,0\right) \equiv a+b_{1}+b_{2}\equiv 0$ we
get
\begin{equation*}
B(g\otimes gx_{1}x_{2};GX_{2},gx_{1})GX_{2}\otimes gx_{1}\otimes x_{2}.
\end{equation*}%
Second summand of the left side of the equality gives us

\begin{eqnarray*}
l_{1} &=&u_{1}=0,l_{2}+u_{2}=1 \\
a &=&1,b_{1}=0,b_{2}-l_{2}=1\Rightarrow b_{2}=1,l_{2}=0,u_{2}=1 \\
d &=&1,e_{1}=1,e_{2}=u_{2}=1.
\end{eqnarray*}%
Since $\alpha \left( 1_{H};0,0,0,1\right) \equiv a+b_{1}+b_{2}\equiv 0$ we
get
\begin{equation*}
+B(x_{2}\otimes gx_{1}x_{2};GX_{2},gx_{1}x_{2})GX_{2}\otimes gx_{1}\otimes
x_{2}.
\end{equation*}

By considering also the right side of the equality we get

\begin{equation*}
-B(x_{2}\otimes gx_{1};GX_{2},gx_{1})+B(g\otimes
gx_{1}x_{2};GX_{2},gx_{1})+B(x_{2}\otimes gx_{1}x_{2};GX_{2},gx_{1}x_{2})=0
\end{equation*}%
which holds in view of the form of the elements.

\subsubsection{Case $GX_{2}\otimes g\otimes gx_{1}x_{2}$}

First summand of the left side of the equality gives us%
\begin{eqnarray*}
l_{1}+u_{1} &=&1,l_{2}=u_{2}=0 \\
a &=&1,b_{1}=l_{1},b_{2}=1 \\
d &=&1,e_{1}=u_{1},e_{2}=0
\end{eqnarray*}%
Since $\alpha \left( x_{2};0,0,1,0\right) \equiv e_{2}=0$ and $\alpha \left(
x_{2};1,0,0,0\right) \equiv a+b_{1}\equiv 0$ we obtain
\begin{equation*}
\left[ +B(g\otimes gx_{1}x_{2};GX_{2},gx_{1})+B(g\otimes
gx_{1}x_{2};GX_{1}X_{2},g)\right] GX_{2}\otimes g\otimes gx_{1}x_{2}.
\end{equation*}%
Second summand of the left side of the equality gives us%
\begin{eqnarray*}
l_{1}+u_{1} &=&1,l_{2}+u_{2}=1 \\
a &=&1,b_{1}=l_{1},b_{2}-l_{2}=1\Rightarrow b_{2}=1,l_{2}=0,u_{2}=1 \\
d &=&1,e_{1}=u_{1},e_{2}=u_{2}=1
\end{eqnarray*}%
Since $\alpha \left( 1_{H};0,1,1,1\right) \equiv 1+e_{2}\equiv 0$ and $%
\alpha \left( 1_{H};1,1,0,1\right) \equiv a+b_{1}\equiv 0$ we obtain%
\begin{equation*}
\left[ B(x_{2}\otimes gx_{1}x_{2};GX_{2},gx_{1}x_{2})+B(x_{2}\otimes
gx_{1}x_{2};GX_{1}X_{2},gx_{2})\right] GX_{2}\otimes g\otimes gx_{1}x_{2}.
\end{equation*}

By considering also the right side of the equality we get%
\begin{gather*}
-B(x_{2}\otimes g;GX_{2},g)+B(g\otimes gx_{1}x_{2};GX_{2},gx_{1})+B(g\otimes
gx_{1}x_{2};GX_{1}X_{2},g)+ \\
+B(x_{2}\otimes gx_{1}x_{2};GX_{2},gx_{1}x_{2})+B(x_{2}\otimes
gx_{1}x_{2};GX_{1}X_{2},gx_{2})=0
\end{gather*}%
which holds in view of the form of the elements.

\subsection{$B\left( x_{2}\otimes gx_{1}x_{2};GX_{1}X_{2},gx_{2}\right) $}

We get%
\begin{eqnarray*}
a &=&b_{1}=b_{2}=1 \\
d &=&e_{2}=1,e_{1}=0
\end{eqnarray*}%
so that we obtain%
\begin{equation*}
GX_{1}^{1_{1}-l_{1}}X_{2}^{1_{2}-l_{2}}\otimes gx_{2}^{1-u_{2}}\otimes
g^{1+l_{1}+l_{2}+u_{2}}x_{1}^{l_{1}}x_{2}^{l_{2}+u_{2}}
\end{equation*}%
\begin{gather*}
\left( -1\right) ^{\alpha \left( 1_{H};0,0,0,0\right) }B(x_{2}\otimes
gx_{1}x_{2};GX_{1}X_{2},gx_{1})GX_{1}X_{2}\otimes gx_{2}\otimes g+ \\
\left( -1\right) ^{\alpha \left( 1_{H};1,0,0,0\right) }B(x_{2}\otimes
gx_{1}x_{2};GX_{1}X_{2},gx_{1})GX_{2}\otimes gx_{2}\otimes x_{1}+ \\
\left( -1\right) ^{\alpha \left( 1_{H};0,1,0,0\right) }B(x_{2}\otimes
gx_{1}x_{2};GX_{1}X_{2},gx_{1})GX_{1}\otimes gx_{2}\otimes x_{2}+ \\
\left( -1\right) ^{\alpha \left( 1_{H};1,1,0,0\right) }B(x_{2}\otimes
gx_{1}x_{2};GX_{1}X_{2},gx_{1})G\otimes gx_{1}\otimes gx_{1}x_{2}+vB\left(
x_{2}\otimes gx_{1}x_{2};GX_{1},gx_{1}x_{2}\right) \\
\left( -1\right) ^{\alpha \left( 1_{H};0,0,0,1\right) }B(x_{2}\otimes
gx_{1}x_{2};GX_{1}X_{2},gx_{1})GX_{1}X_{2}\otimes g\otimes x_{2}+vB\left(
x_{2}\otimes gx_{1}x_{2};GX_{1}X_{2},gx_{2}\right) \\
\left( -1\right) ^{\alpha \left( 1_{H};1,0,0,1\right) }B(x_{2}\otimes
gx_{1}x_{2};GX_{1}X_{2},gx_{1})GX_{2}\otimes g\otimes gx_{1}x_{2}vB\left(
x_{2}\otimes gx_{1}x_{2};GX_{2},gx_{1}x_{2}\right) \\
\left( -1\right) ^{\alpha \left( 1_{H};0,1,0,1\right) }B(x_{2}\otimes
gx_{1}x_{2};GX_{1}X_{2},gx_{1})GX_{1}\otimes g\otimes gx_{1}x_{2}^{1+1}=0 \\
\left( -1\right) ^{\alpha \left( 1_{H};1,1,0,1\right) }B(x_{2}\otimes
gx_{1}x_{2};GX_{1}X_{2},gx_{1})G\otimes g\otimes x_{1}x_{2}^{1+1}=0
\end{gather*}

\subsubsection{Case $GX_{1}\otimes gx_{2}\otimes x_{2}$}

This case was already considered in subsection $B\left( x_{2}\otimes
gx_{1}x_{2};GX_{2},gx_{1}x_{2}\right) .$

\subsubsection{Case $GX_{1}\otimes gx_{2}\otimes x_{2}$}

First summand of the left side of the equality gives us%
\begin{eqnarray*}
l_{1} &=&u_{1}=0,l_{2}=u_{2}=0 \\
a &=&b_{1}=1,b_{2}=0 \\
d &=&e_{2}=1,e_{1}=0
\end{eqnarray*}%
Since $\alpha \left( x_{2};0,0,0,0\right) \equiv a+b_{1}+b_{2}\equiv 0$ we
get
\begin{equation*}
B(g\otimes gx_{1}x_{2};GX_{1},gx_{2})GX_{1}\otimes gx_{2}\otimes x_{2}.
\end{equation*}%
Second summand of the left side of the equality gives us

\begin{eqnarray*}
l_{1} &=&u_{1}=0,l_{2}+u_{2}=1 \\
a &=&b_{1}=1,0,b_{2}=l_{2} \\
d &=&1,e_{1}=0,e_{2}-u_{2}=1\Rightarrow e_{2}=1,u_{2}=0,b_{2}=l_{2}=1.
\end{eqnarray*}%
Since $\alpha \left( 1_{H};0,1,0,0\right) \equiv 0$ we get
\begin{equation*}
+B(x_{2}\otimes gx_{1}x_{2};GX_{1}X_{2},gx_{2})GX_{1}\otimes gx_{2}\otimes
x_{2}.
\end{equation*}

By considering also the right side of the equality we get

\begin{equation*}
-B(x_{2}\otimes gx_{1};GX_{1},gx_{2})+B(g\otimes
gx_{1}x_{2};GX_{1},gx_{2})+B(x_{2}\otimes gx_{1}x_{2};GX_{1}X_{2},gx_{2})=0
\end{equation*}%
which holds in view of the form of the elements.

\subsubsection{Case $G\otimes gx_{1}\otimes gx_{1}x_{2}$}

This case was already considered in subsection $B\left( x_{2}\otimes
gx_{1}x_{2};GX_{1},gx_{1}x_{2}\right) .$

\subsubsection{Case $GX_{1}X_{2}\otimes g\otimes x_{2}$}

This case was already considered in subsection $B\left( x_{2}\otimes
gx_{1}x_{2};GX_{1}X_{2},gx_{2}\right) .$

\subsubsection{Case $GX_{2}\otimes g\otimes gx_{1}x_{2}$}

This case was already considered in subsection $B\left( x_{2}\otimes
gx_{1}x_{2};GX_{2},gx_{1}x_{2}\right) .$

\section{$B\left( x_{1}x_{2}\otimes x_{1}\right) $}

By $\left( \ref{simplx}\right) $ we have

\begin{equation}
B(x_{1}x_{2}\otimes x_{1})=B(x_{1}x_{2}\otimes 1_{H})(1_{A}\otimes
x_{1})-(1_{A}\otimes gx_{1})B(x_{1}x_{2}\otimes 1_{H})(1_{A}\otimes g)
\label{form x1x2otx1}
\end{equation}%
and we obtain%
\begin{eqnarray*}
&&B(x_{1}x_{2}\otimes x_{1})= \\
&&-2B(x_{1}x_{2}\otimes 1_{H};1_{A},gx_{2})1_{A}\otimes gx_{1}x_{2}+ \\
&&-2B(x_{1}x_{2}\otimes 1_{H};G,g)G\otimes gx_{1}+ \\
&&-2B(x_{1}x_{2}\otimes 1_{H};X_{1},g)X_{1}\otimes gx_{1}+ \\
&&-2B(x_{1}x_{2}\otimes 1_{H};X_{2},g)X_{2}\otimes gx_{1}+ \\
&&-2B(x_{1}x_{2}\otimes 1_{H};X_{2},gx_{1}x_{2})X_{1}X_{2}\otimes gx_{1}x_{2}
\\
&&-2B(x_{1}x_{2}\otimes 1_{H};GX_{1},gx_{2})GX_{1}\otimes gx_{1}x_{2}+ \\
&&-2B(x_{1}x_{2}\otimes 1_{H};GX_{2},gx_{2})GX_{2}\otimes gx_{1}x_{2}+ \\
&&+2\left[ +B(x_{1}x_{2}\otimes 1_{H};G,gx_{1}x_{2})-B(x_{1}x_{2}\otimes
1_{H};GX_{2},gx_{1})+B(x_{1}x_{2}\otimes 1_{H};GX_{1},gx_{2})\right]
GX_{1}X_{2}\otimes gx_{1}
\end{eqnarray*}%
We write Casimir condition for $B(x_{1}x_{2}\otimes x_{1})$%
\begin{eqnarray*}
&&\sum_{w_{1}=0}^{1}\sum_{w_{2}=0}^{1}\left( -1\right) ^{\left(
1+w_{2}\right)
}\sum_{a,b_{1},b_{2},d,e_{1},e_{2}=0}^{1}\sum_{l_{1}=0}^{b_{1}}%
\sum_{l_{2}=0}^{b_{2}}\sum_{u_{1}=0}^{e_{1}}\sum_{u_{2}=0}^{e_{2}} \\
&&\left( -1\right) ^{\alpha \left(
x_{1}^{1-w_{1}}x_{2}^{1-w_{2}};l_{1},l_{2},u_{1},u_{2}\right)
}B(g^{w_{1}+w_{2}}x_{1}^{w_{1}}x_{2}^{w_{2}}\otimes
x_{1};G^{a}X_{1}^{b_{1}}X_{2}^{b_{2}},g^{d}x_{1}^{e_{1}}x_{2}^{e_{2}}) \\
&&G^{a}X_{1}^{b_{1}-l_{1}}X_{2}^{b_{2}-l_{2}}\otimes
g^{d}x_{1}^{e_{1}-u_{1}}x_{2}^{e_{2}-u_{2}}\otimes \\
&&g^{a+b_{1}+b_{2}+l_{1}+l_{2}+d+e_{1}+e_{2}+u_{1}+u_{2}}x_{1}^{l_{1}+u_{1}+1-w_{1}}x_{2}^{l_{2}+u_{2}+1-w_{2}}
\\
&=&\sum_{\omega _{1}=0}^{1}B^{A}(x_{1}x_{2}\otimes x_{1}^{1-\omega
_{1}})\otimes B^{H}(x_{1}x_{2}\otimes x_{1}^{1-\omega _{1}})\otimes
g^{1+\omega _{1}}x_{1}^{\omega _{1}}
\end{eqnarray*}%
Thus we get%
\begin{eqnarray*}
&&\sum_{a,b_{1},b_{2},d,e_{1},e_{2}=0}^{1}\sum_{l_{1}=0}^{b_{1}}%
\sum_{l_{2}=0}^{b_{2}}\sum_{u_{1}=0}^{e_{1}}\sum_{u_{2}=0}^{e_{2}}\left(
-1\right) ^{\alpha \left( x_{1};l_{1},l_{2},u_{1},u_{2}\right) } \\
&&B(gx_{2}\otimes
x_{1};G^{a}X_{1}^{b_{1}}X_{2}^{b_{2}},g^{d}x_{1}^{e_{1}}x_{2}^{e_{2}})G^{a}X_{1}^{b_{1}-l_{1}}X_{2}^{b_{2}-l_{2}}\otimes
\\
&&g^{d}x_{1}^{e_{1}-u_{1}}x_{2}^{e_{2}-u_{2}}\otimes
g^{a+b_{1}+b_{2}+l_{1}+l_{2}+d+e_{1}+e_{2}+u_{1}+u_{2}}x_{1}^{l_{1}+u_{1}+1}x_{2}^{l_{2}+u_{2}}
\\
&&-\sum_{a,b_{1},b_{2},d,e_{1},e_{2}=0}^{1}\sum_{l_{1}=0}^{b_{1}}%
\sum_{l_{2}=0}^{b_{2}}\sum_{u_{1}=0}^{e_{1}}\sum_{u_{2}=0}^{e_{2}}\left(
-1\right) ^{\alpha \left( x_{2};l_{1},l_{2},u_{1},u_{2}\right) } \\
&&B(gx_{1}\otimes
x_{1};G^{a}X_{1}^{b_{1}}X_{2}^{b_{2}},g^{d}x_{1}^{e_{1}}x_{2}^{e_{2}})G^{a}X_{1}^{b_{1}-l_{1}}X_{2}^{b_{2}-l_{2}}
\\
&&\otimes g^{d}x_{1}^{e_{1}-u_{1}}x_{2}^{e_{2}-u_{2}}\otimes
g^{a+b_{1}+b_{2}+l_{1}+l_{2}+d+e_{1}+e_{2}+u_{1}+u_{2}}x_{1}^{l_{1}+u_{1}}x_{2}^{l_{2}+u_{2}+1}
\\
&&\text{ \ \ \ \ }+\sum_{a,b_{1},b_{2},d,e_{1},e_{2}=0}^{1}%
\sum_{l_{1}=0}^{b_{1}}\sum_{l_{2}=0}^{b_{2}}\sum_{u_{1}=0}^{e_{1}}%
\sum_{u_{2}=0}^{e_{2}}\left( -1\right) ^{\alpha \left(
x_{1}x_{2};l_{1},l_{2},u_{1},u_{2}\right) } \\
&&B(1_{H}\otimes
x_{1};G^{a}X_{1}^{b_{1}}X_{2}^{b_{2}},g^{d}x_{1}^{e_{1}}x_{2}^{e_{2}})G^{a}X_{1}^{b_{1}-l_{1}}X_{2}^{b_{2}-l_{2}}\otimes
\\
&&g^{d}x_{1}^{e_{1}-u_{1}}x_{2}^{e_{2}-u_{2}}\otimes
g^{a+b_{1}+b_{2}+l_{1}+l_{2}+d+e_{1}+e_{2}+u_{1}+u_{2}}x_{1}^{l_{1}+u_{1}+1}x_{2}^{l_{2}+u_{2}+1}+
\\
&&\text{4}+\sum_{a,b_{1},b_{2},d,e_{1},e_{2}=0}^{1}\sum_{l_{1}=0}^{b_{1}}%
\sum_{l_{2}=0}^{b_{2}}\sum_{u_{1}=0}^{e_{1}}\sum_{u_{2}=0}^{e_{2}}\left(
-1\right) ^{\alpha \left( 1_{H};l_{1},l_{2},u_{1},u_{2}\right) } \\
&&B(x_{1}x_{2}\otimes
x_{1};G^{a}X_{1}^{b_{1}}X_{2}^{b_{2}},g^{d}x_{1}^{e_{1}}x_{2}^{e_{2}})G^{a}X_{1}^{b_{1}-l_{1}}X_{2}^{b_{2}-l_{2}}
\\
&&\otimes g^{d}x_{1}^{e_{1}-u_{1}}x_{2}^{e_{2}-u_{2}}\otimes
g^{a+b_{1}+b_{2}+l_{1}+l_{2}+d+e_{1}+e_{2}+u_{1}+u_{2}}x_{1}^{l_{1}+u_{1}}x_{2}^{l_{2}+u_{2}}+
\\
&=&B^{A}(x_{1}x_{2}\otimes x_{1})\otimes B^{H}(x_{1}x_{2}\otimes
x_{1})\otimes g+ \\
&&B^{A}(x_{1}x_{2}\otimes 1_{H})\otimes B^{H}(x_{1}x_{2}\otimes
1_{H})\otimes x_{1}
\end{eqnarray*}

\subsection{$B(x_{1}x_{2}\otimes x_{1};1_{A},gx_{1}x_{2})$}

We deduce that%
\begin{eqnarray*}
a &=&b_{1}=b_{2}=0 \\
d &=&e_{1}=e_{2}=1
\end{eqnarray*}%
and get%
\begin{eqnarray*}
&&\sum_{u_{1}=0}^{e_{1}}\sum_{u_{2}=0}^{e_{2}}\left( -1\right) ^{\alpha
\left( 1_{H};0,0,u_{1},u_{2}\right) }B(x_{1}x_{2}\otimes
x_{1};1_{A},gx_{1}x_{2})1_{A}\otimes gx_{1}^{1-u_{1}}x_{2}^{1-u_{2}}\otimes
g^{u_{1}+u_{2}+1}x_{1}^{u_{1}}x_{2}^{u_{2}} \\
&=&\left( -1\right) ^{\alpha \left( 1_{H};0,0,0,0\right)
}B(x_{1}x_{2}\otimes x_{1};1_{A},gx_{1}x_{2})1_{A}\otimes gx_{1}x_{2}\otimes
g+ \\
&&\left( -1\right) ^{\alpha \left( 1_{H};0,0,0,1\right) }B(x_{1}x_{2}\otimes
x_{1};1_{A},gx_{1}x_{2})1_{A}\otimes gx_{1}\otimes x_{2}+\text{ } \\
&&\left( -1\right) ^{\alpha \left( 1_{H};0,0,1,0\right) }B(x_{1}x_{2}\otimes
x_{1};1_{A},gx_{1}x_{2})1_{A}\otimes gx_{2}\otimes x_{1}+ \\
&&\left( -1\right) ^{\alpha \left( 1_{H};0,0,1,1\right) }B(x_{1}x_{2}\otimes
x_{1};1_{A},gx_{1}x_{2})1_{A}\otimes g\otimes gx_{1}x_{2}+
\end{eqnarray*}

\subsubsection{Case $1_{A}\otimes gx_{1}\otimes x_{2}$}

We have to consider only the second and the fourth summands of the left side
of the equality.Second summand gives us%
\begin{eqnarray*}
l_{1} &=&u_{1}=l_{2}=u_{2}=0, \\
a &=&b_{1}=b_{2}=0, \\
d &=&e_{1}=1,e_{2}=0.
\end{eqnarray*}%
Since $\alpha \left( x_{2};0,0,0,0\right) \equiv a+b_{1}+b_{2}=0$ we get%
\begin{equation*}
-B(gx_{1}\otimes x_{1};1_{A},gx_{1})1_{A}\otimes gx_{1}\otimes x_{2}.
\end{equation*}%
Fourth summand gives us%
\begin{eqnarray*}
l_{1} &=&u_{1}=0,l_{2}+u_{2}=1 \\
a &=&b_{1}=0,b_{2}=l_{2} \\
d &=&e_{1}=1,e_{2}=l_{2}
\end{eqnarray*}%
Since $\alpha \left( 1_{H};0,0,0,1\right) \equiv a+b_{1}+b_{2}\equiv 0$ and $%
\alpha \left( 1_{H};0,1,0,0\right) \equiv 0,$ we get%
\begin{equation*}
\left[ -B(x_{1}x_{2}\otimes x_{1};1_{A},gx_{1}x_{2})+B(x_{1}x_{2}\otimes
x_{1};X_{2},gx_{1})\right] 1_{A}\otimes gx_{1}\otimes x_{2}.
\end{equation*}%
Since there is nothing in the right side, we get%
\begin{equation*}
B(gx_{1}\otimes x_{1};1_{A},gx_{1})+B(x_{1}x_{2}\otimes
x_{1};1_{A},gx_{1}x_{2})+B(x_{1}x_{2}\otimes x_{1};X_{2},gx_{1})
\end{equation*}%
which holds in view of the form of the elements.

\subsubsection{Case $1_{A}\otimes gx_{2}\otimes x_{1}$}

We have to consider only the first and the fourth summand of the left side
of the equation. First summand gives us%
\begin{eqnarray*}
l_{1} &=&u_{1}=l_{2}=u_{2}=0 \\
a &=&b_{1}=b_{2}=0 \\
d &=&e_{2}=1,e_{1}=0.
\end{eqnarray*}%
Since $\alpha \left( x_{1};0,0,0,0\right) \equiv a+b_{1}+b_{2}=0,$ we get%
\begin{equation*}
B(gx_{2}\otimes x_{1};1_{A},gx_{2})1_{A}\otimes gx_{2}\otimes x_{1}.
\end{equation*}%
Fourth summand is giving us%
\begin{eqnarray*}
l_{1}+u_{1} &=&1,l_{2}=u_{2}=0, \\
a &=&b_{2}=0,b_{1}=l_{1}, \\
d &=&e_{2}=1,e_{1}=u_{1}.
\end{eqnarray*}%
Since $\alpha \left( 1_{H};0,0,1,0\right) \equiv e_{2}+\left(
a+b_{1}+b_{2}\right) =1$ and $\alpha \left( 1_{H};1,0,0,0\right) \equiv
b_{2}=0,$ we get%
\begin{equation*}
\left[ -B(x_{1}x_{2}\otimes x_{1};1_{A},gx_{1}x_{2})+B(x_{1}x_{2}\otimes
x_{1};X_{1},gx_{2})\right] .
\end{equation*}%
By considering also the right side, we obtain%
\begin{gather*}
B(gx_{2}\otimes x_{1};1_{A},gx_{2})-B(x_{1}x_{2}\otimes
x_{1};1_{A},gx_{1}x_{2}) \\
B(x_{1}x_{2}\otimes x_{1};X_{1},gx_{2})-B(x_{1}x_{2}\otimes
1_{H};1_{A},gx_{2})=0
\end{gather*}%
which holds in view of the form of the elements.

\subsubsection{Case $1_{A}\otimes g\otimes gx_{1}x_{2}$}

First summand of the left side of the equality gives us%
\begin{eqnarray*}
l_{1} &=&u_{1}=0,l_{2}+u_{2}=1 \\
a &=&b_{1}=0,b_{2}=l_{2}, \\
d &=&1,e_{1}=0,e_{2}=u_{2}.
\end{eqnarray*}%
Since $\alpha \left( x_{1};0,0,0,1\right) \equiv 1$ and $\alpha \left(
x_{1};0,1,0,0\right) \equiv a+b_{1}+b_{2}+1\equiv 0,$ we get%
\begin{equation*}
\left[ -B(gx_{2}\otimes x_{1};1_{A},gx_{2})+B(gx_{2}\otimes x_{1};X_{2},g)%
\right] 1_{A}\otimes g\otimes gx_{1}x_{2}.
\end{equation*}%
Second summand of the left side of the equation gives us%
\begin{eqnarray*}
l_{2} &=&u_{2}=0,l_{1}+u_{1}=1 \\
a &=&b_{2}=0,b_{1}=l_{1}, \\
d &=&1,e_{2}=0,e_{1}=u_{1}.
\end{eqnarray*}%
Since $\alpha \left( x_{2};0,0,1,0\right) \equiv e_{2}=0$ and $\alpha \left(
x_{2};1,0,0,0\right) \equiv a+b_{1}=1,$ we get%
\begin{equation*}
\left[ -B(gx_{1}\otimes x_{1};1_{A},gx_{1})+B(gx_{1}\otimes x_{1};X_{1},g)%
\right] 1_{A}\otimes g\otimes gx_{1}x_{2}.
\end{equation*}

Third summand of the left side of the equation gives us%
\begin{eqnarray*}
l_{1} &=&u_{1}=l_{2}=u_{2}=0, \\
a &=&b_{1}=b_{2}=0, \\
d &=&1,e_{1}=e_{2}=0.
\end{eqnarray*}%
Since $\alpha \left( x_{1}x_{2};0,0,0,0\right) \equiv 0,$ we get%
\begin{equation*}
B(1_{H}\otimes x_{1};1_{A},g)1_{A}\otimes g\otimes gx_{1}x_{2}.
\end{equation*}%
Fourth summand of the left side of the equation gives us%
\begin{eqnarray*}
l_{1}+u_{1} &=&1,l_{2}+u_{2}=1 \\
a &=&0,b_{1}=l_{1},b_{2}=l_{2}, \\
d &=&1,e_{1}=u_{1},e_{2}=u_{2}.
\end{eqnarray*}%
Since%
\begin{eqnarray*}
&&\alpha \left( 1_{H};0,0,1,1\right) \equiv 1+e_{2}\equiv 0 \\
&&\alpha \left( 1_{H};0,1,1,0\right) \equiv e_{2}+a+b_{1}+b_{2}+1\equiv 0 \\
&&\alpha \left( 1_{H};1,0,0,1\right) \equiv a+b_{1}\equiv 1 \\
&&\alpha \left( 1_{H};1,1,0,0\right) \equiv 1+b_{2}\equiv 0
\end{eqnarray*}%
we get%
\begin{eqnarray*}
&&\lbrack B(x_{1}x_{2}\otimes x_{1};1_{A},gx_{1}x_{2})+B(x_{1}x_{2}\otimes
x_{1};X_{2},gx_{1})+ \\
&&-B(x_{1}x_{2}\otimes x_{1};X_{1},gx_{2})+B(x_{1}x_{2}\otimes
x_{1};X_{1}X_{2},1_{H})]1_{A}\otimes g\otimes gx_{1}x_{2}.
\end{eqnarray*}%
Since there is nothing on the right side, we get%
\begin{eqnarray*}
&&-B(gx_{2}\otimes x_{1};1_{A},gx_{2})+B(gx_{2}\otimes x_{1};X_{2},g)+ \\
&&-B(gx_{1}\otimes x_{1};1_{A},gx_{1})+B(gx_{1}\otimes x_{1};X_{1},g)+ \\
&&+B(1_{H}\otimes x_{1};1_{A},g)+B(x_{1}x_{2}\otimes
x_{1};1_{A},gx_{1}x_{2})+B(x_{1}x_{2}\otimes x_{1};X_{2},gx_{1})+ \\
&&-B(x_{1}x_{2}\otimes x_{1};X_{1},gx_{2})+B(x_{1}x_{2}\otimes
x_{1};X_{1}X_{2},1_{H})=0
\end{eqnarray*}%
which holds in view of the form of the elements.

\subsection{$B\left( x_{1}x_{2}\otimes x_{1};G,gx_{1}\right) $}

We deduce that%
\begin{equation*}
a=1,d=e_{1}=1
\end{equation*}%
and we get%
\begin{eqnarray*}
&&\left( -1\right) ^{\alpha \left( 1_{H};0,0,0,0\right) }B(x_{1}x_{2}\otimes
x_{1};G,gx_{1})G\otimes gx_{1}\otimes g+ \\
&&\left( -1\right) ^{\alpha \left( 1_{H};0,0,u_{1},0\right)
}B(x_{1}x_{2}\otimes x_{1};G,gx_{1})G\otimes g\otimes x_{1}.
\end{eqnarray*}

\subsubsection{Case $G\otimes g\otimes x_{1}$}

We have to consider only the first and the fourth summand of the left side
of the equality

First summand gives us%
\begin{eqnarray*}
l_{1} &=&u_{1}=0,l_{2}=u_{2}=0, \\
a &=&1,b_{1}=b_{2}=0, \\
d &=&1,e_{1}=e_{2}=0.
\end{eqnarray*}%
Since $\alpha \left( x_{1};0,0,0,0\right) \equiv a+b_{1}+b_{2}=1,$ we get
\begin{equation*}
-B(gx_{2}\otimes x_{1};G,g)G\otimes g\otimes x_{1}
\end{equation*}%
Fourth summand gives us%
\begin{eqnarray*}
l_{1}+u_{1} &=&1,l_{2}=u_{2}=0 \\
a &=&1,b_{1}=l_{1},b_{2}=0 \\
d &=&1,e_{1}=u_{1},e_{2}=0.
\end{eqnarray*}%
Since $\alpha \left( 1_{H};0,0,1,0\right) \equiv e_{2}+\left(
a+b_{1}+b_{2}\right) \equiv 1$ and $\alpha \left( 1_{H};1,0,0,0\right)
\equiv b_{2}\equiv 0,$ we get
\begin{equation*}
\left[ -B(x_{1}x_{2}\otimes x_{1};G,gx_{1})+B(x_{1}x_{2}\otimes
x_{1};GX_{1},g)\right] G\otimes g\otimes x_{1}
\end{equation*}%
By considering also the right side we get%
\begin{equation*}
-B(x_{1}x_{2}\otimes 1_{H};G,g)-B(x_{1}x_{2}\otimes
x_{1};G,gx_{1})+B(x_{1}x_{2}\otimes x_{1};GX_{1},g)-B(gx_{2}\otimes
x_{1};G,g)=0
\end{equation*}%
which holds in view of the form of the elements.

\subsection{$B\left( x_{1}x_{2}\otimes x_{1};X_{1},gx_{1}\right) $ \ \ \ \ }

We deduce that%
\begin{equation*}
b_{1}=1,d=e_{1}=1
\end{equation*}%
and we get%
\begin{eqnarray*}
&&+\sum_{l_{1}=0}^{1}\sum_{u_{1}=0}^{1}\left( -1\right) ^{\alpha \left(
1_{H};l_{1},0,u_{1},0\right) }B(x_{1}x_{2}\otimes
x_{1};X_{1},gx_{1})X_{1}^{1-l_{1}}\otimes gx_{1}^{1-u_{1}}\otimes
g^{l_{1}+u_{1}+1}x_{1}^{l_{1}+u_{1}} \\
&=&\left( -1\right) ^{\alpha \left( 1_{H};0,0,1,0\right)
}B(x_{1}x_{2}\otimes x_{1};X_{1},gx_{1})X_{1}\otimes g\otimes x_{1}+ \\
&&\left( -1\right) ^{\alpha \left( 1_{H};1,0,0,0\right) }B(x_{1}x_{2}\otimes
x_{1};X_{1},gx_{1})1_{A}\otimes gx_{1}\otimes x_{1}
\end{eqnarray*}

\subsubsection{Case $X_{1}\otimes g\otimes x_{1}$}

We have to consider only the first and the fourth summand of the left side
of the equality.

First summand gives us%
\begin{eqnarray*}
l_{1} &=&u_{1}=0,l_{2}=u_{2}=0, \\
a &=&b_{2}=0,b_{1}=1 \\
d &=&1,e_{1}=e_{2}=0.
\end{eqnarray*}%
Since $\alpha \left( x_{1};0,0,0,0\right) \equiv a+b_{1}+b_{2}=1$ we get
\begin{equation*}
-B(gx_{2}\otimes x_{1};X_{1},g)X_{1}\otimes g\otimes x_{1}.
\end{equation*}%
Fourth summand gives us%
\begin{eqnarray*}
l_{1}+u_{1} &=&1,l_{2}=u_{2}=0 \\
a &=&b_{2}=0,b_{1}-l_{1}=1\Rightarrow b_{1}=1,l_{1}=0,e_{1}=u_{1}=1 \\
d &=&1,e_{1}=u_{1}=1,e_{2}=0.
\end{eqnarray*}%
Since $\alpha \left( 1_{H};0,0,1,0\right) \equiv e_{2}+\left(
a+b_{1}+b_{2}\right) \equiv 1,$ we get
\begin{equation*}
-B(x_{1}x_{2}\otimes x_{1};X_{1},gx_{1})X_{1}\otimes g\otimes x_{1}.
\end{equation*}%
By considering also the right side we get%
\begin{equation*}
-B(x_{1}x_{2}\otimes 1;X_{1},g)-B(x_{1}x_{2}\otimes
x_{1};X_{1},gx_{1})-B(gx_{2}\otimes x_{1};X_{1},g)=0
\end{equation*}

which holds in view of the form of the elements.

\subsubsection{Case $1_{A}\otimes gx_{1}\otimes x_{1}$}

We have to consider only the first and the fourth summand of the left side
of the equality.

First summand gives us
\begin{eqnarray*}
l_{1} &=&u_{1}=0,l_{2}=u_{2}=0, \\
a &=&b_{1}=b_{2}=0, \\
d &=&e_{1}=1,e_{2}=0.
\end{eqnarray*}%
Since $\alpha \left( x_{1};0,0,0,0\right) \equiv a+b_{1}+b_{2}=0$ we get%
\begin{equation*}
B(gx_{2}\otimes x_{1};1_{A},gx_{1})1_{A}\otimes gx_{1}\otimes x_{1}
\end{equation*}%
Fourth summand gives us%
\begin{eqnarray*}
l_{1}+u_{1} &=&1,l_{2}=u_{2}=0 \\
a &=&b_{2}=0,b_{1}=l_{1}=1 \\
d &=&1,e_{1}-u_{1}=1\Rightarrow e_{1}=1,u_{1}=0,l_{1}=1,e_{2}=0
\end{eqnarray*}%
Since $\alpha \left( 1_{H};1,0,0,0\right) \equiv b_{2}=0,$ we get
\begin{equation*}
B(x_{1}x_{2}\otimes x_{1};X_{1},gx_{1})1_{A}\otimes gx_{1}\otimes x_{1}.
\end{equation*}%
By considering also the right side we get%
\begin{equation*}
-B(x_{1}x_{2}\otimes 1_{H};1_{A},gx_{1})+B(x_{1}x_{2}\otimes
x_{1};X_{1},gx_{1})+B(gx_{2}\otimes x_{1};1_{A},gx_{1})=0
\end{equation*}

which holds in view of the form of the elements.

\subsection{$B\left( x_{1}x_{2}\otimes x_{1};X_{2},gx_{1}\right) $}

We deduce that
\begin{equation*}
b_{2}=1,d=e_{1}=1
\end{equation*}%
and we get%
\begin{eqnarray*}
&&\left( -1\right) ^{\alpha \left( 1_{H};0,0,0,0\right) }B\left(
x_{1}x_{2}\otimes x_{1};X_{2},gx_{1}\right) X_{2}\otimes gx_{1}\otimes g+ \\
&&+\left( -1\right) ^{\alpha \left( 1_{H};0,1,0,0\right) }B\left(
x_{1}x_{2}\otimes x_{1};X_{2},gx_{1}\right) 1_{A}\otimes gx_{1}\otimes x_{2}+
\\
&&+\left( -1\right) ^{\alpha \left( 1_{H};0,0,1,0\right) }B\left(
x_{1}x_{2}\otimes x_{1};X_{2},gx_{1}\right) X_{2}\otimes g\otimes x_{1}+ \\
&&+\left( -1\right) ^{\alpha \left( 1_{H};0,1,1,0\right) }B\left(
x_{1}x_{2}\otimes x_{1};X_{2},gx_{1}\right) 1_{A}\otimes g\otimes gx_{1}x_{2}
\end{eqnarray*}

\subsubsection{Case $1_{A}\otimes gx_{1}\otimes x_{2}$}

This case was already considered in subsection $B(x_{1}x_{2}\otimes
x_{1};1_{A},gx_{1}x_{2}).$

\subsubsection{Case $X_{2}\otimes g\otimes x_{1}$}

We have to consider only the first and the fourth summand of the left side
of the equation. First summand gives us%
\begin{eqnarray*}
l_{1} &=&u_{1}=l_{2}=u_{2}=0 \\
a &=&b_{1}=0,b_{2}=1 \\
d &=&1,e_{1}=e_{2}=0.
\end{eqnarray*}%
Since $\alpha \left( x_{1};0,0,0,0\right) \equiv a+b_{1}+b_{2}=1,$ we get%
\begin{equation*}
-B(gx_{2}\otimes x_{1};X_{2},g)X_{2}\otimes g\otimes x_{1}.
\end{equation*}%
Fourth summand is giving us%
\begin{eqnarray*}
l_{1}+u_{1} &=&1,l_{2}=u_{2}=0, \\
a &=&0,b_{1}=l_{1},b_{2}=1, \\
d &=&1,e_{1}=u_{1},e_{2}=0
\end{eqnarray*}%
Since $\alpha \left( 1_{H};0,0,1,0\right) \equiv e_{2}+\left(
a+b_{1}+b_{2}\right) =1$ and $\alpha \left( 1_{H};1,0,0,0\right) \equiv
b_{2}=1,$ we get%
\begin{equation*}
\left[ -B(x_{1}x_{2}\otimes x_{1};X_{2},gx_{1})-B(x_{1}x_{2}\otimes
x_{1};X_{1}X_{2},g)\right] X_{2}\otimes g\otimes x_{1}.
\end{equation*}%
By considering also the right side, we obtain%
\begin{gather*}
-B(gx_{2}\otimes x_{1};X_{2},g)-B(x_{1}x_{2}\otimes x_{1};X_{2},gx_{1}) \\
-B(x_{1}x_{2}\otimes x_{1};X_{1}X_{2},g)-B(x_{1}x_{2}\otimes 1_{H};X_{2},g)=0
\end{gather*}%
which holds in view of the form of the elements.

\subsubsection{Case $1_{A}\otimes g\otimes gx_{1}x_{2}$}

This case was already considered in subsection $B(x_{1}x_{2}\otimes
x_{1};1_{A},gx_{1}x_{2}).$

\subsection{$B\left( x_{1}x_{2}\otimes x_{1};X_{1}X_{2},gx_{1}x_{2}\right) $}

We deduce that%
\begin{equation*}
a=0,b_{1}=b_{2}=1,d=e_{1}=e_{2}=1
\end{equation*}%
\begin{equation*}
X_{1}^{1-l_{1}}X_{2}^{1-l_{2}}\otimes gx_{1}^{1-u_{1}}x_{2}^{1-u}\otimes
g^{1+l_{1}+l_{2}+u_{1}+u_{2}}x_{1}^{l_{1}+u_{1}}x_{2}^{l_{2}+u_{2}}
\end{equation*}%
and we get

\begin{gather*}
\left( -1\right) ^{\alpha \left( 1_{H};0,0,1,0\right) }B(x_{1}x_{2}\otimes
x_{1};X_{1}X_{2},gx_{1}x_{2})X_{1}X_{2}\otimes gx_{2}\otimes x_{1}+ \\
\left( -1\right) ^{\alpha \left( 1_{H};0,1,1,0\right) }B(x_{1}x_{2}\otimes
x_{1};X_{1}X_{2},gx_{1}x_{2})X_{1}\otimes gx_{2}\otimes gx_{1}x_{2} \\
+\left( -1\right) ^{\alpha \left( 1_{H};1,0,1,0\right) }B(x_{1}x_{2}\otimes
x_{1};X_{1}X_{2},gx_{1}x_{2})X_{1}^{1-l_{1}}X_{2}^{1-l_{2}}\otimes
gx_{2}\otimes g^{+l_{1}+l_{2}}x_{1}^{1+1}x_{2}^{l_{2}}=0 \\
+\left( -1\right) ^{\alpha \left( 1_{H};1,1,1,0\right) }B(x_{1}x_{2}\otimes
x_{1};X_{1}X_{2},gx_{1}x_{2})X_{1}^{1-l_{1}}X_{2}^{1-l_{2}}\otimes
gx_{2}\otimes g^{+l_{1}+l_{2}}x_{1}^{1+1}x_{2}^{l_{2}}=0 \\
+\left( -1\right) ^{\alpha \left( 1_{H};0,0,0,0\right) }B(x_{1}x_{2}\otimes
x_{1};X_{1}X_{2},gx_{1}x_{2})X_{1}X_{2}\otimes gx_{1}x_{2}\otimes g+ \\
+\left( -1\right) ^{\alpha \left( 1_{H};0,1,0,0\right) }B(x_{1}x_{2}\otimes
x_{1};X_{1}X_{2},gx_{1}x_{2})X_{1}\otimes gx_{1}x_{2}\otimes x_{2}+ \\
+\left( -1\right) ^{\alpha \left( 1_{H};1,0,0,0\right) }B(x_{1}x_{2}\otimes
x_{1};X_{1}X_{2},gx_{1}x_{2})X_{2}\otimes gx_{1}x_{2}\otimes x_{1}+ \\
+\left( -1\right) ^{\alpha \left( 1_{H};1,1,0,0\right) }B(x_{1}x_{2}\otimes
x_{1};X_{1}X_{2},gx_{1}x_{2})1_{A}\otimes gx_{1}x_{2}\otimes gx_{1}x_{2}+ \\
+\left( -1\right) ^{\alpha \left( 1_{H};0,0,0,1\right) }B(x_{1}x_{2}\otimes
x_{1};X_{1}X_{2},gx_{1}x_{2})X_{1}X_{2}\otimes gx_{1}\otimes x_{2}+ \\
+\left( -1\right) ^{\alpha \left( 1_{H};0,1,0,1\right) }B(x_{1}x_{2}\otimes
x_{1};X_{1}X_{2},gx_{1}x_{2})X_{1}^{1-l_{1}}X_{2}^{1-l_{2}}\otimes
gx_{1}\otimes g^{1+l_{1}+l_{2}+1}x_{1}^{l_{1}}x_{2}^{1+1}=0 \\
+\left( -1\right) ^{\alpha \left( 1_{H};1,0,0,1\right) }B(x_{1}x_{2}\otimes
x_{1};X_{1}X_{2},gx_{1}x_{2})X_{2}\otimes gx_{1}\otimes gx_{1}x_{2}+ \\
+\left( -1\right) ^{\alpha \left( 1_{H};1,1,0,1\right) }B(x_{1}x_{2}\otimes
x_{1};X_{1}X_{2},gx_{1}x_{2})X_{1}^{1-l_{1}}X_{2}^{1-l_{2}}\otimes
gx_{1}\otimes g^{1+l_{1}+l_{2}+1}x_{1}^{l_{1}}x_{2}^{1+1}=0 \\
+\left( -1\right) ^{\alpha \left( 1_{H};0,0,1,1\right) }B(x_{1}x_{2}\otimes
x_{1};X_{1}X_{2},gx_{1}x_{2})X_{1}X_{2}\otimes g\otimes gx_{1}x_{2}+ \\
+\left( -1\right) ^{\alpha \left( 1_{H};0,1,1,1\right) }B(x_{1}x_{2}\otimes
x_{1};X_{1}X_{2},gx_{1}x_{2})X_{2}^{1-l_{2}}\otimes g\otimes
g^{1+l_{1}+l_{2}}x_{1}^{l_{1}+1}x_{2}^{1+1}=0 \\
+\left( -1\right) ^{\alpha \left( 1_{H};1,0,1,1\right) }B(x_{1}x_{2}\otimes
x_{1};X_{1}X_{2},gx_{1}x_{2})X_{2}^{1-l_{2}}\otimes g\otimes
g^{1+l_{1}+l_{2}}x_{1}^{1+1}x_{2}^{l_{2}+1}=0 \\
+\left( -1\right) ^{\alpha \left( 1_{H};1,1,1,1\right) }B(x_{1}x_{2}\otimes
x_{1};X_{1}X_{2},gx_{1}x_{2})X_{2}^{1-l_{2}}\otimes g\otimes
g^{1+l_{1}+l_{2}}x_{1}^{1+1}x_{2}^{1+1}=0
\end{gather*}

\subsubsection{Case $X_{1}X_{2}\otimes gx_{2}\otimes x_{1}$}

We have to consider only the first and the fourth summand of the left side
of the equation. First summand gives us%
\begin{eqnarray*}
l_{1} &=&u_{1}=l_{2}=u_{2}=0 \\
a &=&0,b_{1}=b_{2}=1 \\
d &=&e_{2}=1,e_{1}=0.
\end{eqnarray*}%
Since $\alpha \left( x_{1};0,0,0,0\right) \equiv a+b_{1}+b_{2}\equiv 0,$ we
get%
\begin{equation*}
+B(gx_{2}\otimes x_{1};X_{1}X_{2},gx_{2})X_{1}X_{2}\otimes gx_{2}\otimes
x_{1}.
\end{equation*}%
Fourth summand is giving us%
\begin{eqnarray*}
l_{1}+u_{1} &=&1,l_{2}=u_{2}=0, \\
a &=&0,b_{1}-l_{1}=1\Rightarrow b_{1}=1,l_{1}=0,u_{1}=1,b_{2}=1, \\
d &=&1,e_{1}=u_{1}=1,e_{2}=1.
\end{eqnarray*}%
Since $\alpha \left( 1_{H};0,0,1,0\right) \equiv e_{2}+\left(
a+b_{1}+b_{2}\right) \equiv 1$ $,$ we get%
\begin{equation*}
-B(x_{1}x_{2}\otimes x_{1};X_{1}X_{2},gx_{1}x_{2})X_{1}X_{2}\otimes
gx_{2}\otimes x_{1}.
\end{equation*}%
By considering also the right side, we obtain%
\begin{equation*}
+B(gx_{2}\otimes x_{1};X_{1}X_{2},gx_{2})-B(x_{1}x_{2}\otimes
x_{1};X_{1}X_{2},gx_{1}x_{2})-B(x_{1}x_{2}\otimes 1_{H};X_{1}X_{2},gx_{2})=0
\end{equation*}%
which holds in view of the form of the elements.

\subsubsection{Case $X_{1}\otimes gx_{2}\otimes gx_{1}x_{2}$}

First summand of the left side of the equality gives us%
\begin{eqnarray*}
l_{1} &=&u_{1}=0,l_{2}+u_{2}=1 \\
a &=&0,b_{1}=1,b_{2}=l_{2}, \\
d &=&1,e_{1}=0,e_{2}-u_{2}=1\Rightarrow e_{2}=1,u_{2}=0,b_{2}=l_{2}=1.
\end{eqnarray*}%
Since $\alpha \left( x_{1};0,1,0,0\right) \equiv a+b_{1}+b_{2}+1\equiv 1,$
we get%
\begin{equation*}
-B(gx_{2}\otimes x_{1};X_{1}X_{2},gx_{2})X_{1}\otimes gx_{2}\otimes
gx_{1}x_{2}.
\end{equation*}%
Second summand of the left side of the equation gives us%
\begin{eqnarray*}
l_{2} &=&u_{2}=0,l_{1}+u_{1}=1 \\
a &=&b_{2}=0,b_{1}-l_{1}=1\Rightarrow b_{1}=1,l_{1}=0,u_{1}=1 \\
d &=&1,e_{2}=1,e_{1}=u_{1}=1.
\end{eqnarray*}%
Since $\alpha \left( x_{2};0,0,1,0\right) \equiv e_{2}=1$ $,$ we get%
\begin{equation*}
+B(gx_{1}\otimes x_{1};X_{1},gx_{1}x_{2})X_{1}\otimes gx_{2}\otimes
gx_{1}x_{2}.
\end{equation*}

Third summand of the left side of the equation gives us%
\begin{eqnarray*}
l_{1} &=&u_{1}=l_{2}=u_{2}=0, \\
a &=&b_{2}=0,b_{1}=1 \\
d &=&e_{2}=1,e_{1}=0.
\end{eqnarray*}%
Since $\alpha \left( x_{1}x_{2};0,0,0,0\right) \equiv 0,$ we get%
\begin{equation*}
B(1_{H}\otimes x_{1};X_{1},gx_{2})X_{1}\otimes gx_{2}\otimes gx_{1}x_{2}.
\end{equation*}%
Fourth summand of the left side of the equation gives us%
\begin{eqnarray*}
l_{1}+u_{1} &=&1,l_{2}+u_{2}=1 \\
a &=&0,b_{1}-l_{1}=1\Rightarrow b_{1}=1,l_{1}=0,u_{1}=1,b_{2}=l_{2}, \\
d &=&1,e_{1}=u_{1}=1,e_{2}-u_{2}=1\Rightarrow e_{2}=1,u_{2}=0,b_{2}=l_{2}=1.
\end{eqnarray*}%
Since $\alpha \left( 1_{H};0,1,1,0\right) \equiv e_{2}+a+b_{1}+b_{2}+1\equiv
0$ we get%
\begin{equation*}
B\left( x_{1}x_{2}\otimes x_{1};X_{1}X_{2},gx_{1}x_{2}\right) X_{1}\otimes
gx_{2}\otimes gx_{1}x_{2}.
\end{equation*}%
Since there is nothing on the right side, we get%
\begin{gather*}
-B(gx_{2}\otimes x_{1};X_{1}X_{2},gx_{2})+B(gx_{1}\otimes
x_{1};X_{1},gx_{1}x_{2})+ \\
+B(1_{H}\otimes x_{1};X_{1},gx_{2})+B\left( x_{1}x_{2}\otimes
x_{1};X_{1}X_{2},gx_{1}x_{2}\right) =0
\end{gather*}%
which holds in view of the form of the elements.

\subsubsection{Case $X_{1}\otimes gx_{1}x_{2}\otimes x_{2}$}

We have to consider only the second and the fourth summands of the left side
of the equality. Second summand gives us%
\begin{eqnarray*}
l_{1} &=&u_{1}=l_{2}=u_{2}=0, \\
a &=&b_{2}=0,b_{1}=1, \\
d &=&e_{1}=e_{2}=1.
\end{eqnarray*}%
Since $\alpha \left( x_{2};0,0,0,0\right) \equiv a+b_{1}+b_{2}=1$ we get%
\begin{equation*}
B(gx_{1}\otimes x_{1};X_{1},gx_{1}x_{2})X_{1}\otimes gx_{1}x_{2}\otimes
x_{2}.
\end{equation*}%
Fourth summand gives us%
\begin{eqnarray*}
l_{1} &=&u_{1}=0,l_{2}+u_{2}=1 \\
a &=&0,b_{1}=1,b_{2}=l_{2} \\
d &=&e_{1}=1,e_{2}-u_{2}=1\Rightarrow e_{2}=1,u_{2}=0,b_{2}=l_{2}=1
\end{eqnarray*}%
Since $\alpha \left( 1_{H};0,1,0,0\right) \equiv 0,$ we get%
\begin{equation*}
B(x_{1}x_{2}\otimes x_{1};X_{1}X_{2},gx_{1}x_{2})X_{1}\otimes
gx_{1}x_{2}\otimes x_{2}.
\end{equation*}%
Since there is nothing in the right side, we get%
\begin{equation*}
B(gx_{1}\otimes x_{1};X_{1},gx_{1}x_{2})+B(x_{1}x_{2}\otimes
x_{1};X_{1}X_{2},gx_{1}x_{2})=0
\end{equation*}%
which holds in view of the form of the elements.

\subsubsection{Case $X_{2}\otimes gx_{1}x_{2}\otimes x_{1}$}

We have to consider only the first and the fourth summand of the left side
of the equation. First summand gives us%
\begin{eqnarray*}
l_{1} &=&u_{1}=l_{2}=u_{2}=0 \\
a &=&b_{1}=0,b_{2}=1 \\
d &=&e_{1}=e_{2}=1.
\end{eqnarray*}%
Since $\alpha \left( x_{1};0,0,0,0\right) \equiv a+b_{1}+b_{2}\equiv 1,$ we
get%
\begin{equation*}
-B(gx_{2}\otimes x_{1};X_{2},gx_{1}x_{2})X_{2}\otimes gx_{1}x_{2}\otimes
x_{1}.
\end{equation*}%
Fourth summand is giving us%
\begin{eqnarray*}
l_{1}+u_{1} &=&1,l_{2}=u_{2}=0, \\
a &=&0,b_{1}=l_{1},b_{2}=1, \\
d &=&1,e_{1}-u_{1}=1\Rightarrow e_{1}=1,u_{1}=0,b_{1}=l_{1}=1,e_{2}=1.
\end{eqnarray*}%
Since $\alpha \left( 1_{H};1,0,0,0\right) \equiv b_{2}=1$ $,$ we get%
\begin{equation*}
-B(x_{1}x_{2}\otimes x_{1};X_{1}X_{2},gx_{1}x_{2})X_{2}\otimes
gx_{1}x_{2}\otimes x_{1}.
\end{equation*}%
By considering also the right side, we obtain%
\begin{equation*}
-B(gx_{2}\otimes x_{1};X_{2},gx_{1}x_{2})-B(x_{1}x_{2}\otimes
x_{1};X_{1}X_{2},gx_{1}x_{2})-B(x_{1}x_{2}\otimes 1_{H};X_{2},gx_{1}x_{2})=0
\end{equation*}%
which holds in view of the form of the elements.

\subsubsection{Case $1_{A}\otimes gx_{1}x_{2}\otimes gx_{1}x_{2}$}

First summand of the left side of the equality gives us%
\begin{eqnarray*}
l_{1} &=&u_{1}=0,l_{2}+u_{2}=1 \\
a &=&b_{1}=0,b_{2}=l_{2}, \\
d &=&1,e_{1}=1 \\
e_{2}-u_{2} &=&1\Rightarrow e_{2}=1,u_{2}=0,b_{2}=l_{2}=1,
\end{eqnarray*}%
Since $\alpha \left( x_{1};0,1,0,0\right) \equiv a+b_{1}+b_{2}+1\equiv 0,$
we get%
\begin{equation*}
+B(gx_{2}\otimes x_{1};X_{2},gx_{1}x_{2})1_{A}\otimes gx_{1}x_{2}\otimes
gx_{1}x_{2}.
\end{equation*}%
Second summand of the left side of the equation gives us%
\begin{eqnarray*}
l_{2} &=&u_{2}=0,l_{1}+u_{1}=1 \\
a &=&b_{2}=0,b_{1}=l_{1}=1 \\
d &=&1,e_{2}=1,e_{1}-u_{1}=1\Rightarrow e_{1}=1,u_{1}=0,b_{1}=l_{1}=1,.
\end{eqnarray*}%
Since $\alpha \left( x_{2};1,0,0,0\right) \equiv a+b_{1}=1$ $,$ we get%
\begin{equation*}
+B(gx_{1}\otimes x_{1};X_{1},gx_{1}x_{2})1_{A}\otimes gx_{1}x_{2}\otimes
gx_{1}x_{2}.
\end{equation*}

Third summand of the left side of the equation gives us%
\begin{eqnarray*}
l_{1} &=&u_{1}=l_{2}=u_{2}=0, \\
a &=&b_{1}=b_{2}=0, \\
d &=&e_{1}=e_{2}=1.
\end{eqnarray*}%
Since $\alpha \left( x_{1}x_{2};0,0,0,0\right) \equiv 0,$ we get%
\begin{equation*}
B(1_{H}\otimes x_{1};1_{A},gx_{1}x_{2})1_{A}\otimes gx_{1}x_{2}\otimes
gx_{1}x_{2}.
\end{equation*}%
Fourth summand of the left side of the equation gives us%
\begin{eqnarray*}
l_{1}+u_{1} &=&1,l_{2}+u_{2}=1 \\
a &=&0,b_{1}=l_{1},b_{2}=l_{2}, \\
d &=&1,e_{1}-u_{1}=1\Rightarrow e_{1}=1,u_{1}=0,b_{1}=l_{1}=1, \\
e_{2}-u_{2} &=&1\Rightarrow e_{2}=1,u_{2}=0,b_{2}=l_{2}=1
\end{eqnarray*}%
Since $\alpha \left( 1_{H};1,1,0,0\right) \equiv 1+b_{2}\equiv 0$ we get%
\begin{equation*}
B\left( x_{1}x_{2}\otimes x_{1};X_{1}X_{2},gx_{1}x_{2}\right) X_{1}\otimes
gx_{2}\otimes gx_{1}x_{2}.
\end{equation*}%
Since there is nothing on the right side, we get%
\begin{eqnarray*}
&&+B(gx_{2}\otimes x_{1};X_{2},gx_{1}x_{2})+B(gx_{1}\otimes
x_{1};X_{1},gx_{1}x_{2})+ \\
&&B(1_{H}\otimes x_{1};1_{A},gx_{1}x_{2})+B\left( x_{1}x_{2}\otimes
x_{1};X_{1}X_{2},gx_{1}x_{2}\right) =0
\end{eqnarray*}

which holds in view of the form of the elements.

\subsubsection{Case $X_{1}X_{2}\otimes gx_{1}\otimes x_{2}$}

We have to consider only the second and the fourth summands of the left side
of the equality. Second summand gives us%
\begin{eqnarray*}
l_{1} &=&u_{1}=l_{2}=u_{2}=0, \\
a &=&0,b_{1}=b_{2}=1 \\
d &=&e_{1}=1,e_{2}=0.
\end{eqnarray*}%
Since $\alpha \left( x_{2};0,0,0,0\right) \equiv a+b_{1}+b_{2}=0$ we get%
\begin{equation*}
-B(gx_{1}\otimes x_{1};X_{1}X_{2},gx_{1})X_{1}X_{2}\otimes gx_{1}\otimes
x_{2}.
\end{equation*}%
Fourth summand gives us%
\begin{eqnarray*}
l_{1} &=&u_{1}=0,l_{2}+u_{2}=1 \\
a &=&0,b_{1}=1,b_{2}-l_{2}=1\Rightarrow b_{2}=1,l_{2}=0,u_{2}=1 \\
d &=&e_{1}=1,e_{2}=u_{2}=1
\end{eqnarray*}%
Since $\alpha \left( 1_{H};0,0,0,1\right) \equiv a+b_{1}+b_{2}\equiv 0,$ we
get%
\begin{equation*}
+B(x_{1}x_{2}\otimes x_{1};X_{1}X_{2},gx_{1}x_{2})X_{1}X_{2}\otimes
gx_{1}\otimes x_{2}.
\end{equation*}%
Since there is nothing in the right side, we get%
\begin{equation*}
-B(gx_{1}\otimes x_{1};X_{1}X_{2},gx_{1})+B(x_{1}x_{2}\otimes
x_{1};X_{1}X_{2},gx_{1}x_{2})=0
\end{equation*}%
which holds in view of the form of the elements.

\subsubsection{Case $X_{2}\otimes gx_{1}\otimes gx_{1}x_{2}$}

First summand of the left side of the equality gives us%
\begin{eqnarray*}
l_{1} &=&u_{1}=0,l_{2}+u_{2}=1 \\
a &=&b_{1}=0,b_{2}-l_{2}=1\Rightarrow b_{2}=1,l_{2}=0,u_{2}=1 \\
d &=&1,e_{1}=1,e_{2}=u_{2}=1.
\end{eqnarray*}%
Since $\alpha \left( x_{1};0,0,0,1\right) \equiv 1,$ we get%
\begin{equation*}
-B(gx_{2}\otimes x_{1};X_{2},gx_{1}x_{2})X_{2}\otimes gx_{1}\otimes
gx_{1}x_{2}.
\end{equation*}%
Second summand of the left side of the equation gives us%
\begin{eqnarray*}
l_{2} &=&u_{2}=0,l_{1}+u_{1}=1 \\
a &=&0,b_{2}=1,b_{1}=l_{1} \\
d &=&1,e_{2}=0,e_{1}-u_{1}=1\Rightarrow e_{1}=1,u_{1}=0,b_{1}=l_{1}=1,.
\end{eqnarray*}%
Since $\alpha \left( x_{2};1,0,0,0\right) \equiv a+b_{1}\equiv 1$ $,$ we get%
\begin{equation*}
+B(gx_{1}\otimes x_{1};X_{1}X_{2},gx_{1})X_{2}\otimes gx_{1}\otimes
gx_{1}x_{2}.
\end{equation*}

Third summand of the left side of the equation gives us%
\begin{eqnarray*}
l_{1} &=&u_{1}=l_{2}=u_{2}=0, \\
a &=&b_{1}=0,b_{2}=1 \\
d &=&e_{1}=1,e_{2}=0.
\end{eqnarray*}%
Since $\alpha \left( x_{1}x_{2};0,0,0,0\right) \equiv 0,$ we get%
\begin{equation*}
B(1_{H}\otimes x_{1};X_{2},gx_{1})X_{2}\otimes gx_{1}\otimes gx_{1}x_{2}.
\end{equation*}%
Fourth summand of the left side of the equation gives us%
\begin{eqnarray*}
l_{1}+u_{1} &=&1,l_{2}+u_{2}=1 \\
a &=&0,b_{1}=l_{1},b_{2}-l_{2}=1\Rightarrow b_{2}=1,l_{2}=0,u_{2}=1 \\
d &=&1,e_{1}-u_{1}=1\Rightarrow e_{1}=1,u_{1}=0,b_{1}=l_{1}=1,e_{2}=u_{2}=1.
\end{eqnarray*}%
Since $\alpha \left( 1_{H};1,0,0,1\right) \equiv a+b_{1}\equiv 1$ we get%
\begin{equation*}
-B\left( x_{1}x_{2}\otimes x_{1};X_{1}X_{2},gx_{1}x_{2}\right) X_{2}\otimes
gx_{1}\otimes gx_{1}x_{2}.
\end{equation*}%
Since there is nothing on the right side, we get%
\begin{gather*}
-B(gx_{2}\otimes x_{1};X_{2},gx_{1}x_{2})+B(gx_{1}\otimes
x_{1};X_{1}X_{2},gx_{1})+ \\
+B(1_{H}\otimes x_{1};X_{2},gx_{1})-B\left( x_{1}x_{2}\otimes
x_{1};X_{1}X_{2},gx_{1}x_{2}\right) =0
\end{gather*}%
which holds in view of the form of the elements.

\subsubsection{Case $X_{1}X_{2}\otimes g\otimes gx_{1}x_{2}$}

First summand of the left side of the equality gives us%
\begin{eqnarray*}
l_{1} &=&u_{1}=0,l_{2}+u_{2}=1 \\
a &=&0,b_{1}=1,b_{2}-l_{2}=1\Rightarrow b_{2}=1,l_{2}=0,u_{2}=1 \\
d &=&1,e_{1}=0,e_{2}=u_{2}=1.
\end{eqnarray*}%
Since $\alpha \left( x_{1};0,0,0,1\right) \equiv 1,$ we get%
\begin{equation*}
-B(gx_{2}\otimes x_{1};X_{1}X_{2},gx_{2})X_{1}X_{2}\otimes g\otimes
gx_{1}x_{2}.
\end{equation*}%
\newline
Second summand of the left side of the equation gives us%
\begin{eqnarray*}
l_{2} &=&u_{2}=0,l_{1}+u_{1}=1 \\
a &=&0,b_{2}=1,b_{1}-l_{1}=1\Rightarrow b_{1}=1,l_{1}=0,u_{1}=1 \\
d &=&1,e_{2}=0,e_{1}=u_{1}=1.
\end{eqnarray*}%
Since $\alpha \left( x_{2};0,0,1,0\right) \equiv e_{2}\equiv 0$ $,$ we get%
\begin{equation*}
-B(gx_{1}\otimes x_{1};X_{1}X_{2},gx_{1})X_{1}X_{2}\otimes g\otimes
gx_{1}x_{2}.
\end{equation*}

Third summand of the left side of the equation gives us%
\begin{eqnarray*}
l_{1} &=&u_{1}=l_{2}=u_{2}=0, \\
a &=&0,b_{1}=b_{2}=1 \\
d &=&1,e_{1}=e_{2}=0.
\end{eqnarray*}%
Since $\alpha \left( x_{1}x_{2};0,0,0,0\right) \equiv 0,$ we get%
\begin{equation*}
B(1_{H}\otimes x_{1};X_{1}X_{2}\otimes g)X_{1}X_{2}\otimes g\otimes
gx_{1}x_{2}.
\end{equation*}%
Fourth summand of the left side of the equation gives us%
\begin{eqnarray*}
l_{1}+u_{1} &=&1,l_{2}+u_{2}=1 \\
a &=&0,b_{1}-l_{1}=1\Rightarrow b_{1}=1,l_{1}=0,u_{1}=1 \\
b_{2}-l_{2} &=&1\Rightarrow b_{2}=1,l_{2}=0,u_{2}=1 \\
d &=&1,e_{1}=u_{1}=1,e_{2}=u_{2}=1.
\end{eqnarray*}%
Since $\alpha \left( 1_{H};0,0,1,1\right) \equiv 1+e_{2}\equiv 0$ we get%
\begin{equation*}
+B\left( x_{1}x_{2}\otimes x_{1};X_{1}X_{2},gx_{1}x_{2}\right)
X_{1}X_{2}\otimes g\otimes gx_{1}x_{2}..
\end{equation*}%
Since there is nothing on the right side, we get%
\begin{gather*}
-B(gx_{2}\otimes x_{1};X_{1}X_{2},gx_{2})-B(gx_{1}\otimes
x_{1};X_{1}X_{2},gx_{1})+ \\
+B(1_{H}\otimes x_{1};X_{1}X_{2}\otimes g)+B\left( x_{1}x_{2}\otimes
x_{1};X_{1}X_{2},gx_{1}x_{2}\right) =0
\end{gather*}%
which holds in view of the form of the elements.

\subsection{$B\left( x_{1}x_{2}\otimes x_{1};GX_{1},gx_{1}x_{2}\right) $}

We deduce that%
\begin{eqnarray*}
a &=&b_{1}=1, \\
d &=&e_{1}=e_{2}=1
\end{eqnarray*}%
and we get%
\begin{gather*}
\left( -1\right) ^{\alpha \left( 1_{H};0,0,0,0\right) }B(x_{1}x_{2}\otimes
x_{1};GX_{1},gx_{1}x_{2})GX_{1}\otimes gx_{1}x_{2}\otimes g \\
\left( -1\right) ^{\alpha \left( 1_{H};1,0,0,0\right) }B(x_{1}x_{2}\otimes
x_{1};GX_{1},gx_{1}x_{2})G\otimes gx_{1}x_{2}\otimes x_{1} \\
\left( -1\right) ^{\alpha \left( 1_{H};0,0,1,0\right) }B(x_{1}x_{2}\otimes
x_{1};GX_{1},gx_{1}x_{2})GX_{1}\otimes gx_{2}\otimes x_{1} \\
\left( -1\right) ^{\alpha \left( 1_{H};1,0,1,0\right) }B(x_{1}x_{2}\otimes
x_{1};GX_{1},gx_{1}x_{2})GX_{1}^{1-l_{1}}\otimes gx_{2}\otimes
g^{+l_{11}}x_{1}^{1+1}=0 \\
\left( -1\right) ^{\alpha \left( 1_{H};0,0,0,1\right) }B(x_{1}x_{2}\otimes
x_{1};GX_{1},gx_{1}x_{2})GX_{1}\otimes gx_{1}\otimes x_{2} \\
\left( -1\right) ^{\alpha \left( 1_{H};1,0,0,1\right) }B(x_{1}x_{2}\otimes
x_{1};GX_{1},gx_{1}x_{2})G\otimes gx_{1}\otimes gx_{1}x_{2} \\
\left( -1\right) ^{\alpha \left( 1_{H};0,0,1,1\right) }B(x_{1}x_{2}\otimes
x_{1};GX_{1},gx_{1}x_{2})GX_{1}\otimes g\otimes gx_{1}x_{2} \\
\left( -1\right) ^{\alpha \left( 1_{H};1,0,1,1\right) }B(x_{1}x_{2}\otimes
x_{1};GX_{1},gx_{1}x_{2})GX_{1}^{1-l_{1}}\otimes g\otimes
g^{+l_{1}+1}x_{1}^{1+1}x_{2}=0
\end{gather*}

\subsubsection{Case $G\otimes gx_{1}x_{2}\otimes x_{1}$}

We have to consider only the first and the fourth summand of the left side
of the equation. First summand gives us%
\begin{eqnarray*}
l_{1} &=&u_{1}=l_{2}=u_{2}=0 \\
a &=&1,b_{1}=b_{2}=0, \\
d &=&e_{1}=e_{2}=1.
\end{eqnarray*}%
Since $\alpha \left( x_{1};0,0,0,0\right) \equiv a+b_{1}+b_{2}\equiv 1,$ we
get%
\begin{equation*}
-B(gx_{2}\otimes x_{1};G,gx_{1}x_{2})G\otimes gx_{1}x_{2}\otimes x_{1}.
\end{equation*}%
Fourth summand is giving us%
\begin{eqnarray*}
l_{1}+u_{1} &=&1,l_{2}=u_{2}=0, \\
a &=&1,b_{2}=0,b_{1}=l_{1}, \\
d &=&1,e_{1}-u_{1}=1\Rightarrow e_{1}=1,u_{1}=0,b_{1}=l_{1}=1,e_{2}=1.
\end{eqnarray*}%
Since $\alpha \left( 1_{H};1,0,0,0\right) \equiv b_{2}=0$ $,$ we get%
\begin{equation*}
B(x_{1}x_{2}\otimes x_{1};GX_{1},gx_{1}x_{2})G\otimes gx_{1}x_{2}\otimes
x_{1}.
\end{equation*}%
By considering also the right side, we obtain%
\begin{equation*}
-B(gx_{2}\otimes x_{1};G,gx_{1}x_{2})+B(x_{1}x_{2}\otimes
x_{1};GX_{1},gx_{1}x_{2})-B(x_{1}x_{2}\otimes 1_{H};G,gx_{1}x_{2})=0
\end{equation*}%
which holds in view of the form of the elements.

\subsubsection{Case $GX_{1}\otimes gx_{2}\otimes x_{1}$}

We have to consider only the first and the fourth summand of the left side
of the equation. First summand gives us%
\begin{eqnarray*}
l_{1} &=&u_{1}=l_{2}=u_{2}=0 \\
a &=&b_{1}=1,b_{2}=0, \\
d &=&e_{1}=e_{2}=1.
\end{eqnarray*}%
Since $\alpha \left( x_{1};0,0,0,0\right) \equiv a+b_{1}+b_{2}\equiv 0,$ we
get%
\begin{equation*}
B(gx_{2}\otimes x_{1};GX_{1},gx_{2})GX_{1}\otimes gx_{2}\otimes x_{1}.
\end{equation*}%
Fourth summand is giving us%
\begin{eqnarray*}
l_{1}+u_{1} &=&1,l_{2}=u_{2}=0, \\
a &=&1,b_{2}=0,b_{1}-l_{1}=1\Rightarrow b_{1}=1,l_{1}=0,u_{1}=1 \\
d &=&1,e_{1}=u_{1}=1,e_{2}=1.
\end{eqnarray*}%
Since $\alpha \left( 1_{H};0,0,1,0\right) \equiv e_{2}+\left(
a+b_{1}+b_{2}\right) \equiv 1$ $,$ we get%
\begin{equation*}
-B(x_{1}x_{2}\otimes x_{1};GX_{1},gx_{1}x_{2})GX_{1}\otimes gx_{2}\otimes
x_{1}.
\end{equation*}%
By considering also the right side, we obtain%
\begin{equation*}
B(gx_{2}\otimes x_{1};GX_{1},gx_{2})-B(x_{1}x_{2}\otimes
x_{1};GX_{1},gx_{1}x_{2})-B(x_{1}x_{2}\otimes 1_{H};GX_{1},gx_{2})=0
\end{equation*}%
which holds in view of the form of the elements.

\subsubsection{Case $GX_{1}\otimes gx_{1}\otimes x_{2}$}

We have to consider only the second and the fourth summands of the left side
of the equality. Second summand gives us%
\begin{eqnarray*}
l_{1} &=&u_{1}=l_{2}=u_{2}=0, \\
a &=&b_{1}=1,b_{2}=0, \\
d &=&e_{1}=1,e_{2}=0.
\end{eqnarray*}%
Since $\alpha \left( x_{2};0,0,0,0\right) \equiv a+b_{1}+b_{2}=0$ we get%
\begin{equation*}
-B(gx_{1}\otimes x_{1};GX_{1},gx_{1})GX_{1}\otimes gx_{1}\otimes x_{2}.
\end{equation*}%
Fourth summand gives us%
\begin{eqnarray*}
l_{1} &=&u_{1}=0,l_{2}+u_{2}=1 \\
a &=&b_{1}=1,b_{2}=l_{2} \\
d &=&e_{1}=1,e_{2}=u_{2}
\end{eqnarray*}%
Since $\alpha \left( 1_{H};0,0,0,1\right) \equiv a+b_{1}+b_{2}\equiv 0$ and $%
\alpha \left( 1_{H};0,1,0,0\right) \equiv 0,$ we get%
\begin{equation*}
\left[ B(x_{1}x_{2}\otimes x_{1};GX_{1},gx_{1}x_{2})+B(x_{1}x_{2}\otimes
x_{1};GX_{1}X_{2},gx_{1})\right] GX_{1}\otimes gx_{1}\otimes x_{2}.
\end{equation*}%
Since there is nothing in the right side, we get%
\begin{equation*}
-B(gx_{1}\otimes x_{1};GX_{1},gx_{1})+B(x_{1}x_{2}\otimes
x_{1};GX_{1},gx_{1}x_{2})+B(x_{1}x_{2}\otimes x_{1};GX_{1}X_{2},gx_{1})=0
\end{equation*}%
which holds in view of the form of the elements.

\subsubsection{Case $G\otimes gx_{1}\otimes gx_{1}x_{2}$}

First summand of the left side of the equality gives us%
\begin{eqnarray*}
l_{1} &=&u_{1}=0,l_{2}+u_{2}=1 \\
a &=&1,b_{1}=0,b_{2}=l_{2}, \\
d &=&1,e_{1}=1,e_{2}=u_{2}.
\end{eqnarray*}
Since $\alpha \left( x_{1};0,0,0,1\right) \equiv 1$ and $\alpha \left(
x_{1};0,1,0,0\right) \equiv a+b_{1}+b_{2}+1\equiv 1,$ we get%
\begin{equation*}
\left[ -B(gx_{2}\otimes x_{1};G,gx_{1}x_{2})-B(gx_{2}\otimes
x_{1};GX_{2},gx_{1})\right] G\otimes gx_{1}\otimes gx_{1}x_{2}.
\end{equation*}%
Second summand of the left side of the equation gives us%
\begin{eqnarray*}
l_{2} &=&u_{2}=0,l_{1}+u_{1}=1 \\
a &=&1,b_{2}=0,b_{1}=l_{1} \\
d &=&1,e_{2}=0,e_{1}-u_{1}=1\Rightarrow e_{1}=1,u_{1}=0,b_{1}=l_{1}=1,.
\end{eqnarray*}%
Since $\alpha \left( x_{2};1,0,0,0\right) \equiv a+b_{1}=0,$ we get%
\begin{equation*}
-B(gx_{1}\otimes x_{1};GX_{1},gx_{1})G\otimes gx_{1}\otimes gx_{1}x_{2}.
\end{equation*}

Third summand of the left side of the equation gives us%
\begin{eqnarray*}
l_{1} &=&u_{1}=l_{2}=u_{2}=0, \\
a &=&1,b_{1}=b_{2}=0, \\
d &=&e_{1}=1,e_{2}=0.
\end{eqnarray*}%
Since $\alpha \left( x_{1}x_{2};0,0,0,0\right) \equiv 0,$ we get%
\begin{equation*}
B(1_{H}\otimes x_{1};G,gx_{1})G\otimes gx_{1}\otimes gx_{1}x_{2}.
\end{equation*}%
Fourth summand of the left side of the equation gives us%
\begin{eqnarray*}
l_{1}+u_{1} &=&1,l_{2}+u_{2}=1 \\
a &=&1,b_{1}=l_{1},b_{2}=l_{2}, \\
d &=&1,e_{1}-u_{1}=1\Rightarrow e_{1}=1,u_{1}=0,b_{1}=l_{1}=1,e_{2}=u_{2}.
\end{eqnarray*}%
Since $\alpha \left( 1_{H};1,0,0,1\right) \equiv $ $a+b_{1}\equiv 0$ and $%
\alpha \left( 1_{H};1,1,0,0\right) $ $\equiv 1+b_{2}\equiv 0$ we get%
\begin{equation*}
\left[ B\left( x_{1}x_{2}\otimes x_{1};GX_{1},gx_{1}x_{2}\right) +B\left(
x_{1}x_{2}\otimes x_{1};GX_{1}X_{2},gx_{1}\right) \right] G\otimes
gx_{1}\otimes gx_{1}x_{2}.
\end{equation*}%
Since there is nothing on the right side, we get%
\begin{gather*}
-B(gx_{2}\otimes x_{1};G,gx_{1}x_{2})-B(gx_{2}\otimes
x_{1};GX_{2},gx_{1})-B(gx_{1}\otimes x_{1};GX_{1},gx_{1}) \\
+B(1_{H}\otimes x_{1};G,gx_{1})+B\left( x_{1}x_{2}\otimes
x_{1};GX_{1},gx_{1}x_{2}\right) +B\left( x_{1}x_{2}\otimes
x_{1};GX_{1}X_{2},gx_{1}\right) =0
\end{gather*}%
which holds in view of the form of the elements.%
\begin{gather*}
-B(gx_{2}\otimes x_{1};G,gx_{1}x_{2})=\left[ +B(x_{1}x_{2}\otimes
1_{H};G,gx_{1}x_{2})+2B(x_{1}x_{2}\otimes 1_{H};GX_{1},gx_{2})\right] \\
-B(gx_{2}\otimes x_{1};GX_{2},gx_{1})=\left[
\begin{array}{c}
+B(x_{1}x_{2}\otimes 1_{H};GX_{2},gx_{1}) \\
-2B(x_{1}x_{2}\otimes 1_{H};G,gx_{1}x_{2})-2B(x_{1}x_{2}\otimes
1_{H};GX_{1},gx_{2})%
\end{array}%
\right] \\
-B(gx_{1}\otimes x_{1};GX_{1},gx_{1})=\left[ -2B(x_{1}x_{2}\otimes
1_{H};G,gx_{1}x_{2})+2B(x_{1}x_{2}\otimes 1_{H};GX_{2},gx_{1})\right] \\
+B(1_{H}\otimes x_{1};G,gx_{1})=+\left[ B(x_{1}x_{2}\otimes
1_{H};G,gx_{1}x_{2})-B(x_{1}x_{2}\otimes 1_{H};GX_{2},gx_{1})\right] \\
+B\left( x_{1}x_{2}\otimes x_{1};GX_{1},gx_{1}x_{2}\right)
=-2B(x_{1}x_{2}\otimes 1_{H};GX_{1},gx_{2}) \\
+B\left( x_{1}x_{2}\otimes x_{1};GX_{1}X_{2},gx_{1}\right) =+\left[
\begin{array}{c}
+2B(x_{1}x_{2}\otimes 1_{H};G,gx_{1}x_{2}) \\
-2B(x_{1}x_{2}\otimes 1_{H};GX_{2},gx_{1})+2B(x_{1}x_{2}\otimes
1_{H};GX_{1},gx_{2})%
\end{array}%
\right]
\end{gather*}

\subsubsection{Case $GX_{1}\otimes g\otimes gx_{1}x_{2}$}

First summand of the left side of the equality gives us%
\begin{eqnarray*}
l_{1} &=&u_{1}=0,l_{2}+u_{2}=1 \\
a &=&1,b_{1}=1,b_{2}=l_{2}, \\
d &=&1,e_{1}=0,e_{2}=u_{2}.
\end{eqnarray*}%
Since $\alpha \left( x_{1};0,0,0,1\right) \equiv 1$ and $\alpha \left(
x_{1};0,1,0,0\right) \equiv a+b_{1}+b_{2}+1\equiv 0,$ we get%
\begin{equation*}
\left[ -B(gx_{2}\otimes x_{1};GX_{1},gx_{2})+B(gx_{2}\otimes
x_{1};GX_{1}X_{2},g)\right] GX_{1}\otimes g\otimes gx_{1}x_{2}.
\end{equation*}%
Second summand of the left side of the equation gives us%
\begin{eqnarray*}
l_{2} &=&u_{2}=0,l_{1}+u_{1}=1 \\
a &=&1,b_{2}=0,b_{1}-l_{1}=1\Rightarrow b_{1}=1,l_{1}=0,u_{1}=1 \\
d &=&1,e_{2}=0,e_{1}=u_{1}=1.
\end{eqnarray*}%
Since $\alpha \left( x_{2};0,0,1,0\right) \equiv e_{2}=0,$ we get%
\begin{equation*}
-B(gx_{1}\otimes x_{1};GX_{1},gx_{1})GX_{1}\otimes g\otimes gx_{1}x_{2}.
\end{equation*}

Third summand of the left side of the equation gives us%
\begin{eqnarray*}
l_{1} &=&u_{1}=l_{2}=u_{2}=0, \\
a &=&b_{1}=1,b_{2}=0, \\
d &=&1,e_{1}=e_{2}=0.
\end{eqnarray*}%
Since $\alpha \left( x_{1}x_{2};0,0,0,0\right) \equiv 0,$ we get%
\begin{equation*}
B(1_{H}\otimes x_{1};GX_{1},g)GX_{1}\otimes g\otimes gx_{1}x_{2}.
\end{equation*}%
Fourth summand of the left side of the equation gives us%
\begin{eqnarray*}
l_{1}+u_{1} &=&1,l_{2}+u_{2}=1 \\
a &=&1,b_{1}-l_{1}=1\Rightarrow b_{1}=1,l_{1}=0,u_{1}=1,b_{2}=l_{2}, \\
d &=&1,e_{1}=u_{1}=1,e_{2}=u_{2}.
\end{eqnarray*}%
Since $\alpha \left( 1_{H};0,0,1,1\right) \equiv $ $1+e_{2}\equiv 0$ and $%
\alpha \left( 1_{H};0,1,1,0\right) $ $\equiv 1+b_{2}\equiv 0$ we get%
\begin{equation*}
\left[ B\left( x_{1}x_{2}\otimes x_{1};GX_{1},gx_{1}x_{2}\right) +B\left(
x_{1}x_{2}\otimes x_{1};GX_{1}X_{2},gx_{1}\right) \right] GX_{1}\otimes
g\otimes gx_{1}x_{2}.
\end{equation*}%
Since there is nothing on the right side, we get%
\begin{gather*}
-B(gx_{2}\otimes x_{1};GX_{1},gx_{2})+B(gx_{2}\otimes
x_{1};GX_{1}X_{2},g)-B(gx_{1}\otimes x_{1};GX_{1},gx_{1})+ \\
+B(1_{H}\otimes x_{1};GX_{1},g)+B\left( x_{1}x_{2}\otimes
x_{1};GX_{1},gx_{1}x_{2}\right) +B\left( x_{1}x_{2}\otimes
x_{1};GX_{1}X_{2},gx_{1}\right) =0
\end{gather*}%
which holds in view of the form of the elements.%
\begin{gather*}
-B(gx_{2}\otimes x_{1};GX_{1},gx_{2})=B(x_{1}x_{2}\otimes
1_{H};GX_{1},gx_{2}) \\
+B(gx_{2}\otimes x_{1};GX_{1}X_{2},g)=\left[
\begin{array}{c}
-B(x_{1}x_{2}\otimes 1_{H};G,gx_{1}x_{2}) \\
+B(x_{1}x_{2}\otimes 1_{H};GX_{2},gx_{1})-B(x_{1}x_{2}\otimes
1_{H};GX_{1},gx_{2})%
\end{array}%
\right] \\
-B(gx_{1}\otimes x_{1};GX_{1},gx_{1})=\left[ -2B(x_{1}x_{2}\otimes
1_{H};G,gx_{1}x_{2})+2B(x_{1}x_{2}\otimes 1_{H};GX_{2},gx_{1})\right] \\
+B(1_{H}\otimes x_{1};GX_{1},g)=\left[ +B(x_{1}x_{2}\otimes
1_{H};G,gx_{1}x_{2})-B(x_{1}x_{2}\otimes 1_{H};GX_{2},gx_{1})\right] \\
+B\left( x_{1}x_{2}\otimes x_{1};GX_{1},gx_{1}x_{2}\right)
=-2B(x_{1}x_{2}\otimes 1_{H};GX_{1},gx_{2}) \\
+B\left( x_{1}x_{2}\otimes x_{1};GX_{1}X_{2},gx_{1}\right) =\left[
\begin{array}{c}
+2B(x_{1}x_{2}\otimes 1_{H};G,gx_{1}x_{2}) \\
-2B(x_{1}x_{2}\otimes 1_{H};GX_{2},gx_{1})+2B(x_{1}x_{2}\otimes
1_{H};GX_{1},gx_{2})%
\end{array}%
\right]
\end{gather*}

\subsection{$B\left( x_{1}x_{2}\otimes x_{1};GX_{2},gx_{1}x_{2}\right) $}

We deduce that%
\begin{eqnarray*}
a &=&b_{2}=1, \\
d &=&e_{1}=e_{2}=1
\end{eqnarray*}%
and we get%
\begin{gather*}
\left( -1\right) ^{\alpha \left( 1_{H};0,0,0,0\right) }B\left(
x_{1}x_{2}\otimes x_{1};GX_{2},gx_{1}x_{2}\,\right) GX_{2}\otimes
gx_{1}x_{2}\otimes g+ \\
\left( -1\right) ^{\alpha \left( 1_{H};0,1,0,0\right) }B\left(
x_{1}x_{2}\otimes x_{1};GX_{2},gx_{1}x_{2}\,\right) G\otimes
gx_{1}x_{2}\otimes x_{2}+ \\
\left( -1\right) ^{\alpha \left( 1_{H};0,0,1,0\right) }B\left(
x_{1}x_{2}\otimes x_{1};GX_{2},gx_{1}x_{2}\,\right) GX_{2}\otimes
gx_{2}\otimes x_{1}+ \\
\left( -1\right) ^{\alpha \left( 1_{H};0,1,1,0\right) }B\left(
x_{1}x_{2}\otimes x_{1};GX_{2},gx_{1}x_{2}\,\right) G\otimes gx_{2}\otimes
gx_{1}x_{2}+ \\
\left( -1\right) ^{\alpha \left( 1_{H};0,0,0,1\right) }B\left(
x_{1}x_{2}\otimes x_{1};GX_{2},gx_{1}x_{2}\,\right) GX_{2}\otimes
gx_{1}\otimes x_{2}+ \\
\left( -1\right) ^{\alpha \left( 1_{H};0,1,0,1\right) }B\left(
x_{1}x_{2}\otimes x_{1};GX_{2},gx_{1}x_{2}\,\right) G\otimes gx_{1}\otimes
g^{l_{2}}x_{2}^{1+1}=0 \\
\left( -1\right) ^{\alpha \left( 1_{H};0,0,1,1\right) }B\left(
x_{1}x_{2}\otimes x_{1};GX_{2},gx_{1}x_{2}\,\right) GX_{2}\otimes g\otimes
gx_{1}x_{2} \\
\left( -1\right) ^{\alpha \left( 1_{H};0,1,1,1\right) }B\left(
x_{1}x_{2}\otimes x_{1};GX_{2},gx_{1}x_{2}\,\right) GX_{2}^{1-l_{2}}\otimes
g\otimes g^{l_{2}+1}x_{1}x_{2}^{1+1}=0
\end{gather*}

\subsubsection{Case $G\otimes gx_{1}x_{2}\otimes x_{2}$}

We have to consider only the second and the fourth summands of the left side
of the equality. Second summand gives us%
\begin{eqnarray*}
l_{1} &=&u_{1}=l_{2}=u_{2}=0, \\
a &=&1,b_{1}=b_{2}=0, \\
d &=&e_{1}=e_{2}=1.
\end{eqnarray*}%
Since $\alpha \left( x_{2};0,0,0,0\right) \equiv a+b_{1}+b_{2}=1$ we get%
\begin{equation*}
B(gx_{1}\otimes x_{1};G,gx_{1}x_{2})G\otimes gx_{1}x_{2}\otimes x_{2}.
\end{equation*}%
Fourth summand gives us%
\begin{eqnarray*}
l_{1} &=&u_{1}=0,l_{2}+u_{2}=1 \\
a &=&1,b_{1}=0,b_{2}=l_{2} \\
d &=&e_{1}=1,e_{2}-u_{2}=1\Rightarrow e_{2}=1,u_{2}=0,b_{2}=l_{2}=1
\end{eqnarray*}%
Since $\alpha \left( 1_{H};0,1,0,0\right) \equiv 0,$ we get%
\begin{equation*}
+B(x_{1}x_{2}\otimes x_{1};GX_{2},gx_{1}x_{2})G\otimes gx_{1}x_{2}\otimes
x_{2}.
\end{equation*}%
Since there is nothing in the right side, we get%
\begin{equation*}
B(gx_{1}\otimes x_{1};G,gx_{1}x_{2})+B(x_{1}x_{2}\otimes
x_{1};GX_{2},gx_{1}x_{2})=0
\end{equation*}%
which holds in view of the form of the elements.

\subsubsection{Case $GX_{2}\otimes gx_{2}\otimes x_{1}$}

We have to consider only the first and the fourth summand of the left side
of the equation. First summand gives us%
\begin{eqnarray*}
l_{1} &=&u_{1}=l_{2}=u_{2}=0 \\
a &=&b_{2}=1,b_{1}=0, \\
d &=&e_{2}=1,e_{1}=0.
\end{eqnarray*}%
Since $\alpha \left( x_{1};0,0,0,0\right) \equiv a+b_{1}+b_{2}\equiv 0,$ we
get%
\begin{equation*}
B(gx_{2}\otimes x_{1};GX_{2},gx_{2})GX_{2}\otimes gx_{2}\otimes x_{1}.
\end{equation*}%
Fourth summand is giving us%
\begin{eqnarray*}
l_{1}+u_{1} &=&1,l_{2}=u_{2}=0, \\
a &=&1,b_{2}=1,b_{1}=l_{1} \\
d &=&e_{2}=1,e_{1}=u_{1}.
\end{eqnarray*}%
Since $\alpha \left( 1_{H};0,0,1,0\right) \equiv e_{2}+\left(
a+b_{1}+b_{2}\right) \equiv 1$ and $\alpha \left( 1_{H};1,0,0,0\right)
\equiv b_{2}\equiv 1$ we get%
\begin{equation*}
\left[ -B(x_{1}x_{2}\otimes x_{1};GX_{2},gx_{1}x_{2})-B(x_{1}x_{2}\otimes
x_{1};GX_{1}X_{2},gx_{2})\right] GX_{2}\otimes gx_{2}\otimes x_{1}.
\end{equation*}%
By considering also the right side, we obtain%
\begin{gather*}
B(gx_{2}\otimes x_{1};GX_{2},gx_{2})-B(x_{1}x_{2}\otimes
x_{1};GX_{2},gx_{1}x_{2}) \\
-B(x_{1}x_{2}\otimes x_{1};GX_{1}X_{2},gx_{2})-B(x_{1}x_{2}\otimes
1_{H};GX_{2},gx_{2})=0
\end{gather*}%
which holds in view of the form of the elements.

\subsubsection{Case $G\otimes gx_{2}\otimes gx_{1}x_{2}$}

First summand of the left side of the equality gives us%
\begin{eqnarray*}
l_{1} &=&u_{1}=0,l_{2}+u_{2}=1 \\
a &=&1,b_{1}=0,b_{2}=l_{2} \\
d &=&1,e_{1}=0,e_{2}-u_{2}=1\Rightarrow e_{2}=1,u_{2}=0,b_{2}=l_{2}=1.
\end{eqnarray*}%
Since $\alpha \left( x_{1};0,1,0,0\right) \equiv a+b_{1}+b_{2}+1\equiv 1,$
we get%
\begin{equation*}
-B(gx_{2}\otimes x_{1};GX_{2},gx_{2})G\otimes gx_{2}\otimes gx_{1}x_{2}.
\end{equation*}%
\newline
Second summand of the left side of the equation gives us%
\begin{eqnarray*}
l_{2} &=&u_{2}=0,l_{1}+u_{1}=1 \\
a &=&1,b_{2}=0,b_{1}=l_{1} \\
d &=&1,e_{2}=1,e_{1}=u_{1}.
\end{eqnarray*}%
Since $\alpha \left( x_{2};0,0,1,0\right) \equiv e_{2}\equiv 1$ and $\alpha
\left( x_{2};1,0,0,0\right) \equiv a+b_{1}\equiv 0,$ we get%
\begin{equation*}
\left[ +B(gx_{1}\otimes x_{1};G,gx_{1}x_{2})-B(gx_{1}\otimes
x_{1};GX_{1},gx_{2})\right] G\otimes gx_{2}\otimes gx_{1}x_{2}.
\end{equation*}

Third summand of the left side of the equation gives us%
\begin{eqnarray*}
l_{1} &=&u_{1}=l_{2}=u_{2}=0, \\
a &=&1,b_{1}=b_{2}=0 \\
d &=&e_{2}=1,e_{1}=0.
\end{eqnarray*}%
Since $\alpha \left( x_{1}x_{2};0,0,0,0\right) \equiv 0,$ we get%
\begin{equation*}
B(1_{H}\otimes x_{1};G\otimes gx_{2})G\otimes gx_{2}\otimes gx_{1}x_{2}.
\end{equation*}%
Fourth summand of the left side of the equation gives us%
\begin{eqnarray*}
l_{1}+u_{1} &=&1,l_{2}+u_{2}=1 \\
a &=&1,b_{1}=l_{1},b_{2}=l_{2} \\
d &=&1,e_{1}=u_{1},e_{2}-u_{2}=1\Rightarrow e_{2}=1,u_{2}=0,b_{2}=l_{2}=1.
\end{eqnarray*}%
Since $\alpha \left( 1_{H};0,1,1,0\right) \equiv e_{2}+a+b_{1}+b_{2}+1\equiv
0$ and $\alpha \left( 1_{H};1,1,0,0\right) \equiv 1+b_{2}\equiv 0$ we get%
\begin{equation*}
\left[ +B\left( x_{1}x_{2}\otimes x_{1};GX_{2},gx_{1}x_{2}\right) +B\left(
x_{1}x_{2}\otimes x_{1};GX_{1}X_{2},gx_{2}\right) \right] G\otimes
gx_{2}\otimes gx_{1}x_{2}.
\end{equation*}%
Since there is nothing on the right side, we get%
\begin{gather*}
-B(gx_{2}\otimes x_{1};GX_{2},gx_{2})+B(gx_{1}\otimes
x_{1};G,gx_{1}x_{2})-B(gx_{1}\otimes x_{1};GX_{1},gx_{2}) \\
+B(1_{H}\otimes x_{1};G\otimes gx_{2})+B\left( x_{1}x_{2}\otimes
x_{1};GX_{2},gx_{1}x_{2}\right) +B\left( x_{1}x_{2}\otimes
x_{1};GX_{1}X_{2},gx_{2}\right) =0
\end{gather*}%
which holds in view of the form of the elements.%
\begin{gather*}
B(gx_{2}\otimes x_{1};GX_{2},gx_{2})=+B(x_{1}x_{2}\otimes
1_{H};GX_{2},gx_{2}) \\
+B(gx_{1}\otimes x_{1};G,gx_{1}x_{2})=+2B(x_{1}x_{2}\otimes
1_{H};GX_{2},gx_{2}) \\
-B(gx_{1}\otimes x_{1};GX_{1},gx_{2})= \\
+B(1_{H}\otimes x_{1};G\otimes gx_{2})=-B(x_{1}x_{2}\otimes
1_{H};GX_{2},gx_{2}) \\
+B\left( x_{1}x_{2}\otimes x_{1};GX_{2},gx_{1}x_{2}\right)
=-2B(x_{1}x_{2}\otimes 1_{H};GX_{2},gx_{2}) \\
+B\left( x_{1}x_{2}\otimes x_{1};GX_{1}X_{2},gx_{2}\right) =0
\end{gather*}

\subsubsection{Case $GX_{2}\otimes gx_{1}\otimes x_{2}$}

We have to consider only the second and the fourth summands of the left side
of the equality. Second summand gives us%
\begin{eqnarray*}
l_{1} &=&u_{1}=l_{2}=u_{2}=0, \\
a &=&b_{2}=1,b_{1}=0, \\
d &=&e_{1}=1,e_{2}=0.
\end{eqnarray*}%
Since $\alpha \left( x_{2};0,0,0,0\right) \equiv a+b_{1}+b_{2}=0$ we get%
\begin{equation*}
-B(gx_{1}\otimes x_{1};GX_{2},gx_{1})GX_{2}\otimes gx_{1}\otimes x_{2}.
\end{equation*}%
Fourth summand gives us%
\begin{eqnarray*}
l_{1} &=&u_{1}=0,l_{2}+u_{2}=1 \\
a &=&1,b_{1}=0,b_{2}-l_{2}=1\Rightarrow b_{2}=1,l_{2}=0,u_{2}=1 \\
d &=&e_{1}=1,e_{2}=u_{2}=1.
\end{eqnarray*}%
Since $\alpha \left( 1_{H};0,0,0,1\right) \equiv a+b_{1}+b_{2}\equiv 0,$ we
get%
\begin{equation*}
+B(x_{1}x_{2}\otimes x_{1};GX_{2},gx_{1}x_{2})GX_{2}\otimes gx_{1}\otimes
x_{2}.
\end{equation*}%
Since there is nothing in the right side, we get%
\begin{equation*}
-B(gx_{1}\otimes x_{1};GX_{2},gx_{1})+B(x_{1}x_{2}\otimes
x_{1};GX_{2},gx_{1}x_{2})=0
\end{equation*}%
which holds in view of the form of the elements.

\subsubsection{Case $GX_{2}\otimes g\otimes gx_{1}x_{2}$}

First summand of the left side of the equality gives us%
\begin{eqnarray*}
l_{1} &=&u_{1}=0,l_{2}+u_{2}=1 \\
a &=&1,b_{1}=0,b_{2}-l_{2}=1\Rightarrow b_{2}=1,l_{2}=0,u_{2}=1 \\
d &=&1,e_{1}=0,e_{2}=u_{2}=1.
\end{eqnarray*}%
Since $\alpha \left( x_{1};0,0,0,1\right) \equiv 1,$ we get%
\begin{equation*}
-B(gx_{2}\otimes x_{1};GX_{2},gx_{2})GX_{2}\otimes g\otimes gx_{1}x_{2}.
\end{equation*}%
\newline
Second summand of the left side of the equation gives us%
\begin{eqnarray*}
l_{2} &=&u_{2}=0,l_{1}+u_{1}=1 \\
a &=&b_{2}=1,0,b_{1}=l_{1} \\
d &=&1,e_{2}=0,e_{1}=u_{1}.
\end{eqnarray*}%
Since $\alpha \left( x_{2};0,0,1,0\right) \equiv e_{2}\equiv 0$ and $\alpha
\left( x_{2};1,0,0,0\right) \equiv a+b_{1}\equiv 0,$ we get%
\begin{equation*}
\left[ -B(gx_{1}\otimes x_{1};GX_{2},gx_{1})-B(gx_{1}\otimes
x_{1};GX_{1}X_{2},gx_{2})\right] GX_{2}\otimes g\otimes gx_{1}x_{2}.
\end{equation*}

Third summand of the left side of the equation gives us%
\begin{eqnarray*}
l_{1} &=&u_{1}=l_{2}=u_{2}=0, \\
a &=&1,b_{1}=b_{2}=0 \\
d &=&e_{2}=1,e_{1}=0.
\end{eqnarray*}%
Since $\alpha \left( x_{1}x_{2};0,0,0,0\right) \equiv 0,$ we get%
\begin{equation*}
B(1_{H}\otimes x_{1};GX_{2}\otimes g)GX_{2}\otimes g\otimes gx_{1}x_{2}.
\end{equation*}%
Fourth summand of the left side of the equation gives us%
\begin{eqnarray*}
l_{1}+u_{1} &=&1,l_{2}+u_{2}=1 \\
a &=&1,b_{1}=l_{1},b_{2}-l_{2}=1\Rightarrow b_{2}=1,l_{2}=0,u_{2}=1, \\
d &=&1,e_{1}=u_{1},e_{2}=u_{2}=1.
\end{eqnarray*}%
Since $\alpha \left( 1_{H};0,0,1,1\right) \equiv 1+e_{2}\equiv 0$ and $%
\alpha \left( 1_{H};1,0,0,1\right) \equiv a+b_{1}\equiv 0$ we get%
\begin{equation*}
\left[ +B\left( x_{1}x_{2}\otimes x_{1};GX_{2},gx_{1}x_{2}\right) +B\left(
x_{1}x_{2}\otimes x_{1};GX_{1}X_{2},gx_{2}\right) \right] GX_{2}\otimes
g\otimes gx_{1}x_{2}.
\end{equation*}%
Since there is nothing on the right side, we get%
\begin{gather*}
-B(gx_{2}\otimes x_{1};GX_{2},gx_{2})-B(gx_{1}\otimes
x_{1};GX_{2},gx_{1})-B(gx_{1}\otimes x_{1};GX_{1}X_{2},gx_{2}) \\
+B(1_{H}\otimes x_{1};GX_{2}\otimes g)+B\left( x_{1}x_{2}\otimes
x_{1};GX_{2},gx_{1}x_{2}\right) +B\left( x_{1}x_{2}\otimes
x_{1};GX_{1}X_{2},gx_{2}\right) =0
\end{gather*}%
which holds in view of the form of the elements.

\subsection{$B\left( x_{1}x_{2}\otimes x_{1};GX_{1}X_{2},gx_{1}\right) $}

We get%
\begin{eqnarray*}
a &=&b_{1}=b_{2}=1 \\
d &=&e_{1}=1,e_{2}=0
\end{eqnarray*}%
so that we obtain%
\begin{gather*}
\left( -1\right) ^{\alpha \left( 1_{H};0,0,0,0\right) }B\left(
x_{1}x_{2}\otimes x_{1};GX_{1}X_{2},gx_{1}\right) GX_{1}X_{2}\otimes
gx_{1}\otimes g+ \\
\left( -1\right) ^{\alpha \left( 1_{H};1,0,0,0\right) }B\left(
x_{1}x_{2}\otimes x_{1};GX_{1}X_{2},gx_{1}\right) GX_{2}\otimes
gx_{1}\otimes x_{1}+ \\
\left( -1\right) ^{\alpha \left( 1_{H};0,1,0,0\right) }B\left(
x_{1}x_{2}\otimes x_{1};GX_{1}X_{2},gx_{1}\right) GX_{1}\otimes
gx_{1}\otimes x_{2}+ \\
\left( -1\right) ^{\alpha \left( 1_{H};1,1,0,0\right) }B\left(
x_{1}x_{2}\otimes x_{1};GX_{1}X_{2},gx_{1}\right) G\otimes gx_{1}\otimes
gx_{1}x_{2}+ \\
\left( -1\right) ^{\alpha \left( 1_{H};0,0,1,0\right) }B\left(
x_{1}x_{2}\otimes x_{1};GX_{1}X_{2},gx_{1}\right) GX_{1}X_{2}\otimes
g\otimes x_{1}+ \\
\left( -1\right) ^{\alpha \left( 1_{H};1,0,1,0\right) }B\left(
x_{1}x_{2}\otimes x_{1};GX_{1}X_{2},gx_{1}\right)
GX_{1}^{1-l_{1}}X_{2}\otimes g\otimes g^{l_{1}}x_{1}^{1+1}=0 \\
\left( -1\right) ^{\alpha \left( 1_{H};0,1,1,0\right) }B\left(
x_{1}x_{2}\otimes x_{1};GX_{1}X_{2},gx_{1}\right) GX_{1}\otimes g\otimes
gx_{1}x_{2}+ \\
\left( -1\right) ^{\alpha \left( 1_{H};1,1,1,0\right) }B\left(
x_{1}x_{2}\otimes x_{1};GX_{1}X_{2},gx_{1}\right) GX_{1}^{1-l_{1}}\otimes
g\otimes g^{l_{1}+1}x_{1}^{1+1}x_{2}=0
\end{gather*}

\subsubsection{Case $GX_{2}\otimes gx_{1}\otimes x_{1}$}

We have to consider only the first and the fourth summand of the left side
of the equation. First summand gives us%
\begin{eqnarray*}
l_{1} &=&u_{1}=l_{2}=u_{2}=0 \\
a &=&b_{2}=1,b_{1}=0, \\
d &=&e_{1}=1,e_{2}=0.
\end{eqnarray*}%
Since $\alpha \left( x_{1};0,0,0,0\right) \equiv a+b_{1}+b_{2}\equiv 0,$ we
get%
\begin{equation*}
B(gx_{2}\otimes x_{1};GX_{2},gx_{1})GX_{2}\otimes gx_{1}\otimes x_{1}.
\end{equation*}%
Fourth summand is giving us%
\begin{eqnarray*}
l_{1}+u_{1} &=&1,l_{2}=u_{2}=0, \\
a &=&1,b_{2}=1,b_{1}=l_{1} \\
d &=&1,e_{1}-u_{1}=1\Rightarrow e_{1}=1,u_{1}=0,b_{1}=l_{1}=1,e_{2}=0.
\end{eqnarray*}%
Since $\alpha \left( 1_{H};1,0,0,0\right) \equiv b_{2}\equiv 1$ we get%
\begin{equation*}
-B(x_{1}x_{2}\otimes x_{1};GX_{1}X_{2},gx_{1})GX_{2}\otimes gx_{1}\otimes
x_{1}.
\end{equation*}%
By considering also the right side, we obtain%
\begin{gather*}
B(gx_{2}\otimes x_{1};GX_{2},gx_{1})-B(x_{1}x_{2}\otimes
x_{1};GX_{1}X_{2},gx_{1})=0 \\
-B(x_{1}x_{2}\otimes 1_{H};GX_{2},gx_{1})=0
\end{gather*}%
which holds in view of the form of the elements.

\subsubsection{Case $GX_{1}\otimes gx_{1}\otimes x_{2}$}

This case was already considered in subsection $B\left( x_{1}x_{2}\otimes
x_{1};GX_{1},gx_{1}x_{2}\right) .$

\subsubsection{Case $G\otimes gx_{1}\otimes gx_{1}x_{2}$}

This case was already considered in subsection $B\left( x_{1}x_{2}\otimes
x_{1};GX_{1},gx_{1}x_{2}\right) .$

\subsubsection{Case $GX_{1}X_{2}\otimes g\otimes x_{1}$}

We have to consider only the first and the fourth summand of the left side
of the equation. First summand gives us%
\begin{eqnarray*}
l_{1} &=&u_{1}=l_{2}=u_{2}=0 \\
a &=&b_{1}=b_{2}=1, \\
d &=&1,e_{1}=e_{2}=0.
\end{eqnarray*}%
Since $\alpha \left( x_{1};0,0,0,0\right) \equiv a+b_{1}+b_{2}\equiv 1,$ we
get%
\begin{equation*}
-B(gx_{2}\otimes x_{1};GX_{1}X_{2},g)GX_{1}X_{2}\otimes g\otimes x_{1}.
\end{equation*}%
Fourth summand is giving us%
\begin{eqnarray*}
l_{1}+u_{1} &=&1,l_{2}=u_{2}=0, \\
a &=&1,b_{2}=1,b_{1}-l_{1}=1\Rightarrow b_{1}=1,l_{1}=0,u_{1}=1 \\
d &=&1,e_{2}=0,e_{1}=u_{1}=1.
\end{eqnarray*}%
Since $\alpha \left( 1_{H};0,0,1,0\right) \equiv e_{2}+\left(
a+b_{1}+b_{2}\right) \equiv 1$ we get%
\begin{equation*}
-B(x_{1}x_{2}\otimes x_{1};GX_{1}X_{2},gx_{1})GX_{1}X_{2}\otimes g\otimes
x_{1}.
\end{equation*}%
By considering also the right side, we obtain%
\begin{equation*}
-B(gx_{2}\otimes x_{1};GX_{1}X_{2},g)-B(x_{1}x_{2}\otimes
x_{1};GX_{1}X_{2},gx_{1})-B(x_{1}x_{2}\otimes 1_{H};GX_{1}X_{2},g)=0
\end{equation*}%
which holds in view of the form of the elements.

\subsubsection{Case $GX_{1}\otimes g\otimes gx_{1}x_{2}$}

This case was already considered in subsection $B\left( x_{1}x_{2}\otimes
x_{1};GX_{1},gx_{1}x_{2}\right) .$

{\Huge \ }

\section{$B\left( x_{1}x_{2}\otimes x_{2}\right) $}

By $\left( \ref{simplx}\right) $ we have%
\begin{equation}
B(x_{1}x_{2}\otimes x_{2})=B(x_{1}x_{2}\otimes 1_{H})(1_{A}\otimes
x_{2})-(1_{A}\otimes gx_{2})B(x_{1}x_{2}\otimes 1_{H})(1_{A}\otimes g)
\label{x1x2otx2}
\end{equation}

and we get
\begin{eqnarray*}
B(x_{1}x_{2}\otimes x_{2}) &=&+2B(x_{1}x_{2}\otimes
1_{H};1_{A},gx_{1})1_{A}\otimes gx_{1}x_{2}+ \\
&&-2B(x_{1}x_{2}\otimes 1_{H};G,g)G\otimes gx_{2}+ \\
&&-2B(x_{1}x_{2}\otimes 1_{H};X_{1},g)X_{1}\otimes gx_{2}+ \\
&&-2B(x_{1}x_{2}\otimes 1_{H};X_{2},g)X_{2}\otimes gx_{2}+ \\
&&+2B(x_{1}x_{2}\otimes 1_{H};X_{1},gx_{1}x_{2})X_{1}X_{2}\otimes gx_{1}x_{2}
\\
&&+2B(x_{1}x_{2}\otimes 1_{H};GX_{1},gx_{1})GX_{1}\otimes gx_{1}x_{2}\text{ }
\\
&&+2B(x_{1}x_{2}\otimes 1_{H};GX_{2},gx_{1})GX_{2}\otimes gx_{1}x_{2}+ \\
&&+\left[
\begin{array}{c}
2B(x_{1}x_{2}\otimes 1_{H};G,gx_{1}x_{2})-2B(x_{1}x_{2}\otimes
1_{H};GX_{2},gx_{1}) \\
+2B(x_{1}x_{2}\otimes 1_{H};GX_{1},gx_{2})%
\end{array}%
\right] GX_{1}X_{2}\otimes gx_{2}
\end{eqnarray*}

We write Casimir condition for $B(x_{1}x_{2}\otimes x_{2}).$%
\begin{eqnarray*}
&&\text{1}\sum_{a,b_{1},b_{2},d,e_{1},e_{2}=0}^{1}\sum_{l_{1}=0}^{b_{1}}%
\sum_{l_{2}=0}^{b_{2}}\sum_{u_{1}=0}^{e_{1}}\sum_{u_{2}=0}^{e_{2}}\left(
-1\right) ^{\alpha \left( x_{1}x_{2};l_{1},l_{2},u_{1},u_{2}\right) } \\
&&B(1_{H}\otimes
x_{2};G^{a}X_{1}^{b_{1}}X_{2}^{b_{2}},g^{d}x_{1}^{e_{1}}x_{2}^{e_{2}}) \\
&&G^{a}X_{1}^{b_{1}-l_{1}}X_{2}^{b_{2}-l_{2}}\otimes
g^{d}x_{1}^{e_{1}-u_{1}}x_{2}^{e_{2}-u_{2}}\otimes
g^{a+b_{1}+b_{2}+l_{1}+l_{2}+d+e_{1}+e_{2}+u_{1}+u_{2}}x_{1}^{l_{1}+u_{1}+1}x_{2}^{l_{2}+u_{2}+1}
\\
&&\text{2}\sum_{a,b_{1},b_{2},d,e_{1},e_{2}=0}^{1}\sum_{l_{1}=0}^{b_{1}}%
\sum_{l_{2}=0}^{b_{2}}\sum_{u_{1}=0}^{e_{1}}\sum_{u_{2}=0}^{e_{2}}\left(
-1\right) ^{\alpha \left( x_{1};l_{1},l_{2},u_{1},u_{2}\right) } \\
&&B(gx_{2}\otimes
x_{2};G^{a}X_{1}^{b_{1}}X_{2}^{b_{2}},g^{d}x_{1}^{e_{1}}x_{2}^{e_{2}}) \\
&&G^{a}X_{1}^{b_{1}-l_{1}}X_{2}^{b_{2}-l_{2}}\otimes
g^{d}x_{1}^{e_{1}-u_{1}}x_{2}^{e_{2}-u_{2}}\otimes
g^{a+b_{1}+b_{2}+l_{1}+l_{2}+d+e_{1}+e_{2}+u_{1}+u_{2}}x_{1}^{l_{1}+u_{1}+1}x_{2}^{l_{2}+u_{2}}
\\
&&\text{3 }\left( -1\right)
\sum_{a,b_{1},b_{2},d,e_{1},e_{2}=0}^{1}\sum_{l_{1}=0}^{b_{1}}%
\sum_{l_{2}=0}^{b_{2}}\sum_{u_{1}=0}^{e_{1}}\sum_{u_{2}=0}^{e_{2}}\left(
-1\right) ^{\alpha \left( x_{2};l_{1},l_{2},u_{1},u_{2}\right) } \\
&&B(gx_{1}\otimes
x_{2};G^{a}X_{1}^{b_{1}}X_{2}^{b_{2}},g^{d}x_{1}^{e_{1}}x_{2}^{e_{2}}) \\
&&G^{a}X_{1}^{b_{1}-l_{1}}X_{2}^{b_{2}-l_{2}}\otimes
g^{d}x_{1}^{e_{1}-u_{1}}x_{2}^{e_{2}-u_{2}}\otimes
g^{a+b_{1}+b_{2}+l_{1}+l_{2}+d+e_{1}+e_{2}+u_{1}+u_{2}}x_{1}^{l_{1}+u_{1}}x_{2}^{l_{2}+u_{2}+1}
\\
&&\text{4}\sum_{a,b_{1},b_{2},d,e_{1},e_{2}=0}^{1}\sum_{l_{1}=0}^{b_{1}}%
\sum_{l_{2}=0}^{b_{2}}\sum_{u_{1}=0}^{e_{1}}\sum_{u_{2}=0}^{e_{2}}\left(
-1\right) ^{\alpha \left( 1_{H};l_{1},l_{2},u_{1},u_{2}\right) } \\
&&B(x_{1}x_{2}\otimes
x_{2};G^{a}X_{1}^{b_{1}}X_{2}^{b_{2}},g^{d}x_{1}^{e_{1}}x_{2}^{e_{2}}) \\
&&G^{a}X_{1}^{b_{1}-l_{1}}X_{2}^{b_{2}-l_{2}}\otimes
g^{d}x_{1}^{e_{1}-u_{1}}x_{2}^{e_{2}-u_{2}}\otimes
g^{a+b_{1}+b_{2}+l_{1}+l_{2}+d+e_{1}+e_{2}+u_{1}+u_{2}}x_{1}^{l_{1}+u_{1}}x_{2}^{l_{2}+u_{2}}
\\
&=&B^{A}(x_{1}x_{2}\otimes x_{2})\otimes B^{H}(x_{1}x_{2}\otimes
x_{2})\otimes g+ \\
&&B^{A}(x_{1}x_{2}\otimes 1_{H})\otimes B^{H}(x_{1}x_{2}\otimes
1_{H})\otimes x_{2}+
\end{eqnarray*}

\subsection{$B(x_{1}x_{2}\otimes x_{2};1_{A},gx_{1}x_{2})$}

We deduce that%
\begin{eqnarray*}
a &=&b_{1}=b_{2}=0 \\
d &=&e_{1}=e_{2}=1 \\
a+b_{1}+b_{2}+l_{1}+l_{2}+d+e_{1}+e_{2}+u_{1}+u_{2} &\equiv &u_{1}+u_{2}+1
\end{eqnarray*}%
and we get%
\begin{eqnarray*}
&&\left( -1\right) ^{\alpha \left( 1_{H};0,0,0,0\right) }B(x_{1}x_{2}\otimes
x_{2};1_{H},gx_{1}x_{2})1_{H}\otimes gx_{1}x_{2}\otimes g \\
&&\left( -1\right) ^{\alpha \left( 1_{H};0,0,0,1\right) }B(x_{1}x_{2}\otimes
x_{2};1_{H},gx_{1}x_{2})1_{H}\otimes gx_{1}\otimes x_{2} \\
&&\left( -1\right) ^{\alpha \left( 1_{H};0,0,1,0\right) }B(x_{1}x_{2}\otimes
x_{2};1_{H},gx_{1}x_{2})1_{H}\otimes gx_{2}\otimes x_{1} \\
&&\left( -1\right) ^{\alpha \left( 1_{H};0,0,1,1\right) }B(x_{1}x_{2}\otimes
x_{2};1_{H},gx_{1}x_{2})1_{H}\otimes g\otimes gx_{1}x_{2}\text{ .}
\end{eqnarray*}

\subsubsection{Case $1_{H}\otimes gx_{1}\otimes x_{2}$}

We have to consider only the third and the fourth summand of the left side
of the equality. Third summand gives us%
\begin{eqnarray*}
l_{1} &=&u_{1}=0,l_{2}=u_{2}=0, \\
a &=&b_{1}=b_{2}=0, \\
d &=&e_{1}=1,e_{2}=0.
\end{eqnarray*}%
Since $\alpha \left( x_{2};0,0,0,0\right) \equiv a+b_{1}+b_{2}=0,$ we get
\begin{equation*}
-B(gx_{1}\otimes x_{2};1_{H},gx_{1})1_{A}\otimes gx_{1}x_{2}\otimes x_{2}.
\end{equation*}%
Fourth summand gives us%
\begin{eqnarray*}
l_{1} &=&u_{1}=0,l_{2}+u_{2}=1, \\
a &=&b_{1}=0,b_{2}=l_{2} \\
d &=&e_{1}=1,e_{2}=u_{2}.
\end{eqnarray*}%
Since $\alpha \left( 1_{H};0,0,1,1\right) \equiv 1+e_{2}\equiv 0$ and $%
\alpha \left( 1_{H};0,1,1,0\right) \equiv e_{2}+\left(
a+b_{1}+b_{2}+1\right) \equiv 0$ we get
\begin{equation*}
\left[ B(x_{1}x_{2}\otimes x_{2};1_{A},gx_{1}x_{2})+B(x_{1}x_{2}\otimes
x_{2};X_{2},gx_{1})\right] 1_{H}\otimes gx_{1}\otimes x_{2}.
\end{equation*}%
By considering also the right side we get%
\begin{equation*}
-B(x_{1}x_{2}\otimes 1_{H};1_{H},gx_{1})+B(x_{1}x_{2}\otimes
x_{2};1_{A},gx_{1}x_{2})+B(x_{1}x_{2}\otimes
x_{2};X_{2},gx_{1})-B(gx_{1}\otimes x_{2};1_{H},gx_{1})=0
\end{equation*}%
which holds in view of the form of the elements.

\subsubsection{Case $1_{A}\otimes gx_{2}\otimes x_{1}$}

We have to consider only the second and the fourth summand of the left side
of the equality.

Second summand gives us%
\begin{eqnarray*}
l_{1} &=&u_{1}=l_{2}=u_{2}=0, \\
a &=&b_{1}=b_{2}=0, \\
d &=&e_{2}=1,e_{1}=0.
\end{eqnarray*}%
Since $\alpha \left( x_{1};0,0,0,0\right) \equiv a+b_{1}+b_{2}=0$ we obtain
\begin{equation*}
B(gx_{2}\otimes x_{2};1_{A},gx_{2})1_{H}\otimes gx_{2}\otimes x_{1}.
\end{equation*}%
Forth summand gives us%
\begin{eqnarray*}
l_{1}+u_{1} &=&1,l_{2}=u_{2}=0 \\
a &=&b_{2}=0,b_{1}=l_{1}, \\
d &=&e_{2}=1,e_{1}=u_{1}.
\end{eqnarray*}%
Since $\alpha \left( 1_{H};0,0,1,0\right) \equiv e_{2}+\left(
a+b_{1}+b_{2}\right) \equiv 1$ and $\alpha \left( 1_{H};1,0,0,0\right)
\equiv b_{2}=0,$ we obtain
\begin{equation*}
\left[ -B(x_{1}x_{2}\otimes x_{2};1_{A},gx_{1}x_{2})+B(x_{1}x_{2}\otimes
x_{2};X_{1},gx_{2})\right] 1_{H}\otimes gx_{2}\otimes x_{1}.
\end{equation*}%
Since there is no term in the right side, we get%
\begin{equation*}
-B(x_{1}x_{2}\otimes x_{2};1_{A},gx_{1}x_{2})+B(x_{1}x_{2}\otimes
x_{2};X_{1},gx_{2})+B(gx_{2}\otimes x_{2};1_{A},gx_{2})=0
\end{equation*}%
which holds in view of the form of the elements.

\subsubsection{Case $1_{A}\otimes g\otimes gx_{1}x_{2}$}

First summand of the left side gives us

\begin{eqnarray*}
l_{1} &=&u_{1}=l_{2}=u_{2}=0 \\
a &=&b_{1}=b_{2}=0 \\
d &=&1,e_{1}=e_{2}=0.
\end{eqnarray*}%
Since $\alpha \left( x_{1}x_{2};0,0,0,0\right) \equiv 0$ we get
\begin{equation*}
B(1_{H}\otimes x_{2};1_{A},g)1_{A}\otimes g\otimes gx_{1}x_{2.}
\end{equation*}%
Second summand of the left side gives us%
\begin{eqnarray*}
l_{1} &=&u_{1}=0,l_{2}+u_{2}=1, \\
a &=&b_{1}=0,b_{2}=l_{2}, \\
d &=&1,e_{1}=0,e_{2}=u_{2}.
\end{eqnarray*}%
Since $\alpha \left( x_{1};0,0,0,1\right) \equiv 1$ and $\alpha \left(
x_{1};0,1,0,0\right) \equiv a+b_{1}+b_{2}+1\equiv 0,$ we obtain
\begin{equation*}
\left[ -B(gx_{2}\otimes x_{2};1_{A},gx_{2})+B(gx_{2}\otimes x_{2};X_{2},g)%
\right] 1_{A}\otimes g\otimes gx_{1}x_{2}.
\end{equation*}%
Third summand of the left side gives us%
\begin{eqnarray*}
l_{1}+u_{1} &=&1,l_{2}=u_{2}=0, \\
a &=&b_{2}=0,b_{1}=l_{1}, \\
d &=&1,e_{1}=u_{1},e_{2}=0.
\end{eqnarray*}%
Since $\alpha \left( x_{2};0,0,1,0\right) \equiv e_{2}=0$ and $\alpha \left(
x_{2};1,0,0,0\right) \equiv a+b_{1}=1,$ we obtain
\begin{equation*}
\left[ -B(gx_{1}\otimes x_{2};1_{A},gx_{1})+B(gx_{1}\otimes x_{2};X_{1},g%
\right] 1_{A}\otimes g\otimes gx_{1}x_{2}
\end{equation*}%
The fourth summand of the left side gives us%
\begin{eqnarray*}
l_{1}+u_{1} &=&1,l_{2}+u_{2}=1, \\
a &=&0,b_{1}=l_{1},b_{2}=l_{2}, \\
d &=&1,e_{1}=u_{1},e_{2}=u_{2}.
\end{eqnarray*}%
Since
\begin{eqnarray*}
\alpha \left( 1_{H};0,0,1,1\right) &\equiv &1+e_{2}\equiv 0 \\
\alpha \left( 1_{H};0,1,1,0\right) &\equiv &e_{2}+\left(
a+b_{1}+b_{2}+1\right) \equiv 0 \\
\alpha \left( 1_{H};1,0,0,1\right) &\equiv &a+b_{1}\equiv 1 \\
\alpha \left( 1_{H};1,1,0,0\right) &\equiv &1+b_{2}\equiv 0
\end{eqnarray*}%
we obtain
\begin{equation*}
\left[
\begin{array}{c}
B(x_{1}x_{2}\otimes x_{2};1_{A},gx_{1}x_{2})+B(x_{1}x_{2}\otimes
x_{2};X_{2},gx_{1}) \\
-B(x_{1}x_{2}\otimes x_{2};X_{1},gx_{2})+B(x_{1}x_{2}\otimes
x_{2};X_{1}X_{2},g)%
\end{array}%
\right] 1_{A}\otimes g\otimes gx_{1}x_{2}
\end{equation*}%
Since there is nothing in the right side we get%
\begin{gather*}
B(1_{H}\otimes x_{2};1_{A},g)+ \\
-B(gx_{2}\otimes x_{2};1_{A},gx_{2})+B(gx_{2}\otimes x_{2};X_{2},g) \\
-B(gx_{1}\otimes x_{2};1_{A},gx_{1})+B(gx_{1}\otimes x_{2};X_{1},g) \\
B(x_{1}x_{2}\otimes x_{2};1_{A},gx_{1}x_{2})+B(x_{1}x_{2}\otimes
x_{2};X_{2},gx_{1})+ \\
-B(x_{1}x_{2}\otimes x_{2};X_{1},gx_{2})+B(x_{1}x_{2}\otimes
x_{2};X_{1}X_{2},g)=0
\end{gather*}%
which holds in view of the form of the elements.

\subsection{$B\left( x_{1}x_{2}\otimes x_{2};G,gx_{2}\right) $}

We deduce that
\begin{equation*}
a=1,d=e_{2}=1
\end{equation*}%
and we get%
\begin{eqnarray*}
&&\left( -1\right) ^{\alpha \left( 1_{H};0,0,0,0\right) }B(x_{1}x_{2}\otimes
x_{1};G,gx_{2})G\otimes gx_{2}\otimes g+ \\
&&\left( -1\right) ^{\alpha \left( 1_{H};0,0,0,1\right) }B(x_{1}x_{2}\otimes
x_{1};G,gx_{2})G\otimes g\otimes x_{2}\text{ .}
\end{eqnarray*}

\subsubsection{Case $G\otimes g\otimes x_{2}$}

We have to consider only the third and the fourth summand of the left side
of the equality. Third summand gives us%
\begin{eqnarray*}
l_{1} &=&u_{1}=0,l_{2}=u_{2}=0 \\
a &=&1,b_{1}=b_{2}=0, \\
d &=&1,e_{1}=e_{2}=0.
\end{eqnarray*}%
Since $\alpha \left( x_{2};0,0,0,0\right) \equiv a+b_{1}+b_{2}=1,$ we get%
\begin{equation*}
B(gx_{1}\otimes x_{2};G,g)G\otimes g\otimes x_{2}.
\end{equation*}%
Fourth summand gives us%
\begin{eqnarray*}
l_{1} &=&u_{1}=0,l_{2}+u_{2}=1 \\
a &=&1,b_{1}=0,b_{2}=l_{2} \\
d &=&1,e_{1}=0,e_{2}=u_{2}.
\end{eqnarray*}%
Since $\alpha \left( 1_{H};0,0,0,1\right) \equiv a+b_{1}+b_{2}\equiv -1$ and
$\alpha \left( 1_{H};0,1,0,0\right) \equiv 0,$ we get%
\begin{equation*}
\left[ -B(x_{1}x_{2}\otimes x_{2};G,g_{1}x_{2})+B(x_{1}x_{2}\otimes
x_{2};GX_{2},g)\right] G\otimes g\otimes x_{2}.
\end{equation*}%
By considering also the right side of the equality we get%
\begin{gather*}
-B(x_{1}x_{2}\otimes 1_{H};G,g)-B(x_{1}x_{2}\otimes
x_{2};G,g_{1}x_{2})+B(x_{1}x_{2}\otimes x_{2};GX_{2},g) \\
+B(gx_{1}\otimes x_{2};G,g)=0
\end{gather*}%
which holds in view of the form of the elements.

\subsection{$B\left( x_{1}x_{2}\otimes x_{2};X_{1},gx_{2}\right) $}

We deduce that%
\begin{equation*}
b_{1}=1,d=e_{2}=1
\end{equation*}%
and we get%
\begin{eqnarray*}
&&\left( -1\right) ^{\alpha \left( 1_{H};0,0,0,0\right) }B(x_{1}x_{2}\otimes
x_{2};X_{1},gx_{2})X_{1}\otimes gx_{2}\otimes g+ \\
&&\left( -1\right) ^{\alpha \left( 1_{H};0,0,0,1\right) }B(x_{1}x_{2}\otimes
x_{2};X_{1},gx_{2})X_{1}\otimes g\otimes x_{2}+\text{ } \\
&&+\left( -1\right) ^{\alpha \left( 1_{H};1,0,0,0\right)
}B(x_{1}x_{2}\otimes x_{2};X_{1},gx_{2})1_{A}\otimes gx_{2}\otimes x_{1}+ \\
&&+\left( -1\right) ^{\alpha \left( 1_{H};1,0,0,1\right)
}B(x_{1}x_{2}\otimes x_{2};X_{1},gx_{2})1_{A}\otimes g\otimes gx_{1}x_{2}%
\text{ .}
\end{eqnarray*}

\subsubsection{Case $X_{1}\otimes g\otimes x_{2}$}

We have to consider only the third and the fourth summand of the left side
of the equality. Third summand gives us%
\begin{eqnarray*}
l_{1} &=&u_{1}=0,l_{2}=u_{2}=0 \\
a &=&b_{2}=0,b_{1}=1, \\
d &=&1,e_{1}=e_{2}=0.
\end{eqnarray*}%
Since $\alpha \left( x_{2};0,0,0,0\right) \equiv a+b_{1}+b_{2}=1,$ we get
\begin{equation*}
B(gx_{1}\otimes x_{2};X_{1},g)X_{1}\otimes g\otimes x_{2}.
\end{equation*}%
Fourth summand gives us%
\begin{eqnarray*}
l_{1} &=&u_{1}=0,l_{2}+u_{2}=1, \\
a &=&0,b_{1}=1,b_{2}=l_{2}, \\
d &=&1,e_{1}=0,e_{2}=u_{2}.
\end{eqnarray*}%
Since $\alpha \left( 1_{H};0,0,0,1\right) \equiv a+b_{1}+b_{2}=1$ and $%
\alpha \left( 1_{H};0,1,0,0\right) \equiv 0,$ we get%
\begin{equation*}
\left[ -B(x_{1}x_{2}\otimes x_{2};X_{1},gx_{2})+B(x_{1}x_{2}\otimes
x_{2};X_{1}X_{2},g)\right] X_{1}\otimes g\otimes x_{2}
\end{equation*}%
By considering also the right side we get%
\begin{gather*}
-B(x_{1}x_{2}\otimes 1_{H};X_{1},g)-B(x_{1}x_{2}\otimes x_{2};X_{1},gx_{2})
\\
+B(x_{1}x_{2}\otimes x_{2};X_{1}X_{2},g)+B(gx_{1}\otimes x_{2};X_{1},g)=0
\end{gather*}%
which holds in view of the form of the elements.

\subsubsection{Case $1_{A}\otimes gx_{2}\otimes x_{1}$}

This case was already done in subsection $B(x_{1}x_{2}\otimes
x_{2};1_{A},gx_{1}x_{2}).$

\subsubsection{Case $1_{A}\otimes g\otimes gx_{1}x_{2}$}

This case was already done in subsection $B(x_{1}x_{2}\otimes
x_{2};1_{A},gx_{1}x_{2}).$

\subsection{$B\left( x_{1}x_{2}\otimes x_{2};X_{2},gx_{2}\right) $}

We deduce that
\begin{equation*}
b_{2}=1,d=e_{2}=1
\end{equation*}%
and we get%
\begin{gather*}
\left( -1\right) ^{\alpha \left( 1_{H};0,0,0,0\right) }B(x_{1}x_{2}\otimes
x_{2};X_{2},gx_{2})X_{2}\otimes gx_{2}\otimes g+ \\
+\left( -1\right) ^{\alpha \left( 1_{H};0,1,0,0\right) }B(x_{1}x_{2}\otimes
x_{2};X_{2},gx_{2})1_{A}\otimes gx_{2}\otimes x_{2}+ \\
+\left( -1\right) ^{\alpha \left( 1_{H};0,0,0,1\right) }B(x_{1}x_{2}\otimes
x_{2};X_{2},gx_{2})X_{2}\otimes g\otimes x_{2}+ \\
++\left( -1\right) ^{\alpha \left( 1_{H};0,1,0,1\right) }B(x_{1}x_{2}\otimes
x_{2};X_{2},gx_{2})X_{2}^{1-l_{2}}\otimes g\otimes g^{l_{2}}x_{2}^{1+1}=0.
\end{gather*}


\subsubsection{Case $1_{A}\otimes gx_{2}\otimes x_{2}$}


We have to consider only the third and the fourth summand of the left side.
Third summand gives us%
\begin{eqnarray*}
l_{1} &=&u_{1}=0,l_{2}=u_{2}=0, \\
a &=&b_{1}=0,b_{2}=1, \\
d &=&1,e_{1}=e_{2}=0.
\end{eqnarray*}%
Since $\alpha \left( x_{2};0,0,0,0\right) \equiv a+b_{1}+b_{2}=1,$ we get
\begin{equation*}
-B(gx_{1}\otimes x_{2};1_{H},gx_{2})1_{A}\otimes gx_{2}\otimes x_{2}.
\end{equation*}%
Fourth summand gives us%
\begin{eqnarray*}
l_{1} &=&u_{1}=0,l_{2}+u_{2}=1, \\
a &=&b_{1}=0,b_{2}=l_{2} \\
d &=&1,e_{1}=0,e_{2}-u_{2}=1\Rightarrow e_{2}=1,u_{2}=0,b_{2}=l_{2}=1.
\end{eqnarray*}%
Since $\alpha \left( 1_{H};0,1,0,0\right) \equiv 0$ we get
\begin{equation*}
B(x_{1}x_{2}\otimes x_{2};X_{2},gx_{2})1_{H}\otimes gx_{2}\otimes x_{2}.
\end{equation*}%
By considering also the right side we get%
\begin{equation*}
-B(x_{1}x_{2}\otimes 1_{H};1_{H},gx_{2})+B(x_{1}x_{2}\otimes
x_{2};X_{2},gx_{2})-B(gx_{1}\otimes x_{2};1_{H},gx_{2})=0
\end{equation*}%
which holds in view of the form of the elements.

\subsubsection{Case $X_{2}\otimes g\otimes x_{2}$}

We have to consider only the third and the fourth summand of the left side.
Third summand gives us%
\begin{eqnarray*}
l_{1} &=&u_{1}=0,l_{2}=u_{2}=0, \\
a &=&b_{1}=0,b_{2}=1 \\
d &=&e_{2}=1,e_{1}=0.
\end{eqnarray*}%
Since $\alpha \left( x_{2};0,0,0,0\right) \equiv a+b_{1}+b_{2}=1,$ we get
\begin{equation*}
B(gx_{1}\otimes x_{2};X_{2},g)X_{2}\otimes g\otimes x_{2}.
\end{equation*}%
Fourth summand gives us%
\begin{eqnarray*}
l_{1} &=&u_{1}=0,l_{2}+u_{2}=1, \\
a &=&b_{1}=0,b_{2}-l_{2}=1\Rightarrow b_{2}=1,l_{2}=0,u_{2}=1 \\
d &=&1,e_{1}=0,e_{2}=u_{2}=1.
\end{eqnarray*}%
Since $\alpha \left( 1_{H};0,0,0,1\right) \equiv a+b_{1}+b_{2}\equiv 1$ we
get
\begin{equation*}
-B(x_{1}x_{2}\otimes x_{2};X_{2},gx_{2})X_{2}\otimes g\otimes x_{2}.
\end{equation*}%
By considering also the right side we get%
\begin{equation*}
-B(x_{1}x_{2}\otimes 1_{H};X_{2},g)-B(x_{1}x_{2}\otimes
x_{2};X_{2},gx_{2})+B(gx_{1}\otimes x_{2};X_{2},g)=0
\end{equation*}%
which holds in view of the form of the elements.

\subsection{$B\left( x_{1}x_{2}\otimes x_{2};X_{1}X_{2},gx_{1}x_{2}\right) $}

We deduce that
\begin{equation*}
a=0,b_{1}=b_{2}=1,d=e_{1}=e_{2}=1
\end{equation*}%
and we get%
\begin{gather*}
\left( -1\right) ^{\alpha \left( 1_{H};0,0,0,0\right) }B\left(
x_{1}x_{2}\otimes x_{2};X_{1}X_{2},gx_{1}x_{2}\right) X_{1}X_{2}\otimes
gx_{1}x_{2}\otimes g \\
+\left( -1\right) ^{\alpha \left( 1_{H};1,0,0,0\right) }B\left(
x_{1}x_{2}\otimes x_{2};X_{1}X_{2},gx_{1}x_{2}\right) X_{2}\otimes
gx_{1}x_{2}\otimes x_{1} \\
+\left( -1\right) ^{\alpha \left( 1_{H};0,1,0,0\right) }B\left(
x_{1}x_{2}\otimes x_{2};X_{1}X_{2},gx_{1}x_{2}\right) X_{1}\otimes
gx_{1}x_{2}\otimes x_{2} \\
+\left( -1\right) ^{\alpha \left( 1_{H};1,1,0,0\right) }B\left(
x_{1}x_{2}\otimes x_{2};X_{1}X_{2},gx_{1}x_{2}\right) 1_{A}\otimes
gx_{1}x_{2}\otimes gx_{1}x_{2} \\
+\left( -1\right) ^{\alpha \left( 1_{H};0,0,1,0\right) }B\left(
x_{1}x_{2}\otimes x_{2};X_{1}X_{2},gx_{1}x_{2}\right) X_{1}X_{2}\otimes
gx_{2}\otimes x_{1} \\
+\left( -1\right) ^{\alpha \left( 1_{H};1,0,1,0\right) }B\left( x_{2}\otimes
x_{1};X_{1}X_{2},gx_{1}x_{2}\right) X_{1}^{1-l_{1}}X_{2}\otimes
gx_{2}\otimes g^{l_{1}}x_{1}^{1+1}=0 \\
+\left( -1\right) ^{\alpha \left( 1_{H};0,1,1,0\right) }B\left(
x_{1}x_{2}\otimes x_{2};X_{1}X_{2},gx_{1}x_{2}\right) X_{1}\otimes
gx_{2}\otimes gx_{1}x_{2} \\
+\left( -1\right) ^{\alpha \left( 1_{H};1,1,1,0\right) }B\left(
x_{1}x_{2}\otimes x_{2};X_{1}X_{2},gx_{1}x_{2}\right) X_{1}^{1-l_{1}}\otimes
gx_{2}\otimes g^{l_{1}+1}x_{1}^{1+1}x_{2}=0 \\
+\left( -1\right) ^{\alpha \left( 1_{H};0,0,0,1\right) }B\left(
x_{1}x_{2}\otimes x_{2};X_{1}X_{2},gx_{1}x_{2}\right) X_{1}X_{2}\otimes
gx_{1}\otimes x_{2} \\
+\left( -1\right) ^{\alpha \left( 1_{H};1,0,0,1\right) }B\left(
x_{1}x_{2}\otimes x_{2};X_{1}X_{2},gx_{1}x_{2}\right) X_{2}\otimes
gx_{1}\otimes gx_{1}x_{2} \\
+\left( -1\right) ^{\alpha \left( 1_{H};l_{1},1,0,1\right) }B\left(
x_{1}x_{2}\otimes x_{2};X_{1}X_{2},gx_{1}x_{2}\right)
X_{1}^{1-l_{1}}X_{2}^{1-l_{2}}\otimes gx_{1}\otimes
g^{l_{1}+l_{2}}x_{1}^{l_{1}}x_{2}^{1+1}=0 \\
+\left( -1\right) ^{\alpha \left( 1_{H};0,0,1,1\right) }B\left(
x_{1}x_{2}\otimes x_{2};X_{1}X_{2},gx_{1}x_{2}\right) X_{1}X_{2}\otimes
g\otimes gx_{1}x_{2} \\
+\left( -1\right) ^{\alpha \left( 1_{H};1,0,1,1\right) }B\left(
x_{1}x_{2}\otimes x_{2};X_{1}X_{2},gx_{1}x_{2}\right)
X_{1}^{1-l_{1}}X_{2}\otimes g\otimes g^{l_{1}+1}x_{1}^{1+1}x_{2}=0 \\
+\left( -1\right) ^{\alpha \left( 1_{H};1,1,1,1\right) }B\left(
x_{1}x_{2}\otimes x_{2};X_{1}X_{2},gx_{1}x_{2}\right)
X_{1}^{1-l_{1}}X_{2}^{1-l_{2}}\otimes g\otimes
g^{l_{1}+l_{2}+1}x_{1}^{l_{1}+1}x_{2}^{1+1}=0.
\end{gather*}

\subsubsection{Case $X_{2}\otimes gx_{1}x_{2}\otimes x_{1}$}

We have to consider only the second and the fourth summand of the left side
of the equality.

We have to consider only the second and the fourth summand of the left side
of the equality. Second summand gives us%
\begin{eqnarray*}
l_{1} &=&u_{1}=l_{2}=u_{2}=0, \\
a &=&b_{1}=0,b_{2}=1, \\
d &=&e_{1}=e_{2}=1.
\end{eqnarray*}%
Since $\alpha \left( x_{1};0,0,0,0\right) \equiv a+b_{1}+b_{2}=1$ we obtain
\begin{equation*}
-B(gx_{2}\otimes x_{2};X_{2},gx_{1}x_{2})X_{2}\otimes gx_{1}x_{2}\otimes
x_{1}.
\end{equation*}%
Forth summand gives us%
\begin{eqnarray*}
l_{1}+u_{1} &=&1,l_{2}=u_{2}=0 \\
a &=&0,b_{2}=1,b_{1}=l_{1}, \\
d &=&e_{2}=1,e_{1}-u_{1}=1\Rightarrow e_{1}=1,u_{1}=0,b_{1}=l_{1}=1
\end{eqnarray*}%
Since $\alpha \left( 1_{H};1,0,0,0\right) \equiv b_{2}=1,$ we obtain
\begin{equation*}
-B(x_{1}x_{2}\otimes x_{2};X_{1}X_{2},gx_{1}x_{2})X_{2}\otimes
gx_{1}x_{2}\otimes x_{1}.
\end{equation*}%
Since there is no term in the right side, we get%
\begin{equation*}
-B(x_{1}x_{2}\otimes x_{2};X_{1}X_{2},gx_{1}x_{2})-B(gx_{2}\otimes
x_{2};X_{2},gx_{1}x_{2})=0
\end{equation*}%
which holds in view of the form of the elements.

\subsubsection{Case $X_{1}\otimes gx_{1}x_{2}\otimes x_{2}$}

We have to consider only the third and the fourth summand of the left side
of the equality. Third summand gives us%
\begin{eqnarray*}
l_{1} &=&u_{1}=0,l_{2}=u_{2}=0, \\
a &=&b_{2}=0,b_{1}=1 \\
d &=&e_{1}=e_{2}=1.
\end{eqnarray*}%
Since $\alpha \left( x_{2};0,0,0,0\right) \equiv a+b_{1}+b_{2}=1,$ we get
\begin{equation*}
B(gx_{1}\otimes x_{2};X_{1},gx_{1}x_{2})X_{1}\otimes gx_{1}x_{2}\otimes
x_{2}.
\end{equation*}%
Fourth summand gives us%
\begin{eqnarray*}
l_{1} &=&u_{1}=0,l_{2}+u_{2}=1, \\
a &=&0,b_{1}=1,b_{2}=l_{2}, \\
d &=&1,e_{1}=1,e_{2}-u_{2}=1\Rightarrow e_{2}=1,u_{2}=0,b_{2}=l_{2}=1.
\end{eqnarray*}%
Since $\alpha \left( 1_{H};0,1,0,0\right) \equiv 0$ we get
\begin{equation*}
B(x_{1}x_{2}\otimes x_{2};X_{1}X_{2},gx_{1}x_{2})X_{1}\otimes
gx_{1}x_{2}\otimes x_{2}.
\end{equation*}%
By considering also the right side we get%
\begin{equation*}
-B(x_{1}x_{2}\otimes 1_{H};X_{1},gx_{1}x_{2})+B(x_{1}x_{2}\otimes
x_{2};X_{1}X_{2},gx_{1}x_{2})+B(gx_{1}\otimes x_{2};X_{1},gx_{1}x_{2})=0
\end{equation*}%
which holds in view of the form of the elements.

\subsubsection{Case $1_{A}\otimes gx_{1}x_{2}\otimes gx_{1}x_{2}$}

First summand of the left side gives us

\begin{eqnarray*}
l_{1} &=&u_{1}=l_{2}=u_{2}=0 \\
a &=&b_{1}=b_{2}=0 \\
d &=&e_{1}=e_{2}=1.
\end{eqnarray*}%
Since $\alpha \left( x_{1}x_{2};0,0,0,0\right) \equiv 0$ we get
\begin{equation*}
B(1_{H}\otimes x_{2};1_{A},gx_{1}x_{2})1_{A}\otimes gx_{1}x_{2}\otimes
gx_{1}x_{2}
\end{equation*}%
Second summand of the left side gives us%
\begin{eqnarray*}
l_{1} &=&u_{1}=0,l_{2}+u_{2}=1, \\
a &=&b_{1}=0,b_{2}=l_{2}, \\
d &=&e_{1}=1,e_{2}-u_{2}=1\Rightarrow e_{2}=1,u_{2}=0,b_{2}=l_{2}=1.
\end{eqnarray*}%
Since $\alpha \left( x_{1};0,1,0,0\right) \equiv a+b_{1}+b_{2}+1\equiv 0,$
we obtain
\begin{equation*}
B(gx_{2}\otimes x_{2};X_{2},gx_{1}x_{2})1_{A}\otimes gx_{1}x_{2}\otimes
gx_{1}x_{2}.
\end{equation*}%
Third summand of the left side gives us%
\begin{eqnarray*}
l_{1}+u_{1} &=&1,l_{2}=u_{2}=0, \\
a &=&b_{2}=0,b_{1}=l_{1}, \\
d &=&e_{2}=1,e_{1}-u_{1}=1\Rightarrow e_{1}=1,u_{1}=0,b_{1}=l_{1}=1.
\end{eqnarray*}%
Since $\alpha \left( x_{2};1,0,0,0\right) \equiv a+b_{1}=1,$ we obtain
\begin{equation*}
+B(gx_{1}\otimes x_{2};X_{1},gx_{1}x_{2})1_{A}\otimes gx_{1}x_{2}\otimes
gx_{1}x_{2}
\end{equation*}%
The fourth summand of the left side gives us%
\begin{eqnarray*}
l_{1}+u_{1} &=&1,l_{2}+u_{2}=1, \\
a &=&0,b_{1}=l_{1},b_{2}=l_{2}, \\
d &=&1,e_{1}-u_{1}=1\Rightarrow e_{1}=1,u_{1}=0,b_{1}=l_{1}=1, \\
e_{2}-u_{2} &=&1\Rightarrow e_{2}=1,u_{2}=0,b_{2}=l_{2}=1
\end{eqnarray*}%
Since $\alpha \left( 1_{H};1,1,0,0\right) \equiv 1+b_{2}\equiv 0$ we obtain
\begin{equation*}
B(x_{1}x_{2}\otimes x_{2};X_{1}X_{2},gx_{1}x_{2})1_{A}\otimes
gx_{1}x_{2}\otimes gx_{1}x_{2}
\end{equation*}%
Since there is nothing in the right side we get%
\begin{gather*}
B(1_{H}\otimes x_{2};1_{A},gx_{1}x_{2})+B(gx_{2}\otimes
x_{2};X_{2},gx_{1}x_{2}) \\
+B(gx_{1}\otimes x_{2};X_{1},gx_{1}x_{2})+B(x_{1}x_{2}\otimes
x_{2};X_{1}X_{2},gx_{1}x_{2})=0
\end{gather*}%
which holds in view of the form of the elements.

\subsubsection{Case $X_{1}X_{2}\otimes gx_{2}\otimes x_{1}$}

We have to consider only the second and the fourth summand of the left side
of the equality.

Second summand gives us%
\begin{eqnarray*}
l_{1} &=&u_{1}=l_{2}=u_{2}=0, \\
a &=&0,b_{1}=b_{2}=1, \\
d &=&e_{2}=1,e_{1}=0.
\end{eqnarray*}%
Since $\alpha \left( x_{1};0,0,0,0\right) \equiv a+b_{1}+b_{2}=0$ we obtain
\begin{equation*}
+B(gx_{2}\otimes x_{2};X_{1}X_{2},gx_{2})X_{1}X_{2}\otimes gx_{2}\otimes
x_{1}.
\end{equation*}%
Forth summand gives us%
\begin{eqnarray*}
l_{1}+u_{1} &=&1,l_{2}=u_{2}=0 \\
a &=&0,b_{2}=1,b_{1}-l_{1}=1\Rightarrow b_{1}=1,l_{1}=0,u_{1}=1, \\
d &=&e_{2}=1,e_{1}=u_{1}=1.
\end{eqnarray*}%
Since $\alpha \left( 1_{H};0,0,1,0\right) \equiv e_{2}+\left(
a+b_{1}+b_{2}\right) \equiv 1,$ we obtain
\begin{equation*}
-B(x_{1}x_{2}\otimes x_{2};X_{1}X_{2},gx_{1}x_{2})X_{1}X_{2}\otimes
gx_{2}\otimes x_{1}.
\end{equation*}%
Since there is no term in the right side, we get%
\begin{equation*}
-B(x_{1}x_{2}\otimes x_{2};X_{1}X_{2},gx_{1}x_{2})+B(gx_{2}\otimes
x_{2};X_{1}X_{2},gx_{2})=0
\end{equation*}%
which holds in view of the form of the elements.

\subsubsection{Case $X_{1}\otimes gx_{2}\otimes gx_{1}x_{2}$}

First summand of the left side gives us

\begin{eqnarray*}
l_{1} &=&u_{1}=l_{2}=u_{2}=0 \\
a &=&b_{2}=0,b_{1}=1 \\
d &=&e_{2}=1,e_{1}=0.
\end{eqnarray*}%
Since $\alpha \left( x_{1}x_{2};0,0,0,0\right) \equiv 0$ we get
\begin{equation*}
B(1_{H}\otimes x_{2};X_{1},gx_{2})X_{1}\otimes gx_{2}\otimes gx_{1}x_{2}
\end{equation*}%
Second summand of the left side gives us%
\begin{eqnarray*}
l_{1} &=&u_{1}=0,l_{2}+u_{2}=1, \\
a &=&0,b_{1}=1,b_{2}=l_{2}, \\
d &=&1,e_{1}=0,e_{2}-u_{2}=1\Rightarrow e_{2}=1,u_{2}=0,b_{2}=l_{2}=1.
\end{eqnarray*}%
Since $\alpha \left( x_{1};0,1,0,0\right) \equiv a+b_{1}+b_{2}+1\equiv 1,$
we obtain
\begin{equation*}
-B(gx_{2}\otimes x_{2};X_{1}X_{2},gx_{2})X_{1}\otimes gx_{2}\otimes
gx_{1}x_{2}.
\end{equation*}%
Third summand of the left side gives us%
\begin{eqnarray*}
l_{1}+u_{1} &=&1,l_{2}=u_{2}=0, \\
a &=&b_{2}=0,b_{1}-l_{1}=1\Rightarrow b_{1}=1,l_{1}=0,u_{1}=1, \\
d &=&e_{2}=1,e_{1}=u_{1}=1.
\end{eqnarray*}%
Since $\alpha \left( x_{2};0,0,1,0\right) \equiv e_{2}=1,$ we obtain
\begin{equation*}
+B(gx_{1}\otimes x_{2};X_{1},gx_{1}x_{2})X_{1}\otimes gx_{2}\otimes
gx_{1}x_{2}
\end{equation*}%
The fourth summand of the left side gives us%
\begin{eqnarray*}
l_{1}+u_{1} &=&1,l_{2}+u_{2}=1, \\
a &=&0,b_{2}=l_{2},b_{1}-l_{1}=1\Rightarrow b_{1}=1,l_{1}=0,u_{1}=1, \\
d &=&1,e_{1}=u_{1}=1, \\
e_{2}-u_{2} &=&1\Rightarrow e_{2}=1,u_{2}=0,b_{2}=l_{2}=1
\end{eqnarray*}%
Since $\alpha \left( 1_{H};0,1,1,0\right) \equiv e_{2}+a+b_{1}+b_{2}+1\equiv
0$ we obtain
\begin{equation*}
B(x_{1}x_{2}\otimes x_{2};X_{1}X_{2},gx_{1}x_{2})X_{1}\otimes gx_{2}\otimes
gx_{1}x_{2}
\end{equation*}%
Since there is nothing in the right side we get%
\begin{gather*}
B(1_{H}\otimes x_{2};X_{1},gx_{2}))-B(gx_{2}\otimes x_{2};X_{1}X_{2},gx_{2})
\\
+B(gx_{1}\otimes x_{2};X_{1},gx_{1}x_{2})+B(x_{1}x_{2}\otimes
x_{2};X_{1}X_{2},gx_{1}x_{2})=0
\end{gather*}%
which holds in view of the form of the elements.

\subsubsection{Case $X_{1}X_{2}\otimes gx_{1}\otimes x_{2}$}

We have to consider only the third and the fourth summand of the left side.
Third summand gives us%
\begin{eqnarray*}
l_{1} &=&u_{1}=0,l_{2}=u_{2}=0, \\
a &=&0,b_{1}=b_{2}=1, \\
d &=&e_{1}=1,e_{2}=0.
\end{eqnarray*}%
Since $\alpha \left( x_{2};0,0,0,0\right) \equiv a+b_{1}+b_{2}=0,$ we get
\begin{equation*}
-B(gx_{1}\otimes x_{2};X_{1}X_{2},gx_{1})X_{1}X_{2}\otimes gx_{1}\otimes
x_{2}.
\end{equation*}%
Fourth summand gives us%
\begin{eqnarray*}
l_{1} &=&u_{1}=0,l_{2}+u_{2}=1, \\
a &=&0,b_{1}=1,b_{2}-l_{2}=1\Rightarrow b_{2}=1,l_{2}=0,u_{2}=1 \\
d &=&e_{1}=1,0,e_{2}=u_{2}=1.
\end{eqnarray*}%
Since $\alpha \left( 1_{H};0,0,0,1\right) \equiv a+b_{1}+b_{2}\equiv 0$ we
get
\begin{equation*}
B(x_{1}x_{2}\otimes x_{2};X_{1}X_{2},gx_{1}x_{2})X_{1}X_{2}\otimes
gx_{1}\otimes x_{2}.
\end{equation*}%
By considering also the right side we get%
\begin{equation*}
-B(x_{1}x_{2}\otimes 1_{H};X_{1}X_{2},gx_{1})+B(x_{1}x_{2}\otimes
x_{2};X_{1}X_{2},gx_{1}x_{2})-B(gx_{1}\otimes x_{2};X_{1}X_{2},gx_{1})=0
\end{equation*}%
which holds in view of the form of the elements.

\subsubsection{Case $X_{2}\otimes gx_{1}\otimes gx_{1}x_{2}$}

First summand of the left side gives us

\begin{eqnarray*}
l_{1} &=&u_{1}=l_{2}=u_{2}=0 \\
a &=&b_{1}=0,b_{2}=1 \\
d &=&e_{1}=1,e_{2}=0.
\end{eqnarray*}%
Since $\alpha \left( x_{1}x_{2};0,0,0,0\right) \equiv 0$ we get
\begin{equation*}
B(1_{H}\otimes x_{2};X_{2},gx_{1})X_{2}\otimes gx_{1}\otimes gx_{1}x_{2}
\end{equation*}%
Second summand of the left side gives us%
\begin{eqnarray*}
l_{1} &=&u_{1}=0,l_{2}+u_{2}=1, \\
a &=&b_{1}=0,b_{2}-l_{2}=1\Rightarrow b_{2}=1,l_{2}=0,u_{2}=1 \\
d &=&e_{1}=1,e_{2}=u_{2}=1.
\end{eqnarray*}%
Since $\alpha \left( x_{1};0,0,0,1\right) \equiv 1,$ we obtain
\begin{equation*}
-B(gx_{2}\otimes x_{2};X_{2},gx_{1}x_{2})X_{2}\otimes gx_{1}\otimes
gx_{1}x_{2}.
\end{equation*}%
Third summand of the left side gives us%
\begin{eqnarray*}
l_{1}+u_{1} &=&1,l_{2}=u_{2}=0, \\
a &=&0,b_{1}=l_{1},b_{2}=1, \\
d &=&1,e_{2}=0,e_{1}-u_{1}=1\Rightarrow e_{1}=1,u_{1}=0,b_{1}=l_{1}=1.
\end{eqnarray*}%
Since $\alpha \left( x_{2};1,0,0,0\right) \equiv a+b_{1}=1,$ we obtain
\begin{equation*}
+B(gx_{1}\otimes x_{2};X_{1}X_{2},gx_{1})X_{2}\otimes gx_{1}\otimes
gx_{1}x_{2}
\end{equation*}%
The fourth summand of the left side gives us%
\begin{eqnarray*}
l_{1}+u_{1} &=&1,l_{2}+u_{2}=1, \\
a &=&0,b_{1}=l_{1},b_{2}-l_{2}=1\Rightarrow b_{2}=1,l_{2}=0,u_{2}=1 \\
d &=&1,e_{1}-u_{1}=1\Rightarrow e_{1}=1,u_{1}=0,b_{1}=l_{1}=1, \\
e_{2} &=&u_{2}=1.
\end{eqnarray*}%
Since $\alpha \left( 1_{H};1,0,0,1\right) \equiv a+b_{1}\equiv 1$ we obtain
\begin{equation*}
-B(x_{1}x_{2}\otimes x_{2};X_{1}X_{2},gx_{1}x_{2})X_{2}\otimes gx_{1}\otimes
gx_{1}x_{2}
\end{equation*}%
Since there is nothing in the right side we get%
\begin{equation*}
B(1_{H}\otimes x_{2};X_{2},gx_{1})-B(gx_{2}\otimes
x_{2};X_{2},gx_{1}x_{2})+B(gx_{1}\otimes
x_{2};X_{1}X_{2},gx_{1})-B(x_{1}x_{2}\otimes x_{2};X_{1}X_{2},gx_{1}x_{2})=0
\end{equation*}%
which holds in view of the form of the elements.

\subsubsection{Case $X_{1}X_{2}\otimes g\otimes gx_{1}x_{2}$}

First summand of the left side gives us

\begin{eqnarray*}
l_{1} &=&u_{1}=l_{2}=u_{2}=0 \\
a &=&0,b_{1}=b_{2}=1, \\
d &=&1,e_{1}=e_{2}=0.
\end{eqnarray*}%
Since $\alpha \left( x_{1}x_{2};0,0,0,0\right) \equiv 0$ we get
\begin{equation*}
B(1_{H}\otimes x_{2};X_{1}X_{2},g)X_{1}X_{2}\otimes g\otimes gx_{1}x_{2}
\end{equation*}%
Second summand of the left side gives us%
\begin{eqnarray*}
l_{1} &=&u_{1}=0,l_{2}+u_{2}=1, \\
a &=&0,b_{1}=1,b_{2}-l_{2}=1\Rightarrow b_{2}=1,l_{2}=0,u_{2}=1 \\
d &=&1,e_{1}=0,e_{2}=u_{2}=1.
\end{eqnarray*}%
Since $\alpha \left( x_{1};0,0,0,1\right) \equiv 1,$ we obtain
\begin{equation*}
-B(gx_{2}\otimes x_{2};X_{1}X_{2},gx_{2})X_{1}X_{2}\otimes g\otimes
gx_{1}x_{2}.
\end{equation*}%
Third summand of the left side gives us%
\begin{eqnarray*}
l_{1}+u_{1} &=&1,l_{2}=u_{2}=0, \\
a &=&0,b_{1}-l_{1}=1\Rightarrow b_{1}=1,l_{1}=0,u_{1}=1,b_{2}=1, \\
d &=&1,e_{1}=u_{1}=1,e_{2}=0.
\end{eqnarray*}%
Since $\alpha \left( x_{2};0,0,1,0\right) \equiv e_{2}=0,$ we obtain
\begin{equation*}
-B(gx_{1}\otimes x_{2};X_{1}X_{2},gx_{1})X_{1}X_{2}\otimes g\otimes
gx_{1}x_{2}
\end{equation*}%
The fourth summand of the left side gives us%
\begin{eqnarray*}
l_{1}+u_{1} &=&1,l_{2}+u_{2}=1, \\
a &=&0,b_{1}-l_{1}=1\Rightarrow b_{1}=1,l_{1}=0,u_{1}=1, \\
b_{2}-l_{2} &=&1\Rightarrow b_{2}=1,l_{2}=0,u_{2}=1 \\
d &=&1,e_{1}=u_{1}=1,e_{2}=u_{2}=1.
\end{eqnarray*}%
Since $\alpha \left( 1_{H};0,0,1,1\right) \equiv 1+e_{2}\equiv 0$ we obtain
\begin{equation*}
B(x_{1}x_{2}\otimes x_{2};X_{1}X_{2},gx_{1}x_{2})X_{1}X_{2}\otimes g\otimes
gx_{1}x_{2}
\end{equation*}%
Since there is nothing in the right side we get%
\begin{gather*}
B(1_{H}\otimes x_{2};X_{1}X_{2},g)-B(gx_{2}\otimes x_{2};X_{1}X_{2},gx_{2})
\\
-B(gx_{1}\otimes x_{2};X_{1}X_{2},gx_{1})+B(x_{1}x_{2}\otimes
x_{2};X_{1}X_{2},gx_{1}x_{2})=0
\end{gather*}%
which holds in view of the form of the elements.

\subsection{$B\left( x_{1}x_{2}\otimes x_{2};GX_{1},gx_{1}x_{2}\right) $}

We deduce that%
\begin{eqnarray*}
a &=&1,b_{1}=1,b_{2}=0 \\
d &=&e_{1}=e_{2}=1
\end{eqnarray*}%
and we get%
\begin{gather*}
\left( -1\right) ^{\alpha \left( 1_{H};0,0,0,0\right) }B(x_{1}x_{2}\otimes
x_{2};GX_{1},gx_{1}x_{2})GX_{1}\otimes gx_{1}x_{2}\otimes g \\
\left( -1\right) ^{\alpha \left( 1_{H};1,0,0,0\right) }B(x_{1}x_{2}\otimes
x_{2};GX_{1},gx_{1}x_{2})G\otimes gx_{1}x_{2}\otimes x_{1} \\
\left( -1\right) ^{\alpha \left( 1_{H};0,0,1,0\right) }B(x_{1}x_{2}\otimes
x_{2};GX_{1},gx_{1}x_{2})GX_{1}\otimes gx_{2}\otimes x_{1} \\
\left( -1\right) ^{\alpha \left( 1_{H};1,0,1,0\right) }B(x_{1}x_{2}\otimes
x_{2};GX_{1},gx_{1}x_{2})GX_{1}^{1-l_{1}}\otimes gx_{2}\otimes
g^{+l_{11}}x_{1}^{1+1}=0 \\
\left( -1\right) ^{\alpha \left( 1_{H};0,0,0,1\right) }B(x_{1}x_{2}\otimes
x_{2};GX_{1},gx_{1}x_{2})GX_{1}\otimes gx_{1}\otimes x_{2} \\
\left( -1\right) ^{\alpha \left( 1_{H};1,0,0,1\right) }B(x_{1}x_{2}\otimes
x_{2};GX_{1},gx_{1}x_{2})G\otimes gx_{1}\otimes gx_{1}x_{2} \\
\left( -1\right) ^{\alpha \left( 1_{H};0,0,1,1\right) }B(x_{1}x_{2}\otimes
x_{2};GX_{1},gx_{1}x_{2})GX_{1}\otimes g\otimes gx_{1}x_{2} \\
\left( -1\right) ^{\alpha \left( 1_{H};1,0,1,1\right) }B(x_{1}x_{2}\otimes
x_{2};GX_{1},gx_{1}x_{2})GX_{1}^{1-l_{1}}\otimes g\otimes
g^{+l_{1}+1}x_{1}^{1+1}x_{2}=0
\end{gather*}

\subsubsection{Case $G\otimes gx_{1}x_{2}\otimes x_{1}$}

We have to consider only the second and the fourth summand of the left side
of the equality.

Second summand gives us%
\begin{eqnarray*}
l_{1} &=&u_{1}=l_{2}=u_{2}=0, \\
a &=&1,b_{1}=b_{2}=0, \\
d &=&e_{1}=e_{2}=1.
\end{eqnarray*}%
Since $\alpha \left( x_{1};0,0,0,0\right) \equiv a+b_{1}+b_{2}=1$ we obtain
\begin{equation*}
-B(gx_{2}\otimes x_{2};G,gx_{1}x_{2})G\otimes gx_{1}x_{2}\otimes x_{1}.
\end{equation*}%
Forth summand gives us%
\begin{eqnarray*}
l_{1}+u_{1} &=&1,l_{2}=u_{2}=0 \\
a &=&1,b_{1}=l_{1},b_{2}=0, \\
d &=&e_{2}=1,e_{1}-u_{1}=1\Rightarrow e_{1}=1,u_{1}=0,b_{1}=l_{1}=1.
\end{eqnarray*}%
Since $\alpha \left( 1_{H};1,0,0,0\right) \equiv b_{2}=0,$ we obtain
\begin{equation*}
B(x_{1}x_{2}\otimes x_{2};GX_{1},gx_{1}x_{2})G\otimes gx_{1}x_{2}\otimes
x_{1}.
\end{equation*}%
Since there is no term in the right side, we get%
\begin{equation*}
-B(gx_{2}\otimes x_{2};G,gx_{1}x_{2})+B(x_{1}x_{2}\otimes
x_{2};GX_{1},gx_{1}x_{2})=0
\end{equation*}%
which holds in view of the form of the elements.

\subsubsection{Case $GX_{1}\otimes gx_{2}\otimes x_{1}$}

We have to consider only the second and the fourth summand of the left side
of the equality.

Second summand gives us%
\begin{eqnarray*}
l_{1} &=&u_{1}=l_{2}=u_{2}=0, \\
a &=&b_{1}=1,b_{2}=0, \\
d &=&e_{2}=1,e_{1}=0.
\end{eqnarray*}%
Since $\alpha \left( x_{1};0,0,0,0\right) \equiv a+b_{1}+b_{2}\equiv 0$ we
obtain
\begin{equation*}
B(gx_{2}\otimes x_{2};GX_{1},gx_{2})GX_{1}\otimes gx_{2}\otimes x_{1}.
\end{equation*}%
Forth summand gives us%
\begin{eqnarray*}
l_{1}+u_{1} &=&1,l_{2}=u_{2}=0 \\
a &=&1,b_{2}=0,b_{1}-l_{1}=1\Rightarrow b_{1}=1,l_{1}=0,u_{1}=1, \\
d &=&e_{2}=1,e_{1}=u_{1}=1.
\end{eqnarray*}%
Since $\alpha \left( 1_{H};0,0,1,0\right) \equiv e_{2}+\left(
a+b_{1}+b_{2}\right) \equiv 1,$ we obtain
\begin{equation*}
-B(x_{1}x_{2}\otimes x_{2};GX_{1},gx_{1}x_{2})GX_{1}\otimes gx_{2}\otimes
x_{1}.
\end{equation*}%
Since there is no term in the right side, we get%
\begin{equation*}
B(gx_{2}\otimes x_{2};GX_{1},gx_{2})-B(x_{1}x_{2}\otimes
x_{2};GX_{1},gx_{1}x_{2})=0
\end{equation*}%
which holds in view of the form of the elements.

\subsubsection{Case $GX_{1}\otimes gx_{1}\otimes x_{2}$}

We have to consider only the third and the fourth summand of the left side.
Third summand gives us%
\begin{eqnarray*}
l_{1} &=&u_{1}=0,l_{2}=u_{2}=0, \\
a &=&b_{1}=1,b_{2}=0, \\
d &=&e_{1}=1,e_{2}=0.
\end{eqnarray*}%
Since $\alpha \left( x_{2};0,0,0,0\right) \equiv a+b_{1}+b_{2}\equiv 0,$ we
get
\begin{equation*}
-B(gx_{1}\otimes x_{2};GX_{1},gx_{1})GX_{1}\otimes gx_{1}\otimes x_{2}.
\end{equation*}%
Fourth summand gives us%
\begin{eqnarray*}
l_{1} &=&u_{1}=0,l_{2}+u_{2}=1, \\
a &=&b_{1}=1,b_{2}=l_{2}, \\
d &=&e_{1}=1,e_{2}=u_{2}.
\end{eqnarray*}%
Since $\alpha \left( 1_{H};0,0,0,1\right) \equiv a+b_{1}+b_{2}\equiv 0$ and $%
\alpha \left( 1_{H};0,1,0,0\right) \equiv 0$ we get
\begin{equation*}
\left[ B(x_{1}x_{2}\otimes x_{2};GX_{1},gx_{1}x_{2})+B(x_{1}x_{2}\otimes
x_{2};GX_{1}X_{2},gx_{1})\right] GX_{1}\otimes gx_{1}\otimes x_{2}.
\end{equation*}%
By considering also the right side we get%
\begin{gather*}
-B(x_{1}x_{2}\otimes 1_{H};GX_{1},gx_{1})+B(x_{1}x_{2}\otimes
x_{2};GX_{1},gx_{1}x_{2})+ \\
B(x_{1}x_{2}\otimes x_{2};GX_{1}X_{2},gx_{1})-B(gx_{1}\otimes
x_{2};GX_{1},gx_{1})=0
\end{gather*}%
which holds in view of the form of the elements.

\subsubsection{Case $G\otimes gx_{1}\otimes gx_{1}x_{2}$}

First summand of the left side gives us

\begin{eqnarray*}
l_{1} &=&u_{1}=l_{2}=u_{2}=0 \\
a &=&1,b_{1}=b_{2}=0, \\
d &=&e_{1}=1,e_{2}=0.
\end{eqnarray*}%
Since $\alpha \left( x_{1}x_{2};0,0,0,0\right) \equiv 0$ we get
\begin{equation*}
B(1_{H}\otimes x_{2};G,gx_{1})G\otimes gx_{1}\otimes gx_{1}x_{2}
\end{equation*}%
Second summand of the left side gives us%
\begin{eqnarray*}
l_{1} &=&u_{1}=0,l_{2}+u_{2}=1, \\
a &=&1,b_{1}=0,b_{2}=l_{2} \\
d &=&e_{1}=1,e_{2}=u_{2}.
\end{eqnarray*}%
Since $\alpha \left( x_{1};0,0,0,1\right) \equiv 1$ and $\alpha \left(
x_{1};0,1,0,0\right) \equiv a+b_{1}+b_{2}+1\equiv 1$ we obtain
\begin{equation*}
\left[ -B(gx_{2}\otimes x_{2};G,gx_{1}x_{2})-B(gx_{2}\otimes
x_{2};GX_{2},gx_{1})\right] G\otimes gx_{1}\otimes gx_{1}x_{2}.
\end{equation*}%
Third summand of the left side gives us%
\begin{eqnarray*}
l_{1}+u_{1} &=&1,l_{2}=u_{2}=0, \\
a &=&1,b_{1}=l_{1},b_{2}=0, \\
d &=&1,e_{2}=0,e_{1}-u_{1}=1\Rightarrow e_{1}=1,u_{1}=0,b_{1}=l_{1}=1.
\end{eqnarray*}%
Since $\alpha \left( x_{2};1,0,0,0\right) \equiv a+b_{1}\equiv 0,$ we obtain
\begin{equation*}
-B(gx_{1}\otimes x_{2};GX_{1},gx_{1})G\otimes gx_{1}\otimes gx_{1}x_{2}.
\end{equation*}%
The fourth summand of the left side gives us%
\begin{eqnarray*}
l_{1}+u_{1} &=&1,l_{2}+u_{2}=1, \\
a &=&1,b_{1}=l_{1}=1,b_{2}=l_{2} \\
d &=&1,e_{2}=u_{2},e_{1}-u_{1}=1\Rightarrow e_{1}=1,u_{1}=0,b_{1}=l_{1}=1.
\end{eqnarray*}%
Since $\alpha \left( 1_{H};1,0,0,1\right) \equiv a+b_{1}\equiv 0$ and $%
\alpha \left( 1_{H};1,1,0,0\right) \equiv 1+b_{2}\equiv 0$ we obtain
\begin{equation*}
\left[ B(x_{1}x_{2}\otimes x_{2};GX_{1},gx_{1}x_{2})+B(x_{1}x_{2}\otimes
x_{2};GX_{1}X_{2},gx_{1})\right] G\otimes gx_{1}\otimes gx_{1}x_{2}.
\end{equation*}%
Since there is nothing in the right side we get%
\begin{gather*}
B(1_{H}\otimes x_{2};G,gx_{1})-B(gx_{2}\otimes
x_{2};G,gx_{1}x_{2})-B(gx_{2}\otimes x_{2};GX_{2},gx_{1}) \\
-B(gx_{1}\otimes x_{2};GX_{1},gx_{1})+B(x_{1}x_{2}\otimes
x_{2};GX_{1},gx_{1}x_{2})+B(x_{1}x_{2}\otimes x_{2};GX_{1}X_{2},gx_{1})=0
\end{gather*}%
which holds in view of the form of the elements.

\subsubsection{Case $GX_{1}\otimes g\otimes gx_{1}x_{2}$}

First summand of the left side gives us

\begin{eqnarray*}
l_{1} &=&u_{1}=l_{2}=u_{2}=0 \\
a &=&b_{1}=1,b_{2}=0, \\
d &=&1,e_{1}=e_{2}=0.
\end{eqnarray*}%
Since $\alpha \left( x_{1}x_{2};0,0,0,0\right) \equiv 0$ we get
\begin{equation*}
B(1_{H}\otimes x_{2};GX_{1},g)GX_{1}\otimes g\otimes gx_{1}x_{2}.
\end{equation*}%
Second summand of the left side gives us%
\begin{eqnarray*}
l_{1} &=&u_{1}=0,l_{2}+u_{2}=1, \\
a &=&b_{1}=1,b_{2}=l_{2} \\
d &=&1,e_{1}=0,e_{2}=u_{2}.
\end{eqnarray*}%
Since $\alpha \left( x_{1};0,0,0,1\right) \equiv 1$ and $\alpha \left(
x_{1};0,1,0,0\right) \equiv a+b_{1}+b_{2}+1\equiv 0$ we obtain
\begin{equation*}
\left[ -B(gx_{2}\otimes x_{2};GX_{1},gx_{2})+B(gx_{2}\otimes
x_{2};GX_{1}X_{2},g)\right] GX_{1}\otimes g\otimes gx_{1}x_{2}.
\end{equation*}%
Third summand of the left side gives us%
\begin{eqnarray*}
l_{1}+u_{1} &=&1,l_{2}=u_{2}=0, \\
a &=&1,b_{2}=0,b_{1}-l_{1}=1\Rightarrow b_{1}=1,l_{1}=0,u_{1}=1, \\
d &=&1,e_{2}=0,e_{1}=u_{1}=1.
\end{eqnarray*}%
Since$\alpha \left( x_{2};0,0,1,0\right) \equiv e\equiv 0,$ we obtain
\begin{equation*}
-B(gx_{1}\otimes x_{2};GX_{1},gx_{1})GX_{1}\otimes g\otimes gx_{1}x_{2}.
\end{equation*}%
The fourth summand of the left side gives us%
\begin{eqnarray*}
l_{1}+u_{1} &=&1,l_{2}+u_{2}=1, \\
a &=&1,b_{2}=l_{2},b_{1}-l_{1}=1\Rightarrow b_{1}=1,l_{1}=0,u_{1}=1 \\
d &=&1,e_{2}=u_{2},e_{1}=u_{1}=1.
\end{eqnarray*}%
Since $\alpha \left( 1_{H};0,0,1,1\right) \equiv 1+e_{2}\equiv 0$ and $%
\alpha \left( 1_{H};0,1,1,0\right) \equiv e_{2}+a+b_{1}+b_{2}+1\equiv 0$ we
obtain
\begin{equation*}
\left[ B(x_{1}x_{2}\otimes x_{2};GX_{1},gx_{1}x_{2})+B(x_{1}x_{2}\otimes
x_{2};GX_{1}X_{2},gx_{1})\right] GX_{1}\otimes g\otimes gx_{1}x_{2}.
\end{equation*}%
Since there is nothing in the right side we get%
\begin{gather*}
B(1_{H}\otimes x_{2};GX_{1},g)-B(gx_{2}\otimes
x_{2};GX_{1},gx_{2})+B(gx_{2}\otimes x_{2};GX_{1}X_{2},g)+ \\
-B(gx_{1}\otimes x_{2};GX_{1},gx_{1})+B(x_{1}x_{2}\otimes
x_{2};GX_{1},gx_{1}x_{2}) \\
+B(x_{1}x_{2}\otimes x_{2};GX_{1}X_{2},gx_{1})=0
\end{gather*}%
which holds in view of the form of the elements.

\subsection{$B\left( x_{1}x_{2}\otimes x_{2};GX_{2},gx_{1}x_{2}\,\right) $}

We have%
\begin{eqnarray*}
a &=&1,b_{1}=0,b_{2}=1 \\
d &=&e_{1}=e_{2}=1
\end{eqnarray*}%
and we get%
\begin{gather*}
\left( -1\right) ^{\alpha \left( 1_{H};0,0,0,0\right) }B\left(
x_{1}x_{2}\otimes x_{2};GX_{2},gx_{1}x_{2}\,\right) GX_{2}\otimes
gx_{1}x_{2}\otimes g+ \\
\left( -1\right) ^{\alpha \left( 1_{H};0,1,0,0\right) }B\left(
x_{1}x_{2}\otimes x_{2};GX_{2},gx_{1}x_{2}\,\right) G\otimes
gx_{1}x_{2}\otimes x_{2}+ \\
\left( -1\right) ^{\alpha \left( 1_{H};0,0,1,0\right) }B\left(
x_{1}x_{2}\otimes x_{2};GX_{2},gx_{1}x_{2}\,\right) GX_{2}\otimes
gx_{2}\otimes x_{1}+ \\
\left( -1\right) ^{\alpha \left( 1_{H};0,1,1,0\right) }B\left(
x_{1}x_{2}\otimes x_{2};GX_{2},gx_{1}x_{2}\,\right) G\otimes gx_{2}\otimes
gx_{1}x_{2}+ \\
\left( -1\right) ^{\alpha \left( 1_{H};0,0,0,1\right) }B\left(
x_{1}x_{2}\otimes x_{2};GX_{2},gx_{1}x_{2}\,\right) GX_{2}\otimes
gx_{1}\otimes x_{2}+ \\
\left( -1\right) ^{\alpha \left( 1_{H};0,1,0,1\right) }B\left(
x_{1}x_{2}\otimes x_{2};GX_{2},gx_{1}x_{2}\,\right) G\otimes gx_{1}\otimes
g^{l_{2}}x_{2}^{1+1}=0 \\
\left( -1\right) ^{\alpha \left( 1_{H};0,0,1,1\right) }B\left(
x_{1}x_{2}\otimes x_{2};GX_{2},gx_{1}x_{2}\,\right) GX_{2}\otimes g\otimes
gx_{1}x_{2} \\
\left( -1\right) ^{\alpha \left( 1_{H};0,1,1,1\right) }B\left(
x_{1}x_{2}\otimes x_{2};GX_{2},gx_{1}x_{2}\,\right) GX_{2}^{1-l_{2}}\otimes
g\otimes g^{l_{2}+1}x_{1}x_{2}^{1+1}=0.
\end{gather*}

\subsubsection{Case $G\otimes gx_{1}x_{2}\otimes x_{2}$}

We have to consider only the third and the fourth summand of the left side.
Third summand gives us%
\begin{eqnarray*}
l_{1} &=&u_{1}=0,l_{2}=u_{2}=0, \\
a &=&1,b_{1}=b_{2}=0, \\
d &=&e_{1}=e_{2}=1.
\end{eqnarray*}%
Since $\alpha \left( x_{2};0,0,0,0\right) \equiv a+b_{1}+b_{2}\equiv 1,$ we
get
\begin{equation*}
B(gx_{1}\otimes x_{2};G,gx_{1}x_{2})G\otimes gx_{1}x_{2}\otimes x_{2}.
\end{equation*}%
Fourth summand gives us%
\begin{eqnarray*}
l_{1} &=&u_{1}=0,l_{2}+u_{2}=1, \\
a &=&1,b_{1}=0,b_{2}=l_{2}, \\
d &=&e_{1}=1,e_{2}-u_{2}=1\Rightarrow e_{2}=1,u_{2}=0,b_{2}=l_{2}=1.
\end{eqnarray*}%
Since $\alpha \left( 1_{H};0,1,0,0\right) \equiv 0$ we get
\begin{equation*}
B(x_{1}x_{2}\otimes x_{2};GX_{2},gx_{1}x_{2})G\otimes gx_{1}x_{2}\otimes
x_{2}.
\end{equation*}%
By considering also the right side we get%
\begin{equation*}
-B(x_{1}x_{2}\otimes 1_{H};G,gx_{1}x_{2})+B(x_{1}x_{2}\otimes
x_{2};GX_{2},gx_{1}x_{2})+B(gx_{1}\otimes x_{2};G,gx_{1}x_{2})=0
\end{equation*}%
which holds in view of the form of the elements.

\subsubsection{Case $GX_{2}\otimes gx_{2}\otimes x_{1}$}

We have to consider only the second and the fourth summand of the left side
of the equality.

Second summand gives us%
\begin{eqnarray*}
l_{1} &=&u_{1}=l_{2}=u_{2}=0, \\
a &=&b_{2}=1,b_{1}=0, \\
d &=&e_{2}=1,e_{1}=0.
\end{eqnarray*}%
Since $\alpha \left( x_{1};0,0,0,0\right) \equiv a+b_{1}+b_{2}\equiv 0$ we
obtain
\begin{equation*}
B(gx_{2}\otimes x_{2};GX_{2},gx_{2})GX_{2}\otimes gx_{2}\otimes x_{1}.
\end{equation*}%
Forth summand gives us%
\begin{eqnarray*}
l_{1}+u_{1} &=&1,l_{2}=u_{2}=0 \\
a &=&b_{2}=1,b_{1}=l_{1}, \\
d &=&e_{2}=1,e_{1}=u_{1}.
\end{eqnarray*}%
Since $\alpha \left( 1_{H};0,0,1,0\right) \equiv e_{2}+\left(
a+b_{1}+b_{2}\right) \equiv 1$ and $\alpha \left( 1_{H};1,0,0,0\right)
\equiv b_{2}=1,$ we obtain
\begin{equation*}
\left[ -B(x_{1}x_{2}\otimes x_{2};GX_{2},gx_{1}x_{2})-B(x_{1}x_{2}\otimes
x_{2};GX_{1}X_{2},gx_{2})\right] GX_{2}\otimes gx_{2}\otimes x_{1}.
\end{equation*}%
Since there is no term in the right side, we get%
\begin{equation*}
B(gx_{2}\otimes x_{2};GX_{2},gx_{2})-B(x_{1}x_{2}\otimes
x_{2};GX_{2},gx_{1}x_{2})-B(x_{1}x_{2}\otimes x_{2};GX_{1}X_{2},gx_{2})=0
\end{equation*}%
which holds in view of the form of the elements.

\subsubsection{Case $G\otimes gx_{2}\otimes gx_{1}x_{2}$}

First summand of the left side gives us

\begin{eqnarray*}
l_{1} &=&u_{1}=l_{2}=u_{2}=0 \\
a &=&1,b_{1}=b_{2}=0, \\
d &=&e_{2}=1,e_{1}=0.
\end{eqnarray*}%
Since $\alpha \left( x_{1}x_{2};0,0,0,0\right) \equiv 0$ we get
\begin{equation*}
B(1_{H}\otimes x_{2};G,gx_{2})G\otimes gx_{2}\otimes gx_{1}x_{2}
\end{equation*}%
Second summand of the left side gives us%
\begin{eqnarray*}
l_{1} &=&u_{1}=0,l_{2}+u_{2}=1, \\
a &=&1,b_{1}=0,b_{2}=l_{2} \\
d &=&1,e_{1}=0,e_{2}-u_{2}=1\Rightarrow e_{2}=1,u_{2}=0,b_{2}=l_{2}=1.
\end{eqnarray*}%
Since $\alpha \left( x_{1};0,1,0,0\right) \equiv a+b_{1}+b_{2}+1\equiv 1$ we
obtain
\begin{equation*}
-B(gx_{2}\otimes x_{2};GX_{2},gx_{2})G\otimes gx_{2}\otimes gx_{1}x_{2}.
\end{equation*}%
Third summand of the left side gives us%
\begin{eqnarray*}
l_{1}+u_{1} &=&1,l_{2}=u_{2}=0, \\
a &=&1,b_{1}=l_{1},b_{2}=0, \\
d &=&e_{2}=1,0,e_{1}=u_{1}.
\end{eqnarray*}%
Since $\alpha \left( x_{2};0,0,1,0\right) \equiv e_{2}=1$ and $\alpha \left(
x_{2};1,0,0,0\right) \equiv a+b_{1}\equiv 0,$ we obtain
\begin{equation*}
\left[ +B(gx_{1}\otimes x_{2};G,gx_{1}x_{2})-B(gx_{1}\otimes
x_{2};GX_{1},gx_{2})\right] G\otimes gx_{2}\otimes gx_{1}x_{2}.
\end{equation*}%
The fourth summand of the left side gives us%
\begin{eqnarray*}
l_{1}+u_{1} &=&1,l_{2}+u_{2}=1, \\
a &=&1,b_{1}=l_{1},b_{2}=l_{2} \\
d &=&1,e_{1}=u_{1},e_{2}-u_{2}=1\Rightarrow e_{2}=1,u_{2}=0,b_{2}=l_{2}=1.
\end{eqnarray*}%
Since $\alpha \left( 1_{H};0,1,1,0\right) \equiv e_{2}+a+b_{1}+b_{2}+1\equiv
0$ and $\alpha \left( 1_{H};1,1,0,0\right) \equiv 1+b_{2}\equiv 0$ we obtain
\begin{equation*}
\left[ B(x_{1}x_{2}\otimes x_{2};GX_{2},gx_{1}x_{2})+B(x_{1}x_{2}\otimes
x_{2};GX_{1}X_{2},gx_{2})\right] G\otimes gx_{2}\otimes gx_{1}x_{2}.
\end{equation*}%
Since there is nothing in the right side we get%
\begin{gather*}
B(1_{H}\otimes x_{2};G,gx_{2})-B(gx_{2}\otimes x_{2};GX_{2},gx_{2})+ \\
+B(gx_{1}\otimes x_{2};G,gx_{1}x_{2})-B(gx_{1}\otimes x_{2};GX_{1},gx_{2}) \\
+B(x_{1}x_{2}\otimes x_{2};GX_{2},gx_{1}x_{2})+B(x_{1}x_{2}\otimes
x_{2};GX_{1}X_{2},gx_{2})=0
\end{gather*}%
which holds in view of the form of the elements.

\subsubsection{Case $GX_{2}\otimes gx_{1}\otimes x_{2}$}

We have to consider only the third and the fourth summand of the left side.
Third summand gives us%
\begin{eqnarray*}
l_{1} &=&u_{1}=0,l_{2}=u_{2}=0, \\
a &=&b_{2}=1,b_{1}=0, \\
d &=&e_{1}=1,e_{2}=0.
\end{eqnarray*}%
Since $\alpha \left( x_{2};0,0,0,0\right) \equiv a+b_{1}+b_{2}\equiv 0,$ we
get
\begin{equation*}
-B(gx_{1}\otimes x_{2};GX_{2},gx_{1})GX_{2}\otimes gx_{1}\otimes x_{2}.
\end{equation*}%
Fourth summand gives us%
\begin{eqnarray*}
l_{1} &=&u_{1}=0,l_{2}+u_{2}=1, \\
a &=&1,b_{1}=0,b_{2}-l_{2}=1\Rightarrow b_{2}=1,l_{2}=0,u_{2}=1 \\
d &=&e_{1}=1,e_{2}=u_{2}=1.
\end{eqnarray*}%
Since $\alpha \left( 1_{H};0,0,0,1\right) \equiv a+b_{1}+b_{2}\equiv 0$ we
get
\begin{equation*}
B(x_{1}x_{2}\otimes x_{2};GX_{2},gx_{1}x_{2})GX_{2}\otimes gx_{1}\otimes
x_{2}.
\end{equation*}%
By considering also the right side we get%
\begin{equation*}
-B(x_{1}x_{2}\otimes 1_{H};GX_{2},gx_{1})+B(x_{1}x_{2}\otimes
x_{2};GX_{2},gx_{1}x_{2})-B(gx_{1}\otimes x_{2};GX_{2},gx_{1})=0
\end{equation*}%
which holds in view of the form of the elements.

\subsubsection{Case $GX_{2}\otimes g\otimes gx_{1}x_{2}$}

First summand of the left side gives us

\begin{eqnarray*}
l_{1} &=&u_{1}=l_{2}=u_{2}=0 \\
a &=&b_{2}=1,b_{1}=0, \\
d &=&1,e_{1}=e_{2}=0.
\end{eqnarray*}%
Since $\alpha \left( x_{1}x_{2};0,0,0,0\right) \equiv 0$ we get
\begin{equation*}
B(1_{H}\otimes x_{2};GX_{2},g)GX_{2}\otimes g\otimes gx_{1}x_{2}.
\end{equation*}%
Second summand of the left side gives us%
\begin{eqnarray*}
l_{1} &=&u_{1}=0,l_{2}+u_{2}=1, \\
a &=&1,b_{1}=0,b_{2}-l_{2}=1\Rightarrow b_{2}=1,l_{2}=0,u_{2}=1 \\
d &=&1,e_{1}=0,e_{2}=u_{2}=1.
\end{eqnarray*}%
Since $\alpha \left( x_{1};0,0,0,1\right) \equiv 1$ we obtain
\begin{equation*}
-B(gx_{2}\otimes x_{2};GX_{2},gx_{2})GX_{2}\otimes g\otimes gx_{1}x_{2}.
\end{equation*}%
Third summand of the left side gives us%
\begin{eqnarray*}
l_{1}+u_{1} &=&1,l_{2}=u_{2}=0, \\
a &=&b_{2}=1,b_{1}=l_{1}, \\
d &=&1,e_{2}=0,e_{1}=u_{1}.
\end{eqnarray*}%
Since $\alpha \left( x_{2};0,0,1,0\right) \equiv e_{2}=0$ and $\alpha \left(
x_{2};1,0,0,0\right) \equiv a+b_{1}\equiv 0,$ we obtain
\begin{equation*}
\left[ -B(gx_{1}\otimes x_{2};GX_{2},gx_{1})-B(gx_{1}\otimes
x_{2};GX_{1}X_{2},g)\right] GX_{2}\otimes g\otimes gx_{1}x_{2}.
\end{equation*}%
The fourth summand of the left side gives us%
\begin{eqnarray*}
l_{1}+u_{1} &=&1,l_{2}+u_{2}=1, \\
a &=&1,b_{1}=l_{1},b_{2}=l_{2} \\
d &=&1,e_{1}=u_{1},e_{2}-u_{2}=1\Rightarrow e_{2}=1,u_{2}=0,b_{2}=l_{2}=1.
\end{eqnarray*}%
Since $\alpha \left( 1_{H};0,1,1,0\right) \equiv e_{2}+a+b_{1}+b_{2}+1\equiv
0$ and $\alpha \left( 1_{H};1,1,0,0\right) \equiv 1+b_{2}\equiv 0$ we obtain
\begin{equation*}
\left[ B(x_{1}x_{2}\otimes x_{2};GX_{2},gx_{1}x_{2})+B(x_{1}x_{2}\otimes
x_{2};GX_{1}X_{2},gx_{2})\right] G\otimes gx_{2}\otimes gx_{1}x_{2}.
\end{equation*}%
Since there is nothing in the right side we get%
\begin{gather*}
B(1_{H}\otimes x_{2};GX_{2},g)-B(gx_{2}\otimes x_{2};GX_{2},gx_{2})+ \\
-B(gx_{1}\otimes x_{2};GX_{2},gx_{1})-B(gx_{1}\otimes x_{2};GX_{1}X_{2},g) \\
+B(x_{1}x_{2}\otimes x_{2};GX_{2},gx_{1}x_{2})+B(x_{1}x_{2}\otimes
x_{2};GX_{1}X_{2},gx_{2})=0
\end{gather*}%
which holds in view of the form of the elements.

\subsection{$B(x_{1}x_{2}\otimes x_{2};GX_{1}X_{2},gx_{2})$}

We deduce that
\begin{equation*}
a=b_{1}=b_{2}=1,d=e_{2}=1
\end{equation*}%
and we get%
\begin{gather*}
\left( -1\right) ^{\alpha \left( 1_{H};0,0,0,0\right) }B(x_{1}x_{2}\otimes
x_{2};GX_{1}X_{2},gx_{2})GX_{1}X_{2}\otimes gx_{2}\otimes g \\
\left( -1\right) ^{\alpha \left( 1_{H};1,0,0,0\right) }B(x_{1}x_{2}\otimes
x_{2};GX_{1}X_{2},gx_{2})GX_{2}\otimes gx_{2}\otimes x_{1} \\
\left( -1\right) ^{\alpha \left( 1_{H};0,1,0,0\right) }B(x_{1}x_{2}\otimes
x_{2};GX_{1}X_{2},gx_{2})GX_{1}\otimes gx_{2}\otimes x_{2} \\
\left( -1\right) ^{\alpha \left( 1_{H};1,1,0,0\right) }B(x_{1}x_{2}\otimes
x_{2};GX_{1}X_{2},gx_{2})G\otimes gx_{2}\otimes gx_{1}x_{2} \\
\left( -1\right) ^{\alpha \left( 1_{H};0,0,0,1\right) }B(x_{1}x_{2}\otimes
x_{2};GX_{1}X_{2},gx_{2})GX_{1}X_{2}\otimes g\otimes x_{2} \\
\left( -1\right) ^{\alpha \left( 1_{H};1,0,0,1\right) }B(x_{1}x_{2}\otimes
x_{2};GX_{1}X_{2},gx_{2})GX_{2}\otimes g\otimes gx_{1}x_{2} \\
\left( -1\right) ^{\alpha \left( 1_{H};0,1,0,1\right) }B(x_{2}\otimes
x_{1};GX_{1}X_{2},gx_{2})GX_{1}\otimes g\otimes gx_{2}^{1+1}=0 \\
\left( -1\right) ^{\alpha \left( 1_{H};1,1,0,1\right) }B(x_{2}\otimes
x_{1};GX_{1}X_{2},gx_{2})GX_{1}^{1-l_{1}}X_{2}^{1-l_{2}}\otimes
gx_{2}^{1-u_{2}}\otimes g^{1+l_{1}+l_{2}+u_{2}}x_{1}x_{2}^{1+1}=0.
\end{gather*}

\subsubsection{Case $GX_{2}\otimes gx_{2}\otimes x_{1}$}

This case was already considered in subsection $B\left( x_{1}x_{2}\otimes
x_{2};GX_{2},gx_{1}x_{2}\,\right) .$

\subsubsection{Case $GX_{1}\otimes gx_{2}\otimes x_{2}$}

We have to consider only the third and the fourth summand of the left side.
Third summand gives us%
\begin{eqnarray*}
l_{1} &=&u_{1}=0,l_{2}=u_{2}=0, \\
a &=&b_{1}=1,b_{2}=0, \\
d &=&e_{2}=1,e_{1}=0.
\end{eqnarray*}%
Since $\alpha \left( x_{2};0,0,0,0\right) \equiv a+b_{1}+b_{2}\equiv 0,$ we
get
\begin{equation*}
-B(gx_{1}\otimes x_{2};GX_{1},gx_{2})GX_{1}\otimes gx_{2}\otimes x_{2}.
\end{equation*}%
Fourth summand gives us%
\begin{eqnarray*}
l_{1} &=&u_{1}=0,l_{2}+u_{2}=1, \\
a &=&b_{1}=1,0,b_{2}=l_{2} \\
d &=&1,e_{1}=0,e_{2}-u_{2}=1\Rightarrow e_{2}=1,u_{2}=0,b_{2}=l_{2}=1.
\end{eqnarray*}%
Since $\alpha \left( 1_{H};0,1,0,0\right) \equiv 0$ we get
\begin{equation*}
B(x_{1}x_{2}\otimes x_{2};GX_{1}X_{2},gx_{2})GX_{1}\otimes gx_{2}\otimes
x_{2}.
\end{equation*}%
By considering also the right side we get%
\begin{equation*}
-B(x_{1}x_{2}\otimes 1;GX_{1},gx_{2})+B(x_{1}x_{2}\otimes
x_{2};GX_{1}X_{2},gx_{2})-B(gx_{1}\otimes x_{2};GX_{1},gx_{2})=0
\end{equation*}%
which holds in view of the form of the elements.

\subsubsection{Case $G\otimes gx_{2}\otimes gx_{1}x_{2}$}

This case was already considered in subsection $B\left( x_{1}x_{2}\otimes
x_{2};GX_{2},gx_{1}x_{2}\,\right) .$

\subsubsection{Case $GX_{1}X_{2}\otimes g\otimes x_{2}$}

We have to consider only the third and the fourth summand of the left side.
Third summand gives us%
\begin{eqnarray*}
l_{1} &=&u_{1}=0,l_{2}=u_{2}=0, \\
a &=&b_{1}=b_{2}=1, \\
d &=&1,e_{1}=e_{2}=0.
\end{eqnarray*}%
Since $\alpha \left( x_{2};0,0,0,0\right) \equiv a+b_{1}+b_{2}\equiv 1,$ we
get
\begin{equation*}
+B(gx_{1}\otimes x_{2};GX_{1}X_{2},g)GX_{1}X_{2}\otimes g\otimes x_{2}.
\end{equation*}%
Fourth summand gives us%
\begin{eqnarray*}
l_{1} &=&u_{1}=0,l_{2}+u_{2}=1, \\
a &=&b_{1}=1,0,b_{2}-l_{2}=1\Rightarrow b_{2}=1,l_{2}=0,u_{2}=1 \\
d &=&1,e_{1}=0,e_{2}=u_{2}=1.
\end{eqnarray*}%
Since $\alpha \left( 1_{H};0,0,0,1\right) \equiv a+b_{1}+b_{2}\equiv 1$ we
get
\begin{equation*}
-B(x_{1}x_{2}\otimes x_{2};GX_{1}X_{2},gx_{2})GX_{1}X_{2}\otimes g\otimes
x_{2}.
\end{equation*}%
By considering also the right side we get%
\begin{equation*}
-B(x_{1}x_{2}\otimes 1;GX_{1}X_{2},g)-B(x_{1}x_{2}\otimes
x_{2};GX_{1}X_{2},gx_{2})+B(gx_{1}\otimes x_{2};GX_{1}X_{2},g)=0
\end{equation*}%
which holds in view of the form of the elements.

\subsubsection{Case $GX_{2}\otimes g\otimes gx_{1}x_{2}$}

This case was already considered in subsection $B\left( x_{1}x_{2}\otimes
x_{2};GX_{2},gx_{1}x_{2}\,\right) .$

\section{$B(x_{1}x_{2}\otimes x_{1}x_{2})$}

By using formula $\left( \ref{1x}\right) $ for $h=x_{1}x_{2},h^{\prime
}=x_{2},i=1$ we get%
\begin{equation*}
B(x_{1}x_{2}\otimes x_{1}x_{2})=B(x_{1}x_{2}\otimes x_{2})(1_{A}\otimes
x_{1})-(1_{A}\otimes gx_{1})B(x_{1}x_{2}\otimes x_{2})(1_{A}\otimes g).
\end{equation*}%
Since%
\begin{equation*}
B(x_{1}x_{2}\otimes x_{2})=B(x_{1}x_{2}\otimes 1_{H})(1_{A}\otimes
x_{2})-(1_{A}\otimes gx_{2})B(x_{1}x_{2}\otimes 1_{H})(1_{A}\otimes g)\left( %
\ref{x1x2otx2}\right)
\end{equation*}%
we obtain%
\begin{eqnarray*}
B(x_{1}x_{2}\otimes x_{1}x_{2}) &=&B(x_{1}x_{2}\otimes 1_{H})(1_{A}\otimes
x_{2})(1_{A}\otimes x_{1})-(1_{A}\otimes gx_{2})B(x_{1}x_{2}\otimes
1_{H})(1_{A}\otimes g)(1_{A}\otimes x_{1}) \\
&&-\left[
\begin{array}{c}
(1_{A}\otimes gx_{1})B(x_{1}x_{2}\otimes 1_{H})(1_{A}\otimes
x_{2})(1_{A}\otimes g) \\
-(1_{A}\otimes gx_{1})(1_{A}\otimes gx_{2})B(x_{1}x_{2}\otimes
1_{H})(1_{A}\otimes g)(1_{A}\otimes g)%
\end{array}%
\right]
\end{eqnarray*}%
i.e.%
\begin{eqnarray}
B(x_{1}x_{2}\otimes x_{1}x_{2}) &=&B(x_{1}x_{2}\otimes 1_{H})(1_{A}\otimes
x_{1}x_{2})+(1_{A}\otimes gx_{2})B(x_{1}x_{2}\otimes 1_{H})(1_{A}\otimes
gx_{1})  \label{form x1x2otx1x2} \\
&&-(1_{A}\otimes gx_{1})B(x_{1}x_{2}\otimes 1_{H})(1_{A}\otimes gx_{2})+
\notag \\
&&+(1_{A}\otimes x_{1}x_{2})B(x_{1}x_{2}\otimes 1_{H})+  \notag
\end{eqnarray}

\begin{gather*}
B(x_{1}x_{2}\otimes x_{1}x_{2})=+4B(x_{1}x_{2}\otimes 1_{H};G,g)G\otimes
gx_{1}x_{2}+ \\
+4B(x_{1}x_{2}\otimes 1_{H};X_{1},g)X_{1}\otimes gx_{1}x_{2}+ \\
+4B(x_{1}x_{2}\otimes 1_{H};X_{2},g)X_{2}\otimes gx_{1}x_{2}+ \\
+\left[
\begin{array}{c}
-4B(x_{1}x_{2}\otimes 1_{H};G,gx_{1}x_{2})+ \\
4B(x_{1}x_{2}\otimes 1_{H};GX_{2},gx_{1})-4B(x_{1}x_{2}\otimes
1_{H};GX_{1},gx_{2})%
\end{array}%
\right] GX_{1}X_{2}\otimes gx_{1}x_{2}
\end{gather*}

We write the Casimir condition equality for $B\left( x_{1}x_{2}\otimes
x_{1}x_{2}\right) .$%
\begin{eqnarray*}
&&\text{1}\sum_{a,b_{1},b_{2},d,e_{1},e_{2}=0}^{1}\sum_{l_{1}=0}^{b_{1}}%
\sum_{l_{2}=0}^{b_{2}}\sum_{u_{1}=0}^{e_{1}}\sum_{u_{2}=0}^{e_{2}}\left(
-1\right) ^{\alpha \left( x_{1}x_{2};l_{1},l_{2},u_{1},u_{2}\right) } \\
&&B(1_{H}\otimes
x_{1}x_{2};G^{a}X_{1}^{b_{1}}X_{2}^{b_{2}},g^{d}x_{1}^{e_{1}}x_{2}^{e_{2}})
\\
&&G^{a}X_{1}^{b_{1}-l_{1}}X_{2}^{b_{2}-l_{2}}\otimes
g^{d}x_{1}^{e_{1}-u_{1}}x_{2}^{e_{2}-u_{2}}\otimes
g^{a+b_{1}+b_{2}+l_{1}+l_{2}+d+e_{1}+e_{2}+u_{1}+u_{2}}x_{1}^{l_{1}+u_{1}+1}x_{2}^{l_{2}+u_{2}+1}+
\\
&&\text{2}+\sum_{a,b_{1},b_{2},d,e_{1},e_{2}=0}^{1}\sum_{l_{1}=0}^{b_{1}}%
\sum_{l_{2}=0}^{b_{2}}\sum_{u_{1}=0}^{e_{1}}\sum_{u_{2}=0}^{e_{2}}\left(
-1\right) ^{\alpha \left( x_{1};l_{1},l_{2},u_{1},u_{2}\right) } \\
&&B(gx_{2}\otimes
x_{1}x_{2};G^{a}X_{1}^{b_{1}}X_{2}^{b_{2}},g^{d}x_{1}^{e_{1}}x_{2}^{e_{2}})
\\
&&G^{a}X_{1}^{b_{1}-l_{1}}X_{2}^{b_{2}-l_{2}}\otimes
g^{d}x_{1}^{e_{1}-u_{1}}x_{2}^{e_{2}-u_{2}}\otimes
g^{a+b_{1}+b_{2}+l_{1}+l_{2}+d+e_{1}+e_{2}+u_{1}+u_{2}}x_{1}^{l_{1}+u_{1}+1}x_{2}^{l_{2}+u_{2}}+
\\
&&\text{3}\left( -1\right)
\sum_{a,b_{1},b_{2},d,e_{1},e_{2}=0}^{1}\sum_{l_{1}=0}^{b_{1}}%
\sum_{l_{2}=0}^{b_{2}}\sum_{u_{1}=0}^{e_{1}}\sum_{u_{2}=0}^{e_{2}}\left(
-1\right) ^{\alpha \left( x_{2};l_{1},l_{2},u_{1},u_{2}\right) } \\
&&B(gx_{1}\otimes
x_{1}x_{2};G^{a}X_{1}^{b_{1}}X_{2}^{b_{2}},g^{d}x_{1}^{e_{1}}x_{2}^{e_{2}})
\\
&&G^{a}X_{1}^{b_{1}-l_{1}}X_{2}^{b_{2}-l_{2}}\otimes
g^{d}x_{1}^{e_{1}-u_{1}}x_{2}^{e_{2}-u_{2}}\otimes
g^{a+b_{1}+b_{2}+l_{1}+l_{2}+d+e_{1}+e_{2}+u_{1}+u_{2}}x_{1}^{l_{1}+u_{1}}x_{2}^{l_{2}+u_{2}+1}+
\\
&&\text{4}\sum_{a,b_{1},b_{2},d,e_{1},e_{2}=0}^{1}\sum_{l_{1}=0}^{b_{1}}%
\sum_{l_{2}=0}^{b_{2}}\sum_{u_{1}=0}^{e_{1}}\sum_{u_{2}=0}^{e_{2}}\left(
-1\right) ^{\alpha \left( 1_{H};l_{1},l_{2},u_{1},u_{2}\right) } \\
&&B(x_{1}x_{2}\otimes
x_{1}x_{2};G^{a}X_{1}^{b_{1}}X_{2}^{b_{2}},g^{d}x_{1}^{e_{1}}x_{2}^{e_{2}})
\\
&&G^{a}X_{1}^{b_{1}-l_{1}}X_{2}^{b_{2}-l_{2}}\otimes
g^{d}x_{1}^{e_{1}-u_{1}}x_{2}^{e_{2}-u_{2}}\otimes
g^{a+b_{1}+b_{2}+l_{1}+l_{2}+d+e_{1}+e_{2}+u_{1}+u_{2}}x_{1}^{l_{1}+u_{1}}x_{2}^{l_{2}+u_{2}}
\\
&=&B^{A}(x_{1}x_{2}\otimes x_{1}x_{2})\otimes B^{H}(x_{1}x_{2}\otimes
x_{1}x_{2})\otimes 1_{H} \\
&&B^{A}(x_{1}x_{2}\otimes x_{1})\otimes B^{H}(x_{1}x_{2}\otimes
x_{1})\otimes gx_{2} \\
&&-B^{A}(x_{1}x_{2}\otimes x_{2})\otimes B^{H}(x_{1}x_{2}\otimes
x_{2})\otimes gx_{1} \\
&&B^{A}(x_{1}x_{2}\otimes 1_{H})\otimes B^{H}(x_{1}x_{2}\otimes
1_{H})\otimes x_{1}x_{2}
\end{eqnarray*}

\subsection{$B\left( x_{1}x_{2}\otimes x_{1}x_{2};G\otimes
gx_{1}x_{2}\right) $}

We deduce that%
\begin{eqnarray*}
a &=&1 \\
b_{1} &=&b_{2}=0 \\
d &=&e_{1}=e_{2}=1 \\
a+b_{1}+b_{2}+l_{1}+l_{2}+d+e_{1}+e_{2}+u_{1}+u_{2} &=&u_{1}+u_{2}
\end{eqnarray*}

and we get%
\begin{eqnarray*}
&&\sum_{u_{1}=0}^{1}\sum_{u_{2}=0}^{1}\left( -1\right) ^{\alpha \left(
1_{H};0,0,u_{1},u_{2}\right) }B(x_{1}x_{2}\otimes
x_{1}x_{2};G^{a}X_{1}^{b_{1}}X_{2}^{b_{2}},g^{d}x_{1}^{e_{1}}x_{2}^{e_{2}})
\\
&&G\otimes gx_{1}^{1-u_{1}}x_{2}^{1-u_{2}}\otimes
g^{u_{1}+u_{2}}x_{1}^{u_{1}}x_{2}^{u_{2}} \\
&=&\left( -1\right) ^{\alpha \left( 1_{H};0,0,0,0\right)
}B(x_{1}x_{2}\otimes
x_{1}x_{2};G^{a}X_{1}^{b_{1}}X_{2}^{b_{2}},g^{d}x_{1}^{e_{1}}x_{2}^{e_{2}})G%
\otimes gx_{1}x_{2}\otimes 1_{H}\text{ } \\
&&\left( -1\right) ^{\alpha \left( 1_{H};0,0,0,1\right) }B(x_{1}x_{2}\otimes
x_{1}x_{2};G^{a}X_{1}^{b_{1}}X_{2}^{b_{2}},g^{d}x_{1}^{e_{1}}x_{2}^{e_{2}})G%
\otimes gx_{1}\otimes gx_{2}\text{ } \\
&&\left( -1\right) ^{\alpha \left( 1_{H};0,0,1,0\right) }B(x_{1}x_{2}\otimes
x_{1}x_{2};G^{a}X_{1}^{b_{1}}X_{2}^{b_{2}},g^{d}x_{1}^{e_{1}}x_{2}^{e_{2}})G%
\otimes gx_{2}\otimes gx_{1} \\
&&\left( -1\right) ^{\alpha \left( 1_{H};0,0,1,1\right) }B(x_{1}x_{2}\otimes
x_{1}x_{2};G^{a}X_{1}^{b_{1}}X_{2}^{b_{2}},g^{d}x_{1}^{e_{1}}x_{2}^{e_{2}})G%
\otimes g\otimes x_{1}x_{2}.
\end{eqnarray*}

\subsubsection{Case $G\otimes gx_{1}\otimes gx_{2}$}

We have only to consider the third and the fourth summands of the left side
of the equality.

Third summand gives us

\begin{eqnarray*}
l_{1} &=&u_{1}=0,l_{2}=u_{2}=0 \\
a &=&1,b_{1}=b_{2}=0
\end{eqnarray*}%
Since $\alpha \left( x_{2};0,0,0,0\right) \equiv a+b_{1}+b_{2}\equiv 1,$ we
get
\begin{equation*}
B(gx_{1}\otimes x_{1}x_{2};G,gx_{1})G\otimes gx_{1}\otimes gx_{2}.
\end{equation*}%
Fourth summand gives us%
\begin{eqnarray*}
l_{1} &=&u_{1}=0,l_{2}+u_{2}=1 \\
a &=&1,b_{1}=0,b_{2}=l_{2}, \\
d &=&1,e_{1}=1,e_{2}=u_{2}.
\end{eqnarray*}%
Since $\alpha \left( 1_{H};0,0,0,1\right) \equiv a+b_{1}+b_{2}=1$ and $%
\alpha \left( 1_{H};0,1,0,0\right) \equiv 0,$ we get%
\begin{equation*}
\left[ -B(x_{1}x_{2}\otimes x_{1}x_{2};G,gx_{1}x_{2})+B(x_{1}x_{2}\otimes
x_{1}x_{2};GX_{2},gx)\right] G\otimes gx_{1}\otimes gx_{2}.
\end{equation*}%
By considering also the right side we obtain%
\begin{gather*}
-B(x_{1}x_{2}\otimes x_{1};G,gx_{1})-B(x_{1}x_{2}\otimes
x_{1}x_{2};G,gx_{1}x_{2}) \\
+B(x_{1}x_{2}\otimes x_{1}x_{2};GX_{2},gx)+B(gx_{1}\otimes
x_{1}x_{2};G,gx_{1})=0
\end{gather*}%
which holds in view of the form of the elements.

\subsubsection{Case $G\otimes gx_{2}\otimes gx_{1}$}

We have only to consider the second and the forth summands of the left side
of the equality.

Second summand gives us

\begin{eqnarray*}
l_{1} &=&u_{1}=0,l_{2}=u_{2}=0, \\
a &=&1,b_{1}=b_{2}=0, \\
d &=&e_{2}=1,e_{1}=0
\end{eqnarray*}%
Since $\alpha \left( x_{1};0,0,0,0\right) \equiv a+b_{1}+b_{2}\equiv 1,$ we
get%
\begin{equation*}
-B(gx_{2}\otimes x_{1}x_{2};G,gx_{2})G\otimes gx_{2}\otimes gx_{1}.
\end{equation*}

Fourth summand gives us%
\begin{eqnarray*}
l_{1}+u_{1} &=&1,l_{2}=u_{2}=0 \\
a &=&1,b_{1}=l_{1},b_{2}=0 \\
d &=&e_{2}=1,e_{1}=u_{1}
\end{eqnarray*}%
Since $\alpha \left( 1_{H};0,0,1,0\right) \equiv e_{2}+\left(
a+b_{1}+b_{2}\right) \equiv 0$ and $\alpha \left( 1_{H};1,0,0,0\right)
\equiv b_{2}=0,$ we obtain%
\begin{equation*}
\left[ B(x_{1}x_{2}\otimes x_{1}x_{2};G,gx_{1}x_{2})+B(x_{1}x_{2}\otimes
x_{1}x_{2};GX_{1},gx_{2})\right] G\otimes gx_{2}\otimes gx_{1}
\end{equation*}%
By considering also the right side, we get%
\begin{gather*}
B(x_{1}x_{2}\otimes x_{2};G,gx_{2})+B(x_{1}x_{2}\otimes
x_{1}x_{2};G,gx_{1}x_{2}) \\
+B(x_{1}x_{2}\otimes x_{1}x_{2};GX_{1},gx_{2})-B(gx_{2}\otimes
x_{1}x_{2};G,gx_{2})=0
\end{gather*}%
which holds in view of the form of the elements.

\subsubsection{Case $G\otimes g\otimes x_{1}x_{2}$}

First summand of the left side of the equality gives%
\begin{eqnarray*}
l_{1} &=&u_{1}=0,l_{2}=u_{2}=0 \\
a &=&1,b_{1}=b_{2}=0 \\
d &=&1,e_{1}=e_{2}=0.
\end{eqnarray*}%
Since $\alpha \left( x_{1}x_{2};0,0,0,0\right) \equiv 0,$ we get%
\begin{equation*}
B(1_{H}\otimes x_{1}x_{2};G,g)G\otimes g\otimes x_{1}x_{2}.
\end{equation*}%
Second summand of the left side of the equality gives us%
\begin{eqnarray*}
l_{1} &=&u_{1}=0,l_{2}+u_{2}=1 \\
a &=&1,b_{1}=0,b_{2}=l_{2} \\
d &=&1,e_{1}=0,e_{2}=u_{2}.
\end{eqnarray*}%
Since $\alpha \left( x_{1};0,0,0,1\right) \equiv 1$ and $\alpha \left(
x_{1};0,1,0,0\right) \equiv a+b_{1}+b_{2}+1\equiv 1$, we get%
\begin{equation*}
\left[ -B(gx_{2}\otimes x_{1}x_{2};G,gx_{2})-B(gx_{2}\otimes
x_{1}x_{2};GX_{2},g)\right] G\otimes g\otimes x_{1}x_{2}.
\end{equation*}%
Third summand of the left side of the equality gives us

\begin{eqnarray*}
l_{1}+u_{1} &=&1,l_{2}=u_{2}=0 \\
a &=&1,b_{1}=l_{1},b_{2}=0 \\
d &=&1,e_{1}=u_{1},e_{2}=0,
\end{eqnarray*}%
Since $\alpha \left( x_{2};0,0,1,0\right) \equiv e_{2}=0$ and $\alpha \left(
x_{2};1,0,0,0\right) \equiv a+b_{1}\equiv 0,$ we get%
\begin{equation*}
\left[ -B(gx_{1}\otimes x_{1}x_{2};G,gx_{1})-B(gx_{1}\otimes
x_{1}x_{2};GX_{1},g)\right] G\otimes g\otimes x_{1}x_{2}.
\end{equation*}%
Fourth summand gives us%
\begin{eqnarray*}
l_{1}+u_{1} &=&1,l_{2}+u_{2}=1, \\
a &=&1,b_{1}=l_{1},b_{2}=l_{2}, \\
d &=&1,e_{1}=u_{1},e_{2}=u_{2}.
\end{eqnarray*}%
Since
\begin{eqnarray*}
\alpha \left( 1_{H};0,0,1,1\right) &\equiv &1+e_{2}\equiv 0 \\
\alpha \left( 1_{H};0,1,1,0\right) &\equiv &e_{2}+\left(
a+b_{1}+b_{2}+1\right) \equiv 1 \\
\alpha \left( 1_{H};1,0,0,1\right) &\equiv &a+b_{1}\equiv 0 \\
\alpha \left( 1_{H};1,1,0,0\right) &\equiv &1+b_{2}\equiv 0
\end{eqnarray*}%
we obtain
\begin{equation*}
\left[
\begin{array}{c}
B(x_{1}x_{2}\otimes x_{1}x_{2};G,gx_{1}x_{2})-B(x_{1}x_{2}\otimes
x_{1}x_{2};GX_{2},gx_{1})+ \\
+B(x_{1}x_{2}\otimes x_{1}x_{2};GX_{1},gx_{2})+B(x_{1}x_{2}\otimes
x_{1}x_{2};GX_{1}X_{2},g)%
\end{array}%
\right] G\otimes g\otimes x_{1}x_{2}.
\end{equation*}%
By taking in account also the right side of the equality we get%
\begin{gather*}
-B(x_{1}x_{2}\otimes 1_{H};G,g)+B(x_{1}x_{2}\otimes
x_{1}x_{2};G,gx_{1}x_{2})-B(x_{1}x_{2}\otimes x_{1}x_{2};GX_{2},gx_{1})+ \\
+B(x_{1}x_{2}\otimes x_{1}x_{2};GX_{1},gx_{2})+B(x_{1}x_{2}\otimes
x_{1}x_{2};GX_{1}X_{2},g)-B(gx_{1}\otimes x_{1}x_{2};G,gx_{1}) \\
-B(gx_{1}\otimes x_{1}x_{2};GX_{1},g)-B(gx_{2}\otimes
x_{1}x_{2};G,gx_{2})-B(gx_{2}\otimes x_{1}x_{2};GX_{2},g)+ \\
B(1_{H}\otimes x_{1}x_{2};G,g)=0
\end{gather*}

which holds in view of the form of the elements.

\subsection{$B\left( x_{1}x_{2}\otimes x_{1}x_{2};X_{1},gx_{1}x_{2}\right) $}

We deduce that%
\begin{eqnarray*}
d &=&e_{1}=e_{2}=1 \\
a &=&b_{2}=0 \\
b_{1} &=&1 \\
a+b_{1}+b_{2}+l_{1}+l_{2}+d+e_{1}+e_{2}+u_{1}+u_{2} &=&l_{1}+u_{1}+u_{2}
\end{eqnarray*}%
and we get

\begin{eqnarray*}
&&\sum_{l_{1}=0}^{1}\sum_{u_{1}=0}^{1}\sum_{u_{2}=0}^{1}\left( -1\right)
^{\alpha \left( 1_{H};l_{1},0,u_{1},u_{2}\right) }B(x_{1}x_{2}\otimes
x_{1}x_{2};GX_{1},gx_{1}x_{2}) \\
&&X_{1}^{1-l_{1}}\otimes gx_{1}^{1-u_{1}}x_{2}^{1-u_{2}}\otimes
g^{l_{1}+u_{1}+u_{2}}x_{1}^{l_{1}+u_{1}}x_{2}^{u_{2}} \\
&=&\left( -1\right) ^{\alpha \left( 1_{H};0,0,0,0\right)
}B(x_{1}x_{2}\otimes x_{1}x_{2};X_{1},gx_{1}x_{2})X_{1}\otimes
gx_{1}x_{2}\otimes 1_{H}+\text{ } \\
&&\left( -1\right) ^{\alpha \left( 1_{H};0,0,0,1\right) }B(x_{1}x_{2}\otimes
x_{1}x_{2};X_{1},gx_{1}x_{2})X_{1}\otimes gx_{1}\otimes gx_{2}+\text{ } \\
&&\left( -1\right) ^{\alpha \left( 1_{H};0,0,1,0\right) }Bx_{1}x_{2}\otimes
x_{1}x_{2};X_{1},gx_{1}x_{2})X_{1}\otimes gx_{2}\otimes gx_{1}+ \\
&&\left( -1\right) ^{\alpha \left( 1_{H};0,0,1,1\right) }B(x_{1}x_{2}\otimes
x_{1}x_{2};X_{1},gx_{1}x_{2})X_{1}\otimes g\otimes x_{1}x_{2}+ \\
&&\left( -1\right) ^{\alpha \left( 1_{H};1,0,0,0\right) }B(x_{1}x_{2}\otimes
x_{1}x_{2};X_{1},gx_{1}x_{2})1_{A}\otimes gx_{1}x_{2}\otimes gx_{1}+ \\
&&\left( -1\right) ^{\alpha \left( 1_{H};1,0,0,1\right) }B(x_{1}x_{2}\otimes
x_{1}x_{2};X_{1},gx_{1}x_{2})1_{A}\otimes gx_{1}\otimes x_{1}x_{2}+ \\
&&\left( -1\right) ^{\alpha \left( 1_{H};1,0,1,0\right) }B(x_{1}x_{2}\otimes
x_{1}x_{2};X_{1},gx_{1}x_{2})1_{A}\otimes
gx_{1}^{1-u_{1}}x_{2}^{1-u_{2}}\otimes x_{1}^{1+1}x_{2}^{u_{2}}+ \\
&&\left( -1\right) ^{\alpha \left( 1_{H};1,0,1,1\right) }B(x_{1}x_{2}\otimes
x_{1}x_{2};X_{1},gx_{1}x_{2})1_{A}\otimes
gx_{1}^{1-u_{1}}x_{2}^{1-u_{2}}\otimes gx_{1}^{1+1}x_{2}^{u_{2}}.
\end{eqnarray*}

\subsubsection{Case $X_{1}\otimes gx_{1}\otimes gx_{2}$}

We have to consider only the third and the fourth summand of the left side
of the equality. Third summand gives us

\begin{eqnarray*}
l_{1} &=&u_{1}=0,l_{2}=u_{2}=0 \\
a &=&b_{2}=0,b_{1}=1 \\
d &=&e_{1}=1,e_{2}=0.
\end{eqnarray*}%
Since $\alpha \left( x_{2};0,0,0,0\right) \equiv a+b_{1}+b_{2}\equiv 1,$ we
get
\begin{equation*}
B(gx_{1}\otimes x_{1}x_{2};X_{1},gx_{1})X_{1}\otimes gx_{1}\otimes gx_{2}
\end{equation*}%
Fourth summand gives us%
\begin{eqnarray*}
l_{1} &=&u_{1}=0,l_{2}+u_{2}=1 \\
a &=&0,b_{1}=1,b_{2}=l_{2} \\
d &=&e_{1}=1,e_{2}=u_{2}.
\end{eqnarray*}%
Since $\alpha \left( 1_{H};0,0,0,1\right) \equiv a+b_{1}+b_{2}=1$ and $%
\alpha \left( 1_{H};0,1,0,0\right) \equiv 0,$ we get
\begin{equation*}
\left[ -B(x_{1}x_{2}\otimes
x_{1}x_{2};X_{1},gx_{1}x_{2})+B(x_{1}x_{2}\otimes
x_{1}x_{2};X_{1}X_{2},gx_{1})\right] X_{1}\otimes gx_{1}\otimes gx_{2}.
\end{equation*}%
By considering also the right side we get%
\begin{gather*}
-B(x_{1}x_{2}\otimes x_{1};X_{1},gx_{1})-B(x_{1}x_{2}\otimes
x_{1}x_{2};X_{1},gx_{1}x_{2}) \\
+B(x_{1}x_{2}\otimes x_{1}x_{2};X_{1}X_{2},gx_{1})+B(gx_{1}\otimes
x_{1}x_{2};X_{1},gx_{1})=0
\end{gather*}%
which holds in view of the form of the elements.

\subsubsection{Case $X_{1}\otimes gx_{2}\otimes gx_{1}$}

We have to consider only the second and the fourth summand of the left side
of the equality. Second summand gives us

\begin{eqnarray*}
l_{1} &=&u_{1}=0,l_{2}=u_{2}=0 \\
a &=&b_{2}=0,b_{1}=1 \\
d &=&e_{2}=1,e_{1}=0.
\end{eqnarray*}%
Since $\alpha \left( x_{1};0,0,0,0\right) \equiv a+b_{1}+b_{2}\equiv 1$ we
get
\begin{equation*}
-B(gx_{2}\otimes x_{1}x_{2};X_{1},gx_{2})X_{1}\otimes gx_{2}\otimes gx_{1}
\end{equation*}

Fourth summand gives us%
\begin{eqnarray*}
l_{1}+u_{1} &=&1,l_{2}=u_{2}=0 \\
a &=&b_{2}=0,b_{1}-l_{1}=1\Rightarrow b_{1}=1,l_{1}=0,u_{1}=1 \\
d &=&e_{2}=1,e_{1}=u_{1}=1.
\end{eqnarray*}%
Since $\alpha \left( 1_{H};0,0,1,0\right) \equiv e_{2}+\left(
a+b_{1}+b_{2}\right) \equiv 0,$ we get
\begin{equation*}
B(x_{1}x_{2}\otimes x_{1}x_{2};X_{1},gx_{1}x_{2})X_{1}\otimes gx_{2}\otimes
gx_{1}
\end{equation*}%
By considering also the right side of the equality, we get%
\begin{equation*}
B(x_{1}x_{2}\otimes x_{2};X_{1},gx_{2})+B(x_{1}x_{2}\otimes
x_{1}x_{2};X_{1},gx_{1}x_{2})-B(gx_{2}\otimes x_{1}x_{2};X_{1},gx_{2})=0
\end{equation*}%
which holds in view of the form of the elements.

\subsubsection{Case $X_{1}\otimes g\otimes x_{1}x_{2}$}

First summand of the left side of the equality gives us%
\begin{eqnarray*}
l_{1} &=&u_{1}=0,l_{2}=u_{2}=0, \\
a &=&b_{2}=0,b_{1}=1 \\
d &=&1,e_{1}=e_{2}=0.
\end{eqnarray*}%
Since $\alpha \left( x_{1}x_{2};0,0,0,0\right) \equiv 0$ we get
\begin{equation*}
B(1_{H}\otimes x_{1}x_{2};X_{1},g)X_{1}\otimes g\otimes x_{1}x_{2}.
\end{equation*}%
Second summand of the left side of the equality gives us%
\begin{eqnarray*}
l_{1} &=&u_{1}=0,l_{2}+u_{2}=1 \\
a &=&b_{1}=0,b_{2}=l_{2} \\
d &=&1,e_{1}=0,e_{2}=u_{2}
\end{eqnarray*}%
Since $\alpha \left( x_{1};0,0,0,1\right) \equiv 1$ and $\alpha \left(
x_{1};0,1,0,0\right) \equiv a+b_{1}+b_{2}+1\equiv 1$, we get%
\begin{equation*}
\left[ -B(gx_{2}\otimes x_{1}x_{2};X_{1},gx_{2})-B(gx_{2}\otimes
x_{1}x_{2};X_{1}X_{2},g)\right] X_{1}\otimes g\otimes x_{1}x_{2}.
\end{equation*}%
Third summand of the left side of the equality gives us

\begin{eqnarray*}
l_{1}+u_{1} &=&1,l_{2}=u_{2}=0 \\
a &=&0,b_{2}=0,b_{1}-1\Rightarrow b_{1}=1,l_{1}=0,u_{1}=1 \\
d &=&1,e_{1}=u_{1}=1,e_{2}=0
\end{eqnarray*}%
Since $\alpha \left( x_{2};0,0,1,0\right) \equiv e_{2}\equiv 0,$ we get%
\begin{equation*}
B(gx_{1}\otimes x_{1}x_{2};X_{1},gx_{1}).
\end{equation*}%
Fourth summand of the left side of the equality gives us%
\begin{eqnarray*}
X_{1}\otimes g\otimes x_{1}x_{2}l_{1}+u_{1} &=&1 \\
l_{2}+u_{2} &=&1 \\
a+b_{1}+b_{2}+d+e_{1}+e_{2} &\equiv &0 \\
d &=&1, \\
e_{1} &=&u_{1} \\
e_{2} &=&u_{2} \\
a &=&0 \\
b_{1}-1 &\Rightarrow &b_{1}=1,l_{1}=0,e_{1}=u_{1}=1 \\
b_{2} &=&l_{2} \\
\alpha \left( 1_{H};0,b_{2},1,e_{2}\right) &\equiv & \\
\alpha \left( 1_{H};0,0,1,1\right) &\equiv &1+e_{2}\equiv 0 \\
\alpha \left( 1_{H};0,1,1,0\right) &\equiv &e_{2}+\left(
a+b_{1}+b_{2}+1\right) \equiv 1
\end{eqnarray*}%
\begin{equation*}
\left[ B(x_{1}x_{2}\otimes x_{1}x_{2};X_{1},gx_{1}x_{2})-B(x_{1}x_{2}\otimes
x_{1}x_{2};X_{1}X_{2},gx_{1})+\right] X_{1}\otimes g\otimes x_{1}x_{2}
\end{equation*}%
By taking in account also the right side we get%
\begin{eqnarray*}
&&-B(x_{1}x_{2}\otimes 1_{H};X_{1},g)+B(x_{1}x_{2}\otimes
x_{1}x_{2};X_{1},gx_{1}x_{2})-B(x_{1}x_{2}\otimes
x_{1}x_{2};X_{1}X_{2},gx_{1}) \\
&&+B(gx_{1}\otimes x_{1}x_{2};X_{1},gx_{1})-B(gx_{2}\otimes
x_{1}x_{2};X_{1},gx_{2})-B(gx_{2}\otimes x_{1}x_{2};X_{1}X_{2},g) \\
B(1_{H}\otimes x_{1}x_{2};X_{1},g) &=&0
\end{eqnarray*}%
Since%
\begin{gather*}
-B(x_{1}x_{2}\otimes 1_{H};X_{1},g)=B(x_{1}x_{2}\otimes 1_{H};1_{A},gx_{1})
\\
+B(x_{1}x_{2}\otimes x_{1}x_{2};X_{1},gx_{1}x_{2})=-4B(x_{1}x_{2}\otimes
1_{H};1_{A},gx_{1}) \\
-B(x_{1}x_{2}\otimes x_{1}x_{2};X_{1}X_{2},gx_{1})=0 \\
+B(gx_{1}\otimes x_{1}x_{2};X_{1},gx_{1})=0 \\
-B(gx_{2}\otimes x_{1}x_{2};X_{1},gx_{2})=2B(x_{1}x_{2}\otimes
1_{H};1_{A},gx_{1}) \\
-B(gx_{2}\otimes x_{1}x_{2};X_{1}X_{2},g)=0 \\
B(1_{H}\otimes x_{1}x_{2};X_{1},g)=-B(x_{1}x_{2}\otimes 1_{H};1_{A},gx_{1})
\end{gather*}%
we get%
\begin{equation*}
2B(x_{1}x_{2}\otimes 1_{H};1_{A},gx_{1})=0
\end{equation*}%
already got.

\subsubsection{Case $1_{A}\otimes gx_{1}x_{2}\otimes gx_{1}$}

Nothing from the first summand of the left side. Second summand

\begin{eqnarray*}
l_{1} &=&u_{1}=0 \\
l_{2} &=&u_{2}=0 \\
a+b_{1}+b_{2}+d+e_{1}+e_{2} &\equiv &1 \\
d &=&1, \\
e_{1} &=&1 \\
e_{2} &=&1 \\
a &=&0 \\
b_{1} &=&0 \\
b_{2} &=&0 \\
\alpha \left( x_{1};0,0,0,0\right) &\equiv &a+b_{1}+b_{2}\equiv 0
\end{eqnarray*}%
\begin{equation*}
B(gx_{2}\otimes x_{1}x_{2};1_{A},gx_{1}x_{2})1_{A}\otimes gx_{1}x_{2}\otimes
gx_{1}
\end{equation*}

Nothing from the third summand. Fourth summand%
\begin{eqnarray*}
l_{1}+u_{1} &=&1 \\
l_{2} &=&u_{2}=0 \\
a+b_{1}+b_{2}+d+e_{1}+e_{2} &\equiv &0 \\
d &=&1, \\
e_{1}-u_{1} &=&1\Rightarrow e_{1}=1,u_{1}=0,l_{1}=1 \\
e_{2} &=&1 \\
a &=&0 \\
b_{1} &=&l_{1}=1 \\
b_{2} &=&0 \\
\alpha \left( 1_{H};1,0,0,0\right) &\equiv &b_{2}=0
\end{eqnarray*}%
\begin{equation*}
B(x_{1}x_{2}\otimes x_{1}x_{2};X_{1},gx_{1}x_{2})1_{A}\otimes
gx_{1}x_{2}\otimes gx_{1}
\end{equation*}%
By considering also the right side we get%
\begin{equation*}
B(x_{1}x_{2}\otimes x_{2};1_{A},gx_{1}x_{2})+B(x_{1}x_{2}\otimes
x_{1}x_{2};X_{1},gx_{1}x_{2})+B(gx_{2}\otimes x_{1}x_{2};1_{A},gx_{1}x_{2})=0
\end{equation*}%
Since%
\begin{eqnarray*}
B(x_{1}x_{2}\otimes x_{2};1_{A},gx_{1}x_{2}) &=&2B(x_{1}x_{2}\otimes
1_{H};1_{A},gx_{1}) \\
+B(x_{1}x_{2}\otimes x_{1}x_{2};X_{1},gx_{1}x_{2}) &=&-4B(x_{1}x_{2}\otimes
1_{H};1_{A},gx_{1}) \\
+B(gx_{2}\otimes x_{1}x_{2};1_{A},gx_{1}x_{2}) &=&2B(x_{1}x_{2}\otimes
1_{H};1_{A},gx_{1})
\end{eqnarray*}%
we get nothing new.

\subsubsection{Case $1_{A}\otimes gx_{1}\otimes x_{1}x_{2}$}

First summand of the left side%
\begin{eqnarray*}
l_{1} &=&u_{1}=0 \\
l_{2} &=&u_{2}=0 \\
a+b_{1}+b_{2}+d+e_{1}+e_{2} &\equiv &0 \\
d &=&1, \\
e_{1} &=&1 \\
e_{2} &=&0 \\
a &=&0 \\
b_{1} &=&0 \\
b_{2} &=&0 \\
\alpha \left( x_{1}x_{2};0,0,0,0\right) &\equiv &0
\end{eqnarray*}%
\begin{equation*}
B(1_{H}\otimes x_{1}x_{2};1_{A},gx_{1})1_{A}\otimes gx_{1}\otimes x_{1}x_{2}
\end{equation*}%
Second summand%
\begin{eqnarray*}
l_{1} &=&u_{1}=0 \\
l_{2}+u_{2} &=&1 \\
a+b_{1}+b_{2}+d+e_{1}+e_{2} &\equiv &1 \\
d &=&1, \\
e_{1} &=&1 \\
e_{2} &=&u_{2} \\
a &=&0 \\
b_{1} &=&0 \\
b_{2} &=&l_{2} \\
\alpha \left( x_{1};0,b_{2},0,e_{2}\right) &\equiv & \\
\alpha \left( x_{1};0,0,0,1\right) &\equiv &1 \\
\alpha \left( x_{1};0,1,0,0\right) &\equiv &a+b_{1}+b_{2}+1\equiv 0
\end{eqnarray*}%
\begin{equation*}
\left[ -B(gx_{2}\otimes x_{1}x_{2};1_{A},gx_{1}x_{2})+B(gx_{2}\otimes
x_{1}x_{2};x_{2},gx_{1})\right] 1_{A}\otimes gx_{1}\otimes x_{1}x_{2}
\end{equation*}%
Third summand

\begin{eqnarray*}
l_{1}+u_{1} &=&1 \\
l_{2} &=&u_{2}=0 \\
a+b_{1}+b_{2}+d+e_{1}+e_{2} &\equiv &1 \\
d &=&1, \\
e_{1}-u_{1} &=&1\Rightarrow e_{1}=1,u_{1}=0,l_{1}=1 \\
e_{2} &=&0 \\
a &=&0 \\
b_{1} &=&l_{1}=1 \\
b_{2} &=&0 \\
\alpha \left( x_{2};1,0,0,0\right) &\equiv &a+b_{1}=1
\end{eqnarray*}%
\begin{equation*}
B(gx_{1}\otimes x_{1}x_{2};X_{1},gx_{1})1_{A}\otimes gx_{1}\otimes x_{1}x_{2}
\end{equation*}%
Fourth summand%
\begin{eqnarray*}
l_{1}+u_{1} &=&1 \\
l_{2}+u_{2} &=&1 \\
a+b_{1}+b_{2}+d+e_{1}+e_{2} &\equiv &0 \\
d &=&1, \\
e_{1}-u_{1} &=&1\Rightarrow e_{1}=1,u_{1}=0,l_{1}=1 \\
e_{2} &=&u_{2} \\
a &=&0 \\
b_{1} &=&l_{1}=1 \\
b_{2} &=&l_{2} \\
\alpha \left( 1_{H};1,b_{2},0,e_{2}\right) &\equiv & \\
\alpha \left( 1_{H};1,0,0,1\right) &\equiv &a+b_{1}=1 \\
\alpha \left( 1_{H};1,1,0,0\right) &\equiv &1+b_{2}\equiv 0
\end{eqnarray*}%
\begin{equation*}
\left[ -B(x_{1}x_{2}\otimes
x_{1}x_{2};X_{1},gx_{1}x_{2})+B(x_{1}x_{2}\otimes
x_{1}x_{2};X_{1}X_{2},gx_{1})+\right] 1_{A}\otimes gx_{1}\otimes x_{1}x_{2}
\end{equation*}%
By taking in account also the right side we get%
\begin{gather*}
-B(x_{1}x_{2}\otimes 1_{H};1_{A},gx_{1})-B(x_{1}x_{2}\otimes
x_{1}x_{2};X_{1},gx_{1}x_{2})+B(x_{1}x_{2}\otimes
x_{1}x_{2};X_{1}X_{2},gx_{1}) \\
B(gx_{1}\otimes x_{1}x_{2};X_{1},gx_{1})-B(gx_{2}\otimes
x_{1}x_{2};1_{A},gx_{1}x_{2})+B(gx_{2}\otimes x_{1}x_{2};x_{2},gx_{1}) \\
+B(1_{H}\otimes x_{1}x_{2};1_{A},gx_{1})=0
\end{gather*}%
Since%
\begin{gather*}
-B(x_{1}x_{2}\otimes 1_{H};1_{A},gx_{1})=-B(x_{1}x_{2}\otimes
1_{H};1_{A},gx_{1}) \\
-B(x_{1}x_{2}\otimes x_{1}x_{2};X_{1},gx_{1}x_{2})=4B(x_{1}x_{2}\otimes
1_{H};1_{A},gx_{1}) \\
+B(x_{1}x_{2}\otimes x_{1}x_{2};X_{1}X_{2},gx_{1})=0 \\
B(gx_{1}\otimes x_{1}x_{2};X_{1},gx_{1})=0 \\
-B(gx_{2}\otimes x_{1}x_{2};1_{A},gx_{1}x_{2})=-2B(x_{1}x_{2}\otimes
1_{H};1_{A},gx_{1}) \\
+B(gx_{2}\otimes x_{1}x_{2};x_{2},gx_{1})=0 \\
+B(1_{H}\otimes x_{1}x_{2};1_{A},gx_{1})=B(x_{1}x_{2}\otimes
1_{H};1_{A},gx_{1})+2B(gx_{2}\otimes 1_{H};1_{A},g)
\end{gather*}%
Since we know that $B(gx_{2}\otimes 1_{H};1_{A},g)=0$ we get
\begin{equation*}
2B(x_{1}x_{2}\otimes 1_{H};1_{A},gx_{1})=0
\end{equation*}%
already got.

\subsection{$B(x_{1}x_{2}\otimes x_{1}x_{2};X_{2},gx_{1}x_{2})$}

\begin{equation*}
a=0,b_{1}=0,b_{2}=1,d=1,e_{1}=1,e_{2}=1
\end{equation*}

and we get%
\begin{multline*}
\sum_{l_{2}=0}^{1}\sum_{u_{1}=0}^{1}\sum_{u_{2}=0}^{1}\left( -1\right)
^{\alpha \left( 1_{H};0,l_{2},u_{1},u_{2}\right) }B(x_{1}x_{2}\otimes
x_{1}x_{2};X_{2},gx_{1}x_{2}) \\
X_{2}^{1-l_{2}}\otimes gx_{1}^{1-u_{1}}x_{2}^{1-u_{2}}\otimes
g^{l_{2}+u_{1}+u_{2}}x_{1}^{u_{1}}x_{2}^{l_{2}+u_{2}}= \\
+\left( -1\right) ^{\alpha \left( 1_{H};0,0,0,0\right) }B(x_{1}x_{2}\otimes
x_{1}x_{2};X_{2},gx_{1}x_{2})X_{2}\otimes gx_{1}x_{2}\otimes 1_{H}+ \\
+\left( -1\right) ^{\alpha \left( 1_{H};0,0,0,1\right) }B(x_{1}x_{2}\otimes
x_{1}x_{2};X_{2},gx_{1}x_{2})X_{2}\otimes gx_{1}\otimes gx_{2}+ \\
+\left( -1\right) ^{\alpha \left( 1_{H};0,0,1,0\right) }B(x_{1}x_{2}\otimes
x_{1}x_{2};X_{2},gx_{1}x_{2})X_{2}\otimes gx_{2}\otimes gx_{1}+ \\
+\left( -1\right) ^{\alpha \left( 1_{H};0,0,1,1\right) }B(x_{1}x_{2}\otimes
x_{1}x_{2};X_{2},gx_{1}x_{2})X_{2}\otimes g\otimes x_{1}x_{2}+ \\
+\left( -1\right) ^{\alpha \left( 1_{H};0,1,0,0\right) }B(x_{1}x_{2}\otimes
x_{1}x_{2};X_{2},gx_{1}x_{2})1_{A}\otimes gx_{1}x_{2}\otimes gx_{2}+ \\
+\left( -1\right) ^{\alpha \left( 1_{H};0,1,0,1\right) }B(x_{1}x_{2}\otimes
x_{1}x_{2};X_{2},gx_{1}x_{2})1_{A}\otimes
gx_{1}^{1-u_{1}}x_{2}^{1-u_{2}}\otimes
g^{l_{2}+u_{1}+u_{2}}x_{1}^{u_{1}}x_{2}^{1+1}+=0 \\
+\left( -1\right) ^{\alpha \left( 1_{H};0,1,1,0\right) }B(x_{1}x_{2}\otimes
x_{1}x_{2};X_{2},gx_{1}x_{2})1_{A}\otimes gx_{2}\otimes x_{1}x_{2}+ \\
+\left( -1\right) ^{\alpha \left( 1_{H};0,11,1\right) }B(x_{1}x_{2}\otimes
x_{1}x_{2};X_{2},gx_{1}x_{2})1_{A}\otimes
gx_{1}^{1-u_{1}}x_{2}^{1-u_{2}}\otimes
g^{l_{2}+u_{1}+u_{2}}x_{1}^{u_{1}}x_{2}^{1+1}=0
\end{multline*}

\subsubsection{Case $X_{2}\otimes gx_{1}\otimes gx_{2}$}

We have to consider only the third and the fourth summand of the left side
of the equality. Third summand gives us

\begin{eqnarray*}
l_{1} &=&u_{1}=0,l_{2}=u_{2}=0 \\
a &=&b_{1}=0,b_{2}=1 \\
d &=&e_{1}=1,e_{2}=0.
\end{eqnarray*}%
Since $\alpha \left( x_{2};0,0,0,0\right) \equiv a+b_{1}+b_{2}\equiv 1,$ we
get
\begin{equation*}
B(gx_{1}\otimes x_{1}x_{2};X_{2},gx_{1})X_{2}\otimes gx_{1}\otimes gx_{2}
\end{equation*}%
Fourth summand gives us%
\begin{eqnarray*}
l_{1} &=&u_{1}=0,l_{2}+u_{2}=1 \\
a &=&b_{1}=0,b_{2}-l_{2}=1\Rightarrow b_{2}=1,l_{2}=0,u_{2}=1 \\
d &=&e_{1}=1,e_{2}=u_{2}=1.
\end{eqnarray*}%
Since $\alpha \left( 1_{H};0,0,0,1\right) \equiv a+b_{1}+b_{2}=1,$ we get
\begin{equation*}
-B(x_{1}x_{2}\otimes x_{1}x_{2};X_{2},gx_{1}x_{2})X_{2}\otimes gx_{1}\otimes
gx_{2}.
\end{equation*}%
By considering also the right side we get%
\begin{equation*}
-B(x_{1}x_{2}\otimes x_{1};X_{2},gx_{1})-B(x_{1}x_{2}\otimes
x_{1}x_{2};X_{2},gx_{1}x_{2})+B(gx_{1}\otimes x_{1}x_{2};X_{2},gx_{1})=0
\end{equation*}%
which holds in view of the form of the elements.

\subsubsection{Case $X_{2}\otimes gx_{2}\otimes gx_{1}$}

We have to consider only the second and the fourth summand of the left side
of the equality. Second summand gives us

\begin{eqnarray*}
l_{1} &=&u_{1}=0,l_{2}=u_{2}=0 \\
a &=&b_{1}=0,b_{2}=1 \\
d &=&e_{2}=1,e_{1}=0.
\end{eqnarray*}%
Since $\alpha \left( x_{1};0,0,0,0\right) \equiv a+b_{1}+b_{2}\equiv 1$ we
get
\begin{equation*}
-B(gx_{2}\otimes x_{1}x_{2};X_{2},gx_{2})X_{2}\otimes gx_{2}\otimes gx_{1}
\end{equation*}

Fourth summand gives us%
\begin{eqnarray*}
l_{1}+u_{1} &=&1,l_{2}=u_{2}=0 \\
a &=&0,b_{1}\newline
=l_{1},b_{2}=1 \\
d &=&e_{2}=1,e_{1}=u_{1}.
\end{eqnarray*}%
Since $\alpha \left( 1_{H};0,0,1,0\right) \equiv e_{2}+\left(
a+b_{1}+b_{2}\right) \equiv 0$ and $\alpha \left( 1_{H};1,0,0,0\right)
\equiv b_{2}=1,$we get
\begin{equation*}
\left[ B(x_{1}x_{2}\otimes x_{1}x_{2};X_{2},gx_{1}x_{2})-B(x_{1}x_{2}\otimes
x_{1}x_{2};X_{1}X_{2},gx_{2})\right] X_{2}\otimes gx_{2}\otimes gx_{1}
\end{equation*}%
By considering also the right side of the equality, we get%
\begin{gather*}
B(x_{1}x_{2}\otimes x_{2};X_{2},gx_{2})+B(x_{1}x_{2}\otimes
x_{1}x_{2};X_{2},gx_{1}x_{2}) \\
-B(x_{1}x_{2}\otimes x_{1}x_{2};X_{1}X_{2},gx_{2})-B(gx_{2}\otimes
x_{1}x_{2};X_{2},gx_{2})=0
\end{gather*}%
which holds in view of the form of the elements.

\subsubsection{Case $X_{2}\otimes g\otimes x_{1}x_{2}$}

First summand of the left side of the equality gives%
\begin{eqnarray*}
l_{1} &=&u_{1}=0,l_{2}=u_{2}=0 \\
a &=&b_{1}=0,b_{2}=1, \\
d &=&1,e_{1}=e_{2}=0.
\end{eqnarray*}%
Since $\alpha \left( x_{1}x_{2};0,0,0,0\right) \equiv 0,$ we get%
\begin{equation*}
B(1_{H}\otimes x_{1}x_{2};X_{2},g)X_{2}\otimes g\otimes x_{1}x_{2}.
\end{equation*}%
Second summand of the left side of the equality gives us%
\begin{eqnarray*}
l_{1} &=&u_{1}=0,l_{2}+u_{2}=1 \\
a &=&b_{1}=0,b_{2}-l_{2}=1\Rightarrow b_{2}=1,l_{2}=0,u_{2}=1 \\
d &=&1,e_{1}=0,e_{2}=u_{2}=1.
\end{eqnarray*}%
Since $\alpha \left( x_{1};0,0,0,1\right) \equiv 1$, we get%
\begin{equation*}
-B(gx_{2}\otimes x_{1}x_{2};X_{2},gx_{2})X_{2}\otimes g\otimes x_{1}x_{2}.
\end{equation*}%
Third summand of the left side of the equality gives us

\begin{eqnarray*}
l_{1}+u_{1} &=&1,l_{2}=u_{2}=0 \\
a &=&0,b_{1}=l_{1},b_{2}=1 \\
d &=&1,e_{1}=u_{1},e_{2}=0,
\end{eqnarray*}%
Since $\alpha \left( x_{2};0,0,1,0\right) \equiv e_{2}=0$ and $\alpha \left(
x_{2};1,0,0,0\right) \equiv a+b_{1}\equiv 1,$ we get%
\begin{equation*}
\left[ -B(gx_{1}\otimes x_{1}x_{2};X_{2},gx_{1})+B(gx_{1}\otimes
x_{1}x_{2};X_{1}X_{2},g)\right] X_{2}\otimes g\otimes x_{1}x_{2}.
\end{equation*}%
Fourth summand gives us%
\begin{eqnarray*}
l_{1}+u_{1} &=&1,l_{2}+u_{2}=1, \\
a &=&0,b_{1}=l_{1},b_{2}-l_{2}=1\Rightarrow b_{2}=1,l_{2}=0,u_{2}=1, \\
d &=&1,e_{1}=u_{1},e_{2}=u_{2}=1.
\end{eqnarray*}%
Since $\alpha \left( 1_{H};0,0,1,1\right) \equiv 1+e_{2}\equiv 0$ and $%
\alpha \left( 1_{H};1,0,0,1\right) \equiv a+b_{1}\equiv 1,$ we obtain
\begin{equation*}
\left[ +B(x_{1}x_{2}\otimes
x_{1}x_{2};X_{2},gx_{1}x_{2})-B(x_{1}x_{2}\otimes
x_{1}x_{2};X_{1}X_{2},gx_{2})\right] X_{2}\otimes g\otimes x_{1}x_{2}.
\end{equation*}%
By taking in account also the right side of the equality we get%
\begin{gather*}
-B(x_{1}x_{2}\otimes 1_{H};X_{2},g)+B(x_{1}x_{2}\otimes
x_{1}x_{2};X_{2},gx_{1}x_{2})-B(x_{1}x_{2}\otimes
x_{1}x_{2};X_{1}X_{2},gx_{2}) \\
-B(gx_{1}\otimes x_{1}x_{2};X_{2},gx_{1})+B(gx_{1}\otimes
x_{1}x_{2};X_{1}X_{2},g)-B(gx_{2}\otimes x_{1}x_{2};X_{2},gx_{2})+ \\
+B(1_{H}\otimes x_{1}x_{2};X_{2},g)=0
\end{gather*}

which holds in view of the form of the elements.

\subsubsection{Case $1_{A}\otimes gx_{1}x_{2}\otimes gx_{2}$}

We have to consider only the third and the fourth summand of the left side
of the equality. Third summand gives us

\begin{eqnarray*}
l_{1} &=&u_{1}=0,l_{2}=u_{2}=0 \\
a &=&b_{1}=b_{2}=0, \\
d &=&e_{1}=e_{2}=1.
\end{eqnarray*}%
Since $\alpha \left( x_{2};0,0,0,0\right) \equiv a+b_{1}+b_{2}\equiv 0,$ we
get
\begin{equation*}
-B(gx_{1}\otimes x_{1}x_{2};1_{A},gx_{1}x_{2})1_{A}\otimes
gx_{1}x_{2}\otimes gx_{2}.
\end{equation*}%
Fourth summand gives us%
\begin{eqnarray*}
l_{1} &=&u_{1}=0,l_{2}+u_{2}=1 \\
a &=&b_{1}=0,b_{2}=l_{2} \\
d &=&e_{1}=1,e_{2}-u_{2}=1\Rightarrow e_{2}=1,u_{2}=0,b_{2}=l_{2}=1.
\end{eqnarray*}%
Since $\alpha \left( 1_{H};0,1,0,0\right) \equiv 0,$ we get
\begin{equation*}
B(x_{1}x_{2}\otimes x_{1}x_{2};X_{2},gx_{1}x_{2})1_{A}\otimes
gx_{1}x_{2}\otimes gx_{2}.
\end{equation*}%
By considering also the right side we get%
\begin{equation*}
-B(x_{1}x_{2}\otimes x_{1};1_{A},gx_{1}x_{2})+B(x_{1}x_{2}\otimes
x_{1}x_{2};X_{2},gx_{1}x_{2})-B(gx_{1}\otimes x_{1}x_{2};1_{A},gx_{1}x_{2})=0
\end{equation*}%
which holds in view of the form of the elements.

\subsubsection{Case $1_{A}\otimes gx_{2}\otimes x_{1}x_{2}$}

First summand of the left side of the equality gives%
\begin{eqnarray*}
l_{1} &=&u_{1}=0,l_{2}=u_{2}=0 \\
a &=&b_{1}=b_{2}=0, \\
d &=&e_{2}=1,e_{1}=0.
\end{eqnarray*}%
Since $\alpha \left( x_{1}x_{2};0,0,0,0\right) \equiv 0,$ we get%
\begin{equation*}
B(1_{H}\otimes x_{1}x_{2};1_{A},gx_{2})1_{A}\otimes gx_{2}\otimes x_{1}x_{2}.
\end{equation*}%
Second summand of the left side of the equality gives us%
\begin{eqnarray*}
l_{1} &=&u_{1}=0,l_{2}+u_{2}=1 \\
a &=&b_{1}=0,b_{2}=l_{2}, \\
d &=&1,e_{1}=0,e_{2}-u_{2}=1\Rightarrow e_{2}=1,u_{2}=0,b_{2}=l_{2}=1.
\end{eqnarray*}%
Since $\alpha \left( x_{1};0,1,0,0\right) \equiv a+b_{1}+b_{2}+1\equiv 0$,
we get%
\begin{equation*}
B(gx_{2}\otimes x_{1}x_{2};X_{2},gx_{2})1_{A}\otimes gx_{2}\otimes
x_{1}x_{2}.
\end{equation*}%
Third summand of the left side of the equality gives us

\begin{eqnarray*}
l_{1}+u_{1} &=&1,l_{2}=u_{2}=0 \\
a &=&b_{2}=0,b_{1}=l_{1}, \\
d &=&1,e_{1}=u_{1},e_{2}=1,
\end{eqnarray*}%
Since $\alpha \left( x_{2};0,0,1,0\right) \equiv e_{2}=1$ and $\alpha \left(
x_{2};1,0,0,0\right) \equiv a+b_{1}\equiv 1,$ we get%
\begin{equation*}
\left[ +B(gx_{1}\otimes x_{1}x_{2};1_{A},gx_{1}x_{2})+B(gx_{1}\otimes
x_{1}x_{2};X_{1},gx_{2})\right] 1_{A}\otimes gx_{2}\otimes x_{1}x_{2}.
\end{equation*}%
Fourth summand gives us%
\begin{eqnarray*}
l_{1}+u_{1} &=&1,l_{2}+u_{2}=1, \\
a &=&0,b_{1}=l_{1},b_{2}=l_{2}, \\
d &=&1,e_{1}=u_{1},e_{2}-u_{2}=1\Rightarrow e_{2}=1,u_{2}=0,b_{2}=l_{2}=1.
\end{eqnarray*}%
Since $\alpha \left( 1_{H};0,1,1,0\right) \equiv e_{2}+a+b_{1}+b_{2}+1\equiv
1$ and $\alpha \left( 1_{H};1,1,0,0\right) \equiv 1+b_{2}\equiv 0,$ we
obtain
\begin{equation*}
\left[ -B(x_{1}x_{2}\otimes
x_{1}x_{2};X_{2},gx_{1}x_{2})+B(x_{1}x_{2}\otimes
x_{1}x_{2};X_{1}X_{2},gx_{2})\right] 1_{A}\otimes gx_{2}\otimes x_{1}x_{2}.
\end{equation*}%
By taking in account also the right side of the equality we get%
\begin{gather*}
-B(x_{1}x_{2}\otimes 1_{H};1_{A},gx_{2})-B(x_{1}x_{2}\otimes
x_{1}x_{2};X_{2},gx_{1}x_{2})+B(x_{1}x_{2}\otimes
x_{1}x_{2};X_{1}X_{2},gx_{2}) \\
+B(gx_{1}\otimes x_{1}x_{2};1_{A},gx_{1}x_{2})+B(gx_{1}\otimes
x_{1}x_{2};X_{1},gx_{2})+B(gx_{2}\otimes x_{1}x_{2};X_{2},gx_{2})+ \\
+B(1_{H}\otimes x_{1}x_{2};1_{A},gx_{2})=0
\end{gather*}

which holds in view of the form of the elements.

\subsection{$B(x_{1}x_{2}\otimes x_{1}x_{2};GX_{1}X_{2},gx_{1}x_{2})$}

We deduce that%
\begin{equation*}
a=b_{1}=b_{2}=d=e_{1}=e_{2}=1
\end{equation*}%
and we obtain%
\begin{eqnarray*}
&&\sum_{l_{1}=0}^{1}\sum_{l_{2}=0}^{1}\sum_{u_{1}=0}^{1}\sum_{u_{2}=0}^{1}%
\left( -1\right) ^{\alpha \left( 1_{H};l_{1},l_{2},u_{1},u_{2}\right)
}B(x_{1}\otimes x_{1}x_{2};GX_{1}X_{2},gx_{1}x_{2}) \\
&&G^{a}X_{1}^{1-l_{1}}X_{2}^{1-l_{2}}\otimes
gx_{1}^{1-u_{1}}x_{2}^{1-u_{2}}\otimes
g^{l_{1}+l_{2}+u_{1}+u_{2}}x_{1}^{l_{1}+u_{1}}x_{2}^{l_{2}+u_{2}}
\end{eqnarray*}%
\begin{gather*}
\left( -1\right) ^{\alpha \left( 1_{H};0,0,0,0\right) }B(x_{1}x_{2}\otimes
x_{1}x_{2};GX_{1}X_{2},gx_{1}x_{2})GX_{1}X_{2}\otimes gx_{1}x_{2}\otimes
1_{H}+ \\
+\left( -1\right) ^{\alpha \left( 1_{H};0,0,0,1\right) }B(x_{1}x_{2}\otimes
x_{1}x_{2};GX_{1}X_{2},gx_{1}x_{2})GX_{1}X_{2}\otimes gx_{1}\otimes gx_{2}%
\text{ } \\
+\left( -1\right) ^{\alpha \left( 1_{H};0,0,1,0\right) }B(x_{1}x_{2}\otimes
x_{1}x_{2};GX_{1}X_{2},gx_{1}x_{2})GX_{1}X_{2}\otimes gx_{2}\otimes gx_{1}+
\\
+\left( -1\right) ^{\alpha \left( 1_{H};0,0,1,1\right) }B(x_{1}x_{2}\otimes
x_{1}x_{2};GX_{1}X_{2},gx_{1}x_{2})GX_{1}X_{2}\otimes g\otimes x_{1}x_{2}+ \\
+\left( -1\right) ^{\alpha \left( 1_{H};0,1,0,0\right) }B(x_{1}x_{2}\otimes
x_{1}x_{2};GX_{1}X_{2},gx_{1}x_{2})GX_{1}\otimes gx_{1}x_{2}\otimes gx_{2}+
\\
+\left( -1\right) ^{\alpha \left( 1_{H};0,1,0,1\right) }B(x_{1}x_{2}\otimes
x_{1}x_{2};GX_{1}X_{2},gx_{1}x_{2})GX_{1}\otimes
gx_{1}x_{2}^{1-u_{2}}\otimes g^{u_{2}}x_{2}^{1+1}+=0 \\
+\left( -1\right) ^{\alpha \left( 1_{H};0,1,1,0\right) }B(x_{1}x_{2}\otimes
x_{1}x_{2};GX_{1}X_{2},gx_{1}x_{2})GX_{1}\otimes gx_{2}\otimes x_{1}x_{2}+ \\
+\left( -1\right) ^{\alpha \left( 1_{H};0,1,1,1\right) }B(x_{1}x_{2}\otimes
x_{1}x_{2};GX_{1}X_{2},gx_{1}x_{2})GX_{1}\otimes gx_{2}^{1-u_{2}}\otimes
g^{1+u_{2}}x_{1}x_{2}^{1+1}+=0 \\
+\left( -1\right) ^{\alpha \left( 1_{H};1,0,0,0\right) }B(x_{1}x_{2}\otimes
x_{1}x_{2};GX_{1}X_{2},gx_{1}x_{2})GX_{2}\otimes gx_{1}x_{2}\otimes gx_{1}+
\\
+\left( -1\right) ^{\alpha \left( 1_{H};1,0,0,1\right) }B(x_{1}x_{2}\otimes
x_{1}x_{2};GX_{1}X_{2},gx_{1}x_{2})GX_{2}\otimes
gx_{1}x_{2}^{1-u_{2}}\otimes g^{1+u_{2}}x_{1}^{1+0}x_{2}^{1+1}+=0 \\
+\left( -1\right) ^{\alpha \left( 1_{H};1,0,1,u_{2}\right)
}B(x_{1}x_{2}\otimes x_{1}x_{2};GX_{1}X_{2},gx_{1}x_{2})GX_{2}\otimes
gx_{1}^{1-u_{1}}x_{2}^{1-u_{2}}\otimes
g^{1+u_{1}+u_{2}}x_{1}^{1+1}x_{2}^{1+u_{2}}+=0 \\
+\left( -1\right) ^{\alpha \left( 1_{H};1,1,0,0\right) }B(x_{1}x_{2}\otimes
x_{1}x_{2};GX_{1}X_{2},gx_{1}x_{2})G\otimes gx_{1}x_{2}\otimes x_{1}x_{2}+ \\
+\left( -1\right) ^{\alpha \left( 1_{H};1,1,0,1\right) }B(x_{1}x_{2}\otimes
x_{1}x_{2};GX_{1}X_{2},gx_{1}x_{2})G\otimes gx_{1}x_{2}^{1-u_{2}}\otimes
g^{u_{2}}x_{1}x_{2}^{1+1}+=0 \\
+\left( -1\right) ^{\alpha \left( 1_{H};1,1,1,u_{2}\right)
}B(x_{1}x_{2}\otimes x_{1}x_{2};GX_{1}X_{2},gx_{1}x_{2})G\otimes
gx_{1}^{1-u_{1}}x_{2}^{1-u_{2}}\otimes
g^{u_{1}+u_{2}}x_{1}^{1+1}x_{2}^{1+u_{2}}=0.
\end{gather*}

\subsubsection{Case $GX_{1}X_{2}\otimes gx_{1}\otimes gx_{2}$}

We have to consider only the third and the fourth summand of the left side
of the equality. Third summand gives us

\begin{eqnarray*}
l_{1} &=&u_{1}=0,l_{2}=u_{2}=0 \\
a &=&b_{1}=b_{2}=1, \\
d &=&e_{1}=1,e_{2}=0.
\end{eqnarray*}%
Since $\alpha \left( x_{2};0,0,0,0\right) \equiv a+b_{1}+b_{2}\equiv 1,$ we
get
\begin{equation*}
B(gx_{1}\otimes x_{1}x_{2};GX_{1}X_{2},gx_{1})GX_{1}X_{2}\otimes
gx_{1}\otimes gx_{2}.
\end{equation*}%
Fourth summand gives us%
\begin{eqnarray*}
l_{1} &=&u_{1}=0,l_{2}+u_{2}=1 \\
a &=&b_{1}=1,b_{2}-l_{2}=1\Rightarrow b_{2}=1,l_{2}=0,u_{2}=1 \\
d &=&e_{1}=1,e_{2}=u_{2}=1.
\end{eqnarray*}%
Since $\alpha \left( 1_{H};0,0,0,1\right) \equiv a+b_{1}+b_{2}\equiv 1,$ we
get
\begin{equation*}
-B(x_{1}x_{2}\otimes x_{1}x_{2};GX_{1}X_{2},gx_{1}x_{2})GX_{1}X_{2}\otimes
gx_{1}\otimes gx_{2}.
\end{equation*}%
By considering also the right side we get%
\begin{equation*}
-B(x_{1}x_{2}\otimes x_{1};GX_{1}X_{2},gx_{1})-B(x_{1}x_{2}\otimes
x_{1}x_{2};GX_{1}X_{2},gx_{1}x_{2})+B(gx_{1}\otimes
x_{1}x_{2};GX_{1}X_{2},gx_{1})=0
\end{equation*}%
which holds in view of the form of the elements.

\subsubsection{Case $GX_{1}X_{2}\otimes gx_{2}\otimes gx_{1}$}

We have to consider only the second and the fourth summand of the left side
of the equality. Second summand gives us

\begin{eqnarray*}
l_{1} &=&u_{1}=0,l_{2}=u_{2}=0 \\
a &=&b_{1}=b_{2}=1 \\
d &=&e_{2}=1,e_{1}=0.
\end{eqnarray*}%
Since $\alpha \left( x_{1};0,0,0,0\right) \equiv a+b_{1}+b_{2}\equiv 1$ we
get
\begin{equation*}
-B(gx_{2}\otimes x_{1}x_{2};GX_{1}X_{2},gx_{2})GX_{1}X_{2}\otimes
gx_{2}\otimes gx_{1}.
\end{equation*}

Fourth summand gives us%
\begin{eqnarray*}
l_{1}+u_{1} &=&1,l_{2}=u_{2}=0 \\
a &=&b_{2}=1,b_{1}-l_{1}=1\Rightarrow b_{1}=1,l_{1}=0,u_{1}=1, \\
d &=&e_{2}=1,e_{1}=u_{1}=1.
\end{eqnarray*}%
Since $\alpha \left( 1_{H};0,0,1,0\right) \equiv e_{2}+\left(
a+b_{1}+b_{2}\right) \equiv 0,$we get
\begin{equation*}
B(x_{1}x_{2}\otimes x_{1}x_{2};GX_{1}X_{2},gx_{1}x_{2})GX_{1}X_{2}\otimes
gx_{2}\otimes gx_{1}.
\end{equation*}%
By considering also the right side of the equality, we get%
\begin{equation*}
B(x_{1}x_{2}\otimes x_{2};GX_{1}X_{2},gx_{2})+B(x_{1}x_{2}\otimes
x_{1}x_{2};GX_{1}X_{2},gx_{1}x_{2})-B(gx_{2}\otimes
x_{1}x_{2};GX_{1}X_{2},gx_{2})=0
\end{equation*}%
which holds in view of the form of the elements.

\subsubsection{Case $GX_{1}X_{2}\otimes g\otimes x_{1}x_{2}$}

First summand of the left side of the equality gives%
\begin{eqnarray*}
l_{1} &=&u_{1}=0,l_{2}=u_{2}=0 \\
a &=&b_{1}=b_{2}=1, \\
d &=&1,e_{1}=e_{2}=0.
\end{eqnarray*}%
Since $\alpha \left( x_{1}x_{2};0,0,0,0\right) \equiv 0,$ we get%
\begin{equation*}
B(1_{H}\otimes x_{1}x_{2};GX_{1}X_{2},g)GX_{1}X_{2}\otimes g\otimes
x_{1}x_{2}.
\end{equation*}%
Second summand of the left side of the equality gives us%
\begin{eqnarray*}
l_{1} &=&u_{1}=0,l_{2}+u_{2}=1 \\
a &=&b_{1}=1,b_{2}-l_{2}=1\Rightarrow b_{2}=1,l_{2}=0,u_{2}=1 \\
d &=&1,e_{1}=0,e_{2}=u_{2}=1.
\end{eqnarray*}%
Since $\alpha \left( x_{1};0,0,0,1\right) \equiv 1$, we get%
\begin{equation*}
-B(gx_{2}\otimes x_{1}x_{2};GX_{1}X_{2},gx_{2})GX_{1}X_{2}\otimes g\otimes
x_{1}x_{2}.
\end{equation*}%
Third summand of the left side of the equality gives us

\begin{eqnarray*}
l_{1}+u_{1} &=&1,l_{2}=u_{2}=0 \\
a &=&b_{2}=1,b_{1}-l_{1}=1\Rightarrow b_{1}=1,l_{1}=0,u_{1}=1, \\
d &=&1,e_{1}=u_{1}=1,e_{2}=0,
\end{eqnarray*}%
Since $\alpha \left( x_{2};0,0,1,0\right) \equiv e_{2}=0,$ we get%
\begin{equation*}
-B(gx_{1}\otimes x_{1}x_{2};GX_{1}X_{2},gx_{1})GX_{1}X_{2}\otimes g\otimes
x_{1}x_{2}.
\end{equation*}%
Fourth summand gives us%
\begin{eqnarray*}
l_{1}+u_{1} &=&1,l_{2}+u_{2}=1, \\
a &=&1,b_{1}-l_{1}=1\Rightarrow b_{1}=1,l_{1}=0,u_{1}=1 \\
b_{2}-l_{2} &=&1\Rightarrow b_{2}=1,l_{2}=0,u_{2}=1, \\
d &=&1,e_{1}=u_{1}=1,e_{2}=u_{2}=1.
\end{eqnarray*}%
Since $\alpha \left( 1_{H};0,0,1,1\right) \equiv 1+e_{2}\equiv 0,$ we obtain
\begin{equation*}
+B(x_{1}x_{2}\otimes x_{1}x_{2};GX_{1}X_{2},gx_{1}x_{2})GX_{1}X_{2}\otimes
g\otimes x_{1}x_{2}.
\end{equation*}%
By taking in account also the right side of the equality we get%
\begin{gather*}
-B(x_{1}x_{2}\otimes 1_{H};GX_{1}X_{2},g)+B(x_{1}x_{2}\otimes
x_{1}x_{2};GX_{1}X_{2},gx_{1}x_{2}) \\
-B(gx_{1}\otimes x_{1}x_{2};GX_{1}X_{2},gx_{1})-B(gx_{2}\otimes
x_{1}x_{2};GX_{1}X_{2},gx_{2}) \\
+B(1_{H}\otimes x_{1}x_{2};GX_{1}X_{2},g)=0
\end{gather*}

which holds in view of the form of the elements.

\subsubsection{Case $GX_{1}\otimes gx_{1}x_{2}\otimes gx_{2}$}

We have to consider only the third and the fourth summand of the left side
of the equality. Third summand gives us

\begin{eqnarray*}
l_{1} &=&u_{1}=0,l_{2}=u_{2}=0 \\
a &=&b_{1}=1,b_{2}=0 \\
d &=&e_{1}=e_{2}=1.
\end{eqnarray*}%
Since $\alpha \left( x_{2};0,0,0,0\right) \equiv a+b_{1}+b_{2}\equiv 0,$ we
get
\begin{equation*}
-B(gx_{1}\otimes x_{1}x_{2};GX_{1},gx_{1}x_{2})GX_{1}\otimes
gx_{1}x_{2}\otimes gx_{2}.
\end{equation*}%
Fourth summand gives us%
\begin{eqnarray*}
l_{1} &=&u_{1}=0,l_{2}+u_{2}=1 \\
a &=&b_{1}=1,b_{2}=l_{2}, \\
d &=&e_{1}=1,e_{2}-u_{2}=1\Rightarrow e_{2}=1,u_{2}=0,b_{2}=l_{2}=1.
\end{eqnarray*}%
Since $\alpha \left( 1_{H};0,1,0,0\right) \equiv 0,$ we get
\begin{equation*}
+B(x_{1}x_{2}\otimes x_{1}x_{2};GX_{1}X_{2},gx_{1}x_{2})GX_{1}\otimes
gx_{1}x_{2}\otimes gx_{2}.
\end{equation*}%
By considering also the right side we get%
\begin{equation*}
-B(x_{1}x_{2}\otimes x_{1};GX_{1},gx_{1}x_{2})+B(x_{1}x_{2}\otimes
x_{1}x_{2};GX_{1}X_{2},gx_{1}x_{2})-B(gx_{1}\otimes
x_{1}x_{2};GX_{1},gx_{1}x_{2})=0
\end{equation*}%
which holds in view of the form of the elements.

\subsubsection{Case $GX_{1}\otimes gx_{2}\otimes x_{1}x_{2}$}

First summand of the left side of the equality gives%
\begin{eqnarray*}
l_{1} &=&u_{1}=0,l_{2}=u_{2}=0 \\
a &=&b_{1}=1,b_{2}=0 \\
d &=&e_{2}=1,e_{1}=0.
\end{eqnarray*}%
Since $\alpha \left( x_{1}x_{2};0,0,0,0\right) \equiv 0,$ we get%
\begin{equation*}
B(1_{H}\otimes x_{1}x_{2};GX_{1}X_{2},g)GX_{1}\otimes gx_{2}\otimes
x_{1}x_{2}.
\end{equation*}%
Second summand of the left side of the equality gives us%
\begin{eqnarray*}
l_{1} &=&u_{1}=0,l_{2}+u_{2}=1 \\
a &=&b_{1}=1,b_{2}=l_{2} \\
d &=&1,e_{1}=0,e_{2}-u_{2}=1\Rightarrow e_{2}=1,u_{2}=0,b_{2}=l_{2}=1.
\end{eqnarray*}%
Since $\alpha \left( x_{1};0,1,0,0\right) \equiv a+b_{1}+b_{2}+1\equiv 1$,
we get%
\begin{equation*}
-B(gx_{2}\otimes x_{1}x_{2};GX_{1}X_{2},gx_{2})GX_{1}\otimes gx_{2}\otimes
x_{1}x_{2}.
\end{equation*}%
Third summand of the left side of the equality gives us

\begin{eqnarray*}
l_{1}+u_{1} &=&1,l_{2}=u_{2}=0 \\
a &=&1,b_{2}=0,b_{1}-l_{1}=1\Rightarrow b_{1}=1,l_{1}=0,u_{1}=1, \\
d &=&e_{2}=1,e_{1}=u_{1}=1,
\end{eqnarray*}%
Since $\alpha \left( x_{2};0,0,1,0\right) \equiv e_{2}=1,$ we get%
\begin{equation*}
B(gx_{1}\otimes x_{1}x_{2};GX_{1},gx_{1}x_{2})GX_{1}\otimes gx_{2}\otimes
x_{1}x_{2}.
\end{equation*}%
Fourth summand gives us%
\begin{eqnarray*}
l_{1}+u_{1} &=&1,l_{2}+u_{2}=1, \\
a &=&1,b_{1}-l_{1}=1\Rightarrow b_{1}=1,l_{1}=0,u_{1}=1,b_{2}=l_{2}, \\
d &=&1,e_{1}=u_{1}=1,e_{2}-u_{2}=1\Rightarrow e_{2}=1,u_{2}=0,b_{2}=l_{2}=1.
\end{eqnarray*}%
Since $\alpha \left( 1_{H};0,1,1,0\right) =e_{2}+a+b_{1}+b_{2}+1\equiv 1,$
we obtain
\begin{equation*}
-B(x_{1}x_{2}\otimes x_{1}x_{2};GX_{1}X_{2},gx_{1}x_{2})GX_{1}\otimes
gx_{2}\otimes x_{1}x_{2}.
\end{equation*}%
By taking in account also the right side of the equality we get%
\begin{gather*}
-B(x_{1}x_{2}\otimes 1_{H};GX_{1},gx_{2})-B(x_{1}x_{2}\otimes
x_{1}x_{2};GX_{1}X_{2},gx_{1}x_{2}) \\
+B(gx_{1}\otimes x_{1}x_{2};GX_{1},gx_{1}x_{2})-B(gx_{2}\otimes
x_{1}x_{2};GX_{1}X_{2},gx_{2}) \\
+B(1_{H}\otimes x_{1}x_{2};GX_{1}X_{2},g)=0
\end{gather*}

which we already got above in Case $GX_{1}X_{2}\otimes g\otimes x_{1}x_{2}.$

\subsubsection{Case $GX_{2}\otimes gx_{1}x_{2}\otimes gx_{1}$}

We have to consider only the second and the fourth summand of the left side
of the equality. Second summand gives us

\begin{eqnarray*}
l_{1} &=&u_{1}=0,l_{2}=u_{2}=0 \\
a &=&b_{2}=1,b_{1}=0 \\
d &=&e_{1}=e_{2}=1.
\end{eqnarray*}%
Since $\alpha \left( x_{1};0,0,0,0\right) \equiv a+b_{1}+b_{2}\equiv 0$ we
get
\begin{equation*}
B(gx_{2}\otimes x_{1}x_{2};GX_{2},gx_{1}x_{2})GX_{2}\otimes
gx_{1}x_{2}\otimes gx_{1}.
\end{equation*}

Fourth summand gives us%
\begin{eqnarray*}
l_{1}+u_{1} &=&1,l_{2}=u_{2}=0 \\
a &=&b_{2}=1,b_{1}=l_{1}, \\
d &=&e_{2}=1,e_{1}-u_{1}=1\Rightarrow e_{1}=1,u_{1}=0,b_{1}=l_{1}=1.
\end{eqnarray*}%
Since $\alpha \left( 1_{H};1,0,0,0\right) \equiv b_{2}\equiv 1,$we get
\begin{equation*}
-B(x_{1}x_{2}\otimes x_{1}x_{2};GX_{1}X_{2},gx_{1}x_{2})GX_{2}\otimes
gx_{1}x_{2}\otimes gx_{1}.
\end{equation*}%
By considering also the right side of the equality, we get%
\begin{equation*}
B(x_{1}x_{2}\otimes x_{2};GX_{2},gx_{1}x_{2})-B(x_{1}x_{2}\otimes
x_{1}x_{2};GX_{1}X_{2},gx_{1}x_{2})+B(gx_{2}\otimes
x_{1}x_{2};GX_{2},gx_{1}x_{2})=0
\end{equation*}%
which holds in view of the form of the elements.

\subsubsection{Case $G\otimes gx_{1}x_{2}\otimes x_{1}x_{2}$}

First summand of the left side of the equality gives%
\begin{eqnarray*}
l_{1} &=&u_{1}=0,l_{2}=u_{2}=0 \\
a &=&1,b_{1}=b_{2}=0 \\
d &=&e_{1}=e_{2}=1.
\end{eqnarray*}%
Since $\alpha \left( x_{1}x_{2};0,0,0,0\right) \equiv 0,$ we get%
\begin{equation*}
B(1_{H}\otimes x_{1}x_{2};G,gx_{1}x_{2})G\otimes gx_{1}x_{2}\otimes
x_{1}x_{2}.
\end{equation*}%
Second summand of the left side of the equality gives us%
\begin{eqnarray*}
l_{1} &=&u_{1}=0,l_{2}+u_{2}=1 \\
a &=&1,b_{1}=0,b_{2}=l_{2}, \\
d &=&e_{1}=1,e_{2}-u_{2}=1\Rightarrow e_{2}=1,u_{2}=0,b_{2}=l_{2}=1.
\end{eqnarray*}%
Since $\alpha \left( x_{1};0,1,0,0\right) \equiv a+b_{1}+b_{2}+1\equiv 1$,
we get%
\begin{equation*}
-B(gx_{2}\otimes x_{1}x_{2};GX_{2},gx_{1}x_{2})G\otimes gx_{1}x_{2}\otimes
x_{1}x_{2}.
\end{equation*}%
Third summand of the left side of the equality gives us

\begin{eqnarray*}
l_{1}+u_{1} &=&1,l_{2}=u_{2}=0 \\
a &=&1,b_{1}=l_{1},b_{2}=0, \\
d &=&1,e_{1}-u_{1}=1\Rightarrow e_{1}=1,u_{1}=0,b_{1}=l_{1}=1,e_{2}=1.
\end{eqnarray*}%
Since $\alpha \left( x_{2};1,0,0,0\right) \equiv a+b_{1}\equiv 0,$ we get%
\begin{equation*}
-B(gx_{1}\otimes x_{1}x_{2};GX_{1},gx_{1}x_{2})G\otimes gx_{1}x_{2}\otimes
x_{1}x_{2}.
\end{equation*}%
Fourth summand gives us%
\begin{eqnarray*}
l_{1}+u_{1} &=&1,l_{2}+u_{2}=1, \\
a &=&1,b_{1}=l_{1},b_{2}=l_{2}, \\
d &=&1,e_{1}-u_{1}=1\Rightarrow e_{1}=1,u_{1}=0,b_{1}=l_{1}=1, \\
e_{2}-u_{2} &=&1\Rightarrow e_{2}=1,u_{2}=0,b_{2}=l_{2}=1.
\end{eqnarray*}%
Since $\alpha \left( 1_{H};1,1,0,0\right) \equiv 1+b_{2}\equiv 0,$ we obtain
\begin{equation*}
+B(x_{1}x_{2}\otimes x_{1}x_{2};GX_{1}X_{2},gx_{1}x_{2})G\otimes
gx_{1}x_{2}\otimes x_{1}x_{2}.
\end{equation*}%
By taking in account also the right side of the equality we get%
\begin{gather*}
-B(x_{1}x_{2}\otimes 1_{H};G,gx_{1}x_{2})+B(x_{1}x_{2}\otimes
x_{1}x_{2};GX_{1}X_{2},gx_{1}x_{2}) \\
-B(gx_{1}\otimes x_{1}x_{2};GX_{1},gx_{1}x_{2})-B(gx_{2}\otimes
x_{1}x_{2};GX_{2},gx_{1}x_{2}) \\
+B(1_{H}\otimes x_{1}x_{2};G,gx_{1}x_{2})=0
\end{gather*}

which holds in view of the form of the elements.

\section{$B(x_{1}x_{2}\otimes gx_{1})$}

By using$\left( \ref{simplgx}\right) ,$ we get%
\begin{eqnarray}
&&B(x_{1}x_{2}\otimes gx_{1})  \label{form x1x2otgx1} \\
&=&(1_{A}\otimes g)B(gx_{1}x_{2}\otimes 1_{H})(1_{A}\otimes gx_{1})  \notag
\\
&&+(1_{A}\otimes x_{1})B(gx_{1}x_{2}\otimes 1_{H})  \notag \\
&&-B(gx_{1}gx_{1}x_{2}\otimes 1_{H})  \notag
\end{eqnarray}%
and we obtain%
\begin{eqnarray*}
B\left( x_{1}x_{2}\otimes gx_{1}\right) &=&2B(gx_{1}x_{2}\otimes
1_{H};1_{A},g)1_{A}\otimes gx_{1}+ \\
&&-2B(gx_{1}x_{2}\otimes 1_{H};G,gx_{2})G\otimes gx_{1}x_{2}+ \\
&&+\left[ +2B(x_{2}\otimes 1_{H};1_{A},gx_{2})+2B(gx_{1}x_{2}\otimes
1_{H};1_{A},gx_{1}x_{2})\right] X_{1}\otimes gx_{1}x_{2}+ \\
&&-2B(x_{1}\otimes 1_{H};1_{A},gx_{2})X_{2}\otimes gx_{1}x_{2}+ \\
&&+\left[
\begin{array}{c}
2B(g\otimes 1_{H};1_{A},g)+2B(x_{2}\otimes \ 1_{H};1_{A},gx_{2}) \\
+2B(x_{1}\otimes 1_{H};1_{A},gx_{1})+2B(gx_{1}x_{2}\otimes
1_{H};1_{A},gx_{1}x_{2})%
\end{array}%
\right] X_{1}X_{2}\otimes gx_{1} \\
&&+\left[ 2B(x_{2}\otimes 1_{H};G,g)-2B(gx_{1}x_{2}\otimes 1_{H};G,gx_{1})%
\right] GX_{1}\otimes gx_{1}+ \\
&&+\left[ -2B(x_{1}\otimes 1_{H};G,g)-2B(gx_{1}x_{2}\otimes 1_{H};G,gx_{2})%
\right] GX_{2}\otimes gx_{1}+ \\
&&+\left[ -2B(g\otimes 1_{H};G,gx_{2})+2B(x_{1}\otimes 1_{H};G,gx_{1}x_{2})%
\right] GX_{1}X_{2}\otimes gx_{1}x_{2}+
\end{eqnarray*}%
We write Casimir formula for $B\left( x_{1}x_{2}\otimes gx_{1}\right) $%
\begin{eqnarray*}
&&\text{1}\sum_{a,b_{1},b_{2},d,e_{1},e_{2}=0}^{1}\sum_{l_{1}=0}^{b_{1}}%
\sum_{l_{2}=0}^{b_{2}}\sum_{u_{1}=0}^{e_{1}}\sum_{u_{2}=0}^{e_{2}}\left(
-1\right) ^{\alpha \left( x_{1}x_{2};l_{1},l_{2},u_{1},u_{2}\right) } \\
&&B(1_{H}\otimes
gx_{1};G^{a}X_{1}^{b_{1}}X_{2}^{b_{2}},g^{d}x_{1}^{e_{1}}x_{2}^{e_{2}}) \\
&&G^{a}X_{1}^{b_{1}-l_{1}}X_{2}^{b_{2}-l_{2}}\otimes
g^{d}x_{1}^{e_{1}-u_{1}}x_{2}^{e_{2}-u_{2}}\otimes
g^{a+b_{1}+b_{2}+l_{1}+l_{2}+d+e_{1}+e_{2}+u_{1}+u_{2}}x_{1}^{l_{1}+u_{1}+1}x_{2}^{l_{2}+u_{2}+1}+
\\
&&\text{2}+\sum_{a,b_{1},b_{2},d,e_{1},e_{2}=0}^{1}\sum_{l_{1}=0}^{b_{1}}%
\sum_{l_{2}=0}^{b_{2}}\sum_{u_{1}=0}^{e_{1}}\sum_{u_{2}=0}^{e_{2}}\left(
-1\right) ^{\alpha \left( x_{1};l_{1},l_{2},u_{1},u_{2}\right) } \\
&&B(gx_{2}\otimes
gx_{1};G^{a}X_{1}^{b_{1}}X_{2}^{b_{2}},g^{d}x_{1}^{e_{1}}x_{2}^{e_{2}}) \\
&&G^{a}X_{1}^{b_{1}-l_{1}}X_{2}^{b_{2}-l_{2}}\otimes
g^{d}x_{1}^{e_{1}-u_{1}}x_{2}^{e_{2}-u_{2}}\otimes
g^{a+b_{1}+b_{2}+l_{1}+l_{2}+d+e_{1}+e_{2}+u_{1}+u_{2}}x_{1}^{l_{1}+u_{1}+1}x_{2}^{l_{2}+u_{2}}+
\\
&&\text{3}\left( -1\right)
\sum_{a,b_{1},b_{2},d,e_{1},e_{2}=0}^{1}\sum_{l_{1}=0}^{b_{1}}%
\sum_{l_{2}=0}^{b_{2}}\sum_{u_{1}=0}^{e_{1}}\sum_{u_{2}=0}^{e_{2}}\left(
-1\right) ^{\alpha \left( x_{2};l_{1},l_{2},u_{1},u_{2}\right) } \\
&&B(gx_{1}\otimes
gx_{1};G^{a}X_{1}^{b_{1}}X_{2}^{b_{2}},g^{d}x_{1}^{e_{1}}x_{2}^{e_{2}}) \\
&&G^{a}X_{1}^{b_{1}-l_{1}}X_{2}^{b_{2}-l_{2}}\otimes
g^{d}x_{1}^{e_{1}-u_{1}}x_{2}^{e_{2}-u_{2}}\otimes
g^{a+b_{1}+b_{2}+l_{1}+l_{2}+d+e_{1}+e_{2}+u_{1}+u_{2}}x_{1}^{l_{1}+u_{1}}x_{2}^{l_{2}+u_{2}+1}+
\\
&&\text{4}\sum_{a,b_{1},b_{2},d,e_{1},e_{2}=0}^{1}\sum_{l_{1}=0}^{b_{1}}%
\sum_{l_{2}=0}^{b_{2}}\sum_{u_{1}=0}^{e_{1}}\sum_{u_{2}=0}^{e_{2}}\left(
-1\right) ^{\alpha \left( 1_{H};l_{1},l_{2},u_{1},u_{2}\right) } \\
&&B(x_{1}x_{2}\otimes
gx_{1};G^{a}X_{1}^{b_{1}}X_{2}^{b_{2}},g^{d}x_{1}^{e_{1}}x_{2}^{e_{2}}) \\
&&G^{a}X_{1}^{b_{1}-l_{1}}X_{2}^{b_{2}-l_{2}}\otimes
g^{d}x_{1}^{e_{1}-u_{1}}x_{2}^{e_{2}-u_{2}}\otimes
g^{a+b_{1}+b_{2}+l_{1}+l_{2}+d+e_{1}+e_{2}+u_{1}+u_{2}}x_{1}^{l_{1}+u_{1}}x_{2}^{l_{2}+u_{2}}
\\
&=&B^{A}(x_{1}x_{2}\otimes gx_{1})\otimes B^{H}(x_{1}x_{2}\otimes
gx_{1})\otimes 1_{H}+ \\
&&B^{A}(x_{1}x_{2}\otimes g)\otimes B^{H}(x_{1}x_{2}\otimes g)\otimes gx_{1}
\end{eqnarray*}

\subsection{$B\left( x_{1}x_{2}\otimes gx_{1};1_{A},gx_{1}\right) $ \ }

We deduce that%
\begin{eqnarray*}
a &=&b_{1}=b_{2}=0 \\
d &=&e_{1}=1,e_{2}=0
\end{eqnarray*}%
and we get
\begin{gather*}
+\left( -1\right) ^{\alpha \left( 1_{H};0,0,0,0\right) }B(x_{1}x_{2}\otimes
gx_{1};1_{A},gx_{1})1_{A}\otimes gx_{1}\otimes 1_{H} \\
+\left( -1\right) ^{\alpha \left( 1_{H};0,0,1,0\right) }B(x_{1}x_{2}\otimes
gx_{1};1_{A},gx_{1})1_{A}\otimes g\otimes gx_{1}
\end{gather*}

\subsubsection{Case $1_{A}\otimes g\otimes gx_{1}$}

We have to consider only the second and the fourth summand of the left side
of the equality. Second summand gives us

\begin{eqnarray*}
l_{1} &=&u_{1}=0,l_{2}=u_{2}=0 \\
a &=&b_{1}=b_{2}=0 \\
d &=&1,e_{1}=e_{2}=0.
\end{eqnarray*}%
Since $\alpha \left( x_{1};0,0,0,0\right) \equiv a+b_{1}+b_{2}\equiv 0,$ we
get%
\begin{equation*}
+B(gx_{2}\otimes gx_{1};1_{A},g)1_{A}\otimes g\otimes gx_{1}.
\end{equation*}

Fourth summand gives us%
\begin{eqnarray*}
l_{1}+u_{1} &=&1,l_{2}=u_{2}=0 \\
a &=&b_{2}=0,b_{1}=l_{1}, \\
d &=&1,e_{2}=0,e_{1}=u_{1}
\end{eqnarray*}%
Since $\alpha \left( 1_{H};0,0,1,0\right) \equiv e_{2}+\left(
a+b_{1}+b_{2}\right) \equiv 0$ and $\alpha \left( 1_{H};1,0,0,0\right)
\equiv b_{2}=0,$ we get
\begin{equation*}
\left[ B(x_{1}x_{2}\otimes gx_{1};1_{A},gx_{1})+B(x_{1}x_{2}\otimes
gx_{1};X_{1},g)\right] 1_{A}\otimes g\otimes gx_{1}.
\end{equation*}%
By considering also the right side, we obtain%
\begin{gather*}
B(x_{1}x_{2}\otimes gx_{1};1_{A},gx_{1})+B(x_{1}x_{2}\otimes gx_{1};X_{1},g)
\\
-B(x_{1}x_{2}\otimes g;1_{A},g)+B(gx_{2}\otimes gx_{1};1_{A},g)=0
\end{gather*}%
which holds in view of the form of the elements.

\subsection{$B\left( x_{1}x_{2}\otimes gx_{1};G,gx_{1}x_{2}\right) $}

We deduce that%
\begin{equation*}
a=1,b_{1}=b_{2}=0,d=e_{1}=e_{2}=1
\end{equation*}%
and we get%
\begin{eqnarray*}
&&\sum_{u_{1}=0}^{1}\sum_{u_{2}=0}^{1}\left( -1\right) ^{\alpha \left(
1_{H};0,0,u_{1},u_{2}\right) }B(x_{1}x_{2}\otimes
gx_{1};G,gx_{1}x_{2})G\otimes gx_{1}^{1-u_{1}}x_{2}^{1-u_{2}}\otimes
g^{u_{1}+u_{2}}x_{1}^{u_{1}}x_{2}^{u_{2}} \\
&=&\left( -1\right) ^{\alpha \left( 1_{H};0,0,0,0\right)
}B(x_{1}x_{2}\otimes gx_{1};G,gx_{1}x_{2})G\otimes gx_{1}x_{2}\otimes 1_{H}+
\\
&&\left( -1\right) ^{\alpha \left( 1_{H};0,0,0,1\right) }B(x_{1}x_{2}\otimes
gx_{1};G,gx_{1}x_{2})G\otimes gx_{1}\otimes gx_{2}+\text{ } \\
&&\left( -1\right) ^{\alpha \left( 1_{H};0,0,1,0\right) }B(x_{1}x_{2}\otimes
gx_{1};G,gx_{1}x_{2})G\otimes gx_{2}\otimes gx_{1} \\
&&\left( -1\right) ^{\alpha \left( 1_{H};0,0,1,1\right) }B(x_{1}x_{2}\otimes
gx_{1};G,gx_{1}x_{2})G\otimes g\otimes x_{1}x_{2}.
\end{eqnarray*}

\subsubsection{Case $G\otimes gx_{1}\otimes gx_{2}$}

We have only to consider the third and the forth summands of the left side
of the equality. Third summand gives us%
\begin{eqnarray*}
l_{1} &=&u_{1}=l_{2}=u_{2}=0, \\
a &=&1,b_{1}=b_{2}=0, \\
d &=&e_{1}=1,e_{2}=0,
\end{eqnarray*}%
Since $\alpha \left( x_{2};0,0,0,0\right) \equiv a+b_{1}+b_{2}=1$ we get
\begin{equation*}
B(gx_{1}\otimes gx_{1};G,gx_{1})G\otimes gx_{1}\otimes gx_{2}.
\end{equation*}%
Fourth summand gives us%
\begin{eqnarray*}
l_{1} &=&u_{1}=0,l_{2}+u_{2}=1 \\
a &=&1,b_{1}=0,b_{2}=l_{2}, \\
d &=&e_{1}=1,e_{2}=u_{2}.
\end{eqnarray*}%
Since $\alpha \left( 1_{H};0,0,0,1\right) \equiv a+b_{1}+b_{2}\equiv 1$ and $%
\alpha \left( 1_{H};0,1,0,0\right) \equiv 0,$ we get
\begin{equation*}
\left[ -B(x_{1}x_{2}\otimes gx_{1};G,gx_{1}x_{2})+B(x_{1}x_{2}\otimes
gx_{1};GX_{2},gx_{1})\right] G\otimes gx_{1}\otimes gx_{2}
\end{equation*}%
Since there is no term to consider in the right side of the equality we
obtain%
\begin{equation*}
-B(x_{1}x_{2}\otimes gx_{1};G,gx_{1}x_{2})+B(x_{1}x_{2}\otimes
gx_{1};GX_{2},gx_{1})+B(gx_{1}\otimes gx_{1};G,gx_{1})=0
\end{equation*}%
which holds in view of the form of the elements.

\subsubsection{Case $G\otimes gx_{2}\otimes gx_{1}$}

We have to consider only the second and the fourth summand of the left side
of the equality. Second summand gives us

\begin{eqnarray*}
l_{1} &=&u_{1}=0,l_{2}=u_{2}=0 \\
a &=&1,b_{1}=b_{2}=0 \\
d &=&e_{2}=1,e_{1}=0.
\end{eqnarray*}%
Since $\alpha \left( x_{1};0,0,0,0\right) \equiv a+b_{1}+b_{2}\equiv 1,$ we
get%
\begin{equation*}
-B(gx_{2}\otimes gx_{1};G_{1},gx_{2})G\otimes gx_{2}\otimes gx_{1}.
\end{equation*}

Fourth summand gives us%
\begin{eqnarray*}
l_{1}+u_{1} &=&1,l_{2}=u_{2}=0 \\
a &=&1,b_{1}=l_{1},b_{2}=0, \\
d &=&e_{2}=1,e_{1}=u_{1}
\end{eqnarray*}%
Since $\alpha \left( 1_{H};0,0,1,0\right) \equiv e_{2}+\left(
a+b_{1}+b_{2}\right) \equiv 0$ and $\alpha \left( 1_{H};1,0,0,0\right)
\equiv b_{2}=0,$ we get
\begin{equation*}
\left[ B(x_{1}x_{2}\otimes gx_{1};G,gx_{1}x_{2})+B(x_{1}x_{2}\otimes
gx_{1};GX_{1},gx_{2})\right] G\otimes gx_{2}\otimes gx_{1}.
\end{equation*}%
By considering also the right side, we obtain%
\begin{gather*}
-B(x_{1}x_{2}\otimes g;G,gx_{2})+B(x_{1}x_{2}\otimes gx_{1};G,gx_{1}x_{2}) \\
+B(x_{1}x_{2}\otimes gx_{1};GX_{1},gx_{2})-B(gx_{2}\otimes
gx_{1};G_{1},gx_{2})=0
\end{gather*}%
which holds in view of the form of the elements.

\subsubsection{Case $G\otimes g\otimes x_{1}x_{2}$}

First summand of the left side of the equality gives us

\begin{eqnarray*}
l_{1} &=&u_{1}=0,l_{2}=u_{2}=0 \\
a &=&1,b_{1}=b_{2}=0, \\
d &=&1,e_{1}=e_{2}=0.
\end{eqnarray*}%
Since $\alpha \left( x_{1}x_{2};0,0,0,0\right) \equiv 0,$ we get
\begin{equation*}
B(1_{H}\otimes gx_{1};G,g)G\otimes g\otimes x_{1}x_{2}
\end{equation*}%
Second summand of the left side gives us%
\begin{eqnarray*}
l_{1} &=&u_{1}=0,l_{2}+u_{2}=1 \\
a &=&1,b_{1}=0,b_{2}=l_{2}, \\
d &=&1,e_{1}=0,e_{2}=u_{2}.
\end{eqnarray*}%
Since $\alpha \left( x_{1};0,0,0,1\right) \equiv 1$ and $\alpha \left(
x_{1};0,1,0,0\right) \equiv a+b_{1}+b_{2}+1\equiv 1,$ we obtain
\begin{equation*}
\left[ -B(gx_{2}\otimes gx_{1};G,gx_{2})-B(gx_{2}\otimes gx_{1};GX_{2},g)%
\right] G\otimes g\otimes x_{1}x_{2}.
\end{equation*}%
Third summand of the left side gives us%
\begin{eqnarray*}
l_{1}+u_{1} &=&1,l_{2}=u_{2}=0, \\
a &=&1,b_{1}=l_{1},b_{2}=0, \\
d &=&1,e_{1}=u_{1},e_{2}=0.
\end{eqnarray*}%
Since $\alpha \left( x_{2};0,0,1,0\right) \equiv e_{2}=0$ and $\alpha \left(
x_{2};1,0,0,0\right) \equiv a+b_{1}\equiv 0,$ we get%
\begin{equation*}
\left[ -B(gx_{1}\otimes gx_{1};G,gx_{1})-B(gx_{1}\otimes gx_{1};GX_{1},g%
\right] G\otimes g\otimes x_{1}x_{2}
\end{equation*}%
The fourth summand of the left side gives us%
\begin{eqnarray*}
l_{1}+u_{1} &=&1,l_{2}+u_{2}=1 \\
a &=&1,b_{1}=l_{1},b_{2}=l_{2}, \\
d &=&1,e_{1}=u_{1},e_{2}=u_{2}.
\end{eqnarray*}%
Since
\begin{eqnarray*}
\alpha \left( 1_{H};0,0,1,1\right) &\equiv &1+e_{2}\equiv 0 \\
\alpha \left( 1_{H};0,1,1,0\right) &\equiv &e_{2}+\left(
a+b_{1}+b_{2}+1\right) \equiv 1 \\
\alpha \left( 1_{H};1,0,0,1\right) &\equiv &a+b_{1}\equiv 0 \\
\alpha \left( 1_{H};1,1,0,0\right) &\equiv &1+b_{2}\equiv 0
\end{eqnarray*}%
we get%
\begin{equation*}
\left[
\begin{array}{c}
B(x_{1}x_{2}\otimes gx_{1};G,gx_{1}x_{2})-B(x_{1}x_{2}\otimes
gx_{1};GX_{2},gx_{1}) \\
+B(x_{1}x_{2}\otimes gx_{1};GX_{1},gx_{2})+B(x_{1}x_{2}\otimes
gx_{1};GX_{1}X_{2},g)%
\end{array}%
\right] G\otimes g\otimes x_{1}x_{2}.
\end{equation*}%
Since there is no term to consider in the right side of the equality we
obtain

\begin{gather*}
B(1_{H}\otimes gx_{1};G,g)-B(gx_{2}\otimes gx_{1};G,gx_{2}) \\
-B(gx_{2}\otimes gx_{1};GX_{2},g)-B(gx_{1}\otimes gx_{1};G,gx_{1}) \\
-B(gx_{1}\otimes gx_{1};GX_{1},g)B(x_{1}x_{2}\otimes gx_{1};G,gx_{1}x_{2}) \\
-B(x_{1}x_{2}\otimes gx_{1};GX_{2},gx_{1})+B(x_{1}x_{2}\otimes
gx_{1};GX_{1},gx_{2}) \\
+B(x_{1}x_{2}\otimes gx_{1};GX_{1}X_{2},g)=0
\end{gather*}%
which holds in view of the form of the elements.

\subsection{$B\left( x_{1}x_{2}\otimes gx_{1};X_{1},gx_{1}x_{2}\right) $ \ }

We deduce that%
\begin{eqnarray*}
a &=&b_{2}=0 \\
b_{1} &=&d=e_{1}=e_{2}=1
\end{eqnarray*}%
and we get%
\begin{eqnarray*}
&&\sum_{u_{1}=0}^{1}\sum_{u_{2}=0}^{1}\left( -1\right) ^{\alpha \left(
1_{H};0,0,u_{1},u_{2}\right) }B(x_{1}x_{2}\otimes
gx_{1};X_{1},gx_{1}x_{2})G\otimes gx_{1}^{1-u_{1}}x_{2}^{1-u_{2}}\otimes
g^{u_{1}+u_{2}}x_{1}^{u_{1}}x_{2}^{u_{2}} \\
&=&\left( -1\right) ^{\alpha \left( 1_{H};0,0,0,0\right)
}B(x_{1}x_{2}\otimes gx_{1};X_{1},gx_{1}x_{2})X_{1}\otimes
gx_{1}x_{2}\otimes 1_{H}+ \\
&&\left( -1\right) ^{\alpha \left( 1_{H};0,0,0,1\right) }B(x_{1}x_{2}\otimes
gx_{1};X_{1},gx_{1}x_{2})X_{1}\otimes gx_{1}\otimes gx_{2}+ \\
&&\left( -1\right) ^{\alpha \left( 1_{H};0,0,1,0\right) }B(x_{1}x_{2}\otimes
gx_{1};X_{1},gx_{1}x_{2})X_{1}\otimes gx_{2}\otimes gx_{1}\text{ } \\
&&\left( -1\right) ^{\alpha \left( 1_{H};0,0,1,1\right) }B(x_{1}x_{2}\otimes
gx_{1};X_{1},gx_{1}x_{2})X_{1}\otimes g\otimes x_{1}x_{2}\text{ }
\end{eqnarray*}

\subsubsection{Case $X_{1}\otimes gx_{1}\otimes gx_{2}$}

We have only to consider the third and the forth summands of the left side
of the equality. Third summand gives us%
\begin{eqnarray*}
l_{1} &=&u_{1}=l_{2}=u_{2}=0, \\
a &=&b_{2}=0,b_{1}=1, \\
d &=&e_{1}=1,e_{2}=0.
\end{eqnarray*}%
Since $\alpha \left( x_{2};0,0,0,0\right) \equiv a+b_{1}+b_{2}=1$ we get%
\begin{equation*}
B(gx_{1}\otimes gx_{1};X_{1},gx_{1})X_{1}\otimes gx_{1}\otimes gx_{2}
\end{equation*}%
Fourth summand gives us%
\begin{eqnarray*}
l_{1} &=&u_{1}=0,l_{2}+u_{2}=1 \\
a &=&0,b_{1}=1,b_{2}=l_{2}, \\
d &=&e_{1}=1,e_{2}=u_{2}.
\end{eqnarray*}%
Since $\alpha \left( 1_{H};0,0,0,1\right) \equiv a+b_{1}+b_{2}\equiv 1$ and $%
\alpha \left( 1_{H};0,1,0,0\right) \equiv 0,$ we get%
\begin{equation*}
\left[ -B(x_{1}x_{2}\otimes gx_{1};X_{1},gx_{1}x_{2})+B(x_{1}x_{2}\otimes
gx_{1};X_{1}X_{2},gx_{1})\right] X_{1}\otimes gx_{1}\otimes gx_{2}.
\end{equation*}%
Since there is no term to consider in the right side of the equality we
obtain%
\begin{equation*}
-B(x_{1}x_{2}\otimes gx_{1};X_{1},gx_{1}x_{2})+B(x_{1}x_{2}\otimes
gx_{1};X_{1}X_{2},gx_{1})+B(gx_{1}\otimes gx_{1};x_{1},gx_{1})=0
\end{equation*}%
which holds in view of the form of the elements.

\subsubsection{Case $X_{1}\otimes gx_{2}\otimes gx_{1}$}

We have to consider only the second and the fourth summand of the left side
of the equality. Second summand gives us

\begin{eqnarray*}
l_{1} &=&u_{1}=0,l_{2}=u_{2}=0 \\
a &=&b_{2}=0,b_{1}=1 \\
d &=&e_{2}=1,e_{1}=0.
\end{eqnarray*}%
Since $\alpha \left( x_{1};0,0,0,0\right) \equiv a+b_{1}+b_{2}\equiv 1,$ we
get%
\begin{equation*}
-B(gx_{2}\otimes gx_{1};X_{1},gx_{2})X_{1}\otimes gx_{2}\otimes gx_{1}
\end{equation*}

Fourth summand gives us%
\begin{eqnarray*}
l_{1}+u_{1} &=&1,l_{2}=u_{2}=0 \\
a &=&1,b_{2}=0,b_{1}-l_{1}=1\Rightarrow b_{1}=1,l_{1}=0,u_{1} \\
d &=&e_{2}=1,e_{1}=u_{1}=1.
\end{eqnarray*}%
Since $\alpha \left( 1_{H};0,0,1,0\right) \equiv e_{2}+\left(
a+b_{1}+b_{2}\right) \equiv 0,$ we get
\begin{equation*}
B(x_{1}x_{2}\otimes gx_{1};X_{1},gx_{1}x_{2})X_{1}\otimes gx_{2}\otimes
gx_{1}
\end{equation*}%
By considering also the right side, we obtain%
\begin{equation*}
-B(x_{1}x_{2}\otimes g;X_{1},gx_{2})+B(x_{1}x_{2}\otimes
gx_{1};X_{1},gx_{1}x_{2})-B(gx_{2}\otimes gx_{1};X_{1},gx_{2})=0
\end{equation*}%
which holds in view of the form of the elements.

\subsubsection{Case $X_{1}\otimes g\otimes x_{1}x_{2}$}

First summand of the left side of the equality gives us%
\begin{eqnarray*}
l_{1} &=&u_{1}=0,l_{2}=u_{2}=0 \\
a &=&b_{2}=0,b_{1}=1, \\
d &=&1,e_{1}=e_{2}=0.
\end{eqnarray*}

Since $\alpha \left( x_{1}x_{2};0,0,0,0\right) \equiv 0,$ we get
\begin{equation*}
B(1_{H}\otimes gx_{1};X_{1},g)X_{1}\otimes g\otimes x_{1}x_{2}
\end{equation*}%
Second summand of the left side gives us%
\begin{eqnarray*}
l_{1} &=&u_{1}=0,l_{2}+u_{2}=1 \\
a &=&0,b_{1}=1,b_{2}=l_{2}, \\
d &=&1,e_{1}=0,e_{2}=u_{2}.
\end{eqnarray*}%
Since $\alpha \left( x_{1};0,0,0,1\right) \equiv 1$ and $\alpha \left(
x_{1};0,1,0,0\right) \equiv a+b_{1}+b_{2}+1\equiv 1,$ we obtain
\begin{equation*}
\left[ -B(gx_{2}\otimes gx_{1};X_{1},gx_{2})-B(gx_{2}\otimes
gx_{1};X_{1}X_{2},g)\right] X_{1}\otimes g\otimes x_{1}x_{2}
\end{equation*}%
Third summand of the left side gives us%
\begin{eqnarray*}
l_{1}+u_{1} &=&1,l_{2}=u_{2}=0, \\
a &=&0,b_{2}=0,b_{1}=l_{1},b_{1}-l_{1}=1\Rightarrow b_{1}=1,l_{1}=0,u_{1}=1
\\
d &=&1,e_{1}=u_{1}=1,e_{2}=0.
\end{eqnarray*}%
Since $\alpha \left( x_{2};0,0,1,0\right) \equiv e_{2}=0$ $,$ we get%
\begin{equation*}
\left[ -B(gx_{1}\otimes gx_{1};X_{1},gx_{1})\right] X_{1}\otimes g\otimes
x_{1}x_{2}
\end{equation*}%
The fourth summand of the left side gives us%
\begin{eqnarray*}
l_{1}+u_{1} &=&1,l_{2}+u_{2}=1 \\
a &=&0,b_{2}=l_{2}, \\
b_{1}-l_{1} &=&1\Rightarrow b_{1}=1,l_{1}=0,u_{1} \\
d &=&1,e_{1}=u_{1}=1,e_{2}=u_{2}.
\end{eqnarray*}%
Since $\alpha \left( 1_{H};0,0,1,1\right) \equiv 1+e_{2}\equiv 0$ and $%
\alpha \left( 1_{H};0,1,1,0\right) \equiv e_{2}+\left(
a+b_{1}+b_{2}+1\right) \equiv 1$ we get%
\begin{equation*}
\left[ B(x_{1}x_{2}\otimes gx_{1};X_{1},gx_{1}x_{2})-B(x_{1}x_{2}\otimes
gx_{1};X_{1}X_{2},gx_{1})\right] X_{1}\otimes g\otimes x_{1}x_{2}
\end{equation*}%
Since there is no term to consider in the right side of the equality we
obtain

\begin{eqnarray*}
&&-B(gx_{1}\otimes gx_{1};X_{1},gx_{1})-B(gx_{2}\otimes
gx_{1};X_{1},gx_{2})-B(gx_{2}\otimes gx_{1};X_{1}X_{2},g) \\
&&B(x_{1}x_{2}\otimes gx_{1};X_{1},gx_{1}x_{2})-B(x_{1}x_{2}\otimes
gx_{1};X_{1}X_{2},gx_{1})+B(1_{H}\otimes gx_{1};X_{1},g) \\
&=&0
\end{eqnarray*}%
which holds in view of the form of the elements.

\subsection{$B\left( x_{1}x_{2}\otimes gx_{1};X_{2},gx_{1}x_{2}\right) $}

We deduce that%
\begin{eqnarray*}
a &=&0,b_{1}=0,b_{2}=1 \\
d &=&e_{1}=e_{2}=1
\end{eqnarray*}%
and we get%
\begin{eqnarray*}
&&\left( -1\right) ^{\alpha \left( 1_{H};0,0,0,0\right) }B(x_{1}x_{2}\otimes
gx_{1};X_{2},gx_{1}x_{2})X_{2}\otimes gx_{1}x_{2}\otimes 1_{H} \\
&&\left( -1\right) ^{\alpha \left( 1_{H};0,0,0,1\right) }B(x_{1}x_{2}\otimes
gx_{1};X_{2},gx_{1}x_{2})X_{2}\otimes gx_{1}\otimes gx_{2} \\
&&\left( -1\right) ^{\alpha \left( 1_{H};0,1,0,0\right) }B(x_{1}x_{2}\otimes
gx_{1};X_{2},gx_{1}x_{2})1_{A}\otimes gx_{1}x_{2}\otimes gx_{2} \\
&&\left( -1\right) ^{\alpha \left( 1_{H};0,1,0,1\right) }B(x_{1}x_{2}\otimes
gx_{1};X_{2},gx_{1}x_{2})X_{2}\otimes gx_{1}\otimes
g^{l_{2}+u_{1}+u_{2}}x_{1}^{u_{1}}x_{2}^{1+1}=0 \\
&&\left( -1\right) ^{\alpha \left( 1_{H};0,0,1,0\right) }B(x_{1}x_{2}\otimes
gx_{1};X_{2},gx_{1}x_{2})X_{2}\otimes gx_{2}\otimes gx_{1} \\
&&\left( -1\right) ^{\alpha \left( 1_{H};0,0,1,1\right) }B(x_{1}x_{2}\otimes
gx_{1};X_{2},gx_{1}x_{2})X_{2}\otimes g\otimes x_{1}x_{2} \\
&&\left( -1\right) ^{\alpha \left( 1_{H};0,1,1,0\right) }B(x_{1}x_{2}\otimes
gx_{1};X_{2},gx_{1}x_{2})1_{A}\otimes gx_{2}\otimes x_{1}x_{2} \\
&&\left( -1\right) ^{\alpha \left( 1_{H};0,1,1,1\right) }B(x_{1}x_{2}\otimes
gx_{1};X_{2},gx_{1}x_{2})X_{2}^{1-l_{2}}\otimes
gx_{1}^{1-u_{1}}x_{2}^{1-u_{2}}\otimes
g^{l_{2}+u_{1}+u_{2}}x_{1}^{u_{1}}x_{2}^{1+1}=0
\end{eqnarray*}

\subsubsection{Case $X_{2}\otimes gx_{1}\otimes gx_{2}$}

We have only to consider the third and the forth summands of the left side
of the equality. Third summand gives us%
\begin{eqnarray*}
l_{1} &=&u_{1}=l_{2}=u_{2}=0, \\
a &=&b_{1}=0,b_{2}=1, \\
d &=&e_{1}=1,e_{2}=0,
\end{eqnarray*}%
Since $\alpha \left( x_{2};0,0,0,0\right) \equiv a+b_{1}+b_{2}=1$ we get
\begin{equation*}
B(gx_{1}\otimes gx_{1};X_{2},gx_{1})X_{2}\otimes gx_{1}\otimes gx_{2}.
\end{equation*}%
Fourth summand gives us%
\begin{eqnarray*}
l_{1} &=&u_{1}=0,l_{2}+u_{2}=1 \\
a &=&b_{1}=0,b_{2}-l_{2}=1\Rightarrow b_{2}=1,l_{2}=0,u_{2}=1 \\
d &=&e_{1}=1,e_{2}=u_{2}=1.
\end{eqnarray*}%
Since $\alpha \left( 1_{H};0,0,0,1\right) \equiv a+b_{1}+b_{2}\equiv 1,$ we
get
\begin{equation*}
-B(x_{1}x_{2}\otimes gx_{1};X_{2},gx_{1}x_{2})X_{2}\otimes gx_{1}\otimes
gx_{2}
\end{equation*}%
Since there is no term to consider in the right side of the equality we
obtain%
\begin{equation*}
B(gx_{1}\otimes gx_{1};X_{2},gx_{1})-B(x_{1}x_{2}\otimes
gx_{1};X_{2},gx_{1}x_{2})=0
\end{equation*}%
which holds in view of the form of the elements.

\subsubsection{Case $1_{A}\otimes gx_{1}x_{2}\otimes gx_{2}$}

We have only to consider the third and the forth summands of the left side
of the equality. Third summand gives us%
\begin{eqnarray*}
l_{1} &=&u_{1}=l_{2}=u_{2}=0, \\
a &=&b_{1}=b_{2}=0, \\
d &=&e_{1}=e_{2}=1,
\end{eqnarray*}%
Since $\alpha \left( x_{2};0,0,0,0\right) \equiv a+b_{1}+b_{2}=0$ we get
\begin{equation*}
-B(gx_{1}\otimes gx_{1};1_{A},gx_{1}x_{2})1_{A}\otimes gx_{1}x_{2}\otimes
gx_{2}.
\end{equation*}%
Fourth summand gives us%
\begin{eqnarray*}
l_{1} &=&u_{1}=0,l_{2}+u_{2}=1 \\
a &=&b_{1}=0,b_{2}=l_{2} \\
d &=&e_{1}=1,e_{2}-u_{2}=1\Rightarrow e_{2}=1,u_{2}=0,b_{2}=l_{2}=1.
\end{eqnarray*}%
Since $\alpha \left( 1_{H};0,1,0,0\right) =0,$ we get
\begin{equation*}
B(x_{1}x_{2}\otimes gx_{1};X_{2},gx_{1}x_{2})1_{A}\otimes gx_{1}x_{2}\otimes
gx_{2}
\end{equation*}%
Since there is no term to consider in the right side of the equality we
obtain%
\begin{equation*}
-B(gx_{1}\otimes gx_{1};1_{A},gx_{1}x_{2})+B(x_{1}x_{2}\otimes
gx_{1};X_{2},gx_{1}x_{2})=0
\end{equation*}%
which holds in view of the form of the elements.

\subsubsection{Case $X_{2}\otimes gx_{2}\otimes gx_{1}$}

We have to consider only the second and the fourth summand of the left side
of the equality. Second summand gives us

\begin{eqnarray*}
l_{1} &=&u_{1}=0,l_{2}=u_{2}=0 \\
a &=&b_{1}=0,b_{2}=1 \\
d &=&e_{2}=1,e_{1}=0.
\end{eqnarray*}%
Since $\alpha \left( x_{1};0,0,0,0\right) \equiv a+b_{1}+b_{2}\equiv 1,$ we
get%
\begin{equation*}
-B(gx_{2}\otimes gx_{1};X_{2},gx_{2})X_{2}\otimes gx_{2}\otimes gx_{1}
\end{equation*}

Fourth summand gives us%
\begin{eqnarray*}
l_{1}+u_{1} &=&1,l_{2}=u_{2}=0 \\
a &=&0,b_{2}=1,b_{1}=l_{1} \\
d &=&e_{2}=1,e_{1}=u_{1}.
\end{eqnarray*}%
Since $\alpha \left( 1_{H};0,0,1,0\right) \equiv e_{2}+\left(
a+b_{1}+b_{2}\right) \equiv 0$ and $\alpha \left( 1_{H};1,0,0,0\right)
=b_{2}=1$ we get
\begin{equation*}
\left[ B(x_{1}x_{2}\otimes gx_{1};X_{2},gx_{1}x_{2})-B(x_{1}x_{2}\otimes
gx_{1};X_{1}X_{2},gx_{2})\right] X_{2}\otimes gx_{2}\otimes gx_{1}
\end{equation*}%
By considering also the right side, we obtain%
\begin{equation*}
-B(x_{1}x_{2}\otimes g;X_{2},gx_{2})+B(x_{1}x_{2}\otimes
gx_{1};X_{2},gx_{1}x_{2})-B(x_{1}x_{2}\otimes
gx_{1};X_{1}X_{2},gx_{2})-B(gx_{2}\otimes gx_{1};X_{2},gx_{2})=0
\end{equation*}%
which holds in view of the form of the elements.

\subsubsection{Case $X_{2}\otimes g\otimes x_{1}x_{2}$}

First summand of the left side of the equality gives us%
\begin{eqnarray*}
l_{1} &=&u_{1}=0,l_{2}=u_{2}=0 \\
a &=&b_{1}=0,b_{2}=1, \\
d &=&1,e_{1}=e_{2}=0.
\end{eqnarray*}

Since $\alpha \left( x_{1}x_{2};0,0,0,0\right) \equiv 0,$ we get
\begin{equation*}
B(1_{H}\otimes gx_{1};X_{2},g)X_{2}\otimes g\otimes x_{1}x_{2}
\end{equation*}%
Second summand of the left side gives us%
\begin{eqnarray*}
l_{1} &=&u_{1}=0,l_{2}+u_{2}=1 \\
a &=&b_{1}=0,1,b_{2}-l_{2}=1\Rightarrow b_{2}=1,l_{2}=0,u_{2}=1, \\
d &=&1,e_{1}=0,e_{2}=u_{2}=1.
\end{eqnarray*}%
Since $\alpha \left( x_{1};0,0,0,1\right) \equiv 1$ $,$ we obtain
\begin{equation*}
-B(gx_{2}\otimes gx_{1};X_{2},gx_{2})X_{2}\otimes g\otimes x_{1}x_{2}
\end{equation*}%
Third summand of the left side gives us%
\begin{eqnarray*}
l_{1}+u_{1} &=&1,l_{2}=u_{2}=0, \\
a &=&0,b_{2}=1,b_{1}=l_{1} \\
d &=&1,e_{1}=u_{1},e_{2}=0.
\end{eqnarray*}%
Since $\alpha \left( x_{2};0,0,1,0\right) \equiv e_{2}=0$ and $,$ $\alpha
\left( x_{2};1,0,0,0\right) =a+b_{1}\equiv 1$we get%
\begin{equation*}
\left[ -B(gx_{1}\otimes gx_{1};X_{2},gx_{1})+B(gx_{1}\otimes
gx_{1};X_{1}X_{2},g)\right] X_{2}\otimes g\otimes x_{1}x_{2}
\end{equation*}%
The fourth summand of the left side gives us%
\begin{eqnarray*}
l_{1}+u_{1} &=&1,l_{2}+u_{2}=1 \\
a &=&0,b_{1}=l_{1}, \\
b_{2}-l_{2} &=&1\Rightarrow b_{2}=1,l_{2}=0,u_{2}=1, \\
d &=&1,e_{1}=u_{1},e_{2}=u_{2}=1.
\end{eqnarray*}%
Since $\alpha \left( 1_{H};0,0,1,1\right) =1+e_{2}\equiv 0$ and $\alpha
\left( 1_{H};1,0,0,1\right) \equiv a+b_{1}=1,$ we get%
\begin{equation*}
\left[ B(x_{1}x_{2}\otimes gx_{1};X_{2},gx_{1}x_{2})-B(x_{1}x_{2}\otimes
gx_{1};X_{1}X_{2},gx_{2})\right] X_{2}\otimes g\otimes x_{1}x_{2}
\end{equation*}%
Since there is no term to consider in the right side of the equality we
obtain

\begin{eqnarray*}
&&-B(gx_{1}\otimes gx_{1};X_{2},gx_{1})+B(gx_{1}\otimes
gx_{1};X_{1}X_{2},g)-B(gx_{2}\otimes gx_{1};X_{2},gx_{2}) \\
&&+B(x_{1}x_{2}\otimes gx_{1};X_{2},gx_{1}x_{2})-B(x_{1}x_{2}\otimes
gx_{1};X_{1}X_{2},gx_{2})+B(1_{H}\otimes gx_{1};X_{2},g) \\
&=&0
\end{eqnarray*}%
which holds in view of the form of the elements.

\subsubsection{Case $1_{A}\otimes gx_{2}\otimes x_{1}x_{2}$}

First summand of the left side of the equality gives us%
\begin{eqnarray*}
l_{1} &=&u_{1}=0,l_{2}=u_{2}=0 \\
a &=&b_{1}=b_{2}=0, \\
d &=&e_{2}=1,e_{1}=0.
\end{eqnarray*}

Since $\alpha \left( x_{1}x_{2};0,0,0,0\right) \equiv 0,$ we get
\begin{equation*}
B(1_{H}\otimes gx_{1};1_{A},gx_{2})1_{A}\otimes gx_{2}\otimes x_{1}x_{2}
\end{equation*}%
Second summand of the left side gives us%
\begin{eqnarray*}
l_{1} &=&u_{1}=0,l_{2}+u_{2}=1 \\
a &=&b_{1}=0,b_{2}=l_{2}, \\
d &=&1,e_{1}=0,e_{2}-u_{2}=1\Rightarrow e_{2}=1,u_{2}=0,b_{2}=l_{2}=1.
\end{eqnarray*}%
Since $\alpha \left( x_{1};0,1,0,0\right) \equiv a+b_{1}+b_{2}+1\equiv 0$ $,$
we obtain
\begin{equation*}
B(gx_{2}\otimes gx_{1};X_{2},gx_{2})1_{A}\otimes gx_{2}\otimes x_{1}x_{2}
\end{equation*}%
Third summand of the left side gives us%
\begin{eqnarray*}
l_{1}+u_{1} &=&1,l_{2}=u_{2}=0, \\
a &=&b_{2}=0,b_{1}=l_{1} \\
d &=&1,e_{1}=u_{1},e_{2}=1.
\end{eqnarray*}%
Since $\alpha \left( x_{2};0,0,1,0\right) \equiv e_{2}=1$ and $,$ $\alpha
\left( x_{2};1,0,0,0\right) =a+b_{1}\equiv 1$we get%
\begin{equation*}
\left[ +B(gx_{1}\otimes gx_{1};1_{A},gx_{1}x_{2})++B(gx_{1}\otimes
gx_{1};X_{1},gx_{2})\right] 1_{A}\otimes gx_{2}\otimes x_{1}x_{2}
\end{equation*}%
The fourth summand of the left side gives us%
\begin{eqnarray*}
l_{1}+u_{1} &=&1,l_{2}+u_{2}=1 \\
a &=&0,b_{1}=l_{1},b_{2}-l_{2}=0 \\
d &=&1,e_{1}=u_{1},e_{2}-u_{2}=1\Rightarrow e_{2}=1,u_{2}=0,b_{2}=l_{2}=1.
\end{eqnarray*}%
Since $\alpha \left( 1_{H};0,1,1,0\right) =e_{2}+a+b_{1}+b_{2}+1\equiv 1$
and $\alpha \left( 1_{H};1,1,0,0\right) \equiv 1+b_{2}\equiv 0,$ we get%
\begin{equation*}
\left[ -B(x_{1}x_{2}\otimes gx_{1};X_{2},gx_{1}x_{2})+B(x_{1}x_{2}\otimes
gx_{1};X_{1}X_{2},gx_{2})\right] 1_{A}\otimes gx_{2}\otimes x_{1}x_{2}
\end{equation*}%
Since there is no term to consider in the right side of the equality we
obtain

\begin{eqnarray*}
&&+B(gx_{1}\otimes gx_{1};1_{A},gx_{1}x_{2})++B(gx_{1}\otimes
gx_{1};X_{1},gx_{2})+B(gx_{2}\otimes gx_{1};X_{2},gx_{2}) \\
&&-B(x_{1}x_{2}\otimes gx_{1};X_{2},gx_{1}x_{2})+B(x_{1}x_{2}\otimes
gx_{1};X_{1}X_{2},gx_{2})+B(1_{H}\otimes gx_{1};1_{A},gx_{2}) \\
&=&0
\end{eqnarray*}%
which holds in view of the form of the elements.

\subsection{$B\left( x_{1}x_{2}\otimes gx_{1};X_{1}X_{2},gx_{1}\right) $}

We deduce that%
\begin{equation*}
a=e_{2}=0,b_{1}=b_{2}=1,d=e_{1}=1.
\end{equation*}%
and we get%
\begin{eqnarray*}
&&\left( -1\right) ^{\alpha \left( 1_{H};0,0,0,0\right) }B(x_{1}x_{2}\otimes
gx_{1};X_{1}X_{2},gx_{1})X_{1}X_{2}\otimes gx_{1}\otimes 1_{H}+\text{ } \\
&&\left( -1\right) ^{\alpha \left( 1_{H};0,0,1,0\right) }B(x_{1}x_{2}\otimes
gx_{1};X_{1}X_{2},gx_{1})X_{1}X_{2}\otimes g\otimes gx_{1} \\
&&\left( -1\right) ^{\alpha \left( 1_{H};0,1,0,0\right) }B(x_{1}x_{2}\otimes
gx_{1};X_{1}X_{2},gx_{1})X_{1}\otimes gx_{1}\otimes gx_{2}+ \\
&&\left( -1\right) ^{\alpha \left( 1_{H};0,1,1,0\right) }B(x_{1}x_{2}\otimes
gx_{1};X_{1}X_{2},gx_{1})X_{1}\otimes g\otimes x_{1}x_{2}+ \\
&&\left( -1\right) ^{\alpha \left( 1_{H};1,0,0,0\right) }B(x_{1}x_{2}\otimes
gx_{1};X_{1}X_{2},gx_{1})X_{2}\otimes gx_{1}\otimes gx_{1}+ \\
&&\left( -1\right) ^{\alpha \left( 1_{H};1,1,0,0\right) }B(x_{1}x_{2}\otimes
gx_{1};X_{1}X_{2},gx_{1})1_{H}\otimes gx_{1}\otimes x_{1}x_{2}+ \\
&&\left( -1\right) ^{\alpha \left( 1_{H};1,0,1,0\right) }B(x_{1}x_{2}\otimes
gx_{1};X_{1}X_{2},gx_{1})X_{2}^{1-l_{2}}\otimes gx_{1}^{1-u_{1}}\otimes
g^{l_{1}+l_{2}+u_{1}}x_{1}^{1+1}x_{2}^{l_{2}}+=0 \\
&&\left( -1\right) ^{\alpha \left( 1_{H};1,1,1,0\right) }B(x_{1}x_{2}\otimes
gx_{1};X_{1}X_{2},gx_{1})X_{2}^{1-l_{2}}\otimes gx_{1}^{1-u_{1}}\otimes
g^{l_{1}+l_{2}+u_{1}}x_{1}^{1+1}x_{2}^{l_{2}}=0
\end{eqnarray*}

\subsubsection{Case $X_{1}X_{2}\otimes g\otimes gx_{1}$}

We have to consider only the second and the fourth summand of the left side
of the equality. Second summand gives us

\begin{eqnarray*}
l_{1} &=&u_{1}=0,l_{2}=u_{2}=0 \\
a &=&0,b_{1}=b_{2}=1
\end{eqnarray*}%
Since $\alpha \left( x_{1};0,0,0,0\right) \equiv a+b_{1}+b_{2}\equiv 0$, we
get
\begin{equation*}
B(gx_{2}\otimes gx_{1};X_{1}X_{2},g)X_{1}X_{2}\otimes g\otimes gx_{1}.
\end{equation*}

Fourth summand gives us%
\begin{eqnarray*}
l_{1}+u_{1} &=&1,l_{2}=u_{2}=0 \\
a &=&0,b_{1}-l_{1}=1\Rightarrow b_{1}=1,l_{1}=0,u_{1}=1,b_{2}=1 \\
d &=&1,e_{1}=u_{1}=1,e_{2}=0.
\end{eqnarray*}%
Since $\alpha \left( 1_{H};0,0,1,0\right) \equiv e_{2}+\left(
a+b_{1}+b_{2}\right) \equiv 0$ we get%
\begin{equation*}
B(x_{1}x_{2}\otimes gx_{1};X_{1}X_{2},gx_{1})X_{1}X_{2}\otimes g\otimes
gx_{1}
\end{equation*}%
By considering also the right side, we get%
\begin{equation*}
-B(x_{1}x_{2}\otimes g;X_{1}X_{2},g)+B(x_{1}x_{2}\otimes
gx_{1};X_{1}X_{2},gx_{1})+B(gx_{2}\otimes gx_{1};X_{1}X_{2},g)=0
\end{equation*}%
which holds in view of the form of the elements.

\subsubsection{Case $X_{1}\otimes gx_{1}\otimes gx_{2}$}

This case was already considered in subsection $B\left( x_{1}x_{2}\otimes
gx_{1};X_{1},gx_{1}x_{2}\right) .$

\subsubsection{Case $X_{1}\otimes g\otimes x_{1}x_{2}$}

This case was already considered in subsection $B\left( x_{1}x_{2}\otimes
gx_{1};X_{1},gx_{1}x_{2}\right) $

\subsubsection{Case $X_{2}\otimes gx_{1}\otimes gx_{1}$}

We have to consider only the second and the fourth summand of the left side
of the equality. Second summand gives us

\begin{eqnarray*}
l_{1} &=&u_{1}=0,l_{2}=u_{2}=0 \\
a &=&b_{1}=0,b_{2}=1, \\
d &=&e_{1}=1,e_{2}=0
\end{eqnarray*}

Since $\alpha \left( x_{1};0,0,0,0\right) \equiv a+b_{1}+b_{2}\equiv 1$, we
get%
\begin{equation*}
-B(gx_{2}\otimes gx_{1};X_{2},gx_{1})X_{2}\otimes gx_{1}\otimes gx_{1}
\end{equation*}

The fourth summand of the left side gives us%
\begin{eqnarray*}
l_{1}+u_{1} &=&1,l_{2}=u_{2}=0 \\
a &=&0,b_{2}=1,b_{1}=l_{1} \\
d &=&1,e_{1}-u_{1}=1\Rightarrow e_{1}=1,u_{1}=0,b_{1}=l_{1}=1 \\
e_{2} &=&0
\end{eqnarray*}%
Since $\alpha \left( 1_{H};1,0,0,0\right) \equiv b_{2}=1$, we get%
\begin{equation*}
-B(x_{1}x_{2}\otimes gx_{1};X_{1}X_{2},gx_{1})X_{2}\otimes gx_{1}\otimes
gx_{1}
\end{equation*}%
By considering also the right side, we get%
\begin{equation*}
-B(x_{1}x_{2}\otimes g;X_{2},gx_{1})-B(x_{1}x_{2}\otimes
gx_{1};X_{1}X_{2},gx_{1})-B(gx_{2}\otimes gx_{1};X_{2},gx_{1})=0
\end{equation*}%
which holds in view of the form of the elements.

\subsubsection{Case $1_{H}\otimes gx_{1}\otimes x_{1}x_{2}$}

This case was already considered in subsection $B\left( x_{1}x_{2}\otimes
x_{1}x_{2};X_{1},gx_{1}x_{2}\right) .$

\subsection{$B\left( x_{1}x_{2}\otimes gx_{1};GX_{1},gx_{1}\right) $}

We deduce that%
\begin{equation*}
a=b_{1}=1,b_{2}=0,d=e_{1}=1,e_{2}=0
\end{equation*}%
and we get%
\begin{eqnarray*}
&&\sum_{l_{1}=0}^{b_{1}}\sum_{u_{1}=0}^{e_{1}}\left( -1\right) ^{\alpha
\left( 1_{H};l_{1},0,u_{1},0\right) }B(x_{1}x_{2}\otimes
gx_{1};GX_{1},gx_{1})GX_{1}^{1-l_{1}}\otimes gx_{1}^{1-u_{1}}\otimes
g^{l_{1}+u_{1}}x_{1}^{l_{1}+u_{1}} \\
&=&\left( -1\right) ^{\alpha \left( 1_{H};0,0,0,0\right)
}B(x_{1}x_{2}\otimes gx_{1};GX_{1},gx_{1})GX_{1}\otimes gx_{1}\otimes 1_{H}+
\\
&&\left( -1\right) ^{\alpha \left( 1_{H};0,0,1,0\right) }B(x_{1}x_{2}\otimes
gx_{1};GX_{1},gx_{1})GX_{1}\otimes g\otimes gx_{1}+ \\
&&\left( -1\right) ^{\alpha \left( 1_{H};1,0,0,0\right) }B(x_{1}x_{2}\otimes
gx_{1};GX_{1},gx_{1})G\otimes gx_{1}\otimes gx_{1}+ \\
&&\left( -1\right) ^{\alpha \left( 1_{H};1,0,1,0\right) }B(x_{1}x_{2}\otimes
gx_{1};GX_{1},gx_{1})G\otimes gx_{1}^{1-u_{1}}\otimes
g^{1+u_{1}}x_{1}^{1+1}=0
\end{eqnarray*}

\subsubsection{Case $GX_{1}\otimes g\otimes gx_{1}$}

We have to consider only the second and the fourth summand of the left side
of the equality. Second summand gives us

\begin{eqnarray*}
l_{1} &=&u_{1}=0,l_{2}=u_{2}=0 \\
a &=&b_{1}=1,b_{2}=0, \\
d &=&1,e_{1}=e_{2}=0
\end{eqnarray*}%
Since $\alpha \left( x_{1};0,0,0,0\right) \equiv a+b_{1}+b_{2}\equiv 0$, we
get

\begin{equation*}
B(gx_{2}\otimes gx_{1};GX_{1},g)GX_{1}\otimes g\otimes gx_{1}.
\end{equation*}

The fourth summand of the left side gives us%
\begin{eqnarray*}
l_{1}+u_{1} &=&1,l_{2}=u_{2}=0 \\
a &=&1,b_{2}=0,b_{1}-l_{1}=1\Rightarrow b_{1}=1,l_{1}=0,u_{1}=1 \\
d &=&1,e_{1}=u_{1}=1,e_{2}=0
\end{eqnarray*}%
Since $\alpha \left( 1_{H};0,0,1,0\right) \equiv e_{2}+\left(
a+b_{1}+b_{2}\right) \equiv 0$, we get%
\begin{equation*}
B(x_{1}x_{2}\otimes gx_{1};GX_{1},gx_{1})GX_{1}\otimes g\otimes gx_{1}
\end{equation*}%
By considering also the right side, we obtain%
\begin{equation*}
-B(x_{1}x_{2}\otimes g;GX_{1},g)+B(x_{1}x_{2}\otimes
gx_{1};GX_{1},gx_{1})+B(gx_{2}\otimes gx_{1};GX_{1},g)=0
\end{equation*}%
which holds in view of the form of the elements.

\subsubsection{Case $G\otimes gx_{1}\otimes gx_{1}$}

We have to consider only the second and the fourth summand of the left side
of the equality. Second summand gives us

\begin{eqnarray*}
l_{1} &=&u_{1}=0,l_{2}=u_{2}=0 \\
a &=&1,b_{1}=b_{2}=0, \\
d &=&e_{1}=1,e_{2}=0
\end{eqnarray*}%
Since $\alpha \left( x_{1};0,0,0,0\right) \equiv a+b_{1}+b_{2}\equiv 1$, we
get

\begin{equation*}
-B(gx_{2}\otimes gx_{1};G,gx_{1})G\otimes gx_{1}\otimes gx_{1}.
\end{equation*}

The fourth summand of the left side gives us%
\begin{eqnarray*}
l_{1}+u_{1} &=&1,l_{2}=u_{2}=0 \\
a &=&1,b_{2}=0,b_{1}=l_{1} \\
d &=&1,e_{1}-u_{1}=1\Rightarrow e_{1}=1,u_{1}=0,b_{1}=l_{1}=1,e_{2}=0
\end{eqnarray*}%
Since $\alpha \left( 1_{H};1,0,0,0\right) \equiv b_{2}=0$, we get
\begin{equation*}
B(x_{1}x_{2}\otimes gx_{1};GX_{1},gx_{1})G\otimes gx_{1}\otimes gx_{1}
\end{equation*}%
By considering also the right side, we obtain%
\begin{equation*}
-B(x_{1}x_{2}\otimes g;G,gx_{1})+B(x_{1}x_{2}\otimes
gx_{1};GX_{1},gx_{1})-B(gx_{2}\otimes gx_{1};G,gx_{1})=0
\end{equation*}%
which holds in view of the form of the elements.

\subsection{$B\left( x_{1}x_{2}\otimes gx_{1};GX_{2},gx_{1}\right) $}

We deduce that%
\begin{equation*}
a=b_{2}=1,b_{1}=0,d=e_{1}=1,e_{2}=0
\end{equation*}%
and we get\

\begin{eqnarray*}
&&\sum_{l_{2}=0}^{1}\sum_{u_{1}=0}^{1}\left( -1\right) ^{\alpha \left(
1_{H};0,l_{2},u_{1},0\right) }B(x_{1}x_{2}\otimes
gx_{1};GX_{2},gx_{1})GX_{2}^{1-l_{2}}\otimes gx_{1}^{1-u_{1}}\otimes
g^{l_{2}+u_{1}}x_{1}^{u_{1}}x_{2}^{l_{2}} \\
&=&\left( -1\right) ^{\alpha \left( 1_{H};0,0,0,0\right)
}B(x_{1}x_{2}\otimes gx_{1};GX_{2},gx_{1})GX_{2}\otimes gx_{1}\otimes 1_{H}+
\\
&&\left( -1\right) ^{\alpha \left( 1_{H};0,0,1,0\right) }B(x_{1}x_{2}\otimes
gx_{1};GX_{2},gx_{1})GX_{2}\otimes g\otimes gx_{1}+ \\
&&\left( -1\right) ^{\alpha \left( 1_{H};0,1,0,0\right) }B(x_{1}x_{2}\otimes
gx_{1};GX_{2},gx_{1})G\otimes gx_{1}\otimes gx_{2}+ \\
&&\left( -1\right) ^{\alpha \left( 1_{H};0,1,1,0\right) }B(x_{1}x_{2}\otimes
gx_{1};GX_{2},gx_{1})
\end{eqnarray*}

\subsubsection{Case $GX_{2}\otimes g\otimes gx_{1}$}

We have to consider only the second and the fourth summand of the left side
of the equality. Second summand gives us

\begin{eqnarray*}
l_{1} &=&u_{1}=0,l_{2}=u_{2}=0 \\
a &=&b_{2}=1,b_{1}=0, \\
d &=&1,e_{1}=e_{2}=0
\end{eqnarray*}

Since $\alpha \left( x_{1};0,0,0,0\right) \equiv a+b_{1}+b_{2}\equiv 0,$ we
get%
\begin{equation*}
B(gx_{2}\otimes gx_{1};GX_{2},g)GX_{2}\otimes g\otimes gx_{1}.
\end{equation*}

The fourth summand of the left side gives us%
\begin{eqnarray*}
l_{1}+u_{1} &=&1,l_{2}=u_{2}=0 \\
a &=&b_{2}=1,0,b_{1}=l_{1} \\
d &=&1,e_{1}=u_{1},e_{2}=0
\end{eqnarray*}%
Since $\alpha \left( 1_{H};0,0,1,0\right) \equiv e_{2}+\left(
a+b_{1}+b_{2}\right) \equiv 0$ and $\alpha \left( 1_{H};1,0,0,0\right)
\equiv b_{2}=1,$ we obtain
\begin{equation*}
\left[ +B(x_{1}x_{2}\otimes gx_{1};GX_{2},gx_{1})-B(x_{1}x_{2}\otimes
gx_{1};GX_{1}X_{2},g)\right] GX_{2}\otimes g\otimes gx_{1}.
\end{equation*}%
By considering also the right side, we get%
\begin{gather*}
-B(x_{1}x_{2}\otimes g;GX_{2},g)+B(x_{1}x_{2}\otimes gx_{1};GX_{2},gx_{1}) \\
-B(x_{1}x_{2}\otimes gx_{1};GX_{1}X_{2},g)+B(gx_{2}\otimes gx_{1};GX_{2},g)=0
\end{gather*}%
which holds in view of the form of the elements.

\subsubsection{Case $G\otimes gx_{1}\otimes gx_{2}$}

This case was already considered in subsection $B\left( x_{1}x_{2}\otimes
gx_{1};G,gx_{1}x_{2}\right) .$

\subsubsection{Case $G\otimes g\otimes x_{1}x_{2}$}

This case was already considered in subsection $B\left( x_{1}x_{2}\otimes
gx_{1};G,gx_{1}x_{2}\right) .$

\subsection{$B\left( x_{1}x_{2}\otimes gx_{1};GX_{1}X_{2},gx_{1}x_{2}\right)
$}

We deduce that%
\begin{equation*}
a=b_{1}=b_{2}=1,d=e_{1}=e_{2}=1
\end{equation*}%
and we get%
\begin{gather*}
\left( -1\right) ^{\alpha \left( 1_{H};0,0,0,0\right) }B(x_{1}x_{2}\otimes
gx_{1};GX_{1}X_{2},gx_{1}x_{2})GX_{1}X_{2}\otimes gx_{1}x_{2}\otimes 1_{H} \\
+\left( -1\right) ^{\alpha \left( 1_{H};1,0,0,0\right) }B(x_{1}x_{2}\otimes
gx_{1};GX_{1}X_{2},gx_{1}x_{2})GX_{2}\otimes gx_{1}x_{2}\otimes gx_{1} \\
+\left( -1\right) ^{\alpha \left( 1_{H};0,1,0,0\right) }B(x_{1}x_{2}\otimes
gx_{1};GX_{1}X_{2},gx_{1}x_{2})GX_{1}\otimes gx_{1}x_{2}\otimes gx_{2}\text{
} \\
+\left( -1\right) ^{\alpha \left( 1_{H};1,1,0,0\right) }B(x_{1}x_{2}\otimes
gx_{1};GX_{1}X_{2},gx_{1}x_{2})G\otimes gx_{1}x_{2}\otimes x_{1}x_{2} \\
+\left( -1\right) ^{\alpha \left( 1_{H};0,0,1,0\right) }B(x_{1}x_{2}\otimes
gx_{1};GX_{1}X_{2},gx_{1}x_{2})GX_{1}X_{2}\otimes gx_{2}\otimes gx_{1} \\
+\left( -1\right) ^{\alpha \left( 1_{H};1,0,1,0\right) }B(x_{1}x_{2}\otimes
gx_{1};GX_{1}X_{2},gx_{1}x_{2})GX_{1}^{1-l_{1}}X_{2}\otimes gx_{2}\otimes
g^{l_{1}+1}x_{1}^{1+1}=0 \\
+\left( -1\right) ^{\alpha \left( 1_{H};0,1,1,0\right) }B(x_{1}x_{2}\otimes
gx_{1};GX_{1}X_{2},gx_{1}x_{2})GX_{1}\otimes gx_{2}\otimes x_{1}x_{2} \\
+\left( -1\right) ^{\alpha \left( 1_{H};1,1,1,0\right) }B(x_{1}x_{2}\otimes
gx_{1};GX_{1}X_{2},gx_{1}x_{2})GX_{1}^{1-l_{1}}\otimes gx_{2}\otimes
g^{l_{1}}x_{1}^{1+1}x_{2}=0 \\
+\left( -1\right) ^{\alpha \left( 1_{H};0,0,0,1\right) }B(x_{1}x_{2}\otimes
gx_{1};GX_{1}X_{2},gx_{1}x_{2})GX_{1}X_{2}\otimes gx_{1}\otimes gx_{2} \\
+\left( -1\right) ^{\alpha \left( 1_{H};1,0,0,1\right) }B(x_{1}x_{2}\otimes
gx_{1};GX_{1}X_{2},gx_{1}x_{2})GX_{2}\otimes gx_{1}\otimes x_{1}x_{2} \\
+\left( -1\right) ^{\alpha \left( 1_{H};l_{1},1,0,1\right)
}B(x_{1}x_{2}\otimes
gx_{1};GX_{1}X_{2},gx_{1}x_{2})GX_{1}^{1-l_{1}}X_{2}^{1-l_{2}}\otimes
gx_{1}\otimes g^{l_{1}+l_{2}+1}x_{1}^{l_{1}}x_{2}^{1+1}=0 \\
+\left( -1\right) ^{\alpha \left( 1_{H};0,0,1,1\right) }B(x_{1}x_{2}\otimes
gx_{1};GX_{1}X_{2},gx_{1}x_{2})GX_{1}X_{2}\otimes g\otimes x_{1}x_{2} \\
+\left( -1\right) ^{\alpha \left( 1_{H};1,0,1,1\right) }B(x_{1}x_{2}\otimes
gx_{1};GX_{1}X_{2},gx_{1}x_{2})GX_{1}^{1-l_{1}}X_{2}\otimes g\otimes
g^{l_{1}+}x_{1}^{1+1}x_{2}=0 \\
+\left( -1\right) ^{\alpha \left( 1_{H};1,1,1,1\right) }B(x_{1}x_{2}\otimes
gx_{1};GX_{1}X_{2},gx_{1}x_{2})GX_{1}^{1-l_{1}}X_{2}^{1-l_{2}}\otimes
g\otimes g^{l_{1}+l_{2}+}x_{1}^{l_{1}+1}x_{2}^{1+1}=0
\end{gather*}

\subsubsection{Case $GX_{2}\otimes gx_{1}x_{2}\otimes gx_{1}$}

We have to consider only the second and the fourth summand of the left side
of the equality. Second summand gives us

\begin{eqnarray*}
l_{1} &=&u_{1}=0,l_{2}=u_{2}=0 \\
a &=&b_{2}=1,b_{1}=0 \\
d &=&e_{1}=e_{2}=1.
\end{eqnarray*}%
Since $\alpha \left( x_{1};0,0,0,0\right) \equiv a+b_{1}+b_{2}\equiv 0,$ we
get%
\begin{equation*}
+B(gx_{2}\otimes gx_{1};GX_{2},gx_{1}x_{2})GX_{2}\otimes gx_{1}x_{2}\otimes
gx_{1}
\end{equation*}

Fourth summand gives us%
\begin{eqnarray*}
l_{1}+u_{1} &=&1,l_{2}=u_{2}=0 \\
a &=&b_{2}=1,b_{1}=l_{1} \\
d &=&e_{2}=1,e_{1}-u_{1}=1\Rightarrow e_{1}=1,u_{1}=0,b_{1}=l_{1}=1.
\end{eqnarray*}%
Since $\alpha \left( 1_{H};1,0,0,0\right) =b_{2}=1$ we get
\begin{equation*}
-B(x_{1}x_{2}\otimes gx_{1};GX_{1}X_{2},gx_{1}x_{2})GX_{2}\otimes
gx_{1}x_{2}\otimes gx_{1}
\end{equation*}%
By considering also the right side, we obtain%
\begin{equation*}
-B(x_{1}x_{2}\otimes g;GX_{2},gx_{1}x_{2})-B(x_{1}x_{2}\otimes
gx_{1};GX_{1}X_{2},gx_{1}x_{2})+B(gx_{2}\otimes gx_{1};GX_{2},gx_{1}x_{2})=0
\end{equation*}%
which holds in view of the form of the elements.

\subsubsection{Case $GX_{1}\otimes gx_{1}x_{2}\otimes gx_{2}$}

We have only to consider the third and the forth summands of the left side
of the equality. Third summand gives us%
\begin{eqnarray*}
l_{1} &=&u_{1}=l_{2}=u_{2}=0, \\
a &=&b_{1}=1,b_{2}=0, \\
d &=&e_{1}=e_{2}=1,
\end{eqnarray*}%
Since $\alpha \left( x_{2};0,0,0,0\right) \equiv a+b_{1}+b_{2}=0$ we get
\begin{equation*}
-B(gx_{1}\otimes gx_{1};GX_{1},gx_{1}x_{2})GX_{1}\otimes gx_{1}x_{2}\otimes
gx_{2}.
\end{equation*}%
Fourth summand gives us%
\begin{eqnarray*}
l_{1} &=&u_{1}=0,l_{2}+u_{2}=1 \\
a &=&b_{1}=1,b_{2}=l_{2} \\
d &=&e_{1}=1,e_{2}-u_{2}=1\Rightarrow e_{2}=1,u_{2}=0,b_{2}=l_{2}=1.
\end{eqnarray*}%
Since $\alpha \left( 1_{H};0,1,0,0\right) =0,$ we get
\begin{equation*}
B(x_{1}x_{2}\otimes gx_{1};GX_{1}X_{2},gx_{1}x_{2})GX_{1}\otimes
gx_{1}x_{2}\otimes gx_{2}
\end{equation*}%
Since there is no term to consider in the right side of the equality we
obtain%
\begin{equation*}
-B(gx_{1}\otimes gx_{1};GX_{1},gx_{1}x_{2})+B(x_{1}x_{2}\otimes
gx_{1};GX_{1}X_{2},gx_{1}x_{2})=0
\end{equation*}%
which holds in view of the form of the elements.

\subsubsection{Case $G\otimes gx_{1}x_{2}\otimes x_{1}x_{2}$}

First summand of the left side of the equality gives us%
\begin{eqnarray*}
l_{1} &=&u_{1}=0,l_{2}=u_{2}=0 \\
a &=&1,b_{1}=b_{2}=0, \\
d &=&e_{1}=e_{2}=1.
\end{eqnarray*}

Since $\alpha \left( x_{1}x_{2};0,0,0,0\right) \equiv 0,$ we get
\begin{equation*}
B(1_{H}\otimes gx_{1};G,gx_{1}x_{2})G\otimes gx_{1}x_{2}\otimes x_{1}x_{2}
\end{equation*}%
Second summand of the left side gives us%
\begin{eqnarray*}
l_{1} &=&u_{1}=0,l_{2}+u_{2}=1 \\
a &=&1,b_{1}=0,b_{2}=l_{2}, \\
d &=&e_{1}=1,e_{2}-u_{2}=1\Rightarrow e_{2}=1,u_{2}=0,b_{2}=l_{2}=1.
\end{eqnarray*}%
Since $\alpha \left( x_{1};0,1,0,0\right) \equiv a+b_{1}+b_{2}+1\equiv 1$ $,$
we obtain
\begin{equation*}
-B(gx_{2}\otimes gx_{1};GX_{2},gx_{1}x_{2})G\otimes gx_{1}x_{2}\otimes
x_{1}x_{2}
\end{equation*}%
Third summand of the left side gives us%
\begin{eqnarray*}
l_{1}+u_{1} &=&1,l_{2}=u_{2}=0, \\
a &=&1,b_{2}=0,b_{1}=l_{1} \\
d &=&1,e_{2}=1. \\
e_{1}-u_{1} &=&1\Rightarrow e_{1}=1,u_{1}=0,b_{1}=l_{1}=1
\end{eqnarray*}%
Since $\alpha \left( x_{2};1,0,0,0\right) =a+b_{1}\equiv 0$ we get%
\begin{equation*}
-B(gx_{1}\otimes gx_{1};GX_{1},gx_{1}x_{2})G\otimes gx_{1}x_{2}\otimes
x_{1}x_{2}
\end{equation*}%
The fourth summand of the left side gives us%
\begin{eqnarray*}
l_{1}+u_{1} &=&1,l_{2}+u_{2}=1 \\
a &=&1,b_{1}=l_{1},b_{2}=l_{2} \\
d &=&1,e_{1}-u_{1}=1\Rightarrow e_{1}=1,u_{1}=0,b_{1}=l_{1}=1, \\
e_{2}-u_{2} &=&1\Rightarrow e_{2}=1,u_{2}=0,b_{2}=l_{2}=1.
\end{eqnarray*}%
Since $\alpha \left( 1_{H};1,1,0,0\right) \equiv 1+b_{2}\equiv 0,$ we get%
\begin{equation*}
B(x_{1}x_{2}\otimes gx_{1};GX_{1}X_{2},gx_{1}x_{2})G\otimes
gx_{1}x_{2}\otimes x_{1}x_{2}
\end{equation*}%
Since there is no term to consider in the right side of the equality we
obtain

\begin{eqnarray*}
&&-B(gx_{1}\otimes gx_{1};GX_{1},gx_{1}x_{2})-B(gx_{2}\otimes
gx_{1};GX_{2},gx_{1}x_{2}) \\
&&+B(x_{1}x_{2}\otimes gx_{1};GX_{1}X_{2},gx_{1}x_{2})+B(1_{H}\otimes
gx_{1};G,gx_{1}x_{2}) \\
&=&0
\end{eqnarray*}%
which holds in view of the form of the elements.

\subsubsection{Case $GX_{1}X_{2}\otimes gx_{2}\otimes gx_{1}$}

We have to consider only the second and the fourth summand of the left side
of the equality. Second summand gives us

\begin{eqnarray*}
l_{1} &=&u_{1}=0,l_{2}=u_{2}=0 \\
a &=&b_{1}=b_{2}=1, \\
d &=&e_{2}=1,e_{1}=0
\end{eqnarray*}

Since $\alpha \left( x_{1};0,0,0,0\right) \equiv a+b_{1}+b_{2}\equiv 1,$ we
get%
\begin{equation*}
-B(gx_{2}\otimes gx_{1};GX_{1}X_{2},gx_{2})GX_{1}X_{2}\otimes gx_{2}\otimes
gx_{1}.
\end{equation*}

The fourth summand of the left side gives us%
\begin{eqnarray*}
l_{1}+u_{1} &=&1,l_{2}=u_{2}=0 \\
a &=&b_{2}=1,0,b_{1}-l_{1}=1\Rightarrow b_{1}=1,l_{1}=0,u_{1}=1 \\
d &=&e_{2}=1,e_{1}=u_{1}=1.
\end{eqnarray*}%
Since $\alpha \left( 1_{H};0,0,1,0\right) \equiv e_{2}+\left(
a+b_{1}+b_{2}\right) \equiv 0,$ we obtain
\begin{equation*}
B(x_{1}x_{2}\otimes gx_{1};GX_{1}X_{2},gx_{1}x_{2})GX_{1}X_{2}\otimes
gx_{2}\otimes gx_{1}.
\end{equation*}%
By considering also the right side, we get%
\begin{gather*}
B(x_{1}x_{2}\otimes gx_{1};GX_{1}X_{2},gx_{1}x_{2}) \\
-B(x_{1}x_{2}\otimes g;GX_{1}X_{2},gx_{2})-B(gx_{2}\otimes
gx_{1};GX_{1}X_{2},gx_{2})=0
\end{gather*}%
which holds in view of the form of the elements.

\subsubsection{Case $GX_{1}\otimes gx_{2}\otimes x_{1}x_{2}$}

First summand of the left side of the equality gives us%
\begin{eqnarray*}
l_{1} &=&u_{1}=0,l_{2}=u_{2}=0 \\
a &=&b_{1}=1,b_{2}=0, \\
d &=&e_{2}=1,e_{1}=0.
\end{eqnarray*}

Since $\alpha \left( x_{1}x_{2};0,0,0,0\right) \equiv 0,$ we get
\begin{equation*}
B(1_{H}\otimes gx_{1};GX_{1},gx_{2})GX_{1}\otimes gx_{2}\otimes x_{1}x_{2}
\end{equation*}%
Second summand of the left side gives us%
\begin{eqnarray*}
l_{1} &=&u_{1}=0,l_{2}+u_{2}=1 \\
a &=&b_{1}=1,0,b_{2}=l_{2}, \\
d &=&1,e_{2}-u_{2}=1\Rightarrow e_{2}=1,u_{2}=0,b_{2}=l_{2}=1,e_{1}=0.
\end{eqnarray*}%
Since $\alpha \left( x_{1};0,1,0,0\right) \equiv a+b_{1}+b_{2}+1\equiv 0$ $,$
we obtain
\begin{equation*}
B(gx_{2}\otimes gx_{1};GX_{1}X_{2},gx_{2})GX_{1}\otimes gx_{2}\otimes
x_{1}x_{2}
\end{equation*}%
Third summand of the left side gives us%
\begin{eqnarray*}
l_{1}+u_{1} &=&1,l_{2}=u_{2}=0, \\
a &=&1,b_{2}=0,b_{1}-l_{1}=1\Rightarrow b_{1}=1,l_{1}=0,u_{1}=1 \\
d &=&1,e_{2}=1,e_{1}=u_{1}=1.
\end{eqnarray*}%
Since $\alpha \left( x_{2};0,0,1,0\right) \equiv e_{2}\equiv 1$ we get%
\begin{equation*}
+B(gx_{1}\otimes gx_{1};GX_{1},gx_{1}x_{2})GX_{1}\otimes gx_{2}\otimes
x_{1}x_{2}
\end{equation*}%
The fourth summand of the left side gives us%
\begin{eqnarray*}
l_{1}+u_{1} &=&1,l_{2}+u_{2}=1 \\
a &=&1,b_{1}-l_{1}=1\Rightarrow b_{1}=1,l_{1}=0,u_{1}=1,b_{2}=l_{2} \\
d &=&1,e_{1}=u_{1}=1, \\
e_{2}-u_{2} &=&1\Rightarrow e_{2}=1,u_{2}=0,b_{2}=l_{2}=1.
\end{eqnarray*}%
Since $\alpha \left( 1_{H};0,1,1,0\right) \equiv e_{2}+a+b_{1}+b_{2}+1\equiv
1,$ we get%
\begin{equation*}
-B(x_{1}x_{2}\otimes gx_{1};GX_{1}X_{2},gx_{1}x_{2})GX_{1}\otimes
gx_{2}\otimes x_{1}x_{2}
\end{equation*}%
Since there is no term to consider in the right side of the equality we
obtain

\begin{eqnarray*}
&&+B(gx_{1}\otimes gx_{1};GX_{1},gx_{1}x_{2})+B(gx_{2}\otimes
gx_{1};GX_{1}X_{2},gx_{2}) \\
&&-B(x_{1}x_{2}\otimes gx_{1};GX_{1}X_{2},gx_{1}x_{2})+B(1_{H}\otimes
gx_{1};GX_{1},gx_{2}) \\
&=&0
\end{eqnarray*}%
which holds in view of the form of the elements.

\subsubsection{Case $GX_{1}X_{2}\otimes gx_{1}\otimes gx_{2}$}

We have only to consider the third and the forth summands of the left side
of the equality. Third summand gives us%
\begin{eqnarray*}
l_{1} &=&u_{1}=l_{2}=u_{2}=0, \\
a &=&b_{1}=b_{2}=1, \\
d &=&e_{1}=1,e_{2}=0
\end{eqnarray*}%
Since $\alpha \left( x_{2};0,0,0,0\right) \equiv a+b_{1}+b_{2}\equiv 1$ we
get
\begin{equation*}
B(gx_{1}\otimes gx_{1};GX_{1}X_{2},gx_{1})GX_{1}X_{2}\otimes gx_{1}\otimes
gx_{2}.
\end{equation*}%
Fourth summand gives us%
\begin{eqnarray*}
l_{1} &=&u_{1}=0,l_{2}+u_{2}=1 \\
a &=&b_{1}=1,b_{2}-l_{2}=1\Rightarrow b_{2}=1,l_{2}=0,u_{2}=1, \\
d &=&e_{1}=1,e_{2}=u_{2}=1.
\end{eqnarray*}%
Since $\alpha \left( 1_{H};0,0,0,1\right) \equiv a+b_{1}+b_{2}\equiv 1,$ we
get
\begin{equation*}
-B(x_{1}x_{2}\otimes gx_{1};GX_{1}X_{2},gx_{1}x_{2})GX_{1}X_{2}\otimes
gx_{1}\otimes gx_{2}.
\end{equation*}%
Since there is no term to consider in the right side of the equality we
obtain%
\begin{equation*}
B(gx_{1}\otimes gx_{1};GX_{1}X_{2},gx_{1})-B(x_{1}x_{2}\otimes
gx_{1};GX_{1}X_{2},gx_{1}x_{2})=0
\end{equation*}%
which holds in view of the form of the elements.

\subsubsection{Case $GX_{2}\otimes gx_{1}\otimes x_{1}x_{2}$}

First summand of the left side of the equality gives us%
\begin{eqnarray*}
l_{1} &=&u_{1}=0,l_{2}=u_{2}=0 \\
a &=&b_{2}=1,b_{1}=0, \\
d &=&e_{1}=1,e_{2}=0.
\end{eqnarray*}

Since $\alpha \left( x_{1}x_{2};0,0,0,0\right) \equiv 0,$ we get
\begin{equation*}
B(1_{H}\otimes gx_{1};GX_{2},gx_{1})GX_{2}\otimes gx_{1}\otimes x_{1}x_{2}
\end{equation*}%
Second summand of the left side gives us%
\begin{eqnarray*}
l_{1} &=&u_{1}=0,l_{2}+u_{2}=1 \\
a &=&1,b_{1}=0,b_{2}-l_{2}=1\Rightarrow b_{2}=1,l_{2}=0,u_{2}=1, \\
d &=&e_{1}=1,e_{2}=u_{2}=1.
\end{eqnarray*}%
Since $\alpha \left( x_{1};0,0,0,1\right) \equiv 1$ $,$ we obtain
\begin{equation*}
-B(gx_{2}\otimes gx_{1};GX_{2},gx_{1}x_{2})GX_{2}\otimes gx_{1}\otimes
x_{1}x_{2}
\end{equation*}%
Third summand of the left side gives us%
\begin{eqnarray*}
l_{1}+u_{1} &=&1,l_{2}=u_{2}=0, \\
a &=&b_{2}=1,b_{1}=l_{1} \\
d &=&1,e_{2}=0, \\
e_{1}-u_{1} &=&1\Rightarrow e_{1}=1,u_{1}=0,b_{1}=l_{1}=1
\end{eqnarray*}%
Since $\alpha \left( x_{2};1,0,0,0\right) =a+b_{1}\equiv 0$ we get%
\begin{equation*}
-B(gx_{1}\otimes gx_{1};GX_{1}X_{2},gx_{1})GX_{2}\otimes gx_{1}\otimes
x_{1}x_{2}
\end{equation*}%
The fourth summand of the left side gives us%
\begin{eqnarray*}
l_{1}+u_{1} &=&1,l_{2}+u_{2}=1 \\
a &=&1,b_{1}=l_{1},b_{2}-l_{2}=1\Rightarrow b_{2}=1,l_{2}=0,u_{2}=1 \\
d &=&1,e_{1}-u_{1}=1\Rightarrow e_{1}=1,u_{1}=0,b_{1}=l_{1}=1, \\
e_{2} &=&u_{2}=1.
\end{eqnarray*}%
Since $\alpha \left( 1_{H};1,0,0,1\right) \equiv a+b_{1}\equiv 0,$ we get%
\begin{equation*}
B(x_{1}x_{2}\otimes gx_{1};GX_{1}X_{2},gx_{1}x_{2})GX_{2}\otimes
gx_{1}\otimes x_{1}x_{2}
\end{equation*}%
Since there is no term to consider in the right side of the equality we
obtain

\begin{eqnarray*}
&&-B(gx_{1}\otimes gx_{1};GX_{1}X_{2},gx_{1})-B(gx_{2}\otimes
gx_{1};GX_{2},gx_{1}x_{2}) \\
&&+B(x_{1}x_{2}\otimes gx_{1};GX_{1}X_{2},gx_{1}x_{2})+B(1_{H}\otimes
gx_{1};GX_{2},gx_{1}) \\
&=&0
\end{eqnarray*}%
which holds in view of the form of the elements.

\subsubsection{Case $GX_{1}X_{2}\otimes g\otimes x_{1}x_{2}$}

First summand of the left side of the equality gives us%
\begin{eqnarray*}
l_{1} &=&u_{1}=0,l_{2}=u_{2}=0 \\
a &=&b_{1}=b_{2}=1, \\
d &=&1,e_{1}=e_{2}=0.
\end{eqnarray*}

Since $\alpha \left( x_{1}x_{2};0,0,0,0\right) \equiv 0,$ we get
\begin{equation*}
B(1_{H}\otimes gx_{1};GX_{1}X_{2},g)GX_{1}X_{2}\otimes g\otimes x_{1}x_{2}
\end{equation*}%
Second summand of the left side gives us%
\begin{eqnarray*}
l_{1} &=&u_{1}=0,l_{2}+u_{2}=1 \\
a &=&b_{1}=1,0,b_{2}-l_{2}=1\Rightarrow b_{2}=1,l_{2}=0,u_{2}=1, \\
d &=&1,e_{2}=u_{2}=1,e_{1}=0.
\end{eqnarray*}%
Since $\alpha \left( x_{1};0,0,0,1\right) \equiv 1$ $,$ we obtain
\begin{equation*}
-B(gx_{2}\otimes gx_{1};GX_{1}X_{2},gx_{2})GX_{1}X_{2}\otimes g\otimes
x_{1}x_{2}.
\end{equation*}%
Third summand of the left side gives us%
\begin{eqnarray*}
l_{1}+u_{1} &=&1,l_{2}=u_{2}=0, \\
a &=&b_{2}=1,b_{1}-l_{1}=1\Rightarrow b_{1}=1,l_{1}=0,u_{1}=1 \\
d &=&1,e_{2}=0,e_{1}=u_{1}=1.
\end{eqnarray*}%
Since $\alpha \left( x_{2};0,0,1,0\right) \equiv e_{2}\equiv 0$ we get%
\begin{equation*}
-B(gx_{1}\otimes gx_{1};GX_{1}X_{2},gx_{1})GX_{1}X_{2}\otimes g\otimes
x_{1}x_{2}.
\end{equation*}%
The fourth summand of the left side gives us%
\begin{eqnarray*}
l_{1}+u_{1} &=&1,l_{2}+u_{2}=1 \\
a &=&1,b_{1}-l_{1}=1\Rightarrow b_{1}=1,l_{1}=0,u_{1}=1, \\
b_{2}-l_{2} &=&1\Rightarrow b_{2}=1,l_{2}=0,u_{2}=1 \\
d &=&1,e_{1}=u_{1}=1, \\
e_{2} &=&u_{2}=1
\end{eqnarray*}%
Since $\alpha \left( 1_{H};0,0,1,1\right) \equiv 1+e_{2}\equiv 0,$ we get%
\begin{equation*}
B(x_{1}x_{2}\otimes gx_{1};GX_{1}X_{2},gx_{1}x_{2})GX_{1}X_{2}\otimes
g\otimes x_{1}x_{2}.
\end{equation*}%
Since there is no term to consider in the right side of the equality we
obtain

\begin{eqnarray*}
&&-B(gx_{1}\otimes gx_{1};GX_{1}X_{2},gx_{1})-B(gx_{2}\otimes
gx_{1};GX_{1}X_{2},gx_{2}) \\
&&+B(x_{1}x_{2}\otimes gx_{1};GX_{1}X_{2},gx_{1}x_{2})+B(1_{H}\otimes
gx_{1};GX_{1}X_{2},g) \\
&=&0
\end{eqnarray*}%
which holds in view of the form of the elements.

\section{$B(x_{1}x_{2}\otimes gx_{2})$}

By using$\left( \ref{simplgx}\right) ,$ we get%
\begin{eqnarray}
&&B(x_{1}x_{2}\otimes gx_{2})  \label{form x1x2otgx2} \\
&=&(1_{A}\otimes g)B(gx_{1}x_{2}\otimes 1_{H})(1_{A}\otimes gx_{2})  \notag
\\
&&+(1_{A}\otimes x_{2})B(gx_{1}x_{2}\otimes 1_{H})  \notag
\end{eqnarray}

so that we obtain
\begin{eqnarray*}
B\left( x_{1}x_{2}\otimes gx_{2}\right) &=&2B(gx_{1}x_{2}\otimes
1_{H};1_{A},g)1_{A}\otimes gx_{2}+ \\
&&+2B(gx_{1}x_{2}\otimes 1_{H};G,gx_{1})G\otimes gx_{1}x_{2}+ \\
&&-2B(x_{2}\otimes 1_{H};1_{A},gx_{1})X_{1}\otimes gx_{1}x_{2}+ \\
&&+\left[ 2B(x_{1}\otimes 1_{H};1_{A},gx_{1})+2B(gx_{1}x_{2}\otimes
1_{H};1_{A},gx_{1}x_{2})\right] X_{2}\otimes gx_{1}x_{2}+ \\
&&+\left[
\begin{array}{c}
2B(g\otimes 1_{H};1_{A},g)+2B(x_{2}\otimes \ 1_{H};1_{A},gx_{2}) \\
+2B(x_{1}\otimes 1_{H};1_{A},gx_{1})+2B(gx_{1}x_{2}\otimes
1_{H};1_{A},gx_{1}x_{2})%
\end{array}%
\right] X_{1}X_{2}\otimes gx_{2}+ \\
&&+\left[ 2B(x_{2}\otimes 1_{H};G,g)-2B(gx_{1}x_{2}\otimes 1_{H};G,gx_{1})%
\right] GX_{1}\otimes gx_{2}+ \\
&&+\left[ -2B(x_{1}\otimes 1_{H};G,g)-2B(gx_{1}x_{2}\otimes 1_{H};G,gx_{2})%
\right] GX_{2}\otimes gx_{2}+ \\
&&+\left[ 2B(g\otimes 1_{H};G,gx_{1})+2B(x_{2}\otimes \ 1_{H};G,gx_{1}x_{2})%
\right] GX_{1}X_{2}\otimes gx_{1}x_{2}.
\end{eqnarray*}%
We write the Casimir equality for $B\left( x_{1}x_{2}\otimes gx_{2}\right) .$%
\begin{eqnarray*}
&&\sum_{w_{1}=0}^{1}\sum_{w_{2}=0}^{1}\left( -1\right) ^{\left(
1+w_{2}\right)
w_{1}}\sum_{a,b_{1},b_{2},d,e_{1},e_{2}=0}^{1}\sum_{l_{1}=0}^{b_{1}}%
\sum_{l_{2}=0}^{b_{2}}\sum_{u_{1}=0}^{e_{1}}\sum_{u_{2}=0}^{e_{2}} \\
&&\left( -1\right) ^{\alpha \left(
x_{1}^{1-w_{1}}x_{2}^{1-w_{2}};l_{1},l_{2},u_{1},u_{2}\right)
}B(g^{w_{1}+w_{2}}x_{1}^{w_{1}}x_{2}^{w_{2}}\otimes
gx_{2};G^{a}X_{1}^{b_{1}}X_{2}^{b_{2}},g^{d}x_{1}^{e_{1}}x_{2}^{e_{2}}) \\
&&G^{a}X_{1}^{b_{1}-l_{1}}X_{2}^{b_{2}-l_{2}}\otimes
g^{d}x_{1}^{e_{1}-u_{1}}x_{2}^{e_{2}-u_{2}}\otimes \\
&&g^{a+b_{1}+b_{2}+l_{1}+l_{2}+d+e_{1}+e_{2}+u_{1}+u_{2}}x_{1}^{l_{1}+u_{1}+1-w_{1}}x_{2}^{l_{2}+u_{2}+1-w_{2}}
\\
&=&\sum_{\omega _{2}=0}^{1}B(x_{1}x_{2}\otimes gx_{2}^{1-\omega
_{2}})\otimes g^{\omega _{2}}x_{2}^{\omega _{2}}
\end{eqnarray*}%
i.e.%
\begin{eqnarray*}
&&\sum_{a,b_{1},b_{2},d,e_{1},e_{2}=0}^{1}\sum_{l_{1}=0}^{b_{1}}%
\sum_{l_{2}=0}^{b_{2}}\sum_{u_{1}=0}^{e_{1}}\sum_{u_{2}=0}^{e_{2}}\left(
-1\right) ^{\alpha \left( x_{1}x_{2};l_{1},l_{2},u_{1},u_{2}\right) } \\
&&B(1_{H}\otimes
gx_{2};G^{a}X_{1}^{b_{1}}X_{2}^{b_{2}},g^{d}x_{1}^{e_{1}}x_{2}^{e_{2}}) \\
&&G^{a}X_{1}^{b_{1}-l_{1}}X_{2}^{b_{2}-l_{2}}\otimes
g^{d}x_{1}^{e_{1}-u_{1}}x_{2}^{e_{2}-u_{2}}\otimes
g^{a+b_{1}+b_{2}+l_{1}+l_{2}+d+e_{1}+e_{2}+u_{1}+u_{2}}x_{1}^{l_{1}+u_{1}+1}x_{2}^{l_{2}+u_{2}+1}+
\\
&&+\sum_{a,b_{1},b_{2},d,e_{1},e_{2}=0}^{1}\sum_{l_{1}=0}^{b_{1}}%
\sum_{l_{2}=0}^{b_{2}}\sum_{u_{1}=0}^{e_{1}}\sum_{u_{2}=0}^{e_{2}}\left(
-1\right) ^{\alpha \left( x_{1};l_{1},l_{2},u_{1},u_{2}\right) } \\
&&B(gx_{2}\otimes
gx_{2};G^{a}X_{1}^{b_{1}}X_{2}^{b_{2}},g^{d}x_{1}^{e_{1}}x_{2}^{e_{2}}) \\
&&G^{a}X_{1}^{b_{1}-l_{1}}X_{2}^{b_{2}-l_{2}}\otimes
g^{d}x_{1}^{e_{1}-u_{1}}x_{2}^{e_{2}-u_{2}}\otimes
g^{a+b_{1}+b_{2}+l_{1}+l_{2}+d+e_{1}+e_{2}+u_{1}+u_{2}}x_{1}^{l_{1}+u_{1}+1}x_{2}^{l_{2}+u_{2}}+
\\
&&+\left( -1\right)
\sum_{a,b_{1},b_{2},d,e_{1},e_{2}=0}^{1}\sum_{l_{1}=0}^{b_{1}}%
\sum_{l_{2}=0}^{b_{2}}\sum_{u_{1}=0}^{e_{1}}\sum_{u_{2}=0}^{e_{2}}\left(
-1\right) ^{\alpha \left( x_{2};l_{1},l_{2},u_{1},u_{2}\right) } \\
&&B(gx_{1}\otimes
gx_{2};G^{a}X_{1}^{b_{1}}X_{2}^{b_{2}},g^{d}x_{1}^{e_{1}}x_{2}^{e_{2}}) \\
&&G^{a}X_{1}^{b_{1}-l_{1}}X_{2}^{b_{2}-l_{2}}\otimes
g^{d}x_{1}^{e_{1}-u_{1}}x_{2}^{e_{2}-u_{2}}\otimes
g^{a+b_{1}+b_{2}+l_{1}+l_{2}+d+e_{1}+e_{2}+u_{1}+u_{2}}x_{1}^{l_{1}+u_{1}}x_{2}^{l_{2}+u_{2}+1}+
\\
&&+\sum_{a,b_{1},b_{2},d,e_{1},e_{2}=0}^{1}\sum_{l_{1}=0}^{b_{1}}%
\sum_{l_{2}=0}^{b_{2}}\sum_{u_{1}=0}^{e_{1}}\sum_{u_{2}=0}^{e_{2}}\left(
-1\right) ^{\alpha \left( 1_{H};l_{1},l_{2},u_{1},u_{2}\right) } \\
&&B(x_{1}x_{2}\otimes
gx_{2};G^{a}X_{1}^{b_{1}}X_{2}^{b_{2}},g^{d}x_{1}^{e_{1}}x_{2}^{e_{2}}) \\
&&G^{a}X_{1}^{b_{1}-l_{1}}X_{2}^{b_{2}-l_{2}}\otimes
g^{d}x_{1}^{e_{1}-u_{1}}x_{2}^{e_{2}-u_{2}}\otimes
g^{a+b_{1}+b_{2}+l_{1}+l_{2}+d+e_{1}+e_{2}+u_{1}+u_{2}}x_{1}^{l_{1}+u_{1}}x_{2}^{l_{2}+u_{2}}+
\\
&=&B(x_{1}x_{2}\otimes gx_{2})\otimes 1_{H}+B(x_{1}x_{2}\otimes g)\otimes
gx_{2}.
\end{eqnarray*}

\subsection{$B\left( x_{1}x_{2}\otimes gx_{2};1_{A},gx_{2}\right) $}

We deduce that%
\begin{eqnarray*}
a &=&b_{1}=b_{2}=0 \\
d &=&e_{2}=1,e_{1}=0
\end{eqnarray*}%
and we get%
\begin{eqnarray*}
&&\left( -1\right) ^{\alpha \left( 1_{H};0,0,0,0\right) }B(x_{1}x_{2}\otimes
gx_{2};1_{A},gx_{2})1_{A}\otimes gx_{2}\otimes 1_{H} \\
&&+\left( -1\right) ^{\alpha \left( 1_{H};0,0,0,1\right)
}B(x_{1}x_{2}\otimes gx_{2};1_{A},gx_{2})1_{A}\otimes g\otimes gx_{2}.
\end{eqnarray*}

\subsubsection{Case $1_{A}\otimes g\otimes gx_{2}$}

We have only to consider the third and the fourth summands of the left side
of the equality. Third summand gives us%
\begin{eqnarray*}
l_{1} &=&u_{1}=l_{2}=u_{2}=0, \\
a &=&b_{1}=b_{2}=0, \\
d &=&1,e_{1}=e_{2}=0
\end{eqnarray*}%
Since $\alpha \left( x_{2};0,0,0,0\right) \equiv a+b_{1}+b_{2}\equiv 0$ we
get
\begin{equation*}
-B(gx_{1}\otimes gx_{2};1_{A},g)1_{A}\otimes g\otimes gx_{2}.
\end{equation*}%
Fourth summand gives us%
\begin{eqnarray*}
l_{1} &=&u_{1}=0,l_{2}+u_{2}=1 \\
a &=&b_{1}=0,b_{2}=l_{2}, \\
d &=&1,e_{1}=0,e_{2}=u_{2}.
\end{eqnarray*}%
Since $\alpha \left( 1_{H};0,0,0,1\right) \equiv a+b_{1}+b_{2}\equiv 0$ and $%
\alpha \left( 1_{H};0,1,0,0\right) \equiv 0,$ we get
\begin{equation*}
\left[ B(x_{1}x_{2}\otimes gx_{2};1_{A},gx_{2})+B(x_{1}x_{2}\otimes
gx_{2};X_{2},g)\right] 1_{A}\otimes g\otimes gx_{2}.
\end{equation*}%
By considering also the right side of the equality we obtain%
\begin{equation*}
-B(x_{1}x_{2}\otimes g;1_{A},g)-B(gx_{1}\otimes
gx_{2};1_{A},g)+B(x_{1}x_{2}\otimes gx_{2};1_{A},gx_{2})+B(x_{1}x_{2}\otimes
gx_{2};X_{2},g)=0
\end{equation*}%
which holds in view of the form of the elements.

\subsection{$B\left( x_{1}x_{2}\otimes gx_{2};G,gx_{1}x_{2}\right) $}

We deduce that%
\begin{equation*}
a=1,b_{1}=b_{2}=0,d=e_{1}=e_{2}=1
\end{equation*}%
and we get%
\begin{eqnarray*}
&&+\sum_{u_{1}=0}^{1}\sum_{u_{2}=0}^{1}\left( -1\right) ^{\alpha \left(
1_{H};0,0,u_{1},u_{2}\right) }B(x_{1}\otimes
x_{1}x_{2};G,gx_{1}x_{2})G\otimes gx_{1}^{1-u_{1}}x_{2}^{1-u_{2}}\otimes
g^{u_{1}+u_{2}}x_{1}^{u_{1}}x_{2}^{u_{2}} \\
&=&\left( -1\right) ^{\alpha \left( 1_{H};0,0,0,0\right)
}B(x_{1}x_{2}\otimes gx_{2};G,gx_{1}x_{2})G\otimes gx_{1}x_{2}\otimes 1_{H}+
\\
&&+\left( -1\right) ^{\alpha \left( 1_{H};0,0,0,1\right)
}B(x_{1}x_{2}\otimes gx_{2};G,gx_{1}x_{2})G\otimes gx_{1}\otimes gx_{2}+ \\
&&+\left( -1\right) ^{\alpha \left( 1_{H};0,0,1,0\right)
}B(x_{1}x_{2}\otimes gx_{2};G,gx_{1}x_{2})G\otimes gx_{2}\otimes gx_{1}+ \\
&&+\left( -1\right) ^{\alpha \left( 1_{H};0,0,1,1\right)
}B(x_{1}x_{2}\otimes gx_{2};G,gx_{1}x_{2})G\otimes g\otimes x_{1}x_{2}
\end{eqnarray*}

\subsubsection{Case $G\otimes gx_{1}\otimes gx_{2}$}

We have only to consider the third and the fourth summands of the left side
of the equality. Third summand gives us%
\begin{eqnarray*}
l_{1} &=&u_{1}=l_{2}=u_{2}=0, \\
a &=&1,b_{1}=b_{2}=0, \\
d &=&e_{1}=1,e_{2}=0
\end{eqnarray*}%
Since $\alpha \left( x_{2};0,0,0,0\right) \equiv a+b_{1}+b_{2}\equiv 1$ we
get
\begin{equation*}
B(gx_{1}\otimes gx_{2};G,gx_{1})G\otimes gx_{1}\otimes gx_{2}.
\end{equation*}%
Fourth summand gives us%
\begin{eqnarray*}
l_{1} &=&u_{1}=0,l_{2}+u_{2}=1 \\
a &=&1,b_{1}=0,b_{2}=l_{2}, \\
d &=&e_{1}=1,0,e_{2}=u_{2}.
\end{eqnarray*}%
Since $\alpha \left( 1_{H};0,0,0,1\right) \equiv a+b_{1}+b_{2}\equiv 1$ and $%
\alpha \left( 1_{H};0,1,0,0\right) \equiv 0,$ we get
\begin{equation*}
\left[ -B(x_{1}x_{2}\otimes gx_{2};G,gx_{1}x_{2})+B(x_{1}x_{2}\otimes
gx_{2};GX_{2},gx_{1})\right] G\otimes gx_{1}\otimes gx_{2}.
\end{equation*}%
By considering also the right side of the equality we obtain%
\begin{equation*}
-B(x_{1}x_{2}\otimes g;G,gx_{1})+B(gx_{1}\otimes
gx_{2};G,gx_{1})-B(x_{1}x_{2}\otimes
gx_{2};G,gx_{1}x_{2})+B(x_{1}x_{2}\otimes gx_{2};GX_{2},gx_{1})=0
\end{equation*}%
which holds in view of the form of the elements.

\subsubsection{Case $G\otimes gx_{2}\otimes gx_{1}$}

We have to consider only the second and the fourth summand of the left side
of the equality. Second summand gives us

\begin{eqnarray*}
l_{1} &=&u_{1}=0,l_{2}=u_{2}=0 \\
a &=&1,b_{1}=b_{2}=0, \\
d &=&e_{2}=1,e_{1}=0
\end{eqnarray*}

Since $\alpha \left( x_{1};0,0,0,0\right) \equiv a+b_{1}+b_{2}\equiv 1,$ we
get%
\begin{equation*}
-B(gx_{2}\otimes gx_{2};G,gx_{2})G\otimes gx_{2}\otimes gx_{1}.
\end{equation*}

The fourth summand of the left side gives us%
\begin{eqnarray*}
l_{1}+u_{1} &=&1,l_{2}=u_{2}=0 \\
a &=&1,b_{2}=0,b_{1}=l_{1} \\
d &=&e_{2}=1,e_{1}=u_{1}.
\end{eqnarray*}%
Since $\alpha \left( 1_{H};0,0,1,0\right) \equiv e_{2}+\left(
a+b_{1}+b_{2}\right) \equiv 0$ and $\alpha \left( 1_{H};0,1,0,0\right) =0,$
we obtain
\begin{equation*}
\left[ B(x_{1}x_{2}\otimes gx_{2};G,gx_{1}x_{2})+B(x_{1}x_{2}\otimes
gx_{2};GX_{1},gx_{2})\right] G\otimes gx_{2}\otimes gx_{1}.
\end{equation*}%
Since there is no term in the right side, we get%
\begin{gather*}
-B(gx_{2}\otimes gx_{2};G,gx_{2})+ \\
B(x_{1}x_{2}\otimes gx_{2};G,gx_{1}x_{2})+B(x_{1}x_{2}\otimes
gx_{2};GX_{1},gx_{2})=0
\end{gather*}%
which holds in view of the form of the elements.

\subsubsection{Case $G\otimes g\otimes x_{1}x_{2}$}

First summand of the left side of the equality gives us%
\begin{eqnarray*}
l_{1} &=&u_{1}=0,l_{2}=u_{2}=0 \\
a &=&1,b_{1}=b_{2}=0, \\
d &=&1,e_{1}=e_{2}=0.
\end{eqnarray*}

Since $\alpha \left( x_{1}x_{2};0,0,0,0\right) \equiv 0,$ we get
\begin{equation*}
B(1_{H}\otimes gx_{2};G,g)G\otimes g\otimes x_{1}x_{2}
\end{equation*}%
Second summand of the left side gives us%
\begin{eqnarray*}
l_{1} &=&u_{1}=0,l_{2}+u_{2}=1 \\
a &=&1,b_{1}=0,b_{2}=l_{2}, \\
d &=&1,e_{1}=0,e_{2}=u_{2}.
\end{eqnarray*}%
Since $\alpha \left( x_{1};0,0,0,1\right) \equiv 1$ and $\alpha \left(
x_{1};0,1,0,0\right) \equiv a+b_{1}+b_{2}+1\equiv 1,$ we obtain
\begin{equation*}
\left[ -B(gx_{2}\otimes gx_{2};G,gx_{2})-B(gx_{2}\otimes gx_{2};GX_{2},g)%
\right] G\otimes g\otimes x_{1}x_{2}
\end{equation*}%
Third summand of the left side gives us%
\begin{eqnarray*}
l_{1}+u_{1} &=&1,l_{2}=u_{2}=0, \\
a &=&1,b_{2}=0,b_{1}=l_{1} \\
d &=&1,e_{2}=0,e_{1}=u_{1}
\end{eqnarray*}%
Since $\alpha \left( x_{2};0,0,1,0\right) \equiv e_{2}=0$ and $\alpha \left(
x_{2};1,0,0,0\right) =a+b_{1}\equiv 0$ we get%
\begin{equation*}
\left[ -B(gx_{1}\otimes gx_{2};G,gx_{1})-B(gx_{1}\otimes gx_{2};GX_{1},g)%
\right] G\otimes g\otimes x_{1}x_{2}
\end{equation*}%
The fourth summand of the left side gives us%
\begin{eqnarray*}
l_{1}+u_{1} &=&1,l_{2}+u_{2}=1 \\
a &=&1,b_{1}=l_{1},b_{2}=l_{2} \\
d &=&1,e_{1}=u_{1},e_{2}=u_{2}.
\end{eqnarray*}%
Since%
\begin{eqnarray*}
\alpha \left( 1_{H};0,0,1,1\right) &\equiv &1+e_{2}\equiv 0\text{, } \\
\alpha \left( 1_{H};0,1,1,0\right) &\equiv &e_{2}+a+b_{1}+b_{2}+1\equiv 1, \\
\alpha \left( 1_{H};1,0,0,1\right) &\equiv &a+b_{1}\equiv 0 \\
\text{and }\alpha \left( 1_{H};1,1,0,0\right) &\equiv &1+b_{2}\equiv 0
\end{eqnarray*}
we get%
\begin{equation*}
\left[
\begin{array}{c}
B(x_{1}x_{2}\otimes gx_{2};G,gx_{1}x_{2})-B(x_{1}x_{2}\otimes
gx_{2};GX_{2},gx_{1}) \\
+B(x_{1}x_{2}\otimes gx_{2};GX_{1},gx_{2})+B(x_{1}x_{2}\otimes
gx_{2};GX_{1}X_{2},g)%
\end{array}%
\right] G\otimes g\otimes x_{1}x_{2}
\end{equation*}%
Since there is no term to consider in the right side of the equality we
obtain

\begin{gather*}
-B(gx_{1}\otimes gx_{2};G,gx_{1})-B(gx_{1}\otimes
gx_{2};GX_{1},g)-B(gx_{2}\otimes gx_{2};G,gx_{2})-B(gx_{2}\otimes
gx_{2};GX_{2},g) \\
B(x_{1}x_{2}\otimes gx_{2};G,gx_{1}x_{2})-B(x_{1}x_{2}\otimes
gx_{2};GX_{2},gx_{1})+B(1_{H}\otimes gx_{2};G,g)+ \\
+B(x_{1}x_{2}\otimes gx_{2};GX_{1},gx_{2})+B(x_{1}x_{2}\otimes
gx_{2};GX_{1}X_{2},g)=0
\end{gather*}%
which holds in view of the form of the elements.

\subsection{$B\left( x_{1}x_{2}\otimes gx_{2};X_{1},gx_{1}x_{2}\right) $}

We deduce that%
\begin{eqnarray*}
a &=&0,b_{1}=1,b_{2}=0,d=1,e_{1}=e_{2}=1. \\
a+b_{1}+b_{2}+l_{1}+l_{2}+d+e_{1}+e_{2}+u_{1}+u_{2} &\equiv
&l_{1}+u_{1}+u_{2}
\end{eqnarray*}%
and we get%
\begin{eqnarray*}
&&\left( -1\right) ^{\alpha \left( 1_{H};0,0,0,0\right) }B(x_{1}x_{2}\otimes
gx_{2};X_{1},gx_{1}x_{2})X_{1}\otimes gx_{1}x_{2}\otimes 1_{H}+ \\
&&\left( -1\right) ^{\alpha \left( 1_{H};0,0,0,1\right) }B(x_{1}x_{2}\otimes
gx_{2};X_{1},gx_{1}x_{2})X_{1}\otimes gx_{1}\otimes gx_{2}+ \\
&&\left( -1\right) ^{\alpha \left( 1_{H};0,0,1,0\right) }B(x_{1}x_{2}\otimes
gx_{2};X_{1},gx_{1}x_{2})X_{1}\otimes gx_{2}\otimes gx_{1}+ \\
&&\left( -1\right) ^{\alpha \left( 1_{H};0,0,1,1\right) }Bx_{1}x_{2}\otimes
gx_{2};X_{1},gx_{1}x_{2})X_{1}\otimes g\otimes x_{1}x_{2}+ \\
&&\left( -1\right) ^{\alpha \left( 1_{H};1,0,0,0\right) }B(x_{1}x_{2}\otimes
gx_{2};X_{1},gx_{1}x_{2})1_{A}\otimes gx_{1}x_{2}\otimes gx_{1}+ \\
&&\left( -1\right) ^{\alpha \left( 1_{H};1,0,0,1\right) }B(x_{1}x_{2}\otimes
gx_{2};X_{1},gx_{1}x_{2})1_{A}\otimes gx_{1}\otimes x_{1}x_{2}+ \\
&&\left( -1\right) ^{\alpha \left( 1_{H};1,0,1,0\right) }B(x_{1}x_{2}\otimes
gx_{2};X_{1},gx_{1}x_{2})1_{A}\otimes gx_{2}^{1-u_{2}}\otimes 0+ \\
&&\left( -1\right) ^{\alpha \left( 1_{H};1,0,1,1\right) }B(x_{1}x_{2}\otimes
gx_{2};X_{1},gx_{1}x_{2})1_{A}\otimes gx_{1}^{1-u_{1}}x_{2}^{1-u_{2}}\otimes
0+
\end{eqnarray*}

\subsubsection{Case $X_{1}\otimes gx_{1}\otimes gx_{2}$}

We have only to consider the third and the fourth summands of the left side
of the equality. Third summand gives us%
\begin{eqnarray*}
l_{1} &=&u_{1}=l_{2}=u_{2}=0, \\
a &=&b_{2}=0,b_{1}=1, \\
d &=&e_{1}=1,e_{2}=0
\end{eqnarray*}%
Since $\alpha \left( x_{2};0,0,0,0\right) \equiv a+b_{1}+b_{2}\equiv 1$ we
get
\begin{equation*}
B(gx_{1}\otimes gx_{2};X_{1},gx_{1})X_{1}\otimes gx_{1}\otimes gx_{2}.
\end{equation*}%
Fourth summand gives us%
\begin{eqnarray*}
l_{1} &=&u_{1}=0,l_{2}+u_{2}=1 \\
a &=&0,b_{1}=1,b_{2}=l_{2}, \\
d &=&e_{1}=1,0,e_{2}=u_{2}.
\end{eqnarray*}%
Since $\alpha \left( 1_{H};0,0,0,1\right) \equiv a+b_{1}+b_{2}\equiv 1$ and $%
\alpha \left( 1_{H};0,1,0,0\right) \equiv 0,$ we get
\begin{equation*}
\left[ -B(x_{1}x_{2}\otimes gx_{2};X_{1},gx_{1}x_{2})+B(x_{1}x_{2}\otimes
gx_{2};X_{1}X_{2},gx_{1})\right] X_{1}\otimes gx_{1}\otimes gx_{2}.
\end{equation*}%
By considering also the right side of the equality we obtain%
\begin{gather*}
-B(x_{1}x_{2}\otimes g;X_{1},gx_{1})+B(gx_{1}\otimes gx_{2};X_{1},gx_{1}) \\
-B(x_{1}x_{2}\otimes gx_{2};X_{1},gx_{1}x_{2})+B(x_{1}x_{2}\otimes
gx_{2};X_{1}X_{2},gx_{1})=0
\end{gather*}%
which holds in view of the form of the elements.

\subsubsection{Case $X_{1}\otimes gx_{2}\otimes gx_{1}$}

We have to consider only the second and the fourth summand of the left side
of the equality. Second summand gives us

\begin{eqnarray*}
l_{1} &=&u_{1}=0,l_{2}=u_{2}=0 \\
a &=&b_{2}=0,b_{1}=1, \\
d &=&e_{2}=1,e_{1}=0
\end{eqnarray*}

Since $\alpha \left( x_{1};0,0,0,0\right) \equiv a+b_{1}+b_{2}\equiv 1,$ we
get%
\begin{equation*}
-B(gx_{2}\otimes gx_{2};X_{1},gx_{2})X_{1}\otimes gx_{2}\otimes gx_{1}.
\end{equation*}

The fourth summand of the left side gives us%
\begin{eqnarray*}
l_{1}+u_{1} &=&1,l_{2}=u_{2}=0 \\
a &=&b_{2}=0,b_{1}-l_{1}=1\Rightarrow b_{1}=1,l_{1}=0,u_{1}=1, \\
d &=&e_{2}=1,e_{1}=u_{1}=1.
\end{eqnarray*}%
Since $\alpha \left( 1_{H};0,0,1,0\right) \equiv e_{2}+\left(
a+b_{1}+b_{2}\right) \equiv 0$ $,$ we obtain
\begin{equation*}
+B(x_{1}x_{2}\otimes gx_{2};X_{1},gx_{1}x_{2})X_{1}\otimes gx_{2}\otimes
gx_{1}.
\end{equation*}%
Since there is no term in the right side, we get%
\begin{equation*}
-B(gx_{2}\otimes gx_{2};X_{1},gx_{2})+B(x_{1}x_{2}\otimes
gx_{2};X_{1},gx_{1}x_{2})=0
\end{equation*}%
which holds in view of the form of the elements.

\subsubsection{Case $X_{1}\otimes g\otimes x_{1}x_{2}$}

First summand of the left side of the equality gives us%
\begin{eqnarray*}
l_{1} &=&u_{1}=0,l_{2}=u_{2}=0 \\
a &=&b_{2}=0,b_{1}=1, \\
d &=&1,e_{1}=e_{2}=0.
\end{eqnarray*}

Since $\alpha \left( x_{1}x_{2};0,0,0,0\right) \equiv 0,$ we get
\begin{equation*}
B(1_{H}\otimes gx_{2};X_{1},g)X_{1}\otimes g\otimes x_{1}x_{2}
\end{equation*}%
Second summand of the left side gives us%
\begin{eqnarray*}
l_{1} &=&u_{1}=0,l_{2}+u_{2}=1 \\
a &=&0,b_{1}=1,b_{2}=l_{2} \\
d &=&1,e_{1}=0,e_{2}=u_{2}.
\end{eqnarray*}%
Since $\alpha \left( x_{1};0,0,0,1\right) \equiv 1$ and $\alpha \left(
x_{1};0,1,0,0\right) \equiv a+b_{1}+b_{2}+1\equiv 1,$ we obtain
\begin{equation*}
\left[ -B(gx_{2}\otimes gx_{2};X_{1},gx_{2})-B(gx_{2}\otimes
gx_{2};X_{1}X_{2},g)\right] X_{1}\otimes g\otimes x_{1}x_{2}
\end{equation*}%
Third summand of the left side gives us%
\begin{eqnarray*}
l_{1}+u_{1} &=&1,l_{2}=u_{2}=0, \\
a &=&b_{2}=0,b_{1}-l_{1}=1\Rightarrow b_{1}=1,l_{1}=0,u_{1}=1 \\
d &=&1,e_{2}=0,e_{1}=u_{1}=1
\end{eqnarray*}%
Since $\alpha \left( x_{2};0,0,1,0\right) \equiv e_{2}=0$ we get%
\begin{equation*}
-B(gx_{1}\otimes gx_{2};X_{1},gx_{1})X_{1}\otimes g\otimes x_{1}x_{2}
\end{equation*}%
The fourth summand of the left side gives us%
\begin{eqnarray*}
l_{1}+u_{1} &=&1,l_{2}+u_{2}=1 \\
a &=&0,b_{1}-l_{1}=1\Rightarrow b_{1}=1,l_{1}=0,u_{1}=1,b_{2}=l_{2} \\
d &=&1,e_{1}=u_{1}=1,e_{2}=u_{2}.
\end{eqnarray*}%
Since $\alpha \left( 1_{H};0,0,1,1\right) \equiv 1+e_{2}\equiv 0$ and $%
\alpha \left( 1_{H};0,1,1,0\right) \equiv e_{2}+a+b_{1}+b_{2}+1\equiv 1,$ we
get%
\begin{equation*}
\left[ B(x_{1}x_{2}\otimes gx_{2};X_{1},gx_{1}x_{2})-B(x_{1}x_{2}\otimes
gx_{2};X_{1}X_{2},gx_{1})\right] X_{1}\otimes g\otimes x_{1}x_{2}
\end{equation*}%
Since there is no term to consider in the right side of the equality we
obtain

\begin{gather*}
-B(gx_{1}\otimes gx_{2};X_{1},gx_{1})-B(gx_{2}\otimes
gx_{2};X_{1},gx_{2})-B(gx_{2}\otimes gx_{2};X_{1}X_{2},g) \\
B(x_{1}x_{2}\otimes gx_{2};X_{1},gx_{1}x_{2})-B(x_{1}x_{2}\otimes
gx_{2};X_{1}X_{2},gx_{1})+B(1_{H}\otimes gx_{2};X_{1},g)+
\end{gather*}%
which holds in view of the form of the elements.

\subsubsection{Case $1_{A}\otimes gx_{1}x_{2}\otimes gx_{1}$}

We have to consider only the second and the fourth summand of the left side
of the equality. Second summand gives us

\begin{eqnarray*}
l_{1} &=&u_{1}=0,l_{2}=u_{2}=0 \\
a &=&b_{2}=b_{1}=0, \\
d &=&e_{1}=e_{2}=1.
\end{eqnarray*}

Since $\alpha \left( x_{1};0,0,0,0\right) \equiv a+b_{1}+b_{2}\equiv 0,$ we
get%
\begin{equation*}
B(gx_{2}\otimes gx_{2};1_{A},gx_{1}x_{2})1_{A}\otimes gx_{1}x_{2}\otimes
gx_{1}.
\end{equation*}

The fourth summand of the left side gives us%
\begin{eqnarray*}
l_{1}+u_{1} &=&1,l_{2}=u_{2}=0 \\
a &=&b_{2}=0,b_{1}-l_{1}=0, \\
d &=&e_{2}=1,e_{1}-u_{1}=1\Rightarrow e_{1}=1,u_{1}=0,b_{1}=l_{1}=1.
\end{eqnarray*}%
Since $\alpha \left( 1_{H};1,0,0,0\right) \equiv b_{2}\equiv 0$ $,$ we
obtain
\begin{equation*}
+B(x_{1}x_{2}\otimes gx_{2};X_{1},gx_{1}x_{2})1_{A}\otimes
gx_{1}x_{2}\otimes gx_{1}.
\end{equation*}%
Since there is no term in the right side, we get%
\begin{equation*}
B(gx_{2}\otimes gx_{2};1_{A},gx_{1}x_{2})+B(x_{1}x_{2}\otimes
gx_{2};X_{1},gx_{1}x_{2})=0
\end{equation*}%
which holds in view of the form of the elements.

\subsubsection{Case $1_{A}\otimes gx_{1}\otimes x_{1}x_{2}$}

First summand of the left side of the equality gives us%
\begin{eqnarray*}
l_{1} &=&u_{1}=0,l_{2}=u_{2}=0 \\
a &=&b_{1}=b_{2}=0, \\
d &=&e_{1}=1,e_{2}=0.
\end{eqnarray*}

Since $\alpha \left( x_{1}x_{2};0,0,0,0\right) \equiv 0,$ we get
\begin{equation*}
B(1_{H}\otimes gx_{2};1_{A},gx_{1})1_{A}\otimes gx_{1}\otimes x_{1}x_{2}
\end{equation*}%
Second summand of the left side gives us%
\begin{eqnarray*}
l_{1} &=&u_{1}=0,l_{2}+u_{2}=1 \\
a &=&b_{1}=0,b_{2}=l_{2} \\
d &=&e_{1}=1,e_{2}=u_{2}.
\end{eqnarray*}%
Since $\alpha \left( x_{1};0,0,0,1\right) \equiv 1$ and $\alpha \left(
x_{1};0,1,0,0\right) \equiv a+b_{1}+b_{2}+1\equiv 0,$ we obtain
\begin{equation*}
\left[ -B(gx_{2}\otimes gx_{2};1_{A},gx_{1}x_{2})+B(gx_{2}\otimes
gx_{2};X_{2},gx_{1})\right] 1_{A}\otimes gx_{1}\otimes x_{1}x_{2}
\end{equation*}%
Third summand of the left side gives us%
\begin{eqnarray*}
l_{1}+u_{1} &=&1,l_{2}=u_{2}=0, \\
a &=&b_{2}=0,b_{1}=l_{1}, \\
d &=&1,e_{2}=0,e_{1}-u_{1}=1\Rightarrow e_{1}=1,u_{1}=0,b_{1}=l_{1}=1.
\end{eqnarray*}%
Since $\alpha \left( x_{2};1,0,0,0\right) \equiv a+b_{1}\equiv 1$ we get%
\begin{equation*}
B(gx_{1}\otimes gx_{2};X_{1},gx_{1})1_{A}\otimes gx_{1}\otimes x_{1}x_{2}
\end{equation*}%
The fourth summand of the left side gives us%
\begin{eqnarray*}
l_{1}+u_{1} &=&1,l_{2}+u_{2}=1 \\
a &=&0,b_{1}=l_{1},b_{2}=l_{2} \\
d &=&1,e_{1}-u_{1}=1\Rightarrow e_{1}=1,u_{1}=0,b_{1}=l_{1}=1,e_{2}=u_{2}.
\end{eqnarray*}%
Since $\alpha \left( 1_{H};1,0,0,1\right) \equiv a+b_{1}\equiv 1$ and $%
\alpha \left( 1_{H};1,1,0,0\right) \equiv 1+b_{2}\equiv 0,$ we get%
\begin{equation*}
\left[ -B(x_{1}x_{2}\otimes gx_{2};X_{1},gx_{1}x_{2})+B(x_{1}x_{2}\otimes
gx_{2};X_{1}X_{2},gx_{1})\right] 1_{A}\otimes gx_{1}\otimes x_{1}x_{2}
\end{equation*}%
Since there is no term to consider in the right side of the equality we
obtain

\begin{gather*}
B(gx_{1}\otimes gx_{2};X_{1},gx_{1})-B(gx_{2}\otimes
gx_{2};1_{A},gx_{1}x_{2})+B(gx_{2}\otimes gx_{2};X_{2},gx_{1}) \\
-B(x_{1}x_{2}\otimes gx_{2};X_{1},gx_{1}x_{2})+B(x_{1}x_{2}\otimes
gx_{2};X_{1}X_{2},gx_{1})+B(1_{H}\otimes gx_{2};1_{A},gx_{1})=0
\end{gather*}%
which holds in view of the form of the elements.

\subsection{$B\left( x_{1}x_{2}\otimes gx_{2};X_{2},gx_{1}x_{2}\right) $}

We deduce that%
\begin{equation*}
a=0,b_{1}=0,b_{2}=1,d=1,e_{1}=1,e_{2}=1
\end{equation*}

and we get%
\begin{multline*}
+\left( -1\right) ^{\alpha \left( 1_{H};0,0,0,0\right) }B(x_{1}x_{2}\otimes
gx_{2};X_{2},gx_{1}x_{2})X_{2}\otimes gx_{1}x_{2}\otimes 1_{H}+ \\
+\left( -1\right) ^{\alpha \left( 1_{H};0,0,0,1\right) }B(x_{1}x_{2}\otimes
gx_{2};X_{2},gx_{1}x_{2})X_{2}\otimes gx_{1}\otimes gx_{2}+ \\
+\left( -1\right) ^{\alpha \left( 1_{H};0,0,1,0\right) }B(x_{1}x_{2}\otimes
gx_{2};X_{2},gx_{1}x_{2})X_{2}\otimes gx_{2}\otimes gx_{1}+ \\
+\left( -1\right) ^{\alpha \left( 1_{H};0,0,1,1\right) }B(x_{1}x_{2}\otimes
gx_{2};X_{2},gx_{1}x_{2})X_{2}\otimes g\otimes x_{1}x_{2}+ \\
+\left( -1\right) ^{\alpha \left( 1_{H};0,1,0,0\right) }B(x_{1}x_{2}\otimes
gx_{2};X_{2},gx_{1}x_{2})1_{A}\otimes gx_{1}x_{2}\otimes gx_{2}+ \\
+\left( -1\right) ^{\alpha \left( 1_{H};0,1,0,1\right) }B(x_{1}x_{2}\otimes
gx_{2};X_{2},gx_{1}x_{2})1_{A}\otimes gx_{1}^{1-u_{1}}x_{2}^{1-u_{2}}\otimes
g^{l_{2}+u_{1}+u_{2}}x_{1}^{u_{1}}x_{2}^{1+1}+=0 \\
+\left( -1\right) ^{\alpha \left( 1_{H};0,1,1,0\right) }B(x_{1}x_{2}\otimes
gx_{2};X_{2},gx_{1}x_{2})1_{A}\otimes gx_{2}\otimes x_{1}x_{2}+ \\
+\left( -1\right) ^{\alpha \left( 1_{H};0,11,1\right) }B(x_{1}x_{2}\otimes
gx_{2};X_{2},gx_{1}x_{2})1_{A}\otimes gx_{1}^{1-u_{1}}x_{2}^{1-u_{2}}\otimes
g^{l_{2}+u_{1}+u_{2}}x_{1}^{u_{1}}x_{2}^{1+1}=0
\end{multline*}

\subsubsection{Case $X_{2}\otimes gx_{1}\otimes gx_{2}$}

We have only to consider the third and the fourth summands of the left side
of the equality. Third summand gives us%
\begin{eqnarray*}
l_{1} &=&u_{1}=l_{2}=u_{2}=0, \\
a &=&b_{1}=0,b_{2}=1, \\
d &=&e_{1}=1,e_{2}=0.
\end{eqnarray*}%
Since $\alpha \left( x_{2};0,0,0,0\right) \equiv a+b_{1}+b_{2}\equiv 1$ we
get
\begin{equation*}
B(gx_{1}\otimes gx_{2};X_{2},gx_{1})X_{2}\otimes gx_{1}\otimes gx_{2}.
\end{equation*}%
Fourth summand gives us%
\begin{eqnarray*}
l_{1} &=&u_{1}=0,l_{2}+u_{2}=1 \\
a &=&b_{1}=0,b_{2}-l_{2}=1\Rightarrow b_{2}=1,l_{2}=0,u_{2}=1, \\
d &=&e_{1}=1,e_{2}=u_{2}=1.
\end{eqnarray*}%
Since $\alpha \left( 1_{H};0,0,0,1\right) \equiv a+b_{1}+b_{2}\equiv 1,$ we
get
\begin{equation*}
-B(x_{1}x_{2}\otimes gx_{2};X_{2},gx_{1}x_{2})X_{2}\otimes gx_{1}\otimes
gx_{2}.
\end{equation*}%
By considering also the right side of the equality we obtain%
\begin{equation*}
-B(x_{1}x_{2}\otimes g;X_{2},gx_{1})+B(gx_{1}\otimes
gx_{2};X_{2},gx_{1})-B(x_{1}x_{2}\otimes gx_{2};X_{2},gx_{1}x_{2})=0
\end{equation*}%
which holds in view of the form of the elements.

\subsubsection{Case $X_{2}\otimes gx_{2}\otimes gx_{1}$}

We have to consider only the second and the fourth summand of the left side
of the equality. Second summand gives us

\begin{eqnarray*}
l_{1} &=&u_{1}=0,l_{2}=u_{2}=0 \\
a &=&b_{1}=0,b_{2}=1, \\
d &=&e_{2}=1,e_{1}=0.
\end{eqnarray*}

Since $\alpha \left( x_{1};0,0,0,0\right) \equiv a+b_{1}+b_{2}\equiv 1,$ we
get%
\begin{equation*}
-B(gx_{2}\otimes gx_{2};X_{2},gx_{2})X_{2}\otimes gx_{2}\otimes gx_{1}.
\end{equation*}

The fourth summand of the left side gives us%
\begin{eqnarray*}
l_{1}+u_{1} &=&1,l_{2}=u_{2}=0 \\
a &=&0,b_{1}=l_{1},b_{2}=1 \\
d &=&e_{2}=1,e_{1}=u_{1}.
\end{eqnarray*}%
Since $\alpha \left( 1_{H};0,0,1,0\right) \equiv e_{2}+\left(
a+b_{1}+b_{2}\right) \equiv 0$ and $\alpha \left( 1_{H};1,0,0,0\right)
\equiv b_{2}\equiv 1,$ we obtain
\begin{equation*}
\left[ +B(x_{1}x_{2}\otimes gx_{2};X_{2},gx_{1}x_{2})-B(x_{1}x_{2}\otimes
gx_{2};X_{1}X_{2},gx_{2})\right] X_{2}\otimes gx_{2}\otimes gx_{1}.
\end{equation*}%
Since there is no term in the right side, we get%
\begin{equation*}
-B(gx_{2}\otimes gx_{2};X_{2},gx_{2})+B(x_{1}x_{2}\otimes
gx_{2};X_{2},gx_{1}x_{2})-B(x_{1}x_{2}\otimes gx_{2};X_{1}X_{2},gx_{2})=0
\end{equation*}%
which holds in view of the form of the elements.

\subsubsection{Case $X_{2}\otimes g\otimes x_{1}x_{2}$}

First summand of the left side of the equality gives us%
\begin{eqnarray*}
l_{1} &=&u_{1}=0,l_{2}=u_{2}=0 \\
a &=&b_{1}=0,b_{2}=1 \\
d &=&1,e_{1}=e_{2}=0.
\end{eqnarray*}

Since $\alpha \left( x_{1}x_{2};0,0,0,0\right) \equiv 0,$ we get
\begin{equation*}
B(1_{H}\otimes gx_{2};X_{2},g)X_{2}\otimes g\otimes x_{1}x_{2}.
\end{equation*}%
Second summand of the left side gives us%
\begin{eqnarray*}
l_{1} &=&u_{1}=0,l_{2}+u_{2}=1 \\
a &=&b_{1}=0,b_{2}-l_{2}=1\Rightarrow b_{2}=1,l_{2}=0,u_{2}=1, \\
d &=&1,e_{1}=0,e_{2}=u_{2}=1.
\end{eqnarray*}%
Since $\alpha \left( x_{1};0,0,0,1\right) \equiv 1,$ we obtain
\begin{equation*}
-B(gx_{2}\otimes gx_{2};X_{2},gx_{2})X_{2}\otimes g\otimes x_{1}x_{2}.
\end{equation*}%
Third summand of the left side gives us%
\begin{eqnarray*}
l_{1}+u_{1} &=&1,l_{2}=u_{2}=0, \\
a &=&0,b_{2}=1,b_{1}=l_{1}, \\
d &=&1,e_{2}=0,e_{1}=u_{1}.
\end{eqnarray*}%
Since $\alpha \left( x_{2};0,0,1,0\right) \equiv e_{2}=0$ $\alpha \left(
x_{2};1,0,0,0\right) \equiv a+b_{1}\equiv 1$ we get%
\begin{equation*}
\left[ -B(gx_{1}\otimes gx_{2};X_{2},gx_{1})+B(gx_{1}\otimes
gx_{2};X_{1}X_{2},g)\right] X_{2}\otimes g\otimes x_{1}x_{2}.
\end{equation*}%
The fourth summand of the left side gives us%
\begin{eqnarray*}
l_{1}+u_{1} &=&1,l_{2}+u_{2}=1 \\
a &=&0,b_{1}=l_{1},b_{2}-l_{2}=1\Rightarrow b_{2}=1,l_{2}=0,u_{2}=1 \\
d &=&1,e_{1}=u_{1},e_{2}=u_{2}=1.
\end{eqnarray*}%
Since $\alpha \left( 1_{H};0,0,1,1\right) \equiv 1+e_{2}\equiv 0$ and $%
\alpha \left( 1_{H};1,0,0,1\right) \equiv a+b_{1}\equiv 1,$ we get%
\begin{equation*}
\left[ +B(x_{1}x_{2}\otimes gx_{2};X_{2},gx_{1}x_{2})-B(x_{1}x_{2}\otimes
gx_{2};X_{1}X_{2},gx_{2})\right] X_{2}\otimes g\otimes x_{1}x_{2}.
\end{equation*}%
Since there is no term to consider in the right side of the equality we
obtain%
\begin{gather*}
-B(gx_{1}\otimes gx_{2};X_{2},gx_{1})+B(gx_{1}\otimes
gx_{2};X_{1}X_{2},g)-B(gx_{2}\otimes gx_{2};X_{2},gx_{2})+ \\
+B(x_{1}x_{2}\otimes gx_{2};X_{2},gx_{1}x_{2})-B(x_{1}x_{2}\otimes
gx_{2};X_{1}X_{2},gx_{2})+B(1_{H}\otimes gx_{2};X_{2},g)=0
\end{gather*}

which holds in view of the form of the elements.

\subsubsection{Case $1_{A}\otimes gx_{1}x_{2}\otimes gx_{2}$}

We have only to consider the third and the fourth summands of the left side
of the equality. Third summand gives us%
\begin{eqnarray*}
l_{1} &=&u_{1}=l_{2}=u_{2}=0, \\
a &=&b_{1}=b_{2}=0, \\
d &=&e_{1}=e_{2}=1.
\end{eqnarray*}%
Since $\alpha \left( x_{2};0,0,0,0\right) \equiv a+b_{1}+b_{2}\equiv 0$ we
get
\begin{equation*}
-B(gx_{1}\otimes gx_{2};1_{A},gx_{1}x_{2})1_{A}\otimes gx_{1}x_{2}\otimes
gx_{2}.
\end{equation*}%
Fourth summand gives us%
\begin{eqnarray*}
l_{1} &=&u_{1}=0,l_{2}+u_{2}=1 \\
a &=&b_{1}=0,b_{2}=l_{2}=1, \\
d &=&e_{1}=1,e_{2}-u_{2}=1\Rightarrow e_{2}=1,u_{2}=0,b_{2}=l_{2}=1.
\end{eqnarray*}%
Since $\alpha \left( 1_{H};0,1,0,0\right) \equiv 0,$ we get
\begin{equation*}
+B(x_{1}x_{2}\otimes gx_{2};X_{2},gx_{1}x_{2})1_{A}\otimes
gx_{1}x_{2}\otimes gx_{2}.
\end{equation*}%
By considering also the right side of the equality we obtain%
\begin{equation*}
-B(x_{1}x_{2}\otimes g;1_{A},gx_{1}x_{2})-B(gx_{1}\otimes
gx_{2};1_{A},gx_{1}x_{2})+B(x_{1}x_{2}\otimes gx_{2};X_{2},gx_{1}x_{2})=0
\end{equation*}%
which holds in view of the form of the elements.

\subsubsection{Case $1_{A}\otimes gx_{2}\otimes x_{1}x_{2}$}

First summand of the left side of the equality gives us%
\begin{eqnarray*}
l_{1} &=&u_{1}=0,l_{2}=u_{2}=0 \\
a &=&b_{1}=b_{2}=0, \\
d &=&e_{2}=1,e_{1}=0.
\end{eqnarray*}

Since $\alpha \left( x_{1}x_{2};0,0,0,0\right) \equiv 0,$ we get
\begin{equation*}
B(1_{H}\otimes gx_{2};1_{A},gx_{2})1_{A}\otimes gx_{2}\otimes x_{1}x_{2}.
\end{equation*}%
Second summand of the left side gives us%
\begin{eqnarray*}
l_{1} &=&u_{1}=0,l_{2}+u_{2}=1 \\
a &=&b_{1}=0,b_{2}=l_{2}, \\
d &=&1,e_{1}=0,e_{2}-u_{2}=1\Rightarrow e_{2}=1,u_{2}=0,b_{2}=l_{2}=1.
\end{eqnarray*}%
Since $\alpha \left( x_{1};0,1,0,0\right) \equiv a+b_{1}+b_{2}+1\equiv 0,$
we obtain
\begin{equation*}
B(gx_{2}\otimes gx_{2};X_{2},gx_{2})1_{A}\otimes gx_{2}\otimes x_{1}x_{2}.
\end{equation*}%
Third summand of the left side gives us%
\begin{eqnarray*}
l_{1}+u_{1} &=&1,l_{2}=u_{2}=0, \\
a &=&b_{2}=0,b_{1}=l_{1}, \\
d &=&1,e_{2}=1,e_{1}=u_{1}.
\end{eqnarray*}%
Since $\alpha \left( x_{2};0,0,1,0\right) \equiv e_{2}=1$ and $\alpha \left(
x_{2};1,0,0,0\right) \equiv a+b_{1}\equiv 1$ we get%
\begin{equation*}
\left[ +B(gx_{1}\otimes gx_{2};1_{A},gx_{1}x_{2})+B(gx_{1}\otimes
gx_{2};X_{1},gx_{2})\right] 1_{A}\otimes gx_{2}\otimes x_{1}x_{2}.
\end{equation*}%
The fourth summand of the left side gives us%
\begin{eqnarray*}
l_{1}+u_{1} &=&1,l_{2}+u_{2}=1 \\
a &=&0,b_{1}=l_{1},b_{2}=l_{2}, \\
d &=&1,e_{1}=u_{1},e_{2}-u_{2}=1\Rightarrow e_{2}=1,u_{2}=0,b_{2}=l_{2}=1.
\end{eqnarray*}%
Since $\alpha \left( 1_{H};0,1,1,0\right) \equiv e_{2}+a+b_{1}+b_{2}+1\equiv
1$ and $\alpha \left( 1_{H};1,1,0,0\right) \equiv 1+b_{2}\equiv 0$ $,$ we get%
\begin{equation*}
\left[ -B(x_{1}x_{2}\otimes gx_{2};X_{2},gx_{1}x_{2})+B(x_{1}x_{2}\otimes
gx_{2};X_{1}X_{2},gx_{2})\right] 1_{A}\otimes gx_{2}\otimes x_{1}x_{2}.
\end{equation*}%
Since there is no term to consider in the right side of the equality we
obtain%
\begin{gather*}
+B(gx_{1}\otimes gx_{2};1_{A},gx_{1}x_{2})+B(gx_{1}\otimes
gx_{2};X_{1},gx_{2})+B(gx_{2}\otimes gx_{2};X_{2},gx_{2})+ \\
-B(x_{1}x_{2}\otimes gx_{2};X_{2},gx_{1}x_{2})+B(x_{1}x_{2}\otimes
gx_{2};X_{1}X_{2},gx_{2})+B(1_{H}\otimes gx_{2};1_{A},gx_{2})=0
\end{gather*}

which holds in view of the form of the elements.

\subsection{$B\left( x_{1}x_{2}\otimes gx_{2};X_{1}X_{2},gx_{2}\right) $}

We have%
\begin{equation*}
a=0,b_{1}=b_{2}=1,d=e_{2}=1
\end{equation*}%
and we get%
\begin{gather*}
+\left( -1\right) ^{\alpha \left( 1_{H};0,0,0,0\right) }B(x_{1}x_{2}\otimes
gx_{2};X_{1}X_{2},gx_{2})X_{1}X_{2}\otimes gx_{2}\otimes 1_{A}+ \\
+\left( -1\right) ^{\alpha \left( 1_{H};0,1,0,0\right) }B(x_{1}x_{2}\otimes
gx_{2};X_{1}X_{2},gx_{2})X_{1}\otimes gx_{2}\otimes gx_{2}+ \\
+\left( -1\right) ^{\alpha \left( 1_{H};1,0,0,0\right) }B(x_{1}x_{2}\otimes
gx_{2};X_{1}X_{2},gx_{2})X_{2}\otimes gx_{2}\otimes gx_{1}+ \\
+\left( -1\right) ^{\alpha \left( 1_{H};1,1,0,0\right) }B(x_{1}x_{2}\otimes
gx_{2};X_{1}X_{2},gx_{2})1_{A}\otimes gx_{2}\otimes x_{1}x_{2}+ \\
+\left( -1\right) ^{\alpha \left( 1_{H};0,0,0,1\right) }B(x_{1}x_{2}\otimes
gx_{2};X_{1}X_{2},gx_{2})X_{1}X_{2}\otimes g\otimes gx_{2}+ \\
+\left( -1\right) ^{\alpha \left( 1_{H};0,1,0,1\right) }B(x_{1}x_{2}\otimes
gx_{2};X_{1}X_{2},gx_{2})X_{1}^{1-l_{1}}X_{2}^{1-l_{2}}\otimes g\otimes
g^{l_{1}+l_{2}+1}x_{1}^{l_{1}}x_{2}^{1+1}=0 \\
+\left( -1\right) ^{\alpha \left( 1_{H};1,0,0,1\right) }B(x_{1}x_{2}\otimes
gx_{2};X_{1}X_{2},gx_{2})X_{2}\otimes g\otimes x_{1}x_{2}+ \\
+\left( -1\right) ^{\alpha \left( 1_{H};1,1,0,1\right) }B(x_{1}x_{2}\otimes
gx_{2};X_{1}X_{2},gx_{2})X_{1}^{1-l_{1}}X_{2}^{1-l_{2}}\otimes g\otimes
g^{1+l_{1}+l_{2}}x_{1}^{l_{1}}x_{2}^{1+1}=0.
\end{gather*}

\subsubsection{Case $X_{1}\otimes gx_{2}\otimes gx_{2}$}

We have only to consider the third and the fourth summands of the left side
of the equality. Third summand gives us%
\begin{eqnarray*}
l_{1} &=&u_{1}=l_{2}=u_{2}=0, \\
a &=&b_{2}=0,b_{1}=1, \\
d &=&e_{1}=1,e_{2}=0.
\end{eqnarray*}%
Since $\alpha \left( x_{2};0,0,0,0\right) \equiv a+b_{1}+b_{2}\equiv 1$ we
get
\begin{equation*}
+B(gx_{1}\otimes gx_{2};X_{1},gx_{2})X_{1}\otimes gx_{2}\otimes gx_{2}.
\end{equation*}%
Fourth summand gives us%
\begin{eqnarray*}
l_{1} &=&u_{1}=0,l_{2}+u_{2}=1 \\
a &=&0,b_{1}=1,b_{2}=l_{2}, \\
d &=&1,e_{1}=0,e_{2}-u_{2}=1\Rightarrow e_{2}=1,u_{2}=0,b_{2}=l_{2}=1.
\end{eqnarray*}%
Since $\alpha \left( 1_{H};0,1,0,0\right) \equiv 0,$ we get
\begin{equation*}
+B(x_{1}x_{2}\otimes gx_{2};X_{1}X_{2},gx_{2})X_{1}\otimes gx_{2}\otimes
gx_{2}.
\end{equation*}%
By considering also the right side of the equality we obtain%
\begin{equation*}
-B(x_{1}x_{2}\otimes g;X_{1},gx_{2})+B(gx_{1}\otimes
gx_{2};X_{1},gx_{2})+B(x_{1}x_{2}\otimes gx_{2};X_{1}X_{2},gx_{2})=0
\end{equation*}%
which holds in view of the form of the elements.

\subsubsection{Case $X_{2}\otimes gx_{2}\otimes gx_{1}$}

This case was already considered in subsection $B\left( x_{1}x_{2}\otimes
gx_{2};X_{2},gx_{1}x_{2}\right) .$

\subsubsection{Case $1_{A}\otimes gx_{2}\otimes x_{1}x_{2}$}

This case was already considered in subsection $B\left( x_{1}x_{2}\otimes
gx_{2};X_{2},gx_{1}x_{2}\right) .$

\subsubsection{Case $X_{1}X_{2}\otimes g\otimes gx_{2}$}

We have only to consider the third and the fourth summands of the left side
of the equality. Third summand gives us%
\begin{eqnarray*}
l_{1} &=&u_{1}=l_{2}=u_{2}=0, \\
a &=&0,b_{1}=b_{2}=1, \\
d &=&1,e_{1}=e_{2}=0.
\end{eqnarray*}%
Since $\alpha \left( x_{2};0,0,0,0\right) \equiv a+b_{1}+b_{2}\equiv 0$ we
get
\begin{equation*}
-B(gx_{1}\otimes gx_{2};X_{1}X_{2},g)X_{1}X_{2}\otimes g\otimes gx_{2}.
\end{equation*}%
Fourth summand gives us%
\begin{eqnarray*}
l_{1} &=&u_{1}=0,l_{2}+u_{2}=1 \\
a &=&0,b_{1}=1,b_{2}-l_{2}=1\Rightarrow b_{2}=1,l_{2}=0,u_{2}=1 \\
d &=&1,e_{1}=0,e_{2}=u_{2}=1.
\end{eqnarray*}%
Since $\alpha \left( 1_{H};0,0,0,1\right) \equiv a+b_{1}+b_{2}\equiv 0,$ we
get
\begin{equation*}
+B(x_{1}x_{2}\otimes gx_{2};X_{1}X_{2},gx_{2})X_{1}X_{2}\otimes g\otimes
gx_{2}.
\end{equation*}%
By considering also the right side of the equality we obtain%
\begin{equation*}
-B(x_{1}x_{2}\otimes g;X_{1}X_{2},g)-B(gx_{1}\otimes
gx_{2};X_{1}X_{2},g)+B(x_{1}x_{2}\otimes gx_{2};X_{1}X_{2},gx_{2})=0
\end{equation*}%
which holds in view of the form of the elements.

\subsubsection{Case $X_{2}\otimes g\otimes x_{1}x_{2}$}

This case was already considered in subsection $B\left( x_{1}x_{2}\otimes
gx_{2};X_{2},gx_{1}x_{2}\right) .$

\subsection{$B\left( x_{1}x_{2}\otimes gx_{2};GX_{1},gx_{2}\right) $}

We deduce that
\begin{equation*}
a=b_{1}=1,d=1,e_{2}=1
\end{equation*}%
and we get%
\begin{eqnarray*}
&&\left( -1\right) ^{\alpha \left( 1_{H};0,0,0,0\right) }B(x_{1}x_{2}\otimes
gx_{2};X_{1},gx_{2})GX_{1}\otimes gx_{2}\otimes 1_{H}+ \\
&&\left( -1\right) ^{\alpha \left( 1_{H};1,0,0,0\right) }B(x_{1}x_{2}\otimes
gx_{2};X_{1},gx_{2})G\otimes gx_{2}\otimes gx_{1} \\
&&\left( -1\right) ^{\alpha \left( 1_{H};0,0,0,1\right) }B(x_{1}x_{2}\otimes
gx_{2};X_{1},gx_{2})GX_{1}\otimes g\otimes gx_{2} \\
&&\left( -1\right) ^{\alpha \left( 1_{H};1,0,0,1\right) }B(x_{1}x_{2}\otimes
gx_{2};X_{1},gx_{2})G\otimes g\otimes x_{1}x_{2}.
\end{eqnarray*}

\subsubsection{Case $G\otimes gx_{2}\otimes gx_{1}$}

This case was already considered in subsection $B\left( x_{1}x_{2}\otimes
gx_{2};G,gx_{1}x_{2}\right) .$

\subsubsection{Case $GX_{1}\otimes g\otimes gx_{2}$}

We have only to consider the third and the fourth summands of the left side
of the equality. Third summand gives us%
\begin{eqnarray*}
l_{1} &=&u_{1}=l_{2}=u_{2}=0, \\
a &=&b_{1}=1,b_{2}=0, \\
d &=&1,e_{1}=e_{2}=0.
\end{eqnarray*}%
Since $\alpha \left( x_{2};0,0,0,0\right) \equiv a+b_{1}+b_{2}\equiv 0$ we
get
\begin{equation*}
-B(gx_{1}\otimes gx_{2};GX_{1},g)GX_{1}\otimes g\otimes gx_{2}.
\end{equation*}%
Fourth summand gives us%
\begin{eqnarray*}
l_{1} &=&u_{1}=0,l_{2}+u_{2}=1 \\
a &=&b_{1}=1,b_{2}=l_{2} \\
d &=&1,e_{1}=0,e_{2}=u_{2}.
\end{eqnarray*}%
Since $\alpha \left( 1_{H};0,0,0,1\right) \equiv a+b_{1}+b_{2}\equiv 0$ and $%
\alpha \left( 1_{H};0,1,0,0\right) \equiv 0,$ we get
\begin{equation*}
\left[ +B(x_{1}x_{2}\otimes gx_{2};GX_{1},gx_{2})+B(x_{1}x_{2}\otimes
gx_{2};GX_{1}X_{2},g)\right] GX_{1}\otimes g\otimes gx_{2}.
\end{equation*}%
By considering also the right side of the equality we obtain%
\begin{gather*}
-B(x_{1}x_{2}\otimes g;GX_{1},g)-B(gx_{1}\otimes gx_{2};GX_{1},g) \\
+B(x_{1}x_{2}\otimes gx_{2};GX_{1},gx_{2})+B(x_{1}x_{2}\otimes
gx_{2};GX_{1}X_{2},g)=0
\end{gather*}%
which holds in view of the form of the elements.

\subsubsection{Case $G\otimes g\otimes x_{1}x_{2}$}

This case was already considered in subsection $B\left( x_{1}x_{2}\otimes
gx_{2};G,gx_{1}x_{2}\right) .$

\subsection{$B\left( x_{1}x_{2}\otimes gx_{2};GX_{2},gx_{2}\right) $}

We deduce that%
\begin{eqnarray*}
a &=&1,b_{1}=0,b_{2}=1 \\
d &=&1,e_{1}=0,e_{2}=1
\end{eqnarray*}%
and we get%
\begin{eqnarray*}
&&\left( -1\right) ^{\alpha \left( 1_{H};0,0,0,0\right) }B(x_{1}x_{2}\otimes
gx_{2};GX_{2},gx_{2})GX_{2}\otimes gx_{2}\otimes 1_{H}+ \\
&&\left( -1\right) ^{\alpha \left( 1_{H};0,0,0,1\right) }Bx_{1}x_{2}\otimes
gx_{2};GX_{2},gx_{2})GX_{2}\otimes g\otimes gx_{2}+ \\
&&\left( -1\right) ^{\alpha \left( 1_{H};0,1,0,0\right) }B(x_{1}x_{2}\otimes
gx_{2};GX_{2},gx_{2})G\otimes gx_{2}\otimes gx_{2}+ \\
&&\left( -1\right) ^{\alpha \left( 1_{H};0,1,0,1\right) }B(x_{1}x_{2}\otimes
gx_{2};GX_{2},gx_{2})GX_{2}^{1-l_{2}}\otimes gx_{2}^{1-u_{2}}\otimes
g^{l_{2}+u_{2}}x_{2}^{1+1}=0
\end{eqnarray*}

\subsubsection{Case $GX_{2}\otimes g\otimes gx_{2}$}

We have only to consider the third and the fourth summands of the left side
of the equality. Third summand gives us%
\begin{eqnarray*}
l_{1} &=&u_{1}=l_{2}=u_{2}=0, \\
a &=&b_{2}=1,b_{1}=0, \\
d &=&1,e_{1}=e_{2}=0.
\end{eqnarray*}%
Since $\alpha \left( x_{2};0,0,0,0\right) \equiv a+b_{1}+b_{2}\equiv 0$ we
get
\begin{equation*}
-B(gx_{1}\otimes gx_{2};GX_{2},g)GX_{2}\otimes g\otimes gx_{2}.
\end{equation*}%
Fourth summand gives us%
\begin{eqnarray*}
l_{1} &=&u_{1}=0,l_{2}+u_{2}=1 \\
a &=&1,b_{1}=0,b_{2}-l_{2}=1\Rightarrow b_{2}=1,l_{2}=0,u_{2}=1, \\
d &=&1,e_{1}=0,e_{2}=u_{2}=1.
\end{eqnarray*}%
Since $\alpha \left( 1_{H};0,0,0,1\right) \equiv a+b_{1}+b_{2}\equiv 0$ $,$
we get
\begin{equation*}
+B(x_{1}x_{2}\otimes gx_{2};GX_{2},gx_{2})GX_{2}\otimes g\otimes gx_{2}.
\end{equation*}%
By considering also the right side of the equality we obtain%
\begin{equation*}
-B(x_{1}x_{2}\otimes g;GX_{2},g)-B(gx_{1}\otimes
gx_{2};GX_{2},g)+B(x_{1}x_{2}\otimes gx_{2};GX_{2},gx_{2})=0
\end{equation*}%
which holds in view of the form of the elements.

\subsubsection{Case $G\otimes gx_{2}\otimes gx_{2}$}

We have only to consider the third and the fourth summands of the left side
of the equality. Third summand gives us%
\begin{eqnarray*}
l_{1} &=&u_{1}=l_{2}=u_{2}=0, \\
a &=&1,b_{1}=b_{2}=0, \\
d &=&e_{2}=1,e_{1}=0.
\end{eqnarray*}%
Since $\alpha \left( x_{2};0,0,0,0\right) \equiv a+b_{1}+b_{2}\equiv 1$ we
get
\begin{equation*}
+B(gx_{1}\otimes gx_{2};G,gx_{2})G\otimes gx_{2}\otimes gx_{2}.
\end{equation*}%
Fourth summand gives us%
\begin{eqnarray*}
l_{1} &=&u_{1}=0,l_{2}+u_{2}=1 \\
a &=&1,b_{1}=0,b_{2}=l_{2}, \\
d &=&1,e_{1}=0,e_{2}-u_{2}=1\Rightarrow e_{2}=1,u_{2}=0,b_{2}=l_{2}=1.
\end{eqnarray*}%
Since $\alpha \left( 1_{H};0,1,0,0\right) \equiv 0$ $,$ we get
\begin{equation*}
+B(x_{1}x_{2}\otimes gx_{2};GX_{2},gx_{2})G\otimes gx_{2}\otimes gx_{2}.
\end{equation*}%
By considering also the right side of the equality we obtain%
\begin{equation*}
-B(x_{1}x_{2}\otimes g;G,gx_{2})+B(gx_{1}\otimes
gx_{2};G,gx_{2})+B(x_{1}x_{2}\otimes gx_{2};GX_{2},gx_{2})=0
\end{equation*}%
which holds in view of the form of the elements.

\subsection{$B\left( x_{1}x_{2}\otimes gx_{2};GX_{1}X_{2},gx_{1}x_{2}\right)
$}

We deduce that%
\begin{equation*}
a=b_{1}=b_{2}=d=e_{1}=e_{2}=1
\end{equation*}%
and we obtain%
\begin{gather*}
\left( -1\right) ^{\alpha \left( 1_{H};0,0,0,0\right) }B(x_{1}x_{2}\otimes
gx_{2};GX_{1}X_{2},gx_{1}x_{2})GX_{1}X_{2}\otimes gx_{1}x_{2}\otimes 1_{H}+
\\
+\left( -1\right) ^{\alpha \left( 1_{H};0,0,0,1\right) }B(x_{1}x_{2}\otimes
gx_{2};GX_{1}X_{2},gx_{1}x_{2})GX_{1}X_{2}\otimes gx_{1}\otimes gx_{2}\text{
\ } \\
+\left( -1\right) ^{\alpha \left( 1_{H};0,0,1,0\right) }B(x_{1}x_{2}\otimes
gx_{2};GX_{1}X_{2},gx_{1}x_{2})GX_{1}X_{2}\otimes gx_{2}\otimes gx_{1}+ \\
+\left( -1\right) ^{\alpha \left( 1_{H};0,0,1,1\right) }B(x_{1}x_{2}\otimes
x_{1}x_{2};GX_{1}X_{2},gx_{1}x_{2})GX_{1}X_{2}\otimes g\otimes x_{1}x_{2}+ \\
+\left( -1\right) ^{\alpha \left( 1_{H};0,1,0,0\right) }B(x_{1}x_{2}\otimes
gx_{2};GX_{1}X_{2},gx_{1}x_{2})GX_{1}\otimes gx_{1}x_{2}\otimes gx_{2}+ \\
+\left( -1\right) ^{\alpha \left( 1_{H};0,1,0,1\right) }B(x_{1}x_{2}\otimes
gx_{2};GX_{1}X_{2},gx_{1}x_{2})GX_{1}\otimes gx_{1}x_{2}^{1-u_{2}}\otimes
g^{u_{2}}x_{2}^{1+1}=0 \\
+\left( -1\right) ^{\alpha \left( 1_{H};0,1,1,0\right) }B(x_{1}x_{2}\otimes
gx_{2};GX_{1}X_{2},gx_{1}x_{2})GX_{1}\otimes gx_{2}\otimes x_{1}x_{2}+ \\
+\left( -1\right) ^{\alpha \left( 1_{H};0,1,1,1\right) }B(x_{1}x_{2}\otimes
gx_{2};GX_{1}X_{2},gx_{1}x_{2})GX_{1}\otimes gx_{2}^{1-u_{2}}\otimes
g^{1+u_{2}}x_{1}x_{2}^{1+1}=0 \\
+\left( -1\right) ^{\alpha \left( 1_{H};1,0,0,0\right) }B(x_{1}x_{2}\otimes
gx_{2};GX_{1}X_{2},gx_{1}x_{2})GX_{2}\otimes gx_{1}x_{2}\otimes gx_{1}+ \\
+\left( -1\right) ^{\alpha \left( 1_{H};1,0,0,1\right) }B(x_{1}x_{2}\otimes
gx_{2};GX_{1}X_{2},gx_{1}x_{2})GX_{2}\otimes gx_{1}x_{2}^{1-u_{2}}\otimes
g^{1+u_{2}}x_{1}^{1+0}x_{2}^{1+1}=0 \\
+\left( -1\right) ^{\alpha \left( 1_{H};1,0,1,u_{2}\right)
}B(x_{1}x_{2}\otimes gx_{2};GX_{1}X_{2},gx_{1}x_{2})GX_{2}\otimes
gx_{1}^{1-u_{1}}x_{2}^{1-u_{2}}\otimes
g^{1+u_{1}+u_{2}}x_{1}^{1+1}x_{2}^{1+u_{2}}=0 \\
+\left( -1\right) ^{\alpha \left( 1_{H};1,1,0,0\right) }B(x_{1}x_{2}\otimes
gx_{2};GX_{1}X_{2},gx_{1}x_{2})G\otimes gx_{1}x_{2}\otimes x_{1}x_{2}+ \\
+\left( -1\right) ^{\alpha \left( 1_{H};1,1,0,1\right) }B(x_{1}x_{2}\otimes
gx_{2};GX_{1}X_{2},gx_{1}x_{2})G\otimes gx_{1}x_{2}^{1-u_{2}}\otimes
g^{u_{2}}x_{1}x_{2}^{1+1}+=0 \\
+\left( -1\right) ^{\alpha \left( 1_{H};1,1,1,u_{2}\right)
}B(x_{1}x_{2}\otimes gx_{2};GX_{1}X_{2},gx_{1}x_{2})G\otimes
gx_{1}^{1-u_{1}}x_{2}^{1-u_{2}}\otimes
g^{u_{1}+u_{2}}x_{1}^{1+1}x_{2}^{1+u_{2}}=0.
\end{gather*}

\subsubsection{Case $GX_{1}X_{2}\otimes gx_{1}\otimes gx_{2}$}

We have only to consider the third and the fourth summands of the left side
of the equality. Third summand gives us%
\begin{eqnarray*}
l_{1} &=&u_{1}=l_{2}=u_{2}=0, \\
a &=&b_{1}=b_{2}=1, \\
d &=&e_{1}=1,e_{2}=0.
\end{eqnarray*}%
Since $\alpha \left( x_{2};0,0,0,0\right) \equiv a+b_{1}+b_{2}\equiv 1$ we
get
\begin{equation*}
+B(gx_{1}\otimes gx_{2};GX_{1}X_{2},gx_{1})GX_{1}X_{2}\otimes gx_{1}\otimes
gx_{2}.
\end{equation*}%
Fourth summand gives us%
\begin{eqnarray*}
l_{1} &=&u_{1}=0,l_{2}+u_{2}=1 \\
a &=&b_{1}=1,b_{2}-l_{2}=1\Rightarrow b_{2}=1,l_{2}=0,u_{2}=1 \\
d &=&e_{1}=1,0,e_{2}=u_{2}=1.
\end{eqnarray*}%
Since $\alpha \left( 1_{H};0,0,0,1\right) \equiv a+b_{1}+b_{2}\equiv 1$ $,$
we get
\begin{equation*}
-B(x_{1}x_{2}\otimes gx_{2};GX_{1}X_{2},gx_{1}x_{2})GX_{1}X_{2}\otimes
gx_{1}\otimes gx_{2}.
\end{equation*}%
By considering also the right side of the equality we obtain%
\begin{equation*}
-B(x_{1}x_{2}\otimes g;GX_{1}X_{2},gx_{1})+B(gx_{1}\otimes
gx_{2};GX_{1}X_{2},gx_{1})-B(x_{1}x_{2}\otimes
gx_{2};GX_{1}X_{2},gx_{1}x_{2})=0
\end{equation*}%
which holds in view of the form of the elements.

\subsubsection{Case $GX_{1}X_{2}\otimes gx_{2}\otimes gx_{1}$}

We have to consider only the second and the fourth summand of the left side
of the equality. Second summand gives us

\begin{eqnarray*}
l_{1} &=&u_{1}=0,l_{2}=u_{2}=0 \\
a &=&b_{1}=b_{2}=1, \\
d &=&e_{2}=1,e_{1}=0.
\end{eqnarray*}

Since $\alpha \left( x_{1};0,0,0,0\right) \equiv a+b_{1}+b_{2}\equiv 1,$ we
get%
\begin{equation*}
-B(gx_{2}\otimes gx_{2};GX_{1}X_{2},gx_{2})GX_{1}X_{2}\otimes gx_{2}\otimes
gx_{1}.
\end{equation*}

The fourth summand of the left side gives us%
\begin{eqnarray*}
l_{1}+u_{1} &=&1,l_{2}=u_{2}=0 \\
a &=&1,b_{1}-l_{1}=1\Rightarrow b_{1}=1,l_{1}=0,u_{1}=1,b_{2}=1 \\
d &=&e_{2}=1,e_{1}=u_{1}=1.
\end{eqnarray*}%
Since $\alpha \left( 1_{H};0,0,1,0\right) \equiv e_{2}+\left(
a+b_{1}+b_{2}\right) \equiv 1$ $,$ we obtain
\begin{equation*}
B(x_{1}x_{2}\otimes gx_{2};GX_{1}X_{2},gx_{1}x_{2})GX_{1}X_{2}\otimes
gx_{2}\otimes gx_{1}.
\end{equation*}%
Since there is no term in the right side, we get%
\begin{equation*}
-B(gx_{2}\otimes gx_{2};GX_{1}X_{2},gx_{2})+B(x_{1}x_{2}\otimes
gx_{2};GX_{1}X_{2},gx_{1}x_{2})=0
\end{equation*}%
which holds in view of the form of the elements.

\subsubsection{Case $GX_{1}X_{2}\otimes g\otimes x_{1}x_{2}$}

First summand of the left side of the equality gives us%
\begin{eqnarray*}
l_{1} &=&u_{1}=0,l_{2}=u_{2}=0 \\
a &=&b_{1}=b_{2}=1, \\
d &=&1,e_{1}=e_{2}=0.
\end{eqnarray*}

Since $\alpha \left( x_{1}x_{2};0,0,0,0\right) \equiv 0,$ we get
\begin{equation*}
B(1_{H}\otimes gx_{2};GX_{1}X_{2},g)GX_{1}X_{2}\otimes g\otimes x_{1}x_{2}.
\end{equation*}%
Second summand of the left side gives us%
\begin{eqnarray*}
l_{1} &=&u_{1}=0,l_{2}+u_{2}=1 \\
a &=&b_{1}=1,b_{2}-l_{2}=1\Rightarrow b_{2}=1,l_{2}=0,u_{2}=1, \\
d &=&1,e_{1}=0,e_{2}=u_{2}=1.
\end{eqnarray*}%
Since $\alpha \left( x_{1};0,0,0,1\right) \equiv 1,$ we obtain
\begin{equation*}
-B(gx_{2}\otimes gx_{2};GX_{1}X_{2},gx_{2})GX_{1}X_{2}\otimes g\otimes
x_{1}x_{2}.
\end{equation*}%
Third summand of the left side gives us%
\begin{eqnarray*}
l_{1}+u_{1} &=&1,l_{2}=u_{2}=0, \\
a &=&b_{2}=1,b_{1}-l_{1}=1\Rightarrow b_{1}=1,l_{1}=0,u_{1}=1 \\
d &=&1,e_{2}=0,e_{1}=u_{1}=1.
\end{eqnarray*}%
Since $\alpha \left( x_{2};0,0,1,0\right) \equiv e_{2}=0$ we get%
\begin{equation*}
-B(gx_{1}\otimes gx_{2};GX_{1}X_{2},gx_{1})GX_{1}X_{2}\otimes g\otimes
x_{1}x_{2}.
\end{equation*}%
The fourth summand of the left side gives us%
\begin{eqnarray*}
l_{1}+u_{1} &=&1,l_{2}+u_{2}=1 \\
a &=&1,b_{1}-l_{1}=1\Rightarrow b_{1}=1,l_{1}=0,u_{1}=1, \\
b_{2}-l_{2} &=&1\Rightarrow b_{2}=1,l_{2}=0,u_{2}=1 \\
d &=&1,e_{1}=u_{1}=1,e_{2}=u_{2}=1.
\end{eqnarray*}%
Since $\alpha \left( 1_{H};0,0,1,1\right) \equiv 1+e_{2}\equiv 0,$ we get%
\begin{equation*}
+B(x_{1}x_{2}\otimes gx_{2};GX_{1}X_{2},gx_{1}x_{2})GX_{1}X_{2}\otimes
g\otimes x_{1}x_{2}.
\end{equation*}%
Since there is no term to consider in the right side of the equality we
obtain%
\begin{gather*}
-B(gx_{1}\otimes gx_{2};GX_{1}X_{2},gx_{1})-B(gx_{2}\otimes
gx_{2};GX_{1}X_{2},gx_{2})+ \\
+B(x_{1}x_{2}\otimes gx_{2};GX_{1}X_{2},gx_{1}x_{2})+B(1_{H}\otimes
gx_{2};GX_{1}X_{2},g)=0
\end{gather*}

which holds in view of the form of the elements.

\subsubsection{Case $GX_{1}\otimes gx_{1}x_{2}\otimes gx_{2}$}

We have only to consider the third and the fourth summands of the left side
of the equality. Third summand gives us%
\begin{eqnarray*}
l_{1} &=&u_{1}=l_{2}=u_{2}=0, \\
a &=&b_{1}=1,b_{2}=0, \\
d &=&e_{1}=e_{2}=1.
\end{eqnarray*}%
Since $\alpha \left( x_{2};0,0,0,0\right) \equiv a+b_{1}+b_{2}\equiv 0$ we
get
\begin{equation*}
-B(gx_{1}\otimes gx_{2};GX_{1},gx_{1}x_{2})GX_{1}\otimes gx_{1}x_{2}\otimes
gx_{2}.
\end{equation*}%
Fourth summand gives us%
\begin{eqnarray*}
l_{1} &=&u_{1}=0,l_{2}+u_{2}=1 \\
a &=&b_{1}=1,b_{2}=l_{2}, \\
d &=&e_{1}=1,e_{2}-u_{2}=1\Rightarrow e_{2}=1,u_{2}=0,b_{2}=l_{2}=1.
\end{eqnarray*}%
Since $\alpha \left( 1_{H};0,1,0,0\right) \equiv 0$ $,$ we get
\begin{equation*}
+B(x_{1}x_{2}\otimes gx_{2};GX_{1}X_{2},gx_{1}x_{2})GX_{1}\otimes
gx_{1}x_{2}\otimes gx_{2}.
\end{equation*}%
By considering also the right side of the equality we obtain%
\begin{equation*}
-B(x_{1}x_{2}\otimes g;GX_{1},gx_{1}x_{2})-B(gx_{1}\otimes
gx_{2};GX_{1},gx_{1}x_{2})+B(x_{1}x_{2}\otimes
gx_{2};GX_{1}X_{2},gx_{1}x_{2})=0
\end{equation*}%
which holds in view of the form of the elements.

\subsubsection{Case $GX_{1}\otimes gx_{2}\otimes x_{1}x_{2}$}

First summand of the left side of the equality gives us%
\begin{eqnarray*}
l_{1} &=&u_{1}=0,l_{2}=u_{2}=0 \\
a &=&b_{1}=1,b_{2}=0 \\
d &=&e_{2}=1,e_{1}=0.
\end{eqnarray*}

Since $\alpha \left( x_{1}x_{2};0,0,0,0\right) \equiv 0,$ we get
\begin{equation*}
B(1_{H}\otimes gx_{2};GX_{1},gx_{2})GX_{1}\otimes gx_{2}\otimes x_{1}x_{2}.
\end{equation*}%
Second summand of the left side gives us%
\begin{eqnarray*}
l_{1} &=&u_{1}=0,l_{2}+u_{2}=1 \\
a &=&b_{1}=1,b_{2}=l_{2}, \\
d &=&1,e_{1}=0,e_{2}-u_{2}=1\Rightarrow e_{2}=1,u_{2}=0,b_{2}=l_{2}=1
\end{eqnarray*}%
Since $\alpha \left( x_{1};0,1,0,0\right) \equiv a+b_{1}+b_{2}+1\equiv 0,$
we obtain
\begin{equation*}
+B(gx_{2}\otimes gx_{2};GX_{1}X_{2},gx_{2})GX_{1}\otimes gx_{2}\otimes
x_{1}x_{2}.
\end{equation*}%
Third summand of the left side gives us%
\begin{eqnarray*}
l_{1}+u_{1} &=&1,l_{2}=u_{2}=0, \\
a &=&1,b_{2}=0,b_{1}-l_{1}=1\Rightarrow b_{1}=1,l_{1}=0,u_{1}=1 \\
d &=&e_{2}=1,e_{1}=u_{1}=1.
\end{eqnarray*}%
Since $\alpha \left( x_{2};0,0,1,0\right) \equiv e_{2}=1$ we get%
\begin{equation*}
+B(gx_{1}\otimes gx_{2};GX_{1},gx_{1}x_{2})GX_{1}\otimes gx_{2}\otimes
x_{1}x_{2}.
\end{equation*}%
The fourth summand of the left side gives us%
\begin{eqnarray*}
l_{1}+u_{1} &=&1,l_{2}+u_{2}=1 \\
a &=&1,b_{1}-l_{1}=1\Rightarrow b_{1}=1,l_{1}=0,u_{1}=1,b_{2}=l_{2}, \\
d &=&1,e_{1}=u_{1}=1,e_{2}-u_{2}=1\Rightarrow e_{2}=1,u_{2}=0,b_{2}=l_{2}=1.
\end{eqnarray*}%
Since $\alpha \left( 1_{H};0,1,1,0\right) \equiv e_{2}+a+b_{1}+b_{2}+1\equiv
1,$ we get%
\begin{equation*}
-B(x_{1}x_{2}\otimes gx_{2};GX_{1}X_{2},gx_{1}x_{2})GX_{1}\otimes
gx_{2}\otimes x_{1}x_{2}.
\end{equation*}%
Since there is no term to consider in the right side of the equality we
obtain%
\begin{gather*}
+B(gx_{1}\otimes gx_{2};GX_{1},gx_{1}x_{2})+B(gx_{2}\otimes
gx_{2};GX_{1}X_{2},gx_{2})+ \\
-B(x_{1}x_{2}\otimes gx_{2};GX_{1}X_{2},gx_{1}x_{2})+B(1_{H}\otimes
gx_{2};GX_{1},gx_{2})=0
\end{gather*}

which holds in view of the form of the elements.

\subsubsection{Case $GX_{2}\otimes gx_{1}x_{2}\otimes gx_{1}$}

We have to consider only the second and the fourth summand of the left side
of the equality. Second summand gives us

\begin{eqnarray*}
l_{1} &=&u_{1}=0,l_{2}=u_{2}=0 \\
a &=&b_{2}=1,b_{1}=0, \\
d &=&e_{2}=e_{1}=1.
\end{eqnarray*}

Since $\alpha \left( x_{1};0,0,0,0\right) \equiv a+b_{1}+b_{2}\equiv 0,$ we
get%
\begin{equation*}
B(gx_{2}\otimes gx_{2};GX_{2},gx_{1}x_{2})GX_{2}\otimes gx_{1}x_{2}\otimes
gx_{1}.
\end{equation*}

The fourth summand of the left side gives us%
\begin{eqnarray*}
l_{1}+u_{1} &=&1,l_{2}=u_{2}=0 \\
a &=&b_{2}=1,b_{1}=l_{1}, \\
d &=&e_{2}=1,e_{1}-u_{1}=1\Rightarrow e_{1}=1,u_{1}=0,b_{1}=l_{1}=1.
\end{eqnarray*}%
Since $\alpha \left( 1_{H};1,0,0,0\right) \equiv b_{2}\equiv 1,$ we obtain
\begin{equation*}
-B(x_{1}x_{2}\otimes gx_{2};GX_{1}X_{2},gx_{1}x_{2})GX_{2}\otimes
gx_{1}x_{2}\otimes gx_{1}.
\end{equation*}%
Since there is no term in the right side, we get%
\begin{equation*}
B(gx_{2}\otimes gx_{2};GX_{2},gx_{1}x_{2})-B(x_{1}x_{2}\otimes
gx_{2};GX_{1}X_{2},gx_{1}x_{2})=0
\end{equation*}%
which holds in view of the form of the elements.

\subsubsection{Case $G\otimes gx_{1}x_{2}\otimes x_{1}x_{2}$}

First summand of the left side of the equality gives us%
\begin{eqnarray*}
l_{1} &=&u_{1}=0,l_{2}=u_{2}=0 \\
a &=&1,b_{1}=b_{2}=0, \\
d &=&e_{1}=e_{2}=1.
\end{eqnarray*}

Since $\alpha \left( x_{1}x_{2};0,0,0,0\right) \equiv 0,$ we get
\begin{equation*}
B(1_{H}\otimes gx_{2};G,gx_{1}x_{2})G\otimes gx_{1}x_{2}\otimes x_{1}x_{2}.
\end{equation*}%
Second summand of the left side gives us%
\begin{eqnarray*}
l_{1} &=&u_{1}=0,l_{2}+u_{2}=1 \\
a &=&1,b_{1}=0,b_{2}=l_{2}, \\
d &=&e_{1}=1,e_{2}-u_{2}=1\Rightarrow e_{2}=1,u_{2}=0,b_{2}=l_{2}=1.
\end{eqnarray*}%
Since $\alpha \left( x_{1};0,1,0,0\right) \equiv a+b_{1}+b_{2}+1\equiv 1,$
we obtain
\begin{equation*}
-B(gx_{2}\otimes gx_{2};GX_{2},gx_{1}x_{2})G\otimes gx_{1}x_{2}\otimes
x_{1}x_{2}.
\end{equation*}%
Third summand of the left side gives us%
\begin{eqnarray*}
l_{1}+u_{1} &=&1,l_{2}=u_{2}=0, \\
a &=&1,b_{2}=0,b_{1}=l_{1}, \\
d &=&e_{2}=1,e_{1}-u_{1}=1\Rightarrow e_{1}=1,u_{1}=0,b_{1}=l_{1}=1.
\end{eqnarray*}%
Since $\alpha \left( x_{2};1,0,0,0\right) \equiv a+b_{1}\equiv 0$ we get%
\begin{equation*}
-B(gx_{1}\otimes gx_{2};GX_{1},gx_{1}x_{2})G\otimes gx_{1}x_{2}\otimes
x_{1}x_{2}.
\end{equation*}%
The fourth summand of the left side gives us%
\begin{eqnarray*}
l_{1}+u_{1} &=&1,l_{2}+u_{2}=1 \\
a &=&1,b_{1}=l_{1},b_{2}=l_{2}, \\
d &=&1,e_{1}-u_{1}=1\Rightarrow e_{1}=1,u_{1}=0,b_{1}=l_{1}=1 \\
e_{2}-u_{2} &=&1\Rightarrow e_{2}=1,u_{2}=0,b_{2}=l_{2}=1.
\end{eqnarray*}%
Since $\alpha \left( 1_{H};1,1,0,0\right) \equiv 1+b_{2}\equiv 0$ $,$ we get%
\begin{equation*}
+B(x_{1}x_{2}\otimes gx_{2};GX_{1}X_{2},gx_{1}x_{2})G\otimes
gx_{1}x_{2}\otimes x_{1}x_{2}.
\end{equation*}%
Since there is no term to consider in the right side of the equality we
obtain%
\begin{gather*}
-B(gx_{1}\otimes gx_{2};GX_{1},gx_{1}x_{2})-B(gx_{2}\otimes
gx_{2};GX_{2},gx_{1}x_{2})+ \\
+B(x_{1}x_{2}\otimes gx_{2};GX_{1}X_{2},gx_{1}x_{2})+B(1_{H}\otimes
gx_{2};G,gx_{1}x_{2})=0
\end{gather*}

which holds in view of the form of the elements.

\section{$B\left( x_{1}x_{2}\otimes gx_{1}x_{2}\right) $}

By using $\ref{simplgxx}$%
\begin{eqnarray}
B(x_{1}x_{2}\otimes gx_{1}x_{2}) &=&(1_{A}\otimes g)B(gx_{1}x_{2}\otimes
1_{H})(1_{A}\otimes gx_{1}x_{2})  \label{form x1x2otgx1gx2} \\
&&-(1_{A}\otimes x_{2})B(gx_{1}x_{2}\otimes 1_{H})(1_{A}\otimes x_{1})
\notag \\
&&+(1_{A}\otimes x_{1})B(gx_{1}x_{2}\otimes 1_{H})(1_{A}\otimes x_{2})
\notag \\
&&+(1_{A}\otimes gx_{1}x_{2})B(gx_{1}x_{2}\otimes 1_{H})(1_{A}\otimes g)
\notag
\end{eqnarray}%
and we obtain%
\begin{eqnarray*}
&&B(x_{1}x_{2}\otimes gx_{1}x_{2}) \\
&=&4B(gx_{1}x_{2}\otimes 1_{H};1_{A},g)1_{A}\otimes gx_{1}x_{2}+ \\
&&+\left[
\begin{array}{c}
4B(g\otimes 1_{H};1_{A},g)+4B(x_{2}\otimes \ 1_{H};1_{A},gx_{2})+ \\
4B(x_{1}\otimes 1_{H};1_{A},gx_{1})+4B(gx_{1}x_{2}\otimes
1_{H};1_{A},gx_{1}x_{2})%
\end{array}%
\right] X_{1}X_{2}\otimes gx_{1}x_{2}+ \\
&&+\left[ 4B(x_{2}\otimes 1_{H};G,g)-4B(gx_{1}x_{2}\otimes 1_{H};G,gx_{1})%
\right] GX_{1}\otimes gx_{1}x_{2}+ \\
&&+\left[ -4B(x_{1}\otimes 1_{H};G,g)-4B(gx_{1}x_{2}\otimes 1_{H};G,gx_{2})%
\right] GX_{2}\otimes gx_{1}x_{2}+
\end{eqnarray*}%
We write the Casimir condition for $B(x_{1}x_{2}\otimes
gx_{1}x_{2}).1_{A}\otimes g\otimes gx_{1}x_{2}$%
\begin{eqnarray*}
&&\sum_{a,b_{1},b_{2},d,e_{1},e_{2}=0}^{1}\sum_{l_{1}=0}^{b_{1}}%
\sum_{l_{2}=0}^{b_{2}}\sum_{u_{1}=0}^{e_{1}}\sum_{u_{2}=0}^{e_{2}}\left(
-1\right) ^{\alpha \left( x_{1}x_{2};l_{1},l_{2},u_{1},u_{2}\right) } \\
&&B(1_{H}\otimes
gx_{1}x_{2};G^{a}X_{1}^{b_{1}}X_{2}^{b_{2}},g^{d}x_{1}^{e_{1}}x_{2}^{e_{2}})
\\
&&G^{a}X_{1}^{b_{1}-l_{1}}X_{2}^{b_{2}-l_{2}}\otimes
g^{d}x_{1}^{e_{1}-u_{1}}x_{2}^{e_{2}-u_{2}}\otimes
g^{a+b_{1}+b_{2}+l_{1}+l_{2}+d+e_{1}+e_{2}+u_{1}+u_{2}}x_{1}^{l_{1}+u_{1}+1}x_{2}^{l_{2}+u_{2}+1}+
\\
&&+\sum_{a,b_{1},b_{2},d,e_{1},e_{2}=0}^{1}\sum_{l_{1}=0}^{b_{1}}%
\sum_{l_{2}=0}^{b_{2}}\sum_{u_{1}=0}^{e_{1}}\sum_{u_{2}=0}^{e_{2}}\left(
-1\right) ^{\alpha \left( x_{1};l_{1},l_{2},u_{1},u_{2}\right) } \\
&&B(gx_{2}\otimes
gx_{1}x_{2};G^{a}X_{1}^{b_{1}}X_{2}^{b_{2}},g^{d}x_{1}^{e_{1}}x_{2}^{e_{2}})
\\
&&G^{a}X_{1}^{b_{1}-l_{1}}X_{2}^{b_{2}-l_{2}}\otimes
g^{d}x_{1}^{e_{1}-u_{1}}x_{2}^{e_{2}-u_{2}}\otimes
g^{a+b_{1}+b_{2}+l_{1}+l_{2}+d+e_{1}+e_{2}+u_{1}+u_{2}}x_{1}^{l_{1}+u_{1}+1}x_{2}^{l_{2}+u_{2}}+
\\
&&+\left( -1\right)
\sum_{a,b_{1},b_{2},d,e_{1},e_{2}=0}^{1}\sum_{l_{1}=0}^{b_{1}}%
\sum_{l_{2}=0}^{b_{2}}\sum_{u_{1}=0}^{e_{1}}\sum_{u_{2}=0}^{e_{2}}\left(
-1\right) ^{\alpha \left( x_{2};l_{1},l_{2},u_{1},u_{2}\right) } \\
&&B(gx_{1}\otimes
gx_{1}x_{2};G^{a}X_{1}^{b_{1}}X_{2}^{b_{2}},g^{d}x_{1}^{e_{1}}x_{2}^{e_{2}})
\\
&&G^{a}X_{1}^{b_{1}-l_{1}}X_{2}^{b_{2}-l_{2}}\otimes
g^{d}x_{1}^{e_{1}-u_{1}}x_{2}^{e_{2}-u_{2}}\otimes
g^{a+b_{1}+b_{2}+l_{1}+l_{2}+d+e_{1}+e_{2}+u_{1}+u_{2}}x_{1}^{l_{1}+u_{1}}x_{2}^{l_{2}+u_{2}+1}+
\\
&&+\sum_{a,b_{1},b_{2},d,e_{1},e_{2}=0}^{1}\sum_{l_{1}=0}^{b_{1}}%
\sum_{l_{2}=0}^{b_{2}}\sum_{u_{1}=0}^{e_{1}}\sum_{u_{2}=0}^{e_{2}}\left(
-1\right) ^{\alpha \left( 1_{H};l_{1},l_{2},u_{1},u_{2}\right) } \\
&&B(x_{1}x_{2}\otimes
gx_{1}x_{2};G^{a}X_{1}^{b_{1}}X_{2}^{b_{2}},g^{d}x_{1}^{e_{1}}x_{2}^{e_{2}})
\\
&&G^{a}X_{1}^{b_{1}-l_{1}}X_{2}^{b_{2}-l_{2}}\otimes
g^{d}x_{1}^{e_{1}-u_{1}}x_{2}^{e_{2}-u_{2}}\otimes
g^{a+b_{1}+b_{2}+l_{1}+l_{2}+d+e_{1}+e_{2}+u_{1}+u_{2}}x_{1}^{l_{1}+u_{1}}x_{2}^{l_{2}+u_{2}}+
\\
&=&B(x_{1}x_{2}\otimes gx_{1}x_{2})\otimes g \\
&&+\text{ }B(x_{1}x_{2}\otimes gx_{1})\otimes x_{2}\text{ } \\
&&-B(x_{1}x_{2}\otimes gx_{2})\otimes x_{1} \\
&&+B(x_{1}x_{2}\otimes g)\otimes gx_{1}x_{2}
\end{eqnarray*}

\subsection{$B(x_{1}x_{2}\otimes gx_{1}x_{2};1_{A},gx_{1}x_{2})$}

We get%
\begin{eqnarray*}
a &=&b_{1}=b_{2}=0 \\
d &=&e_{1}=e_{2}=1
\end{eqnarray*}%
and we obtain%
\begin{eqnarray*}
&&\left( -1\right) ^{\alpha \left( 1_{H};0,0,0,0\right) }B(x_{1}x_{2}\otimes
gx_{1}x_{2};1_{A},gx_{1}x_{2})1_{A}\otimes gx_{1}x_{2}\otimes g+ \\
&&+\left( -1\right) ^{\alpha \left( 1_{H};0,0,0,1\right)
}B(x_{1}x_{2}\otimes gx_{1}x_{2};1_{A},gx_{1}x_{2})1_{A}\otimes
gx_{1}\otimes x_{2}+ \\
&&+\left( -1\right) ^{\alpha \left( 1_{H};0,0,1,0\right)
}B(x_{1}x_{2}\otimes gx_{1}x_{2};1_{A},gx_{1}x_{2})1_{A}\otimes
gx_{2}\otimes x_{1}+ \\
&&+\left( -1\right) ^{\alpha \left( 1_{H};0,0,1,1\right)
}B(x_{1}x_{2}\otimes gx_{1}x_{2};1_{A},gx_{1}x_{2})1_{A}\otimes g\otimes
gx_{1}x_{2}
\end{eqnarray*}

\subsubsection{Case $1_{A}\otimes gx_{1}\otimes x_{2}$}

We have to consider only the third and the fourth summand of the left side
of the equality. Third summand gives us%
\begin{eqnarray*}
l_{1} &=&u_{1}=0,l_{2}=u_{2}=0, \\
a &=&b_{1}=b_{2}=0, \\
d &=&e_{1}=1,e_{2}=0.
\end{eqnarray*}%
Since $\alpha \left( x_{2};0,0,0,0\right) \equiv a+b_{1}+b_{2}=0,$ we get
\begin{equation*}
-B(gx_{1}\otimes gx_{1}x_{2};1_{A},gx_{1})1_{A}\otimes gx_{1}\otimes x_{2}.
\end{equation*}%
Fourth summand gives us%
\begin{eqnarray*}
l_{1} &=&u_{1}=0,l_{2}+u_{2}=1, \\
a &=&b_{1}=0,b_{2}=l_{2}, \\
d &=&e_{1}=1,e_{2}=u_{2}.
\end{eqnarray*}%
Since $\alpha \left( 1_{H};0,0,0,1\right) \equiv a+b_{1}+b_{2}\equiv 0$ and $%
\alpha \left( 1_{H};0,1,0,0\right) \equiv 0$ we get
\begin{equation*}
\left[ B(x_{1}x_{2}\otimes
gx_{1}x_{2};1_{A},gx_{1}x_{2})+B(x_{1}x_{2}\otimes gx_{1}x_{2};X_{2},gx_{1})%
\right] 1_{A}\otimes gx_{1}\otimes x_{2}.
\end{equation*}%
By considering also the right side we get%
\begin{gather*}
B(x_{1}x_{2}\otimes gx_{1}x_{2};1_{A},gx_{1}x_{2})+B(x_{1}x_{2}\otimes
gx_{1}x_{2};X_{2},gx_{1}) \\
-\text{ }B(x_{1}x_{2}\otimes gx_{1}1_{A},gx_{1})-B(gx_{1}\otimes
gx_{1}x_{2};1_{A},gx_{1})=0
\end{gather*}%
which holds in view of the form of the elements.

\subsubsection{Case $1_{A}\otimes gx_{2}\otimes x_{1}$}

We have to consider only the second and the fourth summand of the left side
of the equality.

Second summand gives us%
\begin{eqnarray*}
l_{1} &=&u_{1}=l_{2}=u_{2}=0, \\
a &=&b_{1}=b_{2}=0, \\
d &=&e_{2}=1,e_{1}=0.
\end{eqnarray*}%
Since $\alpha \left( x_{1};0,0,0,0\right) \equiv a+b_{1}+b_{2}\equiv 0$ we
obtain
\begin{equation*}
B(gx_{2}\otimes gx_{1}x_{2};1_{A},gx_{2})1_{A}\otimes gx_{2}\otimes x_{1}.
\end{equation*}%
Forth summand gives us%
\begin{eqnarray*}
l_{1}+u_{1} &=&1,l_{2}=u_{2}=0 \\
a &=&b_{2}=0,b_{1}=l_{1}, \\
d &=&e_{2}=1,e_{1}=u_{1}.
\end{eqnarray*}%
Since $\alpha \left( 1_{H};0,0,1,0\right) \equiv e_{2}+\left(
a+b_{1}+b_{2}\right) \equiv 1$ and $\alpha \left( 1_{H};1,0,0,0\right)
\equiv b_{2}=0,$ we obtain
\begin{equation*}
\left[ -B(x_{1}x_{2}\otimes
gx_{1}x_{2};1_{A},gx_{1}x_{2})+B(x_{1}x_{2}\otimes gx_{1}x_{2};X_{1},gx_{2})%
\right] 1_{A}\otimes gx_{2}\otimes x_{1}.
\end{equation*}%
By considering the right side, we get%
\begin{gather*}
-B(x_{1}x_{2}\otimes gx_{1}x_{2};1_{A},gx_{1}x_{2})+B(x_{1}x_{2}\otimes
gx_{1}x_{2};X_{1},gx_{2})+ \\
+B(gx_{2}\otimes gx_{1}x_{2};1_{A},gx_{2})+B(x_{1}x_{2}\otimes
gx_{2};1_{A},gx_{2})=0
\end{gather*}%
which holds in view of the form of the elements.

\subsubsection{Case $1_{A}\otimes g\otimes gx_{1}x_{2}$}

First summand of the left side of the equality gives us%
\begin{eqnarray*}
l_{1} &=&u_{1}=0,l_{2}=u_{2}=0 \\
a &=&b_{1}=b_{2}=0, \\
d &=&1,e_{1}=e_{2}=0,.
\end{eqnarray*}%
Since $\alpha \left( x_{1}x_{2};0,0,0,0\right) \equiv 0,$ we get%
\begin{equation*}
B(1_{H}\otimes gx_{1}x_{2};1_{A},g)1_{A}\otimes g\otimes gx_{1}x_{2}.
\end{equation*}%
Second summand of the left side of the equation gives us%
\begin{eqnarray*}
l_{1} &=&u_{1}=0,l_{2}+u_{2}=1 \\
a &=&b_{1}=0,b_{2}=l_{2}, \\
d &=&1,e_{1}=0,e_{2}=u_{2}.
\end{eqnarray*}%
Since $\alpha \left( x_{1};0,0,0,1\right) \equiv 1$ and $\alpha \left(
x_{1};0,1,0,0\right) \equiv a+b_{1}+b_{2}+1\equiv 0,$ we get%
\begin{equation*}
\left[ -B(gx_{2}\otimes gx_{1}x_{2};1_{A},gx_{2})+B(gx_{2}\otimes
gx_{1}x_{2};X_{2},g)\right] 1_{A}\otimes g\otimes gx_{1}x_{2}.
\end{equation*}

Third summand of the left side of the equation gives us%
\begin{eqnarray*}
l_{1}+u_{1} &=&1,l_{2}=u_{2}=0, \\
a &=&b_{2}=0,b_{1}=l_{1} \\
d &=&1,e_{2}=0,e_{1}=u_{1}.
\end{eqnarray*}%
Since $\alpha \left( x_{2};0,0,1,0\right) \equiv e_{2}=0$ and $\alpha \left(
x_{2};1,0,0,0\right) \equiv a+b_{1}=1,$ we get%
\begin{equation*}
\left[ -B(gx_{1}\otimes gx_{1}x_{2};1_{A},gx_{1})+B(gx_{1}\otimes
gx_{1}x_{2};X_{1},g)\right] 1_{A}\otimes g\otimes gx_{1}x_{2}.
\end{equation*}%
Fourth summand of the left side of the equation gives us%
\begin{eqnarray*}
l_{1}+u_{1} &=&1,l_{2}+u_{2}=1 \\
a &=&0,b_{1}=l_{1},b_{2}=l_{2}, \\
d &=&1,e_{1}=u_{1},e_{2}=u_{2}.
\end{eqnarray*}%
Since
\begin{eqnarray*}
\alpha \left( 1_{H};0,0,1,1\right) &=&1+e_{2}\text{, } \\
\alpha \left( 1_{H};0,1,1,0\right) &=&e_{2}+a+b_{1}+b_{2}+1, \\
\alpha \left( 1_{H};1,0,0,1\right) &\equiv &a+b_{1}\text{and } \\
\alpha \left( 1_{H};1,1,0,0\right) &\equiv &1+b_{2}
\end{eqnarray*}%
we get%
\begin{equation*}
\left[
\begin{array}{c}
B\left( x_{1}x_{2}\otimes gx_{1}x_{2};1_{A},gx_{1}x_{2}\right) +B\left(
x_{1}x_{2}\otimes gx_{1}x_{2};X_{2},gx_{1}\right) + \\
-B\left( x_{1}x_{2}\otimes gx_{1}x_{2};X_{1},gx_{2}\right) +B\left(
x_{1}x_{2}\otimes gx_{1}x_{2};X_{1}X_{2},g\right)%
\end{array}%
\right] 1_{A}\otimes g\otimes gx_{1}x_{2}.
\end{equation*}%
By considering also the right side, we get%
\begin{gather*}
B(1_{H}\otimes gx_{1}x_{2};1_{A},g)-B(gx_{2}\otimes
gx_{1}x_{2};1_{A},gx_{2})+B(gx_{2}\otimes gx_{1}x_{2};X_{2},g) \\
-B(gx_{1}\otimes gx_{1}x_{2};1_{A},gx_{1})+B(gx_{1}\otimes
gx_{1}x_{2};X_{1},g)+ \\
B\left( x_{1}x_{2}\otimes gx_{1}x_{2};1_{A},gx_{1}x_{2}\right) +B\left(
x_{1}x_{2}\otimes gx_{1}x_{2};X_{2},gx_{1}\right) + \\
-B\left( x_{1}x_{2}\otimes gx_{1}x_{2};X_{1},gx_{2}\right) +B\left(
x_{1}x_{2}\otimes gx_{1}x_{2};X_{1}X_{2},g\right) -B(x_{1}x_{2}\otimes
g;1_{A},g)=0
\end{gather*}%
which holds in view of the form of the elements.

\subsection{$B(x_{1}x_{2}\otimes gx_{1}x_{2};X_{1}X_{2},gx_{1}x_{2})$}

We have%
\begin{equation*}
a=0,b_{1}=b_{2}=1,d=e_{1}=e_{2}=1
\end{equation*}%
and we deduce%
\begin{gather*}
\left( -1\right) ^{\alpha \left( 1_{H};0,0,1,0\right) }B(x_{1}x_{2}\otimes
gx_{1}x_{2};X_{1}X_{2},gx_{1}x_{2})X_{1}X_{2}\otimes gx_{2}\otimes x_{1}+ \\
\left( -1\right) ^{\alpha \left( 1_{H};0,1,1,0\right) }B(x_{1}x_{2}\otimes
gx_{1}x_{2};X_{1}X_{2},gx_{1}x_{2})X_{1}\otimes gx_{2}\otimes gx_{1}x_{2} \\
+\left( -1\right) ^{\alpha \left( 1_{H};1,0,1,0\right) }B(x_{1}x_{2}\otimes
gx_{1}x_{2};X_{1}X_{2},gx_{1}x_{2})X_{1}^{1-l_{1}}X_{2}^{1-l_{2}}\otimes
gx_{2}\otimes g^{+l_{1}+l_{2}}x_{1}^{1+1}x_{2}^{l_{2}}=0 \\
+\left( -1\right) ^{\alpha \left( 1_{H};1,1,1,0\right) }B(x_{1}x_{2}\otimes
gx_{1}x_{2};X_{1}X_{2},gx_{1}x_{2})X_{1}^{1-l_{1}}X_{2}^{1-l_{2}}\otimes
gx_{2}\otimes g^{+l_{1}+l_{2}}x_{1}^{1+1}x_{2}^{l_{2}}=0 \\
+\left( -1\right) ^{\alpha \left( 1_{H};0,0,0,0\right) }B(x_{1}x_{2}\otimes
gx_{1}x_{2};X_{1}X_{2},gx_{1}x_{2})X_{1}X_{2}\otimes gx_{1}x_{2}\otimes g+ \\
+\left( -1\right) ^{\alpha \left( 1_{H};0,1,0,0\right) }B(x_{1}x_{2}\otimes
gx_{1}x_{2};X_{1}X_{2},gx_{1}x_{2})X_{1}\otimes gx_{1}x_{2}\otimes x_{2}+ \\
+\left( -1\right) ^{\alpha \left( 1_{H};1,0,0,0\right) }B(x_{1}x_{2}\otimes
gx_{1}x_{2};X_{1}X_{2},gx_{1}x_{2})X_{2}\otimes gx_{1}x_{2}\otimes x_{1}+ \\
+\left( -1\right) ^{\alpha \left( 1_{H};1,1,0,0\right) }B(x_{1}x_{2}\otimes
gx_{1}x_{2};X_{1}X_{2},gx_{1}x_{2})1_{A}\otimes gx_{1}x_{2}\otimes
gx_{1}x_{2}+ \\
+\left( -1\right) ^{\alpha \left( 1_{H};0,0,0,1\right) }B(x_{1}x_{2}\otimes
gx_{1}x_{2};X_{1}X_{2},gx_{1}x_{2})X_{1}X_{2}\otimes gx_{1}\otimes x_{2}+ \\
+\left( -1\right) ^{\alpha \left( 1_{H};0,1,0,1\right) }B(x_{1}x_{2}\otimes
gx_{1}x_{2};X_{1}X_{2},gx_{1}x_{2})X_{1}^{1-l_{1}}X_{2}^{1-l_{2}}\otimes
gx_{1}\otimes g^{1+l_{1}+l_{2}+1}x_{1}^{l_{1}}x_{2}^{1+1}=0 \\
+\left( -1\right) ^{\alpha \left( 1_{H};1,0,0,1\right) }B(x_{1}x_{2}\otimes
gx_{1}x_{2};X_{1}X_{2},gx_{1}x_{2})X_{2}\otimes gx_{1}\otimes gx_{1}x_{2}+ \\
+\left( -1\right) ^{\alpha \left( 1_{H};1,1,0,1\right) }B(x_{1}x_{2}\otimes
gx_{1}x_{2};X_{1}X_{2},gx_{1}x_{2})X_{1}^{1-l_{1}}X_{2}^{1-l_{2}}\otimes
gx_{1}\otimes g^{1+l_{1}+l_{2}+1}x_{1}^{l_{1}}x_{2}^{1+1}=0 \\
+\left( -1\right) ^{\alpha \left( 1_{H};0,0,1,1\right) }B(x_{1}x_{2}\otimes
gx_{1}x_{2};X_{1}X_{2},gx_{1}x_{2})X_{1}X_{2}\otimes g\otimes gx_{1}x_{2}+ \\
+\left( -1\right) ^{\alpha \left( 1_{H};0,1,1,1\right) }B(x_{1}x_{2}\otimes
gx_{1}x_{2};X_{1}X_{2},gx_{1}x_{2})X_{2}^{1-l_{2}}\otimes g\otimes
g^{1+l_{1}+l_{2}}x_{1}^{l_{1}+1}x_{2}^{1+1}=0 \\
+\left( -1\right) ^{\alpha \left( 1_{H};1,0,1,1\right) }B(x_{1}x_{2}\otimes
gx_{1}x_{2};X_{1}X_{2},gx_{1}x_{2})X_{2}^{1-l_{2}}\otimes g\otimes
g^{1+l_{1}+l_{2}}x_{1}^{1+1}x_{2}^{l_{2}+1}=0 \\
+\left( -1\right) ^{\alpha \left( 1_{H};1,1,1,1\right) }B(x_{1}x_{2}\otimes
gx_{1}x_{2};X_{1}X_{2},gx_{1}x_{2})X_{2}^{1-l_{2}}\otimes g\otimes
g^{1+l_{1}+l_{2}}x_{1}^{1+1}x_{2}^{1+1}=0
\end{gather*}

\subsubsection{$X_{1}X_{2}\otimes gx_{2}\otimes x_{1}$}

We have to consider only the second and the fourth summand of the left side
of the equality.

Second summand gives us%
\begin{eqnarray*}
l_{1} &=&u_{1}=l_{2}=u_{2}=0, \\
a &=&0,b_{1}=b_{2}=1 \\
d &=&e_{2}=1,e_{1}=0.
\end{eqnarray*}%
Since $\alpha \left( x_{1};0,0,0,0\right) \equiv a+b_{1}+b_{2}\equiv 0$ we
obtain
\begin{equation*}
B(gx_{2}\otimes gx_{1}x_{2};X_{1}X_{2},gx_{2})X_{1}X_{2}\otimes
gx_{2}\otimes x_{1}.
\end{equation*}%
Forth summand gives us%
\begin{eqnarray*}
l_{1}+u_{1} &=&1,l_{2}=u_{2}=0 \\
a &=&0,b_{2}=1,b_{1}-l_{1}=1\Rightarrow b_{1}=1,l_{1}=0,u_{1}=1, \\
d &=&e_{2}=1,e_{1}=u_{1}=1.
\end{eqnarray*}%
Since $\alpha \left( 1_{H};0,0,1,0\right) \equiv e_{2}+\left(
a+b_{1}+b_{2}\right) \equiv 1,$ we obtain
\begin{equation*}
-B(x_{1}x_{2}\otimes gx_{1}x_{2};X_{1}X_{1},gx_{1}x_{2})X_{1}X_{2}\otimes
gx_{2}\otimes x_{1}.
\end{equation*}%
By considering the right side, we get,
\begin{gather*}
B(x_{1}x_{2}\otimes gx_{2};X_{1}X_{2},gx_{2})-B(x_{1}x_{2}\otimes
gx_{1}x_{2};X_{1}X_{1},gx_{1}x_{2})+ \\
+B(gx_{2}\otimes gx_{1}x_{2};X_{1}X_{2},gx_{2})=0
\end{gather*}%
which holds in view of the form of the elements.

\subsubsection{Case $X_{1}\otimes gx_{2}\otimes gx_{1}x_{2}$}

First summand of the left side of the equality gives us%
\begin{eqnarray*}
l_{1} &=&u_{1}=0,l_{2}=u_{2}=0 \\
a &=&b_{2}=0,b_{1}=1 \\
d &=&e_{2}=1,e_{1}=0,.
\end{eqnarray*}%
Since $\alpha \left( x_{1}x_{2};0,0,0,0\right) \equiv 0,$ we get%
\begin{equation*}
B(1_{H}\otimes gx_{1}x_{2};X_{1},gx_{2})X_{1}\otimes gx_{2}\otimes
gx_{1}x_{2}.
\end{equation*}%
Second summand of the left side of the equation gives us%
\begin{eqnarray*}
l_{1} &=&u_{1}=0,l_{2}+u_{2}=1 \\
a &=&0,b_{1}=1,b_{2}=l_{2}, \\
d &=&1,e_{1}=0,e_{2}-u_{2}=1\Rightarrow e_{2}=1,u_{2}=0,b_{2}=l_{2}=1.
\end{eqnarray*}%
Since $\alpha \left( x_{1};0,1,0,0\right) \equiv a+b_{1}+b_{2}+1\equiv 1,$
we get%
\begin{equation*}
-B(gx_{2}\otimes gx_{1}x_{2};X_{1}X_{2},gx_{2})X_{1}\otimes gx_{2}\otimes
gx_{1}x_{2}.
\end{equation*}

Third summand of the left side of the equation gives us%
\begin{eqnarray*}
l_{1}+u_{1} &=&1,l_{2}=u_{2}=0, \\
a &=&b_{2}=0,b_{1}-l_{1}=1\Rightarrow b_{1}=1,l_{1}=0,u_{1}=1 \\
d &=&1,e_{2}=1,e_{1}=u_{1}=1.
\end{eqnarray*}%
Since $\alpha \left( x_{2};0,0,1,0\right) \equiv e_{2}=1$ $,$ we get%
\begin{equation*}
+B(gx_{1}\otimes gx_{1}x_{2};X_{1},gx_{1}x_{2})X_{1}\otimes gx_{2}\otimes
gx_{1}x_{2}.
\end{equation*}%
Fourth summand of the left side of the equation gives us%
\begin{eqnarray*}
l_{1}+u_{1} &=&1,l_{2}+u_{2}=1 \\
a &=&0,b_{2}=l_{2}, \\
b_{1}-l_{1} &=&1\Rightarrow b_{1}=1,l_{1}=0,u_{1}=1 \\
d &=&1,e_{1}=u_{1}=1, \\
e_{2}-u_{2} &=&1\Rightarrow e_{2}=1,u_{2}=0,b_{2}=l_{2}=1.
\end{eqnarray*}%
Since $\alpha \left( 1_{H};0,1,1,0\right) =e_{2}+a+b_{1}+b_{2}+1\equiv 0,$
we get%
\begin{equation*}
B\left( x_{1}x_{2}\otimes gx_{1}x_{2};X_{1}X_{2},gx_{1}x_{2}\right)
X_{1}\otimes gx_{2}\otimes gx_{1}x_{2}.
\end{equation*}%
By considering also the right side, we get%
\begin{gather*}
B(1_{H}\otimes gx_{1}x_{2};X_{1},gx_{2})-B(gx_{2}\otimes
gx_{1}x_{2};X_{1}X_{2},gx_{2})+ \\
+B(gx_{1}\otimes gx_{1}x_{2};X_{1},gx_{1}x_{2})+B\left( x_{1}x_{2}\otimes
gx_{1}x_{2};X_{1}X_{2},gx_{1}x_{2}\right) \\
-B(x_{1}x_{2}\otimes g;X_{1},gx_{2})=0
\end{gather*}%
which holds in view of the form of the elements.

\subsubsection{Case $X_{1}\otimes gx_{1}x_{2}\otimes x_{2}$}

We have to consider only the third and the fourth summand of the left side
of the equality. Third summand gives us%
\begin{eqnarray*}
l_{1} &=&u_{1}=0,l_{2}=u_{2}=0, \\
a &=&b_{2}=0,b_{1}=1 \\
d &=&e_{1}=e_{2}=1.
\end{eqnarray*}%
Since $\alpha \left( x_{2};0,0,0,0\right) \equiv a+b_{1}+b_{2}=1,$ we get
\begin{equation*}
B(gx_{1}\otimes gx_{1}x_{2};X_{1},gx_{1}x_{2})X_{1}\otimes
gx_{1}x_{2}\otimes x_{2}.
\end{equation*}%
Fourth summand gives us%
\begin{eqnarray*}
l_{1} &=&u_{1}=0,l_{2}+u_{2}=1, \\
a &=&0,b_{1}=1,b_{2}=l_{2}, \\
d &=&e_{1}=1,e_{2}-u_{2}=1\Rightarrow e_{2}=1,u_{2}=0,b_{2}=l_{2}=1.
\end{eqnarray*}%
Since $\alpha \left( 1_{H};0,1,0,0\right) \equiv 0$ we get
\begin{equation*}
+B(x_{1}x_{2}\otimes gx_{1}x_{2};X_{1}X_{2},gx_{1}x_{2})X_{1}\otimes
gx_{1}x_{2}\otimes x_{2}.
\end{equation*}%
By considering also the right side we get%
\begin{gather*}
+B(x_{1}x_{2}\otimes gx_{1}x_{2};X_{1}X_{2},gx_{1}x_{2})-\text{ }%
B(x_{1}x_{2}\otimes gx_{1};X_{1},gx_{1}x_{2})+ \\
+B(gx_{1}\otimes gx_{1}x_{2};X_{1},gx_{1}x_{2})=0
\end{gather*}%
which holds in view of the form of the elements.

\subsubsection{Case $X_{2}\otimes gx_{1}x_{2}\otimes x_{1}$}

We have to consider only the second and the fourth summand of the left side
of the equality.

Second summand gives us%
\begin{eqnarray*}
l_{1} &=&u_{1}=l_{2}=u_{2}=0, \\
a &=&b_{1}=0,b_{2}=1 \\
d &=&e_{1}=e_{2}=1.
\end{eqnarray*}%
Since $\alpha \left( x_{1};0,0,0,0\right) \equiv a+b_{1}+b_{2}\equiv 1$ we
obtain
\begin{equation*}
-B(gx_{2}\otimes gx_{1}x_{2};X_{2},gx_{1}x_{2})X_{2}\otimes
gx_{1}x_{2}\otimes x_{1}.
\end{equation*}%
Forth summand gives us%
\begin{eqnarray*}
l_{1}+u_{1} &=&1,l_{2}=u_{2}=0 \\
a &=&0,b_{1}=l_{1},b_{2}=1, \\
d &=&e_{2}=1,e_{1}-u_{1}=1\Rightarrow e_{1}=1,u_{1}=0,b_{1}=l_{1}=1.
\end{eqnarray*}%
Since $\alpha \left( 1_{H};1,0,0,0\right) \equiv b_{2}\equiv 1,$ we obtain
\begin{equation*}
-B(x_{1}x_{2}\otimes gx_{1}x_{2};X_{1}X_{2},gx_{1}x_{2})X_{2}\otimes
gx_{1}x_{2}\otimes x_{1}.
\end{equation*}%
By considering the right side, we get,
\begin{gather*}
B(x_{1}x_{2}\otimes gx_{2};X_{2},gx_{1}x_{2})-B(x_{1}x_{2}\otimes
gx_{1}x_{2};X_{1}X_{2},gx_{1}x_{2})+ \\
-B(gx_{2}\otimes gx_{1}x_{2};X_{2},gx_{1}x_{2})=0
\end{gather*}%
which holds in view of the form of the elements.

\subsubsection{Case $1_{A}\otimes gx_{1}x_{2}\otimes gx_{1}x_{2}$}

First summand of the left side of the equality gives us%
\begin{eqnarray*}
l_{1} &=&u_{1}=0,l_{2}=u_{2}=0 \\
a &=&b_{1}=b_{2}=0, \\
d &=&e_{1}=e_{2}=1.
\end{eqnarray*}%
Since $\alpha \left( x_{1}x_{2};0,0,0,0\right) \equiv 0,$ we get%
\begin{equation*}
B(1_{H}\otimes gx_{1}x_{2};1_{A},gx_{1}x_{2})1_{A}\otimes gx_{1}x_{2}\otimes
gx_{1}x_{2}.
\end{equation*}%
Second summand of the left side of the equation gives us%
\begin{eqnarray*}
l_{1} &=&u_{1}=0,l_{2}+u_{2}=1 \\
a &=&b_{1}=0,b_{2}=l_{2}, \\
d &=&e_{1}=1,e_{2}-u_{2}=1\Rightarrow e_{2}=1,u_{2}=0,b_{2}=l_{2}=1.
\end{eqnarray*}%
Since $\alpha \left( x_{1};0,1,0,0\right) \equiv a+b_{1}+b_{2}+1\equiv 0,$
we get%
\begin{equation*}
+B(gx_{2}\otimes gx_{1}x_{2};X_{2},gx_{1}x_{2})1_{A}\otimes
gx_{1}x_{2}\otimes gx_{1}x_{2}.
\end{equation*}

Third summand of the left side of the equation gives us%
\begin{eqnarray*}
l_{1}+u_{1} &=&1,l_{2}=u_{2}=0, \\
a &=&b_{2}=0,b_{1}=l_{1}=1 \\
d &=&e_{2}=1,e_{1}-u_{1}=1\Rightarrow e_{1}=1,u_{1}=0,b_{1}=l_{1}=1.
\end{eqnarray*}%
Since $\alpha \left( x_{2};1,0,0,0\right) \equiv a+b_{1}\equiv 1$ $,$ we get%
\begin{equation*}
+B(gx_{1}\otimes gx_{1}x_{2};X_{1},gx_{1}x_{2})1_{A}\otimes
gx_{1}x_{2}\otimes gx_{1}x_{2}.
\end{equation*}%
Fourth summand of the left side of the equation gives us%
\begin{eqnarray*}
l_{1}+u_{1} &=&1,l_{2}+u_{2}=1 \\
a &=&0,b_{2}=l_{2},b_{1}=l_{1} \\
d &=&1,e_{1}-u_{1}=1\Rightarrow e_{1}=1,u_{1}=0,b_{1}=l_{1}=1, \\
e_{2}-u_{2} &=&1\Rightarrow e_{2}=1,u_{2}=0,b_{2}=l_{2}=1.
\end{eqnarray*}%
Since $\alpha \left( 1_{H};1,1,0,0\right) \equiv 1+b_{2}\equiv 0,$ we get%
\begin{equation*}
B\left( x_{1}x_{2}\otimes gx_{1}x_{2};X_{1}X_{2},gx_{1}x_{2}\right)
1_{A}\otimes gx_{1}x_{2}\otimes gx_{1}x_{2}.
\end{equation*}%
By considering also the right side, we get%
\begin{gather*}
B(1_{H}\otimes gx_{1}x_{2};1_{A},gx_{1}x_{2})+B(gx_{2}\otimes
gx_{1}x_{2};X_{2},gx_{1}x_{2})+ \\
+B(gx_{1}\otimes gx_{1}x_{2};X_{1},gx_{1}x_{2})+B\left( x_{1}x_{2}\otimes
gx_{1}x_{2};X_{1}X_{2},gx_{1}x_{2}\right) \\
-B(x_{1}x_{2}\otimes g;1_{A},gx_{1}x_{2})=0
\end{gather*}%
which holds in view of the form of the elements.

\subsubsection{Case $X_{1}X_{2}\otimes gx_{1}\otimes x_{2}$}

We have to consider only the third and the fourth summand of the left side
of the equality. Third summand gives us%
\begin{eqnarray*}
l_{1} &=&u_{1}=0,l_{2}=u_{2}=0, \\
a &=&0,b_{1}=b_{2}=1 \\
d &=&e_{1}=1,e_{2}=0.
\end{eqnarray*}%
Since $\alpha \left( x_{2};0,0,0,0\right) \equiv a+b_{1}+b_{2}=0,$ we get
\begin{equation*}
-B(gx_{1}\otimes gx_{1}x_{2};X_{1}X_{2},gx_{1})X_{1}X_{2}\otimes
gx_{1}\otimes x_{2}.
\end{equation*}%
Fourth summand gives us%
\begin{eqnarray*}
l_{1} &=&u_{1}=0,l_{2}+u_{2}=1, \\
a &=&0,b_{1}=1,b_{2}-l_{2}=1\Rightarrow b_{2}=1,l_{2}=0,u_{2}=1, \\
d &=&e_{1}=1,e_{2}=u_{2}=1.
\end{eqnarray*}%
Since $\alpha \left( 1_{H};0,0,0,1\right) \equiv a+b_{1}+b_{2}\equiv 0$ we
get
\begin{equation*}
+B(x_{1}x_{2}\otimes gx_{1}x_{2};X_{1}X_{2},gx_{1}x_{2})X_{1}X_{2}\otimes
gx_{1}\otimes x_{2}.
\end{equation*}%
By considering also the right side we get%
\begin{gather*}
+B(x_{1}x_{2}\otimes gx_{1}x_{2};X_{1}X_{2},gx_{1}x_{2})-B(x_{1}x_{2}\otimes
gx_{1};X_{1}X_{2},gx_{1})+ \\
-B(gx_{1}\otimes gx_{1}x_{2};X_{1}X_{2},gx_{1})=0
\end{gather*}%
which holds in view of the form of the elements.

\subsubsection{Case $X_{2}\otimes gx_{1}\otimes gx_{1}x_{2}$}

First summand of the left side of the equality gives us%
\begin{eqnarray*}
l_{1} &=&u_{1}=0,l_{2}=u_{2}=0 \\
a &=&b_{1}=0,b_{2}=1, \\
d &=&e_{1}=1,e_{2}=0.
\end{eqnarray*}%
Since $\alpha \left( x_{1}x_{2};0,0,0,0\right) \equiv 0,$ we get%
\begin{equation*}
B(1_{H}\otimes gx_{1}x_{2};X_{2},gx_{1})X_{2}\otimes gx_{1}\otimes
gx_{1}x_{2}.
\end{equation*}%
Second summand of the left side of the equation gives us%
\begin{eqnarray*}
l_{1} &=&u_{1}=0,l_{2}+u_{2}=1 \\
a &=&b_{1}=0,b_{2}-l_{2}=1\Rightarrow b_{2}=1,l_{2}=0,u_{2}=1, \\
d &=&e_{1}=1,e_{2}=u_{2}=1.
\end{eqnarray*}%
Since $\alpha \left( x_{1};0,0,0,1\right) \equiv 1,$ we get%
\begin{equation*}
-B(gx_{2}\otimes gx_{1}x_{2};X_{2},gx_{1}x_{2})X_{2}\otimes gx_{1}\otimes
gx_{1}x_{2}.
\end{equation*}

Third summand of the left side of the equation gives us%
\begin{eqnarray*}
l_{1}+u_{1} &=&1,l_{2}=u_{2}=0, \\
a &=&0,b_{2}=1,b_{1}=l_{1}, \\
d &=&1,e_{2}=0,e_{1}-u_{1}=1\Rightarrow e_{1}=1,u_{1}=0,b_{1}=l_{1}=1.
\end{eqnarray*}%
Since $\alpha \left( x_{2};1,0,0,0\right) \equiv a+b_{1}\equiv 1$ $,$ we get%
\begin{equation*}
+B(gx_{1}\otimes gx_{1}x_{2};X_{1}X_{2},gx_{1})X_{2}\otimes gx_{1}\otimes
gx_{1}x_{2}.
\end{equation*}%
Fourth summand of the left side of the equation gives us%
\begin{eqnarray*}
l_{1}+u_{1} &=&1,l_{2}+u_{2}=1 \\
a &=&0,b_{1}=l_{1},b_{2}-l_{2}=1\Rightarrow b_{2}=1,l_{2}=0,u_{2}=1, \\
d &=&1,e_{1}-u_{1}=1\Rightarrow e_{1}=1,u_{1}=0,b_{1}=l_{1}=1, \\
e_{2} &=&u_{2}=1.
\end{eqnarray*}%
Since $\alpha \left( 1_{H};1,0,0,1\right) \equiv a+b_{1}\equiv 1,$ we get%
\begin{equation*}
-B\left( x_{1}x_{2}\otimes gx_{1}x_{2};X_{1}X_{2},gx_{1}x_{2}\right)
X_{2}\otimes gx_{1}\otimes gx_{1}x_{2}.
\end{equation*}%
By considering also the right side, we get%
\begin{gather*}
B(1_{H}\otimes gx_{1}x_{2};X_{2},gx_{1})-B(gx_{2}\otimes
gx_{1}x_{2};X_{2},gx_{1}x_{2})+ \\
+B(gx_{1}\otimes gx_{1}x_{2};X_{1}X_{2},gx_{1})-B\left( x_{1}x_{2}\otimes
gx_{1}x_{2};X_{1}X_{2},gx_{1}x_{2}\right) \\
-B(x_{1}x_{2}\otimes g;X_{2},gx_{1})=0
\end{gather*}%
which holds in view of the form of the elements.

\subsubsection{Case $X_{1}X_{2}\otimes g\otimes gx_{1}x_{2}$}

First summand of the left side of the equality gives us%
\begin{eqnarray*}
l_{1} &=&u_{1}=0,l_{2}=u_{2}=0 \\
a &=&0,b_{1}=b_{2}=1, \\
d &=&1,e_{1}=e_{2}=0.
\end{eqnarray*}%
Since $\alpha \left( x_{1}x_{2};0,0,0,0\right) \equiv 0,$ we get%
\begin{equation*}
B(1_{H}\otimes gx_{1}x_{2};X_{1}X_{2},g)X_{1}X_{2}\otimes g\otimes
gx_{1}x_{2}.
\end{equation*}%
Second summand of the left side of the equation gives us%
\begin{eqnarray*}
l_{1} &=&u_{1}=0,l_{2}+u_{2}=1 \\
a &=&0,b_{1}=1,b_{2}-l_{2}=1\Rightarrow b_{2}=1,l_{2}=0,u_{2}=1, \\
d &=&1,e_{1}=0,e_{2}=u_{2}=1.
\end{eqnarray*}%
Since $\alpha \left( x_{1};0,0,0,1\right) \equiv 1,$ we get%
\begin{equation*}
-B(gx_{2}\otimes gx_{1}x_{2};X_{1}X_{2},gx_{2})X_{1}X_{2}\otimes g\otimes
gx_{1}x_{2}.
\end{equation*}

Third summand of the left side of the equation gives us%
\begin{eqnarray*}
l_{1}+u_{1} &=&1,l_{2}=u_{2}=0, \\
a &=&0,b_{2}=1,b_{1}-l_{1}=1\Rightarrow b_{1}=1,l_{1}=0,u_{1}=1, \\
d &=&1,e_{2}=0,e_{1}=u_{1}=1.
\end{eqnarray*}%
Since $\alpha \left( x_{2};0,0,1,0\right) \equiv e_{2}\equiv 0$ $,$ we get%
\begin{equation*}
-B(gx_{1}\otimes gx_{1}x_{2};X_{1}X_{2},gx_{1})X_{1}X_{2}\otimes g\otimes
gx_{1}x_{2}.
\end{equation*}%
Fourth summand of the left side of the equation gives us%
\begin{eqnarray*}
l_{1}+u_{1} &=&1,l_{2}+u_{2}=1 \\
a &=&0,b_{1}-l_{1}=1\Rightarrow b_{1}=1,l_{1}=0,u_{1}=1, \\
b_{2}-l_{2} &=&1\Rightarrow b_{2}=1,l_{2}=0,u_{2}=1, \\
d &=&1,e_{1}=u_{1}=1,e_{2}=u_{2}=1.
\end{eqnarray*}%
Since $\alpha \left( 1_{H};0,0,1,1\right) =1+e_{2}\equiv 0,$ we get%
\begin{equation*}
+B\left( x_{1}x_{2}\otimes gx_{1}x_{2};X_{1}X_{2},gx_{1}x_{2}\right)
X_{1}X_{2}\otimes g\otimes gx_{1}x_{2}.
\end{equation*}%
By considering also the right side, we get%
\begin{gather*}
B(1_{H}\otimes gx_{1}x_{2};X_{1}X_{2},g)-B(gx_{2}\otimes
gx_{1}x_{2};X_{1}X_{2},gx_{2})+ \\
-B(gx_{1}\otimes gx_{1}x_{2};X_{1}X_{2},gx_{1})+B\left( x_{1}x_{2}\otimes
gx_{1}x_{2};X_{1}X_{2},gx_{1}x_{2}\right) \\
-B(x_{1}x_{2}\otimes g;X_{1}X_{2},g)=0
\end{gather*}%
which holds in view of the form of the elements.

\subsection{$B(x_{1}x_{2}\otimes gx_{1}x_{2};GX_{1},gx_{1}x_{2})$}

We deduce that%
\begin{eqnarray*}
a &=&1,b_{1}=1,b_{2}=0 \\
d &=&e_{1}=e_{2}=1
\end{eqnarray*}%
and we get%
\begin{gather*}
\left( -1\right) ^{\alpha \left( 1_{H};0,0,0,0\right) }B(x_{1}x_{2}\otimes
gx_{1}x_{2};GX_{1},gx_{1}x_{2})GX_{1}\otimes gx_{1}x_{2}\otimes g \\
\left( -1\right) ^{\alpha \left( 1_{H};1,0,0,0\right) }B(x_{1}x_{2}\otimes
gx_{1}x_{2};GX_{1},gx_{1}x_{2})G\otimes gx_{1}x_{2}\otimes x_{1} \\
\left( -1\right) ^{\alpha \left( 1_{H};0,0,1,0\right) }B(x_{1}x_{2}\otimes
gx_{1}x_{2};GX_{1},gx_{1}x_{2})GX_{1}\otimes gx_{2}\otimes x_{1} \\
\left( -1\right) ^{\alpha \left( 1_{H};1,0,1,0\right) }B(x_{1}x_{2}\otimes
gx_{1}x_{2};GX_{1},gx_{1}x_{2})GX_{1}^{1-l_{1}}\otimes gx_{2}\otimes
g^{+l_{11}}x_{1}^{1+1}=0 \\
\left( -1\right) ^{\alpha \left( 1_{H};0,0,0,1\right) }B(x_{1}x_{2}\otimes
gx_{1}x_{2};GX_{1},gx_{1}x_{2})GX_{1}\otimes gx_{1}\otimes x_{2} \\
\left( -1\right) ^{\alpha \left( 1_{H};1,0,0,1\right) }B(x_{1}x_{2}\otimes
gx_{1}x_{2};GX_{1},gx_{1}x_{2})G\otimes gx_{1}\otimes gx_{1}x_{2} \\
\left( -1\right) ^{\alpha \left( 1_{H};0,0,1,1\right) }B(x_{1}x_{2}\otimes
gx_{1}x_{2};GX_{1},gx_{1}x_{2})GX_{1}\otimes g\otimes gx_{1}x_{2} \\
\left( -1\right) ^{\alpha \left( 1_{H};1,0,1,1\right) }B(x_{1}x_{2}\otimes
gx_{1}x_{2};GX_{1},gx_{1}x_{2})GX_{1}^{1-l_{1}}\otimes g\otimes
g^{+l_{1}+1}x_{1}^{1+1}x_{2}=0
\end{gather*}

\subsubsection{Case $G\otimes gx_{1}x_{2}\otimes x_{1}$}

We have to consider only the second and the fourth summand of the left side
of the equality.

Second summand gives us%
\begin{eqnarray*}
l_{1} &=&u_{1}=l_{2}=u_{2}=0, \\
a &=&1,b_{1}=b_{2}=0, \\
d &=&e_{1}=e_{2}=1.
\end{eqnarray*}%
Since $\alpha \left( x_{1};0,0,0,0\right) \equiv a+b_{1}+b_{2}\equiv 1$ we
obtain
\begin{equation*}
-B(gx_{2}\otimes gx_{1}x_{2};G,gx_{1}x_{2})G\otimes gx_{1}x_{2}\otimes x_{1}.
\end{equation*}%
Forth summand gives us%
\begin{eqnarray*}
l_{1}+u_{1} &=&1,l_{2}=u_{2}=0 \\
a &=&1,b_{1}=l_{1},b_{2}=0, \\
d &=&e_{2}=1,e_{1}-u_{1}=1\Rightarrow e_{1}=1,u_{1}=0,b_{1}=l_{1}=1.
\end{eqnarray*}%
Since $\alpha \left( 1_{H};1,0,0,0\right) \equiv b_{2}\equiv 0,$ we obtain
\begin{equation*}
+B(x_{1}x_{2}\otimes gx_{1}x_{2};GX_{1},gx_{1}x_{2})G\otimes
gx_{1}x_{2}\otimes x_{1}.
\end{equation*}%
By considering the right side, we get,
\begin{gather*}
B(x_{1}x_{2}\otimes gx_{2};G,gx_{1}x_{2})+B(x_{1}x_{2}\otimes
gx_{1}x_{2};GX_{1},gx_{1}x_{2})+ \\
-B(gx_{2}\otimes gx_{1}x_{2};G,gx_{1}x_{2})=0
\end{gather*}%
which holds in view of the form of the elements.

\subsubsection{Case $GX_{1}\otimes gx_{2}\otimes x_{1}$}

We have to consider only the second and the fourth summand of the left side
of the equality.

Second summand gives us%
\begin{eqnarray*}
l_{1} &=&u_{1}=l_{2}=u_{2}=0, \\
a &=&b_{1}=1,b_{2}=0, \\
d &=&e_{2}=1,e_{1}=0.
\end{eqnarray*}%
Since $\alpha \left( x_{1};0,0,0,0\right) \equiv a+b_{1}+b_{2}\equiv 0$ we
obtain
\begin{equation*}
+B(gx_{2}\otimes gx_{1}x_{2};GX_{1},gx_{2})GX_{1}\otimes gx_{2}\otimes x_{1}.
\end{equation*}%
Forth summand gives us%
\begin{eqnarray*}
l_{1}+u_{1} &=&1,l_{2}=u_{2}=0 \\
a &=&1,b_{1}-l_{1}=1\Rightarrow b_{1}=1,l_{1}=0,u_{1}=1,b_{2}=0, \\
d &=&e_{2}=1,e_{1}=u_{1}=1.
\end{eqnarray*}%
Since $\alpha \left( 1_{H};0,0,1,0\right) \equiv e_{2}+\left(
a+b_{1}+b_{2}\right) \equiv 1,$ we obtain
\begin{equation*}
-B(x_{1}x_{2}\otimes gx_{1}x_{2};GX_{1},gx_{1}x_{2})GX_{1}\otimes
gx_{2}\otimes x_{1}.
\end{equation*}%
By considering the right side, we get,
\begin{gather*}
B(x_{1}x_{2}\otimes gx_{2};GX_{1},gx_{2})-B(x_{1}x_{2}\otimes
gx_{1}x_{2};GX_{1},gx_{1}x_{2})+ \\
+B(gx_{2}\otimes gx_{1}x_{2};GX_{1},gx_{2})=0
\end{gather*}%
which holds in view of the form of the elements.

\subsubsection{Case $GX_{1}\otimes gx_{1}\otimes x_{2}$}

We have to consider only the third and the fourth summand of the left side
of the equality. Third summand gives us%
\begin{eqnarray*}
l_{1} &=&u_{1}=0,l_{2}=u_{2}=0, \\
a &=&b_{1}=1,b_{2}=0 \\
d &=&e_{1}=1,e_{2}=0.
\end{eqnarray*}%
Since $\alpha \left( x_{2};0,0,0,0\right) \equiv a+b_{1}+b_{2}=0,$ we get
\begin{equation*}
-B(gx_{1}\otimes gx_{1}x_{2};GX_{1},gx_{1})GX_{1}\otimes gx_{1}\otimes x_{2}.
\end{equation*}%
Fourth summand gives us%
\begin{eqnarray*}
l_{1} &=&u_{1}=0,l_{2}+u_{2}=1, \\
a &=&1,b_{1}=1,b_{2}=l_{2}, \\
d &=&e_{1}=1,e_{2}=u_{2}.
\end{eqnarray*}%
Since $\alpha \left( 1_{H};0,0,0,1\right) \equiv a+b_{1}+b_{2}\equiv 0$ and $%
\alpha \left( 1_{H};0,1,0,0\right) \equiv 0$ we get
\begin{equation*}
\left[ -B(x_{1}x_{2}\otimes
gx_{1}x_{2};GX_{1},gx_{1}x_{2})+B(x_{1}x_{2}\otimes
gx_{1}x_{2};GX_{1}X_{2},gx_{1})\right] GX_{1}\otimes gx_{1}\otimes x_{2}.
\end{equation*}%
By considering also the right side we get%
\begin{gather*}
-B(x_{1}x_{2}\otimes gx_{1};GX_{1},gx_{1})+B(x_{1}x_{2}\otimes
gx_{1}x_{2};GX_{1},gx_{1}x_{2})+ \\
+B(x_{1}x_{2}\otimes gx_{1}x_{2};GX_{1}X_{2},gx_{1})-B(gx_{1}\otimes
gx_{1}x_{2};GX_{1},gx_{1})=0
\end{gather*}%
which holds in view of the form of the elements.

\subsubsection{Case $G\otimes gx_{1}\otimes gx_{1}x_{2}$}

First summand of the left side of the equality gives us%
\begin{eqnarray*}
l_{1} &=&u_{1}=0,l_{2}=u_{2}=0 \\
a &=&1,b_{1}=b_{2}=0, \\
d &=&e_{1}=1,e_{2}=0.
\end{eqnarray*}%
Since $\alpha \left( x_{1}x_{2};0,0,0,0\right) \equiv 0,$ we get%
\begin{equation*}
B(1_{H}\otimes gx_{1}x_{2};G,gx_{1})G\otimes gx_{1}\otimes gx_{1}x_{2}.
\end{equation*}%
Second summand of the left side of the equation gives us%
\begin{eqnarray*}
l_{1} &=&u_{1}=0,l_{2}+u_{2}=1 \\
a &=&1,b_{1}=0,b_{2}=l_{2}, \\
d &=&e_{1}=1,0,e_{2}=u_{2}.
\end{eqnarray*}%
Since $\alpha \left( x_{1};0,0,0,1\right) \equiv 1$ and $\alpha \left(
x_{1};0,1,0,0\right) \equiv a+b_{1}+b_{2}+1\equiv 1,$ we get%
\begin{equation*}
\left[ -B(gx_{2}\otimes gx_{1}x_{2};G,gx_{1}x_{2})-B(gx_{2}\otimes
gx_{1}x_{2};GX_{2},gx_{1})\right] G\otimes gx_{1}\otimes gx_{1}x_{2}.
\end{equation*}

Third summand of the left side of the equation gives us%
\begin{eqnarray*}
l_{1}+u_{1} &=&1,l_{2}=u_{2}=0, \\
a &=&1,b_{2}=0,b_{1}=l_{1}, \\
d &=&1,e_{2}=0, \\
e_{1}-u_{1} &=&1\Rightarrow e_{1}=1,u_{1}=0,b_{1}=l_{1}=1.
\end{eqnarray*}%
Since $\alpha \left( x_{2};1,0,0,0\right) \equiv a+b_{1}\equiv 0$ $,$ we get%
\begin{equation*}
-B(gx_{1}\otimes gx_{1}x_{2};GX_{1},gx_{1})G\otimes gx_{1}\otimes
gx_{1}x_{2}.
\end{equation*}%
Fourth summand of the left side of the equation gives us%
\begin{eqnarray*}
l_{1}+u_{1} &=&1,l_{2}+u_{2}=1 \\
a &=&1,b_{1}=l_{1},b_{2}=l_{2},, \\
d &=&1,e_{2}=u_{2}, \\
e_{1}-u_{1} &=&1\Rightarrow e_{1}=1,u_{1}=0,b_{1}=l_{1}=1.
\end{eqnarray*}%
Since $\alpha \left( 1_{H};1,0,0,1\right) \equiv a+b_{1}\equiv 0$ and $%
\alpha \left( 1_{H};1,1,0,0\right) \equiv 1+b_{2}\equiv 0,$ we get%
\begin{equation*}
\left[ +B\left( x_{1}x_{2}\otimes gx_{1}x_{2};GX_{1},gx_{1}x_{2}\right)
+B\left( x_{1}x_{2}\otimes gx_{1}x_{2};GX_{1}X_{2},gx_{1}\right) \right]
G\otimes gx_{1}\otimes gx_{1}x_{2}.
\end{equation*}%
By considering also the right side, we get%
\begin{gather*}
B(1_{H}\otimes gx_{1}x_{2};G,gx_{1})-B(gx_{2}\otimes
gx_{1}x_{2};G,gx_{1}x_{2}) \\
-B(gx_{2}\otimes gx_{1}x_{2};GX_{2},gx_{1})-B(gx_{1}\otimes
gx_{1}x_{2};GX_{1},gx_{1})+ \\
+B\left( x_{1}x_{2}\otimes gx_{1}x_{2};GX_{1},gx_{1}x_{2}\right) +B\left(
x_{1}x_{2}\otimes gx_{1}x_{2};GX_{1}X_{2},gx_{1}\right) \\
-B(x_{1}x_{2}\otimes g;G,gx_{1})=0
\end{gather*}%
which holds in view of the form of the elements.

\subsubsection{Case $GX_{1}\otimes g\otimes gx_{1}x_{2}$}

First summand of the left side of the equality gives us%
\begin{eqnarray*}
l_{1} &=&u_{1}=0,l_{2}=u_{2}=0 \\
a &=&b_{1}=1,b_{2}=0, \\
d &=&1,e_{1}=e_{2}=0.
\end{eqnarray*}%
Since $\alpha \left( x_{1}x_{2};0,0,0,0\right) \equiv 0,$ we get%
\begin{equation*}
B(1_{H}\otimes gx_{1}x_{2};GX_{1},g)GX_{1}\otimes g\otimes gx_{1}x_{2}.
\end{equation*}%
Second summand of the left side of the equation gives us%
\begin{eqnarray*}
l_{1} &=&u_{1}=0,l_{2}+u_{2}=1 \\
a &=&b_{1}=1,b_{2}=l_{2}, \\
d &=&1,e_{1}=0,e_{2}=u_{2}.
\end{eqnarray*}%
Since $\alpha \left( x_{1};0,0,0,1\right) \equiv 1$ and $\alpha \left(
x_{1};0,1,0,0\right) \equiv a+b_{1}+b_{2}+1\equiv 0,$ we get%
\begin{equation*}
\left[ -B(gx_{2}\otimes gx_{1}x_{2};GX_{1},gx_{2})+B(gx_{2}\otimes
gx_{1}x_{2};GX_{1}X_{2},g)\right] GX_{1}\otimes g\otimes gx_{1}x_{2}.
\end{equation*}

Third summand of the left side of the equation gives us%
\begin{eqnarray*}
l_{1}+u_{1} &=&1,l_{2}=u_{2}=0, \\
a &=&1,b_{2}=0,b_{1}-l_{1}=1\Rightarrow b_{1}=1,l_{1}=0,u_{1}=1 \\
d &=&1,e_{2}=0,e_{1}=u_{1}=1
\end{eqnarray*}%
Since $\alpha \left( x_{2};0,0,1,0\right) \equiv e_{2}\equiv 0$ $,$ we get%
\begin{equation*}
-B(gx_{1}\otimes gx_{1}x_{2};GX_{1},gx_{1})GX_{1}\otimes g\otimes
gx_{1}x_{2}.
\end{equation*}%
Fourth summand of the left side of the equation gives us%
\begin{eqnarray*}
l_{1}+u_{1} &=&1,l_{2}+u_{2}=1 \\
a &=&1,b_{2}=l_{2}, \\
b_{1}-l_{1} &=&1\Rightarrow b_{1}=1,l_{1}=0,u_{1}=1 \\
d &=&1,e_{2}=u_{2},e_{1}=u_{1}=1
\end{eqnarray*}%
Since $\alpha \left( 1_{H};0,0,1,1\right) \equiv 1+e_{2}\equiv 0$ and $%
\alpha \left( 1_{H};0,1,1,0\right) \equiv e_{2}+a+b_{1}+b_{2}+1\equiv 0,$ we
get%
\begin{equation*}
\left[ +B\left( x_{1}x_{2}\otimes gx_{1}x_{2};GX_{1},gx_{1}x_{2}\right)
+B\left( x_{1}x_{2}\otimes gx_{1}x_{2};GX_{1}X_{2},gx_{1}\right) \right]
GX_{1}\otimes g\otimes gx_{1}x_{2}.
\end{equation*}%
By considering also the right side, we get%
\begin{gather*}
B(1_{H}\otimes gx_{1}x_{2};GX_{1},g)-B(gx_{2}\otimes
gx_{1}x_{2};GX_{1},gx_{2})+ \\
+B(gx_{2}\otimes gx_{1}x_{2};GX_{1}X_{2},g)-B(gx_{1}\otimes
gx_{1}x_{2};GX_{1},gx_{1})+ \\
+B\left( x_{1}x_{2}\otimes gx_{1}x_{2};GX_{1},gx_{1}x_{2}\right) +B\left(
x_{1}x_{2}\otimes gx_{1}x_{2};GX_{1}X_{2},gx_{1}\right) \\
-B(x_{1}x_{2}\otimes g;GX_{1},g)=0
\end{gather*}%
which holds in view of the form of the elements.

\subsection{$B(x_{1}x_{2}\otimes gx_{1}x_{2};GX_{2},gx_{1}x_{2})$}

We have%
\begin{eqnarray*}
a &=&1,b_{1}=0,b_{2}=1 \\
d &=&e_{1}=e_{2}=1
\end{eqnarray*}%
and we get%
\begin{gather*}
\left( -1\right) ^{\alpha \left( 1_{H};0,0,0,0\right) }B\left(
x_{1}x_{2}\otimes gx_{1}x_{2};GX_{2},gx_{1}x_{2}\,\right) GX_{2}\otimes
gx_{1}x_{2}\otimes g+ \\
\left( -1\right) ^{\alpha \left( 1_{H};0,1,0,0\right) }B\left(
x_{1}x_{2}\otimes gx_{1}x_{2};GX_{2},gx_{1}x_{2}\,\right) G\otimes
gx_{1}x_{2}\otimes x_{2}+ \\
\left( -1\right) ^{\alpha \left( 1_{H};0,0,1,0\right) }B\left(
x_{1}x_{2}\otimes gx_{1}x_{2};GX_{2},gx_{1}x_{2}\,\right) GX_{2}\otimes
gx_{2}\otimes x_{1}+ \\
\left( -1\right) ^{\alpha \left( 1_{H};0,1,1,0\right) }B\left(
x_{1}x_{2}\otimes gx_{1}x_{2};GX_{2},gx_{1}x_{2}\,\right) G\otimes
gx_{2}\otimes gx_{1}x_{2}+ \\
\left( -1\right) ^{\alpha \left( 1_{H};0,0,0,1\right) }B\left(
x_{1}x_{2}\otimes gx_{1}x_{2};GX_{2},gx_{1}x_{2}\,\right) GX_{2}\otimes
gx_{1}\otimes x_{2}+ \\
\left( -1\right) ^{\alpha \left( 1_{H};0,1,0,1\right) }B\left(
x_{1}x_{2}\otimes gx_{1}x_{2};GX_{2},gx_{1}x_{2}\,\right) G\otimes
gx_{1}\otimes g^{l_{2}}x_{2}^{1+1}=0 \\
\left( -1\right) ^{\alpha \left( 1_{H};0,0,1,1\right) }B\left(
x_{1}x_{2}\otimes gx_{1}x_{2};GX_{2},gx_{1}x_{2}\,\right) GX_{2}\otimes
g\otimes gx_{1}x_{2} \\
\left( -1\right) ^{\alpha \left( 1_{H};0,1,1,1\right) }B\left(
x_{1}x_{2}\otimes gx_{1}x_{2};GX_{2},gx_{1}x_{2}\,\right)
GX_{2}^{1-l_{2}}\otimes g\otimes g^{l_{2}+1}x_{1}x_{2}^{1+1}=0.
\end{gather*}

\subsubsection{Case $G\otimes gx_{1}x_{2}\otimes x_{2}$}

We have to consider only the third and the fourth summand of the left side
of the equality. Third summand gives us%
\begin{eqnarray*}
l_{1} &=&u_{1}=0,l_{2}=u_{2}=0, \\
a &=&1,b_{1}=b_{2}=0 \\
d &=&e_{1}=e_{2}=1.
\end{eqnarray*}%
Since $\alpha \left( x_{2};0,0,0,0\right) \equiv a+b_{1}+b_{2}=1,$ we get
\begin{equation*}
+B(gx_{1}\otimes gx_{1}x_{2};G,gx_{1}x_{2})G\otimes gx_{1}x_{2}\otimes x_{2}.
\end{equation*}%
Fourth summand gives us%
\begin{eqnarray*}
l_{1} &=&u_{1}=0,l_{2}+u_{2}=1, \\
a &=&1,b_{1}=0,b_{2}=l_{2}, \\
d &=&e_{1}=1, \\
e_{2}-u_{2} &=&1\Rightarrow e_{2}=1,u_{2}=0,b_{2}=l_{2}=1.
\end{eqnarray*}%
Since $\alpha \left( 1_{H};0,1,0,0\right) \equiv 0$ we get
\begin{equation*}
+B(x_{1}x_{2}\otimes gx_{1}x_{2};GX_{2},gx_{1}x_{2})G\otimes
gx_{1}x_{2}\otimes x_{2}.
\end{equation*}%
By considering also the right side we get%
\begin{gather*}
-B(x_{1}x_{2}\otimes gx_{1};G,gx_{1}x_{2})+B(x_{1}x_{2}\otimes
gx_{1}x_{2};GX_{2},gx_{1}x_{2})+ \\
+B(gx_{1}\otimes gx_{1}x_{2};G,gx_{1}x_{2})=0
\end{gather*}%
which holds in view of the form of the elements.

\subsubsection{Case $GX_{2}\otimes gx_{2}\otimes x_{1}$}

We have to consider only the second and the fourth summand of the left side
of the equality.

Second summand gives us%
\begin{eqnarray*}
l_{1} &=&u_{1}=l_{2}=u_{2}=0, \\
a &=&b_{2}=1,b_{1}=0, \\
d &=&e_{2}=1,e_{1}=0.
\end{eqnarray*}%
Since $\alpha \left( x_{1};0,0,0,0\right) \equiv a+b_{1}+b_{2}\equiv 0$ we
obtain
\begin{equation*}
+B(gx_{2}\otimes gx_{1}x_{2};GX_{2},gx_{2})GX_{2}\otimes gx_{2}\otimes x_{1}.
\end{equation*}%
Forth summand gives us%
\begin{eqnarray*}
l_{1}+u_{1} &=&1,l_{2}=u_{2}=0 \\
a &=&b_{2}=1,b_{1}=l_{1}, \\
d &=&e_{2}=1,e_{1}=u_{1}.
\end{eqnarray*}%
Since $\alpha \left( 1_{H};0,0,1,0\right) \equiv e_{2}+\left(
a+b_{1}+b_{2}\right) \equiv 1$ and $\alpha \left( 1_{H};1,0,0,0\right)
\equiv b_{2}=1,$ we obtain
\begin{equation*}
\left[ -B(x_{1}x_{2}\otimes
gx_{1}x_{2};GX_{2},gx_{1}x_{2})-B(x_{1}x_{2}\otimes
gx_{1}x_{2};GX_{1}X_{2},gx_{2})\right] GX_{2}\otimes gx_{2}\otimes x_{1}.
\end{equation*}%
By considering the right side, we get,
\begin{gather*}
B(x_{1}x_{2}\otimes gx_{2};GX_{2},gx_{2})-B(x_{1}x_{2}\otimes
gx_{1}x_{2};GX_{2},gx_{1}x_{2})+ \\
-B(x_{1}x_{2}\otimes gx_{1}x_{2};GX_{1}X_{2},gx_{2})+B(gx_{2}\otimes
gx_{1}x_{2};GX_{2},gx_{2})=0
\end{gather*}%
which holds in view of the form of the elements.

\subsubsection{Case $G\otimes gx_{2}\otimes gx_{1}x_{2}$}

First summand of the left side of the equality gives us%
\begin{eqnarray*}
l_{1} &=&u_{1}=0,l_{2}=u_{2}=0 \\
a &=&1,b_{1}=b_{2}=0, \\
d &=&e_{2}=1,e_{1}=0.
\end{eqnarray*}%
Since $\alpha \left( x_{1}x_{2};0,0,0,0\right) \equiv 0,$ we get%
\begin{equation*}
B(1_{H}\otimes gx_{1}x_{2};G,gx_{2})G\otimes gx_{2}\otimes gx_{1}x_{2}.
\end{equation*}%
Second summand of the left side of the equation gives us%
\begin{eqnarray*}
l_{1} &=&u_{1}=0,l_{2}+u_{2}=1 \\
a &=&1,b_{1}=0,b_{2}=l_{2}, \\
d &=&1,e_{1}=0, \\
e_{2}-u_{2} &=&1\Rightarrow e_{2}=1,u_{2}=0,b_{2}=l_{2}=1.
\end{eqnarray*}%
Since $\alpha \left( x_{1};0,1,0,0\right) \equiv a+b_{1}+b_{2}+1\equiv 1,$
we get%
\begin{equation*}
-B(gx_{2}\otimes gx_{1}x_{2};GX_{2},gx_{2})G\otimes gx_{2}\otimes
gx_{1}x_{2}.
\end{equation*}

Third summand of the left side of the equation gives us%
\begin{eqnarray*}
l_{1}+u_{1} &=&1,l_{2}=u_{2}=0, \\
a &=&1,b_{2}=0,b_{1}=l_{1} \\
d &=&e_{2}=1,e_{1}=u_{1}
\end{eqnarray*}%
Since $\alpha \left( x_{2};0,0,1,0\right) \equiv e_{2}\equiv 1$ $\alpha
\left( x_{2};1,0,0,0\right) \equiv a+b_{1}\equiv 0$ and $,$ we get%
\begin{equation*}
\left[ +B(gx_{1}\otimes gx_{1}x_{2};G,gx_{1}x_{2})-B(gx_{1}\otimes
gx_{1}x_{2};GX_{1},gx_{2})\right] G\otimes gx_{2}\otimes gx_{1}x_{2}.
\end{equation*}%
Fourth summand of the left side of the equation gives us%
\begin{eqnarray*}
l_{1}+u_{1} &=&1,l_{2}+u_{2}=1 \\
a &=&1,b_{2}=l_{2},b_{1}=l_{1} \\
d &=&1,e_{1}=u_{1} \\
e_{2}-u_{2} &=&1\Rightarrow e_{2}=1,u_{2}=0,b_{2}=l_{2}=1
\end{eqnarray*}%
Since $\alpha \left( 1_{H};0,1,1,0\right) \equiv e_{2}+a+b_{1}+b_{2}+1\equiv
0$ and $\alpha \left( 1_{H};1,1,0,0\right) \equiv 1+b_{2}\equiv 0,$ we get%
\begin{equation*}
\left[ +B\left( x_{1}x_{2}\otimes gx_{1}x_{2};GX_{2},gx_{1}x_{2}\right)
+B\left( x_{1}x_{2}\otimes gx_{1}x_{2};GX_{1}X_{2},gx_{2}\right) \right]
G\otimes gx_{2}\otimes gx_{1}x_{2}.
\end{equation*}%
By considering also the right side, we get%
\begin{gather*}
B(1_{H}\otimes gx_{1}x_{2};G,gx_{2})+B(gx_{1}\otimes
gx_{1}x_{2};G,gx_{1}x_{2}) \\
-B(gx_{1}\otimes gx_{1}x_{2};GX_{1},gx_{2})-B(gx_{2}\otimes
gx_{1}x_{2};GX_{2},gx_{2})+ \\
+B\left( x_{1}x_{2}\otimes gx_{1}x_{2};GX_{2},gx_{1}x_{2}\right) +B\left(
x_{1}x_{2}\otimes gx_{1}x_{2};GX_{1}X_{2},gx_{2}\right) \\
-B(x_{1}x_{2}\otimes g;G,gx_{2})=0
\end{gather*}%
which holds in view of the form of the elements.

\subsubsection{Case $GX_{2}\otimes gx_{1}\otimes x_{2}$}

We have to consider only the third and the fourth summand of the left side
of the equality. Third summand gives us%
\begin{eqnarray*}
l_{1} &=&u_{1}=0,l_{2}=u_{2}=0, \\
a &=&b_{2}=1,b_{1}=0 \\
d &=&e_{1}=1,e_{2}=0.
\end{eqnarray*}%
Since $\alpha \left( x_{2};0,0,0,0\right) \equiv a+b_{1}+b_{2}\equiv 0,$ we
get
\begin{equation*}
-B(gx_{1}\otimes gx_{1}x_{2};GX_{2},gx_{1})GX_{2}\otimes gx_{1}\otimes x_{2}.
\end{equation*}%
Fourth summand gives us%
\begin{eqnarray*}
l_{1} &=&u_{1}=0,l_{2}+u_{2}=1, \\
a &=&1,b_{1}=0, \\
b_{2}-l_{2} &=&1\Rightarrow b_{2}=1,l_{2}=0,u_{2}=1 \\
d &=&e_{1}=1,e_{2}=u_{2}=1.
\end{eqnarray*}%
Since $\alpha \left( 1_{H};0,0,0,1\right) \equiv a+b_{1}+b_{2}\equiv 0$ we
get
\begin{equation*}
+B(x_{1}x_{2}\otimes gx_{1}x_{2};GX_{2},gx_{1}x_{2})GX_{2}\otimes
gx_{1}\otimes x_{2}.
\end{equation*}%
By considering also the right side we get%
\begin{gather*}
-B(x_{1}x_{2}\otimes gx_{1};GX_{2},gx_{1})+B(x_{1}x_{2}\otimes
gx_{1}x_{2};GX_{2},gx_{1}x_{2})+ \\
-B(gx_{1}\otimes gx_{1}x_{2};GX_{2},gx_{1})=0
\end{gather*}%
which holds in view of the form of the elements.

\subsubsection{Case $GX_{2}\otimes g\otimes gx_{1}x_{2}$}

First summand of the left side of the equality gives us%
\begin{eqnarray*}
l_{1} &=&u_{1}=0,l_{2}=u_{2}=0 \\
a &=&b_{2}=1,b_{1}=0, \\
d &=&1,e_{1}=e_{2}=0.
\end{eqnarray*}%
Since $\alpha \left( x_{1}x_{2};0,0,0,0\right) \equiv 0,$ we get%
\begin{equation*}
B(1_{H}\otimes gx_{1}x_{2};GX_{2},g)GX_{2}\otimes g\otimes gx_{1}x_{2}.
\end{equation*}%
Second summand of the left side of the equation gives us%
\begin{eqnarray*}
l_{1} &=&u_{1}=0,l_{2}+u_{2}=1 \\
a &=&1,b_{1}=0, \\
b_{2}-l_{2} &=&1\Rightarrow b_{2}=1,l_{2}=0,u_{2}=1 \\
d &=&1,e_{1}=0,e_{2}=u_{2}=1
\end{eqnarray*}%
Since $\alpha \left( x_{1};0,0,0,1\right) \equiv 1,$ we get%
\begin{equation*}
-B(gx_{2}\otimes gx_{1}x_{2};GX_{2},gx_{2})GX_{2}\otimes g\otimes
gx_{1}x_{2}.
\end{equation*}

Third summand of the left side of the equation gives us%
\begin{eqnarray*}
l_{1}+u_{1} &=&1,l_{2}=u_{2}=0, \\
a &=&b_{2}=1,b_{1}=l_{1} \\
d &=&1,e_{2}=0,e_{1}=u_{1}
\end{eqnarray*}%
Since $\alpha \left( x_{2};0,0,1,0\right) \equiv e_{2}\equiv 0$ $\alpha
\left( x_{2};1,0,0,0\right) \equiv a+b_{1}\equiv 0$ and $,$ we get%
\begin{equation*}
\left[ -B(gx_{1}\otimes gx_{1}x_{2};GX_{2},gx_{1})-B(gx_{1}\otimes
gx_{1}x_{2};GX_{1}X_{2},g)\right] GX_{2}\otimes g\otimes gx_{1}x_{2}.
\end{equation*}%
Fourth summand of the left side of the equation gives us%
\begin{eqnarray*}
l_{1}+u_{1} &=&1,l_{2}+u_{2}=1 \\
a &=&1,b_{1}=l_{1}, \\
b_{2}-l_{2} &=&1\Rightarrow b_{2}=1,l_{2}=0,u_{2}=1 \\
d &=&1,e_{1}=u_{1},e_{2}=u_{2}=1
\end{eqnarray*}%
Since $\alpha \left( 1_{H};0,0,1,1\right) \equiv 1+e_{2}\equiv 0$ and $%
\alpha \left( 1_{H};1,0,0,1\right) \equiv a+b_{1}\equiv 0,$ we get%
\begin{equation*}
\left[ +B\left( x_{1}x_{2}\otimes gx_{1}x_{2};GX_{2},gx_{1}x_{2}\right)
+B\left( x_{1}x_{2}\otimes gx_{1}x_{2};GX_{1}X_{2},gx_{2}\right) \right]
GX_{2}\otimes g\otimes gx_{1}x_{2}.
\end{equation*}%
By considering also the right side, we get%
\begin{gather*}
B(1_{H}\otimes gx_{1}x_{2};GX_{2},g)-B(gx_{1}\otimes
gx_{1}x_{2};GX_{2},gx_{1}) \\
-B(gx_{1}\otimes gx_{1}x_{2};GX_{1}X_{2},g)-B(gx_{2}\otimes
gx_{1}x_{2};GX_{2},gx_{2})+ \\
+B\left( x_{1}x_{2}\otimes gx_{1}x_{2};GX_{2},gx_{1}x_{2}\right) +B\left(
x_{1}x_{2}\otimes gx_{1}x_{2};GX_{1}X_{2},gx_{2}\right) \\
-B(x_{1}x_{2}\otimes g;GX_{2},g)=0
\end{gather*}%
which holds in view of the form of the elements.

\part{Condition$\left(\protect\ref{eq.h}\right) $ and condition$\left(
\protect\ref{eq.a}\right) $}

In this part we will explore condition$\left( \ref{eq.h}\right) $ and
condition$\left( \ref{eq.a}\right) .$At the end we will have all the
material needed to prove our fundamental result Theorem $\left( \ref{Theo
sep}\right) $

\section{Equation $\left( \protect\ref{eq.h}\right) $}

In this section we will prove that, in view of our choice to write all $%
B\left( h\otimes h^{\prime }\right) $ by using $B\left( h\otimes
1_{H}\right) $ where $h\otimes 1_{H}$ runs through the seven elements, this
condition $\left( \ref{eq.h}\right) $ will be automatically satisfied.

As we remarked in section $\left( \ref{mor cond}\right) $, we only need to
compute equality $\left( \ref{eq.h}\right) $ for $1_{A}\otimes g$, $%
1_{A}\otimes x_{1}$ and $1_{A}\otimes x_{2}$.

\subsection{Case $1_{A}\otimes g$}

We know that $\left( \ref{eq.h}\right) $ in the case of $1_{A}\otimes g$
assume the form of $\left( \ref{eq.10}\right) $%
\begin{equation*}
B(gh\otimes gh^{\prime })=(1_{A}\otimes g)B(h\otimes h^{\prime
})(1_{A}\otimes g).
\end{equation*}%
Now we use this equality to define the remaining $B(h\otimes h^{\prime })$
using the 34 ones of our list in the Casimir condition. Thus this equality
is automatically satisfied.

\subsection{Case $1_{A}\otimes x_{1}$}

\begin{proposition}
Assume that equation $\left( \ref{eq.h}\right) $ holds for $1_{A}\otimes
x_{1}$ and $h\otimes h^{\prime }$%
\begin{equation}
B(h\otimes h^{\prime })(1_{A}\otimes x_{1})=(1_{A}\otimes x_{1})B(gh\otimes
gh^{\prime })+B(x_{1}h\otimes gh^{\prime })+B(h\otimes x_{1}h^{\prime })
\label{1x1}
\end{equation}%
then it holds for $gh\otimes gh^{\prime }$ i.e.%
\begin{equation*}
B(gh\otimes gh^{\prime })(1_{A}\otimes x_{1})=(1_{A}\otimes
x_{1})B(ggh\otimes ggh^{\prime })+B(x_{1}gh\otimes ggh^{\prime
})+B(gh\otimes x_{1}gh^{\prime })
\end{equation*}
\end{proposition}

\begin{proof}
We conjugate equality $\left( \ref{1x1}\right) $ by $\left( 1_{A}\otimes
g\right) $on both sides and we get%
\begin{gather*}
-\left( 1_{A}\otimes g\right) B(h\otimes h^{\prime })\left( 1_{A}\otimes
g\right) (1_{A}\otimes x_{1})=-(1_{A}\otimes x_{1})\left( 1_{A}\otimes
g\right) B(gh\otimes gh^{\prime })\left( 1_{A}\otimes g\right) \\
+\left( 1_{A}\otimes g\right) B(x_{1}h\otimes gh^{\prime })\left(
1_{A}\otimes g\right) +\left( 1_{A}\otimes g\right) B(h\otimes
x_{1}h^{\prime })\left( 1_{A}\otimes g\right) .
\end{gather*}%
From his we deduce that%
\begin{gather*}
B(gh\otimes gh^{\prime })(1_{A}\otimes x_{1})=(1_{A}\otimes x_{1})B(h\otimes
h^{\prime }) \\
-B(gx_{1}h\otimes h^{\prime })-B(gh\otimes gx_{1}h^{\prime })
\end{gather*}
\end{proof}

\subsubsection{$g\otimes 1_{H}$}

\begin{equation*}
B(g\otimes 1_{H})(1_{A}\otimes x_{1})\overset{?}{=}(1_{A}\otimes
x_{1})B(gg\otimes g)+B(x_{1}g\otimes g)+B(g\otimes x_{1})
\end{equation*}%
\begin{eqnarray*}
&&B(g\otimes 1_{H})(1_{A}\otimes x_{1})\overset{?}{=}(1_{A}\otimes
x_{1})\left( 1_{A}\otimes g\right) B(g\otimes 1_{H})\left( 1_{A}\otimes
g\right) \\
&&-\left( 1_{A}\otimes g\right) B(x_{1}\otimes 1_{H})\left( 1_{A}\otimes
g\right) +\left( 1_{A}\otimes g\right) B(1_{H}\otimes gx_{1})\left(
1_{A}\otimes g\right)
\end{eqnarray*}%
\begin{eqnarray*}
\left( 1_{A}\otimes g\right) &&B(1_{H}\otimes gx_{1})\left( 1_{A}\otimes
g\right) \\
&&\overset{\left( \ref{form 1otgx1}\right) }{=}B(g\otimes
1_{H})(1_{A}\otimes x_{1}) \\
&&+\left( 1_{A}\otimes g\right) (1_{A}\otimes x_{1})B(g\otimes 1_{H})\left(
1_{A}\otimes g\right) \\
&&+\left( 1_{A}\otimes g\right) B(x_{1}\otimes 1_{H})\left( 1_{A}\otimes
g\right)
\end{eqnarray*}%
\begin{eqnarray*}
&&(1_{A}\otimes x_{1})\left( 1_{A}\otimes g\right) B(g\otimes 1_{H})\left(
1_{A}\otimes g\right) -\left( 1_{A}\otimes g\right) B(x_{1}\otimes
1_{H})\left( 1_{A}\otimes g\right) \\
&&+\left( 1_{A}\otimes g\right) B(1_{H}\otimes gx_{1})\left( 1_{A}\otimes
g\right) \\
&=&\left( 1_{A}\otimes gx_{1}\right) B(g\otimes 1_{H})\left( 1_{A}\otimes
g\right) -\left( 1_{A}\otimes g\right) B(x_{1}\otimes 1_{H})\left(
1_{A}\otimes g\right) \\
&&+B(g\otimes 1_{H})(1_{A}\otimes x_{1})-(1_{A}\otimes gx_{1})B(g\otimes
1_{H})\left( 1_{A}\otimes g\right) \\
&&+\left( 1_{A}\otimes g\right) B(x_{1}\otimes 1_{H})\left( 1_{A}\otimes
g\right) \\
&=&B(g\otimes 1_{H})(1_{A}\otimes x_{1})
\end{eqnarray*}

\subsubsection{$x_{1}\otimes 1_{H}$}

\begin{equation*}
B(x_{1}\otimes 1_{H})(1_{A}\otimes x_{1})\overset{?}{=}(1_{A}\otimes
x_{1})B(gx_{1}\otimes g)+B(x_{1}x_{1}\otimes g)+B(x_{1}\otimes x_{1})
\end{equation*}%
We have%
\begin{eqnarray*}
&&(1_{A}\otimes x_{1})B(gx_{1}\otimes g)+B(x_{1}x_{1}\otimes
g)+B(x_{1}\otimes x_{1}) \\
&=&(1_{A}\otimes x_{1})\left( 1_{A}\otimes g\right) B(x_{1}\otimes
1_{H})\left( 1_{A}\otimes g\right) +B(x_{1}\otimes x_{1})
\end{eqnarray*}%
Since%
\begin{equation*}
B(x_{1}\otimes x_{1})\overset{\left( \ref{form x1otx1}\right) }{=}%
B(x_{1}\otimes 1_{H})(1_{A}\otimes x_{1})-(1_{A}\otimes
gx_{1})B(x_{1}\otimes 1_{H})(1_{A}\otimes g)
\end{equation*}%
we obtain%
\begin{gather*}
(1_{A}\otimes x_{1})\left( 1_{A}\otimes g\right) B(x_{1}\otimes 1_{H})\left(
1_{A}\otimes g\right) +B(x_{1}\otimes x_{1}) \\
=(1_{A}\otimes gx_{1})B(x_{1}\otimes 1_{H})\left( 1_{A}\otimes g\right)
+B(x_{1}\otimes 1_{H})(1_{A}\otimes x_{1}) \\
-(1_{A}\otimes gx_{1})B(x_{1}\otimes 1_{H})(1_{A}\otimes g)=B(x_{1}\otimes
1_{H})(1_{A}\otimes x_{1})
\end{gather*}

and we conclude.

\subsubsection{$x_{2}\otimes 1_{H}$}

\begin{equation*}
B(x_{2}\otimes 1_{H})(1_{A}\otimes x_{1})\overset{?}{=}\left( 1_{A}\otimes
gx_{1}\right) B(x_{2}\otimes 1_{H})\left( 1_{A}\otimes g\right)
+B(x_{1}x_{2}\otimes g)+B(x_{2}\otimes x_{1})
\end{equation*}%
Since%
\begin{eqnarray*}
&&\left( 1_{A}\otimes gx_{1}\right) B(x_{2}\otimes 1_{H})\left( 1_{A}\otimes
g\right) +B(x_{1}x_{2}\otimes g)+B(x_{2}\otimes x_{1}) \\
&=&\left( 1_{A}\otimes gx_{1}\right) B(x_{2}\otimes 1_{H})\left(
1_{A}\otimes g\right) +\left( 1_{A}\otimes g\right) B(gx_{1}x_{2}\otimes
1_{H})\left( 1_{A}\otimes g\right) +B(x_{2}\otimes x_{1})
\end{eqnarray*}%
and%
\begin{eqnarray*}
&&B(x_{2}\otimes x_{1})\overset{\left( \ref{form x2otx1}\right) }{=}%
B(x_{2}\otimes 1_{H})(1_{A}\otimes x_{1}) \\
&&-(1_{A}\otimes gx_{1})B(x_{2}\otimes 1_{H})(1_{A}\otimes g)-(1_{A}\otimes
g)B(gx_{1}x_{2}\otimes 1_{H})(1_{A}\otimes g)
\end{eqnarray*}%
we get%
\begin{eqnarray*}
\text{right side} &&\left( 1_{A}\otimes gx_{1}\right) B(x_{2}\otimes
1_{H})\left( 1_{A}\otimes g\right) +\left( 1_{A}\otimes g\right)
B(gx_{1}x_{2}\otimes 1_{H})\left( 1_{A}\otimes g\right) + \\
&&+B(x_{2}\otimes 1_{H})(1_{A}\otimes x_{1}) \\
&&-(1_{A}\otimes gx_{1})B(x_{2}\otimes 1_{H})(1_{A}\otimes g)-(1_{A}\otimes
g)B(gx_{1}x_{2}\otimes 1_{H})(1_{A}\otimes g) \\
&=&B(x_{2}\otimes 1_{H})(1_{A}\otimes x_{1})
\end{eqnarray*}%
and we conclude.

\subsubsection{$x_{1}x_{2}\otimes 1_{H}$}

\begin{equation*}
B(x_{1}x_{2}\otimes 1_{H})(1_{A}\otimes x_{1})\overset{?}{=}\left(
1_{A}\otimes gx_{1}\right) B(x_{1}x_{2}\otimes 1_{H})\left( 1_{A}\otimes
g\right) +B(x_{1}x_{1}x_{2}\otimes g)+B(x_{1}x_{2}\otimes x_{1})
\end{equation*}%
We get%
\begin{eqnarray*}
&&\left( 1_{A}\otimes gx_{1}\right) B(x_{1}x_{2}\otimes 1_{H})\left(
1_{A}\otimes g\right) +B(x_{1}x_{1}x_{2}\otimes g)+B(x_{1}x_{2}\otimes x_{1})
\\
&=&\left( 1_{A}\otimes gx_{1}\right) B(x_{1}x_{2}\otimes 1_{H})\left(
1_{A}\otimes g\right) +B(x_{1}x_{2}\otimes x_{1})
\end{eqnarray*}%
and we conclude.

\subsubsection{$gx_{1}\otimes 1_{H}$}

\begin{equation*}
B\left( gx_{1}\otimes 1_{H}\right) (1_{A}\otimes x_{1})\overset{?}{=}%
(1_{A}\otimes x_{1})B(ggx_{1}\otimes g)+B(x_{1}gx_{1}\otimes
g)+B(gx_{1}\otimes x_{1})
\end{equation*}

We consider the right side%
\begin{eqnarray*}
&&(1_{A}\otimes x_{1})B(x_{1}\otimes g)+B(gx_{1}\otimes x_{1}) \\
&=&(1_{A}\otimes x_{1})\left( 1_{A}\otimes g\right) B(gx_{1}\otimes
1_{H})\left( 1_{A}\otimes g\right) +\left( 1_{A}\otimes g\right)
B(x_{1}\otimes gx_{1})\left( 1_{A}\otimes g\right) \\
&&\overset{\left( \ref{form x1otgx1}\right) }{=}(1_{A}\otimes x_{1})\left(
1_{A}\otimes g\right) B(gx_{1}\otimes 1_{H})\left( 1_{A}\otimes g\right) + \\
&&\left( 1_{A}\otimes g\right) (1_{A}\otimes g)B(gx_{1}\otimes
1_{H})(1_{A}\otimes gx_{1})\left( 1_{A}\otimes g\right) + \\
&&+\left( 1_{A}\otimes g\right) (1_{A}\otimes x_{1})B(gx_{1}\otimes
1_{H})\left( 1_{A}\otimes g\right) \\
&=&\left( 1_{A}\otimes gx_{1}\right) B(gx_{1}\otimes 1_{H})\left(
1_{A}\otimes g\right) + \\
&&++B(gx_{1}\otimes 1_{H})(1_{A}\otimes x_{1}) \\
&&-(1_{A}\otimes gx_{1})B(gx_{1}\otimes 1_{H})\left( 1_{A}\otimes g\right) \\
&=&B(gx_{1}\otimes 1_{H})(1_{A}\otimes x_{1})
\end{eqnarray*}

and we conclude.

\subsubsection{$gx_{2}\otimes 1_{H}$}

\begin{equation*}
B(gx_{2}\otimes 1_{H})(1_{A}\otimes x_{1})\overset{?}{=}(1_{A}\otimes
x_{1})B(ggx_{2}\otimes g)+B(x_{1}gx_{2}\otimes g)+B(gx_{2}\otimes x_{1})
\end{equation*}%
We consider the right side%
\begin{eqnarray*}
&&(1_{A}\otimes x_{1})B(ggx_{2}\otimes g)+B(x_{1}gx_{2}\otimes
g)+B(gx_{2}\otimes x_{1}) \\
&=&(1_{A}\otimes x_{1})\left( 1_{A}\otimes g\right) B(gx_{2}\otimes
1_{H})\left( 1_{A}\otimes g\right) - \\
&&\left( 1_{A}\otimes g\right) B(x_{1}x_{2}\otimes 1_{H})\left( 1_{A}\otimes
g\right) +\left( 1_{A}\otimes g\right) B(x_{2}\otimes gx_{1})\left(
1_{A}\otimes g\right) \\
&&\overset{\left( \ref{form x2otgx1}\right) }{=}(1_{A}\otimes x_{1})\left(
1_{A}\otimes g\right) B(gx_{2}\otimes 1_{H})\left( 1_{A}\otimes g\right) \\
&&-\left( 1_{A}\otimes g\right) B(x_{1}x_{2}\otimes 1_{H})\left(
1_{A}\otimes g\right) + \\
&&+\left( 1_{A}\otimes g\right) (1_{A}\otimes g)B(gx_{2}\otimes
1_{H})(1_{A}\otimes gx_{1})\left( 1_{A}\otimes g\right) \\
&&+\left( 1_{A}\otimes g\right) (1_{A}\otimes x_{1})B(gx_{2}\otimes
1_{H})\left( 1_{A}\otimes g\right) \\
&&+\left( 1_{A}\otimes g\right) B(x_{1}x_{2}\otimes 1_{H})\left(
1_{A}\otimes g\right) \\
&=&\left( 1_{A}\otimes gx_{1}\right) B(gx_{2}\otimes 1_{H})\left(
1_{A}\otimes g\right) \\
&&-\left( 1_{A}\otimes g\right) B(x_{1}x_{2}\otimes 1_{H})\left(
1_{A}\otimes g\right) + \\
&&+B(gx_{2}\otimes 1_{H})(1_{A}\otimes x_{1})-(1_{A}\otimes
gx_{1})B(gx_{2}\otimes 1_{H})\left( 1_{A}\otimes g\right) \\
&=&-\left( 1_{A}\otimes g\right) B(x_{1}x_{2}\otimes 1_{H})\left(
1_{A}\otimes g\right) +B(gx_{2}\otimes 1_{H})(1_{A}\otimes x_{1}) \\
&&+\left( 1_{A}\otimes g\right) B(x_{1}x_{2}\otimes 1_{H})\left(
1_{A}\otimes g\right) \\
&=&B(gx_{2}\otimes 1_{H})(1_{A}\otimes x_{1})
\end{eqnarray*}%
and so we are done.

\subsubsection{$gx_{1}x_{2}\otimes 1_{H}$}

\begin{equation*}
B(gx_{1}x_{2}\otimes 1_{H})(1_{A}\otimes x_{1})\overset{?}{=}(1_{A}\otimes
x_{1})B(ggx_{1}x_{2}\otimes g)+B(x_{1}gx_{1}x_{2}\otimes
g)+B(gx_{1}x_{2}\otimes x_{1})
\end{equation*}%
We consider the right side.%
\begin{eqnarray*}
&&(1_{A}\otimes x_{1})B(ggx_{1}x_{2}\otimes g)+B(x_{1}gx_{1}x_{2}\otimes
g)+B(gx_{1}x_{2}\otimes x_{1}) \\
&=&(1_{A}\otimes x_{1})B(x_{1}x_{2}\otimes g)+B(gx_{1}x_{2}\otimes x_{1}) \\
&=&(1_{A}\otimes x_{1})\left( 1_{A}\otimes g\right) B(gx_{1}x_{2}\otimes
1_{H})\left( 1_{A}\otimes g\right) + \\
&&+\left( 1_{A}\otimes g\right) B(x_{1}x_{2}\otimes gx_{1})\left(
1_{A}\otimes g\right) \\
&&\overset{\left( \ref{form x1x2otgx1}\right) }{=}(1_{A}\otimes x_{1})\left(
1_{A}\otimes g\right) B(gx_{1}x_{2}\otimes 1_{H})\left( 1_{A}\otimes
g\right) + \\
&&+\left( 1_{A}\otimes g\right) (1_{A}\otimes g)B(gx_{1}x_{2}\otimes
1_{H})(1_{A}\otimes gx_{1})\left( 1_{A}\otimes g\right) \\
&&+\left( 1_{A}\otimes g\right) (1_{A}\otimes x_{1})B(gx_{1}x_{2}\otimes
1_{H})\left( 1_{A}\otimes g\right) \\
&=&\left( 1_{A}\otimes gx_{1}\right) B(gx_{1}x_{2}\otimes 1_{H})\left(
1_{A}\otimes g\right) + \\
&&B(gx_{1}x_{2}\otimes 1_{H})(1_{A}\otimes x_{1})+ \\
&&-\left( 1_{A}\otimes gx_{1}\right) B(gx_{1}x_{2}\otimes 1_{H})\left(
1_{A}\otimes g\right)
\end{eqnarray*}%
and we are done.%
\begin{eqnarray*}
&&B(x_{1}x_{2}\otimes gx_{1})\left( \ref{form x1x2otgx1}\right) \\
&=&(1_{A}\otimes g)B(gx_{1}x_{2}\otimes 1_{H})(1_{A}\otimes gx_{1}) \\
&&+(1_{A}\otimes x_{1})B(gx_{1}x_{2}\otimes 1_{H})
\end{eqnarray*}

\subsubsection{$1_{H}\otimes x_{1}$}

\begin{equation*}
B(1_{H}\otimes x_{1})(1_{H}\otimes x_{1})\overset{?}{=}(1_{H}\otimes
x_{1})B(g\otimes gx_{1})+B(x_{1}\otimes gx_{1})+B(1_{H}\otimes x_{1}x_{1})
\end{equation*}

We consider the right side%
\begin{eqnarray*}
&&(1_{H}\otimes x_{1})B(g\otimes gx_{1})+B(x_{1}\otimes
gx_{1})+B(1_{H}\otimes x_{1}x_{1}) \\
&=&(1_{H}\otimes x_{1})(1_{H}\otimes g)B(1_{H}\otimes x_{1})(1_{H}\otimes
g)+B(x_{1}\otimes gx_{1}) \\
&&\overset{\left( \ref{form x1otgx1},\left( \ref{form:1otxi}\right) \right) }%
{=}-(1_{H}\otimes x_{1})(1_{H}\otimes g)(1_{H}\otimes g)B(gx_{1}\otimes
1_{H})(1_{H}\otimes g)(1_{H}\otimes g)+ \\
&&(1_{H}\otimes g)B(gx_{1}\otimes 1_{H})(1_{H}\otimes gx_{1})+(1_{H}\otimes
x_{1})B(gx_{1}\otimes 1_{H}) \\
&=&-(1_{H}\otimes x_{1})B(gx_{1}\otimes 1_{H})+(1_{H}\otimes
g)B(gx_{1}\otimes 1_{H})(1_{H}\otimes gx_{1})+(1_{H}\otimes
x_{1})B(gx_{1}\otimes 1_{H}) \\
&=&(1_{H}\otimes g)B(gx_{1}\otimes 1_{H})(1_{H}\otimes gx_{1})
\end{eqnarray*}%
We consider the left side%
\begin{equation*}
B(1_{H}\otimes x_{1})(1_{H}\otimes x_{1})\overset{\left( \ref{form:1otxi}%
\right) }{=}+(1_{H}\otimes g)B(gx_{1}\otimes 1_{H})(1_{H}\otimes gx_{1})
\end{equation*}%
and we conclude.

\subsubsection{$1_{H}\otimes x_{2}$}

\begin{equation*}
B(1_{H}\otimes x_{2})(1_{A}\otimes x_{1})\overset{?}{=}(1_{A}\otimes
x_{1})B(g\otimes gx_{2})+B(x_{1}\otimes gx_{2})+B(1_{H}\otimes x_{1}x_{2})
\end{equation*}

We consider the left side%
\begin{equation*}
B(1_{H}\otimes x_{2})(1_{A}\otimes x_{1})\overset{\left( \ref{form:1otxi}%
\right) }{=}-(1_{A}\otimes g)B(gx_{2}\otimes 1_{H})(1_{A}\otimes
g)(1_{A}\otimes x_{1})
\end{equation*}

We consider the right side%
\begin{equation*}
(1_{A}\otimes x_{1})B(g\otimes gx_{2})+B(x_{1}\otimes gx_{2})+B(1_{H}\otimes
x_{1}x_{2})
\end{equation*}

\begin{eqnarray*}
(1_{A}\otimes x_{1})B(g\otimes gx_{2}) &=&(1_{A}\otimes x_{1})(1_{H}\otimes
g)B(1_{A}\otimes x_{2})(1_{H}\otimes g) \\
&&\overset{\left( \ref{form:1otxi}\right) }{=}-(1_{A}\otimes
x_{1})(1_{A}\otimes g)(1_{A}\otimes g)B(gx_{2}\otimes 1_{H})(1_{A}\otimes
g)(1_{A}\otimes g) \\
&=&-(1_{A}\otimes x_{1})B(gx_{2}\otimes 1_{H})
\end{eqnarray*}%
\begin{eqnarray*}
&&B(x_{1}\otimes gx_{2})\left( \ref{form x1otgx2}\right) \\
&=&(1_{A}\otimes g)B(gx_{1}\otimes 1_{H})(1_{A}\otimes gx_{2}) \\
&&+(1_{A}\otimes x_{2})B(gx_{1}\otimes 1_{H}) \\
&&-B(gx_{2}gx_{1}\otimes 1_{H})
\end{eqnarray*}%
\begin{equation*}
B(x_{1}\otimes gx_{2})\overset{\left( \ref{form x1otgx2}\right) }{=}%
(1_{A}\otimes g)B(gx_{1}\otimes 1_{H})(1_{A}\otimes gx_{2})+(1_{A}\otimes
x_{2})B(gx_{1}\otimes 1_{H})-B(x_{1}x_{2}\otimes 1_{H})
\end{equation*}%
\begin{eqnarray*}
&&B(1_{H}\otimes x_{1}x_{2})\overset{\left( \ref{form 1otx1x2}\right) }{=} \\
&=&-(1_{A}\otimes g)B(gx_{2}\otimes 1_{H})(1_{A}\otimes g)(1_{A}\otimes
x_{1}) \\
&&+(1_{A}\otimes gx_{1})(1_{A}\otimes g)B(gx_{2}\otimes 1_{H}) \\
&&-(1_{A}\otimes g)B(gx_{1}\otimes 1_{H})(1_{A}\otimes x_{2})(1_{A}\otimes g)
\\
&&+(1_{A}\otimes g)(1_{A}\otimes gx_{2})B(gx_{1}\otimes 1_{H}) \\
&&+B(x_{1}x_{2}\otimes 1_{H})
\end{eqnarray*}%
Therefore we get%
\begin{eqnarray*}
&&(1_{A}\otimes x_{1})B(g\otimes gx_{2})+B(x_{1}\otimes
gx_{2})+B(1_{H}\otimes x_{1}x_{2})\overset{\left( \ref{form:1otxi}\right)
,\left( \ref{form x1otgx2}\right) \left( \ref{form 1otx1x2}\right) }{=} \\
&&-(1_{A}\otimes x_{1})B(gx_{2}\otimes 1_{H})+ \\
&&(1_{A}\otimes g)B(gx_{1}\otimes 1_{H})(1_{A}\otimes gx_{2})+(1_{A}\otimes
x_{2})B(gx_{1}\otimes 1_{H})-B(x_{1}x_{2}\otimes 1_{H})+ \\
&&-(1_{A}\otimes g)B(gx_{2}\otimes 1_{H})(1_{A}\otimes g)(1_{A}\otimes
x_{1})+ \\
&&+(1_{A}\otimes gx_{1})(1_{A}\otimes g)B(gx_{2}\otimes 1_{H})+ \\
&&-(1_{A}\otimes g)B(gx_{1}\otimes 1_{H})(1_{A}\otimes x_{2})(1_{A}\otimes
g)+ \\
&&+(1_{A}\otimes g)(1_{A}\otimes gx_{2})B(gx_{1}\otimes
1_{H})+B(x_{1}x_{2}\otimes 1_{H}) \\
&=&-(1_{A}\otimes g)B(gx_{2}\otimes 1_{H})(1_{A}\otimes g)(1_{A}\otimes
x_{1})
\end{eqnarray*}%
and we conclude.

\subsubsection{$1_{H}\otimes x_{1}x_{2}$}

\begin{equation*}
B(1_{H}\otimes x_{1}x_{2})(1_{A}\otimes x_{1})\overset{?}{=}(1_{A}\otimes
x_{1})B(g\otimes gx_{1}x_{2})+B(x_{1}\otimes gx_{1}x_{2})+B(1_{H}\otimes
x_{1}x_{1}x_{2})
\end{equation*}%
We write the left side%
\begin{eqnarray*}
\text{left side} &&B(1_{H}\otimes x_{1}x_{2})(1_{A}\otimes x_{1})\overset{%
\left( \ref{form 1otx1x2}\right) }{=} \\
&=&+(1_{A}\otimes g)B(gx_{2}\otimes 1_{H})(1_{A}\otimes gx_{1})(1_{A}\otimes
x_{1}) \\
&&+(1_{A}\otimes x_{1})B(gx_{2}\otimes 1_{H})(1_{A}\otimes x_{1}) \\
&&-(1_{A}\otimes g)B(gx_{1}\otimes 1_{H})(1_{A}\otimes gx_{2})(1_{A}\otimes
x_{1}) \\
&&-(1_{A}\otimes x_{2})B(gx_{1}\otimes 1_{H})(1_{A}\otimes x_{1}) \\
&&+B(x_{1}x_{2}\otimes 1_{H})(1_{A}\otimes x_{1})
\end{eqnarray*}%
and we get%
\begin{eqnarray*}
\text{left side } &&B(1_{H}\otimes x_{1}x_{2})(1_{A}\otimes x_{1})= \\
&&+(1_{A}\otimes g)B(gx_{1}\otimes 1_{H})(1_{A}\otimes gx_{1}x_{2}) \\
&&-(1_{A}\otimes x_{2})B(gx_{1}\otimes 1_{H})(1_{A}\otimes x_{1}) \\
&&+B(x_{1}x_{2}\otimes 1_{H})(1_{A}\otimes x_{1})
\end{eqnarray*}%
We write the right side%
\begin{equation*}
\text{right side }(1_{A}\otimes x_{1})B(g\otimes gx_{1}x_{2})+B(x_{1}\otimes
gx_{1}x_{2})
\end{equation*}%
\begin{equation*}
\text{first summand of the right side}(1_{A}\otimes x_{1})B(g\otimes
gx_{1}x_{2})
\end{equation*}%
\begin{eqnarray*}
&&\text{first summand of the right side} \\
(1_{A}\otimes x_{1})B(g\otimes gx_{1}x_{2}) &=&(1_{A}\otimes
x_{1})(1_{A}\otimes g)B(1_{H}\otimes x_{1}x_{2})(1_{A}\otimes g)
\end{eqnarray*}%
\begin{eqnarray*}
&&\text{first summand of the right side} \\
&&(1_{A}\otimes x_{1})(1_{A}\otimes g)B(1_{H}\otimes
x_{1}x_{2})(1_{A}\otimes g)\overset{\left( \ref{form 1otx1x2}\right) }{=} \\
&=&+(1_{A}\otimes x_{1})(1_{A}\otimes g)((1_{A}\otimes x_{1})(1_{A}\otimes
g)1_{A}\otimes g)B(gx_{2}\otimes 1_{H})(1_{A}\otimes gx_{1})(1_{A}\otimes g)
\\
&&+(1_{A}\otimes x_{1})(1_{A}\otimes g)(1_{A}\otimes x_{1})B(gx_{2}\otimes
1_{H})(1_{A}\otimes g) \\
&&-(1_{A}\otimes x_{1})(1_{A}\otimes g)(1_{A}\otimes g)B(gx_{1}\otimes
1_{H})(1_{A}\otimes gx_{2})(1_{A}\otimes g) \\
&&-(1_{A}\otimes x_{1})(1_{A}\otimes g)(1_{A}\otimes x_{2})B(gx_{1}\otimes
1_{H})(1_{A}\otimes g) \\
&&+(1_{A}\otimes x_{1})(1_{A}\otimes g)B(x_{1}x_{2}\otimes
1_{H})(1_{A}\otimes g)
\end{eqnarray*}%
\begin{eqnarray*}
&&\text{Clean first summand of the right side} \\
&&(1_{A}\otimes x_{1})(1_{A}\otimes g)B(1_{H}\otimes
x_{1}x_{2})(1_{A}\otimes g)\overset{\left( \ref{form 1otx1x2}\right) }{=} \\
&=&-(1_{A}\otimes x_{1})B(gx_{1}\otimes 1_{H})(1_{A}\otimes x_{2}) \\
&&-(1_{A}\otimes gx_{1}x_{2})B(gx_{1}\otimes 1_{H})(1_{A}\otimes g) \\
&&+(1_{A}\otimes gx_{1})B(x_{1}x_{2}\otimes 1_{H})(1_{A}\otimes g)
\end{eqnarray*}%
\begin{eqnarray*}
&&\text{second summand of the right side} \\
&&B(x_{1}\otimes gx_{1}x_{2})\overset{\left( \ref{form x1otgx1x2}\right) }{=}%
-B(x_{1}\otimes gx_{2})(1_{A}\otimes x_{1})+(1_{A}\otimes
gx_{1})B(x_{1}\otimes gx_{2})(1_{A}\otimes g)
\end{eqnarray*}%
\begin{eqnarray*}
&&\text{first summand of the second summand of the right side} \\
&&\text{ }B(x_{1}\otimes gx_{2})(1_{A}\otimes x_{1})\overset{\left( \ref%
{form x1otgx2}\right) }{=} \\
&=&-(1_{A}\otimes g)B(gx_{1}\otimes 1_{H})(1_{A}\otimes gx_{2})(1_{A}\otimes
x_{1}) \\
&&-(1_{A}\otimes x_{2})B(gx_{1}\otimes 1_{H})(1_{A}\otimes x_{1}) \\
&&+B(x_{1}x_{2}\otimes 1_{H})(1_{A}\otimes x_{1})
\end{eqnarray*}%
\begin{eqnarray*}
&&\text{second summand of second summand of the right side} \\
&&(1_{A}\otimes gx_{1})B(x_{1}\otimes gx_{2})(1_{A}\otimes g)\overset{\left( %
\ref{form x1otgx2}\right) }{=} \\
&=&(1_{A}\otimes gx_{1})(1_{A}\otimes g)B(gx_{1}\otimes 1_{H})(1_{A}\otimes
gx_{2})(1_{A}\otimes g) \\
&&+(1_{A}\otimes gx_{1})(1_{A}\otimes x_{2})B(gx_{1}\otimes
1_{H})(1_{A}\otimes g) \\
&&-(1_{A}\otimes gx_{1})B(gx_{2}gx_{1}\otimes 1_{H})(1_{A}\otimes g)
\end{eqnarray*}%
\begin{eqnarray*}
&&\text{second summand of second summand of the right side} \\
&&(1_{A}\otimes gx_{1})B(x_{1}\otimes gx_{2})(1_{A}\otimes g)\overset{\left( %
\ref{form x1otgx2}\right) }{=} \\
&=&(1_{A}\otimes x_{1})B(gx_{1}\otimes 1_{H})(1_{A}\otimes x_{2}) \\
&&+(1_{A}\otimes gx_{1}x_{2})B(gx_{1}\otimes 1_{H})(1_{A}\otimes g) \\
&&-(1_{A}\otimes gx_{1})B(x_{1}x_{2}\otimes 1_{H})(1_{A}\otimes g)
\end{eqnarray*}%
\begin{eqnarray*}
&&\text{second summand of the right side} \\
B(x_{1}\otimes gx_{1}x_{2}) &=&-(1_{A}\otimes g)B(gx_{1}\otimes
1_{H})(1_{A}\otimes gx_{2})(1_{A}\otimes x_{1}) \\
&&-(1_{A}\otimes x_{2})B(gx_{1}\otimes 1_{H})(1_{A}\otimes x_{1}) \\
&&+B(x_{1}x_{2}\otimes 1_{H})(1_{A}\otimes x_{1})+ \\
&&(1_{A}\otimes x_{1})B(gx_{1}\otimes 1_{H})(1_{A}\otimes x_{2})+ \\
&&+(1_{A}\otimes gx_{1}x_{2})B(gx_{1}\otimes 1_{H})(1_{A}\otimes g)+ \\
&&-(1_{A}\otimes gx_{1})B(x_{1}x_{2}\otimes 1_{H})(1_{A}\otimes g)
\end{eqnarray*}%
Thus the right side becomes%
\begin{eqnarray*}
&&\text{right side}(1_{A}\otimes x_{1})B(g\otimes
gx_{1}x_{2})+B(x_{1}\otimes gx_{1}x_{2}) \\
&=&-(1_{A}\otimes x_{1})B(gx_{1}\otimes 1_{H})(1_{A}\otimes x_{2}) \\
&&-(1_{A}\otimes gx_{1}x_{2})B(gx_{1}\otimes 1_{H})(1_{A}\otimes g) \\
&&+(1_{A}\otimes gx_{1})B(x_{1}x_{2}\otimes 1_{H})(1_{A}\otimes g) \\
&&-(1_{A}\otimes g)B(gx_{1}\otimes 1_{H})(1_{A}\otimes gx_{2})(1_{A}\otimes
x_{1}) \\
&&-(1_{A}\otimes x_{2})B(gx_{1}\otimes 1_{H})(1_{A}\otimes x_{1}) \\
&&+B(x_{1}x_{2}\otimes 1_{H})(1_{A}\otimes x_{1})+ \\
&&(1_{A}\otimes x_{1})B(gx_{1}\otimes 1_{H})(1_{A}\otimes x_{2})+ \\
&&+(1_{A}\otimes gx_{1}x_{2})B(gx_{1}\otimes 1_{H})(1_{A}\otimes g)+ \\
&&-(1_{A}\otimes gx_{1})B(x_{1}x_{2}\otimes 1_{H})(1_{A}\otimes g)
\end{eqnarray*}%
\begin{eqnarray*}
&&\text{simplifying right side}(1_{A}\otimes x_{1})B(g\otimes
gx_{1}x_{2})+B(x_{1}\otimes gx_{1}x_{2}) \\
&=&+(1_{A}\otimes g)B(gx_{1}\otimes 1_{H})(1_{A}\otimes gx_{1}x_{2}) \\
&&-(1_{A}\otimes x_{2})B(gx_{1}\otimes 1_{H})(1_{A}\otimes x_{1}) \\
&&+B(x_{1}x_{2}\otimes 1_{H})(1_{A}\otimes x_{1})+
\end{eqnarray*}

and we conclude.

\subsubsection{$1_{H}\otimes gx_{1}$}

\begin{equation*}
B(1_{H}\otimes gx_{1})(1_{A}\otimes x_{1})\overset{?}{=}(1_{A}\otimes
x_{1})B(g\otimes ggx_{1})+B(x_{1}\otimes ggx_{1})+B(1_{H}\otimes x_{1}gx_{1})
\end{equation*}%
i.e.%
\begin{equation*}
B(1_{H}\otimes gx_{1})(1_{A}\otimes x_{1})\overset{?}{=}(1_{A}\otimes
x_{1})B(g\otimes x_{1})+B(x_{1}\otimes x_{1})
\end{equation*}%
We consider the left side%
\begin{equation*}
\text{left side }B(1_{H}\otimes gx_{1})(1_{A}\otimes x_{1})
\end{equation*}%
\begin{eqnarray*}
\text{left side} &&B(1_{H}\otimes gx_{1})(1_{A}\otimes x_{1})\overset{\left( %
\ref{form 1otgx1}\right) }{=} \\
&=&(1_{A}\otimes g)B(g\otimes 1_{H})(1_{A}\otimes gx_{1})(1_{A}\otimes x_{1})
\\
&&+(1_{A}\otimes x_{1})B(g\otimes 1_{H})(1_{A}\otimes x_{1})+B(x_{1}\otimes
1_{H})(1_{A}\otimes x_{1})
\end{eqnarray*}%
Thus we get%
\begin{equation*}
\text{left side }B(1_{H}\otimes gx_{1})(1_{A}\otimes x_{1})=B(x_{1}\otimes
1_{H})(1_{A}\otimes x_{1})
\end{equation*}

We consider now the right side%
\begin{eqnarray*}
&&\text{right side }(1_{A}\otimes x_{1})B(g\otimes x_{1})+B(x_{1}\otimes
x_{1}) \\
&=&(1_{A}\otimes x_{1})(1_{A}\otimes g)B(1_{H}\otimes gx_{1})(1_{A}\otimes
g)+B(x_{1}\otimes x_{1}) \\
&=&(1_{A}\otimes gx_{1})B(1_{H}\otimes gx_{1})(1_{A}\otimes
g)+B(x_{1}\otimes x_{1})
\end{eqnarray*}%
\begin{eqnarray*}
&&(1_{A}\otimes gx_{1})B(1_{H}\otimes gx_{1})(1_{A}\otimes g)\overset{\left( %
\ref{form 1otgx1}\right) }{=} \\
&=&(1_{A}\otimes gx_{1})(1_{A}\otimes g)B(g\otimes 1_{H})(1_{A}\otimes
gx_{1})(1_{A}\otimes g) \\
&&+(1_{A}\otimes gx_{1})(1_{A}\otimes x_{1})B(g\otimes 1_{H})(1_{A}\otimes g)
\\
&&+(1_{A}\otimes gx_{1})B(x_{1}\otimes 1_{H})(1_{A}\otimes g)
\end{eqnarray*}%
Thus%
\begin{eqnarray*}
&&(1_{A}\otimes gx_{1})B(1_{H}\otimes gx_{1})(1_{A}\otimes g)\overset{\left( %
\ref{form 1otgx1}\right) }{=} \\
&=&(1_{A}\otimes gx_{1})B(x_{1}\otimes 1_{H})(1_{A}\otimes g)
\end{eqnarray*}%
\begin{equation*}
B(x_{1}\otimes x_{1})\overset{\left( \ref{form x1otx1}\right) }{=}%
B(x_{1}\otimes 1_{H})(1_{A}\otimes x_{1})-(1_{A}\otimes
gx_{1})B(x_{1}\otimes 1_{H})(1_{A}\otimes g)
\end{equation*}%
Thus we get%
\begin{eqnarray*}
\text{right side}(1_{A}\otimes x_{1})B(g\otimes x_{1})+B(x_{1}\otimes x_{1})
&=&(1_{A}\otimes gx_{1})B(x_{1}\otimes 1_{H})(1_{A}\otimes g)+ \\
&&+B(x_{1}\otimes 1_{H})(1_{A}\otimes x_{1})+ \\
&&-(1_{A}\otimes gx_{1})B(x_{1}\otimes 1_{H})(1_{A}\otimes g) \\
&=&B(x_{1}\otimes 1_{H})(1_{A}\otimes x_{1})
\end{eqnarray*}%
and we conclude.

\begin{equation*}
B(h\otimes h^{\prime })(1_{A}\otimes x_{1})\overset{?}{=}(1_{A}\otimes
x_{1})B(gh\otimes gh^{\prime })+B(x_{1}h\otimes gh^{\prime })+B(h\otimes
x_{1}h^{\prime })
\end{equation*}

\subsubsection{$1_{H}\otimes gx_{2}$}

\begin{equation*}
B(1_{H}\otimes gx_{2})(1_{A}\otimes x_{1})\overset{?}{=}(1_{A}\otimes
x_{1})B(g\otimes x_{2})+B(x_{1}\otimes x_{2})-B(1_{H}\otimes gx_{1}x_{2})
\end{equation*}

We consider the left side%
\begin{eqnarray*}
\text{left side} &&B(1_{H}\otimes gx_{2})(1_{A}\otimes x_{1})\overset{\left( %
\ref{form 1otgx2}\right) }{=} \\
&=&-(1_{A}\otimes g)B(g\otimes 1_{H})(1_{A}\otimes gx_{1}x_{2}) \\
&&+(1_{A}\otimes x_{2})B(g\otimes 1_{H})(1_{A}\otimes x_{1}) \\
&&+B(x_{2}\otimes 1_{H})(1_{A}\otimes x_{1})
\end{eqnarray*}%
\begin{equation*}
\text{right side }(1_{A}\otimes x_{1})B(g\otimes x_{2})+B(x_{1}\otimes
x_{2})-B(1_{H}\otimes gx_{1}x_{2})
\end{equation*}%
\begin{eqnarray*}
(1_{A}\otimes x_{1})B(g\otimes x_{2}) &=&(1_{A}\otimes x_{1})(1_{A}\otimes
g)B(1_{H}\otimes gx_{2})(1_{A}\otimes g) \\
&=&(1_{A}\otimes gx_{1})B(1_{H}\otimes gx_{2})(1_{A}\otimes g)
\end{eqnarray*}%
\begin{eqnarray*}
&&(1_{A}\otimes x_{1})B(g\otimes x_{2})=(1_{A}\otimes gx_{1})B(1_{H}\otimes
gx_{2})(1_{A}\otimes g)\overset{\left( \ref{form 1otgx2}\right) }{=} \\
&=&(1_{A}\otimes gx_{1})(1_{A}\otimes g)B(g\otimes 1_{H})(1_{A}\otimes
gx_{2})(1_{A}\otimes g) \\
&&+(1_{A}\otimes gx_{1})(1_{A}\otimes x_{2})B(g\otimes 1_{H})(1_{A}\otimes g)
\\
&&+(1_{A}\otimes gx_{1})B(x_{2}\otimes 1_{H})(1_{A}\otimes g)
\end{eqnarray*}%
\begin{eqnarray*}
&&(1_{A}\otimes x_{1})B(g\otimes x_{2})=(1_{A}\otimes gx_{1})B(1_{H}\otimes
gx_{2})(1_{A}\otimes g)\overset{\left( \ref{form 1otgx2}\right) }{=} \\
&=&(1_{A}\otimes x_{1})B(g\otimes 1_{H})(1_{A}\otimes x_{2}) \\
&&+(1_{A}\otimes gx_{1}x_{2})B(g\otimes 1_{H})(1_{A}\otimes g) \\
&&+(1_{A}\otimes gx_{1})B(x_{2}\otimes 1_{H})(1_{A}\otimes g)
\end{eqnarray*}%
\begin{eqnarray*}
&&-B(1_{H}\otimes gx_{1}x_{2})\overset{\left( \ref{form 1otgx1x2}\right) }{=}%
-(1_{A}\otimes gx_{1})B(1_{H}\otimes gx_{2})\left( 1_{A}\otimes g\right) \\
&&-B(x_{1}\otimes x_{2})+B(1_{H}\otimes gx_{2})(1_{A}\otimes x_{1})
\end{eqnarray*}%
\begin{eqnarray*}
-(1_{A}\otimes gx_{1}) &&B(1_{H}\otimes gx_{2})\left( 1_{A}\otimes g\right)
\\
&&\overset{\left( \ref{form 1otgx2}\right) }{=}-(1_{A}\otimes
gx_{1})(1_{A}\otimes g)B(g\otimes 1_{H})(1_{A}\otimes gx_{2})\left(
1_{A}\otimes g\right) \\
&&+-(1_{A}\otimes gx_{1})(1_{A}\otimes x_{2})B(g\otimes 1_{H})\left(
1_{A}\otimes g\right) \\
&&+-(1_{A}\otimes gx_{1})B(x_{2}\otimes 1_{H})\left( 1_{A}\otimes g\right)
\end{eqnarray*}%
\begin{eqnarray*}
-(1_{A}\otimes gx_{1}) &&B(1_{H}\otimes gx_{2})\left( 1_{A}\otimes g\right)
\\
&&\overset{\left( \ref{form 1otgx2}\right) }{=}-(1_{A}\otimes
x_{1})B(g\otimes 1_{H})(1_{A}\otimes x_{2}) \\
&&-(1_{A}\otimes gx_{1}x_{2})B(g\otimes 1_{H})\left( 1_{A}\otimes g\right) \\
&&-(1_{A}\otimes gx_{1})B(x_{2}\otimes 1_{H})\left( 1_{A}\otimes g\right)
\end{eqnarray*}%
\begin{eqnarray*}
&&B(1_{H}\otimes gx_{2})(1_{A}\otimes x_{1}) \\
&&\overset{\left( \ref{form 1otgx2}\right) }{=}-(1_{A}\otimes g)B(g\otimes
1_{H})(1_{A}\otimes gx_{1}x_{2}) \\
+ &&(1_{A}\otimes x_{2})B(g\otimes 1_{H})(1_{A}\otimes x_{1}) \\
&&+B(x_{2}\otimes 1_{H})(1_{A}\otimes x_{1})
\end{eqnarray*}%
\begin{eqnarray*}
&&-B(1_{H}\otimes gx_{1}x_{2})\overset{\left( \ref{form 1otgx1x2}\right) }{=}%
-(1_{A}\otimes gx_{1})B(1_{H}\otimes gx_{2})\left( 1_{A}\otimes g\right) \\
&&-B(x_{1}\otimes x_{2})+B(1_{H}\otimes gx_{2})(1_{A}\otimes x_{1}) \\
&=&-(1_{A}\otimes x_{1})B(g\otimes 1_{H})(1_{A}\otimes x_{2}) \\
&&-(1_{A}\otimes gx_{1}x_{2})B(g\otimes 1_{H})\left( 1_{A}\otimes g\right) \\
&&-(1_{A}\otimes gx_{1})B(x_{2}\otimes 1_{H})\left( 1_{A}\otimes g\right) \\
&&-B(x_{1}\otimes x_{2})+ \\
&&-(1_{A}\otimes g)B(g\otimes 1_{H})(1_{A}\otimes gx_{1}x_{2})+ \\
&&+(1_{A}\otimes x_{2})B(g\otimes 1_{H})(1_{A}\otimes x_{1})+ \\
&&+B(x_{2}\otimes 1_{H})(1_{A}\otimes x_{1})
\end{eqnarray*}%
\begin{eqnarray*}
&&\text{right side }(1_{A}\otimes x_{1})B(g\otimes x_{2})+B(x_{1}\otimes
x_{2})-B(1_{H}\otimes gx_{1}x_{2}) \\
&=&(1_{A}\otimes x_{1})B(g\otimes 1_{H})(1_{A}\otimes x_{2})+ \\
&&+(1_{A}\otimes gx_{1}x_{2})B(g\otimes 1_{H})(1_{A}\otimes g) \\
&&+(1_{A}\otimes gx_{1})B(x_{2}\otimes 1_{H})(1_{A}\otimes g)+ \\
&&+B(x_{1}\otimes x_{2})+ \\
&&-(1_{A}\otimes x_{1})B(g\otimes 1_{H})(1_{A}\otimes x_{2}) \\
&&-(1_{A}\otimes gx_{1}x_{2})B(g\otimes 1_{H})\left( 1_{A}\otimes g\right) \\
&&-(1_{A}\otimes gx_{1})B(x_{2}\otimes 1_{H})\left( 1_{A}\otimes g\right) \\
&&-B(x_{1}\otimes x_{2})+ \\
&&-(1_{A}\otimes g)B(g\otimes 1_{H})(1_{A}\otimes gx_{1}x_{2})+ \\
&&+(1_{A}\otimes x_{2})B(g\otimes 1_{H})(1_{A}\otimes x_{1})+ \\
&&+B(x_{2}\otimes 1_{H})(1_{A}\otimes x_{1})
\end{eqnarray*}%
\begin{eqnarray*}
\text{clean} &&\text{right side }(1_{A}\otimes x_{1})B(g\otimes
x_{2})+B(x_{1}\otimes x_{2})-B(1_{H}\otimes gx_{1}x_{2}) \\
&=&-(1_{A}\otimes g)B(g\otimes 1_{H})(1_{A}\otimes gx_{1}x_{2})+ \\
&&+(1_{A}\otimes x_{2})B(g\otimes 1_{H})(1_{A}\otimes x_{1})+ \\
&&+B(x_{2}\otimes 1_{H})(1_{A}\otimes x_{1})
\end{eqnarray*}%
and we conclude.

\subsubsection{$1_{H}\otimes gx_{1}x_{2}$}

\begin{equation*}
B(1_{H}\otimes gx_{1}x_{2})(1_{A}\otimes x_{1})\overset{?}{=}(1_{A}\otimes
x_{1})B(g\otimes ggx_{1}x_{2})+B(x_{1}\otimes ggx_{1}x_{2})+B(1_{H}\otimes
x_{1}gx_{1}x_{2})
\end{equation*}%
\begin{equation*}
B(1_{H}\otimes gx_{1}x_{2})(1_{A}\otimes x_{1})\overset{?}{=}(1_{A}\otimes
x_{1})B(g\otimes x_{1}x_{2})+B(x_{1}\otimes x_{1}x_{2})
\end{equation*}

\begin{eqnarray*}
&&\text{left side }B(1_{H}\otimes gx_{1}x_{2})(1_{A}\otimes x_{1})\overset{%
\left( \ref{form 1otgx1x2}\right) }{=} \\
&&(1_{A}\otimes gx_{1})B(1_{H}\otimes gx_{2})\left( 1_{A}\otimes g\right)
(1_{A}\otimes x_{1}) \\
&&-B(1_{H}\otimes gx_{2})(1_{A}\otimes x_{1})(1_{A}\otimes x_{1}) \\
&=&-(1_{A}\otimes gx_{1})B(1_{H}\otimes gx_{2})\left( 1_{A}\otimes
gx_{1}\right) +B(x_{1}\otimes x_{2})(1_{A}\otimes x_{1})
\end{eqnarray*}%
\begin{eqnarray*}
\text{left side} &&B(1_{H}\otimes gx_{1}x_{2})(1_{A}\otimes x_{1})\overset{%
\left( \ref{form 1otgx1x2}\right) }{=} \\
&&(1_{A}\otimes gx_{1})B(1_{H}\otimes gx_{2})\left( 1_{A}\otimes g\right)
(1_{A}\otimes x_{1}) \\
&&+B(x_{1}\otimes x_{2})(1_{A}\otimes x_{1})-B(1_{H}\otimes
gx_{2})(1_{A}\otimes x_{1})(1_{A}\otimes x_{1}) \\
&=&+(1_{A}\otimes x_{1})B(g\otimes 1_{H})(1_{A}\otimes x_{1}x_{2}) \\
&&-(1_{A}\otimes gx_{1}x_{2})B(g\otimes 1_{H})\left( 1_{A}\otimes
gx_{1}\right) \\
&&-(1_{A}\otimes gx_{1})B(x_{2}\otimes 1_{H})\left( 1_{A}\otimes
gx_{1}\right) \\
&&+B(x_{1}\otimes x_{2})(1_{A}\otimes x_{1})+
\end{eqnarray*}%
\begin{gather*}
B(g\otimes x_{1}x_{2})=\left( 1_{A}\otimes g\right) B(1_{H}\otimes
gx_{1}x_{2})\left( 1_{A}\otimes g\right) \overset{\left( \ref{form 1otgx1x2}%
\right) }{=} \\
=+\left( 1_{A}\otimes g\right) (1_{A}\otimes x_{1})B(g\otimes
1_{H})(1_{A}\otimes x_{2})\left( 1_{A}\otimes g\right) \\
+\left( 1_{A}\otimes g\right) (1_{A}\otimes gx_{1}x_{2})B(g\otimes
1_{H})\left( 1_{A}\otimes g\right) \left( 1_{A}\otimes g\right) \\
+\left( 1_{A}\otimes g\right) (1_{A}\otimes gx_{1})B(x_{2}\otimes
1_{H})\left( 1_{A}\otimes g\right) \left( 1_{A}\otimes g\right) \\
+\left( 1_{A}\otimes g\right) B(x_{1}\otimes x_{2})\left( 1_{A}\otimes
g\right) + \\
+\left( 1_{A}\otimes g\right) (1_{A}\otimes g)B(g\otimes 1_{H})(1_{A}\otimes
gx_{1}x_{2})\left( 1_{A}\otimes g\right) + \\
-\left( 1_{A}\otimes g\right) (1_{A}\otimes x_{2})B(g\otimes
1_{H})(1_{A}\otimes x_{1})\left( 1_{A}\otimes g\right) + \\
-\left( 1_{A}\otimes g\right) B(x_{2}\otimes 1_{H})(1_{A}\otimes
x_{1})\left( 1_{A}\otimes g\right)
\end{gather*}%
\begin{gather*}
(1_{A}\otimes x_{1})B(g\otimes x_{1}x_{2})=(1_{A}\otimes x_{1})\left(
1_{A}\otimes g\right) B(1_{H}\otimes gx_{1}x_{2})\left( 1_{A}\otimes
g\right) \overset{\left( \ref{form 1otgx1x2}\right) }{=} \\
=+(1_{A}\otimes x_{1})\left( 1_{A}\otimes g\right) (1_{A}\otimes
x_{1})B(g\otimes 1_{H})(1_{A}\otimes x_{2})\left( 1_{A}\otimes g\right) \\
+(1_{A}\otimes x_{1})\left( 1_{A}\otimes g\right) (1_{A}\otimes
gx_{1}x_{2})B(g\otimes 1_{H})\left( 1_{A}\otimes g\right) \left(
1_{A}\otimes g\right) \\
+(1_{A}\otimes x_{1})\left( 1_{A}\otimes g\right) (1_{A}\otimes
gx_{1})B(x_{2}\otimes 1_{H})\left( 1_{A}\otimes g\right) \left( 1_{A}\otimes
g\right) \\
+(1_{A}\otimes x_{1})\left( 1_{A}\otimes g\right) B(x_{1}\otimes
x_{2})\left( 1_{A}\otimes g\right) + \\
+(1_{A}\otimes x_{1})\left( 1_{A}\otimes g\right) (1_{A}\otimes g)B(g\otimes
1_{H})(1_{A}\otimes gx_{1}x_{2})\left( 1_{A}\otimes g\right) + \\
-(1_{A}\otimes x_{1})\left( 1_{A}\otimes g\right) (1_{A}\otimes
x_{2})B(g\otimes 1_{H})(1_{A}\otimes x_{1})\left( 1_{A}\otimes g\right) + \\
-(1_{A}\otimes x_{1})\left( 1_{A}\otimes g\right) B(x_{2}\otimes
1_{H})(1_{A}\otimes x_{1})\left( 1_{A}\otimes g\right)
\end{gather*}%
\begin{gather*}
(1_{A}\otimes x_{1})B(g\otimes x_{1}x_{2})=(1_{A}\otimes x_{1})\left(
1_{A}\otimes g\right) B(1_{H}\otimes gx_{1}x_{2})\left( 1_{A}\otimes
g\right) \overset{\left( \ref{form 1otgx1x2}\right) }{=} \\
+(1_{A}\otimes gx_{1})B(x_{1}\otimes x_{2})\left( 1_{A}\otimes g\right) + \\
+(1_{A}\otimes x_{1})B(g\otimes 1_{H})(1_{A}\otimes x_{1}x_{2})+ \\
-(1_{A}\otimes gx_{1}x_{2})B(g\otimes 1_{H})(1_{A}\otimes gx_{1})+ \\
-(1_{A}\otimes gx_{1})B(x_{2}\otimes 1_{H})(1_{A}\otimes gx_{1})
\end{gather*}%
\begin{equation*}
B(x_{1}\otimes x_{1}x_{2})\overset{\left( \ref{form x1otx1x2}\right) }{=}%
B(x_{1}\otimes x_{2})(1_{A}\otimes x_{1})-(1_{A}\otimes
gx_{1})B(x_{1}\otimes x_{2})(1_{A}\otimes g)
\end{equation*}%
\begin{gather*}
\text{right side }(1_{A}\otimes x_{1})B(g\otimes x_{1}x_{2})+B(x_{1}\otimes
x_{1}x_{2}) \\
=+(1_{A}\otimes gx_{1})B(x_{1}\otimes x_{2})\left( 1_{A}\otimes g\right) + \\
+(1_{A}\otimes x_{1})B(g\otimes 1_{H})(1_{A}\otimes x_{1}x_{2})+ \\
-(1_{A}\otimes gx_{1}x_{2})B(g\otimes 1_{H})(1_{A}\otimes gx_{1})+ \\
-(1_{A}\otimes gx_{1})B(x_{2}\otimes 1_{H})(1_{A}\otimes gx_{1})+ \\
+B(x_{1}\otimes x_{2})(1_{A}\otimes x_{1})-(1_{A}\otimes
gx_{1})B(x_{1}\otimes x_{2})(1_{A}\otimes g)
\end{gather*}%
\begin{gather*}
\text{clean right side }(1_{A}\otimes x_{1})B(g\otimes
x_{1}x_{2})+B(x_{1}\otimes x_{1}x_{2}) \\
=+(1_{A}\otimes x_{1})B(g\otimes 1_{H})(1_{A}\otimes x_{1}x_{2})+ \\
-(1_{A}\otimes gx_{1}x_{2})B(g\otimes 1_{H})(1_{A}\otimes gx_{1})+ \\
-(1_{A}\otimes gx_{1})B(x_{2}\otimes 1_{H})(1_{A}\otimes gx_{1})+ \\
+B(x_{1}\otimes x_{2})(1_{A}\otimes x_{1})
\end{gather*}%
and we conclude.

\subsubsection{$x_{1}\otimes x_{1}$}

\begin{equation*}
B(x_{1}\otimes x_{1})(1_{A}\otimes x_{1})\overset{?}{=}(1_{A}\otimes
x_{1})B(gx_{1}\otimes gx_{1})+B(x_{1}x_{1}\otimes gx_{1})+B(x_{1}\otimes
x_{1}x_{1})
\end{equation*}%
\begin{equation*}
B(x_{1}\otimes x_{1})(1_{A}\otimes x_{1})\overset{?}{=}(1_{A}\otimes
x_{1})(1_{A}\otimes g)B(x_{1}\otimes x_{1})(1_{A}\otimes g)+
\end{equation*}%
\begin{equation*}
B(x_{1}\otimes x_{1})\overset{\left( \ref{form x1otx1}\right) }{=}%
B(x_{1}\otimes 1_{H})(1_{A}\otimes x_{1})-(1_{A}\otimes
gx_{1})B(x_{1}\otimes 1_{H})(1_{A}\otimes g)
\end{equation*}%
\begin{eqnarray*}
&&\text{left side }B(x_{1}\otimes x_{1})(1_{A}\otimes x_{1}) \\
&=&B(x_{1}\otimes 1_{H})(1_{A}\otimes x_{1})(1_{A}\otimes
x_{1})-(1_{A}\otimes gx_{1})B(x_{1}\otimes 1_{H})(1_{A}\otimes
g)(1_{A}\otimes x_{1}) \\
&=&(1_{A}\otimes gx_{1})B(x_{1}\otimes 1_{H})(1_{A}\otimes gx_{1})
\end{eqnarray*}%
\begin{eqnarray*}
&&\text{right side }(1_{A}\otimes x_{1})(1_{A}\otimes g)B(x_{1}\otimes
x_{1})(1_{A}\otimes g) \\
&=&(1_{A}\otimes x_{1})(1_{A}\otimes g)B(x_{1}\otimes 1_{H})(1_{A}\otimes
x_{1})(1_{A}\otimes g) \\
&&-(1_{A}\otimes x_{1})(1_{A}\otimes g)(1_{A}\otimes gx_{1})B(x_{1}\otimes
1_{H})(1_{A}\otimes g)(1_{A}\otimes g) \\
&=&(1_{A}\otimes x_{1})(1_{A}\otimes g)B(x_{1}\otimes 1_{H})(1_{A}\otimes
x_{1})(1_{A}\otimes g)=(1_{A}\otimes gx_{1})B(x_{1}\otimes
1_{H})(1_{A}\otimes gx_{1})
\end{eqnarray*}%
and we conclude.

\subsubsection{$x_{1}\otimes x_{2}$}

\begin{equation*}
B(x_{1}\otimes x_{2})(1_{A}\otimes x_{1})\overset{?}{=}(1_{A}\otimes
x_{1})B(gx_{1}\otimes gx_{2})+B(x_{1}x_{1}\otimes gx_{2})+B(x_{1}\otimes
x_{1}x_{2})
\end{equation*}%
\begin{equation*}
B(x_{1}\otimes x_{2})(1_{A}\otimes x_{1})\overset{?}{=}(1_{A}\otimes
x_{1})B(gx_{1}\otimes gx_{2})+B(x_{1}\otimes x_{1}x_{2})
\end{equation*}%
\begin{equation*}
B(x_{1}\otimes x_{1}x_{2})\overset{\left( \ref{form x1otx1x2}\right) }{=}%
B(x_{1}\otimes x_{2})(1_{A}\otimes x_{1})-(1_{A}\otimes
gx_{1})B(x_{1}\otimes x_{2})(1_{A}\otimes g)
\end{equation*}%
$B(x_{1}\otimes x_{1}x_{2})=B(x_{1}\otimes x_{2})(1_{A}\otimes
x_{1})-(1_{A}\otimes gx_{1})B(x_{1}\otimes x_{2})(1_{A}\otimes g)$%
\begin{eqnarray*}
&&\text{right side }(1_{A}\otimes x_{1})B(gx_{1}\otimes
gx_{2})+B(x_{1}\otimes x_{1}x_{2})\overset{\left( \ref{form x1otx1x2}\right)
}{=} \\
&&(1_{A}\otimes gx_{1})B(x_{1}\otimes x_{2})(1_{A}\otimes g)+ \\
+ &&B(x_{1}\otimes x_{2})(1_{A}\otimes x_{1})-(1_{A}\otimes
gx_{1})B(x_{1}\otimes x_{2})(1_{A}\otimes g) \\
&=&B(x_{1}\otimes x_{2})(1_{A}\otimes x_{1})
\end{eqnarray*}%
and we conclude.

\subsubsection{$x_{1}\otimes x_{1}x_{2}$}

\begin{equation*}
B(x_{1}\otimes x_{1}x_{2})(1_{A}\otimes x_{1})\overset{?}{=}(1_{A}\otimes
x_{1})B(gx_{1}\otimes gx_{1}x_{2})+B(x_{1}x_{1}\otimes
gx_{1}x_{2})+B(x_{1}\otimes x_{1}x_{1}x_{2})
\end{equation*}%
\begin{equation*}
B(x_{1}\otimes x_{1}x_{2})(1_{A}\otimes x_{1})\overset{?}{=}(1_{A}\otimes
x_{1})B(gx_{1}\otimes gx_{1}x_{2})
\end{equation*}%
\begin{equation*}
B(x_{1}\otimes x_{1}x_{2})\overset{\left( \ref{form x1otx1x2}\right) }{=}%
B(x_{1}\otimes x_{2})(1_{A}\otimes x_{1})-(1_{A}\otimes
gx_{1})B(x_{1}\otimes x_{2})(1_{A}\otimes g)
\end{equation*}%
\begin{equation*}
\text{left side }B(x_{1}\otimes x_{1}x_{2})(1_{A}\otimes
x_{1})=B(x_{1}\otimes x_{2})(1_{A}\otimes x_{1})(1_{A}\otimes
x_{1})-(1_{A}\otimes gx_{1})B(x_{1}\otimes x_{2})(1_{A}\otimes x_{1})
\end{equation*}%
and hence we get%
\begin{equation*}
\text{left side }B(x_{1}\otimes x_{1}x_{2})(1_{A}\otimes
x_{1})=(1_{A}\otimes gx_{1})B(x_{1}\otimes x_{2})(1_{A}\otimes gx_{1})
\end{equation*}%
\begin{gather*}
\text{right side }B(gx_{1}\otimes gx_{1}x_{2})=(1_{A}\otimes
x_{1})(1_{A}\otimes g)B(x_{1}\otimes x_{1}x_{2})(1_{A}\otimes g)= \\
=(1_{A}\otimes x_{1})(1_{A}\otimes g)B(x_{1}\otimes x_{2})(1_{A}\otimes
x_{1})(1_{A}\otimes g) \\
-(1_{A}\otimes x_{1})(1_{A}\otimes g)(1_{A}\otimes gx_{1})B(x_{1}\otimes
x_{2})(1_{A}\otimes g)(1_{A}\otimes g)= \\
=(1_{A}\otimes gx_{1})B(x_{1}\otimes x_{2})(1_{A}\otimes gx_{1})
\end{gather*}%
and we conclude.

\subsubsection{$x_{1}\otimes gx_{1}$}

\begin{equation*}
B(x_{1}\otimes gx_{1})(1_{A}\otimes x_{1})\overset{?}{=}(1_{A}\otimes
x_{1})B(gx_{1}\otimes ggx_{1})+B(x_{1}x_{1}\otimes ggx_{1})+B(x_{1}\otimes
x_{1}gx_{1})
\end{equation*}

i.e.%
\begin{equation*}
B(x_{1}\otimes gx_{1})(1_{A}\otimes x_{1})\overset{?}{=}(1_{A}\otimes
x_{1})B(gx_{1}\otimes x_{1})
\end{equation*}

\begin{equation*}
B(x_{1}\otimes gx_{1})\overset{\left( \ref{form x1otgx1}\right) }{=}%
(1_{A}\otimes g)B(gx_{1}\otimes 1_{H})(1_{A}\otimes gx_{1})+(1_{A}\otimes
x_{1})B(gx_{1}\otimes 1_{H})
\end{equation*}%
\begin{equation*}
\text{left side }B(x_{1}\otimes gx_{1})(1_{A}\otimes x_{1})=(1_{A}\otimes
g)B(gx_{1}\otimes 1_{H})(1_{A}\otimes gx_{1})(1_{A}\otimes
x_{1})+(1_{A}\otimes x_{1})B(gx_{1}\otimes 1_{H})(1_{A}\otimes x_{1})
\end{equation*}%
i.e.%
\begin{equation*}
\text{left side }B(x_{1}\otimes gx_{1})(1_{A}\otimes x_{1})=(1_{A}\otimes
x_{1})B(gx_{1}\otimes 1_{H})(1_{A}\otimes x_{1})
\end{equation*}%
\begin{eqnarray*}
&&B(x_{1}\otimes gx_{1})\overset{\left( \ref{form x1otgx1}\right) }{=}%
(1_{A}\otimes g)B(gx_{1}\otimes 1_{H})(1_{A}\otimes gx_{1})+ \\
&&+(1_{A}\otimes x_{1})B(gx_{1}\otimes 1_{H})
\end{eqnarray*}%
\begin{eqnarray*}
\text{right side }(1_{A}\otimes x_{1})B(gx_{1}\otimes x_{1})
&=&(1_{A}\otimes x_{1})(1_{A}\otimes g)B(x_{1}\otimes gx_{1})(1_{A}\otimes g)
\\
&=&(1_{A}\otimes x_{1})(1_{A}\otimes g)(1_{A}\otimes g)B(gx_{1}\otimes
1_{H})(1_{A}\otimes gx_{1})(1_{A}\otimes g) \\
&&+(1_{A}\otimes x_{1})(1_{A}\otimes g)(1_{A}\otimes x_{1})B(gx_{1}\otimes
1_{H})(1_{A}\otimes g) \\
&=&(1_{A}\otimes x_{1})B(gx_{1}\otimes 1_{H})(1_{A}\otimes x_{1})
\end{eqnarray*}%
and we conclude.

\subsubsection{$x_{1}\otimes gx_{2}$}

\begin{equation*}
B(x_{1}\otimes gx_{2})(1_{A}\otimes x_{1})\overset{?}{=}(1_{A}\otimes
x_{1})B(gx_{1}\otimes ggx_{2})+B(x_{1}x_{1}\otimes gh^{\prime
})+B(x_{1}\otimes x_{1}gx_{2})
\end{equation*}%
\begin{equation*}
B(x_{1}\otimes gx_{2})(1_{A}\otimes x_{1})\overset{?}{=}(1_{A}\otimes
x_{1})B(gx_{1}\otimes x_{2})-B(x_{1}\otimes gx_{1}x_{2})
\end{equation*}%
\begin{equation*}
B(x_{1}\otimes gx_{1}x_{2})\overset{\left( \ref{form x1otgx1x2}\right) }{=}%
-B(x_{1}\otimes gx_{2})(1_{A}\otimes x_{1})+(1_{A}\otimes
gx_{1})B(x_{1}\otimes gx_{2})(1_{A}\otimes g)
\end{equation*}%
\begin{eqnarray*}
&&\text{the right side}(1_{A}\otimes x_{1})B(gx_{1}\otimes
x_{2})-B(x_{1}\otimes gx_{1}x_{2}) \\
&=&(1_{A}\otimes x_{1})B(gx_{1}\otimes x_{2})+B(x_{1}\otimes
gx_{2})(1_{A}\otimes x_{1})-(1_{A}\otimes gx_{1})B(x_{1}\otimes
gx_{2})(1_{A}\otimes g) \\
&=&(1_{A}\otimes x_{1})(1_{A}\otimes g)B(x_{1}\otimes gx_{2})(1_{A}\otimes
g)+B(x_{1}\otimes gx_{2})(1_{A}\otimes x_{1})-(1_{A}\otimes
gx_{1})B(x_{1}\otimes gx_{2})(1_{A}\otimes g)
\end{eqnarray*}%
\begin{eqnarray*}
\text{clean} &&\text{the right side}(1_{A}\otimes x_{1})B(gx_{1}\otimes
x_{2})-B(x_{1}\otimes gx_{1}x_{2}) \\
&=&(1_{A}\otimes x_{1})B(gx_{1}\otimes x_{2})+B(x_{1}\otimes
gx_{2})(1_{A}\otimes x_{1})-(1_{A}\otimes gx_{1})B(x_{1}\otimes
gx_{2})(1_{A}\otimes g) \\
&=&+B(x_{1}\otimes gx_{2})(1_{A}\otimes x_{1})
\end{eqnarray*}%
and we conclude

\subsubsection{$x_{1}\otimes gx_{1}x_{2}$}

\begin{equation*}
B(x_{1}\otimes gx_{1}x_{2})(1_{A}\otimes x_{1})\overset{?}{=}(1_{A}\otimes
x_{1})B(gx_{1}\otimes ggx_{1}x_{2})+B(x_{1}x_{1}\otimes
ggx_{1}x_{2})+B(x_{1}\otimes x_{1}gx_{1}x_{2})
\end{equation*}%
i.e.%
\begin{equation*}
B(x_{1}\otimes gx_{1}x_{2})(1_{A}\otimes x_{1})\overset{?}{=}(1_{A}\otimes
x_{1})(1_{A}\otimes g)B(x_{1}\otimes gx_{1}x_{2})(1_{A}\otimes g)
\end{equation*}

\begin{equation*}
B(x_{1}\otimes gx_{1}x_{2})\overset{\left( \ref{form x1otgx1x2}\right) }{=}%
-B(x_{1}\otimes gx_{2})(1_{A}\otimes x_{1})+(1_{A}\otimes
gx_{1})B(x_{1}\otimes gx_{2})(1_{A}\otimes g)
\end{equation*}%
\begin{eqnarray*}
&&\text{the left side }B(x_{1}\otimes gx_{1}x_{2})(1_{A}\otimes x_{1}) \\
&=&-B(x_{1}\otimes gx_{2})(1_{A}\otimes x_{1})(1_{A}\otimes
x_{1})+(1_{A}\otimes gx_{1})B(x_{1}\otimes gx_{2})(1_{A}\otimes
g)(1_{A}\otimes x_{1}) \\
&=&-(1_{A}\otimes gx_{1})B(x_{1}\otimes gx_{2})(1_{A}\otimes gx_{1})
\end{eqnarray*}%
\begin{eqnarray*}
\text{the right side } &=&(1_{A}\otimes gx_{1})B(x_{1}\otimes
gx_{1}x_{2})(1_{A}\otimes g) \\
&=&-(1_{A}\otimes gx_{1})B(x_{1}\otimes gx_{2})(1_{A}\otimes
x_{1})(1_{A}\otimes g) \\
&&+(1_{A}\otimes gx_{1})(1_{A}\otimes gx_{1})B(x_{1}\otimes
gx_{2})(1_{A}\otimes g)(1_{A}\otimes g) \\
&=&-(1_{A}\otimes gx_{1})B(x_{1}\otimes gx_{2})(1_{A}\otimes gx_{1})
\end{eqnarray*}%
and we conclude.

\subsubsection{$x_{2}\otimes x_{1}$}

\begin{equation*}
B(x_{2}\otimes x_{1})(1_{A}\otimes x_{1})\overset{?}{=}(1_{A}\otimes
x_{1})B(gx_{2}\otimes gx_{1})+B(x_{1}x_{2}\otimes gx_{1})+B(x_{2}\otimes
x_{1}x_{1})
\end{equation*}

\begin{equation*}
B(x_{2}\otimes x_{1})(1_{A}\otimes x_{1})\overset{?}{=}(1_{A}\otimes
x_{1})B(gx_{2}\otimes gx_{1})+B(x_{1}x_{2}\otimes gx_{1})
\end{equation*}

\begin{equation*}
B(x_{2}\otimes x_{1})\overset{\left( \ref{form x2otx1}\right) }{=}%
B(x_{2}\otimes 1_{H})(1_{A}\otimes x_{1})-(1_{A}\otimes
gx_{1})B(x_{2}\otimes 1_{H})(1_{A}\otimes g)-(1_{A}\otimes
g)B(gx_{1}x_{2}\otimes 1_{H})(1_{A}\otimes g)
\end{equation*}

\begin{eqnarray*}
&&\text{the right side }B(x_{2}\otimes x_{1})(1_{A}\otimes x_{1}) \\
&=&B(x_{2}\otimes 1_{H})(1_{A}\otimes x_{1})(1_{A}\otimes
x_{1})-(1_{A}\otimes gx_{1})B(x_{2}\otimes 1_{H})(1_{A}\otimes
g)(1_{A}\otimes x_{1}) \\
&&-(1_{A}\otimes g)B(gx_{1}x_{2}\otimes 1_{H})(1_{A}\otimes g)(1_{A}\otimes
x_{1})
\end{eqnarray*}%
i.e.%
\begin{eqnarray*}
\text{ clean } &&\text{ right side }B(x_{2}\otimes x_{1})(1_{A}\otimes x_{1})
\\
&=&(1_{A}\otimes gx_{1})B(x_{2}\otimes 1_{H})(1_{A}\otimes
gx_{1})+(1_{A}\otimes g)B(gx_{1}x_{2}\otimes 1_{H})(1_{A}\otimes gx_{1})
\end{eqnarray*}%
\begin{gather*}
(1_{A}\otimes x_{1})B(gx_{2}\otimes gx_{1})=(1_{A}\otimes
x_{1})(1_{A}\otimes g)B(x_{2}\otimes x_{1})(1_{A}\otimes g) \\
\overset{\left( \ref{form x2otx1}\right) }{=}(1_{A}\otimes gx_{1})\left[
\begin{array}{c}
B(x_{2}\otimes 1_{H})(1_{A}\otimes x_{1})-(1_{A}\otimes
gx_{1})B(x_{2}\otimes 1_{H})(1_{A}\otimes g) \\
-(1_{A}\otimes g)B(gx_{1}x_{2}\otimes 1_{H})(1_{A}\otimes g)%
\end{array}%
\right] (1_{A}\otimes g) \\
\left[
\begin{array}{c}
(1_{A}\otimes gx_{1})B(x_{2}\otimes 1_{H})(1_{A}\otimes x_{1})(1_{A}\otimes
g) \\
-(1_{A}\otimes gx_{1})(1_{A}\otimes gx_{1})B(x_{2}\otimes
1_{H})(1_{A}\otimes g)(1_{A}\otimes g) \\
-(1_{A}\otimes gx_{1})(1_{A}\otimes g)B(gx_{1}x_{2}\otimes
1_{H})(1_{A}\otimes g)(1_{A}\otimes g)%
\end{array}%
\right] \\
=(1_{A}\otimes gx_{1})B(x_{2}\otimes 1_{H})(1_{A}\otimes
gx_{1})-(1_{A}\otimes x_{1})B(gx_{1}x_{2}\otimes 1_{H})
\end{gather*}%
\begin{eqnarray*}
&&B(x_{1}x_{2}\otimes gx_{1}) \\
&&\overset{\left( \ref{form x1x2otgx1}\right) }{=}(1_{A}\otimes
g)B(gx_{1}x_{2}\otimes 1_{H})(1_{A}\otimes gx_{1})+(1_{A}\otimes
x_{1})B(gx_{1}x_{2}\otimes 1_{H})+
\end{eqnarray*}%
\begin{eqnarray*}
&&\text{right side}(1_{A}\otimes x_{1})B(gx_{2}\otimes
gx_{1})+B(x_{1}x_{2}\otimes gx_{1}) \\
&=&(1_{A}\otimes gx_{1})B(x_{2}\otimes 1_{H})(1_{A}\otimes
gx_{1})-(1_{A}\otimes x_{1})B(gx_{1}x_{2}\otimes 1_{H})+ \\
&&+(1_{A}\otimes g)B(gx_{1}x_{2}\otimes 1_{H})(1_{A}\otimes
gx_{1})+(1_{A}\otimes x_{1})B(gx_{1}x_{2}\otimes 1_{H}) \\
&=&(1_{A}\otimes gx_{1})B(x_{2}\otimes 1_{H})(1_{A}\otimes
gx_{1})+(1_{A}\otimes g)B(gx_{1}x_{2}\otimes 1_{H})(1_{A}\otimes gx_{1})
\end{eqnarray*}%
and we conclude.

\subsubsection{$x_{2}\otimes x_{2}$}

\begin{equation*}
B(x_{2}\otimes x_{2})(1_{A}\otimes x_{1})\overset{?}{=}(1_{A}\otimes
x_{1})B(gx_{2}\otimes gx_{2})+B(x_{1}x_{2}\otimes gx_{2})+B(x_{2}\otimes
x_{1}x_{2})
\end{equation*}

\begin{equation*}
B(x_{2}\otimes x_{2})\overset{\left( \ref{form x2otx2}\right) }{=}%
B(x_{2}\otimes 1_{H})(1_{A}\otimes x_{2})-(1_{A}\otimes
gx_{2})B(x_{2}\otimes 1_{H})(1_{A}\otimes g).
\end{equation*}%
\begin{eqnarray*}
&&\text{left side }B(x_{2}\otimes x_{2})(1_{A}\otimes x_{1})\overset{\left( %
\ref{form x2otx2}\right) }{=} \\
&&B(x_{2}\otimes 1_{H})(1_{A}\otimes x_{2})(1_{A}\otimes
x_{1})-(1_{A}\otimes gx_{2})B(x_{2}\otimes 1_{H})(1_{A}\otimes
g)(1_{A}\otimes x_{1})
\end{eqnarray*}%
i.e.%
\begin{eqnarray*}
&&\text{left side }B(x_{2}\otimes x_{2})(1_{A}\otimes x_{1})\overset{\left( %
\ref{form x2otx2}\right) }{=} \\
&&B(x_{2}\otimes 1_{H})(1_{A}\otimes x_{1}x_{2})+(1_{A}\otimes
gx_{2})B(x_{2}\otimes 1_{H})(1_{A}\otimes gx_{1})
\end{eqnarray*}%
\begin{eqnarray*}
&&B(x_{1}x_{2}\otimes gx_{2}) \\
&&\overset{\left( \ref{form x1x2otgx2}\right) }{=}(1_{A}\otimes
g)B(gx_{1}x_{2}\otimes 1_{H})(1_{A}\otimes gx_{2})+(1_{A}\otimes
x_{2})B(gx_{1}x_{2}\otimes 1_{H})
\end{eqnarray*}%
\begin{eqnarray*}
&&B(x_{2}\otimes x_{1}x_{2})\overset{\left( \ref{form x2otx1x2}\right) }{=}%
B(x_{2}\otimes 1_{H})(1_{A}\otimes x_{1}x_{2}) \\
&&+(1_{A}\otimes gx_{2})B(x_{2}\otimes 1_{H})(1_{A}\otimes gx_{1}) \\
&&-(1_{A}\otimes gx_{1})B(x_{2}\otimes 1_{H})(1_{A}\otimes gx_{2}) \\
&&+(1_{A}\otimes x_{1}x_{2})B(x_{2}\otimes 1_{H}) \\
&&-(1_{A}\otimes g)B(gx_{1}x_{2}\otimes 1_{H})(1_{A}\otimes gx_{2}) \\
&&-(1_{A}\otimes x_{2})B(gx_{1}x_{2}\otimes 1_{H}).
\end{eqnarray*}%
\begin{eqnarray*}
(1_{A}\otimes x_{1})B(gx_{2}\otimes gx_{2}) &=&(1_{A}\otimes
x_{1})(1_{A}\otimes g)B(x_{2}\otimes x_{2})(1_{A}\otimes g) \\
&=&(1_{A}\otimes x_{1})(1_{A}\otimes g)\left[
\begin{array}{c}
B(x_{2}\otimes 1_{H})(1_{A}\otimes x_{2}) \\
-(1_{A}\otimes gx_{2})B(x_{2}\otimes 1_{H})(1_{A}\otimes g)%
\end{array}%
\right] (1_{A}\otimes g) \\
&=&\left[
\begin{array}{c}
(1_{A}\otimes gx_{1})B(x_{2}\otimes 1_{H})(1_{A}\otimes x_{2})(1_{A}\otimes
g) \\
-(1_{A}\otimes gx_{1})(1_{A}\otimes gx_{2})B(x_{2}\otimes
1_{H})(1_{A}\otimes g)(1_{A}\otimes g)%
\end{array}%
\right] \\
&=&(1_{A}\otimes gx_{1})B(x_{2}\otimes 1_{H})(1_{A}\otimes
gx_{2})-(1_{A}\otimes x_{1}x_{2})B(x_{2}\otimes 1_{H})
\end{eqnarray*}%
\begin{eqnarray*}
&&\text{right side }(1_{A}\otimes x_{1})B(gx_{2}\otimes
gx_{2})+B(x_{1}x_{2}\otimes gx_{2})+B(x_{2}\otimes x_{1}x_{2}) \\
&=&(1_{A}\otimes gx_{1})B(x_{2}\otimes 1_{H})(1_{A}\otimes
gx_{2})-(1_{A}\otimes x_{1}x_{2})B(x_{2}\otimes 1_{H})+ \\
&&(1_{A}\otimes g)B(gx_{1}x_{2}\otimes 1_{H})(1_{A}\otimes
gx_{2})+(1_{A}\otimes x_{2})B(gx_{1}x_{2}\otimes 1_{H})+ \\
&&+B(x_{2}\otimes 1_{H})(1_{A}\otimes x_{1}x_{2})+ \\
&&+(1_{A}\otimes gx_{2})B(x_{2}\otimes 1_{H})(1_{A}\otimes gx_{1})+ \\
&&-(1_{A}\otimes gx_{1})B(x_{2}\otimes 1_{H})(1_{A}\otimes gx_{2})+ \\
&&+(1_{A}\otimes x_{1}x_{2})B(x_{2}\otimes 1_{H})+ \\
&&-(1_{A}\otimes g)B(gx_{1}x_{2}\otimes 1_{H})(1_{A}\otimes gx_{2})+ \\
&&-(1_{A}\otimes x_{2})B(gx_{1}x_{2}\otimes 1_{H})
\end{eqnarray*}%
\begin{eqnarray*}
\text{clean} &&\text{right side }(1_{A}\otimes x_{1})B(gx_{2}\otimes
gx_{2})+B(x_{1}x_{2}\otimes gx_{2})+B(x_{2}\otimes x_{1}x_{2}) \\
&=&B(x_{2}\otimes 1_{H})(1_{A}\otimes x_{1}x_{2})+(1_{A}\otimes
gx_{2})B(x_{2}\otimes 1_{H})(1_{A}\otimes gx_{1})
\end{eqnarray*}%
and we conclude.

\subsubsection{$x_{2}\otimes x_{1}x_{2}$}

\begin{equation*}
B(x_{2}\otimes x_{1}x_{2})(1_{A}\otimes x_{1})\overset{?}{=}(1_{A}\otimes
x_{1})B(gx_{2}\otimes gx_{1}x_{2})+B(x_{1}x_{2}\otimes
gx_{1}x_{2})+B(x_{2}\otimes x_{1}x_{1}x_{2})
\end{equation*}%
i.e.%
\begin{equation*}
B(x_{2}\otimes x_{1}x_{2})(1_{A}\otimes x_{1})\overset{?}{=}(1_{A}\otimes
gx_{1})B(x_{2}\otimes x_{1}x_{2})(1_{A}\otimes g)+B(x_{1}x_{2}\otimes
gx_{1}x_{2})
\end{equation*}%
\begin{eqnarray*}
&&B(x_{2}\otimes x_{1}x_{2})\overset{\left( \ref{form x2otx1x2}\right) }{=}%
B(x_{2}\otimes 1_{H})(1_{A}\otimes x_{1}x_{2}) \\
&&+(1_{A}\otimes gx_{2})B(x_{2}\otimes 1_{H})(1_{A}\otimes gx_{1}) \\
&&-(1_{A}\otimes gx_{1})B(x_{2}\otimes 1_{H})(1_{A}\otimes gx_{2}) \\
&&+(1_{A}\otimes x_{1}x_{2})B(x_{2}\otimes 1_{H}) \\
&&-(1_{A}\otimes g)B(gx_{1}x_{2}\otimes 1_{H})(1_{A}\otimes gx_{2}) \\
&&-(1_{A}\otimes x_{2})B(gx_{1}x_{2}\otimes 1_{H}).
\end{eqnarray*}%
\begin{eqnarray*}
\text{left side } &&B(x_{2}\otimes x_{1}x_{2})(1_{A}\otimes x_{1})\overset{%
\left( \ref{form x2otx1x2}\right) }{=}B(x_{2}\otimes 1_{H})(1_{A}\otimes
x_{1}x_{2})(1_{A}\otimes x_{1}) \\
&&+(1_{A}\otimes gx_{2})B(x_{2}\otimes 1_{H})(1_{A}\otimes
gx_{1})(1_{A}\otimes x_{1}) \\
&&-(1_{A}\otimes gx_{1})B(x_{2}\otimes 1_{H})(1_{A}\otimes
gx_{2})(1_{A}\otimes x_{1}) \\
&&+(1_{A}\otimes x_{1}x_{2})B(x_{2}\otimes 1_{H})(1_{A}\otimes x_{1}) \\
&&-(1_{A}\otimes g)B(gx_{1}x_{2}\otimes 1_{H})(1_{A}\otimes
gx_{2})(1_{A}\otimes x_{1}) \\
&&-(1_{A}\otimes x_{2})B(gx_{1}x_{2}\otimes 1_{H})(1_{A}\otimes x_{1})
\end{eqnarray*}%
i.e.%
\begin{eqnarray*}
\text{left side } &&B(x_{2}\otimes x_{1}x_{2})(1_{A}\otimes x_{1})\overset{%
\left( \ref{form x2otx1x2}\right) }{=} \\
&&+(1_{A}\otimes gx_{1})B(x_{2}\otimes 1_{H})(1_{A}\otimes gx_{1}x_{2}) \\
&&+(1_{A}\otimes x_{1}x_{2})B(x_{2}\otimes 1_{H})(1_{A}\otimes x_{1}) \\
&&+(1_{A}\otimes g)B(gx_{1}x_{2}\otimes 1_{H})(1_{A}\otimes gx_{1}x_{2}) \\
&&-(1_{A}\otimes x_{2})B(gx_{1}x_{2}\otimes 1_{H})(1_{A}\otimes x_{1})
\end{eqnarray*}%
\begin{eqnarray*}
&&(1_{A}\otimes gx_{1})B(x_{2}\otimes x_{1}x_{2})(1_{A}\otimes g)\overset{%
\left( \ref{form x2otx1x2}\right) }{=}(1_{A}\otimes gx_{1})B(x_{2}\otimes
1_{H})(1_{A}\otimes x_{1}x_{2})(1_{A}\otimes g) \\
&&+(1_{A}\otimes gx_{1})(1_{A}\otimes gx_{2})B(x_{2}\otimes
1_{H})(1_{A}\otimes gx_{1})(1_{A}\otimes g) \\
&&-(1_{A}\otimes gx_{1})(1_{A}\otimes gx_{1})B(x_{2}\otimes
1_{H})(1_{A}\otimes gx_{2})(1_{A}\otimes g) \\
&&+(1_{A}\otimes gx_{1})(1_{A}\otimes x_{1}x_{2})B(x_{2}\otimes
1_{H})(1_{A}\otimes g) \\
&&-(1_{A}\otimes gx_{1})(1_{A}\otimes g)B(gx_{1}x_{2}\otimes
1_{H})(1_{A}\otimes gx_{2})(1_{A}\otimes g) \\
&&-(1_{A}\otimes gx_{1})(1_{A}\otimes x_{2})B(gx_{1}x_{2}\otimes
1_{H})(1_{A}\otimes g)
\end{eqnarray*}%
i.e.%
\begin{eqnarray*}
&&(1_{A}\otimes gx_{1})B(x_{2}\otimes x_{1}x_{2})(1_{A}\otimes g)\overset{%
\left( \ref{form x2otx1x2}\right) }{=}(1_{A}\otimes gx_{1})B(x_{2}\otimes
1_{H})(1_{A}\otimes gx_{1}x_{2}) \\
&&+(1_{A}\otimes x_{1}x_{2})B(x_{2}\otimes 1_{H})(1_{A}\otimes x_{1}) \\
&&-(1_{A}\otimes x_{1})B(gx_{1}x_{2}\otimes 1_{H})(1_{A}\otimes x_{2}) \\
&&-(1_{A}\otimes gx_{1}x_{2})B(gx_{1}x_{2}\otimes 1_{H})(1_{A}\otimes g)
\end{eqnarray*}%
\begin{eqnarray*}
&&B(x_{1}x_{2}\otimes gx_{1}x_{2})\overset{\left( \ref{form x1x2otgx1gx2}%
\right) }{=}(1_{A}\otimes g)B(gx_{1}x_{2}\otimes 1_{H})(1_{A}\otimes
gx_{1}x_{2}) \\
&&-(1_{A}\otimes x_{2})B(gx_{1}x_{2}\otimes 1_{H})(1_{A}\otimes x_{1}) \\
&&+(1_{A}\otimes x_{1})B(gx_{1}x_{2}\otimes 1_{H})(1_{A}\otimes x_{2}) \\
&&+(1_{A}\otimes gx_{1}x_{2})B(gx_{1}x_{2}\otimes 1_{H})(1_{A}\otimes g)
\end{eqnarray*}%
\begin{eqnarray*}
&&\text{right side }(1_{A}\otimes gx_{1})B(x_{2}\otimes
x_{1}x_{2})(1_{A}\otimes g)++B(x_{1}x_{2}\otimes gx_{1}x_{2})= \\
&&(1_{A}\otimes gx_{1})B(x_{2}\otimes 1_{H})(1_{A}\otimes gx_{1}x_{2})+ \\
&&+(1_{A}\otimes x_{1}x_{2})B(x_{2}\otimes 1_{H})(1_{A}\otimes x_{1})+ \\
&&-(1_{A}\otimes x_{1})B(gx_{1}x_{2}\otimes 1_{H})(1_{A}\otimes x_{2})+ \\
&&-(1_{A}\otimes gx_{1}x_{2})B(gx_{1}x_{2}\otimes 1_{H})(1_{A}\otimes g)+ \\
&&(1_{A}\otimes g)B(gx_{1}x_{2}\otimes 1_{H})(1_{A}\otimes gx_{1}x_{2})+ \\
&&-(1_{A}\otimes x_{2})B(gx_{1}x_{2}\otimes 1_{H})(1_{A}\otimes x_{1})+ \\
&&+(1_{A}\otimes x_{1})B(gx_{1}x_{2}\otimes 1_{H})(1_{A}\otimes x_{2})+ \\
&&+(1_{A}\otimes gx_{1}x_{2})B(gx_{1}x_{2}\otimes 1_{H})(1_{A}\otimes g)+
\end{eqnarray*}%
\begin{eqnarray*}
\text{clean } &&\text{right side }(1_{A}\otimes gx_{1})B(x_{2}\otimes
x_{1}x_{2})(1_{A}\otimes g)+B(x_{1}x_{2}\otimes gx_{1}x_{2})= \\
&&(1_{A}\otimes gx_{1})B(x_{2}\otimes 1_{H})(1_{A}\otimes gx_{1}x_{2})+ \\
&&+(1_{A}\otimes x_{1}x_{2})B(x_{2}\otimes 1_{H})(1_{A}\otimes x_{1})+ \\
&&(1_{A}\otimes g)B(gx_{1}x_{2}\otimes 1_{H})(1_{A}\otimes gx_{1}x_{2})+ \\
&&-(1_{A}\otimes x_{2})B(gx_{1}x_{2}\otimes 1_{H})(1_{A}\otimes x_{1})
\end{eqnarray*}%
and we conclude.

\subsubsection{$x_{2}\otimes gx_{1}$}

\begin{equation*}
B(x_{2}\otimes gx_{1})(1_{A}\otimes x_{1})\overset{?}{=}(1_{A}\otimes
x_{1})B(gx_{2}\otimes ggx_{1})+B(x_{1}x_{2}\otimes ggx_{1})+B(x_{2}\otimes
x_{1}gx_{1})
\end{equation*}%
i.e.%
\begin{equation*}
B(x_{2}\otimes gx_{1})(1_{A}\otimes x_{1})\overset{?}{=}(1_{A}\otimes
x_{1})(1_{A}\otimes g)B(x_{2}\otimes gx_{1})(1_{A}\otimes
g)+B(x_{1}x_{2}\otimes x_{1})
\end{equation*}

\begin{eqnarray*}
&&B(x_{2}\otimes gx_{1}) \\
&&\overset{\left( \ref{form x2otgx1}\right) }{=}(1_{A}\otimes
g)B(gx_{2}\otimes 1_{H})(1_{A}\otimes gx_{1})+(1_{A}\otimes
x_{1})B(gx_{2}\otimes 1_{H})+B(x_{1}x_{2}\otimes 1_{H}).
\end{eqnarray*}%
\begin{eqnarray*}
&&\text{left side }B(x_{2}\otimes gx_{1})(1_{A}\otimes x_{1}) \\
&=&(1_{A}\otimes g)B(gx_{2}\otimes 1_{H})(1_{A}\otimes gx_{1})(1_{A}\otimes
x_{1})+(1_{A}\otimes x_{1})B(gx_{2}\otimes 1_{H})(1_{A}\otimes x_{1}) \\
&&+B(x_{1}x_{2}\otimes 1_{H})(1_{A}\otimes x_{1})
\end{eqnarray*}%
thus we get%
\begin{eqnarray*}
&&\text{left side }B(x_{2}\otimes gx_{1})(1_{A}\otimes x_{1}) \\
&=&(1_{A}\otimes x_{1})B(gx_{2}\otimes 1_{H})(1_{A}\otimes
x_{1})+B(x_{1}x_{2}\otimes 1_{H})(1_{A}\otimes x_{1})
\end{eqnarray*}%
Since
\begin{equation*}
B(x_{1}x_{2}\otimes x_{1})\overset{\left( \ref{form x1x2otx1}\right) }{=}%
B(x_{1}x_{2}\otimes 1_{H})(1_{A}\otimes x_{1})-(1_{A}\otimes
gx_{1})B(x_{1}x_{2}\otimes 1_{H})(1_{A}\otimes g)
\end{equation*}%
we get%
\begin{eqnarray*}
\text{right side} &&(1_{A}\otimes x_{1})(1_{A}\otimes g)B(x_{2}\otimes
gx_{1})(1_{A}\otimes g)+B(x_{1}x_{2}\otimes x_{1}) \\
&=&(1_{A}\otimes x_{1})(1_{A}\otimes g)(1_{A}\otimes g)B(gx_{2}\otimes
1_{H})(1_{A}\otimes gx_{1})(1_{A}\otimes g) \\
&&+(1_{A}\otimes x_{1})(1_{A}\otimes g)(1_{A}\otimes x_{1})B(gx_{2}\otimes
1_{H})(1_{A}\otimes g) \\
&&+(1_{A}\otimes x_{1})(1_{A}\otimes g)B(x_{1}x_{2}\otimes
1_{H})(1_{A}\otimes g)+ \\
&&B(x_{1}x_{2}\otimes 1_{H})(1_{A}\otimes x_{1})-(1_{A}\otimes
gx_{1})B(x_{1}x_{2}\otimes 1_{H})(1_{A}\otimes g)
\end{eqnarray*}%
i.e.%
\begin{eqnarray*}
\text{right side} &&(1_{A}\otimes x_{1})(1_{A}\otimes g)B(x_{2}\otimes
gx_{1})(1_{A}\otimes g)+B(x_{1}x_{2}\otimes x_{1}) \\
&=&(1_{A}\otimes x_{1})B(gx_{2}\otimes 1_{H})(1_{A}\otimes
x_{1})+B(x_{1}x_{2}\otimes 1_{H})(1_{A}\otimes x_{1})
\end{eqnarray*}

and we conclude.

\subsubsection{$x_{2}\otimes gx_{2}$}

\begin{equation*}
B(x_{2}\otimes gx_{2})(1_{A}\otimes x_{1})\overset{?}{=}(1_{A}\otimes
x_{1})B(gx_{2}\otimes ggx_{2})+B(x_{1}x_{2}\otimes ggx_{2})+B(x_{2}\otimes
x_{1}gx_{2})
\end{equation*}%
i.e.%
\begin{equation*}
B(x_{2}\otimes gx_{2})(1_{A}\otimes x_{1})\overset{?}{=}(1_{A}\otimes
x_{1})(1_{A}\otimes g)B(x_{2}\otimes gx_{2})(1_{A}\otimes
g)+B(x_{1}x_{2}\otimes x_{2})-B(x_{2}\otimes gx_{1}x_{2})
\end{equation*}%
\begin{equation*}
B(x_{2}\otimes gx_{2})\overset{\left( \ref{form x2otgx2}\right) }{=}%
(1_{A}\otimes g)B(gx_{2}\otimes 1_{H})(1_{A}\otimes gx_{2})+(1_{A}\otimes
x_{2})B(gx_{2}\otimes 1_{H})
\end{equation*}%
\begin{eqnarray*}
\text{left side }B(x_{2}\otimes gx_{2})(1_{A}\otimes x_{1}) &=&(1_{A}\otimes
g)B(gx_{2}\otimes 1_{H})(1_{A}\otimes gx_{2})(1_{A}\otimes x_{1}) \\
&&+(1_{A}\otimes x_{2})B(gx_{2}\otimes 1_{H})\text{ }(1_{A}\otimes x_{1})
\end{eqnarray*}%
i.e.%
\begin{eqnarray*}
\text{left side }B(x_{2}\otimes gx_{2})(1_{A}\otimes x_{1})
&=&-(1_{A}\otimes g)B(gx_{2}\otimes 1_{H})(1_{A}\otimes gx_{1}x_{2}) \\
&&+(1_{A}\otimes x_{2})B(gx_{2}\otimes 1_{H})\text{ }(1_{A}\otimes x_{1})
\end{eqnarray*}%
\begin{eqnarray*}
(1_{A}\otimes x_{1})(1_{A}\otimes g)B(x_{2}\otimes gx_{2})(1_{A}\otimes g)
&=&(1_{A}\otimes x_{1})(1_{A}\otimes g)(1_{A}\otimes g)B(gx_{2}\otimes
1_{H})(1_{A}\otimes gx_{2})(1_{A}\otimes g) \\
&&+(1_{A}\otimes x_{1})(1_{A}\otimes g)(1_{A}\otimes x_{2})B(gx_{2}\otimes
1_{H})(1_{A}\otimes g)
\end{eqnarray*}%
i.e.%
\begin{eqnarray*}
(1_{A}\otimes x_{1})(1_{A}\otimes g)B(x_{2}\otimes gx_{2})(1_{A}\otimes g)
&=&(1_{A}\otimes x_{1})B(gx_{2}\otimes 1_{H})(1_{A}\otimes x_{2}) \\
&&+(1_{A}\otimes gx_{1}x_{2})B(gx_{2}\otimes 1_{H})(1_{A}\otimes g)
\end{eqnarray*}%
Since%
\begin{equation*}
B(x_{1}x_{2}\otimes x_{2})\overset{\left( \ref{x1x2otx2}\right) }{=}%
B(x_{1}x_{2}\otimes 1_{H})(1_{A}\otimes x_{2})-(1_{A}\otimes
gx_{2})B(x_{1}x_{2}\otimes 1_{H})(1_{A}\otimes g)
\end{equation*}%
we get%
\begin{eqnarray*}
&&(1_{A}\otimes x_{1})(1_{A}\otimes g)B(x_{2}\otimes gx_{2})(1_{A}\otimes
g)+B(x_{1}x_{2}\otimes x_{2}) \\
&=&(1_{A}\otimes x_{1})B(gx_{2}\otimes 1_{H})(1_{A}\otimes x_{2})+ \\
&&+(1_{A}\otimes gx_{1}x_{2})B(gx_{2}\otimes 1_{H})(1_{A}\otimes g)+ \\
&&+B(x_{1}x_{2}\otimes 1_{H})(1_{A}\otimes x_{2})+ \\
&&-(1_{A}\otimes gx_{2})B(x_{1}x_{2}\otimes 1_{H})(1_{A}\otimes g)
\end{eqnarray*}%
Since%
\begin{eqnarray*}
&&-B(x_{2}\otimes gx_{1}x_{2})\overset{\left( \ref{form x2otgx1x2}\right) }{=%
}-(1_{A}\otimes g)B(gx_{2}\otimes 1_{H})(1_{A}\otimes gx_{1}x_{2}) \\
&&+(1_{A}\otimes x_{2})B(gx_{2}\otimes 1_{H})(1_{A}\otimes x_{1}) \\
&&-(1_{A}\otimes x_{1})B(gx_{2}\otimes 1_{H})(1_{A}\otimes x_{2}) \\
&&-(1_{A}\otimes gx_{1}x_{2})B(gx_{2}\otimes 1_{H})(1_{A}\otimes g) \\
&&-B(x_{1}x_{2}\otimes 1_{H})(1_{A}\otimes x_{2}) \\
&&+(1_{A}\otimes gx_{2})B(x_{1}x_{2}\otimes 1_{H})(1_{A}\otimes g)
\end{eqnarray*}%
we get that%
\begin{eqnarray*}
&&\text{right side }(1_{A}\otimes x_{1})(1_{A}\otimes g)B(x_{2}\otimes
gx_{2})(1_{A}\otimes g)+B(x_{1}x_{2}\otimes x_{2})-B(x_{2}\otimes
gx_{1}x_{2}) \\
&=&(1_{A}\otimes x_{1})B(gx_{2}\otimes 1_{H})(1_{A}\otimes x_{2})+ \\
&&+(1_{A}\otimes gx_{1}x_{2})B(gx_{2}\otimes 1_{H})(1_{A}\otimes g)+ \\
&&+B(x_{1}x_{2}\otimes 1_{H})(1_{A}\otimes x_{2})+ \\
&&-(1_{A}\otimes gx_{2})B(x_{1}x_{2}\otimes 1_{H})(1_{A}\otimes g) \\
&&-(1_{A}\otimes g)B(gx_{2}\otimes 1_{H})(1_{A}\otimes gx_{1}x_{2})+ \\
&&+(1_{A}\otimes x_{2})B(gx_{2}\otimes 1_{H})(1_{A}\otimes x_{1})+ \\
&&-(1_{A}\otimes x_{1})B(gx_{2}\otimes 1_{H})(1_{A}\otimes x_{2}) \\
&&-(1_{A}\otimes gx_{1}x_{2})B(gx_{2}\otimes 1_{H})(1_{A}\otimes g) \\
&&-B(x_{1}x_{2}\otimes 1_{H})(1_{A}\otimes x_{2})+ \\
&&+(1_{A}\otimes gx_{2})B(x_{1}x_{2}\otimes 1_{H})(1_{A}\otimes g)
\end{eqnarray*}%
hence we obtain%
\begin{eqnarray*}
&&\text{clean right side }(1_{A}\otimes x_{1})(1_{A}\otimes g)B(x_{2}\otimes
gx_{2})(1_{A}\otimes g)+B(x_{1}x_{2}\otimes x_{2})-B(x_{2}\otimes
gx_{1}x_{2}) \\
&=&-(1_{A}\otimes g)B(gx_{2}\otimes 1_{H})(1_{A}\otimes
gx_{1}x_{2})+(1_{A}\otimes x_{2})B(gx_{2}\otimes 1_{H})(1_{A}\otimes x_{1})
\end{eqnarray*}%
and we conclude.

\subsubsection{$x_{2}\otimes gx_{1}x_{2}$}

\begin{equation*}
B(x_{2}\otimes gx_{1}x_{2})(1_{A}\otimes x_{1})\overset{?}{=}(1_{A}\otimes
x_{1})B(gx_{2}\otimes ggx_{1}x_{2})+B(x_{1}x_{2}\otimes
ggx_{1}x_{2})+B(x_{2}\otimes x_{1}gx_{1}x_{2})
\end{equation*}%
i.e.%
\begin{equation*}
B(x_{2}\otimes gx_{1}x_{2})(1_{A}\otimes x_{1})\overset{?}{=}(1_{A}\otimes
x_{1})(1_{A}\otimes g)B(x_{2}\otimes gx_{1}x_{2})(1_{A}\otimes
g)+B(x_{1}x_{2}\otimes x_{1}x_{2})
\end{equation*}%
Since%
\begin{eqnarray*}
&&B(x_{2}\otimes gx_{1}x_{2})\overset{\left( \ref{form x2otgx1x2}\right) }{=}%
(1_{A}\otimes g)B(gx_{2}\otimes 1_{H})(1_{A}\otimes gx_{1}x_{2}) \\
&&-(1_{A}\otimes x_{2})B(gx_{2}\otimes 1_{H})(1_{A}\otimes x_{1}) \\
&&+(1_{A}\otimes x_{1})B(gx_{2}\otimes 1_{H})(1_{A}\otimes x_{2}) \\
&&+(1_{A}\otimes gx_{1}x_{2})B(gx_{2}\otimes 1_{H})(1_{A}\otimes g) \\
&&+B(x_{1}x_{2}\otimes 1_{H})(1_{A}\otimes x_{2}) \\
&&-(1_{A}\otimes gx_{2})B(x_{1}x_{2}\otimes 1_{H})(1_{A}\otimes g)
\end{eqnarray*}%
we get%
\begin{eqnarray*}
&&\text{left side }B(x_{2}\otimes gx_{1}x_{2})(1_{A}\otimes
x_{1})=(1_{A}\otimes g)B(gx_{2}\otimes 1_{H})(1_{A}\otimes
gx_{1}x_{2})(1_{A}\otimes x_{1}) \\
&&-(1_{A}\otimes x_{2})B(gx_{2}\otimes 1_{H})(1_{A}\otimes
x_{1})(1_{A}\otimes x_{1}) \\
&&+(1_{A}\otimes x_{1})B(gx_{2}\otimes 1_{H})(1_{A}\otimes
x_{2})(1_{A}\otimes x_{1}) \\
&&+(1_{A}\otimes gx_{1}x_{2})B(gx_{2}\otimes 1_{H})(1_{A}\otimes
g)(1_{A}\otimes x_{1}) \\
&&+B(x_{1}x_{2}\otimes 1_{H})(1_{A}\otimes x_{2})(1_{A}\otimes x_{1}) \\
&&-(1_{A}\otimes gx_{2})B(x_{1}x_{2}\otimes 1_{H})(1_{A}\otimes
g)(1_{A}\otimes x_{1})
\end{eqnarray*}%
i.e.%
\begin{eqnarray*}
&&\text{left side }B(x_{2}\otimes gx_{1}x_{2})(1_{A}\otimes x_{1})= \\
&&+(1_{A}\otimes x_{1})B(gx_{2}\otimes 1_{H})(1_{A}\otimes x_{1}x_{2}) \\
&&-(1_{A}\otimes gx_{1}x_{2})B(gx_{2}\otimes 1_{H})(1_{A}\otimes gx_{1}) \\
&&+B(x_{1}x_{2}\otimes 1_{H})(1_{A}\otimes x_{1}x_{2}) \\
&&+(1_{A}\otimes gx_{2})B(x_{1}x_{2}\otimes 1_{H})(1_{A}\otimes gx_{1})
\end{eqnarray*}%
We have%
\begin{eqnarray*}
&&(1_{A}\otimes x_{1})(1_{A}\otimes g)B(x_{2}\otimes
gx_{1}x_{2})(1_{A}\otimes g)\overset{\left( \ref{form x2otgx1x2}\right) }{=}
\\
&&(1_{A}\otimes x_{1})(1_{A}\otimes g)(1_{A}\otimes g)B(gx_{2}\otimes
1_{H})(1_{A}\otimes gx_{1}x_{2})(1_{A}\otimes g) \\
&&-(1_{A}\otimes x_{1})(1_{A}\otimes g)(1_{A}\otimes x_{2})B(gx_{2}\otimes
1_{H})(1_{A}\otimes x_{1})(1_{A}\otimes g) \\
&&+(1_{A}\otimes x_{1})(1_{A}\otimes g)(1_{A}\otimes x_{1})B(gx_{2}\otimes
1_{H})(1_{A}\otimes x_{2})(1_{A}\otimes g) \\
&&+(1_{A}\otimes x_{1})(1_{A}\otimes g)(1_{A}\otimes
gx_{1}x_{2})B(gx_{2}\otimes 1_{H})(1_{A}\otimes g)(1_{A}\otimes g) \\
&&+(1_{A}\otimes x_{1})(1_{A}\otimes g)B(x_{1}x_{2}\otimes
1_{H})(1_{A}\otimes x_{2})(1_{A}\otimes g) \\
&&-(1_{A}\otimes x_{1})(1_{A}\otimes g)(1_{A}\otimes
gx_{2})B(x_{1}x_{2}\otimes 1_{H})(1_{A}\otimes g)(1_{A}\otimes g)
\end{eqnarray*}%
i.e.%
\begin{eqnarray*}
&&(1_{A}\otimes x_{1})(1_{A}\otimes g)B(x_{2}\otimes
gx_{1}x_{2})(1_{A}\otimes g)\overset{\left( \ref{form x2otgx1x2}\right) }{=}
\\
&&(1_{A}\otimes x_{1})B(gx_{2}\otimes 1_{H})(1\otimes x_{1}x_{2}) \\
&&-(1_{A}\otimes gx_{1}x_{2})(1_{A}\otimes )B(gx_{2}\otimes
1_{H})(1_{A}\otimes gx_{1}) \\
&&+(1_{A}\otimes gx_{1})B(x_{1}x_{2}\otimes 1_{H})(1_{A}\otimes gx_{2}) \\
&&-(1_{A}\otimes x_{1}x_{2})B(x_{1}x_{2}\otimes 1_{H})
\end{eqnarray*}%
Since%
\begin{eqnarray*}
&&B(x_{1}x_{2}\otimes x_{1}x_{2})\overset{\left( \ref{form x1x2otx1x2}%
\right) }{=} \\
&&B(x_{1}x_{2}\otimes 1_{H})(1_{A}\otimes x_{1}x_{2})+ \\
&&+(1_{A}\otimes gx_{2})B(x_{1}x_{2}\otimes 1_{H})(1_{A}\otimes gx_{1}) \\
&&-(1_{A}\otimes gx_{1})B(x_{1}x_{2}\otimes 1_{H})(1_{A}\otimes gx_{2})+ \\
&&+(1_{A}\otimes x_{1}x_{2})B(x_{1}x_{2}\otimes 1_{H})+
\end{eqnarray*}%
we get%
\begin{eqnarray*}
\text{right side} &&(1_{A}\otimes x_{1})(1_{A}\otimes g)B(x_{2}\otimes
gx_{1}x_{2})(1_{A}\otimes g)+B(x_{1}x_{2}\otimes x_{1}x_{2})= \\
&&+(1_{A}\otimes x_{1})B(gx_{2}\otimes 1_{H})(1\otimes x_{1}x_{2})+ \\
&&-(1_{A}\otimes gx_{1}x_{2})B(gx_{2}\otimes 1_{H})(1_{A}\otimes gx_{1})+ \\
&&+(1_{A}\otimes gx_{1})B(x_{1}x_{2}\otimes 1_{H})(1_{A}\otimes gx_{2})+ \\
&&-(1_{A}\otimes x_{1}x_{2})B(x_{1}x_{2}\otimes 1_{H})+ \\
&&+B(x_{1}x_{2}\otimes 1_{H})(1_{A}\otimes x_{1}x_{2})+ \\
&&+(1_{A}\otimes gx_{2})B(x_{1}x_{2}\otimes 1_{H})(1_{A}\otimes gx_{1})+ \\
&&-(1_{A}\otimes gx_{1})B(x_{1}x_{2}\otimes 1_{H})(1_{A}\otimes gx_{2})+ \\
&&+(1_{A}\otimes x_{1}x_{2})B(x_{1}x_{2}\otimes 1_{H})+
\end{eqnarray*}%
\begin{eqnarray*}
\text{clean right side} &&(1_{A}\otimes x_{1})(1_{A}\otimes g)B(x_{2}\otimes
gx_{1}x_{2})(1_{A}\otimes g)+B(x_{1}x_{2}\otimes x_{1}x_{2})= \\
&&+(1_{A}\otimes x_{1})B(gx_{2}\otimes 1_{H})(1\otimes x_{1}x_{2})+ \\
&&-(1_{A}\otimes gx_{1}x_{2})B(gx_{2}\otimes 1_{H})(1_{A}\otimes gx_{1})+ \\
&&+B(x_{1}x_{2}\otimes 1_{H})(1_{A}\otimes x_{1}x_{2})+ \\
&&+(1_{A}\otimes gx_{2})B(x_{1}x_{2}\otimes 1_{H})(1_{A}\otimes gx_{1})+
\end{eqnarray*}%
and we conclude.

\subsubsection{$x_{1}x_{2}\otimes x_{1}$}

\begin{equation*}
B(x_{1}x_{2}\otimes x_{1})(1_{A}\otimes x_{1})\overset{?}{=}(1_{A}\otimes
x_{1})B(gx_{1}x_{2}\otimes gx_{1})+B(x_{1}x_{1}x_{2}\otimes
gx_{1})+B(x_{1}x_{2}\otimes x_{1}x_{1})
\end{equation*}%
i.e.%
\begin{equation*}
B(x_{1}x_{2}\otimes x_{1})(1_{A}\otimes x_{1})\overset{?}{=}(1_{A}\otimes
gx_{1})B(x_{1}x_{2}\otimes x_{1})(1_{A}\otimes g)
\end{equation*}%
Since%
\begin{equation*}
B(x_{1}x_{2}\otimes x_{1})\overset{\left( \ref{form x1x2otx1}\right) }{=}%
B(x_{1}x_{2}\otimes 1_{H})(1_{A}\otimes x_{1})-(1_{A}\otimes
gx_{1})B(x_{1}x_{2}\otimes 1_{H})(1_{A}\otimes g)
\end{equation*}%
we get%
\begin{eqnarray*}
\text{left side }B(x_{1}x_{2}\otimes x_{1})(1_{A}\otimes x_{1})
&=&B(x_{1}x_{2}\otimes 1_{H})(1_{A}\otimes x_{1})(1_{A}\otimes x_{1}) \\
&&-(1_{A}\otimes gx_{1})B(x_{1}x_{2}\otimes 1_{H})(1_{A}\otimes
g)(1_{A}\otimes x_{1})\text{ }
\end{eqnarray*}%
i.e.%
\begin{equation*}
\text{left side }B(x_{1}x_{2}\otimes x_{1})(1_{A}\otimes
x_{1})=(1_{A}\otimes gx_{1})B(x_{1}x_{2}\otimes 1_{H})(1_{A}\otimes gx_{1})%
\text{ }
\end{equation*}%
and%
\begin{eqnarray*}
&&\text{right side }B(x_{1}x_{2}\otimes x_{1}) \\
&=&(1_{A}\otimes gx_{1})B(x_{1}x_{2}\otimes 1_{H})(1_{A}\otimes
x_{1})(1_{A}\otimes g) \\
&&-(1_{A}\otimes gx_{1})(1_{A}\otimes gx_{1})B(x_{1}x_{2}\otimes
1_{H})(1_{A}\otimes g)(1_{A}\otimes g)
\end{eqnarray*}%
i.e.%
\begin{equation*}
\text{right side }B(x_{1}x_{2}\otimes x_{1})=(1_{A}\otimes
gx_{1})B(x_{1}x_{2}\otimes 1_{H})(1_{A}\otimes gx_{1})
\end{equation*}%
and we conclude.

\subsubsection{$x_{1}x_{2}\otimes x_{2}$}

\begin{equation*}
B(x_{1}x_{2}\otimes x_{2})(1_{A}\otimes x_{1})\overset{?}{=}(1_{A}\otimes
x_{1})B(gx_{1}x_{2}\otimes gx_{2})+B(x_{1}x_{1}x_{2}\otimes
gx_{2})+B(x_{1}x_{2}\otimes x_{1}x_{2})
\end{equation*}%
i.e.%
\begin{equation*}
B(x_{1}x_{2}\otimes x_{2})(1_{A}\otimes x_{1})\overset{?}{=}(1_{A}\otimes
x_{1})(1_{A}\otimes g)B(x_{1}x_{2}\otimes x_{2})(1_{A}\otimes
g)+B(x_{1}x_{2}\otimes x_{1}x_{2})
\end{equation*}%
Since%
\begin{equation*}
B(x_{1}x_{2}\otimes x_{2})\overset{\left( \ref{x1x2otx2}\right) }{=}%
B(x_{1}x_{2}\otimes 1_{H})(1_{A}\otimes x_{2})-(1_{A}\otimes
gx_{2})B(x_{1}x_{2}\otimes 1_{H})(1_{A}\otimes g)
\end{equation*}%
we get%
\begin{eqnarray*}
&&\text{right side }B(x_{1}x_{2}\otimes x_{2}) \\
&=&B(x_{1}x_{2}\otimes 1_{H})(1_{A}\otimes x_{2})(1_{A}\otimes x_{1}) \\
&&-(1_{A}\otimes gx_{2})B(x_{1}x_{2}\otimes 1_{H})(1_{A}\otimes
g)(1_{A}\otimes x_{1})
\end{eqnarray*}%
i.e.%
\begin{eqnarray*}
&&\text{right side }B(x_{1}x_{2}\otimes x_{2}) \\
&=&B(x_{1}x_{2}\otimes 1_{H})(1_{A}\otimes x_{1}x_{2}) \\
&&+(1_{A}\otimes gx_{2})B(x_{1}x_{2}\otimes 1_{H})(1_{A}\otimes gx_{1})
\end{eqnarray*}%
We have%
\begin{eqnarray*}
&&(1_{A}\otimes x_{1})(1_{A}\otimes g)B(x_{1}x_{2}\otimes x_{2}) \\
&=&(1_{A}\otimes x_{1})(1_{A}\otimes g)B(x_{1}x_{2}\otimes
1_{H})(1_{A}\otimes x_{2})(1_{A}\otimes g) \\
&&-(1_{A}\otimes x_{1})(1_{A}\otimes g)(1_{A}\otimes
gx_{2})B(x_{1}x_{2}\otimes 1_{H})(1_{A}\otimes g)(1_{A}\otimes g)
\end{eqnarray*}%
i.e.%
\begin{eqnarray*}
&&(1_{A}\otimes x_{1})(1_{A}\otimes g)B(x_{1}x_{2}\otimes
x_{2})(1_{A}\otimes g) \\
&=&(1_{A}\otimes gx_{1})B(x_{1}x_{2}\otimes 1_{H})(1_{A}\otimes gx_{2}) \\
&&-(1_{A}\otimes x_{1}x_{2})B(x_{1}x_{2}\otimes 1_{H})
\end{eqnarray*}%
Since%
\begin{eqnarray*}
&&B(x_{1}x_{2}\otimes x_{1}x_{2}) \\
&&\overset{\left( \ref{form x1x2otx1x2}\right) }{=}B(x_{1}x_{2}\otimes
1_{H})(1_{A}\otimes x_{1}x_{2}) \\
&&+(1_{A}\otimes gx_{2})B(x_{1}x_{2}\otimes 1_{H})(1_{A}\otimes gx_{1}) \\
&&-(1_{A}\otimes gx_{1})B(x_{1}x_{2}\otimes 1_{H})(1_{A}\otimes gx_{2})+ \\
&&+(1_{A}\otimes x_{1}x_{2})B(x_{1}x_{2}\otimes 1_{H})+
\end{eqnarray*}%
we get%
\begin{eqnarray*}
&&\text{right side }(1_{A}\otimes x_{1})(1_{A}\otimes g)B(x_{1}x_{2}\otimes
x_{2})+B(x_{1}x_{2}\otimes x_{1}x_{2}) \\
&=&(1_{A}\otimes gx_{1})B(x_{1}x_{2}\otimes 1_{H})(1_{A}\otimes gx_{2}) \\
&&-(1_{A}\otimes x_{1}x_{2})B(x_{1}x_{2}\otimes 1_{H})+ \\
&&+B(x_{1}x_{2}\otimes 1_{H})(1_{A}\otimes x_{1}x_{2})+ \\
&&+(1_{A}\otimes gx_{2})B(x_{1}x_{2}\otimes 1_{H})(1_{A}\otimes gx_{1})+ \\
&&-(1_{A}\otimes gx_{1})B(x_{1}x_{2}\otimes 1_{H})(1_{A}\otimes gx_{2})+ \\
&&+(1_{A}\otimes x_{1}x_{2})B(x_{1}x_{2}\otimes 1_{H})
\end{eqnarray*}%
so that%
\begin{eqnarray*}
&&\text{clean right side }(1_{A}\otimes x_{1})(1_{A}\otimes
g)B(x_{1}x_{2}\otimes x_{2})+B(x_{1}x_{2}\otimes x_{1}x_{2}) \\
&=&+B(x_{1}x_{2}\otimes 1_{H})(1_{A}\otimes x_{1}x_{2})+ \\
&&+(1_{A}\otimes gx_{2})B(x_{1}x_{2}\otimes 1_{H})(1_{A}\otimes gx_{1})
\end{eqnarray*}%
and we conclude.

\subsubsection{$x_{1}x_{2}\otimes x_{1}x_{2}$}

\begin{equation*}
B(x_{1}x_{2}\otimes x_{1}x_{2})(1_{A}\otimes x_{1})\overset{?}{=}%
(1_{A}\otimes x_{1})B(gx_{1}x_{2}\otimes
gx_{1}x_{2})+B(x_{1}x_{1}x_{2}\otimes gx_{1}x_{2})+B(x_{1}x_{2}\otimes
x_{1}x_{1}x_{2})
\end{equation*}%
i.e.%
\begin{equation*}
B(x_{1}x_{2}\otimes x_{1}x_{2})(1_{A}\otimes x_{1})\overset{?}{=}%
(1_{A}\otimes gx_{1})B(x_{1}x_{2}\otimes x_{1}x_{2})(1_{A}\otimes g).
\end{equation*}%
Since%
\begin{eqnarray*}
&&B(x_{1}x_{2}\otimes x_{1}x_{2})\overset{\left( \ref{form x1x2otx1x2}%
\right) }{=} \\
&&B(x_{1}x_{2}\otimes 1_{H})(1_{A}\otimes x_{1}x_{2}) \\
&&+(1_{A}\otimes gx_{2})B(x_{1}x_{2}\otimes 1_{H})(1_{A}\otimes gx_{1}) \\
&&-(1_{A}\otimes gx_{1})B(x_{1}x_{2}\otimes 1_{H})(1_{A}\otimes gx_{2})+ \\
&&+(1_{A}\otimes x_{1}x_{2})B(x_{1}x_{2}\otimes 1_{H})+
\end{eqnarray*}%
we get%
\begin{eqnarray*}
&&\text{left side }B(x_{1}x_{2}\otimes x_{1}x_{2})(1_{A}\otimes x_{1})= \\
&&B(x_{1}x_{2}\otimes 1_{H})(1_{A}\otimes x_{1}x_{2})(1_{A}\otimes x_{1}) \\
&&+(1_{A}\otimes gx_{2})B(x_{1}x_{2}\otimes 1_{H})(1_{A}\otimes
gx_{1})(1_{A}\otimes x_{1}) \\
&&-(1_{A}\otimes gx_{1})B(x_{1}x_{2}\otimes 1_{H})(1_{A}\otimes
gx_{2})(1_{A}\otimes x_{1})+ \\
&&+(1_{A}\otimes x_{1}x_{2})B(x_{1}x_{2}\otimes 1_{H})(1_{A}\otimes x_{1})
\end{eqnarray*}%
i.e.%
\begin{eqnarray*}
&&\text{left side }B(x_{1}x_{2}\otimes x_{1}x_{2})(1_{A}\otimes x_{1})= \\
&&+(1_{A}\otimes gx_{1})B(x_{1}x_{2}\otimes 1_{H})(1_{A}\otimes gx_{1}x_{2})
\\
&&+(1_{A}\otimes x_{1}x_{2})B(x_{1}x_{2}\otimes 1_{H})(1_{A}\otimes x_{1})
\end{eqnarray*}%
and%
\begin{eqnarray*}
&&\text{right side }(1_{A}\otimes gx_{1})B(x_{1}x_{2}\otimes
x_{1}x_{2})(1_{A}\otimes g) \\
&&=(1_{A}\otimes gx_{1})B(x_{1}x_{2}\otimes 1_{H})(1_{A}\otimes
x_{1}x_{2})(1_{A}\otimes g) \\
&&+(1_{A}\otimes gx_{1})(1_{A}\otimes gx_{2})B(x_{1}x_{2}\otimes
1_{H})(1_{A}\otimes gx_{1})(1_{A}\otimes g) \\
&&-(1_{A}\otimes gx_{1})(1_{A}\otimes gx_{1})B(x_{1}x_{2}\otimes
1_{H})(1_{A}\otimes gx_{2})(1_{A}\otimes g)+ \\
&&+(1_{A}\otimes gx_{1})(1_{A}\otimes x_{1}x_{2})B(x_{1}x_{2}\otimes 1_{H})+
\end{eqnarray*}%
\qquad i.e.%
\begin{eqnarray*}
&&\text{right side }(1_{A}\otimes gx_{1})B(x_{1}x_{2}\otimes
x_{1}x_{2})(1_{A}\otimes g) \\
&=&(1_{A}\otimes gx_{1})B(x_{1}x_{2}\otimes 1_{H})(1_{A}\otimes gx_{1}x_{2})
\\
&&+(1_{A}\otimes x_{1}x_{2})B(x_{1}x_{2}\otimes 1_{H})(1_{A}\otimes x_{1})
\end{eqnarray*}%
and we conclude.

\subsubsection{$x_{1}x_{2}\otimes gx_{1}$}

\begin{eqnarray*}
&&B(x_{1}x_{2}\otimes gx_{1})(1_{A}\otimes x_{1}) \\
&&\overset{?}{=}(1_{A}\otimes x_{1})B(gx_{1}x_{2}\otimes
ggx_{1})+B(x_{1}x_{1}x_{2}\otimes ggx_{1}) \\
&&+B(x_{1}x_{2}\otimes x_{1}gx_{1})
\end{eqnarray*}%
i.e.%
\begin{equation*}
B(x_{1}x_{2}\otimes gx_{1})(1_{A}\otimes x_{1})\overset{?}{=}(1_{A}\otimes
gx_{1})B(x_{1}x_{2}\otimes gx_{1})(1_{A}\otimes g)
\end{equation*}%
Since%
\begin{equation*}
B(x_{1}x_{2}\otimes gx_{1})\overset{\left( \ref{form x1x2otgx1}\right) }{=}%
(1_{A}\otimes g)B(gx_{1}x_{2}\otimes 1_{H})(1_{A}\otimes
gx_{1})+(1_{A}\otimes x_{1})B(gx_{1}x_{2}\otimes 1_{H})
\end{equation*}%
we get%
\begin{eqnarray*}
\text{left side }B(x_{1}x_{2}\otimes gx_{1})(1_{A}\otimes x_{1})
&=&(1_{A}\otimes g)B(gx_{1}x_{2}\otimes 1_{H})(1_{A}\otimes
gx_{1})(1_{A}\otimes x_{1}) \\
&&+(1_{A}\otimes x_{1})B(gx_{1}x_{2}\otimes 1_{H})(1_{A}\otimes x_{1})
\end{eqnarray*}%
i.e.

\begin{equation*}
\text{left side }B(x_{1}x_{2}\otimes gx_{1})(1_{A}\otimes
x_{1})=(1_{A}\otimes x_{1})B(gx_{1}x_{2}\otimes 1_{H})(1_{A}\otimes x_{1})
\end{equation*}%
and
\begin{eqnarray*}
&&\text{right side }(1_{A}\otimes gx_{1})B(x_{1}x_{2}\otimes
gx_{1})(1_{A}\otimes g) \\
&=&(1_{A}\otimes gx_{1})(1_{A}\otimes g)B(gx_{1}x_{2}\otimes
1_{H})(1_{A}\otimes gx_{1})(1_{A}\otimes g) \\
&&+(1_{A}\otimes gx_{1})(1_{A}\otimes x_{1})B(gx_{1}x_{2}\otimes
1_{H})(1_{A}\otimes g)
\end{eqnarray*}%
i.e.%
\begin{equation*}
\text{right side }(1_{A}\otimes gx_{1})B(x_{1}x_{2}\otimes
gx_{1})(1_{A}\otimes g)=(1_{A}\otimes x_{1})B(gx_{1}x_{2}\otimes
1_{H})(1_{A}\otimes x_{1})
\end{equation*}%
and we conclude.

\subsubsection{$x_{1}x_{2}\otimes gx_{2}$}

\begin{equation*}
B(x_{1}x_{2}\otimes gx_{2})(1_{A}\otimes x_{1})\overset{?}{=}(1_{A}\otimes
x_{1})B(gx_{1}x_{2}\otimes ggx_{2})+B(x_{1}x_{1}x_{2}\otimes
ggx_{2})+B(x_{1}x_{2}\otimes x_{1}gx_{2})
\end{equation*}%
i.e.%
\begin{equation*}
B(x_{1}x_{2}\otimes gx_{2})(1_{A}\otimes x_{1})\overset{?}{=}(1_{A}\otimes
x_{1})(1_{A}\otimes g)B(x_{1}x_{2}\otimes x_{2})(1_{A}\otimes
g)-B(x_{1}x_{2}\otimes gx_{1}x_{2})
\end{equation*}%
Since%
\begin{equation*}
B(x_{1}x_{2}\otimes gx_{2})\overset{\left( \ref{form x1x2otgx2}\right) }{=}%
(1_{A}\otimes g)B(gx_{1}x_{2}\otimes 1_{H})(1_{A}\otimes
gx_{2})+(1_{A}\otimes x_{2})B(gx_{1}x_{2}\otimes 1_{H})
\end{equation*}%
we get%
\begin{eqnarray*}
&&\text{left side }B(x_{1}x_{2}\otimes gx_{2})(1_{A}\otimes x_{1}) \\
&&\overset{\left( \ref{form x1x2otgx2}\right) }{=}-(1_{A}\otimes
g)B(gx_{1}x_{2}\otimes 1_{H})(1_{A}\otimes gx_{1}x_{2})+(1_{A}\otimes
x_{2})B(gx_{1}x_{2}\otimes 1_{H})(1_{A}\otimes x_{1}).
\end{eqnarray*}%
Since
\begin{eqnarray*}
&&B(x_{1}x_{2}\otimes gx_{1}x_{2})\overset{\left( \ref{form x1x2otgx1gx2}%
\right) }{=}(1_{A}\otimes g)B(gx_{1}x_{2}\otimes 1_{H})(1_{A}\otimes
gx_{1}x_{2})+ \\
&&-(1_{A}\otimes x_{2})B(gx_{1}x_{2}\otimes 1_{H})(1_{A}\otimes x_{1}) \\
&&+(1_{A}\otimes x_{1})B(gx_{1}x_{2}\otimes 1_{H})(1_{A}\otimes x_{2}) \\
&&+(1_{A}\otimes gx_{1}x_{2})B(gx_{1}x_{2}\otimes 1_{H})(1_{A}\otimes g)
\end{eqnarray*}%
we get%
\begin{eqnarray*}
&&\text{right side }(1_{A}\otimes x_{1})(1_{A}\otimes g)B(x_{1}x_{2}\otimes
x_{2})(1_{A}\otimes g)-B(x_{1}x_{2}\otimes gx_{1}x_{2}) \\
&=&(1_{A}\otimes x_{1})(1_{A}\otimes g)(1_{A}\otimes g)B(gx_{1}x_{2}\otimes
1_{H})(1_{A}\otimes gx_{2})(1_{A}\otimes g) \\
&&+(1_{A}\otimes x_{1})(1_{A}\otimes g)(1_{A}\otimes
x_{2})B(gx_{1}x_{2}\otimes 1_{H})(1_{A}\otimes g)+ \\
&&-(1_{A}\otimes g)B(gx_{1}x_{2}\otimes 1_{H})(1_{A}\otimes gx_{1}x_{2})+ \\
&&+(1_{A}\otimes x_{2})B(gx_{1}x_{2}\otimes 1_{H})(1_{A}\otimes x_{1})+ \\
&&-(1_{A}\otimes x_{1})B(gx_{1}x_{2}\otimes 1_{H})(1_{A}\otimes x_{2})+ \\
&&-(1_{A}\otimes gx_{1}x_{2})B(gx_{1}x_{2}\otimes 1_{H})(1_{A}\otimes g)
\end{eqnarray*}%
\begin{eqnarray*}
&&\text{clean right side }(1_{A}\otimes x_{1})(1_{A}\otimes
g)B(x_{1}x_{2}\otimes x_{2})(1_{A}\otimes g)-B(x_{1}x_{2}\otimes gx_{1}x_{2})
\\
&=&-(1_{A}\otimes g)B(gx_{1}x_{2}\otimes 1_{H})(1_{A}\otimes
gx_{1}x_{2})+(1_{A}\otimes x_{2})B(gx_{1}x_{2}\otimes 1_{H})(1_{A}\otimes
x_{1})
\end{eqnarray*}%
and we conclude.

\subsubsection{$x_{1}x_{2}\otimes gx_{1}x_{2}$}

\begin{eqnarray*}
&&B(x_{1}x_{2}\otimes gx_{1}x_{2})(1_{A}\otimes x_{1}) \\
&&\overset{?}{=}(1_{A}\otimes x_{1})B(gx_{1}x_{2}\otimes ggx_{1}x_{2}) \\
&&+B(x_{1}x_{1}x_{2}\otimes ggx_{1}x_{2})+B(x_{1}x_{2}\otimes
x_{1}gx_{1}x_{2})
\end{eqnarray*}%
i.e.%
\begin{equation*}
B(x_{1}x_{2}\otimes gx_{1}x_{2})(1_{A}\otimes x_{1})\overset{?}{=}%
(1_{A}\otimes gx_{1})B(x_{1}x_{2}\otimes gx_{1}x_{2})(1_{A}\otimes g)
\end{equation*}%
Since

\begin{eqnarray*}
&&B(x_{1}x_{2}\otimes gx_{1}x_{2})\overset{\left( \ref{form x1x2otgx1gx2}%
\right) }{=}(1_{A}\otimes g)B(gx_{1}x_{2}\otimes 1_{H})(1_{A}\otimes
gx_{1}x_{2}) \\
&&-(1_{A}\otimes x_{2})B(gx_{1}x_{2}\otimes 1_{H})(1_{A}\otimes x_{1}) \\
&&+(1_{A}\otimes x_{1})B(gx_{1}x_{2}\otimes 1_{H})(1_{A}\otimes x_{2}) \\
&&+(1_{A}\otimes gx_{1}x_{2})B(gx_{1}x_{2}\otimes 1_{H})(1_{A}\otimes g)
\end{eqnarray*}%
we get%
\begin{eqnarray*}
\text{left side } &&B(x_{1}x_{2}\otimes gx_{1}x_{2})(1_{A}\otimes
x_{1})=(1_{A}\otimes g)B(gx_{1}x_{2}\otimes 1_{H})(1_{A}\otimes
gx_{1}x_{2})(1_{A}\otimes x_{1}) \\
&&-(1_{A}\otimes x_{2})B(gx_{1}x_{2}\otimes 1_{H})(1_{A}\otimes
x_{1})(1_{A}\otimes x_{1}) \\
&&+(1_{A}\otimes x_{1})B(gx_{1}x_{2}\otimes 1_{H})(1_{A}\otimes
x_{2})(1_{A}\otimes x_{1}) \\
&&+(1_{A}\otimes gx_{1}x_{2})B(gx_{1}x_{2}\otimes 1_{H})(1_{A}\otimes
g)(1_{A}\otimes x_{1})
\end{eqnarray*}%
i.e.%
\begin{eqnarray*}
\text{left side } &&B(x_{1}x_{2}\otimes gx_{1}x_{2})(1_{A}\otimes x_{1})= \\
&&+(1_{A}\otimes x_{1})B(gx_{1}x_{2}\otimes 1_{H})(1_{A}\otimes x_{1}x_{2})
\\
&&-(1_{A}\otimes gx_{1}x_{2})B(gx_{1}x_{2}\otimes 1_{H})(1_{A}\otimes gx_{1})
\end{eqnarray*}%
and%
\begin{eqnarray*}
&&\text{right side}(1_{A}\otimes gx_{1})B(x_{1}x_{2}\otimes
gx_{1}x_{2})(1_{A}\otimes g)= \\
&&(1_{A}\otimes gx_{1})(1_{A}\otimes g)B(gx_{1}x_{2}\otimes
1_{H})(1_{A}\otimes gx_{1}x_{2})(1_{A}\otimes g) \\
&&-(1_{A}\otimes gx_{1})(1_{A}\otimes x_{2})B(gx_{1}x_{2}\otimes
1_{H})(1_{A}\otimes x_{1})(1_{A}\otimes g) \\
&&+(1_{A}\otimes gx_{1})(1_{A}\otimes x_{1})B(gx_{1}x_{2}\otimes
1_{H})(1_{A}\otimes x_{2})(1_{A}\otimes g) \\
&&+(1_{A}\otimes gx_{1})(1_{A}\otimes gx_{1}x_{2})B(gx_{1}x_{2}\otimes
1_{H})(1_{A}\otimes g)(1_{A}\otimes g)
\end{eqnarray*}%
i.e.%
\begin{eqnarray*}
&&\text{right side}(1_{A}\otimes gx_{1})B(x_{1}x_{2}\otimes
gx_{1}x_{2})(1_{A}\otimes g)= \\
&&(1_{A}\otimes x_{1})B(gx_{1}x_{2}\otimes 1_{H})(1_{A}\otimes x_{1}x_{2}) \\
&&-(1_{A}\otimes gx_{1}x_{2})(1_{A}\otimes )B(gx_{1}x_{2}\otimes
1_{H})(1_{A}\otimes gx_{1})
\end{eqnarray*}%
and we conclude.

\subsection{Case $1_{A}\otimes x_{2}$}

\begin{proposition}
Assume that equation $\left( \ref{eq.h}\right) $ holds for $1_{A}\otimes
x_{2}$ and $h\otimes h^{\prime }$%
\begin{equation}
B(h\otimes h^{\prime })(1_{A}\otimes x_{2})=(1_{A}\otimes x_{2})B(gh\otimes
gh^{\prime })+B(x_{2}h\otimes gh^{\prime })+B(h\otimes x_{2}h^{\prime })
\label{1x2}
\end{equation}%
then it holds for $gh\otimes gh^{\prime }$ i.e.%
\begin{equation*}
B(gh\otimes gh^{\prime })(1_{A}\otimes x_{2})=(1_{A}\otimes
x_{2})B(ggh\otimes ggh^{\prime })+B(x_{2}gh\otimes ggh^{\prime
})+B(gh\otimes x_{2}gh^{\prime })
\end{equation*}
\end{proposition}

\begin{proof}
We conjugate equality $\left( \ref{1x2}\right) $ by $\left( 1_{A}\otimes
g\right) $ on both sides and we get%
\begin{eqnarray*}
-\left( 1_{A}\otimes g\right) B(h\otimes h^{\prime })\left( 1_{A}\otimes
g\right) (1_{A}\otimes x_{2}) &=&-(1_{A}\otimes x_{2})\left( 1_{A}\otimes
g\right) B(gh\otimes gh^{\prime })\left( 1_{A}\otimes g\right) \\
&&+\left( 1_{A}\otimes g\right) B(x_{2}h\otimes gh^{\prime })\left(
1_{A}\otimes g\right) \\
&&+\left( 1_{A}\otimes g\right) B(h\otimes x_{2}h^{\prime })\left(
1_{A}\otimes g\right) .
\end{eqnarray*}%
From this we deduce that%
\begin{eqnarray*}
B(gh\otimes gh^{\prime })(1_{A}\otimes x_{2}) &=&(1_{A}\otimes
x_{2})B(h\otimes h^{\prime }) \\
&&-B(gx_{2}h\otimes h^{\prime })-B(gh\otimes gx_{2}h^{\prime }).
\end{eqnarray*}
\end{proof}

\subsubsection{$g\otimes 1_{H}$}

\begin{equation*}
B(g\otimes 1_{H})(1_{A}\otimes x_{2})\overset{?}{=}(1_{A}\otimes
x_{2})B(gg\otimes g)-B(gx_{2}\otimes g)+B(g\otimes x_{2})
\end{equation*}%
i.e.%
\begin{eqnarray*}
&&B(g\otimes 1_{H})(1_{A}\otimes x_{2})\overset{?}{=}(1_{A}\otimes
x_{2})\left( 1_{A}\otimes g\right) B(g\otimes 1_{H})\left( 1_{A}\otimes
g\right) -\left( 1_{A}\otimes g\right) B(x_{2}\otimes 1_{H})\left(
1_{A}\otimes g\right) \\
&&+\left( 1_{A}\otimes g\right) B(1_{H}\otimes gx_{2})\left( 1_{A}\otimes
g\right)
\end{eqnarray*}%
Since%
\begin{eqnarray*}
&&B(1_{H}\otimes gx_{2})\overset{\left( \ref{form 1otgx2}\right) }{=}%
(1_{A}\otimes g)B(g\otimes 1_{H})(1_{A}\otimes gx_{2}) \\
&&+(1_{A}\otimes x_{2})B(g\otimes 1_{H})+B(x_{2}\otimes 1_{H})
\end{eqnarray*}%
we get
\begin{eqnarray*}
&&\text{right side }(1_{A}\otimes x_{2})\left( 1_{A}\otimes g\right)
B(g\otimes 1_{H})\left( 1_{A}\otimes g\right) -\left( 1_{A}\otimes g\right)
B(x_{2}\otimes 1_{H})\left( 1_{A}\otimes g\right) + \\
&&+\left( 1_{A}\otimes g\right) (1_{A}\otimes g)B(g\otimes
1_{H})(1_{A}\otimes gx_{2})\left( 1_{A}\otimes g\right) + \\
&&+\left( 1_{A}\otimes g\right) (1_{A}\otimes x_{2})B(g\otimes 1_{H})\left(
1_{A}\otimes g\right) +\left( 1_{A}\otimes g\right) B(x_{2}\otimes
1_{H})\left( 1_{A}\otimes g\right)
\end{eqnarray*}%
i.e.%
\begin{equation*}
\text{clean }\text{right side }B(g\otimes 1_{H})(1_{A}\otimes x_{2})
\end{equation*}%
and we conclude.

\subsubsection{$x_{1}\otimes 1_{H}$}

.%
\begin{equation*}
B(x_{1}\otimes 1_{H})(1_{A}\otimes x_{2})\overset{?}{=}(1_{A}\otimes
x_{2})B(gx_{1}\otimes g)+B(x_{2}x_{1}\otimes g)+B(x_{1}\otimes x_{2})
\end{equation*}%
i.e%
\begin{eqnarray*}
&&B(x_{1}\otimes 1_{H})(1_{A}\otimes x_{2})\overset{?}{=}(1_{A}\otimes
x_{2})\left( 1_{A}\otimes g\right) B(x_{1}\otimes 1_{H})\left( 1_{A}\otimes
g\right) \\
&&-\left( 1_{A}\otimes g\right) B(gx_{1}x_{2}\otimes 1_{H})\left(
1_{A}\otimes g\right) +B(x_{1}\otimes x_{2})
\end{eqnarray*}%
Since%
\begin{equation*}
B(x_{1}\otimes x_{2})\overset{\left( \ref{form x1otx2}\right) }{=}%
B(x_{1}\otimes 1_{H})(1_{A}\otimes x_{2})-(1_{A}\otimes
gx_{2})B(x_{1}\otimes 1_{H})(1_{A}\otimes g)
\end{equation*}%
we obtain%
\begin{eqnarray*}
&&\text{right side }(1_{A}\otimes gx_{2})B(x_{1}\otimes 1_{H})\left(
1_{A}\otimes g\right) -\left( 1_{A}\otimes g\right) B(gx_{1}x_{2}\otimes
1_{H})\left( 1_{A}\otimes g\right) \\
&&+B(x_{2}\otimes x_{2}) \\
&=&(1_{A}\otimes x_{2})\left( 1_{A}\otimes g\right) B(x_{1}\otimes
1_{H})\left( 1_{A}\otimes g\right) -\left( 1_{A}\otimes g\right)
B(gx_{1}x_{2}\otimes 1_{H})\left( 1_{A}\otimes g\right) + \\
&&+B(x_{1}\otimes 1_{H})(1_{A}\otimes x_{2})-(1_{A}\otimes
gx_{2})B(x_{1}\otimes 1_{H})(1_{A}\otimes g)+ \\
&&+(1_{A}\otimes g)B(gx_{1}x_{2}\otimes 1_{H})(1_{A}\otimes g)
\end{eqnarray*}%
i.e.%
\begin{eqnarray*}
&&\text{clean right side }(1_{A}\otimes gx_{2})B(x_{1}\otimes 1_{H})\left(
1_{A}\otimes g\right) \\
&&-\left( 1_{A}\otimes g\right) B(gx_{1}x_{2}\otimes 1_{H})\left(
1_{A}\otimes g\right) +B(x_{2}\otimes x_{2}) \\
&=&B(x_{1}\otimes 1_{H})(1_{A}\otimes x_{2})
\end{eqnarray*}%
and we conclude.

\subsubsection{$x_{2}\otimes 1_{H}$}

\begin{equation*}
B(x_{2}\otimes 1_{H})(1_{A}\otimes x_{2})\overset{?}{=}(1_{A}\otimes
x_{2})B(gx_{2}\otimes g)+B(x_{2}x_{2}\otimes g)+B(x_{2}\otimes x_{2})
\end{equation*}%
i.e.%
\begin{equation*}
B(x_{2}\otimes 1_{H})(1_{A}\otimes x_{2})\overset{?}{=}(1_{A}\otimes
x_{2})\left( 1_{A}\otimes g\right) B(x_{2}\otimes 1_{H})\left( 1_{A}\otimes
g\right) +B(x_{2}\otimes x_{2})
\end{equation*}%
Since%
\begin{equation*}
B(x_{2}\otimes x_{2})\overset{\left( \ref{form x2otx2}\right) }{=}%
B(x_{2}\otimes 1_{H})(1_{A}\otimes x_{2})-(1_{A}\otimes
gx_{2})B(x_{2}\otimes 1_{H})(1_{A}\otimes g)
\end{equation*}%
we get%
\begin{eqnarray*}
&&\text{right side }(1_{A}\otimes x_{2})\left( 1_{A}\otimes g\right)
B(x_{2}\otimes 1_{H})\left( 1_{A}\otimes g\right) +B(x_{2}\otimes x_{2}) \\
&=&(1_{A}\otimes x_{2})\left( 1_{A}\otimes g\right) B(x_{2}\otimes
1_{H})\left( 1_{A}\otimes g\right) + \\
&&B(x_{2}\otimes 1_{H})(1_{A}\otimes x_{2})-(1_{A}\otimes
gx_{2})B(x_{2}\otimes 1_{H})(1_{A}\otimes g) \\
&=&B(x_{2}\otimes 1_{H})(1_{A}\otimes x_{2})
\end{eqnarray*}%
and we conclude.

\subsubsection{$x_{1}x_{2}\otimes 1_{H}$}

\begin{equation*}
B(x_{1}x_{2}\otimes 1_{H})(1_{A}\otimes x_{2})\overset{?}{=}(1_{A}\otimes
x_{2})B(gx_{1}x_{2}\otimes g)+B(x_{2}x_{1}x_{2}\otimes
g)+B(x_{1}x_{2}\otimes x_{2})
\end{equation*}%
i.e.%
\begin{equation*}
B(x_{1}x_{2}\otimes 1_{H})(1_{A}\otimes x_{2})\overset{?}{=}(1_{A}\otimes
x_{2})\left( 1_{A}\otimes g\right) B(x_{1}x_{2}\otimes 1_{H})\left(
1_{A}\otimes g\right) +B(x_{1}x_{2}\otimes x_{2}).
\end{equation*}%
Since
\begin{equation*}
B(x_{1}x_{2}\otimes x_{2})\overset{\left( \ref{x1x2otx2}\right) }{=}%
B(x_{1}x_{2}\otimes 1_{H})(1_{A}\otimes x_{2})-(1_{A}\otimes
gx_{2})B(x_{1}x_{2}\otimes 1_{H})(1_{A}\otimes g)
\end{equation*}%
we get%
\begin{eqnarray*}
&&\text{right side }(1_{A}\otimes x_{2})\left( 1_{A}\otimes g\right)
B(x_{1}x_{2}\otimes 1_{H})\left( 1_{A}\otimes g\right) +B(x_{1}x_{2}\otimes
x_{2}) \\
&=&(1_{A}\otimes x_{2})\left( 1_{A}\otimes g\right) B(x_{1}x_{2}\otimes
1_{H})\left( 1_{A}\otimes g\right) + \\
&&B(x_{1}x_{2}\otimes 1_{H})(1_{A}\otimes x_{2})-(1_{A}\otimes
gx_{2})B(x_{1}x_{2}\otimes 1_{H})(1_{A}\otimes g) \\
&=&B(x_{1}x_{2}\otimes 1_{H})(1_{A}\otimes x_{2})
\end{eqnarray*}%
and we conclude.

\subsubsection{$gx_{1}\otimes 1_{H}$}

\begin{equation*}
B(gx_{1}\otimes 1_{H})(1_{A}\otimes x_{2})\overset{?}{=}(1_{A}\otimes
x_{2})\left( 1_{A}\otimes g\right) B(gx_{1}\otimes 1_{H})\left( 1_{A}\otimes
g\right) +B(x_{2}gx_{1}\otimes g)+B(gx_{1}\otimes x_{2})
\end{equation*}%
i.e.%
\begin{eqnarray*}
&&B(gx_{1}\otimes 1_{H})(1_{A}\otimes x_{2})\overset{?}{=}(1_{A}\otimes
x_{2})\left( 1_{A}\otimes g\right) B(gx_{1}\otimes 1_{H})\left( 1_{A}\otimes
g\right) \\
&&+\left( 1_{A}\otimes g\right) B(x_{1}x_{2}\otimes 1_{H})\left(
1_{A}\otimes g\right) +\left( 1_{A}\otimes g\right) B(x_{1}\otimes
gx_{2})1_{A}\otimes g
\end{eqnarray*}%
Since%
\begin{equation*}
B(x_{1}\otimes gx_{2})\overset{\left( \ref{form x1otgx2}\right) }{=}%
(1_{A}\otimes g)B(gx_{1}\otimes 1_{H})(1_{A}\otimes gx_{2})+(1_{A}\otimes
x_{2})B(gx_{1}\otimes 1_{H})-B(x_{1}x_{2}\otimes 1_{H})
\end{equation*}%
we get%
\begin{eqnarray*}
&&\text{right side }(1_{A}\otimes x_{2})\left( 1_{A}\otimes g\right)
B(gx_{1}\otimes 1_{H})\left( 1_{A}\otimes g\right) + \\
&&+\left( 1_{A}\otimes g\right) B(x_{1}x_{2}\otimes 1_{H})\left(
1_{A}\otimes g\right) + \\
&&+(1_{A}\otimes g)(1_{A}\otimes g)B(gx_{1}\otimes 1_{H})(1_{A}\otimes
gx_{2})(1_{A}\otimes g)+ \\
&&+(1_{A}\otimes g)(1_{A}\otimes x_{2})B(gx_{1}\otimes 1_{H})(1_{A}\otimes
g)+ \\
&&-(1_{A}\otimes g)B(x_{1}x_{2}\otimes 1_{H})(1_{A}\otimes g)
\end{eqnarray*}%
i.e.%
\begin{equation*}
\text{right side }B(gx_{1}\otimes 1_{H})(1_{A}\otimes x_{2})
\end{equation*}

and we conclude.

\subsubsection{$gx_{2}\otimes 1_{H}$}

\begin{equation*}
B(gx_{2}\otimes 1_{H})(1_{A}\otimes x_{2})\overset{?}{=}(1_{A}\otimes
x_{2})\left( 1_{A}\otimes g\right) B(gx_{2}\otimes 1_{H})\left( 1_{A}\otimes
g\right) +B(x_{2}gx_{2}\otimes g)+B(gx_{2}\otimes x_{2})
\end{equation*}%
i.e.%
\begin{equation*}
B(gx_{2}\otimes 1_{H})(1_{A}\otimes x_{2})\overset{?}{=}(1_{A}\otimes
x_{2})\left( 1_{A}\otimes g\right) B(gx_{2}\otimes 1_{H})\left( 1_{A}\otimes
g\right) +\left( 1_{A}\otimes g\right) B(x_{2}\otimes gx_{2})\left(
1_{A}\otimes g\right)
\end{equation*}%
Since
\begin{equation*}
B(x_{2}\otimes gx_{2})\overset{\left( \ref{form x2otgx2}\right) }{=}%
(1_{A}\otimes g)B(gx_{2}\otimes 1_{H})(1_{A}\otimes gx_{2})+(1_{A}\otimes
x_{2})B(gx_{2}\otimes 1_{H})
\end{equation*}%
we get%
\begin{eqnarray*}
&&\text{right side }(1_{A}\otimes x_{2})\left( 1_{A}\otimes g\right)
B(gx_{2}\otimes 1_{H})\left( 1_{A}\otimes g\right) +\left( 1_{A}\otimes
g\right) B(x_{2}\otimes gx_{2})\left( 1_{A}\otimes g\right) \\
&=&(1_{A}\otimes x_{2})\left( 1_{A}\otimes g\right) B(gx_{2}\otimes
1_{H})\left( 1_{A}\otimes g\right) + \\
&&\left( 1_{A}\otimes g\right) (1_{A}\otimes g)B(gx_{2}\otimes
1_{H})(1_{A}\otimes gx_{2})\left( 1_{A}\otimes g\right) \\
&&+\left( 1_{A}\otimes g\right) (1_{A}\otimes x_{2})B(gx_{2}\otimes
1_{H})\left( 1_{A}\otimes g\right)
\end{eqnarray*}%
i.e.%
\begin{eqnarray*}
&&\text{clean right side }(1_{A}\otimes x_{2})\left( 1_{A}\otimes g\right)
B(gx_{2}\otimes 1_{H})\left( 1_{A}\otimes g\right) +\left( 1_{A}\otimes
g\right) B(x_{2}\otimes gx_{2})\left( 1_{A}\otimes g\right) \\
&=&B(gx_{2}\otimes 1_{H})(1_{A}\otimes x_{2})
\end{eqnarray*}%
and we conclude.

\subsubsection{$gx_{1}x_{2}\otimes 1_{H}$}

\begin{equation*}
B(gx_{1}x_{2}\otimes 1_{H})(1_{A}\otimes x_{2})\overset{?}{=}(1_{A}\otimes
x_{2})\left( 1_{A}\otimes g\right) B(gx_{1}x_{2}\otimes 1_{H})\left(
1_{A}\otimes g\right) +B(x_{2}gx_{1}x_{2}\otimes g)+B(gx_{1}x_{2}\otimes
x_{2})
\end{equation*}%
i.e.%
\begin{eqnarray*}
&&B(gx_{1}x_{2}\otimes 1_{H})(1_{A}\otimes x_{2})\overset{?}{=}(1_{A}\otimes
x_{2})\left( 1_{A}\otimes g\right) B(gx_{1}x_{2}\otimes 1_{H})\left(
1_{A}\otimes g\right) + \\
&&+\left( 1_{A}\otimes g\right) B(x_{1}x_{2}\otimes gx_{2})\left(
1_{A}\otimes g\right) .
\end{eqnarray*}%
Since
\begin{equation*}
B(x_{1}x_{2}\otimes gx_{2})\overset{\left( \ref{form x1x2otgx2}\right) }{=}%
(1_{A}\otimes g)B(gx_{1}x_{2}\otimes 1_{H})(1_{A}\otimes
gx_{2})+(1_{A}\otimes x_{2})B(gx_{1}x_{2}\otimes 1_{H})
\end{equation*}%
we get%
\begin{eqnarray*}
&&\text{right side }(1_{A}\otimes x_{2})\left( 1_{A}\otimes g\right)
B(gx_{1}x_{2}\otimes 1_{H})\left( 1_{A}\otimes g\right) + \\
&&+\left( 1_{A}\otimes g\right) B(x_{1}x_{2}\otimes gx_{2})\left(
1_{A}\otimes g\right) \\
&=&(1_{A}\otimes x_{2})\left( 1_{A}\otimes g\right) B(gx_{1}x_{2}\otimes
1_{H})\left( 1_{A}\otimes g\right) + \\
&&\left( 1_{A}\otimes g\right) (1_{A}\otimes g)B(gx_{1}x_{2}\otimes
1_{H})(1_{A}\otimes gx_{2})\left( 1_{A}\otimes g\right) + \\
&&\left( 1_{A}\otimes g\right) (1_{A}\otimes x_{2})B(gx_{1}x_{2}\otimes
1_{H})\left( 1_{A}\otimes g\right)
\end{eqnarray*}%
i.e.%
\begin{eqnarray*}
&&\text{right side }(1_{A}\otimes x_{2})\left( 1_{A}\otimes g\right)
B(gx_{1}x_{2}\otimes 1_{H})\left( 1_{A}\otimes g\right) + \\
&&+\left( 1_{A}\otimes g\right) B(x_{1}x_{2}\otimes gx_{2})\left(
1_{A}\otimes g\right) \\
&=&B(gx_{1}x_{2}\otimes 1_{H})(1_{A}\otimes x_{2})
\end{eqnarray*}%
and we conclude.

\subsubsection{$1_{H}\otimes x_{1}$}

\begin{equation*}
B(1\otimes x_{1})(1_{A}\otimes x_{2})\overset{?}{=}(1_{A}\otimes
x_{2})\left( 1_{A}\otimes g\right) B(1_{H}\otimes x_{1})\left( 1_{A}\otimes
g\right) +B(x_{2}\otimes gx_{1})-B(1_{H}\otimes x_{1}x_{2})
\end{equation*}%
We consider the left side%
\begin{equation*}
B(1_{H}\otimes x_{1})(1_{H}\otimes x_{2})\overset{\left( \ref{form:1otxi}%
\right) }{=}(1_{H}\otimes g)B(gx_{1}\otimes 1_{H})(1_{H}\otimes gx_{2})
\end{equation*}%
We have%
\begin{eqnarray*}
&&(1_{A}\otimes x_{2})\left( 1_{A}\otimes g\right) B(1_{H}\otimes
x_{1})\left( 1_{A}\otimes g\right) \overset{\left( \ref{form:1otxi}\right) }{%
=} \\
&=&-(1_{A}\otimes x_{2})\left( 1_{A}\otimes g\right) (1_{H}\otimes
g)B(gx_{1}\otimes 1_{H})(1_{H}\otimes g)\left( 1_{A}\otimes g\right) \\
&=&-(1_{A}\otimes x_{2})B(gx_{1}\otimes 1_{H}).
\end{eqnarray*}%
Since%
\begin{equation*}
B(x_{2}\otimes gx_{1})\overset{\left( \ref{form x2otgx1}\right) }{=}%
(1_{A}\otimes g)B(gx_{2}\otimes 1_{H})(1_{A}\otimes gx_{1})+(1_{A}\otimes
x_{1})B(gx_{2}\otimes 1_{H})+B(x_{1}x_{2}\otimes 1_{H}).
\end{equation*}%
and%
\begin{eqnarray*}
&&B(1_{H}\otimes x_{1}x_{2}) \\
&&\overset{\left( \ref{form 1otx1x2}\right) }{=}+(1_{A}\otimes
g)B(gx_{2}\otimes 1_{H})(1_{A}\otimes gx_{1}) \\
&&+(1_{A}\otimes x_{1})B(gx_{2}\otimes 1_{H}) \\
&&-(1_{A}\otimes g)B(gx_{1}\otimes 1_{H})(1_{A}\otimes gx_{2}) \\
&&-(1_{A}\otimes x_{2})B(gx_{1}\otimes 1_{H}) \\
&&+B(x_{1}x_{2}\otimes 1_{H})
\end{eqnarray*}%
We get%
\begin{eqnarray*}
&&\text{right side }(1_{A}\otimes x_{2})\left( 1_{A}\otimes g\right)
B(1_{H}\otimes x_{1})\left( 1_{A}\otimes g\right) +B(x_{2}\otimes
gx_{1})-B(1_{H}\otimes x_{1}x_{2}) \\
&=&-(1_{A}\otimes x_{2})B(gx_{1}\otimes 1_{H})+ \\
&&+(1_{A}\otimes g)B(gx_{2}\otimes 1_{H})(1_{A}\otimes gx_{1})+ \\
&&+(1_{A}\otimes x_{1})B(gx_{2}\otimes 1_{H}) \\
&&+B(x_{1}x_{2}\otimes 1_{H})+ \\
&&-(1_{A}\otimes g)B(gx_{2}\otimes 1_{H})(1_{A}\otimes gx_{1})+ \\
&&-(1_{A}\otimes x_{1})B(gx_{2}\otimes 1_{H})+ \\
&&+(1_{A}\otimes g)B(gx_{1}\otimes 1_{H})(1_{A}\otimes gx_{2})+ \\
&&+(1_{A}\otimes x_{2})B(gx_{1}\otimes 1_{H})+ \\
&&-B(x_{1}x_{2}\otimes 1_{H})
\end{eqnarray*}%
and hence we have%
\begin{eqnarray*}
&&\text{clean right side }(1_{A}\otimes x_{2})\left( 1_{A}\otimes g\right)
B(1_{H}\otimes x_{1})\left( 1_{A}\otimes g\right) +B(x_{2}\otimes
gx_{1})-B(1_{H}\otimes x_{1}x_{2}) \\
&=&(1_{A}\otimes g)B(gx_{1}\otimes 1_{H})(1_{A}\otimes gx_{2})
\end{eqnarray*}%
and we conclude.

\subsubsection{$1_{H}\otimes x_{2}$}

\begin{equation*}
B(1_{H}\otimes x_{2})(1_{A}\otimes x_{2})\overset{?}{=}(1_{A}\otimes
x_{2})\left( 1_{A}\otimes g\right) B(1_{H}\otimes x_{2})\left( 1_{A}\otimes
g\right) +B(x_{2}\otimes gx_{2})+B(1_{H}\otimes x_{2}x_{2})
\end{equation*}%
i.e.%
\begin{equation*}
B(1_{H}\otimes x_{2})(1_{A}\otimes x_{2})\overset{?}{=}(1_{A}\otimes
x_{2})\left( 1_{A}\otimes g\right) B(1_{H}\otimes x_{2})\left( 1_{A}\otimes
g\right) +B(x_{2}\otimes gx_{2})
\end{equation*}%
Since%
\begin{equation*}
B(1_{H}\otimes x_{2})\overset{\left( \ref{form:1otxi}\right) }{=}%
-(1_{A}\otimes g)B(gx_{2}\otimes 1_{H})(1_{A}\otimes g)
\end{equation*}%
we get%
\begin{eqnarray*}
&&\text{left side }B(1_{H}\otimes x_{2})(1_{A}\otimes x_{2})\overset{\left( %
\ref{form:1otxi}\right) }{=} \\
&=&-(1_{A}\otimes g)B(gx_{2}\otimes 1_{H})(1_{A}\otimes g)(1_{A}\otimes
x_{2}) \\
&=&(1_{A}\otimes g)B(gx_{2}\otimes 1_{H})(1_{A}\otimes gx_{2})
\end{eqnarray*}%
We have%
\begin{eqnarray*}
&&(1_{A}\otimes x_{2})\left( 1_{A}\otimes g\right) B(1_{H}\otimes
x_{2})\left( 1_{A}\otimes g\right) \overset{\left( \ref{form:1otxi}\right) }{%
=} \\
&=&-(1_{A}\otimes x_{2})\left( 1_{A}\otimes g\right) (1_{A}\otimes
g)B(gx_{2}\otimes 1_{H})(1_{A}\otimes g)(1_{A}\otimes g) \\
&=&-(1_{A}\otimes x_{2})B(gx_{2}\otimes 1_{H})
\end{eqnarray*}%
Since%
\begin{equation*}
B(x_{2}\otimes gx_{2})\overset{\left( \ref{form x2otgx2}\right) }{=}%
(1_{A}\otimes g)B(gx_{2}\otimes 1_{H})(1_{A}\otimes gx_{2})+(1_{A}\otimes
x_{2})B(gx_{2}\otimes 1_{H})
\end{equation*}%
we get%
\begin{eqnarray*}
&&\text{right side }(1_{A}\otimes x_{2})\left( 1_{A}\otimes g\right)
B(1_{H}\otimes x_{2})\left( 1_{A}\otimes g\right) +B(x_{2}\otimes gx_{2}) \\
&=&-(1_{A}\otimes x_{2})B(gx_{2}\otimes 1_{H})+(1_{A}\otimes
g)B(gx_{2}\otimes 1_{H})(1_{A}\otimes gx_{2})+(1_{A}\otimes
x_{2})B(gx_{2}\otimes 1_{H}) \\
&=&(1_{A}\otimes g)B(gx_{2}\otimes 1_{H})(1_{A}\otimes gx_{2})
\end{eqnarray*}%
and we conclude.

\subsubsection{$1_{H}\otimes x_{1}x_{2}$}

\begin{equation*}
B(1_{H}\otimes x_{1}x_{2})(1_{A}\otimes x_{2})\overset{?}{=}(1_{A}\otimes
x_{2})\left( 1_{A}\otimes g\right) B(1_{H}\otimes x_{1}x_{2})\left(
1_{A}\otimes g\right) +B(x_{2}\otimes gx_{1}x_{2})+B(1_{H}\otimes
x_{2}x_{1}x_{2})
\end{equation*}%
i.e.%
\begin{equation*}
B(1_{H}\otimes x_{1}x_{2})(1_{A}\otimes x_{2})\overset{?}{=}(1_{A}\otimes
x_{2})\left( 1_{A}\otimes g\right) B(1_{H}\otimes x_{1}x_{2})\left(
1_{A}\otimes g\right) +B(x_{2}\otimes gx_{1}x_{2}).
\end{equation*}

Since%
\begin{eqnarray*}
&&B(1_{H}\otimes x_{1}x_{2}) \\
&&\overset{\left( \ref{form 1otx1x2}\right) }{=}+(1_{A}\otimes
g)B(gx_{2}\otimes 1_{H})(1_{A}\otimes gx_{1}) \\
&&+(1_{A}\otimes x_{1})B(gx_{2}\otimes 1_{H}) \\
&&-(1_{A}\otimes g)B(gx_{1}\otimes 1_{H})(1_{A}\otimes gx_{2}) \\
&&-(1_{A}\otimes x_{2})B(gx_{1}\otimes 1_{H}) \\
&&+B(x_{1}x_{2}\otimes 1_{H})
\end{eqnarray*}%
we get%
\begin{eqnarray*}
&&\text{left side }B(1_{H}\otimes x_{1}x_{2})(1_{A}\otimes x_{2}) \\
&&\overset{\left( \ref{form 1otx1x2}\right) }{=}+(1_{A}\otimes
g)B(gx_{2}\otimes 1_{H})(1_{A}\otimes gx_{1})(1_{A}\otimes x_{2}) \\
&&+(1_{A}\otimes x_{1})B(gx_{2}\otimes 1_{H})(1_{A}\otimes x_{2}) \\
&&-(1_{A}\otimes g)B(gx_{1}\otimes 1_{H})(1_{A}\otimes gx_{2})(1_{A}\otimes
x_{2}) \\
&&-(1_{A}\otimes x_{2})B(gx_{1}\otimes 1_{H})(1_{A}\otimes x_{2}) \\
&&+B(x_{1}x_{2}\otimes 1_{H})(1_{A}\otimes x_{2})
\end{eqnarray*}%
i.e.%
\begin{eqnarray*}
&&\text{left side }B(1_{H}\otimes x_{1}x_{2})(1_{A}\otimes x_{2}) \\
&&=+(1_{A}\otimes g)B(gx_{2}\otimes 1_{H})(1_{A}\otimes gx_{1}x_{2}) \\
&&+(1_{A}\otimes x_{1})B(gx_{2}\otimes 1_{H})(1_{A}\otimes x_{2}) \\
&&-(1_{A}\otimes x_{2})B(gx_{1}\otimes 1_{H})(1_{A}\otimes x_{2}) \\
&&+B(x_{1}x_{2}\otimes 1_{H})(1_{A}\otimes x_{2})
\end{eqnarray*}%
Since%
\begin{eqnarray*}
&&B(x_{2}\otimes gx_{1}x_{2})\overset{\left( \ref{form x2otgx1x2}\right) }{=}%
(1_{A}\otimes g)B(gx_{2}\otimes 1_{H})(1_{A}\otimes gx_{1}x_{2}) \\
&&-(1_{A}\otimes x_{2})B(gx_{2}\otimes 1_{H})(1_{A}\otimes x_{1}) \\
&&+(1_{A}\otimes x_{1})B(gx_{2}\otimes 1_{H})(1_{A}\otimes x_{2}) \\
&&+(1_{A}\otimes gx_{1}x_{2})B(gx_{2}\otimes 1_{H})(1_{A}\otimes g) \\
&&+B(x_{1}x_{2}\otimes 1_{H})(1_{A}\otimes x_{2}) \\
&&-(1_{A}\otimes gx_{2})B(x_{1}x_{2}\otimes 1_{H})(1_{A}\otimes g)
\end{eqnarray*}%
we get%
\begin{equation*}
\text{right side }(1_{A}\otimes x_{2})\left( 1_{A}\otimes g\right)
B(1_{H}\otimes x_{1}x_{2})\left( 1_{A}\otimes g\right) +B(x_{2}\otimes
gx_{1}x_{2})
\end{equation*}%
\begin{gather*}
\text{right side }(1_{A}\otimes x_{2})\left( 1_{A}\otimes g\right)
B(1_{H}\otimes x_{1}x_{2})\left( 1_{A}\otimes g\right) +B(x_{2}\otimes
gx_{1}x_{2})= \\
+(1_{A}\otimes x_{2})\left( 1_{A}\otimes g\right) (1_{A}\otimes
g)B(gx_{2}\otimes 1_{H})(1_{A}\otimes gx_{1})\left( 1_{A}\otimes g\right) \\
+(1_{A}\otimes x_{2})\left( 1_{A}\otimes g\right) (1_{A}\otimes
x_{1})B(gx_{2}\otimes 1_{H})\left( 1_{A}\otimes g\right) \\
-(1_{A}\otimes x_{2})\left( 1_{A}\otimes g\right) (1_{A}\otimes
g)B(gx_{1}\otimes 1_{H})(1_{A}\otimes gx_{2})\left( 1_{A}\otimes g\right) \\
-(1_{A}\otimes x_{2})\left( 1_{A}\otimes g\right) (1_{A}\otimes
x_{2})B(gx_{1}\otimes 1_{H})\left( 1_{A}\otimes g\right) \\
+(1_{A}\otimes x_{2})\left( 1_{A}\otimes g\right) B(x_{1}x_{2}\otimes
1_{H})\left( 1_{A}\otimes g\right) \\
(1_{A}\otimes g)B(gx_{2}\otimes 1_{H})(1_{A}\otimes gx_{1}x_{2})+ \\
-(1_{A}\otimes x_{2})B(gx_{2}\otimes 1_{H})(1_{A}\otimes x_{1})+ \\
+(1_{A}\otimes x_{1})B(gx_{2}\otimes 1_{H})(1_{A}\otimes x_{2})+ \\
+(1_{A}\otimes gx_{1}x_{2})B(gx_{2}\otimes 1_{H})(1_{A}\otimes g)+ \\
+B(x_{1}x_{2}\otimes 1_{H})(1_{A}\otimes x_{2})+ \\
-(1_{A}\otimes gx_{2})B(x_{1}x_{2}\otimes 1_{H})(1_{A}\otimes g)
\end{gather*}%
\begin{gather*}
\text{clean right side }(1_{A}\otimes x_{2})\left( 1_{A}\otimes g\right)
B(1_{H}\otimes x_{1}x_{2})\left( 1_{A}\otimes g\right) +B(x_{2}\otimes
gx_{1}x_{2})= \\
-(1_{A}\otimes x_{2})B(gx_{1}\otimes 1_{H})(1_{A}\otimes x_{2}) \\
(1_{A}\otimes g)B(gx_{2}\otimes 1_{H})(1_{A}\otimes gx_{1}x_{2})+ \\
+(1_{A}\otimes x_{1})B(gx_{2}\otimes 1_{H})(1_{A}\otimes x_{2})+ \\
+B(x_{1}x_{2}\otimes 1_{H})(1_{A}\otimes x_{2})+
\end{gather*}%
and we conclude.

\subsubsection{$1_{H}\otimes gx_{1}$}

\begin{equation*}
B(1_{H}\otimes gx_{1})(1_{A}\otimes x_{2})\overset{?}{=}(1_{A}\otimes
x_{2})\left( 1_{A}\otimes g\right) B(1_{H}\otimes gx_{1})\left( 1_{A}\otimes
g\right) +B(x_{2}\otimes x_{1})+B(1_{H}\otimes gx_{1}x_{2})
\end{equation*}%
Since%
\begin{eqnarray*}
&&B(1_{H}\otimes gx_{1}) \\
&&\overset{\left( \ref{form 1otgx1}\right) }{=}(1_{A}\otimes g)B(g\otimes
1_{H})(1_{A}\otimes gx_{1}) \\
&&+(1_{A}\otimes x_{1})B(g\otimes 1_{H})+B(x_{1}\otimes 1_{H})
\end{eqnarray*}%
we get%
\begin{eqnarray*}
&&\text{left side }B(1_{H}\otimes gx_{1})(1_{A}\otimes x_{2}) \\
&=&(1_{A}\otimes g)B(g\otimes 1_{H})(1_{A}\otimes gx_{1}x_{2}) \\
&&+(1_{A}\otimes x_{1})B(g\otimes 1_{H})(1_{A}\otimes x_{2}) \\
&&+B(x_{1}\otimes 1_{H})(1_{A}\otimes x_{2})
\end{eqnarray*}%
and we have%
\begin{eqnarray*}
&&(1_{A}\otimes x_{2})\left( 1_{A}\otimes g\right) B(1_{H}\otimes
gx_{1})\left( 1_{A}\otimes g\right) \\
&=&(1_{A}\otimes x_{2})B(g\otimes 1_{H})(1_{A}\otimes x_{1}) \\
&&-\left( 1_{A}\otimes gx_{1}x_{2}\right) B(g\otimes 1_{H})\left(
1_{A}\otimes g\right) \\
&&+(1_{A}\otimes x_{2})\left( 1_{A}\otimes g\right) B(x_{1}\otimes
1_{H})\left( 1_{A}\otimes g\right) .
\end{eqnarray*}%
Since
\begin{eqnarray*}
&&B(x_{2}\otimes x_{1})\overset{\left( \ref{form x2otx1}\right) }{=}%
B(x_{2}\otimes 1_{H})(1_{A}\otimes x_{1}) \\
&&-(1_{A}\otimes gx_{1})B(x_{2}\otimes 1_{H})(1_{A}\otimes g) \\
&&-(1_{A}\otimes g)B(gx_{1}x_{2}\otimes 1_{H})(1_{A}\otimes g)
\end{eqnarray*}%
and%
\begin{eqnarray*}
&&B(1_{H}\otimes gx_{1}x_{2})\overset{\left( \ref{form 1otgx1x2}\right) }{=}%
(1_{A}\otimes gx_{1})B(1_{H}\otimes gx_{2})\left( 1_{A}\otimes g\right) \\
&&+B(x_{1}\otimes x_{2})-B(1_{H}\otimes gx_{2})(1_{A}\otimes x_{1}) \\
&=&+(1_{A}\otimes x_{1})B(g\otimes 1_{H})(1_{A}\otimes x_{2}) \\
&&+(1_{A}\otimes gx_{1}x_{2})B(g\otimes 1_{H})\left( 1_{A}\otimes g\right) \\
&&+(1_{A}\otimes gx_{1})B(x_{2}\otimes 1_{H})\left( 1_{A}\otimes g\right) \\
&&+B(x_{1}\otimes 1_{H})(1_{A}\otimes x_{2}) \\
&&-(1_{A}\otimes gx_{2})B(x_{1}\otimes 1_{H})(1_{A}\otimes g) \\
&&+(1_{A}\otimes g)B(gx_{1}x_{2}\otimes 1_{H})(1_{A}\otimes g) \\
&&+(1_{A}\otimes g)B(g\otimes 1_{H})(1_{A}\otimes gx_{1}x_{2})+ \\
&&-(1_{A}\otimes x_{2})B(g\otimes 1_{H})(1_{A}\otimes x_{1})+ \\
&&-B(x_{2}\otimes 1_{H})(1_{A}\otimes x_{1})
\end{eqnarray*}%
we get%
\begin{eqnarray*}
&&\text{right side}(1_{A}\otimes x_{2})\left( 1_{A}\otimes g\right)
B(1_{H}\otimes gx_{1})\left( 1_{A}\otimes g\right) +B(x_{2}\otimes
x_{1})+B(1_{H}\otimes gx_{1}x_{2}) \\
&=&(1_{A}\otimes x_{2})B(g\otimes 1_{H})(1_{A}\otimes x_{1}) \\
&&-\left( 1_{A}\otimes gx_{1}x_{2}\right) B(g\otimes 1_{H})\left(
1_{A}\otimes g\right) \\
&&+(1_{A}\otimes x_{2})\left( 1_{A}\otimes g\right) B(x_{1}\otimes
1_{H})\left( 1_{A}\otimes g\right) . \\
&&B(x_{2}\otimes 1_{H})(1_{A}\otimes x_{1}) \\
&&-(1_{A}\otimes gx_{1})B(x_{2}\otimes 1_{H})(1_{A}\otimes g) \\
&&-(1_{A}\otimes g)B(gx_{1}x_{2}\otimes 1_{H})(1_{A}\otimes g) \\
&&+(1_{A}\otimes x_{1})B(g\otimes 1_{H})(1_{A}\otimes x_{2}) \\
&&+(1_{A}\otimes gx_{1}x_{2})B(g\otimes 1_{H})\left( 1_{A}\otimes g\right) \\
&&+(1_{A}\otimes gx_{1})B(x_{2}\otimes 1_{H})\left( 1_{A}\otimes g\right) \\
&&+B(x_{1}\otimes 1_{H})(1_{A}\otimes x_{2}) \\
&&-(1_{A}\otimes gx_{2})B(x_{1}\otimes 1_{H})(1_{A}\otimes g) \\
&&+(1_{A}\otimes g)B(gx_{1}x_{2}\otimes 1_{H})(1_{A}\otimes g) \\
&&+(1_{A}\otimes g)B(g\otimes 1_{H})(1_{A}\otimes gx_{1}x_{2})+ \\
&&-(1_{A}\otimes x_{2})B(g\otimes 1_{H})(1_{A}\otimes x_{1})+ \\
&&-B(x_{2}\otimes 1_{H})(1_{A}\otimes x_{1})
\end{eqnarray*}%
and hence%
\begin{eqnarray*}
\text{clean } &&\text{right side}(1_{A}\otimes x_{2})\left( 1_{A}\otimes
g\right) B(1_{H}\otimes gx_{1})\left( 1_{A}\otimes g\right) +B(x_{2}\otimes
x_{1})+B(1_{H}\otimes gx_{1}x_{2}) \\
&&+(1_{A}\otimes x_{1})B(g\otimes 1_{H})(1_{A}\otimes x_{2}) \\
&&+B(x_{1}\otimes 1_{H})(1_{A}\otimes x_{2}) \\
&&+(1_{A}\otimes g)B(g\otimes 1_{H})(1_{A}\otimes gx_{1}x_{2})+
\end{eqnarray*}%
and we conclude.

\subsubsection{$1_{H}\otimes gx_{2}$}

\begin{equation*}
B(1_{H}\otimes gx_{2})(1_{A}\otimes x_{2})\overset{?}{=}(1_{A}\otimes
x_{2})\left( 1_{A}\otimes g\right) B(1_{H}\otimes gx_{2})\left( 1_{A}\otimes
g\right) +B(x_{2}\otimes x_{2})+B(1_{H}\otimes x_{2}gx_{2})
\end{equation*}%
i.e.%
\begin{equation*}
B(1_{H}\otimes gx_{2})(1_{A}\otimes x_{2})\overset{?}{=}(1_{A}\otimes
x_{2})\left( 1_{A}\otimes g\right) B(1_{H}\otimes gx_{2})\left( 1_{A}\otimes
g\right) +B(x_{2}\otimes x_{2}).
\end{equation*}%
Since%
\begin{eqnarray*}
&&B(1_{H}\otimes gx_{2}) \\
&&\overset{\left( \ref{form 1otgx2}\right) }{=}(1_{A}\otimes g)B(g\otimes
1_{H})(1_{A}\otimes gx_{2}) \\
&&+(1_{A}\otimes x_{2})B(g\otimes 1_{H}) \\
&&+B(x_{2}\otimes 1_{H})
\end{eqnarray*}%
we get%
\begin{eqnarray*}
&&\text{left side }B(1_{H}\otimes gx_{2})(1_{A}\otimes x_{2}) \\
&&=(1_{A}\otimes g)B(g\otimes 1_{H})(1_{A}\otimes gx_{2})(1_{A}\otimes x_{2})
\\
&&+(1_{A}\otimes x_{2})B(g\otimes 1_{H})(1_{A}\otimes x_{2}) \\
&&+B(x_{2}\otimes 1_{H})(1_{A}\otimes x_{2})
\end{eqnarray*}%
i.e.%
\begin{eqnarray*}
&&\text{left side }B(1_{H}\otimes gx_{2})(1_{A}\otimes x_{2}) \\
&&=(1_{A}\otimes x_{2})B(g\otimes 1_{H})(1_{A}\otimes x_{2}) \\
&&+B(x_{2}\otimes 1_{H})(1_{A}\otimes x_{2})
\end{eqnarray*}%
and%
\begin{eqnarray*}
&&(1_{A}\otimes x_{2})\left( 1_{A}\otimes g\right) B(1_{H}\otimes
gx_{2})\left( 1_{A}\otimes g\right) \\
&=&(1_{A}\otimes x_{2})\left( 1_{A}\otimes g\right) (1_{A}\otimes
g)B(g\otimes 1_{H})(1_{A}\otimes gx_{2})\left( 1_{A}\otimes g\right) \\
&&+(1_{A}\otimes x_{2})\left( 1_{A}\otimes g\right) (1_{A}\otimes
x_{2})B(g\otimes 1_{H})\left( 1_{A}\otimes g\right) \\
&&+(1_{A}\otimes x_{2})\left( 1_{A}\otimes g\right) B(x_{2}\otimes
1_{H})\left( 1_{A}\otimes g\right)
\end{eqnarray*}%
i.e.%
\begin{eqnarray*}
&&(1_{A}\otimes x_{2})\left( 1_{A}\otimes g\right) B(1_{H}\otimes
gx_{2})\left( 1_{A}\otimes g\right) \\
&=&(1_{A}\otimes x_{2})B(g\otimes 1_{H})(1_{A}\otimes x_{2}) \\
&&+(1_{A}\otimes gx_{2})B(x_{2}\otimes 1_{H})\left( 1_{A}\otimes g\right)
\end{eqnarray*}%
Since%
\begin{equation*}
B(x_{2}\otimes x_{2})\overset{\left( \ref{form x2otx2}\right) }{=}%
B(x_{2}\otimes 1_{H})(1_{A}\otimes x_{2})-(1_{A}\otimes
gx_{2})B(x_{2}\otimes 1_{H})(1_{A}\otimes g).
\end{equation*}%
we get%
\begin{eqnarray*}
&&\text{right side }(1_{A}\otimes x_{2})\left( 1_{A}\otimes g\right)
B(1_{H}\otimes gx_{2})\left( 1_{A}\otimes g\right) +B(x_{2}\otimes gx_{2}) \\
&=&(1_{A}\otimes x_{2})B(g\otimes 1_{H})(1_{A}\otimes x_{2}) \\
&&+(1_{A}\otimes gx_{2})B(x_{2}\otimes 1_{H})\left( 1_{A}\otimes g\right) \\
&&B(x_{2}\otimes 1_{H})(1_{A}\otimes x_{2}) \\
&&-(1_{A}\otimes gx_{2})B(x_{2}\otimes 1_{H})(1_{A}\otimes g)
\end{eqnarray*}%
i.e.%
\begin{eqnarray*}
&&\text{clean right side }(1_{A}\otimes x_{2})\left( 1_{A}\otimes g\right)
B(1_{H}\otimes gx_{2})\left( 1_{A}\otimes g\right) +B(x_{2}\otimes gx_{2}) \\
&=&(1_{A}\otimes x_{2})B(g\otimes 1_{H})(1_{A}\otimes x_{2}) \\
&&B(x_{2}\otimes 1_{H})(1_{A}\otimes x_{2})
\end{eqnarray*}%
and we conclude.

\subsubsection{$1_{H}\otimes gx_{1}x_{2}$}

\begin{eqnarray*}
&&B(1_{H}\otimes gx_{1}x_{2})(1_{A}\otimes x_{2})\overset{?}{=}(1_{A}\otimes
x_{2})\left( 1_{A}\otimes g\right) B(1_{H}\otimes gx_{1}x_{2})\left(
1_{A}\otimes g\right) \\
&&+B(x_{2}\otimes ggx_{1}x_{2})+B(1_{H}\otimes x_{2}gx_{1}x_{2})
\end{eqnarray*}%
i.e.%
\begin{equation*}
B(1_{H}\otimes gx_{1}x_{2})(1_{A}\otimes x_{2})\overset{?}{=}(1_{A}\otimes
x_{2})\left( 1_{A}\otimes g\right) B(1_{H}\otimes gx_{1}x_{2})\left(
1_{A}\otimes g\right) +B(x_{2}\otimes x_{1}x_{2}).
\end{equation*}%
Since%
\begin{eqnarray*}
&&B(1_{H}\otimes gx_{1}x_{2})\overset{\left( \ref{form 1otgx1x2}\right) }{=}
\\
&=&+(1_{A}\otimes x_{1})B(g\otimes 1_{H})(1_{A}\otimes x_{2}) \\
&&+(1_{A}\otimes gx_{1}x_{2})B(g\otimes 1_{H})\left( 1_{A}\otimes g\right) \\
&&+(1_{A}\otimes gx_{1})B(x_{2}\otimes 1_{H})\left( 1_{A}\otimes g\right) \\
&&+B(x_{1}\otimes 1_{H})(1_{A}\otimes x_{2}) \\
&&-(1_{A}\otimes gx_{2})B(x_{1}\otimes 1_{H})(1_{A}\otimes g) \\
&&+(1_{A}\otimes g)B(gx_{1}x_{2}\otimes 1_{H})(1_{A}\otimes g) \\
&&+(1_{A}\otimes g)B(g\otimes 1_{H})(1_{A}\otimes gx_{1}x_{2})+ \\
&&-(1_{A}\otimes x_{2})B(g\otimes 1_{H})(1_{A}\otimes x_{1})+ \\
&&-B(x_{2}\otimes 1_{H})(1_{A}\otimes x_{1})
\end{eqnarray*}%
we get%
\begin{eqnarray*}
\text{left side} &&B(1_{H}\otimes gx_{1}x_{2})(1_{A}\otimes x_{2})= \\
&=&+(1_{A}\otimes x_{1})B(g\otimes 1_{H})(1_{A}\otimes x_{2})(1_{A}\otimes
x_{2}) \\
&&+(1_{A}\otimes gx_{1}x_{2})B(g\otimes 1_{H})\left( 1_{A}\otimes g\right)
(1_{A}\otimes x_{2}) \\
&&+(1_{A}\otimes gx_{1})B(x_{2}\otimes 1_{H})\left( 1_{A}\otimes g\right)
(1_{A}\otimes x_{2}) \\
&&+B(x_{1}\otimes 1_{H})(1_{A}\otimes x_{2})(1_{A}\otimes x_{2}) \\
&&-(1_{A}\otimes gx_{2})B(x_{1}\otimes 1_{H})(1_{A}\otimes g)(1_{A}\otimes
x_{2}) \\
&&+(1_{A}\otimes g)B(gx_{1}x_{2}\otimes 1_{H})(1_{A}\otimes g)(1_{A}\otimes
x_{2}) \\
&&+(1_{A}\otimes g)B(g\otimes 1_{H})(1_{A}\otimes gx_{1}x_{2})(1_{A}\otimes
x_{2})+ \\
&&-(1_{A}\otimes x_{2})B(g\otimes 1_{H})(1_{A}\otimes x_{1})(1_{A}\otimes
x_{2})+ \\
&&-B(x_{2}\otimes 1_{H})(1_{A}\otimes x_{1})(1_{A}\otimes x_{2})
\end{eqnarray*}%
i.e.%
\begin{eqnarray*}
\text{left side } &&B(1_{H}\otimes gx_{1}x_{2})(1_{A}\otimes x_{2})= \\
&&-(1_{A}\otimes gx_{1}x_{2})B(g\otimes 1_{H})\left( 1_{A}\otimes
gx_{2}\right) \\
&&-(1_{A}\otimes gx_{1})B(x_{2}\otimes 1_{H})\left( 1_{A}\otimes
gx_{2}\right) \\
&&+(1_{A}\otimes gx_{2})B(x_{1}\otimes 1_{H})(1_{A}\otimes gx_{2}) \\
&&-(1_{A}\otimes g)B(gx_{1}x_{2}\otimes 1_{H})(1_{A}\otimes gx_{2}) \\
&&+(1_{A}\otimes x_{2})B(g\otimes 1_{H})(1_{A}\otimes x_{1}x_{2}) \\
&&+B(x_{2}\otimes 1_{H})(1_{A}\otimes x_{1}x_{2})
\end{eqnarray*}%
and
\begin{eqnarray*}
&&(1_{A}\otimes x_{2})\left( 1_{A}\otimes g\right) B(1_{H}\otimes
gx_{1}x_{2})\left( 1_{A}\otimes g\right) \\
&=&(1_{A}\otimes x_{2})\left( 1_{A}\otimes g\right) (1_{A}\otimes
x_{1})B(g\otimes 1_{H})(1_{A}\otimes x_{2})\left( 1_{A}\otimes g\right) + \\
&&+(1_{A}\otimes x_{2})\left( 1_{A}\otimes g\right) (1_{A}\otimes
gx_{1}x_{2})B(g\otimes 1_{H})\left( 1_{A}\otimes g\right) \left(
1_{A}\otimes g\right) + \\
&&+(1_{A}\otimes x_{2})\left( 1_{A}\otimes g\right) (1_{A}\otimes
gx_{1})B(x_{2}\otimes 1_{H})\left( 1_{A}\otimes g\right) \left( 1_{A}\otimes
g\right) + \\
&&+(1_{A}\otimes x_{2})\left( 1_{A}\otimes g\right) B(x_{1}\otimes
1_{H})(1_{A}\otimes x_{2})\left( 1_{A}\otimes g\right) + \\
&&-(1_{A}\otimes x_{2})\left( 1_{A}\otimes g\right) (1_{A}\otimes
gx_{2})B(x_{1}\otimes 1_{H})(1_{A}\otimes g)\left( 1_{A}\otimes g\right) + \\
&&+(1_{A}\otimes x_{2})\left( 1_{A}\otimes g\right) (1_{A}\otimes
g)B(gx_{1}x_{2}\otimes 1_{H})(1_{A}\otimes g)\left( 1_{A}\otimes g\right) \\
&&+(1_{A}\otimes x_{2})\left( 1_{A}\otimes g\right) (1_{A}\otimes
g)B(g\otimes 1_{H})(1_{A}\otimes gx_{1}x_{2})\left( 1_{A}\otimes g\right) +
\\
&&-(1_{A}\otimes x_{2})\left( 1_{A}\otimes g\right) (1_{A}\otimes
x_{2})B(g\otimes 1_{H})(1_{A}\otimes x_{1})\left( 1_{A}\otimes g\right) + \\
&&-(1_{A}\otimes x_{2})\left( 1_{A}\otimes g\right) B(x_{2}\otimes
1_{H})(1_{A}\otimes x_{1})\left( 1_{A}\otimes g\right)
\end{eqnarray*}%
i.e.%
\begin{eqnarray*}
&&(1_{A}\otimes x_{2})\left( 1_{A}\otimes g\right) B(1_{H}\otimes
gx_{1}x_{2})\left( 1_{A}\otimes g\right) \\
&=&-(1_{A}\otimes gx_{1}x_{2})B(g\otimes 1_{H})(1_{A}\otimes gx_{2})+ \\
&&-(1_{A}\otimes x_{1}x_{2})B(x_{2}\otimes 1_{H})+ \\
&&+(1_{A}\otimes gx_{2})B(x_{1}\otimes 1_{H})(1_{A}\otimes gx_{2})\left(
1_{A}\otimes \right) + \\
&&+(1_{A}\otimes x_{2})B(gx_{1}x_{2}\otimes 1_{H}) \\
&&+(1_{A}\otimes x_{2})B(g\otimes 1_{H})(1_{A}\otimes x_{1}x_{2})+ \\
&&-(1_{A}\otimes gx_{2})B(x_{2}\otimes 1_{H})(1_{A}\otimes gx_{1}).
\end{eqnarray*}%
Since%
\begin{eqnarray*}
&&B(x_{2}\otimes x_{1}x_{2})\overset{\left( \ref{form x2otx1x2}\right) }{=}%
B(x_{2}\otimes 1_{H})(1_{A}\otimes x_{1}x_{2}) \\
&&+(1_{A}\otimes gx_{2})B(x_{2}\otimes 1_{H})(1_{A}\otimes gx_{1}) \\
&&-(1_{A}\otimes gx_{1})B(x_{2}\otimes 1_{H})(1_{A}\otimes gx_{2}) \\
&&+(1_{A}\otimes x_{1}x_{2})B(x_{2}\otimes 1_{H}) \\
&&-(1_{A}\otimes g)B(gx_{1}x_{2}\otimes 1_{H})(1_{A}\otimes gx_{2}) \\
&&-(1_{A}\otimes x_{2})B(gx_{1}x_{2}\otimes 1_{H}).
\end{eqnarray*}%
we get%
\begin{eqnarray*}
\text{right side} &&(1_{A}\otimes x_{2})\left( 1_{A}\otimes g\right)
B(1_{H}\otimes gx_{1}x_{2})\left( 1_{A}\otimes g\right) +B(x_{2}\otimes
x_{1}x_{2}) \\
&=&-(1_{A}\otimes gx_{1}x_{2})B(g\otimes 1_{H})(1_{A}\otimes gx_{2})+ \\
&&-(1_{A}\otimes x_{1}x_{2})B(x_{2}\otimes 1_{H})+ \\
&&+(1_{A}\otimes gx_{2})B(x_{1}\otimes 1_{H})(1_{A}\otimes gx_{2})+ \\
&&+(1_{A}\otimes x_{2})B(gx_{1}x_{2}\otimes 1_{H}) \\
&&+(1_{A}\otimes x_{2})B(g\otimes 1_{H})(1_{A}\otimes x_{1}x_{2})+ \\
&&-(1_{A}\otimes gx_{2})B(x_{2}\otimes 1_{H})(1_{A}\otimes gx_{1}) \\
&&+B(x_{2}\otimes 1_{H})(1_{A}\otimes x_{1}x_{2}) \\
&&+B(x_{2}\otimes 1_{H})(1_{A}\otimes x_{1}x_{2}) \\
&&+(1_{A}\otimes gx_{2})B(x_{2}\otimes 1_{H})(1_{A}\otimes gx_{1}) \\
&&-(1_{A}\otimes gx_{1})B(x_{2}\otimes 1_{H})(1_{A}\otimes gx_{2}) \\
&&+(1_{A}\otimes x_{2})B(g\otimes 1_{H})(1_{A}\otimes x_{1}x_{2}) \\
&&+(1_{A}\otimes x_{1}x_{2})B(x_{2}\otimes 1_{H}) \\
&&-(1_{A}\otimes g)B(gx_{1}x_{2}\otimes 1_{H})(1_{A}\otimes gx_{2}) \\
&&-(1_{A}\otimes x_{2})B(gx_{1}x_{2}\otimes 1_{H})
\end{eqnarray*}%
so that we get%
\begin{eqnarray*}
\text{clean right side} &&(1_{A}\otimes x_{2})\left( 1_{A}\otimes g\right)
B(1_{H}\otimes gx_{1}x_{2})\left( 1_{A}\otimes g\right) +B(x_{2}\otimes
x_{1}x_{2}) \\
&=&-(1_{A}\otimes gx_{1}x_{2})B(g\otimes 1_{H})(1_{A}\otimes gx_{2})+ \\
&&+(1_{A}\otimes gx_{2})B(x_{1}\otimes 1_{H})(1_{A}\otimes gx_{2})+ \\
&&+(1_{A}\otimes x_{2})B(g\otimes 1_{H})(1_{A}\otimes x_{1}x_{2})+ \\
&&+B(x_{2}\otimes 1_{H})(1_{A}\otimes x_{1}x_{2}) \\
&&-(1_{A}\otimes gx_{1})B(x_{2}\otimes 1_{H})(1_{A}\otimes gx_{2}) \\
&&-(1_{A}\otimes g)B(gx_{1}x_{2}\otimes 1_{H})(1_{A}\otimes gx_{2})
\end{eqnarray*}%
and we conclude.

\subsubsection{$x_{1}\otimes x_{1}$}

\begin{equation*}
B(x_{1}\otimes x_{1})(1_{A}\otimes x_{2})\overset{?}{=}(1_{A}\otimes
x_{2})\left( 1_{A}\otimes g\right) B(x_{1}\otimes x_{1})\left( 1_{A}\otimes
g\right) +B(x_{2}x_{1}\otimes gx_{1})+B(x_{1}\otimes x_{2}x_{1})
\end{equation*}%
i.e%
\begin{equation*}
B(x_{1}\otimes x_{1})(1_{A}\otimes x_{2})\overset{?}{=}(1_{A}\otimes
x_{2})\left( 1_{A}\otimes g\right) B(x_{1}\otimes x_{1})\left( 1_{A}\otimes
g\right) -B(x_{1}x_{2}\otimes gx_{1})-B(x_{1}\otimes x_{1}x_{2})
\end{equation*}%
Since%
\begin{equation*}
B(x_{1}\otimes x_{1})\overset{\left( \ref{form x1otx1}\right) }{=}%
B(x_{1}\otimes 1_{H})(1_{A}\otimes x_{1})-(1_{A}\otimes
gx_{1})B(x_{1}\otimes 1_{H})(1_{A}\otimes g)
\end{equation*}%
we get%
\begin{eqnarray*}
\text{left side }B(x_{1}\otimes x_{1})(1_{A}\otimes x_{2}) &=&B(x_{1}\otimes
1_{H})(1_{A}\otimes x_{1})(1_{A}\otimes x_{2}) \\
&&-(1_{A}\otimes gx_{1})B(x_{1}\otimes 1_{H})(1_{A}\otimes g)(1_{A}\otimes
x_{2})
\end{eqnarray*}%
i.e.%
\begin{equation*}
\text{left side }B(x_{1}\otimes x_{1})(1_{A}\otimes x_{2})=-B(x_{1}\otimes
1_{H})(1_{A}\otimes x_{1}x_{2})+(1_{A}\otimes gx_{1})B(x_{1}\otimes
1_{H})(1_{A}\otimes gx_{2})
\end{equation*}%
and%
\begin{eqnarray*}
&&(1_{A}\otimes x_{2})\left( 1_{A}\otimes g\right) B(x_{1}\otimes
x_{1})\left( 1_{A}\otimes g\right) \\
&=&(1_{A}\otimes x_{2})\left( 1_{A}\otimes g\right) B(x_{1}\otimes
1_{H})(1_{A}\otimes x_{1})\left( 1_{A}\otimes g\right) \\
&&-(1_{A}\otimes x_{2})\left( 1_{A}\otimes g\right) (1_{A}\otimes
gx_{1})B(x_{1}\otimes 1_{H})(1_{A}\otimes g)\left( 1_{A}\otimes g\right)
\end{eqnarray*}%
i.e.%
\begin{eqnarray*}
&&(1_{A}\otimes x_{2})\left( 1_{A}\otimes g\right) B(x_{1}\otimes
x_{1})\left( 1_{A}\otimes g\right) \\
&=&(1_{A}\otimes gx_{2})B(x_{1}\otimes 1_{H})(1_{A}\otimes gx_{1}) \\
&&+(1_{A}\otimes x_{1}x_{2})B(x_{1}\otimes 1_{H}).
\end{eqnarray*}%
Since%
\begin{eqnarray*}
&&-B(x_{1}x_{2}\otimes gx_{1}) \\
&&\overset{\left( \ref{form x1x2otgx1}\right) }{=}-(1_{A}\otimes
g)B(gx_{1}x_{2}\otimes 1_{H})(1_{A}\otimes gx_{1}) \\
&&-(1_{A}\otimes x_{1})B(gx_{1}x_{2}\otimes 1_{H})
\end{eqnarray*}%
and%
\begin{eqnarray*}
- &&B(x_{1}\otimes x_{1}x_{2})\overset{\left( \ref{form x1otx1x2}\right) }{=}%
-B(x_{1}\otimes 1_{H})(1_{A}\otimes x_{1}x_{2}) \\
&&-(1_{A}\otimes gx_{2})B(x_{1}\otimes 1_{H})(1_{A}\otimes gx_{1}) \\
&&+(1_{A}\otimes g)B(gx_{1}x_{2}\otimes 1_{H})(1_{A}\otimes gx_{1}) \\
&&+(1_{A}\otimes gx_{1})B(x_{1}\otimes 1_{H})(1_{A}\otimes gx_{2}) \\
&&-(1_{A}\otimes x_{1}x_{2})B(x_{1}\otimes 1_{H}) \\
&&+(1_{A}\otimes x_{1})B(gx_{1}x_{2}\otimes 1_{H})
\end{eqnarray*}%
we get%
\begin{gather*}
\text{right side }(1_{A}\otimes x_{2})\left( 1_{A}\otimes g\right)
B(x_{1}\otimes x_{1})\left( 1_{A}\otimes g\right) -B(x_{1}x_{2}\otimes
gx_{1})-B(x_{1}\otimes x_{1}x_{2}) \\
=(1_{A}\otimes gx_{2})B(x_{1}\otimes 1_{H})(1_{A}\otimes gx_{1}) \\
(1_{A}\otimes x_{1}x_{2})B(x_{1}\otimes 1_{H}) \\
-(1_{A}\otimes g)B(gx_{1}x_{2}\otimes 1_{H})(1_{A}\otimes gx_{1}) \\
-(1_{A}\otimes x_{1})B(gx_{1}x_{2}\otimes 1_{H}) \\
-B(x_{1}\otimes 1_{H})(1_{A}\otimes x_{1}x_{2}) \\
-(1_{A}\otimes gx_{2})B(x_{1}\otimes 1_{H})(1_{A}\otimes gx_{1}) \\
+(1_{A}\otimes g)B(gx_{1}x_{2}\otimes 1_{H})(1_{A}\otimes gx_{1}) \\
+(1_{A}\otimes gx_{1})B(x_{1}\otimes 1_{H})(1_{A}\otimes gx_{2}) \\
-(1_{A}\otimes x_{1}x_{2})B(x_{1}\otimes 1_{H}) \\
+(1_{A}\otimes x_{1})B(gx_{1}x_{2}\otimes 1_{H})
\end{gather*}%
\begin{gather*}
\text{clean right side }(1_{A}\otimes x_{2})\left( 1_{A}\otimes g\right)
B(x_{1}\otimes x_{1})\left( 1_{A}\otimes g\right) -B(x_{1}x_{2}\otimes
gx_{1})-B(x_{1}\otimes x_{1}x_{2}) \\
=-B(x_{1}\otimes 1_{H})(1_{A}\otimes x_{1}x_{2}) \\
+(1_{A}\otimes gx_{1})B(x_{1}\otimes 1_{H})(1_{A}\otimes gx_{2})
\end{gather*}%
and we conclude.

\subsubsection{$x_{1}\otimes x_{2}$}

\begin{equation*}
B(x_{1}\otimes x_{2})(1_{A}\otimes x_{2})\overset{?}{=}(1_{A}\otimes
x_{2})\left( 1_{A}\otimes g\right) B(x_{1}\otimes x_{2})\left( 1_{A}\otimes
g\right) +B(x_{2}x_{1}\otimes gx_{2})+B(x_{1}\otimes x_{2}x_{2})
\end{equation*}%
i.e.%
\begin{equation*}
B(x_{1}\otimes x_{2})(1_{A}\otimes x_{2})\overset{?}{=}(1_{A}\otimes
x_{2})\left( 1_{A}\otimes g\right) B(x_{1}\otimes x_{2})\left( 1_{A}\otimes
g\right) -B(x_{1}x_{2}\otimes gx_{2})
\end{equation*}%
Since%
\begin{eqnarray*}
&&B(x_{1}\otimes x_{2})\overset{\left( \ref{form x1otx2}\right) }{=}%
B(x_{1}\otimes 1_{H})(1_{A}\otimes x_{2})-(1_{A}\otimes
gx_{2})B(x_{1}\otimes 1_{H})(1_{A}\otimes g) \\
&&+(1_{A}\otimes g)B(gx_{1}x_{2}\otimes 1_{H})(1_{A}\otimes g)
\end{eqnarray*}%
we get%
\begin{eqnarray*}
&&\text{left side }B(x_{1}\otimes x_{2})(1_{A}\otimes x_{2}) \\
&=&B(x_{1}\otimes 1_{H})(1_{A}\otimes x_{2})(1_{A}\otimes x_{2}) \\
&&-(1_{A}\otimes gx_{2})B(x_{1}\otimes 1_{H})(1_{A}\otimes g)(1_{A}\otimes
x_{2}) \\
&&+(1_{A}\otimes g)B(gx_{1}x_{2}\otimes 1_{H})(1_{A}\otimes g)(1_{A}\otimes
x_{2})
\end{eqnarray*}%
i.e.%
\begin{eqnarray*}
&&\text{left side }B(x_{1}\otimes x_{2})(1_{A}\otimes x_{2}) \\
&=&(1_{A}\otimes gx_{2})B(x_{1}\otimes 1_{H})(1_{A}\otimes gx_{2}) \\
&&-(1_{A}\otimes g)B(gx_{1}x_{2}\otimes 1_{H})(1_{A}\otimes gx_{2})
\end{eqnarray*}%
and%
\begin{eqnarray*}
(1_{A}\otimes x_{2})\left( 1_{A}\otimes g\right) B(x_{1}\otimes x_{2})\left(
1_{A}\otimes g\right) &=&(1_{A}\otimes x_{2})\left( 1_{A}\otimes g\right)
B(x_{1}\otimes 1_{H})(1_{A}\otimes x_{2})\left( 1_{A}\otimes g\right) \\
&&-(1_{A}\otimes x_{2})\left( 1_{A}\otimes g\right) (1_{A}\otimes
gx_{2})B(x_{1}\otimes 1_{H})(1_{A}\otimes g)\left( 1_{A}\otimes g\right) \\
&&+(1_{A}\otimes x_{2})\left( 1_{A}\otimes g\right) (1_{A}\otimes
g)B(gx_{1}x_{2}\otimes 1_{H})(1_{A}\otimes g)\left( 1_{A}\otimes g\right)
\end{eqnarray*}%
i.e.%
\begin{eqnarray*}
(1_{A}\otimes x_{2})\left( 1_{A}\otimes g\right) B(x_{1}\otimes x_{2})\left(
1_{A}\otimes g\right) &=&(1_{A}\otimes gx_{2})B(x_{1}\otimes
1_{H})(1_{A}\otimes gx_{2}) \\
&&+(1_{A}\otimes x_{2})B(gx_{1}x_{2}\otimes 1_{H})
\end{eqnarray*}%
Since%
\begin{equation*}
-B(x_{1}x_{2}\otimes gx_{2})\overset{\left( \ref{form x1x2otgx2}\right) }{=}%
-(1_{A}\otimes g)B(gx_{1}x_{2}\otimes 1_{H})(1_{A}\otimes
gx_{2})-(1_{A}\otimes x_{2})B(gx_{1}x_{2}\otimes 1_{H})
\end{equation*}%
we get%
\begin{eqnarray*}
&&\text{right side }(1_{A}\otimes x_{2})\left( 1_{A}\otimes g\right)
B(x_{1}\otimes x_{2})\left( 1_{A}\otimes g\right) -B(x_{1}x_{2}\otimes
gx_{2}) \\
&=&(1_{A}\otimes gx_{2})B(x_{1}\otimes 1_{H})(1_{A}\otimes gx_{2}) \\
&&+(1_{A}\otimes x_{2})B(gx_{1}x_{2}\otimes 1_{H}) \\
&&-(1_{A}\otimes g)B(gx_{1}x_{2}\otimes 1_{H})(1_{A}\otimes gx_{2}) \\
&&-(1_{A}\otimes x_{2})B(gx_{1}x_{2}\otimes 1_{H})
\end{eqnarray*}%
i.e.%
\begin{eqnarray*}
&&\text{clean right side }(1_{A}\otimes x_{2})\left( 1_{A}\otimes g\right)
B(x_{1}\otimes x_{2})\left( 1_{A}\otimes g\right) -B(x_{1}x_{2}\otimes
gx_{2}) \\
&=&(1_{A}\otimes gx_{2})B(x_{1}\otimes 1_{H})(1_{A}\otimes gx_{2}) \\
&&-(1_{A}\otimes g)B(gx_{1}x_{2}\otimes 1_{H})(1_{A}\otimes gx_{2})
\end{eqnarray*}%
and we conclude.

\subsubsection{$x_{1}\otimes x_{1}x_{2}$}

\begin{equation*}
B(x_{1}\otimes x_{1}x_{2})(1_{A}\otimes x_{2})\overset{?}{=}(1_{A}\otimes
x_{2})\left( 1_{A}\otimes g\right) B(x_{1}\otimes x_{1}x_{2})\left(
1_{A}\otimes g\right) +B(x_{2}x_{1}\otimes gx_{1}x_{2})+B(x_{1}\otimes
x_{2}x_{1}x_{2})
\end{equation*}%
i.e.%
\begin{equation*}
B(x_{1}\otimes x_{1}x_{2})(1_{A}\otimes x_{2})\overset{?}{=}(1_{A}\otimes
x_{2})\left( 1_{A}\otimes g\right) B(x_{1}\otimes x_{1}x_{2})\left(
1_{A}\otimes g\right) -B(x_{1}x_{2}\otimes gx_{1}x_{2}).
\end{equation*}%
Since
\begin{eqnarray*}
&&B(x_{1}\otimes x_{1}x_{2})\overset{\left( \ref{form x1otx1x2}\right) }{=}%
B(x_{1}\otimes 1_{H})(1_{A}\otimes x_{1}x_{2}) \\
&&+(1_{A}\otimes gx_{2})B(x_{1}\otimes 1_{H})(1_{A}\otimes gx_{1}) \\
&&-(1_{A}\otimes g)B(gx_{1}x_{2}\otimes 1_{H})(1_{A}\otimes gx_{1}) \\
&&-(1_{A}\otimes gx_{1})B(x_{1}\otimes 1_{H})(1_{A}\otimes gx_{2}) \\
&&+(1_{A}\otimes x_{1}x_{2})B(x_{1}\otimes 1_{H}) \\
&&-(1_{A}\otimes x_{1})B(gx_{1}x_{2}\otimes 1_{H})
\end{eqnarray*}%
we get%
\begin{eqnarray*}
&&\text{left side }B(x_{1}\otimes x_{1}x_{2})(1_{A}\otimes
x_{2})=B(x_{1}\otimes 1_{H})(1_{A}\otimes x_{1}x_{2})(1_{A}\otimes x_{2}) \\
&&+(1_{A}\otimes gx_{2})B(x_{1}\otimes 1_{H})(1_{A}\otimes
gx_{1})(1_{A}\otimes x_{2}) \\
&&-(1_{A}\otimes g)B(gx_{1}x_{2}\otimes 1_{H})(1_{A}\otimes
gx_{1})(1_{A}\otimes x_{2}) \\
&&-(1_{A}\otimes gx_{1})B(x_{1}\otimes 1_{H})(1_{A}\otimes
gx_{2})(1_{A}\otimes x_{2}) \\
&&+(1_{A}\otimes x_{1}x_{2})B(x_{1}\otimes 1_{H})(1_{A}\otimes x_{2}) \\
&&-(1_{A}\otimes x_{1})B(gx_{1}x_{2}\otimes 1_{H})(1_{A}\otimes x_{2})
\end{eqnarray*}%
i.e.%
\begin{eqnarray*}
&&\text{left side }B(x_{1}\otimes x_{1}x_{2})(1_{A}\otimes x_{2})= \\
&&+(1_{A}\otimes gx_{2})B(x_{1}\otimes 1_{H})(1_{A}\otimes gx_{1}x_{2}) \\
&&-(1_{A}\otimes g)B(gx_{1}x_{2}\otimes 1_{H})(1_{A}\otimes gx_{1}x_{2}) \\
&&+(1_{A}\otimes x_{1}x_{2})B(x_{1}\otimes 1_{H})(1_{A}\otimes x_{2}) \\
&&-(1_{A}\otimes x_{1})B(gx_{1}x_{2}\otimes 1_{H})(1_{A}\otimes x_{2})
\end{eqnarray*}%
and

\begin{eqnarray*}
&&(1_{A}\otimes x_{2})\left( 1_{A}\otimes g\right) B(x_{1}\otimes
x_{1}x_{2})\left( 1_{A}\otimes g\right) \\
&=&(1_{A}\otimes gx_{2})B(x_{1}\otimes 1_{H})(1_{A}\otimes gx_{1}x_{2}) \\
&&-(1_{A}\otimes x_{2})B(gx_{1}x_{2}\otimes 1_{H})(1_{A}\otimes x_{1}) \\
&&+(1_{A}\otimes x_{1}x_{2})B(x_{1}\otimes 1_{H})(1_{A}\otimes x_{2}) \\
&&+(1_{A}\otimes gx_{1}x_{2})B(gx_{1}x_{2}\otimes 1_{H})\left( 1_{A}\otimes
g\right) .
\end{eqnarray*}%
Since%
\begin{eqnarray*}
- &&B(x_{1}x_{2}\otimes gx_{1}x_{2})\overset{\left( \ref{form x1x2otgx1gx2}%
\right) }{=}-(1_{A}\otimes g)B(gx_{1}x_{2}\otimes 1_{H})(1_{A}\otimes
gx_{1}x_{2}) \\
&&+(1_{A}\otimes x_{2})B(gx_{1}x_{2}\otimes 1_{H})(1_{A}\otimes x_{1}) \\
&&-(1_{A}\otimes x_{1})B(gx_{1}x_{2}\otimes 1_{H})(1_{A}\otimes x_{2}) \\
&&-(1_{A}\otimes gx_{1}x_{2})B(gx_{1}x_{2}\otimes 1_{H})(1_{A}\otimes g)
\end{eqnarray*}%
we get%
\begin{eqnarray*}
\text{right side} &&(1_{A}\otimes x_{2})\left( 1_{A}\otimes g\right)
B(x_{1}\otimes x_{1}x_{2})\left( 1_{A}\otimes g\right) -B(x_{1}x_{2}\otimes
gx_{1}x_{2}) \\
&=&(1_{A}\otimes gx_{2})B(x_{1}\otimes 1_{H})(1_{A}\otimes x_{1}x_{2})\left(
1_{A}\otimes g\right) \\
&&-(1_{A}\otimes x_{2})B(gx_{1}x_{2}\otimes 1_{H})(1_{A}\otimes x_{1}) \\
&&+(1_{A}\otimes x_{1}x_{2})B(x_{1}\otimes 1_{H})(1_{A}\otimes x_{2}) \\
&&+(1_{A}\otimes gx_{1}x_{2})B(gx_{1}x_{2}\otimes 1_{H})\left( 1_{A}\otimes
g\right) \\
&&-(1_{A}\otimes g)B(gx_{1}x_{2}\otimes 1_{H})(1_{A}\otimes gx_{1}x_{2}) \\
&&+(1_{A}\otimes x_{2})B(gx_{1}x_{2}\otimes 1_{H})(1_{A}\otimes x_{1}) \\
&&-(1_{A}\otimes x_{1})B(gx_{1}x_{2}\otimes 1_{H})(1_{A}\otimes x_{2}) \\
&&-(1_{A}\otimes gx_{1}x_{2})B(gx_{1}x_{2}\otimes 1_{H})(1_{A}\otimes g)
\end{eqnarray*}%
\begin{eqnarray*}
\text{clean right side} &&(1_{A}\otimes x_{2})\left( 1_{A}\otimes g\right)
B(x_{1}\otimes x_{1}x_{2})\left( 1_{A}\otimes g\right) -B(x_{1}x_{2}\otimes
gx_{1}x_{2}) \\
&=&(1_{A}\otimes gx_{2})B(x_{1}\otimes 1_{H})(1_{A}\otimes gx_{1}x_{2}) \\
&&+(1_{A}\otimes x_{1}x_{2})B(x_{1}\otimes 1_{H})(1_{A}\otimes x_{2}) \\
&&-(1_{A}\otimes g)B(gx_{1}x_{2}\otimes 1_{H})(1_{A}\otimes gx_{1}x_{2}) \\
&&-(1_{A}\otimes x_{1})B(gx_{1}x_{2}\otimes 1_{H})(1_{A}\otimes x_{2})
\end{eqnarray*}%
and we conclude.

\subsubsection{$x_{1}\otimes gx_{1}$}

\begin{eqnarray*}
&&B(x_{1}\otimes gx_{1})(1_{A}\otimes x_{2})\overset{?}{=}(1_{A}\otimes
x_{2})\left( 1_{A}\otimes g\right) B(x_{1}\otimes gx_{1})\left( 1_{A}\otimes
g\right) \\
&&+B(x_{2}x_{1}\otimes ggx_{1})+B(x_{1}\otimes x_{2}gx_{1})
\end{eqnarray*}%
i.e.%
\begin{eqnarray*}
&&B(x_{1}\otimes gx_{1})(1_{A}\otimes x_{2})\overset{?}{=}(1_{A}\otimes
x_{2})\left( 1_{A}\otimes g\right) B(x_{1}\otimes gx_{1})\left( 1_{A}\otimes
g\right) \\
&&-B(x_{1}x_{2}\otimes x_{1})+B(x_{1}\otimes gx_{1}x_{2})
\end{eqnarray*}%
Since%
\begin{equation*}
B(x_{1}\otimes gx_{1})\overset{\left( \ref{form x1otgx1}\right) }{=}%
(1_{A}\otimes g)B(gx_{1}\otimes 1_{H})(1_{A}\otimes gx_{1})+(1_{A}\otimes
x_{1})B(gx_{1}\otimes 1_{H})
\end{equation*}%
we get%
\begin{eqnarray*}
&&\text{left side }B(x_{1}\otimes gx_{1})(1_{A}\otimes x_{2}) \\
&=&(1_{A}\otimes g)B(gx_{1}\otimes 1_{H})(1_{A}\otimes
gx_{1}x_{2})+(1_{A}\otimes x_{1})B(gx_{1}\otimes 1_{H})(1_{A}\otimes x_{2})
\end{eqnarray*}%
and%
\begin{eqnarray*}
&&(1_{A}\otimes x_{2})\left( 1_{A}\otimes g\right) B(x_{1}\otimes
gx_{1})\left( 1_{A}\otimes g\right) \\
&=&(1_{A}\otimes x_{2})\left( 1_{A}\otimes g\right) (1_{A}\otimes
g)B(gx_{1}\otimes 1_{H})(1_{A}\otimes gx_{1})\left( 1_{A}\otimes g\right) \\
&&+(1_{A}\otimes x_{2})\left( 1_{A}\otimes g\right) (1_{A}\otimes
x_{1})B(gx_{1}\otimes 1_{H})\left( 1_{A}\otimes g\right)
\end{eqnarray*}%
i.e.%
\begin{eqnarray*}
&&(1_{A}\otimes gx_{2})B(x_{1}\otimes gx_{1})\left( 1_{A}\otimes g\right) \\
&=&(1_{A}\otimes x_{2})B(gx_{1}\otimes 1_{H})(1_{A}\otimes x_{1}) \\
&&-(1_{A}\otimes gx_{1}x_{2})B(gx_{1}\otimes 1_{H})\left( 1_{A}\otimes
g\right) .
\end{eqnarray*}%
Since%
\begin{equation*}
B(x_{1}x_{2}\otimes x_{1})\overset{\left( \ref{form x1x2otx1}\right) }{=}%
B(x_{1}x_{2}\otimes 1_{H})(1_{A}\otimes x_{1})-(1_{A}\otimes
gx_{1})B(x_{1}x_{2}\otimes 1_{H})(1_{A}\otimes g)
\end{equation*}%
and%
\begin{eqnarray*}
&&B(x_{1}\otimes gx_{1}x_{2})\overset{\left( \ref{form x1otgx1x2}\right) }{=}%
+(1_{A}\otimes g)B(gx_{1}\otimes 1_{H})(1_{A}\otimes gx_{1}x_{2}) \\
&&-(1_{A}\otimes x_{2})B(gx_{1}\otimes 1_{H})(1_{A}\otimes x_{1}) \\
&&+B(x_{1}x_{2}\otimes 1_{H})(1_{A}\otimes x_{1}) \\
&&+(1_{A}\otimes x_{1})B(gx_{1}\otimes 1_{H})(1_{A}\otimes x_{2}) \\
&&+(1_{A}\otimes gx_{1}x_{2})B(gx_{1}\otimes 1_{H})(1_{A}\otimes g) \\
&&-(1_{A}\otimes gx_{1})B(x_{1}x_{2}\otimes 1_{H})(1_{A}\otimes g)
\end{eqnarray*}%
we get%
\begin{eqnarray*}
\text{right side } &&(1_{A}\otimes gx_{2})B(x_{1}\otimes gx_{1})\left(
1_{A}\otimes g\right) -B(x_{1}x_{2}\otimes x_{1})+B(x_{1}\otimes
gx_{1}x_{2})= \\
&&(1_{A}\otimes x_{2})B(gx_{1}\otimes 1_{H})(1_{A}\otimes x_{1})+ \\
&&-(1_{A}\otimes gx_{1}x_{2})B(gx_{1}\otimes 1_{H})\left( 1_{A}\otimes
g\right) + \\
&&-B(x_{1}x_{2}\otimes 1_{H})(1_{A}\otimes x_{1})+ \\
&&+(1_{A}\otimes gx_{1})B(x_{1}x_{2}\otimes 1_{H})(1_{A}\otimes g)+ \\
&&+(1_{A}\otimes g)B(gx_{1}\otimes 1_{H})(1_{A}\otimes gx_{1}x_{2})+ \\
&&-(1_{A}\otimes x_{2})B(gx_{1}\otimes 1_{H})(1_{A}\otimes x_{1})+ \\
&&+B(x_{1}x_{2}\otimes 1_{H})(1_{A}\otimes x_{1})+ \\
&&+(1_{A}\otimes x_{1})B(gx_{1}\otimes 1_{H})(1_{A}\otimes x_{2})+ \\
&&+(1_{A}\otimes gx_{1}x_{2})B(gx_{1}\otimes 1_{H})(1_{A}\otimes g)+ \\
&&-(1_{A}\otimes gx_{1})B(x_{1}x_{2}\otimes 1_{H})(1_{A}\otimes g)
\end{eqnarray*}%
i.e.%
\begin{eqnarray*}
&&\text{clean right side} \\
\text{ } &&(1_{A}\otimes gx_{2})B(x_{1}\otimes gx_{1})\left( 1_{A}\otimes
g\right) +B(x_{1}x_{2}\otimes x_{1})+B(x_{1}\otimes gx_{1}x_{2})= \\
&&(1_{A}\otimes g)B(gx_{1}\otimes 1_{H})(1_{A}\otimes gx_{1}x_{2}) \\
&&+(1_{A}\otimes x_{1})B(gx_{1}\otimes 1_{H})(1_{A}\otimes x_{2})
\end{eqnarray*}%
and we conclude.

\subsubsection{$x_{1}\otimes gx_{2}$}

\begin{eqnarray*}
&&B(x_{1}\otimes gx_{2})(1_{A}\otimes x_{2})\overset{?}{=}(1_{A}\otimes
x_{2})\left( 1_{A}\otimes g\right) B(x_{1}\otimes gx_{2})\left( 1_{A}\otimes
g\right) \\
&&+B(x_{2}x_{1}\otimes ggx_{2})+B(x_{1}\otimes x_{2}gx_{2})
\end{eqnarray*}%
i.e.%
\begin{equation*}
B(x_{1}\otimes gx_{2})(1_{A}\otimes x_{2})\overset{?}{=}(1_{A}\otimes
x_{2})\left( 1_{A}\otimes g\right) B(x_{1}\otimes gx_{2})\left( 1_{A}\otimes
g\right) -B(x_{1}x_{2}\otimes x_{2}).
\end{equation*}%
Since%
\begin{equation*}
B(x_{1}\otimes gx_{2})\overset{\left( \ref{form x1otgx2}\right) }{=}%
(1_{A}\otimes g)B(gx_{1}\otimes 1_{H})(1_{A}\otimes gx_{2})+(1_{A}\otimes
x_{2})B(gx_{1}\otimes 1_{H})-B(x_{1}x_{2}\otimes 1_{H})
\end{equation*}%
we get%
\begin{eqnarray*}
&&\text{left side }B(x_{1}\otimes gx_{2})(1_{A}\otimes x_{2}) \\
&=&(1_{A}\otimes g)B(gx_{1}\otimes 1_{H})(1_{A}\otimes gx_{2})(1_{A}\otimes
x_{2}) \\
&&+(1_{A}\otimes x_{2})B(gx_{1}\otimes 1_{H})(1_{A}\otimes
x_{2})-B(x_{1}x_{2}\otimes 1_{H})(1_{A}\otimes x_{2})
\end{eqnarray*}%
i.e.%
\begin{eqnarray*}
&&\text{left side }B(x_{1}\otimes gx_{2})(1_{A}\otimes x_{2}) \\
&=&(1_{A}\otimes x_{2})B(gx_{1}\otimes 1_{H})(1_{A}\otimes
x_{2})-B(x_{1}x_{2}\otimes 1_{H})(1_{A}\otimes x_{2})
\end{eqnarray*}%
and%
\begin{eqnarray*}
&&(1_{A}\otimes x_{2})\left( 1_{A}\otimes g\right) B(x_{1}\otimes
gx_{2})\left( 1_{A}\otimes g\right) \\
&=&(1_{A}\otimes x_{2})\left( 1_{A}\otimes g\right) (1_{A}\otimes
g)B(gx_{1}\otimes 1_{H})(1_{A}\otimes gx_{2})\left( 1_{A}\otimes g\right) \\
&&+(1_{A}\otimes x_{2})\left( 1_{A}\otimes g\right) (1_{A}\otimes
x_{2})B(gx_{1}\otimes 1_{H})\left( 1_{A}\otimes g\right) \\
&&-(1_{A}\otimes x_{2})\left( 1_{A}\otimes g\right) B(x_{1}x_{2}\otimes
1_{H})\left( 1_{A}\otimes g\right)
\end{eqnarray*}%
i.e.%
\begin{eqnarray*}
&&(1_{A}\otimes x_{2})\left( 1_{A}\otimes g\right) B(x_{1}\otimes
gx_{2})\left( 1_{A}\otimes g\right) \\
&=&(1_{A}\otimes x_{2})B(gx_{1}\otimes 1_{H})(1_{A}\otimes x_{2}) \\
&&-(1_{A}\otimes gx_{2})B(x_{1}x_{2}\otimes 1_{H})\left( 1_{A}\otimes
g\right) .
\end{eqnarray*}%
Since%
\begin{equation*}
-B(x_{1}x_{2}\otimes x_{2})\overset{\left( \ref{x1x2otx2}\right) }{=}%
-B(x_{1}x_{2}\otimes 1_{H})(1_{A}\otimes x_{2})+(1_{A}\otimes
gx_{2})B(x_{1}x_{2}\otimes 1_{H})(1_{A}\otimes g)
\end{equation*}%
we get%
\begin{eqnarray*}
&&\text{right side }(1_{A}\otimes x_{2})\left( 1_{A}\otimes g\right)
B(x_{1}\otimes gx_{2})\left( 1_{A}\otimes g\right) -B(x_{1}x_{2}\otimes
x_{2}) \\
&=&(1_{A}\otimes x_{2})B(gx_{1}\otimes 1_{H})(1_{A}\otimes x_{2}) \\
&&-(1_{A}\otimes gx_{2})B(x_{1}x_{2}\otimes 1_{H})\left( 1_{A}\otimes
g\right) + \\
&&-B(x_{1}x_{2}\otimes 1_{H})(1_{A}\otimes x_{2})+ \\
&&+(1_{A}\otimes gx_{2})B(x_{1}x_{2}\otimes 1_{H})(1_{A}\otimes g)
\end{eqnarray*}%
i.e.%
\begin{eqnarray*}
&&\text{right side }(1_{A}\otimes x_{2})\left( 1_{A}\otimes g\right)
B(x_{1}\otimes gx_{2})\left( 1_{A}\otimes g\right) -B(x_{1}x_{2}\otimes
x_{2}) \\
&=&(1_{A}\otimes x_{2})B(gx_{1}\otimes 1_{H})(1_{A}\otimes
x_{2})-B(x_{1}x_{2}\otimes 1_{H})(1_{A}\otimes x_{2})
\end{eqnarray*}%
and we conclude.

\subsubsection{$x_{1}\otimes gx_{1}x_{2}$}

\begin{eqnarray*}
&&B(x_{1}\otimes gx_{1}x_{2})(1_{A}\otimes x_{2}) \\
&&\overset{?}{=}(1_{A}\otimes x_{2})\left( 1_{A}\otimes g\right)
B(x_{1}\otimes gx_{1}x_{2})\left( 1_{A}\otimes g\right) +B(x_{2}x_{1}\otimes
ggx_{1}x_{2}) \\
&&+B(x_{1}\otimes x_{2}gx_{1}x_{2})
\end{eqnarray*}%
i.e.%
\begin{equation*}
B(x_{1}\otimes gx_{1}x_{2})(1_{A}\otimes x_{2})\overset{?}{=}(1_{A}\otimes
x_{2})\left( 1_{A}\otimes g\right) B(x_{1}\otimes gx_{1}x_{2})\left(
1_{A}\otimes g\right) -B(x_{1}x_{2}\otimes x_{1}x_{2}).
\end{equation*}%
Since%
\begin{eqnarray*}
&&B(x_{1}\otimes gx_{1}x_{2})\overset{\left( \ref{form x1otgx1x2}\right) }{=}%
+(1_{A}\otimes g)B(gx_{1}\otimes 1_{H})(1_{A}\otimes gx_{1}x_{2}) \\
&&-(1_{A}\otimes x_{2})B(gx_{1}\otimes 1_{H})(1_{A}\otimes x_{1})+ \\
&&+B(x_{1}x_{2}\otimes 1_{H})(1_{A}\otimes x_{1})+ \\
&&+(1_{A}\otimes x_{1})B(gx_{1}\otimes 1_{H})(1_{A}\otimes x_{2})+ \\
&&+(1_{A}\otimes gx_{1}x_{2})B(gx_{1}\otimes 1_{H})(1_{A}\otimes g)+ \\
&&-(1_{A}\otimes gx_{1})B(x_{1}x_{2}\otimes 1_{H})(1_{A}\otimes g)+
\end{eqnarray*}%
we get%
\begin{eqnarray*}
&&\text{left side }B(x_{1}\otimes gx_{1}x_{2})(1_{A}\otimes x_{2}) \\
&=&(1_{A}\otimes g)B(gx_{1}\otimes 1_{H})(1_{A}\otimes
gx_{1}x_{2})(1_{A}\otimes x_{2})+ \\
&&-(1_{A}\otimes x_{2})B(gx_{1}\otimes 1_{H})(1_{A}\otimes
x_{1})(1_{A}\otimes x_{2})+ \\
&&+B(x_{1}x_{2}\otimes 1_{H})(1_{A}\otimes x_{1})(1_{A}\otimes x_{2})+ \\
&&+(1_{A}\otimes x_{1})B(gx_{1}\otimes 1_{H})(1_{A}\otimes
x_{2})(1_{A}\otimes x_{2})+ \\
&&+(1_{A}\otimes gx_{1}x_{2})B(gx_{1}\otimes 1_{H})(1_{A}\otimes
g)(1_{A}\otimes x_{2})+ \\
&&-(1_{A}\otimes gx_{1})B(x_{1}x_{2}\otimes 1_{H})(1_{A}\otimes
g)(1_{A}\otimes x_{2})+
\end{eqnarray*}%
i.e.%
\begin{eqnarray*}
&&\text{left side }B(x_{1}\otimes gx_{1}x_{2})(1_{A}\otimes x_{2}) \\
&=&(1_{A}\otimes x_{2})B(gx_{1}\otimes 1_{H})(1_{A}\otimes x_{1}x_{2})+ \\
&&-B(x_{1}x_{2}\otimes 1_{H})(1_{A}\otimes x_{1}x_{2})+ \\
&&-(1_{A}\otimes gx_{1}x_{2})B(gx_{1}\otimes 1_{H})(1_{A}\otimes gx_{2})+ \\
&&+(1_{A}\otimes gx_{1})B(x_{1}x_{2}\otimes 1_{H})(1_{A}\otimes gx_{2})+
\end{eqnarray*}%
and%
\begin{eqnarray*}
&&(1_{A}\otimes x_{2})\left( 1_{A}\otimes g\right) B(x_{1}\otimes
gx_{1}x_{2})\left( 1_{A}\otimes g\right) \\
&&=(1_{A}\otimes x_{2})\left( 1_{A}\otimes g\right) (1_{A}\otimes
g)B(gx_{1}\otimes 1_{H})(1_{A}\otimes gx_{1}x_{2})\left( 1_{A}\otimes
g\right) \\
&&-(1_{A}\otimes x_{2})\left( 1_{A}\otimes g\right) (1_{A}\otimes
x_{2})B(gx_{1}\otimes 1_{H})(1_{A}\otimes x_{1})\left( 1_{A}\otimes g\right)
\\
&&+(1_{A}\otimes x_{2})\left( 1_{A}\otimes g\right) B(x_{1}x_{2}\otimes
1_{H})(1_{A}\otimes x_{1})\left( 1_{A}\otimes g\right) \\
&&+(1_{A}\otimes x_{2})\left( 1_{A}\otimes g\right) (1_{A}\otimes
x_{1})B(gx_{1}\otimes 1_{H})(1_{A}\otimes x_{2})\left( 1_{A}\otimes g\right)
\\
&&+(1_{A}\otimes x_{2})\left( 1_{A}\otimes g\right) (1_{A}\otimes
gx_{1}x_{2})B(gx_{1}\otimes 1_{H})(1_{A}\otimes g)\left( 1_{A}\otimes
g\right) \\
&&-(1_{A}\otimes x_{2})\left( 1_{A}\otimes g\right) (1_{A}\otimes
gx_{1})B(x_{1}x_{2}\otimes 1_{H})(1_{A}\otimes g)\left( 1_{A}\otimes g\right)
\end{eqnarray*}%
i.e.%
\begin{eqnarray*}
&&(1_{A}\otimes x_{2})\left( 1_{A}\otimes g\right) B(x_{1}\otimes
gx_{1}x_{2})\left( 1_{A}\otimes g\right) \\
= &&+(1_{A}\otimes x_{2})B(gx_{1}\otimes 1_{H})(1_{A}\otimes x_{1}x_{2}) \\
&&+(1_{A}\otimes gx_{2})B(x_{1}x_{2}\otimes 1_{H})(1_{A}\otimes gx_{1}) \\
&&-(1_{A}\otimes gx_{1}x_{2})\left( 1_{A}\otimes \right) (1_{A}\otimes
)B(gx_{1}\otimes 1_{H})(1_{A}\otimes gx_{2}) \\
&&+(1_{A}\otimes x_{1}x_{2})B(x_{1}x_{2}\otimes 1_{H}).
\end{eqnarray*}%
Since%
\begin{eqnarray*}
- &&B(x_{1}x_{2}\otimes x_{1}x_{2})\overset{\left( \ref{form x1x2otx1x2}%
\right) }{=} \\
&&-B(x_{1}x_{2}\otimes 1_{H})(1_{A}\otimes x_{1}x_{2}) \\
&&-(1_{A}\otimes gx_{2})B(x_{1}x_{2}\otimes 1_{H})(1_{A}\otimes gx_{1}) \\
&&+(1_{A}\otimes gx_{1})B(x_{1}x_{2}\otimes 1_{H})(1_{A}\otimes gx_{2})+ \\
&&-(1_{A}\otimes x_{1}x_{2})B(x_{1}x_{2}\otimes 1_{H})
\end{eqnarray*}%
we get%
\begin{eqnarray*}
\text{right side }(1_{A}\otimes x_{2})\left( 1_{A}\otimes g\right)
&&B(x_{1}\otimes gx_{1}x_{2})\left( 1_{A}\otimes g\right)
-B(x_{1}x_{2}\otimes x_{1}x_{2})= \\
&&+(1_{A}\otimes x_{2})B(gx_{1}\otimes 1_{H})(1_{A}\otimes x_{1}x_{2}) \\
&&+(1_{A}\otimes gx_{2})B(x_{1}x_{2}\otimes 1_{H})(1_{A}\otimes gx_{1}) \\
&&-(1_{A}\otimes gx_{1}x_{2})B(gx_{1}\otimes 1_{H})(1_{A}\otimes gx_{2}) \\
&&+(1_{A}\otimes x_{1}x_{2})B(x_{1}x_{2}\otimes 1_{H}) \\
&&-B(x_{1}x_{2}\otimes 1_{H})(1_{A}\otimes x_{1}x_{2}) \\
&&-(1_{A}\otimes gx_{2})B(x_{1}x_{2}\otimes 1_{H})(1_{A}\otimes gx_{1}) \\
&&+(1_{A}\otimes gx_{1})B(x_{1}x_{2}\otimes 1_{H})(1_{A}\otimes gx_{2})+ \\
&&-(1_{A}\otimes x_{1}x_{2})B(x_{1}x_{2}\otimes 1_{H})
\end{eqnarray*}%
i.e.%
\begin{eqnarray*}
\text{clean right side }(1_{A}\otimes x_{2})\left( 1_{A}\otimes g\right)
&&B(x_{1}\otimes gx_{1}x_{2})\left( 1_{A}\otimes g\right)
-B(x_{1}x_{2}\otimes x_{1}x_{2})= \\
&&+(1_{A}\otimes x_{2})B(gx_{1}\otimes 1_{H})(1_{A}\otimes x_{1}x_{2}) \\
&&-(1_{A}\otimes gx_{1}x_{2})B(gx_{1}\otimes 1_{H})(1_{A}\otimes gx_{2}) \\
&&-B(x_{1}x_{2}\otimes 1_{H})(1_{A}\otimes x_{1}x_{2}) \\
&&+(1_{A}\otimes gx_{1})B(x_{1}x_{2}\otimes 1_{H})(1_{A}\otimes gx_{2})+
\end{eqnarray*}%
and we conclude.

\subsubsection{$x_{2}\otimes x_{1}$}

\begin{eqnarray*}
&&B(x_{2}\otimes x_{1})(1_{A}\otimes x_{2})\overset{?}{=}(1_{A}\otimes
x_{2})\left( 1_{A}\otimes g\right) B(x_{2}\otimes x_{1})\left( 1_{A}\otimes
g\right) \\
&&+B(x_{2}x_{2}\otimes gx_{1})+B(x_{2}\otimes x_{2}x_{1})
\end{eqnarray*}%
i.e.%
\begin{equation*}
B(x_{2}\otimes x_{1})(1_{A}\otimes x_{2})\overset{?}{=}(1_{A}\otimes
x_{2})\left( 1_{A}\otimes g\right) B(x_{2}\otimes x_{1})\left( 1_{A}\otimes
g\right) -B(x_{2}\otimes x_{1}x_{2}).
\end{equation*}%
Since
\begin{eqnarray*}
&&B(x_{2}\otimes x_{1})\overset{\left( \ref{form x2otx1}\right) }{=}%
B(x_{2}\otimes 1_{H})(1_{A}\otimes x_{1}) \\
&&-(1_{A}\otimes gx_{1})B(x_{2}\otimes 1_{H})(1_{A}\otimes g) \\
&&-(1_{A}\otimes g)B(gx_{1}x_{2}\otimes 1_{H})(1_{A}\otimes g)
\end{eqnarray*}%
we get%
\begin{eqnarray*}
\text{left side } &&B(x_{2}\otimes x_{1})(1_{A}\otimes )=-B(x_{2}\otimes
1_{H})(1_{A}\otimes x_{1}x_{2}) \\
&&+(1_{A}\otimes gx_{1})B(x_{2}\otimes 1_{H})(1_{A}\otimes gx_{2}) \\
&&+(1_{A}\otimes g)B(gx_{1}x_{2}\otimes 1_{H})(1_{A}\otimes gx_{2})
\end{eqnarray*}%
and%
\begin{eqnarray*}
&&(1_{A}\otimes x_{2})\left( 1_{A}\otimes g\right) B(x_{2}\otimes
x_{1})\left( 1_{A}\otimes g\right) \\
&=&(1_{A}\otimes x_{2})\left( 1_{A}\otimes g\right) B(x_{2}\otimes
1_{H})(1_{A}\otimes x_{1})\left( 1_{A}\otimes g\right) \\
&&-(1_{A}\otimes x_{2})\left( 1_{A}\otimes g\right) (1_{A}\otimes
gx_{1})B(x_{2}\otimes 1_{H})(1_{A}\otimes g)\left( 1_{A}\otimes g\right) \\
&&-(1_{A}\otimes x_{2})\left( 1_{A}\otimes g\right) (1_{A}\otimes
g)B(gx_{1}x_{2}\otimes 1_{H})(1_{A}\otimes g)\left( 1_{A}\otimes g\right)
\end{eqnarray*}%
i.e.%
\begin{eqnarray*}
&&(1_{A}\otimes x_{2})\left( 1_{A}\otimes g\right) B(x_{2}\otimes
x_{1})\left( 1_{A}\otimes g\right) \\
&=&(1_{A}\otimes gx_{2})B(x_{2}\otimes 1_{H})(1_{A}\otimes gx_{1}) \\
&&+(1_{A}\otimes x_{1}x_{2})\left( 1_{A}\otimes \right) (1_{A}\otimes
)B(x_{2}\otimes 1_{H}) \\
&&-(1_{A}\otimes x_{2})B(gx_{1}x_{2}\otimes 1_{H}).
\end{eqnarray*}%
Since%
\begin{eqnarray*}
- &&B(x_{2}\otimes x_{1}x_{2})\overset{\left( \ref{form x2otx1x2}\right) }{=}
\\
&&-B(x_{2}\otimes 1_{H})(1_{A}\otimes x_{1}x_{2}) \\
&&-(1_{A}\otimes gx_{2})B(x_{2}\otimes 1_{H})(1_{A}\otimes gx_{1}) \\
&&+(1_{A}\otimes gx_{1})B(x_{2}\otimes 1_{H})(1_{A}\otimes gx_{2}) \\
&&-(1_{A}\otimes x_{1}x_{2})B(x_{2}\otimes 1_{H}) \\
&&+(1_{A}\otimes g)B(gx_{1}x_{2}\otimes 1_{H})(1_{A}\otimes gx_{2}) \\
&&+(1_{A}\otimes x_{2})B(gx_{1}x_{2}\otimes 1_{H}).
\end{eqnarray*}%
we get%
\begin{gather*}
\text{right side }(1_{A}\otimes x_{2})\left( 1_{A}\otimes g\right)
B(x_{2}\otimes x_{1})\left( 1_{A}\otimes g\right) -B(x_{2}\otimes
x_{1}x_{2})= \\
(1_{A}\otimes gx_{2})B(x_{2}\otimes 1_{H})(1_{A}\otimes gx_{1}) \\
+(1_{A}\otimes x_{1}x_{2})B(x_{2}\otimes 1_{H}) \\
-(1_{A}\otimes x_{2})B(gx_{1}x_{2}\otimes 1_{H}) \\
-B(x_{2}\otimes 1_{H})(1_{A}\otimes x_{1}x_{2}) \\
-(1_{A}\otimes gx_{2})B(x_{2}\otimes 1_{H})(1_{A}\otimes gx_{1}) \\
+(1_{A}\otimes gx_{1})B(x_{2}\otimes 1_{H})(1_{A}\otimes gx_{2}) \\
-(1_{A}\otimes x_{1}x_{2})B(x_{2}\otimes 1_{H}) \\
+(1_{A}\otimes g)B(gx_{1}x_{2}\otimes 1_{H})(1_{A}\otimes gx_{2}) \\
+(1_{A}\otimes x_{2})B(gx_{1}x_{2}\otimes 1_{H}).
\end{gather*}%
i.e.%
\begin{gather*}
\text{clean right side }(1_{A}\otimes x_{2})\left( 1_{A}\otimes g\right)
B(x_{2}\otimes x_{1})\left( 1_{A}\otimes g\right) -B(x_{2}\otimes
x_{1}x_{2})= \\
-B(x_{2}\otimes 1_{H})(1_{A}\otimes x_{1}x_{2}) \\
+(1_{A}\otimes gx_{1})B(x_{2}\otimes 1_{H})(1_{A}\otimes gx_{2}) \\
+(1_{A}\otimes g)B(gx_{1}x_{2}\otimes 1_{H})(1_{A}\otimes gx_{2})
\end{gather*}%
and we conclude.

\subsubsection{$x_{2}\otimes x_{2}$}

\begin{eqnarray*}
&&B(x_{2}\otimes x_{2})(1_{A}\otimes x_{2})\overset{?}{=}(1_{A}\otimes
x_{2})\left( 1_{A}\otimes g\right) B(x_{2}\otimes x_{2})\left( 1_{A}\otimes
g\right) \\
&&+B(x_{2}x_{2}\otimes gx_{2})+B(x_{2}\otimes x_{2}x_{2})
\end{eqnarray*}%
i.e.%
\begin{equation*}
B(x_{2}\otimes x_{2})(1_{A}\otimes x_{2})\overset{?}{=}(1_{A}\otimes
x_{2})\left( 1_{A}\otimes g\right) B(x_{2}\otimes x_{2})\left( 1_{A}\otimes
g\right) .
\end{equation*}%
Since%
\begin{equation*}
B(x_{2}\otimes x_{2})\overset{\left( \ref{form x2otx2}\right) }{=}%
B(x_{2}\otimes 1_{H})(1_{A}\otimes x_{2})-(1_{A}\otimes
gx_{2})B(x_{2}\otimes 1_{H})(1_{A}\otimes g).
\end{equation*}%
we get%
\begin{eqnarray*}
&&\text{left side }B(x_{2}\otimes x_{2})(1_{A}\otimes x_{2}) \\
&=&B(x_{2}\otimes 1_{H})(1_{A}\otimes x_{2})(1_{A}\otimes x_{2})+ \\
&&-(1_{A}\otimes gx_{2})B(x_{2}\otimes 1_{H})(1_{A}\otimes g)(1_{A}\otimes
x_{2})
\end{eqnarray*}%
i.e.%
\begin{equation*}
\text{left side }B(x_{2}\otimes x_{2})(1_{A}\otimes x_{2})=+(1_{A}\otimes
gx_{2})B(x_{2}\otimes 1_{H})(1_{A}\otimes gx_{2})
\end{equation*}%
and%
\begin{eqnarray*}
&&\text{right side }(1_{A}\otimes x_{2})\left( 1_{A}\otimes g\right)
B(x_{2}\otimes x_{2})\left( 1_{A}\otimes g\right) \\
&=&(1_{A}\otimes x_{2})\left( 1_{A}\otimes g\right) B(x_{2}\otimes
1_{H})(1_{A}\otimes x_{2})\left( 1_{A}\otimes g\right) \\
&&-(1_{A}\otimes x_{2})\left( 1_{A}\otimes g\right) (1_{A}\otimes
gx_{2})B(x_{2}\otimes 1_{H})(1_{A}\otimes g)\left( 1_{A}\otimes g\right)
\end{eqnarray*}%
i.e.%
\begin{equation*}
\text{right side }(1_{A}\otimes x_{2})\left( 1_{A}\otimes g\right)
B(x_{2}\otimes x_{2})\left( 1_{A}\otimes g\right) =(1_{A}\otimes
gx_{2})B(x_{2}\otimes 1_{H})(1_{A}\otimes xg_{2})
\end{equation*}%
and we conclude.

\subsubsection{$x_{2}\otimes x_{1}x_{2}$}

\begin{eqnarray*}
&&B(x_{2}\otimes x_{1}x_{2})(1_{A}\otimes x_{2})\overset{?}{=}(1_{A}\otimes
x_{2})\left( 1_{A}\otimes g\right) B(x_{2}\otimes x_{1}x_{2})\left(
1_{A}\otimes g\right) \\
&&+B(x_{2}x_{2}\otimes gx_{1}x_{2})+B(x_{2}\otimes x_{2}x_{1}x_{2})
\end{eqnarray*}%
i.e.%
\begin{equation*}
B(x_{2}\otimes x_{1}x_{2})(1_{A}\otimes x_{2})\overset{?}{=}(1_{A}\otimes
x_{2})\left( 1_{A}\otimes g\right) B(x_{2}\otimes x_{1}x_{2})\left(
1_{A}\otimes g\right)
\end{equation*}%
Since%
\begin{eqnarray*}
&&B(x_{2}\otimes x_{1}x_{2})\overset{\left( \ref{form x2otx1x2}\right) }{=}%
B(x_{2}\otimes 1_{H})(1_{A}\otimes x_{1}x_{2}) \\
&&+(1_{A}\otimes gx_{2})B(x_{2}\otimes 1_{H})(1_{A}\otimes gx_{1}) \\
&&-(1_{A}\otimes gx_{1})B(x_{2}\otimes 1_{H})(1_{A}\otimes gx_{2}) \\
&&+(1_{A}\otimes x_{1}x_{2})B(x_{2}\otimes 1_{H}) \\
&&-(1_{A}\otimes g)B(gx_{1}x_{2}\otimes 1_{H})(1_{A}\otimes gx_{2}) \\
&&-(1_{A}\otimes x_{2})B(gx_{1}x_{2}\otimes 1_{H}).
\end{eqnarray*}%
we get%
\begin{eqnarray*}
&&\text{left side }B(x_{2}\otimes x_{1}x_{2})(1_{A}\otimes
x_{2})=B(x_{2}\otimes 1_{H})(1_{A}\otimes x_{1}x_{2})(1_{A}\otimes x_{2}) \\
&&+(1_{A}\otimes gx_{2})B(x_{2}\otimes 1_{H})(1_{A}\otimes
gx_{1})(1_{A}\otimes x_{2}) \\
&&-(1_{A}\otimes gx_{1})B(x_{2}\otimes 1_{H})(1_{A}\otimes
gx_{2})(1_{A}\otimes x_{2}) \\
&&+(1_{A}\otimes x_{1}x_{2})B(x_{2}\otimes 1_{H})(1_{A}\otimes x_{2}) \\
&&-(1_{A}\otimes g)B(gx_{1}x_{2}\otimes 1_{H})(1_{A}\otimes
gx_{2})(1_{A}\otimes x_{2}) \\
&&-(1_{A}\otimes x_{2})B(gx_{1}x_{2}\otimes 1_{H})(1_{A}\otimes x_{2})
\end{eqnarray*}%
i.e.%
\begin{eqnarray*}
&&\text{left side }B(x_{2}\otimes x_{1}x_{2})(1_{A}\otimes x_{2})= \\
&&+(1_{A}\otimes gx_{2})B(x_{2}\otimes 1_{H})(1_{A}\otimes gx_{1}x_{2}) \\
&&+(1_{A}\otimes x_{1}x_{2})B(x_{2}\otimes 1_{H})(1_{A}\otimes x_{2}) \\
&&-(1_{A}\otimes x_{2})B(gx_{1}x_{2}\otimes 1_{H})(1_{A}\otimes x_{2})
\end{eqnarray*}%
and%
\begin{eqnarray*}
\text{right side } &&(1_{A}\otimes x_{2})\left( 1_{A}\otimes g\right)
B(x_{2}\otimes x_{1}x_{2})\left( 1_{A}\otimes g\right) = \\
&&(1_{A}\otimes x_{2})\left( 1_{A}\otimes g\right) B(x_{2}\otimes
1_{H})(1_{A}\otimes x_{1}x_{2})\left( 1_{A}\otimes g\right) \\
&&+(1_{A}\otimes x_{2})\left( 1_{A}\otimes g\right) (1_{A}\otimes
gx_{2})B(x_{2}\otimes 1_{H})(1_{A}\otimes gx_{1})\left( 1_{A}\otimes g\right)
\\
&&-(1_{A}\otimes x_{2})\left( 1_{A}\otimes g\right) (1_{A}\otimes
gx_{1})B(x_{2}\otimes 1_{H})(1_{A}\otimes gx_{2})\left( 1_{A}\otimes g\right)
\\
&&+(1_{A}\otimes x_{2})\left( 1_{A}\otimes g\right) (1_{A}\otimes
x_{1}x_{2})B(x_{2}\otimes 1_{H})\left( 1_{A}\otimes g\right) \\
&&-(1_{A}\otimes x_{2})\left( 1_{A}\otimes g\right) (1_{A}\otimes
g)B(gx_{1}x_{2}\otimes 1_{H})(1_{A}\otimes gx_{2})\left( 1_{A}\otimes
g\right) \\
&&-(1_{A}\otimes x_{2})\left( 1_{A}\otimes g\right) (1_{A}\otimes
x_{2})B(gx_{1}x_{2}\otimes 1_{H})\left( 1_{A}\otimes g\right)
\end{eqnarray*}%
i.e.%
\begin{eqnarray*}
\text{right side } &&(1_{A}\otimes x_{2})\left( 1_{A}\otimes g\right)
B(x_{2}\otimes x_{1}x_{2})\left( 1_{A}\otimes g\right) = \\
&&(1_{A}\otimes x_{2})\left( 1_{A}\otimes g\right) B(x_{2}\otimes
1_{H})(1_{A}\otimes gx_{1}x_{2})+ \\
&&+(1_{A}\otimes x_{1}x_{2})B(x_{2}\otimes 1_{H})(1_{A}\otimes x_{2})+ \\
&&-(1_{A}\otimes x_{2})B(gx_{1}x_{2}\otimes 1_{H})(1_{A}\otimes x_{2})
\end{eqnarray*}%
and we conclude.

\subsubsection{$x_{2}\otimes gx_{1}$}

\begin{eqnarray*}
&&B(x_{2}\otimes gx_{1})(1_{A}\otimes x_{2})\overset{?}{=}(1_{A}\otimes
x_{2})\left( 1_{A}\otimes g\right) B(x_{2}\otimes gx_{1})\left( 1_{A}\otimes
g\right) + \\
&&+B(x_{2}x_{2}\otimes ggx_{1})+B(x_{2}\otimes x_{2}gx_{1})
\end{eqnarray*}%
i.e%
\begin{equation*}
B(x_{2}\otimes gx_{1})(1_{A}\otimes x_{2})\overset{?}{=}(1_{A}\otimes
x_{2})\left( 1_{A}\otimes g\right) B(x_{2}\otimes gx_{1})\left( 1_{A}\otimes
g\right) +B(x_{2}\otimes gx_{1}x_{2}).
\end{equation*}%
Since%
\begin{eqnarray*}
&&B(x_{2}\otimes gx_{1})\overset{\left( \ref{form x2otgx1}\right) }{=}%
(1_{A}\otimes g)B(gx_{2}\otimes 1_{H})(1_{A}\otimes gx_{1})+ \\
&&+(1_{A}\otimes x_{1})B(gx_{2}\otimes 1_{H})+B(x_{1}x_{2}\otimes 1_{H})
\end{eqnarray*}%
we get%
\begin{eqnarray*}
\text{left side }B(x_{2}\otimes gx_{1}) &=&(1_{A}\otimes g)B(gx_{2}\otimes
1_{H})(1_{A}\otimes gx_{1}x_{2})+ \\
&&+(1_{A}\otimes x_{1})B(gx_{2}\otimes 1_{H})(1_{A}\otimes x_{2})+ \\
&&+B(x_{1}x_{2}\otimes 1_{H})(1_{A}\otimes x_{2})
\end{eqnarray*}%
and%
\begin{eqnarray*}
&&(1_{A}\otimes x_{2})\left( 1_{A}\otimes g\right) B(x_{2}\otimes
gx_{1})\left( 1_{A}\otimes g\right) \\
&=&(1_{A}\otimes x_{2})\left( 1_{A}\otimes g\right) (1_{A}\otimes
g)B(gx_{2}\otimes 1_{H})(1_{A}\otimes gx_{1})\left( 1_{A}\otimes g\right) +
\\
&&+(1_{A}\otimes x_{2})\left( 1_{A}\otimes g\right) (1_{A}\otimes
x_{1})B(gx_{2}\otimes 1_{H})\left( 1_{A}\otimes g\right) + \\
&&+(1_{A}\otimes x_{2})\left( 1_{A}\otimes g\right) B(x_{1}x_{2}\otimes
1_{H})\left( 1_{A}\otimes g\right)
\end{eqnarray*}%
i.e.%
\begin{eqnarray*}
&&(1_{A}\otimes x_{2})\left( 1_{A}\otimes g\right) B(x_{2}\otimes
gx_{1})\left( 1_{A}\otimes g\right) \\
&=&(1_{A}\otimes x_{2})B(gx_{2}\otimes 1_{H})(1_{A}\otimes x_{1}) \\
&&-(1_{A}\otimes gx_{1}x_{2})B(gx_{2}\otimes 1_{H})\left( 1_{A}\otimes
g\right) \\
&&+(1_{A}\otimes gx_{2})B(x_{1}x_{2}\otimes 1_{H})\left( 1_{A}\otimes
g\right) .
\end{eqnarray*}%
Since%
\begin{eqnarray*}
&&B(x_{2}\otimes gx_{1}x_{2})\overset{\left( \ref{form x2otgx1x2}\right) }{=}%
(1_{A}\otimes g)B(gx_{2}\otimes 1_{H})(1_{A}\otimes gx_{1}x_{2}) \\
&&-(1_{A}\otimes x_{2})B(gx_{2}\otimes 1_{H})(1_{A}\otimes x_{1}) \\
&&+(1_{A}\otimes x_{1})B(gx_{2}\otimes 1_{H})(1_{A}\otimes x_{2}) \\
&&+(1_{A}\otimes gx_{1}x_{2})B(gx_{2}\otimes 1_{H})(1_{A}\otimes g) \\
&&+B(x_{1}x_{2}\otimes 1_{H})(1_{A}\otimes x_{2}) \\
&&-(1_{A}\otimes gx_{2})B(x_{1}x_{2}\otimes 1_{H})(1_{A}\otimes g)
\end{eqnarray*}%
we get%
\begin{eqnarray*}
\text{right side } &&(1_{A}\otimes x_{2})\left( 1_{A}\otimes g\right)
B(x_{2}\otimes gx_{1})\left( 1_{A}\otimes g\right) +B(x_{2}\otimes
gx_{1}x_{2})= \\
&&(1_{A}\otimes x_{2})B(gx_{2}\otimes 1_{H})(1_{A}\otimes x_{1})+ \\
&&-(1_{A}\otimes gx_{1}x_{2})B(gx_{2}\otimes 1_{H})\left( 1_{A}\otimes
g\right) \\
&&+(1_{A}\otimes gx_{2})B(x_{1}x_{2}\otimes 1_{H})\left( 1_{A}\otimes
g\right) \\
&&(1_{A}\otimes g)B(gx_{2}\otimes 1_{H})(1_{A}\otimes gx_{1}x_{2}) \\
&&-(1_{A}\otimes x_{2})B(gx_{2}\otimes 1_{H})(1_{A}\otimes x_{1}) \\
&&+(1_{A}\otimes x_{1})B(gx_{2}\otimes 1_{H})(1_{A}\otimes x_{2}) \\
&&+(1_{A}\otimes gx_{1}x_{2})B(gx_{2}\otimes 1_{H})(1_{A}\otimes g) \\
&&+B(x_{1}x_{2}\otimes 1_{H})(1_{A}\otimes x_{2}) \\
&&-(1_{A}\otimes gx_{2})B(x_{1}x_{2}\otimes 1_{H})(1_{A}\otimes g)
\end{eqnarray*}%
i.e.%
\begin{eqnarray*}
\text{clean right side } &&(1_{A}\otimes x_{2})\left( 1_{A}\otimes g\right)
B(x_{2}\otimes gx_{1})\left( 1_{A}\otimes g\right) +B(x_{2}\otimes
gx_{1}x_{2})= \\
&&(1_{A}\otimes g)B(gx_{2}\otimes 1_{H})(1_{A}\otimes gx_{1}x_{2}) \\
&&+(1_{A}\otimes x_{1})B(gx_{2}\otimes 1_{H})(1_{A}\otimes x_{2}) \\
&&+B(x_{1}x_{2}\otimes 1_{H})(1_{A}\otimes x_{2})
\end{eqnarray*}%
and we conclude.

\subsubsection{$x_{2}\otimes gx_{2}$}

\begin{eqnarray*}
&&B(x_{2}\otimes gx_{2})(1_{A}\otimes x_{2})\overset{?}{=}(1_{A}\otimes
x_{2})\left( 1_{A}\otimes g\right) B(x_{2}\otimes gx_{2})\left( 1_{A}\otimes
g\right) + \\
&&+B(x_{2}x_{2}\otimes ggx_{2})+B(x_{2}\otimes x_{2}gx_{2})
\end{eqnarray*}%
i.e.%
\begin{equation*}
B(x_{2}\otimes gx_{2})(1_{A}\otimes x_{2})\overset{?}{=}(1_{A}\otimes
x_{2})\left( 1_{A}\otimes g\right) B(x_{2}\otimes gx_{2})\left( 1_{A}\otimes
g\right)
\end{equation*}%
Since%
\begin{equation*}
B(x_{2}\otimes gx_{2})\overset{\left( \ref{form x2otgx2}\right) }{=}%
(1_{A}\otimes g)B(gx_{2}\otimes 1_{H})(1_{A}\otimes gx_{2})+(1_{A}\otimes
x_{2})B(gx_{2}\otimes 1_{H})
\end{equation*}%
we get%
\begin{eqnarray*}
\text{left side }B(x_{2}\otimes gx_{2})(1_{A}\otimes x_{2}) &=&(1_{A}\otimes
g)B(gx_{2}\otimes 1_{H})(1_{A}\otimes gx_{2})(1_{A}\otimes x_{2}) \\
&&+(1_{A}\otimes x_{2})B(gx_{2}\otimes 1_{H})(1_{A}\otimes x_{2})
\end{eqnarray*}%
i.e.%
\begin{equation*}
\text{left side }B(x_{2}\otimes gx_{2})(1_{A}\otimes x_{2})=(1_{A}\otimes
x_{2})B(gx_{2}\otimes 1_{H})(1_{A}\otimes x_{2})
\end{equation*}%
and%
\begin{eqnarray*}
\text{right side } &=&(1_{A}\otimes x_{2})\left( 1_{A}\otimes g\right)
B(x_{2}\otimes gx_{2})\left( 1_{A}\otimes g\right) \\
&=&(1_{A}\otimes x_{2})\left( 1_{A}\otimes g\right) (1_{A}\otimes
g)B(gx_{2}\otimes 1_{H})(1_{A}\otimes gx_{2})\left( 1_{A}\otimes g\right) \\
&&+(1_{A}\otimes x_{2})\left( 1_{A}\otimes g\right) (1_{A}\otimes
x_{2})B(gx_{2}\otimes 1_{H})\left( 1_{A}\otimes g\right)
\end{eqnarray*}%
i.e.%
\begin{equation*}
\text{right side }=(1_{A}\otimes x_{2})\left( 1_{A}\otimes g\right)
B(x_{2}\otimes gx_{2})\left( 1_{A}\otimes g\right) =(1_{A}\otimes
x_{2})B(gx_{2}\otimes 1_{H})(1_{A}\otimes x_{2})
\end{equation*}%
and we conclude.

\subsubsection{$x_{2}\otimes gx_{1}x_{2}$}

\begin{eqnarray*}
&&B(x_{2}\otimes gx_{1}x_{2})(1_{A}\otimes x_{2})\overset{?}{=}(1_{A}\otimes
x_{2})\left( 1_{A}\otimes g\right) B(x_{2}\otimes gx_{1}x_{2})\left(
1_{A}\otimes g\right) + \\
&&+B(x_{2}x_{2}\otimes ggx_{1}x_{2})+B(x_{2}\otimes x_{2}gx_{1}x_{2})
\end{eqnarray*}%
i.e.

\begin{equation*}
B(x_{2}\otimes gx_{1}x_{2})(1_{A}\otimes x_{2})\overset{?}{=}(1_{A}\otimes
x_{2})\left( 1_{A}\otimes g\right) B(x_{2}\otimes gx_{1}x_{2})\left(
1_{A}\otimes g\right)
\end{equation*}%
Since%
\begin{eqnarray*}
&&B(x_{2}\otimes gx_{1}x_{2})\overset{\left( \ref{form x2otgx1x2}\right) }{=}%
(1_{A}\otimes g)B(gx_{2}\otimes 1_{H})(1_{A}\otimes gx_{1}x_{2}) \\
&&-(1_{A}\otimes x_{2})B(gx_{2}\otimes 1_{H})(1_{A}\otimes x_{1}) \\
&&+(1_{A}\otimes x_{1})B(gx_{2}\otimes 1_{H})(1_{A}\otimes x_{2}) \\
&&+(1_{A}\otimes gx_{1}x_{2})B(gx_{2}\otimes 1_{H})(1_{A}\otimes g) \\
&&+B(x_{1}x_{2}\otimes 1_{H})(1_{A}\otimes x_{2}) \\
&&-(1_{A}\otimes gx_{2})B(x_{1}x_{2}\otimes 1_{H})(1_{A}\otimes g)
\end{eqnarray*}%
we get%
\begin{eqnarray*}
\text{left side } &&B(x_{2}\otimes gx_{1}x_{2})(1_{A}\otimes
x_{2})=(1_{A}\otimes g)B(gx_{2}\otimes 1_{H})(1_{A}\otimes
gx_{1}x_{2})(1_{A}\otimes x_{2}) \\
&&-(1_{A}\otimes x_{2})B(gx_{2}\otimes 1_{H})(1_{A}\otimes
x_{1})(1_{A}\otimes x_{2}) \\
&&+(1_{A}\otimes x_{1})B(gx_{2}\otimes 1_{H})(1_{A}\otimes
x_{2})(1_{A}\otimes x_{2}) \\
&&+(1_{A}\otimes gx_{1}x_{2})B(gx_{2}\otimes 1_{H})(1_{A}\otimes
g)(1_{A}\otimes x_{2}) \\
&&+B(x_{1}x_{2}\otimes 1_{H})(1_{A}\otimes x_{2})(1_{A}\otimes x_{2}) \\
&&-(1_{A}\otimes gx_{2})B(x_{1}x_{2}\otimes 1_{H})(1_{A}\otimes
g)(1_{A}\otimes x_{2})
\end{eqnarray*}%
i.e.%
\begin{eqnarray*}
\text{left side } &&B(x_{2}\otimes gx_{1}x_{2})(1_{A}\otimes x_{2})= \\
&&+(1_{A}\otimes x_{2})B(gx_{2}\otimes 1_{H})(1_{A}\otimes x_{1}x_{2}) \\
&&-(1_{A}\otimes gx_{1}x_{2})B(gx_{2}\otimes 1_{H})(1_{A}\otimes gx_{2}) \\
&&+(1_{A}\otimes gx_{2})B(x_{1}x_{2}\otimes 1_{H})(1_{A}\otimes gx_{2})
\end{eqnarray*}%
and%
\begin{eqnarray*}
&&\text{right side }(1_{A}\otimes x_{2})\left( 1_{A}\otimes g\right)
B(x_{2}\otimes gx_{1}x_{2})\left( 1_{A}\otimes g\right) = \\
&&(1_{A}\otimes x_{2})\left( 1_{A}\otimes g\right) (1_{A}\otimes
g)B(gx_{2}\otimes 1_{H})(1_{A}\otimes gx_{1}x_{2})\left( 1_{A}\otimes
g\right) \\
&&-(1_{A}\otimes x_{2})\left( 1_{A}\otimes g\right) (1_{A}\otimes
x_{2})B(gx_{2}\otimes 1_{H})(1_{A}\otimes x_{1})\left( 1_{A}\otimes g\right)
\\
&&+(1_{A}\otimes x_{2})\left( 1_{A}\otimes g\right) (1_{A}\otimes
x_{1})B(gx_{2}\otimes 1_{H})(1_{A}\otimes x_{2})\left( 1_{A}\otimes g\right)
\\
&&+(1_{A}\otimes x_{2})\left( 1_{A}\otimes g\right) (1_{A}\otimes
gx_{1}x_{2})B(gx_{2}\otimes 1_{H})(1_{A}\otimes g)\left( 1_{A}\otimes
g\right) \\
&&+(1_{A}\otimes x_{2})\left( 1_{A}\otimes g\right) B(x_{1}x_{2}\otimes
1_{H})(1_{A}\otimes x_{2})\left( 1_{A}\otimes g\right) \\
&&-(1_{A}\otimes x_{2})\left( 1_{A}\otimes g\right) (1_{A}\otimes
gx_{2})B(x_{1}x_{2}\otimes 1_{H})(1_{A}\otimes g)\left( 1_{A}\otimes g\right)
\end{eqnarray*}%
i.e.%
\begin{eqnarray*}
&&\text{right side }(1_{A}\otimes x_{2})\left( 1_{A}\otimes g\right)
B(x_{2}\otimes gx_{1}x_{2})\left( 1_{A}\otimes g\right) = \\
&&(1_{A}\otimes x_{2})B(gx_{2}\otimes 1_{H})(1_{A}\otimes x_{1}x_{2}) \\
&&-(1_{A}\otimes gx_{1}x_{2})B(gx_{2}\otimes 1_{H})(1_{A}\otimes gx_{2}) \\
&&+(1_{A}\otimes gx_{2})B(x_{1}x_{2}\otimes 1_{H})(1_{A}\otimes gx_{2})
\end{eqnarray*}%
and we conclude.

\subsubsection{$x_{1}x_{2}\otimes x_{1}$}

\begin{eqnarray*}
&&B(x_{1}x_{2}\otimes x_{1})(1_{A}\otimes x_{2})\overset{?}{=}(1_{A}\otimes
x_{2})\left( 1_{A}\otimes g\right) B(x_{1}x_{2}\otimes x_{1})\left(
1_{A}\otimes g\right) + \\
&&+B(x_{2}x_{1}x_{2}\otimes gx_{1})+B(x_{1}x_{2}\otimes x_{2}x_{1})
\end{eqnarray*}%
i.e.%
\begin{equation*}
B(x_{1}x_{2}\otimes x_{1})(1_{A}\otimes x_{2})\overset{?}{=}(1_{A}\otimes
x_{2})\left( 1_{A}\otimes g\right) B(x_{1}x_{2}\otimes x_{1})\left(
1_{A}\otimes g\right) -B(x_{1}x_{2}\otimes x_{1}x_{2})
\end{equation*}%
Since%
\begin{equation*}
B(x_{1}x_{2}\otimes x_{1})\overset{\left( \ref{form x1x2otx1}\right) }{=}%
B(x_{1}x_{2}\otimes 1_{H})(1_{A}\otimes x_{1})-(1_{A}\otimes
gx_{1})B(x_{1}x_{2}\otimes 1_{H})(1_{A}\otimes g)
\end{equation*}%
we have%
\begin{eqnarray*}
&&\text{left side }B(x_{1}x_{2}\otimes x_{1})(1_{A}\otimes x_{2}) \\
&=&B(x_{1}x_{2}\otimes 1_{H})(1_{A}\otimes x_{1})(1_{A}\otimes x_{2}) \\
&&-(1_{A}\otimes gx_{1})B(x_{1}x_{2}\otimes 1_{H})(1_{A}\otimes
g)(1_{A}\otimes x_{2})
\end{eqnarray*}%
i.e.%
\begin{eqnarray*}
&&\text{left side }B(x_{1}x_{2}\otimes x_{1})(1_{A}\otimes x_{2}) \\
&=&-B(x_{1}x_{2}\otimes 1_{H})(1_{A}\otimes x_{1}x_{2})+(1_{A}\otimes
gx_{1})B(x_{1}x_{2}\otimes 1_{H})(1_{A}\otimes gx_{2})
\end{eqnarray*}%
and%
\begin{eqnarray*}
&&(1_{A}\otimes x_{2})\left( 1_{A}\otimes g\right) B(x_{1}x_{2}\otimes
x_{1})\left( 1_{A}\otimes g\right) \\
&=&(1_{A}\otimes x_{2})\left( 1_{A}\otimes g\right) B(x_{1}x_{2}\otimes
1_{H})(1_{A}\otimes x_{1})\left( 1_{A}\otimes g\right) + \\
&&-(1_{A}\otimes x_{2})\left( 1_{A}\otimes g\right) (1_{A}\otimes
gx_{1})B(x_{1}x_{2}\otimes 1_{H})(1_{A}\otimes g)\left( 1_{A}\otimes g\right)
\end{eqnarray*}%
i.e.%
\begin{eqnarray*}
&&(1_{A}\otimes x_{2})\left( 1_{A}\otimes g\right) B(x_{1}x_{2}\otimes
x_{1})\left( 1_{A}\otimes g\right) \\
&=&(1_{A}\otimes gx_{2})B(x_{1}x_{2}\otimes 1_{H})(1_{A}\otimes gx_{1})+ \\
&&+(1_{A}\otimes x_{1}x_{2})B(x_{1}x_{2}\otimes 1_{H}).
\end{eqnarray*}%
Since%
\begin{eqnarray*}
&&B(x_{1}x_{2}\otimes x_{1}x_{2})\overset{\left( \ref{form x1x2otx1x2}%
\right) }{=}B(x_{1}x_{2}\otimes 1_{H})(1_{A}\otimes x_{1}x_{2}) \\
&&+(1_{A}\otimes gx_{2})B(x_{1}x_{2}\otimes 1_{H})(1_{A}\otimes gx_{1}) \\
&&-(1_{A}\otimes gx_{1})B(x_{1}x_{2}\otimes 1_{H})(1_{A}\otimes gx_{2})+ \\
&&+(1_{A}\otimes x_{1}x_{2})B(x_{1}x_{2}\otimes 1_{H})
\end{eqnarray*}%
we get%
\begin{eqnarray*}
&&\text{right side }(1_{A}\otimes x_{2})\left( 1_{A}\otimes g\right)
B(x_{1}x_{2}\otimes x_{1})\left( 1_{A}\otimes g\right) -B(x_{1}x_{2}\otimes
x_{1}x_{2})= \\
&&(1_{A}\otimes gx_{2})B(x_{1}x_{2}\otimes 1_{H})(1_{A}\otimes gx_{1})+ \\
&&+(1_{A}\otimes x_{1}x_{2})B(x_{1}x_{2}\otimes 1_{H}) \\
&&-B(x_{1}x_{2}\otimes 1_{H})(1_{A}\otimes x_{1}x_{2}) \\
&&-(1_{A}\otimes gx_{2})B(x_{1}x_{2}\otimes 1_{H})(1_{A}\otimes gx_{1}) \\
&&+(1_{A}\otimes gx_{1})B(x_{1}x_{2}\otimes 1_{H})(1_{A}\otimes gx_{2})+ \\
&&-(1_{A}\otimes x_{1}x_{2})B(x_{1}x_{2}\otimes 1_{H})+
\end{eqnarray*}%
i.e.%
\begin{eqnarray*}
&&\text{clean right side }(1_{A}\otimes x_{2})\left( 1_{A}\otimes g\right)
B(x_{1}x_{2}\otimes x_{1})\left( 1_{A}\otimes g\right) -B(x_{1}x_{2}\otimes
x_{1}x_{2})= \\
&=&B(x_{1}x_{2}\otimes 1_{H})(1_{A}\otimes x_{1}x_{2})+(1_{A}\otimes
gx_{1})B(x_{1}x_{2}\otimes 1_{H})(1_{A}\otimes gx_{2})
\end{eqnarray*}%
and we conclude.

\subsubsection{$x_{1}x_{2}\otimes x_{2}$}

\begin{eqnarray*}
&&B(x_{1}x_{2}\otimes x_{2})(1_{A}\otimes x_{2})\overset{?}{=}(1_{A}\otimes
x_{2})\left( 1_{A}\otimes g\right) B(x_{1}x_{2}\otimes x_{2})\left(
1_{A}\otimes g\right) \\
&&+B(x_{2}x_{1}x_{2}\otimes gx_{2})+B(x_{1}x_{2}\otimes x_{2}x_{2})
\end{eqnarray*}%
i.e.%
\begin{equation*}
B(x_{1}x_{2}\otimes x_{2})(1_{A}\otimes x_{2})\overset{?}{=}(1_{A}\otimes
x_{2})\left( 1_{A}\otimes g\right) B(x_{1}x_{2}\otimes x_{2})\left(
1_{A}\otimes g\right)
\end{equation*}%
Since%
\begin{equation*}
B(x_{1}x_{2}\otimes x_{2})\overset{\left( \ref{x1x2otx2}\right) }{=}%
B(x_{1}x_{2}\otimes 1_{H})(1_{A}\otimes x_{2})-(1_{A}\otimes
gx_{2})B(x_{1}x_{2}\otimes 1_{H})(1_{A}\otimes g)
\end{equation*}%
we get%
\begin{eqnarray*}
&&\text{left side }B(x_{1}x_{2}\otimes x_{2})(1_{A}\otimes x_{2}) \\
&=&B(x_{1}x_{2}\otimes 1_{H})(1_{A}\otimes x_{2})(1_{A}\otimes
x_{2})-(1_{A}\otimes gx_{2})B(x_{1}x_{2}\otimes 1_{H})(1_{A}\otimes
g)(1_{A}\otimes x_{2})
\end{eqnarray*}%
i.e.%
\begin{equation*}
\text{left side }B(x_{1}x_{2}\otimes x_{2})(1_{A}\otimes
x_{2})=(1_{A}\otimes gx_{2})B(x_{1}x_{2}\otimes 1_{H})(1_{A}\otimes gx_{2})
\end{equation*}%
and%
\begin{eqnarray*}
&&\text{right side }(1_{A}\otimes x_{2})\left( 1_{A}\otimes g\right)
B(x_{1}x_{2}\otimes x_{2})\left( 1_{A}\otimes g\right) \\
&=&(1_{A}\otimes x_{2})\left( 1_{A}\otimes g\right) B(x_{1}x_{2}\otimes
1_{H})(1_{A}\otimes x_{2})\left( 1_{A}\otimes g\right) + \\
&&-(1_{A}\otimes x_{2})\left( 1_{A}\otimes g\right) (1_{A}\otimes
gx_{2})B(x_{1}x_{2}\otimes 1_{H})(1_{A}\otimes g)\left( 1_{A}\otimes g\right)
\end{eqnarray*}%
i.e.%
\begin{equation*}
\text{right side }(1_{A}\otimes x_{2})\left( 1_{A}\otimes g\right)
B(x_{1}x_{2}\otimes x_{2})\left( 1_{A}\otimes g\right) =(1_{A}\otimes
gx_{2})B(x_{1}x_{2}\otimes 1_{H})(1_{A}\otimes gx_{2})
\end{equation*}%
and we conclude.

\subsubsection{$x_{1}x_{2}\otimes x_{1}x_{2}$}

\begin{eqnarray*}
&&B(x_{1}x_{2}\otimes x_{1}x_{2})(1_{A}\otimes x_{2})\overset{?}{=}%
(1_{A}\otimes x_{2})\left( 1_{A}\otimes g\right) B(x_{1}x_{2}\otimes
x_{1}x_{2})\left( 1_{A}\otimes g\right) \\
&&+B(x_{2}x_{1}x_{2}\otimes gx_{1}x_{2})+B(x_{1}x_{2}\otimes x_{2}x_{1}x_{2})
\end{eqnarray*}%
i.e.%
\begin{equation*}
B(x_{1}x_{2}\otimes x_{1}x_{2})(1_{A}\otimes x_{2})\overset{?}{=}%
(1_{A}\otimes x_{2})\left( 1_{A}\otimes g\right) B(x_{1}x_{2}\otimes
x_{1}x_{2})\left( 1_{A}\otimes g\right)
\end{equation*}%
Since%
\begin{eqnarray*}
&&B(x_{1}x_{2}\otimes x_{1}x_{2})\overset{\left( \ref{form x1x2otx1x2}%
\right) }{=}B(x_{1}x_{2}\otimes 1_{H})(1_{A}\otimes x_{1}x_{2}) \\
&&+(1_{A}\otimes gx_{2})B(x_{1}x_{2}\otimes 1_{H})(1_{A}\otimes gx_{1}) \\
&&-(1_{A}\otimes gx_{1})B(x_{1}x_{2}\otimes 1_{H})(1_{A}\otimes gx_{2})+ \\
&&+(1_{A}\otimes x_{1}x_{2})B(x_{1}x_{2}\otimes 1_{H})+
\end{eqnarray*}%
we get%
\begin{eqnarray*}
&&\text{left side }B(x_{1}x_{2}\otimes x_{1}x_{2})(1_{A}\otimes
x_{2})=B(x_{1}x_{2}\otimes 1_{H})(1_{A}\otimes x_{1}x_{2})(1_{A}\otimes
x_{2}) \\
&&+(1_{A}\otimes gx_{2})B(x_{1}x_{2}\otimes 1_{H})(1_{A}\otimes
gx_{1})(1_{A}\otimes x_{2})+ \\
&&-(1_{A}\otimes gx_{1})B(x_{1}x_{2}\otimes 1_{H})(1_{A}\otimes
gx_{2})(1_{A}\otimes x_{2})+ \\
&&+(1_{A}\otimes x_{1}x_{2})B(x_{1}x_{2}\otimes 1_{H})(1_{A}\otimes x_{2})
\end{eqnarray*}%
i.e.%
\begin{eqnarray*}
&&\text{left side }B(x_{1}x_{2}\otimes x_{1}x_{2})(1_{A}\otimes x_{2})= \\
&&+(1_{A}\otimes gx_{2})B(x_{1}x_{2}\otimes 1_{H})(1_{A}\otimes gx_{1}x_{2})+
\\
&&+(1_{A}\otimes x_{1}x_{2})B(x_{1}x_{2}\otimes 1_{H})(1_{A}\otimes x_{2})
\end{eqnarray*}%
and%
\begin{eqnarray*}
&&\text{right side }(1_{A}\otimes x_{2})\left( 1_{A}\otimes g\right)
B(x_{1}x_{2}\otimes x_{1}x_{2})\left( 1_{A}\otimes g\right) \\
&=&(1_{A}\otimes x_{2})\left( 1_{A}\otimes g\right) B(x_{1}x_{2}\otimes
1_{H})(1_{A}\otimes x_{1}x_{2})\left( 1_{A}\otimes g\right) \\
&&+(1_{A}\otimes x_{2})\left( 1_{A}\otimes g\right) (1_{A}\otimes
gx_{2})B(x_{1}x_{2}\otimes 1_{H})(1_{A}\otimes gx_{1})\left( 1_{A}\otimes
g\right) \\
&&-(1_{A}\otimes x_{2})\left( 1_{A}\otimes g\right) (1_{A}\otimes
gx_{1})B(x_{1}x_{2}\otimes 1_{H})(1_{A}\otimes gx_{2})\left( 1_{A}\otimes
g\right) + \\
&&+(1_{A}\otimes x_{2})\left( 1_{A}\otimes g\right) (1_{A}\otimes
x_{1}x_{2})B(x_{1}x_{2}\otimes 1_{H})\left( 1_{A}\otimes g\right) +
\end{eqnarray*}%
i.e.%
\begin{eqnarray*}
&&\text{right side }(1_{A}\otimes x_{2})\left( 1_{A}\otimes g\right)
B(x_{1}x_{2}\otimes x_{1}x_{2})\left( 1_{A}\otimes g\right) \\
&=&(1_{A}\otimes x_{2})\left( 1_{A}\otimes g\right) B(x_{1}x_{2}\otimes
1_{H})(1_{A}\otimes gx_{1}x_{2}) \\
&&+(1_{A}\otimes x_{1}x_{2})\left( 1_{A}\otimes gg\right) (1_{A}\otimes
)B(x_{1}x_{2}\otimes 1_{H})(1_{A}\otimes x_{2})
\end{eqnarray*}%
and we conclude.

\subsubsection{$x_{1}x_{2}\otimes gx_{1}$}

\begin{eqnarray*}
&&B((x_{1}x_{2}\otimes gx_{1})(1_{A}\otimes x_{2})\overset{?}{=}%
(1_{A}\otimes x_{2})\left( 1_{A}\otimes g\right) B(x_{1}x_{2}\otimes
gx_{1})\left( 1_{A}\otimes g\right) \\
&&+B(x_{2}x_{1}x_{2}\otimes ggx_{1})+B(x_{1}x_{2}\otimes x_{2}gx_{1})
\end{eqnarray*}%
i.e.%
\begin{equation*}
B(x_{1}x_{2}\otimes gx_{1})(1_{A}\otimes x_{2})\overset{?}{=}(1_{A}\otimes
x_{2})\left( 1_{A}\otimes g\right) B(x_{1}x_{2}\otimes gx_{1})\left(
1_{A}\otimes g\right) +B(x_{1}x_{2}\otimes gx_{1}x_{2})
\end{equation*}%
Since%
\begin{eqnarray*}
B(x_{1}x_{2}\otimes gx_{1}) &&\overset{\left( \ref{form x1x2otgx1}\right) }{=%
}(1_{A}\otimes g)B(gx_{1}x_{2}\otimes 1_{H})(1_{A}\otimes gx_{1}) \\
&&+(1_{A}\otimes x_{1})B(gx_{1}x_{2}\otimes 1_{H})
\end{eqnarray*}%
we get%
\begin{eqnarray*}
\text{left side }B(x_{1}x_{2}\otimes gx_{1}) &&(1_{A}\otimes
x_{2})=(1_{A}\otimes g)B(gx_{1}x_{2}\otimes 1_{H})(1_{A}\otimes
gx_{1})(1_{A}\otimes x_{2}) \\
&&+(1_{A}\otimes x_{1})B(gx_{1}x_{2}\otimes 1_{H})(1_{A}\otimes x_{2})
\end{eqnarray*}%
i.e%
\begin{eqnarray*}
\text{left side }B(x_{1}x_{2}\otimes gx_{1}) &&(1_{A}\otimes
x_{2})=(1_{A}\otimes g)B(gx_{1}x_{2}\otimes 1_{H})(1_{A}\otimes gx_{1}x_{2})
\\
&&+(1_{A}\otimes x_{1})B(gx_{1}x_{2}\otimes 1_{H})(1_{A}\otimes x_{2})
\end{eqnarray*}%
and%
\begin{eqnarray*}
&&(1_{A}\otimes x_{2})\left( 1_{A}\otimes g\right) B(x_{1}x_{2}\otimes
gx_{1})\left( 1_{A}\otimes g\right) \\
&=&(1_{A}\otimes x_{2})\left( 1_{A}\otimes g\right) (1_{A}\otimes
g)B(gx_{1}x_{2}\otimes 1_{H})(1_{A}\otimes gx_{1})\left( 1_{A}\otimes
g\right) + \\
&&+(1_{A}\otimes x_{2})\left( 1_{A}\otimes g\right) (1_{A}\otimes
x_{1})B(gx_{1}x_{2}\otimes 1_{H})\left( 1_{A}\otimes g\right)
\end{eqnarray*}%
i.e.%
\begin{eqnarray*}
&&(1_{A}\otimes x_{2})\left( 1_{A}\otimes g\right) B(x_{1}x_{2}\otimes
gx_{1})\left( 1_{A}\otimes g\right) \\
&=&(1_{A}\otimes x_{2})B(gx_{1}x_{2}\otimes 1_{H})(1_{A}\otimes x_{1})+ \\
&&-(1_{A}\otimes gx_{1}x_{2})B(gx_{1}x_{2}\otimes 1_{H})\left( 1_{A}\otimes
g\right)
\end{eqnarray*}%
Since
\begin{eqnarray*}
&&B(x_{1}x_{2}\otimes gx_{1}x_{2})\overset{\left( \ref{form x1x2otgx1gx2}%
\right) }{=}(1_{A}\otimes g)B(gx_{1}x_{2}\otimes 1_{H})(1_{A}\otimes
gx_{1}x_{2}) \\
&&-(1_{A}\otimes x_{2})B(gx_{1}x_{2}\otimes 1_{H})(1_{A}\otimes x_{1})+ \\
&&+(1_{A}\otimes x_{1})B(gx_{1}x_{2}\otimes 1_{H})(1_{A}\otimes x_{2})+ \\
&&+(1_{A}\otimes gx_{1}x_{2})B(gx_{1}x_{2}\otimes 1_{H})(1_{A}\otimes g)+
\end{eqnarray*}%
we get%
\begin{eqnarray*}
&&\text{right side }(1_{A}\otimes x_{2})\left( 1_{A}\otimes g\right)
B(x_{1}x_{2}\otimes gx_{1})\left( 1_{A}\otimes g\right) +B(x_{1}x_{2}\otimes
gx_{1}x_{2})= \\
&&(1_{A}\otimes x_{2})B(gx_{1}x_{2}\otimes 1_{H})(1_{A}\otimes x_{1})+ \\
&&-(1_{A}\otimes gx_{1}x_{2})B(gx_{1}x_{2}\otimes 1_{H})\left( 1_{A}\otimes
g\right) \\
&&(1_{A}\otimes g)B(gx_{1}x_{2}\otimes 1_{H})(1_{A}\otimes gx_{1}x_{2}) \\
&&-(1_{A}\otimes x_{2})B(gx_{1}x_{2}\otimes 1_{H})(1_{A}\otimes x_{1}) \\
&&+(1_{A}\otimes x_{1})B(gx_{1}x_{2}\otimes 1_{H})(1_{A}\otimes x_{2}) \\
&&+(1_{A}\otimes gx_{1}x_{2})B(gx_{1}x_{2}\otimes 1_{H})(1_{A}\otimes g)
\end{eqnarray*}%
and we get%
\begin{eqnarray*}
&&\text{clean right side }(1_{A}\otimes x_{2})\left( 1_{A}\otimes g\right)
B(x_{1}x_{2}\otimes gx_{1})\left( 1_{A}\otimes g\right) +B(x_{1}x_{2}\otimes
gx_{1}x_{2})= \\
&&(1_{A}\otimes g)B(gx_{1}x_{2}\otimes 1_{H})(1_{A}\otimes gx_{1}x_{2}) \\
&&+(1_{A}\otimes x_{1})B(gx_{1}x_{2}\otimes 1_{H})(1_{A}\otimes x_{2})
\end{eqnarray*}%
and we conclude.

\subsubsection{$x_{1}x_{2}\otimes gx_{2}$}

\begin{eqnarray*}
&&B(x_{1}x_{2}\otimes gx_{2})(1_{A}\otimes x_{2})\overset{?}{=}(1_{A}\otimes
x_{2})\left( 1_{A}\otimes g\right) B(x_{1}x_{2}\otimes gx_{2})\left(
1_{A}\otimes g\right) \\
&&+B(x_{2}x_{1}x_{2}\otimes ggx_{2})+B(x_{1}x_{2}\otimes x_{2}gx_{2})
\end{eqnarray*}%
i.e.%
\begin{equation*}
B(x_{1}x_{2}\otimes gx_{2})(1_{A}\otimes x_{2})\overset{?}{=}(1_{A}\otimes
x_{2})\left( 1_{A}\otimes g\right) B(x_{1}x_{2}\otimes gx_{2})\left(
1_{A}\otimes g\right) .
\end{equation*}%
Since%
\begin{equation*}
B(x_{1}x_{2}\otimes gx_{2})\overset{\left( \ref{form x1x2otgx2}\right) }{=}%
(1_{A}\otimes g)B(gx_{1}x_{2}\otimes 1_{H})(1_{A}\otimes
gx_{2})+(1_{A}\otimes x_{2})B(gx_{1}x_{2}\otimes 1_{H})
\end{equation*}%
we get%
\begin{eqnarray*}
&&\text{left side }B(x_{1}x_{2}\otimes gx_{2})(1_{A}\otimes x_{2}) \\
&=&(1_{A}\otimes g)B(gx_{1}x_{2}\otimes 1_{H})(1_{A}\otimes
gx_{2})(1_{A}\otimes x_{2}) \\
&&+(1_{A}\otimes x_{2})B(gx_{1}x_{2}\otimes 1_{H})(1_{A}\otimes x_{2})
\end{eqnarray*}%
i.e.%
\begin{equation*}
\text{left side }B(x_{1}x_{2}\otimes gx_{2})(1_{A}\otimes
x_{2})=(1_{A}\otimes x_{2})B(gx_{1}x_{2}\otimes 1_{H})(1_{A}\otimes x_{2})
\end{equation*}%
and%
\begin{eqnarray*}
&&\text{right side }(1_{A}\otimes x_{2})\left( 1_{A}\otimes g\right)
B(x_{1}x_{2}\otimes gx_{2})\left( 1_{A}\otimes g\right) \\
&=&(1_{A}\otimes x_{2})\left( 1_{A}\otimes g\right) (1_{A}\otimes
g)B(gx_{1}x_{2}\otimes 1_{H})(1_{A}\otimes gx_{2})\left( 1_{A}\otimes
g\right) \\
&&+(1_{A}\otimes x_{2})\left( 1_{A}\otimes g\right) (1_{A}\otimes
x_{2})B(gx_{1}x_{2}\otimes 1_{H})\left( 1_{A}\otimes g\right)
\end{eqnarray*}%
i.e.%
\begin{eqnarray*}
&&\text{right side }(1_{A}\otimes x_{2})\left( 1_{A}\otimes g\right)
B(x_{1}x_{2}\otimes gx_{2})\left( 1_{A}\otimes g\right) \\
&=&(1_{A}\otimes x_{2})B(gx_{1}x_{2}\otimes 1_{H})(1_{A}\otimes x_{2})
\end{eqnarray*}%
and we conclude.

\subsubsection{$x_{1}x_{2}\otimes gx_{1}x_{2}$}

\begin{eqnarray*}
&&B(x_{1}x_{2}\otimes gx_{1}x_{2})(1_{A}\otimes x_{2})\overset{?}{=}%
(1_{A}\otimes x_{2})\left( 1_{A}\otimes g\right) B(x_{1}x_{2}\otimes
gx_{1}x_{2})\left( 1_{A}\otimes g\right) + \\
&&+B(x_{2}x_{1}x_{2}\otimes ggx_{1}x_{2})+B(x_{1}x_{2}\otimes
x_{2}gx_{1}x_{2})
\end{eqnarray*}%
i.e.%
\begin{equation*}
B(x_{1}x_{2}\otimes gx_{1}x_{2})(1_{A}\otimes x_{2})\overset{?}{=}%
(1_{A}\otimes x_{2})\left( 1_{A}\otimes g\right) B(x_{1}x_{2}\otimes
gx_{1}x_{2})\left( 1_{A}\otimes g\right)
\end{equation*}%
Since%
\begin{eqnarray*}
&&B(x_{1}x_{2}\otimes gx_{1}x_{2})\overset{\left( \ref{form x1x2otgx1gx2}%
\right) }{=}(1_{A}\otimes g)B(gx_{1}x_{2}\otimes 1_{H})(1_{A}\otimes
gx_{1}x_{2}) \\
&&-(1_{A}\otimes x_{2})B(gx_{1}x_{2}\otimes 1_{H})(1_{A}\otimes x_{1}) \\
&&+(1_{A}\otimes x_{1})B(gx_{1}x_{2}\otimes 1_{H})(1_{A}\otimes x_{2}) \\
&&+(1_{A}\otimes gx_{1}x_{2})B(gx_{1}x_{2}\otimes 1_{H})(1_{A}\otimes g)
\end{eqnarray*}%
we get%
\begin{eqnarray*}
\text{left side } &&B(x_{1}x_{2}\otimes gx_{1}x_{2})(1_{A}\otimes
x_{2})=(1_{A}\otimes g)B(gx_{1}x_{2}\otimes 1_{H})(1_{A}\otimes
gx_{1}x_{2})(1_{A}\otimes x_{2}) \\
&&-(1_{A}\otimes x_{2})B(gx_{1}x_{2}\otimes 1_{H})(1_{A}\otimes
x_{1})(1_{A}\otimes x_{2}) \\
&&+(1_{A}\otimes x_{1})B(gx_{1}x_{2}\otimes 1_{H})(1_{A}\otimes
x_{2})(1_{A}\otimes x_{2}) \\
&&+(1_{A}\otimes gx_{1}x_{2})B(gx_{1}x_{2}\otimes 1_{H})(1_{A}\otimes
g)(1_{A}\otimes x_{2})
\end{eqnarray*}%
i.e.%
\begin{eqnarray*}
\text{left side } &&B(x_{1}x_{2}\otimes gx_{1}x_{2})(1_{A}\otimes x_{2})= \\
&&+(1_{A}\otimes x_{2})B(gx_{1}x_{2}\otimes 1_{H})(1_{A}\otimes x_{1}x_{2})
\\
&&-(1_{A}\otimes gx_{1}x_{2})B(gx_{1}x_{2}\otimes 1_{H})(1_{A}\otimes gx_{2})
\end{eqnarray*}%
and%
\begin{eqnarray*}
\text{right side } &&(1_{A}\otimes x_{2})\left( 1_{A}\otimes g\right)
B(x_{1}x_{2}\otimes gx_{1}x_{2})\left( 1_{A}\otimes g\right) \\
&=&\text{ }(1_{A}\otimes x_{2})\left( 1_{A}\otimes g\right) (1_{A}\otimes
g)B(gx_{1}x_{2}\otimes 1_{H})(1_{A}\otimes gx_{1}x_{2})\left( 1_{A}\otimes
g\right) + \\
&&-\text{ }(1_{A}\otimes x_{2})\left( 1_{A}\otimes g\right) (1_{A}\otimes
x_{2})B(gx_{1}x_{2}\otimes 1_{H})(1_{A}\otimes x_{1})\left( 1_{A}\otimes
g\right) + \\
&&+\text{ }(1_{A}\otimes x_{2})\left( 1_{A}\otimes g\right) (1_{A}\otimes
x_{1})B(gx_{1}x_{2}\otimes 1_{H})(1_{A}\otimes x_{2})\left( 1_{A}\otimes
g\right) + \\
&&+\text{ }(1_{A}\otimes x_{2})\left( 1_{A}\otimes g\right) (1_{A}\otimes
gx_{1}x_{2})B(gx_{1}x_{2}\otimes 1_{H})(1_{A}\otimes g)\left( 1_{A}\otimes
g\right)
\end{eqnarray*}%
i.e.%
\begin{eqnarray*}
\text{clean right side } &&(1_{A}\otimes x_{2})\left( 1_{A}\otimes g\right)
B(x_{1}x_{2}\otimes gx_{1}x_{2})\left( 1_{A}\otimes g\right) \\
&=&\text{ }(1_{A}\otimes x_{2})B(gx_{1}x_{2}\otimes 1_{H})(1_{A}\otimes
x_{1}x_{2})+ \\
&&-\text{ }(1_{A}\otimes gx_{1}x_{2})B(gx_{1}x_{2}\otimes
1_{H})(1_{A}\otimes gx_{2})
\end{eqnarray*}%
and we conclude.

\section{ Equation $\left( \protect\ref{eq.a}\right) $}

In this section we will explore which conditions should be satisfied in
order that $\left( \ref{eq.a}\right) $ holds. The conditions that we will
get will be, together with the form of $B\left( h\otimes 1_{H}\right) $
where $h$ runs through our seven elements, plus the obvious condition $%
B(1_{H}\otimes 1_{H})=1_{A}\otimes 1_{H},$the ones that will appear in our
main final Theorem \ref{Theo sep}.

We recall $\left( \ref{eq.a}\right) $%
\begin{equation*}
B(h\otimes h^{\prime })(d\otimes 1_{H})=(d_{0}\otimes 1_{H})B(hd_{1}\otimes
h^{\prime }d_{2}).
\end{equation*}%
As we remarked in section $\left( \ref{mor cond}\right) $, we only need to
compute equality $\left( \ref{eq.a}\right) $ for $G\otimes 1_{H}$, $%
X_{1}\otimes 1_{H}$ and $X_{2}\otimes 1_{H}$. Moreover, as we proved in $%
\left( \ref{REa}\right) ,$assuming that equation $\left( \ref{eq.h}\right) $
is satisfied, then if equation $\left( \ref{eq.a}\right) $ holds for $%
h\otimes 1_{H}$ then it holds for any $h\otimes h^{\prime }.$

\subsection{ $d=G$}

\begin{equation}
B(h\otimes h^{\prime })(G\otimes 1_{H})=(G\otimes 1_{H})B(hg\otimes
h^{\prime }g).  \label{eq. dG}
\end{equation}

For $h^{\prime }=1_{H}$ we get%
\begin{equation}
B(h\otimes 1_{H})(G\otimes 1_{H})=(G\otimes 1_{H})B(hg\otimes g).
\label{eq. d1}
\end{equation}%
For $h^{\prime }=g$ we get%
\begin{equation}
B(h\otimes g)(G\otimes 1_{H})=(G\otimes 1_{H})B(hg\otimes 1_{H})
\label{eq.d2}
\end{equation}%
Assume that $\left( \ref{eq. d1}\right) $ and $\left( \ref{eq.10}\right) $
hold and let us prove that $\left( \ref{eq.d2}\right) $ hold.
\begin{eqnarray*}
&&B(h\otimes g)(G\otimes 1_{H})\overset{\left( \ref{eq.10}\right) }{=}\left(
1_{A}\otimes g\right) B(gh\otimes 1_{H})\left( 1_{A}\otimes g\right)
(G\otimes 1_{H}) \\
&=&\left( 1_{A}\otimes g\right) B(gh\otimes 1_{H})(G\otimes 1_{H})\left(
1_{A}\otimes g\right) \\
&&\overset{\left( \ref{eq. d1}\right) }{=}\left( 1_{A}\otimes g\right)
(G\otimes 1_{H})B(ghg\otimes g)\left( 1_{A}\otimes g\right) \\
&=&(G\otimes 1_{H})\left( 1_{A}\otimes g\right) B(ghg\otimes g)\left(
1_{A}\otimes g\right) \\
&&\overset{\left( \ref{eq.10}\right) }{=}(G\otimes 1_{H})B(hg\otimes 1_{H})
\end{eqnarray*}

\begin{eqnarray*}
B(h\otimes h^{\prime }) &=&\sum_{a,b_{1},b_{2}=0, f}^{1}B(h\otimes h^{\prime
};G^{a}X_{1}^{b_{1}}X_{2}^{b_{2}},f)G^{a}X_{1}^{b_{1}}X_{2}^{b_{2}}\otimes f
\\
&=&\sum_{a,d,e_{1},e_{2}=0}^{1}B(h\otimes h^{\prime };G^{a},f)G^{a}\otimes f
\\
&&+\sum_{a,d,e_{1},e_{2}=0}^{1}B(h\otimes h^{\prime
};G^{a}X_{1},f)G^{a}X_{1}\otimes f \\
&&+\sum_{a,d,e_{1},e_{2}=0}^{1}B(h\otimes h^{\prime
};G^{a}X_{2},f)G^{a}X_{2}\otimes f \\
&&+\sum_{a,d,e_{1},e_{2}=0}^{1}B(h\otimes h^{\prime
};G^{a}X_{1}X_{2},f)G^{a}X_{1}X_{2}\otimes f
\end{eqnarray*}

\begin{equation*}
X_{i}G+GX_{i}=\gamma _{i}
\end{equation*}
and%
\begin{equation*}
X_{i}G=\gamma _{i}-GX_{i}
\end{equation*}%
\begin{equation*}
X_{2}G=\gamma _{2}-GX_{2}
\end{equation*}%
\begin{equation*}
X_{1}G=\gamma _{1}-GX_{1}
\end{equation*}

\begin{equation*}
X_{1}X_{2}G=\gamma _{2}X_{1}-X_{1}GX_{2}=\gamma _{2}X_{1}-\gamma
_{1}X_{2}+GX_{1}X_{2}
\end{equation*}%
\begin{equation*}
X_{1}X_{2}+X_{2}X_{1}=\lambda
\end{equation*}%
We compute the left side%
\begin{eqnarray*}
&&B(h\otimes h^{\prime })\left( G\otimes 1_{H}\right) \\
&=&\sum_{a,d,e_{1},e_{2}=0}^{1}B(h\otimes h^{\prime };G^{a},f)G^{a+1}\otimes
f \\
&&+\sum_{a,d,e_{1},e_{2}=0}^{1}B(h\otimes h^{\prime
};G^{a}X_{1},f)G^{a}\left( \gamma _{1}-GX_{1}\right) \otimes f \\
&&+\sum_{a,d,e_{1},e_{2}=0}^{1}B(h\otimes h^{\prime
};G^{a}X_{2},f)G^{a}\left( \gamma _{2}-GX_{2}\right) \otimes f \\
&&+\sum_{a,d,e_{1},e_{2}=0}^{1}B(h\otimes h^{\prime
};G^{a}X_{1}X_{2},f)G^{a}\left( \gamma _{2}X_{1}-\gamma
_{1}X_{2}+GX_{1}X_{2}\right) \otimes f
\end{eqnarray*}%
so that%
\begin{eqnarray*}
&&B(h\otimes h^{\prime })\left( G\otimes 1_{H}\right) \\
&=&\sum_{a,d,e_{1},e_{2}=0}^{1}\left[
\begin{array}{c}
B(h\otimes h^{\prime };G^{a},f)G^{a+1} \\
\left[ +\gamma _{1}B(h\otimes h^{\prime };G^{a}X_{1},f)+\gamma
_{2}B(h\otimes h^{\prime };G^{a}X_{2},f)\right] G^{a}%
\end{array}%
\right] \otimes f \\
&&+\sum_{a,d,e_{1},e_{2}=0}^{1}B(h\otimes h^{\prime };G^{a}X_{1},f)\left[
\gamma _{2}G^{a}-G^{a+1}\right] X_{1}\otimes f \\
&&+\sum_{a,d,e_{1},e_{2}=0}^{1}B(h\otimes h^{\prime };G^{a}X_{2},f)\left[
-\gamma _{1}G^{a}-G^{a+1}\right] X_{2}\otimes f \\
&&+\sum_{a,d,e_{1},e_{2}=0}^{1}B(h\otimes h^{\prime
};G^{a}X_{1}X_{2},f)G^{a+1}X_{1}X_{2}\otimes f
\end{eqnarray*}%
and hence%
\begin{eqnarray*}
&&B(h\otimes h^{\prime })\left( G\otimes 1_{H}\right) \\
&=&\sum_{d,e_{1},e_{2}=0}^{1}\left[
\begin{array}{c}
B(h\otimes h^{\prime };1_{A},f)G \\
\left[ +\gamma _{1}B(h\otimes h^{\prime };X_{1},f)+\gamma _{2}B(h\otimes
h^{\prime };X_{2},f)\right] 1_{A}%
\end{array}%
\right] \otimes f \\
&&\sum_{d,e_{1},e_{2}=0}^{1}\left[
\begin{array}{c}
\alpha B(h\otimes h^{\prime };G,f)1_{A}+ \\
\left[ +\gamma _{1}B(h\otimes h^{\prime };GX_{1},f)+\gamma _{2}B(h\otimes
h^{\prime };GX_{2},f)\right] G%
\end{array}%
\right] \otimes f \\
&&+\sum_{d,e_{1},e_{2}=0}^{1}B(h\otimes h^{\prime };X_{1},f)\left[ \gamma
_{2}-G\right] X_{1}\otimes f \\
&&+\sum_{d,e_{1},e_{2}=0}^{1}B(h\otimes h^{\prime };GX_{1},f)\left[ \gamma
_{2}G-\alpha \right] X_{1}\otimes f \\
&&+\sum_{d,e_{1},e_{2}=0}^{1}B(h\otimes h^{\prime };X_{2},f)\left[ -\gamma
_{1}-G\right] X_{2}\otimes f \\
&&+\sum_{d,e_{1},e_{2}=0}^{1}B(h\otimes h^{\prime };GX_{2},f)\left[ -\gamma
_{1}G-\alpha \right] X_{2}\otimes f \\
&&+\sum_{d,e_{1},e_{2}=0}^{1}B(h\otimes h^{\prime
};X_{1}X_{2},f)GX_{1}X_{2}\otimes f \\
&&+\sum_{d,e_{1},e_{2}=0}^{1}B(h\otimes h^{\prime };GX_{1}X_{2},f)\alpha
X_{1}X_{2}\otimes f
\end{eqnarray*}

We compute the right side%
\begin{eqnarray*}
(G\otimes 1_{H})B(hg\otimes h^{\prime }g)
&=&\sum_{a,b_{1},b_{2},d,e_{1},e_{2}=0}^{1}B(hg\otimes h^{\prime
}g;G^{a}X_{1}^{b_{1}}X_{2}^{b_{2}},f)G^{a+1}X_{1}^{b_{1}}X_{2}^{b_{2}}%
\otimes f \\
&=&\sum_{b_{1},b_{2},d,e_{1},e_{2}=0}^{1}B(hg\otimes h^{\prime
}g;X_{1}^{b_{1}}X_{2}^{b_{2}},f)GX_{1}^{b_{1}}X_{2}^{b_{2}}\otimes f+ \\
&&\sum_{b_{1},b_{2},d,e_{1},e_{2}=0}^{1}B(hg\otimes h^{\prime
}g;GX_{1}^{b_{1}}X_{2}^{b_{2}},f)\alpha X_{1}^{b_{1}}X_{2}^{b_{2}}\otimes f
\end{eqnarray*}%
so that%
\begin{eqnarray*}
(G\otimes 1_{H})B(hg\otimes h^{\prime }g)
&=&\sum_{d,e_{1},e_{2}=0}^{1}B(hg\otimes h^{\prime }g;1_{A},f)G\otimes f \\
&&\sum_{d,e_{1},e_{2}=0}^{1}B(hg\otimes h^{\prime }g;X_{1},f)GX_{1}\otimes f+
\\
&&\sum_{d,e_{1},e_{2}=0}^{1}B(hg\otimes h^{\prime }g;X_{2},f)GX_{2}\otimes f+
\\
&&\sum_{d,e_{1},e_{2}=0}^{1}B(hg\otimes h^{\prime
}g;X_{1}X_{2},f)GX_{1}X_{2}\otimes f+ \\
&&\sum_{d,e_{1},e_{2}=0}^{1}B(hg\otimes h^{\prime }g;G,f)\alpha \otimes f+ \\
&&\sum_{d,e_{1},e_{2}=0}^{1}B(hg\otimes h^{\prime }g;GX_{1},f)\alpha
X_{1}\otimes f+ \\
&&\sum_{d,e_{1},e_{2}=0}^{1}B(hg\otimes h^{\prime }g;GX_{2},f)\alpha
X_{2}\otimes f+ \\
&&\sum_{d,e_{1},e_{2}=0}^{1}B(hg\otimes h^{\prime }g;GX_{1}X_{2},f)\alpha
X_{1}X_{2}\otimes f
\end{eqnarray*}%
By comparing right and left side of the equation, we get%
\begin{eqnarray*}
&&\text{1}\sum_{d,e_{1},e_{2}=0}^{1}\left[
\begin{array}{c}
\alpha B(h\otimes h^{\prime };G,f) \\
+\gamma _{1}B(h\otimes h^{\prime };X_{1},f)+\gamma _{2}B(h\otimes h^{\prime
};X_{2},f)%
\end{array}%
\right] 1_{A}\otimes f \\
&=&\sum_{d,e_{1},e_{2}=0}^{1}B(hg\otimes h^{\prime }g;G,f)\alpha
1_{A}\otimes f
\end{eqnarray*}%
\begin{eqnarray*}
&&\text{2}\sum_{d,e_{1},e_{2}=0}^{1}\left[
\begin{array}{c}
B(h\otimes h^{\prime };1_{A},f) \\
+\gamma _{1}B(h\otimes h^{\prime };GX_{1},f)+\gamma _{2}B(h\otimes h^{\prime
};GX_{2},f)%
\end{array}%
\right] G\otimes f \\
&=&\sum_{d,e_{1},e_{2}=0}^{1}B(hg\otimes h^{\prime }g;1_{A},f)G\otimes f
\end{eqnarray*}%
\begin{eqnarray*}
&&\text{3}\sum_{d,e_{1},e_{2}=0}^{1}\left[ \gamma _{2}B(h\otimes h^{\prime
};X_{1}X_2,f)-\alpha B(h\otimes h^{\prime };GX_{1},f)\right] X_{1}\otimes f
\\
&=&\sum_{d,e_{1},e_{2}=0}^{1}B(hg\otimes h^{\prime }g;GX_{1},f)\alpha
X_{1}\otimes f
\end{eqnarray*}%
\begin{eqnarray*}
&&\text{4}\sum_{d,e_{1},e_{2}=0}^{1}\left[ \left( -\gamma _{1}\right)
B(h\otimes h^{\prime };X_1X_{2},f)+\left( -\alpha \right) B(h\otimes
h^{\prime };GX_{2},f)\right] X_{2}\otimes f \\
&&\sum_{d,e_{1},e_{2}=0}^{1}B(hg\otimes h^{\prime }g;GX_{2},f)\alpha
X_{2}\otimes f
\end{eqnarray*}%
\begin{eqnarray*}
&&\text{5}\sum_{d,e_{1},e_{2}=0}^{1}B(h\otimes h^{\prime
};GX_{1}X_{2},f)\alpha X_{1}X_{2}\otimes f \\
&=&\sum_{d,e_{1},e_{2}=0}^{1}B(hg\otimes h^{\prime }g;GX_{1}X_{2},f)\alpha
X_{1}X_{2}\otimes f
\end{eqnarray*}%
\begin{eqnarray*}
\text{6} &&\sum_{d,e_{1},e_{2}=0}^{1}\left[ -B(h\otimes h^{\prime
};X_{1},f)+ \gamma _{2} B(h\otimes h^{\prime };GX_{1}X_2,f)\right]
GX_{1}\otimes f \\
&=&\sum_{d,e_{1},e_{2}=0}^{1}B(hg\otimes h^{\prime }g;X_{1},f)GX_{1}\otimes f
\end{eqnarray*}%
\begin{eqnarray*}
\text{7} &&+\sum_{d,e_{1},e_{2}=0}^{1}\left[ - B(h\otimes h^{\prime
};X_{2},f) -\gamma _{1}B(h\otimes h^{\prime };GX_1X_{2},f)\right]
GX_{2}\otimes f \\
&=&\sum_{d,e_{1},e_{2}=0}^{1}B(hg\otimes h^{\prime }g;X_{2},f)GX_{2}\otimes f
\end{eqnarray*}%
\begin{eqnarray*}
&&\text{8}\sum_{d,e_{1},e_{2}=0}^{1}B(h\otimes h^{\prime
};X_{1}X_{2},f)GX_{1}X_{2}\otimes f \\
&=&\sum_{d,e_{1},e_{2}=0}^{1}B(hg\otimes h^{\prime
}g;X_{1}X_{2},f)GX_{1}X_{2}\otimes f
\end{eqnarray*}

Therefore $\left( \ref{eq.a}\right) $ for $d=G$ is equivalent to the
following equality to be satisfied for every $f\in H$.

\begin{eqnarray}
\text{1}%
\begin{array}{c}
\alpha B(h\otimes h^{\prime };G,f)+ \\
+\gamma _{1}B(h\otimes h^{\prime };X_{1},f)+\gamma _{2}B(h\otimes h^{\prime
};X_{2},f)%
\end{array}
&=&\alpha B(hg\otimes h^{\prime }g;G,f)  \label{GF1} \\
\text{2}%
\begin{array}{c}
B(h\otimes h^{\prime };1_{A},f) \\
+\gamma _{1}B(h\otimes h^{\prime };GX_{1},f)+\gamma _{2}B(h\otimes h^{\prime
};GX_{2},f)%
\end{array}
&=&B(hg\otimes h^{\prime }g;1_{A},f)  \label{GF2} \\
\text{3}\gamma _{2}B(h\otimes h^{\prime };X_{1}X_{2},f)-\alpha B(h\otimes
h^{\prime };GX_{1},f) &=&\alpha B(hg\otimes h^{\prime }g;GX_{1},f)
\label{GF3} \\
\text{4}-\gamma _{1}B(h\otimes h^{\prime };X_{1}X_{2},f)-\alpha B(h\otimes
h^{\prime };GX_{2},f) &=&\alpha B(hg\otimes h^{\prime }g;GX_{2},f)
\label{GF4} \\
\text{5}\alpha B(h\otimes h^{\prime };GX_{1}X_{2},f) &=&\alpha B(hg\otimes
h^{\prime }g;GX_{1}X_{2},f)  \label{GF5} \\
\text{6}-B(h\otimes h^{\prime };X_{1},f)+\gamma _{2}B(h\otimes h^{\prime
};GX_{1}X_{2},f) &=&B(hg\otimes h^{\prime }g;X_{1},f)  \label{GF6} \\
\text{7}-B(h\otimes h^{\prime };X_{2},f)-\gamma _{1}B(h\otimes h^{\prime
};GX_{1}X_{2},f) &=&B(hg\otimes h^{\prime }g;X_{2},f)  \label{GF7} \\
\text{8}B(h\otimes h^{\prime };X_{1}X_{2},f) &=&B(hg\otimes h^{\prime
}g;X_{1}X_{2},f)  \label{GF8}
\end{eqnarray}

\subsubsection{Case $g\otimes 1_{H}$}

\paragraph{\textbf{Equality }$\left( \mathbf{\protect\ref{GF1}}\right) $}

rewrites as%
\begin{equation*}
\begin{array}{c}
\alpha B(g\otimes 1_{H};G,f)+ \\
+\gamma _{1}B(g\otimes 1_{H};X_{1},f)+\gamma _{2}B(g\otimes 1_{H};X_{2},f)%
\end{array}%
=\alpha B(1_{H}\otimes g;G,f).
\end{equation*}

\subparagraph{Case $f=1_{H}$}

We get%
\begin{equation*}
\begin{array}{c}
\alpha B(g\otimes 1_{H};G,1_{H})+ \\
+\gamma _{1}B(g\otimes 1_{H};X_{1},1_{H})+\gamma _{2}B(g\otimes
1_{H};X_{2},1_{H})%
\end{array}%
=\alpha B(1_{H}\otimes g;G,1_{H}).
\end{equation*}

Taking in account $\left( \ref{eq.10}\right) $ this rewrites as%
\begin{equation*}
\begin{array}{c}
\alpha B(g\otimes 1_{H};G,1_{H})+ \\
+\gamma _{1}B(g\otimes 1_{H};X_{1},1_{H})+\gamma _{2}B(g\otimes
1_{H};X_{2},1_{H})%
\end{array}%
=\alpha B(g\otimes 1_{H};G,1_{H}).
\end{equation*}%
i.e.%
\begin{equation*}
\gamma _{1}B(g\otimes 1_{H};X_{1},1_{H})+\gamma _{2}B(g\otimes
1_{H};X_{2},1_{H})=0.
\end{equation*}%
By taking in account the form of the element, we get%
\begin{equation}
\gamma _{1}B\left( g\otimes 1_{H};1_{A},x_{1}\right) +\gamma _{2}B\left(
g\otimes 1_{H};1_{A},x_{2}\right) =0.  \label{G,g,GF1,1H}
\end{equation}

\subparagraph{Case $f=x_{1}x_{2}$}

\begin{equation*}
\begin{array}{c}
\alpha B(g\otimes 1_{H};G,x_{1}x_{2})+ \\
+\gamma _{1}B(g\otimes 1_{H};X_{1},x_{1}x_{2})+\gamma _{2}B(g\otimes
1_{H};X_{2},x_{1}x_{2})%
\end{array}%
=\alpha B(1_{H}\otimes g;G,x_{1}x_{2}).
\end{equation*}

Taking in account $\left( \ref{eq.10}\right) $ this rewrites as%
\begin{equation*}
\begin{array}{c}
\alpha B(g\otimes 1_{H};G,x_{1}x_{2})+ \\
+\gamma _{1}B(g\otimes 1_{H};X_{1},x_{1}x_{2})+\gamma _{2}B(g\otimes
1_{H};X_{2},x_{1}x_{2})%
\end{array}%
=\alpha B(g\otimes 1_{H};G,x_{1}x_{2}).
\end{equation*}%
i.e.

\begin{equation*}
\gamma _{1}B(g\otimes 1_{H};X_{1},x_{1}x_{2})+\gamma _{2}B(g\otimes
1_{H};X_{2},x_{1}x_{2})=0.
\end{equation*}

By taking in account the form of the element, we get a trivial equality.

\subparagraph{Case $f=gx_{1}$}

\begin{equation*}
\begin{array}{c}
\alpha B(g\otimes 1_{H};G,gx_{1})+ \\
+\gamma _{1}B(g\otimes 1_{H};X_{1},gx_{1})+\gamma _{2}B(g\otimes
1_{H};X_{2},gx_{1})%
\end{array}%
=\alpha B(1_{H}\otimes g;G,gx_{1}).
\end{equation*}

By applying $\left( \ref{eq.10}\right) $ this rewrites as%
\begin{equation*}
\begin{array}{c}
2\alpha B(g\otimes 1_{H};G,gx_{1})+ \\
+\gamma _{1}B(g\otimes 1_{H};X_{1},gx_{1})+\gamma _{2}B(g\otimes
1_{H};X_{2},gx_{1})%
\end{array}%
=0
\end{equation*}%
By taking in account the form of the element, we get
\begin{equation}
2\alpha B(g\otimes 1_{H};G,gx_{1})+\gamma _{2}B(g\otimes
1_{H};1_{A},gx_{1}x_{2})=0  \label{G,g, GF1,gx1}
\end{equation}

\subparagraph{Case $f=gx_{2}$}

\begin{equation*}
\begin{array}{c}
\alpha B(g\otimes 1_{H};G,gx_{2})+ \\
+\gamma _{1}B(g\otimes 1_{H};X_{1},gx_{2})+\gamma _{2}B(g\otimes
1_{H};X_{2},gx_{2})%
\end{array}%
=\alpha B(1_{H}\otimes g;G,gx_{2}).
\end{equation*}

By applying $\left( \ref{eq.10}\right) $ this rewrites as%
\begin{equation*}
\begin{array}{c}
2\alpha B(g\otimes 1_{H};G,gx_{2})+ \\
+\gamma _{1}B(g\otimes 1_{H};X_{1},gx_{2})+\gamma _{2}B(g\otimes
1_{H};X_{2},gx_{2})%
\end{array}%
=0.
\end{equation*}%
By taking in account the form of the element, we get%
\begin{equation}
2\alpha B(g\otimes 1_{H};G,gx_{2})-\gamma _{1}B\left( g\otimes
1_{H};1_{A},gx_{1}x_{2}\right) =0.  \label{G,g, GF1,gx2}
\end{equation}

\paragraph{\textbf{Equality }$\left( \mathbf{\protect\ref{GF2}}\right) $}

rewrites as%
\begin{equation*}
\begin{array}{c}
B(g\otimes 1_{H};1_{A},f) \\
+\gamma _{1}B(g\otimes 1_{H};GX_{1},f)+\gamma _{2}B(g\otimes 1_{H};GX_{2},f)%
\end{array}%
=B(1_{H}\otimes g;1_{A},f)
\end{equation*}

\subparagraph{Case $f=g$%
\protect\begin{equation*}
\protect\begin{array}{c}
B(g\otimes 1_{H};1_{A},g)\protect \\
+\protect\gamma _{1}B(g\otimes 1_{H};GX_{1},g)+\protect\gamma _{2}B(g\otimes
1_{H};GX_{2},g)%
\protect\end{array}%
=B(1_{H}\otimes g;1_{A},g)
\protect\end{equation*}%
}

By applying $\left( \ref{eq.10}\right) $ this rewrites as%
\begin{equation*}
\gamma _{1}B(g\otimes 1_{H};GX_{1},g)+\gamma _{2}B(g\otimes
1_{H};GX_{2},g)=0.
\end{equation*}%
By taking in account the form of the element, we get%
\begin{equation}
\gamma _{1}B\left( g\otimes 1_{H};G,gx_{1}\right) +\gamma _{2}B\left(
g\otimes 1_{H};G,gx_{2}\right) =0.  \label{G,g, GF2,g}
\end{equation}

\subparagraph{Case $f=x_{1}$}

\begin{equation*}
\begin{array}{c}
B(g\otimes 1_{H};1_{A},x_{1}) \\
+\gamma _{1}B(g\otimes 1_{H};GX_{1},x_{1})+\gamma _{2}B(g\otimes
1_{H};GX_{2},x_{1})%
\end{array}%
=B(1_{H}\otimes g;1_{A},x_{1})
\end{equation*}

By applying $\left( \ref{eq.10}\right) $ this rewrites as%
\begin{equation*}
\begin{array}{c}
2B(g\otimes 1_{H};1_{A},x_{1}) \\
+\gamma _{1}B(g\otimes 1_{H};GX_{1},x_{1})+\gamma _{2}B(g\otimes
1_{H};GX_{2},x_{1})%
\end{array}%
=0.
\end{equation*}%
By taking in account the form of the element, we get%
\begin{equation}
2B(g\otimes 1_{H};1_{A},x_{1})-\gamma _{2}B\left( g\otimes
1_{H};G,x_{1}x_{2}\right) =0.  \label{G,g, GF2,x1}
\end{equation}

\subparagraph{Case $f=x_{2}$}

\begin{equation*}
\begin{array}{c}
B(g\otimes 1_{H};1_{A},x_{2}) \\
+\gamma _{1}B(g\otimes 1_{H};GX_{1},x_{2})+\gamma _{2}B(g\otimes
1_{H};GX_{2},x_{2})%
\end{array}%
=B(1_{H}\otimes g;1_{A},x_{2})
\end{equation*}

By applying $\left( \ref{eq.10}\right) $ this rewrites as%
\begin{equation*}
\begin{array}{c}
2B(g\otimes 1_{H};1_{A},x_{2}) \\
+\gamma _{1}B(g\otimes 1_{H};GX_{1},x_{2})+\gamma _{2}B(g\otimes
1_{H};GX_{2},x_{2})%
\end{array}%
=0.
\end{equation*}%
By taking in account the form of the element, we get%
\begin{equation}
2B(g\otimes 1_{H};1_{A},x_{2})+\gamma _{1}B\left( g\otimes
1_{H};G,x_{1}x_{2}\right) =0.  \label{G,g, GF2,x2}
\end{equation}

\subparagraph{Case $f=gx_{1}x_{2}$}

\begin{equation*}
\begin{array}{c}
B(g\otimes 1_{H};1_{A},gx_{1}x_{2}) \\
+\gamma _{1}B(g\otimes 1_{H};GX_{1},gx_{1}x_{2})+\gamma _{2}B(g\otimes
1_{H};GX_{2},gx_{1}x_{2})%
\end{array}%
=B(1_{H}\otimes g;1_{A},gx_{1}x_{2})
\end{equation*}%
By applying $\left( \ref{eq.10}\right) $ this rewrites as%
\begin{equation*}
\gamma _{1}B(g\otimes 1_{H};GX_{1},gx_{1}x_{2})+\gamma _{2}B(g\otimes
1_{H};GX_{2},gx_{1}x_{2})=0.
\end{equation*}%
By taking in account the form of the element, we get a trivial equality.


\paragraph{\textbf{Equality }$\left( \mathbf{\protect\ref{GF3}}\right) $}

rewrites as

\begin{equation*}
\gamma _{2}B(g\otimes 1_{H};X_{1}X_2,f)-\alpha B(g\otimes
1_{H};GX_{1},f)=\alpha B(1_{H}\otimes g;GX_{1},f)
\end{equation*}

\subparagraph{Case $f=g$}

\begin{equation*}
\gamma _{2}B(g\otimes 1_{H};X_{1}X_2,g)-\alpha B(g\otimes
1_{H};GX_{1},g)=\alpha B(1_{H}\otimes g;GX_{1},g)
\end{equation*}%
By applying $\left( \ref{eq.10}\right) $ this rewrites as%
\begin{equation*}
\gamma _{2}B(g\otimes 1_{H};X_{1}X_2,g)-2\alpha B(g\otimes 1_{H};GX_{1},g)=0.
\end{equation*}%
By taking in account the form of the element, we get

\begin{equation*}
\gamma _{2}B(g\otimes 1_{H};1_H,gx_1x_2)+2\alpha B(g\otimes 1_{H};G,gx_1)=0.
\end{equation*}

This have been obtain in Case $gx_1$ of \ref{GF1}.

\subparagraph{Case $f=x_{2}$}

\begin{equation*}
\gamma _{2}B(g\otimes 1_{H};X_{1}X_2,x_{2})-\alpha B(g\otimes
1_{H};GX_{1},x_{2})=\alpha B(1_{H}\otimes g;GX_{1},x_{2})
\end{equation*}

By applying $\left( \ref{eq.10}\right) $ this rewrites as%
\begin{equation*}
\gamma _{2}B(g\otimes 1_{H};X_{1}X_{2},x_{2})=0.
\end{equation*}%
By taking in account the form of the element, we get nothing new.

\paragraph{\textbf{Equality }$\left( \mathbf{\protect\ref{GF4}}\right) $}

rewrites as
\begin{equation*}
-\gamma _{1}B(g\otimes 1_{H};X_{1}X_{2},f)-\alpha B(g\otimes
1_{H};GX_{2},f)=\alpha B(1_{H}\otimes g;GX_{2},f)
\end{equation*}

\subparagraph{Case $f=g$}

\begin{equation*}
-\gamma _{1}B(g\otimes 1_{H};X_1X_{2},g)-\alpha B(g\otimes
1_{H};GX_{2},g)=\alpha B(1_{H}\otimes g;GX_{2},g)
\end{equation*}%
By applying $\left( \ref{eq.10}\right) $ this rewrites as%
\begin{equation*}
-\gamma _{1}B(g\otimes 1_{H};X_1X_{2},g)-2\alpha B(g\otimes 1_{H};GX_{2},g)=0
\end{equation*}

By taking in account the form of the element, we get

\begin{equation*}
-\gamma _{1}B(g\otimes 1_{H};1_H,gx_1x_2)+2\alpha B(g\otimes
1_{H};g,GX_{2})=0
\end{equation*}%
which is $(\ref{G,g, GF1,gx2} ).$

\subparagraph{Case $f=x_{1}$}

\begin{equation*}
-\gamma _{1}B(g\otimes 1_{H};X_{2},x_{1})-\alpha B(g\otimes
1_{H};GX_{2},x_{1})=\alpha B(1_{H}\otimes g;GX_{2},x_{1})
\end{equation*}%
By applying $\left( \ref{eq.10}\right) $ this rewrites as%
\begin{equation*}
\gamma _{1}B(g\otimes 1_{H};X_{1}X_{2},x_{1})=0
\end{equation*}%
Since $B(g\otimes 1_{H};X_{1}X_{2},x_{1})=0,$ we get nothing new.

\paragraph{\textbf{Equality}$\left( \mathbf{\ \protect\ref{GF5}}\right) $}

rewrites as
\begin{equation*}
\alpha B(g\otimes 1_{H};GX_{1}X_{2},f)=\alpha B(1_{H}\otimes
g;GX_{1}X_{2},f).
\end{equation*}

\subparagraph{Case $f=1_{H}$%
\protect\begin{equation*}
\protect\alpha B(g\otimes 1_{H};GX_{1}X_{2},1_{H})=\protect\alpha %
B(1_{H}\otimes g;GX_{1}X_{2},1_{H}).
\protect\end{equation*}%
}

By applying $\left( \ref{eq.10}\right) $ this equality is trivial.

\paragraph{\textbf{Equality }$\left( \mathbf{\protect\ref{GF6}}\right) $%
\textbf{\ \ }}

rewrites as

\begin{equation*}
- B(g\otimes 1_{H};X_{1},f)+ \gamma _{2} B(g\otimes
1_{H};GX_{1}X_2,f)=B(1_{H}\otimes g;X_{1},f)
\end{equation*}

\subparagraph{Case $f=1_{H}$%
\protect\begin{equation*}
- B(g\otimes 1_{H};X_{1},1_{H})+\protect\gamma _{2} B(g\otimes
1_{H};GX_{1}X_2,1_{H})=B(1_{H}\otimes g;X_{1},1_{H})
\protect\end{equation*}%
}

By applying $\left( \ref{eq.10}\right) $ this rewrites as%
\begin{equation*}
- 2B(g\otimes 1_{H};X_{1},1_{H})+\gamma _{2} B(g\otimes
1_{H};GX_{1}X_2,1_{H})=0
\end{equation*}%
By taking in account the form of the element, we get%
\begin{equation*}
-2B(g\otimes 1_{H};1_A,x_1)+ \gamma _{2} B(g\otimes 1_{H};G,x_1x_2)=0
\end{equation*}

This is $(\ref{G,g, GF2,x1})$.

\subparagraph{Case $f=gx_{2}$}

\begin{equation*}
-B(g\otimes 1_{H};X_{1},gx_{2})+\gamma _{2}B(g\otimes
1_{H};GX_{1}X_{2},gx_{2})=B(1_{H}\otimes g;X_{1},gx_{2})
\end{equation*}%
By applying $\left( \ref{eq.10}\right) $ this rewrites as%
\begin{equation*}
\gamma _{2}B(g\otimes 1_{H};GX_{1}X_{2},gx_{2})=0
\end{equation*}%
By taking in account the form of the element, we get nothing new.

\paragraph{\textbf{Equality }$\left( \mathbf{\protect\ref{GF7}}\right) $}

rewrites as

\begin{equation*}
-B(g\otimes 1_{H};X_{2},f)-\gamma _{1}B(g\otimes
1_{H};GX_{1}X_{2},f)=B(1_{H}\otimes g;X_{2},f)
\end{equation*}

\subparagraph{Case $f=1_{H}$}

\begin{equation*}
-B(g\otimes 1_{H};X_{2},1_{H}) -\gamma _{1} B(g\otimes
1_{H};GX_1X_{2},1_{H})=B(1_{H}\otimes g;X_{2},1_{H})
\end{equation*}%
By applying $\left( \ref{eq.10}\right) $ this rewrites as%
\begin{equation*}
- 2B(g\otimes 1_{H};X_{2},1_{H}) -\gamma_{1} B(g\otimes
1_{H};GX_1X_{2},1_{H})=0.
\end{equation*}%
By taking in account the form of the element, we get%
\begin{equation*}
- 2B(g\otimes 1_{H};1_A,x_2) -\gamma_{1} B(g\otimes 1_{H};G,x_1x_2)=0.
\end{equation*}
This is $(\ref{G,g, GF2,x2})$.

\subparagraph{Case $f=gx_{1}$}

\begin{equation*}
-B(g\otimes 1_{H};X_{2},gx_{1}) -\gamma _{1} B(g\otimes
1_{H};GX_{2},gx_{1})=B(1_{H}\otimes g;X_{2},gx_{1})
\end{equation*}%
By applying $\left( \ref{eq.10}\right) $ this rewrites as%
\begin{equation*}
\gamma _{1}B(g\otimes 1_{H};GX_1X_{2},gx_{1})=0
\end{equation*}%
By taking in account the form of the element, we get nothing new.

\paragraph{\textbf{Equality }$\left( \mathbf{\protect\ref{GF8}}\right) $}

\begin{equation*}
B(g\otimes 1_{H};X_{1}X_{2},f)=B(1_{H}\otimes g;X_{1}X_{2},f)
\end{equation*}%
By taking in account the form of the elements we obtain%
\begin{equation*}
B(g\otimes 1_{H};X_{1}X_{2},g)=B(1_{H}\otimes g;X_{1}X_{2},g)
\end{equation*}%
By applying $\left( \ref{eq.10}\right) $ we get nothing new.

\subsubsection{Case $x_{1}\otimes 1_{H}$}

\paragraph{\textbf{Equality }$\left( \mathbf{\protect\ref{GF1}}\right) $}

rewrites as

\begin{equation*}
\begin{array}{c}
\alpha B(x_{1}\otimes 1_{H};G,f)+ \\
+\gamma _{1}B(x_{1}\otimes 1_{H};X_{1},f)+\gamma _{2}B(x_{1}\otimes
1_{H};X_{2},f)%
\end{array}%
=\alpha B(x_{1}g\otimes g;G,f)
\end{equation*}%
i.e.%
\begin{equation*}
\begin{array}{c}
\alpha B(x_{1}\otimes 1_{H};G,f)+ \\
+\gamma _{1}B(x_{1}\otimes 1_{H};X_{1},f)+\gamma _{2}B(x_{1}\otimes
1_{H};X_{2},f)%
\end{array}%
=-\alpha B(gx_{1}\otimes g;G,f)
\end{equation*}

\subparagraph{Case $f=g$%
\protect\begin{equation*}
\protect\begin{array}{c}
\protect\alpha B(x_{1}\otimes 1_{H};G,g)+\protect \\
+\protect\gamma _{1}B(x_{1}\otimes 1_{H};X_{1},g)+\protect\gamma %
_{2}B(x_{1}\otimes 1_{H};X_{2},g)%
\protect\end{array}%
=-\protect\alpha B(gx_{1}\otimes g;G,g).
\protect\end{equation*}%
}

By applying $\left( \ref{eq.10}\right) $ this rewrites as%
\begin{equation*}
\begin{array}{c}
\alpha B(x_{1}\otimes 1_{H};G,g)+ \\
+\gamma _{1}B(x_{1}\otimes 1_{H};X_{1},g)+\gamma _{2}B(x_{1}\otimes
1_{H};X_{2},g)%
\end{array}%
=-\alpha B(x_{1}\otimes 1_{H};G,g).
\end{equation*}%
By taking in account the form of the elements we obtain%
\begin{equation}
\begin{array}{c}
2\alpha B(x_{1}\otimes 1_{H};G,g)+ \\
+\gamma _{1}\left[ -B(g\otimes 1_{H};1_{A},g)-B(x_{1}\otimes
1_{H};1_{A},gx_{1})\right] -\gamma _{2}B(x_{1}\otimes 1_{H};1_{A},gx_{2})%
\end{array}%
=0.  \label{G,x1, GF1,g}
\end{equation}

\subparagraph{Case $f=x_{1}$}

\begin{equation*}
\begin{array}{c}
\alpha B(x_{1}\otimes 1_{H};G,x_{1})+ \\
+\gamma _{1}B(x_{1}\otimes 1_{H};X_{1},x_{1})+\gamma _{2}B(x_{1}\otimes
1_{H};X_{2},x_{1})%
\end{array}%
=-\alpha B(gx_{1}\otimes g;G,x_{1})
\end{equation*}%
By applying $\left( \ref{eq.10}\right) $ this rewrites as%
\begin{equation*}
\gamma _{1}B(x_{1}\otimes 1_{H};X_{1},x_{1})+\gamma _{2}B(x_{1}\otimes
1_{H};X_{2},x_{1})=0.
\end{equation*}%
By taking in account the form of the elements we obtain%
\begin{equation}
\gamma _{1}B(g\otimes 1_{H};1_{A},x_{1})+\gamma _{2}B(x_{1}\otimes
1_{H};1_{A},x_{1}x_{2})=0.  \label{G,x1, GF1,x1}
\end{equation}

\subparagraph{Case $f=x_{2}$}

\begin{equation*}
\begin{array}{c}
\alpha B(x_{1}\otimes 1_{H};G,x_{2})+ \\
+\gamma _{1}B(x_{1}\otimes 1_{H};X_{1},x_{2})+\gamma _{2}B(x_{1}\otimes
1_{H};X_{2},x_{2})%
\end{array}%
=-\alpha B(gx_{1}\otimes g;G,x_{2}).
\end{equation*}%
By applying $\left( \ref{eq.10}\right) $ this rewrites as%
\begin{equation*}
\gamma _{1}B(x_{1}\otimes 1_{H};X_{1},x_{2})+\gamma _{2}B(x_{1}\otimes
1_{H};X_{2},x_{2})=0.
\end{equation*}%
By taking in account the form of the element we obtain%
\begin{equation}
\gamma _{1}[-B(g\otimes 1_{H};1_{A},x_{2})+B(x_{1}\otimes
1_{H};1_{A},x_{1}x_{2})]=0.  \label{G,x1, GF1,x2}
\end{equation}

\subparagraph{Case $f=gx_{1}x_{2}$}

\begin{equation*}
\begin{array}{c}
\alpha B(x_{1}\otimes 1_{H};G,gx_{1}x_{2})+ \\
+\gamma _{1}B(x_{1}\otimes 1_{H};X_{1},gx_{1}x_{2})+\gamma
_{2}B(x_{1}\otimes 1_{H};X_{2},gx_{1}x_{2})%
\end{array}%
=\alpha B(x_{1}g\otimes g;G,gx_{1}x_{2}).
\end{equation*}%
By applying $\left( \ref{eq.10}\right) $ this rewrites as%
\begin{equation*}
\begin{array}{c}
2\alpha B(x_{1}\otimes 1_{H};G,gx_{1}x_{2})+ \\
+\gamma _{1}B(x_{1}\otimes 1_{H};X_{1},gx_{1}x_{2})+\gamma
_{2}B(x_{1}\otimes 1_{H};X_{2},gx_{1}x_{2})%
\end{array}%
=0.
\end{equation*}%
By taking in account the form of the element we obtain%
\begin{equation}
-\gamma _{1}B(g\otimes 1_{H};1_{A},gx_{1}x_{2})+2\alpha B(x_{1}\otimes
1_{H};G,gx_{1}x_{2})=0.  \label{G,x1, GF1,gx1x2}
\end{equation}

\paragraph{\textbf{Equality (\protect\ref{GF2})}}

rewrites as
\begin{equation*}
\begin{array}{c}
B(x_{1}\otimes 1_{H};1_{A},f) \\
+\gamma _{1}B(x_{1}\otimes 1_{H};GX_{1},f)+\gamma _{2}B(x_{1}\otimes
1_{H};GX_{2},f)%
\end{array}%
=B(x_{1}g\otimes g;1_{A},f)
\end{equation*}%
i.e.%
\begin{equation*}
\begin{array}{c}
B(x_{1}\otimes 1_{H};1_{A},f) \\
+\gamma _{1}B(x_{1}\otimes 1_{H};GX_{1},f)+\gamma _{2}B(x_{1}\otimes
1_{H};GX_{2},f)%
\end{array}%
=-B(gx_{1}\otimes g;1_{A},f).
\end{equation*}

\subparagraph{Case $f=1_{H}$}

\begin{equation*}
\begin{array}{c}
B(x_{1}\otimes 1_{H};1_{A},1_{H}) \\
+\gamma _{1}B(x_{1}\otimes 1_{H};GX_{1},1_{H})+\gamma _{2}B(x_{1}\otimes
1_{H};GX_{2},1_{H})%
\end{array}%
=-B(gx_{1}\otimes g;1_{A},1_{H}).
\end{equation*}%
By applying $\left( \ref{eq.10}\right) $ this rewrites as%
\begin{equation*}
\begin{array}{c}
2B(x_{1}\otimes 1_{H};1_{A},1_{H}) \\
+\gamma _{1}B(x_{1}\otimes 1_{H};GX_{1},1_{H})+\gamma _{2}B(x_{1}\otimes
1_{H};GX_{2},1_{H})%
\end{array}%
=0.
\end{equation*}%
By taking in account the form of the element we obtain%
\begin{equation}
\begin{array}{c}
2B(x_{1}\otimes 1_{H};1_{A},1_{H})+ \\
+\gamma _{1}\left[ B(g\otimes 1_{H};G,1_{H})+B(x_{1}\otimes 1_{H};G,x_{1})%
\right] +\gamma _{2}B(x_{1}\otimes 1_{H};G,x_{2})%
\end{array}%
=0.  \label{G,x1, GF2,1H}
\end{equation}

\subparagraph{Case $f=x_{1}x_{2}$}

\begin{equation*}
\begin{array}{c}
B(x_{1}\otimes 1_{H};1_{A},x_{1}x_{2}) \\
+\gamma _{1}B(x_{1}\otimes 1_{H};GX_{1},x_{1}x_{2})+\gamma
_{2}B(x_{1}\otimes 1_{H};GX_{2},x_{1}x_{2})%
\end{array}%
=-B(gx_{1}\otimes g;1_{A},x_{1}x_{2}).
\end{equation*}%
By applying $\left( \ref{eq.10}\right) $ this rewrites as%
\begin{equation*}
\begin{array}{c}
2B(x_{1}\otimes 1_{H};1_{A},x_{1}x_{2}) \\
+\gamma _{1}B(x_{1}\otimes 1_{H};GX_{1},x_{1}x_{2})+\gamma
_{2}B(x_{1}\otimes 1_{H};GX_{2},x_{1}x_{2})%
\end{array}%
=0.
\end{equation*}%
By taking in account the form of the element we obtain%
\begin{equation}
2B(x_{1}\otimes 1_{H};1_{A},x_{1}x_{2})+\gamma _{1}B(g\otimes
1_{H};G,x_{1}x_{2})=0.  \label{G,x1, GF2,x1x2}
\end{equation}

\subparagraph{Case $f=gx_{1}$}

\begin{equation*}
\begin{array}{c}
B(x_{1}\otimes 1_{H};1_{A},gx_{1}) \\
+\gamma _{1}B(x_{1}\otimes 1_{H};GX_{1},gx_{1})+\gamma _{2}B(x_{1}\otimes
1_{H};GX_{2},gx_{1})%
\end{array}%
=-B(gx_{1}\otimes g;1_{A},gx_{1})
\end{equation*}%
By applying $\left( \ref{eq.10}\right) $ this rewrites as%
\begin{equation*}
\gamma _{1}B(x_{1}\otimes 1_{H};GX_{1},gx_{1})+\gamma _{2}B(x_{1}\otimes
1_{H};GX_{2},gx_{1})=0.
\end{equation*}%
By taking in account the form of the elements we obtain%
\begin{equation}
\gamma _{1}B\left( g\otimes 1_{H};G,gx_{1}\right) +\gamma _{2}B(x_{1}\otimes
1_{H};G,gx_{1}x_{2})=0.  \label{G,x1, GF2,gx1}
\end{equation}

\subparagraph{Case $f=gx_{2}$}

\begin{equation*}
\begin{array}{c}
B(x_{1}\otimes 1_{H};1_{A},gx_{2}) \\
+\gamma _{1}B(x_{1}\otimes 1_{H};GX_{1},gx_{2})+\gamma _{2}B(x_{1}\otimes
1_{H};GX_{2},gx_{2})%
\end{array}%
=-B(gx_{1}\otimes g;1_{A},gx_{2}).
\end{equation*}%
By applying $\left( \ref{eq.10}\right) $ this rewrites as%
\begin{equation*}
\gamma _{1}B(x_{1}\otimes 1_{H};GX_{1},gx_{2})+\gamma _{2}B(x_{1}\otimes
1_{H};GX_{2},gx_{2})=0.
\end{equation*}%
By taking in account the form of the elements we obtain%
\begin{equation}
\gamma _{1}[B\left( g\otimes 1_{H};G,gx_{2}\right) -B(x_{1}\otimes
1_{H};G,gx_{1}x_{2})]=0.  \label{G,x1,GF2,gx2}
\end{equation}

\paragraph{\textbf{Equality (\protect\ref{GF3})}}

\begin{equation*}
\gamma _{2}B(x_{1}\otimes 1_{H};X_{1}X_{2},f)-\alpha B(x_{1}\otimes
1_{H};GX_{1},f)=\alpha B(x_{1}g\otimes g;GX_{1},f)
\end{equation*}%
i.e.%
\begin{equation*}
\gamma _{2}B(x_{1}\otimes 1_{H};X_{1}X_{2},f)-\alpha B(x_{1}\otimes
1_{H};GX_{1},f)=-\alpha B(gx_{1}\otimes g;GX_{1},f).
\end{equation*}

\subparagraph{Case $f=1_{H}$}

\begin{equation*}
\gamma _{2}B(x_{1}\otimes 1_{H};X_{1}X_{2},1_{H})-\alpha B(x_{1}\otimes
1_{H};GX_{1},1_{H})=-\alpha B(gx_{1}\otimes g;GX_{1},1_{H}).
\end{equation*}%
By applying $\left( \ref{eq.10}\right) $ this rewrites as%
\begin{equation*}
\gamma _{2}B(x_{1}\otimes 1_{H};X_{1}X_{2},1_{H})-\alpha B(x_{1}\otimes
1_{H};GX_{1},1_{H})=-\alpha B(x_{1}\otimes 1_{H};GX_{1},1_{H}).
\end{equation*}%
By taking in account the form of the elements we obtain
\begin{equation}
\gamma _{2}\left[ -B(g\otimes 1_{H};X_{2},1_{H})+B(x_{1}\otimes
1_{H};1_{A},x_{1}x_{2})\right] =0  \label{G,x1,GF3,1H}
\end{equation}

\subparagraph{Case $f=x_1x_2$}

\begin{equation*}
\gamma _{2}B(x_{1}\otimes 1_{H};X_{1}X_2,x_1x_2)-\alpha B(x_{1}\otimes
1_{H};GX_{1},x_1x_2)=-\alpha B(gx_{1}\otimes g;GX_{1},x_1x_2).
\end{equation*}%
By applying $\left( \ref{eq.10}\right) $ this rewrites as%
\begin{equation*}
\gamma _{2}B(g\otimes 1_{H};X_{1}X_2,x_1x_2)=0.
\end{equation*}%
By taking in account the form of the element we get nothing new.

\subparagraph{Case $f=gx_{1}$}

\begin{equation*}
\gamma _{2}B(x_{1}\otimes 1_{H};X_{1}X_{2},gx_{1})-\alpha B(x_{1}\otimes
1_{H};GX_{1},gx_{1})=-\alpha B(gx_{1}\otimes g;GX_{1},gx_{1}).
\end{equation*}%
By applying $\left( \ref{eq.10}\right) $ this rewrites as%
\begin{equation*}
\gamma _{2}B(x_{1}\otimes 1_{H};X_{1}X_{2},gx_{1})-2\alpha B(x_{1}\otimes
1_{H};GX_{1},gx_{1})=0.
\end{equation*}%
By taking in account the form of the element we get
\begin{equation}
-\gamma _{2}B(g\otimes 1_{H};1_{A},gx_{1}x_{2})-2\alpha B(g\otimes
1_{H};G,gx_{1})=0.  \label{G,x1,GF3,gx1}
\end{equation}

\subparagraph{Case $f=gx_{2}$}

\begin{equation*}
\gamma _{2}B(x_{1}\otimes 1_{H};X_{1}X_{2},gx_{2})-\alpha B(x_{1}\otimes
1_{H};GX_{1},gx_{2})=-\alpha B(gx_{1}\otimes g;GX_{1},gx_{2}).
\end{equation*}%
By applying $\left( \ref{eq.10}\right) $ this rewrites as%
\begin{equation*}
\gamma _{2}B(x_{1}\otimes 1_{H};X_{1}X_{2},gx_{2})-2\alpha B(x_{1}\otimes
1_{H};GX_{1},gx_{2})=0.
\end{equation*}%
By taking in account the form of the element we get%
\begin{equation}
\alpha \left[ -B(g\otimes 1_{H};G,gx_{2})+B(x_{1}\otimes 1_{H};G,gx_{1}x_{2})%
\right] =0.  \label{G,x1, GF3,x2}
\end{equation}

\paragraph{\textbf{Equality (\protect\ref{GF4})}}

\begin{equation*}
-\gamma _{1}B(x_{1}\otimes 1_{H};X_{1}X_{2},f)-\alpha B(x_{1}\otimes
1_{H};GX_{2},f)=\alpha B(x_{1}g\otimes g;GX_{2},f)
\end{equation*}%
i.e.%
\begin{equation*}
-\gamma _{1}B(x_{1}\otimes 1_{H};X_{1}X_{2},f)-\alpha B(x_{1}\otimes
1_{H};GX_{2},f)=-\alpha B(gx_{1}\otimes g;GX_{2},f).
\end{equation*}

\subparagraph{Case $f=1_{H}$}

\begin{equation*}
-\gamma _{1}B(x_{1}\otimes 1_{H};X_{1}X_{2},1_{H})-\alpha B(x_{1}\otimes
1_{H};GX_{2},1_{H})=-\alpha B(gx_{1}\otimes g;GX_{2},1_{H}).
\end{equation*}%
By applying $\left( \ref{eq.10}\right) $ this rewrites as%
\begin{equation*}
\gamma _{1}B(x_{1}\otimes 1_{H};X_{1}X_{2},1_{H})=0
\end{equation*}

In view of the form of the element we get%
\begin{equation}
\gamma _{1}[B(g\otimes 1_{H};1_{A},x_{2})+B(x_{1}\otimes
1_{H};1_{A},x_{1}x_{2})]=0.  \label{G,x1, GF4,1H}
\end{equation}

\subparagraph{Case $f=gx_{1}$}

\begin{equation*}
-\gamma _{1}B(x_{1}\otimes 1_{H};X_1X_{2},gx_{1})-\alpha B(x_{1}\otimes
1_{H};GX_{2},gx_{1})=-\alpha B(gx_{1}\otimes g;GX_{2},gx_{1}).
\end{equation*}%
By applying $\left( \ref{eq.10}\right) $ this rewrites as%
\begin{equation*}
\gamma _{1}B(x_{1}\otimes 1_{H};X_1X_{2},gx_{1})+2\alpha B(x_{1}\otimes
1_{H};GX_{2},gx_{1})=0.
\end{equation*}%
By taking in account the form of the element we get%
\begin{equation}
-\gamma _{1} B(g\otimes 1_H; 1_A,gx_1x_2)+2\alpha B(x_{1}\otimes
1_{H};G,gx_{1}x_{2})=0.  \label{G,x1, GF4,gx1}
\end{equation}

\paragraph{\textbf{Equality (\protect\ref{GF5})}}

rewrites as
\begin{equation*}
\alpha B(x_{1}\otimes 1_{H};GX_{1}X_{2},f)=\alpha B(x_{1}g\otimes
g;GX_{1}X_{2},f)
\end{equation*}%
i.e.%
\begin{equation*}
\alpha B(x_{1}\otimes 1_{H};GX_{1}X_{2},f)=-\alpha B(gx_{1}\otimes
g;GX_{1}X_{2},f).
\end{equation*}

\subparagraph{Case $f=g$}

\begin{equation*}
\alpha B(x_{1}\otimes 1_{H};GX_{1}X_{2},g)=-\alpha B(gx_{1}\otimes
g;GX_{1}X_{2},g).
\end{equation*}%
By applying $\left( \ref{eq.10}\right) $ this rewrites as%
\begin{equation*}
\alpha B(x_{1}\otimes 1_{H};GX_{1}X_{2},g)=-\alpha B(x_{1}\otimes
1_{H};GX_{1}X_{2},g)
\end{equation*}%
i.e.%
\begin{equation*}
\alpha B(x_{1}\otimes 1_{H};GX_{1}X_{2},g)=0
\end{equation*}

By taking account the form of the element we obtain
\begin{equation}
\alpha \lbrack -B(g\otimes 1_{H};G,gx_{2})+B(x_{1}\otimes
1_{H};G,gx_{1}x_{2})]=0.  \label{G,x1, GF5,g}
\end{equation}

\subparagraph{Case $f=x_{1}$}

By applying $\left( \ref{eq.10}\right) $ this rewrites as

\begin{equation*}
\alpha B(x_{1}\otimes 1_{H};GX_{1}X_{2},x_{1})=-\alpha B(gx_{1}\otimes
g;GX_{1}X_{2},x_{1}).
\end{equation*}%
By applying $\left( \ref{eq.10}\right) $ we get nothing new.

\paragraph{\textbf{Equality (\protect\ref{GF6})}}

rewrites as
\begin{equation*}
-B(x_{1}\otimes 1_{x_{1}};X_{1},f)+\gamma _{2}B(x_{1}\otimes
1_{H};GX_{1}X_{2},f)=B(x_{1}g\otimes g;X_{1},f)
\end{equation*}%
i.e.%
\begin{equation*}
-B(x_{1}\otimes 1_{H};X_{1},f)+\gamma _{2}B(x_{1}\otimes
1_{H};GX_{1}X_{2},f)=-B(gx_{1}\otimes g;X_{1},f).
\end{equation*}

\subparagraph{Case $f=g$}

\begin{equation*}
-B(x_{1}\otimes 1_{H};X_{1},g)+\gamma _{2}B(x_{1}\otimes
1_{H};GX_{1}X_2,g)=-B(gx_{1}\otimes g;X_{1},g).
\end{equation*}%
By applying $\left( \ref{eq.10}\right) $ this rewrites as%
\begin{equation*}
-B(x_{1}\otimes 1_{H};X_{1},g)+\gamma _{2}B(x_{1}\otimes
1_{H};GX_{1}X_2,g)=-B(x_{1}\otimes g;X_{1},g)
\end{equation*}%
i.e.%
\begin{equation*}
\gamma _{2}B(x_{1}\otimes 1_{H};GX_{1}X_2,g)=0.
\end{equation*}

This is $\ref{G,x1, GF5,g}$

\subparagraph{Case $f=x_{1}$}

\begin{equation*}
-B(x_{1}\otimes 1_{H};X_{1},x_{1})+\gamma _{2}B(x_{1}\otimes
1_{H};GX_{1}X_2,x_{1})=-B(gx_{1}\otimes g;X_{1},x_{1}).
\end{equation*}%
By applying $\left( \ref{eq.10}\right) $ this rewrites as%
\begin{equation*}
-B(x_{1}\otimes 1_{H};X_{1},x_{1})+\gamma _{2}B(x_{1}\otimes
1_{H};GX_{1}X_2,x_{1})=B(x_{1}\otimes 1_{H};X_{1},x_{1}).
\end{equation*}%
By taking in account the form of the element we get

\begin{equation}
2B(g\otimes 1_{H};1_{A},x_{1}) - \gamma _{2}B(g\otimes 1_{H};G,x_{1}x_2)=0.
\label{G,x1,GF6,x1}
\end{equation}

\subparagraph{Case $f=x_{2}$}

\begin{equation*}
-B(x_{1}\otimes 1_{H};X_{1},x_{2})+\gamma _{2}B(x_{1}\otimes
1_{H};GX_{1}X_2,x_{2})=-B(gx_{1}\otimes g;X_{1},x_{2}).
\end{equation*}%
By applying $\left( \ref{eq.10}\right) $ this rewrites as%
\begin{equation*}
-B(x_{1}\otimes 1_{H};X_{1},x_{2})+\gamma _{2}B(x_{1}\otimes
1_{H};GX_{1}X_2,x_{2})=B(x_{1}\otimes 1_{H};X_{1},x_{2}).
\end{equation*}%
By taking in account the form of the element we get

\begin{equation}
-B(g\otimes 1_{H};1_{A},x_{2})+B(x_1\otimes 1_{H};1_A,x_{1}x_{2})=0.
\label{G,x1,GF6,x2}
\end{equation}

\subparagraph{Case $f=gx_{1}x_{2}$}

\begin{equation*}
-B(x_{1}\otimes 1_{H};X_{1},gx_{1}x_{2})+\gamma _{2}B(x_{1}\otimes
1_{H};GX_{1}X_2,gx_{1}x_{2})=-B(gx_{1}\otimes g;X_{1},gx_{1}x_{2}).
\end{equation*}

By applying $\left( \ref{eq.10}\right) $ this rewrites as%
\begin{equation*}
-B(x_{1}\otimes 1_{H};X_{1},gx_{1}x_{2})+\gamma _{2}B(x_{1}\otimes
1_{H};GX_{1}X_2,gx_{1}x_{2})=-B(x_{1}\otimes 1_{H};X_{1},gx_{1}x_{2})
\end{equation*}%
By taking in account the form of the elements we get nothing new.

\paragraph{\textbf{Equality (\protect\ref{GF7})}}

rewrites as%
\begin{equation*}
-B(x_{1}\otimes 1_{H};X_{2},f)-\gamma _{1}B(x_{1}\otimes
1_{H};GX_{1}X_{2},f)=-B(gx_{1}\otimes g;X_{2},f).
\end{equation*}

\subparagraph{Case $f=g$}

\begin{equation*}
-B(x_{1}\otimes 1_{H};X_{2},g)-\gamma _{1}B(x_{1}\otimes
1_{H};GX_1X_{2},g)=-B(gx_{1}\otimes g;X_{2},g).
\end{equation*}%
By applying $\left( \ref{eq.10}\right) $ this rewrites as%
\begin{equation*}
-B(x_{1}\otimes 1_{H};X_{2},g)-\gamma _{1}B(x_{1}\otimes
1_{H};GX_1X_{2},g)=-B(x_{1}\otimes 1_{H};X_{2},g)
\end{equation*}%
i.e.%
\begin{equation*}
\gamma _{1}B(x_{1}\otimes 1_{H};GX_1X_{2},g)=0.
\end{equation*}%
In view of the form of the element we get
\begin{equation*}
\gamma _{1}[-B(g\otimes 1_H, g,gx_2) + B(x_{1}\otimes 1_{H};G,gx_1x_2)]=0.
\end{equation*}
This is \ref{G,x1,GF2,gx2}.

\subparagraph{Case $f=x_{1}$%
\protect\begin{equation*}
-B(x_{1}\otimes 1_{H};X_{2},x_{1})-\protect\gamma _{1}B(x_{1}\otimes
1_{H};GX_1X_{2},x_{1})=-B(gx_{1}\otimes g;X_{2},x_{1}).
\protect\end{equation*}%
}

By applying $\left( \ref{eq.10}\right) $ this rewrites as%
\begin{equation*}
-B(x_{1}\otimes 1_{H};X_{2},x_{1})-\gamma _{1}B(x_{1}\otimes
1_{H};GX_{2},x_{1})=+B(x_{1}\otimes 1_{H};X_{2},x_{1})
\end{equation*}%
i.e.%
\begin{equation*}
2B(x_{1}\otimes 1_{H};X_{2},x_{1})+\gamma _{1}B(x_{1}\otimes
1_{H};GX_1X_{2},x_{1})=0
\end{equation*}

In view of the form of the element we get%
\begin{equation*}
2B(x_{1}\otimes 1_{H};1_{A},x_{1}x_{2})+\gamma _{1}B(g\otimes
1_{H};G,x_{1}x_{2})=0.
\end{equation*}%
This is $\left( \ref{G,x1, GF2,x1x2}\right) .$

\paragraph{\textbf{Equality (\protect\ref{GF8})}}

rewrites as

\begin{equation*}
B(x_{1}\otimes 1_{H};X_{1}X_{2},f)=-B(gx_{1}\otimes g;X_{1}X_{2},f).
\end{equation*}

\subparagraph{Case $f=1_{H}$}

\begin{equation*}
B(x_{1}\otimes 1_{H};X_{1}X_{2},1_{H})=-B(gx_{1}\otimes g;X_{1}X_{2},1_{H}).
\end{equation*}%
By applying $\left( \ref{eq.10}\right) $ this rewrites as%
\begin{equation*}
B(x_{1}\otimes 1_{H};X_{1}X_{2},1_{H})=-B(x_{1}\otimes
1_{H};X_{1}X_{2},1_{H}).
\end{equation*}%
In view of the form of the elements we get%
\begin{equation*}
B\left( g\otimes 1_{H};1_{A},x_{2}\right) - B(x_{1}\otimes
1_{H};1_{A},x_{1}x_{2})=0
\end{equation*}%
This is (\ref{G,x1,GF6,x2}).

\subparagraph{Case $f=gx_{1}$}

\begin{equation*}
B(x_{1}\otimes 1_{H};X_{1}X_{2},gx_{1})=-B(gx_{1}\otimes
g;X_{1}X_{2},gx_{1}).
\end{equation*}%
By applying $\left( \ref{eq.10}\right) $ this rewrites as%
\begin{equation*}
B(x_{1}\otimes 1_{H};X_{1}X_{2},gx_{1})=B(x_{1}\otimes
1_{H};X_{1}X_{2},gx_{1}).
\end{equation*}%
and we get nothing new.

\subsubsection{Case $x_{2}\otimes 1_{H}$}

\paragraph{\textbf{Equality }$\left( \protect\ref{GF1}\right) $}

rewrites as
\begin{equation*}
\begin{array}{c}
\alpha B(x_{2}\otimes 1_{H};G,f)+ \\
+\gamma _{1}B(x_{2}\otimes 1_{H};X_{1},f)+\gamma _{2}B(x_{2}\otimes
1_{H};X_{2},f)%
\end{array}%
=\alpha B(x_{2}g\otimes g;G,f)
\end{equation*}%
i.e.%
\begin{equation*}
\begin{array}{c}
\alpha B(x_{2}\otimes 1_{H};G,f)+ \\
+\gamma _{1}B(x_{2}\otimes 1_{H};X_{1},f)+\gamma _{2}B(x_{2}\otimes
1_{H};X_{2},f)%
\end{array}%
=-\alpha B(gx_{2}\otimes g;G,f).
\end{equation*}

\subparagraph{Case $f=g$}

\begin{equation*}
\begin{array}{c}
\alpha B(x_{2}\otimes 1_{H};G,g)+ \\
+\gamma _{1}B(x_{2}\otimes 1_{H};X_{1},g)+\gamma _{2}B(x_{2}\otimes
1_{H};X_{2},g)%
\end{array}%
=-\alpha B(gx_{2}\otimes g;G,g).
\end{equation*}

By applying $\left( \ref{eq.10}\right) $ this rewrites as%
\begin{equation*}
\begin{array}{c}
\alpha B(x_{2}\otimes 1_{H};G,g)+ \\
+\gamma _{1}B(x_{2}\otimes 1_{H};X_{1},g)+\gamma _{2}B(x_{2}\otimes
1_{H};X_{2},g)%
\end{array}%
=-\alpha B(x_{2}\otimes 1_{H};G,g).
\end{equation*}%
In view of the form of the element we get%
\begin{equation}
\begin{array}{c}
2\alpha B(x_{2}\otimes 1_{H};G,g)+ \\
-\gamma _{1}B(x_{2}\otimes 1_{H};1_{A},gx_{1})+\gamma _{2}\left[ -B(g\otimes
1_{H};1_{A},g)-B(x_{2}\otimes \ 1_{H};1_{A},gx_{2})\right]%
\end{array}%
=0.  \label{G,x2, GF1,g}
\end{equation}

\subparagraph{Case $f=x_{1}$}

\begin{equation*}
\begin{array}{c}
\alpha B(x_{2}\otimes 1_{H};G,x_{1})+ \\
+\gamma _{1}B(x_{2}\otimes 1_{H};X_{1},x_{1})+\gamma _{2}B(x_{2}\otimes
1_{H};X_{2},x_{1})%
\end{array}%
=-\alpha B(gx_{2}\otimes g;G,x_{1}).
\end{equation*}%
By applying $\left( \ref{eq.10}\right) $ this rewrites as%
\begin{equation*}
\begin{array}{c}
\alpha B(x_{2}\otimes 1_{H};G,x_{1})+ \\
+\gamma _{1}B(x_{2}\otimes 1_{H};X_{1},x_{1})+\gamma _{2}B(x_{2}\otimes
1_{H};X_{2},x_{1})%
\end{array}%
=\alpha B(gx_{2}\otimes g;G,x_{1}).
\end{equation*}%
and we obtain%
\begin{equation*}
\gamma _{1}B(x_{2}\otimes 1_{H};X_{1},x_{1})+\gamma _{2}B(x_{2}\otimes
1_{H};X_{2},x_{1})=0
\end{equation*}

In view of the form of the element we get%
\begin{equation}
\gamma _{2}\left[ -B(g\otimes 1_{H};1_{A},x_{1})-B(x_{2}\otimes \
1_{H};1_{A},x_{1}x_{2})\right] =0  \label{G,x2, GF1,x1}
\end{equation}

\subparagraph{Case $f=x_{2}$}

\begin{equation*}
\begin{array}{c}
\alpha B(x_{2}\otimes 1_{H};G,x_{2})+ \\
+\gamma _{1}B(x_{2}\otimes 1_{H};X_{1},x_{2})+\gamma _{2}B(x_{2}\otimes
1_{H};X_{2},x_{2})%
\end{array}%
=-\alpha B(gx_{2}\otimes g;G,x_{2}).
\end{equation*}%
By applying $\left( \ref{eq.10}\right) $ this rewrites as

\begin{equation*}
\begin{array}{c}
\alpha B(x_{2}\otimes 1_{H};G,x_{2})+ \\
+\gamma _{1}B(x_{2}\otimes 1_{H};X_{1},x_{2})+\gamma _{2}B(x_{2}\otimes
1_{H};X_{2},x_{2})%
\end{array}%
=\alpha B(gx_{2}\otimes g;G,x_{2}).
\end{equation*}%
and we obtain%
\begin{equation*}
\gamma _{1}B(x_{2}\otimes 1_{H};X_{1},x_{2})+\gamma _{2}B(x_{2}\otimes
1_{H};X_{2},x_{2})=0
\end{equation*}%
In view of the form of the element we get%
\begin{equation}
\gamma _{1}B(x_{2}\otimes 1_{H};1_{A},x_{1}x_{2})-\gamma _{2}B(g\otimes
1_{H};1_{A},x_{2})=0  \label{G,x2, GF1,x2}
\end{equation}

\subparagraph{Case $f=gx_{1}x_{2}$}

\begin{equation*}
\begin{array}{c}
\alpha B(x_{2}\otimes 1_{H};G,gx_{1}x_{2})+ \\
+\gamma _{1}B(x_{2}\otimes 1_{H};X_{1},gx_{1}x_{2})+\gamma
_{2}B(x_{2}\otimes 1_{H};X_{2},gx_{1}x_{2})%
\end{array}%
=-\alpha B(gx_{2}\otimes g;G,gx_{1}x_{2}).
\end{equation*}%
By applying $\left( \ref{eq.10}\right) $ this rewrites as%
\begin{equation*}
\begin{array}{c}
\alpha B(x_{2}\otimes 1_{H};G,gx_{1}x_{2})+ \\
+\gamma _{1}B(x_{2}\otimes 1_{H};X_{1},gx_{1}x_{2})+\gamma
_{2}B(x_{2}\otimes 1_{H};X_{2},gx_{1}x_{2})%
\end{array}%
=-\alpha B(x_{2}\otimes 1_{H};G,gx_{1}x_{2})
\end{equation*}%
i.e.%
\begin{equation*}
\begin{array}{c}
2\alpha B(x_{2}\otimes 1_{H};G,gx_{1}x_{2})+ \\
+\gamma _{1}B(x_{2}\otimes 1_{H};X_{1},gx_{1}x_{2})+\gamma
_{2}B(x_{2}\otimes 1_{H};X_{2},gx_{1}x_{2})%
\end{array}%
=0.
\end{equation*}%
In view of the form of the elements we get%
\begin{equation}
2\alpha B(x_{2}\otimes 1_{H};G,gx_{1}x_{2})-\gamma _{2}B\left( g\otimes
1_{H};1_{A},gx_{1}x_{2}\right) =0.  \label{G,x2, GF1,gx1x2}
\end{equation}

\paragraph{\textbf{Equality }$\left( \protect\ref{GF2}\right) $}

rewrites as
\begin{equation*}
\begin{array}{c}
B(x_{2}\otimes 1_{H};1_{A},f) \\
+\gamma _{1}B(x_{2}\otimes 1_{H};GX_{1},f)+\gamma _{2}B(x_{2}\otimes
1_{H};GX_{2},f)%
\end{array}%
=B(x_{2}g\otimes g;1_{A},f)
\end{equation*}%
i.e.%
\begin{equation*}
\begin{array}{c}
B(x_{2}\otimes 1_{H};1_{A},f) \\
+\gamma _{1}B(x_{2}\otimes 1_{H};GX_{1},f)+\gamma _{2}B(x_{2}\otimes
1_{H};GX_{2},f)%
\end{array}%
=-B(gx_{2}\otimes g;1_{A},f).
\end{equation*}

\subparagraph{Case $f=1_{H}$}

\begin{equation*}
\begin{array}{c}
B(x_{2}\otimes 1_{H};1_{A},1_{H}) \\
+\gamma _{1}B(x_{2}\otimes 1_{H};GX_{1},1_{H})+\gamma _{2}B(x_{2}\otimes
1_{H};GX_{2},1_{H})%
\end{array}%
=-B(gx_{2}\otimes g;1_{A},1_{H}).
\end{equation*}%
By applying $\left( \ref{eq.10}\right) $ this rewrites as%
\begin{equation*}
\begin{array}{c}
B(x_{2}\otimes 1_{H};1_{A},1_{H}) \\
+\gamma _{1}B(x_{2}\otimes 1_{H};GX_{1},1_{H})+\gamma _{2}B(x_{2}\otimes
1_{H};GX_{2},1_{H})%
\end{array}%
=-B(x_{2}\otimes 1_{H};1_{A},1_{H})
\end{equation*}%
i.e.%
\begin{equation*}
\begin{array}{c}
2B(x_{2}\otimes 1_{H};1_{A},1_{H}) \\
+\gamma _{1}B(x_{2}\otimes 1_{H};GX_{1},1_{H})+\gamma _{2}B(x_{2}\otimes
1_{H};GX_{2},1_{H})%
\end{array}%
=0.
\end{equation*}%
In view of the form of the elements we get%
\begin{equation}
\begin{array}{c}
2B(x_{2}\otimes 1_{H};1_{A},1_{H}) \\
+\gamma _{1}B(x_{2}\otimes 1_{H};G,x_{1})+\gamma _{2}\left[ B(g\otimes
1_{H};G,1_{H})+B(x_{2}\otimes \ 1_{H};G,x_{2}\right]%
\end{array}%
=0.  \label{G,x2, GF2,1H}
\end{equation}

\subparagraph{Case $f=x_{1}x_{2}$}

\begin{equation*}
\begin{array}{c}
B(x_{2}\otimes 1_{H};1_{A},x_{1}x_{2}) \\
+\gamma _{1}B(x_{2}\otimes 1_{H};GX_{1},x_{1}x_{2})+\gamma
_{2}B(x_{2}\otimes 1_{H};GX_{2},x_{1}x_{2})%
\end{array}%
=-B(gx_{2}\otimes g;1_{A},x_{1}x_{2}).
\end{equation*}%
By applying $\left( \ref{eq.10}\right) $ this rewrites as%
\begin{equation*}
\begin{array}{c}
B(x_{2}\otimes 1_{H};1_{A},x_{1}x_{2}) \\
+\gamma _{1}B(x_{2}\otimes 1_{H};GX_{1},x_{1}x_{2})+\gamma
_{2}B(x_{2}\otimes 1_{H};GX_{2},x_{1}x_{2})%
\end{array}%
=-B(x_{2}\otimes 1_{H};1_{A},x_{1}x_{2})
\end{equation*}%
i.e.%
\begin{equation*}
\begin{array}{c}
2B(x_{2}\otimes 1_{H};1_{A},x_{1}x_{2}) \\
+\gamma _{1}B(x_{2}\otimes 1_{H};GX_{1},x_{1}x_{2})+\gamma
_{2}B(x_{2}\otimes 1_{H};GX_{2},x_{1}x_{2})%
\end{array}%
=0.
\end{equation*}%
In view of the form of the element we get%
\begin{equation}
2B(x_{2}\otimes 1_{H};1_{A},x_{1}x_{2})+\gamma _{2}B\left( g\otimes
1_{H};G,x_{1}x_{2}\right) =0  \label{G,x2, GF2,x1x2}
\end{equation}

\subparagraph{Case $f=gx_{1}$}

\begin{equation*}
\begin{array}{c}
B(x_{2}\otimes 1_{H};1_{A},gx_{1}) \\
+\gamma _{1}B(x_{2}\otimes 1_{H};GX_{1},gx_{1})+\gamma _{2}B(x_{2}\otimes
1_{H};GX_{2},gx_{1})%
\end{array}%
=-B(gx_{2}\otimes g;1_{A},gx_{1}).
\end{equation*}%
By applying $\left( \ref{eq.10}\right) $ this rewrites as
\begin{equation*}
\begin{array}{c}
B(x_{2}\otimes 1_{H};1_{A},gx_{1}) \\
+\gamma _{1}B(x_{2}\otimes 1_{H};GX_{1},gx_{1})+\gamma _{2}B(x_{2}\otimes
1_{H};GX_{2},gx_{1})%
\end{array}%
=B(x_{2}\otimes 1_{H};1_{A},gx_{1}).
\end{equation*}%
i.e.%
\begin{equation*}
\gamma _{1}B(x_{2}\otimes 1_{H};GX_{1},gx_{1})+\gamma _{2}B(x_{2}\otimes
1_{H};GX_{2},gx_{1})=0.
\end{equation*}%
In view of the form of the element we get%
\begin{equation}
\gamma _{2}\left[ B(g\otimes 1_{H};G,gx_{1})+B(x_{2}\otimes \
1_{H};G,gx_{1}x_{2})\right] =0.  \label{G,x2, GF2,gx1}
\end{equation}

\subparagraph{Case $f=gx_{2}$}

\begin{equation*}
\begin{array}{c}
B(x_{2}\otimes 1_{H};1_{A},gx_{2}) \\
+\gamma _{1}B(x_{2}\otimes 1_{H};GX_{1},gx_{2})+\gamma _{2}B(x_{2}\otimes
1_{H};GX_{2},gx_{2})%
\end{array}%
=-B(gx_{2}\otimes g;1_{A},gx_{2}).
\end{equation*}%
By applying $\left( \ref{eq.10}\right) $ this rewrites as
\begin{equation*}
\begin{array}{c}
B(x_{2}\otimes 1_{H};1_{A},gx_{2}) \\
+\gamma _{1}B(x_{2}\otimes 1_{H};GX_{1},gx_{2})+\gamma _{2}B(x_{2}\otimes
1_{H};GX_{2},gx_{2})%
\end{array}%
=B(x_{2}\otimes 1_{H};1_{A},gx_{2}).
\end{equation*}%
i.e.%
\begin{equation*}
\gamma _{1}B(x_{2}\otimes 1_{H};GX_{1},gx_{2})+\gamma _{2}B(x_{2}\otimes
1_{H};GX_{2},gx_{2})=0.
\end{equation*}%
In view of the form of the element we get%
\begin{equation}
-\gamma _{1}B(x_{2}\otimes 1_{H};G,gx_{1}x_{2})+\gamma _{2}B(g\otimes
1_{H};G,gx_{2})=0.  \label{G,x2, GF2,gx2}
\end{equation}

\paragraph{\textbf{Equality }$\left( \protect\ref{GF3}\right) $}

rewrites as
\begin{equation*}
\gamma _{2}B(x_{2}\otimes 1_{H};X_{1}X_{2},f)-\alpha B(x_{2}\otimes
1_{H};GX_{1},f)=\alpha B(x_{2}g\otimes g;GX_{1},f)
\end{equation*}%
i.e.%
\begin{equation*}
\gamma _{2}B(x_{2}\otimes 1_{H};X_{1}X_{2},f)-\alpha B(x_{2}\otimes
1_{H};GX_{1},f)=-\alpha B(gx_{2}\otimes g;GX_{1},f).
\end{equation*}

\subparagraph{Case $f=1_{H}$}

\begin{equation*}
\gamma _{2}B(x_{2}\otimes 1_{H};X_{1}X_{2},1_{H})-\alpha B(x_{2}\otimes
1_{H};GX_{1},1_{H})=-\alpha B(gx_{2}\otimes g;GX_{1},1_{H}).
\end{equation*}%
By applying $\left( \ref{eq.10}\right) $ this rewrites as%
\begin{equation*}
\gamma _{2}B(x_{2}\otimes 1_{H};X_{1}X_{2},1_{H})-\alpha B(x_{2}\otimes
1_{H};GX_{1},1_{H})=-\alpha B(x_{2}\otimes 1_{H};GX_{1},1_{H}).
\end{equation*}%
In view of the form of the element we get
\begin{equation}
\gamma _{2}[B(x_{2}\otimes 1_{H};1_{A},x_{1}x_{2})+B(g\otimes
1_{H};1_{A},x_{1})]=0  \label{G,x2, FG3,1H}
\end{equation}

This is \ref{G,x2, GF1,x1}

\subparagraph{Case $f=gx_{2}$}

\begin{equation*}
\gamma _{2}B(x_{2}\otimes 1_{H};X_{1}X_2,gx_{2})-\alpha B(x_{2}\otimes
1_{H};GX_{1},gx_{2})=-\alpha B(gx_{2}\otimes g;GX_{1},gx_{2}).
\end{equation*}%
By applying $\left( \ref{eq.10}\right) $ this rewrites as%
\begin{equation*}
\gamma _{2}B(x_{2}\otimes 1_{H};X_{1}X_2,gx_{2})-\alpha B(x_{2}\otimes
1_{H};GX_{1},gx_{2})=\alpha B(x_{2}\otimes 1_{H};GX_{1},gx_{2}).
\end{equation*}%
i.e.%
\begin{equation*}
\gamma _{2}B(x_{2}\otimes 1_{H};X_{1}X_2,gx_{2})+2\alpha B(x_{2}\otimes
1_{H};GX_{1},gx_{2})=0.
\end{equation*}%
In view of the form of the element we get

\begin{equation}
-\gamma _{2}B(g\otimes 1_{H};1_{A},gx_{1}x_{2})-2\alpha B(x_{2}\otimes
1_{H};G,gx_{1}x_{2})=0.  \label{G,x2, GF3,gx2}
\end{equation}

\paragraph{\textbf{Equality }$\left( \protect\ref{GF4}\right) $}

rewrites as
\begin{equation*}
-\gamma _{1}B(x_{2}\otimes 1_{H};X_{1}X_{2},f)-\alpha B(x_{2}\otimes
1_{H};GX_{2},f)=\alpha B(x_{2}g\otimes g;GX_{2},f)
\end{equation*}%
i.e.%
\begin{equation*}
-\gamma _{1}B(x_{2}\otimes 1_{H};X_{1}X_{2},f)-\alpha B(x_{2}\otimes
1_{H};GX_{2},f)=-\alpha B(gx_{2}\otimes g;GX_{2},f).
\end{equation*}

\subparagraph{Case $f=1_{H}$}

\begin{equation*}
-\gamma _{1}B(x_{2}\otimes 1_{H};X_{1}X_{2},1_{H})-\alpha B(x_{2}\otimes
1_{H};GX_{2},1_{H})=-\alpha B(gx_{2}\otimes g;GX_{2},1_{H}).
\end{equation*}%
By applying $\left( \ref{eq.10}\right) $ this rewrites as%
\begin{equation*}
-\gamma _{1}B(x_{2}\otimes 1_{H};X_{1}X_{2},1_{H})-\alpha B(x_{2}\otimes
1_{H};GX_{2},1_{H})=-\alpha B(x_{2}\otimes 1_{H};GX_{2},1_{H})
\end{equation*}%
i.e.%
\begin{equation}
\gamma _{1}\left[ B(g\otimes 1_{H};1_{A},x_{1})+B(x_{2}\otimes \
1_{H};1_{A},x_{1}x_{2})\right] =0.  \label{G,x2, GF4,1H}
\end{equation}

\subparagraph{Case $f=x_{1}x_{2}$}

\begin{equation*}
-\gamma _{1}B(x_{2}\otimes 1_{H};X_1X_{2},x_{1}x_{2})-\alpha B(x_{2}\otimes
1_{H};GX_{2},x_{1}x_{2})=-\alpha B(gx_{2}\otimes g;GX_{2},x_{1}x_{2}).
\end{equation*}%
By applying $\left( \ref{eq.10}\right) $ this rewrites as

\begin{equation*}
-\gamma _{1}B(x_{2}\otimes 1_{H};X_1X_{2},x_{1}x_{2})-\alpha B(x_{2}\otimes
1_{H};GX_{2},x_{1}x_{2})=-\alpha B(x_{2}\otimes 1_{H};GX_{2},x_{1}x_{2})
\end{equation*}%
i.e.%
\begin{equation*}
\gamma _{1}B(x_{2}\otimes 1_{H};X_1X_{2},x_{1}x_{2})=0.
\end{equation*}%
In view of the form of the element we get nothing new.

\subparagraph{Case $f=gx_{1}$}

\begin{equation*}
-\gamma _{1}B(x_{2}\otimes 1_{H};X_1X_{2},gx_{1})-\alpha B(x_{2}\otimes
1_{H};GX_{2},gx_{1})=-\alpha B(gx_{2}\otimes g;GX_{2},gx_{1}).
\end{equation*}%
By applying $\left( \ref{eq.10}\right) $ this rewrites as

\begin{equation*}
-\gamma _{1}B(x_{2}\otimes 1_{H};X_{1}X_{2},gx_{1})-\alpha B(x_{2}\otimes
1_{H};GX_{2},x_{1})=\alpha B(x_{2}\otimes 1_{H};GX_{2},gx_{1}).
\end{equation*}%
i.e.%
\begin{equation*}
-\gamma _{1}B(x_{2}\otimes 1_{H};X_{1}X_{2},gx_{1})-2\alpha B(x_{2}\otimes
1_{H};GX_{2},gx_{1})=0.
\end{equation*}%
In view of the form of the element we get%
\begin{equation}
\alpha \left[ B(g\otimes 1_{H};G,gx_{1})+B(x_{2}\otimes \
1_{H};G,gx_{1}x_{2})\right] =0  \label{G,x2, GF4,gx1}
\end{equation}

\subparagraph{Case $f=gx_{2}$}

\begin{equation*}
-\gamma _{1}B(x_{2}\otimes 1_{H};X_1X_{2},gx_{2})-\alpha B(x_{2}\otimes
1_{H};GX_{2},gx_{2})=-\alpha B(gx_{2}\otimes g;GX_{2},gx_{2}).
\end{equation*}%
By applying $\left( \ref{eq.10}\right) $ this rewrites as

\begin{equation*}
-\gamma _{1}B(x_{2}\otimes 1_{H};X_{1}X_{2},gx_{2})-\alpha B(x_{2}\otimes
1_{H};GX_{2},x_{2})=\alpha B(x_{2}\otimes 1_{H};GX_{2},gx_{2}).
\end{equation*}%
i.e.%
\begin{equation*}
-\gamma _{1}B(x_{2}\otimes 1_{H};X_{1}X_{2},gx_{2})-2\alpha B(x_{2}\otimes
1_{H};GX_{2},gx_{2})=0.
\end{equation*}%
In view of the form of the element we get%
\begin{equation}
\gamma _{1}B(g\otimes 1_{H};1_{A},gx_{1}x_{2})-2\alpha B(g\otimes
1_{H};G,gx_{2})=0.  \label{G,x2, GF4,gx2}
\end{equation}

\paragraph{\textbf{Equality }$\left( \protect\ref{GF5}\right) $}

rewrites as
\begin{equation*}
\alpha B(x_{2}\otimes 1_{H};GX_{1}X_{2},f)=\alpha B(x_{2}g\otimes
g;GX_{1}X_{2},f)
\end{equation*}%
i.e.%
\begin{equation*}
\alpha B(x_{2}\otimes 1_{H};GX_{1}X_{2},f)=-\alpha B(gx_{2}\otimes
g;GX_{1}X_{2},f).
\end{equation*}

\subparagraph{Case $f=g$}

\begin{equation*}
\alpha B(x_{2}\otimes 1_{H};GX_{1}X_{2},g)=-\alpha B(gx_{2}\otimes
g;GX_{1}X_{2},g).
\end{equation*}%
By applying $\left( \ref{eq.10}\right) $ this rewrites as%
\begin{equation*}
\alpha B(x_{2}\otimes 1_{H};GX_{1}X_{2},g)=-\alpha B(x_{2}\otimes
1_{H};GX_{1}X_{2},g)
\end{equation*}%
i.e.%
\begin{equation*}
\alpha B(x_{2}\otimes 1_{H};GX_{1}X_{2},g)=0.
\end{equation*}%
In view of the form of the element we get%
\begin{equation}
\alpha \left[ B(g\otimes 1_{H};G,gx_{1})+B(x_{2}\otimes \
1_{H};G,gx_{1}x_{2})\right] =0  \label{G,x2, GF5,g}
\end{equation}

\subparagraph{Case $f=x_{2}$}

\begin{equation*}
\alpha B(x_{2}\otimes 1_{H};GX_{1}X_{2},x_{2})=-\alpha B(gx_{2}\otimes
g;GX_{1}X_{2},x_{2}).
\end{equation*}%
By applying $\left( \ref{eq.10}\right) $ this rewrites as%
\begin{equation*}
\alpha B(x_{2}\otimes 1_{H};GX_{1}X_{2},x_{2})=\alpha B(x_{2}\otimes
1_{H};GX_{1}X_{2},x_{2})
\end{equation*}%
which is obvious.

\paragraph{\textbf{Equality }$\left( \protect\ref{GF6}\right) $}

rewrites as
\begin{equation*}
-B(x_{2}\otimes 1_{H};X_{1},f)+\gamma _{2}B(x_{2}\otimes
1_{H};GX_{1}X_{2},f)=B(x_{2}g\otimes g;X_{1},f)
\end{equation*}%
i.e.%
\begin{equation*}
-B(x_{2}\otimes 1_{H};X_{1},f)+\gamma _{2}B(x_{2}\otimes
1_{H};GX_{1},f)=-B(gx_{2}\otimes g;X_{1},f).
\end{equation*}

\subparagraph{Case $f=g$}

\begin{equation*}
-B(x_{2}\otimes 1_{H};X_{1},g)+\gamma _{2}B(x_{2}\otimes
1_{H};GX_{1}X_{2},g)=-B(gx_{2}\otimes g;X_{1},g).
\end{equation*}%
By applying $\left( \ref{eq.10}\right) $ this rewrites as%
\begin{equation*}
-B(x_{2}\otimes 1_{H};X_{1},g)+\gamma _{2}B(x_{2}\otimes
1_{H};GX_{1}X_{2},g)=-B(x_{2}\otimes 1_{H};X_{1},g)
\end{equation*}%
i.e.%
\begin{equation*}
\gamma _{2}B(x_{2}\otimes 1_{H};GX_{1},g)=0.
\end{equation*}%
In view of the form of the element we get
\begin{equation}
\gamma _{2}[B(g\otimes 1_{H},G,gx_{1})+B(x_{2}\otimes
1_{H};G,gx_{1}x_{2})]=0.  \label{G,x2, GF6,g}
\end{equation}

\subparagraph{Case $f=x_{2}$}

\begin{equation*}
-B(x_{2}\otimes 1_{H};X_{1},x_{2})+\gamma _{2}B(x_{2}\otimes
1_{H};GX_{1}X_{2},x_{2})=-B(gx_{2}\otimes g;X_{1},x_{2}).
\end{equation*}%
By applying $\left( \ref{eq.10}\right) $ this rewrites as%
\begin{equation*}
-B(x_{2}\otimes 1_{H};X_{1},x_{2})+\gamma _{2}B(x_{2}\otimes
1_{H};GX_{1}X_{2},x_{2})=B(x_{2}\otimes 1_{H};X_{1},x_{2})
\end{equation*}%
i.e.%
\begin{equation*}
-2B(x_{2}\otimes 1_{H};X_{1},x_{2})+\gamma _{2}B(x_{2}\otimes
1_{H};GX_{1}X_{2},x_{2})=0.
\end{equation*}%
In view of the form of the element we get%
\begin{equation}
2B(x_{2}\otimes 1_{H};1_{A},x_{1}x_{2})+\gamma _{2}B(g\otimes
1_{H};G,x_{1}x_{2})=0.  \label{G,x2; GF6,x2}
\end{equation}

\paragraph{\textbf{Equality} $\left( \protect\ref{GF7}\right) $}

rewrites as%
\begin{equation*}
-B(x_{2}\otimes 1_{H};X_{2},f)-\gamma _{1}B(x_{2}\otimes
1_{H};GX_1X_{2},f)=B(x_{2}g\otimes g;X_{2},f)
\end{equation*}%
i.e.%
\begin{equation*}
-B(x_{2}\otimes 1_{H};X_{2},f)-\gamma _{1}B(x_{2}\otimes
1_{H};GX_1X_{2},f)=-B(gx_{2}\otimes g;X_{2},f).
\end{equation*}

\subparagraph{Case $f=g$}

\begin{equation*}
-B(x_{2}\otimes 1_{H};X_{2},g)-\gamma _{1}B(x_{2}\otimes
1_{H};GX_1X_{2},g)=-B(gx_{2}\otimes g;X_{2},g).
\end{equation*}%
By applying $\left( \ref{eq.10}\right) $ this rewrites as%
\begin{equation*}
-B(x_{2}\otimes 1_{H};X_{2},g)-\gamma _{1}B(x_{2}\otimes
1_{H};GX_1X_{2},g)=-B(x_{2}\otimes 1_{H};X_{2},g)
\end{equation*}%
i.e.%
\begin{equation*}
\gamma _{1}B(x_{2}\otimes 1_{H};GX_1X_{2},g)=0.
\end{equation*}%
In view of the form of the element we get
\begin{equation}
\gamma _{1}[B(g\otimes 1_H;G,gx_1)+B(x_{2}\otimes 1_{H};G,gx_1x_2)]=0.
\label{G,x2; GF7,g}
\end{equation}

\subparagraph{Case $f=x_{1}$}

\begin{equation*}
-B(x_{2}\otimes 1_{H};X_{2},x_{1})-\gamma _{1}B(x_{2}\otimes
1_{H};GX_1X_{2},x_{1})=-B(gx_{2}\otimes g;X_{2},x_{1}).
\end{equation*}

By applying $\left( \ref{eq.10}\right) $ this rewrites as%
\begin{equation*}
-B(x_{2}\otimes 1_{H};X_{2},x_{1})-\gamma _{1}B(x_{2}\otimes
1_{H};GX_1X_{2},x_{1})=B(x_{2}\otimes 1_{H};X_{2},x_{1})
\end{equation*}%
i.e.%
\begin{equation*}
-2B(x_{2}\otimes 1_{H};X_{2},x_{1})-\gamma _{1}B(x_{2}\otimes
1_{H};GX_1X_{2},x_{1})=0.
\end{equation*}%
In view of the form of the element we get

\begin{equation}
-B(g\otimes 1_{H};1_{A},x_{1})-B(x_{2}\otimes \ 1_{H};1_{A},x_{1}x_{2})=0.
\label{G,x2; GF7,x1}
\end{equation}

\subparagraph{Case $f=x_{2}$}

\begin{equation*}
-B(x_{2}\otimes 1_{H};X_{2},x_{2})-\gamma _{1}B(x_{2}\otimes
1_{H};GX_1X_{2},x_{2})=-B(gx_{2}\otimes g;X_{2},x_{2}).
\end{equation*}

By applying $\left( \ref{eq.10}\right) $ this rewrites as%
\begin{equation*}
-B(x_{2}\otimes 1_{H};X_{2},x_{2})-\gamma _{1}B(x_{2}\otimes
1_{H};GX_{1}X_{2},x_{2})=B(x_{2}\otimes 1_{H};X_{2},x_{2})
\end{equation*}%
i.e.%
\begin{equation*}
2B(x_{2}\otimes 1_{H};X_{2},x_{2})+\gamma _{1}B(x_{2}\otimes
1_{H};GX_{1}X_{2},x_{2})=0.
\end{equation*}%
In view of the form of the element we get%
\begin{equation*}
2B(g\otimes 1_{H};1_{A},x_{2})+\gamma _{1}B(g\otimes 1_{H};G,x_{1}x_{2})=0
\end{equation*}%
which is $\left( \ref{G,g, GF2,x2}\right) .$

\subparagraph{Case $f=gx_{1}x_{2}$}

\begin{equation*}
-B(x_{2}\otimes 1_{H};X_{2},gx_{1}x_{2})-\gamma _{1}B(x_{2}\otimes
1_{H};GX_1X_{2},gx_{1}x_{2})=-B(gx_{2}\otimes g;X_{2},gx_{1}x_{2}).
\end{equation*}

By applying $\left( \ref{eq.10}\right) $ this rewrites as%
\begin{equation*}
-B(x_{2}\otimes 1_{H};X_{2},gx_{1}x_{2})-\gamma _{1}B(x_{2}\otimes
1_{H};GX_1X_{2},gx_{1}x_{2})=-B(x_{2}\otimes 1_{H};X_{2},gx_{1}x_{2})
\end{equation*}%
i.e.%
\begin{equation*}
\gamma _{1}B(x_{2}\otimes 1_{H};GX_1X_{2},gx_{1}x_{2})=0
\end{equation*}%
In view of the form of the elements we get nothing new.

\paragraph{\textbf{Equality }$\left( \protect\ref{GF8}\right) $}

rewrites as%
\begin{equation*}
B(x_{2}\otimes 1_{H};X_{1}X_{2},f)=B(x_{2}g\otimes g;X_{1}X_{2},f)
\end{equation*}%
i.e.%
\begin{equation*}
B(x_{2}\otimes 1_{H};X_{1}X_{2},f)=-B(gx_{2}\otimes g;X_{1}X_{2},f).
\end{equation*}

\subparagraph{Case $f=1_{H}$}

\begin{equation*}
B(x_{2}\otimes 1_{H};X_{1}X_{2},1_{H})=-B(gx_{2}\otimes g;X_{1}X_{2},1_{H}).
\end{equation*}%
By applying $\left( \ref{eq.10}\right) $ this rewrites as%
\begin{equation*}
B(x_{2}\otimes 1_{H};X_{1}X_{2},1_{H})=-B(x_{2}\otimes
1_{H};X_{1}X_{2},1_{H})
\end{equation*}%
i.e.%
\begin{equation*}
B(x_{2}\otimes 1_{H};X_{1}X_{2},1_{H})=0.
\end{equation*}%
In view of the form of the element we get%
\begin{equation*}
B(g\otimes 1_{H};1_{A},x_{1})+B(x_{2}\otimes \ 1_{H};1_{A},x_{1}x_{2})=0.
\end{equation*}

This is (\ref{G,x2; GF7,x1}).

\subparagraph{Case $f=gx_2$}

\begin{equation*}
B(x_{2}\otimes 1_{H};X_{1}X_{2},gx_2)=-B(gx_{2}\otimes g;X_{1}X_{2},gx_2).
\end{equation*}

By applying $\left( \ref{eq.10}\right) $ this rewrites as

\begin{equation*}
B(x_{2}\otimes 1_{H};X_{1}X_{2},gx_2)=B(gx_{2}\otimes g;X_{1}X_{2},gx_2).
\end{equation*}%
which is trivial.

\subsubsection{Case $x_{1}x_{2}\otimes 1_{H}$}

\paragraph{\textbf{Equality (\protect\ref{GF1})}}

rewrites as%
\begin{equation*}
\begin{array}{c}
\alpha B(x_{1}x_{2}\otimes 1_{H};G,f)+ \\
+\gamma _{1}B(x_{1}x_{2}\otimes 1_{H};X_{1},f)+\gamma
_{2}B(x_{1}x_{2}\otimes 1_{H};X_{2},f)%
\end{array}%
=\alpha B(x_{1}x_{2}g\otimes g;G,f)
\end{equation*}%
i.e.%
\begin{equation*}
\begin{array}{c}
\alpha B(x_{1}x_{2}\otimes 1_{H};G,f)+ \\
+\gamma _{1}B(x_{1}x_{2}\otimes 1_{H};X_{1},f)+\gamma
_{2}B(x_{1}x_{2}\otimes 1_{H};X_{2},f)%
\end{array}%
=\alpha B(gx_{1}x_{2}\otimes g;G,f).
\end{equation*}

\subparagraph{Case $f=g$}

\begin{equation*}
\begin{array}{c}
\alpha B(x_{1}x_{2}\otimes 1_{H};G,g)+ \\
+\gamma _{1}B(x_{1}x_{2}\otimes 1_{H};X_{1},g)+\gamma
_{2}B(x_{1}x_{2}\otimes 1_{H};X_{2},g)%
\end{array}%
=\alpha B(gx_{1}x_{2}\otimes g;G,g).
\end{equation*}%
By applying $\left( \ref{eq.10}\right) $ this rewrites as%
\begin{equation*}
\begin{array}{c}
\alpha B(x_{1}x_{2}\otimes 1_{H};G,g)+ \\
+\gamma _{1}B(x_{1}x_{2}\otimes 1_{H};X_{1},g)+\gamma
_{2}B(x_{1}x_{2}\otimes 1_{H};X_{2},g)%
\end{array}%
=\alpha B(x_{1}x_{2}\otimes 1_{H};G,g).
\end{equation*}%
i.e.%
\begin{equation*}
\gamma _{1}B(x_{1}x_{2}\otimes 1_{H};X_{1},g)+\gamma _{2}B(x_{1}x_{2}\otimes
1_{H};X_{2},g)=0.
\end{equation*}%
In view of the form of the element we get%
\begin{equation}
\gamma _{1}B(x_{1}x_{2}\otimes 1_{H};X_{1},g)+\gamma _{2}B(x_{1}x_{2}\otimes
1_{H};X_{2},g)=0.  \label{G,x1x2, GF1,g}
\end{equation}

\subparagraph{Case $f=x_{1}$}

\begin{equation*}
\begin{array}{c}
\alpha B(x_{1}x_{2}\otimes 1_{H};G,x_{1})+ \\
+\gamma _{1}B(x_{1}x_{2}\otimes 1_{H};X_{1},x_{1})+\gamma
_{2}B(x_{1}x_{2}\otimes 1_{H};X_{2},x_{1})%
\end{array}%
=\alpha B(gx_{1}x_{2}\otimes g;G,x_{1}).
\end{equation*}%
By applying $\left( \ref{eq.10}\right) $ this rewrites as%
\begin{equation*}
\begin{array}{c}
\alpha B(x_{1}x_{2}\otimes 1_{H};G,x_{1})+ \\
+\gamma _{1}B(x_{1}x_{2}\otimes 1_{H};X_{1},x_{1})+\gamma
_{2}B(x_{1}x_{2}\otimes 1_{H};X_{2},x_{1})%
\end{array}%
=-\alpha B(x_{1}x_{2}\otimes 1_{H};G,x_{1})
\end{equation*}%
i.e.%
\begin{equation*}
\begin{array}{c}
2\alpha B(x_{1}x_{2}\otimes 1_{H};G,x_{1})+ \\
+\gamma _{1}B(x_{1}x_{2}\otimes 1_{H};X_{1},x_{1})+\gamma
_{2}B(x_{1}x_{2}\otimes 1_{H};X_{2},x_{1})%
\end{array}%
=0.
\end{equation*}%
In view of the form of the element we get%
\begin{equation}
\begin{array}{c}
2\alpha B(x_{1}x_{2}\otimes 1_{H};G,x_{1})+ \\
+\gamma _{1}B(x_{1}x_{2}\otimes 1_{H};X_{1},x_{1})+\gamma
_{2}B(x_{1}x_{2}\otimes 1_{H};X_{2},x_{1})%
\end{array}%
=0.  \label{G,x1x2, GF1,x1}
\end{equation}

\subparagraph{Case $f=x_{2}$}

\begin{equation*}
\begin{array}{c}
\alpha B(x_{1}x_{2}\otimes 1_{H};G,x_{2})+ \\
+\gamma _{1}B(x_{1}x_{2}\otimes 1_{H};X_{1},x_{2})+\gamma
_{2}B(x_{1}x_{2}\otimes 1_{H};X_{2},x_{2})%
\end{array}%
=\alpha B(gx_{1}x_{2}\otimes g;G,x_{2}).
\end{equation*}%
By applying $\left( \ref{eq.10}\right) $ this rewrites as%
\begin{equation*}
\begin{array}{c}
\alpha B(x_{1}x_{2}\otimes 1_{H};G,x_{2})+ \\
+\gamma _{1}B(x_{1}x_{2}\otimes 1_{H};X_{1},x_{2})+\gamma
_{2}B(x_{1}x_{2}\otimes 1_{H};X_{2},x_{2})%
\end{array}%
=-\alpha B(x_{1}x_{2}\otimes 1_{H};G,x_{2})
\end{equation*}%
i.e.%
\begin{equation*}
\begin{array}{c}
2\alpha B(x_{1}x_{2}\otimes 1_{H};G,x_{2})+ \\
+\gamma _{1}B(x_{1}x_{2}\otimes 1_{H};X_{1},x_{2})+\gamma
_{2}B(x_{1}x_{2}\otimes 1_{H};X_{2},x_{2})%
\end{array}%
=0.
\end{equation*}%
In view of the form of the element we get%
\begin{equation}
\begin{array}{c}
2\alpha B(x_{1}x_{2}\otimes 1_{H};G,x_{2})+ \\
+\gamma _{1}B(x_{1}x_{2}\otimes 1_{H};X_{1},x_{2})+\gamma
_{2}B(x_{1}x_{2}\otimes 1_{H};X_{2},x_{2})%
\end{array}%
=0.  \label{G,x1x2, GF1,x2}
\end{equation}

\subparagraph{Case $f=gx_{1}x_{2}$}

\begin{equation*}
\begin{array}{c}
\alpha B(x_{1}x_{2}\otimes 1_{H};G,gx_{1}x_{2})+ \\
+\gamma _{1}B(x_{1}x_{2}\otimes 1_{H};X_{1},gx_{1}x_{2})+\gamma
_{2}B(x_{1}x_{2}\otimes 1_{H};X_{2},gx_{1}x_{2})%
\end{array}%
=\alpha B(gx_{1}x_{2}\otimes g;G,gx_{1}x_{2}).
\end{equation*}%
By applying $\left( \ref{eq.10}\right) $ this rewrites as%
\begin{equation*}
\begin{array}{c}
\alpha B(x_{1}x_{2}\otimes 1_{H};G,gx_{1}x_{2})+ \\
+\gamma _{1}B(x_{1}x_{2}\otimes 1_{H};X_{1},gx_{1}x_{2})+\gamma
_{2}B(x_{1}x_{2}\otimes 1_{H};X_{2},gx_{1}x_{2})%
\end{array}%
=\alpha B(x_{1}x_{2}\otimes 1_{H};G,gx_{1}x_{2})
\end{equation*}%
i.e.%
\begin{equation}
\gamma _{1}B(x_{1}x_{2}\otimes 1_{H};X_{1},gx_{1}x_{2})+\gamma
_{2}B(x_{1}x_{2}\otimes 1_{H};X_{2},gx_{1}x_{2})=0  \label{G,x1x2, GF1,gx1x2}
\end{equation}

\paragraph{\textbf{Equality (\protect\ref{GF2})}}

rewrites as%
\begin{equation*}
\begin{array}{c}
B(x_{1}x_{2}\otimes 1_{H};1_{A},f) \\
+\gamma _{1}B(x_{1}x_{2}\otimes 1_{H};GX_{1},f)+\gamma
_{2}B(x_{1}x_{2}\otimes 1_{H};GX_{2},f)%
\end{array}%
=B(x_{1}x_{2}g\otimes g;1_{A},f)
\end{equation*}%
i.e.%
\begin{equation*}
\begin{array}{c}
B(x_{1}x_{2}\otimes 1_{H};1_{A},f) \\
+\gamma _{1}B(x_{1}x_{2}\otimes 1_{H};GX_{1},f)+\gamma
_{2}B(x_{1}x_{2}\otimes 1_{H};GX_{2},f)%
\end{array}%
=B(gx_{1}x_{2}\otimes g;1_{A},f)
\end{equation*}

\subparagraph{Case $f=1_{H}$}

\begin{equation*}
\begin{array}{c}
B(x_{1}x_{2}\otimes 1_{H};1_{A},1_{H}) \\
+\gamma _{1}B(x_{1}x_{2}\otimes 1_{H};GX_{1},1_{H})+\gamma
_{2}B(x_{1}x_{2}\otimes 1_{H};GX_{2},1_{H})%
\end{array}%
=B(gx_{1}x_{2}\otimes g;1_{A},1_{H})
\end{equation*}%
By applying $\left( \ref{eq.10}\right) $ this rewrites as%
\begin{equation*}
\begin{array}{c}
B(x_{1}x_{2}\otimes 1_{H};1_{A},1_{H}) \\
+\gamma _{1}B(x_{1}x_{2}\otimes 1_{H};GX_{1},1_{H})+\gamma
_{2}B(x_{1}x_{2}\otimes 1_{H};GX_{2},1_{H})%
\end{array}%
=B(x_{1}x_{2}\otimes 1_{H};1_{A},1_{H})
\end{equation*}%
i.e.%
\begin{equation}
\gamma _{1}B(x_{1}x_{2}\otimes 1_{H};GX_{1},1_{H})+\gamma
_{2}B(x_{1}x_{2}\otimes 1_{H};GX_{2},1_{H})=0  \label{G,x1x2, GF2,1H}
\end{equation}

\subparagraph{Case $f=x_{1}x_{2}$}

\begin{equation*}
\begin{array}{c}
B(x_{1}x_{2}\otimes 1_{H};1_{A},x_{1}x_{2}) \\
+\gamma _{1}B(x_{1}x_{2}\otimes 1_{H};GX_{1},x_{1}x_{2})+\gamma
_{2}B(x_{1}x_{2}\otimes 1_{H};GX_{2},x_{1}x_{2})%
\end{array}%
=B(gx_{1}x_{2}\otimes g;1_{A},x_{1}x_{2})
\end{equation*}%
By applying $\left( \ref{eq.10}\right) $ this rewrites as%
\begin{equation*}
\begin{array}{c}
B(x_{1}x_{2}\otimes 1_{H};1_{A},x_{1}x_{2}) \\
+\gamma _{1}B(x_{1}x_{2}\otimes 1_{H};GX_{1},x_{1}x_{2})+\gamma
_{2}B(x_{1}x_{2}\otimes 1_{H};GX_{2},x_{1}x_{2})%
\end{array}%
=B(x_{1}x_{2}\otimes 1_{H};1_{A},x_{1}x_{2})
\end{equation*}%
i.e.%
\begin{equation}
\gamma _{1}B(x_{1}x_{2}\otimes 1_{H};GX_{1},x_{1}x_{2})+\gamma
_{2}B(x_{1}x_{2}\otimes 1_{H};GX_{2},x_{1}x_{2})=0  \label{G,x1x2, GF2,x1x2}
\end{equation}

\subparagraph{Case $f=gx_{1}$}

\begin{equation*}
\begin{array}{c}
B(x_{1}x_{2}\otimes 1_{H};1_{A},gx_{1}) \\
+\gamma _{1}B(x_{1}x_{2}\otimes 1_{H};GX_{1},gx_{1})+\gamma
_{2}B(x_{1}x_{2}\otimes 1_{H};GX_{2},gx_{1})%
\end{array}%
=B(gx_{1}x_{2}\otimes g;1_{A},gx_{1}).
\end{equation*}%
By applying $\left( \ref{eq.10}\right) $ this rewrites as%
\begin{equation*}
\begin{array}{c}
B(x_{1}x_{2}\otimes 1_{H};1_{A},gx_{1}) \\
+\gamma _{1}B(x_{1}x_{2}\otimes 1_{H};GX_{1},gx_{1})+\gamma
_{2}B(x_{1}x_{2}\otimes 1_{H};GX_{2},gx_{1})%
\end{array}%
=-B(x_{1}x_{2}\otimes 1_{H};1_{A},gx_{1}).
\end{equation*}%
i.e.%
\begin{equation}
\begin{array}{c}
2B(x_{1}x_{2}\otimes 1_{H};1_{A},gx_{1}) \\
+\gamma _{1}B(x_{1}x_{2}\otimes 1_{H};GX_{1},gx_{1})+\gamma
_{2}B(x_{1}x_{2}\otimes 1_{H};GX_{2},gx_{1})%
\end{array}%
=0.  \label{G,x1x2, GF2,gx1}
\end{equation}

\subparagraph{Case $f=gx_{2}$}

\begin{equation*}
\begin{array}{c}
B(x_{1}x_{2}\otimes 1_{H};1_{A},gx_{2}) \\
+\gamma _{1}B(x_{1}x_{2}\otimes 1_{H};GX_{1},gx_{2})+\gamma
_{2}B(x_{1}x_{2}\otimes 1_{H};GX_{2},gx_{2})%
\end{array}%
=B(gx_{1}x_{2}\otimes g;1_{A},gx_{2}).
\end{equation*}%
By applying $\left( \ref{eq.10}\right) $ this rewrites as%
\begin{equation*}
\begin{array}{c}
B(x_{1}x_{2}\otimes 1_{H};1_{A},gx_{2}) \\
+\gamma _{1}B(x_{1}x_{2}\otimes 1_{H};GX_{1},gx_{2})+\gamma
_{2}B(x_{1}x_{2}\otimes 1_{H};GX_{2},gx_{2})%
\end{array}%
=-B(x_{1}x_{2}\otimes 1_{H};1_{A},gx_{2})
\end{equation*}%
i.e.%
\begin{equation}
\begin{array}{c}
2B(x_{1}x_{2}\otimes 1_{H};1_{A},gx_{2}) \\
+\gamma _{1}B(x_{1}x_{2}\otimes 1_{H};GX_{1},gx_{2})+\gamma
_{2}B(x_{1}x_{2}\otimes 1_{H};GX_{2},gx_{2})%
\end{array}%
=0  \label{G,x1x2, GF2,gx2}
\end{equation}

\paragraph{\textbf{Equality (\protect\ref{GF3})}}

rewrites as%
\begin{equation*}
\gamma _{2}B(x_{1}x_{2}\otimes 1_{H};X_{1}X_{2},f)-\alpha
B(x_{1}x_{2}\otimes 1_{H};GX_{1},f)=\alpha B(x_{1}x_{2}g\otimes g;GX_{1},f)
\end{equation*}%
i.e.%
\begin{equation*}
\gamma _{2}B(x_{1}x_{2}\otimes 1_{H};X_{1}X_{2},f)-\alpha
B(x_{1}x_{2}\otimes 1_{H};GX_{1},f)=\alpha B(gx_{1}x_{2}\otimes g;GX_{1},f).
\end{equation*}

\subparagraph{Case $f=1_{H}$}

\begin{equation*}
\gamma _{2}B(x_{1}x_{2}\otimes 1_{H};X_{1},1_{H})-\alpha B(x_{1}x_{2}\otimes
1_{H};GX_{1},1_{H})=\alpha B(gx_{1}x_{2}\otimes g;GX_{1},1_{H}).
\end{equation*}%
By applying $\left( \ref{eq.10}\right) $ this rewrites as%
\begin{equation*}
\gamma _{2}B(x_{1}x_{2}\otimes 1_{H};X_{1},1_{H})-\alpha B(x_{1}x_{2}\otimes
1_{H};GX_{1},1_{H})=\alpha B(x_{1}x_{2}\otimes 1_{H};GX_{1},1_{H}).
\end{equation*}%
i.e.%
\begin{equation*}
\gamma _{2}B(x_{1}x_{2}\otimes 1_{H};X_{1},1_{H})-2\alpha
B(x_{1}x_{2}\otimes 1_{H};GX_{1},1_{H})=0.
\end{equation*}%
In view of the form of the element we get%
\begin{equation}
2\alpha B(x_{1}x_{2}\otimes 1_{H};GX_{1},1_{H})=0.  \label{G,x1x2, GF3,1H}
\end{equation}

\subparagraph{Case $f=x_{1}x_{2}$}

\begin{equation*}
\gamma _{2}B(x_{1}x_{2}\otimes 1_{H};X_{1}X_2,x_{1}x_{2})-\alpha
B(x_{1}x_{2}\otimes 1_{H};GX_{1},x_{1}x_{2})=\alpha B(gx_{1}x_{2}\otimes
g;GX_{1},x_{1}x_{2}).
\end{equation*}

By applying $\left( \ref{eq.10}\right) $ this rewrites as%
\begin{equation*}
\gamma _{2}B(x_{1}x_{2}\otimes 1_{H};X_{1}X_{2},x_{1}x_{2})-\alpha
B(x_{1}x_{2}\otimes 1_{H};GX_{1},x_{1}x_{2})=\alpha B(x_{1}x_{2}\otimes
1_{H};GX_{1},x_{1}x_{2})
\end{equation*}%
i.e.%
\begin{equation*}
\gamma _{2}B(x_{1}x_{2}\otimes 1_{H};X_{1}X_{2},x_{1}x_{2})-2\alpha
B(x_{1}x_{2}\otimes 1_{H};GX_{1},x_{1}x_{2})=0.
\end{equation*}%
In view of the form of the element we get%
\begin{equation}
-2\alpha B(x_{1}x_{2}\otimes 1_{H};GX_{1},x_{1}x_{2})=0.
\label{G,x1x2, GF3,x1x2}
\end{equation}

\subparagraph{Case $f=gx_{1}$}

\begin{equation*}
\gamma _{2}B(x_{1}x_{2}\otimes 1_{H};X_{1}X_2,gx_{1})-\alpha
B(x_{1}x_{2}\otimes 1_{H};GX_{1},gx_{1})=\alpha B(gx_{1}x_{2}\otimes
g;GX_{1},gx_{1}).
\end{equation*}

By applying $\left( \ref{eq.10}\right) $ this rewrites as%
\begin{equation*}
\gamma _{2}B(x_{1}x_{2}\otimes 1_{H};X_{1}X_2,gx_{1})-\alpha
B(x_{1}x_{2}\otimes 1_{H};GX_{1},gx_{1})=-\alpha B(x_{1}x_{2}\otimes
1_{H};GX_{1},gx_{1})
\end{equation*}%
i.e.%
\begin{equation*}
\gamma _{2}B(x_{1}x_{2}\otimes 1_{H};X_{1}X_2,gx_{1})=0
\end{equation*}

in view of the form of the element we get

\begin{equation}
\gamma _{2}B(x_{1}x_{2}\otimes 1_{H};X_{1},gx_{1}x_{2})=0
\label{G,x1x2, GF3,gx1}
\end{equation}

\subparagraph{Case $f=gx_{2}$}

\begin{equation*}
\gamma _{2}B(x_{1}x_{2}\otimes 1_{H};X_{1}X_2,gx_{2})-\alpha
B(x_{1}x_{2}\otimes 1_{H};GX_{1},gx_{2})=\alpha B(gx_{1}x_{2}\otimes
g;GX_{1},gx_{2}).
\end{equation*}

By applying $\left( \ref{eq.10}\right) $ this rewrites as%
\begin{equation*}
\gamma _{2}B(x_{1}x_{2}\otimes 1_{H};X_{1}X_2,gx_{2})-\alpha
B(x_{1}x_{2}\otimes 1_{H};GX_{1},gx_{2})=-\alpha B(x_{1}x_{2}\otimes
1_{H};GX_{1},gx_{2})
\end{equation*}%
i.e.%
\begin{equation*}
\gamma _{2}B(x_{1}x_{2}\otimes 1_{H};X_{1}X_2,gx_{2})=0
\end{equation*}%
in view of the form of the element we get

\begin{equation}
\gamma _{2}B(x_{1}x_{2}\otimes 1_{H};X_{2},gx_{1}x_{2})=0
\label{G,x1x2, GF3,gx2}
\end{equation}

\paragraph{\textbf{Equality (\protect\ref{GF4})}}

rewrites as%
\begin{equation*}
-\gamma _{1}B(x_{1}x_{2}\otimes 1_{H};X_1X_{2},f)-\alpha B(x_{1}x_{2}\otimes
1_{H};GX_{2},f)=\alpha B(x_{1}x_{2}g\otimes g;GX_{2},f)
\end{equation*}%
i.e.%
\begin{equation*}
-\gamma _{1}B(x_{1}x_{2}\otimes 1_{H};X_1X_{2},f)-\alpha B(x_{1}x_{2}\otimes
1_{H};GX_{2},f)=\alpha B(gx_{1}x_{2}\otimes g;GX_{2},f).
\end{equation*}

\subparagraph{Case $f=1_{H}$}

\begin{equation*}
-\gamma _{1}B(x_{1}x_{2}\otimes 1_{H};X_1X_{2},1_{H})-\alpha
B(x_{1}x_{2}\otimes 1_{H};GX_{2},1_{H})=\alpha B(gx_{1}x_{2}\otimes
g;GX_{2},1_{H}).
\end{equation*}

By applying $\left( \ref{eq.10}\right) $ this rewrites as%
\begin{equation*}
-\gamma _{1}B(x_{1}x_{2}\otimes 1_{H};X_1X_{2},1_{H})-\alpha
B(x_{1}x_{2}\otimes 1_{H};GX_{2},1_{H})=\alpha B(x_{1}x_{2}\otimes
1_{H};GX_{2},1_{H})
\end{equation*}%
i.e.%
\begin{equation*}
-\gamma _{1}B(x_{1}x_{2}\otimes 1_{H};X_1X_{2},1_{H})-2\alpha
B(x_{1}x_{2}\otimes 1_{H};GX_{2},1_{H})=0.
\end{equation*}%
In view of the form of the element we get

\begin{gather}
-\gamma _{1}[+1-B(x_{1}x_{2}\otimes
1_{H};1_{A},x_{1}x_{2})-B(x_{1}x_{2}\otimes 1_{H};X_{2},x_{1})
\label{G,x1x2, GF4,1H} \\
+B(x_{1}x_{2}\otimes 1_{H};X_{1},x_{2})]-2\alpha B(x_{1}x_{2}\otimes
1_{H};GX_{2},1_{H})=0.  \notag
\end{gather}

\subparagraph{Case $f=x_{1}x_{2}$}

\begin{equation*}
-\gamma _{1}B(x_{1}x_{2}\otimes 1_{H};X_1X_{2},x_{1}x_{2})-\alpha
B(x_{1}x_{2}\otimes 1_{H};GX_{2},x_{1}x_{2})=\alpha B(gx_{1}x_{2}\otimes
g;GX_{2},x_{1}x_{2}).
\end{equation*}

By applying $\left( \ref{eq.10}\right) $ this rewrites as%
\begin{equation*}
-\gamma _{1}B(x_{1}x_{2}\otimes 1_{H};X_{1}X_{2},x_{1}x_{2})-\alpha
B(x_{1}x_{2}\otimes 1_{H};GX_{2},x_{1}x_{2})=\alpha B(x_{1}x_{2}\otimes
1_{H};GX_{2},x_{1}x_{2})
\end{equation*}%
i.e.%
\begin{equation*}
-\gamma _{1}B(x_{1}x_{2}\otimes 1_{H};X_{2},x_{1}x_{2})-2\alpha
B(x_{1}x_{2}\otimes 1_{H};GX_{2},x_{1}x_{2})=0.
\end{equation*}%
In view of the form of the element we get%
\begin{equation}
\alpha B(x_{1}x_{2}\otimes 1_{H};GX_{2},x_{1}x_{2})=0.
\label{G,x1x2, GF4,x1x2}
\end{equation}

\subparagraph{Case $f=gx_{1}$}

\begin{equation*}
-\gamma _{1}B(x_{1}x_{2}\otimes 1_{H};X_1X_{2},gx_{1})-\alpha
B(x_{1}x_{2}\otimes 1_{H};GX_{2},gx_{1})=\alpha B(gx_{1}x_{2}\otimes
g;GX_{2},gx_{1}).
\end{equation*}

By applying $\left( \ref{eq.10}\right) $ this rewrites as%
\begin{equation*}
-\gamma _{1}B(x_{1}x_{2}\otimes 1_{H};X_1X_{2},gx_{1})-\alpha
B(x_{1}x_{2}\otimes 1_{H};GX_{2},gx_{1})=-\alpha B(gx_{1}x_{2}\otimes
g;GX_{2},gx_{1}).
\end{equation*}%
i.e.%
\begin{equation*}
\gamma _{1}B(x_{1}x_{2}\otimes 1_{H};X_1X_{2},gx_{1})=0
\end{equation*}
by the form of the element we get

\begin{equation}
\gamma _{1}B(x_{1}x_{2}\otimes 1_{H};X_{1},gx_{1}x_{2})=0
\label{G,x1x2, GF4,gx1}
\end{equation}

\subparagraph{Case $f=gx_{2}$}

\begin{equation*}
-\gamma _{1}B(x_{1}x_{2}\otimes 1_{H};X_1X_{2},gx_{2})-\alpha
B(x_{1}x_{2}\otimes 1_{H};GX_{2},gx_{2})=\alpha B(gx_{1}x_{2}\otimes
g;GX_{2},gx_{2}).
\end{equation*}

By applying $\left( \ref{eq.10}\right) $ this rewrites as%
\begin{equation*}
-\gamma _{1}B(x_{1}x_{2}\otimes 1_{H};X_1X_{2},gx_{2})-\alpha
B(x_{1}x_{2}\otimes 1_{H};GX_{2},gx_{2})=-\alpha B(gx_{1}x_{2}\otimes
g;GX_{2},gx_{2}).
\end{equation*}%
i.e.%
\begin{equation*}
\gamma _{1}B(x_{1}x_{2}\otimes 1_{H};X_1X_{2},gx_{2})=0
\end{equation*}

by the form of the element we get

\begin{equation}
\gamma _{1}B(x_{1}x_{2}\otimes 1_{H};X_{2},gx_{1}x_{2})=0.
\label{G,x1x2, GF4,gx2}
\end{equation}

\paragraph{\textbf{Equality (\protect\ref{GF5})}}

rewrites as%
\begin{equation*}
\alpha B(x_{1}x_{2}\otimes 1_{H};GX_{1}X_{2},f)=\alpha B(x_{1}x_{2}g\otimes
g;GX_{1}X_{2},f)
\end{equation*}%
i.e.%
\begin{equation*}
\alpha B(x_{1}x_{2}\otimes 1_{H};GX_{1}X_{2},f)=\alpha B(gx_{1}x_{2}\otimes
g;GX_{1}X_{2},f).
\end{equation*}

\subparagraph{Case $f=g$}

\begin{equation*}
\alpha B(x_{1}x_{2}\otimes 1_{H};GX_{1}X_{2},g)=\alpha B(gx_{1}x_{2}\otimes
g;GX_{1}X_{2},g).
\end{equation*}%
By applying $\left( \ref{eq.10}\right) $ this rewrites as%
\begin{equation*}
\alpha B(x_{1}x_{2}\otimes 1_{H};GX_{1}X_{2},g)=\alpha B(x_{1}x_{2}\otimes
1_{H};GX_{1}X_{2},g).
\end{equation*}%
which is trivial.

\subparagraph{Case $f=x_{1}$}

\begin{equation*}
\alpha B(x_{1}x_{2}\otimes 1_{H};GX_{1}X_{2},x_{1})=\alpha
B(gx_{1}x_{2}\otimes g;GX_{1}X_{2},x_{1}).
\end{equation*}%
By applying $\left( \ref{eq.10}\right) $ this rewrites as%
\begin{equation*}
\alpha B(x_{1}x_{2}\otimes 1_{H};GX_{1}X_{2},x_{1})=-\alpha
B(x_{1}x_{2}\otimes 1_{H};GX_{1}X_{2},x_{1})
\end{equation*}%
i.e.%
\begin{equation*}
\alpha B(x_{1}x_{2}\otimes 1_{H};GX_{1}X_{2},x_{1})=0.
\end{equation*}%
In view of the form of the element we get%
\begin{equation*}
\alpha B(x_{1}x_{2}\otimes 1_{H};GX_{1},x_{1}x_{2})=0
\end{equation*}%
which is $\left( \ref{G,x1x2, GF3,x1x2}\right) .$

\subparagraph{Case $f=x_{2}$}

\begin{equation*}
\alpha B(x_{1}x_{2}\otimes 1_{H};GX_{1}X_{2},x_{2})=\alpha
B(gx_{1}x_{2}\otimes g;GX_{1}X_{2},x_{21}).
\end{equation*}%
By applying $\left( \ref{eq.10}\right) $ this rewrites as%
\begin{equation*}
\alpha B(x_{1}x_{2}\otimes 1_{H};GX_{1}X_{2},x_{2})=-\alpha
B(x_{1}x_{2}\otimes 1_{H};GX_{1}X_{2},x_{2})
\end{equation*}%
i.e.%
\begin{equation*}
\alpha B(x_{1}x_{2}\otimes 1_{H};GX_{1}X_{2},x_{2})=0.
\end{equation*}%
In view of the form of the element we get%
\begin{equation*}
\alpha B(x_{1}x_{2}\otimes 1_{H};GX_{2},x_{1}x_{2})=0
\end{equation*}%
which is $\left( \ref{G,x1x2, GF4,x1x2}\right) .$

\paragraph{\textbf{Equality \protect\ref{GF6}}}

rewrites as%
\begin{equation*}
-B(x_{1}x_{2}\otimes 1_{H};X_{1},f)+\gamma _{2}B(x_{1}x_{2}\otimes
1_{H};GX_{1}X_{2},f)=B(x_{1}x_{2}g\otimes g;X_{1},f)
\end{equation*}%
i.e.%
\begin{equation*}
-B(x_{1}x_{2}\otimes 1_{H};X_{1},f)+\gamma _{2}B(x_{1}x_{2}\otimes
1_{H};GX_{1}X_{2},f)=B(gx_{1}x_{2}\otimes g;X_{1},f).
\end{equation*}

\subparagraph{Case $f=g$}

\begin{equation*}
-B(x_{1}x_{2}\otimes 1_{H};X_{1},g)+\gamma _{2}B(x_{1}x_{2}\otimes
1_{H};GX_{1}X_{2},g)=B(gx_{1}x_{2}\otimes g;X_{1},g).
\end{equation*}%
By applying $\left( \ref{eq.10}\right) $ this rewrites as%
\begin{equation*}
-B(x_{1}x_{2}\otimes 1_{H};X_{1},g)+\gamma _{2}B(x_{1}x_{2}\otimes
1_{H};GX_{1}X_{2},g)=B(x_{1}x_{2}\otimes 1_{H};X_{1},g)
\end{equation*}%
i.e.%
\begin{equation*}
-2B(x_{1}x_{2}\otimes 1_{H};X_{1},g)+\gamma _{2}B(x_{1}x_{2}\otimes
1_{H};GX_{1}X_{2},g)=0.
\end{equation*}%
In view of the form of the element we get%
\begin{gather}
-2B(x_{1}x_{2}\otimes 1_{H};X_{1},g)+  \label{G,x1x2, GF6,g} \\
+\gamma _{2}[-B(x_{1}x_{2}\otimes 1_{H};g,gx_{1}x_{2})-B(x_{1}x_{2}\otimes
1_{H};GX_{1},gx_{2})+B(x_{1}x_{2}\otimes 1_{H};GX_{2},gx_{1})]=0  \notag
\end{gather}

\subparagraph{Case $f=x_{1}$}

\begin{equation*}
-B(x_{1}x_{2}\otimes 1_{H};X_{1},x_{1})+\gamma _{2}B(x_{1}x_{2}\otimes
1_{H};GX_{1}X_{2},x_{1})=B(gx_{1}x_{2}\otimes g;X_{1},x_{1}).
\end{equation*}%
By applying $\left( \ref{eq.10}\right) $ this rewrites as%
\begin{equation*}
-B(x_{1}x_{2}\otimes 1_{H};X_{1},x_{1})+\gamma _{2}B(x_{1}x_{2}\otimes
1_{H};GX_{1}X_{2},x_{1})=-B(x_{1}x_{2}\otimes 1_{H};X_{1},x_{1})
\end{equation*}%
i.e.%
\begin{equation}
\gamma _{2}B(x_{1}x_{2}\otimes 1_{H};GX_{1},x_{1}x_{2})=0
\label{G,x1x2, GF6,x1}
\end{equation}

\subparagraph{Case $f=x_{2}$}

\begin{equation*}
-B(x_{1}x_{2}\otimes 1_{H};X_{1},x_{2})+\gamma _{2}B(x_{1}x_{2}\otimes
1_{H};GX_{1}X_{2},x_{2})=B(gx_{1}x_{2}\otimes g;X_{1},x_{2}).
\end{equation*}%
By applying $\left( \ref{eq.10}\right) $ this rewrites as%
\begin{equation*}
-B(x_{1}x_{2}\otimes 1_{H};X_{1},x_{2})+\gamma _{2}B(x_{1}x_{2}\otimes
1_{H};GX_{1}X_{2},x_{2})=-B(x_{1}x_{2}\otimes 1_{H};X_{1},x_{2})
\end{equation*}%
i.e.%
\begin{equation}
\gamma _{2}B(x_{1}x_{2}\otimes 1_{H};GX_{2},x_{1}x_{2})=0.
\label{G,x1x2, GF6,x2}
\end{equation}

\subparagraph{Case $f=gx_{1}x_{2}$}

\begin{equation*}
-B(x_{1}x_{2}\otimes 1_{H};X_{1},gx_{1}x_{2})+\gamma _{2}B(x_{1}x_{2}\otimes
1_{H};GX_{1}X_{2},gx_{1}x_{2})=B(gx_{1}x_{2}\otimes g;X_{1},gx_{1}x_{2}).
\end{equation*}%
By applying $\left( \ref{eq.10}\right) $ this rewrites as%
\begin{equation*}
-B(x_{1}x_{2}\otimes 1_{H};X_{1},gx_{1}x_{2})+\gamma _{2}B(x_{1}x_{2}\otimes
1_{H};GX_{1}X_{2},gx_{1}x_{2})=B(x_{1}x_{2}\otimes 1_{H};X_{1},gx_{1}x_{2})
\end{equation*}%
i.e.%
\begin{equation*}
-2B(x_{1}x_{2}\otimes 1_{H};X_{1},gx_{1}x_{2})+\gamma
_{2}B(x_{1}x_{2}\otimes 1_{H};GX_{1}X_{2},gx_{1}x_{2})=0.
\end{equation*}%
In view of the form of the element we get%
\begin{equation}
B(x_{1}x_{2}\otimes 1_{H};X_{1},gx_{1}x_{2})=0  \label{G,x1x2, GF6,gx1x2}
\end{equation}

\paragraph{\textbf{Equality \protect\ref{GF7}}}

rewrites as
\begin{equation*}
-B(x_{1}x_{2}\otimes 1_{H};X_{2},f)-\gamma _{1}B(x_{1}x_{2}\otimes
1_{H};GX_{1}X_{2},f)=B(x_{1}x_{2}g\otimes g;X_{2},f)
\end{equation*}%
i.e.%
\begin{equation*}
-B(x_{1}x_{2}\otimes 1_{H};X_{2},f)-\gamma _{1}B(x_{1}x_{2}\otimes
1_{H};GX_{1}X_{2},f)=B(gx_{1}x_{2}\otimes g;X_{2},f).
\end{equation*}

\subparagraph{Case $f=g$}

\begin{equation*}
-B(x_{1}x_{2}\otimes 1_{H};X_{2},g)-\gamma _{1}B(x_{1}x_{2}\otimes
1_{H};GX_{1}X_{2},g)=B(gx_{1}x_{2}\otimes g;X_{2},g).
\end{equation*}%
By applying $\left( \ref{eq.10}\right) $ this rewrites as%
\begin{equation*}
-B(x_{1}x_{2}\otimes 1_{H};X_{2},g)-\gamma _{1}B(x_{1}x_{2}\otimes
1_{H};GX_{1}X_{2},g)=B(x_{1}x_{2}\otimes 1_{H};X_{2},g)
\end{equation*}%
i.e.%
\begin{equation*}
-2B(x_{1}x_{2}\otimes 1_{H};X_{2},g)-\gamma _{1}B(x_{1}x_{2}\otimes
1_{H};GX_{1}X_{2},g)=0.
\end{equation*}%
In view of the form of the element we get%
\begin{gather}
-2B(x_{1}x_{2}\otimes 1_{H};X_{2},g)+  \label{G,x1x2, GF7,g} \\
-\gamma _{1}[-B(x_{1}x_{2}\otimes 1_{H};g,gx_{1}x_{2})-B(x_{1}x_{2}\otimes
1_{H};GX_{1},gx_{2})+B(x_{1}x_{2}\otimes 1_{H};GX_{2},gx_{1})]=0  \notag
\end{gather}

\subparagraph{Case $f=x_{1}$}

\begin{equation*}
-B(x_{1}x_{2}\otimes 1_{H};X_{2},x_{1})-\gamma _{1}B(x_{1}x_{2}\otimes
1_{H};GX_1X_{2},x_{1})=B(gx_{1}x_{2}\otimes g;X_{2},x_{1}).
\end{equation*}%
By applying $\left( \ref{eq.10}\right) $ this rewrites as%
\begin{equation*}
-B(x_{1}x_{2}\otimes 1_{H};X_{2},x_{1})-\gamma _{1}B(x_{1}x_{2}\otimes
1_{H};GX_1X_{2},x_{1})=-B(x_{1}x_{2}\otimes 1_{H};X_{2},x_{1})
\end{equation*}%
i.e.%
\begin{equation*}
\gamma _{1}B(x_{1}x_{2}\otimes 1_{H};GX_1X_{2},x_{1})=0
\end{equation*}

In view of the form of the element we get

\begin{equation}
\gamma _{1}B(x_{1}x_{2}\otimes 1_{H};GX_{1},x_{1}x_{2})=0
\label{G,x1x2, GF7,x1}
\end{equation}

\subparagraph{Case $f=x_{2}$}

\begin{equation*}
-B(x_{1}x_{2}\otimes 1_{H};X_{2},x_{2})-\gamma _{1}B(x_{1}x_{2}\otimes
1_{H};GX_{1}X_{2},x_{2})=B(gx_{1}x_{2}\otimes g;X_{2},x_{2}).
\end{equation*}%
By applying $\left( \ref{eq.10}\right) $ this rewrites as%
\begin{equation*}
-B(x_{1}x_{2}\otimes 1_{H};X_{2},x_{2})-\gamma _{1}B(x_{1}x_{2}\otimes
1_{H};GX_{1}X_{2},x_{2})=-B(x_{1}x_{2}\otimes 1_{H};X_{2},x_{2})
\end{equation*}%
i.e.%
\begin{equation}
\gamma _{1}B(x_{1}x_{2}\otimes 1_{H};GX_{2},x_{1}x_{2})=0
\label{G,x1x2, GF7,x2}
\end{equation}

\subparagraph{Case $f=gx_{1}x_{2}$}

\begin{equation*}
-B(x_{1}x_{2}\otimes 1_{H};X_{2},gx_{1}x_{2})-\gamma _{1}B(x_{1}x_{2}\otimes
1_{H};GX_{1}X_{2},gx_{1}x_{2})=B(gx_{1}x_{2}\otimes g;X_{2},gx_{1}x_{2}).
\end{equation*}%
By applying $\left( \ref{eq.10}\right) $ this rewrites as%
\begin{equation*}
-B(x_{1}x_{2}\otimes 1_{H};X_{2},gx_{1}x_{2})-\gamma _{1}B(x_{1}x_{2}\otimes
1_{H};GX_{1}X_{2},gx_{1}x_{2})=B(x_{1}x_{2}\otimes 1_{H};X_{2},gx_{1}x_{2})
\end{equation*}%
i.e.%
\begin{equation*}
-2B(x_{1}x_{2}\otimes 1_{H};X_{2},gx_{1}x_{2})-\gamma
_{1}B(x_{1}x_{2}\otimes 1_{H};GX_{1}X_{2},gx_{1}x_{2})=0.
\end{equation*}%
In view of the form of the element we get%
\begin{equation}
B(x_{1}x_{2}\otimes 1_{H};X_{2},gx_{1}x_{2})=0  \label{G,x1x2, GF7,gx1x2}
\end{equation}

\paragraph{\textbf{Equality} $\left( \protect\ref{GF8}\right) $}

rewrites as%
\begin{equation*}
B(x_{1}x_{2}\otimes 1_{H};X_{1}X_{2},f)=B(x_{1}x_{2}g\otimes g;X_{1}X_{2},f)
\end{equation*}%
i.e.%
\begin{equation*}
B(x_{1}x_{2}\otimes 1_{H};X_{1}X_{2},f)=B(gx_{1}x_{2}\otimes g;X_{1}X_{2},f)
\end{equation*}

\subparagraph{Case $f=1_{H}$}

\begin{equation*}
B(x_{1}x_{2}\otimes 1_{H};X_{1}X_{2},1_{H})=B(gx_{1}x_{2}\otimes
g;X_{1}X_{2},1_{H}).
\end{equation*}%
By applying $\left( \ref{eq.10}\right) $ this rewrites as%
\begin{equation*}
B(x_{1}x_{2}\otimes 1_{H};X_{1}X_{2},1_{H})=B(x_{1}x_{2}\otimes
1_{H};X_{1}X_{2},1_{H}).
\end{equation*}%
which is trivial.

\subparagraph{Case $f=gx_{1}$}

\begin{equation*}
B(x_{1}x_{2}\otimes 1_{H};X_{1}X_{2},gx_{1})=B(gx_{1}x_{2}\otimes
g;X_{1}X_{2},gx_{1})
\end{equation*}%
By applying $\left( \ref{eq.10}\right) $ this rewrites as%
\begin{equation*}
B(x_{1}x_{2}\otimes 1_{H};X_{1}X_{2},gx_{1})=-B(x_{1}x_{2}\otimes
1_{H};X_{1}X_{2},gx_{1})
\end{equation*}%
i.e.%
\begin{equation*}
2B(x_{1}x_{2}\otimes 1_{H};X_{1}X_{2},gx_{1})=0.
\end{equation*}%
In view of the form of the element we get%
\begin{equation*}
B(x_{1}x_{2}\otimes 1_{H};X_{1},gx_{1}x_{2})=0
\end{equation*}%
which is $\left( \ref{G,x1x2, GF6,gx1x2}\right) .$

\subparagraph{Case $f=gx_{2}$}

\begin{equation*}
B(x_{1}x_{2}\otimes 1_{H};X_{1}X_{2},gx_{2})=B(gx_{1}x_{2}\otimes
g;X_{1}X_{2},gx_{2})
\end{equation*}%
By applying $\left( \ref{eq.10}\right) $ this rewrites as%
\begin{equation*}
B(x_{1}x_{2}\otimes 1_{H};X_{1}X_{2},gx_{2})=-B(x_{1}x_{2}\otimes
1_{H};X_{1}X_{2},gx_{2})
\end{equation*}%
i.e.%
\begin{equation*}
2B(x_{1}x_{2}\otimes 1_{H};X_{1}X_{2},gx_{2})=0.
\end{equation*}%
In view of the form of the element we get%
\begin{equation*}
B(x_{1}x_{2}\otimes 1_{H};X_{2},gx_{1}x_{2})=0
\end{equation*}%
which is $\left( \ref{G,x1x2, GF7,gx1x2}\right) .$

\subsubsection{Case $gx_{1}\otimes 1_{H}$}

\paragraph{\textbf{Equality} (\protect\ref{GF1})}

rewrites as%
\begin{equation*}
\begin{array}{c}
\alpha B(gx_{1}\otimes 1_{H};G,f)+ \\
+\gamma _{1}B(gx_{1}\otimes 1_{H};X_{1},f)+\gamma _{2}B(gx_{1}\otimes
1_{H};X_{2},f)%
\end{array}%
=\alpha B(gx_{1}g\otimes g;G,f)
\end{equation*}%
i.e.%
\begin{equation*}
\begin{array}{c}
\alpha B(gx_{1}\otimes 1_{H};G,f)+ \\
+\gamma _{1}B(gx_{1}\otimes 1_{H};X_{1},f)+\gamma _{2}B(gx_{1}\otimes
1_{H};X_{2},f)%
\end{array}%
=-\alpha B(x_{1}\otimes g;G,f).
\end{equation*}

\subparagraph{Case $f=1_H$}

\begin{equation*}
\begin{array}{c}
\alpha B(gx_{1}\otimes 1_H;G,1_H)+ \\
+\gamma _{1}B(gx_{1}\otimes 1_{H};X_{1},1_H)+\gamma _{2}B(gx_{1}\otimes
1_{H};X_{2},1_H)%
\end{array}%
=-\alpha B(x_{1}\otimes g;G,1_H).
\end{equation*}%
By applying $\left( \ref{eq.10}\right) $ this rewrites as%
\begin{equation*}
\begin{array}{c}
\alpha B(gx_{1}\otimes 1_{H};G,1_H)+ \\
+\gamma _{1}B(gx_{1}\otimes 1_{H};X_{1},1_H)+\gamma _{2}B(gx_{1}\otimes
1_{H};X_{2},1_H)%
\end{array}%
=-\alpha B(gx_{1}\otimes 1_{H};G,1_H)
\end{equation*}%
i.e.%
\begin{equation*}
\begin{array}{c}
2\alpha B(gx_{1}\otimes 1_{H};G,1_H)+ \\
+\gamma _{1}B(gx_{1}\otimes 1_{H};X_{1},1_H)+\gamma _{2}B(gx_{1}\otimes
1_{H};X_{2},1_H)%
\end{array}%
=0.
\end{equation*}%
In view of the form of the element we get%
\begin{equation}
\begin{array}{c}
2\alpha \left[ B(x_{1}x_{2}\otimes 1_{H};G,x_{2})-B(x_{1}x_{2}\otimes
1_{H};GX_{2},1_{H})\right] + \\
+\gamma _{1}\left[ -1+B(x_{1}x_{2}\otimes
1_{H};1_{A},x_{1}x_{2})+B(x_{1}x_{2}\otimes 1_{H};X_{2},x_{1})\right]
+\gamma _{2}B(x_{1}x_{2}\otimes 1_{H};X_{2},x_{2})%
\end{array}%
=0.  \label{G,gx1, GF1,1H}
\end{equation}

\subparagraph{Case $f=x_{1}x_{2}$}

\begin{equation*}
\begin{array}{c}
\alpha B(gx_{1}\otimes 1_{H};G,x_{1}x_{2})+ \\
+\gamma _{1}B(gx_{1}\otimes 1_{H};X_{1},x_{1}x_{2})+\gamma
_{2}B(gx_{1}\otimes 1_{H};X_{2},x_{1}x_{2})%
\end{array}%
=-\alpha B(x_{1}\otimes g;G,x_{1}x_{2}).
\end{equation*}%
By applying $\left( \ref{eq.10}\right) $ this rewrites as%
\begin{equation*}
\begin{array}{c}
\alpha B(gx_{1}\otimes 1_{H};G,x_{1}x_{2})+ \\
+\gamma _{1}B(gx_{1}\otimes 1_{H};X_{1},x_{1}x_{2})+\gamma
_{2}B(gx_{1}\otimes 1_{H};X_{2},x_{1}x_{2})%
\end{array}%
=-\alpha B(gx_{1}\otimes 1_{H};G,x_{1}x_{2}).
\end{equation*}%
i.e.%
\begin{equation*}
\begin{array}{c}
2\alpha B(gx_{1}\otimes 1_{H};G,x_{1}x_{2})+ \\
+\gamma _{1}B(gx_{1}\otimes 1_{H};X_{1},x_{1}x_{2})+\gamma
_{2}B(gx_{1}\otimes 1_{H};X_{2},x_{1}x_{2})%
\end{array}%
=0.
\end{equation*}%
In view of the form of the element we get%
\begin{equation*}
-2\alpha B(x_{1}x_{2}\otimes 1_{H};GX_{2},x_{1}x_{2})=0
\end{equation*}%
i.e%
\begin{equation*}
\alpha B(x_{1}x_{2}\otimes 1_{H};GX_{2},x_{1}x_{2})=0
\end{equation*}%
This is $\left( \ref{G,x1x2, GF4,x1x2}\right) .$

\subparagraph{Case $f=gx_{1}$}

\begin{equation*}
\begin{array}{c}
\alpha B(gx_{1}\otimes 1_{H};G,gx_{1})+ \\
+\gamma _{1}B(gx_{1}\otimes 1_{H};X_{1},gx_{1})+\gamma _{2}B(gx_{1}\otimes
1_{H};X_{2},gx_{1})%
\end{array}%
=-\alpha B(x_{1}\otimes g;G,gx_{1}).
\end{equation*}%
By applying $\left( \ref{eq.10}\right) $ this rewrites as%
\begin{equation*}
\begin{array}{c}
\alpha B(gx_{1}\otimes 1_{H};G,gx_{1})+ \\
+\gamma _{1}B(gx_{1}\otimes 1_{H};X_{1},gx_{1})+\gamma _{2}B(gx_{1}\otimes
1_{H};X_{2},gx_{1})%
\end{array}%
=\alpha B(gx_{1}\otimes 1_{H};G,gx_{1})
\end{equation*}%
i.e.%
\begin{equation*}
\gamma _{1}B(gx_{1}\otimes 1_{H};X_{1},gx_{1})+\gamma _{2}B(gx_{1}\otimes
1_{H};X_{2},gx_{1})=0.
\end{equation*}%
In view of the form of the element we get

\begin{equation*}
\gamma _{2}B(x_{1}x_{2}\otimes 1_{H};X_{2},gx_{1}x_{2})=0
\end{equation*}%
which follows from $\left( \ref{G,x1x2, GF3,gx2}\right) .$

\subparagraph{Case $f=gx_{2}$}

\begin{equation*}
\begin{array}{c}
\alpha B(gx_{1}\otimes 1_{H};G,gx_{2})+ \\
+\gamma _{1}B(gx_{1}\otimes 1_{H};X_{1},gx_{2})+\gamma _{2}B(gx_{1}\otimes
1_{H};X_{2},gx_{2})%
\end{array}%
=-\alpha B(x_{1}\otimes g;G,gx_{2}).
\end{equation*}%
By applying $\left( \ref{eq.10}\right) $ this rewrites as%
\begin{equation*}
\begin{array}{c}
\alpha B(gx_{1}\otimes 1_{H};G,gx_{2})+ \\
+\gamma _{1}B(gx_{1}\otimes 1_{H};X_{1},gx_{2})+\gamma _{2}B(gx_{1}\otimes
1_{H};X_{2},gx_{2})%
\end{array}%
=\alpha B(gx_{1}\otimes 1_{H};G,gx_{2})
\end{equation*}%
i.e.%
\begin{equation*}
\gamma _{1}B(gx_{1}\otimes 1_{H};X_{1},gx_{2})+\gamma _{2}B(gx_{1}\otimes
1_{H};X_{2},gx_{2})=0.
\end{equation*}%
In view of the form of the element we get

\begin{equation*}
\gamma _{1}B(x_{1}x_{2}\otimes 1_{H};X_{2},gx_{1}x_{2})=0
\end{equation*}%
which follows from $\left( \ref{G,x1x2, GF4,gx2}\right) .$

\paragraph{\textbf{Equality (\protect\ref{GF2})}}

rewrites as%
\begin{equation*}
\begin{array}{c}
B(gx_{1}\otimes 1_{H};1_{A},f) \\
+\gamma _{1}B(gx_{1}\otimes 1_{H};GX_{1},f)+\gamma _{2}B(gx_{1}\otimes
1_{H};GX_{2},f)%
\end{array}%
=B(gx_{1}g\otimes g;1_{A},f)
\end{equation*}%
i.e.%
\begin{equation*}
\begin{array}{c}
B(gx_{1}\otimes 1_{H};1_{A},f) \\
+\gamma _{1}B(gx_{1}\otimes 1_{H};GX_{1},f)+\gamma _{2}B(gx_{1}\otimes
1_{H};GX_{2},f)%
\end{array}%
=-B(x_{1}\otimes g;1_{A},f).
\end{equation*}

\subparagraph{Case $f=g$}

\begin{equation*}
\begin{array}{c}
B(gx_{1}\otimes 1_{H};1_{A},g) \\
+\gamma _{1}B(gx_{1}\otimes 1_{H};GX_{1},g)+\gamma _{2}B(gx_{1}\otimes
1_{H};GX_{2},g)%
\end{array}%
=-B(x_{1}\otimes g;1_{A},g).
\end{equation*}%
By applying $\left( \ref{eq.10}\right) $ this rewrites as%
\begin{equation*}
\begin{array}{c}
B(gx_{1}\otimes 1_{H};1_{A},g) \\
+\gamma _{1}B(gx_{1}\otimes 1_{H};GX_{1},g)+\gamma _{2}B(gx_{1}\otimes
1_{H};GX_{2},g)%
\end{array}%
=-B(gx_{1}\otimes 1_{H};1_{A},g)
\end{equation*}%
i.e.In view of the form of the element we get%
\begin{equation}
\begin{array}{c}
2\left[ B(x_{1}x_{2}\otimes 1_{H};1_{A},gx_{2})+B(x_{1}x_{2}\otimes
1_{H};X_{2},g)\right] \\
\gamma _{1}\left[ -B(x_{1}x_{2}\otimes
1_{H};G,gx_{1}x_{2})+B(x_{1}x_{2}\otimes 1_{H};GX_{2},gx_{1})\right] +\gamma
_{2}B(x_{1}x_{2}\otimes 1_{H};GX_{2},gx_{2})%
\end{array}%
=0  \label{G,gx1, GF2,g}
\end{equation}

\subparagraph{Case $f=x_{1}$}

\begin{equation*}
\begin{array}{c}
B(gx_{1}\otimes 1_{H};1_{A},x_{1}) \\
+\gamma _{1}B(gx_{1}\otimes 1_{H};GX_{1},x_{1})+\gamma _{2}B(gx_{1}\otimes
1_{H};GX_{2},x_{1})%
\end{array}%
=-B(x_{1}\otimes g;1_{A},x_{1}).
\end{equation*}%
By applying $\left( \ref{eq.10}\right) $ this rewrites as

\begin{equation*}
\begin{array}{c}
B(gx_{1}\otimes 1_{H};1_{A},x_{1}) \\
+\gamma _{1}B(gx_{1}\otimes 1_{H};GX_{1},x_{1})+\gamma _{2}B(gx_{1}\otimes
1_{H};GX_{2},x_{1})%
\end{array}%
=B(gx_{1}\otimes 1_{H};1_{A},x_{1})
\end{equation*}%
i.e.%
\begin{equation*}
\gamma _{1}B(gx_{1}\otimes 1_{H};GX_{1},x_{1})+\gamma _{2}B(gx_{1}\otimes
1_{H};GX_{2},x_{1})=0.
\end{equation*}%
In view of the form of the element we get%
\begin{equation*}
\gamma _{2}B(x_{1}x_{2}\otimes 1_{H};GX_{2},x_{1}x_{2})=0
\end{equation*}%
which follows by $\left( \ref{G,x1x2, GF6,x2}\right) .$

\subparagraph{Case $f=x_{2}$}

\begin{equation*}
\begin{array}{c}
B(gx_{1}\otimes 1_{H};1_{A},x_{2}) \\
+\gamma _{1}B(gx_{1}\otimes 1_{H};GX_{1},x_{2})+\gamma _{2}B(gx_{1}\otimes
1_{H};GX_{2},x_{2})%
\end{array}%
=-B(x_{1}\otimes g;1_{A},x_{2}).
\end{equation*}%
By applying $\left( \ref{eq.10}\right) $ this rewrites as

\begin{equation*}
\begin{array}{c}
B(gx_{1}\otimes 1_{H};1_{A},x_{2}) \\
+\gamma _{1}B(gx_{1}\otimes 1_{H};GX_{1},x_{2})+\gamma _{2}B(gx_{1}\otimes
1_{H};GX_{2},x_{2})%
\end{array}%
=B(gx_{1}\otimes 1_{H};1_{A},x_{2})
\end{equation*}%
i.e.%
\begin{equation*}
\gamma _{1}B(gx_{1}\otimes 1_{H};GX_{1},x_{2})+\gamma _{2}B(gx_{1}\otimes
1_{H};GX_{2},x_{2})=0.
\end{equation*}%
In view of the form of the element we get%
\begin{equation*}
\gamma _{1}B(x_{1}x_{2}\otimes 1_{H};GX_{2},x_{1}x_{2})=0.
\end{equation*}%
which follows by $\left( \ref{G,x1x2, GF7,x2}\right) .$

\paragraph{\textbf{Equality (\protect\ref{GF3})}}

rewrites as%
\begin{equation*}
\gamma _{2}B(gx_{1}\otimes 1_{H};X_{1}X_{2},f)-\alpha B(gx_{1}\otimes
1_{H};GX_{1},f)=\alpha B(gx_{1}g\otimes g;GX_{1},f)
\end{equation*}%
i.e.%
\begin{equation*}
\gamma _{2}B(gx_{1}\otimes 1_{H};X_{1}X_{2},f)-\alpha B(gx_{1}\otimes
1_{H};GX_{1},f)=-\alpha B(x_{1}\otimes g;GX_{1},f).
\end{equation*}

\subparagraph{Case $f=g$}

\begin{equation*}
\gamma _{2}B(gx_{1}\otimes 1_{H};X_{1}X_2,g)-\alpha B(gx_{1}\otimes
1_{H};GX_{1},g)=-\alpha B(x_{1}\otimes g;GX_{1},g).
\end{equation*}%
By applying $\left( \ref{eq.10}\right) $ this rewrites as%
\begin{equation*}
\gamma _{2}B(gx_{1}\otimes 1_{H};X_{1}X_2,g)-\alpha B(gx_{1}\otimes
1_{H};GX_{1},g)=-\alpha B(gx_{1}\otimes 1_{H};GX_{1},g)
\end{equation*}%
i.e.%
\begin{equation*}
\gamma _{2}B(gx_{1}\otimes 1_{H};X_{1}X_2,g)=0.
\end{equation*}%
In view of the form of the element we get

\begin{equation*}
\gamma _{2}B(x_{1}x_2\otimes 1_{H};X_2,gx_1x_2)=0.
\end{equation*}%
This is (\ref{G,x1x2, GF3,gx2})

\subparagraph{Case $f=x_{2}$}

\begin{equation*}
\gamma _{2}B(gx_{1}\otimes 1_{H};X_{1}X_{2},x_{2})-\alpha B(gx_{1}\otimes
1_{H};GX_{1},x_{2})=-\alpha B(x_{1}\otimes g;GX_{1},x_{2}).
\end{equation*}%
By applying $\left( \ref{eq.10}\right) $ this rewrites as%
\begin{equation*}
\gamma _{2}B(gx_{1}\otimes 1_{H};X_{1}X_{2},x_{2})-\alpha B(gx_{1}\otimes
1_{H};GX_{1},x_{2})=\alpha B(gx_{1}\otimes 1_{H};GX_{1},x_{2})
\end{equation*}%
i.e.%
\begin{equation*}
\gamma _{2}B(gx_{1}\otimes 1_{H};X_{1}X_{2},x_{2})-2\alpha B(gx_{1}\otimes
1_{H};GX_{1},x_{2})=0.
\end{equation*}%
In view of the form of the element we get%
\begin{equation*}
\ -2\alpha B(x_{1}x_{2}\otimes 1_{H};GX_{2},x_{1}x_{2})=0.
\end{equation*}%
This is $\left( \ref{G,x1x2, GF4,x1x2}\right) $

\subparagraph{Case $f=gx_{2}$}

\begin{equation*}
\gamma _{2}B(gx_{1}\otimes 1_{H};X_{1}X_2,gx_{2})-\alpha B(gx_{1}\otimes
1_{H};GX_{1},gx_{2})=-\alpha B(x_{1}\otimes g;GX_{1},gx_{2}).
\end{equation*}%
By applying $\left( \ref{eq.10}\right) $ this rewrites as%
\begin{equation*}
\gamma _{2}B(gx_{1}\otimes 1_{H};X_{1}X_2,gx_{2})-\alpha B(gx_{1}\otimes
1_{H};GX_{1},gx_{2})=\alpha B(gx_{1}\otimes 1_{H};GX_{1},gx_{2})
\end{equation*}%
i.e.%
\begin{equation*}
\gamma _{2}B(gx_{1}\otimes 1_{H};X_{1}X_2,gx_{2})-2\alpha B(gx_{1}\otimes
1_{H};GX_{1},gx_{2})=0.
\end{equation*}%
In view of the form of the element we get nothing new.

\paragraph{\textbf{Equality (\protect\ref{GF4})}}

rewrites as%
\begin{equation*}
-\gamma _{1}B(gx_{1}\otimes 1_{H};X_{1}X_{2},f)-\alpha B(gx_{1}\otimes
1_{H};GX_{2},f)=\alpha B(gx_{1}g\otimes g;GX_{2},f)
\end{equation*}%
i.e.

\begin{equation*}
-\gamma _{1}B(gx_{1}\otimes 1_{H};X_1X_{2},f)-\alpha B(gx_{1}\otimes
1_{H};GX_{2},f)=-\alpha B(x_{1}\otimes g;GX_{2},f)
\end{equation*}

\subparagraph{Case $f=g$%
\protect\begin{equation*}
-\protect\gamma _{1}B(gx_{1}\otimes 1_{H};X_1X_{2},g)-\protect\alpha %
B(gx_{1}\otimes 1_{H};GX_{2},g)=-\protect\alpha B(x_{1}\otimes g;GX_{2},g).
\protect\end{equation*}%
}

By applying $\left( \ref{eq.10}\right) $ this rewrites as%
\begin{equation*}
-\gamma _{1}B(gx_{1}\otimes 1_{H};X_1X_{2},g)-\alpha B(gx_{1}\otimes
1_{H};GX_{2},g)=-\alpha B(gx_{1}\otimes 1_{H};GX_{2},g)
\end{equation*}%
i.e.%
\begin{equation*}
\gamma _{1}B(gx_{1}\otimes 1_{H};X_1X_{2},g)=0.
\end{equation*}%
In view of the form of the element we get

\begin{equation*}
\gamma _{1}B(x_{1}x_1\otimes 1_{H};X_{2},gx_1x_2)=0.
\end{equation*}%
This is (\ref{G,x1x2, GF4,gx2})

\subparagraph{Case $f=x_{1}$%
\protect\begin{equation*}
-\protect\gamma _{1}B(gx_{1}\otimes 1_{H};X_1X_{2},x_{1})-\protect\alpha %
B(gx_{1}\otimes 1_{H};GX_{2},x_{1})=-\protect\alpha B(x_{1}\otimes
g;GX_{2},x_{1}).
\protect\end{equation*}%
}

By applying $\left( \ref{eq.10}\right) $ this rewrites as%
\begin{equation*}
-\gamma _{1}B(gx_{1}\otimes 1_{H};X_1X_{2},x_{1})-\alpha B(gx_{1}\otimes
1_{H};GX_{2},x_{1})=+\alpha B(gx_{1}\otimes 1_{H};GX_{2},x_{1})
\end{equation*}%
i.e.%
\begin{equation*}
-\gamma _{1}B(gx_{1}\otimes 1_{H};X_1X_{2},x_{1})-2\alpha B(gx_{1}\otimes
1_{H};GX_{2},x_{1})=0
\end{equation*}%
In view of the form of the element we get%
\begin{equation*}
2\alpha B(x_{1}x_{2}\otimes 1_{H};GX_{2},x_{1}x_{2})=0
\end{equation*}%
which is $\left( \ref{G,x1x2, GF4,x1x2}\right) .$

\paragraph{\textbf{Equality (\protect\ref{GF5})}}

rewrites as%
\begin{equation*}
\alpha B(gx_{1}\otimes 1_{H};GX_{1}X_{2},f)=\alpha B(gx_{1}g\otimes
g;GX_{1}X_{2},f)
\end{equation*}%
i.e.%
\begin{equation*}
\alpha B(gx_{1}\otimes 1_{H};GX_{1}X_{2},f)=-\alpha B(x_{1}\otimes
g;GX_{1}X_{2},f).
\end{equation*}

\subparagraph{Case $f=1_{H}$}

\begin{equation*}
\alpha B(gx_{1}\otimes 1_{H};GX_{1}X_{2},1_{H})=-\alpha B(x_{1}\otimes
g;GX_{1}X_{2},1_{H}).
\end{equation*}%
By applying $\left( \ref{eq.10}\right) $ this rewrites as%
\begin{equation*}
\alpha B(gx_{1}\otimes 1_{H};GX_{1}X_{2},1_{H})=-\alpha B(gx_{1}\otimes
1_{H};GX_{1}X_{2},1_{H})
\end{equation*}%
i.e.%
\begin{equation*}
2\alpha B(gx_{1}\otimes 1_{H};GX_{1}X_{2},1_{H})=0.
\end{equation*}%
In view of the form of the element we get

\begin{equation*}
\alpha B(x_{1}x_{2}\otimes 1_{H};GX_{2},x_{1}x_{2})=0
\end{equation*}%
which follows from $\left( \ref{G,x1x2, GF4,x1x2}\right) .$

\paragraph{\textbf{Equality (\protect\ref{GF6})}}

rewrites as%
\begin{equation*}
-B(gx_{1}\otimes 1_{H};X_{1},f)+\gamma _{2}B(gx_{1}\otimes
1_{H};GX_{1}X_{2},f)=B(gx_{1}g\otimes g;X_{1},f)
\end{equation*}%
i.e.%
\begin{equation*}
-B(gx_{1}\otimes 1_{H};X_{1},f)+\gamma _{2}B(gx_{1}\otimes
1_{H};GX_{1}X_{2},f)=-B(x_{1}\otimes g;X_{1},f).
\end{equation*}

\subparagraph{Case $f=1_{H}$}

\begin{equation*}
-B(gx_{1}\otimes 1_{H};X_{1},1_{H})+\gamma _{2}B(gx_{1}\otimes
1_{H};GX_{1}X_2,1_{H})=-B(x_{1}\otimes g;X_{1},1_{H}).
\end{equation*}%
By applying $\left( \ref{eq.10}\right) $ this rewrites as%
\begin{equation*}
-B(gx_{1}\otimes 1_{H};X_{1},1_{H})+\gamma _{2}B(gx_{1}\otimes
1_{H};GX_{1}X_2,1_{H})=-B(gx_{1}\otimes 1_{H};X_{1},1_{H})
\end{equation*}%
i.e.%
\begin{equation*}
\gamma _{2}B(gx_{1}\otimes 1_{H};GX_{1}X_2,1_{H})=0.
\end{equation*}%
In view of the form of the element we get
\begin{equation*}
\gamma _{2}B(x_{1}x_2\otimes 1_{H};GX_2,x_1x_2)=0.
\end{equation*}%
This is (\ref{G,x1x2, GF6,x2}).

\subparagraph{Case $f=gx_{2}$}

\begin{equation*}
-B(gx_{1}\otimes 1_{H};X_{1},gx_{2})+\gamma _{2}B(gx_{1}\otimes
1_{H};GX_{1}X_{2},gx_{2})=-B(x_{1}\otimes g;X_{1},gx_{2}).
\end{equation*}%
By applying $\left( \ref{eq.10}\right) $ this rewrites as%
\begin{equation*}
-B(gx_{1}\otimes 1_{H};X_{1},gx_{2})+\gamma _{2}B(gx_{1}\otimes
1_{H};GX_{1}X_{2},gx_{2})=+B(gx_{1}\otimes 1_{H};X_{1},gx_{2})
\end{equation*}%
i.e.%
\begin{equation*}
-2B(gx_{1}\otimes 1_{H};X_{1},gx_{2})+\gamma _{2}B(gx_{1}\otimes
1_{H};GX_{1}X_{2},gx_{2})=0.
\end{equation*}%
In view of the form of the element we get%
\begin{equation*}
B(x_{1}x_{2}\otimes 1_{H};X_{2},gx_{1}x_{2})=0
\end{equation*}%
which is $\left( \ref{G,x1x2, GF7,gx1x2}\right) .$

\paragraph{\textbf{Equality (\protect\ref{GF7})}}

rewrites as%
\begin{equation*}
-B(gx_{1}\otimes 1_{H};X_{2},f)-\gamma _{1}B(gx_{1}\otimes
1_{H};GX_{1}X_{2},f)=B(gx_{1}g\otimes g;X_{2},f)
\end{equation*}%
i.e.%
\begin{equation*}
-B(gx_{1}\otimes 1_{H};X_{2},f)-\gamma _{1}B(gx_{1}\otimes
1_{H};GX_{1}X_{2},f)=-B(x_{1}\otimes g;X_{2},f).
\end{equation*}

\subparagraph{Case $f=1_{H}$}

\begin{equation*}
-B(gx_{1}\otimes 1_{H};X_{2},1_{H})-\gamma _{1}B(gx_{1}\otimes
1_{H};GX_1X_{2},1_{H})=-B(x_{1}\otimes g;X_{2},1_{H}).
\end{equation*}%
By applying $\left( \ref{eq.10}\right) $ this rewrites as%
\begin{equation*}
-B(gx_{1}\otimes 1_{H};X_{2},1_{H})-\gamma _{1}B(gx_{1}\otimes
1_{H};GX_1X_{2},1_{H})=-B(gx_{1}\otimes 1_{H};X_{2},1_{H})
\end{equation*}%
i.e.%
\begin{equation*}
\gamma _{1}B(gx_{1}\otimes 1_{H};GX_1X_{2},1_{H})=0.
\end{equation*}%
In view of the form of the element we get

\begin{equation*}
\gamma _{1}B(x_{1}x_2\otimes 1_{H};GX_{2},x_1x_2)=0.
\end{equation*}%
This is (\ref{G,x1x2, GF7,x2}).

\subparagraph{Case $f=gx_{1}$}

\begin{equation*}
-B(gx_{1}\otimes 1_{H};X_{2},gx_{1})-\gamma _{1}B(gx_{1}\otimes
1_{H};GX_{1}X_{2},gx_{1})=-B(x_{1}\otimes g;X_{2},gx_{1}).
\end{equation*}%
By applying $\left( \ref{eq.10}\right) $ this rewrites as%
\begin{equation*}
-B(gx_{1}\otimes 1_{H};X_{2},gx_{1})-\gamma _{1}B(gx_{1}\otimes
1_{H};GX_{1}X_{2},gx_{1})=B(gx_{1}\otimes 1_{H};X_{2},gx_{1})
\end{equation*}%
i.e.%
\begin{equation*}
-2B(gx_{1}\otimes 1_{H};X_{2},gx_{1})-\gamma _{1}B(gx_{1}\otimes
1_{H};GX_{1}X_{2},gx_{1})=0.
\end{equation*}%
In view of the form of the element we get%
\begin{equation*}
-2B(x_{1}x_{2}\otimes 1_{H};X_{2},gx_{1}x_{2})=0.
\end{equation*}%
This follows from $\left( \ref{G,x1x2, GF7,gx1x2}\right) .$

\paragraph{\textbf{Equality (\protect\ref{GF8})}}

rewrites as%
\begin{equation*}
B(gx_{1}\otimes 1_{H};X_{1}X_{2},f)=B(gx_{1}g\otimes g;X_{1}X_{2},f)
\end{equation*}%
i.e.%
\begin{equation*}
B(gx_{1}\otimes 1_{H};X_{1}X_{2},f)=-B(x_{1}\otimes g;X_{1}X_{2},f).
\end{equation*}

\subparagraph{Case $f=g$}

\begin{equation*}
B(gx_{1}\otimes 1_{H};X_{1}X_{2},g)=-B(x_{1}\otimes g;X_{1}X_{2},g).
\end{equation*}%
By applying $\left( \ref{eq.10}\right) $ this rewrites as%
\begin{equation*}
2B(gx_{1}\otimes 1_{H};X_{1}X_{2},g)=0.
\end{equation*}%
In view of the form of the element we get%
\begin{equation*}
B(x_{1}x_{2}\otimes 1_{H};X_{2},gx_{1}x_{2})=0
\end{equation*}%
which is $\left( \ref{G,x1x2, GF7,gx1x2}\right) .$

\subsubsection{Case $gx_{2}\otimes 1_{H}$}

\paragraph{\textbf{Equality }$\left( \protect\ref{GF1}\right) $}

rewrites as

\begin{equation*}
\begin{array}{c}
\alpha B(gx_{2}\otimes 1_{H};G,f)+ \\
+\gamma _{1}B(gx_{2}\otimes 1_{H};X_{1},f)+\gamma _{2}B(gx_{2}\otimes
1_{H};X_{2},f)%
\end{array}%
=\alpha B(gx_{2}g\otimes g;G,f)
\end{equation*}%
i.e.%
\begin{equation*}
\begin{array}{c}
\alpha B(gx_{2}\otimes 1_{H};G,f)+ \\
+\gamma _{1}B(gx_{2}\otimes 1_{H};X_{1},f)+\gamma _{2}B(gx_{2}\otimes
1_{H};X_{2},f)%
\end{array}%
=-\alpha B(x_{2}\otimes g;G,f)
\end{equation*}

\subparagraph{Case $f=1_{H}$}

\begin{equation*}
\begin{array}{c}
\alpha B(gx_{2}\otimes 1_{H};G,1_{H})+ \\
+\gamma _{1}B(gx_{2}\otimes 1_{H};X_{1},1_{H})+\gamma _{2}B(gx_{2}\otimes
1_{H};X_{2},1_{H})%
\end{array}%
=-\alpha B(x_{2}\otimes g;G,1_{H})
\end{equation*}

By applying $\left( \ref{eq.10}\right) $ this rewrites as%
\begin{equation*}
\begin{array}{c}
\alpha B(gx_{2}\otimes 1_{H};G,1_{H})+ \\
+\gamma _{1}B(gx_{2}\otimes 1_{H};X_{1},1_{H})+\gamma _{2}B(gx_{2}\otimes
1_{H};X_{2},1_{H})%
\end{array}%
=-\alpha B(gx_{2}\otimes 1_{H};G,1_{H})
\end{equation*}%
In view of the form of the element we get%
\begin{equation}
\begin{array}{c}
2\alpha \left[ -B(x_{1}x_{2}\otimes 1_{H};G,x_{1})+B(x_{1}x_{2}\otimes
1_{H};GX_{1},1_{H})\right] + \\
-\gamma _{1}B(x_{1}x_{2}\otimes 1_{H};X_{1},x_{1}) \\
+\gamma _{2}\left[ -1+B(x_{1}x_{2}\otimes
1_{H};1_{A},x_{1}x_{2})-B(x_{1}x_{2}\otimes 1_{H};X_{1},x_{2})\right]%
\end{array}%
=0  \label{G,gx2, GF1,1H}
\end{equation}

\subparagraph{Case $f=x_{1}x_{2}$}

\begin{equation*}
\begin{array}{c}
\alpha B(gx_{2}\otimes 1_{H};G,x_{1}x_{2})+ \\
+\gamma _{1}B(gx_{2}\otimes 1_{H};X_{1},x_{1}x_{2})+\gamma
_{2}B(gx_{2}\otimes 1_{H};X_{2},x_{1}x_{2})%
\end{array}%
=-\alpha B(x_{2}\otimes g;G,x_{1}x_{2})
\end{equation*}%
By applying $\left( \ref{eq.10}\right) $ this rewrites as

\begin{equation*}
\begin{array}{c}
\alpha B(gx_{2}\otimes 1_{H};G,x_{1}x_{2})+ \\
+\gamma _{1}B(gx_{2}\otimes 1_{H};X_{1},x_{1}x_{2})+\gamma
_{2}B(gx_{2}\otimes 1_{H};X_{2},x_{1}x_{2})%
\end{array}%
=-\alpha B(gx_{2}\otimes 1_{H};G,x_{1}x_{2})
\end{equation*}%
i.e.%
\begin{equation*}
\begin{array}{c}
2\alpha B(gx_{2}\otimes 1_{H};G,x_{1}x_{2})+ \\
+\gamma _{1}B(gx_{2}\otimes 1_{H};X_{1},x_{1}x_{2})+\gamma
_{2}B(gx_{2}\otimes 1_{H};X_{2},x_{1}x_{2})%
\end{array}%
=0.
\end{equation*}%
In view of the form of the element we get

\begin{equation*}
2\alpha B(x_{1}x_{2}\otimes 1_{H};GX_{1},x_{1}x_{2})=0
\end{equation*}%
which is $\left( \ref{G,x1x2, GF3,x1x2}\right) .$

\subparagraph{Case $f=gx_{1}$}

\begin{equation*}
\begin{array}{c}
\alpha B(gx_{2}\otimes 1_{H};G,gx_{1})+ \\
+\gamma _{1}B(gx_{2}\otimes 1_{H};X_{1},gx_{1})+\gamma _{2}B(gx_{2}\otimes
1_{H};X_{2},gx_{1})%
\end{array}%
=-\alpha B(x_{2}\otimes g;G,gx_{1})
\end{equation*}%
By applying $\left( \ref{eq.10}\right) $ this rewrites as%
\begin{equation*}
\begin{array}{c}
\alpha B(gx_{2}\otimes 1_{H};G,gx_{1})+ \\
+\gamma _{1}B(gx_{2}\otimes 1_{H};X_{1},gx_{1})+\gamma _{2}B(gx_{2}\otimes
1_{H};X_{2},gx_{1})%
\end{array}%
=\alpha B(gx_{2}\otimes 1_{H};G,gx_{1})
\end{equation*}%
i.e.%
\begin{equation*}
\gamma _{1}B(gx_{2}\otimes 1_{H};X_{1},gx_{1})+\gamma _{2}B(gx_{2}\otimes
1_{H};X_{2},gx_{1})=0.
\end{equation*}%
In view of the form of the element we get%
\begin{equation*}
\gamma _{2}B(x_{1}x_{2}\otimes 1_{H};X_{1},gx_{1}x_{2})=0.
\end{equation*}%
This is $\left( \ref{G,x1x2, GF3,gx1}\right) .$

\subparagraph{Case $f=gx_{2}$}

\begin{equation*}
\begin{array}{c}
\alpha B(gx_{2}\otimes 1_{H};G,gx_{2})+ \\
+\gamma _{1}B(gx_{2}\otimes 1_{H};X_{1},gx_{2})+\gamma _{2}B(gx_{2}\otimes
1_{H};X_{2},gx_{2})%
\end{array}%
=-\alpha B(x_{2}\otimes g;G,gx_{2})
\end{equation*}%
By applying $\left( \ref{eq.10}\right) $ this rewrites as%
\begin{equation*}
\begin{array}{c}
\alpha B(gx_{2}\otimes 1_{H};G,gx_{2})+ \\
+\gamma _{1}B(gx_{2}\otimes 1_{H};X_{1},gx_{2})+\gamma _{2}B(gx_{2}\otimes
1_{H};X_{2},gx_{2})%
\end{array}%
=\alpha B(gx_{2}\otimes 1_{H};G,gx_{2})
\end{equation*}%
i.e.%
\begin{equation*}
\gamma _{1}B(gx_{2}\otimes 1_{H};X_{1},gx_{2})+\gamma _{2}B(gx_{2}\otimes
1_{H};X_{2},gx_{2})=0.
\end{equation*}%
In view of the form of the element we get%
\begin{equation*}
\gamma _{1}B(x_{1}x_{2}\otimes 1_{H};X_{1},gx_{1}x_{2})=0
\end{equation*}%
which follows from $\left( \ref{G,x1x2, GF4,gx1}\right) .$

\paragraph{\textbf{Equality }$\left( \mathbf{\protect\ref{GF2}}\right) $}

rewrites as%
\begin{equation*}
\begin{array}{c}
B(gx_{2}\otimes 1_{H};1_{A},f) \\
+\gamma _{1}B(gx_{2}\otimes 1_{H};GX_{1},f)+\gamma _{2}B(gx_{2}\otimes
1_{H};GX_{2},f)%
\end{array}%
=B(gx_{2}g\otimes g;1_{A},f)
\end{equation*}%
i.e.

\begin{equation*}
\begin{array}{c}
B(gx_{2}\otimes 1_{H};1_{A},f) \\
+\gamma _{1}B(gx_{2}\otimes 1_{H};GX_{1},f)+\gamma _{2}B(gx_{2}\otimes
1_{H};GX_{2},f)%
\end{array}%
=-B(x_{2}\otimes g;1_{A},f).
\end{equation*}

\subparagraph{Case $f=g$}

\begin{equation*}
\begin{array}{c}
B(gx_{2}\otimes 1_{H};1_{A},g) \\
+\gamma _{1}B(gx_{2}\otimes 1_{H};GX_{1},g)+\gamma _{2}B(gx_{2}\otimes
1_{H};GX_{2},g)%
\end{array}%
=-B(x_{2}\otimes g;1_{A},g).
\end{equation*}%
By applying $\left( \ref{eq.10}\right) $ this rewrites as%
\begin{equation*}
\begin{array}{c}
B(gx_{2}\otimes 1_{H};1_{A},g) \\
+\gamma _{1}B(gx_{2}\otimes 1_{H};GX_{1},g)+\gamma _{2}B(gx_{2}\otimes
1_{H};GX_{2},g)%
\end{array}%
=-B(gx_{2}\otimes 1_{H};1_{A},g)
\end{equation*}%
i.e.%
\begin{equation*}
\begin{array}{c}
2B(gx_{2}\otimes 1_{H};1_{A},g) \\
+\gamma _{1}B(gx_{2}\otimes 1_{H};GX_{1},g)+\gamma _{2}B(gx_{2}\otimes
1_{H};GX_{2},g)%
\end{array}%
=0.
\end{equation*}%
In view of the form of the element we get%
\begin{equation}
\begin{array}{c}
2\left[ -B(x_{1}x_{2}\otimes 1_{H};1_{A},gx_{1})-B(x_{1}x_{2}\otimes
1_{H};X_{1},g)\right] \\
-\gamma _{1}B(x_{1}x_{2}\otimes 1_{H};GX_{1},gx_{1}) \\
+\gamma _{2}\left[ -B(x_{1}x_{2}\otimes
1_{H};G,gx_{1}x_{2})-B(x_{1}x_{2}\otimes 1_{H};GX_{1},gx_{2})\right]%
\end{array}%
=0.  \label{G,gx2, GF2,g}
\end{equation}

\subparagraph{Case $f=x_{1}$}

\begin{equation*}
\begin{array}{c}
B(gx_{2}\otimes 1_{H};1_{A},x_{1}) \\
+\gamma _{1}B(gx_{2}\otimes 1_{H};GX_{1},x_{1})+\gamma _{2}B(gx_{2}\otimes
1_{H};GX_{2},x_{1})%
\end{array}%
=-B(x_{2}\otimes g;1_{A},x_{1}).
\end{equation*}%
By applying $\left( \ref{eq.10}\right) $ this rewrites as%
\begin{equation*}
\begin{array}{c}
B(gx_{2}\otimes 1_{H};1_{A},x_{1}) \\
+\gamma _{1}B(gx_{2}\otimes 1_{H};GX_{1},x_{1})+\gamma _{2}B(gx_{2}\otimes
1_{H};GX_{2},x_{1})%
\end{array}%
=B(gx_{2}\otimes 1_{H};1_{A},x_{1})
\end{equation*}%
i.e.%
\begin{equation*}
\gamma _{1}B(gx_{2}\otimes 1_{H};GX_{1},x_{1})+\gamma _{2}B(gx_{2}\otimes
1_{H};GX_{2},x_{1})=0
\end{equation*}%
In view of the form of the element we get%
\begin{equation*}
\gamma _{2}B(x_{1}x_{2}\otimes 1_{H};GX_{1},x_{1}x_{2})=0
\end{equation*}%
which is $\left( \ref{G,x1x2, GF6,x1}\right) $

\subparagraph{Case $f=x_{2}$}

\begin{equation*}
\begin{array}{c}
B(gx_{2}\otimes 1_{H};1_{A},x_{2}) \\
+\gamma _{1}B(gx_{2}\otimes 1_{H};GX_{1},x_{2})+\gamma _{2}B(gx_{2}\otimes
1_{H};GX_{2},x_{2})%
\end{array}%
=-B(x_{2}\otimes g;1_{A},x_{2}).
\end{equation*}%
By applying $\left( \ref{eq.10}\right) $ this rewrites as%
\begin{equation*}
\begin{array}{c}
B(gx_{2}\otimes 1_{H};1_{A},x_{2}) \\
+\gamma _{1}B(gx_{2}\otimes 1_{H};GX_{1},x_{2})+\gamma _{2}B(gx_{2}\otimes
1_{H};GX_{2},x_{2})%
\end{array}%
=B(gx_{2}\otimes 1_{H};1_{A},x_{2})
\end{equation*}%
i.e.%
\begin{equation*}
\gamma _{1}B(gx_{2}\otimes 1_{H};GX_{1},x_{2})+\gamma _{2}B(gx_{2}\otimes
1_{H};GX_{2},x_{2})=0
\end{equation*}%
In view of the form of the element we get%
\begin{equation}
\gamma _{1}B(x_{1}x_{2}\otimes 1_{H};GX_{1},x_{1}x_{2})=0
\label{G,gx2, GF2,x2}
\end{equation}

\subparagraph{Case $f=gx_{1}x_{2}$}

\begin{equation*}
\begin{array}{c}
B(gx_{2}\otimes 1_{H};1_{A},gx_{1}x_{2}) \\
+\gamma _{1}B(gx_{2}\otimes 1_{H};GX_{1},gx_{1}x_{2})+\gamma
_{2}B(gx_{2}\otimes 1_{H};GX_{2},gx_{1}x_{2})%
\end{array}%
=-B(x_{2}\otimes g;1_{A},gx_{1}x_{2}).
\end{equation*}%
By applying $\left( \ref{eq.10}\right) $ this rewrites as%
\begin{equation*}
\begin{array}{c}
B(gx_{2}\otimes 1_{H};1_{A},gx_{1}x_{2}) \\
+\gamma _{1}B(gx_{2}\otimes 1_{H};GX_{1},gx_{1}x_{2})+\gamma
_{2}B(gx_{2}\otimes 1_{H};GX_{2},gx_{1}x_{2})%
\end{array}%
=-B(gx_{2}\otimes 1_{H};1_{A},gx_{1}x_{2}).
\end{equation*}%
i.e.%
\begin{equation*}
\begin{array}{c}
2B(gx_{2}\otimes 1_{H};1_{A},gx_{1}x_{2}) \\
+\gamma _{1}B(gx_{2}\otimes 1_{H};GX_{1},gx_{1}x_{2})+\gamma
_{2}B(gx_{2}\otimes 1_{H};GX_{2},gx_{1}x_{2})%
\end{array}%
=0
\end{equation*}%
In view of the form of the element we get%
\begin{equation*}
B(x_{1}x_{2}\otimes 1_{H};X_{1},gx_{1}x_{2})=0
\end{equation*}%
which is $\left( \ref{G,x1x2, GF6,gx1x2}\right) .$

\paragraph{\textbf{Equality }$\left( \protect\ref{GF3}\right) $}

rewrites as%
\begin{equation*}
\gamma _{2}B(gx_{2}\otimes 1_{H};X_{1}X_{2},f)-\alpha B(gx_{2}\otimes
1_{H};GX_{1},f)=\alpha B(gx_{2}g\otimes g;GX_{1},f)
\end{equation*}%
i.e.%
\begin{equation*}
\gamma _{2}B(gx_{2}\otimes 1_{H};X_{1}X_{2},f)-\alpha B(gx_{2}\otimes
1_{H};GX_{1},f)=-\alpha B(x_{2}\otimes g;GX_{1},f).
\end{equation*}

\subparagraph{Case $f=g$}

\begin{equation*}
\gamma _{2}B(gx_{2}\otimes 1_{H};X_{1}X_{2},g)-\alpha B(gx_{2}\otimes
1_{H};GX_{1},g)=-\alpha B(x_{2}\otimes g;GX_{1},g).
\end{equation*}%
By applying $\left( \ref{eq.10}\right) $ this rewrites as%
\begin{equation*}
\gamma _{2}B(gx_{2}\otimes 1_{H};X_{1}X_{2},g)-\alpha B(gx_{2}\otimes
1_{H};GX_{1},g)=-\alpha B(gx_{2}\otimes 1_{H};GX_{1},g)
\end{equation*}%
i.e.%
\begin{equation*}
\gamma _{2}B(gx_{2}\otimes 1_{H};X_{1}X_{2},g)=0.
\end{equation*}%
In view of the form of the element we get%
\begin{equation*}
\gamma _{2}B(x_{1}x_{2}\otimes 1_{H};X_{1},gx_{1}x_{2})=0
\end{equation*}%
which is $\left( \ref{G,x1x2, GF3,gx1}\right) $

\subparagraph{Case $f=x_{2}$}

\begin{equation*}
\gamma _{2}B(gx_{2}\otimes 1_{H};X_{1}X_{2},x_{2})-\alpha B(gx_{2}\otimes
1_{H};GX_{1},x_{2})=-\alpha B(x_{2}\otimes g;GX_{1},x_{2}).
\end{equation*}%
By applying $\left( \ref{eq.10}\right) $ this rewrites as%
\begin{equation*}
\gamma _{2}B(gx_{2}\otimes 1_{H};X_{1}X_{2},x_{2})-\alpha B(gx_{2}\otimes
1_{H};GX_{1},x_{2})=\alpha B(gx_{2}\otimes 1_{H};GX_{1},x_{2})
\end{equation*}%
i.e.%
\begin{equation*}
\gamma _{2}B(gx_{2}\otimes 1_{H};X_{1}X_{2},x_{2})-2\alpha B(gx_{2}\otimes
1_{H};GX_{1},x_{2})=0.
\end{equation*}%
In view of the form of the element we get%
\begin{equation*}
\alpha B(x_{1}x_{2}\otimes 1_{H};GX_{1},x_{1}x_{2})=0
\end{equation*}%
which is $\left( \ref{G,x1x2, GF3,x1x2}\right) .$

\paragraph{\textbf{Equality }$\left( \protect\ref{GF4}\right) $}

rewrites as%
\begin{equation*}
-\gamma _{1}B(gx_{2}\otimes 1_{H};X_{1}X_{2},f)-\alpha B(gx_{2}\otimes
1_{H};GX_{2},f)=\alpha B(gx_{2}g\otimes g;GX_{2},f)
\end{equation*}%
i.e.%
\begin{equation*}
-\gamma _{1}B(gx_{2}\otimes 1_{H};X_{1}X_{2},f)-\alpha B(gx_{2}\otimes
1_{H};GX_{2},f)=-\alpha B(x_{2}\otimes g;GX_{2},f).
\end{equation*}

\subparagraph{Case $f=g$}

\begin{equation*}
-\gamma _{1}B(gx_{2}\otimes 1_{H};X_{1}X_{2},g)-\alpha B(gx_{2}\otimes
1_{H};GX_{2},g)=-\alpha B(x_{2}\otimes g;GX_{2},g).
\end{equation*}%
By applying $\left( \ref{eq.10}\right) $ this rewrites as%
\begin{equation*}
-\gamma _{1}B(gx_{2}\otimes 1_{H};X_{1}X_{2},g)-\alpha B(gx_{2}\otimes
1_{H};GX_{2},g)=-\alpha B(gx_{2}\otimes 1_{H};GX_{2},g)
\end{equation*}%
i.e.%
\begin{equation*}
\gamma _{1}B(gx_{2}\otimes 1_{H};X_{1}X_{2},g)=0.
\end{equation*}%
By the form of the element we get%
\begin{equation*}
\gamma _{1}B(x_{1}x_{2}\otimes 1_{H};X_{1},gx_{1}x_{2})=0
\end{equation*}%
which follows from $\left( \ref{G,x1x2, GF6,gx1x2}\right) $.

\subparagraph{Case $f=x_{1}$}

\begin{equation*}
-\gamma _{1}B(gx_{2}\otimes 1_{H};X_{1}X_{2},x_{1})-\alpha B(gx_{2}\otimes
1_{H};GX_{2},x_{1})=-\alpha B(x_{2}\otimes g;GX_{2},x_{1}).
\end{equation*}%
By applying $\left( \ref{eq.10}\right) $ this rewrites as%
\begin{equation*}
-\gamma _{1}B(gx_{2}\otimes 1_{H};X_{1}X_{2},x_{1})-\alpha B(gx_{2}\otimes
1_{H};GX_{2},x_{1})=\alpha B(gx_{2}\otimes 1_{H};GX_{2},x_{1})
\end{equation*}%
i.e.%
\begin{equation*}
-\gamma _{1}B(gx_{2}\otimes 1_{H};X_{1}X_{2},x_{1})-2\alpha B(gx_{2}\otimes
1_{H};GX_{2},x_{1})=0.
\end{equation*}%
In view of the form of the element we get%
\begin{equation*}
\alpha B(x_{1}x_{2}\otimes 1_{H};GX_{1},x_{1}x_{2})=0
\end{equation*}%
which is $\left( \ref{G,x1x2, GF3,x1x2}\right) .$

\paragraph{\textbf{Equality }$\left( \protect\ref{GF5}\right) $}

rewrites as%
\begin{equation*}
\alpha B(gx_{2}\otimes 1_{H};GX_{1}X_{2},f)=\alpha B(gx_{2}g\otimes
g;GX_{1}X_{2},f)
\end{equation*}%
i.e.%
\begin{equation*}
\alpha B(gx_{2}\otimes 1_{H};GX_{1}X_{2},f)=-\alpha B(x_{2}\otimes
g;GX_{1}X_{2},f).
\end{equation*}

\subparagraph{Case $f=1_{H}$}

\begin{equation*}
\alpha B(gx_{2}\otimes 1_{H};GX_{1}X_{2},1_{H})=-\alpha B(x_{2}\otimes
g;GX_{1}X_{2},1_{H}).
\end{equation*}%
By applying $\left( \ref{eq.10}\right) $ this rewrites as%
\begin{equation*}
\alpha B(gx_{2}\otimes 1_{H};GX_{1}X_{2},1_{H})=-\alpha B(gx_{2}\otimes
1_{H};GX_{1}X_{2},1_{H})
\end{equation*}%
i.e.%
\begin{equation*}
2\alpha B(gx_{2}\otimes 1_{H};GX_{1}X_{2},1_{H})=0
\end{equation*}%
In view of the form of the element we get%
\begin{equation*}
\alpha B(x_{1}x_{2}\otimes 1_{H};GX_{1},x_{1}x_{2})=0
\end{equation*}%
which is $\left( \ref{G,x1x2, GF3,x1x2}\right) .$

\paragraph{\textbf{Equality }$\left( \protect\ref{GF6}\right) $}

rewrites as%
\begin{equation*}
-B(gx_{2}\otimes 1_{H};X_{1},f)+\gamma _{2}B(gx_{2}\otimes
1_{H};GX_{1}X_{2},f)=B(gx_{2}g\otimes g;X_{1},f)
\end{equation*}%
i.e.%
\begin{equation*}
-B(gx_{2}\otimes 1_{H};X_{1},f)+\gamma _{2}B(gx_{2}\otimes
1_{H};GX_{1}X_{2},f)=-B(x_{2}\otimes g;X_{1},f).
\end{equation*}

\subparagraph{Case $f=1_{H}$}

\begin{equation*}
-B(gx_{2}\otimes 1_{H};X_{1},1_{H})+\gamma _{2}B(gx_{2}\otimes
1_{H};GX_{1}X_{2},1_{H})=-B(x_{2}\otimes g;X_{1},1_{H}).
\end{equation*}%
By applying $\left( \ref{eq.10}\right) $ this rewrites as%
\begin{equation*}
-B(gx_{2}\otimes 1_{H};X_{1},1_{H})+\gamma _{2}B(gx_{2}\otimes
1_{H};GX_{1}X_{2},1_{H})=-B(gx_{2}\otimes 1_{H};X_{1},1_{H})
\end{equation*}%
i.e.%
\begin{equation*}
\gamma _{2}B(gx_{2}\otimes 1_{H};GX_{1}X_{2},1_{H})=0.
\end{equation*}%
In view of the form of the element we get%
\begin{equation*}
\gamma _{2}B(x_{1}x_{2}\otimes 1_{H};GX_{1},x_{1}x_{2})=0
\end{equation*}%
which is $\left( \ref{G,x1x2, GF6,x1}\right) .$

\subparagraph{Case $f=gx_{2}$}

\begin{equation*}
-B(gx_{2}\otimes 1_{H};X_{1},gx_{2})+\gamma _{2}B(gx_{2}\otimes
1_{H};GX_{1}X_{2},gx_{2})=-B(x_{2}\otimes g;X_{1},gx_{2}).
\end{equation*}%
By applying $\left( \ref{eq.10}\right) $ this rewrites as%
\begin{equation*}
-B(gx_{2}\otimes 1_{H};X_{1},gx_{2})+\gamma _{2}B(gx_{2}\otimes
1_{H};GX_{1}X_{2},gx_{2})=B(gx_{2}\otimes 1_{H};X_{1},gx_{2})
\end{equation*}%
i.e.%
\begin{equation*}
-2B(gx_{2}\otimes 1_{H};X_{1},gx_{2})+\gamma _{2}B(gx_{2}\otimes
1_{H};GX_{1}X_{2},gx_{2})=0.
\end{equation*}%
In view of the form of the element we get%
\begin{equation*}
B(x_{1}x_{2}\otimes 1_{H};X_{1},gx_{1}x_{2})=0
\end{equation*}%
which is $\left( \ref{G,x1x2, GF6,gx1x2}\right) .$

\paragraph{\textbf{Equality }$\left( \protect\ref{GF7}\right) $}

rewrites as%
\begin{equation*}
-B(gx_{2}\otimes 1_{H};X_{2},f)-\gamma _{1}B(gx_{2}\otimes
1_{H};GX_{1}X_{2},f)=B(gx_{2}g\otimes g;X_{2},f)
\end{equation*}%
i.e.%
\begin{equation*}
-B(gx_{2}\otimes 1_{H};X_{2},f)-\gamma _{1}B(gx_{2}\otimes
1_{H};GX_{1}X_{2},f)=-B(x_{2}\otimes g;X_{2},f).
\end{equation*}

\subparagraph{Case $f=1_{H}$}

\begin{equation*}
-B(gx_{2}\otimes 1_{H};X_{2},1_{H})-\gamma _{1}B(gx_{2}\otimes
1_{H};GX_{1}X_{2},1_{H})=-B(x_{2}\otimes g;X_{2},1_{H}).
\end{equation*}%
By applying $\left( \ref{eq.10}\right) $ this rewrites as%
\begin{equation*}
-B(gx_{2}\otimes 1_{H};X_{2},1_{H})-\gamma _{1}B(gx_{2}\otimes
1_{H};GX_{1}X_{2},1_{H})=-B(gx_{2}\otimes 1_{H};X_{2},1_{H})
\end{equation*}%
i.e.%
\begin{equation*}
B(gx_{2}\otimes 1_{H};GX_{1}X_{2},1_{H})=0.
\end{equation*}%
In view of the form of the element we get%
\begin{equation*}
\gamma _{1}B(x_{1}x_{2}\otimes 1_{H};GX_{1},x_{1}x_{2})=0
\end{equation*}%
which is $\left( \ref{G,x1x2, GF7,x1}\right) .$

\subparagraph{Case $f=gx_{1}$}

\begin{equation*}
-B(gx_{2}\otimes 1_{H};X_{2},gx_{1})-\gamma _{1}B(gx_{2}\otimes
1_{H};GX_{1}X_{2},gx_{1})=-B(x_{2}\otimes g;X_{2},gx_{1}).
\end{equation*}%
By applying $\left( \ref{eq.10}\right) $ this rewrites as%
\begin{equation*}
-B(gx_{2}\otimes 1_{H};X_{2},gx_{1})-\gamma _{1}B(gx_{2}\otimes
1_{H};GX_{1}X_{2},gx_{1})=B(gx_{2}\otimes 1_{H};X_{2},gx_{1})
\end{equation*}%
i.e.%
\begin{equation*}
-2B(gx_{2}\otimes 1_{H};X_{2},gx_{1})-\gamma _{1}B(gx_{2}\otimes
1_{H};GX_{1}X_{2},gx_{1})=0.
\end{equation*}%
In view of the form of the element we get%
\begin{equation*}
B(x_{1}x_{2}\otimes 1_{H};X_{1},gx_{1}x_{2})=0
\end{equation*}%
which is $\left( \ref{G,x1x2, GF6,gx1x2}\right) .$

\paragraph{\textbf{Equality }$\left( \protect\ref{GF8}\right) $}

rewrites as%
\begin{equation*}
B(gx_{2}\otimes 1_{H};X_{1}X_{2},f)=B(gx_{2}g\otimes g;X_{1}X_{2},f)
\end{equation*}%
i.e.%
\begin{equation*}
B(gx_{2}\otimes 1_{H};X_{1}X_{2},f)=-B(x_{2}\otimes g;X_{1}X_{2},f).
\end{equation*}

\subparagraph{Case $f=g$}

\begin{equation*}
B(gx_{2}\otimes 1_{H};X_{1}X_{2},g)=-B(x_{2}\otimes g;X_{1}X_{2},g).
\end{equation*}%
By applying $\left( \ref{eq.10}\right) $ this rewrites as%
\begin{equation*}
B(gx_{2}\otimes 1_{H};X_{1}X_{2},g)=-B(gx_{2}\otimes 1_{H};X_{1}X_{2},g)
\end{equation*}%
i.e.%
\begin{equation*}
2B(gx_{2}\otimes 1_{H};X_{1}X_{2},g)=0.
\end{equation*}

In view of the form of the element we get%
\begin{equation*}
B(x_{1}x_{2}\otimes 1_{H};X_{1},gx_{1}x_{2})=0
\end{equation*}%
which is $\left( \ref{G,x1x2, GF6,gx1x2}\right) .$

\subsubsection{Case $gx_{1}x_{2}\otimes 1_{H}$}

\paragraph{\textbf{Equality }$\left( \protect\ref{GF1}\right) $}

rewrites as%
\begin{equation*}
\begin{array}{c}
\alpha B(gx_{1}x_{2}\otimes 1_{H};G,f)+ \\
+\gamma _{1}B(gx_{1}x_{2}\otimes 1_{H};X_{1},f)+\gamma
_{2}B(gx_{1}x_{2}\otimes 1_{H};X_{2},f)%
\end{array}%
=\alpha B(gx_{1}x_{2}g\otimes g;G,f)
\end{equation*}%
i.e.%
\begin{equation*}
\begin{array}{c}
\alpha B(gx_{1}x_{2}\otimes 1_{H};G,f)+ \\
+\gamma _{1}B(gx_{1}x_{2}\otimes 1_{H};X_{1},f)+\gamma
_{2}B(gx_{1}x_{2}\otimes 1_{H};X_{2},f)%
\end{array}%
=\alpha B(x_{1}x_{2}\otimes g;G,f).
\end{equation*}

\subparagraph{Case $f=1_{H}$}

\begin{equation*}
\begin{array}{c}
\alpha B(gx_{1}x_{2}\otimes 1_{H};G,1_{H})+ \\
+\gamma _{1}B(gx_{1}x_{2}\otimes 1_{H};X_{1},1_{H})+\gamma
_{2}B(gx_{1}x_{2}\otimes 1_{H};X_{2},1_{H})%
\end{array}%
=\alpha B(x_{1}x_{2}\otimes g;G,1_{H}).
\end{equation*}%
By applying $\left( \ref{eq.10}\right) $ this rewrites as%
\begin{equation*}
\begin{array}{c}
\alpha B(gx_{1}x_{2}\otimes 1_{H};G,1_{H})+ \\
+\gamma _{1}B(gx_{1}x_{2}\otimes 1_{H};X_{1},1_{H})+\gamma
_{2}B(gx_{1}x_{2}\otimes 1_{H};X_{2},1_{H})%
\end{array}%
=\alpha B(gx_{1}x_{2}\otimes 1_{H};G,1_{H})
\end{equation*}%
i.e.%
\begin{equation*}
\gamma _{1}B(gx_{1}x_{2}\otimes 1_{H};X_{1},1_{H})+\gamma
_{2}B(gx_{1}x_{2}\otimes 1_{H};X_{2},1_{H})=0.
\end{equation*}%
In view of the form of the element we get%
\begin{gather}
\gamma _{1}\left[ -B(x_{2}\otimes 1_{H};1_{A},1_{H})+B(gx_{1}x_{2}\otimes
1_{H};1_{A},x_{1})\right]  \label{G,gx1x2, GF1,1H} \\
+\gamma _{2}\left[ B(x_{1}\otimes 1_{H};1_{A},1_{H})+B(gx_{1}x_{2}\otimes
1_{H};1_{A},x_{2})\right] =0.  \notag
\end{gather}

\subparagraph{Case $f=x_{1}x_{2}$}

\begin{equation*}
\begin{array}{c}
\alpha B(gx_{1}x_{2}\otimes 1_{H};G,x_{1}x_{2})+ \\
+\gamma _{1}B(gx_{1}x_{2}\otimes 1_{H};X_{1},x_{1}x_{2})+\gamma
_{2}B(gx_{1}x_{2}\otimes 1_{H};X_{2},x_{1}x_{2})%
\end{array}%
=\alpha B(x_{1}x_{2}\otimes g;G,x_{1}x_{2}).
\end{equation*}%
By applying $\left( \ref{eq.10}\right) $ this rewrites as%
\begin{equation*}
\begin{array}{c}
\alpha B(gx_{1}x_{2}\otimes 1_{H};G,x_{1}x_{2})+ \\
+\gamma _{1}B(gx_{1}x_{2}\otimes 1_{H};X_{1},x_{1}x_{2})+\gamma
_{2}B(gx_{1}x_{2}\otimes 1_{H};X_{2},x_{1}x_{2})%
\end{array}%
=\alpha B(gx_{1}x_{2}\otimes 1_{H};G,x_{1}x_{2})
\end{equation*}%
i.e.%
\begin{equation*}
\gamma _{1}B(gx_{1}x_{2}\otimes 1_{H};X_{1},x_{1}x_{2})+\gamma
_{2}B(gx_{1}x_{2}\otimes 1_{H};X_{2},x_{1}x_{2})=0.
\end{equation*}%
In view of the form of the element we get%
\begin{equation}
-\gamma _{1}B(x_{2}\otimes 1_{H};1_{A},x_{1}x_{2})+\gamma _{2}B(x_{1}\otimes
1_{H};1_{A},x_{1}x_{2})=0.  \label{G,gx1x2, GF1,x1x2}
\end{equation}

\subparagraph{Case $f=gx_{1}$}

\begin{equation*}
\begin{array}{c}
\alpha B(gx_{1}x_{2}\otimes 1_{H};G,gx_{1})+ \\
+\gamma _{1}B(gx_{1}x_{2}\otimes 1_{H};X_{1},gx_{1})+\gamma
_{2}B(gx_{1}x_{2}\otimes 1_{H};X_{2},gx_{1})%
\end{array}%
=\alpha B(x_{1}x_{2}\otimes g;G,gx_{1}).
\end{equation*}%
By applying $\left( \ref{eq.10}\right) $ this rewrites as%
\begin{equation*}
\begin{array}{c}
\alpha B(gx_{1}x_{2}\otimes 1_{H};G,gx_{1})+ \\
+\gamma _{1}B(gx_{1}x_{2}\otimes 1_{H};X_{1},gx_{1})+\gamma
_{2}B(gx_{1}x_{2}\otimes 1_{H};X_{2},gx_{1})%
\end{array}%
=-\alpha B(gx_{1}x_{2}\otimes 1_{H};G,gx_{1})
\end{equation*}%
i.e.%
\begin{equation*}
\begin{array}{c}
2\alpha B(gx_{1}x_{2}\otimes 1_{H};G,gx_{1})+ \\
+\gamma _{1}B(gx_{1}x_{2}\otimes 1_{H};X_{1},gx_{1})+\gamma
_{2}B(gx_{1}x_{2}\otimes 1_{H};X_{2},gx_{1})%
\end{array}%
=0.
\end{equation*}%
In view of the form of the element we get%
\begin{equation}
\begin{array}{c}
2\alpha B(gx_{1}x_{2}\otimes 1_{H};G,gx_{1})+ \\
-\gamma _{1}B(x_{2}\otimes 1_{H};1_{A},gx_{1})+\gamma _{2}\left[
B(x_{1}\otimes 1_{H};1_{A},gx_{1})+B(gx_{1}x_{2}\otimes
1_{H};1_{A},gx_{1}x_{2})\right]%
\end{array}%
=0.  \label{G,gx1x2, GF1,gx1}
\end{equation}

\subparagraph{Case $f=gx_{1}$}

\begin{equation*}
\begin{array}{c}
\alpha B(gx_{1}x_{2}\otimes 1_{H};G,gx_{2})+ \\
+\gamma _{1}B(gx_{1}x_{2}\otimes 1_{H};X_{1},gx_{2})+\gamma
_{2}B(gx_{1}x_{2}\otimes 1_{H};X_{2},gx_{2})%
\end{array}%
=\alpha B(x_{1}x_{2}\otimes g;G,gx_{2}).
\end{equation*}%
By applying $\left( \ref{eq.10}\right) $ this rewrites as%
\begin{equation*}
\begin{array}{c}
\alpha B(gx_{1}x_{2}\otimes 1_{H};G,gx_{2})+ \\
+\gamma _{1}B(gx_{1}x_{2}\otimes 1_{H};X_{1},gx_{2})+\gamma
_{2}B(gx_{1}x_{2}\otimes 1_{H};X_{2},gx_{2})%
\end{array}%
=-\alpha B(gx_{1}x_{2}\otimes 1_{H};G,gx_{2})
\end{equation*}%
i.e.%
\begin{equation*}
\begin{array}{c}
2\alpha B(gx_{1}x_{2}\otimes 1_{H};G,gx_{2})+ \\
+\gamma _{1}B(gx_{1}x_{2}\otimes 1_{H};X_{1},gx_{2})+\gamma
_{2}B(gx_{1}x_{2}\otimes 1_{H};X_{2},gx_{2})%
\end{array}%
=0.
\end{equation*}%
In view of the form of the element we get%
\begin{equation}
\begin{array}{c}
2\alpha B(gx_{1}x_{2}\otimes 1_{H};G,gx_{2})+ \\
+\gamma _{1}\left[ -B(x_{2}\otimes 1_{H};1_{A},gx_{2})-B(gx_{1}x_{2}\otimes
1_{H};1_{A},gx_{1}x_{2})\right] +\gamma _{2}B(x_{1}\otimes
1_{H};1_{A},gx_{2})%
\end{array}%
=0.  \label{G,gx1x2, GF1,gx2}
\end{equation}

\paragraph{\textbf{Equality }$\left( \protect\ref{GF2}\right) $}

rewrites as%
\begin{equation*}
\begin{array}{c}
B(gx_{1}x_{2}\otimes 1_{H};1_{A},f) \\
+\gamma _{1}B(gx_{1}x_{2}\otimes 1_{H};GX_{1},f)+\gamma
_{2}B(gx_{1}x_{2}\otimes 1_{H};GX_{2},f)%
\end{array}%
=B(gx_{1}x_{2}g\otimes g;1_{A},f)
\end{equation*}%
i.e.%
\begin{equation*}
\begin{array}{c}
B(gx_{1}x_{2}\otimes 1_{H};1_{A},f) \\
+\gamma _{1}B(gx_{1}x_{2}\otimes 1_{H};GX_{1},f)+\gamma
_{2}B(gx_{1}x_{2}\otimes 1_{H};GX_{2},f)%
\end{array}%
=B(x_{1}x_{2}\otimes g;1_{A},f).
\end{equation*}

\subparagraph{Case $f=g$}

\begin{equation*}
\begin{array}{c}
B(gx_{1}x_{2}\otimes 1_{H};1_{A},g) \\
+\gamma _{1}B(gx_{1}x_{2}\otimes 1_{H};GX_{1},g)+\gamma
_{2}B(gx_{1}x_{2}\otimes 1_{H};GX_{2},g)%
\end{array}%
=B(x_{1}x_{2}\otimes g;1_{A},g).
\end{equation*}%
By applying $\left( \ref{eq.10}\right) $ this rewrites as%
\begin{equation*}
\begin{array}{c}
B(gx_{1}x_{2}\otimes 1_{H};1_{A},g) \\
+\gamma _{1}B(gx_{1}x_{2}\otimes 1_{H};GX_{1},g)+\gamma
_{2}B(gx_{1}x_{2}\otimes 1_{H};GX_{2},g)%
\end{array}%
=B(gx_{1}x_{2}\otimes 1_{H};1_{A},g)
\end{equation*}%
i.e.%
\begin{equation*}
\gamma _{1}B(gx_{1}x_{2}\otimes 1_{H};GX_{1},g)+\gamma
_{2}B(gx_{1}x_{2}\otimes 1_{H};GX_{2},g)=0.
\end{equation*}%
In view of the form of the element we get%
\begin{gather}
\gamma _{1}\left[ B(x_{2}\otimes 1_{H};G,g)-B(gx_{1}x_{2}\otimes
1_{H};G,gx_{1})\right]  \label{G,gx1x2, GF2,g} \\
+\gamma _{2}\left[ -B(x_{1}\otimes 1_{H};G,g)-B(gx_{1}x_{2}\otimes
1_{H};G,gx_{2})\right] =0.  \notag
\end{gather}

\subparagraph{Case $f=x_{1}$}

\begin{equation*}
\begin{array}{c}
B(gx_{1}x_{2}\otimes 1_{H};1_{A},x_{1}) \\
+\gamma _{1}B(gx_{1}x_{2}\otimes 1_{H};GX_{1},x_{1})+\gamma
_{2}B(gx_{1}x_{2}\otimes 1_{H};GX_{2},x_{1})%
\end{array}%
=B(x_{1}x_{2}\otimes g;1_{A},x_{1}).
\end{equation*}%
By applying $\left( \ref{eq.10}\right) $ this rewrites as%
\begin{equation*}
\begin{array}{c}
B(gx_{1}x_{2}\otimes 1_{H};1_{A},x_{1}) \\
+\gamma _{1}B(gx_{1}x_{2}\otimes 1_{H};GX_{1},x_{1})+\gamma
_{2}B(gx_{1}x_{2}\otimes 1_{H};GX_{2},x_{1})%
\end{array}%
=-B(gx_{1}x_{2}\otimes 1_{h};1_{A},x_{1})
\end{equation*}%
i.e.%
\begin{equation*}
\begin{array}{c}
2B(gx_{1}x_{2}\otimes 1_{H};1_{A},x_{1}) \\
+\gamma _{1}B(gx_{1}x_{2}\otimes 1_{H};GX_{1},x_{1})+\gamma
_{2}B(gx_{1}x_{2}\otimes 1_{H};GX_{2},x_{1})%
\end{array}%
=0.
\end{equation*}%
In view of the form of the element we get%
\begin{equation}
\begin{array}{c}
2B(gx_{1}x_{2}\otimes 1_{H};1_{A},x_{1}) \\
+\gamma _{1}B(x_{2}\otimes 1_{H};G,x_{1})+\gamma _{2}\left[ -B(x_{1}\otimes
1_{H};G,x_{1})-B(gx_{1}x_{2}\otimes 1_{H};G,x_{1}x_{2})\right]%
\end{array}%
=0.  \label{G,gx1x2, GF2,x1}
\end{equation}

\subparagraph{Case $f=x_{2}$}

\begin{equation*}
\begin{array}{c}
B(gx_{1}x_{2}\otimes 1_{H};1_{A},x_{2}) \\
+\gamma _{1}B(gx_{1}x_{2}\otimes 1_{H};GX_{1},x_{2})+\gamma
_{2}B(gx_{1}x_{2}\otimes 1_{H};GX_{2},x_{2})%
\end{array}%
=B(x_{1}x_{2}\otimes g;1_{A},x_{2}).
\end{equation*}%
By applying $\left( \ref{eq.10}\right) $ this rewrites as%
\begin{equation*}
\begin{array}{c}
B(gx_{1}x_{2}\otimes 1_{H};1_{A},x_{2}) \\
+\gamma _{1}B(gx_{1}x_{2}\otimes 1_{H};GX_{1},x_{2})+\gamma
_{2}B(gx_{1}x_{2}\otimes 1_{H};GX_{2},x_{2})%
\end{array}%
=-B(gx_{1}x_{2}\otimes 1_{h};1_{A},x_{2})
\end{equation*}%
i.e.%
\begin{equation*}
\begin{array}{c}
2B(gx_{1}x_{2}\otimes 1_{H};1_{A},x_{2}) \\
+\gamma _{1}B(gx_{1}x_{2}\otimes 1_{H};GX_{1},x_{2})+\gamma
_{2}B(gx_{1}x_{2}\otimes 1_{H};GX_{2},x_{2})%
\end{array}%
=0.
\end{equation*}%
In view of the form of the element we get%
\begin{equation}
\begin{array}{c}
2B(gx_{1}x_{2}\otimes 1_{H};1_{A},x_{2}) \\
+\gamma _{1}\left[ B(x_{2}\otimes 1_{H};G,x_{2})+B(gx_{1}x_{2}\otimes
1_{H};G,x_{1}x_{2})\right] -\gamma _{2}B(x_{1}\otimes 1_{H};G,x_{2})%
\end{array}%
=0.  \label{G,gx1x2, GF2,x2}
\end{equation}

\subparagraph{Case $f=gx_{1}x_{2}$}

\begin{equation*}
\begin{array}{c}
B(gx_{1}x_{2}\otimes 1_{H};1_{A},gx_{1}x_{2}) \\
+\gamma _{1}B(gx_{1}x_{2}\otimes 1_{H};GX_{1},gx_{1}x_{2})+\gamma
_{2}B(gx_{1}x_{2}\otimes 1_{H};GX_{2},gx_{1}x_{2})%
\end{array}%
=B(x_{1}x_{2}\otimes g;1_{A},gx_{1}x_{2}).
\end{equation*}%
By applying $\left( \ref{eq.10}\right) $ this rewrites as%
\begin{equation*}
\begin{array}{c}
B(gx_{1}x_{2}\otimes 1_{H};1_{A},gx_{1}x_{2}) \\
+\gamma _{1}B(gx_{1}x_{2}\otimes 1_{H};GX_{1},gx_{1}x_{2})+\gamma
_{2}B(gx_{1}x_{2}\otimes 1_{H};GX_{2},gx_{1}x_{2})%
\end{array}%
=B(gx_{1}x_{2}\otimes 1_{H};1_{A},gx_{1}x_{2})
\end{equation*}%
i.e.%
\begin{equation*}
\gamma _{1}B(gx_{1}x_{2}\otimes 1_{H};GX_{1},gx_{1}x_{2})+\gamma
_{2}B(gx_{1}x_{2}\otimes 1_{H};GX_{2},gx_{1}x_{2})=0.
\end{equation*}%
In view of the form of the element we get%
\begin{equation}
\gamma _{1}B(x_{2}\otimes 1_{H};G,gx_{1}x_{2})-\gamma _{2}B(x_{1}\otimes
1_{H};G,gx_{1}x_{2})=0.  \label{G,gx1x2, GF2,gx1x2}
\end{equation}

\paragraph{\textbf{Equality }$\left( \protect\ref{GF3}\right) $}

rewrites as%
\begin{equation*}
\gamma _{2}B(gx_{1}x_{2}\otimes 1_{H};X_{1}X_{2},f)-\alpha
B(gx_{1}x_{2}\otimes 1_{H};GX_{1},f)=\alpha B(gx_{1}x_{2}g\otimes g;GX_{1},f)
\end{equation*}%
i.e.%
\begin{equation*}
\gamma _{2}B(gx_{1}x_{2}\otimes 1_{H};X_{1}X_{2},f)-\alpha
B(gx_{1}x_{2}\otimes 1_{H};GX_{1},f)=\alpha B(x_{1}x_{2}\otimes g;GX_{1},f).
\end{equation*}

\subparagraph{Case $f=g$}

\begin{equation*}
\gamma _{2}B(gx_{1}x_{2}\otimes 1_{H};X_{1}X_{2},g)-\alpha
B(gx_{1}x_{2}\otimes 1_{H};GX_{1},g)=\alpha B(x_{1}x_{2}\otimes g;GX_{1},g).
\end{equation*}%
By applying $\left( \ref{eq.10}\right) $ this rewrites as%
\begin{equation*}
\gamma _{2}B(gx_{1}x_{2}\otimes 1_{H};X_{1}X_{2},g)-\alpha
B(gx_{1}x_{2}\otimes 1_{H};GX_{1},g)=\alpha B(gx_{1}x_{2}\otimes
1_{H};GX_{1},g)
\end{equation*}%
i.e.%
\begin{equation*}
\gamma _{2}B(gx_{1}x_{2}\otimes 1_{H};X_{1}X_{2},g)-2\alpha
B(gx_{1}x_{2}\otimes 1_{H};GX_{1},g)=0.
\end{equation*}%
In view of the form of the element we get%
\begin{gather}
\gamma _{2}\left[ B(g\otimes 1_{H};1_{A},g)+B(x_{2}\otimes \
1_{H};1_{A},gx_{2})+B(x_{1}\otimes 1_{H};1_{A},gx_{1})+B(gx_{1}x_{2}\otimes
1_{H};1_{A},gx_{1}x_{2})\right] +  \label{G,gx1x2, GF3,g} \\
2\alpha \left[ -B(x_{2}\otimes 1_{H};G,g)+B(gx_{1}x_{2}\otimes
1_{H};G,gx_{1})\right] =0.  \notag
\end{gather}

\subparagraph{Case $f=x_{1}$}

\begin{equation*}
\gamma _{2}B(gx_{1}x_{2}\otimes 1_{H};X_{1}X_{2},x_{1})-\alpha
B(gx_{1}x_{2}\otimes 1_{H};GX_{1},x_{1})=\alpha B(x_{1}x_{2}\otimes
g;GX_{1},x_{1}).
\end{equation*}%
By applying $\left( \ref{eq.10}\right) $ this rewrites as%
\begin{equation*}
\gamma _{2}B(gx_{1}x_{2}\otimes 1_{H};X_{1}X_{2},x_{1})-\alpha
B(gx_{1}x_{2}\otimes 1_{H};GX_{1},x_{1})=-\alpha B(gx_{1}x_{2}\otimes
1_{H};GX_{1},x_{1})
\end{equation*}%
i.e.%
\begin{equation*}
B(gx_{1}x_{2}\otimes 1_{H};X_{1}X_{2},x_{1})=0.
\end{equation*}%
In view of the form of the element we get%
\begin{equation}
\gamma _{2}\left[ B(g\otimes 1_{H};1_{A},x_{1})+B(x_{2}\otimes \
1_{H};1_{A},x_{1}x_{2})\right] =0.  \label{G,gx1x2, GF3,x1}
\end{equation}

\subparagraph{Case $f=x_{2}$}

\begin{equation*}
\gamma _{2}B(gx_{1}x_{2}\otimes 1_{H};X_{1}X_{2},x_{2})-\alpha
B(gx_{1}x_{2}\otimes 1_{H};GX_{1},x_{2})=\alpha B(x_{1}x_{2}\otimes
g;GX_{1},x_{2}).
\end{equation*}%
By applying $\left( \ref{eq.10}\right) $ this rewrites as%
\begin{equation*}
\gamma _{2}B(gx_{1}x_{2}\otimes 1_{H};X_{1}X_{2},x_{2})-\alpha
B(gx_{1}x_{2}\otimes 1_{H};GX_{1},x_{2})=-\alpha B(gx_{1}x_{2}\otimes
1_{H};GX_{1},x_{2})
\end{equation*}%
i.e.%
\begin{equation*}
\gamma _{2}B(gx_{1}x_{2}\otimes 1_{H};X_{1}X_{2},x_{2})=0
\end{equation*}%
In view of the form of the element we get%
\begin{equation}
\gamma _{2}\left[ B(g\otimes 1_{H};X_{2},1_{H})-B(x_{1}\otimes
1_{H};1_{A},x_{1}x_{2})\right] =0  \label{G,gx1x2, GF3,x2}
\end{equation}

\subparagraph{Case $f=gx_{1}x_{2}$}

\begin{eqnarray*}
&&\gamma _{2}B(gx_{1}x_{2}\otimes 1_{H};X_{1}X_{2},gx_{1}x_{2})-\alpha
B(gx_{1}x_{2}\otimes 1_{H};GX_{1},gx_{1}x_{2}) \\
&=&\alpha B(x_{1}x_{2}\otimes g;GX_{1},gx_{1}x_{2}).
\end{eqnarray*}%
By applying $\left( \ref{eq.10}\right) $ this rewrites as%
\begin{eqnarray*}
&&\gamma _{2}B(gx_{1}x_{2}\otimes 1_{H};X_{1}X_{2},gx_{1}x_{2})-\alpha
B(gx_{1}x_{2}\otimes 1_{H};GX_{1},gx_{1}x_{2}) \\
&=&\alpha B(gx_{1}x_{2}\otimes 1_{H};GX_{1},gx_{1}x_{2})
\end{eqnarray*}%
i.e.%
\begin{equation*}
\gamma _{2}B(gx_{1}x_{2}\otimes 1_{H};X_{1}X_{2},gx_{1}x_{2})-2\alpha
B(gx_{1}x_{2}\otimes 1_{H};GX_{1},gx_{1}x_{2})=0.
\end{equation*}%
In view of the form of the element we get%
\begin{equation}
\gamma _{2}B\left( g\otimes 1_{H};1_{A},gx_{1}x_{2}\right) -2\alpha
B(x_{2}\otimes 1_{H};G,gx_{1}x_{2})=0.  \label{G,gx1x2, GF3,gx1x2}
\end{equation}

\paragraph{\textbf{Equality }$\left( \protect\ref{GF4}\right) $}

rewrites as%
\begin{equation*}
-\gamma _{1}B(gx_{1}x_{2}\otimes 1_{H};X_{1}X_{2},f)-\alpha
B(gx_{1}x_{2}\otimes 1_{H};GX_{2},f)=\alpha B(gx_{1}x_{2}g\otimes g;GX_{2},f)
\end{equation*}%
i.e.%
\begin{equation*}
-\gamma _{1}B(gx_{1}x_{2}\otimes 1_{H};X_{1}X_{2},f)-\alpha
B(gx_{1}x_{2}\otimes 1_{H};GX_{2},f)=\alpha B(x_{1}x_{2}\otimes g;GX_{2},f).
\end{equation*}

\subparagraph{Case $f=g$}

\begin{equation*}
-\gamma _{1}B(gx_{1}x_{2}\otimes 1_{H};X_{1}X_{2},g)-\alpha
B(gx_{1}x_{2}\otimes 1_{H};GX_{2},g)=\alpha B(x_{1}x_{2}\otimes g;GX_{2},g).
\end{equation*}%
By applying $\left( \ref{eq.10}\right) $ this rewrites as%
\begin{equation*}
-\gamma _{1}B(gx_{1}x_{2}\otimes 1_{H};X_{1}X_{2},g)-\alpha
B(gx_{1}x_{2}\otimes 1_{H};GX_{2},g)=\alpha B(gx_{1}x_{2}\otimes
1_{H};GX_{2},g)
\end{equation*}%
i.e.%
\begin{equation*}
-\gamma _{1}B(gx_{1}x_{2}\otimes 1_{H};X_{1}X_{2},g)-2\alpha
B(gx_{1}x_{2}\otimes 1_{H};GX_{2},g)=0
\end{equation*}%
In view of the form of the element we get%
\begin{gather}
-\gamma _{1}\left[
\begin{array}{c}
B(g\otimes 1_{H};1_{A},g)+B(x_{2}\otimes \ 1_{H};1_{A},gx_{2}) \\
+B(x_{1}\otimes 1_{H};1_{A},gx_{1})+B(gx_{1}x_{2}\otimes
1_{H};1_{A},gx_{1}x_{2})%
\end{array}%
\right] +  \label{G,gx1x2, GF4,g} \\
-2\alpha \left[ -B(x_{1}\otimes 1_{H};G,g)-B(gx_{1}x_{2}\otimes
1_{H};G,gx_{2})\right] =0  \notag
\end{gather}

\subparagraph{Case $f=x_{1}$%
\protect\begin{equation*}
-\protect\gamma _{1}B(gx_{1}x_{2}\otimes 1_{H};X_{1}X_{2},x_{1})-\protect%
\alpha B(gx_{1}x_{2}\otimes 1_{H};GX_{2},x_{1})=\protect\alpha %
B(x_{1}x_{2}\otimes g;GX_{2},x_{1}).
\protect\end{equation*}%
}

By applying $\left( \ref{eq.10}\right) $ this rewrites as%
\begin{equation*}
-\gamma _{1}B(gx_{1}x_{2}\otimes 1_{H};X_{1}X_{2},x_{1})-\alpha
B(gx_{1}x_{2}\otimes 1_{H};GX_{2},x_{1})=-\alpha B(gx_{1}x_{2}\otimes
1_{H};GX_{2},x_{1})
\end{equation*}%
i.e.%
\begin{equation*}
-\gamma _{1}B(gx_{1}x_{2}\otimes 1_{H};X_{1}X_{2},x_{1})=0
\end{equation*}%
In view of the form of the element we get%
\begin{equation}
\gamma _{1}\left[ B(g\otimes 1_{H};1_{A},x_{1})+B(x_{2}\otimes \
1_{H};1_{A},x_{1}x_{2})\right] =0  \label{G,gx1x2, GF4,x1}
\end{equation}

\subparagraph{Case $f=x_{2}$%
\protect\begin{eqnarray*}
&&-\protect\gamma _{1}B(gx_{1}x_{2}\otimes 1_{H};X_{1}X_{2},x_{2})-\protect%
\alpha B(gx_{1}x_{2}\otimes 1_{H};GX_{2},x_{2}) \\
W &=&\protect\alpha B(x_{1}x_{2}\otimes g;GX_{2},x_{2}).
\protect\end{eqnarray*}%
}

By applying $\left( \ref{eq.10}\right) $ this rewrites as%
\begin{eqnarray*}
&&-\gamma _{1}B(gx_{1}x_{2}\otimes 1_{H};X_{1}X_{2},x_{2})-\alpha
B(gx_{1}x_{2}\otimes 1_{H};GX_{2},x_{2}) \\
&=&-\alpha B(gx_{1}x_{2}\otimes 1_{H};GX_{2},x_{2})
\end{eqnarray*}%
In view of the form of the element we get%
\begin{equation}
-\gamma _{1}\left[ B(g\otimes 1_{H};1_{A},x_{2})-B(x_{1}\otimes
1_{H};1_{A},x_{1}x_{2})\right] =0.  \label{G,gx1x2, GF4,x2}
\end{equation}

\subparagraph{Case $f=gx_{1}x_{2}$}

\begin{eqnarray*}
&&-\gamma _{1}B(gx_{1}x_{2}\otimes 1_{H};X_{1}X_{2},gx_{1}x_{2})-\alpha
B(gx_{1}x_{2}\otimes 1_{H};GX_{2},gx_{1}x_{2}) \\
&=&\alpha B(x_{1}x_{2}\otimes g;GX_{2},gx_{1}x_{2}).
\end{eqnarray*}%
By applying $\left( \ref{eq.10}\right) $ this rewrites as%
\begin{eqnarray*}
&&-\gamma _{1}B(gx_{1}x_{2}\otimes 1_{H};X_{1}X_{2},gx_{1}x_{2})-\alpha
B(gx_{1}x_{2}\otimes 1_{H};GX_{2},gx_{1}x_{2}) \\
&=&\alpha B(gx_{1}x_{2}\otimes 1_{H};GX_{2},gx_{1}x_{2})
\end{eqnarray*}%
i.e.%
\begin{equation*}
-\gamma _{1}B(gx_{1}x_{2}\otimes 1_{H};X_{1}X_{2},gx_{1}x_{2})-2\alpha
B(gx_{1}x_{2}\otimes 1_{H};GX_{2},gx_{1}x_{2})=0
\end{equation*}%
In view of the form of the element we get%
\begin{equation}
-\gamma _{1}B\left( g\otimes 1_{H};1_{A},gx_{1}x_{2}\right) +2\alpha
B(x_{1}\otimes 1_{H};G,gx_{1}x_{2})=0.  \label{G,gx1x2, GF4,gx1x2}
\end{equation}

\paragraph{\textbf{Equality }$\left( \protect\ref{GF5}\right) $}

rewrites as%
\begin{equation*}
\alpha B(gx_{1}x_{2}\otimes 1_{H};GX_{1}X_{2},f)=\alpha
B(gx_{1}x_{2}g\otimes g;GX_{1}X_{2},f)
\end{equation*}%
i.e.%
\begin{equation*}
\alpha B(gx_{1}x_{2}\otimes 1_{H};GX_{1}X_{2},f)=\alpha B(x_{1}x_{2}\otimes
g;GX_{1}X_{2},f).
\end{equation*}

\subparagraph{Case $f=1_{H}$}

\begin{equation*}
\alpha B(gx_{1}x_{2}\otimes 1_{H};GX_{1}X_{2},1_{H})=\alpha
B(x_{1}x_{2}\otimes g;GX_{1}X_{2},1_{H}).
\end{equation*}%
By applying $\left( \ref{eq.10}\right) $ this rewrites as%
\begin{equation*}
\alpha B(gx_{1}x_{2}\otimes 1_{H};GX_{1}X_{2},1_{H})=\alpha
B(gx_{1}x_{2}\otimes 1_{H};GX_{1}X_{2},1_{H})
\end{equation*}%
which is trivial.

\subparagraph{Case $f=x_{1}x_{2}$}

\begin{equation*}
\alpha B(gx_{1}x_{2}\otimes 1_{H};GX_{1}X_{2},x_{1}x_{2})=\alpha
B(x_{1}x_{2}\otimes g;GX_{1}X_{2},x_{1}x_{2}).
\end{equation*}%
By applying $\left( \ref{eq.10}\right) $ this rewrites as%
\begin{equation*}
\alpha B(gx_{1}x_{2}\otimes 1_{H};GX_{1}X_{2},x_{1}x_{2})=\alpha
B(gx_{1}x_{2}\otimes 1_{H};GX_{1}X_{2},x_{1}x_{2})
\end{equation*}%
which is trivial.

\subparagraph{Case $f=gx_{1}$}

\begin{equation*}
\alpha B(gx_{1}x_{2}\otimes 1_{H};GX_{1}X_{2},gx_{1})=\alpha
B(x_{1}x_{2}\otimes g;GX_{1}X_{2},gx_{1}).
\end{equation*}%
By applying $\left( \ref{eq.10}\right) $ this rewrites as%
\begin{equation*}
\alpha B(gx_{1}x_{2}\otimes 1_{H};GX_{1}X_{2},gx_{1})=-\alpha
B(gx_{1}x_{2}\otimes 1_{H};GX_{1}X_{2},gx_{1})
\end{equation*}%
i.e.%
\begin{equation*}
2\alpha B(gx_{1}x_{2}\otimes 1_{H};GX_{1}X_{2},gx_{1})=0.
\end{equation*}%
In view of the form of the element we get%
\begin{equation*}
\alpha \left[ B(g\otimes 1_{H};G,gx_{1})+B(x_{2}\otimes \
1_{H};G,gx_{1}x_{2})\right] =0.
\end{equation*}%
this is $\left( \ref{G,x2, GF5,g}\right) .$

\subparagraph{Case $f=gx_{2}$}

\begin{equation*}
\alpha B(gx_{1}x_{2}\otimes 1_{H};GX_{1}X_{2},gx_{2})=\alpha
B(x_{1}x_{2}\otimes g;GX_{1}X_{2},gx_{2}).
\end{equation*}%
By applying $\left( \ref{eq.10}\right) $ this rewrites as%
\begin{equation*}
\alpha B(gx_{1}x_{2}\otimes 1_{H};GX_{1}X_{2},gx_{2})=-\alpha
B(gx_{1}x_{2}\otimes 1_{H};GX_{1}X_{2},gx_{2})
\end{equation*}%
i.e.%
\begin{equation*}
2\alpha B(gx_{1}x_{2}\otimes 1_{H};GX_{1}X_{2},gx_{2})=0.
\end{equation*}%
In view of the form of the element we get%
\begin{equation}
\alpha \left[ B(g\otimes 1_{H};G,gx_{2})-B(x_{1}\otimes 1_{H};G,gx_{1}x_{2})%
\right] =0.  \label{G,gx1x2, GF5,gx2}
\end{equation}

\paragraph{\textbf{Equality }$\left( \protect\ref{GF6}\right) $}

rewrites as%
\begin{equation*}
-B(gx_{1}x_{2}\otimes 1_{H};X_{1},f)+\gamma _{2}B(gx_{1}x_{2}\otimes
1_{H};GX_{1}X_{2},f)=B(gx_{1}x_{2}g\otimes g;X_{1},f)
\end{equation*}%
i.e.%
\begin{equation*}
-B(gx_{1}x_{2}\otimes 1_{H};X_{1},f)+\gamma _{2}B(gx_{1}x_{2}\otimes
1_{H};GX_{1}X_{2},f)=B(x_{1}x_{2}\otimes g;X_{1},f).
\end{equation*}

\subparagraph{Case $f=1_{H}$}

\begin{equation*}
-B(gx_{1}x_{2}\otimes 1_{H};X_{1},1_{H})+\gamma _{2}B(gx_{1}x_{2}\otimes
1_{H};GX_{1}X_{2},1_{H})=B(x_{1}x_{2}\otimes g;X_{1},1_{H}).
\end{equation*}

By applying $\left( \ref{eq.10}\right) $ this rewrites as%
\begin{equation*}
-B(gx_{1}x_{2}\otimes 1_{H};X_{1},1_{H})+\gamma _{2}B(gx_{1}x_{2}\otimes
1_{H};GX_{1}X_{2},1_{H})=B(gx_{1}x_{2}\otimes 1_{H};X_{1},1_{H})
\end{equation*}%
i.e.%
\begin{equation*}
-2B(gx_{1}x_{2}\otimes 1_{H};X_{1},1_{H})+\gamma _{2}B(gx_{1}x_{2}\otimes
1_{H};GX_{1}X_{2},1_{H})=0.
\end{equation*}%
In view of the form of the element we get%
\begin{gather}
-2\left[ -B(x_{2}\otimes 1_{H};1_{A},1_{H})+B(gx_{1}x_{2}\otimes
1_{H};1_{A},x_{1})\right] +  \label{G,gx1x2, GF6,1H} \\
\gamma _{2}\left[ B(g\otimes 1_{H};G,1_{H})+B(x_{2}\otimes \
1_{H};G,x_{2})+B(x_{1}\otimes 1_{H};G,x_{1})+B(gx_{1}x_{2}\otimes
1_{H};G,x_{1}x_{2})\right]  \notag \\
=0  \notag
\end{gather}

\subparagraph{Case $f=x_{1}x_{2}$}

\begin{equation*}
-B(gx_{1}x_{2}\otimes 1_{H};X_{1},x_{1}x_{2})+\gamma
_{2}B(gx_{1}x_{2}\otimes 1_{H};GX_{1}X_{2},x_{1}x_{2})=B(x_{1}x_{2}\otimes
g;X_{1},x_{1}x_{2}).
\end{equation*}

By applying $\left( \ref{eq.10}\right) $ this rewrites as%
\begin{equation*}
-B(gx_{1}x_{2}\otimes 1_{H};X_{1},x_{1}x_{2})+\gamma
_{2}B(gx_{1}x_{2}\otimes 1_{H};GX_{1}X_{2},x_{1}x_{2})=B(gx_{1}x_{2}\otimes
1_{H};X_{1},x_{1}x_{2})
\end{equation*}%
i.e.%
\begin{equation*}
-2B(gx_{1}x_{2}\otimes 1_{H};X_{1},x_{1}x_{2})+\gamma
_{2}B(gx_{1}x_{2}\otimes 1_{H};GX_{1}X_{2},x_{1}x_{2})=0.
\end{equation*}%
In view of the form of the element we get%
\begin{equation*}
\gamma _{2}B(g\otimes 1_{H};G,x_{1}x_{2})+2B(x_{2}\otimes
1_{H};1_{A},x_{1}x_{2})=0
\end{equation*}%
which is $\left( \ref{G,x2, GF2,x1x2}\right) .$

\subparagraph{Case $f=gx_{1}$}

\begin{equation*}
-B(gx_{1}x_{2}\otimes 1_{H};X_{1},gx_{1})+\gamma _{2}B(gx_{1}x_{2}\otimes
1_{H};GX_{1}X_{2},gx_{1})=B(x_{1}x_{2}\otimes g;X_{1},gx_{1}).
\end{equation*}%
By applying $\left( \ref{eq.10}\right) $ this rewrites as%
\begin{equation*}
-B(gx_{1}x_{2}\otimes 1_{H};X_{1},gx_{1})+\gamma _{2}B(gx_{1}x_{2}\otimes
1_{H};GX_{1}X_{2},gx_{1})=-B(gx_{1}x_{2}\otimes 1_{H};X_{1},gx_{1})
\end{equation*}%
i.e.%
\begin{equation*}
\gamma _{2}B(gx_{1}x_{2}\otimes 1_{H};GX_{1}X_{2},gx_{1})=0.
\end{equation*}%
In view of the form of the element we get%
\begin{equation}
\gamma _{2}\left[ B(g\otimes 1_{H};G,gx_{1})+B(x_{2}\otimes \
1_{H};G,gx_{1}x_{2})\right] =0.  \label{G,gx1x2, GF6,gx1}
\end{equation}%
.

\subparagraph{Case $f=gx_{2}$}

\begin{equation*}
-B(gx_{1}x_{2}\otimes 1_{H};X_{1},gx_{2})+\gamma _{2}B(gx_{1}x_{2}\otimes
1_{H};GX_{1}X_{2},gx_{2})=B(x_{1}x_{2}\otimes g;X_{1},gx_{2}).
\end{equation*}%
By applying $\left( \ref{eq.10}\right) $ this rewrites as%
\begin{equation*}
-B(gx_{1}x_{2}\otimes 1_{H};X_{1},gx_{2})+\gamma _{2}B(gx_{1}x_{2}\otimes
1_{H};GX_{1}X_{2},gx_{2})=-B(gx_{1}x_{2}\otimes 1_{H};X_{1},gx_{2})
\end{equation*}%
i.e.%
\begin{equation*}
\gamma _{2}B(gx_{1}x_{2}\otimes 1_{H};GX_{1}X_{2},gx_{2})=0.
\end{equation*}%
In view of the form of the element we get%
\begin{equation}
\gamma _{2}\left[ B(g\otimes 1_{H};G,gx_{2})-B(x_{1}\otimes
1_{H};G,gx_{1}x_{2})\right] =0.  \label{G,gx1x2, GF6,gx2}
\end{equation}%
.

\paragraph{\textbf{Equality }$\left( \protect\ref{GF7}\right) $}

rewrites as%
\begin{equation*}
-B(gx_{1}x_{2}\otimes 1_{H};X_{2},f)-\gamma _{1}B(gx_{1}x_{2}\otimes
1_{H};GX_{1}X_{2},f)=B(gx_{1}x_{2}g\otimes g;X_{2},f)
\end{equation*}%
i.e.%
\begin{equation*}
-B(gx_{1}x_{2}\otimes 1_{H};X_{2},f)-\gamma _{1}B(gx_{1}x_{2}\otimes
1_{H};GX_{1}X_{2},f)=B(x_{1}x_{2}\otimes g;X_{2},f).
\end{equation*}

\subparagraph{Case $f=1_{H}$}

\begin{equation*}
-B(gx_{1}x_{2}\otimes 1_{H};X_{2},1_{H})-\gamma _{1}B(gx_{1}x_{2}\otimes
1_{H};GX_{1}X_{2},1_{H})=B(x_{1}x_{2}\otimes g;X_{2},1_{H}).
\end{equation*}%
By applying $\left( \ref{eq.10}\right) $ this rewrites as%
\begin{equation*}
-B(gx_{1}x_{2}\otimes 1_{H};X_{2},1_{H})-\gamma _{1}B(gx_{1}x_{2}\otimes
1_{H};GX_{1}X_{2},1_{H})=B(gx_{1}x_{2}\otimes 1_{H};X_{2},1_{H})
\end{equation*}%
i.e.%
\begin{equation*}
-2B(gx_{1}x_{2}\otimes 1_{H};X_{2},1_{H})-\gamma _{1}B(gx_{1}x_{2}\otimes
1_{H};GX_{1}X_{2},1_{H})=0.
\end{equation*}%
In view of the form of the element we get%
\begin{gather}
-\gamma _{1}\left[ B(g\otimes 1_{H};G,1_{H})+B(x_{2}\otimes \
1_{H};G,x_{2})+B(x_{1}\otimes 1_{H};G,x_{1})+B(gx_{1}x_{2}\otimes
1_{H};G,x_{1}x_{2})\right]  \label{G,gx1x2, GF7,1H} \\
-2\left[ B(x_{1}\otimes 1_{H};1_{A},1_{H})+B(gx_{1}x_{2}\otimes
1_{H};1_{A},x_{2})\right] =0  \notag
\end{gather}

\subparagraph{Case $f=x_{1}x_{2}$}

\begin{equation*}
-B(gx_{1}x_{2}\otimes 1_{H};X_{2},x_{1}x_{2})-\gamma
_{1}B(gx_{1}x_{2}\otimes 1_{H};GX_{1}X_{2},x_{1}x_{2})=B(x_{1}x_{2}\otimes
g;X_{2},x_{1}x_{2}).
\end{equation*}%
By applying $\left( \ref{eq.10}\right) $ this rewrites as%
\begin{equation*}
-B(gx_{1}x_{2}\otimes 1_{H};X_{2},x_{1}x_{2})-\gamma
_{1}B(gx_{1}x_{2}\otimes 1_{H};GX_{1}X_{2},x_{1}x_{2})=B(gx_{1}x_{2}\otimes
1_{H};X_{2},x_{1}x_{2})
\end{equation*}%
i.e.%
\begin{equation*}
-2B(gx_{1}x_{2}\otimes 1_{H};X_{2},x_{1}x_{2})-\gamma
_{1}B(gx_{1}x_{2}\otimes 1_{H};GX_{1}X_{2},x_{1}x_{2})=0.
\end{equation*}%
In view of the form of the element we get%
\begin{equation}
-\gamma _{1}B(g\otimes 1_{H};G,x_{1}x_{2})-2B(x_{1}\otimes
1_{H};1_{A},x_{1}x_{2})=0.  \label{G,gx1x2, GF7,x1x2}
\end{equation}

\subparagraph{Case $f=gx_{1}$}

\begin{equation*}
-B(gx_{1}x_{2}\otimes 1_{H};X_{2},gx_{1})-\gamma _{1}B(gx_{1}x_{2}\otimes
1_{H};GX_{1}X_{2},gx_{1})=B(x_{1}x_{2}\otimes g;X_{2},gx_{1}).
\end{equation*}%
By applying $\left( \ref{eq.10}\right) $ this rewrites as%
\begin{equation*}
-B(gx_{1}x_{2}\otimes 1_{H};X_{2},gx_{1})-\gamma _{1}B(gx_{1}x_{2}\otimes
1_{H};GX_{1}X_{2},gx_{1})=-B(gx_{1}x_{2}\otimes 1_{H};X_{2},gx_{1})
\end{equation*}%
i.e.%
\begin{equation*}
\gamma _{1}B(gx_{1}x_{2}\otimes 1_{H};GX_{1}X_{2},gx_{1})=0.
\end{equation*}%
In view of the form of the element we get%
\begin{equation}
\gamma _{1}\left[ B(g\otimes 1_{H};G,gx_{1})+B(x_{2}\otimes \
1_{H};G,gx_{1}x_{2})\right] =0.  \label{G,gx1x2, GF7,gx1}
\end{equation}

\subparagraph{Case $f=gx_{2}$}

\begin{equation*}
-B(gx_{1}x_{2}\otimes 1_{H};X_{2},gx_{2})-\gamma _{1}B(gx_{1}x_{2}\otimes
1_{H};GX_{1}X_{2},gx_{2})=B(x_{1}x_{2}\otimes g;X_{2},gx_{2}).
\end{equation*}%
By applying $\left( \ref{eq.10}\right) $ this rewrites as%
\begin{equation*}
-B(gx_{1}x_{2}\otimes 1_{H};X_{2},gx_{2})-\gamma _{1}B(gx_{1}x_{2}\otimes
1_{H};GX_{1}X_{2},gx_{2})=-B(gx_{1}x_{2}\otimes 1_{H};X_{2},gx_{2})
\end{equation*}%
i.e.%
\begin{equation*}
\gamma _{1}B(gx_{1}x_{2}\otimes 1_{H};GX_{1}X_{2},gx_{2})=0.
\end{equation*}%
In view of the form of the element we get%
\begin{equation}
\gamma _{1}\left[ B(g\otimes 1_{H};G,gx_{2})-B(x_{1}\otimes
1_{H};G,gx_{1}x_{2})\right] =0.  \label{G,gx1x2, GF7,gx2}
\end{equation}

\paragraph{\textbf{Equality }$\left( \protect\ref{GF8}\right) $}

rewrites as%
\begin{equation*}
B(gx_{1}x_{2}\otimes 1_{H};X_{1}X_{2},f)=B(gx_{1}x_{2}g\otimes
g;X_{1}X_{2},f)
\end{equation*}%
i.e.%
\begin{equation*}
B(gx_{1}x_{2}\otimes 1_{H};X_{1}X_{2},f)=B(x_{1}x_{2}\otimes g;X_{1}X_{2},f).
\end{equation*}

\subparagraph{Case $f=g$}

\begin{equation*}
B(gx_{1}x_{2}\otimes 1_{H};X_{1}X_{2},g)=B(x_{1}x_{2}\otimes g;X_{1}X_{2},g).
\end{equation*}%
By applying $\left( \ref{eq.10}\right) $ this rewrites as%
\begin{equation*}
B(gx_{1}x_{2}\otimes 1_{H};X_{1}X_{2},g)=B(gx_{1}x_{2}\otimes
1_{H};X_{1}X_{2},g)
\end{equation*}%
which is trivial.

\subparagraph{Case $f=x_{1}$}

\begin{equation*}
B(gx_{1}x_{2}\otimes 1_{H};X_{1}X_{2},x_{1})=B(x_{1}x_{2}\otimes
g;X_{1}X_{2},x_{1}).
\end{equation*}%
By applying $\left( \ref{eq.10}\right) $ this rewrites as%
\begin{equation*}
2B(gx_{1}x_{2}\otimes 1_{H};X_{1}X_{2},x_{1})=0.
\end{equation*}%
In view of the form of the element we get%
\begin{equation*}
B(g\otimes 1_{H};1_{A},x_{1})+B(x_{2}\otimes \ 1_{H};1_{A},x_{1}x_{2})=0
\end{equation*}%
which is $\left( \ref{G,x2; GF7,x1}\right) .$

\subparagraph{Case $f=x_{2}$}

\begin{equation*}
B(gx_{1}x_{2}\otimes 1_{H};X_{1}X_{2},x_{2})=B(x_{1}x_{2}\otimes
g;X_{1}X_{2},x_{2}).
\end{equation*}%
By applying $\left( \ref{eq.10}\right) $ this rewrites as%
\begin{equation*}
2B(gx_{1}x_{2}\otimes 1_{H};X_{1}X_{2},x_{2})=0.
\end{equation*}%
In view of the form of the element we get%
\begin{equation}
B(g\otimes 1_{H};1_{A},x_{2})-B(x_{1}\otimes 1_{H};1_{A},x_{1}x_{2})=0.
\label{G,gx1x2, GF8,x2}
\end{equation}

\subparagraph{Case $f=gx_{1}x_{2}$}

\begin{equation*}
B(gx_{1}x_{2}\otimes 1_{H};X_{1}X_{2},gx_{1}x_{2})=B(x_{1}x_{2}\otimes
g;X_{1}X_{2},gx_{1}x_{2}).
\end{equation*}%
By applying $\left( \ref{eq.10}\right) $ this is trivial.

\subsection{$d=X_{i}$}

We write $\left( \ref{eq.a}\right) $ for $d=X_{i}$%
\begin{equation*}
X_{i}\otimes g\otimes g+1_{A}\otimes x_{i}\otimes g+1_{A}\otimes
1_{H}\otimes x_{i}
\end{equation*}%
\begin{eqnarray*}
B(h\otimes h^{\prime })(X_{i}\otimes 1_{H}) &=&(X_{i}\otimes
1_{H})B(hg\otimes h^{\prime }g)+ \\
&&+(1_{A}\otimes 1_{H})B(hx_{i}\otimes h^{\prime }g)+ \\
&&+(1_{A}\otimes 1_{H})B(h\otimes h^{\prime }x_{i})
\end{eqnarray*}

\subsection{$d=X_{1}$}

We write

\begin{eqnarray*}
B(h\otimes h^{\prime }) &=&\sum_{f}[B(h\otimes h^{\prime
};1_{A},f)1_{A}\otimes f+B(h\otimes h^{\prime };G,f)G\otimes f+ \\
&&B(h\otimes h^{\prime };X_{1},f)X_{1}\otimes f+B(h\otimes h^{\prime
};X_{2},f)X_{2}\otimes f. \\
&&+B(h\otimes h^{\prime };X_{1}X_{2},f)X_{1}X_{2}\otimes f+B(h\otimes
h^{\prime };GX_{1},f)GX_{1}\otimes f+ \\
&&B(h\otimes h^{\prime };GX_{2},f)GX_{2}\otimes f+B(h\otimes h^{\prime
};GX_{1}X_{2},f)GX_{1}X_{2}\otimes ]f.
\end{eqnarray*}

We recall that%
\begin{eqnarray*}
1_{A}X_{1} &=&X_{1} \\
GX_{1} &=&GX_{1} \\
X_{1}X_{1} &=&\beta _{1} \\
X_{2}X_{1} &=&\lambda -X_{1}X_{2} \\
X_{1}X_{2}X_{1} &=&X_{1}\left( \lambda -X_{1}X_{2}\right) =\lambda
X_{1}-\beta _{1}X_{2} \\
GX_{1}X_{1} &=&\beta _{1}G \\
GX_{2}X_{1} &=&G\left( \lambda -X_{1}X_{2}\right) =\lambda G-GX_{1}X_{2} \\
GX_{1}X_{2}X_{1} &=&G\left( \lambda X_{1}-\beta _{1}X_{2}\right) =\lambda
GX_{1}-\beta _{1}GX_{2}
\end{eqnarray*}

We write the left side%
\begin{eqnarray*}
B(h\otimes h^{\prime })\left( X_{1}\otimes 1_{H}\right) = &&\sum_f\left[
\beta _{1}B(h\otimes h^{\prime };X_{1},f)+\lambda B(h\otimes h^{\prime
};X_{2},f)\right] 1_{A}\otimes f+ \\
&&\sum_f\left[ \beta _{1}B(h\otimes h^{\prime };GX_{1},f)+\lambda B(h\otimes
h^{\prime };GX_{2},f)\right] G\otimes f+ \\
&&\sum_f\left[ B(h\otimes h^{\prime };1_{A},f)+\lambda B(h\otimes h^{\prime
};X_{1}X_{2},f)\right] X_{1}\otimes f+ \\
&&\sum_f -\beta _{1}B(h\otimes h^{\prime };X_{1}X_{2},f)X_{2}\otimes f+ \\
&&\sum_f -B(h\otimes h^{\prime };X_{2},f)X_{1}X_{2}\otimes f+ \\
&&\sum_f\left[ B(h\otimes h^{\prime };G,f)+\lambda B(h\otimes h^{\prime
};GX_{1}X_{2},f)\right] GX_{1}\otimes f+ \\
&& \sum_f -\beta _1B(h\otimes h^{\prime };GX_{1}X_{2},f)GX_{2}\otimes f \\
&&\sum_f -B(h\otimes h^{\prime };GX_{2},f)GX_{1}X_{2}\otimes f+
\end{eqnarray*}

We recall that%
\begin{eqnarray*}
X_{1}1_{A} &=&X_{1} \\
X_{1}G &=&-GX_{1}+\gamma _{1} \\
X_{1}X_{1} &=&\beta _{1} \\
X_{1}X_{2} &=&X_{1}X_{2} \\
X_{1}X_{1}X_{2} &=&\beta _{1}X_{2} \\
X_{1}GX_{1} &=&-\beta _{1}G+\gamma _{1}X_{1} \\
X_{1}GX_{2} &=&-GX_{1}X_{2}+\gamma _{1}X_{2} \\
X_{1}GX_{1}X_{2} &.=&-\beta _{1}GX_{2}+\gamma _{1}X_{1}X_{2}
\end{eqnarray*}

We write the right side%
\begin{eqnarray*}
&&1\left( X_{1}\otimes 1_{H}\right) B(hg\otimes h^{\prime }g)= \\
&&[\gamma _{1}B(hg\otimes h^{\prime }g;G,f)+\beta _{1}\sum_{f}B(hg\otimes
h^{\prime }g;X_{1},f)]1_{A}\otimes f+ \\
&&-\beta _{1}\sum_{f}B(hg\otimes h^{\prime }g;GX_{1},f)G\otimes f+ \\
&&\sum_{f}[B(hg\otimes h^{\prime }g;1_{A},f)+B(hg\otimes h^{\prime
}g;GX_{1},f)]X_{1}\otimes f+ \\
&&\sum_{f}[\beta _{1}B(hg\otimes h^{\prime }g;X_{1}X_{2},f)+\gamma
_{1}B(hg\otimes h^{\prime }g;GX_{2},f)]X_{2}\otimes f+ \\
&&\sum_{f}[B(hg\otimes h^{\prime }g;X_{2},f)+\gamma _{1}B(hg\otimes
h^{\prime }g;GX_{1}X_{2},f)X_{1}X_{2}\otimes f+ \\
&&-\sum_{d,e_{1},e_{2}=0}^{1}B(hg\otimes h^{\prime }g;G,f)GX_{1}\otimes f+ \\
&&-\beta _{1}\sum_{f}B(hg\otimes h^{\prime }g;GX_{1}X_{2},f)GX_{2}\otimes f
\\
&&-\sum_{f}B(hg\otimes h^{\prime }g;GX_{2},f)GX_{1}X_{2}\otimes f+ \\
&&
\end{eqnarray*}%
\begin{eqnarray*}
\text{2}B(hx_{1}\otimes h^{\prime }g)
&=&\sum_{a,b_{1},b_{2},d,e_{1},e_{2}=0}^{1}B(hx_{1}\otimes h^{\prime
}g;G^{a}X_{1}^{b_{1}}X_{2}^{b_{2}},f)G^{a}X_{1}^{b_{1}}X_{2}^{b_{2}}\otimes f
\\
&=&\sum_{d,e_{1},e_{2}=0}^{1}B(hx_{1}\otimes h^{\prime
}g;1_{A},f)1_{A}\otimes f \\
&&\sum_{d,e_{1},e_{2}=0}^{1}B(hx_{1}\otimes h^{\prime }g;G,f)G\otimes f \\
&&+\sum_{d,e_{1},e_{2}=0}^{1}B(hx_{1}\otimes h^{\prime
}g;X_{1},f)X_{1}\otimes f \\
&&+\sum_{d,e_{1},e_{2}=0}^{1}B(hx_{1}\otimes h^{\prime
}g;GX_{1},f)GX_{1}\otimes f \\
&&+\sum_{d,e_{1},e_{2}=0}^{1}B(hx_{1}\otimes h^{\prime
}g;X_{2},f)X_{2}\otimes f \\
&&+\sum_{d,e_{1},e_{2}=0}^{1}B(hx_{1}\otimes h^{\prime
}g;GX_{2},f)GX_{2}\otimes f \\
&&+\sum_{d,e_{1},e_{2}=0}^{1}B(hx_{1}\otimes h^{\prime
}g;X_{1}X_{2},f)X_{1}X_{2}\otimes f \\
&&+\sum_{d,e_{1},e_{2}=0}^{1}B(hx_{1}\otimes h^{\prime
}g;GX_{1}X_{2},f)GX_{1}X_{2}\otimes f
\end{eqnarray*}%
\begin{eqnarray*}
\text{3}B(h\otimes h^{\prime }x_{1})
&=&\sum_{a,b_{1},b_{2},d,e_{1},e_{2}=0}^{1}B(h\otimes h^{\prime
}x_{1};G^{a}X_{1}^{b_{1}}X_{2}^{b_{2}},f)G^{a}X_{1}^{b_{1}}X_{2}^{b_{2}}%
\otimes f \\
&=&\sum_{d,e_{1},e_{2}=0}^{1}B(h\otimes h^{\prime
}x_{1};1_{A},f)1_{A}\otimes f \\
&&\sum_{d,e_{1},e_{2}=0}^{1}B(h\otimes h^{\prime }x_{1};G,f)G\otimes f \\
&&+\sum_{d,e_{1},e_{2}=0}^{1}B(h\otimes h^{\prime
}x_{1};X_{1},f)X_{1}\otimes f \\
&&+\sum_{d,e_{1},e_{2}=0}^{1}B(h\otimes h^{\prime
}x_{1};GX_{1},f)GX_{1}\otimes f \\
&&+\sum_{d,e_{1},e_{2}=0}^{1}B(h\otimes h^{\prime
}x_{1};X_{2},f)X_{2}\otimes f \\
&&+\sum_{d,e_{1},e_{2}=0}^{1}B(h\otimes h^{\prime
}x_{1};GX_{2},f)GX_{2}\otimes f \\
&&+\sum_{d,e_{1},e_{2}=0}^{1}B(h\otimes h^{\prime
}x_{1};X_{1}X_{2},f)X_{1}X_{2}\otimes f \\
&&+\sum_{d,e_{1},e_{2}=0}^{1}B(h\otimes h^{\prime
}x_{1};GX_{1}X_{2},f)GX_{1}X_{2}\otimes f
\end{eqnarray*}%
so that we deduce that%
\begin{gather}
\beta _{1}B(h\otimes h^{\prime };X_{1},f)+\lambda B(h\otimes h^{\prime
};X_{2},f)=\gamma _{1}B(hg\otimes h^{\prime }g;G,f)  \notag \\
+\beta _{1}B(hg\otimes h^{\prime }g;X_{1},f)+B(hx_{1}\otimes h^{\prime
}g;1_{A},f)+B(h\otimes h^{\prime }x_{1};1_{A},f)  \label{X1F11}
\end{gather}%
\begin{eqnarray}
&&\beta _{1}B(h\otimes h^{\prime };GX_{1},f)+\lambda B(h\otimes h^{\prime
};GX_{2},f)  \label{X1F21} \\
&=&-\beta _{1}B(hg\otimes h^{\prime }g;GX_{1},f)+B(hx_{1}\otimes h^{\prime
}g;G,f)+B(h\otimes h^{\prime }x_{1};G,f)  \notag
\end{eqnarray}

\begin{gather}
B(h\otimes h^{\prime };1_{A},f)+\lambda B(h\otimes h^{\prime };X_{1}X_{2},f)=
\notag \\
B(hg\otimes h^{\prime }g;1_{A},f)+\gamma _{1}B(hg\otimes h^{\prime
}g;GX_{1},f)  \label{X1F31} \\
+B(hx_{1}\otimes h^{\prime }g;X_{1},f)+B(h\otimes h^{\prime }x_{1};X_{1},f)
\notag
\end{gather}%
\begin{gather}
-\beta _{1}B(h\otimes h^{\prime };X_{1}X_{2},f)=\beta _{1}B(hg\otimes
h^{\prime }g;X_{1}X_{2},f)  \label{X1F41} \\
+\gamma _{1}B(hg\otimes h^{\prime }g;GX_{2},f)+B(hx_{1}\otimes h^{\prime
}g;X_{2},f)+B(h\otimes h^{\prime }x_{1};X_{2},f)  \notag
\end{gather}%
\begin{eqnarray}
-B(h\otimes h^{\prime };X_{2},f) &=&B(hg\otimes h^{\prime }g;X_{2},f)+\gamma
_{1}B(hg\otimes h^{\prime }g;GX_{1}X_{2},f)  \label{X1F51} \\
&&+B(hx_{1}\otimes h^{\prime }g;X_{1}X_{2},f)+B(h\otimes h^{\prime
}x_{1};X_{1}X_{2},f)  \notag
\end{eqnarray}%
\begin{eqnarray}
&&B(h\otimes h^{\prime };G,f)+\lambda B(h\otimes h^{\prime };GX_{1}X_{2},f)
\label{X1F61} \\
&=&-B(hg\otimes h^{\prime }g;G,f)+B(hx_{1}\otimes h^{\prime
}g;GX_{1},f)+B(h\otimes h^{\prime }x_{1};GX_{1},f)  \notag
\end{eqnarray}%
\begin{gather}
-\beta _{1}B(h\otimes h^{\prime };GX_{1}X_{2},f)=-\beta _{1}B(hg\otimes
h^{\prime }g;GX_{1}X_{2},f)  \label{X1F71} \\
+B(hx_{1}\otimes h^{\prime }g;GX_{2},f)+B(h\otimes h^{\prime }x_{1};GX_{2},f)
\notag
\end{gather}%
\begin{gather}
-B(h\otimes h^{\prime };GX_{2},f)=-B(hg\otimes h^{\prime }g;GX_{2},f)
\label{X1F81} \\
+B(hx_{1}\otimes h^{\prime }g;GX_{1}X_{2},f)+B(h\otimes h^{\prime
}x_{1};GX_{1}X_{2},f)  \notag
\end{gather}

\subsubsection{Case $g\otimes 1_{H}$}

\paragraph{Equality $\left( \protect\ref{X1F11}\right) $}

rewrites as%
\begin{gather*}
\beta _{1}B(g\otimes 1_{H};X_{1},f)+\lambda B(g\otimes 1_{H};X_{2},f)=\gamma
_{1}B(1_{H}\otimes g;G,f) \\
+\beta _{1}B(1_{H}\otimes g;X_{1},f)+B(gx_{1}\otimes g;1_{A},f)+B(g\otimes
x_{1};1_{A},f)
\end{gather*}

\subparagraph{Case $f=1_{H}$}

\begin{gather*}
\beta _{1}B(g\otimes 1_{H};X_{1},1_{H})+\lambda B(g\otimes
1_{H};X_{2},1_{H})=\gamma _{1}B(1_{H}\otimes g;G,1_{H}) \\
\beta _{1}B(1_{H}\otimes g;X_{1},1_{H})+B(gx_{1}\otimes
g;1_{A},1_{H})+B(g\otimes x_{1};1_{A},1_{H}).
\end{gather*}%
By applying $\left( \ref{eq.10}\right) $ this rewrites as%
\begin{gather*}
\beta _{1}B(g\otimes 1_{H};X_{1},1_{H})+\lambda B(g\otimes
1_{H};X_{2},1_{H})-\gamma _{1}B(g\otimes 1_{H};G,1_{H}) \\
-\beta _{1}B(g\otimes 1_{H};X_{1},1_{H})-B(x_{1}\otimes
1_{H};1_{A},1_{H})-B(1_{H}\otimes gx_{1};1_{A},1_{H})=0.
\end{gather*}%
i.e.%
\begin{equation*}
\lambda B(g\otimes 1_{H};X_{2},1_{H})-\gamma _{1}B(g\otimes
1_{H};G,1_{H})-B(x_{1}\otimes 1_{H};1_{A},1_{H})-B(1_{H}\otimes
gx_{1};1_{A},1_{H})=0
\end{equation*}%
In view of the form of the elements we get%
\begin{equation}
\lambda B\left( g\otimes 1_{H};1_{A},x_{2}\right) -\gamma _{1}B(g\otimes
1_{H};G,1_{H})-2B(x_{1}\otimes 1_{H};1_{A},1_{H})=0.  \label{X1,g, X1F11,1H}
\end{equation}

\subparagraph{Case $f=g$}

\begin{eqnarray*}
&&\beta _{1}B(g\otimes 1_{H};X_{1},g)+\lambda B(g\otimes 1_{H};X_{2},g) \\
&=&\gamma _{1}B(1_{H}\otimes g;G,g)+\beta _{1}B(1_{H}\otimes
g;X_{1},g)+B(gx_{1}\otimes g;1_{A},g)+B(g\otimes x_{1};1_{A},g).
\end{eqnarray*}%
By applying $\left( \ref{eq.10}\right) $ this rewrites as%
\begin{eqnarray*}
&&\beta _{1}B(g\otimes 1_{H};X_{1},g)+\lambda B(g\otimes 1_{H};X_{2},g) \\
&=&\gamma _{1}B(g\otimes 1_{H};G,g)+\beta _{1}B(g\otimes
1_{H};X_{1},g)+B(x_{1}\otimes 1_{H};1_{A},g)+B(1_{H}\otimes gx_{1};1_{A},g)
\end{eqnarray*}%
i.e.%
\begin{equation*}
\lambda B(g\otimes 1_{H};X_{2},g)-\gamma _{1}B(g\otimes
1_{H};G,g)-B(x_{1}\otimes 1_{H};1_{A},g)-B(1_{H}\otimes gx_{1};1_{A},g)=0.
\end{equation*}%
In view of the form of the elements this is trivial.

\subparagraph{Case $f=x_{1}x_{2}$}

\begin{eqnarray*}
&&\beta _{1}B(g\otimes 1_{H};X_{1},x_{1}x_{2})+\lambda B(g\otimes
1_{H};X_{2},x_{1}x_{2}) \\
&=&\beta _{1}B(1_{H}\otimes g;X_{1},x_{1}x_{2})+B(gx_{1}\otimes
g;1_{A},x_{1}x_{2})+B(g\otimes x_{1};1_{A},x_{1}x_{2})
\end{eqnarray*}

By applying $\left( \ref{eq.10}\right) $ this rewrites as%
\begin{eqnarray*}
&&\lambda B(g\otimes 1_{H};X_{2},x_{1}x_{2}) \\
&=&\gamma _{1}B(g\otimes 1_{H};G,x_{1}x_{2})+B(x_{1}\otimes
1_{H};1_{A},x_{1}x_{2})+B(1_{H}\otimes gx_{1};1_{A},x_{1}x_{2}).
\end{eqnarray*}%
In view of the form of the elements we get%
\begin{equation*}
\gamma _{1}B(g\otimes 1_{H};G,x_{1}x_{2})+2B\left( x_{1}\otimes
1_{H};1_{A},x_{1}x_{2}\right) =0.
\end{equation*}%
which is $\left( \ref{G,x1, GF2,x1x2}\right) .$

\subparagraph{Case $f=gx_{1}$%
\protect\begin{gather*}
\protect\beta _{1}B(g\otimes 1_{H};X_{1},gx_{1})+\protect\lambda B(g\otimes
1_{H};X_{2},gx_{1})=\protect\gamma _{1}B(1_{H}\otimes g;G,gx_{1}) \\
+\protect\beta _{1}B(1_{H}\otimes g;X_{1},gx_{1})+B(gx_{1}\otimes
g;1_{A},gx_{1})+B(g\otimes x_{1};1_{A},gx_{1})
\protect\end{gather*}%
}

By applying $\left( \ref{eq.10}\right) $ this rewrites as%
\begin{gather*}
\beta _{1}B(g\otimes 1_{H};X_{1},gx_{1})+\lambda B(g\otimes
1_{H};X_{2},gx_{1})=-\gamma _{1}B(g\otimes 1_{H};G,gx_{1}) \\
-\beta _{1}B(g\otimes 1_{H};X_{1},gx_{1})-B(x_{1}\otimes
1_{H};1_{A},gx_{1})-B(1_{H}\otimes gx_{1};1_{A},gx_{1})
\end{gather*}%
In view of the form of the elements we get%
\begin{eqnarray}
&&+\gamma _{1}B(g\otimes 1_{H};G,gx_{1})+\lambda B(g\otimes
1_{H};X_{2},gx_{1})+2B(x_{1}\otimes 1_{H};1_{A},gx_{1})+
\label{X1,g, X1F11,gx1} \\
&&+2B\left( g\otimes 1_{H};1_{A},g\right) =0  \notag
\end{eqnarray}

\subparagraph{Case $f=gx_{2}$%
\protect\begin{eqnarray*}
&&\protect\beta _{1}B(g\otimes 1_{H};X_{1},gx_{2})+\protect\lambda %
B(g\otimes 1_{H};X_{2},gx_{2}) \\
&=&\protect\gamma _{1}B(1_{H}\otimes g;G,gx_{2})+\protect\beta %
_{1}B(1_{H}\otimes g;X_{1},gx_{2})+ \\
&&+B(gx_{1}\otimes g;1_{A},gx_{2})+B(g\otimes x_{1};1_{A},gx_{2}).
\protect\end{eqnarray*}%
}

By applying $\left( \ref{eq.10}\right) $ this rewrites as In view of the
form of the elements we get%
\begin{gather*}
\beta _{1}B(g\otimes 1_{H};X_{1},gx_{2})+\lambda B(g\otimes
1_{H};X_{2},gx_{2})= \\
-\gamma _{1}B(g\otimes 1_{H};G,gx_{2})-\beta _{1}B(g\otimes
1_{H};X_{1},gx_{2}) \\
-B(x_{1}\otimes 1_{H};1_{A},gx_{2})-B(1_{H}\otimes gx_{1};1_{A},gx_{2}).
\end{gather*}%
In view of the form of the elements, we get%
\begin{equation}
\gamma _{1}B(g\otimes 1_{H};G,gx_{2})+2B(x_{1}\otimes
1_{H};1_{A},gx_{2})-2\beta _{1}B\left( g\otimes
1_{H};1_{A},gx_{1}x_{2}\right) =0.  \label{X1,g, X1F11,gx2}
\end{equation}

\paragraph{Equality $\left( \protect\ref{X1F21}\right) $}

rewrites as%
\begin{eqnarray*}
&&\beta _{1}B(g\otimes 1_{H};GX_{1},f)+\lambda B(g\otimes 1_{H};GX_{2},f) \\
&=&-\beta _{1}B(1_{H}\otimes g;GX_{1},f)+B(gx_{1}\otimes g;G,f)+B(g\otimes
x_{1};G,f).
\end{eqnarray*}

\subparagraph{Case $f=g$}

\begin{eqnarray*}
&&\beta _{1}B(g\otimes 1_{H};GX_{1},g)+\lambda B(g\otimes 1_{H};GX_{2},g) \\
&=&-\beta _{1}B(1_{H}\otimes g;GX_{1},g)+B(gx_{1}\otimes g;G,g)+B(g\otimes
x_{1};G,g).
\end{eqnarray*}%
By applying $\left( \ref{eq.10}\right) $ this rewrites as%
\begin{eqnarray*}
&&\beta _{1}B(g\otimes 1_{H};GX_{1},g)+\lambda B(g\otimes 1_{H};GX_{2},g) \\
&=&-\beta _{1}B(g\otimes 1_{H};GX_{1},g)+B(x_{1}\otimes
1_{H};G,g)+B(1_{H}\otimes gx_{1};G,g).
\end{eqnarray*}%
In view of the form of the elements we get

\begin{equation}
2\beta _{1}B\left( g\otimes 1_{H};G,gx_{1}\right) +\lambda B\left( g\otimes
1_{H};G,gx_{2}\right) +2B(x_{1}\otimes 1_{H};G,g)=0.  \label{X1,g,X1F21,g}
\end{equation}

\subparagraph{Case $f=x_{1}$}

\begin{eqnarray*}
&&\beta _{1}B(g\otimes 1_{H};GX_{1},x_{1})+\lambda B(g\otimes
1_{H};GX_{2},x_{1}) \\
&=&-\beta _{1}B(1_{H}\otimes g;GX_{1},x_{1})+B(gx_{1}\otimes
g;G,x_{1})+B(g\otimes x_{1};G,x_{1}).
\end{eqnarray*}%
By applying $\left( \ref{eq.10}\right) $ this rewrites as%
\begin{eqnarray*}
&&\beta _{1}B(g\otimes 1_{H};GX_{1},x_{1})+\lambda B(g\otimes
1_{H};GX_{2},x_{1}) \\
&=&\beta _{1}B(g\otimes 1_{H};GX_{1},x_{1})-B(x_{1}\otimes
1_{H};G,x_{1})-B(1\otimes gx_{1};G,x_{1}).
\end{eqnarray*}

In view of the form of the elements we get%
\begin{equation}
-\lambda B\left( g\otimes 1_{H};G,x_{1}x_{2}\right) +2B(x_{1}\otimes
1_{H};G,x_{1})=0.  \label{X1,g,X1F21,x1}
\end{equation}

\subparagraph{Case $f=x_{2}$}

\begin{eqnarray*}
&&\beta _{1}B(g\otimes 1_{H};GX_{1},x_{2})+\lambda B(g\otimes
1_{H};GX_{2},x_{2}) \\
&=&-\beta _{1}B(1_{H}\otimes g;GX_{1},x_{2})+B(gx_{1}\otimes
g;G,x_{2})+B(g\otimes x_{1};G,x_{2}).
\end{eqnarray*}%
By applying $\left( \ref{eq.10}\right) $ this rewrites as%
\begin{eqnarray*}
&&\beta _{1}B(g\otimes 1_{H};GX_{1},x_{2})+\lambda B(g\otimes
1_{H};GX_{2},x_{2}) \\
&=&\beta _{1}B(g\otimes 1_{H};GX_{1},x_{2})-B(x_{1}\otimes
1_{H};G,x_{2})-B(1\otimes gx_{1};G,x_{2}).
\end{eqnarray*}%
In view of the form of the elements we get%
\begin{equation}
B(x_{1}\otimes 1_{H};G,x_{2})=0.  \label{X1,g,X1F21,x2}
\end{equation}

\subparagraph{Case $f=gx_{1}x_{2}$}

\begin{eqnarray*}
&&\beta _{1}B(g\otimes 1_{H};GX_{1},gx_{1}x_{2})+\lambda B(g\otimes
1_{H};GX_{2},gx_{1}x_{2}) \\
&=&-\beta _{1}B(1_{H}\otimes g;GX_{1},gx_{1}x_{2})+B(gx_{1}\otimes
g;G,gx_{1}x_{2})+B(g\otimes x_{1};G,gx_{1}x_{2}).
\end{eqnarray*}%
By applying $\left( \ref{eq.10}\right) $ this rewrites as%
\begin{eqnarray*}
&&\beta _{1}B(g\otimes 1_{H};GX_{1},gx_{1}x_{2})+\lambda B(g\otimes
1_{H};GX_{2},gx_{1}x_{2}) \\
&=&-\beta _{1}B(g\otimes 1_{H};GX_{1},gx_{1}x_{2})+B(x_{1}\otimes
1_{H};G,gx_{1}x_{2})+B(1_{H}\otimes gx_{1};G,gx_{1}x_{2}).
\end{eqnarray*}%
By the form of the element we get%
\begin{equation}
B\left( g\otimes 1_{H};G,gx_{2}\right) -B\left( x_{1}\otimes
1_{H};G,gx_{1}x_{2}\right) =0.  \label{X1,g,X1F21,gx1x2}
\end{equation}

\paragraph{Equality $\left( \protect\ref{X1F31}\right) $}

rewrites as%
\begin{eqnarray*}
&&B(g\otimes 1_{H};1_{A},f)+\lambda B(g\otimes 1_{H};X_{1}X_{2},f) \\
&=&B(1_{H}\otimes g;1_{A},f)+\gamma _{1}B(1_{H}\otimes g;GX_{1},f)+ \\
&&+B(gx_{1}\otimes g;X_{1},f)+B(g\otimes x_{1};X_{1},f).
\end{eqnarray*}

\subparagraph{Case $f=g$}

\begin{eqnarray*}
&&B(g\otimes 1_{H};1_{A},g)+\lambda B(g\otimes 1_{H};X_{1}X_{2},g) \\
&=&B(1_{H}\otimes g;1_{A},g)+B(gx_{1}\otimes g;X_{1},g)+B(g\otimes
x_{1};X_{1},g).
\end{eqnarray*}%
By applying $\left( \ref{eq.10}\right) $ this rewrites as%
\begin{eqnarray*}
&&B(g\otimes 1_{H};1_{A},g)+\lambda B(g\otimes 1_{H};X_{1}X_{2},g) \\
&=&B(g\otimes 1_{H};1_{A},g)+\gamma _{1}B(g\otimes 1_{H};GX_{1},g) \\
&&+B(x_{1}\otimes 1_{H};X_{1},g)+B(1_{H}\otimes gx_{1};X_{1},g).
\end{eqnarray*}%
In view of the form of the elements we get%
\begin{gather}
+\lambda B\left( g\otimes 1_{H};1_{A},gx_{1}x_{2}\right) +\gamma _{1}B\left(
g\otimes 1_{H};G,gx_{1}\right) +  \label{X1,g, X1F31,g} \\
+2B(g\otimes 1_{H};1_{A},g)+2B(x_{1}\otimes 1_{H};1_{A},gx_{1})=0.  \notag
\end{gather}

\subparagraph{Case $f=x_{1}$}

\begin{eqnarray*}
&&B(g\otimes 1_{H};1_{A},x_{1})+\lambda B(g\otimes 1_{H};X_{1}X_{2},x_{1}) \\
&=&B(1_{H}\otimes g;1_{A},x_{1})+\gamma _{1}B(1_{H}\otimes g;GX_{1},x_{1}) \\
&&+B(gx_{1}\otimes g;X_{1},x_{1})+B(g\otimes x_{1};X_{1},x_{1}).
\end{eqnarray*}%
By applying $\left( \ref{eq.10}\right) $ this rewrites as%
\begin{eqnarray*}
&&2B(g\otimes 1_{H};1_{A},x_{1})+\lambda B(g\otimes 1_{H};X_{1}X_{2},x_{1})
\\
&=&-B(x_{1}\otimes 1_{H};X_{1},x_{1})-\gamma _{1}B(g\otimes
1_{H};GX_{1},x_{1}) \\
&&-B(1_{H}\otimes gx_{1};X_{1},x_{1}).
\end{eqnarray*}%
In view of the form of the elements we get%
\begin{equation*}
2B(g\otimes 1_{H};1_{A},x_{1})=2B(g\otimes 1_{H};1_{A},x_{1})
\end{equation*}%
which is trivial.

\subparagraph{Case $f=x_{2}$}

\begin{eqnarray*}
&&B(g\otimes 1_{H};1_{A},x_{2})+\lambda B(g\otimes 1_{H};X_{1}X_{2},x_{2}) \\
&=&B(1_{H}\otimes g;1_{A},x_{2})+\gamma _{1}B(1_{H}\otimes g;GX_{1},x_{2}) \\
&&+B(gx_{1}\otimes g;X_{1},x_{2})+B(g\otimes x_{1};X_{1},x_{2}).
\end{eqnarray*}%
By applying $\left( \ref{eq.10}\right) $ this rewrites as%
\begin{eqnarray*}
&&2B(g\otimes 1_{H};1_{A},x_{2})+\lambda B(g\otimes 1_{H};X_{1}X_{2},x_{2})
\\
&=&-B(x_{1}\otimes 1_{H};X_{1},x_{2})-\gamma _{1}B(g\otimes
1_{H};GX_{1},x_{2})-B(1_{H}\otimes gx_{1};X_{1},x_{2}).
\end{eqnarray*}%
In view of the form of the elements we get%
\begin{equation*}
2B(x_{1}\otimes 1_{H};1_{A},x_{1}x_{2})+\gamma _{1}B\left( g\otimes
1_{H};G,x_{1}x_{2}\right) =0.
\end{equation*}%
which is $\left( \ref{G,g, GF2,x2}\right) .$

\subparagraph{Case $f=gx_{1}x_{2}$}

\begin{eqnarray*}
&&B(g\otimes 1_{H};1_{A},gx_{1}x_{2})+\lambda B(g\otimes
1_{H};X_{1}X_{2},gx_{1}x_{2}) \\
&=&B(1_{H}\otimes g;1_{A},gx_{1}x_{2})+\gamma _{1}B(1_{H}\otimes
g;GX_{1},gx_{1}x_{2}) \\
&&+B(gx_{1}\otimes g;X_{1},gx_{1}x_{2})+B(g\otimes x_{1};X_{1},gx_{1}x_{2}).
\end{eqnarray*}%
By applying $\left( \ref{eq.10}\right) $ this rewrites as%
\begin{eqnarray*}
&&B(g\otimes 1_{H};1_{A},gx_{1}x_{2})+\lambda B(g\otimes
1_{H};X_{1}X_{2},gx_{1}x_{2}) \\
&=&B(g\otimes 1_{H};1_{A},gx_{1}x_{2})+\gamma _{1}B(g\otimes
1_{H};GX_{1},gx_{1}x_{2}) \\
&&+B(x_{1}\otimes 1_{H};X_{1},gx_{1}x_{2})+B(1_{H}\otimes
gx_{1};X_{1},gx_{1}x_{2}).
\end{eqnarray*}%
In view of the form of the elements we get%
\begin{equation*}
-B(g\otimes 1_{H};1_{A},gx_{1}x_{2})+B\left( g\otimes
1_{H};1_{A},gx_{1}x_{2}\right) =0.
\end{equation*}%
which is trivial.

\paragraph{Equality $\left( \protect\ref{X1F41}\right) $}

rewrites as%
\begin{eqnarray*}
&&-\beta _{1}B(g\otimes 1_{H};X_{1}X_{2},f) \\
&=&\beta _{1}B(1_{H}\otimes g;X_{1}X_{2},f)+\gamma _{1}B(1_{H}\otimes
g;GX_{2},f) \\
&&+B(gx_{1}\otimes g;X_{2},f)+B(g\otimes x_{1};X_{2},f).
\end{eqnarray*}

\subparagraph{Case $f=g$}

\begin{eqnarray*}
&&-\beta _{1}B(g\otimes 1_{H};X_{1}X_{2},g) \\
&=&\beta _{1}B(1_{H}\otimes g;X_{1}X_{2},g)+ \\
&&+\gamma _{1}B(1_{H}\otimes g;GX_{2},g)+ \\
&&+B(gx_{1}\otimes g;X_{2},g)+B(g\otimes x_{1};X_{2},g).
\end{eqnarray*}%
By applying $\left( \ref{eq.10}\right) $ this rewrites as%
\begin{eqnarray*}
&&-2\beta _{1}B(g\otimes 1_{H};X_{1}X_{2},g) \\
&=&B(x_{1}\otimes 1_{H};X_{2},g)+\gamma _{1}B(g\otimes
1_{H};GX_{2},g)+B(1_{H}\otimes gx_{1};X_{2},g).
\end{eqnarray*}%
In view of the form of the elements we get%
\begin{equation}
-2\beta _{1}B\left( g\otimes 1_{H};1_{A},gx_{1}x_{2}\right) +2B(x_{1}\otimes
1_{H};1_{A},gx_{2})+\gamma _{1}B\left( g\otimes 1_{H};G,gx_{2}\right) =0.
\label{X1,g, X1F41,g}
\end{equation}

$.$

\subparagraph{Case $f=x_{1}$}

\begin{eqnarray*}
&&-\beta _{1}B(g\otimes 1_{H};X_{1}X_{2},x_{1}) \\
&=&\beta _{1}B(1_{H}\otimes g;X_{1}X_{2},x_{1})+\gamma _{1}B(1_{H}\otimes
g;GX_{2},x_{1}) \\
&&+B(gx_{1}\otimes g;X_{2},x_{1})+B(g\otimes x_{1};X_{2},x_{1}).
\end{eqnarray*}%
By applying $\left( \ref{eq.10}\right) $ this rewrites as%
\begin{equation*}
0=-B(x_{1}\otimes 1_{H};X_{2},x_{1})-\gamma _{1}B(g\otimes
1_{H};GX_{2},x_{1})-B(1_{H}\otimes gx_{1};X_{2},x_{1}).
\end{equation*}%
In view of the form of the elements we get%
\begin{equation}
2B(x_{1}\otimes 1_{H};1_{A},x_{1}x_{2})+\gamma _{1}B\left( g\otimes
1_{H};G,x_{1}x_{2}\right) =0.  \label{X1,g, X1F41,x1}
\end{equation}

\paragraph{Equality $\left( \protect\ref{X1F51}\right) $}

rewrites as%
\begin{eqnarray*}
&&-B(g\otimes 1_{H};X_{2},f) \\
&=&B(1_{H}\otimes g;X_{2},f)+\gamma _{1}B(1_{H}\otimes g;GX_{1}X_{2},f) \\
&&+B(gx_{1}\otimes g;X_{1}X_{2},f)+B(g\otimes x_{1};X_{1}X_{2},f).
\end{eqnarray*}

\subparagraph{Case $f=1_{H}$}

\begin{eqnarray*}
&&-B(g\otimes 1_{H};X_{2},1_{H}) \\
&=&B(1_{H}\otimes g;X_{2},1_{H})+\gamma _{1}B(1_{H}\otimes
g;GX_{1}X_{2},1_{H})+ \\
&&+B(gx_{1}\otimes g;X_{1}X_{2},1_{H})+B(g\otimes x_{1};X_{1}X_{2},1_{H}).
\end{eqnarray*}%
By applying $\left( \ref{eq.10}\right) $ this rewrites as%
\begin{eqnarray*}
&&-2B(g\otimes 1_{H};X_{2},1_{H}) \\
&=&\gamma _{1}B(g\otimes 1_{H};GX_{1}X_{2},1_{H})+ \\
&&+B(x_{1}\otimes 1_{H};X_{1}X_{2},1_{H})+B(1_{H}\otimes
gx_{1};X_{1}X_{2},1_{H}).
\end{eqnarray*}%
In view of the form of the elements we get%
\begin{eqnarray}
&&-2B\left( g\otimes 1_{H};1_{A},x_{2}\right)  \label{X1,g,X1F51,1H} \\
&=&\gamma _{1}B\left( g\otimes 1_{H};G,x_{1}x_{2}\right) -2B\left( g\otimes
1_{H};1_{A},x_{2}\right) +2B(x_{1}\otimes 1_{H};1_{A},x_{1}x_{2}).  \notag
\end{eqnarray}

\subparagraph{Case $f=gx_{1}$}

\begin{eqnarray*}
&&-B(g\otimes 1_{H};X_{2},gx_{1}) \\
&=&B(1_{H}\otimes g;X_{2},gx_{1})+\gamma _{1}B(1_{H}\otimes
g;GX_{1}X_{2},gx_{1}) \\
&&+B(gx_{1}\otimes g;X_{1}X_{2},gx_{1})+B(g\otimes x_{1};X_{1}X_{2},gx_{1}).
\end{eqnarray*}%
By applying $\left( \ref{eq.10}\right) $ this rewrites as%
\begin{eqnarray*}
&&-B(g\otimes 1_{H};X_{2},gx_{1}) \\
&=&-B(g\otimes 1_{H};X_{2},gx_{1})-\gamma _{1}B(g\otimes
1_{H};GX_{1}X_{2},gx_{1}) \\
&&-B(x_{1}\otimes 1_{H};X_{1}X_{2},gx_{1})-B(1_{H}\otimes
gx_{1};X_{1}X_{2},gx_{1})
\end{eqnarray*}%
i.e.%
\begin{equation*}
0=+B(x_{1}\otimes 1_{H};X_{1}X_{2},gx_{1})+B(1_{H}\otimes
gx_{1};X_{1}X_{2},gx_{1})
\end{equation*}%
In view of the form of the elements we get%
\begin{equation*}
0=-B\left( g\otimes 1_{H};1_{A},gx_{1}x_{2}\right) +B\left( g\otimes
1_{H};1_{A},gx_{1}x_{2}\right)
\end{equation*}%
which is trivial.

\paragraph{Equality $\left( \protect\ref{X1F61}\right) $}

rewrites as%
\begin{eqnarray*}
&&B(g\otimes 1_{H};G,f)+\lambda B(g\otimes 1_{H};GX_{1}X_{2},f) \\
&=&-B(1_{H}\otimes g;G,f)+B(gx_{1}\otimes g;GX_{1},f)+B(g\otimes
x_{1};GX_{1},f).
\end{eqnarray*}

\subparagraph{Case $f=1_{H}$%
\protect\begin{eqnarray*}
&&B(g\otimes 1_{H};G,1_{H})+\protect\lambda B(g\otimes
1_{H};GX_{1}X_{2},1_{H}) \\
&=&-B(1_{H}\otimes g;G,1_{H})+B(gx_{1}\otimes g;GX_{1},1_{H})+B(g\otimes
x_{1};GX_{1},1_{H}).
\protect\end{eqnarray*}%
}

By applying $\left( \ref{eq.10}\right) $ this rewrites as%
\begin{eqnarray*}
&&2B(g\otimes 1_{H};G,1_{H})+\lambda B(g\otimes 1_{H};GX_{1}X_{2},1_{H}) \\
&=&B(x_{1}\otimes 1_{H};GX_{1},1_{H})+B(1_{H}\otimes gx_{1};GX_{1},1_{H}).
\end{eqnarray*}%
In view of the form of the elements we get

\begin{eqnarray*}
&&2B(g\otimes 1_{H};G,1_{H})+\lambda B\left( g\otimes
1_{H};G,x_{1}x_{2}\right) \\
&=&2B(g\otimes 1_{H};G,1_{H})+2B(x_{1}\otimes 1_{H};G,x_{1})
\end{eqnarray*}%
i.e.%
\begin{equation*}
\lambda B\left( g\otimes 1_{H};G,x_{1}x_{2}\right) =2B(x_{1}\otimes
1_{H};G,x_{1})
\end{equation*}%
which is $\left( \ref{X1,g,X1F21,x1}\right) .$

\subparagraph{Case $f=x_{1}x_{2}$%
\protect\begin{eqnarray*}
&&B(g\otimes 1_{H};G,x_{1}x_{2})+\protect\lambda B(g\otimes
1_{H};GX_{1}X_{2},x_{1}x_{2}) \\
&=&-B(1_{H}\otimes g;G,x_{1}x_{2})+B(gx_{1}\otimes g;GX_{1},x_{1}x_{2}) \\
&&+B(g\otimes x_{1};GX_{1},x_{1}x_{2}).
\protect\end{eqnarray*}%
}

By applying $\left( \ref{eq.10}\right) $ this rewrites as%
\begin{eqnarray*}
&&2B(g\otimes 1_{H};G,x_{1}x_{2})+\lambda B(g\otimes
1_{H};GX_{1}X_{2},x_{1}x_{2}) \\
&=&B(x_{1}\otimes 1_{H};GX_{1},x_{1}x_{2})+B(1_{H}\otimes
gx_{1};GX_{1},x_{1}x_{2}).
\end{eqnarray*}%
In view of the form of the elements we get%
\begin{equation*}
2B(g\otimes 1_{H};G,x_{1}x_{2})=2B\left( g\otimes 1_{H};G,x_{1}x_{2}\right) .
\end{equation*}

\subparagraph{Case $f=gx_{1}$}

\begin{eqnarray*}
&&B(g\otimes 1_{H};G,gx_{1})+\lambda B(g\otimes 1_{H};GX_{1}X_{2},gx_{1}) \\
&=&-B(1_{H}\otimes g;G,gx_{1})+B(gx_{1}\otimes g;GX_{1},gx_{1})+B(g\otimes
x_{1};GX_{1},gx_{1}).
\end{eqnarray*}%
By applying $\left( \ref{eq.10}\right) $ this rewrites as%
\begin{eqnarray*}
&&B(g\otimes 1_{H};G,gx_{1})+\lambda B(g\otimes 1_{H};GX_{1}X_{2},gx_{1}) \\
&=&B(g\otimes 1_{H};G,gx_{1})-B(x_{1}\otimes
1_{H};GX_{1},gx_{1})-B(1_{H}\otimes gx_{1};GX_{1},gx_{1})
\end{eqnarray*}%
i.e.%
\begin{equation*}
\lambda B(g\otimes 1_{H};GX_{1}X_{2},gx_{1})=-B(x_{1}\otimes
1_{H};GX_{1},gx_{1})-B(1_{H}\otimes gx_{1};GX_{1},gx_{1})
\end{equation*}%
In view of the form of the elements we get%
\begin{equation*}
0=-B\left( g\otimes 1_{H};G,gx_{1}\right) +B\left( g\otimes
1_{H};G,gx_{1}\right) .
\end{equation*}

\subparagraph{Case $f=gx_{2}$}

\begin{eqnarray*}
&&B(g\otimes 1_{H};G,gx_{2})+\lambda B(g\otimes 1_{H};GX_{1}X_{2},gx_{2}) \\
&=&-B(1_{H}\otimes g;G,gx_{2})+B(gx_{1}\otimes g;GX_{1},gx_{2})+B(g\otimes
x_{1};GX_{1},gx_{2}).
\end{eqnarray*}%
By applying $\left( \ref{eq.10}\right) $ this rewrites as%
\begin{eqnarray*}
&&B(g\otimes 1_{H};G,gx_{2})+\lambda B(g\otimes 1_{H};GX_{1}X_{2},gx_{2}) \\
&=&B(g\otimes 1_{H};G,gx_{2})-B(x_{1}\otimes
1_{H};GX_{1},gx_{2})-B(1_{H}\otimes gx_{1};GX_{1},gx_{2})
\end{eqnarray*}%
i.e.%
\begin{equation*}
\lambda B(g\otimes 1_{H};GX_{1}X_{2},gx_{2})=-B(x_{1}\otimes
1_{H};GX_{1},gx_{2})-B(1_{H}\otimes gx_{1};GX_{1},gx_{2})
\end{equation*}%
In view of the form of the elements we get%
\begin{eqnarray*}
0 &=&-\left[ -B(g\otimes 1_{H};GX_{2},g)-B(x_{1}\otimes 1_{H};G,gx_{1}x_{2})%
\right] \\
&&-\left[ B(g\otimes 1_{H};G,gx_{2})-B(x_{1}\otimes 1_{H};G,gx_{1}x_{2})%
\right]
\end{eqnarray*}%
\begin{eqnarray*}
0 &=&B(g\otimes 1_{H};GX_{2},g)+B(x_{1}\otimes 1_{H};G,gx_{1}x_{2}) \\
&&-B(g\otimes 1_{H};G,gx_{2})+B(x_{1}\otimes 1_{H};G,gx_{1}x_{2})
\end{eqnarray*}%
\begin{eqnarray*}
0 &=&-B\left( g\otimes 1_{H};G,gx_{2}\right) +B(x_{1}\otimes
1_{H};G,gx_{1}x_{2}) \\
&&-B(g\otimes 1_{H};G,gx_{2})+B(x_{1}\otimes 1_{H};G,gx_{1}x_{2})
\end{eqnarray*}%
i.e.%
\begin{equation}
B(x_{1}\otimes 1_{H};G,gx_{1}x_{2})-B(g\otimes 1_{H};G,gx_{2})=0.
\label{X1,x2,X1F71,gx1x2}
\end{equation}

\paragraph{Equality $\left( \protect\ref{X1F71}\right) $}

rewrites as%
\begin{eqnarray*}
&&-\beta _{1}B(g\otimes 1_{H};GX_{1}X_{2},f) \\
&=&-\beta _{1}B(1_{H}\otimes g;GX_{1}X_{2},f)+ \\
&&+B(gx_{1}\otimes g;GX_{2},f)+B(g\otimes x_{1};GX_{2},f).
\end{eqnarray*}

\subparagraph{Case $f=1_{H}$}

\begin{eqnarray*}
&&-\beta _{1}B(g\otimes 1_{H};GX_{1}X_{2},1_{H}) \\
&=&-\beta _{1}B(1_{H}\otimes g;GX_{1}X_{2},1_{H})+ \\
&&+B(gx_{1}\otimes g;GX_{2},1_{H})+B(g\otimes x_{1};GX_{2},1_{H}).
\end{eqnarray*}%
By applying $\left( \ref{eq.10}\right) $ this rewrites as%
\begin{eqnarray*}
&&-\beta _{1}B(g\otimes 1_{H};GX_{1}X_{2},1_{H}) \\
&=&-\beta _{1}B(g\otimes 1_{H};GX_{1}X_{2},1_{H})+ \\
&&+B(x_{1}\otimes 1_{H};GX_{2},1_{H})+B(1_{H}\otimes gx_{1};GX_{2},1_{H})
\end{eqnarray*}%
i.e.%
\begin{equation*}
0=B(x_{1}\otimes 1_{H};GX_{2},1_{H})+B(1_{H}\otimes gx_{1};GX_{2},1_{H})
\end{equation*}%
In view of the form of the elements we get%
\begin{equation*}
0=B(x_{1}\otimes 1_{H};G,x_{2})
\end{equation*}%
which is $\left( \ref{X1,g,X1F21,x2}\right) .$

\subparagraph{Case $f=gx_{1}$}

\begin{eqnarray*}
&&-\beta _{1}B(g\otimes 1_{H};GX_{1}X_{2},gx_{1}) \\
&=&-\beta _{1}B(1_{H}\otimes g;GX_{1}X_{2},gx_{1})+ \\
&&+B(gx_{1}\otimes g;GX_{2},gx_{1})+B(g\otimes x_{1};GX_{2},gx_{1}).
\end{eqnarray*}%
By applying $\left( \ref{eq.10}\right) $ this rewrites as

\begin{eqnarray*}
&&2\beta _{1}B(g\otimes 1_{H};GX_{1}X_{2},gx_{1}) \\
&=&B(x_{1}\otimes 1_{H};GX_{2},gx_{1})+B(1_{H}\otimes gx_{1};GX_{2},gx_{1}).
\end{eqnarray*}%
In view of the form of the elements we get%
\begin{gather*}
0=B(x_{1}\otimes 1_{H};G,gx_{1}x_{2})+ \\
+\left[ -2B\left( g\otimes 1_{H};G,gx_{2}\right) +B(x_{1}\otimes
1_{H};G,gx_{1}x_{2})\right]
\end{gather*}%
i.e.%
\begin{equation*}
B(g\otimes 1_{H};G,gx_{2})-B(x_{1}\otimes 1_{H};G,gx_{1}x_{2})=0
\end{equation*}%
which is $\left( \ref{X1,g,X1F21,gx1x2}\right) .$

\paragraph{Equality $\left( \protect\ref{X1F81}\right) $}

rewrites as%
\begin{eqnarray*}
&&-B(g\otimes 1_{H};GX_{2},f) \\
&=&-B(1_{H}\otimes g;GX_{2},f)+ \\
&&+B(gx_{1}\otimes g;GX_{1}X_{2},f)+B(g\otimes x_{1};GX_{1}X_{2},f).
\end{eqnarray*}

\subparagraph{Case $f=g$}

\begin{eqnarray*}
&&-B(g\otimes 1_{H};GX_{2},g) \\
&=&-B(1_{H}\otimes g;GX_{2},g)+B(gx_{1}\otimes g;GX_{1}X_{2},g)+ \\
&&+B(g\otimes x_{1};GX_{1}X_{2},g).
\end{eqnarray*}%
By applying $\left( \ref{eq.10}\right) $ this rewrites as%
\begin{eqnarray*}
&&-B(g\otimes 1_{H};GX_{2},g) \\
&=&-B(g\otimes 1_{H};GX_{2},g)+B(x_{1}\otimes 1_{H};GX_{1}X_{2},g)+ \\
&&+B(1_{H}\otimes gx_{1};GX_{1}X_{2},g)
\end{eqnarray*}%
i.e.%
\begin{equation*}
0=B(x_{1}\otimes 1_{H};GX_{1}X_{2},g)+B(1_{H}\otimes gx_{1};GX_{1}X_{2},g).
\end{equation*}%
In view of the form of the elements we get%
\begin{eqnarray*}
0 &=&\left[ -B(g\otimes 1_{H};G,gx_{2})+B(x_{1}\otimes 1_{H};G,gx_{1}x_{2})%
\right] \\
&&+\left[ -B(g\otimes 1_{H};G,gx_{2})+B(x_{1}\otimes 1_{H};G,gx_{1}x_{2})%
\right]
\end{eqnarray*}%
i.e.%
\begin{equation*}
B(g\otimes 1_{H};G,gx_{2})-B(x_{1}\otimes 1_{H};G,gx_{1}x_{2})=0
\end{equation*}%
which is $\left( \ref{X1,g,X1F21,gx1x2}\right) .$

\subparagraph{Case $f=x_{1}$}

\begin{eqnarray*}
&&-B(g\otimes 1_{H};GX_{2},x_{1}) \\
&=&-B(1_{H}\otimes g;GX_{2},x_{1})+B(gx_{1}\otimes g;GX_{1}X_{2},x_{1}) \\
&&+B(g\otimes x_{1};GX_{1}X_{2},x_{1}).
\end{eqnarray*}%
By applying $\left( \ref{eq.10}\right) $ this rewrites as%
\begin{eqnarray*}
&&2B(g\otimes 1_{H};GX_{2},x_{1}) \\
&=&B(x_{1}\otimes 1_{H};GX_{1}X_{2},x_{1})+B(x_{1}\otimes
1_{H};GX_{1}X_{2},x_{1}).
\end{eqnarray*}%
In view of the form of the elements we get

\begin{eqnarray*}
&&-2B\left( g\otimes 1_{H};G,x_{1}x_{2}\right) \\
&=&-B(g\otimes 1_{H};G,x_{1}x_{2})-B\left( g\otimes
1_{H};G,x_{1}x_{2}\right) .
\end{eqnarray*}

\subsubsection{Case $x_{1}\otimes 1_{H}$}

\paragraph{Equality $\left( \protect\ref{X1F11}\right) $}

rewrites as%
\begin{eqnarray*}
&&\beta _{1}B(x_{1}\otimes 1_{H};X_{1},f)+\lambda B(x_{1}\otimes
1_{H};X_{2},f) \\
&=&\gamma _{1}B(x_{1}g\otimes g;G,f)+\beta _{1}B(x_{1}g\otimes g;X_{1},f)+ \\
&&+B(x_{1}x_{1}\otimes g;1_{A},f)+B(x_{1}\otimes x_{1};1_{A},f)
\end{eqnarray*}%
i.e,%
\begin{eqnarray*}
&&\beta _{1}B(x_{1}\otimes 1_{H};X_{1},f)+\lambda B(x_{1}\otimes
1_{H};X_{2},f) \\
&=&\gamma _{1}B(x_{1}g\otimes g;G,f)-\beta _{1}B(gx_{1}\otimes g;X_{1},f)+ \\
&&+B(x_{1}\otimes x_{1};1_{A},f).
\end{eqnarray*}

\subparagraph{Case $f=g$}

\begin{eqnarray*}
&&\beta _{1}B(x_{1}\otimes 1_{H};X_{1},g)+\lambda B(x_{1}\otimes
1_{H};X_{2},g) \\
&=&\gamma _{1}B(x_{1}g\otimes g;G,g)-\beta _{1}B(gx_{1}\otimes g;X_{1},g)+ \\
&&+B(x_{1}\otimes x_{1};1_{A},g).
\end{eqnarray*}%
By applying $\left( \ref{eq.10}\right) $ this rewrites as%
\begin{eqnarray*}
&&2\beta _{1}B(x_{1}\otimes 1_{H};X_{1},g)+\lambda B(x_{1}\otimes
1_{H};X_{2},g) \\
&=&\gamma _{1}B(x_{1}\otimes 1_{H};G,g)+B(x_{1}\otimes x_{1};1_{H},g).
\end{eqnarray*}%
In view of the form of the elements we get

\begin{eqnarray}
&&2\beta _{1}\left[ -B(g\otimes 1_{H};1_{A},g)-B(x_{1}\otimes
1_{H};1_{A},gx_{1})\right] +  \label{X1,x1,X1F11,g} \\
&&-\lambda B(x_{1}\otimes 1_{H};1_{A},gx_{2})  \notag \\
&=&\gamma _{1}B(x_{1}\otimes 1_{H};G,g)  \notag
\end{eqnarray}

\subparagraph{Case $f=x_{1}$}

\begin{eqnarray*}
&&\beta _{1}B(x_{1}\otimes 1_{H};X_{1},x_{1})+\lambda B(x_{1}\otimes
1_{H};X_{2},x_{1}) \\
&=&\gamma _{1}B(x_{1}g\otimes g;G,x_{1})-\beta _{1}B(gx_{1}\otimes
g;X_{1},x_{1})+ \\
&&+B(x_{1}\otimes x_{1};1_{A},x_{1}).
\end{eqnarray*}%
By applying $\left( \ref{eq.10}\right) $ this rewrites as%
\begin{eqnarray*}
&&\beta _{1}B(x_{1}\otimes 1_{H};X_{1},x_{1})+\lambda B(x_{1}\otimes
1_{H};X_{2},x_{1}) \\
&=&\beta _{1}B(x_{1}\otimes 1_{H};X_{1},x_{1})+B(x_{1}\otimes
x_{1};1_{A},x_{1}).
\end{eqnarray*}%
i.e.%
\begin{equation*}
\lambda B(x_{1}\otimes 1_{H};X_{2},x_{1})=-\gamma _{1}B(x_{1}\otimes
1_{H};G,x_{1})+B(x_{1}\otimes x_{1};1_{A},x_{1}).
\end{equation*}%
In view of the form of the elements we get%
\begin{equation}
-\lambda B(x_{1}\otimes 1_{H};1_{A},x_{1}x_{2})+\gamma _{1}B(x_{1}\otimes
1_{H};G,x_{1})=0.  \label{X1,x1,X1F11,x1}
\end{equation}

\subparagraph{Case $f=x_{2}$}

\begin{eqnarray*}
&&\beta _{1}B(x_{1}\otimes 1_{H};X_{1},x_{2})+\lambda B(x_{1}\otimes
1_{H};X_{2},x_{2}) \\
&=&\gamma _{1}B(x_{1}g\otimes g;G,x_{2})-\beta _{1}B(gx_{1}\otimes
g;X_{1},x_{2})+ \\
&&+B(x_{1}\otimes x_{1};1_{A},x_{2}).
\end{eqnarray*}%
By applying $\left( \ref{eq.10}\right) $ this rewrites as%
\begin{eqnarray*}
&&\beta _{1}B(x_{1}\otimes 1_{H};X_{1},x_{2})+\lambda B(x_{1}\otimes
1_{H};X_{2},x_{2}) \\
&=&-\gamma _{1}B(x_{1}\otimes 1_{H};G,x_{2})+\beta _{1}B(x_{1}\otimes
1_{H};X_{1},x_{2})+ \\
&&+B(x_{1}\otimes x_{1};1_{A},x_{2}).
\end{eqnarray*}%
i.e.%
\begin{equation*}
\lambda B(x_{1}\otimes 1_{H};X_{2},x_{2})=+B(x_{1}\otimes x_{1};1_{A},x_{2}).
\end{equation*}%
In view of the form of the elements we get%
\begin{equation}
0=\gamma _{1}B(x_{1}\otimes 1_{H};G,x_{2}).  \label{X1,x1,X1F11,x2}
\end{equation}

\subparagraph{Case $f=gx_{1}x_{2}$}

\begin{eqnarray*}
&&\beta _{1}B(x_{1}\otimes 1_{H};X_{1},gx_{1}x_{2})+\lambda B(x_{1}\otimes
1_{H};X_{2},gx_{1}x_{2}) \\
&=&-\beta _{1}B(gx_{1}\otimes g;X_{1},gx_{1}x_{2})+B(x_{1}\otimes
x_{1};1_{A},gx_{1}x_{2}).
\end{eqnarray*}%
By applying $\left( \ref{eq.10}\right) $ this rewrites as%
\begin{eqnarray*}
&&2\beta _{1}B(x_{1}\otimes 1_{H};X_{1},gx_{1}x_{2})+\lambda B(x_{1}\otimes
1_{H};X_{2},gx_{1}x_{2}) \\
&=&B(x_{1}\otimes x_{1};1_{H},gx_{1}x_{2}).
\end{eqnarray*}%
In view of the form of the elements we get%
\begin{eqnarray}
&&-2\beta _{1}B(g\otimes 1_{H};1_{A},gx_{1}x_{2})  \label{X1,x1,X1F11,gx1x2}
\\
&=&\gamma _{1}B(x_{1}\otimes 1_{H};G,gx_{1}x_{2})-2B\left( x_{1}\otimes
1_{H};1_{A},gx_{2}\right) .  \notag
\end{eqnarray}

\paragraph{Equality $\left( \protect\ref{X1F21}\right) $}

rewrites as%
\begin{eqnarray*}
&&\beta _{1}B(x_{1}\otimes 1_{H};GX_{1},f)+\lambda B(x_{1}\otimes
1_{H};GX_{2},f) \\
&=&-\beta _{1}B(x_{1}g\otimes g;GX_{1},f)+B(x_{1}x_{1}\otimes h^{\prime
}g;G,f) \\
&&+B(x_{1}\otimes x_{1};G,f)
\end{eqnarray*}%
i.e.%
\begin{eqnarray*}
&&\beta _{1}B(x_{1}\otimes 1_{H};GX_{1},f)+\lambda B(x_{1}\otimes
1_{H};GX_{2},f) \\
&=&-\beta _{1}B(x_{1}g\otimes g;GX_{1},f)+B(x_{1}\otimes x_{1};G,f).
\end{eqnarray*}

\subparagraph{Case $f=1_{H}$}

\begin{eqnarray*}
&&\beta _{1}B(x_{1}\otimes 1_{H};GX_{1},1_{H})+\lambda B(x_{1}\otimes
1_{H};GX_{2},1_{H}) \\
&=&-\beta _{1}B(x_{1}g\otimes g;GX_{1},1_{H})+B(x_{1}\otimes x_{1};G,1_{H}).
\end{eqnarray*}%
By applying $\left( \ref{eq.10}\right) $ this rewrites as%
\begin{equation*}
2\beta _{1}B(x_{1}\otimes 1_{H};GX_{1},1_{H})+\lambda B(x_{1}\otimes
1_{H};GX_{2},1_{H})=B(x_{1}\otimes x_{1};G,1_{H}).
\end{equation*}%
In view of the form of the elements we get%
\begin{equation*}
2\beta _{1}\left[ B(g\otimes 1_{H};G,1_{H})+B(x_{1}\otimes 1_{H};G,x_{1})%
\right] +\lambda B(x_{1}\otimes 1_{H};G,x_{2})=0.
\end{equation*}%
By $\left( \ref{X1,g,X1F21,x2}\right) $ we obtain%
\begin{equation}
\beta _{1}\left[ B(g\otimes 1_{H};G,1_{H})+B(x_{1}\otimes 1_{H};G,x_{1})%
\right] =0  \label{X1,x1,X1F21,1}
\end{equation}

\subparagraph{Case $f=x_{1}x_{2}$}

\begin{eqnarray*}
&&\beta _{1}B(x_{1}\otimes 1_{H};GX_{1},x_{1}x_{2})+\lambda B(x_{1}\otimes
1_{H};GX_{2},x_{1}x_{2}) \\
&=&-\beta _{1}B(x_{1}g\otimes g;GX_{1},1_{H})+B(x_{1}\otimes
x_{1};G,x_{1}x_{2}).
\end{eqnarray*}%
By applying $\left( \ref{eq.10}\right) $ this rewrites as%
\begin{equation*}
2\beta _{1}B(x_{1}\otimes 1_{H};GX_{1},x_{1}x_{2})+\lambda B(x_{1}\otimes
1_{H};GX_{2},x_{1}x_{2})=B(x_{1}\otimes x_{1};G,x_{1}x_{2}).
\end{equation*}%
In view of the form of the elements we get%
\begin{equation}
\beta _{1}B\left( g\otimes 1_{H};G,x_{1}x_{2}\right) =0.
\label{X1,x1,X1F21,x1x2}
\end{equation}

\subparagraph{Case $f=gx_{1}$}

\begin{eqnarray*}
&&\beta _{1}B(x_{1}\otimes 1_{H};GX_{1},gx_{1})+\lambda B(x_{1}\otimes
1_{H};GX_{2},gx_{1}) \\
&=&-\beta _{1}B(x_{1}g\otimes g;GX_{1},gx_{1})+B(x_{1}\otimes
x_{1};G,gx_{1}).
\end{eqnarray*}%
By applying $\left( \ref{eq.10}\right) $ this rewrites as%
\begin{eqnarray*}
&&\beta _{1}B(x_{1}\otimes 1_{H};GX_{1},gx_{1})+\lambda B(x_{1}\otimes
1_{H};GX_{2},gx_{1}) \\
&=&\beta _{1}B(x_{1}\otimes 1_{H};GX_{1},gx_{1})+B(x_{1}\otimes
x_{1};G,gx_{1})
\end{eqnarray*}%
i.e.%
\begin{equation*}
\lambda B(x_{1}\otimes 1_{H};GX_{2},gx_{1})=B(x_{1}\otimes x_{1};G,gx_{1}).
\end{equation*}%
In view of the form of the elements we get%
\begin{equation}
\lambda B(x_{1}\otimes 1_{H};G,gx_{1}x_{2})=-2B\left( x_{1}\otimes
1_{H};G,g\right) .  \label{X1,x1,X1F21,gx1}
\end{equation}

\subparagraph{Case $f=gx_{2}$}

\begin{eqnarray*}
&&\beta _{1}B(x_{1}\otimes 1_{H};GX_{1},gx_{2})+\lambda B(x_{1}\otimes
1_{H};GX_{2},gx_{2}) \\
&=&-\beta _{1}B(x_{1}g\otimes g;GX_{1},gx_{2})+B(x_{1}\otimes
x_{1};G,gx_{2}).
\end{eqnarray*}%
By applying $\left( \ref{eq.10}\right) $ this rewrites as%
\begin{eqnarray*}
&&\beta _{1}B(x_{1}\otimes 1_{H};GX_{1},gx_{2})+\lambda B(x_{1}\otimes
1_{H};GX_{2},gx_{2}) \\
&=&\beta _{1}B(x_{1}\otimes 1_{H};GX_{1},gx_{2})+B(x_{1}\otimes
x_{1};G,gx_{2})
\end{eqnarray*}%
i.e.%
\begin{equation*}
\lambda B(x_{1}\otimes 1_{H};GX_{2},gx_{2})=B(x_{1}\otimes x_{1};G,gx_{2}).
\end{equation*}%
In view of the form of the elements we get nothing new.

\paragraph{Equality $\left( \protect\ref{X1F31}\right) $}

rewrites as%
\begin{eqnarray*}
&&B(x_{1}\otimes 1_{H};1_{A},f)+\lambda B(x_{1}\otimes 1_{H};X_{1}X_{2},f) \\
&=&\gamma _{1}B(x_{1}g\otimes g;GX_{1},f)+B(x_{1}g\otimes g;1_{A},f)+ \\
&&+B(x_{1}\otimes x_{1};X_{1},f).
\end{eqnarray*}

\subparagraph{Case $f=1_{H}$%
\protect\begin{eqnarray*}
&&B(x_{1}\otimes 1_{H};1_{A},1_{H})+\protect\lambda B(x_{1}\otimes
1_{H};X_{1}X_{2},1_{H}) \\
&=&\protect\gamma _{1}B(x_{1}g\otimes g;GX_{1},1_{H})+B(x_{1}g\otimes
g;1_{A},1_{H})+ \\
&&+B(x_{1}\otimes x_{1};X_{1},1_{H}).
\protect\end{eqnarray*}%
}

By applying $\left( \ref{eq.10}\right) $ this rewrites as%
\begin{eqnarray*}
&&2B(x_{1}\otimes 1_{H};1_{A},1_{H})+\lambda B(x_{1}\otimes
1_{H};X_{1}X_{2},1_{H}) \\
&=&-\gamma _{1}B(x_{1}\otimes 1_{H};GX_{1},1_{H})+B(x_{1}\otimes
x_{1};X_{1},1_{H}).
\end{eqnarray*}%
In view of the form of the elements we get%
\begin{gather}
2B(x_{1}\otimes 1_{H};1_{A},1_{H})+  \label{X1,x1,X1F31,1H} \\
+\lambda \left[ -B\left( g\otimes 1_{H};1_{A},x_{2}\right) +B(x_{1}\otimes
1_{H};1_{A},x_{1}x_{2})\right]  \notag \\
+\gamma _{1}\left[ B(g\otimes 1_{H};G,1_{H})+B(x_{1}\otimes 1_{H};G,x_{1})%
\right] =0.  \notag
\end{gather}

\subparagraph{Case $f=x_{1}x_{2}$%
\protect\begin{eqnarray*}
&&B(x_{1}\otimes 1_{H};1_{A},x_{1}x_{2})+\protect\lambda B(x_{1}\otimes
1_{H};X_{1}X_{2},x_{1}x_{2}) \\
&=&\protect\gamma _{1}B(x_{1}g\otimes g;GX_{1},x_{1}x_{2})+B(x_{1}g\otimes
g;1_{A},x_{1}x_{2})+ \\
&&+B(x_{1}\otimes x_{1};X_{1},x_{1}x_{2}).
\protect\end{eqnarray*}%
}

By applying $\left( \ref{eq.10}\right) $ this rewrites as%
\begin{eqnarray*}
&&2B(x_{1}\otimes 1_{H};1_{A},x_{1}x_{2})+\lambda B(x_{1}\otimes
1_{H};X_{1}X_{2},x_{1}x_{2}) \\
&=&-\gamma _{1}B(x_{1}\otimes 1_{H};GX_{1},x_{1}x_{2})+B(x_{1}\otimes
x_{1};X_{1},x_{1}x_{2}).
\end{eqnarray*}%
In view of the form of the elements we get%
\begin{equation}
2B\left( x_{1}\otimes 1_{H};1_{A},x_{1}x_{2}\right) +\gamma _{1}B\left(
g\otimes 1_{H};G,x_{1}x_{2}\right) =0  \label{X1,x1,X1F31,x1x2}
\end{equation}

\subparagraph{Case $f=gx_{1}$}

\begin{eqnarray*}
&&B(x_{1}\otimes 1_{H};1_{A},gx_{1})+\lambda B(x_{1}\otimes
1_{H};X_{1}X_{2},gx_{1}) \\
&=&\gamma _{1}B(x_{1}g\otimes g;GX_{1},gx_{1})+B(x_{1}g\otimes
g;1_{A},gx_{1})+ \\
&&+B(x_{1}\otimes x_{1};X_{1},gx_{1}).
\end{eqnarray*}

By applying $\left( \ref{eq.10}\right) $ this rewrites as
\begin{eqnarray*}
&&B(x_{1}\otimes 1_{H};1_{A},gx_{1})+\lambda B(x_{1}\otimes
1_{H};X_{1}X_{2},gx_{1}) \\
&=&\gamma _{1}B(x_{1}\otimes 1_{H};GX_{1},gx_{1})+B(x_{1}\otimes
1_{H};1_{A},gx_{1})+ \\
&&+B(x_{1}\otimes x_{1};X_{1},gx_{1})
\end{eqnarray*}%
i.e.%
\begin{equation*}
\lambda B(x_{1}\otimes 1_{H};X_{1}X_{2},gx_{1})=-\gamma _{1}B(g\otimes
1_{H};GX_{1},g)+B(x_{1}\otimes x_{1};X_{1},gx_{1}).
\end{equation*}%
In view of the form of the elements we get%
\begin{equation*}
-\lambda B(g\otimes 1_{H};X_{1}X_{2},g)-\gamma _{1}B\left( g\otimes
1_{H};G,gx_{1}\right) -B(x_{1}\otimes x_{1};X_{1},gx_{1})=0
\end{equation*}%
\begin{gather}
-\lambda B\left( g\otimes 1_{H};1_{A},gx_{1}x_{2}\right) -\gamma _{1}B\left(
g\otimes 1_{H};G,gx_{1}\right)  \label{X1,x1,X1F31,gx1} \\
-2\left[ B(g\otimes 1_{H};1_{A},g)+B(x_{1}\otimes 1_{H};1_{A},gx_{1})\right]
=0  \notag
\end{gather}

\subparagraph{Case $f=gx_{2}$}

\begin{eqnarray*}
&&B(x_{1}\otimes 1_{H};1_{A},gx_{2})+\lambda B(x_{1}\otimes
1_{H};X_{1}X_{2},gx_{2}) \\
&=&\gamma _{1}B(x_{1}\otimes 1_{H};GX_{1},gx_{2})+B(x_{1}g\otimes
g;1_{A},gx_{2}) \\
&&+B(x_{1}\otimes x_{1};X_{1},gx_{2}).
\end{eqnarray*}

By applying $\left( \ref{eq.10}\right) $ this rewrites as
\begin{eqnarray*}
&&B(x_{1}\otimes 1_{H};1_{A},gx_{2})+\lambda B(x_{1}\otimes
1_{H};X_{1}X_{2},gx_{2})+ \\
&&-\gamma _{1}B\left( x_{1}\otimes 1_{H};GX_{1},gx_{2}\right) \\
&=&B(x_{1}\otimes 1_{H};1_{A},gx_{2})+B(x_{1}\otimes x_{1};X_{1},gx_{2})
\end{eqnarray*}%
i.e.%
\begin{equation*}
+\lambda B(x_{1}\otimes 1_{H};X_{1}X_{2},gx_{2})-\gamma _{1}B\left(
x_{1}\otimes 1_{H};GX_{1},gx_{2}\right) =+B(x_{1}\otimes x_{1};X_{1},gx_{2})
\end{equation*}%
In view of the form of the elements we get
\begin{equation}
\gamma _{1}\left[ -B\left( g\otimes 1_{H};G,gx_{2}\right) +B(x_{1}\otimes
1_{H};G,gx_{1}x_{2})\right] =0  \label{X1,x1,X1F31,gx2}
\end{equation}

\paragraph{Equality $\left( \protect\ref{X1F41}\right) $}

rewrites as

\begin{eqnarray*}
&&-\beta _{1}B(x_{1}\otimes 1_{H};X_{1}X_{2},f) \\
&=&-\beta _{1}B(gx_{1}\otimes g;X_{1}X_{2},f)+ \\
&&-\gamma _{1}B(gx_{1}\otimes g;GX_{2},f)+B(x_{1}\otimes x_{1};X_{2},f)
\end{eqnarray*}

\subparagraph{Case $f=1_{H}$}

\begin{eqnarray*}
&&-\beta _{1}B(x_{1}\otimes 1_{H};X_{1}X_{2},1_{H}) \\
&=&-\beta _{1}B(gx_{1}\otimes g;X_{1}X_{2},1_{H})-\gamma _{1}B(gx_{1}\otimes
g;GX_{2},1_{H})+ \\
&&+B(x_{1}\otimes x_{1};X_{2},1_{H})
\end{eqnarray*}%
By applying $\left( \ref{eq.10}\right) $ this rewrites as%
\begin{equation*}
0=-\gamma _{1}B(x_{1}\otimes 1_{H};GX_{2},1_{H})+B(x_{1}\otimes
x_{1};X_{2},1_{H}).
\end{equation*}%
In view of the form of the elements we get

\begin{equation}
0=\gamma _{1}B(x_{1}\otimes 1_{H};G,x_{2}).  \label{X1,x1,X1F41,1H}
\end{equation}

\subparagraph{Case $f=gx_{1}$}

\begin{eqnarray*}
&&2\beta _{1}B(gx_{1}\otimes g;X_{1}X_{2},gx_{1}) \\
&=&-\gamma _{1}B(gx_{1}\otimes g;GX_{2},gx_{1})+B(x_{1}\otimes
x_{1};X_{2},gx_{1}).
\end{eqnarray*}%
By applying $\left( \ref{eq.10}\right) $ this rewrites as%
\begin{eqnarray*}
&&-2\beta _{1}B(x_{1}\otimes 1_{H};X_{1}X_{2},gx_{1}) \\
&=&\gamma _{1}B(x_{1}\otimes 1_{H};GX_{2},gx_{1})+B(x_{1}\otimes
x_{1};X_{2},gx_{1}).
\end{eqnarray*}%
In view of the form of the elements we get%
\begin{equation}
2\beta _{1}B\left( g\otimes 1_{H};1_{A},gx_{1}x_{2}\right) -\gamma
_{1}B(x_{1}\otimes 1_{H};G,gx_{1}x_{2})=0.  \label{X1,x1,X1F41,gx1}
\end{equation}

\paragraph{Equality $\left( \protect\ref{X1F51}\right) $}

rewrites as%
\begin{eqnarray*}
&&-B(x_{1}\otimes 1_{H};X_{2},f) \\
&=&-B(gx_{1}\otimes g;X_{2},f)-\gamma _{1}B(gx_{1}\otimes g;GX_{1}X_{2},f)+
\\
&&+B(x_{1}\otimes x_{1};X_{1}X_{2},f).
\end{eqnarray*}

\subparagraph{Case $f=g$}

\begin{eqnarray*}
&&-B(x_{1}\otimes 1_{H};X_{2},g) \\
&=&-B(gx_{1}\otimes g;X_{2},g)-\gamma _{1}B(gx_{1}\otimes g;GX_{1}X_{2},g)+
\\
&&+B(x_{1}\otimes x_{1};X_{1}X_{2},g).
\end{eqnarray*}%
By applying $\left( \ref{eq.10}\right) $ this rewrites as%
\begin{equation*}
-\gamma _{1}B(x_{1}\otimes 1_{H};GX_{1}X_{2},g)+B(x_{1}\otimes
x_{1};X_{1}X_{2},g)=0.
\end{equation*}%
In view of the form of the elements we get%
\begin{equation*}
\gamma _{1}\left[ -B(g\otimes 1_{H};G,gx_{2})+B(x_{1}\otimes
1_{H};G,gx_{1}x_{2})\right] =0.
\end{equation*}%
which follows from $\left( \ref{X1,g,X1F21,gx1x2}\right) .$

\subparagraph{Case $f=gx_{1}$}

\begin{eqnarray*}
&&-B(x_{1}\otimes 1_{H};X_{2},gx_{1}) \\
&=&-B(gx_{1}\otimes g;X_{2},gx_{1})-\gamma _{1}B(gx_{1}\otimes
g;GX_{1}X_{2},gx_{1}) \\
&&+B(x_{1}\otimes x_{1};X_{1}X_{2},gx_{1}).
\end{eqnarray*}%
By applying $\left( \ref{eq.10}\right) $ this rewrites as%
\begin{eqnarray*}
&&-2B(x_{1}\otimes 1_{H};X_{2},gx_{1}) \\
&=&\gamma _{1}B(x_{1}\otimes 1_{H};GX_{1}X_{2},gx_{1})+B(x_{1}\otimes
x_{1};X_{1}X_{2},gx_{1}).
\end{eqnarray*}%
In view of the form of the elements we get nothing new.

\paragraph{Equality $\left( \protect\ref{X1F61}\right) $}

rewrites as%
\begin{eqnarray*}
&&B(x_{1}\otimes 1_{H};G,f)+\lambda B(x_{1}\otimes 1_{H};GX_{1}X_{2},f \\
&=&B(gx_{1}\otimes g;G,f)+B(x_{1}\otimes x_{1};GX_{1},f).
\end{eqnarray*}

\subparagraph{Case $f=g$}

\begin{eqnarray*}
&&B(x_{1}\otimes 1_{H};G,g)+\lambda B(x_{1}\otimes 1_{H};GX_{1}X_{2},g) \\
&=&B(gx_{1}\otimes g;G,g)+B(x_{1}\otimes x_{1};GX_{1},g).
\end{eqnarray*}%
By applying $\left( \ref{eq.10}\right) $ this rewrites as%
\begin{equation*}
\lambda B(x_{1}\otimes 1_{H};GX_{1}X_{2},g)=B(x_{1}\otimes x_{1};GX_{1},g).
\end{equation*}%
In view of the form of the elements we get%
\begin{equation*}
\lambda \left[ -B(g\otimes 1_{H};G,gx_{2})+B(x_{1}\otimes
1_{H};G,gx_{1}x_{2})\right] =0
\end{equation*}%
which follows from $\left( \ref{X1,g,X1F21,gx1x2}\right) .$

\subparagraph{Case $f=x_{1}$}

\begin{eqnarray*}
&&B(x_{1}\otimes 1_{H};G,x_{1})+\lambda B(x_{1}\otimes
1_{H};GX_{1}X_{2},x_{1}) \\
&=&B(gx_{1}\otimes g;G,x_{1})+B(x_{1}\otimes x_{1};GX_{1},x_{1}).
\end{eqnarray*}%
By applying $\left( \ref{eq.10}\right) $ this rewrites as%
\begin{equation*}
2B(x_{1}\otimes 1_{H};G,x_{1})+\lambda B(x_{1}\otimes
1_{H};GX_{1}X_{2},x_{1})=B(x_{1}\otimes x_{1};GX_{1},x_{1}).
\end{equation*}%
In view of the form of the elements we get%
\begin{equation*}
2B(x_{1}\otimes 1_{H};G,x_{1})-\lambda B(g\otimes 1_{H};G,x_{1}x_{2})=0
\end{equation*}%
which is $\left( \ref{X1,g,X1F21,x1}\right) .$

\subparagraph{Case $f=x_{2}$}

\begin{equation*}
B(x_{1}\otimes 1_{H};G,x_{2})+\lambda B(x_{1}\otimes
1_{H};GX_{1}X_{2},x_{2})=B(gx_{1}\otimes g;G,x_{2})+B(x_{1}\otimes
x_{1};GX_{1},x_{2}).
\end{equation*}%
By applying $\left( \ref{eq.10}\right) $ this rewrites as%
\begin{equation*}
2B(x_{1}\otimes 1_{H};G,x_{2})+\lambda B(x_{1}\otimes
1_{H};GX_{1}X_{2},x_{2})=B(x_{1}\otimes x_{1};GX_{1},x_{2}).
\end{equation*}%
In view of the form of the elements we get%
\begin{equation*}
B(x_{1}\otimes 1_{H};G,x_{2})=0
\end{equation*}%
which is $\left( \ref{X1,g,X1F21,x2}\right) .$

\subparagraph{Case $f=gx_{1}x_{2}$}

\begin{equation*}
B(x_{1}\otimes 1_{H};G,gx_{1}x_{2})+\lambda B(x_{1}\otimes
1_{H};GX_{1}X_{2},gx_{1}x_{2})=B(gx_{1}\otimes
g;G,gx_{1}x_{2})+B(x_{1}\otimes x_{1};GX_{1},gx_{1}x_{2}).
\end{equation*}%
By applying $\left( \ref{eq.10}\right) $ this rewrites as%
\begin{equation*}
\lambda B(x_{1}\otimes 1_{H};GX_{1}X_{2},gx_{1}x_{2})=B(x_{1}\otimes
x_{1};GX_{1},gx_{1}x_{2}).
\end{equation*}%
In view of the form of the elements we get

\begin{equation*}
0=-B\left( g\otimes 1_{H};G,gx_{2}\right) +B(x_{1}\otimes
1_{H};G,gx_{1}x_{2}).
\end{equation*}%
which is $\left( \ref{X1,g,X1F21,gx1x2}\right) .$

\paragraph{Equality $\left( \protect\ref{X1F71}\right) $}

rewrites as%
\begin{equation*}
-\beta _{1}B(x_{1}\otimes 1_{H};GX_{1}X_{2},f)=\beta _{1}B(gx_{1}\otimes
g;GX_{1}X_{2},f)+B(x_{1}\otimes x_{1};GX_{2},f)
\end{equation*}

\subparagraph{Case $f=g$}

\begin{equation*}
-\beta _{1}B(x_{1}\otimes 1_{H};GX_{1}X_{2},g)=\beta _{1}B(gx_{1}\otimes
g;GX_{1}X_{2},g)+B(x_{1}\otimes x_{1};GX_{2},g)
\end{equation*}%
By applying $\left( \ref{eq.10}\right) $ this rewrites as%
\begin{equation*}
-2\beta _{1}B(x_{1}\otimes 1_{H};GX_{1}X_{2},g)=B(x_{1}\otimes
x_{1};GX_{2},g).
\end{equation*}%
In view of the form of the elements we get%
\begin{equation*}
\beta _{1}\left[ -B(g\otimes 1_{H};G,gx_{2})+B(x_{1}\otimes
1_{H};G,gx_{1}x_{2})\right] =0.
\end{equation*}

which follows from $\left( \ref{X1,g,X1F21,gx1x2}\right) .$

\subparagraph{Case $f=x_{1}$}

\begin{equation*}
-\beta _{1}B(x_{1}\otimes 1_{H};GX_{1}X_{2},x_{1})=\beta _{1}B(gx_{1}\otimes
g;GX_{1}X_{2},x_{1})+B(x_{1}\otimes x_{1};GX_{2},x_{1})
\end{equation*}%
By applying $\left( \ref{eq.10}\right) $ this rewrites as%
\begin{equation*}
0=B(x_{1}\otimes x_{1};GX_{2},x_{1})
\end{equation*}%
which is already known.

\paragraph{Equality $\left( \protect\ref{X1F81}\right) $}

rewrites as%
\begin{equation*}
-B(x_{1}\otimes 1_{H};GX_{2},f)=B(gx_{1}\otimes g;GX_{2},f)+B(x_{1}\otimes
x_{1};GX_{1}X_{2},f).
\end{equation*}

\subparagraph{Case $f=1_{H}$}

\begin{equation*}
-B(x_{1}\otimes 1_{H};GX_{2},1_{H})=B(gx_{1}\otimes
g;GX_{2},1_{H})+B(x_{1}\otimes x_{1};GX_{1}X_{2},1_{H}).
\end{equation*}%
By applying $\left( \ref{eq.10}\right) $ this rewrites as%
\begin{equation*}
-2B(x_{1}\otimes 1_{H};GX_{2},1_{H})=B(x_{1}\otimes x_{1};GX_{1}X_{2},1_{H}).
\end{equation*}%
In view of the form of the elements we get%
\begin{equation*}
B(x_{1}\otimes 1_{H};G,x_{2})=0
\end{equation*}

which is $\left( \ref{X1,g,X1F21,x2}\right) .$

\subparagraph{Case $f=gx_{1}$}

\begin{equation*}
-B(x_{1}\otimes 1_{H};GX_{2},gx_{1})=B(gx_{1}\otimes
g;GX_{2},gx_{1})+B(x_{1}\otimes x_{1};GX_{1}X_{2},gx_{1}).
\end{equation*}%
By applying $\left( \ref{eq.10}\right) $ this rewrites as%
\begin{equation*}
0=B(g\otimes 1_{H};G,gx_{2})-B(x_{1}\otimes 1_{H};G,gx_{1}x_{2})
\end{equation*}%
which is $\left( \ref{X1,g,X1F21,gx1x2}\right) .$

\subsubsection{Case $x_{2}\otimes 1_{H}$}

\paragraph{Equality $\left( \protect\ref{X1F11}\right) $}

rewrites as%
\begin{eqnarray*}
&&\beta _{1}B(x_{2}\otimes 1_{H};X_{1},f)+\lambda B(x_{2}\otimes
1_{H};X_{2},f) \\
&=&-\gamma _{1}B(gx_{2}\otimes g;G,f)-\beta _{1}B(gx_{2}\otimes g;X_{1},f)+
\\
&&-B(x_{1}x_{2}\otimes g;1_{A},f)+B(x_{2}\otimes x_{1};1_{A},f).
\end{eqnarray*}

\subparagraph{Case $f=g$}

\begin{eqnarray*}
&&\beta _{1}B(x_{2}\otimes 1_{H};X_{1},g)+\lambda B(x_{2}\otimes
1_{H};X_{2},g) \\
&=&-\gamma _{1}B(gx_{2}\otimes g;G,g)-\beta _{1}B(gx_{2}\otimes g;X_{1},g)+
\\
&&-B(x_{1}x_{2}\otimes g;1_{A},g)+B(x_{2}\otimes x_{1};1_{A},g).
\end{eqnarray*}%
By applying $\left( \ref{eq.10}\right) $ this rewrites as%
\begin{eqnarray*}
&&\beta _{1}B(x_{2}\otimes 1_{H};X_{1},g)+\lambda B(x_{2}\otimes
1_{H};X_{2},g) \\
&=&-\gamma _{1}B(x_{2}\otimes 1_{H};G,g)-\beta _{1}B(x_{2}\otimes
1_{H};X_{1},g)+ \\
&&-B(gx_{1}x_{2}\otimes 1_{H};1_{A},g)+B(x_{2}\otimes x_{1};1_{A},g).
\end{eqnarray*}%
In view of the form of the elements we get

\begin{eqnarray}
&&-2\beta _{1}B(x_{2}\otimes 1_{H};1_{A},gx_{1})+\lambda \left[ -B(g\otimes
1_{H};1_{A},g)-B(x_{2}\otimes \ 1_{H};1_{A},gx_{2})\right]
\label{X1,x2,X1F11,g} \\
&=&-\gamma _{1}B(x_{2}\otimes 1_{H};G,g)-2B(gx_{1}x_{2}\otimes
1_{H};1_{A},g).  \notag
\end{eqnarray}

\subparagraph{Case $f=x_{1}$}

\begin{eqnarray*}
&&\beta _{1}B(x_{2}\otimes 1_{H};X_{1},x_{1})+\lambda B(x_{2}\otimes
1_{H};X_{2},x_{1}) \\
&=&-\gamma _{1}B(gx_{2}\otimes g;G,x_{1})-\beta _{1}B(gx_{2}\otimes
g;X_{1},x_{1})+ \\
&&-B(x_{1}x_{2}\otimes g;1_{A},x_{1})+B(x_{2}\otimes x_{1};1_{A},x_{1}).
\end{eqnarray*}%
By applying $\left( \ref{eq.10}\right) $ this rewrites as%
\begin{equation*}
\lambda B(x_{2}\otimes 1_{H};X_{2},x_{1})=\gamma _{1}B(x_{2}\otimes
1_{H};G,x_{1})+B(gx_{1}x_{2}\otimes 1_{H};1_{A},x_{1})+B(x_{2}\otimes
x_{1};1_{A},x_{1}).
\end{equation*}%
In view of the form of the elements we get

\begin{equation}
\lambda \left[ -B(g\otimes 1_{H};1_{A},x_{1})-B(x_{2}\otimes \
1_{H};1_{A},x_{1}x_{2})\right] =\gamma _{1}B(x_{2}\otimes
1_{H};G,x_{1})+2B(gx_{1}x_{2}\otimes 1_{H};1_{A},x_{1}).
\label{X1,x2,X1F11,x1}
\end{equation}

\subparagraph{Case $f=x_{2}$}

\begin{eqnarray*}
\beta _{1}B(x_{2}\otimes 1_{H};X_{1},x_{2})+\lambda B(x_{2}\otimes
1_{H};X_{2},x_{2}) &=&-\gamma _{1}B(gx_{2}\otimes g;G,x_{2})-\beta
_{1}B(gx_{2}\otimes g;X_{1},x_{2}) \\
&&-B(x_{1}x_{2}\otimes g;1_{A},x_{2})+B(x_{2}\otimes x_{1};1_{A},x_{2}).
\end{eqnarray*}%
By applying $\left( \ref{eq.10}\right) $ this rewrites as%
\begin{equation*}
\lambda B(x_{2}\otimes 1_{H};X_{2},x_{2})=\gamma _{1}B(x_{2}\otimes
1_{H};G,x_{2})+B(gx_{1}x_{2}\otimes 1_{H};1_{A},x_{2})+B(x_{2}\otimes
x_{1};1_{A},x_{2}).
\end{equation*}%
In view of the form of the elements we get%
\begin{equation}
-\lambda B(g\otimes 1_{H};1_{A},x_{2})=\gamma _{1}B(x_{2}\otimes
1_{H};G,x_{2})+2B(gx_{1}x_{2}\otimes 1_{H};1_{A},x_{2}).
\label{X1,x2,X1F11,x2}
\end{equation}

\subparagraph{Case $f=gx_{1}x_{2}$}

\begin{eqnarray*}
&&\beta _{1}B(x_{2}\otimes 1_{H};X_{1},gx_{1}x_{2})+\lambda B(x_{2}\otimes
1_{H};X_{2},gx_{1}x_{2}) \\
&=&-\gamma _{1}B(gx_{2}\otimes g;G,gx_{1}x_{2})-\beta _{1}B(gx_{2}\otimes
g;X_{1},gx_{1}x_{2})+ \\
&&-B(x_{1}x_{2}\otimes g;1_{A},gx_{1}x_{2})+B(x_{2}\otimes
x_{1};1_{A},gx_{1}x_{2}).
\end{eqnarray*}%
By applying $\left( \ref{eq.10}\right) $ this rewrites as%
\begin{eqnarray*}
&&2\beta _{1}B(x_{2}\otimes 1_{H};X_{1},gx_{1}x_{2})+\lambda B(x_{2}\otimes
1_{H};X_{2},gx_{1}x_{2}) \\
&=&-\gamma _{1}B(x_{2}\otimes 1_{H};G,gx_{1}x_{2})-B(gx_{1}x_{2}\otimes
1_{H};1_{A},gx_{1}x_{2})+B(x_{2}\otimes x_{1};1_{A},gx_{1}x_{2}).
\end{eqnarray*}%
In view of the form of the elements we get

\begin{gather}
\lambda B\left( g\otimes 1_{H};1_{A},gx_{1}x_{2}\right) -\gamma
_{1}B(x_{2}\otimes 1_{H};G,gx_{1}x_{2})  \label{X1,x2,X1F11,gx1x2} \\
-2B\left( x_{2}\otimes 1_{H};1_{A},gx_{2}\right) -2B(gx_{1}x_{2}\otimes
1_{H};1_{A},gx_{1}x_{2})=0.  \notag
\end{gather}

\paragraph{Equality $\left( \protect\ref{X1F21}\right) $}

rewrites as%
\begin{eqnarray*}
&&\beta _{1}B(x_{2}\otimes 1_{H};GX_{1},f)+\lambda B(x_{2}\otimes
1_{H};GX_{2},f) \\
&=&\beta _{1}B(gx_{1}\otimes g;GX_{1},f)-B(x_{1}x_{2}\otimes
g;G,f)+B(x_{2}\otimes x_{1};G,f).
\end{eqnarray*}

\subparagraph{Case $f=1_{H}$}

\begin{eqnarray*}
&&\beta _{1}B(x_{2}\otimes 1_{H};GX_{1},1_{H})+\lambda B(x_{2}\otimes
1_{H};GX_{2},1_{H}) \\
&=&\beta _{1}B(gx_{1}\otimes g;GX_{1},1_{H})-B(x_{1}x_{2}\otimes
g;G,1_{H})+B(x_{2}\otimes x_{1};G,1_{H}).
\end{eqnarray*}%
By applying $\left( \ref{eq.10}\right) $ this rewrites as%
\begin{equation*}
\lambda B(x_{2}\otimes 1_{H};GX_{2},1_{H})=-B(gx_{1}x_{2}\otimes
1_{H};G,1_{H})+B(x_{2}\otimes x_{1};G,1_{H}).
\end{equation*}%
In view of the form of the elements we get%
\begin{equation}
\lambda \left[ B(g\otimes 1_{H};G,1_{H})+B(x_{2}\otimes \ 1_{H};G,x_{2}%
\right] =-2B(gx_{1}x_{2}\otimes 1_{H};G,1_{H}).  \label{X1,x2,X1F21,1H}
\end{equation}

\subparagraph{Case $f=x_{1}x_{2}$}

\begin{eqnarray*}
&&\beta _{1}B(x_{2}\otimes 1_{H};GX_{1},x_{1}x_{2})+\lambda B(x_{2}\otimes
1_{H};GX_{2},x_{1}x_{2}) \\
&=&\beta _{1}B(gx_{1}\otimes g;GX_{1},x_{1}x_{2})-B(x_{1}x_{2}\otimes
g;G,x_{1}x_{2})+ \\
&&+B(x_{2}\otimes x_{1};G,x_{1}x_{2}).
\end{eqnarray*}%
By applying $\left( \ref{eq.10}\right) $ this rewrites as%
\begin{equation*}
\lambda B(x_{2}\otimes 1_{H};GX_{2},x_{1}x_{2})=-B(gx_{1}x_{2}\otimes
1_{H};G,x_{1}x_{2})+B(x_{2}\otimes x_{1};G,x_{1}x_{2}).
\end{equation*}%
In view of the form of the elements we get

\begin{equation}
\lambda B\left( g\otimes 1_{H};G,x_{1}x_{2}\right) +2B(gx_{1}x_{2}\otimes
1_{H};G,x_{1}x_{2})=0.  \label{X1,x2,X1F21,x1x2}
\end{equation}

\subparagraph{Case $f=gx_{1}$}

\begin{eqnarray*}
&&\beta _{1}B(x_{2}\otimes 1_{H};GX_{1},gx_{1})+\lambda B(x_{2}\otimes
1_{H};GX_{2},gx_{1}) \\
&=&\beta _{1}B(gx_{1}\otimes g;GX_{1},gx_{1})-B(x_{1}x_{2}\otimes
g;G,gx_{1})+B(x_{2}\otimes x_{1};G,gx_{1}).
\end{eqnarray*}%
By applying $\left( \ref{eq.10}\right) $ this rewrites as%
\begin{equation*}
2\beta _{1}B(x_{2}\otimes 1_{H};GX_{1},gx_{1})+\lambda B(x_{2}\otimes
1_{H};GX_{2},gx_{1})=+B(gx_{1}x_{2}\otimes 1_{H};G,gx_{1})+B(x_{2}\otimes
x_{1};G,gx_{1}).
\end{equation*}%
In view of the form of the elements we get%
\begin{equation}
\lambda \left[ B(g\otimes 1_{H};G,gx_{1})+B(x_{2}\otimes \
1_{H};G,gx_{1}x_{2})\right] =2B(gx_{1}x_{2}\otimes 1_{H};G,gx_{1})-2B\left(
x_{2}\otimes 1_{H};G,g\right) .  \label{X1,x2,X1F21,gx1}
\end{equation}

\subparagraph{Case $f=gx_{2}$}

\begin{eqnarray*}
&&\beta _{1}B(x_{2}\otimes 1_{H};GX_{1},gx_{2})+\lambda B(x_{2}\otimes
1_{H};GX_{2},gx_{2}) \\
&=&\beta _{1}B(gx_{1}\otimes g;GX_{1},gx_{2})-B(x_{1}x_{2}\otimes
g;G,gx_{2})+B(x_{2}\otimes x_{1};G,gx_{2}).
\end{eqnarray*}%
By applying $\left( \ref{eq.10}\right) $ this rewrites as%
\begin{eqnarray*}
&&2\beta _{1}B(x_{2}\otimes 1_{H};GX_{1},gx_{2})+\lambda B(x_{2}\otimes
1_{H};GX_{2},gx_{2}) \\
&=&B(gx_{1}x_{2}\otimes 1_{H};G,gx_{2})+B(x_{2}\otimes x_{1};G,gx_{2}).
\end{eqnarray*}%
In view of the form of the elements we get%
\begin{equation}
-2\beta _{1}B(x_{2}\otimes 1_{H};G,gx_{1}x_{2})+\lambda B(g\otimes
1_{H};G,gx_{2})-2B(gx_{1}x_{2}\otimes 1_{H};G,gx_{2})=0.
\label{X1,x2,X1F21,gx2}
\end{equation}

\paragraph{Equality $\left( \protect\ref{X1F31}\right) $}

rewrites as%
\begin{eqnarray*}
&&B(x_{2}\otimes 1_{H};1_{A},f)+\lambda B(x_{2}\otimes 1_{H};X_{1}X_{2},f) \\
&=&-\gamma _{1}B(gx_{2}\otimes g;GX_{1},f)-B(gx_{2}\otimes g;1_{A},f)+ \\
&&-B(x_{1}x_{2}\otimes g;X_{1},f)+B(x_{2}\otimes x_{1};X_{1},f).
\end{eqnarray*}

\subparagraph{Case $f=1_{H}$}

\begin{eqnarray*}
&&B(x_{2}\otimes 1_{H};1_{A},1_{H})+\lambda B(x_{2}\otimes
1_{H};X_{1}X_{2},1_{H}) \\
&=&-\gamma _{1}B(gx_{2}\otimes g;GX_{1},g)-B(gx_{2}\otimes
g;1_{A},1_{H})-B(x_{1}x_{2}\otimes g;X_{1},1_{H})+B(x_{2}\otimes
x_{1};X_{1},1_{H}).
\end{eqnarray*}%
By applying $\left( \ref{eq.10}\right) $ this rewrites as%
\begin{eqnarray*}
&&2B(x_{2}\otimes 1_{H};1_{A},1_{H})+\lambda B(x_{2}\otimes
1_{H};X_{1}X_{2},1_{H}) \\
&=&-\gamma _{1}B(x_{2}\otimes 1_{H};GX_{1},1_{H})-B(gx_{1}x_{2}\otimes
1_{H};X_{1},1_{H})+B(x_{2}\otimes x_{1};X_{1},1_{H}).
\end{eqnarray*}%
In view of the form of the elements we get

\begin{eqnarray*}
&&\lambda \left[ B(g\otimes 1_{H};1_{A},x_{1})+B(x_{2}\otimes \
1_{H};1_{A},x_{1}x_{2})\right] \\
&=&-\gamma _{1}B(x_{2}\otimes 1_{H};G,x_{1})-2B(gx_{1}x_{2}\otimes
1_{H};1_{A},x_{1})
\end{eqnarray*}%
which is $\left( \ref{X1,x2,X1F11,x1}\right) $.

\subparagraph{Case $f=x_{1}x_{2}$}

\begin{eqnarray*}
&&B(x_{2}\otimes 1_{H};1_{A},x_{1}x_{2})+\lambda B(x_{2}\otimes
1_{H};X_{1}X_{2},x_{1}x_{2}) \\
&=&-\gamma _{1}B(gx_{2}\otimes g;GX_{1},x_{1}x_{2})-B(gx_{2}\otimes
g;1_{A},x_{1}x_{2})+ \\
&&-B(x_{1}x_{2}\otimes g;X_{1},x_{1}x_{2})+B(x_{2}\otimes
x_{1};X_{1},x_{1}x_{2}).
\end{eqnarray*}%
By applying $\left( \ref{eq.10}\right) $ this rewrites as%
\begin{eqnarray*}
&&2B(x_{2}\otimes 1_{H};1_{A},x_{1}x_{2})+\lambda B(x_{2}\otimes
1_{H};X_{1}X_{2},x_{1}x_{2}) \\
&=&-\gamma _{1}B(x_{2}\otimes 1_{H};GX_{1},x_{1}x_{2})-B(gx_{1}x_{2}\otimes
1_{H};X_{1},x_{1}x_{2})+B(x_{2}\otimes x_{1};X_{1},x_{1}x_{2}).
\end{eqnarray*}%
In view of the form of the elements we get

\begin{equation*}
2B(x_{2}\otimes 1_{H};1_{A},x_{1}x_{2})=2B(x_{2}\otimes
1_{H};1_{A},x_{1}x_{2})
\end{equation*}%
which is trivial.

\subparagraph{Case $f=gx_{1}$}

\begin{eqnarray*}
&&B(x_{2}\otimes 1_{H};1_{A},gx_{1})+\lambda B(x_{2}\otimes
1_{H};X_{1}X_{2},gx_{1}) \\
&=&-\gamma _{1}B(gx_{2}\otimes g;GX_{1},gx_{1})-B(gx_{2}\otimes
g;1_{A},gx_{1})+ \\
&&-B(x_{1}x_{2}\otimes g;X_{1},gx_{1})+B(x_{2}\otimes x_{1};X_{1},gx_{1}).
\end{eqnarray*}%
By applying $\left( \ref{eq.10}\right) $ this rewrites as%
\begin{eqnarray*}
&&\lambda B(x_{2}\otimes 1_{H};X_{1}X_{2},gx_{1}) \\
&=&\gamma _{1}B(x_{2}\otimes 1_{H};GX_{1},gx_{1})+B(gx_{1}x_{2}\otimes
1_{H};X_{1},gx_{1})+B(x_{2}\otimes x_{1};X_{1},gx_{1}).
\end{eqnarray*}%
In view of the form of the elements we get%
\begin{equation*}
0=-B(x_{2}\otimes 1_{H};1_{A},gx_{1})+B(x_{2}\otimes 1_{H};1_{A},gx_{1})
\end{equation*}%
which is trivial.

\subparagraph{Case $f=gx_{2}$}

\begin{eqnarray*}
&&B(x_{2}\otimes 1_{H};1_{A},gx_{2})+\lambda B(x_{2}\otimes
1_{H};X_{1}X_{2},gx_{2}) \\
&=&-\gamma _{1}B(gx_{2}\otimes g;GX_{1},gx_{2})-B(gx_{2}\otimes
g;1_{A},gx_{2})+ \\
&&-B(x_{1}x_{2}\otimes g;X_{1},gx_{2})+B(x_{2}\otimes x_{1};X_{1},gx_{2}).
\end{eqnarray*}%
By applying $\left( \ref{eq.10}\right) $ this rewrites as%
\begin{equation*}
\lambda B(x_{2}\otimes 1_{H};X_{1}X_{2},gx_{2})=\gamma _{1}B(x_{2}\otimes
1_{H};GX_{1},gx_{2})+B(gx_{1}x_{2}\otimes 1_{H};X_{1},gx_{2})+B(x_{2}\otimes
x_{1};X_{1},gx_{2}).
\end{equation*}%
In view of the form of the elements we get%
\begin{eqnarray}
&&-\lambda B\left( g\otimes 1_{H};1_{A},gx_{1}x_{2}\right)
\label{X1,x2,X1F31,gx2} \\
&=&-\gamma _{1}B(x_{2}\otimes 1_{H};G,gx_{1}x_{2})-2B(x_{2}\otimes
1_{H};1_{A},gx_{2})-2B(gx_{1}x_{2}\otimes 1_{H};1_{A},gx_{1}x_{2}).  \notag
\end{eqnarray}

\paragraph{Equality $\left( \protect\ref{X1F41}\right) $}

rewrites as%
\begin{eqnarray*}
&&-\beta _{1}B(x_{2}\otimes 1_{H};X_{1}X_{2},f) \\
&=&-\beta _{1}B(gx_{2}\otimes g;X_{1}X_{2},f)-\gamma _{1}B(gx_{2}\otimes
g;GX_{2},f)+ \\
&&-B(x_{1}x_{2}\otimes g;X_{2},f)+B(x_{2}\otimes x_{1};X_{2},f).
\end{eqnarray*}

\subparagraph{Case $f=1_{H}$}

\begin{equation*}
0=-\gamma _{1}B(gx_{2}\otimes g;GX_{2},1_{H})-B(x_{1}x_{2}\otimes
g;X_{2},1_{H})+B(x_{2}\otimes x_{1};X_{2},1_{H}).
\end{equation*}%
By applying $\left( \ref{eq.10}\right) $ this rewrites as%
\begin{equation*}
0=-\gamma _{1}B(x_{2}\otimes 1_{H};GX_{2},1_{H})-B(gx_{1}x_{2}\otimes
1_{H};X_{2},1_{H})+B(x_{2}\otimes x_{1};X_{2},1_{H}).
\end{equation*}%
In view of the form of the elements we get%
\begin{eqnarray}
0 &=&-\gamma _{1}\left[ B(g\otimes 1_{H};G,1_{H})+B(x_{2}\otimes \
1_{H};G,x_{2}\right] +  \label{X1,x2,X1F41,1H} \\
&&-2B(x_{1}\otimes 1_{H};1_{A},1_{H})-2B(gx_{1}x_{2}\otimes
1_{H};1_{A},x_{2}).  \notag
\end{eqnarray}

\subparagraph{Case $f=x_{1}x_{2}$%
\protect\begin{gather*}
-\protect\beta _{1}B(x_{2}\otimes 1_{H};X_{1}X_{2},x_{1}x_{2}) \\
=-\protect\gamma _{1}B(gx_{2}\otimes g;GX_{2},x_{1}x_{2})-\protect\beta %
_{1}B(gx_{2}\otimes g;X_{1}X_{2},x_{1}x_{2}) \\
-B(x_{1}x_{2}\otimes g;X_{2},x_{1}x_{2})+B(x_{2}\otimes
x_{1};X_{2},x_{1}x_{2}).
\protect\end{gather*}%
}

By applying $\left( \ref{eq.10}\right) $ this rewrites as%
\begin{equation*}
0=-\gamma _{1}B(x_{2}\otimes 1_{H};GX_{2},x_{1}x_{2})-B(gx_{1}x_{2}\otimes
1_{H};X_{2},x_{1}x_{2})+B(x_{2}\otimes x_{1};X_{2},x_{1}x_{2}).
\end{equation*}%
In view of the form of the elements we get%
\begin{equation}
0=-\gamma _{1}B\left( g\otimes 1_{H};G,x_{1}x_{2}\right) -2B(x_{1}\otimes
1_{H};1_{A},x_{1}x_{2}).  \label{X1,x2,X1F41,x1x2}
\end{equation}

\subparagraph{Case $f=gx_{1}$}

\begin{equation*}
2\beta _{1}B(gx_{2}\otimes g;X_{1}X_{2},gx_{1})=-B(x_{1}x_{2}\otimes
g;X_{2},gx_{1})+B(x_{2}\otimes x_{1};X_{2},gx_{1}).
\end{equation*}%
By applying $\left( \ref{eq.10}\right) $ this rewrites as%
\begin{equation*}
-2\beta _{1}B(x_{2}\otimes 1_{H};X_{1}X_{2},gx_{1})=B(gx_{1}x_{2}\otimes
1_{H};X_{2},gx_{1})+B(x_{2}\otimes x_{1};X_{2},gx_{1}).
\end{equation*}%
In view of the form of the elements we get%
\begin{gather}
B(x_{1}\otimes 1_{H};1_{A},gx_{1})+B(gx_{1}x_{2}\otimes
1_{H};1_{A},gx_{1}x_{2})  \label{X1,x2,X1F41,gx1} \\
+B(g\otimes 1_{H};1_{A},g)+B(x_{2}\otimes \ 1_{H};1_{A},gx_{2})=0.  \notag
\end{gather}

\subparagraph{Case $f=gx_{2}$}

\begin{eqnarray*}
&&2\beta _{1}B(gx_{2}\otimes g;X_{1}X_{2},gx_{2}) \\
&=&-\gamma _{1}B(gx_{2}\otimes g;GX_{2},gx_{2})-B(x_{1}x_{2}\otimes
g;X_{2},gx_{2})+B(x_{2}\otimes x_{1};X_{2},gx_{2}).
\end{eqnarray*}%
By applying $\left( \ref{eq.10}\right) $ this rewrites as%
\begin{eqnarray*}
&&-2\beta _{1}B(x_{2}\otimes 1_{H};X_{1}X_{2},gx_{2}) \\
&=&\gamma _{1}B(x_{2}\otimes 1_{H};GX_{2},gx_{2})+B(gx_{1}x_{2}\otimes
1_{H};X_{2},gx_{2})+B(x_{2}\otimes x_{1};X_{2},gx_{2}).
\end{eqnarray*}%
In view of the form of the elements we get%
\begin{equation}
2\beta _{1}B\left( g\otimes 1_{H};1_{A},gx_{1}x_{2}\right) =\gamma
_{1}B(g\otimes 1_{H};G,gx_{2})+2B(x_{1}\otimes 1_{H};1_{A},gx_{2}).
\label{X1,x2,X1F41,gx2}
\end{equation}

\paragraph{Equality $\left( \protect\ref{X1F51}\right) $}

rewrites as%
\begin{eqnarray*}
-B(x_{2}\otimes 1_{H};X_{2},f) &=&-\gamma _{1}B(gx_{2}\otimes
g;GX_{1}X_{2},f)-B(gx_{2}\otimes g;X_{2},f)+ \\
&&-B(x_{1}x_{2}\otimes g;X_{1}X_{2},f)+B(x_{2}\otimes x_{1};X_{1}X_{2},f).
\end{eqnarray*}

\subparagraph{Case $f=g$}

\begin{eqnarray*}
&&-B(x_{2}\otimes 1_{H};X_{2},g) \\
&=&-\gamma _{1}B(gx_{2}\otimes g;GX_{1}X_{2},g)-B(gx_{2}\otimes g;X_{2},g)+
\\
&&-B(x_{1}x_{2}\otimes g;X_{1}X_{2},g)+B(x_{2}\otimes x_{1};X_{1}X_{2},g).
\end{eqnarray*}%
By applying $\left( \ref{eq.10}\right) $ this rewrites as

\begin{equation*}
0=-\gamma _{1}B(x_{2}\otimes 1_{H};GX_{1}X_{2},g)-B(gx_{1}x_{2}\otimes
1_{H};X_{1}X_{2},g)+B(x_{2}\otimes x_{1};X_{1}X_{2},g).
\end{equation*}%
In view of the form of the elements we get

\begin{eqnarray}
0 &=&\gamma _{1}\left[ B(g\otimes 1_{H};G,gx_{1})+B(x_{2}\otimes \
1_{H};G,gx_{1}x_{2})\right]  \label{X1,x2,X1F51,g} \\
&&+\left[
\begin{array}{c}
2B(g\otimes 1_{H};1_{A},g)+2B(x_{2}\otimes \ 1_{H};1_{A},gx_{2}) \\
+2B(x_{1}\otimes 1_{H};1_{A},gx_{1})+2B(gx_{1}x_{2}\otimes
1_{H};1_{A},gx_{1}x_{2})%
\end{array}%
\right]  \notag
\end{eqnarray}

\subparagraph{Case $f=x_{1}$}

\begin{eqnarray*}
&&-B(x_{2}\otimes 1_{H};X_{2},x_{1}) \\
&=&-\gamma _{1}B(gx_{2}\otimes g;GX_{1}X_{2},x_{1})-B(gx_{2}\otimes
g;X_{2},x_{1})+ \\
&&-B(x_{1}x_{2}\otimes g;X_{1}X_{2},x_{1})+B(x_{2}\otimes
x_{1};X_{1}X_{2},x_{1}).
\end{eqnarray*}%
By applying $\left( \ref{eq.10}\right) $ this rewrites as%
\begin{eqnarray*}
&&-2B(x_{2}\otimes 1_{H};X_{2},x_{1}) \\
&=&\gamma _{1}B(x_{2}\otimes 1_{H};GX_{1}X_{2},x_{1})+B(gx_{1}x_{2}\otimes
1_{H};X_{1}X_{2},x_{1})+B(x_{2}\otimes x_{1};X_{1}X_{2},x_{1}).
\end{eqnarray*}%
In view of the form of the elements we get%
\begin{equation*}
B(g\otimes 1_{H};1_{A},x_{1})+B(x_{2}\otimes \
1_{H};1_{A},x_{1}x_{2})=B(g\otimes 1_{H};1_{A},x_{1})+B(x_{2}\otimes \
1_{H};1_{A},x_{1}x_{2})
\end{equation*}%
which is trivial.

\subparagraph{Case $f=x_{2}$}

\begin{eqnarray*}
&&-B(x_{2}\otimes 1_{H};X_{2},x_{2}) \\
&=&-\gamma _{1}B(gx_{2}\otimes g;GX_{1}X_{2},x_{2})-B(gx_{2}\otimes
g;X_{2},x_{2})+ \\
&&-B(x_{1}x_{2}\otimes g;X_{1}X_{2},x_{2})+B(x_{2}\otimes
x_{1};X_{1}X_{2},x_{2}).
\end{eqnarray*}%
By applying $\left( \ref{eq.10}\right) $ this rewrites as%
\begin{eqnarray*}
&&-2B(x_{2}\otimes 1_{H};X_{2},x_{2}) \\
&=&\gamma _{1}B(x_{2}\otimes 1_{H};GX_{1}X_{2},x_{2})+B(gx_{1}x_{2}\otimes
1_{H};X_{1}X_{2},x_{2})+B(x_{2}\otimes x_{1};X_{1}X_{2},x_{2}).
\end{eqnarray*}%
In view of the form of the elements we get%
\begin{equation*}
B(g\otimes 1_{H};1_{A},x_{2})=-\gamma _{1}B\left( g\otimes
1_{H};G,x_{1}x_{2}\right) -B(x_{1}\otimes
1_{H};1_{A},x_{1}x_{2})+B(x_{2}\otimes x_{1};X_{1}X_{2},x_{2}).
\end{equation*}%
\begin{equation}
0=\gamma _{1}B\left( g\otimes 1_{H};G,x_{1}x_{2}\right) +2B(x_{1}\otimes
1_{H};1_{A},x_{1}x_{2}).  \label{X1,x2,X1F51,x2}
\end{equation}

\subparagraph{Case $f=gx_{1}x_{2}$}

\begin{eqnarray*}
&&-B(x_{2}\otimes 1_{H};X_{2},gx_{1}x_{2}) \\
&=&-\gamma _{1}B(gx_{2}\otimes g;GX_{1}X_{2},gx_{1}x_{2})-B(gx_{2}\otimes
g;X_{2},gx_{1}x_{2}) \\
&&-B(x_{1}x_{2}\otimes g;X_{1}X_{2},gx_{1}x_{2})+B(x_{2}\otimes
x_{1};X_{1}X,gx_{1}x_{2}).
\end{eqnarray*}%
By applying $\left( \ref{eq.10}\right) $ this rewrites as

\begin{eqnarray*}
0 &=&-\gamma _{1}B(x_{2}\otimes 1_{H};GX_{1}X_{2},gx_{1}x_{2})+ \\
&&-B(gx_{1}x_{2}\otimes 1_{H};X_{1}X_{2},gx_{1}x_{2})+B(x_{2}\otimes
x_{1};X_{1}X_{2},gx_{1}x_{2}).
\end{eqnarray*}%
In view of the form of the elements we get%
\begin{equation*}
0=-B\left( g\otimes 1_{H};1_{A},gx_{1}x_{2}\right) +B\left( g\otimes
1_{H};1_{A},gx_{1}x_{2}\right) .
\end{equation*}%
which is trivial.

\paragraph{Equality $\left( \protect\ref{X1F61}\right) $}

rewrites as%
\begin{eqnarray*}
&&B(x_{2}\otimes 1_{H};G,f)+\lambda B(x_{2}\otimes 1_{H};GX_{1}X_{2},f) \\
&=&B(gx_{2}\otimes g;G,f)-B(x_{1}x_{2}\otimes g;GX_{1},f)+B(x_{2}\otimes
x_{1};GX_{1},f).
\end{eqnarray*}

\subparagraph{Case $f=g$%
\protect\begin{eqnarray*}
&&B(x_{2}\otimes 1_{H};G,g)+\protect\lambda B(x_{2}\otimes
1_{H};GX_{1}X_{2},g) \\
&=&B(gx_{2}\otimes g;G,g)-B(x_{1}x_{2}\otimes g;GX_{1},g)+B(x_{2}\otimes
x_{1};GX_{1},g).
\protect\end{eqnarray*}%
}

By applying $\left( \ref{eq.10}\right) $ this rewrites as%
\begin{equation*}
\lambda B(x_{2}\otimes 1_{H};GX_{1}X_{2},g)=-B(gx_{1}x_{2}\otimes
1_{H};GX_{1},g)+B(x_{2}\otimes x_{1};GX_{1},g).
\end{equation*}

In view of the form of the elements we get%
\begin{eqnarray}
&&\lambda \left[ B(g\otimes 1_{H};G,gx_{1})+B(x_{2}\otimes \
1_{H};G,gx_{1}x_{2})\right]  \label{X1,x2,X1F61,g} \\
&=&-2\left[ B(x_{2}\otimes 1_{H};G,g)-B(gx_{1}x_{2}\otimes 1_{H};G,gx_{1})%
\right] .  \notag
\end{eqnarray}

\subparagraph{Case $f=x_{1}$}

\begin{eqnarray*}
&&B(x_{2}\otimes 1_{H};G,x_{1})+\lambda B(x_{2}\otimes
1_{H};GX_{1}X_{2},x_{1}) \\
&=&B(gx_{2}\otimes g;G,x_{1})-B(x_{1}x_{2}\otimes
g;GX_{1},x_{1})+B(x_{2}\otimes x_{1};GX_{1},x_{1}).
\end{eqnarray*}%
By applying $\left( \ref{eq.10}\right) $ this rewrites as%
\begin{eqnarray*}
&&2B(x_{2}\otimes 1_{H};G,x_{1})+\lambda B(x_{2}\otimes
1_{H};GX_{1}X_{2},x_{1}) \\
&=&B(gx_{1}x_{2}\otimes 1_{H};GX_{1},x_{1})+B(x_{2}\otimes
x_{1};GX_{1},x_{1}).
\end{eqnarray*}%
In view of the form of the elements we get%
\begin{equation*}
2B(x_{2}\otimes 1_{H};G,x_{1})=B(x_{2}\otimes 1_{H};G,x_{1})+B(x_{2}\otimes
1_{H};G,x_{1})
\end{equation*}%
which is trivial.

\subparagraph{Case $f=x_{2}$}

\begin{eqnarray*}
&&B(x_{2}\otimes 1_{H};G,x_{2})+\lambda B(x_{2}\otimes
1_{H};GX_{1}X_{2},x_{2}) \\
&=&B(gx_{2}\otimes g;G,x_{2})-B(x_{1}x_{2}\otimes
g;GX_{1},x_{2})+B(x_{2}\otimes x_{1};GX_{1},x_{2}).
\end{eqnarray*}%
By applying $\left( \ref{eq.10}\right) $ this rewrites as%
\begin{eqnarray*}
&&2B(x_{2}\otimes 1_{H};G,x_{2})+\lambda B(x_{2}\otimes
1_{H};GX_{1}X_{2},x_{2}) \\
&=&B(gx_{1}x_{2}\otimes 1_{H};GX_{1},x_{2})+B(x_{2}\otimes
x_{1};GX_{1},x_{2}).
\end{eqnarray*}%
In view of the form of the elements we get%
\begin{eqnarray*}
&&2B(x_{2}\otimes 1_{H};G,x_{2})-\lambda B\left( g\otimes
1_{H};G,x_{1}x_{2}\right) \\
&=&2\left[ B(x_{2}\otimes 1_{H};G,x_{2})+B(gx_{1}x_{2}\otimes
1_{H};G,x_{1}x_{2})\right]
\end{eqnarray*}%
i.e.%
\begin{equation*}
\lambda B\left( g\otimes 1_{H};G,x_{1}x_{2}\right) +2B(gx_{1}x_{2}\otimes
1_{H};G,x_{1}x_{2})=0
\end{equation*}%
which is $\left( \ref{X1,x2,X1F21,x1x2}\right) .$

\subparagraph{Case $f=gx_{1}x_{2}$%
\protect\begin{eqnarray*}
&&B(x_{2}\otimes 1_{H};G,gx_{1}x_{2})+\protect\lambda B(x_{2}\otimes
1_{H};GX_{1}X_{2},gx_{1}x_{2}) \\
&=&B(gx_{2}\otimes g;G,gx_{1}x_{2})-B(x_{1}x_{2}\otimes
g;GX_{1},gx_{1}x_{2})+B(x_{2}\otimes x_{1};GX_{1},gx_{1}x_{2}).
\protect\end{eqnarray*}%
}

By applying $\left( \ref{eq.10}\right) $ this rewrites as%
\begin{equation*}
\lambda B(x_{2}\otimes 1_{H};GX_{1}X_{2},gx_{1}x_{2})=-B(gx_{1}x_{2}\otimes
1_{H};GX_{1},gx_{1}x_{2})+B(x_{2}\otimes x_{1};GX_{1},gx_{1}x_{2}).
\end{equation*}

In view of the form of the elements we get%
\begin{equation*}
0=-B(x_{2}\otimes 1_{H};G,gx_{1}x_{2})+B(x_{2}\otimes 1_{H};G,gx_{1}x_{2})
\end{equation*}%
which is trivial.

\paragraph{Equality $\left( \protect\ref{X1F71}\right) $}

rewrites as%
\begin{eqnarray*}
&&-\beta _{1}B(x_{2}\otimes 1_{H};GX_{1}X_{2},f) \\
&=&\beta _{1}B(gx_{2}\otimes g;GX_{1}X_{2},f)-B(x_{1}x_{2}\otimes
g;GX_{2},f)+B(x_{2}\otimes x_{1};GX_{2},f).
\end{eqnarray*}

\subparagraph{Case $f=g$}

\begin{eqnarray*}
&&-\beta _{1}B(x_{2}\otimes 1_{H};GX_{1}X_{2},g) \\
&=&\beta _{1}B(gx_{2}\otimes g;GX_{1}X_{2},g)-B(x_{1}x_{2}\otimes
g;GX_{2},g)+B(x_{2}\otimes x_{1};GX_{2},g).
\end{eqnarray*}%
By applying $\left( \ref{eq.10}\right) $ this rewrites as%
\begin{equation*}
-2\beta _{1}B(x_{2}\otimes 1_{H};GX_{1}X_{2},g)=-B(gx_{1}x_{2}\otimes
1_{H};GX_{2},g)+B(x_{2}\otimes x_{1};GX_{2},g).
\end{equation*}%
In view of the form of the elements we get%
\begin{eqnarray*}
&&-2\beta _{1}\left[ B(g\otimes 1_{H};G,gx_{1})+B(x_{2}\otimes \
1_{H};G,gx_{1}x_{2})\right] \\
&=&\left[ B(x_{1}\otimes 1_{H};G,g)+B(gx_{1}x_{2}\otimes 1_{H};G,gx_{2})%
\right] +B(x_{2}\otimes x_{1};GX_{2},g)
\end{eqnarray*}%
\begin{equation}
\beta _{1}\left[ B(g\otimes 1_{H};G,gx_{1})+B(x_{2}\otimes \
1_{H};G,gx_{1}x_{2})\right] =0  \label{X1,x2,X1F71,g}
\end{equation}

\subparagraph{Case $f=x_{1}$}

\begin{eqnarray*}
&&-\beta _{1}B(x_{2}\otimes 1_{H};GX_{1}X_{2},x_{1}) \\
&=&\beta _{1}B(gx_{2}\otimes g;GX_{1}X_{2},x_{1})-B(x_{1}x_{2}\otimes
g;GX_{2},x_{1})+B(x_{2}\otimes x_{1};GX_{2},x_{1}).
\end{eqnarray*}%
By applying $\left( \ref{eq.10}\right) $ this rewrites as%
\begin{equation*}
+B(gx_{1}x_{2}\otimes 1_{H};GX_{2},x_{1})+B(x_{2}\otimes
x_{1};GX_{2},x_{1})=0.
\end{equation*}%
In view of the form of the elements we get%
\begin{gather*}
+\left[ -B(x_{1}\otimes 1_{H};G,x_{1})-B(gx_{1}x_{2}\otimes
1_{H};G,x_{1}x_{2})\right] \\
+\left[ -B(x_{1}\otimes 1_{H};G,x_{1})-B(gx_{1}x_{2}\otimes
1_{H};G,x_{1}x_{2})\right] =0.
\end{gather*}%
\begin{equation}
B(x_{1}\otimes 1_{H};G,x_{1})+B(gx_{1}x_{2}\otimes 1_{H};G,x_{1}x_{2})=0
\label{X1,x2,X1F71,x1}
\end{equation}

\subparagraph{Case $f=x_{2}$}

\begin{eqnarray*}
&&-\beta _{1}B(x_{2}\otimes 1_{H};GX_{1}X_{2},x_{2}) \\
&=&\beta _{1}B(gx_{2}\otimes g;GX_{1}X_{2},x_{2})-B(x_{1}x_{2}\otimes
g;GX_{2},x_{2})+B(x_{2}\otimes x_{1};GX_{2},x_{2}).
\end{eqnarray*}%
By applying $\left( \ref{eq.10}\right) $ this rewrites as%
\begin{equation*}
+B(gx_{1}x_{2}\otimes 1_{H};GX_{2},x_{2})+B(x_{2}\otimes
x_{1};GX_{2},x_{2})=0.
\end{equation*}%
In view of the form of the elements we get%
\begin{equation*}
B(x_{1}\otimes 1_{H};G,x_{2})=0
\end{equation*}%
which is $\left( \ref{X1,g,X1F21,x2}\right) .$

\subparagraph{Case $f=gx_{1}x_{2}$}

\begin{eqnarray*}
&&-\beta _{1}B(x_{2}\otimes 1_{H};GX_{1}X_{2},gx_{1}x_{2}) \\
&=&\beta _{1}B(gx_{2}\otimes g;GX_{1}X_{2},gx_{1}x_{2})-B(x_{1}x_{2}\otimes
g;GX_{2},gx_{1}x_{2})+ \\
&&+B(x_{2}\otimes x_{1};GX_{2},gx_{1}x_{2}).
\end{eqnarray*}%
By applying $\left( \ref{eq.10}\right) $ this rewrites as%
\begin{equation*}
-2\beta _{1}B(x_{2}\otimes
1_{H};GX_{1}X_{2},gx_{1}x_{2})=-B(gx_{1}x_{2}\otimes
1_{H};GX_{2},gx_{1}x_{2})+B(x_{2}\otimes x_{1};GX_{2},gx_{1}x_{2}).
\end{equation*}%
In view of the form of the elements we get%
\begin{equation*}
B(x_{1}\otimes 1_{H};G,gx_{1}x_{2})+\left[ -2B(g\otimes
1_{H};G,gx_{2})+B(x_{1}\otimes 1_{H};G,gx_{1}x_{2})\right] =0
\end{equation*}%
i.e.%
\begin{equation*}
B(x_{1}\otimes 1_{H};G,gx_{1}x_{2})-B(g\otimes 1_{H};G,gx_{2})=0
\end{equation*}%
which is $\left( \ref{X1,g,X1F21,gx1x2}\right) .$

\paragraph{Equality $\left( \protect\ref{X1F81}\right) $}

rewrites as%
\begin{eqnarray*}
&&-B(x_{2}\otimes 1_{H};GX_{2},f) \\
&=&B(gx_{2}\otimes g;GX_{2},f)-B(x_{1}x_{2}\otimes
g;GX_{1}X_{2},f)+B(x_{2}\otimes x_{1};GX_{1}X_{2},f).
\end{eqnarray*}

\subparagraph{Case $f=1_{H}$}

\begin{eqnarray*}
&&-B(x_{2}\otimes 1_{H};GX_{2},1_{H}) \\
&=&B(gx_{2}\otimes g;GX_{2},1_{H})-B(x_{1}x_{2}\otimes
g;GX_{1}X_{2},1_{H})+B(x_{2}\otimes x_{1};GX_{1}X_{2},1_{H}).
\end{eqnarray*}%
By applying $\left( \ref{eq.10}\right) $ this rewrites as%
\begin{eqnarray*}
&&-2B(x_{2}\otimes 1_{H};GX_{2},1_{H}) \\
&=&-B(gx_{1}x_{2}\otimes 1_{H};GX_{1}X_{2},1_{H})+B(x_{2}\otimes
x_{1};GX_{1}X_{2},1_{H}).
\end{eqnarray*}%
In view of the form of the elements we get%
\begin{gather*}
B(g\otimes 1_{H};G,1_{H})+B(x_{2}\otimes \ 1_{H};G,x_{2})=B(g\otimes
1_{H};G,1_{H})+ \\
+B(x_{2}\otimes \ 1_{H};G,x_{2})+B(x_{1}\otimes
1_{H};G,x_{1})+B(gx_{1}x_{2}\otimes 1_{H};G,x_{1}x_{2})
\end{gather*}%
i.e.%
\begin{equation*}
B(x_{1}\otimes 1_{H};G,x_{1})+B(gx_{1}x_{2}\otimes 1_{H};G,x_{1}x_{2})=0
\end{equation*}%
which is $\left( \ref{X1,x2,X1F71,x1}\right) .$

\subparagraph{Case $f=x_{1}x_{2}$}

\begin{eqnarray*}
&&-B(x_{2}\otimes 1_{H};GX_{2},x_{1}x_{2}) \\
&=&B(gx_{2}\otimes g;GX_{2},x_{1}x_{2})-B(x_{1}x_{2}\otimes
g;GX_{1}X_{2},x_{1}x_{2})+B(x_{2}\otimes x_{1};GX_{1}X_{2},x_{1}x_{2}).
\end{eqnarray*}%
By applying $\left( \ref{eq.10}\right) $ this rewrites as%
\begin{equation*}
-2B(x_{2}\otimes 1_{H};GX_{2},x_{1}x_{2})=-B(gx_{1}x_{2}\otimes
1_{H};GX_{1}X_{2},x_{1}x_{2})+B(x_{2}\otimes x_{1};GX_{1}X_{2},x_{1}x_{2}).
\end{equation*}%
In view of the form of the elements we get%
\begin{equation*}
-2B\left( g\otimes 1_{H};G,x_{1}x_{2}\right) =-2B(g\otimes
1_{H};G,x_{1}x_{2})
\end{equation*}%
which is trivial.

\subparagraph{Case $f=gx_{1}$}

\begin{eqnarray*}
&&-B(x_{2}\otimes 1_{H};GX_{2},gx_{1}) \\
&=&B(gx_{2}\otimes g;GX_{2},gx_{1})-B(x_{1}x_{2}\otimes
g;GX_{1}X_{2},gx_{1})+B(x_{2}\otimes x_{1};GX_{1}X_{2},gx_{1}).
\end{eqnarray*}%
By applying $\left( \ref{eq.10}\right) $ this rewrites as

\begin{equation*}
B(gx_{1}x_{2}\otimes 1_{H};GX_{1}X_{2},gx_{1})+B(x_{2}\otimes
x_{1};GX_{1}X_{2},gx_{1})=0.
\end{equation*}%
In view of the form of the elements we get%
\begin{gather*}
\left[ B(g\otimes 1_{H};G,gx_{1})+B(x_{2}\otimes \ 1_{H};G,gx_{1}x_{2})%
\right] \\
+\left[ -B(g\otimes 1_{H};G,gx_{1})-B(x_{2}\otimes \ 1_{H};G,gx_{1}x_{2})%
\right] =0
\end{gather*}%
which is trivial.

\subparagraph{Case $f=gx_{2}$}

\begin{eqnarray*}
&&-B(x_{2}\otimes 1_{H};GX_{2},gx_{2}) \\
&=&B(gx_{2}\otimes g;GX_{2},gx_{2})-B(x_{1}x_{2}\otimes
g;GX_{1}X_{2},gx_{2})+B(x_{2}\otimes x_{1};GX_{1}X_{2},gx_{2}).
\end{eqnarray*}%
By applying $\left( \ref{eq.10}\right) $ this rewrites as

\begin{equation*}
B(gx_{1}x_{2}\otimes 1_{H};GX_{1}X_{2},gx_{2})+B(x_{2}\otimes
x_{1};GX_{1}X_{2},gx_{2})=0.
\end{equation*}%
In view of the form of the elements we get%
\begin{equation*}
B(g\otimes 1_{H};G,gx_{2})-B(x_{1}\otimes 1_{H};G,gx_{1}x_{2})=0
\end{equation*}%
which is $\left( \ref{X1,x2,X1F71,gx1x2}\right) .$

\subsubsection{Case $x_{1}x_{2}\otimes 1_{H}$}

\paragraph{Equality $\left( \protect\ref{X1F11}\right) $}

rewrites as%
\begin{eqnarray*}
&&\beta _{1}B(x_{1}x_{2}\otimes 1_{H};X_{1},f)+\lambda B(x_{1}x_{2}\otimes
1_{H};X_{2},f) \\
&=&\gamma _{1}B(gx_{1}x_{2}\otimes g;G,f)+\beta _{1}B(gx_{1}x_{2}\otimes
g;X_{1},f)+B(x_{1}x_{2}\otimes x_{1};1_{A},f).
\end{eqnarray*}

\subparagraph{Case $f=g$}

\begin{eqnarray*}
&&\beta _{1}B(x_{1}x_{2}\otimes 1_{H};X_{1},g)+\lambda B(x_{1}x_{2}\otimes
1_{H};X_{2},g) \\
&=&\gamma _{1}B(gx_{1}x_{2}\otimes g;G,g)+\beta _{1}B(gx_{1}x_{2}\otimes
g;X_{1},g)+B(x_{1}x_{2}\otimes x_{1};1_{A},g).
\end{eqnarray*}%
By applying $\left( \ref{eq.10}\right) $ this rewrites as%
\begin{equation*}
\lambda B(x_{1}x_{2}\otimes 1_{H};X_{2},g)=\gamma _{1}B(x_{1}x_{2}\otimes
1_{H};G,g)+B(x_{1}x_{2}\otimes x_{1};1_{A},g).
\end{equation*}%
In view of the form of the elements we get%
\begin{equation}
\lambda B(x_{1}x_{2}\otimes 1_{H};X_{2},g)-\gamma _{1}B(x_{1}x_{2}\otimes
1_{H};G,g)=0  \label{X1,x1x2,X1F11,g}
\end{equation}

\subparagraph{Case $f=x_{1}$}

\begin{gather*}
\beta _{1}B(x_{1}x_{2}\otimes 1_{H};X_{1},x_{1})+\lambda B(x_{1}x_{2}\otimes
1_{H};X_{2},x_{1})= \\
\gamma _{1}B(gx_{1}x_{2}\otimes g;G,x_{1})+\beta _{1}B(gx_{1}x_{2}\otimes
g;X_{1},x_{1})+B(x_{1}x_{2}\otimes x_{1};1_{A},x_{1}).
\end{gather*}%
By applying $\left( \ref{eq.10}\right) $ this rewrites as%
\begin{eqnarray*}
&&2\beta _{1}B(x_{1}x_{2}\otimes 1_{H};X_{1},x_{1})+\lambda
B(x_{1}x_{2}\otimes 1_{H};X_{2},x_{1}) \\
&=&-\gamma _{1}B(x_{1}x_{2}\otimes 1_{H};G,x_{1})+B(x_{1}x_{2}\otimes
x_{1};1_{A},x_{1}).
\end{eqnarray*}%
In view of the form of the elements we get%
\begin{equation}
2\beta _{1}B(x_{1}x_{2}\otimes 1_{H};X_{1},x_{1})+\lambda
B(x_{1}x_{2}\otimes 1_{H};X_{2},x_{1})+\gamma _{1}B(x_{1}x_{2}\otimes
1_{H};G,x_{1})=0.  \label{X1,x1x2,X1F11,x1}
\end{equation}

\subparagraph{Case $f=x_{2}$}

\begin{eqnarray*}
&&\beta _{1}B(x_{1}x_{2}\otimes 1_{H};X_{1},x_{2})+\lambda
B(x_{1}x_{2}\otimes 1_{H};X_{2},x_{2}) \\
&=&\gamma _{1}B(gx_{1}x_{2}\otimes g;G,x_{2})+\beta _{1}B(gx_{1}x_{2}\otimes
g;X_{1},x_{2}) \\
&&+B(x_{1}x_{2}\otimes x_{1};1_{A},x_{2}).
\end{eqnarray*}%
By applying $\left( \ref{eq.10}\right) $ this rewrites as%
\begin{eqnarray*}
&&2\beta _{1}B(x_{1}x_{2}\otimes 1_{H};X_{1},x_{2})+\lambda
B(x_{1}x_{2}\otimes 1_{H};X_{2},x_{2}) \\
&=&-\gamma _{1}B(x_{1}x_{2}\otimes 1_{H};G,x_{2})+B(x_{1}x_{2}\otimes
x_{1};1_{A},x_{2}).
\end{eqnarray*}%
In view of the form of the elements we get%
\begin{equation}
2\beta _{1}B(x_{1}x_{2}\otimes 1_{H};X_{1},x_{2})+\lambda
B(x_{1}x_{2}\otimes 1_{H};X_{2},x_{2})+\gamma _{1}B(x_{1}x_{2}\otimes
1_{H};G,x_{2})=0.  \label{X1,x1x2,X1F11,x2}
\end{equation}

\subparagraph{Case $f=gx_{1}x_{2}$}

\begin{eqnarray*}
&&\beta _{1}B(x_{1}x_{2}\otimes 1_{H};X_{1},gx_{1}x_{2})+\lambda
B(x_{1}x_{2}\otimes 1_{H};X_{2},gx_{1}x_{2}) \\
&=&\gamma _{1}B(gx_{1}x_{2}\otimes g;G,gx_{1}x_{2})+\beta
_{1}B(gx_{1}x_{2}\otimes g;X_{1},gx_{1}x_{2})+B(x_{1}x_{2}\otimes
x_{1};1_{A},gx_{1}x_{2}).
\end{eqnarray*}%
By applying $\left( \ref{eq.10}\right) $ this rewrites as%
\begin{equation*}
+\lambda B(x_{1}x_{2}\otimes 1_{H};X_{2},gx_{1}x_{2})=\gamma
_{1}B(x_{1}x_{2}\otimes 1_{H};G,gx_{1}x_{2})+B(x_{1}x_{2}\otimes
x_{1};1_{A},gx_{1}x_{2}).
\end{equation*}%
In view of the form of the elements we get%
\begin{equation*}
+\lambda B(x_{1}x_{2}\otimes 1_{H};X_{2},gx_{1}x_{2})=\gamma
_{1}B(x_{1}x_{2}\otimes 1_{H};G,gx_{1}x_{2})-2B(x_{1}x_{2}\otimes
1_{H};1_{A},gx_{2}).
\end{equation*}%
By $\left( \ref{G,x1x2, GF7,gx1x2}\right) $ we obtain

\begin{equation}
\gamma _{1}B(x_{1}x_{2}\otimes 1_{H};G,gx_{1}x_{2})-2B(x_{1}x_{2}\otimes
1_{H};1_{A},gx_{2})=0.  \label{X1,x1x2,X1F11,gx1x2}
\end{equation}

\paragraph{Equality $\left( \protect\ref{X1F21}\right) $}

rewrites as%
\begin{eqnarray*}
&&\beta _{1}B(x_{1}x_{2}\otimes 1_{H};GX_{1},f)+\lambda B(x_{1}x_{2}\otimes
1_{H};GX_{2},f) \\
&=&-\beta _{1}B(gx_{1}x_{2}\otimes g;GX_{1},f)+B(x_{1}x_{2}\otimes
x_{1};G,f).
\end{eqnarray*}

\subparagraph{Case $f=1_{H}$%
\protect\begin{eqnarray*}
&&\protect\beta _{1}B(x_{1}x_{2}\otimes 1_{H};GX_{1},1_{H})+\protect\lambda %
B(x_{1}x_{2}\otimes 1_{H};GX_{2},1_{H}) \\
&=&-\protect\beta _{1}B(gx_{1}x_{2}\otimes
g;GX_{1},1_{H})+B(x_{1}x_{2}\otimes x_{1};G,1_{H}).
\protect\end{eqnarray*}%
}

By applying $\left( \ref{eq.10}\right) $ this rewrites as%
\begin{equation*}
2\beta _{1}B(x_{1}x_{2}\otimes 1_{H};GX_{1},1_{H})+\lambda
B(x_{1}x_{2}\otimes 1_{H};GX_{2},1_{H})=B(x_{1}x_{2}\otimes x_{1};G,1_{H}).
\end{equation*}%
In view of the form of the elements we get

\begin{equation}
2\beta _{1}B(x_{1}x_{2}\otimes 1_{H};GX_{1},1_{H})+\lambda
B(x_{1}x_{2}\otimes 1_{H};GX_{2},1_{H})=0.  \label{X1,x1x2,X1F21,1H}
\end{equation}

\subparagraph{Case $f=gx_{1}$}

\begin{eqnarray*}
&&\beta _{1}B(x_{1}x_{2}\otimes 1_{H};GX_{1},gx_{1})+\lambda
B(x_{1}x_{2}\otimes 1_{H};GX_{2},gx_{1}) \\
&=&-\beta _{1}B(gx_{1}x_{2}\otimes g;GX_{1},gx_{1})+B(x_{1}x_{2}\otimes
x_{1};G,gx_{1}).
\end{eqnarray*}%
By applying $\left( \ref{eq.10}\right) $ this rewrites as%
\begin{equation*}
\lambda B(x_{1}x_{2}\otimes 1_{H};GX_{2},gx_{1})=B(x_{1}x_{2}\otimes
x_{1};G,gx_{1}).
\end{equation*}%
In view of the form of the elements we get%
\begin{equation}
\lambda B(x_{1}x_{2}\otimes 1_{H};GX_{2},gx_{1})+2B(x_{1}x_{2}\otimes
1_{H};G,g)=0.  \label{X1,x1x2,X1F21,gx1}
\end{equation}

\subparagraph{Case $f=gx_{2}$}

\begin{eqnarray*}
&&\beta _{1}B(x_{1}x_{2}\otimes 1_{H};GX_{1},gx_{2})+\lambda
B(x_{1}x_{2}\otimes 1_{H};GX_{2},gx_{2}) \\
&=&-\beta _{1}B(gx_{1}x_{2}\otimes g;GX_{1},gx_{2})+B(x_{1}x_{2}\otimes
x_{1};G,gx_{2}).
\end{eqnarray*}%
By applying $\left( \ref{eq.10}\right) $ this rewrites as%
\begin{equation*}
\lambda B(x_{1}x_{2}\otimes 1_{H};GX_{2},gx_{2})=B(x_{1}x_{2}\otimes
x_{1};G,gx_{2}).
\end{equation*}%
In view of the form of the elements we get%
\begin{equation}
\lambda B(x_{1}x_{2}\otimes 1_{H};GX_{2},gx_{2})=0.
\label{X1,x1x2,X1F21,gx2}
\end{equation}

\paragraph{Equality $\left( \protect\ref{X1F31}\right) $}

rewrites as%
\begin{eqnarray*}
&&B(x_{1}x_{2}\otimes 1_{H};1_{A},f)+\lambda B(x_{1}x_{2}\otimes
1_{H};X_{1}X_{2},f) \\
&=&\gamma _{1}B(gx_{1}x_{2}\otimes g;GX_{1},f)+B(gx_{1}x_{2}\otimes
g;1_{A},f)+ \\
&&+B(x_{1}x_{2}\otimes x_{1};X_{1},f)
\end{eqnarray*}

\subparagraph{Case $f=1_{H}$}

\begin{eqnarray*}
&&B(x_{1}x_{2}\otimes 1_{H};1_{A},1_{H})+\lambda B(x_{1}x_{2}\otimes
1_{H};X_{1}X_{2},1_{H}) \\
&=&\gamma _{1}B(gx_{1}x_{2}\otimes g;GX_{1},1_{H})+B(gx_{1}x_{2}\otimes
g;1_{A},1_{H})+ \\
&&+B(x_{1}x_{2}\otimes x_{1};X_{1},1_{H})
\end{eqnarray*}%
By applying $\left( \ref{eq.10}\right) $ this rewrites as%
\begin{equation*}
\lambda B(x_{1}x_{2}\otimes 1_{H};X_{1}X_{2},1_{H})=\gamma
_{1}B(x_{1}x_{2}\otimes 1_{H};GX_{1},1_{H})+B(x_{1}x_{2}\otimes
x_{1};X_{1},1_{H})
\end{equation*}%
In view of the form of the elements we get%
\begin{gather}
\lambda \left[
\begin{array}{c}
+1-B(x_{1}x_{2}\otimes 1_{H};1_{A},x_{1}x_{2}) \\
-B(x_{1}x_{2}\otimes 1_{H};X_{2},x_{1})+B(x_{1}x_{2}\otimes
1_{H};X_{1},x_{2})%
\end{array}%
\right] +  \label{X1,x1x2,X1F31,1H} \\
-\gamma _{1}B(x_{1}x_{2}\otimes 1_{H};GX_{1},1_{H})=0  \notag
\end{gather}

\subparagraph{Case $f=x_{1}x_{2}$}

\begin{eqnarray*}
&&B(x_{1}x_{2}\otimes 1_{H};1_{A},x_{1}x_{2})+\lambda B(x_{1}x_{2}\otimes
1_{H};X_{1}X_{2},x_{1}x_{2}) \\
&=&\gamma _{1}B(gx_{1}x_{2}\otimes g;GX_{1},x_{1}x_{2})+B(gx_{1}x_{2}\otimes
g;1_{A},x_{1}x_{2})+B(x_{1}x_{2}\otimes x_{1};X_{1},x_{1}x_{2})
\end{eqnarray*}%
By applying $\left( \ref{eq.10}\right) $ this rewrites as%
\begin{equation*}
\lambda B(x_{1}x_{2}\otimes 1_{H};X_{1}X_{2},x_{1}x_{2})=\gamma
_{1}B(x_{1}x_{2}\otimes 1_{H};GX_{1},x_{1}x_{2})+B(x_{1}x_{2}\otimes
x_{1};X_{1},x_{1}x_{2})
\end{equation*}%
In view of the form of the elements we get%
\begin{equation}
\gamma _{1}B(x_{1}x_{2}\otimes 1_{H};GX_{1},x_{1}x_{2})=0
\label{X1,x1x2,X1F31,x1x2}
\end{equation}%
.

\subparagraph{Case $f=gx_{1}$}

\begin{eqnarray*}
&&B(x_{1}x_{2}\otimes 1_{H};1_{A},gx_{1})+\lambda B(x_{1}x_{2}\otimes
1_{H};X_{1}X_{2},gx_{1}) \\
&=&\gamma _{1}B(gx_{1}x_{2}\otimes g;GX_{1},gx_{1})+B(gx_{1}x_{2}\otimes
g;1_{A},gx_{1})+B(x_{1}x_{2}\otimes x_{1};X_{1},gx_{1})
\end{eqnarray*}%
By applying $\left( \ref{eq.10}\right) $ this rewrites as%
\begin{eqnarray*}
&&2B(x_{1}x_{2}\otimes 1_{H};1_{A},gx_{1})+\lambda B(x_{1}x_{2}\otimes
1_{H};X_{1}X_{2},gx_{1}) \\
&=&-\gamma _{1}B(x_{1}x_{2}\otimes 1_{H};GX_{1},gx_{1})+B(x_{1}x_{2}\otimes
x_{1};X_{1},gx_{1}).
\end{eqnarray*}%
In view of the form of the elements we get%
\begin{eqnarray*}
&&2B(x_{1}x_{2}\otimes 1_{H};1_{A},gx_{1})+\lambda B(x_{1}x_{2}\otimes
1_{H};X_{1},gx_{1}x_{2}) \\
&=&-\gamma _{1}B(x_{1}x_{2}\otimes 1_{H};GX_{1},gx_{1})-2B(x_{1}x_{2}\otimes
1_{H};X_{1},g).
\end{eqnarray*}

By $\left( \ref{G,x1x2, GF6,gx1x2}\right) $ we obtain

\begin{equation}
2B(x_{1}x_{2}\otimes 1_{H};1_{A},gx_{1})+\gamma _{1}B(x_{1}x_{2}\otimes
1_{H};GX_{1},gx_{1})+2B(x_{1}x_{2}\otimes 1_{H};X_{1},g)=0
\label{X1,x1x2,X1F31,gx1}
\end{equation}

\subparagraph{Case $f=gx_{2}$}

\begin{eqnarray*}
&&B(x_{1}x_{2}\otimes 1_{H};1_{A},gx_{2})+\lambda B(x_{1}x_{2}\otimes
1_{H};X_{1}X_{2},gx_{2}) \\
&=&\gamma _{1}B(gx_{1}x_{2}\otimes g;GX_{1},gx_{2})+B(gx_{1}x_{2}\otimes
g;1_{A},gx_{2})+ \\
&&+B(x_{1}x_{2}\otimes x_{1};X_{1},gx_{2})
\end{eqnarray*}%
By applying $\left( \ref{eq.10}\right) $ this rewrites as%
\begin{eqnarray*}
&&2B(x_{1}x_{2}\otimes 1_{H};1_{A},gx_{2})+\lambda B(x_{1}x_{2}\otimes
1_{H};X_{1}X_{2},gx_{2}) \\
&=&-\gamma _{1}B(x_{1}x_{2}\otimes 1_{H};GX_{1},gx_{2})+B(x_{1}x_{2}\otimes
x_{1};X_{1},gx_{2}).
\end{eqnarray*}%
In view of the form of the elements we get

\begin{eqnarray*}
&&2B(x_{1}x_{2}\otimes 1_{H};1_{A},gx_{2})+\lambda B(x_{1}x_{2}\otimes
1_{H};X_{2},gx_{1}x_{2}) \\
&=&-\gamma _{1}B(x_{1}x_{2}\otimes 1_{H};GX_{1},gx_{2})+B(x_{1}x_{2}\otimes
x_{1};X_{1},gx_{2}).
\end{eqnarray*}%
In view of the form of the elements we get%
\begin{gather*}
2B(x_{1}x_{2}\otimes 1_{H};1_{A},gx_{2})+\lambda B(x_{1}x_{2}\otimes
1_{H};X_{2},gx_{1}x_{2}) \\
+\gamma _{1}B(x_{1}x_{2}\otimes 1_{H};GX_{1},gx_{2})=0.
\end{gather*}%
By $\left( \ref{G,x1x2, GF7,gx1x2}\right) $ we get%
\begin{equation}
2B(x_{1}x_{2}\otimes 1_{H};1_{A},gx_{2})+\gamma _{1}B(x_{1}x_{2}\otimes
1_{H};GX_{1},gx_{2})=0.  \label{X1,x1x2,X1F31,gx2}
\end{equation}

\paragraph{Equality $\left( \protect\ref{X1F41}\right) $}

rewrites as%
\begin{gather*}
-\beta _{1}B(x_{1}x_{2}\otimes 1_{H};X_{1}X_{2},f)=\beta
_{1}B(gx_{1}x_{2}\otimes g;X_{1}X_{2},f)+ \\
+\gamma _{1}B(gx_{1}x_{2}\otimes g;GX_{2},f)+B(x_{1}x_{2}\otimes
x_{1};X_{2},f)
\end{gather*}

\subparagraph{Case $f=1_{H}$}

\begin{gather*}
-\beta _{1}B(x_{1}x_{2}\otimes 1_{H};X_{1}X_{2},1_{H})=\beta
_{1}B(gx_{1}x_{2}\otimes g;X_{1}X_{2},1_{H})+ \\
+\gamma _{1}B(gx_{1}x_{2}\otimes g;GX_{2},1_{H})+B(x_{1}x_{2}\otimes
x_{1};X_{2},1_{H})
\end{gather*}

By applying $\left( \ref{eq.10}\right) $ this rewrites as%
\begin{gather*}
2\beta _{1}B(x_{1}x_{2}\otimes 1_{H};X_{1}X_{2},1_{H})+\gamma
_{1}B(x_{1}x_{2}\otimes 1_{H};GX_{2},1_{H}) \\
+B(x_{1}x_{2}\otimes x_{1};X_{2},1_{H})=0.
\end{gather*}%
In view of the form of the elements we get

\begin{gather}
2\beta _{1}\left[
\begin{array}{c}
+1-B(x_{1}x_{2}\otimes 1_{H};1_{A},x_{1}x_{2})-B(x_{1}x_{2}\otimes
1_{H};X_{2},x_{1}) \\
+B(x_{1}x_{2}\otimes 1_{H};X_{1},x_{2})%
\end{array}%
\right]  \label{X1,x1x2,X1F41,1H} \\
+\gamma _{1}B(x_{1}x_{2}\otimes 1_{H};GX_{2},1_{H})=0.  \notag
\end{gather}

\subparagraph{Case $f=gx_{1}$}

\begin{gather*}
-\beta _{1}B(x_{1}x_{2}\otimes 1_{H};X_{1}X_{2},gx_{1})=\beta
_{1}B(gx_{1}x_{2}\otimes g;X_{1}X_{2},gx_{1}) \\
+\gamma _{1}B(gx_{1}x_{2}\otimes g;GX_{2},gx_{1})+B(x_{1}x_{2}\otimes
x_{1};X_{2},gx_{1})
\end{gather*}%
By applying $\left( \ref{eq.10}\right) $ this rewrites as%
\begin{equation*}
-\gamma _{1}B(x_{1}x_{2}\otimes 1_{H};GX_{2},gx_{1})+B(x_{1}x_{2}\otimes
x_{1};X_{2},gx_{1})=0
\end{equation*}%
In view of the form of the elements we get%
\begin{equation}
\gamma _{1}B(x_{1}x_{2}\otimes 1_{H};GX_{2},gx_{1})+2B(x_{1}x_{2}\otimes
1_{H};X_{2},g)=0.  \label{X1,x1x2,X1F41,gx1}
\end{equation}

\subparagraph{Case $f=gx_{2}$}

\begin{gather*}
-\beta _{1}B(x_{1}x_{2}\otimes 1_{H};X_{1}X_{2},gx_{2})=\beta
_{1}B(gx_{1}x_{2}\otimes g;X_{1}X_{2},gx_{2}) \\
+\gamma _{1}B(gx_{1}x_{2}\otimes g;GX_{2},gx_{2})+B(x_{1}x_{2}\otimes
x_{1};X_{2},gx_{2})
\end{gather*}%
By applying $\left( \ref{eq.10}\right) $ this rewrites as%
\begin{equation*}
-\gamma _{1}B(x_{1}x_{2}\otimes 1_{H};GX_{2},gx_{2})+B(x_{1}x_{2}\otimes
x_{1};X_{2},gx_{2})=0
\end{equation*}%
In view of the form of the elements we get%
\begin{equation}
\gamma _{1}B(x_{1}x_{2}\otimes 1_{H};GX_{2},gx_{2})=0
\label{X1,x1x2,X1F41,gx2}
\end{equation}

\paragraph{Equality $\left( \protect\ref{X1F51}\right) $}

rewrites as%
\begin{eqnarray*}
&&-B(x_{1}x_{2}\otimes 1_{H};X_{2},f) \\
&=&B(gx_{1}x_{2}\otimes g;X_{2},f)+\gamma _{1}B(gx_{1}x_{2}\otimes
g;GX_{1}X_{2},f)+ \\
&&+B(x_{1}x_{2}\otimes x_{1};X_{1}X_{2},f)
\end{eqnarray*}

\subparagraph{Case $f=g$}

\begin{gather*}
-B(x_{1}x_{2}\otimes 1_{H};X_{2},g)=B(gx_{1}x_{2}\otimes g;X_{2},g)+ \\
+\gamma _{1}B(gx_{1}x_{2}\otimes g;GX_{1}X_{2},g)+B(x_{1}x_{2}\otimes
x_{1};X_{1}X_{2},g)
\end{gather*}%
By applying $\left( \ref{eq.10}\right) $ this rewrites as%
\begin{equation*}
2B(x_{1}x_{2}\otimes 1_{H};X_{2},g)+\gamma _{1}B(x_{1}x_{2}\otimes
1_{H};GX_{1}X_{2},g)+B(x_{1}x_{2}\otimes x_{1};X_{1}X_{2},g)=0
\end{equation*}%
In view of the form of the elements we get%
\begin{gather}
2B(x_{1}x_{2}\otimes 1_{H};X_{2},g)+  \label{X1,x1x2,X1F51,g} \\
+\gamma _{1}\left[
\begin{array}{c}
-B(x_{1}x_{2}\otimes 1_{H};G,gx_{1}x_{2})+ \\
B(x_{1}x_{2}\otimes 1_{H};GX_{2},gx_{1})+ \\
-B(x_{1}x_{2}\otimes 1_{H};GX_{1},gx_{2})%
\end{array}%
\right] =0  \notag
\end{gather}

\subparagraph{Case $f=x_{1}$%
\protect\begin{eqnarray*}
&&-B(x_{1}x_{2}\otimes 1_{H};X_{2},x_{1}) \\
&=&B(gx_{1}x_{2}\otimes g;X_{2},x_{1})+\protect\gamma _{1}B(gx_{1}x_{2}%
\otimes g;GX_{1}X_{2},x_{1})+B(x_{1}x_{2}\otimes x_{1};X_{1}X_{2},x_{1})
\protect\end{eqnarray*}%
}

By applying $\left( \ref{eq.10}\right) $ this rewrites as%
\begin{equation*}
-\gamma _{1}B(x_{1}x_{2}\otimes 1_{H};GX_{1},x_{1}x_{2})+B(x_{1}x_{2}\otimes
x_{1};X_{1}X_{2},x_{1})=0
\end{equation*}%
In view of the form of the elements we get%
\begin{equation*}
\gamma _{1}B(x_{1}x_{2}\otimes 1_{H};GX_{1},x_{1}x_{2})=0
\end{equation*}

which is $\left( \ref{X1,x1x2,X1F31,x1x2}\right) .$

\subparagraph{Case $f=x_{2}$%
\protect\begin{eqnarray*}
&&-B(x_{1}x_{2}\otimes 1_{H};X_{2},x_{2}) \\
&=&B(gx_{1}x_{2}\otimes g;X_{2},x_{2})+\protect\gamma _{1}B(gx_{1}x_{2}%
\otimes g;GX_{1}X_{2},x_{2})+B(x_{1}x_{2}\otimes x_{1};X_{1}X_{2},x_{2})
\protect\end{eqnarray*}%
}

By applying $\left( \ref{eq.10}\right) $ this rewrites as%
\begin{equation*}
-\gamma _{1}B(x_{1}x_{2}\otimes 1_{H};GX_{1}X_{2},x_{2})+B(x_{1}x_{2}\otimes
x_{1};X_{1}X_{2},x_{2})=0
\end{equation*}%
In view of the form of the elements we get%
\begin{equation*}
\gamma _{1}B(x_{1}x_{2}\otimes 1_{H};GX_{2},x_{1}x_{2})=0
\end{equation*}

which is $\left( \ref{G,x1x2, GF7,x2}\right) .$

\paragraph{Equality $\left( \protect\ref{X1F61}\right) $}

rewrites as%
\begin{eqnarray*}
&&B(x_{1}x_{2}\otimes 1_{H};G,f)+\lambda B(x_{1}x_{2}\otimes
1_{H};GX_{1}X_{2},f) \\
&=&-B(gx_{1}x_{2}\otimes g;G,f)+B(x_{1}x_{2}\otimes x_{1};GX_{1},f).
\end{eqnarray*}

\subparagraph{Case $f=g$}

\begin{eqnarray*}
&&B(x_{1}x_{2}\otimes 1_{H};G,g)+\lambda B(x_{1}x_{2}\otimes
1_{H};GX_{1}X_{2},g) \\
&=&-B(gx_{1}x_{2}\otimes g;G,g)+B(x_{1}x_{2}\otimes x_{1};GX_{1},g)
\end{eqnarray*}%
By applying $\left( \ref{eq.10}\right) $ this rewrites as%
\begin{equation*}
2B(x_{1}x_{2}\otimes 1_{H};G,g)+\lambda B(x_{1}x_{2}\otimes
1_{H};GX_{1}X_{2},g)=B(x_{1}x_{2}\otimes x_{1};GX_{1},g)
\end{equation*}%
In view of the form of the elements we get

\begin{equation}
2B(x_{1}x_{2}\otimes 1_{H};G,g)+\lambda \left[
\begin{array}{c}
-B(x_{1}x_{2}\otimes 1_{H};G,gx_{1}x_{2}) \\
+B(x_{1}x_{2}\otimes 1_{H};GX_{2},gx_{1})+ \\
-B(x_{1}x_{2}\otimes 1_{H};GX_{1},gx_{2})%
\end{array}%
\right] =0  \label{X1,x1x2,X1F61,g}
\end{equation}

\subparagraph{Case $f=x_{1}$}

\begin{eqnarray*}
&&B(x_{1}x_{2}\otimes 1_{H};G,x_{1})+\lambda B(x_{1}x_{2}\otimes
1_{H};GX_{1}X_{2},x_{1}) \\
&=&-B(gx_{1}x_{2}\otimes g;G,x_{1})+B(x_{1}x_{2}\otimes x_{1};GX_{1},x_{1})
\end{eqnarray*}%
By applying $\left( \ref{eq.10}\right) $ this rewrites as%
\begin{equation*}
\lambda B(x_{1}x_{2}\otimes 1_{H};GX_{1},x_{1}x_{2})+B(x_{1}x_{2}\otimes
x_{1};GX_{1},x_{1})=0.
\end{equation*}%
In view of the form of the elements we get

\begin{equation}
\lambda B(x_{1}x_{2}\otimes 1_{H};GX_{1},x_{1}x_{2})=0.
\label{X1,x1x2,X1F61,x1}
\end{equation}

\subparagraph{Case $f=x_{2}$}

\begin{eqnarray*}
&&B(x_{1}x_{2}\otimes 1_{H};G,x_{2})+\lambda B(x_{1}x_{2}\otimes
1_{H};GX_{1}X_{2},x_{2}) \\
&=&-B(gx_{1}x_{2}\otimes g;G,x_{2})+B(x_{1}x_{2}\otimes x_{1};GX_{1},x_{2})
\end{eqnarray*}

By applying $\left( \ref{eq.10}\right) $ this rewrites as%
\begin{equation*}
-\lambda B(x_{1}x_{2}\otimes 1_{H};GX_{2},x_{1}x_{2})=+B(x_{1}x_{2}\otimes
x_{1};GX_{1},x_{2})
\end{equation*}%
In view of the form of the elements we get%
\begin{equation}
\lambda B(x_{1}x_{2}\otimes 1_{H};GX_{2},x_{1}x_{2})=0
\label{X1,x1x2,X1F61,x2}
\end{equation}

\subparagraph{Case $f=gx_{1}x_{2}$}

\begin{eqnarray*}
&&B(x_{1}x_{2}\otimes 1_{H};G,gx_{1}x_{2})+\lambda B(x_{1}x_{2}\otimes
1_{H};GX_{1}X_{2},gx_{1}x_{2}) \\
&=&-B(gx_{1}x_{2}\otimes g;G,gx_{1}x_{2})+B(x_{1}x_{2}\otimes
x_{1};GX_{1},gx_{1}x_{2})
\end{eqnarray*}%
By applying $\left( \ref{eq.10}\right) $ this rewrites as%
\begin{eqnarray*}
&&2B(x_{1}x_{2}\otimes 1_{H};G,gx_{1}x_{2})+\lambda B(x_{1}x_{2}\otimes
1_{H};GX_{1}X_{2},gx_{1}x_{2}) \\
&=&B(x_{1}x_{2}\otimes x_{1};GX_{1},gx_{1}x_{2})
\end{eqnarray*}%
In view of the form of the elements we get%
\begin{equation}
B(x_{1}x_{2}\otimes 1_{H};G,gx_{1}x_{2})+B(x_{1}x_{2}\otimes
1_{H};GX_{1},gx_{2})=0  \label{X1,x1x2,X1F61,gx1x2}
\end{equation}

\paragraph{Equality $\left( \protect\ref{X1F71}\right) $}

rewrites as%
\begin{eqnarray*}
&&-\beta _{1}B(x_{1}x_{2}\otimes 1_{H};GX_{1}X_{2},f) \\
&=&-\beta _{1}B(gx_{1}x_{2}\otimes g;GX_{1}X_{2},f)+B(x_{1}x_{2}\otimes
x_{1};GX_{2},f)
\end{eqnarray*}

\subparagraph{Case $f=g$}

\begin{eqnarray*}
&&-\beta _{1}B(x_{1}x_{2}\otimes 1_{H};GX_{1}X_{2},g) \\
&=&-\beta _{1}B(gx_{1}x_{2}\otimes g;GX_{1}X_{2},g)+B(x_{1}x_{2}\otimes
x_{1};GX_{2},g)
\end{eqnarray*}%
By applying $\left( \ref{eq.10}\right) $ this rewrites as%
\begin{equation*}
B(x_{1}x_{2}\otimes x_{1};GX_{2},g)=0.
\end{equation*}%
In view of the form of the element, this is trivial.

\subparagraph{Case $f=x_{1}$}

\begin{equation*}
-\beta _{1}B(x_{1}x_{2}\otimes 1_{H};GX_{1}X_{2},x_{1})=-\beta
_{1}B(gx_{1}x_{2}\otimes g;GX_{1}X_{2},x_{1})+B(x_{1}x_{2}\otimes
x_{1};GX_{2},x_{1})
\end{equation*}%
By applying $\left( \ref{eq.10}\right) $ this rewrites as%
\begin{equation*}
+2\beta _{1}B(x_{1}x_{2}\otimes 1_{H};GX_{1}X_{2},x_{1})+B(x_{1}x_{2}\otimes
x_{1};GX_{2},x_{1})=0
\end{equation*}%
In view of the form of the elements we get

\begin{equation}
\beta _{1}B(x_{1}x_{2}\otimes 1_{H};GX_{1},x_{1}x_{2})=0
\label{X1,x1x2,X1F71,x1}
\end{equation}

\subparagraph{Case $f=x_{2}$}

\begin{equation*}
-\beta _{1}B(x_{1}x_{2}\otimes 1_{H};GX_{1}X_{2},x_{2})=-\beta
_{1}B(gx_{1}x_{2}\otimes g;GX_{1}X_{2},x_{2})+B(x_{1}x_{2}\otimes
x_{1};GX_{2},x_{2})
\end{equation*}%
By applying $\left( \ref{eq.10}\right) $ this rewrites as%
\begin{equation*}
+2\beta _{1}B(x_{1}x_{2}\otimes 1_{H};GX_{1}X_{2},x_{2})+B(x_{1}x_{2}\otimes
x_{1};GX_{2},x_{2})=0
\end{equation*}

In view of the form of the elements we get%
\begin{equation}
\beta _{1}B(x_{1}x_{2}\otimes 1_{H};GX_{2},x_{1}x_{2})=0
\label{X1,x1x2,X1F71,x2}
\end{equation}

\subparagraph{Case $f=gx_{1}x_{2}$}

\begin{eqnarray*}
&&-\beta _{1}B(x_{1}x_{2}\otimes 1_{H};GX_{1}X_{2},gx_{1}x_{2}) \\
&=&-\beta _{1}B(gx_{1}x_{2}\otimes
g;GX_{1}X_{2},gx_{1}x_{2})+B(x_{1}x_{2}\otimes x_{1};GX_{2},gx_{1}x_{2})
\end{eqnarray*}%
By applying $\left( \ref{eq.10}\right) $ this rewrites as%
\begin{equation*}
B(x_{1}x_{2}\otimes x_{1};GX_{2},gx_{1}x_{2})=0.
\end{equation*}%
In view of the form of the element we get

\begin{equation}
B(x_{1}x_{2}\otimes 1_{H};GX_{2},gx_{2})=0.  \label{X1,x1x2,X1F71,gx1x2}
\end{equation}

\paragraph{Equality $\left( \protect\ref{X1F81}\right) $}

rewrites as

\begin{equation*}
-B(x_{1}x_{2}\otimes 1_{H};GX_{2},f)=-B(gx_{1}x_{2}\otimes
g;GX_{2},f)+B(x_{1}x_{2}\otimes x_{1};GX_{1}X_{2},f).
\end{equation*}

\subparagraph{Case $f=1_{H}$}

\begin{equation*}
-B(x_{1}x_{2}\otimes 1_{H};GX_{2},1_{H})=-B(gx_{1}x_{2}\otimes
g;GX_{2},1_{H})+B(x_{1}x_{2}\otimes x_{1};GX_{1}X_{2},1_{H}).
\end{equation*}%
By applying $\left( \ref{eq.10}\right) $ this rewrites as%
\begin{equation*}
B(x_{1}x_{2}\otimes x_{1};GX_{1}X_{2},1_{H})=0
\end{equation*}%
which, in view of the form of the element, is trivial.

\paragraph{Case $f=x_{1}x_{2}$}

\begin{equation*}
-B(x_{1}x_{2}\otimes 1_{H};GX_{2},x_{1}x_{2})=-B(gx_{1}x_{2}\otimes
g;GX_{2},x_{1}x_{2})+B(x_{1}x_{2}\otimes x_{1};GX_{1}X_{2},x_{1}x_{2}).
\end{equation*}%
By applying $\left( \ref{eq.10}\right) $ this rewrites as%
\begin{equation*}
B(x_{1}x_{2}\otimes x_{1};GX_{1}X_{2},x_{1}x_{2})=0
\end{equation*}%
which, in view of the form of the element, is trivial.

\paragraph{Case $f=gx_{1}$}

\begin{equation*}
-B(x_{1}x_{2}\otimes 1_{H};GX_{2},gx_{1})=-B(gx_{1}x_{2}\otimes
g;GX_{2},gx_{1})+B(x_{1}x_{2}\otimes x_{1};GX_{1}X_{2},gx_{1}).
\end{equation*}%
By applying $\left( \ref{eq.10}\right) $ this rewrites as%
\begin{equation*}
2B(x_{1}x_{2}\otimes 1_{H};GX_{2},gx_{1})+B(x_{1}x_{2}\otimes
x_{1};GX_{1}X_{2},gx_{1})=0.
\end{equation*}%
In view of the form of the elements we get

\begin{equation}
+B(x_{1}x_{2}\otimes 1_{H};G,gx_{1}x_{2})+B(x_{1}x_{2}\otimes
1_{H};GX_{1},gx_{2})=0.  \label{X1,x1x2,X1F81,gx1}
\end{equation}

\paragraph{Case $f=gx_{2}$}

\begin{equation*}
-B(x_{1}x_{2}\otimes 1_{H};GX_{2},gx_{2})=-B(gx_{1}x_{2}\otimes
g;GX_{2},gx_{2})+B(x_{1}x_{2}\otimes x_{1};GX_{1}X_{2},gx_{2}).
\end{equation*}%
By applying $\left( \ref{eq.10}\right) $ this rewrites as%
\begin{equation*}
2B(x_{1}x_{2}\otimes 1_{H};GX_{2},gx_{2})+B(x_{1}x_{2}\otimes
x_{1};GX_{1}X_{2},gx_{2}).
\end{equation*}%
In view of the form of the elements we get%
\begin{equation*}
B(x_{1}x_{2}\otimes 1_{H};GX_{2},gx_{2})=0
\end{equation*}%
which is $\left( \ref{X1,x1x2,X1F71,gx1x2}\right) .$

\subsubsection{Case $gx_{1}\otimes 1_{H}$}

\paragraph{Equality $\left( \protect\ref{X1F11}\right) $}

rewrites as%
\begin{gather*}
\beta _{1}B(gx_{1}\otimes 1_{H};X_{1},f)+\lambda B(gx_{1}\otimes
1_{H};X_{2},f)=-\gamma _{1}B(x_{1}\otimes g;G,f) \\
-\beta _{1}B(x_{1}\otimes g;X_{1},f)+B(gx_{1}\otimes x_{1};1_{A},f)
\end{gather*}

\subparagraph{Case $f=1_{H}$}

\begin{gather*}
\beta _{1}B(gx_{1}\otimes 1_{H};X_{1},1_{H})+\lambda B(gx_{1}\otimes
1_{H};X_{2},1_{H})=-\gamma _{1}B(x_{1}\otimes g;G,1_{H}) \\
-\beta _{1}B(x_{1}\otimes g;X_{1},1_{H})+B(gx_{1}\otimes x_{1};1_{A},1_{H})
\end{gather*}%
By applying $\left( \ref{eq.10}\right) $ this rewrites as%
\begin{gather*}
2\beta _{1}B(gx_{1}\otimes 1_{H};X_{1},1_{H})+\lambda B(gx_{1}\otimes
1_{H};X_{2},1_{H})+\gamma _{1}B(gx_{1}\otimes 1_{H};G,1_{H}) \\
-B(gx_{1}\otimes x_{1};1_{A},1_{H})=0
\end{gather*}%
In view of the form of the elements we get%
\begin{gather}
2\beta _{1}\left[ -1+B(x_{1}x_{2}\otimes
1_{H};1_{A},x_{1}x_{2})+B(x_{1}x_{2}\otimes 1_{H};X_{2},x_{1})\right]
+\lambda B(x_{1}x_{2}\otimes 1_{H};X_{2},x_{2})  \label{X1,gx1,X1F11,1H} \\
=-\gamma _{1}\left[ B(x_{1}x_{2}\otimes 1_{H};G,x_{2})-B(x_{1}x_{2}\otimes
1_{H};GX_{2},1_{H})\right] +  \notag \\
-2\left[ B(x_{1}x_{2}\otimes 1_{H};1_{A},gx_{2})+B(x_{1}x_{2}\otimes
1_{H};X_{2},g)\right] .  \notag
\end{gather}

\subparagraph{Case $f=x_{1}x_{2}$}

\begin{gather*}
\beta _{1}B(gx_{1}\otimes 1_{H};X_{1},x_{1}x_{2})+\lambda B(gx_{1}\otimes
1_{H};X_{2},x_{1}x_{2})=-\gamma _{1}B(x_{1}\otimes g;G,x_{1}x_{2}) \\
-\beta _{1}B(x_{1}\otimes g;X_{1},x_{1}x_{2})+B(gx_{1}\otimes
x_{1};1_{A},x_{1}x_{2}).
\end{gather*}%
By applying $\left( \ref{eq.10}\right) $ this rewrites as%
\begin{gather*}
\beta _{1}B(gx_{1}\otimes 1_{H};X_{1},x_{1}x_{2})+\lambda B(gx_{1}\otimes
1_{H};X_{2},x_{1}x_{2})=-\gamma _{1}B(gx_{1}\otimes 1_{H};G,x_{1}x_{2}) \\
-\beta _{1}B(gx_{1}\otimes 1_{H};X_{1},x_{1}x_{2})+B(gx_{1}\otimes
x_{1};1_{A},x_{1}x_{2}).
\end{gather*}%
In view of the form of the elements we get%
\begin{equation*}
\gamma _{1}B(x_{1}x_{2}\otimes 1_{H};GX_{2},x_{1}x_{2})=0
\end{equation*}%
which is $\left( \ref{G,x1x2, GF7,x2}\right) .$

\subparagraph{Case $f=gx_{1}$}

\begin{gather*}
\beta _{1}B(gx_{1}\otimes 1_{H};X_{1},gx_{1})+\lambda B(gx_{1}\otimes
1_{H};X_{2},gx_{1})=-\gamma _{1}B(x_{1}\otimes g;G,gx_{1}) \\
-\beta _{1}B(x_{1}\otimes g;X_{1},gx_{1})+B(gx_{1}\otimes x_{1};1_{A},gx_{1})
\end{gather*}%
By applying $\left( \ref{eq.10}\right) $ this rewrites as%
\begin{equation*}
\lambda B(gx_{1}\otimes 1_{H};X_{2},gx_{1})-\gamma _{1}B(gx_{1}\otimes
1_{H};G,gx_{1})-B(gx_{1}\otimes x_{1};1_{A},gx_{1})=0.
\end{equation*}%
In view of the form of the elements we get%
\begin{gather}
\lambda B(x_{1}x_{2}\otimes 1_{H};X_{2},gx_{1}x_{2})-\gamma _{1}\left[
\begin{array}{c}
B(x_{1}x_{2}\otimes 1_{H};G,gx_{1}x_{2})+ \\
-B(x_{1}x_{2}\otimes 1_{H};GX_{2},gx_{1})%
\end{array}%
\right]  \label{X1,gx1,X1F11,gx1} \\
+2\left[ B(x_{1}x_{2}\otimes 1_{H};1_{A},gx_{2})+B(x_{1}x_{2}\otimes
1_{H};X_{2},g)\right] =0.  \notag
\end{gather}

\subparagraph{Case $f=gx_{2}$}

\begin{gather*}
\beta _{1}B(gx_{1}\otimes 1_{H};X_{1},gx_{2})+\lambda B(gx_{1}\otimes
1_{H};X_{2},gx_{2})=-\gamma _{1}B(x_{1}\otimes g;G,gx_{2}) \\
-\beta _{1}B(x_{1}\otimes g;X_{1},gx_{2})+B(gx_{1}\otimes x_{1};1_{A},gx_{2})
\end{gather*}%
By applying $\left( \ref{eq.10}\right) $ this rewrites as%
\begin{equation*}
\lambda B(gx_{1}\otimes 1_{H};X_{2},gx_{2})=\gamma _{1}B(gx_{1}\otimes
1_{H};G,gx_{2})+B(gx_{1}\otimes x_{1};1_{A},gx_{2}).
\end{equation*}%
In view of the form of the elements we get%
\begin{equation*}
\gamma _{1}B(x_{1}x_{2}\otimes 1_{H};GX_{2},gx_{2})=0
\end{equation*}

which is $\left( \ref{X1,x1x2,X1F41,gx2}\right) .$

\paragraph{Equality $\left( \protect\ref{X1F21}\right) $}

rewrites as%
\begin{eqnarray*}
&&\beta _{1}B(gx_{1}\otimes 1_{H};GX_{1},f)+\lambda B(gx_{1}\otimes
1_{H};GX_{2},f) \\
&=&\beta _{1}B(x_{1}\otimes g;GX_{1},f)+B(gx_{1}\otimes x_{1};G,f).
\end{eqnarray*}

\subparagraph{Case $f=g$}

\begin{eqnarray*}
&&\beta _{1}B(gx_{1}\otimes 1_{H};GX_{1},g)+\lambda B(gx_{1}\otimes
1_{H};GX_{2},g) \\
&=&\beta _{1}B(x_{1}\otimes g;GX_{1},g)+B(gx_{1}\otimes x_{1};G,g).
\end{eqnarray*}%
By applying $\left( \ref{eq.10}\right) $ this rewrites as%
\begin{equation*}
\lambda B(x_{1}x_{2}\otimes 1_{H};GX_{2},gx_{2})=0)
\end{equation*}%
which is $\left( \ref{X1,x1x2,X1F21,gx2}\right) .$

\subparagraph{Case $f=x_{1}$}

\begin{eqnarray*}
&&\beta _{1}B(gx_{1}\otimes 1_{H};GX_{1},x_{1})+\lambda B(gx_{1}\otimes
1_{H};GX_{2},x_{1}) \\
&=&\beta _{1}B(x_{1}\otimes g;GX_{1},x_{1})+B(gx_{1}\otimes x_{1};G,x_{1}).
\end{eqnarray*}%
By applying $\left( \ref{eq.10}\right) $ this rewrites as%
\begin{eqnarray*}
&&2\beta _{1}B(gx_{1}\otimes 1_{H};GX_{1},x_{1})+\lambda B(gx_{1}\otimes
1_{H};GX_{2},x_{1}) \\
&=&B(gx_{1}\otimes x_{1};G,x_{1}).
\end{eqnarray*}%
In view of the form of the elements we get%
\begin{equation*}
\lambda B(x_{1}x_{2}\otimes 1_{H};GX_{2},x_{1}x_{2})=0
\end{equation*}%
which is $\left( \ref{X1,x1x2,X1F61,x2}\right) .$

\subparagraph{Case $f=x_{2}$}

\begin{eqnarray*}
&&\beta _{1}B(gx_{1}\otimes 1_{H};GX_{1},x_{2})+\lambda B(gx_{1}\otimes
1_{H};GX_{2},x_{2}) \\
&=&\beta _{1}B(x_{1}\otimes g;GX_{1},x_{2})+B(gx_{1}\otimes x_{1};G,x_{2}).
\end{eqnarray*}%
By applying $\left( \ref{eq.10}\right) $ this rewrites as%
\begin{equation*}
\beta _{1}B(x_{1}x_{2}\otimes 1_{H};GX_{2},x_{1}x_{2})=0
\end{equation*}%
which is $\left( \ref{X1,x1x2,X1F71,x2}\right) .$

\subparagraph{Case $f=gx_{1}x_{2}$}

\begin{eqnarray*}
&&\beta _{1}B(gx_{1}\otimes 1_{H};GX_{1},gx_{1}x_{2})+\lambda
B(gx_{1}\otimes 1_{H};GX_{2},gx_{1}x_{2}) \\
&=&\beta _{1}B(x_{1}\otimes g;GX_{1},gx_{1}x_{2})+B(gx_{1}\otimes
x_{1};G,gx_{1}x_{2}).
\end{eqnarray*}%
By applying $\left( \ref{eq.10}\right) $ this rewrites as%
\begin{equation*}
+\lambda B(gx_{1}\otimes 1_{H};GX_{2},gx_{1}x_{2})-B(gx_{1}\otimes
x_{1};G,gx_{1}x_{2})=0.
\end{equation*}%
In view of the form of the elements we get%
\begin{equation*}
B(x_{1}x_{2}\otimes 1_{H};GX_{2},gx_{2})=0
\end{equation*}%
which is $\left( \ref{X1,x1x2,X1F71,gx1x2}\right) .$

\paragraph{Equality $\left( \protect\ref{X1F31}\right) $}

rewrites as

\begin{gather*}
B(gx_{1}\otimes 1_{H};1_{A},f)+\lambda B(gx_{1}\otimes 1_{H};X_{1}X_{2},f)=
\\
-B(x_{1}\otimes g;1_{A},f)-\gamma _{1}B(x_{1}\otimes
g;GX_{1},f)+B(gx_{1}\otimes x_{1};X_{1},f).
\end{gather*}

\subparagraph{Case $f=g$}

\begin{gather*}
B(gx_{1}\otimes 1_{H};1_{A},g)+\lambda B(gx_{1}\otimes 1_{H};X_{1}X_{2},g)=
\\
-B(x_{1}\otimes g;1_{A},g)-\gamma _{1}B(x_{1}\otimes
g;GX_{1},g)+B(gx_{1}\otimes x_{1};X_{1},g).
\end{gather*}%
By applying $\left( \ref{eq.10}\right) $ this rewrites as%
\begin{gather*}
2B(gx_{1}\otimes 1_{H};1_{A},g)+\lambda B(gx_{1}\otimes 1_{H};X_{1}X_{2},g)=
\\
-\gamma _{1}B(gx_{1}\otimes 1_{H};GX_{1},g)+B(gx_{1}\otimes x_{1};X_{1},g).
\end{gather*}%
In view of the form of the elements we get%
\begin{gather}
2\left[ B(x_{1}x_{2}\otimes 1_{H};1_{A},gx_{2})+B(x_{1}x_{2}\otimes
1_{H};X_{2},g)\right] +  \label{X1,gx1,X1F31,g} \\
+\lambda B(x_{1}x_{2}\otimes 1_{H};X_{2},gx_{1}x_{2})  \notag \\
+\gamma _{1}\left[ -B(x_{1}x_{2}\otimes
1_{H};G,gx_{1}x_{2})+B(x_{1}x_{2}\otimes 1_{H};GX_{2},gx_{1})\right] =0.
\notag
\end{gather}

\subparagraph{Case $f=x_{1}$}

\begin{gather*}
B(gx_{1}\otimes 1_{H};1_{A},x_{1})+\lambda B(gx_{1}\otimes
1_{H};X_{1}X_{2},x_{1})= \\
-B(x_{1}\otimes g;1_{A},x_{1})-\gamma _{1}B(x_{1}\otimes
g;GX_{1},x_{1})+B(gx_{1}\otimes x_{1};X_{1},x_{1}).
\end{gather*}%
By applying $\left( \ref{eq.10}\right) $ this rewrites as%
\begin{gather*}
+\lambda B(gx_{1}\otimes 1_{H};X_{1}X_{2},x_{1})= \\
-\gamma _{1}B(x_{1}\otimes g;GX_{1},x_{1})+B(gx_{1}\otimes
x_{1};X_{1},x_{1}).
\end{gather*}%
In view of the form of the elements this is trivial.

\subparagraph{Case $f=x_{2}$}

\begin{gather*}
B(gx_{1}\otimes 1_{H};1_{A},x_{2})+\lambda B(gx_{1}\otimes
1_{H};X_{1}X_{2},x_{2})= \\
-B(x_{1}\otimes g;1_{A},x_{2})-\gamma _{1}B(x_{1}\otimes
g;GX_{1},x_{2})+B(gx_{1}\otimes x_{1};X_{1},x_{2}).
\end{gather*}%
By applying $\left( \ref{eq.10}\right) $ this rewrites as%
\begin{gather*}
+\lambda B(gx_{1}\otimes 1_{H};X_{1}X_{2},x_{2})= \\
-\gamma _{1}B(x_{1}\otimes g;GX_{1},x_{2})+B(gx_{1}\otimes
x_{1};X_{1},x_{2}).
\end{gather*}%
In view of the form of the elements we get%
\begin{equation*}
\gamma _{1}B(x_{1}x_{2}\otimes 1_{H};GX_{2},x_{1}x_{2})=0
\end{equation*}%
which is $\left( \ref{G,x1x2, GF7,x2}\right) .$

\subparagraph{Case $f=gx_{1}$}

\begin{gather*}
B(gx_{1}\otimes 1_{H};1_{A},gx_{1})+\lambda B(gx_{1}\otimes
1_{H};X_{1}X_{2},gx_{1})= \\
-B(x_{1}\otimes g;1_{A},gx_{1})-\gamma _{1}B(x_{1}\otimes
g;GX_{1},gx_{1})+B(gx_{1}\otimes x_{1};X_{1},gx_{1}).
\end{gather*}%
By applying $\left( \ref{eq.10}\right) $ this rewrites as%
\begin{gather*}
+\lambda B(gx_{1}\otimes 1_{H};X_{1}X_{2},gx_{1})= \\
-\gamma _{1}B(x_{1}\otimes g;GX_{1},gx_{1})+B(gx_{1}\otimes
x_{1};X_{1},gx_{1}).
\end{gather*}%
In view of the form of the elements this is trivial.

\subparagraph{Case $f=gx_{1}x_{2}$}

\begin{gather*}
B(gx_{1}\otimes 1_{H};1_{A},gx_{1}x_{2})+\lambda B(gx_{1}\otimes
1_{H};X_{1}X_{2},gx_{1}x_{2})= \\
-B(x_{1}\otimes g;1_{A},gx_{1}x_{2})-\gamma _{1}B(x_{1}\otimes
g;GX_{1},gx_{1}x_{2})+B(gx_{1}\otimes x_{1};X_{1},gx_{1}x_{2}).
\end{gather*}%
By applying $\left( \ref{eq.10}\right) $ this rewrites as%
\begin{gather*}
2B(gx_{1}\otimes 1_{H};1_{A},gx_{1}x_{2})+\lambda B(gx_{1}\otimes
1_{H};X_{1}X_{2},gx_{1}x_{2})= \\
-\gamma _{1}B(x_{1}\otimes g;GX_{1},gx_{1}x_{2})+.
\end{gather*}%
In view of the form of the elements this is trivial.

\paragraph{Equality $\left( \protect\ref{X1F41}\right) $}

rewrites as%
\begin{eqnarray*}
&&-\beta _{1}B(gx_{1}\otimes 1_{H};X_{1}X_{2},f) \\
&=&-\beta _{1}B(x_{1}\otimes g;X_{1}X_{2},f)-\gamma _{1}B(x_{1}\otimes
g;GX_{2},f)+B(gx_{1}\otimes x_{1};X_{2},f).
\end{eqnarray*}

\subparagraph{Case $f=g$}

\begin{eqnarray*}
&&-\beta _{1}B(gx_{1}\otimes 1_{H};X_{1}X_{2},g) \\
&=&-\beta _{1}B(x_{1}\otimes g;X_{1}X_{2},g)-\gamma _{1}B(x_{1}\otimes
g;GX_{2},g)+B(gx_{1}\otimes x_{1};X_{2},g).
\end{eqnarray*}%
By applying $\left( \ref{eq.10}\right) $ this rewrites as%
\begin{equation*}
-\gamma _{1}B(gx_{1}\otimes 1_{H},GX_{2},g)+B(gx_{1}\otimes x_{1};X_{2},g)=0.
\end{equation*}%
In view of the form of the elements we get%
\begin{equation*}
\gamma _{1}B(x_{1}x_{2}\otimes 1_{H};GX_{2},gx_{2})=0
\end{equation*}%
which is $\left( \ref{X1,x1x2,X1F41,gx2}\right) .$

\subparagraph{Case $f=x_{1}$}

\begin{eqnarray*}
&&-\beta _{1}B(gx_{1}\otimes 1_{H};X_{1}X_{2},x_{1}) \\
&=&-\beta _{1}B(x_{1}\otimes g;X_{1}X_{2},x_{1})-\gamma _{1}B(x_{1}\otimes
g;GX_{2},x_{1})+B(gx_{1}\otimes x_{1};X_{2},x_{1}).
\end{eqnarray*}%
By applying $\left( \ref{eq.10}\right) $ this rewrites as

\subparagraph{%
\protect\begin{equation*}
-2\protect\beta _{1}B(gx_{1}\otimes 1_{H};X_{1}X_{2},x_{1})=\protect\gamma %
_{1}B(gx_{1}\otimes 1_{H};GX_{2},x_{1})+B(gx_{1}\otimes
x_{1};X_{2},x_{1})=0.
\protect\end{equation*}%
}

In view of the form of the elements we get%
\begin{equation*}
\gamma _{1}B(x_{1}x_{2}\otimes 1_{H};GX_{2},x_{1}x_{2})=0
\end{equation*}%
which is $\left( \ref{G,x1x2, GF7,x2}\right) .$

\paragraph{Equality $\left( \protect\ref{X1F51}\right) $}

rewrites as%
\begin{eqnarray*}
-B(gx_{1}\otimes 1_{H};X_{2},f) &=&-B(x_{1}\otimes g;X_{2},f)-\gamma
_{1}B(x_{1}\otimes g;GX_{1}X_{2},f) \\
&&+B(gx_{1}\otimes x_{1};X_{1}X_{2},f).
\end{eqnarray*}

\subparagraph{Case $f=1_{H}$}

\begin{eqnarray*}
-B(gx_{1}\otimes 1_{H};X_{2},1_{H}) &=&-B(x_{1}\otimes g;X_{2},1_{H})-\gamma
_{1}B(x_{1}\otimes g;GX_{1}X_{2},1_{H}) \\
&&+B(gx_{1}\otimes x_{1};X_{1}X_{2},1_{H}).
\end{eqnarray*}%
By applying $\left( \ref{eq.10}\right) $ this rewrites as%
\begin{equation*}
-\gamma _{1}B(gx_{1}\otimes 1_{H};GX_{1}X_{2},1_{H})+B(gx_{1}\otimes
x_{1};X_{1}X_{2},1_{H})=0.
\end{equation*}%
In view of the form of the elements we get%
\begin{equation*}
\gamma _{1}B(x_{1}x_{2}\otimes 1_{H};GX_{2},x_{1}x_{2})=0
\end{equation*}%
which is $\left( \ref{G,x1x2, GF7,x2}\right) .$

\subparagraph{Case $f=gx_{1}$}

\begin{eqnarray*}
-B(gx_{1}\otimes 1_{H};X_{2},gx_{1}) &=&-B(x_{1}\otimes
g;X_{2},gx_{1})-\gamma _{1}B(x_{1}\otimes g;GX_{1}X_{2},gx_{1}) \\
&&+B(gx_{1}\otimes x_{1};X_{1}X_{2},gx_{1}).
\end{eqnarray*}%
By applying $\left( \ref{eq.10}\right) $ this rewrites as%
\begin{equation*}
2B(gx_{1}\otimes 1_{H};X_{2},gx_{1})+\gamma _{1}B(gx_{1}\otimes
1_{H};GX_{1}X_{2},gx_{1})+B(gx_{1}\otimes x_{1};X_{1}X_{2},gx_{1})=0.
\end{equation*}%
In view of the form of the elements we get%
\begin{equation*}
2B(x_{1}x_{2}\otimes 1_{H};X_{2},gx_{1}x_{2})-2B(x_{1}x_{2}\otimes
1_{H};X_{2},gx_{1}x_{2})=0
\end{equation*}%
which is trivial.

\paragraph{Equality $\left( \protect\ref{X1F61}\right) $}

rewrites as%
\begin{eqnarray*}
&&B(gx_{1}\otimes 1_{H};G,f)+\lambda B(gx_{1}\otimes 1_{H};GX_{1}X_{2},f) \\
&=&B(x_{1}\otimes g;G,f)+B(gx_{1}\otimes x_{1};GX_{1},f).
\end{eqnarray*}

\subparagraph{Case $f=1_{H}$}

\begin{eqnarray*}
&&B(gx_{1}\otimes 1_{H};G,1_{H})+\lambda B(gx_{1}\otimes
1_{H};GX_{1}X_{2},1_{H}) \\
&=&B(x_{1}\otimes g;G,1_{H})+B(gx_{1}\otimes x_{1};GX_{1},1_{H}).
\end{eqnarray*}%
By applying $\left( \ref{eq.10}\right) $ this rewrites as%
\begin{equation*}
-\lambda B(gx_{1}\otimes 1_{H};GX_{1}X_{2},1_{H})+B(gx_{1}\otimes
x_{1};GX_{1},1_{H})=0.
\end{equation*}%
In view of the form of the elements we get%
\begin{equation*}
\lambda B(x_{1}x_{2}\otimes 1_{H};GX_{2},x_{1}x_{2})=0
\end{equation*}%
which is $\left( \ref{X1,x1x2,X1F61,x2}\right) .$

\subparagraph{Case $f=x_{1}x_{2}$}

\begin{eqnarray*}
&&B(gx_{1}\otimes 1_{H};G,x_{1}x_{2})+\lambda B(gx_{1}\otimes
1_{H};GX_{1}X_{2},x_{1}x_{2}) \\
&=&B(x_{1}\otimes g;G,x_{1}x_{2})+B(gx_{1}\otimes x_{1};GX_{1},x_{1}x_{2}).
\end{eqnarray*}%
By applying $\left( \ref{eq.10}\right) $ this rewrites as%
\begin{equation*}
-\lambda B(gx_{1}\otimes 1_{H};GX_{1}X_{2},x_{1}x_{2})+B(gx_{1}\otimes
x_{1};GX_{1},x_{1}x_{2})=0.
\end{equation*}%
In view of the form of the elements this is trivial.

\subparagraph{Case $f=gx_{1}$}

\begin{eqnarray*}
&&B(gx_{1}\otimes 1_{H};G,gx_{1})+\lambda B(gx_{1}\otimes
1_{H};GX_{1}X_{2},gx_{1}) \\
&=&B(x_{1}\otimes g;G,gx_{1})+B(gx_{1}\otimes x_{1};GX_{1},gx_{1}).
\end{eqnarray*}%
By applying $\left( \ref{eq.10}\right) $ this rewrites as%
\begin{eqnarray*}
&&2B(gx_{1}\otimes 1_{H};G,gx_{1})+\lambda B(gx_{1}\otimes
1_{H};GX_{1}X_{2},gx_{1}) \\
&=&+B(gx_{1}\otimes x_{1};GX_{1},gx_{1}).
\end{eqnarray*}%
In view of the form of the elements we get%
\begin{gather*}
-\left[ B(x_{1}x_{2}\otimes 1_{H};G,gx_{1}x_{2})-B(x_{1}x_{2}\otimes
1_{H};GX_{2},gx_{1})\right] \\
+\left[ +B(x_{1}x_{2}\otimes 1_{H};G,gx_{1}x_{2})-B(x_{1}x_{2}\otimes
1_{H};GX_{2},gx_{1})\right] =0
\end{gather*}%
which is trivial.

\subparagraph{Case $f=gx_{2}$}

\begin{eqnarray*}
&&B(gx_{1}\otimes 1_{H};G,gx_{2})+\lambda B(gx_{1}\otimes
1_{H};GX_{1}X_{2},gx_{2}) \\
&=&B(x_{1}\otimes g;G,gx_{2})+B(gx_{1}\otimes x_{1};GX_{1},gx_{2}).
\end{eqnarray*}%
By applying $\left( \ref{eq.10}\right) $ this rewrites as%
\begin{eqnarray*}
&&2B(gx_{1}\otimes 1_{H};G,gx_{2})+\lambda B(gx_{1}\otimes
1_{H};GX_{1}X_{2},gx_{2}) \\
&=&+B(gx_{1}\otimes x_{1};GX_{1},gx_{2}).
\end{eqnarray*}%
In view of the form of the elements we get%
\begin{equation*}
B(x_{1}x_{2}\otimes 1_{H};GX_{2},gx_{2})=0
\end{equation*}%
which is $\left( \ref{X1,x1x2,X1F71,gx1x2}\right) .$

\paragraph{Equality $\left( \protect\ref{X1F71}\right) $}

rewrites as%
\begin{gather*}
-\beta _{1}B(gx_{1}\otimes 1_{H};GX_{1}X_{2},f)=+\beta _{1}B(x_{1}\otimes
g;GX_{1}X_{2},f) \\
+B(gx_{1}\otimes x_{1};GX_{2},f).
\end{gather*}%
Case $f=1_{H}$

\begin{gather*}
-\beta _{1}B(gx_{1}\otimes 1_{H};GX_{1}X_{2},1_{H})=+\beta
_{1}B(x_{1}\otimes g;GX_{1}X_{2},1_{H}) \\
+B(gx_{1}\otimes x_{1};GX_{2},1_{H}).
\end{gather*}%
By applying $\left( \ref{eq.10}\right) $ this rewrites as%
\begin{equation*}
2\beta _{1}B(gx_{1}\otimes 1_{H};GX_{1}X_{2},1_{H})+B(gx_{1}\otimes
x_{1};GX_{2},1_{H})=0.
\end{equation*}%
In view of the form of the elements we get

\begin{equation*}
\beta _{1}B(x_{1}x_{2}\otimes 1_{H};GX_{2},x_{1}x_{2})=0
\end{equation*}%
which is $\left( \ref{X1,x1x2,X1F71,x2}\right) .$

\subparagraph{Case $f=gx_{1}$}

\begin{gather*}
-\beta _{1}B(gx_{1}\otimes 1_{H};GX_{1}X_{2},gx_{1})=+\beta
_{1}B(x_{1}\otimes g;GX_{1}X_{2},gx_{1}) \\
+B(gx_{1}\otimes x_{1};GX_{2},gx_{1}).
\end{gather*}%
By applying $\left( \ref{eq.10}\right) $ this rewrites as%
\begin{equation*}
B(gx_{1}\otimes x_{1};GX_{2},gx_{1})=0
\end{equation*}%
In view of the form of the element we get%
\begin{equation*}
B(x_{1}x_{2}\otimes 1_{H};GX_{2},gx_{2})=0
\end{equation*}%
which is $\left( \ref{X1,x1x2,X1F71,gx1x2}\right) .$

\paragraph{Equality $\left( \protect\ref{X1F81}\right) $}

rewrites as%
\begin{gather*}
-B(gx_{1}\otimes 1_{H};GX_{2},f)=+B(x_{1}\otimes g;GX_{2},f) \\
+B(gx_{1}\otimes x_{1};GX_{1}X_{2},f)
\end{gather*}

\subparagraph{Case $f=g$}

\begin{gather*}
-B(gx_{1}\otimes 1_{H};GX_{2},g)=+B(x_{1}\otimes g;GX_{2},g) \\
+B(gx_{1}\otimes x_{1};GX_{1}X_{2},g)
\end{gather*}%
By applying $\left( \ref{eq.10}\right) $ this rewrites as%
\begin{equation*}
2B(gx_{1}\otimes 1_{H};GX_{2},g)+B(gx_{1}\otimes x_{1};GX_{1}X_{2},g)=0.
\end{equation*}%
In view of the form of the elements we get%
\begin{equation*}
B(x_{1}x_{2}\otimes 1_{H};GX_{2},gx_{2})=0
\end{equation*}%
which is $\left( \ref{X1,x1x2,X1F71,gx1x2}\right) .$

\subparagraph{Case $f=$ $x_{1}$}

\begin{gather*}
-B(gx_{1}\otimes 1_{H};GX_{2},x_{1})=+B(x_{1}\otimes g;GX_{2},x_{1}) \\
+B(gx_{1}\otimes x_{1};GX_{1}X_{2},x_{1}).
\end{gather*}%
By applying $\left( \ref{eq.10}\right) $ this rewrites as%
\begin{equation*}
B(gx_{1}\otimes x_{1};GX_{1}X_{2},x_{1})=0
\end{equation*}%
which is trivial in view of the form of the element.

\subsubsection{Case $gx_{2}\otimes 1_{H}$}

\paragraph{Equality $\left( \protect\ref{X1F11}\right) $}

rewrites as%
\begin{gather*}
\beta _{1}B(gx_{2}\otimes 1_{H};X_{1},f)+\lambda B(gx_{2}\otimes
1_{H};X_{2},f)=-\gamma _{1}B(x_{2}\otimes g;G,f) \\
-\beta _{1}B(x_{2}\otimes g;X_{1},f)-B(gx_{1}x_{2}\otimes
g;1_{A},f)+B(gx_{2}\otimes x_{1};1_{A},f).
\end{gather*}

\subparagraph{Case $f=1_{H}$}

\begin{gather*}
\beta _{1}B(gx_{2}\otimes 1_{H};X_{1},1_{H})+\lambda B(gx_{2}\otimes
1_{H};X_{2},1_{H})=-\gamma _{1}B(x_{2}\otimes g;G,1_{H}) \\
-\beta _{1}B(x_{2}\otimes g;X_{1},1_{H})-B(gx_{1}x_{2}\otimes
g;1_{A},1_{H})+B(gx_{2}\otimes x_{1};1_{A},1_{H}).
\end{gather*}%
By applying $\left( \ref{eq.10}\right) $ this rewrites as%
\begin{gather*}
2\beta _{1}B(gx_{2}\otimes 1_{H};X_{1},1_{H})+\lambda B(gx_{2}\otimes
1_{H};X_{2},1_{H})=-\gamma _{1}B(gx_{2}\otimes 1_{H};G,1_{H}) \\
-B(x_{1}x_{2}\otimes 1_{H};1_{A},1_{H})+B(gx_{2}\otimes x_{1};1_{A},1_{H}).
\end{gather*}%
In view of the form of the elements we get%
\begin{gather}
-2\beta _{1}B(x_{1}x_{2}\otimes 1_{H};X_{1},x_{1})+\lambda \left[
-1+B(x_{1}x_{2}\otimes 1_{H};1_{A},x_{1}x_{2})-B(x_{1}x_{2}\otimes
1_{H};X_{1},x_{2})\right]  \label{X1,gx2,X1F11,1H} \\
+\gamma _{1}\left[ -B(x_{1}x_{2}\otimes 1_{H};G,x_{1})+B(x_{1}x_{2}\otimes
1_{H};GX_{1},1_{H})\right] =0.  \notag
\end{gather}

\subparagraph{Case $f=g$}

\begin{gather*}
\beta _{1}B(gx_{2}\otimes 1_{H};X_{1},g)+\lambda B(gx_{2}\otimes
1_{H};X_{2},g)=-\gamma _{1}B(x_{2}\otimes g;G,g) \\
-\beta _{1}B(x_{2}\otimes g;X_{1},g)-B(gx_{1}x_{2}\otimes
g;1_{A},g)+B(gx_{2}\otimes x_{1};1_{A},g).
\end{gather*}%
By applying $\left( \ref{eq.10}\right) $ this rewrites as%
\begin{gather*}
\beta _{1}B(gx_{2}\otimes 1_{H};X_{1},g)+\lambda B(gx_{2}\otimes
1_{H};X_{2},g)=-\gamma _{1}B(gx_{2}\otimes 1_{H};G,g) \\
-\beta _{1}B(gx_{2}\otimes 1_{H};X_{1},g)-B(x_{1}x_{2}\otimes
1_{H};1_{A},g)+B(gx_{2}\otimes x_{1};1_{A},g).
\end{gather*}%
In view of the form of the elements we get a trivial equality.

\subparagraph{Case $f=x_{1}x_{2}$}

\begin{gather*}
\beta _{1}B(gx_{2}\otimes 1_{H};X_{1},x_{1}x_{2})+\lambda B(gx_{2}\otimes
1_{H};X_{2},x_{1}x_{2})=-\gamma _{1}B(x_{2}\otimes g;G,x_{1}x_{2}) \\
-\beta _{1}B(x_{2}\otimes g;X_{1},x_{1}x_{2})-B(gx_{1}x_{2}\otimes
g;1_{A},x_{1}x_{2})+B(gx_{2}\otimes x_{1};1_{A},x_{1}x_{2}).
\end{gather*}%
By applying $\left( \ref{eq.10}\right) $ this rewrites as%
\begin{gather*}
2\beta _{1}B(gx_{2}\otimes 1_{H};X_{1},x_{1}x_{2})+\lambda B(gx_{2}\otimes
1_{H};X_{2},x_{1}x_{2})=-\gamma _{1}B(gx_{2}\otimes 1_{H};G,x_{1}x_{2}) \\
-B(x_{1}x_{2}\otimes 1_{H};1_{A},x_{1}x_{2})+B(gx_{2}\otimes
x_{1};1_{A},x_{1}x_{2}).
\end{gather*}%
In view of the form of the elements we get%
\begin{equation*}
\gamma _{1}B(x_{1}x_{2}\otimes 1_{H};GX_{1},x_{1}x_{2})=0
\end{equation*}%
which is $\left( \ref{G,x1x2, GF7,x1}\right) .$

\subparagraph{Case $f=gx_{1}$}

\begin{gather*}
\beta _{1}B(gx_{2}\otimes 1_{H};X_{1},gx_{1})+\lambda B(gx_{2}\otimes
1_{H};X_{2},gx_{1})=-\gamma _{1}B(x_{2}\otimes g;G,gx_{1}) \\
-\beta _{1}B(x_{2}\otimes g;X_{1},gx_{1})-B(gx_{1}x_{2}\otimes
g;1_{A},gx_{1})+B(gx_{2}\otimes x_{1};1_{A},gx_{1}).
\end{gather*}%
By applying $\left( \ref{eq.10}\right) $ this rewrites as%
\begin{gather*}
+\lambda B(gx_{2}\otimes 1_{H};X_{2},gx_{1})=+\gamma _{1}B(gx_{2}\otimes
1_{H};G,gx_{1}) \\
+B(x_{1}x_{2}\otimes 1_{H};1_{A},gx_{1})+B(gx_{2}\otimes x_{1};1_{A},gx_{1}).
\end{gather*}%
In view of the form of the elements we get%
\begin{gather}
\lambda B(x_{1}x_{2}\otimes 1_{H};X_{1},gx_{1}x_{2})+\gamma
_{1}B(x_{1}x_{2}\otimes 1_{H};GX_{1},gx_{1})  \label{X1,gx2,X1F11,gx1} \\
+2B(x_{1}x_{2}\otimes 1_{H};1_{A},gx_{1})+2B(x_{1}x_{2}\otimes
1_{H};X_{1},g)=0.  \notag
\end{gather}

\subparagraph{Case $f=gx_{2}$}

\begin{gather*}
\beta _{1}B(gx_{2}\otimes 1_{H};X_{1},gx_{2})+\lambda B(gx_{2}\otimes
1_{H};X_{2},gx_{2})=-\gamma _{1}B(x_{2}\otimes g;G,gx_{2}) \\
-\beta _{1}B(x_{2}\otimes g;X_{1},gx_{2})-B(gx_{1}x_{2}\otimes
g;1_{A},gx_{2})+B(gx_{2}\otimes x_{1};1_{A},gx_{2}).
\end{gather*}%
By applying $\left( \ref{eq.10}\right) $ this rewrites as%
\begin{gather*}
+\lambda B(gx_{2}\otimes 1_{H};X_{2},gx_{2})=+\gamma _{1}B(gx_{2}\otimes
1_{H};G,gx_{2}) \\
+B(x_{1}x_{2}\otimes 1_{H};1_{A},gx_{2})+B(gx_{2}\otimes x_{1};1_{A},gx_{2}).
\end{gather*}%
In view of the form of the elements we get%
\begin{equation}
+\gamma _{1}\left[ B(x_{1}x_{2}\otimes
1_{H};G,gx_{1}x_{2})+B(x_{1}x_{2}\otimes 1_{H};GX_{1},gx_{2})\right] =0
\label{X1,gx2,X1F11,gx2}
\end{equation}

\paragraph{Equality $\left( \protect\ref{X1F21}\right) $}

rewrites as%
\begin{eqnarray*}
&&\beta _{1}B(gx_{2}\otimes 1_{H};GX_{1},f)+\lambda B(gx_{2}\otimes
1_{H};GX_{2},f) \\
&=&+\beta _{1}B(x_{2}\otimes g;GX_{1},f)-B(gx_{1}x_{2}\otimes
g;G,f)+B(gx_{2}\otimes x_{1};G,f)
\end{eqnarray*}

\subparagraph{Case $f=g$}

\begin{eqnarray*}
&&\beta _{1}B(gx_{2}\otimes 1_{H};GX_{1},g)+\lambda B(gx_{2}\otimes
1_{H};GX_{2},g) \\
&=&\beta _{1}B(x_{2}\otimes g;GX_{1},g)-B(gx_{1}x_{2}\otimes
g;G,g)+B(gx_{2}\otimes x_{1};G,g)
\end{eqnarray*}%
By applying $\left( \ref{eq.10}\right) $ this rewrites as%
\begin{eqnarray*}
&&+\lambda B(gx_{2}\otimes 1_{H};GX_{2},g) \\
&=&-B(x_{1}x_{2}\otimes 1_{H};G,g)+B(gx_{2}\otimes x_{1};G,g).
\end{eqnarray*}%
In view of the form of the elements we get

\begin{equation*}
\lambda \left[ B(x_{1}x_{2}\otimes 1_{H};G,gx_{1}x_{2})+B(x_{1}x_{2}\otimes
1_{H};GX_{1},gx_{2})\right] =0
\end{equation*}%
which follows from $\left( \ref{X1,x1x2,X1F81,gx1}\right) .$

\subparagraph{Case $f=x_{1}$}

\begin{eqnarray*}
&&\beta _{1}B(gx_{2}\otimes 1_{H};GX_{1},x_{1})+\lambda B(gx_{2}\otimes
1_{H};GX_{2},x_{1}) \\
&=&\beta _{1}B(x_{2}\otimes g;GX_{1},x_{1})-B(gx_{1}x_{2}\otimes
g;G,x_{1})+B(gx_{2}\otimes x_{1};G,x_{1})
\end{eqnarray*}%
By applying $\left( \ref{eq.10}\right) $ this rewrites as%
\begin{eqnarray*}
&&2\beta _{1}B(gx_{2}\otimes 1_{H};GX_{1},x_{1})+\lambda B(gx_{2}\otimes
1_{H};GX_{2},x_{1}) \\
&=&B(x_{1}x_{2}\otimes 1_{H};G,x_{1})+B(gx_{2}\otimes x_{1};G,x_{1})
\end{eqnarray*}%
In view of the form of the elements we get%
\begin{equation*}
\lambda B(x_{1}x_{2}\otimes 1_{H};GX_{1},x_{1}x_{2})=0
\end{equation*}%
which is $\left( \ref{X1,x1x2,X1F61,x1}\right) .$

\subparagraph{Case $f=x_{2}$}

\begin{eqnarray*}
&&\beta _{1}B(gx_{2}\otimes 1_{H};GX_{1},x_{2})+\lambda B(gx_{2}\otimes
1_{H};GX_{2},x_{2}) \\
&=&\beta _{1}B(x_{2}\otimes g;GX_{1},x_{2})-B(gx_{1}x_{2}\otimes
g;G,x_{2})+B(gx_{2}\otimes x_{1};G,x_{2})
\end{eqnarray*}%
By applying $\left( \ref{eq.10}\right) $ this rewrites as%
\begin{eqnarray*}
&&2\beta _{1}B(gx_{2}\otimes 1_{H};GX_{1},x_{2})+\lambda B(gx_{2}\otimes
1_{H};GX_{2},x_{2}) \\
&=&+B(x_{1}x_{2}\otimes 1_{H};G,x_{2})+B(gx_{2}\otimes x_{1};G,x_{2})
\end{eqnarray*}%
In view of the form of the elements we get%
\begin{equation*}
\beta _{1}B(x_{1}x_{2}\otimes 1_{H};GX_{1},x_{1}x_{2})=0
\end{equation*}%
which is $\left( \ref{X1,x1x2,X1F71,x1}\right) .$

\subparagraph{Case $f=gx_{1}x_{2}$}

\begin{eqnarray*}
&&\beta _{1}B(gx_{2}\otimes 1_{H};GX_{1},gx_{1}x_{2})+\lambda
B(gx_{2}\otimes 1_{H};GX_{2},gx_{1}x_{2}) \\
&=&\beta _{1}B(x_{2}\otimes g;GX_{1},gx_{1}x_{2})-B(gx_{1}x_{2}\otimes
g;G,gx_{1}x_{2})+B(gx_{2}\otimes x_{1};G,gx_{1}x_{2})
\end{eqnarray*}%
By applying $\left( \ref{eq.10}\right) $ this rewrites as%
\begin{equation*}
-\lambda B(gx_{2}\otimes 1_{H};GX_{2},gx_{1}x_{2})-B(x_{1}x_{2}\otimes
1_{H};G,gx_{1}x_{2})+B(gx_{2}\otimes x_{1};G,gx_{1}x_{2})=0
\end{equation*}%
In view of the form of the elements we get%
\begin{equation*}
B(x_{1}x_{2}\otimes 1_{H};G,gx_{1}x_{2})+B(x_{1}x_{2}\otimes
1_{H};GX_{1},gx_{2})=0
\end{equation*}%
which is $\left( \ref{X1,x1x2,X1F81,gx1}\right) .$

\paragraph{Equality $\left( \protect\ref{X1F31}\right) $}

rewrites as%
\begin{gather*}
B(gx_{2}\otimes 1_{H};1_{A},f)+\lambda B(gx_{2}\otimes 1_{H};X_{1}X_{2},f)=
\\
-B(x_{2}\otimes g;1_{A},f)-\gamma _{1}B(x_{2}\otimes
g;GX_{1},f)-B(gx_{1}x_{2}\otimes g;X_{1},f)+B(gx_{2}\otimes x_{1};X_{1},f)
\end{gather*}

\subparagraph{Case $f=g$}

\begin{gather*}
B(gx_{2}\otimes 1_{H};1_{A},g)+\lambda B(gx_{2}\otimes 1_{H};X_{1}X_{2},g)=
\\
-B(x_{2}\otimes g;1_{A},g)-\gamma _{1}B(x_{2}\otimes
g;GX_{1},g)-B(gx_{1}x_{2}\otimes g;X_{1},g)+B(gx_{2}\otimes x_{1};X_{1},g)
\end{gather*}%
By applying $\left( \ref{eq.10}\right) $ this rewrites as

\begin{gather*}
2B(gx_{2}\otimes 1_{H};1_{A},g)+\lambda B(gx_{2}\otimes 1_{H};X_{1}X_{2},g)=
\\
-\gamma _{1}B(gx_{2}\otimes 1_{H};GX_{1},g)-B(x_{1}x_{2}\otimes
1_{H};X_{1},g)+B(gx_{2}\otimes x_{1};X_{1},g)
\end{gather*}%
In view of the form of the elements we get%
\begin{gather}
2\left[ B(x_{1}x_{2}\otimes 1_{H};1_{A},gx_{1})+B(x_{1}x_{2}\otimes
1_{H};X_{1},g)\right]  \label{X1,gx2,X131,g} \\
+\lambda B(x_{1}x_{2}\otimes 1_{H};X_{1},gx_{1}x_{2})+\gamma
_{1}B(x_{1}x_{2}\otimes 1_{H};GX_{1},gx_{1})=0  \notag
\end{gather}

\subparagraph{Case $f=x_{1}$}

\begin{gather*}
B(gx_{2}\otimes 1_{H};1_{A},x_{1})+\lambda B(gx_{2}\otimes
1_{H};X_{1}X_{2},x_{1})= \\
-B(x_{2}\otimes g;1_{A},x_{1})-\gamma _{1}B(x_{2}\otimes g;GX_{1},x_{1})+ \\
-B(gx_{1}x_{2}\otimes g;X_{1},x_{1})+B(gx_{2}\otimes x_{1};X_{1},x_{1})
\end{gather*}%
By applying $\left( \ref{eq.10}\right) $ this rewrites and the form of the
elements we get%
\begin{equation*}
B(x_{1}x_{2}\otimes 1_{H};X_{1},x_{1})-B(x_{1}x_{2}\otimes
1_{H};X_{1},x_{1})=0
\end{equation*}%
which is trivial.

\subparagraph{Case $f=x_{2}$}

\begin{gather*}
B(gx_{2}\otimes 1_{H};1_{A},x_{2})+\lambda B(gx_{2}\otimes
1_{H};X_{1}X_{2},x_{2})= \\
-B(x_{2}\otimes g;1_{A},x_{2})-\gamma _{1}B(x_{2}\otimes g;GX_{1},x_{2})+ \\
-B(gx_{1}x_{2}\otimes g;X_{1},x_{2})+B(gx_{2}\otimes x_{1};X_{1},x_{2})
\end{gather*}%
By applying $\left( \ref{eq.10}\right) $ this rewrites as%
\begin{gather*}
\lambda B(gx_{2}\otimes 1_{H};X_{1}X_{2},x_{2})= \\
+\gamma _{1}B(gx_{2}\otimes 1_{H};GX_{1},x_{2})+B(x_{1}x_{2}\otimes
1_{H};X_{1},x_{2})+B(gx_{2}\otimes x_{1};X_{1},x_{2})
\end{gather*}%
In view of the form of the elements we get%
\begin{equation*}
\gamma _{1}B(x_{1}x_{2}\otimes 1_{H};GX_{1},x_{1}x_{2})=0
\end{equation*}%
which is $\left( \ref{G,x1x2, GF7,x1}\right) .$

\subparagraph{Case $f=gx_{1}x_{2}$}

\begin{gather*}
B(gx_{2}\otimes 1_{H};1_{A},gx_{1}x_{2})+\lambda B(gx_{2}\otimes
1_{H};X_{1}X_{2},gx_{1}x_{2})= \\
-B(x_{2}\otimes g;1_{A},gx_{1}x_{2})-\gamma _{1}B(x_{2}\otimes
g;GX_{1},gx_{1}x_{2})+ \\
-B(gx_{1}x_{2}\otimes g;X_{1},gx_{1}x_{2})+B(gx_{2}\otimes
x_{1};X_{1},gx_{1}x_{2})
\end{gather*}%
By applying $\left( \ref{eq.10}\right) $ this rewrites as%
\begin{gather*}
2B(gx_{2}\otimes 1_{H};1_{A},gx_{1}x_{2})+\lambda B(gx_{2}\otimes
1_{H};X_{1}X_{2},gx_{1}x_{2})= \\
-\gamma _{1}B(gx_{2}\otimes 1_{H};GX_{1},gx_{1}x_{2})-B(x_{1}x_{2}\otimes
1_{H};X_{1},gx_{1}x_{2})+B(gx_{2}\otimes x_{1};X_{1},gx_{1}x_{2})
\end{gather*}%
In view of the form of the elements we get

\begin{equation*}
B(x_{1}x_{2}\otimes 1_{H};X_{1},gx_{1}x_{2})=0
\end{equation*}%
which is $\left( \ref{G,x1x2, GF6,gx1x2}\right) .$

\paragraph{Equality $\left( \protect\ref{X1F41}\right) $}

rewrites as%
\begin{eqnarray*}
&&-\beta _{1}B(gx_{2}\otimes 1_{H};X_{1}X_{2},f) \\
&=&-\beta _{1}B(x_{2}\otimes g;X_{1}X_{2},f)-\gamma _{1}B(x_{2}\otimes
g;GX_{2},f)+ \\
&&-B(gx_{1}x_{2}\otimes g;X_{2},f)+B(gx_{2}\otimes x_{1};X_{2},f)
\end{eqnarray*}

\subparagraph{Case $f=g$}

By applying $\left( \ref{eq.10}\right) $ this rewrites as%
\begin{equation*}
-\gamma _{1}B(gx_{2}\otimes 1_{H};GX_{2},g)-B(x_{1}x_{2}\otimes
1_{H};X_{2},g)-B(gx_{2}\otimes x_{1};X_{2},g)=0
\end{equation*}%
In view of the form of the elements we get%
\begin{equation}
\gamma _{1}\left[ -B(x_{1}x_{2}\otimes
1_{H};G,gx_{1}x_{2})+B(x_{1}x_{2}\otimes 1_{H};GX_{1},gx_{2})\right] =0
\label{X1,gx2,X141,g}
\end{equation}

\subparagraph{Case $f=x_{1}$}

\begin{eqnarray*}
&&-\beta _{1}B(gx_{2}\otimes 1_{H};X_{1}X_{2},x_{1}) \\
&=&-\beta _{1}B(x_{2}\otimes g;X_{1}X_{2},x_{1})-\gamma _{1}B(x_{2}\otimes
g;GX_{2},x_{1})+ \\
&&-B(gx_{1}x_{2}\otimes g;X_{2},x_{1})+B(gx_{2}\otimes x_{1};X_{2},x_{1}).
\end{eqnarray*}%
By applying $\left( \ref{eq.10}\right) $ this rewrites as%
\begin{gather*}
2\beta _{1}B(gx_{2}\otimes 1_{H};X_{1}X_{2},x_{1})+\gamma
_{1}B(gx_{2}\otimes 1_{H};GX_{2},x_{1}) \\
+B(x_{1}x_{2}\otimes 1_{H};X_{2},x_{1})+B(gx_{2}\otimes x_{1};X_{2},x_{1})=0.
\end{gather*}%
In view of the form of the elements we get%
\begin{equation*}
\gamma _{1}B(x_{1}x_{2}\otimes 1_{H};GX_{1},x_{1}x_{2})=0
\end{equation*}%
which is $\left( \ref{G,x1x2, GF7,x1}\right) .$

\subparagraph{Case $f=x_{2}$}

\begin{eqnarray*}
&&-\beta _{1}B(gx_{2}\otimes 1_{H};X_{1}X_{2},x_{2}) \\
&=&-\beta _{1}B(x_{2}\otimes g;X_{1}X_{2},x_{2})-\gamma _{1}B(x_{2}\otimes
g;GX_{2},x_{2})+ \\
&&-B(gx_{1}x_{2}\otimes g;X_{2},x_{2})+B(gx_{2}\otimes x_{1};X_{2},x_{2})
\end{eqnarray*}%
By applying $\left( \ref{eq.10}\right) $ this rewrites as%
\begin{eqnarray*}
&&2\beta _{1}B(gx_{2}\otimes 1_{H};X_{1}X_{2},x_{2})+\gamma
_{1}B(gx_{2}\otimes 1_{H};GX_{2},x_{2}) \\
&&+B(x_{1}x_{2}\otimes 1_{H};X_{2},x_{2})+B(gx_{2}\otimes x_{1};X_{2},x_{2}).
\end{eqnarray*}%
In view of the form of the elements we get%
\begin{equation*}
B(x_{1}x_{2}\otimes 1_{H};X_{2},x_{2})-B(x_{1}x_{2}\otimes
1_{H};X_{2},x_{2})=0
\end{equation*}%
which is trivial.

\subparagraph{Case $f=gx_{1}x_{2}$}

\begin{eqnarray*}
&&-\beta _{1}B(gx_{2}\otimes 1_{H};X_{1}X_{2},gx_{1}x_{2}) \\
&=&-\beta _{1}B(x_{2}\otimes g;X_{1}X_{2},gx_{1}x_{2})-\gamma
_{1}B(x_{2}\otimes g;GX_{2},gx_{1}x_{2})+ \\
&&-B(gx_{1}x_{2}\otimes g;X_{2},gx_{1}x_{2})+B(gx_{2}\otimes
x_{1};X_{2},gx_{1}x_{2})
\end{eqnarray*}%
By applying $\left( \ref{eq.10}\right) $ this rewrites as%
\begin{equation*}
-\gamma _{1}B(gx_{2}\otimes 1_{H};GX_{2},gx_{1}x_{2})-B(x_{1}x_{2}\otimes
1_{H};X_{2},gx_{1}x_{2})+B(gx_{2}\otimes x_{1};X_{2},gx_{1}x_{2})=0.
\end{equation*}%
In view of the form of the elements we get%
\begin{equation*}
-B(x_{1}x_{2}\otimes 1_{H};X_{2},gx_{1}x_{2})+B(x_{1}x_{2}\otimes
1_{H};X_{2},gx_{1}x_{2})=0.
\end{equation*}%
which is trivial.

\paragraph{Equality $\left( \protect\ref{X1F51}\right) $}

rewrites as%
\begin{eqnarray*}
-B(gx_{2}\otimes 1_{H};X_{2},f) &=&-B(x_{2}\otimes g;X_{2},f)-\gamma
_{1}B(x_{2}\otimes g;GX_{1}X_{2},f) \\
&&-B(gx_{1}x_{2}\otimes g;X_{1}X_{2},f)+B(gx_{2}\otimes x_{1};X_{1}X_{2},f).
\end{eqnarray*}

\subparagraph{Case $f=1_{H}$}

\begin{eqnarray*}
-B(gx_{2}\otimes 1_{H};X_{2},1_{H}) &=&-B(x_{2}\otimes g;X_{2},1_{H})-\gamma
_{1}B(x_{2}\otimes g;GX_{1}X_{2},1_{H}) \\
&&-B(gx_{1}x_{2}\otimes g;X_{1}X_{2},1_{H})+B(gx_{2}\otimes
x_{1};X_{1}X_{2},1_{H}).
\end{eqnarray*}%
By applying $\left( \ref{eq.10}\right) $ this rewrites as%
\begin{equation*}
\gamma _{1}B(gx_{2}\otimes 1_{H};GX_{1}X_{2},1_{H})-B(x_{1}x_{2}\otimes
1_{H};X_{1}X_{2},1_{H})+B(gx_{2}\otimes x_{1};X_{1}X_{2},1_{H}).
\end{equation*}%
In view of the form of the elements we get%
\begin{equation*}
\gamma _{1}B(x_{1}x_{2}\otimes 1_{H};GX_{1},x_{1}x_{2})=0
\end{equation*}

which is $\left( \ref{G,x1x2, GF7,x1}\right) .$

\subparagraph{Case $f=gx_{1}$}

\begin{eqnarray*}
-B(gx_{2}\otimes 1_{H};X_{2},gx_{1}) &=&-B(x_{2}\otimes
g;X_{2},gx_{1})-\gamma _{1}B(x_{2}\otimes g;GX_{1}X_{2},gx_{1}) \\
&&-B(gx_{1}x_{2}\otimes g;X_{1}X_{2},gx_{1})+B(gx_{2}\otimes
x_{1};X_{1}X_{2},gx_{1}).
\end{eqnarray*}%
By applying $\left( \ref{eq.10}\right) $ this rewrites as%
\begin{eqnarray*}
&&2B(gx_{2}\otimes 1_{H};X_{2},gx_{1})+\gamma _{1}B(gx_{2}\otimes
1_{H};GX_{1}X_{2},gx_{1}) \\
&&+B(x_{1}x_{2}\otimes 1_{H};X_{1}X_{2},gx_{1})+B(gx_{2}\otimes
x_{1};X_{1}X_{2},gx_{1}).
\end{eqnarray*}%
In view of the form of the elements we get%
\begin{equation*}
B(x_{1}x_{2}\otimes 1_{H};X_{1},gx_{1}x_{2})=0
\end{equation*}%
which is $\left( \ref{G,x1x2, GF6,gx1x2}\right) .$

\subparagraph{Case $f=gx_{2}$}

\begin{eqnarray*}
-B(gx_{2}\otimes 1_{H};X_{2},gx_{2}) &=&-B(x_{2}\otimes
g;X_{2},gx_{2})-\gamma _{1}B(x_{2}\otimes g;GX_{1}X_{2},gx_{2}) \\
&&-B(gx_{1}x_{2}\otimes g;X_{1}X_{2},gx_{2})+B(gx_{2}\otimes
x_{1};X_{1}X_{2},gx_{2}).
\end{eqnarray*}

By applying $\left( \ref{eq.10}\right) $ this rewrites as%
\begin{eqnarray*}
&&2B(gx_{2}\otimes 1_{H};X_{2},gx_{2})+\gamma _{1}B(x_{2}\otimes
g;GX_{1}X_{2},gx_{2}) \\
&&+B(x_{1}x_{2}\otimes 1_{H};X_{1}X_{2},gx_{2})+B(gx_{2}\otimes
x_{1};X_{1}X_{2},gx_{2}).
\end{eqnarray*}%
In view of the form of the elements we get

\begin{equation*}
B(x_{1}x_{2}\otimes 1_{H};X_{2},gx_{1}x_{2})-B(x_{1}x_{2}\otimes
1_{H};X_{2},gx_{1}x_{2})=0
\end{equation*}%
which is trivial.

\paragraph{Equality $\left( \protect\ref{X1F61}\right) $}

rewrites as%
\begin{eqnarray*}
&&B(gx_{2}\otimes 1_{H};G,f)+\lambda B(gx_{2}\otimes 1_{H};GX_{1}X_{2},f) \\
&=&+B(x_{2}\otimes g;G,f)-B(gx_{1}x_{2}\otimes g;GX_{1},f)+B(gx_{2}\otimes
x_{1};GX_{1},f).
\end{eqnarray*}

\subparagraph{Case $f=1_{H}$}

\begin{eqnarray*}
&&B(gx_{2}\otimes 1_{H};G,1_{H})+\lambda B(gx_{2}\otimes
1_{H};GX_{1}X_{2},1_{H}) \\
&=&+B(x_{2}\otimes g;G,1_{H})-B(gx_{1}x_{2}\otimes
g;GX_{1},1_{H})+B(gx_{2}\otimes x_{1};GX_{1},1_{H}).
\end{eqnarray*}%
By applying $\left( \ref{eq.10}\right) $ this rewrites as%
\begin{equation*}
\lambda B(gx_{2}\otimes 1_{H};GX_{1}X_{2},1_{H})+B(x_{1}x_{2}\otimes
1_{H};GX_{1},1_{H})-B(gx_{2}\otimes x_{1};GX_{1},1_{H})
\end{equation*}

In view of the form of the elements we get%
\begin{equation*}
\lambda B(x_{1}x_{2}\otimes 1_{H};GX_{1},x_{1}x_{2})=0
\end{equation*}%
which is $\left( \ref{X1,x1x2,X1F61,x1}\right) .$

\subparagraph{Case $f=x_{1}x_{2}$}

\begin{eqnarray*}
&&B(gx_{2}\otimes 1_{H};G,x_{1}x_{2})+\lambda B(gx_{2}\otimes
1_{H};GX_{1}X_{2},x_{1}x_{2}) \\
&=&+B(x_{2}\otimes g;G,x_{1}x_{2})-B(gx_{1}x_{2}\otimes
g;GX_{1},x_{1}x_{2})+B(gx_{2}\otimes x_{1};GX_{1},x_{1}x_{2}).
\end{eqnarray*}%
By applying $\left( \ref{eq.10}\right) $ this rewrites as

\begin{eqnarray*}
&&+\lambda B(gx_{2}\otimes 1_{H};GX_{1}X_{2},x_{1}x_{2}) \\
&=&-B(x_{1}x_{2}\otimes 1_{H};GX_{1},x_{1}x_{2})+B(gx_{2}\otimes
x_{1};GX_{1},x_{1}x_{2}).
\end{eqnarray*}%
In view of the form of the elements we get%
\begin{equation*}
-B(x_{1}x_{2}\otimes 1_{H};GX_{1},x_{1}x_{2})+B(x_{1}x_{2}\otimes
1_{H};GX_{1},x_{1}x_{2})=0
\end{equation*}%
which is trivial

\subparagraph{Case $f=gx_{1}$}

\begin{eqnarray*}
&&B(gx_{2}\otimes 1_{H};G,gx_{1})+\lambda B(gx_{2}\otimes
1_{H};GX_{1}X_{2},gx_{1}) \\
&=&+B(x_{2}\otimes g;G,gx_{1})-B(gx_{1}x_{2}\otimes
g;GX_{1},gx_{1})+B(gx_{2}\otimes x_{1};GX_{1},gx_{1}).
\end{eqnarray*}%
By applying $\left( \ref{eq.10}\right) $ this rewrites as

\begin{eqnarray*}
&&+\lambda B(gx_{2}\otimes 1_{H};GX_{1}X_{2},gx_{1}) \\
&=&-2B(gx_{2}\otimes 1_{H};G,gx_{1})+B(x_{1}x_{2}\otimes
1_{H};GX_{1},gx_{1})+B(gx_{2}\otimes x_{1};GX_{1},gx_{1}).
\end{eqnarray*}%
In view of the form of the elements we get%
\begin{equation*}
-B(x_{1}x_{2}\otimes 1_{H};GX_{1},gx_{1})+B(x_{1}x_{2}\otimes
1_{H};GX_{1},gx_{1})=0
\end{equation*}%
which is trivial

\subparagraph{Case $f=gx_{2}$}

\begin{eqnarray*}
&&B(gx_{2}\otimes 1_{H};G,gx_{2})+\lambda B(gx_{2}\otimes
1_{H};GX_{1}X_{2},gx_{2}) \\
&=&+B(x_{2}\otimes g;G,gx_{2})-B(gx_{1}x_{2}\otimes
g;GX_{1},gx_{2})+B(gx_{2}\otimes x_{1};GX_{1},gx_{2}).
\end{eqnarray*}%
By applying $\left( \ref{eq.10}\right) $ this rewrites as%
\begin{equation*}
2B(gx_{2}\otimes 1_{H};G,gx_{2})-B(x_{1}x_{2}\otimes
1_{H};GX_{1},gx_{2})-B(gx_{2}\otimes x_{1};GX_{1},gx_{2})=0.
\end{equation*}%
In view of the form of the elements we get%
\begin{equation*}
B(x_{1}x_{2}\otimes 1_{H};G,gx_{1}x_{2})+B(x_{1}x_{2}\otimes
1_{H};GX_{1},gx_{2})=0
\end{equation*}%
which is $\left( \ref{X1,x1x2,X1F81,gx1}\right) .$

\paragraph{Equality $\left( \protect\ref{X1F71}\right) $}

rewrites as

\begin{gather*}
-\beta _{1}B(gx_{2}\otimes 1_{H};GX_{1}X_{2},f)=+\beta _{1}B(x_{2}\otimes
g;GX_{1}X_{2},f) \\
-B(gx_{1}x_{2}\otimes g;GX_{2},f)+B(gx_{2}\otimes x_{1};GX_{2},f)
\end{gather*}

\subparagraph{Case $f=1_{H}$}

\begin{gather*}
-\beta _{1}B(gx_{2}\otimes 1_{H};GX_{1}X_{2},1_{H})=+\beta
_{1}B(x_{2}\otimes g;GX_{1}X_{2},1_{H}) \\
-B(gx_{1}x_{2}\otimes g;GX_{2},1_{H})+B(gx_{2}\otimes x_{1};GX_{2},1_{H})
\end{gather*}%
By applying $\left( \ref{eq.10}\right) $ this rewrites as

\begin{gather*}
2\beta _{1}B(x_{1}x_{2}\otimes 1_{H};GX_{1},1_{H}) \\
-B(x_{1}x_{2}\otimes 1_{H};GX_{2},1_{H})+B(gx_{2}\otimes x_{1};GX_{2},1_{H}).
\end{gather*}%
In view of the form of the elements we get%
\begin{equation}
\beta _{1}B(gx_{2}\otimes 1_{H};GX_{1}X_{2},1_{H})=0.  \label{X1,gx2,X171,1H}
\end{equation}

\subparagraph{Case $f=x_{1}x_{2}$}

\begin{gather*}
-\beta _{1}B(gx_{2}\otimes 1_{H};GX_{1}X_{2},x_{1}x_{2})=+\beta
_{1}B(x_{2}\otimes g;GX_{1}X_{2},x_{1}x_{2}) \\
-B(gx_{1}x_{2}\otimes g;GX_{2},x_{1}x_{2})+B(gx_{2}\otimes
x_{1};GX_{2},x_{1}x_{2})
\end{gather*}%
By applying $\left( \ref{eq.10}\right) $ this rewrites as%
\begin{gather*}
2\beta _{1}B(x_{2}\otimes g;GX_{1}X_{2},x_{1}x_{2}) \\
-B(x_{1}x_{2}\otimes 1_{H};GX_{2},x_{1}x_{2})+B(gx_{2}\otimes
x_{1};GX_{2},x_{1}x_{2})=0.
\end{gather*}%
In view of the form of the elements we get%
\begin{equation*}
-B(x_{1}x_{2}\otimes 1_{H};GX_{2},x_{1}x_{2})+B(x_{1}x_{2}\otimes
1_{H};GX_{2},x_{1}x_{2})=0
\end{equation*}%
which is trivial.

\subparagraph{Case $f=gx_{1}$%
\protect\begin{gather*}
-\protect\beta _{1}B(gx_{2}\otimes 1_{H};GX_{1}X_{2},gx_{1})=+\protect\beta %
_{1}B(x_{2}\otimes g;GX_{1}X_{2},gx_{1}) \\
-B(gx_{1}x_{2}\otimes g;GX_{2},gx_{1})+B(gx_{2}\otimes x_{1};GX_{2},gx_{1})
\protect\end{gather*}%
}

By applying $\left( \ref{eq.10}\right) $ this rewrites as%
\begin{equation*}
+B(x_{1}x_{2}\otimes 1_{H};GX_{2},gx_{1})+B(gx_{2}\otimes
x_{1};GX_{2},gx_{1})
\end{equation*}%
In view of the form of the elements we get%
\begin{equation*}
+B(x_{1}x_{2}\otimes 1_{H};G,gx_{1}x_{2})+B(x_{1}x_{2}\otimes
1_{H};GX_{1},gx_{2})=0.
\end{equation*}%
which is $\left( \ref{X1,x1x2,X1F81,gx1}\right) .$

\subparagraph{Case $f=gx_{2}$}

\begin{gather*}
-\beta _{1}B(gx_{2}\otimes 1_{H};GX_{1}X_{2},gx_{2})=+\beta
_{1}B(x_{2}\otimes g;GX_{1}X_{2},gx_{2}) \\
-B(gx_{1}x_{2}\otimes g;GX_{2},gx_{2})+B(gx_{2}\otimes x_{1};GX_{2},gx_{2})
\end{gather*}%
By applying $\left( \ref{eq.10}\right) $ this rewrites and in view of the
form of the elements we get%
\begin{equation*}
B(x_{1}x_{2}\otimes 1_{H};GX_{2},gx_{2})-B(x_{1}x_{2}\otimes
1_{H};GX_{2},gx_{2})=0.
\end{equation*}%
which is trivial.

\paragraph{Equality $\left( \protect\ref{X1F81}\right) $}

rewrites as

\begin{gather*}
-B(gx_{2}\otimes 1_{H};GX_{2},f)=B(x_{2}\otimes g;GX_{2},f) \\
-B(gx_{1}x_{2}\otimes g;GX_{1}X_{2},f)+B(gx_{2}\otimes x_{1};GX_{1}X_{2},f).
\end{gather*}

\subparagraph{Case $f=g$}

\begin{gather*}
-B(gx_{2}\otimes 1_{H};GX_{2},g)=B(x_{2}\otimes g;GX_{2},g) \\
-B(gx_{1}x_{2}\otimes g;GX_{1}X_{2},g)+B(gx_{2}\otimes x_{1};GX_{1}X_{2},g).
\end{gather*}

By applying $\left( \ref{eq.10}\right) $ this rewrites as%
\begin{gather*}
2B(gx_{2}\otimes 1_{H};GX_{2},g) \\
-B(x_{1}x_{2}\otimes 1_{H};GX_{1}X_{2},g)+B(gx_{2}\otimes
x_{1};GX_{1}X_{2},g).
\end{gather*}%
In view of the form of the elements we get%
\begin{equation*}
B(x_{1}x_{2}\otimes 1_{H};G,gx_{1}x_{2})+B(x_{1}x_{2}\otimes
1_{H};GX_{1},gx_{2})=0
\end{equation*}%
which is $\left( \ref{X1,x1x2,X1F81,gx1}\right) .$

\subparagraph{Case $f=x_{1}$}

\begin{gather*}
-B(gx_{2}\otimes 1_{H};GX_{2},x_{1})=B(x_{2}\otimes g;GX_{2},x_{1}) \\
-B(gx_{1}x_{2}\otimes g;GX_{1}X_{2},x_{1})+B(gx_{2}\otimes
x_{1};GX_{1}X_{2},x_{1}).
\end{gather*}%
By applying $\left( \ref{eq.10}\right) $ this rewrites as%
\begin{equation*}
-B(x_{1}x_{2}\otimes 1_{H};GX_{1}X_{2},x_{1})+B(gx_{2}\otimes
x_{1};GX_{1}X_{2},x_{1})=0.
\end{equation*}

\begin{equation*}
-B(x_{1}x_{2}\otimes 1_{H};GX_{1},x_{1}x_{2})+B(x_{1}x_{2}\otimes
1_{H};GX_{1},x_{1}x_{2})=0
\end{equation*}%
which is trivial.

\subparagraph{Case $f=x_{2}$}

\begin{gather*}
-B(gx_{2}\otimes 1_{H};GX_{2},x_{2})=B(x_{2}\otimes g;GX_{2},x_{2}) \\
-B(gx_{1}x_{2}\otimes g;GX_{1}X_{2},x_{2})+B(gx_{2}\otimes
x_{1};GX_{1}X_{2},x_{2}).
\end{gather*}%
By applying $\left( \ref{eq.10}\right) $ this rewrites as%
\begin{equation*}
+B(x_{1}x_{2}\otimes 1_{H};GX_{1}X_{2},x_{2})+B(gx_{2}\otimes
x_{1};GX_{1}X_{2},x_{2})=0.
\end{equation*}%
In view of the form of the elements we get%
\begin{equation*}
-B(x_{1}x_{2}\otimes 1_{H};GX_{2},x_{1}x_{2})+B(x_{1}x_{2}\otimes
1_{H};GX_{2},x_{1}x_{2})=0
\end{equation*}%
which is trivial.

\subsubsection{Case $gx_{1}x_{2}\otimes 1_{H}$}

\paragraph{Equality $\left( \protect\ref{X1F11}\right) $}

rewrites as%
\begin{gather*}
\beta _{1}B(gx_{1}x_{2}\otimes 1_{H};X_{1},f)+\lambda B(gx_{1}x_{2}\otimes
1_{H};X_{2},f)=\gamma _{1}B(x_{1}x_{2}\otimes g;G,f) \\
+\beta _{1}B(x_{1}x_{2}\otimes g;X_{1},f)+B(gx_{1}x_{2}\otimes
x_{1};1_{A},f).
\end{gather*}%
Case $f=1_{H}$

\begin{gather*}
\beta _{1}B(gx_{1}x_{2}\otimes 1_{H};X_{1},1_{H})+\lambda
B(gx_{1}x_{2}\otimes 1_{H};X_{2},1_{H})=\gamma _{1}B(x_{1}x_{2}\otimes
g;G,1_{H}) \\
+\beta _{1}B(x_{1}x_{2}\otimes g;X_{1},1_{H})+B(gx_{1}x_{2}\otimes
x_{1};1_{A},1_{H}).
\end{gather*}

By applying $\left( \ref{eq.10}\right) $ this rewrites as%
\begin{gather*}
+\lambda B(gx_{1}x_{2}\otimes 1_{H};X_{2},1_{H})=\gamma
_{1}B(gx_{1}x_{2}\otimes 1_{H};G,1_{H}) \\
+B(gx_{1}x_{2}\otimes x_{1};1_{A},1_{H}).
\end{gather*}%
In view of the form of the elements we get%
\begin{equation}
+\lambda \left[ B(x_{1}\otimes 1_{H};1_{A},1_{H})+B(gx_{1}x_{2}\otimes
1_{H};1_{A},x_{2})\right] -\gamma _{1}B(gx_{1}x_{2}\otimes 1_{H};G,1_{H})=0
\label{X1,gx1x2,X1F11,1H}
\end{equation}

\subparagraph{Case $f=x_{1}x_{2}$}

\begin{gather*}
\beta _{1}B(gx_{1}x_{2}\otimes 1_{H};X_{1},x_{1}x_{2})+\lambda
B(gx_{1}x_{2}\otimes 1_{H};X_{2},x_{1}x_{2})=\gamma _{1}B(x_{1}x_{2}\otimes
g;G,x_{1}x_{2}) \\
+\beta _{1}B(x_{1}x_{2}\otimes g;X_{1},x_{1}x_{2})+B(gx_{1}x_{2}\otimes
x_{1};1_{A},x_{1}x_{2}).
\end{gather*}%
By applying $\left( \ref{eq.10}\right) $ this rewrites as%
\begin{gather*}
+\lambda B(gx_{1}x_{2}\otimes 1_{H};X_{2},x_{1}x_{2})=\gamma
_{1}B(gx_{1}x_{2}\otimes 1_{H};G,x_{1}x_{2}) \\
+B(gx_{1}x_{2}\otimes x_{1};1_{A},x_{1}x_{2}).
\end{gather*}%
In view of the form of the elements we get%
\begin{equation}
+\lambda B(x_{1}\otimes 1_{H};1_{A},x_{1}x_{2})-\gamma
_{1}B(gx_{1}x_{2}\otimes 1_{H};G,x_{1}x_{2})=0  \label{X1,gx1x2,X1F11,x1x2}
\end{equation}

\subparagraph{Case $f=gx_{1}$}

\begin{gather*}
\beta _{1}B(gx_{1}x_{2}\otimes 1_{H};X_{1},gx_{1})+\lambda
B(gx_{1}x_{2}\otimes 1_{H};X_{2},gx_{1})=\gamma _{1}B(x_{1}x_{2}\otimes
g;G,gx_{1}) \\
+\beta _{1}B(x_{1}x_{2}\otimes g;X_{1},gx_{1})+B(gx_{1}x_{2}\otimes
x_{1};1_{A},gx_{1}).
\end{gather*}%
By applying $\left( \ref{eq.10}\right) $ this rewrites as%
\begin{gather*}
2\beta _{1}B(gx_{1}x_{2}\otimes 1_{H};X_{1},gx_{1})+\lambda
B(gx_{1}x_{2}\otimes 1_{H};X_{2},gx_{1})+\gamma _{1}B(gx_{1}x_{2}\otimes
1_{H};G,gx_{1}) \\
=B(gx_{1}x_{2}\otimes x_{1};1_{A},gx_{1}).
\end{gather*}%
In view of the form of the elements we get%
\begin{gather}
-2\beta _{1}B(x_{2}\otimes 1_{H};1_{A},gx_{1})+\lambda \left[ B(x_{1}\otimes
1_{H};1_{A},gx_{1})+B(gx_{1}x_{2}\otimes 1_{H};1_{A},gx_{1}x_{2})\right]
\label{X1,gx1x2,X1F11,gx1} \\
+\gamma _{1}B(gx_{1}x_{2}\otimes 1_{H};G,gx_{1})+2B(gx_{1}x_{2}\otimes
1_{H};1_{A},g)=0.  \notag
\end{gather}

\subparagraph{Case $f=gx_{2}$}

\begin{gather*}
\beta _{1}B(gx_{1}x_{2}\otimes 1_{H};X_{1},gx_{2})+\lambda
B(gx_{1}x_{2}\otimes 1_{H};X_{2},gx_{2})=\gamma _{1}B(x_{1}x_{2}\otimes
g;G,gx_{2}) \\
+\beta _{1}B(x_{1}x_{2}\otimes g;X_{1},gx_{2})+B(gx_{1}x_{2}\otimes
x_{1};1_{A},gx_{2}).
\end{gather*}%
By applying $\left( \ref{eq.10}\right) $ this rewrites as%
\begin{gather*}
2\beta _{1}B(gx_{1}x_{2}\otimes 1_{H};X_{1},gx_{2})+\lambda
B(gx_{1}x_{2}\otimes 1_{H};X_{2},gx_{2})+\gamma _{1}B(gx_{1}x_{2}\otimes
1_{H};G,gx_{2})= \\
+B(gx_{1}x_{2}\otimes x_{1};1_{A},gx_{2}).
\end{gather*}%
In view of the form of the elements we get%
\begin{gather}
2\beta _{1}\left[ -B(x_{2}\otimes 1_{H};1_{A},gx_{2})-B(gx_{1}x_{2}\otimes
1_{H};1_{A},gx_{1}x_{2})\right]  \label{X1,gx1x2,X171,gx2} \\
+\lambda B(x_{1}\otimes 1_{H};1_{A},gx_{2})+\gamma _{1}B(gx_{1}x_{2}\otimes
1_{H};G,gx_{2})=0.  \notag
\end{gather}

\paragraph{Equality $\left( \protect\ref{X1F21}\right) $}

rewrites as%
\begin{eqnarray*}
&&\beta _{1}B(gx_{1}x_{2}\otimes 1_{H};GX_{1},f)+\lambda
B(gx_{1}x_{2}\otimes 1_{H};GX_{2},f) \\
&=&-\beta _{1}B(x_{1}x_{2}\otimes g;GX_{1},f)+B(gx_{1}x_{2}\otimes x_{1};G,f)
\end{eqnarray*}

\subparagraph{Case $f=g$}

\begin{eqnarray*}
&&\beta _{1}B(gx_{1}x_{2}\otimes 1_{H};GX_{1},g)+\lambda
B(gx_{1}x_{2}\otimes 1_{H};GX_{2},g) \\
&=&-\beta _{1}B(x_{1}x_{2}\otimes g;GX_{1},g)+B(gx_{1}x_{2}\otimes x_{1};G,g)
\end{eqnarray*}%
By applying $\left( \ref{eq.10}\right) $ this rewrites as%
\begin{eqnarray*}
&&2\beta _{1}B(gx_{1}x_{2}\otimes 1_{H};GX_{1},g)+\lambda
B(gx_{1}x_{2}\otimes 1_{H};GX_{2},g) \\
&=&+B(gx_{1}x_{2}\otimes x_{1};G,g)
\end{eqnarray*}%
In view of the form of the elements we get%
\begin{gather}
2\beta _{1}\left[ B(x_{2}\otimes 1_{H};G,g)-B(gx_{1}x_{2}\otimes
1_{H};G,gx_{1})\right]  \label{X1,gx1x2,XF21,g} \\
+\lambda \left[ -B(x_{1}\otimes 1_{H};G,g)-B(gx_{1}x_{2}\otimes
1_{H};G,gx_{2})\right] =0  \notag
\end{gather}

\subparagraph{Case $f=x_{1}$}

\begin{eqnarray*}
&&\beta _{1}B(gx_{1}x_{2}\otimes 1_{H};GX_{1},x_{1})+\lambda
B(gx_{1}x_{2}\otimes 1_{H};GX_{2},x_{1}) \\
&=&-\beta _{1}B(x_{1}x_{2}\otimes g;GX_{1},x_{1})+B(x_{1}x_{2}\otimes
g;G,x_{1})+B(gx_{1}x_{2}\otimes x_{1};G,x_{1})
\end{eqnarray*}%
By applying $\left( \ref{eq.10}\right) $ this rewrites as%
\begin{equation*}
+\lambda B(gx_{1}x_{2}\otimes 1_{H};GX_{2},x_{1})+B(gx_{1}x_{2}\otimes
1_{H};G,x_{1})-B(gx_{1}x_{2}\otimes x_{1};G,x_{1})=0
\end{equation*}%
In view of the form of the elements we get%
\begin{equation*}
+\lambda \left[ B(x_{1}\otimes 1_{H};G,x_{1})+B(gx_{1}x_{2}\otimes
1_{H};G,x_{1}x_{2})\right] =0
\end{equation*}%
which follows from $\left( \ref{X1,x2,X1F71,x1}\right) $

\subparagraph{Case $f=x_{2}$}

\begin{eqnarray*}
&&\beta _{1}B(gx_{1}x_{2}\otimes 1_{H};GX_{1},x_{2})+\lambda
B(gx_{1}x_{2}\otimes 1_{H};GX_{2},x_{2}) \\
&=&-\beta _{1}B(x_{1}x_{2}\otimes g;GX_{1},x_{2})+B(x_{1}x_{2}\otimes
g;G,x_{2})+B(gx_{1}x_{2}\otimes x_{1};G,x_{2})
\end{eqnarray*}%
By applying $\left( \ref{eq.10}\right) $ this rewrites as%
\begin{equation*}
+\lambda B(gx_{1}x_{2}\otimes 1_{H};GX_{2},x_{2})+B(gx_{1}x_{2}\otimes
1_{H};G,x_{2})-B(gx_{1}x_{2}\otimes x_{1};G,x_{2})=0
\end{equation*}%
In view of the form of the elements we get%
\begin{equation*}
\lambda B(x_{1}\otimes 1_{H};G,x_{2})=0
\end{equation*}%
which follows from $\left( \ref{X1,g,X1F21,x2}\right) .$

\subparagraph{Case $f=gx_{1}x_{2}$}

\begin{eqnarray*}
&&\beta _{1}B(gx_{1}x_{2}\otimes 1_{H};GX_{1},gx_{1}x_{2})+\lambda
B(gx_{1}x_{2}\otimes 1_{H};GX_{2},gx_{1}x_{2}) \\
&=&-\beta _{1}B(x_{1}x_{2}\otimes g;GX_{1},gx_{1}x_{2})+B(gx_{1}x_{2}\otimes
x_{1};G,gx_{1}x_{2})
\end{eqnarray*}%
By applying $\left( \ref{eq.10}\right) $ this rewrites as%
\begin{eqnarray*}
&&2\beta _{1}B(gx_{1}x_{2}\otimes 1_{H};GX_{1},gx_{1}x_{2})+\lambda
B(gx_{1}x_{2}\otimes 1_{H};GX_{2},gx_{1}x_{2}) \\
&=&+B(gx_{1}x_{2}\otimes x_{1};G,gx_{1}x_{2})
\end{eqnarray*}%
In view of the form of the elements we get

\begin{equation}
2\beta _{1}B(x_{2}\otimes 1_{H};G,gx_{1}x_{2})-\lambda B(x_{1}\otimes
1_{H};G,gx_{1}x_{2})=0  \label{X1,gx1x2,X1F21,gx1x2}
\end{equation}

\paragraph{Equality $\left( \protect\ref{X1F31}\right) $}

rewrites as%
\begin{gather*}
B(gx_{1}x_{2}\otimes 1_{H};1_{A},f)+\lambda B(gx_{1}x_{2}\otimes
1_{H};X_{1}X_{2},f)= \\
B(x_{1}x_{2}\otimes g;1_{A},f)+\gamma _{1}B(x_{1}x_{2}\otimes
g;GX_{1},f)+B(gx_{1}x_{2}\otimes x_{1};X_{1},f)
\end{gather*}

\subparagraph{Case $f=g$}

\begin{gather*}
B(gx_{1}x_{2}\otimes 1_{H};1_{A},g)+\lambda B(gx_{1}x_{2}\otimes
1_{H};X_{1}X_{2},g)= \\
B(x_{1}x_{2}\otimes g;1_{A},g)+\gamma _{1}B(x_{1}x_{2}\otimes
g;GX_{1},g)+B(gx_{1}x_{2}\otimes x_{1};X_{1},g)
\end{gather*}%
By applying $\left( \ref{eq.10}\right) $ this rewrites as%
\begin{gather*}
+\lambda B(gx_{1}x_{2}\otimes 1_{H};X_{1}X_{2},g)= \\
+\gamma _{1}B(gx_{1}x_{2}\otimes 1_{H};GX_{1},g)+B(gx_{1}x_{2}\otimes
x_{1};X_{1},g)
\end{gather*}%
In view of the form of the elements we get%
\begin{gather}
+\lambda \left[
\begin{array}{c}
B(g\otimes 1_{H};1_{A},g)+B(x_{2}\otimes \ 1_{H};1_{A},gx_{2})+ \\
+B(x_{1}\otimes 1_{H};1_{A},gx_{1})+B(gx_{1}x_{2}\otimes
1_{H};1_{A},gx_{1}x_{2})%
\end{array}%
\right]  \label{X1,gx1x2,X1F31,g} \\
-\gamma _{1}\left[ B(x_{2}\otimes 1_{H};G,g)-B(gx_{1}x_{2}\otimes
1_{H};G,gx_{1})\right] =0  \notag
\end{gather}

\subparagraph{Case $f=x_{1}$}

\begin{gather*}
B(gx_{1}x_{2}\otimes 1_{H};1_{A},x_{1})+\lambda B(gx_{1}x_{2}\otimes
1_{H};X_{1}X_{2},x_{1})= \\
B(x_{1}x_{2}\otimes g;1_{A},x_{1})+\gamma _{1}B(x_{1}x_{2}\otimes
g;GX_{1},x_{1})+B(gx_{1}x_{2}\otimes x_{1};X_{1},x_{1})
\end{gather*}%
By applying $\left( \ref{eq.10}\right) $ this rewrites as%
\begin{gather*}
2B(gx_{1}x_{2}\otimes 1_{H};1_{A},x_{1})+\lambda B(gx_{1}x_{2}\otimes
1_{H};X_{1}X_{2},x_{1})= \\
-\gamma _{1}B(gx_{1}x_{2}\otimes 1_{H};GX_{1},x_{1})+B(gx_{1}x_{2}\otimes
x_{1};X_{1},x_{1})
\end{gather*}%
In view of the form of the elements we get%
\begin{gather}
2B(gx_{1}x_{2}\otimes 1_{H};1_{A},x_{1})+\lambda \left[ B(g\otimes
1_{H};1_{A},x_{1})+B(x_{2}\otimes \ 1_{H};1_{A},x_{1}x_{2})\right]
\label{X1,gx1x2,X1F31,x1} \\
+\gamma _{1}B(x_{2}\otimes 1_{H};G,x_{1})=0  \notag
\end{gather}

\subparagraph{Case $f=x_{2}$}

\begin{gather*}
B(gx_{1}x_{2}\otimes 1_{H};1_{A},x_{2})+\lambda B(gx_{1}x_{2}\otimes
1_{H};X_{1}X_{2},x_{2})= \\
B(x_{1}x_{2}\otimes g;1_{A},x_{2})+\gamma _{1}B(x_{1}x_{2}\otimes
g;GX_{1},x_{2})+B(gx_{1}x_{2}\otimes x_{1};X_{1},x_{2})
\end{gather*}%
By applying $\left( \ref{eq.10}\right) $ this rewrites as%
\begin{gather*}
2B(gx_{1}x_{2}\otimes 1_{H};1_{A},x_{2})+\lambda B(gx_{1}x_{2}\otimes
1_{H};X_{1}X_{2},x_{2})+ \\
+\gamma _{1}B(gx_{1}x_{2}\otimes 1_{H};GX_{1},x_{2})-B(gx_{1}x_{2}\otimes
x_{1};X_{1},x_{2})=0
\end{gather*}%
In view of the form of the elements we get%
\begin{gather}
2B(gx_{1}x_{2}\otimes 1_{H};1_{A},x_{2})+\lambda \left[ B(g\otimes
1_{H};X_{2},1_{H})-B(x_{1}\otimes 1_{H};1_{A},x_{1}x_{2})\right] +
\label{X1,gx1x2,X1F31,x2} \\
+\gamma _{1}\left[ B(x_{2}\otimes 1_{H};G,x_{2})+B(gx_{1}x_{2}\otimes
1_{H};G,x_{1}x_{2})\right] =0  \notag
\end{gather}

\subparagraph{Case $f=gx_{1}x_{2}$}

\begin{gather*}
B(gx_{1}x_{2}\otimes 1_{H};1_{A},gx_{1}x_{2})+\lambda B(gx_{1}x_{2}\otimes
1_{H};X_{1}X_{2},gx_{1}x_{2})= \\
B(x_{1}x_{2}\otimes g;1_{A},gx_{1}x_{2})+\gamma _{1}B(x_{1}x_{2}\otimes
g;GX_{1},gx_{1}x_{2})+B(gx_{1}x_{2}\otimes x_{1};X_{1},gx_{1}x_{2})
\end{gather*}%
By applying $\left( \ref{eq.10}\right) $ this rewrites as%
\begin{gather*}
+\lambda B(gx_{1}x_{2}\otimes 1_{H};X_{1}X_{2},gx_{1}x_{2})+ \\
-\gamma _{1}B(gx_{1}x_{2}\otimes
1_{H};GX_{1},gx_{1}x_{2})-B(gx_{1}x_{2}\otimes x_{1};X_{1},gx_{1}x_{2})=0
\end{gather*}%
In view of the form of the elements we get%
\begin{gather}
+\lambda B\left( g\otimes 1_{H};1_{A},gx_{1}x_{2}\right) +
\label{X1,gx1x2,X1F31,gx1x2} \\
-\gamma _{1}B(x_{2}\otimes 1_{H};G,gx_{1}x_{2})-\left[ +2B(x_{2}\otimes
1_{H};1_{A},gx_{2})+2B(gx_{1}x_{2}\otimes 1_{H};1_{A},gx_{1}x_{2})\right] =0
\notag
\end{gather}

\paragraph{Equality $\left( \protect\ref{X1F41}\right) $}

rewrites as%
\begin{eqnarray*}
-\beta _{1}B(gx_{1}x_{2}\otimes 1_{H};X_{1}X_{2},f) &=&\beta
_{1}B(x_{1}x_{2}\otimes g;X_{1}X_{2},f) \\
&&+\gamma _{1}B(x_{1}x_{2}\otimes g;GX_{2},f)+B(gx_{1}x_{2}\otimes
x_{1};X_{2},f)
\end{eqnarray*}

\subparagraph{Case $f=g$}

\begin{gather*}
-\beta _{1}B(gx_{1}x_{2}\otimes 1_{H};X_{1}X_{2},g)=\beta
_{1}B(x_{1}x_{2}\otimes g;X_{1}X_{2},g) \\
+\gamma _{1}B(x_{1}x_{2}\otimes g;GX_{2},g)+B(gx_{1}x_{2}\otimes
x_{1};X_{2},g)
\end{gather*}%
By applying $\left( \ref{eq.10}\right) $ this rewrites as%
\begin{eqnarray*}
&&2\beta _{1}B(gx_{1}x_{2}\otimes 1_{H};X_{1}X_{2},g) \\
&&+\gamma _{1}B(gx_{1}x_{2}\otimes 1_{H};GX_{2},g)+B(gx_{1}x_{2}\otimes
x_{1};X_{2},g)=0
\end{eqnarray*}%
In view of the form of the elements we get%
\begin{eqnarray}
&&2\beta _{1}\left[ B(g\otimes 1_{H};1_{A},g)+B(x_{2}\otimes \
1_{H};1_{A},gx_{2})+B(x_{1}\otimes 1_{H};1_{A},gx_{1})+B(gx_{1}x_{2}\otimes
1_{H};1_{A},gx_{1}x_{2})\right]  \label{X1,gx1x2,X1F41,g} \\
&&-\gamma _{1}\left[ B(x_{1}\otimes 1_{H};G,g)+B(gx_{1}x_{2}\otimes
1_{H};G,gx_{2})\right] =0  \notag
\end{eqnarray}

\subparagraph{Case $f=x_{1}$}

\begin{gather*}
-\beta _{1}B(gx_{1}x_{2}\otimes 1_{H};X_{1}X_{2},x_{1})=\beta
_{1}B(x_{1}x_{2}\otimes g;X_{1}X_{2},x_{1}) \\
+\gamma _{1}B(x_{1}x_{2}\otimes g;GX_{2},x_{1})+B(gx_{1}x_{2}\otimes
x_{1};X_{2},x_{1})
\end{gather*}%
By applying $\left( \ref{eq.10}\right) $ this rewrites as%
\begin{equation*}
-\gamma _{1}B(gx_{1}x_{2}\otimes 1_{H};GX_{2},x_{1})+B(gx_{1}x_{2}\otimes
x_{1};X_{2},x_{1})=0
\end{equation*}%
In view of the form of the elements we get%
\begin{equation*}
+\gamma _{1}\left[ B(x_{1}\otimes 1_{H};G,x_{1})+B(gx_{1}x_{2}\otimes
1_{H};G,x_{1}x_{2})\right] =0
\end{equation*}%
which follows from $\left( \ref{X1,x2,X1F71,x1}\right) $

\subparagraph{Case $f=x_{2}$}

\begin{gather*}
-\beta _{1}B(gx_{1}x_{2}\otimes 1_{H};X_{1}X_{2},x_{2})=\beta
_{1}B(x_{1}x_{2}\otimes g;X_{1}X_{2},x_{2}) \\
+\gamma _{1}B(x_{1}x_{2}\otimes g;GX_{2},x_{2})+B(x_{1}x_{2}\otimes
x_{1};X_{2},x_{2})
\end{gather*}%
By applying $\left( \ref{eq.10}\right) $ this rewrites as%
\begin{equation*}
-\gamma _{1}B(gx_{1}x_{2}\otimes 1_{H};GX_{2},x_{2})+B(x_{1}x_{2}\otimes
x_{1};X_{2},x_{2})=0
\end{equation*}%
In view of the form of the elements we get%
\begin{equation*}
\gamma _{1}B(x_{1}\otimes 1_{H};G,x_{2})=0
\end{equation*}%
which follows from $\left( \ref{X1,g,X1F21,x2}\right) $

\subparagraph{Case $f=gx_{1}x_{2}$}

\begin{gather*}
-\beta _{1}B(gx_{1}x_{2}\otimes 1_{H};X_{1}X_{2},gx_{1}x_{2})=\beta
_{1}B(x_{1}x_{2}\otimes g;X_{1}X_{2},gx_{1}x_{2}) \\
+\gamma _{1}B(x_{1}x_{2}\otimes g;GX_{2},gx_{1}x_{2})+B(x_{1}x_{2}\otimes
x_{1};X_{2},gx_{1}x_{2})
\end{gather*}%
By applying $\left( \ref{eq.10}\right) $ this rewrites as%
\begin{gather*}
2\beta _{1}B(gx_{1}x_{2}\otimes 1_{H};X_{1}X_{2},gx_{1}x_{2})+ \\
+\gamma _{1}B(gx_{1}x_{2}\otimes
1_{H};GX_{2},gx_{1}x_{2})+B(x_{1}x_{2}\otimes x_{1};X_{2},gx_{1}x_{2})=0
\end{gather*}%
In view of the form of the elements we get%
\begin{eqnarray}
&&2\beta _{1}B\left( g\otimes 1_{H};1_{A},gx_{1}x_{2}\right)
\label{X1,gx1x2,X1F41,gx1x2} \\
&&-\gamma _{1}B(x_{1}\otimes 1_{H};G,gx_{1}x_{2})-2B(x_{1}\otimes
1_{H};1_{A},gx_{2})=0  \notag
\end{eqnarray}

\paragraph{Equality $\left( \protect\ref{X1F51}\right) $}

rewrites as

\begin{eqnarray*}
-B(gx_{1}x_{2}\otimes 1_{H};X_{2},f) &=&B(x_{1}x_{2}\otimes
g;X_{2},f)+\gamma _{1}B(x_{1}x_{2}\otimes g;GX_{1}X_{2},f) \\
&&+B(gx_{1}x_{2}\otimes x_{1};X_{1}X_{2},f)
\end{eqnarray*}

\subparagraph{Case $f=1_{H}$}

\begin{eqnarray*}
-B(gx_{1}x_{2}\otimes 1_{H};X_{2},1_{H}) &=&B(x_{1}x_{2}\otimes
g;X_{2},1_{H})+\gamma _{1}B(x_{1}x_{2}\otimes g;GX_{1}X_{2},1_{H}) \\
&&+B(gx_{1}x_{2}\otimes x_{1};X_{1}X_{2},1_{H})
\end{eqnarray*}%
By applying $\left( \ref{eq.10}\right) $ this rewrites as%
\begin{equation*}
+2B(gx_{1}x_{2}\otimes 1_{H};X_{2},1_{H})+\gamma _{1}B(x_{1}x_{2}\otimes
g;GX_{1}X_{2},1_{H})+B(gx_{1}x_{2}\otimes x_{1};X_{1}X_{2},1_{H})=0
\end{equation*}%
In view of the form of the elements we get%
\begin{gather}
2\left[ B(x_{1}\otimes 1_{H};1_{A},1_{H})+B(gx_{1}x_{2}\otimes
1_{H};1_{A},x_{2})\right]  \label{X1,gx1x2,X1F51,1H} \\
+\gamma _{1}\left[
\begin{array}{c}
B(g\otimes 1_{H};G,1_{H})+B(x_{2}\otimes \ 1_{H};G,x_{2})+ \\
+B(x_{1}\otimes 1_{H};G,x_{1})+B(gx_{1}x_{2}\otimes 1_{H};G,x_{1}x_{2})%
\end{array}%
\right] =0  \notag
\end{gather}

\subparagraph{Case $f=x_{1}x_{2}$}

\begin{eqnarray*}
&&-B(gx_{1}x_{2}\otimes 1_{H};X_{2},x_{1}x_{2}) \\
&=&B(x_{1}x_{2}\otimes g;X_{2},x_{1}x_{2})+\gamma _{1}B(x_{1}x_{2}\otimes
g;GX_{1}X_{2},x_{1}x_{2}) \\
&&+B(gx_{1}x_{2}\otimes x_{1};X_{1}X_{2},x_{1}x_{2})
\end{eqnarray*}%
By applying $\left( \ref{eq.10}\right) $ this rewrites as%
\begin{eqnarray*}
&&2B(gx_{1}x_{2}\otimes 1_{H};X_{2},x_{1}x_{2})+\gamma
_{1}B(gx_{1}x_{2}\otimes 1_{H};GX_{1}X_{2},x_{1}x_{2}) \\
&&+B(gx_{1}x_{2}\otimes x_{1};X_{1}X_{2},x_{1}x_{2})=0
\end{eqnarray*}%
In view of the form of the elements we get%
\begin{equation*}
2B(x_{1}\otimes 1_{H};1_{A},x_{1}x_{2})+\gamma _{1}B(g\otimes
1_{H};G,x_{1}x_{2})=0
\end{equation*}%
which is $\left( \ref{G,x1, GF2,x1x2}\right) .$

\subparagraph{Case $f=gx_{1}$}

\begin{eqnarray*}
-B(gx_{1}x_{2}\otimes 1_{H};X_{2},gx_{1}) &=&B(x_{1}x_{2}\otimes
g;X_{2},gx_{1})+\gamma _{1}B(x_{1}x_{2}\otimes g;GX_{1}X_{2},gx_{1}) \\
&&+B(gx_{1}x_{2}\otimes x_{1};X_{1}X_{2},gx_{1})
\end{eqnarray*}%
By applying $\left( \ref{eq.10}\right) $ this rewrites as%
\begin{equation*}
\gamma _{1}B(gx_{1}x_{2}\otimes
1_{H};GX_{1}X_{2},gx_{1})+B(gx_{1}x_{2}\otimes x_{1};X_{1}X_{2},gx_{1})=0
\end{equation*}%
In view of the form of the elements we get%
\begin{gather}
+\gamma _{1}\left[ B(g\otimes 1_{H};G,gx_{1})+B(x_{2}\otimes \
1_{H};G,gx_{1}x_{2})\right]  \label{X1,gx1x2,X1F51,gx1} \\
-\left[
\begin{array}{c}
2B(g\otimes 1_{H};1_{A},g)+2B(x_{2}\otimes \ 1_{H};1_{A},gx_{2})+ \\
+2B(x_{1}\otimes 1_{H};1_{A},gx_{1})+2B(gx_{1}x_{2}\otimes
1_{H};1_{A},gx_{1}x_{2})%
\end{array}%
\right] =0  \notag
\end{gather}

\subparagraph{Case $f=gx_{2}$}

\begin{eqnarray*}
-B(gx_{1}x_{2}\otimes 1_{H};X_{2},gx_{2}) &=&B(x_{1}x_{2}\otimes
g;X_{2},gx_{2})+\gamma _{1}B(x_{1}x_{2}\otimes g;GX_{1}X_{2},gx_{2}) \\
&&+B(gx_{1}x_{2}\otimes x_{1};X_{1}X_{2},gx_{2})
\end{eqnarray*}%
By applying $\left( \ref{eq.10}\right) $ this rewrites as%
\begin{equation*}
-\gamma _{1}B(gx_{1}x_{2}\otimes
1_{H};GX_{1}X_{2},gx_{2})+B(gx_{1}x_{2}\otimes x_{1};X_{1}X_{2},gx_{2})=0
\end{equation*}%
In view of the form of the elements we get

\begin{equation*}
\gamma _{1}\left[ B(g\otimes 1_{H};G,gx_{2})-B(x_{1}\otimes
1_{H};G,gx_{1}x_{2})\right] =0
\end{equation*}%
which is $\left( \ref{G,gx1x2, GF7,gx2}\right) .$

\paragraph{Equality $\left( \protect\ref{X1F61}\right) $}

rewrites as%
\begin{eqnarray*}
&&B(gx_{1}x_{2}\otimes 1_{H};G,f)+\lambda B(gx_{1}x_{2}\otimes
1_{H};GX_{1}X_{2},f) \\
&=&-B(x_{1}x_{2}\otimes g;G,f)+B(gx_{1}x_{2}\otimes x_{1};GX_{1},f)
\end{eqnarray*}

\subparagraph{Case $f=1_{H}$}

\begin{eqnarray*}
&&B(gx_{1}x_{2}\otimes 1_{H};G,1_{H})+\lambda B(gx_{1}x_{2}\otimes
1_{H};GX_{1}X_{2},1_{H}) \\
&=&-B(x_{1}x_{2}\otimes g;G,1_{H})+B(gx_{1}x_{2}\otimes x_{1};GX_{1},1_{H})
\end{eqnarray*}%
By applying $\left( \ref{eq.10}\right) $ this rewrites as%
\begin{equation*}
2B(gx_{1}x_{2}\otimes 1_{H};G,1_{H})+\lambda B(gx_{1}x_{2}\otimes
1_{H};GX_{1}X_{2},1_{H})-B(gx_{1}x_{2}\otimes x_{1};GX_{1},1_{H})=0
\end{equation*}%
In view of the form of the elements we get%
\begin{gather}
2B(gx_{1}x_{2}\otimes 1_{H};G,1_{H})+  \label{X1,gx1x2,X1F61,1H} \\
+\lambda \left[ B(g\otimes 1_{H};G,1_{H})+B(x_{2}\otimes \
1_{H};G,x_{2})+B(x_{1}\otimes 1_{H};G,x_{1})+B(gx_{1}x_{2}\otimes
1_{H};G,x_{1}x_{2})\right] =0  \notag
\end{gather}

\subparagraph{Case $f=x_{1}x_{2}$}

\begin{eqnarray*}
&&B(gx_{1}x_{2}\otimes 1_{H};G,x_{1}x_{2})+\lambda B(gx_{1}x_{2}\otimes
1_{H};GX_{1}X_{2},x_{1}x_{2}) \\
&=&-B(x_{1}x_{2}\otimes g;G,x_{1}x_{2})+B(gx_{1}x_{2}\otimes
x_{1};GX_{1},x_{1}x_{2})
\end{eqnarray*}%
By applying $\left( \ref{eq.10}\right) $ this rewrites as%
\begin{gather*}
2B(gx_{1}x_{2}\otimes 1_{H};G,x_{1}x_{2})+\lambda B(gx_{1}x_{2}\otimes
1_{H};GX_{1}X_{2},x_{1}x_{2}) \\
-B(gx_{1}x_{2}\otimes x_{1};GX_{1},x_{1}x_{2})=0
\end{gather*}%
In view of the form of the elements we get%
\begin{equation}
2B(gx_{1}x_{2}\otimes 1_{H};G,x_{1}x_{2})+\lambda B(g\otimes
1_{H};G,x_{1}x_{2})=0  \label{X1,gx1x2,X1F61,x1x2}
\end{equation}

\subparagraph{Case $f=gx_{1}$}

\begin{eqnarray*}
&&B(gx_{1}x_{2}\otimes 1_{H};G,gx_{1})+\lambda B(gx_{1}x_{2}\otimes
1_{H};GX_{1}X_{2},gx_{1}) \\
&=&-B(x_{1}x_{2}\otimes g;G,gx_{1})+B(gx_{1}x_{2}\otimes x_{1};GX_{1},gx_{1})
\end{eqnarray*}%
By applying $\left( \ref{eq.10}\right) $ this rewrites as%
\begin{equation*}
\lambda B(gx_{1}x_{2}\otimes 1_{H};GX_{1}X_{2},gx_{1})-B(gx_{1}x_{2}\otimes
x_{1};GX_{1},gx_{1})=0
\end{equation*}%
In view of the form of the elements we get%
\begin{gather}
\lambda \left[ B(g\otimes 1_{H};G,gx_{1})+B(x_{2}\otimes \
1_{H};G,gx_{1}x_{2})\right] +  \label{X1,gx1x2,X1F61,gx1} \\
+\left[ 2B(x_{2}\otimes 1_{H};G,g)-2B(gx_{1}x_{2}\otimes 1_{H};G,gx_{1})%
\right] =0  \notag
\end{gather}

\subparagraph{Case $f=gx_{2}$}

\begin{eqnarray*}
&&B(gx_{1}x_{2}\otimes 1_{H};G,gx_{2})+\lambda B(gx_{1}x_{2}\otimes
1_{H};GX_{1}X_{2},gx_{2}) \\
&=&-B(x_{1}x_{2}\otimes g;G,gx_{2})+B(gx_{1}x_{2}\otimes x_{1};GX_{1},gx_{2})
\end{eqnarray*}%
By applying $\left( \ref{eq.10}\right) $ this rewrites as

\begin{equation*}
+\lambda B(gx_{1}x_{2}\otimes 1_{H};GX_{1}X_{2},gx_{2})-B(gx_{1}x_{2}\otimes
x_{1};GX_{1},gx_{2})=0
\end{equation*}%
In view of the form of the elements we get%
\begin{equation*}
+\lambda \left[ B(g\otimes 1_{H};G,gx_{2})-B(x_{1}\otimes
1_{H};G,gx_{1}x_{2})\right] =0
\end{equation*}%
which follows from $\left( \ref{X1,g,X1F21,gx1x2}\right) .$

\paragraph{Equality $\left( \protect\ref{X1F71}\right) $}

rewrites as

\begin{eqnarray*}
&&-\beta _{1}B(gx_{1}x_{2}\otimes 1_{H};GX_{1}X_{2},f) \\
&=&-\beta _{1}B(x_{1}x_{2}\otimes g;GX_{1}X_{2},f)+B(gx_{1}x_{2}\otimes
x_{1};GX_{2},f).
\end{eqnarray*}

\subparagraph{Case $f=1_{H}$}

\begin{eqnarray*}
&&-\beta _{1}B(gx_{1}x_{2}\otimes 1_{H};GX_{1}X_{2},1_{H}) \\
&=&-\beta _{1}B(x_{1}x_{2}\otimes g;GX_{1}X_{2},1_{H})+B(gx_{1}x_{2}\otimes
x_{1};GX_{2},1_{H}).
\end{eqnarray*}%
By applying $\left( \ref{eq.10}\right) $ this rewrites as%
\begin{equation*}
B(gx_{1}x_{2}\otimes x_{1};GX_{2},1_{H})=0
\end{equation*}%
which is already known.

\subparagraph{Case $f=x_{1}x_{2}$}

\begin{eqnarray*}
&&-\beta _{1}B(gx_{1}x_{2}\otimes 1_{H};GX_{1}X_{2},x_{1}x_{2}) \\
&=&-\beta _{1}B(x_{1}x_{2}\otimes
g;GX_{1}X_{2},x_{1}x_{2})+B(gx_{1}x_{2}\otimes x_{1};GX_{2},x_{1}x_{2}).
\end{eqnarray*}%
By applying $\left( \ref{eq.10}\right) $ this rewrites as%
\begin{equation*}
B(gx_{1}x_{2}\otimes x_{1};GX_{2},x_{1}x_{2})=0
\end{equation*}%
which is already known.

\subparagraph{Case $f=gx_{1}$}

\begin{eqnarray*}
&&-\beta _{1}B(gx_{1}x_{2}\otimes 1_{H};GX_{1}X_{2},gx_{1}) \\
&=&-\beta _{1}B(x_{1}x_{2}\otimes g;GX_{1}X_{2},gx_{2})+B(gx_{1}x_{2}\otimes
x_{1};GX_{2},gx_{1}).
\end{eqnarray*}

By applying $\left( \ref{eq.10}\right) $ this rewrites as%
\begin{equation*}
2\beta _{1}B(gx_{1}x_{2}\otimes
1_{H};GX_{1}X_{2},gx_{1})+B(gx_{1}x_{2}\otimes x_{1};GX_{2},gx_{1})=0.
\end{equation*}%
In view of the form of the elements we get%
\begin{gather}
\beta _{1}\left[ B(g\otimes 1_{H};G,gx_{1})+B(x_{2}\otimes \
1_{H};G,gx_{1}x_{2})\right]  \label{X1,gx1x2,X1F71,gx1} \\
+B(x_{1}\otimes 1_{H};G,g)+B(gx_{1}x_{2}\otimes 1_{H};G,gx_{2})=0.  \notag
\end{gather}

\subparagraph{Case $f=gx_{2}$}

\begin{eqnarray*}
&&-\beta _{1}B(gx_{1}x_{2}\otimes 1_{H};GX_{1}X_{2},gx_{2}) \\
&=&-\beta _{1}B(x_{1}x_{2}\otimes g;GX_{1}X_{2},gx_{2})+B(gx_{1}x_{2}\otimes
x_{1};GX_{2},gx_{2}).
\end{eqnarray*}%
By applying $\left( \ref{eq.10}\right) $ this rewrites as%
\begin{equation*}
2\beta _{1}B(gx_{1}x_{2}\otimes
1_{H};GX_{1}X_{2},gx_{2})+B(gx_{1}x_{2}\otimes x_{1};GX_{2},gx_{2})=0.
\end{equation*}%
In view of the form of the elements we get

\begin{equation*}
2\beta _{1}\left[ B(g\otimes 1_{H};G,gx_{2})-B(x_{1}\otimes
1_{H};G,gx_{1}x_{2})\right] =0.
\end{equation*}%
which follows from $\left( \ref{X1,g,X1F21,gx1x2}\right) .$

\paragraph{Equality $\left( \protect\ref{X1F81}\right) $}

rewrites as

\begin{eqnarray*}
&&-B(gx_{1}x_{2}\otimes 1_{H};GX_{2},f) \\
&=&-B(x_{1}x_{2}\otimes g;GX_{2},f)+B(gx_{1}x_{2}\otimes
x_{1};GX_{1}X_{2},f).
\end{eqnarray*}

\subparagraph{Case $f=g$}

\begin{equation*}
-B(gx_{1}x_{2}\otimes 1_{H};GX_{2},g)=-B(x_{1}x_{2}\otimes
g;GX_{2},g)+B(gx_{1}x_{2}\otimes x_{1};GX_{1}X_{2},g).
\end{equation*}%
By applying $\left( \ref{eq.10}\right) $ this rewrites as%
\begin{equation*}
B(gx_{1}x_{2}\otimes x_{1};GX_{1}X_{2},g)=0
\end{equation*}%
which is already known.

\subparagraph{Case $f=x_{1}$}

\begin{equation*}
-B(gx_{1}x_{2}\otimes 1_{H};GX_{2},x_{1})=-B(x_{1}x_{2}\otimes
g;GX_{2},x_{1})+B(gx_{1}x_{2}\otimes x_{1};GX_{1}X_{2},x_{1}).
\end{equation*}%
By applying $\left( \ref{eq.10}\right) $ this rewrites as%
\begin{equation*}
2B(gx_{1}x_{2}\otimes 1_{H};GX_{2},x_{1})+B(gx_{1}x_{2}\otimes
x_{1};GX_{1}X_{2},x_{1})=0.
\end{equation*}%
In view of the form of the elements we get%
\begin{equation*}
B(x_{1}\otimes 1_{H};G,x_{1})+B(gx_{1}x_{2}\otimes 1_{H};G,x_{1}x_{2})=0
\end{equation*}%
which is $\left( \ref{X1,x2,X1F71,x1}\right) .$

\subparagraph{Case $f=x_{2}$}

\begin{equation*}
-B(gx_{1}x_{2}\otimes 1_{H};GX_{2},x_{2})=-B(x_{1}x_{2}\otimes
g;GX_{2},x_{2})+B(gx_{1}x_{2}\otimes x_{1};GX_{1}X_{2},x_{2}).
\end{equation*}%
By applying $\left( \ref{eq.10}\right) $ this rewrites as

\begin{equation*}
2B(gx_{1}x_{2}\otimes 1_{H};GX_{2},x_{2})+B(gx_{1}x_{2}\otimes
x_{1};GX_{1}X_{2},x_{2})=0.
\end{equation*}%
In view of the form of the elements we get

\begin{equation*}
B(x_{1}\otimes 1_{H};G,x_{2})=0
\end{equation*}%
which is $\left( \ref{X1,g,X1F21,x2}\right) .$

\subparagraph{Case $f=gx_{1}x_{2}$}

\begin{eqnarray*}
&&-B(gx_{1}x_{2}\otimes 1_{H};GX_{2},gx_{1}x_{2}) \\
&=&-B(x_{1}x_{2}\otimes g;GX_{2},gx_{1}x_{2})+B(gx_{1}x_{2}\otimes
x_{1};GX_{1}X_{2},gx_{1}x_{2}).
\end{eqnarray*}%
By applying $\left( \ref{eq.10}\right) $ this rewrites as%
\begin{equation*}
B(gx_{1}x_{2}\otimes x_{1};GX_{1}X_{2},gx_{1}x_{2})=0.
\end{equation*}

In view of the form of the elements we get%
\begin{equation*}
B(g\otimes 1_{H};G,gx_{2})-B(x_{1}\otimes 1_{H};G,gx_{1}x_{2})=0
\end{equation*}%
which is $\left( \ref{X1,g,X1F21,gx1x2}\right) .$

\subsection{$d=X_{2}$}

\begin{equation*}
B(h\otimes h^{\prime })(X_{2}\otimes 1_{H}) =(X_{2}\otimes 1_{H})B(hg\otimes
h^{\prime }g)+ B(hx_{2}\otimes h^{\prime }g)+ B(h\otimes h^{\prime }x_{2})
\end{equation*}

\begin{eqnarray*}
B(h\otimes h^{\prime }) &=&\sum_{f}[B(h\otimes h^{\prime
};1_{A},f)1_{A}\otimes f+B(h\otimes h^{\prime };G,f)G\otimes f+ \\
&&B(h\otimes h^{\prime };X_{1},f)X_{1}\otimes f+B(h\otimes h^{\prime
};X_{2},f)X_{2}\otimes f. \\
&&+B(h\otimes h^{\prime };X_{1}X_{2},f)X_{1}X_{2}\otimes f+B(h\otimes
h^{\prime };GX_{1},f)GX_{1}\otimes f+ \\
&&B(h\otimes h^{\prime };GX_{2},f)GX_{2}\otimes f+B(h\otimes h^{\prime
};GX_{1}X_{2},f)GX_{1}X_{2}\otimes ]f
\end{eqnarray*}

We recall that%
\begin{eqnarray*}
1_{A}X_{2} &=&X_{2} \\
GX_{2} &=&GX_{2} \\
X_{1}X_{2} &=&X_{1}X_{2} \\
X_{2}X_{2} &=&\beta _{2} \\
X_{1}X_{2}X_{2} &=&\beta _{2}X_{1} \\
GX_{1}X_{2} &=&GX_{1}X_{2} \\
GX_{2}X_{2} &=&\beta _{2}G \\
GX_{1}X_{2}X_{2} &=&\beta _{2}GX_{1}
\end{eqnarray*}

We write the left side%
\begin{eqnarray*}
B(h\otimes h^{\prime })\left( X_{2}\otimes 1_{H}\right) &=&\sum_{f}\beta
_{2}B(h\otimes h^{\prime };X_{2},f)1_{A}\otimes f+ \\
&&+\sum_{f}\beta _{2}B(h\otimes h^{\prime };GX_{2},f)G\otimes f+ \\
&&+\sum_{f}\beta _{2}B(h\otimes h^{\prime };X_{1}X_{2},f)X_{1}\otimes f+ \\
&&+\sum_{f}B(h\otimes h^{\prime };1_{H},f)X_{2}\otimes f+ \\
&&+\sum_{f}B(h\otimes h^{\prime };X_{1},f)X_{1}X_{2}\otimes f+ \\
&&+\sum_{f}\beta _{2}B(h\otimes h^{\prime };GX_{1}X_{2},f)GX_{1}\otimes f+ \\
&&+\sum_{f}B(h\otimes h^{\prime };G,f)GX_{2}\otimes f \\
&&+\sum_{f}B(h\otimes h^{\prime };GX_{1},f)GX_{1}X_{2}\otimes f
\end{eqnarray*}

We write the right side

$G^{2}=\alpha $, $X_{i}^{2}=\beta _{i}$ and $X_{i}G+GX_{i}=\gamma _{i}$ and $%
X_{1}X_{2}+X_{2}X_{1}=\lambda $.

\begin{eqnarray*}
X_{2}1_{A} &=&X_{2} \\
X_{2}G& =&-GX_{2}+\gamma_2 \\
X_{2}X_{1} &=&-X_1X_2+\lambda \\
X_{2}X_{2} &=&\beta_{2} \\
X_{2}X_{1}X_{2} &=& -\beta_2X_1 +\lambda X_2 \\
X_{2}GX_{1}&=&-GX_2X_1+\gamma_2X_1 = GX_1X_2 - \lambda G +\gamma_2X_1 \\
X_{2}GX_{2} & = & -\beta_2 G+ \gamma_2X_2 \\
X_{2}GX_{1}X_{2} &. = & \beta_2GX_1 - \lambda GX_2 +\gamma_2X_1X_2
\end{eqnarray*}

\begin{eqnarray*}
&&\left( X_{2}\otimes 1_{H}\right) B(hg\otimes h^{\prime }g)= \\
\sum_{f} &&[\gamma _{2}B(hg\otimes h^{\prime }g;G,f)+\lambda B(hg\otimes
h^{\prime }g;X_{1},f)+\beta _{2}B(hg\otimes h^{\prime
}g;X_{2},f)]1_{A}\otimes f \\
+\sum_{f} &&[-\beta _{2}B(hg\otimes h^{\prime }g;GX_{2},f)-\lambda
B(hg\otimes h^{\prime }g;GX_{1},f)]G\otimes f \\
+\sum_{f} &&[-\beta _{2}B(hg\otimes h^{\prime }g;X_{1}X_{2},f)+\gamma
_{2}B(hg\otimes h^{\prime }g;GX_{1},f)]X_{1}\otimes f \\
+\sum_{f} &&[B(hg\otimes h^{\prime }g;1_{A},f)+\gamma _{2}B(hg\otimes
h^{\prime }g;GX_{2},f)+\lambda B(hg\otimes h^{\prime
}g;X_{1}X_{2},f)]X_{2}\otimes f \\
+\sum_{f} &&[-B(hg\otimes h^{\prime }g;X_{1},f)+\gamma _{2}B(hg\otimes
h^{\prime }g;GX_{1}X_{2},f)]X_{1}X_{2}\otimes f \\
+\sum_{f} &&+\beta _{2}B(hg\otimes h^{\prime }g;GX_{1}X_{2},f)GX_{1}\otimes f
\\
+\sum_{f} &&[-B(hg\otimes h^{\prime }g;G,f)-\lambda B(hg\otimes h^{\prime
}g;GX_{1}X_{2},f)]GX_{2}\otimes f \\
+\sum_{f} &&B(hg\otimes h^{\prime }g;GX_{1},f)GX_{1}X_{2}\otimes f
\end{eqnarray*}

so we obtain

\begin{eqnarray}
\beta _{2}B(h\otimes h^{\prime };X_{2},f) &=&\gamma _{2}B(hg\otimes
h^{\prime }g;G,f)+\lambda B(hg\otimes h^{\prime }g;X_{1},f)+\beta
_{2}B(hg\otimes h^{\prime }g;X_{2},f)  \notag \\
&&+B(hx_{2}\otimes h^{\prime }g;1_{A},f)+B(h\otimes h^{\prime
}x_{2};1_{A},f).  \label{X2F11}
\end{eqnarray}

\begin{eqnarray}
\beta _{2}B(h\otimes h^{\prime };GX_{2},f) &=&-\beta _{2}B(hg\otimes
h^{\prime }g;GX_{2},f)-\lambda B(hg\otimes h^{\prime }g;GX_{1},f)
\label{X2F21} \\
&&+B(hx_{2}\otimes h^{\prime }g;G,f)+B(h\otimes h^{\prime }x_{2};G,f).
\notag
\end{eqnarray}

\begin{eqnarray}
\beta _{2}B(h\otimes h^{\prime };X_{1}X_{2},f) &=&-\beta _{2}B(hg\otimes
h^{\prime }g;X_{1}X_{2},f)+\gamma _{2}B(hg\otimes h^{\prime }g;GX_{1},f)
\label{X2F31} \\
&&+B(hx_{2}\otimes h^{\prime }g;X_{1},f)+B(h\otimes h^{\prime
}x_{2};X_{1},f).  \notag
\end{eqnarray}

\begin{eqnarray}
B(h\otimes h^{\prime };1_{A},f) &=&B(hg\otimes h^{\prime }g;1_{A},f)+\gamma
_{2}B(hg\otimes h^{\prime }g;GX_{2},f)+\lambda B(hg\otimes h^{\prime
}g;X_{1}X_{2},f)  \label{X2F41} \\
&&+B(hx_{2}\otimes h^{\prime }g;X_{2},f)+B(h\otimes h^{\prime
}x_{2};X_{2},f).  \notag
\end{eqnarray}

\begin{eqnarray}
B(h\otimes h^{\prime };X_{1},f) &=&-B(hg\otimes h^{\prime }g;X_{1},f)+\gamma
_{2}B(hg\otimes h^{\prime }g;GX_{1}X_{2},f)  \label{X2F51} \\
&&+B(hx_{2}\otimes h^{\prime }g;X_{1}X_{2},f)+B(h\otimes h^{\prime
}x_{2};X_{1}X_{2},f)  \notag
\end{eqnarray}

\begin{gather}
\beta _{2}B(h\otimes h^{\prime };GX_{1}X_{2},f)=+\beta _{2}B(hg\otimes
h^{\prime }g;GX_{1}X_{2},f)  \notag \\
+B(hx_{2}\otimes h^{\prime }g;GX_{1},f)+B(h\otimes h^{\prime
}x_{2};GX_{1},f).  \label{X2F61}
\end{gather}

\begin{eqnarray}
B(h\otimes h^{\prime };G,f) &=&-B(hg\otimes h^{\prime }g;G,f)-\lambda
B(hg\otimes h^{\prime }g;GX_{1}X_{2},f)  \label{X2F71} \\
&&+B(hx_{2}\otimes h^{\prime }g;GX_{2},f)+B(h\otimes h^{\prime
}x_{2};GX_{2},f).  \notag
\end{eqnarray}

\begin{gather}
B(h\otimes h^{\prime };GX_{1},f)=B(hg\otimes h^{\prime }g;GX_{1},f)
\label{X2F81} \\
+B(hx_{2}\otimes h^{\prime }g;GX_{1}X_{2},f)+B(h\otimes h^{\prime
}x_{2};GX_{1}X_{2},f).  \notag
\end{gather}

\subsubsection{Case $g\otimes 1_{H}$}

\paragraph{\textbf{Equality }$\left( \protect\ref{X2F11}\right) $}

rewrites as

\begin{gather}
\beta _{2}B(g\otimes 1_H );X_{2},f) = \gamma_2 B(1_H\otimes g;G,f) + \lambda
B(1_H\otimes g;X_1,f)  \notag \\
+ \beta_2B(1_H\otimes g;X_2,f) +B(gx_{2}\otimes g;1_{A},f)+B(g\otimes
x_{2};1_{A},f)  \notag
\end{gather}

\subparagraph{Case $f=1_{H}$}

\begin{gather}
\beta _{2}B(g\otimes 1_H );X_{2},1_H) = \gamma_2 B(1_H\otimes g;G,1_H) +
\lambda B(1_H\otimes g;X_1,1_H)  \notag \\
+ \beta_2B(1_H\otimes g;X_2,1_H) +B(gx_{2}\otimes g;1_{A},1_H)+B(g\otimes
x_{2};1_{A},1_H)  \notag
\end{gather}

By applying $\left( \ref{eq.10}\right) $ this rewrites as

\begin{gather}
\gamma_2 B(g\otimes 1_H;G,1_H) + \lambda B(g\otimes 1_H;X_1,1_H)  \notag \\
+B(x_{2}\otimes 1_H;1_{A},1_H) +B(1_H\otimes gx_{2};1_{A},1_H) =0  \notag
\end{gather}

In view of the form of the elements we get

\begin{gather}
\gamma _{2}B(g\otimes 1_{H};G,1_{H})+\lambda B(g\otimes 1_{H};1_{A},x_{1})
\label{X2,g,X2F11,1H} \\
+2B(x_{2}\otimes 1_{H};1_{A},1_{H})=0  \notag
\end{gather}

\subparagraph{Case $f=x_1x_2$}

\begin{gather}
\beta _{2}B(g\otimes 1_H );X_{2},x_1x_2) = \gamma_2 B(1_H\otimes g;G,x_1x_2)
+ \lambda B(1_H\otimes g;X_1,x_1x_2)  \notag \\
+ \beta_2B(1_H\otimes g;X_2,x_1x_2) +B(gx_{2}\otimes
g;1_{A},x_1x_2)+B(g\otimes x_{2};1_{A},x_1x_2)  \notag
\end{gather}

By applying $\left( \ref{eq.10}\right) $ this rewrites as

\begin{gather}
\gamma_2 B(g\otimes 1_H;G,x_1x_2) + \lambda B(g\otimes 1_H;X_1,x_1x_2)
\notag \\
+2B(x_{2}\otimes 1_H;1_{A},x_1x_2)+B(1_H\otimes gx_{2};1_{A},x_1x_2)=0
\notag
\end{gather}

In view of the form of the elements we get

\begin{equation*}
\gamma _{2}B(g\otimes 1_{H};G,x_{1}x_{2})+2B(x_{2}\otimes
1_{H};1_{A},x_{1}x_{2})=0
\end{equation*}%
which is $\left( \ref{G,x2, GF2,x1x2}\right) .$

\subparagraph{Case $f=gx_1$}

\begin{gather}
\beta _{2}B(g\otimes 1_H );X_{2},gx_1) = \gamma_2 B(1_H\otimes g;G,gx_1) +
\lambda B(1_H\otimes g;X_1,gx_1)  \notag \\
+ \beta_2B(1_H\otimes g;X_2,gx_1) +B(gx_{2}\otimes g;1_{A},gx_1)+B(g\otimes
x_{2};1_{A},gx_1)  \notag
\end{gather}

By applying $\left( \ref{eq.10}\right) $ this rewrites as

\begin{gather}
2\beta _{2}B(g\otimes 1_H );X_{2},gx_1) +\gamma_2 B(g\otimes 1_H;G,gx_1) +
\lambda B(g\otimes 1_H;X_1,gx_1)  \notag \\
+B(x_{2}\otimes 1_H;1_{A},gx_1)+B(1_H\otimes gx_{2};1_{A},gx_1) =0  \notag
\end{gather}

In view of the form of the elements we get

\begin{gather}
2\beta _{2}B(g\otimes 1_{H};1_{A},gx_{1}x_{2})+\gamma _{2}B(g\otimes
1_{H};G,gx_{1})  \label{X2,g,X2F11,gx1} \\
+2B(x_{2}\otimes 1_{H};1_{A},gx_{1})=0  \notag
\end{gather}

\subparagraph{Case $f=gx_{2}$}

\begin{gather}
\beta _{2}B(g\otimes 1_{H};X_{2},gx_{2})=\gamma _{2}B(1_{H}\otimes
g;G,gx_{2})+\lambda B(1_{H}\otimes g;X_{1},gx_{2})  \notag \\
+\beta _{2}B(1_{H}\otimes g;X_{2},gx_{2})+B(gx_{2}\otimes
g;1_{A},gx_{2})+B(g\otimes x_{2};1_{A},gx_{2})  \notag
\end{gather}

By applying $\left( \ref{eq.10}\right) $ this rewrites as

\begin{gather}
2\beta _{2}B(g\otimes 1_{H};X_{2},gx_{2})+\gamma _{2}B(g\otimes
1_{H};G,gx_{2})+\lambda B(g\otimes 1_{H};X_{1},gx_{2})  \notag \\
+B(x_{2}\otimes 1_{H};1_{A},gx_{2})+B(1_{H}\otimes gx_{2};1_{A},gx_{2})=0
\notag
\end{gather}

In view of the form of the elements we get

\begin{gather}
\gamma _{2}B(g\otimes 1_{H};G,gx_{2})-\lambda B(g\otimes
1_{H};1_{A},gx_{1}x_{2})  \label{X2,g,X2F11,gx2} \\
+2B(x_{2}\otimes 1_{H};1_{A},gx_{2})+2B(g\otimes 1_{H};1_{A},g)=0.  \notag
\end{gather}

\paragraph{Equality $\left( \protect\ref{X2F21}\right) $}

rewrites as

\begin{eqnarray}
\beta_{2}B(g\otimes 1_H;GX_{2},f) &=& -\beta_2 B(1_H\otimes g;GX_2,f) -
\lambda B(1_H\otimes g;GX_1,f)  \notag \\
&&+B(gx_{2}\otimes g;G,f)+B(g\otimes x_{2};G,f)  \notag
\end{eqnarray}

\subparagraph{Case $f=g$}

\begin{eqnarray}
\beta_{2}B(g\otimes 1_H;GX_{2},g) &=& -\beta_2 B(1_H\otimes g;GX_2,g) -
\lambda B(1_H\otimes g;GX_1,g)  \notag \\
&&+B(gx_{2}\otimes g;G,g)+B(g\otimes x_{2};G,g)  \notag
\end{eqnarray}

By applying $\left( \ref{eq.10}\right) $ this rewrites as

\begin{eqnarray*}
2\beta_{2}B(g\otimes 1_H;GX_{2},g) &=& - \lambda B(g\otimes 1_H;GX_1,g) \\
&&+B(x_{2}\otimes 1_H;G,g)+B(1_H\otimes gx_{2};G,g)
\end{eqnarray*}

In view of the form of the elements we get

\begin{equation}
2\beta _{2}B(g\otimes 1_{H};G,gx_{2})+\lambda B(g\otimes
1_{H};G,gx_{1})+2B(x_{2}\otimes 1_{H};G,g)=0  \label{X2,g,X2F21,g}
\end{equation}

\subparagraph{Case $f=x_1$}

\begin{eqnarray}
\beta_{2}B(g\otimes 1_H;GX_{2},x_1) &=& -\beta_2 B(1_H\otimes g;GX_2,x_1) -
\lambda B(1_H\otimes g;GX_1,x_1)  \notag \\
&&+B(gx_{2}\otimes g;G,x_1)+B(g\otimes x_{2};G,x_1)  \notag
\end{eqnarray}

By applying $\left( \ref{eq.10}\right) $ this rewrites as

\begin{equation}
\lambda B(g\otimes 1_H;GX_1,x_1) -B(x_{2}\otimes 1_H;G,x_1)-B(1_H\otimes
gx_{2};G,x_1) =0  \notag
\end{equation}

In view of the form of the elements we get

\begin{equation}
B(x_{2}\otimes 1_{H};G,x_{1})=0  \label{X2,g,X2F21,x1}
\end{equation}

\subparagraph{Case $f=x_2$}

\begin{eqnarray}
\beta_{2}B(g\otimes 1_H;GX_{2},x_2) &=& -\beta_2 B(1_H\otimes g;GX_2,x_2) -
\lambda B(1_H\otimes g;GX_1,x_2)  \notag \\
&&+B(gx_{2}\otimes g;G,x_2)+B(g\otimes x_{2};G,x_2)  \notag
\end{eqnarray}

By applying $\left( \ref{eq.10}\right) $ this rewrites as

\begin{equation}
\lambda B(g\otimes 1_H;GX_1,x_2) -B(x_{2}\otimes 1_H;G,x_2)-B(1_H\otimes
gx_{2};G,x_2) =0  \notag
\end{equation}

In view of the form of the elements we get

\begin{equation}
\lambda B(g\otimes 1_{H};G,x_{1}x_{2})-2B(x_{2}\otimes 1_{H};G,x_{2})=0
\label{X2,g,X2F21,x2}
\end{equation}

\paragraph{Equality $\left( \protect\ref{X2F31}\right) $}

rewrites as

\begin{eqnarray*}
\beta _{2}B(g\otimes 1_{H};X_{1}X_{2},f) &=&-\beta _{2}B(1_{H}\otimes
g;X_{1}X_{2},f)+\gamma _{2}B(1_{H}\otimes g;GX_{1},f) \\
&&+B(gx_{2}\otimes g;X_{1},f)+B(g\otimes x_{2};X_{1},f)
\end{eqnarray*}

\subparagraph{Case $f=g$}

\begin{eqnarray*}
\beta _{2}B(g\otimes 1_{H};X_{1}X_{2},g) &=&-\beta _{2}B(1_{H}\otimes
g;X_{1}X_{2},g)+\gamma _{2}B(1_{H}\otimes g;GX_{1},g) \\
&&+B(gx_{2}\otimes g;X_{1},g)+B(g\otimes x_{2};X_{1},g)
\end{eqnarray*}

By applying $\left( \ref{eq.10}\right) $ this rewrites as

\begin{gather*}
2\beta _{2}B(g\otimes 1_{H};X_{1}X_{2},g)-\gamma _{2}B(g\otimes
1_{H};GX_{1},g) \\
-B(x_{2}\otimes 1_{H};X_{1},g)-B(1_{H}\otimes gx_{2};X_{1},g)=0
\end{gather*}

In view of the form of the element we get

\begin{equation}
2\beta _{2}B(g\otimes 1_{H};1_{H},gx_{1}x_{2})+\gamma _{2}B(g\otimes
1_{H};G,gx_{1})+2B(x_{2}\otimes 1_{H};1_{H},gx_{1})=0.  \label{X2,g,X2F31,g}
\end{equation}

\subparagraph{Case $f=x_2$}

\begin{eqnarray*}
+\gamma_2 B(1_H\otimes g;GX_1,x_2) +B(gx_{2}\otimes g;X_1,x_2)+B(g\otimes
x_{2};X_1,x_2)& = & 0
\end{eqnarray*}

By applying $\left( \ref{eq.10}\right) $ this rewrites as

\begin{equation*}
-\gamma _{2}B(g\otimes 1_{H};GX_{1},x_{2})-B(x_{2}\otimes
1_{H};X_{1},x_{2})-B(1_{H}\otimes gx_{2};X_{1},x_{2})=0
\end{equation*}

In view of the form of the element we get

\begin{equation}
\gamma _{2}B(g\otimes 1_{H};G,x_{1}x_{2})+2B(x_{2}\otimes
1_{H};1,x_{1}x_{2})=0  \label{X2,g,X2F31,x2}
\end{equation}

\paragraph{Equality $\left( \protect\ref{X2F41}\right) $}

rewrites as

\begin{eqnarray*}
B(g\otimes 1_H;1_A,f) &=& B(1_H\otimes g;1_A,f) +\gamma_2 B(1_H\otimes
g;GX_2,f) +\lambda B(1_H\otimes g;X_1X_2,f) \\
&&+B(gx_{2}\otimes g;X_2,f) +B(g\otimes x_{2};X_2,f)
\end{eqnarray*}

\subparagraph{Case $f=g$}

\begin{eqnarray*}
\gamma_2 B(1_H\otimes g;GX_2,g) +\lambda B(1_H\otimes g;X_1X_2,g)
+B(gx_{2}\otimes g;X_2,g) && \\
+B(g\otimes x_{2};X_2,g)&=&0
\end{eqnarray*}

By applying $\left( \ref{eq.10}\right) $ this rewrites as

\begin{eqnarray*}
\gamma_2 B(g\otimes 1_H;GX_2,g) +\lambda B(g\otimes 1_H;X_1X_2,g)
+B(x_{2}\otimes 1_H;X_2,g) && \\
+B(1\otimes gx_{2};X_2,g)&=&0
\end{eqnarray*}

In view of the form of the element we get

\begin{eqnarray*}
-\gamma _{2}B(g\otimes 1_{H};G,gx_{2})+\lambda B(g\otimes
1_{H};1_{A},gx_{1}x_{2}) && \\
-2B(g\otimes 1_{H};1_{A},g)-2B(x_{2}\otimes 1_{H};1_{A},gx_{2}) &=&0
\end{eqnarray*}%
This is $\left( \ref{X2,g,X2F11,gx2}\right) .$

\subparagraph{Case $f=x_1$}

\begin{eqnarray*}
B(g\otimes 1_H;1_A,x_1)= B(1_H\otimes g;1_A,x_1)+\gamma_2 B(1_H\otimes
g;GX_2,x_1) \\
+\lambda B(1_H\otimes g;X_1X_2,x_1) +B(gx_{2}\otimes g;X_2,x_1) +B(g\otimes
x_{2};X_2,x_1)
\end{eqnarray*}

By applying $\left( \ref{eq.10}\right) $ this rewrites as

\begin{eqnarray*}
2B(g\otimes 1_H;1_A,x_1)+\gamma_2 B(g\otimes 1_H;GX_2,x_1) -\lambda
B(g\otimes 1_H;X_1X_2,x_1) && \\
+B(x_{2}\otimes 1_H;X_2,x_1) -B(1_H\otimes gx_{2};X_2,x_1)&=&0
\end{eqnarray*}

In view of the form of the element we get

\begin{eqnarray}
3B(g\otimes 1_{H};1_{A},x_{1})-\gamma _{2}B(g\otimes 1_{H};G,x_{1}x_{2}) &&
\label{X2,g,X2F41,x1} \\
+B(x_{2}\otimes 1_{H};X_{2},x_{1})+B\left( x_{2}\otimes
1_{H};1_{A},x_{1}x_{2}\right) &=&0  \notag
\end{eqnarray}

\subparagraph{Case $f=x_2$}

\begin{eqnarray*}
B(g\otimes 1_H;1_A,x_2)= B(1_H\otimes g;1_A,x_2)+\gamma_2 B(1_H\otimes
g;GX_2,x_2) \\
+\lambda B(1_H\otimes g;X_1X_2,x_2) +B(gx_{2}\otimes g;X_2,x_2) +B(g\otimes
x_{2};X_2,x_2)
\end{eqnarray*}

By applying $\left( \ref{eq.10}\right) $ this rewrites as

\begin{eqnarray*}
2B(g\otimes 1_H;1_A,x_2)+\gamma_2 B(g\otimes 1_H;GX_2,x_2) -\lambda
B(g\otimes 1_H;X_1X_2,x_2) && \\
+B(x_{2}\otimes 1_H;X_2,x_2) -B(1_H\otimes gx_{2};X_2,x_2)&=&0
\end{eqnarray*}

In view of the form of the elements we get nothing new.

\subparagraph{Case $f=gx_1x_2$}

\begin{eqnarray*}
B(g\otimes 1_H;1_A,gx_1x_2)= B(1_H\otimes g;1_A,gx_1x_2)+\gamma_2
B(1_H\otimes g;GX_2,gx_1x_2) \\
+\lambda B(1_H\otimes g;X_1X_2,gx_1x_2) +B(gx_{2}\otimes g;X_2,gx_1x_2)
+B(g\otimes x_{2};X_2,gx_1x_2)
\end{eqnarray*}

By applying $\left( \ref{eq.10}\right) $ this rewrites as

\begin{eqnarray*}
2B(g\otimes 1_H;1_A,gx_1x_2)+\gamma_2 B(g\otimes 1_H;GX_2,gx_1x_2) -\lambda
B(g\otimes 1_H;X_1X_2,gx_1x_2) && \\
+B(x_{2}\otimes 1_H;X_2,gx_1x_2) -B(1_H\otimes gx_{2};X_2,gx_1x_2)&=&0
\end{eqnarray*}

In view of the form of the elements we get nothing new

\paragraph{Equality $\left( \protect\ref{X2F51}\right) $}

rewrites as

\begin{eqnarray*}
B(g\otimes 1_H;X_{1},f) &=& - B(1_H\otimes g;X_1,f) +\gamma_2 B(1_H\otimes
g;GX_1X_2,f)  \notag \\
&&+B(gx_{2}\otimes g;X_1X_2,f)+B(g\otimes x_{2};X_1X_2,f)
\end{eqnarray*}

\subparagraph{Case $f=1_H$}

\begin{eqnarray*}
B(g\otimes 1_H;X_{1},1_H) &=& - B(1_H\otimes g;X_1,1_H) +\gamma_2
B(1_H\otimes g;GX_1X_2,1_H)  \notag \\
&&+B(gx_{2}\otimes g;X_1X_2,1_H)+B(g\otimes x_{2};X_1X_2,1_H)
\end{eqnarray*}

By applying $\left( \ref{eq.10}\right) $ this rewrites as

\begin{eqnarray*}
2B(g\otimes 1_H;X_{1},1_H) -\gamma_2 B(g\otimes 1_H;GX_1X_2,1_H)  \notag \\
-B(x_{2}\otimes 1_H;X_1X_2,1_H)-B(1_H\otimes gx_{2};X_1X_2,1_H) =0
\end{eqnarray*}

In view of the form of the element we get

\begin{equation*}
\gamma _{2}B(g\otimes 1_{H};G,x_{1}x_{2})+2B(x_{2}\otimes
1_{H};1_{A},x_{1}x_{2})=0
\end{equation*}%
which is $\left( \ref{G,x2, GF2,x1x2}\right) .$

\subparagraph{Case $f=gx_2$}

\begin{eqnarray*}
B(g\otimes 1_H;X_{1},gx_2) &=& -B(1_H\otimes g;X_1,gx_2) +\gamma_2
B(1_H\otimes g;GX_1X_2,gx_2)  \notag \\
&&+B(gx_{2}\otimes g;X_1X_2,gx_2)+B(g\otimes x_{2};X_1X_2,gx_2)
\end{eqnarray*}

By applying $\left( \ref{eq.10}\right) $ this rewrites as

\begin{equation*}
-\gamma _{2}B(g\otimes 1_{H};GX_{1}X_{2},gx_{2})-B(x_{2}\otimes
1_{H};X_{1}X_{2},gx_{2})-B(1_{H}\otimes gx_{2};X_{1}X_{2},gx_{2})=0
\end{equation*}

In view of the form of the elements this is trivial.

\paragraph{Equality $\left( \protect\ref{X2F61}\right) $}

rewrites as

\begin{eqnarray*}
\beta_2 B(g\otimes 1_H;GX_1X_2,f) &=& +\beta_2 B(1_H\otimes g;GX_1X_2,f)
\notag \\
&&+B(gx_{2}\otimes g;GX_1,f)+B(g\otimes x_{2};GX_1,f)
\end{eqnarray*}

\subparagraph{Case $f=1_H$}

\begin{eqnarray*}
\beta_2 B(g\otimes 1_H;GX_1X_2,1_H) &=& +\beta_2 B(1_H\otimes g;GX_1X_2,1_H)
\notag \\
&&+B(gx_{2}\otimes g;GX_1,1_H)+B(g\otimes x_{2};GX_1,1_H)
\end{eqnarray*}

By applying $\left( \ref{eq.10}\right) $ this rewrites as

\begin{eqnarray*}
B(x_{2}\otimes 1_H;GX_1,1_H)+B(1_H\otimes gx_{2};GX_1,1_H)=0
\end{eqnarray*}

In view of the form of the element we get

\begin{equation*}
B(x_{2}\otimes 1_{H};G,x_{1})=0
\end{equation*}%
which is $\left( \ref{X2,g,X2F21,x1}\right) .$

\subparagraph{Case $f=gx_2$}

\begin{eqnarray*}
\beta _{2}B(g\otimes 1_{H};GX_{1}X_{2},gx_{2}) &=&+\beta _{2}B(1_{H}\otimes
g;GX_{1}X_{2},gx_{2}) \\
&&+B(gx_{2}\otimes g;GX_{1},gx_{2})+B(g\otimes x_{2};GX_{1},gx_{2})
\end{eqnarray*}

By applying $\left( \ref{eq.10}\right) $ this rewrites as

\begin{eqnarray*}
+B(x_{2}\otimes 1_H;GX_1,gx_2)+B(1_H\otimes gx_{2};GX_1,gx_2) =0
\end{eqnarray*}

In view of the form of the element we get

\begin{eqnarray*}
-B(x_{2}\otimes 1_H;X_2,gx_1x_2)-2B(g\otimes 1_H,G,gx_1) -B(x_2\otimes
1_H;G,gx_1x_2) =0
\end{eqnarray*}

thus

\begin{equation}
B\left( g\otimes 1_{H};1_{A},gx_{1}x_{2}\right) +B(g\otimes 1_{H},G,gx_{1})=0
\label{X2,g,X2F61,gx2}
\end{equation}

\paragraph{Equality $\left( \protect\ref{X2F71}\right) $}

rewrites as

\begin{eqnarray*}
B(g\otimes 1_H;G,f) &=&-B(1_H\otimes g;G,f)- \lambda B(1_H\otimes
g;GX_1X_2,f)  \notag \\
&&+B(gx_{2}\otimes g;GX_2,f)+B(g\otimes x_{2};GX_2,f)
\end{eqnarray*}

\subparagraph{Case $f=1_H$}

\begin{eqnarray*}
B(g\otimes 1_H;G,1_H) &=&-B(1_H\otimes g;G,1_H)- \lambda B(1_H\otimes
g;GX_1X_2,1_H)  \notag \\
&&+B(gx_{2}\otimes g;GX_2,1_H)+B(g\otimes x_{2};GX_2,1_H)
\end{eqnarray*}

By applying $\left( \ref{eq.10}\right) $ this rewrites as

\begin{eqnarray*}
2 B(g\otimes 1_H;G,1_H) - \lambda B(g\otimes 1_H;GX_1X_2,1_H) \\
+B(x_{2}\otimes 1_H;GX_2,1_H)+B(1_H\otimes gx_{2};GX_2,1_H)=0
\end{eqnarray*}

In view of the form of the element we get

\begin{gather}
2B(g\otimes 1_{H};G,1_{H})-\lambda B(g\otimes 1_{H};G,x_{1}x_{2})
\label{X2,g,X2F71,1H} \\
+2B(g\otimes 1_{H};G,1_{H})+2B(x_{2}\otimes 1_{H};G,x_{2})=0  \notag
\end{gather}

\subparagraph{Case $f=x_1x_2$}

\begin{eqnarray*}
B(g\otimes 1_H;G,x_1x_2) &=&-B(1_H\otimes g;G,x_1x_2)- \lambda B(1_H\otimes
g;GX_1X_2,x_1x_2)  \notag \\
&&+B(gx_{2}\otimes g;GX_2,x_1x_2)+B(g\otimes x_{2};GX_2,x_1x_2)
\end{eqnarray*}

By applying $\left( \ref{eq.10}\right) $ this rewrites as

\begin{eqnarray*}
2B(g\otimes 1_H;G,x_1x_2) +\lambda B(1_H\otimes g;GX_1X_2,x_1x_2) \\
-B(x_{2}\otimes 1_H;GX_2,x_1x_2)-B(1_H\otimes gx_{2};GX_2,x_1x_2)=0
\end{eqnarray*}

In view of the form of the element we get

\begin{equation*}
B(g\otimes 1_{H};G,x_{1}x_{2})-B\left( g\otimes 1_{H};G,x_{1}x_{2}\right) =0
\end{equation*}%
which is trivial.

\subparagraph{Case $f=gx_1$}

\begin{eqnarray*}
B(g\otimes 1_H;G,gx_1) &=&-B(1_H\otimes g;G,gx_1)- \lambda B(1_H\otimes
g;GX_1X_2,gx_1)  \notag \\
&&+B(gx_{2}\otimes g;GX_2,gx_1)+B(g\otimes x_{2};GX_2,gx_1)
\end{eqnarray*}

By applying $\left( \ref{eq.10}\right) $ this rewrites as

\begin{eqnarray*}
- \lambda B(g\otimes 1_H;GX_1X_2,gx_1) +B(x_{2}\otimes
1_H;GX_2,gx_1)+B(1_H\otimes gx_{2};GX_2,gx_1)=0
\end{eqnarray*}

In view of the form of the element we get

\begin{equation}
B(g\otimes 1_{H};G,gx_{1})+B(x_{2}\otimes 1_{H};G,gx_{1}x_{2})=0
\label{X2,g,X2F71,gx1}
\end{equation}

\subparagraph{Case $f=gx_2$}

\begin{eqnarray*}
B(g\otimes 1_H;G,gx_2) &=&-B(1_H\otimes g;G,gx_2)- \lambda B(1_H\otimes
g;GX_1X_2,gx_2)  \notag \\
&&+B(gx_{2}\otimes g;GX_2,gx_2)+B(g\otimes x_{2};GX_2,gx_2)
\end{eqnarray*}

By applying $\left( \ref{eq.10}\right) $ this rewrites as

\begin{eqnarray*}
- \lambda B(g\otimes 1_H;GX_1X_2,gx_2) +B(x_{2}\otimes
1_H;GX_2,gx_2)+B(1_H\otimes gx_{2};GX_2,gx_2)=0
\end{eqnarray*}

In view of the form of the elements this equality is trivial.

\paragraph{Equality $\left( \protect\ref{X2F81}\right) $}

rewrites as

\begin{gather*}
B(g\otimes 1_{H};GX_{1},f)=B(1_{H}\otimes g;GX_{1},f) \\
+B(gx_{2}\otimes g;GX_{1}X_{2},f)+B(g\otimes x_{2};GX_{1}X_{2},f)
\end{gather*}

\subparagraph{Case $f=g$}

\begin{eqnarray*}
B(g\otimes 1_H;GX_{1},g) &=&B(1_H\otimes g;GX_1,g)  \notag \\
&&+B(gx_{2}\otimes g;GX_1X_2,g)+B(g\otimes x_{2};GX_1X_2,g)
\end{eqnarray*}

By applying $\left( \ref{eq.10}\right) $ this rewrites as

\begin{eqnarray*}
B(x_{2}\otimes 1_H;GX_1X_2,g)+B(1_H\otimes gx_{2};GX_1X_2,g) =0
\end{eqnarray*}

In view of the form of the elements we get

\begin{equation*}
B(g\otimes 1_{H};G,gx_{1})+B(x_{2}\otimes 1_{H};G,gx_{1}x_{2})=0
\end{equation*}%
which is $\left( \ref{X2,g,X2F71,gx1}\right) .$

\subparagraph{Case $f=x_2$}

\begin{gather*}
B(g\otimes 1_{H};GX_{1},x_{2})=B(1_{H}\otimes g;GX_{1},x_{2}) \\
+B(gx_{2}\otimes g;GX_{1}X_{2},x_{2})+B(g\otimes x_{2};GX_{1}X_{2},x_{2})
\end{gather*}

By applying $\left( \ref{eq.10}\right) $ this rewrites as
\begin{eqnarray*}
2B(g\otimes 1_H;GX_{1},x_2) -B(gx_{2}\otimes g;GX_1X_2,x_2)-B(g\otimes
x_{2};GX_1X_2,x_2)=0
\end{eqnarray*}

In view of the form of the elements we get%
\begin{equation*}
2B(g\otimes 1_{H};G,x_{1}x_{2})-B\left( g\otimes 1_{H};G,x_{1}x_{2}\right)
-B\left( g\otimes 1_{H};G,x_{1}x_{2}\right) =0
\end{equation*}

which is trivial.

\subsubsection{Case $x_1\otimes 1_{H}$}

\paragraph{Equality $\left( \protect\ref{X2F11}\right) $}

rewrites as

\begin{eqnarray*}
\beta _{2}B(x_{1}\otimes 1_{H};X_{2},f) &=&-\gamma _{2}B(gx_{1}\otimes
g;G,f)-\lambda B(gx_{1}\otimes g;X_{1},f)-\beta _{2}B(gx_{1}\otimes
g;X_{2},f) \\
&&+B(x_{1}x_{2}\otimes g;1_{A},f)+B(x_{1}\otimes x_{2};1_{A},f)
\end{eqnarray*}

\subparagraph{Case $f=g$}

\begin{eqnarray*}
\beta _{2}B(x_{1}\otimes 1_{H};X_{2},g) &=&-\gamma _{2}B(gx_{1}\otimes
g;G,g)-\lambda B(gx_{1}\otimes g;X_{1},g)-\beta _{2}B(gx_{1}\otimes
g;X_{2},g) \\
&&+B(x_{1}x_{2}\otimes g;1_{A},g)+B(x_{1}\otimes x_{2};1_{A},g)
\end{eqnarray*}%
By applying $\left( \ref{eq.10}\right) $ this rewrites as%
\begin{eqnarray*}
2\beta _{2}B(x_{1}\otimes 1_{H};X_{2},g) &=&-\gamma _{2}B(x_{1}\otimes
1_{H};G,g)-\lambda B(x_{1}\otimes 1_{H};X_{1},g)) \\
&&+B(gx_{1}x_{2}\otimes 1_{H};1_{A},g)+B(x_{1}\otimes x_{2};1_{A},g)
\end{eqnarray*}%
In view of the form of the elements we get%
\begin{eqnarray}
&&-2\beta _{2}B(x_{1}\otimes 1_{H};1_{A},gx_{2})+\gamma _{2}B(x_{1}\otimes
1_{H};G,g)  \label{X2,x1,X2F11,g} \\
&&-\lambda \left[ B(g\otimes 1_{H};1_{A},g)+B(x_{1}\otimes
1_{H};1_{A},gx_{1})\right] -2B(gx_{1}x_{2}\otimes 1_{H};1_{A},g)=0  \notag
\end{eqnarray}

\subparagraph{Case $f=x_{1}$}

\begin{eqnarray*}
\beta _{2}B(x_{1}\otimes 1_{H};X_{2},x_{1}) &=&-\gamma _{2}B(gx_{1}\otimes
g;G,x_{1})-\lambda B(gx_{1}\otimes g;X_{1},x_{1})-\beta _{2}B(gx_{1}\otimes
g;X_{2},x_{1}) \\
&&+B(x_{1}x_{2}\otimes g;1_{A},x_{1})+B(x_{1}\otimes x_{2};1_{A},x_{1})
\end{eqnarray*}%
By applying $\left( \ref{eq.10}\right) $ this rewrites as%
\begin{eqnarray*}
&&+\gamma _{2}B(x_{1}\otimes 1_{H};G,x_{1})+\lambda B(x_{1}\otimes
1_{H};X_{1},x_{1}) \\
&&-B(gx_{1}x_{2}\otimes 1_{H};1_{A},x_{1})+B(x_{1}\otimes
x_{2};1_{A},x_{1})=0
\end{eqnarray*}%
In view of the form of the elements we get%
\begin{equation}
\gamma _{2}B(x_{1}\otimes 1_{H};G,x_{1})-\lambda B(g\otimes
1_{H};1_{A},x_{1})-2B(gx_{1}x_{2}\otimes 1_{H};1_{A},x_{1})=0
\label{X2,x1,X2F11,x1}
\end{equation}

\subparagraph{Case $f=x_{2}$}

\begin{eqnarray*}
\beta _{2}B(x_{1}\otimes 1_{H};X_{2},x_{2}) &=&-\gamma _{2}B(gx_{1}\otimes
g;G,x_{2})-\lambda B(gx_{1}\otimes g;X_{1},x_{2})-\beta _{2}B(gx_{1}\otimes
g;X_{2},x_{2}) \\
&&+B(x_{1}x_{2}\otimes g;1_{A},x_{2})+B(x_{1}\otimes x_{2};1_{A},x_{2})
\end{eqnarray*}%
By applying $\left( \ref{eq.10}\right) $ this rewrites as%
\begin{eqnarray*}
&&+\gamma _{2}B(x_{1}\otimes 1_{H};G,x_{2})+\lambda B(x_{1}\otimes
1_{H};X_{1},x_{2}) \\
&&-B(gx_{1}x_{2}\otimes 1_{H};1_{A},x_{2})+B(x_{1}\otimes
x_{2};1_{A},x_{2})=0.
\end{eqnarray*}%
In view of the form of the elements we get%
\begin{eqnarray}
&&-\gamma _{2}B(x_{1}\otimes 1_{H};G,x_{2})+\lambda \left[ B(g\otimes
1_{H};1_{A},x_{2})-B(x_{1}\otimes 1_{H};1_{A},x_{1}x_{2})\right]
\label{X2,x1,X2F11,x2} \\
&&+2B(gx_{1}x_{2}\otimes 1_{H};1_{A},x_{2})=0.  \notag
\end{eqnarray}

\subparagraph{Case $f=gx_{1}x_{2}$}

\begin{eqnarray*}
\beta _{2}B(x_{1}\otimes 1_{H};X_{2},gx_{1}x_{2}) &=&-\gamma
_{2}B(gx_{1}\otimes g;G,gx_{1}x_{2})-\lambda B(gx_{1}\otimes
g;X_{1},gx_{1}x_{2})-\beta _{2}B(gx_{1}\otimes g;X_{2},gx_{1}x_{2}) \\
&&+B(x_{1}x_{2}\otimes g;1_{A},gx_{1}x_{2})+B(x_{1}\otimes
x_{2};1_{A},gx_{1}x_{2})
\end{eqnarray*}%
By applying $\left( \ref{eq.10}\right) $ this rewrites as

\begin{eqnarray*}
&&2\beta _{2}B(x_{1}\otimes 1_{H};X_{2},gx_{1}x_{2})+\gamma
_{2}B(x_{1}\otimes 1_{H};G,gx_{1}x_{2})+\lambda B(x_{1}\otimes
1_{H};X_{1},gx_{1}x_{2}) \\
&&-B(gx_{1}x_{2}\otimes 1_{H};1_{A},gx_{1}x_{2})-B(x_{1}\otimes
x_{2};1_{A},gx_{1}x_{2})=0
\end{eqnarray*}%
In view of the form of the elements we get%
\begin{eqnarray}
&&+\gamma _{2}B(x_{1}\otimes 1_{H};G,gx_{1}x_{2})-\lambda B(g\otimes
1_{H};1_{A},gx_{1}x_{2})  \label{X2,x1,X2F11,gx1x2} \\
&&-2B(gx_{1}x_{2}\otimes 1_{H};1_{A},gx_{1}x_{2})-2B\left( x_{1}\otimes
1_{H};1_{A},gx_{1}\right) =0.  \notag
\end{eqnarray}

\paragraph{Equality $\left( \protect\ref{X2F21}\right) $}

rewrites as%
\begin{eqnarray*}
\beta _{2}B(x_{1}\otimes 1_{H};GX_{2},f) &=&\beta _{2}B(gx_{1}\otimes
g;GX_{2},f)+\lambda B(gx_{1}\otimes g;GX_{1},f) \\
&&+B(x_{1}x_{2}\otimes g;G,f)+B(x_{1}\otimes x_{2};G,f)
\end{eqnarray*}

\subparagraph{Case $f=1_{H}$}

\begin{eqnarray*}
\beta _{2}B(x_{1}\otimes 1_{H};GX_{2},1_{H}) &=&+\beta _{2}B(gx_{1}\otimes
g;GX_{2},1_{H})+\lambda B(gx_{1}\otimes g;GX_{1},1_{H}) \\
&&+B(x_{1}x_{2}\otimes g;G,1_{H})+B(x_{1}\otimes x_{2};G,1_{H}).
\end{eqnarray*}%
By applying $\left( \ref{eq.10}\right) $ this rewrites as%
\begin{equation*}
+\lambda B(x_{1}\otimes 1_{H};GX_{1},1_{H})+B(gx_{1}x_{2}\otimes
1_{H};G,1_{H})+B(x_{1}\otimes x_{2};G,1_{H})=0.
\end{equation*}%
In view of the form of the elements we get%
\begin{equation}
\lambda \left[ B(g\otimes 1_{H};G,1_{H})+B(x_{1}\otimes 1_{H};G,x_{1})\right]
+2B(gx_{1}x_{2}\otimes 1_{H};G,1_{H})=0  \label{X2,x1,X2F21,1H}
\end{equation}

\subparagraph{Case $f=x_{1}x_{2}$}

\begin{eqnarray*}
\beta _{2}B(x_{1}\otimes 1_{H};GX_{2},x_{1}x_{2}) &=&+\beta
_{2}B(gx_{1}\otimes g;GX_{2},x_{1}x_{2})+\lambda B(gx_{1}\otimes
g;GX_{1},x_{1}x_{2}) \\
&&+B(x_{1}x_{2}\otimes g;G,x_{1}x_{2})+B(x_{1}\otimes x_{2};G,x_{1}x_{2}).
\end{eqnarray*}%
By applying $\left( \ref{eq.10}\right) $ this rewrites as%
\begin{equation*}
\lambda B(x_{1}\otimes 1_{H};GX_{1},x_{1}x_{2})+B(gx_{1}x_{2}\otimes
1_{H};G,x_{1}x_{2})+B(x_{1}\otimes x_{2};G,x_{1}x_{2})=0.
\end{equation*}%
In view of the form of the elements we get%
\begin{equation*}
\lambda B\left( g\otimes 1_{H};G,x_{1}x_{2}\right) +2B(gx_{1}x_{2}\otimes
1_{H};G,x_{1}x_{2})=0.
\end{equation*}%
which is $\left( \ref{X1,x2,X1F21,x1x2}\right) .$

\subparagraph{Case $f=gx_{1}$}

\begin{eqnarray*}
\beta _{2}B(x_{1}\otimes 1_{H};GX_{2},gx_{1}) &=&+\beta _{2}B(gx_{1}\otimes
g;GX_{2},gx_{1})+\lambda B(gx_{1}\otimes g;GX_{1},gx_{1}) \\
&&+B(x_{1}x_{2}\otimes g;G,gx_{1})+B(x_{1}\otimes x_{2};G,gx_{1}).
\end{eqnarray*}%
By applying $\left( \ref{eq.10}\right) $ this rewrites as%
\begin{eqnarray*}
&&2\beta _{2}B(x_{1}\otimes 1_{H};GX_{2},gx_{1})+\lambda B(x_{1}\otimes
1_{H};GX_{1},gx_{1}) \\
&&+B(gx_{1}x_{2}\otimes 1_{H};G,gx_{1})-B(x_{1}\otimes x_{2};G,gx_{1})=0.
\end{eqnarray*}%
In view of the form of the elements we get%
\begin{equation}
2\beta _{2}B(x_{1}\otimes 1_{H};G,gx_{1}x_{2})+\lambda B\left( g\otimes
1_{H};G,gx_{1}\right) +2B(gx_{1}x_{2}\otimes 1_{H};G,gx_{1})=0.
\label{X2,x1,X2F21,gx1}
\end{equation}

\subparagraph{Case $f=gx_{2}$}

\begin{eqnarray*}
\beta _{2}B(x_{1}\otimes 1_{H};GX_{2},gx_{2}) &=&+\beta _{2}B(gx_{1}\otimes
g;GX_{2},gx_{2})+\lambda B(gx_{1}\otimes g;GX_{1},gx_{2}) \\
&&+B(x_{1}x_{2}\otimes g;G,gx_{2})+B(x_{1}\otimes x_{2};G,gx_{2}).
\end{eqnarray*}%
By applying $\left( \ref{eq.10}\right) $ this rewrites as%
\begin{eqnarray*}
&&2\beta _{2}B(x_{1}\otimes 1_{H};GX_{2},gx_{2})+\lambda B(x_{1}\otimes
1_{H};GX_{1},gx_{2}) \\
&&+B(gx_{1}x_{2}\otimes 1_{H};G,gx_{2})-B(x_{1}\otimes x_{2};G,gx_{2})=0.
\end{eqnarray*}%
In view of the form of the elements we get%
\begin{eqnarray}
&&+\lambda \left[ B\left( g\otimes 1_{H};G,gx_{2}\right) -B(x_{1}\otimes
1_{H};G,gx_{1}x_{2})\right]  \label{X2,x1,X2F21,gx2} \\
&&+2B(gx_{1}x_{2}\otimes 1_{H};G,gx_{2})+2B\left( x_{1}\otimes
1_{H};G,g\right) =0.  \notag
\end{eqnarray}

\paragraph{Equality $\left( \protect\ref{X2F31}\right) $}

rewrites as

\subparagraph{%
\protect\begin{eqnarray*}
\protect\beta _{2}B(x_{1}\otimes 1_{H};X_{1}X_{2},f) &=&+\protect\beta %
_{2}B(gx_{1}\otimes g;X_{1}X_{2},f)-\protect\gamma _{2}B(gx_{1}\otimes
g;GX_{1},f) \\
&&+B(x_{1}x_{2}\otimes g;X_{1},f)+B(x_{1}\otimes x_{2};X_{1},f).
\protect\end{eqnarray*}%
Case $f=1_{H}$}

\begin{eqnarray*}
\beta _{2}B(x_{1}\otimes 1_{H};X_{1}X_{2},1_{H}) &=&+\beta
_{2}B(gx_{1}\otimes g;X_{1}X_{2},1_{H})-\gamma _{2}B(gx_{1}\otimes
g;GX_{1},1_{H}) \\
&&+B(x_{1}x_{2}\otimes g;X_{1},1_{H})+B(x_{1}\otimes x_{2};X_{1},1_{H}).
\end{eqnarray*}%
By applying $\left( \ref{eq.10}\right) $ this rewrites as%
\begin{equation*}
\gamma _{2}B(x_{1}\otimes 1_{H};GX_{1},1_{H})-B(gx_{1}x_{2}\otimes
1_{H};X_{1},1_{H})-B(x_{1}\otimes x_{2};X_{1},1_{H})=0.
\end{equation*}%
In view of the form of the elements we get%
\begin{gather}
\gamma _{2}\left[ B(g\otimes 1_{H};G,1_{H})+B(x_{1}\otimes 1_{H};G,x_{1})%
\right]  \label{X2,x1,X2F31,1H} \\
+2B(x_{2}\otimes 1_{H};1_{A},1_{H})-2B(gx_{1}x_{2}\otimes
1_{H};1_{A},x_{1})=0.  \notag
\end{gather}

\subparagraph{Case $f=x_{1}x_{2}$}

\begin{eqnarray*}
\beta _{2}B(x_{1}\otimes 1_{H};X_{1}X_{2},x_{1}x_{2}) &=&+\beta
_{2}B(gx_{1}\otimes g;X_{1}X_{2},x_{1}x_{2})-\gamma _{2}B(gx_{1}\otimes
g;GX_{1},x_{1}x_{2}) \\
&&+B(x_{1}x_{2}\otimes g;X_{1},x_{1}x_{2})+B(x_{1}\otimes
x_{2};X_{1},x_{1}x_{2}).
\end{eqnarray*}%
By applying $\left( \ref{eq.10}\right) $ this rewrites as%
\begin{equation*}
\gamma _{2}B(x_{1}\otimes 1_{H};GX_{1},x_{1}x_{2})-B(gx_{1}x_{2}\otimes
1_{H};X_{1},x_{1}x_{2})-B(x_{1}\otimes x_{2};X_{1},x_{1}x_{2})=0.
\end{equation*}%
In view of the form of the elements we get%
\begin{equation*}
\gamma _{2}B\left( g\otimes 1_{H};G,x_{1}x_{2}\right) +2B(x_{2}\otimes
1_{H};1_{A},x_{1}x_{2})=0.
\end{equation*}%
which is $\left( \ref{G,x2, GF2,x1x2}\right) .$

\subparagraph{Case $f=gx_{1}$}

\begin{eqnarray*}
\beta _{2}B(x_{1}\otimes 1_{H};X_{1}X_{2},gx_{1}) &=&+\beta
_{2}B(gx_{1}\otimes g;X_{1}X_{2},gx_{1})-\gamma _{2}B(gx_{1}\otimes
g;GX_{1},gx_{1}) \\
&&+B(x_{1}x_{2}\otimes g;X_{1},gx_{1})+B(x_{1}\otimes x_{2};X_{1},gx_{1}).
\end{eqnarray*}%
By applying $\left( \ref{eq.10}\right) $ this rewrites as%
\begin{gather*}
2\beta _{2}B(x_{1}\otimes 1_{H};X_{1}X_{2},gx_{1})-\gamma _{2}B(x_{1}\otimes
1_{H};GX_{1},gx_{1}) \\
+B(gx_{1}x_{2}\otimes 1_{H};X_{1},gx_{1})-B(x_{1}\otimes
x_{2};X_{1},gx_{1})=0.
\end{gather*}%
In view of the form of the elements we get%
\begin{equation*}
2\beta _{2}B\left( g\otimes 1_{H};1_{A},gx_{1}x_{2}\right) +\gamma
_{2}B\left( g\otimes 1_{H};G,gx_{1}\right) +2B(x_{2}\otimes
1_{H};1_{A},gx_{1})=0
\end{equation*}%
which is $\left( \ref{X2,g,X2F11,gx1}\right) $

\subparagraph{Case $f=gx_{2}$}

\begin{eqnarray*}
\beta _{2}B(x_{1}\otimes 1_{H};X_{1}X_{2},gx_{2}) &=&+\beta
_{2}B(gx_{1}\otimes g;X_{1}X_{2},gx_{2})-\gamma _{2}B(gx_{1}\otimes
g;GX_{1},gx_{2}) \\
&&+B(x_{1}x_{2}\otimes g;X_{1},gx_{2})+B(x_{1}\otimes x_{2};X_{1},gx_{2}).
\end{eqnarray*}%
By applying $\left( \ref{eq.10}\right) $ this rewrites as%
\begin{eqnarray*}
&&2\beta _{2}B(x_{1}\otimes 1_{H};X_{1}X_{2},gx_{2})-\gamma
_{2}B(x_{1}\otimes 1_{H};GX_{1},gx_{2}) \\
&&+B(gx_{1}x_{2}\otimes 1_{H};X_{1},gx_{2})-B(x_{1}\otimes
x_{2};X_{1},gx_{2})=0.
\end{eqnarray*}%
In view of the form of the elements we get%
\begin{gather}
\gamma _{2}\left[ -B(g\otimes 1_{H};G,gx_{2})+B(x_{1}\otimes
1_{H};G,gx_{1}x_{2})\right] +  \label{X2,x1,X2F31,gx2} \\
-2B(x_{2}\otimes 1_{H};1_{A},gx_{2})-2B(gx_{1}x_{2}\otimes
1_{H};1_{A},gx_{1}x_{2})  \notag \\
-2B(g\otimes 1_{H};1_{A},g)-2B(x_{1}\otimes 1_{H};1_{A},gx_{1})=0.  \notag
\end{gather}

\paragraph{Equality $\left( \protect\ref{X2F41}\right) $}

rewrites as%
\begin{eqnarray*}
B(x_{1}\otimes 1_{H};1_{A},f) &=&-B(gx_{1}\otimes g;1_{A},f)-\gamma
_{2}B(gx_{1}\otimes g;GX_{2},f)-\lambda B(gx_{1}\otimes g;X_{1}X_{2},f) \\
&&+B(x_{1}x_{2}\otimes g;X_{2},f)+B(x_{1}\otimes x_{2};X_{2},f)
\end{eqnarray*}

\subparagraph{Case $f=1_{H}$}

\begin{eqnarray*}
B(x_{1}\otimes 1_{H};1_{A},1_{H}) &=&-B(gx_{1}\otimes g;1_{A},1_{H})-\gamma
_{2}B(gx_{1}\otimes g;GX_{2},1_{H})-\lambda B(gx_{1}\otimes
g;X_{1}X_{2},1_{H}) \\
&&+B(x_{1}x_{2}\otimes g;X_{2},1_{H})+B(x_{1}\otimes x_{2};X_{2},1_{H})
\end{eqnarray*}%
By applying $\left( \ref{eq.10}\right) $ this rewrites as%
\begin{eqnarray*}
&&2B(x_{1}\otimes 1_{H};1_{A},1_{H})+\gamma _{2}B(x_{1}\otimes
1_{H};GX_{2},1_{H})+\lambda B(x_{1}\otimes 1_{H};X_{1}X_{2},1_{H}) \\
&&-B(gx_{1}x_{2}\otimes 1_{H};X_{2},1_{H})-B(x_{1}\otimes
x_{2};X_{2},1_{H})=0
\end{eqnarray*}%
In view of the form of the elements we get%
\begin{gather}
+\gamma _{2}B(x_{1}\otimes 1_{H};G,x_{2})  \label{X2,x1,X2F41,1H} \\
+\lambda \left[ -B(g\otimes 1_{H};1_{A},x_{2})+B(x_{1}\otimes
1_{H};1_{A},x_{1}x_{2})\right]  \notag \\
-2B(gx_{1}x_{2}\otimes 1_{H};1_{A},x_{2})=0  \notag
\end{gather}

\subparagraph{Case $f=x_{1}x_{2}$}

\begin{eqnarray*}
&&B(x_{1}\otimes 1_{H};1_{A},x_{1}x_{2}) \\
&=&-B(gx_{1}\otimes g;1_{A},x_{1}x_{2})-\gamma _{2}B(gx_{1}\otimes
g;GX_{2},x_{1}x_{2})-\lambda B(gx_{1}\otimes g;X_{1}X_{2},x_{1}x_{2}) \\
&&+B(x_{1}x_{2}\otimes g;X_{2},x_{1}x_{2})+B(x_{1}\otimes
x_{2};X_{2},x_{1}x_{2})
\end{eqnarray*}%
By applying $\left( \ref{eq.10}\right) $ this rewrites as%
\begin{gather*}
2B(x_{1}\otimes 1_{H};1_{A},x_{1}x_{2})+\gamma _{2}B(x_{1}\otimes
1_{H};GX_{2},x_{1}x_{2})+\lambda B(x_{1}\otimes 1_{H};X_{1}X_{2},x_{1}x_{2})
\\
-B(gx_{1}x_{2}\otimes 1_{H};X_{2},x_{1}x_{2})-B(x_{1}\otimes
x_{2};X_{2},x_{1}x_{2})=0
\end{gather*}%
In view of the form of the elements we get%
\begin{equation*}
B(x_{1}\otimes 1_{H};1_{A},x_{1}x_{2})-B(x_{1}\otimes
1_{H};1_{A},x_{1}x_{2})=0
\end{equation*}%
which is trivial.

\subparagraph{Case $f=gx_{1}$}

\begin{eqnarray*}
&&B(x_{1}\otimes 1_{H};1_{A},gx_{1}) \\
&=&-B(gx_{1}\otimes g;1_{A},gx_{1})-\gamma _{2}B(gx_{1}\otimes
g;GX_{2},gx_{1})-\lambda B(gx_{1}\otimes g;X_{1}X_{2},gx_{1})+ \\
&&+B(x_{1}x_{2}\otimes g;X_{2},gx_{1})+B(x_{1}\otimes x_{2};X_{2},gx_{1})
\end{eqnarray*}%
By applying $\left( \ref{eq.10}\right) $ this rewrites as%
\begin{eqnarray*}
&&-\gamma _{2}B(x_{1}\otimes 1_{H};;GX_{2},gx_{1})-\lambda B(x_{1}\otimes
1_{H};;X_{1}X_{2},gx_{1})+ \\
&&+B(gx_{1}x_{2}\otimes 1_{H};X_{2},gx_{1})-B(x_{1}\otimes
x_{2};X_{2},gx_{1})=0.
\end{eqnarray*}%
In view of the form of the elements we get%
\begin{eqnarray}
&&-\gamma _{2}B(x_{1}\otimes 1_{H};G,gx_{1}x_{2})+\lambda B\left( g\otimes
1_{H};1_{A},gx_{1}x_{2}\right) +  \label{X2,x1,X2F41,gx1} \\
&&+2B(x_{1}\otimes 1_{H};1_{A},gx_{1})+2B(gx_{1}x_{2}\otimes
1_{H};1_{A},gx_{1}x_{2})=0.  \notag
\end{eqnarray}

\subparagraph{Case $f=gx_{2}$}

\begin{eqnarray*}
&&B(x_{1}\otimes 1_{H};1_{A},gx_{2}) \\
&=&-B(gx_{1}\otimes g;1_{A},gx_{2})-\gamma _{2}B(gx_{1}\otimes
g;GX_{2},gx_{2})-\lambda B(gx_{1}\otimes g;X_{1}X_{2},gx_{2})+ \\
&&+B(x_{1}x_{2}\otimes g;X_{2},gx_{2})+B(x_{1}\otimes x_{2};X_{2},gx_{2})
\end{eqnarray*}%
By applying $\left( \ref{eq.10}\right) $ this rewrites as%
\begin{eqnarray*}
&&-\gamma _{2}B(x_{1}\otimes 1_{H};GX_{2},gx_{2})-\lambda B(x_{1}\otimes
1_{H};X_{1}X_{2},gx_{2})+ \\
&&+B(gx_{1}x_{2}\otimes 1_{H};X_{2},gx_{2})-B(x_{1}\otimes
x_{2};X_{2},gx_{2})=0.
\end{eqnarray*}%
In view of the form of the elements we get%
\begin{equation*}
B(x_{1}\otimes 1_{H};1_{A},gx_{2})-B(x_{1}\otimes 1_{H};1_{A},gx_{2})=0
\end{equation*}%
which is trivial.

\paragraph{Equality $\left( \protect\ref{X2F51}\right) $}

rewrites as%
\begin{eqnarray*}
B(x_{1}\otimes 1_{H};X_{1},f) &=&B(fx_{1}\otimes f;X_{1},f)-\gamma
_{2}B(fx_{1}\otimes f;GX_{1}X_{2},f) \\
&&+B(x_{1}x_{2}\otimes f;X_{1}X_{2},f)+B(x_{1}\otimes x_{2};X_{1}X_{2},f)
\end{eqnarray*}

\subparagraph{Case $f=g$}

\begin{eqnarray*}
B(x_{1}\otimes 1_{H};X_{1},g) &=&B(gx_{1}\otimes g;X_{1},g)-\gamma
_{2}B(gx_{1}\otimes g;GX_{1}X_{2},g) \\
&&+B(x_{1}x_{2}\otimes g;X_{1}X_{2},g)+B(x_{1}\otimes x_{2};X_{1}X_{2},g)
\end{eqnarray*}%
By applying $\left( \ref{eq.10}\right) $ this rewrites as%
\begin{eqnarray*}
&&+\gamma _{2}B(x_{1}\otimes 1_{H};GX_{1}X_{2},g) \\
&&-B(gx_{1}x_{2}\otimes 1_{H};X_{1}X_{2},g)-B(x_{1}\otimes
x_{2};X_{1}X_{2},g)=0
\end{eqnarray*}%
In view of the form of the elements we get%
\begin{gather*}
+\gamma _{2}\left[ -B(g\otimes 1_{H};G,gx_{2})+B(x_{1}\otimes
1_{H};G,gx_{1}x_{2})\right] \\
-2B(g\otimes 1_{H};1_{A},g)-2B(x_{2}\otimes \ 1_{H};1_{A},gx_{2}) \\
-2B(x_{1}\otimes 1_{H};1_{A},gx_{1})-2B(gx_{1}x_{2}\otimes
1_{H};1_{A},gx_{1}x_{2})=0
\end{gather*}%
This is $\left( \ref{X2,x1,X2F31,gx2}\right) .$

\subparagraph{Case $f=x_{1}$}

\begin{eqnarray*}
B(x_{1}\otimes 1_{H};X_{1},x_{1}) &=&+B(gx_{1}\otimes g;X_{1},x_{1})-\gamma
_{2}B(gx_{1}\otimes g;GX_{1}X_{2},x_{1}) \\
&&+B(x_{1}x_{2}\otimes g;X_{1}X_{2},x_{1})+B(x_{1}\otimes
x_{2};X_{1}X_{2},x_{1})
\end{eqnarray*}%
By applying $\left( \ref{eq.10}\right) $ this rewrites as%
\begin{eqnarray*}
&&2B(x_{1}\otimes 1_{H};X_{1},x_{1})-\gamma _{2}B(x_{1}\otimes
1_{H};GX_{1}X_{2},x_{1}) \\
&&+B(gx_{1}x_{2}\otimes 1_{H};X_{1}X_{2},x_{1})-B(x_{1}\otimes
x_{2};X_{1}X_{2},x_{1})=0.
\end{eqnarray*}%
In view of the form of the elements we get%
\begin{equation*}
2B(x_{2}\otimes \ 1_{H};1_{A},x_{1}x_{2})+\gamma _{2}B(g\otimes
1_{H};G,x_{1}x_{2})=0
\end{equation*}%
which is $\left( \ref{G,x2; GF6,x2}\right) .$

\subparagraph{Case $f=x_{2}$}

\begin{eqnarray*}
B(x_{1}\otimes 1_{H};X_{1},x_{2}) &=&+B(gx_{1}\otimes g;X_{1},x_{2})-\gamma
_{2}B(gx_{1}\otimes g;GX_{1}X_{2},x_{2}) \\
&&+B(x_{1}x_{2}\otimes g;X_{1}X_{2},x_{2})+B(x_{1}\otimes
x_{2};X_{1}X_{2},x_{2})
\end{eqnarray*}%
By applying $\left( \ref{eq.10}\right) $ this rewrites as%
\begin{eqnarray*}
&&2B(x_{1}\otimes 1_{H};X_{1},x_{2})-\gamma _{2}B(x_{1}\otimes
1_{H};GX_{1}X_{2},x_{2}) \\
&&+B(gx_{1}x_{2}\otimes 1_{H};X_{1}X_{2},x_{2})-B(x_{1}\otimes
x_{2};X_{1}X_{2},x_{2})=0.
\end{eqnarray*}%
In view of the form of the elements we get%
\begin{eqnarray*}
&&-B(g\otimes 1_{H};1_{A},x_{2})+B(x_{1}\otimes 1_{H};1_{A},x_{1}x_{2}) \\
&&+B(g\otimes 1_{H};1_{A},x_{2})-B(x_{1}\otimes 1_{H};1_{A},x_{1}x_{2})=0.
\end{eqnarray*}%
which is trivial.

\subparagraph{Case $f=gx_{1}x_{2}$}

\begin{eqnarray*}
B(x_{1}\otimes 1_{H};X_{1},gx_{1}x_{2}) &=&+B(gx_{1}\otimes
g;X_{1},gx_{1}x_{2})-\gamma _{2}B(gx_{1}\otimes g;GX_{1}X_{2},gx_{1}x_{2}) \\
&&+B(x_{1}x_{2}\otimes g;X_{1}X_{2},gx_{1}x_{2})+B(x_{1}\otimes
x_{2};X_{1}X_{2},gx_{1}x_{2})
\end{eqnarray*}%
By applying $\left( \ref{eq.10}\right) $ this rewrites as%
\begin{eqnarray*}
&&\gamma _{2}B(x_{1}\otimes 1_{H};GX_{1}X_{2},gx_{1}x_{2}) \\
&&-B(gx_{1}x_{2}\otimes 1_{H};X_{1}X_{2},gx_{1}x_{2})-B(x_{1}\otimes
x_{2};X_{1}X_{2},gx_{1}x_{2})=0
\end{eqnarray*}%
In view of the form of the elements we get%
\begin{equation*}
-B\left( g\otimes 1_{H};1_{A},gx_{1}x_{2}\right) +2B\left( g\otimes
1_{H};1_{A},gx_{1}x_{2}\right) -B\left( g\otimes
1_{H};1_{A},gx_{1}x_{2}\right) =0
\end{equation*}%
which is trivial.

\paragraph{Equality $\left( \protect\ref{X2F61}\right) $}

rewrites as%
\begin{eqnarray*}
\beta _{2}B(x_{1}\otimes 1_{H};GX_{1}X_{2},f) &=&-\beta _{2}B(gx_{1}\otimes
g;GX_{1}X_{2},f) \\
&&+B(x_{1}x_{2}\otimes g;GX_{1},f)+B(x_{1}\otimes x_{2};GX_{1},f)
\end{eqnarray*}

\subparagraph{Case $f=g$}

\begin{eqnarray*}
\beta _{2}B(x_{1}\otimes 1_{H};GX_{1}X_{2},g) &=&-\beta _{2}B(gx_{1}\otimes
g;GX_{1}X_{2},g) \\
&&+B(x_{1}x_{2}\otimes g;GX_{1},g)+B(x_{1}\otimes x_{2};GX_{1},g)
\end{eqnarray*}%
By applying $\left( \ref{eq.10}\right) $ this rewrites as

\begin{eqnarray*}
&&2\beta _{2}B(x_{1}\otimes 1_{H};GX_{1}X_{2},g) \\
&&-B(gx_{1}x_{2}\otimes 1_{H};GX_{1},g)-B(x_{1}\otimes x_{2};GX_{1},g)=0.
\end{eqnarray*}%
In view of the form of the elements we get%
\begin{eqnarray}
&&\beta _{2}\left[ -B(g\otimes 1_{H};G,gx_{2})+B(x_{1}\otimes
1_{H};G,gx_{1}x_{2})\right]  \label{X2,x1,X2F61,g} \\
&&-B(x_{2}\otimes 1_{H};G,g)+B(gx_{1}x_{2}\otimes 1_{H};G,gx_{1})=0  \notag
\end{eqnarray}

\subparagraph{Case $f=x_{1}$}

\begin{eqnarray*}
\beta _{2}B(x_{1}\otimes 1_{H};GX_{1}X_{2},x_{1}) &=&-\beta
_{2}B(gx_{1}\otimes g;GX_{1}X_{2},x_{1}) \\
&&+B(x_{1}x_{2}\otimes g;GX_{1},x_{1})+B(x_{1}\otimes x_{2};GX_{1},x_{1})
\end{eqnarray*}%
By applying $\left( \ref{eq.10}\right) $ this rewrites as%
\begin{equation*}
-B(gx_{1}x_{2}\otimes 1_{H};GX_{1},x_{1})+B(x_{1}\otimes
x_{2};GX_{1},x_{1})=0
\end{equation*}%
In view of the form of the elements we get%
\begin{equation*}
B(x_{2}\otimes 1_{H};G,x_{1})=0
\end{equation*}%
which is $\left( \ref{X2,g,X2F21,x1}\right) .$

\subparagraph{Case $f=x_{2}$}

\begin{eqnarray*}
\beta _{2}B(x_{1}\otimes 1_{H};GX_{1}X_{2},x_{2}) &=&-\beta
_{2}B(gx_{1}\otimes g;GX_{1}X_{2},x_{2}) \\
&&+B(x_{1}x_{2}\otimes g;GX_{1},x_{2})+B(x_{1}\otimes x_{2};GX_{1},x_{2})
\end{eqnarray*}%
By applying $\left( \ref{eq.10}\right) $ this rewrites as%
\begin{equation*}
-B(gx_{1}x_{2}\otimes 1_{H};GX_{1},x_{2})+B(x_{1}\otimes x_{2};GX_{1},x_{2}).
\end{equation*}%
In view of the form of the elements we get%
\begin{equation}
B(x_{2}\otimes 1_{H};G,x_{2})+B(gx_{1}x_{2}\otimes 1_{H};G,x_{1}x_{2})=0
\label{X2,x1,X2F61,x2}
\end{equation}

\subparagraph{Case $f=gx_{1}x_{2}$}

\begin{eqnarray*}
\beta _{2}B(x_{1}\otimes 1_{H};GX_{1}X_{2},gx_{1}x_{2}) &=&-\beta
_{2}B(gx_{1}\otimes g;GX_{1}X_{2},gx_{1}x_{2}) \\
&&+B(x_{1}x_{2}\otimes g;GX_{1},gx_{1}x_{2})+B(x_{1}\otimes
x_{2};GX_{1},gx_{1}x_{2})
\end{eqnarray*}%
By applying $\left( \ref{eq.10}\right) $ this rewrites as%
\begin{eqnarray*}
&&2\beta _{2}B(x_{1}\otimes 1_{H};GX_{1}X_{2},gx_{1}x_{2}) \\
&&-B(gx_{1}x_{2}\otimes 1_{H};GX_{1},gx_{1}x_{2})-B(x_{1}\otimes
x_{2};GX_{1},gx_{1}x_{2})=0.
\end{eqnarray*}%
In view of the form of the elements we get%
\begin{equation*}
B\left( g\otimes 1_{H};G,gx_{1}\right) +B(x_{2}\otimes 1_{H};G,gx_{1}x_{2})=0
\end{equation*}%
which is $\left( \ref{X2,g,X2F71,gx1}\right) .$

\paragraph{\textbf{Equality} $\left( \protect\ref{X2F71}\right) $}

rewrites as%
\begin{eqnarray*}
&&B(x_{1}\otimes 1_{H};G,f) \\
&=&B(gx_{1}\otimes g;G,f)+\lambda B(gx_{1}\otimes g;GX_{1}X_{2},f) \\
&&+B(x_{1}x_{2}\otimes g;GX_{2},f)+B(x_{1}\otimes x_{2};GX_{2},f).
\end{eqnarray*}

\subparagraph{Case $f=g$}

\begin{eqnarray*}
B(x_{1}\otimes 1_{H};G,g) &=&B(gx_{1}\otimes g;G,g)+\lambda B(gx_{1}\otimes
g;GX_{1}X_{2},g) \\
&&+B(x_{1}x_{2}\otimes g;GX_{2},g)+B(x_{1}\otimes x_{2};GX_{2},g).
\end{eqnarray*}%
By applying $\left( \ref{eq.10}\right) $ this rewrites as%
\begin{eqnarray*}
&&-\lambda B(x_{1}\otimes 1_{H};GX_{1}X_{2},g) \\
&&-B(gx_{1}x_{2}\otimes 1_{H};GX_{2},g)-B(x_{1}\otimes x_{2};GX_{2},g).
\end{eqnarray*}%
In view of the form of the elements we get%
\begin{eqnarray*}
&&\lambda \left[ B(g\otimes 1_{H};G,gx_{2})-B(x_{1}\otimes
1_{H};G,gx_{1}x_{2})\right] \\
&&+2B(x_{1}\otimes 1_{H};G,g)+2B(gx_{1}x_{2}\otimes 1_{H};G,gx_{2})=0.
\end{eqnarray*}%
By $\left( \ref{X1,g,X1F21,gx1x2}\right) ,$we obtain%
\begin{equation}
B(x_{1}\otimes 1_{H};G,g)+B(gx_{1}x_{2}\otimes 1_{H};G,gx_{2})=0
\label{X2,x1,X2F71,g}
\end{equation}

\subparagraph{Case $f=x_{1}$}

\begin{eqnarray*}
B(x_{1}\otimes 1_{H};G,x_{1}) &=&+B(gx_{1}\otimes g;G,x_{1})+\lambda
B(gx_{1}\otimes g;GX_{1}X_{2},x_{1}) \\
&&+B(x_{1}x_{2}\otimes g;GX_{2},x_{1})+B(x_{1}\otimes x_{2};GX_{2},x_{1}).
\end{eqnarray*}%
By applying $\left( \ref{eq.10}\right) $ this rewrites as%
\begin{eqnarray*}
&&2B(x_{1}\otimes 1_{H};G,x_{1})+\lambda B(x_{1}\otimes
1_{H};GX_{1}X_{2},x_{1}) \\
&&+B(gx_{1}x_{2}\otimes 1_{H};GX_{2},x_{1})-B(x_{1}\otimes
x_{2};GX_{2},x_{1})=0.
\end{eqnarray*}%
In view of the form of the elements we get%
\begin{equation*}
\lambda B(g\otimes 1_{H};G,x_{1}x_{2})+2B(gx_{1}x_{2}\otimes
1_{H};G,x_{1}x_{2})=0
\end{equation*}%
which is $\left( \ref{X1,gx1x2,X1F61,x1x2}\right) .$

\subparagraph{Case $f=x_{2}$}

\begin{eqnarray*}
B(x_{1}\otimes 1_{H};G,x_{2}) &=&+B(gx_{1}\otimes g;G,x_{2})+\lambda
B(gx_{1}\otimes g;GX_{1}X_{2},x_{2}) \\
&&+B(x_{1}x_{2}\otimes g;GX_{2},x_{2})+B(x_{1}\otimes x_{2};GX_{2},x_{2}).
\end{eqnarray*}%
By applying $\left( \ref{eq.10}\right) $ this rewrites as%
\begin{eqnarray*}
&&2B(x_{1}\otimes 1_{H};G,x_{2})+\lambda B(x_{1}\otimes
1_{H};GX_{1}X_{2},x_{2}) \\
&&+B(gx_{1}x_{2}\otimes 1_{H};GX_{2},x_{2})-B(x_{1}\otimes
x_{2};GX_{2},x_{2})=0.
\end{eqnarray*}%
In view of the form of the elements we get%
\begin{equation*}
B(x_{1}\otimes 1_{H};G,x_{2})-B(x_{1}\otimes 1_{H};G,x_{2})=0
\end{equation*}%
which is trivial.

\subparagraph{Case $f=gx_{1}x_{2}$}

\begin{eqnarray*}
B(x_{1}\otimes 1_{H};G,gx_{1}x_{2}) &=&+B(gx_{1}\otimes
g;G,gx_{1}x_{2})+\lambda B(gx_{1}\otimes g;GX_{1}X_{2},gx_{1}x_{2}) \\
&&+B(x_{1}x_{2}\otimes g;GX_{2},gx_{1}x_{2})+B(x_{1}\otimes
x_{2};GX_{2},gx_{1}x_{2}).
\end{eqnarray*}%
By applying $\left( \ref{eq.10}\right) $ this rewrites as%
\begin{gather*}
-\lambda B(x_{1}\otimes 1_{H};GX_{1}X_{2},gx_{1}x_{2}) \\
-B(gx_{1}x_{2}\otimes 1_{H};GX_{2},gx_{1}x_{2})-B(x_{1}\otimes
x_{2};GX_{2},gx_{1}x_{2})=0
\end{gather*}%
In view of the form of the elements we get%
\begin{equation*}
B(x_{1}\otimes 1_{H};G,gx_{1}x_{2})-B(x_{1}\otimes 1_{H};G,gx_{1}x_{2})=0
\end{equation*}%
which is trivial.

\paragraph{\textbf{Equality} $\left( \protect\ref{X2F81}\right) $}

rewrites as%
\begin{eqnarray*}
B(x_{1}\otimes 1_{H};GX_{1},f) &=&-B(gx_{1}\otimes g;GX_{1},f) \\
&&+B(x_{1}x_{2}\otimes g;GX_{1}X_{2},f)+B(x_{1}\otimes x_{2};GX_{1}X_{2},f).
\end{eqnarray*}

\subparagraph{Case $f=1_{H}$}

\begin{eqnarray*}
B(x_{1}\otimes 1_{H};GX_{1},1_{H}) &=&-B(gx_{1}\otimes g;GX_{1},1_{H}) \\
&&+B(x_{1}x_{2}\otimes g;GX_{1}X_{2},1_{H})+B(x_{1}\otimes
x_{2};GX_{1}X_{2},1_{H}).
\end{eqnarray*}%
By applying $\left( \ref{eq.10}\right) $ this rewrites as%
\begin{eqnarray*}
&&2B(x_{1}\otimes 1_{H};GX_{1},1_{H}) \\
&&-B(gx_{1}x_{2}\otimes 1_{H};GX_{1}X_{2},1_{H})-B(x_{1}\otimes
x_{2};GX_{1}X_{2},1_{H})=0.
\end{eqnarray*}

In view of the form of the elements we get%
\begin{equation}
B(x_{2}\otimes \ 1_{H};G,x_{2})+B(gx_{1}x_{2}\otimes 1_{H};G,x_{1}x_{2})=0
\label{X2,x1,X2F81,1H}
\end{equation}

\subparagraph{Case $f=x_{1}x_{2}$}

\begin{eqnarray*}
B(x_{1}\otimes 1_{H};GX_{1},x_{1}x_{2}) &=&-B(gx_{1}\otimes
g;GX_{1},x_{1}x_{2}) \\
&&+B(x_{1}x_{2}\otimes g;GX_{1}X_{2},x_{1}x_{2})+B(x_{1}\otimes
x_{2};GX_{1}X_{2},x_{1}x_{2}).
\end{eqnarray*}%
By applying $\left( \ref{eq.10}\right) $ this rewrites as%
\begin{eqnarray*}
&&2B(x_{1}\otimes 1_{H};GX_{1},x_{1}x_{2}) \\
&&-B(gx_{1}x_{2}\otimes 1_{H};GX_{1}X_{2},x_{1}x_{2})-B(x_{1}\otimes
x_{2};GX_{1}X_{2},x_{1}x_{2})=0.
\end{eqnarray*}%
In view of the form of the elements we get%
\begin{equation*}
B\left( g\otimes 1_{H};G,x_{1}x_{2}\right) -B(g\otimes 1_{H};G,x_{1}x_{2})=0
\end{equation*}%
which is trivial

\subparagraph{Case $f=gx_{1}$}

\begin{eqnarray*}
B(x_{1}\otimes 1_{H};GX_{1},gx_{1}) &=&-B(gx_{1}\otimes g;GX_{1},gx_{1}) \\
&&+B(x_{1}x_{2}\otimes g;GX_{1}X_{2},gx_{1})+B(x_{1}\otimes
x_{2};GX_{1}X_{2},gx_{1}).
\end{eqnarray*}%
By applying $\left( \ref{eq.10}\right) $ this rewrites as%
\begin{equation*}
B(gx_{1}x_{2}\otimes 1_{H};GX_{1}X_{2},gx_{1})-B(x_{1}\otimes
x_{2};GX_{1}X_{2},gx_{1})=0.
\end{equation*}%
In view of the form of the elements we get%
\begin{equation*}
B(g\otimes 1_{H};G,gx_{1})+B(x_{2}\otimes \
1_{H};G,gx_{1}x_{2})-B(x_{1}\otimes x_{2};GX_{1}X_{2},gx_{1})=0.
\end{equation*}

\subparagraph{Case $f=gx_{2}$}

\begin{eqnarray*}
B(x_{1}\otimes 1_{H};GX_{1},gx_{2}) &=&-B(gx_{1}\otimes g;GX_{1},gx_{2}) \\
&&+B(x_{1}x_{2}\otimes g;GX_{1}X_{2},gx_{2})+B(x_{1}\otimes
x_{2};GX_{1}X_{2},gx_{2}).
\end{eqnarray*}

By applying $\left( \ref{eq.10}\right) $ this rewrites as%
\begin{equation*}
B(gx_{1}x_{2}\otimes 1_{H};GX_{1}X_{2},gx_{2})-B(x_{1}\otimes
x_{2};GX_{1}X_{2},gx_{2})=0.
\end{equation*}%
In view of the form of the elements we get%
\begin{gather*}
B(g\otimes 1_{H};G,gx_{2})-B(x_{1}\otimes 1_{H};G,gx_{1}x_{2}) \\
-B(g\otimes 1_{H};G,gx_{2})+B(x_{1}\otimes 1_{H};G,gx_{1}x_{2})=0
\end{gather*}%
which is trivial.

\subsubsection{Case $x_{2}\otimes 1_{H}$}

\paragraph{\textbf{Equality} $\left( \protect\ref{X2F11}\right) $}

rewrites as

\begin{eqnarray*}
\beta _{2}B(x_{2}\otimes 1_{H};X_{2},f) &=&-\gamma _{2}B(gx_{2}\otimes
g;G,f)-\lambda B(gx_{2}\otimes g;X_{1},f) \\
&&-\beta _{2}B(gx_{2}\otimes g;X_{2},f)+B(x_{2}\otimes x_{2};1_{A},f)
\end{eqnarray*}

\subparagraph{Case $f=g$}

\begin{eqnarray*}
\beta _{2}B(x_{2}\otimes 1_{H};X_{2},g) &=&-\gamma _{2}B(gx_{2}\otimes
g;G,g)-\lambda B(gx_{2}\otimes g;X_{1},g) \\
&&-\beta _{2}B(gx_{2}\otimes g;X_{2},g)+B(x_{2}\otimes x_{2};1_{A},g)
\end{eqnarray*}%
By applying $\left( \ref{eq.10}\right) $ this rewrites as%
\begin{eqnarray*}
2\beta _{2}B(x_{2}\otimes 1_{H};X_{2},g)+\gamma _{2}B(x_{2}\otimes 1_{H};G,g)
\\
\lambda B(x_{2}\otimes 1_{H};X_{1},g)-B(x_{2}\otimes x_{2};1_{A},g)=0.
\end{eqnarray*}%
In view of the form of the elements we get%
\begin{eqnarray}
2\beta _{2}\left[ B(g\otimes 1_{H};1_{A},g)+B(x_{2}\otimes \
1_{H};1_{A},gx_{2})\right]  \label{X2,x2,X2F11,g} \\
\gamma _{2}B(x_{2}\otimes 1_{H};G,g)+\lambda B(x_{2}\otimes
1_{H};1_{A},gx_{1})=0.  \notag
\end{eqnarray}

\subparagraph{Case $f=x_{1}$}

\begin{eqnarray*}
\beta _{2}B(x_{2}\otimes 1_{H};X_{2},x_{1}) &=&-\gamma _{2}B(gx_{2}\otimes
g;G,x_{1})-\lambda B(gx_{2}\otimes g;X_{1},x_{1}) \\
&&-\beta _{2}B(gx_{2}\otimes g;X_{2},x_{1})+B(x_{2}\otimes x_{2};1_{A},x_{1})
\end{eqnarray*}%
By applying $\left( \ref{eq.10}\right) $ this rewrites as%
\begin{eqnarray*}
&&-\gamma _{2}B(x_{2}\otimes 1_{H};G,x_{1})-\lambda B(x_{2}\otimes
1_{H};X_{1},x_{1}) \\
&&-B(x_{2}\otimes x_{2};1_{A},x_{1})=0
\end{eqnarray*}%
In view of the form of the elements we get%
\begin{equation*}
\gamma _{2}B(x_{2}\otimes 1_{H};G,x_{1})=0
\end{equation*}%
which follows from $\left( \ref{X2,g,X2F21,x1}\right) .$

\subparagraph{Case $f=x_{2}$}

\begin{eqnarray*}
\beta _{2}B(v;X_{2},x_{2}) &=&-\gamma _{2}B(gx_{2}\otimes g;G,x_{2})-\lambda
B(gx_{2}\otimes g;X_{1},x_{2}) \\
&&-\beta _{2}B(gx_{2}\otimes g;X_{2},x_{2})+B(x_{2}\otimes x_{2};1_{A},x_{2})
\end{eqnarray*}%
By applying $\left( \ref{eq.10}\right) $ this rewrites as%
\begin{eqnarray*}
&&-\gamma _{2}B(x_{2}\otimes 1_{H};G,x_{2})-\lambda B(x_{2}\otimes
1_{H};X_{1},x_{2}) \\
&&-B(x_{2}\otimes x_{2};1_{A},x_{2})=0.
\end{eqnarray*}%
In view of the form of the elements we get%
\begin{equation}
\gamma _{2}B(x_{2}\otimes 1_{H};G,x_{2})+\lambda B(x_{2}\otimes
1_{H};1_{A},x_{1}x_{2})=0.  \label{X2,x2,X2F11,x2}
\end{equation}

\subparagraph{Case $f=gx_{1}x_{2}$}

\begin{eqnarray*}
\beta _{2}B(x_{2}\otimes 1_{H};X_{2},gx_{1}x_{2}) &=&-\gamma
_{2}B(gx_{2}\otimes g;G,gx_{1}x_{2})-\lambda B(gx_{2}\otimes
g;X_{1},gx_{1}x_{2}) \\
&&-\beta _{2}B(gx_{2}\otimes g;X_{2},gx_{1}x_{2})+B(x_{2}\otimes
x_{2};1_{A},gx_{1}x_{2})
\end{eqnarray*}%
By applying $\left( \ref{eq.10}\right) $ this rewrites as%
\begin{eqnarray*}
&&2\beta _{2}B(x_{2}\otimes 1_{H};X_{2},gx_{1}x_{2})+\gamma
_{2}B(x_{2}\otimes 1_{H};G,gx_{1}x_{2}) \\
&&+\lambda B(x_{2}\otimes 1_{H};X_{1},gx_{1}x_{2})-B(x_{2}\otimes
x_{2};1_{A},gx_{1}x_{2})=0
\end{eqnarray*}%
In view of the form of the elements we get%
\begin{eqnarray}
+ &&2\beta _{2}B\left( g\otimes 1_{H};1_{A},gx_{1}x_{2}\right) -\gamma
_{2}B(x_{2}\otimes 1_{H};G,gx_{1}x_{2})  \label{X2,x2,X2F11,gx1x2} \\
&&+2B\left( x_{2}\otimes 1_{H};1_{A},gx_{1}\right) =0  \notag
\end{eqnarray}

\paragraph{\textbf{Equality} $\left( \protect\ref{X2F21}\right) $}

rewrites as%
\begin{eqnarray*}
\beta _{2}B(x_{2}\otimes 1_{H};GX_{2},f) &=&+\beta _{2}B(gx_{2}\otimes
g;GX_{2},f) \\
&&+\lambda B(gx_{2}\otimes g;GX_{1},f)+B(x_{2}\otimes x_{2};G,f)
\end{eqnarray*}

\subparagraph{Case $f=1_{H}$}

\begin{eqnarray*}
\beta _{2}B(x_{2}\otimes 1_{H};GX_{2},1_{H}) &=&+\beta _{2}B(gx_{2}\otimes
g;GX_{2},1_{H}) \\
&&+\lambda B(gx_{2}\otimes g;GX_{1},1_{H})+B(x_{2}\otimes x_{2};G,1_{H})
\end{eqnarray*}

By applying $\left( \ref{eq.10}\right) $ this rewrites as%
\begin{equation*}
+\lambda B(x_{2}\otimes 1_{H};GX_{1},1_{H})+B(x_{2}\otimes x_{2};G,1_{H})=0
\end{equation*}%
In view of the form of the elements we get%
\begin{equation*}
\lambda B(x_{2}\otimes 1_{H};G,x_{1})=0
\end{equation*}%
which holds in view of $\left( \ref{X2,g,X2F21,x1}\right) .$

\subparagraph{Case $f=x_{1}x_{2}$}

\begin{eqnarray*}
\beta _{2}B(x_{2}\otimes 1_{H};GX_{2},x_{1}x_{2}) &=&+\beta
_{2}B(gx_{2}\otimes g;GX_{2},x_{1}x_{2}) \\
&&+\lambda B(gx_{2}\otimes g;GX_{1},x_{1}x_{2})+B(x_{2}\otimes
x_{2};G,x_{1}x_{2})
\end{eqnarray*}%
By applying $\left( \ref{eq.10}\right) $ this rewrites as%
\begin{equation*}
-\lambda B(x_{2}\otimes 1_{H};GX_{1},x_{1}x_{2})-B(x_{2}\otimes
x_{2};G,x_{1}x_{2})=0
\end{equation*}%
In view of the form of the elements we get%
\begin{equation*}
B(x_{2}\otimes x_{2};G,x_{1}x_{2})=0
\end{equation*}%
which is already known.

\subparagraph{Case $f=gx_{1}$}

\begin{eqnarray*}
\beta _{2}B(x_{2}\otimes 1_{H};GX_{2},gx_{1}) &=&+\beta _{2}B(gx_{2}\otimes
g;GX_{2},gx_{1}) \\
&&+\lambda B(gx_{2}\otimes g;GX_{1},gx_{1})+B(x_{2}\otimes x_{2};G,gx_{1})
\end{eqnarray*}

By applying $\left( \ref{eq.10}\right) $ this rewrites as%
\begin{eqnarray*}
&&2\beta _{2}B(x_{2}\otimes 1_{H};GX_{2},gx_{1})+\lambda B(x_{2}\otimes
1_{H};GX_{1},gx_{1}) \\
&&-B(x_{2}\otimes x_{2};G,gx_{1})=0.
\end{eqnarray*}%
In view of the form of the elements we get%
\begin{equation}
\beta _{2}\left[ B(g\otimes 1_{H};G,gx_{1})+B(x_{2}\otimes \
1_{H};G,gx_{1}x_{2})\right] =0  \label{X2,x2,X2F21,gx1}
\end{equation}

\subparagraph{Case $f=gx_{2}$}

\begin{eqnarray*}
\beta _{2}B(x_{2}\otimes 1_{H};GX_{2},gx_{2}) &=&+\beta _{2}B(gx_{2}\otimes
g;GX_{2},gx_{2}) \\
&&+\lambda B(gx_{2}\otimes g;GX_{1},gx_{2})+B(x_{2}\otimes x_{2};G,gx_{2})
\end{eqnarray*}%
By applying $\left( \ref{eq.10}\right) $ this rewrites as%
\begin{eqnarray*}
&&2\beta _{2}B(x_{2}\otimes 1_{H};GX_{2},gx_{2})+\lambda B(x_{2}\otimes
1_{H};GX_{1},gx_{2}) \\
&&-B(x_{2}\otimes x_{2};G,gx_{2})=0
\end{eqnarray*}%
In view of the form of the elements we get%
\begin{eqnarray}
&&2\beta _{2}B(g\otimes 1_{H};G,gx_{2})-\lambda B(x_{2}\otimes
1_{H};G,gx_{1}x_{2})  \label{X2,x2,X2F21,gx2} \\
&&+2B\left( x_{2}\otimes 1_{H};G,g\right) =0  \notag
\end{eqnarray}

\paragraph{\textbf{Equality} $\left( \protect\ref{X2F31}\right) $}

rewrites as%
\begin{eqnarray*}
\beta _{2}B(x_{2}\otimes 1_{H};X_{1}X_{2},f) &=&+\beta _{2}B(gx_{2}\otimes
g;X_{1}X_{2},f) \\
&&-\gamma _{2}B(gx_{2}\otimes g;GX_{1},f)+B(x_{2}\otimes x_{2};X_{1},f).
\end{eqnarray*}

\subparagraph{Case $f=1_{H}$}

\begin{eqnarray*}
\beta _{2}B(x_{2}\otimes 1_{H};X_{1}X_{2},1_{H}) &=&+\beta
_{2}B(gx_{2}\otimes g;X_{1}X_{2},1_{H}) \\
&&-\gamma _{2}B(gx_{2}\otimes g;GX_{1},1_{H})+B(x_{2}\otimes
x_{2};X_{1},1_{H}).
\end{eqnarray*}%
By applying $\left( \ref{eq.10}\right) $ this rewrites as%
\begin{equation*}
\gamma _{2}B(x_{2}\otimes 1_{H};GX_{1},1_{H})-B(x_{2}\otimes
x_{2};X_{1},1_{H})=0.
\end{equation*}%
In view of the form of the elements we get%
\begin{equation*}
\gamma _{2}B(x_{2}\otimes 1_{H};G,x_{1})=0
\end{equation*}%
which follows from $\left( \ref{X2,g,X2F21,x1}\right) .$

\subparagraph{Case $f=x_{1}x_{2}$}

\begin{eqnarray*}
\beta _{2}B(x_{2}\otimes 1_{H};X_{1}X_{2},x_{1}x_{2}) &=&+\beta
_{2}B(gx_{2}\otimes g;X_{1}X_{2},x_{1}x_{2}) \\
&&-\gamma _{2}B(gx_{2}\otimes g;GX_{1},x_{1}x_{2})+B(x_{2}\otimes
x_{2};X_{1},x_{1}x_{2}).
\end{eqnarray*}%
By applying $\left( \ref{eq.10}\right) $ this rewrites as%
\begin{equation*}
\gamma _{2}B(x_{2}\otimes 1_{H};GX_{1},x_{1}x_{2})-B(x_{2}\otimes
x_{2};X_{1},x_{1}x_{2})=0.
\end{equation*}%
In view of the form of the elements we get%
\begin{equation*}
B(x_{2}\otimes x_{2};X_{1},x_{1}x_{2})=0
\end{equation*}%
which is already known.

\subparagraph{Case $f=gx_{1}$}

\begin{eqnarray*}
\beta _{2}B(x_{2}\otimes 1_{H};X_{1}X_{2},gx_{1}) &=&+\beta
_{2}B(gx_{2}\otimes g;X_{1}X_{2},gx_{1}) \\
&&-\gamma _{2}B(gx_{2}\otimes g;GX_{1},gx_{1})+B(x_{2}\otimes
x_{2};X_{1},gx_{1}).
\end{eqnarray*}%
By applying $\left( \ref{eq.10}\right) $ this rewrites as%
\begin{gather*}
2\beta _{2}B(x_{2}\otimes 1_{H};X_{1}X_{2},gx_{1})-\gamma _{2}B(x_{2}\otimes
1_{H};GX_{1},gx_{1}) \\
-B(x_{2}\otimes x_{2};X_{1},gx_{1})=0.
\end{gather*}%
In view of the form of the elements we get%
\begin{gather*}
2\beta _{2}B(x_{2}\otimes 1_{H};X_{1}X_{2},gx_{1})-\gamma _{2}B(x_{2}\otimes
1_{H};GX_{1},gx_{1}) \\
-B(x_{2}\otimes x_{2};X_{1},gx_{1})=0.
\end{gather*}%
which is already known.

\subparagraph{Case $f=gx_{2}$}

\begin{eqnarray*}
\beta _{2}B(x_{2}\otimes 1_{H};X_{1}X_{2},gx_{2}) &=&+\beta
_{2}B(gx_{2}\otimes g;X_{1}X_{2},gx_{2}) \\
&&-\gamma _{2}B(gx_{2}\otimes g;GX_{1},gx_{2})+B(x_{2}\otimes
x_{2};X_{1},gx_{2}).
\end{eqnarray*}%
By applying $\left( \ref{eq.10}\right) $ this rewrites as%
\begin{eqnarray*}
&&2\beta _{2}B(x_{2}\otimes 1_{H};X_{1}X_{2},gx_{2})-\gamma
_{2}B(x_{2}\otimes 1_{H};GX_{1},gx_{2}) \\
&&-B(x_{2}\otimes x_{2};X_{1},gx_{2})=0.
\end{eqnarray*}%
In view of the form of the elements we get

\begin{eqnarray*}
- &&2\beta _{2}B(g\otimes 1_{H};1_{A},gx_{1}x_{2})+\gamma _{2}B(x_{2}\otimes
1_{H};G,gx_{1}x_{2}) \\
&&-2B(x_{2}\otimes 1_{H};1_{A},gx_{1})=0.
\end{eqnarray*}%
This is $\left( \ref{X2,x2,X2F11,gx1x2}\right) .$

\paragraph{\textbf{Equality} $\left( \protect\ref{X2F41}\right) $}

rewrites as%
\begin{eqnarray*}
B(x_{2}\otimes 1_{H};1_{A},f) &=&-B(gx_{2}\otimes g;1_{A},f)-\gamma
_{2}B(gx_{2}\otimes g;GX_{2},f) \\
&&-\lambda B(gx_{2}\otimes g;X_{1}X_{2},f)+B(x_{2}\otimes x_{2};X_{2},f).
\end{eqnarray*}

\subparagraph{Case $f=1_{H}$}

\begin{eqnarray*}
B(x_{2}\otimes 1_{H};1_{A},1_{H}) &=&-B(gx_{2}\otimes g;1_{A},1_{H})-\gamma
_{2}B(gx_{2}\otimes g;GX_{2},1_{H}) \\
&&-\lambda B(gx_{2}\otimes g;X_{1}X_{2},1_{H})+B(x_{2}\otimes
x_{2};X_{2},1_{H}).
\end{eqnarray*}%
By applying $\left( \ref{eq.10}\right) $ this rewrites as%
\begin{eqnarray*}
&&2B(x_{2}\otimes 1_{H};1_{A},1_{H})+\gamma _{2}B(x_{2}\otimes
1_{H};GX_{2},1_{H}) \\
&&+\lambda B(x_{2}\otimes 1_{H};X_{1}X_{2},1_{H})-B(x_{2}\otimes
x_{2};X_{2},1_{H})=0.
\end{eqnarray*}%
In view of the form of the elements we get%
\begin{eqnarray}
&&2B(x_{2}\otimes 1_{H};1_{A},1_{H})+\gamma _{2}\left[ B(g\otimes
1_{H};G,1_{H})+B(x_{2}\otimes \ 1_{H};G,x_{2}\right]  \label{X2,x2,X2F41,1H}
\\
&&+\lambda \left[ B(g\otimes 1_{H};1_{A},x_{1})+B(x_{2}\otimes \
1_{H};1_{A},x_{1}x_{2})\right] =0.  \notag
\end{eqnarray}

\subparagraph{Case $f=x_{1}x_{2}$}

\begin{eqnarray*}
B(x_{2}\otimes 1_{H};1_{A},x_{1}x_{2}) &=&-B(gx_{2}\otimes
g;1_{A},x_{1}x_{2})-\gamma _{2}B(gx_{2}\otimes g;GX_{2},x_{1}x_{2}) \\
&&-\lambda B(gx_{2}\otimes g;X_{1}X_{2},x_{1}x_{2})+B(x_{2}\otimes
x_{2};X_{2},x_{1}x_{2}).
\end{eqnarray*}%
By applying $\left( \ref{eq.10}\right) $ this rewrites as%
\begin{eqnarray*}
&&2B(x_{2}\otimes 1_{H};1_{A},x_{1}x_{2})+\gamma _{2}B(x_{2}\otimes
1_{H};GX_{2},x_{1}x_{2}) \\
&&+\lambda B(x_{2}\otimes 1_{H};X_{1}X_{2},x_{1}x_{2})-B(x_{2}\otimes
x_{2};X_{2},x_{1}x_{2})=0.
\end{eqnarray*}%
In view of the form of the elements we get%
\begin{equation*}
2B(x_{2}\otimes 1_{H};1_{A},x_{1}x_{2})+\gamma _{2}B(g\otimes
1_{H};G,x_{1}x_{2})=0
\end{equation*}%
which is $\left( \ref{G,x2; GF6,x2}\right) .$

\subparagraph{Case $f=gx_{1}$%
\protect\begin{eqnarray*}
B(x_{2}\otimes 1_{H};1_{A},gx_{1}) &=&-B(gx_{2}\otimes g;1_{A},gx_{1})-%
\protect\gamma _{2}B(gx_{2}\otimes g;GX_{2},gx_{1}) \\
&&-\protect\lambda B(gx_{2}\otimes g;X_{1}X_{2},gx_{1})+B(x_{2}\otimes
x_{2};X_{2},gx_{1}).
\protect\end{eqnarray*}%
}

By applying $\left( \ref{eq.10}\right) $ this rewrites as%
\begin{eqnarray*}
&&-\gamma _{2}B(x_{2}\otimes 1_{H};GX_{2},gx_{1})-\lambda B(x_{2}\otimes
1_{H};X_{1}X_{2},gx_{1}) \\
&&-B(x_{2}\otimes x_{2};X_{2},gx_{1})=0.
\end{eqnarray*}%
In view of the form of the elements we get

\begin{equation*}
\gamma _{2}\left[ B(g\otimes 1_{H};G,gx_{1})+B(x_{2}\otimes \
1_{H};G,gx_{1}x_{2})\right] =0
\end{equation*}%
which is $\left( \ref{G,x2, GF2,gx1}\right) .$

\subparagraph{Case $f=gx_{2}$}

\begin{eqnarray*}
B(x_{2}\otimes 1_{H};1_{A},gx_{2}) &=&-B(gx_{2}\otimes
g;1_{A},gx_{2})-\gamma _{2}B(gx_{2}\otimes g;GX_{2},gx_{2}) \\
&&-\lambda B(gx_{2}\otimes g;X_{1}X_{2},gx_{2})+B(x_{2}\otimes
x_{2};X_{2},gx_{2}).
\end{eqnarray*}%
By applying $\left( \ref{eq.10}\right) $ this rewrites as%
\begin{eqnarray*}
&&-\gamma _{2}B(x_{2}\otimes 1_{H};GX_{2},gx_{2})-\lambda B(x_{2}\otimes
1_{H};X_{1}X_{2},gx_{2}) \\
&&-B(x_{2}\otimes x_{2};X_{2},gx_{2})=0.
\end{eqnarray*}%
In view of the form of the elements we get

\begin{eqnarray}
&&\gamma _{2}B(g\otimes 1_{H};G,gx_{2})-\lambda B\left( g\otimes
1_{H};1_{A},gx_{1}x_{2}\right)  \label{X2,x2,X2F41,gx2} \\
&&+2\left[ +B(g\otimes 1_{H};1_{A},g)+B(x_{2}\otimes \ 1_{H};1_{A},gx_{2})%
\right] =0.  \notag
\end{eqnarray}

\paragraph{\textbf{Equality} $\left( \protect\ref{X2F51}\right) $}

rewrites as

\subparagraph{Case $f=g$}

\begin{eqnarray*}
B(x_{2}\otimes 1_{H};X_{1},g) &=&B(gx_{2}\otimes g;X_{1},g) \\
&&-\gamma _{2}B(gx_{2}\otimes g;GX_{1}X_{2},g)+B(x_{2}\otimes
x_{2};X_{1}X_{2},g)
\end{eqnarray*}

By applying $\left( \ref{eq.10}\right) $ this rewrites as
\begin{equation*}
\gamma _{2}B(x_{2}\otimes 1_{H};GX_{1}X_{2},g)-B(x_{2}\otimes
x_{2};X_{1}X_{2},g)=0
\end{equation*}%
In view of the form of the elements we get

\begin{equation*}
\gamma _{2}\left[ B(g\otimes 1_{H};G,gx_{1})+B(x_{2}\otimes \
1_{H};G,gx_{1}x_{2})\right] =0
\end{equation*}%
which is $\left( \ref{G,x2, GF4,gx1}\right) .$

\subparagraph{Case $f=x_{2}$}

\begin{eqnarray*}
B(x_{2}\otimes 1_{H};X_{1},x_{2}) &=&+B(gx_{2}\otimes g;X_{1},x_{2}) \\
&&-\gamma _{2}B(gx_{2}\otimes g;GX_{1}X_{2},x_{2})+B(x_{2}\otimes
x_{2};X_{1}X_{2},x_{2})
\end{eqnarray*}%
By applying $\left( \ref{eq.10}\right) $ this rewrites as

\begin{gather*}
2B(x_{2}\otimes 1_{H};X_{1},x_{2})-\gamma _{2}B(gx_{2}\otimes
g;GX_{1}X_{2},x_{2}) \\
-B(x_{2}\otimes x_{2};X_{1}X_{2},x_{2})=0
\end{gather*}%
In view of the form of the elements we get%
\begin{equation*}
2B(x_{2}\otimes 1_{H};1_{A},x_{1}x_{2})+\gamma _{2}B\left( g\otimes
1_{H};G,x_{1}x_{2}\right) =0
\end{equation*}%
which is $\left( \ref{G,x2, GF2,x1x2}\right) .$

\paragraph{\textbf{Equality} $\left( \protect\ref{X2F61}\right) $}

rewrites as%
\begin{eqnarray*}
\beta _{2}B(x_{2}\otimes 1_{H};GX_{1}X_{2},f) &=&-\beta _{2}B(gx_{2}\otimes
g;GX_{1}X_{2},f) \\
&&+B(x_{2}\otimes x_{2};GX_{1},f)
\end{eqnarray*}

\subparagraph{Case $f=g$}

\begin{eqnarray*}
\beta _{2}B(x_{2}\otimes 1_{H};GX_{1}X_{2},g) &=&-\beta _{2}B(gx_{2}\otimes
g;GX_{1}X_{2},g) \\
&&+B(x_{2}\otimes x_{2};GX_{1},g)
\end{eqnarray*}%
By applying $\left( \ref{eq.10}\right) $ this rewrites as%
\begin{eqnarray*}
&&2\beta _{2}B(x_{2}\otimes 1_{H};GX_{1}X_{2},g) \\
&&-B(x_{2}\otimes x_{2};GX_{1},g)=0
\end{eqnarray*}%
In view of the form of the elements we get%
\begin{equation*}
\beta _{2}\left[ B(g\otimes 1_{H};G,gx_{1})+B(x_{2}\otimes \
1_{H};G,gx_{1}x_{2})\right]
\end{equation*}%
which is $\left( \ref{X2,x2,X2F21,gx1}\right) .$

\subparagraph{Case $f=x_{2}$}

\begin{eqnarray*}
\beta _{2}B(x_{2}\otimes 1_{H};GX_{1}X_{2},x_{2}) &=&-\beta
_{2}B(gx_{2}\otimes g;GX_{1}X_{2},x_{2}) \\
&&+B(x_{2}\otimes x_{2};GX_{1},x_{2}).
\end{eqnarray*}%
By applying $\left( \ref{eq.10}\right) $ this rewrites as%
\begin{equation*}
B(x_{2}\otimes x_{2};GX_{1},x_{2})=0
\end{equation*}%
which is already known.

\paragraph{\textbf{Equality} $\left( \protect\ref{X2F71}\right) $}

rewrites as%
\begin{eqnarray*}
B(x_{2}\otimes 1_{H};G,f) &=&B(gx_{2}\otimes g;G,f)+\lambda B(gx_{2}\otimes
g;GX_{1}X_{2},f) \\
&&+B(x_{2}\otimes x_{2};GX_{2},f)
\end{eqnarray*}

\subparagraph{Case $f=g$}

\begin{eqnarray*}
B(x_{2}\otimes 1_{H};G,g) &=&B(gx_{2}\otimes g;G,g)+\lambda B(gx_{2}\otimes
g;GX_{1}X_{2},g) \\
&&+B(x_{2}\otimes x_{2};GX_{2},g)
\end{eqnarray*}%
By applying $\left( \ref{eq.10}\right) $ this rewrites as
\begin{equation*}
-\lambda B(x_{2}\otimes 1_{H};GX_{1}X_{2},g)-B(x_{2}\otimes x_{2};GX_{2},g)=0
\end{equation*}%
In view of the form of the elements we get%
\begin{equation}
\lambda \left[ B(g\otimes 1_{H};G,gx_{1})+B(x_{2}\otimes \
1_{H};G,gx_{1}x_{2})\right] =0  \label{X2,x2,X2F71,g}
\end{equation}

\subparagraph{Case $f=x_{1}$}

\begin{eqnarray*}
B(x_{2}\otimes 1_{H};G,x_{1}) &=&B(gx_{2}\otimes g;G,x_{1})+\lambda
B(gx_{2}\otimes g;GX_{1}X_{2},x_{1}) \\
&&+B(x_{2}\otimes x_{2};GX_{2},x_{1})
\end{eqnarray*}%
By applying $\left( \ref{eq.10}\right) $ this rewrites as
\begin{eqnarray*}
&&2B(x_{2}\otimes 1_{H};G,x_{1})+\lambda B(x_{2}\otimes
1_{H};GX_{1}X_{2},x_{1}) \\
&&-B(x_{2}\otimes x_{2};GX_{2},x_{1})=0
\end{eqnarray*}%
In view of the form of the elements we get%
\begin{equation*}
B(x_{2}\otimes 1_{H};G,x_{1})=0
\end{equation*}%
which is $\left( \ref{X2,g,X2F21,x1}\right) .$

\subparagraph{Case $f=x_{2}$}

\begin{eqnarray*}
B(x_{2}\otimes 1_{H};G,x_{2}) &=&B(gx_{2}\otimes g;G,x_{2})+\lambda
B(gx_{2}\otimes g;GX_{1}X_{2},x_{2}) \\
&&+B(x_{2}\otimes x_{2};GX_{2},x_{2}).
\end{eqnarray*}%
By applying $\left( \ref{eq.10}\right) $ this rewrites as
\begin{eqnarray*}
&&2B(x_{2}\otimes 1_{H};G,x_{2})+\lambda B(x_{2}\otimes
1_{H};GX_{1}X_{2},x_{2}) \\
&&+B(x_{2}\otimes x_{2};GX_{2},x_{2})=0.
\end{eqnarray*}%
In view of the form of the elements we get%
\begin{equation*}
\lambda B(g\otimes 1_{H};G,x_{1}x_{2})-2B(x_{2}\otimes 1_{H};G,x_{2})=0
\end{equation*}%
which is $\left( \ref{X2,g,X2F21,x2}\right) .$

\subparagraph{Case $f=gx_{1}x_{2}$}

\begin{eqnarray*}
B(x_{2}\otimes 1_{H};G,gx_{1}x_{2}) &=&B(gx_{2}\otimes
g;G,gx_{1}x_{2})+\lambda B(gx_{2}\otimes g;GX_{1}X_{2},gx_{1}x_{2}) \\
&&+B(x_{2}\otimes x_{2};GX_{2},gx_{1}x_{2})
\end{eqnarray*}%
By applying $\left( \ref{eq.10}\right) $ this rewrites as
\begin{equation*}
-\lambda B(x_{2}\otimes 1_{H};GX_{1}X_{2},gx_{1}x_{2})-B(x_{2}\otimes
x_{2};GX_{2},gx_{1}x_{2})=0
\end{equation*}%
In view of the form of the elements we get

\begin{equation*}
B(g\otimes 1_{H};G,gx_{1})+B(x_{2}\otimes \ 1_{H};G,gx_{1}x_{2})=0
\end{equation*}%
which is $\left( \ref{X2,g,X2F71,gx1}\right) .$

\paragraph{\textbf{Equality} $\left( \protect\ref{X2F81}\right) $}

rewrites as%
\begin{eqnarray*}
B(x_{2}\otimes \ 1_{H};GX_{1},f) &=&-B(gx_{2}\otimes g;GX_{1},f) \\
&&+B(x_{2}\otimes x_{2};GX_{1}X_{2},f)
\end{eqnarray*}

\subparagraph{Case $f=1_{H}$}

\begin{eqnarray*}
B(x_{2}\otimes \ 1_{H};GX_{1},1_{H}) &=&-B(gx_{2}\otimes g;GX_{1},1_{H}) \\
&&+B(x_{2}\otimes x_{2};GX_{1}X_{2},1_{H}).
\end{eqnarray*}%
By applying $\left( \ref{eq.10}\right) $ this rewrites as%
\begin{equation*}
2B(x_{2}\otimes \ 1_{H};GX_{1},1_{H})-B(x_{2}\otimes
x_{2};GX_{1}X_{2},1_{H})=0
\end{equation*}%
In view of the form of the elements we get

\begin{equation*}
B(x_{2}\otimes 1_{H};G,x_{1})=0
\end{equation*}%
which is $\left( \ref{X2,g,X2F21,x1}\right) .$

\subparagraph{Case $f=gx_{2}$}

\begin{eqnarray*}
B(x_{2}\otimes \ 1_{H};GX_{1},gx_{2}) &=&-B(gx_{2}\otimes g;GX_{1},gx_{2}) \\
&&+B(x_{2}\otimes x_{2};GX_{1}X_{2},gx_{2}).
\end{eqnarray*}%
By applying $\left( \ref{eq.10}\right) $ this rewrites as%
\begin{equation*}
B(x_{2}\otimes x_{2};GX_{1}X_{2},gx_{2})=0.
\end{equation*}%
In view of the form of the elements we get%
\begin{equation*}
B(g\otimes 1_{H};G,gx_{1})+B(x_{2}\otimes \ 1_{H};G,gx_{1}x_{2})=0
\end{equation*}%
which is $\left( \ref{X2,g,X2F71,gx1}\right) .$

\subsubsection{Case $x_{1}x_{2}\otimes 1_{H}$}

\paragraph{\textbf{Equality} $\left( \protect\ref{X2F11}\right) $}

rewrites as%
\begin{eqnarray*}
\beta _{2}B(x_{1}x_{2}\otimes 1_{H};X_{2},f) &=&\gamma
_{2}B(gx_{1}x_{2}\otimes g;G,f)+\lambda B(gx_{1}x_{2}\otimes
g;X_{1},f)+\beta _{2}B(gx_{1}x_{2}\otimes g;X_{2},f) \\
&&+B(x_{1}x_{2}\otimes x_{2};1_{A},f).
\end{eqnarray*}

\subparagraph{Case $f=g$}

\begin{eqnarray*}
\beta _{2}B(x_{1}x_{2}\otimes 1_{H};X_{2},g) &=&\gamma
_{2}B(gx_{1}x_{2}\otimes g;G,g)+\lambda B(gx_{1}x_{2}\otimes
g;X_{1},g)+\beta _{2}B(gx_{1}x_{2}\otimes g;X_{2},g) \\
&&+B(x_{1}x_{2}\otimes x_{2};1_{A},g).
\end{eqnarray*}%
By applying $\left( \ref{eq.10}\right) $ this rewrites as
\begin{eqnarray*}
&&-\gamma _{2}B(x_{1}x_{2}\otimes 1_{H};G,g)-\lambda B(x_{1}x_{2}\otimes
1_{H};X_{1},g)+ \\
&&-B(x_{1}x_{2}\otimes x_{2};1_{A},g)=0.
\end{eqnarray*}%
In view of the form of the elements we get%
\begin{equation}
\gamma _{2}B(x_{1}x_{2}\otimes 1_{H};G,g)+\lambda B(x_{1}x_{2}\otimes
1_{H};X_{1},g)=0.  \label{X2,x1x2,X2F11,g}
\end{equation}

\subparagraph{Case $f=x_{1}$}

\begin{eqnarray*}
\beta _{2}B(x_{1}x_{2}\otimes 1_{H};X_{2},x_{1}) &=&\gamma
_{2}B(gx_{1}x_{2}\otimes g;G,x_{1})+\lambda B(gx_{1}x_{2}\otimes
g;X_{1},x_{1})+\beta _{2}B(gx_{1}x_{2}\otimes g;X_{2},x_{1}) \\
&&+B(x_{1}x_{2}\otimes x_{2};1_{A},x_{1}).
\end{eqnarray*}%
By applying $\left( \ref{eq.10}\right) $ this rewrites as
\begin{gather*}
2\beta _{2}B(x_{1}x_{2}\otimes 1_{H};X_{2},x_{1})+\gamma
_{2}B(x_{1}x_{2}\otimes 1_{H};G,x_{1}) \\
+\lambda B(x_{1}x_{2}\otimes 1_{H};X_{1},x_{1})-B(x_{1}x_{2}\otimes
x_{2};1_{A},x_{1})=0.
\end{gather*}%
In view of the form of the elements we get%
\begin{gather}
2\beta _{2}B(x_{1}x_{2}\otimes 1_{H};X_{2},x_{1})+\gamma
_{2}B(x_{1}x_{2}\otimes 1_{H};G,x_{1})  \label{X2,x1x2,X2F11,x1} \\
+\lambda B(x_{1}x_{2}\otimes 1_{H};X_{1},x_{1})=0.  \notag
\end{gather}

\subparagraph{Case $f=x_{2}$}

\begin{eqnarray*}
\beta _{2}B(x_{1}x_{2}\otimes 1_{H};X_{2},x_{2}) &=&\gamma
_{2}B(gx_{1}x_{2}\otimes g;G,x_{2})+\lambda B(gx_{1}x_{2}\otimes
g;X_{1},x_{2})+\beta _{2}B(gx_{1}x_{2}\otimes g;X_{2},x_{2}) \\
&&+B(x_{1}x_{2}\otimes x_{2};1_{A},x_{2}).
\end{eqnarray*}%
By applying $\left( \ref{eq.10}\right) $ this rewrites as
\begin{gather*}
2\beta _{2}B(x_{1}x_{2}\otimes 1_{H};X_{2},x_{2})+\gamma
_{2}B(x_{1}x_{2}\otimes 1_{H};G,x_{2}) \\
+\lambda B(x_{1}x_{2}\otimes 1_{H};X_{1},x_{2})+B(x_{1}x_{2}\otimes
x_{2};1_{A},x_{2})=0.
\end{gather*}%
In view of the form of the elements we get%
\begin{gather}
2\beta _{2}B(x_{1}x_{2}\otimes 1_{H};X_{2},x_{2})+\gamma
_{2}B(x_{1}x_{2}\otimes 1_{H};G,x_{2})  \label{X2,x1x2,X2F11,x2} \\
+\lambda B(x_{1}x_{2}\otimes 1_{H};X_{1},x_{2})=0.  \notag
\end{gather}

\subparagraph{Case $f=gx_{1}x_{2}$}

\begin{eqnarray*}
\beta _{2}B(x_{1}x_{2}\otimes 1_{H};X_{2},gx_{1}x_{2}) &=&\gamma
_{2}B(gx_{1}x_{2}\otimes g;G,gx_{1}x_{2})+\lambda B(gx_{1}x_{2}\otimes
g;X_{1},gx_{1}x_{2})+ \\
&&+\beta _{2}B(gx_{1}x_{2}\otimes g;X_{2},gx_{1}x_{2})+B(x_{1}x_{2}\otimes
x_{2};1_{A},gx_{1}x_{2}).
\end{eqnarray*}%
By applying $\left( \ref{eq.10}\right) $ this rewrites as
\begin{eqnarray*}
&&\gamma _{2}B(x_{1}x_{2}\otimes 1_{H};G,gx_{1}x_{2})+\lambda
B(x_{1}x_{2}\otimes 1_{H};X_{1},gx_{1}x_{2}) \\
&&+B(x_{1}x_{2}\otimes x_{2};1_{A},gx_{1}x_{2})=0.
\end{eqnarray*}%
In view of the form of the elements we get%
\begin{eqnarray}
&&\gamma _{2}B(x_{1}x_{2}\otimes 1_{H};G,gx_{1}x_{2})+\lambda
B(x_{1}x_{2}\otimes 1_{H};X_{1},gx_{1}x_{2})  \label{X2,x1x2,X2F11,gx1x2} \\
&&+2B(x_{1}x_{2}\otimes 1_{H};1_{A},gx_{1})=0.  \notag
\end{eqnarray}

\paragraph{\textbf{Equality} $\left( \protect\ref{X2F21}\right) $}

rewrites as

\subparagraph{%
\protect\begin{gather*}
\protect\beta _{2}B(x_{1}x_{2}\otimes 1_{H};GX_{2},f)=-\protect\beta %
_{2}B(gx_{1}x_{2}\otimes g;GX_{2},f)+ \\
-\protect\lambda B(gx_{1}x_{2}\otimes g;GX_{1},f)+B(x_{1}x_{2}\otimes
x_{2};G,f).
\protect\end{gather*}%
Case $f=1_{H}$}

\begin{gather*}
\beta _{2}B(x_{1}x_{2}\otimes 1_{H};GX_{2},1_{H})=-\beta
_{2}B(gx_{1}x_{2}\otimes g;GX_{2},1_{H})+ \\
-\lambda B(gx_{1}x_{2}\otimes g;GX_{1},1_{H})+B(x_{1}x_{2}\otimes
x_{2};G,1_{H}).
\end{gather*}%
By applying $\left( \ref{eq.10}\right) $ this rewrites as%
\begin{eqnarray*}
&&2\beta _{2}B(x_{1}x_{2}\otimes 1_{H};GX_{2},1_{H})+\lambda
B(x_{1}x_{2}\otimes 1_{H};GX_{1},1_{H}) \\
&&-B(x_{1}x_{2}\otimes x_{2};G,1_{H})=0.
\end{eqnarray*}%
In view of the form of the elements we get%
\begin{equation}
2\beta _{2}B(x_{1}x_{2}\otimes 1_{H};GX_{2},1_{H})+\lambda
B(x_{1}x_{2}\otimes 1_{H};GX_{1},1_{H})=0.  \label{X2,x1x2,X2F21,1H}
\end{equation}

\subparagraph{Case $f=x_{1}x_{2}$}

\begin{eqnarray*}
&&\beta _{2}B(x_{1}x_{2}\otimes 1_{H};GX_{2},x_{1}x_{2}) \\
&=&-\beta _{2}B(gx_{1}x_{2}\otimes g;GX_{2},x_{1}x_{2})-\lambda
B(gx_{1}x_{2}\otimes g;GX_{1},x_{1}x_{2})+ \\
&&+B(x_{1}x_{2}\otimes x_{2};G,x_{1}x_{2}).
\end{eqnarray*}%
By applying $\left( \ref{eq.10}\right) $ this rewrites as%
\begin{eqnarray*}
&&2\beta _{2}B(x_{1}x_{2}\otimes 1_{H};GX_{2},x_{1}x_{2})+\lambda
B(x_{1}x_{2}\otimes 1_{H};GX_{1},x_{1}x_{2}) \\
&&-B(x_{1}x_{2}\otimes x_{2};G,x_{1}x_{2})=0.
\end{eqnarray*}%
In view of the form of the elements we get%
\begin{equation*}
2\beta _{2}B(x_{1}x_{2}\otimes 1_{H};GX_{2},x_{1}x_{2})+\lambda
B(x_{1}x_{2}\otimes 1_{H};GX_{1},x_{1}x_{2})=0.
\end{equation*}%
In view of $\left( \ref{X1,x1x2,X1F61,x1}\right) $ we get%
\begin{equation}
\beta _{2}B(x_{1}x_{2}\otimes 1_{H};GX_{2},x_{1}x_{2})=0.
\label{X2,x1x2,X2F21,x1x2}
\end{equation}

\subparagraph{Case $f=gx_{1}$}

\begin{eqnarray*}
&&\beta _{2}B(x_{1}x_{2}\otimes 1_{H};GX_{2},gx_{1}) \\
&=&-\beta _{2}B(gx_{1}x_{2}\otimes g;GX_{2},gx_{1})-\lambda
B(gx_{1}x_{2}\otimes g;GX_{1},gx_{1})+ \\
&&+B(x_{1}x_{2}\otimes x_{2};G,gx_{1}).
\end{eqnarray*}%
By applying $\left( \ref{eq.10}\right) $ this rewrites as%
\begin{equation*}
-\lambda B(x_{1}x_{2}\otimes 1_{H};GX_{1},gx_{1})-B(x_{1}x_{2}\otimes
x_{2};G,gx_{1})=0.
\end{equation*}%
In view of the form of the elements we get%
\begin{equation}
\lambda B(x_{1}x_{2}\otimes 1_{H};GX_{1},gx_{1})=0.
\label{X2,x1x2,X2F21,gx1}
\end{equation}

\subparagraph{Case $f=gx_{2}$}

\begin{eqnarray*}
&&\beta _{2}B(x_{1}x_{2}\otimes 1_{H};GX_{2},gx_{2}) \\
&=&-\beta _{2}B(gx_{1}x_{2}\otimes g;GX_{2},gx_{2})-\lambda
B(gx_{1}x_{2}\otimes g;GX_{1},gx_{2})+ \\
&&+B(x_{1}x_{2}\otimes x_{2};G,gx_{2}).
\end{eqnarray*}%
By applying $\left( \ref{eq.10}\right) $ this rewrites as%
\begin{equation*}
\lambda B(x_{1}x_{2}\otimes 1_{H};GX_{1},gx_{2})-B(x_{1}x_{2}\otimes
x_{2};G,gx_{2})=0.
\end{equation*}%
In view of the form of the elements we get%
\begin{equation}
2B(x_{1}x_{2}\otimes 1_{H};G,g)-\lambda B(x_{1}x_{2}\otimes
1_{H};GX_{1},gx_{2})=0.  \label{X2,x1x2,X2F21,gx2}
\end{equation}

\paragraph{\textbf{Equality} $\left( \protect\ref{X2F31}\right) $}

rewrites as%
\begin{eqnarray*}
\beta _{2}B(x_{1}x_{2}\otimes 1_{H};X_{1}X_{2},f) &=&-\beta
_{2}B(gx_{1}x_{2}\otimes g;X_{1}X_{2},f) \\
&&+\gamma _{2}B(gx_{1}x_{2}\otimes g;GX_{1},f)+B(x_{1}x_{2}\otimes
x_{2};X_{1},f).
\end{eqnarray*}

\subparagraph{Case $f=1_{H}$}

\begin{eqnarray*}
\beta _{2}B(x_{1}x_{2}\otimes 1_{H};X_{1}X_{2},1_{H}) &=&-\beta
_{2}B(gx_{1}x_{2}\otimes g;X_{1}X_{2},1_{H}) \\
&&+\gamma _{2}B(gx_{1}x_{2}\otimes g;GX_{1},1_{H})+B(x_{1}x_{2}\otimes
x_{2};X_{1},1_{H}).
\end{eqnarray*}%
By applying $\left( \ref{eq.10}\right) $ this rewrites as%
\begin{eqnarray*}
&&2\beta _{2}B(x_{1}x_{2}\otimes 1_{H};X_{1}X_{2},1_{H})-\gamma
_{2}B(x_{1}x_{2}\otimes 1_{H};GX_{1},1_{H}) \\
&&-B(x_{1}x_{2}\otimes x_{2};X_{1},1_{H})=0.
\end{eqnarray*}%
In view of the form of the elements we get%
\begin{eqnarray}
&&2\beta _{2}\left[
\begin{array}{c}
+1-B(x_{1}x_{2}\otimes 1_{H};1_{A},x_{1}x_{2}) \\
-B(x_{1}x_{2}\otimes 1_{H};X_{2},x_{1})+B(x_{1}x_{2}\otimes
1_{H};X_{1},x_{2})%
\end{array}%
\right]  \label{X2,x1x2,X2F31,1H} \\
&&-\gamma _{2}B(x_{1}x_{2}\otimes 1_{H};GX_{1},1_{H})=0.  \notag
\end{eqnarray}

\subparagraph{Case $f=x_{1}x_{2}$}

\begin{eqnarray*}
\beta _{2}B(x_{1}x_{2}\otimes 1_{H};X_{1}X_{2},x_{1}x_{2}) &=&-\beta
_{2}B(gx_{1}x_{2}\otimes g;X_{1}X_{2},x_{1}x_{2}) \\
&&+\gamma _{2}B(gx_{1}x_{2}\otimes g;GX_{1},x_{1}x_{2})+B(x_{1}x_{2}\otimes
x_{2};X_{1},x_{1}x_{2}).
\end{eqnarray*}%
By applying $\left( \ref{eq.10}\right) $ this rewrites as%
\begin{eqnarray*}
&&2\beta _{2}B(x_{1}x_{2}\otimes 1_{H};X_{1}X_{2},x_{1}x_{2})-\gamma
_{2}B(x_{1}x_{2}\otimes 1_{H};GX_{1},x_{1}x_{2}) \\
&&-B(x_{1}x_{2}\otimes x_{2};X_{1},x_{1}x_{2})=0.
\end{eqnarray*}%
In view of the form of the elements we get%
\begin{equation*}
\gamma _{2}B(x_{1}x_{2}\otimes 1_{H};GX_{1},x_{1}x_{2})=0.
\end{equation*}%
which is $\left( \ref{G,x1x2, GF6,x1}\right) .$

\subparagraph{Case $f=gx_{1}$}

\begin{eqnarray*}
&&\beta _{2}B(x_{1}x_{2}\otimes 1_{H};X_{1}X_{2},gx_{1}) \\
&&+\gamma _{2}B(gx_{1}x_{2}\otimes g;GX_{1},gx_{1})+B(x_{1}x_{2}\otimes
x_{2};X_{1},gx_{1}).
\end{eqnarray*}%
By applying $\left( \ref{eq.10}\right) $ this rewrites as%
\begin{equation*}
\gamma _{2}B(x_{1}x_{2}\otimes 1_{H};GX_{1},gx_{1})-B(x_{1}x_{2}\otimes
x_{2};X_{1},gx_{1})=0.
\end{equation*}%
In view of the form of the elements we get%
\begin{equation}
\gamma _{2}B(x_{1}x_{2}\otimes 1_{H};GX_{1},gx_{1})=0
\label{X2,x1x2,X2F31,gx1}
\end{equation}

\subparagraph{Case $f=gx_{2}$}

\begin{eqnarray*}
\beta _{2}B(x_{1}x_{2}\otimes 1_{H};X_{1}X_{2},gx_{2}) &=&-\beta
_{2}B(gx_{1}x_{2}\otimes g;X_{1}X_{2},gx_{2}) \\
&&+\gamma _{2}B(gx_{1}x_{2}\otimes g;GX_{1},gx_{2})+B(x_{1}x_{2}\otimes
x_{2};X_{1},gx_{2}).
\end{eqnarray*}%
By applying $\left( \ref{eq.10}\right) $ this rewrites as%
\begin{equation*}
\gamma _{2}B(x_{1}x_{2}\otimes 1_{H};GX_{1},gx_{2})-B(x_{1}x_{2}\otimes
x_{2};X_{1},gx_{2})=0.
\end{equation*}%
In view of the form of the elements we get%
\begin{equation}
\gamma _{2}B(x_{1}x_{2}\otimes 1_{H};GX_{1},gx_{2})+2B(x_{1}x_{2}\otimes
1_{H};X_{1},g)=0.  \label{X2,x1x2,X2F31,gx2}
\end{equation}

\paragraph{\textbf{Equality} $\left( \protect\ref{X2F41}\right) $}

rewrites as

\begin{eqnarray*}
B(x_{1}x_{2}\otimes 1_{H};1_{A},f) &=&B(gx_{1}x_{2}\otimes g;1_{A},f)+\gamma
_{2}B(gx_{1}x_{2}\otimes g;GX_{2},f) \\
+\lambda B(gx_{1}x_{2}\otimes g;X_{1}X_{2},f) &&+B(x_{1}x_{2}\otimes
x_{2};X_{2},f).
\end{eqnarray*}

\subparagraph{Case $f=1_{H}$}

\begin{eqnarray*}
B(x_{1}x_{2}\otimes 1_{H};1_{A},1_{H}) &=&B(gx_{1}x_{2}\otimes
g;1_{A},1_{H})+\gamma _{2}B(gx_{1}x_{2}\otimes g;GX_{2},1_{H}) \\
+\lambda B(gx_{1}x_{2}\otimes g;X_{1}X_{2},1_{H}) &&+B(x_{1}x_{2}\otimes
x_{2};X_{2},1_{H}).
\end{eqnarray*}%
By applying $\left( \ref{eq.10}\right) $ this rewrites as
\begin{eqnarray*}
&&-\gamma _{2}B(x_{1}x_{2}\otimes 1_{H};GX_{2},1_{H})-\lambda
B(x_{1}x_{2}\otimes 1_{H};X_{1}X_{2},1_{H}) \\
&&-B(x_{1}x_{2}\otimes x_{2};X_{2},1_{H})=0.
\end{eqnarray*}%
In view of the form of the elements we get%
\begin{gather}
-\gamma _{2}B(x_{1}x_{2}\otimes 1_{H};GX_{2},1_{H})+
\label{X2,x1x2,X2F41,1H} \\
-\lambda \left[
\begin{array}{c}
+1-B(x_{1}x_{2}\otimes 1_{H};1_{A},x_{1}x_{2}) \\
-B(x_{1}x_{2}\otimes 1_{H};X_{2},x_{1})+B(x_{1}x_{2}\otimes
1_{H};X_{1},x_{2})%
\end{array}%
\right] =0.  \notag
\end{gather}

\subparagraph{Case $f=x_{1}x_{2}$}

\begin{eqnarray*}
B(x_{1}x_{2}\otimes 1_{H};1_{A},x_{1}x_{2}) &=&B(gx_{1}x_{2}\otimes
g;1_{A},x_{1}x_{2})+\gamma _{2}B(gx_{1}x_{2}\otimes g;GX_{2},x_{1}x_{2}) \\
+\lambda B(gx_{1}x_{2}\otimes g;X_{1}X_{2},x_{1}x_{2})
&&+B(x_{1}x_{2}\otimes x_{2};X_{2},x_{1}x_{2}).
\end{eqnarray*}%
By applying $\left( \ref{eq.10}\right) $ this rewrites as
\begin{gather*}
-\gamma _{2}B(x_{1}x_{2}\otimes 1_{H};GX_{2},x_{1}x_{2})-\lambda
B(x_{1}x_{2}\otimes 1_{H};X_{1}X_{2},x_{1}x_{2}) \\
-B(x_{1}x_{2}\otimes x_{2};X_{2},x_{1}x_{2})=0.
\end{gather*}%
In view of the form of the elements we get%
\begin{equation*}
\gamma _{2}B(x_{1}x_{2}\otimes 1_{H};GX_{2},x_{1}x_{2})=0
\end{equation*}%
which is $\left( \ref{G,x1x2, GF6,x2}\right) .$

\subparagraph{Case $f=gx_{1}$}

\begin{eqnarray*}
B(x_{1}x_{2}\otimes 1_{H};1_{A},gx_{1}) &=&B(gx_{1}x_{2}\otimes
g;1_{A},gx_{1})+\gamma _{2}B(gx_{1}x_{2}\otimes g;GX_{2},gx_{1}) \\
+\lambda B(gx_{1}x_{2}\otimes g;X_{1}X_{2},gx_{1}) &&+B(x_{1}x_{2}\otimes
x_{2};X_{2},gx_{1}).
\end{eqnarray*}%
By applying $\left( \ref{eq.10}\right) $ this rewrites as
\begin{eqnarray*}
&&2B(x_{1}x_{2}\otimes 1_{H};1_{A},gx_{1})+\gamma _{2}B(x_{1}x_{2}\otimes
1_{H};GX_{2},gx_{1}) \\
+\lambda B(x_{1}x_{2}\otimes 1_{H};X_{1}X_{2},gx_{1}) &&-B(x_{1}x_{2}\otimes
x_{2};X_{2},gx_{1})=0.
\end{eqnarray*}%
In view of the form of the elements we get%
\begin{gather*}
2B(x_{1}x_{2}\otimes 1_{H};1_{A},gx_{1})+\gamma _{2}B(x_{1}x_{2}\otimes
1_{H};GX_{2},gx_{1}) \\
+\lambda B(x_{1}x_{2}\otimes 1_{H};X_{1},gx_{1}x_{2})=0
\end{gather*}%
In view of $\left( \ref{G,x1x2, GF6,gx1x2}\right) $we obtain%
\begin{equation*}
2B(x_{1}x_{2}\otimes 1_{H};1_{A},gx_{1})+\gamma _{2}B(x_{1}x_{2}\otimes
1_{H};GX_{2},gx_{1})=0
\end{equation*}%
which is $\left( \ref{G,x1x2, GF2,gx1}\right) $

\subparagraph{Case $f=gx_{2}$}

\begin{eqnarray*}
B(x_{1}x_{2}\otimes 1_{H};1_{A},gx_{2}) &=&B(gx_{1}x_{2}\otimes
g;1_{A},gx_{2})+\gamma _{2}B(gx_{1}x_{2}\otimes g;GX_{2},gx_{2}) \\
+\lambda B(gx_{1}x_{2}\otimes g;X_{1}X_{2},gx_{2}) &&+B(x_{1}x_{2}\otimes
x_{2};X_{2},gx_{2}).
\end{eqnarray*}%
By applying $\left( \ref{eq.10}\right) $ this rewrites as
\begin{eqnarray*}
&&2B(x_{1}x_{2}\otimes 1_{H};1_{A},gx_{2})+\gamma _{2}B(x_{1}x_{2}\otimes
1_{H};GX_{2},gx_{2}) \\
+\lambda B(x_{1}x_{2}\otimes 1_{H};X_{1}X_{2},gx_{2}) &&-B(x_{1}x_{2}\otimes
x_{2};X_{2},gx_{2})=0.
\end{eqnarray*}%
In view of the form of the elements we get%
\begin{gather}
2B(x_{1}x_{2}\otimes 1_{H};1_{A},gx_{2})+\gamma _{2}B(x_{1}x_{2}\otimes
1_{H};GX_{2},gx_{2})  \label{X2,x1x2,X2F41,gx2} \\
+\lambda B(x_{1}x_{2}\otimes 1_{H};X_{2},gx_{1}x_{2})+2B(x_{1}x_{2}\otimes
1_{H};X_{2},g)=0.  \notag
\end{gather}

\paragraph{\textbf{Equality} $\left( \protect\ref{X2F51}\right) $}

rewrites as

\begin{eqnarray*}
B(x_{1}x_{2}\otimes 1_{H};X_{1},f) &=&-B(gx_{1}x_{2}\otimes g;X_{1},f)+ \\
&&\gamma _{2}B(gx_{1}x_{2}\otimes g;GX_{1}X_{2},f)+B((x_{1}x_{2}\otimes
x_{2};X_{1}X_{2},f).
\end{eqnarray*}

\subparagraph{Case $f=g$}

\begin{eqnarray*}
B(x_{1}x_{2}\otimes 1_{H};X_{1},g) &=&-B(gx_{1}x_{2}\otimes g;X_{1},g)+ \\
&&\gamma _{2}B(gx_{1}x_{2}\otimes g;GX_{1}X_{2},g)+B(x_{1}x_{2}\otimes
x_{2};X_{1}X_{2},g).
\end{eqnarray*}%
By applying $\left( \ref{eq.10}\right) $ this rewrites as%
\begin{gather*}
2B(x_{1}x_{2}\otimes 1_{H};X_{1},g)-\gamma _{2}B(x_{1}x_{2}\otimes
1_{H};GX_{1}X_{2},g)+ \\
-B(x_{1}x_{2}\otimes x_{2};X_{1}X_{2},g)=0.
\end{gather*}%
In view of the form of the elements we get%
\begin{gather}
2B(x_{1}x_{2}\otimes 1_{H};X_{1},g)+  \label{X2,x1x2,X2F51,g} \\
-\gamma _{2}\left[
\begin{array}{c}
-B(x_{1}x_{2}\otimes 1_{H};G,gx_{1}x_{2})+ \\
B(x_{1}x_{2}\otimes 1_{H};GX_{2},gx_{1})-B(x_{1}x_{2}\otimes
1_{H};GX_{1},gx_{2})%
\end{array}%
\right] =0.  \notag
\end{gather}

\subparagraph{Case $f=x_{1}$}

\begin{eqnarray*}
B(x_{1}x_{2}\otimes 1_{H};X_{1},x_{1}) &=&-B(gx_{1}x_{2}\otimes
g;X_{1},x_{1})+ \\
&&\gamma _{2}B(gx_{1}x_{2}\otimes g;GX_{1}X_{2},x_{1})+B(x_{1}x_{2}\otimes
x_{2};X_{1}X_{2},x_{1}).
\end{eqnarray*}%
By applying $\left( \ref{eq.10}\right) $ this rewrites as%
\begin{equation*}
\gamma _{2}B(x_{1}x_{2}\otimes 1_{H};GX_{1}X_{2},x_{1})-B(x_{1}x_{2}\otimes
x_{2};X_{1}X_{2},x_{1})=0.
\end{equation*}%
In view of the form of the elements we get%
\begin{equation*}
\gamma _{2}B(x_{1}x_{2}\otimes 1_{H};GX_{1},x_{1}x_{2})=0
\end{equation*}%
which is $\left( \ref{G,x1x2, GF6,x1}\right) .$

\subparagraph{Case $f=x_{2}$}

\begin{eqnarray*}
B(x_{1}x_{2}\otimes 1_{H};X_{1},x_{2}) &=&-B(gx_{1}x_{2}\otimes
g;X_{1},x_{2})+ \\
&&\gamma _{2}B(gx_{1}x_{2}\otimes g;GX_{1}X_{2},x_{2})+B(x_{1}x_{2}\otimes
x_{2};X_{1}X_{2},x_{2}).
\end{eqnarray*}%
By applying $\left( \ref{eq.10}\right) $ this rewrites as%
\begin{equation*}
\gamma _{2}B(x_{1}x_{2}\otimes 1_{H};GX_{1}X_{2},x_{2})-B(x_{1}x_{2}\otimes
x_{2};X_{1}X_{2},x_{2})=0.
\end{equation*}%
In view of the form of the elements we get%
\begin{equation*}
\gamma _{2}B(x_{1}x_{2}\otimes 1_{H};GX_{2},x_{1}x_{2})=0
\end{equation*}%
which is $\left( \ref{G,x1x2, GF6,x2}\right) .$

\subparagraph{Case $f=gx_{1}x_{2}$}

\begin{eqnarray*}
&&B(x_{1}x_{2}\otimes 1_{H};X_{1},gx_{1}x_{2}) \\
&=&-B(gx_{1}x_{2}\otimes g;X_{1},gx_{1}x_{2})+ \\
&&+\gamma _{2}B(gx_{1}x_{2}\otimes
g;GX_{1}X_{2},gx_{1}x_{2})+B(x_{1}x_{2}\otimes x_{2};X_{1}X_{2},gx_{1}x_{2}).
\end{eqnarray*}%
By applying $\left( \ref{eq.10}\right) $ this rewrites as%
\begin{eqnarray*}
&&2B(x_{1}x_{2}\otimes 1_{H};X_{1},gx_{1}x_{2})-\gamma
_{2}B(gx_{1}x_{2}\otimes g;GX_{1}X_{2},gx_{1}x_{2}) \\
&&-B(x_{1}x_{2}\otimes x_{2};X_{1}X_{2},gx_{1}x_{2})=0
\end{eqnarray*}%
In view of the form of the elements we get%
\begin{equation*}
2B(x_{1}x_{2}\otimes 1_{H};X_{1},gx_{1}x_{2})-2B(x_{1}x_{2}\otimes
1_{H};X_{1},gx_{1}x_{2})=0
\end{equation*}%
which is trivial.

\paragraph{\textbf{Equality} $\left( \protect\ref{X2F61}\right) $}

rewrites as

\begin{eqnarray*}
\beta _{2}B(x_{1}x_{2}\otimes 1_{H};GX_{1}X_{2},f) &=&+\beta
_{2}B(gx_{1}x_{2}\otimes g;GX_{1}X_{2},f) \\
&&+B(x_{1}x_{2}\otimes x_{2};GX_{1},f).
\end{eqnarray*}

\subparagraph{Case $f=g$}

\begin{eqnarray*}
\beta _{2}B(x_{1}x_{2}\otimes 1_{H};GX_{1}X_{2},g) &=&+\beta
_{2}B(gx_{1}x_{2}\otimes g;GX_{1}X_{2},g) \\
&&+B(x_{1}x_{2}\otimes x_{2};GX_{1},g).
\end{eqnarray*}%
By applying $\left( \ref{eq.10}\right) $ this rewrites as
\begin{equation*}
B(x_{1}x_{2}\otimes x_{2};GX_{1},g)=0
\end{equation*}%
which is already known.

\subparagraph{Case $f=x_{1}$}

\begin{eqnarray*}
\beta _{2}B(x_{1}x_{2}\otimes 1_{H};GX_{1}X_{2},x_{1}) &=&+\beta
_{2}B(gx_{1}x_{2}\otimes g;GX_{1}X_{2},x_{1}) \\
&&+B(x_{1}x_{2}\otimes x_{2};GX_{1},x_{1}).
\end{eqnarray*}

By applying $\left( \ref{eq.10}\right) $ this rewrites as%
\begin{equation*}
2\beta _{2}B(x_{1}x_{2}\otimes 1_{H};GX_{1}X_{2},x_{1})-B(x_{1}x_{2}\otimes
x_{2};GX_{1},x_{1})=0.
\end{equation*}%
In view of the form of the elements we get%
\begin{equation}
\beta _{2}B(x_{1}x_{2}\otimes 1_{H};GX_{1},x_{1}x_{2})=0.
\label{X2,x1x2,X2F61,x1}
\end{equation}

\subparagraph{Case $f=x_{2}$}

\begin{eqnarray*}
\beta _{2}B(x_{1}x_{2}\otimes 1_{H};GX_{1}X_{2},x_{2}) &=&+\beta
_{2}B(gx_{1}x_{2}\otimes g;GX_{1}X_{2},x_{2}) \\
&&+B(x_{1}x_{2}\otimes x_{2};GX_{1},x_{2}).
\end{eqnarray*}%
By applying $\left( \ref{eq.10}\right) $ this rewrites as
\begin{equation*}
2\beta _{2}B(x_{1}x_{2}\otimes 1_{H};GX_{1}X_{2},x_{2})-B(x_{1}x_{2}\otimes
x_{2};GX_{1},x_{2})=0.
\end{equation*}%
In view of the form of the elements we get%
\begin{equation*}
\beta _{2}B(x_{1}x_{2}\otimes 1_{H};GX_{2},x_{1}x_{2})=0
\end{equation*}%
which is $\left( \ref{X2,x1x2,X2F21,x1x2}\right) .$

\subparagraph{Case $f=gx_{1}x_{2}$}

\begin{eqnarray*}
\beta _{2}B(x_{1}x_{2}\otimes 1_{H};GX_{1}X_{2},gx_{1}x_{2}) &=&+\beta
_{2}B(gx_{1}x_{2}\otimes g;GX_{1}X_{2},gx_{1}x_{2}) \\
&&+B(x_{1}x_{2}\otimes x_{2};GX_{1},gx_{1}x_{2}).
\end{eqnarray*}%
By applying $\left( \ref{eq.10}\right) $ this rewrites as
\begin{equation*}
B(x_{1}x_{2}\otimes x_{2};GX_{1},gx_{1}x_{2})=0.
\end{equation*}%
In view of the form of the elements we get%
\begin{equation}
B(x_{1}x_{2}\otimes 1_{H};GX_{1},gx_{1})=0.  \label{X2,x1x2,X2F61,gx1x2}
\end{equation}

\paragraph{\textbf{Equality} $\left( \protect\ref{X2F71}\right) $}

rewrites as

\begin{eqnarray*}
B(x_{1}x_{2}\otimes 1_{H};G,f) &=&-B(gx_{1}x_{2}\otimes g;G,f)-\lambda
B(gx_{1}x_{2}\otimes g;GX_{1}X_{2},f) \\
&&+B(x_{1}x_{2}\otimes x_{2};GX_{2},f).
\end{eqnarray*}

\subparagraph{Case $f=g$}

\begin{eqnarray*}
B(x_{1}x_{2}\otimes 1_{H};G,g) &=&-B(gx_{1}x_{2}\otimes g;G,g)-\lambda
B(gx_{1}x_{2}\otimes g;GX_{1}X_{2},g) \\
&&+B(x_{1}x_{2}\otimes x_{2};GX_{2},g).
\end{eqnarray*}%
By applying $\left( \ref{eq.10}\right) $ this rewrites as%
\begin{eqnarray*}
&&2B(x_{1}x_{2}\otimes 1_{H};G,g)+\lambda B(gx_{1}x_{2}\otimes
g;GX_{1}X_{2},g) \\
&&-B(x_{1}x_{2}\otimes x_{2};GX_{2},g)=0.
\end{eqnarray*}%
In view of the form of the elements we get%
\begin{gather}
2B(x_{1}x_{2}\otimes 1_{H};G,g)  \label{X2,x1x2,X2F71,g} \\
+\lambda \left[
\begin{array}{c}
-B(x_{1}x_{2}\otimes 1_{H};G,gx_{1}x_{2})+ \\
B(x_{1}x_{2}\otimes 1_{H};GX_{2},gx_{1})-B(x_{1}x_{2}\otimes
1_{H};GX_{1},gx_{2})%
\end{array}%
\right] =0.  \notag
\end{gather}

\subparagraph{Case $f=x_{1}$}

\begin{eqnarray*}
B(x_{1}x_{2}\otimes 1_{H};G,x_{1}) &=&-B(gx_{1}x_{2}\otimes
g;G,x_{1})-\lambda B(gx_{1}x_{2}\otimes g;GX_{1}X_{2},x_{1}) \\
&&+B(x_{1}x_{2}\otimes x_{2};GX_{2},x_{1}).
\end{eqnarray*}%
By applying $\left( \ref{eq.10}\right) $ this rewrites as%
\begin{equation*}
-\lambda B(x_{1}x_{2}\otimes 1_{H};GX_{1}X_{2},x_{1})-B(x_{1}x_{2}\otimes
x_{2};GX_{2},x_{1})=0.
\end{equation*}%
In view of the form of the elements we get%
\begin{equation*}
\lambda B(x_{1}x_{2}\otimes 1_{H};GX_{1},x_{1}x_{2})=0
\end{equation*}%
which is $\left( \ref{X1,x1x2,X1F61,x1}\right) .$

\subparagraph{Case $f=x_{2}$}

\begin{eqnarray*}
B(x_{1}x_{2}\otimes 1_{H};G,x_{2}) &=&-B(gx_{1}x_{2}\otimes
g;G,x_{2})-\lambda B(gx_{1}x_{2}\otimes g;GX_{1}X_{2},x_{2}) \\
&&+B(x_{1}x_{2}\otimes x_{2};GX_{2},x_{2}).
\end{eqnarray*}%
By applying $\left( \ref{eq.10}\right) $ this rewrites as%
\begin{equation*}
-\lambda B(x_{1}x_{2}\otimes 1_{H};GX_{1}X_{2},x_{2})-B(x_{1}x_{2}\otimes
x_{2};GX_{2},x_{2})=0.
\end{equation*}%
In view of the form of the elements we get%
\begin{equation*}
\lambda B(x_{1}x_{2}\otimes 1_{H};GX_{2},x_{1}x_{2})=0
\end{equation*}%
which is $\left( \ref{X1,x1x2,X1F61,x2}\right) .$

\subparagraph{Case $f=gx_{1}x_{2}$}

\begin{eqnarray*}
B(x_{1}x_{2}\otimes 1_{H};G,gx_{1}x_{2}) &=&-B(gx_{1}x_{2}\otimes
g;G,gx_{1}x_{2})-\lambda B(gx_{1}x_{2}\otimes g;GX_{1}X_{2},gx_{1}x_{2}) \\
&&+B(x_{1}x_{2}\otimes x_{2};GX_{2},gx_{1}x_{2}).
\end{eqnarray*}%
By applying $\left( \ref{eq.10}\right) $ this rewrites as%
\begin{eqnarray*}
&&2B(x_{1}x_{2}\otimes 1_{H};G,gx_{1}x_{2})+\lambda B(x_{1}x_{2}\otimes
1_{H};GX_{1}X_{2},gx_{1}x_{2}) \\
&&-B(x_{1}x_{2}\otimes x_{2};GX_{2},gx_{1}x_{2})=0.
\end{eqnarray*}%
In view of the form of the elements we get%
\begin{equation}
B(x_{1}x_{2}\otimes 1_{H};G,gx_{1}x_{2})-B(x_{1}x_{2}\otimes
1_{H};GX_{2},gx_{1})=0.  \label{X2,x1x2,X2F71,gx1x2}
\end{equation}

\paragraph{\textbf{Equality} $\left( \protect\ref{X2F81}\right) $}

rewrites as

\begin{eqnarray*}
B(x_{1}x_{2}\otimes 1_{H};GX_{1},f) &=&B(gx_{1}x_{2}\otimes g;GX_{1},f) \\
&&+B(x_{1}x_{2}\otimes x_{2};GX_{1}X_{2},f).
\end{eqnarray*}

\subparagraph{Case $f=1_{H}$}

\begin{eqnarray*}
B(x_{1}x_{2}\otimes 1_{H};GX_{1},1_{H}) &=&B(gx_{1}x_{2}\otimes
g;GX_{1},1_{H}) \\
&&+B(x_{1}x_{2}\otimes x_{2};GX_{1}X_{2},1_{H}).
\end{eqnarray*}%
By applying $\left( \ref{eq.10}\right) $ this rewrites as%
\begin{equation*}
B(x_{1}x_{2}\otimes x_{2};GX_{1}X_{2},1_{H})=0
\end{equation*}%
which is already known.

\subparagraph{Case $f=x_{1}x_{2}$}

\begin{eqnarray*}
B(x_{1}x_{2}\otimes 1_{H};GX_{1},x_{1}x_{2}) &=&B(gx_{1}x_{2}\otimes
g;GX_{1},x_{1}x_{2}) \\
&&+B(x_{1}x_{2}\otimes x_{2};GX_{1}X_{2},x_{1}x_{2}).
\end{eqnarray*}%
By applying $\left( \ref{eq.10}\right) $ this rewrites as%
\begin{equation*}
B(x_{1}x_{2}\otimes x_{2};GX_{1}X_{2},x_{1}x_{2})=0
\end{equation*}%
which is already known.

\subparagraph{Case $f=gx_{1}$}

\begin{eqnarray*}
B(x_{1}x_{2}\otimes 1_{H};GX_{1},gx_{1}) &=&B(gx_{1}x_{2}\otimes
g;GX_{1},gx_{1})+ \\
&&+B(x_{1}x_{2}\otimes x_{2};GX_{1}X_{2},gx_{1}).
\end{eqnarray*}%
By applying $\left( \ref{eq.10}\right) $ this rewrites as%
\begin{equation*}
2B(x_{1}x_{2}\otimes 1_{H};GX_{1},gx_{1})-B(x_{1}x_{2}\otimes
x_{2};GX_{1}X_{2},gx_{1})=0
\end{equation*}%
In view of the form of the elements we get%
\begin{equation*}
B(x_{1}x_{2}\otimes 1_{H};GX_{1},gx_{1})=0
\end{equation*}%
which is $\left( \ref{X2,x1x2,X2F61,gx1x2}\right) .$

\subparagraph{Case $f=gx_{2}$}

\begin{eqnarray*}
B(x_{1}x_{2}\otimes 1_{H};GX_{1},gx_{2}) &=&B(gx_{1}x_{2}\otimes
g;GX_{1},gx_{2})+ \\
&&+B(x_{1}x_{2}\otimes x_{2};GX_{1}X_{2},gx_{2}).
\end{eqnarray*}%
By applying $\left( \ref{eq.10}\right) $ this rewrites as%
\begin{equation*}
2B(x_{1}x_{2}\otimes 1_{H};GX_{1},gx_{2})-B(x_{1}x_{2}\otimes
x_{2};GX_{1}X_{2},gx_{2})=0.
\end{equation*}%
In view of the form of the elements we get%
\begin{equation*}
B(x_{1}x_{2}\otimes 1_{H};G,gx_{1}x_{2})-B(x_{1}x_{2}\otimes
1_{H};GX_{2},gx_{1})=0
\end{equation*}%
which is $\left( \ref{X2,x1x2,X2F71,gx1x2}\right) .$

\subsubsection{Case $gx_{1}\otimes 1_{H}$}

\paragraph{\textbf{Equality} $\left( \protect\ref{X2F11}\right) $}

rewrites as%
\begin{eqnarray*}
&&\beta _{2}B(gx_{1}\otimes 1_{H};X_{2},f) \\
&=&-\gamma _{2}B(x_{1}\otimes g;G,f)-\lambda B(x_{1}\otimes g;X_{1},f)-\beta
_{2}B(x_{1}\otimes g;X_{2},f)+ \\
&&+B(gx_{1}x_{2}\otimes g;1_{A},f)+B(gx_{1}\otimes x_{2};1_{A},f).
\end{eqnarray*}

\subparagraph{Case $f=1_{H}$}

\begin{eqnarray*}
&&\beta _{2}B(gx_{1}\otimes 1_{H};X_{2},1_{H}) \\
&=&-\gamma _{2}B(x_{1}\otimes g;G,1_{H})-\lambda B(x_{1}\otimes
g;X_{1},1_{H})-\beta _{2}B(x_{1}\otimes g;X_{2},1_{H})+ \\
&&+B(gx_{1}x_{2}\otimes g;1_{A},1_{H})+B(gx_{1}\otimes x_{2};1_{A},1_{H}).
\end{eqnarray*}%
By applying $\left( \ref{eq.10}\right) $ this rewrites as
\begin{gather*}
2\beta _{2}B(gx_{1}\otimes 1_{H};X_{2},1_{H})+\gamma _{2}B(gx_{1}\otimes
1_{H};G,1_{H})+\lambda B(gx_{1}\otimes 1_{H};X_{1},1_{H})+ \\
-B(x_{1}x_{2}\otimes 1_{H};1_{A},1_{H})-B(gx_{1}\otimes x_{2};1_{A},1_{H})=0.
\end{gather*}%
In view of the form of the elements we get%
\begin{gather}
2\beta _{2}B(x_{1}x_{2}\otimes 1_{H};X_{2},x_{2})+\gamma _{2}\left[
B(x_{1}x_{2}\otimes 1_{H};G,x_{2})-B(x_{1}x_{2}\otimes 1_{H};GX_{2},1_{H})%
\right] +  \label{X2,gx1,X2F11,1H} \\
+\lambda \left[ -1+B(x_{1}x_{2}\otimes
1_{H};1_{A},x_{1}x_{2})+B(x_{1}x_{2}\otimes 1_{H};X_{2},x_{1})\right] =0.
\notag
\end{gather}

\subparagraph{Case $f=x_{1}x_{2}$}

\begin{eqnarray*}
&&\beta _{2}B(gx_{1}\otimes 1_{H};X_{2},x_{1}x_{2}) \\
&=&-\gamma _{2}B(x_{1}\otimes g;G,x_{1}x_{2})-\lambda B(x_{1}\otimes
g;X_{1},x_{1}x_{2})-\beta _{2}B(x_{1}\otimes g;X_{2},x_{1}x_{2})+ \\
&&+B(gx_{1}x_{2}\otimes g;1_{A},x_{1}x_{2})+B(gx_{1}\otimes
x_{2};1_{A},x_{1}x_{2}).
\end{eqnarray*}%
By applying $\left( \ref{eq.10}\right) $ this rewrites as
\begin{gather*}
2\beta _{2}B(gx_{1}\otimes 1_{H};X_{2},x_{1}x_{2})+\gamma
_{2}B(gx_{1}\otimes 1_{H};G,x_{1}x_{2})+\lambda B(gx_{1}\otimes
1_{H};X_{1},x_{1}x_{2})+ \\
-B(x_{1}x_{2}\otimes 1_{H};1_{A},x_{1}x_{2})-B(gx_{1}\otimes
x_{2};1_{A},x_{1}x_{2})=0.
\end{gather*}%
In view of the form of the elements we get%
\begin{equation*}
\gamma _{2}B(x_{1}x_{2}\otimes 1_{H};GX_{2},x_{1}x_{2})=0
\end{equation*}%
which is $\left( \ref{G,x1x2, GF6,x2}\right) .$

\subparagraph{Case $f=gx_{1}$}

\begin{eqnarray*}
&&\beta _{2}B(gx_{1}\otimes 1_{H};X_{2},gx_{1}) \\
&=&-\gamma _{2}B(x_{1}\otimes g;G,gx_{1})-\lambda B(x_{1}\otimes
g;X_{1},gx_{1})-\beta _{2}B(x_{1}\otimes g;X_{2},gx_{1})+ \\
&&+B(gx_{1}x_{2}\otimes g;1_{A},gx_{1})+B(gx_{1}\otimes x_{2};1_{A},gx_{1}).
\end{eqnarray*}%
By applying $\left( \ref{eq.10}\right) $ this rewrites as
\begin{eqnarray*}
&&-\gamma _{2}B(gx_{1}\otimes 1_{H};G,gx_{1})-\lambda B(gx_{1}\otimes
1_{H};X_{1},gx_{1}) \\
&&+B(x_{1}x_{2}\otimes 1_{H};1_{A},gx_{1})-B(gx_{1}\otimes
x_{2};1_{A},gx_{1})=0.
\end{eqnarray*}%
In view of the form of the elements we get%
\begin{equation*}
\gamma _{2}\left[ B(x_{1}x_{2}\otimes
1_{H};G,gx_{1}x_{2})-B(x_{1}x_{2}\otimes 1_{H};GX_{2},gx_{1})\right] =0
\end{equation*}%
which follows from $\left( \ref{X2,x1x2,X2F71,gx1x2}\right) .$

\subparagraph{Case $f=gx_{2}$}

\begin{eqnarray*}
\beta _{2}B(gx_{1}\otimes 1_{H};X_{2},gx_{2}) &=&-\gamma _{2}B(x_{1}\otimes
g;G,gx_{2})-\lambda B(x_{1}\otimes g;X_{1},gx_{2})-\beta _{2}B(x_{1}\otimes
g;X_{2},gx_{2})+ \\
&&+B(gx_{1}x_{2}\otimes g;1_{A},gx_{2})+B(gx_{1}\otimes x_{2};1_{A},gx_{2}).
\end{eqnarray*}%
By applying $\left( \ref{eq.10}\right) $ this rewrites as
\begin{eqnarray*}
&&-\gamma _{2}B(gx_{1}\otimes 1_{H};G,gx_{2})-\lambda B(gx_{1}\otimes
1_{H};X_{1},gx_{2}) \\
&&+B(x_{1}x_{2}\otimes 1_{H};1_{A},gx_{2})-B(gx_{1}\otimes
x_{2};1_{A},gx_{2}).
\end{eqnarray*}%
In view of the form of the elements we get%
\begin{eqnarray*}
&&\gamma _{2}B(x_{1}x_{2}\otimes 1_{H};GX_{2},gx_{2})+\lambda
B(x_{1}x_{2}\otimes 1_{H};X_{2},gx_{1}x_{2}) \\
&&+2B(x_{1}x_{2}\otimes 1_{H};1_{A},gx_{2})+2B(x_{1}x_{2}\otimes
1_{H};X_{2},g).
\end{eqnarray*}%
In view of $\left( \ref{X1,x1x2,X1F71,gx1x2}\right) ,\left( \ref{G,x1x2,
GF7,gx1x2}\right) $%
\begin{equation}
B(x_{1}x_{2}\otimes 1_{H};1_{A},gx_{2})+B(x_{1}x_{2}\otimes 1_{H};X_{2},g)=0.
\label{X2,gx1,X2F11,gx2}
\end{equation}

\paragraph{\textbf{Equality} $\left( \protect\ref{X2F21}\right) $}

rewrites as%
\begin{eqnarray*}
\beta _{2}B(gx_{1}\otimes 1_{H};GX_{2},f) &=&+\beta _{2}B(x_{1}\otimes
g;GX_{2},f)+\lambda B(x_{1}\otimes g;GX_{1},f) \\
&&+B(gx_{1}x_{2}\otimes g;G,f)+B(gx_{1}\otimes x_{2};G,f).
\end{eqnarray*}

\subparagraph{Case $f=g$}

\begin{eqnarray*}
\beta _{2}B(gx_{1}\otimes 1_{H};GX_{2},g) &=&+\beta _{2}B(x_{1}\otimes
g;GX_{2},g)+\lambda B(x_{1}\otimes g;GX_{1},g) \\
&&+B(gx_{1}x_{2}\otimes g;G,g)+B(gx_{1}\otimes x_{2};G,g).
\end{eqnarray*}%
By applying $\left( \ref{eq.10}\right) $ this rewrites as%
\begin{eqnarray*}
&&-\lambda B(gx_{1}\otimes 1_{H};GX_{1},g) \\
&&-B(x_{1}x_{2}\otimes 1_{H};G,g)-B(gx_{1}\otimes x_{2};G,g)=0.
\end{eqnarray*}%
In view of the form of the elements we get%
\begin{equation}
\lambda \left[ B(x_{1}x_{2}\otimes 1_{H};G,gx_{1}x_{2})-B(x_{1}x_{2}\otimes
1_{H};GX_{2},gx_{1})\right] =0.  \label{X2,gx1,X2F21,g}
\end{equation}

\subparagraph{Case $f=x_{1}$}

\begin{eqnarray*}
\beta _{2}B(gx_{1}\otimes 1_{H};GX_{2},x_{1}) &=&+\beta _{2}B(x_{1}\otimes
g;GX_{2},x_{1})+\lambda B(x_{1}\otimes g;GX_{1},x_{1}) \\
&&+B(gx_{1}x_{2}\otimes g;G,x_{1})+B(gx_{1}\otimes x_{2};G,x_{1}).
\end{eqnarray*}%
By applying $\left( \ref{eq.10}\right) $ this rewrites as%
\begin{eqnarray*}
&&2\beta _{2}B(gx_{1}\otimes 1_{H};GX_{2},x_{1})+\lambda B(gx_{1}\otimes
1_{H};GX_{1},x_{1}) \\
&&+B(x_{1}x_{2}\otimes 1_{H};G,x_{1})-B(gx_{1}\otimes x_{2};G,x_{1})=0.
\end{eqnarray*}%
In view of the form of the elements we get%
\begin{equation*}
\beta _{2}B(x_{1}x_{2}\otimes 1_{H};GX_{2},x_{1}x_{2})=0
\end{equation*}%
which is $\left( \ref{X2,x1x2,X2F21,x1x2}\right) .$

\subparagraph{Case $f=x_{2}$}

\begin{eqnarray*}
\beta _{2}B(gx_{1}\otimes 1_{H};GX_{2},x_{2}) &=&+\beta _{2}B(x_{1}\otimes
g;GX_{2},x_{2})+\lambda B(x_{1}\otimes g;GX_{1},x_{2}) \\
&&+B(gx_{1}x_{2}\otimes g;G,x_{2})+B(gx_{1}\otimes x_{2};G,x_{2}).
\end{eqnarray*}%
By applying $\left( \ref{eq.10}\right) $ this rewrites as%
\begin{eqnarray*}
&&2\beta _{2}B(gx_{1}\otimes 1_{H};GX_{2},x_{2})+\lambda B(gx_{1}\otimes
1_{H};GX_{1},x_{2}) \\
&&+B(x_{1}x_{2}\otimes 1_{H};G,x_{2})-B(gx_{1}\otimes x_{2};G,x_{2})=0.
\end{eqnarray*}%
In view of the form of the elements we get%
\begin{equation*}
\lambda B(x_{1}x_{2}\otimes 1_{H};GX_{2},x_{1}x_{2})=0
\end{equation*}%
which is $\left( \ref{X1,x1x2,X1F61,x2}\right) .$

\subparagraph{Case $f=gx_{1}x_{2}$}

\begin{eqnarray*}
&&\beta _{2}B(gx_{1}\otimes 1_{H};GX_{2},gx_{1}x_{2}) \\
&=&+\beta _{2}B(x_{1}\otimes g;GX_{2},gx_{1}x_{2})+\lambda B(x_{1}\otimes
g;GX_{1},gx_{1}x_{2})+ \\
&&+B(gx_{1}x_{2}\otimes g;G,gx_{1}x_{2})+B(gx_{1}\otimes
x_{2};G,gx_{1}x_{2}).
\end{eqnarray*}%
By applying $\left( \ref{eq.10}\right) $ this rewrites as%
\begin{eqnarray*}
&&-\lambda B(gx_{1}\otimes 1_{H};GX_{1},gx_{1}x_{2}) \\
&&-B(x_{1}x_{2}\otimes 1_{H};G,gx_{1}x_{2})-B(gx_{1}\otimes
x_{2};G,gx_{1}x_{2})=0.
\end{eqnarray*}%
In view of the form of the elements we get

\begin{equation*}
B(x_{1}x_{2}\otimes 1_{H};G,gx_{1}x_{2})-B(x_{1}x_{2}\otimes
1_{H};GX_{2},gx_{1})=0
\end{equation*}%
which is $\left( \ref{X2,x1x2,X2F71,gx1x2}\right) .$

\paragraph{\textbf{Equality} $\left( \protect\ref{X2F31}\right) $}

rewrites as%
\begin{eqnarray*}
\beta _{2}B(gx_{1}\otimes 1_{H};X_{1}X_{2},f) &=&+\beta _{2}B(x_{1}\otimes
g;X_{1}X_{2},f)-\gamma _{2}B(x_{1}\otimes g;GX_{1},f) \\
&&+B(gx_{1}x_{2}\otimes g;X_{1},f)+B(gx_{1}\otimes x_{2};X_{1},f).
\end{eqnarray*}%
Case $f=g$

\begin{eqnarray*}
\beta _{2}B(gx_{1}\otimes 1_{H};X_{1}X_{2},g) &=&+\beta _{2}B(x_{1}\otimes
g;X_{1}X_{2},g)-\gamma _{2}B(x_{1}\otimes g;GX_{1},g) \\
&&+B(gx_{1}x_{2}\otimes g;X_{1},g)+B(gx_{1}\otimes x_{2};X_{1},g).
\end{eqnarray*}%
By applying $\left( \ref{eq.10}\right) $ this rewrites as%
\begin{eqnarray*}
&&+\gamma _{2}B(gx_{1}\otimes 1_{H};GX_{1},g) \\
&&-B(x_{1}x_{2}\otimes 1_{H};X_{1},g)-B(gx_{1}\otimes x_{2};X_{1},g)=0.
\end{eqnarray*}%
In view of the form of the elements we get%
\begin{equation*}
\gamma _{2}\left[ B(x_{1}x_{2}\otimes
1_{H};G,gx_{1}x_{2})-B(x_{1}x_{2}\otimes 1_{H};GX_{2},gx_{1})\right] =0
\end{equation*}%
which follows from $\left( \ref{X2,x1x2,X2F71,gx1x2}\right) .$

\subparagraph{Case $f=x_{1}$}

\begin{eqnarray*}
\beta _{2}B(gx_{1}\otimes 1_{H};X_{1}X_{2},x_{1}) &=&+\beta
_{2}B(x_{1}\otimes g;X_{1}X_{2},x_{1})-\gamma _{2}B(x_{1}\otimes
g;GX_{1},x_{1}) \\
&&+B(gx_{1}x_{2}\otimes g;X_{1},x_{1})+B(gx_{1}\otimes x_{2};X_{1},x_{1}).
\end{eqnarray*}

By applying $\left( \ref{eq.10}\right) $ this rewrites as%
\begin{eqnarray*}
&&2\beta _{2}B(gx_{1}\otimes 1_{H};X_{1}X_{2},x_{1})-\gamma
_{2}B(gx_{1}\otimes 1_{H};GX_{1},x_{1}) \\
&&+B(x_{1}x_{2}\otimes 1_{H};X_{1},x_{1})-B(gx_{1}\otimes
x_{2};X_{1},x_{1})=0.
\end{eqnarray*}%
In view of the form of the elements we get

\begin{equation*}
B(x_{1}x_{2}\otimes 1_{H};X_{1},x_{1})-B(x_{1}x_{2}\otimes
1_{H};X_{1},x_{1})=0
\end{equation*}%
which is trivial.

\subparagraph{Case $f=x_{2}$}

\begin{eqnarray*}
\beta _{2}B(gx_{1}\otimes 1_{H};X_{1}X_{2},x_{2}) &=&+\beta
_{2}B(x_{1}\otimes g;X_{1}X_{2},x_{2})-\gamma _{2}B(x_{1}\otimes
g;GX_{1},x_{2}) \\
&&+B(gx_{1}x_{2}\otimes g;X_{1},x_{2})+B(gx_{1}\otimes x_{2};X_{1},x_{2}).
\end{eqnarray*}

By applying $\left( \ref{eq.10}\right) $ this rewrites as%
\begin{eqnarray*}
&&2\beta _{2}B(gx_{1}\otimes 1_{H};X_{1}X_{2},x_{2})-\gamma
_{2}B(gx_{1}\otimes 1_{H};GX_{1},x_{2}) \\
&&+B(x_{1}x_{2}\otimes 1_{H};X_{1},x_{2})-B(gx_{1}\otimes
x_{2};X_{1},x_{2})=0.
\end{eqnarray*}%
In view of the form of the elements we get%
\begin{equation*}
\gamma _{2}B(x_{1}x_{2}\otimes 1_{H};GX_{2},x_{1}x_{2})=0
\end{equation*}%
which is $\left( \ref{G,x1x2, GF6,x2}\right) .$

\subparagraph{Case $f=gx_{1}x_{2}$}

\begin{eqnarray*}
\beta _{2}B(gx_{1}\otimes 1_{H};X_{1}X_{2},gx_{1}x_{2}) &=&+\beta
_{2}B(x_{1}\otimes g;X_{1}X_{2},gx_{1}x_{2})-\gamma _{2}B(x_{1}\otimes
g;GX_{1},gx_{1}x_{2}) \\
&&+B(gx_{1}x_{2}\otimes g;X_{1},gx_{1}x_{2})+B(gx_{1}\otimes
x_{2};X_{1},gx_{1}x_{2}).
\end{eqnarray*}%
By applying $\left( \ref{eq.10}\right) $ this rewrites as%
\begin{eqnarray*}
&&\gamma _{2}B(gx_{1}\otimes 1_{H};GX_{1},gx_{1}x_{2}) \\
&&-B(x_{1}x_{2}\otimes 1_{H};X_{1},gx_{1}x_{2})-B(gx_{1}\otimes
x_{2};X_{1},gx_{1}x_{2})=0.
\end{eqnarray*}%
In view of the form of the elements we get%
\begin{equation*}
B(x_{1}x_{2}\otimes 1_{H};X_{1},gx_{1}x_{2})-B(x_{1}x_{2}\otimes
1_{H};X_{1},gx_{1}x_{2})=0
\end{equation*}%
which is trivial.

\paragraph{\textbf{Equality} $\left( \protect\ref{X2F41}\right) $}

rewrites as%
\begin{eqnarray*}
B(gx_{1}\otimes 1_{H};1_{A},f) &=&-B(x_{1}\otimes g;1_{A},f)-\gamma
_{2}B(x_{1}\otimes g;GX_{2},f)-\lambda B(x_{1}\otimes g;X_{1}X_{2},f) \\
&&+B(gx_{1}x_{2}\otimes g;X_{2},f)+B(gx_{1}\otimes x_{2};X_{2},f).
\end{eqnarray*}

\subparagraph{Case $f=g$}

\begin{eqnarray*}
B(gx_{1}\otimes 1_{H};1_{A},g) &=&-B(x_{1}\otimes g;1_{A},g)-\gamma
_{2}B(x_{1}\otimes g;GX_{2},g)-\lambda B(x_{1}\otimes g;X_{1}X_{2},g) \\
&&+B(gx_{1}x_{2}\otimes g;X_{2},g)+B(gx_{1}\otimes x_{2};X_{2},g).
\end{eqnarray*}

By applying $\left( \ref{eq.10}\right) $ this rewrites as%
\begin{eqnarray*}
&&2B(gx_{1}\otimes 1_{H};1_{A},g)+\gamma _{2}B(gx_{1}\otimes
1_{H};GX_{2},g)+\lambda B(gx_{1}\otimes 1_{H};X_{1}X_{2},g) \\
&&-B(x_{1}x_{2}\otimes 1_{H};X_{2},g)-B(gx_{1}\otimes x_{2};X_{2},g)=0.
\end{eqnarray*}
In view of the form of the elements we get

\begin{gather*}
2\left[ B(x_{1}x_{2}\otimes 1_{H};1_{A},gx_{2})+B(x_{1}x_{2}\otimes
1_{H};X_{2},g)\right] \\
+\gamma _{2}B(x_{1}x_{2}\otimes 1_{H};GX_{2},gx_{2})+\lambda
B(x_{1}x_{2}\otimes 1_{H};X_{2},gx_{1}x_{2})=0.
\end{gather*}%
This follows from $\left( \ref{X2,gx1,X2F11,gx2}\right) ,\left( \ref%
{X1,x1x2,X1F71,gx1x2}\right) $ and $\left( \ref{G,x1x2, GF7,gx1x2}\right) $

\subparagraph{Case $f=x_{1}$}

\begin{eqnarray*}
B(gx_{1}\otimes 1_{H};1_{A},x_{1}) &=&-B(x_{1}\otimes g;1_{A},x_{1})-\gamma
_{2}B(x_{1}\otimes g;GX_{2},x_{1})-\lambda B(x_{1}\otimes g;X_{1}X_{2},x_{1})
\\
&&+B(gx_{1}x_{2}\otimes g;X_{2},x_{1})+B(gx_{1}\otimes x_{2};X_{2},x_{1}).
\end{eqnarray*}%
By applying $\left( \ref{eq.10}\right) $ this rewrites as
\begin{eqnarray*}
&&-\gamma _{2}B(gx_{1}\otimes 1_{H};GX_{2},x_{1})-\lambda B(gx_{1}\otimes
1_{H};X_{1}X_{2},x_{1}) \\
&&+B(x_{1}x_{2}\otimes 1_{H};X_{2},x_{1})-B(gx_{1}\otimes
x_{2};X_{2},x_{1})=0.
\end{eqnarray*}%
In view of the form of the elements we get%
\begin{equation*}
\gamma _{2}B(x_{1}x_{2}\otimes 1_{H};GX_{2},x_{1}x_{2})=0
\end{equation*}%
which is $\left( \ref{G,x1x2, GF6,x2}\right) .$

\subparagraph{Case $f=x_{2}$}

\begin{eqnarray*}
B(gx_{1}\otimes 1_{H};1_{A},x_{2}) &=&-B(x_{1}\otimes g;1_{A},x_{2})-\gamma
_{2}B(x_{1}\otimes g;GX_{2},x_{2})-\lambda B(x_{1}\otimes g;X_{1}X_{2},x_{2})
\\
&&+B(gx_{1}x_{2}\otimes g;X_{2},x_{2})+B(gx_{1}\otimes x_{2};X_{2},x_{2}).
\end{eqnarray*}%
By applying $\left( \ref{eq.10}\right) $ this rewrites as
\begin{eqnarray*}
&&-\gamma _{2}B(gx_{1}\otimes 1_{H};GX_{2},x_{2})-\lambda B(gx_{1}\otimes
1_{H};X_{1}X_{2},x_{2}) \\
&&+B(x_{1}x_{2}\otimes 1_{H};X_{2},x_{2})-B(gx_{1}\otimes
x_{2};X_{2},x_{2})=0.
\end{eqnarray*}%
In view of the form of the elements we get%
\begin{equation*}
B(x_{1}x_{2}\otimes 1_{H};X_{2},x_{2})-B(x_{1}x_{2}\otimes
1_{H};X_{2},x_{2})=0
\end{equation*}%
which is trivial.

\subparagraph{Case $f=gx_{1}x_{2}$}

\begin{gather*}
B(gx_{1}\otimes 1_{H};1_{A},gx_{1}x_{2})=-B(x_{1}\otimes
g;1_{A},gx_{1}x_{2})+ \\
-\gamma _{2}B(x_{1}\otimes g;GX_{2},gx_{1}x_{2})-\lambda B(x_{1}\otimes
g;X_{1}X_{2},gx_{1}x_{2}) \\
+B(gx_{1}x_{2}\otimes g;X_{2},gx_{1}x_{2})+B(gx_{1}\otimes
x_{2};X_{2},gx_{1}x_{2}).
\end{gather*}%
By applying $\left( \ref{eq.10}\right) $ this rewrites as
\begin{gather*}
2B(gx_{1}\otimes 1_{H};1_{A},gx_{1}x_{2})+\gamma _{2}B(gx_{1}\otimes
1_{H};GX_{2},gx_{1}x_{2})+ \\
+\lambda B(gx_{1}\otimes 1_{H};X_{1}X_{2},gx_{1}x_{2})-B(x_{1}x_{2}\otimes
1_{H};X_{2},gx_{1}x_{2})+ \\
-B(gx_{1}\otimes x_{2};X_{2},gx_{1}x_{2})=0.
\end{gather*}%
In view of the form of the elements we get%
\begin{equation*}
B(x_{1}x_{2}\otimes 1_{H};X_{2},gx_{1}x_{2})-B(x_{1}x_{2}\otimes
1_{H};X_{2},gx_{1}x_{2})=0
\end{equation*}%
which is trivial.

\paragraph{\textbf{Equality} $\left( \protect\ref{X2F51}\right) $}

rewrites as%
\begin{eqnarray*}
B(gx_{1}\otimes 1_{H};X_{1},f) &=&+B(x_{1}\otimes g;X_{1},f)-\gamma
_{2}B(x_{1}\otimes g;GX_{1}X_{2},f) \\
&&+B(gx_{1}x_{2}\otimes g;X_{1}X_{2},f)+B(gx_{1}\otimes x_{2};X_{1}X_{2},f)
\end{eqnarray*}

\subparagraph{Case $f=1_{H}$}

\begin{eqnarray*}
B(gx_{1}\otimes 1_{H};X_{1},1_{H}) &=&+B(x_{1}\otimes g;X_{1},1_{H})-\gamma
_{2}B(x_{1}\otimes g;GX_{1}X_{2},1_{H}) \\
&&+B(gx_{1}x_{2}\otimes g;X_{1}X_{2},1_{H})+B(gx_{1}\otimes
x_{2};X_{1}X_{2},1_{H})
\end{eqnarray*}%
By applying $\left( \ref{eq.10}\right) $ this rewrites as

\begin{eqnarray*}
&&+\gamma _{2}B(gx_{1}\otimes 1_{H};GX_{1}X_{2},1_{H}) \\
&&-B(x_{1}x_{2}\otimes 1_{H};X_{1}X_{2},1_{H})-B(gx_{1}\otimes
x_{2};X_{1}X_{2},1_{H})=0
\end{eqnarray*}%
In view of the form of the elements we get%
\begin{equation*}
\gamma _{2}B(x_{1}x_{2}\otimes 1_{H};GX_{2},x_{1}x_{2})=0
\end{equation*}%
which is $\left( \ref{G,x1x2, GF6,x2}\right) .$

\subparagraph{Case $f=gx_{1}$}

\begin{eqnarray*}
B(gx_{1}\otimes 1_{H};X_{1},gx_{1}) &=&+B(x_{1}\otimes
g;X_{1},gx_{1})-\gamma _{2}B(x_{1}\otimes g;GX_{1}X_{2},gx_{1}) \\
&&+B(gx_{1}x_{2}\otimes g;X_{1}X_{2},gx_{1})+B(gx_{1}\otimes
x_{2};X_{1}X_{2},gx_{1})
\end{eqnarray*}%
By applying $\left( \ref{eq.10}\right) $ this rewrites as

\begin{eqnarray*}
&&2B(gx_{1}\otimes 1_{H};X_{1},gx_{1})-\gamma _{2}B(gx_{1}\otimes
1_{H};GX_{1}X_{2},gx_{1}) \\
&&+B(x_{1}x_{2}\otimes 1_{H};X_{1}X_{2},gx_{1})-B(gx_{1}\otimes
x_{2};X_{1}X_{2},gx_{1})=0
\end{eqnarray*}%
In view of the form of the elements we get%
\begin{equation*}
B(x_{1}x_{2}\otimes 1_{H};X_{1},gx_{1}x_{2})-B(x_{1}x_{2}\otimes
1_{H};X_{1},gx_{1}x_{2})=0
\end{equation*}%
which is trivial.

\subparagraph{Case $f=gx_{2}$}

\begin{eqnarray*}
B(gx_{1}\otimes 1_{H};X_{1},gx_{2}) &=&+B(x_{1}\otimes
g;X_{1},gx_{2})-\gamma _{2}B(x_{1}\otimes g;GX_{1}X_{2},gx_{2}) \\
&&+B(gx_{1}x_{2}\otimes g;X_{1}X_{2},gx_{2})+B(gx_{1}\otimes
x_{2};X_{1}X_{2},gx_{2}).
\end{eqnarray*}%
By applying $\left( \ref{eq.10}\right) $ this rewrites as

\begin{eqnarray*}
&&2B(gx_{1}\otimes 1_{H};X_{1},gx_{2})-\gamma _{2}B(gx_{1}\otimes
1_{H};GX_{1}X_{2},gx_{2}) \\
&&+B(x_{1}x_{2}\otimes 1_{H};X_{1}X_{2},gx_{2})-B(gx_{1}\otimes
x_{2};X_{1}X_{2},gx_{2})=0.
\end{eqnarray*}%
In view of the form of the elements we get%
\begin{equation*}
B(x_{1}x_{2}\otimes 1_{H};X_{2},gx_{1}x_{2})-B(x_{1}x_{2}\otimes
1_{H};X_{2},gx_{1}x_{2})=0
\end{equation*}%
which is trivial.

\paragraph{\textbf{Equality} $\left( \protect\ref{X2F61}\right) $}

rewrites as%
\begin{eqnarray*}
\beta _{2}B(gx_{1}\otimes 1_{H};GX_{1}X_{2},f) &=&-\beta _{2}B(x_{1}\otimes
g;GX_{1}X_{2},f) \\
&&+B(gx_{1}x_{2}\otimes g;GX_{1},f)+B(gx_{1}\otimes x_{2};GX_{1},f).
\end{eqnarray*}

\subparagraph{Case $f=1_{H}$}

\begin{eqnarray*}
\beta _{2}B(gx_{1}\otimes 1_{H};GX_{1}X_{2},1_{H}) &=&-\beta
_{2}B(x_{1}\otimes g;GX_{1}X_{2},1_{H}) \\
&&+B(gx_{1}x_{2}\otimes g;GX_{1},1_{H})+B(gx_{1}\otimes x_{2};GX_{1},1_{H}).
\end{eqnarray*}%
By applying $\left( \ref{eq.10}\right) $ this rewrites as
\begin{eqnarray*}
&&2\beta _{2}B(gx_{1}\otimes 1_{H};GX_{1}X_{2},1_{H}) \\
&&-B(x_{1}x_{2}\otimes 1_{H};GX_{1},1_{H})-B(gx_{1}\otimes
x_{2};GX_{1},1_{H})=0.
\end{eqnarray*}%
In view of the form of the elements we get%
\begin{equation*}
\beta _{2}B(x_{1}x_{2}\otimes 1_{H};GX_{2},x_{1}x_{2})=0
\end{equation*}%
which is $\left( \ref{X2,x1x2,X2F21,x1x2}\right) .$

\subparagraph{Case $f=x_{1}x_{2}$}

\begin{eqnarray*}
\beta _{2}B(gx_{1}\otimes 1_{H};GX_{1}X_{2},x_{1}x_{2}) &=&-\beta
_{2}B(x_{1}\otimes g;GX_{1}X_{2},x_{1}x_{2}) \\
&&+B(gx_{1}x_{2}\otimes g;GX_{1},x_{1}x_{2})+B(gx_{1}\otimes
x_{2};GX_{1},x_{1}x_{2}).
\end{eqnarray*}%
By applying $\left( \ref{eq.10}\right) $ this rewrites as
\begin{eqnarray*}
&&2\beta _{2}B(gx_{1}\otimes 1_{H};GX_{1}X_{2},x_{1}x_{2}) \\
&&-B(x_{1}x_{2}\otimes 1_{H};GX_{1},x_{1}x_{2})-B(gx_{1}\otimes
x_{2};GX_{1},x_{1}x_{2})=0.
\end{eqnarray*}%
In view of the form of the elements we get%
\begin{equation*}
-B(x_{1}x_{2}\otimes 1_{H};GX_{1},x_{1}x_{2})+B(x_{1}x_{2}\otimes
1_{H};GX_{1},x_{1}x_{2})=0
\end{equation*}

which is trivial.

\subparagraph{Case $f=gx_{1}$}

\begin{eqnarray*}
\beta _{2}B(gx_{1}\otimes 1_{H};GX_{1}X_{2},gx_{1}) &=&-\beta
_{2}B(x_{1}\otimes g;GX_{1}X_{2},gx_{1}) \\
&&+B(gx_{1}x_{2}\otimes g;GX_{1},gx_{1})+B(gx_{1}\otimes
x_{2};GX_{1},gx_{1}).
\end{eqnarray*}%
By applying $\left( \ref{eq.10}\right) $ this rewrites as
\begin{equation*}
-B(x_{1}x_{2}\otimes 1_{H};GX_{1},gx_{1})+B(gx_{1}\otimes
x_{2};GX_{1},gx_{1})=0.
\end{equation*}%
In view of the form of the elements we get%
\begin{equation*}
-B(x_{1}x_{2}\otimes 1_{H};GX_{1},gx_{1})+B(x_{1}x_{2}\otimes
1_{H};GX_{1},gx_{1})=0
\end{equation*}%
which is trivial.

\subparagraph{Case $f=gx_{2}$}

\begin{eqnarray*}
\beta _{2}B(gx_{1}\otimes 1_{H};GX_{1}X_{2},gx_{2}) &=&-\beta
_{2}B(x_{1}\otimes g;GX_{1}X_{2},gx_{2}) \\
&&+B(gx_{1}x_{2}\otimes g;GX_{1},gx_{2})+B(gx_{1}\otimes
x_{2};GX_{1},gx_{2}).
\end{eqnarray*}%
By applying $\left( \ref{eq.10}\right) $ this rewrites as
\begin{equation*}
-B(x_{1}x_{2}\otimes 1_{H};GX_{1},gx_{2})+B(gx_{1}\otimes
x_{2};GX_{1},gx_{2})=0.
\end{equation*}%
In view of the form of the elements we get%
\begin{equation*}
B(x_{1}x_{2}\otimes 1_{H};G,gx_{1}x_{2})-B(x_{1}x_{2}\otimes
1_{H};GX_{2},gx_{1})=0
\end{equation*}%
which is $\left( \ref{X2,x1x2,X2F71,gx1x2}\right) .$

\paragraph{\textbf{Equality} $\left( \protect\ref{X2F71}\right) $}

rewrites as%
\begin{eqnarray*}
B(gx_{1}\otimes 1_{H};G,f) &=&+B(x_{1}\otimes g;G,f)+\lambda B(x_{1}\otimes
g;GX_{1}X_{2},f) \\
&&+B(gx_{1}x_{2}\otimes g;GX_{2},f)+B(gx_{1}\otimes x_{2};GX_{2},f).
\end{eqnarray*}

\subparagraph{Case $f=1_{H}$}

\begin{eqnarray*}
B(gx_{1}\otimes 1_{H};G,1_{H}) &=&+B(x_{1}\otimes g;G,1_{H})+\lambda
B(x_{1}\otimes g;GX_{1}X_{2},1_{H}) \\
&&+B(gx_{1}x_{2}\otimes g;GX_{2},1_{H})+B(gx_{1}\otimes x_{2};GX_{2},1_{H}).
\end{eqnarray*}%
By applying $\left( \ref{eq.10}\right) $ this rewrites as
\begin{eqnarray*}
&&-\lambda B(gx_{1}\otimes 1_{H};GX_{1}X_{2},1_{H}) \\
&&-B(x_{1}x_{2}\otimes 1_{H};GX_{2},1_{H})-B(gx_{1}\otimes
x_{2};GX_{2},1_{H})=0.
\end{eqnarray*}%
In view of the form of the elements we get%
\begin{equation*}
\lambda B(x_{1}x_{2}\otimes 1_{H};GX_{2},x_{1}x_{2})=0
\end{equation*}%
which is $\left( \ref{X1,x1x2,X1F61,x2}\right) .$

\subparagraph{Case $f=x_{1}x_{2}$}

\begin{eqnarray*}
B(gx_{1}\otimes 1_{H};G,x_{1}x_{2}) &=&+B(x_{1}\otimes
g;G,x_{1}x_{2})+\lambda B(x_{1}\otimes g;GX_{1}X_{2},x_{1}x_{2}) \\
&&+B(gx_{1}x_{2}\otimes g;GX_{2},x_{1}x_{2})+B(gx_{1}\otimes
x_{2};GX_{2},x_{1}x_{2}).
\end{eqnarray*}

By applying $\left( \ref{eq.10}\right) $ this rewrites as
\begin{eqnarray*}
&&-\lambda B(gx_{1}\otimes 1_{H};GX_{1}X_{2},x_{1}x_{2}) \\
&&-B(x_{1}x_{2}\otimes 1_{H};GX_{2},x_{1}x_{2})-B(gx_{1}\otimes
x_{2};GX_{2},x_{1}x_{2})=0.
\end{eqnarray*}%
In view of the form of the elements we get%
\begin{equation*}
-B(x_{1}x_{2}\otimes 1_{H};GX_{2},x_{1}x_{2})+B(x_{1}x_{2}\otimes
1_{H};GX_{2},x_{1}x_{2})=0
\end{equation*}%
which is trivial.

\subparagraph{Case $f=gx_{1}$}

\begin{eqnarray*}
B(gx_{1}\otimes 1_{H};G,gx_{1}) &=&+B(x_{1}\otimes g;G,gx_{1})+\lambda
B(x_{1}\otimes g;GX_{1}X_{2},gx_{1}) \\
&&+B(gx_{1}x_{2}\otimes g;GX_{2},gx_{1})+B(gx_{1}\otimes
x_{2};GX_{2},gx_{1}).
\end{eqnarray*}%
By applying $\left( \ref{eq.10}\right) $ this rewrites as
\begin{eqnarray*}
&&2B(gx_{1}\otimes 1_{H};G,gx_{1})+\lambda B(gx_{1}\otimes
1_{H};GX_{1}X_{2},gx_{1}) \\
&&+B(x_{1}x_{2}\otimes 1_{H};GX_{2},gx_{1})-B(gx_{1}\otimes
x_{2};GX_{2},gx_{1})=0.
\end{eqnarray*}%
In view of the form of the elements we get%
\begin{equation*}
B(x_{1}x_{2}\otimes 1_{H};G,gx_{1}x_{2})-B(x_{1}x_{2}\otimes
1_{H};GX_{2},gx_{1})=0
\end{equation*}%
which is $\left( \ref{X2,x1x2,X2F71,gx1x2}\right) .$

\subparagraph{Case $f=gx_{2}$}

\begin{eqnarray*}
B(gx_{1}\otimes 1_{H};G,gx_{2}) &=&+B(x_{1}\otimes g;G,gx_{2})+\lambda
B(x_{1}\otimes g;GX_{1}X_{2},gx_{2}) \\
&&+B(gx_{1}x_{2}\otimes g;GX_{2},gx_{2})+B(gx_{1}\otimes
x_{2};GX_{2},gx_{2}).
\end{eqnarray*}%
By applying $\left( \ref{eq.10}\right) $ this rewrites as%
\begin{eqnarray*}
&&2B(gx_{1}\otimes 1_{H};G,gx_{2})+\lambda B(gx_{1}\otimes
1_{H};GX_{1}X_{2},gx_{2}) \\
&&+B(x_{1}x_{2}\otimes 1_{H};GX_{2},gx_{2})-B(gx_{1}\otimes
x_{2};GX_{2},gx_{2})=0.
\end{eqnarray*}
In view of the form of the elements we get%
\begin{equation*}
-B(x_{1}x_{2}\otimes 1_{H};GX_{2},gx_{2})+B(x_{1}x_{2}\otimes
1_{H};GX_{2},gx_{2})=0
\end{equation*}%
which is trivial.

\paragraph{\textbf{Equality} $\left( \protect\ref{X2F81}\right) $}

rewrites as%
\begin{eqnarray*}
B(gx_{1}\otimes 1_{H};GX_{1},f) &=&-B(x_{1}\otimes g;GX_{1},f) \\
&&+B(gx_{1}x_{2}\otimes g;GX_{1}X_{2},f)+B(gx_{1}\otimes
x_{2};GX_{1}X_{2},f).
\end{eqnarray*}

\subparagraph{Case $f=g$}

\begin{eqnarray*}
B(gx_{1}\otimes 1_{H};GX_{1},g) &=&-B(x_{1}\otimes g;GX_{1},g) \\
&&+B(gx_{1}x_{2}\otimes g;GX_{1}X_{2},g)+B(gx_{1}\otimes
x_{2};GX_{1}X_{2},g).
\end{eqnarray*}%
By applying $\left( \ref{eq.10}\right) $ this rewrites as
\begin{eqnarray*}
&&2B(gx_{1}\otimes 1_{H};GX_{1},g) \\
&&-B(x_{1}x_{2}\otimes 1_{H};GX_{1}X_{2},g)-B(gx_{1}\otimes
x_{2};GX_{1}X_{2},g)=0.
\end{eqnarray*}%
In view of the form of the elements we get%
\begin{equation*}
B(x_{1}x_{2}\otimes 1_{H};G,gx_{1}x_{2})-B(x_{1}x_{2}\otimes
1_{H};GX_{2},gx_{1})=0
\end{equation*}%
which is $\left( \ref{X2,x1x2,X2F71,gx1x2}\right) .$

\subparagraph{Case $f=x_{1}$}

\begin{eqnarray*}
B(gx_{1}\otimes 1_{H};GX_{1},x_{1}) &=&-B(x_{1}\otimes g;GX_{1},x_{1}) \\
&&+B(gx_{1}x_{2}\otimes g;GX_{1}X_{2},x_{1})+B(gx_{1}\otimes
x_{2};GX_{1}X_{2},x_{1}).
\end{eqnarray*}%
By applying $\left( \ref{eq.10}\right) $ this rewrites as
\begin{equation*}
B(x_{1}x_{2}\otimes 1_{H};GX_{1}X_{2},x_{1})-B(gx_{1}\otimes
x_{2};GX_{1}X_{2},x_{1})=0.
\end{equation*}%
In view of the form of the elements we get%
\begin{equation*}
-B(x_{1}x_{2}\otimes 1_{H};GX_{1},x_{1}x_{2})+B(x_{1}x_{2}\otimes
1_{H};GX_{1},x_{1}x_{2})=0
\end{equation*}%
which is trivial.

\subparagraph{Case $f=x_{2}$}

\begin{eqnarray*}
B(gx_{1}\otimes 1_{H};GX_{1},x_{2}) &=&-B(x_{1}\otimes g;GX_{1},x_{2}) \\
&&+B(gx_{1}x_{2}\otimes g;GX_{1}X_{2},x_{2})+B(gx_{1}\otimes
x_{2};GX_{1}X_{2},x_{2}).
\end{eqnarray*}%
By applying $\left( \ref{eq.10}\right) $ this rewrites as
\begin{equation*}
B(x_{1}x_{2}\otimes 1_{H};GX_{1}X_{2},x_{2})-B(gx_{1}\otimes
x_{2};GX_{1}X_{2},x_{2}).
\end{equation*}%
In view of the form of the elements we get%
\begin{equation*}
B(x_{1}x_{2}\otimes 1_{H};GX_{2},x_{1}x_{2})-B(x_{1}x_{2}\otimes
1_{H};GX_{2},x_{1}x_{2})=0
\end{equation*}%
which is trivial.

\subsubsection{Case $gx_{2}\otimes 1_{H}$}

\paragraph{\textbf{Equality} $\left( \protect\ref{X2F11}\right) $}

rewrites as%
\begin{eqnarray*}
\beta _{2}B(gx_{2}\otimes 1_{H};X_{2},f) &=&-\gamma _{2}B(x_{2}\otimes
g;G,f)-\lambda B(x_{2}\otimes g;X_{1},f) \\
&&-\beta _{2}B(x_{2}\otimes g;X_{2},f)+B(gx_{2}\otimes x_{2};1_{A},f).
\end{eqnarray*}

\subparagraph{Case $f=1_{H}$}

\begin{eqnarray*}
\beta _{2}B(gx_{2}\otimes 1_{H};X_{2},1_{H}) &=&-\gamma _{2}B(x_{2}\otimes
g;G,1_{H})-\lambda B(x_{2}\otimes g;X_{1},1_{H}) \\
&&-\beta _{2}B(x_{2}\otimes g;X_{2},1_{H})+B(gx_{2}\otimes
x_{2};1_{A},1_{H}).
\end{eqnarray*}%
By applying $\left( \ref{eq.10}\right) $ this rewrites as%
\begin{eqnarray*}
&&2\beta _{2}B(gx_{2}\otimes 1_{H};X_{2},1_{H})+\gamma _{2}B(gx_{2}\otimes
1_{H};G,1_{H})+ \\
\lambda B(gx_{2}\otimes 1_{H};X_{1},1_{H}) &&-B(gx_{2}\otimes
x_{2};1_{A},1_{H})=0.
\end{eqnarray*}%
In view of the form of the elements we get%
\begin{gather}
2\beta _{2}\left[ -1+B(x_{1}x_{2}\otimes
1_{H};1_{A},x_{1}x_{2})-B(x_{1}x_{2}\otimes 1_{H};X_{1},x_{2})\right]
\label{X2,gx2,X2F11,1H} \\
+\gamma _{2}\left[ -B(x_{1}x_{2}\otimes 1_{H};G,x_{1})+B(x_{1}x_{2}\otimes
1_{H};GX_{1},1_{H})\right] +  \notag \\
-\lambda B(x_{1}x_{2}\otimes 1_{H};X_{1},x_{1})=0.  \notag
\end{gather}

\subparagraph{Case $f=x_{1}x_{2}$}

\begin{eqnarray*}
\beta _{2}B(gx_{2}\otimes 1_{H};X_{2},x_{1}x_{2}) &=&-\gamma
_{2}B(x_{2}\otimes g;G,x_{1}x_{2})-\lambda B(x_{2}\otimes g;X_{1},x_{1}x_{2})
\\
&&-\beta _{2}B(x_{2}\otimes g;X_{2},x_{1}x_{2})+B(gx_{2}\otimes
x_{2};1_{A},x_{1}x_{2}).
\end{eqnarray*}%
By applying $\left( \ref{eq.10}\right) $ this rewrites as%
\begin{eqnarray*}
&&2\beta _{2}B(gx_{2}\otimes 1_{H};X_{2},x_{1}x_{2})+\gamma
_{2}B(gx_{2}\otimes 1_{H};G,x_{1}x_{2}) \\
&&+\lambda B(gx_{2}\otimes 1_{H};X_{1},x_{1}x_{2})-B(gx_{2}\otimes
x_{2};1_{A},x_{1}x_{2})=0.
\end{eqnarray*}%
In view of the form of the elements we get%
\begin{equation*}
\gamma _{2}B(x_{1}x_{2}\otimes 1_{H};GX_{1},x_{1}x_{2})=0
\end{equation*}%
which is $\left( \ref{G,x1x2, GF6,x1}\right) .$

\subparagraph{Case $f=gx_{1}$}

\begin{eqnarray*}
\beta _{2}B(gx_{2}\otimes 1_{H};X_{2},gx_{1}) &=&-\gamma _{2}B(x_{2}\otimes
g;G,gx_{1})-\lambda B(x_{2}\otimes g;X_{1},gx_{1}) \\
&&-\beta _{2}B(x_{2}\otimes g;X_{2},gx_{1})+B(gx_{2}\otimes
x_{2};1_{A},gx_{1}).
\end{eqnarray*}%
By applying $\left( \ref{eq.10}\right) $ this rewrites as%
\begin{eqnarray*}
&&-\gamma _{2}B(gx_{2}\otimes 1_{H};G,gx_{1})-\lambda B(gx_{2}\otimes
1_{H};X_{1},gx_{1}) \\
&&-B(gx_{2}\otimes x_{2};1_{A},gx_{1})=0.
\end{eqnarray*}%
In view of the form of the elements we get%
\begin{equation*}
\gamma _{2}B(x_{1}x_{2}\otimes 1_{H};GX_{1},gx_{1})=0
\end{equation*}%
which is $\left( \ref{X2,x1x2,X2F31,gx1}\right) .$

\subparagraph{Case $f=gx_{2}$}

\begin{eqnarray*}
\beta _{2}B(gx_{2}\otimes 1_{H};X_{2},gx_{2}) &=&-\gamma _{2}B(x_{2}\otimes
g;G,gx_{2})-\lambda B(x_{2}\otimes g;X_{1},gx_{2}) \\
&&-\beta _{2}B(x_{2}\otimes g;X_{2},gx_{2})+B(gx_{2}\otimes
x_{2};1_{A},gx_{2}).
\end{eqnarray*}%
By applying $\left( \ref{eq.10}\right) $ this rewrites as%
\begin{eqnarray*}
&&-\gamma _{2}B(gx_{2}\otimes 1_{H};G,gx_{2})-\lambda B(gx_{2}\otimes
1_{H};X_{1},gx_{2}) \\
&&-B(gx_{2}\otimes x_{2};1_{A},gx_{2})=0.
\end{eqnarray*}%
In view of the form of the elements we get

\begin{eqnarray}
&&\gamma _{2}\left[ B(x_{1}x_{2}\otimes
1_{H};G,gx_{1}x_{2})+B(x_{1}x_{2}\otimes 1_{H};GX_{1},gx_{2})\right]
\label{X2,gx2,X2F11,gx2} \\
&&+\lambda B(x_{1}x_{2}\otimes 1_{H};X_{1},gx_{1}x_{2})  \notag \\
&&+\left[ +2B(x_{1}x_{2}\otimes 1_{H};1_{A},gx_{1})+2B(x_{1}x_{2}\otimes
1_{H};X_{1},g)\right] =0.  \notag
\end{eqnarray}

\paragraph{\textbf{Equality} $\left( \protect\ref{X2F21}\right) $}

rewrites as%
\begin{eqnarray*}
\beta _{2}B(gx_{2}\otimes 1_{H};GX_{2},f) &=&+\beta _{2}B(x_{2}\otimes
g;GX_{2},f)+ \\
&&\lambda B(x_{2}\otimes g;GX_{1},f)+B(gx_{2}\otimes x_{2};G,f).
\end{eqnarray*}

\subparagraph{Case $f=g$}

\begin{eqnarray*}
\beta _{2}B(gx_{2}\otimes 1_{H};GX_{2},g) &=&+\beta _{2}B(x_{2}\otimes
g;GX_{2},g)+ \\
&&\lambda B(x_{2}\otimes g;GX_{1},g)+B(gx_{2}\otimes x_{2};G,g).
\end{eqnarray*}%
By applying $\left( \ref{eq.10}\right) $ this rewrites as%
\begin{equation*}
-\lambda B(gx_{2}\otimes 1_{H};GX_{1},g)-B(gx_{2}\otimes x_{2};G,g)=0.
\end{equation*}%
In view of the form of the elements we get%
\begin{equation*}
\lambda B(x_{1}x_{2}\otimes 1_{H};GX_{1},gx_{1})=0
\end{equation*}%
which is $\left( \ref{X2,x1x2,X2F21,gx1}\right) .$

\subparagraph{Case $f=x_{1}$}

\begin{eqnarray*}
\beta _{2}B(gx_{2}\otimes 1_{H};GX_{2},x_{1}) &=&+\beta _{2}B(x_{2}\otimes
g;GX_{2},x_{1})+ \\
&&\lambda B(x_{2}\otimes g;GX_{1},x_{1})+B(gx_{2}\otimes x_{2};G,x_{1}).
\end{eqnarray*}%
By applying $\left( \ref{eq.10}\right) $ this rewrites as%
\begin{eqnarray*}
&&2\beta _{2}B(gx_{2}\otimes 1_{H};GX_{2},x_{1})+ \\
+ &&\lambda B(gx_{2}\otimes 1_{H};GX_{1},x_{1})-B(gx_{2}\otimes
x_{2};G,x_{1})=0.
\end{eqnarray*}%
In view of the form of the elements we get%
\begin{equation*}
\beta _{2}B(x_{1}x_{2}\otimes 1_{H};GX_{1},x_{1}x_{2})=0
\end{equation*}%
which is $\left( \ref{X2,x1x2,X2F61,x1}\right) .$

\subparagraph{Case $f=x_{2}$}

\begin{eqnarray*}
\beta _{2}B(gx_{2}\otimes 1_{H};GX_{2},x_{2}) &=&+\beta _{2}B(x_{2}\otimes
g;GX_{2},x_{2})+ \\
&&\lambda B(x_{2}\otimes g;GX_{1},x_{2})+B(gx_{2}\otimes x_{2};G,x_{2}).
\end{eqnarray*}%
By applying $\left( \ref{eq.10}\right) $ this rewrites as%
\begin{eqnarray*}
&&2\beta _{2}B(gx_{2}\otimes 1_{H};GX_{2},x_{2}) \\
&&\lambda B(gx_{2}\otimes 1_{H};GX_{1},x_{2})-B(gx_{2}\otimes
x_{2};G,x_{2})=0.
\end{eqnarray*}%
In view of the form of the elements we get%
\begin{equation*}
\lambda B(x_{1}x_{2}\otimes 1_{H};GX_{1},x_{1}x_{2})=0
\end{equation*}%
which is $\left( \ref{X1,x1x2,X1F61,x1}\right) .$

\subparagraph{Case $f=gx_{1}x_{2}$}

\begin{eqnarray*}
\beta _{2}B(gx_{2}\otimes 1_{H};GX_{2},gx_{1}x_{2}) &=&+\beta
_{2}B(x_{2}\otimes g;GX_{2},gx_{1}x_{2})+ \\
&&\lambda B(x_{2}\otimes g;GX_{1},gx_{1}x_{2})+B(gx_{2}\otimes
x_{2};G,gx_{1}x_{2}).
\end{eqnarray*}%
By applying $\left( \ref{eq.10}\right) $ this rewrites as%
\begin{equation*}
-\lambda B(gx_{2}\otimes 1_{H};GX_{1},gx_{1}x_{2})-B(gx_{2}\otimes
x_{2};G,gx_{1}x_{2})=0.
\end{equation*}%
In view of the form of the elements we get

\begin{equation*}
B(x_{1}x_{2}\otimes 1_{H};GX_{1},gx_{1})=0
\end{equation*}%
which is $\left( \ref{X2,x1x2,X2F61,gx1x2}\right) .$

\paragraph{\textbf{Equality} $\left( \protect\ref{X2F31}\right) $}

rewrites as%
\begin{eqnarray*}
\beta _{2}B(gx_{2}\otimes 1_{H};X_{1}X_{2},f) &=&+\beta _{2}B(x_{2}\otimes
g;X_{1}X_{2},f) \\
&&-\gamma _{2}B(x_{2}\otimes g;GX_{1},f)+B(gx_{2}\otimes x_{2};X_{1},f).
\end{eqnarray*}

\subparagraph{Case $f=g$}

\begin{eqnarray*}
\beta _{2}B(gx_{2}\otimes 1_{H};X_{1}X_{2},g) &=&+\beta _{2}B(x_{2}\otimes
g;X_{1}X_{2},g) \\
&&-\gamma _{2}B(x_{2}\otimes g;GX_{1},g)+B(gx_{2}\otimes x_{2};X_{1},g).
\end{eqnarray*}%
By applying $\left( \ref{eq.10}\right) $ this rewrites as%
\begin{equation*}
-\gamma _{2}B(gx_{2}\otimes 1_{H};GX_{1},g)+B(gx_{2}\otimes x_{2};X_{1},g)=0.
\end{equation*}%
In view of the form of the elements we get%
\begin{equation*}
\gamma _{2}B(x_{1}x_{2}\otimes 1_{H};GX_{1},gx_{1})=0
\end{equation*}%
which follows from $\left( \ref{X2,x1x2,X2F61,gx1x2}\right) .$

\subparagraph{Case $f=x_{2}$}

\begin{eqnarray*}
\beta _{2}B(gx_{2}\otimes 1_{H};X_{1}X_{2},x_{2}) &=&+\beta
_{2}B(x_{2}\otimes g;X_{1}X_{2},x_{2}) \\
&&-\gamma _{2}B(x_{2}\otimes g;GX_{1},x_{2})+B(gx_{2}\otimes
x_{2};X_{1},x_{2}).
\end{eqnarray*}%
By applying $\left( \ref{eq.10}\right) $ this rewrites as%
\begin{eqnarray*}
&&2\beta _{2}B(gx_{2}\otimes 1_{H};X_{1}X_{2},x_{2}) \\
&&-\gamma _{2}B(gx_{2}\otimes 1_{H};GX_{1},x_{2})-B(gx_{2}\otimes
x_{2};X_{1},x_{2})=0.
\end{eqnarray*}%
In view of the form of the elements we get%
\begin{equation*}
\gamma _{2}B(x_{1}x_{2}\otimes 1_{H};GX_{1},x_{1}x_{2})=0
\end{equation*}%
which is $\left( \ref{G,x1x2, GF6,x1}\right) .$

\paragraph{\textbf{Equality} $\left( \protect\ref{X2F41}\right) $}

rewrites as

\begin{eqnarray*}
B(gx_{2}\otimes 1_{H};1_{A},f) &=&-B(x_{2}\otimes g;1_{A},f)-\gamma
_{2}B(x_{2}\otimes g;GX_{2},f) \\
&&-\lambda B(x_{2}\otimes g;X_{1}X_{2},f)+B(gx_{2}\otimes x_{2};X_{2},f).
\end{eqnarray*}

\subparagraph{Case $f=g$}

\begin{eqnarray*}
B(gx_{2}\otimes 1_{H};1_{A},g) &=&-B(x_{2}\otimes g;1_{A},g)-\gamma
_{2}B(x_{2}\otimes g;GX_{2},g) \\
&&-\lambda B(x_{2}\otimes g;X_{1}X_{2},g)+B(gx_{2}\otimes x_{2};X_{2},g).
\end{eqnarray*}%
By applying $\left( \ref{eq.10}\right) $ this rewrites as%
\begin{eqnarray*}
&&2B(gx_{2}\otimes 1_{H};1_{A},g)+\gamma _{2}B(gx_{2}\otimes 1_{H};GX_{2},g)
\\
&&+\lambda B(gx_{2}\otimes 1_{H};X_{1}X_{2},g)-B(gx_{2}\otimes
x_{2};X_{2},g)=0.
\end{eqnarray*}%
In view of the form of the elements we get%
\begin{eqnarray}
&&2\left[ -B(x_{1}x_{2}\otimes 1_{H};1_{A},gx_{1})-B(x_{1}x_{2}\otimes
1_{H};X_{1},g)\right]  \label{X2,gx2,X2F41,g} \\
&&+\gamma _{2}\left[ -B(x_{1}x_{2}\otimes
1_{H};G,gx_{1}x_{2})-B(x_{1}x_{2}\otimes 1_{H};GX_{1},gx_{2})\right]  \notag
\\
&&-\lambda B(x_{1}x_{2}\otimes 1_{H};X_{1},gx_{1}x_{2})=0.  \notag
\end{eqnarray}

\subparagraph{Case $f=x_{1}$}

\begin{eqnarray*}
B(gx_{2}\otimes 1_{H};1_{A},x_{1}) &=&-B(x_{2}\otimes g;1_{A},x_{1})-\gamma
_{2}B(x_{2}\otimes g;GX_{2},x_{1}) \\
&&-\lambda B(x_{2}\otimes g;X_{1}X_{2},x_{1})+B(gx_{2}\otimes
x_{2};X_{2},x_{1}).
\end{eqnarray*}%
By applying $\left( \ref{eq.10}\right) $ this rewrites as%
\begin{eqnarray*}
&&-\gamma _{2}B(gx_{2}\otimes 1_{H};GX_{2},x_{1}) \\
&&-\lambda B(gx_{2}\otimes 1_{H};X_{1}X_{2},x_{1})-B(gx_{2}\otimes
x_{2};X_{2},x_{1}).
\end{eqnarray*}%
In view of the form of the elements we get%
\begin{equation*}
\gamma _{2}B(x_{1}x_{2}\otimes 1_{H};GX_{1},x_{1}x_{2})=0
\end{equation*}%
which is $\left( \ref{G,x1x2, GF6,x1}\right) .$

\subparagraph{Case $f=x_{2}$}

\begin{eqnarray*}
B(gx_{2}\otimes 1_{H};1_{A},x_{2}) &=&-B(x_{2}\otimes g;1_{A},x_{2})-\gamma
_{2}B(x_{2}\otimes g;GX_{2},x_{2}) \\
&&-\lambda B(x_{2}\otimes g;X_{1}X_{2},x_{2})+B(gx_{2}\otimes
x_{2};X_{2},x_{2}).
\end{eqnarray*}%
By applying $\left( \ref{eq.10}\right) $ this rewrites as%
\begin{eqnarray*}
&&-\gamma _{2}B(gx_{2}\otimes 1_{H};GX_{2},x_{2}) \\
&&-\lambda B(gx_{2}\otimes 1_{H};X_{1}X_{2},x_{2})-B(gx_{2}\otimes
x_{2};X_{2},x_{2})=0.
\end{eqnarray*}%
In view of the form of the elements we get%
\begin{equation*}
B(gx_{2}\otimes x_{2};X_{2},x_{2})=0
\end{equation*}%
which is already known.

\subparagraph{Case $f=gx_{1}x_{2}$}

\begin{eqnarray*}
B(gx_{2}\otimes 1_{H};1_{A},gx_{1}x_{2}) &=&-B(x_{2}\otimes
g;1_{A},gx_{1}x_{2})-\gamma _{2}B(x_{2}\otimes g;GX_{2},gx_{1}x_{2}) \\
&&-\lambda B(x_{2}\otimes g;X_{1}X_{2},gx_{1}x_{2})+B(gx_{2}\otimes
x_{2};X_{2},gx_{1}x_{2}).
\end{eqnarray*}%
By applying $\left( \ref{eq.10}\right) $ this rewrites as%
\begin{eqnarray*}
&&2B(gx_{2}\otimes 1_{H};1_{A},gx_{1}x_{2})+\gamma _{2}B(gx_{2}\otimes
1_{H};GX_{2},gx_{1}x_{2}) \\
&&+\lambda B(gx_{2}\otimes 1_{H};X_{1}X_{2},gx_{1}x_{2})-B(gx_{2}\otimes
x_{2};X_{2},gx_{1}x_{2})=0.
\end{eqnarray*}%
In view of the form of the elements we get%
\begin{equation*}
-2B(x_{1}x_{2}\otimes 1_{H};X_{1},gx_{1}x_{2})+2B(x_{1}x_{2}\otimes
1_{H};X_{1},gx_{1}x_{2})=0
\end{equation*}%
which is trivial.

\paragraph{\textbf{Equality} $\left( \protect\ref{X2F51}\right) $}

rewrites as

\begin{eqnarray*}
B(gx_{2}\otimes 1_{H};X_{1},f) &=&+B(x_{2}\otimes g;X_{1},f)-\gamma
_{2}B(x_{2}\otimes g;GX_{1}X_{2},f) \\
&&+B(gx_{2}\otimes x_{2};X_{1}X_{2},f).
\end{eqnarray*}

\subparagraph{Case $f=1_{H}$}

\begin{eqnarray*}
B(gx_{2}\otimes 1_{H};X_{1},1_{H}) &=&+B(x_{2}\otimes g;X_{1},1_{H})-\gamma
_{2}B(x_{2}\otimes g;GX_{1}X_{2},1_{H}) \\
&&+B(gx_{2}\otimes x_{2};X_{1}X_{2},1_{H}).
\end{eqnarray*}%
By applying $\left( \ref{eq.10}\right) $ this rewrites as%
\begin{equation*}
+\gamma _{2}B(gx_{2}\otimes 1_{H};GX_{1}X_{2},1_{H})-B(gx_{2}\otimes
x_{2};X_{1}X_{2},1_{H})=0.
\end{equation*}%
In view of the form of the elements we get%
\begin{equation*}
\gamma _{2}B(x_{1}x_{2}\otimes 1_{H};GX_{1},x_{1}x_{2})=0
\end{equation*}%
which is $\left( \ref{G,x1x2, GF6,x1}\right) .$

\subparagraph{Case $f=gx_{2}$}

\begin{eqnarray*}
B(gx_{2}\otimes 1_{H};X_{1},gx_{2}) &=&+B(x_{2}\otimes
g;X_{1},gx_{2})-\gamma _{2}B(x_{2}\otimes g;GX_{1}X_{2},gx_{2}) \\
&&+B(gx_{2}\otimes x_{2};X_{1}X_{2},gx_{2}).
\end{eqnarray*}%
By applying $\left( \ref{eq.10}\right) $ this rewrites as%
\begin{eqnarray*}
&&2B(gx_{2}\otimes 1_{H};X_{1},gx_{2})-\gamma _{2}B(gx_{2}\otimes
1_{H};GX_{1}X_{2},gx_{2}) \\
&&-B(gx_{2}\otimes x_{2};X_{1}X_{2},gx_{2})=0.
\end{eqnarray*}%
In view of the form of the elements we get%
\begin{equation*}
2B(x_{1}x_{2}\otimes 1_{H};X_{1},gx_{1}x_{2})-2B(x_{1}x_{2}\otimes
1_{H};X_{1},gx_{1}x_{2})=0
\end{equation*}

which is trivial.

\paragraph{\textbf{Equality} $\left( \protect\ref{X2F61}\right) $}

rewrites as

\subparagraph{%
\protect\begin{eqnarray*}
\protect\beta _{2}B(gx_{2}\otimes 1_{H};GX_{1}X_{2},f) &=&-\protect\beta %
_{2}B(x_{2}\otimes g;GX_{1}X_{2},f) \\
&&+B(gx_{2}\otimes x_{2};GX_{1},f).
\protect\end{eqnarray*}%
}

\subparagraph{Case $f=1_{H}$}

\begin{eqnarray*}
\beta _{2}B(gx_{2}\otimes 1_{H};GX_{1}X_{2},1_{H}) &=&-\beta
_{2}B(x_{2}\otimes g;GX_{1}X_{2},1_{H}) \\
&&+B(gx_{2}\otimes x_{2};GX_{1},1_{H}).
\end{eqnarray*}%
By applying $\left( \ref{eq.10}\right) $ this rewrites as%
\begin{equation*}
2\beta _{2}B(gx_{2}\otimes 1_{H};GX_{1}X_{2},1_{H})-B(gx_{2}\otimes
x_{2};GX_{1},1_{H})=0.
\end{equation*}%
In view of the form of the elements we get%
\begin{equation*}
\beta _{2}B(x_{1}x_{2}\otimes 1_{H};GX_{1},x_{1}x_{2})=0
\end{equation*}%
which is $\left( \ref{X2,x1x2,X2F61,x1}\right) .$

\subparagraph{Case $f=gx_{2}$}

\begin{eqnarray*}
\beta _{2}B(gx_{2}\otimes 1_{H};GX_{1}X_{2},gx_{2}) &=&-\beta
_{2}B(x_{2}\otimes g;GX_{1}X_{2},gx_{2}) \\
&&+B(gx_{2}\otimes x_{2};GX_{1},gx_{2}).
\end{eqnarray*}%
By applying $\left( \ref{eq.10}\right) $ this rewrites as%
\begin{equation*}
B(gx_{2}\otimes x_{2};GX_{1},gx_{2})=0
\end{equation*}%
In view of the form of the element we get%
\begin{equation*}
B(x_{1}x_{2}\otimes 1_{H};GX_{1},gx_{1})=0
\end{equation*}%
which is $\left( \ref{X2,x1x2,X2F61,gx1x2}\right) .$

\paragraph{\textbf{Equality} $\left( \protect\ref{X2F71}\right) $}

rewrites as%
\begin{eqnarray*}
B(gx_{2}\otimes 1_{H};G,f) &=&+B(x_{2}\otimes g;G,f)+\lambda B(x_{2}\otimes
g;GX_{1}X_{2},f) \\
&&+B(gx_{2}\otimes x_{2};GX_{2},f).
\end{eqnarray*}

\subparagraph{Case $f=1_{H}$}

\begin{eqnarray*}
B(gx_{2}\otimes 1_{H};G,1_{H}) &=&+B(x_{2}\otimes g;G,1_{H})+\lambda
B(x_{2}\otimes g;GX_{1}X_{2},1_{H}) \\
&&+B(gx_{2}\otimes x_{2};GX_{2},1_{H}).
\end{eqnarray*}%
By applying $\left( \ref{eq.10}\right) $ this rewrites as%
\begin{equation*}
\lambda B(gx_{2}\otimes 1_{H};GX_{1}X_{2},1_{H})+B(gx_{2}\otimes
x_{2};GX_{2},1_{H})=0.
\end{equation*}%
In view of the form of the elements we get%
\begin{equation*}
\lambda B(x_{1}x_{2}\otimes 1_{H};GX_{1},x_{1}x_{2})=0
\end{equation*}%
which is $\left( \ref{X1,x1x2,X1F61,x1}\right) .$

\subparagraph{Case $f=x_{1}x_{2}$}

\begin{eqnarray*}
B(gx_{2}\otimes 1_{H};G,x_{1}x_{2}) &=&+B(x_{2}\otimes
g;G,x_{1}x_{2})+\lambda B(x_{2}\otimes g;GX_{1}X_{2},x_{1}x_{2}) \\
&&+B(gx_{2}\otimes x_{2};GX_{2},x_{1}x_{2}).
\end{eqnarray*}%
By applying $\left( \ref{eq.10}\right) $ this rewrites as%
\begin{equation*}
\lambda B(gx_{2}\otimes 1_{H};GX_{1}X_{2},x_{1}x_{2})+B(gx_{2}\otimes
x_{2};GX_{2},x_{1}x_{2})=0.
\end{equation*}%
In view of the form of the elements this is trivial.

\subparagraph{Case $f=gx_{1}$}

\begin{eqnarray*}
B(gx_{2}\otimes 1_{H};G,gx_{1}) &=&+B(x_{2}\otimes g;G,gx_{1})+\lambda
B(x_{2}\otimes g;GX_{1}X_{2},gx_{1}) \\
&&+B(gx_{2}\otimes x_{2};GX_{2},gx_{1}).
\end{eqnarray*}%
By applying $\left( \ref{eq.10}\right) $ this rewrites as%
\begin{eqnarray*}
&&2B(gx_{2}\otimes 1_{H};G,gx_{1})+\lambda B(gx_{2}\otimes
1_{H};GX_{1}X_{2},gx_{1}) \\
&&-B(gx_{2}\otimes x_{2};GX_{2},gx_{1})=0.
\end{eqnarray*}%
In view of the form of the elements we get%
\begin{equation*}
B(x_{1}x_{2}\otimes 1_{H};GX_{1},gx_{1})=0
\end{equation*}%
which is $\left( \ref{X2,x1x2,X2F61,gx1x2}\right) .$

\subparagraph{Case $f=gx_{2}$}

\begin{eqnarray*}
B(gx_{2}\otimes 1_{H};G,gx_{2}) &=&+B(x_{2}\otimes g;G,gx_{2})+\lambda
B(x_{2}\otimes g;GX_{1}X_{2},gx_{2}) \\
&&+B(gx_{2}\otimes x_{2};GX_{2},gx_{2}).
\end{eqnarray*}%
By applying $\left( \ref{eq.10}\right) $ this rewrites as%
\begin{eqnarray*}
&&2B(gx_{2}\otimes 1_{H};G,gx_{2})+\lambda B(gx_{2}\otimes
1_{H};GX_{1}X_{2},gx_{2}) \\
&&-B(gx_{2}\otimes x_{2};GX_{2},gx_{2})=0.
\end{eqnarray*}%
In view of the form of the elements we get%
\begin{eqnarray*}
&&2\left[ B(x_{1}x_{2}\otimes 1_{H};G,gx_{1}x_{2})+B(x_{1}x_{2}\otimes
1_{H};GX_{1},gx_{2})\right] \\
&&+\left[ -2B(x_{1}x_{2}\otimes 1_{H};G,gx_{1}x_{2})-2B(x_{1}x_{2}\otimes
1_{H};GX_{1},gx_{2})\right] =0
\end{eqnarray*}%
which is trivial.

\paragraph{\textbf{Equality} $\left( \protect\ref{X2F81}\right) $}

rewrites as%
\begin{equation*}
B(gx_{2}\otimes 1_{H};GX_{1},f)=-B(x_{2}\otimes g;GX_{1},f)+B(gx_{2}\otimes
x_{2};GX_{1}X_{2},f)=0.
\end{equation*}

\subparagraph{Case $f=g$}

\begin{equation*}
B(gx_{2}\otimes 1_{H};GX_{1},g)=-B(x_{2}\otimes g;GX_{1},g)+B(gx_{2}\otimes
x_{2};GX_{1}X_{2},g)=0.
\end{equation*}%
By applying $\left( \ref{eq.10}\right) $ this rewrites as%
\begin{equation*}
2B(gx_{2}\otimes 1_{H};GX_{1},g)-B(gx_{2}\otimes x_{2};GX_{1}X_{2},g)=0.
\end{equation*}%
In view of the form of the elements we get%
\begin{equation*}
B(x_{1}x_{2}\otimes 1_{H};GX_{1},gx_{1})=0
\end{equation*}%
which is $\left( \ref{X2,x1x2,X2F61,gx1x2}\right) .$

\subparagraph{Case $f=x_{2}$}

\begin{equation*}
B(gx_{2}\otimes 1_{H};GX_{1},x_{2})=-B(x_{2}\otimes
g;GX_{1},x_{2})+B(gx_{2}\otimes x_{2};GX_{1}X_{2},x_{2})=0.
\end{equation*}%
By applying $\left( \ref{eq.10}\right) $ this rewrites as%
\begin{equation*}
B(gx_{2}\otimes x_{2};GX_{1}X_{2},x_{2})=0
\end{equation*}%
In view of the form of the element, this is already known.

\subsubsection{Case $gx_{1}x_{2}\otimes 1_{H}$}

\paragraph{\textbf{Equality} $\left( \protect\ref{X2F11}\right) $}

rewrites as%
\begin{eqnarray*}
\beta _{2}B(gx_{1}x_{2}\otimes 1_{H};X_{2},f) &=&\gamma
_{2}B(x_{1}x_{2}\otimes g;G,f)+\lambda B(x_{1}x_{2}\otimes g;X_{1},f) \\
&&+\beta _{2}B(x_{1}x_{2}\otimes g;X_{2},f)+B(gx_{1}x_{2}\otimes
x_{2};1_{A},f).
\end{eqnarray*}

\subparagraph{Case $f=1_{H}$}

\begin{eqnarray*}
\beta _{2}B(gx_{1}x_{2}\otimes 1_{H};X_{2},1_{H}) &=&\gamma
_{2}B(x_{1}x_{2}\otimes g;G,1_{H})+\lambda B(x_{1}x_{2}\otimes g;X_{1},1_{H})
\\
&&+\beta _{2}B(x_{1}x_{2}\otimes g;X_{2},1_{H})+B(gx_{1}x_{2}\otimes
x_{2};1_{A},1_{H}).
\end{eqnarray*}%
By applying $\left( \ref{eq.10}\right) $ this rewrites as%
\begin{eqnarray*}
&&-\gamma _{2}B(gx_{1}x_{2}\otimes 1_{H};G,1_{H})-\lambda
B(gx_{1}x_{2}\otimes 1_{H};X_{1},1_{H}) \\
&&-B(gx_{1}x_{2}\otimes x_{2};1_{A},1_{H})=0.
\end{eqnarray*}%
In view of the form of the elements we get%
\begin{equation}
-\gamma _{2}B(gx_{1}x_{2}\otimes 1_{H};G,1_{H})+\lambda \left[
+B(x_{2}\otimes 1_{H};1_{A},1_{H})-B(gx_{1}x_{2}\otimes 1_{H};1_{A},x_{1})%
\right] =0  \label{X2,gx1x2,X2F11,1H}
\end{equation}

\subparagraph{Case $f=x_{1}x_{2}$}

\begin{eqnarray*}
\beta _{2}B(gx_{1}x_{2}\otimes 1_{H};X_{2},x_{1}x_{2}) &=&\gamma
_{2}B(x_{1}x_{2}\otimes g;G,x_{1}x_{2})+\lambda B(x_{1}x_{2}\otimes
g;X_{1},x_{1}x_{2}) \\
&&+\beta _{2}B(x_{1}x_{2}\otimes g;X_{2},x_{1}x_{2})+B(gx_{1}x_{2}\otimes
x_{2};1_{A},x_{1}x_{2}).
\end{eqnarray*}%
By applying $\left( \ref{eq.10}\right) $ this rewrites as%
\begin{eqnarray*}
&&-\gamma _{2}B(gx_{1}x_{2}\otimes 1_{H};G,x_{1}x_{2})-\lambda
B(gx_{1}x_{2}\otimes 1_{H};X_{1},x_{1}x_{2}) \\
&&-B(gx_{1}x_{2}\otimes x_{2};1_{A},x_{1}x_{2})=0.
\end{eqnarray*}%
In view of the form of the elements we get%
\begin{equation}
\gamma _{2}B(gx_{1}x_{2}\otimes 1_{H};G,x_{1}x_{2})-\lambda B(x_{2}\otimes
1_{H};1_{A},x_{1}x_{2})=0  \label{X2,gx1x2,X2F11,x1x2}
\end{equation}

\subparagraph{Case $f=gx_{1}$}

\begin{eqnarray*}
\beta _{2}B(gx_{1}x_{2}\otimes 1_{H};X_{2},gx_{1}) &=&\gamma
_{2}B(x_{1}x_{2}\otimes g;G,gx_{1})+\lambda B(x_{1}x_{2}\otimes
g;X_{1},gx_{1}) \\
&&+\beta _{2}B(x_{1}x_{2}\otimes g;X_{2},gx_{1})+B(gx_{1}x_{2}\otimes
x_{2};1_{A},gx_{1}).
\end{eqnarray*}%
By applying $\left( \ref{eq.10}\right) $ this rewrites as%
\begin{gather*}
2\beta _{2}B(gx_{1}x_{2}\otimes 1_{H};X_{2},gx_{1})+\gamma
_{2}B(gx_{1}x_{2}\otimes 1_{H};G,gx_{1}) \\
+\lambda B(gx_{1}x_{2}\otimes 1_{H};X_{1},gx_{1})-B(gx_{1}x_{2}\otimes
x_{2};1_{A},gx_{1})=0.
\end{gather*}%
In view of the form of the elements we get%
\begin{gather}
2\beta _{2}\left[ B(x_{1}\otimes 1_{H};1_{A},gx_{1})+B(gx_{1}x_{2}\otimes
1_{H};1_{A},gx_{1}x_{2})\right]  \label{X2,gx1x2,X2F11,gx1} \\
+\gamma _{2}B(gx_{1}x_{2}\otimes 1_{H};G,gx_{1})-\lambda B(x_{2}\otimes
1_{H};1_{A},gx_{1})=0.  \notag
\end{gather}

\subparagraph{Case $f=gx_{2}$}

\begin{eqnarray*}
\beta _{2}B(gx_{1}x_{2}\otimes 1_{H};X_{2},gx_{2}) &=&\gamma
_{2}B(x_{1}x_{2}\otimes g;G,gx_{2})+\lambda B(x_{1}x_{2}\otimes
g;X_{1},gx_{2}) \\
&&+\beta _{2}B(x_{1}x_{2}\otimes g;X_{2},gx_{2})+B(gx_{1}x_{2}\otimes
x_{2};1_{A},gx_{2}).
\end{eqnarray*}%
By applying $\left( \ref{eq.10}\right) $ this rewrites as%
\begin{eqnarray*}
&&2\beta _{2}B(gx_{1}x_{2}\otimes 1_{H};X_{2},gx_{2})+\gamma
_{2}B(gx_{1}x_{2}\otimes 1_{H};G,gx_{2}) \\
&&+\lambda B(gx_{1}x_{2}\otimes 1_{H};X_{1},gx_{2})-B(gx_{1}x_{2}\otimes
x_{2};1_{A},gx_{2})=0.
\end{eqnarray*}%
In view of the form of the elements we get

\begin{gather}
2\beta _{2}B(x_{1}\otimes 1_{H};1_{A},gx_{2})+\gamma
_{2}B(gx_{1}x_{2}\otimes 1_{H};G,gx_{2})  \label{X2,gx1x2,X2F11,gx2} \\
+\lambda \left[ -B(x_{2}\otimes 1_{H};1_{A},gx_{2})-B(gx_{1}x_{2}\otimes
1_{H};1_{A},gx_{1}x_{2})\right]  \notag \\
+2B(gx_{1}x_{2}\otimes 1_{H};1_{A},g)=0.  \notag
\end{gather}

\paragraph{\textbf{Equality} $\left( \protect\ref{X2F21}\right) $}

rewrites as%
\begin{gather*}
\beta _{2}B(gx_{1}x_{2}\otimes 1_{H};GX_{2},f)=-\beta
_{2}B(x_{1}x_{2}\otimes g;GX_{2},f) \\
-\lambda B(x_{1}x_{2}\otimes g;GX_{1},f)+B(gx_{1}x_{2}\otimes x_{2};G,f).
\end{gather*}

\subparagraph{Case $f=g$}

\begin{gather*}
\beta _{2}B(gx_{1}x_{2}\otimes 1_{H};GX_{2},g)=-\beta
_{2}B(x_{1}x_{2}\otimes g;GX_{2},g) \\
-\lambda B(x_{1}x_{2}\otimes g;GX_{1},g)+B(gx_{1}x_{2}\otimes x_{2};G,g).
\end{gather*}%
By applying $\left( \ref{eq.10}\right) $ this rewrites as%
\begin{gather*}
2\beta _{2}B(gx_{1}x_{2}\otimes 1_{H};GX_{2},g) \\
+\lambda B(gx_{1}x_{2}\otimes 1_{H};GX_{1},g)-B(gx_{1}x_{2}\otimes
x_{2};G,g)=0.
\end{gather*}%
In view of the form of the elements we get%
\begin{gather}
2\beta _{2}\left[ -B(x_{1}\otimes 1_{H};G,g)-B(gx_{1}x_{2}\otimes
1_{H};G,gx_{2})\right]  \label{X2,gx1x2,X2F21,g} \\
+\lambda \left[ B(x_{2}\otimes 1_{H};G,g)-B(gx_{1}x_{2}\otimes
1_{H};G,gx_{1})\right] =0.  \notag
\end{gather}

\subparagraph{Case $f=x_{1}$}

\begin{gather*}
\beta _{2}B(gx_{1}x_{2}\otimes 1_{H};GX_{2},x_{1})=-\beta
_{2}B(x_{1}x_{2}\otimes g;GX_{2},x_{1}) \\
-\lambda B(x_{1}x_{2}\otimes g;GX_{1},x_{1})+B(gx_{1}x_{2}\otimes
x_{2};G,x_{1}).
\end{gather*}

By applying $\left( \ref{eq.10}\right) $ this rewrites as%
\begin{equation*}
-\lambda B(gx_{1}x_{2}\otimes 1_{H};GX_{1},x_{1})-B(gx_{1}x_{2}\otimes
x_{2};G,x_{1})=0.
\end{equation*}%
In view of the form of the elements we get%
\begin{equation*}
\lambda B(x_{2}\otimes 1_{H};G,x_{1})=0
\end{equation*}%
which holds in view of $\left( \ref{X2,g,X2F21,x1}\right) .$

\subparagraph{Case $f=x_{2}$}

\begin{gather*}
\beta _{2}B(gx_{1}x_{2}\otimes 1_{H};GX_{2},x_{2})=-\beta
_{2}B(x_{1}x_{2}\otimes g;GX_{2},x_{2}) \\
-\lambda B(x_{1}x_{2}\otimes g;GX_{1},x_{2})+B(gx_{1}x_{2}\otimes
x_{2};G,x_{2}).
\end{gather*}

By applying $\left( \ref{eq.10}\right) $ this rewrites as%
\begin{equation*}
-\lambda B(gx_{1}x_{2}\otimes 1_{H};GX_{1},x_{2})-B(gx_{1}x_{2}\otimes
x_{2};G,x_{2})=0.
\end{equation*}%
In view of the form of the elements we get%
\begin{equation*}
\lambda \left[ B(x_{2}\otimes 1_{H};G,x_{2})+B(gx_{1}x_{2}\otimes
1_{H};G,x_{1}x_{2})\right] =0
\end{equation*}%
which follows from $\left( \ref{X2,x1,X2F61,x2}\right) .$

\subparagraph{Case $f=gx_{1}x_{2}$}

\begin{gather*}
\beta _{2}B(gx_{1}x_{2}\otimes 1_{H};GX_{2},gx_{1}x_{2})=-\beta
_{2}B(x_{1}x_{2}\otimes g;GX_{2},gx_{1}x_{2}) \\
-\lambda B(x_{1}x_{2}\otimes g;GX_{1},gx_{1}x_{2})+B(gx_{1}x_{2}\otimes
x_{2};G,gx_{1}x_{2}).
\end{gather*}%
By applying $\left( \ref{eq.10}\right) $ this rewrites as%
\begin{gather*}
2\beta _{2}B(gx_{1}x_{2}\otimes 1_{H};GX_{2},gx_{1}x_{2}) \\
+\lambda B(gx_{1}x_{2}\otimes 1_{H};GX_{1},gx_{1}x_{2})-B(gx_{1}x_{2}\otimes
x_{2};G,gx_{1}x_{2})=0.
\end{gather*}%
In view of the form of the elements we get

\begin{gather}
-2\beta _{2}B(x_{1}\otimes 1_{H};G,gx_{1}x_{2})  \label{X2,gx1x2,X2F21,gx1x2}
\\
+\lambda B(x_{2}\otimes 1_{H};G,gx_{1}x_{2})-2B(gx_{1}x_{2}\otimes
1_{H};G,gx_{1})=0.  \notag
\end{gather}

\paragraph{\textbf{Equality} $\left( \protect\ref{X2F31}\right) $}

rewrites as%
\begin{eqnarray*}
\beta _{2}B(gx_{1}x_{2}\otimes 1_{H};X_{1}X_{2},f) &=&-\beta
_{2}B(x_{1}x_{2}\otimes g;X_{1}X_{2},f) \\
+\gamma _{2}B(x_{1}x_{2}\otimes g;GX_{1},f) &&+B(gx_{1}x_{2}\otimes
x_{2};X_{1},f).
\end{eqnarray*}

\subparagraph{Case $f=g$}

\begin{eqnarray*}
\beta _{2}B(gx_{1}x_{2}\otimes 1_{H};X_{1}X_{2},g) &=&-\beta
_{2}B(x_{1}x_{2}\otimes g;X_{1}X_{2},g) \\
+\gamma _{2}B(x_{1}x_{2}\otimes g;GX_{1},g) &&+B(gx_{1}x_{2}\otimes
x_{2};X_{1},g).
\end{eqnarray*}%
By applying $\left( \ref{eq.10}\right) $ this rewrites as%
\begin{eqnarray*}
&&2\beta _{2}B(gx_{1}x_{2}\otimes 1_{H};X_{1}X_{2},g) \\
-\gamma _{2}B(gx_{1}x_{2}\otimes 1_{H};GX_{1},g) &&-B(gx_{1}x_{2}\otimes
x_{2};X_{1},g)=0.
\end{eqnarray*}%
In view of the form of the elements we get

\begin{gather}
2\beta _{2}\left[
\begin{array}{c}
B(g\otimes 1_{H};1_{A},g)+B(x_{2}\otimes \ 1_{H};1_{A},gx_{2}) \\
+B(x_{1}\otimes 1_{H};1_{A},gx_{1})+B(gx_{1}x_{2}\otimes
1_{H};1_{A},gx_{1}x_{2})%
\end{array}%
\right]  \label{X2,gx1x2,X2F31,g} \\
-\gamma _{2}\left[ B(x_{2}\otimes 1_{H};G,g)-B(gx_{1}x_{2}\otimes
1_{H};G,gx_{1})\right] =0.  \notag
\end{gather}

\subparagraph{Case $f=x_{1}$}

\begin{eqnarray*}
\beta _{2}B(gx_{1}x_{2}\otimes 1_{H};X_{1}X_{2},x_{1}) &=&-\beta
_{2}B(x_{1}x_{2}\otimes g;X_{1}X_{2},x_{1}) \\
+\gamma _{2}B(x_{1}x_{2}\otimes g;GX_{1},x_{1}) &&+B(gx_{1}x_{2}\otimes
x_{2};X_{1},x_{1}).
\end{eqnarray*}%
By applying $\left( \ref{eq.10}\right) $ this rewrites as%
\begin{equation*}
\gamma _{2}B(gx_{1}x_{2}\otimes 1_{H};GX_{1},x_{1})-B(gx_{1}x_{2}\otimes
x_{2};X_{1},x_{1})=0.
\end{equation*}%
In view of the form of the elements we get%
\begin{equation*}
\gamma _{2}B(x_{2}\otimes 1_{H};G,x_{1})=0
\end{equation*}%
which follows from $\left( \ref{X2,g,X2F21,x1}\right) .$

\subparagraph{Case $f=x_{2}$}

\begin{eqnarray*}
\beta _{2}B(gx_{1}x_{2}\otimes 1_{H};X_{1}X_{2},x_{2}) &=&-\beta
_{2}B(c;X_{1}X_{2},x_{2}) \\
+\gamma _{2}B(x_{1}x_{2}\otimes g;GX_{1},x_{2}) &&+B(gx_{1}x_{2}\otimes
x_{2};X_{1},x_{2}).
\end{eqnarray*}%
By applying $\left( \ref{eq.10}\right) $ this rewrites as%
\begin{equation*}
\gamma _{2}B(gx_{1}x_{2}\otimes 1_{H};GX_{1},x_{2})-B(gx_{1}x_{2}\otimes
x_{2};X_{1},x_{2})=0.
\end{equation*}%
In view of the form of the elements we get%
\begin{equation*}
\gamma _{2}\left[ B(x_{2}\otimes 1_{H};G,x_{2})+B(gx_{1}x_{2}\otimes
1_{H};G,x_{1}x_{2})\right] =0
\end{equation*}%
which follows from $\left( \ref{X2,x1,X2F61,x2}\right) .$

\subparagraph{Case $f=gx_{1}x_{2}$}

\begin{eqnarray*}
\beta _{2}B(gx_{1}x_{2}\otimes 1_{H};X_{1}X_{2},gx_{1}x_{2}) &=&-\beta
_{2}B(c;X_{1}X_{2},gx_{1}x_{2}) \\
+\gamma _{2}B(x_{1}x_{2}\otimes g;GX_{1},gx_{1}x_{2})
&&+B(gx_{1}x_{2}\otimes x_{2};X_{1},gx_{1}x_{2}).
\end{eqnarray*}%
By applying $\left( \ref{eq.10}\right) $ this rewrites as%
\begin{gather*}
2\beta _{2}B(gx_{1}x_{2}\otimes 1_{H};X_{1}X_{2},gx_{1}x_{2})+ \\
-\gamma _{2}B(gx_{1}x_{2}\otimes
1_{H};GX_{1},gx_{1}x_{2})-B(gx_{1}x_{2}\otimes x_{2};X_{1},gx_{1}x_{2})=0.
\end{gather*}%
In view of the form of the elements we get%
\begin{gather*}
2\beta _{2}B\left( g\otimes 1_{H};1_{A},gx_{1}x_{2}\right) + \\
-\gamma _{2}B(x_{2}\otimes 1_{H};G,gx_{1}x_{2})+2B(x_{2}\otimes
1_{H};1_{A},gx_{1})=0.
\end{gather*}

which is $\left( \ref{X2,x2,X2F11,gx1x2}\right) .$

\paragraph{\textbf{Equality} $\left( \protect\ref{X2F41}\right) $}

rewrites as%
\begin{eqnarray*}
B(gx_{1}x_{2}\otimes 1_{H};1_{A},f) &=&B(x_{1}x_{2}\otimes g;1_{A},f)+\gamma
_{2}B(x_{1}x_{2}\otimes g;GX_{2},f) \\
&&+\lambda B(x_{1}x_{2}\otimes g;X_{1}X_{2},f)+B(gx_{1}x_{2}\otimes
x_{2};X_{2},f).
\end{eqnarray*}

\subparagraph{Case $f=g$%
\protect\begin{eqnarray*}
B(gx_{1}x_{2}\otimes 1_{H};1_{A},g) &=&B(x_{1}x_{2}\otimes g;1_{A},g)+%
\protect\gamma _{2}B(x_{1}x_{2}\otimes g;GX_{2},g) \\
&&+\protect\lambda B(x_{1}x_{2}\otimes g;X_{1}X_{2},g)+B(gx_{1}x_{2}\otimes
x_{2};X_{2},g).
\protect\end{eqnarray*}%
}

By applying $\left( \ref{eq.10}\right) $ this rewrites as%
\begin{eqnarray*}
&&-\gamma _{2}B(gx_{1}x_{2}\otimes 1_{H};GX_{2},g) \\
&&-\lambda B(gx_{1}x_{2}\otimes 1_{H};X_{1}X_{2},g)-B(gx_{1}x_{2}\otimes
x_{2};X_{2},g)=0.
\end{eqnarray*}%
In view of the form of the elements we get%
\begin{eqnarray}
&&\gamma _{2}\left[ B(x_{1}\otimes 1_{H};G,g)+B(gx_{1}x_{2}\otimes
1_{H};G,gx_{2})\right]  \label{X2,gx1x2,X2F41,g} \\
&&-\lambda \left[
\begin{array}{c}
B(g\otimes 1_{H};1_{A},g)+B(x_{2}\otimes \ 1_{H};1_{A},gx_{2}) \\
+B(x_{1}\otimes 1_{H};1_{A},gx_{1})+B(gx_{1}x_{2}\otimes
1_{H};1_{A},gx_{1}x_{2})%
\end{array}%
\right] =0.  \notag
\end{eqnarray}

\subparagraph{Case $f=x_{1}$}

\begin{eqnarray*}
B(gx_{1}x_{2}\otimes 1_{H};1_{A},x_{1}) &=&B(x_{1}x_{2}\otimes
g;1_{A},x_{1})+\gamma _{2}B(x_{1}x_{2}\otimes g;GX_{2},x_{1}) \\
&&+\lambda B(x_{1}x_{2}\otimes g;X_{1}X_{2},x_{1})+B(gx_{1}x_{2}\otimes
x_{2};X_{2},x_{1}).
\end{eqnarray*}

By applying $\left( \ref{eq.10}\right) $ this rewrites as%
\begin{eqnarray*}
&&2B(gx_{1}x_{2}\otimes 1_{H};1_{A},x_{1})+\gamma _{2}B(gx_{1}x_{2}\otimes
1_{H};GX_{2},x_{1}) \\
&&+\lambda B(gx_{1}x_{2}\otimes 1_{H};X_{1}X_{2},x_{1})-B(gx_{1}x_{2}\otimes
x_{2};X_{2},x_{1})=0.
\end{eqnarray*}%
In view of the form of the elements we get%
\begin{eqnarray}
&&2B(gx_{1}x_{2}\otimes 1_{H};1_{A},x_{1})+\gamma _{2}\left[ -B(x_{1}\otimes
1_{H};G,x_{1})-B(gx_{1}x_{2}\otimes 1_{H};G,x_{1}x_{2})\right]
\label{X2,gx1x2,X2F41,x1} \\
&&+\lambda \left[ B(g\otimes 1_{H};1_{A},x_{1})+B(x_{2}\otimes \
1_{H};1_{A},x_{1}x_{2})\right] =0.  \notag
\end{eqnarray}

\subparagraph{Case $f=x_{2}$}

\subparagraph{%
\protect\begin{eqnarray*}
B(gx_{1}x_{2}\otimes 1_{H};1_{A},x_{2}) &=&B(x_{1}x_{2}\otimes
g;1_{A},x_{2})+\protect\gamma _{2}B(x_{1}x_{2}\otimes g;GX_{2},x_{2}) \\
&&+\protect\lambda B(x_{1}x_{2}\otimes
g;X_{1}X_{2},x_{2})+B(gx_{1}x_{2}\otimes x_{2};X_{2},x_{2}).
\protect\end{eqnarray*}%
}

By applying $\left( \ref{eq.10}\right) $ this rewrites as%
\begin{eqnarray*}
&&2B(gx_{1}x_{2}\otimes 1_{H};1_{A},x_{2})+\gamma _{2}B(gx_{1}x_{2}\otimes
1_{H};GX_{2},x_{2}) \\
&&+\lambda B(gx_{1}x_{2}\otimes 1_{H};X_{1}X_{2},x_{2})-B(gx_{1}x_{2}\otimes
x_{2};X_{2},x_{2})=0.
\end{eqnarray*}%
In view of the form of the elements we get%
\begin{eqnarray*}
&&2B(gx_{1}x_{2}\otimes 1_{H};1_{A},x_{2})-\gamma _{2}B(x_{1}\otimes
1_{H};G,x_{2}) \\
&&+\lambda \left[ B(g\otimes 1_{H};X_{2},1_{H})-B(x_{1}\otimes
1_{H};1_{A},x_{1}x_{2})\right] =0.
\end{eqnarray*}%
which is $\left( \ref{X2,x1,X2F11,x2}\right) .$

\subparagraph{Case $f=gx_{1}x_{2}$}

\subparagraph{%
\protect\begin{eqnarray*}
B(gx_{1}x_{2}\otimes 1_{H};1_{A},gx_{1}x_{2}) &=&B(x_{1}x_{2}\otimes
g;1_{A},gx_{1}x_{2})+\protect\gamma _{2}B(x_{1}x_{2}\otimes
g;GX_{2},gx_{1}x_{2}) \\
&&+\protect\lambda B(x_{1}x_{2}\otimes
g;X_{1}X_{2},gx_{1}x_{2})+B(gx_{1}x_{2}\otimes x_{2};X_{2},gx_{1}x_{2}).
\protect\end{eqnarray*}%
}

By applying $\left( \ref{eq.10}\right) $ this rewrites as%
\begin{eqnarray*}
&&-\gamma _{2}B(gx_{1}x_{2}\otimes 1_{H};GX_{2},gx_{1}x_{2}) \\
&&-\lambda B(gx_{1}x_{2}\otimes
1_{H};X_{1}X_{2},gx_{1}x_{2})-B(gx_{1}x_{2}\otimes
x_{2};X_{2},gx_{1}x_{2})=0.
\end{eqnarray*}%
In view of the form of the elements we get

\begin{gather*}
\gamma _{2}B(x_{1}\otimes 1_{H};G,gx_{1}x_{2})-\lambda B\left( g\otimes
1_{H};1_{A},gx_{1}x_{2}\right) \\
-\left[ 2B(x_{1}\otimes 1_{H};1_{A},gx_{1})+2B(gx_{1}x_{2}\otimes
1_{H};1_{A},gx_{1}x_{2})\right] =0.
\end{gather*}%
which is $\left( \ref{X2,x1,X2F11,gx1x2}\right) .$

\paragraph{\textbf{Equality} $\left( \protect\ref{X2F51}\right) $}

rewrites as%
\begin{gather*}
B(gx_{1}x_{2}\otimes 1_{H};X_{1},f)=-B(x_{1}x_{2}\otimes g;X_{1},f) \\
+\gamma _{2}B(x_{1}x_{2}\otimes g;GX_{1}X_{2},f)+B(gx_{1}x_{2}\otimes
x_{2};X_{1}X_{2},f)
\end{gather*}

\subparagraph{Case $f=1_{H}$}

\begin{gather*}
B(gx_{1}x_{2}\otimes 1_{H};X_{1},1_{H})=-B(x_{1}x_{2}\otimes g;X_{1},1_{H})
\\
+\gamma _{2}B(x_{1}x_{2}\otimes g;GX_{1}X_{2},1_{H})+B(gx_{1}x_{2}\otimes
x_{2};X_{1}X_{2},1_{H}).
\end{gather*}%
By applying $\left( \ref{eq.10}\right) $ this rewrites as%
\begin{gather*}
2B(gx_{1}x_{2}\otimes 1_{H};X_{1},1_{H})+ \\
-\gamma _{2}B(gx_{1}x_{2}\otimes
1_{H};GX_{1}X_{2},1_{H})-B(gx_{1}x_{2}\otimes x_{2};X_{1}X_{2},1_{H})=0.
\end{gather*}%
In view of the form of the elements we get

\begin{gather}
2\left[ -B(x_{2}\otimes 1_{H};1_{A},1_{H})+B(gx_{1}x_{2}\otimes
1_{H};1_{A},x_{1})\right] +  \label{X2,gx1x2,X2F51,1H} \\
-\gamma _{2}\left[
\begin{array}{c}
B(g\otimes 1_{H};G,1_{H})+B(x_{2}\otimes \ 1_{H};G,x_{2}) \\
+B(x_{1}\otimes 1_{H};G,x_{1})+B(gx_{1}x_{2}\otimes 1_{H};G,x_{1}x_{2})%
\end{array}%
\right] =0.  \notag
\end{gather}

\subparagraph{Case $f=x_{1}x_{2}$}

\begin{gather*}
B(gx_{1}x_{2}\otimes 1_{H};X_{1},x_{1}x_{2})=-B(x_{1}x_{2}\otimes
g;X_{1},x_{1}x_{2}) \\
+\gamma _{2}B(x_{1}x_{2}\otimes
g;GX_{1}X_{2},x_{1}x_{2})+B(gx_{1}x_{2}\otimes x_{2};X_{1}X_{2},x_{1}x_{2}).
\end{gather*}%
By applying $\left( \ref{eq.10}\right) $ this rewrites as%
\begin{gather*}
2B(gx_{1}x_{2}\otimes 1_{H};X_{1},x_{1}x_{2}) \\
-\gamma _{2}B(gx_{1}x_{2}\otimes
1_{H};GX_{1}X_{2},x_{1}x_{2})-B(gx_{1}x_{2}\otimes
x_{2};X_{1}X_{2},x_{1}x_{2})=0.
\end{gather*}%
In view of the form of the elements we get%
\begin{equation*}
2B(x_{2}\otimes 1_{H};1_{A},x_{1}x_{2})+\gamma _{2}B(g\otimes
1_{H};G,x_{1}x_{2})=0
\end{equation*}%
which is $\left( \ref{G,x2; GF6,x2}\right) $.

\subparagraph{Case $f=gx_{1}$}

\begin{gather*}
B(gx_{1}x_{2}\otimes 1_{H};X_{1},gx_{1})=-B(x_{1}x_{2}\otimes g;X_{1},gx_{1})
\\
+\gamma _{2}B(x_{1}x_{2}\otimes g;GX_{1}X_{2},gx_{1})+B(gx_{1}x_{2}\otimes
x_{2};X_{1}X_{2},gx_{1}).
\end{gather*}%
By applying $\left( \ref{eq.10}\right) $ this rewrites as%
\begin{equation*}
\gamma _{2}B(gx_{1}x_{2}\otimes
1_{H};GX_{1}X_{2},gx_{1})-B(gx_{1}x_{2}\otimes x_{2};X_{1}X_{2},gx_{1})=0.
\end{equation*}%
In view of the form of the elements we get%
\begin{equation*}
\gamma _{2}\left[ B(g\otimes 1_{H};G,gx_{1})+B(x_{2}\otimes \
1_{H};G,gx_{1}x_{2})\right] =0
\end{equation*}%
which is $\left( \ref{G,x2, GF2,gx1}\right) .$

\subparagraph{Case $f=gx_{2}$}

\begin{gather*}
B(gx_{1}x_{2}\otimes 1_{H};X_{1},f)=-B(x_{1}x_{2}\otimes g;X_{1},f) \\
+\gamma _{2}B(x_{1}x_{2}\otimes g;GX_{1}X_{2},f)+B(gx_{1}x_{2}\otimes
x_{2};X_{1}X_{2},f)
\end{gather*}%
By applying $\left( \ref{eq.10}\right) $ this rewrites as%
\begin{equation*}
\gamma _{2}B(gx_{1}x_{2}\otimes
1_{H};GX_{1}X_{2},gx_{2})-B(gx_{1}x_{2}\otimes x_{2};X_{1}X_{2},gx_{2})=0.
\end{equation*}%
In view of the form of the elements we get%
\begin{gather*}
\gamma _{2}\left[ B(g\otimes 1_{H};G,gx_{2})-B(x_{1}\otimes
1_{H};G,gx_{1}x_{2})\right] + \\
-\left[
\begin{array}{c}
2B(g\otimes 1_{H};1_{A},g)+2B(x_{2}\otimes \ 1_{H};1_{A},gx_{2}) \\
+2B(x_{1}\otimes 1_{H};1_{A},gx_{1})+2B(gx_{1}x_{2}\otimes
1_{H};1_{A},gx_{1}x_{2})%
\end{array}%
\right] =0.
\end{gather*}%
in view of $\left( \ref{G,gx1x2, GF6,gx2}\right) ,$we get%
\begin{gather*}
B(g\otimes 1_{H};1_{A},g)+B(x_{2}\otimes \ 1_{H};1_{A},gx_{2})+ \\
+B(x_{1}\otimes 1_{H};1_{A},gx_{1})+B(gx_{1}x_{2}\otimes
1_{H};1_{A},gx_{1}x_{2})=0
\end{gather*}%
which is $\left( \ref{X1,x2,X1F41,gx1}\right) .$

\paragraph{\textbf{Equality} $\left( \protect\ref{X2F61}\right) $}

rewrites as%
\begin{gather*}
\beta _{2}B(gx_{1}x_{2}\otimes 1_{H};GX_{1}X_{2},f)=+\beta
_{2}B(x_{1}x_{2}\otimes g;GX_{1}X_{2},f) \\
+B(gx_{1}x_{2}\otimes x_{2};GX_{1},f).
\end{gather*}

\subparagraph{Case $f=1_{H}$}

\begin{gather*}
\beta _{2}B(gx_{1}x_{2}\otimes 1_{H};GX_{1}X_{2},1_{H})=+\beta
_{2}B(x_{1}x_{2}\otimes g;GX_{1}X_{2},1_{H}) \\
+B(gx_{1}x_{2}\otimes x_{2};GX_{1},1_{H}).
\end{gather*}%
In view of the form of the elements we get%
\begin{equation*}
B(gx_{1}x_{2}\otimes x_{2};GX_{1},1_{H})=0.
\end{equation*}%
In view of the form of the element, this is already known.

\subparagraph{Case $f=x_{1}x_{2}$}

\begin{gather*}
\beta _{2}B(gx_{1}x_{2}\otimes 1_{H};GX_{1}X_{2},x_{1}x_{2})=+\beta
_{2}B(x_{1}x_{2}\otimes g;GX_{1}X_{2},x_{1}x_{2}) \\
+B(gx_{1}x_{2}\otimes x_{2};GX_{1},x_{1}x_{2}).
\end{gather*}%
In view of the form of the elements we get%
\begin{equation*}
B(gx_{1}x_{2}\otimes x_{2};GX_{1},x_{1}x_{2})=0.
\end{equation*}%
In view of the form of the element, this is already known.

\subparagraph{Case $f=gx_{1}$}

\begin{gather*}
\beta _{2}B(gx_{1}x_{2}\otimes 1_{H};GX_{1}X_{2},gx_{1})=+\beta
_{2}B(x_{1}x_{2}\otimes g;GX_{1}X_{2},gx_{1}) \\
+B(gx_{1}x_{2}\otimes x_{2};GX_{1},gx_{1}).
\end{gather*}%
In view of the form of the elements we get%
\begin{equation*}
2\beta _{2}B(gx_{1}x_{2}\otimes
1_{H};GX_{1}X_{2},gx_{1})-B(gx_{1}x_{2}\otimes x_{2};GX_{1},gx_{1})=0
\end{equation*}%
In view of the form of the elements we get%
\begin{equation*}
\beta _{2}\left[ B(g\otimes 1_{H};G,gx_{1})+B(x_{2}\otimes \
1_{H};G,gx_{1}x_{2})\right] =0
\end{equation*}%
which is $\left( \ref{X2,x2,X2F21,gx1}\right) .$

\subparagraph{Case $f=gx_{2}$}

\begin{gather*}
\beta _{2}B(gx_{1}x_{2}\otimes 1_{H};GX_{1}X_{2},gx_{2})=+\beta
_{2}B(x_{1}x_{2}\otimes g;GX_{1}X_{2},gx_{2}) \\
+B(gx_{1}x_{2}\otimes x_{2};GX_{1},gx_{2}).
\end{gather*}%
In view of the form of the elements we get%
\begin{equation*}
2\beta _{2}B(gx_{1}x_{2}\otimes
1_{H};GX_{1}X_{2},gx_{2})-B(gx_{1}x_{2}\otimes x_{2};GX_{1},gx_{2})=0
\end{equation*}%
In view of the form of the elements we get%
\begin{gather*}
2\beta _{2}\left[ B(g\otimes 1_{H};G,gx_{2})-B(x_{1}\otimes
1_{H};G,gx_{1}x_{2})\right] \\
+\left[ 2B(x_{2}\otimes 1_{H};G,g)-2B(gx_{1}x_{2}\otimes 1_{H};G,gx_{1})%
\right] =0.
\end{gather*}%
In view of $\left( \ref{X1,g,X1F21,gx1x2}\right) ,$we get%
\begin{equation}
B(x_{2}\otimes 1_{H};G,g)-B(gx_{1}x_{2}\otimes 1_{H};G,gx_{1})=0.
\label{X2,gx1x2,X2F61,gx2}
\end{equation}

\paragraph{\textbf{Equality} $\left( \protect\ref{X2F71}\right) $}

rewrites as%
\begin{gather*}
B(gx_{1}x_{2}\otimes 1_{H};G,f)=-B(x_{1}x_{2}\otimes g;G,f) \\
-\lambda B(x_{1}x_{2}\otimes g;GX_{1}X_{2},f)+B(gx_{1}x_{2}\otimes
x_{2};GX_{2},f).
\end{gather*}

\subparagraph{Case $f=1_{H}$}

\begin{gather*}
B(gx_{1}x_{2}\otimes 1_{H};G,1_{H})=-B(x_{1}x_{2}\otimes g;G,1_{H}) \\
-\lambda B(x_{1}x_{2}\otimes g;GX_{1}X_{2},1_{H})+B(gx_{1}x_{2}\otimes
x_{2};GX_{2},1_{H}).
\end{gather*}%
By applying $\left( \ref{eq.10}\right) $ this rewrites as%
\begin{gather*}
2B(gx_{1}x_{2}\otimes 1_{H};G,1_{H}) \\
+\lambda B(gx_{1}x_{2}\otimes 1_{H};GX_{1}X_{2},1_{H})-B(gx_{1}x_{2}\otimes
x_{2};GX_{2},1_{H})=0.
\end{gather*}%
In view of the form of the elements we get%
\begin{gather}
2B(gx_{1}x_{2}\otimes 1_{H};G,1_{H})  \label{X2,gx1x2,X2F71,1H} \\
+\lambda \left[
\begin{array}{c}
B(g\otimes 1_{H};G,1_{H})+B(x_{2}\otimes \ 1_{H};G,x_{2}) \\
+B(x_{1}\otimes 1_{H};G,x_{1})+B(gx_{1}x_{2}\otimes 1_{H};G,x_{1}x_{2})%
\end{array}%
\right] =0.  \notag
\end{gather}

\subparagraph{Case $f=x_{1}x_{2}$}

\begin{gather*}
B(gx_{1}x_{2}\otimes 1_{H};G,x_{1}x_{2})=-B(x_{1}x_{2}\otimes g;G,x_{1}x_{2})
\\
+\lambda B(x_{1}x_{2}\otimes g;GX_{1}X_{2},x_{1}x_{2})+B(gx_{1}x_{2}\otimes
x_{2};GX_{2},x_{1}x_{2}).
\end{gather*}%
By applying $\left( \ref{eq.10}\right) $ this rewrites as%
\begin{gather*}
2B(gx_{1}x_{2}\otimes 1_{H};G,x_{1}x_{2}) \\
+\lambda B(gx_{1}x_{2}\otimes
1_{H};GX_{1}X_{2},x_{1}x_{2})-B(gx_{1}x_{2}\otimes
x_{2};GX_{2},x_{1}x_{2})=0.
\end{gather*}%
In view of the form of the elements we get%
\begin{equation*}
2B(gx_{1}x_{2}\otimes 1_{H};G,x_{1}x_{2})+\lambda B(g\otimes
1_{H};G,x_{1}x_{2})=0
\end{equation*}%
which is $\left( \ref{X1,gx1x2,X1F61,x1x2}\right) .$

\subparagraph{Case $f=gx_{1}$}

\begin{gather*}
B(gx_{1}x_{2}\otimes 1_{H};G,gx_{1})=-B(x_{1}x_{2}\otimes g;G,gx_{1}) \\
-\lambda B(x_{1}x_{2}\otimes g;GX_{1}X_{2},gx_{1})+B(gx_{1}x_{2}\otimes
x_{2};GX_{2},gx_{1}).
\end{gather*}%
By applying $\left( \ref{eq.10}\right) $ this rewrites as%
\begin{equation*}
-\lambda B(gx_{1}x_{2}\otimes 1_{H};GX_{1}X_{2},gx_{1})-B(gx_{1}x_{2}\otimes
x_{2};GX_{2},gx_{1})=0.
\end{equation*}%
In view of the form of the elements we get%
\begin{equation*}
\lambda \left[ B(g\otimes 1_{H};G,gx_{1})+B(x_{2}\otimes \
1_{H};G,gx_{1}x_{2})\right] =0
\end{equation*}%
which is $\left( \ref{X2,x2,X2F71,g}\right) .$

\subparagraph{Case $f=gx_{2}$}

\begin{gather*}
B(gx_{1}x_{2}\otimes 1_{H};G,gx_{2})=-B(x_{1}x_{2}\otimes g;G,gx_{2}) \\
-\lambda B(x_{1}x_{2}\otimes g;GX_{1}X_{2},gx_{2})+B(gx_{1}x_{2}\otimes
x_{2};GX_{2},gx_{2}).
\end{gather*}%
By applying $\left( \ref{eq.10}\right) $ this rewrites as%
\begin{equation*}
-\lambda B(gx_{1}x_{2}\otimes 1_{H};GX_{1}X_{2},gx_{2})-B(gx_{1}x_{2}\otimes
x_{2};GX_{2},gx_{2})=0.
\end{equation*}%
In view of the form of the elements we get%
\begin{gather*}
\lambda \left[ B(g\otimes 1_{H};G,gx_{2})-B(x_{1}\otimes 1_{H};G,gx_{1}x_{2})%
\right] + \\
2\left[ B(x_{1}\otimes 1_{H};G,g)+B(gx_{1}x_{2}\otimes 1_{H};G,gx_{2})\right]
=0.
\end{gather*}%
By means of $\left( \ref{X1,g,X1F21,gx1x2}\right) ,$we get%
\begin{equation*}
B(x_{1}\otimes 1_{H};G,g)+B(gx_{1}x_{2}\otimes 1_{H};G,gx_{2})=0
\end{equation*}%
which is $\left( \ref{X2,x1,X2F71,g}\right) .$

\paragraph{\textbf{Equality} $\left( \protect\ref{X2F81}\right) $}

rewrites as%
\begin{gather*}
B(gx_{1}x_{2}\otimes 1_{H};GX_{1},f)=B(x_{1}x_{2}\otimes g;GX_{1},f) \\
+B(gx_{1}x_{2}\otimes x_{2};GX_{1}X_{2},f).
\end{gather*}

\subparagraph{Case $f=g$}

\begin{gather*}
B(gx_{1}x_{2}\otimes 1_{H};GX_{1},g)=B(x_{1}x_{2}\otimes g;GX_{1},g) \\
+B(gx_{1}x_{2}\otimes x_{2};GX_{1}X_{2},g).
\end{gather*}%
By applying $\left( \ref{eq.10}\right) $ this rewrites as%
\begin{equation*}
B(gx_{1}x_{2}\otimes x_{2};GX_{1}X_{2},g)=0
\end{equation*}%
which is known in view of the form of the element.

\subparagraph{Case $f=x_{1}$}

\begin{gather*}
B(gx_{1}x_{2}\otimes 1_{H};GX_{1},x_{1})=B(x_{1}x_{2}\otimes g;GX_{1},x_{1})
\\
+B(gx_{1}x_{2}\otimes x_{2};GX_{1}X_{2},x_{1}).
\end{gather*}%
By applying $\left( \ref{eq.10}\right) $ this rewrites as%
\begin{equation*}
2B(gx_{1}x_{2}\otimes 1_{H};GX_{1},x_{1})-B(gx_{1}x_{2}\otimes
x_{2};GX_{1}X_{2},x_{1})=0
\end{equation*}%
In view of the form of the elements we get%
\begin{equation*}
B(x_{2}\otimes 1_{H};G,x_{1})=0
\end{equation*}%
which is $\left( \ref{X2,g,X2F21,x1}\right) .$

\subparagraph{Case $f=x_{2}$}

\begin{gather*}
B(gx_{1}x_{2}\otimes 1_{H};GX_{1},x_{2})=B(x_{1}x_{2}\otimes g;GX_{1},x_{2})
\\
+B(gx_{1}x_{2}\otimes x_{2};GX_{1}X_{2},x_{2}).
\end{gather*}%
By applying $\left( \ref{eq.10}\right) $ this rewrites as%
\begin{equation*}
2B(gx_{1}x_{2}\otimes 1_{H};GX_{1},x_{2})-B(gx_{1}x_{2}\otimes
x_{2};GX_{1}X_{2},x_{2})=0.
\end{equation*}%
In view of the form of the elements we get%
\begin{equation*}
B(x_{2}\otimes 1_{H};G,x_{2})+B(gx_{1}x_{2}\otimes 1_{H};G,x_{1}x_{2})=0
\end{equation*}%
which is $\left( \ref{X2,x1,X2F61,x2}\right) .$

\subparagraph{Case $f=gx_{1}x_{2}$}

\begin{gather*}
B(gx_{1}x_{2}\otimes 1_{H};GX_{1},gx_{1}x_{2})=B(x_{1}x_{2}\otimes
g;GX_{1},gx_{1}x_{2}) \\
+B(gx_{1}x_{2}\otimes x_{2};GX_{1}X_{2},gx_{1}x_{2}).
\end{gather*}%
By applying $\left( \ref{eq.10}\right) $ this rewrites as%
\begin{equation*}
B(gx_{1}x_{2}\otimes x_{2};GX_{1}X_{2},gx_{1}x_{2})=0.
\end{equation*}%
In view of the form of the elements we get%
\begin{equation*}
B(g\otimes 1_{H};G,gx_{1})+B(x_{2}\otimes \ 1_{H};G,gx_{1}x_{2})=0
\end{equation*}%
which is $\left( \ref{X2,g,X2F71,gx1}\right) .$

\section{LIST\ OF\ ALL EQUALITIES GOT FROM $\left( \protect\ref{eq.h}\right)
$\label{LAE}}

\subsection{$G$}

\begin{equation*}
\gamma _{1}B\left( g\otimes 1_{H};1_{A},x_{1}\right) +\gamma _{2}B\left(
g\otimes 1_{H};1_{A},x_{2}\right) =0\text{ }\left( \ref{G,g,GF1,1H}\right)
\end{equation*}

\begin{equation*}
2\alpha B(g\otimes 1_{H};G,gx_{1})+\gamma _{2}B(g\otimes
1_{H};1_{A},gx_{1}x_{2})=0\left( \ref{G,g, GF1,gx1}\right)
\end{equation*}

\begin{equation*}
2\alpha B(g\otimes 1_{H};G,gx_{2})-\gamma _{1}B\left( g\otimes
1_{H};1_{A},gx_{1}x_{2}\right) =0.\left( \ref{G,g, GF1,gx2}\right)
\end{equation*}

\begin{equation*}
\gamma _{1}B\left( g\otimes 1_{H};G,gx_{1}\right) +\gamma _{2}B\left(
g\otimes 1_{H};G,gx_{2}\right) =0.\left( \ref{G,g, GF2,g}\right)
\end{equation*}%
\begin{equation*}
2B(g\otimes 1_{H};1_{A},x_{1})-\gamma _{2}B\left( g\otimes
1_{H};G,x_{1}x_{2}\right) =0.\left( \ref{G,g, GF2,x1}\right)
\end{equation*}%
\begin{equation*}
2B(g\otimes 1_{H};1_{A},x_{2})+\gamma _{1}B(g\otimes
1_{H};G,x_{1}x_{2})=0\left( \ref{G,g, GF2,x2}\right)
\end{equation*}%
\begin{equation*}
\begin{array}{c}
2\alpha B(x_{1}\otimes 1_{H};G,g)+ \\
+\gamma _{1}\left[ -B(g\otimes 1_{H};1_{A},g)-B(x_{1}\otimes
1_{H};1_{A},gx_{1})\right] -\gamma _{2}B(x_{1}\otimes 1_{H};1_{A},gx_{2})%
\end{array}%
=0.\left( \ref{G,x1, GF1,g}\right)
\end{equation*}%
\begin{equation*}
\gamma _{1}B(g\otimes 1_{H};1_{A},x_{1})+\gamma _{2}B(g\otimes
1_{H};1_{A},x_{2})=0.\left( \ref{G,x1, GF1,x1}\right)
\end{equation*}%
This is $\left( \ref{G,g,GF1,1H}\right) $%
\begin{equation*}
\gamma _{1}[-B(g\otimes 1_{H};1_{A},x_{2})+B(x_{1}\otimes
1_{H};1_{A},x_{1}x_{2})]=0.\left( \ref{G,x1, GF1,x2}\right)
\end{equation*}

\begin{equation*}
-\gamma _{1}B(g\otimes 1_{H};1_{A},gx_{1}x_{2})+2\alpha B(x_{1}\otimes
1_{H};G,gx_{1}x_{2})=0.\left( \ref{G,x1, GF1,gx1x2}\right)
\end{equation*}%
\begin{equation*}
\begin{array}{c}
2B(x_{1}\otimes 1_{H};1_{A},1_{H})+ \\
+\gamma _{1}\left[ B(g\otimes 1_{H};G,1_{H})+B(x_{1}\otimes 1_{H};G,x_{1})%
\right] +\gamma _{2}B(x_{1}\otimes 1_{H};G,x_{2})%
\end{array}%
=0.\left( \ref{G,x1, GF2,1H}\right)
\end{equation*}%
\begin{equation*}
2B(x_{1}\otimes 1_{H};1_{A},x_{1}x_{2})+\gamma _{1}B(g\otimes
1_{H};G,x_{1}x_{2})=0\text{ }\left( \ref{G,x1, GF2,x1x2}\right)
\end{equation*}%
\begin{equation*}
\gamma _{1}B\left( g\otimes 1_{H};G,gx_{1}\right) +\gamma _{2}B(x_{1}\otimes
1_{H};G,gx_{1}x_{2})=0.\left( \ref{G,x1, GF2,gx1}\right)
\end{equation*}

\begin{equation*}
\gamma _{1}[B\left( g\otimes 1_{H};G,gx_{2}\right) -B(x_{1}\otimes
1_{H};G,gx_{1}x_{2})]=0.\left( \ref{G,x1,GF2,gx2}\right)
\end{equation*}

\begin{equation*}
\gamma _{2}\left[ -B(g\otimes 1_{H};X_{2},1_{H})+B(x_{1}\otimes
1_{H};1_{A},x_{1}x_{2})\right] =0\left( \ref{G,x1,GF3,1H}\right)
\end{equation*}%
\begin{equation*}
-\gamma _{2}B(g\otimes 1_{H};1_{A},gx_{1}x_{2})-2\alpha B(g\otimes
1_{H};G,gx_{1})=0.\left( \ref{G,x1,GF3,gx1}\right)
\end{equation*}

\begin{equation*}
\alpha \left[ -B(g\otimes 1_{H};g,gx_{2})+B(x_{1}\otimes 1_{H};g,gx_{1}x_{2})%
\right] =0.\left( \ref{G,x1, GF3,x2}\right)
\end{equation*}%
\begin{equation*}
\gamma _{1}[B(g\otimes 1_{H};1_{A},x_{2})+B(x_{1}\otimes
1_{H};1_{A},x_{1}x_{2})]=0.\left( \ref{G,x1, GF4,1H}\right)
\end{equation*}

\begin{equation*}
-\gamma _{1}B(g\otimes 1_{H};1_{A},gx_{1}x_{2})+2\alpha B(x_{1}\otimes
1_{H};G,gx_{1}x_{2})=0.\left( \ref{G,x1, GF4,gx1}\right)
\end{equation*}

\begin{equation*}
\alpha \lbrack -B(g\otimes 1_{H};g,gx_{2})+B(x_{1}\otimes
1_{H};G,gx_{1}x_{2})]=0.\left( \ref{G,x1, GF5,g}\right)
\end{equation*}%
\begin{equation*}
2B(g\otimes 1_{H};1_{A},x_{1})-\gamma _{2}B(g\otimes
1_{H};G,x_{1}x_{2})=0.\left( \ref{G,x1,GF6,x1}\right)
\end{equation*}%
\begin{equation*}
B(g\otimes 1_{H};1_{A},x_{2})=B(x_{1}\otimes 1_{H};1_{A},x_{1}x_{2}).\left( %
\ref{G,x1,GF6,x2}\right)
\end{equation*}%
\begin{equation*}
\begin{array}{c}
2\alpha B(x_{2}\otimes 1_{H};G,g)+ \\
-\gamma _{1}B(x_{2}\otimes 1_{H};1_{A},gx_{1})+\gamma _{2}\left[ -B(g\otimes
1_{H};1_{A},g)-B(x_{2}\otimes \ 1_{H};1_{A},gx_{2})\right]%
\end{array}%
=0.\left( \ref{G,x2, GF1,g}\right)
\end{equation*}%
\begin{equation*}
\gamma _{2}\left[ -B(g\otimes 1_{H};1_{A},x_{1})-B(x_{2}\otimes \
1_{H};1_{A},x_{1}x_{2})\right] =0\left( \ref{G,x2, GF1,x1}\right)
\end{equation*}

\begin{equation*}
\gamma _{1}B(x_{2}\otimes 1_{H};1_{A},x_{1}x_{2})-\gamma _{2}B(g\otimes
1_{H};1_{A},x_{2})=0\left( \ref{G,x2, GF1,x2}\right)
\end{equation*}

\begin{equation*}
2\alpha B(x_{2}\otimes 1_{H};G,gx_{1}x_{2})-\gamma _{2}B\left( g\otimes
1_{H};1_{A},gx_{1}x_{2}\right) =0.\left( \ref{G,x2, GF1,gx1x2}\right)
\end{equation*}%
\begin{equation*}
\begin{array}{c}
2B(x_{2}\otimes 1_{H};1_{A},1_{H}) \\
+\gamma _{1}B(x_{2}\otimes 1_{H};G,x_{1})+\gamma _{2}\left[ B(g\otimes
1_{H};G,1_{H})+B(x_{2}\otimes \ 1_{H};G,x_{2}\right]%
\end{array}%
=0.\left( \ref{G,x2, GF2,1H}\right)
\end{equation*}%
\begin{equation*}
2B(x_{2}\otimes 1_{H};1_{A},x_{1}x_{2})+\gamma _{2}B\left( g\otimes
1_{H};G,x_{1}x_{2}\right) =0=0\left( \ref{G,x2, GF2,x1x2}\right)
\end{equation*}%
\begin{equation*}
\gamma _{2}\left[ B(g\otimes 1_{H};G,gx_{1})+B(x_{2}\otimes \
1_{H};G,gx_{1}x_{2})\right] =0.\left( \ref{G,x2, GF2,gx1}\right)
\end{equation*}

\begin{equation*}
-\gamma _{1}B(x_{2}\otimes 1_{H};G,gx_{1}x_{2})+\gamma _{2}B(g\otimes
1_{H};G,gx_{2})=0.\left( \ref{G,x2, GF2,gx2}\right)
\end{equation*}%
\begin{equation*}
\gamma _{2}[B(x_{2}\otimes 1_{H};1_{A},x_{1}x_{2})+B(g\otimes
1_{H};1_{A},x_{1})]=0\left( \ref{G,x2, FG3,1H}\right)
\end{equation*}

\begin{equation*}
-\gamma _{2}B(g\otimes 1_{H};1_{A},gx_{1}x_{2})-2\alpha B(x_{2}\otimes
1_{H};G,gx_{1}x_{2})=0.\left( \ref{G,x2, GF3,gx2}\right)
\end{equation*}%
\begin{equation*}
\gamma _{1}\left[ B(g\otimes 1_{H};1_{A},x_{1})+B(x_{2}\otimes \
1_{H};1_{A},x_{1}x_{2})\right] =0.\left( \ref{G,x2, GF4,1H}\right)
\end{equation*}%
\begin{equation*}
\alpha \left[ B(g\otimes 1_{H};G,gx_{1})+B(x_{2}\otimes \
1_{H};G,gx_{1}x_{2})\right] =0\left( \ref{G,x2, GF4,gx1}\right)
\end{equation*}

\begin{equation*}
\gamma _{1}B(g\otimes 1_{H};1_{A},gx_{1}x_{2})-2\alpha B(g\otimes
1_{H};G,gx_{2})=0.\left( \ref{G,x2, GF4,gx2}\right)
\end{equation*}%
\begin{equation*}
\alpha \left[ B(g\otimes 1_{H};G,gx_{1})+B(x_{2}\otimes \
1_{H};G,gx_{1}x_{2})\right] =0\left( \ref{G,x2, GF5,g}\right)
\end{equation*}

\begin{equation*}
\gamma _{2}[B(g\otimes 1_{H},G,gx_{1})+B(x_{2}\otimes
1_{H};G,gx_{1}x_{2})=0.=0.\left( \ref{G,x2, GF6,g}\right)
\end{equation*}%
\begin{equation*}
2B(x_{2}\otimes 1_{H};1_{A},x_{1}x_{2})+\gamma _{2}B(g\otimes
1_{H};G,x_{1}x_{2})=0.\left( \ref{G,x2; GF6,x2}\right)
\end{equation*}

\begin{equation*}
\gamma _{1}[B(g\otimes 1_{H};G,gx_{1})+B(x_{2}\otimes
1_{H};G,gx_{1}x_{2})]=0.\left( \ref{G,x2; GF7,g}\right)
\end{equation*}%
\begin{equation*}
+B(g\otimes 1_{H};1_{A},x_{1})+B(x_{2}\otimes \
1_{H};1_{A},x_{1}x_{2})=0.\left( \ref{G,x2; GF7,x1}\right)
\end{equation*}%
\begin{equation*}
\gamma _{1}B(x_{1}x_{2}\otimes 1_{H};X_{1},g)+\gamma _{2}B(x_{1}x_{2}\otimes
1_{H};X_{2},g)=0.\left( \ref{G,x1x2, GF1,g}\right)
\end{equation*}

\begin{equation*}
\begin{array}{c}
2\alpha B(x_{1}x_{2}\otimes 1_{H};G,x_{1})+ \\
+\gamma _{1}B(x_{1}x_{2}\otimes 1_{H};X_{1},x_{1})+\gamma
_{2}B(x_{1}x_{2}\otimes 1_{H};X_{2},x_{1})%
\end{array}%
=0.\left( \ref{G,x1x2, GF1,x1}\right)
\end{equation*}%
\begin{equation*}
\begin{array}{c}
2\alpha B(x_{1}x_{2}\otimes 1_{H};G,x_{2})+ \\
+\gamma _{1}B(x_{1}x_{2}\otimes 1_{H};X_{1},x_{2})+\gamma
_{2}B(x_{1}x_{2}\otimes 1_{H};X_{2},x_{2})%
\end{array}%
=0.\left( \ref{G,x1x2, GF1,x2}\right)
\end{equation*}%
\begin{equation*}
\gamma _{1}B(x_{1}x_{2}\otimes 1_{H};X_{1},gx_{1}x_{2})+\gamma
_{2}B(x_{1}x_{2}\otimes 1_{H};X_{2},gx_{1}x_{2})=0\text{ }\left( \ref%
{G,x1x2, GF1,gx1x2}\right)
\end{equation*}

\begin{equation*}
\gamma _{1}B(x_{1}x_{2}\otimes 1_{H};GX_{1},1_{H})+\gamma
_{2}B(x_{1}x_{2}\otimes 1_{H};GX_{2},1_{H})=0\left( \ref{G,x1x2, GF2,1H}%
\right)
\end{equation*}%
\begin{equation*}
\gamma _{1}B(x_{1}x_{2}\otimes 1_{H};GX_{1},x_{1}x_{2})+\gamma
_{2}B(x_{1}x_{2}\otimes 1_{H};GX_{2},x_{1}x_{2})=0\left( \ref{G,x1x2,
GF2,x1x2}\right)
\end{equation*}%
\begin{equation*}
\begin{array}{c}
2B(x_{1}x_{2}\otimes 1_{H};1_{A},gx_{1}) \\
+\gamma _{1}B(x_{1}x_{2}\otimes 1_{H};GX_{1},gx_{1})+\gamma
_{2}B(x_{1}x_{2}\otimes 1_{H};GX_{2},gx_{1})%
\end{array}%
=0.\left( \ref{G,x1x2, GF2,gx1}\right)
\end{equation*}%
\begin{equation*}
\begin{array}{c}
2B(x_{1}x_{2}\otimes 1_{H};1_{A},gx_{2}) \\
+\gamma _{1}B(x_{1}x_{2}\otimes 1_{H};GX_{1},gx_{2})+\gamma
_{2}B(x_{1}x_{2}\otimes 1_{H};GX_{2},gx_{2})%
\end{array}%
=0\left( \ref{G,x1x2, GF2,gx2}\right)
\end{equation*}

\begin{equation*}
2\alpha B(x_{1}x_{2}\otimes 1_{H};GX_{1},1_{H})=0.\left( \ref{G,x1x2, GF3,1H}%
\right)
\end{equation*}

\begin{equation*}
\alpha B(x_{1}x_{2}\otimes 1_{H};GX_{1},x_{1}x_{2})=0.\left( \ref{G,x1x2,
GF3,x1x2}\right)
\end{equation*}%
\begin{equation*}
-2\alpha B(x_{1}x_{2}\otimes 1_{H};GX_{1},x_{1}x_{2})=0\text{ }\left( \ref%
{G,x1x2, GF3,gx1}\right)
\end{equation*}%
\begin{equation*}
\gamma _{2}B(x_{1}x_{2}\otimes 1_{H};X_{2},gx_{1}x_{2})=0\text{ }\left( \ref%
{G,x1x2, GF3,gx2}\right)
\end{equation*}%
\begin{equation*}
-\gamma _{1}[1-B(x_{1}x_{2}\otimes
1_{H};1_{A},x_{1}x_{2})+B(x_{1}x_{2}\otimes
1_{H};x_{1},x_{2})-B(x_{1}x_{2}\otimes 1_{H};x_{2},x_{1})]-2\alpha
B(x_{1}x_{2}\otimes 1_{H};GX_{2},1_{H}=0.\left( \ref{G,x1x2, GF4,1H}\right)
\end{equation*}

\begin{equation*}
\alpha B(x_{1}x_{2}\otimes 1_{H};GX_{2},x_{1}x_{2})=0.\left( \ref{G,x1x2,
GF4,x1x2}\right)
\end{equation*}%
\begin{equation*}
\gamma _{1}B(x_{1}x_{2}\otimes 1_{H};X_{1},gx_{1}x_{2})=0\text{ }.\left( \ref%
{G,x1x2, GF4,gx1}\right)
\end{equation*}%
\begin{equation*}
\gamma _{1}B(x_{1}x_{2}\otimes 1_{H};X_{2},gx_{1}x_{2})=0.\left( \ref%
{G,x1x2, GF4,gx2}\right)
\end{equation*}%
\begin{equation*}
-2B(x_{1}x_{2}\otimes 1_{H};X_{1},g)+\gamma _{2}[-B(x_{1}x_{2}\otimes
1_{H};g,gx_{1}x_{2})-B(x_{1}x_{2}\otimes
1_{H};GX_{1},gx_{2})+B(x_{1}x_{2}\otimes 1_{H};GX_{2},gx_{1})]=0\left( \ref%
{G,x1x2, GF6,g}\right)
\end{equation*}%
\begin{equation*}
\gamma _{2}B(x_{1}x_{2}\otimes 1_{H};GX_{1},x_{1}x_{2})=0.\left( \ref%
{G,x1x2, GF6,x1}\right)
\end{equation*}%
\begin{equation*}
\gamma _{2}B(x_{1}x_{2}\otimes 1_{H};GX_{2},x_{1}x_{2})=0.\left( \ref%
{G,x1x2, GF6,x2}\right)
\end{equation*}%
\begin{equation*}
B(x_{1}x_{2}\otimes 1_{H};X_{1},gx_{1}x_{2})=0\text{ }\left( \ref{G,x1x2,
GF6,gx1x2}\right)
\end{equation*}%
\begin{gather*}
-2B(x_{1}x_{2}\otimes 1_{H};X_{2},g)+ \\
-\gamma _{1}[-B(x_{1}x_{2}\otimes 1_{H};g,gx_{1}x_{2})-B(x_{1}x_{2}\otimes
1_{H};GX_{1},gx_{2})+B(x_{1}x_{2}\otimes 1_{H};GX_{2},gx_{1})]=0\left( \ref%
{G,x1x2, GF7,g}\right)
\end{gather*}%
\begin{equation*}
\gamma _{1}B(x_{1}x_{2}\otimes 1_{H};GX_{1},x_{1}x_{2})=0\left( \ref{G,x1x2,
GF7,x1}\right)
\end{equation*}

\begin{equation*}
\gamma _{1}B(x_{1}x_{2}\otimes 1_{H};GX_{2},x_{1}x_{2})=0\left( \ref{G,x1x2,
GF7,x2}\right)
\end{equation*}%
\begin{equation*}
B(x_{1}x_{2}\otimes 1_{H};X_{2},gx_{1}x_{2})=0\left( \ref{G,x1x2, GF7,gx1x2}%
\right)
\end{equation*}%
\begin{equation*}
\begin{array}{c}
2\alpha \left[ B(x_{1}x_{2}\otimes 1_{H};G,x_{2})-B(x_{1}x_{2}\otimes
1_{H};GX_{2},1_{H})\right] + \\
+\gamma _{1}\left[ -1+B(x_{1}x_{2}\otimes
1_{H};1_{A},x_{1}x_{2})+B(x_{1}x_{2}\otimes 1_{H};X_{2},x_{1})\right]
+\gamma _{2}B(x_{1}x_{2}\otimes 1_{H};X_{2},x_{2})%
\end{array}%
=0.\left( \ref{G,gx1, GF1,1H}\right)
\end{equation*}

\begin{equation*}
\begin{array}{c}
2\left[ B(x_{1}x_{2}\otimes 1_{H};1_{A},gx_{2})+B(x_{1}x_{2}\otimes
1_{H};X_{2},g)\right] \\
\gamma _{1}\left[ -B(x_{1}x_{2}\otimes
1_{H};G,gx_{1}x_{2})+B(x_{1}x_{2}\otimes 1_{H};GX_{2},gx_{1})\right] +\gamma
_{2}B(x_{1}x_{2}\otimes 1_{H};GX_{2},gx_{2})%
\end{array}%
\left( \ref{G,gx1, GF2,g}\right)
\end{equation*}%
\begin{equation*}
\begin{array}{c}
2\alpha \left[ -B(x_{1}x_{2}\otimes 1_{H};G,x_{1})+B(x_{1}x_{2}\otimes
1_{H};GX_{1},1_{H})\right] + \\
-\gamma _{1}B(x_{1}x_{2}\otimes 1_{H};X_{1},x_{1})+\gamma _{2}\left[
-1+B(x_{1}x_{2}\otimes 1_{H};1_{A},x_{1}x_{2})-B(x_{1}x_{2}\otimes
1_{H};X_{1},x_{2})\right]%
\end{array}%
=0\left( \ref{G,gx2, GF1,1H}\right)
\end{equation*}%
\begin{equation*}
\begin{array}{c}
2\left[ -B(x_{1}x_{2}\otimes 1_{H};1_{A},gx_{1})-B(x_{1}x_{2}\otimes
1_{H};X_{1},g)\right] \\
-\gamma _{1}B(x_{1}x_{2}\otimes 1_{H};GX_{1},gx_{1})+\gamma _{2}\left[
-B(x_{1}x_{2}\otimes 1_{H};G,gx_{1}x_{2})-B(x_{1}x_{2}\otimes
1_{H};GX_{1},gx_{2})\right]%
\end{array}%
=0.\left( \ref{G,gx2, GF2,g}\right)
\end{equation*}%
\begin{equation*}
\gamma _{1}B(x_{1}x_{2}\otimes 1_{H};GX_{1},x_{1}x_{2})=0\left( \ref{G,gx2,
GF2,x2}\right)
\end{equation*}

\begin{equation*}
\gamma _{1}\left[ -B(x_{2}\otimes 1_{H};1_{A},1_{H})+B(gx_{1}x_{2}\otimes
1_{H};1_{A},x_{1})\right] +\gamma _{2}\left[ B(x_{1}\otimes
1_{H};1_{A},1_{H})+B(gx_{1}x_{2}\otimes 1_{H};1_{A},x_{2})\right] =0\left( %
\ref{G,gx1x2, GF1,1H}\right)
\end{equation*}%
\begin{equation*}
-\gamma _{1}B(x_{2}\otimes 1_{H};1_{A},x_{1}x_{2})+\gamma _{2}B(x_{1}\otimes
1_{H};1_{A},x_{1}x_{2})=0.\left( \ref{G,gx1x2, GF1,x1x2}\right)
\end{equation*}%
\begin{equation*}
\begin{array}{c}
2\alpha B(x_{2}\otimes 1_{H};G,g)+ \\
-\gamma _{1}B(x_{2}\otimes 1_{H};1_{A},gx_{1})+\gamma _{2}\left[
B(x_{1}\otimes 1_{H};1_{A},gx_{1})+B(gx_{1}x_{2}\otimes
1_{H};1_{A},gx_{1}x_{2})\right]%
\end{array}%
=0.\left( \ref{G,gx1x2, GF1,gx1}\right)
\end{equation*}

\begin{equation*}
\begin{array}{c}
2\alpha B(gx_{1}x_{2}\otimes 1_{H};G,gx_{2})+ \\
+\gamma _{1}\left[ -B(x_{2}\otimes 1_{H};1_{A},gx_{2})-B(gx_{1}x_{2}\otimes
1_{H};1_{A},gx_{1}x_{2})\right] +\gamma _{2}B(x_{1}\otimes
1_{H};1_{A},gx_{2})%
\end{array}%
=0.\left( \ref{G,gx1x2, GF1,gx2}\right)
\end{equation*}%
\begin{gather*}
\gamma _{1}\left[ B(x_{2}\otimes 1_{H};G,g)-B(gx_{1}x_{2}\otimes
1_{H};G,gx_{1})\right] \left( \ref{G,gx1x2, GF2,g}\right) \\
+\gamma _{2}\left[ -B(x_{1}\otimes 1_{H};G,g)-B(gx_{1}x_{2}\otimes
1_{H};G,gx_{2})\right] =0.
\end{gather*}%
\begin{equation*}
\begin{array}{c}
2B(gx_{1}x_{2}\otimes 1_{H};1_{A},x_{1}) \\
+\gamma _{1}B(x_{2}\otimes 1_{H};G,x_{1})+\gamma _{2}\left[ -B(x_{1}\otimes
1_{H};G,x_{1})-B(gx_{1}x_{2}\otimes 1_{H};G,x_{1}x_{2})\right]%
\end{array}%
=0.\left( \ref{G,gx1x2, GF2,x1}\right)
\end{equation*}%
\begin{equation*}
\begin{array}{c}
2B(gx_{1}x_{2}\otimes 1_{H};1_{A},x_{2}) \\
+\gamma _{1}\left[ B(x_{2}\otimes 1_{H};G,x_{2})+B(gx_{1}x_{2}\otimes
1_{H};G,x_{1}x_{2})\right] -\gamma _{2}B(x_{1}\otimes 1_{H};G,x_{2})%
\end{array}%
=0.\left( \ref{G,gx1x2, GF2,x2}\right)
\end{equation*}%
\begin{equation*}
\gamma _{1}B(x_{2}\otimes 1_{H};G,gx_{1}x_{2})-\gamma _{2}B(x_{1}\otimes
1_{H};G,gx_{1}x_{2})=0.\left( \ref{G,gx1x2, GF2,gx1x2}\right)
\end{equation*}%
\begin{gather*}
\gamma _{2}\left[ B(g\otimes 1_{H};1_{A},g)+B(x_{2}\otimes \
1_{H};1_{A},gx_{2})+B(x_{1}\otimes 1_{H};1_{A},gx_{1})+B(gx_{1}x_{2}\otimes
1_{H};1_{A},gx_{1}x_{2})\right] +\left( \ref{G,gx1x2, GF3,g}\right) \\
2\alpha \left[ -B(x_{2}\otimes 1_{H};G,g)+B(x_{2}\otimes 1_{H};G,g)\right]
=0.
\end{gather*}

\begin{equation*}
\gamma _{2}\left[ B(g\otimes 1_{H};1_{A},x_{1})+B(x_{2}\otimes \
1_{H};1_{A},x_{1}x_{2})\right] =0.\left( \ref{G,gx1x2, GF3,x1}\right)
\end{equation*}

\begin{equation*}
\gamma _{2}\left[ B(g\otimes 1_{H};X_{2},1_{H})-B(x_{1}\otimes
1_{H};1_{A},x_{1}x_{2})\right] =0\left( \ref{G,gx1x2, GF3,x2}\right)
\end{equation*}

\begin{equation*}
\gamma _{2}B\left( g\otimes 1_{H};1_{A},gx_{1}x_{2}\right) -2\alpha
B(x_{2}\otimes 1_{H};G,gx_{1}x_{2})=0.\left( \ref{G,gx1x2, GF3,gx1x2}\right)
\end{equation*}

\begin{gather*}
-\gamma _{1}\left[ B(g\otimes 1_{H};1_{A},g)+B(x_{2}\otimes \
1_{H};1_{A},gx_{2})+B(x_{1}\otimes 1_{H};1_{A},gx_{1})+B(gx_{1}x_{2}\otimes
1_{H};1_{A},gx_{1}x_{2})\right] +\left( \ref{G,gx1x2, GF4,g}\right) \\
-2\alpha \left[ -B(x_{1}\otimes 1_{H};G,g)-B(gx_{1}x_{2}\otimes
1_{H};G,gx_{2})\right] =0
\end{gather*}%
\begin{equation*}
\gamma _{1}\left[ B(g\otimes 1_{H};1_{A},x_{1})+B(x_{2}\otimes \
1_{H};1_{A},x_{1}x_{2})\right] =0\left( \ref{G,gx1x2, GF4,x1}\right)
\end{equation*}%
\begin{equation*}
-\gamma _{1}\left[ B(g\otimes 1_{H};1_{A},x_{2})-B(x_{1}\otimes
1_{H};1_{A},x_{1}x_{2})\right] =0.\left( \ref{G,gx1x2, GF4,x2}\right)
\end{equation*}

\begin{equation*}
-\gamma _{1}B\left( g\otimes 1_{H};1_{A},gx_{1}x_{2}\right) +2\alpha
B(x_{1}\otimes 1_{H};G,gx_{1}x_{2})=0.\left( \ref{G,gx1x2, GF4,gx1x2}\right)
\end{equation*}

\begin{equation*}
\alpha \left[ B(g\otimes 1_{H};G,gx_{2})-B(x_{1}\otimes 1_{H};G,gx_{1}x_{2})%
\right] =0.\left( \ref{G,gx1x2, GF5,gx2}\right)
\end{equation*}

\begin{gather*}
-2\left[ -B(x_{2}\otimes 1_{H};1_{A},1_{H})+B(gx_{1}x_{2}\otimes
1_{H};1_{A},x_{1})\right] +\left( \ref{G,gx1x2, GF6,1H}\right) \\
\gamma _{2}\left[ B(g\otimes 1_{H};G,1_{H})+B(x_{2}\otimes \
1_{H};G,x_{2})+B(x_{1}\otimes 1_{H};G,x_{1})+B(gx_{1}x_{2}\otimes
1_{H};G,x_{1}x_{2})\right] \\
=0
\end{gather*}%
\begin{equation*}
\gamma _{2}\left[ B(g\otimes 1_{H};G,gx_{2})-B(x_{1}\otimes
1_{H};G,gx_{1}x_{2})\right] =0.\left( \ref{G,gx1x2, GF6,gx1}\right)
\end{equation*}%
\begin{equation*}
\gamma _{2}\left[ B(g\otimes 1_{H};G,gx_{2})-B(x_{1}\otimes
1_{H};G,gx_{1}x_{2})\right] =0.\left( \ref{G,gx1x2, GF6,gx2}\right)
\end{equation*}%
\begin{gather*}
-\gamma _{1}\left[ B(g\otimes 1_{H};G,1_{H})+B(x_{2}\otimes \
1_{H};G,x_{2})+B(x_{1}\otimes 1_{H};G,x_{1})+B(gx_{1}x_{2}\otimes
1_{H};G,x_{1}x_{2})\right] \\
-2\left[ B(x_{1}\otimes 1_{H};1_{A},1_{H})+B(gx_{1}x_{2}\otimes
1_{H};1_{A},x_{2})\right] =0\left( \ref{G,gx1x2, GF7,1H}\right)
\end{gather*}%
\begin{equation*}
-\gamma _{1}B(g\otimes 1_{H};G,x_{1}x_{2})-2B(x_{1}\otimes
1_{H};1_{A},x_{1}x_{2})=0.\left( \ref{G,gx1x2, GF7,x1x2}\right)
\end{equation*}%
\begin{equation*}
\gamma _{1}\left[ B(g\otimes 1_{H};G,gx_{1})+B(x_{2}\otimes \
1_{H};G,gx_{1}x_{2})\right] =0.\left( \ref{G,gx1x2, GF7,gx1}\right)
\end{equation*}%
\begin{equation*}
\gamma _{1}\left[ B(g\otimes 1_{H};G,gx_{2})-B(x_{1}\otimes
1_{H};G,gx_{1}x_{2})\right] =0.\left( \ref{G,gx1x2, GF7,gx2}\right)
\end{equation*}

\subsection{$X_{1}$}

\begin{equation*}
\lambda B\left( g\otimes 1_{H};1_{A},x_{2}\right) -\gamma _{1}B(g\otimes
1_{H};G,1_{H})-2B(x_{1}\otimes 1_{H};1_{A},1_{H})=0.\left( \ref{X1,g,
X1F11,1H}\right)
\end{equation*}%
\begin{eqnarray*}
&&+\gamma _{1}B(g\otimes 1_{H};G,gx_{1})+\lambda B(g\otimes
1_{H};X_{2},gx_{1})+B(x_{1}\otimes 1_{H};1_{A},gx_{1})+\left( \ref{X1,g,
X1F11,gx1}\right) \\
&&+2B\left( g\otimes 1_{H};1_{A},g\right) +B\left( x_{1}\otimes
1_{H};1_{A},gx_{1}\right) =0
\end{eqnarray*}

\begin{equation*}
\gamma _{1}B(g\otimes 1_{H};G,gx_{2})+2B(x_{1}\otimes
1_{H};1_{A},gx_{2})-2\beta _{1}B\left( g\otimes
1_{H};1_{A},gx_{1}x_{2}\right) =0.\left( \ref{X1,g, X1F11,gx2}\right)
\end{equation*}

\begin{equation*}
2\beta _{1}B\left( g\otimes 1_{H};G,gx_{1}\right) +\lambda B\left( g\otimes
1_{H};G,gx_{2}\right) +2B(x_{1}\otimes 1_{H};G,g)=0.\left( \ref{X1,g,X1F21,g}%
\right)
\end{equation*}%
\begin{equation*}
-\lambda B\left( g\otimes 1_{H};G,x_{1}x_{2}\right) +2B(x_{1}\otimes
1_{H};G,x_{1})=0.\left( \ref{X1,g,X1F21,x1}\right)
\end{equation*}%
\begin{equation*}
B(x_{1}\otimes 1_{H};G,x_{2})=0.\left( \ref{X1,g,X1F21,x2}\right)
\end{equation*}

\begin{equation*}
B\left( g\otimes 1_{H};G,gx_{2}\right) -B\left( x_{1}\otimes
1_{H};G,gx_{1}x_{2}\right) =0\left( \ref{X1,g,X1F21,gx1x2}\right)
\end{equation*}%
\begin{gather*}
+\lambda B\left( g\otimes 1_{H};1_{A},gx_{1}x_{2}\right) +\gamma _{1}B\left(
g\otimes 1_{H};G,gx_{1}\right) +\left( \ref{X1,g, X1F31,g}\right) \\
+2B(g\otimes 1_{H};1_{A},g)+2B(x_{1}\otimes 1_{H};1_{A},gx_{1})=0.
\end{gather*}

\begin{equation*}
-2\beta _{1}B\left( g\otimes 1_{H};1_{A},gx_{1}x_{2}\right) +2B(x_{1}\otimes
1_{H};1_{A},gx_{2})+\gamma _{1}B\left( g\otimes 1_{H};G,gx_{2}\right)
=0.\left( \ref{X1,g, X1F41,g}\right)
\end{equation*}%
\begin{equation*}
2B(x_{1}\otimes 1_{H};1_{A},x_{1}x_{2})+\gamma _{1}B\left( g\otimes
1_{H};G,x_{1}x_{2}\right) =0.\left( \ref{X1,g, X1F41,x1}\right)
\end{equation*}

\begin{equation*}
-2B\left( g\otimes 1_{H};1_{A},x_{2}\right) =\gamma _{1}B\left( g\otimes
1_{H};G,x_{1}x_{2}\right) -2B\left( g\otimes 1_{H};1_{A},x_{2}\right)
+2B(x_{1}\otimes 1_{H};1_{A},x_{1}x_{2})\left( \ref{X1,g,X1F51,1H}\right)
\end{equation*}%
\begin{equation*}
B(x_{1}\otimes 1_{H};G,gx_{1}x_{2})-B(g\otimes 1_{H};G,gx_{2})=0.\left( \ref%
{X1,x2,X1F71,gx1x2}\right)
\end{equation*}%
\begin{equation*}
2\beta _{1}\left[ -B(g\otimes 1_{H};1_{A},g)-B(x_{1}\otimes
1_{H};1_{A},gx_{1})\right] -\lambda B(x_{1}\otimes
1_{H};1_{A},gx_{2})=\gamma _{1}B(x_{1}\otimes 1_{H};G,g)\left( \ref%
{X1,x1,X1F11,g}\right)
\end{equation*}%
\begin{equation*}
-\lambda B(x_{1}\otimes 1_{H};1_{A},x_{1}x_{2})+\gamma _{1}B(x_{1}\otimes
1_{H};G,x_{1})=0.\left( \ref{X1,x1,X1F11,x1}\right)
\end{equation*}%
\begin{equation*}
0=\gamma _{1}B(x_{1}\otimes 1_{H};G,x_{2}).\left( \ref{X1,x1,X1F11,x2}\right)
\end{equation*}

\begin{equation*}
-2\beta _{1}B(g\otimes 1_{H};1_{A},gx_{1}x_{2})=\gamma _{1}B(x_{1}\otimes
1_{H};G,gx_{1}x_{2})-2B\left( x_{1}\otimes 1_{H};1_{A},gx_{2}\right) .\left( %
\ref{X1,x1,X1F11,gx1x2}\right)
\end{equation*}%
\begin{equation*}
\beta _{1}\left[ B(g\otimes 1_{H};G,1_{H})+B(x_{1}\otimes 1_{H};G,x_{1})%
\right] =0\left( \ref{X1,x1,X1F21,1}\right)
\end{equation*}%
\begin{equation*}
\beta _{1}B\left( g\otimes 1_{H};G,x_{1}x_{2}\right) =0.\left( \ref%
{X1,x1,X1F21,x1x2}\right)
\end{equation*}%
\begin{equation*}
\lambda B(x_{1}\otimes 1_{H};G,gx_{1}x_{2})=-2B\left( x_{1}\otimes
1_{H};G,g\right) .\left( \ref{X1,x1,X1F21,gx1}\right)
\end{equation*}%
\begin{gather*}
2B(x_{1}\otimes 1_{H};1_{A},1_{H})+\lambda \left[ -B\left( g\otimes
1_{H};1_{A},x_{2}\right) +B(x_{1}\otimes 1_{H};1_{A},x_{1}x_{2})\right]
\left( \ref{X1,x1,X1F31,1H}\right) \\
+\gamma _{1}\left[ B(g\otimes 1_{H};G,1_{H})+B(x_{1}\otimes 1_{H};G,x_{1})%
\right] =0.
\end{gather*}%
\begin{equation*}
2B\left( x_{1}\otimes 1_{H};1_{A},x_{1}x_{2}\right) +\gamma _{1}B\left(
g\otimes 1_{H};G,x_{1}x_{2}\right) =0\left( \ref{X1,x1,X1F31,x1x2}\right)
\end{equation*}

\begin{equation*}
-\lambda B\left( g\otimes 1_{H};1_{A},gx_{1}x_{2}\right) -\gamma _{1}B\left(
g\otimes 1_{H};G,gx_{1}\right) -2\left[ B(g\otimes
1_{H};1_{A},g)+B(x_{1}\otimes 1_{H};1_{A},gx_{1})\right] =0\left( \ref%
{X1,x1,X1F31,gx1}\right)
\end{equation*}

\begin{equation*}
\gamma _{1}\left[ -B\left( g\otimes 1_{H};G,gx_{2}\right) +B(x_{1}\otimes
1_{H};G,gx_{1}x_{2})\right] =0\left( \ref{X1,x1,X1F31,gx2}\right)
\end{equation*}%
\begin{equation*}
0=\gamma _{1}B(x_{1}\otimes 1_{H};G,x_{2}).\left( \ref{X1,x1,X1F41,1H}\right)
\end{equation*}

\begin{equation*}
2\beta _{1}B\left( g\otimes 1_{H};1_{A},gx_{1}x_{2}\right) -\gamma
_{1}B(x_{1}\otimes 1_{H};G,gx_{1}x_{2})-B(x_{1}\otimes
x_{1};X_{2},gx_{1})=0.\left( \ref{X1,x1,X1F41,gx1}\right)
\end{equation*}

\begin{eqnarray*}
&&-2\beta _{1}B(x_{2}\otimes 1_{H};1_{A},gx_{1})+\lambda \left[ -B(g\otimes
1_{H};1_{A},g)-B(x_{2}\otimes \ 1_{H};1_{A},gx_{2})\right] \\
&=&-\gamma _{1}B(x_{2}\otimes 1_{H};G,g)-2B(gx_{1}x_{2}\otimes
1_{H};1_{A},g).\left( \ref{X1,x2,X1F11,g}\right)
\end{eqnarray*}%
\begin{equation*}
\lambda \left[ -B(g\otimes 1_{H};1_{A},x_{1})-B(x_{2}\otimes \
1_{H};1_{A},x_{1}x_{2})\right] =\gamma _{1}B(x_{2}\otimes
1_{H};G,x_{1})+2B(gx_{1}x_{2}\otimes 1_{H};1_{A},x_{1}).\left( \ref%
{X1,x2,X1F11,x1}\right)
\end{equation*}%
\begin{equation*}
-\lambda B(g\otimes 1_{H};1_{A},x_{2})=\gamma _{1}B(x_{2}\otimes
1_{H};G,x_{2})+2B(gx_{1}x_{2}\otimes 1_{H};1_{A},x_{2}).\left( \ref%
{X1,x2,X1F11,x2}\right)
\end{equation*}

\begin{gather*}
\lambda B\left( g\otimes 1_{H};1_{A},gx_{1}x_{2}\right) -\gamma
_{1}B(x_{2}\otimes 1_{H};G,gx_{1}x_{2})\left( \ref{X1,x2,X1F11,gx1x2}\right)
\\
-2B\left( x_{2}\otimes 1_{H};1_{A},gx_{2}\right) -2B(gx_{1}x_{2}\otimes
1_{H};1_{A},gx_{1}x_{2})=0.
\end{gather*}%
\begin{equation*}
\lambda \left[ B(g\otimes 1_{H};G,1_{H})+B(x_{2}\otimes \ 1_{H};G,x_{2}%
\right] =-2B(gx_{1}x_{2}\otimes 1_{H};G,1_{H}).\left( \ref{X1,x2,X1F21,1H}%
\right)
\end{equation*}%
\begin{equation*}
\lambda B\left( g\otimes 1_{H};G,x_{1}x_{2}\right) +2B(gx_{1}x_{2}\otimes
1_{H};G,x_{1}x_{2})=0.\left( \ref{X1,x2,X1F21,x1x2}\right)
\end{equation*}%
\begin{equation*}
\lambda \left[ B(g\otimes 1_{H};G,gx_{1})+B(x_{2}\otimes \
1_{H};G,gx_{1}x_{2})\right] =2B(gx_{1}x_{2}\otimes 1_{H};G,gx_{1})-2B\left(
x_{2}\otimes 1_{H};G,g\right) .\left( \ref{X1,x2,X1F21,gx1}\right)
\end{equation*}

\begin{equation*}
-2\beta _{1}B(x_{2}\otimes 1_{H};G,gx_{1}x_{2})+\lambda B(g\otimes
1_{H};G,gx_{2})-2B(gx_{1}x_{2}\otimes 1_{H};G,gx_{2})=0.\left( \ref%
{X1,x2,X1F21,gx2}\right)
\end{equation*}%
\begin{equation*}
-\lambda B\left( g\otimes 1_{H};1_{A},gx_{1}x_{2}\right) =-\gamma
_{1}B(x_{2}\otimes 1_{H};G,gx_{1}x_{2})-2B(x_{2}\otimes
1_{H};1_{A},gx_{2})-2B(gx_{1}x_{2}\otimes 1_{H};1_{A},gx_{1}x_{2}).\left( %
\ref{X1,x2,X1F31,gx2}\right)
\end{equation*}

\begin{equation*}
0=-\gamma _{1}\left[ B(g\otimes 1_{H};G,1_{H})+B(x_{2}\otimes \ 1_{H};G,x_{2}%
\right] -2B(x_{1}\otimes 1_{H};1_{A},1_{H})-2B(gx_{1}x_{2}\otimes
1_{H};1_{A},x_{2}).\left( \ref{X1,x2,X1F41,1H}\right)
\end{equation*}%
\begin{equation*}
0=-\gamma _{1}B\left( g\otimes 1_{H};G,x_{1}x_{2}\right) -2B(x_{1}\otimes
1_{H};1_{A},x_{1}x_{2}).\left( \ref{X1,x2,X1F41,x1x2}\right)
\end{equation*}%
\begin{gather*}
B(x_{1}\otimes 1_{H};1_{A},gx_{1})+B(gx_{1}x_{2}\otimes
1_{H};1_{A},gx_{1}x_{2})\left( \ref{X1,x2,X1F41,gx1}\right) \\
+B(g\otimes 1_{H};1_{A},g)+B(x_{2}\otimes \ 1_{H};1_{A},gx_{2})=0.
\end{gather*}

\begin{equation*}
2\beta _{1}B\left( g\otimes 1_{H};1_{A},gx_{1}x_{2}\right) =\gamma
_{1}B(g\otimes 1_{H};G,gx_{2})+2B(x_{1}\otimes 1_{H};1_{A},gx_{2}).\left( %
\ref{X1,x2,X1F41,gx2}\right)
\end{equation*}

\begin{eqnarray*}
0 &=&\gamma _{1}\left[ B(g\otimes 1_{H};G,gx_{1})+B(x_{2}\otimes \
1_{H};G,gx_{1}x_{2})\right] \left( \ref{X1,x2,X1F51,g}\right) \\
&&+\left[ 2B(g\otimes 1_{H};1_{A},g)+2B(x_{2}\otimes \
1_{H};1_{A},gx_{2})+2B(x_{1}\otimes
1_{H};1_{A},gx_{1})+2B(gx_{1}x_{2}\otimes 1_{H};1_{A},gx_{1}x_{2})\right]
\end{eqnarray*}

\begin{equation*}
0=\gamma _{1}B\left( g\otimes 1_{H};G,x_{1}x_{2}\right) +2B(x_{1}\otimes
1_{H};1_{A},x_{1}x_{2}).\left( \ref{X1,x2,X1F51,x2}\right)
\end{equation*}%
\begin{equation*}
\lambda \left[ B(g\otimes 1_{H};G,gx_{1})+B(x_{2}\otimes \
1_{H};G,gx_{1}x_{2})\right] =-2\left[ B(x_{2}\otimes
1_{H};G,g)-B(gx_{1}x_{2}\otimes 1_{H};G,gx_{1})\right] \left( \ref%
{X1,x2,X1F61,g}\right)
\end{equation*}%
\begin{equation*}
\beta _{1}\left[ B(g\otimes 1_{H};G,gx_{1})+B(x_{2}\otimes \
1_{H};G,gx_{1}x_{2})\right] =0\left( \ref{X1,x2,X1F71,g}\right)
\end{equation*}%
\begin{equation*}
B(x_{1}\otimes 1_{H};G,x_{1})+B(gx_{1}x_{2}\otimes
1_{H};G,x_{1}x_{2})=0\left( \ref{X1,x2,X1F71,x1}\right)
\end{equation*}%
\begin{equation*}
\lambda B(x_{1}x_{2}\otimes 1_{H};X_{2},g)-\gamma _{1}B(x_{1}x_{2}\otimes
1_{H};G,g)=0\left( \ref{X1,x1x2,X1F11,g}\right)
\end{equation*}%
\begin{equation*}
2\beta _{1}B(x_{1}x_{2}\otimes 1_{H};X_{1},x_{1})+\lambda
B(x_{1}x_{2}\otimes 1_{H};X_{2},x_{1})+\gamma _{1}B(x_{1}x_{2}\otimes
1_{H};G,x_{1})=0\left( \ref{X1,x1x2,X1F11,x1}\right)
\end{equation*}

\begin{equation*}
2\beta _{1}B(x_{1}x_{2}\otimes 1_{H};X_{1},x_{2})+\lambda
B(x_{1}x_{2}\otimes 1_{H};X_{2},x_{2})+\gamma _{1}B(x_{1}x_{2}\otimes
1_{H};G,x_{2})=0.\left( \ref{X1,x1x2,X1F11,x2}\right)
\end{equation*}

\begin{equation*}
\gamma _{1}B(x_{1}x_{2}\otimes 1_{H};G,gx_{1}x_{2})-2B(x_{1}x_{2}\otimes
1_{H};1_{A},gx_{2})=0.\left( \ref{X1,x1x2,X1F11,gx1x2}\right)
\end{equation*}%
\begin{equation*}
2\beta _{1}B(x_{1}x_{2}\otimes 1_{H};GX_{1},1_{H})+\lambda
B(x_{1}x_{2}\otimes 1_{H};GX_{2},1_{H})=0.\left( \ref{X1,x1x2,X1F21,1H}%
\right)
\end{equation*}%
\begin{equation*}
\lambda B(x_{1}x_{2}\otimes 1_{H};GX_{2},gx_{1})+2B(x_{1}x_{2}\otimes
1_{H};G,g)=0.\left( \ref{X1,x1x2,X1F21,gx1}\right)
\end{equation*}%
\begin{equation*}
\lambda B(x_{1}x_{2}\otimes 1_{H};GX_{2},gx_{2})=0.\left( \ref%
{X1,x1x2,X1F21,gx2}\right)
\end{equation*}

\begin{gather*}
\lambda \left[ +1-B(x_{1}x_{2}\otimes
1_{H};1_{A},x_{1}x_{2})-B(x_{1}x_{2}\otimes
1_{H};X_{2},x_{1})+B(x_{1}x_{2}\otimes 1_{H};X_{1},x_{2})\right] \left( \ref%
{X1,x1x2,X1F31,1H}\right) \\
-\gamma _{1}B(x_{1}x_{2}\otimes 1_{H};GX_{1},1_{H})=0
\end{gather*}%
\begin{equation*}
\gamma _{1}B(x_{1}x_{2}\otimes 1_{H};GX_{1},x_{1}x_{2})=0\left( \ref%
{X1,x1x2,X1F31,x1x2}\right)
\end{equation*}%
\begin{equation*}
2B(x_{1}x_{2}\otimes 1_{H};1_{A},gx_{1})+\gamma _{1}B(x_{1}x_{2}\otimes
1_{H};GX_{1},gx_{1})+2B(x_{1}x_{2}\otimes 1_{H};X_{1},g)=0\left( \ref%
{X1,x1x2,X1F31,gx1}\right)
\end{equation*}%
\begin{equation*}
2B(x_{1}x_{2}\otimes 1_{H};1_{A},gx_{2})+\gamma _{1}B(x_{1}x_{2}\otimes
1_{H};GX_{1},gx_{2})=0.\left( \ref{X1,x1x2,X1F31,gx2}\right)
\end{equation*}

\begin{equation*}
2\beta _{1}\left[
\begin{array}{c}
+1-B(x_{1}x_{2}\otimes 1_{H};1_{A},x_{1}x_{2})-B(x_{1}x_{2}\otimes
1_{H};X_{2},x_{1}) \\
+B(x_{1}x_{2}\otimes 1_{H};X_{1},x_{2})%
\end{array}%
\right] +\gamma _{1}B(x_{1}x_{2}\otimes 1_{H};GX_{2},1_{H})=0.\left( \ref%
{X1,x1x2,X1F41,1H}\right)
\end{equation*}

\begin{equation*}
\gamma _{1}B(x_{1}x_{2}\otimes 1_{H};GX_{2},gx_{1})+2B(x_{1}x_{2}\otimes
1_{H};X_{2},g)=0.\left( \ref{X1,x1x2,X1F41,gx1}\right)
\end{equation*}%
\begin{equation*}
\gamma _{1}B(x_{1}x_{2}\otimes 1_{H};GX_{2},gx_{2})=0\left( \ref%
{X1,x1x2,X1F41,gx2}\right)
\end{equation*}%
\begin{gather*}
2B(x_{1}x_{2}\otimes 1_{H};X_{2},g)+\left( \ref{X1,x1x2,X1F51,g}\right) \\
+\gamma _{1}\left[ -B(x_{1}x_{2}\otimes
1_{H};G,gx_{1}x_{2})+B(x_{1}x_{2}\otimes
1_{H};GX_{2},gx_{1})-B(x_{1}x_{2}\otimes 1_{H};GX_{1},gx_{2})\right] =0
\end{gather*}%
\begin{equation*}
2B(x_{1}x_{2}\otimes 1_{H};G,g)+\lambda \left[
\begin{array}{c}
-B(x_{1}x_{2}\otimes 1_{H};G,gx_{1}x_{2}) \\
+B(x_{1}x_{2}\otimes 1_{H};GX_{2},gx_{1})-B(x_{1}x_{2}\otimes
1_{H};GX_{1},gx_{2})%
\end{array}%
\right] =0\left( \ref{X1,x1x2,X1F61,g}\right)
\end{equation*}%
\begin{equation*}
\lambda B(x_{1}x_{2}\otimes 1_{H};GX_{1},x_{1}x_{2})=0.\left( \ref%
{X1,x1x2,X1F61,x1}\right)
\end{equation*}%
\begin{equation*}
\lambda B(x_{1}x_{2}\otimes 1_{H};GX_{2},x_{1}x_{2})=0\left( \ref%
{X1,x1x2,X1F61,x2}\right)
\end{equation*}%
\begin{equation*}
B(x_{1}x_{2}\otimes 1_{H};G,gx_{1}x_{2})+B(x_{1}x_{2}\otimes
1_{H};GX_{1},gx_{2})=0\left( \ref{X1,x1x2,X1F61,gx1x2}\right)
\end{equation*}%
\begin{equation*}
\beta _{1}B(x_{1}x_{2}\otimes 1_{H};GX_{1},x_{1}x_{2})=0\left( \ref%
{X1,x1x2,X1F71,x1}\right)
\end{equation*}%
\begin{equation*}
\beta _{1}B(x_{1}x_{2}\otimes 1_{H};GX_{2},x_{1}x_{2})=0\left( \ref%
{X1,x1x2,X1F71,x2}\right)
\end{equation*}%
\begin{equation*}
B(x_{1}x_{2}\otimes 1_{H};GX_{2},gx_{2})=0.\left( \ref{X1,x1x2,X1F71,gx1x2}%
\right)
\end{equation*}%
\begin{equation*}
+B(x_{1}x_{2}\otimes 1_{H};G,gx_{1}x_{2})+B(x_{1}x_{2}\otimes
1_{H};GX_{1},gx_{2})=0.\left( \ref{X1,x1x2,X1F81,gx1}\right)
\end{equation*}%
\begin{gather*}
2\beta _{1}\left[ -1+B(x_{1}x_{2}\otimes
1_{H};1_{A},x_{1}x_{2})+B(x_{1}x_{2}\otimes 1_{H};X_{2},x_{1})\right]
+\lambda B(x_{1}x_{2}\otimes 1_{H};X_{2},x_{2})\left( \ref{X1,gx1,X1F11,1H}%
\right) \\
=-\gamma _{1}\left[ B(x_{1}x_{2}\otimes 1_{H};G,x_{2})-B(x_{1}x_{2}\otimes
1_{H};GX_{2},1_{H})\right] + \\
-2\left[ B(x_{1}x_{2}\otimes 1_{H};1_{A},gx_{2})+B(x_{1}x_{2}\otimes
1_{H};X_{2},g)\right] .
\end{gather*}

\begin{gather*}
\lambda B(x_{1}x_{2}\otimes 1_{H};X_{2},gx_{1}x_{2})-\gamma _{1}\left[
B(x_{1}x_{2}\otimes 1_{H};G,gx_{1}x_{2})-B(x_{1}x_{2}\otimes
1_{H};GX_{2},gx_{1})\right] \left( \ref{X1,gx1,X1F11,gx1}\right) \\
+2\left[ B(x_{1}x_{2}\otimes 1_{H};1_{A},gx_{2})+B(x_{1}x_{2}\otimes
1_{H};X_{2},g)\right] =0.
\end{gather*}%
\begin{gather*}
2\left[ B(x_{1}x_{2}\otimes 1_{H};1_{A},gx_{2})+B(x_{1}x_{2}\otimes
1_{H};X_{2},g)\right] +\left( \ref{X1,gx1,X1F31,g}\right) \\
+\lambda B(x_{1}x_{2}\otimes 1_{H};X_{2},gx_{1}x_{2}) \\
+\gamma _{1}\left[ -B(x_{1}x_{2}\otimes
1_{H};G,gx_{1}x_{2})+B(x_{1}x_{2}\otimes 1_{H};GX_{2},gx_{1})\right] =0.
\end{gather*}%
\begin{gather*}
-2\beta _{1}B(x_{1}x_{2}\otimes 1_{H};X_{1},x_{1})+\lambda \left[
-1+B(x_{1}x_{2}\otimes 1_{H};1_{A},x_{1}x_{2})-B(x_{1}x_{2}\otimes
1_{H};X_{1},x_{2})\right] \left( \ref{X1,gx2,X1F11,1H}\right) \\
+\gamma _{1}\left[ -B(x_{1}x_{2}\otimes 1_{H};G,x_{1})+B(x_{1}x_{2}\otimes
1_{H};GX_{1},1_{H})\right] =0.
\end{gather*}%
\begin{gather*}
\lambda B(x_{1}x_{2}\otimes 1_{H};X_{1},gx_{1}x_{2})+\gamma
_{1}B(x_{1}x_{2}\otimes 1_{H};GX_{1},gx_{1})\left( \ref{X1,gx2,X1F11,gx1}%
\right) \\
+2B(x_{1}x_{2}\otimes 1_{H};1_{A},gx_{1})+2B(x_{1}x_{2}\otimes
1_{H};X_{1},g)=0.
\end{gather*}%
\begin{equation*}
+\gamma _{1}\left[ B(x_{1}x_{2}\otimes
1_{H};G,gx_{1}x_{2})+B(x_{1}x_{2}\otimes 1_{H};GX_{1},gx_{2})\right]
=0\left( \ref{X1,gx2,X1F11,gx2}\right)
\end{equation*}%
\begin{gather*}
2\left[ B(x_{1}x_{2}\otimes 1_{H};1_{A},gx_{1})+B(x_{1}x_{2}\otimes
1_{H};X_{1},g)\right] \left( \ref{X1,gx2,X131,g}\right) \\
+\lambda B(x_{1}x_{2}\otimes 1_{H};X_{1},gx_{1}x_{2})+\gamma
_{1}B(x_{1}x_{2}\otimes 1_{H};GX_{1},gx_{1})=0
\end{gather*}%
\begin{equation*}
\gamma _{1}\left[ -B(x_{1}x_{2}\otimes
1_{H};G,gx_{1}x_{2})+B(x_{1}x_{2}\otimes 1_{H};GX_{1},gx_{2})\right]
=0\left( \ref{X1,gx2,X141,g}\right)
\end{equation*}%
\begin{equation*}
\beta _{1}B(gx_{2}\otimes 1_{H};GX_{1}X_{2},1_{H})=0.\left( \ref%
{X1,gx2,X171,1H}\right)
\end{equation*}%
\begin{equation*}
+\lambda \left[ B(x_{1}\otimes 1_{H};1_{A},1_{H})+B(gx_{1}x_{2}\otimes
1_{H};1_{A},x_{2})\right] -\gamma _{1}B(gx_{1}x_{2}\otimes
1_{H};G,1_{H})=0\left( \ref{X1,gx1x2,X1F11,1H}\right)
\end{equation*}%
\begin{equation*}
+\lambda B(x_{1}\otimes 1_{H};1_{A},x_{1}x_{2})-\gamma
_{1}B(gx_{1}x_{2}\otimes 1_{H};G,x_{1}x_{2})=0\left( \ref%
{X1,gx1x2,X1F11,x1x2}\right)
\end{equation*}%
\begin{gather*}
-2\beta _{1}B(x_{2}\otimes 1_{H};1_{A},gx_{1})+\lambda \left[ B(x_{1}\otimes
1_{H};1_{A},gx_{1})+B(gx_{1}x_{2}\otimes 1_{H};1_{A},gx_{1}x_{2})\right]
\left( \ref{X1,gx1x2,X1F11,gx1}\right) \\
+\gamma _{1}B(gx_{1}x_{2}\otimes 1_{H};G,gx_{1})+2B(gx_{1}x_{2}\otimes
1_{H};1_{A},g)=0.
\end{gather*}%
\begin{gather*}
2\beta _{1}\left[ -B(x_{2}\otimes 1_{H};1_{A},gx_{2})-B(gx_{1}x_{2}\otimes
1_{H};1_{A},gx_{1}x_{2})\right] \left( \ref{X1,gx1x2,X171,gx2}\right) \\
+\lambda B(x_{1}\otimes 1_{H};1_{A},gx_{2})+\gamma _{1}B(gx_{1}x_{2}\otimes
1_{H};G,gx_{2})=0.
\end{gather*}%
\begin{gather*}
2\beta _{1}\left[ B(x_{2}\otimes 1_{H};G,g)-B(gx_{1}x_{2}\otimes
1_{H};G,gx_{1})\right] \left( \ref{X1,gx1x2,XF21,g}\right) \\
+\lambda \left[ -B(x_{1}\otimes 1_{H};G,g)-B(gx_{1}x_{2}\otimes
1_{H};G,gx_{2})\right] =0
\end{gather*}

\begin{equation*}
2\beta _{1}B(x_{2}\otimes 1_{H};G,gx_{1}x_{2})-\lambda B(x_{1}\otimes
1_{H};G,gx_{1}x_{2})=0\left( \ref{X1,gx1x2,X1F21,gx1x2}\right)
\end{equation*}%
\begin{gather*}
+\lambda \left[ B(g\otimes 1_{H};1_{A},g)+B(x_{2}\otimes \
1_{H};1_{A},gx_{2})+B(x_{1}\otimes 1_{H};1_{A},gx_{1})+B(gx_{1}x_{2}\otimes
1_{H};1_{A},gx_{1}x_{2})\right] \\
-\gamma _{1}\left[ B(x_{2}\otimes 1_{H};G,g)-B(gx_{1}x_{2}\otimes
1_{H};G,gx_{1})\right] =0\left( \ref{X1,gx1x2,X1F31,g}\right)
\end{gather*}%
\begin{gather*}
2B(gx_{1}x_{2}\otimes 1_{H};1_{A},x_{1})+\lambda \left[ B(g\otimes
1_{H};1_{A},x_{1})+B(x_{2}\otimes \ 1_{H};1_{A},x_{1}x_{2})\right] \\
+\gamma _{1}B(x_{2}\otimes 1_{H};G,x_{1})=0\left( \ref{X1,gx1x2,X1F31,x1}%
\right)
\end{gather*}

\begin{gather*}
2B(gx_{1}x_{2}\otimes 1_{H};1_{A},x_{2})+\lambda \left[ B(g\otimes
1_{H};X_{2},1_{H})-B(x_{1}\otimes 1_{H};1_{A},x_{1}x_{2})\right] +\left( \ref%
{X1,gx1x2,X1F31,x2}\right) \\
+\gamma _{1}\left[ B(x_{2}\otimes 1_{H};G,x_{2})+B(gx_{1}x_{2}\otimes
1_{H};G,x_{1}x_{2})\right] =0
\end{gather*}

\begin{gather*}
+\lambda B\left( g\otimes 1_{H};1_{A},gx_{1}x_{2}\right) +\left( \ref%
{X1,gx1x2,X1F31,gx1x2}\right) \\
-\gamma _{1}B(x_{2}\otimes 1_{H};G,gx_{1}x_{2})-\left[ +2B(x_{2}\otimes
1_{H};1_{A},gx_{2})+2B(gx_{1}x_{2}\otimes 1_{H};1_{A},gx_{1}x_{2})\right] =0
\end{gather*}

\begin{eqnarray*}
&&2\beta _{1}\left[ B(g\otimes 1_{H};1_{A},g)+B(x_{2}\otimes \
1_{H};1_{A},gx_{2})+B(x_{1}\otimes 1_{H};1_{A},gx_{1})+B(gx_{1}x_{2}\otimes
1_{H};1_{A},gx_{1}x_{2})\right] \left( \ref{X1,gx1x2,X1F41,g}\right) \\
&&-\gamma _{1}\left[ B(x_{1}\otimes 1_{H};G,g)+B(gx_{1}x_{2}\otimes
1_{H};G,gx_{2})\right] =0
\end{eqnarray*}%
\begin{eqnarray*}
&&2\beta _{1}B\left( g\otimes 1_{H};1_{A},gx_{1}x_{2}\right) \left( \ref%
{X1,gx1x2,X1F41,gx1x2}\right) \\
&&-\gamma _{1}B(x_{1}\otimes 1_{H};G,gx_{1}x_{2})-2B(x_{1}\otimes
1_{H};1_{A},gx_{2})=0
\end{eqnarray*}%
\begin{gather*}
2\left[ B(x_{1}\otimes 1_{H};1_{A},1_{H})+B(gx_{1}x_{2}\otimes
1_{H};1_{A},x_{2})\right] \left( \ref{X1,gx1x2,X1F51,1H}\right) \\
+\gamma _{1}\left[ B(g\otimes 1_{H};G,1_{H})+B(x_{2}\otimes \
1_{H};G,x_{2})+B(x_{1}\otimes 1_{H};G,x_{1})+B(gx_{1}x_{2}\otimes
1_{H};G,x_{1}x_{2})\right] =0
\end{gather*}%
\begin{gather*}
+\gamma _{1}\left[ B(g\otimes 1_{H};G,gx_{1})+B(x_{2}\otimes \
1_{H};G,gx_{1}x_{2})\right] \left( \ref{X1,gx1x2,X1F51,gx1}\right) \\
+-\left[ 2B(g\otimes 1_{H};1_{A},g)+2B(x_{2}\otimes \
1_{H};1_{A},gx_{2})+2B(x_{1}\otimes
1_{H};1_{A},gx_{1})+2B(gx_{1}x_{2}\otimes 1_{H};1_{A},gx_{1}x_{2})\right] =0
\end{gather*}%
\begin{gather*}
2B(gx_{1}x_{2}\otimes 1_{H};G,1_{H})\left( \ref{X1,gx1x2,X1F61,1H}\right) \\
+\lambda \left[ B(g\otimes 1_{H};G,1_{H})+B(x_{2}\otimes \
1_{H};G,x_{2})+B(x_{1}\otimes 1_{H};G,x_{1})+B(gx_{1}x_{2}\otimes
1_{H};G,x_{1}x_{2})\right] =0
\end{gather*}%
\begin{equation*}
2B(gx_{1}x_{2}\otimes 1_{H};G,x_{1}x_{2})+\lambda B(g\otimes
1_{H};G,x_{1}x_{2})=0\left( \ref{X1,gx1x2,X1F61,x1x2}\right)
\end{equation*}%
\begin{equation*}
\lambda \left[ B(g\otimes 1_{H};G,gx_{1})+B(x_{2}\otimes \
1_{H};G,gx_{1}x_{2})\right] +\left[ 2B(x_{2}\otimes
1_{H};G,g)-2B(gx_{1}x_{2}\otimes 1_{H};G,gx_{1})\right] =0\left( \ref%
{X1,gx1x2,X1F61,gx1}\right)
\end{equation*}%
\begin{gather*}
\beta _{1}\left[ B(g\otimes 1_{H};G,gx_{1})+B(x_{2}\otimes \
1_{H};G,gx_{1}x_{2})\right] \left( \ref{X1,gx1x2,X1F71,gx1}\right) \\
+B(x_{1}\otimes 1_{H};G,g)+B(gx_{1}x_{2}\otimes 1_{H};G,gx_{2})=0.
\end{gather*}

\subsection{$X_{2}$}

\begin{gather*}
\gamma _{2}B(g\otimes 1_{H};G,1_{H})+\lambda B(g\otimes
1_{H};1_{A},x_{1})\left( \ref{X2,g,X2F11,1H}\right) \\
+2B(x_{2}\otimes 1_{H};1_{A},1_{H})=0
\end{gather*}%
\begin{gather*}
2\beta _{2}B(g\otimes 1_{H};1_{A},gx_{1}x_{2})+\gamma _{2}B(g\otimes
1_{H};G,gx_{1})\left( \ref{X2,g,X2F11,gx1}\right) \\
+2B(x_{2}\otimes 1_{H};1_{A},gx_{1})=0
\end{gather*}

\begin{gather*}
\gamma _{2}B(g\otimes 1_{H};G,gx_{2})+\lambda B(g\otimes
1_{H};1_{A},gx_{1}x_{2})\left( \ref{X2,g,X2F11,gx2}\right) \\
+2B(x_{2}\otimes 1_{H};1_{A},gx_{2})+2B(g\otimes 1_{H};1_{A},g)=0.
\end{gather*}

\begin{equation*}
2\beta _{2}B(g\otimes 1_{H};G,gx_{2})+\lambda B(g\otimes
1_{H};G,gx_{1})+2B(x_{2}\otimes 1_{H};G,g)=0\left( \ref{X2,g,X2F21,g}\right)
\end{equation*}%
\begin{equation*}
B(x_{2}\otimes 1_{H};G,x_{1})=0\left( \ref{X2,g,X2F21,x1}\right)
\end{equation*}%
\begin{equation*}
\lambda B(g\otimes 1_{H};G,x_{1}x_{2})-2B(x_{2}\otimes 1_{H};G,x_{2})=0\ref%
{X2,g,X2F21,x2}
\end{equation*}%
\begin{equation*}
2\beta _{2}B(g\otimes 1_{H};1_{H},gx_{1}x_{2})+\gamma _{2}B(g\otimes
1_{H};G,gx_{1})+2B(x_{2}\otimes 1_{H};1_{H},gx_{1})=0.\left( \ref%
{X2,g,X2F31,g}\right)
\end{equation*}%
\begin{equation*}
\gamma _{2}B(g\otimes 1_{H};G,x_{1}x_{2})+2B(x_{2}\otimes
1_{H};1,x_{1}x_{2})=0\left( \ref{X2,g,X2F31,x2}\right)
\end{equation*}%
\begin{eqnarray*}
3B(g\otimes 1_{H};1_{A},x_{1})-\gamma _{2}B(g\otimes
1_{H};G,x_{1}x_{2})\left( \ref{X2,g,X2F41,x1}\right) && \\
+B(x_{2}\otimes 1_{H};X_{2},x_{1})+B\left( x_{2}\otimes
1_{H};1_{A},x_{1}x_{2}\right) &=&0
\end{eqnarray*}%
\begin{equation*}
B\left( g\otimes 1_{H};1_{A},gx_{1}x_{2}\right) +B(g\otimes
1_{H},G,gx_{1})=0\left( \ref{X2,g,X2F61,gx2}\right)
\end{equation*}%
\begin{gather*}
2B(g\otimes 1_{H};G,1_{H})-\lambda B(g\otimes 1_{H};G,x_{1}x_{2})\left( \ref%
{X2,g,X2F71,1H}\right) \\
+2B(g\otimes 1_{H};G,1_{H})+2B(x_{2}\otimes 1_{H};G,x_{2})=0
\end{gather*}%
\begin{equation*}
B(g\otimes 1_{H};G,gx_{1})+B(x_{2}\otimes 1_{H};G,gx_{1}x_{2})=0\left( \ref%
{X2,g,X2F71,gx1}\right)
\end{equation*}

\begin{eqnarray*}
&&-2\beta _{2}B(x_{1}\otimes 1_{H};1_{A},gx_{2})+\gamma _{2}B(x_{1}\otimes
1_{H};G,g)\left( \ref{X2,x1,X2F11,g}\right) \\
&&-\lambda \left[ B(g\otimes 1_{H};1_{A},g)+B(x_{1}\otimes
1_{H};1_{A},gx_{1})\right] -2B(gx_{1}x_{2}\otimes 1_{H};1_{A},g)=0
\end{eqnarray*}

\begin{equation*}
\gamma _{2}B(x_{1}\otimes 1_{H};G,x_{1})-\lambda B(g\otimes
1_{H};1_{A},x_{1})-2B(gx_{1}x_{2}\otimes 1_{H};1_{A},x_{1})=0\left( \ref%
{X2,x1,X2F11,x1}\right)
\end{equation*}%
\begin{eqnarray*}
&&-\gamma _{2}B(x_{1}\otimes 1_{H};G,x_{2})+\lambda \left[ B(g\otimes
1_{H};1_{A},x_{2})-B(x_{1}\otimes 1_{H};1_{A},x_{1}x_{2})\right] \left( \ref%
{X2,x1,X2F11,x2}\right) \\
&&+2B(gx_{1}x_{2}\otimes 1_{H};1_{A},x_{2})=0.
\end{eqnarray*}%
\begin{eqnarray*}
&&+\gamma _{2}B(x_{1}\otimes 1_{H};G,gx_{1}x_{2})-\lambda B(g\otimes
1_{H};1_{A},gx_{1}x_{2})\left( \ref{X2,x1,X2F11,gx1x2}\right) \\
&&-2B(gx_{1}x_{2}\otimes 1_{H};1_{A},gx_{1}x_{2})-2B\left( x_{1}\otimes
1_{H};1_{A},gx_{1}\right) =0
\end{eqnarray*}%
\begin{equation*}
\lambda \left[ B(g\otimes 1_{H};G,1_{H})+B(x_{1}\otimes 1_{H};G,x_{1})\right]
+2B(gx_{1}x_{2}\otimes 1_{H};G,1_{H})=0\left( \ref{X2,x1,X2F21,1H}\right)
\end{equation*}

\begin{equation*}
2\beta _{2}B(x_{1}\otimes 1_{H};G,gx_{1}x_{2})+\lambda B\left( g\otimes
1_{H};G,gx_{1}\right) +2B(x_{2}\otimes 1_{H};G,g)=0.\left( \ref%
{X2,x1,X2F21,gx1}\right)
\end{equation*}

\begin{eqnarray*}
&&\lambda \left[ B\left( g\otimes 1_{H};G,gx_{2}\right) -B(x_{1}\otimes
1_{H};G,gx_{1}x_{2})\right] \left( \ref{X2,x1,X2F21,gx2}\right) \\
&&+2-B(x_{1}\otimes 1_{H};G,g_{1})+2B\left( x_{1}\otimes 1_{H};G,g\right) =0.
\end{eqnarray*}

\begin{gather*}
\gamma _{2}\left[ B(g\otimes 1_{H};G,1_{H})+B(x_{1}\otimes 1_{H};G,x_{1})%
\right] \left( \ref{X2,x1,X2F31,1H}\right) \\
+2B(x_{2}\otimes 1_{H};1_{A},1_{H})-2B(gx_{1}x_{2}\otimes
1_{H};1_{A},x_{1})=0.
\end{gather*}%
\begin{gather*}
\gamma _{2}\left[ -B(g\otimes 1_{H};G,gx_{2})+B(x_{1}\otimes
1_{H};G,gx_{1}x_{2})\right] +\left( \ref{X2,x1,X2F31,gx2}\right) \\
-2B(x_{2}\otimes 1_{H};1_{A},gx_{2})-2B(gx_{1}x_{2}\otimes
1_{H};1_{A},gx_{1}x_{2}) \\
-2B(g\otimes 1_{H};1_{A},g)-2B(x_{1}\otimes 1_{H};1_{A},gx_{1})=0.
\end{gather*}%
\begin{gather*}
+\gamma _{2}B(x_{1}\otimes 1_{H};G,x_{2})\left( \ref{X2,x1,X2F41,1H}\right)
\\
+\lambda \left[ -B(g\otimes 1_{H};1_{A},x_{2})+B(x_{1}\otimes
1_{H};1_{A},x_{1}x_{2})\right] \\
-2B(gx_{1}x_{2}\otimes 1_{H};1_{A},x_{2})=0
\end{gather*}

\begin{eqnarray*}
&&-\gamma _{2}B(x_{1}\otimes 1_{H};G,gx_{1}x_{2})+\lambda B\left( g\otimes
1_{H};1_{A},gx_{1}x_{2}\right) \left( \ref{X2,x1,X2F41,gx1}\right) \\
&&+2B(x_{1}\otimes 1_{H};1_{A},gx_{1})+2B(gx_{1}x_{2}\otimes
1_{H};1_{A},gx_{1}x_{2})=0.
\end{eqnarray*}

\begin{eqnarray*}
&&\beta _{2}\left[ -B(g\otimes 1_{H};G,gx_{2})+B(x_{1}\otimes
1_{H};G,gx_{1}x_{2})\right] \left( \ref{X2,x1,X2F61,g}\right) \\
&&-B(x_{2}\otimes 1_{H};G,g)+B(gx_{1}x_{2}\otimes 1_{H};G,gx_{1})=0
\end{eqnarray*}%
\begin{equation*}
B(x_{2}\otimes 1_{H};G,x_{2})+B(gx_{1}x_{2}\otimes
1_{H};G,x_{1}x_{2})=0\left( \ref{X2,x1,X2F61,x2}\right)
\end{equation*}

\begin{equation*}
B(x_{1}\otimes 1_{H};G,g)+B(gx_{1}x_{2}\otimes 1_{H};G,gx_{2})=0\left( \ref%
{X2,x1,X2F71,g}\right)
\end{equation*}%
\begin{equation*}
B(x_{2}\otimes \ 1_{H};G,x_{2})+B(gx_{1}x_{2}\otimes
1_{H};G,x_{1}x_{2})=0\left( \ref{X2,x1,X2F81,1H}\right)
\end{equation*}%
\begin{eqnarray*}
2 &&\beta _{2}\left[ B(g\otimes 1_{H};1_{A},g)+B(x_{2}\otimes \
1_{H};1_{A},gx_{2})\right] \left( \ref{X2,x2,X2F11,g}\right) \\
&&-\gamma _{2}B(x_{2}\otimes 1_{H};G,g)+\lambda B(x_{2}\otimes
1_{H};1_{A},gx_{1})=0.
\end{eqnarray*}%
\begin{equation*}
\gamma _{2}B(x_{2}\otimes 1_{H};G,x_{2})+\lambda B(x_{2}\otimes
1_{H};1_{A},x_{1}x_{2})=0.\left( \ref{X2,x2,X2F11,x2}\right)
\end{equation*}

\begin{eqnarray*}
+ &&2\beta _{2}B\left( g\otimes 1_{H};1_{A},gx_{1}x_{2}\right) -\gamma
_{2}B(x_{2}\otimes 1_{H};G,gx_{1}x_{2})\left( \ref{X2,x2,X2F11,gx1x2}\right)
\\
&&+2B\left( x_{2}\otimes 1_{H};1_{A},gx_{1}\right) =0
\end{eqnarray*}

\begin{equation*}
\beta _{2}\left[ B(g\otimes 1_{H};G,gx_{1})+B(x_{2}\otimes \
1_{H};G,gx_{1}x_{2})\right] =0\left( \ref{X2,x2,X2F21,gx1}\right)
\end{equation*}%
\begin{eqnarray*}
&&2\beta _{2}B(g\otimes 1_{H};G,gx_{2})-\lambda B(x_{2}\otimes
1_{H};G,gx_{1}x_{2})\left( \ref{X2,x2,X2F21,gx2}\right) \\
&&+2B\left( x_{2}\otimes 1_{H};G,g\right) =0
\end{eqnarray*}

\begin{eqnarray*}
&&2B(x_{2}\otimes 1_{H};1_{A},1_{H})+\gamma _{2}\left[ B(g\otimes
1_{H};G,1_{H})+B(x_{2}\otimes \ 1_{H};G,x_{2}\right] \left( \ref%
{X2,x2,X2F41,1H}\right) \\
&&+\lambda \left[ B(g\otimes 1_{H};1_{A},x_{1})+B(x_{2}\otimes \
1_{H};1_{A},x_{1}x_{2})\right] =0.
\end{eqnarray*}

\begin{eqnarray*}
&&\gamma _{2}B(g\otimes 1_{H};G,gx_{2})-\lambda B\left( g\otimes
1_{H};1_{A},gx_{1}x_{2}\right) \left( \ref{X2,x2,X2F41,gx2}\right) \\
&&+2\left[ +B(g\otimes 1_{H};1_{A},g)+B(x_{2}\otimes \ 1_{H};1_{A},gx_{2})%
\right] =0.
\end{eqnarray*}

\begin{equation*}
\lambda \left[ B(g\otimes 1_{H};G,gx_{1})+B(x_{2}\otimes \
1_{H};G,gx_{1}x_{2})\right] =0\left( \ref{X2,x2,X2F71,g}\right)
\end{equation*}%
\begin{equation*}
\gamma _{2}B(x_{1}x_{2}\otimes 1_{H};G,g)+\lambda B(x_{1}x_{2}\otimes
1_{H};X_{1},g)=0.\left( \ref{X2,x1x2,X2F11,g}\right)
\end{equation*}

\begin{gather*}
2\beta _{2}B(x_{1}x_{2}\otimes 1_{H};X_{2},x_{1})+\gamma
_{2}B(x_{1}x_{2}\otimes 1_{H};G,x_{1})\left( \ref{X2,x1x2,X2F11,x1}\right) \\
+\lambda B(x_{1}x_{2}\otimes 1_{H};X_{1},x_{1})=0.
\end{gather*}

\begin{gather*}
2\beta _{2}B(x_{1}x_{2}\otimes 1_{H};X_{2},x_{2})+\gamma
_{2}B(x_{1}x_{2}\otimes 1_{H};G,x_{2})\left( \ref{X2,x1x2,X2F11,x2}\right) \\
+\lambda B(x_{1}x_{2}\otimes 1_{H};X_{1},x_{2})=0.
\end{gather*}%
\begin{eqnarray*}
&&\gamma _{2}B(x_{1}x_{2}\otimes 1_{H};G,gx_{1}x_{2})+\lambda
B(x_{1}x_{2}\otimes 1_{H};X_{1},gx_{1}x_{2})\left( \ref{X2,x1x2,X2F11,gx1x2}%
\right) \\
&&+2B(x_{1}x_{2}\otimes 1_{H};1_{A},gx_{1})=0.
\end{eqnarray*}

\begin{equation*}
2\beta _{2}B(x_{1}x_{2}\otimes 1_{H};GX_{2},1_{H})+\lambda
B(x_{1}x_{2}\otimes 1_{H};GX_{1},1_{H})=0.\left( \ref{X2,x1x2,X2F21,1H}%
\right)
\end{equation*}%
\begin{equation*}
\beta _{2}B(x_{1}x_{2}\otimes 1_{H};GX_{2},x_{1}x_{2})=0.\left( \ref%
{X2,x1x2,X2F21,x1x2}\right)
\end{equation*}

\begin{equation*}
\lambda B(x_{1}x_{2}\otimes 1_{H};GX_{1},gx_{1})=0.\left( \ref%
{X2,x1x2,X2F21,gx1}\right)
\end{equation*}%
\begin{equation*}
2B(x_{1}x_{2}\otimes 1_{H};G,g)-\lambda B(x_{1}x_{2}\otimes
1_{H};GX_{1},gx_{2})=0.\left( \ref{X2,x1x2,X2F21,gx2}\right)
\end{equation*}%
\begin{eqnarray*}
&&2\beta _{2}\left[
\begin{array}{c}
+1-B(x_{1}x_{2}\otimes 1_{H};1_{A},x_{1}x_{2}) \\
-B(x_{1}x_{2}\otimes 1_{H};X_{2},x_{1})+B(x_{1}x_{2}\otimes
1_{H};X_{1},x_{2})%
\end{array}%
\right] \left( \ref{X2,x1x2,X2F31,1H}\right) \\
&&-\gamma _{2}B(x_{1}x_{2}\otimes 1_{H};GX_{1},1_{H})=0.
\end{eqnarray*}%
\begin{equation*}
\gamma _{2}B(x_{1}x_{2}\otimes 1_{H};GX_{1},gx_{1})=0\left( \ref%
{X2,x1x2,X2F31,gx1}\right)
\end{equation*}%
\begin{equation*}
\gamma _{2}B(x_{1}x_{2}\otimes 1_{H};GX_{1},gx_{2})+2B(x_{1}x_{2}\otimes
1_{H};X_{1},g)=0.\left( \ref{X2,x1x2,X2F31,gx2}\right)
\end{equation*}

\begin{equation*}
-\gamma _{2}B(x_{1}x_{2}\otimes 1_{H};GX_{2},1_{H})-\lambda \left[
\begin{array}{c}
+1-B(x_{1}x_{2}\otimes 1_{H};1_{A},x_{1}x_{2}) \\
-B(x_{1}x_{2}\otimes 1_{H};X_{2},x_{1})+B(x_{1}x_{2}\otimes
1_{H};X_{1},x_{2})%
\end{array}%
\right] =0.\left( \ref{X2,x1x2,X2F41,1H}\right)
\end{equation*}%
\begin{gather*}
2B(x_{1}x_{2}\otimes 1_{H};1_{A},gx_{2})+\gamma _{2}B(x_{1}x_{2}\otimes
1_{H};GX_{2},gx_{2})\left( \ref{X2,x1x2,X2F41,gx2}\right) \\
+\lambda B(x_{1}x_{2}\otimes 1_{H};X_{2},gx_{1}x_{2})+2B(x_{1}x_{2}\otimes
1_{H};X_{2},g)=0.
\end{gather*}

\begin{gather*}
2B(x_{1}x_{2}\otimes 1_{H};X_{1},g)-\left( \ref{X2,x1x2,X2F51,g}\right) \\
\gamma _{2}\left[
\begin{array}{c}
-B(x_{1}x_{2}\otimes 1_{H};G,gx_{1}x_{2})+ \\
B(x_{1}x_{2}\otimes 1_{H};GX_{2},gx_{1})-B(x_{1}x_{2}\otimes
1_{H};GX_{1},gx_{2})%
\end{array}%
\right] =0.
\end{gather*}%
\begin{equation*}
\beta _{2}B(x_{1}x_{2}\otimes 1_{H};GX_{1},x_{1}x_{2})=0.\left( \ref%
{X2,x1x2,X2F61,x1}\right)
\end{equation*}%
\begin{equation*}
B(x_{1}x_{2}\otimes 1_{H};GX_{1},gx_{1})=0.\left( \ref{X2,x1x2,X2F61,gx1x2}%
\right)
\end{equation*}%
\begin{gather*}
2B(x_{1}x_{2}\otimes 1_{H};G,g)\left( \ref{X2,x1x2,X2F71,g}\right) \\
+\lambda \left[
\begin{array}{c}
-B(x_{1}x_{2}\otimes 1_{H};G,gx_{1}x_{2})+ \\
B(x_{1}x_{2}\otimes 1_{H};GX_{2},gx_{1})-B(x_{1}x_{2}\otimes
1_{H};GX_{1},gx_{2})%
\end{array}%
\right] =0.
\end{gather*}%
\begin{equation*}
B(x_{1}x_{2}\otimes 1_{H};G,gx_{1}x_{2})-B(x_{1}x_{2}\otimes
1_{H};GX_{2},gx_{1})=0.\left( \ref{X2,x1x2,X2F71,gx1x2}\right)
\end{equation*}%
\begin{gather*}
2\beta _{2}B(x_{1}x_{2}\otimes 1_{H};X_{2},x_{2})+\gamma _{2}\left[
B(x_{1}x_{2}\otimes 1_{H};G,x_{2})-B(x_{1}x_{2}\otimes 1_{H};GX_{2},1_{H})%
\right] \left( \ref{X2,gx1,X2F11,1H}\right) \\
+\lambda \left[ -1+B(x_{1}x_{2}\otimes
1_{H};1_{A},x_{1}x_{2})+B(x_{1}x_{2}\otimes 1_{H};X_{2},x_{1})\right] =0.
\end{gather*}

\begin{equation*}
B(x_{1}x_{2}\otimes 1_{H};1_{A},gx_{2})+B(x_{1}x_{2}\otimes
1_{H};X_{2},g)=0.\left( \ref{X2,gx1,X2F11,gx2}\right)
\end{equation*}%
\begin{equation*}
\lambda \left[ B(x_{1}x_{2}\otimes 1_{H};G,gx_{1}x_{2})-B(x_{1}x_{2}\otimes
1_{H};GX_{2},gx_{1})\right] =0.\left( \ref{X2,gx1,X2F21,g}\right)
\end{equation*}%
\begin{gather*}
2\beta _{2}\left[ -1+B(x_{1}x_{2}\otimes
1_{H};1_{A},x_{1}x_{2})-B(x_{1}x_{2}\otimes 1_{H};X_{1},x_{2})\right] \left( %
\ref{X2,gx2,X2F11,1H}\right) \\
+\gamma _{2}\left[ -B(x_{1}x_{2}\otimes 1_{H};G,x_{1})+B(x_{1}x_{2}\otimes
1_{H};GX_{1},1_{H})\right] + \\
-\lambda B(x_{1}x_{2}\otimes 1_{H};X_{1},x_{1})=0.
\end{gather*}

\begin{eqnarray*}
&&\gamma _{2}\left[ B(x_{1}x_{2}\otimes
1_{H};G,gx_{1}x_{2})+B(x_{1}x_{2}\otimes 1_{H};GX_{1},gx_{2})\right] \left( %
\ref{X2,gx2,X2F11,gx2}\right) \\
&&+\lambda B(x_{1}x_{2}\otimes 1_{H};X_{1},gx_{1}x_{2}) \\
&&+\left[ +2B(x_{1}x_{2}\otimes 1_{H};1_{A},gx_{1})+2B(x_{1}x_{2}\otimes
1_{H};X_{1},g)\right] =0.
\end{eqnarray*}%
\begin{eqnarray*}
&&2\left[ -B(x_{1}x_{2}\otimes 1_{H};1_{A},gx_{1})-B(x_{1}x_{2}\otimes
1_{H};X_{1},g)\right] \left( \ref{X2,gx2,X2F41,g}\right) \\
&&+\gamma _{2}\left[ -B(x_{1}x_{2}\otimes
1_{H};G,gx_{1}x_{2})-B(x_{1}x_{2}\otimes 1_{H};GX_{1},gx_{2})\right] \\
&&-\lambda B(x_{1}x_{2}\otimes 1_{H};X_{1},gx_{1}x_{2})=0.
\end{eqnarray*}%
\begin{equation*}
-\gamma _{2}B(gx_{1}x_{2}\otimes 1_{H};G,1_{H})+\lambda \left[
+B(x_{2}\otimes 1_{H};1_{A},1_{H})-B(gx_{1}x_{2}\otimes 1_{H};1_{A},x_{1})%
\right] =0\left( \ref{X2,gx1x2,X2F11,1H}\right)
\end{equation*}%
\begin{equation*}
\gamma _{2}B(gx_{1}x_{2}\otimes 1_{H};G,x_{1}x_{2})-\lambda B(x_{2}\otimes
1_{H};1_{A},x_{1}x_{2})=0\left( \ref{X2,gx1x2,X2F11,x1x2}\right)
\end{equation*}%
\begin{gather*}
2\beta _{2}\left[ B(x_{1}\otimes 1_{H};1_{A},gx_{1})+B(gx_{1}x_{2}\otimes
1_{H};1_{A},gx_{1}x_{2})\right] \left( \ref{X2,gx1x2,X2F11,gx1}\right) \\
+\gamma _{2}B(gx_{1}x_{2}\otimes 1_{H};G,gx_{1})-\lambda B(x_{2}\otimes
1_{H};1_{A},gx_{1})=0.
\end{gather*}

\begin{gather*}
2\beta _{2}B(x_{1}\otimes 1_{H};1_{A},gx_{2})+\gamma
_{2}B(gx_{1}x_{2}\otimes 1_{H};G,gx_{2})\left( \ref{X2,gx1x2,X2F11,gx2}%
\right) \\
+\lambda \left[ -B(x_{2}\otimes 1_{H};1_{A},gx_{2})-B(gx_{1}x_{2}\otimes
1_{H};1_{A},gx_{1}x_{2})\right] \\
+2B(gx_{1}x_{2}\otimes 1_{H};1_{A},g)=0.
\end{gather*}

\begin{gather*}
2\beta _{2}\left[ -B(x_{1}\otimes 1_{H};G,g)-B(gx_{1}x_{2}\otimes
1_{H};G,gx_{2})\right] \left( \ref{X2,gx1x2,X2F21,g}\right) \\
+\lambda \left[ B(x_{2}\otimes 1_{H};G,g)-B(gx_{1}x_{2}\otimes
1_{H};G,gx_{1})\right] =0.
\end{gather*}%
\begin{gather*}
-2\beta _{2}B(x_{1}\otimes 1_{H};G,gx_{1}x_{2})\left( \ref%
{X2,gx1x2,X2F21,gx1x2}\right) \\
+\lambda B(x_{2}\otimes 1_{H};G,gx_{1}x_{2})-2B(gx_{1}x_{2}\otimes
1_{H};G,gx_{1})=0.
\end{gather*}%
\begin{gather*}
2\beta _{2}\left[
\begin{array}{c}
B(g\otimes 1_{H};1_{A},g)+B(x_{2}\otimes \ 1_{H};1_{A},gx_{2}) \\
+B(x_{1}\otimes 1_{H};1_{A},gx_{1})+B(gx_{1}x_{2}\otimes
1_{H};1_{A},gx_{1}x_{2})%
\end{array}%
\right] \left( \ref{X2,gx1x2,X2F31,g}\right) \\
-\gamma _{2}\left[ B(x_{2}\otimes 1_{H};G,g)-B(gx_{1}x_{2}\otimes
1_{H};G,gx_{1})\right] =0.
\end{gather*}

\begin{eqnarray*}
&&\gamma _{2}\left[ B(x_{1}\otimes 1_{H};G,g)+B(gx_{1}x_{2}\otimes
1_{H};G,gx_{2})\right] \left( \ref{X2,gx1x2,X2F41,g}\right) \\
&&-\lambda \left[
\begin{array}{c}
B(g\otimes 1_{H};1_{A},g)+B(x_{2}\otimes \ 1_{H};1_{A},gx_{2}) \\
+B(x_{1}\otimes 1_{H};1_{A},gx_{1})+B(gx_{1}x_{2}\otimes
1_{H};1_{A},gx_{1}x_{2})%
\end{array}%
\right] =0.
\end{eqnarray*}

\begin{eqnarray*}
&&2B(gx_{1}x_{2}\otimes 1_{H};1_{A},x_{1})+\gamma _{2}\left[ -B(x_{1}\otimes
1_{H};G,x_{1})-B(gx_{1}x_{2}\otimes 1_{H};G,x_{1}x_{2})\right] \left( \ref%
{X2,gx1x2,X2F41,x1}\right) \\
&&+\lambda \left[ B(g\otimes 1_{H};1_{A},x_{1})+B(x_{2}\otimes \
1_{H};1_{A},x_{1}x_{2})\right] =0.
\end{eqnarray*}

\begin{gather*}
2\left[ -B(x_{2}\otimes 1_{H};1_{A},1_{H})+B(gx_{1}x_{2}\otimes
1_{H};1_{A},x_{1})\right] +\left( \ref{X2,gx1x2,X2F51,1H}\right) \\
-\gamma _{2}\left[
\begin{array}{c}
B(g\otimes 1_{H};G,1_{H})+B(x_{2}\otimes \ 1_{H};G,x_{2}) \\
+B(x_{1}\otimes 1_{H};G,x_{1})+B(gx_{1}x_{2}\otimes 1_{H};G,x_{1}x_{2})%
\end{array}%
\right] =0.
\end{gather*}

\begin{equation*}
B(x_{2}\otimes 1_{H};G,g)-B(gx_{1}x_{2}\otimes 1_{H};G,gx_{1})=0.\left( \ref%
{X2,gx1x2,X2F61,gx2}\right)
\end{equation*}%
\begin{gather*}
2B(gx_{1}x_{2}\otimes 1_{H};G,1_{H})\left( \ref{X2,gx1x2,X2F71,1H}\right) \\
+\lambda \left[
\begin{array}{c}
B(g\otimes 1_{H};G,1_{H})+B(x_{2}\otimes \ 1_{H};G,x_{2}) \\
+B(x_{1}\otimes 1_{H};G,x_{1})+B(gx_{1}x_{2}\otimes 1_{H};G,x_{1}x_{2})%
\end{array}%
\right] =0.
\end{gather*}

\section{Non degenerate case}

From now on we will assume that all our constants $\alpha ,\beta ,\gamma
_{1},\gamma _{2}$ and $\lambda $ are different from zero.

\subsection{LIST\ OF\ ALL\ MONOMIAL\ EQUALITIES 0\label{LME0}}

We list all monomial equalities we got in $\left( \ref{LAE}\right) .$In
doing this, because of our previous assumptions, we will get rid of all the
constants and keep same label.
\begin{equation*}
B\left( g\otimes 1_{H};G,x_{1}x_{2}\right) =0.\left( \ref{X1,x1,X1F21,x1x2}%
\right)
\end{equation*}%
\begin{equation*}
B(x_{1}\otimes 1_{H};G,x_{2})=0.\left( \ref{X1,g,X1F21,x2}\right)
\end{equation*}%
\begin{equation*}
B(x_{1}\otimes 1_{H};G,x_{2})=0.\left( \ref{X1,x1,X1F11,x2}\right)
\end{equation*}%
\begin{equation*}
B(x_{1}\otimes 1_{H};G,x_{2})=0.\left( \ref{X1,x1,X1F41,1H}\right)
\end{equation*}%
\begin{equation*}
B(x_{1}x_{2}\otimes 1_{H};X_{1},gx_{1}x_{2})=0\text{ }.\left( \ref{G,x1x2,
GF4,gx1}\right)
\end{equation*}%
\begin{equation*}
B(x_{1}x_{2}\otimes 1_{H};X_{1},gx_{1}x_{2})=0\text{ }\left( \ref{G,x1x2,
GF6,gx1x2}\right)
\end{equation*}%
\begin{equation*}
B(x_{1}x_{2}\otimes 1_{H};X_{2},gx_{1}x_{2})=0\text{ }\left( \ref{G,x1x2,
GF3,gx2}\right)
\end{equation*}%
\begin{equation*}
B(x_{1}x_{2}\otimes 1_{H};X_{2},gx_{1}x_{2})=0.\left( \ref{G,x1x2, GF4,gx2}%
\right)
\end{equation*}%
\begin{equation*}
B(x_{1}x_{2}\otimes 1_{H};X_{2},gx_{1}x_{2})=0.\left( \ref{G,x1x2, GF4,gx2}%
\right)
\end{equation*}%
\begin{equation*}
B(x_{1}x_{2}\otimes 1_{H};X_{2},gx_{1}x_{2})=0\left( \ref{G,x1x2, GF7,gx1x2}%
\right)
\end{equation*}%
\begin{equation*}
B(x_{1}x_{2}\otimes 1_{H};GX_{1},1_{H})=0.\left( \ref{G,x1x2, GF3,1H}\right)
\end{equation*}

\begin{equation*}
B(x_{1}x_{2}\otimes 1_{H};GX_{1},x_{1}x_{2})=0.\left( \ref{G,x1x2, GF3,x1x2}%
\right)
\end{equation*}%
\begin{equation*}
B(x_{1}x_{2}\otimes 1_{H};GX_{1},x_{1}x_{2})=0\text{ }\left( \ref{G,x1x2,
GF3,gx1}\right)
\end{equation*}%
\begin{equation*}
B(x_{1}x_{2}\otimes 1_{H};GX_{1},x_{1}x_{2})=0.\left( \ref{G,x1x2, GF6,x1}%
\right)
\end{equation*}%
\begin{equation*}
B(x_{1}x_{2}\otimes 1_{H};GX_{1},x_{1}x_{2})=0\left( \ref{G,x1x2, GF7,x1}%
\right)
\end{equation*}%
\begin{equation*}
B(x_{1}x_{2}\otimes 1_{H};GX_{1},x_{1}x_{2})=0\left( \ref{G,gx2, GF2,x2}%
\right)
\end{equation*}%
\begin{equation*}
B(x_{1}x_{2}\otimes 1_{H};GX_{1},x_{1}x_{2})=0\left( \ref{X1,x1x2,X1F31,x1x2}%
\right)
\end{equation*}%
\begin{equation*}
B(x_{1}x_{2}\otimes 1_{H};GX_{2},x_{1}x_{2})=0.\left( \ref{G,x1x2, GF4,x1x2}%
\right)
\end{equation*}%
\begin{equation*}
B(x_{1}x_{2}\otimes 1_{H};GX_{2},x_{1}x_{2})=0.\left( \ref{G,x1x2, GF6,x2}%
\right)
\end{equation*}

\begin{equation*}
B(x_{1}x_{2}\otimes 1_{H};GX_{2},x_{1}x_{2})=0\left( \ref{G,x1x2, GF7,x2}%
\right)
\end{equation*}%
\begin{equation*}
B(x_{1}x_{2}\otimes 1_{H};GX_{2},gx_{2})=0.\left( \ref{X1,x1x2,X1F21,gx2}%
\right)
\end{equation*}%
\begin{equation*}
B(x_{1}x_{2}\otimes 1_{H};GX_{1},x_{1}x_{2})=0\left( \ref{X1,x1x2,X1F71,x1}%
\right)
\end{equation*}%
\begin{equation*}
B(x_{1}x_{2}\otimes 1_{H};GX_{1},x_{1}x_{2})=0.\left( \ref{X1,x1x2,X1F61,x1}%
\right)
\end{equation*}%
\begin{equation*}
B(x_{1}x_{2}\otimes 1_{H};GX_{1},x_{1}x_{2})=0.\left( \ref{X2,x1x2,X2F61,x1}%
\right)
\end{equation*}%
\begin{equation*}
B(x_{1}x_{2}\otimes 1_{H};GX_{1},gx_{1})=0.\left( \ref{X2,x1x2,X2F61,gx1x2}%
\right)
\end{equation*}%
\begin{equation*}
B(x_{1}x_{2}\otimes 1_{H};GX_{1},gx_{1})=0\left( \ref{X2,x1x2,X2F31,gx1}%
\right)
\end{equation*}%
\begin{equation*}
B(x_{1}x_{2}\otimes 1_{H};GX_{1},gx_{1})=0.\left( \ref{X2,x1x2,X2F21,gx1}%
\right)
\end{equation*}%
\begin{equation*}
B(x_{1}x_{2}\otimes 1_{H};GX_{1},gx_{1})=0.\left( \ref{X2,x1x2,X2F61,gx1x2}%
\right)
\end{equation*}%
\begin{equation*}
B(x_{1}x_{2}\otimes 1_{H};GX_{1},gx_{1})=0\left( \ref{X2,x1x2,X2F31,gx1}%
\right)
\end{equation*}%
\begin{equation*}
B(x_{1}x_{2}\otimes 1_{H};GX_{2},x_{1}x_{2})=0.\left( \ref%
{X2,x1x2,X2F21,x1x2}\right)
\end{equation*}%
\begin{equation*}
B(x_{1}x_{2}\otimes 1_{H};GX_{2},x_{1}x_{2})=0\left( \ref{X1,x1x2,X1F71,x2}%
\right)
\end{equation*}%
\begin{equation*}
B(x_{1}x_{2}\otimes 1_{H};GX_{2},x_{1}x_{2})=0\left( \ref{X1,x1x2,X1F61,x2}%
\right)
\end{equation*}%
\begin{equation*}
B(x_{1}x_{2}\otimes 1_{H};GX_{2},gx_{2})=0.\left( \ref{X1,x1x2,X1F71,gx1x2}%
\right)
\end{equation*}%
\begin{equation*}
B(x_{1}x_{2}\otimes 1_{H};GX_{2},gx_{2})=0\left( \ref{X1,x1x2,X1F41,gx2}%
\right)
\end{equation*}%
\begin{equation*}
B(gx_{2}\otimes 1_{H};GX_{1}X_{2},1_{H})=0.\left( \ref{X1,gx2,X171,1H}\right)
\end{equation*}

We make now a list without repetitions. When an equality is repeated, we
will take the one which appears as first in the paper and relabel it.%
\begin{equation}
B\left( g\otimes 1_{H};G,x_{1}x_{2}\right) =0.\left( \ref{X1,x1,X1F21,x1x2}%
\right)  \label{got1,G,x1x2}
\end{equation}%
\begin{equation}
B(x_{1}\otimes 1_{H};G,x_{2})=0.\left( \ref{X1,g,X1F21,x2}\right)
\label{x1ot1,G,x2}
\end{equation}%
\begin{equation}
B(x_{1}x_{2}\otimes 1_{H};X_{1},gx_{1}x_{2})=0\text{ }.\left( \ref{G,x1x2,
GF4,gx1}\right)  \label{x1x2ot1,X1,gx1x2}
\end{equation}%
\begin{equation}
B(x_{1}x_{2}\otimes 1_{H};X_{2},gx_{1}x_{2})=0\text{ }\left( \ref{G,x1x2,
GF3,gx2}\right)  \label{x1x2ot1,X2,gx1x2}
\end{equation}%
\begin{equation}
B(x_{1}x_{2}\otimes 1_{H};GX_{1},1_{H})=0.\left( \ref{G,x1x2, GF3,1H}\right)
\label{x1x2ot1,GX1,1}
\end{equation}%
\begin{equation}
B(x_{1}x_{2}\otimes 1_{H};GX_{1},x_{1}x_{2})=0\text{ }\left( \ref{G,x1x2,
GF3,gx1}\right)  \label{x1x2ot1,GX1,x1x2}
\end{equation}%
\begin{equation}
B(x_{1}x_{2}\otimes 1_{H};GX_{1},gx_{1})=0.\left( \ref{X2,x1x2,X2F21,gx1}%
\right)  \label{x1x2ot1,GX1,gx1}
\end{equation}%
\begin{equation}
B(x_{1}x_{2}\otimes 1_{H};GX_{2},x_{1}x_{2})=0.\left( \ref%
{X2,x1x2,X2F21,x1x2}\right)  \label{x1x2ot1,GX2,x1x2}
\end{equation}%
\begin{equation}
B(x_{1}x_{2}\otimes 1_{H};GX_{2},gx_{2})=0\left( \ref{X1,x1x2,X1F41,gx2}%
\right)  \label{x1x2ot1;GX2,gx2}
\end{equation}%
\begin{equation}
B(gx_{2}\otimes 1_{H};GX_{1}X_{2},1_{H})=0\left( \ref{X1,gx2,X171,1H}\right)
\label{gx2ot1,GX1X2,1}
\end{equation}%
Now we take the complete list of equalities $\left( \ref{LAE}\right) $ and
cancel there the terms above.

\subsection{LIST\ OF\ ALL\ EQUALITIES 1\label{LAE1}}

\subsubsection{$G$}

\begin{equation*}
\gamma _{1}B\left( g\otimes 1_{H};1_{A},x_{1}\right) +\gamma _{2}B\left(
g\otimes 1_{H};1_{A},x_{2}\right) =0\text{ }\left( \ref{G,g,GF1,1H}\right)
\end{equation*}

\begin{equation*}
2\alpha B(g\otimes 1_{H};G,gx_{1})+\gamma _{2}B(g\otimes
1_{H};1_{A},gx_{1}x_{2})=0\left( \ref{G,g, GF1,gx1}\right)
\end{equation*}

\begin{equation*}
2\alpha B(g\otimes 1_{H};G,gx_{2})-\gamma _{1}B\left( g\otimes
1_{H};1_{A},gx_{1}x_{2}\right) =0.\left( \ref{G,g, GF1,gx2}\right)
\end{equation*}

\begin{equation*}
\gamma _{1}B\left( g\otimes 1_{H};G,gx_{1}\right) +\gamma _{2}B\left(
g\otimes 1_{H};G,gx_{2}\right) =0.\left( \ref{G,g, GF2,g}\right)
\end{equation*}%
\begin{equation*}
B(g\otimes 1_{H};1_{A},x_{1})=0.\left( \ref{G,g, GF2,x1}\right)
\end{equation*}%
\begin{equation*}
B(g\otimes 1_{H};1_{A},x_{2})=0\left( \ref{G,g, GF2,x2}\right)
\end{equation*}%
\begin{equation*}
\begin{array}{c}
2\alpha B(x_{1}\otimes 1_{H};G,g)+ \\
+\gamma _{1}\left[ -B(g\otimes 1_{H};1_{A},g)-B(x_{1}\otimes
1_{H};1_{A},gx_{1})\right] -\gamma _{2}B(x_{1}\otimes 1_{H};1_{A},gx_{2})%
\end{array}%
=0.\left( \ref{G,x1, GF1,g}\right)
\end{equation*}%
\begin{equation*}
\gamma _{1}B(g\otimes 1_{H};1_{A},x_{1})+\gamma _{2}B(g\otimes
1_{H};1_{A},x_{2})=0.\left( \ref{G,x1, GF1,x1}\right)
\end{equation*}%
\begin{equation*}
\gamma _{1}[-B(g\otimes 1_{H};1_{A},x_{2})+B(x_{1}\otimes
1_{H};1_{A},x_{1}x_{2})]=0.\left( \ref{G,x1, GF1,x2}\right)
\end{equation*}

\begin{equation*}
-\gamma _{1}B(g\otimes 1_{H};1_{A},gx_{1}x_{2})+2\alpha B(x_{1}\otimes
1_{H};G,gx_{1}x_{2})=0.\left( \ref{G,x1, GF1,gx1x2}\right)
\end{equation*}%
\begin{equation*}
\begin{array}{c}
2B(x_{1}\otimes 1_{H};1_{A},1_{H})+ \\
+\gamma _{1}\left[ B(g\otimes 1_{H};G,1_{H})+B(x_{1}\otimes 1_{H};G,x_{1})%
\right]%
\end{array}%
=0.\left( \ref{G,x1, GF2,1H}\right)
\end{equation*}%
\begin{equation*}
2B(x_{1}\otimes 1_{H};1_{A},x_{1}x_{2})+\gamma _{1}B(g\otimes
1_{H};G,x_{1}x_{2})=0\text{ }\left( \ref{G,x1, GF2,x1x2}\right)
\end{equation*}%
\begin{equation*}
\gamma _{1}B\left( g\otimes 1_{H};G,gx_{1}\right) +\gamma _{2}B(x_{1}\otimes
1_{H};G,gx_{1}x_{2})=0.\left( \ref{G,x1, GF2,gx1}\right)
\end{equation*}

\begin{equation*}
\gamma _{1}[B\left( g\otimes 1_{H};G,gx_{2}\right) -B(x_{1}\otimes
1_{H};G,gx_{1}x_{2})]=0.\left( \ref{G,x1,GF2,gx2}\right)
\end{equation*}

\begin{equation*}
\gamma _{2}\left[ -B(g\otimes 1_{H};X_{2},1_{H})+B(x_{1}\otimes
1_{H};1_{A},x_{1}x_{2})\right] =0\left( \ref{G,x1,GF3,1H}\right)
\end{equation*}%
\begin{equation*}
-\gamma _{2}B(g\otimes 1_{H};1_{A},gx_{1}x_{2})-2\alpha B(g\otimes
1_{H};G,gx_{1})=0.\left( \ref{G,x1,GF3,gx1}\right)
\end{equation*}

\begin{equation*}
\alpha \left[ -B(g\otimes 1_{H};g,gx_{2})+B(x_{1}\otimes 1_{H};g,gx_{1}x_{2})%
\right] =0.\left( \ref{G,x1, GF3,x2}\right)
\end{equation*}%
\begin{equation*}
\gamma _{1}[B(g\otimes 1_{H};1_{A},x_{2})+B(x_{1}\otimes
1_{H};1_{A},x_{1}x_{2})]=0.\left( \ref{G,x1, GF4,1H}\right)
\end{equation*}

\begin{equation*}
-\gamma _{1}B(g\otimes 1_{H};1_{A},gx_{1}x_{2})+2\alpha B(x_{1}\otimes
1_{H};G,gx_{1}x_{2})=0.\left( \ref{G,x1, GF4,gx1}\right)
\end{equation*}

\begin{equation*}
\alpha \lbrack -B(g\otimes 1_{H};g,gx_{2})+B(x_{1}\otimes
1_{H};G,gx_{1}x_{2})]=0.\left( \ref{G,x1, GF5,g}\right)
\end{equation*}%
\begin{equation*}
2B(g\otimes 1_{H};1_{A},x_{1})-\gamma _{2}B(g\otimes
1_{H};G,x_{1}x_{2})=0.\left( \ref{G,x1,GF6,x1}\right)
\end{equation*}%
\begin{equation*}
B(g\otimes 1_{H};1_{A},x_{2})=B(x_{1}\otimes 1_{H};1_{A},x_{1}x_{2}).\left( %
\ref{G,x1,GF6,x2}\right)
\end{equation*}%
\begin{equation*}
\begin{array}{c}
2\alpha B(x_{2}\otimes 1_{H};G,g)+ \\
-\gamma _{1}B(x_{2}\otimes 1_{H};1_{A},gx_{1})+\gamma _{2}\left[ -B(g\otimes
1_{H};1_{A},g)-B(x_{2}\otimes \ 1_{H};1_{A},gx_{2})\right]%
\end{array}%
=0.\left( \ref{G,x2, GF1,g}\right)
\end{equation*}%
\begin{equation*}
\gamma _{2}\left[ -B(g\otimes 1_{H};1_{A},x_{1})-B(x_{2}\otimes \
1_{H};1_{A},x_{1}x_{2})\right] =0\left( \ref{G,x2, GF1,x1}\right)
\end{equation*}

\begin{equation*}
\gamma _{1}B(x_{2}\otimes 1_{H};1_{A},x_{1}x_{2})-\gamma _{2}B(g\otimes
1_{H};1_{A},x_{2})=0\left( \ref{G,x2, GF1,x2}\right)
\end{equation*}

\begin{equation*}
2\alpha B(x_{2}\otimes 1_{H};G,gx_{1}x_{2})-\gamma _{2}B\left( g\otimes
1_{H};1_{A},gx_{1}x_{2}\right) =0.\left( \ref{G,x2, GF1,gx1x2}\right)
\end{equation*}%
\begin{equation*}
\begin{array}{c}
2B(x_{2}\otimes 1_{H};1_{A},1_{H}) \\
+\gamma _{1}B(x_{2}\otimes 1_{H};G,x_{1})+\gamma _{2}\left[ B(g\otimes
1_{H};G,1_{H})+B(x_{2}\otimes \ 1_{H};G,x_{2}\right]%
\end{array}%
=0.\left( \ref{G,x2, GF2,1H}\right)
\end{equation*}%
\begin{equation*}
2B(x_{2}\otimes 1_{H};1_{A},x_{1}x_{2})=0=0\left( \ref{G,x2, GF2,x1x2}\right)
\end{equation*}%
\begin{equation*}
\gamma _{2}\left[ B(g\otimes 1_{H};G,gx_{1})+B(x_{2}\otimes \
1_{H};G,gx_{1}x_{2})\right] =0.\left( \ref{G,x2, GF2,gx1}\right)
\end{equation*}

\begin{equation*}
-\gamma _{1}B(x_{2}\otimes 1_{H};G,gx_{1}x_{2})+\gamma _{2}B(g\otimes
1_{H};G,gx_{2})=0.\left( \ref{G,x2, GF2,gx2}\right)
\end{equation*}%
\begin{equation*}
\gamma _{2}[B(x_{2}\otimes 1_{H};1_{A},x_{1}x_{2})+B(g\otimes
1_{H};1_{A},x_{1})]=0\left( \ref{G,x2, FG3,1H}\right)
\end{equation*}

\begin{equation*}
-\gamma _{2}B(g\otimes 1_{H};1_{A},gx_{1}x_{2})-2\alpha B(x_{2}\otimes
1_{H};G,gx_{1}x_{2})=0.\left( \ref{G,x2, GF3,gx2}\right)
\end{equation*}%
\begin{equation*}
\gamma _{1}\left[ B(g\otimes 1_{H};1_{A},x_{1})+B(x_{2}\otimes \
1_{H};1_{A},x_{1}x_{2})\right] =0.\left( \ref{G,x2, GF4,1H}\right)
\end{equation*}%
\begin{equation*}
\alpha \left[ B(g\otimes 1_{H};G,gx_{1})+B(x_{2}\otimes \
1_{H};G,gx_{1}x_{2})\right] =0\left( \ref{G,x2, GF4,gx1}\right)
\end{equation*}

\begin{equation*}
\gamma _{1}B(g\otimes 1_{H};1_{A},gx_{1}x_{2})-2\alpha B(g\otimes
1_{H};G,gx_{2})=0.\left( \ref{G,x2, GF4,gx2}\right)
\end{equation*}%
\begin{equation*}
\alpha \left[ B(g\otimes 1_{H};G,gx_{1})+B(x_{2}\otimes \
1_{H};G,gx_{1}x_{2})\right] =0\left( \ref{G,x2, GF5,g}\right)
\end{equation*}

\begin{equation*}
\gamma _{2}[B(g\otimes 1_{H},G,gx_{1})+B(x_{2}\otimes
1_{H};G,gx_{1}x_{2})=0.=0.\left( \ref{G,x2, GF6,g}\right)
\end{equation*}%
\begin{equation*}
2B(x_{2}\otimes 1_{H};1_{A},x_{1}x_{2})+\gamma _{2}B(g\otimes
1_{H};G,x_{1}x_{2})=0.\left( \ref{G,x2; GF6,x2}\right)
\end{equation*}

\begin{equation*}
\gamma _{1}[B(g\otimes 1_{H};G,gx_{1})+B(x_{2}\otimes
1_{H};G,gx_{1}x_{2})]=0.\left( \ref{G,x2; GF7,g}\right)
\end{equation*}%
\begin{equation*}
+B(g\otimes 1_{H};1_{A},x_{1})+B(x_{2}\otimes \
1_{H};1_{A},x_{1}x_{2})=0.\left( \ref{G,x2; GF7,x1}\right)
\end{equation*}%
\begin{equation*}
\gamma _{1}B(x_{1}x_{2}\otimes 1_{H};X_{1},g)+\gamma _{2}B(x_{1}x_{2}\otimes
1_{H};X_{2},g)=0.\left( \ref{G,x1x2, GF1,g}\right)
\end{equation*}

\begin{equation*}
\begin{array}{c}
2\alpha B(x_{1}x_{2}\otimes 1_{H};G,x_{1})+ \\
+\gamma _{1}B(x_{1}x_{2}\otimes 1_{H};X_{1},x_{1})+\gamma
_{2}B(x_{1}x_{2}\otimes 1_{H};X_{2},x_{1})%
\end{array}%
=0.\left( \ref{G,x1x2, GF1,x1}\right)
\end{equation*}%
\begin{equation*}
\begin{array}{c}
2\alpha B(x_{1}x_{2}\otimes 1_{H};G,x_{2})+ \\
+\gamma _{1}B(x_{1}x_{2}\otimes 1_{H};X_{1},x_{2})+\gamma
_{2}B(x_{1}x_{2}\otimes 1_{H};X_{2},x_{2})%
\end{array}%
=0.\left( \ref{G,x1x2, GF1,x2}\right)
\end{equation*}

\begin{equation*}
+\gamma _{2}B(x_{1}x_{2}\otimes 1_{H};GX_{2},1_{H})=0\left( \ref{G,x1x2,
GF2,1H}\right)
\end{equation*}%
\begin{equation*}
\begin{array}{c}
2B(x_{1}x_{2}\otimes 1_{H};1_{A},gx_{1}) \\
+\gamma _{2}B(x_{1}x_{2}\otimes 1_{H};GX_{2},gx_{1})%
\end{array}%
=0.\left( \ref{G,x1x2, GF2,gx1}\right)
\end{equation*}%
\begin{equation*}
\begin{array}{c}
2B(x_{1}x_{2}\otimes 1_{H};1_{A},gx_{2}) \\
+\gamma _{1}B(x_{1}x_{2}\otimes 1_{H};GX_{1},gx_{2})%
\end{array}%
=0\left( \ref{G,x1x2, GF2,gx2}\right)
\end{equation*}

\begin{equation*}
-2\alpha B(x_{1}x_{2}\otimes 1_{H};GX_{1},x_{1}x_{2})=0\text{ }\left( \ref%
{G,x1x2, GF3,gx1}\right)
\end{equation*}%
\begin{gather*}
-\gamma _{1}[1-B(x_{1}x_{2}\otimes
1_{H};1_{A},x_{1}x_{2})+B(x_{1}x_{2}\otimes
1_{H};x_{1},x_{2})-B(x_{1}x_{2}\otimes 1_{H};x_{2},x_{1})] \\
-2\alpha B(x_{1}x_{2}\otimes 1_{H};GX_{2},1_{H}=0.\left( \ref{G,x1x2, GF4,1H}%
\right)
\end{gather*}

\begin{gather*}
-2B(x_{1}x_{2}\otimes 1_{H};X_{1},g)+\gamma _{2}[-B(x_{1}x_{2}\otimes
1_{H};g,gx_{1}x_{2}) \\
-B(x_{1}x_{2}\otimes 1_{H};GX_{1},gx_{2})+B(x_{1}x_{2}\otimes
1_{H};GX_{2},gx_{1})]=0\left( \ref{G,x1x2, GF6,g}\right)
\end{gather*}%
\begin{gather*}
-2B(x_{1}x_{2}\otimes 1_{H};X_{2},g)+ \\
-\gamma _{1}[-B(x_{1}x_{2}\otimes 1_{H};g,gx_{1}x_{2})-B(x_{1}x_{2}\otimes
1_{H};GX_{1},gx_{2})+B(x_{1}x_{2}\otimes 1_{H};GX_{2},gx_{1})]=0\left( \ref%
{G,x1x2, GF7,g}\right)
\end{gather*}

\begin{equation*}
\begin{array}{c}
2\alpha \left[ B(x_{1}x_{2}\otimes 1_{H};G,x_{2})-B(x_{1}x_{2}\otimes
1_{H};GX_{2},1_{H})\right] + \\
+\gamma _{1}\left[ -1+B(x_{1}x_{2}\otimes
1_{H};1_{A},x_{1}x_{2})+B(x_{1}x_{2}\otimes 1_{H};X_{2},x_{1})\right]
+\gamma _{2}B(x_{1}x_{2}\otimes 1_{H};X_{2},x_{2})%
\end{array}%
=0.\left( \ref{G,gx1, GF1,1H}\right)
\end{equation*}

\begin{equation*}
\begin{array}{c}
2\left[ B(x_{1}x_{2}\otimes 1_{H};1_{A},gx_{2})+B(x_{1}x_{2}\otimes
1_{H};X_{2},g)\right] \\
\gamma _{1}\left[ -B(x_{1}x_{2}\otimes
1_{H};G,gx_{1}x_{2})+B(x_{1}x_{2}\otimes 1_{H};GX_{2},gx_{1})\right]%
\end{array}%
\left( \ref{G,gx1, GF2,g}\right)
\end{equation*}%
\begin{equation*}
\begin{array}{c}
2\alpha \left[ -B(x_{1}x_{2}\otimes 1_{H};G,x_{1})\right] + \\
-\gamma _{1}B(x_{1}x_{2}\otimes 1_{H};X_{1},x_{1})+\gamma _{2}\left[
-1+B(x_{1}x_{2}\otimes 1_{H};1_{A},x_{1}x_{2})-B(x_{1}x_{2}\otimes
1_{H};X_{1},x_{2})\right]%
\end{array}%
=0\left( \ref{G,gx2, GF1,1H}\right)
\end{equation*}%
\begin{equation*}
\begin{array}{c}
2\left[ -B(x_{1}x_{2}\otimes 1_{H};1_{A},gx_{1})-B(x_{1}x_{2}\otimes
1_{H};X_{1},g)\right] \\
+\gamma _{2}\left[ -B(x_{1}x_{2}\otimes
1_{H};G,gx_{1}x_{2})-B(x_{1}x_{2}\otimes 1_{H};GX_{1},gx_{2})\right]%
\end{array}%
=0.\left( \ref{G,gx2, GF2,g}\right)
\end{equation*}

\begin{gather*}
\gamma _{1}\left[ -B(x_{2}\otimes 1_{H};1_{A},1_{H})+B(gx_{1}x_{2}\otimes
1_{H};1_{A},x_{1})\right] \\
+\gamma _{2}\left[ B(x_{1}\otimes 1_{H};1_{A},1_{H})+B(gx_{1}x_{2}\otimes
1_{H};1_{A},x_{2})\right] =0\left( \ref{G,gx1x2, GF1,1H}\right)
\end{gather*}%
\begin{equation*}
-\gamma _{1}B(x_{2}\otimes 1_{H};1_{A},x_{1}x_{2})+\gamma _{2}B(x_{1}\otimes
1_{H};1_{A},x_{1}x_{2})=0.\left( \ref{G,gx1x2, GF1,x1x2}\right)
\end{equation*}%
\begin{equation*}
\begin{array}{c}
2\alpha B(x_{2}\otimes 1_{H};G,g)+ \\
-\gamma _{1}B(x_{2}\otimes 1_{H};1_{A},gx_{1})+\gamma _{2}\left[
B(x_{1}\otimes 1_{H};1_{A},gx_{1})+B(gx_{1}x_{2}\otimes
1_{H};1_{A},gx_{1}x_{2})\right]%
\end{array}%
=0.\left( \ref{G,gx1x2, GF1,gx1}\right)
\end{equation*}

\begin{equation*}
\begin{array}{c}
2\alpha B(gx_{1}x_{2}\otimes 1_{H};G,gx_{2})+ \\
+\gamma _{1}\left[ -B(x_{2}\otimes 1_{H};1_{A},gx_{2})-B(gx_{1}x_{2}\otimes
1_{H};1_{A},gx_{1}x_{2})\right] +\gamma _{2}B(x_{1}\otimes
1_{H};1_{A},gx_{2})%
\end{array}%
=0.\left( \ref{G,gx1x2, GF1,gx2}\right)
\end{equation*}%
\begin{gather*}
\gamma _{1}\left[ B(x_{2}\otimes 1_{H};G,g)-B(gx_{1}x_{2}\otimes
1_{H};G,gx_{1})\right] \left( \ref{G,gx1x2, GF2,g}\right) \\
+\gamma _{2}\left[ -B(x_{1}\otimes 1_{H};G,g)-B(gx_{1}x_{2}\otimes
1_{H};G,gx_{2})\right] =0.
\end{gather*}%
\begin{equation*}
\begin{array}{c}
2B(gx_{1}x_{2}\otimes 1_{H};1_{A},x_{1}) \\
+\gamma _{1}B(x_{2}\otimes 1_{H};G,x_{1})+\gamma _{2}\left[ -B(x_{1}\otimes
1_{H};G,x_{1})-B(gx_{1}x_{2}\otimes 1_{H};G,x_{1}x_{2})\right]%
\end{array}%
=0.\left( \ref{G,gx1x2, GF2,x1}\right)
\end{equation*}%
\begin{equation*}
\begin{array}{c}
2B(gx_{1}x_{2}\otimes 1_{H};1_{A},x_{2}) \\
+\gamma _{1}\left[ B(x_{2}\otimes 1_{H};G,x_{2})+B(gx_{1}x_{2}\otimes
1_{H};G,x_{1}x_{2})\right]%
\end{array}%
=0.\left( \ref{G,gx1x2, GF2,x2}\right)
\end{equation*}%
\begin{equation*}
\gamma _{1}B(x_{2}\otimes 1_{H};G,gx_{1}x_{2})-\gamma _{2}B(x_{1}\otimes
1_{H};G,gx_{1}x_{2})=0.\left( \ref{G,gx1x2, GF2,gx1x2}\right)
\end{equation*}%
\begin{gather*}
\gamma _{2}\left[
\begin{array}{c}
B(g\otimes 1_{H};1_{A},g)+B(x_{2}\otimes \ 1_{H};1_{A},gx_{2})+ \\
+B(x_{1}\otimes 1_{H};1_{A},gx_{1})+B(gx_{1}x_{2}\otimes
1_{H};1_{A},gx_{1}x_{2})%
\end{array}%
\right] +\left( \ref{G,gx1x2, GF3,g}\right) \\
2\alpha \left[ -B(x_{2}\otimes 1_{H};G,g)+B(x_{2}\otimes 1_{H};G,g)\right]
=0.
\end{gather*}

\begin{equation*}
\gamma _{2}\left[ B(g\otimes 1_{H};1_{A},x_{1})+B(x_{2}\otimes \
1_{H};1_{A},x_{1}x_{2})\right] =0.\left( \ref{G,gx1x2, GF3,x1}\right)
\end{equation*}

\begin{equation*}
\gamma _{2}\left[ B(g\otimes 1_{H};X_{2},1_{H})-B(x_{1}\otimes
1_{H};1_{A},x_{1}x_{2})\right] =0\left( \ref{G,gx1x2, GF3,x2}\right)
\end{equation*}

\begin{equation*}
\gamma _{2}B\left( g\otimes 1_{H};1_{A},gx_{1}x_{2}\right) -2\alpha
B(x_{2}\otimes 1_{H};G,gx_{1}x_{2})=0.\left( \ref{G,gx1x2, GF3,gx1x2}\right)
\end{equation*}

\begin{gather*}
-\gamma _{1}\left[
\begin{array}{c}
B(g\otimes 1_{H};1_{A},g)+B(x_{2}\otimes \ 1_{H};1_{A},gx_{2})+ \\
+B(x_{1}\otimes 1_{H};1_{A},gx_{1})+B(gx_{1}x_{2}\otimes
1_{H};1_{A},gx_{1}x_{2})%
\end{array}%
\right] +\left( \ref{G,gx1x2, GF4,g}\right) \\
-2\alpha \left[ -B(x_{1}\otimes 1_{H};G,g)-B(gx_{1}x_{2}\otimes
1_{H};G,gx_{2})\right] =0
\end{gather*}%
\begin{equation*}
\gamma _{1}\left[ B(g\otimes 1_{H};1_{A},x_{1})+B(x_{2}\otimes \
1_{H};1_{A},x_{1}x_{2})\right] =0\left( \ref{G,gx1x2, GF4,x1}\right)
\end{equation*}%
\begin{equation*}
-\gamma _{1}\left[ B(g\otimes 1_{H};1_{A},x_{2})-B(x_{1}\otimes
1_{H};1_{A},x_{1}x_{2})\right] =0.\left( \ref{G,gx1x2, GF4,x2}\right)
\end{equation*}

\begin{equation*}
-\gamma _{1}B\left( g\otimes 1_{H};1_{A},gx_{1}x_{2}\right) +2\alpha
B(x_{1}\otimes 1_{H};G,gx_{1}x_{2})=0.\left( \ref{G,gx1x2, GF4,gx1x2}\right)
\end{equation*}

\begin{equation*}
\alpha \left[ B(g\otimes 1_{H};G,gx_{2})-B(x_{1}\otimes 1_{H};G,gx_{1}x_{2})%
\right] =0.\left( \ref{G,gx1x2, GF5,gx2}\right)
\end{equation*}

\begin{gather*}
-2\left[ -B(x_{2}\otimes 1_{H};1_{A},1_{H})+B(gx_{1}x_{2}\otimes
1_{H};1_{A},x_{1})\right] +\left( \ref{G,gx1x2, GF6,1H}\right) \\
+\gamma _{2}\left[
\begin{array}{c}
B(g\otimes 1_{H};G,1_{H})+B(x_{2}\otimes \ 1_{H};G,x_{2})+ \\
+B(x_{1}\otimes 1_{H};G,x_{1})+B(gx_{1}x_{2}\otimes 1_{H};G,x_{1}x_{2})%
\end{array}%
\right] \\
=0
\end{gather*}%
\begin{equation*}
\gamma _{2}\left[ B(g\otimes 1_{H};G,gx_{2})-B(x_{1}\otimes
1_{H};G,gx_{1}x_{2})\right] =0.\left( \ref{G,gx1x2, GF6,gx1}\right)
\end{equation*}%
\begin{equation*}
\gamma _{2}\left[ B(g\otimes 1_{H};G,gx_{2})-B(x_{1}\otimes
1_{H};G,gx_{1}x_{2})\right] =0.\left( \ref{G,gx1x2, GF6,gx2}\right)
\end{equation*}%
\begin{gather*}
-\gamma _{1}\left[
\begin{array}{c}
B(g\otimes 1_{H};G,1_{H})+B(x_{2}\otimes \ 1_{H};G,x_{2}) \\
+B(x_{1}\otimes 1_{H};G,x_{1})+B(gx_{1}x_{2}\otimes 1_{H};G,x_{1}x_{2})%
\end{array}%
\right] + \\
-2\left[ B(x_{1}\otimes 1_{H};1_{A},1_{H})+B(gx_{1}x_{2}\otimes
1_{H};1_{A},x_{2})\right] =0\left( \ref{G,gx1x2, GF7,1H}\right)
\end{gather*}%
\begin{equation*}
-\gamma _{1}B(g\otimes 1_{H};G,x_{1}x_{2})-2B(x_{1}\otimes
1_{H};1_{A},x_{1}x_{2})=0.\left( \ref{G,gx1x2, GF7,x1x2}\right)
\end{equation*}%
\begin{equation*}
\gamma _{1}\left[ B(g\otimes 1_{H};G,gx_{1})+B(x_{2}\otimes \
1_{H};G,gx_{1}x_{2})\right] =0.\left( \ref{G,gx1x2, GF7,gx1}\right)
\end{equation*}%
\begin{equation*}
\gamma _{1}\left[ B(g\otimes 1_{H};G,gx_{2})-B(x_{1}\otimes
1_{H};G,gx_{1}x_{2})\right] =0.\left( \ref{G,gx1x2, GF7,gx2}\right)
\end{equation*}

\subsubsection{$X_{1}$}

\begin{equation*}
\lambda B\left( g\otimes 1_{H};1_{A},x_{2}\right) -\gamma _{1}B(g\otimes
1_{H};G,1_{H})-2B(x_{1}\otimes 1_{H};1_{A},1_{H})=0.\left( \ref{X1,g,
X1F11,1H}\right)
\end{equation*}%
\begin{eqnarray*}
&&+\gamma _{1}B(g\otimes 1_{H};G,gx_{1})+\lambda B(g\otimes
1_{H};X_{2},gx_{1})+B(x_{1}\otimes 1_{H};1_{A},gx_{1})+\left( \ref{X1,g,
X1F11,gx1}\right) \\
&&+2B\left( g\otimes 1_{H};1_{A},g\right) +B\left( x_{1}\otimes
1_{H};1_{A},gx_{1}\right) =0
\end{eqnarray*}

\begin{equation*}
\gamma _{1}B(g\otimes 1_{H};G,gx_{2})+2B(x_{1}\otimes
1_{H};1_{A},gx_{2})-2\beta _{1}B\left( g\otimes
1_{H};1_{A},gx_{1}x_{2}\right) =0.\left( \ref{X1,g, X1F11,gx2}\right)
\end{equation*}

\begin{equation*}
2\beta _{1}B\left( g\otimes 1_{H};G,gx_{1}\right) +\lambda B\left( g\otimes
1_{H};G,gx_{2}\right) +2B(x_{1}\otimes 1_{H};G,g)=0.\left( \ref{X1,g,X1F21,g}%
\right)
\end{equation*}%
\begin{equation*}
+2B(x_{1}\otimes 1_{H};G,x_{1})=0.\left( \ref{X1,g,X1F21,x1}\right)
\end{equation*}

\begin{equation*}
B\left( g\otimes 1_{H};G,gx_{2}\right) -B\left( x_{1}\otimes
1_{H};G,gx_{1}x_{2}\right) =0\left( \ref{X1,g,X1F21,gx1x2}\right)
\end{equation*}%
\begin{gather*}
+\lambda B\left( g\otimes 1_{H};1_{A},gx_{1}x_{2}\right) +\gamma _{1}B\left(
g\otimes 1_{H};G,gx_{1}\right) +\left( \ref{X1,g, X1F31,g}\right) \\
+2B(g\otimes 1_{H};1_{A},g)+2B(x_{1}\otimes 1_{H};1_{A},gx_{1})=0.
\end{gather*}

\begin{equation*}
-2\beta _{1}B\left( g\otimes 1_{H};1_{A},gx_{1}x_{2}\right) +2B(x_{1}\otimes
1_{H};1_{A},gx_{2})+\gamma _{1}B\left( g\otimes 1_{H};G,gx_{2}\right)
=0.\left( \ref{X1,g, X1F41,g}\right)
\end{equation*}%
\begin{equation*}
2B(x_{1}\otimes 1_{H};1_{A},x_{1}x_{2})=0.\left( \ref{X1,g, X1F41,x1}\right)
\end{equation*}

\begin{equation*}
-2B\left( g\otimes 1_{H};1_{A},x_{2}\right) =-2B\left( g\otimes
1_{H};1_{A},x_{2}\right) +2B(x_{1}\otimes 1_{H};1_{A},x_{1}x_{2})\left( \ref%
{X1,g,X1F51,1H}\right)
\end{equation*}%
\begin{equation*}
B(x_{1}\otimes 1_{H};G,gx_{1}x_{2})-B(g\otimes 1_{H};G,gx_{2})=0.\left( \ref%
{X1,x2,X1F71,gx1x2}\right)
\end{equation*}%
\begin{gather*}
2\beta _{1}\left[ -B(g\otimes 1_{H};1_{A},g)-B(x_{1}\otimes
1_{H};1_{A},gx_{1})\right] + \\
-\lambda B(x_{1}\otimes 1_{H};1_{A},gx_{2})=\gamma _{1}B(x_{1}\otimes
1_{H};G,g)\left( \ref{X1,x1,X1F11,g}\right)
\end{gather*}%
\begin{equation*}
-\lambda B(x_{1}\otimes 1_{H};1_{A},x_{1}x_{2})+\gamma _{1}B(x_{1}\otimes
1_{H};G,x_{1})=0.\left( \ref{X1,x1,X1F11,x1}\right)
\end{equation*}

\begin{equation*}
-2\beta _{1}B(g\otimes 1_{H};1_{A},gx_{1}x_{2})=\gamma _{1}B(x_{1}\otimes
1_{H};G,gx_{1}x_{2})-2B\left( x_{1}\otimes 1_{H};1_{A},gx_{2}\right) .\left( %
\ref{X1,x1,X1F11,gx1x2}\right)
\end{equation*}%
\begin{equation*}
\beta _{1}\left[ B(g\otimes 1_{H};G,1_{H})+B(x_{1}\otimes 1_{H};G,x_{1})%
\right] =0\left( \ref{X1,x1,X1F21,1}\right)
\end{equation*}%
\begin{equation*}
\lambda B(x_{1}\otimes 1_{H};G,gx_{1}x_{2})=-2B\left( x_{1}\otimes
1_{H};G,g\right) .\left( \ref{X1,x1,X1F21,gx1}\right)
\end{equation*}%
\begin{gather*}
2B(x_{1}\otimes 1_{H};1_{A},1_{H})+\lambda \left[ -B\left( g\otimes
1_{H};1_{A},x_{2}\right) +B(x_{1}\otimes 1_{H};1_{A},x_{1}x_{2})\right]
\left( \ref{X1,x1,X1F31,1H}\right) \\
+\gamma _{1}\left[ B(g\otimes 1_{H};G,1_{H})+B(x_{1}\otimes 1_{H};G,x_{1})%
\right] =0.
\end{gather*}%
\begin{equation*}
2B\left( x_{1}\otimes 1_{H};1_{A},x_{1}x_{2}\right) =0\left( \ref%
{X1,x1,X1F31,x1x2}\right)
\end{equation*}

\begin{gather*}
-\lambda B\left( g\otimes 1_{H};1_{A},gx_{1}x_{2}\right) -\gamma _{1}B\left(
g\otimes 1_{H};G,gx_{1}\right) + \\
-2\left[ B(g\otimes 1_{H};1_{A},g)+B(x_{1}\otimes 1_{H};1_{A},gx_{1})\right]
=0\left( \ref{X1,x1,X1F31,gx1}\right)
\end{gather*}

\begin{equation*}
\gamma _{1}\left[ -B\left( g\otimes 1_{H};G,gx_{2}\right) +B(x_{1}\otimes
1_{H};G,gx_{1}x_{2})\right] =0\left( \ref{X1,x1,X1F31,gx2}\right)
\end{equation*}

\begin{equation*}
2\beta _{1}B\left( g\otimes 1_{H};1_{A},gx_{1}x_{2}\right) -\gamma
_{1}B(x_{1}\otimes 1_{H};G,gx_{1}x_{2})-B(x_{1}\otimes
x_{1};X_{2},gx_{1})=0.\left( \ref{X1,x1,X1F41,gx1}\right)
\end{equation*}

\begin{eqnarray*}
&&-2\beta _{1}B(x_{2}\otimes 1_{H};1_{A},gx_{1})+\lambda \left[ -B(g\otimes
1_{H};1_{A},g)-B(x_{2}\otimes \ 1_{H};1_{A},gx_{2})\right] \\
&=&-\gamma _{1}B(x_{2}\otimes 1_{H};G,g)-2B(gx_{1}x_{2}\otimes
1_{H};1_{A},g).\left( \ref{X1,x2,X1F11,g}\right)
\end{eqnarray*}%
\begin{eqnarray*}
&&\lambda \left[ -B(g\otimes 1_{H};1_{A},x_{1})-B(x_{2}\otimes \
1_{H};1_{A},x_{1}x_{2})\right] \\
&=&\gamma _{1}B(x_{2}\otimes 1_{H};G,x_{1})+2B(gx_{1}x_{2}\otimes
1_{H};1_{A},x_{1}).\left( \ref{X1,x2,X1F11,x1}\right)
\end{eqnarray*}%
\begin{equation*}
-\lambda B(g\otimes 1_{H};1_{A},x_{2})=\gamma _{1}B(x_{2}\otimes
1_{H};G,x_{2})+2B(gx_{1}x_{2}\otimes 1_{H};1_{A},x_{2}).\left( \ref%
{X1,x2,X1F11,x2}\right)
\end{equation*}

\begin{gather*}
\lambda B\left( g\otimes 1_{H};1_{A},gx_{1}x_{2}\right) -\gamma
_{1}B(x_{2}\otimes 1_{H};G,gx_{1}x_{2})\left( \ref{X1,x2,X1F11,gx1x2}\right)
\\
-2B\left( x_{2}\otimes 1_{H};1_{A},gx_{2}\right) -2B(gx_{1}x_{2}\otimes
1_{H};1_{A},gx_{1}x_{2})=0.
\end{gather*}%
\begin{equation*}
\lambda \left[ B(g\otimes 1_{H};G,1_{H})+B(x_{2}\otimes \ 1_{H};G,x_{2}%
\right] =-2B(gx_{1}x_{2}\otimes 1_{H};G,1_{H}).\left( \ref{X1,x2,X1F21,1H}%
\right)
\end{equation*}%
\begin{equation*}
+2B(gx_{1}x_{2}\otimes 1_{H};G,x_{1}x_{2})=0.\left( \ref{X1,x2,X1F21,x1x2}%
\right)
\end{equation*}%
\begin{eqnarray*}
&&\lambda \left[ B(g\otimes 1_{H};G,gx_{1})+B(x_{2}\otimes \
1_{H};G,gx_{1}x_{2})\right] \\
&=&2B(gx_{1}x_{2}\otimes 1_{H};G,gx_{1})-2B\left( x_{2}\otimes
1_{H};G,g\right) .\left( \ref{X1,x2,X1F21,gx1}\right)
\end{eqnarray*}

\begin{equation*}
-2\beta _{1}B(x_{2}\otimes 1_{H};G,gx_{1}x_{2})+\lambda B(g\otimes
1_{H};G,gx_{2})-2B(gx_{1}x_{2}\otimes 1_{H};G,gx_{2})=0.\left( \ref%
{X1,x2,X1F21,gx2}\right)
\end{equation*}%
\begin{gather*}
-\lambda B\left( g\otimes 1_{H};1_{A},gx_{1}x_{2}\right) =-\gamma
_{1}B(x_{2}\otimes 1_{H};G,gx_{1}x_{2})+ \\
-2B(x_{2}\otimes 1_{H};1_{A},gx_{2})-2B(gx_{1}x_{2}\otimes
1_{H};1_{A},gx_{1}x_{2}).\left( \ref{X1,x2,X1F31,gx2}\right)
\end{gather*}

\begin{equation*}
0=-\gamma _{1}\left[ B(g\otimes 1_{H};G,1_{H})+B(x_{2}\otimes \ 1_{H};G,x_{2}%
\right] -2B(x_{1}\otimes 1_{H};1_{A},1_{H})-2B(gx_{1}x_{2}\otimes
1_{H};1_{A},x_{2}).\left( \ref{X1,x2,X1F41,1H}\right)
\end{equation*}%
\begin{equation*}
0=-2B(x_{1}\otimes 1_{H};1_{A},x_{1}x_{2}).\left( \ref{X1,x2,X1F41,x1x2}%
\right)
\end{equation*}%
\begin{gather*}
B(x_{1}\otimes 1_{H};1_{A},gx_{1})+B(gx_{1}x_{2}\otimes
1_{H};1_{A},gx_{1}x_{2})\left( \ref{X1,x2,X1F41,gx1}\right) \\
+B(g\otimes 1_{H};1_{A},g)+B(x_{2}\otimes \ 1_{H};1_{A},gx_{2})=0.
\end{gather*}

\begin{equation*}
2\beta _{1}B\left( g\otimes 1_{H};1_{A},gx_{1}x_{2}\right) =\gamma
_{1}B(g\otimes 1_{H};G,gx_{2})+2B(x_{1}\otimes 1_{H};1_{A},gx_{2}).\left( %
\ref{X1,x2,X1F41,gx2}\right)
\end{equation*}

\begin{eqnarray*}
0 &=&\gamma _{1}\left[ B(g\otimes 1_{H};G,gx_{1})+B(x_{2}\otimes \
1_{H};G,gx_{1}x_{2})\right] \left( \ref{X1,x2,X1F51,g}\right) \\
&&+\left[
\begin{array}{c}
2B(g\otimes 1_{H};1_{A},g)+2B(x_{2}\otimes \ 1_{H};1_{A},gx_{2})+ \\
+2B(x_{1}\otimes 1_{H};1_{A},gx_{1})+2B(gx_{1}x_{2}\otimes
1_{H};1_{A},gx_{1}x_{2})%
\end{array}%
\right]
\end{eqnarray*}

\begin{equation*}
0=+2B(x_{1}\otimes 1_{H};1_{A},x_{1}x_{2}).\left( \ref{X1,x2,X1F51,x2}\right)
\end{equation*}%
\begin{equation*}
\lambda \left[ B(g\otimes 1_{H};G,gx_{1})+B(x_{2}\otimes \
1_{H};G,gx_{1}x_{2})\right] =-2\left[ B(x_{2}\otimes
1_{H};G,g)-B(gx_{1}x_{2}\otimes 1_{H};G,gx_{1})\right] \left( \ref%
{X1,x2,X1F61,g}\right)
\end{equation*}%
\begin{equation*}
\beta _{1}\left[ B(g\otimes 1_{H};G,gx_{1})+B(x_{2}\otimes \
1_{H};G,gx_{1}x_{2})\right] =0\left( \ref{X1,x2,X1F71,g}\right)
\end{equation*}%
\begin{equation*}
B(x_{1}\otimes 1_{H};G,x_{1})+B(gx_{1}x_{2}\otimes
1_{H};G,x_{1}x_{2})=0\left( \ref{X1,x2,X1F71,x1}\right)
\end{equation*}%
\begin{equation*}
\lambda B(x_{1}x_{2}\otimes 1_{H};X_{2},g)-\gamma _{1}B(x_{1}x_{2}\otimes
1_{H};G,g)=0\left( \ref{X1,x1x2,X1F11,g}\right)
\end{equation*}%
\begin{equation*}
2\beta _{1}B(x_{1}x_{2}\otimes 1_{H};X_{1},x_{1})+\lambda
B(x_{1}x_{2}\otimes 1_{H};X_{2},x_{1})+\gamma _{1}B(x_{1}x_{2}\otimes
1_{H};G,x_{1})=0\left( \ref{X1,x1x2,X1F11,x1}\right)
\end{equation*}

\begin{equation*}
2\beta _{1}B(x_{1}x_{2}\otimes 1_{H};X_{1},x_{2})+\lambda
B(x_{1}x_{2}\otimes 1_{H};X_{2},x_{2})+\gamma _{1}B(x_{1}x_{2}\otimes
1_{H};G,x_{2})=0.\left( \ref{X1,x1x2,X1F11,x2}\right)
\end{equation*}

\begin{equation*}
\gamma _{1}B(x_{1}x_{2}\otimes 1_{H};G,gx_{1}x_{2})-2B(x_{1}x_{2}\otimes
1_{H};1_{A},gx_{2})=0.\left( \ref{X1,x1x2,X1F11,gx1x2}\right)
\end{equation*}%
\begin{equation*}
+\lambda B(x_{1}x_{2}\otimes 1_{H};GX_{2},1_{H})=0.\left( \ref%
{X1,x1x2,X1F21,1H}\right)
\end{equation*}%
\begin{equation*}
\lambda B(x_{1}x_{2}\otimes 1_{H};GX_{2},gx_{1})+2B(x_{1}x_{2}\otimes
1_{H};G,g)=0.\left( \ref{X1,x1x2,X1F21,gx1}\right)
\end{equation*}

\begin{equation*}
\lambda \left[
\begin{array}{c}
+1-B(x_{1}x_{2}\otimes 1_{H};1_{A},x_{1}x_{2}) \\
-B(x_{1}x_{2}\otimes 1_{H};X_{2},x_{1})+B(x_{1}x_{2}\otimes
1_{H};X_{1},x_{2})%
\end{array}%
\right] =0\left( \ref{X1,x1x2,X1F31,1H}\right)
\end{equation*}%
\begin{equation*}
2B(x_{1}x_{2}\otimes 1_{H};1_{A},gx_{1})+2B(x_{1}x_{2}\otimes
1_{H};X_{1},g)=0\left( \ref{X1,x1x2,X1F31,gx1}\right)
\end{equation*}%
\begin{equation*}
2B(x_{1}x_{2}\otimes 1_{H};1_{A},gx_{2})+\gamma _{1}B(x_{1}x_{2}\otimes
1_{H};GX_{1},gx_{2})=0.\left( \ref{X1,x1x2,X1F31,gx2}\right)
\end{equation*}

\begin{gather*}
2\beta _{1}\left[
\begin{array}{c}
+1-B(x_{1}x_{2}\otimes 1_{H};1_{A},x_{1}x_{2})-B(x_{1}x_{2}\otimes
1_{H};X_{2},x_{1}) \\
+B(x_{1}x_{2}\otimes 1_{H};X_{1},x_{2})%
\end{array}%
\right] \\
+\gamma _{1}B(x_{1}x_{2}\otimes 1_{H};GX_{2},1_{H})=0.\left( \ref%
{X1,x1x2,X1F41,1H}\right)
\end{gather*}

\begin{equation*}
\gamma _{1}B(x_{1}x_{2}\otimes 1_{H};GX_{2},gx_{1})+2B(x_{1}x_{2}\otimes
1_{H};X_{2},g)=0.\left( \ref{X1,x1x2,X1F41,gx1}\right)
\end{equation*}%
\begin{gather*}
2B(x_{1}x_{2}\otimes 1_{H};X_{2},g)+\left( \ref{X1,x1x2,X1F51,g}\right) \\
+\gamma _{1}\left[ -B(x_{1}x_{2}\otimes
1_{H};G,gx_{1}x_{2})+B(x_{1}x_{2}\otimes
1_{H};GX_{2},gx_{1})-B(x_{1}x_{2}\otimes 1_{H};GX_{1},gx_{2})\right] =0
\end{gather*}%
\begin{equation*}
2B(x_{1}x_{2}\otimes 1_{H};G,g)+\lambda \left[
\begin{array}{c}
-B(x_{1}x_{2}\otimes 1_{H};G,gx_{1}x_{2}) \\
+B(x_{1}x_{2}\otimes 1_{H};GX_{2},gx_{1})-B(x_{1}x_{2}\otimes
1_{H};GX_{1},gx_{2})%
\end{array}%
\right] =0\left( \ref{X1,x1x2,X1F61,g}\right)
\end{equation*}%
\begin{equation*}
B(x_{1}x_{2}\otimes 1_{H};G,gx_{1}x_{2})+B(x_{1}x_{2}\otimes
1_{H};GX_{1},gx_{2})=0\left( \ref{X1,x1x2,X1F61,gx1x2}\right)
\end{equation*}%
\begin{equation*}
+B(x_{1}x_{2}\otimes 1_{H};G,gx_{1}x_{2})+B(x_{1}x_{2}\otimes
1_{H};GX_{1},gx_{2})=0.\left( \ref{X1,x1x2,X1F81,gx1}\right)
\end{equation*}%
\begin{gather*}
2\beta _{1}\left[ -1+B(x_{1}x_{2}\otimes
1_{H};1_{A},x_{1}x_{2})+B(x_{1}x_{2}\otimes 1_{H};X_{2},x_{1})\right]
+\lambda B(x_{1}x_{2}\otimes 1_{H};X_{2},x_{2})\left( \ref{X1,gx1,X1F11,1H}%
\right) \\
=-\gamma _{1}\left[ B(x_{1}x_{2}\otimes 1_{H};G,x_{2})-B(x_{1}x_{2}\otimes
1_{H};GX_{2},1_{H})\right] + \\
-2\left[ B(x_{1}x_{2}\otimes 1_{H};1_{A},gx_{2})+B(x_{1}x_{2}\otimes
1_{H};X_{2},g)\right] .
\end{gather*}

\begin{gather*}
-\gamma _{1}\left[ B(x_{1}x_{2}\otimes
1_{H};G,gx_{1}x_{2})-B(x_{1}x_{2}\otimes 1_{H};GX_{2},gx_{1})\right] \left( %
\ref{X1,gx1,X1F11,gx1}\right) \\
+2\left[ B(x_{1}x_{2}\otimes 1_{H};1_{A},gx_{2})+B(x_{1}x_{2}\otimes
1_{H};X_{2},g)\right] =0.
\end{gather*}%
\begin{gather*}
2\left[ B(x_{1}x_{2}\otimes 1_{H};1_{A},gx_{2})+B(x_{1}x_{2}\otimes
1_{H};X_{2},g)\right] +\left( \ref{X1,gx1,X1F31,g}\right) \\
+\gamma _{1}\left[ -B(x_{1}x_{2}\otimes
1_{H};G,gx_{1}x_{2})+B(x_{1}x_{2}\otimes 1_{H};GX_{2},gx_{1})\right] =0.
\end{gather*}%
\begin{gather*}
-2\beta _{1}B(x_{1}x_{2}\otimes 1_{H};X_{1},x_{1})+ \\
+\lambda \left[ -1+B(x_{1}x_{2}\otimes
1_{H};1_{A},x_{1}x_{2})-B(x_{1}x_{2}\otimes 1_{H};X_{1},x_{2})\right] \left( %
\ref{X1,gx2,X1F11,1H}\right) \\
+\gamma _{1}\left[ -B(x_{1}x_{2}\otimes 1_{H};G,x_{1})\right] =0.
\end{gather*}%
\begin{gather*}
\left( \ref{X1,gx2,X1F11,gx1}\right) \\
+2B(x_{1}x_{2}\otimes 1_{H};1_{A},gx_{1})+2B(x_{1}x_{2}\otimes
1_{H};X_{1},g)=0.
\end{gather*}%
\begin{equation*}
+\gamma _{1}\left[ B(x_{1}x_{2}\otimes
1_{H};G,gx_{1}x_{2})+B(x_{1}x_{2}\otimes 1_{H};GX_{1},gx_{2})\right]
=0\left( \ref{X1,gx2,X1F11,gx2}\right)
\end{equation*}%
\begin{equation*}
2\left[ B(x_{1}x_{2}\otimes 1_{H};1_{A},gx_{1})+B(x_{1}x_{2}\otimes
1_{H};X_{1},g)\right] =0\left( \ref{X1,gx2,X131,g}\right)
\end{equation*}%
\begin{equation*}
\gamma _{1}\left[ -B(x_{1}x_{2}\otimes
1_{H};G,gx_{1}x_{2})+B(x_{1}x_{2}\otimes 1_{H};GX_{1},gx_{2})\right]
=0\left( \ref{X1,gx2,X141,g}\right)
\end{equation*}%
\begin{equation*}
+\lambda \left[ B(x_{1}\otimes 1_{H};1_{A},1_{H})+B(gx_{1}x_{2}\otimes
1_{H};1_{A},x_{2})\right] -\gamma _{1}B(gx_{1}x_{2}\otimes
1_{H};G,1_{H})=0\left( \ref{X1,gx1x2,X1F11,1H}\right)
\end{equation*}%
\begin{equation*}
+\lambda B(x_{1}\otimes 1_{H};1_{A},x_{1}x_{2})-\gamma
_{1}B(gx_{1}x_{2}\otimes 1_{H};G,x_{1}x_{2})=0\left( \ref%
{X1,gx1x2,X1F11,x1x2}\right)
\end{equation*}%
\begin{gather*}
-2\beta _{1}B(x_{2}\otimes 1_{H};1_{A},gx_{1})+\lambda \left[ B(x_{1}\otimes
1_{H};1_{A},gx_{1})+B(gx_{1}x_{2}\otimes 1_{H};1_{A},gx_{1}x_{2})\right]
\left( \ref{X1,gx1x2,X1F11,gx1}\right) \\
+\gamma _{1}B(gx_{1}x_{2}\otimes 1_{H};G,gx_{1})+2B(gx_{1}x_{2}\otimes
1_{H};1_{A},g)=0.
\end{gather*}%
\begin{gather*}
2\beta _{1}\left[ -B(x_{2}\otimes 1_{H};1_{A},gx_{2})-B(gx_{1}x_{2}\otimes
1_{H};1_{A},gx_{1}x_{2})\right] \left( \ref{X1,gx1x2,X171,gx2}\right) \\
+\lambda B(x_{1}\otimes 1_{H};1_{A},gx_{2})+\gamma _{1}B(gx_{1}x_{2}\otimes
1_{H};G,gx_{2})=0.
\end{gather*}%
\begin{gather*}
2\beta _{1}\left[ B(x_{2}\otimes 1_{H};G,g)-B(gx_{1}x_{2}\otimes
1_{H};G,gx_{1})\right] \left( \ref{X1,gx1x2,XF21,g}\right) \\
+\lambda \left[ -B(x_{1}\otimes 1_{H};G,g)-B(gx_{1}x_{2}\otimes
1_{H};G,gx_{2})\right] =0
\end{gather*}

\begin{equation*}
2\beta _{1}B(x_{2}\otimes 1_{H};G,gx_{1}x_{2})-\lambda B(x_{1}\otimes
1_{H};G,gx_{1}x_{2})=0\left( \ref{X1,gx1x2,X1F21,gx1x2}\right)
\end{equation*}%
\begin{gather*}
+\lambda \left[
\begin{array}{c}
B(g\otimes 1_{H};1_{A},g)+B(x_{2}\otimes \ 1_{H};1_{A},gx_{2}) \\
+B(x_{1}\otimes 1_{H};1_{A},gx_{1})+B(gx_{1}x_{2}\otimes
1_{H};1_{A},gx_{1}x_{2})%
\end{array}%
\right] \\
-\gamma _{1}\left[ B(x_{2}\otimes 1_{H};G,g)-B(gx_{1}x_{2}\otimes
1_{H};G,gx_{1})\right] =0\left( \ref{X1,gx1x2,X1F31,g}\right)
\end{gather*}%
\begin{gather*}
2B(gx_{1}x_{2}\otimes 1_{H};1_{A},x_{1})+\lambda \left[ B(g\otimes
1_{H};1_{A},x_{1})+B(x_{2}\otimes \ 1_{H};1_{A},x_{1}x_{2})\right] \\
+\gamma _{1}B(x_{2}\otimes 1_{H};G,x_{1})=0\left( \ref{X1,gx1x2,X1F31,x1}%
\right)
\end{gather*}

\begin{gather*}
2B(gx_{1}x_{2}\otimes 1_{H};1_{A},x_{2})+\lambda \left[ B(g\otimes
1_{H};X_{2},1_{H})-B(x_{1}\otimes 1_{H};1_{A},x_{1}x_{2})\right] +\left( \ref%
{X1,gx1x2,X1F31,x2}\right) \\
+\gamma _{1}\left[ B(x_{2}\otimes 1_{H};G,x_{2})+B(gx_{1}x_{2}\otimes
1_{H};G,x_{1}x_{2})\right] =0
\end{gather*}

\begin{gather*}
+\lambda B\left( g\otimes 1_{H};1_{A},gx_{1}x_{2}\right) +\left( \ref%
{X1,gx1x2,X1F31,gx1x2}\right) \\
-\gamma _{1}B(x_{2}\otimes 1_{H};G,gx_{1}x_{2})-\left[ +2B(x_{2}\otimes
1_{H};1_{A},gx_{2})+2B(gx_{1}x_{2}\otimes 1_{H};1_{A},gx_{1}x_{2})\right] =0
\end{gather*}

\begin{eqnarray*}
&&2\beta _{1}\left[
\begin{array}{c}
B(g\otimes 1_{H};1_{A},g)+B(x_{2}\otimes \ 1_{H};1_{A},gx_{2}) \\
+B(x_{1}\otimes 1_{H};1_{A},gx_{1})+B(gx_{1}x_{2}\otimes
1_{H};1_{A},gx_{1}x_{2})%
\end{array}%
\right] \left( \ref{X1,gx1x2,X1F41,g}\right) \\
&&-\gamma _{1}\left[ B(x_{1}\otimes 1_{H};G,g)+B(gx_{1}x_{2}\otimes
1_{H};G,gx_{2})\right] =0
\end{eqnarray*}%
\begin{eqnarray*}
&&2\beta _{1}B\left( g\otimes 1_{H};1_{A},gx_{1}x_{2}\right) \left( \ref%
{X1,gx1x2,X1F41,gx1x2}\right) \\
&&-\gamma _{1}B(x_{1}\otimes 1_{H};G,gx_{1}x_{2})-2B(x_{1}\otimes
1_{H};1_{A},gx_{2})=0
\end{eqnarray*}%
\begin{gather*}
2\left[ B(x_{1}\otimes 1_{H};1_{A},1_{H})+B(gx_{1}x_{2}\otimes
1_{H};1_{A},x_{2})\right] \left( \ref{X1,gx1x2,X1F51,1H}\right) \\
+\gamma _{1}\left[ B(g\otimes 1_{H};G,1_{H})+B(x_{2}\otimes \
1_{H};G,x_{2})+B(x_{1}\otimes 1_{H};G,x_{1})+B(gx_{1}x_{2}\otimes
1_{H};G,x_{1}x_{2})\right] =0
\end{gather*}%
\begin{gather*}
+\gamma _{1}\left[ B(g\otimes 1_{H};G,gx_{1})+B(x_{2}\otimes \
1_{H};G,gx_{1}x_{2})\right] +\left( \ref{X1,gx1x2,X1F51,gx1}\right) \\
-\left[
\begin{array}{c}
2B(g\otimes 1_{H};1_{A},g)+2B(x_{2}\otimes \ 1_{H};1_{A},gx_{2}) \\
+2B(x_{1}\otimes 1_{H};1_{A},gx_{1})+2B(gx_{1}x_{2}\otimes
1_{H};1_{A},gx_{1}x_{2})%
\end{array}%
\right] =0
\end{gather*}%
\begin{gather*}
2B(gx_{1}x_{2}\otimes 1_{H};G,1_{H})\left( \ref{X1,gx1x2,X1F61,1H}\right) \\
+\lambda \left[ B(g\otimes 1_{H};G,1_{H})+B(x_{2}\otimes \
1_{H};G,x_{2})+B(x_{1}\otimes 1_{H};G,x_{1})+B(gx_{1}x_{2}\otimes
1_{H};G,x_{1}x_{2})\right] =0
\end{gather*}%
\begin{equation*}
2B(gx_{1}x_{2}\otimes 1_{H};G,x_{1}x_{2})+\lambda B(g\otimes
1_{H};G,x_{1}x_{2})=0\left( \ref{X1,gx1x2,X1F61,x1x2}\right)
\end{equation*}%
\begin{gather*}
\lambda \left[ B(g\otimes 1_{H};G,gx_{1})+B(x_{2}\otimes \
1_{H};G,gx_{1}x_{2})\right] \\
+\left[ 2B(x_{2}\otimes 1_{H};G,g)-2B(gx_{1}x_{2}\otimes 1_{H};G,gx_{1})%
\right] =0\left( \ref{X1,gx1x2,X1F61,gx1}\right)
\end{gather*}%
\begin{gather*}
\beta _{1}\left[ B(g\otimes 1_{H};G,gx_{1})+B(x_{2}\otimes \
1_{H};G,gx_{1}x_{2})\right] \left( \ref{X1,gx1x2,X1F71,gx1}\right) \\
+B(x_{1}\otimes 1_{H};G,g)+B(gx_{1}x_{2}\otimes 1_{H};G,gx_{2})=0.
\end{gather*}

\subsubsection{$X_{2}$}

\begin{gather*}
\gamma _{2}B(g\otimes 1_{H};G,1_{H})+\lambda B(g\otimes
1_{H};1_{A},x_{1})\left( \ref{X2,g,X2F11,1H}\right) \\
+2B(x_{2}\otimes 1_{H};1_{A},1_{H})=0
\end{gather*}%
\begin{gather*}
2\beta _{2}B(g\otimes 1_{H};1_{A},gx_{1}x_{2})+\gamma _{2}B(g\otimes
1_{H};G,gx_{1})\left( \ref{X2,g,X2F11,gx1}\right) \\
+2B(x_{2}\otimes 1_{H};1_{A},gx_{1})=0
\end{gather*}

\begin{gather*}
\gamma _{2}B(g\otimes 1_{H};G,gx_{2})+\lambda B(g\otimes
1_{H};1_{A},gx_{1}x_{2})\left( \ref{X2,g,X2F11,gx2}\right) \\
+2B(x_{2}\otimes 1_{H};1_{A},gx_{2})+2B(g\otimes 1_{H};1_{A},g)=0.
\end{gather*}

\begin{equation*}
2\beta _{2}B(g\otimes 1_{H};G,gx_{2})+\lambda B(g\otimes
1_{H};G,gx_{1})+2B(x_{2}\otimes 1_{H};G,g)=0\left( \ref{X2,g,X2F21,g}\right)
\end{equation*}%
\begin{equation*}
B(x_{2}\otimes 1_{H};G,x_{1})=0\left( \ref{X2,g,X2F21,x1}\right)
\end{equation*}%
\begin{equation*}
\lambda B(g\otimes 1_{H};G,x_{1}x_{2})-2B(x_{2}\otimes 1_{H};G,x_{2})=0\ref%
{X2,g,X2F21,x2}
\end{equation*}%
\begin{equation*}
2\beta _{2}B(g\otimes 1_{H};1_{H},gx_{1}x_{2})+\gamma _{2}B(g\otimes
1_{H};G,gx_{1})+2B(x_{2}\otimes 1_{H};1_{H},gx_{1})=0.\left( \ref%
{X2,g,X2F31,g}\right)
\end{equation*}%
\begin{equation*}
\gamma _{2}B(g\otimes 1_{H};G,x_{1}x_{2})+2B(x_{2}\otimes
1_{H};1,x_{1}x_{2})=0\left( \ref{X2,g,X2F31,x2}\right)
\end{equation*}%
\begin{eqnarray*}
3B(g\otimes 1_{H};1_{A},x_{1})-\gamma _{2}B(g\otimes
1_{H};G,x_{1}x_{2})\left( \ref{X2,g,X2F41,x1}\right) && \\
+B(x_{2}\otimes 1_{H};X_{2},x_{1})+B\left( x_{2}\otimes
1_{H};1_{A},x_{1}x_{2}\right) &=&0
\end{eqnarray*}%
\begin{equation*}
B\left( g\otimes 1_{H};1_{A},gx_{1}x_{2}\right) +B(g\otimes
1_{H},G,gx_{1})=0\left( \ref{X2,g,X2F61,gx2}\right)
\end{equation*}%
\begin{gather*}
2B(g\otimes 1_{H};G,1_{H})-\lambda B(g\otimes 1_{H};G,x_{1}x_{2})\left( \ref%
{X2,g,X2F71,1H}\right) \\
+2B(g\otimes 1_{H};G,1_{H})+2B(x_{2}\otimes 1_{H};G,x_{2})=0
\end{gather*}%
\begin{equation*}
B(g\otimes 1_{H};G,gx_{1})+B(x_{2}\otimes 1_{H};G,gx_{1}x_{2})=0\left( \ref%
{X2,g,X2F71,gx1}\right)
\end{equation*}

\begin{eqnarray*}
&&-2\beta _{2}B(x_{1}\otimes 1_{H};1_{A},gx_{2})+\gamma _{2}B(x_{1}\otimes
1_{H};G,g)\left( \ref{X2,x1,X2F11,g}\right) \\
&&-\lambda \left[ B(g\otimes 1_{H};1_{A},g)+B(x_{1}\otimes
1_{H};1_{A},gx_{1})\right] -2B(gx_{1}x_{2}\otimes 1_{H};1_{A},g)=0
\end{eqnarray*}

\begin{equation*}
\gamma _{2}B(x_{1}\otimes 1_{H};G,x_{1})-\lambda B(g\otimes
1_{H};1_{A},x_{1})-2B(gx_{1}x_{2}\otimes 1_{H};1_{A},x_{1})=0\left( \ref%
{X2,x1,X2F11,x1}\right)
\end{equation*}%
\begin{eqnarray*}
&&+\lambda \left[ B(g\otimes 1_{H};1_{A},x_{2})-B(x_{1}\otimes
1_{H};1_{A},x_{1}x_{2})\right] \left( \ref{X2,x1,X2F11,x2}\right) \\
&&+2B(gx_{1}x_{2}\otimes 1_{H};1_{A},x_{2})=0.
\end{eqnarray*}%
\begin{eqnarray*}
&&+\gamma _{2}B(x_{1}\otimes 1_{H};G,gx_{1}x_{2})-\lambda B(g\otimes
1_{H};1_{A},gx_{1}x_{2})\left( \ref{X2,x1,X2F11,gx1x2}\right) \\
&&-2B(gx_{1}x_{2}\otimes 1_{H};1_{A},gx_{1}x_{2})-2B\left( x_{1}\otimes
1_{H};1_{A},gx_{1}\right) =0
\end{eqnarray*}%
\begin{equation*}
\lambda \left[ B(g\otimes 1_{H};G,1_{H})+B(x_{1}\otimes 1_{H};G,x_{1})\right]
+2B(gx_{1}x_{2}\otimes 1_{H};G,1_{H})=0\left( \ref{X2,x1,X2F21,1H}\right)
\end{equation*}

\begin{equation*}
2\beta _{2}B(x_{1}\otimes 1_{H};G,gx_{1}x_{2})+\lambda B\left( g\otimes
1_{H};G,gx_{1}\right) +2B(x_{2}\otimes 1_{H};G,g)=0.\left( \ref%
{X2,x1,X2F21,gx1}\right)
\end{equation*}

\begin{eqnarray*}
&&\lambda \left[ B\left( g\otimes 1_{H};G,gx_{2}\right) -B(x_{1}\otimes
1_{H};G,gx_{1}x_{2})\right] \left( \ref{X2,x1,X2F21,gx2}\right) \\
&&+2-B(x_{1}\otimes 1_{H};G,g_{1})+2B\left( x_{1}\otimes 1_{H};G,g\right) =0.
\end{eqnarray*}

\begin{gather*}
\gamma _{2}\left[ B(g\otimes 1_{H};G,1_{H})+B(x_{1}\otimes 1_{H};G,x_{1})%
\right] \left( \ref{X2,x1,X2F31,1H}\right) \\
+2B(x_{2}\otimes 1_{H};1_{A},1_{H})-2B(gx_{1}x_{2}\otimes
1_{H};1_{A},x_{1})=0.
\end{gather*}%
\begin{gather*}
\gamma _{2}\left[ -B(g\otimes 1_{H};G,gx_{2})+B(x_{1}\otimes
1_{H};G,gx_{1}x_{2})\right] +\left( \ref{X2,x1,X2F31,gx2}\right) \\
-2B(x_{2}\otimes 1_{H};1_{A},gx_{2})-2B(gx_{1}x_{2}\otimes
1_{H};1_{A},gx_{1}x_{2}) \\
-2B(g\otimes 1_{H};1_{A},g)-2B(x_{1}\otimes 1_{H};1_{A},gx_{1})=0.
\end{gather*}%
\begin{gather*}
\left( \ref{X2,x1,X2F41,1H}\right) \\
+\lambda \left[ -B(g\otimes 1_{H};1_{A},x_{2})+B(x_{1}\otimes
1_{H};1_{A},x_{1}x_{2})\right] \\
-2B(gx_{1}x_{2}\otimes 1_{H};1_{A},x_{2})=0
\end{gather*}

\begin{eqnarray*}
&&-\gamma _{2}B(x_{1}\otimes 1_{H};G,gx_{1}x_{2})+\lambda B\left( g\otimes
1_{H};1_{A},gx_{1}x_{2}\right) \left( \ref{X2,x1,X2F41,gx1}\right) \\
&&+2B(x_{1}\otimes 1_{H};1_{A},gx_{1})+2B(gx_{1}x_{2}\otimes
1_{H};1_{A},gx_{1}x_{2})=0.
\end{eqnarray*}

\begin{eqnarray*}
&&\beta _{2}\left[ -B(g\otimes 1_{H};G,gx_{2})+B(x_{1}\otimes
1_{H};G,gx_{1}x_{2})\right] \left( \ref{X2,x1,X2F61,g}\right) \\
&&-B(x_{2}\otimes 1_{H};G,g)+B(gx_{1}x_{2}\otimes 1_{H};G,gx_{1})=0
\end{eqnarray*}%
\begin{equation*}
B(x_{2}\otimes 1_{H};G,x_{2})+B(gx_{1}x_{2}\otimes
1_{H};G,x_{1}x_{2})=0\left( \ref{X2,x1,X2F61,x2}\right)
\end{equation*}

\begin{equation*}
B(x_{1}\otimes 1_{H};G,g)+B(gx_{1}x_{2}\otimes 1_{H};G,gx_{2})=0\left( \ref%
{X2,x1,X2F71,g}\right)
\end{equation*}%
\begin{equation*}
B(x_{2}\otimes \ 1_{H};G,x_{2})+B(gx_{1}x_{2}\otimes
1_{H};G,x_{1}x_{2})=0\left( \ref{X2,x1,X2F81,1H}\right)
\end{equation*}%
\begin{eqnarray*}
2 &&\beta _{2}\left[ B(g\otimes 1_{H};1_{A},g)+B(x_{2}\otimes \
1_{H};1_{A},gx_{2})\right] \left( \ref{X2,x2,X2F11,g}\right) \\
&&-\gamma _{2}B(x_{2}\otimes 1_{H};G,g)+\lambda B(x_{2}\otimes
1_{H};1_{A},gx_{1})=0.
\end{eqnarray*}%
\begin{equation*}
\gamma _{2}B(x_{2}\otimes 1_{H};G,x_{2})+\lambda B(x_{2}\otimes
1_{H};1_{A},x_{1}x_{2})=0.\left( \ref{X2,x2,X2F11,x2}\right)
\end{equation*}

\begin{eqnarray*}
+ &&2\beta _{2}B\left( g\otimes 1_{H};1_{A},gx_{1}x_{2}\right) -\gamma
_{2}B(x_{2}\otimes 1_{H};G,gx_{1}x_{2})\left( \ref{X2,x2,X2F11,gx1x2}\right)
\\
&&+2B\left( x_{2}\otimes 1_{H};1_{A},gx_{1}\right) =0
\end{eqnarray*}

\begin{equation*}
\beta _{2}\left[ B(g\otimes 1_{H};G,gx_{1})+B(x_{2}\otimes \
1_{H};G,gx_{1}x_{2})\right] =0\left( \ref{X2,x2,X2F21,gx1}\right)
\end{equation*}%
\begin{eqnarray*}
&&2\beta _{2}B(g\otimes 1_{H};G,gx_{2})-\lambda B(x_{2}\otimes
1_{H};G,gx_{1}x_{2})\left( \ref{X2,x2,X2F21,gx2}\right) \\
&&+2B\left( x_{2}\otimes 1_{H};G,g\right) =0
\end{eqnarray*}

\begin{eqnarray*}
&&2B(x_{2}\otimes 1_{H};1_{A},1_{H})+\gamma _{2}\left[ B(g\otimes
1_{H};G,1_{H})+B(x_{2}\otimes \ 1_{H};G,x_{2}\right] \left( \ref%
{X2,x2,X2F41,1H}\right) \\
&&+\lambda \left[ B(g\otimes 1_{H};1_{A},x_{1})+B(x_{2}\otimes \
1_{H};1_{A},x_{1}x_{2})\right] =0.
\end{eqnarray*}

\begin{eqnarray*}
&&\gamma _{2}B(g\otimes 1_{H};G,gx_{2})-\lambda B\left( g\otimes
1_{H};1_{A},gx_{1}x_{2}\right) \left( \ref{X2,x2,X2F41,gx2}\right) \\
&&+2\left[ +B(g\otimes 1_{H};1_{A},g)+B(x_{2}\otimes \ 1_{H};1_{A},gx_{2})%
\right] =0.
\end{eqnarray*}

\begin{equation*}
\lambda \left[ B(g\otimes 1_{H};G,gx_{1})+B(x_{2}\otimes \
1_{H};G,gx_{1}x_{2})\right] =0\left( \ref{X2,x2,X2F71,g}\right)
\end{equation*}%
\begin{equation*}
\gamma _{2}B(x_{1}x_{2}\otimes 1_{H};G,g)+\lambda B(x_{1}x_{2}\otimes
1_{H};X_{1},g)=0.\left( \ref{X2,x1x2,X2F11,g}\right)
\end{equation*}

\begin{gather*}
2\beta _{2}B(x_{1}x_{2}\otimes 1_{H};X_{2},x_{1})+\gamma
_{2}B(x_{1}x_{2}\otimes 1_{H};G,x_{1})\left( \ref{X2,x1x2,X2F11,x1}\right) \\
+\lambda B(x_{1}x_{2}\otimes 1_{H};X_{1},x_{1})=0.
\end{gather*}

\begin{gather*}
2\beta _{2}B(x_{1}x_{2}\otimes 1_{H};X_{2},x_{2})+\gamma
_{2}B(x_{1}x_{2}\otimes 1_{H};G,x_{2})\left( \ref{X2,x1x2,X2F11,x2}\right) \\
+\lambda B(x_{1}x_{2}\otimes 1_{H};X_{1},x_{2})=0.
\end{gather*}%
\begin{eqnarray*}
&&\gamma _{2}B(x_{1}x_{2}\otimes 1_{H};G,gx_{1}x_{2})+\left( \ref%
{X2,x1x2,X2F11,gx1x2}\right) \\
&&+2B(x_{1}x_{2}\otimes 1_{H};1_{A},gx_{1})=0.
\end{eqnarray*}

\begin{equation*}
2\beta _{2}B(x_{1}x_{2}\otimes 1_{H};GX_{2},1_{H})+=0.\left( \ref%
{X2,x1x2,X2F21,1H}\right)
\end{equation*}

\begin{equation*}
2B(x_{1}x_{2}\otimes 1_{H};G,g)-\lambda B(x_{1}x_{2}\otimes
1_{H};GX_{1},gx_{2})=0.\left( \ref{X2,x1x2,X2F21,gx2}\right)
\end{equation*}%
\begin{eqnarray*}
&&2\beta _{2}\left[
\begin{array}{c}
+1-B(x_{1}x_{2}\otimes 1_{H};1_{A},x_{1}x_{2}) \\
-B(x_{1}x_{2}\otimes 1_{H};X_{2},x_{1})+B(x_{1}x_{2}\otimes
1_{H};X_{1},x_{2})%
\end{array}%
\right] \left( \ref{X2,x1x2,X2F31,1H}\right) \\
&=&0.
\end{eqnarray*}%
\begin{equation*}
\gamma _{2}B(x_{1}x_{2}\otimes 1_{H};GX_{1},gx_{2})+2B(x_{1}x_{2}\otimes
1_{H};X_{1},g)=0.\left( \ref{X2,x1x2,X2F31,gx2}\right)
\end{equation*}

\begin{eqnarray*}
&&-\gamma _{2}B(x_{1}x_{2}\otimes 1_{H};GX_{2},1_{H})+ \\
-\lambda \left[
\begin{array}{c}
+1-B(x_{1}x_{2}\otimes 1_{H};1_{A},x_{1}x_{2}) \\
-B(x_{1}x_{2}\otimes 1_{H};X_{2},x_{1})+B(x_{1}x_{2}\otimes
1_{H};X_{1},x_{2})%
\end{array}%
\right] &=&0.\left( \ref{X2,x1x2,X2F41,1H}\right)
\end{eqnarray*}%
\begin{gather*}
2B(x_{1}x_{2}\otimes 1_{H};1_{A},gx_{2})+\left( \ref{X2,x1x2,X2F41,gx2}%
\right) \\
+2B(x_{1}x_{2}\otimes 1_{H};X_{2},g)=0.
\end{gather*}

\begin{gather*}
2B(x_{1}x_{2}\otimes 1_{H};X_{1},g)-\left( \ref{X2,x1x2,X2F51,g}\right) \\
\gamma _{2}\left[
\begin{array}{c}
-B(x_{1}x_{2}\otimes 1_{H};G,gx_{1}x_{2})+ \\
B(x_{1}x_{2}\otimes 1_{H};GX_{2},gx_{1})-B(x_{1}x_{2}\otimes
1_{H};GX_{1},gx_{2})%
\end{array}%
\right] =0.
\end{gather*}%
\begin{gather*}
2B(x_{1}x_{2}\otimes 1_{H};G,g)\left( \ref{X2,x1x2,X2F71,g}\right) \\
+\lambda \left[
\begin{array}{c}
-B(x_{1}x_{2}\otimes 1_{H};G,gx_{1}x_{2})+ \\
B(x_{1}x_{2}\otimes 1_{H};GX_{2},gx_{1})-B(x_{1}x_{2}\otimes
1_{H};GX_{1},gx_{2})%
\end{array}%
\right] =0.
\end{gather*}%
\begin{equation*}
B(x_{1}x_{2}\otimes 1_{H};G,gx_{1}x_{2})-B(x_{1}x_{2}\otimes
1_{H};GX_{2},gx_{1})=0.\left( \ref{X2,x1x2,X2F71,gx1x2}\right)
\end{equation*}%
\begin{gather*}
2\beta _{2}B(x_{1}x_{2}\otimes 1_{H};X_{2},x_{2})+\gamma _{2}\left[
B(x_{1}x_{2}\otimes 1_{H};G,x_{2})-B(x_{1}x_{2}\otimes 1_{H};GX_{2},1_{H})%
\right] \left( \ref{X2,gx1,X2F11,1H}\right) \\
+\lambda \left[ -1+B(x_{1}x_{2}\otimes
1_{H};1_{A},x_{1}x_{2})+B(x_{1}x_{2}\otimes 1_{H};X_{2},x_{1})\right] =0.
\end{gather*}

\begin{equation*}
B(x_{1}x_{2}\otimes 1_{H};1_{A},gx_{2})+B(x_{1}x_{2}\otimes
1_{H};X_{2},g)=0.\left( \ref{X2,gx1,X2F11,gx2}\right)
\end{equation*}%
\begin{equation*}
\lambda \left[ B(x_{1}x_{2}\otimes 1_{H};G,gx_{1}x_{2})-B(x_{1}x_{2}\otimes
1_{H};GX_{2},gx_{1})\right] =0.\left( \ref{X2,gx1,X2F21,g}\right)
\end{equation*}%
\begin{gather*}
2\beta _{2}\left[ -1+B(x_{1}x_{2}\otimes
1_{H};1_{A},x_{1}x_{2})-B(x_{1}x_{2}\otimes 1_{H};X_{1},x_{2})\right] \left( %
\ref{X2,gx2,X2F11,1H}\right) \\
+\gamma _{2}\left[ -B(x_{1}x_{2}\otimes 1_{H};G,x_{1})\right] + \\
-\lambda B(x_{1}x_{2}\otimes 1_{H};X_{1},x_{1})=0.
\end{gather*}

\begin{eqnarray*}
&&\gamma _{2}\left[ B(x_{1}x_{2}\otimes
1_{H};G,gx_{1}x_{2})+B(x_{1}x_{2}\otimes 1_{H};GX_{1},gx_{2})\right] \left( %
\ref{X2,gx2,X2F11,gx2}\right) \\
&&+ \\
&&+\left[ +2B(x_{1}x_{2}\otimes 1_{H};1_{A},gx_{1})+2B(x_{1}x_{2}\otimes
1_{H};X_{1},g)\right] =0.
\end{eqnarray*}%
\begin{eqnarray*}
&&2\left[ -B(x_{1}x_{2}\otimes 1_{H};1_{A},gx_{1})-B(x_{1}x_{2}\otimes
1_{H};X_{1},g)\right] \left( \ref{X2,gx2,X2F41,g}\right) \\
&&+\gamma _{2}\left[ -B(x_{1}x_{2}\otimes
1_{H};G,gx_{1}x_{2})-B(x_{1}x_{2}\otimes 1_{H};GX_{1},gx_{2})\right] \\
&&=0.
\end{eqnarray*}%
\begin{equation*}
-\gamma _{2}B(gx_{1}x_{2}\otimes 1_{H};G,1_{H})+\lambda \left[
+B(x_{2}\otimes 1_{H};1_{A},1_{H})-B(gx_{1}x_{2}\otimes 1_{H};1_{A},x_{1})%
\right] =0\left( \ref{X2,gx1x2,X2F11,1H}\right)
\end{equation*}%
\begin{equation*}
\gamma _{2}B(gx_{1}x_{2}\otimes 1_{H};G,x_{1}x_{2})-\lambda B(x_{2}\otimes
1_{H};1_{A},x_{1}x_{2})=0\left( \ref{X2,gx1x2,X2F11,x1x2}\right)
\end{equation*}%
\begin{gather*}
2\beta _{2}\left[ B(x_{1}\otimes 1_{H};1_{A},gx_{1})+B(gx_{1}x_{2}\otimes
1_{H};1_{A},gx_{1}x_{2})\right] \left( \ref{X2,gx1x2,X2F11,gx1}\right) \\
+\gamma _{2}B(gx_{1}x_{2}\otimes 1_{H};G,gx_{1})-\lambda B(x_{2}\otimes
1_{H};1_{A},gx_{1})=0.
\end{gather*}

\begin{gather*}
2\beta _{2}B(x_{1}\otimes 1_{H};1_{A},gx_{2})+\gamma
_{2}B(gx_{1}x_{2}\otimes 1_{H};G,gx_{2})\left( \ref{X2,gx1x2,X2F11,gx2}%
\right) \\
+\lambda \left[ -B(x_{2}\otimes 1_{H};1_{A},gx_{2})-B(gx_{1}x_{2}\otimes
1_{H};1_{A},gx_{1}x_{2})\right] \\
+2\ =0.
\end{gather*}

\begin{gather*}
2\beta _{2}\left[ -B(x_{1}\otimes 1_{H};G,g)-B(gx_{1}x_{2}\otimes
1_{H};G,gx_{2})\right] \left( \ref{X2,gx1x2,X2F21,g}\right) \\
+\lambda \left[ B(x_{2}\otimes 1_{H};G,g)-B(gx_{1}x_{2}\otimes
1_{H};G,gx_{1})\right] =0.
\end{gather*}%
\begin{gather*}
-2\beta _{2}B(x_{1}\otimes 1_{H};G,gx_{1}x_{2})\left( \ref%
{X2,gx1x2,X2F21,gx1x2}\right) \\
+\lambda B(x_{2}\otimes 1_{H};G,gx_{1}x_{2})-2B(gx_{1}x_{2}\otimes
1_{H};G,gx_{1})=0.
\end{gather*}%
\begin{gather*}
2\beta _{2}\left[
\begin{array}{c}
B(g\otimes 1_{H};1_{A},g)+B(x_{2}\otimes \ 1_{H};1_{A},gx_{2}) \\
+B(x_{1}\otimes 1_{H};1_{A},gx_{1})+B(gx_{1}x_{2}\otimes
1_{H};1_{A},gx_{1}x_{2})%
\end{array}%
\right] \left( \ref{X2,gx1x2,X2F31,g}\right) \\
-\gamma _{2}\left[ B(x_{2}\otimes 1_{H};G,g)-B(gx_{1}x_{2}\otimes
1_{H};G,gx_{1})\right] =0.
\end{gather*}

\begin{eqnarray*}
&&\gamma _{2}\left[ B(x_{1}\otimes 1_{H};G,g)+B(gx_{1}x_{2}\otimes
1_{H};G,gx_{2})\right] \left( \ref{X2,gx1x2,X2F41,g}\right) \\
&&-\lambda \left[
\begin{array}{c}
B(g\otimes 1_{H};1_{A},g)+B(x_{2}\otimes \ 1_{H};1_{A},gx_{2}) \\
+B(x_{1}\otimes 1_{H};1_{A},gx_{1})+B(gx_{1}x_{2}\otimes
1_{H};1_{A},gx_{1}x_{2})%
\end{array}%
\right] =0.
\end{eqnarray*}

\begin{eqnarray*}
&&2B(gx_{1}x_{2}\otimes 1_{H};1_{A},x_{1})+\gamma _{2}\left[
\begin{array}{c}
-B(x_{1}\otimes 1_{H};G,x_{1})+ \\
-B(gx_{1}x_{2}\otimes 1_{H};G,x_{1}x_{2})%
\end{array}%
\right] \left( \ref{X2,gx1x2,X2F41,x1}\right) \\
&&+\lambda \left[ B(g\otimes 1_{H};1_{A},x_{1})+B(x_{2}\otimes \
1_{H};1_{A},x_{1}x_{2})\right] =0.
\end{eqnarray*}

\begin{gather*}
2\left[ -B(x_{2}\otimes 1_{H};1_{A},1_{H})+B(gx_{1}x_{2}\otimes
1_{H};1_{A},x_{1})\right] +\left( \ref{X2,gx1x2,X2F51,1H}\right) \\
-\gamma _{2}\left[
\begin{array}{c}
B(g\otimes 1_{H};G,1_{H})+B(x_{2}\otimes \ 1_{H};G,x_{2}) \\
+B(x_{1}\otimes 1_{H};G,x_{1})+B(gx_{1}x_{2}\otimes 1_{H};G,x_{1}x_{2})%
\end{array}%
\right] =0.
\end{gather*}

\begin{equation*}
B(x_{2}\otimes 1_{H};G,g)-B(gx_{1}x_{2}\otimes 1_{H};G,gx_{1})=0.\left( \ref%
{X2,gx1x2,X2F61,gx2}\right)
\end{equation*}%
\begin{gather*}
2B(gx_{1}x_{2}\otimes 1_{H};G,1_{H})\left( \ref{X2,gx1x2,X2F71,1H}\right) \\
+\lambda \left[
\begin{array}{c}
B(g\otimes 1_{H};G,1_{H})+B(x_{2}\otimes \ 1_{H};G,x_{2}) \\
+B(x_{1}\otimes 1_{H};G,x_{1})+B(gx_{1}x_{2}\otimes 1_{H};G,x_{1}x_{2})%
\end{array}%
\right] =0.
\end{gather*}

\subsection{LIST\ OF\ MONOMIAL\ EQUALITIES 1\label{LME1}}

We list the new monomial equalities we obtained so far and take out
constants.%
\begin{equation*}
B(g\otimes 1_{H};1_{A},x_{1})=0.\left( \ref{G,g, GF2,x1}\right)
\end{equation*}%
\begin{equation*}
B(g\otimes 1_{H};1_{A},x_{2})=0\left( \ref{G,g, GF2,x2}\right)
\end{equation*}%
\begin{equation*}
B(x_{1}\otimes 1_{H};1_{A},x_{1}x_{2})=0.\left( \ref{X1,g, X1F41,x1}\right)
\end{equation*}%
\begin{equation*}
B(x_{1}\otimes 1_{H};1_{A},x_{1}x_{2})=0.\left( \ref{X1,x2,X1F41,x1x2}\right)
\end{equation*}%
\begin{equation*}
B(x_{1}\otimes 1_{H};1_{A},x_{1}x_{2})=0.\left( \ref{X1,x2,X1F51,x2}\right)
\end{equation*}%
\begin{equation*}
B\left( x_{1}\otimes 1_{H};1_{A},x_{1}x_{2}\right) =0\left( \ref%
{X1,x1,X1F31,x1x2}\right)
\end{equation*}%
\begin{equation*}
B(x_{1}\otimes 1_{H};1_{A},x_{1}x_{2})=0.\left( \ref{X1,x2,X1F41,x1x2}\right)
\end{equation*}%
\begin{equation*}
B(x_{1}\otimes 1_{H};1_{A},x_{1}x_{2})=0.\left( \ref{X1,x2,X1F51,x2}\right)
\end{equation*}%
\begin{equation*}
B(x_{1}\otimes 1_{H};G,x_{1})=0.\left( \ref{X1,g,X1F21,x1}\right)
\end{equation*}%
\begin{equation*}
B(x_{2}\otimes 1_{H};1_{A},x_{1}x_{2})=0\left( \ref{G,x2, GF2,x1x2}\right)
\end{equation*}%
\begin{equation*}
B(x_{2}\otimes 1_{H};G,x_{1})=0\left( \ref{X2,g,X2F21,x1}\right)
\end{equation*}%
\begin{equation*}
B(x_{1}x_{2}\otimes 1_{H};GX_{1},x_{1}x_{2})=0\text{ }\left( \ref{G,x1x2,
GF3,gx1}\right) \text{ this is \ref{x1x2ot1,GX1,x1x2}}
\end{equation*}%
\begin{equation*}
B(x_{1}x_{2}\otimes 1_{H};GX_{2},1_{H})=0.\left( \ref{X1,x1x2,X1F21,1H}%
\right)
\end{equation*}%
\begin{equation*}
B(x_{1}x_{2}\otimes 1_{H};GX_{2},1_{H})=0\left( \ref{G,x1x2, GF2,1H}\right)
\end{equation*}%
\begin{equation*}
B(x_{1}x_{2}\otimes 1_{H};GX_{2},1_{H})=0.\left( \ref{X2,x1x2,X2F21,1H}%
\right)
\end{equation*}%
\begin{equation*}
B(gx_{1}x_{2}\otimes 1_{H};G,x_{1}x_{2})=0.\left( \ref{X1,x2,X1F21,x1x2}%
\right)
\end{equation*}%
We take one of the repetition as before and relabel them.

\begin{equation}
B(g\otimes 1_{H};1_{A},x_{1})=0.\left( \ref{G,g, GF2,x1}\right)
\label{got1,1,x1}
\end{equation}%
\begin{equation}
B(g\otimes 1_{H};1_{A},x_{2})=0\left( \ref{G,g, GF2,x2}\right)
\label{got1,1,x2}
\end{equation}%
\begin{equation}
B(x_{1}\otimes 1_{H};1_{A},x_{1}x_{2})=0.\left( \ref{X1,g, X1F41,x1}\right)
\label{x1ot1,1,x1x2}
\end{equation}%
\begin{equation}
B(x_{1}\otimes 1_{H};G,x_{1})=0.\left( \ref{X1,g,X1F21,x1}\right)
\label{x1ot1,G,x1}
\end{equation}%
\begin{equation}
B(x_{2}\otimes 1_{H};1_{A},x_{1}x_{2})=0\left( \ref{G,x2, GF2,x1x2}\right)
\label{x2ot1,1,x1x2}
\end{equation}%
\begin{equation}
B(x_{2}\otimes 1_{H};G,x_{1})=0\left( \ref{X2,g,X2F21,x1}\right)
\label{x2ot1,G,x1}
\end{equation}%
\begin{equation}
B(x_{1}x_{2}\otimes 1_{H};GX_{2},1_{H})=0.\left( \ref{X1,x1x2,X1F21,1H}%
\right)  \label{x1x2ot1,GX2,1}
\end{equation}%
\begin{equation*}
B(gx_{1}x_{2}\otimes 1_{H};G,x_{1}x_{2})=0.\text{we do not consider }\left( %
\ref{X2,x1x2,X2F21,x1x2}\right) \text{this as it is}\left( \ref%
{X1,x2,X1F21,x1x2}\right) \text{ }
\end{equation*}

Now we take the complete list of equalities $\left( \ref{LAE1}\right) $ and
cancel there the terms above.

\subsection{LIST\ OF\ ALL\ EQUALITIES 2\label{LAE2}}

\subsubsection{$G$}

\begin{equation*}
\gamma _{1}B\left( g\otimes 1_{H};1_{A},x_{1}\right) +\gamma _{2}B\left(
g\otimes 1_{H};1_{A},x_{2}\right) =0\text{ }\left( \ref{G,g,GF1,1H}\right)
\end{equation*}

\begin{equation*}
2\alpha B(g\otimes 1_{H};G,gx_{1})+\gamma _{2}B(g\otimes
1_{H};1_{A},gx_{1}x_{2})=0\left( \ref{G,g, GF1,gx1}\right)
\end{equation*}

\begin{equation*}
2\alpha B(g\otimes 1_{H};G,gx_{2})-\gamma _{1}B\left( g\otimes
1_{H};1_{A},gx_{1}x_{2}\right) =0.\left( \ref{G,g, GF1,gx2}\right)
\end{equation*}

\begin{equation*}
\gamma _{1}B\left( g\otimes 1_{H};G,gx_{1}\right) +\gamma _{2}B\left(
g\otimes 1_{H};G,gx_{2}\right) =0.\left( \ref{G,g, GF2,g}\right)
\end{equation*}%
\begin{equation*}
\begin{array}{c}
2\alpha B(x_{1}\otimes 1_{H};G,g)+ \\
+\gamma _{1}\left[ -B(g\otimes 1_{H};1_{A},g)-B(x_{1}\otimes
1_{H};1_{A},gx_{1})\right] -\gamma _{2}B(x_{1}\otimes 1_{H};1_{A},gx_{2})%
\end{array}%
=0.\left( \ref{G,x1, GF1,g}\right)
\end{equation*}

\begin{equation*}
-\gamma _{1}B(g\otimes 1_{H};1_{A},gx_{1}x_{2})+2\alpha B(x_{1}\otimes
1_{H};G,gx_{1}x_{2})=0.\left( \ref{G,x1, GF1,gx1x2}\right)
\end{equation*}%
\begin{equation*}
2B(x_{1}\otimes 1_{H};1_{A},1_{H})+\gamma _{1}B(g\otimes
1_{H};G,1_{H})=0.\left( \ref{G,x1, GF2,1H}\right)
\end{equation*}%
\begin{equation*}
2B(x_{1}\otimes 1_{H};1_{A},x_{1}x_{2})+\gamma _{1}B(g\otimes
1_{H};G,x_{1}x_{2})=0\text{ }\left( \ref{G,x1, GF2,x1x2}\right)
\end{equation*}%
\begin{equation*}
\gamma _{1}B\left( g\otimes 1_{H};G,gx_{1}\right) +\gamma _{2}B(x_{1}\otimes
1_{H};G,gx_{1}x_{2})=0.\left( \ref{G,x1, GF2,gx1}\right)
\end{equation*}

\begin{equation*}
\gamma _{1}[B\left( g\otimes 1_{H};G,gx_{2}\right) -B(x_{1}\otimes
1_{H};G,gx_{1}x_{2})]=0.\left( \ref{G,x1,GF2,gx2}\right)
\end{equation*}

\begin{equation*}
B(g\otimes 1_{H};X_{2},1_{H})=0\left( \ref{G,x1,GF3,1H}\right)
\end{equation*}%
\begin{equation*}
-\gamma _{2}B(g\otimes 1_{H};1_{A},gx_{1}x_{2})-2\alpha B(g\otimes
1_{H};G,gx_{1})=0.\left( \ref{G,x1,GF3,gx1}\right)
\end{equation*}

\begin{equation*}
\alpha \left[ -B(g\otimes 1_{H};g,gx_{2})+B(x_{1}\otimes 1_{H};g,gx_{1}x_{2})%
\right] =0.\left( \ref{G,x1, GF3,x2}\right)
\end{equation*}%
\begin{equation*}
0=0.\left( \ref{G,x1, GF4,1H}\right)
\end{equation*}

\begin{equation*}
-\gamma _{1}B(g\otimes 1_{H};1_{A},gx_{1}x_{2})+2\alpha B(x_{1}\otimes
1_{H};G,gx_{1}x_{2})=0.\left( \ref{G,x1, GF4,gx1}\right)
\end{equation*}

\begin{equation*}
\alpha \lbrack -B(g\otimes 1_{H};g,gx_{2})+B(x_{1}\otimes
1_{H};G,gx_{1}x_{2})]=0.\left( \ref{G,x1, GF5,g}\right)
\end{equation*}%
\begin{equation*}
B(g\otimes 1_{H};G,x_{1}x_{2})=0.\left( \ref{G,x1,GF6,x1}\right)
\end{equation*}%
\begin{equation*}
0=0.\left( \ref{G,x1,GF6,x2}\right)
\end{equation*}%
\begin{equation*}
\begin{array}{c}
2\alpha B(x_{2}\otimes 1_{H};G,g)+ \\
-\gamma _{1}B(x_{2}\otimes 1_{H};1_{A},gx_{1})+\gamma _{2}\left[ -B(g\otimes
1_{H};1_{A},g)-B(x_{2}\otimes \ 1_{H};1_{A},gx_{2})\right]%
\end{array}%
=0.\left( \ref{G,x2, GF1,g}\right)
\end{equation*}%
\begin{equation*}
B(x_{2}\otimes \ 1_{H};1_{A},x_{1}x_{2})=0\left( \ref{G,x2, GF1,x1}\right)
\end{equation*}

\begin{equation*}
2\alpha B(x_{2}\otimes 1_{H};G,gx_{1}x_{2})-\gamma _{2}B\left( g\otimes
1_{H};1_{A},gx_{1}x_{2}\right) =0.\left( \ref{G,x2, GF1,gx1x2}\right)
\end{equation*}%
\begin{equation*}
\begin{array}{c}
2B(x_{2}\otimes 1_{H};1_{A},1_{H}) \\
+\gamma _{2}\left[ B(g\otimes 1_{H};G,1_{H})+B(x_{2}\otimes \ 1_{H};G,x_{2}%
\right]%
\end{array}%
=0.\left( \ref{G,x2, GF2,1H}\right)
\end{equation*}%
\begin{equation*}
0=0\left( \ref{G,x2, GF2,x1x2}\right)
\end{equation*}%
\begin{equation*}
\gamma _{2}\left[ B(g\otimes 1_{H};G,gx_{1})+B(x_{2}\otimes \
1_{H};G,gx_{1}x_{2})\right] =0.\left( \ref{G,x2, GF2,gx1}\right)
\end{equation*}

\begin{equation*}
-\gamma _{1}B(x_{2}\otimes 1_{H};G,gx_{1}x_{2})+\gamma _{2}B(g\otimes
1_{H};G,gx_{2})=0.\left( \ref{G,x2, GF2,gx2}\right)
\end{equation*}%
\begin{equation*}
0=0\left( \ref{G,x2, FG3,1H}\right)
\end{equation*}

\begin{equation*}
-\gamma _{2}B(g\otimes 1_{H};1_{A},gx_{1}x_{2})-2\alpha B(x_{2}\otimes
1_{H};G,gx_{1}x_{2})=0.\left( \ref{G,x2, GF3,gx2}\right)
\end{equation*}%
\begin{equation*}
B(x_{2}\otimes \ 1_{H};1_{A},x_{1}x_{2})=0.\left( \ref{G,x2, GF4,1H}\right)
\end{equation*}%
\begin{equation*}
\alpha \left[ B(g\otimes 1_{H};G,gx_{1})+B(x_{2}\otimes \
1_{H};G,gx_{1}x_{2})\right] =0\left( \ref{G,x2, GF4,gx1}\right)
\end{equation*}

\begin{equation*}
\gamma _{1}B(g\otimes 1_{H};1_{A},gx_{1}x_{2})-2\alpha B(g\otimes
1_{H};G,gx_{2})=0.\left( \ref{G,x2, GF4,gx2}\right)
\end{equation*}%
\begin{equation*}
\alpha \left[ B(g\otimes 1_{H};G,gx_{1})+B(x_{2}\otimes \
1_{H};G,gx_{1}x_{2})\right] =0\left( \ref{G,x2, GF5,g}\right)
\end{equation*}

\begin{equation*}
\gamma _{2}[B(g\otimes 1_{H},G,gx_{1})+B(x_{2}\otimes
1_{H};G,gx_{1}x_{2})=0.=0.\left( \ref{G,x2, GF6,g}\right)
\end{equation*}%
\begin{equation*}
B(g\otimes 1_{H};G,x_{1}x_{2})=0.\left( \ref{G,x2; GF6,x2}\right)
\end{equation*}

\begin{equation*}
\gamma _{1}[B(g\otimes 1_{H};G,gx_{1})+B(x_{2}\otimes
1_{H};G,gx_{1}x_{2})]=0.\left( \ref{G,x2; GF7,g}\right)
\end{equation*}%
\begin{equation*}
B(x_{2}\otimes \ 1_{H};1_{A},x_{1}x_{2})=0.\left( \ref{G,x2; GF7,x1}\right)
\end{equation*}%
\begin{equation*}
\gamma _{1}B(x_{1}x_{2}\otimes 1_{H};X_{1},g)+\gamma _{2}B(x_{1}x_{2}\otimes
1_{H};X_{2},g)=0.\left( \ref{G,x1x2, GF1,g}\right)
\end{equation*}

\begin{equation*}
\begin{array}{c}
2\alpha B(x_{1}x_{2}\otimes 1_{H};G,x_{1})+ \\
+\gamma _{1}B(x_{1}x_{2}\otimes 1_{H};X_{1},x_{1})+\gamma
_{2}B(x_{1}x_{2}\otimes 1_{H};X_{2},x_{1})%
\end{array}%
=0.\left( \ref{G,x1x2, GF1,x1}\right)
\end{equation*}%
\begin{equation*}
\begin{array}{c}
2\alpha B(x_{1}x_{2}\otimes 1_{H};G,x_{2})+ \\
+\gamma _{1}B(x_{1}x_{2}\otimes 1_{H};X_{1},x_{2})+\gamma
_{2}B(x_{1}x_{2}\otimes 1_{H};X_{2},x_{2})%
\end{array}%
=0.\left( \ref{G,x1x2, GF1,x2}\right)
\end{equation*}

\begin{equation*}
0=0\left( \ref{G,x1x2, GF2,1H}\right)
\end{equation*}%
\begin{equation*}
\begin{array}{c}
2B(x_{1}x_{2}\otimes 1_{H};1_{A},gx_{1}) \\
+\gamma _{2}B(x_{1}x_{2}\otimes 1_{H};GX_{2},gx_{1})%
\end{array}%
=0.\left( \ref{G,x1x2, GF2,gx1}\right)
\end{equation*}%
\begin{equation*}
\begin{array}{c}
2B(x_{1}x_{2}\otimes 1_{H};1_{A},gx_{2}) \\
+\gamma _{1}B(x_{1}x_{2}\otimes 1_{H};GX_{1},gx_{2})%
\end{array}%
=0\left( \ref{G,x1x2, GF2,gx2}\right)
\end{equation*}

\begin{equation*}
0=0\text{ }\left( \ref{G,x1x2, GF3,gx1}\right)
\end{equation*}%
\begin{equation*}
-\gamma _{1}[1-B(x_{1}x_{2}\otimes
1_{H};1_{A},x_{1}x_{2})+B(x_{1}x_{2}\otimes
1_{H};x_{1},x_{2})-B(x_{1}x_{2}\otimes 1_{H};x_{2},x_{1})]-2\alpha
B(x_{1}x_{2}\otimes 1_{H};GX_{2},1_{H}=0.\left( \ref{G,x1x2, GF4,1H}\right)
\end{equation*}

\begin{equation*}
-2B(x_{1}x_{2}\otimes 1_{H};X_{1},g)+\gamma _{2}[-B(x_{1}x_{2}\otimes
1_{H};g,gx_{1}x_{2})-B(x_{1}x_{2}\otimes
1_{H};GX_{1},gx_{2})+B(x_{1}x_{2}\otimes 1_{H};GX_{2},gx_{1})]=0\left( \ref%
{G,x1x2, GF6,g}\right)
\end{equation*}%
\begin{gather*}
-2B(x_{1}x_{2}\otimes 1_{H};X_{2},g)+ \\
-\gamma _{1}[-B(x_{1}x_{2}\otimes 1_{H};g,gx_{1}x_{2})-B(x_{1}x_{2}\otimes
1_{H};GX_{1},gx_{2})+B(x_{1}x_{2}\otimes 1_{H};GX_{2},gx_{1})]=0\left( \ref%
{G,x1x2, GF7,g}\right)
\end{gather*}

\begin{equation*}
\begin{array}{c}
2\alpha B(x_{1}x_{2}\otimes 1_{H};G,x_{2})+ \\
+\gamma _{1}\left[ -1+B(x_{1}x_{2}\otimes
1_{H};1_{A},x_{1}x_{2})+B(x_{1}x_{2}\otimes 1_{H};X_{2},x_{1})\right]
+\gamma _{2}B(x_{1}x_{2}\otimes 1_{H};X_{2},x_{2})%
\end{array}%
=0.\left( \ref{G,gx1, GF1,1H}\right)
\end{equation*}

\begin{equation*}
\begin{array}{c}
2\left[ B(x_{1}x_{2}\otimes 1_{H};1_{A},gx_{2})+B(x_{1}x_{2}\otimes
1_{H};X_{2},g)\right] \\
\gamma _{1}\left[ -B(x_{1}x_{2}\otimes
1_{H};G,gx_{1}x_{2})+B(x_{1}x_{2}\otimes 1_{H};GX_{2},gx_{1})\right]%
\end{array}%
\left( \ref{G,gx1, GF2,g}\right)
\end{equation*}%
\begin{equation*}
\begin{array}{c}
2\alpha \left[ -B(x_{1}x_{2}\otimes 1_{H};G,x_{1})\right] + \\
-\gamma _{1}B(x_{1}x_{2}\otimes 1_{H};X_{1},x_{1})+\gamma _{2}\left[
-1+B(x_{1}x_{2}\otimes 1_{H};1_{A},x_{1}x_{2})-B(x_{1}x_{2}\otimes
1_{H};X_{1},x_{2})\right]%
\end{array}%
=0\left( \ref{G,gx2, GF1,1H}\right)
\end{equation*}%
\begin{equation*}
\begin{array}{c}
2\left[ -B(x_{1}x_{2}\otimes 1_{H};1_{A},gx_{1})-B(x_{1}x_{2}\otimes
1_{H};X_{1},g)\right] \\
+\gamma _{2}\left[ -B(x_{1}x_{2}\otimes
1_{H};G,gx_{1}x_{2})-B(x_{1}x_{2}\otimes 1_{H};GX_{1},gx_{2})\right]%
\end{array}%
=0.\left( \ref{G,gx2, GF2,g}\right)
\end{equation*}

\begin{gather*}
\gamma _{1}\left[ -B(x_{2}\otimes 1_{H};1_{A},1_{H})+B(gx_{1}x_{2}\otimes
1_{H};1_{A},x_{1})\right] + \\
+\gamma _{2}\left[ B(x_{1}\otimes 1_{H};1_{A},1_{H})+B(gx_{1}x_{2}\otimes
1_{H};1_{A},x_{2})\right] =0\left( \ref{G,gx1x2, GF1,1H}\right)
\end{gather*}%
\begin{equation*}
0=0.\left( \ref{G,gx1x2, GF1,x1x2}\right)
\end{equation*}%
\begin{equation*}
\begin{array}{c}
2\alpha B(x_{2}\otimes 1_{H};G,g)+ \\
-\gamma _{1}B(x_{2}\otimes 1_{H};1_{A},gx_{1})+\gamma _{2}\left[
B(x_{1}\otimes 1_{H};1_{A},gx_{1})+B(gx_{1}x_{2}\otimes
1_{H};1_{A},gx_{1}x_{2})\right]%
\end{array}%
=0.\left( \ref{G,gx1x2, GF1,gx1}\right)
\end{equation*}

\begin{equation*}
\begin{array}{c}
2\alpha B(gx_{1}x_{2}\otimes 1_{H};G,gx_{2})+ \\
+\gamma _{1}\left[ -B(x_{2}\otimes 1_{H};1_{A},gx_{2})-B(gx_{1}x_{2}\otimes
1_{H};1_{A},gx_{1}x_{2})\right] +\gamma _{2}B(x_{1}\otimes
1_{H};1_{A},gx_{2})%
\end{array}%
=0.\left( \ref{G,gx1x2, GF1,gx2}\right)
\end{equation*}%
\begin{gather*}
\gamma _{1}\left[ B(x_{2}\otimes 1_{H};G,g)-B(gx_{1}x_{2}\otimes
1_{H};G,gx_{1})\right] \left( \ref{G,gx1x2, GF2,g}\right) \\
+\gamma _{2}\left[ -B(x_{1}\otimes 1_{H};G,g)-B(gx_{1}x_{2}\otimes
1_{H};G,gx_{2})\right] =0.
\end{gather*}%
\begin{equation*}
2B(gx_{1}x_{2}\otimes 1_{H};1_{A},x_{1})-\gamma _{2}B(gx_{1}x_{2}\otimes
1_{H};G,x_{1}x_{2})=0.\left( \ref{G,gx1x2, GF2,x1}\right)
\end{equation*}%
\begin{equation*}
\begin{array}{c}
2B(gx_{1}x_{2}\otimes 1_{H};1_{A},x_{2}) \\
+\gamma _{1}\left[ B(x_{2}\otimes 1_{H};G,x_{2})+B(gx_{1}x_{2}\otimes
1_{H};G,x_{1}x_{2})\right]%
\end{array}%
=0.\left( \ref{G,gx1x2, GF2,x2}\right)
\end{equation*}%
\begin{equation*}
\gamma _{1}B(x_{2}\otimes 1_{H};G,gx_{1}x_{2})-\gamma _{2}B(x_{1}\otimes
1_{H};G,gx_{1}x_{2})=0.\left( \ref{G,gx1x2, GF2,gx1x2}\right)
\end{equation*}%
\begin{gather*}
\gamma _{2}\left[
\begin{array}{c}
B(g\otimes 1_{H};1_{A},g)+B(x_{2}\otimes \ 1_{H};1_{A},gx_{2})+ \\
+B(x_{1}\otimes 1_{H};1_{A},gx_{1})+B(gx_{1}x_{2}\otimes
1_{H};1_{A},gx_{1}x_{2})%
\end{array}%
\right] +\left( \ref{G,gx1x2, GF3,g}\right) \\
+2\alpha \left[ -B(x_{2}\otimes 1_{H};G,g)+B(x_{2}\otimes 1_{H};G,g)\right]
=0.
\end{gather*}

\begin{equation*}
B(x_{2}\otimes \ 1_{H};1_{A},x_{1}x_{2})=0.\left( \ref{G,gx1x2, GF3,x1}%
\right)
\end{equation*}

\begin{equation*}
B(g\otimes 1_{H};X_{2},1_{H})=0\left( \ref{G,gx1x2, GF3,x2}\right)
\end{equation*}

\begin{equation*}
\gamma _{2}B\left( g\otimes 1_{H};1_{A},gx_{1}x_{2}\right) -2\alpha
B(x_{2}\otimes 1_{H};G,gx_{1}x_{2})=0.\left( \ref{G,gx1x2, GF3,gx1x2}\right)
\end{equation*}

\begin{gather*}
-\gamma _{1}\left[ B(g\otimes 1_{H};1_{A},g)+B(x_{2}\otimes \
1_{H};1_{A},gx_{2})+B(x_{1}\otimes 1_{H};1_{A},gx_{1})+B(gx_{1}x_{2}\otimes
1_{H};1_{A},gx_{1}x_{2})\right] +\left( \ref{G,gx1x2, GF4,g}\right) \\
-2\alpha \left[ -B(x_{1}\otimes 1_{H};G,g)-B(gx_{1}x_{2}\otimes
1_{H};G,gx_{2})\right] =0
\end{gather*}%
\begin{equation*}
B(x_{2}\otimes \ 1_{H};1_{A},x_{1}x_{2})=0\left( \ref{G,gx1x2, GF4,x1}\right)
\end{equation*}%
\begin{equation*}
0=0.\left( \ref{G,gx1x2, GF4,x2}\right)
\end{equation*}

\begin{equation*}
-\gamma _{1}B\left( g\otimes 1_{H};1_{A},gx_{1}x_{2}\right) +2\alpha
B(x_{1}\otimes 1_{H};G,gx_{1}x_{2})=0.\left( \ref{G,gx1x2, GF4,gx1x2}\right)
\end{equation*}

\begin{equation*}
\alpha \left[ B(g\otimes 1_{H};G,gx_{2})-B(x_{1}\otimes 1_{H};G,gx_{1}x_{2})%
\right] =0.\left( \ref{G,gx1x2, GF5,gx2}\right)
\end{equation*}

\begin{gather*}
-2\left[ -B(x_{2}\otimes 1_{H};1_{A},1_{H})+B(gx_{1}x_{2}\otimes
1_{H};1_{A},x_{1})\right] +\left( \ref{G,gx1x2, GF6,1H}\right) \\
\gamma _{2}\left[ B(g\otimes 1_{H};G,1_{H})+B(x_{2}\otimes \
1_{H};G,x_{2})+B(gx_{1}x_{2}\otimes 1_{H};G,x_{1}x_{2})\right] \\
=0
\end{gather*}%
\begin{equation*}
\gamma _{2}\left[ B(g\otimes 1_{H};G,gx_{2})-B(x_{1}\otimes
1_{H};G,gx_{1}x_{2})\right] =0.\left( \ref{G,gx1x2, GF6,gx1}\right)
\end{equation*}%
\begin{equation*}
\gamma _{2}\left[ B(g\otimes 1_{H};G,gx_{2})-B(x_{1}\otimes
1_{H};G,gx_{1}x_{2})\right] =0.\left( \ref{G,gx1x2, GF6,gx2}\right)
\end{equation*}%
\begin{gather*}
-\gamma _{1}\left[ B(g\otimes 1_{H};G,1_{H})+B(x_{2}\otimes \
1_{H};G,x_{2})+B(gx_{1}x_{2}\otimes 1_{H};G,x_{1}x_{2})\right] \\
-2\left[ B(x_{1}\otimes 1_{H};1_{A},1_{H})+B(gx_{1}x_{2}\otimes
1_{H};1_{A},x_{2})\right] =0\left( \ref{G,gx1x2, GF7,1H}\right)
\end{gather*}%
\begin{equation*}
B(g\otimes 1_{H};G,x_{1}x_{2})=0.\left( \ref{G,gx1x2, GF7,x1x2}\right)
\end{equation*}%
\begin{equation*}
\gamma _{1}\left[ B(g\otimes 1_{H};G,gx_{1})+B(x_{2}\otimes \
1_{H};G,gx_{1}x_{2})\right] =0.\left( \ref{G,gx1x2, GF7,gx1}\right)
\end{equation*}%
\begin{equation*}
\gamma _{1}\left[ B(g\otimes 1_{H};G,gx_{2})-B(x_{1}\otimes
1_{H};G,gx_{1}x_{2})\right] =0.\left( \ref{G,gx1x2, GF7,gx2}\right)
\end{equation*}

\subsubsection{$X_{1}$}

\begin{equation*}
\lambda B\left( g\otimes 1_{H};1_{A},x_{2}\right) -\gamma _{1}B(g\otimes
1_{H};G,1_{H})-2B(x_{1}\otimes 1_{H};1_{A},1_{H})=0.\left( \ref{X1,g,
X1F11,1H}\right)
\end{equation*}%
\begin{eqnarray*}
&&+\gamma _{1}B(g\otimes 1_{H};G,gx_{1})+\lambda B(g\otimes
1_{H};X_{2},gx_{1})+B(x_{1}\otimes 1_{H};1_{A},gx_{1})+\left( \ref{X1,g,
X1F11,gx1}\right) \\
&&+2B\left( g\otimes 1_{H};1_{A},g\right) +B\left( x_{1}\otimes
1_{H};1_{A},gx_{1}\right) =0
\end{eqnarray*}

\begin{equation*}
\gamma _{1}B(g\otimes 1_{H};G,gx_{2})+2B(x_{1}\otimes
1_{H};1_{A},gx_{2})-2\beta _{1}B\left( g\otimes
1_{H};1_{A},gx_{1}x_{2}\right) =0.\left( \ref{X1,g, X1F11,gx2}\right)
\end{equation*}

\begin{equation*}
2\beta _{1}B\left( g\otimes 1_{H};G,gx_{1}\right) +\lambda B\left( g\otimes
1_{H};G,gx_{2}\right) +2B(x_{1}\otimes 1_{H};G,g)=0.\left( \ref{X1,g,X1F21,g}%
\right)
\end{equation*}%
\begin{equation*}
0=0.\left( \ref{X1,g,X1F21,x1}\right)
\end{equation*}

\begin{equation*}
B\left( g\otimes 1_{H};G,gx_{2}\right) -B\left( x_{1}\otimes
1_{H};G,gx_{1}x_{2}\right) =0\left( \ref{X1,g,X1F21,gx1x2}\right)
\end{equation*}%
\begin{gather*}
+\lambda B\left( g\otimes 1_{H};1_{A},gx_{1}x_{2}\right) +\gamma _{1}B\left(
g\otimes 1_{H};G,gx_{1}\right) +\left( \ref{X1,g, X1F31,g}\right) \\
+2B(g\otimes 1_{H};1_{A},g)+2B(x_{1}\otimes 1_{H};1_{A},gx_{1})=0.
\end{gather*}

\begin{equation*}
-2\beta _{1}B\left( g\otimes 1_{H};1_{A},gx_{1}x_{2}\right) +2B(x_{1}\otimes
1_{H};1_{A},gx_{2})+\gamma _{1}B\left( g\otimes 1_{H};G,gx_{2}\right)
=0.\left( \ref{X1,g, X1F41,g}\right)
\end{equation*}%
\begin{equation*}
0=0.\left( \ref{X1,g, X1F41,x1}\right)
\end{equation*}

\begin{equation*}
0=0\left( \ref{X1,g,X1F51,1H}\right)
\end{equation*}%
\begin{equation*}
B(x_{1}\otimes 1_{H};G,gx_{1}x_{2})-B(g\otimes 1_{H};G,gx_{2})=0.\left( \ref%
{X1,x2,X1F71,gx1x2}\right)
\end{equation*}%
\begin{gather*}
2\beta _{1}\left[ -B(g\otimes 1_{H};1_{A},g)-B(x_{1}\otimes
1_{H};1_{A},gx_{1})\right] + \\
-\lambda B(x_{1}\otimes 1_{H};1_{A},gx_{2})=\gamma _{1}B(x_{1}\otimes
1_{H};G,g)\left( \ref{X1,x1,X1F11,g}\right)
\end{gather*}%
\begin{equation*}
0=0.\left( \ref{X1,x1,X1F11,x1}\right)
\end{equation*}

\begin{equation*}
-2\beta _{1}B(g\otimes 1_{H};1_{A},gx_{1}x_{2})=\gamma _{1}B(x_{1}\otimes
1_{H};G,gx_{1}x_{2})-2B\left( x_{1}\otimes 1_{H};1_{A},gx_{2}\right) .\left( %
\ref{X1,x1,X1F11,gx1x2}\right)
\end{equation*}%
\begin{equation*}
B(g\otimes 1_{H};G,1_{H})=0\left( \ref{X1,x1,X1F21,1}\right)
\end{equation*}%
\begin{equation*}
\lambda B(x_{1}\otimes 1_{H};G,gx_{1}x_{2})=-2B\left( x_{1}\otimes
1_{H};G,g\right) .\left( \ref{X1,x1,X1F21,gx1}\right)
\end{equation*}%
\begin{gather*}
2B(x_{1}\otimes 1_{H};1_{A},1_{H})-\lambda B\left( g\otimes
1_{H};1_{A},x_{2}\right) +\left( \ref{X1,x1,X1F31,1H}\right) \\
+\gamma _{1}B(g\otimes 1_{H};G,1_{H})=0.
\end{gather*}%
\begin{equation*}
2B\left( x_{1}\otimes 1_{H};1_{A},x_{1}x_{2}\right) =0\left( \ref%
{X1,x1,X1F31,x1x2}\right)
\end{equation*}

\begin{gather*}
-\lambda B\left( g\otimes 1_{H};1_{A},gx_{1}x_{2}\right) -\gamma _{1}B\left(
g\otimes 1_{H};G,gx_{1}\right) + \\
-2\left[ B(g\otimes 1_{H};1_{A},g)+B(x_{1}\otimes 1_{H};1_{A},gx_{1})\right]
=0\left( \ref{X1,x1,X1F31,gx1}\right)
\end{gather*}

\begin{equation*}
\gamma _{1}\left[ -B\left( g\otimes 1_{H};G,gx_{2}\right) +B(x_{1}\otimes
1_{H};G,gx_{1}x_{2})\right] =0\left( \ref{X1,x1,X1F31,gx2}\right)
\end{equation*}

\begin{gather*}
2\beta _{1}B\left( g\otimes 1_{H};1_{A},gx_{1}x_{2}\right) + \\
-\gamma _{1}B(x_{1}\otimes 1_{H};G,gx_{1}x_{2})-B(x_{1}\otimes
x_{1};X_{2},gx_{1})=0.\left( \ref{X1,x1,X1F41,gx1}\right)
\end{gather*}

\begin{eqnarray*}
&&-2\beta _{1}B(x_{2}\otimes 1_{H};1_{A},gx_{1})+\lambda \left[ -B(g\otimes
1_{H};1_{A},g)-B(x_{2}\otimes \ 1_{H};1_{A},gx_{2})\right] \\
&=&-\gamma _{1}B(x_{2}\otimes 1_{H};G,g)-2B(gx_{1}x_{2}\otimes
1_{H};1_{A},g).\left( \ref{X1,x2,X1F11,g}\right)
\end{eqnarray*}%
\begin{equation*}
\lambda \left[ -B(x_{2}\otimes \ 1_{H};1_{A},x_{1}x_{2})\right]
=+2B(gx_{1}x_{2}\otimes 1_{H};1_{A},x_{1}).\left( \ref{X1,x2,X1F11,x1}\right)
\end{equation*}%
\begin{equation*}
\gamma _{1}B(x_{2}\otimes 1_{H};G,x_{2})+2B(gx_{1}x_{2}\otimes
1_{H};1_{A},x_{2})=0.\left( \ref{X1,x2,X1F11,x2}\right)
\end{equation*}

\begin{gather*}
\lambda B\left( g\otimes 1_{H};1_{A},gx_{1}x_{2}\right) -\gamma
_{1}B(x_{2}\otimes 1_{H};G,gx_{1}x_{2})\left( \ref{X1,x2,X1F11,gx1x2}\right)
\\
-2B\left( x_{2}\otimes 1_{H};1_{A},gx_{2}\right) -2B(gx_{1}x_{2}\otimes
1_{H};1_{A},gx_{1}x_{2})=0.
\end{gather*}%
\begin{eqnarray*}
&&\lambda \left[ B(g\otimes 1_{H};G,1_{H})+B(x_{2}\otimes \ 1_{H};G,x_{2}%
\right] \\
&=&-2B(gx_{1}x_{2}\otimes 1_{H};G,1_{H}).\left( \ref{X1,x2,X1F21,1H}\right)
\end{eqnarray*}%
\begin{equation*}
B(gx_{1}x_{2}\otimes 1_{H};G,x_{1}x_{2})=0.\left( \ref{X1,x2,X1F21,x1x2}%
\right)
\end{equation*}%
\begin{eqnarray*}
&&\lambda \left[ B(g\otimes 1_{H};G,gx_{1})+B(x_{2}\otimes \
1_{H};G,gx_{1}x_{2})\right] \\
&=&2B(gx_{1}x_{2}\otimes 1_{H};G,gx_{1})-2B\left( x_{2}\otimes
1_{H};G,g\right) .\left( \ref{X1,x2,X1F21,gx1}\right)
\end{eqnarray*}

\begin{equation*}
-2\beta _{1}B(x_{2}\otimes 1_{H};G,gx_{1}x_{2})+\lambda B(g\otimes
1_{H};G,gx_{2})-2B(gx_{1}x_{2}\otimes 1_{H};G,gx_{2})=0.\left( \ref%
{X1,x2,X1F21,gx2}\right)
\end{equation*}%
\begin{gather*}
-\lambda B\left( g\otimes 1_{H};1_{A},gx_{1}x_{2}\right) =-\gamma
_{1}B(x_{2}\otimes 1_{H};G,gx_{1}x_{2})+ \\
-2B(x_{2}\otimes 1_{H};1_{A},gx_{2})-2B(gx_{1}x_{2}\otimes
1_{H};1_{A},gx_{1}x_{2}).\left( \ref{X1,x2,X1F31,gx2}\right)
\end{gather*}

\begin{gather*}
0=-\gamma _{1}\left[ B(g\otimes 1_{H};G,1_{H})+B(x_{2}\otimes \ 1_{H};G,x_{2}%
\right] + \\
-2B(x_{1}\otimes 1_{H};1_{A},1_{H})-2B(gx_{1}x_{2}\otimes
1_{H};1_{A},x_{2}).\left( \ref{X1,x2,X1F41,1H}\right)
\end{gather*}%
\begin{equation*}
0=0.\left( \ref{X1,x2,X1F41,x1x2}\right)
\end{equation*}%
\begin{gather*}
B(x_{1}\otimes 1_{H};1_{A},gx_{1})+B(gx_{1}x_{2}\otimes
1_{H};1_{A},gx_{1}x_{2})\left( \ref{X1,x2,X1F41,gx1}\right) \\
+B(g\otimes 1_{H};1_{A},g)+B(x_{2}\otimes \ 1_{H};1_{A},gx_{2})=0.
\end{gather*}

\begin{equation*}
2\beta _{1}B\left( g\otimes 1_{H};1_{A},gx_{1}x_{2}\right) =\gamma
_{1}B(g\otimes 1_{H};G,gx_{2})+2B(x_{1}\otimes 1_{H};1_{A},gx_{2}).\left( %
\ref{X1,x2,X1F41,gx2}\right)
\end{equation*}

\begin{eqnarray*}
0 &=&\gamma _{1}\left[ B(g\otimes 1_{H};G,gx_{1})+B(x_{2}\otimes \
1_{H};G,gx_{1}x_{2})\right] \left( \ref{X1,x2,X1F51,g}\right) \\
&&+\left[
\begin{array}{c}
2B(g\otimes 1_{H};1_{A},g)+2B(x_{2}\otimes \ 1_{H};1_{A},gx_{2})+ \\
+2B(x_{1}\otimes 1_{H};1_{A},gx_{1})+2B(gx_{1}x_{2}\otimes
1_{H};1_{A},gx_{1}x_{2})%
\end{array}%
\right]
\end{eqnarray*}

\begin{equation*}
0=0.\left( \ref{X1,x2,X1F51,x2}\right)
\end{equation*}%
\begin{eqnarray*}
&&\lambda \left[ B(g\otimes 1_{H};G,gx_{1})+B(x_{2}\otimes \
1_{H};G,gx_{1}x_{2})\right] \\
&=&-2\left[ B(x_{2}\otimes 1_{H};G,g)-B(gx_{1}x_{2}\otimes 1_{H};G,gx_{1})%
\right] \left( \ref{X1,x2,X1F61,g}\right)
\end{eqnarray*}%
\begin{equation*}
\beta _{1}\left[ B(g\otimes 1_{H};G,gx_{1})+B(x_{2}\otimes \
1_{H};G,gx_{1}x_{2})\right] =0\left( \ref{X1,x2,X1F71,g}\right)
\end{equation*}%
\begin{equation*}
B(gx_{1}x_{2}\otimes 1_{H};G,x_{1}x_{2})=0\left( \ref{X1,x2,X1F71,x1}\right)
\end{equation*}%
\begin{equation*}
\lambda B(x_{1}x_{2}\otimes 1_{H};X_{2},g)-\gamma _{1}B(x_{1}x_{2}\otimes
1_{H};G,g)=0\left( \ref{X1,x1x2,X1F11,g}\right)
\end{equation*}%
\begin{equation*}
2\beta _{1}B(x_{1}x_{2}\otimes 1_{H};X_{1},x_{1})+\lambda
B(x_{1}x_{2}\otimes 1_{H};X_{2},x_{1})+\gamma _{1}B(x_{1}x_{2}\otimes
1_{H};G,x_{1})=0\left( \ref{X1,x1x2,X1F11,x1}\right)
\end{equation*}

\begin{equation*}
2\beta _{1}B(x_{1}x_{2}\otimes 1_{H};X_{1},x_{2})+\lambda
B(x_{1}x_{2}\otimes 1_{H};X_{2},x_{2})+\gamma _{1}B(x_{1}x_{2}\otimes
1_{H};G,x_{2})=0.\left( \ref{X1,x1x2,X1F11,x2}\right)
\end{equation*}

\begin{equation*}
\gamma _{1}B(x_{1}x_{2}\otimes 1_{H};G,gx_{1}x_{2})-2B(x_{1}x_{2}\otimes
1_{H};1_{A},gx_{2})=0.\left( \ref{X1,x1x2,X1F11,gx1x2}\right)
\end{equation*}%
\begin{equation*}
0=0.\left( \ref{X1,x1x2,X1F21,1H}\right)
\end{equation*}%
\begin{equation*}
\lambda B(x_{1}x_{2}\otimes 1_{H};GX_{2},gx_{1})+2B(x_{1}x_{2}\otimes
1_{H};G,g)=0.\left( \ref{X1,x1x2,X1F21,gx1}\right)
\end{equation*}

\begin{gather*}
\lambda \left[ +1-B(x_{1}x_{2}\otimes
1_{H};1_{A},x_{1}x_{2})-B(x_{1}x_{2}\otimes
1_{H};X_{2},x_{1})+B(x_{1}x_{2}\otimes 1_{H};X_{1},x_{2})\right] \left( \ref%
{X1,x1x2,X1F31,1H}\right) \\
=0
\end{gather*}%
\begin{equation*}
2B(x_{1}x_{2}\otimes 1_{H};1_{A},gx_{1})+2B(x_{1}x_{2}\otimes
1_{H};X_{1},g)=0\left( \ref{X1,x1x2,X1F31,gx1}\right)
\end{equation*}%
\begin{equation*}
2B(x_{1}x_{2}\otimes 1_{H};1_{A},gx_{2})+\gamma _{1}B(x_{1}x_{2}\otimes
1_{H};GX_{1},gx_{2})=0.\left( \ref{X1,x1x2,X1F31,gx2}\right)
\end{equation*}

\begin{equation*}
\left[
\begin{array}{c}
+1-B(x_{1}x_{2}\otimes 1_{H};1_{A},x_{1}x_{2})-B(x_{1}x_{2}\otimes
1_{H};X_{2},x_{1}) \\
+B(x_{1}x_{2}\otimes 1_{H};X_{1},x_{2})%
\end{array}%
\right] =0.\left( \ref{X1,x1x2,X1F41,1H}\right)
\end{equation*}

\begin{equation*}
\gamma _{1}B(x_{1}x_{2}\otimes 1_{H};GX_{2},gx_{1})+2B(x_{1}x_{2}\otimes
1_{H};X_{2},g)=0.\left( \ref{X1,x1x2,X1F41,gx1}\right)
\end{equation*}%
\begin{gather*}
2B(x_{1}x_{2}\otimes 1_{H};X_{2},g)+\left( \ref{X1,x1x2,X1F51,g}\right) \\
+\gamma _{1}\left[ -B(x_{1}x_{2}\otimes
1_{H};G,gx_{1}x_{2})+B(x_{1}x_{2}\otimes
1_{H};GX_{2},gx_{1})-B(x_{1}x_{2}\otimes 1_{H};GX_{1},gx_{2})\right] =0
\end{gather*}%
\begin{equation*}
2B(x_{1}x_{2}\otimes 1_{H};G,g)+\lambda \left[
\begin{array}{c}
-B(x_{1}x_{2}\otimes 1_{H};G,gx_{1}x_{2}) \\
+B(x_{1}x_{2}\otimes 1_{H};GX_{2},gx_{1})-B(x_{1}x_{2}\otimes
1_{H};GX_{1},gx_{2})%
\end{array}%
\right] =0\left( \ref{X1,x1x2,X1F61,g}\right)
\end{equation*}%
\begin{equation*}
B(x_{1}x_{2}\otimes 1_{H};G,gx_{1}x_{2})+B(x_{1}x_{2}\otimes
1_{H};GX_{1},gx_{2})=0\left( \ref{X1,x1x2,X1F61,gx1x2}\right)
\end{equation*}%
\begin{equation*}
+B(x_{1}x_{2}\otimes 1_{H};G,gx_{1}x_{2})+B(x_{1}x_{2}\otimes
1_{H};GX_{1},gx_{2})=0.\left( \ref{X1,x1x2,X1F81,gx1}\right)
\end{equation*}%
\begin{gather*}
2\beta _{1}\left[ -1+B(x_{1}x_{2}\otimes
1_{H};1_{A},x_{1}x_{2})+B(x_{1}x_{2}\otimes 1_{H};X_{2},x_{1})\right]
+\lambda B(x_{1}x_{2}\otimes 1_{H};X_{2},x_{2})\left( \ref{X1,gx1,X1F11,1H}%
\right) \\
=-\gamma _{1}B(x_{1}x_{2}\otimes 1_{H};G,x_{2})+ \\
-2\left[ B(x_{1}x_{2}\otimes 1_{H};1_{A},gx_{2})+B(x_{1}x_{2}\otimes
1_{H};X_{2},g)\right] .
\end{gather*}

\begin{gather*}
-\gamma _{1}\left[ B(x_{1}x_{2}\otimes
1_{H};G,gx_{1}x_{2})-B(x_{1}x_{2}\otimes 1_{H};GX_{2},gx_{1})\right] \left( %
\ref{X1,gx1,X1F11,gx1}\right) \\
+2\left[ B(x_{1}x_{2}\otimes 1_{H};1_{A},gx_{2})+B(x_{1}x_{2}\otimes
1_{H};X_{2},g)\right] =0.
\end{gather*}%
\begin{gather*}
2\left[ B(x_{1}x_{2}\otimes 1_{H};1_{A},gx_{2})+B(x_{1}x_{2}\otimes
1_{H};X_{2},g)\right] +\left( \ref{X1,gx1,X1F31,g}\right) \\
+\gamma _{1}\left[ -B(x_{1}x_{2}\otimes
1_{H};G,gx_{1}x_{2})+B(x_{1}x_{2}\otimes 1_{H};GX_{2},gx_{1})\right] =0.
\end{gather*}%
\begin{gather*}
-2\beta _{1}B(x_{1}x_{2}\otimes 1_{H};X_{1},x_{1})+\lambda \left[
-1+B(x_{1}x_{2}\otimes 1_{H};1_{A},x_{1}x_{2})-B(x_{1}x_{2}\otimes
1_{H};X_{1},x_{2})\right] +\left( \ref{X1,gx2,X1F11,1H}\right) \\
+\gamma _{1}\left[ -B(x_{1}x_{2}\otimes 1_{H};G,x_{1})\right] =0.
\end{gather*}%
\begin{gather*}
\left( \ref{X1,gx2,X1F11,gx1}\right) \\
+2B(x_{1}x_{2}\otimes 1_{H};1_{A},gx_{1})+2B(x_{1}x_{2}\otimes
1_{H};X_{1},g)=0.
\end{gather*}%
\begin{equation*}
+\gamma _{1}\left[ B(x_{1}x_{2}\otimes
1_{H};G,gx_{1}x_{2})+B(x_{1}x_{2}\otimes 1_{H};GX_{1},gx_{2})\right]
=0\left( \ref{X1,gx2,X1F11,gx2}\right)
\end{equation*}%
\begin{gather*}
2\left[ B(x_{1}x_{2}\otimes 1_{H};1_{A},gx_{1})+B(x_{1}x_{2}\otimes
1_{H};X_{1},g)\right] \left( \ref{X1,gx2,X131,g}\right) \\
=0
\end{gather*}%
\begin{equation*}
\gamma _{1}\left[ -B(x_{1}x_{2}\otimes
1_{H};G,gx_{1}x_{2})+B(x_{1}x_{2}\otimes 1_{H};GX_{1},gx_{2})\right]
=0\left( \ref{X1,gx2,X141,g}\right)
\end{equation*}%
\begin{equation*}
+\lambda \left[ B(x_{1}\otimes 1_{H};1_{A},1_{H})+B(gx_{1}x_{2}\otimes
1_{H};1_{A},x_{2})\right] -\gamma _{1}B(gx_{1}x_{2}\otimes
1_{H};G,1_{H})=0\left( \ref{X1,gx1x2,X1F11,1H}\right)
\end{equation*}%
\begin{equation*}
B(gx_{1}x_{2}\otimes 1_{H};G,x_{1}x_{2})=0\left( \ref{X1,gx1x2,X1F11,x1x2}%
\right)
\end{equation*}%
\begin{gather*}
-2\beta _{1}B(x_{2}\otimes 1_{H};1_{A},gx_{1})+\lambda \left[ B(x_{1}\otimes
1_{H};1_{A},gx_{1})+B(gx_{1}x_{2}\otimes 1_{H};1_{A},gx_{1}x_{2})\right]
\left( \ref{X1,gx1x2,X1F11,gx1}\right) \\
+\gamma _{1}B(gx_{1}x_{2}\otimes 1_{H};G,gx_{1})+2B(gx_{1}x_{2}\otimes
1_{H};1_{A},g)=0.
\end{gather*}%
\begin{gather*}
2\beta _{1}\left[ -B(x_{2}\otimes 1_{H};1_{A},gx_{2})-B(gx_{1}x_{2}\otimes
1_{H};1_{A},gx_{1}x_{2})\right] \left( \ref{X1,gx1x2,X171,gx2}\right) \\
+\lambda B(x_{1}\otimes 1_{H};1_{A},gx_{2})+\gamma _{1}B(gx_{1}x_{2}\otimes
1_{H};G,gx_{2})=0.
\end{gather*}%
\begin{gather*}
2\beta _{1}\left[ B(x_{2}\otimes 1_{H};G,g)-B(gx_{1}x_{2}\otimes
1_{H};G,gx_{1})\right] \left( \ref{X1,gx1x2,XF21,g}\right) \\
+\lambda \left[ -B(x_{1}\otimes 1_{H};G,g)-B(gx_{1}x_{2}\otimes
1_{H};G,gx_{2})\right] =0
\end{gather*}

\begin{equation*}
2\beta _{1}B(x_{2}\otimes 1_{H};G,gx_{1}x_{2})-\lambda B(x_{1}\otimes
1_{H};G,gx_{1}x_{2})=0\left( \ref{X1,gx1x2,X1F21,gx1x2}\right)
\end{equation*}%
\begin{gather*}
+\lambda \left[
\begin{array}{c}
B(g\otimes 1_{H};1_{A},g)+B(x_{2}\otimes \ 1_{H};1_{A},gx_{2})+ \\
+B(x_{1}\otimes 1_{H};1_{A},gx_{1})+B(gx_{1}x_{2}\otimes
1_{H};1_{A},gx_{1}x_{2})%
\end{array}%
\right] \\
-\gamma _{1}\left[ B(x_{2}\otimes 1_{H};G,g)-B(gx_{1}x_{2}\otimes
1_{H};G,gx_{1})\right] =0\left( \ref{X1,gx1x2,X1F31,g}\right)
\end{gather*}%
\begin{gather*}
2B(gx_{1}x_{2}\otimes 1_{H};1_{A},x_{1})+\lambda \left[ B(x_{2}\otimes \
1_{H};1_{A},x_{1}x_{2})\right] \\
=0\left( \ref{X1,gx1x2,X1F31,x1}\right)
\end{gather*}

\begin{gather*}
2B(gx_{1}x_{2}\otimes 1_{H};1_{A},x_{2})+\lambda B(g\otimes
1_{H};X_{2},1_{H})+\left( \ref{X1,gx1x2,X1F31,x2}\right) \\
+\gamma _{1}\left[ B(x_{2}\otimes 1_{H};G,x_{2})+B(gx_{1}x_{2}\otimes
1_{H};G,x_{1}x_{2})\right] =0
\end{gather*}

\begin{gather*}
+\lambda B\left( g\otimes 1_{H};1_{A},gx_{1}x_{2}\right) +\left( \ref%
{X1,gx1x2,X1F31,gx1x2}\right) \\
-\gamma _{1}B(x_{2}\otimes 1_{H};G,gx_{1}x_{2})-\left[ +2B(x_{2}\otimes
1_{H};1_{A},gx_{2})+2B(gx_{1}x_{2}\otimes 1_{H};1_{A},gx_{1}x_{2})\right] =0
\end{gather*}

\begin{eqnarray*}
&&2\beta _{1}\left[
\begin{array}{c}
B(g\otimes 1_{H};1_{A},g)+B(x_{2}\otimes \ 1_{H};1_{A},gx_{2})+ \\
+B(x_{1}\otimes 1_{H};1_{A},gx_{1})+B(gx_{1}x_{2}\otimes
1_{H};1_{A},gx_{1}x_{2})%
\end{array}%
\right] \left( \ref{X1,gx1x2,X1F41,g}\right) \\
&&-\gamma _{1}\left[ B(x_{1}\otimes 1_{H};G,g)+B(gx_{1}x_{2}\otimes
1_{H};G,gx_{2})\right] =0
\end{eqnarray*}%
\begin{eqnarray*}
&&2\beta _{1}B\left( g\otimes 1_{H};1_{A},gx_{1}x_{2}\right) \left( \ref%
{X1,gx1x2,X1F41,gx1x2}\right) \\
&&-\gamma _{1}B(x_{1}\otimes 1_{H};G,gx_{1}x_{2})-2B(x_{1}\otimes
1_{H};1_{A},gx_{2})=0
\end{eqnarray*}%
\begin{gather*}
2\left[ B(x_{1}\otimes 1_{H};1_{A},1_{H})+B(gx_{1}x_{2}\otimes
1_{H};1_{A},x_{2})\right] \left( \ref{X1,gx1x2,X1F51,1H}\right) \\
+\gamma _{1}\left[ B(g\otimes 1_{H};G,1_{H})+B(x_{2}\otimes \
1_{H};G,x_{2})+B(gx_{1}x_{2}\otimes 1_{H};G,x_{1}x_{2})\right] =0
\end{gather*}%
\begin{gather*}
+\gamma _{1}\left[ B(g\otimes 1_{H};G,gx_{1})+B(x_{2}\otimes \
1_{H};G,gx_{1}x_{2})\right] +\left( \ref{X1,gx1x2,X1F51,gx1}\right) \\
-\left[
\begin{array}{c}
2B(g\otimes 1_{H};1_{A},g)+2B(x_{2}\otimes \ 1_{H};1_{A},gx_{2})+ \\
+2B(x_{1}\otimes 1_{H};1_{A},gx_{1})+2B(gx_{1}x_{2}\otimes
1_{H};1_{A},gx_{1}x_{2})%
\end{array}%
\right] =0
\end{gather*}%
\begin{gather*}
2B(gx_{1}x_{2}\otimes 1_{H};G,1_{H})+\left( \ref{X1,gx1x2,X1F61,1H}\right) \\
+\lambda \left[ B(g\otimes 1_{H};G,1_{H})+B(x_{2}\otimes \
1_{H};G,x_{2})+B(gx_{1}x_{2}\otimes 1_{H};G,x_{1}x_{2})\right] =0
\end{gather*}%
\begin{equation*}
2B(gx_{1}x_{2}\otimes 1_{H};G,x_{1}x_{2})+\lambda B(g\otimes
1_{H};G,x_{1}x_{2})=0\left( \ref{X1,gx1x2,X1F61,x1x2}\right)
\end{equation*}%
\begin{gather*}
\lambda \left[ B(g\otimes 1_{H};G,gx_{1})+B(x_{2}\otimes \
1_{H};G,gx_{1}x_{2})\right] + \\
+\left[ 2B(x_{2}\otimes 1_{H};G,g)-2B(gx_{1}x_{2}\otimes 1_{H};G,gx_{1})%
\right] =0\left( \ref{X1,gx1x2,X1F61,gx1}\right)
\end{gather*}%
\begin{gather*}
\beta _{1}\left[ B(g\otimes 1_{H};G,gx_{1})+B(x_{2}\otimes \
1_{H};G,gx_{1}x_{2})\right] +\left( \ref{X1,gx1x2,X1F71,gx1}\right) \\
+B(x_{1}\otimes 1_{H};G,g)+B(gx_{1}x_{2}\otimes 1_{H};G,gx_{2})=0.
\end{gather*}

\subsubsection{$X_{2}$}

\begin{gather*}
\gamma _{2}B(g\otimes 1_{H};G,1_{H})+\left( \ref{X2,g,X2F11,1H}\right) \\
+2B(x_{2}\otimes 1_{H};1_{A},1_{H})=0
\end{gather*}%
\begin{gather*}
2\beta _{2}B(g\otimes 1_{H};1_{A},gx_{1}x_{2})+\gamma _{2}B(g\otimes
1_{H};G,gx_{1})\left( \ref{X2,g,X2F11,gx1}\right) \\
+2B(x_{2}\otimes 1_{H};1_{A},gx_{1})=0
\end{gather*}

\begin{gather*}
\gamma _{2}B(g\otimes 1_{H};G,gx_{2})+\lambda B(g\otimes
1_{H};1_{A},gx_{1}x_{2})\left( \ref{X2,g,X2F11,gx2}\right) \\
+2B(x_{2}\otimes 1_{H};1_{A},gx_{2})+2B(g\otimes 1_{H};1_{A},g)=0.
\end{gather*}

\begin{equation*}
2\beta _{2}B(g\otimes 1_{H};G,gx_{2})+\lambda B(g\otimes
1_{H};G,gx_{1})+2B(x_{2}\otimes 1_{H};G,g)=0\left( \ref{X2,g,X2F21,g}\right)
\end{equation*}%
\begin{equation*}
0=0\left( \ref{X2,g,X2F21,x1}\right)
\end{equation*}%
\begin{equation*}
\lambda B(g\otimes 1_{H};G,x_{1}x_{2})-2B(x_{2}\otimes 1_{H};G,x_{2})=0\ref%
{X2,g,X2F21,x2}
\end{equation*}%
\begin{equation*}
2\beta _{2}B(g\otimes 1_{H};1_{H},gx_{1}x_{2})+\gamma _{2}B(g\otimes
1_{H};G,gx_{1})+2B(x_{2}\otimes 1_{H};1_{H},gx_{1})=0.\left( \ref%
{X2,g,X2F31,g}\right)
\end{equation*}%
\begin{equation*}
\gamma _{2}B(g\otimes 1_{H};G,x_{1}x_{2})+2B(x_{2}\otimes
1_{H};1,x_{1}x_{2})=0\left( \ref{X2,g,X2F31,x2}\right)
\end{equation*}%
\begin{eqnarray*}
-\gamma _{2}B(g\otimes 1_{H};G,x_{1}x_{2})\left( \ref{X2,g,X2F41,x1}\right)
&& \\
+B(x_{2}\otimes 1_{H};X_{2},x_{1})+B\left( x_{2}\otimes
1_{H};1_{A},x_{1}x_{2}\right) &=&0
\end{eqnarray*}%
\begin{equation*}
B\left( g\otimes 1_{H};1_{A},gx_{1}x_{2}\right) +B(g\otimes
1_{H},G,gx_{1})=0\left( \ref{X2,g,X2F61,gx2}\right)
\end{equation*}%
\begin{gather*}
2B(g\otimes 1_{H};G,1_{H})-\lambda B(g\otimes 1_{H};G,x_{1}x_{2})\left( \ref%
{X2,g,X2F71,1H}\right) \\
+2B(g\otimes 1_{H};G,1_{H})+2B(x_{2}\otimes 1_{H};G,x_{2})=0
\end{gather*}%
\begin{equation*}
B(g\otimes 1_{H};G,gx_{1})+B(x_{2}\otimes 1_{H};G,gx_{1}x_{2})=0\left( \ref%
{X2,g,X2F71,gx1}\right)
\end{equation*}

\begin{eqnarray*}
&&-2\beta _{2}B(x_{1}\otimes 1_{H};1_{A},gx_{2})+\gamma _{2}B(x_{1}\otimes
1_{H};G,g)\left( \ref{X2,x1,X2F11,g}\right) \\
&&-\lambda \left[ B(g\otimes 1_{H};1_{A},g)+B(x_{1}\otimes
1_{H};1_{A},gx_{1})\right] -2B(gx_{1}x_{2}\otimes 1_{H};1_{A},g)=0
\end{eqnarray*}

\begin{equation*}
B(gx_{1}x_{2}\otimes 1_{H};1_{A},x_{1})=0\left( \ref{X2,x1,X2F11,x1}\right)
\end{equation*}%
\begin{eqnarray*}
&&\left( \ref{X2,x1,X2F11,x2}\right) \\
&&+2B(gx_{1}x_{2}\otimes 1_{H};1_{A},x_{2})=0.
\end{eqnarray*}%
\begin{eqnarray*}
&&+\gamma _{2}B(x_{1}\otimes 1_{H};G,gx_{1}x_{2})-\lambda B(g\otimes
1_{H};1_{A},gx_{1}x_{2})\left( \ref{X2,x1,X2F11,gx1x2}\right) \\
&&-2B(gx_{1}x_{2}\otimes 1_{H};1_{A},gx_{1}x_{2})-2B\left( x_{1}\otimes
1_{H};1_{A},gx_{1}\right) =0
\end{eqnarray*}%
\begin{equation*}
\lambda B(g\otimes 1_{H};G,1_{H})+2B(gx_{1}x_{2}\otimes
1_{H};G,1_{H})=0\left( \ref{X2,x1,X2F21,1H}\right)
\end{equation*}

\begin{equation*}
2\beta _{2}B(x_{1}\otimes 1_{H};G,gx_{1}x_{2})+\lambda B\left( g\otimes
1_{H};G,gx_{1}\right) +2B(x_{2}\otimes 1_{H};G,g)=0.\left( \ref%
{X2,x1,X2F21,gx1}\right)
\end{equation*}

\begin{eqnarray*}
&&\lambda \left[ B\left( g\otimes 1_{H};G,gx_{2}\right) -B(x_{1}\otimes
1_{H};G,gx_{1}x_{2})\right] \left( \ref{X2,x1,X2F21,gx2}\right) \\
&&+2-B(x_{1}\otimes 1_{H};G,g_{1})+2B\left( x_{1}\otimes 1_{H};G,g\right) =0.
\end{eqnarray*}

\begin{gather*}
\gamma _{2}B(g\otimes 1_{H};G,1_{H})+\left( \ref{X2,x1,X2F31,1H}\right) \\
+2B(x_{2}\otimes 1_{H};1_{A},1_{H})-2B(gx_{1}x_{2}\otimes
1_{H};1_{A},x_{1})=0.
\end{gather*}%
\begin{gather*}
\gamma _{2}\left[ -B(g\otimes 1_{H};G,gx_{2})+B(x_{1}\otimes
1_{H};G,gx_{1}x_{2})\right] +\left( \ref{X2,x1,X2F31,gx2}\right) \\
-2B(x_{2}\otimes 1_{H};1_{A},gx_{2})-2B(gx_{1}x_{2}\otimes
1_{H};1_{A},gx_{1}x_{2}) \\
-2B(g\otimes 1_{H};1_{A},g)-2B(x_{1}\otimes 1_{H};1_{A},gx_{1})=0.
\end{gather*}%
\begin{gather*}
\left( \ref{X2,x1,X2F41,1H}\right) \\
\\
B(gx_{1}x_{2}\otimes 1_{H};1_{A},x_{2})=0
\end{gather*}

\begin{eqnarray*}
&&-\gamma _{2}B(x_{1}\otimes 1_{H};G,gx_{1}x_{2})+\lambda B\left( g\otimes
1_{H};1_{A},gx_{1}x_{2}\right) \left( \ref{X2,x1,X2F41,gx1}\right) \\
&&+2B(x_{1}\otimes 1_{H};1_{A},gx_{1})+2B(gx_{1}x_{2}\otimes
1_{H};1_{A},gx_{1}x_{2})=0.
\end{eqnarray*}

\begin{eqnarray*}
&&\beta _{2}\left[ -B(g\otimes 1_{H};G,gx_{2})+B(x_{1}\otimes
1_{H};G,gx_{1}x_{2})\right] \left( \ref{X2,x1,X2F61,g}\right) \\
&&-B(x_{2}\otimes 1_{H};G,g)+B(gx_{1}x_{2}\otimes 1_{H};G,gx_{1})=0
\end{eqnarray*}%
\begin{equation*}
B(x_{2}\otimes 1_{H};G,x_{2})+B(gx_{1}x_{2}\otimes
1_{H};G,x_{1}x_{2})=0\left( \ref{X2,x1,X2F61,x2}\right)
\end{equation*}

\begin{equation*}
B(x_{1}\otimes 1_{H};G,g)+B(gx_{1}x_{2}\otimes 1_{H};G,gx_{2})=0\left( \ref%
{X2,x1,X2F71,g}\right)
\end{equation*}%
\begin{equation*}
B(x_{2}\otimes \ 1_{H};G,x_{2})+B(gx_{1}x_{2}\otimes
1_{H};G,x_{1}x_{2})=0\left( \ref{X2,x1,X2F81,1H}\right)
\end{equation*}%
\begin{eqnarray*}
2 &&\beta _{2}\left[ B(g\otimes 1_{H};1_{A},g)+B(x_{2}\otimes \
1_{H};1_{A},gx_{2})\right] \left( \ref{X2,x2,X2F11,g}\right) \\
&&-\gamma _{2}B(x_{2}\otimes 1_{H};G,g)+\lambda B(x_{2}\otimes
1_{H};1_{A},gx_{1})=0.
\end{eqnarray*}%
\begin{equation*}
B(x_{2}\otimes 1_{H};G,x_{2})=0.\left( \ref{X2,x2,X2F11,x2}\right)
\end{equation*}

\begin{eqnarray*}
+ &&2\beta _{2}B\left( g\otimes 1_{H};1_{A},gx_{1}x_{2}\right) -\gamma
_{2}B(x_{2}\otimes 1_{H};G,gx_{1}x_{2})\left( \ref{X2,x2,X2F11,gx1x2}\right)
\\
&&+2B\left( x_{2}\otimes 1_{H};1_{A},gx_{1}\right) =0
\end{eqnarray*}

\begin{equation*}
\beta _{2}\left[ B(g\otimes 1_{H};G,gx_{1})+B(x_{2}\otimes \
1_{H};G,gx_{1}x_{2})\right] =0\left( \ref{X2,x2,X2F21,gx1}\right)
\end{equation*}%
\begin{eqnarray*}
&&2\beta _{2}B(g\otimes 1_{H};G,gx_{2})-\lambda B(x_{2}\otimes
1_{H};G,gx_{1}x_{2})\left( \ref{X2,x2,X2F21,gx2}\right) \\
&&+2B\left( x_{2}\otimes 1_{H};G,g\right) =0
\end{eqnarray*}

\begin{eqnarray*}
&&2B(x_{2}\otimes 1_{H};1_{A},1_{H})+\gamma _{2}\left[ B(g\otimes
1_{H};G,1_{H})+B(x_{2}\otimes \ 1_{H};G,x_{2}\right] \left( \ref%
{X2,x2,X2F41,1H}\right) \\
&&+\lambda B(x_{2}\otimes \ 1_{H};1_{A},x_{1}x_{2})=0.
\end{eqnarray*}

\begin{eqnarray*}
&&\gamma _{2}B(g\otimes 1_{H};G,gx_{2})-\lambda B\left( g\otimes
1_{H};1_{A},gx_{1}x_{2}\right) \left( \ref{X2,x2,X2F41,gx2}\right) \\
&&+2\left[ +B(g\otimes 1_{H};1_{A},g)+B(x_{2}\otimes \ 1_{H};1_{A},gx_{2})%
\right] =0.
\end{eqnarray*}

\begin{equation*}
\lambda \left[ B(g\otimes 1_{H};G,gx_{1})+B(x_{2}\otimes \
1_{H};G,gx_{1}x_{2})\right] =0\left( \ref{X2,x2,X2F71,g}\right)
\end{equation*}%
\begin{equation*}
\gamma _{2}B(x_{1}x_{2}\otimes 1_{H};G,g)+\lambda B(x_{1}x_{2}\otimes
1_{H};X_{1},g)=0.\left( \ref{X2,x1x2,X2F11,g}\right)
\end{equation*}

\begin{gather*}
2\beta _{2}B(x_{1}x_{2}\otimes 1_{H};X_{2},x_{1})+\gamma
_{2}B(x_{1}x_{2}\otimes 1_{H};G,x_{1})\left( \ref{X2,x1x2,X2F11,x1}\right) \\
+\lambda B(x_{1}x_{2}\otimes 1_{H};X_{1},x_{1})=0.
\end{gather*}

\begin{gather*}
2\beta _{2}B(x_{1}x_{2}\otimes 1_{H};X_{2},x_{2})+\gamma
_{2}B(x_{1}x_{2}\otimes 1_{H};G,x_{2})\left( \ref{X2,x1x2,X2F11,x2}\right) \\
+\lambda B(x_{1}x_{2}\otimes 1_{H};X_{1},x_{2})=0.
\end{gather*}%
\begin{eqnarray*}
&&\gamma _{2}B(x_{1}x_{2}\otimes 1_{H};G,gx_{1}x_{2})+\left( \ref%
{X2,x1x2,X2F11,gx1x2}\right) \\
&&+2B(x_{1}x_{2}\otimes 1_{H};1_{A},gx_{1})=0.
\end{eqnarray*}

\begin{equation*}
0=0.\left( \ref{X2,x1x2,X2F21,1H}\right)
\end{equation*}

\begin{equation*}
2B(x_{1}x_{2}\otimes 1_{H};G,g)-\lambda B(x_{1}x_{2}\otimes
1_{H};GX_{1},gx_{2})=0.\left( \ref{X2,x1x2,X2F21,gx2}\right)
\end{equation*}%
\begin{eqnarray*}
&&2\beta _{2}\left[
\begin{array}{c}
+1-B(x_{1}x_{2}\otimes 1_{H};1_{A},x_{1}x_{2}) \\
-B(x_{1}x_{2}\otimes 1_{H};X_{2},x_{1})+B(x_{1}x_{2}\otimes
1_{H};X_{1},x_{2})%
\end{array}%
\right] \left( \ref{X2,x1x2,X2F31,1H}\right) \\
&=&0.
\end{eqnarray*}%
\begin{equation*}
\gamma _{2}B(x_{1}x_{2}\otimes 1_{H};GX_{1},gx_{2})+2B(x_{1}x_{2}\otimes
1_{H};X_{1},g)=0.\left( \ref{X2,x1x2,X2F31,gx2}\right)
\end{equation*}

\begin{equation*}
\left[
\begin{array}{c}
+1-B(x_{1}x_{2}\otimes 1_{H};1_{A},x_{1}x_{2}) \\
-B(x_{1}x_{2}\otimes 1_{H};X_{2},x_{1})+B(x_{1}x_{2}\otimes
1_{H};X_{1},x_{2})%
\end{array}%
\right] =0.\left( \ref{X2,x1x2,X2F41,1H}\right)
\end{equation*}%
\begin{gather*}
2B(x_{1}x_{2}\otimes 1_{H};1_{A},gx_{2})+\left( \ref{X2,x1x2,X2F41,gx2}%
\right) \\
+2B(x_{1}x_{2}\otimes 1_{H};X_{2},g)=0.
\end{gather*}

\begin{gather*}
2B(x_{1}x_{2}\otimes 1_{H};X_{1},g)-\left( \ref{X2,x1x2,X2F51,g}\right) \\
\gamma _{2}\left[
\begin{array}{c}
-B(x_{1}x_{2}\otimes 1_{H};G,gx_{1}x_{2})+ \\
B(x_{1}x_{2}\otimes 1_{H};GX_{2},gx_{1})-B(x_{1}x_{2}\otimes
1_{H};GX_{1},gx_{2})%
\end{array}%
\right] =0.
\end{gather*}%
\begin{gather*}
2B(x_{1}x_{2}\otimes 1_{H};G,g)\left( \ref{X2,x1x2,X2F71,g}\right) \\
+\lambda \left[
\begin{array}{c}
-B(x_{1}x_{2}\otimes 1_{H};G,gx_{1}x_{2})+ \\
B(x_{1}x_{2}\otimes 1_{H};GX_{2},gx_{1})-B(x_{1}x_{2}\otimes
1_{H};GX_{1},gx_{2})%
\end{array}%
\right] =0.
\end{gather*}%
\begin{equation*}
B(x_{1}x_{2}\otimes 1_{H};G,gx_{1}x_{2})-B(x_{1}x_{2}\otimes
1_{H};GX_{2},gx_{1})=0.\left( \ref{X2,x1x2,X2F71,gx1x2}\right)
\end{equation*}%
\begin{gather*}
2\beta _{2}B(x_{1}x_{2}\otimes 1_{H};X_{2},x_{2})+\gamma
_{2}B(x_{1}x_{2}\otimes 1_{H};G,x_{2})\left( \ref{X2,gx1,X2F11,1H}\right) \\
+\lambda \left[ -1+B(x_{1}x_{2}\otimes
1_{H};1_{A},x_{1}x_{2})+B(x_{1}x_{2}\otimes 1_{H};X_{2},x_{1})\right] =0.
\end{gather*}

\begin{equation*}
B(x_{1}x_{2}\otimes 1_{H};1_{A},gx_{2})+B(x_{1}x_{2}\otimes
1_{H};X_{2},g)=0.\left( \ref{X2,gx1,X2F11,gx2}\right)
\end{equation*}%
\begin{equation*}
\lambda \left[ B(x_{1}x_{2}\otimes 1_{H};G,gx_{1}x_{2})-B(x_{1}x_{2}\otimes
1_{H};GX_{2},gx_{1})\right] =0.\left( \ref{X2,gx1,X2F21,g}\right)
\end{equation*}%
\begin{gather*}
2\beta _{2}\left[ -1+B(x_{1}x_{2}\otimes
1_{H};1_{A},x_{1}x_{2})-B(x_{1}x_{2}\otimes 1_{H};X_{1},x_{2})\right] \left( %
\ref{X2,gx2,X2F11,1H}\right) \\
+\gamma _{2}\left[ -B(x_{1}x_{2}\otimes 1_{H};G,x_{1})\right] + \\
-\lambda B(x_{1}x_{2}\otimes 1_{H};X_{1},x_{1})=0.
\end{gather*}

\begin{eqnarray*}
&&\gamma _{2}\left[ B(x_{1}x_{2}\otimes
1_{H};G,gx_{1}x_{2})+B(x_{1}x_{2}\otimes 1_{H};GX_{1},gx_{2})\right] \left( %
\ref{X2,gx2,X2F11,gx2}\right) \\
&&+ \\
&&+\left[ +2B(x_{1}x_{2}\otimes 1_{H};1_{A},gx_{1})+2B(x_{1}x_{2}\otimes
1_{H};X_{1},g)\right] =0.
\end{eqnarray*}%
\begin{eqnarray*}
&&2\left[ -B(x_{1}x_{2}\otimes 1_{H};1_{A},gx_{1})-B(x_{1}x_{2}\otimes
1_{H};X_{1},g)\right] \left( \ref{X2,gx2,X2F41,g}\right) \\
&&+\gamma _{2}\left[ -B(x_{1}x_{2}\otimes
1_{H};G,gx_{1}x_{2})-B(x_{1}x_{2}\otimes 1_{H};GX_{1},gx_{2})\right] \\
&=&0.
\end{eqnarray*}%
\begin{equation*}
-\gamma _{2}B(gx_{1}x_{2}\otimes 1_{H};G,1_{H})+\lambda \left[
+B(x_{2}\otimes 1_{H};1_{A},1_{H})-B(gx_{1}x_{2}\otimes 1_{H};1_{A},x_{1})%
\right] =0\left( \ref{X2,gx1x2,X2F11,1H}\right)
\end{equation*}%
\begin{equation*}
B(gx_{1}x_{2}\otimes 1_{H};G,x_{1}x_{2})=0\left( \ref{X2,gx1x2,X2F11,x1x2}%
\right)
\end{equation*}%
\begin{gather*}
2\beta _{2}\left[ B(x_{1}\otimes 1_{H};1_{A},gx_{1})+B(gx_{1}x_{2}\otimes
1_{H};1_{A},gx_{1}x_{2})\right] \left( \ref{X2,gx1x2,X2F11,gx1}\right) \\
+\gamma _{2}B(gx_{1}x_{2}\otimes 1_{H};G,gx_{1})-\lambda B(x_{2}\otimes
1_{H};1_{A},gx_{1})=0.
\end{gather*}

\begin{gather*}
2\beta _{2}B(x_{1}\otimes 1_{H};1_{A},gx_{2})+\gamma
_{2}B(gx_{1}x_{2}\otimes 1_{H};G,gx_{2})\left( \ref{X2,gx1x2,X2F11,gx2}%
\right) \\
+\lambda \left[ -B(x_{2}\otimes 1_{H};1_{A},gx_{2})-B(gx_{1}x_{2}\otimes
1_{H};1_{A},gx_{1}x_{2})\right] \\
+2B(gx_{1}x_{2}\otimes 1_{H};1_{A},g)=0.
\end{gather*}

\begin{gather*}
2\beta _{2}\left[ -B(x_{1}\otimes 1_{H};G,g)-B(gx_{1}x_{2}\otimes
1_{H};G,gx_{2})\right] \left( \ref{X2,gx1x2,X2F21,g}\right) \\
+\lambda \left[ B(x_{2}\otimes 1_{H};G,g)-B(gx_{1}x_{2}\otimes
1_{H};G,gx_{1})\right] =0.
\end{gather*}%
\begin{gather*}
-2\beta _{2}B(x_{1}\otimes 1_{H};G,gx_{1}x_{2})\left( \ref%
{X2,gx1x2,X2F21,gx1x2}\right) \\
+\lambda B(x_{2}\otimes 1_{H};G,gx_{1}x_{2})-2B(gx_{1}x_{2}\otimes
1_{H};G,gx_{1})=0.
\end{gather*}%
\begin{gather*}
2\beta _{2}\left[
\begin{array}{c}
B(g\otimes 1_{H};1_{A},g)+B(x_{2}\otimes \ 1_{H};1_{A},gx_{2}) \\
+B(x_{1}\otimes 1_{H};1_{A},gx_{1})+B(gx_{1}x_{2}\otimes
1_{H};1_{A},gx_{1}x_{2})%
\end{array}%
\right] \left( \ref{X2,gx1x2,X2F31,g}\right) \\
-\gamma _{2}\left[ B(x_{2}\otimes 1_{H};G,g)-B(gx_{1}x_{2}\otimes
1_{H};G,gx_{1})\right] =0.
\end{gather*}

\begin{eqnarray*}
&&\gamma _{2}\left[ B(x_{1}\otimes 1_{H};G,g)+B(gx_{1}x_{2}\otimes
1_{H};G,gx_{2})\right] \left( \ref{X2,gx1x2,X2F41,g}\right) \\
&&-\lambda \left[
\begin{array}{c}
B(g\otimes 1_{H};1_{A},g)+B(x_{2}\otimes \ 1_{H};1_{A},gx_{2}) \\
+B(x_{1}\otimes 1_{H};1_{A},gx_{1})+B(gx_{1}x_{2}\otimes
1_{H};1_{A},gx_{1}x_{2})%
\end{array}%
\right] =0.
\end{eqnarray*}

\begin{eqnarray*}
&&2B(gx_{1}x_{2}\otimes 1_{H};1_{A},x_{1})-\gamma _{2}B(gx_{1}x_{2}\otimes
1_{H};G,x_{1}x_{2})\left( \ref{X2,gx1x2,X2F41,x1}\right) \\
&&+\lambda B(x_{2}\otimes \ 1_{H};1_{A},x_{1}x_{2})=0.
\end{eqnarray*}

\begin{gather*}
2\left[ -B(x_{2}\otimes 1_{H};1_{A},1_{H})+B(gx_{1}x_{2}\otimes
1_{H};1_{A},x_{1})\right] +\left( \ref{X2,gx1x2,X2F51,1H}\right) \\
-\gamma _{2}\left[
\begin{array}{c}
B(g\otimes 1_{H};G,1_{H})+B(x_{2}\otimes \ 1_{H};G,x_{2}) \\
+B(gx_{1}x_{2}\otimes 1_{H};G,x_{1}x_{2})%
\end{array}%
\right] =0.
\end{gather*}

\begin{equation*}
B(x_{2}\otimes 1_{H};G,g)-B(gx_{1}x_{2}\otimes 1_{H};G,gx_{1})=0.\left( \ref%
{X2,gx1x2,X2F61,gx2}\right)
\end{equation*}%
\begin{gather*}
2B(gx_{1}x_{2}\otimes 1_{H};G,1_{H})\left( \ref{X2,gx1x2,X2F71,1H}\right) \\
+\lambda \left[
\begin{array}{c}
B(g\otimes 1_{H};G,1_{H})+B(x_{2}\otimes \ 1_{H};G,x_{2}) \\
+B(gx_{1}x_{2}\otimes 1_{H};G,x_{1}x_{2})%
\end{array}%
\right] =0.
\end{gather*}

\subsection{LIST\ OF\ MONOMIAL\ EQUALITIES 2 \label{LME2}}

We list the new monomial equalities we obtained so far and take out constants%
\begin{equation*}
B(g\otimes 1_{H};G,1_{H})=0\left( \ref{X1,x1,X1F21,1}\right)
\end{equation*}%
\begin{equation*}
B(g\otimes 1_{H};X_{2},1_{H})=0\left( \ref{G,x1,GF3,1H}\right)
\end{equation*}%
\begin{equation*}
B(g\otimes 1_{H};G,x_{1}x_{2})=0.\left( \ref{G,x1,GF6,x1}\right) \text{this
is already in }\ref{LME0}
\end{equation*}%
\begin{equation*}
B(g\otimes 1_{H};G,x_{1}x_{2})=0.\left( \ref{G,x2; GF6,x2}\right)
\end{equation*}%
\begin{equation*}
B(g\otimes 1_{H};G,x_{1}x_{2})=0.\left( \ref{G,gx1x2, GF7,x1x2}\right)
\end{equation*}%
\begin{equation*}
B(g\otimes 1_{H};X_{2},1_{H})=0\left( \ref{G,gx1x2, GF3,x2}\right)
\end{equation*}%
\begin{equation*}
B\left( x_{1}\otimes 1_{H};1_{A},x_{1}x_{2}\right) =0\left( \ref%
{X1,x1,X1F31,x1x2}\right)
\end{equation*}%
\begin{equation*}
B(x_{2}\otimes \ 1_{H};1_{A},x_{1}x_{2})=0\left( \ref{G,x2, GF1,x1}\right)
\end{equation*}%
\begin{equation*}
B(x_{2}\otimes \ 1_{H};1_{A},x_{1}x_{2})=0.\left( \ref{G,x2, GF4,1H}\right)
\end{equation*}%
\begin{equation*}
B(x_{2}\otimes \ 1_{H};1_{A},x_{1}x_{2})=0.\left( \ref{G,x2; GF7,x1}\right)
\end{equation*}%
\begin{equation*}
B(x_{2}\otimes \ 1_{H};1_{A},x_{1}x_{2})=0.\left( \ref{G,gx1x2, GF3,x1}%
\right)
\end{equation*}%
\begin{equation*}
B(x_{2}\otimes \ 1_{H};1_{A},x_{1}x_{2})=0\left( \ref{G,gx1x2, GF4,x1}\right)
\end{equation*}%
\begin{equation*}
B(x_{2}\otimes 1_{H};G,x_{2})=0.\left( \ref{X2,x2,X2F11,x2}\right)
\end{equation*}

\begin{equation*}
B(gx_{1}x_{2}\otimes 1_{H};1_{A},x_{1})=0\left( \ref{X2,x1,X2F11,x1}\right)
\end{equation*}%
\begin{equation*}
B(gx_{1}x_{2}\otimes 1_{H};1_{A},x_{1})=0\left( \ref{X2,x1,X2F11,x1}\right)
\end{equation*}%
\begin{equation*}
B(gx_{1}x_{2}\otimes 1_{H};1_{A},x_{2})=0.\left( \ref{X2,x1,X2F11,x2}\right)
\end{equation*}%
\begin{equation*}
B(gx_{1}x_{2}\otimes 1_{H};1_{A},x_{2})=0\left( \ref{X2,x1,X2F41,1H}\right)
\end{equation*}%
\begin{equation*}
B(gx_{1}x_{2}\otimes 1_{H};G,x_{1}x_{2})=0.\left( \ref{X1,x2,X1F21,x1x2}%
\right)
\end{equation*}%
\begin{equation*}
B(gx_{1}x_{2}\otimes 1_{H};G,x_{1}x_{2})=0\left( \ref{X1,x2,X1F71,x1}\right)
\end{equation*}%
\begin{equation*}
B(gx_{1}x_{2}\otimes 1_{H};G,x_{1}x_{2})=0\left( \ref{X1,gx1x2,X1F11,x1x2}%
\right)
\end{equation*}%
\begin{equation*}
B(gx_{1}x_{2}\otimes 1_{H};G,x_{1}x_{2})=0\left( \ref{X1,gx1x2,X1F11,x1x2}%
\right)
\end{equation*}%
\begin{equation*}
B(gx_{1}x_{2}\otimes 1_{H};G,x_{1}x_{2})=0\left( \ref{X2,gx1x2,X2F11,x1x2}%
\right)
\end{equation*}%
We take one of the repetition as before and relabel them%
\begin{equation}
B(g\otimes 1_{H};G,1_{H})=0\left( \ref{X1,x1,X1F21,1}\right)
\label{got1,G,1}
\end{equation}%
\begin{equation}
B(g\otimes 1_{H};X_{2},1_{H})=0\left( \ref{G,x1,GF3,1H}\right)
\label{got1,X2,1}
\end{equation}%
\begin{equation*}
B\left( x_{1}\otimes 1_{H};1_{A},x_{1}x_{2}\right) =0\left( \ref%
{X1,x1,X1F31,x1x2}\right) \text{ this was already in }\left( \ref{LME1}%
\right)
\end{equation*}%
\begin{equation*}
B(x_{2}\otimes \ 1_{H};1_{A},x_{1}x_{2})=0.\left( \ref{G,x2, GF4,1H}\right)
\text{ this was already in }\left( \ref{LME1}\right)
\end{equation*}%
\begin{equation}
B(x_{2}\otimes 1_{H};G,x_{2})=0.\left( \ref{X2,x2,X2F11,x2}\right)
\label{x2ot1,G,x2}
\end{equation}

\begin{equation}
B(gx_{1}x_{2}\otimes 1_{H};1_{A},x_{1})=0\left( \ref{X2,x1,X2F11,x1}\right)
\label{gx1x2ot1,1,x1}
\end{equation}%
\begin{equation}
B(gx_{1}x_{2}\otimes 1_{H};1_{A},x_{2})=0.\left( \ref{X2,x1,X2F11,x2}\right)
\label{gx1x2ot1,1,x2}
\end{equation}%
\begin{equation*}
B(gx_{1}x_{2}\otimes 1_{H};G,x_{1}x_{2})=0.\left( \ref{X1,x2,X1F21,x1x2}%
\right) \text{ this was already in }\left( \ref{LME1}\right)
\end{equation*}%
Now we take the complete list of equalities 1 $\left( \ref{LAE2}\right) $
and cancel there the terms above.

\subsection{LIST\ OF\ ALL\ EQUALITIES 3\label{LAE3}}

\subsubsection{$G$}

\begin{equation*}
\gamma _{1}B\left( g\otimes 1_{H};1_{A},x_{1}\right) +\gamma _{2}B\left(
g\otimes 1_{H};1_{A},x_{2}\right) =0\text{ }\left( \ref{G,g,GF1,1H}\right)
\end{equation*}

\begin{equation*}
2\alpha B(g\otimes 1_{H};G,gx_{1})+\gamma _{2}B(g\otimes
1_{H};1_{A},gx_{1}x_{2})=0\left( \ref{G,g, GF1,gx1}\right)
\end{equation*}

\begin{equation*}
2\alpha B(g\otimes 1_{H};G,gx_{2})-\gamma _{1}B\left( g\otimes
1_{H};1_{A},gx_{1}x_{2}\right) =0.\left( \ref{G,g, GF1,gx2}\right)
\end{equation*}

\begin{equation*}
\gamma _{1}B\left( g\otimes 1_{H};G,gx_{1}\right) +\gamma _{2}B\left(
g\otimes 1_{H};G,gx_{2}\right) =0.\left( \ref{G,g, GF2,g}\right)
\end{equation*}%
\begin{equation*}
\begin{array}{c}
2\alpha B(x_{1}\otimes 1_{H};G,g)+ \\
+\gamma _{1}\left[ -B(g\otimes 1_{H};1_{A},g)-B(x_{1}\otimes
1_{H};1_{A},gx_{1})\right] -\gamma _{2}B(x_{1}\otimes 1_{H};1_{A},gx_{2})%
\end{array}%
=0.\left( \ref{G,x1, GF1,g}\right)
\end{equation*}

\begin{equation*}
-\gamma _{1}B(g\otimes 1_{H};1_{A},gx_{1}x_{2})+2\alpha B(x_{1}\otimes
1_{H};G,gx_{1}x_{2})=0.\left( \ref{G,x1, GF1,gx1x2}\right)
\end{equation*}%
\begin{equation*}
B(x_{1}\otimes 1_{H};1_{A},1_{H})=0.\left( \ref{G,x1, GF2,1H}\right)
\end{equation*}%
\begin{equation*}
B(x_{1}\otimes 1_{H};1_{A},x_{1}x_{2})=0\text{ }\left( \ref{G,x1, GF2,x1x2}%
\right)
\end{equation*}%
\begin{equation*}
\gamma _{1}B\left( g\otimes 1_{H};G,gx_{1}\right) +\gamma _{2}B(x_{1}\otimes
1_{H};G,gx_{1}x_{2})=0.\left( \ref{G,x1, GF2,gx1}\right)
\end{equation*}

\begin{equation*}
\gamma _{1}[B\left( g\otimes 1_{H};G,gx_{2}\right) -B(x_{1}\otimes
1_{H};G,gx_{1}x_{2})]=0.\left( \ref{G,x1,GF2,gx2}\right)
\end{equation*}

\begin{equation*}
0=0\left( \ref{G,x1,GF3,1H}\right)
\end{equation*}%
\begin{equation*}
-\gamma _{2}B(g\otimes 1_{H};1_{A},gx_{1}x_{2})-2\alpha B(g\otimes
1_{H};G,gx_{1})=0.\left( \ref{G,x1,GF3,gx1}\right)
\end{equation*}

\begin{equation*}
\alpha \left[ -B(g\otimes 1_{H};g,gx_{2})+B(x_{1}\otimes 1_{H};g,gx_{1}x_{2})%
\right] =0.\left( \ref{G,x1, GF3,x2}\right)
\end{equation*}%
\begin{equation*}
0=0.\left( \ref{G,x1, GF4,1H}\right)
\end{equation*}

\begin{equation*}
-\gamma _{1}B(g\otimes 1_{H};1_{A},gx_{1}x_{2})+2\alpha B(x_{1}\otimes
1_{H};G,gx_{1}x_{2})=0.\left( \ref{G,x1, GF4,gx1}\right)
\end{equation*}

\begin{equation*}
\alpha \lbrack -B(g\otimes 1_{H};g,gx_{2})+B(x_{1}\otimes
1_{H};G,gx_{1}x_{2})]=0.\left( \ref{G,x1, GF5,g}\right)
\end{equation*}%
\begin{equation*}
0=0.\left( \ref{G,x1,GF6,x1}\right)
\end{equation*}%
\begin{equation*}
0=0.\left( \ref{G,x1,GF6,x2}\right)
\end{equation*}%
\begin{equation*}
\begin{array}{c}
2\alpha B(x_{2}\otimes 1_{H};G,g)+ \\
-\gamma _{1}B(x_{2}\otimes 1_{H};1_{A},gx_{1})+\gamma _{2}\left[ -B(g\otimes
1_{H};1_{A},g)-B(x_{2}\otimes \ 1_{H};1_{A},gx_{2})\right]%
\end{array}%
=0.\left( \ref{G,x2, GF1,g}\right)
\end{equation*}%
\begin{equation*}
0=0\left( \ref{G,x2, GF1,x1}\right)
\end{equation*}

\begin{equation*}
2\alpha B(x_{2}\otimes 1_{H};G,gx_{1}x_{2})-\gamma _{2}B\left( g\otimes
1_{H};1_{A},gx_{1}x_{2}\right) =0.\left( \ref{G,x2, GF1,gx1x2}\right)
\end{equation*}%
\begin{equation*}
2B(x_{2}\otimes 1_{H};1_{A},1_{H})+\gamma _{2}B(x_{2}\otimes \
1_{H};G,x_{2})=0.\left( \ref{G,x2, GF2,1H}\right)
\end{equation*}%
\begin{equation*}
0=0\left( \ref{G,x2, GF2,x1x2}\right)
\end{equation*}%
\begin{equation*}
\gamma _{2}\left[ B(g\otimes 1_{H};G,gx_{1})+B(x_{2}\otimes \
1_{H};G,gx_{1}x_{2})\right] =0.\left( \ref{G,x2, GF2,gx1}\right)
\end{equation*}

\begin{equation*}
-\gamma _{1}B(x_{2}\otimes 1_{H};G,gx_{1}x_{2})+\gamma _{2}B(g\otimes
1_{H};G,gx_{2})=0.\left( \ref{G,x2, GF2,gx2}\right)
\end{equation*}%
\begin{equation*}
0=0\left( \ref{G,x2, FG3,1H}\right)
\end{equation*}

\begin{equation*}
-\gamma _{2}B(g\otimes 1_{H};1_{A},gx_{1}x_{2})-2\alpha B(x_{2}\otimes
1_{H};G,gx_{1}x_{2})=0.\left( \ref{G,x2, GF3,gx2}\right)
\end{equation*}%
\begin{equation*}
0=0.\left( \ref{G,x2, GF4,1H}\right)
\end{equation*}%
\begin{equation*}
\alpha \left[ B(g\otimes 1_{H};G,gx_{1})+B(x_{2}\otimes \
1_{H};G,gx_{1}x_{2})\right] =0\left( \ref{G,x2, GF4,gx1}\right)
\end{equation*}

\begin{equation*}
\gamma _{1}B(g\otimes 1_{H};1_{A},gx_{1}x_{2})-2\alpha B(g\otimes
1_{H};G,gx_{2})=0.\left( \ref{G,x2, GF4,gx2}\right)
\end{equation*}%
\begin{equation*}
\alpha \left[ B(g\otimes 1_{H};G,gx_{1})+B(x_{2}\otimes \
1_{H};G,gx_{1}x_{2})\right] =0\left( \ref{G,x2, GF5,g}\right)
\end{equation*}

\begin{equation*}
\gamma _{2}[B(g\otimes 1_{H},G,gx_{1})+B(x_{2}\otimes
1_{H};G,gx_{1}x_{2})=0.=0.\left( \ref{G,x2, GF6,g}\right)
\end{equation*}%
\begin{equation*}
0=0.\left( \ref{G,x2; GF6,x2}\right)
\end{equation*}

\begin{equation*}
\gamma _{1}[B(g\otimes 1_{H};G,gx_{1})+B(x_{2}\otimes
1_{H};G,gx_{1}x_{2})]=0.\left( \ref{G,x2; GF7,g}\right)
\end{equation*}%
\begin{equation*}
0=0.\left( \ref{G,x2; GF7,x1}\right)
\end{equation*}%
\begin{equation*}
\gamma _{1}B(x_{1}x_{2}\otimes 1_{H};X_{1},g)+\gamma _{2}B(x_{1}x_{2}\otimes
1_{H};X_{2},g)=0.\left( \ref{G,x1x2, GF1,g}\right)
\end{equation*}

\begin{equation*}
\begin{array}{c}
2\alpha B(x_{1}x_{2}\otimes 1_{H};G,x_{1})+ \\
+\gamma _{1}B(x_{1}x_{2}\otimes 1_{H};X_{1},x_{1})+\gamma
_{2}B(x_{1}x_{2}\otimes 1_{H};X_{2},x_{1})%
\end{array}%
=0.\left( \ref{G,x1x2, GF1,x1}\right)
\end{equation*}%
\begin{equation*}
\begin{array}{c}
2\alpha B(x_{1}x_{2}\otimes 1_{H};G,x_{2})+ \\
+\gamma _{1}B(x_{1}x_{2}\otimes 1_{H};X_{1},x_{2})+\gamma
_{2}B(x_{1}x_{2}\otimes 1_{H};X_{2},x_{2})%
\end{array}%
=0.\left( \ref{G,x1x2, GF1,x2}\right)
\end{equation*}

\begin{equation*}
0=0\left( \ref{G,x1x2, GF2,1H}\right)
\end{equation*}%
\begin{equation*}
\begin{array}{c}
2B(x_{1}x_{2}\otimes 1_{H};1_{A},gx_{1})+ \\
+\gamma _{2}B(x_{1}x_{2}\otimes 1_{H};GX_{2},gx_{1})%
\end{array}%
=0.\left( \ref{G,x1x2, GF2,gx1}\right)
\end{equation*}%
\begin{equation*}
\begin{array}{c}
2B(x_{1}x_{2}\otimes 1_{H};1_{A},gx_{2})+ \\
+\gamma _{1}B(x_{1}x_{2}\otimes 1_{H};GX_{1},gx_{2})%
\end{array}%
=0\left( \ref{G,x1x2, GF2,gx2}\right)
\end{equation*}

\begin{equation*}
0=0\text{ }\left( \ref{G,x1x2, GF3,gx1}\right)
\end{equation*}%
\begin{gather*}
-\gamma _{1}[1-B(x_{1}x_{2}\otimes 1_{H};1_{A},x_{1}x_{2})+ \\
+B(x_{1}x_{2}\otimes 1_{H};x_{1},x_{2})-B(x_{1}x_{2}\otimes
1_{H};x_{2},x_{1})]+ \\
-2\alpha B(x_{1}x_{2}\otimes 1_{H};GX_{2},1_{H}=0.\left( \ref{G,x1x2, GF4,1H}%
\right)
\end{gather*}

\begin{gather*}
-2B(x_{1}x_{2}\otimes 1_{H};X_{1},g)+\gamma _{2}[-B(x_{1}x_{2}\otimes
1_{H};g,gx_{1}x_{2})+ \\
-B(x_{1}x_{2}\otimes 1_{H};GX_{1},gx_{2})+B(x_{1}x_{2}\otimes
1_{H};GX_{2},gx_{1})]=0\left( \ref{G,x1x2, GF6,g}\right)
\end{gather*}%
\begin{gather*}
-2B(x_{1}x_{2}\otimes 1_{H};X_{2},g)+ \\
-\gamma _{1}[-B(x_{1}x_{2}\otimes 1_{H};g,gx_{1}x_{2})+ \\
-B(x_{1}x_{2}\otimes 1_{H};GX_{1},gx_{2})+B(x_{1}x_{2}\otimes
1_{H};GX_{2},gx_{1})]=0\left( \ref{G,x1x2, GF7,g}\right)
\end{gather*}

\begin{equation*}
\begin{array}{c}
2\alpha B(x_{1}x_{2}\otimes 1_{H};G,x_{2})+ \\
+\gamma _{1}\left[
\begin{array}{c}
-1+B(x_{1}x_{2}\otimes 1_{H};1_{A},x_{1}x_{2}) \\
+B(x_{1}x_{2}\otimes 1_{H};X_{2},x_{1})%
\end{array}%
\right] +\gamma _{2}B(x_{1}x_{2}\otimes 1_{H};X_{2},x_{2})%
\end{array}%
=0.\left( \ref{G,gx1, GF1,1H}\right)
\end{equation*}

\begin{equation*}
\begin{array}{c}
2\left[ B(x_{1}x_{2}\otimes 1_{H};1_{A},gx_{2})+B(x_{1}x_{2}\otimes
1_{H};X_{2},g)\right] \\
\gamma _{1}\left[ -B(x_{1}x_{2}\otimes
1_{H};G,gx_{1}x_{2})+B(x_{1}x_{2}\otimes 1_{H};GX_{2},gx_{1})\right]%
\end{array}%
\left( \ref{G,gx1, GF2,g}\right)
\end{equation*}%
\begin{equation*}
\begin{array}{c}
2\alpha \left[ -B(x_{1}x_{2}\otimes 1_{H};G,x_{1})\right] + \\
-\gamma _{1}B(x_{1}x_{2}\otimes 1_{H};X_{1},x_{1})+\gamma _{2}\left[
-1+B(x_{1}x_{2}\otimes 1_{H};1_{A},x_{1}x_{2})-B(x_{1}x_{2}\otimes
1_{H};X_{1},x_{2})\right]%
\end{array}%
=0\left( \ref{G,gx2, GF1,1H}\right)
\end{equation*}%
\begin{equation*}
\begin{array}{c}
2\left[ -B(x_{1}x_{2}\otimes 1_{H};1_{A},gx_{1})-B(x_{1}x_{2}\otimes
1_{H};X_{1},g)\right] \\
+\gamma _{2}\left[ -B(x_{1}x_{2}\otimes
1_{H};G,gx_{1}x_{2})-B(x_{1}x_{2}\otimes 1_{H};GX_{1},gx_{2})\right]%
\end{array}%
=0.\left( \ref{G,gx2, GF2,g}\right)
\end{equation*}

\begin{equation*}
-\gamma _{1}B(x_{2}\otimes 1_{H};1_{A},1_{H})+\gamma _{2}B(x_{1}\otimes
1_{H};1_{A},1_{H})=0\left( \ref{G,gx1x2, GF1,1H}\right)
\end{equation*}%
\begin{equation*}
0=0.\left( \ref{G,gx1x2, GF1,x1x2}\right)
\end{equation*}%
\begin{equation*}
\begin{array}{c}
2\alpha B(x_{2}\otimes 1_{H};G,g)+ \\
-\gamma _{1}B(x_{2}\otimes 1_{H};1_{A},gx_{1})+\gamma _{2}\left[
\begin{array}{c}
B(x_{1}\otimes 1_{H};1_{A},gx_{1}) \\
+B(gx_{1}x_{2}\otimes 1_{H};1_{A},gx_{1}x_{2})%
\end{array}%
\right]%
\end{array}%
=0.\left( \ref{G,gx1x2, GF1,gx1}\right)
\end{equation*}

\begin{equation*}
\begin{array}{c}
2\alpha B(gx_{1}x_{2}\otimes 1_{H};G,gx_{2})+ \\
+\gamma _{1}\left[ -B(x_{2}\otimes 1_{H};1_{A},gx_{2})-B(gx_{1}x_{2}\otimes
1_{H};1_{A},gx_{1}x_{2})\right] +\gamma _{2}B(x_{1}\otimes
1_{H};1_{A},gx_{2})%
\end{array}%
=0.\left( \ref{G,gx1x2, GF1,gx2}\right)
\end{equation*}%
\begin{gather*}
\gamma _{1}\left[ B(x_{2}\otimes 1_{H};G,g)-B(gx_{1}x_{2}\otimes
1_{H};G,gx_{1})\right] \left( \ref{G,gx1x2, GF2,g}\right) \\
+\gamma _{2}\left[ -B(x_{1}\otimes 1_{H};G,g)-B(gx_{1}x_{2}\otimes
1_{H};G,gx_{2})\right] =0.
\end{gather*}%
\begin{equation*}
0=0.\left( \ref{G,gx1x2, GF2,x1}\right)
\end{equation*}%
\begin{equation*}
0=0.\left( \ref{G,gx1x2, GF2,x2}\right)
\end{equation*}%
\begin{equation*}
\gamma _{1}B(x_{2}\otimes 1_{H};G,gx_{1}x_{2})-\gamma _{2}B(x_{1}\otimes
1_{H};G,gx_{1}x_{2})=0.\left( \ref{G,gx1x2, GF2,gx1x2}\right)
\end{equation*}%
\begin{gather*}
\gamma _{2}\left[
\begin{array}{c}
B(g\otimes 1_{H};1_{A},g)+B(x_{2}\otimes \ 1_{H};1_{A},gx_{2})+ \\
+B(x_{1}\otimes 1_{H};1_{A},gx_{1})+B(gx_{1}x_{2}\otimes
1_{H};1_{A},gx_{1}x_{2})%
\end{array}%
\right] +\left( \ref{G,gx1x2, GF3,g}\right) \\
2\alpha \left[ -B(x_{2}\otimes 1_{H};G,g)+B(x_{2}\otimes 1_{H};G,g)\right]
=0.
\end{gather*}

\begin{equation*}
0=0.\left( \ref{G,gx1x2, GF3,x1}\right)
\end{equation*}

\begin{equation*}
0=0\left( \ref{G,gx1x2, GF3,x2}\right)
\end{equation*}

\begin{equation*}
\gamma _{2}B\left( g\otimes 1_{H};1_{A},gx_{1}x_{2}\right) -2\alpha
B(x_{2}\otimes 1_{H};G,gx_{1}x_{2})=0.\left( \ref{G,gx1x2, GF3,gx1x2}\right)
\end{equation*}

\begin{gather*}
-\gamma _{1}\left[
\begin{array}{c}
B(g\otimes 1_{H};1_{A},g)+B(x_{2}\otimes \ 1_{H};1_{A},gx_{2})+ \\
+B(x_{1}\otimes 1_{H};1_{A},gx_{1})+B(gx_{1}x_{2}\otimes
1_{H};1_{A},gx_{1}x_{2})%
\end{array}%
\right] +\left( \ref{G,gx1x2, GF4,g}\right) \\
-2\alpha \left[ -B(x_{1}\otimes 1_{H};G,g)-B(gx_{1}x_{2}\otimes
1_{H};G,gx_{2})\right] =0
\end{gather*}%
\begin{equation*}
0=0\left( \ref{G,gx1x2, GF4,x1}\right)
\end{equation*}%
\begin{equation*}
0=0.\left( \ref{G,gx1x2, GF4,x2}\right)
\end{equation*}

\begin{equation*}
-\gamma _{1}B\left( g\otimes 1_{H};1_{A},gx_{1}x_{2}\right) +2\alpha
B(x_{1}\otimes 1_{H};G,gx_{1}x_{2})=0.\left( \ref{G,gx1x2, GF4,gx1x2}\right)
\end{equation*}

\begin{equation*}
\alpha \left[ B(g\otimes 1_{H};G,gx_{2})-B(x_{1}\otimes 1_{H};G,gx_{1}x_{2})%
\right] =0.\left( \ref{G,gx1x2, GF5,gx2}\right)
\end{equation*}

\begin{equation*}
2B(x_{2}\otimes 1_{H};1_{A},1_{H})+\gamma _{2}B(x_{2}\otimes \
1_{H};G,x_{2})=0\left( \ref{G,gx1x2, GF6,1H}\right)
\end{equation*}%
\begin{equation*}
\gamma _{2}\left[ B(g\otimes 1_{H};G,gx_{2})-B(x_{1}\otimes
1_{H};G,gx_{1}x_{2})\right] =0.\left( \ref{G,gx1x2, GF6,gx1}\right)
\end{equation*}%
\begin{equation*}
\gamma _{2}\left[ B(g\otimes 1_{H};G,gx_{2})-B(x_{1}\otimes
1_{H};G,gx_{1}x_{2})\right] =0.\left( \ref{G,gx1x2, GF6,gx2}\right)
\end{equation*}%
\begin{equation*}
-\gamma _{1}B(x_{2}\otimes \ 1_{H};G,x_{2})-2B(x_{1}\otimes
1_{H};1_{A},1_{H})=0\left( \ref{G,gx1x2, GF7,1H}\right)
\end{equation*}%
\begin{equation*}
0=0.\left( \ref{G,gx1x2, GF7,x1x2}\right)
\end{equation*}%
\begin{equation*}
\gamma _{1}\left[ B(g\otimes 1_{H};G,gx_{1})+B(x_{2}\otimes \
1_{H};G,gx_{1}x_{2})\right] =0.\left( \ref{G,gx1x2, GF7,gx1}\right)
\end{equation*}%
\begin{equation*}
\gamma _{1}\left[ B(g\otimes 1_{H};G,gx_{2})-B(x_{1}\otimes
1_{H};G,gx_{1}x_{2})\right] =0.\left( \ref{G,gx1x2, GF7,gx2}\right)
\end{equation*}

\subsubsection{$X_{1}$}

\begin{equation*}
\lambda B\left( g\otimes 1_{H};1_{A},x_{2}\right) -2B(x_{1}\otimes
1_{H};1_{A},1_{H})=0.\left( \ref{X1,g, X1F11,1H}\right)
\end{equation*}%
\begin{eqnarray*}
&&+\gamma _{1}B(g\otimes 1_{H};G,gx_{1})+\lambda B(g\otimes
1_{H};X_{2},gx_{1})+B(x_{1}\otimes 1_{H};1_{A},gx_{1})+\left( \ref{X1,g,
X1F11,gx1}\right) \\
&&+2B\left( g\otimes 1_{H};1_{A},g\right) +B\left( x_{1}\otimes
1_{H};1_{A},gx_{1}\right) =0
\end{eqnarray*}

\begin{equation*}
\gamma _{1}B(g\otimes 1_{H};G,gx_{2})+2B(x_{1}\otimes
1_{H};1_{A},gx_{2})-2\beta _{1}B\left( g\otimes
1_{H};1_{A},gx_{1}x_{2}\right) =0.\left( \ref{X1,g, X1F11,gx2}\right)
\end{equation*}

\begin{equation*}
2\beta _{1}B\left( g\otimes 1_{H};G,gx_{1}\right) +\lambda B\left( g\otimes
1_{H};G,gx_{2}\right) +2B(x_{1}\otimes 1_{H};G,g)=0.\left( \ref{X1,g,X1F21,g}%
\right)
\end{equation*}%
\begin{equation*}
0=0.\left( \ref{X1,g,X1F21,x1}\right)
\end{equation*}

\begin{equation*}
B\left( g\otimes 1_{H};G,gx_{2}\right) -B\left( x_{1}\otimes
1_{H};G,gx_{1}x_{2}\right) =0\left( \ref{X1,g,X1F21,gx1x2}\right)
\end{equation*}%
\begin{gather*}
+\lambda B\left( g\otimes 1_{H};1_{A},gx_{1}x_{2}\right) +\gamma _{1}B\left(
g\otimes 1_{H};G,gx_{1}\right) +\left( \ref{X1,g, X1F31,g}\right) \\
+2B(g\otimes 1_{H};1_{A},g)+2B(x_{1}\otimes 1_{H};1_{A},gx_{1})=0.
\end{gather*}

\begin{equation*}
-2\beta _{1}B\left( g\otimes 1_{H};1_{A},gx_{1}x_{2}\right) +2B(x_{1}\otimes
1_{H};1_{A},gx_{2})+\gamma _{1}B\left( g\otimes 1_{H};G,gx_{2}\right)
=0.\left( \ref{X1,g, X1F41,g}\right)
\end{equation*}%
\begin{equation*}
0=0.\left( \ref{X1,g, X1F41,x1}\right)
\end{equation*}

\begin{equation*}
0=0\left( \ref{X1,g,X1F51,1H}\right)
\end{equation*}%
\begin{equation*}
B(x_{1}\otimes 1_{H};G,gx_{1}x_{2})-B(g\otimes 1_{H};G,gx_{2})=0.\left( \ref%
{X1,x2,X1F71,gx1x2}\right)
\end{equation*}%
\begin{eqnarray*}
&&2\beta _{1}\left[ -B(g\otimes 1_{H};1_{A},g)-B(x_{1}\otimes
1_{H};1_{A},gx_{1})\right] -\lambda B(x_{1}\otimes 1_{H};1_{A},gx_{2})+ \\
&=&\gamma _{1}B(x_{1}\otimes 1_{H};G,g)\left( \ref{X1,x1,X1F11,g}\right)
\end{eqnarray*}%
\begin{equation*}
0=0.\left( \ref{X1,x1,X1F11,x1}\right)
\end{equation*}

\begin{equation*}
-2\beta _{1}B(g\otimes 1_{H};1_{A},gx_{1}x_{2})=\gamma _{1}B(x_{1}\otimes
1_{H};G,gx_{1}x_{2})-2B\left( x_{1}\otimes 1_{H};1_{A},gx_{2}\right) .\left( %
\ref{X1,x1,X1F11,gx1x2}\right)
\end{equation*}%
\begin{equation*}
0=0\left( \ref{X1,x1,X1F21,1}\right)
\end{equation*}%
\begin{equation*}
\lambda B(x_{1}\otimes 1_{H};G,gx_{1}x_{2})=-2B\left( x_{1}\otimes
1_{H};G,g\right) .\left( \ref{X1,x1,X1F21,gx1}\right)
\end{equation*}%
\begin{equation*}
2B(x_{1}\otimes 1_{H};1_{A},1_{H})-\lambda B\left( g\otimes
1_{H};1_{A},x_{2}\right) =0.\left( \ref{X1,x1,X1F31,1H}\right)
\end{equation*}%
\begin{equation*}
0=0\left( \ref{X1,x1,X1F31,x1x2}\right)
\end{equation*}

\begin{gather*}
-\lambda B\left( g\otimes 1_{H};1_{A},gx_{1}x_{2}\right) -\gamma _{1}B\left(
g\otimes 1_{H};G,gx_{1}\right) + \\
-2\left[ B(g\otimes 1_{H};1_{A},g)+B(x_{1}\otimes 1_{H};1_{A},gx_{1})\right]
=0\left( \ref{X1,x1,X1F31,gx1}\right)
\end{gather*}

\begin{equation*}
\gamma _{1}\left[ -B\left( g\otimes 1_{H};G,gx_{2}\right) +B(x_{1}\otimes
1_{H};G,gx_{1}x_{2})\right] =0\left( \ref{X1,x1,X1F31,gx2}\right)
\end{equation*}

\begin{equation*}
2\beta _{1}B\left( g\otimes 1_{H};1_{A},gx_{1}x_{2}\right) -\gamma
_{1}B(x_{1}\otimes 1_{H};G,gx_{1}x_{2})-B(x_{1}\otimes
x_{1};X_{2},gx_{1})=0.\left( \ref{X1,x1,X1F41,gx1}\right)
\end{equation*}

\begin{eqnarray*}
&&-2\beta _{1}B(x_{2}\otimes 1_{H};1_{A},gx_{1})+\lambda \left[ -B(g\otimes
1_{H};1_{A},g)-B(x_{2}\otimes \ 1_{H};1_{A},gx_{2})\right] \\
&=&-\gamma _{1}B(x_{2}\otimes 1_{H};G,g)-2B(gx_{1}x_{2}\otimes
1_{H};1_{A},g).\left( \ref{X1,x2,X1F11,g}\right)
\end{eqnarray*}%
\begin{equation*}
0=0.\left( \ref{X1,x2,X1F11,x1}\right)
\end{equation*}%
\begin{equation*}
0=0.\left( \ref{X1,x2,X1F11,x2}\right)
\end{equation*}

\begin{gather*}
\lambda B\left( g\otimes 1_{H};1_{A},gx_{1}x_{2}\right) -\gamma
_{1}B(x_{2}\otimes 1_{H};G,gx_{1}x_{2})\left( \ref{X1,x2,X1F11,gx1x2}\right)
\\
-2B\left( x_{2}\otimes 1_{H};1_{A},gx_{2}\right) -2B(gx_{1}x_{2}\otimes
1_{H};1_{A},gx_{1}x_{2})=0.
\end{gather*}%
\begin{equation*}
\lambda B(x_{2}\otimes \ 1_{H};G,x_{2}=-2B(gx_{1}x_{2}\otimes
1_{H};G,1_{H}).\left( \ref{X1,x2,X1F21,1H}\right)
\end{equation*}%
\begin{equation*}
0=0.\left( \ref{X1,x2,X1F21,x1x2}\right)
\end{equation*}%
\begin{eqnarray*}
&&\lambda \left[ B(g\otimes 1_{H};G,gx_{1})+B(x_{2}\otimes \
1_{H};G,gx_{1}x_{2})\right] \\
&=&2B(gx_{1}x_{2}\otimes 1_{H};G,gx_{1})-2B\left( x_{2}\otimes
1_{H};G,g\right) .\left( \ref{X1,x2,X1F21,gx1}\right)
\end{eqnarray*}

\begin{equation*}
-2\beta _{1}B(x_{2}\otimes 1_{H};G,gx_{1}x_{2})+\lambda B(g\otimes
1_{H};G,gx_{2})-2B(gx_{1}x_{2}\otimes 1_{H};G,gx_{2})=0.\left( \ref%
{X1,x2,X1F21,gx2}\right)
\end{equation*}%
\begin{gather*}
-\lambda B\left( g\otimes 1_{H};1_{A},gx_{1}x_{2}\right) =-\gamma
_{1}B(x_{2}\otimes 1_{H};G,gx_{1}x_{2})+ \\
-2B(x_{2}\otimes 1_{H};1_{A},gx_{2})-2B(gx_{1}x_{2}\otimes
1_{H};1_{A},gx_{1}x_{2}).\left( \ref{X1,x2,X1F31,gx2}\right)
\end{gather*}

\begin{equation*}
0=-\gamma _{1}B(x_{2}\otimes \ 1_{H};G,x_{2})-2B(x_{1}\otimes
1_{H};1_{A},1_{H}).\left( \ref{X1,x2,X1F41,1H}\right)
\end{equation*}%
\begin{equation*}
0=0.\left( \ref{X1,x2,X1F41,x1x2}\right)
\end{equation*}%
\begin{gather*}
B(x_{1}\otimes 1_{H};1_{A},gx_{1})+B(gx_{1}x_{2}\otimes
1_{H};1_{A},gx_{1}x_{2})\left( \ref{X1,x2,X1F41,gx1}\right) \\
+B(g\otimes 1_{H};1_{A},g)+B(x_{2}\otimes \ 1_{H};1_{A},gx_{2})=0.
\end{gather*}

\begin{equation*}
2\beta _{1}B\left( g\otimes 1_{H};1_{A},gx_{1}x_{2}\right) =\gamma
_{1}B(g\otimes 1_{H};G,gx_{2})+2B(x_{1}\otimes 1_{H};1_{A},gx_{2}).\left( %
\ref{X1,x2,X1F41,gx2}\right)
\end{equation*}

\begin{eqnarray*}
0 &=&\gamma _{1}\left[ B(g\otimes 1_{H};G,gx_{1})+B(x_{2}\otimes \
1_{H};G,gx_{1}x_{2})\right] +\left( \ref{X1,x2,X1F51,g}\right) \\
&&+\left[
\begin{array}{c}
2B(g\otimes 1_{H};1_{A},g)+2B(x_{2}\otimes \ 1_{H};1_{A},gx_{2}) \\
+2B(x_{1}\otimes 1_{H};1_{A},gx_{1})+2B(gx_{1}x_{2}\otimes
1_{H};1_{A},gx_{1}x_{2})%
\end{array}%
\right]
\end{eqnarray*}

\begin{equation*}
0=0.\left( \ref{X1,x2,X1F51,x2}\right)
\end{equation*}%
\begin{eqnarray*}
&&\lambda \left[ B(g\otimes 1_{H};G,gx_{1})+B(x_{2}\otimes \
1_{H};G,gx_{1}x_{2})\right] \\
&=&-2\left[ B(x_{2}\otimes 1_{H};G,g)-B(gx_{1}x_{2}\otimes 1_{H};G,gx_{1})%
\right] \left( \ref{X1,x2,X1F61,g}\right)
\end{eqnarray*}%
\begin{equation*}
\beta _{1}\left[ B(g\otimes 1_{H};G,gx_{1})+B(x_{2}\otimes \
1_{H};G,gx_{1}x_{2})\right] =0\left( \ref{X1,x2,X1F71,g}\right)
\end{equation*}%
\begin{equation*}
0=0\left( \ref{X1,x2,X1F71,x1}\right)
\end{equation*}%
\begin{equation*}
\lambda B(x_{1}x_{2}\otimes 1_{H};X_{2},g)-\gamma _{1}B(x_{1}x_{2}\otimes
1_{H};G,g)=0\left( \ref{X1,x1x2,X1F11,g}\right)
\end{equation*}%
\begin{gather*}
2\beta _{1}B(x_{1}x_{2}\otimes 1_{H};X_{1},x_{1})+\lambda
B(x_{1}x_{2}\otimes 1_{H};X_{2},x_{1})+ \\
+\gamma _{1}B(x_{1}x_{2}\otimes 1_{H};G,x_{1})=0\left( \ref{X1,x1x2,X1F11,x1}%
\right)
\end{gather*}

\begin{gather*}
2\beta _{1}B(x_{1}x_{2}\otimes 1_{H};X_{1},x_{2})+\lambda
B(x_{1}x_{2}\otimes 1_{H};X_{2},x_{2})+ \\
+\gamma _{1}B(x_{1}x_{2}\otimes 1_{H};G,x_{2})=0.\left( \ref%
{X1,x1x2,X1F11,x2}\right)
\end{gather*}

\begin{equation*}
\gamma _{1}B(x_{1}x_{2}\otimes 1_{H};G,gx_{1}x_{2})-2B(x_{1}x_{2}\otimes
1_{H};1_{A},gx_{2})=0.\left( \ref{X1,x1x2,X1F11,gx1x2}\right)
\end{equation*}%
\begin{equation*}
0=0.\left( \ref{X1,x1x2,X1F21,1H}\right)
\end{equation*}%
\begin{equation*}
\lambda B(x_{1}x_{2}\otimes 1_{H};GX_{2},gx_{1})+2B(x_{1}x_{2}\otimes
1_{H};G,g)=0.\left( \ref{X1,x1x2,X1F21,gx1}\right)
\end{equation*}

\begin{equation*}
\lambda \left[
\begin{array}{c}
+1-B(x_{1}x_{2}\otimes 1_{H};1_{A},x_{1}x_{2})-B(x_{1}x_{2}\otimes
1_{H};X_{2},x_{1})+ \\
+B(x_{1}x_{2}\otimes 1_{H};X_{1},x_{2})%
\end{array}%
\right] =0\left( \ref{X1,x1x2,X1F31,1H}\right)
\end{equation*}%
\begin{equation*}
2B(x_{1}x_{2}\otimes 1_{H};1_{A},gx_{1})+2B(x_{1}x_{2}\otimes
1_{H};X_{1},g)=0\left( \ref{X1,x1x2,X1F31,gx1}\right)
\end{equation*}%
\begin{equation*}
2B(x_{1}x_{2}\otimes 1_{H};1_{A},gx_{2})+\gamma _{1}B(x_{1}x_{2}\otimes
1_{H};GX_{1},gx_{2})=0.\left( \ref{X1,x1x2,X1F31,gx2}\right)
\end{equation*}

\begin{equation*}
\left[
\begin{array}{c}
+1-B(x_{1}x_{2}\otimes 1_{H};1_{A},x_{1}x_{2})-B(x_{1}x_{2}\otimes
1_{H};X_{2},x_{1}) \\
+B(x_{1}x_{2}\otimes 1_{H};X_{1},x_{2})%
\end{array}%
\right] =0.\left( \ref{X1,x1x2,X1F41,1H}\right)
\end{equation*}

\begin{equation*}
\gamma _{1}B(x_{1}x_{2}\otimes 1_{H};GX_{2},gx_{1})+2B(x_{1}x_{2}\otimes
1_{H};X_{2},g)=0.\left( \ref{X1,x1x2,X1F41,gx1}\right)
\end{equation*}%
\begin{gather*}
2B(x_{1}x_{2}\otimes 1_{H};X_{2},g)+\left( \ref{X1,x1x2,X1F51,g}\right) \\
+\gamma _{1}\left[ -B(x_{1}x_{2}\otimes
1_{H};G,gx_{1}x_{2})+B(x_{1}x_{2}\otimes
1_{H};GX_{2},gx_{1})-B(x_{1}x_{2}\otimes 1_{H};GX_{1},gx_{2})\right] =0
\end{gather*}%
\begin{gather*}
2B(x_{1}x_{2}\otimes 1_{H};G,g)+ \\
+\lambda \left[
\begin{array}{c}
-B(x_{1}x_{2}\otimes 1_{H};G,gx_{1}x_{2}) \\
+B(x_{1}x_{2}\otimes 1_{H};GX_{2},gx_{1})-B(x_{1}x_{2}\otimes
1_{H};GX_{1},gx_{2})%
\end{array}%
\right] =0\left( \ref{X1,x1x2,X1F61,g}\right)
\end{gather*}%
\begin{equation*}
B(x_{1}x_{2}\otimes 1_{H};G,gx_{1}x_{2})+B(x_{1}x_{2}\otimes
1_{H};GX_{1},gx_{2})=0\left( \ref{X1,x1x2,X1F61,gx1x2}\right)
\end{equation*}%
\begin{equation*}
+B(x_{1}x_{2}\otimes 1_{H};G,gx_{1}x_{2})+B(x_{1}x_{2}\otimes
1_{H};GX_{1},gx_{2})=0.\left( \ref{X1,x1x2,X1F81,gx1}\right)
\end{equation*}%
\begin{gather*}
2\beta _{1}\left[ -1+B(x_{1}x_{2}\otimes
1_{H};1_{A},x_{1}x_{2})+B(x_{1}x_{2}\otimes 1_{H};X_{2},x_{1})\right]
+\lambda B(x_{1}x_{2}\otimes 1_{H};X_{2},x_{2})\left( \ref{X1,gx1,X1F11,1H}%
\right) \\
=-\gamma _{1}B(x_{1}x_{2}\otimes 1_{H};G,x_{2})+ \\
-2\left[ B(x_{1}x_{2}\otimes 1_{H};1_{A},gx_{2})+B(x_{1}x_{2}\otimes
1_{H};X_{2},g)\right] .
\end{gather*}

\begin{gather*}
-\gamma _{1}\left[ B(x_{1}x_{2}\otimes
1_{H};G,gx_{1}x_{2})-B(x_{1}x_{2}\otimes 1_{H};GX_{2},gx_{1})\right] \left( %
\ref{X1,gx1,X1F11,gx1}\right) \\
+2\left[ B(x_{1}x_{2}\otimes 1_{H};1_{A},gx_{2})+B(x_{1}x_{2}\otimes
1_{H};X_{2},g)\right] =0.
\end{gather*}%
\begin{gather*}
2\left[ B(x_{1}x_{2}\otimes 1_{H};1_{A},gx_{2})+B(x_{1}x_{2}\otimes
1_{H};X_{2},g)\right] +\left( \ref{X1,gx1,X1F31,g}\right) \\
+ \\
+\gamma _{1}\left[ -B(x_{1}x_{2}\otimes
1_{H};G,gx_{1}x_{2})+B(x_{1}x_{2}\otimes 1_{H};GX_{2},gx_{1})\right] =0.
\end{gather*}%
\begin{gather*}
-2\beta _{1}B(x_{1}x_{2}\otimes 1_{H};X_{1},x_{1})+ \\
+\lambda \left[ -1+B(x_{1}x_{2}\otimes
1_{H};1_{A},x_{1}x_{2})-B(x_{1}x_{2}\otimes 1_{H};X_{1},x_{2})\right] \left( %
\ref{X1,gx2,X1F11,1H}\right) \\
+\gamma _{1}\left[ -B(x_{1}x_{2}\otimes 1_{H};G,x_{1})\right] =0.
\end{gather*}%
\begin{equation*}
+2B(x_{1}x_{2}\otimes 1_{H};1_{A},gx_{1})+2B(x_{1}x_{2}\otimes
1_{H};X_{1},g)=0.\left( \ref{X1,gx2,X1F11,gx1}\right)
\end{equation*}%
\begin{equation*}
+\gamma _{1}\left[ B(x_{1}x_{2}\otimes
1_{H};G,gx_{1}x_{2})+B(x_{1}x_{2}\otimes 1_{H};GX_{1},gx_{2})\right]
=0\left( \ref{X1,gx2,X1F11,gx2}\right)
\end{equation*}%
\begin{gather*}
2\left[ B(x_{1}x_{2}\otimes 1_{H};1_{A},gx_{1})+B(x_{1}x_{2}\otimes
1_{H};X_{1},g)\right] \left( \ref{X1,gx2,X131,g}\right) \\
=0
\end{gather*}%
\begin{equation*}
\gamma _{1}\left[ -B(x_{1}x_{2}\otimes
1_{H};G,gx_{1}x_{2})+B(x_{1}x_{2}\otimes 1_{H};GX_{1},gx_{2})\right]
=0\left( \ref{X1,gx2,X141,g}\right)
\end{equation*}%
\begin{equation*}
\lambda B(x_{1}\otimes 1_{H};1_{A},1_{H})-\gamma _{1}B(gx_{1}x_{2}\otimes
1_{H};G,1_{H})=0\left( \ref{X1,gx1x2,X1F11,1H}\right)
\end{equation*}%
\begin{equation*}
0=0\left( \ref{X1,gx1x2,X1F11,x1x2}\right)
\end{equation*}%
\begin{gather*}
-2\beta _{1}B(x_{2}\otimes 1_{H};1_{A},gx_{1})+\lambda \left[ B(x_{1}\otimes
1_{H};1_{A},gx_{1})+B(gx_{1}x_{2}\otimes 1_{H};1_{A},gx_{1}x_{2})\right]
\left( \ref{X1,gx1x2,X1F11,gx1}\right) \\
+\gamma _{1}B(gx_{1}x_{2}\otimes 1_{H};G,gx_{1})+2B(gx_{1}x_{2}\otimes
1_{H};1_{A},g)=0.
\end{gather*}%
\begin{gather*}
2\beta _{1}\left[ -B(x_{2}\otimes 1_{H};1_{A},gx_{2})-B(gx_{1}x_{2}\otimes
1_{H};1_{A},gx_{1}x_{2})\right] \left( \ref{X1,gx1x2,X171,gx2}\right) \\
+\lambda B(x_{1}\otimes 1_{H};1_{A},gx_{2})+\gamma _{1}B(gx_{1}x_{2}\otimes
1_{H};G,gx_{2})=0.
\end{gather*}%
\begin{gather*}
2\beta _{1}\left[ B(x_{2}\otimes 1_{H};G,g)-B(gx_{1}x_{2}\otimes
1_{H};G,gx_{1})\right] \left( \ref{X1,gx1x2,XF21,g}\right) \\
+\lambda \left[ -B(x_{1}\otimes 1_{H};G,g)-B(gx_{1}x_{2}\otimes
1_{H};G,gx_{2})\right] =0
\end{gather*}

\begin{equation*}
2\beta _{1}B(x_{2}\otimes 1_{H};G,gx_{1}x_{2})-\lambda B(x_{1}\otimes
1_{H};G,gx_{1}x_{2})=0\left( \ref{X1,gx1x2,X1F21,gx1x2}\right)
\end{equation*}%
\begin{gather*}
+\lambda \left[ B(g\otimes 1_{H};1_{A},g)+B(x_{2}\otimes \
1_{H};1_{A},gx_{2})+B(x_{1}\otimes 1_{H};1_{A},gx_{1})+B(gx_{1}x_{2}\otimes
1_{H};1_{A},gx_{1}x_{2})\right] \\
-\gamma _{1}\left[ B(x_{2}\otimes 1_{H};G,g)-B(gx_{1}x_{2}\otimes
1_{H};G,gx_{1})\right] =0\left( \ref{X1,gx1x2,X1F31,g}\right)
\end{gather*}%
\begin{equation*}
B(x_{2}\otimes \ 1_{H};1_{A},x_{1}x_{2})=0\left( \ref{X1,gx1x2,X1F31,x1}%
\right)
\end{equation*}

\begin{equation*}
0=0\left( \ref{X1,gx1x2,X1F31,x2}\right)
\end{equation*}

\begin{gather*}
+\lambda B\left( g\otimes 1_{H};1_{A},gx_{1}x_{2}\right) +\left( \ref%
{X1,gx1x2,X1F31,gx1x2}\right) \\
-\gamma _{1}B(x_{2}\otimes 1_{H};G,gx_{1}x_{2})-\left[ +2B(x_{2}\otimes
1_{H};1_{A},gx_{2})+2B(gx_{1}x_{2}\otimes 1_{H};1_{A},gx_{1}x_{2})\right] =0
\end{gather*}

\begin{eqnarray*}
&&2\beta _{1}\left[
\begin{array}{c}
B(g\otimes 1_{H};1_{A},g)+B(x_{2}\otimes \ 1_{H};1_{A},gx_{2})+ \\
+B(x_{1}\otimes 1_{H};1_{A},gx_{1})+B(gx_{1}x_{2}\otimes
1_{H};1_{A},gx_{1}x_{2})%
\end{array}%
\right] \left( \ref{X1,gx1x2,X1F41,g}\right) \\
&&-\gamma _{1}\left[ B(x_{1}\otimes 1_{H};G,g)+B(gx_{1}x_{2}\otimes
1_{H};G,gx_{2})\right] =0
\end{eqnarray*}%
\begin{eqnarray*}
&&2\beta _{1}B\left( g\otimes 1_{H};1_{A},gx_{1}x_{2}\right) \left( \ref%
{X1,gx1x2,X1F41,gx1x2}\right) \\
&&-\gamma _{1}B(x_{1}\otimes 1_{H};G,gx_{1}x_{2})-2B(x_{1}\otimes
1_{H};1_{A},gx_{2})=0
\end{eqnarray*}%
\begin{equation*}
2B(x_{1}\otimes 1_{H};1_{A},1_{H})+\gamma _{1}B(x_{2}\otimes \
1_{H};G,x_{2})=0\left( \ref{X1,gx1x2,X1F51,1H}\right)
\end{equation*}%
\begin{gather*}
+\gamma _{1}\left[ B(g\otimes 1_{H};G,gx_{1})+B(x_{2}\otimes \
1_{H};G,gx_{1}x_{2})\right] +\left( \ref{X1,gx1x2,X1F51,gx1}\right) \\
-\left[
\begin{array}{c}
2B(g\otimes 1_{H};1_{A},g)+2B(x_{2}\otimes \ 1_{H};1_{A},gx_{2})+ \\
+2B(x_{1}\otimes 1_{H};1_{A},gx_{1})+2B(gx_{1}x_{2}\otimes
1_{H};1_{A},gx_{1}x_{2})%
\end{array}%
\right] =0
\end{gather*}%
\begin{equation*}
2B(gx_{1}x_{2}\otimes 1_{H};G,1_{H})+\lambda B(x_{2}\otimes \
1_{H};G,x_{2})=0\left( \ref{X1,gx1x2,X1F61,1H}\right)
\end{equation*}%
\begin{equation*}
0=0\left( \ref{X1,gx1x2,X1F61,x1x2}\right)
\end{equation*}%
\begin{gather*}
\lambda \left[ B(g\otimes 1_{H};G,gx_{1})+B(x_{2}\otimes \
1_{H};G,gx_{1}x_{2})\right] + \\
+\left[ 2B(x_{2}\otimes 1_{H};G,g)-2B(gx_{1}x_{2}\otimes 1_{H};G,gx_{1})%
\right] =0\left( \ref{X1,gx1x2,X1F61,gx1}\right)
\end{gather*}%
\begin{gather*}
\beta _{1}\left[ B(g\otimes 1_{H};G,gx_{1})+B(x_{2}\otimes \
1_{H};G,gx_{1}x_{2})\right] \left( \ref{X1,gx1x2,X1F71,gx1}\right) \\
+B(x_{1}\otimes 1_{H};G,g)+B(gx_{1}x_{2}\otimes 1_{H};G,gx_{2})=0.
\end{gather*}

\subsubsection{$X_{2}$}

\begin{equation*}
B(x_{2}\otimes 1_{H};1_{A},1_{H})=0\left( \ref{X2,g,X2F11,1H}\right)
\end{equation*}%
\begin{gather*}
2\beta _{2}B(g\otimes 1_{H};1_{A},gx_{1}x_{2})+\gamma _{2}B(g\otimes
1_{H};G,gx_{1})\left( \ref{X2,g,X2F11,gx1}\right) \\
+2B(x_{2}\otimes 1_{H};1_{A},gx_{1})=0
\end{gather*}

\begin{gather*}
\gamma _{2}B(g\otimes 1_{H};G,gx_{2})+\lambda B(g\otimes
1_{H};1_{A},gx_{1}x_{2})\left( \ref{X2,g,X2F11,gx2}\right) \\
+2B(x_{2}\otimes 1_{H};1_{A},gx_{2})+2B(g\otimes 1_{H};1_{A},g)=0.
\end{gather*}

\begin{equation*}
2\beta _{2}B(g\otimes 1_{H};G,gx_{2})+\lambda B(g\otimes
1_{H};G,gx_{1})+2B(x_{2}\otimes 1_{H};G,g)=0\left( \ref{X2,g,X2F21,g}\right)
\end{equation*}%
\begin{equation*}
0=0\left( \ref{X2,g,X2F21,x1}\right)
\end{equation*}%
\begin{equation*}
0=0\ref{X2,g,X2F21,x2}
\end{equation*}%
\begin{equation*}
2\beta _{2}B(g\otimes 1_{H};1_{H},gx_{1}x_{2})+\gamma _{2}B(g\otimes
1_{H};G,gx_{1})+2B(x_{2}\otimes 1_{H};1_{H},gx_{1})=0.\left( \ref%
{X2,g,X2F31,g}\right)
\end{equation*}%
\begin{equation*}
B(x_{2}\otimes 1_{H};1,x_{1}x_{2})=0\left( \ref{X2,g,X2F31,x2}\right)
\end{equation*}%
\begin{equation*}
B(x_{2}\otimes 1_{H};X_{2},x_{1})+B\left( x_{2}\otimes
1_{H};1_{A},x_{1}x_{2}\right) =0\left( \ref{X2,g,X2F41,x1}\right)
\end{equation*}%
\begin{equation*}
B\left( g\otimes 1_{H};1_{A},gx_{1}x_{2}\right) +B(g\otimes
1_{H},G,gx_{1})=0\left( \ref{X2,g,X2F61,gx2}\right)
\end{equation*}%
\begin{equation*}
B(g\otimes 1_{H};G,x_{1}x_{2})=0\left( \ref{X2,g,X2F71,1H}\right)
\end{equation*}%
\begin{equation*}
B(g\otimes 1_{H};G,gx_{1})+B(x_{2}\otimes 1_{H};G,gx_{1}x_{2})=0\left( \ref%
{X2,g,X2F71,gx1}\right)
\end{equation*}

\begin{eqnarray*}
&&-2\beta _{2}B(x_{1}\otimes 1_{H};1_{A},gx_{2})+\gamma _{2}B(x_{1}\otimes
1_{H};G,g)\left( \ref{X2,x1,X2F11,g}\right) \\
&&-\lambda \left[ B(g\otimes 1_{H};1_{A},g)+B(x_{1}\otimes
1_{H};1_{A},gx_{1})\right] -2B(gx_{1}x_{2}\otimes 1_{H};1_{A},g)=0
\end{eqnarray*}

\begin{equation*}
0=0\left( \ref{X2,x1,X2F11,x1}\right)
\end{equation*}%
\begin{equation*}
0=0.\left( \ref{X2,x1,X2F11,x2}\right)
\end{equation*}%
\begin{eqnarray*}
&&+\gamma _{2}B(x_{1}\otimes 1_{H};G,gx_{1}x_{2})-\lambda B(g\otimes
1_{H};1_{A},gx_{1}x_{2})\left( \ref{X2,x1,X2F11,gx1x2}\right) \\
&&-2B(gx_{1}x_{2}\otimes 1_{H};1_{A},gx_{1}x_{2})-2B\left( x_{1}\otimes
1_{H};1_{A},gx_{1}\right) =0
\end{eqnarray*}%
\begin{equation*}
B(gx_{1}x_{2}\otimes 1_{H};G,1_{H})=0\left( \ref{X2,x1,X2F21,1H}\right)
\end{equation*}

\begin{equation*}
2\beta _{2}B(x_{1}\otimes 1_{H};G,gx_{1}x_{2})+\lambda B\left( g\otimes
1_{H};G,gx_{1}\right) +2B(x_{2}\otimes 1_{H};G,g)=0.\left( \ref%
{X2,x1,X2F21,gx1}\right)
\end{equation*}

\begin{eqnarray*}
&&\lambda \left[ B\left( g\otimes 1_{H};G,gx_{2}\right) -B(x_{1}\otimes
1_{H};G,gx_{1}x_{2})\right] \left( \ref{X2,x1,X2F21,gx2}\right) \\
&&+2-B(x_{1}\otimes 1_{H};G,g_{1})+2B\left( x_{1}\otimes 1_{H};G,g\right) =0.
\end{eqnarray*}

\begin{equation*}
B(x_{2}\otimes 1_{H};1_{A},1_{H})=0.\left( \ref{X2,x1,X2F31,1H}\right)
\end{equation*}%
\begin{gather*}
\gamma _{2}\left[ -B(g\otimes 1_{H};G,gx_{2})+B(x_{1}\otimes
1_{H};G,gx_{1}x_{2})\right] +\left( \ref{X2,x1,X2F31,gx2}\right) \\
-2B(x_{2}\otimes 1_{H};1_{A},gx_{2})-2B(gx_{1}x_{2}\otimes
1_{H};1_{A},gx_{1}x_{2}) \\
-2B(g\otimes 1_{H};1_{A},g)-2B(x_{1}\otimes 1_{H};1_{A},gx_{1})=0.
\end{gather*}%
\begin{equation*}
0=0.\left( \ref{X2,x1,X2F41,1H}\right)
\end{equation*}

\begin{eqnarray*}
&&-\gamma _{2}B(x_{1}\otimes 1_{H};G,gx_{1}x_{2})+\lambda B\left( g\otimes
1_{H};1_{A},gx_{1}x_{2}\right) \left( \ref{X2,x1,X2F41,gx1}\right) \\
&&+2B(x_{1}\otimes 1_{H};1_{A},gx_{1})+2B(gx_{1}x_{2}\otimes
1_{H};1_{A},gx_{1}x_{2})=0.
\end{eqnarray*}

\begin{eqnarray*}
&&\beta _{2}\left[ -B(g\otimes 1_{H};G,gx_{2})+B(x_{1}\otimes
1_{H};G,gx_{1}x_{2})\right] \left( \ref{X2,x1,X2F61,g}\right) \\
&&-B(x_{2}\otimes 1_{H};G,g)+B(gx_{1}x_{2}\otimes 1_{H};G,gx_{1})=0
\end{eqnarray*}%
\begin{equation*}
0=0\left( \ref{X2,x1,X2F61,x2}\right)
\end{equation*}

\begin{equation*}
B(x_{1}\otimes 1_{H};G,g)+B(gx_{1}x_{2}\otimes 1_{H};G,gx_{2})=0\left( \ref%
{X2,x1,X2F71,g}\right)
\end{equation*}%
\begin{equation*}
B(x_{2}\otimes \ 1_{H};G,x_{2})=0\left( \ref{X2,x1,X2F81,1H}\right)
\end{equation*}%
\begin{eqnarray*}
2 &&\beta _{2}\left[ B(g\otimes 1_{H};1_{A},g)+B(x_{2}\otimes \
1_{H};1_{A},gx_{2})\right] \left( \ref{X2,x2,X2F11,g}\right) \\
&&-\gamma _{2}B(x_{2}\otimes 1_{H};G,g)+\lambda B(x_{2}\otimes
1_{H};1_{A},gx_{1})=0.
\end{eqnarray*}%
\begin{equation*}
0=0.\left( \ref{X2,x2,X2F11,x2}\right)
\end{equation*}

\begin{eqnarray*}
+ &&2\beta _{2}B\left( g\otimes 1_{H};1_{A},gx_{1}x_{2}\right) -\gamma
_{2}B(x_{2}\otimes 1_{H};G,gx_{1}x_{2})\left( \ref{X2,x2,X2F11,gx1x2}\right)
\\
&&+2B\left( x_{2}\otimes 1_{H};1_{A},gx_{1}\right) =0
\end{eqnarray*}

\begin{equation*}
\beta _{2}\left[ B(g\otimes 1_{H};G,gx_{1})+B(x_{2}\otimes \
1_{H};G,gx_{1}x_{2})\right] =0\left( \ref{X2,x2,X2F21,gx1}\right)
\end{equation*}%
\begin{eqnarray*}
&&2\beta _{2}B(g\otimes 1_{H};G,gx_{2})-\lambda B(x_{2}\otimes
1_{H};G,gx_{1}x_{2})\left( \ref{X2,x2,X2F21,gx2}\right) \\
&&+2B\left( x_{2}\otimes 1_{H};G,g\right) =0
\end{eqnarray*}

\begin{equation*}
2B(x_{2}\otimes 1_{H};1_{A},1_{H})+\gamma _{2}B(x_{2}\otimes \
1_{H};G,x_{2})=0.\left( \ref{X2,x2,X2F41,1H}\right)
\end{equation*}

\begin{eqnarray*}
&&\gamma _{2}B(g\otimes 1_{H};G,gx_{2})-\lambda B\left( g\otimes
1_{H};1_{A},gx_{1}x_{2}\right) \left( \ref{X2,x2,X2F41,gx2}\right) \\
&&+2\left[ +B(g\otimes 1_{H};1_{A},g)+B(x_{2}\otimes \ 1_{H};1_{A},gx_{2})%
\right] =0.
\end{eqnarray*}

\begin{equation*}
\lambda \left[ B(g\otimes 1_{H};G,gx_{1})+B(x_{2}\otimes \
1_{H};G,gx_{1}x_{2})\right] =0\left( \ref{X2,x2,X2F71,g}\right)
\end{equation*}%
\begin{equation*}
\gamma _{2}B(x_{1}x_{2}\otimes 1_{H};G,g)+\lambda B(x_{1}x_{2}\otimes
1_{H};X_{1},g)=0.\left( \ref{X2,x1x2,X2F11,g}\right)
\end{equation*}

\begin{gather*}
2\beta _{2}B(x_{1}x_{2}\otimes 1_{H};X_{2},x_{1})+\gamma
_{2}B(x_{1}x_{2}\otimes 1_{H};G,x_{1})\left( \ref{X2,x1x2,X2F11,x1}\right) \\
+\lambda B(x_{1}x_{2}\otimes 1_{H};X_{1},x_{1})=0.
\end{gather*}

\begin{gather*}
2\beta _{2}B(x_{1}x_{2}\otimes 1_{H};X_{2},x_{2})+\gamma
_{2}B(x_{1}x_{2}\otimes 1_{H};G,x_{2})\left( \ref{X2,x1x2,X2F11,x2}\right) \\
+\lambda B(x_{1}x_{2}\otimes 1_{H};X_{1},x_{2})=0.
\end{gather*}%
\begin{eqnarray*}
&&\gamma _{2}B(x_{1}x_{2}\otimes 1_{H};G,gx_{1}x_{2})+\left( \ref%
{X2,x1x2,X2F11,gx1x2}\right) \\
&&+2B(x_{1}x_{2}\otimes 1_{H};1_{A},gx_{1})=0.
\end{eqnarray*}

\begin{equation*}
0=0.\left( \ref{X2,x1x2,X2F21,1H}\right)
\end{equation*}

\begin{equation*}
2B(x_{1}x_{2}\otimes 1_{H};G,g)-\lambda B(x_{1}x_{2}\otimes
1_{H};GX_{1},gx_{2})=0.\left( \ref{X2,x1x2,X2F21,gx2}\right)
\end{equation*}%
\begin{eqnarray*}
&&2\beta _{2}\left[
\begin{array}{c}
+1-B(x_{1}x_{2}\otimes 1_{H};1_{A},x_{1}x_{2}) \\
-B(x_{1}x_{2}\otimes 1_{H};X_{2},x_{1})+B(x_{1}x_{2}\otimes
1_{H};X_{1},x_{2})%
\end{array}%
\right] \left( \ref{X2,x1x2,X2F31,1H}\right) \\
&=&0.
\end{eqnarray*}%
\begin{equation*}
\gamma _{2}B(x_{1}x_{2}\otimes 1_{H};GX_{1},gx_{2})+2B(x_{1}x_{2}\otimes
1_{H};X_{1},g)=0.\left( \ref{X2,x1x2,X2F31,gx2}\right)
\end{equation*}

\begin{equation*}
\left[
\begin{array}{c}
+1-B(x_{1}x_{2}\otimes 1_{H};1_{A},x_{1}x_{2}) \\
-B(x_{1}x_{2}\otimes 1_{H};X_{2},x_{1})+B(x_{1}x_{2}\otimes
1_{H};X_{1},x_{2})%
\end{array}%
\right] =0.\left( \ref{X2,x1x2,X2F41,1H}\right)
\end{equation*}%
\begin{gather*}
2B(x_{1}x_{2}\otimes 1_{H};1_{A},gx_{2})+\left( \ref{X2,x1x2,X2F41,gx2}%
\right) \\
+2B(x_{1}x_{2}\otimes 1_{H};X_{2},g)=0.
\end{gather*}

\begin{gather*}
2B(x_{1}x_{2}\otimes 1_{H};X_{1},g)-\left( \ref{X2,x1x2,X2F51,g}\right) \\
\gamma _{2}\left[
\begin{array}{c}
-B(x_{1}x_{2}\otimes 1_{H};G,gx_{1}x_{2})+ \\
B(x_{1}x_{2}\otimes 1_{H};GX_{2},gx_{1})-B(x_{1}x_{2}\otimes
1_{H};GX_{1},gx_{2})%
\end{array}%
\right] =0.
\end{gather*}%
\begin{gather*}
2B(x_{1}x_{2}\otimes 1_{H};G,g)\left( \ref{X2,x1x2,X2F71,g}\right) \\
+\lambda \left[
\begin{array}{c}
-B(x_{1}x_{2}\otimes 1_{H};G,gx_{1}x_{2})+ \\
B(x_{1}x_{2}\otimes 1_{H};GX_{2},gx_{1})-B(x_{1}x_{2}\otimes
1_{H};GX_{1},gx_{2})%
\end{array}%
\right] =0.
\end{gather*}%
\begin{equation*}
B(x_{1}x_{2}\otimes 1_{H};G,gx_{1}x_{2})-B(x_{1}x_{2}\otimes
1_{H};GX_{2},gx_{1})=0.\left( \ref{X2,x1x2,X2F71,gx1x2}\right)
\end{equation*}%
\begin{gather*}
2\beta _{2}B(x_{1}x_{2}\otimes 1_{H};X_{2},x_{2})+\gamma
_{2}B(x_{1}x_{2}\otimes 1_{H};G,x_{2})\left( \ref{X2,gx1,X2F11,1H}\right) \\
+\lambda \left[ -1+B(x_{1}x_{2}\otimes
1_{H};1_{A},x_{1}x_{2})+B(x_{1}x_{2}\otimes 1_{H};X_{2},x_{1})\right] =0.
\end{gather*}

\begin{equation*}
B(x_{1}x_{2}\otimes 1_{H};1_{A},gx_{2})+B(x_{1}x_{2}\otimes
1_{H};X_{2},g)=0.\left( \ref{X2,gx1,X2F11,gx2}\right)
\end{equation*}%
\begin{equation*}
\lambda \left[ B(x_{1}x_{2}\otimes 1_{H};G,gx_{1}x_{2})-B(x_{1}x_{2}\otimes
1_{H};GX_{2},gx_{1})\right] =0.\left( \ref{X2,gx1,X2F21,g}\right)
\end{equation*}%
\begin{gather*}
2\beta _{2}\left[ -1+B(x_{1}x_{2}\otimes
1_{H};1_{A},x_{1}x_{2})-B(x_{1}x_{2}\otimes 1_{H};X_{1},x_{2})\right] \left( %
\ref{X2,gx2,X2F11,1H}\right) \\
+\gamma _{2}\left[ -B(x_{1}x_{2}\otimes 1_{H};G,x_{1})\right] + \\
-\lambda B(x_{1}x_{2}\otimes 1_{H};X_{1},x_{1})=0.
\end{gather*}

\begin{eqnarray*}
&&\gamma _{2}\left[ B(x_{1}x_{2}\otimes
1_{H};G,gx_{1}x_{2})+B(x_{1}x_{2}\otimes 1_{H};GX_{1},gx_{2})\right] \left( %
\ref{X2,gx2,X2F11,gx2}\right) \\
&&+ \\
&&+\left[ +2B(x_{1}x_{2}\otimes 1_{H};1_{A},gx_{1})+2B(x_{1}x_{2}\otimes
1_{H};X_{1},g)\right] =0.
\end{eqnarray*}%
\begin{eqnarray*}
&&2\left[ -B(x_{1}x_{2}\otimes 1_{H};1_{A},gx_{1})-B(x_{1}x_{2}\otimes
1_{H};X_{1},g)\right] \left( \ref{X2,gx2,X2F41,g}\right) \\
&&+\gamma _{2}\left[ -B(x_{1}x_{2}\otimes
1_{H};G,gx_{1}x_{2})-B(x_{1}x_{2}\otimes 1_{H};GX_{1},gx_{2})\right] \\
&=&0.
\end{eqnarray*}%
\begin{equation*}
-\gamma _{2}B(gx_{1}x_{2}\otimes 1_{H};G,1_{H})+\lambda B(x_{2}\otimes
1_{H};1_{A},1_{H})=0\left( \ref{X2,gx1x2,X2F11,1H}\right)
\end{equation*}%
\begin{equation*}
0=0\left( \ref{X2,gx1x2,X2F11,x1x2}\right)
\end{equation*}%
\begin{gather*}
2\beta _{2}\left[ B(x_{1}\otimes 1_{H};1_{A},gx_{1})+B(gx_{1}x_{2}\otimes
1_{H};1_{A},gx_{1}x_{2})\right] \left( \ref{X2,gx1x2,X2F11,gx1}\right) \\
+\gamma _{2}B(gx_{1}x_{2}\otimes 1_{H};G,gx_{1})-\lambda B(x_{2}\otimes
1_{H};1_{A},gx_{1})=0.
\end{gather*}

\begin{gather*}
2\beta _{2}B(x_{1}\otimes 1_{H};1_{A},gx_{2})+\gamma
_{2}B(gx_{1}x_{2}\otimes 1_{H};G,gx_{2})\left( \ref{X2,gx1x2,X2F11,gx2}%
\right) \\
+\lambda \left[ -B(x_{2}\otimes 1_{H};1_{A},gx_{2})-B(gx_{1}x_{2}\otimes
1_{H};1_{A},gx_{1}x_{2})\right] \\
+2B(gx_{1}x_{2}\otimes 1_{H};1_{A},g)=0.
\end{gather*}

\begin{gather*}
2\beta _{2}\left[ -B(x_{1}\otimes 1_{H};G,g)-B(gx_{1}x_{2}\otimes
1_{H};G,gx_{2})\right] \left( \ref{X2,gx1x2,X2F21,g}\right) \\
+\lambda \left[ B(x_{2}\otimes 1_{H};G,g)-B(gx_{1}x_{2}\otimes
1_{H};G,gx_{1})\right] =0.
\end{gather*}%
\begin{gather*}
-2\beta _{2}B(x_{1}\otimes 1_{H};G,gx_{1}x_{2})\left( \ref%
{X2,gx1x2,X2F21,gx1x2}\right) \\
+\lambda B(x_{2}\otimes 1_{H};G,gx_{1}x_{2})-2B(gx_{1}x_{2}\otimes
1_{H};G,gx_{1})=0.
\end{gather*}%
\begin{gather*}
2\beta _{2}\left[
\begin{array}{c}
B(g\otimes 1_{H};1_{A},g)+B(x_{2}\otimes \ 1_{H};1_{A},gx_{2}) \\
+B(x_{1}\otimes 1_{H};1_{A},gx_{1})+B(gx_{1}x_{2}\otimes
1_{H};1_{A},gx_{1}x_{2})%
\end{array}%
\right] \left( \ref{X2,gx1x2,X2F31,g}\right) \\
-\gamma _{2}\left[ B(x_{2}\otimes 1_{H};G,g)-B(gx_{1}x_{2}\otimes
1_{H};G,gx_{1})\right] =0.
\end{gather*}

\begin{eqnarray*}
&&\gamma _{2}\left[ B(x_{1}\otimes 1_{H};G,g)+B(gx_{1}x_{2}\otimes
1_{H};G,gx_{2})\right] \left( \ref{X2,gx1x2,X2F41,g}\right) \\
&&-\lambda \left[
\begin{array}{c}
B(g\otimes 1_{H};1_{A},g)+B(x_{2}\otimes \ 1_{H};1_{A},gx_{2}) \\
+B(x_{1}\otimes 1_{H};1_{A},gx_{1})+B(gx_{1}x_{2}\otimes
1_{H};1_{A},gx_{1}x_{2})%
\end{array}%
\right] =0.
\end{eqnarray*}

\begin{equation*}
0=0.\left( \ref{X2,gx1x2,X2F41,x1}\right)
\end{equation*}

\begin{equation*}
-2B(x_{2}\otimes 1_{H};1_{A},1_{H})-\gamma _{2}B(x_{2}\otimes \
1_{H};G,x_{2})=0.\left( \ref{X2,gx1x2,X2F51,1H}\right)
\end{equation*}

\begin{equation*}
B(x_{2}\otimes 1_{H};G,g)-B(gx_{1}x_{2}\otimes 1_{H};G,gx_{1})=0.\left( \ref%
{X2,gx1x2,X2F61,gx2}\right)
\end{equation*}%
\begin{equation*}
2B(gx_{1}x_{2}\otimes 1_{H};G,1_{H})+\lambda B(x_{2}\otimes \
1_{H};G,x_{2})=0.\left( \ref{X2,gx1x2,X2F71,1H}\right)
\end{equation*}

\subsection{LIST\ OF\ MONOMIAL\ EQUALITIES 3\label{LME3}}

We list the new monoidal equalities we obtained above taking out constants%
\begin{equation*}
B(g\otimes 1_{H};G,x_{1}x_{2})=0\left( \ref{X2,g,X2F71,1H}\right)
\end{equation*}%
\begin{equation*}
B(x_{1}\otimes 1_{H};1_{A},1_{H})=0.\left( \ref{G,x1, GF2,1H}\right)
\end{equation*}%
\begin{equation*}
B(x_{1}\otimes 1_{H};1_{A},x_{1}x_{2})=0\text{ }\left( \ref{G,x1, GF2,x1x2}%
\right)
\end{equation*}%
\begin{equation*}
B(x_{2}\otimes 1_{H};1_{A},1_{H})=0\left( \ref{X2,g,X2F11,1H}\right)
\end{equation*}%
\begin{equation*}
B(x_{2}\otimes 1_{H};1_{A},1_{H})=0.\left( \ref{X2,x1,X2F31,1H}\right)
\end{equation*}%
\begin{equation*}
B(x_{2}\otimes \ 1_{H};1_{A},x_{1}x_{2})=0\left( \ref{X1,gx1x2,X1F31,x1}%
\right)
\end{equation*}%
\begin{equation*}
B(x_{2}\otimes 1_{H};1,x_{1}x_{2})=0\left( \ref{X2,g,X2F31,x2}\right)
\end{equation*}%
\begin{equation*}
B(x_{2}\otimes \ 1_{H};G,x_{2})=0\left( \ref{X2,x1,X2F81,1H}\right)
\end{equation*}%
\begin{equation*}
B(gx_{1}x_{2}\otimes 1_{H};G,1_{H})=0\left( \ref{X2,x1,X2F21,1H}\right)
\end{equation*}%
We take one of the repetition as before and relabel them

\begin{equation*}
B(g\otimes 1_{H};G,x_{1}x_{2})=0\left( \ref{X2,g,X2F71,1H}\right) \text{
this is already in \ref{LME0}}
\end{equation*}%
\begin{equation}
B(x_{1}\otimes 1_{H};1_{A},1_{H})=0.\left( \ref{G,x1, GF2,1H}\right)
\label{x1ot1,1,1}
\end{equation}%
\begin{equation*}
B(x_{1}\otimes 1_{H};1_{A},x_{1}x_{2})=0\text{ }\left( \ref{G,x1, GF2,x1x2}%
\right) \text{ this is already in \ref{LME1}}
\end{equation*}%
\begin{equation}
B(x_{2}\otimes 1_{H};1_{A},1_{H})=0\left( \ref{X2,g,X2F11,1H}\right)
\label{x2ot1,1,1}
\end{equation}%
\begin{equation*}
B(x_{2}\otimes \ 1_{H};1_{A},x_{1}x_{2})=0\left( \ref{X1,gx1x2,X1F31,x1}%
\right) \text{ this was already in }\left( \ref{LME1}\right)
\end{equation*}%
\begin{equation*}
B(x_{2}\otimes \ 1_{H};G,x_{2})=0\left( \ref{X2,x1,X2F81,1H}\right) \text{
this was already in \ref{LME2}}
\end{equation*}%
\begin{equation}
B(gx_{1}x_{2}\otimes 1_{H};G,1_{H})=0\left( \ref{X2,x1,X2F21,1H}\right)
\label{gx1x2ot1,G,1}
\end{equation}

Now we take the complete list of equalities $\left( \ref{LAE3}\right) $ and
cancel there the terms above.

\subsection{LIST\ OF\ ALL\ EQUALITIES 4\label{LAE4}}

\subsubsection{$G$}

\begin{equation*}
\gamma _{1}B\left( g\otimes 1_{H};1_{A},x_{1}\right) +\gamma _{2}B\left(
g\otimes 1_{H};1_{A},x_{2}\right) =0\text{ }\left( \ref{G,g,GF1,1H}\right)
\end{equation*}

\begin{equation*}
2\alpha B(g\otimes 1_{H};G,gx_{1})+\gamma _{2}B(g\otimes
1_{H};1_{A},gx_{1}x_{2})=0\left( \ref{G,g, GF1,gx1}\right)
\end{equation*}

\begin{equation*}
2\alpha B(g\otimes 1_{H};G,gx_{2})-\gamma _{1}B\left( g\otimes
1_{H};1_{A},gx_{1}x_{2}\right) =0.\left( \ref{G,g, GF1,gx2}\right)
\end{equation*}

\begin{equation*}
\gamma _{1}B\left( g\otimes 1_{H};G,gx_{1}\right) +\gamma _{2}B\left(
g\otimes 1_{H};G,gx_{2}\right) =0.\left( \ref{G,g, GF2,g}\right)
\end{equation*}%
\begin{equation*}
\begin{array}{c}
2\alpha B(x_{1}\otimes 1_{H};G,g)+ \\
+\gamma _{1}\left[ -B(g\otimes 1_{H};1_{A},g)-B(x_{1}\otimes
1_{H};1_{A},gx_{1})\right] -\gamma _{2}B(x_{1}\otimes 1_{H};1_{A},gx_{2})%
\end{array}%
=0.\left( \ref{G,x1, GF1,g}\right)
\end{equation*}

\begin{equation*}
-\gamma _{1}B(g\otimes 1_{H};1_{A},gx_{1}x_{2})+2\alpha B(x_{1}\otimes
1_{H};G,gx_{1}x_{2})=0.\left( \ref{G,x1, GF1,gx1x2}\right)
\end{equation*}%
\begin{equation*}
0=0.\left( \ref{G,x1, GF2,1H}\right)
\end{equation*}%
\begin{equation*}
0=0\text{ }\left( \ref{G,x1, GF2,x1x2}\right)
\end{equation*}%
\begin{equation*}
\gamma _{1}B\left( g\otimes 1_{H};G,gx_{1}\right) +\gamma _{2}B(x_{1}\otimes
1_{H};G,gx_{1}x_{2})=0.\left( \ref{G,x1, GF2,gx1}\right)
\end{equation*}

\begin{equation*}
\gamma _{1}[B\left( g\otimes 1_{H};G,gx_{2}\right) -B(x_{1}\otimes
1_{H};G,gx_{1}x_{2})]=0.\left( \ref{G,x1,GF2,gx2}\right)
\end{equation*}

\begin{equation*}
0=0\left( \ref{G,x1,GF3,1H}\right)
\end{equation*}%
\begin{equation*}
-\gamma _{2}B(g\otimes 1_{H};1_{A},gx_{1}x_{2})-2\alpha B(g\otimes
1_{H};G,gx_{1})=0.\left( \ref{G,x1,GF3,gx1}\right)
\end{equation*}

\begin{equation*}
\alpha \left[ -B(g\otimes 1_{H};g,gx_{2})+B(x_{1}\otimes 1_{H};g,gx_{1}x_{2})%
\right] =0.\left( \ref{G,x1, GF3,x2}\right)
\end{equation*}%
\begin{equation*}
0=0.\left( \ref{G,x1, GF4,1H}\right)
\end{equation*}

\begin{equation*}
-\gamma _{1}B(g\otimes 1_{H};1_{A},gx_{1}x_{2})+2\alpha B(x_{1}\otimes
1_{H};G,gx_{1}x_{2})=0.\left( \ref{G,x1, GF4,gx1}\right)
\end{equation*}

\begin{equation*}
\alpha \lbrack -B(g\otimes 1_{H};g,gx_{2})+B(x_{1}\otimes
1_{H};G,gx_{1}x_{2})]=0.\left( \ref{G,x1, GF5,g}\right)
\end{equation*}%
\begin{equation*}
0=0.\left( \ref{G,x1,GF6,x1}\right)
\end{equation*}%
\begin{equation*}
0=0.\left( \ref{G,x1,GF6,x2}\right)
\end{equation*}%
\begin{equation*}
\begin{array}{c}
2\alpha B(x_{2}\otimes 1_{H};G,g)+ \\
-\gamma _{1}B(x_{2}\otimes 1_{H};1_{A},gx_{1})+\gamma _{2}\left[ -B(g\otimes
1_{H};1_{A},g)-B(x_{2}\otimes \ 1_{H};1_{A},gx_{2})\right]%
\end{array}%
=0.\left( \ref{G,x2, GF1,g}\right)
\end{equation*}%
\begin{equation*}
0=0\left( \ref{G,x2, GF1,x1}\right)
\end{equation*}

\begin{equation*}
2\alpha B(x_{2}\otimes 1_{H};G,gx_{1}x_{2})-\gamma _{2}B\left( g\otimes
1_{H};1_{A},gx_{1}x_{2}\right) =0.\left( \ref{G,x2, GF1,gx1x2}\right)
\end{equation*}%
\begin{equation*}
0=0.\left( \ref{G,x2, GF2,1H}\right)
\end{equation*}%
\begin{equation*}
0=0\left( \ref{G,x2, GF2,x1x2}\right)
\end{equation*}%
\begin{equation*}
\gamma _{2}\left[ B(g\otimes 1_{H};G,gx_{1})+B(x_{2}\otimes \
1_{H};G,gx_{1}x_{2})\right] =0.\left( \ref{G,x2, GF2,gx1}\right)
\end{equation*}

\begin{equation*}
-\gamma _{1}B(x_{2}\otimes 1_{H};G,gx_{1}x_{2})+\gamma _{2}B(g\otimes
1_{H};G,gx_{2})=0.\left( \ref{G,x2, GF2,gx2}\right)
\end{equation*}%
\begin{equation*}
0=0\left( \ref{G,x2, FG3,1H}\right)
\end{equation*}

\begin{equation*}
-\gamma _{2}B(g\otimes 1_{H};1_{A},gx_{1}x_{2})-2\alpha B(x_{2}\otimes
1_{H};G,gx_{1}x_{2})=0.\left( \ref{G,x2, GF3,gx2}\right)
\end{equation*}%
\begin{equation*}
0=0.\left( \ref{G,x2, GF4,1H}\right)
\end{equation*}%
\begin{equation*}
\alpha \left[ B(g\otimes 1_{H};G,gx_{1})+B(x_{2}\otimes \
1_{H};G,gx_{1}x_{2})\right] =0\left( \ref{G,x2, GF4,gx1}\right)
\end{equation*}

\begin{equation*}
\gamma _{1}B(g\otimes 1_{H};1_{A},gx_{1}x_{2})-2\alpha B(g\otimes
1_{H};G,gx_{2})=0.\left( \ref{G,x2, GF4,gx2}\right)
\end{equation*}%
\begin{equation*}
\alpha \left[ B(g\otimes 1_{H};G,gx_{1})+B(x_{2}\otimes \
1_{H};G,gx_{1}x_{2})\right] =0\left( \ref{G,x2, GF5,g}\right)
\end{equation*}

\begin{equation*}
\gamma _{2}[B(g\otimes 1_{H},G,gx_{1})+B(x_{2}\otimes
1_{H};G,gx_{1}x_{2})=0.=0.\left( \ref{G,x2, GF6,g}\right)
\end{equation*}%
\begin{equation*}
0=0.\left( \ref{G,x2; GF6,x2}\right)
\end{equation*}

\begin{equation*}
\gamma _{1}[B(g\otimes 1_{H};G,gx_{1})+B(x_{2}\otimes
1_{H};G,gx_{1}x_{2})]=0.\left( \ref{G,x2; GF7,g}\right)
\end{equation*}%
\begin{equation*}
0=0.\left( \ref{G,x2; GF7,x1}\right)
\end{equation*}%
\begin{equation*}
\gamma _{1}B(x_{1}x_{2}\otimes 1_{H};X_{1},g)+\gamma _{2}B(x_{1}x_{2}\otimes
1_{H};X_{2},g)=0.\left( \ref{G,x1x2, GF1,g}\right)
\end{equation*}

\begin{equation*}
\begin{array}{c}
2\alpha B(x_{1}x_{2}\otimes 1_{H};G,x_{1})+ \\
+\gamma _{1}B(x_{1}x_{2}\otimes 1_{H};X_{1},x_{1})+\gamma
_{2}B(x_{1}x_{2}\otimes 1_{H};X_{2},x_{1})%
\end{array}%
=0.\left( \ref{G,x1x2, GF1,x1}\right)
\end{equation*}%
\begin{equation*}
\begin{array}{c}
2\alpha B(x_{1}x_{2}\otimes 1_{H};G,x_{2})+ \\
+\gamma _{1}B(x_{1}x_{2}\otimes 1_{H};X_{1},x_{2})+\gamma
_{2}B(x_{1}x_{2}\otimes 1_{H};X_{2},x_{2})%
\end{array}%
=0.\left( \ref{G,x1x2, GF1,x2}\right)
\end{equation*}

\begin{equation*}
0=0\left( \ref{G,x1x2, GF2,1H}\right)
\end{equation*}%
\begin{equation*}
\begin{array}{c}
2B(x_{1}x_{2}\otimes 1_{H};1_{A},gx_{1}) \\
+\gamma _{2}B(x_{1}x_{2}\otimes 1_{H};GX_{2},gx_{1})%
\end{array}%
=0.\left( \ref{G,x1x2, GF2,gx1}\right)
\end{equation*}%
\begin{equation*}
\begin{array}{c}
2B(x_{1}x_{2}\otimes 1_{H};1_{A},gx_{2}) \\
+\gamma _{1}B(x_{1}x_{2}\otimes 1_{H};GX_{1},gx_{2})%
\end{array}%
=0\left( \ref{G,x1x2, GF2,gx2}\right)
\end{equation*}

\begin{equation*}
0=0\text{ }\left( \ref{G,x1x2, GF3,gx1}\right)
\end{equation*}%
\begin{gather*}
-\gamma _{1}[1-B(x_{1}x_{2}\otimes
1_{H};1_{A},x_{1}x_{2})+B(x_{1}x_{2}\otimes 1_{H};x_{1},x_{2})+ \\
-B(x_{1}x_{2}\otimes 1_{H};x_{2},x_{1})]-2\alpha B(x_{1}x_{2}\otimes
1_{H};GX_{2},1_{H}=0.\left( \ref{G,x1x2, GF4,1H}\right)
\end{gather*}

\begin{gather*}
-2B(x_{1}x_{2}\otimes 1_{H};X_{1},g)+\gamma _{2}[-B(x_{1}x_{2}\otimes
1_{H};g,gx_{1}x_{2})+ \\
-B(x_{1}x_{2}\otimes 1_{H};GX_{1},gx_{2})+B(x_{1}x_{2}\otimes
1_{H};GX_{2},gx_{1})]=0\left( \ref{G,x1x2, GF6,g}\right)
\end{gather*}%
\begin{gather*}
-2B(x_{1}x_{2}\otimes 1_{H};X_{2},g)+ \\
-\gamma _{1}[-B(x_{1}x_{2}\otimes 1_{H};g,gx_{1}x_{2})+ \\
-B(x_{1}x_{2}\otimes 1_{H};GX_{1},gx_{2})+B(x_{1}x_{2}\otimes
1_{H};GX_{2},gx_{1})]=0\left( \ref{G,x1x2, GF7,g}\right)
\end{gather*}

\begin{equation*}
\begin{array}{c}
2\alpha B(x_{1}x_{2}\otimes 1_{H};G,x_{2})+ \\
+\gamma _{1}\left[
\begin{array}{c}
-1+B(x_{1}x_{2}\otimes 1_{H};1_{A},x_{1}x_{2})+ \\
+B(x_{1}x_{2}\otimes 1_{H};X_{2},x_{1})%
\end{array}%
\right] +\gamma _{2}B(x_{1}x_{2}\otimes 1_{H};X_{2},x_{2})%
\end{array}%
=0.\left( \ref{G,gx1, GF1,1H}\right)
\end{equation*}

\begin{equation*}
\begin{array}{c}
2\left[ B(x_{1}x_{2}\otimes 1_{H};1_{A},gx_{2})+B(x_{1}x_{2}\otimes
1_{H};X_{2},g)\right] \\
\gamma _{1}\left[ -B(x_{1}x_{2}\otimes
1_{H};G,gx_{1}x_{2})+B(x_{1}x_{2}\otimes 1_{H};GX_{2},gx_{1})\right]%
\end{array}%
\left( \ref{G,gx1, GF2,g}\right)
\end{equation*}%
\begin{equation*}
\begin{array}{c}
2\alpha \left[ -B(x_{1}x_{2}\otimes 1_{H};G,x_{1})\right] + \\
-\gamma _{1}B(x_{1}x_{2}\otimes 1_{H};X_{1},x_{1})+\gamma _{2}\left[
-1+B(x_{1}x_{2}\otimes 1_{H};1_{A},x_{1}x_{2})-B(x_{1}x_{2}\otimes
1_{H};X_{1},x_{2})\right]%
\end{array}%
=0\left( \ref{G,gx2, GF1,1H}\right)
\end{equation*}%
\begin{equation*}
\begin{array}{c}
2\left[ -B(x_{1}x_{2}\otimes 1_{H};1_{A},gx_{1})-B(x_{1}x_{2}\otimes
1_{H};X_{1},g)\right] \\
+\gamma _{2}\left[ -B(x_{1}x_{2}\otimes
1_{H};G,gx_{1}x_{2})-B(x_{1}x_{2}\otimes 1_{H};GX_{1},gx_{2})\right]%
\end{array}%
=0.\left( \ref{G,gx2, GF2,g}\right)
\end{equation*}

\begin{equation*}
0=0\left( \ref{G,gx1x2, GF1,1H}\right)
\end{equation*}%
\begin{equation*}
0=0.\left( \ref{G,gx1x2, GF1,x1x2}\right)
\end{equation*}%
\begin{equation*}
\begin{array}{c}
2\alpha B(x_{2}\otimes 1_{H};G,g)+ \\
-\gamma _{1}B(x_{2}\otimes 1_{H};1_{A},gx_{1})+\gamma _{2}\left[
B(x_{1}\otimes 1_{H};1_{A},gx_{1})+B(gx_{1}x_{2}\otimes
1_{H};1_{A},gx_{1}x_{2})\right]%
\end{array}%
=0.\left( \ref{G,gx1x2, GF1,gx1}\right)
\end{equation*}

\begin{equation*}
\begin{array}{c}
2\alpha B(gx_{1}x_{2}\otimes 1_{H};G,gx_{2})+ \\
+\gamma _{1}\left[ -B(x_{2}\otimes 1_{H};1_{A},gx_{2})-B(gx_{1}x_{2}\otimes
1_{H};1_{A},gx_{1}x_{2})\right] +\gamma _{2}B(x_{1}\otimes
1_{H};1_{A},gx_{2})%
\end{array}%
=0.\left( \ref{G,gx1x2, GF1,gx2}\right)
\end{equation*}%
\begin{gather*}
\gamma _{1}\left[ B(x_{2}\otimes 1_{H};G,g)-B(gx_{1}x_{2}\otimes
1_{H};G,gx_{1})\right] \left( \ref{G,gx1x2, GF2,g}\right) \\
+\gamma _{2}\left[ -B(x_{1}\otimes 1_{H};G,g)-B(gx_{1}x_{2}\otimes
1_{H};G,gx_{2})\right] =0.
\end{gather*}%
\begin{equation*}
0=0.\left( \ref{G,gx1x2, GF2,x1}\right)
\end{equation*}%
\begin{equation*}
0=0.\left( \ref{G,gx1x2, GF2,x2}\right)
\end{equation*}%
\begin{equation*}
\gamma _{1}B(x_{2}\otimes 1_{H};G,gx_{1}x_{2})-\gamma _{2}B(x_{1}\otimes
1_{H};G,gx_{1}x_{2})=0.\left( \ref{G,gx1x2, GF2,gx1x2}\right)
\end{equation*}%
\begin{gather*}
\gamma _{2}\left[
\begin{array}{c}
B(g\otimes 1_{H};1_{A},g)+B(x_{2}\otimes \ 1_{H};1_{A},gx_{2})+ \\
+B(x_{1}\otimes 1_{H};1_{A},gx_{1})+B(gx_{1}x_{2}\otimes
1_{H};1_{A},gx_{1}x_{2})%
\end{array}%
\right] +\left( \ref{G,gx1x2, GF3,g}\right) \\
+2\alpha \left[ -B(x_{2}\otimes 1_{H};G,g)+B(x_{2}\otimes 1_{H};G,g)\right]
=0.
\end{gather*}

\begin{equation*}
0=0.\left( \ref{G,gx1x2, GF3,x1}\right)
\end{equation*}

\begin{equation*}
0=0\left( \ref{G,gx1x2, GF3,x2}\right)
\end{equation*}

\begin{equation*}
\gamma _{2}B\left( g\otimes 1_{H};1_{A},gx_{1}x_{2}\right) -2\alpha
B(x_{2}\otimes 1_{H};G,gx_{1}x_{2})=0.\left( \ref{G,gx1x2, GF3,gx1x2}\right)
\end{equation*}

\begin{gather*}
-\gamma _{1}\left[
\begin{array}{c}
B(g\otimes 1_{H};1_{A},g)+B(x_{2}\otimes \ 1_{H};1_{A},gx_{2})+ \\
+B(x_{1}\otimes 1_{H};1_{A},gx_{1})+B(gx_{1}x_{2}\otimes
1_{H};1_{A},gx_{1}x_{2})%
\end{array}%
\right] +\left( \ref{G,gx1x2, GF4,g}\right) \\
-2\alpha \left[ -B(x_{1}\otimes 1_{H};G,g)-B(gx_{1}x_{2}\otimes
1_{H};G,gx_{2})\right] =0
\end{gather*}%
\begin{equation*}
0=0\left( \ref{G,gx1x2, GF4,x1}\right)
\end{equation*}%
\begin{equation*}
0=0.\left( \ref{G,gx1x2, GF4,x2}\right)
\end{equation*}

\begin{equation*}
-\gamma _{1}B\left( g\otimes 1_{H};1_{A},gx_{1}x_{2}\right) +2\alpha
B(x_{1}\otimes 1_{H};G,gx_{1}x_{2})=0.\left( \ref{G,gx1x2, GF4,gx1x2}\right)
\end{equation*}

\begin{equation*}
\alpha \left[ B(g\otimes 1_{H};G,gx_{2})-B(x_{1}\otimes 1_{H};G,gx_{1}x_{2})%
\right] =0.\left( \ref{G,gx1x2, GF5,gx2}\right)
\end{equation*}

\begin{equation*}
0=0\left( \ref{G,gx1x2, GF6,1H}\right)
\end{equation*}%
\begin{equation*}
\gamma _{2}\left[ B(g\otimes 1_{H};G,gx_{2})-B(x_{1}\otimes
1_{H};G,gx_{1}x_{2})\right] =0.\left( \ref{G,gx1x2, GF6,gx1}\right)
\end{equation*}%
\begin{equation*}
\gamma _{2}\left[ B(g\otimes 1_{H};G,gx_{2})-B(x_{1}\otimes
1_{H};G,gx_{1}x_{2})\right] =0.\left( \ref{G,gx1x2, GF6,gx2}\right)
\end{equation*}%
\begin{equation*}
0=0\left( \ref{G,gx1x2, GF7,1H}\right)
\end{equation*}%
\begin{equation*}
0=0.\left( \ref{G,gx1x2, GF7,x1x2}\right)
\end{equation*}%
\begin{equation*}
\gamma _{1}\left[ B(g\otimes 1_{H};G,gx_{1})+B(x_{2}\otimes \
1_{H};G,gx_{1}x_{2})\right] =0.\left( \ref{G,gx1x2, GF7,gx1}\right)
\end{equation*}%
\begin{equation*}
\gamma _{1}\left[ B(g\otimes 1_{H};G,gx_{2})-B(x_{1}\otimes
1_{H};G,gx_{1}x_{2})\right] =0.\left( \ref{G,gx1x2, GF7,gx2}\right)
\end{equation*}

\subsubsection{$X_{1}$}

\begin{equation*}
B\left( g\otimes 1_{H};1_{A},x_{2}\right) =0.\left( \ref{X1,g, X1F11,1H}%
\right)
\end{equation*}%
\begin{eqnarray*}
&&+\gamma _{1}B(g\otimes 1_{H};G,gx_{1})+\lambda B(g\otimes
1_{H};X_{2},gx_{1})+B(x_{1}\otimes 1_{H};1_{A},gx_{1})+\left( \ref{X1,g,
X1F11,gx1}\right) \\
&&+2B\left( g\otimes 1_{H};1_{A},g\right) +B\left( x_{1}\otimes
1_{H};1_{A},gx_{1}\right) =0
\end{eqnarray*}

\begin{equation*}
\gamma _{1}B(g\otimes 1_{H};G,gx_{2})+2B(x_{1}\otimes
1_{H};1_{A},gx_{2})-2\beta _{1}B\left( g\otimes
1_{H};1_{A},gx_{1}x_{2}\right) =0.\left( \ref{X1,g, X1F11,gx2}\right)
\end{equation*}

\begin{equation*}
2\beta _{1}B\left( g\otimes 1_{H};G,gx_{1}\right) +\lambda B\left( g\otimes
1_{H};G,gx_{2}\right) +2B(x_{1}\otimes 1_{H};G,g)=0.\left( \ref{X1,g,X1F21,g}%
\right)
\end{equation*}%
\begin{equation*}
0=0.\left( \ref{X1,g,X1F21,x1}\right)
\end{equation*}

\begin{equation*}
B\left( g\otimes 1_{H};G,gx_{2}\right) -B\left( x_{1}\otimes
1_{H};G,gx_{1}x_{2}\right) =0\left( \ref{X1,g,X1F21,gx1x2}\right)
\end{equation*}%
\begin{gather*}
+\lambda B\left( g\otimes 1_{H};1_{A},gx_{1}x_{2}\right) +\gamma _{1}B\left(
g\otimes 1_{H};G,gx_{1}\right) +\left( \ref{X1,g, X1F31,g}\right) \\
+2B(g\otimes 1_{H};1_{A},g)+2B(x_{1}\otimes 1_{H};1_{A},gx_{1})=0.
\end{gather*}

\begin{equation*}
-2\beta _{1}B\left( g\otimes 1_{H};1_{A},gx_{1}x_{2}\right) +2B(x_{1}\otimes
1_{H};1_{A},gx_{2})+\gamma _{1}B\left( g\otimes 1_{H};G,gx_{2}\right)
=0.\left( \ref{X1,g, X1F41,g}\right)
\end{equation*}%
\begin{equation*}
0=0.\left( \ref{X1,g, X1F41,x1}\right)
\end{equation*}

\begin{equation*}
0=0\left( \ref{X1,g,X1F51,1H}\right)
\end{equation*}%
\begin{equation*}
B(x_{1}\otimes 1_{H};G,gx_{1}x_{2})-B(g\otimes 1_{H};G,gx_{2})=0.\left( \ref%
{X1,x2,X1F71,gx1x2}\right)
\end{equation*}%
\begin{equation*}
2\beta _{1}\left[ -B(g\otimes 1_{H};1_{A},g)-B(x_{1}\otimes
1_{H};1_{A},gx_{1})\right] -\lambda B(x_{1}\otimes
1_{H};1_{A},gx_{2})=\gamma _{1}B(x_{1}\otimes 1_{H};G,g)\left( \ref%
{X1,x1,X1F11,g}\right)
\end{equation*}%
\begin{equation*}
0=0.\left( \ref{X1,x1,X1F11,x1}\right)
\end{equation*}

\begin{equation*}
-2\beta _{1}B(g\otimes 1_{H};1_{A},gx_{1}x_{2})=\gamma _{1}B(x_{1}\otimes
1_{H};G,gx_{1}x_{2})-2B\left( x_{1}\otimes 1_{H};1_{A},gx_{2}\right) .\left( %
\ref{X1,x1,X1F11,gx1x2}\right)
\end{equation*}%
\begin{equation*}
0=0\left( \ref{X1,x1,X1F21,1}\right)
\end{equation*}%
\begin{equation*}
\lambda B(x_{1}\otimes 1_{H};G,gx_{1}x_{2})=-2B\left( x_{1}\otimes
1_{H};G,g\right) .\left( \ref{X1,x1,X1F21,gx1}\right)
\end{equation*}%
\begin{equation*}
B\left( g\otimes 1_{H};1_{A},x_{2}\right) =0.\left( \ref{X1,x1,X1F31,1H}%
\right)
\end{equation*}%
\begin{equation*}
0=0\left( \ref{X1,x1,X1F31,x1x2}\right)
\end{equation*}

\begin{gather*}
-\lambda B\left( g\otimes 1_{H};1_{A},gx_{1}x_{2}\right) -\gamma _{1}B\left(
g\otimes 1_{H};G,gx_{1}\right) + \\
-2\left[ B(g\otimes 1_{H};1_{A},g)+B(x_{1}\otimes 1_{H};1_{A},gx_{1})\right]
=0\left( \ref{X1,x1,X1F31,gx1}\right)
\end{gather*}

\begin{equation*}
\gamma _{1}\left[ -B\left( g\otimes 1_{H};G,gx_{2}\right) +B(x_{1}\otimes
1_{H};G,gx_{1}x_{2})\right] =0\left( \ref{X1,x1,X1F31,gx2}\right)
\end{equation*}

\begin{equation*}
2\beta _{1}B\left( g\otimes 1_{H};1_{A},gx_{1}x_{2}\right) -\gamma
_{1}B(x_{1}\otimes 1_{H};G,gx_{1}x_{2})-B(x_{1}\otimes
x_{1};X_{2},gx_{1})=0.\left( \ref{X1,x1,X1F41,gx1}\right)
\end{equation*}

\begin{eqnarray*}
&&-2\beta _{1}B(x_{2}\otimes 1_{H};1_{A},gx_{1})+\lambda \left[ -B(g\otimes
1_{H};1_{A},g)-B(x_{2}\otimes \ 1_{H};1_{A},gx_{2})\right] \\
&=&-\gamma _{1}B(x_{2}\otimes 1_{H};G,g)-2B(gx_{1}x_{2}\otimes
1_{H};1_{A},g).\left( \ref{X1,x2,X1F11,g}\right)
\end{eqnarray*}%
\begin{equation*}
0=0.\left( \ref{X1,x2,X1F11,x1}\right)
\end{equation*}%
\begin{equation*}
0=0.\left( \ref{X1,x2,X1F11,x2}\right)
\end{equation*}

\begin{gather*}
\lambda B\left( g\otimes 1_{H};1_{A},gx_{1}x_{2}\right) -\gamma
_{1}B(x_{2}\otimes 1_{H};G,gx_{1}x_{2})\left( \ref{X1,x2,X1F11,gx1x2}\right)
\\
-2B\left( x_{2}\otimes 1_{H};1_{A},gx_{2}\right) -2B(gx_{1}x_{2}\otimes
1_{H};1_{A},gx_{1}x_{2})=0.
\end{gather*}%
\begin{equation*}
B(x_{2}\otimes \ 1_{H};G,x_{2})=0.\left( \ref{X1,x2,X1F21,1H}\right)
\end{equation*}%
\begin{equation*}
0=0.\left( \ref{X1,x2,X1F21,x1x2}\right)
\end{equation*}%
\begin{eqnarray*}
&&\lambda \left[ B(g\otimes 1_{H};G,gx_{1})+B(x_{2}\otimes \
1_{H};G,gx_{1}x_{2})\right] \\
&=&2B(gx_{1}x_{2}\otimes 1_{H};G,gx_{1})-2B\left( x_{2}\otimes
1_{H};G,g\right) .\left( \ref{X1,x2,X1F21,gx1}\right)
\end{eqnarray*}

\begin{equation*}
-2\beta _{1}B(x_{2}\otimes 1_{H};G,gx_{1}x_{2})+\lambda B(g\otimes
1_{H};G,gx_{2})-2B(gx_{1}x_{2}\otimes 1_{H};G,gx_{2})=0.\left( \ref%
{X1,x2,X1F21,gx2}\right)
\end{equation*}%
\begin{gather*}
-\lambda B\left( g\otimes 1_{H};1_{A},gx_{1}x_{2}\right) =-\gamma
_{1}B(x_{2}\otimes 1_{H};G,gx_{1}x_{2})+ \\
-2B(x_{2}\otimes 1_{H};1_{A},gx_{2})-2B(gx_{1}x_{2}\otimes
1_{H};1_{A},gx_{1}x_{2}).\left( \ref{X1,x2,X1F31,gx2}\right)
\end{gather*}

\begin{equation*}
0=0.\left( \ref{X1,x2,X1F41,1H}\right)
\end{equation*}%
\begin{equation*}
0=0.\left( \ref{X1,x2,X1F41,x1x2}\right)
\end{equation*}%
\begin{gather*}
B(x_{1}\otimes 1_{H};1_{A},gx_{1})+B(gx_{1}x_{2}\otimes
1_{H};1_{A},gx_{1}x_{2})\left( \ref{X1,x2,X1F41,gx1}\right) \\
+B(g\otimes 1_{H};1_{A},g)+B(x_{2}\otimes \ 1_{H};1_{A},gx_{2})=0.
\end{gather*}

\begin{equation*}
2\beta _{1}B\left( g\otimes 1_{H};1_{A},gx_{1}x_{2}\right) =\gamma
_{1}B(g\otimes 1_{H};G,gx_{2})+2B(x_{1}\otimes 1_{H};1_{A},gx_{2}).\left( %
\ref{X1,x2,X1F41,gx2}\right)
\end{equation*}

\begin{eqnarray*}
0 &=&\gamma _{1}\left[ B(g\otimes 1_{H};G,gx_{1})+B(x_{2}\otimes \
1_{H};G,gx_{1}x_{2})\right] +\left( \ref{X1,x2,X1F51,g}\right) \\
&&+\left[
\begin{array}{c}
2B(g\otimes 1_{H};1_{A},g)+2B(x_{2}\otimes \ 1_{H};1_{A},gx_{2})+ \\
+2B(x_{1}\otimes 1_{H};1_{A},gx_{1})+2B(gx_{1}x_{2}\otimes
1_{H};1_{A},gx_{1}x_{2})%
\end{array}%
\right]
\end{eqnarray*}

\begin{equation*}
0=0.\left( \ref{X1,x2,X1F51,x2}\right)
\end{equation*}%
\begin{eqnarray*}
&&\lambda \left[ B(g\otimes 1_{H};G,gx_{1})+B(x_{2}\otimes \
1_{H};G,gx_{1}x_{2})\right] \\
&=&-2\left[ B(x_{2}\otimes 1_{H};G,g)-B(gx_{1}x_{2}\otimes 1_{H};G,gx_{1})%
\right] \left( \ref{X1,x2,X1F61,g}\right)
\end{eqnarray*}%
\begin{equation*}
\beta _{1}\left[ B(g\otimes 1_{H};G,gx_{1})+B(x_{2}\otimes \
1_{H};G,gx_{1}x_{2})\right] =0\left( \ref{X1,x2,X1F71,g}\right)
\end{equation*}%
\begin{equation*}
0=0\left( \ref{X1,x2,X1F71,x1}\right)
\end{equation*}%
\begin{equation*}
\lambda B(x_{1}x_{2}\otimes 1_{H};X_{2},g)-\gamma _{1}B(x_{1}x_{2}\otimes
1_{H};G,g)=0\left( \ref{X1,x1x2,X1F11,g}\right)
\end{equation*}%
\begin{equation*}
2\beta _{1}B(x_{1}x_{2}\otimes 1_{H};X_{1},x_{1})+\lambda
B(x_{1}x_{2}\otimes 1_{H};X_{2},x_{1})+\gamma _{1}B(x_{1}x_{2}\otimes
1_{H};G,x_{1})=0\left( \ref{X1,x1x2,X1F11,x1}\right)
\end{equation*}

\begin{equation*}
2\beta _{1}B(x_{1}x_{2}\otimes 1_{H};X_{1},x_{2})+\lambda
B(x_{1}x_{2}\otimes 1_{H};X_{2},x_{2})+\gamma _{1}B(x_{1}x_{2}\otimes
1_{H};G,x_{2})=0.\left( \ref{X1,x1x2,X1F11,x2}\right)
\end{equation*}

\begin{equation*}
\gamma _{1}B(x_{1}x_{2}\otimes 1_{H};G,gx_{1}x_{2})-2B(x_{1}x_{2}\otimes
1_{H};1_{A},gx_{2})=0.\left( \ref{X1,x1x2,X1F11,gx1x2}\right)
\end{equation*}%
\begin{equation*}
0=0.\left( \ref{X1,x1x2,X1F21,1H}\right)
\end{equation*}%
\begin{equation*}
\lambda B(x_{1}x_{2}\otimes 1_{H};GX_{2},gx_{1})+2B(x_{1}x_{2}\otimes
1_{H};G,g)=0.\left( \ref{X1,x1x2,X1F21,gx1}\right)
\end{equation*}

\begin{gather*}
\lambda \left[ +1-B(x_{1}x_{2}\otimes
1_{H};1_{A},x_{1}x_{2})-B(x_{1}x_{2}\otimes
1_{H};X_{2},x_{1})+B(x_{1}x_{2}\otimes 1_{H};X_{1},x_{2})\right] \left( \ref%
{X1,x1x2,X1F31,1H}\right) \\
=0
\end{gather*}%
\begin{equation*}
2B(x_{1}x_{2}\otimes 1_{H};1_{A},gx_{1})+2B(x_{1}x_{2}\otimes
1_{H};X_{1},g)=0\left( \ref{X1,x1x2,X1F31,gx1}\right)
\end{equation*}%
\begin{equation*}
2B(x_{1}x_{2}\otimes 1_{H};1_{A},gx_{2})+\gamma _{1}B(x_{1}x_{2}\otimes
1_{H};GX_{1},gx_{2})=0.\left( \ref{X1,x1x2,X1F31,gx2}\right)
\end{equation*}

\begin{equation*}
\left[
\begin{array}{c}
+1-B(x_{1}x_{2}\otimes 1_{H};1_{A},x_{1}x_{2})-B(x_{1}x_{2}\otimes
1_{H};X_{2},x_{1}) \\
+B(x_{1}x_{2}\otimes 1_{H};X_{1},x_{2})%
\end{array}%
\right] =0.\left( \ref{X1,x1x2,X1F41,1H}\right)
\end{equation*}

\begin{equation*}
\gamma _{1}B(x_{1}x_{2}\otimes 1_{H};GX_{2},gx_{1})+2B(x_{1}x_{2}\otimes
1_{H};X_{2},g)=0.\left( \ref{X1,x1x2,X1F41,gx1}\right)
\end{equation*}%
\begin{gather*}
2B(x_{1}x_{2}\otimes 1_{H};X_{2},g)+\left( \ref{X1,x1x2,X1F51,g}\right) \\
+\gamma _{1}\left[ -B(x_{1}x_{2}\otimes
1_{H};G,gx_{1}x_{2})+B(x_{1}x_{2}\otimes
1_{H};GX_{2},gx_{1})-B(x_{1}x_{2}\otimes 1_{H};GX_{1},gx_{2})\right] =0
\end{gather*}%
\begin{equation*}
2B(x_{1}x_{2}\otimes 1_{H};G,g)+\lambda \left[
\begin{array}{c}
-B(x_{1}x_{2}\otimes 1_{H};G,gx_{1}x_{2}) \\
+B(x_{1}x_{2}\otimes 1_{H};GX_{2},gx_{1})-B(x_{1}x_{2}\otimes
1_{H};GX_{1},gx_{2})%
\end{array}%
\right] =0\left( \ref{X1,x1x2,X1F61,g}\right)
\end{equation*}%
\begin{equation*}
B(x_{1}x_{2}\otimes 1_{H};G,gx_{1}x_{2})+B(x_{1}x_{2}\otimes
1_{H};GX_{1},gx_{2})=0\left( \ref{X1,x1x2,X1F61,gx1x2}\right)
\end{equation*}%
\begin{equation*}
+B(x_{1}x_{2}\otimes 1_{H};G,gx_{1}x_{2})+B(x_{1}x_{2}\otimes
1_{H};GX_{1},gx_{2})=0.\left( \ref{X1,x1x2,X1F81,gx1}\right)
\end{equation*}%
\begin{gather*}
2\beta _{1}\left[ -1+B(x_{1}x_{2}\otimes
1_{H};1_{A},x_{1}x_{2})+B(x_{1}x_{2}\otimes 1_{H};X_{2},x_{1})\right]
+\lambda B(x_{1}x_{2}\otimes 1_{H};X_{2},x_{2})\left( \ref{X1,gx1,X1F11,1H}%
\right) \\
=-\gamma _{1}B(x_{1}x_{2}\otimes 1_{H};G,x_{2})+ \\
-2\left[ B(x_{1}x_{2}\otimes 1_{H};1_{A},gx_{2})+B(x_{1}x_{2}\otimes
1_{H};X_{2},g)\right] .
\end{gather*}

\begin{gather*}
-\gamma _{1}\left[ B(x_{1}x_{2}\otimes
1_{H};G,gx_{1}x_{2})-B(x_{1}x_{2}\otimes 1_{H};GX_{2},gx_{1})\right] \left( %
\ref{X1,gx1,X1F11,gx1}\right) \\
+2\left[ B(x_{1}x_{2}\otimes 1_{H};1_{A},gx_{2})+B(x_{1}x_{2}\otimes
1_{H};X_{2},g)\right] =0.
\end{gather*}%
\begin{gather*}
2\left[ B(x_{1}x_{2}\otimes 1_{H};1_{A},gx_{2})+B(x_{1}x_{2}\otimes
1_{H};X_{2},g)\right] +\left( \ref{X1,gx1,X1F31,g}\right) \\
+\gamma _{1}\left[ -B(x_{1}x_{2}\otimes
1_{H};G,gx_{1}x_{2})+B(x_{1}x_{2}\otimes 1_{H};GX_{2},gx_{1})\right] =0.
\end{gather*}%
\begin{gather*}
-2\beta _{1}B(x_{1}x_{2}\otimes 1_{H};X_{1},x_{1})+\lambda \left[
-1+B(x_{1}x_{2}\otimes 1_{H};1_{A},x_{1}x_{2})-B(x_{1}x_{2}\otimes
1_{H};X_{1},x_{2})\right] \left( \ref{X1,gx2,X1F11,1H}\right) \\
+\gamma _{1}\left[ -B(x_{1}x_{2}\otimes 1_{H};G,x_{1})\right] =0.
\end{gather*}%
\begin{gather*}
\left( \ref{X1,gx2,X1F11,gx1}\right) \\
+2B(x_{1}x_{2}\otimes 1_{H};1_{A},gx_{1})+2B(x_{1}x_{2}\otimes
1_{H};X_{1},g)=0.
\end{gather*}%
\begin{equation*}
+\gamma _{1}\left[ B(x_{1}x_{2}\otimes
1_{H};G,gx_{1}x_{2})+B(x_{1}x_{2}\otimes 1_{H};GX_{1},gx_{2})\right]
=0\left( \ref{X1,gx2,X1F11,gx2}\right)
\end{equation*}%
\begin{gather*}
2\left[ B(x_{1}x_{2}\otimes 1_{H};1_{A},gx_{1})+B(x_{1}x_{2}\otimes
1_{H};X_{1},g)\right] \left( \ref{X1,gx2,X131,g}\right) \\
=0
\end{gather*}%
\begin{equation*}
\gamma _{1}\left[ -B(x_{1}x_{2}\otimes
1_{H};G,gx_{1}x_{2})+B(x_{1}x_{2}\otimes 1_{H};GX_{1},gx_{2})\right]
=0\left( \ref{X1,gx2,X141,g}\right)
\end{equation*}%
\begin{equation*}
B(x_{1}\otimes 1_{H};1_{A},1_{H})=0\left( \ref{X1,gx1x2,X1F11,1H}\right)
\end{equation*}%
\begin{equation*}
0=0\left( \ref{X1,gx1x2,X1F11,x1x2}\right)
\end{equation*}%
\begin{gather*}
-2\beta _{1}B(x_{2}\otimes 1_{H};1_{A},gx_{1})+\lambda \left[ B(x_{1}\otimes
1_{H};1_{A},gx_{1})+B(gx_{1}x_{2}\otimes 1_{H};1_{A},gx_{1}x_{2})\right]
\left( \ref{X1,gx1x2,X1F11,gx1}\right) \\
+\gamma _{1}B(gx_{1}x_{2}\otimes 1_{H};G,gx_{1})+2B(gx_{1}x_{2}\otimes
1_{H};1_{A},g)=0.
\end{gather*}%
\begin{gather*}
2\beta _{1}\left[ -B(x_{2}\otimes 1_{H};1_{A},gx_{2})-B(gx_{1}x_{2}\otimes
1_{H};1_{A},gx_{1}x_{2})\right] \left( \ref{X1,gx1x2,X171,gx2}\right) \\
+\lambda B(x_{1}\otimes 1_{H};1_{A},gx_{2})+\gamma _{1}B(gx_{1}x_{2}\otimes
1_{H};G,gx_{2})=0.
\end{gather*}%
\begin{gather*}
2\beta _{1}\left[ B(x_{2}\otimes 1_{H};G,g)-B(gx_{1}x_{2}\otimes
1_{H};G,gx_{1})\right] \left( \ref{X1,gx1x2,XF21,g}\right) \\
+\lambda \left[ -B(x_{1}\otimes 1_{H};G,g)-B(gx_{1}x_{2}\otimes
1_{H};G,gx_{2})\right] =0
\end{gather*}

\begin{equation*}
2\beta _{1}B(x_{2}\otimes 1_{H};G,gx_{1}x_{2})-\lambda B(x_{1}\otimes
1_{H};G,gx_{1}x_{2})=0\left( \ref{X1,gx1x2,X1F21,gx1x2}\right)
\end{equation*}%
\begin{gather*}
\lambda \left[
\begin{array}{c}
B(g\otimes 1_{H};1_{A},g)+B(x_{2}\otimes \ 1_{H};1_{A},gx_{2})+ \\
+B(x_{1}\otimes 1_{H};1_{A},gx_{1})+B(gx_{1}x_{2}\otimes
1_{H};1_{A},gx_{1}x_{2})%
\end{array}%
\right] + \\
-\gamma _{1}\left[ B(x_{2}\otimes 1_{H};G,g)-B(gx_{1}x_{2}\otimes
1_{H};G,gx_{1})\right] =0\left( \ref{X1,gx1x2,X1F31,g}\right)
\end{gather*}%
\begin{equation*}
0=0\left( \ref{X1,gx1x2,X1F31,x1}\right)
\end{equation*}

\begin{equation*}
0=0\left( \ref{X1,gx1x2,X1F31,x2}\right)
\end{equation*}

\begin{gather*}
+\lambda B\left( g\otimes 1_{H};1_{A},gx_{1}x_{2}\right) +\left( \ref%
{X1,gx1x2,X1F31,gx1x2}\right) \\
-\gamma _{1}B(x_{2}\otimes 1_{H};G,gx_{1}x_{2})-\left[ +2B(x_{2}\otimes
1_{H};1_{A},gx_{2})+2B(gx_{1}x_{2}\otimes 1_{H};1_{A},gx_{1}x_{2})\right] =0
\end{gather*}

\begin{eqnarray*}
&&2\beta _{1}\left[ B(g\otimes 1_{H};1_{A},g)+B(x_{2}\otimes \
1_{H};1_{A},gx_{2})+B(x_{1}\otimes 1_{H};1_{A},gx_{1})+B(gx_{1}x_{2}\otimes
1_{H};1_{A},gx_{1}x_{2})\right] \left( \ref{X1,gx1x2,X1F41,g}\right) \\
&&-\gamma _{1}\left[ B(x_{1}\otimes 1_{H};G,g)+B(gx_{1}x_{2}\otimes
1_{H};G,gx_{2})\right] =0
\end{eqnarray*}%
\begin{eqnarray*}
&&2\beta _{1}B\left( g\otimes 1_{H};1_{A},gx_{1}x_{2}\right) \left( \ref%
{X1,gx1x2,X1F41,gx1x2}\right) \\
&&-\gamma _{1}B(x_{1}\otimes 1_{H};G,gx_{1}x_{2})-2B(x_{1}\otimes
1_{H};1_{A},gx_{2})=0
\end{eqnarray*}%
\begin{equation*}
B(x_{1}\otimes 1_{H};1_{A},1_{H})=0\left( \ref{X1,gx1x2,X1F51,1H}\right)
\end{equation*}%
\begin{gather*}
+\gamma _{1}\left[ B(g\otimes 1_{H};G,gx_{1})+B(x_{2}\otimes \
1_{H};G,gx_{1}x_{2})\right] +\left( \ref{X1,gx1x2,X1F51,gx1}\right) \\
-\left[
\begin{array}{c}
2B(g\otimes 1_{H};1_{A},g)+2B(x_{2}\otimes \ 1_{H};1_{A},gx_{2})+ \\
+2B(x_{1}\otimes 1_{H};1_{A},gx_{1})+2B(gx_{1}x_{2}\otimes
1_{H};1_{A},gx_{1}x_{2})%
\end{array}%
\right] =0
\end{gather*}%
\begin{equation*}
0=0\left( \ref{X1,gx1x2,X1F61,1H}\right)
\end{equation*}%
\begin{equation*}
0=0\left( \ref{X1,gx1x2,X1F61,x1x2}\right)
\end{equation*}%
\begin{gather*}
\lambda \left[ B(g\otimes 1_{H};G,gx_{1})+B(x_{2}\otimes \
1_{H};G,gx_{1}x_{2})\right] + \\
+\left[ 2B(x_{2}\otimes 1_{H};G,g)-2B(gx_{1}x_{2}\otimes 1_{H};G,gx_{1})%
\right] =0\left( \ref{X1,gx1x2,X1F61,gx1}\right)
\end{gather*}%
\begin{gather*}
\beta _{1}\left[ B(g\otimes 1_{H};G,gx_{1})+B(x_{2}\otimes \
1_{H};G,gx_{1}x_{2})\right] \left( \ref{X1,gx1x2,X1F71,gx1}\right) \\
+B(x_{1}\otimes 1_{H};G,g)+B(gx_{1}x_{2}\otimes 1_{H};G,gx_{2})=0.
\end{gather*}

\subsubsection{$X_{2}$}

\begin{equation*}
0=0\left( \ref{X2,g,X2F11,1H}\right)
\end{equation*}%
\begin{gather*}
2\beta _{2}B(g\otimes 1_{H};1_{A},gx_{1}x_{2})+\gamma _{2}B(g\otimes
1_{H};G,gx_{1})\left( \ref{X2,g,X2F11,gx1}\right) \\
+2B(x_{2}\otimes 1_{H};1_{A},gx_{1})=0
\end{gather*}

\begin{gather*}
\gamma _{2}B(g\otimes 1_{H};G,gx_{2})+\lambda B(g\otimes
1_{H};1_{A},gx_{1}x_{2})\left( \ref{X2,g,X2F11,gx2}\right) \\
+2B(x_{2}\otimes 1_{H};1_{A},gx_{2})+2B(g\otimes 1_{H};1_{A},g)=0.
\end{gather*}

\begin{equation*}
2\beta _{2}B(g\otimes 1_{H};G,gx_{2})+\lambda B(g\otimes
1_{H};G,gx_{1})+2B(x_{2}\otimes 1_{H};G,g)=0\left( \ref{X2,g,X2F21,g}\right)
\end{equation*}%
\begin{equation*}
0=0\left( \ref{X2,g,X2F21,x1}\right)
\end{equation*}%
\begin{equation*}
0=0\ref{X2,g,X2F21,x2}
\end{equation*}%
\begin{equation*}
2\beta _{2}B(g\otimes 1_{H};1_{H},gx_{1}x_{2})+\gamma _{2}B(g\otimes
1_{H};G,gx_{1})+2B(x_{2}\otimes 1_{H};1_{H},gx_{1})=0.\left( \ref%
{X2,g,X2F31,g}\right)
\end{equation*}%
\begin{equation*}
B(x_{2}\otimes 1_{H};1,x_{1}x_{2})=0\left( \ref{X2,g,X2F31,x2}\right)
\end{equation*}%
\begin{equation*}
B(x_{2}\otimes 1_{H};X_{2},x_{1})+B\left( x_{2}\otimes
1_{H};1_{A},x_{1}x_{2}\right) =0\left( \ref{X2,g,X2F41,x1}\right)
\end{equation*}%
\begin{equation*}
B\left( g\otimes 1_{H};1_{A},gx_{1}x_{2}\right) +B(g\otimes
1_{H},G,gx_{1})=0\left( \ref{X2,g,X2F61,gx2}\right)
\end{equation*}%
\begin{equation*}
0=0\left( \ref{X2,g,X2F71,1H}\right)
\end{equation*}%
\begin{equation*}
B(g\otimes 1_{H};G,gx_{1})+B(x_{2}\otimes 1_{H};G,gx_{1}x_{2})=0\left( \ref%
{X2,g,X2F71,gx1}\right)
\end{equation*}

\begin{eqnarray*}
&&-2\beta _{2}B(x_{1}\otimes 1_{H};1_{A},gx_{2})+\gamma _{2}B(x_{1}\otimes
1_{H};G,g)\left( \ref{X2,x1,X2F11,g}\right) \\
&&-\lambda \left[ B(g\otimes 1_{H};1_{A},g)+B(x_{1}\otimes
1_{H};1_{A},gx_{1})\right] -2B(gx_{1}x_{2}\otimes 1_{H};1_{A},g)=0
\end{eqnarray*}

\begin{equation*}
0=0\left( \ref{X2,x1,X2F11,x1}\right)
\end{equation*}%
\begin{equation*}
0=0.\left( \ref{X2,x1,X2F11,x2}\right)
\end{equation*}%
\begin{eqnarray*}
&&+\gamma _{2}B(x_{1}\otimes 1_{H};G,gx_{1}x_{2})-\lambda B(g\otimes
1_{H};1_{A},gx_{1}x_{2})\left( \ref{X2,x1,X2F11,gx1x2}\right) \\
&&-2B(gx_{1}x_{2}\otimes 1_{H};1_{A},gx_{1}x_{2})-2B\left( x_{1}\otimes
1_{H};1_{A},gx_{1}\right) =0
\end{eqnarray*}%
\begin{equation*}
0=0\left( \ref{X2,x1,X2F21,1H}\right)
\end{equation*}

\begin{equation*}
2\beta _{2}B(x_{1}\otimes 1_{H};G,gx_{1}x_{2})+\lambda B\left( g\otimes
1_{H};G,gx_{1}\right) +2B(x_{2}\otimes 1_{H};G,g)=0.\left( \ref%
{X2,x1,X2F21,gx1}\right)
\end{equation*}

\begin{eqnarray*}
&&\lambda \left[ B\left( g\otimes 1_{H};G,gx_{2}\right) -B(x_{1}\otimes
1_{H};G,gx_{1}x_{2})\right] \left( \ref{X2,x1,X2F21,gx2}\right) \\
&&+2-B(x_{1}\otimes 1_{H};G,g_{1})+2B\left( x_{1}\otimes 1_{H};G,g\right) =0.
\end{eqnarray*}

\begin{equation*}
0=0.\left( \ref{X2,x1,X2F31,1H}\right)
\end{equation*}%
\begin{gather*}
\gamma _{2}\left[ -B(g\otimes 1_{H};G,gx_{2})+B(x_{1}\otimes
1_{H};G,gx_{1}x_{2})\right] +\left( \ref{X2,x1,X2F31,gx2}\right) \\
-2B(x_{2}\otimes 1_{H};1_{A},gx_{2})-2B(gx_{1}x_{2}\otimes
1_{H};1_{A},gx_{1}x_{2}) \\
-2B(g\otimes 1_{H};1_{A},g)-2B(x_{1}\otimes 1_{H};1_{A},gx_{1})=0.
\end{gather*}%
\begin{equation*}
0=0.\left( \ref{X2,x1,X2F41,1H}\right)
\end{equation*}

\begin{eqnarray*}
&&-\gamma _{2}B(x_{1}\otimes 1_{H};G,gx_{1}x_{2})+\lambda B\left( g\otimes
1_{H};1_{A},gx_{1}x_{2}\right) \left( \ref{X2,x1,X2F41,gx1}\right) \\
&&+2B(x_{1}\otimes 1_{H};1_{A},gx_{1})+2B(gx_{1}x_{2}\otimes
1_{H};1_{A},gx_{1}x_{2})=0.
\end{eqnarray*}

\begin{eqnarray*}
&&\beta _{2}\left[ -B(g\otimes 1_{H};G,gx_{2})+B(x_{1}\otimes
1_{H};G,gx_{1}x_{2})\right] \left( \ref{X2,x1,X2F61,g}\right) \\
&&-B(x_{2}\otimes 1_{H};G,g)+B(gx_{1}x_{2}\otimes 1_{H};G,gx_{1})=0
\end{eqnarray*}%
\begin{equation*}
0=0\left( \ref{X2,x1,X2F61,x2}\right)
\end{equation*}

\begin{equation*}
B(x_{1}\otimes 1_{H};G,g)+B(gx_{1}x_{2}\otimes 1_{H};G,gx_{2})=0\left( \ref%
{X2,x1,X2F71,g}\right)
\end{equation*}%
\begin{equation*}
0=0\left( \ref{X2,x1,X2F81,1H}\right)
\end{equation*}%
\begin{eqnarray*}
2 &&\beta _{2}\left[ B(g\otimes 1_{H};1_{A},g)+B(x_{2}\otimes \
1_{H};1_{A},gx_{2})\right] \left( \ref{X2,x2,X2F11,g}\right) \\
&&-\gamma _{2}B(x_{2}\otimes 1_{H};G,g)+\lambda B(x_{2}\otimes
1_{H};1_{A},gx_{1})=0.
\end{eqnarray*}%
\begin{equation*}
0=0.\left( \ref{X2,x2,X2F11,x2}\right)
\end{equation*}

\begin{eqnarray*}
+ &&2\beta _{2}B\left( g\otimes 1_{H};1_{A},gx_{1}x_{2}\right) -\gamma
_{2}B(x_{2}\otimes 1_{H};G,gx_{1}x_{2})\left( \ref{X2,x2,X2F11,gx1x2}\right)
\\
&&+2B\left( x_{2}\otimes 1_{H};1_{A},gx_{1}\right) =0
\end{eqnarray*}

\begin{equation*}
\beta _{2}\left[ B(g\otimes 1_{H};G,gx_{1})+B(x_{2}\otimes \
1_{H};G,gx_{1}x_{2})\right] =0\left( \ref{X2,x2,X2F21,gx1}\right)
\end{equation*}%
\begin{eqnarray*}
&&2\beta _{2}B(g\otimes 1_{H};G,gx_{2})-\lambda B(x_{2}\otimes
1_{H};G,gx_{1}x_{2})\left( \ref{X2,x2,X2F21,gx2}\right) \\
&&+2B\left( x_{2}\otimes 1_{H};G,g\right) =0
\end{eqnarray*}

\begin{equation*}
0=0.\left( \ref{X2,x2,X2F41,1H}\right)
\end{equation*}

\begin{eqnarray*}
&&\gamma _{2}B(g\otimes 1_{H};G,gx_{2})-\lambda B\left( g\otimes
1_{H};1_{A},gx_{1}x_{2}\right) \left( \ref{X2,x2,X2F41,gx2}\right) \\
&&+2\left[ +B(g\otimes 1_{H};1_{A},g)+B(x_{2}\otimes \ 1_{H};1_{A},gx_{2})%
\right] =0.
\end{eqnarray*}

\begin{equation*}
\lambda \left[ B(g\otimes 1_{H};G,gx_{1})+B(x_{2}\otimes \
1_{H};G,gx_{1}x_{2})\right] =0\left( \ref{X2,x2,X2F71,g}\right)
\end{equation*}%
\begin{equation*}
\gamma _{2}B(x_{1}x_{2}\otimes 1_{H};G,g)+\lambda B(x_{1}x_{2}\otimes
1_{H};X_{1},g)=0.\left( \ref{X2,x1x2,X2F11,g}\right)
\end{equation*}

\begin{gather*}
2\beta _{2}B(x_{1}x_{2}\otimes 1_{H};X_{2},x_{1})+\gamma
_{2}B(x_{1}x_{2}\otimes 1_{H};G,x_{1})\left( \ref{X2,x1x2,X2F11,x1}\right) \\
+\lambda B(x_{1}x_{2}\otimes 1_{H};X_{1},x_{1})=0.
\end{gather*}

\begin{gather*}
2\beta _{2}B(x_{1}x_{2}\otimes 1_{H};X_{2},x_{2})+\gamma
_{2}B(x_{1}x_{2}\otimes 1_{H};G,x_{2})\left( \ref{X2,x1x2,X2F11,x2}\right) \\
+\lambda B(x_{1}x_{2}\otimes 1_{H};X_{1},x_{2})=0.
\end{gather*}%
\begin{eqnarray*}
&&\gamma _{2}B(x_{1}x_{2}\otimes 1_{H};G,gx_{1}x_{2})+\left( \ref%
{X2,x1x2,X2F11,gx1x2}\right) \\
&&+2B(x_{1}x_{2}\otimes 1_{H};1_{A},gx_{1})=0.
\end{eqnarray*}

\begin{equation*}
0=0.\left( \ref{X2,x1x2,X2F21,1H}\right)
\end{equation*}

\begin{equation*}
2B(x_{1}x_{2}\otimes 1_{H};G,g)-\lambda B(x_{1}x_{2}\otimes
1_{H};GX_{1},gx_{2})=0.\left( \ref{X2,x1x2,X2F21,gx2}\right)
\end{equation*}%
\begin{eqnarray*}
&&2\beta _{2}\left[
\begin{array}{c}
+1-B(x_{1}x_{2}\otimes 1_{H};1_{A},x_{1}x_{2}) \\
-B(x_{1}x_{2}\otimes 1_{H};X_{2},x_{1})+B(x_{1}x_{2}\otimes
1_{H};X_{1},x_{2})%
\end{array}%
\right] \left( \ref{X2,x1x2,X2F31,1H}\right) \\
&=&0.
\end{eqnarray*}%
\begin{equation*}
\gamma _{2}B(x_{1}x_{2}\otimes 1_{H};GX_{1},gx_{2})+2B(x_{1}x_{2}\otimes
1_{H};X_{1},g)=0.\left( \ref{X2,x1x2,X2F31,gx2}\right)
\end{equation*}

\begin{equation*}
\left[
\begin{array}{c}
+1-B(x_{1}x_{2}\otimes 1_{H};1_{A},x_{1}x_{2}) \\
-B(x_{1}x_{2}\otimes 1_{H};X_{2},x_{1})+B(x_{1}x_{2}\otimes
1_{H};X_{1},x_{2})%
\end{array}%
\right] =0.\left( \ref{X2,x1x2,X2F41,1H}\right)
\end{equation*}%
\begin{gather*}
2B(x_{1}x_{2}\otimes 1_{H};1_{A},gx_{2})+\left( \ref{X2,x1x2,X2F41,gx2}%
\right) \\
+2B(x_{1}x_{2}\otimes 1_{H};X_{2},g)=0.
\end{gather*}

\begin{gather*}
2B(x_{1}x_{2}\otimes 1_{H};X_{1},g)-\left( \ref{X2,x1x2,X2F51,g}\right) \\
\gamma _{2}\left[
\begin{array}{c}
-B(x_{1}x_{2}\otimes 1_{H};G,gx_{1}x_{2})+ \\
B(x_{1}x_{2}\otimes 1_{H};GX_{2},gx_{1})-B(x_{1}x_{2}\otimes
1_{H};GX_{1},gx_{2})%
\end{array}%
\right] =0.
\end{gather*}%
\begin{gather*}
2B(x_{1}x_{2}\otimes 1_{H};G,g)\left( \ref{X2,x1x2,X2F71,g}\right) \\
+\lambda \left[
\begin{array}{c}
-B(x_{1}x_{2}\otimes 1_{H};G,gx_{1}x_{2})+ \\
B(x_{1}x_{2}\otimes 1_{H};GX_{2},gx_{1})-B(x_{1}x_{2}\otimes
1_{H};GX_{1},gx_{2})%
\end{array}%
\right] =0.
\end{gather*}%
\begin{equation*}
B(x_{1}x_{2}\otimes 1_{H};G,gx_{1}x_{2})-B(x_{1}x_{2}\otimes
1_{H};GX_{2},gx_{1})=0.\left( \ref{X2,x1x2,X2F71,gx1x2}\right)
\end{equation*}%
\begin{gather*}
2\beta _{2}B(x_{1}x_{2}\otimes 1_{H};X_{2},x_{2})+\gamma
_{2}B(x_{1}x_{2}\otimes 1_{H};G,x_{2})\left( \ref{X2,gx1,X2F11,1H}\right) \\
+\lambda \left[ -1+B(x_{1}x_{2}\otimes
1_{H};1_{A},x_{1}x_{2})+B(x_{1}x_{2}\otimes 1_{H};X_{2},x_{1})\right] =0.
\end{gather*}

\begin{equation*}
B(x_{1}x_{2}\otimes 1_{H};1_{A},gx_{2})+B(x_{1}x_{2}\otimes
1_{H};X_{2},g)=0.\left( \ref{X2,gx1,X2F11,gx2}\right)
\end{equation*}%
\begin{equation*}
\lambda \left[ B(x_{1}x_{2}\otimes 1_{H};G,gx_{1}x_{2})-B(x_{1}x_{2}\otimes
1_{H};GX_{2},gx_{1})\right] =0.\left( \ref{X2,gx1,X2F21,g}\right)
\end{equation*}%
\begin{gather*}
2\beta _{2}\left[ -1+B(x_{1}x_{2}\otimes
1_{H};1_{A},x_{1}x_{2})-B(x_{1}x_{2}\otimes 1_{H};X_{1},x_{2})\right] \left( %
\ref{X2,gx2,X2F11,1H}\right) \\
+\gamma _{2}\left[ -B(x_{1}x_{2}\otimes 1_{H};G,x_{1})\right] + \\
-\lambda B(x_{1}x_{2}\otimes 1_{H};X_{1},x_{1})=0.
\end{gather*}

\begin{eqnarray*}
&&\gamma _{2}\left[ B(x_{1}x_{2}\otimes
1_{H};G,gx_{1}x_{2})+B(x_{1}x_{2}\otimes 1_{H};GX_{1},gx_{2})\right] \left( %
\ref{X2,gx2,X2F11,gx2}\right) \\
&&+ \\
&&+\left[ +2B(x_{1}x_{2}\otimes 1_{H};1_{A},gx_{1})+2B(x_{1}x_{2}\otimes
1_{H};X_{1},g)\right] =0.
\end{eqnarray*}%
\begin{eqnarray*}
&&2\left[ -B(x_{1}x_{2}\otimes 1_{H};1_{A},gx_{1})-B(x_{1}x_{2}\otimes
1_{H};X_{1},g)\right] \left( \ref{X2,gx2,X2F41,g}\right) \\
&&+\gamma _{2}\left[ -B(x_{1}x_{2}\otimes
1_{H};G,gx_{1}x_{2})-B(x_{1}x_{2}\otimes 1_{H};GX_{1},gx_{2})\right] \\
&=&0.
\end{eqnarray*}%
\begin{equation*}
0=0\left( \ref{X2,gx1x2,X2F11,1H}\right)
\end{equation*}%
\begin{equation*}
0=0\left( \ref{X2,gx1x2,X2F11,x1x2}\right)
\end{equation*}%
\begin{gather*}
2\beta _{2}\left[ B(x_{1}\otimes 1_{H};1_{A},gx_{1})+B(gx_{1}x_{2}\otimes
1_{H};1_{A},gx_{1}x_{2})\right] \left( \ref{X2,gx1x2,X2F11,gx1}\right) \\
+\gamma _{2}B(gx_{1}x_{2}\otimes 1_{H};G,gx_{1})-\lambda B(x_{2}\otimes
1_{H};1_{A},gx_{1})=0.
\end{gather*}

\begin{gather*}
2\beta _{2}B(x_{1}\otimes 1_{H};1_{A},gx_{2})+\gamma
_{2}B(gx_{1}x_{2}\otimes 1_{H};G,gx_{2})\left( \ref{X2,gx1x2,X2F11,gx2}%
\right) \\
+\lambda \left[ -B(x_{2}\otimes 1_{H};1_{A},gx_{2})-B(gx_{1}x_{2}\otimes
1_{H};1_{A},gx_{1}x_{2})\right] \\
+2B(gx_{1}x_{2}\otimes 1_{H};1_{A},g)=0.
\end{gather*}

\begin{gather*}
2\beta _{2}\left[ -B(x_{1}\otimes 1_{H};G,g)-B(gx_{1}x_{2}\otimes
1_{H};G,gx_{2})\right] \left( \ref{X2,gx1x2,X2F21,g}\right) \\
+\lambda \left[ B(x_{2}\otimes 1_{H};G,g)-B(gx_{1}x_{2}\otimes
1_{H};G,gx_{1})\right] =0.
\end{gather*}%
\begin{gather*}
-2\beta _{2}B(x_{1}\otimes 1_{H};G,gx_{1}x_{2})\left( \ref%
{X2,gx1x2,X2F21,gx1x2}\right) \\
+\lambda B(x_{2}\otimes 1_{H};G,gx_{1}x_{2})-2B(gx_{1}x_{2}\otimes
1_{H};G,gx_{1})=0.
\end{gather*}%
\begin{gather*}
2\beta _{2}\left[
\begin{array}{c}
B(g\otimes 1_{H};1_{A},g)+B(x_{2}\otimes \ 1_{H};1_{A},gx_{2}) \\
+B(x_{1}\otimes 1_{H};1_{A},gx_{1})+B(gx_{1}x_{2}\otimes
1_{H};1_{A},gx_{1}x_{2})%
\end{array}%
\right] \left( \ref{X2,gx1x2,X2F31,g}\right) \\
-\gamma _{2}\left[ B(x_{2}\otimes 1_{H};G,g)-B(gx_{1}x_{2}\otimes
1_{H};G,gx_{1})\right] =0.
\end{gather*}

\begin{eqnarray*}
&&\gamma _{2}\left[ B(x_{1}\otimes 1_{H};G,g)+B(gx_{1}x_{2}\otimes
1_{H};G,gx_{2})\right] \left( \ref{X2,gx1x2,X2F41,g}\right) \\
&&-\lambda \left[
\begin{array}{c}
B(g\otimes 1_{H};1_{A},g)+B(x_{2}\otimes \ 1_{H};1_{A},gx_{2}) \\
+B(x_{1}\otimes 1_{H};1_{A},gx_{1})+B(gx_{1}x_{2}\otimes
1_{H};1_{A},gx_{1}x_{2})%
\end{array}%
\right] =0.
\end{eqnarray*}

\begin{equation*}
0=0.\left( \ref{X2,gx1x2,X2F41,x1}\right)
\end{equation*}

\begin{equation*}
0=0.\left( \ref{X2,gx1x2,X2F51,1H}\right)
\end{equation*}

\begin{equation*}
B(x_{2}\otimes 1_{H};G,g)-B(gx_{1}x_{2}\otimes 1_{H};G,gx_{1})=0.\left( \ref%
{X2,gx1x2,X2F61,gx2}\right)
\end{equation*}%
\begin{equation*}
0=0.\left( \ref{X2,gx1x2,X2F71,1H}\right)
\end{equation*}

\subsection{LIST\ OF\ MONOMIAL\ EQUALITIES 4}

We list the new monomial equalities we obtained above taking out constants
and relabel them.

\begin{equation*}
B\left( g\otimes 1_{H};1_{A},x_{2}\right) =0.\left( \ref{X1,g, X1F11,1H}%
\right)
\end{equation*}%
\begin{equation*}
B\left( g\otimes 1_{H};1_{A},x_{2}\right) =0.\left( \ref{X1,x1,X1F31,1H}%
\right)
\end{equation*}%
\begin{equation*}
B(x_{2}\otimes \ 1_{H};G,x_{2})=0.\left( \ref{X1,x2,X1F21,1H}\right)
\end{equation*}%
\begin{equation*}
B(x_{1}\otimes 1_{H};1_{A},1_{H})=0\left( \ref{X1,gx1x2,X1F11,1H}\right)
\end{equation*}%
\begin{equation*}
B(x_{1}\otimes 1_{H};1_{A},1_{H})=0\left( \ref{X1,gx1x2,X1F51,1H}\right)
\end{equation*}%
\begin{equation*}
B(x_{2}\otimes 1_{H};1,x_{1}x_{2})=0\left( \ref{X2,g,X2F31,x2}\right)
\end{equation*}%
Thus we get%
\begin{equation*}
B\left( g\otimes 1_{H};1_{A},x_{2}\right) =0.\left( \ref{X1,g, X1F11,1H}%
\right) \text{this appeared already in }\left( \ref{LME1}\right) \text{ }
\end{equation*}%
\begin{equation*}
B(x_{1}\otimes 1_{H};1_{A},1_{H})=0\left( \ref{X1,gx1x2,X1F11,1H}\right)
\text{ this appeared already in }\left( \ref{LME3}\right)
\end{equation*}%
\begin{equation*}
B(x_{2}\otimes 1_{H};1,x_{1}x_{2})=0\left( \ref{X2,g,X2F31,x2}\right) \text{%
this appeared already in }\left( \ref{LME1}\right)
\end{equation*}%
\begin{equation*}
B(x_{2}\otimes \ 1_{H};G,x_{2})=0.\left( \ref{X1,x2,X1F21,1H}\right) \text{%
this appeared already in }\left( \ref{LME3}\right)
\end{equation*}%
Now we take the complete list of equalities $\left( \ref{LAE3}\right) $ and
cancel there the terms above

\subsection{LIST\ OF\ ALL\ EQUALITIES 5\label{LAE5}}

\subsubsection{$G$}

\begin{equation*}
0=0\text{ }\left( \ref{G,g,GF1,1H}\right)
\end{equation*}

\begin{equation*}
2\alpha B(g\otimes 1_{H};G,gx_{1})+\gamma _{2}B(g\otimes
1_{H};1_{A},gx_{1}x_{2})=0\left( \ref{G,g, GF1,gx1}\right)
\end{equation*}

\begin{equation*}
2\alpha B(g\otimes 1_{H};G,gx_{2})-\gamma _{1}B\left( g\otimes
1_{H};1_{A},gx_{1}x_{2}\right) =0.\left( \ref{G,g, GF1,gx2}\right)
\end{equation*}

\begin{equation*}
\gamma _{1}B\left( g\otimes 1_{H};G,gx_{1}\right) +\gamma _{2}B\left(
g\otimes 1_{H};G,gx_{2}\right) =0.\left( \ref{G,g, GF2,g}\right)
\end{equation*}%
\begin{equation*}
\begin{array}{c}
2\alpha B(x_{1}\otimes 1_{H};G,g)+ \\
+\gamma _{1}\left[ -B(g\otimes 1_{H};1_{A},g)-B(x_{1}\otimes
1_{H};1_{A},gx_{1})\right] -\gamma _{2}B(x_{1}\otimes 1_{H};1_{A},gx_{2})%
\end{array}%
=0.\left( \ref{G,x1, GF1,g}\right)
\end{equation*}

\begin{equation*}
-\gamma _{1}B(g\otimes 1_{H};1_{A},gx_{1}x_{2})+2\alpha B(x_{1}\otimes
1_{H};G,gx_{1}x_{2})=0.\left( \ref{G,x1, GF1,gx1x2}\right)
\end{equation*}%
\begin{equation*}
0=0.\left( \ref{G,x1, GF2,1H}\right)
\end{equation*}%
\begin{equation*}
0=0\text{ }\left( \ref{G,x1, GF2,x1x2}\right)
\end{equation*}%
\begin{equation*}
\gamma _{1}B\left( g\otimes 1_{H};G,gx_{1}\right) +\gamma _{2}B(x_{1}\otimes
1_{H};G,gx_{1}x_{2})=0.\left( \ref{G,x1, GF2,gx1}\right)
\end{equation*}

\begin{equation*}
\gamma _{1}[B\left( g\otimes 1_{H};G,gx_{2}\right) -B(x_{1}\otimes
1_{H};G,gx_{1}x_{2})]=0.\left( \ref{G,x1,GF2,gx2}\right)
\end{equation*}

\begin{equation*}
0=0\left( \ref{G,x1,GF3,1H}\right)
\end{equation*}%
\begin{equation*}
-\gamma _{2}B(g\otimes 1_{H};1_{A},gx_{1}x_{2})-2\alpha B(g\otimes
1_{H};G,gx_{1})=0.\left( \ref{G,x1,GF3,gx1}\right)
\end{equation*}

\begin{equation*}
\alpha \left[ -B(g\otimes 1_{H};g,gx_{2})+B(x_{1}\otimes 1_{H};g,gx_{1}x_{2})%
\right] =0.\left( \ref{G,x1, GF3,x2}\right)
\end{equation*}%
\begin{equation*}
0=0.\left( \ref{G,x1, GF4,1H}\right)
\end{equation*}

\begin{equation*}
-\gamma _{1}B(g\otimes 1_{H};1_{A},gx_{1}x_{2})+2\alpha B(x_{1}\otimes
1_{H};G,gx_{1}x_{2})=0.\left( \ref{G,x1, GF4,gx1}\right)
\end{equation*}

\begin{equation*}
\alpha \lbrack -B(g\otimes 1_{H};g,gx_{2})+B(x_{1}\otimes
1_{H};G,gx_{1}x_{2})]=0.\left( \ref{G,x1, GF5,g}\right)
\end{equation*}%
\begin{equation*}
0=0.\left( \ref{G,x1,GF6,x1}\right)
\end{equation*}%
\begin{equation*}
0=0.\left( \ref{G,x1,GF6,x2}\right)
\end{equation*}%
\begin{equation*}
\begin{array}{c}
2\alpha B(x_{2}\otimes 1_{H};G,g)+ \\
-\gamma _{1}B(x_{2}\otimes 1_{H};1_{A},gx_{1})+\gamma _{2}\left[ -B(g\otimes
1_{H};1_{A},g)-B(x_{2}\otimes \ 1_{H};1_{A},gx_{2})\right]%
\end{array}%
=0.\left( \ref{G,x2, GF1,g}\right)
\end{equation*}%
\begin{equation*}
0=0\left( \ref{G,x2, GF1,x1}\right)
\end{equation*}

\begin{equation*}
2\alpha B(x_{2}\otimes 1_{H};G,gx_{1}x_{2})-\gamma _{2}B\left( g\otimes
1_{H};1_{A},gx_{1}x_{2}\right) =0.\left( \ref{G,x2, GF1,gx1x2}\right)
\end{equation*}%
\begin{equation*}
0=0.\left( \ref{G,x2, GF2,1H}\right)
\end{equation*}%
\begin{equation*}
0=0\left( \ref{G,x2, GF2,x1x2}\right)
\end{equation*}%
\begin{equation*}
\gamma _{2}\left[ B(g\otimes 1_{H};G,gx_{1})+B(x_{2}\otimes \
1_{H};G,gx_{1}x_{2})\right] =0.\left( \ref{G,x2, GF2,gx1}\right)
\end{equation*}

\begin{equation*}
-\gamma _{1}B(x_{2}\otimes 1_{H};G,gx_{1}x_{2})+\gamma _{2}B(g\otimes
1_{H};G,gx_{2})=0.\left( \ref{G,x2, GF2,gx2}\right)
\end{equation*}%
\begin{equation*}
0=0\left( \ref{G,x2, FG3,1H}\right)
\end{equation*}

\begin{equation*}
-\gamma _{2}B(g\otimes 1_{H};1_{A},gx_{1}x_{2})-2\alpha B(x_{2}\otimes
1_{H};G,gx_{1}x_{2})=0.\left( \ref{G,x2, GF3,gx2}\right)
\end{equation*}%
\begin{equation*}
0=0.\left( \ref{G,x2, GF4,1H}\right)
\end{equation*}%
\begin{equation*}
\alpha \left[ B(g\otimes 1_{H};G,gx_{1})+B(x_{2}\otimes \
1_{H};G,gx_{1}x_{2})\right] =0\left( \ref{G,x2, GF4,gx1}\right)
\end{equation*}

\begin{equation*}
\gamma _{1}B(g\otimes 1_{H};1_{A},gx_{1}x_{2})-2\alpha B(g\otimes
1_{H};G,gx_{2})=0.\left( \ref{G,x2, GF4,gx2}\right)
\end{equation*}%
\begin{equation*}
\alpha \left[ B(g\otimes 1_{H};G,gx_{1})+B(x_{2}\otimes \
1_{H};G,gx_{1}x_{2})\right] =0\left( \ref{G,x2, GF5,g}\right)
\end{equation*}

\begin{equation*}
\gamma _{2}[B(g\otimes 1_{H},G,gx_{1})+B(x_{2}\otimes
1_{H};G,gx_{1}x_{2})]=0.\left( \ref{G,x2, GF6,g}\right)
\end{equation*}%
\begin{equation*}
0=0.\left( \ref{G,x2; GF6,x2}\right)
\end{equation*}

\begin{equation*}
\gamma _{1}[B(g\otimes 1_{H};G,gx_{1})+B(x_{2}\otimes
1_{H};G,gx_{1}x_{2})]=0.\left( \ref{G,x2; GF7,g}\right)
\end{equation*}%
\begin{equation*}
0=0.\left( \ref{G,x2; GF7,x1}\right)
\end{equation*}%
\begin{equation*}
\gamma _{1}B(x_{1}x_{2}\otimes 1_{H};X_{1},g)+\gamma _{2}B(x_{1}x_{2}\otimes
1_{H};X_{2},g)=0.\left( \ref{G,x1x2, GF1,g}\right)
\end{equation*}

\begin{equation*}
\begin{array}{c}
2\alpha B(x_{1}x_{2}\otimes 1_{H};G,x_{1})+ \\
+\gamma _{1}B(x_{1}x_{2}\otimes 1_{H};X_{1},x_{1})+\gamma
_{2}B(x_{1}x_{2}\otimes 1_{H};X_{2},x_{1})%
\end{array}%
=0.\left( \ref{G,x1x2, GF1,x1}\right)
\end{equation*}%
\begin{equation*}
\begin{array}{c}
2\alpha B(x_{1}x_{2}\otimes 1_{H};G,x_{2})+ \\
+\gamma _{1}B(x_{1}x_{2}\otimes 1_{H};X_{1},x_{2})+\gamma
_{2}B(x_{1}x_{2}\otimes 1_{H};X_{2},x_{2})%
\end{array}%
=0.\left( \ref{G,x1x2, GF1,x2}\right)
\end{equation*}

\begin{equation*}
0=0\left( \ref{G,x1x2, GF2,1H}\right)
\end{equation*}%
\begin{equation*}
\begin{array}{c}
2B(x_{1}x_{2}\otimes 1_{H};1_{A},gx_{1}) \\
+\gamma _{2}B(x_{1}x_{2}\otimes 1_{H};GX_{2},gx_{1})%
\end{array}%
=0.\left( \ref{G,x1x2, GF2,gx1}\right)
\end{equation*}%
\begin{equation*}
\begin{array}{c}
2B(x_{1}x_{2}\otimes 1_{H};1_{A},gx_{2}) \\
+\gamma _{1}B(x_{1}x_{2}\otimes 1_{H};GX_{1},gx_{2})%
\end{array}%
=0\left( \ref{G,x1x2, GF2,gx2}\right)
\end{equation*}

\begin{equation*}
0=0\text{ }\left( \ref{G,x1x2, GF3,gx1}\right)
\end{equation*}%
\begin{gather*}
-\gamma _{1}[1-B(x_{1}x_{2}\otimes 1_{H};1_{A},x_{1}x_{2})+ \\
+B(x_{1}x_{2}\otimes 1_{H};x_{1},x_{2})-B(x_{1}x_{2}\otimes
1_{H};x_{2},x_{1})]+ \\
-2\alpha B(x_{1}x_{2}\otimes 1_{H};GX_{2},1_{H}=0.\left( \ref{G,x1x2, GF4,1H}%
\right)
\end{gather*}

\begin{gather*}
-2B(x_{1}x_{2}\otimes 1_{H};X_{1},g)+\gamma _{2}[-B(x_{1}x_{2}\otimes
1_{H};g,gx_{1}x_{2})+ \\
-B(x_{1}x_{2}\otimes 1_{H};GX_{1},gx_{2})+B(x_{1}x_{2}\otimes
1_{H};GX_{2},gx_{1})]=0\left( \ref{G,x1x2, GF6,g}\right)
\end{gather*}%
\begin{gather*}
-2B(x_{1}x_{2}\otimes 1_{H};X_{2},g)+ \\
-\gamma _{1}[-B(x_{1}x_{2}\otimes 1_{H};g,gx_{1}x_{2})-B(x_{1}x_{2}\otimes
1_{H};GX_{1},gx_{2})+B(x_{1}x_{2}\otimes 1_{H};GX_{2},gx_{1})]=0\left( \ref%
{G,x1x2, GF7,g}\right)
\end{gather*}

\begin{equation*}
\begin{array}{c}
2\alpha B(x_{1}x_{2}\otimes 1_{H};G,x_{2})+ \\
+\gamma _{1}\left[ -1+B(x_{1}x_{2}\otimes
1_{H};1_{A},x_{1}x_{2})+B(x_{1}x_{2}\otimes 1_{H};X_{2},x_{1})\right]
+\gamma _{2}B(x_{1}x_{2}\otimes 1_{H};X_{2},x_{2})%
\end{array}%
=0.\left( \ref{G,gx1, GF1,1H}\right)
\end{equation*}

\begin{equation*}
\begin{array}{c}
2\left[ B(x_{1}x_{2}\otimes 1_{H};1_{A},gx_{2})+B(x_{1}x_{2}\otimes
1_{H};X_{2},g)\right] \\
+\gamma _{1}\left[ -B(x_{1}x_{2}\otimes
1_{H};G,gx_{1}x_{2})+B(x_{1}x_{2}\otimes 1_{H};GX_{2},gx_{1})\right]%
\end{array}%
\left( \ref{G,gx1, GF2,g}\right)
\end{equation*}%
\begin{equation*}
\begin{array}{c}
2\alpha \left[ -B(x_{1}x_{2}\otimes 1_{H};G,x_{1})\right] + \\
-\gamma _{1}B(x_{1}x_{2}\otimes 1_{H};X_{1},x_{1})+\gamma _{2}\left[
-1+B(x_{1}x_{2}\otimes 1_{H};1_{A},x_{1}x_{2})-B(x_{1}x_{2}\otimes
1_{H};X_{1},x_{2})\right]%
\end{array}%
=0\left( \ref{G,gx2, GF1,1H}\right)
\end{equation*}%
\begin{equation*}
\begin{array}{c}
2\left[ -B(x_{1}x_{2}\otimes 1_{H};1_{A},gx_{1})-B(x_{1}x_{2}\otimes
1_{H};X_{1},g)\right] \\
+\gamma _{2}\left[ -B(x_{1}x_{2}\otimes
1_{H};G,gx_{1}x_{2})-B(x_{1}x_{2}\otimes 1_{H};GX_{1},gx_{2})\right]%
\end{array}%
=0.\left( \ref{G,gx2, GF2,g}\right)
\end{equation*}

\begin{equation*}
0=0\left( \ref{G,gx1x2, GF1,1H}\right)
\end{equation*}%
\begin{equation*}
0=0.\left( \ref{G,gx1x2, GF1,x1x2}\right)
\end{equation*}%
\begin{equation*}
\begin{array}{c}
2\alpha B(x_{2}\otimes 1_{H};G,g)+ \\
-\gamma _{1}B(x_{2}\otimes 1_{H};1_{A},gx_{1})+\gamma _{2}\left[
B(x_{1}\otimes 1_{H};1_{A},gx_{1})+B(gx_{1}x_{2}\otimes
1_{H};1_{A},gx_{1}x_{2})\right]%
\end{array}%
=0.\left( \ref{G,gx1x2, GF1,gx1}\right)
\end{equation*}

\begin{equation*}
\begin{array}{c}
2\alpha B(gx_{1}x_{2}\otimes 1_{H};G,gx_{2})+ \\
+\gamma _{1}\left[ -B(x_{2}\otimes 1_{H};1_{A},gx_{2})-B(gx_{1}x_{2}\otimes
1_{H};1_{A},gx_{1}x_{2})\right] +\gamma _{2}B(x_{1}\otimes
1_{H};1_{A},gx_{2})%
\end{array}%
=0.\left( \ref{G,gx1x2, GF1,gx2}\right)
\end{equation*}%
\begin{gather*}
\gamma _{1}\left[ B(x_{2}\otimes 1_{H};G,g)-B(gx_{1}x_{2}\otimes
1_{H};G,gx_{1})\right] \left( \ref{G,gx1x2, GF2,g}\right) \\
+\gamma _{2}\left[ -B(x_{1}\otimes 1_{H};G,g)-B(gx_{1}x_{2}\otimes
1_{H};G,gx_{2})\right] =0.
\end{gather*}%
\begin{equation*}
0=0.\left( \ref{G,gx1x2, GF2,x1}\right)
\end{equation*}%
\begin{equation*}
0=0.\left( \ref{G,gx1x2, GF2,x2}\right)
\end{equation*}%
\begin{equation*}
\gamma _{1}B(x_{2}\otimes 1_{H};G,gx_{1}x_{2})-\gamma _{2}B(x_{1}\otimes
1_{H};G,gx_{1}x_{2})=0.\left( \ref{G,gx1x2, GF2,gx1x2}\right)
\end{equation*}%
\begin{gather*}
\gamma _{2}\left[
\begin{array}{c}
B(g\otimes 1_{H};1_{A},g)+B(x_{2}\otimes \ 1_{H};1_{A},gx_{2})+ \\
+B(x_{1}\otimes 1_{H};1_{A},gx_{1})+B(gx_{1}x_{2}\otimes
1_{H};1_{A},gx_{1}x_{2})%
\end{array}%
\right] +\left( \ref{G,gx1x2, GF3,g}\right) \\
+2\alpha \left[ -B(x_{2}\otimes 1_{H};G,g)+B(x_{2}\otimes 1_{H};G,g)\right]
=0.
\end{gather*}

\begin{equation*}
0=0.\left( \ref{G,gx1x2, GF3,x1}\right)
\end{equation*}

\begin{equation*}
0=0\left( \ref{G,gx1x2, GF3,x2}\right)
\end{equation*}

\begin{equation*}
\gamma _{2}B\left( g\otimes 1_{H};1_{A},gx_{1}x_{2}\right) -2\alpha
B(x_{2}\otimes 1_{H};G,gx_{1}x_{2})=0.\left( \ref{G,gx1x2, GF3,gx1x2}\right)
\end{equation*}

\begin{gather*}
-\gamma _{1}\left[
\begin{array}{c}
B(g\otimes 1_{H};1_{A},g)+B(x_{2}\otimes \ 1_{H};1_{A},gx_{2})+ \\
+B(x_{1}\otimes 1_{H};1_{A},gx_{1})+B(gx_{1}x_{2}\otimes
1_{H};1_{A},gx_{1}x_{2})%
\end{array}%
\right] +\left( \ref{G,gx1x2, GF4,g}\right) \\
-2\alpha \left[ -B(x_{1}\otimes 1_{H};G,g)-B(gx_{1}x_{2}\otimes
1_{H};G,gx_{2})\right] =0
\end{gather*}%
\begin{equation*}
0=0\left( \ref{G,gx1x2, GF4,x1}\right)
\end{equation*}%
\begin{equation*}
0=0.\left( \ref{G,gx1x2, GF4,x2}\right)
\end{equation*}

\begin{equation*}
-\gamma _{1}B\left( g\otimes 1_{H};1_{A},gx_{1}x_{2}\right) +2\alpha
B(x_{1}\otimes 1_{H};G,gx_{1}x_{2})=0.\left( \ref{G,gx1x2, GF4,gx1x2}\right)
\end{equation*}

\begin{equation*}
\alpha \left[ B(g\otimes 1_{H};G,gx_{2})-B(x_{1}\otimes 1_{H};G,gx_{1}x_{2})%
\right] =0.\left( \ref{G,gx1x2, GF5,gx2}\right)
\end{equation*}

\begin{equation*}
0=0\left( \ref{G,gx1x2, GF6,1H}\right)
\end{equation*}%
\begin{equation*}
\gamma _{2}\left[ B(g\otimes 1_{H};G,gx_{2})-B(x_{1}\otimes
1_{H};G,gx_{1}x_{2})\right] =0.\left( \ref{G,gx1x2, GF6,gx1}\right)
\end{equation*}%
\begin{equation*}
\gamma _{2}\left[ B(g\otimes 1_{H};G,gx_{2})-B(x_{1}\otimes
1_{H};G,gx_{1}x_{2})\right] =0.\left( \ref{G,gx1x2, GF6,gx2}\right)
\end{equation*}%
\begin{equation*}
0=0\left( \ref{G,gx1x2, GF7,1H}\right)
\end{equation*}%
\begin{equation*}
0=0.\left( \ref{G,gx1x2, GF7,x1x2}\right)
\end{equation*}%
\begin{equation*}
\gamma _{1}\left[ B(g\otimes 1_{H};G,gx_{1})+B(x_{2}\otimes \
1_{H};G,gx_{1}x_{2})\right] =0.\left( \ref{G,gx1x2, GF7,gx1}\right)
\end{equation*}%
\begin{equation*}
\gamma _{1}\left[ B(g\otimes 1_{H};G,gx_{2})-B(x_{1}\otimes
1_{H};G,gx_{1}x_{2})\right] =0.\left( \ref{G,gx1x2, GF7,gx2}\right)
\end{equation*}

\subsubsection{$X_{1}$}

\begin{equation*}
0=0.\left( \ref{X1,g, X1F11,1H}\right)
\end{equation*}%
\begin{eqnarray*}
&&+\gamma _{1}B(g\otimes 1_{H};G,gx_{1})+\lambda B(g\otimes
1_{H};X_{2},gx_{1})+B(x_{1}\otimes 1_{H};1_{A},gx_{1})+\left( \ref{X1,g,
X1F11,gx1}\right) \\
&&+2B\left( g\otimes 1_{H};1_{A},g\right) +B\left( x_{1}\otimes
1_{H};1_{A},gx_{1}\right) =0
\end{eqnarray*}

\begin{equation*}
\gamma _{1}B(g\otimes 1_{H};G,gx_{2})+2B(x_{1}\otimes
1_{H};1_{A},gx_{2})-2\beta _{1}B\left( g\otimes
1_{H};1_{A},gx_{1}x_{2}\right) =0.\left( \ref{X1,g, X1F11,gx2}\right)
\end{equation*}

\begin{equation*}
2\beta _{1}B\left( g\otimes 1_{H};G,gx_{1}\right) +\lambda B\left( g\otimes
1_{H};G,gx_{2}\right) +2B(x_{1}\otimes 1_{H};G,g)=0.\left( \ref{X1,g,X1F21,g}%
\right)
\end{equation*}%
\begin{equation*}
0=0.\left( \ref{X1,g,X1F21,x1}\right)
\end{equation*}

\begin{equation*}
B\left( g\otimes 1_{H};G,gx_{2}\right) -B\left( x_{1}\otimes
1_{H};G,gx_{1}x_{2}\right) =0\left( \ref{X1,g,X1F21,gx1x2}\right)
\end{equation*}%
\begin{gather*}
+\lambda B\left( g\otimes 1_{H};1_{A},gx_{1}x_{2}\right) +\gamma _{1}B\left(
g\otimes 1_{H};G,gx_{1}\right) +\left( \ref{X1,g, X1F31,g}\right) \\
+2B(g\otimes 1_{H};1_{A},g)+2B(x_{1}\otimes 1_{H};1_{A},gx_{1})=0.
\end{gather*}

\begin{equation*}
-2\beta _{1}B\left( g\otimes 1_{H};1_{A},gx_{1}x_{2}\right) +2B(x_{1}\otimes
1_{H};1_{A},gx_{2})+\gamma _{1}B\left( g\otimes 1_{H};G,gx_{2}\right)
=0.\left( \ref{X1,g, X1F41,g}\right)
\end{equation*}%
\begin{equation*}
0=0.\left( \ref{X1,g, X1F41,x1}\right)
\end{equation*}

\begin{equation*}
0=0\left( \ref{X1,g,X1F51,1H}\right)
\end{equation*}%
\begin{equation*}
B(x_{1}\otimes 1_{H};G,gx_{1}x_{2})-B(g\otimes 1_{H};G,gx_{2})=0.\left( \ref%
{X1,x2,X1F71,gx1x2}\right)
\end{equation*}%
\begin{equation*}
2\beta _{1}\left[ -B(g\otimes 1_{H};1_{A},g)-B(x_{1}\otimes
1_{H};1_{A},gx_{1})\right] -\lambda B(x_{1}\otimes
1_{H};1_{A},gx_{2})=\gamma _{1}B(x_{1}\otimes 1_{H};G,g)\left( \ref%
{X1,x1,X1F11,g}\right)
\end{equation*}%
\begin{equation*}
0=0.\left( \ref{X1,x1,X1F11,x1}\right)
\end{equation*}

\begin{equation*}
-2\beta _{1}B(g\otimes 1_{H};1_{A},gx_{1}x_{2})=\gamma _{1}B(x_{1}\otimes
1_{H};G,gx_{1}x_{2})-2B\left( x_{1}\otimes 1_{H};1_{A},gx_{2}\right) .\left( %
\ref{X1,x1,X1F11,gx1x2}\right)
\end{equation*}%
\begin{equation*}
0=0\left( \ref{X1,x1,X1F21,1}\right)
\end{equation*}%
\begin{equation*}
\lambda B(x_{1}\otimes 1_{H};G,gx_{1}x_{2})=-2B\left( x_{1}\otimes
1_{H};G,g\right) .\left( \ref{X1,x1,X1F21,gx1}\right)
\end{equation*}%
\begin{equation*}
0=0.\left( \ref{X1,x1,X1F31,1H}\right)
\end{equation*}%
\begin{equation*}
0=0\left( \ref{X1,x1,X1F31,x1x2}\right)
\end{equation*}

\begin{gather*}
-\lambda B\left( g\otimes 1_{H};1_{A},gx_{1}x_{2}\right) -\gamma _{1}B\left(
g\otimes 1_{H};G,gx_{1}\right) + \\
-2\left[ B(g\otimes 1_{H};1_{A},g)+B(x_{1}\otimes 1_{H};1_{A},gx_{1})\right]
=0\left( \ref{X1,x1,X1F31,gx1}\right)
\end{gather*}

\begin{equation*}
\gamma _{1}\left[ -B\left( g\otimes 1_{H};G,gx_{2}\right) +B(x_{1}\otimes
1_{H};G,gx_{1}x_{2})\right] =0\left( \ref{X1,x1,X1F31,gx2}\right)
\end{equation*}

\begin{equation*}
2\beta _{1}B\left( g\otimes 1_{H};1_{A},gx_{1}x_{2}\right) -\gamma
_{1}B(x_{1}\otimes 1_{H};G,gx_{1}x_{2})-B(x_{1}\otimes
x_{1};X_{2},gx_{1})=0.\left( \ref{X1,x1,X1F41,gx1}\right)
\end{equation*}

\begin{eqnarray*}
&&-2\beta _{1}B(x_{2}\otimes 1_{H};1_{A},gx_{1})+\lambda \left[ -B(g\otimes
1_{H};1_{A},g)-B(x_{2}\otimes \ 1_{H};1_{A},gx_{2})\right] \\
&=&-\gamma _{1}B(x_{2}\otimes 1_{H};G,g)-2B(gx_{1}x_{2}\otimes
1_{H};1_{A},g).\left( \ref{X1,x2,X1F11,g}\right)
\end{eqnarray*}%
\begin{equation*}
0=0.\left( \ref{X1,x2,X1F11,x1}\right)
\end{equation*}%
\begin{equation*}
0=0.\left( \ref{X1,x2,X1F11,x2}\right)
\end{equation*}

\begin{gather*}
\lambda B\left( g\otimes 1_{H};1_{A},gx_{1}x_{2}\right) -\gamma
_{1}B(x_{2}\otimes 1_{H};G,gx_{1}x_{2})\left( \ref{X1,x2,X1F11,gx1x2}\right)
\\
-2B\left( x_{2}\otimes 1_{H};1_{A},gx_{2}\right) -2B(gx_{1}x_{2}\otimes
1_{H};1_{A},gx_{1}x_{2})=0.
\end{gather*}%
\begin{equation*}
0=0.\left( \ref{X1,x2,X1F21,1H}\right)
\end{equation*}%
\begin{equation*}
0=0.\left( \ref{X1,x2,X1F21,x1x2}\right)
\end{equation*}%
\begin{eqnarray*}
&&\lambda \left[ B(g\otimes 1_{H};G,gx_{1})+B(x_{2}\otimes \
1_{H};G,gx_{1}x_{2})\right] \\
&=&2B(gx_{1}x_{2}\otimes 1_{H};G,gx_{1})-2B\left( x_{2}\otimes
1_{H};G,g\right) .\left( \ref{X1,x2,X1F21,gx1}\right)
\end{eqnarray*}

\begin{gather*}
-2\beta _{1}B(x_{2}\otimes 1_{H};G,gx_{1}x_{2})+\lambda B(g\otimes
1_{H};G,gx_{2}) \\
-2B(gx_{1}x_{2}\otimes 1_{H};G,gx_{2})=0.\left( \ref{X1,x2,X1F21,gx2}\right)
\end{gather*}%
\begin{gather*}
-\lambda B\left( g\otimes 1_{H};1_{A},gx_{1}x_{2}\right) =-\gamma
_{1}B(x_{2}\otimes 1_{H};G,gx_{1}x_{2})+ \\
-2B(x_{2}\otimes 1_{H};1_{A},gx_{2})-2B(gx_{1}x_{2}\otimes
1_{H};1_{A},gx_{1}x_{2}).\left( \ref{X1,x2,X1F31,gx2}\right)
\end{gather*}

\begin{equation*}
0=0.\left( \ref{X1,x2,X1F41,1H}\right)
\end{equation*}%
\begin{equation*}
0=0.\left( \ref{X1,x2,X1F41,x1x2}\right)
\end{equation*}%
\begin{gather*}
B(x_{1}\otimes 1_{H};1_{A},gx_{1})+B(gx_{1}x_{2}\otimes
1_{H};1_{A},gx_{1}x_{2})\left( \ref{X1,x2,X1F41,gx1}\right) \\
+B(g\otimes 1_{H};1_{A},g)+B(x_{2}\otimes \ 1_{H};1_{A},gx_{2})=0.
\end{gather*}

\begin{equation*}
2\beta _{1}B\left( g\otimes 1_{H};1_{A},gx_{1}x_{2}\right) =\gamma
_{1}B(g\otimes 1_{H};G,gx_{2})+2B(x_{1}\otimes 1_{H};1_{A},gx_{2}).\left( %
\ref{X1,x2,X1F41,gx2}\right)
\end{equation*}

\begin{eqnarray*}
0 &=&\gamma _{1}\left[ B(g\otimes 1_{H};G,gx_{1})+B(x_{2}\otimes \
1_{H};G,gx_{1}x_{2})\right] \left( \ref{X1,x2,X1F51,g}\right) \\
&&+\left[
\begin{array}{c}
2B(g\otimes 1_{H};1_{A},g)+2B(x_{2}\otimes \ 1_{H};1_{A},gx_{2})+ \\
+2B(x_{1}\otimes 1_{H};1_{A},gx_{1})+2B(gx_{1}x_{2}\otimes
1_{H};1_{A},gx_{1}x_{2})%
\end{array}%
\right]
\end{eqnarray*}

\begin{equation*}
0=0.\left( \ref{X1,x2,X1F51,x2}\right)
\end{equation*}%
\begin{eqnarray*}
&&\lambda \left[ B(g\otimes 1_{H};G,gx_{1})+B(x_{2}\otimes \
1_{H};G,gx_{1}x_{2})\right] \\
&=&-2\left[ B(x_{2}\otimes 1_{H};G,g)-B(gx_{1}x_{2}\otimes 1_{H};G,gx_{1})%
\right] \left( \ref{X1,x2,X1F61,g}\right)
\end{eqnarray*}%
\begin{equation*}
\beta _{1}\left[ B(g\otimes 1_{H};G,gx_{1})+B(x_{2}\otimes \
1_{H};G,gx_{1}x_{2})\right] =0\left( \ref{X1,x2,X1F71,g}\right)
\end{equation*}%
\begin{equation*}
0=0\left( \ref{X1,x2,X1F71,x1}\right)
\end{equation*}%
\begin{equation*}
\lambda B(x_{1}x_{2}\otimes 1_{H};X_{2},g)-\gamma _{1}B(x_{1}x_{2}\otimes
1_{H};G,g)=0\left( \ref{X1,x1x2,X1F11,g}\right)
\end{equation*}%
\begin{equation*}
2\beta _{1}B(x_{1}x_{2}\otimes 1_{H};X_{1},x_{1})+\lambda
B(x_{1}x_{2}\otimes 1_{H};X_{2},x_{1})+\gamma _{1}B(x_{1}x_{2}\otimes
1_{H};G,x_{1})=0\left( \ref{X1,x1x2,X1F11,x1}\right)
\end{equation*}

\begin{equation*}
2\beta _{1}B(x_{1}x_{2}\otimes 1_{H};X_{1},x_{2})+\lambda
B(x_{1}x_{2}\otimes 1_{H};X_{2},x_{2})+\gamma _{1}B(x_{1}x_{2}\otimes
1_{H};G,x_{2})=0.\left( \ref{X1,x1x2,X1F11,x2}\right)
\end{equation*}

\begin{equation*}
\gamma _{1}B(x_{1}x_{2}\otimes 1_{H};G,gx_{1}x_{2})-2B(x_{1}x_{2}\otimes
1_{H};1_{A},gx_{2})=0.\left( \ref{X1,x1x2,X1F11,gx1x2}\right)
\end{equation*}%
\begin{equation*}
0=0.\left( \ref{X1,x1x2,X1F21,1H}\right)
\end{equation*}%
\begin{equation*}
\lambda B(x_{1}x_{2}\otimes 1_{H};GX_{2},gx_{1})+2B(x_{1}x_{2}\otimes
1_{H};G,g)=0.\left( \ref{X1,x1x2,X1F21,gx1}\right)
\end{equation*}

\begin{gather*}
\lambda \left[ +1-B(x_{1}x_{2}\otimes
1_{H};1_{A},x_{1}x_{2})-B(x_{1}x_{2}\otimes
1_{H};X_{2},x_{1})+B(x_{1}x_{2}\otimes 1_{H};X_{1},x_{2})\right] \left( \ref%
{X1,x1x2,X1F31,1H}\right) \\
=0
\end{gather*}%
\begin{equation*}
2B(x_{1}x_{2}\otimes 1_{H};1_{A},gx_{1})+2B(x_{1}x_{2}\otimes
1_{H};X_{1},g)=0\left( \ref{X1,x1x2,X1F31,gx1}\right)
\end{equation*}%
\begin{equation*}
2B(x_{1}x_{2}\otimes 1_{H};1_{A},gx_{2})+\gamma _{1}B(x_{1}x_{2}\otimes
1_{H};GX_{1},gx_{2})=0.\left( \ref{X1,x1x2,X1F31,gx2}\right)
\end{equation*}

\begin{equation*}
\left[
\begin{array}{c}
+1-B(x_{1}x_{2}\otimes 1_{H};1_{A},x_{1}x_{2})-B(x_{1}x_{2}\otimes
1_{H};X_{2},x_{1}) \\
+B(x_{1}x_{2}\otimes 1_{H};X_{1},x_{2})%
\end{array}%
\right] =0.\left( \ref{X1,x1x2,X1F41,1H}\right)
\end{equation*}

\begin{equation*}
\gamma _{1}B(x_{1}x_{2}\otimes 1_{H};GX_{2},gx_{1})+2B(x_{1}x_{2}\otimes
1_{H};X_{2},g)=0.\left( \ref{X1,x1x2,X1F41,gx1}\right)
\end{equation*}%
\begin{gather*}
2B(x_{1}x_{2}\otimes 1_{H};X_{2},g)+\left( \ref{X1,x1x2,X1F51,g}\right) \\
+\gamma _{1}\left[ -B(x_{1}x_{2}\otimes
1_{H};G,gx_{1}x_{2})+B(x_{1}x_{2}\otimes
1_{H};GX_{2},gx_{1})-B(x_{1}x_{2}\otimes 1_{H};GX_{1},gx_{2})\right] =0
\end{gather*}%
\begin{equation*}
2B(x_{1}x_{2}\otimes 1_{H};G,g)+\lambda \left[
\begin{array}{c}
-B(x_{1}x_{2}\otimes 1_{H};G,gx_{1}x_{2}) \\
+B(x_{1}x_{2}\otimes 1_{H};GX_{2},gx_{1})-B(x_{1}x_{2}\otimes
1_{H};GX_{1},gx_{2})%
\end{array}%
\right] =0\left( \ref{X1,x1x2,X1F61,g}\right)
\end{equation*}%
\begin{equation*}
B(x_{1}x_{2}\otimes 1_{H};G,gx_{1}x_{2})+B(x_{1}x_{2}\otimes
1_{H};GX_{1},gx_{2})=0\left( \ref{X1,x1x2,X1F61,gx1x2}\right)
\end{equation*}%
\begin{equation*}
+B(x_{1}x_{2}\otimes 1_{H};G,gx_{1}x_{2})+B(x_{1}x_{2}\otimes
1_{H};GX_{1},gx_{2})=0.\left( \ref{X1,x1x2,X1F81,gx1}\right)
\end{equation*}%
\begin{gather*}
2\beta _{1}\left[ -1+B(x_{1}x_{2}\otimes
1_{H};1_{A},x_{1}x_{2})+B(x_{1}x_{2}\otimes 1_{H};X_{2},x_{1})\right]
+\lambda B(x_{1}x_{2}\otimes 1_{H};X_{2},x_{2})\left( \ref{X1,gx1,X1F11,1H}%
\right) \\
=-\gamma _{1}B(x_{1}x_{2}\otimes 1_{H};G,x_{2})+ \\
-2\left[ B(x_{1}x_{2}\otimes 1_{H};1_{A},gx_{2})+B(x_{1}x_{2}\otimes
1_{H};X_{2},g)\right] .
\end{gather*}

\begin{gather*}
-\gamma _{1}\left[ B(x_{1}x_{2}\otimes
1_{H};G,gx_{1}x_{2})-B(x_{1}x_{2}\otimes 1_{H};GX_{2},gx_{1})\right] \left( %
\ref{X1,gx1,X1F11,gx1}\right) \\
+2\left[ B(x_{1}x_{2}\otimes 1_{H};1_{A},gx_{2})+B(x_{1}x_{2}\otimes
1_{H};X_{2},g)\right] =0.
\end{gather*}%
\begin{gather*}
2\left[ B(x_{1}x_{2}\otimes 1_{H};1_{A},gx_{2})+B(x_{1}x_{2}\otimes
1_{H};X_{2},g)\right] +\left( \ref{X1,gx1,X1F31,g}\right) \\
+ \\
+\gamma _{1}\left[ -B(x_{1}x_{2}\otimes
1_{H};G,gx_{1}x_{2})+B(x_{1}x_{2}\otimes 1_{H};GX_{2},gx_{1})\right] =0.
\end{gather*}%
\begin{gather*}
-2\beta _{1}B(x_{1}x_{2}\otimes 1_{H};X_{1},x_{1})+\lambda \left[
-1+B(x_{1}x_{2}\otimes 1_{H};1_{A},x_{1}x_{2})-B(x_{1}x_{2}\otimes
1_{H};X_{1},x_{2})\right] \left( \ref{X1,gx2,X1F11,1H}\right) \\
+\gamma _{1}\left[ -B(x_{1}x_{2}\otimes 1_{H};G,x_{1})\right] =0.
\end{gather*}%
\begin{equation*}
2B(x_{1}x_{2}\otimes 1_{H};1_{A},gx_{1})+2B(x_{1}x_{2}\otimes
1_{H};X_{1},g)=0\left( \ref{X1,gx2,X1F11,gx1}\right)
\end{equation*}%
\begin{equation*}
+\gamma _{1}\left[ B(x_{1}x_{2}\otimes
1_{H};G,gx_{1}x_{2})+B(x_{1}x_{2}\otimes 1_{H};GX_{1},gx_{2})\right]
=0\left( \ref{X1,gx2,X1F11,gx2}\right)
\end{equation*}%
\begin{equation*}
2\left[ B(x_{1}x_{2}\otimes 1_{H};1_{A},gx_{1})+B(x_{1}x_{2}\otimes
1_{H};X_{1},g)\right] =0\left( \ref{X1,gx2,X131,g}\right)
\end{equation*}%
\begin{equation*}
\gamma _{1}\left[ -B(x_{1}x_{2}\otimes
1_{H};G,gx_{1}x_{2})+B(x_{1}x_{2}\otimes 1_{H};GX_{1},gx_{2})\right]
=0\left( \ref{X1,gx2,X141,g}\right)
\end{equation*}%
\begin{equation*}
0=0\left( \ref{X1,gx1x2,X1F11,1H}\right)
\end{equation*}%
\begin{equation*}
0=0\left( \ref{X1,gx1x2,X1F11,x1x2}\right)
\end{equation*}%
\begin{gather*}
-2\beta _{1}B(x_{2}\otimes 1_{H};1_{A},gx_{1})+\lambda \left[ B(x_{1}\otimes
1_{H};1_{A},gx_{1})+B(gx_{1}x_{2}\otimes 1_{H};1_{A},gx_{1}x_{2})\right]
\left( \ref{X1,gx1x2,X1F11,gx1}\right) \\
+\gamma _{1}B(gx_{1}x_{2}\otimes 1_{H};G,gx_{1})+2B(gx_{1}x_{2}\otimes
1_{H};1_{A},g)=0.
\end{gather*}%
\begin{gather*}
2\beta _{1}\left[ -B(x_{2}\otimes 1_{H};1_{A},gx_{2})-B(gx_{1}x_{2}\otimes
1_{H};1_{A},gx_{1}x_{2})\right] \left( \ref{X1,gx1x2,X171,gx2}\right) \\
+\lambda B(x_{1}\otimes 1_{H};1_{A},gx_{2})+\gamma _{1}B(gx_{1}x_{2}\otimes
1_{H};G,gx_{2})=0.
\end{gather*}%
\begin{gather*}
2\beta _{1}\left[ B(x_{2}\otimes 1_{H};G,g)-B(gx_{1}x_{2}\otimes
1_{H};G,gx_{1})\right] \left( \ref{X1,gx1x2,XF21,g}\right) \\
+\lambda \left[ -B(x_{1}\otimes 1_{H};G,g)-B(gx_{1}x_{2}\otimes
1_{H};G,gx_{2})\right] =0
\end{gather*}

\begin{equation*}
2\beta _{1}B(x_{2}\otimes 1_{H};G,gx_{1}x_{2})-\lambda B(x_{1}\otimes
1_{H};G,gx_{1}x_{2})=0\left( \ref{X1,gx1x2,X1F21,gx1x2}\right)
\end{equation*}%
\begin{gather*}
+\lambda \left[ B(g\otimes 1_{H};1_{A},g)+B(x_{2}\otimes \
1_{H};1_{A},gx_{2})+B(x_{1}\otimes 1_{H};1_{A},gx_{1})+B(gx_{1}x_{2}\otimes
1_{H};1_{A},gx_{1}x_{2})\right] \\
-\gamma _{1}\left[ B(x_{2}\otimes 1_{H};G,g)-B(gx_{1}x_{2}\otimes
1_{H};G,gx_{1})\right] =0\left( \ref{X1,gx1x2,X1F31,g}\right)
\end{gather*}%
\begin{equation*}
0=0\left( \ref{X1,gx1x2,X1F31,x1}\right)
\end{equation*}

\begin{equation*}
0=0\left( \ref{X1,gx1x2,X1F31,x2}\right)
\end{equation*}

\begin{gather*}
\lambda B\left( g\otimes 1_{H};1_{A},gx_{1}x_{2}\right) +\left( \ref%
{X1,gx1x2,X1F31,gx1x2}\right) \\
-\gamma _{1}B(x_{2}\otimes 1_{H};G,gx_{1}x_{2})-\left[ +2B(x_{2}\otimes
1_{H};1_{A},gx_{2})+2B(gx_{1}x_{2}\otimes 1_{H};1_{A},gx_{1}x_{2})\right] =0
\end{gather*}

\begin{eqnarray*}
&&2\beta _{1}\left[
\begin{array}{c}
B(g\otimes 1_{H};1_{A},g)+B(x_{2}\otimes \ 1_{H};1_{A},gx_{2})+ \\
+B(x_{1}\otimes 1_{H};1_{A},gx_{1})+B(gx_{1}x_{2}\otimes
1_{H};1_{A},gx_{1}x_{2})%
\end{array}%
\right] +\left( \ref{X1,gx1x2,X1F41,g}\right) \\
&&-\gamma _{1}\left[ B(x_{1}\otimes 1_{H};G,g)+B(gx_{1}x_{2}\otimes
1_{H};G,gx_{2})\right] =0
\end{eqnarray*}%
\begin{eqnarray*}
&&2\beta _{1}B\left( g\otimes 1_{H};1_{A},gx_{1}x_{2}\right) \left( \ref%
{X1,gx1x2,X1F41,gx1x2}\right) \\
&&-\gamma _{1}B(x_{1}\otimes 1_{H};G,gx_{1}x_{2})-2B(x_{1}\otimes
1_{H};1_{A},gx_{2})=0
\end{eqnarray*}%
\begin{equation*}
0=0\left( \ref{X1,gx1x2,X1F51,1H}\right)
\end{equation*}%
\begin{gather*}
+\gamma _{1}\left[ B(g\otimes 1_{H};G,gx_{1})+B(x_{2}\otimes \
1_{H};G,gx_{1}x_{2})\right] +\left( \ref{X1,gx1x2,X1F51,gx1}\right) \\
-\left[
\begin{array}{c}
2B(g\otimes 1_{H};1_{A},g)+2B(x_{2}\otimes \ 1_{H};1_{A},gx_{2})+ \\
+2B(x_{1}\otimes 1_{H};1_{A},gx_{1})+2B(gx_{1}x_{2}\otimes
1_{H};1_{A},gx_{1}x_{2})%
\end{array}%
\right] =0
\end{gather*}%
\begin{equation*}
0=0\left( \ref{X1,gx1x2,X1F61,1H}\right)
\end{equation*}%
\begin{equation*}
0=0\left( \ref{X1,gx1x2,X1F61,x1x2}\right)
\end{equation*}%
\begin{gather*}
\lambda \left[ B(g\otimes 1_{H};G,gx_{1})+B(x_{2}\otimes \
1_{H};G,gx_{1}x_{2})\right] + \\
+\left[ 2B(x_{2}\otimes 1_{H};G,g)-2B(gx_{1}x_{2}\otimes 1_{H};G,gx_{1})%
\right] =0\left( \ref{X1,gx1x2,X1F61,gx1}\right)
\end{gather*}%
\begin{gather*}
\beta _{1}\left[ B(g\otimes 1_{H};G,gx_{1})+B(x_{2}\otimes \
1_{H};G,gx_{1}x_{2})\right] \left( \ref{X1,gx1x2,X1F71,gx1}\right) \\
+B(x_{1}\otimes 1_{H};G,g)+B(gx_{1}x_{2}\otimes 1_{H};G,gx_{2})=0.
\end{gather*}

\subsubsection{$X_{2}$}

\begin{equation*}
0=0\left( \ref{X2,g,X2F11,1H}\right)
\end{equation*}%
\begin{gather*}
2\beta _{2}B(g\otimes 1_{H};1_{A},gx_{1}x_{2})+\gamma _{2}B(g\otimes
1_{H};G,gx_{1})\left( \ref{X2,g,X2F11,gx1}\right) \\
+2B(x_{2}\otimes 1_{H};1_{A},gx_{1})=0
\end{gather*}

\begin{gather*}
\gamma _{2}B(g\otimes 1_{H};G,gx_{2})+\lambda B(g\otimes
1_{H};1_{A},gx_{1}x_{2})\left( \ref{X2,g,X2F11,gx2}\right) \\
+2B(x_{2}\otimes 1_{H};1_{A},gx_{2})+2B(g\otimes 1_{H};1_{A},g)=0.
\end{gather*}

\begin{equation*}
2\beta _{2}B(g\otimes 1_{H};G,gx_{2})+\lambda B(g\otimes
1_{H};G,gx_{1})+2B(x_{2}\otimes 1_{H};G,g)=0\left( \ref{X2,g,X2F21,g}\right)
\end{equation*}%
\begin{equation*}
0=0\left( \ref{X2,g,X2F21,x1}\right)
\end{equation*}%
\begin{equation*}
0=0\ref{X2,g,X2F21,x2}
\end{equation*}%
\begin{equation*}
2\beta _{2}B(g\otimes 1_{H};1_{H},gx_{1}x_{2})+\gamma _{2}B(g\otimes
1_{H};G,gx_{1})+2B(x_{2}\otimes 1_{H};1_{H},gx_{1})=0.\left( \ref%
{X2,g,X2F31,g}\right)
\end{equation*}%
\begin{equation*}
0=0\left( \ref{X2,g,X2F31,x2}\right)
\end{equation*}%
\begin{equation*}
B(x_{2}\otimes 1_{H};X_{2},x_{1})+B\left( x_{2}\otimes
1_{H};1_{A},x_{1}x_{2}\right) =0\left( \ref{X2,g,X2F41,x1}\right)
\end{equation*}%
\begin{equation*}
B\left( g\otimes 1_{H};1_{A},gx_{1}x_{2}\right) +B(g\otimes
1_{H},G,gx_{1})=0\left( \ref{X2,g,X2F61,gx2}\right)
\end{equation*}%
\begin{equation*}
0=0\left( \ref{X2,g,X2F71,1H}\right)
\end{equation*}%
\begin{equation*}
B(g\otimes 1_{H};G,gx_{1})+B(x_{2}\otimes 1_{H};G,gx_{1}x_{2})=0\left( \ref%
{X2,g,X2F71,gx1}\right)
\end{equation*}

\begin{eqnarray*}
&&-2\beta _{2}B(x_{1}\otimes 1_{H};1_{A},gx_{2})+\gamma _{2}B(x_{1}\otimes
1_{H};G,g)\left( \ref{X2,x1,X2F11,g}\right) \\
&&-\lambda \left[ B(g\otimes 1_{H};1_{A},g)+B(x_{1}\otimes
1_{H};1_{A},gx_{1})\right] -2B(gx_{1}x_{2}\otimes 1_{H};1_{A},g)=0
\end{eqnarray*}

\begin{equation*}
0=0\left( \ref{X2,x1,X2F11,x1}\right)
\end{equation*}%
\begin{equation*}
0=0.\left( \ref{X2,x1,X2F11,x2}\right)
\end{equation*}%
\begin{eqnarray*}
&&+\gamma _{2}B(x_{1}\otimes 1_{H};G,gx_{1}x_{2})-\lambda B(g\otimes
1_{H};1_{A},gx_{1}x_{2})\left( \ref{X2,x1,X2F11,gx1x2}\right) \\
&&-2B(gx_{1}x_{2}\otimes 1_{H};1_{A},gx_{1}x_{2})-2B\left( x_{1}\otimes
1_{H};1_{A},gx_{1}\right) =0
\end{eqnarray*}%
\begin{equation*}
0=0\left( \ref{X2,x1,X2F21,1H}\right)
\end{equation*}

\begin{equation*}
2\beta _{2}B(x_{1}\otimes 1_{H};G,gx_{1}x_{2})+\lambda B\left( g\otimes
1_{H};G,gx_{1}\right) +2B(x_{2}\otimes 1_{H};G,g)=0.\left( \ref%
{X2,x1,X2F21,gx1}\right)
\end{equation*}

\begin{eqnarray*}
&&\lambda \left[ B\left( g\otimes 1_{H};G,gx_{2}\right) -B(x_{1}\otimes
1_{H};G,gx_{1}x_{2})\right] \left( \ref{X2,x1,X2F21,gx2}\right) \\
&&+2-B(x_{1}\otimes 1_{H};G,g_{1})+2B\left( x_{1}\otimes 1_{H};G,g\right) =0.
\end{eqnarray*}

\begin{equation*}
0=0.\left( \ref{X2,x1,X2F31,1H}\right)
\end{equation*}%
\begin{gather*}
\gamma _{2}\left[ -B(g\otimes 1_{H};G,gx_{2})+B(x_{1}\otimes
1_{H};G,gx_{1}x_{2})\right] +\left( \ref{X2,x1,X2F31,gx2}\right) \\
-2B(x_{2}\otimes 1_{H};1_{A},gx_{2})-2B(gx_{1}x_{2}\otimes
1_{H};1_{A},gx_{1}x_{2}) \\
-2B(g\otimes 1_{H};1_{A},g)-2B(x_{1}\otimes 1_{H};1_{A},gx_{1})=0.
\end{gather*}%
\begin{equation*}
0=0.\left( \ref{X2,x1,X2F41,1H}\right)
\end{equation*}

\begin{eqnarray*}
&&-\gamma _{2}B(x_{1}\otimes 1_{H};G,gx_{1}x_{2})+\lambda B\left( g\otimes
1_{H};1_{A},gx_{1}x_{2}\right) \left( \ref{X2,x1,X2F41,gx1}\right) \\
&&+2B(x_{1}\otimes 1_{H};1_{A},gx_{1})+2B(gx_{1}x_{2}\otimes
1_{H};1_{A},gx_{1}x_{2})=0.
\end{eqnarray*}

\begin{eqnarray*}
&&\beta _{2}\left[ -B(g\otimes 1_{H};G,gx_{2})+B(x_{1}\otimes
1_{H};G,gx_{1}x_{2})\right] \left( \ref{X2,x1,X2F61,g}\right) \\
&&-B(x_{2}\otimes 1_{H};G,g)+B(gx_{1}x_{2}\otimes 1_{H};G,gx_{1})=0
\end{eqnarray*}%
\begin{equation*}
0=0\left( \ref{X2,x1,X2F61,x2}\right)
\end{equation*}

\begin{equation*}
B(x_{1}\otimes 1_{H};G,g)+B(gx_{1}x_{2}\otimes 1_{H};G,gx_{2})=0\left( \ref%
{X2,x1,X2F71,g}\right)
\end{equation*}%
\begin{equation*}
0=0\left( \ref{X2,x1,X2F81,1H}\right)
\end{equation*}%
\begin{eqnarray*}
2 &&\beta _{2}\left[ B(g\otimes 1_{H};1_{A},g)+B(x_{2}\otimes \
1_{H};1_{A},gx_{2})\right] \left( \ref{X2,x2,X2F11,g}\right) \\
&&-\gamma _{2}B(x_{2}\otimes 1_{H};G,g)+\lambda B(x_{2}\otimes
1_{H};1_{A},gx_{1})=0.
\end{eqnarray*}%
\begin{equation*}
0=0.\left( \ref{X2,x2,X2F11,x2}\right)
\end{equation*}

\begin{eqnarray*}
+ &&2\beta _{2}B\left( g\otimes 1_{H};1_{A},gx_{1}x_{2}\right) -\gamma
_{2}B(x_{2}\otimes 1_{H};G,gx_{1}x_{2})\left( \ref{X2,x2,X2F11,gx1x2}\right)
\\
&&+2B\left( x_{2}\otimes 1_{H};1_{A},gx_{1}\right) =0
\end{eqnarray*}

\begin{equation*}
\beta _{2}\left[ B(g\otimes 1_{H};G,gx_{1})+B(x_{2}\otimes \
1_{H};G,gx_{1}x_{2})\right] =0\left( \ref{X2,x2,X2F21,gx1}\right)
\end{equation*}%
\begin{eqnarray*}
&&2\beta _{2}B(g\otimes 1_{H};G,gx_{2})-\lambda B(x_{2}\otimes
1_{H};G,gx_{1}x_{2})\left( \ref{X2,x2,X2F21,gx2}\right) \\
&&+2B\left( x_{2}\otimes 1_{H};G,g\right) =0
\end{eqnarray*}

\begin{equation*}
0=0.\left( \ref{X2,x2,X2F41,1H}\right)
\end{equation*}

\begin{eqnarray*}
&&\gamma _{2}B(g\otimes 1_{H};G,gx_{2})-\lambda B\left( g\otimes
1_{H};1_{A},gx_{1}x_{2}\right) \left( \ref{X2,x2,X2F41,gx2}\right) \\
&&+2\left[ +B(g\otimes 1_{H};1_{A},g)+B(x_{2}\otimes \ 1_{H};1_{A},gx_{2})%
\right] =0.
\end{eqnarray*}

\begin{equation*}
\lambda \left[ B(g\otimes 1_{H};G,gx_{1})+B(x_{2}\otimes \
1_{H};G,gx_{1}x_{2})\right] =0\left( \ref{X2,x2,X2F71,g}\right)
\end{equation*}%
\begin{equation*}
\gamma _{2}B(x_{1}x_{2}\otimes 1_{H};G,g)+\lambda B(x_{1}x_{2}\otimes
1_{H};X_{1},g)=0.\left( \ref{X2,x1x2,X2F11,g}\right)
\end{equation*}

\begin{gather*}
2\beta _{2}B(x_{1}x_{2}\otimes 1_{H};X_{2},x_{1})+\gamma
_{2}B(x_{1}x_{2}\otimes 1_{H};G,x_{1})\left( \ref{X2,x1x2,X2F11,x1}\right) \\
+\lambda B(x_{1}x_{2}\otimes 1_{H};X_{1},x_{1})=0.
\end{gather*}

\begin{gather*}
2\beta _{2}B(x_{1}x_{2}\otimes 1_{H};X_{2},x_{2})+\gamma
_{2}B(x_{1}x_{2}\otimes 1_{H};G,x_{2})\left( \ref{X2,x1x2,X2F11,x2}\right) \\
+\lambda B(x_{1}x_{2}\otimes 1_{H};X_{1},x_{2})=0.
\end{gather*}%
\begin{eqnarray*}
&&\gamma _{2}B(x_{1}x_{2}\otimes 1_{H};G,gx_{1}x_{2})+\left( \ref%
{X2,x1x2,X2F11,gx1x2}\right) \\
&&+2B(x_{1}x_{2}\otimes 1_{H};1_{A},gx_{1})=0.
\end{eqnarray*}

\begin{equation*}
0=0.\left( \ref{X2,x1x2,X2F21,1H}\right)
\end{equation*}

\begin{equation*}
2B(x_{1}x_{2}\otimes 1_{H};G,g)-\lambda B(x_{1}x_{2}\otimes
1_{H};GX_{1},gx_{2})=0.\left( \ref{X2,x1x2,X2F21,gx2}\right)
\end{equation*}%
\begin{eqnarray*}
&&2\beta _{2}\left[
\begin{array}{c}
+1-B(x_{1}x_{2}\otimes 1_{H};1_{A},x_{1}x_{2}) \\
-B(x_{1}x_{2}\otimes 1_{H};X_{2},x_{1})+B(x_{1}x_{2}\otimes
1_{H};X_{1},x_{2})%
\end{array}%
\right] \left( \ref{X2,x1x2,X2F31,1H}\right) \\
&=&0.
\end{eqnarray*}%
\begin{equation*}
\gamma _{2}B(x_{1}x_{2}\otimes 1_{H};GX_{1},gx_{2})+2B(x_{1}x_{2}\otimes
1_{H};X_{1},g)=0.\left( \ref{X2,x1x2,X2F31,gx2}\right)
\end{equation*}

\begin{equation*}
\left[
\begin{array}{c}
+1-B(x_{1}x_{2}\otimes 1_{H};1_{A},x_{1}x_{2}) \\
-B(x_{1}x_{2}\otimes 1_{H};X_{2},x_{1})+B(x_{1}x_{2}\otimes
1_{H};X_{1},x_{2})%
\end{array}%
\right] =0.\left( \ref{X2,x1x2,X2F41,1H}\right)
\end{equation*}%
\begin{gather*}
2B(x_{1}x_{2}\otimes 1_{H};1_{A},gx_{2})+\left( \ref{X2,x1x2,X2F41,gx2}%
\right) \\
+2B(x_{1}x_{2}\otimes 1_{H};X_{2},g)=0.
\end{gather*}

\begin{gather*}
2B(x_{1}x_{2}\otimes 1_{H};X_{1},g)-\left( \ref{X2,x1x2,X2F51,g}\right) \\
\gamma _{2}\left[
\begin{array}{c}
-B(x_{1}x_{2}\otimes 1_{H};G,gx_{1}x_{2})+ \\
B(x_{1}x_{2}\otimes 1_{H};GX_{2},gx_{1})-B(x_{1}x_{2}\otimes
1_{H};GX_{1},gx_{2})%
\end{array}%
\right] =0.
\end{gather*}%
\begin{gather*}
2B(x_{1}x_{2}\otimes 1_{H};G,g)\left( \ref{X2,x1x2,X2F71,g}\right) \\
+\lambda \left[
\begin{array}{c}
-B(x_{1}x_{2}\otimes 1_{H};G,gx_{1}x_{2})+ \\
B(x_{1}x_{2}\otimes 1_{H};GX_{2},gx_{1})-B(x_{1}x_{2}\otimes
1_{H};GX_{1},gx_{2})%
\end{array}%
\right] =0.
\end{gather*}%
\begin{equation*}
B(x_{1}x_{2}\otimes 1_{H};G,gx_{1}x_{2})-B(x_{1}x_{2}\otimes
1_{H};GX_{2},gx_{1})=0.\left( \ref{X2,x1x2,X2F71,gx1x2}\right)
\end{equation*}%
\begin{gather*}
2\beta _{2}B(x_{1}x_{2}\otimes 1_{H};X_{2},x_{2})+\gamma
_{2}B(x_{1}x_{2}\otimes 1_{H};G,x_{2})\left( \ref{X2,gx1,X2F11,1H}\right) \\
+\lambda \left[ -1+B(x_{1}x_{2}\otimes
1_{H};1_{A},x_{1}x_{2})+B(x_{1}x_{2}\otimes 1_{H};X_{2},x_{1})\right] =0.
\end{gather*}

\begin{equation*}
B(x_{1}x_{2}\otimes 1_{H};1_{A},gx_{2})+B(x_{1}x_{2}\otimes
1_{H};X_{2},g)=0.\left( \ref{X2,gx1,X2F11,gx2}\right)
\end{equation*}%
\begin{equation*}
\lambda \left[ B(x_{1}x_{2}\otimes 1_{H};G,gx_{1}x_{2})-B(x_{1}x_{2}\otimes
1_{H};GX_{2},gx_{1})\right] =0.\left( \ref{X2,gx1,X2F21,g}\right)
\end{equation*}%
\begin{gather*}
2\beta _{2}\left[ -1+B(x_{1}x_{2}\otimes
1_{H};1_{A},x_{1}x_{2})-B(x_{1}x_{2}\otimes 1_{H};X_{1},x_{2})\right] \left( %
\ref{X2,gx2,X2F11,1H}\right) \\
+\gamma _{2}\left[ -B(x_{1}x_{2}\otimes 1_{H};G,x_{1})\right] + \\
-\lambda B(x_{1}x_{2}\otimes 1_{H};X_{1},x_{1})=0.
\end{gather*}

\begin{eqnarray*}
&&\gamma _{2}\left[ B(x_{1}x_{2}\otimes
1_{H};G,gx_{1}x_{2})+B(x_{1}x_{2}\otimes 1_{H};GX_{1},gx_{2})\right] \left( %
\ref{X2,gx2,X2F11,gx2}\right) \\
&&+ \\
&&+\left[ +2B(x_{1}x_{2}\otimes 1_{H};1_{A},gx_{1})+2B(x_{1}x_{2}\otimes
1_{H};X_{1},g)\right] =0.
\end{eqnarray*}%
\begin{eqnarray*}
&&2\left[ -B(x_{1}x_{2}\otimes 1_{H};1_{A},gx_{1})-B(x_{1}x_{2}\otimes
1_{H};X_{1},g)\right] \left( \ref{X2,gx2,X2F41,g}\right) \\
&&+\gamma _{2}\left[ -B(x_{1}x_{2}\otimes
1_{H};G,gx_{1}x_{2})-B(x_{1}x_{2}\otimes 1_{H};GX_{1},gx_{2})\right] \\
&=&0.
\end{eqnarray*}%
\begin{equation*}
0=0\left( \ref{X2,gx1x2,X2F11,1H}\right)
\end{equation*}%
\begin{equation*}
0=0\left( \ref{X2,gx1x2,X2F11,x1x2}\right)
\end{equation*}%
\begin{gather*}
2\beta _{2}\left[ B(x_{1}\otimes 1_{H};1_{A},gx_{1})+B(gx_{1}x_{2}\otimes
1_{H};1_{A},gx_{1}x_{2})\right] \left( \ref{X2,gx1x2,X2F11,gx1}\right) \\
+\gamma _{2}B(gx_{1}x_{2}\otimes 1_{H};G,gx_{1})-\lambda B(x_{2}\otimes
1_{H};1_{A},gx_{1})=0.
\end{gather*}

\begin{gather*}
2\beta _{2}B(x_{1}\otimes 1_{H};1_{A},gx_{2})+\gamma
_{2}B(gx_{1}x_{2}\otimes 1_{H};G,gx_{2})\left( \ref{X2,gx1x2,X2F11,gx2}%
\right) \\
+\lambda \left[ -B(x_{2}\otimes 1_{H};1_{A},gx_{2})-B(gx_{1}x_{2}\otimes
1_{H};1_{A},gx_{1}x_{2})\right] \\
+2B(gx_{1}x_{2}\otimes 1_{H};1_{A},g)=0.
\end{gather*}

\begin{gather*}
2\beta _{2}\left[ -B(x_{1}\otimes 1_{H};G,g)-B(gx_{1}x_{2}\otimes
1_{H};G,gx_{2})\right] \left( \ref{X2,gx1x2,X2F21,g}\right) \\
+\lambda \left[ B(x_{2}\otimes 1_{H};G,g)-B(gx_{1}x_{2}\otimes
1_{H};G,gx_{1})\right] =0.
\end{gather*}%
\begin{gather*}
-2\beta _{2}B(x_{1}\otimes 1_{H};G,gx_{1}x_{2})\left( \ref%
{X2,gx1x2,X2F21,gx1x2}\right) \\
+\lambda B(x_{2}\otimes 1_{H};G,gx_{1}x_{2})-2B(gx_{1}x_{2}\otimes
1_{H};G,gx_{1})=0.
\end{gather*}%
\begin{gather*}
2\beta _{2}\left[
\begin{array}{c}
B(g\otimes 1_{H};1_{A},g)+B(x_{2}\otimes \ 1_{H};1_{A},gx_{2}) \\
+B(x_{1}\otimes 1_{H};1_{A},gx_{1})+B(gx_{1}x_{2}\otimes
1_{H};1_{A},gx_{1}x_{2})%
\end{array}%
\right] \left( \ref{X2,gx1x2,X2F31,g}\right) \\
-\gamma _{2}\left[ B(x_{2}\otimes 1_{H};G,g)-B(gx_{1}x_{2}\otimes
1_{H};G,gx_{1})\right] =0.
\end{gather*}

\begin{eqnarray*}
&&\gamma _{2}\left[ B(x_{1}\otimes 1_{H};G,g)+B(gx_{1}x_{2}\otimes
1_{H};G,gx_{2})\right] \left( \ref{X2,gx1x2,X2F41,g}\right) \\
&&-\lambda \left[
\begin{array}{c}
B(g\otimes 1_{H};1_{A},g)+B(x_{2}\otimes \ 1_{H};1_{A},gx_{2}) \\
+B(x_{1}\otimes 1_{H};1_{A},gx_{1})+B(gx_{1}x_{2}\otimes
1_{H};1_{A},gx_{1}x_{2})%
\end{array}%
\right] =0.
\end{eqnarray*}

\begin{equation*}
0=0.\left( \ref{X2,gx1x2,X2F41,x1}\right)
\end{equation*}

\begin{equation*}
0=0.\left( \ref{X2,gx1x2,X2F51,1H}\right)
\end{equation*}

\begin{equation*}
B(x_{2}\otimes 1_{H};G,g)-B(gx_{1}x_{2}\otimes 1_{H};G,gx_{1})=0.\left( \ref%
{X2,gx1x2,X2F61,gx2}\right)
\end{equation*}%
\begin{equation*}
0=0.\left( \ref{X2,gx1x2,X2F71,1H}\right)
\end{equation*}

\subsection{LIST\ OF\ ALL\ MONOMIAL\ EQUALITIES}

\begin{equation*}
B(g\otimes 1_{H};1_{A},x_{1})=0.\left( \ref{G,g, GF2,x1}\right)
\end{equation*}%
\begin{equation*}
B(g\otimes 1_{H};1_{A},x_{2})=0.\left( \ref{G,g, GF2,x2}\right)
\end{equation*}%
\begin{equation*}
B(g\otimes 1_{H};G,1_{H})=0\left( \ref{X1,x1,X1F21,1}\right)
\end{equation*}%
\begin{equation*}
B(g\otimes 1_{H};G,x_{1}x_{2})=0.\left( \ref{G,x1,GF6,x1}\right)
\end{equation*}%
\begin{equation*}
B(g\otimes 1_{H};X_{2},1_{H})=0\left( \ref{G,x1,GF3,1H}\right)
\end{equation*}

\begin{equation*}
B(x_{1}\otimes 1_{H};1_{A},1_{H})=0.\left( \ref{G,x1, GF2,1H}\right)
\end{equation*}%
\begin{equation*}
B(x_{1}\otimes 1_{H};1_{A},x_{1}x_{2})=0.\left( \ref{X1,g, X1F41,x1}\right)
\end{equation*}%
\begin{equation*}
B(x_{1}\otimes 1_{H};G,x_{1})=0.\left( \ref{X1,g,X1F21,x1}\right)
\end{equation*}%
\begin{equation*}
B(x_{1}\otimes 1_{H};G,x_{2})=0.\left( \ref{X1,g,X1F21,x2}\right)
\end{equation*}%
\begin{equation*}
B(x_{2}\otimes 1_{H};1_{A},1_{H})=0\left( \ref{X2,g,X2F11,1H}\right)
\end{equation*}%
\begin{equation*}
B(x_{2}\otimes 1_{H};1_{A},x_{1}x_{2})=0\left( \ref{G,x2, GF2,x1x2}\right)
\end{equation*}%
\begin{equation*}
B(x_{2}\otimes 1_{H};G,x_{1})=0\left( \ref{X2,g,X2F21,x1}\right)
\end{equation*}%
\begin{equation*}
B(x_{2}\otimes 1_{H};G,x_{2})=0.\left( \ref{X2,x2,X2F11,x2}\right)
\end{equation*}%
\begin{equation*}
B(x_{1}x_{2}\otimes 1_{H};X_{1},gx_{1}x_{2})=0\text{ }.\left( \ref{G,x1x2,
GF4,gx1}\right)
\end{equation*}%
\begin{equation*}
B(x_{1}x_{2}\otimes 1_{H};X_{2},gx_{1}x_{2})=0\text{ }\left( \ref{G,x1x2,
GF3,gx2}\right)
\end{equation*}%
\begin{equation*}
B(x_{1}x_{2}\otimes 1_{H};GX_{1},1_{H})=0.\left( \ref{G,x1x2, GF3,1H}\right)
\end{equation*}%
\begin{equation*}
B(x_{1}x_{2}\otimes 1_{H};GX_{1},x_{1}x_{2})=0\text{ }\left( \ref{G,x1x2,
GF3,gx1}\right)
\end{equation*}%
\begin{equation*}
B(x_{1}x_{2}\otimes 1_{H};GX_{1},gx_{1})=0.\left( \ref{X2,x1x2,X2F21,gx1}%
\right)
\end{equation*}%
\begin{equation*}
B(x_{1}x_{2}\otimes 1_{H};GX_{2},1_{H})=0.\left( \ref{X1,x1x2,X1F21,1H}%
\right)
\end{equation*}%
\begin{equation*}
B(x_{1}x_{2}\otimes 1_{H};GX_{2},x_{1}x_{2})=0.\left( \ref%
{X2,x1x2,X2F21,x1x2}\right)
\end{equation*}%
\begin{equation*}
B(x_{1}x_{2}\otimes 1_{H};GX_{2},gx_{2})=0\left( \ref{X1,x1x2,X1F41,gx2}%
\right)
\end{equation*}

\begin{equation*}
B(gx_{2}\otimes 1_{H};GX_{1}X_{2},1_{H})=0.\left( \ref{X1,gx2,X171,1H}\right)
\end{equation*}%
\begin{equation*}
B(gx_{1}x_{2}\otimes 1_{H};1_{A},x_{1})=0\left( \ref{X2,x1,X2F11,x1}\right)
\end{equation*}%
\begin{equation*}
B(gx_{1}x_{2}\otimes 1_{H};1_{A},x_{2})=0.\left( \ref{X2,x1,X2F11,x2}\right)
\end{equation*}

\begin{equation*}
B(gx_{1}x_{2}\otimes 1_{H};G,1_{H})=0\left( \ref{X2,x1,X2F21,1H}\right)
\end{equation*}

\subsection{LIST\ OF\ BINOMIAL EQUALITIES\ in $\left( \protect\ref{LAE5}%
\right) $ with coefficients 1 or -1\label{BIN 1-1}}

We will delete the multiplying constants, if any and reorder them

\begin{equation*}
B(g\otimes 1_{H},G,gx_{1})+B(x_{2}\otimes 1_{H};G,gx_{1}x_{2})=0\left( \ref%
{G,x2, GF6,g}\right) .
\end{equation*}%
\begin{equation*}
B\left( g\otimes 1_{H};1_{A},gx_{1}x_{2}\right) +B(g\otimes
1_{H},G,gx_{1})=0\left( \ref{X2,g,X2F61,gx2}\right)
\end{equation*}%
\begin{equation*}
B(g\otimes 1_{H};G,gx_{1})=-B(x_{2}\otimes \ 1_{H};G,gx_{1}x_{2}).\left( \ref%
{G,x2, GF2,gx1}\right)
\end{equation*}%
\begin{equation*}
B(g\otimes 1_{H};G,gx_{2})=B(x_{1}\otimes 1_{H};G,gx_{1}x_{2}).\left( \ref%
{G,x1, GF3,x2}\right)
\end{equation*}

\begin{equation*}
B(x_{1}\otimes 1_{H};G,g)+B(gx_{1}x_{2}\otimes 1_{H};G,gx_{2})=0\left( \ref%
{X2,x1,X2F71,g}\right)
\end{equation*}%
\begin{equation*}
B(x_{2}\otimes 1_{H};G,g)-B(gx_{1}x_{2}\otimes 1_{H};G,gx_{1})=0.\left( \ref%
{X2,gx1x2,X2F61,gx2}\right)
\end{equation*}%
\begin{equation*}
B(x_{2}\otimes 1_{H};X_{2},x_{1})+B\left( x_{2}\otimes
1_{H};1_{A},x_{1}x_{2}\right) =0\left( \ref{X2,g,X2F41,x1}\right)
\end{equation*}%
\begin{equation*}
B(x_{1}x_{2}\otimes 1_{H};1_{A},gx_{1})+B(x_{1}x_{2}\otimes
1_{H};X_{1},g)=0\left( \ref{X1,gx2,X131,g}\right)
\end{equation*}%
\begin{equation*}
B(x_{1}x_{2}\otimes 1_{H};1_{A},gx_{2})+B(x_{1}x_{2}\otimes
1_{H};X_{2},g)=0.\left( \ref{X2,gx1,X2F11,gx2}\right)
\end{equation*}

\begin{equation*}
B(x_{1}x_{2}\otimes 1_{H};G,gx_{1}x_{2})+B(x_{1}x_{2}\otimes
1_{H};GX_{1},gx_{2})=0\left( \ref{X1,gx2,X1F11,gx2}\right)
\end{equation*}%
\begin{equation*}
-B(x_{1}x_{2}\otimes 1_{H};G,gx_{1}x_{2})+B(x_{1}x_{2}\otimes
1_{H};GX_{1},gx_{2})=0\left( \ref{X1,gx2,X141,g}\right)
\end{equation*}

\begin{equation*}
B(x_{1}x_{2}\otimes 1_{H};G,gx_{1}x_{2})-B(x_{1}x_{2}\otimes
1_{H};GX_{2},gx_{1})=0.\left( \ref{X2,x1x2,X2F71,gx1x2}\right)
\end{equation*}
Comparing $\left( \ref{X1,gx2,X1F11,gx2}\right) $ and $\left( \ref%
{X1,gx2,X141,g}\right) $and $\left( \ref{X2,x1x2,X2F71,gx1x2}\right) $we
deduce that%
\begin{eqnarray}
B(x_{1}x_{2}\otimes 1_{H};G,gx_{1}x_{2}) &=&0  \label{B1} \\
B(x_{1}x_{2}\otimes 1_{H};GX_{1},gx_{2}) &=&0  \label{B2} \\
B(x_{1}x_{2}\otimes 1_{H};GX_{2},gx_{1}) &=&0  \label{B3}
\end{eqnarray}

\subsubsection{List of binomial equalities with different coefficient}

These binomial equalities have different coefficient. We delete the $\left( %
\ref{B1}\right) ,\left( \ref{B2}\right) $ and $\left( \ref{B3}\right) $
above in the following equalities
\begin{equation*}
2\alpha B(g\otimes 1_{H};G,gx_{1})+\gamma _{2}B(g\otimes
1_{H};1_{A},gx_{1}x_{2})=0\left( \ref{G,g, GF1,gx1}\right)
\end{equation*}

\begin{equation*}
2\alpha B(g\otimes 1_{H};G,gx_{2})-\gamma _{1}B\left( g\otimes
1_{H};1_{A},gx_{1}x_{2}\right) =0.\left( \ref{G,g, GF1,gx2}\right)
\end{equation*}

\begin{equation*}
\gamma _{1}B\left( g\otimes 1_{H};G,gx_{1}\right) +\gamma _{2}B\left(
g\otimes 1_{H};G,gx_{2}\right) =0.\left( \ref{G,g, GF2,g}\right)
\end{equation*}

\begin{equation*}
-\gamma _{1}B(g\otimes 1_{H};1_{A},gx_{1}x_{2})+2\alpha B(x_{1}\otimes
1_{H};G,gx_{1}x_{2})=0.\left( \ref{G,x1, GF1,gx1x2}\right)
\end{equation*}%
\begin{equation*}
\gamma _{1}B\left( g\otimes 1_{H};G,gx_{1}\right) +\gamma _{2}B(x_{1}\otimes
1_{H};G,gx_{1}x_{2})=0.\left( \ref{G,x1, GF2,gx1}\right)
\end{equation*}

\begin{equation*}
-\gamma _{2}B(g\otimes 1_{H};1_{A},gx_{1}x_{2})-2\alpha B(g\otimes
1_{H};G,gx_{1})=0.\left( \ref{G,x1,GF3,gx1}\right)
\end{equation*}

\begin{equation*}
-\gamma _{1}B(g\otimes 1_{H};1_{A},gx_{1}x_{2})+2\alpha B(x_{1}\otimes
1_{H};G,gx_{1}x_{2})=0.\left( \ref{G,x1, GF4,gx1}\right)
\end{equation*}

\begin{equation*}
2\alpha B(x_{2}\otimes 1_{H};G,gx_{1}x_{2})-\gamma _{2}B\left( g\otimes
1_{H};1_{A},gx_{1}x_{2}\right) =0.\left( \ref{G,x2, GF1,gx1x2}\right)
\end{equation*}

\begin{equation*}
-\gamma _{1}B(x_{2}\otimes 1_{H};G,gx_{1}x_{2})+\gamma _{2}B(g\otimes
1_{H};G,gx_{2})=0.\left( \ref{G,x2, GF2,gx2}\right)
\end{equation*}

\begin{equation*}
-\gamma _{2}B(g\otimes 1_{H};1_{A},gx_{1}x_{2})-2\alpha B(x_{2}\otimes
1_{H};G,gx_{1}x_{2})=0.\left( \ref{G,x2, GF3,gx2}\right)
\end{equation*}

\begin{equation*}
\gamma _{1}B(g\otimes 1_{H};1_{A},gx_{1}x_{2})-2\alpha B(g\otimes
1_{H};G,gx_{2})=0.\left( \ref{G,x2, GF4,gx2}\right)
\end{equation*}

\begin{equation*}
\gamma _{1}B(x_{1}x_{2}\otimes 1_{H};X_{1},g)+\gamma _{2}B(x_{1}x_{2}\otimes
1_{H};X_{2},g)=0.\left( \ref{G,x1x2, GF1,g}\right)
\end{equation*}

\begin{equation*}
B(x_{1}x_{2}\otimes 1_{H};1_{A},gx_{1})=0.\left( \ref{G,x1x2, GF2,gx1}\right)
\end{equation*}%
\begin{equation*}
B(x_{1}x_{2}\otimes 1_{H};1_{A},gx_{2})=0\left( \ref{G,x1x2, GF2,gx2}\right)
\end{equation*}

\begin{equation*}
\gamma _{1}B(x_{2}\otimes 1_{H};G,gx_{1}x_{2})-\gamma _{2}B(x_{1}\otimes
1_{H};G,gx_{1}x_{2})=0.\left( \ref{G,gx1x2, GF2,gx1x2}\right)
\end{equation*}

\begin{equation*}
\gamma _{2}B\left( g\otimes 1_{H};1_{A},gx_{1}x_{2}\right) -2\alpha
B(x_{2}\otimes 1_{H};G,gx_{1}x_{2})=0.\left( \ref{G,gx1x2, GF3,gx1x2}\right)
\end{equation*}

\begin{equation*}
-\gamma _{1}B\left( g\otimes 1_{H};1_{A},gx_{1}x_{2}\right) +2\alpha
B(x_{1}\otimes 1_{H};G,gx_{1}x_{2})=0.\left( \ref{G,gx1x2, GF4,gx1x2}\right)
\end{equation*}

\begin{equation*}
\lambda B(x_{1}\otimes 1_{H};G,gx_{1}x_{2})=-2B\left( x_{1}\otimes
1_{H};G,g\right) .\left( \ref{X1,x1,X1F21,gx1}\right)
\end{equation*}

\begin{equation*}
\lambda B(x_{1}x_{2}\otimes 1_{H};X_{2},g)-\gamma _{1}B(x_{1}x_{2}\otimes
1_{H};G,g)=0\left( \ref{X1,x1x2,X1F11,g}\right)
\end{equation*}

\begin{equation*}
B(x_{1}x_{2}\otimes 1_{H};1_{A},gx_{2})=0.\left( \ref{X1,x1x2,X1F11,gx1x2}%
\right)
\end{equation*}%
\begin{equation*}
B(x_{1}x_{2}\otimes 1_{H};G,g)=0.\left( \ref{X1,x1x2,X1F21,gx1}\right)
\end{equation*}

\begin{equation*}
B(x_{1}x_{2}\otimes 1_{H};1_{A},gx_{2})=0.\left( \ref{X1,x1x2,X1F31,gx2}%
\right)
\end{equation*}

\begin{equation*}
B(x_{1}x_{2}\otimes 1_{H};X_{2},g)=0.\left( \ref{X1,x1x2,X1F41,gx1}\right)
\end{equation*}

\begin{equation*}
2\beta _{1}B(x_{2}\otimes 1_{H};G,gx_{1}x_{2})-\lambda B(x_{1}\otimes
1_{H};G,gx_{1}x_{2})=0\left( \ref{X1,gx1x2,X1F21,gx1x2}\right)
\end{equation*}

\begin{equation*}
\gamma _{2}B(x_{1}x_{2}\otimes 1_{H};G,g)+\lambda B(x_{1}x_{2}\otimes
1_{H};X_{1},g)=0.\left( \ref{X2,x1x2,X2F11,g}\right)
\end{equation*}

\begin{equation*}
B(x_{1}x_{2}\otimes 1_{H};1_{A},gx_{1})=0\left( \ref{X2,x1x2,X2F11,gx1x2}%
\right)
\end{equation*}

\begin{equation*}
2B(x_{1}x_{2}\otimes 1_{H};G,g)=0.\left( \ref{X2,x1x2,X2F21,gx2}\right)
\end{equation*}%
\begin{equation*}
+2B(x_{1}x_{2}\otimes 1_{H};X_{1},g)=0.\left( \ref{X2,x1x2,X2F31,gx2}\right)
\end{equation*}

\subsubsection{REPETITIONS}

\begin{equation*}
\lbrack -B(g\otimes 1_{H};G,gx_{2})+B(x_{1}\otimes
1_{H};G,gx_{1}x_{2})]=0.\left( \ref{G,x1, GF5,g}\right) \text{ this is }%
\left( \ref{G,x1, GF3,x2}\right)
\end{equation*}%
\begin{equation*}
\left[ B(g\otimes 1_{H};G,gx_{2})-B(x_{1}\otimes 1_{H};G,gx_{1}x_{2})\right]
=0.\left( \ref{G,gx1x2, GF5,gx2}\right) \text{ this is }\left( \ref{G,x1,
GF3,x2}\right)
\end{equation*}%
\begin{equation*}
B\left( g\otimes 1_{H};G,gx_{2}\right) -B(x_{1}\otimes
1_{H};G,gx_{1}x_{2})]=0.\left( \ref{G,x1,GF2,gx2}\right) \text{ this is }%
\left( \ref{G,x1, GF3,x2}\right)
\end{equation*}%
\begin{equation*}
B(g\otimes 1_{H};G,gx_{2})-B(x_{1}\otimes 1_{H};G,gx_{1}x_{2})=0.\left( \ref%
{G,gx1x2, GF6,gx1}\right) \text{ this is }\left( \ref{G,x1, GF3,x2}\right)
\end{equation*}%
\begin{equation*}
B(g\otimes 1_{H};G,gx_{2})-B(x_{1}\otimes 1_{H};G,gx_{1}x_{2})=0.\left( \ref%
{G,gx1x2, GF6,gx2}\right) \text{ this is }\left( \ref{G,x1, GF3,x2}\right)
\end{equation*}%
\begin{equation*}
B(g\otimes 1_{H};G,gx_{2})-B(x_{1}\otimes 1_{H};G,gx_{1}x_{2})=0.\left( \ref%
{G,gx1x2, GF7,gx2}\right) \text{ this is }\left( \ref{G,x1, GF3,x2}\right)
\end{equation*}%
\begin{equation*}
B\left( g\otimes 1_{H};G,gx_{2}\right) -B\left( x_{1}\otimes
1_{H};G,gx_{1}x_{2}\right) =0\left( \ref{X1,g,X1F21,gx1x2}\right) \text{
this is }\left( \ref{G,x1, GF3,x2}\right)
\end{equation*}

\begin{equation*}
B(g\otimes 1_{H};G,gx_{1})+B(x_{2}\otimes \ 1_{H};G,gx_{1}x_{2})=0\left( \ref%
{G,x2, GF5,g}\right) \text{ this is }.\left( \ref{G,x2, GF2,gx1}\right)
\end{equation*}%
\begin{equation*}
\left[ B(g\otimes 1_{H};G,gx_{1})+B(x_{2}\otimes \ 1_{H};G,gx_{1}x_{2})%
\right] =0\left( \ref{G,x2, GF4,gx1}\right) \text{ this is }\left( \ref%
{G,x2, GF2,gx1}\right)
\end{equation*}

\begin{equation*}
\lbrack B(g\otimes 1_{H},G,gx_{1})+B(x_{2}\otimes
1_{H};G,gx_{1}x_{2})=0.=0.\left( \ref{G,x2, GF6,g}\right) \text{ this is }%
\left( \ref{G,x2, GF2,gx1}\right)
\end{equation*}%
\begin{equation*}
B(g\otimes 1_{H};G,gx_{1})+B(x_{2}\otimes 1_{H};G,gx_{1}x_{2})=0\left( \ref%
{X2,g,X2F71,gx1}\right) \text{ this is }\left( \ref{G,x2, GF2,gx1}\right)
\end{equation*}%
\begin{equation*}
B(g\otimes 1_{H};G,gx_{1})+B(x_{2}\otimes 1_{H};G,gx_{1}x_{2})]=0.\left( \ref%
{G,x2; GF7,g}\right) \text{ this is }\left( \ref{G,x2, GF2,gx1}\right)
\end{equation*}%
\begin{equation*}
B(g\otimes 1_{H};G,gx_{1})+B(x_{2}\otimes \ 1_{H};G,gx_{1}x_{2})=0\left( \ref%
{X1,x2,X1F71,g}\right) \text{ this is }\left( \ref{G,x2, GF2,gx1}\right)
\end{equation*}%
\begin{equation*}
B(g\otimes 1_{H};G,gx_{1})+B(x_{2}\otimes \ 1_{H};G,gx_{1}x_{2})=0\left( \ref%
{X2,x2,X2F21,gx1}\right) \text{ this is }\left( \ref{G,x2, GF2,gx1}\right)
\end{equation*}%
\begin{equation*}
B(g\otimes 1_{H};G,gx_{1})+B(x_{2}\otimes \ 1_{H};G,gx_{1}x_{2})=0.\left( %
\ref{G,gx1x2, GF7,gx1}\right) \text{this is }\left( \ref{G,x2, GF2,gx1}%
\right)
\end{equation*}

\begin{equation*}
B(g\otimes 1_{H};G,gx_{1})+B(x_{2}\otimes \ 1_{H};G,gx_{1}x_{2})=0\left( \ref%
{X2,x2,X2F71,g}\right) \text{this is }\left( \ref{G,x2, GF2,gx1}\right)
\end{equation*}%
\begin{equation*}
B(g\otimes 1_{H};G,gx_{1})+B(x_{2}\otimes \ 1_{H};G,gx_{1}x_{2})=0\left( \ref%
{X2,x2,X2F21,gx1}\right) \text{this is }\left( \ref{G,x2, GF2,gx1}\right)
\end{equation*}%
\begin{equation*}
B(g\otimes 1_{H};G,gx_{1})+B(x_{2}\otimes \ 1_{H};G,gx_{1}x_{2})=0\left( \ref%
{X2,x2,X2F71,g}\right) \text{this is}\left( \ref{G,x2, GF2,gx1}\right)
\end{equation*}%
\begin{equation*}
\lbrack -B(g\otimes 1_{H};G,gx_{2})+B(x_{1}\otimes
1_{H};G,gx_{1}x_{2})]=0.\left( \ref{G,x1, GF5,g}\right) \text{ this is }%
\left( \ref{G,x1, GF3,x2}\right)
\end{equation*}%
\begin{equation*}
\left[ B(g\otimes 1_{H};G,gx_{2})-B(x_{1}\otimes 1_{H};G,gx_{1}x_{2})\right]
=0.\left( \ref{G,gx1x2, GF5,gx2}\right) \text{ this is }\left( \ref{G,x1,
GF3,x2}\right)
\end{equation*}%
\begin{equation*}
B\left( g\otimes 1_{H};G,gx_{2}\right) -B(x_{1}\otimes
1_{H};G,gx_{1}x_{2})]=0.\left( \ref{G,x1,GF2,gx2}\right) \text{ this is }%
\left( \ref{G,x1, GF3,x2}\right)
\end{equation*}%
\begin{equation*}
B(g\otimes 1_{H};G,gx_{2})-B(x_{1}\otimes 1_{H};G,gx_{1}x_{2})=0.\left( \ref%
{G,gx1x2, GF6,gx1}\right) \text{ this is }\left( \ref{G,x1, GF3,x2}\right)
\end{equation*}%
\begin{equation*}
B(g\otimes 1_{H};G,gx_{2})-B(x_{1}\otimes 1_{H};G,gx_{1}x_{2})=0.\left( \ref%
{G,gx1x2, GF6,gx2}\right) \text{ this is }\left( \ref{G,x1, GF3,x2}\right)
\end{equation*}%
\begin{equation*}
B(g\otimes 1_{H};G,gx_{2})-B(x_{1}\otimes 1_{H};G,gx_{1}x_{2})=0.\left( \ref%
{G,gx1x2, GF7,gx2}\right) \text{ this is }\left( \ref{G,x1, GF3,x2}\right)
\end{equation*}%
\begin{equation*}
B\left( g\otimes 1_{H};G,gx_{2}\right) -B\left( x_{1}\otimes
1_{H};G,gx_{1}x_{2}\right) =0\left( \ref{X1,g,X1F21,gx1x2}\right) \text{
this is }\left( \ref{G,x1, GF3,x2}\right)
\end{equation*}%
\begin{equation*}
B(g\otimes 1_{H};G,gx_{2})-B(x_{1}\otimes 1_{H};G,gx_{1}x_{2})=0.\left( \ref%
{X1,x2,X1F71,gx1x2}\right) \text{ this is }\left( \ref{G,x1, GF3,x2}\right)
\end{equation*}%
\begin{equation*}
-B\left( g\otimes 1_{H};G,gx_{2}\right) +B(x_{1}\otimes
1_{H};G,gx_{1}x_{2})=0\left( \ref{X1,x1,X1F31,gx2}\right) \text{ this is }%
\left( \ref{G,x1, GF3,x2}\right)
\end{equation*}%
\begin{equation*}
B(x_{1}x_{2}\otimes 1_{H};1_{A},gx_{1})+B(x_{1}x_{2}\otimes
1_{H};X_{1},g)=0.\left( \ref{X1,gx2,X1F11,gx1}\right) \text{ this is }\left( %
\ref{X1,gx2,X131,g}\right)
\end{equation*}%
\begin{equation*}
B(x_{1}x_{2}\otimes 1_{H};1_{A},gx_{1})+B(x_{1}x_{2}\otimes
1_{H};X_{1},g)=0\left( \ref{X1,x1x2,X1F31,gx1}\right) \text{ this is }\left( %
\ref{X1,gx2,X131,g}\right)
\end{equation*}%
\begin{equation*}
+B(x_{1}x_{2}\otimes 1_{H};G,gx_{1}x_{2})+B(x_{1}x_{2}\otimes
1_{H};GX_{1},gx_{2})=0.\left( \ref{X1,x1x2,X1F81,gx1}\right) \text{ this is }%
\left( \ref{X1,gx2,X1F11,gx2}\right)
\end{equation*}%
\begin{equation*}
B(x_{1}x_{2}\otimes 1_{H};G,gx_{1}x_{2})-B(x_{1}x_{2}\otimes
1_{H};GX_{2},gx_{1})=0.\left( \ref{X2,gx1,X2F21,g}\right) \text{ this is }%
\left( \ref{X2,x1x2,X2F71,gx1x2}\right)
\end{equation*}%
\begin{equation*}
B(x_{1}x_{2}\otimes 1_{H};G,gx_{1}x_{2})+B(x_{1}x_{2}\otimes
1_{H};GX_{1},gx_{2})=0\left( \ref{X1,x1x2,X1F61,gx1x2}\right) \text{ this is
}\left( \ref{X1,gx2,X1F11,gx2}\right)
\end{equation*}%
\begin{equation*}
B(x_{1}x_{2}\otimes 1_{H};1_{A},gx_{2})+B(x_{1}x_{2}\otimes
1_{H};X_{2},g)=0.\left( \ref{X2,x1x2,X2F41,gx2}\right) \text{ this is }%
\left( \ref{X2,gx1,X2F11,gx2}\right)
\end{equation*}

\subsection{LIST\ OF\ THE\ NEW\ MONOIDAL EQUALITIES\ WE\ GOT\ ABOVE}

\begin{equation*}
B(x_{1}x_{2}\otimes 1_{H};1_{A},gx_{1})=0.\left( \ref{G,x1x2, GF2,gx1}\right)
\end{equation*}%
\begin{equation*}
B(x_{1}x_{2}\otimes 1_{H};1_{A},gx_{1})=0\left( \ref{X2,x1x2,X2F11,gx1x2}%
\right) \text{this is }\left( \ref{G,x1x2, GF2,gx1}\right)
\end{equation*}%
\begin{equation*}
B(x_{1}x_{2}\otimes 1_{H};1_{A},gx_{2})=0\left( \ref{G,x1x2, GF2,gx2}\right)
\end{equation*}%
\begin{equation*}
B(x_{1}x_{2}\otimes 1_{H};1_{A},gx_{2})=0.\left( \ref{X1,x1x2,X1F31,gx2}%
\right) \left( \ref{G,x1x2, GF2,gx2}\right)
\end{equation*}%
\begin{equation*}
B(x_{1}x_{2}\otimes 1_{H};1_{A},gx_{2})=0.\left( \ref{X1,x1x2,X1F11,gx1x2}%
\right) \text{this is}\left( \ref{G,x1x2, GF2,gx2}\right)
\end{equation*}%
\begin{equation*}
B(x_{1}x_{2}\otimes 1_{H};G,g)=0.\left( \ref{X1,x1x2,X1F21,gx1}\right)
\end{equation*}%
\begin{equation*}
B(x_{1}x_{2}\otimes 1_{H};G,g)=0.\left( \ref{X2,x1x2,X2F21,gx2}\right) \text{
this is }\left( \ref{X1,x1x2,X1F21,gx1}\right)
\end{equation*}%
\begin{equation*}
B(x_{1}x_{2}\otimes 1_{H};G,gx_{1}x_{2})=0\left( \ref{B1}\right)
\end{equation*}

\begin{equation*}
B(x_{1}x_{2}\otimes 1_{H};X_{1},g)=0.\left( \ref{X2,x1x2,X2F31,gx2}\right)
\end{equation*}

\begin{equation*}
B(x_{1}x_{2}\otimes 1_{H};X_{2},g)=0.\left( \ref{X1,x1x2,X1F41,gx1}\right)
\end{equation*}%
\begin{equation*}
B(x_{1}x_{2}\otimes 1_{H};GX_{1},gx_{2})=0\left( \ref{B2}\right)
\end{equation*}%
\begin{equation*}
B(x_{1}x_{2}\otimes 1_{H};GX_{2},gx_{1})=0\left( \ref{B3}\right)
\end{equation*}

\subsubsection{LIST\ OF\ EQUALITIES\ ABOVE\ WITH\ NO\ REPETITION}

\begin{equation}
B(x_{1}x_{2}\otimes 1_{H};1_{A},gx_{1})=0.\left( \ref{G,x1x2, GF2,gx1}\right)
\label{x1x2ot1,1,gx1}
\end{equation}%
\begin{equation}
B(x_{1}x_{2}\otimes 1_{H};1_{A},gx_{2})=0\left( \ref{G,x1x2, GF2,gx2}\right)
\label{x1x2ot1,1,gx2}
\end{equation}%
\begin{equation}
B(x_{1}x_{2}\otimes 1_{H};G,g)=0.\left( \ref{X1,x1x2,X1F21,gx1}\right)
\label{x1x2ot1,G,g}
\end{equation}%
\begin{equation}
B(x_{1}x_{2}\otimes 1_{H};G,gx_{1}x_{2})=0\left( \ref{B1}\right)
\label{x1x2ot1,G.gx1x2}
\end{equation}

\begin{equation}
B(x_{1}x_{2}\otimes 1_{H};X_{1},g)=0.\left( \ref{X2,x1x2,X2F31,gx2}\right)
\label{x1x2ot1,X1,g}
\end{equation}

\begin{equation}
B(x_{1}x_{2}\otimes 1_{H};X_{2},g)=0.\left( \ref{X1,x1x2,X1F41,gx1}\right)
\label{x1x2ot1,X2,g}
\end{equation}%
\begin{equation}
B(x_{1}x_{2}\otimes 1_{H};GX_{1},gx_{2})=0\left( \ref{B2}\right)
\label{x1x2ot1,GX1,gx2}
\end{equation}%
\begin{equation}
B(x_{1}x_{2}\otimes 1_{H};GX_{2},gx_{1})=0\left( \ref{B3}\right)
\label{x1x2ot1,GX2,gx1}
\end{equation}

\subsection{Substituting these equalities above in the list of binomial
equalities}

Now we take this list and substitute them in a copy of binomials equalities
down here%
\begin{equation*}
2\alpha B(g\otimes 1_{H};G,gx_{1})+\gamma _{2}B(g\otimes
1_{H};1_{A},gx_{1}x_{2})=0\left( \ref{G,g, GF1,gx1}\right)
\end{equation*}

\begin{equation*}
2\alpha B(g\otimes 1_{H};G,gx_{2})-\gamma _{1}B\left( g\otimes
1_{H};1_{A},gx_{1}x_{2}\right) =0.\left( \ref{G,g, GF1,gx2}\right)
\end{equation*}

\begin{equation*}
\gamma _{1}B\left( g\otimes 1_{H};G,gx_{1}\right) +\gamma _{2}B\left(
g\otimes 1_{H};G,gx_{2}\right) =0.\left( \ref{G,g, GF2,g}\right)
\end{equation*}

\begin{equation*}
-\gamma _{1}B(g\otimes 1_{H};1_{A},gx_{1}x_{2})+2\alpha B(x_{1}\otimes
1_{H};G,gx_{1}x_{2})=0.\left( \ref{G,x1, GF1,gx1x2}\right)
\end{equation*}%
\begin{equation*}
\gamma _{1}B\left( g\otimes 1_{H};G,gx_{1}\right) +\gamma _{2}B(x_{1}\otimes
1_{H};G,gx_{1}x_{2})=0.\left( \ref{G,x1, GF2,gx1}\right)
\end{equation*}

\begin{equation*}
-\gamma _{2}B(g\otimes 1_{H};1_{A},gx_{1}x_{2})-2\alpha B(g\otimes
1_{H};G,gx_{1})=0.\left( \ref{G,x1,GF3,gx1}\right)
\end{equation*}

\begin{equation*}
-\gamma _{1}B(g\otimes 1_{H};1_{A},gx_{1}x_{2})+2\alpha B(x_{1}\otimes
1_{H};G,gx_{1}x_{2})=0.\left( \ref{G,x1, GF4,gx1}\right)
\end{equation*}

\begin{equation*}
2\alpha B(x_{2}\otimes 1_{H};G,gx_{1}x_{2})-\gamma _{2}B\left( g\otimes
1_{H};1_{A},gx_{1}x_{2}\right) =0.\left( \ref{G,x2, GF1,gx1x2}\right)
\end{equation*}

\begin{equation*}
-\gamma _{1}B(x_{2}\otimes 1_{H};G,gx_{1}x_{2})+\gamma _{2}B(g\otimes
1_{H};G,gx_{2})=0.\left( \ref{G,x2, GF2,gx2}\right)
\end{equation*}

\begin{equation*}
-\gamma _{2}B(g\otimes 1_{H};1_{A},gx_{1}x_{2})-2\alpha B(x_{2}\otimes
1_{H};G,gx_{1}x_{2})=0.\left( \ref{G,x2, GF3,gx2}\right)
\end{equation*}

\begin{equation*}
\gamma _{1}B(g\otimes 1_{H};1_{A},gx_{1}x_{2})-2\alpha B(g\otimes
1_{H};G,gx_{2})=0.\left( \ref{G,x2, GF4,gx2}\right)
\end{equation*}

\begin{equation*}
0=0.\left( \ref{G,x1x2, GF1,g}\right)
\end{equation*}

\begin{equation*}
\gamma _{1}B(x_{2}\otimes 1_{H};G,gx_{1}x_{2})-\gamma _{2}B(x_{1}\otimes
1_{H};G,gx_{1}x_{2})=0.\left( \ref{G,gx1x2, GF2,gx1x2}\right)
\end{equation*}

\begin{equation*}
\gamma _{2}B\left( g\otimes 1_{H};1_{A},gx_{1}x_{2}\right) -2\alpha
B(x_{2}\otimes 1_{H};G,gx_{1}x_{2})=0.\left( \ref{G,gx1x2, GF3,gx1x2}\right)
\end{equation*}

\begin{equation*}
-\gamma _{1}B\left( g\otimes 1_{H};1_{A},gx_{1}x_{2}\right) +2\alpha
B(x_{1}\otimes 1_{H};G,gx_{1}x_{2})=0.\left( \ref{G,gx1x2, GF4,gx1x2}\right)
\end{equation*}

\begin{equation*}
\lambda B(x_{1}\otimes 1_{H};G,gx_{1}x_{2})=-2B\left( x_{1}\otimes
1_{H};G,g\right) .\left( \ref{X1,x1,X1F21,gx1}\right)
\end{equation*}

\begin{equation*}
0=0\left( \ref{X1,x1x2,X1F11,g}\right)
\end{equation*}

\begin{equation*}
2\beta _{1}B(x_{2}\otimes 1_{H};G,gx_{1}x_{2})-\lambda B(x_{1}\otimes
1_{H};G,gx_{1}x_{2})=0\left( \ref{X1,gx1x2,X1F21,gx1x2}\right)
\end{equation*}

\begin{equation*}
0=0.\left( \ref{X2,x1x2,X2F11,g}\right)
\end{equation*}

Thus we got no new monomial equalities.

\subsection{COMPLETE\ LIST\ OF\ MONOMIAL\ EQUALITIES\ UP\ TO\ NOW}

\begin{equation*}
B(g\otimes 1_{H};1_{A},x_{1})=0.\left( \ref{G,g, GF2,x1}\right)
\end{equation*}%
\begin{equation*}
B(g\otimes 1_{H};1_{A},x_{2})=0\left( \ref{G,g, GF2,x2}\right)
\end{equation*}%
\begin{equation*}
B\left( g\otimes 1_{H};1_{A},x_{2}\right) =0.\left( \ref{X1,g, X1F11,1H}%
\right) \text{ this is }\left( \ref{G,g, GF2,x2}\right)
\end{equation*}%
\begin{equation*}
B\left( g\otimes 1_{H};1_{A},x_{2}\right) =0.\left( \ref{X1,x1,X1F31,1H}%
\right) \text{ this is }\left( \ref{G,g, GF2,x2}\right)
\end{equation*}%
\begin{equation*}
B(g\otimes 1_{H};G,1_{H})=0\left( \ref{X1,x1,X1F21,1}\right)
\end{equation*}%
\begin{equation*}
B\left( g\otimes 1_{H};G,x_{1}x_{2}\right) =0.\left( \ref{X1,x1,X1F21,x1x2}%
\right)
\end{equation*}%
\begin{equation*}
B(g\otimes 1_{H};G,x_{1}x_{2})=0.\left( \ref{G,x1,GF6,x1}\right)
\end{equation*}%
\begin{equation*}
B(g\otimes 1_{H};G,x_{1}x_{2})=0\left( \ref{X2,g,X2F71,1H}\right) \text{
this is }\left( \ref{G,x1,GF6,x1}\right)
\end{equation*}%
\begin{equation*}
B(g\otimes 1_{H};X_{2},1_{H})=0\left( \ref{G,x1,GF3,1H}\right) \text{this is
}B\left( g\otimes 1_{H};1_{A},x_{2}\right) \text{ and so it is }\left( \ref%
{G,g, GF2,x2}\right)
\end{equation*}%
\begin{equation*}
B(x_{1}\otimes 1_{H};1_{A},1_{H})=0.\left( \ref{G,x1, GF2,1H}\right)
\end{equation*}%
\begin{equation*}
B(x_{1}\otimes 1_{H};1_{A},1_{H})=0\left( \ref{X1,gx1x2,X1F11,1H}\right)
\text{ this is }\left( \ref{G,x1, GF2,1H}\right)
\end{equation*}%
\begin{equation*}
B(x_{1}\otimes 1_{H};1_{A},1_{H})=0\left( \ref{X1,gx1x2,X1F51,1H}\right)
\text{this is }\left( \ref{G,x1, GF2,1H}\right)
\end{equation*}%
\begin{equation*}
B(x_{1}\otimes 1_{H};1_{A},x_{1}x_{2})=0.\left( \ref{X1,g, X1F41,x1}\right)
\end{equation*}%
\begin{equation*}
B(x_{1}\otimes 1_{H};1_{A},x_{1}x_{2})=0\text{ }\left( \ref{G,x1, GF2,x1x2}%
\right) \text{ this is a repetition of }\left( \ref{X1,g, X1F41,x1}\right)
\end{equation*}%
\begin{equation*}
B\left( x_{1}\otimes 1_{H};1_{A},x_{1}x_{2}\right) =0\left( \ref%
{X1,x1,X1F31,x1x2}\right) \text{ this was already in }\left( \ref{LME1}%
\right)
\end{equation*}%
\begin{equation*}
B(x_{1}\otimes 1_{H};G,x_{1})=0.\left( \ref{X1,g,X1F21,x1}\right)
\end{equation*}%
\begin{equation*}
B(x_{1}\otimes 1_{H};G,x_{2})=0.\left( \ref{X1,g,X1F21,x2}\right)
\end{equation*}%
\begin{equation*}
B(x_{2}\otimes 1_{H};1_{A},1_{H})=0\left( \ref{X2,g,X2F11,1H}\right)
\end{equation*}%
\begin{equation*}
B(x_{2}\otimes 1_{H};1_{A},x_{1}x_{2})=0\left( \ref{G,x2, GF2,x1x2}\right)
\end{equation*}%
\begin{equation*}
B(x_{2}\otimes \ 1_{H};1_{A},x_{1}x_{2})=0.\left( \ref{G,x2, GF4,1H}\right)
\text{ this was already in }\left( \ref{LME1}\right)
\end{equation*}%
\begin{equation*}
B(x_{2}\otimes \ 1_{H};1_{A},x_{1}x_{2})=0\left( \ref{X1,gx1x2,X1F31,x1}%
\right) \text{ this was already in }\left( \ref{LME1}\right)
\end{equation*}%
\begin{equation*}
B(x_{2}\otimes 1_{H};1,x_{1}x_{2})=0\left( \ref{X2,g,X2F31,x2}\right) \text{
this is }\left( \ref{G,x2, GF2,x1x2}\right)
\end{equation*}%
\begin{equation*}
B(x_{2}\otimes 1_{H};G,x_{1})=0\left( \ref{X2,g,X2F21,x1}\right)
\end{equation*}%
\begin{equation*}
B(x_{2}\otimes \ 1_{H};G,x_{2})=0.\left( \ref{X1,x2,X1F21,1H}\right)
\end{equation*}%
\begin{equation*}
B(x_{2}\otimes \ 1_{H};G,x_{2})=0\left( \ref{X2,x1,X2F81,1H}\right) \text{
this is }\left( \ref{X1,x2,X1F21,1H}\right)
\end{equation*}%
\begin{equation*}
B(x_{2}\otimes 1_{H};G,x_{2})=0.\left( \ref{X2,x2,X2F11,x2}\right) \text{
this is }\left( \ref{X1,x2,X1F21,1H}\right)
\end{equation*}%
\begin{equation*}
B(x_{1}x_{2}\otimes 1_{H};1_{A},gx_{1})=0.\left( \ref{G,x1x2, GF2,gx1}\right)
\end{equation*}%
\begin{equation*}
B(x_{1}x_{2}\otimes 1_{H};1_{A},gx_{2})=0\left( \ref{G,x1x2, GF2,gx2}\right)
\end{equation*}%
\begin{equation*}
B(x_{1}x_{2}\otimes 1_{H};G,g)=0.\left( \ref{X1,x1x2,X1F21,gx1}\right)
\end{equation*}%
\begin{equation*}
B(x_{1}x_{2}\otimes 1_{H};G,gx_{1}x_{2})=0\left( \ref{B1}\right)
\end{equation*}%
\begin{equation*}
B(x_{1}x_{2}\otimes 1_{H};X_{1},g)=0.\left( \ref{X2,x1x2,X2F31,gx2}\right)
\end{equation*}%
\begin{equation*}
B(x_{1}x_{2}\otimes 1_{H};X_{1},gx_{1}x_{2})=0\text{ }.\left( \ref{G,x1x2,
GF4,gx1}\right)
\end{equation*}%
\begin{equation*}
B(x_{1}x_{2}\otimes 1_{H};X_{2},g)=0.\left( \ref{X1,x1x2,X1F41,gx1}\right)
\end{equation*}%
\begin{equation*}
B(x_{1}x_{2}\otimes 1_{H};X_{2},gx_{1}x_{2})=0\text{ }\left( \ref{G,x1x2,
GF3,gx2}\right)
\end{equation*}%
\begin{equation*}
B(x_{1}x_{2}\otimes 1_{H};GX_{1},1_{H})=0.\left( \ref{G,x1x2, GF3,1H}\right)
\end{equation*}%
\begin{equation*}
B(x_{1}x_{2}\otimes 1_{H};GX_{1},x_{1}x_{2})=0.\left( \ref{G,x1x2, GF3,x1x2}%
\right)
\end{equation*}%
\begin{equation*}
B(x_{1}x_{2}\otimes 1_{H};GX_{1},x_{1}x_{2})=0\text{ }\left( \ref{G,x1x2,
GF3,gx1}\right) \text{ this is a repetition of }\left( \ref{G,x1x2, GF3,x1x2}%
\right)
\end{equation*}%
\begin{equation*}
B(x_{1}x_{2}\otimes 1_{H};GX_{1},gx_{1})=0.\left( \ref{X2,x1x2,X2F21,gx1}%
\right)
\end{equation*}%
\begin{equation*}
B(x_{1}x_{2}\otimes 1_{H};GX_{1},gx_{2})=0\left( \ref{B2}\right)
\end{equation*}%
\begin{equation*}
B(x_{1}x_{2}\otimes 1_{H};GX_{2},1_{H})=0.\left( \ref{X1,x1x2,X1F21,1H}%
\right)
\end{equation*}%
\begin{equation*}
B(x_{1}x_{2}\otimes 1_{H};GX_{2},x_{1}x_{2})=0.\left( \ref%
{X2,x1x2,X2F21,x1x2}\right)
\end{equation*}%
\begin{equation*}
B(x_{1}x_{2}\otimes 1_{H};GX_{2},gx_{1})=0\left( \ref{B3}\right)
\end{equation*}%
\begin{equation*}
B(x_{1}x_{2}\otimes 1_{H};GX_{2},gx_{2})=0\left( \ref{X1,x1x2,X1F41,gx2}%
\right)
\end{equation*}%
\begin{equation*}
B(gx_{2}\otimes 1_{H};GX_{1}X_{2},1_{H})=0.\left( \ref{X1,gx2,X171,1H}\right)
\end{equation*}%
\begin{equation*}
B(gx_{1}x_{2}\otimes 1_{H};1_{A},x_{1})=0\left( \ref{X2,x1,X2F11,x1}\right)
\end{equation*}%
\begin{equation*}
B(gx_{1}x_{2}\otimes 1_{H};1_{A},x_{2})=0.\left( \ref{X2,x1,X2F11,x2}\right)
\end{equation*}%
\begin{equation*}
B(gx_{1}x_{2}\otimes 1_{H};G,1_{H})=0\left( \ref{X2,x1,X2F21,1H}\right)
\end{equation*}%
\begin{equation*}
B(gx_{1}x_{2}\otimes 1_{H};G,x_{1}x_{2})=0.\left( \ref{X1,x2,X1F21,x1x2}%
\right) \text{ this was already in }\left( \ref{LME1}\right)
\end{equation*}

\begin{equation*}
B(gx_{1}x_{2}\otimes 1_{H};G,x_{1}x_{2})=0.\text{we do not consider }\left( %
\ref{X2,x1x2,X2F21,x1x2}\right) \text{this as it is}\left( \ref%
{X1,x2,X1F21,x1x2}\right) \text{ }
\end{equation*}

\subsection{GETTING\ MORE\ MONOMIAL\ EQUALITIES}

We now take all equalities in \ref{BIN 1-1} .

The first one we rewrite as%
\begin{equation*}
B(g\otimes 1_{H},G,gx_{1})=-B\left( g\otimes 1_{H};1_{A},gx_{1}x_{2}\right)
\left( \ref{X2,g,X2F61,gx2}\right)
\end{equation*}%
and now we replace in all other equalities $B(g\otimes 1_{H},G,gx_{1})$ with
$-B\left( g\otimes 1_{H};1_{A},gx_{1}x_{2}\right) .$ Now we replace%
\begin{equation*}
B(x_{2}\otimes \ 1_{H};G,gx_{1}x_{2})=B(g\otimes
1_{H};1_{A},gx_{1}x_{2})\left( \ref{G,x2, GF2,gx1}\right)
\end{equation*}%
we go on this way%
\begin{equation*}
B(x_{1}\otimes 1_{H};G,gx_{1}x_{2})=B(g\otimes 1_{H};G,gx_{2})\left( \ref%
{G,x1, GF3,x2}\right)
\end{equation*}

\begin{equation*}
B(x_{1}\otimes 1_{H};G,g)+B(gx_{1}x_{2}\otimes 1_{H};G,gx_{2})=0\left( \ref%
{X2,x1,X2F71,g}\right)
\end{equation*}%
\begin{equation*}
B(x_{2}\otimes 1_{H};G,g)-B(gx_{1}x_{2}\otimes 1_{H};G,gx_{1})=0.\left( \ref%
{X2,gx1x2,X2F61,gx2}\right)
\end{equation*}%
\begin{equation*}
B(x_{2}\otimes 1_{H};X_{2},x_{1})+B\left( x_{2}\otimes
1_{H};1_{A},x_{1}x_{2}\right) =0\left( \ref{X2,g,X2F41,x1}\right)
\end{equation*}%
\begin{equation*}
B(x_{1}x_{2}\otimes 1_{H};1_{A},gx_{1})+B(x_{1}x_{2}\otimes
1_{H};X_{1},g)=0\left( \ref{X1,gx2,X131,g}\right)
\end{equation*}%
\begin{equation*}
B(x_{1}x_{2}\otimes 1_{H};1_{A},gx_{2})+B(x_{1}x_{2}\otimes
1_{H};X_{2},g)=0.\left( \ref{X2,gx1,X2F11,gx2}\right)
\end{equation*}

\begin{equation*}
B(x_{1}x_{2}\otimes 1_{H};G,gx_{1}x_{2})+B(x_{1}x_{2}\otimes
1_{H};GX_{1},gx_{2})=0\left( \ref{X1,gx2,X1F11,gx2}\right)
\end{equation*}%
\begin{equation*}
-B(x_{1}x_{2}\otimes 1_{H};G,gx_{1}x_{2})+B(x_{1}x_{2}\otimes
1_{H};GX_{1},gx_{2})=0\left( \ref{X1,gx2,X141,g}\right)
\end{equation*}

\begin{equation*}
B(x_{1}x_{2}\otimes 1_{H};G,gx_{1}x_{2})-B(x_{1}x_{2}\otimes
1_{H};GX_{2},gx_{1})=0.\left( \ref{X2,x1x2,X2F71,gx1x2}\right)
\end{equation*}

\begin{equation*}
-2\alpha B(g\otimes 1_{H};1_{A},gx_{1}x_{2})+\gamma _{2}B(g\otimes
1_{H};1_{A},gx_{1}x_{2})=0\left( \ref{G,g, GF1,gx1}\right)
\end{equation*}

\begin{equation*}
2\alpha B(g\otimes 1_{H};G,gx_{2})-\gamma _{1}B\left( g\otimes
1_{H};1_{A},gx_{1}x_{2}\right) =0.\left( \ref{G,g, GF1,gx2}\right)
\end{equation*}

\begin{equation*}
\gamma _{1}B\left( g\otimes 1_{H};G,gx_{1}\right) +\gamma _{2}B\left(
g\otimes 1_{H};G,gx_{2}\right) =0.\left( \ref{G,g, GF2,g}\right)
\end{equation*}

\begin{equation*}
-\gamma _{1}B(g\otimes 1_{H};1_{A},gx_{1}x_{2})+2\alpha B(g\otimes
1_{H};G,gx_{2})=0.\left( \ref{G,x1, GF1,gx1x2}\right) \text{this is }\left( %
\ref{G,g, GF1,gx2}\right)
\end{equation*}%
\begin{equation*}
\gamma _{1}B\left( g\otimes 1_{H};G,gx_{1}\right) +\gamma _{2}B(g\otimes
1_{H};G,gx_{2})=0.\left( \ref{G,x1, GF2,gx1}\right) \text{ this is }\left( %
\ref{G,g, GF2,g}\right)
\end{equation*}

\begin{equation*}
-\gamma _{2}B(g\otimes 1_{H};1_{A},gx_{1}x_{2})+2\alpha B(g\otimes
1_{H};1_{A},gx_{1}x_{2})=0.\left( \ref{G,x1,GF3,gx1}\right) \text{ this is }%
\left( \ref{G,g, GF1,gx1}\right)
\end{equation*}

\begin{equation*}
-\gamma _{1}B(g\otimes 1_{H};1_{A},gx_{1}x_{2})+2\alpha B(g\otimes
1_{H};G,gx_{2})=0.\left( \ref{G,x1, GF4,gx1}\right) \text{this is }\left( %
\ref{G,g, GF1,gx2}\right)
\end{equation*}

\begin{equation*}
2\alpha B(x_{2}\otimes 1_{H};G,gx_{1}x_{2})-\gamma _{2}B\left( g\otimes
1_{H};1_{A},gx_{1}x_{2}\right) =0.\left( \ref{G,x2, GF1,gx1x2}\right)
\end{equation*}

\begin{equation*}
-\gamma _{1}B(x_{2}\otimes 1_{H};G,gx_{1}x_{2})+\gamma _{2}B(g\otimes
1_{H};G,gx_{2})=0.\left( \ref{G,x2, GF2,gx2}\right)
\end{equation*}

From%
\begin{equation*}
-\gamma _{2}B(g\otimes 1_{H};1_{A},gx_{1}x_{2})-2\alpha B(x_{2}\otimes
1_{H};G,gx_{1}x_{2})=0.\left( \ref{G,x2, GF3,gx2}\right)
\end{equation*}%
and from%
\begin{equation*}
\gamma _{2}B\left( g\otimes 1_{H};1_{A},gx_{1}x_{2}\right) -2\alpha
B(x_{2}\otimes 1_{H};G,gx_{1}x_{2})=0.\left( \ref{G,gx1x2, GF3,gx1x2}\right)
\end{equation*}%
we get%
\begin{equation}
B\left( g\otimes 1_{H};1_{A},gx_{1}x_{2}\right) =0  \label{got1,1,gx1x2}
\end{equation}%
and%
\begin{equation}
B(x_{2}\otimes 1_{H};G,gx_{1}x_{2})=0  \label{x2ot1,G,gx1x2}
\end{equation}

Now from%
\begin{equation*}
\gamma _{1}B(g\otimes 1_{H};1_{A},gx_{1}x_{2})-2\alpha B(g\otimes
1_{H};G,gx_{2})=0.\left( \ref{G,x2, GF4,gx2}\right)
\end{equation*}%
we get%
\begin{equation}
B(g\otimes 1_{H};G,gx_{2})=0  \label{got1,G,gx2}
\end{equation}%
and from
\begin{equation*}
\gamma _{1}B\left( g\otimes 1_{H};G,gx_{1}\right) +\gamma _{2}B\left(
g\otimes 1_{H};G,gx_{2}\right) =0.\left( \ref{G,g, GF2,g}\right)
\end{equation*}%
we get%
\begin{equation}
B\left( g\otimes 1_{H};G,gx_{1}\right) =0  \label{got1,G,gx1}
\end{equation}

\begin{equation*}
\gamma _{1}B(x_{1}x_{2}\otimes 1_{H};X_{1},g)+\gamma _{2}B(x_{1}x_{2}\otimes
1_{H};X_{2},g)=0.\left( \ref{G,x1x2, GF1,g}\right)
\end{equation*}

\begin{equation*}
\gamma _{1}B(x_{2}\otimes 1_{H};G,gx_{1}x_{2})-\gamma _{2}B(g\otimes
1_{H};G,gx_{2})=0.\left( \ref{G,gx1x2, GF2,gx1x2}\right) \text{trivial}
\end{equation*}

\begin{equation*}
-\gamma _{1}B\left( g\otimes 1_{H};1_{A},gx_{1}x_{2}\right) +2\alpha
B(g\otimes 1_{H};G,gx_{2})=0.\left( \ref{G,gx1x2, GF4,gx1x2}\right) \text{%
trivial now}
\end{equation*}

From%
\begin{equation*}
\lambda B(g\otimes 1_{H};G,gx_{2})=-2B\left( x_{1}\otimes 1_{H};G,g\right)
.\left( \ref{X1,x1,X1F21,gx1}\right)
\end{equation*}%
we get%
\begin{equation}
B\left( x_{1}\otimes 1_{H};G,g\right) =0  \label{x1ot1,G,g}
\end{equation}

\begin{equation*}
\lambda B(x_{1}x_{2}\otimes 1_{H};X_{2},g)-\gamma _{1}B(x_{1}x_{2}\otimes
1_{H};G,g)=0\left( \ref{X1,x1x2,X1F11,g}\right) \text{trivial}
\end{equation*}%
\begin{equation*}
2\beta _{1}B(x_{2}\otimes 1_{H};G,gx_{1}x_{2})-\lambda B(g\otimes
1_{H};G,gx_{2})=0\left( \ref{X1,gx1x2,X1F21,gx1x2}\right) \text{trivial}
\end{equation*}

\begin{equation*}
\gamma _{2}B(x_{1}x_{2}\otimes 1_{H};G,g)+\lambda B(x_{1}x_{2}\otimes
1_{H};X_{1},g)=0.\left( \ref{X2,x1x2,X2F11,g}\right) \text{trivial\ref%
{x2ot1,G,gx1x2}}
\end{equation*}

From%
\begin{equation*}
\gamma _{1}B(g\otimes 1_{H};G,gx_{2})+2B(x_{1}\otimes
1_{H};1_{A},gx_{2})-2\beta _{1}B\left( g\otimes
1_{H};1_{A},gx_{1}x_{2}\right) =0.\left( \ref{X1,g, X1F11,gx2}\right)
\end{equation*}%
by using $\left( \ref{got1,G,gx2}\right) $ and $\left( \ref{got1,1,gx1x2}%
\right) ,$we get%
\begin{equation}
B(x_{1}\otimes 1_{H};1_{A},gx_{2})=0.  \label{x1ot1,1,gx2}
\end{equation}%
From%
\begin{equation*}
2\beta _{2}B(g\otimes 1_{H};G,gx_{2})+\lambda B(g\otimes
1_{H};G,gx_{1})+2B(x_{2}\otimes 1_{H};G,g)=0\left( \ref{X2,g,X2F21,g}\right)
\end{equation*}%
by using $\left( \ref{got1,G,gx2}\right) $ and $\left( \ref{got1,G,gx1}%
\right) ,$we get%
\begin{equation}
B(x_{2}\otimes 1_{H};G,g)=0  \label{x2ot1,G,g}
\end{equation}%
From%
\begin{eqnarray*}
3B(g\otimes 1_{H};1_{A},x_{1})-\gamma _{2}B(g\otimes
1_{H};G,x_{1}x_{2})\left( \ref{X2,g,X2F41,x1}\right) && \\
+B(x_{2}\otimes 1_{H};X_{2},x_{1})+B\left( x_{2}\otimes
1_{H};1_{A},x_{1}x_{2}\right) &=&0
\end{eqnarray*}%
by using $\left( \ref{G,g, GF2,x1}\right) ,\left( \ref{G,x1,GF6,x1}\right) $
and $\left( \ref{G,x2, GF2,x1x2}\right) ,$ we get
\begin{equation}
B(x_{2}\otimes 1_{H};X_{2},x_{1})=0  \label{x2ot1,X2,x1}
\end{equation}

From
\begin{gather*}
2\beta _{2}B(g\otimes 1_{H};1_{A},gx_{1}x_{2})+\gamma _{2}B(g\otimes
1_{H};G,gx_{1})\left( \ref{X2,g,X2F11,gx1}\right) \\
+2B(x_{2}\otimes 1_{H};1_{A},gx_{1})=0
\end{gather*}%
by using $\left( \ref{got1,1,gx1x2}\right) $ and $\left( \ref{got1,G,gx1}%
\right) ,$we get

\begin{equation}
B(x_{2}\otimes 1_{H};1_{A},gx_{1})=0  \label{x2ot1,1,gx1}
\end{equation}%
From
\begin{equation*}
-\lambda B\left( g\otimes 1_{H};1_{A},gx_{1}x_{2}\right) -\gamma _{1}B\left(
g\otimes 1_{H};G,gx_{1}\right) -2\left[ B(g\otimes
1_{H};1_{A},g)+B(x_{1}\otimes 1_{H};1_{A},gx_{1})\right] =0\left( \ref%
{X1,x1,X1F31,gx1}\right)
\end{equation*}%
By using $\left( \ref{got1,1,gx1x2}\right) $ and$\left( \ref{got1,G,gx1}%
\right) $ we get%
\begin{equation}
B(g\otimes 1_{H};1_{A},g)+B(x_{1}\otimes 1_{H};1_{A},gx_{1})=0  \label{B4}
\end{equation}%
Now from%
\begin{eqnarray*}
&&-2\beta _{2}B(x_{1}\otimes 1_{H};1_{A},gx_{2})+\gamma _{2}B(x_{1}\otimes
1_{H};G,g)\left( \ref{X2,x1,X2F11,g}\right) \\
&&-\lambda \left[ B(g\otimes 1_{H};1_{A},g)+B(x_{1}\otimes
1_{H};1_{A},gx_{1})\right] -2B(gx_{1}x_{2}\otimes 1_{H};1_{A},g)=0
\end{eqnarray*}%
By using $\left( \ref{x1ot1,1,gx2}\right) $, $\left( \ref{x1ot1,G,g}\right) $
and $\left( \ref{B4}\right) $ we get%
\begin{equation}
B(gx_{1}x_{2}\otimes 1_{H};1_{A},g)=0  \label{gx1x2ot1,1,g}
\end{equation}

From%
\begin{gather*}
-2\left[ -B(x_{2}\otimes 1_{H};1_{A},1_{H})+B(gx_{1}x_{2}\otimes
1_{H};1_{A},x_{1})\right] +\left( \ref{G,gx1x2, GF6,1H}\right) \\
\gamma _{2}\left[ B(g\otimes 1_{H};G,1_{H})+B(x_{2}\otimes \
1_{H};G,x_{2})+B(x_{1}\otimes 1_{H};G,x_{1})+B(gx_{1}x_{2}\otimes
1_{H};G,x_{1}x_{2})\right] \\
=0
\end{gather*}%
by using $\left( \ref{x2ot1,1,1}\right) ,\left( \ref{gx1x2ot1,1,x1}\right)
,\left( \ref{got1,G,1}\right) ,\left( \ref{x2ot1,G,x2}\right) $ and $\left( %
\ref{x1ot1,G,x1}\right) $ we get%
\begin{equation}
B(gx_{1}x_{2}\otimes 1_{H};G,x_{1}x_{2})=0  \label{gx1x2ot1,G,x1x2}
\end{equation}

From
\begin{equation*}
\gamma _{1}B\left( g\otimes 1_{H};G,gx_{1}\right) +\gamma _{2}B(x_{1}\otimes
1_{H};G,gx_{1}x_{2})=0\left( \ref{G,x1, GF2,gx1}\right)
\end{equation*}%
by using $\left( \ref{got1,G,gx1}\right) ,$ we get%
\begin{equation}
B(x_{1}\otimes 1_{H};G,gx_{1}x_{2})=0  \label{x1ot1,G,gx1x2}
\end{equation}%
From
\begin{equation*}
B(x_{1}\otimes 1_{H};G,g)+B(gx_{1}x_{2}\otimes 1_{H};G,gx_{2})=0\left( \ref%
{X2,x1,X2F71,g}\right)
\end{equation*}%
by $\left( \ref{x1ot1,G,g}\right) ,$we get%
\begin{equation}
B(gx_{1}x_{2}\otimes 1_{H};G,gx_{2})=0  \label{gx1x2ot1,G,gx2}
\end{equation}%
From
\begin{gather*}
2\beta _{2}\left[ -B(x_{1}\otimes 1_{H};G,g)-B(gx_{1}x_{2}\otimes
1_{H};G,gx_{2})\right] \left( \ref{X2,gx1x2,X2F21,g}\right) \\
+\lambda \left[ B(x_{2}\otimes 1_{H};G,g)-B(gx_{1}x_{2}\otimes
1_{H};G,gx_{1})\right] =0.
\end{gather*}%
by using $\left( \ref{x1ot1,G,g}\right) ,\left( \ref{gx1x2ot1,G,gx2}\right) $
and $\left( \ref{x2ot1,G,g}\right) ,$we get%
\begin{equation}
B(gx_{1}x_{2}\otimes 1_{H};G,gx_{1})=0  \label{gx1x2ot1,G,gx1}
\end{equation}

By using%
\begin{equation*}
B(g\otimes 1_{H};1_{A},g)+B(x_{1}\otimes 1_{H};1_{A},gx_{1})=0.\left( \ref%
{G,x1, GF1,g}\right)
\end{equation*}%
from
\begin{eqnarray*}
&&\lambda B(g\otimes 1_{H};X_{2},gx_{1})+2B(x_{1}\otimes
1_{H};1_{A},gx_{1})+\left( \ref{X1,g, X1F11,gx1}\right) \\
&&+2B\left( g\otimes 1_{H};1_{A},g\right) =0
\end{eqnarray*}%
we get%
\begin{equation}
B(g\otimes 1_{H};X_{2},gx_{1})=0  \label{got1,X2,gx1}
\end{equation}%
and $\left( \ref{X1,g, X1F11,gx1}\right) $ becomes trivial in view of $%
\left( \ref{got1,X2,gx1}\right) $ and $\left( \ref{G,x1, GF1,g}\right) .$

\subsection{LIST\ OF\ ALL\ MONOMIAL\ EQUALITIES\label{LME}}

\begin{equation*}
B(g\otimes 1_{H};1_{A},x_{1})=0\left( \ref{got1,1,x1}\right)
\end{equation*}%
\begin{equation*}
B(g\otimes 1_{H};1_{A},x_{2})=0\left( \ref{got1,1,x2}\right)
\end{equation*}%
\begin{equation*}
B\left( g\otimes 1_{H};1_{A},gx_{1}x_{2}\right) =0\left( \ref{got1,1,gx1x2}%
\right)
\end{equation*}%
\begin{equation*}
B(g\otimes 1_{H};G,1_{H})=0\left( \ref{got1,G,1}\right)
\end{equation*}

\begin{equation*}
B\left( g\otimes 1_{H};G,x_{1}x_{2}\right) =0.\left( \ref{got1,G,x1x2}\right)
\end{equation*}%
\begin{equation*}
B(g\otimes 1_{H};G,gx_{1})=0\left( \ref{got1,G,gx1}\right)
\end{equation*}%
\begin{equation*}
B(g\otimes 1_{H};G,gx_{2})=0\left( \ref{got1,G,gx2}\right)
\end{equation*}%
\begin{equation*}
B(g\otimes 1_{H};X_{2},1_{H})=0\left( \ref{got1,X2,1}\right)
\end{equation*}%
\begin{equation*}
B(g\otimes 1_{H};X_{2},gx_{1})=0\left( \ref{got1,X2,gx1}\right)
\end{equation*}%
\begin{equation*}
B(x_{1}\otimes 1_{H};1_{A},1_{H})=0.\left( \ref{x1ot1,1,1}\right)
\end{equation*}%
\begin{equation*}
B(x_{1}\otimes 1_{H};1_{A},x_{1}x_{2})=0.\left( \ref{x1ot1,1,x1x2}\right)
\end{equation*}%
\begin{equation*}
B(x_{1}\otimes 1_{H};1_{A},gx_{2})=0\left( \ref{x1ot1,1,gx2}\right)
\end{equation*}%
\begin{equation*}
B\left( x_{1}\otimes 1_{H};G,g\right) =0\left( \ref{x1ot1,G,g}\right)
\end{equation*}%
\begin{equation*}
B(x_{1}\otimes 1_{H};G,x_{1})=0.\left( \ref{x1ot1,G,x1}\right)
\end{equation*}%
\begin{equation*}
B(x_{1}\otimes 1_{H};G,x_{2})=0.\left( \ref{x1ot1,G,x2}\right)
\end{equation*}%
\begin{equation*}
B(x_{1}\otimes 1_{H};G,gx_{1}x_{2})=0\left( \ref{x1ot1,G,gx1x2}\right)
\end{equation*}%
\begin{equation*}
B(x_{2}\otimes 1_{H};1_{A},1_{H})=0\left( \ref{x2ot1,1,1}\right)
\end{equation*}%
\begin{equation*}
B(x_{2}\otimes 1_{H};1_{A},gx_{1})=0\left( \ref{x2ot1,1,gx1}\right)
\end{equation*}%
\begin{equation*}
B(x_{2}\otimes 1_{H};1_{A},x_{1}x_{2})=0\left( \ref{x2ot1,1,x1x2}\right)
\end{equation*}%
\begin{equation*}
B(x_{2}\otimes 1_{H};G,g)=0\left( \ref{x2ot1,G,g}\right)
\end{equation*}%
\begin{equation*}
B(x_{2}\otimes 1_{H};G,x_{1})=0\left( \ref{x2ot1,G,x1}\right)
\end{equation*}%
\begin{equation*}
B(x_{2}\otimes 1_{H};G,x_{2})=0\left( \ref{x2ot1,G,x2}\right)
\end{equation*}%
\begin{equation*}
B(x_{2}\otimes 1_{H};G,gx_{1}x_{2})=0\left( \ref{x2ot1,G,gx1x2}\right)
\end{equation*}%
\begin{equation*}
B(x_{2}\otimes 1_{H};X_{2},x_{1})=0\left( \ref{x2ot1,X2,x1}\right)
\end{equation*}%
\begin{equation*}
B(x_{1}x_{2}\otimes 1_{H};1_{A},gx_{1})=0.\left( \ref{x1x2ot1,1,gx1}\right)
\end{equation*}%
\begin{equation*}
B(x_{1}x_{2}\otimes 1_{H};1_{A},gx_{2})=0\left( \ref{x1x2ot1,1,gx2}\right)
\end{equation*}%
\begin{equation*}
B(x_{1}x_{2}\otimes 1_{H};G,g)=0\left( \ref{x1x2ot1,G,g}\right)
\end{equation*}%
\begin{equation*}
B(x_{1}x_{2}\otimes 1_{H};G,gx_{1}x_{2})=0\left( \ref{x1x2ot1,G.gx1x2}\right)
\end{equation*}%
\begin{equation*}
B(x_{1}x_{2}\otimes 1_{H};X_{1},g)=0\left( \ref{x1x2ot1,X1,g}\right)
\end{equation*}%
\begin{equation*}
B(x_{1}x_{2}\otimes 1_{H};X_{1},gx_{1}x_{2})=0.\left( \ref{x1x2ot1,X1,gx1x2}%
\right)
\end{equation*}%
\begin{equation*}
B(x_{1}x_{2}\otimes 1_{H};X_{2},g)=0\left( \ref{x1x2ot1,X2,g}\right)
\end{equation*}%
\begin{equation*}
B(x_{1}x_{2}\otimes 1_{H};X_{2},gx_{1}x_{2})=0\left( \ref{x1x2ot1,X2,gx1x2}%
\right)
\end{equation*}%
\begin{equation*}
B(x_{1}x_{2}\otimes 1_{H};GX_{1},1_{H})=0.\left( \ref{x1x2ot1,GX1,1}\right)
\end{equation*}%
\begin{equation*}
B(x_{1}x_{2}\otimes 1_{H};GX_{1},x_{1}x_{2})=0\left( \ref{x1x2ot1,GX1,x1x2}%
\right)
\end{equation*}%
\begin{equation*}
B(x_{1}x_{2}\otimes 1_{H};GX_{1},gx_{1})=0.\left( \ref{x1x2ot1,GX1,gx1}%
\right)
\end{equation*}%
\begin{equation*}
B(x_{1}x_{2}\otimes 1_{H};GX_{1},gx_{2})=0\left( \ref{x1x2ot1,GX1,gx2}\right)
\end{equation*}%
\begin{equation*}
B(x_{1}x_{2}\otimes 1_{H};GX_{2},1_{H})=0.\left( \ref{x1x2ot1,GX2,1}\right)
\end{equation*}%
\begin{equation*}
B(x_{1}x_{2}\otimes 1_{H};GX_{2},x_{1}x_{2})=0.\left( \ref{x1x2ot1,GX2,x1x2}%
\right)
\end{equation*}%
\begin{equation*}
B(x_{1}x_{2}\otimes 1_{H};GX_{2},gx_{1})=0\left( \ref{x1x2ot1,GX2,gx1}\right)
\end{equation*}%
\begin{equation*}
B(x_{1}x_{2}\otimes 1_{H};GX_{2},gx_{2})=0\left( \ref{x1x2ot1;GX2,gx2}\right)
\end{equation*}%
\begin{equation*}
B(gx_{2}\otimes 1_{H};GX_{1}X_{2},1_{H})=0\left( \ref{gx2ot1,GX1X2,1}\right)
\end{equation*}%
\begin{equation*}
B(gx_{1}x_{2}\otimes 1_{H};1_{A},g)=0\left( \ref{gx1x2ot1,1,g}\right)
\end{equation*}%
\begin{equation*}
B(gx_{1}x_{2}\otimes 1_{H};1_{A},x_{1})=0\left( \ref{gx1x2ot1,1,x1}\right)
\end{equation*}%
\begin{equation*}
B(gx_{1}x_{2}\otimes 1_{H};1_{A},x_{2})=0.\left( \ref{gx1x2ot1,1,x2}\right)
\end{equation*}%
\begin{equation*}
B(gx_{1}x_{2}\otimes 1_{H};G,1_{H})=0\left( \ref{gx1x2ot1,G,1}\right)
\end{equation*}%
\begin{equation*}
B(gx_{1}x_{2}\otimes 1_{H};G,x_{1}x_{2})=0\left( \ref{gx1x2ot1,G,x1x2}\right)
\end{equation*}%
\begin{equation*}
B(gx_{1}x_{2}\otimes 1_{H};G,gx_{1})=0\left( \ref{gx1x2ot1,G,gx1}\right)
\end{equation*}%
\begin{equation*}
B(gx_{1}x_{2}\otimes 1_{H};G,gx_{2})=0\left( \ref{gx1x2ot1,G,gx2}\right)
\end{equation*}

\subsection{LIST\ OF\ ALL\ EQUALITIES 6\label{LAE6}}

Now we take the list of all equalities \ref{LAE5} and substitute inside it
the equalities above and reorder them.

\begin{equation*}
B(g\otimes 1_{H};1_{A},g)+B(x_{1}\otimes 1_{H};1_{A},gx_{1})=0.\left( \ref%
{G,x1, GF1,g}\right)
\end{equation*}

\begin{eqnarray*}
B(g\otimes 1_{H};1_{A},g)+B(x_{1}\otimes 1_{H};1_{A},gx_{1}) &=&0.\left( \ref%
{X1,g, X1F31,g}\right) \\
&&\text{ this is }\left( \ref{G,x1, GF1,g}\right) \text{ }
\end{eqnarray*}%
\begin{eqnarray*}
B(g\otimes 1_{H};1_{A},g)+B(x_{1}\otimes 1_{H};1_{A},gx_{1}) &=&0\left( \ref%
{X1,x1,X1F11,g}\right) \text{ } \\
&&\text{this is }\left( \ref{G,x1, GF1,g}\right)
\end{eqnarray*}%
\begin{eqnarray*}
B(g\otimes 1_{H};1_{A},g)+B(x_{1}\otimes 1_{H};1_{A},gx_{1}) &=&0\left( \ref%
{X1,x1,X1F31,gx1}\right) \\
&&\text{this is }\left( \ref{G,x1, GF1,g}\right)
\end{eqnarray*}

\begin{eqnarray*}
&&B(g\otimes 1_{H};1_{A},g)+B(x_{1}\otimes 1_{H};1_{A},gx_{1})=0\left( \ref%
{X2,x1,X2F11,g}\right) \\
&&\text{this is }\left( \ref{G,x1, GF1,g}\right)
\end{eqnarray*}%
\begin{equation*}
B(g\otimes 1_{H};1_{A},g)+B(x_{2}\otimes \ 1_{H};1_{A},gx_{2})=0.\left( \ref%
{G,x2, GF1,g}\right)
\end{equation*}%
\begin{eqnarray*}
B(g\otimes 1_{H};1_{A},g)+B(x_{2}\otimes \ 1_{H};1_{A},gx_{2}) &=&0\left( %
\ref{X1,x2,X1F11,g}\right) \\
&&\text{ this is }\left( \ref{G,x2, GF1,g}\right)
\end{eqnarray*}%
\begin{gather*}
B(g\otimes 1_{H};1_{A},g)+B(x_{2}\otimes 1_{H};1_{A},gx_{2})=0\left( \ref%
{X2,g,X2F11,gx2}\right) \\
\text{ this is }\left( \ref{G,x2, GF1,g}\right)
\end{gather*}%
\begin{eqnarray*}
B(g\otimes 1_{H};1_{A},g)+B(x_{2}\otimes \ 1_{H};1_{A},gx_{2}) &=&0.\left( %
\ref{X2,x2,X2F11,g}\right) \\
&&\text{ this is }\left( \ref{G,x2, GF1,g}\right)
\end{eqnarray*}%
\begin{eqnarray*}
B(g\otimes 1_{H};1_{A},g)+B(x_{2}\otimes \ 1_{H};1_{A},gx_{2}) &=&0\left( %
\ref{X2,x2,X2F41,gx2}\right) \\
&&\text{ this is }\left( \ref{G,x2, GF1,g}\right)
\end{eqnarray*}%
\begin{equation*}
\beta _{1}B(g\otimes 1_{H};1_{A},gx_{1}x_{2})=B\left( x_{1}\otimes
1_{H};1_{A},gx_{2}\right) .\left( \ref{X1,x1,X1F11,gx1x2}\right)
\end{equation*}%
\begin{equation*}
B(x_{1}\otimes 1_{H};1_{A},gx_{1})+B(gx_{1}x_{2}\otimes
1_{H};1_{A},gx_{1}x_{2})=0.\left( \ref{G,gx1x2, GF1,gx1}\right)
\end{equation*}%
\begin{eqnarray*}
B(x_{1}\otimes 1_{H};1_{A},gx_{1})+B(gx_{1}x_{2}\otimes
1_{H};1_{A},gx_{1}x_{2}) &=&0\left( \ref{X1,gx1x2,X1F11,gx1}\right) \text{ }
\\
&&\text{this is }\left( \ref{G,gx1x2, GF1,gx1}\right)
\end{eqnarray*}%
\begin{eqnarray*}
B\left( x_{1}\otimes 1_{H};1_{A},gx_{1}\right) +B(gx_{1}x_{2}\otimes
1_{H};1_{A},gx_{1}x_{2}) &=&0\left( \ref{X2,x1,X2F11,gx1x2}\right) \\
&&\text{this is }\left( \ref{G,gx1x2, GF1,gx1}\right)
\end{eqnarray*}%
\begin{eqnarray*}
&&B(x_{1}\otimes 1_{H};1_{A},gx_{1})+B(gx_{1}x_{2}\otimes
1_{H};1_{A},gx_{1}x_{2})=0\left( \ref{X2,x1,X2F41,gx1}\right) \\
&&\text{this is }\left( \ref{G,gx1x2, GF1,gx1}\right)
\end{eqnarray*}%
\begin{gather*}
B(x_{1}\otimes 1_{H};1_{A},gx_{1})+B(gx_{1}x_{2}\otimes
1_{H};1_{A},gx_{1}x_{2})=0\left( \ref{X2,gx1x2,X2F11,gx1}\right) \\
\text{this is }\left( \ref{G,gx1x2, GF1,gx1}\right)
\end{gather*}

\begin{equation*}
B(x_{2}\otimes 1_{H};1_{A},gx_{2})+B(gx_{1}x_{2}\otimes
1_{H};1_{A},gx_{1}x_{2})=0.\left( \ref{G,gx1x2, GF1,gx2}\right)
\end{equation*}%
\begin{eqnarray*}
B\left( x_{2}\otimes 1_{H};1_{A},gx_{2}\right) +B(gx_{1}x_{2}\otimes
1_{H};1_{A},gx_{1}x_{2}) &=&0\left( \ref{X1,x2,X1F11,gx1x2}\right) \\
&&\text{this is }\left( \ref{G,gx1x2, GF1,gx2}\right)
\end{eqnarray*}%
\begin{eqnarray*}
B(x_{2}\otimes 1_{H};1_{A},gx_{2})+B(gx_{1}x_{2}\otimes
1_{H};1_{A},gx_{1}x_{2}) &=&0.\left( \ref{X1,x2,X1F31,gx2}\right) \\
&&\text{ this is }\left( \ref{G,gx1x2, GF1,gx2}\right)
\end{eqnarray*}%
\begin{eqnarray*}
B(x_{2}\otimes 1_{H};1_{A},gx_{2})+B(gx_{1}x_{2}\otimes
1_{H};1_{A},gx_{1}x_{2}) &=&0\left( \ref{X1,gx1x2,X171,gx2}\right) \\
&&\text{this is }\left( \ref{G,gx1x2, GF1,gx2}\right)
\end{eqnarray*}%
\begin{eqnarray*}
B(x_{2}\otimes 1_{H};1_{A},gx_{2})+B(gx_{1}x_{2}\otimes
1_{H};1_{A},gx_{1}x_{2}) &=&0\left( \ref{X1,gx1x2,X1F31,gx1x2}\right) \\
&&\text{this is }\left( \ref{G,gx1x2, GF1,gx2}\right)
\end{eqnarray*}%
\begin{eqnarray*}
&&\left[
\begin{array}{c}
B(g\otimes 1_{H};1_{A},g)+B(x_{2}\otimes \ 1_{H};1_{A},gx_{2})+ \\
B(x_{1}\otimes 1_{H};1_{A},gx_{1})+B(gx_{1}x_{2}\otimes
1_{H};1_{A},gx_{1}x_{2})%
\end{array}%
\right] +\left( \ref{G,gx1x2, GF4,g}\right) \\
&&\text{ this is trivial in view of }\left( \ref{G,x2, GF1,g}\right) \text{
and }\left( \ref{G,gx1x2, GF1,gx1}\right)
\end{eqnarray*}%
\begin{eqnarray*}
&&\left[
\begin{array}{c}
B(g\otimes 1_{H};1_{A},g)+B(x_{2}\otimes \ 1_{H};1_{A},gx_{2}) \\
+B(x_{1}\otimes 1_{H};1_{A},gx_{1})+B(gx_{1}x_{2}\otimes
1_{H};1_{A},gx_{1}x_{2})%
\end{array}%
\right] +\left( \ref{G,gx1x2, GF3,g}\right) \\
&&\text{ this is }\left( \ref{G,gx1x2, GF4,g}\right) \text{ and therefore
trivial.}
\end{eqnarray*}%
\begin{eqnarray*}
&&\left[
\begin{array}{c}
B(g\otimes 1_{H};1_{A},g)+B(x_{2}\otimes \ 1_{H};1_{A},gx_{2}) \\
+B(x_{1}\otimes 1_{H};1_{A},gx_{1})+B(gx_{1}x_{2}\otimes
1_{H};1_{A},gx_{1}x_{2})%
\end{array}%
\right] \left( \ref{X1,x2,X1F51,g}\right) \\
&&\text{ this is }\left( \ref{G,gx1x2, GF4,g}\right) \text{ and therefore
trivial.}
\end{eqnarray*}%
\begin{gather*}
\left[
\begin{array}{c}
B(g\otimes 1_{H};1_{A},g)+B(x_{2}\otimes \ 1_{H};1_{A},gx_{2}) \\
+B(x_{1}\otimes 1_{H};1_{A},gx_{1})+B(gx_{1}x_{2}\otimes
1_{H};1_{A},gx_{1}x_{2})%
\end{array}%
\right] =0\left( \ref{X1,gx1x2,X1F31,g}\right) \\
\text{this is }\left( \ref{G,gx1x2, GF4,g}\right) \text{ and therefore
trivial.}
\end{gather*}%
\begin{gather*}
B(x_{1}\otimes 1_{H};1_{A},gx_{1})+B(gx_{1}x_{2}\otimes
1_{H};1_{A},gx_{1}x_{2})\left( \ref{X1,x2,X1F41,gx1}\right) \\
+B(g\otimes 1_{H};1_{A},g)+B(x_{2}\otimes \ 1_{H};1_{A},gx_{2})=0. \\
\text{this is }\left( \ref{G,gx1x2, GF4,g}\right) \text{ and therefore
trivial.}
\end{gather*}%
\begin{gather*}
\left[
\begin{array}{c}
B(g\otimes 1_{H};1_{A},g)+B(x_{2}\otimes \ 1_{H};1_{A},gx_{2}) \\
+B(x_{1}\otimes 1_{H};1_{A},gx_{1})+B(gx_{1}x_{2}\otimes
1_{H};1_{A},gx_{1}x_{2})%
\end{array}%
\right] =0\left( \ref{X1,gx1x2,X1F41,g}\right) \\
\text{this is }\left( \ref{G,gx1x2, GF4,g}\right) \text{ and therefore
trivial.}
\end{gather*}%
\begin{gather*}
\left[
\begin{array}{c}
B(g\otimes 1_{H};1_{A},g)+B(x_{2}\otimes \ 1_{H};1_{A},gx_{2}) \\
+B(x_{1}\otimes 1_{H};1_{A},gx_{1})+B(gx_{1}x_{2}\otimes
1_{H};1_{A},gx_{1}x_{2})%
\end{array}%
\right] =0\left( \ref{X1,gx1x2,X1F51,gx1}\right) \\
\text{this is }\left( \ref{G,gx1x2, GF4,g}\right) \text{ and therefore
trivial.}
\end{gather*}%
\begin{gather*}
B(g\otimes 1_{H};1_{A},g)+B(x_{1}\otimes 1_{H};1_{A},gx_{1})\left( \ref%
{X2,x1,X2F31,gx2}\right) \\
+B(x_{2}\otimes 1_{H};1_{A},gx_{2})+B(gx_{1}x_{2}\otimes
1_{H};1_{A},gx_{1}x_{2})=0 \\
\text{this is }\left( \ref{G,gx1x2, GF4,g}\right) \text{ and therefore
trivial.}
\end{gather*}%
\begin{gather*}
\left[
\begin{array}{c}
B(g\otimes 1_{H};1_{A},g)+B(x_{2}\otimes \ 1_{H};1_{A},gx_{2}) \\
+B(x_{1}\otimes 1_{H};1_{A},gx_{1})+B(gx_{1}x_{2}\otimes
1_{H};1_{A},gx_{1}x_{2})%
\end{array}%
\right] =0.\left( \ref{X2,gx1x2,X2F31,g}\right) \\
\text{this is }\left( \ref{G,gx1x2, GF4,g}\right) \text{ and therefore
trivial.}
\end{gather*}%
\begin{eqnarray*}
&&\left[
\begin{array}{c}
B(g\otimes 1_{H};1_{A},g)+B(x_{2}\otimes \ 1_{H};1_{A},gx_{2}) \\
+B(x_{1}\otimes 1_{H};1_{A},gx_{1})+B(gx_{1}x_{2}\otimes
1_{H};1_{A},gx_{1}x_{2})%
\end{array}%
\right] \left( \ref{X2,gx1x2,X2F41,g}\right) \\
&&\text{this is }\left( \ref{G,gx1x2, GF4,g}\right) \text{ and therefore
trivial.}
\end{eqnarray*}%
\begin{eqnarray*}
&&\lambda B(g\otimes 1_{H};X_{2},gx_{1})+2B(x_{1}\otimes
1_{H};1_{A},gx_{1})+\left( \ref{X1,g, X1F11,gx1}\right) \\
&&+2B\left( g\otimes 1_{H};1_{A},g\right) +=0\text{ this is trivial by }%
\left( \ref{got1,X2,gx1}\right) \text{and}\left( \ref{G,x1, GF1,g}\right)
\end{eqnarray*}%
\begin{equation*}
\begin{array}{c}
2\alpha B(x_{1}x_{2}\otimes 1_{H};G,x_{1})+ \\
+\gamma _{1}B(x_{1}x_{2}\otimes 1_{H};X_{1},x_{1})+\gamma
_{2}B(x_{1}x_{2}\otimes 1_{H};X_{2},x_{1})%
\end{array}%
=0.\left( \ref{G,x1x2, GF1,x1}\right)
\end{equation*}%
\begin{equation*}
\begin{array}{c}
2\alpha B(x_{1}x_{2}\otimes 1_{H};G,x_{2})+ \\
+\gamma _{1}B(x_{1}x_{2}\otimes 1_{H};X_{1},x_{2})+\gamma
_{2}B(x_{1}x_{2}\otimes 1_{H};X_{2},x_{2})%
\end{array}%
=0.\left( \ref{G,x1x2, GF1,x2}\right)
\end{equation*}

\begin{gather*}
1-B(x_{1}x_{2}\otimes 1_{H};1_{A},x_{1}x_{2})+B(x_{1}x_{2}\otimes
1_{H};x_{1},x_{2})-B(x_{1}x_{2}\otimes 1_{H};x_{2},x_{1})=0 \\
\left( \ref{G,x1x2, GF4,1H}\right)
\end{gather*}%
\begin{eqnarray*}
\left[
\begin{array}{c}
+1-B(x_{1}x_{2}\otimes 1_{H};1_{A},x_{1}x_{2})-B(x_{1}x_{2}\otimes
1_{H};X_{2},x_{1}) \\
+B(x_{1}x_{2}\otimes 1_{H};X_{1},x_{2})%
\end{array}%
\right] &=&0\left( \ref{X1,x1x2,X1F31,1H}\right) \\
&&\text{ this is }\left( \ref{G,x1x2, GF4,1H}\right)
\end{eqnarray*}%
\begin{eqnarray*}
\left[
\begin{array}{c}
+1-B(x_{1}x_{2}\otimes 1_{H};1_{A},x_{1}x_{2})-B(x_{1}x_{2}\otimes
1_{H};X_{2},x_{1}) \\
+B(x_{1}x_{2}\otimes 1_{H};X_{1},x_{2})%
\end{array}%
\right] &=&0.\left( \ref{X1,x1x2,X1F41,1H}\right) \\
&&\text{ this is }\left( \ref{G,x1x2, GF4,1H}\right)
\end{eqnarray*}%
\begin{eqnarray*}
\left[
\begin{array}{c}
+1-B(x_{1}x_{2}\otimes 1_{H};1_{A},x_{1}x_{2}) \\
-B(x_{1}x_{2}\otimes 1_{H};X_{2},x_{1})+B(x_{1}x_{2}\otimes
1_{H};X_{1},x_{2})%
\end{array}%
\right] &=&0.\left( \ref{X2,x1x2,X2F31,1H}\right) \\
&&\text{ this is }\left( \ref{X1,x1x2,X1F31,1H}\right)
\end{eqnarray*}%
\begin{eqnarray*}
\left[
\begin{array}{c}
+1-B(x_{1}x_{2}\otimes 1_{H};1_{A},x_{1}x_{2}) \\
-B(x_{1}x_{2}\otimes 1_{H};X_{2},x_{1})+B(x_{1}x_{2}\otimes
1_{H};X_{1},x_{2})%
\end{array}%
\right] &=&0.\left( \ref{X2,x1x2,X2F41,1H}\right) \\
&&\text{ this is }\left( \ref{G,x1x2, GF4,1H}\right)
\end{eqnarray*}

\begin{equation*}
\begin{array}{c}
2\alpha B(x_{1}x_{2}\otimes 1_{H};G,x_{2})+ \\
+\gamma _{1}\left[ -1+B(x_{1}x_{2}\otimes
1_{H};1_{A},x_{1}x_{2})+B(x_{1}x_{2}\otimes 1_{H};X_{2},x_{1})\right]
+\gamma _{2}B(x_{1}x_{2}\otimes 1_{H};X_{2},x_{2})%
\end{array}%
=0.\left( \ref{G,gx1, GF1,1H}\right)
\end{equation*}%
\begin{equation*}
\begin{array}{c}
2\alpha \left[ -B(x_{1}x_{2}\otimes 1_{H};G,x_{1})\right] + \\
-\gamma _{1}B(x_{1}x_{2}\otimes 1_{H};X_{1},x_{1})+\gamma _{2}\left[
-1+B(x_{1}x_{2}\otimes 1_{H};1_{A},x_{1}x_{2})-B(x_{1}x_{2}\otimes
1_{H};X_{1},x_{2})\right]%
\end{array}%
=0\left( \ref{G,gx2, GF1,1H}\right)
\end{equation*}%
\begin{eqnarray*}
&&2\beta _{1}B(x_{1}x_{2}\otimes 1_{H};X_{1},x_{1})+\lambda
B(x_{1}x_{2}\otimes 1_{H};X_{2},x_{1}) \\
+\gamma _{1}B(x_{1}x_{2}\otimes 1_{H};G,x_{1}) &=&0\left( \ref%
{X1,x1x2,X1F11,x1}\right)
\end{eqnarray*}%
\begin{eqnarray*}
&&2\beta _{1}B(x_{1}x_{2}\otimes 1_{H};X_{1},x_{2})+\lambda
B(x_{1}x_{2}\otimes 1_{H};X_{2},x_{2}) \\
+\gamma _{1}B(x_{1}x_{2}\otimes 1_{H};G,x_{2}) &=&0.\left( \ref%
{X1,x1x2,X1F11,x2}\right)
\end{eqnarray*}%
\begin{gather*}
2\beta _{2}B(x_{1}x_{2}\otimes 1_{H};X_{2},x_{1})+\gamma
_{2}B(x_{1}x_{2}\otimes 1_{H};G,x_{1})\left( \ref{X2,x1x2,X2F11,x1}\right) \\
+\lambda B(x_{1}x_{2}\otimes 1_{H};X_{1},x_{1})=0.
\end{gather*}

\begin{gather*}
2\beta _{2}B(x_{1}x_{2}\otimes 1_{H};X_{2},x_{2})+\gamma
_{2}B(x_{1}x_{2}\otimes 1_{H};G,x_{2})\left( \ref{X2,x1x2,X2F11,x2}\right) \\
+\lambda B(x_{1}x_{2}\otimes 1_{H};X_{1},x_{2})=0.
\end{gather*}

\begin{gather*}
-2\beta _{1}B(x_{1}x_{2}\otimes 1_{H};X_{1},x_{1})+\lambda \left[
-1+B(x_{1}x_{2}\otimes 1_{H};1_{A},x_{1}x_{2})-B(x_{1}x_{2}\otimes
1_{H};X_{1},x_{2})\right] \left( \ref{X1,gx2,X1F11,1H}\right) \\
+\gamma _{1}\left[ -B(x_{1}x_{2}\otimes 1_{H};G,x_{1})\right] =0.
\end{gather*}%
\begin{gather*}
2\beta _{1}\left[ -1+B(x_{1}x_{2}\otimes
1_{H};1_{A},x_{1}x_{2})+B(x_{1}x_{2}\otimes 1_{H};X_{2},x_{1})\right] \left( %
\ref{X1,gx1,X1F11,1H}\right) \\
+\lambda B(x_{1}x_{2}\otimes 1_{H};X_{2},x_{2})+\gamma
_{1}B(x_{1}x_{2}\otimes 1_{H};G,x_{2})=0
\end{gather*}%
\begin{gather*}
2\beta _{2}\left[ -1+B(x_{1}x_{2}\otimes
1_{H};1_{A},x_{1}x_{2})-B(x_{1}x_{2}\otimes 1_{H};X_{1},x_{2})\right] \left( %
\ref{X2,gx2,X2F11,1H}\right) \\
+\gamma _{2}\left[ -B(x_{1}x_{2}\otimes 1_{H};G,x_{1})\right] -\lambda
B(x_{1}x_{2}\otimes 1_{H};X_{1},x_{1})=0.
\end{gather*}

\begin{gather*}
2\beta _{2}B(x_{1}x_{2}\otimes 1_{H};X_{2},x_{2})+\gamma
_{2}B(x_{1}x_{2}\otimes 1_{H};G,x_{2})\left( \ref{X2,gx1,X2F11,1H}\right) \\
+\lambda \left[ -1+B(x_{1}x_{2}\otimes
1_{H};1_{A},x_{1}x_{2})+B(x_{1}x_{2}\otimes 1_{H};X_{2},x_{1})\right] =0.
\end{gather*}%
By $\left( \ref{G,gx1x2, GF1,gx1}\right) $ and $\left( \ref{x1ot1,1,gx2}%
\right) $
\begin{gather*}
2\beta _{2}B(x_{1}\otimes 1_{H};1_{A},gx_{2})+\left( \ref{X2,gx1x2,X2F11,gx2}%
\right) \\
+\lambda \left[ -B(x_{2}\otimes 1_{H};1_{A},gx_{2})-B(gx_{1}x_{2}\otimes
1_{H};1_{A},gx_{1}x_{2})\right] =0.
\end{gather*}%
becomes trivial.

By $\left( \ref{G,x1, GF1,g}\right) $ and $\left( \ref{got1,X2,gx1}\right) $
the equality
\begin{eqnarray*}
&&\lambda B(g\otimes 1_{H};X_{2},gx_{1})+2B(x_{1}\otimes
1_{H};1_{A},gx_{1})+\left( \ref{X1,g, X1F11,gx1}\right) \\
&&+2B\left( g\otimes 1_{H};1_{A},g\right) =0
\end{eqnarray*}%
becomes trivial.%
\begin{equation*}
0=0\text{ }\left( \ref{G,g,GF1,1H}\right)
\end{equation*}

\begin{equation*}
2\alpha 0+\gamma _{2}0=0\left( \ref{G,g, GF1,gx1}\right)
\end{equation*}

\begin{equation*}
2\alpha 0-\gamma _{1}0=0.\left( \ref{G,g, GF1,gx2}\right)
\end{equation*}

\begin{equation*}
\gamma _{1}0+\gamma _{2}0=0.\left( \ref{G,g, GF2,g}\right)
\end{equation*}

\begin{equation*}
-\gamma _{1}0+2\alpha 0=0.\left( \ref{G,x1, GF1,gx1x2}\right)
\end{equation*}%
\begin{equation*}
0=0.\left( \ref{G,x1, GF2,1H}\right)
\end{equation*}%
\begin{equation*}
0=0\text{ }\left( \ref{G,x1, GF2,x1x2}\right)
\end{equation*}%
\begin{equation*}
\gamma _{1}0+\gamma _{2}0=0.\left( \ref{G,x1, GF2,gx1}\right)
\end{equation*}

\begin{equation*}
\gamma _{1}[0-0]=0.\left( \ref{G,x1,GF2,gx2}\right)
\end{equation*}

\begin{equation*}
0=0\left( \ref{G,x1,GF3,1H}\right)
\end{equation*}%
\begin{equation*}
-\gamma _{2}0-2\alpha 0=0.\left( \ref{G,x1,GF3,gx1}\right)
\end{equation*}

\begin{equation*}
\alpha \left[ -0+0\right] =0.\left( \ref{G,x1, GF3,x2}\right)
\end{equation*}%
\begin{equation*}
0=0.\left( \ref{G,x1, GF4,1H}\right)
\end{equation*}

\begin{equation*}
-\gamma _{1}0+2\alpha 0=0.\left( \ref{G,x1, GF4,gx1}\right)
\end{equation*}

\begin{equation*}
\alpha \lbrack 0+0]=0.\left( \ref{G,x1, GF5,g}\right)
\end{equation*}%
\begin{equation*}
0=0.\left( \ref{G,x1,GF6,x1}\right)
\end{equation*}%
\begin{equation*}
0=0.\left( \ref{G,x1,GF6,x2}\right)
\end{equation*}%
\begin{equation*}
0=0\left( \ref{G,x2, GF1,x1}\right)
\end{equation*}

\begin{equation*}
2\alpha 0-\gamma _{2}0=0.\left( \ref{G,x2, GF1,gx1x2}\right)
\end{equation*}%
\begin{equation*}
0=0.\left( \ref{G,x2, GF2,1H}\right)
\end{equation*}%
\begin{equation*}
0=0\left( \ref{G,x2, GF2,x1x2}\right)
\end{equation*}%
\begin{equation*}
\gamma _{2}\left[ 0+0\right] =0.\left( \ref{G,x2, GF2,gx1}\right)
\end{equation*}

\begin{equation*}
-\gamma _{1}0+\gamma _{2}0=0.\left( \ref{G,x2, GF2,gx2}\right)
\end{equation*}%
\begin{equation*}
0=0\left( \ref{G,x2, FG3,1H}\right)
\end{equation*}

\begin{equation*}
-\gamma _{2}0-2\alpha 0=0.\left( \ref{G,x2, GF3,gx2}\right)
\end{equation*}%
\begin{equation*}
0=0.\left( \ref{G,x2, GF4,1H}\right)
\end{equation*}%
\begin{equation*}
\alpha \left[ 0+0\right] =0\left( \ref{G,x2, GF4,gx1}\right)
\end{equation*}

\begin{equation*}
\gamma _{1}0-2\alpha 0=0.\left( \ref{G,x2, GF4,gx2}\right)
\end{equation*}%
\begin{equation*}
\alpha \left[ 0+0\right] =0\left( \ref{G,x2, GF5,g}\right)
\end{equation*}

\begin{equation*}
\gamma _{2}[0+0=0.=0.\left( \ref{G,x2, GF6,g}\right)
\end{equation*}%
\begin{equation*}
0=0.\left( \ref{G,x2; GF6,x2}\right)
\end{equation*}

\begin{equation*}
\gamma _{1}[0+0]=0.\left( \ref{G,x2; GF7,g}\right)
\end{equation*}%
\begin{equation*}
0=0.\left( \ref{G,x2; GF7,x1}\right)
\end{equation*}%
\begin{equation*}
\gamma _{1}0+\gamma _{2}0=0.\left( \ref{G,x1x2, GF1,g}\right)
\end{equation*}

\begin{equation*}
0=0\left( \ref{G,x1x2, GF2,1H}\right)
\end{equation*}%
\begin{equation*}
\begin{array}{c}
20 \\
+\gamma _{2}0%
\end{array}%
=0.\left( \ref{G,x1x2, GF2,gx1}\right)
\end{equation*}%
\begin{equation*}
\begin{array}{c}
20 \\
+\gamma _{1}0%
\end{array}%
=0\left( \ref{G,x1x2, GF2,gx2}\right)
\end{equation*}

\begin{equation*}
0=0\text{ }\left( \ref{G,x1x2, GF3,gx1}\right)
\end{equation*}

\begin{equation*}
-20+\gamma _{2}[-0-0+0]=0\left( \ref{G,x1x2, GF6,g}\right)
\end{equation*}%
\begin{gather*}
-20+ \\
-\gamma _{1}[-0-0+0]=0\left( \ref{G,x1x2, GF7,g}\right)
\end{gather*}

\begin{equation*}
\begin{array}{c}
2\left[ 0+0\right] \\
\gamma _{1}\left[ -0+0\right]%
\end{array}%
\left( \ref{G,gx1, GF2,g}\right)
\end{equation*}%
\begin{equation*}
\begin{array}{c}
2\left[ -0-0\right] \\
+\gamma _{2}\left[ -0-0\right]%
\end{array}%
=0.\left( \ref{G,gx2, GF2,g}\right)
\end{equation*}

\begin{equation*}
0=0\left( \ref{G,gx1x2, GF1,1H}\right)
\end{equation*}%
\begin{equation*}
0=0.\left( \ref{G,gx1x2, GF1,x1x2}\right)
\end{equation*}%
\begin{gather*}
\gamma _{1}\left[ 0-0\right] \left( \ref{G,gx1x2, GF2,g}\right) \\
+\gamma _{2}\left[ -0-0\right] =0.
\end{gather*}%
\begin{equation*}
0=0.\left( \ref{G,gx1x2, GF2,x1}\right)
\end{equation*}%
\begin{equation*}
0=0.\left( \ref{G,gx1x2, GF2,x2}\right)
\end{equation*}%
\begin{equation*}
\gamma _{1}0-\gamma _{2}0=0.\left( \ref{G,gx1x2, GF2,gx1x2}\right)
\end{equation*}

\begin{equation*}
0=0.\left( \ref{G,gx1x2, GF3,x1}\right)
\end{equation*}

\begin{equation*}
0=0\left( \ref{G,gx1x2, GF3,x2}\right)
\end{equation*}

\begin{equation*}
\gamma _{2}0-2\alpha 0=0.\left( \ref{G,gx1x2, GF3,gx1x2}\right)
\end{equation*}

\begin{equation*}
0=0\left( \ref{G,gx1x2, GF4,x1}\right)
\end{equation*}%
\begin{equation*}
0=0.\left( \ref{G,gx1x2, GF4,x2}\right)
\end{equation*}

\begin{equation*}
-\gamma _{1}0+2\alpha 0=0.\left( \ref{G,gx1x2, GF4,gx1x2}\right)
\end{equation*}

\begin{equation*}
\alpha \left[ 0-0\right] =0.\left( \ref{G,gx1x2, GF5,gx2}\right)
\end{equation*}

\begin{equation*}
0=0\left( \ref{G,gx1x2, GF6,1H}\right)
\end{equation*}%
\begin{equation*}
\gamma _{2}\left[ 0-0\right] =0.\left( \ref{G,gx1x2, GF6,gx1}\right)
\end{equation*}%
\begin{equation*}
\gamma _{2}\left[ 0-0\right] =0.\left( \ref{G,gx1x2, GF6,gx2}\right)
\end{equation*}%
\begin{equation*}
0=0\left( \ref{G,gx1x2, GF7,1H}\right)
\end{equation*}%
\begin{equation*}
0=0.\left( \ref{G,gx1x2, GF7,x1x2}\right)
\end{equation*}%
\begin{equation*}
\gamma _{1}\left[ 0+0\right] =0.\left( \ref{G,gx1x2, GF7,gx1}\right)
\end{equation*}%
\begin{equation*}
\gamma _{1}\left[ 0-0\right] =0.\left( \ref{G,gx1x2, GF7,gx2}\right)
\end{equation*}

\begin{equation*}
0=0.\left( \ref{X1,g, X1F11,1H}\right)
\end{equation*}

\begin{equation*}
\gamma _{1}0+20-2\beta _{1}0=0.\left( \ref{X1,g, X1F11,gx2}\right)
\end{equation*}

\begin{equation*}
2\beta _{1}0+\lambda 0+20=0.\left( \ref{X1,g,X1F21,g}\right)
\end{equation*}%
\begin{equation*}
0=0.\left( \ref{X1,g,X1F21,x1}\right)
\end{equation*}

\begin{equation*}
0-0=0\left( \ref{X1,g,X1F21,gx1x2}\right)
\end{equation*}

\begin{equation*}
-2\beta _{1}0+20+\gamma _{1}0=0.\left( \ref{X1,g, X1F41,g}\right)
\end{equation*}%
\begin{equation*}
0=0.\left( \ref{X1,g, X1F41,x1}\right)
\end{equation*}

\begin{equation*}
0=0\left( \ref{X1,g,X1F51,1H}\right)
\end{equation*}%
\begin{equation*}
0-0=0.\left( \ref{X1,x2,X1F71,gx1x2}\right)
\end{equation*}%
\begin{equation*}
0=0.\left( \ref{X1,x1,X1F11,x1}\right)
\end{equation*}

\begin{equation*}
0=0\left( \ref{X1,x1,X1F21,1}\right)
\end{equation*}%
\begin{equation*}
\lambda 0=-20.\left( \ref{X1,x1,X1F21,gx1}\right)
\end{equation*}%
\begin{equation*}
0=0.\left( \ref{X1,x1,X1F31,1H}\right)
\end{equation*}%
\begin{equation*}
0=0\left( \ref{X1,x1,X1F31,x1x2}\right)
\end{equation*}

\begin{equation*}
\gamma _{1}\left[ -0+0\right] =0\left( \ref{X1,x1,X1F31,gx2}\right)
\end{equation*}

\begin{equation*}
2\beta _{1}0-\gamma _{1}0=0.\left( \ref{X1,x1,X1F41,gx1}\right)
\end{equation*}

\begin{equation*}
0=0.\left( \ref{X1,x2,X1F11,x1}\right)
\end{equation*}%
\begin{equation*}
0=0.\left( \ref{X1,x2,X1F11,x2}\right)
\end{equation*}

\begin{equation*}
0=0.\left( \ref{X1,x2,X1F21,1H}\right)
\end{equation*}%
\begin{equation*}
0=0.\left( \ref{X1,x2,X1F21,x1x2}\right)
\end{equation*}%
\begin{equation*}
\lambda \left[ 0+0\right] =20-20.\left( \ref{X1,x2,X1F21,gx1}\right)
\end{equation*}

\begin{equation*}
-2\beta _{1}0+\lambda 0-20=0.\left( \ref{X1,x2,X1F21,gx2}\right)
\end{equation*}

\begin{equation*}
0=0.\left( \ref{X1,x2,X1F41,1H}\right)
\end{equation*}%
\begin{equation*}
0=0.\left( \ref{X1,x2,X1F41,x1x2}\right)
\end{equation*}

\begin{equation*}
2\beta _{1}0=\gamma _{1}0+20.\left( \ref{X1,x2,X1F41,gx2}\right)
\end{equation*}

\begin{equation*}
0=0.\left( \ref{X1,x2,X1F51,x2}\right)
\end{equation*}%
\begin{equation*}
\lambda \left[ 0+0\right] =-2\left[ 0-0\right] \left( \ref{X1,x2,X1F61,g}%
\right)
\end{equation*}%
\begin{equation*}
\beta _{1}\left[ 0+0\right] =0\left( \ref{X1,x2,X1F71,g}\right)
\end{equation*}%
\begin{equation*}
0=0\left( \ref{X1,x2,X1F71,x1}\right)
\end{equation*}%
\begin{equation*}
\lambda 0-\gamma _{1}0=0\left( \ref{X1,x1x2,X1F11,g}\right)
\end{equation*}%
\begin{equation*}
\gamma _{1}0-20=0.\left( \ref{X1,x1x2,X1F11,gx1x2}\right)
\end{equation*}%
\begin{equation*}
0=0.\left( \ref{X1,x1x2,X1F21,1H}\right)
\end{equation*}%
\begin{equation*}
\lambda 0+20=0.\left( \ref{X1,x1x2,X1F21,gx1}\right)
\end{equation*}

\begin{equation*}
20+20=0\left( \ref{X1,x1x2,X1F31,gx1}\right)
\end{equation*}%
\begin{equation*}
20+\gamma _{1}0=0.\left( \ref{X1,x1x2,X1F31,gx2}\right)
\end{equation*}

\begin{equation*}
\gamma _{1}0+20=0.\left( \ref{X1,x1x2,X1F41,gx1}\right)
\end{equation*}%
\begin{gather*}
20+\left( \ref{X1,x1x2,X1F51,g}\right) \\
+\gamma _{1}\left[ -0+0-0\right] =0
\end{gather*}%
\begin{equation*}
20+\lambda \left[
\begin{array}{c}
-0 \\
+0-0%
\end{array}%
\right] =0\left( \ref{X1,x1x2,X1F61,g}\right)
\end{equation*}%
\begin{equation*}
0+0=0\left( \ref{X1,x1x2,X1F61,gx1x2}\right)
\end{equation*}%
\begin{equation*}
+0+0=0.\left( \ref{X1,x1x2,X1F81,gx1}\right)
\end{equation*}

\begin{gather*}
-\gamma _{1}\left[ 0-0\right] \left( \ref{X1,gx1,X1F11,gx1}\right) \\
+2\left[ 0+0\right] =0.
\end{gather*}%
\begin{gather*}
2\left[ 0+0\right] +\left( \ref{X1,gx1,X1F31,g}\right) \\
+ \\
+\gamma _{1}\left[ -0+0\right] =0.
\end{gather*}%
\begin{gather*}
\left( \ref{X1,gx2,X1F11,gx1}\right) \\
+20+20=0.
\end{gather*}%
\begin{equation*}
+\gamma _{1}\left[ 0+0\right] =0\left( \ref{X1,gx2,X1F11,gx2}\right)
\end{equation*}%
\begin{gather*}
2\left[ 0+0\right] \left( \ref{X1,gx2,X131,g}\right) \\
=0
\end{gather*}%
\begin{equation*}
\gamma _{1}\left[ -0+0\right] =0\left( \ref{X1,gx2,X141,g}\right)
\end{equation*}%
\begin{equation*}
0=0\left( \ref{X1,gx1x2,X1F11,1H}\right)
\end{equation*}%
\begin{equation*}
0=0\left( \ref{X1,gx1x2,X1F11,x1x2}\right)
\end{equation*}%
\begin{gather*}
2\beta _{1}\left[ 0-0\right] \left( \ref{X1,gx1x2,XF21,g}\right) \\
+\lambda \left[ 0-0\right] =0
\end{gather*}

\begin{equation*}
2\beta _{1}0-\lambda 0=0\left( \ref{X1,gx1x2,X1F21,gx1x2}\right)
\end{equation*}

\begin{equation*}
0=0\left( \ref{X1,gx1x2,X1F31,x1}\right)
\end{equation*}

\begin{equation*}
0=0\left( \ref{X1,gx1x2,X1F31,x2}\right)
\end{equation*}

\begin{eqnarray*}
&&2\beta _{1}0\left( \ref{X1,gx1x2,X1F41,gx1x2}\right) \\
&&-\gamma _{1}0-20=0
\end{eqnarray*}%
\begin{equation*}
0=0\left( \ref{X1,gx1x2,X1F51,1H}\right)
\end{equation*}%
\begin{equation*}
0=0\left( \ref{X1,gx1x2,X1F61,1H}\right)
\end{equation*}%
\begin{equation*}
0=0\left( \ref{X1,gx1x2,X1F61,x1x2}\right)
\end{equation*}

\begin{equation*}
\lambda \left[ 0+0\right] +\left[ 20-20\right] =0\left( \ref%
{X1,gx1x2,X1F61,gx1}\right)
\end{equation*}%
\begin{gather*}
\beta _{1}\left[ 0+0\right] \left( \ref{X1,gx1x2,X1F71,gx1}\right) \\
+0+0=0.
\end{gather*}

\begin{equation*}
0=0\left( \ref{X2,g,X2F11,1H}\right)
\end{equation*}%
\begin{gather*}
2\beta _{2}0+\gamma _{2}0\left( \ref{X2,g,X2F11,gx1}\right) \\
+20=0
\end{gather*}

\begin{equation*}
2\beta _{2}0+\lambda 0+20=0\left( \ref{X2,g,X2F21,g}\right)
\end{equation*}%
\begin{equation*}
0=0\left( \ref{X2,g,X2F21,x1}\right)
\end{equation*}%
\begin{equation*}
0=0\ref{X2,g,X2F21,x2}
\end{equation*}%
\begin{equation*}
2\beta _{2}0+\gamma _{2}0+20=0.\left( \ref{X2,g,X2F31,g}\right)
\end{equation*}%
\begin{equation*}
0=0\left( \ref{X2,g,X2F31,x2}\right)
\end{equation*}%
\begin{equation*}
0+0=0\left( \ref{X2,g,X2F41,x1}\right)
\end{equation*}%
\begin{equation*}
0+0=0\left( \ref{X2,g,X2F61,gx2}\right)
\end{equation*}%
\begin{equation*}
0=0\left( \ref{X2,g,X2F71,1H}\right)
\end{equation*}%
\begin{equation*}
0+0=0\left( \ref{X2,g,X2F71,gx1}\right)
\end{equation*}

\begin{equation*}
0=0\left( \ref{X2,x1,X2F11,x1}\right)
\end{equation*}%
\begin{equation*}
0=0.\left( \ref{X2,x1,X2F11,x2}\right)
\end{equation*}

\begin{equation*}
0=0\left( \ref{X2,x1,X2F21,1H}\right)
\end{equation*}

\begin{equation*}
2\beta _{2}0+\lambda 0+20=0.\left( \ref{X2,x1,X2F21,gx1}\right)
\end{equation*}

\begin{eqnarray*}
&&+\lambda \left[ 0-0\right] \\
&&+2B0+20=0.
\end{eqnarray*}%
\begin{eqnarray*}
&&\lambda \left[ 0-0\right] \left( \ref{X2,x1,X2F21,gx2}\right) \\
&&+2-0+20=0.
\end{eqnarray*}

\begin{equation*}
0=0.\left( \ref{X2,x1,X2F31,1H}\right)
\end{equation*}

\begin{equation*}
0=0.\left( \ref{X2,x1,X2F41,1H}\right)
\end{equation*}

\begin{eqnarray*}
&&\beta _{2}\left[ -0+0\right] \left( \ref{X2,x1,X2F61,g}\right) \\
&&-0+0=0
\end{eqnarray*}%
\begin{equation*}
0=0\left( \ref{X2,x1,X2F61,x2}\right)
\end{equation*}

\begin{equation*}
0+0=0\left( \ref{X2,x1,X2F71,g}\right)
\end{equation*}%
\begin{equation*}
0=0\left( \ref{X2,x1,X2F81,1H}\right)
\end{equation*}%
\begin{equation*}
0=0.\left( \ref{X2,x2,X2F11,x2}\right)
\end{equation*}

\begin{eqnarray*}
+ &&2\beta _{2}0-\gamma _{2}0\left( \ref{X2,x2,X2F11,gx1x2}\right) \\
&&+20=0
\end{eqnarray*}

\begin{equation*}
\beta _{2}\left[ 0+0\right] =0\left( \ref{X2,x2,X2F21,gx1}\right)
\end{equation*}%
\begin{eqnarray*}
&&2\beta _{2}0-\lambda 0\left( \ref{X2,x2,X2F21,gx2}\right) \\
&&+20=0
\end{eqnarray*}

\begin{equation*}
0=0.\left( \ref{X2,x2,X2F41,1H}\right)
\end{equation*}

\begin{equation*}
\lambda \left[ 0+0\right] =0\left( \ref{X2,x2,X2F71,g}\right)
\end{equation*}%
\begin{equation*}
\gamma _{2}0+\lambda 0=0.\left( \ref{X2,x1x2,X2F11,g}\right)
\end{equation*}

\begin{eqnarray*}
&&\gamma _{2}0+\left( \ref{X2,x1x2,X2F11,gx1x2}\right) \\
&&+20=0.
\end{eqnarray*}

\begin{equation*}
0=0.\left( \ref{X2,x1x2,X2F21,1H}\right)
\end{equation*}

\begin{equation*}
20-\lambda 0=0.\left( \ref{X2,x1x2,X2F21,gx2}\right)
\end{equation*}%
\begin{equation*}
\gamma _{2}0+20=0.\left( \ref{X2,x1x2,X2F31,gx2}\right)
\end{equation*}

\begin{gather*}
20+\left( \ref{X2,x1x2,X2F41,gx2}\right) \\
+20=0.
\end{gather*}

\begin{gather*}
20-\left( \ref{X2,x1x2,X2F51,g}\right) \\
\gamma _{2}\left[
\begin{array}{c}
-0+ \\
0-0%
\end{array}%
\right] =0.
\end{gather*}%
\begin{gather*}
20\left( \ref{X2,x1x2,X2F71,g}\right) \\
+\lambda \left[
\begin{array}{c}
-0+ \\
0-0%
\end{array}%
\right] =0.
\end{gather*}%
\begin{equation*}
0-0=0.\left( \ref{X2,x1x2,X2F71,gx1x2}\right)
\end{equation*}

\begin{equation*}
0+0=0.\left( \ref{X2,gx1,X2F11,gx2}\right)
\end{equation*}%
\begin{equation*}
\lambda \left[ 0-0\right] =0.\left( \ref{X2,gx1,X2F21,g}\right)
\end{equation*}

\begin{eqnarray*}
&&\gamma _{2}\left[ 0+0\right] \left( \ref{X2,gx2,X2F11,gx2}\right) \\
&&+ \\
&&+\left[ +20+20\right] =0.
\end{eqnarray*}%
\begin{eqnarray*}
&&2\left[ -0-0\right] \left( \ref{X2,gx2,X2F41,g}\right) \\
&&+\gamma _{2}\left[ -0-0\right] \\
&=&0.
\end{eqnarray*}%
\begin{equation*}
0=0\left( \ref{X2,gx1x2,X2F11,1H}\right)
\end{equation*}%
\begin{equation*}
0=0\left( \ref{X2,gx1x2,X2F11,x1x2}\right)
\end{equation*}

\begin{gather*}
2\beta _{2}\left[ -0-0\right] \left( \ref{X2,gx1x2,X2F21,g}\right) \\
+\lambda \left[ 0-0\right] =0.
\end{gather*}%
\begin{gather*}
-2\beta _{2}0\left( \ref{X2,gx1x2,X2F21,gx1x2}\right) \\
+\lambda 0-20=0.
\end{gather*}

\begin{equation*}
0=0.\left( \ref{X2,gx1x2,X2F41,x1}\right)
\end{equation*}

\begin{equation*}
0=0.\left( \ref{X2,gx1x2,X2F51,1H}\right)
\end{equation*}

\begin{equation*}
0-0=0.\left( \ref{X2,gx1x2,X2F61,gx2}\right)
\end{equation*}%
\begin{equation*}
0=0.\left( \ref{X2,gx1x2,X2F71,1H}\right)
\end{equation*}

\subsection{CLEAN\ LIST\ OF WHAT\ WE\ GOT\ ABOVE}

\subsubsection{Binomial equalities without terms of $B\left(
x_{1}x_{2}\otimes 1_{H}\right) $ \label{BEW}}

\begin{equation*}
B(g\otimes 1_{H};1_{A},g)+B(x_{1}\otimes 1_{H};1_{A},gx_{1})=0.\left( \ref%
{G,x1, GF1,g}\right)
\end{equation*}%
\begin{equation*}
B(g\otimes 1_{H};1_{A},g)+B(x_{2}\otimes \ 1_{H};1_{A},gx_{2})=0.\left( \ref%
{G,x2, GF1,g}\right)
\end{equation*}%
\begin{equation*}
\beta _{1}B(g\otimes 1_{H};1_{A},gx_{1}x_{2})=B\left( x_{1}\otimes
1_{H};1_{A},gx_{2}\right) .\left( \ref{X1,x1,X1F11,gx1x2}\right)
\end{equation*}%
\begin{equation*}
B(x_{1}\otimes 1_{H};1_{A},gx_{1})+B(gx_{1}x_{2}\otimes
1_{H};1_{A},gx_{1}x_{2})=0.\left( \ref{G,gx1x2, GF1,gx1}\right)
\end{equation*}%
\begin{equation*}
B(x_{2}\otimes 1_{H};1_{A},gx_{2})+B(gx_{1}x_{2}\otimes
1_{H};1_{A},gx_{1}x_{2})=0.\left( \ref{G,gx1x2, GF1,gx2}\right)
\end{equation*}

\subsubsection{The final form of our seven elements}

We rewrite the equalities $\left( \ref{got1}\right) ,\left( \ref{x1}\right)
,\left( \ref{x2}\right) ,\left( \ref{x1x2}\right) ,\left( \ref{gx1}\right)
,\left( \ref{gx2}\right) $ and $\left( \ref{gx1x2}\right) $

taking in account all monomial equalities in $\left( \ref{LME}\right) $ and
all equalities in $\left( \ref{BEW}\right) .$

\begin{equation}
B\left( g\otimes 1_{H}\right) =B\left( g\otimes 1_{H};1_{A},g\right)
1_{A}\otimes g  \label{mbg}
\end{equation}%
\begin{equation}
B\left( x_{1}\otimes 1_{H}\right) =-B(g\otimes 1_{H};1_{A},g)1_{A}\otimes
gx_{1}  \label{mbx1}
\end{equation}%
\begin{equation}
B\left( x_{2}\otimes 1_{H}\right) =-B(g\otimes 1_{H};1_{A},g)1_{A}\otimes
gx_{2}  \label{mbx2}
\end{equation}%
\begin{equation}
B\left( gx_{1}x_{2}\otimes 1_{H}\right) =B(g\otimes
1_{H};1_{A},g)1_{A}\otimes gx_{1}x_{2}  \label{mbgx1x2}
\end{equation}%
\begin{eqnarray}
B(x_{1}x_{2}\otimes 1_{H}) &=&B(x_{1}x_{2}\otimes
1_{H};1_{A},1_{H})1_{A}\otimes 1_{H}+  \label{mx1x2} \\
&&+B(x_{1}x_{2}\otimes 1_{H};1_{A},x_{1}x_{2})1_{A}\otimes x_{1}x_{2}+
\notag \\
&&+B(x_{1}x_{2}\otimes 1_{H};G,x_{1})G\otimes x_{1}+  \notag \\
&&+B(x_{1}x_{2}\otimes 1_{H};G,x_{2})G\otimes x_{2}+  \notag \\
&&+B(x_{1}x_{2}\otimes 1_{H};X_{1},x_{1})X_{1}\otimes x_{1}+  \notag \\
&&+B(x_{1}x_{2}\otimes 1_{H};X_{1},x_{2})X_{1}\otimes x_{2}+  \notag \\
&&+B(x_{1}x_{2}\otimes 1_{H};X_{2},x_{1})X_{2}\otimes x_{1}+  \notag \\
&&+B(x_{1}x_{2}\otimes 1_{H};X_{2},x_{2})X_{2}\otimes x_{2}+  \notag \\
&&+\left[ +1-B(x_{1}x_{2}\otimes 1_{H};1_{A},x_{1}x_{2})-B(x_{1}x_{2}\otimes
1_{H};X_{2},x_{1})+B(x_{1}x_{2}\otimes 1_{H};X_{1},x_{2})\right]
X_{1}X_{2}\otimes 1_{H}  \notag
\end{eqnarray}%
\begin{eqnarray}
B\left( gx_{1}\otimes 1_{H}\right) &=&+\left[ B(x_{1}x_{2}\otimes
1_{H};1_{A},x_{1}x_{2})+B(x_{1}x_{2}\otimes 1_{H};X_{2},x_{1})\right]
1_{A}\otimes x_{1}+  \label{mgx1} \\
&&+B(x_{1}x_{2}\otimes 1_{H};X_{2},x_{2})1_{A}\otimes x_{2}+  \notag \\
&&\left[ B(x_{1}x_{2}\otimes 1_{H};G,x_{2})-B(x_{1}x_{2}\otimes
1_{H};GX_{2},1_{H})\right] G\otimes 1_{H}+  \notag \\
&&+\left[ -1+B(x_{1}x_{2}\otimes 1_{H};1_{A},x_{1}x_{2})+B(x_{1}x_{2}\otimes
1_{H};X_{2},x_{1})\right] X_{1}\otimes 1_{H}+  \notag \\
&&+B(x_{1}x_{2}\otimes 1_{H};X_{2},x_{2})X_{2}\otimes 1_{H}  \notag
\end{eqnarray}

\begin{eqnarray}
B\left( gx_{2}\otimes 1_{H}\right) &=&-B(x_{1}x_{2}\otimes
1_{H};X_{1},x_{1})1_{A}\otimes x_{1}+  \label{mgx2} \\
&&+\left[ B(x_{1}x_{2}\otimes 1_{H};1_{A},x_{1}x_{2})-B(x_{1}x_{2}\otimes
1_{H};X_{1},x_{2})\right] 1_{A}\otimes x_{2}+  \notag \\
&&-B(x_{1}x_{2}\otimes 1_{H};G,x_{1})G\otimes 1_{H}+  \notag \\
&&-B(x_{1}x_{2}\otimes 1_{H};X_{1},x_{1})X_{1}\otimes 1_{H}+  \notag \\
&&\left[ -1+B(x_{1}x_{2}\otimes 1_{H};1_{A},x_{1}x_{2})-B(x_{1}x_{2}\otimes
1_{H};X_{1},x_{2})\right] X_{2}\otimes 1  \notag
\end{eqnarray}

\bigskip

\subsubsection{Binomial equalities containing only terms of $B\left(
x_{1}x_{2}\otimes 1_{H}\right) $ \label{BEX}}

\begin{gather*}
1-B(x_{1}x_{2}\otimes 1_{H};1_{A},x_{1}x_{2})+B(x_{1}x_{2}\otimes
1_{H};X_{1},x_{2})-B(x_{1}x_{2}\otimes 1_{H};X_{2},x_{1})=0 \\
\left( \ref{G,x1x2, GF4,1H}\right)
\end{gather*}%
\begin{equation*}
\begin{array}{c}
2\alpha B(x_{1}x_{2}\otimes 1_{H};G,x_{1})+ \\
+\gamma _{1}B(x_{1}x_{2}\otimes 1_{H};X_{1},x_{1})+\gamma
_{2}B(x_{1}x_{2}\otimes 1_{H};X_{2},x_{1})%
\end{array}%
=0.\left( \ref{G,x1x2, GF1,x1}\right)
\end{equation*}%
In view of $\left( \ref{G,x1x2, GF4,1H}\right) ,$%
\begin{equation*}
\begin{array}{c}
2\alpha B(x_{1}x_{2}\otimes 1_{H};G,x_{2})+ \\
+\gamma _{1}\left[ -1+B(x_{1}x_{2}\otimes
1_{H};1_{A},x_{1}x_{2})+B(x_{1}x_{2}\otimes 1_{H};X_{2},x_{1})\right]
+\gamma _{2}B(x_{1}x_{2}\otimes 1_{H};X_{2},x_{2})%
\end{array}%
=0.\left( \ref{G,gx1, GF1,1H}\right)
\end{equation*}%
becomes%
\begin{eqnarray*}
\begin{array}{c}
2\alpha B(x_{1}x_{2}\otimes 1_{H};G,x_{1})+ \\
+\gamma _{1}B(x_{1}x_{2}\otimes 1_{H};X_{1},x_{2})+\gamma
_{2}B(x_{1}x_{2}\otimes 1_{H};X_{2},x_{1})%
\end{array}
&=&0.\left( \ref{G,x1x2, GF1,x1}\right) \\
&&\text{this is }\left( \ref{G,x1x2, GF1,x1}\right)
\end{eqnarray*}%
In view of $\left( \ref{G,x1x2, GF4,1H}\right) $%
\begin{equation*}
\begin{array}{c}
2\alpha \left[ -B(x_{1}x_{2}\otimes 1_{H};G,x_{1})\right] +\left( \ref%
{G,gx2, GF1,1H}\right) \\
-\gamma _{1}B(x_{1}x_{2}\otimes 1_{H};X_{1},x_{1})+\gamma _{2}\left[
-1+B(x_{1}x_{2}\otimes 1_{H};1_{A},x_{1}x_{2})-B(x_{1}x_{2}\otimes
1_{H};X_{1},x_{2})\right]%
\end{array}%
\end{equation*}%
becomes%
\begin{eqnarray*}
\begin{array}{c}
2\alpha \left[ -B(x_{1}x_{2}\otimes 1_{H};G,x_{1})\right] + \\
-\gamma _{1}B(x_{1}x_{2}\otimes 1_{H};X_{1},x_{1})-\gamma
_{2}B(x_{1}x_{2}\otimes 1_{H};X_{2},x_{1})%
\end{array}
&=&0\left( \ref{G,gx2, GF1,1H}\right) \\
&&\text{which is }\left( \ref{G,x1x2, GF1,x1}\right) .
\end{eqnarray*}%
\begin{equation*}
\begin{array}{c}
2\alpha B(x_{1}x_{2}\otimes 1_{H};G,x_{2})+ \\
+\gamma _{1}B(x_{1}x_{2}\otimes 1_{H};X_{1},x_{2})+\gamma
_{2}B(x_{1}x_{2}\otimes 1_{H};X_{2},x_{2})%
\end{array}%
=0.\left( \ref{G,x1x2, GF1,x2}\right)
\end{equation*}%
\begin{eqnarray*}
&&2\beta _{1}B(x_{1}x_{2}\otimes 1_{H};X_{1},x_{1})+\lambda
B(x_{1}x_{2}\otimes 1_{H};X_{2},x_{1}) \\
+\gamma _{1}B(x_{1}x_{2}\otimes 1_{H};G,x_{1}) &=&0\left( \ref%
{X1,x1x2,X1F11,x1}\right)
\end{eqnarray*}%
\begin{eqnarray*}
&&2\beta _{1}B(x_{1}x_{2}\otimes 1_{H};X_{1},x_{2})+\lambda
B(x_{1}x_{2}\otimes 1_{H};X_{2},x_{2}) \\
+\gamma _{1}B(x_{1}x_{2}\otimes 1_{H};G,x_{2}) &=&0.\left( \ref%
{X1,x1x2,X1F11,x2}\right)
\end{eqnarray*}%
\begin{gather*}
2\beta _{2}B(x_{1}x_{2}\otimes 1_{H};X_{2},x_{1})+\gamma
_{2}B(x_{1}x_{2}\otimes 1_{H};G,x_{1})\left( \ref{X2,x1x2,X2F11,x1}\right) \\
+\lambda B(x_{1}x_{2}\otimes 1_{H};X_{1},x_{1})=0.
\end{gather*}

\begin{gather*}
2\beta _{2}B(x_{1}x_{2}\otimes 1_{H};X_{2},x_{2})+\gamma
_{2}B(x_{1}x_{2}\otimes 1_{H};G,x_{2})\left( \ref{X2,x1x2,X2F11,x2}\right) \\
+\lambda B(x_{1}x_{2}\otimes 1_{H};X_{1},x_{2})=0.
\end{gather*}

In view of $\left( \ref{G,x1x2, GF4,1H}\right) $%
\begin{equation*}
1-B(x_{1}x_{2}\otimes 1_{H};1_{A},x_{1}x_{2})+B(x_{1}x_{2}\otimes
1_{H};X_{1},x_{2})-B(x_{1}x_{2}\otimes 1_{H};X_{2},x_{1})=0
\end{equation*}%
\begin{gather*}
-2\beta _{1}B(x_{1}x_{2}\otimes 1_{H};X_{1},x_{1})+\lambda \left[
-1+B(x_{1}x_{2}\otimes 1_{H};1_{A},x_{1}x_{2})-B(x_{1}x_{2}\otimes
1_{H};X_{1},x_{2})\right] \left( \ref{X1,gx2,X1F11,1H}\right) \\
+\gamma _{1}\left[ -B(x_{1}x_{2}\otimes 1_{H};G,x_{1})\right] =0.
\end{gather*}%
becomes%
\begin{gather*}
-2\beta _{1}B(x_{1}x_{2}\otimes 1_{H};X_{1},x_{1})-\lambda
B(x_{1}x_{2}\otimes 1_{H};X_{2},x_{1})\left( \ref{X1,gx2,X1F11,1H}\right) \\
+\gamma _{1}\left[ -B(x_{1}x_{2}\otimes 1_{H};G,x_{1})\right] =0. \\
\text{which is }\left( \ref{X1,x1x2,X1F11,x1}\right)
\end{gather*}%
In view of $\left( \ref{G,x1x2, GF4,1H}\right) $%
\begin{equation*}
1-B(x_{1}x_{2}\otimes 1_{H};1_{A},x_{1}x_{2})+B(x_{1}x_{2}\otimes
1_{H};X_{1},x_{2})-B(x_{1}x_{2}\otimes 1_{H};X_{2},x_{1})=0
\end{equation*}%
\begin{gather*}
2\beta _{1}\left[ -1+B(x_{1}x_{2}\otimes
1_{H};1_{A},x_{1}x_{2})+B(x_{1}x_{2}\otimes 1_{H};X_{2},x_{1})\right] \left( %
\ref{X1,gx1,X1F11,1H}\right) \\
+\lambda B(x_{1}x_{2}\otimes 1_{H};X_{2},x_{2})+\gamma
_{1}B(x_{1}x_{2}\otimes 1_{H};G,x_{2})=0
\end{gather*}%
becomes%
\begin{gather*}
2\beta _{1}B(x_{1}x_{2}\otimes 1_{H};X_{1},x_{2})+\lambda
B(x_{1}x_{2}\otimes 1_{H};X_{2},x_{2})+\gamma _{1}B(x_{1}x_{2}\otimes
1_{H};G,x_{2})=0\left( \ref{X1,gx1,X1F11,1H}\right) \text{ } \\
\text{which is }\left( \ref{X1,x1x2,X1F11,x2}\right) .
\end{gather*}%
In view of $\left( \ref{G,x1x2, GF4,1H}\right) $%
\begin{equation*}
1-B(x_{1}x_{2}\otimes 1_{H};1_{A},x_{1}x_{2})+B(x_{1}x_{2}\otimes
1_{H};X_{1},x_{2})-B(x_{1}x_{2}\otimes 1_{H};X_{2},x_{1})=0
\end{equation*}%
\begin{gather*}
2\beta _{2}\left[ -1+B(x_{1}x_{2}\otimes
1_{H};1_{A},x_{1}x_{2})-B(x_{1}x_{2}\otimes 1_{H};X_{1},x_{2})\right] \left( %
\ref{X2,gx2,X2F11,1H}\right) \\
+\gamma _{2}\left[ -B(x_{1}x_{2}\otimes 1_{H};G,x_{1})\right] -\lambda
B(x_{1}x_{2}\otimes 1_{H};X_{1},x_{1})=0.
\end{gather*}%
becomes%
\begin{gather*}
-2\beta _{2}B(x_{1}x_{2}\otimes 1_{H};X_{2},x_{1})\left( \ref%
{X2,gx2,X2F11,1H}\right) \\
+\gamma _{2}\left[ -B(x_{1}x_{2}\otimes 1_{H};G,x_{1})\right] -\lambda
B(x_{1}x_{2}\otimes 1_{H};X_{1},x_{1})=0. \\
\text{which is }\left( \ref{X2,x1x2,X2F11,x1}\right)
\end{gather*}

In view of $\left( \ref{G,x1x2, GF4,1H}\right) $%
\begin{equation*}
1-B(x_{1}x_{2}\otimes 1_{H};1_{A},x_{1}x_{2})+B(x_{1}x_{2}\otimes
1_{H};X_{1},x_{2})-B(x_{1}x_{2}\otimes 1_{H};X_{2},x_{1})=0
\end{equation*}%
\begin{gather*}
2\beta _{2}B(x_{1}x_{2}\otimes 1_{H};X_{2},x_{2})+\gamma
_{2}B(x_{1}x_{2}\otimes 1_{H};G,x_{2})\left( \ref{X2,gx1,X2F11,1H}\right) \\
+\lambda \left[ -1+B(x_{1}x_{2}\otimes
1_{H};1_{A},x_{1}x_{2})+B(x_{1}x_{2}\otimes 1_{H};X_{2},x_{1})\right] =0.
\end{gather*}%
becomes%
\begin{gather*}
2\beta _{2}B(x_{1}x_{2}\otimes 1_{H};X_{2},x_{2})+\gamma
_{2}B(x_{1}x_{2}\otimes 1_{H};G,x_{2})\left( \ref{X2,gx1,X2F11,1H}\right) \\
+\lambda B(x_{1}x_{2}\otimes 1_{H};X_{2},x_{1})=0. \\
\text{which is }\left( \ref{X2,x1x2,X2F11,x2}\right)
\end{gather*}

\begin{gather*}
1-B(x_{1}x_{2}\otimes 1_{H};1_{A},x_{1}x_{2})+B(x_{1}x_{2}\otimes
1_{H};X_{1},x_{2})-B(x_{1}x_{2}\otimes 1_{H};X_{2},x_{1})=0 \\
\left( \ref{G,x1x2, GF4,1H}\right)
\end{gather*}%
\begin{equation*}
\begin{array}{c}
2\alpha B(x_{1}x_{2}\otimes 1_{H};G,x_{1})+ \\
+\gamma _{1}B(x_{1}x_{2}\otimes 1_{H};X_{1},x_{1})+\gamma
_{2}B(x_{1}x_{2}\otimes 1_{H};X_{2},x_{1})%
\end{array}%
=0.\left( \ref{G,x1x2, GF1,x1}\right)
\end{equation*}%
\begin{equation*}
\begin{array}{c}
2\alpha B(x_{1}x_{2}\otimes 1_{H};G,x_{2})+ \\
+\gamma _{1}B(x_{1}x_{2}\otimes 1_{H};X_{1},x_{2})+\gamma
_{2}B(x_{1}x_{2}\otimes 1_{H};X_{2},x_{2})%
\end{array}%
=0.\left( \ref{G,x1x2, GF1,x2}\right)
\end{equation*}%
\begin{eqnarray*}
&&2\beta _{1}B(x_{1}x_{2}\otimes 1_{H};X_{1},x_{1})+\lambda
B(x_{1}x_{2}\otimes 1_{H};X_{2},x_{1}) \\
+\gamma _{1}B(x_{1}x_{2}\otimes 1_{H};G,x_{1}) &=&0\left( \ref%
{X1,x1x2,X1F11,x1}\right)
\end{eqnarray*}%
\begin{eqnarray*}
&&2\beta _{1}B(x_{1}x_{2}\otimes 1_{H};X_{1},x_{2})+\lambda
B(x_{1}x_{2}\otimes 1_{H};X_{2},x_{2}) \\
+\gamma _{1}B(x_{1}x_{2}\otimes 1_{H};G,x_{2}) &=&0.\left( \ref%
{X1,x1x2,X1F11,x2}\right)
\end{eqnarray*}%
\begin{gather*}
2\beta _{2}B(x_{1}x_{2}\otimes 1_{H};X_{2},x_{1})+\gamma
_{2}B(x_{1}x_{2}\otimes 1_{H};G,x_{1})\left( \ref{X2,x1x2,X2F11,x1}\right) \\
+\lambda B(x_{1}x_{2}\otimes 1_{H};X_{1},x_{1})=0.
\end{gather*}

\begin{gather*}
2\beta _{2}B(x_{1}x_{2}\otimes 1_{H};X_{2},x_{2})+\gamma
_{2}B(x_{1}x_{2}\otimes 1_{H};G,x_{2})\left( \ref{X2,x1x2,X2F11,x2}\right) \\
+\lambda B(x_{1}x_{2}\otimes 1_{H};X_{1},x_{2})=0.
\end{gather*}

Now we set%
\begin{eqnarray}
A &=&B(x_{1}x_{2}\otimes 1_{H};1_{A},x_{1}x_{2})  \label{variables} \\
B &=&B(x_{1}x_{2}\otimes 1_{H};G,x_{1})  \notag \\
C &=&B(x_{1}x_{2}\otimes 1_{H};G,x_{2})  \notag \\
D &=&B(x_{1}x_{2}\otimes 1_{H};X_{1},x_{1})  \notag \\
E &=&B(x_{1}x_{2}\otimes 1_{H};X_{1},x_{2})  \notag \\
F &=&B(x_{1}x_{2}\otimes 1_{H};X_{2},x_{1})  \notag \\
G &=&B(x_{1}x_{2}\otimes 1_{H};X_{2},x_{2})  \notag
\end{eqnarray}%
and rewrite all our equation using this relabeling.%
\begin{equation*}
-A+E-F=-1\left( \ref{G,x1x2, GF4,1H}\right)
\end{equation*}%
\begin{equation*}
2\alpha B+\gamma _{1}D+\gamma _{2}F=0\left( \ref{G,x1x2, GF1,x1}\right)
\end{equation*}%
\begin{equation*}
2\alpha C+\gamma _{1}E+\gamma _{2}G=0\left( \ref{G,x1x2, GF1,x2}\right)
\end{equation*}%
\begin{equation*}
\gamma _{1}B+2\beta _{1}D+\lambda F=0\left( \ref{X1,x1x2,X1F11,x1}\right)
\end{equation*}%
\begin{equation*}
\gamma _{1}C+2\beta _{1}E+\lambda G=0\left( \ref{X1,x1x2,X1F11,x2}\right)
\end{equation*}%
\begin{equation*}
\gamma _{2}B+\lambda D+2\beta _{2}F=0\left( \ref{X2,x1x2,X2F11,x1}\right)
\end{equation*}

\begin{equation*}
\gamma _{2}C+\lambda E+2\beta _{2}G=0\left( \ref{X2,x1x2,X2F11,x2}\right)
\end{equation*}%
We write the matrix of the coefficients:

\begin{equation*}
\begin{array}{ccccccc}
-1 & 0 & 0 & 0 & 1 & -1 & 0 \\
0 & 2\alpha & 0 & \gamma _{1} & 0 & \gamma _{2} & 0 \\
0 & 0 & 2\alpha & 0 & \gamma _{1} & 0 & \gamma _{2} \\
0 & \gamma _{1} & 0 & 2\beta _{1} & 0 & \lambda & 0 \\
0 & 0 & \gamma _{1} & 0 & 2\beta _{1} & 0 & \lambda \\
0 & \gamma _{2} & 0 & \lambda & 0 & 2\beta _{2} & 0 \\
0 & 0 & \gamma _{2} & 0 & \lambda & 0 & 2\beta _{2}%
\end{array}%
\end{equation*}%
and now we find an equivalent matrix by transposition of rows.

\begin{equation}
\begin{array}{ccccccc}
0 & 2\alpha & \gamma _{1} & \gamma _{2} & 0 & 0 & 0 \\
0 & \gamma _{1} & 2\beta _{1} & \lambda & 0 & 0 & 0 \\
0 & \gamma _{2} & \lambda & 2\beta _{2} & 0 & 0 & 0 \\
0 & 0 & 0 & 0 & 2\alpha & \gamma _{1} & \gamma _{2} \\
0 & 0 & 0 & 0 & \gamma _{1} & 2\beta _{1} & \lambda \\
0 & 0 & 0 & 0 & \gamma _{2} & \lambda & 2\beta _{2} \\
-1 & 0 & 0 & -1 & 1 & 0 & 0%
\end{array}
\label{form matrix}
\end{equation}%
We consider the matrix%
\begin{equation*}
\begin{array}{ccc}
2\alpha & \gamma _{1} & \gamma _{2} \\
\gamma _{1} & 2\beta _{1} & \lambda \\
\gamma _{2} & \lambda & 2\beta _{2}%
\end{array}%
\end{equation*}%
with determinant determinant:%
\begin{equation}
\Delta =-2\alpha \lambda ^{2}+2\lambda \gamma _{1}\gamma _{2}-2\beta
_{2}\gamma _{1}^{2}-2\beta _{1}\gamma _{2}^{2}+8\alpha \beta _{1}\allowbreak
\beta _{2}  \label{form;Delta}
\end{equation}

Then the determinant of our matrix is
\begin{equation*}
-\Delta ^{2}.
\end{equation*}%
We consider the system $%
\begin{array}{ccccccc}
0 & 2\alpha & \gamma _{1} & \gamma _{2} & 0 & 0 & 0 \\
0 & \gamma _{1} & 2\beta _{1} & \lambda & 0 & 0 & 0 \\
0 & \gamma _{2} & \lambda & 2\beta _{2} & 0 & 0 & 0 \\
0 & 0 & 0 & 0 & 2\alpha & \gamma _{1} & \gamma _{2} \\
0 & 0 & 0 & 0 & \gamma _{1} & 2\beta _{1} & \lambda \\
0 & 0 & 0 & 0 & \gamma _{2} & \lambda & 2\beta _{2} \\
-1 & 0 & 0 & -1 & 1 & 0 & 0%
\end{array}%
\begin{array}{c}
A \\
B \\
C \\
D \\
E \\
F \\
G%
\end{array}%
=%
\begin{array}{c}
0 \\
0 \\
0 \\
0 \\
0 \\
0 \\
-1%
\end{array}%
$

\begin{itemize}
\item If $\Delta \neq 0$ then we have exactly one solution, namely $%
\begin{array}{c}
1 \\
0 \\
0 \\
0 \\
0 \\
0 \\
0%
\end{array}%
$

\item If $\Delta =0$ then we have infinite solutions.
\end{itemize}

\section{The separability result}

\begin{theorem}
\label{Theo sep}Let $A=Cl(\alpha ,\beta _{i},\gamma _{i},\lambda )$ be the
Clifford algebra considered as a comodulo algebra over $E\left( 2\right) $
the Radford $8$ dimensional Hopf algebra as in $\left( \ref{CLIFF}\right) $.
Assume that $\mathrm{char}\left( \Bbbk \right) \neq 2$ and that all defining
constants of $A$ are nonzero.

Then the cowreath $(A\otimes H^{op},H)$ is separable with respect to any
bilinear form satisfying all monomial equalities in $\left( \ref{LME}\right)
$, all binomial equalities in $\left( \ref{BEW}\right) $ where the six
elements satisfy $\left( \ref{got1}\right) ,\left( \ref{x1}\right) ,\left( %
\ref{x2}\right) ,\left( \ref{x1x2}\right) ,\left( \ref{gx1}\right) ,\left( %
\ref{gx2}\right) $ and $\left( \ref{gx1x2}\right) $ i.e. they satisfy $%
\left( \ref{mbg}\right) \left( \ref{mbx1}\right) \left( \ref{mbx2}\right)
\left( \ref{mx1x2}\right) \left( \ref{mgx1}\right) \left( \ref{mgx2}\right)
\left( \ref{mbgx1x2}\right) .$

Then we consider the determinant $\Delta $ defined in $\left( \ref%
{form;Delta}\right) .$Moreover we have two cases.

\begin{enumerate}
\item Whenever $\Delta \neq 0$ then the following further condition should
be satisfied%
\begin{eqnarray*}
B(x_{1}x_{2}\otimes 1_{H};1_{A},x_{1}x_{2}) &=&1 \\
B(x_{1}x_{2}\otimes 1_{H};G,x_{1}) &=&0 \\
B(x_{1}x_{2}\otimes 1_{H};G,x_{2}) &=&0 \\
B(x_{1}x_{2}\otimes 1_{H};X_{1},x_{1}) &=&0 \\
B(x_{1}x_{2}\otimes 1_{H};X_{1},x_{2}) &=&0 \\
B(x_{1}x_{2}\otimes 1_{H};X_{2},x_{1}) &=&0 \\
B(x_{1}x_{2}\otimes 1_{H};X_{2},x_{2}) &=&0
\end{eqnarray*}

\item Whenever $\Delta =0$ we have infinite solutions..
\end{enumerate}
\end{theorem}

\begin{proof}
In view of \cite[Proposition 7.4]{BCT1}, we have to find a bilinear form
\begin{equation*}
B:H\otimes H\rightarrow A\otimes H^{op}
\end{equation*}%
which satisfies Casimir condition, morphism condition and also the
normalized condition. In paragraph $\left( \ref{MCNC}\right) $ we proved
that if $B$ satisfies $\left( \ref{eq.h}\right) $ of the morphism conditions
and $\left( \ref{1ot1}\right) $ then it satisfies normalized condition.

In part $\left( \ref{PART CASIMIR}\right) $ we consider our eight elements $%
g\otimes 1_{H},x_{1}\otimes 1_{H},x_{1}\otimes 1_{H},x_{2}\otimes
1_{H},x_{1}x_{2}\otimes 1_{H},gx_{1}\otimes 1_{H},gx_{2}\otimes
1_{H},gx_{1}x_{2}\otimes 1_{H}$ and we will assume that Casimir condition \
holds for $B$ with respect to them and deduce the particular form that we
will get in that case namely $\left( \ref{1ot1}\right) ,\left( \ref{got1}%
\right) ,\left( \ref{x1}\right) ,\left( \ref{x2}\right) ,\left( \ref{x1x2}%
\right) ,\left( \ref{gx1}\right) ,\left( \ref{gx2}\right) $ and $\left( \ref%
{gx1x2}\right) .$Then using $\left( \ref{eq.a}\right) $ and $\left( \ref%
{eq.h}\right) $ we could prove that Casimir condition is automatically
satisfied in the case of all other elements. then in section we proved that
in view of our choice to write all $B\left( h\otimes h^{\prime }\right) $ by
using $B\left( h\otimes 1_{H}\right) $ where $h\otimes 1_{H}$ runs through
the seven elements, this condition $\left( \ref{eq.h}\right) $ will be
automatically satisfied. Thus we were left with Equation $\left( \ref{eq.a}%
\right) .$ After getting in $\left( \ref{LAE}\right) $ all the equalities
that should be satisfied in order that this condition hold we choose to
restrict to the non degenerate case i.e. that all defining constants of $A$
are nonzero. Then we were reduced that our seven elements should satisfy all
monomial equalities in $\left( \ref{LME}\right) $ and all binomial
equalities in $\left( \ref{BEW}\right) $ i.e. they satisfy $\left( \ref{1ot1}%
\right) ,\left( \ref{mbg}\right) ,\left( \ref{mbx1}\right) ,\left( \ref{mbx2}%
\right) ,\left( \ref{mx1x2}\right) ,\left( \ref{mgx1}\right) ,\left( \ref%
{mgx2}\right) $ and $\left( \ref{mbgx1x2}\right) .$Then we studied the
remaining equalities $\left( \ref{BEX}\right) .$The study of these
equalities led us to consider a system of linear equalities in the variables
listed in $\left( \ref{variables}\right) .$Then we considered the matrix of
the system we got. Its determinant $\Delta $ gave us the two cases one has
to consider to find the solution of the system.
\end{proof}

\end{document}